%% file: matrixDecom_EN.tex
\documentclass[twoside,11pt,dvipsnames, openleft]{book}
\usepackage{jmlr2e}
\usepackage{bm}
\usepackage{cancel}

\renewcommand{\arraystretch}{1.15}

\usepackage[english]{babel}

\usepackage{imakeidx}
\makeindex[columns=2, title=Alphabetical Index]

\usepackage{stackengine}

\usepackage{dsfont}

\usepackage{extarrows}

\usepackage[framemethod=tikz]{mdframed}

\usepackage{amsmath}
\usepackage{arydshln}

\numberwithin{equation}{chapter}

\usepackage{blkarray}

\usepackage[absolute,overlay]{textpos}

\usepackage[final]{pdfpages}

\usepackage[noframe]{showframe}
\usepackage{framed}
\usepackage{lipsum}

\newenvironment{svgraybox}{%
	\MakeFramed{\advance\hsize-\width \FrameRestore\FrameRestore}}%
{\endMakeFramed}

\usepackage{xcolor}
\newcommand{\cleft}[2][.]{%
	\begingroup\colorlet{savedleftcolor}{.}%
	\color{#1}\left#2\color{savedleftcolor}%
}
\newcommand{\cright}[2][.]{%
	\color{#1}\right#2\endgroup
}

\usepackage{tcolorbox}
\usepackage{graphicx}             
\usepackage[hang]{subfigure}

\usepackage{algorithm}
\usepackage{algpseudocode}
\usepackage{graphics}
\usepackage{epsfig}

\usepackage{blkarray}
\usepackage{listings}

\usepackage[Bjornstrup]{fncychap}

\newcommand{\colortitlechap}{\color[RGB]{60,113,183}} 
\newcommand{\colornumberchap}{\color[RGB]{60,113,183}} 
\newcommand{\colorbackchap}{\colorbox[RGB]{200,200,200}} 

\makeatletter
\renewcommand{\DOCH}{%
		\settowidth{\py}{\CNoV\thechapter}
		\addtolength{\py}{-10pt}%
		\fboxsep=0pt%
		\colorbackchap{\rule{0pt}{40pt}\parbox[b]{\textwidth}{\hfill}}%
		\kern-\py\raise20pt%
		\hbox{\colornumberchap\CNoV\thechapter}\\%
	}

\renewcommand{\DOTI}[1]{%
		\nointerlineskip\raggedright%
		\fboxsep=\myhi%
		\vskip-1ex%
		\colorbackchap{\parbox[t]{\mylen}{\CTV\FmTi{\colortitlechap#1}}}\par\nobreak%
		\vskip 40\p@%
	}

\renewcommand{\DOTIS}[1]{%
		\fboxsep=0pt%
		\colorbackchap{\rule{0pt}{40pt}\parbox[b]{\textwidth}{\hfill}}\\%
		\nointerlineskip\raggedright%
		\fboxsep=\myhi%
		\colorbackchap{\parbox[t]{\mylen}{\CTV\FmTi{\colortitlechap#1}}}\par\nobreak%
		\vskip 40\p@%
	}
\makeatother

\usepackage{fancyhdr, blindtext}

\newcommand{\changefonts}{%
	\fontsize{9}{11}\selectfont
}

\fancypagestyle{plain}{%

	\fancyhf{}  
}

\usepackage[colorlinks]{hyperref}
\usepackage{hyperref}
\hypersetup{
    colorlinks=true, 
    linktoc=all,     
    linkcolor=mydarkblue,  
    anchorcolor=blue,
    citecolor=mydarkgreen,
}
\usepackage[hyperpageref]{backref} 

\usepackage{booktabs}  

\usepackage{bookmark}
\bookmarksetup{
	numbered,
	open
}

\usepackage{tikz}
\usetikzlibrary{mindmap,trees,backgrounds}
\usetikzlibrary{arrows.meta}
\usetikzlibrary{patterns}

\usetikzlibrary{decorations.text}

\usetikzlibrary{calc,backgrounds}

\newcommand*\circled[1]{\tikz[baseline=(char.base)]{
		\node[shape=circle,draw,inner sep=2pt] (char) {#1};}}
	
\usepackage{changepage}                 

\usepackage{hhline}
\newcommand*{\bmboxtimes}{\multicolumn{1}{|c|}{\bm{\boxtimes}}}
\newcommand*{\colorbmboxtimes}{\multicolumn{1}{|c|}{\textcolor{mylightbluetext}{\bm{\boxtimes}}}}

\usetikzlibrary{positioning,shapes,shadows,arrows}
\tikzstyle{condition}=[rectangle, draw=black, rounded corners, fill=colorqr, drop shadow,
text centered, anchor=north, text=black, text width=3cm]
\tikzstyle{abstract}=[rectangle, draw=black, rounded corners, fill=blue!30, drop shadow,
text centered, anchor=north, text=black, text width=3cm]
\tikzstyle{comment}=[rectangle, draw=black, rounded corners, fill=color1, drop shadow,
text centered, anchor=north, text=black, text width=3cm]
\tikzstyle{myarrow}=[->, >=open triangle 90, thick]
\tikzstyle{line}=[-, thick]

\usepackage{tikz-3dplot}
\tikzset{>=latex}
\usetikzlibrary{matrix}

\newcommand{\slicedist}{0.5}
\newcommand{\slicedistsmall}{0.008}

\newcommand{\Depth}{2}
\newcommand{\Height}{2}
\newcommand{\Width}{2}
\newcommand{\xx}{1}
\newcommand{\yy}{1}
\newcommand{\zz}{1}
\newcommand{\xd}{4.5}
\newcommand{\yd}{-0.3}
\newcommand{\ye}{1.1}
\newcommand{\zd}{1}

\usepackage{sidecap}
\sidecaptionvpos{figure}{t}

\usepackage{verbatimbox}

\usepackage[mathscr]{euscript}

\usepackage{stmaryrd}

\newlength{\offsetpage}
\newenvironment{widepage}{\begin{adjustwidth}{-\offsetpage/2}{-\offsetpage/2}
		\addtolength{\textwidth}{\offsetpage}}%
	{\end{adjustwidth}}
\newenvironment{widepage1}{\begin{adjustwidth}{-0.9cm}{-\offsetpage/2}
		\addtolength{\textwidth}{\offsetpage}}%
	{\end{adjustwidth}}

\usepackage{minitoc}   
\let\cleardoublepage\clearpage
\usepackage[titletoc]{appendix}

\usepackage[labelfont=bf]{caption}
\newcommand\myhrulefill[1]{\leavevmode\leaders\hrule height#1\hfill\kern0pt}
\DeclareCaptionFormat{myformat}{
	{\color[RGB]{0,0,0}\myhrulefill{0.12em}}
	\\#1#2#3
}
\usepackage{afterpage}

\makeatletter \renewcommand*\cleardoublepage{
\clearpage
\if@twoside   
	\ifodd\c@page 
		\hbox{}\newpage
		\if@twocolumn\hbox{}   
			\newpage
		\fi
	\fi
\fi
} \makeatother
\let\originalpart=\part
\def\part#1{\cleardoublepage\clearpage \pagecolor{\partcolor} \originalpart{#1}\nopagecolor }

\usepackage{yfonts,color,lettrine}
\definecolor{caligraphcolor}{HTML}{74AECB}

\usepackage{setspace}

\AtBeginDocument{\setlength{\DefaultFindent}{0.5em}}
\def\algoalign#1{\parbox[t]{\dimexpr\linewidth-\algorithmicindent}{#1}}

\jmlrheading{1}{2021-}{1-48}{4/00}{10/00}{Jun Lu}
\firstpageno{1}

\input{symbols.tex}
\input{symbols_theorems}

\begin{document}
\includepdf[pages=-,scale=1.1,offset=0mm -4mm,]{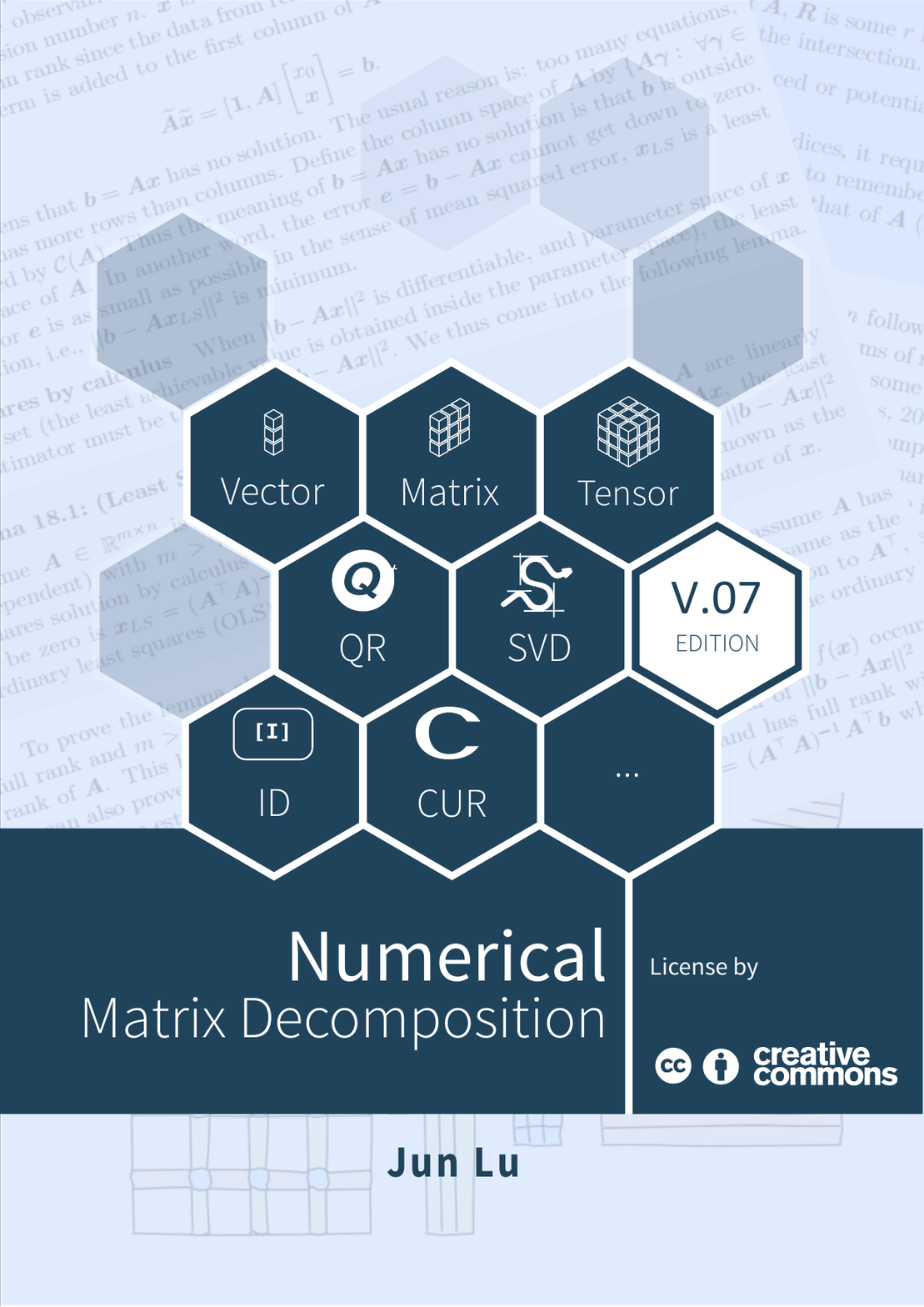} 

\frontmatter
\pagecolor{lowblue}
\newpage
\thispagestyle{empty}
\includegraphics{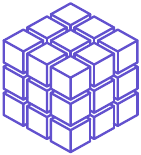}
\begingroup
\hypersetup{
	colorlinks=true, 
	linktoc=all,     
	linkcolor=mydarkblue,  
	anchorcolor=blue,
	citecolor=mydarkgreen,
	urlcolor=mydarkblue,
}
\begin{textblock}{14}(1,3.2)
{\Large \noindent This book gives a comprehensive summary of the important matrix decomposition methods such as QR, ID, CUR, and SVD, as well as their state-of-the-art applications. It's worth noting that the improved programs for most decomposition algorithms were introduced that could potentially reduce their complexity. It is an interesting book for the researchers or engineers in this field.
	
--\href{https://scholar.google.com.hk/citations?user=DpRPZhgAAAAJ&hl=en&oi=ao}{Hua Zhang}, Professor of Environmental Engineering, Tongji University

}
 
\end{textblock}
\endgroup
\clearpage

\clearpage
\nopagecolor 

\thispagestyle{empty}  
\title{Numerical Matrix Decomposition}

\author{
\begin{center}
\name Jun Lu \\ 
\email jun.lu.locky@gmail.com \\
\copyright \,\, 2021$\sim$ \,\, Jun Lu\\
Version 0.7
\end{center}
}

\maketitle

\chapter*{\centering \begin{normalsize}Preface\end{normalsize}}

Matrix decomposition is a fundamental concept in linear algebra that has far-reaching implications in various fields such as engineering, physics, computer science, and data analysis. At its core, matrix decomposition involves breaking down a matrix into a product of simpler matrices, which can be easier to work with or reveal important properties of the original matrix. This book aims to provide a comprehensive and in-depth exploration of matrix decomposition methods, their theoretical foundations, and their practical applications.

In 1954, Alston S. Householder published \textit{Principles of Numerical Analysis}, one of the earliest modern treatments of matrix decomposition approaches that favored a (block) LU decomposition---the factorization of a matrix into the product of lower and upper triangular matrices.
Since then, matrix decomposition has become a core technology in machine learning and data analysis, largely attributed to the development of the backpropagation algorithm for training neural networks and its ability to reduce the dimensionality of data, making it easier for machine learning or optimization algorithms to process.
The primary aim of this book is to provide a self-contained introduction to the concepts and mathematical tools in numerical linear algebra and matrix analysis, enabling a seamless introduction to matrix decomposition techniques and their applications in subsequent sections.
However, we clearly recognize our inability to cover all the useful and interesting results concerning matrix decomposition, given the limited scope to present this discussion, e.g., the separate analysis of Euclidean space, Hermitian space, Hilbert space, and concepts in the complex domain. We refer the reader to the literature in the field of linear algebra for a more detailed introduction to these related fields.

This book is primarily a summary of the purpose and significance of important matrix decomposition methods, such as LU, CUR, QR, and SVD, along with their origins and complexities, which shed light on their modern applications. Most importantly, this book presents improved procedures for some of the calculations involved in decomposition algorithms, potentially reducing their computational complexity. Throughout this work, we maintain a decomposition-centric perspective, introducing the necessary background when needed. While many other linear algebra textbooks treat matrix decomposition methods as secondary, we place a primary focus on these methods, with fundamental concepts serving as essential tools. The only prerequisite for this book is a basic course in linear algebra. Beyond that, we offer a self-contained development with rigorous proofs provided throughout.
The book covers several topics that are treated differently compared to classical numerical linear algebra textbooks. Some examples of these less common issues include:
\begin{itemize}
\item The treatment of numerical methods is comprehensive, e.g., the comparison of computational complexity between Givens rotation and Householder reflector algorithms, as well as a juxtaposition of Bayesian approaches and randomized algorithms.

\item Exploration and utilization of gradient-based methods, such as alternating least squares, alternating direction methods of multipliers, block coordinate descent, and nonnegative matrix factorization, in addition to low-rank approximation via the singular value decomposition (SVD).

\item Establishing connections between different decompositional methods.

\item Applications such as low-rank neural networks, least squares for rank-deficient matrices, image and matrix compression, and low-rank approximation for movie recommender systems.

\item Providing proofs of existence from different perspective, accommodating  different conditions or prerequisites.

\item A unified discussion of matrix decomposition methods from the perspective of biconjugate decomposition.

\item Low-rank approximation using alternating least squares and Newton's method, unified descriptions of nonnegative matrix factorization approaches, and low-rank matrix factorization using special matrix products.
\end{itemize}

This book seeks to serve as a comprehensive and accessible introduction to matrix decomposition, offering readers a bridge between theoretical concepts and practical applications. It is designed for readers with a foundational understanding of linear algebra and aims to achieve the following objectives:
\begin{itemize}
	\item \textit{Explore core principles.} Present the mathematical foundations of matrix decomposition, ranging from basic methods such as LU, Cholesky, and QR decomposition to advanced techniques like SVD, eigenvalue decomposition, and their modern extensions.
	\item \textit{Highlight practical applications.} Demonstrate the relevance of decomposition methods in diverse fields, including optimization, machine learning, neural network compression, and data interpretation.
	\item \textit{Address computational aspects.} Provide insights into the computational complexities, implementation challenges, and cost-reduction strategies associated with various decomposition algorithms.
	\item \textit{Facilitate problem-solving.} Equip readers with the tools to understand and solve problems involving matrices, emphasizing how decomposition can simplify complex tasks and provide deeper insights.
\end{itemize}

\chapter*{\centering \begin{normalsize}Keywords\end{normalsize}}
Existence and computing of matrix decompositions, Complexity, Floating-point operations (flops), Low-rank approximation, Pivot, LU decomposition for nonzero leading principal minors, Data distillation, CR decomposition, CUR/Skeleton decomposition, Interpolative decomposition, Biconjugate decomposition, Coordinate transformation, Hessenberg decomposition, ULV decomposition, URV decomposition, Rank decomposition, Gram-Schmidt process, Householder reflector, Givens rotation, Rank-revealing decomposition, Cholesky decomposition and update/downdate, Eigenvalue problems, Alternating least squares, 
Low-rank Hadamard decomposition, Low-rank Kronecker decomposition, LoRA,
Randomized algorithm, Tensor decomposition, CP decomposition, Tucker decomposition, High-order SVD. 

\vspace{5em}
\noindent
\textit{Acknowledgment: }
We extend our gratitude to Gilbert Strang for raising the question formulated in Corollary~\ref{corollary:invertible-intersection}, reviewing the manuscript, providing valuable ideas and references related to the three factorizations stemming from the elimination steps, and generously sharing the manuscript of \citet{strang2021three}. 
We also appreciate Federico Poloni for his insightful comments on the proof of alternating least squares.
Special thanks are also due to anonymous professors
who have consulted with us and expressed interest in adopting this book as  course material in college settings.
The author would also like to thank Joerg Osterrieder, Christine P. Chai, and Xuanyu Ye for their collaboration  on the Bayesian approach for nonnegative matrix factorization and (intervened) interpolative decomposition, which has illuminated various aspects of the content presented in this book.
The author also wishes to acknowledge Nicolas P. Rougier for open-sourcing the poster design \citep{nicolas2015nn}.


\input{chapter-worldmap.tex}

\newpage
\begingroup
\hypersetup{
linkcolor=structurecolor,
linktoc=page,  
}
\dominitoc
\pdfbookmark{\contentsname}{toc} 
\tableofcontents 
\endgroup

\input{chapter-notation}

\mainmatter
\input{chapter-intro.tex}

\input{chapter-LU.tex}

\input{chapter-cholesky.tex}
\input{chapter-qr.tex} 
\input{chapter-utv.tex}

\input{chapter-cr.tex}
\input{chapter-skeleton.tex}

\input{chapter-id.tex}

\input{chapter-id-bid.tex}

\input{chapter-hessenberg.tex}

\input{chapter-tridiagonal.tex}

\input{chapter-bidiagonal.tex}

\input{chapter-eigenvalue.tex}

\input{chapter-schur.tex}

\input{chapter-spectral.tex}

\input{chapter-svd.tex}

\input{chapter-eigencalc.tex}

\input{chapter-als.tex}

\input{chapter-NMF.tex}
\input{chapter-special.tex}
\input{chapter-biconjugate.tex}

\input{chapter-tensornotation.tex}

\input{chapter-tensorcp.tex}

\input{chapter-tensortucker.tex}

\input{chapter-tensorhosvd.tex}

\input{chapter-tensorttdecom.tex}

\input{chapter-app_proofs.tex}

\input{chapter-app_inequality}

\input{chapter-app_norm}

\input{chapter-app_proj.tex}

\input{chapter-app_pseudo.tex}

\input{chapter-app_condition}


\newpage
\vskip 0.2in
\bibliography{bib}

\clearpage
\printindex
\clearpage

\backmatter
\clearpage
\pagestyle{empty}




\end{document}

%% file: symbols.tex
\newenvironment{sbmatrix}[1]{\def\mysubscript{#1}\mathop\bgroup\begin{bmatrix}}{\end{bmatrix}\egroup_{\textstyle\mathstrut\mysubscript}}

\newcommand{\visavi}{vis-\`a-vis }
\newcommand{\holders}{\text{H{\"o}lder's} }
\newcommand{\cond}{\text{cond} }
\newcommand{\relerror}{\text{RelError} }
\newcommand{\abserror}{\text{AbsError} }
\newcommand{\mybackslash}{\char`\\}

\newcommand{\longeqf}{\scalebox{4}[1]{=}}

\newcommand{\leadto}{\qquad\underrightarrow{ \text{leads to} }\qquad}

\newcommand{\off}{\text{off}}
\newcommand{\gapthree}{\,\,\,}
\newcommand{\gaps}{\,\,\,\,}
\newcommand{\gapforall}{\,\,}
\newcommand{\gap}{\,\,\,\,\,\,\,\,}  
\mathchardef\mhyphen="2D

\newcommand{\exampbar}{\hfill $\square$\par}
\newcommand{\argmax}{\text{arg max}}
\newcommand{\argmin}{\text{arg min}}

\newcommand\mathopmax[1]{\mathop{\max}_{#1}}
\newcommand\mathopmin[1]{\mathop{\min}_{#1}}
\newcommand\mathoplim[1]{\mathop{\lim}_{#1}}

\newcommand\inner[2]{\left\langle#1, #2\right\rangle}
\newcommand\abs[1]{\left\lvert#1\right\rvert}
\newcommand\norm[1]{\left\lVert#1\right\rVert}
\newcommand\normzero[1]{\left\lVert#1\right\rVert_0}
\newcommand\normone[1]{\left\lVert#1\right\rVert_1}
\newcommand\normmone[1]{\left\lVert#1\right\rVert_{m_1}}
\newcommand\normtwo[1]{\left\lVert#1\right\rVert_2}

\newcommand\normf[1]{\left\lVert#1\right\rVert_F}
\newcommand\norminf[1]{\left\lVert#1\right\rVert_{\infty}}
\newcommand\normminf[1]{\left\lVert#1\right\rVert_{m_\infty}}

\newcommand{\trans}[1]{\ensuremath{#1^{ \top}}}
\newcommand{\dotL}{\dot{\bm{L}}}
\newcommand{\dotU}{\dot{\bm{U}}}
\newcommand{\ddotL}{\ddot{\bm{L}}}

\newcommand{\convd}{\breve d}
\newcommand{\concd}{\widehat d}
\newcommand{\cnstd}{\bar d}

\newcommand{\concG}{\widehat G}

\newcommand{\indicator}{\mathds{1}}
\def\1{\bm{1}}

\newcommand{\R}{\mathbb{R}}
\newcommand{\real}{\mathbb{R}}
\newcommand{\complex}{\mathbb{C}}
\newcommand{\symmetric}{\mathbb{S}}
\newcommand{\psd}{\mathbb{S}_{+}}
\newcommand{\pd}{\mathbb{S}_{++}}
\newcommand{\nsd}{\mathbb{S}_{-}}
\newcommand{\nd}{\mathbb{S}_{--}}
\newcommand{\integer}{\mathbb{Z}}
\newcommand{\naturalset}{\mathbb{N}}

\newcommand{\KL}{D_{\mathrm{KL}}}

\definecolor{titlepagecolor}{cmyk}{75,68,67,90}
\definecolor{titlepagecolor2}{rgb}{1.0, 0.08, 0.58}
\definecolor{emerald}{rgb}{0.31, 0.78, 0.47}
\definecolor{deeppink}{HTML}{D14064}
\definecolor{lowpink}{HTML}{ffe6ec}
\newcommand{\partcolor}{gray!65} 
\definecolor{lowblue}{HTML}{E1EBFE}

\let\oldforall\forall
\renewcommand{\forall}{\oldforall\, }

\usepackage{color}   

\newcommand{\mdframecolor}{gray!10}
\newcommand{\mdframehideline}{true}
\newcommand{\mdframehidelineNote}{true}

\definecolor{mylightbluetitle}{RGB}{60,113,183}
\definecolor{mylightbluetext}{rgb}{0,0.08,0.45}
\definecolor{structurecolorblue}{RGB}{60,113,183}
\definecolor{structurecolorgreen}{RGB}{63,145,182}
\colorlet{structurecolor}{structurecolorblue}
\definecolor{structurecolorelegant}{RGB}{60,113,183}
\definecolor{structurecolorlt}{RGB}{31,119,185}

\definecolor{structurecolorHighTheoremBlue}{RGB}{220,227,248}
\definecolor{structurecolorHighTheoremGreen}{RGB}{188,222,231}
\colorlet{structurecolorHighTheorem}{structurecolorHighTheoremBlue}
\newcommand{\mdframecolorTheorem}{gray!35}  

\newcommand{\mdframecolorBox}{gray!15}  

\definecolor{winestain}{rgb}{0.5,0,0}
\definecolor{mydarkblue}{rgb}{0,0.08,0.45}
\definecolor{mydeepblue}{rgb}{0,0.08,0.65}
\definecolor{mydarkred}{rgb}{0.70,0.00,0.00}
\definecolor{mydarkgreen}{rgb}{0.00,0.30,0.00}
\definecolor{mydarkyellow}{RGB}{197,151,13}
\definecolor{mydarkpurple}{RGB}{90,35,140}
\definecolor{mydarkgray}{RGB}{64,64,64}

\definecolor{color0}  {RGB}{174,225,254} 
\definecolor{color1}  {RGB}{220,227,248} 
\definecolor{color2}  {RGB}{28,130,185} 
\definecolor{color3}  {RGB}{255,253,250} 
\definecolor{colormiddleright}  {RGB}{245,253,250} 
\definecolor{colorbottomleft}  {RGB}{255,243,250} 
\definecolor{coloruppermiddle}  {RGB}{255,253,230} 
\definecolor{colormiddleleft}  {RGB}{255,244,237}
\definecolor{colorcr}  {RGB}{249,253,232} 
\definecolor{colorreduction}  {RGB}{255,235,254} 
\definecolor{colorqr}  {RGB}{254,221,199} 
\definecolor{colorbiconjugate}  {RGB}{251,149,161} 
\definecolor{colorsvd}  {RGB}{215,247,235} 
\definecolor{colorupperright}  {RGB}{239,246,251} 
\definecolor{colorspectral}  {RGB}{206,226,243} 
\definecolor{colorbottomright}  {RGB}{220,224,236} 
\definecolor{coloreigenvalue}  {RGB}{197,203,224} 
\definecolor{colorcp} {RGB}{217, 234, 186} 
\definecolor{colorcpborder} {RGB}{233, 243, 216} 
\definecolor{colorupperleft}  {RGB}{235,243,240} 
\definecolor{colorsemidefinite}  {RGB}{217,232,226} 
\definecolor{colormiddle} {RGB}{235, 240,255}
\definecolor{colorlu}  {RGB}{220,227,255} 
\definecolor{colorals}  {RGB}{240,230,255} 
\definecolor{coloralsbkg}  {RGB}{248,243,255} 
\definecolor{canaryyellow}{rgb}{1.0, 0.75, 0.0}
\definecolor{bluepigment}{rgb}{0.0, 0.0, 1.0}
\definecolor{canarypurple}{RGB}{208, 13, 241}
\definecolor{colorGreenOcre}{RGB}{51,102,0} 
\definecolor{colorBlue2}{RGB}{200,207,248}
\definecolor{shadecolor}{gray}{0.75}

\newcommand{\hadaprod}{\circledast}
\newcommand{\khatrirao}{\odot}
\newcommand{\kronecker}{\otimes}

\newcommand{\dist}{\mathrm{dist}}
\newcommand{\Cov}{\mathbb{C}\mathrm{ov}}

\newcommand{\Exp}{\mathbb{E}}
\newcommand{\Var}{\mathbb{V}\mathrm{ar}}

\newcommand{\prob}{\mathbb{P}}

\newcommand{\gammadist}{\mathcal{G}}
\newcommand{\inversegammadist}{\mathcal{G}^{-1}}

\newcommand{\normal}{\mathcal{N}}
\newcommand{\truncatednormal}{\mathcal{TN}}
\newcommand{\generaltruncatednormal}{\mathcal{GTN}}
\newcommand{\gtnsng}{\mathcal{GTNSNG}}

\newcommand{\betadist}{\mathrm{Beta}}

\newcommand{\bernoulli}{\mathrm{Bernoulli}}  

\newcommand{\tr}{\mathrm{tr}}

\newcommand{\cspace}{\mathcal{C}}
\newcommand{\nspace}{\mathcal{N}}
\newcommand{\bzero}{\mathbf{0}}
\newcommand{\bone}{\mathbf{1}}
\newcommand{\diag}{\mathrm{diag}}
\newcommand{\defect}{\mathrm{def}}
\newcommand{\adjugate}{\mathrm{adj}}
\newcommand{\rank}{\mathrm{rank}}
\newcommand{\inertia}{\mathrm{in}}

\newcommand{\trace}{\mathrm{tr}}
\newcommand{\nilp}{\mathrm{nilp}}
\newcommand{\spn}{\mathrm{span}}
\newcommand{\nr}{\mathrm{nr}}
\newcommand{\kyfannorm}{Ky Fan $k$-norm }
\newcommand{\kyfanpknorm}{Ky Fan $p$-$k$-norm }

\DeclareMathOperator{\sign}{sign}   

\newcommand{\bOmega}{\boldsymbol\Omega}

\newcommand{\balpha}{\boldsymbol\alpha}
\newcommand{\bbeta}{\boldsymbol\beta}
\newcommand{\bdelta}{\boldsymbol\delta}

\newcommand{\bPhi}{\boldsymbol\Phi}
\newcommand{\bPsi}{\boldsymbol\Psi}
\newcommand{\bpsi}{\boldsymbol\psi}
\newcommand{\bgamma}{\boldsymbol\gamma}
\newcommand{\bepsilon}{\boldsymbol\epsilon}
\newcommand{\bLambda}{\boldsymbol\Lambda}
\newcommand{\boldeta}{\boldsymbol\eta}
\newcommand{\blambda}{\boldsymbol\lambda}
\newcommand{\bphi}{\boldsymbol\phi}

\newcommand{\bmu}{\boldsymbol\mu}

\newcommand{\bSigma}{\boldsymbol\Sigma}

\newcommand{\abovebX}[1]{\overset{\mathit{(#1)}}{\bX}}

\newcommand{\widebarbW}{\overline{\bm{W}}}

\newcommand{\widebarbZ}{\overline{\bm{Z}}}
\newcommand{\widebarba}{\overline{\bm{a}}}
\newcommand{\widebarbb}{\overline{\bm{b}}}

\newcommand{\widebarbw}{\overline{\bm{w}}}

\newcommand{\widebarbz}{\overline{\bm{z}}}

\newcommand{\mathcalC}{\mathcal{C}}

\newcommand{\mathcalI}{\mathcal{I}}
\newcommand{\mathcalJ}{\mathcal{J}}

\newcommand{\mathcalM}{\mathcal{M}}
\newcommand{\mathcalN}{\mathcal{N}}
\newcommand{\mathcalO}{\mathcal{O}}
\newcommand{\mathcalP}{\mathcal{P}}

\newcommand{\mathcalS}{\mathcal{S}}
\newcommand{\mathcalT}{\mathcal{T}}
\newcommand{\mathcalU}{\mathcal{U}}
\newcommand{\mathcalV}{\mathcal{V}}
\newcommand{\mathcalW}{\mathcal{W}}

\newcommand{\widetildebA}{\widetilde{\bm{A}}}
\newcommand{\widetildebB}{\widetilde{\bm{B}}}
\newcommand{\widetildebC}{\widetilde{\bm{C}}}
\newcommand{\widetildebD}{\widetilde{\bm{D}}}
\newcommand{\widetildebE}{\widetilde{\bm{E}}}

\newcommand{\widetildebG}{\widetilde{\bm{G}}}
\newcommand{\widetildebH}{\widetilde{\bm{H}}}

\newcommand{\widetildebJ}{\widetilde{\bm{J}}}

\newcommand{\widetildebM}{\widetilde{\bm{M}}}

\newcommand{\widetildebQ}{\widetilde{\bm{Q}}}
\newcommand{\widetildebR}{\widetilde{\bm{R}}}
\newcommand{\widetildebS}{\widetilde{\bm{S}}}

\newcommand{\widetildebU}{\widetilde{\bm{U}}}
\newcommand{\widetildebV}{\widetilde{\bm{V}}}
\newcommand{\widetildebW}{\widetilde{\bm{W}}}

\newcommand{\widetildebZ}{\widetilde{\bm{Z}}}
\newcommand{\widetildeba}{\widetilde{\bm{a}}}
\newcommand{\widetildebb}{\widetilde{\bm{b}}}
\newcommand{\widetildebc}{\widetilde{\bm{c}}}

\newcommand{\widetildebl}{\widetilde{\bm{l}}}

\newcommand{\widetildebv}{\widetilde{\bm{v}}}
\newcommand{\widetildebw}{\widetilde{\bm{w}}}
\newcommand{\widetildebx}{\widetilde{\bm{x}}}
\newcommand{\widetildeby}{\widetilde{\bm{y}}}
\newcommand{\widetildebz}{\widetilde{\bm{z}}}

\newcommand{\widetildea}{\widetilde{a}}

\newcommand{\widetildey}{\widetilde{y}}
\newcommand{\widetildez}{\widetilde{z}}

\newcommand{\ba}{\bm{a}}
\newcommand{\bA}{\bm{A}}
\newcommand{\bb}{\bm{b}}
\newcommand{\bB}{\bm{B}}
\newcommand{\bc}{\bm{c}}
\newcommand{\bC}{\bm{C}}  
\newcommand{\bd}{\bm{d}}
\newcommand{\bD}{\bm{D}}
\newcommand{\be}{\bm{e}}
\newcommand{\bE}{\bm{E}}
\newcommand{\bff}{\bm{f}}
\newcommand{\bF}{\bm{F}}
\newcommand{\bg}{\bm{g}}
\newcommand{\bG}{\bm{G}}
\newcommand{\bh}{\bm{h}}
\newcommand{\bH}{\bm{H}}

\newcommand{\bI}{\bm{I}}

\newcommand{\bJ}{\bm{J}}

\newcommand{\bK}{\bm{K}}
\newcommand{\bl}{\bm{l}}
\newcommand{\bL}{\bm{L}}
\newcommand{\bmm}{\bm{m}}
\newcommand{\bM}{\bm{M}}
\newcommand{\bn}{\bm{n}}
\newcommand{\bN}{\bm{N}}
\newcommand{\bo}{\bm{o}}

\newcommand{\bp}{\bm{p}}
\newcommand{\bP}{\bm{P}}
\newcommand{\bq}{\bm{q}}
\newcommand{\bQ}{\bm{Q}}
\newcommand{\br}{\bm{r}}
\newcommand{\bR}{\bm{R}}
\newcommand{\bs}{\bm{s}}
\newcommand{\bS}{\bm{S}}
\newcommand{\bt}{\bm{t}}
\newcommand{\bT}{\bm{T}}
\newcommand{\bu}{\bm{u}}
\newcommand{\bU}{\bm{U}}
\newcommand{\bv}{\bm{v}}
\newcommand{\bV}{\bm{V}}
\newcommand{\bw}{\bm{w}}
\newcommand{\bW}{\bm{W}}
\newcommand{\bx}{\bm{x}}
\newcommand{\bX}{\bm{X}}
\newcommand{\by}{\bm{y}}
\newcommand{\bY}{\bm{Y}}
\newcommand{\bz}{\bm{z}}
\newcommand{\bZ}{\bm{Z}}

\newcommand{\elA}{{\mathscr{A}}}

\newcommand{\eG}{\bm{\mathscr{G}}}
\newcommand{\eA}{\bm{\mathscr{A}}}
\newcommand{\eC}{\bm{\mathscr{C}}}
\newcommand{\eB}{\bm{\mathscr{B}}}
\newcommand{\eY}{\bm{\mathscr{Y}}}
\newcommand{\eX}{\bm{\mathscr{X}}}
\newcommand{\eL}{\bm{\mathscr{L}}}
\newcommand{\weX}{\widehat{\bm{\mathscr{X}}}}

\newcommand{\widehatbA}{\widehat{\bm{A}}}
\newcommand{\widehatbU}{\widehat{\bm{U}}}

\newcommand{\brx}{\bm{R}_\gamma^{(x)}}
\newcommand{\bry}{\bm{R}_\gamma^{(y)}} 

\newcommand{\bryt}{\bm{R}_\gamma^{(y)\top}} 
\newcommand{\bxgamma}{\bX_\gamma}
\newcommand{\bygamma}{\bY_\gamma}
\newcommand{\bugamma}{\bU_\gamma}
\newcommand{\bvgamma}{\bV_\gamma}
\newcommand{\bomegagamma}{\bOmega_\gamma}
\newcommand{\svd}{\mathrm{svd}}

\def\vmu{{\bm{\mu}}}
\def\vtheta{{\bm{\theta}}}
\def\va{{\bm{a}}}

\def\vx{{\bm{x}}}

\def\mA{{\bm{A}}}
\def\mB{{\bm{B}}}

\def\mH{{\bm{H}}}

\def\mJ{{\bm{J}}}

\def\mX{{\bm{X}}}

\def\mSigma{{\bm{\Sigma}}}

\DeclareMathAlphabet{\mathsfit}{\encodingdefault}{\sfdefault}{m}{sl}
\SetMathAlphabet{\mathsfit}{bold}{\encodingdefault}{\sfdefault}{bx}{n}

\def\sA{{\mathbb{A}}}
\def\sB{{\mathbb{B}}}

\def\sD{{\mathbb{D}}}
\def\sF{{\mathbb{F}}}
\def\sG{{\mathbb{G}}}

\def\sK{{\mathbb{K}}}

\def\sO{{\mathbb{O}}}

\def\sS{{\mathbb{S}}}

\def\sV{{\mathbb{V}}}

\def\sX{{\mathbb{X}}}
\def\sY{{\mathbb{Y}}}

\def\ra{{\textnormal{a}}}
\def\rb{{\textnormal{b}}}
\def\rc{{\textnormal{c}}}

\def\rp{{\textnormal{p}}}
\def\rq{{\textnormal{q}}}

\def\rx{{\textnormal{x}}}
\def\ry{{\textnormal{y}}}

\def\rva{{\mathbf{a}}}
\def\rvb{{\mathbf{b}}}

\def\rvv{{\mathbf{v}}}

\def\rvy{{\mathbf{y}}}

\def\rmA{{\mathbf{A}}}
\def\rmB{{\mathbf{B}}}

\def\rmX{{\mathbf{X}}}
\def\rmY{{\mathbf{Y}}}

\def\eva{{a}}

\def\gD{{\mathcal{D}}}
\def\gE{{\mathcal{E}}}

%% file: symbols_theorems.tex
\usepackage[noframe]{showframe}
\usepackage{framed}
\usepackage{lipsum}

\newcommand{\roundcornertheorem}{0pt}
\newcommand{\linewidththeorem}{0.1pt}
\newcommand{\frametitlerulewidththeorem}{0.1pt}
\newcommand{\innerbottommargintheorem}{2pt}
\newcommand{\innerleftmargintheorem}{2pt}
\newcommand{\innerrightmargintheorem}{2pt}
\newcommand{\innertopmargintheorem}{2pt}
\newcommand{\outerlinewidththeorem}{1pt}

\newcommand{\textremark}{3pt}

\newenvironment{theoremHigh}[1][]{%
\refstepcounter{theo}%
\ifstrempty{#1}%
{\newcommand{\theoName}{}}
{\newcommand{\theoName}{:~(#1)}}
\mdfsetup{
	backgroundcolor=structurecolorHighTheorem,
	linecolor=structurecolorHighTheorem,
	frametitlerulewidth=\frametitlerulewidththeorem,
	roundcorner=\roundcornertheorem,
	linewidth=\linewidththeorem,
	innerbottommargin=\innerbottommargintheorem,
	innerleftmargin=\innerleftmargintheorem,
	innerrightmargin=\innerrightmargintheorem,
	innertopmargin=\innertopmargintheorem,
	outerlinewidth=\outerlinewidththeorem,
	topline=false,
	innertopmargin=-5pt,
	innerbottommargin=1pt,
	linewidth=0,
	startinnercode=\paragraph{{\strut Theorem~\thetheo\theoName}}
}
\begin{mdframed}[]\relax%
}{\end{mdframed}}
\newenvironment{corollaryHigh}[1][]{%
\refstepcounter{theo}%
\ifstrempty{#1}%
{\newcommand{\theoName}{}}
{\newcommand{\theoName}{:~(#1)}}
\mdfsetup{
	backgroundcolor=structurecolorHighTheorem,
	linecolor=structurecolorHighTheorem,
	frametitlerulewidth=\frametitlerulewidththeorem,
	roundcorner=\roundcornertheorem,
	linewidth=\linewidththeorem,
	innerbottommargin=\innerbottommargintheorem,
	innerleftmargin=\innerleftmargintheorem,
	innerrightmargin=\innerrightmargintheorem,
	innertopmargin=\innertopmargintheorem,
	outerlinewidth=\outerlinewidththeorem,
	topline=false,
	innertopmargin=-5pt,
	innerbottommargin=1pt,
	linewidth=0,
	startinnercode=\paragraph{{\strut Corollary~\thetheo\theoName}}
}
\begin{mdframed}[]\relax%
}{\end{mdframed}}

\newenvironment{theorem}[1][]{%
	\refstepcounter{theo}%
	\ifstrempty{#1}%
	{\newcommand{\theoName}{}}
	{\newcommand{\theoName}{:~(#1)}}
	\mdfsetup{
		backgroundcolor=\mdframecolorTheorem,
		linecolor=\mdframecolorTheorem,
		frametitlerulewidth=\frametitlerulewidththeorem,
		roundcorner=\roundcornertheorem,
		linewidth=\linewidththeorem,
		innerbottommargin=\innerbottommargintheorem,
		innerleftmargin=\innerleftmargintheorem,
		innerrightmargin=\innerrightmargintheorem,
		innertopmargin=\innertopmargintheorem,
		outerlinewidth=\outerlinewidththeorem,
		topline=false,
		innertopmargin=-5pt,
		innerbottommargin=1pt,
		linewidth=0,
		startinnercode=\paragraph{{\strut Theorem~\thetheo\theoName}}
	}
	\begin{mdframed}[]\relax%
	}{\end{mdframed}}
\newenvironment{corollary}[1][]{%
	\refstepcounter{theo}%
\ifstrempty{#1}%
{\newcommand{\theoName}{}}
{\newcommand{\theoName}{:~(#1)}}
\mdfsetup{
	backgroundcolor=\mdframecolorTheorem,
	linecolor=\mdframecolorTheorem,
	frametitlerulewidth=\frametitlerulewidththeorem,
	roundcorner=\roundcornertheorem,
	linewidth=\linewidththeorem,
	innerbottommargin=\innerbottommargintheorem,
	innerleftmargin=\innerleftmargintheorem,
	innerrightmargin=\innerrightmargintheorem,
	innertopmargin=\innertopmargintheorem,
	outerlinewidth=\outerlinewidththeorem,
	topline=false,
	innertopmargin=-5pt,
	innerbottommargin=1pt,
	linewidth=0,
	startinnercode=\paragraph{{\strut Corollary~\thetheo\theoName}}
}
	\begin{mdframed}[]\relax%
	}{\end{mdframed}}
\newenvironment{lemma}[1][]{%
	\refstepcounter{theo}%
\ifstrempty{#1}%
{\newcommand{\theoName}{}}
{\newcommand{\theoName}{:~(#1)}}
\mdfsetup{
	backgroundcolor=\mdframecolorTheorem,
	linecolor=\mdframecolorTheorem,
	frametitlerulewidth=\frametitlerulewidththeorem,
	roundcorner=\roundcornertheorem,
	linewidth=\linewidththeorem,
	innerbottommargin=\innerbottommargintheorem,
	innerleftmargin=\innerleftmargintheorem,
	innerrightmargin=\innerrightmargintheorem,
	innertopmargin=\innertopmargintheorem,
	outerlinewidth=\outerlinewidththeorem,
	topline=false,
	innertopmargin=-5pt,
	innerbottommargin=1pt,
	linewidth=0,
	startinnercode=\paragraph{{\strut Lemma~\thetheo\theoName}}
}
\begin{mdframed}[]\relax%
}{\end{mdframed}}

\newenvironment{proposition}[1][]{%
\refstepcounter{theo}%
\ifstrempty{#1}%
{\newcommand{\theoName}{}}
{\newcommand{\theoName}{:~(#1)}}
\mdfsetup{
	backgroundcolor=\mdframecolorTheorem,
	linecolor=\mdframecolorTheorem,
	frametitlerulewidth=\frametitlerulewidththeorem,
	roundcorner=\roundcornertheorem,
	linewidth=\linewidththeorem,
	innerbottommargin=\innerbottommargintheorem,
	innerleftmargin=\innerleftmargintheorem,
	innerrightmargin=\innerrightmargintheorem,
	innertopmargin=\innertopmargintheorem,
	outerlinewidth=\outerlinewidththeorem,
	topline=false,
	innertopmargin=-5pt,
	innerbottommargin=1pt,
	linewidth=0,
	startinnercode=\paragraph{{\strut Proposition~\thetheo\theoName}}
}
\begin{mdframed}[]\relax%
}{\end{mdframed}}

\usepackage{amsthm}
\newtheoremstyle{normalfontstyle} 
{3pt}                           
{3pt}                           
{\normalfont}                   
{}                              
{\bfseries}                     
{}                             
{ }                             
{}                              

\usepackage{thmtools}
\declaretheoremstyle[
spaceabove=3pt,
spacebelow=3pt,
headfont=\bfseries,
notefont=\bfseries, 
notebraces={(}{)}, 
bodyfont=\normalfont,
postheadspace=1em, 
]{normalfontboldhead}

\theoremstyle{normalfontstyle}
\newcommand{\BlackBox}{\rule{1.5ex}{1.5ex}}  
\renewenvironment{proof}{\par\noindent{\bf Proof\ }}{\hfill\BlackBox\\[2mm]}

\declaretheorem[style=normalfontboldhead, name=Definition, numberlike=theo]{definitionT}
\newmdenv[skipabove=7pt,
skipbelow=7pt,
rightline=false,
leftline=true,
topline=false,
bottomline=false,
linecolor=mydarkblue,
innerleftmargin=5pt,
innerrightmargin=5pt,
innertopmargin=0pt,
leftmargin=2cm,
rightmargin=0cm,
linewidth=4pt,
innerbottommargin=0pt]{dBox}
\newenvironment{definition}{\begin{dBox}\begin{definitionT}}{\end{definitionT}\end{dBox}}

\declaretheorem[style=normalfontboldhead, name=Exercise, numberlike=theo]{exerciseC}
\newmdenv[skipabove=7pt,
skipbelow=7pt,
rightline=false,
leftline=true,
topline=false,
bottomline=false,
linecolor=mydarkgreen,
innerleftmargin=5pt,
innerrightmargin=5pt,
innertopmargin=0pt,
leftmargin=2cm,
rightmargin=0cm,
linewidth=4pt,
innerbottommargin=0pt]{eBox}
\newenvironment{exercise}{\begin{eBox}\begin{exerciseC}}{\end{exerciseC}\end{eBox}}

\declaretheorem[style=normalfontboldhead, name=Remark, numberlike=theo]{remarekC}
\newmdenv[skipabove=7pt,
skipbelow=7pt,
rightline=false,
leftline=true,
topline=false,
bottomline=false,
linecolor=mydarkpurple,
innerleftmargin=5pt,
innerrightmargin=5pt,
innertopmargin=0pt,
leftmargin=2cm,
rightmargin=0cm,
linewidth=4pt,
innerbottommargin=0pt]{rBox}
\newenvironment{remark}{\begin{rBox}\begin{remarekC}}{\end{remarekC}\end{rBox}}

\declaretheorem[style=normalfontboldhead, name=Example, numberlike=theo]{exampleC}
\newmdenv[skipabove=7pt,
skipbelow=7pt,
rightline=false,
leftline=false,
topline=false,
bottomline=false,
linecolor=mydarkgreen,
innerleftmargin=1pt,
innerrightmargin=5pt,
innertopmargin=0pt,
leftmargin=2cm,
rightmargin=0cm,
linewidth=4pt,
innerbottommargin=0pt]{xBox}
\newenvironment{example}{\begin{xBox}\begin{exampleC}}{\exampbar\end{exampleC}\end{xBox}}

\usepackage{adforn}  
\newcommand{\xchaptertitle}{Chapter~\thechapter~}
\newcommand{\problemname}{Problems}
\newenvironment{problemset}[1][\xchaptertitle~\problemname]{
\vspace*{10pt}
\begin{center}
\phantomsection\addcontentsline{toc}{section}{\texorpdfstring{\xchaptertitle~\problemname}{\problemname}}
\markright{#1}
\textcolor{structurecolor}{\Large\bfseries\adftripleflourishleft~#1~\adftripleflourishright}
\end{center}
\begin{enumerate}[ref=\thechapter.\theenumi]}{
\end{enumerate}}

\usepackage[shortlabels]{enumitem}  
\usepackage{enumitem}
\newcommand*{\eitemi}{\tikz \draw [baseline, ball color=structurecolor,draw=none] circle (2pt);}
\newcommand*{\eitemii}{\tikz \draw [baseline, fill=structurecolor,draw=none,circular drop shadow] circle (2pt);}
\newcommand*{\eitemiii}{\tikz \draw [baseline, fill=structurecolor,draw=none] circle (2pt);}
\setlist[enumerate,1]{label=\color{black}\arabic*.,itemsep=0pt,partopsep=0pt,parsep=\parskip,topsep=3pt}
\setlist[enumerate,2]{label=\color{black}(\alph*).,itemsep=0pt,partopsep=0pt,parsep=\parskip,topsep=3pt}
\setlist[enumerate,3]{label=\color{black}\Roman*.,itemsep=0pt,partopsep=0pt,parsep=\parskip,topsep=3pt}
\setlist[enumerate,4]{label=\color{black}\Alph*.,itemsep=0pt,partopsep=0pt,parsep=\parskip,topsep=3pt}
\setlist[itemize,1]{label={\eitemi},itemsep=0pt,partopsep=0pt,parsep=\parskip,topsep=3pt}
\setlist[itemize,2]{label={\eitemii},itemsep=0pt,partopsep=0pt,parsep=\parskip,topsep=3pt}
\setlist[itemize,3]{label={\eitemiii},itemsep=0pt,partopsep=0pt,parsep=\parskip,topsep=3pt}

%% file: chapter-worldmap.tex
\begin{figure}[htbp]
\centering
\begin{widepage1}
\resizebox{1.0\textwidth}{!}{%
\begin{tikzpicture}[>=latex]
%
%
\tikzstyle{state} = [draw, very thick, fill=white, rectangle, minimum height=3em, minimum width=6em, node distance=8em, font={\sffamily\bfseries}]
\tikzstyle{stateEdgePortion} = [black,thick];
\tikzstyle{stateEdge} = [stateEdgePortion,->];
\tikzstyle{stateEdge2} = [stateEdgePortion,<->];
\tikzstyle{edgeLabel} = [pos=0.5, text centered, font={\sffamily\small}];

%
%
\node[state, name=semidefinite,fill={colorsemidefinite}] {Semidefinite};
\node[state, name=upperleft, left of=semidefinite, xshift=-8em,fill={colorsemidefinite}] {\parbox{6.5em}{Rank-Revealing\\Semidefinite}};
\node[state, name=upperright, right of=semidefinite, xshift=+6em, fill={colorspectral}] {Spectral};
\node[state, name=schur, below of=upperright, xshift=+4em,yshift=3em, fill={colorspectral}] {Schur};
\node[state, name=cr, below of=semidefinite,yshift=-2em,  fill={colorcr}] {CR};
\node[state, name=rank, above of=cr,yshift=-3.2em,  fill={colorcr}] {Rank};
\node[state, name=interpolative, left of=cr, yshift=-2em,  fill={colorcr}] {Interpolative};
\node[state, name=skeleton, right of=cr, yshift=-2em, fill={colorcr}] {\parbox{3.5em}{ Skeleton (CUR)}};
\node[state, name=cholesky, below of=cr, right of=cr, xshift=8em,yshift=1em, fill={colorlu}] {Cholesky};
\node[ellipse, name=bunchkaufman, draw,font={\sffamily\bfseries}, below of=cholesky,xshift=-1.2em,yshift=-1.5em,fill={colorlu}]  {BunchKaufman};
\node[state, name=lu, below of=cr, left of=cr, xshift=-8em, yshift=1em, fill={colorlu}] {LU};
\node[state, name=biconjugate, below of=cr, node distance=14em,yshift=-1em, fill={colorbiconjugate}] {Biconjugate};

\node[state, name=qr, below of=biconjugate, left of=biconjugate, node distance=7em, xshift=-9em, yshift=-1em, fill={colorqr}] {QR (QL)};
\node[state, name=lq, above of=qr, yshift=0.3em, xshift=-2em, yshift=-2.5em, fill={colorqr}] {LQ (RQ)};
\node[state, name=utv, right of=lq, yshift=0em, xshift=2.4em, fill={colorqr}] {UTV};
\node[state, name=twosidedortho, draw, above of=lq,xshift=5em, yshift=-3em,fill={colorqr}]  {\parbox{5em}{Two-Sided \\Orthogonal}};

\node[state, name=svd, below of=biconjugate, right of=biconjugate, node distance=7em, xshift=9em,yshift=-1em, fill={colorsvd}] {SVD};
\node[state, name=closing, below of=biconjugate, node distance=14em,xshift=-8em, fill={colorreduction}] {\parbox{4.3em}{Hessenberg\\\,(HT form)}};
\node[state, name=polar, above of=svd, yshift=-1em, fill={colorsvd}] {Polar};
\node[state, name=bottomright, below of=svd,yshift=1em,xshift=1em,fill={coloreigenvalue}] {Eigenvalue};
\node[state, name=lastAck, below of=bottomright,fill={coloreigenvalue}] {Jordan};

\node[state, name=tridiagonal, below of=closing, fill={colorreduction}] {Tridiagonal};
\node[state, name=bottomleft, below of=qr, node distance=16em, yshift=1em,xshift=-2.5em, fill={colorreduction}] {Bidiagonal};
\node[ellipse, name=als, draw,font={\sffamily\bfseries}, left of=polar,xshift=-5em,yshift=-1em,fill={colorals}]  {ALS};
\node[ellipse, name=nmf, draw,font={\sffamily\bfseries}, above of=als,yshift=1.5em,fill={colorals}]  {NMF};

\node[state, name=cp, right of=tridiagonal,yshift=-1em,xshift=0.85em,fill={colorcp}] {CP};
\node[state, name=tucker, above of=cp,yshift=-3em,fill={colorcp}] {Tucker};
\node[state, name=tensortrain, right of=tucker,xshift=-0.5em,fill={colorcp}] {TensorTrain};
\node[state, name=hosvd, above of=tucker,yshift=-3em,fill={colorcp}] {HOSVD};


\draw ($(tucker.north) + (0,0em)$)
edge[stateEdge] node[edgeLabel, yshift=0em, xshift=0em]{Orthogonal} 
(hosvd.south) ;

\draw ($(svd.west) + (0,-1em)$)
edge[stateEdge] node[edgeLabel, yshift=0em, xshift=-0.7em]{} 
($(tucker.north) + (3em,0em)$) ;
\draw[decoration={text along path,
text={|\sffamily\small|Left Singular Vectors},text align={center}},decorate] ($(tucker.north) + (3em,0.3em)$) to [bend left=0]  ($(svd.west) + (0,-0.7em)$);

\draw ($(svd.south) + (0,0em)$)
edge[stateEdge] node[edgeLabel, yshift=0em, xshift=-0.7em]{} 
(tensortrain.north) ;
\draw[decoration={text along path,
text={|\sffamily\small|TT-SVD Algo.},text align={center}},decorate] ($(tensortrain.north) + (0,0.25em)$) to [bend left=0]  ($(svd.south) + (0,0.25em)$);

\draw ($(svd.west) + (0,-.5em)$)
edge[stateEdge] node[edgeLabel, yshift=0em, xshift=-0.7em]{} 
($(hosvd.north) + (3em,0em)$) ;
\draw[decoration={text along path,
text={|\sffamily\footnotesize |Left Singular Vectors},text align={center}},decorate] ($(hosvd.north) + (-1em,-1em)$) to [bend left=0]  ($(svd.west) + (0,-.5em)$);

\draw (cp.north)
edge[stateEdge] node[edgeLabel, yshift=0em, xshift=0.em]{More General} 
(tucker.south) ;

\draw (cp.east)
edge[stateEdge, bend right=22.5] node[edgeLabel, yshift=0em, xshift=0.em]{Clear Rank} 
(tensortrain.south) ;

\draw[decoration={text along path,
	text={|\sffamily|Low-Rank},text align={center}},decorate] ($(als.south) + (-1em,0em)$)	 to [bend right=22.5]  ($(svd.west) + (-0em,0em)$) ;
\draw ($(svd.west) + (0em,1em)$)
edge[stateEdge2, bend left=22.5] node[edgeLabel, yshift=0em, xshift=-0em]{} 
(als.south) ;

\draw (lq.east)
edge[stateEdge] node[edgeLabel, yshift=0em, xshift=-0.7em]{\parbox{3em}{Rank\\Estimation}} 
(utv.west) ;

\draw (rank.south)
edge[stateEdge] node[edgeLabel, yshift=0em, xshift=-0em]{Special Case} 
(cr.north) ;

\draw (rank.west)
edge[stateEdge, bend left=-30] node[edgeLabel, yshift=0em, xshift=-0em]{} 
($(interpolative.north) + (-2.5em,0em)$) ;
\draw[decoration={text along path,
	text={|\sffamily|Special Case},text align={center}},decorate] ($(interpolative.north) + (-3em,0.3em)$) to [bend left=30]  ($(rank.west) + (-0em,0.em)$) ; 

\draw (rank.east)
edge[stateEdge, bend left=30] node[edgeLabel, yshift=0em, xshift=-0em]{} 
($(skeleton.north) + (2.5em,0em)$) ;
\draw[decoration={text along path,
	text={|\sffamily|``Special" Case},text align={center}},decorate] ($(rank.east) + (0.7em,0.2em)$) to [bend left=30]  ($(skeleton.north) + (2.5em,0.5em)$);

\draw ($(cr.west) + (0,1em)$)
edge[stateEdge, bend left=-12.5] node[edgeLabel, xshift=-0.5em, yshift=0.8em]{} 
(interpolative.north);
\draw[decoration={text along path,
	text={|\sffamily|Independent},text align={center}},decorate] ($(interpolative.north) + (-0.5em,0.1em)$) to [bend left=12.5]  ($(cr.west) + (.5em,1.2em)$);
\draw[decoration={text along path,
	text={|\sffamily|Columns},text align={center}},decorate] ($(interpolative.north) + (-0.8em,-1.1em)$) to [bend left=12.5]  ($(cr.west) + (2.5em,0.3em)$);

\draw ($(qr.north) + (3em,0em)$)
edge[stateEdge, bend left=-22.5] node[edgeLabel, xshift=-1em,yshift=0.5em]{Rank Estimation} 
(utv.south) ;

\draw ($(lq.north) + (0em,0em)$)
edge[stateEdge, bend left=22.5] node[edgeLabel, xshift=1.5em,yshift=-0.5em]{Phase 1} 
(twosidedortho.west) ;
\draw (twosidedortho.east) 
edge[stateEdge, bend left=+22.5] node[edgeLabel, xshift=-1.5em,yshift=-0.5em]{Phase 2} 
($(utv.north) + (0em,0em)$);

\draw ($(qr.north -| lq.south)$)
edge[stateEdge, bend left=0] node[edgeLabel, yshift=0em, xshift=-0.5em]{Row Space} 
(lq.south) ;

\draw ($(semidefinite.west) + (0,.5em)$) 
edge[stateEdge, bend left=1] node[edgeLabel, xshift=0em,yshift=0.5em]{Reveal Rank} 
($(upperleft.east) + (-0em,.5em)$);

\draw  ($(lu.east) + (0,1em)$)
edge[stateEdge, bend left=-22.5] node[edgeLabel, xshift=0em, yshift=0.1em]{RREF via Elimination} 
($(cr.south) + (-1em,0)$);

\draw ($(lu.north) + (.5em,0)$) 
edge[stateEdge, bend left=45] node[edgeLabel,xshift=-0em,yshift=1em]{PSD} 
($(semidefinite.west) + (0,-.5em)$);

\draw ($(cholesky.north) + (-.5em,0)$) 
edge[stateEdge, bend right=45] node[edgeLabel,xshift=-2.3em, yshift=1em]{PSD} 
($(semidefinite.east) + (0,-.5em)$);
\draw ($(semidefinite.east) + (0,.5em)$) 
edge[stateEdge, bend left=45] node[edgeLabel,xshift=0em,yshift=1em]{PD} 
($(cholesky.north) + (.5em,0)$);

\draw (lu.east)
edge[stateEdge] node[edgeLabel, yshift=1em]{Positive Definite} 
(cholesky.west) ;

\draw (cholesky.south)
edge[stateEdge] node[edgeLabel, yshift=-0.1em]{Indefinite} 
(bunchkaufman.north) ;

\draw (bottomright.south)
edge[stateEdge] node[edgeLabel, xshift=0.5em,yshift=0em]{\parbox{6em}{Generalized \\Eigenvector}} 
(lastAck.north);

\draw[decoration={text along path,
text={|\sffamily|(X,Y)=(I,I)},text align={center}},decorate] ($(lu.east) + (2em,-0.5em)$) to [bend left=22.5]  ($(biconjugate.north) + (2em,-2em)$);
\draw ($(biconjugate.north) + (-1em,0)$)
edge[stateEdge, bend left=-22.5] node[edgeLabel, xshift=2em, yshift=-0.5em]{} 
($(lu.east) + (0,-1em)$);

\draw[decoration={text along path,
text={|\sffamily|(X,Y)=(I,I)},text align={center}},decorate] ($(biconjugate.north) + (-2em,-2.2em)$) to [bend left=22.5]  ($(cholesky.west) + (-2em,-0.7em)$);
\draw ($(biconjugate.north) + (1em,0)$)
edge[stateEdge, bend right=-22.5] node[edgeLabel, xshift=-2em, yshift=-0.5em]{} 
($(cholesky.west) + (0,-1em)$);

\draw[decoration={text along path,
text={|\sffamily|(X,Y)=(I,A)$},text align={center}},decorate] 	($(qr.east) + (2em,-0.6em)$) to [bend right=22.5]  ($(biconjugate.south) + (2em,2em)$);
\draw ($(biconjugate.south) + (-1em,0)$) 
edge[stateEdge, bend left=22.5] node[edgeLabel, xshift=2em, yshift=1em]{} 
($(qr.east) + (0,.5em)$);

\draw[decoration={text along path,
text={|\sffamily|(X,Y)=(V,U)$},text align={center}},decorate] 	($(biconjugate.south) + (-2em,1.7em)$)  to [bend right=22.5]  ($(svd.west) + (-2em,-0.3em)$);
\draw ($(biconjugate.south) + (1em,0)$) 
edge[stateEdge, bend right=22.5] node[edgeLabel, xshift=-2em, yshift=1em]{} 
($(svd.west) + (0,.5em)$);

\draw (qr.south) 
edge[stateEdge2, bend left=-22.5] node[edgeLabel, xshift=-0em]{Householder} 
(bottomleft.north);

\draw ($(qr.east) + (0,-.5em)$) 
edge[stateEdge2, bend left=22.5] node[edgeLabel, yshift=-1.3em,xshift=2em]{\parbox{4.3em}{Householder\\ \,Givens}} 
(closing.north);

\draw (qr.south east) 
edge[stateEdge2, bend left=-32.5] node[edgeLabel, xshift=0em, yshift=2em, text width=4em]{Householder}
(tridiagonal.north west);

\draw (bottomleft.east)
edge[stateEdge, bend right=12.5] node[edgeLabel, xshift=-0em, yshift=-1em]{$\bT=\bB^\top\bB$} 
(tridiagonal.west);

\draw ($(cr.east) + (0,1em)$)
edge[stateEdge, bend left=12.5] node[edgeLabel, xshift=0.5em, yshift=0.8em]{} 
(skeleton.north);
\draw[decoration={text along path,
	text={|\sffamily|Same C},text align={center}},decorate] ($(cr.east) + (-1.8em,1.2em)$) to [bend left=20.5]  ($(skeleton.north) + (2.5em,-1.5em)$);

\draw(als.north)
edge[stateEdge] node[edgeLabel, xshift=-0em, yshift=0em]{Nonnegative} 
(nmf.south);

\draw(schur.north)
edge[stateEdge, bend right=12.5] node[edgeLabel, xshift=-0em, yshift=-1em]{Symmetric} 
(upperright.east);

\draw (closing) 
edge[stateEdge] node[edgeLabel, xshift=-0em]{Symmetric} 
(tridiagonal);

\draw (svd) 
edge[stateEdge] node[edgeLabel,xshift=1.6em, yshift=0em]{\parbox{12em}{Same Derivation,\\Different Transform.}}
(polar);

%
%
\coordinate (lastAck2ClosedA) at ($(bottomright.east) + (2em,0)$);
\coordinate (lastAck2ClosedB) at ($(upperright.north -| lastAck.east) + (2em,1em)$);
\coordinate (lastAck2ClosedC) at ($(upperright.north) + (0.5em,1em)$);
\draw (bottomright.east) edge[stateEdgePortion] (lastAck2ClosedA);
\draw (lastAck2ClosedA) edge[stateEdgePortion] (lastAck2ClosedB);
\draw (lastAck2ClosedB) edge[stateEdgePortion] 
node[edgeLabel,xshift=-0.5em, yshift=1.1em]{Diagonalization}(lastAck2ClosedC);
\draw (lastAck2ClosedC) edge[stateEdge] ($(upperright.north) + (0.5em,0)$);

\coordinate (svd2ClosedA) at ($(svd.east) + (0em,0em)$);
\coordinate (svd2ClosedB) at ($(svd.east -| svd.east) + (2em,0em)$);
\coordinate (svd2ClosedC) at ($(skeleton.east) + (10em,0em)$);
\draw (svd.east) edge[stateEdgePortion] (svd2ClosedA);
\draw (svd2ClosedA) edge[stateEdgePortion] (svd2ClosedB);
\draw (svd2ClosedB) edge[stateEdgePortion] 
node[edgeLabel,xshift=-0.5em, yshift=1.1em]{}(svd2ClosedC);
\draw (svd2ClosedC) edge[stateEdge] ($(skeleton.east) + (0em,0)$);
\coordinate (svd2ClosedBB) at ($(svd.east -| svd.east) + (1.3em,0em)$);
\coordinate (svd2ClosedCC) at ($(skeleton.east) + (9.3em,0em)$);
\draw[decoration={text along path,
text={|\sffamily|Pseudoskeleton / Low-Rank Approx.},text align={center}},decorate] (svd2ClosedCC) to [bend left=180]  (svd2ClosedBB);


\coordinate (qr2inter1) at ($(qr.west) + (0,0em)$);
\coordinate (qr2inter2) at ($(qr.west) + (-2.8em,0em)$);
\coordinate (qr2inter3) at ($(interpolative.west -| qr.west) + (-2.8em,0em)$);
\draw (qr2inter1) edge[stateEdgePortion] (qr2inter2);
\draw (qr2inter2) edge[stateEdgePortion] (qr2inter3);
\draw (qr2inter3) edge[stateEdge] 
node[edgeLabel, text width=7.25em, yshift=0.8em]{CPQR or RRQR} 
(interpolative.west);

\coordinate (tri2svdA) at ($(tridiagonal.east)$);
\coordinate (tri2svdB) at ($(tridiagonal.east) + (1.2em,0)$);
\coordinate (tri2svdC) at ($(svd.west -| tri2svdB)$);
\draw (tri2svdA) edge[stateEdgePortion] (tri2svdB);
\draw (tri2svdB) edge[stateEdgePortion] (tri2svdC);
\draw (tri2svdC)
edge[stateEdge] node[edgeLabel, yshift=-0.2em,xshift=-1.6em]{\parbox{13em}{Compute SVD:\\QR Algorithm With Bidiagonal}} 
(svd.west) ;


%
%
\coordinate (bottomleft2ClosedA) at ($(tridiagonal.south) + (0,-2em)$);
\coordinate (bottomleft2ClosedB) at ($(bottomleft.south -| bottomleft.west) + (-1.3em,-1em)$);
\coordinate (bottomleft2ClosedC) at ($(upperright.north -| bottomleft.west) + (-1.3em,1em)$);
\coordinate (bottomleft2ClosedD) at ($(upperright.north) + (-0.5em,1em)$);
\draw (tridiagonal.south) edge[stateEdgePortion] (bottomleft2ClosedA);
\draw (bottomleft2ClosedA) edge[stateEdgePortion] (bottomleft2ClosedB);
\draw (bottomleft2ClosedB) edge[stateEdgePortion] (bottomleft2ClosedC);
\draw (bottomleft2ClosedC) edge[stateEdgePortion] 
node[edgeLabel, text width=7.25em, yshift=1.1em]{QR Algorithm} 
(bottomleft2ClosedD);
\draw (bottomleft2ClosedD) edge[stateEdge] ($(upperright.north) + (-0.5em,0)$);

\begin{pgfonlayer}{background}
\draw [join=round,cyan,dotted,fill={colormiddleright}] ($(svd.north west) + (-1em, -4em)$) rectangle ($(polar.south east) + (1em, 4em)$);
\draw [join=round,cyan,dotted,fill={colorbottomleft}] ($(qr.north west) + (-3em, 1em)$) rectangle ($(tridiagonal.south east) + (1em, -1.5em)$);
\draw [join=round,cyan,dotted,fill={coloruppermiddle}] ($(rank.north west  -| interpolative.north west) + (-0.8em, +0.5em)$) rectangle ($(skeleton.south east) + (0.8em, -0.5em)$);
\draw [join=round,cyan,dotted,fill={colormiddleleft}] ($(qr.south west) + (-2.5em, -0.5em)$) rectangle ($(twosidedortho.north east -| utv.north east) + (0.6em, 0.5em)$);		
\draw [join=round,cyan,dotted,fill={colorbottomright}] ($(lastAck.south west) + (-0.5em, -0.5em)$) rectangle ($(bottomright.north east) + (0.6em, 0.5em)$);
\draw [join=round,cyan,dotted,fill={colorupperright}] ($(schur.south east) + (0.5em, -0.5em)$) rectangle ($(upperright.north west) + (-0.5em, 0.5em)$);
\draw [join=round,cyan,dotted,fill={colorupperleft}] ($(upperleft.south west) + (-0.5em, -0.5em)$) rectangle ($(semidefinite.north east) + (0.5em, 0.5em)$);
\draw [join=round,cyan,dotted,fill={colormiddle}] ($(lu.south west) + (-0.5em, -0.5em)$) rectangle ($(cholesky.north east) + (0.5em, 0.5em)$);
\draw [join=round,cyan,dotted,fill={coloralsbkg}] ($(als.south west) + (-1.2em, -1.2em)$) rectangle ($(nmf.north east) + (1.2em, 1.2em)$);
\draw [join=round,cyan,dotted,fill={colorcpborder}] ($(cp.south west) + (-0.5em, -0.5em)$) rectangle ($(hosvd.north east -|tensortrain.north east) + (0.5em, 0.5em)$);
\end{pgfonlayer}
\end{tikzpicture}
}

\end{widepage1}
\caption{Matrix Decomposition World Map.}
\label{fig:matrix-decom-world-picture}
\end{figure}

\begin{figure}[htbp]
	\centering
	\begin{widepage}
		\resizebox{1.0\textwidth}{!}{%
\begin{tikzpicture}[node distance=2cm]
\node (matrix) [condition, rectangle split, rectangle split parts=2]
{
	\textbf{Matrix $\bA$}
	\nodepart{second}$\bA\in \real^{m\times n}$
};
\node (AuxNode01) [text width=4cm, below=of matrix] {};
\node (square) [condition, rectangle split, rectangle split parts=2, below=of AuxNode01,xshift=-4.3cm,yshift=3cm]
{
	\textbf{Square}
	\nodepart{second}$\bA\in \real^{n\times n}$
};
\node (rectangular) [condition, rectangle split, rectangle split parts=2, below=of AuxNode01,xshift=5.7cm,yshift=3cm]
{
	\textbf{Rectangular}
	\nodepart{second}$m\neq n$
};
\node (AuxNode02) [text width=0.5cm, below=of square] {};
\node (symmetric) [condition, rectangle split, rectangle split parts=2, below=of AuxNode02,xshift=-3cm,yshift=3cm]
{
	\textbf{Symmetric}
	\nodepart{second}$\bA=\bA^\top$
};
\node (asymmetric) [condition, rectangle split, rectangle split parts=2, below=of AuxNode02,xshift=2cm,yshift=3cm]
{
	\textbf{Asymmetric}
	\nodepart{second}$\bA\neq \bA^\top$
};

\node (AuxNode03) [below=of symmetric] {};
\node (pd) [abstract, rectangle split, rectangle split parts=2, below=of AuxNode03, xshift=0cm, yshift=3cm]
{
	\textbf{PD or PSD}
	\nodepart{second}$\bx^\top\bA\bx\geq0$
};
\node (cholesky) [comment, rectangle split, rectangle split parts=2, below=0.2cm of pd, text justified]
{
	\textbf{Cholesky, PD}
	\nodepart{second}$\bA=\bR^\top\bR$ 
};
\node (pcholesky) [comment, rectangle split, rectangle split parts=2, below=0.2cm of cholesky, text justified]
{
	\textbf{Pivoted Cholesky, PD}
	\nodepart{second}$\bP\bA\bP^\top = \bR^\top\bR $
};
\node (semidefinite) [comment, rectangle split, rectangle split parts=2, below=0.2cm of pcholesky, text justified]
{
	\textbf{Semidefi. PSD}
	\nodepart{second}$\bA=\bR^\top\bR$
};
\node (rrsemidefinite) [comment, rectangle split, rectangle split parts=2, below=0.2cm of semidefinite, text justified]
{
	\textbf{RR Semi., PSD}
	\nodepart{second}$\bP^\top\bA\bP=\bR^\top\bR$
};
\node (bunchkaufman) [comment, rectangle split, rectangle split parts=2, below=0.2cm of rrsemidefinite, text justified]
{
	\textbf{BunchKaufman}
	\nodepart{second}$\bP^\top\bA\bP=\bL\bB\bL^\top$
};
\node (eigen) [abstract, rectangle split, rectangle split parts=2, below=1cm of bunchkaufman]
{
	\textbf{Eigenvector}
	\nodepart{second}Orthogonal
};
\node (spectral) [comment, rectangle split, rectangle split parts=2, below=0.2cm of eigen, text justified]
{
	\textbf{Spectral}
	\nodepart{second}$\bA=\bQ\bLambda\bQ^\top$
};
\node (AuxNode04) [below=of asymmetric] {};
\node (generalizedeigen) [abstract, rectangle split, rectangle split parts=2, below=of AuxNode04, xshift=-1.3cm, yshift=3cm]
{
	\textbf{General Eigen}
	\nodepart{second}Generalized 
};
\node (eigenvalue) [comment, rectangle split, rectangle split parts=2, below=0.2cm of generalizedeigen, text justified]
{
	\textbf{Eigenvalue}
	\nodepart{second}$\bA=\bX\bLambda\bX^{-1}$
};
\node (schur) [comment, rectangle split, rectangle split parts=2, below=0.2cm of eigenvalue, text justified]
{
	\textbf{Schur}
	\nodepart{second}$\bA=\bQ\bU\bQ^\top$
};
\node (jordan) [comment, rectangle split, rectangle split parts=2, below=0.2cm of schur, text justified]
{
	\textbf{Jordan}
	\nodepart{second}$\bA=\bX\bJ\bX^{-1}$
};
\node (basis) [abstract, rectangle split, rectangle split parts=2, below=1cm of jordan]
{
	\textbf{Basis}
	\nodepart{second}Column and Row
};
\node (twosided) [comment, rectangle split, rectangle split parts=2, below=0.2cm of basis, text justified]
{
	\textbf{Two-Sided}
	\nodepart{second}$\bA\bP\bA=\bU\bF\bV^\top$
};
\node (reduction) [abstract, rectangle split, rectangle split parts=2, right=0.4cm of generalizedeigen]
{
	\textbf{Reduction}
	\nodepart{second}Ortho. Similar
};
\node (hessenberg) [comment, rectangle split, rectangle split parts=2, below=0.2cm of reduction, text justified]
{
	\textbf{Hessenberg}
	\nodepart{second}$\bA=\bQ\bH\bQ^\top$
};
\node (tridiagonal) [comment, rectangle split, rectangle split parts=2, below=0.2cm of hessenberg, text justified]
{
	\textbf{Tridiagonal}
	\nodepart{second}$\bA=\bQ\bT\bQ^\top$
};
\node (gaussianelimi) [abstract, rectangle split, rectangle split parts=2, below=1cm of tridiagonal]
{
	\textbf{Gaussian Elim.}
	\nodepart{second}Element. Trans.
};
\node (lu) [comment, rectangle split, rectangle split parts=2, below=0.2cm of gaussianelimi, text justified]
{
	\textbf{LU}
	\nodepart{second}$\bA=\bL\bU$\\ $\bA=\bL\bD\bU$\\$\bA=\bP\bL\bU$
};
\node (rrlu) [comment, rectangle split, rectangle split parts=2, below=0.2cm of lu, text justified]
{
	\textbf{RRLU}
	\nodepart{second}$\bP\bA\bQ=\bL\bU$
};
\node (cplu) [comment, rectangle split, rectangle split parts=2, below=0.2cm of rrlu, text justified]
{
	\textbf{Complete \\Pivoting LU}
	\nodepart{second}$\bP\bA\bQ=\bL\bU$
};

\node (AuxNode05) [below=of rectangular] {};
\node (fullrank) [abstract, rectangle split, rectangle split parts=2, left=of AuxNode05, xshift=2cm]
{
	\textbf{Full Rank}
	\nodepart{second}Basis Space
};
\node (general) [abstract, rectangle split, rectangle split parts=2, right=of AuxNode05, xshift=-2cm]
{
	\textbf{General}
	\nodepart{second}$\bA\in \real^{m\times n}$
};
\node (qr) [comment, rectangle split, rectangle split parts=2, below=0.2cm of fullrank, text justified]
{
	\textbf{QR,\\column space}
	\nodepart{second}$\bA=\bQ\bR$
};
\node (lq) [comment, rectangle split, rectangle split parts=2, below=0.2cm of qr, text justified]
{
	\textbf{LQ, row space}
	\nodepart{second}$\bA=\bL\bQ$
};

\node (CoolingLoopInstants) [comment, rectangle split, rectangle split parts=2, below=0.2cm of general, text justified]
{
	\textbf{SVD}
	\nodepart{second}$\bA=\bU\bSigma\bV^\top$
};
\node (polar) [comment, rectangle split, rectangle split parts=2, below=0.2cm of CoolingLoopInstants, text justified]
{
	\textbf{Polar}
	\nodepart{second}$\bA=\bQ_l\bS_l$
};
\node (bidiagonal) [comment, rectangle split, rectangle split parts=2, below=0.2cm of polar, text justified]
{
	\textbf{Bidiagonal}
	\nodepart{second}$\bA=\bU\bB\bV^\top$
};
\node (utv) [comment, rectangle split, rectangle split parts=2, below=0.2cm of bidiagonal, text justified]
{
	\textbf{UTV}
	\nodepart{second}$\bA=\bU\bT\bV$
};
\node (rrqr) [comment, rectangle split, rectangle split parts=2, below=0.2cm of utv, text justified]
{
	\textbf{RRQR}
	\nodepart{second}$\bA=\bQ\bR\bP^\top$
};
\node (tensordecom) [condition, rectangle split,rectangle split parts=2, below=1cm of rrqr, text justified]
{
	\textbf{Tensor Decom.}
	\nodepart{second}$\eX\in \real^{I\times J\times K}$
};
\node (cp) [comment, rectangle split, rectangle split parts=2, below=0.2cm of tensordecom, text justified]
{
	\textbf{CP}
	\nodepart{second}$\sum_{r=1}^{R} \ba_r\circ \bb\circ  \bc$
};
\node (tucker) [comment, rectangle split, rectangle split parts=2, below=0.2cm of cp, text justified]
{
	\textbf{Tucker/HOSVD}
	\nodepart{second}$\llbracket\eG; \bA, \bB, \bC   \rrbracket  $
};
\node (tensortrain) [comment, rectangle split, rectangle split parts=2, below=0.2cm of tucker, text justified]
{
	\textbf{Tensor-Train}
	\nodepart{second}$\bA\boxtimes \bB\boxtimes \bC$
};

\node (datadistill) [abstract, rectangle split, rectangle split parts=2, below=1cm of lq]
{
	\textbf{Data Distill}
	\nodepart{second}Low-Rank App.
};
\node (cr) [comment, rectangle split, rectangle split parts=2, below=0.2cm of datadistill, text justified]
{
	\textbf{CR}
	\nodepart{second}$\bA=\bC\bR$
};
\node (rank) [comment, rectangle split, rectangle split parts=2, below=0.2cm of cr, text justified]
{
	\textbf{Rank}
	\nodepart{second}$\bA=\bD\bF$
};
\node (skeleton) [comment, rectangle split, rectangle split parts=2, below=0.2cm of rank, text justified]
{
	\textbf{Skeleton}
	\nodepart{second}$\bA=\bC\bU^{-1}\bR$
};
\node (interpolative) [comment, rectangle split, rectangle split parts=2, below=0.2cm of skeleton, text justified]
{
	\textbf{Interpolative}
	\nodepart{second}$\bA=\bC\bW$
};
\node (als) [comment, rectangle split, rectangle split parts=2, below=0.2cm of interpolative, text justified]
{
	\textbf{ALS}
	\nodepart{second}$\bA\approx\bW\bZ$
};
\node (nmf) [comment, rectangle split, rectangle split parts=2, below=0.2cm of als, text justified]
{
	\textbf{NMF}
	\nodepart{second}$\bA\approx\bW\bZ$
};
\draw[myarrow] (square.north) -- ++(0,0.8) -| (matrix.south);
\draw[line] (square.north) -- ++(0,0.8) -| (rectangular.north);

\draw[myarrow] (symmetric.north) -- ++(0,0.8) -| (square.south);
\draw[line] (symmetric.north) -- ++(0,0.8) -| (asymmetric.north);

\draw[line] (pd.west) -- ++(-0.2,0);
\draw[myarrow] (eigen.west) -- ++(-0.2,0) -- ([yshift=0.5cm, xshift=-0.2cm] pd.north west) -|
([xshift=-1cm]symmetric.south);

\draw[line] (basis.west) -- ++(-0.2,0);
\draw[line] (generalizedeigen.west) -- ++(-0.2,0);
\draw[myarrow] (basis.west) -- ++(-0.2,0) -- ([yshift=0.5cm, xshift=-0.2cm] generalizedeigen.north west) -|
([xshift=-1cm]asymmetric.south);
([xshift=1cm]asymmetric.south);

\draw[line] (reduction.west) -- ++(-0.2,0);
\draw[line] (gaussianelimi.west) -- ++(-0.2,0);
\draw[myarrow] (gaussianelimi.west) -- ++(-0.2,0) -- ([yshift=0.5cm, xshift=-0.2cm] reduction.north west)  -|
([xshift=1cm]asymmetric.south);

\draw[line] (datadistill.west) -- ++(-0.2,0);
\draw[line] (fullrank.west) -- ++(-0.2,0);
\draw[myarrow] (datadistill.west) -- ++(-0.2,0) -- ([yshift=0.5cm, xshift=-0.2cm] fullrank.north west)  -|
([xshift=-1cm]rectangular.south);

\draw[myarrow] (general.north) -- ++(0,0.8) -| ([xshift=1cm]rectangular.south);


\end{tikzpicture}
}
\end{widepage}
\caption{Matrix Decomposition World Map Under Conditions.}
\label{fig:matrix-decom-world-picture2}
\end{figure}

%% file: chapter-notation.tex
\chapter*{Notation}\label{notation}
\index{Notation}


This section provides a concise reference describing notation used throughout this
book.
If you are unfamiliar with any of the corresponding mathematical concepts,
the book describes most of these ideas in Chapter~\ref{chapter_introduction}.

\vspace{0.4in}
\begin{minipage}{\textwidth}
\centerline{\bf General Notations}
\bgroup
\def\arraystretch{1.5}
\begin{tabular}{cp{4.25in}}
$\displaystyle \triangleq$ & equals by definition \\ 
$\displaystyle :=, \leftarrow$ & equals by assignment \\ 
$\displaystyle \pi$ &	 3.141592$\ldots$ \\
$\displaystyle e, \exp$ &	 2.71828$\ldots$
\end{tabular}
\egroup
\end{minipage}

\vspace{0.4in}
\begin{minipage}{\textwidth}
\centerline{\bf Numbers and Arrays}
\bgroup
\def\arraystretch{1.5}
\begin{tabular}{cp{4.25in}}
$\displaystyle a$   & A scalar (integer or real)\\
$\displaystyle \ba$ & A vector\\
$\displaystyle \bA$ & A matrix\\
$\displaystyle \eA$ & A tensor\\
$\displaystyle \bI_n$ & Identity matrix with $n$ rows and $n$ columns\\
$\displaystyle \bI$   & Identity matrix with dimensionality implied by context\\
$\displaystyle \be_i$ & Standard basis/canonical orthonormal  vector $[0,\dots,0,1,0,\dots,0]$ with a 1 at position $i$\\
$\displaystyle \text{diag}(\va)$ & A square, diagonal matrix with diagonal entries given by $\ba$\\
$\displaystyle \text{diag}(\bA, \bB,\ldots)$&  a block-diagonal matrix=($\bA\oplus \bB\oplus \ldots$), the direct sum notation\\
$\displaystyle \ra$   & A scalar random variable\\
$\displaystyle \rva$  & A vector-valued random variable\\
$\displaystyle \rmA$  & A matrix-valued random variable\\
\end{tabular}
\egroup
\index{Scalar}
\index{Vector}
\index{Matrix}
\index{Tensor}
\end{minipage}

\index{Sets}
\vspace{0.2in}
\begin{minipage}{\textwidth}
\centerline{\bf Sets}
\bgroup
\def\arraystretch{1.5}
\begin{tabular}{cp{4.25in}}
$\displaystyle \sA, \sS$ & A set\\
$\displaystyle \varnothing$ & The null set \\
$\displaystyle \real, \real_+, \real_{++}, \complex$ & The set of real/nonnegative/positive/complex numbers \\
$\displaystyle \naturalset, \integer$ & The set of natural/integer numbers \\
$\displaystyle \sF$ & The set of either real or complex numbers \\
$\displaystyle \{0, 1\}$ & The set containing 0 and 1 \\
$\displaystyle \{0, 1, \dots, n \}$ & The set of all integers between $0$ and $n$\\
$\displaystyle [a, b]$ & The real interval including $a$ and $b$\\
$\displaystyle (a, b]$ & The real interval excluding $a$ but including $b$\\
$\displaystyle \sA \backslash \sB$ & Set subtraction, i.e., the set containing the elements of $\sA$ that are not in $\sB$\\
\end{tabular}
\egroup
\index{Scalar}
\index{Vector}
\index{Matrix}
\index{Tensor}

\index{Set}
\end{minipage}

\index{Matrix indexing}
\vspace{0.2in}
\begin{minipage}{\textwidth}
\centerline{\bf Indexing}
\bgroup
\def\arraystretch{1.5}
\begin{tabular}{cp{4.25in}}
$\displaystyle \eva_i$ & Element $i$ of vector $\va$, with indexing starting at 1 \\
$\displaystyle \ba_{-i}$ & All elements of vector $\va$ except for element $i$ \\
$\displaystyle  \bA[i,j]=a_{ij}$ & Element $(i, j)$ of matrix $\mA$ \\
$\displaystyle \mA_{i, :}=\bA[i,:]$ & Row $i$ of matrix $\mA$ \\
$\displaystyle \mA_{:, i}=\bA[:,i]=\ba_i$ & Column $i$ of matrix $\mA$ \\
$\displaystyle \eA_{ijk}=\elA_{ijk}=a_{ijk}$ & Element $(i, j, k)$ of a 3-D tensor $\eA$\\
$\displaystyle \eA_{:, :, i}=\eA[:,:,i]$ & 2-D slice of a 3-D tensor $\eA$\\
\end{tabular}
\egroup
\end{minipage}

\vspace{0.4in}
\begin{minipage}{\textwidth}
\centerline{\bf Calculus}
\bgroup
\def\arraystretch{1.5}
\begin{tabular}{cp{4.25in}}
$\displaystyle\frac{d y} {d x}$ & Derivative of $y$ with respect to $x$\\ [2ex]
$\displaystyle \frac{\partial y} {\partial x} $ & Partial derivative of $y$ with respect to $x$ \\
$\displaystyle \nabla_{\bx} y $ & Gradient of $y$ with respect to $\bx$ \\
$\displaystyle \nabla_{\bX} y $ & Matrix derivatives of $y$ with respect to $\bX$ \\
$\displaystyle \nabla_{\eX} y $ & Tensor containing derivatives of $y$ with respect to $\eX$ \\
$\displaystyle \frac{\partial f}{\partial \vx} $ & Jacobian matrix $\mJ \in \R^{m\times n}$ of $f: \R^n \rightarrow \R^m$\\
$\displaystyle \nabla_\vx^2 f(\vx)\text{ or }\mH( f)(\vx)$ & The Hessian matrix of $f$ at input point $\vx$\\
$\displaystyle \int f(\vx) d\vx $ & Definite integral over the entire domain of $\vx$ \\
$\displaystyle \int_\sS f(\vx) d\vx$ & Definite integral with respect to $\vx$ over the set $\sS$ \\
\end{tabular}
\egroup
\index{Derivative}
\index{Integral}
\index{Jacobian matrix}
\index{Hessian matrix}
\end{minipage}

\vspace{0.4in}
\begin{minipage}{\textwidth}
\centerline{\bf Probability and Information Theory}
\bgroup
\def\arraystretch{1.5}
\begin{tabular}{cp{4.25in}}
$\displaystyle \ra \bot \rb$ & The random variables $\ra$ and $\rb$ are independent\\
$\displaystyle \ra \bot \rb \mid \rc $ & They are conditionally independent given $\rc$\\
$\displaystyle P(\ra)$ & A probability distribution over a discrete variable\\
$\displaystyle p(\ra)$ & A probability distribution over a continuous variable, or over
a variable whose type has not been specified\\
$\displaystyle \ra \sim P$ & Random variable $\ra$ has distribution $P$\\
$\displaystyle  \Exp_{\rx\sim P} [ f(x) ]\text{ or } \Exp [f(x)]$ & Expectation of $f(x)$ with respect to $P(\rx)$ \\
$\displaystyle \Var[f(x)] $ &  Variance of $f(x)$ under $P(\rx)$ \\
$\displaystyle \Cov[f(x),g(x)] $ & Covariance of $f(x)$ and $g(x)$ under $P(\rx)$\\
$\displaystyle H(\rx) $ & Shannon entropy of the random variable $\rx$\\
$\displaystyle \KL ( P \Vert Q ) $ & Kullback-Leibler divergence of P and Q \\
$\displaystyle \mathcal{N} ( \vx ; \vmu , \mSigma)$ & Gaussian distribution %
over $\vx$ with mean $\vmu$ and covariance $\mSigma$ \\
\end{tabular}
\egroup
\index{Independence}
\index{Conditional independence}
\index{Variance}
\index{Covariance}
\index{Kullback-Leibler divergence}
\index{Shannon entropy}
\end{minipage}

\vspace{0.4in}
\begin{minipage}{\textwidth}
\centerline{\bf Functions}
\bgroup
\def\arraystretch{1.5}
\begin{tabular}{cp{4.25in}}
$\displaystyle f: \sA \rightarrow \sB$ & The function $f$ with domain $\sA$ and range $\sB$\\
$\displaystyle f \circ g $ & Composition of the functions $f$ and $g$ \\
$\displaystyle f(\vx ; \vtheta) $ & A function of $\vx$ parametrized by $\vtheta$.
(Sometimes we write $f(\vx)$ and omit the argument $\vtheta$ to lighten notation) \\
$\displaystyle \log(x), \ln(x)$ & Natural logarithm of $x$ \\
$\displaystyle \sigma(x),\, \text{Sigmoid}(x)$ & Logistic sigmoid, i.e., $\displaystyle \frac{1} {1 + \exp\{-x\}}$ \\
$\displaystyle \zeta(x)$ & Softplus, $\log(1 + \exp\{x\})$ \\
$\displaystyle \norm{\bx}_p $ & $\ell_p$ norm of $\vx$ \\
$\displaystyle \norm{\bx}_1, \norm{\bx}_2, \norm{\bx}_\infty$ & $\ell_1$ norm, $\ell_2$ norm, $\ell_\infty$ norm\\
$\displaystyle \normf{\bA}, \normtwo{\bA}$ & Frobenius and spectral norms \\
$\displaystyle \norm{\bA}_{[k]}, \norm{\bA}_{[k,p]}, \norm{\bA}_{[q,p]} $ & Ky Fan $k$-norm, \kyfanpknorm, Schatten $p$-norm \\
$\displaystyle x^+, [x]_+$ & Positive part of $x$, i.e., $\max(0,x)$\\
$\displaystyle u(x)$ & Step function with value 1 when $x\geq0$ and value 0 otherwise\\
$\displaystyle \text{Re}(x), \text{Im}(x)$ & Real and imaginary part of $x$\\
$\displaystyle \indicator\{\mathrm{condition}\}$ & is 1 if the condition is true, 0 otherwise\\
\end{tabular}
\egroup
\index{Sigmoid}
\index{Softplus}
\index{Norm}
\end{minipage}
Sometimes we use a function $f$ whose argument is a scalar but apply
it to a vector, matrix, or tensor: $f(\vx)$, $f(\mX)$, or $f(\eX)$.
This denotes the application of $f$ to the
array element-wise. For example, if $\eC = \sigma(\eX)$, then $\mathscr{C}_{ijk} = \sigma(\mathscr{X}_{ijk})$
for all valid values of $i$, $j$, and $k$.

\vspace{0.2in}
\begin{minipage}{\textwidth}
\centerline{\bf Linear Algebra Operations}
\bgroup
\def\arraystretch{1.5}
\begin{tabular}{cp{4.25in}}
$\displaystyle \abs{\bA}$ & Componentwise absolute matrix $\mA$ \\
$\displaystyle \frac{[\bA]}{[\bB]}$ &  Componentwise division between two matrices\\
$\displaystyle \bA^\top$ & Transpose of matrix $\mA$ \\
$\displaystyle \bA^{-1}$ & Inverse of $\mA$\\
$\displaystyle \bA^+$ & Moore-Penrose pseudo-inverse of $\mA$\\
$\displaystyle \bA \circledast \bB $ & Element-wise (Hadamard) product of $\mA$ and $\mB$ \\
$\displaystyle \bA \kronecker \bB $ & Kronecker product of $\mA$ and $\mB$ \\
$\displaystyle \bA \khatrirao \bB $ & Khatri-Rao product of $\mA$ and $\mB$ \\
$\displaystyle \eX \circ \eY $ & Tensor outer product of $\eX$ and $\eY$ \\
$\displaystyle \det(\bA)$ & Determinant of $\bA$ \\
$\displaystyle \trace(\bA)$ & Trace of $\bA$ \\
$\displaystyle \mathrm{rref}(\bA)$ & Reduced row echelon form of $\bA$ \\
$\displaystyle \cspace(\bA)$ & Column space of $\bA$ \\
$\displaystyle \nspace(\bA)$ & Null space of $\bA$ \\
$\displaystyle \mathcalV, \mathcalW$ & A general subspace \\
$\displaystyle \dim(\mathcalV)$ & Dimension of the space $\mathcalV$ \\
$\displaystyle \defect(\bA)$ & Defect or nullity of $\bA$ \\
$\displaystyle \rank(\bA)$ & Rank of $\bA$ \\
$\displaystyle \trace(\bA)$ & Trace of $\bA$ \\
$\displaystyle \adjugate(\bA)$ & Adjugate of $\bA$ \\
$\displaystyle \bN_n$ & Standard nilpotent matrix $\in\real^{n\times n}$\\
$\displaystyle \rho(\bA), \Lambda(\bA), \widehat{\Lambda}(\bA)$ & Spectral radius, spectrum, multispectrum of $\bA$\\
$\displaystyle \lambda_{\min}(\bA),\lambda_{\max}(\bA)$ & Smallest/largest eigenvalues of $\bA$\\
$\displaystyle \sigma_{\min}(\bA),\sigma_{\max}(\bA)$ & Smallest/largest singuar values of $\bA$\\
$\displaystyle d_{\min}(\bA),d_{\max}(\bA)$ & Smallest/largest diagonal values of $\bA$\\
$\displaystyle \bA\geq \bzero, \bA>\bzero$ & Nonnegative/positive matrix $\bA$\\
$\displaystyle \bA\succeq \bzero\, (\bA\in\psd^n), \bA\succ \bzero\, (\bA\in\pd^n)$ & Positive semidefinite/definite matrix $\bA$\\
subdiagonal/superdiagonal & Entries below/above the main diagonal
\end{tabular}
\egroup
\index{Transpose}
\index{Element-wise product|see {Hadamard product}}
\index{Hadamard product}
\index{Determinant}
\end{minipage}


\vspace{0.4in}
\begin{minipage}{\textwidth}
\centerline{\bf Abbreviations}
\bgroup
\def\arraystretch{1.5}
\begin{tabular}{cp{4.25in}}
PD & Positive definite  \\
PSD & Positive semidefinite \\
CPQR & Column-pivoted QR decomposition \\
CGS & Classical Gram-Schimdt process\\
MGS & Modified Gram-Schmidt process\\
CR & Column-row decomposition \\
REF & Row echelon form\\
RREF & Reduced row echelon form\\
ID & Interpolative decomposition\\
IID & Intervened interpolative decomposition \\
BID & Bayesian interpolative decomposition \\
MCMC & Markov chain Monte Carlo \\
i.i.d. & Independently and identically distributed \\
p.d.f. & Probability density function \\
p.m.f. & Probability mass function \\
SVD &  Singular value decomposition \\
PCA & Principal component analysis \\
OLS & Ordinary least squares\\
TN & Truncated-normal distribution \\
GTN & General-truncated-normal distribution\\
ALS & Alternating least squares \\
GD & Gradient descent\\
SGD & Stochastic gradient descent \\
MU & Multiplicative update \\
MSE & Mean squared error\\
NMF & Nonnegative matrix factorization \\
ADMM & Alternating direction methods of multipliers\\
HOSVD & High-order SVD \\
CP & Canonical Polyadic Decomposition\\
TT & Tensor-train decomposition\\

\end{tabular}
\egroup
\end{minipage}

\clearpage


%% file: chapter-intro.tex

\newpage
\chapter{Introduction and Background}\label{chapter_introduction}
\begingroup
\hypersetup{
	linkcolor=structurecolor,
	linktoc=page,  
}
\minitoc \newpage
\endgroup
\section*{Introduction and Background}
\addcontentsline{toc}{section}{Introduction and Background}

The world of matrices is as vast as it is complex, with applications ranging from the smallest quantum systems to  large-scale data analytics problems. One of the most powerful tools for understanding and utilizing matrices is decomposition---the process of reducing a matrix to its   constituent parts.  
Matrix decomposition is also a fundamental concept in linear algebra that has wide-ranging applications across various fields, including computer science, engineering, physics, and economics. At its core, this process involves breaking down a complex matrix into a set of simpler matrices that, when combined, can reconstruct the original matrix. This technique is not only theoretically fascinating but also practically invaluable, as it simplifies the solution of systems of linear equations, reduces computational complexity, and provides insights into the underlying structure of data.

The journey of matrix decomposition began with the pioneering work of Alston S. Householder in the mid-20th century, which laid the groundwork for modern numerical analysis \citep{householder1954principles}. Since then, the field has evolved dramatically, with significant developments such as the backpropagation algorithm for neural networks, the increasing need for data dimensionality reduction in machine learning, and the demand for the expressive power of low-rank matrices in natural language processing and large language models.
Today, matrix decomposition has become a core technology in statistics \citep{banerjee2014linear, gentle1998numerical}, optimization \citep{gill2021numerical}, recommender systems \citep{symeonidis2016matrix}, and machine learning \citep{goodfellow2016deep, bishop2006pattern}, largely due to advancements in deep neural networks and transformer modeling.

The primary aim of this book is to provide a self-contained introduction to the concepts and mathematical tools in numerical linear algebra and matrix analysis, 
and based on them to present major matrix decomposition techniques and their applications to the reader.
We clearly realize our inability to cover all the useful and interesting topics concerning matrix decomposition. 
We refer the reader to literature in the field of linear algebra and data analysis for a more detailed introduction to  related areas. Some excellent examples include \citet{householder2006principles, elden2007matrix, trefethen1997numerical, strang1993introduction, stewart2000decompositional, gentle2007matrix, higham2002accuracy, quarteroni2010numerical,
horn2012matrix, golub2013matrix, beck2017first, gallier2017fundamentals, boyd2018introduction, strang2019linear, van2020advanced, strang2021every}. 
Moreover, this book aims to serve as a comprehensive guide for readers who wish to delve into the intricacies of matrix decomposition, grasp its applications in-depth,  and explore the calculation and complexity of decompositional methods, with an emphasis on cost reduction.

A matrix decomposition is a way of reducing a complex matrix into its constituent parts, which are in simpler forms and represent the original matrix as a product of two (or more) \textit{factor matrices}. 
The underlying principle of the decompositional approach to matrix computation is to simplify complex matrix operations  by performing them on the decomposed parts.
Formally, a matrix decomposition task on matrix $\bA$ can be cast as 
\begin{itemize}
\item $\bA=\bQ\bU$: where $\bQ$ is an orthogonal matrix that spans the same column space as $\bA$, and $\bU$ is a relatively simple and sparse matrix to reconstruct $\bA$.
\item $\bA=\bQ\bT\bQ^\top$: where $\bQ$ is orthogonal such that $\bA$ and $\bT$ are \textit{similar matrices}
\footnote{See Definition~\ref{definition:similar-matrices}  for a rigorous definition.} 
that share properties such as same eigenvalues, ranks, and sparsity. Moreover,
operations on $\bT$ are often simplier than those on $\bA$.
\item $\bA=\bU\bT\bV$: where $\bU$ and $\bV$ are orthogonal matrices such that the columns of $\bU$ and the rows of $\bV$ constitute  orthonormal bases for the column  and row spaces of $\bA$, respectively.
\item $\underset{m\times n}{\bA}=\underset{m\times r}{\bB}\gapthree \underset{r\times n}{\bC}$: where $\bB$ and $\bC$ are full-rank matrices that can reduce the memory storage of $\bA$. In practice, a low-rank approximation $\underset{m\times n}{\bA}\approx \underset{m\times k}{\bD}\gapthree \underset{k\times n}{\bF}$ can be employed, where $k<r$ is called the \textit{numerical rank} of the matrix such that  the matrix can be stored much more inexpensively and can be multiplied rapidly with
vectors or other matrices (e.g., the LoRA used in large language models \citep{hu2021lora}).
An approximation of the form $\bA=\bD\bF$ is useful for storing the matrix $\bA$ more frugally (using $k(m+n)$ floating-point numbers instead of $mn$  for  $\bA$), for efficiently computing a matrix-vector product $\bb = \bA\bx$ (via $\bc = \bF\bx$ and $\bb = \bD\bc$), for data interpretation, and much more.

\index{Numerical rank}
\index{Linear system}
\item A matrix decomposition, which is usually expensive to
compute, can be reused to solve new problems involving the
original matrix in different scenarios (e.g., once the factorization of $\bA$ is obtained, it can be reused to solve the set of linear systems $\{\bb_1=\bA\bx_1, \bb_2=\bA\bx_2, \ldots, \bb_k=\bA\bx_k\}$).
\item More generally, a matrix decomposition 
provides a way to comprehend a matrix as a product of multiple smaller and simpler matrices,
each with a clear geometric transformation meaning (see Section~\ref{section:coordinate-transformation}).
\end{itemize}

Matrix decomposition algorithms can be summarized into numerous categories.
\begin{enumerate}
\item Factorizations arise from Gaussian elimination, including the LU decomposition and its positive definite alternative, which is called the Cholesky decomposition;
\item Factorizations obtained when either orthogonalizing the columns or the rows of a matrix such that the data can be explained well on an orthonormal basis;
\item Factorizations where the matrices are skeletoned such that a subset of the columns or the rows can represent the whole data in a small reconstruction error, whilst, the sparsity and nonnegativity of the matrices are kept as they are, the purpose of which is to find a concise but interpretable representation of the original matrix;
\item Reduction to the Hessenberg, tridiagonal, or bidiagonal form, the result of which is that  the properties of the matrices, e.g., ranks, eigenvalues, and so on, can be explored in these reduced matrices.
\item Factorizations result from computations based on the eigenvalues of matrices;
\item Factorizations that can be cast as a special kind of decompositions that involve optimization methods and high-level ideas, the category of which may not be straightforward to determine.
\end{enumerate}

To be more specific, this book delves into  essential matrix decomposition methods such as QR, ID, CUR, and SVD, providing not only the theoretical underpinnings but also the computational complexities associated with each technique. We explore the nuances of algorithms like Givens rotation, Householder reflector, Bayesian approaches, and randomized algorithms.
Matrix decomposition is not just a theoretical exercise; it has practical applications in low-rank neural networks, image compression, and movie recommender systems, to name a few. The book highlights these applications, showcasing how decomposition techniques can be leveraged to solve real-world problems.

This book offers the first systematic comparisons among 
different variants of matrix decomposition. 
While method comparisons have been conducted in fields like linear modeling \citep{searle2016linear},
a comprehensive methodological analysis of matrix decomposition is lacking.
Many papers introducing new matrix decomposition applications fail to offer a thorough comparison with different decomposition approaches, such as those for neural network morphism and network architecture search.
We provide an overview of various approaches in the literature, ranging from fundamental ones covered in linear algebra courses to more complex methods for theoretical analysis. 
The vignettes in this book provide ample evidence regarding the benefits of harnessing the decompositional nature of the acquired matrix.
Figures~\ref{fig:matrix-decom-world-picture} and \ref{fig:matrix-decom-world-picture2}
provide a visual representation of the relationships and distinctions between each decomposition method.
Throughout the book, each chapter includes a section of problem sets that readers may choose to skip. However, there may occasionally be exercises following the introduction of certain definitions. These exercises are closely related to the concepts just introduced and should not be omitted.

\paragraph{Notation and preliminaries.} In the rest of this section, we provide a brief review of fundamental concepts in linear algebra, calculus, and optimization that will help us appreciate the structural properties of matrix factorization problems later in the book.  
For the rest of the important concepts, we define and discuss them as per need for clarity.
This chapter is not intended to be, by any means, a comprehensive treatment of these subjects, and  interested readers who wish to explore these subject in greater depth are encouraged to refer to advanced linear algebra and calculus books.
Throughout this text, unless stated otherwise, we simplify matters by considering only matrices with real entries. 
Additionally, unless otherwise noted, the eigenvalues of the matrices discussed are also real.

\section{Linear Algebra}
\paragraph{The vector space $\real^n$ and matrix space $\real^{m\times n}$.}
The vector space $\real^n$ is the set of $n$-dimensional column vectors with real components. 
Throughout the book, our primary focus will be on problems constrained within the $\real^n$ vector space. However, in a few instances, we will also explore other vector spaces, such as the  nonnegative vector space.
Similarly, the matrix space $\real^{m\times n}$ is the set of all real-valued matrices of order $m\times n$.

In all cases, scalars will be denoted in a non-bold font (script letters), possibly with subscripts (e.g., $a$, $\alpha$, $\alpha_i$). We will use \textbf{boldface} lowercase letters, possibly with subscripts, to denote vectors (e.g., $\bmu$, $\bx$, $\bx_n$, $\bz$), and
\textbf{boldface} uppercase letters, possibly with subscripts, to denote matrices (e.g., $\bA$, $\bL_j$). 
The $i$-th element of a vector $\bz$ will be denoted by $z_i$ in the non-bold font (or in rare cases as $\bz_i$ in bold font).
In the meantime, the \textit{normal fonts or plain typefaces} of scalars denote  \textbf{random variables} (e.g., $\textnormal{a}$ and $\textnormal{b}_1$ are random variables, while italics $a$ and $b_1$ are scalars); 
the normal fonts of \textbf{boldface} lowercase letters, possibly with subscripts, denote \textbf{random vectors} (e.g., $\rva$ and $\rvb_1$ are random vectors, while italics $\ba$ and $\bb_1$ are vectors); 
and the normal fonts of \textbf{boldface} uppercase letters, possibly with subscripts, denote \textbf{random matrices} (e.g., $\rmA$ and $\rmB_1$ are random matrices, while italics $\bA$ and $\bB_1$ are matrices).

The $n$-th element in a sequence is denoted by a superscript in parentheses. 
For instance, $\bA^{(n)}$ represents the $n$-th matrix in a sequence, and $\ba^{(k)}$ represents the $k$-th vector in a sequence.

Subarrays are formed when a subset of the indices is fixed.
\textit{The $i$-th row and $j$-th column value of matrix $\bA$ (referred to as entry ($i,j$) of $\bA$) will be denoted by $\bA_{ij}$ if block submatrices are involved, or as $a_{ij}$ alternatively if block submatrices are not involved} (in this case, $\bA\in\real^{m\times n}$ can be denoted as $\bA=\big\{a_{ij}\big\}_{i,j=1}^{m,n}=[a_{ij}]$). 
Furthermore, it will be helpful to utilize the \textbf{Matlab-style notation}: the submatrix of  matrix $\bA$ spanning from the $i$-th row to the $j$-th row and the $k$-th column to the $m$-th column  will be denoted by $\bA_{i:j,k:m}$ or $\bA[i:j,k:m]$. 
A colon is used to indicate all elements of a dimension, e.g., $\bA_{:,k:m}=\bA[:,k:m]$ denotes the $k$-th column to the $m$-th column of  matrix $\bA$, and $\bA_{:,k}=\bA[:,k]$ denotes the $k$-th column of $\bA$. Alternatively, the $k$-th column of matrix $\bA$ may be denoted more compactly as $\ba_k$. 
\index{Matlab-style notation}

When the index is not continuous, given ordered subindex sets $I$ and $J$, $\bA[I, J]$ denotes the submatrix of $\bA$ obtained by extracting the rows and columns of $\bA$ indexed by $I$ and $J$, respectively; and $\bA[:, J]$ denotes the submatrix of $\bA$ obtained by extracting the columns of $\bA$ indexed by $J$, where again the colon operator implies all indices, and  the $[:, J]$ syntax in this expression selects all rows from $\bA$ and only the columns specified by the indices in $J$.

\begin{definition}[Matlab Notation]\label{definition:matlabnotation}
Let $\bA=[a_{ij}]\in \real^{m\times n}$, and let $I=[i_1, i_2, \ldots, i_k]$ and $J=[j_1, j_2, \ldots, j_l]$ be two index vectors. Then $\bA[I,J]$ denotes the $k\times l$ submatrix
$$
\bA[I,J]=\bA_{I,J}
\triangleq
\begin{bmatrix}
a_{i_1,j_1} & a_{i_1,j_2} &\ldots & a_{i_1,j_l}\\
a_{i_2,j_1} & a_{i_2,j_2} &\ldots & a_{i_2,j_l}\\
\vdots & \vdots&\ddots & \vdots\\
a_{i_k,j_1} & a_{i_k,j_2} &\ldots & a_{i_k,j_l}\\
\end{bmatrix}.
$$
Whilst, $\bA[I,:]=\bA_{I,:}$ denotes a $k\times n$ submatrix, and $\bA[:,J]=\bA_{:,J}$ denotes a $m\times l$ submatrix analogously.

We note that it does not matter whether the index vectors $I$ and $J$ are row vectors or column vectors. 
What's important is which axis they index (either rows  or columns of $\bA$).
Note that the ranges of the indices are given as folows:
$$
\left\{
\begin{aligned}
0&\leq \min(I) \leq \max(I)\leq m;\\
0&\leq \min(J) \leq \max(J)\leq n.
\end{aligned}
\right.
$$
\end{definition}

In all cases, vectors are formulated in a column rather than in a row (except the index vectors $I$ and $J$ introduced above). 
A row vector will be denoted by the transpose of a column vector, such as $\ba^\top$. 
A column vector with specific values is separated by the semicolon symbol $``;"$, for instance, $\bx=[1;2;3]$ is a column vector in $\real^3$. 
Similarly, a  row vector with specific values is split by the comma symbol $``,"$, e.g., $\by=[1,2,3]$ is a row vector with 3 values. 
Additionally, a column vector can also be denoted by the transpose of a row vector, e.g., $\by=[1,2,3]^\top$ is also a column vector.

The transpose of a matrix $\bA$ will be represented as $\bA^\top$, and its inverse will be denoted by $\bA^{-1}$. 
We follow the standard convention by denoting  the $p \times p$ identity matrix as $\bI_p$ (or simply as $\bI$ when the size can be determined from the context). A vector or matrix of all zeros will be denoted by a \textbf{boldface} zero, $\bzero$, whose size should be clear from  context; or we denote $\bzero_p$ as the vector of all zeros with $p$ entries and $\bzero_{m\times n}$ as the matrix of all zeros with dimensions $m\times n$.
Similarly, a vector or matrix of all ones will be denoted by a \textbf{boldface} one $\bone$, whose size is clear from  context; or we denote $\bone_p$ as the vector of all ones with $p$ entries and $\bone_{m\times n}$ as the matrix of all ones with dimensions $m\times n$. The $p\times p$ all-ones matrix can be denoted by $\bA=\bone_p\bone_p^\top$.
We will frequently suppress the subscripts of these matrices when the dimensions are clear from the context.
When it comes to diagonal matrices, we employ $\bD=\diag([d_1,d_2,\ldots,d_n])$ to abbreviate the $n$-by-$n$ diagonal matrix with  elements $d_1,d_2,\ldots,d_n$ on its \textit{main diagonal} (or simply diagonal, i.e., the diagonal running from the top left corner to the bottom right corner of the matrix).
For a \textit{block-diagonal matrix} $\bA=\tiny\begin{bmatrix}
\bA_{11} & & \\
& \ddots & \\
& & \bA_{nn}
\end{bmatrix}$, we denote it by $\bA=\diag(\bA_{11}, \bA_{22}, \ldots,\bA_{nn})$ or $\bA=\bA_{11}\oplus \bA_{22}\oplus \ldots\oplus\bA_{nn}$, where the latter is called the \textit{direct sum} notation.

\index{Direct sum}

We will  use $\be_1, \be_2, \ldots, \be_n$ to represent the \textit{(unit) standard basis} of $\real^n$, where $\be_i$ is the vector whose $i$-th component is one while all the others are zero (a.k.a., the $i$-th \textit{canonical vector} of the Euclidean space $\real^n$).

\begin{definition}[Nonnegative Orthant, Positive Orthant, and Unit-Simplex]\label{definition:simplex}
The \textit{nonnegative orthant} is a subset of $\real^n$ that consists of all vectors in $\real^n$ with nonnegative components and is denoted by $\real_+^n$: 
$$
\real_+^n=\{[x_1, x_2, \ldots, x_n]^\top\in\real^n: x_1, x_2, \ldots, x_n\geq 0\}.
$$
Similarly, the \textit{positive orthant} consists of all vectors in $\real^n$ with positive components and is denoted by $\real_{++}^n$:
$$
\real_{++}^n=\{[x_1, x_2, \ldots, x_n]^\top\in\real^n: x_1, x_2, \ldots, x_n> 0\}.
$$
The \textit{unit-simplex} is a subset of $\real^n$ comprising all nonnegative vectors whose sum is 1:
$$
\Delta_n=\{\bx=[x_1, x_2, \ldots, x_n]^\top\in\real^n:  x_1, x_2, \ldots, x_n\geq  0,
\,\, \sum_{i=1}^{n}x_i=1 \}.
$$
\end{definition}

%

\index{Eigenvalue}
\index{Eigenvector}
\begin{definition}[Eigenvalue, Eigenvector]
Given any vector space $\sF$ and any linear map $\bA: \sF \rightarrow \sF$ (or simply a real matrix $\bA\in\real^{n\times n}$), a scalar $\lambda \in \sK$ is called a \textit{(right) eigenvalue, or proper value, or characteristic value} of $\bA$, if there exists some nonzero vector $\bu \in \sF$ such that
\begin{equation*}
\bA \bu = \lambda \bu.
\end{equation*}
And $\bu$ is called a \textit{(right) eigenvector} of $\bA$ associated with $\lambda$.

On the other hand, $\kappa$ is referred to as a \textit{left eigenvalue} if there exists some nonzero vector $\bv\in \sF$ such that 
$$
\bv^\top\bA = \kappa \bv^\top.
$$
And $\bv$ is called a \textit{(left) eigenvector} of $\bA$ associated with $\kappa$.

When it is clear from the context, we will simply use the term ``eigenvalue/eigenvector" instead of ``right eigenvalue/eigenvector."

\end{definition}
In simple terms, an eigenvector $\bu$ of a matrix $\bA$ represents a direction that remains unchanged when transformed into the coordinate system defined by the columns of $\bA$ (see Section~\ref{section:coordinate-transformation} for more details on coordinate transformations).
In fact, real-valued matrices can have complex eigenvalues. However, all the eigenvalues of symmetric matrices are real (see Theorem~\ref{theorem:spectral_theorem}).

\index{Spectrum}
\index{Spectral radius}
\begin{definition}[Spectrum and Spectral Radius]\label{definition:spectrum}
The set of all eigenvalues of $\bA$ is called the \textit{spectrum} of $\bA$ and is denoted by $\Lambda(\bA)$. 
The set of all eigenvalues of $\bA$ including their algebraic multiplicity (Definition~\ref{definition:eigen_multipli}) is called the \textit{multispectrum} of $\bA$ and is denoted by $\widehat{\Lambda}(\bA)$. It holds that $\Lambda(\bA)\subseteq \widehat{\Lambda}(\bA)$.
The largest magnitude of the eigenvalues is known as the spectral radius $\rho(\bA)$:
$
\rho(\bA) = \mathop{\max}_{\lambda\in \Lambda(\bA)}  \abs{\lambda}.
$
\end{definition}

%

\begin{exercise}[Orthogonality in Eigenvectors]
Suppose $\bA\bu=\lambda\bu$ and $\bv^\top\bA=\kappa\bv^\top$ with $\lambda\neq \kappa$. Show that $\bv^\top\bu=0$.
\end{exercise}

Moreover, the tuple $(\lambda, \bu)$ above is referred to as an \textit{eigenpair}. Intuitively, these definitions imply that multiplying the matrix $\bA$ by the vector $\bu$ results in a new vector that is in the same direction as $\bu$, but only scaled by a factor $\lambda$. For any eigenvector $\bu$, we can scale it by a scalar $s$ such that $s\bu$ remains an eigenvector of $\bA$; and for this reason, we call the eigenvector an eigenvector of $\bA$ associated with the eigenvalue $\lambda$. 
To avoid ambiguity, it is customary to assume that the eigenvector is normalized to have length one and the first entry is positive (or negative) since both $\bu$ and $-\bu$ are eigenvectors.

In linear algebra, every vector space has a basis, and any vector within the space can be expressed as  a linear combination of the basis vectors. 
We then define the span and dimension of a subspace via the basis.

\index{Subspace}
\begin{definition}[Subspace]
	A nonempty subset $\mathcal{V}$ of $\real^n$ is called a subspace if $x\ba+y\ba\in \mathcal{V}$ for every $\ba,\bb\in \mathcal{V}$ and every $x,y\in \real$.
\end{definition}

\index{Span}
\begin{definition}[Span]
	If every vector $\bv$ in the subspace $\mathcal{V}$ can be expressed as a linear combination of $\{\ba_1, \ba_2, \ldots,$ $\ba_m\}$, then $\{\ba_1, \ba_2, \ldots, \ba_m\}$ is said to span $\mathcal{V}$.
\end{definition}

\index{Linearly independent}
In this context, we will often use the idea of the linear independence of a set of vectors. Two equivalent definitions are given as follows.
\begin{definition}[Linearly Independent]
A set of vectors $\{\ba_1, \ba_2, \ldots, \ba_m\}$ is called linearly independent if there is no combination that  can yield $x_1\ba_1+x_2\ba_2+\ldots+x_m\ba_m=\bzero $ unless all $x_i$'s are equal to  zero. An equivalent definition is that $\ba_1\neq \bzero$, and for every $k>1$, the vector $\ba_k$ does not belong to the span of $\{\ba_1, \ba_2, \ldots, \ba_{k-1}\}$.
\end{definition}

\index{Basis}
\index{Dimension}
\begin{definition}[Basis and Dimension]
A set of vectors $\{\ba_1, \ba_2, \ldots, \ba_m\}$ is called a basis of $\mathcal{V}$ if they are linearly independent and span $\mathcal{V}$. Every basis of a given subspace contains  the same number of vectors, and this number is referred to as  the dimension of the subspace $\mathcal{V}$. By convention, the dimension of the  subspace $\{\bzero\}$ is zero. 
Additionally, every subspace with a  nonzero dimension has an orthogonal basis; that is, a basis can be chosen such that all its vectors are mutually orthogonal.
\end{definition}

\index{Column space}
\begin{definition}[Column Space (Range)]
Let $\bA\in\real^{m\times n}$ be any real matrix. Then, the \textit{column space (or range)} of $\bA$ is defined as the set spanned by its columns:
\begin{equation*}
\mathcal{C} (\bA) = \{ \by\in \mathbb{R}^m: \exists\, \bx \in \mathbb{R}^n, \, \by = \bA \bx \}.
\end{equation*}
And the \textit{row space} of $\bA$ is the set spanned by its rows, which is equivalent to the column space of $\bA^\top$:
\begin{equation*}
\mathcal{C} (\bA^\top) = \{ \bx\in \mathbb{R}^n: \exists\, \by \in \mathbb{R}^m, \, \bx = \bA^\top \by \}.
\end{equation*}
\end{definition}

\index{Null space}
\begin{definition}[Null Space (Nullspace, Kernel)]
Let $\bA\in\real^{m\times n}$ be any real matrix. Then, the \textit{null space (or kernel, or nullspace)} of $\bA$ is defined as the  set:
\begin{equation*}
\nspace (\bA) = \{\by \in \mathbb{R}^n:  \, \bA \by = \bzero \}.
\end{equation*}
And the null space of $\bA^\top$ is defined as 	\begin{equation*}
\nspace (\bA^\top) = \{\bx \in \mathbb{R}^m:  \, \bA^\top \bx = \bzero \}.
\end{equation*}
\end{definition}

Both the column space of $\bA$ and the null space of $\bA^\top$ are subspaces of $\real^m$. In fact, every vector in $\nspace(\bA^\top)$ is orthogonal  to the vectors in $\cspace(\bA)$, and vice versa. Similarly, every vector in $\nspace(\bA)$ is also perpendicular to the vectors in $\cspace(\bA^\top)$, and vice versa.

\index{Rank}
The rank of a matrix provides information about the linear independence of the matrix's rows or columns.
\begin{definition}[Rank]
The rank of a matrix $\bA\in \real^{m\times n}$ is the dimension of its column space. 
That is, the rank of $\bA$ is equal to the maximum number of linearly independent columns of $\bA$, and it is also the maximum number of linearly independent rows of $\bA$. The matrix $\bA$ and its transpose $\bA^\top$ have the same rank. 
A matrix $\bA$ is said to have full rank if its rank equals $\min\{m,n\}$.
Specifically, given a vector $\bu \in \real^m$ and a vector $\bv \in \real^n$, then the $m\times n$ matrix $\bu\bv^\top$ obtained by the outer product of these vectors has a  rank of 1. 
In summary, the rank of a matrix is equal to:
\begin{itemize}
\item the number of linearly independent columns;
\item the number of linearly independent rows.
\end{itemize}
And remarkably, these two quantities are always equal (see Appendix~\ref{append:row-equal-column}).
\end{definition}

\index{Direct sum}
\begin{exercise}[Rank of Direct Sum]
Show that the rank of a direct sum is the sum of the ranks of the summands: $\rank(\bA)=\rank(\bA_{11})+\rank(\bA_{22})+\ldots+\rank(\bA_{nn})$  if $\bA=\bA_{11}\oplus \bA_{22}\oplus \ldots\oplus\bA_{nn}$.
\end{exercise}

\index{Sylvester's inequality}
\index{Frobenius' inequality}
We may prove some of the rank properties when they are relevant to the topics discussed in the book. Here are some of these inequalities collected in the following remark.
\begin{remark}[Rank Properties]\label{remark:rank_prop}
We have the following rank properties:
\begin{itemize}
\item $\rank(\bA)=\rank(\bB\bA)=\rank(\bA\bC)=\rank(\bB\bA\bC)$, given $\bA\in\real^{m\times n}$, and nonsingular $\bB\in\real^{m\times m}$ and nonsingular $\bC\in\real^{n\times n}$.
\item  When one or more rows and/or columns are removed from a matrix, the rank of the resulting submatrix will not exceed the rank of the original matrix.
\item \textit{Sylvester's inequality.} $\rank(\bA)+\rank(\bB)-r \leq \rank(\bA\bB)\leq \min\{\rank(\bA), \rank(\bB)\}$ if $\bA\in\real^{m\times r}$ and $\bB\in\real^{r\times n}$.
\item \textit{Frobenius' inequality.} $\rank(\bA\bB)+\rank(\bB\bC)\leq \rank(\bB)+\rank(\bA\bB\bC)$ given appropriate matrices $\bA,\bB$, and $\bC$.
\item $\abs{\rank(\bA) -\rank(\bB)} \leq \rank(\bA+\bB) \leq \rank(\bA)+\rank(\bB)$. 
\end{itemize}
\end{remark}

In certain cases, we may also consider a special type of rank known as the \textit{$k$-Rank or Kruskal rank}.  The $k$-rank appears in the formulation
of the famous Kruskal condition for CANDECOMP-PARAFAC uniqueness \citep{carroll1970analysis,harshman1970foundations, de2008decompositionspart}.
\begin{definition}[$k$-Rank or Kruskal Rank]\label{definition:krus_rk}
The \textit{Kruskal rank or $k$-rank} of a matrix $\bA$, denoted by $\rank_k(\bA)$ or $k_{\bA}$, is the largest  number $r$ such that any set of $r$ columns of $\bA$ is linearly independent.
\end{definition}
When the matrix has a $k$-rank of $\sigma$, this means that no column in any subset of size $\sigma$ can be expressed as a linear combination of the others in that subset.
Clearly, a rank-$r$ matrix can have a $k$-rank of 1 if there are two identical columns in the matrix. 
Thus, the  $k$-rank provides insight into a specific type of structural redundancy within the matrix.
Specifically, 
\begin{itemize}
\item When the matrix has full column rank,  adding columns to the matrix may increase or decrease the $k$-rank.
\item When the matrix does not have full column rank,  adding columns to the matrix can never increase the $k$-rank.
\end{itemize}

\begin{lemma}[Premultiplying with Nonsingular]\label{lemma:left_mul_krank}
Premultiplying a matrix by a nonsingular matrix can be understood as applying the same nonsingular transformation to each column of the original matrix. This transformation preserves the linear independence of the columns; thus, the linear dependence relationships among the columns remain unchanged.
Therefore, premultiplying a matrix with a nonsingular matrix will not change either the rank or $k$-rank of  the result.
\end{lemma}

Similarly, we can also consider the $k^\prime$-rank of partitioned matrices.
\begin{definition}[$k^\prime$-Rank]
The \textit{$k^\prime$-rank} of a  partitioned matrix $\bA$, denoted by $\rank_{k^\prime}(\bA)$ or $k_{\bA}^\prime$, is the largest  number $r$ such that any set of $r$
submatrices of $\bA$ yields a set of linearly independent columns.
\end{definition}

\index{Orthogonal complement}
\begin{definition}[Orthogonal Complement in General]
The orthogonal complement $\mathcalV^\perp\subseteq \real^m$ of a subspace $\mathcalV\subseteq\real^m$ consists of all vectors that are perpendicular to $\mathcalV$. That is,
$$
\mathcalV^\perp = \{\bv\in \real^m : \bv^\top\bu=0, \,\,\, \forall \bu\in \mathcalV  \}.
$$
These two subspaces are mutually exclusive yet collectively span the entire space $\real^m$.
The dimensions of $\mathcalV$ and $\mathcalV^\perp$ sum up to the dimension of the entire space: $\dim(\mathcalV)+\dim(\mathcalV^\perp)=m$. Furthermore, it holds that $(\mathcalV^\perp)^\perp=\mathcalV$.
\end{definition}

\paragraph{Subspace intersection.} Consider general subspaces that are not necessarily disjoint and do not necessarily span the entire space, we have the following result. 
Given a subspace $\mathcalV$ and its two subspaces $\mathcalV_1$ and $\mathcalV_2$, their dimensions follow:
\begin{equation}
\begin{aligned}
\dim(\mathcalV_1&\cap \mathcalV_2) + \dim(\mathcalV_1+ \mathcalV_2)=\dim(\mathcalV_1)+\dim(\mathcalV_2)\\
&\implies \dim(\mathcalV_1\cap \mathcalV_2) \geq \dim(\mathcalV_1)+\dim(\mathcalV_2) - \dim(\mathcalV).
\end{aligned}
\end{equation}
\begin{exercise}
Given two matrices $\bA\in\real^{n\times p}$ and $\bB\in\real^{n\times q}$, show that 
$$
\begin{aligned}
\cspace(\bA)+\cspace(\bB) &= \cspace([\bA,\bB]); 
\qquad
\nspace(\bA)\cap \nspace(\bB) &=\nspace\big(\footnotesize\begin{bmatrix}
	\bA\\
	\bB
\end{bmatrix}\normalsize\big).
\end{aligned}
$$
\end{exercise}

\index{Orthogonal complement}
Consider the preceding definition, for example, we can explicitly define the orthogonal complement of the column space in the following definition.
\begin{definition}[Orthogonal Complement of Column Space]
Let $\bA$ be an $m \times n$ real matrix. Then, the orthogonal complement of $\mathcalC(\bA)$, denoted by $\mathcal{C}^{\bot}(\bA)$, is the subspace defined as:
\begin{equation*}
\begin{aligned}
\mathcal{C}^{\bot}(\bA) &= \{\by\in \mathbb{R}^m: \, \by^\top \bA \bx=\bzero, \, \forall \bx \in \mathbb{R}^n \} 
=\{\by\in \mathbb{R}^m: \, \by^\top \bv = \bzero, \, \forall \bv \in \mathcalC(\bA) \}.
\end{aligned}
\end{equation*}
\end{definition}
We can then identify the four fundamental spaces associated with any matrix $\bA\in \real^{m\times n}$ of rank $r$:
\begin{itemize}
\item $\cspace(\bA)$: Column space of $\bA$, i.e., linear combinations of columns with dimension $r$;

\item $\nspace(\bA)$: Null space of $\bA$, i.e., all $\bx$ satisfying $\bA\bx=\bzero$ with dimension $n-r$;

\item  $\cspace(\bA^\top)$: Row space of $\bA$, i.e., linear combinations of rows with dimension $r$;

\item  $\nspace(\bA^\top)$: Left null space of $\bA$, i.e., all $\by$ satisfying $\bA^\top \by=\bzero$ with dimension $m-r$.
\end{itemize}
Furthermore, $\nspace(\bA)$ is the orthogonal complement of $\cspace(\bA^\top)$, and $\cspace(\bA)$ is the orthogonal complement of $\nspace(\bA^\top)$. The proof can be found in Appendix~\ref{appendix:fundamental-rank-nullity}.
In many text, the dimension of the null space is called the \textit{defect or nullity} for the square matrix  $\bA\in\real^{n\times n}$:
\begin{equation}
\defect(\bA) = \dim(\nspace(\bA)).
\end{equation}
\paragraph{Complementary nullities.}
Consider a row index set $I$ and a column index set $J$ with cardinality $\abs{I}=i$ and $\abs{J}=j$. Their complementary sets are $I^C=\{1,2,\ldots,n\}\backslash I$ and $J^C=\{1,2,\ldots,n\}\backslash J$, respectively.
Then, it follows that 
\begin{equation}
\begin{aligned}
\defect(\bA[I,J]) &= \defect(\bA^{-1}[J^C, I^C]);\\
\rank(\bA[I,J]) &= \rank(\bA^{-1}[J^C, I^C]) +i+j-n.
\end{aligned}
\end{equation}

\paragraph{Fredholm alternative theorem.}
Consider the linear system $\bA\bx=\bb$, where $\bA$ and $\bb$ are given. The linear system has a solution if and only if $\bb$ is in the column space of $\bA$: $\bb\in\cspace(\bA)$. 
This is also equivalent to saying that $\bb$ is orthogonal to every vector $\by\in\nspace(\bA^\top)$: $\bb^\top\by=0$ if $\by\in\nspace(\bA^\top)$.
This is known as the \textit{Fredholm alternative theorem}.

\index{Fredholm alternative theorem}
\index{Fundamental subspace}

\index{Orthogonal matrix}
\begin{definition}[Orthogonal Matrix, Semi-Orthogonal Matrix]
A real square matrix $\bQ\in\real^{n\times n}$ is called an \textit{orthogonal matrix} if the inverse of $\bQ$ is equal to its  transpose, that is, $\bQ^{-1}=\bQ^\top$ and $\bQ\bQ^\top = \bQ^\top\bQ = \bI$.~\footnote{Similarly, $\bU\in\complex^{n\times n}$ is called unitary if $\bU\bU^*=\bU^*\bU=\bI$, i.e., the inverse is equal to its conjugate transpose.} 
In other words, suppose $\bQ=[\bq_1, \bq_2, \ldots, \bq_n]$, where $\bq_i \in \real^n$ for all $i \in \{1, 2, \ldots, n\}$. Then, $\bq_i^\top \bq_j = \delta(i,j)$, where $\delta(i,j)$ is the Kronecker delta function. 
For any vector $\bx$, the orthogonal matrix will preserve the length: $\normtwo{\bQ\bx}= \normtwo{\bx}$ (Exercise~\ref{exercise:orthogo_ell2}).
Note that, since the orthogonal matrix $\bQ$ contains unit-length columns, the columns are mutually orthonormal. However, the term \textit{orthonormal matrix} is \textbf{not} used for historical reasons.

If $\bQ$ contains only $\gamma$ of these columns with $\gamma<n$, then $\bQ^\top\bQ = \bI_\gamma$ stills holds, where $\bI_\gamma$ is the $\gamma\times \gamma$ identity matrix. 
But $\bQ\bQ^\top=\bI$ will not hold. 
In this case, $\bQ$ is called \textit{semi-orthogonal}.
\end{definition}

We should also notice that the rows of an orthogonal matrix are also mutually orthonormal according to the definition provided above.
The transformation of a vector $\bx$ under an orthogonal matrix has the property: $\norm{\bQ\bx}_2=\norm{\bx}_2$ (where $\norm{\cdot}_2$ denotes the norm of a vector; see Definition~\ref{definition:vec_l2_norm}). 
That is, the length remains unchanged under an orthogonal transformation, a property often referred to as  \textit{Euclidean isometry}.

\index{Permutation matrix}
\begin{definition}[Permutation Matrix]\label{definition:permutation-matrix}
A permutation matrix $\bP\in \real^{n\times n}$ is a square binary matrix that contains exactly one entry of 1 in each row and each column, with 0's elsewhere. By definition, $\bP\bP^\top=\bI$ and $\bP$ is orthogonal. 

\paragraph{Row perspective.} The permutation matrix $\bP$ rearranges the rows of the identity matrix $\bI$ in any desired order, and this order determines the sequence of row permutations.
In other words, to permute the rows of matrix $\bA$, we simply multiply on the left by $\bP\bA$. 

\paragraph{Column perspective.} Alternatively, the permutation matrix $\bP$ rearranges the columns of the identity matrix $\bI$ in any desired order, and this order dictates the sequence of column permutations. 
To achieve the column permutation of $\bA$, we multiply on the right by  $\bA\bP$.

\paragraph{Reversal matrix.} A special permutation matrix $\bP=[\be_n, \be_{n-1}, \ldots, \be_1]\in\real^{n\times n}$, where $\be_i$ represents the $i$-th unit basis vector in $\real^n$, is called the \textit{reversal matrix}. The reversal matrix satisfies that $\bP^2=\bI$ and $\bP^\top=\bP$ (orthogonality).
It  can be shown that if $\bR$ is upper triangular, then $\bL=\bP\bR\bP$ is lower triangular; the diagonal values of $\bL$ are those of $\bR$, with the order reversed. 
If $\bQ$ is orthogonal, then both $\bQ\bP$ and $\bP\bQ$ are also orthogonal.
\end{definition}

For example, let $\bx\in\real^n$. Then, $\bP_d=[\be_2, \be_3, \ldots, \be_n, \be_1]$ is called a \textit{downshift permutation} such that $\bP_d\bx=[x_2,x_3,\ldots,x_n,x_1]^\top$; $\bP_e=[\be_n, \be_{n-1}, \ldots, \be_1]$ is called  an \textit{exchange permutation} such that 
$\bP_e\bx=[x_n,x_{n-1},\ldots,x_1]^\top$.
If further $n=pr$, i.e., $n$ can be divided by $p$, the \textit{mod-$p$ perfect shuffle permutation} is defined as
\begin{equation}
\begin{aligned}
\bP_{p,r}&=\bI_n[(1:r:n), (2:r:n), \ldots, (r:r:n)];\\
\bP_{p,r}^\top &= \bI_n[(1:p:n), (2:p:n), \ldots, (p:p:n)],
\end{aligned}
\end{equation}
such  that $\bP_{p,r}\bx=[(x_1, x_{1+r}, \ldots, x_{pr-r}), (x_2, x_{2+r}, \ldots, x_{pr-r+1}), \ldots]$.

The permutation matrix $\bP$ can be more efficiently represented using a vector $J \in \integer_{++}^n$ of indices such that $\bP = \bI[:, J]$, where $\bI$ is the $n\times n$ identity matrix. Notably, the elements in vector $J$ sum to $1+2+\ldots+n= \frac{n^2+n}{2}$.

\begin{example}[Permutation]
Suppose
$$\bA=\begin{bmatrix}
1 & 2&3\\
4&5&6\\
7&8&9
\end{bmatrix}
\qquad \text{and} \qquad
\bP=\begin{bmatrix}
	&1&\\
	&&1\\
	1&&
\end{bmatrix}.
$$
The row  and column permutations are given by 
$$
\bP\bA = \begin{bmatrix}
4&5&6\\
7&8&9\\
1 & 2&3\\
\end{bmatrix}
\qquad \text{and} \qquad
\bA\bP = \begin{bmatrix}
	3 & 1 & 2 \\
	6 & 4 & 5\\
	9 & 7 & 8
\end{bmatrix},
$$
respectively, 
where the order of the rows of $\bA$ appearing in $\bP\bA$ corresponds to the order of the rows of $\bI$ in $\bP$,
and
where the order of the columns of $\bA$ appearing in $\bA\bP$ corresponds to the order of the columns of $\bI$ in $\bP$. 
\end{example}

\index{Generalized permutation matrix}
\paragraph{Generalized permutation matrices.}
Similarly, a matrix $\bP\in\real^{n\times n}$ is called a \textit{generalized permutation matrix}
if all its elements are either 0, 1, or $-1$, and each row and each column has exactly one nonzero element.

\begin{exercise}
Suppose $\bx\in[5,2,3,9]^\top$, and denote the reordered nondecreasing vector of $\bx$ by $\bx^{\uparrow}=[2,3,5,9]^\top$. Show that there exists a permutation matrix $\bP$ such that $\bP\bx=\bx^{\uparrow}$, and there exists a generalized permutation matrix $\bQ$ such that $\bQ\bx=\abs{\bx}^{\uparrow}$.
\end{exercise}

\index{Selection matrix}
\begin{definition}[Selection Matrix]\label{definition:selection-matrix}
A selection matrix $\bS\in \real^{n\times n}$ is a square diagonal matrix with diagonal entries that are either 1 or 0.  
The entries with 1 indicate the rows or columns that will be selected.

\paragraph{Row perspective.} The selection matrix $\bS$ contains the rows of the identity matrix $\bI$ if we intend to select the corresponding rows; otherwise, it masks the rows in the identity matrix $\bI$ with zeros.
To select the rows of matrix $\bA$, we perform a left multiplication by $\bS\bA$. 

\paragraph{Column perspective.} Alternatively, the selection matrix $\bS$ contains the columns of the identity matrix $\bI$ if we want to select the corresponding columns; otherwise, it masks  the columns in the identity matrix $\bI$ with zeros. Therefore, the column selection of $\bA$ is equivalent to postmultiplying by $\bA\bS$.
\end{definition}

\begin{example}[Selection and Permutation]
Suppose, 
$$\bA=\begin{bmatrix}
1 & 2&3\\
4&5&6\\
7&8&9
\end{bmatrix}
\qquad \text{and} \qquad
\bS=\begin{bmatrix}
1&&\\
& 0&\\
&&1
\end{bmatrix}.
$$
The row and columns selections are  given by 
$$
\bS\bA = \begin{bmatrix}
1&2&3\\
0&0&0\\
7 & 8&9\\
\end{bmatrix}
\qquad \text{and} \qquad
\bA\bS = \begin{bmatrix}
	1& 0 & 3 \\
	4 & 0 & 6\\
	7 & 0 & 9
\end{bmatrix},
$$
respectively,
where the rows of $\bA$ appearing in $\bS\bA$ correspond to the row entries of $\bS$, 
and  the columns of $\bA$ appearing in $\bA\bS$ correspond to  the column entries of $\bS$.
Suppose we want to reorder the selected rows or columns into the upper-left corner of the final matrix, we can construct a permutation as follows
$$
\bP=\begin{bmatrix}
1&&\\
& & 1\\
&1&
\end{bmatrix}, 
$$
such that 
$$
\bP\bS\bA = \begin{bmatrix}
1&2&3\\
7 & 8&9\\
0&0&0\\
\end{bmatrix}
\qquad \text{and} \qquad
\bA\bS\bP = \begin{bmatrix}
	1 & 3& 0 \\
	4  & 6& 0\\
	7  & 9& 0
\end{bmatrix}.
$$
This technique is essential for certain mathematical proofs, such as the properties of positive definite matrices in  Lemma~\ref{lemma:pd-more-properties}. 
\end{example}

\index{Normal matrix}
\index{Hermitian matrix}
\index{Orthogonal matrix}
\index{Unitary matrix}
\index{Range-symmetric}
\index{Range-Hermitian}
For real matrices, we recall the following special matrices.
\begin{definition}[Normal and Other Special Real Matrices]\label{definition:speci_mat}
Give a real $m\times n$ matrix $\bA\in \real^{m\times n}$, the matrix is 
\begin{enumerate}
\item \textit{normal} if $\bA^\top\bA=\bA\bA^\top$;
\item \textit{symmetric} if $\bA^\top=\bA$ (with $m=n$), denoted by $\bA\in\sS^n$~\footnote{Similarly, if $\bA$ is positive semidefinite, then $\bA\in\sS_+^n$; and if $\bA$ is positive definite, then $\bA\in\sS_{++}^n$. Obviously, the inclusion $\sS_{++}^n\subseteq \sS_{+}^n \subseteq \sS^n$ holds.};
\item \textit{range-symmetric} if $\cspace(\bA)=\cspace(\bA^\top)$  (with $m=n$);
\item \textit{skew-symmetric} if $\bA^\top=-\bA$ (with $m=n$);
\item \textit{orthogonal} if $\bA^\top\bA  =\bA\bA^\top = \bI_n$ (with $m=n$), denoted by $\bA\in\sO^n$;
\item $\bA$ and $\bB$ are said to commute if $\bA\bB=\bB\bA$ (with $\bA,\bB\in\real^{n\times n}$);
\item $\bA$ and $\bB$ are said to anticommute if $\bA\bB=-\bB\bA$ (with $\bA,\bB\in\real^{n\times n}$).
\end{enumerate}
\end{definition}
While most of our discussions will focus on real matrices, it is also beneficial to be aware of the definitions of special complex matrices.
\begin{definition}[Normal and Other Special Complex Matrices]\label{definition:complex_special}
	Give a complex  $m\times n$ matrix $\bA\in \complex^{m\times n}$, the matrix is 
\begin{enumerate}
\item \textit{normal} if $\bA^\ast\bA=\bA\bA^\ast$ (not to be confused with normal distributions in statistics, $^*$ denotes the conjugate transpose of a complex matrix);
\item \textit{Hermitian} if $\bA^\ast=\bA$ (with $m=n$);
\item \textit{range-Hermitian} if $\cspace(\bA)=\cspace(\bA^\ast)$ (with $m=n$);
\item \textit{skew-Hermitian} if $\bA^\ast=-\bA$ (with $m=n$);
\item \textit{unitary} if $\bA^\ast\bA  =\bA\bA^\ast = \bI_n$ (with $m=n$); 
\item  $\bA$ and $\bB$ are said to commute if $\bA\bB=\bB\bA$ (with $\bA,\bB\in\complex^{n\times n}$);
\item $\bA$ and $\bB$ are said to anticommute if $\bA\bB=-\bB\bA$ (with $\bA,\bB\in\complex^{n\times n}$).
\end{enumerate}
where the \textit{conjugate} $\overline{\bA}$ of $\bA$ is the $m\times n$ matrix $\overline{\bA}=[b_{ij}]$ such that $b_{ij}=\overline{a_{ij}}$ and $\bA^\ast=\overline{(\bA^\top)}$.
\paragraph{Observation.} The set of normal matrices includes all  Hermitian, skew-Hermitian, and unitary matrices. Additionally, it also contains other matrices, such as 
$
\scriptsize
\begin{bmatrix}
1 & -1 \\
1 & 1
\end{bmatrix}.
$
\end{definition}

\begin{exercise}[Left Eigenvector of Normal Matrices]\label{exercise:left_right_normal}
Suppose $\bA\in\complex^{n\times n}$ is normal. Show that $\bA$ and $\bA^*$ have the same eigenvectors ($\bA\bx=\lambda\bx\implies \bA^*\bx=\overline{\lambda}\bx$).
Based on this, show that the right eigenvector of a normal matrix is also a left eigenvector of the matrix corresponding to the same eigenvalue ($\bA\bx=\lambda\bx\implies \bx^*\bA=\lambda\bx^*$).
\end{exercise}

\begin{example}
Given matrices
$$
\bA=\begin{bmatrix}
1 & 0 \\
0 & 1
\end{bmatrix}, 
\gap 
\bB=\begin{bmatrix}
	0 & 1 \\
	1 & 0
\end{bmatrix}, 
\gap 
\bC=\begin{bmatrix}
	1 & 0 \\
	0 & -1
\end{bmatrix}.
$$
Then, $\bA$ and $\bB$ commute, and $\bB$ and $\bC$ anticommute.
\end{example}

\begin{remark}[Commutativity in Block-Diagonal Matrix]\label{remark:commute}
Suppose $\bD$ is a diagonal matrix: $\bD=\{\bD_{ij}\}_{i,j=1}^{p,p}$ where $\bD_{ij}=\bzero $ if $i\neq j$, and $\bD_{ii}=\lambda_i \bI_{n_i}$ for all $i\in\{1,2,\ldots,p\}$ ($\sum_{i=1}^{p} n_i=n$). Partition $\bA$ conformally with $\bD$ such that $\bA=\{\bA_{ij}\}_{i,j=1}^{p,p}$. 
Then, $\bA\bD=\bD\bA$ if and only if $\lambda_{i}\bA_{ij}=\bA_{ij}\lambda_j$ for all $i,j\in\{1,2,\ldots,p\}$. This also indicates that $\bA$ is a block-diagonal matrix: $\bA=\diag(\bA_{11}, \bA_{22}, \ldots, \bA_{pp})$.
\end{remark}
\begin{exercise}[Properties of Commutativity]\label{exercise:prop_comm}
Let $\bB\in\real^{n\times n}$ be a matrix that commutes with $\bA\in\real^{n\times n}$. Show that $\bA^k\bB=\bB\bA^k$ for $k=1,2,\ldots$. 
This indicates if $\bC=p(\bA)=a_0\bI+a_1\bA+a_2\bA^2+\ldots +a_n\bA^n$, i.e., a polynomial of $\bA$, then $\bC$ also commutes with $\bB$ if $\bA\bB=\bB\bA$.
\textit{Hint: prove by induction.}
\end{exercise}

\index{Nonnegative matrix}
\begin{definition}[Positive, Nonnegative Matrix]\label{definition:pox_nonnegamat}
Let $\bA\in \real^{m\times n}$. Then,
\begin{itemize}
\item  $\bA$ is \textit{nonnegative} if all its entries $a_{ij}$ are real numbers and $a_{ij}\geq 0$ for all $i\in\{1,2,\ldots, m\}$ and $j\in\{1,2,\ldots,n\}$. This is denoted by $\bA\geq \bzero_{m,n}$ or $\bzero_{m,n}\leq \bA$ ($\bA\geq \bzero_{n}$ if $n=m$) or simply by $\bA\geq \bzero$.
\item $\bA$ is \textit{positive} if all its entries are real numbers and $a_{ij} >0$ for all $i\in\{1,2,\ldots, m\}$ and $j\in\{1,2,\ldots,n\}$. This is denoted by $\bA> \bzero_{m,n}$ or $ \bzero_{m,n}<\bA$ ($\bA> \bzero_{n}$ if $n=m$) or simply by $\bA>\bzero$. 
\item We denote  $\bA\geq\bB$ if $\bA-\bB\geq \bzero_{m,n}$ ($\bA>\bB$ if $\bA-\bB> \bzero_{m,n}$).~\footnote{Not to be confused with a positive definite matrix $\bA\succ\bzero$ or a positive semidefinite matrix $\bA\succeq \bzero$; see Definition~\ref{definition:psd-pd-defini}.}
\end{itemize}

The negative and nonpositive matrices are defined in a similar way. 
The set of nonnegative (resp. positive, negative, nonpositive) $m\times n$ matrices  is  denoted by $\real^{m\times n}_{\geq 0}$ (resp. $\real^{m\times n}_{>0}$, $\real^{m\times n}_{<0}$, $\real^{m\times n}_{\leq 0}$) or by $\real^{m\times n}_{+}$ (resp. $\real^{m\times n}_{++}$, $\real^{m\times n}_{--}$, $\real^{m\times n}_{-}$). 

Given a general matrix $\bA\in\real^{m\times n}$, the \textit{entrywise absolute matrix} is denoted by $\abs{\bA}\in\real^{m\times n}_{\geq 0}$, which is nonnegative. (If $\bA\in\complex^{m\times n}$ is complex, then the ($i,j$)-th element of $\abs{\bA}$ is $\sqrt{a_{ij}\overline{a}_{ij}}$).

A nonnegative $\bA\in\real^{n\times n}$ is said to be \textit{regular or primitive} if there exists a $k\geq 1$ such that $\bA^k$ is positive.
\end{definition}

\begin{definition}[Reducible and Irreducible Matrix]\label{definition:reduc_irreduc}
Let $\bA\in\real^{n\times n}$. Then, $\bA$ is called \textit{reducible} if there exists a permutation matrix $\bP$ such that 
$$
\bP\bA\bP^\top = \begin{bmatrix}
	\bA_{11} & \bA_{12}\\
	\bzero & \bA_{22}
\end{bmatrix},
$$
where $\bA\in\real^{r\times r}$, $\bA_{12}\in\real^{r\times (n-r)}$, and $\bA_{22}\in\real^{(n-r)\times (n-r)}$. 
Otherwise, $\bA$ is said to be \textit{irreducible}.
\end{definition}

\index{Nilpotent matrix}
\begin{definition}[Nilpotent Matrix]\label{definition:niopotent_mat}
A matrix $\bA\in \real^{n\times n}$ is \textit{nilpotent} if there exists a $k$ such that $\bA^k = \bzero$. 
The \textit{nilpotency} of $\bA$ is the number $\nilp(\bA)=\min\{k\in\integer_+ | \bA^k=\bzero\}$.
The $n\times n$ \textit{standard nilpotent matrix}, which has one's on the superdiagonal values and zero's elsewhere, is denoted by $\bN_n$ ($\bN_1=0$ and $\bN_0=0_{0\times 0}$ by convention). It can be shown that $\rank(\bN_n^k)=n-k$ for all $k\in\{0,1,\ldots,n\}$.
\end{definition}
The matrix $\bA=
\scriptsize
\begin{bmatrix}
0 & 1 \\
0 & 0
\end{bmatrix}$
is nilpotent since $\bA^2=\bzero$.
More generally, any triangular matrix with zero's along the main diagonal is nilpotent. For example, the matrix 
$
\bA=
\scriptsize
\begin{bmatrix}
0 & 1 & 3\\
0 & 0 & 4\\
0 & 0 & 0
\end{bmatrix}
$
is nilpotent since $\bA^3=\bzero$.

\begin{exercise}
Show that the following matrix is nilpotent:
\begin{equation}
\bA = 
\begin{bmatrix}
0 & a_{12} &  \cdots & 0 \\
0 & 0 & \ddots & 0 \\
0 & \cdots  & \ddots & a_{n-1,n} \\
0 &\cdots & \cdots & 0 
\end{bmatrix}.
\end{equation}
\end{exercise}


\index{Idempotent matrix}
\begin{definition}[Idempotent, Involutory, Skew-Involutory Matrix]\label{definition:Idempotent_mat}
A matrix $\bA\in \real^{n\times n}$ is called \textit{idempotent} if $\bA^2 = \bA$. The matrix  $\bA^\perp\triangleq\bI-\bA$ is called the \textit{complementary idempotent} of $\bA$, which is also idempotent:
\begin{itemize}
\item $\bA$ is called the idempotent matrix onto $\mathcalV=\cspace(\bA)$ along $\mathcalV^\perp=\nspace(\bA^\top)$.
\item Similarly, $\bA^\perp$ is called the idempotent matrix onto $\cspace(\bA^\perp)$ along $\nspace(\bA^{\perp\top})$.
\item $2\bA-\bI$ is involutory (see below).
\end{itemize}
A matrix $\bA\in \real^{n\times n}$ is  \textit{involutory} if $\bA^2=\bI$ ($\det(\bA)=\pm 1$); and $\bA$ is \textit{skew-involutory} if $\bA^2=-\bI$.
\end{definition}

\begin{definition}[Projection Matrix]\label{definition:projection_matrix_intro}
A matrix $\bA\in\real^{n\times n}$ is called an \textit{(orthogonal) projector or projection matrix} if it is symmetric and idempotent (Appendix~\ref{section:by-geometry-hat-matrix}); therefore, $\bA$ is positive semidefinite by definition (Definition~\ref{definition:psd-pd-defini}).
\begin{itemize}
\item In this sense, a projector is a normal matrix each of whose eigenvalues is 1 or 0 (Lemma~\ref{lemma:eigenvalues-of-projection}, Proposition~\ref{proposition:eigen-of-projection-matrix}).
\item  An \textit{elementary projector} is a projector exactly one of whose eigenvalues is 0.
\item A key property of a projector $\bA$ is that $\bA\bx=\bx$ if $\bx\in\cspace(\bA)$ (Proposition~\ref{proposition:column-space-of-projection}).
\item If $\bP$ has full column rank, then $\bA\triangleq\bP(\bP^\top\bP)^{-1}\bP^\top$ is a projector (Proposition~\ref{proposition:projection-from-matrix}); this is commonly used in least squares problems (Section~\ref{section:pre_ls}).
\item If $\bA$ is a projector onto a subspace $\mathcalV$, then $\bA^{\perp}\triangleq \bI-\bA$ is also a projector onto the complement subspace $\mathcalV^\perp$, called the \textit{complementary projector} of $\bA$, satisfying $$\cspace(\bA)^{\perp}=\nspace(\bA)=\cspace(\bA^\perp)=\mathcalV^\perp.$$
\end{itemize}

\end{definition}
\begin{exercise}[Projection Matrix]\label{exercise:projection_matrix_intro1}
Let $\bA\in\real^{n\times n}$ be idempotent. Show that the following statements are equivalent:
\begin{itemize}
\item $\bA$ is a projection matrix.
\item $\bA$ is symmetric.
\item $\bA\bA^\top\bA=\bA$.
\item $\bA$ is normal.
\end{itemize}
\end{exercise}
\begin{exercise}[Projection Matrix]\label{exercise:projection_matrix_intro2}
Let $\bx\in\real^n$ be nonzero. Show that $\bA\triangleq \bI-\frac{1}{\bx^\top\bx}\bx\bx^\top$ is an elementary projection matrix satisfying (a). $\rank(\bA)=n-1$;  (b). $\cspace(\bA)=\spn\{\bx\}^\perp$; (c). $\nspace(\bA)=\spn\{\bx\}$.
For the other way around, if $\bA\in\real^{n\times n}$ is a projector with $\rank(\bA)=n-1$, show that there is a nonzero vector $\bx\in\nspace(\bA)$ such that $\bA= \bI-\frac{1}{\bx^\top\bx}\bx\bx^\top$.
\end{exercise}

\begin{definition}[Reflection Matrix]\label{definition:reflection_mat_intro}
A matrix $\bA$ is called a \textit{reflector or reflection matrix} if it is symmetric and orthogonal (equivalently stated as symmetric and involutory; for example, the Householder reflector in Definition~\ref{definition:householder-reflector} is a reflector.)
\begin{itemize}
\item  In this sense, a reflector is a normal matrix each of whose eigenvalues is 1 or $-1$.
\item  An \textit{elementary reflector} is a reflector exactly one of whose eigenvalues is $-1$.
\end{itemize}
\end{definition}
\index{Orthogonal projection}
\index{Projection matrix}
\index{Projector}
\index{Reflection matrix}
\index{Reflector}
\index{Elementary matrix}
An \textit{elementary matrix} is a nonsingular matrix obtained by adding an outer-product matrix to the identity matrix: $\bA\triangleq\bI+\bx\by^\top$~\footnote{$\bA$ is nonsingular if $\bx^\top\by\neq -1$; see determinant results in Exercise~\ref{exercise:det_int}.}. In this context, elementary reflectors are considered elementary matrices. However, elementary projectors are not elementary matrices because they are singular.

\begin{exercise}[Reflection Matrix]\label{exer:reflection_mat_intro0}
Let $\bA\in\real^{n\times n}$. Show that the following statements are equivalent:
\begin{itemize}
\item $\bA$ is a reflection matrix.
\item $\bA=\bA\bA^\top+\bA^\top-\bI$.
\item $\bA=\frac{1}{2}(\bA+\bI)(\bA^\top+\bI)-\bI$.
\end{itemize}
\end{exercise}

\begin{exercise}[Reflection Matrix]\label{exer:reflection_mat_intro1}
Show that $\bB\triangleq 2\bA-\bI$ is a reflector if $\bA$ is a projector.
Conversely, show that $\bA\triangleq \frac{1}{2}(\bB+\bI)$ is a projector if $\bB$ is a reflector.
\end{exercise}
More properties of the reflector are discussed in Remark~\ref{remark:householder-propes} and Exercise~\ref{exercise:household_reflec1} when we discuss the Householder reflector for QR decomposition.

\index{Stochastic matrix}
\begin{definition}[Stochastic Matrix]\label{definition:stochastic_mat}
A nonnegative matrix $\bA\in\real^{n\times n}$ is called \textit{(row) stochastic} if $a_{ij}\geq 0$ for all $1\leq i,j \leq n$ and $\sum_{j=1}^{n} a_{ij}=1$ for all $1\leq i \leq n$; 
that is, $\bA\bone=\bone$, i.e., each row sums to 1. 
Theorem~\ref{theorem:spec_stoc} shows that a stochastic matrix has a spectral radius of one.

Similarly, a nonnegative matrix $\bA$ is \textit{column stochastic} if $\bone^\top\bA=\bone^\top$.
If both $\bA$ and $\bA^\top$ are  row stochastic, then $\bA$ is said to be \textit{doubly stochastic}, i.e., $\bA\bone=\bone$ and $\bone^\top\bA=\bone^\top$.
If $\bA\bone\leq \bone$ and $\bone^\top\bA\leq\bone^\top$, then $\bA$ is called \textit{doubly substochastic}.
\end{definition}
For example, the following $3\times 3$ matrix is a stochastic matrix:
$$
\bA = 
\begin{bmatrix}
	\frac{1}{2}& \frac{1}{2}& 0 \\
	\frac{1}{6}& \frac{1}{2}& \frac{1}{3} \\
	0 & \frac{1}{3}& \frac{2}{3}
\end{bmatrix}.
$$
On the other hand, permutation matrices are a special kind of doubly stochastic matrices.
\begin{exercise}[Doubly Substochastic]\label{exercise:dou_substo_int}
Let $\bQ,\bU\in\real^{n\times n}$ be orthogonal, and let $\bP\in\real^{n\times n}$ whose elements are $p_{ij}=\abs{q_{ij}u_{ji}}$ for $i,j\in\{1,2,\ldots,n\}$. 
Show that $\bP$ is doubly substochastic.
\textit{Hint: Use the Cauchy-Schwarz inequality (Equation~\eqref{equation:gen_ca_sc_1}).}
\end{exercise}

\begin{exercise}[Doubly Substochastic]\label{exercise:dou_substo_int2}
Let $\bA\in\real^{n\times n}$ be doubly substochastic. Show that there exists a doubly stochastic matrix $\bB\in\real^{n\times n}$ such that $\bB-\bA$ is nonnegative.
\textit{Hint: Show that if $\bA$ is doubly substochastic, there exists a doubly substochastic $\bC$ such that $\bC-\bA$ is nonnegative and $\bC$ is closer to a doubly stochastic matrix in the sense that the  entries of $\bA\bone$ and $\bone^\top\bA$  are closer to 1.}
\end{exercise}

\index{Determinant}
\begin{definition}[Determinant: Laplace Expansion by Minors]\label{definition:determinant}
Let $\bA\in\real^{n\times n}$ be any square matrix, and let $\bA_{ij}\in\real^{(n-1)\times (n-1)}$ denote the submatrix of $\bA$ obtained by deleting the $i$-th row and $j$-th column. 
The determinant of $\bA$ can be computed recursively using the following equations:
\begin{equation}\label{equation:def_det}
\det(\bA)=\sum_{k=1}^{n} (-1)^{i+k} a_{ik}\det(\bA_{ik})=\sum_{k=1}^{n}(-1)^{k+j} a_{kj}\det(\bA_{kj}),
\end{equation}
where the first equation is the \textit{Laplace expansion by minors along row $i$}, and the second equation is the \textit{Laplace expansion by minors along column $j$}.
Equivalently, given a cardinality $r$, and consider an index set $J\subseteq\{1,2,\ldots,n\}$ with cardinality $r$ ($\abs{J}=r$) and its complementary set $J^C=\{1,2,\ldots,n\}\backslash J$, we have:
$$
\begin{aligned}
\det(\bA) 
=\sum_{I} (-1)^{\gamma} \det(\bA[I,J])\det(\bA[I^C, J^C])
=\sum_{I} (-1)^{\gamma} \det(\bA[J,I])\det(\bA[J^C, I^C]),
\end{aligned}
$$
where $\gamma=\sum_{i\in I} i +\sum_{j\in J}j$, and the sum is taken over all the index sets $I\subseteq\{1,2,\ldots,n\}$ with cardinality $r$.
When $r=1$, this reduces to \eqref{equation:def_det}.
\end{definition}

\begin{remark}[Determinant: Alternating Sums and Permutations]\label{remark:det_altsumper}
Let the function $p:\{1,2,\ldots,n\}\rightarrow \{1,2,\ldots,n\}$ be a one-to-one function  of permutations, i.e., $p(i)=i$ in the identity case.
Then, there are $n!$ distinct permutations of the set $\{1,2,\ldots,n\}$. For a given permutation, let $\text{sgn}(p)=1$ if the minimum number of transpositions to achieve this permutation is even; and $\text{sgn}(p)=-1$ otherwise.
Then, the determinant can be equivalently defined  as
\begin{equation}
\det(\bA) = \sum_{p} \left( \text{sgn}(p) \prod_{i=1}^{n} a_{i p(i)} \right).
\end{equation}
\end{remark}

\index{Sylvester's identity for determinant}
\begin{theorem}[Sylvester's Identity]
Let $\bA\in\real^{n\times n}$, and  let $\bA_k=\bA[1:k,1:k]$ be the $k\times k$ leading principal submatrix of $\bA$. Let $\bB=[b_{ij}]\in\real^{(n-k)\times (n-k)}$, where $b_{ij}=\det(\bA[I,J])$, with  $I=\{1,2,\ldots,k, k+i\}$ and $J=\{1,2,\ldots,k, k + j\}$ for all $i,j\in\{1,2,\ldots,n-k\}$. Then,
$$
\det(\bB) = (\det(\bA_k))^{n-k-1} \det(\bA)
\quad\stackrel{\text{nonsingular} \bA_k}{\Longrightarrow }\quad
\det(\bA)=\frac{\det(\bB)}{(\det(\bA_k))^{n-k-1}} .
$$
And if $k=n-1$, this implies $\det(\bB)=\det(\bA)$.
\end{theorem}

The determinant of a square matrix maps the matrix into a scalar value.
In the case of a $2\times 2$ matrix, the determinant represents the area of the parallelogram spanned by the column vectors of the matrix. It is positive if the orientation is counterclockwise and negative if clockwise.
For a $3\times 3$ matrix, the determinant corresponds to the volume of the parallelepiped formed by the three column vectors of the matrix. Again, the sign indicates whether the orientation is preserved or reversed.
For a matrix representing a linear transformation in $n$-dimensional space, the absolute value of the determinant gives the factor by which the volume changes under this transformation. If the determinant is positive, the orientation (or handedness) of the basis is preserved; if negative, it is reversed.
A matrix is invertible if and only if its determinant is nonzero. This means that the transformation does not collapse the space into a lower dimension or a single point, which would happen if the determinant were zero.
The determinant of a matrix is equal to the product of its eigenvalues (Theorem~\ref{theorem:eigen_trace}). This means that the determinant reflects the combined effect of all the stretching factors applied by the matrix to the eigenvectors.
We provide a few properties of the determinant.
\begin{remark}[Properties of Determinant]\label{remark:determinant-intermezzo}
The following properties hold for the determinant of matrices:
\begin{itemize}
\item \textit{Multiplicativity.} The determinant of the product of two matrices is given by $\det(\bA\bB)=\det (\bA)\det(\bB)$.

\item \textit{Transpose.} The determinant of the transpose is $\det(\bA^\top) = \det(\bA)$.

\item Suppose matrix $\bA$ has an eigenvalue $\lambda$, then $\det(\bA-\lambda\bI) =0$.

\item The determinant of any singular matrix is 0 (nonzero for only nonsingular matrices).

\item The determinant  of any identity matrix is $1$.

\item The determinant  of an orthogonal matrix $\bQ$: 
$$
\det(\bQ) = \det(\bQ^\top) = \pm 1, \qquad \text{since  } \det(\bQ^\top)\det(\bQ)=\det(\bQ^\top\bQ)=\det(\bI)=1.
$$

\item For any square matrix $\bA$, and given an orthogonal matrix $\bQ$: 
$$
\det(\bA) = \det(\bQ^\top) \det(\bA)\det(\bQ) =\det(\bQ^\top\bA\bQ).
$$

\item Suppose $\bA\in\real^{n\times n}$, then $\det(-\bA) = (-1)^n \det(\bA)$.

\item \textit{Determinant of block matrices.} Let $\bM=\scriptsize\begin{bmatrix}
\bA&\bB\\
\bC &\bD 
\end{bmatrix}$. Then,
\begin{itemize}
\item If $\bA\in\real^{p\times p}$ is nonsingular, then, $\det(\bM)=\det(\bA)\det(\bD-\bC\bA^{-1}\bB)$ and $\rank(\bM)=p+\rank(\bD-\bC\bA^{-1}\bB)$.
\item If $\bD\in\real^{q\times q}$ is nonsingular, then, $\det(\bM)=\det(\bA-\bB\bD^{-1}\bC)\det(\bD)$ and $\rank(\bM)=q+\rank(\bA-\bB\bD^{-1}\bC)$. 
\item If both $\bA$ and $\bD$ are square and nonsingular, then $\bD-\bC\bA^{-1}\bB$ is nonsingular if and only if $\bA-\bB\bD^{-1}\bC$ is nonsingular.
\end{itemize}
\textit{Hint: see Problem~\ref{prob:det_block1}$\sim$\ref{prob:det_block2}.}

\item \textit{Determinant of block matrices.} More generally, let  $\bM\in\real^{n\times n}$. Consider an index set $I\subseteq\{1,2,\ldots,n\}$  and its complementary set $J=\{1,2,\ldots,n\}\backslash I$. If $\bM[I,I]$ is nonsingular, then, 
\begin{equation}\label{equation:det_block_rem}
\det(\bM)=\det(\bM[I,I])\det(\bM[J,J]-\bM[J,I]\bM[I,I]^{-1}\bM[I,J]).
\end{equation}
When $J$ is a single element, it reduces to \eqref{equation:det_block_rem2}.

\index{Cauchy-Binet formula}
\item \textit{Cauchy-Binet formula.} Consider $\bC=\bA\bB$ where $\bA\in\real^{m\times r}$ and $\bB\in\real^{r\times n}$, and let $I$ be an index set of $\{1,2,\ldots,m\}$ and $J$ be an index set of $\{1,2,\ldots,n\}$ with $\abs{I}=\abs{J}=k$, i.e., the same cardinality of $k$: $k\leq \min\{m,n,r\}$. Then, 
$$
\det(\bC[I,J])
=
\sum_{K:\abs{K}=k}
\det(\bA[I,K])\det(\bB[K,J]),
$$
where $K\subseteq\{1,2,\ldots,r\}$ and the sum is taken over all the index sets with cardinality $k$.
The Cauchy-Binet formula mimics the matrix multiplication.

\end{itemize}
More properties are discussed in Exercise~\ref{exercise:det_int}. The determinants of  elementary row transformations are discussed in Problem~\ref{problem:row_det_chg}.
\end{remark}

\begin{exercise}[Determinant Results]\label{exercise:det_int}
We have shown several important results regarding the determinant of a matrix in Remark~\ref{remark:determinant-intermezzo}. Given $\bA\in\real^{n\times n}$, show that
\begin{itemize}
	\item $\det(c\bA)=c^n\det(\bA)$. (\textit{Hint: Use induction.})
	\item $\det(\bA^{-1}) = 1/\det(\bA)$.
	\item $\det(\bA^m) = \det(\bA)^m$.
	\item $\det(\bI-\bu\bv^\top) = 1-\bu^\top\bv$. And $\bA\triangleq\bI-\bu\bv^\top\in\real^{n\times n}$ is nonsingular if and only if $\bu^\top\bv\neq 1$; and $\rank(\bA)=n-1$ if and only if $\bu^\top\bv=1$.
\end{itemize}
\end{exercise}

\index{Adjugate}
\begin{definition}[Adjugate]\label{definition:adjugate}
Let $\bA\in\real^{n\times n}$ be any square matrix. Then, the adjugate of $\bA$ is an $n\times n$ matrix whose $(i,j)$-th element is given by 
\begin{equation}\label{equation:adjug1}
\adjugate(\bA)_{ij} = (-1)^{i+j} \det\left(\bA\big[\{j\}^C, \{i\}^C\big]\right),
\end{equation}
where $\{i\}^C$ is the complementary set of $\{1,2,\ldots,n\}$: $\{i\}^C=\{1,2,\ldots,n\}\backslash i$.
Comparing this with the definition of determinants, we have
\begin{equation}\label{equation:adjug2}
\adjugate(\bA)\bA = \bA \adjugate(\bA) = \det(\bA)\bI.
\end{equation}
This shows that $\adjugate(\bA)$ is nonsingular (resp. upper triangular, diagonal) if $\bA$ is nonsingular (resp. upper triangular, diagonal): 
\begin{equation}\label{equation:adjug3}
\adjugate(\bA)=\det(\bA)\bA^{-1}.
\end{equation}
If $\bA$ is singular, then every column of $\bA$ lies in the null space of $\adjugate(\bA)$: the nullity/defect of $\adjugate(\bA)$ is at least the rank of $\bA$.
\end{definition}
For example, $\adjugate\big(\scriptsize\begin{bmatrix}
	a& b\\
	c& d
\end{bmatrix}\normalsize\big)
=
\scriptsize
\begin{bmatrix}
d & -b \\
-c & a
\end{bmatrix}
$.

\index{Cauchy expansion}
\begin{exercise}[Properties of Adjugate]\label{exercise:pro_adjug}
Let $\bA\in\real^{n\times n}$. Show that 
\begin{itemize}
\item $\adjugate(c\bA)=c^{n-1}\adjugate(\bA)\implies \adjugate(c\bI)=c^{n-1}\bI$.
\item $\det(\adjugate(\bA)) = (\det(\bA))^{n-1}$.
\item If $\bA$ is nonsingular, $\bA^{-1}=(\det(\bA))^{-1}\adjugate(\bA)$.
\item If $\bA$ is nonsingular, $\adjugate(\bA^{-1})=\bA/\det(\bA)$.
\item If $\bA$ is singular and $\rank(\bA)\leq n-2$, $\adjugate(\bA)=\bzero$.
\item If $\bA$ is singular and $\rank(\bA)= n-1$, $\rank(\adjugate(\bA))=1$.
\item If $\bA,\bB$ are nonsingular, $\adjugate(\bA\bB)=\adjugate(\bA)\adjugate(\bB)$. (This actually can be extended to all $\bA,\bB$ due to continuity.)
\item If $\bA$ is nonsingular, $\adjugate(\adjugate(\bA))=(\det(\bA))^{n-2}\bA$. (This actually can be extended to all $\bA$ due to continuity.)
\item If $\bA$ is nonsingular, $\adjugate(\bA^\top)=\adjugate(\bA)^\top$. (This actually can be extended to all $\bA$ due to continuity.) 
\item The adjugate is the transpose of the gradient of $\det(\bA)$: $\adjugate(\bA)_{ij}=\frac{\partial}{\partial a_{ji}}\det(\bA)$.
\item When $J$ in \eqref{equation:det_block_rem} is a single element, \eqref{equation:det_block_rem} reduces to 
\begin{equation}\label{equation:det_block_rem2}
\det(\bM)=\bM[J,J]\det(\bM[I,I]) - \bM[J,I]\adjugate\big(\bM[I,I]\big)\bM[I,J],
\end{equation}
which is valid even if $\bM[I,I]$ is singular (due to continuity).
When $J=\{n\}$ if $\bM\in\real^{n\times n}$, \eqref{equation:det_block_rem2} is called the \textit{Cauchy expansion} for the determinant of a \textit{bordered matrix}.
\end{itemize}
\end{exercise}

\paragraph{Consequence of $\rank(\bA)= n-1$.}
From the above results, if $\bA$ is singular and $\rank(\bA)= n-1$, $\rank(\adjugate(\bA))=1$. Therefore, there exist nonzero vectors $\bu, \bv\in\real^n$ and a nonzero scalar $a$ such that $\adjugate(\bA)=a\bu\bv^\top$. This indicates, by \eqref{equation:adjug2}, 
$
\bA\bu\bv^\top =\bu\bv^\top\bA=\bzero.
$
Hence, $\bA\bu=\bzero$ and $\bv^\top\bA=\bzero$.

From an introductory course on linear algebra, we have the following remark regarding  the equivalent claims of nonsingular matrices.
\begin{remark}[List of Equivalence of Nonsingularity for a Matrix]
Let $\bA\in \real^{n\times n}$ be any square matrix. Then, the following claims are equivalent:
\begin{itemize}
\item $\bA$ is nonsingular;
\footnote{The source of the name is discussed in Section~\ref{section:nonsingul_inver_svd}, which is a result of the singular value decomposition (SVD).}
\item $\bA$ is invertible, i.e., $\bA^{-1}$ exists; 
\item $\bA\bx=\bb$ has a unique solution $\bx = \bA^{-1}\bb$;
\item $\bA\bx = \bzero$ has a unique, trivial solution: $\bx=\bzero$;
\item Columns of $\bA$ are linearly independent;
\item Rows of $\bA$ are linearly independent;
\item $\det(\bA) \neq 0$; 
\item $\dim(\nspace(\bA))=0$;
\item $\nspace(\bA) = \{\bzero\}$, i.e., the null space is trivial;
\item $\cspace(\bA)=\cspace(\bA^\top) = \real^n$, i.e., the column space or row space spans the entire $\real^n$;
\item $\bA$ has full rank $r=n$;
\item The reduced row echelon form is $\bR=\bI$;
\item $\bA^\top\bA$ is symmetric positive definite;
\item $\bA$ has $n$ nonzero (positive) singular values;
\item All eigenvalues are nonzero.
\end{itemize}
\end{remark}

It is important to keep these equivalences in mind. On the other hand, the following remark shows the equivalent claims for singular matrices.
\begin{remark}[List of Equivalence of Singularity for a Matrix]
Let $\bA\in \real^{n\times n}$ be any square matrix with an eigenpair $(\lambda, \bu)$. Then, the following claims are equivalent:
\begin{itemize}
	\item $(\bA-\lambda\bI)$ is singular; \footnote{Again, the source of the name is discussed in Section~\ref{section:nonsingul_inver_svd}, which is a result of the singular value decomposition (SVD).}
	\item $(\bA-\lambda\bI)$ is not invertible;
	\item $(\bA-\lambda\bI)\bx = \bzero$ has nonzero $\bx\neq \bzero$ solutions, and $\bx=\bu$ is one  such solution;
	\item $(\bA-\lambda\bI)$ has linearly dependent columns;
	\item $\det(\bA-\lambda\bI) = 0$; 
	\item $\dim(\nspace(\bA-\lambda\bI))>0$;
	\item Null space of $(\bA-\lambda\bI)$ is nontrivial;
	\item Columns of $(\bA-\lambda\bI)$ are linearly dependent;
	\item Rows of $(\bA-\lambda\bI)$ are linearly dependent;
	\item $(\bA-\lambda\bI)$ has rank $r<n$;
	\item Dimension of column space = dimension of row space = $r<n$;
	\item $(\bA-\lambda\bI)^\top(\bA-\lambda\bI)$ is symmetric semidefinite;
	\item $(\bA-\lambda\bI)$ has $r<n$ nonzero (positive) singular values;
	\item Zero is an eigenvalue of $(\bA-\lambda\bI)$.
\end{itemize}
\end{remark}

\index{Linear system}
\begin{remark}[Linear System Solution]\label{remark:into_lin_syssol}
Given a linear system $\bA\bx = \bb$ with $\bA\in\real^{m\times n}$.
If there is at least one solution, the linear system is called \textit{consistent}; otherwise, it is called \textit{inconsistent}. It can be shown that the system is consistent if and only if $\rank([\bA,\bb])=\rank(\bA)$, i.e., the \textit{augmented matrix} $[\bA,\bb]$ has the same rank as the observed coefficient matrix $\bA$; this also means $\bb$ is in the column space of $\bA$.
Specifically, we can show that 
\begin{itemize}
\item When $\bA$ has full column rank $n$, the linear system has a  unique solution: $\widehat{\bx} = (\bA^\top\bA)^{-1}\bA^\top \bb$. Refer to its description in the left inverse (Theorem~\ref{theorem:unique-linear-system-solution}).
\item When $\bA$ has full row rank $m$, the linear system has at least one solution: $\widehat{\bx} = \bA_R^{-1}\bb$, where $\bA_R^{-1}$ is a right inverse of $\bA$. Refer to its description in the right inverse (Theorem~\ref{theorem:always-have-solution-right-inverse}).
\end{itemize}
\end{remark}

\index{Linear interpolation}
\begin{remark}[Linear Interpolation Problem]
A problem related to the linear system problem is called the \textit{linear interpolation problem}. Let $\bb\in\real^m$ and $\bx\in\real^n$.
Then, there exists a matrix $\bA\in\real^{m\times n}$ such that $\bA\bx=\bb$ if and only if $\bb=\bzero$ or $\bx\neq \bzero$.
(a). If $\bb=\bzero$, then one of such $\bA$ is $\bA=\bzero$. (b) If $\bx\neq \bzero$, then one of such $\bA$ is $\bA=\frac{\bb\bx^\top}{\bx^\top\bx}$.
\end{remark}

\index{Vector norm}
\index{Matrix norm}
\begin{definition}[Vector $\ell_2$ Norm]\label{definition:vec_l2_norm}
For a vector $\bx\in\real^n$, the \textit{$\ell_2$ vector norm} is defined as $\norm{\bx}_2 = \sqrt{x_1^2+x_2^2+\ldots+x_n^2}$.
\end{definition}

For a matrix $\bA\in\real^{m\times n}$, we  define the (matrix) Frobenius norm as follows. A detailed discussion on matrix norms can be found in Appendix~\ref{appendix:matrix-norm}. 
\begin{definition}[Matrix Frobenius Norm\index{Frobenius norm}]\label{definition:frobernius-in-svd}
The \textit{Frobenius norm} of a matrix $\bA\in \real^{m\times n}$ is defined as 
$$
\norm{\bA}_F = \sqrt{\sum_{i=1,j=1}^{m,n} (a_{ij})^2}=\sqrt{\trace(\bA\bA^\top)}=\sqrt{\trace(\bA^\top\bA)} = \sqrt{\sigma_1^2+\sigma_2^2+\ldots+\sigma_r^2}, 
$$
where $\sigma_1, \sigma_2, \ldots, \sigma_r$ are nonzero singular values of $\bA$.
\end{definition}
The Frobenius norm can be viewed as the $\ell_2$ norm applied to the vectorized form of the  matrix.
The following result regarding the Frobenius norm of the difference between two matrices is useful in various applications, such as low-rank matrix approximations.
\begin{lemma}\label{lemma:fro_norm_diff}
Let $\bA,\bB\in\real^{m\times n}$. The Frobenius norm of $\bA-\bB$ is 
$$
\norm{\bA-\bB}_F^2 = \norm{\bA}_F + \norm{\bB}_F - 2\trace(\bA^\top\bB).
$$
\end{lemma}

The spectral norm is defined as follows.
\begin{definition}[Matrix Spectral Norm]\label{definition:spectral_norm}
The \textit{spectral norm} of a matrix $\bA\in \real^{m\times n}$ is defined as 
$$
\norm{\bA}_2 = \mathop{\max}_{\bx\neq\bzero} \frac{\norm{\bA\bx}_2}{\norm{\bx}_2}  =\mathop{\max}_{\bx\in \real^n: \norm{\bx}_2=1}  \norm{\bA\bx}_2 ,
$$
which is also the maximal singular value of $\bA$, i.e., $\norm{\bA}_2 = \sigma_{\max}(\bA)$.
The definition also indicates the inequality: $\normtwo{\bA\bx}\leq \sigma_{\max}(\bA)\normtwo{\bx}$ for any $\bx\in\real^n$.
\end{definition}

We note that the Frobenius norm serves as the matrix counterpart of the vector $\ell_2$ norm.
For simplicity, we will not explicitly denote the full subscript for both the vector $\ell_2$ norm and the matrix Frobenius norm when it is clear from the context which one we are referring to: $\norm{\bA}=\norm{\bA}_F$ and $\norm{\bx}=\norm{\bx}_2$.
This means we treat the vector $\ell_2$ norm and the Frobenius norm as standard vector and mawtrix norms, respectively.
However, for the spectral norm, the subscript $\norm{\bA}_2$ should \textbf{not} be omitted.


The \textit{$\ell_2$ vector inner product} of a real vector space  is a function that associates to each pair of vectors $\bx, \by$ a real number, which is denoted by $\langle\bx,\by\rangle=\bx^\top\by$. 
Similarly, the \textit{Frobenius matrix inner product} of a real matrix space is a function that associates to each pair of matrices $\bA,\bB\in\real^{m\times n}$ a real number, which is denoted by $\langle \bA, \bB\rangle=\sum_{i,j=1}^{m,n}a_{ij}b_{ij}=\trace(\bA\bB^\top)$.

Given a specific norm definition, the concepts of an open ball and a closed ball are introduced as follows:
\begin{definition}[Open Ball, Closed Ball]
Let $\norm{\cdot}: \real^n\rightarrow \real_+$ be a norm function. The \textit{open ball} with center $\bc\in\real^n$ and radius $r$  is defined as 
$$
B(\bc, r) =\{\bx\in\real^n: \norm{\bx-\bc} <r\}.
$$
Similarly, the \textit{closed ball} with center $\bc\in\real^n$ and radius $r$  is defined as 
$$
B[\bc,r] =\{\bx\in\real^n: \norm{\bx-\bc} \leq r\}.
$$
\end{definition}

\section{Decomposition-Related Definitions}
To facilitate our discussion, we will introduce several definitions related to matrix decompositions. For simplicity, these terms are defined in the complex domain, although our primary focus is on the real domain.
\begin{definition}[Equivalent, Biequivalent]\label{definition:biequivalent}
Let $\bA,\bB\in\complex^{m\times n}$. Then, 
\begin{itemize}
\item $\bA$ and $\bB$ are \textit{left equivalent} if there exists a nonsingular $\bZ_1\in\complex^{m\times m}$ such that $\bA=\bZ_1\bB$.
\item $\bA$ and $\bB$ are \textit{right equivalent} if there exists a nonsingular $\bZ_2\in\complex^{n\times n}$ such that $\bA=\bB\bZ_2$.
\item $\bA$ and $\bB$ are \textit{biequivalent} if there exist nonsingular $\bZ_1\in\complex^{m\times m}$ and  $\bZ_2\in\complex^{n\times n}$ such that $\bA=\bZ_1\bB\bZ_2$.
\item $\bA$ and $\bB$ are \textit{unitarily left equivalent (resp. orthogonally left equivalent)}  if there exists a unitary (resp. real orthogonal) $\bZ_1\in\complex^{m\times m}$ such that $\bA=\bZ_1\bB$ (e.g., QR, QL).
\item $\bA$ and $\bB$ are \textit{unitarily right equivalent (resp. orthogonally right equivalent)} if there exists a unitary (resp. real orthogonal) $\bZ_2\in\complex^{n\times n}$ such that $\bA=\bB\bZ_2$ (e.g., LQ, RQ).
\item $\bA$ and $\bB$ are \textit{unitarily biequivalent (or simply unitarily equivalent; resp. orthogonally equivalent)} if there exist unitary (resp. real orthogonal) $\bZ_1\in\complex^{m\times m}$ and  $\bZ_2\in\complex^{n\times n}$ such that $\bA=\bZ_1\bB\bZ_2$ (e.g., SVD, UTV, Hessenberg, tridiagonal, bidiagonal).
\end{itemize}
The corresponding transformations are referred to as  left equivalence, right equivalence, biequivalence, unitary left equivalence, unitary right equivalence, and unitary (resp. orthogonal) biequivalence transformations, respectively.
\end{definition}

We will cover similarity transformations in more detail in Chapter~\ref{chapter:hessenberg}. Let's now briefly examine the related concepts of similarity and congruence.
\begin{definition}[Similar, Congruent]\label{definition:simiar_congru}
Let $\bA,\bB\in\complex^{n\times n}$. Then, 
\begin{itemize}
\item $\bA$ and $\bB$ are \textit{similar} if there exists a nonsingular $\bZ\in\complex^{n\times n}$ such that $\bA=\bZ\bB\bZ^{-1}$. If $\bZ$ is a permutation matrix, then $\bA$ and $\bB$ are called \textit{permutation similar} (e.g., block-diagonalization).
\item $\bA$ and $\bB$ are \textit{congruent (or star-congruent)} if there exists a nonsingular $\bZ\in\complex^{n\times n}$ such that $\bZ\bB\bZ^*$.
\item $\bA$ and $\bB$ are \textit{T-congruent} if there exists a nonsingular $\bZ\in\complex^{n\times n}$  such that $\bZ\bB\bZ^\top$.
\item $\bA$ and $\bB$ are \textit{unitarily similar (resp. orthogonally similar)} if there exists a unitary (resp. orthogonal) $\bZ\in\complex^{n\times n}$ such that $\bA=\bZ\bB\bZ^{-1}=\bZ\bB\bZ^*$ (e.g., Hessenberg, tridiagonal, spectral, Schur, Eigenvalue, Jordan).
\end{itemize}
The corresponding transformations are called similarity, congruence, T-congruence, and unitary (resp. orthogonal) similarity transformations, respectively.
\end{definition}

\begin{exercise}[Properties of Similarity Transformation]
Let $\bA$ and $\bB$ be similar, denoted by $\bA\sim\bB$. Show that the similarity relation is 
\begin{itemize}
\item \textit{Reflexive.} $\bA\sim\bA$.
\item \textit{Symmetric.} $\bA\sim\bB$ whenever $\bB\sim\bA$.
\item \textit{Transitive.} $\bA\sim\bC$ whenever $\bA\sim\bB$ and $\bB\sim\bC$.
\end{itemize}
How about the biequivalence transformation?
\end{exercise}


\section{Differentiability and Differential Calculus}
\index{Continuously differentiability}
\index{Second-order partial derivative}
\begin{definition}[Directional Derivative, Partial Derivative]\label{definition:partial_deri}
Given a function $f$ defined over a set $\sS\subseteq \real^n$ and a nonzero vector $\bd\in\real^n$, the \textit{directional derivative} of $f$ at $\bx$ w.r.t. the direction $\bd$ is given by, if the limit exists, 
$$
\mathop{\lim}_{t\rightarrow 0}
\frac{f(\bx+t\bd) - f(\bx)}{t}.
$$
And it is denoted by $f^\prime(\bx; \bd)$ or $D_{\bd}f(\bx)$. 
The directional derivative is sometimes called the \textit{G\^ateaux derivative}.

For any $i\in\{1,2,\ldots,n\}$, the directional derivative at $\bx$ w.r.t. the direction of the $i$-th standard basis vector  $\be_i$ (if it exists) is called the $i$-th \textit{partial derivative} and is denoted by $\frac{\partial f}{\partial x_i} (\bx)$, $D_{\be_i}f(\bx)$, or $\partial_i f(\bx)$.
\end{definition}

\begin{figure*}[h]
\centering  
\subfigtopskip=2pt 
\subfigbottomskip=9pt 
\subfigcapskip=-5pt 
\includegraphics[width=0.75\textwidth]{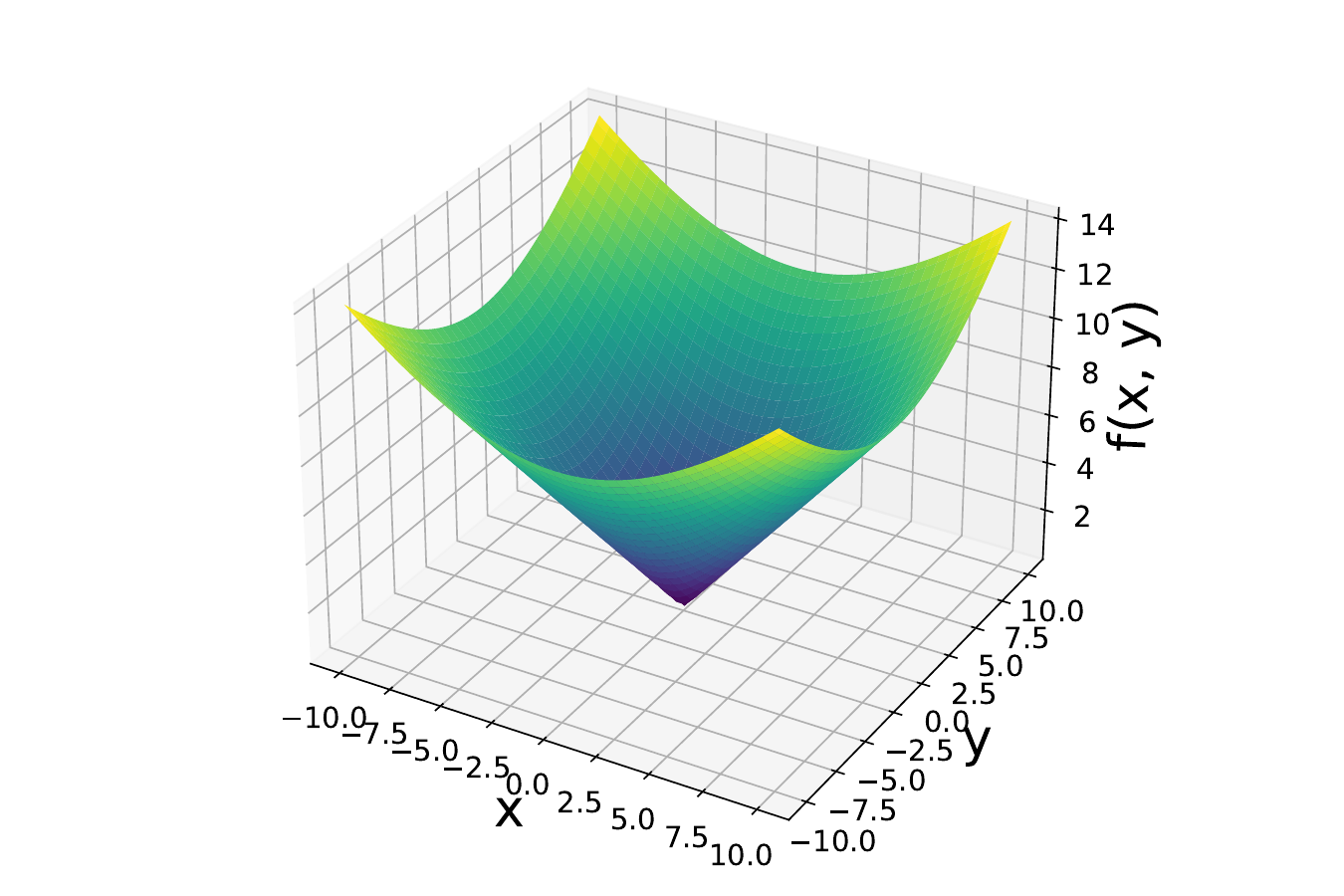}
\caption{Plot for function $f(x, y) = \sqrt{x^2+y^2}$, in which case any directional derivative for the direction $\bd=(a,b)$ with $a\neq 0$ and $b\neq 0$ at point $(0,0)$ exists. However, the partial derivatives at this point do not exist.}
\label{fig:direc_deri_no_partial}
\end{figure*}

Even if a function can have a directional derivative in every direction $\bd$ at some points, the partial derivative may not exist.
For example, given the direction $\bd=(a,b)$ with $a\neq 0$ and $b\neq 0$, the directional derivation of the function $f(x,y)=\sqrt{x^2+y^2}$ (see Figure~\ref{fig:direc_deri_no_partial}) at $(0,0)$ can be obtained by 
$$
f^\prime(0,0; \bd) = \mathop{\lim}_{t\rightarrow 0}
\frac{\sqrt{(0+ta)^2+(0+tb)^2} -\sqrt{0^2+0^2} }{t} = \sqrt{a^2+b^2}.
$$
When $\bd$ is a unit vector, the directional derivative is 1.
However, it can be shown that the partial derivatives are 
$$
\begin{aligned}
\frac{\partial f}{\partial x}(0,0) &= \mathop{\lim}_{h\rightarrow 0} \frac{\sqrt{(0+h)^2+0^2}}{h} = \mathop{\lim}_{h\rightarrow 0} \frac{\abs{h}}{h},\\
\frac{\partial f}{\partial y}(0,0) &= \mathop{\lim}_{h\rightarrow 0} \frac{\sqrt{0^2+(0+h)^2}}{h} = \mathop{\lim}_{h\rightarrow 0} \frac{\abs{h}}{h}.
\end{aligned}
$$
When $h>0$, the partial derivatives are 1; when $h<0$, the partial derivaties are $-1$. Therefore, the partial derivatives do not exist.

If all the partial derivatives of a function $f$ exist at a point $\bx\in\real^n$, then the \textit{gradient} of $f$ at $\bx$ is defined as the column vector containing all the partial derivatives:
$$
\nabla f(\bx)=
\begin{bmatrix}
\frac{\partial f}{\partial x_1} (\bx)\\
\frac{\partial f}{\partial x_2} (\bx)\\
\vdots \\
\frac{\partial f}{\partial x_n} (\bx)	
\end{bmatrix}
\in \real^n.
$$

\begin{exercise}
Demonstrate that a function can have partial derivatives but is not necessarily continuous. Also, show that a function can be continuous but does not necessarily have partial derivatives.
\end{exercise}

A function $f$ defined over an open set $\sS\subseteq \real^n$ is called \textit{differentiable}
if all the partial derivatives exist (i.e., the derivatives exist in univariate cases, and the gradients exist in multivariate cases). This is actually the definition of \textit{Fr\'echet differentiability}.

\begin{definition}[(Fr\'echet) Differentiability]
Given a function $f:\real^n\rightarrow \real$, the function $f$ is said to be differentiable at $\bx$ if there exists a vector $\bg\in\real^n$ such that 
$$
\mathop{\lim}_{\bd\rightarrow\bzero} \frac{f(\bx+\bd)-f(\bx)-\bg^\top\bd}{\normtwo{\bd}}=0.
$$
The unique vector $\bg$ is equal to the gradient $\nabla f(\bx)$.
\end{definition}

Moreover, a function $f$ defined over an open set $\sS\subseteq \real^n$ is called \textit{continuously differentiable} over $\sS$ if all the partial derivatives exist and are continuous on $\sS$.

A differentiable function may not be continuously differentiable. For example, consider  the following function $f$:
$$
f(x) = 
\begin{cases}
x^2 \sin(\frac{1}{x}), & \text{if } x \neq 0, \\
0, & \text{if } x = 0.
\end{cases}
$$
This function is differentiable everywhere, including at $x=0$. To see this, for $x \neq 0$, $f(x)$ is a product of two differentiable functions ($x^2$ and $\sin(1/x)$), hence it is differentiable.
At $x = 0$, we can calculate the derivative using the limit definition:
$$
f^\prime(0) = \mathop{\lim}_{t\rightarrow 0}
\frac{f(t) - f(0) }{t}
=\mathop{\lim}_{t\rightarrow 0} \frac{t^2 \sin(\frac{1}{t}) - 0 }{t}
=t \sin(\frac{1}{t}) 
$$
Since $\abs{\sin(1/t)}\leq 1$ for all $t \neq 0$, the limit exists and is 0. Thus, $f(x)$ is differentiable at $x = 0$ with $f^\prime(0) = 0$.
However,  $f(x)$ is not continuously differentiable at $x = 0$. The derivative of $f(x)$ when $x\neq 0$ is 
$$
f^\prime(x) = 2x \sin(\frac{1}{x}) - \cos(\frac{1}{x}).
$$
The limit $\mathop{\lim}_{x\rightarrow 0}f^\prime(x)$ does not exist since the sine and cosine  functions oscillate as $x$ approaches 0.

In the setting of \textit{continuously differentiability}, the directional derivative and gradient have the following relationship:
\begin{equation}
f^\prime(\bx; \bd) = \nabla f(\bx)^\top \bd, \gap \text{for all }\bx\in\sS \text{ and }\bd\in\real^n.
\end{equation} 
And, to recap, in the setting of  differentiability, we also have
\begin{equation}
\mathop{\lim}_{\bd\rightarrow \bzero}
\frac{f(\bx+\bd) - f(\bx) - \nabla f(\bx)^\top \bd}{\norm{\bd}} = 0\gap 
\text{for all }\bx\in\sS,
\end{equation}
or 
\begin{equation}
f(\by) = f(\bx)+\nabla f(\bx)^\top (\by-\bx) + o(\norm{\by-\bx}),
\end{equation}
where $o(\cdot): \real_+\rightarrow \real$ is a one-dimensional function satisfying $\frac{o(t)}{t}\rightarrow 0$ as $t\rightarrow 0^+$.
Therefore, any differentiable function indicates that $\mathop{\lim}_{\bx\rightarrow \ba} f(\bx) = f(\ba)$. Hence, any differentiable function is continuous. 

On the other hand, in the setting of differentiability, although the partial derivatives may not be continuous, the partial derivatives exist for all differentiable points.

The partial derivative $\frac{\partial f}{\partial x_i} (\bx)$ is also a real-valued function of $\bx\in\sS$ that can be partially differentiated. The $j$-th partial derivative of $\frac{\partial f}{\partial x_i} (\bx)$ is defined as 
$$
\frac{\partial^2 f}{\partial x_j\partial x_i} (\bx)=
\frac{\partial \left(\frac{\partial f}{\partial x_i} (\bx)\right)}{\partial x_j} (\bx).
$$
This is called the ($j,i$)-th \textit{second-order partial derivative} of function $f$.
A function $f$ defined over an open set $\sS\subseteq$ is called \textit{twice continuously differentiable} over $\sS$ if all the second-order partial derivatives exist and are continuous over $\sS$. In the setting of twice continuously differentiability, the second-order partial derivative are symmetric:
$$
\frac{\partial^2 f}{\partial x_j\partial x_i} (\bx)=
\frac{\partial^2 f}{\partial x_i\partial x_j} (\bx).
$$
The \textit{Hessian} of the function $f$ at a point $\bx\in\sS$ is defined as the symmetric $n\times n$ matrix 
$$
\nabla^2f(\bx)=
\begin{bmatrix}
\frac{\partial^2 f}{\partial x_1^2} (\bx) & 
\frac{\partial^2 f}{\partial x_1\partial x_2} (\bx) & \ldots &
\frac{\partial^2 f}{\partial x_1\partial x_n} (\bx)\\
\frac{\partial^2 f}{\partial x_2\partial x_1} (\bx) & 
\frac{\partial^2 f}{\partial x_2\partial x_2} (\bx) & \ldots &
\frac{\partial^2 f}{\partial x_2\partial x_n} (\bx)\\
\vdots & 
\vdots & \ddots &
\vdots\\
\frac{\partial^2 f}{\partial x_n\partial x_1} (\bx) & 
\frac{\partial^2 f}{\partial x_n\partial x_2} (\bx) & \ldots &
\frac{\partial^2 f}{\partial x_n^2} (\bx)
\end{bmatrix}.
$$
We provide a simple proof of Taylor's expansion in Appendix~\ref{appendix:taylor-expansion}  for one-dimensional functions. 
In the case of high-dimensional functions, we have the following three approximation results.

\begin{theorem}[Mean Value Theorem]\label{theorem:mean_approx}
Let $f(\bx):\sS\rightarrow \real$ be a  continuously differentiable function over an open set $\sS\subseteq\real^n$, and given two points $\bx, \by$. Then there exists $\bx^\star\in[\bx,\by]$ such that 
$$
f(\by) = f(\bx)+ \nabla f(\bx^\star)^\top (\by-\bx).
$$ 
\end{theorem}

\begin{theorem}[Linear Approximation Theorem]\label{theorem:linear_approx}
Let $f(\bx):\sS\rightarrow \real$ be a twice continuously differentiable function over an open set $\sS\subseteq\real^n$, and given two points $\bx, \by$. Then there exists $\bx^\star\in[\bx,\by]$ such that 
$$
f(\by) = f(\bx)+ \nabla f(\bx)^\top (\by-\bx) + \frac{1}{2} (\by-\bx)^\top \nabla^2 f(\bx^\star) (\by-\bx).
$$ 
\end{theorem}

\begin{theorem}[Quadratic Approximation Theorem]\label{theorem:quad_app_theo}
Let $f(\bx):\sS\rightarrow \real$ be a twice continuously differentiable function over an open set $\sS\subseteq\real^n$, and given two points $\bx, \by$. Then it follows that 
$$
f(\by) = f(\bx)+ \nabla f(\bx)^\top (\by-\bx) + \frac{1}{2} (\by-\bx)^\top \nabla^2 f(\bx) (\by-\bx)
+
o(\norm{\by-\bx}^2).
$$ 
\end{theorem}

\section{Convex Functions}

We  briefly introduce the notion of convex functions.
First, we define the concepts of open and closed sets.

\begin{definition}[Interior Points, Open Sets, Closed Sets]
Given a set $\sS\subseteq \real^n$, a point $\ba\in\sS$ is called an \textit{interior point} of $\sS$ if there exists a scalar $s>0$ such that the open ball $B(\ba, s)\subseteq\sS$.

A set $\sS_1$ is said to be open if it consists of only interior points. That is, for every $\bx\in\sS$, there exists a scalar $s>0$ such that $B(\bx, s)\subseteq\sS_1$.

A set $\sS_2$ is said to be closed if its complement $\sS_2^c$ is open.
Alternatively, $\sS_2$ is closed if it contains all the limit points of convergent sequence of points in $\sS_2$; that is, for every sequence of points $\{\bx_i\}_{i\geq 1}\subseteq \sS_2$ satisfying $\bx_i\rightarrow \bx^\star$ as $i\rightarrow \infty$, it follows that $\bx^\star\in\sS_2$.
\end{definition}

Given vectors $\bx$ and $\by$,  the point $\lambda\bx+(1-\lambda)\by$ with $\lambda\in[0,1]$ is known as  a convex combination of the two vectors.
A set that is closed under arbitrary convex combinations is referred to as a convex set.
A standard definition is given below.
\begin{definition}[Convex Set]\label{definition:convexset}
A set $\sS\subseteq \real^n$ is called \textit{convex} if, for any $\bx,\by\in\sS$ and $\lambda\in[0,1]$, the point $\lambda\bx+(1-\lambda)\by$ is also in $\sS$.
\end{definition}
Geometrically, convex sets  contain all line segments
that join two points within the set. Consequently, these sets do not feature any concave indentations.

A related concept is that of convex functions, which exhibit specific behavior under convex combinations. We now recall the definition:
\begin{definition}[Convex Functions]\label{definition:convexfuncs}
A function $f:\sS\rightarrow \real$ defined over a convex set $\sS\subseteq \real^n$ is called \textit{convex} if 
$$
f(\lambda\bx+(1-\lambda)\by)
\leq \lambda f(\bx) + (1-\lambda) f(\by), \text{ for any }\bx,\by\in\sS, \lambda\in[0,1].
$$
And the function  $f$ is called \textit{strictly convex} if 
$$
f(\lambda\bx+(1-\lambda)\by)
< \lambda f(\bx) + (1-\lambda) f(\by),\text{ for any }\bx\neq \by\in\sS, \lambda\in(0,1).
$$
\end{definition}
\begin{exercise}[Convex Functions]\label{exercise:convexfuncs}
Let $f:\sS\rightarrow \real$ be a convex function, and let $p\geq 2$ be any integer. Show that 
$
f\big(\sum_{i=1}^{p} \lambda_i\bx_i\big) \leq \sum_{i=1}^{p}\lambda_i f(\bx_i),
$
if $\lambda_i\geq 0$ and $\sum_{i=1}^{p}\lambda_i=1$ (i.e., $\{\lambda_i\}$ belongs to the unit simplex in $\real^p$.)
\end{exercise}

In scenarios involving multiple variables, the input vector  $\bx$ can be partitioned into $\bx=(\ba, \bb)\in\real^n$. 
When discussing properties of functions in this context, convexity can sometimes be specified as \textit{joint convexity}, meaning that the function is convex with respect to all components of the input simultaneously.
Additionally, there are instances where a function exhibits \textit{partial convexity}, often termed \textit{marginal convexity}, indicating that the function is convex when considered with respect to one or more subsets of its input variables while keeping others constant.

\begin{definition}[Marginally Convex]\label{definition:marginal_convexfuncs}
A function $f(\ba,\bb):\sS\rightarrow \real$ (i.e., $(\ba,\bb)\in\real^n$) defined over a convex set $\sS\subseteq \real^n$ is called \textit{marginally convex} if 
$$
f(\lambda\bx+(1-\lambda)\by, \bb)
\leq \lambda f(\bx, \bb) + (1-\lambda) f(\by, \bb), \text{ for any }(\bx,\bb),(\by,\bb)\in\sS, \lambda\in[0,1].
$$
\end{definition}

\begin{exercise}
Show that the norm function $\normtwo{\bx}^2$ is convex.
\end{exercise}

\begin{exercise}
Suppose the function $f$ is convex. Show that the directional derivative $g(\bd)=f'(\bx;\bd)$ is also convex.
\end{exercise}

\paragraph{First-order characterizations of convex functions.}
Convex functions need not be differentiable. However, when they are continuously differentiable, the convex function can be characterized by the gradient inequality.
\begin{theorem}[Gradient Inequality]\label{theorem:conv_gradient_ineq}
If the function $f:\sS\rightarrow \real$ is a  continuously differentiable function defined over a convex set $\sS\subseteq\real^n$, then the function is convex over $\sS$ if and only if 
$$
f(\bx) +\nabla f(\bx)^\top (\by-\bx)\leq f(\by), \text{ for any $\bx,\by\in\sS$}.
$$
Similarly, the function is strictly convex over $\sS$ if and only if 
$$
f(\bx) +\nabla f(\bx)^\top (\by-\bx)< f(\by), \text{ for any $\bx\neq \by\in\sS$}.
$$
This means the graph of a convex function lies above the tangent plane at any point.
For concave or strictly concave functions, the inequality signs are reversed.
\end{theorem}

\paragraph{Second-order characterizations of convex functions.}
When the function is further assumed to be twice continuously differentiable, convexity can be defined by the positive semidefiniteness (see Definition~\ref{definition:psd-pd-defini}) of the Hessian matrix.
\begin{theorem}[PSD Hessian of Convex Functions]
If the function $f:\sS\rightarrow \real$ is a twice continuously differentiable function defined over an open convex set $\sS\subseteq\real^n$, then the function is convex if and only if $\nabla^2 f(\bx)\geq \bzero $ for any $\bx\in\sS$.
Moreover, if $\nabla^2 f(\bx)> \bzero $, then the function is strictly convex over $\sS$~\footnote{The condition is only sufficient but not necessary. For example, $f(x)=x^6$ is strictly convex, but $f^{\prime\prime}(x)=30x^4$ is equal to zero at $x=0$.}. 
\end{theorem}

\begin{definition}[Strongly Convex/Smooth Functions]\label{definition:scss_func}
A continuously differentiable function $f:\real^n\rightarrow \real$ is called \textit{$\sigma$-strongly convex (SC)} and \textit{$L$-strongly smooth (SS)} if, for every $\bx,\by\in\real^n$, it follows that 
$$
\frac{\sigma}{2} \norm{\bx-\by}^2
\leq f(\by)-f(\bx)-f(\bx)^\top (\by-\bx)
\leq \frac{L}{2} \norm{\bx-\by}^2.~\footnote{In many texts, the second inequality is called the \textit{descent lemma} for $L$-strongly smooth functions since $f(\by) \leq f(\bx)+f(\bx)^\top (\by-\bx)+\frac{L}{2} \norm{\bx-\by}^2$, i.e., an update of the function value  from $f(\bx)$ to $f(\by)$ is upper-bounded.}
$$
It is evident that if the function is $L_1$-strongly smooth, it is also $L_2$-strongly smooth for any $L_2\geq L_1$.
A $\sigma$-strongly convex function is necessarily a convex function.
While if the function is $\sigma_1$-strongly convex, it is also $\sigma_2$-strongly convex for any $\sigma_2\in(0, \sigma_1]$.
\end{definition}
Strong convexity and smoothness are important properties in the analysis of gradient descent algorithms (see, for example, Problem~\ref{problem:mono_gd}).
For the strongly convexity, it ensures the function curve upwards sharply and provides a lower bound on the curvature of the function, which prevents the gradient from vanishing too quickly.
The strongly smoothness, on the other hand, provides an upper bound on the curvature of the function and also means $\norm{\nabla f(\bx) - \nabla f(\by)} \leq L \norm{\bx-\by}$, ensuring that the gradient does not change too rapidly. 
That is to say, it controls the magnitude of gradients and prevents  the gradient descent steps from being too aggressive.

\section{Matrix Decomposition in a Nutshell}
We now provide a brief overview of the decompositional methods that will be covered in detail throughout the rest of the book:
\begin{enumerate}[leftmargin=*]
\item \textit{LU.} $\bA=\bL\bU = 
\left(\parbox{8em}{lower triangular $\bL$\\ 1's on the diagonal}\right)
\left(\parbox{9.5em}{upper triangular $\bU$\\pivots on the diagonal}\right)$

\textbf{Requirements:} $\bA$ has nonzero leading principal minors, i.e., no row permutations involve to reduce $\bA$ into upper triangular via  Gaussian elimination (Theorem~\ref{theorem:lu-factorization-without-permutation}, p.~\pageref{theorem:lu-factorization-without-permutation}).

\item \textit{LU.} $\bA=\bL\bD\bU = \left(\parbox{8em}{lower triangular $\bL$\\ 1's on the diagonal}\right)
\left(\parbox{6em}{pivot matrix  \\$\bD$ is diagonal}\right)
\left(\parbox{8.5em}{upper triangular $\bU$\\1's on the diagonal}\right)$

\textbf{Requirements:} $\bA$ has nonzero leading principal minors, i.e., no row permutations involves to reduce $\bA$ into upper triangular via  Gaussian elimination. The matrix $\bD$ contains the pivots such that $\bU$ has 1's on the diagonal. When $\bA$ is symmetric, it follows that $\bA=\bL\bD\bL^\top$ (Corollary~\ref{corollary:ldu-decom}, p.~\pageref{corollary:ldu-decom}).

\item \textit{PLU.} $\bA=\bP\bL\bU =\left(\parbox{7em}{permutation \\matrix  $\bP$ avoids\\ zeros for\\ eliminations}\right) \left(\parbox{8em}{lower triangular $\bL$\\ 1's on the diagonal}\right)\left(\parbox{9.5em}{upper triangular $\bU$\\pivots on the diagonal}\right)$

\textbf{Requirements:} $\bA$ is nonsingular, then $\bP,\bL$, and $\bU$ are nonsingular as well. $\bP^\top$ performs the row permutation in advance such that $\bP^\top\bA$ has nonzero leading principal minors (Theorem~\ref{theorem:lu-factorization-with-permutation}, p.~\pageref{theorem:lu-factorization-with-permutation}).

\item \textit{RRLU.} $\bP\bA\bQ = \begin{bmatrix}
	\bL_{11} & \bzero \\
	\bL_{21}^\top & \bI 
\end{bmatrix}
\begin{bmatrix}
	\bU_{11} & \bU_{12} \\
	\bzero & \bzero 
\end{bmatrix}$

\textbf{Requirements:} Any matrix with rank $r$ such that $\bP$ and $\bQ$ are permutations, $\bL_{11}$ is lower triangular with 1's on the diagonal, $\bU_{11}$ is upper triangular with pivots on the diagonal, where both of them reveal the rank of matrix $\bA$ (Section~\ref{section:rank-reveal-lu-short}, p.~\pageref{section:rank-reveal-lu-short}).

\item \textit{Complete Pivoting LU.} $\bP\bA\bQ = \bL\bU=
\left(\parbox{8em}{lower triangular $\bL$\\ 1's on the diagonal}\right)\left(\parbox{9.5em}{upper triangular $\bU$\\pivots on the diagonal}\right)$

\textbf{Requirements:} $\bA$ is nonsingular such that $\bP$ and $\bQ$ are permutations, $\bL$ is unit lower triangular with 1's on the diagonal, $\bU$ is upper triangular with pivots on the diagonal (Section~\ref{section:complete-pivoting}, p.~\pageref{section:complete-pivoting}).

\item \textit{Cholesky.} $\bA=\bR^\top\bR$=(lower triangular)(upper triangular)

\textbf{Requirements:} $\bA$ is positive definite (symmetric) such that the diagonal of $\bR$ contains the diagonals of $\sqrt{\bD}$, which are positive from the LU decomposition. Trivially, $\bR^\top$ can be set as $\bL\sqrt{\bD}$ from the LU decomposition. When $\bA$ is positive semidefinite, $\bR$ is upper triangular with possible zeros on the diagonal (Theorem~\ref{theorem:cholesky-factor-exist}, p.~\pageref{theorem:cholesky-factor-exist}; Theorem~\ref{theorem:semidefinite-factor-exist}, p.~\pageref{theorem:semidefinite-factor-exist}).

\item \textit{Pivoted Cholesky.} $\bP\bA\bP^\top = \bR^\top\bR $=(lower triangular)(upper triangular)

\textbf{Requirements:} $\bA$ is positive definite (symmetric) such that $\bR$ is upper triangular (Section~\ref{section:volum-picot-cholesjy}, p.~\pageref{section:volum-picot-cholesjy}).

\item \textit{Semidefinite Rank-Revealing.} $\bP^\top \bA\bP  = \bR^\top\bR, \gap \mathrm{with} \gap 
\bR = \begin{bmatrix}
	\bR_{11} & \bR_{12}\\
	\bzero &\bzero 
\end{bmatrix} $

\textbf{Requirements:} $\bA$ is positive semidefinite with rank $r$ such that after permuting by $\bP$, $\bR_{11}$ is upper triangular with rank $r$ (Theorem~\ref{theorem:semidefinite-factor-rank-reveal}, p.~\pageref{theorem:semidefinite-factor-rank-reveal}).

\item \textit{Bunch Kaufman.}  $\bP^\top \bA\bP  = \bL\bB\bL^\top $=(unit lower triangular)(block-diagonal)(unit upper triangular)

\textbf{Requirements:} $\bA$ is symmetric such that $\bB$ is a block-diagonal matrix with each diagonal block of  $\bB$ being either a  $1\times 1$ or a  $2\times 2$ matrix (Theorem~\ref{theorem:Bunch_Kaufman}, p.~\pageref{theorem:Bunch_Kaufman}).

\item \textit{CR.} $\bA=\bC\bR$=(first independent columns of $\bA$)(basis for row space of $\bA$)

\textbf{Requirements:} Any matrix $\bA$ with rank $r$, $\bC$ contains the first $r$ linearly independent columns of $\bA$, $\bR$ is the reduced row echelon form with zero rows removed. $\bR$ contains an $r\times r$ identity submatrix (Theorem~\ref{theorem:cr-decomposition}, p.~\pageref{theorem:cr-decomposition}).

\item \textit{Rank.} $\bA=\bD\bF$=(column basis of $\bA$)(row basis of $\bA$)

\textbf{Requirements:} Any matrix $\bA$ with rank $r$, $\bD$ contains $r$ columns that span the column space of $\bA$, $\bF$ contains $r$ rows that span the row space of $\bA$ (Theorem~\ref{theorem:rank-decomposition}, p.~\pageref{theorem:rank-decomposition}).

\item \textit{Skeleton.} $\bA=\bC\bU^{-1}\bR$=($r$ columns of $\bA$)(intersection of $\bC$ and $\bR$)$^{-1}$($r$ rows of $\bA$)

\textbf{Requirements:} Any matrix $\bA$ with rank $r$ such that $\bC$ and $\bR$ are derived directly from the columns and rows of $\bA$, respectively. $\bU$ is the mixing matrix that is on the intersection of $\bC$ and $\bR$ (Theorem~\ref{theorem:skeleton-decomposition}, p.~\pageref{theorem:skeleton-decomposition}).

\item \textit{Interpolative.} $\bA=\bC\bW$=(independent columns of $\bA$)(columns consisting of an identity submatrix)

\textbf{Requirements:} Any matrix $\bA$ with rank $r$ such that $\bW$ contains an $r\times r$ identity submatrix (in the sense of permutation). The elements of $\bW$ are no greater than 1 in absolute value (Theorem~\ref{theorem:interpolative-decomposition}, p.~\pageref{theorem:interpolative-decomposition}; Theorem~\ref{theorem:interpolative-decomposition-row}, p.~\pageref{theorem:interpolative-decomposition-row}).

\item \textit{QR.} $\bA = \bQ\bR$=(orthonormal columns in $\bQ$)(upper triangular $\bR$)

\textbf{Requirements:} Rectangular matrix $\bA$ has linearly independent columns; otherwise $\bR$ will be singular (diagonals contain at least one zero). In the full QR decomposition, $\bQ$ is orthogonal such that $\bQ^{-1}=\bQ^\top$ (Theorem~\ref{theorem:qr-decomposition}, p.~\pageref{theorem:qr-decomposition}).

\item \textit{QL.} $\bA\bP = \bQ\bL$=(orthonormal columns in $\bQ$)(lower triangular $\bL$)

\textbf{Requirements:} $\bP$ is a reversal matrix (Definition~\ref{definition:permutation-matrix}). Rectangular matrix $\bA$ has linearly independent columns; otherwise $\bL$ will be singular (diagonals contain at least one zero). In the full QL decomposition, $\bQ$ is orthogonal such that $\bQ^{-1}=\bQ^\top$ (Section~\ref{section:lq-decomp}, p.~\pageref{section:lq-decomp}).

\item \textit{CPQR.} 
$\bA = \bQ
\begin{bmatrix}
	\bR_{11} & \bR_{12} \\
	\bzero   & \bzero 
\end{bmatrix}\bP^\top
=\left(\parbox{5.9em}{orthogonal $\bQ$}\right)
\left(\parbox{8.em}{upper triangular $\bR_{11}$ $|$ full matrix $\bR_{12}$}\right)
\left(\parbox{5.2em}{permutation $\bP^\top$}\right)$

\textbf{Requirements:} Any rectangular matrix $\bA$ with rank $r$ such that $\bR_{11}$ is $r\times r$ upper triangular (Theorem~\ref{theorem:rank-revealing-qr-general}, p.~\pageref{theorem:rank-revealing-qr-general}).

\item \textit{RRQR.} 
$\bA = \bQ
\begin{bmatrix}
	\bR_{11} & \bR_{12} \\
	\bzero   & \bR_{22}
\end{bmatrix}\bP^\top
=\left(\parbox{5.9em}{orthogonal $\bQ$}\right)
\left(\parbox{8.em}{upper triangular $\bR_{11}$ $|$ full matrix $\bR_{12}$ $|$ $\bR_{22}$ small in norm}\right)
\left(\parbox{5.2em}{permutation $\bP^\top$}\right)$

\textbf{Requirements:} Any rectangular matrix $\bA$ with rank $r$ such that $\bR_{11}$ is $r\times r$ upper triangular and $\bR_{22}$ is small in norm (Section~\ref{section:rank-r-qr}, p.~\pageref{section:rank-r-qr}).

\item \textit{LQ.} $\bA=\bL\bQ$=(lower triangular $\bL$)(orthonormal rows in $\bQ$)

\textbf{Requirements:} Rectangular matrix $\bA$ has linearly independent rows; otherwise $\bL$ will be singular (diagonals contain at least one zero). In the full LQ decomposition, $\bQ$ is orthogonal such that $\bQ^{-1}=\bQ^\top$ (Theorem~\ref{theorem:lq-decomposition}, p.~\pageref{theorem:lq-decomposition}).

\item \textit{RQ.} $\bP\bA=\bR\bQ$=(upper triangular $\bR$)(orthonormal rows in $\bQ$)

\textbf{Requirements:} $\bP$ is a reversal matrix (Definition~\ref{definition:permutation-matrix}). Rectangular matrix $\bA$ has linearly independent rows; otherwise $\bR$ will be singular (diagonals contain at least one zero). In the full RQ decomposition, $\bQ$ is orthogonal such that $\bQ^{-1}=\bQ^\top$ (Section~\ref{section:lq-decomp}, p.~\pageref{section:lq-decomp}).

\item \textit{RPLQ.} 
$\bA = \bP^\top
\begin{bmatrix}
		\bL_{11} & \bzero \\
\bL_{21} & \bzero 
\end{bmatrix}\bQ
=\left(\text{permutation $\bP^\top$}\right)
\left(\parbox{8.em}{lower triangular $\bL_{11}$ $|$ full matrix $\bL_{21}$}\right)
\left(\parbox{5.9em}{orthogonal $\bQ$}\right)$

\textbf{Requirements:} Any rectangular matrix $\bA$ with rank $r$ such that $\bL_{11}$ is $r\times r$ lower triangular (Section~\ref{section:lq-decomp}, p.~\pageref{section:lq-decomp}).

\item \textit{Two-Sided Orthogonal.} $\bA\bP\bA=\bU\bF\bV^\top=\left(\parbox{5.3em}{orthonormal  \\ basis in $\bU$}\right)
\left(\parbox{4.5em}{upper-left \\$r\times r$\\ submatrix\\ is nonzero}\right)
\left(\parbox{5.3em}{orthonormal \\ basis in $\bV$}\right)$

\textbf{Requirements:} Square matrix with rank $r$ such that the first $r$ columns of $\bU$ span the column space of $\bA$, and the remaining columns span the null space of $\bA^\top$. Whilst, the first $r$ columns of $\bV$ span the row space of $\bA$, and the remaining columns span the null space of $\bA$ (Theorem~\ref{theorem:two-sided-orthogonal}, p.~\pageref{theorem:two-sided-orthogonal}).

\item \textit{UTV.} $\bA = \bU\bT\bV$=(orthogonal $\bU$)
$\left(\parbox{8.4em}{upper triangular $\bR$\\lower triangular $\bL$}\right)$
(orthogonal $\bV$)

\textbf{Requirements:} Any matrix $\bA$ with rank $r$ such that $\bT$ is lower or upper triangular with rank $r$ (Theorem~\ref{theorem:ulv-decomposition}, p.~\pageref{theorem:ulv-decomposition}).

\item \textit{Hessenberg.} $\bA=\bQ\bH\bQ^\top$=(orthogonal matrix $\bQ$)(Hessenberg matrix $\bH$)($\bQ^\top=\bQ^{-1}$)

\textbf{Requirements:} Any square matrix $\bA$ such that $\bA$ and $\bH$ are orthogonally similar matrices with the same rank, trace, and eigenvalues (Theorem~\ref{theorem:hessenberg-decom}, p.~\pageref{theorem:hessenberg-decom}).

\item \textit{Hessenberg-Triangular Form.} $\bA=\bQ\bH\bZ^\top$ and $\bB=\bQ\bR\bZ^\top$, ($\bH$ Hessenberg)($\bR$ upper triangular)

\textbf{Requirements:} Given a matrix pair $(\bA, \bB)$ where $\bA, \bB \in \real^{n \times n}$, there exist  orthogonal matrices $\bQ, \bZ \in \real^{n \times n}$ such that $\bQ^\top \bA \bZ$ is upper Hessenberg while $\bQ^\top \bB \bZ$ is upper triangular (Section~\ref{section:ht_form}, p.~\pageref{section:ht_form}).

\item \textit{Tridiagonal.} $\bA=\bQ\bT\bQ^\top$=(orthogonal $\bQ$)(symmetric tridiagonal $\bT$)($\bQ^\top=\bQ^{-1}$)

\textbf{Requirements:} Any symmetric matrix $\bA$ such that $\bA$ and $\bT$ are orthogonally similar matrices with the same rank, trace, and eigenvalues (Theorem~\ref{definition:tridiagonal-hessenbert}, p.~\pageref{definition:tridiagonal-hessenbert}).

\item \textit{Bidiagonal.} $\bA=\bU\bB\bV^\top$=(orthogonal matrix $\bU$)(bidiagonal $\bB$)(orthogonal matrix $\bV$)

\textbf{Requirements:} Any rectangular matrix $\bA$ such that $\bB$ is upper bidiagonal. $\bT=\bB^\top\bB$ is tridiagonal if $\bB$ is bidiagonal such that $\bA^\top\bA = \bV\bB^\top\bB\bV^\top$ is the tridiagonal decomposition of $\bA^\top\bA$ (Theorem~\ref{theorem:Golub-Kahan-Bidiagonalization-decom}, p.~\pageref{theorem:Golub-Kahan-Bidiagonalization-decom}).

\item \textit{Eigenvalue.} $\bA = \bX\bLambda\bX^{-1}
=\left(\parbox{6.5em}{eigenvectors in\\ columns of $\bX$}\right)
\left(\parbox{6.1em}{eigenvalues in \\diagonal of $\bLambda$}\right)
\left(\parbox{7em}{left eigenvectors \\in $\bX^{-1}$}\right)$

\textbf{Requirements:} Square matrix $\bA$ has $n$ linearly independent eigenvectors (Theorem~\ref{theorem:eigenvalue-decomposition}, p.~\pageref{theorem:eigenvalue-decomposition}).

\item \textit{Schur.} $\bA = \bQ\bU\bQ^\top
=\left(\parbox{6.6em}{orthogonormal \\columns in $\bQ$}\right)
\left(\parbox{8.5em}{upper triangular $\bU$}\right)
\left(\parbox{4.8em}{$\bQ^\top=\bQ^{-1}$}\right)$

\textbf{Requirements:} $\bA$ is any real matrix with real eigenvalues (Theorem~\ref{theorem:schur-decomposition}, p.~\pageref{theorem:schur-decomposition}).

\item \textit{Spectral.} $\bA=\bQ\bLambda\bQ^\top
=\left(\parbox{7.6em}{orthogonormal \\eigenvectors in $\bQ$}\right)
\left(\parbox{7.8em}{real eigenvalues in \\diagonal of $\bLambda$}\right)
\left(\parbox{7em}{left eigenvectors \\in $\bQ^{\top}$}\right)$

\textbf{Requirements:} $\bA$ is real and symmetric, which is a special case of the Schur decomposition (Theorem~\ref{theorem:spectral_theorem}, p.~\pageref{theorem:spectral_theorem}). When $\bA$ is complex or real normal, $\bQ$ is complex unitary and $\bLambda$ contains  complex eigenvalues; when $\bA$ is complex Hermitian, $\bQ$ is complex unitary  and $\bLambda$ is real diagonal (Theorem~\ref{theorem:normal_Complex_spectral_theorem}, p.~\pageref{theorem:normal_Complex_spectral_theorem}).
See also the spectral decomposition for the pseudo-inverses of normal matrices in Theorem~\ref{theorem:norm_pseudo_spec} (p.~\pageref{theorem:norm_pseudo_spec}).

\item \textit{Jordan.} $\bA = \bX\bJ\bX^{-1}
=\left(\parbox{7.3em}{generalized eigenvectors \\in columns of $\bX$}\right)
\left(\parbox{6.4em}{Jordan blocks \\in $\bJ$}\right)
\left(\parbox{8.8em}{left generalized  \\eigenvectors in $\bX^{-1}$}\right)$

\textbf{Requirements:} $\bA$ is any square matrix, $\bJ$ is a Jordan form matrix, containing identical eigenvalues in each block (Theorem~\ref{theorem:jordan-decomposition}, p.~\pageref{theorem:jordan-decomposition}).

\item \textit{SVD.} $\bA=\bU \bSigma \bV^\top
=\left(\parbox{5.5em}{left singular \\vectors in $\bU$}\right)
\left(\parbox{7.6em}{singular values in \\diagonal of $\bSigma$}\right)
\left(\parbox{5.8em}{right singular \\vectors in $\bV$}\right)$

\textbf{Requirements:} Any matrix $\bA$ such that the diagonals of $\bSigma^2$ contain the eigenvalues of $\bA^\top\bA$ or $\bA\bA^\top$. In the full SVD, both $\bU$ and $\bV$ are orthogonal (Theorem~\ref{theorem:full_svd_rectangular}, p.~\pageref{theorem:full_svd_rectangular}). 

\item \textit{Polar.} $\bA=\bQ_l\bS_l
=\left(\parbox{4.6em}{orthogonal matrix $\bQ_l$}\right)
\left(\parbox{6.5em}{positive \\semidefinite $\bS_l$}\right)$

\textbf{Requirements:} $\bA$ is any square matrix. The polar decomposition indicates $\bS_l^2 = \bA^\top\bA$. When $\bA$ is nonsingular, $\bS_l$ is positive definite. Moreover, there exists a right polar decomposition $\bA=\bS_r\bQ_r$ such that $\bS_r^2 = \bA\bA^\top$ (Theorem~\ref{theorem:polar-decomposition}, p.~\pageref{theorem:polar-decomposition}).

\item \textit{ALS.} $\bA \approx \bW\bZ$=$\left(s=\rank(\bW)<\rank(\bA)=r\right)$
$\left(s=\rank(\bZ)<\rank(\bA)=r\right)$

\textbf{Requirements:} For any rectangular matrix $\bA$, $\bW\bZ$ approximates $\bA$ in terms of the Frobenius norm (Chapter~\ref{chapter:als}, p.~\pageref{chapter:als}).

\item \textit{NMF.} $\bA \approx \bW\bZ$=$\left(s=\rank(\bW)<\rank(\bA)=r\right)$
$\left(s=\rank(\bZ)<\rank(\bA)=r\right)$

\textbf{Requirements:} For any rectangular matrix $\bA$, $\bW\bZ$ approximates $\bA$ in terms of the Frobenius norm with all entries of $\bA,\bW$, and $\bZ$ being nonegative (Section~\ref{chapter:nmf}, p.~\pageref{chapter:nmf}).

\item \textit{Biconjugate.} $\bA = \bPhi \bOmega^{-1} \bPsi^\top$

\textbf{Requirements:} Any rectangular matrix $\bA$ with rank $r$ such that $\bOmega$ is $r\times r$ diagonal. The columns of $\bPhi$ and $\bPsi$ come from the Wedderburn sequence (Theorem~\ref{theorem:biconjugate-form1}, p.~\pageref{theorem:biconjugate-form1}; Theorem~\ref{theorem:biconjugate-form2}, p.~\pageref{theorem:biconjugate-form2}).

\item \textit{CP.} $\eX \approx \llbracket \bA^{(1)}, \bA^{(2)}, \ldots, \bA^{(N)}   \rrbracket  
=
\sum_{r=1}^{R} \ba_r^{(1)} \circ \ba_r^{(2)} \circ \ldots \circ \ba_r^{(N)} $

\textbf{Requirements:} Any tensor $\eX \in \real^{I_1\times I_2\times \ldots \times I_N}$ such that the CP decomposition is a low-rank approximation of the original tensor (Theorem~\ref{theorem:cp-decomp}, p.~\pageref{theorem:cp-decomp}).

\item \textit{Tucker.} $\eX \approx \llbracket\eG; \bA^{(1)}, \bA^{(2)}, \ldots, \bA^{(N)}   \rrbracket  
=
\sum_{r_1=1}^{R_1} \sum_{r_2=1}^{R_2} \ldots \sum_{r_N=1}^{R_N}
g_{r_1r_2\ldots r_N}
\ba_r^{(1)} \circ \ba_r^{(2)} \circ \ldots \circ \ba_r^{(N)} $

\textbf{Requirements:} Any tensor $\eX \in \real^{I_1\times I_2\times \ldots \times I_N}$ such that the Tucker decomposition captures the principal components in high-order dimensions (Theorem~\ref{theorem:tucker-decomp}, p.~\pageref{theorem:tucker-decomp}).

\item \textit{HOSVD.} $\eX \approx \llbracket\eG; \bA^{(1)}, \bA^{(2)}, \ldots, \bA^{(N)}   \rrbracket  
=
\sum_{r_1=1}^{R_1} \sum_{r_2=1}^{R_2} \ldots \sum_{r_N=1}^{R_N}
g_{r_1r_2\ldots r_N}
\ba_r^{(1)} \circ \ba_r^{(2)} \circ \ldots \circ \ba_r^{(N)} $

\textbf{Requirements:} Any tensor $\eX \in \real^{I_1\times I_2\times \ldots \times I_N}$ such that the HOSVD captures the principal components in high-order dimensions, where  the slices of $\eG$ are mutually orthogonal and  are ordered in  descending order (Theorem~\ref{theorem:hosvd-decomp}, p.~\pageref{theorem:hosvd-decomp}).

\item \textit{TT.} $\eX \approx \eG^{(1)} \boxtimes \eG^{(2)} \boxtimes \ldots \boxtimes \eG^{(N)}$

\textbf{Requirements:} Any tensor $\eX \in \real^{I_1\times I_2\times \ldots \times I_N}$ such that the TT decomposition factors the tensor into $N$ third-order tensors (Theorem~\ref{theorem:ttrain-decomp}, p.~\pageref{theorem:ttrain-decomp}).
 
\end{enumerate}




%% file: chapter-LU.tex
\part{Gaussian Elimination}

\section*{Introduction}
\lettrine{\color{caligraphcolor}I}
In linear algebra, \textit{Gaussian elimination} is often referred to as the  \textit{row reduction} process, which is an algorithm for solving systems of linear equations such that the original full linear equation is converted into  \textit{row echelon form}   through a series of elementary row operations (or in some cases into an upper triangular one) whilst the computational complexity is reduced.
The elementary row operations used in Gaussian elimination are as follows:
\begin{enumerate}
\item Swapping two rows.
\item Multiplying a row by a nonzero scalar.
\item Adding a scalar multiple of one row to another row.

\end{enumerate}

The row reduction process in Gaussian elimination can be divided into two phases. 
The first phase, also known as \textit{forward elimination}, transforms the system into an upper triangular matrix or, more generally, into the row echelon form. During this transformation, key properties of the system, such as its rank, are preserved. These properties provide insights into the existence, uniqueness, and nature of the solutions.

The second phase involves \textit{back substitution}. Here, the system is further transformed into the reduced row echelon form, from which the solution can be directly read off.

Gaussian elimination can suffer from numerical instability when dealing with ill-conditioned or nearly singular matrices. In such situations, alternative methods, like partial pivoting, are often employed to achieve a more accurate solution.

In this part, we will introduce  LU and Cholesky decompositions. Another pertinent decomposition,  CR decomposition, will be delayed into Part~\ref{part:data-interation}  as it skeletonizes the matrix and compresses the matrix into a thin one whilst sparsity and  nonnegativity of the matrices are preserved.

\chapter{LU Decomposition}
\begingroup
\hypersetup{
	linkcolor=structurecolor,
	linktoc=page,  
}
\minitoc \newpage
\endgroup

\section{LU Decomposition}
\lettrine{\color{caligraphcolor}P}
Perhaps the most well-known and fundamental matrix decomposition we should be aware of is the \textit{LU decomposition}, where \textit{L} stands for lower triangular and \textit{U} stands for upper triangular. 
The decomposition factors a given matrix into a product of a lower triangular matrix and an upper triangular matrix. LU decomposition is closely related to Gaussian elimination, and both methods are widely employed for solving linear systems, computing determinants, and inverting matrices.
We will now present the results in the following theorem, with the proof of its existence to be discussed in subsequent sections.

\index{Decomposition: LU}
\begin{theoremHigh}[LU Decomposition with Permutation]\label{theorem:lu-factorization-with-permutation}
Let $\bA$ be any nonsingular $n\times n$ square matrix. Then, it  can be factored as 
\begin{equation}
	\bA = \bP\bL\bU, \nonumber
\end{equation}
where $\bP$ is a permutation matrix, $\bL$ is a \textit{unit lower triangular matrix} (i.e., a lower triangular matrix with all 1's on the diagonal), and $\bU$ is a \textit{nonsingular} upper triangular matrix. 
\end{theoremHigh}

Henceforth, in this text,  results related to decomposition will be presented in blue boxes, while other claims and theorems will be presented in gray boxes. This convention will be followed throughout the book without any special mention.
The LU decomposition originated from Gauss's elimination algorithm, which he sketched in 1809 \citep{gauss1877theoria} and presented in full in 1810 \citep{gauss1828disquisitio}.

One of the main applications of LU decomposition is in solving linear systems. 
Given a system of linear equations, $\bA\bx = \bb$, where $\bA$ is a square matrix and $\bb$ is a target vector, the LU decomposition can be used to find $\bx$ efficiently. First, the matrix $\bA$ is factored into matrices $\bP$, $\bL$, and $\bU$, and then the system of equations is rewritten as $\bP\bw = \bb$, $\bL\bv=\bw$, and $\bU\bx = \bv$. The system can then be solved using forward and backward substitutions. The LU decomposition is particularly useful when solving systems of equations that have the same coefficient matrix but different right-hand sides, i.e., solving the $n$ linear equations $\bA\bx = \bb_1$, $\bA\bx = \bb_2, \ldots, \bA\bx = \bb_n$ simultaneously. See \S~\ref{section:lu-linear-sistem} for more details.

LU decomposition is also used in matrix inversion, where the inverse of a matrix $\bA$ is calculated as $\bA^{-1} =\bU^{-1} \bL^{-1}\bP^\top$. However, this method is not always the most efficient way to compute the inverse of a matrix, and other methods such as the Gauss-Jordan method or the QR decomposition may be more suitable for some applications \citep{arias2011suitable, ji2012gauss, ma2011qr}. See \S~\ref{section:inverse-by-lu}  for more details.

Another application of LU decomposition is in computing the determinant of a matrix. The determinant of a matrix can be calculated as the product of the diagonal elements of the $\bU$ matrix, multiplied by the product of the diagonal elements of the $\bL$ matrix and a constant ($-1$) raised to the power of the number of row exchanges performed during the LU decomposition. See \S~\ref{sectin:app_lu_det}  for more details.

In addition to its use in solving linear systems, LU decomposition is also used in other numerical methods such as computational fluid dynamics, finite element analysis, and optimization algorithms \citep{shakib1991new, johan1992data, maurer2011parallel, harbrecht2012low}.

\begin{remark}[Decomposition Notation]
The above decomposition applies to any nonsingular matrix $\bA$. We will see that this decomposition arises from  \textit{elimination steps}, in which case row operations of subtraction and exchange of two rows are allowed; the subtractions are recorded in the matrix $\bL$, and the row exchanges are recorded in the matrix $\bP$. To make this row exchange explicit, the common form for the  decomposition above is $\bQ\bA=\bL\bU$, where $\bQ=\bP^\top$ precisely records the  row exchanges in $\bA$. Otherwise, the $\bP$ would record the row exchanges of $\bL\bU$. 
In our case, we will focus on making the decomposition clear for the matrix $\bA$ itself, rather than for $\bQ\bA$. For this reason, we will place the permutation matrix on the right-hand side of the equation for the remainder of the text without specific notice.
\end{remark}

In certain situations, the permutation matrix $\bP$ may not be necessary for the decomposition.
Instead, we can rely on the conditions related to leading principal minors. To illustrate this approach, we provide a rigorous definition that will be crucial for the subsequent discussion.
\index{Principal minors}
\begin{definition}[Principal Minors]\label{definition:principle-minors}
Let $\bA$ be an $n\times n$ square matrix. A $k \times k$ submatrix of $\bA$ obtained by deleting any $n-k$ columns and the same $n-k$ rows from $\bA$ is called a $k$-th order \textit{principal submatrix} of $\bA$. The determinant of a $k \times k$ principal submatrix is called a $k$-th order \textit{principal minor} of $\bA$.
\end{definition}

We may be able to obtain a particular type of principal minors under mild conditions on the selected indices of the submatrix, i.e., the \textit{leading principal minors}.
\begin{definition}[Leading Principal Minors\index{Leading principal minors}]\label{definition:leading-principle-minors}
Let $\bA$ be an $n\times n$ square matrix. A $k \times k$ submatrix of $\bA$ obtained by deleting the \textbf{last} $n-k$ columns and the \textbf{last} $n-k$ rows from $\bA$ is called the $k$-th order \textit{leading principal submatrix} of $\bA$. 
That is, the $k\times k$ submatrix taken from the top-left corner of $\bA$, i.e., $\bA[1:k,1:k]$ using  Matlab-style notation. The determinant of the $k \times k$ leading principal submatrix is called the $k$-th order \textit{leading principal minor} of $\bA$.
\end{definition}

Given an $n\times n$ matrix $\bA$ with $(i,j)$-th entry being $a_{ij}$, let $\bA[1:k,1:k]$ denote the $k\times k$ submatrix taken from the top-left corner of $\bA$. That is, 
$$
\bA_{1:k,1:k} = \bA[1:k,1:k] = 
\begin{bmatrix}
a_{11} & a_{12} & \ldots & a_{1k}\\
a_{21} & a_{22} & \ldots & a_{2k}\\
\vdots & \vdots & \ddots & \vdots \\
a_{k1} & a_{k2} & \ldots & a_{kk}\\
\end{bmatrix}.~\footnote{We may use the two Matlab-style notations interchangeably: $\bA_{1:k,1:k} = \bA[1:k,1:k]$.}
$$
Then $\Delta_k=\det(\bA[1:k,1:k] )$ is the $k$-th order leading principal minor of $\bA$.

When certain conditions are met regarding the leading principal minors of matrix $\bA$, it is possible to perform the LU decomposition without the need for the permutation matrix $\bP$.
\begin{theoremHigh}[LU Decomposition without Permutation]\label{theorem:lu-factorization-without-permutation}
Let  $\bA$ be any $n\times n$ square matrix. If all the leading principal minors are nonzero, i.e., $\det(\bA[1:k,1:k])\neq 0$  for all $k\in \{1,2,\ldots, n\}$, then $\bA$ can be factored as 
\begin{equation}
	\bA = \bL\bU, \nonumber
\end{equation}
where $\bL$ is a unit lower triangular matrix (i.e., a lower triangular matrix with all 1's on the diagonal), and $\bU$ is a \textit{nonsingular} upper triangular matrix. 
Specifically, this decomposition is \textbf{unique}. See Corollary~\ref{corollary:unique-lu-without-permutation}.
\end{theoremHigh}

\begin{remark}[Other Forms of the LU Decomposition without Permutation]
In Theorem~\ref{theorem:lu-factorization-without-permutation}, we assume that the leading principal minors are nonzero; in other words, this implies that both the leading principal submatrices  and the matrix $\bA$ are nonsingular.

\paragraph{Singular $\bA$.} In the preceding theorem, we assume $\bA$ is nonsingular. The LU decomposition also exists for singular matrix $\bA$. However, in such cases, the matrix $\bU$ will also be singular. This will be shown in the following section, where it will be shown that if matrix $\bA$ is singular, certain pivots will be zero, leading to  corresponding zero diagonal values in  $\bU$.

\paragraph{Singular leading principal submatrices.} Even when assuming that matrix $\bA$ is nonsingular, it's important to note that the leading principal submatrices might still be singular. 
If certain leading principal minors are zero, the LU decomposition exists, but in such cases, it is not unique.
\end{remark}

The origin of this decomposition will be addressed in the following section. 
Additionally, there are extensions of the LU decomposition to encompass non-square or singular matrices, such as rank-revealing LU decomposition. 
For a more in-depth exploration, readers are recommended to refer to \citet{stewart2000decompositional, pan2000existence, miranian2003strong, dopico2006multiple}. Alternatively, a brief discussion will be provided in  Section~\ref{section:rank-reveal-lu-short}. 

\index{Linear system}
\index{Backward substitution}
\index{Gaussian elimination}
\section{Relation to Gaussian Elimination}\label{section:gaussian-elimination}
Solving the linear system equation $\bA\bx=\bb$ constitutes a fundamental problem in linear algebra.
\textit{Gaussian elimination} transforms a linear system into an upper triangular form by applying simple \textit{elementary row transformations} on the left of the linear system within $n-1$ stages if $\bA\in \real^{n\times n}$. 
Consequently, solving the system becomes more straightforward through \textit{backward substitution}. The elementary transformation is precisely defined as follows.

\begin{definition}[Elementary Transformation\index{Elementary transformation}]\label{definition:elemen_trans}
Given a  square matrix $\bA$, the following three transformations are referred to as \textit{elementary row (resp. column) transformations}:
\begin{enumerate}
\item   Interchanging two rows (resp. columns) of $\bA$ (using permutation matrices in Definition~\ref{definition:permutation-matrix});
\item  Multiplying all elements of a row (resp. a column) of $\bA$ by a nonzero scalar;
\item   Adding any row (resp. column) of $\bA$ multiplied by a nonzero number to any other row (or column).
\end{enumerate}
The three types of transformation will modify the determinant of a matrix; see Problem~\ref{problem:row_det_chg}.
\end{definition}
Specifically, the elementary row transformations of $\bA$ are represented by unit lower triangular matrices, multiplying on the left of $\bA$  (e.g., $\bE\bA$), and the elementary column transformations of $\bA$ are represented by unit upper triangular matrices, multiplying on the right  of  $\bA$ (e.g., $\bA\bE$). For example, 
the following is the type-2 elementary transformation matrix:
\begin{equation}\label{equation:ele_exam_type2}
\textbf{(Type-2):}\quad \bE_{i;\alpha}\triangleq  \diag(1,\ldots,1,\alpha,1,\ldots,1) = \bI + (\alpha-1)\be_i\be_i^\top\in\real^{n\times n}.
\end{equation} 
$\bE_{i;\alpha}\bA$ multiplies the $i$-th row of $\bA$ by $\alpha$, and $\bA\bE_{i;\alpha}$ multiplies the $i$-th column of $\bA$ by $\alpha$.
And the following matrix is the type-3 elementary transformation matrix:
\begin{equation}\label{equation:ele_exam_eijal}
\textbf{(Type-3):}\qquad \bE_{i,j;\alpha}\triangleq \bI+\alpha\be_i\be_j^\top\in\real^{n\times n}.
\end{equation} 
When $i>j$, $\bE_{i,j;\alpha}$ is an elementary row transformation such that $\bE_{i,j;\alpha}\bA$ adds $\alpha$ multiples of the $j$-th row to the $i$-th row; 
when $i<j$, $\bE_{i,j;\alpha}$ is an elementary column transformation such that $\bA\bE_{i,j;\alpha}$ adds $\alpha$ multiples of the $i$-th column to the $j$-th column. 

Gaussian elimination is  described by the third type (elementary row transformation) above.
Suppose the upper triangular matrix obtained by Gaussian elimination is given by $\bU = \bE_{n-1}\bE_{n-2}\ldots\bE_1\bA$  ($n-1$ steps); and at the $k$-th stage, the $k$-th column of $\bE_{k-1}\bE_{k-2}\ldots\bE_1\bA$ is represented as $\bx\in \real^n$. Gaussian elimination will introduce zeros below the diagonal of $\bx$ using the transformation
\begin{equation}\label{equation:elimination_mat}
\textbf{(Elimination matrix)}:\qquad \bE_k \triangleq \bI - \bz_k \be_k^\top,
\end{equation} 
where $\be_k \in \real^n$ is the $k$-th unit basis vector, and $\bz_k\in \real^n$ is given by 
$$
\bz_k = [0, \ldots, 0, z_{k+1}, \ldots, z_n]^\top, \qquad z_i= \frac{x_{i}}{x_{k}}, \gap\forall i \in \{k+1,\ldots, n\}.
$$
We realize that $\bE_k$ is a unit lower triangular matrix (with $1$'s on the diagonal) and has only the $k$-th column of the lower submatrix as nonzero,
$$
\bE_k=
\footnotesize
\begin{bmatrix}
1 & \ldots & 0&  0 & \ldots & 0\\
\vdots & \ddots & \vdots & \vdots & \ddots & \vdots \\
0 & \ldots & 1  & 0 & \ldots & 0 \\
0 & \ldots & -z_{k+1}  & 1 & \ldots & 0\\
\vdots & \ddots & \vdots & \vdots & \ddots & \vdots \\
0 & \ldots & -z_n & 0 & \ldots & 1
\end{bmatrix}.
$$
And multiplying on the left by $\bE_k$ will introduce zeros below the diagonal:
$$
\bE_k \bx =
\footnotesize
\begin{bmatrix}
	1 & \ldots & 0&  0 & \ldots & 0\\
	\vdots & \ddots & \vdots & \vdots & \ddots & \vdots \\
	0 & \ldots & 1  & 0 & \ldots & 0 \\
	0 & \ldots & -z_{k+1}  & 1 & \ldots & 0\\
	\vdots & \ddots & \vdots & \vdots & \ddots & \vdots \\
	0 & \ldots & -z_n & 0 & \ldots & 1
\end{bmatrix}
\begin{bmatrix}
	x_1 \\
	\vdots \\
	x_k\\
	x_{k+1} \\
	\vdots \\
	x_n
\end{bmatrix}=
\begin{bmatrix}
	x_1 \\
	\vdots \\
	x_k\\
0 \\
	\vdots \\
0
\end{bmatrix}.
$$
If we let $\bA_1 = \bA$ and $\bA_k$ containing elements $\{a_{ij}^{(k)}\}$ be the matrix obtained after $k-1$ elimination steps ($1\leq k<n$), then $\bA_k$ has the following form
$$
\bA_k=
\footnotesize
\begin{bmatrix}
a_{11}^{(k)} & \ldots & a_{1k}^{(k)} & \ldots & a_{1n}^{(k)}\\
\vdots & \ddots & \vdots & \ddots  & \vdots \\
0 & \ldots & a_{kk}^{(k)}  & \ldots &   a_{kn}^{(k)}\\
\vdots & \ddots & \vdots & \ddots  & \vdots \\
0 & \ldots & a_{nk}^{(k)} & \ldots & a_{nn}^{(k)}
\end{bmatrix}.
$$
In each stage, $\det(\bA_k) = \pm\det(\bA)$. And therefore, $\bA_k$ is nonsingular if $\bA$ is nonsingular.
\footnote{To prove that $\bA_k$ is nonsingular, we can demonstrate it recursively by showing that $\bA_k = \bE_{k-1}\bA_{k-1}$, where $\bE_{k-1}$ and $\bA_{k-1}$ are both nonsingular.} 

For example, we write out the Gaussian elimination steps for a $4\times 4$ matrix. For simplicity, we assume there are no row permutations. And in the following matrix, $\boxtimes$ represents a value that is not necessarily zero, and \textbf{boldface} indicates the value has just been changed. 
\begin{tcolorbox}[title={A Trivial Gaussian Elimination For a $4\times 4$ Matrix},colback=\mdframecolorTheorem]
\begin{equation}\label{equation:elmination-steps}
\footnotesize
\begin{sbmatrix}{\bA}
	\boxtimes & \boxtimes & \boxtimes & \boxtimes \\
	\boxtimes & \boxtimes & \boxtimes & \boxtimes \\
	\boxtimes & \boxtimes & \boxtimes & \boxtimes \\
	\boxtimes & \boxtimes & \boxtimes & \boxtimes
\end{sbmatrix}
\stackrel{\bE_1}{\longrightarrow}
\begin{sbmatrix}{\bE_1\bA}
	\boxtimes & \boxtimes & \boxtimes & \boxtimes \\
	\bm{0} & \textcolor{mylightbluetext}{\bm{\boxtimes}} & \bm{\boxtimes} & \bm{\boxtimes} \\
	\bm{0} & \bm{\boxtimes} & \bm{\boxtimes} & \bm{\boxtimes} \\
	\bm{0} & \bm{\boxtimes} & \bm{\boxtimes} & \bm{\boxtimes}
\end{sbmatrix}
\stackrel{\bE_2}{\longrightarrow}
\begin{sbmatrix}{\bE_2\bE_1\bA}
	\boxtimes & \boxtimes & \boxtimes & \boxtimes \\
	 0 & \textcolor{mylightbluetext}{\boxtimes} & \boxtimes & \boxtimes \\
	 0  & \bm{0} & \textcolor{mylightbluetext}{\bm{\boxtimes}} & \bm{\boxtimes} \\
      0 & \bm{0} & \bm{\boxtimes} & \bm{\boxtimes}
\end{sbmatrix}
\stackrel{\bE_3}{\longrightarrow}
\begin{sbmatrix}{\bE_3\bE_2\bE_1\bA}
	\boxtimes & \boxtimes & \boxtimes & \boxtimes \\
	0 &  \textcolor{mylightbluetext}{\boxtimes} & \boxtimes & \boxtimes \\
    0 & 0  & \textcolor{mylightbluetext}{\boxtimes} & \boxtimes \\
	0 & 0  & \bm{0} & \textcolor{mylightbluetext}{\bm{\boxtimes}}
\end{sbmatrix},
\end{equation}
\end{tcolorbox}
\noindent where $\bE_1, \bE_2$, and $\bE_3$ are lower triangular matrices. Specifically, as discussed above, Gaussian transformation matrices $\bE_i$'s are unit lower triangular matrices with $1$'s on the diagonal. This can also be explained as follows: for the $k$-th transformation $\bE_k$, applied to the matrix $\bE_{k-1}\ldots\bE_1\bA$, the transformation subtracts multiples of the $k$-th row from rows $\{k+1, k+2, \ldots, n\}$ to obtain zeros below the diagonal in the $k$-th column of the matrix. And  rows $\{1, 2, \ldots, k-1\}$ are not used in this process. 

To illustrate further, consider the transformation example at stage $1$. 
We left-multiply by $\bE_1$, subtracting multiples of the $1$-st row from rows $2, 3,$ and $4$, and setting the first entries of rows $2, 3,$ and $4$ to zero.  
Similar procedures apply for steps 2 and 3. 
By setting $\bL\triangleq\bE_1^{-1}\bE_2^{-1}\bE_3^{-1}$ and letting the matrix after elimination be $\bU$, \footnote{The inverses of unit lower triangular matrices are also unit lower triangular matrices. And the products of unit lower triangular matrices are also unit lower triangular matrices.} we obtain $\bA=\bL\bU$. Thus, an LU decomposition is derived for this $4\times 4$ matrix $\bA$.

\index{Pivot}
\index{Pivot variable}
\index{Free variable}
\begin{definition}[Pivot, Pivot Variables, Free Variables]\label{definition:pivot}
The first nonzero entry in the row after each elimination step is called a \textit{pivot}. For example, the \textcolor{mylightbluetext}{blue} crosses in Equation~\eqref{equation:elmination-steps} are pivots.

In the context of solving systems of linear equations, variables associated with pivot columns (columns containing pivot elements) are called ``\textit{pivot variables}." 
While a \textit{free variable} is a variable that does not correspond to a pivot column.
\end{definition}

\index{Echelon form}
\index{Row echelon form}
\index{Column echelon form}
\index{Reduced row echelon form}
\begin{definition}[Echelon Forms]\label{definition:row_echwlo_lu}
When the leading coefficient (the first nonzero entry) of each row of  a matrix is to the right of the leading coefficient of the row above it, we call the matrix is in \textit{row echelon form}.
This form makes it easy to solve the system by back substitution.
Similarly, a matrix is in \textit{column echelon form} if its transpose is in row echelon form.
Moreover, if a matrix is in row echelon form, then the first nonzero entry of each row is called a \textit{pivot} (Definition~\ref{definition:pivot}), and the columns containing these pivots are called \textit{pivot columns}.

In addition, a matrix is in \textit{reduced row echelon form} (RREF, a.k.a. \textit{row canonical form}) if it satisfies the following conditions:
\begin{enumerate}
\item It is in row echelon form;
\item The leading coefficient in each nonzero row is a 1 (pivots are 1, or called a leading 1);
\item Each column containing a learning 1 has zeros in all  other entries in that column;
\item All nonzero rows are placed at the bottom of the matrix.
\end{enumerate}

An important property of  echelon forms is that  all row echelon forms and the reduced row echelon form have the same number of zero rows, and the pivots are located in the same indices (albeit the row echelon form is not unique).
\end{definition}


However, there is a possibility that $a_{11}$ (the entry (1,1) of matrix $\bA$) may take on a value of zero at times. 
In such cases, no matrix $\bE_1$ can ensure the success of the subsequent elimination step.
Thus, it becomes necessary to interchange the first row and the second row (or some other rows) using a permutation matrix $\bP_1$. 
This process is commonly known as  \textit{pivoting} or simply \textit{permutation}:
\begin{tcolorbox}[title={Gaussian Elimination With a Permutation in the Beginning},colback=\mdframecolorTheorem]
$$
\scriptsize
\begin{aligned}
\begin{sbmatrix}{\bA}
	0 & \boxtimes & \boxtimes & \boxtimes \\
	\boxtimes & \boxtimes & \boxtimes & \boxtimes \\
	\boxtimes & \boxtimes & \boxtimes & \boxtimes \\
	\boxtimes & \boxtimes & \boxtimes & \boxtimes
\end{sbmatrix}
&\stackrel{\bP_1}{\longrightarrow}
\begin{sbmatrix}{\bP_1\bA}
	\bm{\boxtimes} & \bm{\boxtimes} & \bm{\boxtimes} & \bm{\boxtimes} \\
	\bm{0} & \bm{\boxtimes} & \bm{\boxtimes} & \bm{\boxtimes} \\
	\boxtimes & \boxtimes & \boxtimes & \boxtimes \\
	\boxtimes & \boxtimes & \boxtimes & \boxtimes
\end{sbmatrix}
\stackrel{\bE_1}{\longrightarrow}
\begin{sbmatrix}{\bE_1\bP_1\bA}
	\boxtimes & \boxtimes & \boxtimes & \boxtimes \\
	\bm{0} & \textcolor{mylightbluetext}{\bm{\boxtimes}} & \bm{\boxtimes} & \bm{\boxtimes} \\
	\bm{0} & \bm{\boxtimes} & \bm{\boxtimes} & \bm{\boxtimes} \\
	\bm{0} & \bm{\boxtimes} & \bm{\boxtimes} & \bm{\boxtimes}
\end{sbmatrix}
\stackrel{\bE_2}{\longrightarrow}
\begin{sbmatrix}{\bE_2\bE_1\bP_1\bA}
	\boxtimes & \boxtimes & \boxtimes & \boxtimes \\
	0 & \textcolor{mylightbluetext}{\boxtimes} & \boxtimes & \boxtimes \\
	0  & \bm{0} & \textcolor{mylightbluetext}{\bm{\boxtimes}} & \bm{\boxtimes} \\
	0 & \bm{0} & \bm{\boxtimes} & \bm{\boxtimes}
\end{sbmatrix}
\stackrel{\bE_3}{\longrightarrow}
\begin{sbmatrix}{\bE_3\bE_2\bE_1\bP_1\bA}
	\boxtimes & \boxtimes & \boxtimes & \boxtimes \\
	0 &  \textcolor{mylightbluetext}{\boxtimes} & \boxtimes & \boxtimes \\
	0 & 0  & \textcolor{mylightbluetext}{\boxtimes} & \boxtimes \\
	0 & 0  & \bm{0} & \textcolor{mylightbluetext}{\bm{\boxtimes}}
\end{sbmatrix}.
\end{aligned}
$$
\end{tcolorbox}
\noindent By setting $\bL\triangleq\bE_1^{-1}\bE_2^{-1}\bE_3^{-1}$ and $\bP\triangleq\bP_1^{-1}$, we obtain $\bA=\bP\bL\bU$. Therefore, we obtain a full LU decomposition with permutation for this $4\times 4$ matrix $\bA$.

In certain scenarios, additional permutation matrices, such as $\bP_2, \bP_3, \ldots$, may be interspersed among the lower triangular  $\bE_i$'s. 
An example is illustrated as follows:
\begin{tcolorbox}[title={Gaussian Elimination With a Permutation in Between},colback=\mdframecolorTheorem]
$$
\footnotesize
\begin{sbmatrix}{\bA}
	\boxtimes & \boxtimes & \boxtimes & \boxtimes \\
	\boxtimes & \boxtimes & \boxtimes & \boxtimes \\
	\boxtimes & \boxtimes & \boxtimes & \boxtimes \\
	\boxtimes & \boxtimes & \boxtimes & \boxtimes \\
\end{sbmatrix}
\stackrel{\bE_1}{\longrightarrow}
\begin{sbmatrix}{\bE_1\bA}
	\boxtimes & \boxtimes & \boxtimes & \boxtimes \\
	\bm{0} & \bm{0} & \bm{\boxtimes} & \bm{\boxtimes} \\
	\bm{0} & \bm{\boxtimes} & \bm{\boxtimes} & \bm{\boxtimes} \\
	\bm{0} & \bm{\boxtimes} & \bm{\boxtimes} & \bm{\boxtimes}
\end{sbmatrix}
\stackrel{\bP_1}{\longrightarrow}
\begin{sbmatrix}{\bP_1\bE_1\bA}
	\boxtimes & \boxtimes & \boxtimes & \boxtimes \\
	\bm{0} & \textcolor{mylightbluetext}{\bm{\boxtimes}} & \bm{\boxtimes} & \bm{\boxtimes} \\
	\bm{0}  & \bm{0} & \textcolor{mylightbluetext}{\bm{\boxtimes}} & \bm{\boxtimes} \\
	0 & \boxtimes & \boxtimes & \boxtimes
\end{sbmatrix}
\stackrel{\bE_2}{\longrightarrow}
\begin{sbmatrix}{\bE_2\bP_1\bE_1\bA}
	\boxtimes & \boxtimes & \boxtimes & \boxtimes \\
	0 &  \textcolor{mylightbluetext}{\boxtimes} & \boxtimes & \boxtimes \\
	0 & 0  & \textcolor{mylightbluetext}{\boxtimes} & \boxtimes \\
	0 & \bm{0}  & \bm{0} & \textcolor{mylightbluetext}{\bm{\boxtimes}}
\end{sbmatrix}.
$$
\end{tcolorbox}
\noindent In this case, we find $\bU=\bE_2\bP_1\bE_1\bA$. In Sections~\ref{section:lu-perm},~\ref{sec:compute-lu-pivoting},  or~\ref{section:partial-pivot-lu}, we will demonstrate that the permutations in-between will result in the same representation $\bA=\bP\bL\bU$, where $\bP$ takes account of all the permutations.

The above examples can be easily extended to any $n\times n$ matrix if we assume there are no row permutations in the process. And we will have $n-1$ such lower triangular transformations.
The $k$-th transformation $\bE_k$ introduces zeros below the diagonal in the $k$-th column of $\bA$ by subtracting multiples of the $k$-th row from rows $\{k+1, k+2, \ldots, n\}$.
Finally, by setting $\bL\triangleq\bE_1^{-1}\bE_2^{-1}\ldots \bE_{n-1}^{-1}$, we obtain the LU decomposition $\bA=\bL\bU$ (without permutation).

\subsection*{Properties of Elementary Transformation}
From the above examples for row elementary operations, we conclude the following result.
\begin{proposition}[Row Space after Row Operations]\label{proposition:rowspa_rowele}
Given a matrix $\bA\in\real^{m\times n}$, operated by a sequence of row transformations $\bE_1, \bE_2, \ldots, \bE_k$, and let $\bE = \bE_k\bE_{k-1}\ldots \bE_1$. Then, the row space of $\bB=\bE\bA$ is identical to the row space of $\bA$. 
\end{proposition}
\begin{proof}[of Proposition~\ref{proposition:rowspa_rowele}]
	Since the rows of $\bB$ are linear combinations of the rows of $\bA$, it follows that $\cspace(\bB^\top)\subseteq \cspace(\bA^\top)$. 
Moreover, since  the row transformations are invertible, $\bE= \bE_k\bE_{k-1}\ldots \bE_1$ is also invertible. Therefore, $\bA = {\bE^{-1}}\bB$, which implies that the rows of $\bA$ are also  linear combinations of the rows of $\bB$: $\cspace(\bA^\top)\subseteq \cspace(\bB^\top)$. 
Thus, $\cspace(\bA^\top)= \cspace(\bB^\top)$.
\end{proof}
Note, however, the column space of $\bA$ and $\bB$ may be different. But since the dimension of the row space is equal to that of the column space (Theorem~\ref{theorem:equal-dimension-rank}, proof in Appendix~\ref{append:row-equal-column}), the dimensions of the column spaces of $\bA$ and $\bB$ are the same.

\section{Existence of  LU Decomposition without Permutation}\label{section:exist-lu-without-perm}
Gaussian elimination, or Gaussian transformation, reveals the origin of the LU decomposition. Subsequently, we rigorously establish the validity of Theorem~\ref{theorem:lu-factorization-without-permutation} through an inductive proof, affirming the existence of the LU decomposition without permutation.

\begin{proof}[{of Theorem~\ref{theorem:lu-factorization-without-permutation}: LU Decomposition without Permutation}]
We will establish, through induction, that any square matrix $\bA$ of order $n$ with nonzero leading principal minors admits the LU decomposition $\bA=\bL\bU$.
The case for $1\times 1$ matrices is trivial, as we set $L=1$ and $U=A$, yielding $A=LU$.

Suppose, the inductive step, that any $k\times k$ matrix $\bA_k$ with all the leading principal minors being nonzero has an LU decomposition without permutation. If we can prove that any $(k+1)\times(k+1)$ matrix $\bA_{k+1}$ can also be factored as this LU decomposition without permutation, then we complete the proof.

Consider any $(k+1)\times(k+1)$ matrix $\bA_{k+1}$. Suppose the $k$-th order leading principal submatrix of $\bA_{k+1}$ is $\bA_k$ with size $k\times k$. Then, by the inductive hypothesis, $\bA_k$ can be factored as $\bA_k =  \bL_k\bU_k$, where $\bL_k$ is a unit lower triangular matrix and $\bU_k$ is a nonsingular upper triangular matrix based on the assumption. We can express $\bA_{k+1}$ as
$
\bA_{k+1} = \begin{bmatrix}
	\bA_k & \bb \\
	\bc^\top & d
\end{bmatrix}.
$
Then, it can be factored as:
$$
\bA_{k+1} = \begin{bmatrix}
	\bA_k & \bb \\
	\bc^\top & d
\end{bmatrix}
=
\begin{bmatrix}
	\bL_k &\bzero \\
	\bx^\top  & 1 
\end{bmatrix}
\begin{bmatrix}
	\bU_k & \by\\
	\bzero & z
\end{bmatrix} = \bL_{k+1}\bU_{k+1},
$$
where $\bb = \bL_k\by$, $\bc^\top = \bx^\top\bU_k$, $d = \bx^\top\by + z$, $\bL_{k+1}=\begin{bmatrix}
	\bL_k &\bzero \\
	\bx^\top  & 1 
\end{bmatrix}$, and $\bU_{k+1}=\begin{bmatrix}
\bU_k & \by\\
\bzero & z
\end{bmatrix}$. 
From the assumption, $\bL_k$ and $\bU_k$ are nonsingular. Therefore,
$$
\by = \bL_k^{-1}\bb, \qquad \bx^\top=\bc^\top\bU_k^{-1}, \qquad  z=d - \bx^\top\by.
$$
To complete the proof, we need to show that $z$ is nonzero such that $\bU_{k+1}$ is nonsingular.

Given that all the leading principal minors of $\bA_{k+1}$ are nonzero, 
we have $\det(\bA_{k+1})=$
\footnote{By the fact that if matrix $\bM$ has a block formulation: $\bM=\begin{bmatrix}
		\bA & \bB \\
		\bC & \bD 
	\end{bmatrix}$, then $\det(\bM) = \det(\bA)\det(\bD-\bC\bA^{-1}\bB)$.}
$\det(\bA_k)\cdot$ $\det(d-\bc^\top\bA_k^{-1}\bb)  
=\det(\bA_k)\cdot(d-\bc^\top\bA_k^{-1}\bb) \neq 0$ since $d-\bc^\top\bA_k^{-1}\bb$ is a scalar.
Given  $\det(\bA_k)\neq 0$ from the assumption, we obtain $d-\bc^\top\bA_k^{-1}\bb \neq 0$. Substituting $\bb = \bL_k\by$ and $\bc^\top = \bx^\top\bU_k$ into the formula, we obtain $d-\bx^\top\bU_k\bA_k^{-1}\bL_k\by =d-\bx^\top\bU_k(\bL_k\bU_k)^{-1}\bL_k\by =d-\bx^\top\by \neq 0$, which is exactly the form of $z\neq 0$. Thus, we find $\bL_{k+1}$ with all diagonal values being 1, and $\bU_{k+1}$ with all the values on the diagonal being nonzero, signifying that $\bL_{k+1}$ and $\bU_{k+1}$ are nonsingular. \footnote{A triangular matrix (upper or lower) is nonsingular if and only if all the entries on its main diagonal are nonzero.}
This completes the proof.
\end{proof}

We further establish that the LU decomposition is unique when no permutations are involved.
\begin{corollary}[Uniqueness of the LU Decomposition without Permutation]\label{corollary:unique-lu-without-permutation}
Let $\bA$ be any $n\times n$ square matrix. If $\bA$ has nonzero leaning principal minors, then the LU decomposition is unique.
\end{corollary}
\begin{proof}[of Corollary~\ref{corollary:unique-lu-without-permutation}]
Suppose the LU decomposition is not unique, then we can find two decompositions such that $\bA=\bL_1\bU_1 = \bL_2\bU_2$, which implies $\bL_2^{-1}\bL_1=\bU_2\bU_1^{-1}$. 
The left side of the equation, $\bL_2^{-1}\bL_1$,  is a unit lower triangular matrix, and the right side of the equation, $\bU_2\bU_1^{-1}$, is an upper triangular matrix. 
This implies that both sides of the above equation must be diagonal matrices. 
Since the inverse of a unit lower triangular matrix is also a unit lower triangular matrix, and the product of unit lower triangular matrices is also a unit lower triangular matrix, we conclude that $\bL_2^{-1}\bL_1 = \bI$.
The equality implies that both sides are identity matrices, leading to $\bL_1=\bL_2$ and $\bU_1=\bU_2$, which leads to a contradiction. 
\end{proof}

In the proof of Theorem~\ref{theorem:lu-factorization-without-permutation}, we demonstrated that the upper triangular matrix has nonzero  diagonal values  when all the leading principal minors of $\bA$ are  nonzero. 
An alternative formulation arises by dividing each row of $\bU$ by its corresponding diagonal value, resulting in the \textit{LDU decomposition}.
\begin{corollaryHigh}[LDU Decomposition]\label{corollary:ldu-decom}
Let $\bA$ be  any $n\times n$ square matrix. If all the leading principal minors are nonzero, i.e., $\det(\bA[1:k,1:k])\neq 0$, for all $k\in \{1,2,\ldots, n\}$, then $\bA$ can be \textbf{uniquely} factored as 
\begin{equation}
	\bA = \bL\bD\bU, \nonumber
\end{equation}
where $\bL$ is a unit lower triangular matrix, $\bU$ is a \textbf{unit} upper triangular matrix, and $\bD$ is a diagonal matrix. 
\end{corollaryHigh} 
The proof is straightforward. Given the LU decomposition of  $\bA=\bL\bR$, we can identify a diagonal matrix $\bD=\diag(r_{11}, r_{22}, \ldots, r_{nn})$ such that $\bD^{-1}\bR = \bU$ forms a unit upper triangular matrix. And the uniqueness follows from the uniqueness of the LU decomposition.

\index{Schur complement}
\section{Schur Complement and Matrix Inverse Decomposition}\label{section:schur-complement}
Before delving into the existence of LU decomposition with permutation, we first discuss the \textit{Schur complement} of a matrix and its decompositions.
Let $\bM$ be an $n\times n$ square matrix with the following $2\times 2$ block matrix format:
\begin{equation}\label{equation:schu_set}
\bM = \begin{bmatrix}
\underset{p\times p}{\bA} & \underset{p\times q}{\bB} \\
\underset{q\times p}{\bC} & \underset{q\times q}{\bD}
\end{bmatrix},
\end{equation}
where $\bA$ is a $p\times p$ matrix and $\bD$ is a $q\times q$ matrix with $n=p+q$. And it is clear that $\bB$ is a $p\times q$ matrix and $\bC$ is a $q\times p$ matrix. We have  the following \textit{Schur complement}:
\begin{itemize}
\item If $\bD$ is invertible, $\Delta_{\bD} = \bA -\bB \bD^{-1}\bC$ is called the Schur complement of $\bD$ in $\bM$;~\footnote{When $q=1$, the determinant of $\bM$ can be obtained by Cauchy expansion; see~\eqref{equation:det_block_rem2} and Problem~\ref{prob:cauch_exp}.}

\item If $\bA$ is invertible, $\Delta_{\bA} = \bD - \bC\bA^{-1}\bB$ is called the Schur complement of $\bA$ in $\bM$.
\end{itemize}
More generally, consider an index set $I\subseteq \{1,2,\ldots,n\}$ and its complementary set $J=\{1,2,\ldots,n\}/I$. Then, 
\begin{itemize}
\item If $\bM[I,I]$ is invertible, $\Delta_{I}=\bM[J,J] - \bM[J,I]\bM[I,I]^{-1}\bA[I,J]$ is called the Schur complement of $\bM[I,I]$ in $\bM$.
\end{itemize}

Two special Schur complements will be used extensively in the following sections.
\begin{remark}[Schur Complement of $m_{11}$ and $m_{nn}$]
For any matrix $\bM\in\real^{n\times n}$, 
let the entry $(1,1)$ of $\bM$ be $m_{11}$ and assume $m_{11}$ is not zero, then 
$$
\Delta_{m_{11}} = \bM_{2:n,2:n} -\frac{1}{m_{11}} \bM_{2:n,1}\bM_{1,2:n} \in \real^{(n-1)\times (n-1)}
$$ 
is called the Schur complement of $m_{11}$ in $\bM$. 
Similarly, when $m_{nn}\neq 0$, the Schur complement of  $m_{nn}$ in $\bM$ is given by 
$$
\Delta_{m_{nn}} =\bM_{1:n-1,1:n-1} - \frac{1}{m_{nn}}\bM_{1:n-1,n} \bM_{n,1:n-1}  \in \real^{(n-1)\times (n-1)}.
$$
\end{remark}

When certain blocks in a matrix are nonsingular, the block matrix can be decomposed using Schur complements.
\begin{lemma}[Factorization in Schur Complements]\label{lemma:fac_schur_a}
Consider the setting in \eqref{equation:schu_set}. If $\bA$ is invertible, then $\bM$ can be factored as 
$$
\footnotesize
\begin{aligned}
\bM 
=
\begin{bmatrix}
\bA &\bB \\
\bC & \bD
\end{bmatrix}
&=
\begin{bmatrix}
\bI & \bzero \\
-\bC\bA^{-1} & \bI 
\end{bmatrix}^{-1}
\begin{bmatrix}
\bA & \bzero \\
\bzero & \Delta_{\bA}
\end{bmatrix}
\begin{bmatrix}
\bI &-\bA^{-1}\bB \\
\bzero & \bI
\end{bmatrix}^{-1}
=
\begin{bmatrix}
\bI & \bzero \\
\bC\bA^{-1} & \bI 
\end{bmatrix}
\begin{bmatrix}
\bA & \bzero \\
\bzero & \Delta_{\bA}
\end{bmatrix}
\begin{bmatrix}
\bI &\bA^{-1}\bB \\
\bzero & \bI
\end{bmatrix},
\end{aligned}
$$
which indicates that $\bM$ can be written as a product of a (block) upper triangular, a (block) diagonal matrix, and a (block) lower triangular matrix. 
Similarly, if  $\bD$  is invertible, then $\bM$ can be factored as 
$$
\setlength{\arraycolsep}{4.5pt}
\footnotesize
\begin{aligned}
\bM 
=
\begin{bmatrix}
\bA &\bB \\
\bC & \bD
\end{bmatrix}
&=
\begin{bmatrix}
\bI & -\bB\bD^{-1} \\
\bzero & \bI 
\end{bmatrix}^{-1}
\begin{bmatrix}
\Delta_{\bD} & \bzero \\
\bzero & \bD
\end{bmatrix}
\begin{bmatrix}
\bI & \bzero \\
-\bD^{-1}\bC & \bI
\end{bmatrix}^{-1}
=
\begin{bmatrix}
\bI & \bB\bD^{-1} \\
\bzero & \bI
\end{bmatrix}
\begin{bmatrix}
\Delta_{\bD} & \bzero \\
\bzero & \bD
\end{bmatrix}
\begin{bmatrix}
\bI & \bzero \\
\bD^{-1}\bC & \bI 
\end{bmatrix}
.
\end{aligned}
$$
\end{lemma}

\begin{proof}[of Lemma~\ref{lemma:fac_schur_a}]
If $\bA$ is invertible, by mimicking the process of Gaussian elimination (Section~\ref{section:gaussian-elimination}) on the block matrix, we can lower triangularize and upper triangularize $\bM$ as follows (called \textit{block Gaussian elimination}, which is unique; see Problem~\ref{problem:block_gaus_uni}):
\begin{equation}\label{equation:block_gauss}
\footnotesize
\begin{aligned}
\begin{bmatrix}
\bI & \bzero \\
-\bC\bA^{-1} & \bI 
\end{bmatrix}
\begin{bmatrix}
\bA &\bB \\
\bC & \bD
\end{bmatrix}
&=
\begin{bmatrix}
\bA &\bB \\
\bzero & \Delta_{\bA}
\end{bmatrix}
\text{ and }
\begin{bmatrix}
\bA &\bB \\
\bC & \bD
\end{bmatrix}
\begin{bmatrix}
\bI &-\bA^{-1}\bB \\
\bzero & \bI
\end{bmatrix}
&=
\begin{bmatrix}
\bA &\bzero \\
\bC & \Delta_{\bA}
\end{bmatrix}.
\end{aligned}
\end{equation}
Then, multiplying $\bM$ from the left by the (block) lower triangular matrix and from the right by the (block) upper triangular matrix, we get:
$$
\footnotesize
\begin{bmatrix}
\bI & \bzero \\
-\bC\bA^{-1} & \bI 
\end{bmatrix}
\begin{bmatrix}
\bA &\bB \\
\bC & \bD
\end{bmatrix}
\begin{bmatrix}
\bI &-\bA^{-1}\bB \\
\bzero & \bI
\end{bmatrix}
=
\begin{bmatrix}
\bA & \bzero \\
\bzero & \Delta_{\bA}
\end{bmatrix}.
$$
Using the following fact, we obtain the result
$$
\footnotesize
\begin{bmatrix}
\bI & \bzero \\
-\bC\bA^{-1} & \bI 
\end{bmatrix}
\begin{bmatrix}
\bI & \bzero \\
\bC\bA^{-1} & \bI 
\end{bmatrix} = \bI
\qquad\text{and} \qquad 
\begin{bmatrix}
\bI &-\bA^{-1}\bB \\
\bzero & \bI
\end{bmatrix}
\begin{bmatrix}
\bI &\bA^{-1}\bB \\
\bzero & \bI
\end{bmatrix}
=\bI.
$$
The second part can be proved similarly.
\end{proof}

The determinant of block matrices (Remark~\ref{remark:determinant-intermezzo}) shows that $\det(\Delta_{\bA}) = \det(\bM)/\det(\bA)$.
And the rank of $\Delta_{\bA}$ can be determined as follows.
\begin{exercise}[Rank of $\Delta_{\bA}$]
Show that $\rank(\bM) = \rank(\bA) +\rank(\Delta_{\bA})$.
Thus, $\rank(\bM) = \rank(\bA) $ if and only if $\Delta_{\bA}=\bzero$, and $\bM$ is nonsingular if and only if $\Delta_{\bA}$ is nonsingular (when $\bA$ is assumed nonsingular).
\end{exercise}

It is then straightforward to compute the inverse of $\bM$ since the inversse of the (block) upper triangular, the (block) diagonal matrix, and the (block) lower triangular matrix are relatively easy. If we assume $\Delta_{\bA}$ is invertible and all the relevant inverses exist, then we have 
\begin{equation}\label{equation:schur_inv1}
\setlength{\arraycolsep}{4pt}
\footnotesize
\begin{aligned}
\bM^{-1}
&
=
\begin{bmatrix}
\bI &-\bA^{-1}\bB \\
\bzero & \bI
\end{bmatrix}
\begin{bmatrix}
\bA^{-1} & \bzero \\
\bzero & \Delta_{\bA}^{-1}
\end{bmatrix}
\begin{bmatrix}
\bI & \bzero \\
-\bC\bA^{-1} & \bI 
\end{bmatrix}
=
\begin{bmatrix}
\bA^{-1} + \bA^{-1}\bB\Delta_{\bA}^{-1}\bC\bA^{-1}& -\bA^{-1}\bB\Delta_{\bA}^{-1}\\
-\Delta_{\bA}^{-1}\bC\bA^{-1} & \Delta_{\bA}^{-1}
\end{bmatrix}.
\end{aligned}
\end{equation}
And if $\Delta_{\bD}$ is invertible and all the relevant inverses exist, the inverse of $\bM$ is
\noindent
\begin{equation}\label{equation:schur_inv2}
\setlength{\arraycolsep}{3.55pt}
\footnotesize
\begin{aligned}
\bM^{-1}
&=
\begin{bmatrix}
\bI & \bzero \\
-\bD^{-1}\bC & \bI 
\end{bmatrix}
\begin{bmatrix}
\Delta_{\bD}^{-1} & \bzero \\
\bzero & \bD^{-1}
\end{bmatrix}
\begin{bmatrix}
\bI & -\bB\bD^{-1} \\
\bzero & \bI
\end{bmatrix}
=
\begin{bmatrix}
\Delta_{\bD}^{-1} & -\Delta_{\bD}^{-1}\bB\bD^{-1} \\
-\bD^{-1}\bC \Delta_{\bD}^{-1} & \bD^{-1}+\bD^{-1}\bC \Delta_{\bD}^{-1}\bB\bD^{-1}
\end{bmatrix}.
\end{aligned}
\end{equation}
\begin{exercise}[Invertibility of $\bM$, $\bA$, $\bD$, $\Delta_{\bA}$, and $\Delta_{\bD}$]\label{exercise:inv_m}
Equation~\eqref{equation:schur_inv1} and \eqref{equation:schur_inv2} present the formulas for the inverse of $\bM$ if the relevant inverses exist.
To be more specific,
\begin{itemize}
\item Suppose $\bA$ is invertible, show that  $\bM$ is invertible if and only if  $\Delta_{\bA}$ is invertible.
\item Suppose $\bD$ is invertible, show that  $\bM$ is invertible if and only if $\Delta_{\bD}$ is invertible.
\end{itemize}
This also implies if both $\bA$ and $\bD$ are invertible, then if any of $(\bM, \Delta_{\bA}, \Delta_{\bD})$ is invertible, all the three matrices are invertible.
\end{exercise}

\begin{remark}[Invertibility of $\bM$]
Exercise~\ref{exercise:inv_m} shows $\bM$ is invertible if $\bA$ or $\bD$ is invertible with some mild conditions.
However, the invertibility of $\bA$ or $\bD$ is not a necessary condition. For example, $\bM=\bM^{-1}=\scriptsize\begin{bmatrix}
0 & 1 \\
1 & 0
\end{bmatrix}$
is invertible, while $\bA$ and $\bD$ are not invertible in this case.
\end{remark}

A more general block of the inverse of $\bM$ is discussed in Problem~\ref{prob:mt_inv_gen}.
Equation~\eqref{equation:schur_inv1} also shows (by Remark~\ref{remark:determinant-intermezzo}):
\begin{equation}\label{equation:schur_inv1_det}
\bM^{-1}[p+1:n, p+1:n] = \Delta_{\bA}^{-1}
\implies
\det\big(\bM^{-1}[p+1:n, p+1:n]\big) = \frac{\det(\bA)}{\det(\bM)},
\end{equation}
where the latter is a special form of \textit{Jacobi's equality} $\det\big(\bM^{-1}[J,J]\big)=\det(\bM[I,I])/\det(\bM)$ (Problem~\ref{prob:cramer_adj_5}, where we prove it using the definitions of determinants and adjugate).
\begin{exercise}[Jacobi's Equality]\label{exercise:jacob_eq}
Prove  Jacobi's equality (see Problem~\ref{prob:cramer_adj_5})  using the Schur complement of $\bM[I,I]$ in $\bM$.
\end{exercise}

%
%

Moreover, by writing out the two inverse formulas using  $\Delta_{\bA}$ and $\Delta_{\bD}$, setting $\bD = \bI$, and changing $\bB$ to $-\bB$, we get: 
\begin{equation}
(\bA +\bB\bC)^{-1} = \bA^{-1} - \bA^{-1} \bB (\bI-\bC\bA^{-1}\bB)^{-1}\bC\bA^{-1},
\end{equation}
which is known as the \textit{matrix inversion lemma} \citep{boyd2004convex, gallier2010schur}.
More generally, we have the following results.
\begin{lemma}[Matrix Inversion Lemma]\label{lemma:matr_inv_ap}
Consider the setting in \eqref{equation:schu_set}. 
Suppose both $\bA$ and $\bD$ are invertible.
Since the inverse of a matrix is unique, 
combining Equation~\eqref{equation:schur_inv1} and \eqref{equation:schur_inv2}, if $\bM$ is invertible, we conclude the following identities
$$
\begin{aligned}
\Delta_{\bD}^{-1} &= \bA^{-1} + \bA^{-1}\bB\Delta_{\bA}^{-1}\bC\bA^{-1}
\gap \text{and} \gap 
\Delta_{\bA}^{-1} &=\bD^{-1}+\bD^{-1}\bC \Delta_{\bD}^{-1}\bB\bD^{-1}.
\end{aligned}
$$

\end{lemma}

Various identities can be demonstrated by applying the matrix inversion lemma.

\paragraph{Woodbury matrix identity.} Based on Lemma~\ref{lemma:matr_inv_ap}, if we further change $\bC$ to $-\bC$ and substitute $\bD$ with $\bD^{-1}$, we obtain the \textit{Woodbury matrix identity}:
\begin{equation}\label{equation:schur_inv1_det_sour}
(\bA+\bB\bD\bC)^{-1} = \bA^{-1} - \bA^{-1} \bB(\bD^{-1} + \bC\bA^{-1}\bB)^{-1}\bC\bA^{-1}.
\end{equation}
Therefore, if $\bA\in\real^{p\times p}$, $\bB\in\real^{p\times q}$, $\bC\in\real^{q\times p}$, and $\bD\in\real^{q\times q}$ with $p>q$. $\bB\bD\bC$ is thus a low-rank perturbation \visavi $\bA$. Thus, considerable speedup can be achieved to obtain the inverse of the perturbed matrix.
Similarly, for the determinants:
\begin{equation}\label{equation:schur_inv1_det_res}
\abs{\bA+\bB\bD\bC}=
\abs{\bA}\cdot \abs{\bD}\cdot\abs{\bD^{-1} + \bC\bA^{-1}\bB}.
\end{equation}

\paragraph{Sherman-Morrison formula.}
Similarly, based on Lemma~\ref{lemma:matr_inv_ap}, if we further set $q=1$, i.e., $\bA\in\real^{p\times p}$, $\bB=\bb\in\real^{p\times 1}$, $\bC=\bc^\top\in\real^{1\times p}$, and $\bD=-1$, we obtain the \textit{Sherman-Morrison formula}:
\begin{equation}
(\bA+\bb\bc^\top)^{-1} = \bA^{-1} - \frac{\bA^{-1}\bb\bc^\top\bA^{-1}}{1+\bc^\top\bA^{-1}\bb},
\end{equation}
if $(\bA+\bb\bc^\top)$ is invertible.

\section{Existence of  LU Decomposition with Permutation}\label{section:lu-perm}
In Theorem~\ref{theorem:lu-factorization-without-permutation}, we require that the matrix $\bA$ must have nonzero leading principal minors. 
However, this condition is not obligatory. 
Even if the leading principal minors are zero, nonsingular matrices can still undergo an LU decomposition, albeit with an additional permutation. The proof remains grounded in the principle of induction.
\begin{proof}[{of Theorem~\ref{theorem:lu-factorization-with-permutation}: LU Decomposition with Permutation}]
We start by noting that any $1\times 1$ nonsingular matrix has a full LU decomposition (i.e., with permutation), expressed as $A=PLU$ by simply setting $P=1$, $L=1$, and $U=A$.
We will demonstrate that if every $(n-1)\times (n-1)$ nonsingular matrix has a full LU decomposition, then the same holds for every $n\times n$ nonsingular matrix. 
By induction, we establish that every nonsingular matrix has a full LU decomposition.

The proof will be structured as follows: Firstly, we establish that \textcolor{black}{If $\bA$ is nonsingular, then its row-permuted matrix $\bB$ is also nonsingular}. 
Subsequently, we demonstrate that the \textit{Schur complement} of $b_{11}$ in $\bB$ is also nonsingular. 
Finally, we leverage this property to  formulate the decomposition of $\bA$ in terms of $\bB$.

We notice that at least one element in the first column of $\bA$ must be nonzero, otherwise $\bA$ will be singular. 
Consequently, we can employ a row permutation to ensure that the element in entry $(1,1)$ becomes nonzero. That is, there exists a permutation $\bP_1$ such that $\bB = \bP_1 \bA$ and $b_{11} \neq 0$. Since $\bA$ and $\bP_1$ are both nonsingular, and the product of nonsingular matrices is also nonsingular, we conclude that $\bB$ is also nonsingular.

\begin{mdframed}[hidealllines=\mdframehidelineNote,backgroundcolor=\mdframecolor,frametitle={Schur complement of $\bB$ is also nonsingular:}]
Now, we consider the Schur complement of $b_{11}$ in $\bB$ with size $(n-1)\times (n-1)$:
$$
\widehat{\bB} = \bB_{2:n,2:n} -\frac{1}{b_{11}} \bB_{2:n,1} \bB_{1,2:n}.
$$
Suppose there is an $(n-1)$-vector $\bx$ that satisfies
\begin{equation}\label{equ:lu-pivot1}
	\widehat{\bB} \bx = \bzero.
\end{equation}
Then, $\bx$ and $y=-\frac{1}{b_{11}}\bB_{1,2:n} \cdot \bx  $ satisfy
$$
\bB 
\left[
\begin{matrix}
	y \\
	 \bx
\end{matrix}
\right]
=
\left[
\begin{matrix}
	b_{11} & \bB_{1,2:n} \\
	\bB_{2:n,1} & \bB_{2:n,2:n}
\end{matrix}
\right]
\left[
\begin{matrix}
	y \\
	 \bx
\end{matrix}
\right]
=
\left[
\begin{matrix}
	0 \\
	\bzero
\end{matrix}
\right].
$$
Since $\bB$ is nonsingular, $\bx$ and $y$ must be zero (its null space is of dimension 0). 
Therefore, Equation~\eqref{equ:lu-pivot1} holds only when $\bx=\bzero$, which means that the null space of $\widehat{\bB}$ is of dimension 0. Consequently, $\widehat{\bB}$ is nonsingular with size $(n-1)\times(n-1)$.
\end{mdframed}

By the induction assumption,  any $(n-1)\times(n-1)$ nonsingular matrix can be factorized as the full LU decomposition form:
$
\widehat{\bB} = \bP_2\bL_2\bU_2.
$
We can then factor $\bA$ as
\begin{equation*}
\begin{aligned}
\bA &= \bP_1^\top  
\left[
\begin{matrix}
b_{11} & \bB_{1,2:n} \\
\bB_{2:n,1} & \bB_{2:n,2:n}
\end{matrix}
\right]
= \bP_1^\top  
\left[
\begin{matrix}
1 & \bzero \\
\bzero & \bP_2
\end{matrix}
\right] 
\left[
\begin{matrix}
b_{11} & \bB_{1,2:n} \\
\bP_2^\top \bB_{2:n,1} &\bP_2^\top \bB_{2:n,2:n}
\end{matrix}
\right]\\
&= \bP_1^\top  
\left[
\begin{matrix}
1 & \bzero \\
\bzero & \bP_2
\end{matrix}
\right] 
\left[
\begin{matrix}
b_{11} & \bB_{1,2:n} \\
\bP_2^\top \bB_{2:n,1} & \textcolor{mylightbluetext}{\bL_2\bU_2}+\bP_2^\top \textcolor{mylightbluetext}{\frac{1}{b_{11}} \bB_{2:n,1} \bB_{1,2:n}}
\end{matrix}
\right]\\
&= \bP_1^\top  
\left[
\begin{matrix}
1 & \bzero \\
\bzero & \bP_2
\end{matrix}
\right] 
\left[
\begin{matrix}
1 & \bzero \\
\frac{1}{b_{11}}\bP_2^\top \bB_{2:n,1} & \bL_2
\end{matrix}
\right] 
\left[
\begin{matrix}
b_{11} & \bB_{1,2:n} \\
\bzero & \bU_2
\end{matrix}
\right].\\
\end{aligned}
\end{equation*}
Therefore, we find the full LU decomposition of $\bA=\bP\bL\bU$ by defining
$$
\bP \triangleq \bP_1^\top  
\left[
\begin{matrix}
	1 & 0 \\
	0 & \bP_2
\end{matrix}
\right], \qquad
\bL\triangleq\left[
\begin{matrix}
	1 & 0 \\
	\frac{1}{b_{11}}\bP_2^\top \bB_{2:n,1} & \bL_2
\end{matrix}
\right], \qquad
\bU\triangleq
\left[
\begin{matrix}
	b_{11} & \bB_{1,2:n} \\
	\bzero & \bU_2
\end{matrix}
\right],
$$
from which the result follows. We will formulate this process into Algorithm~\ref{alg:lu-with-pivoting} to compute this decomposition.
\end{proof}

\section{Computing  LU without Pivoting Recursively: A=LU}\label{section:compute-lu-without-pivot}
As a starting point, we will consider the most common and simplest case, which does not involve  row exchanges:
$
\bA = \bL\bU, 
$
where $\bL$ is a unit lower triangular matrix and $\bU$ is a nonsingular upper triangular matrix. We refer to this decomposition as the LU decomposition without pivoting (or without permutation).

In Section~\ref{section:gaussian-elimination}, we discussed the connection between  LU decomposition and  Gaussian elimination. 
The LU decomposition of a matrix can be found by first applying Gaussian elimination to $\bA$ to obtain $\bU$, and then examining the multipliers in the Gaussian elimination process to determine the entries below the main diagonal of $\bL$. We here provide another method for finding the LU decomposition without going through the process of  Gaussian elimination.

Again, we define $\bA_{i:j,m:n}$ as the $(j-i+1)\times(n-m+1)$ submatrix of $\bA$ with rows $i, i+1, \ldots, j$ and columns $m, m+1, \ldots, n$ of $\bA$. And simply, $a_{ij}$ as the ($i,j$)-th entry of $\bA$. 

Assume $\bA$ has the form $	\bA = \bL\bU$ of the LU decomposition. We will see, this assumption is equivalent to stating that $a_{11}$ is not zero, which implies the leading principal minors are nonzero.

From the properties of lower   and upper triangular matrices, we suppose that $\bA$ can be factored as 
$$
\bA=\left[
\begin{matrix}
	a_{11} & \bA_{1,2:n} \\
	\bA_{2:n,1} & \bA_{2:n,2:n}
\end{matrix}
\right] 
=\left[
\begin{matrix}
	1 & 0 \\
	\bL_{2:n,1} & \bL_{2:n,2:n}
\end{matrix}
\right] 
\left[
\begin{matrix}
	u_{11} & \bU_{1,2:n} \\
	0 & \bU_{2:n,2:n}
\end{matrix}
\right] = \bL\bU.
$$
Expanding  the product on the right-hand side of the above equation, we obtain 
$$
\bA=\left[
\begin{matrix}
	a_{11} & \bA_{1,2:n} \\
	\bA_{2:n,1} & \bA_{2:n,2:n}
\end{matrix}
\right] 
=\left[
\begin{matrix}
	u_{11} & \bU_{1,2:n} \\
	u_{11} \bL_{2:n,1} & \bL_{2:n,1}\bU_{1,2:n} +\bL_{2:n,2:n} \bU_{2:n,2:n}
\end{matrix}
\right] ,
$$
which allows us to determine the values of $\bL$ and $\bU$ by
\begin{equation}\label{equation:lu1_eq1}
\begin{aligned}
	u_{11} &= a_{11}   \\
	\bU_{1,2:n} &= \bA_{1,2:n}
\end{aligned}
\bigg\}
\qquad \mathrm{i.e.,} \qquad \bU_{1,1:n} = \bA_{1,1:n},
\end{equation}
\begin{equation}\label{equation:lu1_eq2}
\bL_{2:n,1} = \frac{1}{a_{11}} \bA_{2:n,1},
\end{equation}
and 
\begin{equation}
\bL_{2:n,2:n}\bU_{2:n,2:n} = \bA_{2:n,2:n} - \bL_{2:n,1}\bU_{1,2:n} = \bA_{2:n,2:n} -\frac{1}{a_{11}} \bA_{2:n,1}\bA_{1,2:n}.
\end{equation}
As $\bL_{2:n,2:n}\in \real^{(n-1)\times(n-1)}$ is also a unit lower triangular matrix and $\bU_{2:n,2:n}\in \real^{(n-1)\times(n-1)}$ is a nonsingular upper triangular matrix, both of which are of size $(n-1)\times (n-1)$ from the context. Let $\bA_2 = \bA_{2:n,2:n} -\frac{1}{a_{11}} \bA_{2:n,1}\bA_{1,2:n} $. So we can calculate $\bL_{2:n,2:n}$ and $\bU_{2:n,2:n}$ by factoring $\bA_2$ as 
$$
\bA_2 = \bL_{2:n,2:n}\bU_{2:n,2:n},
$$
which is an LU decomposition of a matrix with size $(n-1)\times(n-1)$. This suggests a recursive algorithm: to factorize a matrix of size $n\times n $, we calculate the first column of $\bL$ (with the first column being implicitly determined by Equation~\eqref{equation:lu1_eq2}) and the first row of $\bU$ (with the first row being implicitly determined by Equation~\eqref{equation:lu1_eq1}) leaving the other $n-1$ columns and $n-1$ rows of them to the next round. Continuing recursively, we arrive at a decomposition of a $1\times 1$ matrix.

\paragraph{A word on the leading principal minors.} In this process, we only assume the elements in entry $(1,1)$ of $\bA, \bA_2, \bA_3, \ldots, \bA_n$ are nonzero. This is actually the same assumption as the leading principal minors of $\bA$ are all nonzero. The recursive process is outlined in Algorithm~\ref{alg:lu-without-pivoting}.

\begin{algorithm}[H] 
\caption{LU Decomposition without Pivoting Recursively} 
\label{alg:lu-without-pivoting} 
\begin{algorithmic}[1] 
\Require 
Matrix $\bA$ is nonsingular and square with size $n\times n $; 
\State Calculate the first row of $\bU$: $\bU_{1,1:n} \leftarrow \bA_{1,1:n}$; \Comment{0 flops}
\State Calculate the first column of $\bL$: $l_{11}\leftarrow 1$ and $\bL_{2:n,1} \leftarrow \frac{1}{a_{11}} \bA_{2:n,1}$; \Comment{$n-1$ flops}
\State Calculate the LU decomposition 
$$\bA_2\leftarrow\bA_{2:n,2:n} -\frac{1}{a_{11}} \bA_{2:n,1}\bA_{1,2:n} = \bL_{2:n,2:n}\bU_{2:n,2:n};$$ \Comment{$2(n-1)^2$ flops}
\end{algorithmic} 
\end{algorithm}

\index{Floating-point operations}
\index{Flops}
\paragraph{Operation count.} The LU decomposition algorithm without pivoting is the first algorithm we have presented in this book. It is essential to evaluate its computational cost. 
To do so, we follow the standard route and count the number of \textit{floating-point operations (flops)} required by the algorithm \citep{golub2013matrix}. Each \textit{addition (+), subtraction ($-$), multiplication ($\times$), division ($\div$), square root ($\sqrt{\cdot}$), taking square ($a^2$), taking absolute value ($|a|$), and comparing two floating-point numbers ($a>b$)} operation counts as one flop. Note that we have the convention that an \textit{assignment operation} (e.g., set $a=b$) does not count as one flop.
It is also customary to count only the leading terms of the complexity. 
However, we should also note that the flop counts are often very crude measures of processing time and they may even be deceptive in some cases.
There are circumstances where algorithms with the same number of flops can execute with an extremely different amount of time due to the access patterns and memory hierarchy of modern computers \citep{elden2007matrix}.

\begin{theorem}[Algorithm Complexity: LU without Pivoting Recursively]\label{theorem:lu-complexity}
	Algorithm~\ref{alg:lu-without-pivoting} requires $\sim (2/3)n^3$ flops to compute the LU decomposition of an $n\times n$ matrix. Note that the theorem expresses only the leading term of the flop count. And the symbol ``$\sim$" has the usual asymptotic meaning
	\begin{equation*}
	\lim_{n \to +\infty} \frac{\mathrm{number\, of\, flops}}{(2/3)n^3} = 1.
	\end{equation*}
\end{theorem}

\begin{proof}[of Theorem~\ref{theorem:lu-complexity}]
The step 1 in Algorithm~\ref{alg:lu-without-pivoting} costs 0 flops and step 2 involves $(n-1)$ divisions. 

In step 3, we can compute $\frac{1}{a_{11}} \bA_{2:n,1}$ firstly, which costs $0$ flops as it has been calculated in step 2, and the outer product with $\bA_{1,2:n}$ costs $(n-1)^2$ multiplications/flops. So the computation of $\frac{1}{a_{11}} \bA_{2:n,1}\bA_{1,2:n}$ involves $(n-1)^2$ multiplications.
If we calculate $\bA_{2:n,1}\bA_{1,2:n}$ firstly, then the costs of $\frac{1}{a_{11}} \bA_{2:n,1}\bA_{1,2:n}$ is $2(n-1)^2$ totally. 
So we choose the first way to do the computation. Furthermore, The subtraction of the two matrices requires $(n-1)^2$ flops. As a result, the total cost of step 3 is $2(n-1)^2$ flops. 

It can be shown that the costs of the first recursive loop is $2(n-1)^2 + (n-1) = 2n^2-3n+1$ flops. Let $f(n)=2n^2-3n+1$, the overall cost can then be calculated as 
$$
\mathrm{cost}=f(n)+f(n-1)+\ldots+f(1).
$$
Simple calculation~\footnote{By the fact that $1^2+2^2+\ldots+n^2 = \frac{2n^3+3n^2+n}{6}$ and $1+2+\ldots+n=\frac{n(n+1)}{2}$.} can show that the complexity is $(2/3)n^3-(1/2)n^2-(1/6)n$ flops, or $(2/3)n^3$ flops if we keep only the leading term.
\end{proof}

\subsection{Complexity of Matrix and Vector Operations}\label{section:compl_mvops}
The calculation of the complexity extensively relies on the complexity of the multiplication of two matrices so that we formulate the finding in the following lemma.
\begin{lemma}[Vector Inner Product Complexity]
Let $\bv,\bw\in \real^{n}$ be two vectors. The  inner product of the two vectors $\bv^\top\bw$ is given by $\bv^\top\bw=v_1w_1+v_2w_2+\ldots v_nw_n$, which involves $n$ scalar multiplications and $n-1$ scalar additions. Therefore, the complexity for the inner product is $2n-1$ flops.
\end{lemma}

The matrix multiplication thus relies on the complexity of the inner product.
\begin{lemma}[Matrix Multiplication Complexity]\label{lemma:matrix-multi-complexity}
Let $\bA\in\real^{m\times n}$ and $\bB\in \real^{n\times k}$ be two matrices. Then, the complexity of the multiplication $\bC=\bA\bB$ is $mk(2n-1)$ flops.
\end{lemma}
\begin{proof}[of Lemma~\ref{lemma:matrix-multi-complexity}]
We notice that each entry of $\bC$ involves a vector inner product, which requires $n$ multiplications and $n-1$ additions. And there are $mk$ such entries, which leads to the conclusion.
\end{proof}

\section{Computing  LU without Pivoting Element-Wise: A=LU}\label{section:compute-lu-without-pivot-doolittle}
We observe that computing the LU decomposition is equivalent to solving the following equations
$$
a_{ij} = \sum_{s=1}^{\min(i,j)} l_{is} u_{sj}, \qquad \forall i,j\in \{1,2,\ldots, n\}.
$$
Furthermore, when $i\leq j$, the above equation can be decomposed into 
$$
a_{ij} = \sum_{s=1}^{i-1} l_{is} u_{sj} + u_{ij}, \qquad \text{since $l_{ii}=1$}.
$$
Suppose we know the first $k-1$ columns of $\bL$ and the first $k-1$ rows of $\bU$, we have the following observations:
$$
\begin{aligned}
a_{kj} &= \sum_{s=1}^{k-1} l_{ks} u_{sj} + u_{kj}, \qquad &\text{for all $j\in \{k,k+1,\ldots, n\}$,} \qquad \text{since $k\leq j$};\\
a_{ik} &= \sum_{s=1}^{k-1} l_{is} u_{sk}+l_{ik} u_{kk}, \qquad &\text{for all $i\in \{k+1,k+2,\ldots,n\}$},\qquad \text{since $i\geq k$}.
\end{aligned}
$$
Therefore, the $k$-th row of $\bU$ and $k$-th column of $\bL$ can be obtained by 
$$
\begin{aligned}
u_{kj} &= a_{kj} - \sum_{s=1}^{k-1} l_{ks} u_{sj} , \qquad &\text{for all $j\in \{k,k+1,\ldots, n\}$,} \qquad \text{since $k\leq j$};\\
l_{ik}  &=(a_{ik} - \sum_{s=1}^{k-1} l_{is} u_{sk})/u_{kk}, \qquad &\text{for all $i\in \{k+1,k+2,\ldots,n\}$},\qquad \text{since $i\geq k$}.
\end{aligned}
$$
This method is known as  \textit{Doolittle's method}, the values of $\bU$ and $\bL$ can be computed element-wise, and the process is formulated in Algorithm~\ref{alg:compute-lu-element-level}. We notice that, mathematically,  Doolittle's method is equivalent to the recursive algorithm  but from a different perspective.

\begin{algorithm}[h] 
\caption{LU Decomposition without Pivoting Element-Wise by Doolittle's Method} 
\label{alg:compute-lu-element-level} 
\begin{algorithmic}[1] 
\Require 
Matrix $\bA$ with size $n\times n$; 
\For{$k=1$ to $n$}
\State //i.e., the $k$-th row of $\bU$ and $k$-th column of $\bL$;
\For{$j=k$ to $n$}
\State $u_{kj} \leftarrow a_{kj} - \sum_{s=1}^{k-1} l_{ks} u_{sj}$;
\EndFor
\For{$i=k+1$ to $n$}
\State$l_{ik} \leftarrow (a_{ik} - \sum_{s=1}^{k-1} l_{is} u_{sk})/u_{kk}$;
\EndFor
\EndFor
\end{algorithmic} 
\end{algorithm}

	\begin{theorem}[Algorithm Complexity: LU without Pivoting Element-Wise]\label{theorem:lu-complexity-doolittle}
	Algorithm~\ref{alg:compute-lu-element-level} requires $\sim (2/3)n^3$ flops to compute the LU decomposition of an $n\times n$ matrix. 
	\end{theorem}

\begin{proof}[of Theorem~\ref{theorem:lu-complexity-doolittle}]
The step 4 in Algorithm~\ref{alg:compute-lu-element-level} requires $(k-1)$ multiplications, $(k-2)$ additions, and 1 subtractions for each loop $(k,j)$. And there are $n-k+1$ such loops for  $j$, resulting in a total of $(2k-2)(n-k+1)=-2k^2+(2n+4)k-2(n+1)$ flops  from step 4 for each loop $k$.

Similarly, the step 7 requires $(k-1)$ multiplications, $(k-2)$ additions, $1$ subtraction, and 1 division for each loop $(k,i)$. And there are $n-k$ such loops for  $i$, resulting a total of $(2k-1)(n-k)=-2k^2 + (2n+1)k-n$ flops  from step 7 for each loop $k$. 

Thus, for each loop $k$, the total complexity is $(2k-2)(n-k+1)+(2k-1)(n-k) = -4k^2 + (4n+5)k -(3n+2)$ flops. Let $f(k)=-4k^2 + (4n+5)k -(3n+2)$, the final complexity can be calculated by 
$$
\mathrm{cost=} f(1)+f(2) +\ldots +f(n).
$$
A simple calculation  shows that the complexity is $(2/3)n^3$ flops if we keep only the leading term.
\end{proof}

\index{Thin matrix}
\subsection{Extension to Thin Matrices}
We notice that the complexity of Algorithm~\ref{alg:compute-lu-element-level}  is the same as that of Algorithm~\ref{alg:lu-without-pivoting}. Doolittle's method is mathematically equivalent to the recursive algorithm. 
However,  Doolittle's method can be extended to compute the LU decomposition for $\bA\in\real^{m\times n}$ with $m\geq n$, which is commonly referred to as a  \textit{thin matrix}. 
The LU decomposition is given by $\bA=\bL\bU$, where $\bL\in\real^{m\times n}$ is a thin matrix and $\bU\in \real^{n\times n}$ is upper triangular. That is, $\bL$ is \textit{trapezoidal}: $l_{ij}=0$ for $i<j$. The procedure is similar and is outlined in Algorithm~\ref{alg:compute-lu-element-level-thin}, where the difference is illustrated in blue text. 

\begin{algorithm}[h] 
\caption{Thin LU Decomposition without Pivoting Element-Wise} 
\label{alg:compute-lu-element-level-thin} 
\begin{algorithmic}[1] 
\Require 
Matrix $\bA$ with size $\textcolor{mylightbluetext}{m}\times n$; 
\For{$k=1$ to $n$}
\State //i.e., the $k$-th row of $\bU$ and $k$-th column of $\bL$;
\For{$j=k$ to $n$}
\State $u_{kj} \leftarrow a_{kj} - \sum_{s=1}^{k-1} l_{ks} u_{sj}$;
\EndFor
\For{$i=k+1$ to $\textcolor{mylightbluetext}{m}$}
\State$l_{ik}  \leftarrow (a_{ik} - \sum_{s=1}^{k-1} l_{is} u_{sk})/u_{kk}$;
\EndFor
\EndFor
\end{algorithmic} 
\end{algorithm}

\begin{theorem}[Algorithm Complexity: LU Thin Matrix]\label{theorem:lu-complexity-doolittle-thin}
Algorithm~\ref{alg:compute-lu-element-level-thin} requires $\sim n^2(m-n/3)$ flops to compute the LU decomposition of an $m\times n$ matrix.  When $m=n$, the complexity of Algorithm~\ref{alg:compute-lu-element-level-thin} is $(2/3)n^3$ flops, which is the same as that of Algorithm~\ref{alg:compute-lu-element-level}.
\end{theorem}

\begin{proof}[of Theorem~\ref{theorem:lu-complexity-doolittle-thin}]
The complexity of step 4 is the same as that in Algorithm~\ref{alg:compute-lu-element-level}, which is $(2k-2)(n-k+1)=-2k^2+(2n+4)k-2(n+1)$ flops totally from step 4 for each loop $k$.
	
The complexity of step 7 is slightly different, where we replace $n$ by $m$, and it requires $(2k-1)(m-k)=-2k^2 + (2m+1)k-m$ flops totally from step 7 for each loop $k$. 

Thus, for each loop $k$, the total complexity is $(2k-2)(n-k+1)+(2k-1)(m-k) = -4k^2 + (2m+2n+5)k -(2n+m+2)$ flops. Let $f(k)=-4k^2 + (2m+2n+5)k -(2n+m+2)$, the final complexity can be calculated by 
$$
\mathrm{cost=} f(1)+f(2) +\ldots+ f(n).
$$
A simple calculation can show that the complexity is $n^2(m-n/3)$ flops if we keep only the leading term.
\end{proof}

\section{Computing  LU with Pivoting: A=PLU}\label{sec:compute-lu-pivoting}
Furthermore, we extend Algorithm~\ref{alg:lu-without-pivoting} to perform a  full LU decomposition with $\bA=\bP\bL\bU$. Note that we assume $a_{11}$ is nonzero in Algorithm~\ref{alg:lu-without-pivoting}. This is not necessarily true. We will avoid this assumption by introducing a permutation matrix. The following algorithm is derived from the proof of Theorem~\ref{theorem:lu-factorization-with-permutation} 
in Section~\ref{section:lu-perm}.

\begin{algorithm}[H] 
\caption{LU Decomposition with Pivoting} 
\label{alg:lu-with-pivoting} 
\begin{algorithmic}[1] 
\Require 
Matrix $\bA$ is nonsingular and square with size $n\times n$; 
\State Choose permutation matrix $\bP_1$ such that $\bB = \bP_1 \bA$ and $b_{11}\neq 0$; \Comment{0 flops}
\State Calculate the $\widehat{\bB}$ for next round: $\widehat{\bB}\leftarrow \bB_{2:n,2:n} -\frac{1}{b_{11}} \bB_{2:n,1} \bB_{1,2:n} \equiv \bP_2\bL_2 \bU_2$; \Comment{$2(n-1)^2+(n-1)$ flops}
\State Calculate the full LU decomposition of $\bA=\bP\bL\bU$ with 
$$
\bP \triangleq \bP_1^\top  
\left[
\begin{matrix}
1 & 0 \\
0 & \bP_2
\end{matrix}
\right], \qquad
\bL\triangleq
\left[
\begin{matrix}
1 & 0 \\
\frac{1}{b_{11}}\bP_2^\top \bB_{2:n,1} & \bL_2
\end{matrix}
\right], \qquad
\bU\triangleq
\left[
\begin{matrix}
b_{11} & \bB_{1,2:n} \\
\bzero & \bU_2
\end{matrix}
\right]
$$
\gap \gap \Comment{$n-1$ flops}
\end{algorithmic} 
\end{algorithm}

\begin{theorem}[Algorithm Complexity: LU with Pivoting]\label{theorem:lu-complexity-with-pivoting}
Algorithm~\ref{alg:lu-with-pivoting} requires $\sim (2/3)n^3$ flops to compute a full LU decomposition of an $n\times n$ nonsingular matrix.
\end{theorem}

\begin{proof}[of Theorem~\ref{theorem:lu-complexity-with-pivoting}]
Step 1 costs 0 flops as it only involves assignment operations (although finding $b_{11}\neq 0$ needs compare operations, we will ignore for brevity) and step 2 involves $(n-1)^2$ multiplications, $(n-1)$ divisions, and $(n-1)^2$ subtractions, which costs $2(n-1)^2+(n-1)$ flops to compute $\widehat{\bB}=\bB_{2:n,2:n} -\frac{1}{b_{11}} \bB_{2:n,1} \bB_{1,2:n}$ as shown in the proof of Theorem~\ref{theorem:lu-complexity}.

The computation of step 3 results from $\frac{1}{b_{11}}\bP_2^\top \bB_{2:n,1}$, which costs $n-1$ flops as the permutation operation does not count towards the flop count. 

So it costs $(2n^2-2n)$ flops in the first recursive loop. Let $f(n) = 2n^2-2n$, the final complexity can be calculated by 
$$
\mathrm{cost=} f(n)+f(n-1)+\ldots+f(1).
$$
Simple calculations can show that the complexity is $(2/3)n^3-(2/3)n$ flops, or $(2/3)n^3$ flops if we keep only the leading term.
\end{proof}

\index{Bandwidth}
\section{Bandwidth Preservation in  LU Decomposition without Permutation}\label{section:band_lu_wop}

For any matrix, its bandwidth can be defined as follows.
\begin{definition}[Matrix Bandwidth\index{Matrix bandwidth}]\label{defin:matrix-bandwidth}
For any matrix $\bA\in \real^{n\times n}$ with entry $(i,j)$ denoted by $a_{ij}$, the matrix $\bA$ has an \textit{upper bandwidth $q$} if $a_{ij} =0$ for $j>i+q$, and a \textit{lower bandwidth $p$} if $a_{ij}=0$ for $i>j+p$. 

An example of a $7\times 7$ matrix with an upper bandwidth $2$ and a lower bandwidth $3$ is shown below:
$$
\footnotesize
\begin{bmatrix}
\boxtimes & \boxtimes & \boxtimes & 0& 0 & 0 & 0\\
\boxtimes & \boxtimes & \boxtimes & \boxtimes& 0 & 0 & 0\\
\boxtimes & \boxtimes & \boxtimes & \boxtimes& \boxtimes & 0 & 0\\
\boxtimes & \boxtimes & \boxtimes & \boxtimes& \boxtimes & \boxtimes & 0\\
0 & \boxtimes & \boxtimes & \boxtimes& \boxtimes & \boxtimes & \boxtimes\\
0 & 0 & \boxtimes & \boxtimes& \boxtimes & \boxtimes & \boxtimes\\
0 & 0 & 0 & \boxtimes& \boxtimes & \boxtimes & \boxtimes\\
\end{bmatrix}.
$$
\end{definition}

It's important to note that the bandwidth is zero for diagonal matrices and increases as the nonzero entries spread further away from the diagonal.

In numerical analysis, the bandwidth of a matrix is crucial for understanding the computational cost of algorithms, particularly when solving systems of linear equations; while in  signal processing, the bandwidth of a matrix may relate to the frequencies represented in a discrete signal.
The concept of bandwidth extends to rectangular matrices, where it characterizes the ``spread" of nonzero elements along the rows or columns. For rectangular matrices, the bandwidth is defined independently for rows and columns: the \textit{column bandwidth} is defined as the maximum distance between the diagonal and the nonzero entries in any column; and the \textit{row bandwidth} is defined as the maximum distance between the diagonal and the nonzero entries in any row. Sometimes, the \textit{total bandwidth} is considered, which is the larger of the column bandwidth and the row bandwidth.

A matrix with a small bandwidth is considered ``narrow" because its nonzero entries are concentrated near the diagonal. Narrow matrices can often be processed more efficiently in algorithms, such as in solving systems of linear equations or in fast Fourier transform (FFT) algorithms, because of the reduced number of operations required.

\index{Bandwidth preservation}
Next, we demonstrate that the bandwidth of a matrix is preserved after performing the LU decomposition without permutation; that is, the factored components maintain the same bandwidth as the original matrix.
\begin{lemma}[Bandwidth Preservation]\label{lemma:lu-bandwidth-presev}
Let $\bA\in \real^{n\times n}$ be any square matrix with an upper bandwidth $q$ and a lower bandwidth $p$. If $\bA$ admits the LU decomposition $\bA=\bL\bU$, then the upper triangular factor $\bU$ maintains an upper bandwidth of $q$, and  lower triangular factor $\bL$ exhibits a lower bandwidth of $p$.
\end{lemma}
\begin{proof}[of Lemma~\ref{lemma:lu-bandwidth-presev}]
Building on  the computation of the LU decomposition without permutation (as described in Section~\ref{section:compute-lu-without-pivot}), and
leveraging  the properties of lower and upper triangular matrices, we can express the decomposition of $\bA$ as follows:
$$
\begin{aligned}
\bA&=\left[
\begin{matrix}
	a_{11} & \bA_{1,2:n} \\
	\bA_{2:n,1} & \bA_{2:n,2:n}
\end{matrix}
\right] 
=\left[
\begin{matrix}
	1 & \bzero \\
	\frac{1}{a_{11}} \bA_{2:n,1}  & \bI_{n-1}
\end{matrix}
\right] 
\left[
\begin{matrix}
	a_{11} & \bA_{1,2:n}\\
	\bzero & \bS
\end{matrix}
\right]
= \bL_1 \bU_1,
\end{aligned}
$$
where $\bS =\bA_{2:n,2:n} - \frac{1}{a_{11}}\bA_{2:n,1}\bA_{1,2:n}$ is the Schur complement of $a_{11}$ in $\bA$.
For simplicity, we can name this decomposition of $\bA$ as the $s$-decomposition of $\bA$.
The first column of $\bL_1$ and the first row of $\bU_1$ exhibit the desired structure, with bandwidth $p$ and $q$, respectively. Additionally,  the Schur complement $\bS$ of $a_{11}$ has an upper bandwidth $q-1$ and a lower bandwidth $p-1$, respectively. The overall result is derived through induction on the $s$-decomposition of $\bS$. 
\end{proof}

\index{Tridiagonal matrix}
\subsection*{LU Decomposition of Tridiagonal Matrices}
We consider a tridiagonal matrix of the following form \footnote{We will extensively discuss the tridiagonal decomposition for symmetric matrices in Section~\ref{section:tridiagonal_decom}, in which case a symmetric matrix can be factored into a symmetric tridiagonal matrix. }
$$
\bA = 
\footnotesize
\begin{bmatrix}
b_1 & c_1 & & & & &\\
a_2 & b_2 & c_2 & & & & \\
&     a_3 & b_3 & c_3 & & & \\
& &         \ddots & \ddots  & \ddots & & \\
& & &                a_{n-2}& b_{n-2} & c_{n-2}  & \\
& & & &                       a_{n-1} & b_{n-1} & c_{n-1}  \\
& & & & &                               a_n     & b_n 
\end{bmatrix}
\normalsize
\in \real^{n\times n}.
$$
That is, the tridiagonal matrix $\bA$ has (upper and lower) bandwidth 1. Then the leading principal minors (Definition~\ref{definition:leading-principle-minors}) of $\bA$ are 
$$
\Delta_0=1, \gap 
\Delta_1 = b_1, \gap 
\Delta_k = b_k\Delta_{k-1} - a_kc_{k-1}\Delta_{k-2}, \gap k\in\{2,3,\ldots,n\},
$$
where $\Delta_k = \det(\bA[1:k,1:k])$ with the convention $\Delta_0=1$.
\begin{exercise}
Prove the above result. \textit{Hint: Use induction}. 
\end{exercise}
Then the LU decomposition of the matrix $\bA$ is 
$$
\stackrel{\bL\bU=}{\bA} = 
\footnotesize
\begin{bmatrix}
1&&&&& \\
a_2 \frac{\Delta_0}{\Delta_1}&1&&&& \\
&a_3 \frac{\Delta_1}{\Delta_2}&1&&& \\
&&\ddots &\ddots&& \\
&&&a_{n-1} \frac{\Delta_{n-3}}{\Delta_{n-2}}&1& \\
&&&&a_{n} \frac{\Delta_{n-2}}{\Delta_{n-1}}&1 \\
\end{bmatrix}
\begin{bmatrix}
\frac{\Delta_1}{\Delta_0}& c_1&&&& \\
&\frac{\Delta_2}{\Delta_1}& c_2&&& \\
&&\frac{\Delta_3}{\Delta_2}& c_3&& \\
&& &\ddots&\ddots& \\
&&& &\frac{\Delta_{n-1}}{\Delta_{n-2}}& c_{n-1} \\
&&&& &\frac{\Delta_n}{\Delta_{n-1}} \\
\end{bmatrix}.
$$
Both the lower  and upper triangular matrices preserve the bandwidth value of 1.

\index{Linear system}
\subsection*{Solving Linear System of Tridiagonal Matrices}
Consider the linear system $\bA\bx=\bz$ with tridiagonal $\bA\in\real^{n\times n}$ and $\Delta_k\neq 0$ for $k\in\{1,2,\ldots, n\}$. The linear system can be efficiently solved using the LU decomposition. Let the diagonal matrix 
$\bD = \diag(\frac{\Delta_1}{\Delta_0}, \frac{\Delta_2}{\Delta_1}, \ldots, \frac{\Delta_{n}}{\Delta_{n-1}})$ and 
$$
d_1=\frac{c_1}{b_1},\gap
d_k = c_k\frac{\Delta_{k-1}}{\Delta_k},\,\, k\in\{2,3,\ldots,n-1\}, \gap 
d_n = \frac{\Delta_n}{\Delta_{n-1}}=b_n-a_n d_{n-1}.
$$
Then $\bA$ can be factored as $\bA=(\bL\bD)(\bD^{-1}\bU)$:
$$
\bA = 
\footnotesize
\begin{bmatrix}
\frac{c_1}{d_1}&&&&& \\
a_2 & \frac{c_2}{d_2}&&&& \\
&a_3 &\frac{c_3}{d_3}&&& \\
&&\ddots &\ddots&& \\
&&&a_{n-1} &\frac{c_{n-1}}{d_{n-1}}& \\
&&&&a_{n} & d_n\\
\end{bmatrix}
\begin{bmatrix}
1& d_1&&&& \\
&1& d_2&&& \\
&&1& d_3&& \\
&& &\ddots&\ddots& \\
&&& &1& d_{n-1} \\
&&&& &1\\
\end{bmatrix}.
$$
Then the linear system $\bA\bx=\bz$ can be solved firstly by $(\bL\bD)\bw = \bz$, and then by $(\bD^{-1}\bU)\bx = \bw$. 
The first system is obtained by
$$
w_1 = z_1 \frac{d_1}{c_1},\gap 
w_k = (z_k - a_kw_{k-1}) \frac{d_k}{c_k}, \,\, k\in\{2,3,\ldots,n-1\},\gap
w_n = (z_n-a_n w_{n-1}) \frac{1}{d_n}. 
$$
The second system is obtained by 
$$
x_n=w_n, \gap 
x_k = (w_{k-1} - x_{k-1})/d_{k-1}, \,\,k\in\{1,2,\ldots, n-1\}.
$$

\section{Block LU Decomposition}
Another form of the LU decomposition involves factoring the matrix into  block triangular matrices: the original matrix is factored into the product of a block lower triangular matrix and a block upper triangular matrix.

\begin{theoremHigh}[Block LU Decomposition without Permutation]\label{theorem:block-lu-factorization-without-permutation}
Let $\bA$ be  any $n\times n$ square matrix. If the first $m$ \textit{leading principal block submatrices} are nonsingular, then $\bA$ can be factored as 
\begin{equation}
\bA = \bL\bU
=
\begin{bmatrix}
\bI  & & & \\
\bL_{21} & \bI & & \\
\vdots & & \ddots & \\
\bL_{m1} & \ldots & \bL_{m,m-1} & \bI 
\end{bmatrix}
\begin{bmatrix}
\bU_{11}  &\bU_{12} & \ldots & \bU_{1m}\\
  & \bU_{22} & &  \vdots \\
  & & \ddots & \bU_{m-1,m}\\
  & & & \bU_{mm}\\
\end{bmatrix}
, \nonumber
\end{equation}
where $\bL_{i,j}$'s and $\bU_{ij}$'s are some block matrices. 
Specifically, this decomposition is unique.
\end{theoremHigh}
Note that the matrix $\bU$ in the above theorem is not necessarily upper triangular. An illustrative  example can be shown as follows:
$$
\bA = 
\left[\begin{array}{cc;{2pt/2pt}cc}
0& 1 &   1 & 1\\
-1& 2 &   -1 & 2\\\hdashline[2pt/2pt]
2& 1 &   4 & 2\\	
1& 2 &   3 & 3\\
\end{array}\right]
=
\left[\begin{array}{cc;{2pt/2pt}cc}
	1& 0 &   0& 0\\
	0& 1 &   0 & 0\\\hdashline[2pt/2pt]
	5& -2 &   1 & 0\\	
	4& -1 &   0 & 1\\
\end{array}\right]
\left[\begin{array}{cc;{2pt/2pt}cc}
	0& 1 &   1 & 1\\
	-1& 2 &   -1 & 2\\\hdashline[2pt/2pt]
	0& 0 &   -3 & 1\\	
	0& 0 &   -2 & 1\\
\end{array}\right].
$$
The standard non-block LU decomposition fails on $\bA$ since the entry $(1,1)$ is zero. However, the block LU decomposition exists.

\index{Linear system}
\section{Application: Linear System via  LU Decomposition}\label{section:lu-linear-sistem}
Consider a well-determined linear system $\bA\bx = \bb$, where $\bA$  is an  $n\times n $ nonsingular matrix. 
Instead of solving the system by computing the inverse of $\bA$, we solve linear equations by the LU decomposition. Suppose $\bA$ admits the LU decomposition $\bA = \bP\bL\bU$; the solution of the linear system can be obtained using the following algorithm.

\begin{algorithm}[H] 
\caption{Solving Linear Equations Using LU Decomposition} 
\label{alg:linear-equation-by-LU} 
\begin{algorithmic}[1] 
\Require 
Matrix $\bA$ is nonsingular and square with size $n\times n $, solve $\bA\bx=\bb$; 
\State LU Decomposition: factor $\bA$ as $\bA=\bP\bL\bU$; \Comment{(2/3)$n^3$ flops}
\State Permutation: $\bw = \bP^\top\bb$; \Comment{0 flops }
\State Forward substitution: solve $\bL\bv = \bw$; \Comment{$1+3+\ldots + (2n-1)=n^2$ flops}
\State Backward substitution: solve $\bU\bx= \bv$; \Comment{$1+3+\ldots + (2n-1)=n^2$ flops}
\end{algorithmic} 
\end{algorithm}

The complexity of the decomposition step is $(2/3)n^3$ flops (Theorem~\ref{theorem:lu-complexity-with-pivoting}). Both the backward and forward substitution steps each incur a cost of $1+3+\ldots + (2n-1)=n^2$ flops. Therefore, the total cost for computing the linear system using  LU factorization is $(2/3)n^3 + 2n^2$ flops. If we keep only the leading term,  Algorithm~\ref{alg:linear-equation-by-LU} costs approximately $(2/3)n^3$ flops, where the most significant complexity comes from the LU decomposition.

\paragraph{Linear system via the block LU decomposition.} In the context of a block LU decomposition of $\bA=\bL\bU$, we need to solve $\bL\bv = \bw$ and $\bU\bx = \bv$. However, the latter system is not triangular and requires some extra computations.

\section{Application: Computing  Inverses of Nonsingular Matrices}\label{section:inverse-by-lu}
According to Theorem~\ref{theorem:lu-factorization-with-permutation}, given any nonsingular matrix $\bA\in \real^{n\times n}$, a full LU factorization exists, represented as $\bA=\bP\bL\bU$. Then the inverse can be obtained by solving the matrix equation
$$
\bA\bX = \bI,
$$
which involves a computation of  $n$ linear systems: $\bA\bx_i = \be_i$ for all $i \in \{1, 2, \ldots, n\}$, where $\bx_i$ is the $i$-the column of $\bX$ and $\be_i$ is the $i$-th column of $\bI$ (i.e., the $i$-th unit basis vector). 

\index{Linear system}
\begin{theorem}[Inverse of Nonsingular Matrix by Linear System]\label{theorem:inverse_vanilla}
Computing the inverse of a nonsingular matrix $\bA \in \real^{n\times n}$ by solving $n$ linear systems requires $\sim (2/3)n^3 + n(2n^2)=(8/3)n^3$ flops, where the term $(2/3)n^3$ comes from the computation of the LU decomposition of $\bA$, as detailed in  Theorem~\ref{theorem:lu-complexity-with-pivoting}.
\end{theorem}
The proof is straightforward  using Algorithm~\ref{alg:linear-equation-by-LU}.

\index{Row partitions}
\index{Column partitions}

However, the complexity can be reduced by taking the advantage of the structures of $\bU$ and $\bL$.
We find that the inverse of the nonsingular matrix is $\bA^{-1} = \bU^{-1}\bL^{-1}\bP^{-1}=\bU^{-1}\bL^{-1}\bP^\top$.
\begin{theorem}[Inverse of Nonsingular Matrix by LU Factorization]\label{theorem:inverse-by-lu2}
Computing the inverse of a nonsingular matrix $\bA \in \real^{n\times n}$ by $\bA^{-1} = \bU^{-1}\bL^{-1}\bP^\top$ needs $\sim (2/3)n^3 + (4/3)n^3=2n^3$ flops, where $(2/3)n^3$ comes from the computation of the LU decomposition of $\bA$ by Theorem~\ref{theorem:lu-complexity-with-pivoting}.
\end{theorem}
\begin{proof}[of Theorem~\ref{theorem:inverse-by-lu2}]
We note the computation of $\bU^{-1}\bL^{-1}\bP^\top$ comes from $\bY=\bU^{-1}\bL^{-1}$. Here, $\bU^{-1}$ is an upper triangular, $\bL^{-1}$ is a unit lower triangular matrix.
Suppose 
$$\bU^{-1} = \bZ=\begin{bmatrix}
	-\bz_1^\top-\\
	-\bz_2^\top-\\
	\vdots \\
	-\bz_n^\top- 
\end{bmatrix}\qquad
\text{ and }\qquad
\bU = [\bu_1, \bu_2, \ldots, \bu_n], 
$$
are the row partitions and column partitions of $\bU^{-1}$ and $\bU$, respectively.
Since both $\bZ$ and $\bU$ are upper triangular matrices. Thus we have
$$
\bI=
\footnotesize
\begin{bmatrix}
	\bz_1^\top\bu_1=1 & \bz_1^\top\bu_2=0 & \bz_1^\top\bu_3=0 & \bz_1^\top\bu_4=0 & \ldots  & \bz_1^\top\bu_n=0 \\
	\bz_2^\top\bu_1=0 & \bz_2^\top\bu_2=1 & \bz_2^\top\bu_3=0            & \bz_2^\top\bu_4=0 & \ldots  & \bz_2^\top\bu_n=0 \\
	\bz_3^\top\bu_1=0 & \bz_3^\top\bu_2=0 & \bz_3^\top\bu_3=1 &\bz_3^\top\bu_4=0  & \ldots & \bz_3^\top\bu_n=0 \\
	\bz_4^\top\bu_1=0 & \bz_2^\top\bu_2=0 & \bz_4^\top\bu_3=0 &\bz_4^\top\bu_4=1  & \ldots & \bz_4^\top\bu_n=0 \\
	\vdots  &\vdots  & \vdots &\vdots  & \ldots & \vdots \\
	\bz_n^\top\bu_1=0 & \bz_n^\top\bu_2=0 & \bz_n^\top\bu_3=0&\bz_n^\top\bu_4=0  & \ldots & \bz_n^\top\bu_n=1
\end{bmatrix}.
$$
By $\bz_1^\top\bu_1=1$, we can compute the first component of $\bz_1$ with 1 flop; By $\bz_1^\top\bu_2=0$, we can compute the second component of $\bz_1$ with 3 flops, as the first component is already calculated; $\ldots$. 
Then we list the complexity of each inner product in an $n\times n$ matrix, where each entry ($i,j$) denotes the cost to calculate the $(i,j)$-th element of $\bU^{-1}$:
$$
\mathrm{cost}=
\footnotesize
\begin{bmatrix}
\bz_1^\top\bu_1=1 & \bz_1^\top\bu_2=3           & \bz_1^\top\bu_3=5 & \bz_1^\top\bu_4=7 & \ldots  & \bz_1^\top\bu_n=2n-1 \\
0 & \bz_2^\top\bu_2=1 & \bz_2^\top\bu_3=3            & \bz_2^\top\bu_4=5 & \ldots  & \bz_2^\top\bu_n=2n-3 \\
0 & 0 & \bz_3^\top\bu_3=1 &\bz_3^\top\bu_4=3  & \ldots & \bz_3^\top\bu_n=2n-5 \\
0 & 0 & 0 &\bz_4^\top\bu_4=1  & \ldots & \bz_4^\top\bu_n=2n-7 \\
\vdots  &\vdots  & \vdots &\vdots  & \ldots & \vdots \\
0 & 0 & 0& 0 & \ldots & \bz_n^\top\bu_n=1
\end{bmatrix}=
\begin{bmatrix}
	n^2 \\
	(n-1)^2 \\
	(n-2)^2 \\
	(n-3)^2 \\
	\vdots \\
	1 
\end{bmatrix},
$$ 
which is a sum of $n$ sets of \textit{arithmetic sequences} and the sum of each sequence is displayed  in the last equality above. Thus the total cost to compute $\bU^{-1}$ is \underline{$(2n^3+3n^2+n)/6$} flops. 
 Similarly, the computation of  $\bL^{-1}$ also requires \underline{$(2n^3+3n^2+n)/6$} flops.

A moment of reflexion on the multiplication of $\bY=\bU^{-1}\bL^{-1}$ would reveal that:
\begin{itemize}
\item The entry $(1,1)$ of $y_{11}$ involves the computation of an inner product with $n$-dimension, which takes $2n-1$ flops (i.e., $n$ multiplications and $n-1$ additions).

\item The entry $(1,2)$ of $y_{12}$ involves the computation of an inner product with $(n-1)$-dimension, which takes $2n-3$ flops. 

\item The process can go on, we write the flops in an $n\times n$ matrix with each entry $(i,j)$-th meaning the number of flops to calculate the $(i,j)$-th element of $\bY$:
$$
\text{costs of $\bY=\bU^{-1}\bL^{-1}$}=
\footnotesize
\begin{bmatrix}
	\{ \textcolor{red}{2n-1}             & \textcolor{red}{2n-3} & \textcolor{red}{2n-5} & \textcolor{red}{2n-7} & \textcolor{red}{\ldots} & \textcolor{red}{1}\} \\
	\overbrace{\textcolor{cyan}{2n-3}}	 & \{ \textcolor{green}{2n-3} & \textcolor{green}{2n-5} & \textcolor{green}{2n-7} & \textcolor{green}{\ldots} & \textcolor{green}{1}\} \\
	\textcolor{cyan}{2n-5}        & \overbrace{\textcolor{magenta}{2n-5}} & \{ \textcolor{mylightbluetext}{2n-5} & \textcolor{mylightbluetext}{2n-7} & \textcolor{mylightbluetext}{\ldots} & \textcolor{mylightbluetext}{1}\} \\
	\textcolor{cyan}{2n-7}        & \textcolor{magenta}{2n-7}  &  \overbrace{\textcolor{olive}{2n-7}} & \{\textcolor{orange}{2n-7} &\textcolor{orange}{\ldots}  & \textcolor{orange}{1}\} \\
	\textcolor{cyan}{\vdots}      & \textcolor{magenta}{\vdots} & \textcolor{olive}{\vdots} & \textcolor{purple}{\vdots} & \ldots & \vdots  \\
	\underbrace{\textcolor{cyan}{1}}  & \underbrace{\textcolor{magenta}{1}}          & \underbrace{\textcolor{olive}{1}}      & \underbrace{\textcolor{purple}{1}}  & \vdots & \textcolor{teal}{1} \\
\end{bmatrix}.
$$
The matrix is symmetric, and the overall complexity can be determined by summing several arithmetic sequences, with each sequence represented by a different color. To enhance clarity, each arithmetic sequence is highlighted using braces in the matrix above. Straightforward calculations reveal that the complexity amounts to \underline{$(2/3)n^3 + (1/3)n$}.
\end{itemize}

As a result, the total cost is then $(2/3)n^3 + (1/3)n^3 + (1/3)n^3 + (2/3)n^3=2n^3$ flops if we keep only the leading term, where the first $(2/3)n^3$ comes from the computation of the LU decomposition of $\bA$ (rather than $8/3n^3$ flops in Theorem~\ref{theorem:inverse_vanilla}).
\end{proof}
\section{Application: Computing the Determinant}\label{sectin:app_lu_det}
We can easily find the determinant of a matrix   using the LU decomposition.
If $\bA=\bL\bU$, then $\det(\bA) = \det(\bL\bU) = \det(\bL)\det(\bU) = u_{11}u_{22}\ldots u_{nn}$, where $u_{ii}$ denotes the $i$-th diagonal of $\bU$ for $i\in \{1,2,\ldots,n\}$. \footnote{The determinant of a lower triangular matrix (or an upper triangular matrix) is the product of its diagonal entries.}

Moreover, for the LU decomposition with permutation $\bA=\bP\bL\bU$,  we have $\det(\bA) = \det(\bP\bL\bU) = \det(\bP)u_{11}u_{22}\ldots u_{nn}$.
The determinant of a permutation matrix is either 1 or –1, as 
changing rows around (which alters the sign of the determinant \footnote{The determinant changes sign when two rows are exchanged (sign reversal).}) transforms a permutation matrix
into the identity matrix $\bI$, whose determinant is one.

\section{Pivoting}
\subsection{Partial Pivoting}\label{section:partial-pivot-lu}
In practical applications, the practice of pivoting is considered beneficial even when technically unnecessary. 
When dealing with a linear system via the LU decomposition as shown in Algorithm~\ref{alg:linear-equation-by-LU}, if the diagonal entries of $\bU$ are small, it can lead to inaccurate solutions for the linear solution. 
Thus, it is common practice to select the largest entry as the pivot. 
This approach, known as the \textit{partial pivoting}, is commonly employed to enhance numerical stability. 
For example, consider the partial pivoting for a $4\times 4$ matrix:
\begin{tcolorbox}[title={Partial Pivoting For a $4\times 4$ Matrix},colback=\mdframecolorTheorem]
\footnotesize
\begin{equation}\label{equation:elmination-steps2}
\begin{sbmatrix}{\bA}
\boxtimes & \boxtimes & \boxtimes & \boxtimes \\
\boxtimes & \boxtimes & \boxtimes & \boxtimes \\
\boxtimes & \boxtimes & \boxtimes & \boxtimes \\
\boxtimes & \boxtimes & \boxtimes & \boxtimes
\end{sbmatrix}
\stackrel{\bE_1}{\longrightarrow}
\begin{sbmatrix}{\bE_1\bA}
\boxtimes & \boxtimes & \boxtimes & \boxtimes \\
0 & \bm{2} & \bm{\boxtimes} & \bm{\boxtimes} \\
0 & \bm{5} & \bm{\boxtimes} & \bm{\boxtimes} \\
0 & \bm{7} & \bm{\boxtimes} & \bm{\boxtimes}
\end{sbmatrix}
\stackrel{\bP_1}{\longrightarrow}
\begin{sbmatrix}{\bP_1\bE_1\bA}
\boxtimes & \boxtimes & \boxtimes & \boxtimes \\
0 & \textcolor{mylightbluetext}{7} & \textcolor{mylightbluetext}{\boxtimes} & \textcolor{mylightbluetext}{\boxtimes} \\
0  & 5 & \boxtimes & \boxtimes \\
0 & \textcolor{mylightbluetext}{2} & \textcolor{mylightbluetext}{\boxtimes} & \textcolor{mylightbluetext}{\boxtimes}
\end{sbmatrix}
\stackrel{\bE_2}{\longrightarrow}
\begin{sbmatrix}{\bE_2\bP_1\bE_1\bA}
\boxtimes & \boxtimes & \boxtimes & \boxtimes \\
0 &  7 & \boxtimes & \boxtimes \\
0 & \bm{0}  & \textcolor{mylightbluetext}{\bm{\boxtimes}} & \bm{\boxtimes} \\
0 & \bm{0} & \bm{0} & \textcolor{mylightbluetext}{\bm{\boxtimes}}
\end{sbmatrix},
\end{equation}
\end{tcolorbox}
\noindent in which case, we pick 7 as the pivot after the transformation by $\bE_1$, even if such a choice may seem unnecessary. This interchange permutation ensures that no multiplier is greater than 1 in absolute value during the Gaussian elimination process. 
A concrete example is presented below.
\begin{example}[Partial Pivoting]
Suppose 
$$
\bA= 
\begin{bmatrix}
2 & 10 & 5\\
1 & 4  & -2 \\
6 & 8  & 4
\end{bmatrix}.
$$
To get the smallest possible multipliers in the first Gaussian transformation, we need to interchange the largest value in the first column to entry $(1,1)$. The permutation matrix $\bP_1$ is doing so
$$
\bP_1 =
\begin{bmatrix}
&& 1\\
&1& \\
1& &
\end{bmatrix}
, \qquad \text{such that} \qquad 
\bU=\bP_1\bA = 
\begin{bmatrix}
6 & 8  & 4\\
1 & 4  & -2 \\
2 & 10 & 5
\end{bmatrix}.
$$
Then it follows that $\bE_1$ will introduce zero below the entry $(1,1)$, 
$$
\bE_1 = 
\begin{bmatrix}
1 & &\\
-1/6 & 1 \\
-1/3 & 0 & 1
\end{bmatrix}, 
\qquad \text{and} \qquad 
\bU=\bE_1\bP_1\bA=
\begin{bmatrix}
	6 & 8  & 4\\
    0  & 8/3  & -8/3 \\
    0  &  22/3  & 11/3
\end{bmatrix}.
$$
Now, we pivot 22/3 before 8/3, and the permutation is given by 
$$
\bP_2 =
\begin{bmatrix}
1	&& \\
	& &1 \\
	&1 &
\end{bmatrix}
, \qquad \text{such that} \qquad 
\bU=\bP_2\bE_1\bP_1\bA = 
\begin{bmatrix}
	6 & 8  & 4\\
	0  &  22/3  & 11/3\\
		0  & 8/3  & -8/3 \\
\end{bmatrix}.
$$
Finally, the Gaussian transformation to introduce zero below entry ($2,2$) is given by 
$$
\bE_2 =
\begin{bmatrix}
1 & & \\
& 1 & \\
& -4/11 & 1
\end{bmatrix}
\qquad \text{and} \qquad 
\bU=\bE_2 \bP_2\bE_1\bP_1\bA=
\begin{bmatrix}
	6 & 8  & 4\\
	0  &  22/3  & 11/3\\
	0  & 0  & -4\\
\end{bmatrix}.
$$
This results in the final  $\bU$.
\end{example}
As discussed in \S~\ref{section:gaussian-elimination}, the Gaussian transformation $\bE_k$ in the $k$-th step of the partial pivoting is given by 
$
\bE_k = \bI - \bz_k \be_k^\top,
$
where $\be_k \in \real^n$ is the $k$-th unit basis vector, and $\bz_k\in \real^n$ is given by 
$$
\bz_k = [0, \ldots, 0, z_{k+1}, \ldots, z_n]^\top, \qquad z_i= \frac{u_{ik}}{u_{kk}},\,\, \forall i \in \{k+1,\ldots, n\}.
$$
We realize that $\bE_k$ is a unit lower triangular matrix with only the $k$-th column of the lower submatrix being nonzero:
$$
\begin{blockarray}{ccccccc}
	\begin{block}{c[ccccc]c}
		&	\bI & \bzero & \bzero & \bzero & \bzero &    \\
		& \bzero & \ddots & \bzero & \bzero	 & \bzero  &   \\
		\bE_k=	& \bzero&\bzero & 1 &\bzero	 &  \bzero  &  \\
		& \bzero & \bzero & \boxtimes & \bI	 & \bzero  &  \\
		& \bzero & \bzero & \boxtimes & \bzero  & \bI &  \\
	\end{block}
	&  & &  k&  &  &  \\
\end{blockarray}.
$$
More generally, the procedure for computing the LU decomposition with partial pivoting of $\bA\in\real^{n\times n}$ is outlined in Algorithm~\ref{alg:lu-partial-pivot}.

\begin{algorithm}[h] 
\caption{LU Decomposition with Partial Pivoting} 
\label{alg:lu-partial-pivot} 
\begin{algorithmic}[1] 
\Require 
Matrix $\bA$ with size $n\times n$;  
\State Initialize $\bU \leftarrow \bA$;
\For{$k=1$ to $n-1$} \Comment{i.e., get the $k$-th column of $\bU$}
\State \parbox[t]{\dimexpr\linewidth-\algorithmicindent}{
Find a row permutation matrix $\bP_k$ that swaps $u_{kk}$ with the largest element in $|\bU[k:n,k]|$;
}
\State $\bU \leftarrow \bP_k\bU$;
\State \parbox[t]{\dimexpr\linewidth-\algorithmicindent}{
Determine the Gaussian transformation $\bE_k $ to introduce zeros below the diagonal in the $k$-th column of $\bU$;}
\State $\bU \leftarrow \bE_k\bU$;
\EndFor
\State Output $\bU$;
\end{algorithmic} 
\end{algorithm}
The algorithm requires $\sim2/3(n^3)$ flops and $(n-1)+(n-2)+\ldots + 1 \sim \mathcalO(n^2)$ comparisons, which result from the pivoting procedure.
Upon completion, the upper triangular matrix $\bU$ is given by 
\begin{equation}\label{equation:partial_pivot_u1}
\bU = \bE_{n-1}\bP_{n-1} \ldots \bE_2\bP_2\bE_1\bP_1\bA.
\end{equation}

\paragraph{Computing the final $\bL$.} Here, we demonstrate that Algorithm~\ref{alg:lu-partial-pivot} computes the LU decomposition in the following form
$$
\bA = \bP\bL\bU,
$$
where $\bP\triangleq\bP_1 \bP_2\ldots\bP_{n-1}$ takes account of all the interchanges, $\bU$ is the upper triangular matrix obtained directly from the algorithm, and $\bL$ is unit lower triangular with $|l_{ij}|\leq 1$ for all $1 \leq i,j\leq n$. 
We note that the submatrix $\bL[k+1:n,k]$ is a permuted version of $\bE_k$'s multipliers. To see this, we notice that the permutation matrices used in the algorithm fall into a special kind of permutation matrix since we only interchange two rows of the matrix. \textit{This implies the $\bP_k$'s are symmetric and $\bP_k^2 = \bI$ for $k\in \{1,2,\ldots,n-1\}$}.
Let
$$
\bM_k \triangleq (\bP_{n-1} \ldots \bP_{k+1}) \bE_k (\bP_{k+1} \ldots \bP_{n-1}).
$$
Then, substituting into Equation~\eqref{equation:partial_pivot_u1}, $\bU$ can be written as 
$$
\bU = \bM_{n-1}\ldots \bM_2\bM_1 \bP^\top \bA.
$$
To analyze $\bM_k$, we realize that $\bP_{k+1}$ is a permutation with the upper-left $k\times k$ block being an identity matrix and $\bE_k=\bI_n-\bz_k\be_k^\top$ by \eqref{equation:elimination_mat}. Thus:
$$
\begin{aligned}
\bM_k &= (\bP_{n-1} \ldots \bP_{k+1}) (\bI_n-\bz_k\be_k^\top) (\bP_{k+1} \ldots \bP_{n-1})\\
&=\bI_n - (\bP_{n-1} \ldots \bP_{k+1})(\bz_k\be_k^\top)(\bP_{k+1} \ldots \bP_{n-1}) \\
&=\bI_n - (\bP_{n-1} \ldots \bP_{k+1}\bz_k) (\be_k^\top\bP_{k+1} \ldots \bP_{n-1}) \\
&=\bI_n - (\bP_{n-1} \ldots \bP_{k+1}\bz_k)\be_k^\top.  \qquad &\text{(since $\be_k^\top\bP_{k+1} \ldots \bP_{n-1} = \be_k^\top$)}
\end{aligned}
$$
This implies that $\bM_k$ is unit lower triangular with its $k$-th column representing a permuted version of $\bE_k$. Consequently, the ultimate lower triangular $\bL$ is  given by 
$$
\bL \triangleq \bM_1^{-1}\bM_2^{-1} \ldots \bM_{n-1}^{-1}.
$$
We obtain the full LU decomposition $\bA =\bP\bL\bU$.

\subsection{Complete Pivoting}\label{section:complete-pivoting}
In partial pivoting, when introducing zeros below the diagonal of the $k$-th column of $\bU$, the $k$-th pivot is determined by scanning the current subcolumn $\bU[k:n,k]$. In \textit{complete pivoting}, the largest absolute entry in the current submatrix $\bU[k:n,k:n]$ is interchanged into the entry $(k,k)$ of $\bU$. Therefore, an additional \textit{column permutation} $\bQ_k$ is applied at each step. The final upper triangular matrix $\bU$ is obtained by 
$$
\bU = \bE_{n-1}\bP_{n-1}\ldots (\bE_2\bP_2(\bE_1\bP_1\bA\bQ_1)\bQ_2) \ldots \bQ_{n-1}.
$$
Similarly, the complete pivoting algorithm is formulated in Algorithm~\ref{alg:lu-complete-pivot}. 
\begin{algorithm}[H] 
\caption{LU Decomposition with Complete Pivoting} 
\label{alg:lu-complete-pivot} 
\begin{algorithmic}[1] 
\Require 
Matrix $\bA$ with size $n\times n$;  
\State Initialize $\bU \leftarrow \bA$;
\For{$k=1$ to $n-1$} \Comment{the value $k$ is to get the $k$-th column of $\bU$}
\State \algoalign{
Find a row permutation matrix $\bP_k$ and a column permutation $\bQ_k$ that swaps $u_{kk}$ with the largest element in $|\bU[k:n,k:n]|$, say $u_{ab} = \max{|\bU[k:n,k:n]|}$;}
\State $\bU \leftarrow\bP_k\bU\bQ_k$;
\State \algoalign{Determine the Gaussian transformation $\bE_k $ to introduce zeros below the diagonal in the $k$-th column of $\bU$;}
\State $\bU \leftarrow \bE_k\bU$;
\EndFor
\State Output $\bU$;
\end{algorithmic} 
\end{algorithm}
The algorithm requires $~2/3(n^3)$ flops and $(n^2+(n-1)^2+\ldots +1^2)\sim \mathcalO(n^3)$ comparisons resulting from the pivoting procedure. Again, let $\bP\triangleq\bP_1 \bP_2\ldots\bP_{n-1}$, $\bQ \triangleq \bQ_1 \bQ_2\ldots\bQ_{n-1}$,
$$
\bM_k \triangleq (\bP_{n-1} \ldots \bP_{k+1}) \bE_k (\bP_{k+1} \ldots \bP_{n-1}), \qquad \text{for all $k\in \{1,2,\ldots,n-1\}$},
$$
and 
$
\bL \triangleq \bM_1^{-1}\bM_2^{-1} \ldots \bM_{n-1}^{-1}.
$
We have $\bA = \bP\bL\bU \bQ^\top$ or $\bP^\top\bA\bQ = \bL\bU$ as the final decomposition.

\subsection{Rook Pivoting}
\textit{Rook pivoting} provides an alternative to  partial and complete pivoting. Instead of selecting the largest value in $\abs{\bU[k:n,k:n]}$ at the $k$-th stage, it searches for an element in $\abs{\bU[k:n,k:n]}$ that is \textit{maximal in both its row and column}. Apparently, the rook pivoting is not unique such that we could find many entries that satisfy the criteria. For instance, given a submatrix $\bU[k:n,k:n]$ as follows
$$
\bU[k:n,k:n] = 
\begin{bmatrix}
1 & 2 & 3 & 4 \\
2 & 3 & 7 & 3 \\
5 & 2 & 1 & 2 \\
2 & 1 & 2 & 1 \\
\end{bmatrix},
$$
where, in complete pivoting, the element $7$ would be chosen. And one among $5, 4, 7$ will be identified as a rook pivot.

\index{Rank-revealing}\index{Rank-revealing LU}
\section{Rank-Revealing LU Decomposition}\label{section:rank-reveal-lu-short}
In many applications, a factorization produced by Gaussian elimination with pivoting when $\bA$ has rank $r$ will reveal the rank in the following form:
$$
\bP\bA\bQ = 
\begin{bmatrix}
\bL_{11} & \bzero \\
\bL_{21}^\top & \bI 
\end{bmatrix}
\begin{bmatrix}
\bU_{11} & \bU_{12} \\
\bzero & \bzero 
\end{bmatrix},
$$
where $\bL_{11}\in \real^{r\times r}$ and $\bU_{11}\in \real^{r\times r}$ are nonsingular, $\bL_{21}, \bU_{21}\in \real^{r\times (n-r)}$, and $\bP,\bQ$ are permutations. Gaussian elimination with rook  or complete pivoting can result in such decompositions \citep{hwang1992rank, higham2002accuracy}.

\index{Perturbation}
\section{Rate of Change of L and U*}
If $\bA$ has a unique LU decomposition, and a small perturbation $\Delta\bA$ such that the LU decomposition of $\bA+\Delta\bA$ also exists and is given by $\bA+\Delta\bA = (\bL+\Delta\bL)(\bU+\Delta\bU)$, then $\Delta\bA$ can be obtained by (ignoring $\Delta\bL\cdot \Delta\bU$)
$$
\Delta\bA = \Delta\bL \cdot \bU + \bL \cdot  \Delta\bU.
$$
Furthermore, we have 
$$
\bL^{-1}\cdot \Delta\bA \cdot \bU^{-1} = \bL^{-1}\cdot  \Delta\bL + \Delta\bU \cdot \bU^{-1}.
$$
Since both $\bL$ and $(\bL+\Delta\bL)$ are unit lower triangular matrices, $\Delta \bL$ is \textit{strictly lower triangular}, i.e., lower triangular with zeros on the diagonal. 
Similarly, since both $\bU$ and $\bU+\Delta\bU$ are upper triangular matrices, $\Delta\bU$ is therefore upper triangular. 
Thus,
$$
\Delta\bL = \bL \cdot slt(\bL^{-1}\cdot \Delta\bA \cdot \bU^{-1}), 
\qquad 
\Delta\bU = ut(\bL^{-1}\cdot \Delta\bA \cdot \bU^{-1})  \bU,
$$
where $slt(\bB)$ denotes the strictly lower triangular part of $\bB$, and $ut(\bB)$ represents the upper triangular part of $\bB$. Clearly, the \textit{sensitivity}, a.k.a., the \textit{rate of change} of $\bL$ and $\bU$ depends on the inverse of $\bL$ and $\bU$. More generally, we have the following result.

\begin{theorem}[Rate of Change of $\bL$ and $\bU$]\label{theorem:rate-lu-decomposition}
Let $\bA\in \real^{n\times n}$ be given. Suppose  $\bA $ has nonzero leading principal minors, i.e., $\det(\bA[1:k,1:k])\neq 0$, for all $k\in \{1,2,\ldots, n\}$, with LU decomposition $\bA = \bL\bU$. Let $\bG$ be a real $n\times n$ matrix and let $\Delta\bA = \epsilon \bG$, for some $\epsilon \geq 0$. If $\epsilon$ is sufficiently small such that all the leading principal minors of $\bA+t\bG$ are nonzero for all $|t|\leq \epsilon$. Then $\bA+\Delta\bA$ has the LU decomposition 
$$
\bA+\Delta\bA = (\bL+\Delta\bL)(\bU+\Delta\bU),
$$
with $\Delta\bL$ and $\Delta\bU$ satisfying 
$$
\begin{aligned}
\Delta \bL &= \epsilon \dotL(0) + \mathcalO(\epsilon^2), \\
\Delta\bU &= \epsilon \dotU(0) + \mathcalO(\epsilon^2),
\end{aligned}
$$
where $\dotL(0) = \dotL(t)|_{t=0}$ and $\dotU(0) = \dotU(t)|_{t=0}$, and $\dotL(t)$ and $\dotU(t)$ are defined by the unique LU decomposition
\begin{equation}\label{equation:partial-eq}
\bA + t\bG = \bL(t)\bU(t), \qquad |t| \leq \epsilon.
\end{equation}
\end{theorem}

\begin{proof}[of Theorem~\ref{theorem:rate-lu-decomposition}]
We notice that 
\begin{equation}\label{equation:rage-lu-eq1}
\begin{aligned}
	\bL(t) &= \bL +\Delta\bL, \\
	\dotL(t)&=\frac{\partial \bL(t)}{\partial t} .
\end{aligned}
\end{equation}
By Taylor's expansion (Appendix~\ref{appendix:taylor-expansion}), we have 
$$
\dotL(\epsilon) =\frac{\partial \bL(\epsilon)}{\partial \epsilon} \approx \dotL(0) + \ddotL(\epsilon) \epsilon + \mathcalO(\epsilon^2).
$$
Since $\epsilon$ is small enough, Equation~\eqref{equation:rage-lu-eq1} can be written as 
$$
\dotL(\epsilon)=\frac{\partial \bL(\epsilon)}{\partial \epsilon} = \frac{\Delta\bL}{\epsilon},
$$
which results in 
$$
\Delta\bL = \epsilon \dotL(0) + \mathcalO(\epsilon^2).
$$
Similarly, we can also obtain 
$$
\Delta\bU = \epsilon \dotU(0) + \mathcalO(\epsilon^2).
$$
This completes the proof.
\end{proof}
Following from  Theorem~\ref{theorem:rate-lu-decomposition} above, we can also claim that 
\begin{align}
\bL\dotU(0) &+ \dotL(0)\bU = \bG, \label{equation:partial-lu-rate-1}\\
\dotL(0) &= \bL \cdot slt(\bL^{-1}\bG\bU^{-1}), \label{equation:partial-lu-rate-2}\\
\dotU(0) &= ut(\bL^{-1}\bG\bU^{-1}) \bU. \label{equation:partial-lu-rate-3}
\end{align}
Note that $\bL(0) = \bL$, $\bL(\epsilon) =\bL+\Delta\bL$, $\bU(0) = \bU$, and $\bU(\epsilon)=\bU +\Delta\bU$. If we differentiate Equation~\eqref{equation:partial-eq}, and set $t=0$, we can obtain the result in Equation~\eqref{equation:partial-lu-rate-1}. Further, multiplying Equation~\eqref{equation:partial-lu-rate-1} left by $\bL^{-1}$ and right by $\bU^{-1}$, we have 
$$
\dotU(0)\bU^{-1} + \bL^{-1}\dotL(0) =\bL^{-1}\bG\bU^{-1},
$$
where $\dotL(0)$ is strictly lower triangular and $\dotU(0)$ is upper triangular. Thus, the results in Equation~\eqref{equation:partial-lu-rate-2} and \eqref{equation:partial-lu-rate-3} can be proved.

In the above theorem, we prove that the sensitivity of $\Delta\bL$ and $\Delta\bU$ is bounded by $\dotL(0)$ and $\dotU(0)$ to some $\epsilon$ value. More results on the perturbation of the LU decomposition can be found in \citet{baarland1991perturbation,sun1992rounding, sun1992componentwise, stewart1993perturbation, stewart1997perturbation, chang1997pertubation, chang1998sensitivity, higham2002accuracy, bueno2004stability}.

\index{Row equivalent}
\begin{problemset}
\item Use row reduction	to solve the linear system:
$$
\begin{aligned}
2x_1 + 3x_2 + 4x_3 &= 9,\\
x_1 + 2x_2+ 3x_3 &= 5,\\
3x_1 + 4x_2 + 5x_3 &= 7.
\end{aligned}
$$
	
\item Two matrices $\bA$ and $\bB$ are said to be  \textit{row equivalent} (denoted by $\bA\stackrel{r}{\sim}\bB$) if there is a sequence of elementary row transformations  that can be performed on $\bA$ to obtain $\bB$. 
\begin{itemize}
\item Show that $\bA\stackrel{r}{\sim}\bB$ if and only if $\bA=\bP\bB$ for some nonsingular matrix $\bP$.
\item Show that if $\bA\stackrel{r}{\sim}\bC$ and $\bB\stackrel{r}{\sim}\bC$, then $\bA\stackrel{r}{\sim}\bB$.
\item Show that if $\bA\stackrel{r}{\sim}\bB$ and $\bB\stackrel{r}{\sim}\bC$, then $\bA\stackrel{r}{\sim}\bC$.
\item Show that if $\bA\stackrel{r}{\sim}\bB$, then $\bB\stackrel{r}{\sim}\bA$.
\item Show that $\bA\stackrel{r}{\sim}\bB$ if both $\bA$ and $\bB$ are nonsingular.
\end{itemize}

\item Let $\bE$ and $\bF$ be two matrices derived from the identity matrix by adding multiples of row $i$ to rows $j$ and $k$, where $i\neq j$ and $i\neq k$ respectively. Show that $\bE\bF=\bF\bE$.
	
\item Show that  the matrix $
\scriptsize
\begin{bmatrix}
0& 1\\
1 & 0
\end{bmatrix}$ does not admit an LU decomposition.

\item Let $\bL_1$ and $\bL_2$ be nonsingular lower triangular matrices, and $\bU_1$ and $\bU_2$ be nonsingular upper triangular matrices.
Prove that $\bL_1\bU_1=\bL_2\bU_2$ if and only if there exists an nonsingular diagonal matrix $\bD$ such that $\bL_1=\bL_2\bD$ and $\bU_1=\bD^{-1}\bU_2$.

\item We have demonstrated  elementary row transformations of a matrix is achieved by left multiplication of the matrix (Definition~\ref{definition:elemen_trans}). Discuss the specific forms of 1). an elementary row transformation of interchanging two rows of a matrix; 2). the one of multiplying all elements of a row by a scalar; 3). the one of adding a scalar multiple of a  row to another row. How about the three elementary column transformations? \textit{Hint: Equations~\eqref{equation:ele_exam_type2} and \eqref{equation:ele_exam_eijal}.}

\item \label{problem:row_det_chg} Consider the three types of elementary row transformation in Definition~\ref{definition:elemen_trans}. Show that the type-1 transformation on the determinant is to multiply it by $-1$; the type-2 transformation is to multiply it by the nonzero scaling factor; the type-3 transformation leaves the determinant unchanged. \textit{Hint: Use induction.}


\item \label{prob:mt_inv_gen} Following Schur complement and Lemma~\ref{lemma:fac_schur_a}, consider a matrix $\bA\in\real^{n\times n}$, an index set $I$ and its complement $J=\{1,2,\ldots,n\}/ I$. Show that 
\begin{equation}
\begin{aligned}
	\bA^{-1}[I, I] &= \left(\bA[I,I] - \bA[I, J]\bA[J,J]^{-1}\bA[J,I] \right)^{-1};\\
	\bA^{-1}[I,J] &= \bA[I,I]^{-1} \bA[I,J] \left(\bA[J,I]\bA[I,I]^{-1}\bA[I,J]-\bA[J,J]\right)^{-1}\\
	&= \left(\bA[I,J]\bA[J,J]^{-1}\bA[J,I]-\bA[I,I]\right)^{-1} \bA[I,J]\bA[J,J]^{-1},
\end{aligned}
\end{equation}
where $\bA^{-1}[I,J]$ is the submatrix of $\bA^{-1}$, and $\bA[I,I]^{-1}$ is the inverse of $\bA[I,I]$.

\item \label{prob:det_block1} \textbf{Determinant of block matrices.} Let $\bM=\scriptsize\begin{bmatrix}
	\bA& \bB\\
	\bzero & \bD 
\end{bmatrix}$, show that $\det(\bM)=\det(\bA)\det(\bD)$. \textit{Hint: Use the factorization in Schur complements; Lemma~\ref{lemma:fac_schur_a}.}

\item \label{prob:det_block2} \textbf{Determinant of block matrices.} Given the knowledge that $\det(\bM)=\det(\bA)\det(\bD)$ if $\bM=\scriptsize\begin{bmatrix}
\bA& \bB\\
\bzero & \bD 
\end{bmatrix}$.
How about the determinant of $\scriptsize\begin{bmatrix}
	\bA& \bzero\\
	\bC  & \bD 
\end{bmatrix}$?
Prove the determinant of block matrices in Remark~\ref{remark:determinant-intermezzo}. \textit{Hint: Use the factorization in Schur complements; Lemma~\ref{lemma:fac_schur_a}.}

\item \label{prob:det_block3} \textbf{Determinant of block matrices.}
Let $\bA,\bB,\bC,\bD$ be conformable matrices, show that 
$$
\begin{aligned}
&\det\big(
\footnotesize
\begin{bmatrix}
	\bA & \bA\bB \\
	\bC & \bD
\end{bmatrix}
\big)
=
\det(\bA)\det(\bD-\bC\bB);
\quad 
\det\big(
\footnotesize
\begin{bmatrix}
	\bA & \bB \\
	\bC\bA & \bD
\end{bmatrix}
\big)
=
\det(\bA)\det(\bD-\bC\bB);\\
&\det\big(
\footnotesize
\begin{bmatrix}
	\bA & \bB\bD \\
	\bC & \bD
\end{bmatrix}
\big)
=
\det(\bD)\det(\bA-\bB\bC);
\quad 
\det\big(
\footnotesize
\begin{bmatrix}
	\bA & \bB \\
	\bD\bC & \bD
\end{bmatrix}
\big)
=
\det(\bD)\det(\bA-\bB\bC).
\end{aligned}
$$
\textit{Hint: Use the factorization in block forms similar to Problem~\ref{prob:det_block1} and \ref{prob:det_block2}.}

\item \textbf{Rank of block matrices.} Let $\bM=\scriptsize\begin{bmatrix}
\bI_m & \bA\\
\bB & \bI_n 
\end{bmatrix}$, where $\bA\in\real^{m\times n}$ and $\bB\in\real^{n\times m}$. Show that $\rank(\bM)=m+\rank(\bI_n-\bB\bA)=n+\rank(\bI_m-\bA\bB)$. \textit{Hint: Use the factorization in block forms similar to Problem~\ref{prob:det_block1} and \ref{prob:det_block2}.}

\item Prove \eqref{equation:schur_inv1_det_res} based on \eqref{equation:schur_inv1_det_sour} and \eqref{equation:schur_inv1_det}.

\item \label{prob:cauch_exp} \textbf{Cauchy’s formula.} Consider $\bM=\scriptsize\begin{bmatrix}
\bA & \bb\\
\bc^\top & d
\end{bmatrix}$, when $d\neq 0$, show that 
\begin{equation}
d\det(\bA)-\bc^\top \adjugate(\bA)\bb 
=
d\det(\bA-d^{-1}\bb\bc^\top)
\end{equation}
This is known as \textit{Cauchy’s formula for the determinant of a rank-one perturbation} when $d=-1$:
\begin{equation}
\det(\bA+\bb\bc^\top) = \det(\bA) +\bc^\top \adjugate(\bA)\bb.
\end{equation}
\textit{Hint: Use \eqref{equation:det_block_rem2} and determinant of block matrices in Remark~\ref{remark:determinant-intermezzo}.}

\item \label{problem:block_gaus_uni} Prove the block Gaussian elimination in \eqref{equation:block_gauss} is unique if we want to obtain $\Delta_{\bA}$. 

\item \textbf{Recursive Schur complement.} Consider the block matrix $\bM$ in \eqref{equation:schu_set}: $\bM=\scriptsize\begin{bmatrix}
	\bA&\bB\\
	\bC&\bD
\end{bmatrix}$, where $\bA=\scriptsize\begin{bmatrix}
\bA_{11} &\bA_{12} \\
\bA_{21} &\bA_{22}
\end{bmatrix}$, $\bC=\begin{bmatrix}
\bC_{1} & \bC_{2}
\end{bmatrix}$, and 
$\bB=\scriptsize\begin{bmatrix}
\bB_1\\
\bB_2 
\end{bmatrix}$:
$$
\bM=
\begin{bmatrix}
\bA_{11} &\bA_{12} & \bB_1  \\
\bA_{21} &\bA_{22} & \bB_2 \\
\bC_{1} & \bC_{2} & \bD 
\end{bmatrix}.
$$
Show that $\Delta_{\bA} = f(\Delta_{\bA_{11}}, \delta_{\bA_{11}})$, where 
$\Delta_{\bA}$ denotes the Schur complement of $\bA$ in $\bM$, 
$\Delta_{\bA_{11}}$ denotes the Schur complement of $\bA_{11}$ in $\bM$,  
$\delta_{\bA_{11}}$ denotes the Schur complement of $\bA_{11}$ in $\bA$, 
and $f(\Delta_{\bA_{11}}, \delta_{\bA_{11}})$ denotes the Schur complement of  $\delta_{\bA_{11}}$ in $\Delta_{\bA_{11}}$.

\index{Strassen's algorithm}
\item \label{prob:strassen} \textbf{Strassen's algorithm.} Let $\bA,\bB\in\real^{2^n\times 2^n}$, where $N\triangleq 2^n$. Section~\ref{section:compl_mvops} shows that the matrix multiplication $\bA\bB$ costs $\mathcalO(8\cdot 2^{3n-2})=\mathcalO(2\cdot N^3)$ flops. 
Consider their block forms:
$\bA=\scriptsize\begin{bmatrix}
	\bA_{11} & \bA_{12} \\ 
	\bA_{21} & \bA_{22}
\end{bmatrix}$ and  $\bB=\scriptsize\begin{bmatrix}
	\bB_{11} & \bB_{12} \\ 
	\bB_{21} & \bB_{22}
\end{bmatrix}$, where each block is of size $2^{n-1}\times 2^{n-1}$. Then $\bA\bB =\scriptsize \begin{bmatrix}
\bA_{11}\bB_{11}+\bA_{12}\bB_{21} & \bA_{11} \bB_{12} + \bA_{12}\bB_{22} \\
\bA_{21}\bB_{11}+\bA_{22}\bB_{21} & \bA_{21} \bB_{12} + \bA_{22}\bB_{22}
\end{bmatrix}
\triangleq
\begin{bmatrix}
\bC_{11} & \bC_{12}\\
\bC_{21} & \bC_{22}
\end{bmatrix}.
$
Show that the complexity remains the same using this block calculation.
Consider further 
$$
\begin{aligned}
\bS_1 &\triangleq (\bA_{11} + \bA_{22})(\bB_{11} + \bB_{22});  \quad
&\bS_2 &\triangleq (\bA_{21} + \bA_{22}) \bB_{11}; \\
\bS_3 &\triangleq \bA_{11}(\bB_{12} - \bB_{22}) ;  \quad
&\bS_4 &\triangleq \bA_{22}(\bB_{21} - \bB_{11}); \\
\bS_5 &\triangleq (\bA_{11} + \bA_{12}) \bB_{22} ;  \quad
&\bS_6 &\triangleq (\bA_{21} - \bA_{11})(\bB_{11} + \bB_{12}) ;\\
\bS_7 &\triangleq (\bA_{12} - \bA_{22})(\bB_{21} + \bB_{22}).
\end{aligned}
$$
Show that 
$$
\begin{aligned}
\bC_{11} &= \bS_{1} + \bS_{4} - \bS_{5} + \bS_{7};  \quad
&\bC_{12} &= \bS_{3} + \bS_{5}; \\
\bC_{21} &= \bS_{2} + \bS_{4};  \quad
&\bC_{22} &= \bS_{1} - \bS_{2} + \bS_{3} + \bS_{6}.
\end{aligned}
$$
And show that the complexity is $\mathcalO(7\cdot 2^{3n-2})$. 
Let $f(n)$ be the complexity of the matrix multiplication of two $2^n\times 2^n$ matrices. It follows that $f(n) = 7f(n-1)$ if we use this block calculation approach. This shows the complexity reduces to $\mathcalO(2\cdot N^{\log_2^7})$ if we continue to use this block calculation, down to a series of 1-by-1 matrices.

\end{problemset}

%% file: chapter-cholesky.tex
\chapter{Cholesky Decomposition}
\begingroup
\hypersetup{
	linkcolor=structurecolor,
	linktoc=page,  
}
\minitoc \newpage
\endgroup

\section{Cholesky Decomposition}
\lettrine{\color{caligraphcolor}P}
Positive definiteness or positive semidefiniteness is one of the highest accolades to which a matrix can aspire. 
In this chapter, we will introduce decompositional approaches for the two special kinds of matrices and we first illustrate the most famous Cholesky decomposition as follows.

\index{Decomposition: Cholesky}
\begin{theoremHigh}[Cholesky Decomposition]\label{theorem:cholesky-factor-exist}
Let $\bA\in \real^{n\times n}$ be any \textit{positive definite} (PD) matrix. Then, it can be factored as 
$$
\bA = \bR^\top\bR,
$$
where $\bR \in \real^{n\times n}$ is an upper triangular matrix \textbf{with strictly positive diagonal elements}~\footnote{If $\bA$ is positive semidefinite (PSD), then the diagonal entries can be chosen nonnegative.}. This decomposition is known as the \textit{Cholesky decomposition}  of $\bA$, and $\bR$ is called the \textit{Cholesky factor} or \textit{Cholesky triangle} of $\bA$.

Alternatively, $\bA$ can be factored as $\bA=\bL\bL^\top$, where $\bL=\bR^\top$ is a lower triangular matrix \textit{with strictly positive diagonals}.

Specifically, the Cholesky decomposition is \textbf{unique} (Corollary~\ref{corollary:unique-cholesky-main}).
In cases where the diagonal elements of $\bR$ are not restricted to positive values, then the factorization $\bA=\bR^\top\bR$ is \textbf{not unique}.
\footnote{The Cholesky decomposition also holds for \textit{complex Hermitian positive definite} matrix $\bA\in\complex^{n\times n}$. In such cases, $\bA=\bR^\ast\bR$ for a unique upper triangular matrix $\bR$ with positive diagonal entries.}
\end{theoremHigh}

The Cholesky decomposition is named after the French military officer and mathematician, \textit{Andr\'{e}-Louis Cholesky} (1875{\textendash}1918), who developed the Cholesky decomposition in his surveying work. Similar to the LU decomposition for solving linear systems, the Cholesky decomposition is further used primarily for solving positive definite linear systems. The development on the solution is similar to that of the LU decomposition, as detailed in Section~\ref{section:lu-linear-sistem}, and we shall not repeat the details here.

%
%

\section{Existence of  Cholesky Decomposition via Recursive Calculation}

In this section, we will demonstrate  the existence of the Cholesky decomposition through recursive calculations. Additionally, in Section~\ref{section:cholesky-by-qr-spectral}, we will also provide an alternative proof  of the Cholesky decomposition's existence through  QR  and spectral decompositions.

Before we proceed to show the existence of the Cholesky decomposition,  we introduce the necessary definitions and lemmas.
\begin{definition}[Positive Definite and Positive Semidefinite\index{Positive definite}\index{Positive semidefinite}]\label{definition:psd-pd-defini}
A symmetric matrix $\bA\in \real^{n\times n}$ ($\bA\in\symmetric^n$) is considered \textit{positive definite (PD)} if $\bx^\top\bA\bx>0$ for all nonzero $\bx\in \real^n$, denoted by $\bA\succ \bzero$ or $\bA\in\pd^n$.
And a matrix $\bA\in \real^{n\times n}$ is called \textit{positive semidefinite (PSD)} if $\bx^\top\bA\bx \geq 0$ for all $\bx\in \real^n$,  denoted by $\bA\succeq\bzero$ or $\bA\in\psd^n$. 
\footnote{
In this book, a positive definite or a semidefinite matrix is always assumed to be symmetric, i.e., the notion of a positive definite matrix or semidefinite matrix is only interesting for symmetric matrices. And we only consider $\bx\in\real^n$; while the complex case for PD matrices are equivalent; see Problem~\ref{prob:pd_cmequiv}.
}
\footnote{Similarly, a complex matrix $\bA$ is said to be \textit{Hermitian positive definite (HSD)}  if $\bA$ is Hermitian (Definition~\ref{definition:complex_special}) and $\bz^\ast\bA\bz>0$ for all $\bz\in\complex^n$ with $\bz\neq \bzero$. However, the complex analog of  PD/PSD can be defined for non-Hermitian matrices. We will mostly consider PD/PSD for real symmetric matrices in this book.
}
\footnote{A symmetric matrix $\bA\in\real^{n\times n}$ is called \textit{negative definite (ND) } if $\bx^\top\bA\bx<0$ for all nonzero $\bx\in\real^n$, denoted by $\bA\prec\bzero$ or $\bA\in\nd^n$; 
a symmetric matrix $\bA\in\real^{n\times n}$ is called \textit{negative semidefinite (NSD)}  if $\bx^\top\bA\bx\leq 0$ for all $\bx\in\real^n$, denoted by $\bA\preceq\bzero$ or $\bA\in\nsd^n$;
and a symmetric matrix $\bA\in\real^{n\times n}$ is called \textit{indefinite (ID)}  if there exist $\bx$ and $\by\in\real^n$ such that $\bx^\top\bA\bx<0$ and $\by^\top\bA\by>0$.
 }
\end{definition}
For any matrix $\bA\in\real^{n\times n}$ (not necessarily symmetric), let $P(\bA) \triangleq \frac{1}{2} (\bA+\bA^\top)$ (the symmetric part of $\bA$),  we have 
$$
\bx^\top\bA\bx = \bx^\top P(\bA)\bx= \sum_{i,j=1}^{n,n} x_ix_j a_{ij}, \gap \forall \bx\in\real^n.
$$
This explains why the concept of a positive definite matrix is meaningful primarily for symmetric matrices.
\begin{exercise}[Principal Submatrices of PD]
Show that all principal submatrices of a PD (resp. PSD) matrix are also PD (resp. PSD).
\end{exercise}

\begin{exercise}
Let $\bA\in\real^{n\times n}$.
Show that $\bx^\top\bA\bx=0$ for all $\bx\in\real^n$ if and only if $\bA$ is skew-symmetric (Definition~\ref{definition:speci_mat}).
If further $\bA$ is symmetric, show that $\bx^\top\bA\bx=0$ for all $\bx\in\real^n$ if and only if $\bA=\bzero$.
\end{exercise}

One of the prerequisites for the Cholesky decomposition is that the matrix must be positive definite, as defined above.
We now outline several properties of a positive definite matrix as follows.
The first property demonstrates that PD and PSD matrices have special forms of eigenvalues.
\begin{tcolorbox}[title={Positive Definite Matrix Property 1 of 8: An Alternative Definition of PD},colback=\mdframecolorTheorem]

\begin{theorem}[Eigenvalue Characterization Theorem]\label{theorem:eigen_charac}
A matrix $\bA$ is positive definite if and only if it contains only \textit{positive eigenvalues}. Similarly, a matrix $\bA$ is positive semidefinite if and only if it contains only  \textit{nonnegative eigenvalues}.~\footnote{The trace, determinant, and principal minors of a positive (semi)definite matrix is discussed in Problem~\ref{prob:tr_de_pd}.}
Moreover, we have the following implications:
\begin{itemize}
\item $\bA-\gamma\bI\succeq \bzero$ if and only if $\lambda_{\min}(\bA) \geq \gamma$;
\item  $\bA-\gamma\bI\succ \bzero$ if and only if $\lambda_{\min}(\bA) > \gamma$;
\item $\bA-\gamma\bI\preceq \bzero$ if and only if $\lambda_{\max}(\bA) \leq \gamma$;
\item $\bA-\gamma\bI\prec \bzero$ if and only if $\lambda_{\max}(\bA) < \gamma$;
\item $\lambda_{\min}(\bA)\bI\preceq \bA \preceq \lambda_{\max}(\bA)\bI$,
\end{itemize}
where $\lambda_{\min}(\bA)$ and $\lambda_{\max}(\bA)$ represent the minimum and maximum eigenvalues of $\bA$, respectively, and $\bB \prec \bC$ means $\bC-\bB$ is PSD.
\end{theorem}
\end{tcolorbox}
The forward implication can be shown that $\bx^\top\bA\bx=\lambda\bx^\top\bx>0$ such that $ \lambda=(\bx^\top\bA\bx)/(\bx^\top\bx)>0$ (resp. $\geq 0$) if $\bA$ is PD (resp. PSD).
The full proof of this equivalence can  be found in Section~\ref{section:equivalent-pd-psd}, based on the spectral theorem.
This lemma provides an alternative definition of positive definiteness and positive semidefiniteness in terms of the eigenvalues of the matrix, which is a fundamental property for the Cholesky decomposition.
More eigenvalue-related results of PD/PSD matrices are discussed in Section~\ref{section:eig_ine_pdpsd} as a consequence of the spectral decomposition, e.g.,  Weyl's inequality for symmetric matrices in Theorem~\ref{theorem:psd_diff} that relates the eigenvalues of two matrices when their difference is positive semidefinite.

\begin{exercise}[Power of PSD]
Let $\bA$ be PSD. Show that  $\bA^k$ is also PSD for $k=1,2,\ldots$.
\end{exercise}

Although not all components of a positive definite matrix are necessarily positive, the diagonal components of such a matrix are guaranteed to be positive, as stated in the following result.
\begin{tcolorbox}[title={Positive Definite Matrix Property 2 of 8},colback=\mdframecolorTheorem]
\begin{theorem}[Positive Diagonals of Positive Definite Matrices]\label{theorem:positive-in-pd}
The diagonal elements of a positive definite matrix $\bA\in\real^{n\times n}$ are all \textit{positive}. And similarly, the diagonal elements of a positive semidefinite matrix $\bB\in\real^{n\times n}$ are all \textit{nonnegative}.
\end{theorem}
\end{tcolorbox}
\begin{proof}[of Theorem~\ref{theorem:positive-in-pd}]
By the definition of positive definite matrices, we have $\bx^\top\bA \bx >0$ for all nonzero vectors $\bx$. Consider the specific case where $\bx=\be_i$ with $\be_i$ being the $i$-th unit basis vector having the $i$-th entry equal to 1 and all other entries equal to 0. Then,
$$
\be_i^\top\bA \be_i = a_{ii}>0, \gapforall \forall i \in \{1, 2, \ldots, n\},
$$	
where $a_{ii}$ is the $i$-th diagonal component. The proof for the second part follows a similar argument. This completes the proof.
\end{proof}

It can be proposed that if a matrix contains both positive and negative diagonal elements, then the matrix is indefinite.
\begin{proposition}[Diagonals of Indefinite Matrices]\label{proposition:indefin_diag}
Given a symmetric matrix $\bA\in\real^{n\times n}$, if there exist both positive and negative elements in its diagonal, then $\bA$ is indefinite.
\end{proposition}
\begin{proof}[of Proposition~\ref{proposition:indefin_diag}]
Suppose $a_{ii}>0$ and $a_{jj}<0$ for $i\neq j, i,j\in\{1,2,\ldots, n\}$. Then it follows that 
$$
\be_i^\top\bA\be_i =a_{ii}>0
\gap \text{and}\gap 
\be_j^\top\bA\be_j=a_{jj}<0,
$$
from which the result follows.
\end{proof}
However, the reverse claim of Proposition~\ref{proposition:indefin_diag} is not true; see Problem~\ref{problem:indefinite_reverse}.

The complete pivoting in the Cholesky decomposition is simpler compared to that in  the LU decomposition (Section~\ref{section:complete-pivoting}). This simplicity arises because the maximum value in a positive definite matrix is located on its diagonal.
\begin{tcolorbox}[title={Positive Definite Matrix Property 3 of 8},colback=\mdframecolorTheorem]
\begin{theorem}[Maximal Value of Positive Definite Matrices]\label{theorem:positive-maximal-diagonal}
The largest element in a positive definite matrix lies on the diagonal. A similar argument applies to positive semidefinite matrices.
\end{theorem}
\end{tcolorbox}
\begin{proof}[of Theorem~\ref{theorem:positive-maximal-diagonal}]
Let $\be_i$ and $\be_j$ represent the $i$-th and $j$-th unit vectors, respectively. Moreover, let $a_{ij}$ denote the entry at position $(i,j)$ in the positive definite matrix $\bA$. Then, it follows that 
$$
(\be_i-\be_j)^\top\bA(\be_i-\be_j) = a_{ii}+a_{jj}-2a_{ij}>0.
$$
Thus, either $a_{ii}$ or $a_{jj}$ must be greater than $a_{ij}$. 
By iterating through all entries, the conclusion is reached. Regarding positive semidefinite matrices, it can be deduced that the largest element appears on the diagonal, and there may be non-diagonal elements equal to the maximum value.
\end{proof}

\begin{exercise}
Given an $n\times n$ positive semidefinite matrix $\bA$, show that:
\begin{itemize}
\item For any $i\neq j$, it holds that $a_{ii}a_{jj} \geq a_{ij}^2$.
\item If $a_{ii}=0$ for some $i\in\{1,2,\ldots,n\}$, then the $i$-th row of $\bA$ is a zero vector.
\end{itemize}
Hint: use the result of Theorem~\ref{theorem:positive-maximal-diagonal}  and the AM-GM inequality (Equation~\eqref{equation:amgm_ineq}).
\end{exercise}

\begin{tcolorbox}[title={Positive Definite Matrix Property 4 of 8},colback=\mdframecolorTheorem]
Any positive definite matrix $\bA$ is nonsingular. This can be shown that if $\bA\bx=\bzero$, then $\bx^\top\bA\bx=0$. This implies $\bx$ must be $\bzero$. 
Therefore, $\bA$ is of full rank.
\end{tcolorbox}
\begin{exercise}[Singular PSD]\label{exercise:sing_psd}
Similarly, let $\bA$ be positive semidefinite. Show that $\bx^\top\bA\bx=0$ if and only if $\bA\bx=\bzero$.
Show that a positive semidefinite $\bA$ is positive definite if and only if it is nonsingular. See Problem~\ref{prob:quad_nsdpsd} for more implication.
\end{exercise}

Similar to the LU decomposition, the existence of the Cholesky decomposition also depends on the Schur complement.
\begin{tcolorbox}[title={Positive Definite Matrix Property 5 of 8},colback=\mdframecolorTheorem]
\begin{theorem}[Schur Complement of Positive Definite Matrices\index{Schur complement}]\label{theorem:pd-of-schur}
Let $\bA\in \real^{n\times n}$ be  any positive definite (resp. positive semidefinite) matrix. Then, its Schur complement of $a_{11}$ is given by $\bS_{n-1}=\bA_{2:n,2:n}-\frac{1}{a_{11}} \bA_{2:n,1}\bA_{2:n,1}^\top$, which is also positive definite (resp. positive semidefinite). 

\paragraph{Notation clarification.} Note that the subscript $n-1$ of $\bS_{n-1}$ indicates it is of size $(n-1)\times (n-1)$, implying that it is a Schur complement derived from an $n\times n$ positive definite matrix. 
This notation will be consistently utilized in subsequent sections.
\end{theorem}	
\end{tcolorbox}

\begin{proof}[of Theorem~\ref{theorem:pd-of-schur}]
For any nonzero vector $\bv\in \real^{n-1}$, we can construct a vector $\bx\in \real^n$ as follows:
$$
\bx = 
\begin{bmatrix}
	-\frac{1}{a_{11}} \bA_{2:n,1}^\top  \bv \\
	\bv
\end{bmatrix},
$$
which is nonzero. Then we have
$$
\setlength{\arraycolsep}{3pt}
\footnotesize
\begin{aligned}
\bx^\top\bA\bx 
&= 
\left[-\frac{\bv^\top \bA_{2:n,1}}{a_{11}} , \bv^\top\right]
\begin{bmatrix}
a_{11} & \bA_{2:n,1}^\top \\
\bA_{2:n,1} & \bA_{2:n,2:n}
\end{bmatrix}
\begin{bmatrix}
-\frac{1}{a_{11}} \bA_{2:n,1}^\top  \bv \\
\bv
\end{bmatrix} 
= \left[-\frac{\bv^\top \bA_{2:n,1}}{a_{11}} , \bv^\top\right]
\begin{bmatrix}
0 \\
\bS_{n-1}\bv
\end{bmatrix} 
= \bv^\top\bS_{n-1}\bv.
\end{aligned}
$$
Since $\bA$ is positive definite, the condition $\bx^\top\bA\bx = \bv^\top\bS_{n-1}\bv >0$ holds for all nonzero vectors $\bv$. 
Consequently, the Schur complement $\bS_{n-1}$ is also positive definite.
The proof for the positive semidefinite part follows similarly.
\end{proof}


\index{Schur complement}
\paragraph{A word on the Schur complement.} In the proof of Theorem~\ref{theorem:lu-factorization-with-permutation}, we demonstrated that the Schur complement $\bS_{n-1}=\bA_{2:n,2:n}-\frac{1}{a_{11}} \bA_{2:n,1}\bA_{2:n,1}^\top$ remains nonsingular if $\bA$ is nonsngular and $a_{11}\neq 0$. Similarly, the Schur complement of $a_{nn}$ in $\bA$ is given by $\bar{\bS}_{n-1} =\bA_{1:n-1,1:n-1} - \frac{1}{a_{nn}}\bA_{1:n-1,n} \bA_{1:n-1,n}^\top$, which is also positive definite if $\bA$ is positive definite. This property can help prove the fact that the leading principal minors of positive definite matrices are all positive. See Section~\ref{appendix:leading-minors-pd}  for more details. 
\begin{exercise}[Schur Complement of Positive Definite Matrices]\label{exercise:schu_block}
More generally, let 
\begin{equation}
\bM = \begin{bmatrix}
	\underset{p\times p}{\bA} & \underset{p\times q}{\bB} \\
	\underset{q\times p}{\bC} & \underset{q\times q}{\bD}
\end{bmatrix}.
\end{equation}
If $\bD$ is invertible, $\Delta_{\bD} = \bA -\bB \bD^{-1}\bC$ is called the Schur complement of $\bD$ in $\bM$; if $\bA$ is invertible, $\Delta_{\bA} = \bD - \bC\bA^{-1}\bB$ is called the Schur complement of $\bA$ in $\bM$ (see Section~\ref{section:schur-complement}).
Show that 
\begin{itemize}
\item If $\bA \succ \bzero$, then $\bA,\bD\succ \bzero$; if $\bA \succeq \bzero$, then $\bA,\bD\succeq \bzero$.
\item $\bM \succ \bzero \iff \bA, \Delta_{\bA} \succ \bzero \iff \bD, \Delta_{\bD} \succ \bzero.$
\item If $\bA \succ \bzero$, then $\bM \succeq \bzero \iff \Delta_{\bA} \succeq \bzero.$
\item If $\bD \succ \bzero$, then $\bM \succeq \bzero \iff \Delta_{\bD} \succeq \bzero.$
\end{itemize}
\end{exercise}

\begin{tcolorbox}[title={Positive Definite Matrix Property 6 of 8},colback=\mdframecolorTheorem]
The set of $n\times n$  positive definite matrices is convex. That is, for any $\lambda\in\real$ such that $0\leq \lambda \leq 1$, the matrix $\lambda \bA + (1-\lambda)\bB$ is positive definite if $\bA,\bB\in\real^{n\times n}$ are both  positive definite. This can be seen that 
$$
\bx^\top\left(\lambda \bA + (1-\lambda)\bB\right)\bx
=\lambda \bx^\top\bA\bx + (1-\lambda ) \bx^\top\bB\bx 
 >0
$$
since one of $\lambda, (1-\lambda)$ is positive, and $\bx^\top\bA\bx , \bx^\top\bB\bx>0$.
\end{tcolorbox}

\index{Recursive algorithm}
We proceed to demonstrate the existence of the Cholesky decomposition using these lemmas.

\begin{proof}[{of Theorem~\ref{theorem:cholesky-factor-exist}: Existence of Cholesky Decomposition Recursively}]
Given a positive definite matrix $\bA$, we can express it explicitly (since $a_{11}$ is positive according to  Theorem~\ref{theorem:positive-in-pd}):
$$
\footnotesize
\begin{aligned}
\bA &= 
\begin{bmatrix}
	a_{11} & \bA_{2:n,1}^\top \\
	\bA_{2:n,1} & \bA_{2:n,2:n}
\end{bmatrix} 
=\begin{bmatrix}
	\sqrt{a_{11}} &\bzero\\
	\frac{1}{\sqrt{a_{11}}} \bA_{2:n,1} &\bI 
\end{bmatrix}
\begin{bmatrix}
	\sqrt{a_{11}} & \frac{1}{\sqrt{a_{11}}}\bA_{2:n,1}^\top \\
	\bzero & \bA_{2:n,2:n}-\frac{1}{a_{11}} \bA_{2:n,1}\bA_{2:n,1}^\top
\end{bmatrix}\\
&=\begin{bmatrix}
	\sqrt{a_{11}} &\bzero\\
	\frac{1}{\sqrt{a_{11}}} \bA_{2:n,1} &\bI 
\end{bmatrix}
\begin{bmatrix}
	1 & \bzero \\
	\bzero & \bA_{2:n,2:n}-\frac{1}{a_{11}} \bA_{2:n,1}\bA_{2:n,1}^\top
\end{bmatrix}
\begin{bmatrix}
	\sqrt{a_{11}} & \frac{1}{\sqrt{a_{11}}}\bA_{2:n,1}^\top \\
	\bzero & \bI
\end{bmatrix}
=\bR_1^\top
\begin{bmatrix}
	1 & \bzero \\
	\bzero & \bS_{n-1}
\end{bmatrix}
\bR_1,
\end{aligned}
$$
where  
$\bR_1 \triangleq 
\footnotesize
\begin{bmatrix}
\sqrt{a_{11}} & \frac{1}{\sqrt{a_{11}}}\bA_{2:n,1}^\top \\
\bzero & \bI
\end{bmatrix}.
$
Having established the positive definiteness of the Schur complement $\bS_{n-1}$ in Theorem~\ref{theorem:pd-of-schur}, we can  factorize it in the same way as
$$
\bS_{n-1}=
\widehat{\bR}_2^\top
\begin{bmatrix}
1 & \bzero \\
\bzero & \bS_{n-2}
\end{bmatrix}
\widehat{\bR}_2,
\quad\text{with positive definite $\bS_{n-2}$}.
$$
Therefore, we have
$$
\setlength{\arraycolsep}{1.2pt}
\footnotesize
\begin{aligned}
\bA &= \bR_1^\top
\begin{bmatrix}
	1 & \bzero \\
	\bzero & \widehat{\bR}_2^\top
	\scriptsize
	\begin{bmatrix}
		1 & \bzero \\
		\bzero & \bS_{n-2}
	\end{bmatrix}
	\footnotesize
	\widehat{\bR}_2
\end{bmatrix}
\bR_1
=
\bR_1^\top
\begin{bmatrix}
	1 &\bzero \\
	\bzero &\widehat{\bR}_2^\top
\end{bmatrix}
\begin{bmatrix}
	1 &\bzero \\
	\bzero &
	\scriptsize\begin{bmatrix}
		1 & \bzero \\
		\bzero & \bS_{n-2}
	\end{bmatrix}
	\footnotesize
\end{bmatrix}
\begin{bmatrix}
	1 &\bzero \\
	\bzero &\widehat{\bR}_2
\end{bmatrix}
\bR_1
=
\bR_1^\top \bR_2^\top
\begin{bmatrix}
	1 &\bzero \\
	\bzero &\scriptsize\begin{bmatrix}
		1 & \bzero \\
		\bzero & \bS_{n-2}
	\end{bmatrix}
	\footnotesize
\end{bmatrix}
\bR_2 \bR_1.
\end{aligned}
$$
The same formula can be recursively applied. This process gradually continues down to the bottom right corner, yielding the decomposition
$$
\begin{aligned}
\bA &= \bR_1^\top\bR_2^\top\ldots \bR_n^\top \bR_n\ldots \bR_2\bR_1
= \bR^\top \bR,
\end{aligned}
$$
where $\bR_1, \bR_2, \ldots, \bR_n$ are upper triangular matrices with positive diagonal elements. Their product, $\bR\triangleq\bR_1\bR_2\ldots\bR_n$, is also an upper triangular matrix with positive diagonal elements, from which the result follows.
\end{proof}
The process in the proof can also be used to calculate the Cholesky decomposition and to determine the complexity of the algorithm. In Section~\ref{section:compute-cholesky}, we will undertake the computation in a similar manner but from a different perspective.  

The reverse of the Cholesky decomposition also holds.
\begin{lemma}[$\bR^\top\bR$ is PD]\label{lemma:r-to-pd}
Let $\bR$ be  any upper triangular matrix  with positive diagonal elements. Then, 
$
\bA = \bR^\top\bR
$
is positive definite.
\end{lemma}
\begin{proof}[of Lemma~\ref{lemma:r-to-pd}]
If an upper triangular matrix $\bR$ has positive diagonals, it possesses  full column rank. And the null space of $\bR$ is of dimension 0 by the fundamental theorem of linear algebra (discussed in Appendix~\ref{appendix:fundamental-rank-nullity}). As a result, $\bR\bx \neq \bzero$ for any nonzero vector $\bx$. Thus, $\bx^\top\bA\bx = \normtwo{\bR\bx}^2 >0$ for any nonzero vector $\bx$.
\end{proof}
This lemma   applies not only to  upper triangular matrices $\bR$ but also to any matrix $\bR$ with linearly independent columns.

\paragraph{A word on the two claims.} Combining Theorem~\ref{theorem:cholesky-factor-exist}  and Lemma~\ref{lemma:r-to-pd}, we can claim that matrix $\bA$ is positive definite if and only if $\bA$ can be factored as $\bA=\bR^\top\bR$, where $\bR$ is an upper triangular matrix with positive diagonals.

\section{Sylvester's Criterion: Leading Principal Minors of PD Matrices}\label{appendix:leading-minors-pd}

In Theorem~\ref{theorem:pd-of-schur}, we established that for any positive definite matrix $\bA\in \real^{n\times n}$, its Schur complement of $a_{11}$ is $\bS_{n-1}=\bA_{2:n,2:n}-\frac{1}{a_{11}} \bA_{2:n,1}\bA_{2:n,1}^\top$ and it is also positive definite.
This property extends to  its Schur complement of $a_{nn}$: $\bS_{n-1}^\prime = \bA_{1:n-1,1:n-1} -\frac{1}{a_{nn}} \bA_{1:n-1,n}\bA_{1:n-1,n}^\top$  also preserves positive definiteness.

We then claim that all the leading principal minors (Definition~\ref{definition:leading-principle-minors}) of a positive definite matrix $\bA \in \real^{n\times n}$ are positive. This assertion is also known as \textit{Sylvester's criterion} \citep{swamy1973sylvester, gilbert1991positive}. Recall that these positive leading principal minors imply the existence of the LU decomposition for positive definite matrix $\bA$ by Theorem~\ref{theorem:lu-factorization-without-permutation}.

In order to demonstrate Sylvester's criterion, the following lemma is required.
\begin{tcolorbox}[title={Positive Definite Matrix Property 7 of 8},colback=\mdframecolorTheorem]
\begin{theorem}[Quadratic PD]\label{theorem:quadratic-pd}
Let $\bA$ be given, and  let $\bE$ be any    invertible matrix. Then,  the matrix $\bA$ is positive definite (resp. positive semidefinite) if and only if $\bE^\top\bA\bE$ is also positive definite (resp. positive semidefinite).
\end{theorem}
\end{tcolorbox}
\begin{proof}[of Theorem~\ref{theorem:quadratic-pd}]
Suppose $\bA$ is positive definite; then, for any nonzero vector $\bx$, we have $\bx^\top \bE^\top\bA\bE \bx = \by^\top\bA\by > 0$ since $\bE$ is invertible such that $\bE\bx$ is nonzero.~\footnote{Since the null space of $\bE$ is of dimension 0 and the only solution for $\bE\bx=\bzero$ is the trivial solution $\bx=\bzero$.} This implies $\bE^\top\bA\bE$ is PD.
For the converse claim, suppose $\bE^\top\bA\bE$ is positive definite; then, for any nonzero $\bx$, we have $\bx^\top \bE^\top\bA\bE \bx>0$. For any nonzero $\by$, there exists a nonzero $\bx$ such that $\by =\bE\bx$ since $\bE$ is invertible. This implies $\bA$ is PD as well.
The PSD case can be derived similarly.
\end{proof}
More generally, we can show the following results.
\begin{exercise}[PSD/PD with Congruence]\label{exer:congru_pd}
Let $\bA\in\real^{n\times n}$ and $\bS\in\real^{n\times m}$. Show that 
\begin{itemize}
\item If $\bA$ is PSD, then $\bS^\top\bA\bS$ is PSD,  $\rank(\bS^\top\bA\bS)=\rank(\bA\bS)$, and $\nspace(\bS^\top\bA\bS)=\nspace(\bA\bS)$.
\item If $\bA$ is PD, then $\rank(\bS^\top\bA\bS)=\rank(\bS)$, and $\bS^\top\bA\bS$ is PD if and only if $\rank(\bS)=m$.
\end{itemize}
\end{exercise}

We now provide a rigorous proof for Sylvester's criterion.
\begin{tcolorbox}[title={Positive Definite Matrix Property 8 of 8},colback=\mdframecolorTheorem]
\begin{theorem}[Sylvester's Criterion\index{Sylvester's criterion}]\label{theorem:sylvester-criterion}
Let $\bA\in \real^{n\times n}$ be any  real symmetric matrix. Then, $\bA$ is positive definite if and only if all of its leading principal minors  are positive. 
\end{theorem}
\end{tcolorbox}

\begin{proof}[of Theorem~\ref{theorem:sylvester-criterion}]
We will prove this by forward and reverse implication separately as follows.
\paragraph{Forward implication.} We will prove  the forward implication by induction. Suppose $\bA$ is positive definite. Since all the   diagonal components of positive definite matrices are  positive (Theorem~\ref{theorem:positive-in-pd}), the base case for $n=1$ is trivial that $\det(\bA)> 0$ if $\bA$ is a scalar.

Suppose that all the leading principal minors of $k\times k$ PD matrices are  positive. We need to show that this holds for $(k+1)\times (k+1)$ PD matrices.

For a $(k+1)\times (k+1)$ matrix with the block form $\bM=\scriptsize\begin{bmatrix}
\bA & \bb\\
\bb^\top & d
\end{bmatrix}$, where $\bA$ is a $k\times k$ submatrix. Then its Schur complement of $d$, denoted as $\bS_{k} = \bA - \frac{1}{d} \bb\bb^\top$, is also positive definite and its determinant is positive from the assumption. Therefore, $\det(\bM) = \det(d)\det( \bA - \frac{1}{d} \bb\bb^\top) $=  
$d\cdot \det( \bA - \frac{1}{d} \bb\bb^\top)>0$ (determinant of block matrices, Remark~\ref{remark:determinant-intermezzo}), which completes the proof.

\paragraph{Reverse implication.} Conversely, suppose all the leading principal minors of $\bA\in \real^{n\times n}$ are positive, i.e., all the leading principal submatrices are nonsingular. 
Further, denote the $(i,j)$-th entry of $\bA$ by $a_{ij}$. We realize that $a_{11}>0$ by the assumption.
We can subtract multiples of the first row of $\bA$ to zero out the entries in the first column of $\bA$ below the first diagonal $a_{11}$. That is,
$$
\bA = 
\begin{bmatrix}
a_{11} & a_{12} & \ldots &a_{1n}\\
a_{21} & a_{22} & \ldots &a_{2n}\\
\vdots & \vdots & \ddots &\vdots\\
a_{n1} & a_{n2} & \ldots &a_{nn}\\
\end{bmatrix}
\stackrel{\bE_1 \bA}{\longrightarrow}
\begin{bmatrix}
a_{11} & a_{12} & \ldots &a_{1n}\\
0 & a_{22} & \ldots &a_{2n}\\
\vdots & \vdots & \ddots &\vdots\\
0 & a_{n2} & \ldots &a_{nn}\\
\end{bmatrix}.
$$

Next, subtract multiples of the first column of $\bE_1\bA$  from the other columns of $\bE_1\bA$ to zero out the entries in the first row of $\bE_1\bA$ to the right of the first column. Since $\bA$ is symmetric, right-multiplying by $\bE_1^\top$ accomplishes this:
$$
\bA = 
\begin{bmatrix}
a_{11} & a_{12} & \ldots &a_{1n}\\
a_{21} & a_{22} & \ldots &a_{2n}\\
\vdots & \vdots & \ddots &\vdots\\
a_{n1} & a_{n2} & \ldots &a_{nn}\\
\end{bmatrix}
\stackrel{\bE_1 \bA}{\longrightarrow}
\begin{bmatrix}
a_{11} & a_{12} & \ldots &a_{1n}\\
0 & a_{22} & \ldots &a_{2n}\\
\vdots & \vdots & \ddots &\vdots\\
0 & a_{n2} & \ldots &a_{nn}\\
\end{bmatrix}
\stackrel{\bE_1 \bA\bE_1^\top}{\longrightarrow}
\begin{bmatrix}
a_{11} & 0 & \ldots &0\\
0 & a_{22} & \ldots &a_{2n}\\
\vdots & \vdots & \ddots &\vdots\\
0 & a_{n2} & \ldots &a_{nn}\\
\end{bmatrix}.
$$
This operation ensures that the principal minors of $\bA$ remain unchanged. 
The leading principal minors of $\bE_1 \bA\bE_1^\top$  are identical to those of  $\bA$.

By continuing  this process, we  transform $\bA$ into a diagonal matrix $\bE_n \ldots \bE_1 \bA\bE_1^\top\ldots\bE_n^\top$ with diagonal values identical to those of $\bA$ and positively valued. 
Let $\bE \triangleq \bE_n \ldots \bE_1$, which is an invertible matrix. Apparently, $\bE \bA \bE^\top$ is PD, implying that $\bA$ is PD as well from Theorem~\ref{theorem:quadratic-pd}.
\end{proof}
Alternatively, Sylvester's criterion can be proved using characteristic polynomial; see, for example, \citet{horn2012matrix}.

\section{Existence of  Cholesky Decomposition via  LU Decomposition}
According to  Theorem~\ref{theorem:sylvester-criterion} regarding Sylvester's criterion and  Theorem~\ref{theorem:lu-factorization-without-permutation} stating the existence of LU decomposition without permutation,
a positive definite matrix $\bA$ can be uniquely decomposed as $\bA=\bL\bU_0$, where $\bL$ is a unit lower triangular matrix and $\bU_0$ is an upper triangular matrix.
Additionally, we also have the fact that \textit{the signs of the pivots of a symmetric matrix are the same as the signs of the eigenvalues} \citep{strang1993introduction}:
$$
\text{number of positive pivots = number of positive eigenvalues. \index{Pivot} }
$$
We can express $\bA = \bL\bU_0$ explicitly as follows
$$
\begin{aligned}
\bA = \bL\bU_0 &=
\footnotesize
\begin{bmatrix}
1 & 0 & \ldots & 0 \\
l_{21} & 1 & \ldots & 0\\
\vdots & \vdots & \ddots & \vdots\\
l_{n1} & l_{n2} & \ldots & 1
\end{bmatrix}
\begin{bmatrix}
u_{11} & u_{12} & \ldots & u_{1n} \\
0 & u_{22} & \ldots & u_{2n}\\
\vdots & \vdots & \ddots & \vdots\\
0 & 0 & \ldots & u_{nn}
\end{bmatrix}.\\
\end{aligned}
$$
This indicates that the diagonal entries of $\bU_0$ correspond to the pivots of $\bA$.
Since all the eigenvalues of PD matrices are positive (see Theorem~\ref{theorem:eigen_charac}, which is a consequence of the spectral decomposition), the diagonals of $\bU_0$ are positive. 

By extracting the diagonal elements of $\bU_0$ into a diagonal matrix $\bD$, we can express $\bU_0=\bD\bU$ using the following equation:
$$
\begin{aligned}
\bA = \bL\bU_0 =
\footnotesize
\begin{bmatrix}
1 & 0 & \ldots & 0 \\
l_{21} & 1 & \ldots & 0\\
\vdots & \vdots & \ddots & \vdots\\
l_{n1} & l_{n2} & \ldots & 1
\end{bmatrix}
\begin{bmatrix}
u_{11} & 0 & \ldots & 0 \\
0 & u_{22} & \ldots & 0\\
\vdots & \vdots & \ddots & \vdots\\
0 & 0 & \ldots & u_{nn}
\end{bmatrix}
\begin{bmatrix}
1 & u_{12}/u_{11} & \ldots & u_{1n}/u_{11} \\
0 & 1 & \ldots & u_{2n}/u_{22}\\
\vdots & \vdots & \ddots & \vdots\\
0 & 0 & \ldots & 1
\end{bmatrix}
\normalsize
\triangleq\bL\bD\bU,
\end{aligned}
$$
where $\bU$ is a \textit{unit} upper triangular matrix.
By the uniqueness of the LU decomposition without permutation in Corollary~\ref{corollary:unique-lu-without-permutation} and the symmetry of $\bA$, it follows that $\bU = \bL^\top$ and $\bA = \bL\bD\bL^\top$. Since the diagonals of $\bD$ are positive, we can set $\bR \triangleq \bD^{1/2}\bL^\top$, where $\bD^{1/2}=\diag(\sqrt{u_{11}}, \sqrt{u_{22}}, \ldots, \sqrt{u_{nn}})$ such that $\bA = \bR^\top\bR$ is the Cholesky decomposition of $\bA$, and $\bR$ is upper triangular with positive diagonals.

\subsection{Diagonal Values of the Upper Triangular Matrix}\label{section:cholesky-diagonals}

Let $\bA$ be a PD  matrix and consider $\bA$ as a block matrix $\bA = \scriptsize\begin{bmatrix}
	\bA_{k} & \bA_{12} \\
	\bA_{21} & \bA_{22}
\end{bmatrix}$, where $\bA_{k}\in \real^{k\times k}$. The block LU decomposition of $\bA$ is given by 
$$
\begin{aligned}
\bA &= \begin{bmatrix}
\bA_{k} & \bA_{12} \\
\bA_{21} & \bA_{22}
\end{bmatrix}
=\bL\bU_0=
\begin{bmatrix}
\bL_{k} & \bzero \\
\bL_{21} & \bL_{22}
\end{bmatrix}
\begin{bmatrix}
\bU_{k} & \bU_{12} \\
\bzero & \bU_{22}
\end{bmatrix} 
=\begin{bmatrix}
\bL_{k}\bU_{k} & \bL_{k}\bU_{12} \\
\bL_{21}\bU_{k}  & \bL_{21}\bU_{12}+\bL_{22}\bU_{22}
\end{bmatrix}.
\end{aligned}
$$
Then the $k$-th order leading principal minor (Definition~\ref{definition:leading-principle-minors}), $\Delta_k=\det(\bA[1:k,1:k] ) = \det(\bA_{k})$,  is given by 
$$
\Delta_k = \det(\bA_{k}) = \det(\bL_{k}\bU_{k} ) = \det(\bL_{k} )\det(\bU_{k}).
$$
We notice that $\bL_{k}$ is a unit lower triangular matrix and $\bU_{k}$ is an upper triangular matrix. 
Due to the property that the determinant of a lower triangular matrix (or an upper triangular matrix) is the product of its diagonal entries, we obtain 
$$
\Delta_k = \det(\bU_{k})= u_{11} u_{22}\ldots u_{kk},
$$
i.e., the $k$-th order leading principal minor of $\bA$ is equal to the determinant of the $k\times k$ submatrix of $\bU_0$. 
This quantity is also equivalent to the product of the first $k$ diagonals of $\bD$ ($\bD$ is the matrix from $\bA = \bL\bD\bL^\top$). Let $\bD \triangleq \diag(d_1, d_2, \ldots, d_n)$; therefore, we have 
$$
\Delta_k = d_1 d_2\ldots d_k = \Delta_{k-1}d_k.
$$
This provides us with an alternative representation of $\bD$, namely, the \textbf{squared} diagonal values of $\bR$ ($\bR$ is the Cholesky factor from $\bA=\bR^\top\bR$), which can be expressed as
$$
\bD = \diag\left(\Delta_1, \frac{\Delta_2}{\Delta_1}, \ldots, \frac{\Delta_n}{\Delta_{n-1}}\right),
$$
where $\Delta_k$ is the $k$-th  order leading principal minor of $\bA$, for all $k\in \{1,2,\ldots, n\}$. That is, the diagonal values of $\bR$ are given by 
$$
\diag\left(\sqrt{\Delta_1}, \sqrt{\frac{\Delta_2}{\Delta_1}}, \ldots, \sqrt{\frac{\Delta_n}{\Delta_{n-1}}}\right).
$$

\subsection{Block Cholesky Decomposition}
Following from the previous section, suppose $\bA\in\real^{n\times n}$ is a PD matrix and take $\bA$ as a block matrix $\bA = \scriptsize\begin{bmatrix}
\bA_k & \bA_{12} \\
\bA_{21} & \bA_{22}
\end{bmatrix}$, where $\bA_k\in \real^{k\times k}$. Its block LU decomposition is given by 
$$
\begin{aligned}
\bA &= \begin{bmatrix}
	\bA_k & \bA_{12} \\
	\bA_{21} & \bA_{22}
\end{bmatrix}
=\bL\bU_0=
\begin{bmatrix}
	\bL_k & \bzero \\
	\bL_{21} & \bL_{22}
\end{bmatrix}
\begin{bmatrix}
	\bU_k & \bU_{12} \\
	\bzero & \bU_{22}
\end{bmatrix} 
=\begin{bmatrix}
	\bL_k\bU_k & \bL_{k}\bU_{12} \\
	\bL_{21}\bU_{k}  & \bL_{21}\bU_{12}+\bL_{22}\bU_{22}
\end{bmatrix},\\
\end{aligned}
$$
where the $k$-th order leading principal submatrix $\bA_k$ of $\bA$ also possesses an LU decomposition, given by $\bA_k=\bL_k\bU_k$. 
It is evident that the Cholesky decomposition of an $n\times n$ matrix includes $n-1$ additional Cholesky decompositions within it: $\bA_k = \bR_k^\top\bR_k$, for all $k\in \{1,2,\ldots, n-1\}$. This is particularly true because any leading principal submatrix $\bA_k$ of a positive definite matrix $\bA$ is also positive definite. 
To illustrate, for a positive definite matrix $\bA_{k+1} \in \real^{(k+1)\times(k+1)}$, and any nonzero vector $\bx_k\in\real^k$ appended by a zero element $\bx_{k+1} = \footnotesize\begin{bmatrix}
	\bx_k\\
	0
\end{bmatrix}$, it follows that
$$
\bx_k^\top\bA_k\bx_k = \bx_{k+1}^\top\bA_{k+1}\bx_{k+1} >0,
$$
which implies $\bA_k$ is positive definite. If we start with $\bA\in \real^{n\times n}$, we will recursively conclude that $\bA_{n-1}$ is PD, $\bA_{n-2}$ is PD, $\ldots$ (i.e., consistent with  Sylvester's criterion). And each of these matrices admits a Cholesky decomposition.

\section{Existence of  Cholesky Decomposition via Induction}
In the previous section, we established the existence of the Cholesky decomposition using LU decomposition without permutation. Following the proof of the LU decomposition in Section~\ref{section:exist-lu-without-perm}, it becomes evident that the existence of the Cholesky decomposition can be directly inferred through induction as well.
\begin{proof}[{of Theorem~\ref{theorem:cholesky-factor-exist}: Existence of Cholesky Decomposition by Induction\index{Induction}}]
We will prove by induction that every $n\times n$ positive definite matrix $\bA$ can be decomposed as $\bA=\bR^\top\bR$. The $1\times 1$ base case is trivial by setting $R=\sqrt{A}$, thus, $A=R^2$. 

Suppose that for any $k\times k$ PD matrix $\bA_k$, it has a Cholesky decomposition. If we can prove any $(k+1)\times(k+1)$ PD matrix $\bA_{k+1}$ can also be factored as this Cholesky decomposition, then we complete the proof.

Consider any $(k+1)\times(k+1)$ PD matrix $\bA_{k+1}$, which can be expressed as
$
\bA_{k+1} = 
\scriptsize\begin{bmatrix}
\bA_k & \bb \\
\bb^\top & d
\end{bmatrix}.
$
We note that $\bA_k$ is PD. 
Based on the inductive hypothesis, $\bA_k$ can be decomposed as $\bA_k = \bR_k^\top\bR_k$ using the Cholesky decomposition.
We can construct the upper triangular matrix 
$$
\bR_{k+1}
\triangleq
\begin{bmatrix}
\bR_k & \br\\
0 & s
\end{bmatrix}
\qquad\text{such that}\qquad 
\bR_{k+1}^\top\bR_{k+1} = 
\begin{bmatrix}
\bR_k^\top\bR_k & \bR_k^\top \br\\
\br^\top \bR_k & \br^\top\br+s^2
\end{bmatrix}.
$$
To complete the proof, we need to show that  $\bR_{k+1}^\top \bR_{k+1} = \bA_{k+1}$ is the Cholesky decomposition of $\bA_{k+1}$ (which requires the value $s$ to be positive). That is, we need to prove
$$
\begin{aligned}
\bb &= \bR_k^\top \br
\qquad \text{and}\qquad 
d = \br^\top\br+s^2.
\end{aligned}
$$
As $\bR_k$ is nonsingular, there exists a unique solution  for $\br$ and $s$ such that
$$
\begin{aligned}
\br &= \bR_k^{-\top}\bb
\qquad \text{and}\qquad 
s = \sqrt{d - \br^\top\br} = \sqrt{d - \bb^\top\bA_k^{-1}\bb},
\end{aligned}
$$
where we assume $s$ is nonnegative. However, we need to further show  that $s$ is not only nonnegative but also positive. Since $\bA_k$ is PD, from Sylvester's criterion, and the determinant of a block matrix (Remark~\ref{remark:determinant-intermezzo}). We have
$$
\det(\bA_{k+1}) = \det(\bA_k)\det(d- \bb^\top\bA_k^{-1}\bb) =  \det(\bA_k)(d- \bb^\top\bA_k^{-1}\bb)>0.
$$
Since $ \det(\bA_k)>0$, we then obtain that $(d- \bb^\top\bA_k^{-1}\bb)>0$, and this implies $s>0$.
We complete the proof.
\end{proof}

\section{Uniqueness of  Cholesky Decomposition}
The uniqueness of the Cholesky decomposition can be an immediate consequence of the uniqueness of the LU decomposition without permutation.  Alternatively, a rigorous proof is provided as follows.
\begin{corollary}[Uniqueness of Cholesky Decomposition\index{Uniqueness}]\label{corollary:unique-cholesky-main}
	The Cholesky decomposition $\bA=\bR^\top\bR$ for any positive definite matrix $\bA\in \real^{n\times n}$ is unique.
\end{corollary}
\begin{proof}[of Corollary~\ref{corollary:unique-cholesky-main}]
Suppose the Cholesky decomposition is not unique. Then, there exist two different decompositions such that $\bA=\bR_1^\top\bR_1 = \bR_2^\top\bR_2$, which implies
$$
\bR_1\bR_2^{-1}= \bR_1^{-\top} \bR_2^\top.
$$
Since the inverse of an upper triangular matrix is also an upper triangular matrix, and the product of two upper triangular matrices is also an upper triangular matrix,~\footnote{Same for lower triangular matrices: the inverse of a lower triangular matrix is also a lower triangular matrix, and the product of two lower triangular matrices is also a lower triangular matrix.} we realize that the left-hand side of the above equation is an upper triangular matrix, and its right-hand side  is a lower triangular matrix. This implies $\bR_1\bR_2^{-1}= \bR_1^{-\top} \bR_2^\top$ must be a diagonal matrix, and $\bR_1^{-\top} \bR_2^\top= (\bR_1^{-\top} \bR_2^\top)^\top = \bR_2\bR_1^{-1}$.
Let $\bLambda \triangleq \bR_1\bR_2^{-1}= \bR_2\bR_1^{-1}$ be the diagonal matrix. We notice that the diagonal value of $\bLambda$ is the product of the corresponding diagonal values of $\bR_1$ and $\bR_2^{-1}$ (or $\bR_2$ and $\bR_1^{-1}$). That is, for 
$$
\bR_1=
\footnotesize
\begin{bmatrix}
r_{11} & r_{12} & \ldots & r_{1n} \\
0 & r_{22} & \ldots & r_{2n}\\
\vdots & \vdots & \ddots & \vdots\\
0 & 0 & \ldots & r_{nn}
\end{bmatrix}
\gap\text{and}\gap
\normalsize
\bR_2=
\footnotesize
\begin{bmatrix}
s_{11} & s_{12} & \ldots & s_{1n} \\
0 & s_{22} & \ldots & s_{2n}\\
\vdots & \vdots & \ddots & \vdots\\
0 & 0 & \ldots & s_{nn}
\end{bmatrix},
$$
we have,
$$
\begin{aligned}
\bR_1\bR_2^{-1}=
\footnotesize
\begin{bmatrix}
\frac{r_{11}}{s_{11}} & 0 & \ldots & 0 \\
0 & \frac{r_{22}}{s_{22}} & \ldots & 0\\
\vdots & \vdots & \ddots & \vdots\\
0 & 0 & \ldots & \frac{r_{nn}}{s_{nn}}
\end{bmatrix}
\normalsize
=
\footnotesize
\begin{bmatrix}
\frac{s_{11}}{r_{11}} & 0 & \ldots & 0 \\
0 & \frac{s_{22}}{r_{22}} & \ldots & 0\\
\vdots & \vdots & \ddots & \vdots\\
0 & 0 & \ldots & \frac{s_{nn}}{r_{nn}}
\end{bmatrix}
\normalsize
=\bR_2\bR_1^{-1}.
\end{aligned}
$$ 
Since both $\bR_1$ and $\bR_2$ have positive diagonals, we can deduce that  $r_{11}=s_{11}, r_{22}=s_{22}, \ldots, r_{nn}=s_{nn}$. And $\bLambda = \bR_1\bR_2^{-1}= \bR_2\bR_1^{-1}  =\bI$.
That is, $\bR_1=\bR_2$, which results in a contradiction. Therefore, the Cholesky decomposition is  unique. 
\end{proof}
Based on the preceding discussion, when the diagonal elements of $\bR_1$ and $\bR_2$ are not restricted to positive values, then $r_{ii} = \pm s_{ii}$ for $i\in\{1,2,\ldots, n\}$. In this case, the factorization $\bA=\bR^\top\bR$ is \textit{not unique}.

\section{Computing  Cholesky Decomposition Recursively}\label{section:compute-cholesky}
Similar to computing the LU decomposition,
to compute the Cholesky decomposition, we start by writing out the equality $\bA=\bR^\top\bR$:
$$
\footnotesize
\setlength{\arraycolsep}{2pt}
\begin{aligned}
\bA=
\begin{bmatrix}
	a_{11} & \bA_{1,2:n} \\
	\bA_{2:n,1} & \bA_{2:n,2:n}
\end{bmatrix}
&=
\begin{bmatrix}
	r_{11} & 0 \\
	\bR_{1,2:n}^\top & \bR_{2:n,2:n}^\top
\end{bmatrix}
\begin{bmatrix}
	r_{11} & \bR_{1,2:n} \\
	0 & \bR_{2:n,2:n}
\end{bmatrix}
=
\begin{bmatrix}
r_{11}^2 & r_{11}\bR_{1,2:n} \\
r_{11}\bR_{1,2:n}^\top & \bR_{1,2:n}^\top\bR_{1,2:n} + \bR_{2:n,2:n}^\top\bR_{2:n,2:n}
\end{bmatrix},  
\end{aligned}
$$
which allows us to determine the first row of $\bR$ by 
$$
r_{11} = \sqrt{a_{11}}
\gap \text{and}\gap 
 \bR_{1,2:n} = \frac{1}{r_{11}}\bA_{1,2:n}.
$$
Let $\bA_2=\bR_{2:n,2:n}^\top\bR_{2:n,2:n}$. The equality $\bA_{2:n,2:n} = \bR_{1,2:n}^\top\bR_{1,2:n} + \bR_{2:n,2:n}^\top\bR_{2:n,2:n}$ and the symmetry of $\bA$ indicate
$$
\begin{aligned}
\bA_2=\bR_{2:n,2:n}^\top\bR_{2:n,2:n} &= \bA_{2:n,2:n} - \bR_{1,2:n}^\top\bR_{1,2:n} 
= \bA_{2:n,2:n} - \frac{1}{a_{11}} \bA_{2:n,1}\bA_{1,2:n},
\end{aligned}
$$
where $\bA_2$ is the Schur complement of $a_{11}$ in $\bA$ and has a  size of $(n-1)\times (n-1)$. To obtain $\bR_{2:n,2:n}$, we must compute the Cholesky decomposition of the matrix $\bA_2$ of shape $(n-1)\times (n-1)$. Again, this is a recursive algorithm and is outlined in Algorithm~\ref{alg:compute-choklesky}.

\begin{algorithm}[H] 
\caption{Cholesky Decomposition via Recursive Algorithm} 
\label{alg:compute-choklesky} 
\begin{algorithmic}[1] 
\Require 
Positive definite matrix $\bA$ with size $n\times n$; 
\State Calculate first row of $\bR$ by $r_{11} \leftarrow \sqrt{a_{11}}, \bR_{1,2:n} \leftarrow \frac{1}{r_{11}}\bA_{1,2:n}$; \Comment{$n$ flops}
\State Compute the Cholesky decomposition of the $(n-1)\times (n-1)$ matrix
$$
\bA_2\leftarrow\bR_{2:n,2:n}^\top\bR_{2:n,2:n}=\bA_{2:n,2:n} - \frac{1}{a_{11}} \bA_{2:n,1}\bA_{1,2:n};
$$
\Comment{$n^2-n$ flops}

\end{algorithmic} 
\end{algorithm}


\begin{theorem}[Algorithm Complexity: Cholesky Recursively]\label{theorem:cholesky-complexity}
Algorithm~\ref{alg:compute-choklesky} requires $\sim(1/3)n^3$ flops to compute the Cholesky decomposition of an $n\times n$ positive definite matrix.
\end{theorem}

\begin{proof}[of Theorem~\ref{theorem:cholesky-complexity}]
Step 1 involves 1 square root and $(n-1)$ divisions, which take $n$ flops totally.  

For step 2, note that $\frac{1}{a_{11}} \bA_{2:n,1}\bA_{1,2:n} = \left(\frac{1}{\sqrt{a_{11}}} \bA_{2:n,1}\right)\left(\frac{1}{\sqrt{a_{11}}}\bA_{1,2:n}\right) = \bR_{1,2:n}^\top \bR_{1,2:n}$. If we calculate the complexity directly from the equation in step 2, we will get the same complexity as the LU decomposition. But the symmetry of $\bR_{1,2:n}^\top \bR_{1,2:n}$ can be adopted, the complexity of $\bR_{1,2:n}^\top \bR_{1,2:n}$ reduces from $(n-1)\times(n-1)$ multiplications to $1+2+\ldots+(n-1)=\frac{n^2-n}{2}$ multiplications, which is almost half of the original complexity. The cost of matrix division reduces from $(n-1)\times(n-1)$ to $1+2+\ldots+(n-1)=\frac{n^2-n}{2}$ as well. So the total cost for step 2 is $n^2-n$ flops.

Let $f(n) = n^2-n + n = n^2$, the total complexity is 
$$
\mathrm{cost} = f(n)+f(n-1)+\ldots +f(1).
$$
Simple calculations  show that  the total complexity for all the recursive steps is $\frac{2n^3+3n^2+n}{6}$ flops, which is $(1/3)n^3$ flops if we keep only the leading term.
\end{proof}

The Cholesky decomposition computation mentioned above has an important application in testing the positive definiteness of a symmetric matrix. To perform the test, one can apply the algorithm mentioned above and declare the matrix as positive definite if the algorithm completes without encountering any negative or zero pivots (as described in step 1 above). Otherwise, if the algorithm encounters such pivots, the matrix is deemed not positive definite.


To end up this section, we provide the full pseudo code for Algorithm~\ref{alg:compute-choklesky} as shown in Algorithm~\ref{alg:compute-choklesky11} (compare the two algorithms).  
\begin{algorithm}[H] 
\caption{Cholesky Decomposition via Recursive Algorithm: Full Pseudo Code} 
\label{alg:compute-choklesky11} 
\begin{algorithmic}[1] 
\Require 
Positive definite matrix $\bA$ with size $n\times n$; 
\For{$k=1$ to $n$} \Comment{compute the $k$-th row of $\bR$}
\State $r_{kk} \leftarrow \sqrt{a_{kk}}$; \Comment{first element of $k$-th row, 1 flop}
\State $\bR_{k,k+1:n} \leftarrow \frac{1}{r_{kk}} \bA_{k,k+1:n}$; \Comment{the rest elements of $k$-th row, $n-k$ flops}
\State $\bA_{k+1:n,k+1:n} \leftarrow \bA_{k+1:n,k+1:n} - \bR_{k,k+1:n}^\top\bR_{k,k+1:n}$; \Comment{$2(1+2+\ldots+(n-k))$ flops}
\EndFor


\end{algorithmic} 
\end{algorithm}
\section{Computing  Cholesky Decomposition Element-Wise}

It is also common to compute the Cholesky decomposition via element-level equations derived directly from solving the matrix equation $\bA=\bR^\top\bR$. We notice that the entry $(i,j)$ of $\bA$ is $a_{ij} = \bR_{:,i}^\top \bR_{:,j} = \sum_{k=1}^{i} r_{ki}r_{kj}$ if $i<j$. This further implies, if $i<j$, we have
$$
\begin{aligned}
a_{ij} &= \bR_{:,i}^\top \bR_{:,j} = \sum_{k=1}^{i} r_{ki}r_{kj} 
= \sum_{k=1}^{i-1} r_{ki}r_{kj} + r_{ii}r_{ij}
\implies
r_{ij} = (a_{ij} - \sum_{k=1}^{i-1} r_{ki}r_{kj})/r_{ii},
\gap 
\text{if }i<j.
\end{aligned}
$$
On the other hand, if $i=j$, we have 
$$
\begin{aligned}
a_{jj} &= \sum_{k=1}^{j} r_{kj}^2=\sum_{k=1}^{j-1} r_{kj}^2 + r_{jj}^2
&\implies
r_{jj} = \sqrt{a_{jj} - \sum_{k=1}^{j-1} r_{kj}^2}.
\end{aligned}
$$
If we equate the elements of $\bR$ by taking a column at a time and start with $r_{11} = \sqrt{a_{11}}$, the element-level algorithm is formulated in Algorithm~\ref{alg:compute-choklesky-element-level}.

\begin{algorithm}[H] 
\caption{Cholesky Decomposition Element-Wise: $\bA=\bR^\top\bR$} 
\label{alg:compute-choklesky-element-level} 
\begin{algorithmic}[1] 
\Require 
Positive definite matrix $\bA$ with size $n\times n$; 
\State Calculate first element of $\bR$ by $r_{11} \leftarrow \sqrt{a_{11}}$; 
\For{$j=1$ to $n$} \Comment{Compute the $j$-th column of $\bR$}
\For{$i=1$ to $j-1$} 
\State $r_{ij} \leftarrow (a_{ij} - \sum_{k=1}^{i-1} r_{ki}r_{kj})/r_{ii}$, since $i<j$;
\EndFor
\State $r_{jj} \leftarrow \sqrt{a_{jj}- \sum_{k=1}^{j-1}r_{kj}^2}$;
\EndFor
\State Output $\bA=\bR^\top\bR$.
\end{algorithmic} 
\end{algorithm}

\begin{theorem}[Algorithm Complexity: Cholesky Element-Wise]\label{theorem:cholesky-complexity-element}
Algorithm~\ref{alg:compute-choklesky-element-level} requires $\sim(1/3)n^3$ flops to compute the Cholesky decomposition of an $n\times n$ positive definite matrix.
\footnote{The Cholesky method requires half of the flops required by  Doolittle's method for LU decomposition in Theorem~\ref{theorem:lu-complexity-doolittle}.}
\end{theorem}

\begin{proof}[of Theorem~\ref{theorem:cholesky-complexity-element}]
For step 4 in Algorithm~\ref{alg:compute-choklesky-element-level}, that is, for each $(j, i)$ loop, step 4 involves $(i-1)$ multiplications, $(i-2)$ additions, $1$ subtraction, and 1 division, which is $2i-1$ flops. Let $f(k)=2k-1$, the total flops required from step 4 for each loop $j$ is given by 
$$
f(1)+f(2)+\ldots +f(j-1) = j^2-2j+1 \,\,\,\,\mathrm{flops}.
$$	
Further, for each loop $j$, step 6 requires $j-1$ multiplications, $j-2$ additions, $1$ subtraction, and $1$ square root. That is, step 6 involves $2j-1$ flops for each loop $j$. Combining the flops involved in step 4 and step 6, it requires $j^2$ flops for each loop $j$ totally. Let $g(k)=k^2$, the total complexity is thus given by 
$$
g(1)+g(2)+\ldots +g(n) = \frac{2n^3+3n^2+n}{6} \,\,\,\,\mathrm{flops},
$$	
which is $(1/3)n^3$ flops if we keep only the leading term.
\end{proof}
The complexity of Algorithm~\ref{alg:compute-choklesky-element-level} is equivalent to that of Algorithm~\ref{alg:compute-choklesky}. 
Indeed, the ideas behind both algorithms are similar.

On the other hand, Algorithm~\ref{alg:compute-choklesky-element-level}  can be modified to compute the Cholesky decomposition in the form $\bA=\bL\bD\bL^\top$, where $\bL$ is unit lower triangular and $\bD$ is diagonal, as outlined in Algorithm~\ref{alg:compute-choklesky-_ldl}, whose Step 3 and Step 5 are derived from (since $l_{ii}=1, \forall i\in\{1,2,\ldots,n\}$):
$$
\begin{aligned}
a_{jj}&=\sum_{k=1}^{j-1}d_{kk} l_{jk}^2 + d_{jj};\\
a_{ij}&= d_{jj} l_{ij}+ \sum_{k=1}^{j-1} d_{kk} l_{ik}l_{jk}, \gap \text{if }i>j.
\end{aligned}
$$
\begin{exercise}
	Derive the complexity of Algorithm~\ref{alg:compute-choklesky-_ldl}.
\end{exercise}
This form of Cholesky decomposition is useful for determining the condition number (see Section~\ref{section:qr_condition} and Appendix~\ref{appendix:condition_number}) of a PD matrix. 
In essence, the condition number of a function measures the sensitivity of the output value to small changes in the input; a smaller condition number indicates better numerical stability. 
For positive definite linear systems, the condition number is defined as the ratio of the largest eigenvalue to the smallest eigenvalue.
The condition number of a positive definite matrix is lower bounded by the diagonal matrix in the Cholesky decomposition (see Problem~\ref{problem:cond_pd}):
\begin{equation}\label{equation:cond_pd_ineq}
\cond(\bA) \geq \cond(\bD).
\end{equation}
This can be proven by showing that $\lambda_{\max}\geq d_{\max}$ and $\lambda_{\min}\leq d_{\min}$, where $\lambda_{\max}$ and $\lambda_{\min}$ are the largest and smallest eigenvalue of $\bA$, and $d_{\max}$ and $d_{\min}$ are the largest and smallest diagonals of $\bD$.
Therefore, this form of the Cholesky decomposition can be utilized  to modify Newton's method; see \S~\ref{section:app_cho_md_newton}.

\begin{algorithm}[h] 
\caption{Cholesky Decomposition Element-Wise: $\bA=\bL\bD\bL^\top$}  
\label{alg:compute-choklesky-_ldl} 
\begin{algorithmic}[1] 
\Require 
Positive definite matrix $\bA$ with size $n\times n$; 
\For{$j=1$ to $n$} \Comment{Compute the $j$-th column of $\bL$}
\State $l_{jj}\leftarrow1$;
\State $c_{jj}\leftarrow a_{jj}-\sum_{k=1}^{j-1}d_{kk} l_{jk}^2$;
\State $d_{jj}\leftarrow c_{jj}$
\For{$i=j+1$ to $n$} 
\State $c_{ij}\leftarrow a_{ij}-\sum_{k=1}^{j-1}d_{kk} l_{ik}l_{jk}$, since $i>j$;
\State $l_{ij}\leftarrow \frac{c_{ij}}{d_{jj}}$;
\EndFor
\EndFor
\State Output $\bA=\bL\bD\bL^\top$, where $\bD=\diag(d_{11}, d_{22},\ldots,d_{nn})$.
\end{algorithmic} 
\end{algorithm}

\section{More Properties of Positive Definite Matrices}
\begin{lemma}[PD Properties\index{Positive definite}]\label{lemma:pd-more-properties}
Given a positive definite matrix $\bA\in \real^{n\times n}$, the following properties hold:
\begin{enumerate}
\item \textit{Principal minors criterion.} Any principal minors are positive (see Definition~\ref{definition:principle-minors}, not necessarily  the leading principal minors);
\item Suppose the diagonal values of $\bA$ are $a_{ii}$, for all $i\in \{1,2,\ldots, n\}$, then $\det(\bA) \leq \prod_{i=1}^{n}a_{ii}$. 
Equality holds when $\bA$ is a diagonal matrix.
\end{enumerate}
\end{lemma}
\begin{proof}[of Lemma~\ref{lemma:pd-more-properties}]
For 1). Similar to the permutation matrix (Definition~\ref{definition:permutation-matrix}), we can define a selection matrix. For a row selection matrix $\bP$ (Definition~\ref{definition:selection-matrix}), if we select a row of the matrix $\bA$ by $\bP\bA$, the corresponding diagonal of $\bP$ will be 1, and 0 otherwise. In this sense,  $\bP^\top=\bP$ acts as a column selection matrix such that $\bA\bP^\top$ selects the corresponding columns. Then  any $k\times k$ submatrix of $\bA$ can be obtained by $\bP\bA\bP^\top$. For any vector $\bx\in \real^n$, where not all the corresponding $k$ entries are zero (such that $\bx^\top\bP$ is nonzero), we have 
$$
\bx^\top \bP\bA\bP^\top \bx >0.
$$
Since $(n-k)$ rows of $\bP$ are zero and the corresponding $(n-k)$ columns of $\bP$ are zero as well, these rows and columns will not count any value for the above equation, and they can be removed.
This implies that the $k\times k$ submatrix $\bA_k \in \real^{k\times k}$ of $\bA$ and its corresponding $k\times k$ sub-selection matrix $\bP_k \in \real^{k\times k}$ (which is an identity matrix here) satisfy
$$
\bx_k\bP_k\bA_k\bP_k^\top \bx_k >0,
$$
where $\bx_k \in \real^{k}$. Since the vector $\bx$ is any nonzero vector in $\real^n$ (such that $\bx^\top\bP$ is nonzero), the $\bx_k$ is thus also any nonzero vector in $\real^k$.
This results in that $\bA_k$ is PD.

For 2). Consider the LU decomposition of $\bA$ as
$\bA = \bL\bU_0$, where 
$$
\begin{aligned}
\bA = \bL\bU_0 &=
\footnotesize
\begin{bmatrix}
1 & 0 & \ldots & 0 \\
l_{21} & 1 & \ldots & 0\\
\vdots & \vdots & \ddots & \vdots\\
l_{n1} & l_{n2} & \ldots & 1
\end{bmatrix}
\begin{bmatrix}
u_{11} & u_{12} & \ldots & u_{1n} \\
0 & u_{22} & \ldots & u_{2n}\\
\vdots & \vdots & \ddots & \vdots\\
0 & 0 & \ldots & u_{nn}
\end{bmatrix}\\
&=
\footnotesize
\begin{bmatrix}
1 & 0 & \ldots & 0 \\
l_{21} & 1 & \ldots & 0\\
\vdots & \vdots & \ddots & \vdots\\
l_{n1} & l_{n2} & \ldots & 1
\end{bmatrix}
\begin{bmatrix}
u_{11} & 0 & \ldots & 0 \\
0 & u_{22} & \ldots & 0\\
\vdots & \vdots & \ddots & \vdots\\
0 & 0 & \ldots & u_{nn}
\end{bmatrix}
\begin{bmatrix}
1 & u_{12}/u_{11} & \ldots & u_{1n}/u_{11} \\
0 & 1 & \ldots & u_{2n}/u_{22}\\
\vdots & \vdots & \ddots & \vdots\\
0 & 0 & \ldots & 1
\end{bmatrix}\\
&=
\footnotesize
\begin{bmatrix}
1 & 0 & \ldots & 0 \\
l_{21} & 1 & \ldots & 0\\
\vdots & \vdots & \ddots & \vdots\\
l_{n1} & l_{n2} & \ldots & 1
\end{bmatrix}
\begin{bmatrix}
u_{11} & 0 & \ldots & 0 \\
0 & u_{22} & \ldots & 0\\
\vdots & \vdots & \ddots & \vdots\\
0 & 0 & \ldots & u_{nn}
\end{bmatrix}
\begin{bmatrix}
1 & l_{21} & \ldots & l_{n1} \\
0 & 1 & \ldots & l_{n2}\\
\vdots & \vdots & \ddots & \vdots\\
0 & 0 & \ldots & 1
\end{bmatrix}
\normalsize
=\bL\bD\bL^\top.
\end{aligned}
$$
That is, for $j<i$, we have $u_{ji}/u_{jj} = l_{ij}$, which leads to $u_{ji} = l_{ij}u_{jj}$.
As discussed in Section~\ref{section:cholesky-diagonals}  that $\det(\bA) = \det(\bL)\det(\bU_0) = \prod_{i=1}^{n}u_{ii}\leq \prod_{i=1}^{n}a_{ii}$, where the last inequality arises from the symmetry of $\bA$ such that 
$$
a_{ii} = \left(\sum_{j=1}^{i-1} l_{ij}u_{ji}\right) +u_{ii} = \left(\sum_{j=1}^{i-1} l_{ij}^2\cdot u_{jj}\right) +u_{ii} \geq u_{ii}.
$$ 
And the equality can be achieved when $\bA$ is a diagonal matrix.
\end{proof}

As a matter of fact, the principal minors criterion states that a symmetric matrix is positive definite if and only if all its principal minors are positive. The reverse implication is left as an exercise.
Note that while the principal minors criterion is useful for detecting positive definiteness of a matrix, it is not suitable for detecting positive semidefiniteness.

\paragraph{Quadratic form of positive definite matrices.}
\begin{figure}[h]
	\centering  
	\vspace{-0.35cm} 
	\subfigtopskip=2pt 
	\subfigbottomskip=2pt 
	\subfigcapskip=-5pt 
	\subfigure[Positive definite matrix: $\bA = \begin{bmatrix}
		200 & 0 \\ 0 & 200
	\end{bmatrix}$.]{\label{fig:quadratic_PD}
		\includegraphics[width=0.485\linewidth]{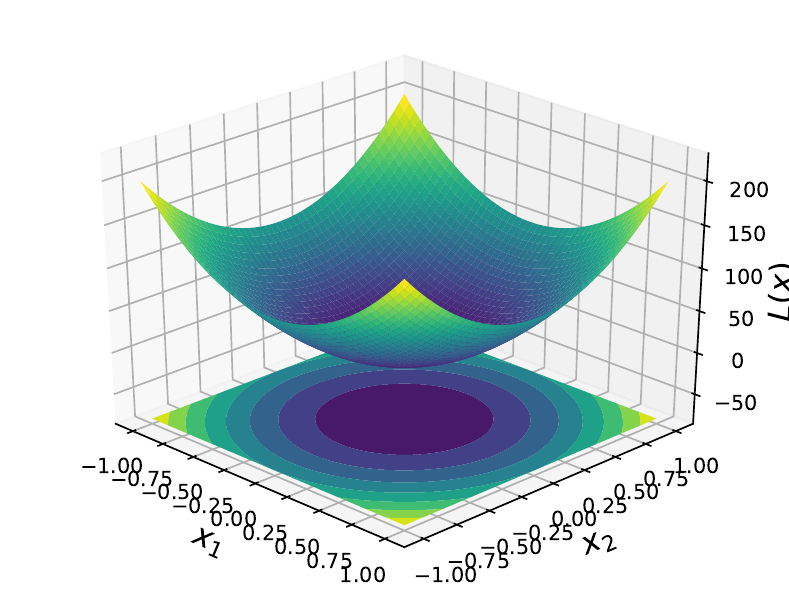}}
	\subfigure[Negative definite matrix: $\bA = \begin{bmatrix}
		-200 & 0 \\ 0 & -200
	\end{bmatrix}$.]{\label{fig:quadratic_ND}
		\includegraphics[width=0.485\linewidth]{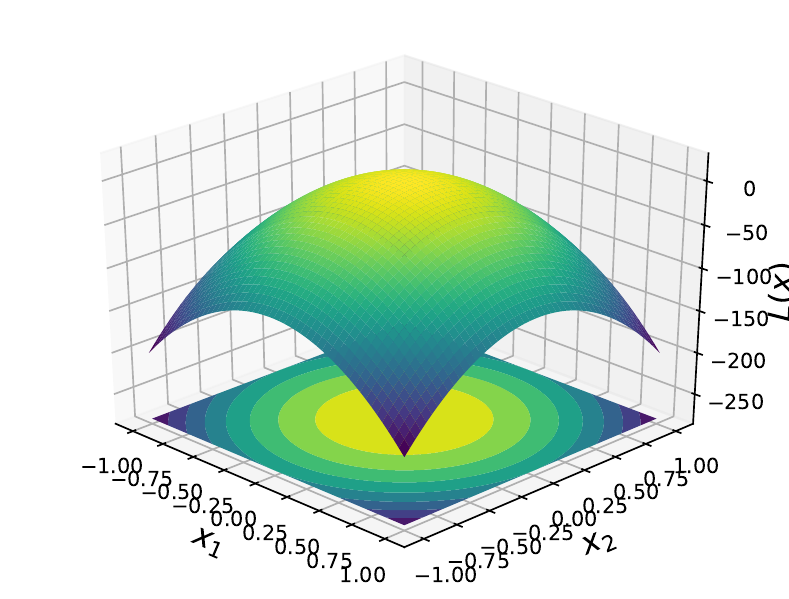}}
	\subfigure[Semidefinite matrix: $\bA = \begin{bmatrix}
		200 & 0 \\ 0 & 0
	\end{bmatrix}$. A line runs through the bottom of the valley is the set of solutions.]{\label{fig:quadratic_singular}
		\includegraphics[width=0.485\linewidth]{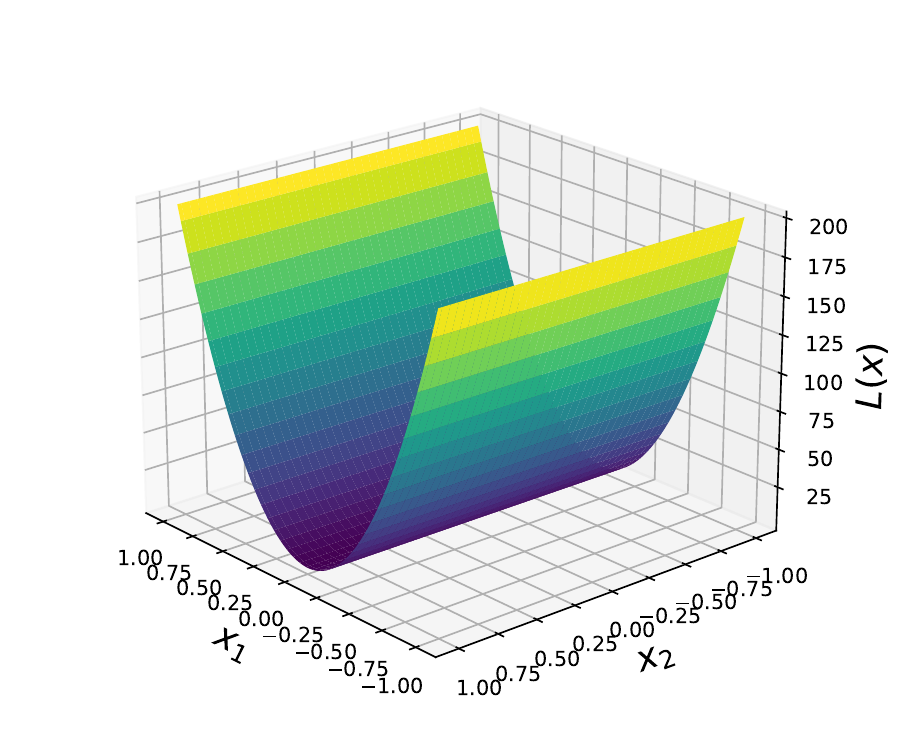}}
	\subfigure[Indefinte matrix: $\bA = \begin{bmatrix}
		200 & 0 \\ 0 & -200
	\end{bmatrix}$.]{\label{fig:quadratic_saddle}
		\includegraphics[width=0.485\linewidth]{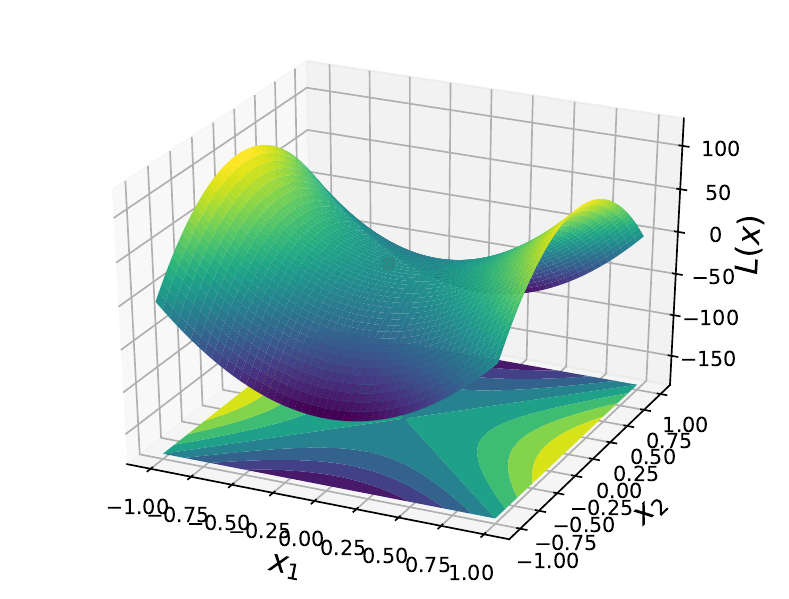}}
\caption{Loss surfaces for different quadratic forms, providing the surface plots and contour plots (\textcolor{mylightbluetext}{blue}=low,
\textcolor{mydarkyellow}{yellow}=high), where the upper graphs are the surface plots, and the lower ones are their projection (i.e., contours).}
\label{fig:different_quadratics}
\end{figure}

We  further  discuss  linear systems with different types of matrices, the convex quadratic:
\begin{equation}\label{equation:quadratic-form-general-form}
	L(\bx) = \frac{1}{2} \bx^\top \bA \bx - \bb^\top \bx + c, \gap \bx\in \real^d,
\end{equation}
where $\bA\in \real^{d\times d}$, $\bb \in \real^d$, and $c$ is a scalar constant. Though the quadratic form in Equation~\eqref{equation:quadratic-form-general-form} is an extremely simple model, it is rich enough to approximate many other functions, e.g., the Fisher information matrix \citep{amari1998natural}, and it captures key features of pathological curvature. The gradient of $L(\bx)$ at point $\bx$ is given by 
$$
	\nabla L(\bx) = \frac{1}{2} (\bA^\top +\bA) \bx - \bb.
$$
The unique minimum of the function is the solution to the linear system $\frac{1}{2} (\bA^\top +\bA) \bx=  \bb $:
$$
	\bx_* = 2(\bA^\top +\bA)^{-1}\bb.
$$
If $\bA$ is symmetric, the equation reduces to 
$
\nabla L(\bx) = \bA \bx - \bb.
$
Then the unique minimum of the function is the solution of the linear system $\bA\bx=\bb$, where $\bA$ and $\bb$ are known matrix or vector, and $\bx$ is an unknown vector; and the optimal point of $\bx$ is thus given by 
$$
\bx_*  = \bA^{-1}\bb.
$$
For different types of matrix $\bA$, the loss surface of $L(\bx)$ will be different as shown in Figure~\ref{fig:different_quadratics}. When $\bA$ is positive definite, the surface is a \textit{convex bowl}; when $\bA$ is negative definite, on the contrary, the surface is a \textit{concave bowl}. $\bA$ also could be singular, in which case $\bA\bx-\bb=\bzero$ has more than one solution, and the set of solutions is a line (in the two-dimensional case) or a hyperplane (in the high-dimensional case).
This situation is similar to the case of a semidefinite quadratic form, as shown in Figure~\ref{fig:quadratic_singular}.
Moreover, $\bA$ could be none of the above, then there exists a \textit{saddle point} (see Figure~\ref{fig:quadratic_saddle}), where the gradient descent may fail \citep{lu2022gradient}. In such senses, other methods, e.g., perturbed gradient descent \citep{jin2017escape}, can be applied to escape  saddle points. 

\paragraph{Diagonally dominant matrices.}
A specific form of diagonally dominant matrices is a significant subset of positive semidefinite matrices.
\begin{definition}[Diagonally Dominant Matrices]
Given a symmetric matrix $\bA\in\real^{n\times n}$, then $\bA$ is called diagonally dominant if 
$$
\abs{a_{ii}} \geq \sum_{j\neq i} \abs{a_{ij}}, \qquad \forall i\in\{1,2,\ldots, n\};
$$
and $\bA$ is called strictly diagonally dominant if 
$$
\abs{a_{ii}} > \sum_{j\neq i} \abs{a_{ij}}, \qquad \forall i\in\{1,2,\ldots, n\}.
$$
\end{definition}

We  now show that \textit{diagonally dominant matrices} with nonnegative diagonal elements are positive semidefinite and that \textit{strictly diagonally dominant matrices} with positive diagonal elements are positive definite. 

\begin{theorem}[Positive Definiteness of Diagonally Dominant Matrices]\label{theorem:pd_diag_domi}
Given a symmetric matrix  $\bA\in\real^{n\times n}$, 
\begin{enumerate}
\item If $\bA$ is diagonally dominant with nonnegative diagonals, then $\bA$ is positive semidefinite; 
\item If $\bA$ is strictly diagonally dominant with positive diagonals, then $\bA$ is positive definite.
\end{enumerate}

\end{theorem}
\begin{proof}[of Theorem~\ref{theorem:pd_diag_domi}]
For 1). Suppose $\bA$ is not positive semidefinite with a negative eigenvalue $\lambda$ associated with an eigenvector $\bv$ such that $\bA\bv=\lambda\bv$. 
Let $v_i$ be the element of $\bv$ with largest magnitude.
Consider the $i$-th element of $(\bA-\lambda\bI)\bv=\bzero$, we have
$$
\abs{a_{ii} - \lambda}\cdot \abs{v_i}
=
\abs{\sum_{j\neq i} a_{ij} v_j}
\leq 
\left( \sum_{ j\neq i} \abs{a_{ij}} \right) \abs{v_i}
\leq \abs{a_{ij}} \cdot \abs{v_i}.
$$
This implies $\abs{a_{ii} - \lambda} \leq \abs{a_{ij}}$ and leads to a contradiction.

For 2). From part (1), we known that $\bA$ is positive semidefinite. Suppose $\bA$ is not positive definite with a zero eigenvalue $0$ associated with a nonzero eigenvector $\bv$ such that $\bA\bv=\bzero$.   Similarly, we have 
$$
\abs{a_{ii}}\cdot \abs{v_i}
=
\abs{\sum_{j\neq i} a_{ij} v_j}
\leq 
\left( \sum_{ j\neq i} \abs{a_{ij}} \right) \abs{v_i}
< \abs{a_{ij}} \cdot \abs{v_i},
$$
which is impossible and the result follows.
\end{proof}

\section{Concluding Remarks on Positive Definite Matrices}\label{section:conc_pd}
In Section~\ref{section:equivalent-pd-psd}, we will prove that a matrix $\bA$ is PD if and only if $\bA$ can be factored as $\bA=\bP^\top\bP$, where $\bP$ is nonsingular. Additionally, in Section~\ref{section:unique-posere-pd}, we will show  that a PD matrix $\bA$ can be uniquely factored as $\bA =\bB^k$ for $k=\{2,3,\ldots\}$, where $\bB$ is also PD. Both results are consequences of the spectral decomposition of PD matrices.

In conclusion, for a positive definite matrix $\bA$, four distinct factorizations are feasible. We can factor it into: (a). $\bA=\bR^\top\bR$, where $\bR$ is an upper triangular matrix with positive diagonals, as shown in Theorem~\ref{theorem:cholesky-factor-exist} by Cholesky decomposition; (b). $\bA = \bP^\top\bP$, where $\bP$ is nonsingular in Theorem~\ref{theorem:nonsingular-factor-of-PD}; (c). or $\bA=\bS^\top\bS$, where $\bS$ has full column rank; (d). and $\bA = \bB^k$ for $k=\{2,3,\ldots\}$, where $\bB$ is PD in Theorem~\ref{theorem:unique-factor-pd}. For clarity, the different factorizations of a positive definite matrix $\bA$ are summarized in Figure~\ref{fig:pd-summary}.

\begin{figure}[htbp]
\centering
\centering
\resizebox{0.9\textwidth}{!}{%
\begin{tikzpicture}[>=latex]

\tikzstyle{state} = [draw, very thick, fill=white, rectangle, minimum height=3em, minimum width=6em, node distance=8em, font={\sffamily\bfseries}]
\tikzstyle{stateEdgePortion} = [black,thick];
\tikzstyle{stateEdge} = [stateEdgePortion,->];
\tikzstyle{stateEdge2} = [stateEdgePortion,<->];
\tikzstyle{edgeLabel} = [pos=0.5, text centered, font={\sffamily\small}];

\node[ellipse, name=pdmatrix, draw,font={\sffamily\bfseries},  node distance=7em, xshift=-9em, yshift=-1em,fill={colorals}]  {PD Matrix $\bA$};
\node[state, name=bsqure, below of=pdmatrix, xshift=-5em, yshift=1em, fill={colorlu}] {$\bB^k$};
\node[state, name=ptp, right of=bsqure, xshift=14em, fill={colorlu}] {$\bP^\top\bP$};
\node[state, name=rsqure, left of=bsqure, xshift=-4em, fill={colorlu}] {$\bR^\top\bR$};
\node[state, name=sts, right of=bsqure, xshift=3em,  fill={colorlu}] {$\bS^\top\bS$};
\node[ellipse, name=utv, below of=pdmatrix,draw,  node distance=7em, xshift=-5em, yshift=-4em,font={\tiny},fill={coloruppermiddle}]  {PD $\bB$};
\node[ellipse, name=upperr, left of=utv, draw, node distance=8em, xshift=-4em,font={\tiny},fill={coloruppermiddle}]  {\parbox{6em}{Upper \\Triangular $\bR$}};
\node[ellipse, name=nonp, right of=utv,draw,  node distance=8em, xshift=14em, font={\tiny},fill={coloruppermiddle}]  {\parbox{7em}{Nonsingular $\bP$}};
\node[ellipse, name=nons, right of=utv,draw,  node distance=8em, xshift=3em, font={\tiny},fill={coloruppermiddle}]  {\parbox{6em}{Full column\\ rank $\bS$}};

\coordinate (lq2inter3) at ($(pdmatrix.east -| ptp.north) + (-0em,0em)$);
\draw (pdmatrix.east) edge[stateEdgePortion] (lq2inter3);
\draw (lq2inter3) edge[stateEdge] 
node[edgeLabel, text width=7.25em, yshift=0.8em]{\parbox{5em}{Spectral\\Decomposition}} (ptp.north);

\coordinate (rqr2inter1) at ($(pdmatrix.west) + (0,0em)$);
\coordinate (rqr2inter3) at ($(rqr2inter1-| rsqure.north) + (-0em,0em)$);
\draw (rqr2inter1) edge[stateEdgePortion] (rqr2inter3);
\draw (rqr2inter3) edge[stateEdge] 
node[edgeLabel, text width=8em, yshift=0.8em]{\parbox{2em}{LU/\\Spectral/\\Recursive}} (rsqure.north);

\draw (pdmatrix.south)
edge[stateEdge] node[edgeLabel, xshift=-1.0em,yshift=0.5em]{\parbox{5em}{Spectral\\Decomposition} } 
(bsqure.north);

\draw (pdmatrix.south)
edge[stateEdge] node[edgeLabel,xshift=1.0em, yshift=0.5em]{\parbox{5em}{Spectral\\Decomposition} } 
(sts.north);

\draw (upperr.north)
edge[stateEdge] node[edgeLabel,yshift=0.5em]{} 
(rsqure.south);

\draw (utv.north)
edge[stateEdge] node[edgeLabel,yshift=0.5em]{} 
(bsqure.south);

\draw (nonp.north)
edge[stateEdge] node[edgeLabel,yshift=0.5em]{} 
(ptp.south);

\draw (nons.north)
edge[stateEdge] node[edgeLabel,yshift=0.5em]{} 
(sts.south);

\begin{pgfonlayer}{background}
\draw [join=round,cyan,dotted,fill={colormiddle}] ($(upperr.south west) + (-1.6em, -1em)$) rectangle ($( pdmatrix.east-|ptp.north east) + (1.6em, +1.8em)$);
\end{pgfonlayer}

\end{tikzpicture}
}
\caption{Demonstration of different factorizations for a positive definite matrix $\bA$.}
\label{fig:pd-summary}
\end{figure}

\section{Pivoted Cholesky Decomposition}\label{section:volum-picot-cholesjy}
If $\bP$ is a permutation matrix and $\bA$ is positive definite, then $\bP^\top\bA\bP$ is said to be a \textit{diagonal permutation} of $\bA$ (among other things, it permutes the diagonals of $\bA$). Any such diagonal permutation of $\bA$
is positive definite and has a Cholesky factor. Such a factorization is called a \textit{pivoted Cholesky factorization}. There are many ways to pivot a
Cholesky decomposition, but the most common one is the complete pivoting (see Section~\ref{section:complete-pivoting}) such that 
$$
\bP\bA\bP^\top = \bR^\top\bR 
$$
is the column-pivoted Cholesky decomposition of $\bA$, where $\bP$ is a permutation matrix, and $\bR$ is upper triangular.

Following the Cholesky decomposition via recursive calculation in Algorithm~\ref{alg:compute-choklesky11}, we notice from Theorem~\ref{theorem:positive-maximal-diagonal} that the maximal element of a PD matrix lies on the diagonal. Therefore, the complete pivoting algorithm for Cholesky decomposition can only search in the diagonals. The procedure is shown in Algorithm~\ref{alg:cholesky-complete-pivot22}.
\begin{algorithm}[h] 
\caption{Cholesky Decomposition via Recursive Algorithm: Complete Pivoting} 
\label{alg:cholesky-complete-pivot22}  
\begin{algorithmic}[1] 
\Require 
Positive definite matrix $\bA$ with size $n\times n$; 
\State $\bR \in \real^{n\times n}$ is initialized with all zeros;
\For{$k=1$ to $n$} \Comment{compute the $k$-th row of $\bR$}
\State Search in the diagonals of $\bA_{k:n,k:n}$ such that $a_{vv} = \max \bA_{k:n,k:n}$;
\State Swap the $k$-th and $v$-th column of $\bA$: $\bA_{:,k}\leftrightarrow\bA_{:,v}$ by column permutation;
\State Swap the $k$-th and $v$-th column of $\bR$: $\bR_{:,k}\leftrightarrow\bR_{:,v}$ by column permutation;
\State Swap the $k$-th and $v$-th row of $\bA$: $\bA_{k,:}\leftrightarrow\bA_{v,:}$ by row permutation;
\State $r_{kk} \leftarrow \sqrt{a_{kk}}$; \Comment{first element of $k$-th row, 1 flop}
\State $\bR_{k,k+1:n} \leftarrow \frac{1}{r_{kk}} \bA_{k,k+1:n}$; \Comment{the rest elements of $k$-th row, $n-k$ flops}
\State $\bA_{k+1:n,k+1:n} \leftarrow \bA_{k+1:n,k+1:n} - \bR_{k,k+1:n}^\top\bR_{k,k+1:n}$; \Comment{$2(1+2+\ldots+(n-k))$ flops}
\EndFor
\end{algorithmic} 
\end{algorithm}


\section{Decomposition for Semidefinite Matrices}
For positive semidefinite matrices, the Cholesky-like decomposition also exists with a slight modification.

\begin{theoremHigh}[Semidefinite Decomposition\index{Positive semidefinite}]\label{theorem:semidefinite-factor-exist}
Let $\bA\in \real^{n\times n}$ be any  \textit{positive semidefinite} matrix. Then, it can be factored as 
$$
\bA = \bR^\top\bR,
$$
where $\bR \in \real^{n\times n}$ is an upper triangular matrix with possible \textbf{zero} diagonal elements, and the factorization is \textbf{not unique} in general. 
\end{theoremHigh}
For such a decomposition, the diagonal of $\bR$ may not display the rank of $\bA$ \citep{higham2009cholesky}. 
\begin{example}[\citep{higham2009cholesky}]
Consider the matrix 
$$
\bA = 
\begin{bmatrix}
	1 & -1 & 1 \\
	-1 & 1 & -1 \\
	1 & -1 & 2 
\end{bmatrix}.
$$
The semidefinite decomposition is given by 
$$
\bA = 
\begin{bmatrix}
	1 & 0 & 0 \\
	-1 & 0 & 0 \\
	1 & 1 & 0 
\end{bmatrix}
\begin{bmatrix}
	1 & -1 & 1 \\
	0 & 0 & 1 \\
	0 & 0 & 0
\end{bmatrix}
=\bR^\top\bR.
$$
The matrix $\bA$ has a rank of 2, while $\bR$ has only one nonzero element on its diagonal.
\end{example}

We notice that all PD matrices have full rank, and this characteristic profoundly influences many of our proofs discussed above. 
This property can be established by Sylvester's criterion (Theorem~\ref{theorem:sylvester-criterion}) that all the leading principal minors of PD matrices are positive. Alternatively, we can simply prove that if a PD matrix $\bA$ were rank-deficient, this would imply that the null space of $\bA$ has a positive dimension, and there would exist a nonzero vector $\bx$ in the null space of $\bA$ such that $\bA\bx=\bzero$, i.e., $\bx^\top\bA\bx=0$. This leads to a contradiction to the definition of  PD matrices.

However, this is not necessarily true for PSD matrices, where the dimension of the null space can be greater than zero.
Therefore, and more generally, a rank-revealing decomposition for semidefinite decomposition is provided as follows.\index{Rank-revealing}\index{Semidefinite rank-revealing}
\begin{theoremHigh}[Semidefinite Rank-Revealing Decomposition\index{Rank-revealing}]\label{theorem:semidefinite-factor-rank-reveal}
Let $\bA\in \real^{n\times n}$ be any  positive semidefinite matrix with rank $r$. Then, it can be factored as 
$$
\bP^\top \bA\bP  = \bR^\top\bR, \qquad \mathrm{with} \qquad 
\bR = \begin{bmatrix}
	\bR_{11} & \bR_{12}\\
	\bzero &\bzero 
\end{bmatrix} \in \real^{n\times n},
$$
where $\bP$ is a permutation matrix, $\bR_{11} \in \real^{r\times r}$ is an upper triangular matrix with positive diagonal elements, and $\bR_{12}\in \real^{r\times (n-r)}$. 
\end{theoremHigh}
The proof for the existence of the above rank-revealing decomposition for semidefinite matrices is delayed to Section~\ref{section:semi-rank-reveal-proof} as a consequence of the spectral decomposition (Theorem~\ref{theorem:spectral_theorem}) and the column-pivoted QR decomposition (Theorem~\ref{theorem:rank-revealing-qr-general}). Whereas, the rigorous proof for the trivial Semidefinite Decomposition Theorem~\ref{theorem:semidefinite-factor-exist} can be a direct result of the spectral decomposition and the standard QR decomposition (Theorem~\ref{theorem:qr-decomposition}).

\index{Bunch-Kaufman decomposition}
When the matrix $\bA$ is symmetric and indefinite, we can use a \textit{symmetric indefinite  decomposition}, also known as the \textit{Bunch Kaufman decomposition} \citep{bunch1977some}. 
\begin{theoremHigh}[Bunch Kaufman Decomposition]\label{theorem:Bunch_Kaufman}
Let $\bA\in \real^{n\times n}$ be any symmetric (indefinite) matrix. Then, it can be factored as 
$$
\bP^\top \bA\bP  = \bL\bB\bL^\top, 
$$
where $\bP$ is a permutation matrix, $\bL$ is an unit lower triangular matrix, and $\bB$ is a block-diagonal matrix with each diagonal block of  $\bB$ being either a  $1\times 1$ or a  $2\times 2$ matrix.
\end{theoremHigh}
This type of decomposition is sometimes referred to as an 
$\bL\bB\bL^\top$ decomposition.
It is particularly useful in practical applications for solving linear systems and computing eigenvalues of matrices, especially when a direct Cholesky decomposition is not applicable (e.g., when the matrix is not positive definite) \citep{dumas2018symmetric}.

\index{Newton's method}
\index{Modified Newton's method}
\section{Application: Modified Newton's Method}\label{section:app_cho_md_newton}
When optimizing or minimizing a function $f(\bx)$ over $\bx$, 
the plain Newton's method~\footnote{see, for example, \citet{lu2022gradient}.} updates the estimate at the  $t$-th iteration as
$$
\bx^{(t+1)} \leftarrow \bx^{(t)} + \bd^{(t)},
$$
where $(\nabla^2 f(\bx^{(t)}) )\bd^{(t)} = -\nabla f(\bx^{(t)})$ determines the ``candidate" descent direction $\bd^{(t)}$.
The vector $\bd^{(t)}$ is  a descent direction only when the Hessian $(\nabla^2 f(\bx^{(t)}) )$ is PD, which is not always the case.

The modified Newton's method addresses this issue by approximating the Hessian with $\bH^{(t)} = \nabla^2 f(\bx^{(t)}) +\bE^{(t)}$, ensuring that $\bH^{(t)}$ is PD \citep{gill2019practical}.
Given the Cholesky decomposition in the form $\nabla^2 f(\bx^{(t)})=\bL\bD\bL^\top=\bR^\top\bR$ (where $\bR=\bD^{1/2}\bL^\top$) and the condition number inequality $\cond(\nabla^2 f(\bx^{(t)}))\geq \cond(\bD)$ (see Equation~\eqref{equation:cond_pd_ineq}), the goal of modified Newton's method can be approximately achieved by adjusting the diagonals of $\bD$.
To be more specific, when computing the Cholesky decomposition using Algorithm~\ref{alg:compute-choklesky-_ldl}, the modified Newton's method imposes bounds on the diagonal $d_{jj}$, given two parameters $\alpha$ and $\beta$, such that 
$$
d_{jj} \geq \alpha,  \gap  l_{ij}\sqrt{d_{jj}}\leq \beta, \,i=\{j+1,j+2,\ldots,n\}.
$$
The latter constraint serves  to upper-bound each row of $\bR$, since $\bR=\bD^{1/2}\bL^\top$. And this is equivalent to updating each $d_{jj}$ in Algorithm~\ref{alg:compute-choklesky-_ldl} by
$$
d_{jj} = \max\left\{ \abs{c_{jj}},\, \beta,\, \mathop{\max}_{i>j}\abs{c_{ij}}\right\}.
$$

\section{Application: Rank-One Update/Downdate}\label{section:cholesky-rank-one-update}
Updating linear systems after low-rank modifications of the system matrix is a  common practice in machine learning, statistics, and many other fields. 
However, it is widely recognized that such updates can introduce significant instabilities due to round-off errors \citep{seeger2004low}. 
When the system matrix is positive definite, employing a representation based on Cholesky decomposition is usually preferable as it provides improved numerical stability. 
In this section, we will present a proof for the rank-one update/downdate using Cholesky decomposition.
\subsection{Rank-One Update}\index{Rank-one update}
A rank-one update $\bA^\prime$ of a matrix $\bA$ by a vector $\bv$ is of the form \citep{gill1974methods, bojanczyk1987note,  chang1997pertubation, davis1999modifying, seeger2004low, chen2008algorithm, davis2008user, higham2009cholesky}:
\begin{equation*}
\begin{aligned}
	\bA^\prime &= \bA + \bv \bv^\top;\\
		\downarrow &\gap  \downarrow\\
	\bR^{\prime\top}\bR^\prime &= \bR^\top\bR + \bv \bv^\top.
\end{aligned}
\end{equation*}
If the Cholesky factor $\bR$ of $\bA \in \real^{n\times n}$ has already been computed, we can efficiently obtain the Cholesky factor $\bR^\prime$ of $\bA^\prime$.
Note that $\bA^\prime$ differs from $\bA$ only by the symmetric rank-one matrix. 
Therefore, we can compute $\bR^\prime$ from $\bR$ using the rank-one Cholesky update, which takes $\mathcalO(n^2)$ operations, each saving from $\mathcalO(n^3)$ complexity if we were to recompute the Cholesky decomposition of $\bA^\prime$ from scratch,
given that we  know $\bR$, the Cholesky decomposition of $\bA$ up front, i.e., we want to compute the Cholesky decomposition of $\bA^\prime$ via that of $\bA$. 
To see this,
suppose there exists a set of orthogonal matrices $\bQ_n \bQ_{n-1}\ldots \bQ_1$ such that  
$$
\bQ_n \bQ_{n-1}\ldots \bQ_1 
\begin{bmatrix}
\bv^\top \\
\bR
\end{bmatrix}
=
\begin{bmatrix}
\bzero \\
\bR^\prime
\end{bmatrix}.
$$
Then we find out the expression for the Cholesky factor of $\bA^\prime$ by $\bR^\prime$. 
Specifically, multiplying the left-hand side of the above equation by its transpose yields
$$
\begin{bmatrix}
\bv & \bR^\top
\end{bmatrix}
\bQ_1^\top \ldots \bQ_{n-1}^\top\bQ_n^\top
\bQ_n \bQ_{n-1}\ldots \bQ_1 
\begin{bmatrix}
\bv^\top \\
\bR
\end{bmatrix}
= \bR^\top\bR + \bv \bv^\top.
$$ 
Similarly, multiplying the right-hand side by its transpose gives
$$
\begin{bmatrix}
\bzero & \bR^{\prime\top}
\end{bmatrix}
\begin{bmatrix}
\bzero \\
\bR^\prime
\end{bmatrix}=\bR^{\prime\top}\bR^\prime,
$$
which matches the left-hand side equation. \textit{Givens rotations} are examples of such orthogonal matrices that can transfer $\bR$ and $\bv$ into $\bR^\prime$. We will discuss the fundamental meaning of Givens rotations shortly to prove the existence of the QR decomposition in Section~\ref{section:qr-givens}. 
In this section, we will only introduce their definition and present the results directly. Feel free to skip this section for a first reading.

\begin{definition}[$n$-th Order Givens Rotation]\label{definition:givens-rotation-in-qr}
A Givens rotation is represented by a matrix of the following form
$$
\bG_{kl}=
\begin{bmatrix}
	1 &          &   &  &   &   & &  & &\\
	& \ddots  &  &  &  & && & &\\
	&      & 1 &  & & &  && &\\
	&      &  & c &  &  &  & s & &\\
	&& &   & 1 & & && &\\
	&& &   &   &\ddots &  && &\\
	&& &  &   &  & 1&& &\\
	&& & -s &  &  & &c& &\\
	&& & &  &  & & &1 & \\
	&& & &  &  & & & &\ddots
\end{bmatrix}_{n\times n},
$$
where the $(k,k), (k,l), (l,k), (l,l)$ entries are $c, s, -s, c$ respectively, and $s = \sin \theta$ and $c=\cos \theta$ for some $\theta$.

Let $\bdelta_k \in \real^n$ be the zero vector except that the entry $k$ is 1 (the $k$-th unit basis vector). Then, mathematically, the Givens rotation defined above can be denoted by 
$$
\bG_{kl}\triangleq\bG_{kl}(\theta)= \bI + (c-1)(\bdelta_k\bdelta_k^\top + \bdelta_l\bdelta_l^\top) + s(\bdelta_k\bdelta_l^\top -\bdelta_l\bdelta_k^\top ),
$$
where the subscripts $k$ and $l$ indicate that the rotation occurs \textbf{in plane $k$ and $l$}.

Specifically, one can also define the $n$-th order Givens rotation, where $(k,k),$ $(k,l),$ $(l,k),$ $(l,l)$ entries are $c, \textcolor{mylightbluetext}{-s, s}, c$ respectively (note the difference in the sign of $s$). The underlying principles remain the same.
\end{definition}

\begin{exercise}
Show that $\bG_{kl}(-\theta)^{-1}=\bG_{kl}(\theta)$. \textit{Hint: Use the orthogonality of $\bG_{kl}(\theta)$.}
\end{exercise}

It can be easily verified that an $n$-th order Givens rotation is an orthogonal matrix, and its determinant is 1. For any vector $\bx =[x_1, x_2, \ldots, x_n]^\top \in \real^n$, the result of applying a Givens rotation $\bG_{kl}$ to $\bx$ is $\by = \bG_{kl}\bx$, where
$$ 
\left\{
\begin{aligned}
&y_k = c \cdot x_k + s\cdot x_l,   \\
&y_l = -s\cdot x_k +c\cdot x_l,  \\
&y_j = x_j . &  (j\neq k,l) 
\end{aligned}
\right.
$$
That is, a Givens rotation applied to $\bx$ rotates two components of $\bx$ by some angle $\theta$, while keeping all other components unchanged.

Now, let's consider a Givens rotation of order $(n+1)$, where the rotation is indexed from $0$ to $n$. This rotation can be expressed as 
$$
\bG_k \triangleq \bI + (c_k-1)(\bdelta_0\bdelta_0^\top + \bdelta_k\bdelta_k^\top) + s_k(\bdelta_0\bdelta_k^\top -\bdelta_k\bdelta_0^\top ),
$$
where $c_k = \cos \theta_k, s_k=\sin\theta_k$ for some $\theta_k$, $\bG_k \in \real^{(n+1)\times (n+1)}$, and $\bdelta_k\in \real^{n+1}$ is a zero vector except that the $(k+1)$-th entry is 1.

\begin{mdframed}[hidealllines=\mdframehideline,backgroundcolor=\mdframecolor]
Taking out the $k$-th column of the following equation 
$$
\begin{bmatrix}
	\bv^\top \\
	\bR
\end{bmatrix}
\longrightarrow 
\begin{bmatrix}
	\bzero \\
	\bR^\prime
\end{bmatrix},
$$
where we let the $k$-th element of $\bv$ be $v_k$, and the $k$-th diagonal of $\bR$ be $r_{kk}$.
We realize that $\sqrt{v_k^2 + r_{kk}^2} \neq 0$,
and let $c_k \triangleq \frac{r_{kk}}{\sqrt{v_k^2 + r_{kk}^2}}$, $s_k\triangleq-\frac{v_k}{\sqrt{v_k^2 + r_{kk}^2}}$. Then we have
$$ 
\left\{
\begin{aligned}
	&v_k \rightarrow c_kv_k+s_kr_{kk}=0;   \\
	&r_{kk}\rightarrow -s_k v_k +c_kr_{kk}= \sqrt{v_k^2 + r_{kk}^2} = r^\prime_{kk} .  \\
\end{aligned}
\right.
$$
That is, $\bG_k$ will introduce a zero value to the $k$-th element of $\bv$ and a nonzero value to $r_{kk}$.
\end{mdframed}
This finding above is crucial for the rank-one update. And we obtain 
$$
\bG_n \bG_{n-1}\ldots \bG_1 
\begin{bmatrix}
\bv^\top \\
\bR
\end{bmatrix}
=
\begin{bmatrix}
\bzero \\
\bR^\prime
\end{bmatrix}.
$$
For each Givens rotation, it takes $6n$ flops. And there are $n$ such rotations, which requires $6n^2$ flops if keeping only the leading term. The complexity to calculate the Cholesky factor of $\bA^\prime$ is thus reduced from $\frac{1}{3} n^3$ to $6n^2$ flops using the rank-one update, provided that we already know the Cholesky factor of $\bA$. 
The above algorithm is essential for reducing the complexity in the posterior calculation of Bayesian inference for Gaussian mixture model \citep{lu2021bayes}. At each stage, $k$ new samples are added or removed from an existing cluster, which corresponds to applying $k$ rank-one updates.

\subsection{Rank-One Downdate}\index{Rank-one downdate}
Let us consider the scenario where we have computed the Cholesky factor of $\bA$, and $\bA^\prime$ is the rank-one downdate of $\bA$ given by the following expression:
\begin{equation*}
\begin{aligned}
	\bA^\prime &= \bA - \bv \bv^\top;\\
		\downarrow &\gap  \downarrow\\
	\bR^{\prime\top}\bR^\prime &= \bR^\top\bR - \bv \bv^\top. 
\end{aligned}
\end{equation*}
The algorithm is similar by proceeding as follows:
\begin{equation}\label{equation:rank-one-downdate}
\bG_1 \bG_{2}\ldots \bG_n
\begin{bmatrix}
	\bzero \\
	\bR
\end{bmatrix}
=
\begin{bmatrix}
	\bv^\top \\
	\bR^\prime
\end{bmatrix}.
\end{equation}
Again, the set of Givens rotations $
\bG_k = \bI + (c_k-1)(\bdelta_0\bdelta_0^\top + \bdelta_k\bdelta_k^\top) + s_k(\bdelta_0\bdelta_k^\top -\bdelta_k\bdelta_0^\top )
$
for $k\in\{1,2,\ldots,n\}$
can be constructed as follows:
\begin{mdframed}[hidealllines=\mdframehideline,backgroundcolor=\mdframecolor]
Taking out the $k$-th column of the following equation 
$$
\begin{bmatrix}
	\bzero \\
	\bR
\end{bmatrix}
\longrightarrow
\begin{bmatrix}
	\bv^\top \\
	\bR^\prime
\end{bmatrix}.
$$
We realize that $r_{kk} \neq 0$,
and let $c_k\triangleq\frac{\sqrt{r_{kk}^2 - v_k^2}}{r_{kk}}$, $s_k \triangleq \frac{v_k}{r_{kk}}$. Then, we have
$$ 
\left\{
\begin{aligned}
	& 0 \rightarrow s_kr_{kk}=v_k;   \\
	&r_{kk}\rightarrow c_k r_{kk}= \sqrt{r_{kk}^2-v_k^2  }=r^\prime_{kk} .  \\
\end{aligned}
\right.
$$
This requires $r^2_{kk} > v_k^2$ to make $\bA^\prime$ to be positive definite. Otherwise, $c_k$ as defined above will not exist.
\end{mdframed}
Again, one can verify that, multiplying the left-hand side of Equation~\eqref{equation:rank-one-downdate} by its transpose, we have 
$$
\begin{bmatrix}
\bzero & \bR^\top
\end{bmatrix}
\bG_n^\top\ldots \bG_{2}^\top\bG_1^\top
\bG_1 \bG_{2}\ldots \bG_n
\begin{bmatrix}
\bzero \\
\bR
\end{bmatrix} =\bR^\top\bR.
$$
And multiplying the r.h.s. by its transpose, we have 
$$
\begin{bmatrix}
\bv & \bR^{\prime\top} 
\end{bmatrix}
\begin{bmatrix}
\bv^\top \\
\bR^\prime
\end{bmatrix}=\bv\bv^\top + \bR^{\prime\top}\bR^\prime.
$$
This results in $\bR^{\prime\top}\bR^\prime = \bR^\top\bR - \bv \bv^\top$.

\section{Application: Indefinite Rank-Two Update}\index{Rank-two update}

Let $\bA = \bR^\top\bR$ be the Cholesky decomposition of $\bA$, \citet{goldfarb1976factorized, seeger2004low} introduced a stable method for the indefinite rank-two update of the form 
$$
\bA^\prime = (\bI+\bv\bu^\top)\bA(\bI+\bu\bv^\top).
$$
Let
$$
\bigg\{ 
\begin{aligned}
\bz &\triangleq \bR^{-\top}\bv,   \\
\bw &\triangleq \bR\bu,
\end{aligned}
\qquad 
\implies
\qquad  
\bigg\{
\begin{aligned}
\bv &= \bR^{\top}\bz,   \\
\bu &= \bR^{-1}\bw.
\end{aligned}
$$
And suppose the LQ decomposition \footnote{We will shortly introduce in Theorem~\ref{theorem:lq-decomposition}.} of $\bI+\bz\bw^\top$ is given by $\bI+\bz\bw^\top =\bL\bQ$, where $\bL$ is lower triangular and $\bQ$ is orthogonal. Thus, we have 
$$
\begin{aligned}
\bA^\prime &= (\bI+\bv\bu^\top)\bA(\bI+\bu\bv^\top)
= (\bI+\bR^{\top}\bz   \bw^\top \bR^{-\top })\bA(\bI+\bR^{-1}\bw \bz^\top\bR)\\
&= \bR^\top (\bI+\bz\bw^\top)(\bI+\bw\bz^\top)\bR
= \bR^\top\bL \bQ\bQ^\top \bL^\top\bR \\
&= \bR^\top\bL \bL^\top\bR.
\end{aligned}
$$
Let $\bR^\prime \triangleq \bR^\top\bL$, which is lower triangular. This establishes the Cholesky decomposition of $\bA^\prime$.


%

\begin{problemset}
\item \label{prob:pd_cmequiv} \textbf{PD.} Let $\bA\in\real^{n\times n}$. Show that the following two statements are equivalent:
\begin{itemize}
	\item $\bx^*\bA\bx>0$ for all nonzero $\bx\in\complex^n$.
	\item $\bx^\top\bA\bx>0$ for all nonzero $\bx\in\real^n$.
\end{itemize}

\item \label{prob:tr_de_pd} \textbf{Trace, det of PD/PSD/ND matrices.} Let $\bA$ be positive definite (resp. positive semidefinite), show that $\trace(\bA), \det(\bA)$, and the principal minors of $\bA$ are all positive (resp. nonnegative). Moreover, $\trace(\bA)=0$ if and only if $\bA=\bzero$. 
Let $\bB\in\real^{n\times n}$ be negative definite. Show that $\trace(\bB)$ is negative; $\det(\bB)$ is negative for odd $n$ and positive for even $n$.
\textit{Hint: Use Theorem~\ref{theorem:eigen_charac}, Theorem~\ref{theorem:eigen_trace}, and Theorem~\ref{theorem:eigen_trace2}.}

\item Show that the following matrix is positive definite and find its Cholesky decomposition:
$$
\bA = 
\begin{bmatrix}
5 & -1 & 0 \\
-1 & 4 & 2 \\
0 & 2 & 8
\end{bmatrix}.
$$
\item Let $\bA,\bB\in\real^{n\times n}$ be  positive semidefinite matrices. Show that $\bA+\bB$ is also positive semidefinite.
\item Let $\bA\in\real^{n\times n}$ and $\bB\in\real^{m\times m}$ be  two symmetric matrices. Show that the following two claims are equivalent:
\begin{enumerate}
\item $\bA$ and $\bB$ are positive semidefinite.
\item $\scriptsize\begin{bmatrix}
\bA& \bzero \\
\bzero & \bB
\end{bmatrix}$ is positive semidefinite.
\end{enumerate}
\item 
Given a matrix $\bB\in\real^{n\times k}$ and let $\bA=\bB\bB^\top$. Show that 
\begin{enumerate}
\item $\bA$ is positive semidefinite.
\item $\bA$ is positive definite if and only if $\bB$ has full row rank.
\end{enumerate}
\item Show that if $\bA$ is positive semidefinite and nonsingular, then $\bA^{-1}$ is positive definite.

\item \label{problem:indefinite_reverse} We have shown in Proposition~\ref{proposition:indefin_diag} that if a matrix contains both positive and negative elements, then the matrix is indefinite. Show that the reverse claim is not true. \textit{Hint: Find a matrix that contains only positive or negative elements on  its diagonal; however, this matrix has both positive and negative eigenvalues.}


\item \textbf{Quadratic form.} Consider the quadratic form $L(\bx) = \frac{1}{2} \bx^\top \bA \bx - \bb^\top \bx + c$ in Equation~\eqref{equation:quadratic-form-general-form}, where $\bA\in\real^{d\times d}$ is symmetric, $ \bx\in \real^d$, and $c\in\real$. Show that the following two claims are equivalent~\footnote{This result finds its application in non-convex quadratic constrained quadratic problems (QCQPs) \citep{beck2014introduction}.}:
\begin{itemize}
\item $L(\bx) \geq 0$ for all $\bx\in\real^d$.
\item $\scriptsize\begin{bmatrix}
\bA & -\bb\\
-\bb^\top & 2c
\end{bmatrix}\geq 0$.
\end{itemize} 
\textit{Hint: Apply eigenvalue characterization theorem on $\bA$.}

\item \textbf{Quadratic form.} Consider again  the quadratic form $L(\bx) = \frac{1}{2} \bx^\top \bA \bx - \bb^\top \bx + c$ in Equation~\eqref{equation:quadratic-form-general-form}, where $\bA\in\real^{d\times d}$, $ \bx\in \real^d$, and $c\in\real$. Suppose $\bA$ is positive semidefinite. 
Show that $L(\bx)$ is bounded below over $\real^n$ if and only if $\bb$ is in the column space of $\bA$.

\item \textbf{Quadratic form.} Consider again   the quadratic form $L(\bx) = \frac{1}{2} \bx^\top \bA \bx - \bb^\top \bx + c$ in Equation~\eqref{equation:quadratic-form-general-form}. Show that $L(\bx)$ is coercive~\footnote{A function $f(\bx):\real^n\rightarrow \real$ is called coercive if $\mathop{\lim}_{\bx\rightarrow \infty} f(x)=\infty$.} if and only if $\bA$ is PD.

\item \label{prob:sam_quad} \textbf{Quadratic form.} Let $\bA\in\real^{n\times n}$ be a general square matrix (not necessarily symmetric). Show that $\bx^\top\bA\bx=\bx^\top[\frac{1}{2}(\bA+\bA^\top)]\bx$. The latter quadratic form is induced from a symmetric matrix.

\item \textbf{Constructing symmetric form.} Let $P(\bA)\triangleq\frac{1}{2}(\bA+\bA^\top)$ if $\bA\in\real^{n\times n}$. Show that 
\begin{itemize}
\item \textit{Null space.} $\nspace(\bA)\subset \nspace(P(\bA))$ and $\nspace(\bA^\top)\subset \nspace(P(\bA))$ such that $\rank(P(\bA))\leq \rank(\bA)$.
\item When $\rank(P(\bA))= \rank(\bA)$, it follows that $\bA$, $\bA^\top$, and $P(\bA)$ have the same null space.
\end{itemize}
\textit{Hint: Consider the quadratic form $\bx^\top\bA\bx$ and $\bx^\top\bA^\top\bx$, and  use Exercise~\ref{exercise:sing_psd}.}

\item \label{prob:quad_nsdpsd} Let $\bA\in\real^{n\times n}$. It can be easily obtained that $\nspace(\bA)\subseteq\{\bx\in\real^n: \bx^\top\bA\bx=0\}$, i.e., the null space of $\bA$ is a subset of the set making the quadratic form vanish. Show that the equality is attained if and only if $\bA$ is PSD or NSD. \textit{Hint: Exercise~\ref{exercise:sing_psd}.}

\item \label{problem:cond_pd} Given the Cholesky decomposition of a PD matrix: $\bA=\bL\bD\bL^\top$, show that $\cond(\bA)\geq \cond(\bD)$.

\end{problemset}

%% file: chapter-qr.tex
\newpage 
\part{Triangularization, Orthogonalization, and Gram-Schmidt Process}
\section*{Introduction}
\lettrine{\color{caligraphcolor}G}
Let  $\bA=[\ba_1, \ba_2, \ldots, \ba_l]$ be an $m\times l$ matrix with $m\geq l$, the process of \textit{orthonormalization} produces an $m\times l$ matrix $\bQ_l=[\bq_1, \bq_2, \ldots, \bq_l]$ such that
$$
\spn([\bq_1, \bq_2, \ldots, \bq_k]) = \spn([\bq_1, \ba_2, \ldots, \ba_k]), \gap \text{for all } k\in \{1,2,\ldots, l\}.
$$ 
Additionally, the columns of $\bQ_l$ are \textit{mutually orthonormal}:
$$
\bq_i^\top\bq_j = 
\left\{
\begin{aligned}
&1, \qquad \text{if $i=j$};\\
&0, \qquad \text{if $i\neq j$}.
\end{aligned}
\right\} \implies  \bQ_l^\top\bQ_l = \bI_l,
$$
where $\bI_l$ is the $l\times l$ identity matrix. When we complete $\bQ_l$ into $m$ mutually orthonormal columns $\bQ=[\bq_1, \bq_2, \ldots, \bQ_m]$, $\bQ$ becomes a square matrix, and it follows that $\bQ\bQ^\top=\bQ^\top\bQ=\bI_m$. This extension, often referred to as  \textit{column completion}, can be achieved using the \textit{Gram-Schmidt process}, which will be introduced in the sequel. 

For the rest of this part, we will explore several factorization methods that focus on orthogonalization or orthonormalization, each with different emphases, such as on the column space, row space, or both. An overview of these concepts is provided in Figure~\ref{fig:orthgonal-world-picture}.

\begin{figure}[htbp]
\centering
\begin{widepage}
\centering
\resizebox{0.6\textwidth}{!}{%
\begin{tikzpicture}[>=latex]

\tikzstyle{state} = [draw, very thick, fill=white, rectangle, minimum height=3em, minimum width=6em, node distance=8em, font={\sffamily\bfseries}]
\tikzstyle{stateEdgePortion} = [black,thick];
\tikzstyle{stateEdge} = [stateEdgePortion,->];
\tikzstyle{stateEdge2} = [stateEdgePortion,<->];
\tikzstyle{edgeLabel} = [pos=0.5, text centered, font={\sffamily\small}];

\node[state, name=qr, node distance=7em, xshift=-9em, yshift=-1em, fill={colorqr}] {QR};
\node[state, name=lq, above of=qr, xshift=-10em, yshift=-4em, fill={colorqr}] {LQ};
\node[state, name=utv, above of=qr, xshift=10em, yshift=-4em, fill={colorqr}] {UTV};
\node[state, name=twosidedortho, draw, above of=qr,xshift=0em, yshift=0em,fill={colorqr}]  {\parbox{5em}{Two-Sided \\Orthogonal}};
	
\draw (lq.east)
edge[stateEdge] node[edgeLabel, yshift=0em, xshift=-0.7em]{\parbox{3em}{Rank\\Estimation}} 
(utv.west) ;

\draw ($(qr.east)$)
edge[stateEdge, bend left=-12.5] node[edgeLabel, xshift=-0.3em,yshift=0.5em]{Rank Estimation} 
(utv.south) ;

\draw ($(lq.north) + (0em,0em)$)
edge[stateEdge, bend left=12.5] node[edgeLabel, xshift=1.5em,yshift=-0.5em]{Phase 1} 
(twosidedortho.west) ;
\draw (twosidedortho.east) 
edge[stateEdge, bend left=+12.5] node[edgeLabel, xshift=-1.5em,yshift=-0.5em]{Phase 2} 
($(utv.north) + (0em,0em)$);

\draw ($(qr.west)$)
edge[stateEdge, bend left=+12.5] node[edgeLabel, yshift=0.5em, xshift=1em]{Row Space} 
(lq.south) ;

\begin{pgfonlayer}{background}
\draw [join=round,cyan,dotted,fill={colormiddleleft}] ($(qr.south west -| lq.west) + (-0.5em, -0.5em)$) rectangle ($(twosidedortho.north east -| utv.north east) + (0.6em, 0.5em)$);		
\end{pgfonlayer}

\end{tikzpicture}
}
\end{widepage}
\caption{Orthogonalization World Map. Also, consider its position within the matrix decomposition world map, as illustrated in Figure~\ref{fig:matrix-decom-world-picture}.}
\label{fig:orthgonal-world-picture}
\end{figure}

\newpage
\chapter{QR Decomposition and Orthogonality}
\begingroup
\hypersetup{
	linkcolor=structurecolor,
	linktoc=page,  
}
\minitoc \newpage
\endgroup
\section{QR Decomposition}
\lettrine{\color{caligraphcolor}I}
In many applications, we are interested in the column space of a matrix $\bA=[\ba_1, \ba_2, \ldots, \ba_n] \in \real^{m\times n}$. The successive subspaces spanned by the columns $\ba_1, \ba_2, \ldots$ of $\bA$ are crucial for understanding the inherent structure and properties of the matrix:
$$
\cspace([\ba_1])\,\,\,\, \subseteq\,\,\,\, \cspace([\ba_1, \ba_2]) \,\,\,\,\subseteq\,\,\,\, \cspace([\ba_1, \ba_2, \ba_3])\,\,\,\, \subseteq\,\,\,\, \ldots,
$$
where $\cspace([\ldots])$ represents the subspace spanned by the vectors included in the brackets (alternatively, it can be denoted by $\cspace([\ldots])=\spn\{\ldots\}$).
Moreover, the notion of orthogonal or orthonormal bases within the column space plays a crucial role in various algorithms, allowing for efficient computations and interpretations.
QR factorization is an approach employed to analyze and decompose matrices in a way that exposes their column space structure.
The idea of QR decomposition involves  the construction of a sequence of orthonormal vectors $\bq_1, \bq_2, \ldots$ that span the same successive subspaces:
$$
\big\{\cspace([\bq_1])=\cspace([\ba_1])\big\}\subseteq 
\big\{\cspace([\bq_1, \bq_2])=\cspace([\ba_1, \ba_2])\big\}\subseteq
\big\{\cspace([\bq_1, \bq_2, \bq_3])=\cspace([\ba_1, \ba_2, \ba_3])\big\} 
\subseteq \ldots.
$$

We provide the result of QR decomposition in the following theorem and we delay the discussion of its existence in  subsequent sections.
\index{Decomposition: QR}
\begin{theoremHigh}[QR Decomposition]\label{theorem:qr-decomposition}
Let $\bA=[\ba_1, \ba_2, \ldots, \ba_n] \in\real^{m\times n}$ be any matrix (regardless of whether its columns are linearly independent or dependent) with $m\geq n$. Then, it can be factored as 
$$
\bA = \bQ\bR,~\footnote{If $\bA$ is complex, then $\bQ$ is unitary (or semi-unitary with orthonormal columns) and $\bR$ is complex upper triangular.}
$$
where 
\begin{enumerate}
\item  \textbf{Reduced}: $\bQ$ is $m\times n$ with orthonormal columns, and $\bR$ is an $n\times n$ upper triangular matrix, known as the \textbf{reduced QR decomposition};
\item \textbf{Full}: $\bQ$ is $m\times m$ with orthonormal columns, and $\bR$ is an $m\times n$ upper triangular matrix, known as the \textbf{full QR decomposition}. If we further restrict the upper triangular matrix to be a square matrix, the full QR decomposition can be denoted by 
$$
\bA = \bQ\begin{bmatrix}
	\bR_0\\
	\bzero
\end{bmatrix},
$$
where $\bR_0$ is an $n\times n$ upper triangular matrix. 
\end{enumerate}

Specifically, when $\bA$ has full rank, i.e., $\bA$  has linearly independent columns, $\bR$ also exhibits linearly independent columns, and $\bR$ is nonsingular in the \textit{reduced} case. 
This implies that the diagonals of $\bR$ are nonzero. Under this condition, when we further restrict elements on the diagonal of $\bR$ to be  \textbf{positive} and $\rank(\bA)=n$, the \textit{reduced} QR decomposition is \textbf{unique} (Corollary~\ref{corollary:unique-qr}). The \textit{full} QR decomposition is normally not unique since the right-most $(m-n)$ columns in $\bQ$ can be arranged in any order.
\end{theoremHigh}

The QR decomposition is a method of factorizing a matrix into an orthogonal matrix  $\bQ$ and an upper triangular matrix $\bR$. This technique is particularly useful in linear algebra and has a variety of significant applications. For example, QR decomposition can be employed to solve linear systems of equations $\bA\bx=\bb$. Since $\bQ$ is orthogonal ($\bQ\bQ^\top=\bI$), the system can be transformed into $\bR\bx=\bQ^\top\bb$, which is easier to solve due to the upper triangular nature of $\bR$.

In the field of statistics and data fitting, QR decomposition is used to find the least squares solution for an overdetermined system, where there are more equations than unknowns. This provides a method to minimize the sum of the squares of the residuals. Moreover, QR decomposition can be utilized to compute the eigenvalues and eigenvectors of a matrix, especially within algorithms like the \textit{QR algorithm} (Chapter~\ref{section:eigenvalue-problem}).

If a matrix $\bA$ can be decomposed into $\bQ$ and $\bR$, then its inverse $\bA^{-1}$ can be calculated as $\bR^{-1}\bQ^\top$. This approach is more stable and efficient than direct inversion for certain types of matrices.  The QR decomposition is also numerically stable, meaning it is less susceptible to round-off errors compared to some other matrix decompositions, making it ideal for computational tasks. In certain applications, the structure of the matrix, like sparsity, can be exploited during the QR decomposition process to reduce computational complexity.



\section{Projecting a Vector Onto Another Vector}\label{section:project-onto-a-vector}

\begin{figure}[h!]
\centering             
\vspace{-0.35cm} 
\subfigtopskip=2pt 
\subfigbottomskip=2pt 
\subfigcapskip=-5pt 
\subfigure[Project onto a line.]{\label{fig:project-line}
	\includegraphics[width=0.4\linewidth]{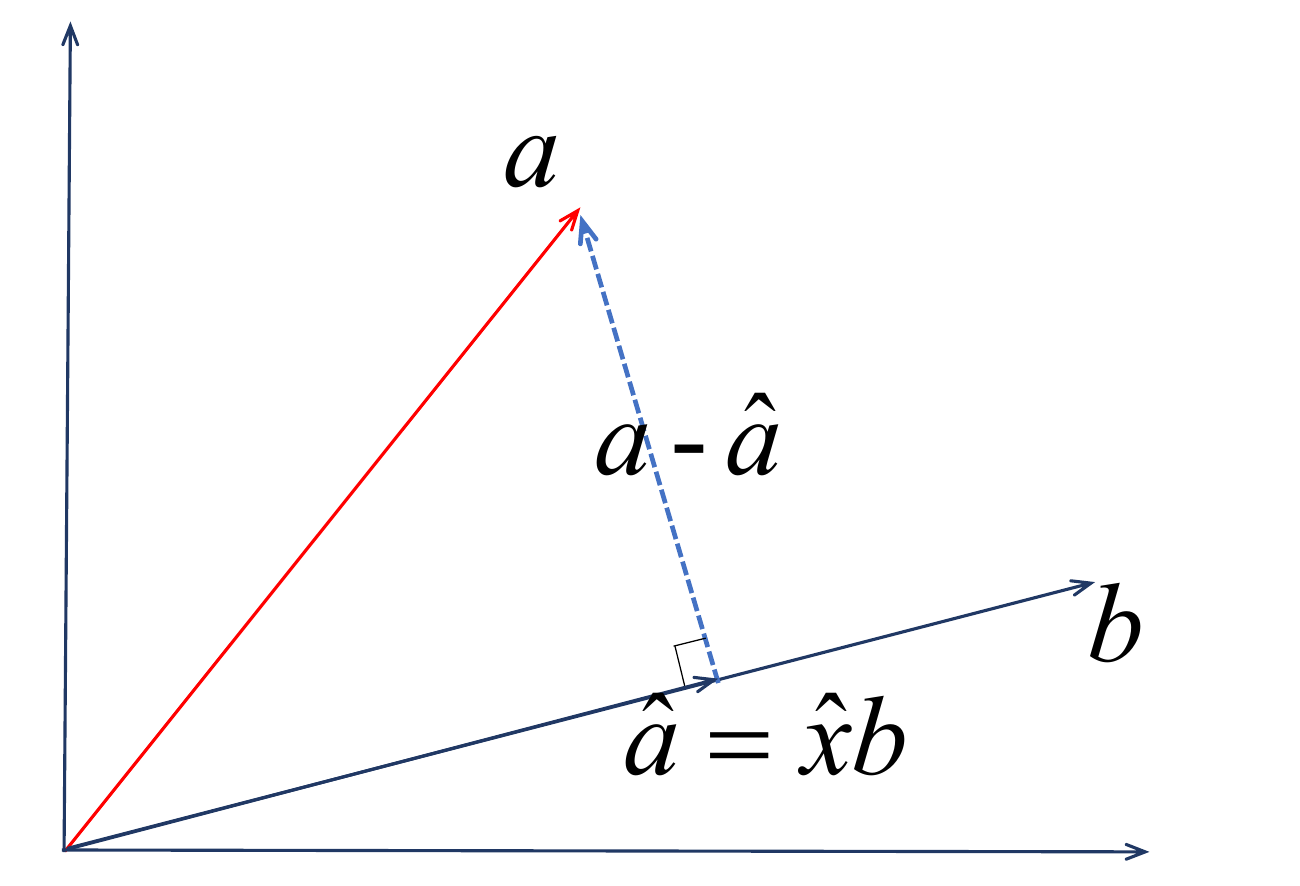}}
\quad 
\subfigure[Project onto a space.]{\label{fig:project-space}
	\includegraphics[width=0.4\linewidth]{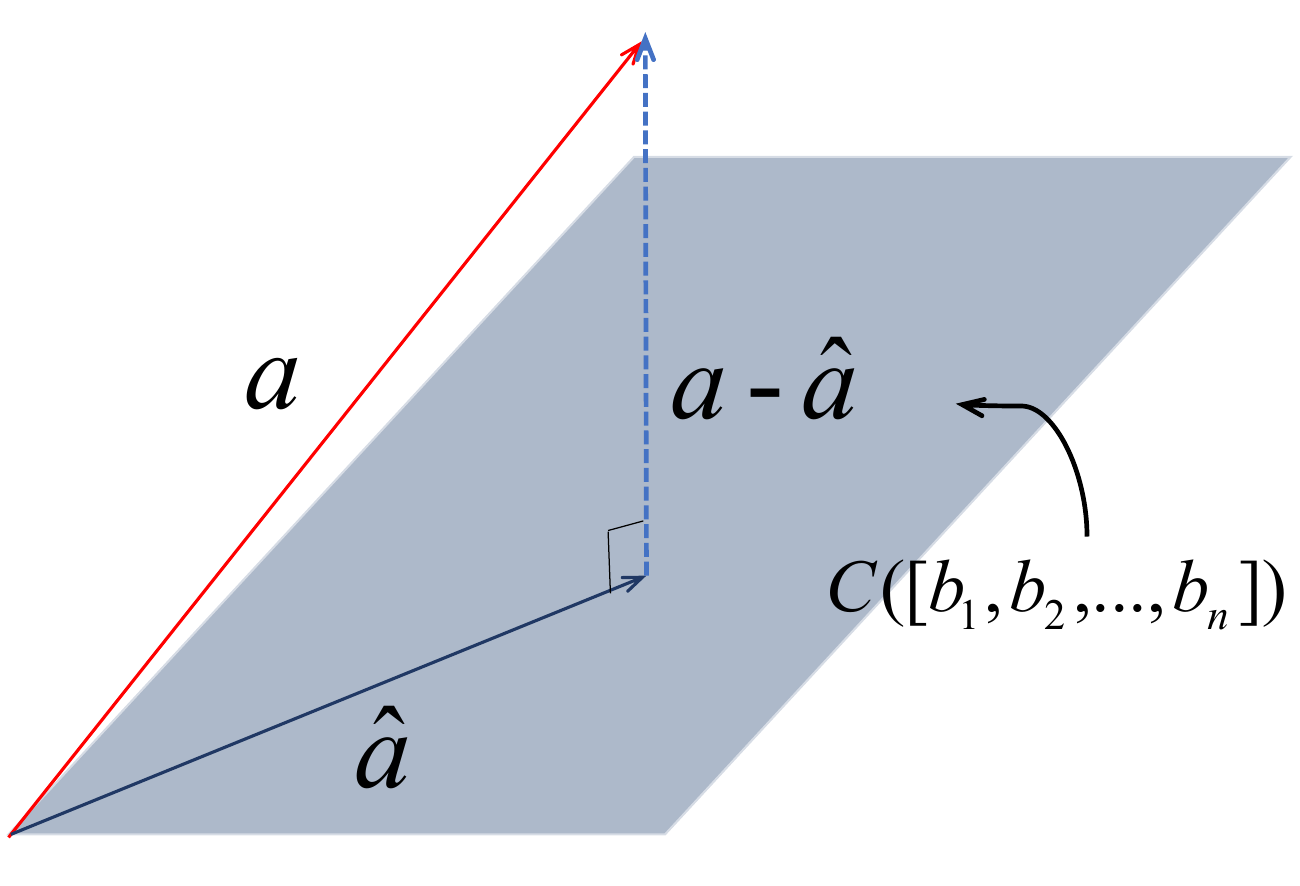}}
\caption{Project a vector onto a line and a space.}
\label{fig:projection-qr}
\end{figure}

The process of projecting  a vector $\ba$ onto a vector $\bb$ involves finding the vector closest to $\ba$ along the line determined by $\bb$. 
The projected vector, denoted as $\widehat{\ba}$, is a scalar multiple of $\bb$. 
Let $\widehat{\ba} \triangleq \widehat{x} \bb$, then, by construction, $\ba-\widehat{\ba}$ is perpendicular to $\bb$, as illustrated in Figure~\ref{fig:project-line}. 
Consequently, we obtain the following result:
\begin{tcolorbox}[title={Project Vector $\ba$ Onto Vector $\bb$}]
$\ba^\perp \triangleq\ba-\widehat{\ba}$ is perpendicular to $\bb$, so $(\ba-\widehat{x}\bb)^\top\bb=0$: $\widehat{x}$ = $\frac{\ba^\top\bb}{\bb^\top\bb}$ and $\widehat{\ba} = \frac{\ba^\top\bb}{\bb^\top\bb}\bb = \frac{\bb\bb^\top}{\bb^\top\bb}\ba$.~\footnote{The result can be generalized to a norm derived from an inner product (Definition~\ref{definition:inner_prod}, Exercise~\ref{exercise:norm_innr_pro}). Let $\norm{\cdot}$ be a norm on $\real^n$ that is derived from an inner product. Then, $\widehat{x}$ that minimizes the value $\norm{\ba-\widehat{x}\bb}$ is $\widehat{x}=\frac{\inner{\ba}{\bb}}{\norm{\bb}^2}$. When the inner product is the Euclidean inner product, it reduces to the main result shown above.}
\end{tcolorbox}

The above discussion leads to the following lemma.
\begin{lemma}[Finding an Orthogonal Vector]
Given two unit vectors $\bu$ and $\bv$ (i.e., $\normtwo{\bu}=\normtwo{\bv}=1$), then $\bw=\bv-\bu\bu^\top\bv$ is orthogonal to $\bu$.
If $\bu$ and $\bv$ are not unit, we can normalize them by setting $\bu:=\frac{\bu}{\normtwo{\bu}}$ and $\bv:=\frac{\bv}{\normtwo{\bv}}$ to achieve the same result.
\end{lemma}

\section{Projecting a Vector Onto a Plane}\label{section:project-onto-a-plane}
Similarly, the process of projecting  a vector $\ba$ onto a space spanned by vectors $\bb_1, \bb_2, \ldots, \bb_n$ involves finding the vector that is closest to $\ba$  within the column space of $[\bb_1, \bb_2, \ldots, \bb_n]$. 
The projection vector $\widehat{\ba}$ can be expressed as a linear combination of $\bb_1, \bb_2, \ldots, \bb_n$: $\widehat{\ba} \triangleq \widehat{x}_1\bb_1+ \widehat{x}_2\bb_2+\ldots+\widehat{x}_n\bb_n$. 
This problem can be formulated as a least squares problem. 
The projection can be obtained by solving the \textit{normal equation} $\bB^\top\bB\widehat{\bx} = \bB^\top\ba$, where $\bB=[\bb_1, \bb_2, \ldots, \bb_n]$ and $\widehat{\bx}=[\widehat{x}_1, \widehat{x}_2, \ldots, \widehat{x}_n]$. We refer the details of this projection view in the least squares to \citet{strang1993introduction, trefethen1997numerical, yang2000matrix, golub2013matrix, lu2021rigorous} as it is not the main interest of this book. For each vector $\bb_i$, the projection of $\ba$ in the direction of $\bb_i$ can be analogously obtained by 
$$
\widehat{\ba}_i = \frac{\bb_i\bb_i^\top}{\bb_i^\top\bb_i}\ba, \gap \forall i \in \{1,2,\ldots, n\}.
$$
Let $\widehat{\ba}\triangleq\sum_{i=1}^{n}\widehat{\ba}_i$, this results in 
$$
\ba^\perp \triangleq (\ba-\widehat{\ba}) \perp \cspace(\bB),
$$
i.e., $(\ba-\widehat{\ba})$ is perpendicular to the column space of $\bB=[\bb_1, \bb_2, \ldots, \bb_n]$, as illustrated in Figure~\ref{fig:project-space}.

\section{Existence of  QR Decomposition via  Gram-Schmidt Process}\label{section:gram-schmidt-process}
\subsubsection*{\textbf{First View by Projection Directly}}
Given three linearly independent vectors $\{\ba_1, \ba_2, \ba_3\}$ and the space spanned by the three linearly independent vectors $\cspace{([\ba_1, \ba_2, \ba_3])}$, i.e., the column space of the matrix $[\ba_1, \ba_2, \ba_3]$. We intend to construct three orthogonal vectors $\{\bb_1, \bb_2, \bb_3\}$ such that $\cspace{([\bb_1, \bb_2, \bb_3])}$ = $\cspace{([\ba_1, \ba_2, \ba_3])}$. Subsequently, we normalize the orthogonal vectors by dividing each by its length. This process produces three mutually orthonormal vectors $\bq_1 \triangleq \frac{\bb_1}{\normtwo{\bb_1}}$, $\bq_2 \triangleq \frac{\bb_2}{\normtwo{\bb_2}}$, and $\bq_2 \triangleq \frac{\bb_2}{\normtwo{\bb_2}}$.

\begin{figure}[H]
\centering  
\vspace{-0.35cm} 
\subfigtopskip=2pt 
\subfigbottomskip=2pt 
\subfigcapskip=-5pt 
\subfigure[Project $\ba_2$ onto the space perpendicular to $\bb_1$.]{\label{fig:gram-schmidt1}
	\includegraphics[width=0.4\linewidth]{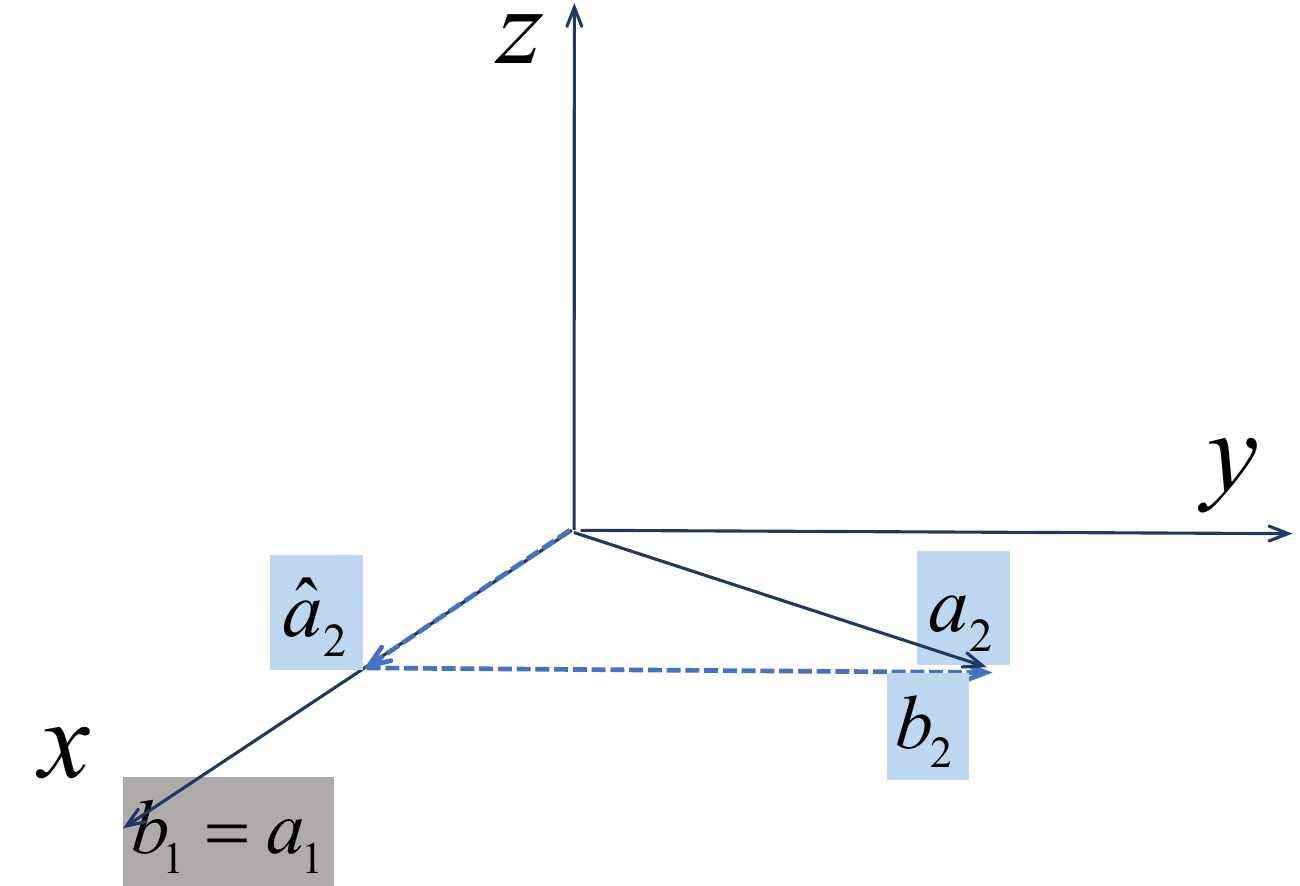}}
\quad 
\subfigure[Project $\ba_3$ onto the space perpendicular to $\bb_1, \bb_2$.]{\label{fig:gram-schmidt2}
	\includegraphics[width=0.4\linewidth]{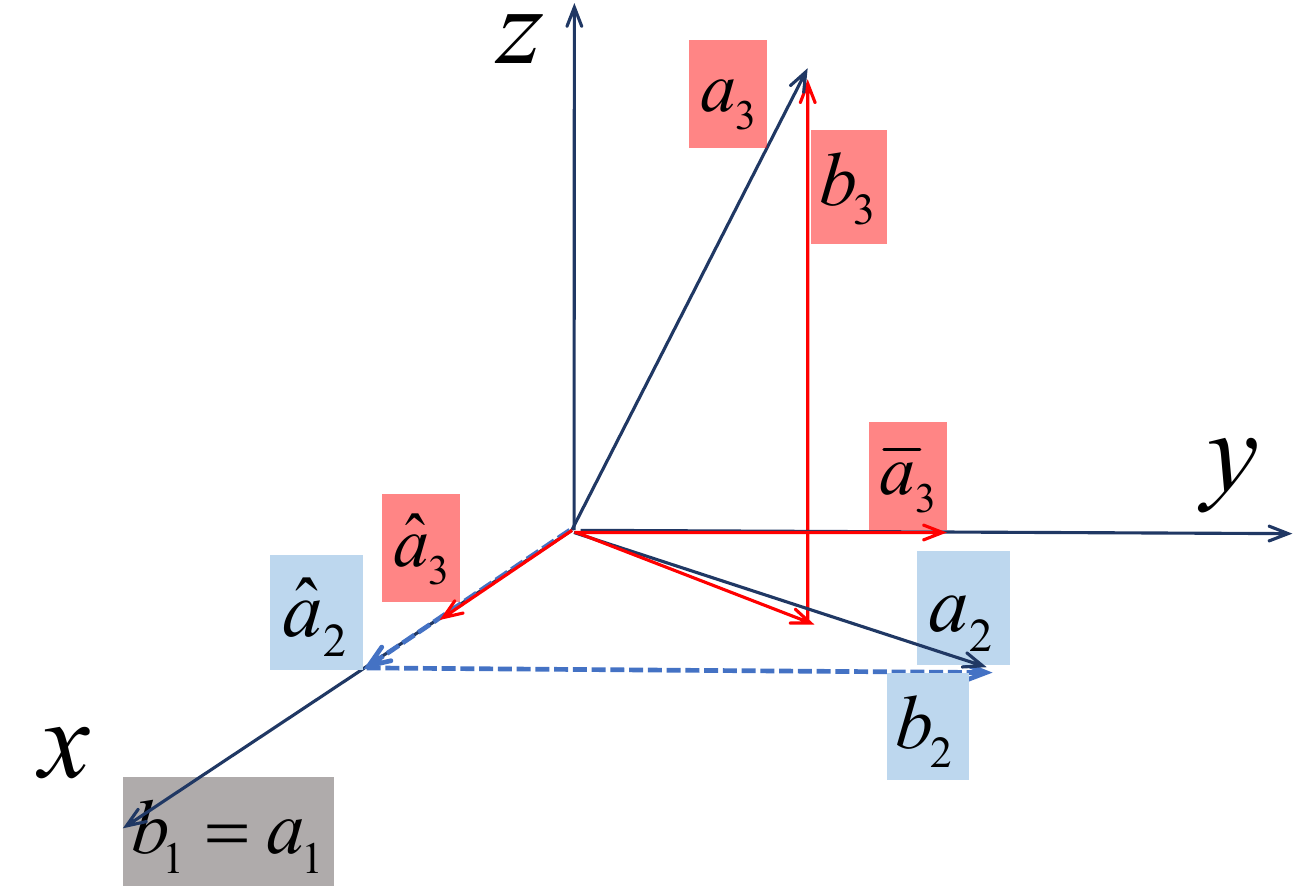}}
\caption{The Gram-Schmidt process.}
\label{fig:gram-schmidt-12}
\end{figure}

For the first vector, we simply set $\bb_1 = \ba_1$. The second vector $\bb_2$ must be perpendicular to the first one. This is achieved by considering the vector $\ba_2$ and subtracting its projection along $\bb_1$:
\begin{equation}
\begin{aligned}
\bb_2 &= \ba_2- \frac{\bb_1 \bb_1^\top}{\bb_1^\top\bb_1} \ba_2 = \left(\bI- \frac{\bb_1 \bb_1^\top}{\bb_1^\top\bb_1} \right)\ba_2   \qquad &(\text{Projection view})\\
&= \ba_2-  \underbrace{\frac{ \bb_1^\top \ba_2}{\bb_1^\top\bb_1} \bb_1}_{\widehat{\ba}_2}, \qquad &(\text{Combination view}) \nonumber
\end{aligned}
\end{equation}
where the first equation shows that $\bb_2$ results from the multiplication of the matrix $\big(\bI- \frac{\bb_1 \bb_1^\top}{\bb_1^\top\bb_1} \big)$ and the vector $\ba_2$, signifying the projection of $\ba_2$ onto the orthogonal complement space of $\cspace{([\bb_1])}$. The second equality in the above equation shows that $\ba_2$ is a linear combination of $\bb_1$ and $\bb_2$.
Clearly, the space spanned by $\bb_1, \bb_2$ coincides with the space spanned by $\ba_1, \ba_2$. The situation is illustrated in Figure~\ref{fig:gram-schmidt1}, where we designate \textbf{the direction of $\bb_1$ as the $x$-axis in the Cartesian coordinate system}, and $\widehat{\ba}_2$ represents the projection of $\ba_2$ onto the line defined by $\bb_1$. 
The figure clearly demonstrates that  the component of $\ba_2$ that is perpendicular to $\bb_1$ is $\bb_2 = \ba_2 - \widehat{\ba}_2$ (Pythagoras' theorem).

For the third vector $\bb_3$, it must be perpendicular to both  $\bb_1$ and $\bb_2$, which corresponds to the vector $\ba_3$ subtracting its projection along the plane spanned by $\bb_1$ and $\bb_2$:
\begin{equation}\label{equation:gram-schdt-eq2}
\begin{aligned}
\bb_3 &= \ba_3- \frac{\bb_1 \bb_1^\top}{\bb_1^\top\bb_1} \ba_3 - \frac{\bb_2 \bb_2^\top}{\bb_2^\top\bb_2} \ba_3 = \left(\bI- \frac{\bb_1 \bb_1^\top}{\bb_1^\top\bb_1}  - \frac{\bb_2 \bb_2^\top}{\bb_2^\top\bb_2} \right)\ba_3   \qquad &(\text{Projection view})\\
&= \ba_3- \underbrace{\frac{ \bb_1^\top\ba_3}{\bb_1^\top\bb_1} \bb_1}_{\widehat{\ba}_3} - \underbrace{\frac{ \bb_2^\top\ba_3}{\bb_2^\top\bb_2}  \bb_2}_{\bar{\ba}_3},    \qquad &(\text{Combination view})
\end{aligned}
\end{equation}
where the first equation shows that $\bb_3$ is a multiplication of the matrix $\left(\bI- \frac{\bb_1 \bb_1^\top}{\bb_1^\top\bb_1}  - \frac{\bb_2 \bb_2^\top}{\bb_2^\top\bb_2} \right)$ and the vector $\ba_3$, signifying the projection of $\ba_3$ onto the orthogonal complement space of $\cspace{([\bb_1, \bb_2])}$. The second equality in the above equation shows that $\ba_3$ is a linear combination of $\bb_1, \bb_2, $ and $\bb_3$. We will see this property is essential in the idea of the QR decomposition.
Once again, it can be shown that the space spanned by $\bb_1, \bb_2, \bb_3$ is identical to the space spanned by $\ba_1, \ba_2, \ba_3$. 
The scenario is illustrated in Figure~\ref{fig:gram-schmidt2}, where we designate \textbf{the direction of $\bb_2$ as the $y$-axis in the Cartesian coordinate system}. Here, $\widehat{\ba}_3$ is the projection of $\ba_3$ onto line $\bb_1$, while $\bar{\ba}_3$ represents the projection of $\ba_3$ onto line $\bb_2$. 
It can be shown that the component of $\ba_3$ perpendicular to both $\bb_1$ and $\bb_2$ is $\bb_3=\ba_3-\widehat{\ba}_3-\bar{\ba}_3$ from the figure.

Finally, we normalize each vector by dividing its  length, resulting in three orthonormal vectors $\bq_1 \triangleq \frac{\bb_1}{\normtwo{\bb_1}}$, $\bq_2 \triangleq \frac{\bb_2}{\normtwo{\bb_2}}$, and $\bq_2 \triangleq \frac{\bb_2}{\normtwo{\bb_2}}$.

This procedure can be extended to a set of vectors, not limited to just three, and is known as the \textit{Gram-Schmidt process}. After applying this process, the matrix $\bA$ will be triangularized. The method is named after \textit{Jørgen Pedersen Gram} and \textit{Erhard Schmidt}, but it appeared earlier in the work of \textit{Pierre-Simon Laplace} within the theory of Lie group decomposition.

As mentioned earlier, the concept behind QR decomposition involves constructing a sequence of orthonormal vectors $\bq_1, \bq_2, \ldots$, which collectively span the same successive subspaces:
$$
\big\{\cspace([\bq_1])=\cspace([\ba_1])\big\} \subseteq
\big\{\cspace([\bq_1, \bq_2])=\cspace([\ba_1, \ba_2])\big\} \subseteq
\big\{\cspace([\bq_1, \bq_2, \bq_3])=\cspace([\ba_1, \ba_2, \ba_3])\big\} 
\subseteq \ldots.
$$
This implies that  any vector $\ba_k$ resides in the space spanned by $\cspace([\bq_1, \bq_2, \ldots, \bq_k])$.~\footnote{Moreover, any vector $\bq_k$ resides in the space spanned by $\cspace([\ba_1, \ba_2, \ldots, \ba_k])$.} Once these orthonormal vectors are identified, the reconstruction of $\ba_i$'s from the orthogonal matrix $\bQ=[\bq_1, \bq_2, \ldots, \bq_n]$ requires an upper triangular matrix $\bR$  such that $\bA = \bQ\bR$.

\index{Gram–Schmidt}
While the Gram-Schmidt process is a well-known algorithm for QR decomposition, it is not the only one. Several alternative QR decomposition algorithms, such as \textit{Householder reflections} and \textit{Givens rotations}, are available. These methods are often more robust in the presence of round-off errors and may alter the order in which the columns of $\bA$ are processed.

\subsubsection*{\textbf{Another View by Inner Product (with Projection Implicitly)}}
In the direct projection view mentioned above, we first find orthogonal vectors and then normalize them to have a unit length. 
This projection is based on the result derived in Section~\ref{section:project-onto-a-vector}, which states that the projection of the vector $\ba$ onto the vector $\bb$ is given by $\widehat{\ba} = \frac{\ba^\top\bb}{\bb^\top\bb}\bb = \frac{\bb\bb^\top}{\bb^\top\bb}\ba$. 
If $\bb$ has a unit length, this simplifies to:
\begin{equation}\label{equation:qr-projection-unit}
\widehat{\ba} = (\ba^\top\bb) \bb.
\end{equation}

Once again, considering three linearly independent vectors $\{\ba_1, \ba_2, \ba_3\}$ and the space spanned by these vectors, denoted by $\cspace{([\ba_1, \ba_2, \ba_3])}$. 
We will employ a slight abuse of notation in this context, where we will replace $\bb_k$ in the first view with $\ba_k^\perp$  to emphasize that $\ba_k^\perp$ is orthogonal to $\{\bq_1, \ldots, \bq_{k-1}\}$. The process proceeds as follows:
\begin{itemize}
\item Compute $\bq_1$ of unit length so that $\cspace([\bq_1]) = \cspace([\ba_1])$:
\begin{itemize}
\item Compute the length of $\ba_1$: $r_{11} = \normtwo{\ba_1}$;
\item Sets $\bq_1$ to a unit vector in the direction of $\ba_1$: $\bq_1=\frac{\ba_1}{\normtwo{\ba_1}}=\frac{\ba_1}{r_{11}}$;
\end{itemize}
\item Compute $\bq_2$ of unit length so that $\cspace([\bq_1,\bq_2]) = \cspace([\ba_1,\ba_2])$:
\begin{itemize}
\item Compute $r_{12}$ so that $r_{12}\bq_1 = (\ba_2^\top\bq_1) \bq_1$ equals the component of $\ba_2$ in the direction of $\bq_1$ by Equation~\eqref{equation:qr-projection-unit};
\item Compute the component of $\ba_2$ that is orthogonal to $\bq_1$: $\ba_2^\perp = \ba_2 - r_{12}\bq_1$;
\item Compute the length of vector $\ba_2^\perp$: $r_{22}=\normtwo{\ba_2^\perp}$;
\item Set $\bq_2$ to a unit vector in the direction of ${\ba_2^\perp}$: $\bq_2 = \frac{\ba_2^\perp}{\normtwo{\ba_2^\perp}}=\frac{\ba_2^\perp}{r_{22}}$;
\item This results in:
$$
\begin{bmatrix}
	\ba_1 & \ba_2 
\end{bmatrix}=
\begin{bmatrix}
\bq_1 & \bq_2 
\end{bmatrix}
\begin{bmatrix}
r_{11}& r_{12}\\
0 & r_{22} 
\end{bmatrix};
$$
\end{itemize}
\item Compute $\bq_3$ of unit length so that $\cspace([\bq_1,\bq_2,\bq_3]) = \cspace([\ba_1,\ba_2,\ba_3])$:
\begin{itemize}
\item Compute $r_{13}$ so that $r_{13}\bq_1 = (\ba_3^\top\bq_1) \bq_1$ equals the component of $\ba_3$ in the direction of $\bq_1$ by Equation~\eqref{equation:qr-projection-unit};
\item Compute $r_{23}$ so that $r_{23}\bq_2 = (\ba_3^\top\bq_2) \bq_2$ equals the component of $\ba_3$ in the direction of $\bq_2$ by Equation~\eqref{equation:qr-projection-unit};
\item Compute the component of $\ba_3$ that is orthogonal to $\bq_1$ and $\bq_2$: $\ba_3^\perp = \ba_3 - r_{13}\bq_1 - r_{23}\bq_2$;
\item Compute the length of vector $\ba_3^\perp$: $r_{33}=\normtwo{\ba_3^\perp}$;
\item Set $\bq_3$ to a unit vector in the direction of ${\ba_3^\perp}$: $\bq_3 = \frac{\ba_3^\perp}{\normtwo{\ba_3^\perp}}=\frac{\ba_3^\perp}{r_{33}}$;
\item This results in:
$$
\begin{bmatrix}
\ba_1 & \ba_2 & \ba_3  
\end{bmatrix}=
\begin{bmatrix}
\bq_1 & \bq_2 & \bq_3 
\end{bmatrix}
\begin{bmatrix}
r_{11}& r_{12} & r_{13}\\
0 & r_{22}  & r_{23} \\
0 & 0 & r_{33}
\end{bmatrix}.
$$
\end{itemize}
\end{itemize}
Again, this idea can be extended to a set of vectors rather than only three. The process described above reveals the interpretation of $r_{ij}$ in the triangular matrix $\bR$, where each $r_{ij}$ represents the component of $\ba_j$ in the direction of $\bq_i$. 
This aligns with the outcome of matrix multiplication:
$
\ba_j = \sum_{i=1}^{j} r_{ij} \bq_j.
$

\subsubsection*{\textbf{Main Proof}}
Although the existence of the QR decomposition is conceptually intuitive based on the two views of the Gram-Schmidt process discussed above, a formal proof is clunky and requires the use of an inductive hypothesis. We will now provide a rigorous proof.

\begin{proof}[{of Theorem~\ref{theorem:qr-decomposition}: QR Decomposition}]
We will prove by induction that every $m\times n$ matrix $\bA$ with linearly independent columns admits a \textit{reduced} QR decomposition (the dependent case will be discussed in the sequel). The \textit{full} QR decomposition can be obtained by completing the orthonormal columns in $\bQ$. The case for a $1\times 1$ matrix is trivial by setting $Q=1, R=A$. Thus, $A=QR$. 

Assume that for any $m\times k$ matrix $\bA_k$ with linearly independent columns, the reduced QR decomposition exists. 
To complete the proof, we need to show that any $m\times (k+1)$ matrix $\bA_{k+1}$ can also be factorized using this reduced QR decomposition. 
Suppose $\bA_{k+1}=[\bA_k, \ba_{k+1}]$, where $\bA_k$ admits the reduced QR decomposition by inductive hypothesis
$$
\bA_k = \bQ_k\bR_k,
$$ 
where $\bQ_k$ consists of orthonormal columns $\bQ_k^\top\bQ_k=\bI_k$ and $\bR_k \in \real^{k\times k}$ is upper triangular. 
Moreover, according to the induction hypothesis,
if the values on the diagonal of $\bR_k$ are chosen to be positive, then the reduced QR decomposition of
$\bA_k = \bQ_k\bR_k$ is \textbf{unique}.

Suppose further that $\bA_{k+1}$ can be factored as
\begin{equation}\label{equation:qr-rigorous-proof}
\bA_{k+1} = 
\begin{bmatrix}
	\bA_k& \ba_{k+1}
\end{bmatrix}
=
\begin{bmatrix}
	\widetildebQ_k& \bq_{k+1}
\end{bmatrix}
\begin{bmatrix}
	\widetildebR_k& 
	\begin{matrix}
		\br_{k}\\
		r_{k+1}
	\end{matrix}
\end{bmatrix}
,
\end{equation}
where it is apparent that $\bA_k = \widetildebQ_k \widetildebR_k$ is a reduced QR decomposition of $\bA_k$.
Furthermore, if restrict the diagonal values of $\widetildebR_k$ to be positive, the factorization is \textbf{unique} implying that  $\widetildebQ_k = \bQ_k$ and $\widetildebR_k=\bR_k$. Moreover, Equation~\eqref{equation:qr-rigorous-proof} implies 
$$
\begin{aligned}
\ba_{k+1} = \bQ_k\br_k + \bq_{k+1}r_{k+1} \\
\implies  &\bQ_k^\top \ba_{k+1} = \bQ_k^\top ( \bQ_k\br_k + \bq_{k+1}r_{k+1})=\br_k + \bQ_k^\top\bq_{k+1}r_{k+1}.
\end{aligned}
$$
Since we assume that the columns of $\bQ_k$ are orthonormal to $\bq_{k+1}$, it follows that $\bQ_k^\top\bq_{k+1}=\bzero$ and $\br_k=\bQ_k^\top \ba_{k+1}$. When $\bQ_k$ is fixed, the value of $\br_k$ is \textbf{uniquely} determined.

Let $\ba_{k+1}^\perp = \ba_{k+1}-\bQ_k\br_k$, then $\ba_{k+1}^\perp$ is orthogonal to the columns of $\bQ_k$. To see this, $\bQ_k^\top \ba_{k+1}^\perp = \bQ_k^\top(\ba_{k+1}-\bQ_k\br_k) =\bzero$ as we construct $\br_k$ by $\br_k=\bQ_k^\top \ba_{k+1}$. 
Since $\ba_{k+1}$ is linearly independent from the columns of $\bA_k$, so it is also linearly independent from the columns of $\bQ_k$, $\ba_{k+1}^\perp$ is thus nonzero. Therefore, let $r_{k+1}=\normtwo{\ba_{k+1}^\perp}$ and $\bq_{k+1} = \ba_{k+1}^\perp / r_{k+1}$, we find the \textbf{unique} reduced QR decomposition of $\bA_{k+1}$ with positive diagonals in the upper triangular matrix. This completes the proof.
\end{proof}

\index{Semi-orthogonal}
\section{Orthogonal vs Orthonormal}\index{Orthogonal}\index{Orthonormal}\label{section:orthogonal-orthonormal-qr}
Vectors $\bq_1, \bq_2, \ldots, \bq_n\in \real^m$ are \textit{mutually orthogonal} if their dot products $\bq_i^\top\bq_j$ are zero whenever $i \neq j$. 
Upon normalizing each vector by its length, the resulting vectors become \textit{orthogonal unit vectors}. Then the vectors $\bq_1, \bq_2, \ldots, \bq_n$ are considered \textit{mutually orthonormal}, and these orthonormal vectors are typically arranged into a matrix $\bQ$:
\begin{itemize}
\item When $m\neq n$: the matrix $\bQ$ is easy to work with because $\bQ^\top\bQ=\bI \in \real^{n\times n}$. Such a matrix $\bQ$ with $m\neq n$ is sometimes referred to as a \textit{semi-orthogonal} matrix.

\item When $m= n$: the matrix $\bQ$ is square, the condition $\bQ^\top\bQ=\bI$ implies that $\bQ^\top=\bQ^{-1}$, i.e., the transpose of $\bQ$ serves as its inverse. 
Then we also have $\bQ\bQ^\top=\bI$, i.e., $\bQ^\top$ serves as the \textit{two-sided inverse} of $\bQ$. This type of $\bQ$ is referred to as an \textit{orthogonal matrix}. \footnote{Note that we use the term \textit{orthogonal matrix} to refer to the matrix $\bQ$ with orthonormal columns. The term \textit{orthonormal matrix} is \textbf{not} used for historical reasons.}
\end{itemize}
To see this, we have
$$
\begin{bmatrix}
\bq_1^\top \\
\bq_2^\top\\
\vdots \\
\bq_n^\top
\end{bmatrix}
\begin{bmatrix}
\bq_1 &\bq_2 & \ldots & \bq_n
\end{bmatrix}
= 
\begin{bmatrix}
1 & & & \\
& 1 & & \\
& & \ddots & \\
& & & 1
\end{bmatrix}.
$$
In other words, $\bq_i^\top \bq_j = \delta_{ij}$, where $\delta_{ij}$ represents the \textit{Kronecker delta}. The columns of an orthogonal matrix $\bQ\in \real^{n\times n}$ form an \textit{orthonormal basis} for $\real^n$. \footnote{Notice again that the \textit{orthogonal matrix} $\bQ$ contains an \textit{orthonormal basis}, \textbf{not} an orthogonal basis.}

Orthogonal matrices can be regarded as matrices that transform the basis of other matrices. Hence they preserve the angle (inner product) between  vectors:
$$
\text{(inner product): } \qquad \bu^\top \bv = (\bQ\bu)^\top(\bQ\bv).
$$
The above invariance of the inner products of angles between the vectors also relies on the invariance of their lengths:
$$
\text{(length): } \qquad \normtwo{\bQ\bu} = \normtwo{\bu}.
$$
In real-valued cases, when multiplied by an orthogonal matrix $\bQ$, the original vector space undergoes \textit{rotation} (if $\det(\bQ) = 1$) or \textit{reflection} (if $\det(\bQ) = -1$). 
Many decomposition algorithms will result in two orthogonal matrices, leading to such rotations or reflections occurring twice. 
See Section~\ref{section:coordinate-transformation} for a discussion on  coordinate transformations in matrix decomposition.

\begin{exercise}[Orthogonal]
Prove that the following statements about a square matrix $\bQ\in\real^{n\times n}$ are equivalent:
\begin{itemize}
\item $\bQ$ is orthogonal.
\item $\bQ^\top$ is orthogonal.
\item $\bQ\bQ^\top=\bQ^\top\bQ=\bI$.
\item $\bQ$ is nonsingular and $\bQ^\top=\bQ^{-1}$.
\item The rows of $\bQ$ are orthogonormal.
\item The columns of $\bQ$ are orthonormal.
\item For all $\bx\in\real^n$, it follows that $\normtwo{\bx}=\normtwo{\bQ\bx}$.
\end{itemize}
\paragraph{Unitary.$^*$} Let $\bU\in\complex^{n\times n}$ be unitary. Show that $\overline{\bU}$ , $\overline{\bU}^\top$ , and $\overline{\bU}^*$ are all unitary.
\end{exercise}
Further properties related to the preservation of orthogonal transformations are explored in Problems~\ref{prob:ortho_prese1}$\sim$\ref{prob:ortho_presen}.

\begin{figure}[H]
\centering   
\vspace{-0.35cm}  
\subfigtopskip=2pt  
\subfigbottomskip=2pt  
\subfigcapskip=-5pt  
\subfigure[$\bQ_1$ rotates $\bv$ along the $\bzero$ point.]{\label{fig:orthogonal-rotate1}
	\includegraphics[width=0.4\linewidth]{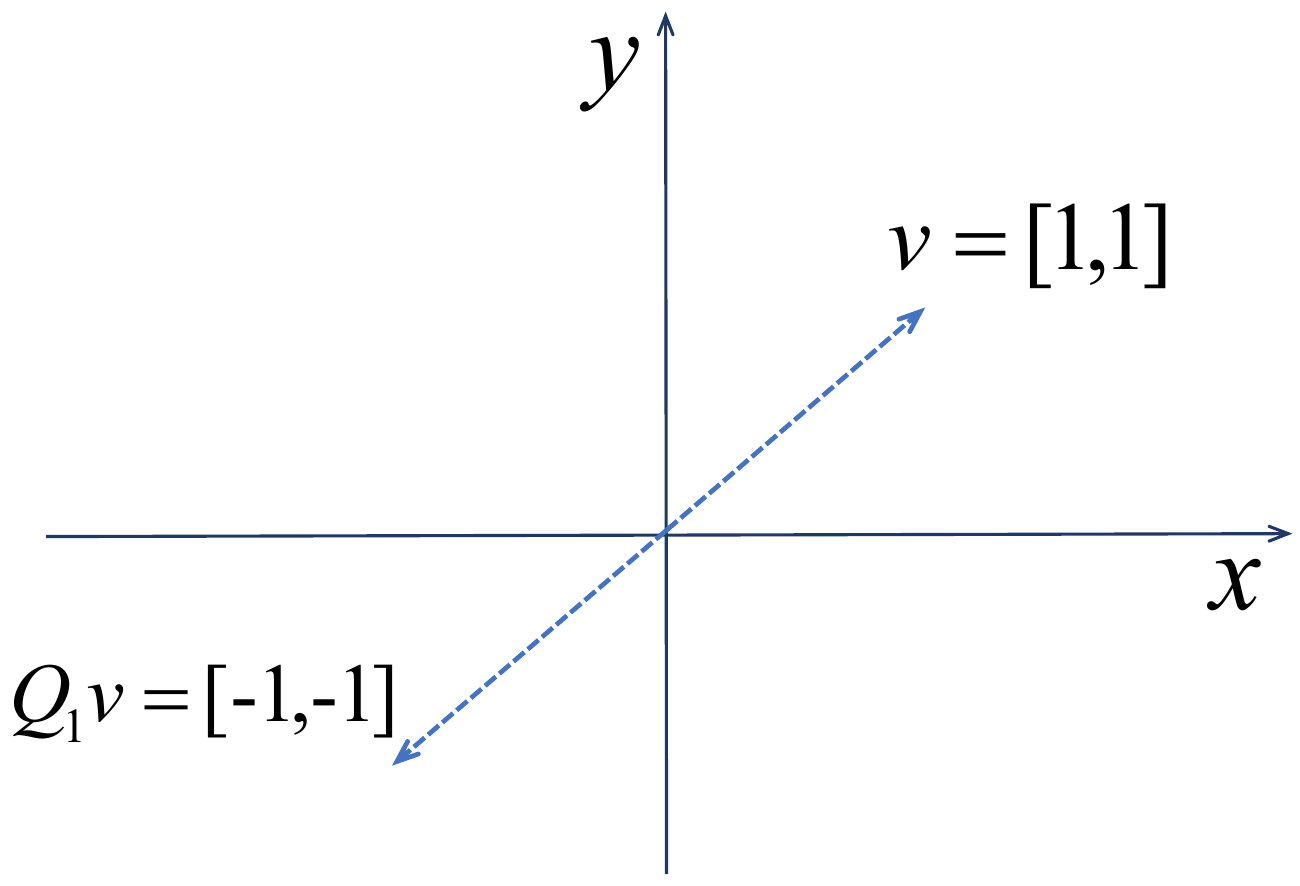}}
\quad 
\subfigure[$\bQ_2$ reflects $\bv$ across the $x$-axis.]{\label{fig:orthogonal-rotate2}
	\includegraphics[width=0.4\linewidth]{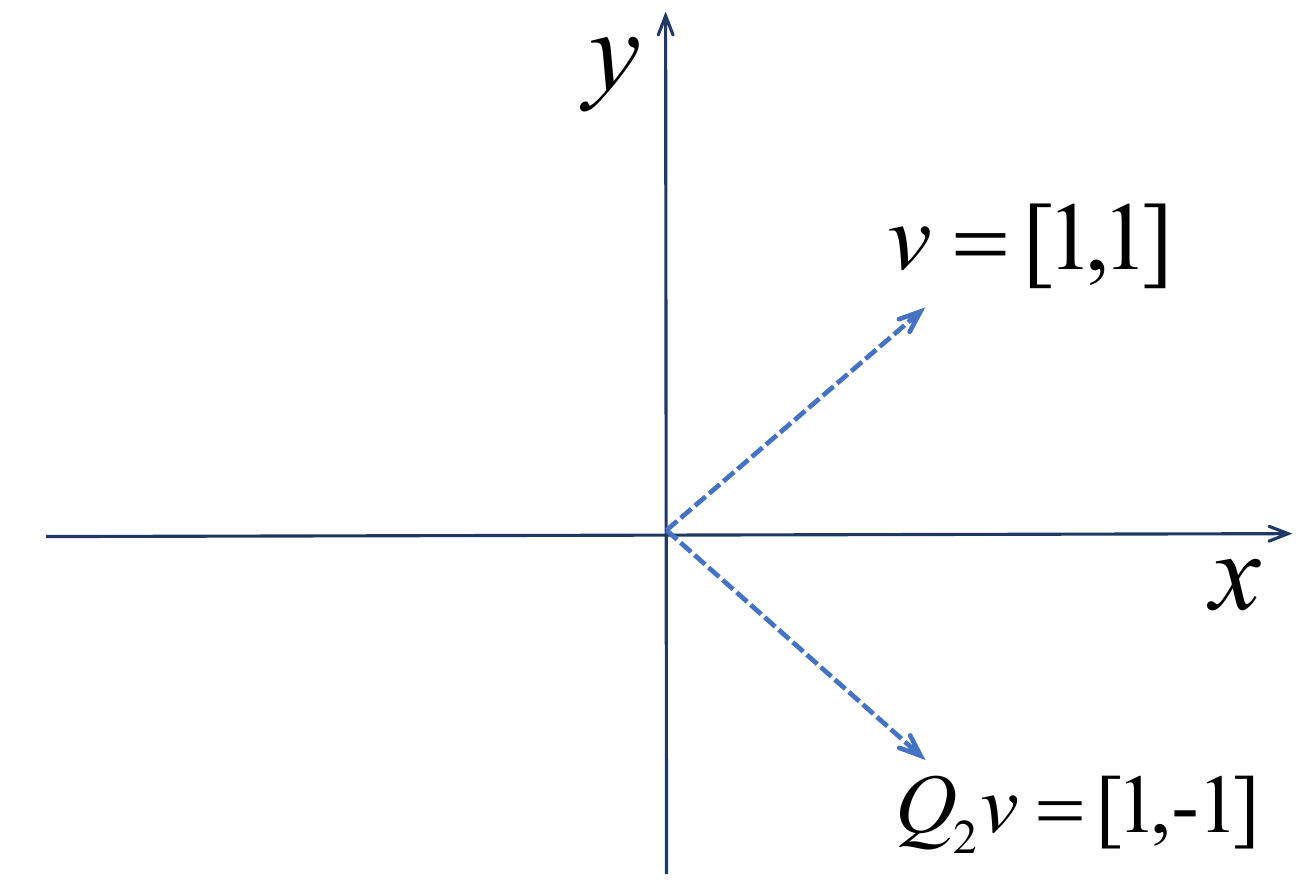}}
\caption{Rotation and reflection in orthogonal matrices.}
\label{fig:orthogonal-rotate-reflect}
\end{figure}

\begin{example}[Rotation and Reflection in Orthogonal Matrices]
To illustrate rotation and reflection in orthogonal matrices, suppose 
$$
\bQ_1 = \begin{bmatrix}
-1 & \\
& -1
\end{bmatrix}
\qquad 
\text{and}
\qquad 
\bQ_2 = \begin{bmatrix}
1 & \\
& -1
\end{bmatrix},
$$
where $\det(\bQ_1)=1$ and $\det(\bQ_2)=-1$. For the vector $\bv=[1,1]^\top$, we have 
$$
\bQ_1\bv = \begin{bmatrix}
-1 \\
-1
\end{bmatrix}
\qquad 
\text{and}
\qquad 
\bQ_2\bv = \begin{bmatrix}
1 \\
-1
\end{bmatrix}.
$$
Thus, $\bQ_1$ rotates $\bv$ along the point $\bzero$, and $\bQ_2$ reflects $\bv$ across the $x$-axis.
The illustration of the rotation and reflection is shown in Figure~\ref{fig:orthogonal-rotate-reflect}.
\end{example}

\section{Properties of  QR Decomposition}
\subsection{Orthonormal Basis via QR Decomposition}
For any matrix $\bA$, we have the property: $\nspace(\bA^\top)$ is the orthogonal complement of the column space $\cspace(\bA)$ in $\real^m$: $\dim(\nspace(\bA^\top))+\dim(\cspace(\bA))=m$.
This is known as the \textit{rank-nullity theorem}. And a proof of this theorem can be found in Appendix~\ref{appendix:fundamental-rank-nullity}.
The QR decomposition allows us to find a basis for the corresponding subspaces.  
In singular value decomposition (SVD, Chapter~\ref{chapter:SVD}), we will also find  orthonormal bases for $\nspace(\bA)$ and $\cspace(\bA^\top)$.

\begin{lemma}[Orthonormal Basis in $\real^m$]\label{lemma:qr-four-orthonormal-Basis}
For the full QR decomposition of a matrix $\bA\in \real^{m\times n}$ with full rank $n$ and $m\geq n$, the following properties hold:
\begin{itemize}
\item $\{\bq_1,\bq_2, \ldots,\bq_n\}$ is an orthonormal basis for $\cspace(\bA)$;
\item $\{\bq_{n+1},\bq_{n+2}, \ldots,\bq_m\}$ is an orthonormal basis for $\nspace(\bA^\top)$. 
\end{itemize}
\end{lemma}
\begin{proof}[of Lemma~\ref{lemma:qr-four-orthonormal-Basis}]
According to the Gram-Schmidt process, it can be easily shown that $\spn\{\ba_1, \ba_2, \ldots, \ba_k\}$~\footnote{Note that we use the notation such that $\spn\{\ba_1, \ba_2, \ldots\}$ is equal to the column space of the matrix column space $\cspace([\ba_1, \ba_2, \ldots])$.}
is equal to $\spn\{\bq_1, \bq_2, \ldots, \bq_k\}$ for all $k\in\{1, 2, \ldots, n\}$. 
Therefore, we can conclude that $\cspace(\bA) = \spn\{\ba_1, \ba_2, \ldots, \ba_n\}$ =$\spn\{\bq_1, \bq_2, \ldots, \bq_n\}$, and $\{\bq_1, \bq_2, \ldots, \bq_n\}$ serves as an orthonormal basis for the column space of $\bA$. 
Additionally, since $\nspace(\bA^\top)$ and $\cspace(\bA)$ are orthogonal complement subspaces, denoted by $ \nspace(\bA^\top) \bot\cspace(\bA)$, 
we then have 
$\dim(\nspace(\bA^\top))$$=m-\dim(\cspace(\bA))=m-n$. 
And the space spanned by $\{\bq_{n+1},\bq_{n+2}, \ldots,\bq_m\}$ is also $\bot\cspace(\bA)$ with a dimension of $m-n$. Thus, $\{\bq_{n+1},\bq_{n+2}, \ldots,\bq_m\}$ forms an orthonormal basis for $\nspace(\bA^\top)$.
\end{proof}

\subsection{Hadamard Inequalities}\label{section:hadama_ineq}
The QR decomposition reveals the famous determinantal inequality by \textit{Jacques Hadamard} in 1893.
\begin{theorem}[Hadamard Inequality]\label{theorem:hadamard_ineq1}
Let $\bA\in\real^{n\times n}$ be any matrix. Then, 
$$
|\det(\bA)| \leq \prod_{i=1}^{n}\Vert\ba_i\Vert = \prod_{i=1}^{n}\left(\sum_{j=1}^{n}a_{ij}^2\right)^{1/2}.
$$
The equality holds when $\bA$ has orthogonal columns.
Applying the results to $\bA^\top$, we also obtain 
$$
|\det(\bA)| \leq  \prod_{j=1}^{n}\left(\sum_{i=1}^{n}a_{ij}^2\right)^{1/2}.
$$
\end{theorem}
\begin{proof}[of Theorem~\ref{theorem:hadamard_ineq1}]
If matrix  $\bA$ is singular, then $\det(\bA)=0$, and the result is trivial. If $\bA$ is nonsingular, then $\bA$ admits the unique QR decomposition $\bA=\bQ\bR$, where $\bR$ is upper triangular with positive diagonals. Since $\det(\bQ)=\pm1$, we obtain 
$$
\begin{aligned}
|\det(\bA)|  &= \vert\det(\bQ)\vert \cdot  \vert\det(\bR)\vert = \prod_{j=1}^n r_{jj};\\
\Vert\ba_j\Vert &= \Vert\bQ\br_j\Vert = \Vert\br_j\Vert = \sqrt{\sum_{i=1}^{n} r_{ij}^2} \geq r_{jj}, \gap
j\in\{1,2,\ldots, n\}.
\end{aligned}
$$
Therefore, 
$$
|\det(\bA)| = \prod_{j=1}^n r_{jj}
\leq \prod_{j=1}^n\Vert\ba_j\Vert.
$$
This completes the proof.
\end{proof}

The Hadamard inequality holds similarly for positive semidefinite matrices.
\begin{theorem}[Hadamard Inequality for PSD]\label{theorem:hadamard_ineq_psd}
Let  $\bA\in\real^{n\times n}$ be any  positive semidefinite matrix. Then,
$$
|\det(\bA)| \leq \prod_{i=1}^{n}a_{ii}.
$$
The equality holds when $\bA$ is a diagonal matrix.
\end{theorem}
\begin{proof}[of Theorem~\ref{theorem:hadamard_ineq_psd}]
If $\det(\bA)=0$, the inequality is trivial. 
If $\det(\bA)\neq 0$, then $\bA$ is positive definite and admits the unique Cholesky decomposition $\bA = \bR^\top\bR$, where $\bR$ is upper triangular with positive diagonals. Then we have 
$$
\det(\bA) = \det(\bR^\top\bR) =\det(\bR)^2.
$$
By Theorem~\ref{theorem:hadamard_ineq1}, we have 
$$
\begin{aligned}
\det(\bR) \leq \prod_{j=1}^{n}\left(\sum_{i=1}^{n}r_{ij}^2\right)^{1/2}
=\prod_{j=1}^{n} (a_{jj})^{1/2}
\implies
|\det(\bA)| =\det(\bR)^2
\leq \prod_{j=1}^{n} a_{jj}.
\end{aligned}
$$
The equality holds when $\bR$ has orthogonal columns, i.e., $\bA$ is a diagonal matrix.
\end{proof}

\subsection{Cholesky for QR}
We have demonstrated the existence of the QR decomposition for rectangular matrices. When $\bA$ is a nonsingular square matrix in $\real^{n\times n}$, it is evident that $\bA^\top\bA$ is symmetric positive definite.
Then, the matrix $\bA^\top\bA$ admits the Cholesky decomposition: $\bA^\top\bA=\bR^\top\bR$, where $\bR$ is nonsingular and upper triangular (Theorem~\ref{theorem:cholesky-factor-exist}).
Furthermore, we have $\bQ=\bA\bR^{-1}$ as an orthogonal matrix, such that $\bA=\bQ\bR$ represents the QR decomposition of $\bA$.
\begin{proof}
Since $\bA^\top\bA$ is nonsingular, by Theorem~\ref{theorem:cholesky-factor-exist}, $\bR$ is also nonsingular.
We then have 
$$
\bQ\bQ^\top = \bA\bR^{-1} \bR^{-\top}\bA^\top
=\bA (\bR^\top\bR)^{-1} \bA^\top
= \bA (\bA^\top\bA)^{-1} \bA^\top
=\bA\bA^{-1}\bA^{-\top} \bA^\top
=\bI.
$$
Therefore, $\bQ$ is orthogonal.
\end{proof}

\section{Computing  Reduced QR Decomposition via  Gram-Schmidt Process}\label{section:qr-gram-compute}

Expressing the reduced QR decomposition in the form $\bA = \bQ\bR$, where $\bQ\in \real^{m\times n}$ and $\bR\in \real^{n\times n}$, we have:
\begin{equation}
\bA=\left[
\begin{matrix}
\ba_1 & \ba_2 & \ldots & \ba_n
\end{matrix}
\right] 
=\left[
\begin{matrix}
\bq_1 & \bq_2 & \ldots & \bq_n
\end{matrix}
\right] 
\begin{bmatrix}
r_{11} & r_{12}& \dots & r_{1n}\\
    & r_{22}& \dots & r_{2n}\\
    &       &    \ddots  & \vdots \\
\multicolumn{2}{c}{\raisebox{1.3ex}[0pt]{\Huge0}} & & r_{nn} \nonumber
\end{bmatrix}.
\end{equation}
The orthogonal matrix $\bQ$ can be easily calculated using the Gram-Schmidt process. 
To understand the origin of the upper triangular matrix $\bR$, let us explicitly express these equations:
\begin{equation*}
\begin{aligned}
\ba_1 & = r_{11}\bq_1 &= \sum_{i=1}^{1} r_{i1}\bq_1, \\
\ba_2 & = r_{12}\bq_1 + r_{22}\bq_2&= \sum_{i=1}^{2} r_{i2}\bq_2, \\
&\vdots& \\
\ba_k &= r_{1k}\bq_1 + r_{2k}\bq_2 + \ldots + r_{kk}\bq_k  &= \sum_{i=1}^{k} r_{ik} \bq_k,\\ 
&\vdots& \\
\end{aligned}
\end{equation*}
which corresponds to the second equation of Equation~\eqref{equation:gram-schdt-eq2} and conforms to the structure of an upper triangular matrix $\bR$. Extending the concept from  Equation~\eqref{equation:gram-schdt-eq2} to the $k$-th term, we obtain
$$
\begin{aligned}
\ba_k &= \sum_{i=1}^{k-1}(\bq_i^\top\ba_k)\bq_i + \ba_k^\perp 
= \sum_{i=1}^{k-1}(\bq_i^\top\ba_k)\bq_i + \normtwo{\ba_k^\perp}\cdot \bq_k, 
\end{aligned}
$$
indicating that we can gradually orthonormalize $\bA$ to obtain an orthonormal set $\bQ=[\bq_1, \bq_2, \ldots, \bq_n]$ by 
\begin{equation}\label{equation:qr-gsp-equation}
\left\{
\begin{aligned}
r_{ik} &= \bq_i^\top\ba_k, \,\,\,\,\forall i \in \{1,2,\ldots, k-1\};\\ 
\ba_k^\perp&\triangleq \ba_k-\sum_{i=1}^{k-1}r_{ik}\bq_i;\\
r_{kk} &= \normtwo{\ba_k^\perp};\\
\bq_k &= \ba_k^\perp/r_{kk}.
\end{aligned}
\right.
\end{equation}
The procedure is formulated in Algorithm~\ref{alg:reduced-qr}.
\begin{algorithm}[h] 
\caption{Reduced QR Decomposition via Gram-Schmidt Process} 
\label{alg:reduced-qr} 
\begin{algorithmic}[1] 
\Require Matrix $\bA$ has linearly independent columns with size $m\times n $ and $m\geq n$; 
\For{$k=1$ to $n$} \Comment{compute the $k$-th column of $\bQ,\bR$}
\For{$i=1$ to $k-1$}
\State $r_{ik} \leftarrow\bq_i^\top\ba_k$; \Comment{entry ($i,k$) of $\bR$, $2m-1$ flops}
\EndFor  \Comment{all $k-1$ iterations: $(k-1)(2m-1)$ flops}
\State $\ba_k^\perp\leftarrow  \ba_k-\sum_{i=1}^{k-1}r_{ik}\bq_i$;\Comment{$2m(k-1)$ flops}
\State $r_{kk} \leftarrow \normtwo{\ba_k^\perp}$; \Comment{main diagonal of $\bR$, $2m$ flops}
\State $\bq_k \leftarrow \ba_k^\perp/r_{kk}$; \Comment{$m$ flops}
\EndFor
\State Output $\bQ=[\bq_1, \ldots, \bq_n]$ and $\bR$ with entry $(i,k)$ being $r_{ik}$.
\end{algorithmic} 
\end{algorithm}

\begin{theorem}[Algorithm Complexity: Reduced QR via Gram-Schmidt]\label{theorem:qr-reduced}
Algorithm~\ref{alg:reduced-qr} requires $\sim 2mn^2$ flops to compute the reduced QR decomposition of an $m\times n$ matrix with linearly independent columns and $m\geq n$.
\end{theorem}

\begin{proof}[of Theorem~\ref{theorem:qr-reduced}]
In step 3, the computation of $r_{ik}$ involves a vector inner product, which requires $m$ multiplications and $m-1$ additions, resulting in  $\underline{(k-1)(2m-1)}$ flops for all the $k-1$ iterations.

In step 5, the computation of $r_{ik}\bq_i$ requires $m$ flops and there are $k-1$ such scalar-vector multiplications, leading to a total of $m(k-1)$ flops. 
Additionally, the vector subtraction and additions require another $m(k-1)$ flops. Thus, step $5$ incurs a total cost of $\underline{2m(k-1)}$ flops. 

In step 6, computing the vector norm requires a vector inner product followed by a square root operation, which takes $\underline{2m}$ flops. 

Step 7 incurs a cost of $\underline{m}$ flops for the divisions.

Therefore, to compute the $k$-th column of $\bQ$ and $\bR$, it requires $(k-1)(2m-1)+ 2m(k-1)+2m+m=4mk-m-k+1$ flops. Let $f(k) = 4mk-m-k+1$, the total complexity can be obtained by 
$
\mathrm{cost}= f(n) +f(n-1) + \ldots +f(1).
$
Simple calculations  show the total complexity is $2mn^2+mn-3m-\frac{n^2-n}{2}$ flops, or $2mn^2$ flops if we keep only the leading term.
\end{proof}

\index{Orthogonal projection}
\index{Projection matrix}
\index{Projector}
\subsubsection*{\textbf{Orthogonal Projection: Preliminary for MGS}}
Upon revisiting  Equation~\eqref{equation:qr-gsp-equation}, i.e., step 2 to step 6 in Algorithm~\ref{alg:reduced-qr}, we observe that the first two equalities imply that
\begin{equation}\label{equation:qr-gsp-equation2}
\left.
\begin{aligned}
r_{ik} &= \bq_i^\top\ba_k, \,\forall i \in \{1,2,\ldots, k-1\}\\ 
\ba_k^\perp&\triangleq \ba_k-\sum_{i=1}^{k-1}r_{ik}\bq_i\\
\end{aligned}
\right\}
\rightarrow 
\ba_k^\perp= \ba_k- \bQ_{k-1}\bQ_{k-1}^\top \ba_k=(\bI-\bQ_{k-1}\bQ_{k-1}^\top )\ba_k,
\end{equation}
where $\bQ_{k-1}\triangleq[\bq_1,\bq_2,\ldots, \bq_{k-1}]$. This implies that $\bq_k$ can be obtained by 
$$
\bq_k = \frac{\ba_k^\perp}{\normtwo{\ba_k^\perp}} = \frac{(\bI-\bQ_{k-1}\bQ_{k-1}^\top )\ba_k}{\normtwo{(\bI-\bQ_{k-1}\bQ_{k-1}^\top )\ba_k}}.
$$
The matrix $(\bI-\bQ_{k-1}\bQ_{k-1}^\top )$ in the above equation is known as an \textit{orthogonal projection matrix} \footnote{More details can be referred to Appendix~\ref{section:by-geometry-hat-matrix}, Section~\ref{section:spec_app_eigproj}, or Definition~\ref{definition:projection_matrix_intro}.} that  projects $\ba_k$ \textbf{along} the column space of $\bQ_{k-1}$,
i.e., it projects a vector so that the projected vector is perpendicular to the column space of $\bQ_{k-1}$. 
Consequently, $\ba_k^\perp$ or $\bq_k$ calculated in this way will be orthogonal to $\cspace(\bQ_{k-1})$, residing in the null space of $\bQ_{k-1}^\top$, i.e., the space of $\nspace(\bQ_{k-1}^\top)$ by the fundamental theorem of linear algebra (Theorem~\ref{theorem:fundamental-linear-algebra}). 
\index{Fundamental theorem}

Let $\bP_1\triangleq(\bI-\bQ_{k-1}\bQ_{k-1}^\top )$. It can be shown  that $\bP_1=(\bI-\bQ_{k-1}\bQ_{k-1}^\top )$ is an orthogonal projection matrix such that $\bP_1\bv$ will project the vector $\bv$ onto the null space of $\bQ_{k-1}^\top$. 
Additionally, let $\bP_2\triangleq\bQ_{k-1}\bQ_{k-1}^\top$; then, $\bP_2$ is also an orthogonal projection matrix such that $\bP_2\bv$ will project the vector $\bv$ onto the column space of $\bQ_{k-1}$  ($\bP_2$ is called a \textit{complementary projector} of $\bP_1$, vice versa; Definition~\ref{definition:projection_matrix_intro}).

\begin{figure}[h!]
\centering  
\vspace{-0.35cm} 
\subfigtopskip=2pt 
\subfigbottomskip=2pt 
\subfigcapskip=-5pt 
\subfigure[Orthogonal projection.]{\label{fig:project-oblique}
	\includegraphics[width=0.4\linewidth]{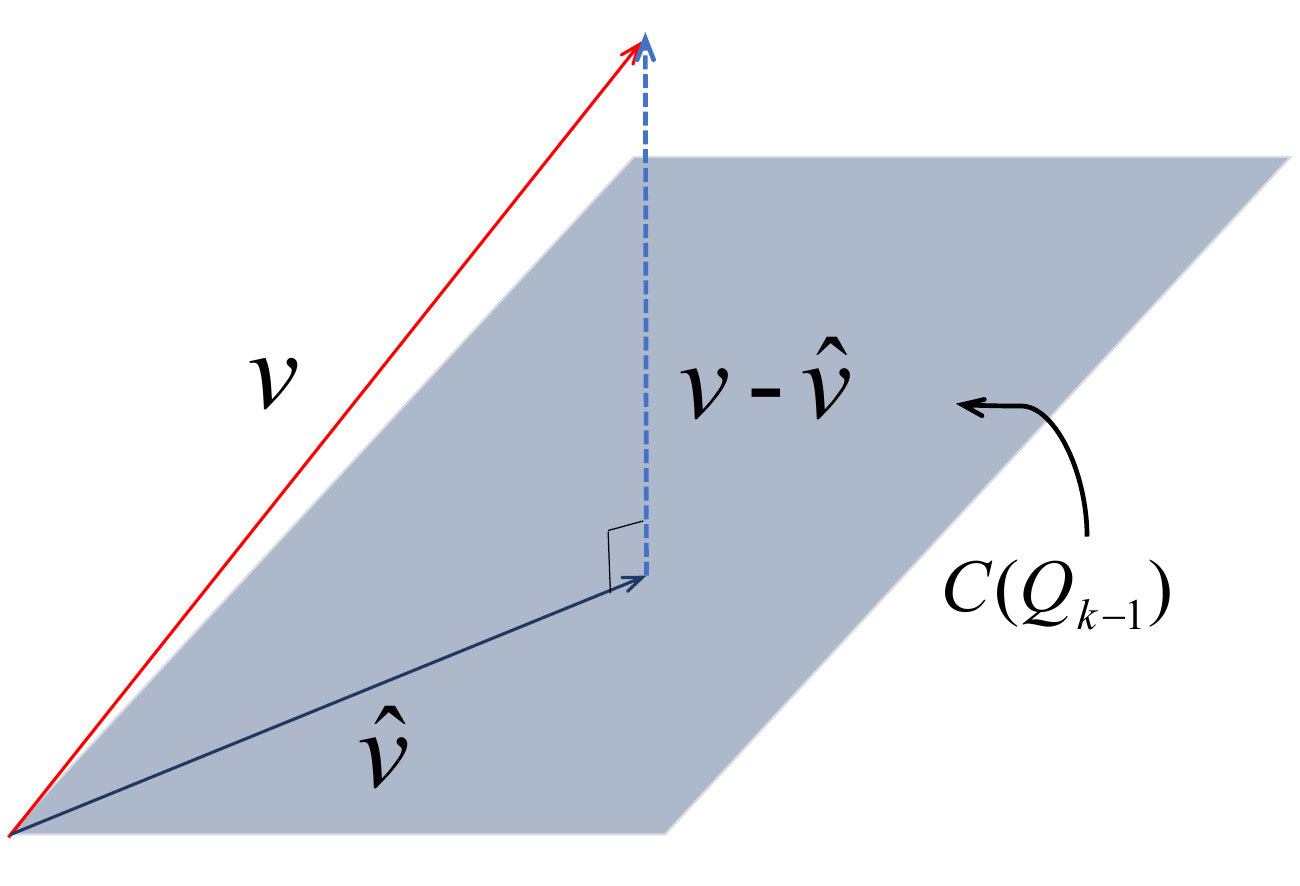}}
\quad 
\subfigure[Oblique projection.]{\label{fig:project-orthogonal}
	\includegraphics[width=0.4\linewidth]{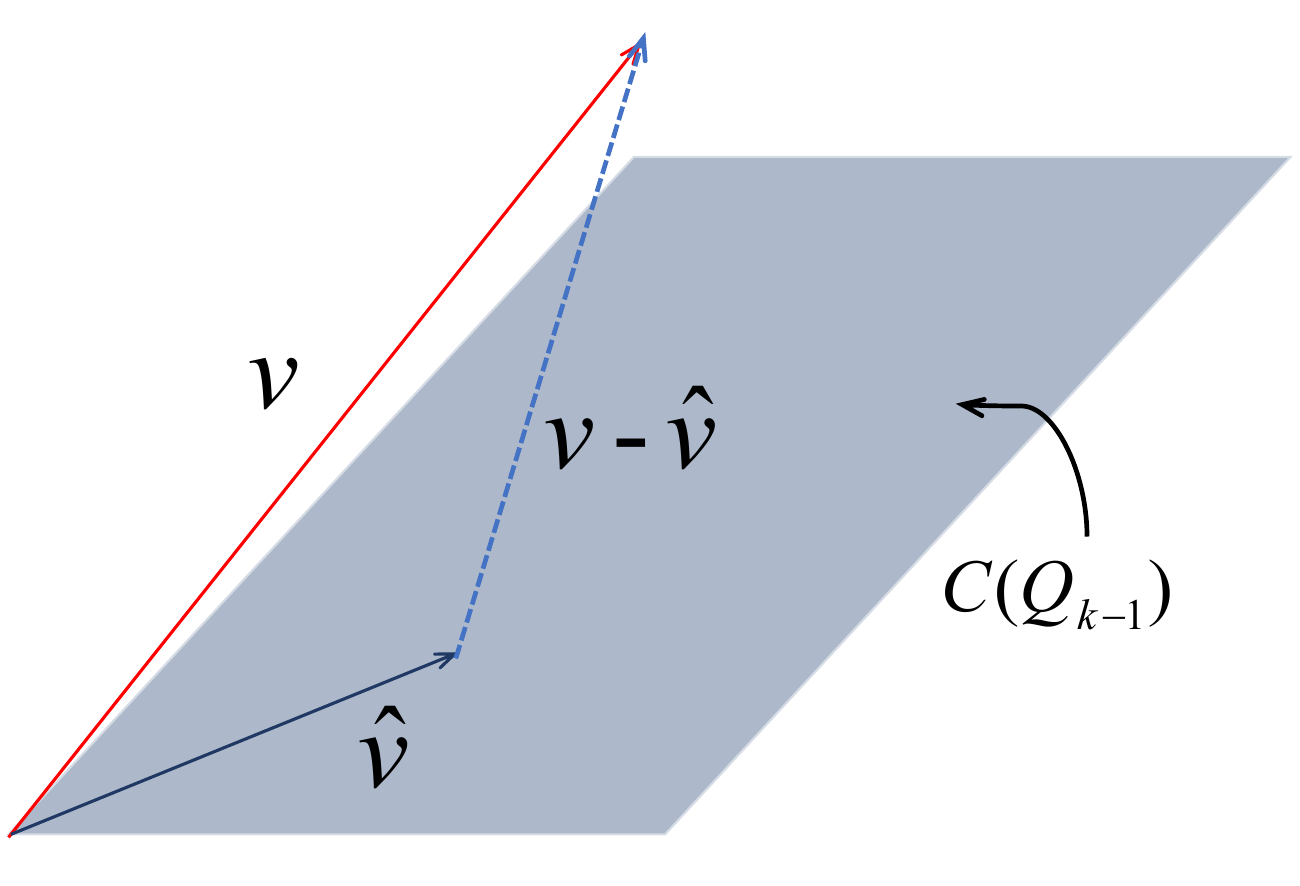}}
\caption{Demonstration of the difference between orthogonal projection and oblique projection.}
\label{fig:projection-oblie-and-orthogonal}
\end{figure}

But why can the matrices $\bP_1$ and $\bP_2$ effectively  project a vector onto the corresponding subspaces? We will show in Lemma~\ref{lemma:rank-of-ttt} that the column space of $\bQ_{k-1}$ is equal to the column space of $\bQ_{k-1}\bQ_{k-1}^\top$: 
$$
\cspace(\bQ_{k-1})=\cspace(\bQ_{k-1}\bQ_{k-1}^\top)=\cspace(\bP_2).
$$
Therefore, the result of $\bP_2\bv$ is a linear combination of the columns of $\bP_2$, which lies in the column space of $\bP_2$ or the column space of $\bQ_{k-1}$. The formal definition of a \textit{projection matrix} $\bP$ states that it is idempotent, satisfying $\bP^2=\bP$. This means that projecting a vector twice using $\bP$ yields the same result as projecting it once \footnote{Refer to Definition~\ref{definition:projection-matrix}.}. What makes the above $\bP_2=\bQ_{k-1}\bQ_{k-1}^\top $ special  is that the projection $\widehat{\bv}$ of any vector $\bv$ is perpendicular to the residual vector $\bv-\widehat{\bv}$:
$$
(\widehat{\bv}=\bP_2\bv) \perp (\bv-\widehat{\bv}).
$$
This goes to the original definition we provided earlier: the \textit{orthogonal projection matrix} \footnote{Refer to Definition~\ref{definition:orthogonal-projection-matrix}.}. 
To avoid confusion, the term \textit{oblique projection matrix} may be used in the nonorthogonal case, where the difference is illustrated in Figure~\ref{fig:projection-oblie-and-orthogonal}. 

When $\bP_2$ is an orthogonal projection matrix, $\bP_1\triangleq\bI-\bP_2$ also serves as an orthogonal projection matrix that can project any vector onto the space perpendicular to  $\cspace(\bQ_{k-1})$, i.e., onto $\nspace(\bQ_{k-1}^\top)$. Therefore, we conclude the presence of two orthogonal projections:
$$
\left\{
\begin{aligned}
\bP_1&=\bI-\bP_2: &\gap& \text{project onto $\nspace(\bQ_{k-1}^\top)$, \textbf{along} the column space of $\bQ_{k-1}$;} \\
\bP_2&=\bQ_{k-1}\bQ_{k-1}^\top: &\gap& \text{project onto $\cspace(\bQ_{k-1})$, \textbf{onto} the column space of $\bQ_{k-1}$} .
\end{aligned}
\right.
$$

An additional result to note arises when the columns of $\bQ_{k-1}$ are mutually orthonormal. In this case, we can observe the following decomposition:
\begin{equation}\label{equation:qr-orthogonal-equality}
\boxed{\bP_1 = \bI - \bQ_{k-1}\bQ_{k-1}^\top = (\bI-\bq_1\bq_1^\top)(\bI-\bq_2\bq_2^\top)\ldots (\bI-\bq_{k-1}\bq_{k-1}^\top),}
\end{equation}
where $\bQ_{k-1}=[\bq_1,\bq_2,\ldots, \bq_{k-1}]$ and each $(\bI-\bq_i\bq_i^\top)$ serves to project a vector onto the perpendicular space of $\bq_i$.

\subsubsection*{\textbf{Modified Gram-Schmidt (MGS) Process}}
To emphasize the modified Gram-Schmidt process and to make a connection to the equivalent projection in Equation~\eqref{equation:qr-orthogonal-equality}, we begin by illustrating a lemma that presents an alternative approach for obtaining the entries in the upper triangular matrix $\bR$ of the QR decomposition.
\begin{lemma}[Modified Gram-Schmidt]
Given a vector set $[\ba_1, \ba_2, \ldots,\ba_{k-1}, \ba_k$], where the first $k-1$ column are spanned by $k-1$ orthonormal vectors $[\bq_1, \bq_2, \ldots, \bq_{k-1}]$:
$$
\cspace([\ba_1, \ba_2, \ldots,\ba_{i}]) = \cspace([\bq_1, \bq_2, \ldots, \bq_{i}]), \gapforall \forall i\in \{1,2,\ldots, k-1\}. 
$$
Therefore, $r_{ik} = \bq_i^\top\ba_k$ represents the magnitude of the projection of $\ba_k$ on the vector $\bq_i$ (since $\bq_i$ is of unit length). Then it follows that 
$$
\begin{aligned}
\bq_i^\top\ba_k &= \bq_i^\top (\ba_k \underbrace{- r_{1k}\bq_1 - r_{2k}\bq_2 - \ldots - r_{i-1,k}\bq_{i-1}}_{\text{orthogonal to $\bq_i$}})\\
&= \bq_i^\top \big(\ba_k - \sum_{j=1}^{i-1}r_{jk}\bq_j \big), \gap \forall i\in \{1,2,\ldots, k-1\}.
\end{aligned}
$$
This can be easily verified since $\bq_i$ is orthonormal to $\{\bq_1, \bq_2, \ldots, \bq_{i-1}\}$. This observation implies another update for the $k$-th column of $\bR$. 
\end{lemma}
The lemma above reveals a second algorithm to compute the reduced QR decomposition of a matrix, as shown in Algorithm~\ref{alg:qr-mgs-right} of which the algorithm on the left is identical to Algorithm~\ref{alg:reduced-qr} (with slight modifications) to highlight the differences.

\noindent
\begin{minipage}[t]{0.495\linewidth}
\begin{algorithm}[H] 
\caption{CGS (=Algorithm~\ref{alg:reduced-qr})} 
\label{alg:qr-mgs-left}
\begin{algorithmic}[1] 
\Require $\bA\in \real^{m\times n}$ with full column rank;
\For{$k=1$ to $n$} 
\State $\ba_k^\perp\leftarrow\ba_k$;
\For{$i=1$ to $k-1$}
\State $\boxed{r_{ik} \leftarrow\bq_i^\top\ba_k}$;
\State $\ba_k^\perp\leftarrow \ba_k^\perp-r_{ik}\bq_i$; \,\,($\dagger$)
\EndFor 
\State $r_{kk} \leftarrow \normtwo{\ba_k^\perp}$; 
\State $\bq_k \leftarrow \ba_k^\perp/r_{kk}$; 
\EndFor 
\end{algorithmic} 
\end{algorithm}
\end{minipage}%
\hfil 
\begin{minipage}[t]{0.495\linewidth}
\begin{algorithm}[H] 
\caption{MGS}
\label{alg:qr-mgs-right}
\begin{algorithmic}[1] 
\Require $\bA\in \real^{m\times n}$ with full column rank;
\For{$k=1$ to $n$} 
\State $\ba_k^\perp\leftarrow\ba_k$;
\For{$i=1$ to $k-1$}
\State $\boxed{r_{ik} \leftarrow\bq_i^\top\textcolor{mylightbluetext}{\ba_k^\perp}}$;
\State $\ba_k^\perp\leftarrow \ba_k^\perp-r_{ik}\bq_i$; \,\,($\ast$)
\EndFor 
\State $r_{kk} \leftarrow \normtwo{\ba_k^\perp}$; 
\State $\bq_k \leftarrow \ba_k^\perp/r_{kk}$; 
\EndFor 
\end{algorithmic} 
\end{algorithm}
\end{minipage}
\index{Modified Gram–Schmidt}

The process described above is referred to as the \textit{modified Gram-Schmidt (MGS) process}, whereas the previous one is also known as the \textit{classical Gram-Schmidt (CGS) process}. In theory, both CGS and MGS are equivalent in the
sense that they compute exactly the same QR decompositions when exact arithmetic is employed. However, in practice, with the presence of round-off errors, the orthonormal columns of $\bQ$ computed by MGS tend to be ``more orthonormal" than those computed by CGS. 

To see the equivalence of the above two algorithms,
we note that the equality ($\dagger$) in Algorithm~\ref{alg:qr-mgs-left} is equivalent to 
\begin{equation}\label{equation:qr-orthogonal-equality2}
\ba_k^\perp = \ba_k - (\bq_1^\top\ba_k)\bq_1 - (\bq_2^\top\ba_k)\bq_2-\ldots -(\bq_{k-1}^\top\ba_k)\bq_{k-1} = (\bI - \bQ_{k-1}\bQ_{k-1}^\top) \ba_k.
\end{equation}
And
the equality ($\ast$) in Algorithm~\ref{alg:qr-mgs-right} can be reformulated as (via the step 4 and step 5 of the algorithm)
$$
\begin{aligned}
\ba_k^\perp &:= \ba_k^\perp-r_{ik}\bq_i
= \ba_k^\perp-(\bq_i^\top\ba_k^\perp)\bq_i
=\ba_k^\perp-\bq_i\bq_i^\top\ba_k^\perp
=(\bI-\bq_i\bq_i^\top)\ba_k^\perp.
\end{aligned}
$$
That is, $\ba_k^\perp$ will be updated by 
\begin{equation}\label{equation:qr-orthogonal-equality3}
\ba_k^\perp=
\left\{(\bI-\bq_{k-1}\bq_{k-1}^\top)\ldots\left[(\bI-\bq_2\bq_2^\top)\left((\bI-\bq_1\bq_1^\top) \ba_k\right)\right]\right\},
\end{equation}
where the nested parentheses in MGS denote the computation order.
The comparison of \eqref{equation:qr-orthogonal-equality2} and \eqref{equation:qr-orthogonal-equality3} matches the orthogonal projection matrix equality in Equation~\eqref{equation:qr-orthogonal-equality} that
$$
\begin{aligned}
\bP_1 &= \bI-\bQ_{k-1}\bQ_{k-1}^\top
=(\bI-\bq_1\bq_1^\top)(\bI-\bq_2\bq_2^\top)\ldots (\bI-\bq_{k-1}\bq_{k-1}^\top)
=\prod_{i=1}^{k-1}(\bI-\bq_i\bq_i^\top),
\end{aligned}
$$
where $\bQ_{k-1}=[\bq_1,\bq_2,\ldots, \bq_{k-1}]$.

\begin{figure}[H]
\centering  
\vspace{-0.35cm} 
\subfigtopskip=2pt 
\subfigbottomskip=2pt 
\subfigcapskip=-5pt 
\subfigure[CGS, step 1: \textcolor{mylightbluetext}{blue} vector; step 2: \textcolor{mydarkgreen}{green} vector; step 3: \textcolor{mydarkpurple}{purple} vector.]{\label{fig:projection-mgs-demons-cgs}
	\includegraphics[width=0.4\linewidth]{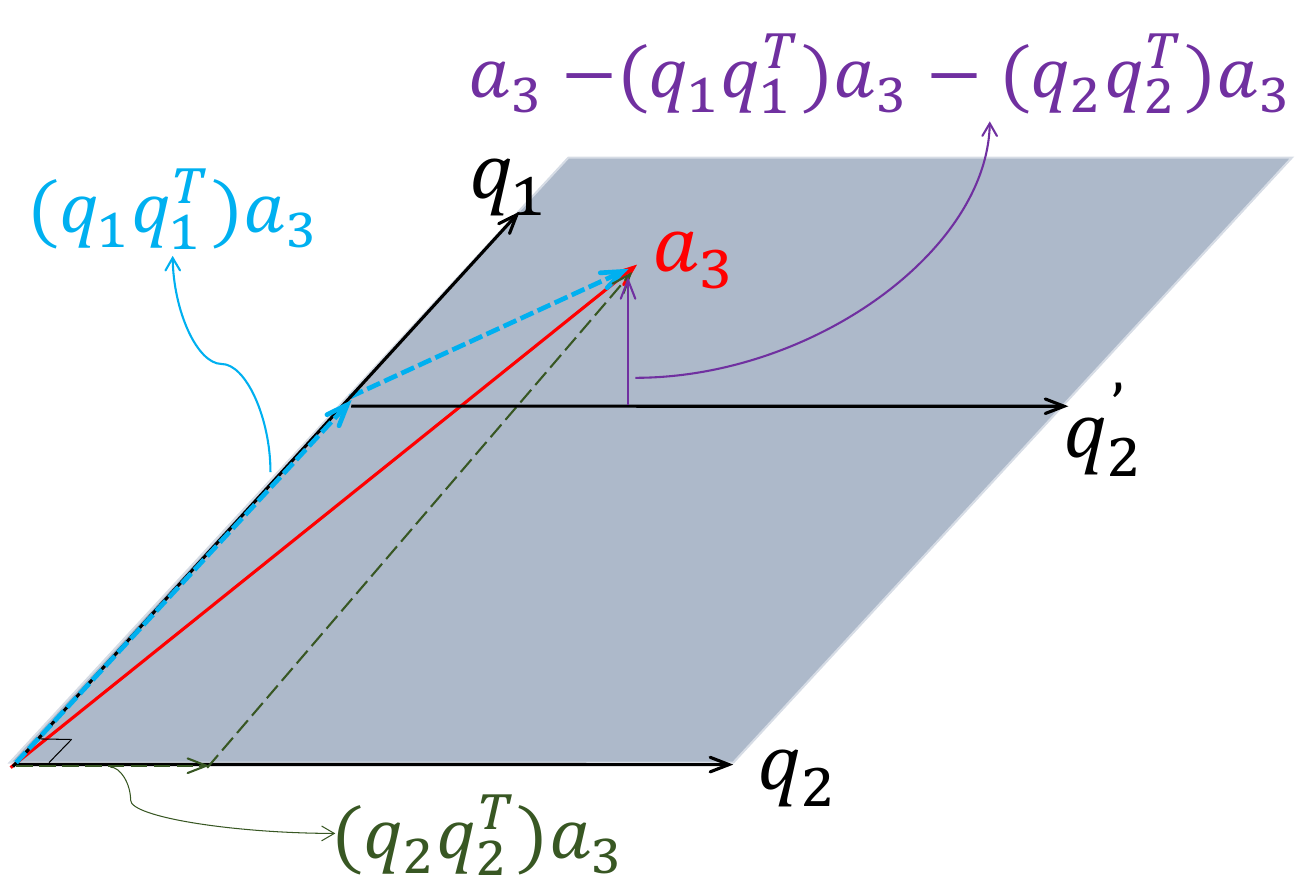}}
\quad 
\subfigure[MGS, step 1: \textcolor{mylightbluetext}{blue} vector; step 2: \textcolor{mydarkpurple}{purple} vector.]{\label{fig:projection-mgs-demons-mgs}
	\includegraphics[width=0.4\linewidth]{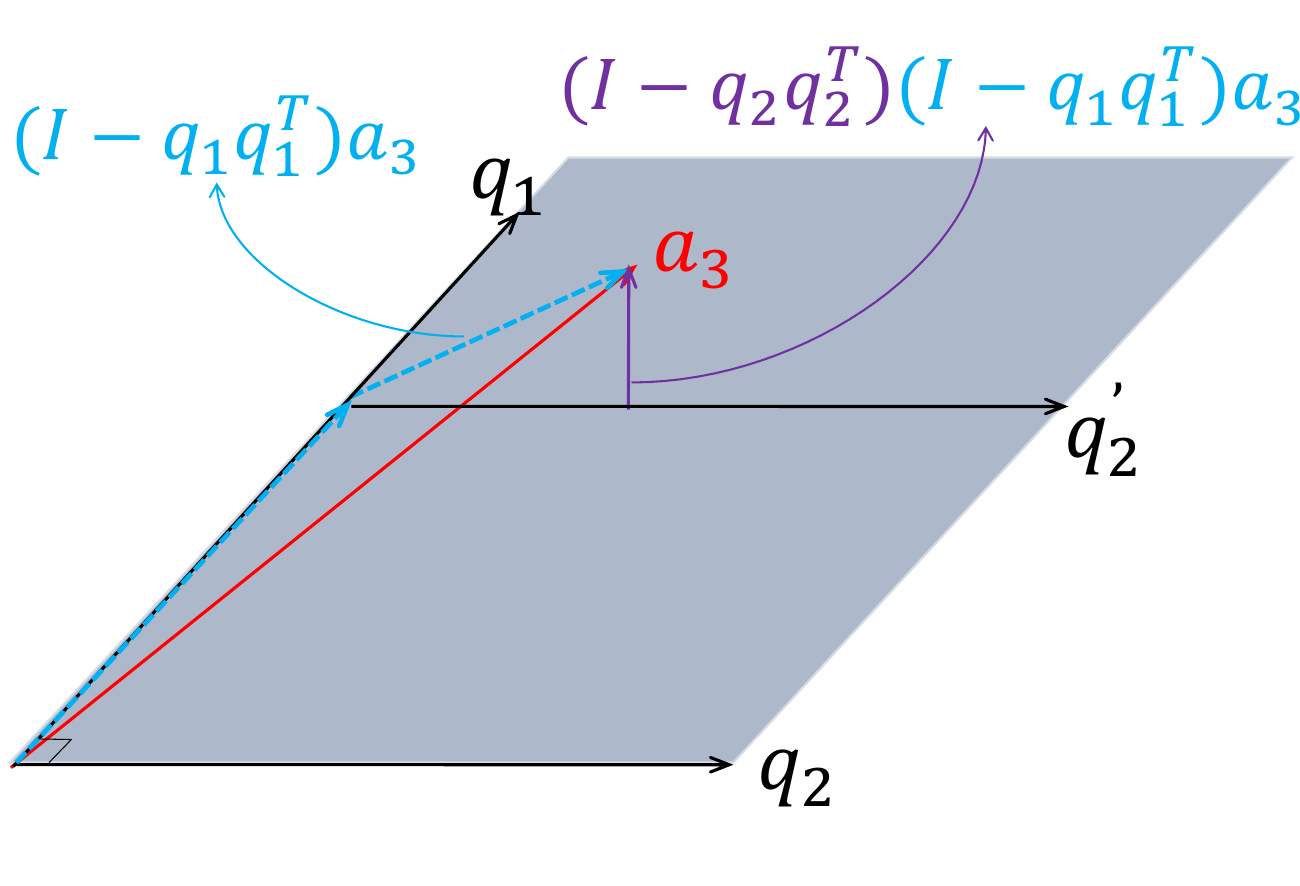}}
\caption{CGS vs MGS in three-dimensional space, where $\bq_2^\prime$ is parallel to $\bq_2$ so that projecting onto $\bq_2$ is equivalent to projecting onto $\bq_2^\prime$.}
\label{fig:projection-mgs-demons-3d}
\end{figure}

\paragraph{What's the difference?}
Taking a three-column matrix $\bA=[\ba_1, \ba_2, \ba_3]$ as an example. Suppose we have already computed $\{\bq_1, \bq_2\}$ such that $\spn\{\bq_1, \bq_2\}=\spn\{\ba_1, \ba_2\}$, and we want to proceed to compute  $\bq_3$.

In the CGS algorithm, the orthogonalization of column $\ba_3$ against columns $\{\bq_1, \bq_2\}$ is achieved by projecting the original column $\ba_3$ of $\bA$ onto $\{\bq_1, \bq_2\}$, respectively, followed by subtracting these projections at once:
\begin{equation}\label{equation:cgs-3d-exmp}
\left\{
\begin{aligned}
\ba_3^\perp &= \ba_3 - (\bq_1^\top\ba_3)\bq_1 - (\bq_2^\top\ba_3)\bq_2\\
&= \ba_3 - (\bq_1\bq_1^\top)\ba_3 - \boxed{(\bq_2\bq_2^\top)\textcolor{mylightbluetext}{\ba_3}};\\
\bq_3 &=  {\ba_3^\perp}/{\normtwo{\ba_3^\perp}},
\end{aligned}
\right.
\end{equation}
as shown in Figure~\ref{fig:projection-mgs-demons-cgs}.

In the MGS algorithm, on the other hand, the components along each $\{\bq_1, \bq_2\}$ are immediately subtracted out of the rest of the column $\ba_3$ as soon as the vectors  $\{\bq_1, \bq_2\}$ are computed. 
Therefore, the orthogonalization of column $\ba_3$ against $\{\bq_1, \bq_2\}$ is not performed by projecting the original column $\ba_3$ against $\{\bq_1, \bq_2\}$ as it is in CGS, but rather against a vector obtained by subtracting from that column $\ba_3$ of $\bA$ the components in the direction of $\bq_1, \bq_2$ successively. This is important because the error components of $\bq_3$ in $\spn\{\bq_1, \bq_2\}$ will be smaller (we will  discuss this further in the following paragraphs).

More precisely, in the MGS algorithm, the orthogonalization of column $\ba_3$ against $\bq_1$ is performed by subtracting the component of $\bq_1$ from the vector $\ba_3$:
$$
\ba_3^{(1) }=  (\bI-\bq_1\bq_1^\top)\ba_3 = \ba_3 - (\bq_1\bq_1^\top)\ba_3,
$$
where $\ba_3^{(1) }$ represents the component of $\ba_3$ that lies in a space perpendicular to $\bq_1$. And the further step is performed by 
\begin{equation}\label{equation:mgs-3d-exmp}
\begin{aligned}
\ba_3^{(2) }=  (\bI-\bq_2\bq_2^\top)\ba_3^{(1) }&=\ba_3^{(1) }-(\bq_2\bq_2^\top)\ba_3^{(1) }\\
&=\ba_3 - (\bq_1\bq_1^\top)\ba_3-\boxed{(\bq_2\bq_2^\top)\textcolor{mylightbluetext}{\ba_3^{(1) }}},
\end{aligned}
\end{equation}
where $\ba_3^{(2) }$ represents the component of $\ba_3^{(1) }$ that lies in a space perpendicular to $\bq_2$. And we highlight the difference from the CGS in Equation~\eqref{equation:cgs-3d-exmp} using \textcolor{mylightbluetext}{\boxed{blue}} text. 
As a result, $\ba_3^{(2) }$ corresponds to the component of $\ba_3$ that lies in the space perpendicular to $\{\bq_1, \bq_2\}$, as shown in Figure~\ref{fig:projection-mgs-demons-mgs}. 

\index{Cancellation}
\index{Catastrophic cancellation}
\subsubsection*{\textbf{Main Difference and Catastrophic Cancellation}}
The key difference between the CGS and MGS processes lies in the fact that $\ba_3$ can generally have large components in $\spn\{\bq_1, \bq_2\}$. 
In such cases, one starts with large
values and ends up with small values that yields significant relative errors in them. 
This phenomenon is commonly referred to as \textit{catastrophic cancellation}. 
Whereas, $\ba_3^{(1) }$ lies in the direction perpendicular to $\bq_1$ and carries only a small ``error" component in the direction of $\bq_1$. Comparing the \fbox{boxed} terms in Equations~\eqref{equation:cgs-3d-exmp} and \eqref{equation:mgs-3d-exmp}, it is not hard to see that $(\bq_2\bq_2^\top)\ba_3^{(1) }$ in Equation~\eqref{equation:mgs-3d-exmp} is more accurate based on the above argument. And thus, because of the much smaller error in this projection factor, the MGS introduces smaller orthogonalization error at each subtraction step compared to the CGS method. In fact, it can be shown that the final $\bQ$ obtained in the CGS satisfies
$$
\normtwo{\bI-\bQ\bQ^\top} \leq \mathcalO(\epsilon \kappa^2(\bA)),
$$
where $\kappa(\bA)$ is a value larger than 1 determined by $\bA$.
Whereas, in MGS, the error satisfies
$$
\normtwo{\bI-\bQ\bQ^\top} \leq \mathcalO(\epsilon \kappa(\bA)).
$$
That is, the $\bQ$ obtained via MGS is ``more orthogonal."
Therefore, we summarize the difference between the CGS and MGS processes for obtaining $\bq_k$ from the $k$-th column $\ba_k$ of $\bA$, given the orthonormalized vectors $\{\bq_1, \bq_2, \ldots, \bq_{k-1}\}$:
$$
\begin{aligned}
\text{(CGS)}: &\, \text{obtain $\bq_k$ by normalizing $\ba_k^\perp=(\bI-\bQ_{k-1}\bQ_{k-1}^\top)\ba_k$;} \\
\text{(MGS)}: &\, \text{obtain $\bq_k$ by normalizing  $\ba_k^\perp=\big\{(\bI-\bq_{k-1}\bq_{k-1}^\top)\ldots\big[(\bI-\bq_2\bq_2^\top)\big((\bI-\bq_1\bq_1^\top) \ba_k\big)\big]\big\}$.} 
\end{aligned}
$$

\subsubsection*{\textbf{Triangular Orthogonalization in CGS and MGS}}
We now illustrate that in the CGS or MGS processes, the orthogonal matrix $\bQ$ is obtained using a set of triangular matrices. For simplicity, we only discuss the situation in the CGS and follow up the three-column example, where $\bA=[\ba_1, \ba_2, \ba_3]\in \real^{3\times 3}$. The above discussion shows that the mutually orthonormal vectors $\{\bq_1, \bq_2, \bq_3\}$ can be obtained as follows:
$$
\left\{
\begin{aligned}
\bq_1 &= {\ba_1}/{r_{11}};\\
\bq_2 &= {(\ba_2 - r_{12}\bq_1)}/{r_{22}};\\
\bq_3 &= {(\ba_3-r_{13}\bq_1-r_{23}\bq_2)}/{r_{33}}.
\end{aligned}
\right.
$$
Whilst, the three mutually orthonormal vectors can be equivalently obtained by 
$$
\bQ\bR_3\bR_2\bR_1= \bA \leadto \bQ=\bA\bR_1^{-1}\bR_2^{-1}\bR_3^{-1},
$$
where 
$$
\bR_3=
\begin{bmatrix}
1 & 0 & r_{13} \\
0 & 1 & r_{23}\\
0 & 0 & r_{33}
\end{bmatrix},
\gap
\bR_2=
\begin{bmatrix}
1 & r_{12} & 0 \\
0 & r_{22} &0\\
0 & 0 & 1
\end{bmatrix},
\gap
\bR_1=
\begin{bmatrix}
r_{11} & 0 & 0 \\
0 & 1 &0\\
0 & 0 & 1
\end{bmatrix},
$$
such that 
$$
\bR_3\bR_2\bR_1 = \bR=
\begin{bmatrix}
r_{11}& r_{12} & r_{13}\\
0 & r_{22} & r_{13}\\
0 & 0 & r_{33}
\end{bmatrix}.
$$
The  procedure  $\bA\bR_1^{-1}\bR_2^{-1}\bR_3^{-1}$ obtains the $\{\bq_1, \bq_2, \bq_3\}$ in a successive manner, where $\bA\bR_1^{-1}$ places $\bq_1$ into the first column of $\bQ$; $(\bA\bR_1^{-1})\bR_2^{-1}$ places $\bq_2$ into the second column of $\bQ$; and $(\bA\bR_1^{-1}\bR_2^{-1})\bR_3^{-1}$ places $\bq_3$ into the third column of $\bQ$. 
This procedure is called  \textit{triangular orthogonalization} in the Gram-Schmidt process. The triangular orthogonalization is problematic in the sense that the condition number of a triangular matrix ($\bR_1^{-1}, \bR_2^{-1}, \bR_3^{-1}$ in the above three-column example) can be anything (i.e., can be arbitrarily large). And the Gram-Schmidt process contains a series of them, where the condition number (see Appendix~\ref{appendix:condition_number} for more details on this topic) can grow very large, leading to numerical instability in the orthogonalization process.

\subsubsection*{\textbf{More to go, preliminaries for Householder and Givens methods}}
Although both methods have their advantages, in practice,  MGS  usually outperforms  CGS.  
An example will be provided in the sequel to illustrate this. 
However,  MGS can still fall victim to the \textit{catastrophic cancellation} problem. Suppose, in iteration $k$ of the MGS Algorithm~\ref{alg:qr-mgs-right}, $\ba_k$ is almost in the span of $\{\bq_1, \bq_2, \ldots, \bq_{k-1}\}$. This will result in that $\ba_k^\perp$ has only a small component that is perpendicular to $\spn\{\bq_1, \bq_2, \ldots, \bq_{k-1}\}$, whereas the ``error" component in the $\spn\{\bq_1, \bq_2, \ldots, \bq_{k-1}\}$ will be amplified, resulting in $\bQ$ being less orthonormal. 
As mentioned earlier, both the CGS and MGS methods suffer from the same main disadvantage: they obtain the orthogonal matrix $\bQ$ through the upper triangular matrix $\bR$. Specifically, for $\bA\in \real^{m\times n}$, $\bQ$ can be obtained by using the following equation:
$$
\bQ=\bA\underbrace{\bR_1^{-1}\bR_2^{-1}\ldots \bR_n^{-1}}_{\bR^{-1}}.
$$

In this case, if we can find a successive set of orthogonal matrices $\{\bQ_1, \bQ_2, \ldots, \bQ_l\}$ such that $\bQ_l\ldots\bQ_2\bQ_1\bA$ is triangularized, then $\bQ=(\bQ_l\ldots\bQ_2\bQ_1)^\top$ will be ``more" orthogonal than that in the CGS or MGS method since the condition numbers for the orthogonal matrices are all 1. We will discuss this method in Section~\ref{section:qr-via-householder}  and \ref{section:qr-givens}  using Householder reflectors and Givens rotations.

\subsubsection*{\textbf{Example for MGS vs CGS}}
To catch a glimpse of the differences of the MGS and CGS approaches we have introduced, let us start by trying it out on an example with a $4\times 3$ \textit{Lauchli matrix}.
The Lauchli matrix is a rectangular matrix of size $(n+1) \times n$ that has ones on the top row and the parameter $\epsilon=\sqrt{\epsilon_{mach}}$ on the diagonal starting at entry $(2,1)$, i.e., on the lower subdiagonal.
The value $\epsilon_{mach}$ represents the machine precision.
\index{Lauchli matrix}
\begin{example}[MGS vs CGS]
Let $\epsilon = \sqrt{\epsilon_{mach}}$ and we will examine the QR decomposition of the following matrix using CGS and MGS methods:
$$
\bA = 
\begin{bmatrix}
1& 1 & 1 \\
\epsilon & 0  & 0 \\
0 & \epsilon  & 0 \\
0 & 0 & \epsilon
\end{bmatrix}
=
\begin{bmatrix}
\ba_1 & \ba_2 & \ba_3 
\end{bmatrix}.
$$
Note that during the calculation, $1 + \epsilon_{\text{mach}}$ will be rounded to 1.
The CGS method follows the steps outlined below:
\begin{itemize}
\item Compute the unit-length vector $\bq_1$ such that $\cspace([\bq_1]) = \cspace([\ba_1])$:
\begin{itemize}
\item Compute $r_{11}$: $r_{11} = \normtwo{\ba_1}=\sqrt{1+\epsilon_{mach}}\approx1$;
\item Compute $\bq_1$: $\bq_1=\frac{\ba_1}{\normtwo{\ba_1}}=\ba_1$;
\end{itemize}
\item Compute the unit-length vector $\bq_2$ such that $\cspace([\bq_1,\bq_2]) = \cspace([\ba_1,\ba_2])$:
\begin{itemize}
\item Compute $r_{12}$: $r_{12} = \ba_2^\top\bq_1=1$; 
\item Compute $\ba_2^\perp$: $\ba_2^\perp = \ba_2 - r_{12}\bq_1=\ba_2-\ba_1 = [0,-\epsilon, \epsilon, 0]^\top$;
\item Compute $r_{22}$: $r_{22}=\normtwo{\ba_2^\perp}=\sqrt{2\epsilon_{mach}}=\sqrt{2}\epsilon$;
\item Compute $\bq_2$: $\bq_2 = \frac{\ba_2^\perp}{\normtwo{\ba_2^\perp}}=\frac{\ba_2^\perp}{r_{22}}=[0,-\frac{1}{\sqrt{2}}, \frac{1}{\sqrt{2}}, 0]^\top$;
\end{itemize}
\item Compute the unit-length vector $\bq_3$ such that $\cspace([\bq_1,\bq_2,\bq_3]) = \cspace([\ba_1,\ba_2,\ba_3])$:
\begin{itemize}
\item Compute $r_{13}$: $r_{13}= \ba_3^\top\bq_1=\ba_3^\top\ba_1=1$;
\item Compute $r_{23}$: $r_{23}= \ba_3^\top\bq_2=0$;
\item Compute $\ba_3^\perp$: $\ba_3^\perp = \ba_3 - r_{13}\bq_1 - r_{23}\bq_2=[0,-\epsilon, 0, \epsilon]^\top$;
\item Compute $r_{33}$: $r_{33}=\normtwo{\ba_3^\perp}=\sqrt{2\epsilon_{mach}}=\sqrt{2}\epsilon$;
\item Compute $\bq_3$: $\bq_3 = \frac{\ba_3^\perp}{\normtwo{\ba_3^\perp}}=\frac{\ba_3^\perp}{r_{33}}=[0,-\frac{1}{\sqrt{2}},0, \frac{1}{\sqrt{2}}]^\top$;
\end{itemize}
\item This results in 
$$
\bA=
\begin{bmatrix}
1& 0& 0 \\
\epsilon & -\frac{1}{\sqrt{2}}  & -\frac{1}{\sqrt{2}} \\
0 & \frac{1}{\sqrt{2}}  & 0 \\
0 & 0 & \frac{1}{\sqrt{2}}
\end{bmatrix}
\begin{bmatrix}
1& 1 & 1 \\
0 & \sqrt{2}\epsilon  & 0 \\
0 & 0  & \sqrt{2}\epsilon 
\end{bmatrix}=
\bQ_1\bR_1
$$
\end{itemize}

Whilst, the MGS proceeds as follows:
\begin{itemize}
\item Compute the unit-length vector $\bq_1$ such that $\cspace([\bq_1]) = \cspace([\ba_1])$:
\begin{itemize}
\item Compute $r_{11}$: $r_{11} = \normtwo{\ba_1}=\sqrt{1+\epsilon_{mach}}\approx1$;
\item Compute $\bq_1$: $\bq_1=\frac{\ba_1}{\normtwo{\ba_1}}=\ba_1$;
\end{itemize}
\item Compute the unit-length vector $\bq_2$ such that $\cspace([\bq_1,\bq_2]) = \cspace([\ba_1,\ba_2])$:
\begin{itemize}
\item Compute $r_{12}$: $r_{12} = \ba_2^\top\bq_1=1$; 
\item Compute $\ba_2^\perp$: $\ba_2^\perp = \ba_2 - r_{12}\bq_1=\ba_2-\ba_1 = [0,-\epsilon, \epsilon, 0]^\top$;
\item Compute $r_{22}$: $r_{22}=\normtwo{\ba_2^\perp}=\sqrt{2\epsilon_{mach}}=\sqrt{2}\epsilon$;
\item Compute $\bq_2$: $\bq_2 = \frac{\ba_2^\perp}{\normtwo{\ba_2^\perp}}=\frac{\ba_2^\perp}{r_{22}}=[0,-\frac{1}{\sqrt{2}}, \frac{1}{\sqrt{2}}, 0]^\top$;
\item \textcolor{mylightbluetext}{Up until now, we remain identical to the CGS};
\end{itemize}
\item Compute the unit-length vector $\bq_3$ such that $\cspace([\bq_1,\bq_2,\bq_3]) = \cspace([\ba_1,\ba_2,\ba_3])$:
\begin{itemize}
\item Compute $r_{13}$: $r_{13}= \ba_3^\top\bq_1=\ba_3^\top\ba_1=1$;
\item Compute temporary $\ba_3^\perp=\ba_3 - r_{13}\bq_1=[0,-\epsilon, 0, \epsilon]^\top$;
\item Compute $r_{23}$: $r_{23}= \bq_2^\top \textcolor{mylightbluetext}{\ba_3^\perp}=\textcolor{mylightbluetext}{\frac{\epsilon}{\sqrt{2}}}$;
\item Compute final $\ba_3^\perp$: $\ba_3^\perp = \underbrace{\ba_3 - r_{13}\bq_1}_{\text{the old $\ba_3^\perp$}} - r_{23}\bq_2=\textcolor{mylightbluetext}{[0,-\epsilon/2, -\epsilon/2, \epsilon]^\top}$;
\item Compute $r_{33}$: $r_{33}=\normtwo{\ba_3^\perp}=\sqrt{2\epsilon_{mach}}=\textcolor{mylightbluetext}{\frac{\sqrt{6}}{2}\epsilon}$;
\item Compute $\bq_3$: $\bq_3 = \frac{\ba_3^\perp}{\normtwo{\ba_3^\perp}}=\frac{\ba_3^\perp}{r_{33}}=\textcolor{mylightbluetext}{[0,-\frac{1}{\sqrt{6}}, -\frac{1}{\sqrt{6}}, \frac{2}{\sqrt{6}}]^\top}$;
\end{itemize}
\item This results in 
$$
\bA=
\begin{bmatrix}
1& 0& 0 \\
\epsilon & -\frac{1}{\sqrt{2}}  & \textcolor{mylightbluetext}{-\frac{1}{\sqrt{6}}}\\
0 & \frac{1}{\sqrt{2}}  &  \textcolor{mylightbluetext}{-\frac{1}{\sqrt{6}}}\\
0 & 0 &  \textcolor{mylightbluetext}{\frac{2}{\sqrt{6}}}\\
\end{bmatrix}
\begin{bmatrix}
1& 1 & 1 \\
0 & \sqrt{2}\epsilon  & \textcolor{mylightbluetext}{\frac{\epsilon}{2}} \\
0 & 0  & \frac{\sqrt{6}}{2}\epsilon
\end{bmatrix}=
\bQ_2\bR_2.
$$

\end{itemize}
We notice that 
$$
\bQ_1^\top\bQ_1=
\begin{bmatrix}
1+\epsilon_{mach} & -\frac{1}{\sqrt{2}}\epsilon&-\frac{1}{\sqrt{2}}\epsilon\\
-\frac{1}{\sqrt{2}}\epsilon & 1 & \frac{1}{2}\\
-\frac{1}{\sqrt{2}}\epsilon & \frac{1}{2} & 1
\end{bmatrix}
\gap \text{and}\gap
\bQ_2^\top\bQ_2=
\begin{bmatrix}
1+\epsilon_{mach} & -\frac{1}{\sqrt{2}}\epsilon & -\frac{1}{\sqrt{6}}\epsilon\\
-\frac{1}{\sqrt{2}}\epsilon & 1 & 0\\
-\frac{1}{\sqrt{6}}\epsilon & 0 & 1
\end{bmatrix},
$$
which shows that $\bQ_2$ is better in the sense of orthogonality.
\end{example}

\subsubsection*{\textbf{Row-Wise MGS, Recursive Algorithm}}
The algorithms presented  in Algorithm~\ref{alg:qr-mgs-left} and \ref{alg:qr-mgs-right} are used to calculate the entries of the upper triangular matrix $\bR$ in an element-wise and column-by-column manner. 
Suppose $\bA$ has column partition $\bA=[\ba_1, \bA_2]$, where $\bA_2=[\ba_2, \ba_3, \ldots, \ba_n]\in \real^{m\times (n-1)}$. Notice in the CGS Algorithm~\ref{alg:qr-mgs-left}, the first row of $\bR$ can be obtained by 
$$
\left.
\begin{aligned}
r_{11} &= \normtwo{\ba_1}\\
r_{1k} & = \bq_1^\top\ba_k, \,\,\, \forall k\in\{2,3,\ldots, n\}
\end{aligned}\right\}
\leadto
\left\{
\begin{aligned}
r_{11} &= \normtwo{\ba_1}\\
\br_{12}^\top & = \bq_1^\top\bA_2, \,\,\, \br_{12}=[r_{12}, r_{13}, \ldots, r_{1n}].
\end{aligned}\right.
$$
Therefore, the QR decomposition of $\bA$ is given by 
$$
\bA = 
\begin{bmatrix}
\ba_1 & \bA_2
\end{bmatrix}=
\begin{bmatrix}
\bq_1 & \bQ_2
\end{bmatrix}
\begin{bmatrix}
r_{11} & \br_{12}^\top\\
\bzero & \bR_{22}
\end{bmatrix}
=
\begin{bmatrix}
r_{11}\bq_1 & \bq_1\br_{12}^\top +\bQ_2\bR_{22}
\end{bmatrix},
$$
where the matrix $\bQ_2\in \real^{m\times (n-1)}$ consists of mutually orthonormal columns   and $\bR_{22}\in \real^{(n-1)\times (n-1)}$ is upper triangular. 
Consequently, $\bQ_2\bR_{22} $ represents the reduced QR decomposition of $\bA_2-\bq_1\br_{12}^\top$, which reveals a recursive algorithm for the reduced QR decomposition of $\bA$. 
This approach is equivalent to the MGS method that subtracts each component in the span of $\{\bq_1, \bq_2, \ldots, \bq_{k-1}\}$ when computing column $k$ of $\bQ$ (i.e., equality ($*$) in Algorithm~\ref{alg:qr-mgs-right}). The process is described in Algorithm~\ref{alg:qr-mgs-fulll-rowwise-recursive}.
\begin{algorithm}[H] 
\caption{MGS (\textcolor{mylightbluetext}{Row-Wise and Recursively})=Algorithm~\ref{alg:qr-mgs-right}}
\label{alg:qr-mgs-fulll-rowwise-recursive}
\begin{algorithmic}[1] 
\Require $\bA\in \real^{m\times n}$ with full column rank;
\For{$k=1$ to $n$}  \Comment{i.e., compute $k$-th column of $\bQ$ and $k$-th row of $\bR$}
\State $\ba_1\leftarrow\bA[:,1]$; \Comment{$1$-st column of $\bA\in \real^{m\times (n-k+1)}$}
\State $r_{kk}\leftarrow\normtwo{\ba_1}$;\Comment{$\ba_1\in \real^{m\times 1}$}
\State $\bq_k \leftarrow \ba_1/r_{kk}$;
\State $\bA_2\leftarrow\bA[:,2:n]\in \real^{m\times (n-k)}$; \Comment{2-nd to $n$-th column of $\bA$}
\State $\br_{k2}^\top\leftarrow\bq_k^\top\bA_2$; \Comment{$\br_{k2}^\top\in \real^{1\times (n-k)}$}
\State $\bA\leftarrow\bA_2-\bq_k\br_{k2}^\top$; \Comment{$\bA \in \real^{m\times (n-k)}$}
\EndFor 
\State Output $\bQ=[\bq_1, \ldots, \bq_n]$ and $\bR$ with entry $(i,k)$ being $r_{ik}$.
\end{algorithmic} 
\end{algorithm}

More compactly, Algorithm~\ref{alg:qr-mgs-fulll-rowwise-recursive} can be equivalently stated as Algorithm~\ref{alg:mgs_rowwise-recursive_comp}.
\begin{algorithm}[H] 
\caption{MGS (\textcolor{mylightbluetext}{Row-Wise and Recursively})=Algorithm~\ref{alg:qr-mgs-right}=Algorithm~\ref{alg:qr-mgs-fulll-rowwise-recursive}}
\label{alg:mgs_rowwise-recursive_comp}
\begin{algorithmic}[1] 
\Require $\bA\in \real^{m\times n}$ with full column rank;
\For{$k=1$ to $n$}  \Comment{i.e., compute $k$-th column of $\bQ$ and $k$-th row of $\bR$}
\State $\bq_k\leftarrow\bA[:,k]/\normtwo{\bA[:,k]}$; \Comment{Normalize $k$-th column of $\bA\in \real^{m\times n}$}
\State $\br_k^\top \leftarrow \bq_k^\top \bA$; \Comment{$\br_i^\top\in\real^{1\times n}$, $k$-th row of $\bR$}
\State $\bA \leftarrow\bA-\bq\br^\top$; \Comment{MGS step}
\EndFor 
\State Output $\bQ=[\bq_1, \ldots, \bq_n]$ and $\bR$ with entry $k$-th row being $\br_k^\top$.
\end{algorithmic} 
\end{algorithm}

To enhance the orthogonality of the ${\bq_i}$'s and improve numerical accuracy, an additional re-orthonormalization step can be performed. This step becomes necessary as the basis vectors generated tend to lose their orthonormality during the process. 
The re-orthonormalization steps are highlighted in blue in Algorithm~\ref{alg:mgs_recursive_comp_moreortho}.
\begin{algorithm}[H] 
\caption{MGS (\textcolor{mylightbluetext}{Row-Wise, Re-Orthonormalization} based on Algorithm~\ref{alg:mgs_rowwise-recursive_comp})}
\label{alg:mgs_recursive_comp_moreortho}
\begin{algorithmic}[1] 
\Require $\bA\in \real^{m\times n}$ with full column rank;
\For{$k=1$ to $n$}  \Comment{i.e., compute $k$-th column of $\bQ$ and $k$-th row of $\bR$}
\State $\bq_k\leftarrow\bA[:,k]/\normtwo{\bA[:,k]}$; \Comment{Normalize $k$-th column of $\bA\in \real^{m\times n}$}
\State \textcolor{mylightbluetext}{$\bq\leftarrow\bq - \bQ_{k-1}\bQ_{k-1}^\top \bq$}; \Comment{$\bQ_{k-1} =[\bq_1,\bq_2, \ldots, \bq_{k-1}]\in\real^{m\times (k-1)} $}
\State \textcolor{mylightbluetext}{$\bq \leftarrow \bq/\normtwo{\bq}$};
\State $\br_k^\top \leftarrow \bq_k^\top \bA$; \Comment{$\br_i^\top\in\real^{1\times n}$, $k$-th row of $\bR$}
\State $\bA \leftarrow\bA-\bq\br^\top$; \Comment{MGS step}
\EndFor 
\State Output $\bQ=[\bq_1, \ldots, \bq_n]$ and $\bR$ with entry $k$-th row being $\br_k^\top$.
\end{algorithmic} 
\end{algorithm}

\section{Computing  Full QR Decomposition via  Gram-Schmidt Process}\label{section:silentcolu_qrdecomp}
\begin{figure}[h]
\centering  
\vspace{-0.35cm} 
\subfigtopskip=2pt 
\subfigbottomskip=2pt
\subfigcapskip=-5pt 
\subfigure[Reduced QR decomposition.]{\label{fig:gphalf}
	\includegraphics[width=0.47\linewidth]{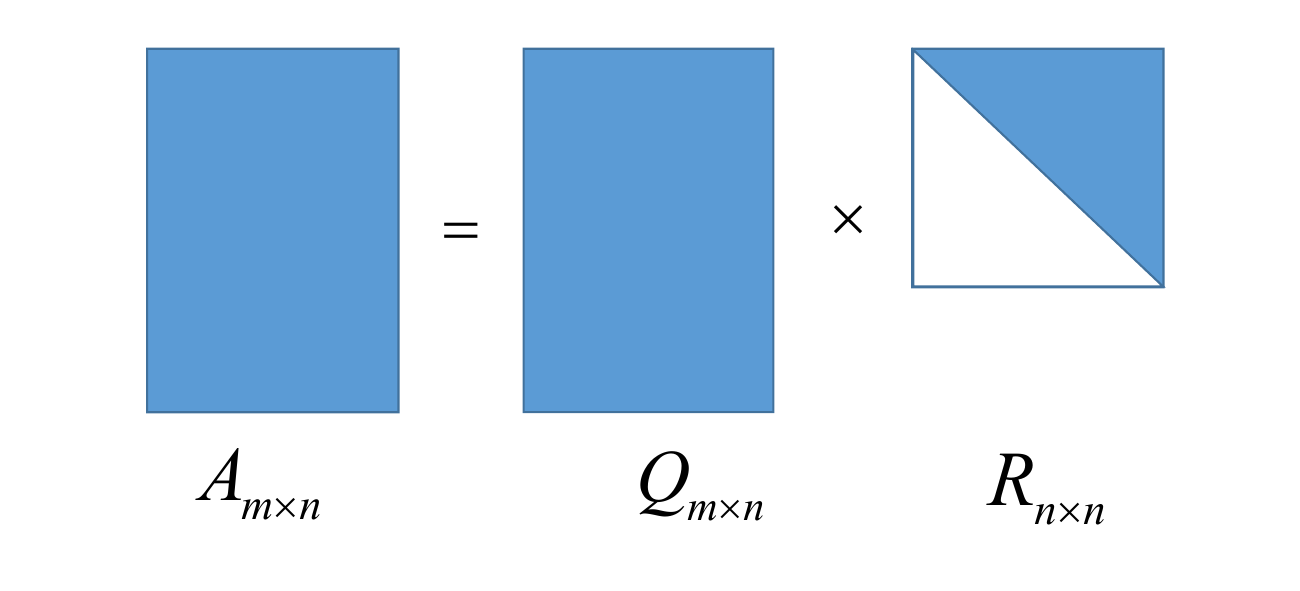}}
\quad 
\subfigure[Full QR decomposition.]{\label{fig:gpall}
	\includegraphics[width=0.47\linewidth]{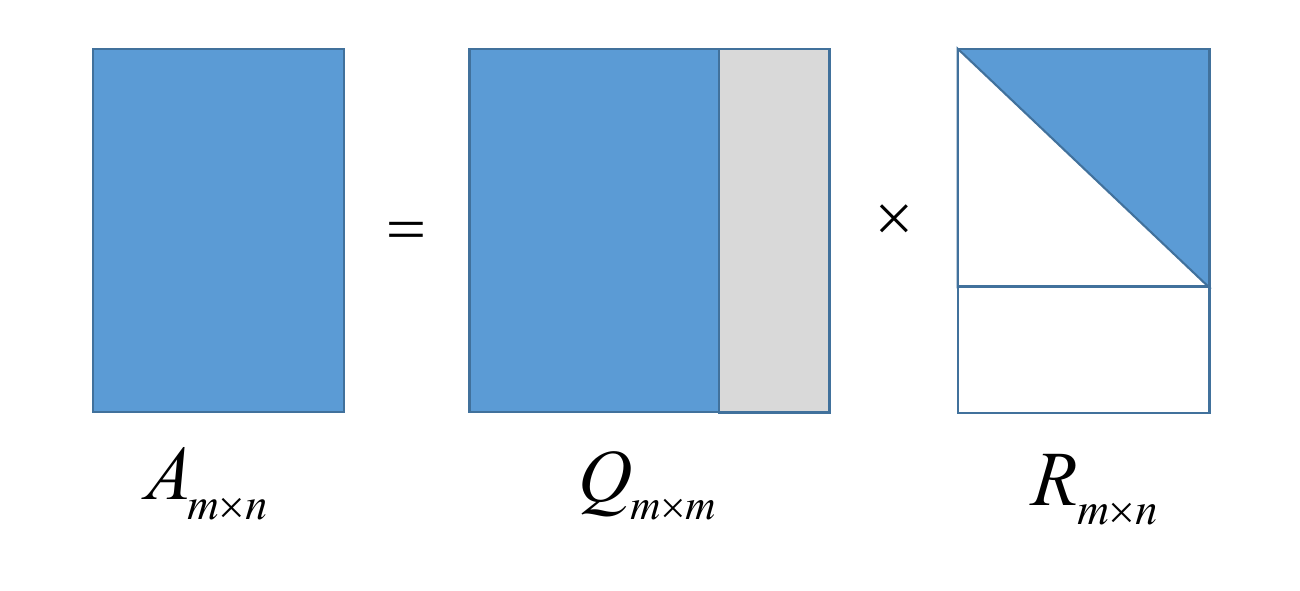}}
\caption{Comparison between the reduced and full QR decompositions. White entries are zero, and \textcolor{mylightbluetext}{blue} entries are not necessarily zero. \textcolor{mydarkgray}{Gray} columns denote silent columns.}
\label{fig:qr-comparison}
\end{figure}

A full QR decomposition of an $m\times n$ matrix with linearly independent columns extends the process by appending additional $m-n$ orthonormal columns to $\bQ$, thereby transforming it into an $m\times m$ orthogonal matrix. Simultaneously, $\bR$ is augmented with rows of zeros to attain an $m\times n$ upper triangular matrix.
We refer to the additional columns in $\bQ$ as \textbf{silent columns} and the additional rows in $\bR$ as \textbf{silent rows}. The comparison between the reduced   and the full QR decompositions is shown in Figure~\ref{fig:qr-comparison}, where silent columns in $\bQ$ are represented in \textcolor{mydarkgray}{gray}, blank entries denote zero elements, and \textcolor{mylightbluetext}{blue} entries  indicate elements that may not necessarily be zero.

\section{Dependent Columns}\label{section:dependent-gram-schmidt-process}
In our earlier discussions, we assumed that the matrix $\bA$ has linearly independent columns. While this condition simplifies our analysis, it is not strictly required for all scenarios. 
Consider step $k$ Algorithm~\ref{alg:reduced-qr}, where $\ba_k$ lies in the plane spanned by $\bq_1, \bq_2, \ldots, \bq_{k-1}$ (which is equivalent to the space spanned by $\ba_1, \ba_2, \ldots, \ba_{k-1}$), indicating that the vectors $\ba_1, \ba_2, \ldots, \ba_k$ are dependent. 
Then $r_{kk}$ will be zero and $\bq_k$ cannot be determined  due to division by zero.
In such cases, we can arbitrarily choose $\bq_k$ to be any normalized vector that is orthogonal to $\cspace([\bq_1, \bq_2, \ldots, \bq_{k-1}])$ and proceed with the Gram-Schmidt process. 
Again, when dealing with a matrix $\bA$ that has dependent columns, we have both reduced and full QR decomposition algorithms. We reformulate the step $k$ in the algorithm as follows:
$$
\bq_k=\left\{
\begin{aligned}
&\big(\ba_k-\sum_{i=1}^{k-1}r_{ik}\bq_i\big)/r_{kk}, \qquad r_{ik}=\bq_i^\top\ba_k, r_{kk}=\normtwo{\ba_k-\sum_{i=1}^{k-1}r_{ik}\bq_i}, &\mathrm{if\,} r_{kk}\neq0, \\
&\text{pick  one in\,}\cspace^{\bot}([\bq_1, \bq_2, \ldots, \bq_{k-1}]),\text{ and normalize},\qquad &\mathrm{if\,} r_{kk}=0.
\end{aligned}
\right.
$$

This idea can be further extended such that,  when $\bq_k$ does not exist, we simply skip the current step and add the silent columns at the end of the process.  Consequently, the QR decomposition of a matrix with dependent columns is not unique. However, by adhering to a systematic process or a methodical procedure, the QR decomposition of any matrix remains unique.

Moreover, this insight also aids in determining the linear independence of a set of vectors. 
Whenever $r_{kk}$ in Algorithm~\ref{alg:reduced-qr} becomes zero, we report the vectors $\ba_1, \ba_2, \ldots, \ba_k$ are dependent and terminate the algorithm for ``\textit{independence checking}."

\section{QR with Column Pivoting: Column-Pivoted QR (CPQR)}\label{section:cpqr}

Suppose $\bA$ has linearly dependent columns; a column-pivoted QR (CPQR) decomposition can be found as follows.
\begin{theoremHigh}[Column-Pivoted QR Decomposition\index{Column-pivoted QR (CPQR)}]\label{theorem:rank-revealing-qr-general}
Let $\bA=[\ba_1, \ba_2, \ldots, \ba_n]$ be an  $m\times n$ matrix, with $m\geq n$ and rank $r$. Then, it can be factored as 
$$
\bA\bP = \bQ
\begin{bmatrix}
\bR_{11} & \bR_{12} \\
\bzero   & \bzero 
\end{bmatrix},
$$
where $\bR_{11} \in \real^{r\times r}$ is upper triangular, $\bR_{12} \in \real^{r\times (n-r)}$, $\bQ\in \real^{m\times m}$ is an orthogonal matrix, and $\bP$ is a permutation matrix. This is also known as the \textbf{full} CPQR decomposition. Similarly, the \textbf{reduced} version is given by 
$$
\bA\bP = \bQ_r
\begin{bmatrix}
\bR_{11} & \bR_{12} \\ 
\end{bmatrix},
$$
where $\bR_{11} \in \real^{r\times r}$ is upper triangular, $\bR_{12} \in \real^{r\times (n-r)}$, $\bQ_r\in \real^{m\times r}$ contains orthonormal columns, and $\bP$ is a permutation matrix.
\end{theoremHigh}
\subsection{A Simple CPQR via CGS}
\begin{figure}[H]
\centering  
\vspace{-0.35cm} 
\subfigtopskip=2pt 
\subfigbottomskip=2pt 
\subfigcapskip=-5pt 
\subfigure[Reduced CPQR decomposition.]{\label{fig:gphalf-rank-reveal}
	\includegraphics[width=0.475\linewidth]{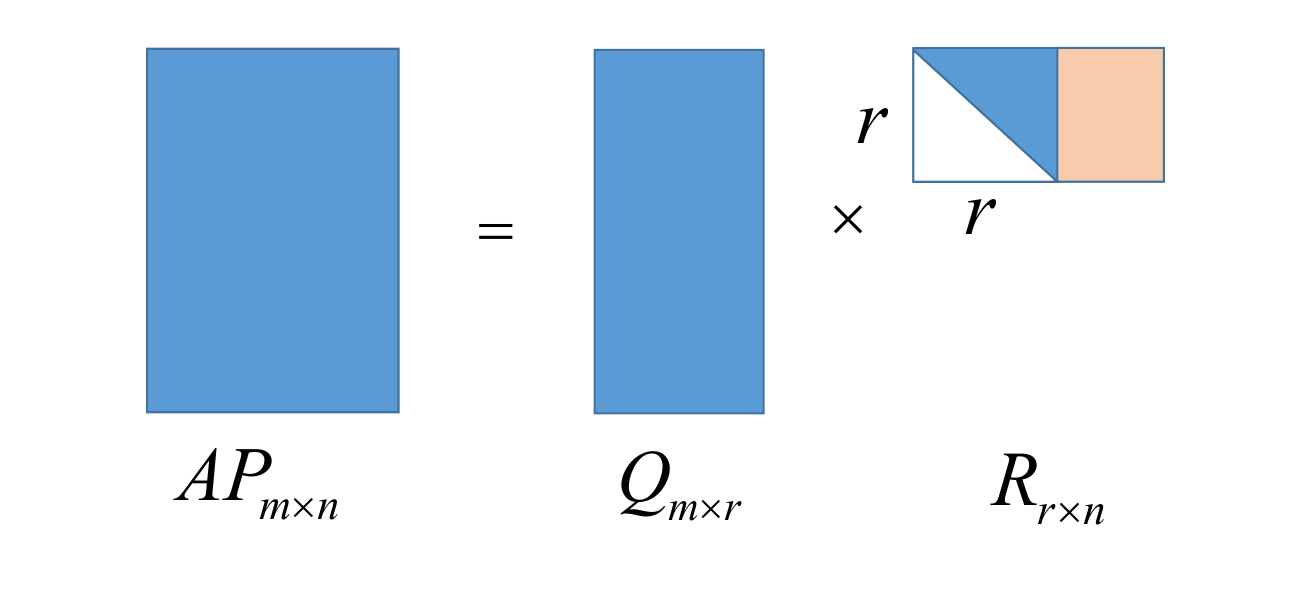}}
\quad 
\subfigure[Full CPQR decomposition.]{\label{fig:gpall-rank-reveal}
	\includegraphics[width=0.475\linewidth]{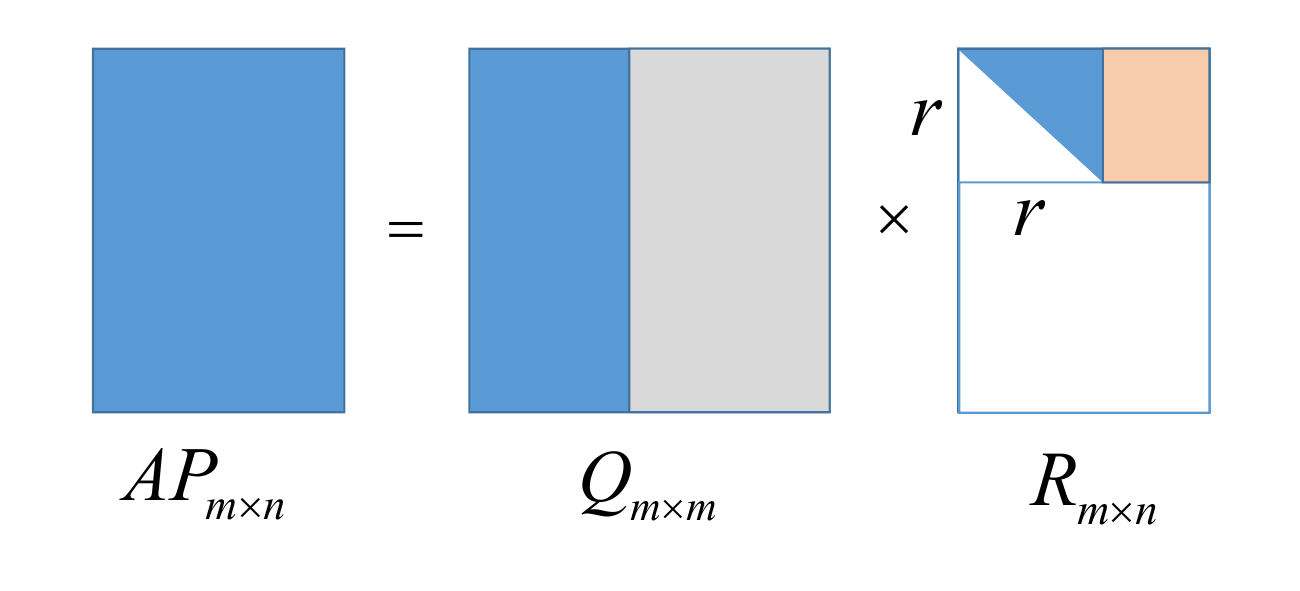}}
\caption{Comparison between the reduced and full CPQR decompositions.}
\label{fig:qr-comparison-rank-reveal}
\end{figure}

The CPQR decomposition can be computed using the classical Gram-Schmidt process.
Following the QR decomposition for dependent columns, when $r_{kk}=0$, it indicates that the column $k$ of $\bA$ is dependent on the previous $k-1$ columns. 
In such situations, we perform a column permutation, moving this column to the last position, and then continue with the Gram-Schmidt process.
We notice that $\bP$ represents the permutation matrix that interchanges the dependent columns into the last $n-r$ columns. 
If the first $r$ columns of $\bA\bP$ are denoted as  $[\widehat{\ba}_1, \widehat{\ba}_2, \ldots, \widehat{\ba}_r]$,  the span of these columns is equivalent to the span of $\bQ_r$ (in the reduced version) or the span of $\bQ[:,1:r]$ (in the full version):
$$
\cspace([\widehat{\ba}_1, \widehat{\ba}_2, \ldots, \widehat{\ba}_r]) = \cspace(\bQ_r) = \cspace(\bQ_{:,1:r}).
$$
And $\bR_{12}$ represents a matrix that recovers the dependent $n-r$ columns from the column space of $\bQ_r$ or the column space of $\bQ[:,1:r]$. The comparison between the reduced and full CPQR decompositions is illustrated in Figure~\ref{fig:qr-comparison-rank-reveal}, where silent columns in $\bQ$ are denoted in \textcolor{mydarkgray}{gray}, blank entries are zero, and \textcolor{mylightbluetext}{blue}/\textcolor{orange}{orange} entries are elements that are not necessarily zero.

\begin{algorithm}[h] 
\caption{\textcolor{mylightbluetext}{Simple} Reduced CPQR Decomposition via CGS} 
\label{alg:reduced-qr-rank-revealing} 
\begin{algorithmic}[1] 
	\Require 
	Matrix $\bA$ with size $m\times n$ and $m\geq n$;  	
	\State $cnt = 0$; \Comment{i.e., the count for the permutations}
	\State $\bq_1   \leftarrow \ba_1/r_{11},  r_{11}=\normtwo{\ba_1}$; \Comment{i.e., the first column of $\bR_{11}$}
	\For{$k=2$ to $n$}   \Comment{i.e., compute column $k$ of $\bR_{11}$}
	\State Set the initial value $r_{kk}=0$;
	\While{$r_{kk} ==0$} \Comment{the column is dependent if $r_{kk}$ is equal to 0}
	\State $r_{ik}\leftarrow\bq_i^\top\ba_k$, $\forall i \in \{1, 2, \ldots, k-1\}$;\Comment{first $k-1$ elements in column $k$ of $\bR_{11}$}
	\State $r_{kk}\leftarrow\normtwo{\ba_k-\sum_{i=1}^{k-1}r_{ik}\bq_i}$;\Comment{$k$-th element in column $k$ of $\bR_{11}$}
	\State $\bq_k \leftarrow (\ba_k-\sum_{i=1}^{k-1}r_{ik}\bq_i)/r_{kk}$ if $r_{kk}\neq 0$;
	\If{$r_{kk}==0$}
	\State $cnt = cnt+1$;
	\State Permute the column $k$ to the last column;
	\State i.e., $[\ba_k, \ba_{k+1}, \ldots, \ba_n] \rightarrow [\ba_{k+1}, \ba_{k+2}, \ldots, \ba_n, \ba_k]$;
	\EndIf
	\EndWhile
	\If{$k+cnt == n$}
	\State rank $r\leftarrow k$;
	\For{$k=r+1$ to $n$} \Comment{i.e., compute column $k$ of $\bR_{12}$}
	\State $r_{ik}\leftarrow\bq_i^\top\ba_k$, $\forall i \in \{1, 2, \ldots, \textcolor{mylightbluetext}{r}\}$;
	\EndFor
	\EndIf
	\State Output rank $r$, output $\bR_{11}, \bR_{12}$, $\bQ_r=[\bq_1, \bq_2, \ldots, \bq_r]$. And exit the loop.
	\EndFor
\end{algorithmic} 
\end{algorithm}

The reduced CPQR algorithm is formulated in Algorithm~\ref{alg:reduced-qr-rank-revealing}.  Additionally, it is straightforward to determine the last $n-r$ columns of $\bQ$ in the orthogonal complement of $\cspace(\bQ_r)$.
Note that the step 6 of Algorithm~\ref{alg:reduced-qr-rank-revealing} can be rewritten as $\bQ_{k-1}^\top \ba_k$, where $\bQ_{k-1}=[\bq_1, \ldots, \bq_{k-1}]$.

\begin{theorem}[Algorithm Complexity: Reduced CPQR]\label{theorem:qr-reduced-rank-revealing}
Algorithm~\ref{alg:reduced-qr-rank-revealing} requires $\sim 6mnr-4mr^2$ flops to compute a CPQR decomposition of an $m\times n$ matrix with $m\geq n$ and rank $r$.
\end{theorem}

\begin{proof}[of Theorem~\ref{theorem:qr-reduced-rank-revealing}]
From Theorem~\ref{theorem:qr-reduced}, to compute $\bR_{11}$, we just need to replace $n$ in Theorem~\ref{theorem:qr-reduced} by $r$ such that the complexity of computing $\bR_{11}$ is \underline{$2mr^2$} flops if keep only the leading term.

To compute $\bR_{12}$, there are $r\times (n-r)$ values, each taking $2m-1$ flops ($m$ multiplications and $m-1$ additions) according to step 18. That is \underline{$r(n-r)(2m-1)$} flops. 

The upper bound on steps 6, 7, and 8 is to set $k-1=r$ in these steps. 
This results in \underline{$(2m-1)r$} flops for step 6 (involving $mr$ multiplications and $(m-1)r$ additions), \underline{$2m(r+1)$} flops for step 7 (involving $mr$ multiplications, $(r-1)m$ additions, $m$ subtractions, and $2m$ for the norm), and \underline{$m$} flops for step 8. The total complexity of steps 6, 7, and 8 is thus \underline{$4mr+3m-r$} flops. And there are $n-r$ such iterations, the total number of flops required to find the dependent columns is \underline{$(4mr+3m-r)(n-r)$}.

Therefore, the final complexity is 
$$
2mr^2 + r(n-r)(2m-1)+ (4mr+3m-r)(n-r),
$$ 
which simplifies to $6mnr-4mr^2$ flops if we keep only the leading term.
\end{proof}
We notice that when $r=n$, the complexity of Algorithm~\ref{alg:reduced-qr-rank-revealing} is $2mn^2$, which matches  the complexity of Algorithm~\ref{alg:reduced-qr}.

\subsection{A Practical CPQR via CGS}\label{section:practical-cpqr-cgs}
We notice that the simple CPQR algorithm pivots the first $r$ independent columns to the first $r$ columns of $\bA\bP$. Let $\bA_1$ denote the first $r$ columns of $\bA\bP$, and let $\bA_2$ denote the remaining columns. Then, from the full CPQR, we have
$$
\bA\bP=
[\bA_1, \bA_2] = 
\bQ \begin{bmatrix}
\bR_{11} & \bR_{12} \\
\bzero   & \bzero 
\end{bmatrix}
=\left[
\bQ \begin{bmatrix}
\bR_{11}  \\
\bzero   
\end{bmatrix}
,
\bQ \begin{bmatrix}
\bR_{12}  \\
\bzero   
\end{bmatrix}
\right]
.
$$ 
It is not hard to see that 
$$
\normf{\bA_2} = \normf{\bQ \begin{bmatrix}
		\bR_{12}  \\
		\bzero   
\end{bmatrix}}
=
\normf{\begin{bmatrix}
		\bR_{12}  \\
		\bzero   
\end{bmatrix}}
=
\normf{\bR_{12}},
$$
where the penultimate equality is derived from the preservation of the matrix norm under orthogonal transformations (Proposition~\ref{proposition:frobenius-orthogonal-equi}). Therefore, the norm of $\bR_{12}$ is determined by the norm of $\bA_2$. When favoring a well-conditioned CPQR, it is desirable for $\bR_{12}$ to have a small norm. A practical approach for CPQR decomposition is to initially permute the columns of the matrix $\bA$ such that they are ordered in decreasing vector norm:
$$
\widetildebA = \bA\bP_0 = [\ba_{j_1}, \ba_{j_2}, \ldots, \ba_{j_n}],
$$
where $\{j_1, j_2, \ldots, j_n\}$ represents a permuted index set of $\{1,2,\ldots, n\}$, and 
$$
\normtwo{\ba_{j_1}}\geq  \normtwo{\ba_{j_2}}\geq  \ldots\geq  \normtwo{\ba_{j_n}}.
$$
Then, we apply the ``simple" reduced CPQR decomposition to $\widetildebA$ such that $\widetildebA \bP_1= \bQ_r[\bR_{11}, \bR_{12}]$. The ``practical" reduced CPQR of $\bA$ is then recovered as
$$
\bA\underbrace{\bP_0\bP_1}_{\bP} =\bQ_r[\bR_{11}, \bR_{12}].
$$
When computing the $\ell_2$ vector norm (i.e., inner product, Appendix~\ref{section:vector-norm}), extra $n(2m-1)$ flops are required to compute the norms of the  $n$ column vectors and $\frac{n(n-1)}{2}$ comparisons are needed to establish the order of the norms.  

\subsection{A Practical CPQR via MGS}
Based on the recursive MGS in Algorithm~\ref{alg:qr-mgs-fulll-rowwise-recursive}, we can also develop a practical CPQR. The algorithm is presented  in Algorithm~\ref{alg:cpqr-partial-fact-cpqr}, where the sole difference from Algorithm~\ref{alg:qr-mgs-fulll-rowwise-recursive} is highlighted in the \textcolor{mylightbluetext}{blue} text: we permute the column with the largest norm into the first column.

\begin{algorithm}[htp] 
\caption{\textcolor{mylightbluetext}{Practical} CPQR via MGS (\textcolor{mylightbluetext}{Row-Wise and Recursively}). 
The algorithm is derived from Algorithm~\ref{alg:qr-mgs-fulll-rowwise-recursive} and a similar procedure can be derived based on Algorithm~\ref{alg:mgs_rowwise-recursive_comp} and Algorithm~\ref{alg:mgs_recursive_comp_moreortho}.
}
\label{alg:cpqr-partial-fact-cpqr}
\begin{algorithmic}[1] 
	\Require $\bA\in \real^{m\times n}$ with \textcolor{mylightbluetext}{exact} rank $r$;
	\For{$k=1$ to $n$}  \Comment{i.e., compute $k$-th column of $\bQ$ and $k$-th row of $\bR$}
	\State \textcolor{mylightbluetext}{Find the column with largest norm in $\bA$, and permute to first column;}
	\State $\ba_1\leftarrow\bA[:,1]$; \Comment{$1$-st column of $\bA\in \real^{m\times (n-k+1)}$}
	\State $r_{kk}\leftarrow\normtwo{\ba_1}$;\Comment{$\ba_1\in \real^{m\times 1}$}
	\State $\bq_k \leftarrow \ba_1/r_{kk}$;
	\State $\br_{k2}^\top\leftarrow\bq_k^\top\bA_2$; \Comment{$\bA_2=\bA[:,2:n]\in \real^{m\times (n-k)}$, $\br_{k2}^\top\in \real^{1\times (n-k)}$}
	\State $\bA\leftarrow\bA_2-\bq_k\br_{k2}^\top$; \Comment{$\bA \in \real^{m\times (n-k)}$}
	\State \textcolor{mylightbluetext}{Exit when $r_{kk}=0$;}
	\EndFor 
	\State Output $\bQ=[\bq_1, \ldots, \bq_n]$ and $\bR$ with entry $(i,k)$ being $r_{ik}$.
\end{algorithmic} 
\end{algorithm}

The main difference is in that, in each iteration, we need to calculate  the norms of all the (remaining) columns of $\bA$ rather than  computing the norms all at once as done in CGS. At iteration $k$, we need to compute the reduced QR decomposition of a matrix of size $m\times (n-k+1)$ if the original matrix $\bA$ is of size $m\times n$. That is, extra $(n-k+1)(2m-1)$ flops are required to proceed with the CPQR via MGS. Let $f(k)=(n-k+1)(2m-1)$; simple calculation can show that additional complexity for CPQR via MGS is:
\begin{equation}\label{equation:mgs-cpqr-extra1}
\text{extra cost = }f(1)+f(2)+\ldots +f(n)\sim mn^2 \text{  flops},
\end{equation}
if only keep the leading term. This costs more than $n(2m-1)$ flops in the ``practical" CPQR via CGS.

However, this additional cost in CPQR via MGS can be mitigated to some extent.
Suppose the column partition of $\bA\in \real^{m\times n}$ is $\bA=[\ba_1, \ba_2, \ldots, \ba_n]$, and let the squared norm of each column be given in the vector 
$$
\bl_a = 
\begin{bmatrix}
l_1 \\
l_2 \\
\vdots \\
l_n
\end{bmatrix}=
\begin{bmatrix}
\normtwo{\ba_1}^2\\
\normtwo{\ba_2}^2\\
\vdots \\
\normtwo{\ba_n}^2\\
\end{bmatrix}.
$$
Suppose further that $\bq\in \real^m$ is a unit-length vector such that $\bq^\top\bq=1$, and $\br\in \real^n$ is a vector given by 
$$
\br = \bA^\top\bq
=
\begin{bmatrix}
r_1 \\
r_2 \\
\vdots \\
r_n
\end{bmatrix}
. \gap \text{(similar to the step 6 of  Algorithm~\ref{alg:cpqr-partial-fact-cpqr})}
$$
Let further $\bB=\bA-\bq\br^\top=[\bb_1, \bb_2, \ldots, \bb_n]$ (similar to the step 7 of  Algorithm~\ref{alg:cpqr-partial-fact-cpqr}). 
The vector representing the squared lengths of $\bB$ is given by
$$
\bl_b = 
\begin{bmatrix}
s_1 \\
s_2 \\
\vdots \\
s_n
\end{bmatrix}=
\begin{bmatrix}
\normtwo{\bb_1}^2\\
\normtwo{\bb_2}^2\\
\vdots \\
\normtwo{\bb_n}^2\\
\end{bmatrix}
=
\begin{bmatrix}
l_1 -r_1^2\\
l_2 -r_2^2\\
\vdots \\
l_n-r_n^2
\end{bmatrix}.
$$
This can be easily verified since $\bb_i = \ba_i -r_{i}\bq =\ba_i-(\ba_i^\top\bq) \bq$ such that 
$$
\normtwo{\bb_i}^2 = \normtwo{\ba_i -r_{i}\bq}^2=( \ba_i -r_{i}\bq)^\top( \ba_i -r_{i}\bq)=l_i-r_i^2.
$$
Come back to step 2 of Algorithm~\ref{alg:cpqr-partial-fact-cpqr}, suppose we have computed the squared norms of the columns from the original matrix $\bA\in \real^{m\times n}$ (which requires $n(2m-1)$, the same as that in the ``practical" CPQR via CGS). The squared norms of the columns from $\bA_2-\bq_1\br_{12}^\top$ (suppose $k=1$ in step 7 of Algorithm~\ref{alg:cpqr-partial-fact-cpqr}) can be obtained with an additional $2(n-1)$ flops. 
Over the course of the $n$ iterations, the total cost is $2(n-1)+2(n-2)+\ldots +2(1)=n^2-n$ flops. This is significantly less than the complexity of $\sim mn^2$ in Equation~\eqref{equation:mgs-cpqr-extra1}.

\index{Partial factorization}
\index{Rank-revealing QR decomposition (RRQR)}
\subsection{Partial Factorization for CPQR: Extra Bonus of CPQR via MGS}\label{section:partial-cpqr-mgs}
An additional advantage of the CPQR via MGS is the ability to perform \textit{partial factorization} at a certain stage.  
We notice that in step 2 of Algorithm~\ref{alg:cpqr-partial-fact-cpqr}, we permute the column with the largest norm into the first column; and step 4 of Algorithm~\ref{alg:cpqr-partial-fact-cpqr} involves computing the norm and placing it along the main diagonal of the upper triangular matrix $\bR$. 
In the case where $\bA$ has full rank $n$, the values of $r_{kk}$ in all iterations will be positive. 
In the scenario where $\bA$ has an ``exact" rank $r$, the algorithm will stop after iteration $r$. 
However, when $\bA$ has an ``effective" rank $r$ with rank-deficiency
\footnote{\textit{Effective rank}, or also known as the \textit{numerical rank}. Assume the $i$-th largest singular value of $\bA$ is denoted by $\sigma_i(\bA)$. Then if $\sigma_r(\bA)\gg \sigma_{r+1}(\bA)\approx 0$, $r$ is known as the numerical rank of $\bA$. The singular value of matrix $\bA$ will be introduced in the SVD section (Chapter~\ref{chapter:SVD}). Whereas, when $\sigma_i(\bA)>\sigma_{r+1}(\bA)=0$, it is known as having \textit{exact rank} $r$ as we have used in most of our discussions.}, the algorithm can proceed well in the first $r$ iterations since $\{r_{11}, r_{22}, \ldots, r_{rr}\}$ are relatively large in value that are far from 0. When it comes to iteration $k=r+1, \ldots$,  the value of $r_{kk}$ becomes small, which means the column $k$ of $\bA\bP$ has a small component in the direction of $\bq_k$ and is ``almost" dependent on the previous $(k-1)$ orthonormal basis $\{\bq_1, \bq_2, \ldots, \bq_{k-1}\}$.  This situation holds for the remaining columns $\{k+1, k+2, \ldots\}$ columns since $r_{kk}\geq r_{k+1, k+1}\geq \ldots r_{nn}$. The partial factorization CPQR via MGS is detailed  in Algorithm~\ref{alg:cpqr-partial-fact-cpqr-partial}.
This concept is closely related to the \textit{rank-revealing QR decomposition (RRQR)}~\footnote{The rigorous definition of a RRQR is discussed in Definition~\ref{definition:rrqr_in_svd} using the concept of numerical rank of a matrix.} that we will introduce in the subsequent sections. The algorithm will result in the factorization 
\begin{equation}
\bA\bP = 
\bQ\bR=
\bQ
\begin{bmatrix}
	\bR_{11} & \bR_{12} \\
	\bzero & \bR_{22}
\end{bmatrix},
\end{equation}
where $\bR_{22}$ is small in norm.

However, this partial factorization does not work for the CPQR via CGS. In the ``practical" CPQR with CGS (Algorithm~\ref{alg:reduced-qr-rank-revealing}, where $\bA$ is permuted all at once at the beginning of the procedure), it becomes inefficient when $r \ll \min(m, n)$. 
This is because we only permute the column with the largest norm to the beginning of $\bA\bP$. A value $r_{kk}$ close to 0 only means that column $k$ of $\bA\bP$ is almost dependent on $\{\bq_1, \bq_2, \ldots, \bq_{k-1}\}$, but it does not imply the remaining columns $\{k+1, k+2, \ldots\}$ of $\bA\bP$ are also almost dependent on $\{\bq_1, \bq_2, \ldots, \bq_k\}$.

\begin{algorithm}[htp] 
\caption{\textcolor{mylightbluetext}{Practical} and \textcolor{mylightbluetext}{Partial} CPQR via MGS (\textcolor{mylightbluetext}{Row-Wise and Recursively}).
The algorithm is based on Algorithm~\ref{alg:qr-mgs-fulll-rowwise-recursive} and similar procedure can be derived based on Algorithm~\ref{alg:mgs_rowwise-recursive_comp} and Algorithm~\ref{alg:mgs_recursive_comp_moreortho}.
}
\label{alg:cpqr-partial-fact-cpqr-partial}
\begin{algorithmic}[1] 
\Require $\bA\in \real^{m\times n}$ with rank-deficiency;
\State Select a stopping criteria $\delta$;
\For{$k=1$ to $n$}  \Comment{i.e., compute $k$-th column of $\bQ$ and $k$-th row of $\bR$}
\State \textcolor{mylightbluetext}{Find the column with largest norm in $\bA$, and permute to first column;}
\State $\ba_1\leftarrow \bA[:,1]$; \Comment{$1$-st column of $\bA\in \real^{m\times (n-k+1)}$}
\State $r_{kk}\leftarrow\normtwo{\ba_1}$;\Comment{$\ba_1\in \real^{m}$}
\State $\bq_k \leftarrow \ba_1/r_{kk}$;
\State $\br_{k2}^\top\leftarrow\bq_k^\top\bA_2$; \Comment{$\bA_2=\bA[:,2:n]\in \real^{m\times (n-k)}$, $\br_{k2}^\top\in \real^{1\times (n-k)}$}
\State $\bA\leftarrow\bA_2-\bq_k\br_{k2}^\top$; \Comment{$\bA \in \real^{m\times (n-k)}$}
\State \textcolor{mylightbluetext}{Exit when $r_{kk}<\delta$, set effective rank $r=k$;}
\EndFor 
\State Output $\bQ=[\bq_1, \ldots, \bq_n]$, $\bR$ with entry $(i,k)$ being $r_{ik}$, and \textcolor{mylightbluetext}{effective rank $r$}.
\end{algorithmic} 
\end{algorithm}

\index{Low-rank approximation}
\subsection{Low-Rank Approximation via MGS}
For  low-rank approximation using the Gram-Schmidt process, given a rank-$r$ matrix $\bA\in\real^{m\times n}$, the goal is to find an orthonormal set of vectors $\{\bq_i\}_{i=1}^k$ that form a set of basis vectors for expressing the columns of $\bA$, where $k<r$. That is, let $\bQ_k = [\bq_1, \bq_2, \ldots, \bq_k]$, we want to find the $\bQ_k$ such that 
$$
\text{LR Approx.: }\gap \normtwo{\bA - \bQ_k\bQ_k^\top\bA}=
\min\{\normtwo{\bA-\bB}:\rank(\bB)=k\},
$$
where $\normtwo{\cdot}$ represents the spectral norm of a matrix, which we have introduced in Definition~\ref{definition:spectral_norm} or will be further discussed in Appendix~\ref{appendix:matrix-norm-sect2}.

Based on Algorithm~\ref{alg:cpqr-partial-fact-cpqr-partial}, we can obtain the decomposition at $k$-th iteration:
\begin{equation}\label{equation:cpqr_los_rank}
\bA\bP_k = \bQ_k\bR_k + [\bzero,  \widetildebA_k],
\end{equation}
where $\bP_k\in\real^{n\times n}$, $\bR_k\in\real^{k\times n}$, $\widetildebA_k\in\real^{m\times (n-k)}$, and $\bQ_k=[\bq_1, \bq_2, \ldots, \bq_k]\in\real^{m\times k}$. 
Specifically, $\bP_k$ accounts for all permutations up to the $k$-th iteration, $\bR_k$ consists of the first $k$ rows of the upper triangular matrix, $\bQ_k$ comprises $k$ orthonormal basis vectors, and $\widetildebA_k$ contains the remaining columns that have not yet been chosen as pivot columns, i.e., $\widetildebA_k = \bA_{:,k+1:n} $ at the $k$-th iteration.
Equation~\eqref{equation:cpqr_los_rank} implies 
$$
\bA = \bQ_k\bR_k\bP_k^\top +[\bzero, \widetildebA_k]\bP_k^\top.
$$
Suppose our objective is to compute a low-rank approximation of $\bA$ that achieves a precision of $\delta$ in terms of the spectral norm. 
Then after the $k$-th step, we set $\bB=\bQ_k\bR_k\bP_k^\top$ 
\footnote{It can be shown that $\rank(\bB)=k$.}.
By evaluating the value of $\normtwo{\bA-\bB}$, we can determine whether to terminate the iteration when this value falls below $\delta$.

\paragraph{Low-rank approximation via SVD.}
In Section~\ref{section:svd-low-rank-approxi}, we will introduce  low-rank approximation using SVD. (Feel free to skip this section for a first reading.)
Suppose the rank-$r$ matrix $\bA$ admits \textit{reduced SVD} $\underset{m \times n}{\bA}=\underset{m \times r}{\bU} \,\,\, \underset{r \times r}{\bSigma}\,\,\,\underset{r \times n}{\bV^\top}$.
From the discussion in Section~\ref{section:svd-low-rank-approxi}, we can set $\bQ_k = \bU_k = [\bu_1, \bu_2, \ldots, \bu_k]$, i.e., containing the left singular vectors of $\bA$, such that 
$$
\normtwo{\bA - \bU_k\bU_k^\top\bA}=
\min\{\normtwo{\bA-\bB}:\rank(\bB)=k\}=\sigma_{k+1}.
$$ 
This result shows that the lower bound for the norm is the  $(k+1)$-th largest singular value $\sigma_{k+1}$ of $\bA$.
And the low-rank approximation then relies on the calculation of the SVD, which is a computationally expensive method to calculate. However, the low-rank approximation using MGS provides a faster algorithm to find a close approximation to the optimal basis vectors.

\index{Column pivoting}
\index{Revealing rank-one deficiency}
\index{Rank-revealing}
\section{QR with Column Pivoting: Revealing Rank-One Deficiency}\label{section:rank-one-qr-revealing}
We notice that Algorithm~\ref{alg:reduced-qr-rank-revealing} is just one of such methods used to find the column permutation when $\bA$ is rank-deficient, where we interchange the first $r$ linearly independent  columns of the matrix $\bA$ into the first $r$ columns of  $\bA\bP$. 
In cases where $\bA$ is nearly rank-one deficient, we aim to find a column permutation of $\bA$ such that the resulting pivotal element $r_{nn}$ of the QR decomposition is small. This problem is commonly referred to as the \textit{revealing rank-one deficiency} problem and falls under the broader category of \textit{rank-revealing QR decomposition (RRQR)} problems, which in general is a QR decomposition that exposes the rank deficiency of a matrix.

The RRQR problem is particularly useful in the sense that it allows us to infer the numerical rank of a matrix without explicitly computing its singular value decomposition (SVD), which can be a significant advantage in terms of computational time and resources.
In least squares problems, where one seeks the best approximate solution to an overdetermined system of equations, the rank of the coefficient matrix plays a crucial role. An RRQR factorization can help identify the effective rank and thus the number of linearly independent equations, which is essential for solving such problems accurately.
On the other hand, in statistical modeling and machine learning, selecting a subset of regressors that best explains the variability in the response variable is a common task. RRQR factorization can assist in identifying the most relevant subset of variables by revealing the rank structure of the matrix formed by these variables, e.g., finding independent and significant alphas signals for quantitative strategies \citep{lu2022feature}.

\begin{theoremHigh}[Revealing Rank-One Deficiency \citep{chan1987rank}]\label{theorem:finding-good-qr-ordering}
Let $\bA\in \real^{m\times n}$, and let $\bv\in \real^n$ be a unit $\ell_2$ norm vector (i.e., $\normtwo{\bv}=1$). Then there exists a permutation $\bP$ such that the reduced QR decomposition
$$
\bA\bP = \bQ\bR
$$ 
satisfies $r_{nn} \leq \sqrt{n} \epsilon$, where $\epsilon = \normtwo{\bA\bv}$, and $r_{nn}$ represents the $n$-th diagonal of $\bR$. 
Note the sizes of $\bQ\in \real^{m\times n}$ and $\bR\in \real^{n\times n}$ in the reduced QR decomposition.
\end{theoremHigh}
\begin{proof}[of Theorem~\ref{theorem:finding-good-qr-ordering}]
Suppose $\bP\in \real^{n\times n}$ is a permutation matrix and let $\bw=\bP^\top\bv$, where
$$
|w_n| = \max |v_i|,  \,\,\,\, \forall i \in \{1,2,\ldots,n\}.
$$
That is, we interchange the entry with largest magnitude into the last entry
such that the last component of $\bw$ is equal to the maximal component of $\bv$ in absolute value. 
Consequently, we have $|w_n| \geq 1/\sqrt{n}$. Suppose the QR decomposition of $\bA\bP$ is $\bA\bP = \bQ\bR$, then 
$$
\epsilon = \normtwo{\bA\bv} = \normtwo{(\bQ^\top\bA\bP) (\bP^\top\bv)} = \normtwo{\bR\bw} = 
\normtwo{\begin{bmatrix}
		\vdots \\
		r_{nn} w_n
\end{bmatrix}}
\geq |r_{nn} w_n| \geq |r_{nn}|/\sqrt{n},
$$
where the second equality  is a result of length preservation under orthogonal transformations, and $\bP$ is orthogonal satisfying $\bP\bP^\top=\bI$. This completes the proof.
\end{proof}
The following discussion relies on the existence of the SVD, which will be introduced in Chapter~\ref{chapter:SVD}. Feel free to skip at a first reading. Suppose the SVD of $\bA$ is given by $\bA = \sum_{i=1}^{n} \sigma_i \bu_i\bv_i^\top$, where $\sigma_i$'s are the singular values with $\sigma_1 \geq \sigma_2 \geq \ldots \geq \sigma_n$, i.e., $\sigma_n$ is the smallest singular value (suppose $m\geq n$), and $\bu_i$'s and $\bv_i$'s are the left and right singular vectors, respectively.
Then, if we let $\bv = \bv_n$ such that $\bA\bv_n = \sigma_n \bu_n$, \footnote{We will prove that the right singular vector of $\bA$ is equal to the right singular vector of $\bR$ if the matrix $\bA$ has the QR decomposition $\bA=\bQ\bR$, as stated in Lemma~\ref{lemma:svd-for-qr}. The claim can also applies to the singular values. Therefore, $\bv_n$ here is also the right singular vector of $\bR$.} we obtain 
$$
\normtwo{\bA\bv} = \sigma_n. 
$$ 
By constructing a permutation matrix $\bP$ satisfying
$$
|\bP^\top \bv|_n = \max |v_i|,  \,\,\,\, \forall i \in \{1,2,\ldots,n\},
$$
we can find a QR decomposition of $\bA\bP=\bQ\bR$ with a pivot $r_{nn}$ smaller than $\sqrt{n}\sigma_n$. 
In the scenario where $\bA$ is rank-one deficient, $\sigma_n$ will approach zero, resulting in a pivot $r_{nn}$ that is bounded to a small magnitude close to zero.

\section{QR with Column Pivoting: Revealing Rank-r Deficiency*}\label{section:rank-r-qr}\index{Rank-revealing QR}
Building upon the previous section, let us consider the computation of the reduced QR decomposition when $\bA\in \real^{m\times n}$ is approximately rank-$r$ deficient \footnote{Note that rank-$r$ here does not mean the matrix has a rank of $r$, but rather its  rank is $(\min\{m,n\}-r)$.} with $r>1$. Our goal now is to identify a permutation $\bP$ satisfying
\begin{equation}\label{equation:rankr-reval-qr}
\bA\bP = 
\bQ\bR=
\bQ
\begin{bmatrix}
\bL & \bM \\
\bzero & \bN
\end{bmatrix},
\end{equation}\footnote{To abuse the notation, we use the notation $\bL,\bM,\bN$ for clarity in the derivation. However, to maintain consistency with other contexts, it is preferable to replace $\bL,\bM,\bN$ with $\bR_{11}, \bR_{12}, \bR_{22}$.}
where $\bN \in \real^{r\times r}$, and $\norm{\bN}$ is small in some norm (and $\bL\in \real^{(n-r)\times (n-r)}, \bM\in \real^{(n-r)\times r}$ that can be inferred from context).

\begin{algorithm}[h] 
\caption{Reveal Rank-$r$ Deficiency} 
\label{alg:qr-reveal-rank-r} 
\begin{algorithmic}[1] 
\Require 
Matrix $\bA$ with size $m\times n$ and $m\geq n$, and effective $\rank(\bA)=n-r$; 
\State Initialize $\bW \in \real^{n\times r}$ to zero; \Comment{store the singular vectors}
\State Initial QR decomposition by $\bA = \bQ\bR$;
\For{$i=n$ to $n-r+1$}
\State $\bL \leftarrow $ leading $i\times i$ block of $\bR$;
\State Compute the singular vector $\bv\in \real^i$ corresponding to the min singular value of $\bL$;
\State Compute a permutation $\widetilde{\Pi} \in \real^{i\times i}$ such that $|\widetilde{\Pi}^\top \bv|_i = \max |\bv_j|, \forall j \in \{1,2,\ldots, i\}$;
\State Assign $\footnotesize\begin{bmatrix}
\bv \\
\bzero 
\end{bmatrix}$ to the $i$-th column of $\bW$;
\State Compute $\bW \leftarrow \Pi^\top \bW$, where $\Pi = \footnotesize\begin{bmatrix}
\widetilde{\Pi} & \bzero \\
\bzero & \bI 
\end{bmatrix}$;
\State Compute the QR decomposition: $\bL\widetilde{\Pi} = \widetilde{\bQ} \widetilde{\bL}$;
\State $\bP \leftarrow \bP \Pi$;
\State $\bQ \leftarrow \bQ \footnotesize\begin{bmatrix}
\widetilde{\bQ} & \bzero  \\
\bzero & \bI 
\end{bmatrix}$;
\State $\bR \leftarrow  \footnotesize\begin{bmatrix}
\widetilde{\bL} & \widetilde{\bQ}^\top \bB \\
\bzero & \bC 
\end{bmatrix}$, where $\footnotesize\bB =\bR[1:n-i, n-i+1:n]$, $\bC = \bR[n-i+1:n, n-i+1:n]$;
\EndFor
\end{algorithmic} 
\end{algorithm}

A recursive algorithm can be applied to achieve this. Suppose we have already isolated a small $k\times k$ block $\bN_k$, based on which, if we can isolate a small $(k+1)\times (k+1)$ block $\bN_{k+1}$, then we can find the permutation matrix recursively. To repeat, suppose we have the permutation $\bP_k$ such that the bottom-right matrix $\bN_k \in \real^{k\times k}$ has a small norm:
$$
\bA\bP_k = \bQ_k \bR_k=
\bQ_k 
\begin{bmatrix}
\bL_k & \bM_k \\
\bzero & \bN_k
\end{bmatrix}.
$$
We want to determine a permutation $\bP_{k+1}$ such that $\bN_{k+1} \in \real^{(k+1)\times (k+1)}$ also possesses a small norm:
$$
\textbf{(Goal): }\gap
\bA\bP_{k+1} = \bQ_{k+1} \bR_{k+1}=
\bQ_{k+1}
\begin{bmatrix}
\bL_{k+1} & \bM_{k+1} \\
\bzero & \bN_{k+1}
\end{bmatrix}.
$$
To achieve this, based on the algorithm presented in the previous section, there is an $(n-k)\times (n-k)$ permutation matrix $\widetilde{\bP}_{k+1}$ such that $\bL_k \in \real^{(n-k)\times (n-k)}$ admits the QR decomposition $\bL_k \widetilde{\bP}_{k+1} = \widetilde{\bQ}_{k+1}\widetilde{\bL}_k$, where the entry $(n-k, n-k)$ of $\widetilde{\bL}_k$ is small. By constructing 
$$
\bP_{k+1} \triangleq \bP_k
\begin{bmatrix}
\widetilde{\bP}_{k+1} & \bzero \\
\bzero & \bI 
\end{bmatrix}
\qquad\text{and}\qquad 
\bQ_{k+1} \triangleq \bQ_k 
\begin{bmatrix}
\widetilde{\bQ}_{k+1} & \bzero \\
\bzero & \bI 
\end{bmatrix},
$$
we have 
$$
\textbf{(Achieved): }\gap
\bA \bP_{k+1} = \bQ_{k+1}
\begin{bmatrix}
\widetilde{\bL}_k & \widetilde{\bQ}_{k+1}^\top  \bM_k \\
\bzero & \bN_k  
\end{bmatrix}.
$$
We know that entry $(n-k, n-k)$ of $\widetilde{\bL}_k$ is small. 
If we can prove the last row of $\widetilde{\bQ}_{k+1}^\top  \bM_k$ is small in norm, then we find the QR decomposition revealing rank-$(k+1)$ deficiency (refer to \citet{chan1987rank} for a proof). And the procedure is formulated in Algorithm~\ref{alg:qr-reveal-rank-r}.

\section{Existence of  QR Decomposition via  Householder Reflector}\label{section:qr-via-householder}
\textit{Householder matrices}, also known as \textit{Householder reflectors}, which can reflect vectors, play an crucial role in numerical linear algebra for tasks such as solving linear systems, addressing least squares problems, and deriving Hessenberg forms. In this section, we present how Householder reflectors can be used to prove the existence of the QR decomposition.
\subsubsection*{\textbf{Householder Reflectors}}
Let's begin with the formal definition of a Householder reflector, exploring its properties thereafter.
\begin{definition}[Householder Reflector\index{Householder reflector}]\label{definition:householder-reflector}
Let $\bu \in \real^n$ be a unit vector ($\normtwo{\bu}=1$). 
The matrix $\bH = \bI - 2\bu\bu^\top$ is referred to as a \textit{Householder reflector}, a.k.a., a \textit{Householder transformation}. We call this $\bH$ the Householder reflector associated with the unit vector $\bu$, where the unit vector $\bu$ is also known as the \textit{Householder vector}. 
When a vector $\bx$ is multiplied by $\bH$, it undergoes reflection with respect to the hyperplane $\spn\{\bu\}^\perp$.

Note that if $\normtwo{\bu} \neq 1$, we can define the Householder reflector $\bH$ as $\bH = \bI - 2  \frac{\bu\bu^\top}{\bu^\top\bu}$.
\end{definition}
From the definition of the Householder reflector, we can derive the following corollary, which states that certain vectors remain unchanged under the action of the Householder reflector.
\begin{corollary}[Unreflected by Householder]
Let $\bu\in\real^n$ be given with $\normtwo{\bu}=1$, and  define the Householder reflector as $\bH=\bI-2\bu\bu^\top$. 
The Householder reflector leaves  any vector $\bv$  perpendicular to $\bu$ unchanged; that is, $\bH\bv=\bv$ if $\bu^\top\bv=0$. 
\end{corollary}
The proof is straightforward since $(\bI - 2\bu\bu^\top)\bv = \bv - 2\bu\bu^\top\bv=\bv$.

\begin{SCfigure}
\centering
\includegraphics[width=0.5\textwidth]{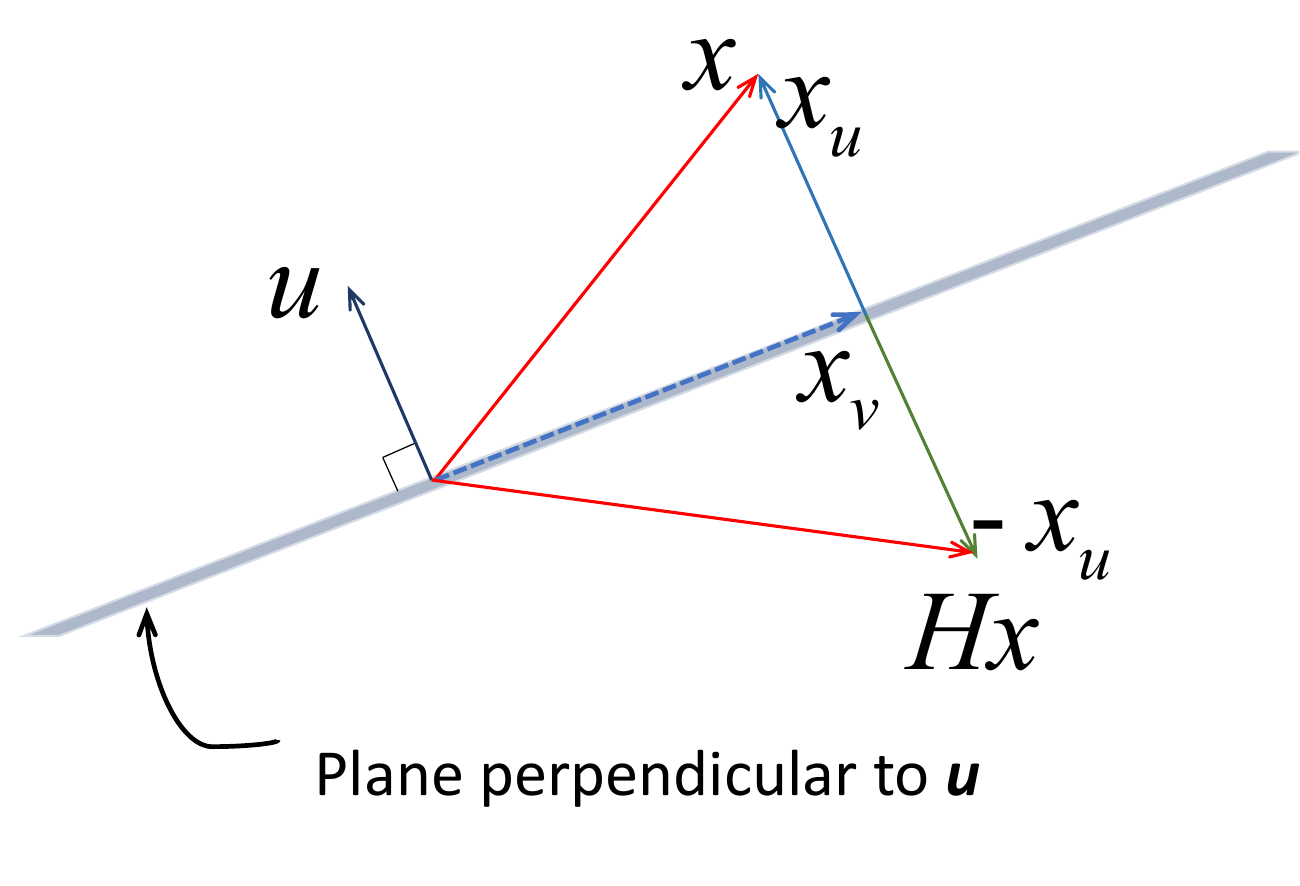}
\caption{Demonstration of the Householder reflector. The Householder reflector obtained by $\bH=\bI-2\bu\bu^\top$, where $\normtwo{\bu}=1$, will reflect a vector $\bx$ along the plane perpendicular to $\bu$. Specifically, it transforms $\bx=\bx_{\bv} + \bx_{\bu}$ into $\bx_{\bv} - \bx_{\bu}$.}
\label{fig:householder}
\end{SCfigure}

Suppose $\bu$ is a unit vector with $\normtwo{\bu}=1$, and let $\bv$ be a vector perpendicular to $\bu$. Then any vector $\bx$ in the plane can be decomposed into two components:
$$
\bx = \bx_{\bu} + \bx_{\bv},
$$ 
where the first component $\bx_{\bu}$ is parallel to $\bu$ and the second one $\bx_{\bv}$ is perpendicular to $\bu$ (i.e., parallel to $\bv$). 
Referring to Section~\ref{section:project-onto-a-vector} on vector projection, $\bx_{\bu}$ can be computed as $\bx_{\bu} = \frac{\bu\bu^\top}{\bu^\top\bu} \bx = \bu\bu^\top\bx$, representing  the projection of $\bx$ onto  $\bu$. 
Applying the Householder reflector associated with $\bu$ to the vector $\bx$, we obtain:
$$\bH\bx = (\bI - 2\bu\bu^\top)(\bx_{\bv} + \bx_{\bu}) = \bx_{\bv} -\bu\bu^\top \bx = \bx_{\bv} - \bx_{\bu},
$$ 
which means the Householder reflector transforms $\bx_{\bv} + \bx_{\bu}$ into $\bx_{\bv} - \bx_{\bu}$.
In other words, the space perpendicular to $\bu$ acts as a mirror, and any vector $\bx$ is reflected by the Householder reflector associated with $\bu$ (i.e., reflected by the hyperplane $\spn\{\bu\}^\perp$). 
The situation is illustrated in Figure~\ref{fig:householder}.

The previous discussion explains how to determine the reflected vector using the Householder reflector. 
However, an additional question arises: If we know in advance that two vectors are reflections of each other, how can we find the corresponding Householder reflector? 
The property is crucial for computing the QR decomposition, where we aim to reflect/transform a column into a specific form.
\begin{corollary}[Reflection Theorem: Finding the Householder Reflector]\label{corollary:householder-reflect-finding}
Suppose $\bx$ is reflected to $\by$ by a Householder reflector, with $\normtwo{\bx} = \normtwo{\by}$. The (unique) Householder reflector can be obtained by
$$
\bH = \bI - 2 \bu\bu^\top, \text{ where } \bu = \frac{\bx-\by}{\normtwo{\bx-\by}}.
$$
\end{corollary}
\begin{proof}[of Corollary~\ref{corollary:householder-reflect-finding}]
Write out the equation, we have 
$$
\begin{aligned}
\bH\bx &= \bx - 2 \bu\bu^\top\bx =\bx - 2\frac{(\bx-\by)(\bx^\top-\by^\top)}{(\bx-\by)^\top(\bx-\by)} \bx
= \bx - (\bx-\by) = \by.
\end{aligned}
$$
Note that the condition $\normtwo{\bx} = \normtwo{\by}$ is required to prove the result.
\end{proof}

Householder reflectors are useful for setting a block of components of a given vector to zero. Particularly, it is often desirable to set the vector $\ba\in \real^n$ to  zero, except for the $i$-th element. 
In such cases, the Householder vector can be chosen as:
$$
\bu = \frac{\ba - r\be_i}{\normtwo{\ba - r\be_i}}, \qquad \text{where } r = \pm\normtwo{\ba},
$$
which is a reasonable Householder vector since $\normtwo{\ba} = \normtwo{r\be_i} = \abs{r}$. We carefully notice that when $r=\normtwo{\ba}$, $\ba$ is reflected to $\normtwo{\ba}\be_i$ via the Householder reflector $\bH = \bI - 2 \bu\bu^\top$; conversely, when $r=-\normtwo{\ba}$, $\ba$ is reflected to $-\normtwo{\ba}\be_i$.

Recalling from Section~\ref{section:qr-gram-compute}, we claim the Householder or Givens method utilizes a set of orthogonal matrices to triangularize the matrix, thereby obtaining the QR decomposition and achieving a higher level of orthogonality in this context. The Householder reflector serves as one such orthogonal matrix for this purpose. 
In the following remark, we present additional properties of the Householder reflector.
\begin{remark}[Householder Properties]\label{remark:householder-propes}
A Householder reflector $\bH$ possesses the following properties: \footnote{Compare the properties of the Householder reflector with the properties of a projection matrix in Section~\ref{section:qr-gram-compute}}:
\begin{itemize}
\item $\bH\bH = \bI$, i.e., reflecting a vector twice is equivalent to not reflecting it at all.

\item Symmetry: $\bH = \bH^\top$.

\item Orthogonality: $\bH^\top\bH = \bH\bH^\top = \bI$ such that the Householder reflector is an orthogonal matrix.

\item $\bH\bu = -\bu$, if $\bH = \bI - 2 \bu\bu^\top$.
\item Unit eigenvalues: the eigenvalue of $\bH$ is either $1$ or $-1$. Given an eigenpair $(\lambda, \bx)$ of $\bH$, it follows that $\normtwo{\bH\bx}=\normtwo{\lambda\bx}=\normtwo{\bx}$. Therefore, $\lambda=\pm 1$.
\item The determinant of a Householder reflector is $-1$ (see Problem~\ref{problem:eig_det_house}).
\end{itemize}
\end{remark}

\begin{exercise}[Reflection Matrix]\label{exercise:household_reflec1}
Show that the Householder reflector $\bH = \bI - 2  \frac{\bu\bu^\top}{\bu^\top\bu}$ is an elementary reflector (Definition~\ref{definition:reflection_mat_intro}) satisfying $\det(\bH)=-1$, $\bH\bu=-\bu$, and $\frac{1}{2}(\bH+\bI)$ is an  elementary projector (Definition~\ref{definition:projection_matrix_intro}).
\end{exercise}

\subsubsection*{\textbf{Orthogonal Triangularization}}
To reiterate, as discussed in the Gram-Schmidt section, QR decomposition involves using a triangular matrix to orthogonalize a matrix $\bA$. 
Building upon this concept, if we have a sequence of orthogonal matrices that can successively transform $\bA$ into an upper triangular form, we can also construct the QR decomposition.
In particular, suppose we have an orthogonal matrix $\bQ_1$ that introduces zeros into all entries of the first column of $\bA$ except for the element at position (1,1); 
and similarly, an orthogonal matrix $\bQ_2$ that introduces zeros into the second column except for the elements at positions (2,1) and (2,2); 
$\ldots$. Then we can also find the QR decomposition.
To achieve this zero introduction, we could reflect the columns of the matrix to a basis vector $\be_1$ whose entries are all zero except the first entry.

Let $\bA=[\ba_1, \ba_2, \ldots, \ba_n]\in \real^{m\times n}$ be the column partition of $\bA$, and define further
\begin{equation}\label{equation:qr-householder-to-chooose-r-numeraically}
r_1 \triangleq \normtwo{\ba_1},\qquad \bu_1 \triangleq \frac{\ba_1 - r_1 \be_1}{\normtwo{\ba_1 - r_1 \be_1}}, \qquad  \text{and}\qquad \bH_1 \triangleq \bI - 2\bu_1\bu_1^\top,
\end{equation}
where $\be_1$ here is the first unit basis in $\real^m$, i.e., $\be_1 = [1;0;0;\ldots;0]\in \real^m$.
Then,
\begin{equation}\label{equation:householder-qr-projection-step1}
\bH_1\bA = [\bH_1\ba_1, \bH_1\ba_2, \ldots, \bH_1\ba_n] 
\triangleq
\begin{bmatrix}
r_1 & \bR_{1,2:n} \\
\bzero&  \bB_2
\end{bmatrix},
\end{equation}
which reflects $\ba_1$ to $r_1\be_1$ and introduces zeros below the diagonal in the first column. We observe that the entries below $r_1$ become zero after this specific reflection. Notice that we reflect $\ba_1$ to $\normtwo{\ba_1}\be_1$, where both vectors have the same length, rather than reflect $\ba_1$ to $\be_1$ directly to  ensure \textbf{numerical stability} and to match the requirement stated in Corollary~\ref{corollary:householder-reflect-finding}. 

\textbf{Choice of $r_1$:}  It should be noted that the choice of $r_1$ is \textbf{not unique}. To ensure \textbf{numerical stability}, it is desirable to choose $r_1 =-\text{sign}(a_{11}) \normtwo{\ba_1}$, where $a_{11}$ represents the first component of $\ba_{1}$. 
Alternatively, one can also choose $r_1 =\text{sign}(a_{11}) \normtwo{\ba_1}$, as long as $\normtwo{\ba_1}$ is equal to $\normtwo{r_1\be_1}$. However, a detailed discussion of this topic is beyond the scope of this section.

Next, we can apply this process to $\bB_2$ in Equation~\eqref{equation:householder-qr-projection-step1} to transform the entries  below the element (2,2) into zeros. 
Note that we do not apply this process to the entire matrix $\bH_1\bA$ but only to the submatrix $\bB_2$ in it because we have already introduced zeros in the first column, and reflecting again will reintroduce nonzero values back and destroy what have accomplished.

Suppose $\bB_2 = [\bb_2, \bb_3, \ldots, \bb_n]$ is the column partition of $\bB_2$, and define 
$$
r_2 \triangleq \normtwo{\bb_2},\qquad \bu_2 \triangleq \frac{\bb_2 - r_2 \be_1}{\normtwo{\bb_2 - r_2 \be_1}},   \qquad \widetilde{\bH}_2 \triangleq \bI - 2\bu_2\bu_2^\top, 
\qquad \text{and}\qquad  
\bH_2 \triangleq
\begin{bmatrix}
1 & \bzero \\
\bzero & \widetilde{\bH}_2
\end{bmatrix}.
$$
In this context, $\be_1$ denotes the first unit basis in $\real^{m-1}$, and $\bH_2$ is also an orthogonal matrix since $\widetilde{\bH}_2$ is orthogonal.
Applying $\widetilde{\bH_2}$ or $\bH_2$ yields
$$
\widetilde{\bH_2}\bB_2 = [\bH_2\bb_2, \bH_2\bb_3, \ldots, \bH_2\bb_n]
\triangleq
\begin{bmatrix}
r_2 & \bR_{2,3:n} \\
\bzero &\bC_3
\end{bmatrix},
$$
and
$$
\bH_2\bH_1\bA = [\bH_2\bH_1\ba_1, \bH_2\bH_1\ba_2, \ldots, \bH_2\bH_1\ba_n]
\triangleq
\begin{bmatrix}
r_1 & r_{12} & \bR_{1,3:n} \\
0 & r_2 & \bR_{2,3:n} \\
\bzero &  \bzero &\bC_3
\end{bmatrix}.
$$

The same process can go on. And if $\bA\in \real^{m\times n}$, after $n$ stages, we will finally triangularize $\bA = (\bH_n \bH_{n-1}\ldots\bH_1)^{-1} \bR = \bQ\bR$. Since the $\bH_i$'s are symmetric and orthogonal (Remark~\ref{remark:householder-propes}), we have orthogonal $\bQ=(\bH_n \bH_{n-1}\ldots\bH_1)^{-1} = \bH_1 \bH_2\ldots\bH_n$.

An example of a $5\times 4$ matrix is shown as follows, where $\boxtimes$ represents a value that is not necessarily zero, and \textbf{boldface} indicates the value has just been changed:
$$
\footnotesize
\begin{aligned}
\begin{sbmatrix}{\bA}
\boxtimes & \boxtimes & \boxtimes & \boxtimes \\
\boxtimes & \boxtimes & \boxtimes & \boxtimes \\
\boxtimes & \boxtimes & \boxtimes & \boxtimes \\
\boxtimes & \boxtimes & \boxtimes & \boxtimes \\
\boxtimes & \boxtimes & \boxtimes & \boxtimes
\end{sbmatrix}
&\stackrel{\bH_1}{\rightarrow}
\begin{sbmatrix}{\bH_1\bA}
\bm{\boxtimes} & \bm{\boxtimes} & \bm{\boxtimes} & \bm{\boxtimes} \\
\bm{0} & \bm{\boxtimes} & \bm{\boxtimes} & \bm{\boxtimes} \\
\bm{0} & \bm{\boxtimes} & \bm{\boxtimes} & \bm{\boxtimes} \\
\bm{0} & \bm{\boxtimes} & \bm{\boxtimes} & \bm{\boxtimes} \\
\bm{0} & \bm{\boxtimes} & \bm{\boxtimes} & \bm{\boxtimes}
\end{sbmatrix}
\stackrel{\bH_2}{\rightarrow}
\begin{sbmatrix}{\bH_2\bH_1\bA}
\boxtimes & \boxtimes & \boxtimes & \boxtimes \\
0 & \bm{\boxtimes} & \bm{\boxtimes} & \bm{\boxtimes} \\
0 & \bm{0} & \bm{\boxtimes} & \bm{\boxtimes} \\
0 & \bm{0} & \bm{\boxtimes} & \bm{\boxtimes} \\
0 & \bm{0} & \bm{\boxtimes} & \bm{\boxtimes}
\end{sbmatrix}
\stackrel{\bH_3}{\rightarrow}
\begin{sbmatrix}{\bH_3\bH_2\bH_1\bA}
\boxtimes & \boxtimes & \boxtimes & \boxtimes \\
0 & \boxtimes & \boxtimes & \boxtimes \\
0 & 0 & \bm{\boxtimes} & \bm{\boxtimes} \\
0 & 0 & \bm{0} & \bm{\boxtimes} \\
0 & 0 & \bm{0} & \bm{\boxtimes}
\end{sbmatrix}
\stackrel{\bH_4}{\rightarrow}
\begin{sbmatrix}{\bH_4\bH_3\bH_2\bH_1\bA}
\boxtimes & \boxtimes & \boxtimes & \boxtimes \\
0 & \boxtimes & \boxtimes & \boxtimes \\
0 & 0 & \boxtimes & \boxtimes \\
0 & 0 & 0 & \bm{\boxtimes} \\
0 & 0 & 0 & \bm{0}
\end{sbmatrix}.
\end{aligned}
$$

\paragraph{A closer look at the QR factorization.} The Householder algorithm is a procedure  that transforms a matrix into  triangular form through a sequence of orthogonal matrix operations. In the Gram-Schmidt process (both CGS and MGS), we use a triangular matrix to orthogonalize the matrix. 
In contrast, the Householder algorithm employs orthogonal matrices for matrix triangularization. 
The difference between these two approaches can be summarized as follows:
\begin{itemize}
\item Gram-Schmidt algorithm: triangular orthogonalization; 
\item Householder algorithm: orthogonal triangularization.
\end{itemize}

We further notice that, in the Householder algorithm or the Givens algorithm (which we will discuss shortly), a set of orthogonal matrices are applied so that the QR decomposition obtained is a \textit{full} QR decomposition. Whereas, the direct QR decomposition obtained by CGS or MGS is a \textit{reduced} decomposition (although the silent columns or rows can be further added to find the full version).

\section{Computing  Full QR Decomposition via  Householder Reflector}\label{section:householder_qr_cp}
Since $\bA$ has $n$ columns, in each step $i\in \{1, 2, \ldots, n\}$, where we introduce zeros below the diagonal in the $i$-th column, we perform operations on a submatrix of size  $(m-i+1)\times(n-i+1)$. 
To compute the upper triangular matrix $\bR = \bH_n\bH_{n-1}\ldots\bH_1 \bA$, we notice that 
$$
\begin{aligned}
\bR 
&= \left(\bH_n\ldots\left(\bH_3(\bH_2(\bH_1 \bA))\right)\right)\\
&= 
\begin{bmatrix}
\bI_{n-1}& \bzero \\
\bzero& \bI - 2\bu_n \bu_n^\top 
\end{bmatrix}
\ldots
\begin{bmatrix}
\bI_2& \bzero \\
\bzero& \bI - 2\bu_3\bu_3^\top 
\end{bmatrix}
\begin{bmatrix}
\bI_1& \bzero \\
\bzero& \bI - 2\bu_2\bu_2^\top 
\end{bmatrix}
\begin{bmatrix}
\bI - 2\bu_1\bu_1^\top
\end{bmatrix}\bA,
\end{aligned}
$$
where the parentheses indicate the order of the computation, the upper-left of $\bH_2$ is a $1\times 1$ identity matrix, and it will not modify the \textbf{first row} and \textbf{first column}
\footnote{Since the first column of $\bH_1\bA$ contains zeros below the main diagonal. Similarly for the remainders.} 
of $\bH_1\bA$ from the ``triangular property"; and the upper-left of $\bH_3$ is a $2\times 2$ identity matrix, which will not modify the \textbf{first 2 rows} and \textbf{first 2 columns} of $\bH_2\bH_2\bA$; $\ldots$. This property yields the step 8 in Algorithm~\ref{alg:qr-decomposition-householder}, which operates only on the $i:m$ rows and $i+1:n$ columns of $\bR$ in iteration $i$ (since the $i$-th column takes 1 flop explicitly in step 7, though this hardly reduces the complexity).
After the Householder transformation, we output the final triangular matrix $\bR$, and the process is shown in Algorithm~\ref{alg:qr-decomposition-householder}. 

Furthermore, in Algorithm~\ref{alg:qr-decomposition-householder}, in order to obtain the final orthogonal matrix $\bQ=\bH_1 \bH_2\ldots\bH_n$, we notice that
$$
\begin{aligned}
\bQ
&=(((\bH_1 \bH_2)\bH_3)\ldots\bH_n) \\
&= \begin{bmatrix}
	\bI - 2\bu_1\bu_1^\top
\end{bmatrix}
\begin{bmatrix}
	\bI_1& \bzero \\
	\bzero& \bI - 2\bu_2\bu_2^\top 
\end{bmatrix}
\begin{bmatrix}
	\bI_2& \bzero \\
	\bzero& \bI - 2\bu_3\bu_3^\top 
\end{bmatrix}
\ldots
\begin{bmatrix}
	\bI_{n-1}& \bzero \\
	\bzero& \bI - 2\bu_n\bu_n^\top 
\end{bmatrix},
\end{aligned}
$$
where, once again, the parentheses indicate the order of computation. The upper-left of $\bH_2$ is a $1\times 1$ identity matrix, and it will not modify the \textbf{first column} of $\bH_1$; and the upper-left of $\bH_3$ is a $2\times 2$ identity matrix, which will not modify the \textbf{first 2 columns} of $\bH_1\bH_2$; $\ldots$. This property yields the step 14 in the algorithm.
\begin{algorithm}[h] 
\caption{Full QR Decomposition via the Householder Reflector} 
\label{alg:qr-decomposition-householder} 
\begin{algorithmic}[1] 
\Require Matrix $\bA$ with size $m\times n $ and $m\geq n$; 
\Statex \textbf{Stage A: Get triangular matrix $\bR$}
\State Initially set $\bR = \bA$; 
\For{$i=1$ to $n$} 
\State $\ba \leftarrow \bR_{i:m,i}$, i.e., first column of $\bR_{i:m,i:n} \in \real^{(m-i+1)\times(n-i+1)}$;
\State $r \leftarrow \normtwo{\ba}$; \Comment{$2(m-i+1)$ flops}
\State $\bu_i \leftarrow \ba-r\be_1 \in \real^{m-i+1}$; \Comment{ $1$ flop}
\State $\bu_i \leftarrow \bu_i / \normtwo{\bu_i}$ ; \Comment{$3(m-i+1)$ flops}
\State $r_{ii} \leftarrow r$, $\bR_{i+1:m,i}\leftarrow\bzero$; \Comment{0 flops, update first column of $\bR_{i:m,i:n}$}
\State $\bR_{i:m,i+1:n} \leftarrow \bR_{i:m,i+1:n} - 2\bu_i (\bu_i^\top \bR_{i:m,i+1:n})$; \Comment{update $i+1:n$ columns of $\bR_{i:m,i:n}$}
\EndFor
\State Output $\bR$ as the triangular matrix;
\Statex \textbf{Stage B: Get orthogonal matrix $\bQ$}
\State Get $\bQ\triangleq\bH_1 \bH_2\ldots\bH_n$;
\State Initially set $\bQ \leftarrow \bH_1$;
\For{$i=1$ to $n-1$} 
\State $\bQ_{1:m,i+1:m} \leftarrow \bQ_{1:m,i+1:m}(\bI - 2\bu_{i+1}\bu_{i+1}^\top)=\bQ_{1:m,i+1:m}-\bQ_{1:m,i+1:m}2\bu_{i+1}\bu_{i+1}^\top$;
\EndFor
\State Output $\bQ$ as the orthogonal matrix;
\end{algorithmic} 
\end{algorithm}

\begin{theorem}[Algorithm Complexity: QR via Householder]\label{theorem:qr-full-householder}
Algorithm~\ref{alg:qr-decomposition-householder} requires $\sim 2mn^2-\frac{2}{3}n^3$ flops to compute a full QR decomposition of an $m\times n$ matrix with linearly independent columns and $m\geq n$. Further, if $\bQ$ is needed explicitly~\footnote{In many problems, $\bQ$ may not be needed, which is called \textit{Q-less QR decomposition}. For example, when solving the linear system $\bA\bx=\bb$, we may construct $\bQ^\top\bb=\bR^{-\top}\bA^\top\bb$ and compute $\bQ^\top\bb$ rather than $\bQ$.}, an additional $\sim 4m^2n-2mn^2$ flops are required.
\end{theorem}

\begin{proof}[of Theorem~\ref{theorem:qr-full-householder}]
For the $i$-th iteration of the loop, the submatrix $\bA_{i:m,i:n}$ is of size $(m-i+1)\times(n-i+1)$. Thus, $\ba_1$ is in $\real^{m-i+1}$. 

In step 4, the computation of $r = \normtwo{\ba}$ involves $m-i+1$ multiplications, $m-i$ additions, and 1 square root operation, resulting in a total of \underline{$2(m-i+1)$} flops.

In step 5, the computation of $\bu_i = \ba-r\be_1$ involves 1 subtraction, which requires \underline{$1$} flop due to the special structure of $\be_1$;

In step 6, similar to step 4, it requires $2(m-i+1)$ flops ($m-i+1$ multiplications, $m-i$ additions, and 1 square root) to compute the norm $\normtwo{\bu_i}$ and $m-i+1$ additional divisions, resulting in a total of \underline{$3(m-i+1)$} flops.


In step 8, consider loop $i=1$. the computation of $\bu_1^\top \bR_{1:m,2:n}$ requires $(n-1)\times$  ($m$ multiplications and $m-1$ additions),  resulting in \underline{$(n-1)(2m-1)$} flops. 
The computation of $2\bu_1$ requires \underline{$m$} multiplications. 
Additionally, the product of $2\bu_1 (\bu_1^\top \bR_{1:m,2:n})$ requires \underline{$m(n-1)$} multiplications to obtain an $m\times (n-1)$ matrix. 
The final matrix subtraction needs \underline{$m(n-1)$} subtractions. 
Therefore, the total complexity for step 8 when $i=1$ is \underline{$4m(n-1)+m-(n-1)$} flops. This analysis can be applied to any loop $i$, and the complexity of step 8 for loop $i$ can be obtained by \underline{$4(m-i+1)(n-i)+m-n+i$} flops.

Thus, for loop $i$, the total complexity from step 3 to step 8 can be defined as $f(i)$ flops.
To compute $\bR$, the final complexity is 
$$
\mathrm{cost}=f(1)+f(2)+\ldots +f(n).
$$
Simple calculation shows that the sum of $n$ loops  amounts to \underline{$2mn^2-\frac{2}{3}n^3$} flops if we keep only the leading terms.

To obtain the final orthogonal matrix $\bQ$, since we have already computed $2\bu_{i+1}$ in step 8, this will not incur additional costs. The computation of $\bQ_{1:m,i+1:m}2\bu_{i+1}$ involves $m\times$  ($m-i$ multiplications and $m-i-1$ additions), which is \underline{$m(2(m-i)-1)$} flops. Right multiplying by $\bu_{i+1}^\top$ takes \underline{$m(m-i)$} multiplications. Further, the final matrix subtraction requires \underline{$m(m-i)$} subtractions. So, in loop $i$, the complexity of step 14 is $g(i)=4m(m-i)-m = 4m^2-4mi-m$ flops. To compute $\bQ$, the final complexity is 
$$
\mathrm{cost}=g(1)+g(2)+\ldots +g(n-1).
$$
Simple calculation  shows the sum of $n-1$ loops is $4m^2n-2mn^2-4m^2+mn+m$ flops, or \underline{$4m^2n-2mn^2$} flops if we keep only the leading terms.
\end{proof}

After computing the full QR decomposition using the Householder algorithm, it is straightforward to recover the reduced QR decomposition by  removing the silent columns in $\bQ$ and silent rows in $\bR$. However, there does not exist a direct way to compute the reduced QR decomposition without the full decompositional form.

In \citet{golub2013matrix}, a Householder method for rank-revealing QR decomposition is discussed. The complexity is $4mnr-2r^2(m+n)+4r^3/3$ flops for rank-$r$ matrix $\bA$. 
When $r=n$, the complexity matches the result stated in Theorem~\ref{theorem:qr-full-householder}.

\section{Existence of  QR Decomposition via  Givens Rotation}\label{section:qr-givens}
We have seen that Givens rotations can be utilized to find the rank-one update/downdate of the Cholesky decomposition, as discussed in  Section~\ref{section:cholesky-rank-one-update}. Now let's take a look at what  the Givens rotations accomplishes through some specific examples. Consider the following $2\times 2$ orthogonal matrices
\begin{equation}\label{equation:hou_gi_exa}
\bF = 
\begin{bmatrix}
	-c & s\\
	s & c
\end{bmatrix}, 
\qquad 
\bJ=
\begin{bmatrix}
	c & -s \\
	s & c
\end{bmatrix},
\qquad \text{and}\qquad 
\bG=
\begin{bmatrix}
	c & s \\
	-s & c
\end{bmatrix},
\end{equation}
where $s = \sin(\theta)$ and $c=\cos(\theta)$ for some $\theta$. The first matrix has $\det(\bF)=-1$ and is a special case of a Householder reflector in  two-dimensional space, such that $\bF=\bI-2\bu\bu^\top$, where $\bu=\begin{bmatrix}
\sqrt{\frac{1+c}{2}}, &\sqrt{\frac{1-c}{2}}
\end{bmatrix}^\top$ or $\bu=\begin{bmatrix}
-\sqrt{\frac{1+c}{2}}, &-\sqrt{\frac{1-c}{2}}
\end{bmatrix}^\top$. The latter two matrices have $\det(\bJ)=\det(\bG)=1$ and induce rotations instead of reflections. 
These matrices are referred to as {\textit{Givens rotations}}.

\begin{figure}[H]
\centering  
\vspace{-0.35cm} 
\subfigtopskip=2pt 
\subfigbottomskip=2pt 
\subfigcapskip=-5pt 
\subfigure[$\by = \bJ\bx$, counter-clockwise rotation.]{\label{fig:rotation1}
	\includegraphics[width=0.4\linewidth]{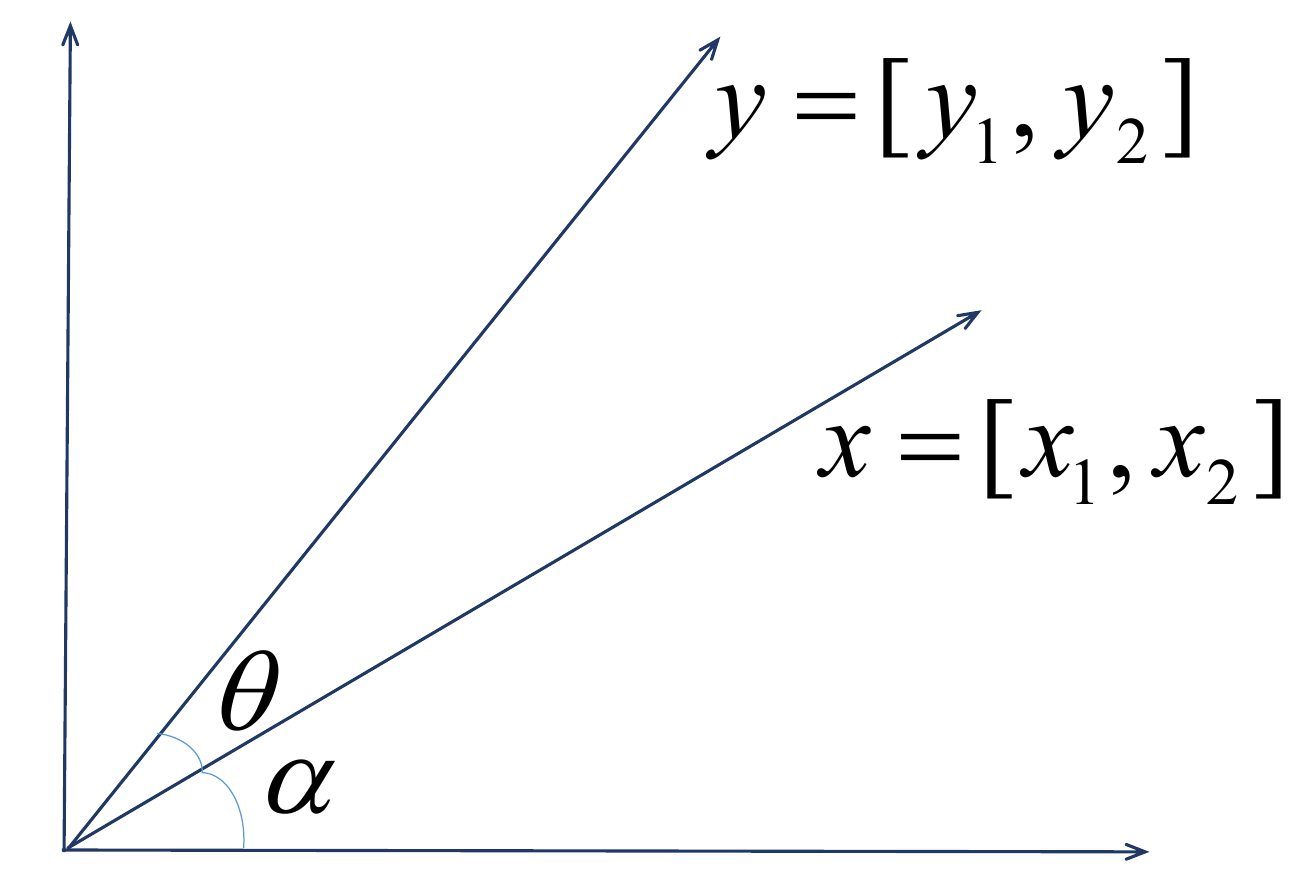}}
\quad \quad 
\subfigure[$\by = \bG\bx$, clockwise rotation.]{\label{fig:rotation2}
	\includegraphics[width=0.4\linewidth]{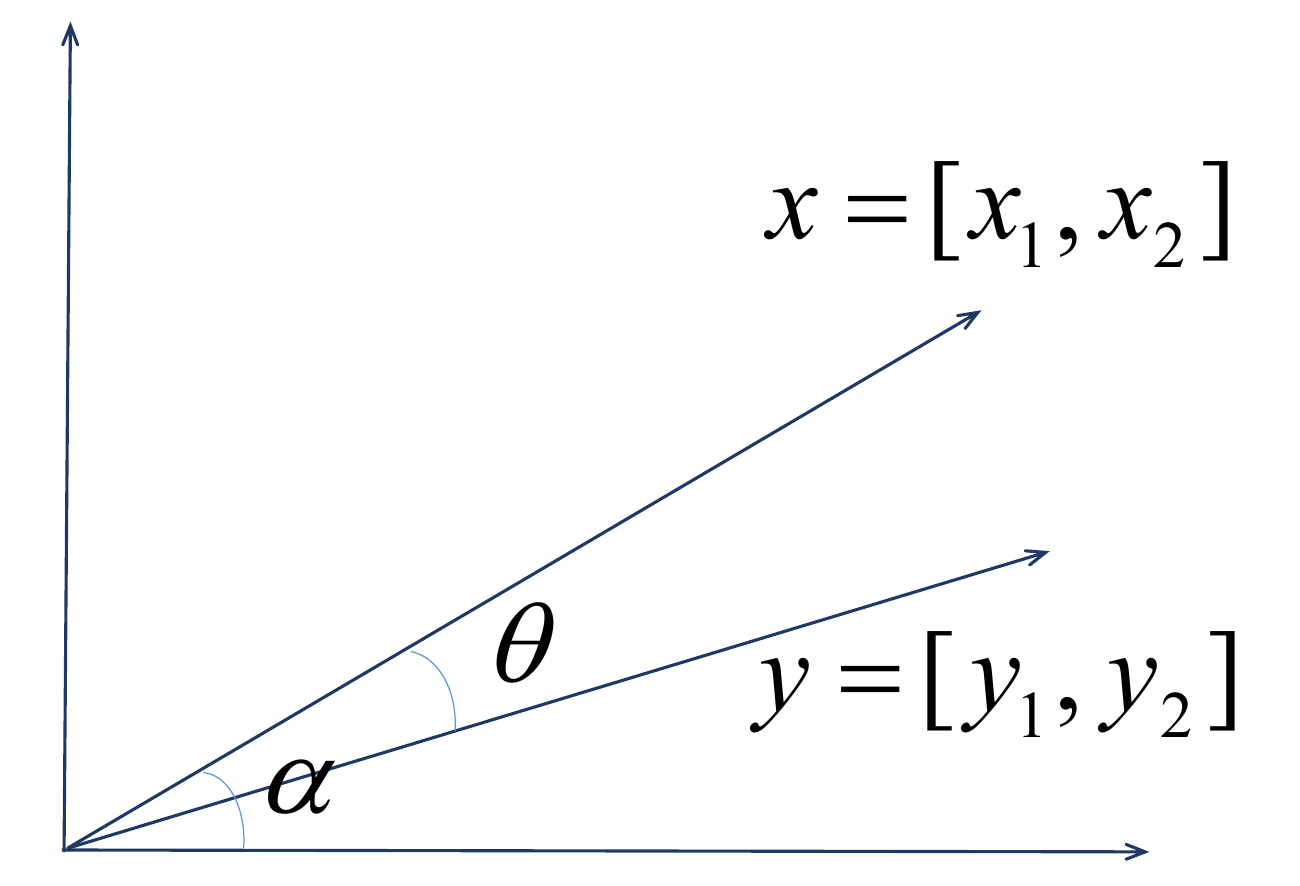}}
\caption{Illustration  of two specific Givens rotations.}
\label{fig:rotation}
\end{figure}

Figure~\ref{fig:rotation1} demonstrate the rotation of a vector $\bx$ under the action of matrix $\bJ$, resulting in $\by = \bJ\bx$. 
Specifically,
$
\left\{
\begin{aligned}
&y_1 = c\cdot x_1 - s\cdot x_2;   \\
&y_2 = s \cdot x_1 + c\cdot x_2.
\end{aligned}
\right.
$
We aim to confirm that the angle between $\bx$ and $\by $ is actually $\theta$ (representing a counter-clockwise rotation) under the Givens rotation $\bJ$. Firstly, we have 
$$ 
\left\{
\begin{aligned}
&\cos(\alpha) =\frac{x_1}{\sqrt{x_1^2+x_2^2}};   \\
&\sin (\alpha) =\frac{x_2}{\sqrt{x_1^2+x_2^2}},
\end{aligned}
\right.
\qquad 
\text{and }\qquad 
\left\{
\begin{aligned}
&\cos(\theta) =c;   \\
&\sin (\theta) =s.
\end{aligned}
\right.
$$
This implies that $\cos(\theta+\alpha) = \cos(\theta)\cos(\alpha)-\sin(\theta)\sin(\alpha)$.
If we can demonstrate that $\cos(\theta+\alpha) = \cos(\theta)\cos(\alpha)-\sin(\theta)\sin(\alpha)$ is equal to $\frac{y_1}{\sqrt{y_1^2+y_2^2}}$, then we complete the proof.
For the former one, $\cos(\theta+\alpha) = \cos(\theta)\cos(\alpha)-\sin(\theta)\sin(\alpha)=\frac{c\cdot x_1 - s\cdot x_2}{\sqrt{x_1^2+x_2^2}}$. For the latter one, it can be verified that $\sqrt{y_1^2+y_2^2}=\sqrt{x_1^2+x_2^2}$, and thus $\frac{y_1}{\sqrt{y_1^2+y_2^2}} = \frac{c\cdot x_1 - s\cdot x_2}{\sqrt{x_1^2+x_2^2}}$. This completes the proof. Similarly, we can also show that the angle between vectors $\by=\bG\bx$ and $\bx$ is also $\theta$ in Figure~\ref{fig:rotation2}, and the rotation is clockwise.

\begin{exercise}[Decomposition of Rotation Matrices]
Consider the matrix $\bJ$ or $\bG$ in \eqref{equation:hou_gi_exa}. Show that when $s\neq 0$ (i.e., $\theta\neq k\pi$ for integer $k$), $\bJ$ or $\bG$ can be factored as $\bJ/\bG=\bR\bL\bR$, where 
$$
\bR = \begin{bmatrix}
1 & r \\
0 & 1
\end{bmatrix}
\quad\text{and}\quad
\bL= \begin{bmatrix}
	1 & 0 \\
	l & 1
\end{bmatrix}.
$$
Therefore, the plane rotation can be expressed as three transformations.
\end{exercise}

\index{Givens plane}
\index{Givens rotation}
More generally, we have defined the $n$-th order Givens rotation in Definition~\ref{definition:givens-rotation-in-qr}.
It can be easily verified that the $n$-th order Givens rotation is an orthogonal matrix with a determinant of 1. 
For any vector $\bx =[x_1, x_2, \ldots, x_n]^\top \in \real^n$, we can express the transformed vector $\by$ as $\by = \bG_{kl}\bx$ (the subscripts $k,l$ indicate the rotations occur \textbf{in plane $k$ and $l$}), where
$$ 
\by=\left\{
\begin{aligned}
&y_k = c \cdot x_k + s\cdot x_l,   \\
&y_l = -s\cdot x_k +c\cdot x_l,  \\
&y_j = x_j . &  (j\neq k,l) 
\end{aligned}
\right.
$$
That is, a Givens rotation applied to $\bx$ rotates two components of $\bx$ by an angle $\theta$, while keeping all other components unchanged.
When $\sqrt{x_k^2 + x_l^2} \neq 0$,
we can define $c \triangleq \frac{x_k}{\sqrt{x_k^2 + x_l^2}}$ and $s\triangleq\frac{x_l}{\sqrt{x_k^2 + x_l^2}}$. Then, 
$$ 
\by=\left\{
\begin{aligned}
&y_k = \sqrt{x_k^2 + x_l^2},   \\
&y_l = 0,  \\
&y_j = x_j . &  (j\neq k,l) 
\end{aligned}
\right.
$$
This  finding is crucial for performing the QR decomposition using Givens rotations.

\begin{corollary}[Basis From Givens Rotations Forwards]\label{corollary:basis-from-givens}
Let $\bx \in \real^n$ be any vector. Then, there exists a set of Givens rotations $\{\bG_{12}, \bG_{13}, \ldots, \bG_{1n}\}$ such that the transformation $\bG_{1n}\ldots \bG_{13}\bG_{12}\bx = \normtwo{\bx}\be_1$, where $\be_1\in \real^n$ represents the first unit basis in $\real^n$.
\end{corollary}
\begin{proof}[of Corollary~\ref{corollary:basis-from-givens}]
Based on the aforementioned finding, we can identify a set of Givens rotations $\bG_{12}, \bG_{13}$, and $\bG_{14}$ such that 
$$
\begin{aligned}
\bG_{12}\bx &= \left[\sqrt{x_1^2 + x_2^2}, 0, x_3, \ldots, x_n \right]^\top;\\
\bG_{13}\bG_{12}\bx &= \left[\sqrt{x_1^2 + x_2^2+x_3^2}, 0, 0, x_4, \ldots, x_n \right]^\top;\\
\bG_{14}\bG_{13}\bG_{12}\bx &= \left[\sqrt{x_1^2 + x_2^2+x_3^2+x_4^2},0, 0, 0, x_5, \ldots, x_n \right]^\top.
\end{aligned}
$$
Continuing this process, we obtain $\bG_{1n}\ldots \bG_{13}\bG_{12} = \normtwo{\bx}\be_1$.
\end{proof}

\begin{remark}[Basis From Givens Rotations Backwards]\label{remark:basis-from-givens2}
In Corollary~\ref{corollary:basis-from-givens}, we derive the Givens rotations that introduce zeros from the second entry to the $n$-th entry of a vector (i.e., in forward order). 
However, there are cases where we desire the reverse order, i.e., introducing zeros from the $n$-th entry to the second entry such that $\bG_{12}\bG_{13}\ldots  \bG_{1n}\bx = \normtwo{\bx}\be_1$, where $\be_1\in \real^n$ represents the first unit basis in $\real^n$.

The procedure is similar. We can find $\bG_{1n},\bG_{1,(n-1)}, \bG_{1,(n-2)}$ using a similar approach, such that 
$$
\bG_{1n}\bx = \left[\sqrt{x_1^2 + x_n^2}, x_2, x_3, \ldots, x_{n-1}, 0 \right]^\top,
$$
$$
\bG_{1,(n-1)}\bG_{1n}\bx = \left[\sqrt{x_1^2 + x_{n-1}^2+x_n^2}, x_2, x_3, \ldots, x_{n-2}, 0, 0  \right]^\top,
$$
and
$$
\bG_{1,(n-2)}\bG_{1,(n-1)}\bG_{1n}\bx = \left[\sqrt{x_1^2 + x_{n-2}^2+x_{n-1}^2+x_n^2}, x_2, x_3,  \ldots,x_{n-3},0, 0, 0  \right]^\top.
$$
Continuing  this process, we obtain $\bG_{12}\bG_{13}\ldots  \bG_{1n}\bx = \normtwo{\bx}\be_1$.

\paragraph{An alternative form.} Alternatively, there are rotations $\{\bG_{12}, \bG_{23}, \ldots,\bG_{(n-1),n}\}$ such that $\bG_{12} \bG_{23} \ldots \bG_{(n-1),n}\bx = \normtwo{\bx}\be_1$, where
$$
\bG_{(n-1),n}\bx = \left[x_1, x_2, \ldots, x_{n-2},\sqrt{x_{n-1}^2 + x_n^2}, 0 \right]^\top,
$$
$$
\bG_{(n-2),(n-1)}\bG_{(n-1),n}\bx = \left[x_1, x_2, \ldots,x_{n-3}, \sqrt{x_{n-2}^2+x_{n-1}^2 + x_n^2}, 0, 0  \right]^\top,
$$
and
$$\footnotesize
\begin{aligned}
\bG_{(n-3),(n-2)}\bG_{(n-2),(n-1)}\bG_{(n-1),n}\bx &= 
\bigg[x_1,  x_2, \ldots ,x_{n-4}, \sqrt{x_{n-3}^2+x_{n-2}^2+x_{n-1}^2 + x_n^2},0, 0, 0  \bigg]^\top.
\end{aligned}
$$
Continuing  this process, we obtain $\bG_{12} \bG_{23} \ldots \bG_{(n-1),n}\bx = \normtwo{\bx}\be_1$.

The  backward Givens rotation basis update discussed above  will be proved useful in the context of rank-one changes of the QR decomposition (Section~\ref{section:qr-rank-one-changes}).
\end{remark}

From Corollary~\ref{corollary:basis-from-givens} mentioned earlier, to introduce zeros, we could \textbf{rotate} the columns of the matrix to align with basis vector $\be_1$, where all entries except the first entry are zero.
Let $\bA=[\ba_1, \ba_2, \ldots, \ba_n] \in \real^{m\times n}$ be the column partition of matrix $\bA$, and let 
\begin{equation}\label{equation:qr-rotate-to-chooose-r-numeraically}
\bG_1 \triangleq \bG_{1m}\ldots \bG_{13}\bG_{12}.
\end{equation}
Then,
\begin{equation}\label{equation:rotate-qr-projection-step1}
\begin{aligned}
\bG_1\bA &= [\bG_1\ba_1, \bG_1\ba_2, \ldots, \bG_1\ba_n] 
\triangleq
\begin{bmatrix}
\normtwo{\ba_1} & \bR_{1,2:n} \\
\bzero&  \bB_2
\end{bmatrix},
\end{aligned}
\end{equation}
which rotates $\ba_1$ to $\normtwo{\ba_1}\be_1$ and introduces zeros below the diagonal in the $1$-st column of $\bA$, where $\be_1$ represents the first unit basis in $\real^m$, i.e., $\be_1 = [1;0;0;\ldots;0]\in \real^m$. It is important to note that the Givens rotation $\bG_1$ will sequentially affects the pair of elements ($1$-st, $2$-nd), ($1$-st, $3$-rd), ($1$-st, $4$-th), \ldots, ($1$-st, $m$-th) for any vector $\bv\in \real^m$.

We can then apply this process to $\bB_2$ in Equation~\eqref{equation:rotate-qr-projection-step1} to set all entries below the (2,2)-th entry to zero.
Let $\bB_2 = [\bb_2, \bb_3, \ldots, \bb_n]$ be the column partition of $\bB_2$, and consider the Givens rotation sequence 
$$
\bG_2 \triangleq \bG_{2m}\ldots\bG_{24}\bG_{23},
$$
where the rotations $\bG_{2m}, \ldots, \bG_{24}, \bG_{23}$ can be inferred from the context.
Then, 
$$
\begin{aligned}
\bG_2
\begin{bmatrix}
\bR_{1,2:n} \\
\bB_2
\end{bmatrix}
&\triangleq 
\begin{bmatrix}
r_{12} & \bR_{1,3:n} \\
\normtwo{\bb_2} & \bR_{2,3:n} \\
\bzero &\bC_3
\end{bmatrix}
\end{aligned}
\quad\text{and}\quad
\bG_2\bG_1\bA 
\triangleq
\begin{bmatrix}
	\normtwo{\ba_1} & r_{12} & \bR_{1,3:n} \\
	0 & \normtwo{\bb_2} & \bR_{2,3:n} \\
	\bzero &  \bzero &\bC_3
\end{bmatrix},
$$
where the zeros in the first columns are kept since rotations of zero values are still zero.

The same process can continue, ultimately resulting in the triangularization of $\bA = (\bG_n \bG_{n-1}\ldots\bG_1)^{-1} \bR = \bQ\bR$. 
Since the matrices $\bG_i$'s are orthogonal, the orthogonal $\bQ$ can also be obtained as $\bQ=(\bG_n \bG_{n-1}\ldots\bG_1)^{-1} = \bG_1^\top \bG_2^\top\ldots\bG_n^\top =(\bG_n \ldots \bG_2 \bG_1)^\top $, and 
\begin{equation}\label{equation:givens-q}
\begin{aligned}
\bQ=(\bG_n \ldots \bG_2 \bG_1)^\top 
=\left\{(\bG_{nm} \ldots \bG_{n,(n+1)}) \ldots (\bG_{2m}\ldots \bG_{23})  ( \bG_{1m} \ldots \bG_{12} )\right\}^\top .
\end{aligned}
\end{equation}

\paragraph{When will the Givens method work better?} In practice, the Givens rotation algorithm performs better than the Householder algorithm when $\bA$ already contains many zeros below the main diagonal. 
Therefore,  Givens rotations are particularly suitable for rank-one changes in the QR decomposition, as such changes introduce only a small number of nonzero values (Section~\ref{section:qr-rank-one-changes}).

An illustrative  example of a $5\times 4$ matrix is shown below, where $\boxtimes$ represents a value that is not necessarily zero, and \textbf{boldface} indicates the value has just been changed. 
\paragraph{Givens rotations in $\bG_1$.} For a $5\times 4$ example, we realize that $\bG_1 = \bG_{15}\bG_{14}\bG_{13}\bG_{12}$. And the process is shown as follows:
$$
\footnotesize
\begin{aligned}
\begin{sbmatrix}{\bA}
	\boxtimes & \boxtimes & \boxtimes & \boxtimes \\
	\boxtimes & \boxtimes & \boxtimes & \boxtimes \\
	\boxtimes & \boxtimes & \boxtimes & \boxtimes \\
	\boxtimes & \boxtimes & \boxtimes & \boxtimes \\
	\boxtimes & \boxtimes & \boxtimes & \boxtimes
\end{sbmatrix}
\stackrel{\bG_{12}}{\rightarrow}
\begin{sbmatrix}{\bG_{12}\bA}
\bm{\boxtimes} & \bm{\boxtimes} & \bm{\boxtimes} & \bm{\boxtimes} \\
\bm{0} & \bm{\boxtimes} & \bm{\boxtimes} & \bm{\boxtimes} \\
\boxtimes & \boxtimes & \boxtimes & \boxtimes \\
\boxtimes & \boxtimes & \boxtimes & \boxtimes \\
\boxtimes & \boxtimes & \boxtimes & \boxtimes
\end{sbmatrix}
\stackrel{\bG_{13}}{\rightarrow}
\begin{sbmatrix}{\bG_{13}\bG_{12}\bA}
	\bm{\boxtimes} & \bm{\boxtimes} & \bm{\boxtimes} & \bm{\boxtimes} \\
	0 & \boxtimes & \boxtimes &\boxtimes \\
	\bm{0} & \bm{\boxtimes} & \bm{\boxtimes} & \bm{\boxtimes} \\
	\boxtimes & \boxtimes & \boxtimes & \boxtimes \\
\boxtimes & \boxtimes & \boxtimes & \boxtimes
\end{sbmatrix}
\stackrel{\bG_{14}}{\rightarrow}
\begin{sbmatrix}{\bG_{14}\bG_{13}\bG_{12}\bA}
\bm{\boxtimes} & \bm{\boxtimes} & \bm{\boxtimes} & \bm{\boxtimes} \\
0 & \boxtimes & \boxtimes &\boxtimes \\
0 & \boxtimes & \boxtimes & \boxtimes \\
\bm{0} & \bm{\boxtimes} & \bm{\boxtimes} & \bm{\boxtimes} \\
\boxtimes & \boxtimes & \boxtimes & \boxtimes
\end{sbmatrix}
\stackrel{\bG_{15}}{\rightarrow}
\begin{sbmatrix}{\bG_{15}\bG_{14}\bG_{13}\bG_{12}\bA}
\bm{\boxtimes} & \bm{\boxtimes} & \bm{\boxtimes} & \bm{\boxtimes} \\
0 & \boxtimes & \boxtimes &\boxtimes \\
0 & \boxtimes & \boxtimes & \boxtimes \\
0 & \boxtimes & \boxtimes & \boxtimes \\
\bm{0} & \bm{\boxtimes} & \bm{\boxtimes} & \bm{\boxtimes} \\
\end{sbmatrix}.
\end{aligned}
$$

\paragraph{Givens rotation as a big picture.} Considering $\bG_1, \bG_2, \bG_3, \bG_4$ as a single matrix, we obtain
$$
\footnotesize
\begin{aligned}
\begin{sbmatrix}{\bA}
	\boxtimes & \boxtimes & \boxtimes & \boxtimes \\
	\boxtimes & \boxtimes & \boxtimes & \boxtimes \\
	\boxtimes & \boxtimes & \boxtimes & \boxtimes \\
	\boxtimes & \boxtimes & \boxtimes & \boxtimes \\
	\boxtimes & \boxtimes & \boxtimes & \boxtimes
\end{sbmatrix}
&\stackrel{\bG_1}{\rightarrow}
\begin{sbmatrix}{\bG_1\bA}
	\bm{\boxtimes} & \bm{\boxtimes} & \bm{\boxtimes} & \bm{\boxtimes} \\
	\bm{0} & \bm{\boxtimes} & \bm{\boxtimes} & \bm{\boxtimes} \\
	\bm{0} & \bm{\boxtimes} & \bm{\boxtimes} & \bm{\boxtimes} \\
	\bm{0} & \bm{\boxtimes} & \bm{\boxtimes} & \bm{\boxtimes} \\
	\bm{0} & \bm{\boxtimes} & \bm{\boxtimes} & \bm{\boxtimes}
\end{sbmatrix}
\stackrel{\bG_2}{\rightarrow}
\begin{sbmatrix}{\bG_2\bG_1\bA}
	\boxtimes & \boxtimes & \boxtimes & \boxtimes \\
	0 & \bm{\boxtimes} & \bm{\boxtimes} & \bm{\boxtimes} \\
	0 & \bm{0} & \bm{\boxtimes} & \bm{\boxtimes} \\
	0 & \bm{0} & \bm{\boxtimes} & \bm{\boxtimes} \\
	0 & \bm{0} & \bm{\boxtimes} & \bm{\boxtimes}
\end{sbmatrix}
\stackrel{\bG_3}{\rightarrow}
\begin{sbmatrix}{\bG_3\bG_2\bG_1\bA}
	\boxtimes & \boxtimes & \boxtimes & \boxtimes \\
	0 & \boxtimes & \boxtimes & \boxtimes \\
	0 & 0 & \bm{\boxtimes} & \bm{\boxtimes} \\
	0 & 0 & \bm{0} & \bm{\boxtimes} \\
	0 & 0 & \bm{0} & \bm{\boxtimes}
\end{sbmatrix}
\stackrel{\bG_4}{\rightarrow}
\begin{sbmatrix}{\bG_4\bG_3\bG_2\bG_1\bA}
	\boxtimes & \boxtimes & \boxtimes & \boxtimes \\
	0 & \boxtimes & \boxtimes & \boxtimes \\
	0 & 0 & \boxtimes & \boxtimes \\
	0 & 0 & 0 & \bm{\boxtimes} \\
	0 & 0 & 0 & \bm{0}
\end{sbmatrix}.
\end{aligned}
$$

\paragraph{Order of introducing zeros.} When using Givens rotations to compute the QR decomposition, there is flexibility in choosing the order in which the zeros are introduced in $\bR$. In our approach, we introduce zeros column by column. Alternatively, it is also possible to introduce zeros row by row.

\section{Computing  Full QR Decomposition via  Givens Rotation}
The algorithm for computing the full QR decomposition using Givens rotations is straightforward, as demonstrated in the example above. This process is illustrated in Algorithm~\ref{alg:qr-decomposition-givens}.
\begin{algorithm}[h] 
\caption{Full QR Decomposition via the Givens Rotation} 
\label{alg:qr-decomposition-givens} 
\begin{algorithmic}[1] 
\Require matrix $\bA$ with size $m\times n $ and $m\geq n$; 
\State Initially set $\bR \leftarrow \bA$, $\bQ\leftarrow\bI$;
\Statex \textbf{Stage A: Get triangular matrix $\bR$}
\For{$i=1$ to $n$} 
\For{$j=i+1$ to $m$}
\State Get Givens rotation $\bG_{i,j}$ with the following parameters $c, s$:
\State $c \leftarrow \frac{x_k}{\sqrt{x_k^2 + x_l^2}}$, $s\leftarrow\frac{x_l}{\sqrt{x_k^2 + x_l^2}}$, where $x_k = r_{ii}$, $x_l = r_{ji}$;
\State Calculate $\bR \leftarrow \bG_{i,j}\bR $ in following two steps:
\State $i$-th row: $\bR_{i,:}\leftarrow c\cdot \bR_{i,:} + s \bR_{j,:} $;
\State $j$-th row: $\bR_{j,:} \leftarrow -s\cdot \bR_{i,:} + c \bR_{j,:} $;
\EndFor
\EndFor
\State Output $\bR$ as the triangular matrix;
\Statex \textbf{Stage B: Get orthogonal matrix $\bQ$}
\For{$i=1$ to $n$} 
\For{$j=i+1$ to $m$}
\State Calculate $\bQ \leftarrow \bG_{i,j}\bQ $ in following two steps:
\State $i$-th row: $\bQ_{i,:} \leftarrow c\cdot \bQ_{i,:} + s \bQ_{j,:} $;
\State $j$-th row: $\bQ_{j,:} \leftarrow -s\cdot \bQ_{i,:} + c \bQ_{j,:} $;
\EndFor
\EndFor
\State Output $\bQ\leftarrow\bQ^\top$ from Equation~\eqref{equation:givens-q};
\end{algorithmic} 
\end{algorithm}

\begin{theorem}[Algorithm Complexity: QR via Givens]\label{theorem:qr-full-givens}
Algorithm~\ref{alg:qr-decomposition-givens}  requires $\sim 3mn^2-n^3$ flops to compute a full QR decomposition of an $m\times n$ matrix with linearly independent columns and $m\geq n$. Further, if $\bQ$ is needed explicitly, an additional $\sim 3mn^2-n^3$ flops required.
\end{theorem}

\begin{proof}[of Theorem~\ref{theorem:qr-full-givens}]
For step 5, each iteration $i,j$ requires $6$ flops (including 2 square operations, 1 addition, 1 square root, and 2 divisions). And there are $(m-i)$ iterations for each $i$, which means a total of $(m-1)+(m-2)+\ldots+(m-n)=(mn-\frac{n^2+n}{2})$ iterations. 
Therefore, the overall complexity for all the step 5's is \underline{$6(mn-\frac{n^2+n}{2})$} flops.

For each iteration $i$, steps 7 and 8 involve operating on two vectors of length $(n-i+1)$. 
The two steps take $6(n-i+1)$ flops for each iteration $i$ (consisting of $4(n-i+1)$ multiplications and $2(n-i+1)$ additions). 
And for each $i$, there are $(m-i)$ such operations resulting in a total of  $(m-i)\times 6(n-i+1)$ flops for each $i$. Let $f(i) = (m-i)\times 6(n-i+1) = m(n+1)-(m+n+1)i +i^2$. The overall complexity of the two steps is equal to 
$$
\mathrm{cost} =f(1)+ f(2) +\ldots +f(n),
$$
or \underline{$3mn^2-n^3$} flops if keeping only the leading terms.

Similarly, the complexity of steps 15 and 16 can be determined as \underline{$3mn^2-n^3$} flops when considering only the leading terms.
\end{proof}

Similar to the Householder algorithm, once the full QR decomposition has been computed using the Givens algorithm, it is straightforward to obtain the reduced QR decomposition by eliminating the silent columns in $\bQ$ and silent rows in $\bR$.

\index{Uniqueness}
\section{Uniqueness of  QR Decomposition}\label{section:nonunique-qr}
The QR decomposition results obtained from the Gram-Schmidt process, the Householder algorithm, and the Givens algorithm can differ. 
Even within the Householder algorithm, there are different methods to choose the sign of $r_1$ in Equation~\eqref{equation:qr-householder-to-chooose-r-numeraically}. 
Consequently, from this point, the QR decomposition of a matrix is not unique. The non-uniqueness of the QR decomposition can be demonstrated by the following example.

\begin{example}[Non-Uniqueness of the QR Decomposition\index{Uniqueness}]\label{example:qr_unique}
Let $\bA$ be
$
\bA = 
\scriptsize
\begin{bmatrix}
4 & 1 \\
3 & 2
\end{bmatrix}.
$
The QR decomposition of $\bA$ can be obtained by 
$$
\begin{aligned}
\bA 
&= \bQ_1\bR_1\equiv
\begin{bmatrix}
	0.8 & -0.6 \\
	0.6 & 0.8
\end{bmatrix}
\begin{bmatrix}
	5 & 2 \\
	0 & 1
\end{bmatrix}
&=& \bQ_2\bR_2\equiv
\begin{bmatrix}
0.8 & 0.6 \\
0.6 & -0.8
\end{bmatrix}
\begin{bmatrix}
5 & 2 \\
0 & -1
\end{bmatrix}\\
&=\bQ_3\bR_3\equiv
\begin{bmatrix}
-0.8 & -0.6 \\
-0.6 & 0.8
\end{bmatrix}
\begin{bmatrix}
-5 & -2 \\
0 & 1
\end{bmatrix}
&=& \bQ_4\bR_4\equiv
\begin{bmatrix}
-0.8 & 0.6 \\
-0.6 & -0.8
\end{bmatrix}
\begin{bmatrix}
-5 & -2 \\
0 & -1
\end{bmatrix}.
\end{aligned}
$$
Thus, the QR decomposition of $\bA$ is not unique.
\end{example}

The key observation of the above example can be stated in the following lemma, which highlights the non-uniqueness of the QR decomposition for a matrix.
\begin{lemma}[Non-uniqueness QR Decompositions]
Let $\bR$ be nonsingular and upper triangular. Then there exists a diagonal matrix $\bD$ with diagonal entries $\pm 1$ such that $\widehat{\bR}=\bD\bR$ is also nonsingular upper triangular and contains only positive diagonal entries.
\end{lemma}
\begin{exercise}
Find the diagonal matrices $\bD_2$, $\bD_3$, and $\bD_4$ that can be used to recover $\bR_1$ from $\bR_2, \bR_3$, and $\bR_4$ in Example~\ref{example:qr_unique}.
\end{exercise}

However, if we use the procedure described in the Gram-Schmidt process or systematically choose the sign in the Householder algorithm,  the QR decomposition becomes unique.
Moreover, the uniqueness of the \textit{reduced} QR decomposition for a matrix $\bA$ with full column rank  is guaranteed when $\bR$ has positive diagonals as demonstrated in the ``main proof" of Section~\ref{section:gram-schmidt-process} through inductive analysis.
We present an alternative proof for the uniqueness of the \textit{reduced} QR decomposition for matrices if the diagonal values of $\bR$ are positive. 
This proof provides insights into the \textit{implicit Q theorem} in the Hessenberg decomposition (Section~\ref{section:hessenberg-decomposition}) or tridiagonal decomposition (Section~\ref{section:tridiagonal-decomposition}).
\begin{corollary}[Uniqueness of the reduced QR Decomposition]\label{corollary:unique-qr}
Let $\bA$ be an $m\times n$ matrix with full column rank $n$ and $m\geq n$. The uniqueness of the \textit{reduced} QR decomposition is guaranteed when the main diagonal values of $\bR$ are positive.
\end{corollary}
\begin{proof}[of Corollary~\ref{corollary:unique-qr}]
Suppose the \textit{reduced} QR decomposition is not unique, we can complete it into a \textit{full} QR decomposition.
In this case, we can find two such full decompositions such that $\bA=\bQ_1\bR_1 = \bQ_2\bR_2$.
This implies $\bR_1 = \bQ_1^{-1}\bQ_2\bR_2 \triangleq \bV \bR_2$, where $\bV\triangleq \bQ_1^{-1}\bQ_2$ is an orthogonal matrix. 
By expanding the equation, we have:
$$
\begin{aligned}
	\bR_1 &= 
	\begin{bmatrix}
		r_{11} & r_{12}& \dots & r_{1n}\\
		& r_{22}& \dots & r_{2n}\\
		&       &    \ddots  & \vdots \\
		\multicolumn{2}{c}{\raisebox{1.3ex}[0pt]{\Huge0}} & & r_{nn} \\
		\bzero & \bzero &\ldots & \bzero
	\end{bmatrix}=
	\begin{bmatrix}
		v_{11}& v_{12} & \ldots & v_{1m}\\
		v_{21} & v_{22} & \ldots & v_{2m}\\
		\vdots & \vdots & \ddots & \vdots\\
		v_{m1} & v_{m2} & \ldots & v_{\textcolor{black}{mm}}
	\end{bmatrix}
	\begin{bmatrix}
		s_{11} & s_{12}& \dots & s_{1n}\\
		& s_{22}& \dots & s_{2n}\\
		&       &    \ddots  & \vdots \\
		\multicolumn{2}{c}{\raisebox{1.3ex}[0pt]{\Huge0}} & & s_{nn} \\
		\bzero & \bzero &\ldots & \bzero
	\end{bmatrix}= \bV\bR_2,
\end{aligned}
$$
This implies 
$$
r_{11} = v_{11} s_{11}, \qquad v_{21}=v_{31}=v_{41}=\ldots=v_{m1}=0.
$$
Since $\bV$ contains mutually orthonormal columns, and the first column of $\bV$ has a  norm of 1, we can conclude that $v_{11} = \pm 1$.
Given that $r_{ii}> 0$ and $s_{ii}> 0$ for $i\in \{1,2,\ldots,n\}$ based on the assumption, we have $r_{11}> 0$ and $s_{11}> 0$, implying that $v_{11}$ can only be positive 1. Additionally, since $\bV$ is an orthogonal matrix, we also have 
$$
v_{12}=v_{13}=v_{14}=\ldots=v_{1m}=0. 
$$
Applying this process to the submatrices of $\bR_1, \bV, $ and $\bR_2$, we will find the upper-left submatrix of $\bV$ is an identity matrix: $\bV[1:n,1:n]=\bI_n$, which means $\bR_1=\bR_2$. This, in turn, implies $\bQ_1[:,1:n]=\bQ_2[:,1:n]$ and leads to a contradiction such that the reduced QR decomposition is unique.
\end{proof}
We observe that the uniqueness of the aforementioned reduced QR decomposition is derived from the positivity of the diagonal elements in $\bR$. 
Consequently, if we impose the condition that the diagonal values of $\bR$ are positive in the Householder or Givens algorithms, the decomposition will be unique. 
In fact, the Gram-Schmidt process outlined in Algorithm~\ref{alg:reduced-qr} constitutes such a decomposition, where the diagonal elements of $\bR$ will be positive when $\bA$ has full column rank. 
However, if $\bA$ possesses dependent columns, the diagonal entries of $\bR$ can only be \textbf{nonnegative}, and the factorization may not be unique.

\paragraph{Uniqueness up to the signs of rows of $\bR$.}  When the restriction on the diagonal values of $\bR$ to be positive is lifted, the (reduced) QR decomposition loses its uniqueness. However, the uniqueness remains valid up to the signs of the rows of $\bR$, as illustrated in Example~\ref{example:qr_unique}.

\section{LQ, QL, RQ Decompositions}\label{section:lq-decomp}
We have previously established the existence of the QR decomposition through the Gram-Schmidt process.
In this case, our focus is on the column space of a given matrix $\bA=[\ba_1, \ba_2, \ldots, \ba_n] \in \real^{m\times n}$. The successive subspaces spanned by the columns $\ba_1, \ba_2, \ldots$ of $\bA$ are
$$
\cspace([\ba_1])\,\,\,\, \subseteq\,\,\,\, \cspace([\ba_1, \ba_2]) \,\,\,\,\subseteq\,\,\,\, \cspace([\ba_1, \ba_2, \ba_3])\,\,\,\, \subseteq\,\,\,\, \ldots,
$$
The concept behind the QR decomposition involves generating a sequence of orthonormal vectors $\bq_1, \bq_2, \ldots$ that span the same successive subspaces:
$$
\left\{\cspace([\bq_1])=\cspace([\ba_1]) \right\}\,\,\,\, \subseteq\,\,\,\, \{\cspace([\bq_1, \bq_2])=\cspace([\ba_1, \ba_2])\} \,\,\,\, \subseteq\,\,\,\, \ldots,
$$
However, in many applications (see \citet{schilders2009solution} and Appendix~\ref{appendix:KKT_proj}), there is also a need to consider the row space of a matrix $\bB=[\bb_1^\top; \bb_2^\top; \ldots;\bb_m^\top] \in \real^{m\times n}$, where $\bb_i$ is the $i$-th row of $\bB$. The successive subspaces spanned by the rows $\bb_1, \bb_2, \ldots$ of $\bB$ are
$$
\cspace([\bb_1])\,\,\,\, \subseteq\,\,\,\, \cspace([\bb_1, \bb_2]) \,\,\,\,\subseteq\,\,\,\, \cspace([\bb_1, \bb_2, \bb_3])\,\,\,\, \subseteq\,\,\,\, \ldots.
$$
The QR decomposition thus has its sibling, which finds the orthogonal row space.
By applying the QR decomposition to $\bB^\top = \bQ_0\bR$, we recover the LQ decomposition of the matrix $\bB = \bL \bQ$, where $\bQ \triangleq \bQ_0^\top$ and $\bL \triangleq \bR^\top$.
The LQ decomposition is helpful in demonstrating the existence of the UTV decomposition in the following chapter.

\begin{figure}[H]
\centering   
\vspace{-0.35cm}  
\subfigtopskip=2pt  
\subfigbottomskip=2pt  
\subfigcapskip=-5pt  
\subfigure[Reduced LQ decomposition.]{\label{fig:lqhalf}
	\includegraphics[width=0.47\linewidth]{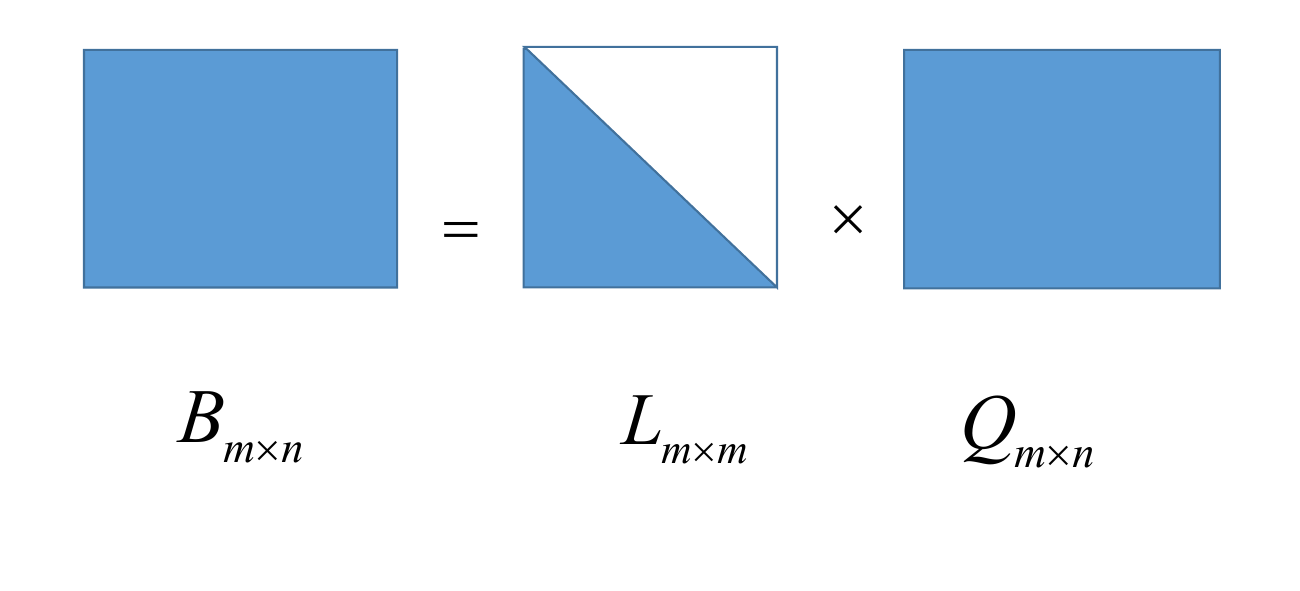}}
\quad 
\subfigure[Full LQ decomposition.]{\label{fig:lqall}
	\includegraphics[width=0.47\linewidth]{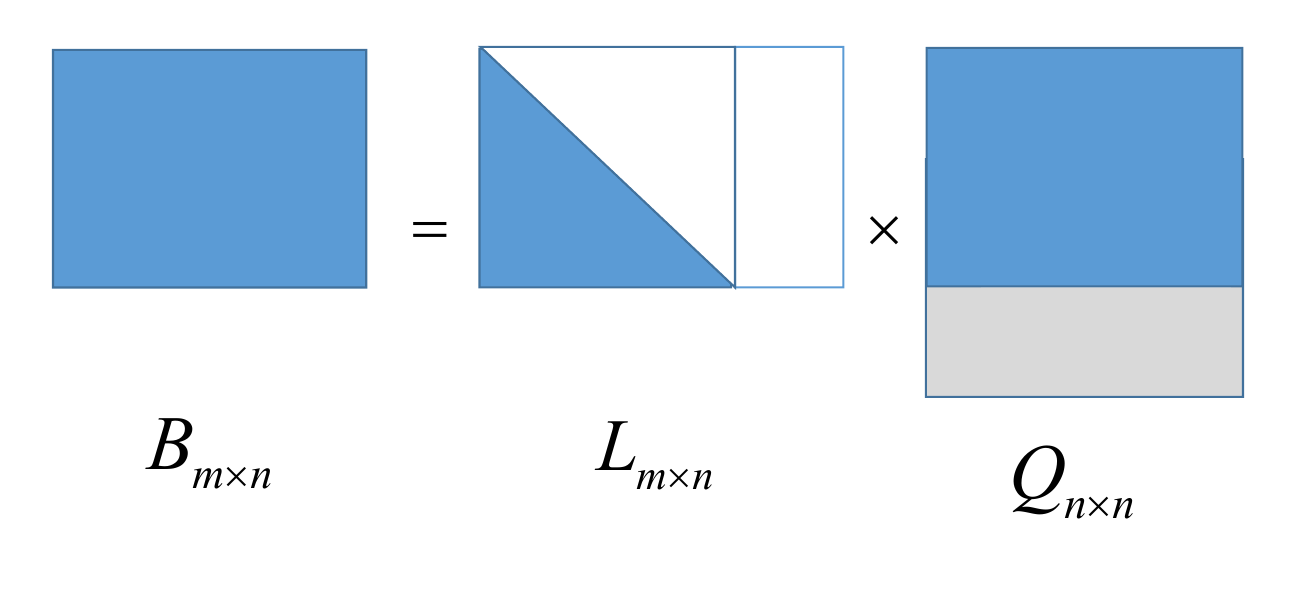}}
\caption{Comparison between the reduced and full LQ decomposition. White entries are zero, and blue entries are not necessarily zero. Gray columns denote silent rows.}
\label{fig:lq-comparison}
\end{figure}

\begin{theoremHigh}[LQ Decomposition]\label{theorem:lq-decomposition}
Let $\bB$ be any $m\times n$ matrix (regardless of whether its rows are linearly independent or dependent) with $n\geq m$. Then, it can be factored as 
$$
\bB = \bL\bQ,
$$
where 
\begin{enumerate}
\item \textbf{Reduced}: $\bL$ is an $m\times m$ lower triangular matrix, and $\bQ$ is $m\times n$ with orthonormal rows, known as the \textbf{reduced LQ decomposition};

\item \textbf{Full}: $\bL$ is an $m\times n$ lower triangular matrix, and $\bQ$ is $n\times n$ with orthonormal rows, known as the \textbf{full LQ decomposition}. If we further restrict the lower triangular matrix to be a square matrix, the full LQ decomposition can be denoted by 
$$
\bB = \begin{bmatrix}
	\bL_0 & \bzero
\end{bmatrix}\bQ,
$$
where $\bL_0$ is an $m\times m$ square lower triangular matrix.
\end{enumerate}
\end{theoremHigh}
Similarly, a comparison between the reduced and full LQ decomposition is shown in Figure~\ref{fig:lq-comparison}, where white entries are zero,  blue entries are not necessarily zero, and gray columns denote silent rows.

\paragraph{Row-pivoted LQ (RPLQ)\index{Row-pivoted}\index{RPLQ}.} Analogous to the column-pivoted QR decomposition discussed in Section~\ref{section:cpqr}, we can also define a \textit{row-pivoted LQ (RPLQ)} decomposition for a rank-$r$ matrix $\bB\in\real^{m\times n}$: 
$$
\left\{
\begin{aligned}
\text{Reduced RPLQ: }&\qquad 
\bP\bB &=& 
\underbrace{\begin{bmatrix}
	\bL_{11} \\
	\bL_{21}
\end{bmatrix}}_{m\times r}
\underbrace{\bQ_r }_{r\times n};\\
\text{Full RPLQ: }&\qquad 
\bP\bB &=& 
\underbrace{\begin{bmatrix}
	\bL_{11} & \bzero \\
	\bL_{21} & \bzero 
\end{bmatrix}}_{m\times m}
\underbrace{\bQ }_{m\times n},\\
\end{aligned}
\right.
$$
where $\bL_{11}\in \real^{r\times r}$ is lower triangular, $\bQ_r$ or $\bQ[1:r,:]$ spans the same row space as $\bB$, and $\bP$ is a permutation matrix that  moves independent rows to the top.

\index{QL decomposition}
\index{RQ decomposition}
\paragraph{QL and RQ decompositions.}
Consider the reduce QR decomposition for $\bA=\bQ\bR$, where $\bQ\in\real^{m\times n}$ is semi-orthogonal and $\bR\in\real^{n\times n}$, and let $\bP\in\real^{n\times n}$ be an $n$-th order reversal matrix (Definition~\ref{definition:permutation-matrix}).
Then, it follows that 
$$
\bA\bP = \underbrace{\bQ\bP}_{\bU} \underbrace{\bP\bR\bP}_{\bL},
$$ 
where $\bA\bP=[\ba_n, \ba_{n-1}, \ldots,\ba_1]$ indicates the column-reversed matrix of $\bA$, $\bL=\bP\bR\bP$ is a lower triangular matrix, $\bU=\bQ\bP=[\bq_n, \bq_{n-1}, \ldots, \bq_1]$ is the column-reversed matrix of $\bQ$ (which is also semi-orthogonal), and $m\geq n$. This is called the \textit{reduced QL} decomposition of $\bA$, and its full version can be identified similarly. 

Analogously, when $n\geq m$, we may find the RQ decomposition based on the LQ decomposition.

\section{Two-Sided Orthogonal Decomposition}
To this end, a direct implication of the CPQR and RPLQ decompositions is the \textit{two-sided orthogonal decomposition}, which simultaneously determines the orthonormal row basis and orthonormal column basis.
\begin{theoremHigh}[Two-Sided Orthogonal Decomposition]\label{theorem:two-sided-orthogonal}
Let $\bA\in \real^{n\times n}$ be a square matrix with rank $r$. Then, the full CPQR and RPLQ decompositions of $\bA$ are given by 
$$\bA\bP_1=\bQ_1
\begin{bmatrix}
\bR_{11} & \bR_{12}\\
\bzero & \bzero
\end{bmatrix}
\quad \text{and}\quad
\bP_2\bA=
\begin{bmatrix}
\bL_{11} & \bzero \\
\bL_{21} & \bzero 
\end{bmatrix}
\bQ_2,
$$ 
respectively. Then we would find out 
$$
\bA\bP\bA = \bQ_1
\underbrace{\begin{bmatrix}
	\bR_{11}\bL_{11}+\bR_{12}\bL_{21} & \bzero \\
	\bzero & \bzero 
\end{bmatrix}}_{\text{rank $r$}}
\bQ_2,
$$
where the first $r$ columns of $\bQ_1$ span the same column space as $\bA$, the first $r$ rows of $\bQ_2$ span the same row space as $\bA$, and $\bP$ is a permutation matrix. We refer to this decomposition as the \textit{two-sided orthogonal decomposition}.
\end{theoremHigh}

This decomposition shares a similar property with the singular value decomposition (SVD): $\bA=\bU\bSigma\bV^\top$, where the first $r$ columns of $\bU$ span the same column space as $\bA$ and the first $r$ columns of $\bV$ span the same row space as $\bA$ (as detailed in  Proposition~\ref{proposition:svd-four-orthonormal-Basis}). Therefore, the two-sided orthogonal decomposition can be regarded as an inexpensive alternative in this regard.

\begin{proposition}[Four Orthonormal Basis]
Given the two-sided orthogonal decomposition of matrix $\bA\in \real^{n\times n}$ with rank $r$: $\bA\bP\bA = \bU \bF\bV^\top$, where $\bU=[\bu_1, \bu_2, \ldots,\bu_n]$ and $\bV=[\bv_1, \bv_2, \ldots, \bv_n]$ are the column partitions of $\bU$ and $\bV$, respectively, we can observe the following property:
\begin{itemize}
\item $\{\bv_1, \bv_2, \ldots, \bv_r\} $ is an orthonormal basis of $\cspace(\bA^\top)$;

\item $\{\bv_{r+1},\bv_{r+2}, \ldots, \bv_n\}$ is an orthonormal basis of $\nspace(\bA)$;

\item $\{\bu_1,\bu_2, \ldots,\bu_r\}$ is an orthonormal basis of $\cspace(\bA)$;

\item $\{\bu_{r+1}, \bu_{r+2},\ldots,\bu_n\}$ is an orthonormal basis of $\nspace(\bA^\top)$. 
\end{itemize}
\end{proposition}

\section{Applications}
\index{Least squares}
\subsection{Application: Least Squares via the Full QR Decomposition}\label{section:application-ls-qr}
Let us consider the overdetermined system 
\begin{equation}\label{equation:qr_lr_1}
\bA\bx = \bb,
\end{equation}
where $\bA\in \real^{m\times n}$ is the data matrix, and $\bb\in \real^m$ with $m\geq n$ is the observation matrix (when $n>m$, any linear
model is over-parametrized, and regularization is needed to achieve a stable
fit; $\ell_1$ or even $\ell_0$ regularizations are necessary to find a sparse solution) \footnote{
In this system, $m$ represents the number of equations and $n$ represents the number of unknowns.
When $m>n$, we say the linear system has more equations than unknowns.
This scenario attracts more attention because repeated measurements are taken to minimize errors in the system.
While in the setting of \textit{compressed sensing}, we usually receive the vector $\bb$ and want to reconstruct the original signal $\bx$, which is assumed to be sparse and $n>m$, i.e., reconstructing the signal from a low-dimensional vector.
It is also worth mentioning that when $n>m$, we may have many solutions to the system, and the \textit{Hessian} of the least squares problem is singular. Therefore, regularization is required to obtain a meaningful solution.
}. 
Since the number of equations $m$ exceed the number of unknowns $n$, the linear system may not have a solution; alternatively, due to measurement errors and other factors, the equality relations may not hold precisely.
However, $\bA$ will typically have full column rank since the data from real-world applications often have a low chance of being dependent (or become independent after post-processing).
And the least squares (LS) solution for this system can be obtained using $\bx_{LS} = (\bA^\top\bA)^{-1}\bA^\top\bb$ to minimize the squared $\ell_2$ norm of the residual, which is denoted by 
\begin{equation}\label{equation:qr_lr_2}
\normtwo{\bA\bx-\bb}^2.~\footnote{The least squares method uses the $\ell_2$ norm to measure the size of the error because the $\ell_2$ norm is smooth and differentiable, and corresponds to Gaussian additive noise; see Section~\ref{section:more_err_sta_als}.}
\end{equation}
The idea of least squares is very intuitive: if system~\eqref{equation:qr_lr_1} has a solution, then finding the global optimal solution to problem~\eqref{equation:qr_lr_2} is equivalent to finding the solution to the system. 
When the system does not have a solution, problem~\eqref{equation:qr_lr_2} actually provides a solution with the smallest possible error in some sense.
Here, $\bA^\top\bA$ is invertible due to the full column rank of $\bA$ and the fact that $\rank(\bA^\top\bA)=\rank(\bA)$.

However, computing the inverse of a matrix is a nontrivial task. Instead, we can employ the QR decomposition method to find the least squares solution, as demonstrated in the following theorem.
\begin{theorem}[LS via QR for Full Column Rank Matrix]\label{theorem:qr-for-ls}
Let $\bA\in \real^{m\times n}$, and let $\bA=\bQ\bR$ be its full QR decomposition, where  $\bQ\in\real^{m\times m}$ is an orthogonal matrix,  $\bR\in \real^{m\times n}$ is an upper triangular matrix appended by additional $m-n$ zero rows, and $\bA$ has full column rank with $m\geq n$.~\footnote{In Section~\ref{section:ls-utv} and~\ref{section:application-ls-svd}, we briefly discuss how to use UTV decomposition and SVD to address the rank-deficient least squares problems.
See Problem~\ref{prob:als_pseudo1}$\sim$\ref{prob:als_pseudon} for more results.}
Suppose $\bR = 
\footnotesize
\begin{bmatrix}
\bR_1 \\
\bzero
\end{bmatrix}$, where $\bR_1 \in \real^{n\times n}$ is the square upper triangular part of $\bR$, and $\bb\in \real^m$. Then, the LS solution to $\bA\bx=\bb$ is given by 
$$
\bx_{LS} = \bR_1^{-1}\bc, 
$$
where $\bc$ is the first $n$ components of $\bQ^\top\bb$.
\end{theorem}

\begin{proof}[of Theorem~\ref{theorem:qr-for-ls}]
Given that $\bA=\bQ\bR$ represents the full QR decomposition of $\bA$ and $m\geq n$, the last $m-n$ rows of $\bR$ are zero, as shown in Figure~\ref{fig:qr-comparison}. Hence, $\bR_1 \in \real^{n\times n}$ is the square upper triangular part of $\bR$, and 
$
\bQ^\top \bA = \bR = 
\footnotesize\begin{bmatrix}
\bR_1 \\
\bzero
\end{bmatrix}.
$
Thus, we have
$$
\begin{aligned}
\normtwo{\bA\bx-\bb}^2 &= (\bA\bx-\bb)^\top(\bA\bx-\bb)\stackrel{*}{=}(\bA\bx-\bb)^\top\bQ\bQ^\top (\bA\bx-\bb) \\
&\stackrel{+}{=}\normtwo{\bQ^\top \bA \bx-\bQ^\top\bb}^2 
=\normtwo{\begin{bmatrix}
		\bR_1 \\
		\bzero
	\end{bmatrix} \bx-\bQ^\top\bb}^2
=\normtwo{\bR_1\bx - \bc}^2+\normtwo{\bd}^2,
\end{aligned}
$$ 
\footnote{equality ($*$) follows from the orthogonality of $\bQ$; equality (+) follows from the invariance of the norm under orthogonal transformations.}~where $\bc$ represents the first $n$ components of $\bQ^\top\bb$, and $\bd$ represents the last $m-n$ components of $\bQ^\top\bb$. Then the LS solution can be obtained by performing back substitution on the upper triangular system $\bR_1\bx = \bc$, which can be expressed as $\bx_{LS} = \bR_1^{-1}\bc$.
\end{proof}

The linear system $\bR_1\bx=\bc$ has an exact solution since $\bR_1$ has full rank (invertible). Therefore, the \textit{least squares residual} of $\normtwo{\bA\bx-\bb}$ is 
$$
\mathop{\min}_{\bx} \normtwo{\bA\bx-\bb}^2 =
\mathop{\min}_{\bx} \normtwo{\bR_1\bx - \bc}^2+\normtwo{\bd}^2
=
\normtwo{\bd}^2.
$$

To verify Theorem~\ref{theorem:qr-for-ls}, we consider the full QR decomposition of $\bA = \bQ\bR$, where $\bQ\in \real^{m\times m}$ and $\bR\in \real^{m\times n}$. 
Along with the LS solution $\bx_{LS} = (\bA^\top\bA)^{-1}\bA^\top\bb$, we obtain
\begin{equation}\label{equation:qr-for-ls-1}
\begin{aligned}
\bx_{LS} 
&= (\bA^\top\bA)^{-1}\bA^\top\bb 
= (\bR^\top\bQ^\top\bQ\bR)^{-1} \bR^\top\bQ^\top \bb
= (\bR^\top\bR)^{-1} \bR^\top\bQ^\top \bb \\
&= (\bR_1^\top\bR_1)^{-1} \bR^\top\bQ^\top \bb 
=\bR_1^{-1} \bR_1^{-\top} \bR^\top\bQ^\top \bb
= \bR_1^{-1} \bR_1^{-\top} \bR_1^\top\bQ_1^\top \bb 
=\bR_1^{-1} \bQ_1^\top \bb,
\end{aligned}
\end{equation}
where $\bR = 
\footnotesize
\begin{bmatrix}
\bR_1 \\
\bzero
\end{bmatrix}$ and $\bR_1\in \real^{n\times n}$  is an upper triangular matrix.
Additionally,  $\bQ_1 =\bQ_{:,1:n}\in \real^{m\times n}$ denotes the first $n$ columns of $\bQ$ (i.e., $\bQ_1\bR_1$ is the reduced QR decomposition of $\bA$). Then the result of Equation~\eqref{equation:qr-for-ls-1} agrees with Theorem~\ref{theorem:qr-for-ls}.

In conclusion, utilizing the QR decomposition allows us to directly derive the least squares solution, which aligns with the proposition stated in Theorem~\ref{theorem:qr-for-ls}.
Additionally, we indirectly validate the least squares result (based on calculus) through the QR decomposition, confirming the consistency  between the two approaches. 
For readers interested in the linear algebra aspect of least squares, 
Appendix~\ref{appendix:ls-fundation-theorem} provides a visual representation of the least squares method for a full column rank matrix $\bA$ based on the fundamental theorem of linear algebra.
Alternatively, a comprehensive discussion can be found in \citet{lu2021revisit}.

\index{Fundamental theorem}

\paragraph{Complexity.} 
In the least squares solution using QR decomposition, it is unnecessary to compute the orthogonal matrix $\bQ$ entirely. This can be observed by performing the Householder algorithm (Algorithm~\ref{alg:qr-decomposition-householder}) on the augmented matrix $[\bA~|~\bb]$:
\begin{equation}\label{equation:ls_hourse}
[\bA~|~\bb] \rightarrow \bQ^\top[\bA~|~\bb] 
=
\begin{bmatrix}
	\bQ_1^\top \\
	\bQ_2^\top
\end{bmatrix}
[\bA~|~\bb]
=
\begin{bmatrix}
	\bR_1 & \bQ_1^\top\bb \\
	\bzero & \bQ_2^\top\bb
\end{bmatrix},
\end{equation}
where $\bQ=[\bQ_1, \bQ_2]\in\real^{m\times m}$ with $\bQ_1\in\real^{m\times n}$ and $\bQ_2 \in\real^{m\times (m-n)}$. 
Applying the Householder algorithm allows us to directly compute $\bR_1$ and $\bQ_1^\top\bb$.
Hence, the complexity of the least squares method using QR decomposition corresponds to the first part of Theorem~\ref{theorem:qr-full-householder}, which is approximately $\sim 2mn^2-\frac{2}{3}n^3$ flops. 
The computation cost of obtaining $\bQ_1^\top\bb$ is not included in the final flops count, as it requires $\mathcalO(mn)$ flops (which is not significant compared to the leading term).

\paragraph{Update of least squares problem.}
In the context of least squares, each row of $\bA$ and $\bb$ is  referred to as an \textit{observation}.
In real-world application, new observations may be received. 
When performing the optimization process from scratch, obtaining the solution of the least squares problem would require approximately $\sim 2(m+1)n^2-\frac{2}{3}n^3$ flops. 
Building upon Equation~\eqref{equation:ls_hourse}, let's consider a new observation $[\ba^\top,\beta]$ (more generally, one can consider a set of new observations rather than just one), which leads us to the following reduction:
\begin{equation}\label{equation:qr_update}
\left[
\begin{array}{c|c}
	\bA & \bb  \\
	\ba^\top & \beta
\end{array}
\right]
\rightarrow 
\begin{bmatrix}
	\bQ^\top & \bzero \\
	\bzero & 1
\end{bmatrix}
\left[
\begin{array}{c|c}
	\bA & \bb  \\
	\ba^\top & \beta
\end{array}
\right]
=
\underbrace{\begin{bmatrix}
		\bR_1 & \bQ_1^\top\bb \\
		\bzero & \bQ_2^\top\bb \\
		\ba^\top & \beta
\end{bmatrix}}_{\triangleq\bZ}
\end{equation}
Therefore, the updated least squares solution is obtained by transforming $\bZ\in\real^{(m+1)\times (n+1)}$ into an upper triangular matrix (actually, we transform the left $n$ columns of $\bZ$ into an upper triangular matrix). This can be done by applying a set of rotations in the $(1,m+1)$ plane, $(2,m+1)$ plane, $\ldots$, $(n,m+1)$ plane that introduce zero in the $(m+1, 1), (m+1,2),\ldots,(m+1,n)$-th entries of $\bZ$, respectively. 
The computational cost for this operation is $\mathcalO(mn)$ flops. See Section~\ref{section:qr_adddelrow} for more details.

\index{Perturbation}
\index{Least squares}
\paragraph{Perturbation theory of LS via Householder QR.}
Following \citet{higham2002accuracy, elden2007matrix}, it has been established that the least squares solution of a perturbed linear system remains stable when computed using the Householder QR algorithm.
\begin{theorem}[Perturbed LS via Householder QR]
Let $\bA\in\real^{m\times n}$ with $m\geq n$ and full rank. Then, the solution $\bx_{LS}$ computed using Householder QR decomposition as stated in Theorem~\ref{theorem:qr-for-ls} is the exact least squares solution of 
$$
\mathop{\min}_{\bx} 
\normtwo{
(\bA+\delta\bA)\bx - (\bb+\delta\bb)
},
$$
where 
$$
\normf{\delta\bA} \leq c_1 mn \mu \normf{\bA} + \mathcalO(\mu^2), \gap 
\normtwo{\delta\bb} \leq c_2 mn\normtwo{\bb} + \mathcalO(\mu^2),
$$
and $c_1$ and $c_2$ are small constants.
\end{theorem}

\begin{exercise}
Examine the least squares problem using column-pivoted QR decomposition (Theorem~\ref{theorem:rank-revealing-qr-general}) for rank-deficient matrices.
\end{exercise}

\index{Gauss-Newton method}
\index{Nonlinear least squares}
\subsection{Application: Gauss-Newton and Levenberg-Marquardt Method}
In the previous section, we considered the least squares problem for linear systems, which is often referred to as linear least squares.
When the residual $\br(\bx)$ in Equation~\eqref{equation:qr_lr_2} is not linear, we have the nonlinear least squares problem~\footnote{More details can refer to, for example, \citet{madsen2004methods}.}:
$$
\bx^* = \mathop{\argmax}_{\bx} \left\{f(\bx) = \frac{1}{2} \normtwo{\br(\bx)}^2\right\}, 
\gap \br(\bx)\in\real^m, \,\bx\in\real^n, \, m\geq n.
$$
When $\br(\bx)=\bA\bx-\bb$, this reduces to the linear least squares problem in Equation~\eqref{equation:qr_lr_2}.
The gradient and Hessian of $f(\bx)$ are 
\begin{equation}\label{equation:gaus_new_jaco}
\begin{aligned}
\nabla f(\bx) &= \bJ(\bx)^\top \br(\bx)
\quad\text{and} \quad
\nabla^2 f(\bx) = \bJ(\bx)^\top\bJ(\bx) + \sum_{i=1}^{m} r_i(\bx)\nabla^2r_{i}(\bx),
\end{aligned}
\end{equation}
where $\bJ(\bx)\in\real^{m\times n}$ is the Jacobian matrix (see Problem~\ref{problem:gaus_new_jco}).
The plain Newton's method~\footnote{See, for example, \citet{lu2022gradient}.} is an iterative update method that considers the update at the $t$-th iteration by 
$$
\bx^{(t+1)} \leftarrow \bx^{(t)} + \bd^{(t)},
$$
where $(\nabla^2 f(\bx^{(t)}) )\bd^{(t)} = -\nabla f(\bx^{(t)})$ determines  the ``candidate" descent direction $\bd^{(t)}$.
For brevity, we omit the superscript $t$ and 
use the linear Taylor's approximation (Theorem~\ref{theorem:linear_approx}), we have 
$$
\nabla f(\bx + \bd) 
\approx
\nabla f(\bx ) +\nabla^2 f(\bx )^\top \bd.
$$
Therefore, Newton's method can be understood as finding $\bd$ making $\nabla f(\bx + \bd) =\nabla f(\bx ) +\nabla^2 f(\bx )^\top \bd$ towards $\bzero$ (i.e., a stationary point).
Taking the quadratic Taylor's approximation (Theorem~\ref{theorem:quad_app_theo}), we have 
\begin{equation}\label{equation:ng_secondapp}
f(\bx + \bd) 
= 
f(\bx ) +\nabla f(\bx )^\top \bd + \frac{1}{2}\bd^\top \nabla^2f(\bx) \bd + o(\normtwo{\bd}^2).
\end{equation}
When $\bx$ is a stationary point, $\nabla f(\bx )=\bzero$. Suppose further that $f(\bx)$ has a positive definite Hessian: $\nabla^2f(\bx)\succ 0$; this implies  that $\nabla^2f(\bx)$ has a minimum eigenvalue $\lambda_{\min}>0$ (Theorem~\ref{theorem:eigen_charac}) and $\bd^\top \nabla^2f(\bx) \bd \geq \lambda \normtwo{\bd}^2$ for all $\lambda_{\min}>\lambda>0$. This again implies that the third term in \eqref{equation:ng_secondapp} dominates the fourth term. Therefore, $\bx$ is a \textit{local minimizer} (a minimum point around the neighborhood of $\bx$ with some radius $r$, $B(\bx,r)$) when $\bx$ is a stationary point and $\nabla^2f(\bx)$ is positive definite.

\paragraph{Gauss-Newtom method.}
However, since  the Hessian $\nabla^2r_{i}(\bx)$ can be intractable, the Gauss-Newton method considers to use $\bJ(\bx)^\top\bJ(\bx)$ to approximate the Hessian $\nabla^2 f(\bx)$ and obtains the ``candidate" descent direction
$$
\bJ(\bx^{(t)})^\top\bJ(\bx^{(t)})\bd^{(t)} = -\bJ(\bx^{(t)})^\top \br(\bx^{(t)}).
$$
The ``candidate" descent direction can be equivalently obtained by 
\begin{equation}\label{equation:gasnewton}
\textbf{(Gauss-Newton):}\gap \bd^{(t)} = \mathop{\argmin}_{\bd} \normtwo{\bJ(\bx^{(t)})\bd+\br(\bx^{(t)})}^2,
\end{equation}
which is a linear least squares problem and can be solved using Theorem~\ref{theorem:qr-for-ls} (when $\bJ(\bx^{(t)})$ has full rank).
Let $\bJ(\bx^{(t)})$ admit the reduced QR decomposition $\bJ(\bx^{(t)})=\bQ^{(t)}\bR^{(t)}$, the ``candidate" descent direction can be obtained by 
$$
\bd^{(t)} \leftarrow (\bR^{(t)})^{-1}(\bQ^{(t)})^\top\br(\bx^{(t)}).
$$
Thus, there is no need to compute the inverse of $\bJ(\bx^{(t)})^\top\bJ(\bx^{(t)})$.
When $(\bd^{(t)})^\top\nabla f(\bx^{(t)})\leq 0$, the  direction $\bd^{(t)}$ is called a ``true" descent direction (as opposed to the ``candidate" descent direction we used previously).
We can confirm that when $\bJ(\bx^{(t)})$ has full rank $n$ (since $m\geq n$), 
$$
(\bd^{(t)})^\top\nabla f(\bx^{(t)}) = (\bd^{(t)})^\top \bJ(\bx^{(t)})^\top \br(\bx^{(t)}) = -\normtwo{\bJ(\bx^{(t)})\bd^{(t)} }^2\leq 0.
$$
Therefore, the obtained direction $\bd^{(t)}$ is a ``true" descent direction.

\index{Levenberg-Marquardt method}
\index{Trust region method}
\index{KKT condition}
\paragraph{Levenberg-Marquardt (LM) method.}
The Levenberg-Marquardt method also considers the same problem in \eqref{equation:gasnewton}, but with an additional constraint $\normtwo{\bd}\leq \Delta^{(t)}$ \citep{levenberg1944method, marquardt1963algorithm, wright1985inexact}:
\begin{equation}\label{equation:lmmethod}
\textbf{(LM-1):}\gap \bd^{(t)} = \mathop{\argmin}_{\bd} \normtwo{\bJ(\bx^{(t)})\bd+\br(\bx^{(t)})}^2, \gap \text{s.t.}\gap \normtwo{\bd}\leq \Delta^{(t)}.
\end{equation}
This is equivalently to, using Lagrange  multiplier, the following problem
\begin{equation}\label{equation:lmmethod_lag}
\begin{aligned}
\textbf{(LM-2):}\gap 
\bd^{(t)} &= \mathop{\argmin}_{\bd} \normtwo{\bJ(\bx^{(t)})\bd+\br(\bx^{(t)})}^2 + \lambda\normtwo{\bd}^2\\
&= \mathop{\argmin}_{\bd} \normtwo{\begin{bmatrix}
	\bJ(\bx^{(t)})\\
	\sqrt{\lambda}\bI 
\end{bmatrix}\bd
+
\begin{bmatrix}
\br(\bx^{(t)})\\
\bzero 
\end{bmatrix}
}^2,
\end{aligned}
\end{equation}
where $\lambda$ is a Lagrange multiplier related to $\Delta^{(t)}$. The second form above is an updated least squares problem. Given the knowledge of the QR decomposition of $\bJ(\bx^{(t)}) =\bQ^{(t)}\bR^{(t)}$, the least squares problem can be solved using the \textit{update of least squares problem} (Equation~\eqref{equation:qr_update}).

The update from LM method is a special form of \textit{trust region method}, where the function value at the $t$-th iteration is updated by a quadratic approximation (assume $f(\bx)$ is twice continuously differentiable, Theorem~\ref{theorem:quad_app_theo}):
\begin{equation}\label{equation:trust_reg}
\bd^{(t)} = 
\mathop{\argmin}_{\bd}
\big\{
f(\bx^{(t)}+\bd) 
\approx 
f(\bx^{(t)}) +\nabla f(\bx^{(t)})^\top\bd +\frac{1}{2} \bd^\top \bH^{(t)} \bd
\big\},
\,  \text{s.t.}\, \normtwo{\bd}\leq \Delta^{(t)},
\end{equation}
where the ball $B[\bx^{(t)}, \Delta^{(t)}]=\{\bx\in\real^n: \normtwo{\bx-\bx^{(t)}} \leq \Delta^{(t)}\}$ is called the \textit{trust region}, and $\Delta^{(t)}$ is called the \textit{trust region radius}.

In the trust region method, if the step is deemed good (the gain ratio of the actual reduction to the predicted reduction is  sufficiently large), the current point is updated to the new point, and the trust region may be expanded.
If the step is not good (the ratio is small), the current point remains unchanged, and the trust region is usually shrunk to improve the accuracy of the model.

\begin{theorem}[Trust Region Method]
Consider the problem
$$
\mathop{\argmin}_{\bd}\, g(\bd) = \bff +\bg^\top\bd +\frac{1}{2} \bd^\top\bH\bd \gap\text{s.t.}\gap \normtwo{\bd}\leq \Delta.
$$
Then $\bd^*$ is the optimum if and only if there exists a value $\lambda\geq 0$ such that
$$
\begin{aligned}
(\bH+\lambda\bI)\bd^* &= -\bg;\\
\lambda (\Delta - \normtwo{\bd^*}) &=\bzero; \\
(\bH+\lambda \bI) &\text{ is positive semidefinite}.
\end{aligned}
$$ 
\end{theorem}
The first two conditions are from KKT conditions (see, for example, Appendix~\ref{appendix:KKT_proj}), and the third condition comes from $g(\bd)\geq g(\bd^*)$. The proof is left as an exercise. 

The comparison between \eqref{equation:lmmethod_lag} and \eqref{equation:trust_reg} reveals that the LM method is a special form of trust region method by approximating  $\bH^{(t)}\approx\bJ(\bx^{(t)})^\top\bJ(\bx^{(t)})$. Therefore, one may also consider the following optimization to solve the LM method:
$$
\textbf{(LM-3):}\gap \bd^{(t)} = \mathop{\argmin}_{\bd} \big(\bJ(\bx^{(t)})^\top\bJ(\bx^{(t)}) + \lambda\bI\big)\bd =-\bJ(\bx^{(t)})\br(\bx^{(t)}).
$$

\index{Rank-one update}
\subsection{Application: Rank-One Changes}\label{section:qr-rank-one-changes}
We have previously explored the rank-one update/downdate of the Cholesky decomposition in Section~\ref{section:cholesky-rank-one-update}. 
The rank-one change of the matrix $\bA$ in the QR decomposition, denoted by $\bA^\prime$, is defined in a similar manner:
$$
\begin{aligned}
	\bA^\prime &= \bA + \bu\bv^\top, \\
	\downarrow &\gap  \downarrow\\
	\bQ^\prime\bR^\prime &=\bQ\bR + \bu\bv^\top,
\end{aligned}
$$
where if we set $\bA^\prime = \bA - (-\bu)\bv^\top$, we recover the downdate form such that the update and downdate in the QR decomposition are equivalent. 
To reiterate,
the rank-one update/downdate problem involves finding the QR decomposition of $\bA^\prime$ when the QR decomposition of $\bA$ has already been computed.
Let $\bw \triangleq \bQ^\top\bu$, then we have
$
\bA^\prime = \bQ(\bR + \bw\bv^\top).
$
Based on the second form described in Remark~\ref{remark:basis-from-givens2}  regarding the introduction of zeros in a backward manner, there exists a set of Givens rotations $\bG_{12} \bG_{23} \ldots \bG_{(n-1),n}$ that satisfy:
$$
\bG_{12} \bG_{23} \ldots \bG_{(n-1),n} \bw = \pm \normtwo{\bw} \be_1,
$$
where $\bG_{(k-1),k}$ represents the Givens rotation effecting in the $(k-1)$-th and $k$-th plane, which introduces zero in the $k$-th entry of $\bw$. Applying these rotations to $\bR$, we obtain
$
\bG_{12} \bG_{23} \ldots \bG_{(n-1),n}\bR \triangleq \bH_0 ,
$
where the Givens rotations in this \textit{reverse order} (\textit{backward rotations}) are employed to transform the upper triangular matrix $\bR$ into a ``simple" \textit{upper Hessenberg matrix}, which closely resembles upper triangular matrices (see Definition~\ref{definition:upper-hessenbert} that we will introduce in the Hessenberg decomposition chapter). 
However, if the rotations transform $\bw$ into $\pm \normtwo{\bw}\be_1$ in the \textit{forward order} (\textit{forward rotations}) as demonstrated in Corollary~\ref{corollary:basis-from-givens}, we will not obtain this upper Hessenberg $\bH_0$. To see this, consider the matrix $\bR\in \real^{5\times 5}$, an example is shown as follows, where $\boxtimes$ represents a value that is not necessarily zero, and \textbf{boldface} indicates the value has just been changed. Using backward rotations results in the upper Hessenberg $\bH_0$, which is easier to manage:
$$
\begin{aligned}
\text{\parbox{7em}{Backward\\(Right Way)}: }
\footnotesize
\begin{sbmatrix}{\bR}
\boxtimes & \boxtimes & \boxtimes & \boxtimes& \boxtimes \\
0 & \boxtimes & \boxtimes & \boxtimes& \boxtimes \\
0 & 0 & \boxtimes & \boxtimes& \boxtimes \\
0 & 0 & 0 & \boxtimes& \boxtimes \\
0 & 0 & 0 & 0& \boxtimes
\end{sbmatrix}
&\stackrel{\bG_{45}}{\rightarrow}
\footnotesize
\begin{sbmatrix}{\bG_{45}\bR}
\boxtimes & \boxtimes & \boxtimes & \boxtimes& \boxtimes \\
0 & \boxtimes & \boxtimes & \boxtimes& \boxtimes \\
0 & 0 & \boxtimes & \boxtimes& \boxtimes \\
0 & 0 & 0 & \bm{\boxtimes}& \bm{\boxtimes} \\
0 & 0 & 0 & \bm{\boxtimes}& \bm{\boxtimes}
\end{sbmatrix}
\stackrel{\bG_{34}}{\rightarrow}
\footnotesize
\begin{sbmatrix}{\bG_{34}\bG_{45}\bR}
\boxtimes & \boxtimes & \boxtimes & \boxtimes& \boxtimes \\
0 & \boxtimes & \boxtimes & \boxtimes& \boxtimes \\
0 & 0 & \bm{\boxtimes} & \bm{\boxtimes}& \bm{\boxtimes} \\
0 & 0 & \bm{\boxtimes} & \bm{\boxtimes}& \bm{\boxtimes} \\
0 & 0 & 0 & \boxtimes & \boxtimes
\end{sbmatrix}\\
&\stackrel{\bG_{23}}{\rightarrow}
\footnotesize
\begin{sbmatrix}{\bG_{23}\bG_{34}\bG_{45}\bR}
\boxtimes & \boxtimes & \boxtimes & \boxtimes& \boxtimes \\
0 & \bm{\boxtimes} & \bm{\boxtimes} & \bm{\boxtimes}& \bm{\boxtimes} \\
0 & \bm{\boxtimes} & \bm{\boxtimes} & \bm{\boxtimes}& \bm{\boxtimes} \\
0 & 0 & \boxtimes & \boxtimes& \boxtimes \\
0 & 0 & 0 & \boxtimes& \boxtimes
\end{sbmatrix}
\stackrel{\bG_{12}}{\rightarrow}
\footnotesize
\begin{sbmatrix}{\bG_{12}\bG_{23}\bG_{34}\bG_{45}\bR}
\bm{\boxtimes} & \bm{\boxtimes} & \bm{\boxtimes} & \bm{\boxtimes}& \bm{\boxtimes} \\
\bm{\boxtimes} & \bm{\boxtimes} & \bm{\boxtimes} & \bm{\boxtimes}& \bm{\boxtimes} \\
0 & \boxtimes & \boxtimes & \boxtimes& \boxtimes \\
0 & 0 & \boxtimes & \boxtimes& \boxtimes \\
0 & 0 & 0 & \boxtimes& \boxtimes
\end{sbmatrix}.
\end{aligned}
$$
And the forward rotations result in a \textbf{full (non-sparse) matrix}:
$$
\begin{aligned}
	\text{\parbox{7em}{Forward\\(Wrong Way)}: }
	\footnotesize
	\begin{sbmatrix}{\bR}
		\boxtimes & \boxtimes & \boxtimes & \boxtimes& \boxtimes \\
		0 & \boxtimes & \boxtimes & \boxtimes& \boxtimes \\
		0 & 0 & \boxtimes & \boxtimes& \boxtimes \\
		0 & 0 & 0 & \boxtimes& \boxtimes \\
		0 & 0 & 0 & 0& \boxtimes
	\end{sbmatrix}
	&\stackrel{\bG_{12}}{\rightarrow}
	\footnotesize
	\begin{sbmatrix}{\bG_{12}\bR}
		\bm{\boxtimes} & \bm{\boxtimes} & \bm{\boxtimes} & \bm{\boxtimes}& \bm{\boxtimes} \\
		\bm{\boxtimes} & \bm{\boxtimes} & \bm{\boxtimes} & \bm{\boxtimes}& \bm{\boxtimes} \\
		0 & 0 & \boxtimes & \boxtimes& \boxtimes \\
		0 & 0 & 0 & \boxtimes& \boxtimes \\
		0 & 0 & 0 & 0& \boxtimes
	\end{sbmatrix}
	\stackrel{\bG_{23}}{\rightarrow}
	\footnotesize
	\begin{sbmatrix}{\bG_{23}\bG_{12}\bR}
		\boxtimes & \boxtimes & \boxtimes & \boxtimes& \boxtimes \\
		\bm{\boxtimes} & \bm{\boxtimes} & \bm{\boxtimes} & \bm{\boxtimes}& \bm{\boxtimes} \\
		\bm{\boxtimes} & \bm{\boxtimes} & \bm{\boxtimes} & \bm{\boxtimes}& \bm{\boxtimes} \\
		0 & 0 & 0 & \boxtimes& \boxtimes \\
		0 & 0 & 0 & 0& \boxtimes
	\end{sbmatrix}\\
	&\stackrel{\bG_{34}}{\rightarrow}
	\footnotesize
	\begin{sbmatrix}{\bG_{34}\bG_{23}\bG_{12}\bR}
		\boxtimes & \boxtimes & \boxtimes & \boxtimes& \boxtimes \\
		\boxtimes & \boxtimes & \boxtimes & \boxtimes& \boxtimes \\
		\bm{\boxtimes} & \bm{\boxtimes} & \bm{\boxtimes} & \bm{\boxtimes}& \bm{\boxtimes} \\
		\bm{\boxtimes} & \bm{\boxtimes} & \bm{\boxtimes} & \bm{\boxtimes}& \bm{\boxtimes} \\
		0 & 0 & 0 & 0& \boxtimes
	\end{sbmatrix}
	\stackrel{\bG_{45}}{\rightarrow}
	\footnotesize
	\begin{sbmatrix}{\bG_{45}\bG_{34}\bG_{23}\bG_{12}\bR}
		\boxtimes & \boxtimes & \boxtimes & \boxtimes& \boxtimes \\
		\boxtimes & \boxtimes & \boxtimes & \boxtimes& \boxtimes \\
		\boxtimes & \boxtimes & \boxtimes & \boxtimes& \boxtimes \\
		\bm{\boxtimes} & \bm{\boxtimes} & \bm{\boxtimes} & \bm{\boxtimes}& \bm{\boxtimes} \\
		\bm{\boxtimes} & \bm{\boxtimes} & \bm{\boxtimes} & \bm{\boxtimes}& \bm{\boxtimes} \\
	\end{sbmatrix}.
\end{aligned}
$$
That is,  backward rotations will preserve a lot of the zeros as they are, whereas the forward rotations will eliminate these zeros.
In general,  backward rotations yield
$$
\bG_{12} \bG_{23} \ldots \bG_{(n-1),n} (\bR+\bw\bv^\top) = \bH_0  \pm \normtwo{\bw} \be_1 \bv^\top \triangleq \bH, 
$$
which is  in upper Hessenberg form. Similar to the triangularization process achieved through Givens rotations discussed in Section~\ref{section:qr-givens}, we can find a sequence of rotations $\bJ_{12}, \bJ_{23}, \ldots, \bJ_{(n-1),n}$ such that
$$
\bJ_{(n-1),n} \ldots \bJ_{23}\bJ_{12}\bH \triangleq \bR^\prime
$$
is upper triangular. 
Following the example of a $5\times 5$ matrix, the triangularization process is presented below:
$$
\begin{aligned}
	\underbrace{\bH_0  \pm \normtwo{\bw} \be_1 \bv^\top}_{\bH} =
	\footnotesize
	\begin{sbmatrix}{\bH}
		\boxtimes & \boxtimes & \boxtimes & \boxtimes & \boxtimes \\
		\boxtimes & \boxtimes & \boxtimes & \boxtimes & \boxtimes \\
		0 & \boxtimes & \boxtimes & \boxtimes& \boxtimes \\
		0 & 0 & \boxtimes & \boxtimes& \boxtimes \\
		0 & 0 & 0 & \boxtimes& \boxtimes
	\end{sbmatrix}
	&\stackrel{\bJ_{12}}{\rightarrow}
	\footnotesize
	\begin{sbmatrix}{\bJ_{12}\bH}
		\bm{\boxtimes} & \bm{\boxtimes} & \bm{\boxtimes} & \bm{\boxtimes} & \bm{\boxtimes} \\
		\bm{0} & \bm{\boxtimes} & \bm{\boxtimes} & \bm{\boxtimes} & \bm{\boxtimes} \\
		0 & \boxtimes & \boxtimes & \boxtimes& \boxtimes \\
		0 & 0 & \boxtimes & \boxtimes& \boxtimes \\
		0 & 0 & 0 & \boxtimes& \boxtimes
	\end{sbmatrix}
	\stackrel{\bJ_{23}}{\rightarrow}
	\footnotesize
	\begin{sbmatrix}{\bJ_{23}\bJ_{12}\bH}
		\boxtimes & \boxtimes & \boxtimes & \boxtimes & \boxtimes \\
		0 & \bm{\boxtimes} & \bm{\boxtimes} & \bm{\boxtimes} & \bm{\boxtimes} \\
		0 & \bm{0} & \bm{\boxtimes} & \bm{\boxtimes} & \bm{\boxtimes} \\
		0 & 0 & \boxtimes & \boxtimes& \boxtimes \\
		0 & 0 & 0 & \boxtimes& \boxtimes
	\end{sbmatrix}\\
	&\stackrel{\bJ_{34}}{\rightarrow}
	\footnotesize
	\begin{sbmatrix}{\bJ_{34}\bJ_{23}\bJ_{12}\bH}
		\boxtimes & \boxtimes & \boxtimes & \boxtimes & \boxtimes \\
		0 & \boxtimes & \boxtimes & \boxtimes & \boxtimes \\
		0 & 0 & \bm{\boxtimes} & \bm{\boxtimes} & \bm{\boxtimes} \\
		0 & 0 & \bm{0} & \bm{\boxtimes} & \bm{\boxtimes} \\
		0 & 0 & 0 & \boxtimes& \boxtimes
	\end{sbmatrix}
	\stackrel{\bJ_{45}}{\rightarrow}
	\footnotesize
	\begin{sbmatrix}{\bJ_{45}\bJ_{34}\bJ_{23}\bJ_{12}\bH}
		\boxtimes & \boxtimes & \boxtimes & \boxtimes & \boxtimes \\
		0 & 0 & \boxtimes & \boxtimes & \boxtimes \\
		0 & 0 & \boxtimes & \boxtimes& \boxtimes \\
		0 & 0 & 0 & \bm{\boxtimes} & \bm{\boxtimes} \\
		0 & 0 & 0 & \bm{0} & \bm{\boxtimes} \\
	\end{sbmatrix}.
\end{aligned}
$$
And the QR decomposition of $\bA^\prime$ is then given by 
$$
\bA^\prime = \bQ^\prime \bR^\prime,
$$
where 
\begin{equation}\label{equation:qr-rank-one-update}
\left\{
\begin{aligned}
\bR^\prime &\triangleq(\bJ_{(n-1),n} \ldots \bJ_{23}\bJ_{12}) (\bG_{12} \bG_{23} \ldots \bG_{(n-1),n}) (\bR+\bw\bv^\top);\\
\bQ^\prime &\triangleq \bQ\left\{(\bJ_{(n-1),n} \ldots \bJ_{23}\bJ_{12}) (\bG_{12} \bG_{23} \ldots \bG_{(n-1),n}) \right\}^\top; \\
\text{(or) }\bQ^{\prime\top}&\triangleq \left\{(\bJ_{(n-1),n} \ldots \bJ_{23}\bJ_{12}) (\bG_{12} \bG_{23} \ldots \bG_{(n-1),n}) \right\}\bQ^\top .
\end{aligned}
\right.
\end{equation}
The procedure is outlined in Algorithm~\ref{alg:qr-rankoneChange}.

\begin{algorithm}[htp] 
	\caption{QR Rank-One Changes} 
	\label{alg:qr-rankoneChange} 
	\begin{algorithmic}[1] 
\Require Matrix $\bA \in \real^{n\times n}$ with QR decomposition $\bA=\bQ\bR$, and $\bA^\prime = \bA+\bu\bv^\top$; 
\Statex \textbf{Stage A: Transfer $\bw$ to a first basis vector, $\bw\rightarrow \normtwo{\bw}\be_1$}
\State Calculate $\bw\leftarrow\bQ^\top\bu$; \Comment{$n(2n-1)$ flops}
\State Calculate $\bH \leftarrow \bR$;
\For{$i=n-1$ to $1$} 
\State Get Givens rotation $\bG_{i,i+1}$ with the following parameters $c, s$:
\State $c \leftarrow \frac{x_k}{\sqrt{x_k^2 + x_l^2}}$, $s\leftarrow\frac{x_l}{\sqrt{x_k^2 + x_l^2}}$, where $x_k = \bw_i$, $x_l = \bw_{i+1}$; \Comment{6 flops}
\State Calculate $\bH \leftarrow \bG_{i,i+1}\bH $ in following two steps:
\State $i$-th row: $\bH_{i,:} \leftarrow c\cdot \bH_{i,:} + s\cdot \bH_{j,:} $, where $j=i+1$; \Comment{$3(n-i+1)$ flops}
\State $(i+1)$-th row: $\bH_{i+1,:} \leftarrow -s\cdot \bH_{i,:} + c\cdot \bH_{j,:} $, where $j=i+1$; \Comment{$3(n-i+1)$ flops}
\EndFor

\Statex \textbf{Stage B: Triangularize $\bR^\prime$}
\State Set $\bR^\prime \leftarrow\bH \pm \normtwo{\bw} \be_1 \bv^\top$; \Comment{$\bH, \bR^\prime$ are both upper Hessenberg}
\For{$i=1$ to $n-1$} 
\State Get Givens rotation $\bJ_{i,i+1}$ with the following parameters $c, s$:
\State $c \leftarrow \frac{x_k}{\sqrt{x_k^2 + x_l^2}}$, $s\leftarrow\frac{x_l}{\sqrt{x_k^2 + x_l^2}}$, where $x_k = \bH_{i,i}$, $x_l = \bH_{i+1,i}$;
\State Calculate $\bR^\prime \leftarrow \bJ_{i,i+1}\bR^\prime $ in following two steps:
\State $i$-th row: $\bR^\prime_{i,:} \leftarrow c\cdot \bR^\prime_{i,:} + s\cdot \bR^\prime_{j,:} $, where $j=i+1$;
\State $(i+1)$-th row: $\bR^\prime_{i+1,:} \leftarrow -s\cdot \bR^\prime_{i,:} + c\cdot \bR^\prime_{j,:} $, where $j=i+1$;
\EndFor
\State Output $\bR^\prime$; 

\Statex \textbf{Stage C: Obtain the orthogonal matrix $\bQ^\prime$}
\State Set $\bQ^{\prime\top} = \bQ^\top$;
\For{$i=n-1$ to $1$}    \Comment{The following $c,s$ are from step 5}
\State $i$-th row: $\bQ^{\prime\top}_{i,:} \leftarrow c\cdot \bQ^{\prime\top}_{i,:} + s\cdot \bQ^{\prime\top}_{j,:} $, where $j=i+1$;  \Comment{$6n$ flops}
\State $(i+1)$-th row: $\bQ^{\prime\top}_{i+1,:} \leftarrow -s\cdot \bQ^{\prime\top}_{i,:} + c\cdot \bQ^{\prime\top}_{j,:} $, where $j=i+1$;\Comment{$6n$ flops}
\EndFor
\For{$i=1$ to $n-1$} \Comment{The following $c,s$ are from step 13}
\State $i$-th row: $\bQ^{\prime\top}_{i,:} \leftarrow c\cdot \bQ^{\prime\top}_{i,:} + s\cdot \bQ^{\prime\top}_{j,:} $, where $j=i+1$; \Comment{$6n$ flops}
\State $(i+1)$-th row: $\bQ^{\prime\top}_{i+1,:} \leftarrow -s\cdot \bQ^{\prime\top}_{i,:} + c\cdot \bQ^{\prime\top}_{j,:} $, where $j=i+1$;\Comment{$6n$ flops}
\EndFor
\State Output $\bQ^\prime$;

	\end{algorithmic} 
\end{algorithm}

We state the complexity of the rank-one update in the following theorem.
\begin{theorem}[Algorithm Complexity: QR Rank-One Change]\label{theorem:qr-full-givens-rank1}
Algorithm~\ref{alg:qr-rankoneChange} requires $\sim 8n^2$ flops to compute the full QR decomposition of an $\bA^\prime \in \real^{n\times n}$ matrix with a rank-one change to $\bA$, given that the full QR decomposition of $\bA$ is already known. Further, if $\bQ^\prime$ is needed explicitly, an additional $\sim 12n^2$ flops are required.
\end{theorem}

\begin{proof}[of Theorem~\ref{theorem:qr-full-givens-rank1}]
For step 1, it is straightforward that calculating $\bw=\bQ^\top\bu$ in step 1 requires \underline{(*). $n(2n-1)=2n^2-n$} flops. 
	
In step 5, each iteration $i$ requires $6$ flops (involving 2 square operations, 1 addition, 1 square root, and 2 divisions). And there are $n-1$ such iterations so that the complexity for all the step 5's is \underline{$6(n-1)$} flops.
	
For each iteration $i$, step 7 and step 8 involve operating on two length-$(n-i+1)$ vectors. The two steps require $6(n-i+1)$ flops for each iteration $i$ (which consist of $4(n-i+1)$ multiplications and $2(n-i+1)$ additions). 
Let $f(i) =  6(n-i+1)$, the total complexity for the two steps is equal to 
$$
\mathrm{cost} =f(1)+ f(2) +\ldots +f(n-1) = \underline{3n^2-3n} \,\, \mathrm{flops}.
$$
Therefore, the complexity for step 3 to step 9 is \underline{(*). $6(n-1) + 3n^2-3n = 3n^2+3n-6$} flops. Similarly, the complexity for step 11 to step 16 is also  \underline{(*). $3n^2+3n-6$} flops.
	
However, for each iteration $i$, step 21 and step 22 operate on two length-$n$ vectors. The two steps require  $6n$ flops for each iteration $i$, resulting in a total complexity of  \underline{(*). $6n(n-1) =6n^2-6n$} flops for step 20 to step 23. 
Similarly, step 24 to step 27 also require \underline{(*). $6n^2-6n$} flops. 
	
Therefore, the final complexity can be obtained by summing up the equations marked as (*), resulting in \underline{$20n^2$} flops if keep only the leading terms ($8n^2$ flops for calculating $\bR^\prime$, and $12n^2$ flops for calculating $\bQ^\prime$).
	
Note that in each iteration $i$, step 7 and step 8 involve operations on two length-$(n-i+1)$ vectors since $\bR$ is upper triangular. 
If we do not assume this specific structure, the overall complexity would be $26n^2$ flops, as stated in \citet{golub2013matrix}.
\end{proof}
The algorithm can be easily applied to a rectangular matrix $\bA\in \real^{m\times n}$ or to the sum  $\bA+\bU\bV^\top$, where $\bU\in \real^{m\times k}$ and $\bV\in \real^{n\times k}$.

\subsection{Application: Appending or Deleting a Column}\label{section:append-column-qr}
In certain applications, such as the $F$-test for least squares via QR decomposition \citep{lu2021rigorous}, there arises a need to remove or add a column to the observed matrix. The objective is to efficiently obtain the QR decomposition of the modified matrix.
\paragraph{Deleting a column.}
Suppose the QR decomposition of $\bA\in \real^{m\times n}$ is given by $\bA=\bQ\bR$, where the column partition of $\bA$ is $\bA=[\ba_1,\ba_2,\ldots,\ba_n]$. Now, if we delete the $k$-th column of $\bA$, resulting in $\bA^\prime \triangleq [\ba_1,\ldots,\ba_{k-1},\ba_{k+1},\ldots,\ba_n] \in \real^{m\times (n-1)}$, we want to compute the QR decomposition of $\bA^\prime$ efficiently.
Additionally, let $\bR$ have the following structure:
$$
\begin{aligned}
\begin{blockarray}{ccccc}
\begin{block}{c[ccc]c}
	        &	\bR_{11} & \ba & \bR_{12} & k-1  \\
	\bR \triangleq	&  \bzero    &	r_{kk} & \bb^\top  & 1 \\
         	&	\bzero  &\bzero& \bR_{22} & m-k  \\
\end{block}
& k-1 & 1 &  n-k & \\
\end{blockarray}, 
\end{aligned}
\qquad\text{such that}\qquad
\bQ^\top \bA^\prime = 
\begin{bmatrix}
\bR_{11} &\bR_{12} \\
\bzero & \bb^\top \\
\bzero & \bR_{22}
\end{bmatrix} \triangleq \bH
$$
is upper Hessenberg. 
An example is provided below to illustrate the case of a $6\times 5$ matrix, where $k=3$ and the $k$-th column is deleted:
$$
\begin{aligned}
	\begin{sbmatrix}{\bR = \bQ^\top\bA}
		\boxtimes & \boxtimes & \boxtimes & \boxtimes & \boxtimes \\
		0 & \boxtimes & \boxtimes & \boxtimes & \boxtimes \\
		0 & 0 & \boxtimes & \boxtimes& \boxtimes \\
		0 & 0 & 0 & \boxtimes& \boxtimes \\
		0 & 0 & 0 & 0& \boxtimes \\
		0 & 0 & 0 & 0& 0
	\end{sbmatrix}
	&\quad\implies\quad
	\begin{sbmatrix}{\bH = \bQ^\top\bA^\prime}
		\boxtimes & \boxtimes  & \boxtimes & \boxtimes \\
		0 & \boxtimes  & \boxtimes & \boxtimes \\
		0 & 0 & \boxtimes& \boxtimes \\
		0 & 0  & \boxtimes& \boxtimes \\
		0 & 0  & 0& \boxtimes \\
		0 & 0  & 0& 0
	\end{sbmatrix}.
\end{aligned}
$$

Again, for columns $k$ to $n-1$ of $\bH$, there exists a set of rotations $\bG_{k,k+1}$, $\bG_{k+1,k+2}$, $\ldots$, $\bG_{n-1,n}$ that can introduce zeros in the elements $h_{k+1,k}$, $h_{k+2,k+1}$, $\ldots$, $h_{n,n-1}$ of $\bH$. Then the triangular matrix $\bR^\prime$ is given by 
$$
\bR^\prime \triangleq \bG_{n-1,n}\ldots \bG_{k+1,k+2}\bG_{k,k+1}\bQ^\top \bA^\prime.
$$
The orthogonal matrix can be obtained through the following procedure:
\begin{equation}\label{equation:qr-delete-column-finalq}
\bQ^\prime = (\bG_{n-1,n}\ldots \bG_{k+1,k+2}\bG_{k,k+1}\bQ^\top )^\top = \bQ \bG_{k,k+1}^\top  \bG_{k+1,k+2}^\top \ldots \bG_{n-1,n}^\top,
\end{equation}
such that $\bA^\prime = \bQ^\prime\bR^\prime$. 
The procedure is outlined in Algorithm~\ref{alg:qr-delete-a-column}. And the $6\times 5$ example is shown as follows, where $\boxtimes$ represents a value that is not necessarily zero, and \textbf{boldface} indicates the value has just been changed:
$$
\begin{aligned}
\begin{sbmatrix}{\bR = \bQ^\top\bA}
\boxtimes & \boxtimes & \boxtimes & \boxtimes & \boxtimes \\
0 & \boxtimes & \boxtimes & \boxtimes & \boxtimes \\
0 & 0 & \boxtimes & \boxtimes& \boxtimes \\
0 & 0 & 0 & \boxtimes& \boxtimes \\
0 & 0 & 0 & 0& \boxtimes \\
0 & 0 & 0 & 0& 0
\end{sbmatrix}
&\stackrel{k=3}{\rightarrow}
\begin{sbmatrix}{\bH = \bQ^\top\bA^\prime}
\boxtimes & \boxtimes  & \boxtimes & \boxtimes \\
0 & \boxtimes  & \boxtimes & \boxtimes \\
0 & 0 & \boxtimes& \boxtimes \\
0 & 0  & \boxtimes& \boxtimes \\
0 & 0  & 0& \boxtimes \\
0 & 0  & 0& 0
\end{sbmatrix}
\stackrel{\bG_{34}}{\rightarrow}
\begin{sbmatrix}{\bG_{34}\bH }
	\boxtimes & \boxtimes  & \boxtimes & \boxtimes \\
	0 & \boxtimes  & \boxtimes & \boxtimes \\
	0 & 0 & \bm{\boxtimes}& \bm{\boxtimes} \\
	0 & 0  & \bm{0}& \bm{\boxtimes} \\
	0 & 0  & 0& \boxtimes \\
	0 & 0  & 0& 0
\end{sbmatrix}
\stackrel{\bG_{45}}{\rightarrow}
\begin{sbmatrix}{\bG_{45}\bG_{34}\bH}
	\boxtimes & \boxtimes  & \boxtimes & \boxtimes \\
	0 & \boxtimes  & \boxtimes & \boxtimes \\
	0 & 0 & \boxtimes & \boxtimes \\
	0 & 0  &0 & \bm{\boxtimes} \\
	0 & 0  & 0& \bm{0} \\
	0 & 0  & 0& 0
\end{sbmatrix}.
\end{aligned}
$$

\begin{algorithm}[h] 
\caption{QR Deleting a Column} 
\label{alg:qr-delete-a-column} 
\begin{algorithmic}[1] 
\Require Matrix $\bA \in \real^{m\times n}$ with full QR decomposition $\bA=\bQ\bR$, and $\bA^\prime \in \real^{m\times (n-1)}$ by deleting column $k$ of $\bA$; 
\Statex \textbf{Stage A: Triangularize $\bH$}
\State Obtain $\bH$ by deleting column $k$ of $\bR$, that is, $\bH=\bQ^\top\bA^\prime$;
\For{$i=k$ to $n-1$} 
\State Get Givens rotation $\bG_{i,i+1}$ with the following parameters $c, s$:
\State $c \leftarrow \frac{x_k}{\sqrt{x_k^2 + x_l^2}}$, $s\leftarrow\frac{x_l}{\sqrt{x_k^2 + x_l^2}}$, where $x_k = h_{ii}$, $x_l = h_{i+1,i}$;
\State Calculate $\bH \leftarrow \bG_{i,i+1}\bH $ in following two steps:
\State $i$-th row: $\bH_{i,:} \leftarrow c\cdot \bH_{i,:} + s \bH_{j,:} $, where $j=i+1$;
\State $(i+1)$-th row: $\bH_{i+1,:} \leftarrow -s\cdot \bH_{i,:} + c \bH_{j,:} $, where $j=i+1$;
\EndFor
\State Set $\bR^\prime \leftarrow \bH$ and output $\bR^\prime$;
\Statex \textbf{Stage B: Obtain the orthogonal matrix $\bQ^\prime$}
\State Set $\bQ^\prime \leftarrow \bQ^\top $;
\For{$i=k$ to $n-1$} 
\State $c \leftarrow \frac{x_k}{\sqrt{x_k^2 + x_l^2}}$, $s\leftarrow\frac{x_l}{\sqrt{x_k^2 + x_l^2}}$, where $x_k$, $x_l$ are from step 4;
\State Calculate $\bQ^\prime \leftarrow \bG_{i,i+1}\bQ^\prime $ in following two steps:
\State $i$-th row: $\bQ^\prime_{i,:} \leftarrow c\cdot \bQ^\prime_{i,:} + s \bQ^\prime_{j,:} $, where $j=i+1$;
\State $(i+1)$-th row: $\bQ^\prime_{i+1,:} \leftarrow -s\cdot \bQ^\prime_{i,:} + c \bQ^\prime_{j,:} $, where $j=i+1$;
\EndFor
\State Output $\bQ^\prime \leftarrow \bQ^{\prime\top}$ from Equation~\eqref{equation:qr-delete-column-finalq}; 
\end{algorithmic} 
\end{algorithm}

\begin{theorem}[Algorithm Complexity: QR Deleting Column]\label{theorem:qr-full-givens-delete-column}
Algorithm~\ref{alg:qr-delete-a-column}  requires $\sim 3n^2-6nk+3k^2$ flops to compute a full QR decomposition of matrix $\bA^\prime \in \real^{m\times (n-1)}$.
This matrix is obtained by deleting column $k$ from $\bA\in \real^{m\times n}$, assuming the full QR decomposition of $\bA$ is already known. 
Further, if $\bQ^\prime$ is needed explicitly, an additional $\sim 6m(n-k)$ flops are required.
\end{theorem}

\begin{proof}[of Theorem~\ref{theorem:qr-full-givens-delete-column}]
For step 4, each iteration $i$ requires $6$ flops (consisting of 2 square operations, 1 addition, 1 square root, and 2 divisions). 
And there are $n-k$ iterations in total, resulting in a complexity of \underline{$6(n-k)$} flops for all the step 4 operations.
	
For each iteration $i$, step 6 and step 7 operate on two length-$(n-i)$ vectors since $\bH$ is upper Hessenberg. The two steps require  $6(n-i)$ flops for each iteration $i$ (involving  $4(n-i)$ multiplications and $2(n-i)$ additions).  
Let $f(i) =  6(n-i)$, the total complexity for the two steps is given by
$$
\mathrm{cost} =f(k)+ f(k+1) +\ldots + f(n-1) = \underline{3n^2-6nk+3k^2+3n-3k} \,\, \mathrm{flops}.
$$
Therefore, the complexity for step 2 to step 8 is \underline{(*). $3n^2-6nk+3k^2$} flops if we keep only the leading terms.

However, for each iteration $i$, step 14 and step 15 involve operations on two length-$m$ vectors. The two steps require $6m$ flops for each $i$, and  there are $n-k$ such iterations. Therefore, the total complexity of steps 11 to 16 is \underline{(*). $6m(n-k)$} flops.  
\end{proof}
Note that the value of column $k$ affects the complexity: when $k=n$, the complexity is 0; and when $k=1$, the complexity reaches its maximal value.

\paragraph{Appending a column.}
Similarly, suppose $\widetilde{\bA} \triangleq [\ba_1,\ba_k,\bw,\ba_{k+1},\ldots,\ba_n]$, where we append a vector $\bw$ to the $(k+1)$-th column of $\bA$. Applying the orthogonal transformation $\bQ^\top$, we have
$$
\bQ^\top \widetilde{\bA} = [\bQ^\top\ba_1,\ldots, \bQ^\top\ba_k, \bQ^\top\bw, \bQ^\top\ba_{k+1}, \ldots,\bQ^\top\ba_n] \triangleq \widetilde{\bH}.
$$
A set of Givens rotations $\bJ_{m-1,m}, \bJ_{m-2,m-1}, \ldots, \bJ_{k+1,k+2}$ can introduce zeros in the elements $\widetilde{h}_{m,k+1}$, $\widetilde{h}_{m-1,k+1}$, $\ldots$, $\widetilde{h}_{k+2,k+1}$ of $\widetilde{\bH}$, thereby achieving the desired result of the updated QR decomposition. That is,
$$
\widetilde{\bR} \triangleq\bJ_{k+1,k+2}\ldots \bJ_{m-2,m-1}\bJ_{m-1,m} \bQ^\top \widetilde{\bA},
$$
is upper triangular.
Suppose $\widetilde{\bH}$ is a $6\times 5$ matrix and $k=2$. An example is shown below, where $\boxtimes$ represents a value that is not necessarily zero, and \textbf{boldface} indicates the value has just been changed:
$$
\footnotesize
\begin{aligned}
\begin{sbmatrix}{\widetilde{\bH}}
\boxtimes & \boxtimes & \boxtimes & \boxtimes & \boxtimes \\
0 & \boxtimes & \boxtimes & \boxtimes & \boxtimes \\
0 & 0 & \boxtimes & \boxtimes& \boxtimes \\
0 & 0 & \boxtimes & 0& \boxtimes \\
0 & 0 & \boxtimes & 0& 0 \\
0 & 0 & \boxtimes & 0& 0
\end{sbmatrix}
&\stackrel{\bJ_{56}}{\rightarrow}
\begin{sbmatrix}{\bJ_{56}\widetilde{\bH} \rightarrow \widetilde{h}_{63}=0}
\boxtimes & \boxtimes & \boxtimes & \boxtimes & \boxtimes \\
0 & \boxtimes & \boxtimes & \boxtimes & \boxtimes \\
0 & 0 & \boxtimes & \boxtimes& \boxtimes \\
0 & 0 & \boxtimes & 0& \boxtimes \\
0 & 0 & \bm{\boxtimes} & 0& 0 \\
0 & 0 & \bm{0} & 0& 0
\end{sbmatrix}
\stackrel{\bJ_{45}}{\rightarrow}
\begin{sbmatrix}{\bJ_{45}\bJ_{56}\widetilde{\bH}\rightarrow \widetilde{h}_{53}=0}
\boxtimes & \boxtimes & \boxtimes & \boxtimes & \boxtimes \\
0 & \boxtimes & \boxtimes & \boxtimes & \boxtimes \\
0 & 0 & \boxtimes & \boxtimes& \boxtimes \\
0 & 0 & \bm{\boxtimes} & 0& \bm{\boxtimes} \\
0 & 0 & \bm{0} & 0& \bm{\boxtimes}\\
0 & 0 & 0 & 0& 0
\end{sbmatrix}
\stackrel{\bJ_{34}}{\rightarrow}
\begin{sbmatrix}{\bJ_{34}\bJ_{45}\bJ_{56}\widetilde{\bH}\rightarrow \widetilde{h}_{43}=0}
\boxtimes & \boxtimes & \boxtimes & \boxtimes & \boxtimes \\
0 & \boxtimes & \boxtimes & \boxtimes & \boxtimes \\
0 & 0 & \bm{\boxtimes}  & \bm{\boxtimes} & \bm{\boxtimes}  \\
0 & 0 & \bm{0}  & \bm{\boxtimes} & \bm{\boxtimes}  \\
0 & 0 & 0 & 0& \boxtimes\\
0 & 0 & 0 & 0& 0
\end{sbmatrix} = \widetilde{\bR}.
\end{aligned}
$$ 
Finally, we obtain the orthogonal matrix 
\begin{equation}\label{equation:qr-add-column-finalq}
\widetilde{\bQ} = (\bJ_{k+1,k+2}\ldots \bJ_{m-2,m-1}\bJ_{m-1,m} \bQ^\top )^\top = \bQ \bJ_{m-1,m}^\top  \bJ_{m-2,m-1}^\top \ldots \bJ_{k+1,k+2}^\top,
\end{equation}
such that $\widetilde{\bA} = \widetilde{\bQ}\widetilde{\bR}$. The procedure is formulated in Algorithm~\ref{alg:qr-adding-a-column}.

\begin{algorithm}[h] 
\caption{QR Adding a Column} 
\label{alg:qr-adding-a-column} 
\begin{algorithmic}[1] 
\Require Matrix $\bA \in \real^{m\times n}$ with full QR decomposition $\bA=\bQ\bR$, and $\widetilde{\bA}\in \real^{m\times (n+1)}$ by adding column $\bw$ into $(k+1)$-th column of $\bA$; 
\Statex \textbf{Stage A: Triangularize $\widetilde{\bH}$}
\State Calculate $\bQ^\top\bw$;
\State Obtain $\widetilde{\bH}$ by inserting $\bQ^\top \bw$ into $(k+1)$-th column of $\bR$;
\For{$i=m-1$ to $k+1$} 
\State Get Givens rotation $\bJ_{i,i+1}$ with the following parameters $c, s$:
\State $c \leftarrow \frac{x_k}{\sqrt{x_k^2 + x_l^2}}$, $s\leftarrow\frac{x_l}{\sqrt{x_k^2 + x_l^2}}$, where $x_k = \widetilde{h}_{i,k+1}$, $x_l = \widetilde{h}_{i+1,k+1}$;
\State Calculate $\widetilde{\bH} \leftarrow \bJ_{i,i+1}\bH $ in following two steps:
\State $i$-th row: $\widetilde{\bH}_{i,:} \leftarrow c\cdot \widetilde{\bH}_{i,:} + s \widetilde{\bH}_{j,:} $, where $j=i+1$;
\State $(i+1)$-th row: $\widetilde{\bH}_{i+1,:} \leftarrow -s\cdot \widetilde{\bH}_{i,:} + c \widetilde{\bH}_{j,:} $, where $j=i+1$;
\EndFor
\State Set $\widetilde{\bR} \leftarrow \widetilde{\bH}$ and output $\widetilde{\bR}$;
\Statex \textbf{Stage B: Obtain the orthogonal matrix $\widetilde{\bQ}$}
\State Set $\widetilde{\bQ} \leftarrow \bQ^\top $;
\For{$i=m-1$ to $k+1$} 
\State $c \leftarrow \frac{x_k}{\sqrt{x_k^2 + x_l^2}}$, $s\leftarrow\frac{x_l}{\sqrt{x_k^2 + x_l^2}}$, where $x_k$, $x_l$ are from step 5;
\State Calculate $\widetilde{\bQ} \leftarrow \bJ_{i,i+1}\widetilde{\bQ} $ in following two steps:
\State $i$-th row: $\widetilde{\bQ}_{i,:} \leftarrow c\cdot \widetilde{\bQ}_{i,:} + s \widetilde{\bQ}_{j,:} $, where $j=i+1$;
\State $(i+1)$-th row: $\bQ_{i+1,:} \leftarrow -s\cdot \widetilde{\bQ}_{i,:} + c \widetilde{\bQ}_{j,:} $, where $j=i+1$;
\EndFor
\State Output $\widetilde{\bQ} \leftarrow \widetilde{\bQ}^{\top}$ from Equation~\eqref{equation:qr-add-column-finalq}; 
\end{algorithmic} 
\end{algorithm}
\begin{theorem}[Algorithm Complexity: QR Adding Column]\label{theorem:qr-full-givens-add-column}
Algorithm~\ref{alg:qr-adding-a-column} requires $\sim 2m^2+6(mn+k^2-nk-mk)$ flops to compute a full QR decomposition of matrix $\widetilde{\bA} \in \real^{m\times (n+1)}$, where we add a column to  the $(k+1)$-th column of $\bA\in \real^{m\times n}$ and the full QR decomposition of $\bA$ is known. Further, if $\widetilde{\bQ}$ is needed explicitly, an additional $\sim 6m(m-k)$ flops are required.
\end{theorem}

\begin{proof}[of Theorem~\ref{theorem:qr-full-givens-add-column}]
It is evident that step 1 requires \underline{$m(2m-1)$} flops to calculate $\bQ^\top\bw$.

In step 5, each iteration $i$ requires $6$ flops (including 2 square operations, 1 addition, 1 square root, and 2 divisions). And there are $m-k-1$ such iterations, resulting in a complexity of \underline{$6(m-k-1)$} flops for all step 5 operations.

For each iteration $i$, steps 7 and 8 operate on two vectors of length $(n-k+1)$ because  $\widetilde{\bH}$ is upper Hessenberg. The two steps require $6(n-k+1)$ flops for each iteration $i$ (consisting of $4(n-k+1)$ multiplications and $2(n-k+1)$ additions).  
Since there are $m-k-1$ such iterations, the overall computational complexity for all instances of steps 7 and 8 is \underline{$6(n-k+1)(m-k-1)$} flops.

Therefore, the complexity for steps 1 to 9 is given by
$$
\begin{aligned}
&\gap m(2m-1)+6(m-k-1)+6(n-k+1)(m-k-1)
=\underline{2m^2+6(n-k+2)(m-k-1)}, 
\end{aligned}
$$ 
or \underline{$2m^2+6(mn+k^2-nk-mk)$} flops if we keep only the leading terms.

However, for each iteration $i$, step 15 and step 16 involve the manipulation of two vectors of length $m$. The two steps require $6m$ flops for each iteration $i$. Since there are a total of $(m-k-1)$ such iterations, the overall complexity for steps 12 to 17 amounts to \underline{$6m(m-k-1)$} flops. 
Alternatively, if we consider only the leading terms, the complexity can be simplified to \underline{$6m(m-k)$} flops. 
\end{proof}
Note that the column number $k$ plays a significant role in determining the complexity. 
When $k=n$, the complexity is $2m^2$ flops. On the other hand, when $k=1$, the complexity reaches its maximum value.

\paragraph{Real-world application.} 
The method described above is valuable for efficient variable selection in the least squares problem using QR decomposition. 
In this approach, we iteratively eliminate a column from the data matrix $\bA$ and perform an $F$-test to determine the significance of the corresponding variable. If the variable is deemed insignificant, it is removed, leading to a simpler model. A brief overview is provided below, and for further information, please refer to \citet{lu2021rigorous}. 

Following the setup described in Section~\ref{section:application-ls-qr}, 
let's consider the overdetermined system $\bA\bx = \bb$, where $\bA\in \real^{m\times n}$ represents the data matrix, and $\bb\in \real^m$ with $m\geq n$ is the observation matrix. 
The LS solution is obtained by minimizing $\normtwo{\bA\bx-\bb}^2$, and it can be expressed as $\bx_{LS} = (\bA^\top\bA)^{-1}\bA^\top\bb$.
In the derivation above, since $\bA$ has full column rank and $\rank(\bA^\top\bA)=\rank(\bA)$, $\bA^\top\bA$ is invertible.

Suppose we remove a column from $\bA$ to obtain $\widehat{\bA}$. Consequently, the LS solution changes from $\bx_{LS}$ to $\widehat{\bx}_{LS}$.
Define 
$$
\begin{aligned}
	RSS(\widehat{\bx}_{LS}) &= \normtwo{\bb - \widehat{\bb}_{LS}}^2, \qquad \text{where~ } \widehat{\bb}_{LS} = \widehat{\bA}\widehat{\bx}_{LS},   \\
	RSS(\bx_{LS})&= \normtwo{\bb - \bb_{LS}}^2, \qquad \text{where~ } \bb_{LS} = \bA\bx_{LS},\\
	\bH & = \bA(\bA^\top\bA)^{-1}\bA^\top, \\
	\widehat{\bH} &= \widehat{\bA}(\widehat{\bA}^\top \widehat{\bA})^{-1}\widehat{\bA}^\top. \\
\end{aligned}
$$
Suppose the \textit{reduced} QR decompositions of $\bA$ and $\widehat{\bA}$ are given by $\bA=\bQ\bR$ and $\widehat{\bA}=\widehat{\bQ}\widehat{\bR}$.
Thus $RSS(\bx_{LS}) = \bb^\top (\bI-\bH)\bb = \bb^\top\bb - (\bb^\top\bQ)(\bQ^\top\bb)$
and $RSS(\widehat{\bx}_{LS}) - RSS(\bx_{LS}) = \normtwo{\bb_{LS} - \widehat{\bb}_{LS}}^2 = \bb^\top(\bH-\widehat{\bH})\bb=(\bb^\top\bQ)(\bQ^\top\bb)-(\bb^\top\widehat{\bQ})(\widehat{\bQ}^\top\bb)$, which are the differences of two inner products.
It can be shown that $RSS(\bx_{LS})\sim \sigma^2 \chi^2_{(m-n)}$, which follows a Chi-square distribution, and $\sigma$ is the noise level. Under the hypothesis that the deleted column is not significant, we could conclude that 
$$
T={\frac{1}{n-q}\big(RSS(\widehat{\bx}_{LS}) - RSS(\bx_{LS})\big)  }/{\frac{1}{m-n}RSS(\bx_{LS})} \sim F_{n-q,m-n},
$$
which is the \textbf{test statistic for the $F$-test} with $q=n-1$.
Suppose we have the data set $(\ba_1, b_1)$, $(\ba_2, b_2)$, $\ldots$, $(\ba_n, b_m)$, and we observe $T=t$ for this particular data set. Then 
$$
p=P[T((\ba_1, b_1), (\ba_2, b_2), \ldots, (\ba_n, b_m)) \geq t] = P[F_{n-q,m-n} \geq t].
$$
We reject the hypothesis if $p<\alpha$, for some small $\alpha$, say 0.05. This is known as the \textit{p-value}.

\subsection{Application: Appending or Deleting a Row}\label{section:qr_adddelrow}
Similarly, we may also need to append or delete a row from the observed matrix. The goal now is to determine an efficient method to obtain the QR decomposition for this modified matrix.

\paragraph{Appending a row.}
Suppose the full QR decomposition of $\bA\in \real^{m\times n}$ is given by $\bA= \scriptsize\begin{bmatrix}
	\bA_1 \\
	\bA_2 
\end{bmatrix}=\bQ\bR$, where $\bA_1\in \real^{k\times n}$ and $\bA_2 \in \real^{(m-k)\times n}$. Now, if we append a row such that $\bA^\prime = \scriptsize\begin{bmatrix}
	\bA_1 \\
	\bw^\top\\
	\bA_2 
\end{bmatrix} \in \real^{(m+1)\times n}$, we aim to efficiently obtain the full QR decomposition of $\bA^\prime$. Construct a permutation matrix, denoted by $\bP$, such that
$$
\bP=
\begin{bmatrix}
	\bzero & 1 & \bzero  \\
	\bI_k & \bzero & \bzero \\
	\bzero & \bzero & \bI_{m-k}
\end{bmatrix}
\longrightarrow
\bP 
\begin{bmatrix}
	\bA_1 \\
	\bw^\top\\
	\bA_2 
\end{bmatrix}
=
\begin{bmatrix}
	\bw^\top\\
	\bA_1 \\
	\bA_2 
\end{bmatrix}
\quad 
\implies
\quad
\begin{bmatrix}
	1 & \bzero \\
	\bzero & \bQ^\top 
\end{bmatrix} 
\bP 
\bA^\prime 
=
\begin{bmatrix}
	\bw^\top \\
	\bR 
\end{bmatrix}=\bH 
$$
is upper Hessenberg. Similarly, a set of rotations $\bG_{12}, \bG_{23}, \ldots, \bG_{n,n+1}$ can be applied to introduce zeros in the elements $h_{21}$, $h_{32}$, $\ldots$, $h_{n+1,n}$ of $\bH$. The resulting triangular matrix $\bR^\prime$ is given by 
$$
\bR^\prime \triangleq \bG_{n,n+1}\ldots \bG_{23}\bG_{12}\begin{bmatrix}
	1 & \bzero \\
	\bzero & \bQ^\top 
\end{bmatrix} 
\bP   \bA^\prime.
$$
And the orthogonal matrix  is obtained by
$$
\bQ^\prime \triangleq \left(\bG_{n,n+1}\ldots \bG_{23}\bG_{12}\begin{bmatrix}
	1 & \bzero \\
	\bzero & \bQ^\top 
\end{bmatrix} 
\bP \right)^\top 
= 
\bP^\top 
\begin{bmatrix}
	1 & \bzero \\
	\bzero & \bQ 
\end{bmatrix}  
\bG_{12}^\top  \bG_{23}^\top \ldots \bG_{n,n+1}^\top.
$$
Thus, the QR decomposition of $\bA^\prime$ is  $\bA^\prime = \bQ^\prime\bR^\prime$.

\paragraph{Deleting a row.} Suppose $\bA = 
\scriptsize
\begin{bmatrix}
\bA_1 \\
\bw^\top\\
\bA_2 
\end{bmatrix} \in \real^{m\times n}$, where $\bA_1\in\real^{k\times n}$ and $\bA_2 \in \real^{(m-k-1)\times n}$. The full QR decomposition of $\bA$ is given by $\bA=\bQ\bR$, where $\bQ\in \real^{m\times m}$ and $\bR\in \real^{m\times n}$. We want to compute the full QR decomposition of 
$\widetilde{\bA} =  
\footnotesize
\begin{bmatrix}
\bA_1 \\
\bA_2 
\end{bmatrix}$ efficiently (assuming $m-1\geq n$). Analogously, we can construct a permutation matrix 
$$
\bP = 
\begin{bmatrix}
	\bzero & 1 & \bzero  \\
	\bI_k & \bzero & \bzero \\
	\bzero & \bzero & \bI_{m-k-1}
\end{bmatrix}
\implies 
\bP\bA = 
\begin{bmatrix}
	\bzero & 1 & \bzero  \\
	\bI_k & \bzero & \bzero \\
	\bzero & \bzero & \bI_{m-k-1}
\end{bmatrix}
\begin{bmatrix}
	\bA_1 \\
	\bw^\top\\
	\bA_2 
\end{bmatrix}=
\begin{bmatrix}
	\bw^\top \\
	\bA_1\\
	\bA_2 
\end{bmatrix} = \bP\bQ\bR \triangleq\bM\bR ,
$$
where $\bM \triangleq \bP\bQ$ is an orthogonal matrix. 
Let $\bmm^\top$ be the first row of $\bM$, and construct a set of Givens rotations $\bG_{m-1,m}, \bG_{m-2,m-1}, \ldots, \bG_{1,2}$, which introduce zeros in the elements $m_m, m_{m-1}, \ldots, m_2$ of $\bmm$, respectively. 
By applying these rotations, we can obtain $\bG_{1,2}\ldots \bG_{m-2,m-1}\bG_{m-1,m}\bmm = \alpha \be_1$, where $\alpha = \pm 1$. Therefore, we obtain the following result:
$$
\bG_{1,2}\ldots \bG_{m-2,m-1}\bG_{m-1,m} \bR \triangleq
\begin{blockarray}{cc}
	\begin{block}{[c]c}
		\bv^\top & 1 \\
		\bR_1& m-1  \\
	\end{block}
\end{blockarray},
$$ 
which is upper Hessenberg with $\bR_1\in \real^{(m-1)\times n}$ being upper triangular. And 
$$
\bM   \bG_{m-1,m}^\top \bG_{m-2,m-1}^\top \ldots  \bG_{1,2}^\top  \triangleq
\begin{bmatrix}
	\alpha & \bzero \\
	\bzero & \bQ_1 
\end{bmatrix},
$$
where $\bQ_1\in \real^{(m-1)\times (m-1)}$ is an orthogonal matrix. The bottom-left block of the  matrix above is a zero vector since $\alpha=\pm 1$ and $\bM$ is orthogonal. To see this, let $\bG\triangleq\bG_{m-1,m}^\top \bG_{m-2,m-1}^\top $ $\ldots  \bG_{1,2}^\top $, with its first  column denoted as  $\bg$. Writing $\bM$ as the row partition $\bM = [\bmm^\top; \bmm_2^\top; \bmm_3^\top; \ldots, \bmm_{m}^\top]$,  we have
$$
\begin{aligned}
	\bmm^\top\bg &= \pm 1 \qquad  \rightarrow \qquad  \bg = \pm \bmm, \\
	\bmm_i^\top \bmm &=0,  \qquad \forall i \in \{2,3,\ldots,m\}.
\end{aligned}
$$
This results in 
$$
\begin{aligned}
\bP\bA&=\bM\bR
=(\bM   \bG_{m-1,m}^\top \bG_{m-2,m-1}^\top \ldots  \bG_{1,2}\top ) (\bG_{1,2}\ldots \bG_{m-2,m-1}\bG_{m-1,m} \bR ) \\
&=
\begin{bmatrix}
\alpha & \bzero \\
\bzero & \bQ_1 
\end{bmatrix}
\begin{bmatrix}
\bv^\top \\
\bR_1 
\end{bmatrix} = 
\begin{bmatrix}
\alpha \bv^\top \\
\bQ_1\bR_1 
\end{bmatrix}
=
\begin{bmatrix}
\bw^\top \\
\widetilde{\bA}
\end{bmatrix}
.
\end{aligned}
$$
This implies that $\bQ_1\bR_1$ is the full QR decomposition of $\widetilde{\bA}=\footnotesize\begin{bmatrix}
\bA_1\\
\bA_2 
\end{bmatrix}$.

\index{Condition number}
\index{Ill-condition}
\index{Well-condition}
\subsection{Application: Reducing the Ill-Condition via the QR decomposition}\label{section:qr_condition}

\index{Linear system}
\paragraph{Well-determined linear system.}
Consider the well-determined linear equation $\bA\bx = \bb$, where $\bA\in \real^{n\times n}$ is nonsingular. The solution $\bx = \bA^{-1}\bb$ exists and is unique. Now suppose the vector $\bb$ is perturbed by $\delta\bb$. 
In this case, the solution is given by:
$$
\bx + \delta\bx = \bA^{-1}(\bb+\delta\bb)=\bA^{-1}\bb+\bA^{-1}\delta\bb = \bx + \bA^{-1}\delta\bb.
$$
That is, $\delta\bx = \bA^{-1}\delta\bb $. 
Utilizing the properties of the matrix-vector inequality of the spectral matrix norm (see Definition~\ref{definition:spectral_norm} or  Appendix~\ref{appendix:matrix-norm-sect2}), we obtain 
\begin{equation}\label{equation:qr_well_cond}
\normtwo{\delta\bx} = \normtwo{\bA^{-1}\delta\bb} \leq \normtwo{\bA^{-1} }  \normtwo{\delta\bb}.
\end{equation}
This is referred to as the \textit{absolute error bound} of $\normtwo{\delta\bx}$. 
Furthermore, when $\normtwo{\bA^{-1} }$ is small, small perturbations in $\bb$ (i.e., $\normtwo{\delta\bb}$ is small) will result in small changes in $\bx$. 
However, if $\normtwo{\bA^{-1}}$ is large, the variations in $\bx$ may be substantial.

Now, we divide the  equation above by $\normtwo{\bx}$:
$$
\frac{\normtwo{\delta\bx} }{\normtwo{\bx}} \leq \normtwo{\bA^{-1}}  \frac{\normtwo{\delta\bb}}{\normtwo{\bx}}.
$$
By $ \normtwo{\bb} = \normtwo{\bA\bx} \rightarrow \normtwo{\bx} \geq \frac{\normtwo{\bb}}{\normtwo{\bA}}$, we can conclude:
$$
\frac{\normtwo{\delta\bx} }{\normtwo{\bx}} 
\leq \normtwo{\bA} \normtwo{\bA^{-1}}  \frac{\normtwo{\delta\bb}}{\normtwo{\bb}}.
$$
This is known as the \textit{relative error bound} of $\frac{\normtwo{\delta\bx} }{\normtwo{\bx}}$. The product $\normtwo{\bA} \normtwo{\bA^{-1}} $ is called the \textit{condition number} of $\bA$ and is denoted by $\kappa(\bA)$:
\begin{equation}\label{equation:qr_condition_num}
\kappa(\bA) = \normtwo{\bA} \normtwo{\bA^{-1}}.
\end{equation}
Similarly, utilizing the inequality of the Frobenius norm $\normtwo{\bb} = \normtwo{\bA\bx} \rightarrow \normtwo{\bb} \leq \normf{\bA} \normtwo{\bx}$, we can also define the condition number as
\begin{equation}
\kappa(\bA)= \normf{\bA} \normf{\bA^{-1}}.
\end{equation}
In this context, we will focus on the first definition of the condition number. 
When the relative error in $\bx$ is not significantly greater than the relative error in $\bb$, we classify the matrix as \textit{well-conditioned}. That is, if the condition number $\kappa(\bA)$ is small, the matrix $\bA$ is considered well-conditioned. 
Conversely, if the condition number is large, the matrix is characterized as \textit{ill-conditioned}.

\index{Least squares}
\paragraph{Overdetermined linear system.}
Now, let's delve into the overdetermined linear equation $\bA\bx = \bb$, where $\bA\in \real^{m\times n}$ with $m\geq n$. Assuming $\bA$ has full column rank such that $\bA^\top\bA$ is invertible,~\footnote{In Section~\ref{section:ls-utv} and~\ref{section:application-ls-svd}, we briefly discuss how to use UTV decomposition and SVD to address the rank-deficient least squares problems.
	See Problem~\ref{prob:als_pseudo1}$\sim$\ref{prob:als_pseudon} for more results.} the unique least squares solution is given by 
$
\bx = (\bA^\top\bA)^{-1}\bA^\top\bb, 
$
which results from the normal equation (see Definition~\ref{definition:normal-equation-als}):
$
(\bA^\top\bA)\bx = \bA^\top\bb.
$
It can be shown that the condition number is given by 
$$
\kappa(\bA^\top\bA) = \kappa(\bA)^2.
$$

\begin{example}[Reducing Ill-Condition]
Consider the matrix
$$ 
\bA = 
\begin{bmatrix}
1+\delta & 1\\
1 & 1+\delta
\end{bmatrix},
$$
where $\delta$ is small. The condition number of $\bA$ is of the order $\delta^{-1}$. If we use the QR decomposition $\bA = \bQ\bR$ to solve the linear equation, then 
$$
\kappa(\bQ)=1, \qquad 
\kappa(\bA) \rightarrow \kappa(\bQ^\top\bA) = \kappa(\bR).
$$
If one believes/proves that the condition number of $\bR$ is smaller than $\bA$, then we can overcome the ill-conditioning issue to some extent. Specifically, four QR decomposition results of $\bA$ are given by (suppose the diagonal of the matrix $\bR$ can be negative) 
$$
\begin{aligned}
\bA &= \bQ_1\bR_1
= 
\begin{bmatrix}
q_{11} & q_{12} \\
q_{21} & q_{22}
\end{bmatrix}
\begin{bmatrix}
r_{11} & r_{12} \\
0 & r_{22}
\end{bmatrix}\\
&=
\begin{bmatrix}
\frac{1+\delta}{\sqrt{\delta^2+2\delta+2}} & \frac{1}{\sqrt{\delta^2+2\delta+2}}\\
\frac{1}{\sqrt{\delta^2+2\delta+2}} & -\frac{1+\delta}{\sqrt{\delta^2+2\delta+2}}
\end{bmatrix}
\begin{bmatrix}
\sqrt{\delta^2+2\delta+2} & \frac{\delta^2+\delta}{(1+\delta)\sqrt{\delta^2+2\delta+2}} + \frac{\sqrt{\delta^2+2\delta+2}}{1+\delta}\\
0 & - \frac{\delta^2+2\delta}{\sqrt{\delta^2+2\delta+2}}
\end{bmatrix},
\end{aligned}
$$
or 
$$
\begin{aligned}
\bA &= \bQ_2\bR_2
= 
\begin{bmatrix}
q_{11} & \textcolor{mylightbluetext}{-}q_{12} \\
q_{21} & \textcolor{mylightbluetext}{-}q_{22}
\end{bmatrix}
\begin{bmatrix}
r_{11} & r_{12} \\
0 & \textcolor{mylightbluetext}{-}r_{22}
\end{bmatrix}\\
&= 
\begin{bmatrix}
\frac{1+\delta}{\sqrt{\delta^2+2\delta+2}} & -\frac{1}{\sqrt{\delta^2+2\delta+2}}\\
\frac{1}{\sqrt{\delta^2+2\delta+2}} & \frac{1+\delta}{\sqrt{\delta^2+2\delta+2}}
\end{bmatrix}
\begin{bmatrix}
\sqrt{\delta^2+2\delta+2} & \frac{\delta^2+\delta}{(1+\delta)\sqrt{\delta^2+2\delta+2}} + \frac{\sqrt{\delta^2+2\delta+2}}{1+\delta}\\
0 & \frac{\delta^2+2\delta}{\sqrt{\delta^2+2\delta+2}}
\end{bmatrix},
\end{aligned}
$$
or 
$$
\begin{aligned}
\bA &= \bQ_3\bR_3
= 
\begin{bmatrix}
\textcolor{mylightbluetext}{-}q_{11} & q_{12} \\
\textcolor{mylightbluetext}{-}q_{21} & q_{22}
\end{bmatrix}
\begin{bmatrix}
\textcolor{mylightbluetext}{-}r_{11} & \textcolor{mylightbluetext}{-}r_{12} \\
0 & r_{22}
\end{bmatrix}\\
&=
\begin{bmatrix}
-\frac{1+\delta}{\sqrt{\delta^2+2\delta+2}} & \frac{1}{\sqrt{\delta^2+2\delta+2}}\\
-\frac{1}{\sqrt{\delta^2+2\delta+2}} & -\frac{1+\delta}{\sqrt{\delta^2+2\delta+2}}
\end{bmatrix}
\begin{bmatrix}
-\sqrt{\delta^2+2\delta+2} & -\frac{\delta^2+\delta}{(1+\delta)\sqrt{\delta^2+2\delta+2}} - \frac{\sqrt{\delta^2+2\delta+2}}{1+\delta}\\
0 & - \frac{\delta^2+2\delta}{\sqrt{\delta^2+2\delta+2}}
\end{bmatrix},
\end{aligned}
$$
or 
$$
\begin{aligned}
\bA &= \bQ_4\bR_4
= 
\begin{bmatrix}
\textcolor{mylightbluetext}{-}q_{11} & \textcolor{mylightbluetext}{-}q_{12} \\
\textcolor{mylightbluetext}{-}q_{21} & \textcolor{mylightbluetext}{-}q_{22}
\end{bmatrix}
\begin{bmatrix}
\textcolor{mylightbluetext}{-}r_{11} & \textcolor{mylightbluetext}{-}r_{12} \\
0 & \textcolor{mylightbluetext}{-}r_{22}
\end{bmatrix}\\
&=
\begin{bmatrix}
-\frac{1+\delta}{\sqrt{\delta^2+2\delta+2}} & -\frac{1}{\sqrt{\delta^2+2\delta+2}}\\
-\frac{1}{\sqrt{\delta^2+2\delta+2}} & \frac{1+\delta}{\sqrt{\delta^2+2\delta+2}}
\end{bmatrix}
\begin{bmatrix}
-\sqrt{\delta^2+2\delta+2} & -\frac{\delta^2+\delta}{(1+\delta)\sqrt{\delta^2+2\delta+2}} - \frac{\sqrt{\delta^2+2\delta+2}}{1+\delta}\\
0 & \frac{\delta^2+2\delta}{\sqrt{\delta^2+2\delta+2}}
\end{bmatrix}.
\end{aligned}
$$
Suppose $\delta=0.01$, the condition number of $\bA$ is 200, and the condition number of $\bR_1$, $\bR_2$, $\bR_3$, or $\bR_4$ is 1.
\end{example}

The example presented above, which demonstrates solving linear equations using the QR decomposition, helps mitigate the issue of ill-conditioning. For more detailed information on numerical stability, please refer to \citet{higham2002accuracy, zhang2017matrix, golub2013matrix, boyd2018introduction}, as well as Appendix~\ref{appendix:condition_number}.

\begin{problemset}
	
\item \textbf{Adjugate of orthogonal.} Let $\bQ\in\real^{n\times n}$ be orthogonal. Show that $\adjugate(\bQ) = \det(\bQ)\bQ^\top$ such that $\adjugate(\bQ)$ is also orthogonal (Definition~\ref{definition:adjugate}).
	
\item Prove that if $\bA$ is triangular and orthogonal, then $\bA$ must be diagonal.

\item Let $\bA\in\real^{n\times n}$ be skew-symmetric ($\bA^\top=-\bA$). Show that $(\bI-\bA)^{-1}(\bI+\bA)$ is orthogonal.

\item Let  $\bu$ and $\bv$ be two orthogonal unit vectors. Show that $\bu+\bv$ is orthogonal to $\bu-\bv$.

\item \textbf{Reflector.} Let  $\bu\in\real^n$ and $\bv\in\real^n$ be two orthogonal  vectors (not necessarily unit), where $\bu\in\mathcalV$ and $\bv\in\mathcalV^\perp$. Define $\ba\triangleq\bu+\bv$ and $\bb\triangleq\bu-\bv$. Show that there exists a unique elementary reflector (Householder reflector) $\bH\in\real^{n\times n}$ (Definition~\ref{definition:reflection_mat_intro}) such that $\bH\ba = \bb$.
Moreover, if $\mathcalV=\{\bw\}^\perp$, show that $\bH=\bI-2\frac{\bw\bw^\top}{\bw^\top\bw}$.

\item Let $\bu=[-\sin(\theta), \cos(\theta)]$ be a unit vector. Show that the Householder reflector determined by $\bu$ is $\bH=\scriptsize\begin{bmatrix}
\cos(2\theta) & \sin (2\theta) \\
\sin(2\theta) & -\cos(2\theta) 
\end{bmatrix}$.

\item \label{prob:semi_eq1} Let $\bQ,\bU\in\real^{m\times n}$ be two semi-orthogonal matrices with $m\geq n$. Show that $\bQ$ and $\bU$ have the same column space if and only if there exists an orthogonal matrix $\bP\in\real^{n\times n}$ such that $\bQ=\bU\bP$.

\item  \label{prob:semi_eq2} Let $\bQ,\bU\in\real^{n\times n}$ be orthogonal. Show that there exists an orthogonal matrix $\bP$ such that $\bQ=\bP\bU$.

\item  Let $\bQ,\bU\in\real^{m\times n}$ be two semi-orthogonal matrices with $m\geq n$. Show that there exists an orthogonal matrix $\bP\in\real^{m\times m}$ such that $\bQ=\bP\bU$. Compare this result to Problems~\ref{prob:semi_eq1} and \ref{prob:semi_eq2}. \textit{Hint: Complete the semi-orthogonal matrices into $m\times m$ orthogonal matrices.}

\item Let $\bA$ admit QR decomposition $\bA=\bQ\bR$. Show that $\bA=\bQ\bR$ is normal if and  only if $\bR\bQ$ is normal.

\item \label{problem:part_ortho} Consider the partition of an orthogonal matrix
$
\bQ = 
\scriptsize
\begin{bmatrix}
\underset{p\times p}{\bA} &\underset{p\times q}{\bB} \\
\underset{q\times p}{\bC} & \underset{q\times q}{\bD}
\end{bmatrix}\in\real^{n\times n}.
$
Show that $\rank(\bB)=\rank(\bC)$ and $\rank(\bD)=n+\rank(\bA)-2p$.

\item Consider the rank of matrices:
\begin{itemize}
\item Suppose matrices $\bA$ and $\bB$ have full column ranks. Show that $\bA\bB$ has full column rank.
\item Suppose $\bA\bB$ has full column ranks. Show that $\bB$ also has full column rank. However, it cannot be guaranteed that $\bA$ will have full column rank.
\item Discuss the rank of the upper triangular matrices  obtained from the QR decompositions of $\bA\bB$, $\bA$, and $\bB$ in various cases of the matrices involved.
\end{itemize}

\item We have stated in Theorem~\ref{theorem:qr-decomposition} that
$\bR$ is nonsingular in the reduced QR decomposition when 
 $\bA$ has full column rank $n$.
Suppose $\bA$ does not have full column rank. Examine  the relationship between the rank of $\bA$ and the number of nonzero entries in $\bR$.

\item Use the Gram-Schmidt process, Householder transformation, Givens rotation methods to find the orthonormal basis for the space spanned by the vectors
$$
\begin{aligned}
\bv_1 &= [1,3,7,5]^\top, \gap 
\bv_2 &= [6,3,6,3]^\top,  \gap 
\bv_3 &= [5,2,7,4]^\top.
\end{aligned}
$$

\item \textbf{Distance between a vector and a hyperplane.} Given a nonzero vector $\bzero\neq \ba\in\real^n$ and a scalar $\beta$,  we define  the hyperplane $H(\ba, \beta) = \{\bx\in\real^n:\ba^\top\bx+\beta=0\}$. For any vector $\by\in\real^n$, equipped with the projection along a line in Section~\ref{section:project-onto-a-vector}, show that the distance between the vector $\by$ and the hyperplane $H(\ba, \beta)$ is given by 
\begin{equation}
d(\by, H(\ba, \beta)) = \frac{\abs{\ba^\top\by + \beta}}{\normtwo{\ba}}.
\end{equation}
\textit{Hint: Choose two random points on the plane and first show that $\ba$ is orthogonal to the plane.}

\index{Nonlinear least squares}
\item \label{problem:gaus_new_jco} Prove Equation~\eqref{equation:gaus_new_jaco}, the gradient and Hessian of nonlinear least squares problems. \textit{Hint: Derive element-wise: 
$$
\frac{\partial f(\bx)}{\partial x_j} = \sum_{i=1}^{m} r_i(\bx) \frac{\partial r_i (\bx)}{\partial x_j}, 
\gap
\frac{\partial^2 f(\bx)}{\partial x_j \partial x_k} =\sum_{i=1}^{m} 
\left( 
\frac{\partial r_i(\bx)}{\partial x_j}\frac{\partial r_i(\bx)}{\partial x_k}
+
r_i(\bx)\frac{\partial^2 r_i(\bx)}{\partial x_j\partial x_k}
\right).
$$}

\item Consider $r(\bx)=[x_1, x_1+2x_2^2]^\top$, compute the Jacobian matrix of $r(\bx)$.

\item Although we use the fact that every orthogonal (orthonormal) list of vectors are linearly independent throughout our discussions, prove this rigorously. \textit{Hint: Assume in contrast they are linearly dependent and find a contradiction.}

\item Let $\bA\in\real^{m\times n}$ be given with $m\geq n$. Provide an algorithm that uses Householder reflectors to compute an orthogonal matrix $\bQ\in\real^{m\times m}$ such that $\bA=\bQ\bL$, where  $\bL[1:n, 1:n]$ is  lower triangular and $\bL[n + 1:m, :] = \bzero$.

\index{Range-symmetric}
\item Let $\bA\in\real^{n\times n}$ with rank $r$. Show that $\bA$ is range-symmetric (Definition~\ref{definition:speci_mat}) if and only if there exist a nonsingular matrix $\bS\in\real^{n\times n}$ and a nonsingular matrix $\bM\in\real^{r\times r}$ such that 
$$
\bA=\bS\begin{bmatrix}
\bM & \bzero\\
\bzero & \bzero 
\end{bmatrix}
\bS^\top .
$$
\textit{Hint: Consider the QR decomposition of $\bS=\bQ\bR$.}

\item Let $\bA\in\real^{n\times n}$ be skew-symmetric (Definition~\ref{definition:speci_mat}). Show that $\bI+\bA$ is nonsingular,  $\bB\triangleq(\bI-\bA)(\bI+\bA)^{-1}$ is orthogonal, $\bI+\bB = 2(\bI+\bA)^{-1}$, and $\det(\bB)=1$.

\item \label{prob:ortho_prese1} \textbf{Orthogonal preservation.} Let $\bQ\in\real^{n\times n}$ be orthogonal. Show that $\bx,\by\in\real^n$ are orthogonal if and only if $\bQ\bx$ and $\bQ\by$ are orthogonal.

\item \label{prob:ortho_presen} \textbf{Orthogonal preservation.} Let $\bQ\in\real^{n\times n}$ be orthogonal, and  let $\lambda$ be an  eigenvalue of $\bQ$. Show that $\lambda=\pm 1$, and $\bx\in\real^n$ is a (right) eigenvector of $\bQ$ associated with $\lambda$ if and only if $\bx$ is a left eigenvector
of $\bQ$ associated with $\lambda$.
\end{problemset}

%
%
%

%% file: chapter-utv.tex
\newpage
\chapter{UTV Decomposition: ULV and URV Decomposition}\label{section:ulv-urv-decomposition}
\begingroup
\hypersetup{
	linkcolor=structurecolor,
	linktoc=page,  
}
\minitoc \newpage
\endgroup

\section{UTV Decomposition}
\lettrine{\color{caligraphcolor}T}
The UTV decomposition extends the QR decomposition by factoring the matrix into two orthogonal matrices $\bA=\bU\bT\bV$, where $\bU$ and $\bV$ are orthogonal, and $\bT$ is either upper or lower triangular.~\footnote{These decompositions fall into a category known as the \textit{double-sided orthogonal decomposition}. We will see the UTV decomposition, complete orthogonal decomposition (Theorem~\ref{theorem:complete-orthogonal-decom}), and singular value decomposition are all in this notion.} 
The resulting $\bT$  provides support for rank estimation. 
Additionally, $\bT$ can take the form of a lower triangular matrix, leading to the ULV decomposition, or an upper triangular matrix, resulting in the URV decomposition. 
The UTV framework bears a resemblance to the singular value decomposition (SVD, see Chapter~\ref{chapter:SVD}) and  can be regarded as an inexpensive alternative to the SVD.
While the SVD provides a full diagonalization of the matrix, the UTV decomposition offers a more computationally efficient approach while still retaining many of the useful properties of the SVD. This makes the UTV decomposition particularly valuable in applications where computational resources are limited but accurate rank estimation and other matrix properties are needed.

\index{Decomposition: UTV}
\begin{theoremHigh}[Full ULV Decomposition]\label{theorem:ulv-decomposition}
Let $\bA\in\real^{m\times n}$ be given  with rank $r$. Then it  can be factored as 
$$
\bA = \bU \begin{bmatrix}
\bL & \bzero \\
\bzero & \bzero 
\end{bmatrix}\bV,
$$
where $\bU\in \real^{m\times m}$ and $\bV\in \real^{n\times n}$ are  orthogonal matrices, and $\bL\in \real^{r\times r}$ is a lower triangular matrix  with full rank.
\end{theoremHigh}

The existence of the ULV decomposition is  derived from the QR and LQ decompositions.
\begin{proof}[of Theorem~\ref{theorem:ulv-decomposition}]
For any rank-$r$ matrix $\bA=[\ba_1, \ba_2, \ldots, \ba_n]$, a column permutation matrix $\bP$ (Definition~\ref{definition:permutation-matrix}) can be employed to ensure  that the linearly independent columns of $\bA$ appear in the first $r$ columns of $\bA\bP$. Without loss of generality, we assume that  $\bb_1, \bb_2, \ldots, \bb_r$ represent  the $r$ linearly independent columns of $\bA$ and 
$$
\bA\bP = [\bb_1, \bb_2, \ldots, \bb_r, \bb_{r+1}, \ldots, \bb_n].
$$
Let $\bZ \triangleq [\bb_1, \bb_2, \ldots, \bb_r] \in \real^{m\times r}$. Since each $\bb_i$ lies in the column space of $\bZ$, we can find a matrix $\bE\in \real^{r\times (n-r)}$ such that 
$$
[\bb_{r+1}, \bb_{r+2}, \ldots, \bb_n] = \bZ \bE
\quad\implies\quad
\bA\bP = [\bb_1, \bb_2, \ldots, \bb_r, \bb_{r+1}, \ldots, \bb_n] = \bZ
\begin{bmatrix}
	\bI_r & \bE
\end{bmatrix},
$$
where $\bI_r$ denotes the $r\times r$ identity matrix. Additionally, $\bZ\in \real^{m\times r}$ has full column rank such that its full QR decomposition is given by $\bZ = \bU
\footnotesize
\begin{bmatrix}
\bR \\
\bzero
\end{bmatrix}$, 
where $\bR\in \real^{r\times r}$ is an upper triangular matrix with full rank, and $\bU$ is  orthogonal. This implies 
\begin{equation}\label{equation:ulv-smpl}
\bA\bP = \bZ
\begin{bmatrix}
	\bI_r & \bE
\end{bmatrix}
=
\bU\begin{bmatrix}
	\bR \\
	\bzero
\end{bmatrix}
\begin{bmatrix}
	\bI_r & \bE
\end{bmatrix}
=
\bU\begin{bmatrix}
	\bR & \bR\bE \\
	\bzero & \bzero 
\end{bmatrix}.
\end{equation}
The matrix $\bR$ has full rank, implying that $[\bR,\bR\bE ]$ also has full rank. 
As a result,  the full LQ decomposition of $[\bR,\bR\bE ]$
can be expressed as  $[\bL , \bzero ] \bV_0$, where $\bL\in \real^{r\times r}$ is a lower triangular matrix, and $\bV_0$ is an orthogonal matrix. 
Substituting this into Equation~\eqref{equation:ulv-smpl}, we obtain 
$$
\bA = \bU \begin{bmatrix}
	\bL & \bzero \\
	\bzero & \bzero 
\end{bmatrix}\bV_0 \bP^{-1}.
$$
Let $\bV \triangleq\bV_0 \bP^{-1}$, where $\bV_0$ and $\bP$ are both orthogonal matrices. 
Since $\bV_0 \bP^{-1}$ is a product of two orthogonal matrices, it is also an orthogonal matrix. This completes the proof.
\end{proof}
An alternative perspective on the proof of the ULV decomposition will be presented shortly in the proof of Theorem~\ref{theorem:complete-orthogonal-decom}. This approach involves both the rank-revealing QR decomposition and the standard QR decomposition.

\begin{figure}[H]
\centering
\vspace{-0.35cm}
\subfigtopskip=2pt
\subfigbottomskip=2pt
\subfigcapskip=-5pt
\subfigure[Reduced ULV decomposition.]{\label{fig:ulvhalf}
	\includegraphics[width=0.47\linewidth]{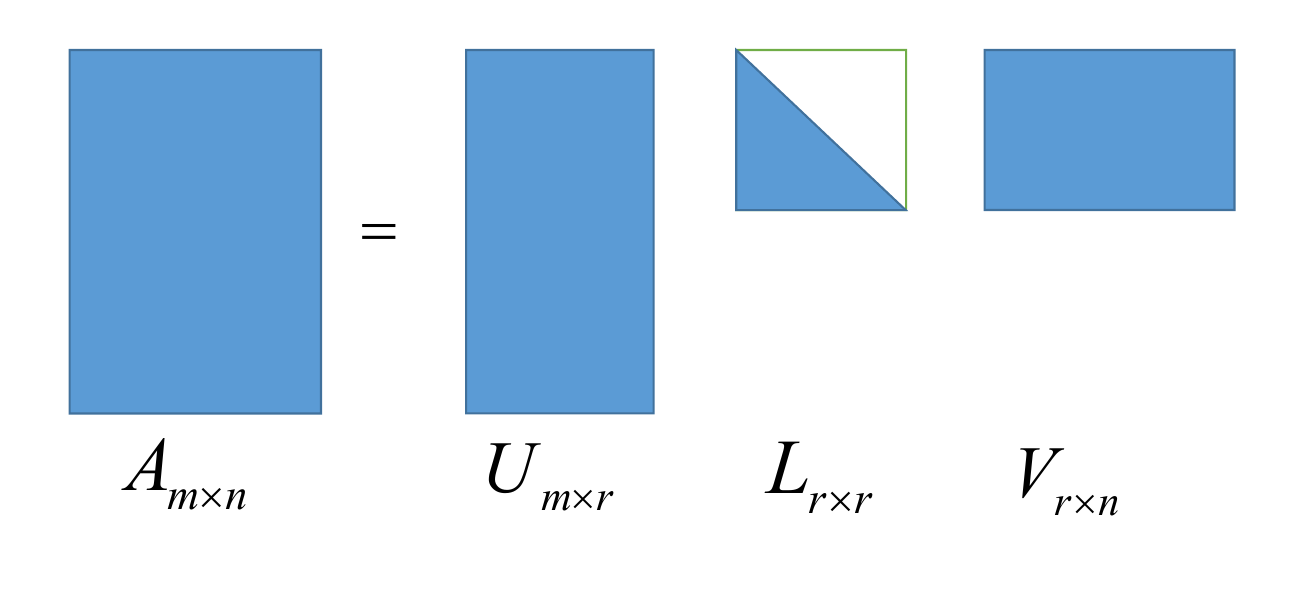}}
\quad 
\subfigure[Full ULV decomposition.]{\label{fig:ulvall}
	\includegraphics[width=0.47\linewidth]{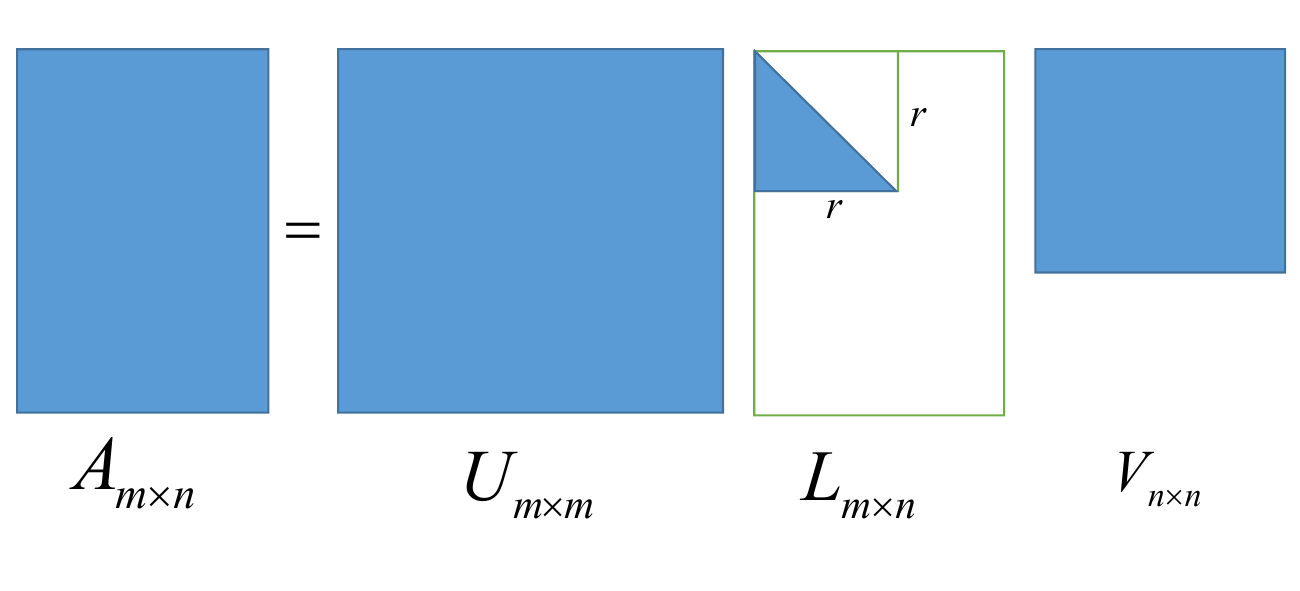}}
\quad
\subfigure[Reduced URV decomposition.]{\label{fig:urvhalf}
	\includegraphics[width=0.47\linewidth]{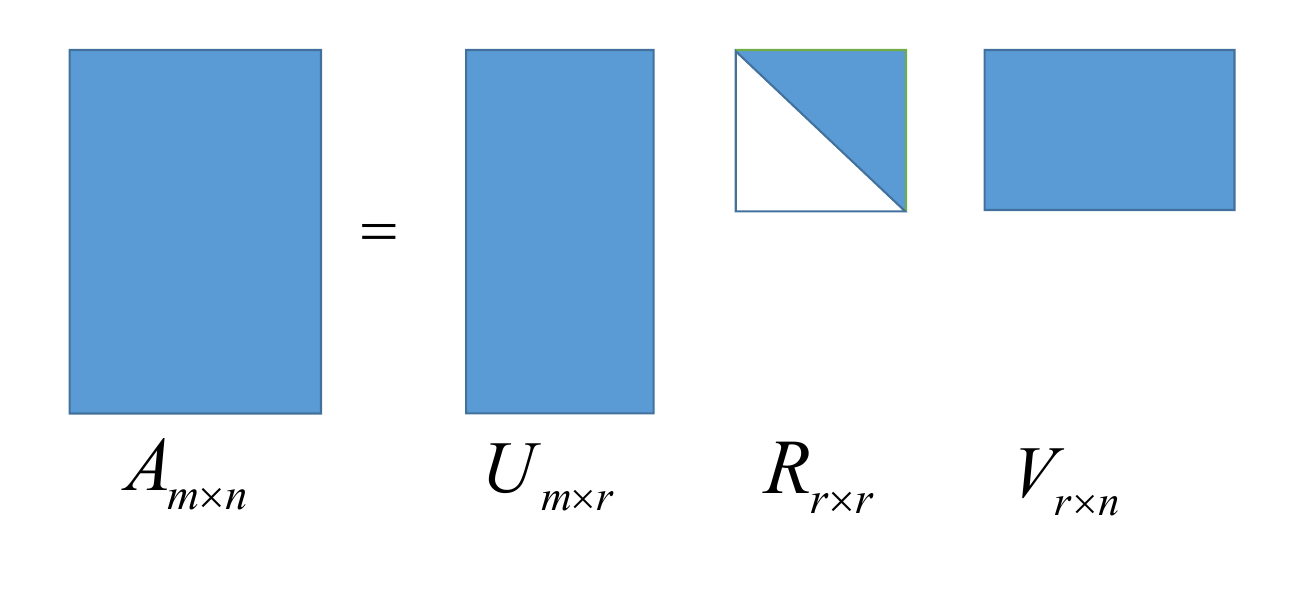}}
\quad 
\subfigure[Full URV decomposition.]{\label{fig:urvall}
	\includegraphics[width=0.47\linewidth]{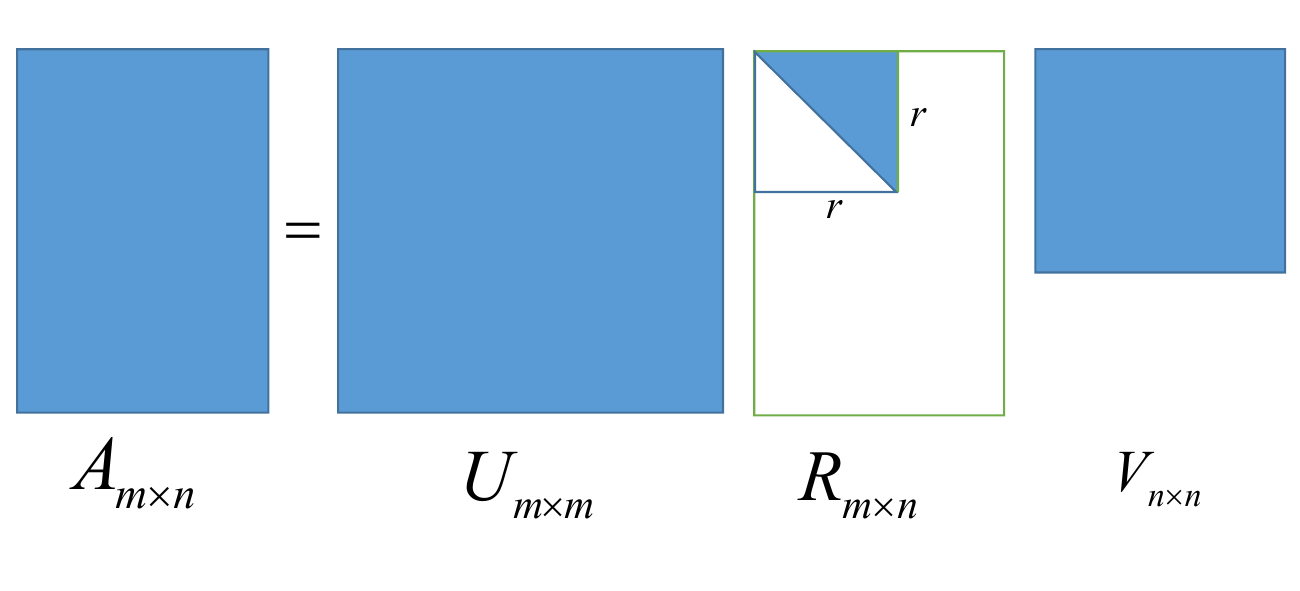}}
\caption{Comparison between the reduced and full ULV, and between the reduced and full URV, where white entries are zero and blue entries are not necessarily zero.}
\label{fig:ulv-comparison}
\end{figure}

Consider the ULV decomposition of the matrix $\bA$ given by
$$
\bA = \bU \begin{bmatrix}
	\bL & \bzero \\
	\bzero & \bzero 
\end{bmatrix}\bV.
$$
Let $\bU_0 \triangleq \bU[:,1:r]$ and $\bV_0 \triangleq \bV[1:r,:]$, i.e., $\bU_0$ contains only the first $r$ columns of $\bU$, and $\bV_0$ contains only the first $r$ rows of $\bV$. 
Consequently, we still obtain $\bA = \bU_0 \bL\bV_0$. 
This representation is known as the \textit{reduced ULV decomposition}.
The comparison between the reduced and the full ULV decompositions is depicted in Figure~\ref{fig:ulv-comparison}, where white entries are zero and blue entries are not necessarily zero.

Similarly, we can also establish  the URV decomposition as follows.
\begin{theoremHigh}[Full URV Decomposition]\label{theorem:urv-decomposition}
Let $\bA\in\real^{m\times n}$ be given  with rank $r$. Then it  can be factored as 
$$
\bA = \bU \begin{bmatrix}
\bR & \bzero \\
\bzero & \bzero 
\end{bmatrix}\bV,
$$
where $\bU\in \real^{m\times m}$ and $\bV\in \real^{n\times n}$ are two orthogonal matrices, and $\bR\in \real^{r\times r}$ is an upper triangular matrix with full rank.
\end{theoremHigh}
\begin{proof}[of Theorem~\ref{theorem:urv-decomposition}]
To abuse the notation, let  $\bA=[\ba_1, \ba_2, \ldots, \ba_m]^\top$ be the row partition of matrix $\bA$ with rank-$r$, where $\ba_1, \ba_2, \ldots, \ba_m\in\real^n$ are the rows of $\bA$, we can  construct a row permutation matrix $\bP$ such that the independent rows of $\bA$ appear in the first $r$ rows of $\bP\bA$. Without loss of generality, let's assume $\bb_1^\top, \bb_2^\top, \ldots, \bb_r^\top$ represent the $r$ linearly independent rows of $\bA$ and 
$$
\bP\bA = 
\footnotesize
\begin{bmatrix}
\bb_1,
\bb_2,
\ldots,
\bb_r,
\bb_{r+1},
\ldots ,
\bb_m
\end{bmatrix}^\top.
$$
Let $\bZ \triangleq [\bb_1, \bb_2, \ldots, \bb_r]^\top \in \real^{r\times n}$. Since each $\bb_i$ lies in the row space of $\bZ$, we can find a matrix $\bE\in \real^{(m-r)\times r}$ such that 
$$
\begin{bmatrix}
\bb_{r+1}^\top\\
\bb_{r+2}^\top \\
\vdots\\
\bb_m^\top
\end{bmatrix} 
\normalsize
= 
\bE\bZ
\gap 
\implies 
\gap 
\bP\bA = 
\footnotesize\begin{bmatrix}
\bb_1^\top \\
\bb_2^\top\\
\vdots\\
\bb_r^\top\\
\bb_{r+1}^\top\\
\vdots \\
\bb_m^\top
\end{bmatrix}
\normalsize
=
\begin{bmatrix}
\bI_r \\
\bE
\end{bmatrix}
\bZ,
$$
where $\bI_r$ denotes an $r\times r$ identity matrix. 
Furthermore, $\bZ\in \real^{r\times n}$ has full row rank such that its full LQ decomposition is given by $\bZ = \begin{bmatrix}
	\bL & \bzero 
\end{bmatrix} \bV$, where $\bL\in \real^{r\times r}$ is a lower triangular matrix with full rank and $\bV$ is an orthogonal matrix. 
This result indicates that
\begin{equation}\label{equation:urv-smpl}
	\bP\bA = \begin{bmatrix}
		\bI_r \\
		\bE
	\end{bmatrix}
	\bZ
	=
	\begin{bmatrix}
		\bI_r \\
		\bE
	\end{bmatrix}
	\begin{bmatrix}
		\bL & \bzero
	\end{bmatrix}
	\bV
	=
	\begin{bmatrix}
		\bL & \bzero \\
		\bE\bL & \bzero 
	\end{bmatrix}\bV.
\end{equation}
Since $\bL$ has full rank, it follows that 
$
\footnotesize
\begin{bmatrix}
	\bL \\
	\bE\bL 
\end{bmatrix}$ also has full rank. 
As a result, its full QR decomposition is given by 
$\bU_0
\footnotesize
\begin{bmatrix}
\bR \\
\bzero 
\end{bmatrix}$, 
where $\bR\in \real^{r\times r}$ is upper triangular and $\bU_0$ is an orthogonal matrix. Substituting this into Equation~\eqref{equation:urv-smpl}, we obtain 
$$
\bA =\bP^{-1} \bU_0 
\begin{bmatrix}
\bR & \bzero \\
\bzero & \bzero 
\end{bmatrix}\bV.
$$
Let $\bU \triangleq\bP^{-1} \bU_0$, which is the product of two orthogonal matrices, and is also an orthogonal matrix, from which the result follows.
\end{proof}

Additionally, there exists a reduced version of the URV decomposition, and the distinction between the full and reduced URV can be inferred from the context, as illustrated in Figure~\ref{fig:ulv-comparison}. 
The ULV and URV decompositions are sometimes collectively referred to as the \textit{UTV decomposition framework} \citep{fierro1997low, golub2013matrix}.
\begin{figure}[htbp]
	\centering
	\centering
	\resizebox{1.0\textwidth}{!}{%
		\begin{tikzpicture}[>=latex]
			
			\tikzstyle{state} = [draw, very thick, fill=white, rectangle, minimum height=3em, minimum width=6em, node distance=8em, font={\sffamily\bfseries}]
			\tikzstyle{stateEdgePortion} = [black,thick];
			\tikzstyle{stateEdge} = [stateEdgePortion,->];
			\tikzstyle{stateEdge2} = [stateEdgePortion,<->];
			\tikzstyle{edgeLabel} = [pos=0.5, text centered, font={\sffamily\small}];

			\node[ellipse, name=subspace, draw,font={\sffamily\bfseries},  node distance=7em, xshift=-9em, yshift=-1em,fill={colorals}]  {Subspaces};
			\node[state, name=qr, below of=subspace, xshift=-9em, yshift=0em, fill={colorlu}] {Full QR};
			\node[state, name=lq, below of=subspace, xshift=9em, yshift=0em, fill={colorlu}] {Full LQ};
			\node[state, name=rqr, left of=qr, xshift=-3em, fill={colorqr}] {Reduced QR};
			\node[state, name=rlq, right of=lq, xshift=3em, fill={colorqr}] {Reduced LQ};
			\node[state, name=utv, below of=subspace, xshift=0em, yshift=-9em, fill={colorlu}] {Full URV/ULV};
			\node[state, name=rutv, draw, right of=utv,xshift=7em, yshift=0em,fill={colorqr}] 
			{Reduced URV/ULV};
			
			\coordinate (qr2inter3) at ($(subspace.west -| qr.north) + (-0em,0em)$);
			\draw (subspace.west) edge[stateEdgePortion] (qr2inter3);
			\draw (qr2inter3) edge[stateEdge] 
			node[edgeLabel, text width=8em, yshift=0.8em]{\parbox{2em}{Column\\Space}} (qr.north);
			
			\coordinate (lq2inter3) at ($(subspace.east -| lq.north) + (-0em,0em)$);
			\draw (subspace.east) edge[stateEdgePortion] (lq2inter3);
			\draw (lq2inter3) edge[stateEdge] 
			node[edgeLabel, text width=7.25em, yshift=0.8em]{\parbox{2em}{Row\\Space}} (lq.north);
			
			\coordinate (rqr2inter1) at ($(subspace.north) + (0,1em)$);
			\coordinate (rqr2inter3) at ($(rqr2inter1-| rqr.north) + (-0em,0em)$);
			\draw (subspace.north) edge[stateEdgePortion] (rqr2inter1);
			\draw (rqr2inter1) edge[stateEdgePortion] (rqr2inter3);
			\draw (rqr2inter3) edge[stateEdge] 
			node[edgeLabel, text width=8em, yshift=0.8em]{\parbox{2em}{Column\\Space}} (rqr.north);

			\coordinate (rlr2inter3) at ($(rqr2inter1-| rlq.north) + (-0em,0em)$);
			\draw (rqr2inter1) edge[stateEdgePortion] (rlr2inter3);
			\draw (rlr2inter3) edge[stateEdge] 
			node[edgeLabel, text width=8em, yshift=0.8em]{\parbox{2em}{Row\\Space}} (rlq.north);

			\draw (rqr.east)
			edge[stateEdge] node[edgeLabel,yshift=0.5em]{Complete} 
			(qr.west) ;
			
			\draw (rlq.west)
			edge[stateEdge] node[edgeLabel,yshift=0.5em]{Complete} 
			(lq.east) ;
			
			\draw (utv.east)
			edge[stateEdge] node[edgeLabel,yshift=0.5em]{Deleting} 
			(rutv.west) ;
			
			\coordinate (qr2utv) at ($(qr.south) + (-0em,-3em)$);
			\coordinate (qr2utv2) at ($(qr2utv -| utv.north) + (-0em,0em)$);
			\coordinate (qr2utv3) at ($(lq.south) + (-0em,-3em)$);
			\draw (qr.south) edge[stateEdgePortion] (qr2utv);
			\draw (qr2utv) edge[stateEdgePortion] (qr2utv2);
			\draw (lq.south) edge[stateEdgePortion] (qr2utv3);
			\draw (qr2utv3) edge[stateEdgePortion] (qr2utv2);
			\draw (qr2utv2) edge[stateEdge] (utv.north);
			
			\begin{pgfonlayer}{background}
				\draw [join=round,cyan,dotted,fill={colormiddleleft}] ($(subspace.north west -|rqr.north west) + (-0.6em, 2em)$) rectangle ($(rutv.east -| rlq.south east) + (0.6em, -2.0em)$);
			\end{pgfonlayer}
		\end{tikzpicture}
	}
	\caption{Derive the ULV/URV by the QR and LQ.}
	\label{fig:decomQR}
\end{figure}

The relationship between the QR decomposition  and ULV/URV decompositions  is illustrated  in Figure~\ref{fig:decomQR}. 
Moreover, the ULV and URV decompositions can be employed to demonstrate a fundamental property in linear algebra, which we will soon observe: the row rank and column rank of any matrix are equal. Whereas, an elementary construction to prove this claim is provided in Appendix~\ref{append:row-equal-column}.

\index{Fundamental theorem}

We will soon observe that the forms of ULV and URV are closely related to the singular value decomposition (SVD). All three factorizations decompose the matrix $\bA$ into two orthogonal matrices.
Specifically, using ULV and URV, we can obtain sets of bases for the four subspaces of $\bA$ as described in the fundamental theorem of linear algebra. 
For instance, considering ULV, the first $r$ columns of $\bU$ form an orthonormal basis for the column space $\cspace(\bA)$, and the last $(m-r)$ columns of $\bU$ form an orthonormal basis for the null space $\nspace(\bA^\top)$.  
Similarly, the first $r$ rows of $\bV$ form an orthonormal basis for the row space $\cspace(\bA^\top)$, and the last $(n-r)$ rows form an orthonormal basis for the null space $\nspace(\bA)$ (similarly to the two-sided orthogonal decomposition, Theorem~\ref{theorem:two-sided-orthogonal}):
\begin{equation}\label{equation:ortho_space_utv}
\begin{aligned}
	\cspace(\bA) &= \spn\{\bu_1, \bu_2, \ldots, \bu_r\}, \qquad
	\nspace(\bA) = \spn\{\bv_{r+1}, \bv_{r+2}, \ldots, \bv_n\}, \\
	\cspace(\bA^\top) &= \spn\{ \bv_1, \bv_2, \ldots,\bv_r \} ,\qquad
	\nspace(\bA^\top) =\spn\{\bu_{r+1}, \bu_{r+2}, \ldots, \bu_m\}.
\end{aligned}
\end{equation}
The SVD extends its reach by establishing a connection between the two pairs of orthonormal bases: the transformation from the column basis to the row basis, and from the left null space basis to the right null space basis. More details about this connection will be presented in the dedicated SVD section.

\section{Complete Orthogonal Decomposition}
The \textit{complete orthogonal decomposition}, which also factors a matrix into two orthogonal matrices, is closely related to the UTV decomposition.
\begin{theoremHigh}[Complete Orthogonal Decomposition]\label{theorem:complete-orthogonal-decom}
Let $\bA\in\real^{m\times n}$ be given  with rank $r$. Then it  can be factored as 
$$
\bA = \bU \begin{bmatrix}
	\bT & \bzero \\
	\bzero & \bzero 
\end{bmatrix}\bV,
$$
where $\bU\in \real^{m\times m}$ and $\bV\in \real^{n\times n}$ are two orthogonal matrices, and $\bT\in \real^{r\times r}$ is a rank-$r$ matrix.
\end{theoremHigh}
\begin{proof}[of Theorem~\ref{theorem:complete-orthogonal-decom}]
By utilizing the column-pivoted QR decomposition (Theorem~\ref{theorem:rank-revealing-qr-general}), $\bA$ can be factored as 
$$
\bQ_1^\top \bA\bP = 
\begin{bmatrix}
	\bR_{11} & \bR_{12} \\
	\bzero   & \bzero 
\end{bmatrix},
$$
where $\bR_{11} \in \real^{r\times r}$ is upper triangular, $\bR_{12} \in \real^{r\times (n-r)}$, $\bQ_1\in \real^{m\times m}$ is an orthogonal matrix, and $\bP$ is a permutation matrix. 
Then it is not difficult to find a decomposition satisfying
\begin{equation}\label{equation:orthogonal-complete-qr-or-not}
\begin{bmatrix}
	\bR_{11}^\top \\
	\bR_{12}^\top
\end{bmatrix}
\triangleq
\bQ_2
\begin{bmatrix}
	\bS \\
	\bzero 
\end{bmatrix},
\end{equation}
where $\bQ_2$ is an orthogonal matrix, and $\bS$ is a rank-$r$ matrix. The decomposition is reasonable because  the matrix
$\footnotesize
\begin{bmatrix}
\bR_{11}^\top \\
\bR_{12}^\top
\end{bmatrix} \in \real^{n\times r}$ has rank $r$ whose columns stay in a subspace of $\real^{n}$. Nevertheless, the columns of $\bQ_2$ span the entire space of $\real^n$, where we can assume that the first $r$ columns of $\bQ_2$ span the same space as that of  $\footnotesize\begin{bmatrix}
\bR_{11}^\top \\
\bR_{12}^\top
\end{bmatrix}$. The matrix $\footnotesize\begin{bmatrix}
\bS \\
\bzero 
\end{bmatrix}$ serves to transform $\bQ_2$ into $\footnotesize\begin{bmatrix}
\bR_{11}^\top \\
\bR_{12}^\top
\end{bmatrix}$.
Then,  it follows that 
$$
\bQ_1^\top \bA\bP \bQ_2 = 
\begin{bmatrix}
	\bS^\top & \bzero \\
	\bzero & \bzero 
\end{bmatrix}.
$$
Let $\bU \triangleq\bQ_1$, $\bV\triangleq\bQ_2^\top\bP^\top$, and $\bT \triangleq \bS^\top$, we complete the proof.
\end{proof}
We observe that if we consider Equation~\eqref{equation:orthogonal-complete-qr-or-not} as the reduced QR decomposition of $
\footnotesize 
\begin{bmatrix}
	\bR_{11}^\top \\
	\bR_{12}^\top
\end{bmatrix}$, then the complete orthogonal decomposition reduces to the ULV decomposition.

\section{Applications}
\index{Least squares}\index{Rank-deficient}\index{Rank-deficiency}
\subsection{Application: Least Squares via UTV for Rank-Deficient Matrices}\label{section:ls-utv}
In Section~\ref{section:application-ls-qr}, we introduced the least squares (LS) via the full QR decomposition for full-rank matrices.
However, if often happens that the matrix may be rank-deficient. If $\bA$ does not have full column rank, $\bA^\top\bA$ is not invertible. 
In such cases, we can  use the ULV/URV decomposition to determine the least squares solution, as illustrated in the following theorem.
\begin{theorem}[LS via ULV/URV for Rank-Deficient Matrix]\label{theorem:qr-for-ls-urv}
Let $\bA\in \real^{m\times n}$ with rank $r$ and $m\geq n$~\footnote{See Problem~\ref{prob:als_pseudo1}$\sim$\ref{prob:als_pseudon} for more results on least squares with rank-deficient data matrices.}. Suppose $\bA=\bU\bT\bV$ is its full ULV/URV decomposition, where  $\bU\in\real^{m\times m}$ and $\bV\in \real^{n\times n}$ are orthogonal matrices, and
$$
\bT = \begin{bmatrix}
\bT_{11} & \bzero \\
\bzero & \bzero
\end{bmatrix}\in \real^{m\times n},
$$
where $\bT_{11} \in \real^{r\times r}$ is either a lower triangular matrix or an upper triangular matrix. 
Given the known vector $\bb\in \real^m$,  the LS solution with the minimal $\ell_2$ norm to $\bA\bx=\bb$ is given by 
$$
\bx_{LS} = \bV^\top 
\begin{bmatrix}
\bT_{11}^{-1}\bc\\
\bzero	
\end{bmatrix},
$$
where $\bc$ is the first $r$ components of $\bU^\top\bb$. 
\end{theorem}
	
\begin{proof}[of Theorem~\ref{theorem:qr-for-ls-urv}]
Since $\bA=\bU\bT\bV$ is the full UTV decomposition of $\bA$ and $m\geq n$, thus we have
$$
\begin{aligned}
\normtwo{\bA\bx-\bb}^2 &= (\bA\bx-\bb)^\top(\bA\bx-\bb)
\stackrel{*}{=}(\bA\bx-\bb)^\top\bU\bU^\top (\bA\bx-\bb) \\
&\stackrel{+}{=}\normtwo{\bU^\top \bA \bx-\bU^\top\bb}^2
=\normtwo{\bU^\top\bU\bT\bV \bx-\bU^\top\bb}^2\\
&=\normtwo{\bT\bV \bx-\bU^\top\bb}^2
=\normtwo{\bT_{11}\be - \bc}^2+\normtwo{\bd}^2,   
\end{aligned}
$$ 
\footnote{Equality ($*$) follows from the orthogonality of $\bU$; equality (+) follows from the invariance of the norm under orthogonal transformations.}
where $\bc$ corresponds to the first $r$ components of $\bU^\top\bb$, $\bd$ corresponds to the remaining  $m-r$ components of $\bU^\top\bb$, $\be$ represents the first $r$ components of $\bV\bx$, and  $\bff$ represents the remaining $n-r$ components of $\bV\bx$, i.e.,
$$
\bU^\top\bb 
= \begin{bmatrix}
\bc \\
\bd 
\end{bmatrix}
\qquad 
\text{and}\qquad
\bV\bx 
= \begin{bmatrix}
\be \\
\bff
\end{bmatrix}.
$$
The least squares solution can be determined through the process of backward or forward substitution in the upper or lower triangular system, denoted by  $\bT_{11}\be = \bc$, i.e., $\be = \bT_{11}^{-1}\bc$. 
To ensure that $\bx$ achieves the minimal $\ell_2$ norm, it is necessary for the vector $\bff$ to be null. 
That is,
$
\bx_{LS} = \bV^\top 
\footnotesize
\begin{bmatrix}
\bT_{11}^{-1}\bc\\
\bzero	
\end{bmatrix}.
$
This completes the proof.
\end{proof}

Moreover, we will shortly find that a similar argument can also be made via the singular value decomposition (SVD) in Section~\ref{section:application-ls-svd}.
This is due to the shared structural resemblance between SVD and ULV/URV, both being flanked by two orthogonal matrices.  
However, SVD extends its capabilities by diagonalizing the matrix. Enthusiastic readers can also obtain the least squares solution through the column-pivoted QR decomposition.

\paragraph{A word on the minimal $\ell_2$ norm LS solution.} 
In the context of the least squares problem, the ensemble of all solutions that minimize the objective function can be represented as follows:
$$
\mathcal{X} = \{\bx\in \real^n: \normtwo{\bA\bx-\bb} =\min\normtwo{\bA\bx-\bb} \},
$$
which is convex. 
Given $\bx_1$ and $\bx_2$ belonging to $\mathcal{X}$, and with $\lambda$ ranging within the interval $[0,1]$, the subsequent relationship holds true: 
$$
\normtwo{\bA(\lambda\bx_1 + (1-\lambda)\bx_2) -\bb} \leq \lambda\normtwo{\bA\bx_1-\bb} +(1-\lambda)\normtwo{\bA\bx_2-\bb} = \mathop{\min}_{\bx\in \real^n} \normtwo{\bA\bx-\bb}. 
$$
Thus, we have $\lambda\bx_1 + (1-\lambda)\bx_2 \in \mathcal{X}$ (Definition~\ref{definition:convexset}). 
In the presented proof, deviating from the setting $\bff=\bzero$ would yield multiple least squares solutions. 
Nonetheless, the unique solution with the minimal $\ell_2$ norm prevails as the optimal least squares solution. 
Within the context of the full-rank scenario detailed in Section~\ref{section:application-ls-qr}, the solution to the least squares problem remains distinct and invariably possesses the minimum $\ell_2$ norm property.
See also \citet{foster2003solving, golub2013matrix} for a comprehensive exploration of this subject.

\subsection{Application: Row Rank equals Column Rank  via UTV}
As previously discussed, the UTV framework can establish a significant theorem in linear algebra, demonstrating the equivalence of the row rank and column rank of a matrix.
In contrast, an elementary demonstration of this assertion is provided in Appendix~\ref{append:row-equal-column}.

Notice that in order to utilize the UTV framework in the proof, a slight adjustment to the claim regarding the existence of the UTV decomposition needs to be taken care of. For example, in Theorem~\ref{theorem:ulv-decomposition}, the condition for matrix $\bA$ is that it possesses a rank of $r$. 
However, since a rank $r$ inherently implies equivalence between row rank and column rank, a more suitable phrasing for this purpose is to state that matrix $\bA$ has a column rank of $r$ in Theorem~\ref{theorem:ulv-decomposition}.
For further elucidation, refer to the work of \citet{lu2021column}.

\begin{theorem}[Row Rank Equals Column Rank\index{Matrix rank}]\label{theorem:equal-dimension-rank}
The dimension of the column space of a matrix $\bA\in \real^{m\times n}$ is equal to the dimension of its
row space, i.e., the row rank and the column rank of a matrix $\bA$ are equal.
\end{theorem}

\begin{proof}[of Theorem~\ref{theorem:equal-dimension-rank}]
Any $m\times n$ matrix $\bA$ with rank $r$ can be factored as 
$$
\bA = \bU_0 \begin{bmatrix}
\bL & \bzero \\
\bzero & \bzero 
\end{bmatrix}\bV_0,
$$
where $\bU_0\in \real^{m\times m}$ and $\bV_0\in \real^{n\times n}$ are two orthogonal matrices, 
while $\bL\in \real^{r\times r}$ is a lower triangular matrix~\footnote{Instead of relying on the ULV decomposition, some texts employ elementary transformations $\bE_1$ and $\bE_2$ to establish the result, with the matrix $\bA$ expressed as 
$
\bA = \bE_1 
\scriptsize
\begin{bmatrix}
\bI_r & \bzero \\
\bzero & \bzero 
\end{bmatrix}
\normalsize\bE_2$.
}. 
Let 
$$\bD \triangleq \begin{bmatrix}
\bL & \bzero \\
\bzero & \bzero 
\end{bmatrix},
$$
the row rank and column rank of $\bD$ are apparently the same. 
To complete the proof, we must establish that the column rank of $\bA$ matches that of $\bD$, and similarly, the row rank of $\bA$ aligns with that of $\bD$.

Let $\bU \triangleq \bU_0^\top$ and $\bV\triangleq\bV_0^\top$, yielding  $\bD = \bU\bA\bV$. 
Breaking down the  concept above into two steps, a moment of reflexion reveals that, if we could first establish that the row rank and column rank of $\bA$ coincide with those of $\bU\bA$, and subsequently, prove that the row rank and column rank of $\bU\bA$ align with those of $\bU\bA\bV$, the proof can be successfully concluded.

\paragraph{Row rank and column rank of $\bA$ are equal to those of $\bU\bA$.} 
Let $\bB \triangleq \bU\bA$, and further consider $\bA=[\ba_1,\ba_2,\ldots, \ba_n]$ and $\bB=[\bb_1,\bb_2,\ldots,\bb_n]$ as the column partitions of $\bA$ and $\bB$, respectively. 
Thus, we have $[\bb_1,\bb_2,\ldots,\bb_n] = [\bU\ba_1,\bU\ba_2,\ldots, \bU\ba_n]$. If a linear combination $x_1\ba_1+x_2\ba_2+\ldots+x_n\ba_n=\bzero$, then premultiplying by $\bU$ gives
$$
\bU(x_1\ba_1+x_2\ba_2+\ldots+x_n\ba_n) = x_1\bb_1+x_2\bb_2+\ldots+x_n\bb_n = \bzero.
$$
Consider distinct indices $j_1, j_2, \ldots, j_r$ ranging from 1 to $n$. If the set $\{\ba_{j_1}, \ba_{j_2}, \ldots, \ba_{j_r}\}$ is linearly independent, then the set $\{\bb_{j_1}, \bb_{j_2}, \ldots, \bb_{j_r}\}$ must also exhibit linearly independence. This observation signifies that 
$$
\dim(\cspace(\bB)) \leq \dim(\cspace(\bA)).
$$
Similarly, employing $\bA = \bU^\top\bB$, we can deduce that
$$
\dim(\cspace(\bA)) \leq \dim(\cspace(\bB)).
$$
This implies 
$$
\dim(\cspace(\bB)) = \dim(\cspace(\bA)).
$$
Apply the process to $\bB^\top$ and $\bA^\top$, we have 
$$
\dim(\cspace(\bB^\top)) = \dim(\cspace(\bA^\top)).
$$
This indicates the row rank and column rank of $\bA$ and $\bB=\bU\bA$ are the same.
\paragraph{Row rank and column rank of $\bU\bA$ are equal to those of $\bU\bA\bV$.}
Similarly, by applying above discussion to $\bU\bA$ and $\bU\bA\bV$, we can also demonstrate that the row rank and column rank of $\bU\bA$ and $\bU\bA\bV$ are the same. This concludes the proof.
\end{proof}

\begin{problemset}
\item Find the UTV decomposition of $\diag(\bA,\bB)$ using the UTV decompositions of $\bA$ and $\bB$.
\item Prove rigorously that the four subspaces in Equation~\eqref{equation:ortho_space_utv} can be spanned using the UTV decomposition.

\item Let $\bA=\bU\bB\bV$ be given, where $\bU$ and $\bV$ are orthogonal. Show that $\sum_{i,j}\abs{a_{ij}}^2=\sum_{i,j} \abs{b_{ij}}^2$.
\textit{Hint: Show that $\trace(\bA^\top\bA)=\trace(\bB^\top\bB)$.}

\item Describe the process of computing the UTV decomposition using either Householder reflectors or Givens rotations. What is the computational complexity associated with each method?

\item Let $\bA,\bB\in\complex^{n\times n}$ be given, where $\bB$ is nonsingular. Show that there exist unitary matrices $\bU,\bV\in\complex^{n\times n}$ such that $\bA=\bU\bT_A\bV$ and $\bB=\bU\bT_B\bV$, where $\bT_A$ and $\bT_B$ are upper triangular.  Furthermore, show that  the main diagonal entries of $\bT_B^{-1}\bT_A$ are the eigenvalues $\bB^{-1}\bA$. 
When these eigenvalues are real, show that all the underlying matrices can be chosen to be real, and $\bU$ and $\bV$ are orthogonal.
\textit{Hint: Consider the Schur decomposition of $\bB^{-1}\bA=\bU\bT\bU^*$ (Theorem~\ref{theorem:schur-decomposition_complex}) and the QR decomposition of $\bB\bU$.}

\item Let $\bA,\bB\in\real^{n\times n}$ be given, where $\bB$ is nonsingular. Show that there are orthogonal matrices $\bU,\bV\in\real^{n\times n}$ such that $\bA=\bU\bT_A\bV$ and $\bB=\bU\bT_B\bV$, where $\bT_A$ and $\bT_B$ are upper quasi-triangular (see Theorem~\ref{theorem:quasi_schur}). \textit{Hint: Consider the quasi-triangular Schur decomposition of $\bB^{-1}\bA=\bU\bT\bU^\top$ (Theorem~\ref{theorem:quasi_schur}) and the QR decomposition of $\bB\bU$.}

\item \textbf{Read Chapter~\ref{chapter:SVD} first.} Let $\bA,\bB\in\real^{m\times n}$ and consider the definitions of left equivalence, right equivalence, and orthogonally biequivalence in Definition~\ref{definition:biequivalent}. Show that 
\begin{itemize}
\item  The matrices $\bA$ and $\bB$ are orthogonally left equivalent if and only if $\bA^\top\bA=\bB^\top\bB$.
\item  The matrices $\bA$ and $\bB$ are orthogonally right equivalent if and only if $\bA\bA^\top =\bB\bB^\top$.
\item The matrices $\bA$ and $\bB$  are orthogonally biequivalent if and only if $\bA$ and $\bB$ have the same singular values with the same multiplicity.
\end{itemize}

\item  Let $\bA,\bB\in\real^{m\times n}$. Show that 
\begin{itemize}
\item  The matrices $\bA$ and $\bB$ are left equivalent if and only if $\nspace(\bA) = \nspace(B)$.
\item  The matrices $\bA$ and $\bB$ are right equivalent if and only $\cspace(\bA) = \cspace(\bB)$.
\item  The matrices $\bA$ and $\bB$ are biequivalent if and only if $\rank(\bA)=\rank(\bB)$.
\end{itemize}
\end{problemset}

%% file: chapter-cr.tex
\newpage 
\part{Data Interpretation and Information Distillation}\label{part:data-interation}
\index{Fundamental theorem}
\section*{Introduction}
\lettrine{\color{caligraphcolor}F}
For a matrix $\bA\in \real^{m\times n}$ with rank $r$, the fundamental theorem of linear algebra (Theorem~\ref{theorem:fundamental-linear-algebra}) ensures the existence of $r$ linearly independent columns and rows.  
Consequently, $\bA$ admits 
$$
\text{(DI1)} \qquad \underset{m\times n}{\bA} = \underset{m\times r}{\bC} \,\,\, \underset{r\times n}{\bF},
$$
where $\bC$ consists of $r$ linearly independent columns of $\bA$ and $\bF$ is used to reconstruct all the columns of $\bA$, since all the columns of $\bC\bF$ are linear combinations of the columns of $\bC$. 
To illustrate  this, if we let $\bF=[\bff_1, \bff_2, \ldots, \bff_n]$ be the column partition of $\bF$, then
$$
\bA = \bC\bF=[\bC\bff_1, \bC\bff_2, \ldots, \bC\bff_n],
$$
where each column $\bC\bff_i$ is formed as a linear combination of the columns in $\bC$. 
The columns of $\bC$ are commonly referred to as the \textit{spanning columns} of $\bA$. 
Alternatively, $\bA$ can be represented as
$$
\text{(DI2)} \qquad \underset{m\times n}{\bA} = \underset{m\times r}{\bD} \,\,\, \underset{r\times n}{\bR},
$$
where $\bR$ comprises  $r$ linearly independent rows of $\bA$ and $\bD$ is utilized to reconstruct all the rows of $\bA$, as the rows of $\bD\bR$ can be formed through linear combinations of the rows in $\bR$. 
To understand this, if we let $\bD=[\bd_1^\top; \bd_2^\top; \ldots; \bd_m^\top]$ be the row partition of $\bD$, then
$$
\bA = \bD\bR=
\begin{bmatrix}
	\bd_1^\top\bR \\
	\bd_2^\top\bR \\
	\vdots \\
	\bd_m^\top\bR 
\end{bmatrix},
$$
where each row $\bd_i^\top\bR$ is a  linear combination of the rows in $\bR$. And the rows of $\bR$ are commonly referred to as the \textit{spanning rows} of $\bA$.

The factorization in (DI1) resembles the QR decomposition, with the first factor spanning the column space of $\bA$ (and being orthogonal in the latter decomposition).
On the other hand, (DI2) is similar to the LQ decomposition.  
Nevertheless, (DI1) holds several advantages when compared to methods such as the QR decomposition:
\begin{itemize}
\item It is sometimes advantageous to work with a basis that consists of a subset of the
columns of $\bA$ itself. 
In such instances, the stipulation of orthogonal or orthogonormal basis vectors is often relinquished;
\item If $\bA$ is sparse and nonnegative, then $\bC$ retains these properties;
\item In general, the (DI1) factorization demands less memory for storage and exhibits greater efficiency in computation compared to the QR decomposition.
The (DI1) is extensively employed as a feature selection tool, facilitating essential extraction and enabling processing of voluminous data that initially exceeds RAM capacity.
From the above representation, it becomes evident that only $r(m+n)$ entries need to be stored, in contrast to the $mn$ entries in the original matrix $\bA$ or the $mn+\frac{(n+1)n}{2}$ entries in the \textit{reduced} QR decomposition (Figure~\ref{fig:qr-comparison});

\item Finding the indices associated with the spanning columns often proves invaluable for data interpretation and analysis. This approach facilitates identifying a subset of columns that encapsulates the matrix's essential information. Furthermore, if the columns of matrix $\bA$ hold specific interpretations, such as transactions in a transaction data set, the same interpretation extends to the columns of $\bC$;
\item Calculating the least squares solution $\bA\bx=\bb$ can be achieved by determining the least squares solution of $\bC\widetildebx=\bb$. The former one projects $\bb$ onto the column space of $\bA$, while the latter projects it onto the column space of $\bC$. This factorization (DI1) can be regarded as the first phase of variable selection in the context of least squares or linear models. A brief review of the variable selection procedure is given in Section~\ref{section:append-column-qr}  as an application of the QR decomposition.

\item The (DI1) often preserves ``the physical essence" of a problem in a way that the QR does not;

\item  We can eliminate  the non-relevant portions of the data, which consist of error and redundant information, via (DI1).
\end{itemize}



For the rest of this part, we will explore various adaptations of the aforementioned factorization, each with its own distinct focus.
These variations encompass aspects such as descriptiveness, interpretability (easy to interpret or not), conditioning (well-conditioned or not), and the way for selecting basis vectors for the column space, for the row space, or for both spaces. An overview of the broader context is illustrated in Figure~\ref{fig:data-interpretation-world-picture}.


\begin{figure}[htbp]
\centering
\begin{widepage}
	\centering
\resizebox{0.65\textwidth}{!}{%
\begin{tikzpicture}[>=latex]
%
%
\tikzstyle{state} = [draw, very thick, fill=white, rectangle, minimum height=3em, minimum width=6em, node distance=8em, font={\sffamily\bfseries}]
\tikzstyle{stateEdgePortion} = [black,thick];
\tikzstyle{stateEdge} = [stateEdgePortion,->];
\tikzstyle{stateEdge2} = [stateEdgePortion,<->];
\tikzstyle{edgeLabel} = [pos=0.5, text centered, font={\sffamily\small}];

\node[state, name=cr, fill={colorcr}] {CR};
\node[state, name=rank, above of=cr,yshift=-3.2em,  fill={colorcr}] {Rank};
\node[state, name=interpolative, left of=cr, yshift=-2em,  fill={colorcr}] {Interpolative};
\node[state, name=skeleton, right of=cr, yshift=-2em, fill={colorcr}] {\parbox{3.5em}{ Skeleton (CUR)}};

\draw (rank.south)
edge[stateEdge] node[edgeLabel, yshift=0em, xshift=-0em]{Special Case} 
(cr.north) ;

\draw ($(cr.east) + (0,1em)$)
edge[stateEdge, bend left=12.5] node[edgeLabel, xshift=0.5em, yshift=0.8em]{} 
(skeleton.north);
\draw[decoration={text along path,
	text={|\sffamily|Same C},text align={center}},decorate] ($(cr.east) + (-1.8em,1.2em)$) to [bend left=20.5]  ($(skeleton.north) + (2.5em,-1.5em)$);

\draw (rank.west)
edge[stateEdge, bend left=-30] node[edgeLabel, yshift=0em, xshift=-0em]{} 
($(interpolative.north) + (-2.5em,0em)$) ;
\draw[decoration={text along path,
text={|\sffamily|Special Case},text align={center}},decorate] ($(interpolative.north) + (-3em,0.3em)$) to [bend left=30]  ($(rank.west) + (-0em,0.em)$) ; 

\draw (rank.east)
edge[stateEdge, bend left=30] node[edgeLabel, yshift=0em, xshift=-0em]{} 
($(skeleton.north) + (2.5em,0em)$) ;
\draw[decoration={text along path,
text={|\sffamily|``Special" Case},text align={center}},decorate] ($(rank.east) + (0.7em,0.2em)$) to [bend left=30]  ($(skeleton.north) + (2.5em,0.5em)$);

\draw ($(cr.west) + (0,1em)$)
edge[stateEdge, bend left=-12.5] node[edgeLabel, xshift=-0.5em, yshift=0.8em]{} 
(interpolative.north);
\draw[decoration={text along path,
text={|\sffamily|Independent},text align={center}},decorate] ($(interpolative.north) + (-0.5em,0.1em)$) to [bend left=12.5]  ($(cr.west) + (.5em,1.2em)$);
\draw[decoration={text along path,
text={|\sffamily|Columns},text align={center}},decorate] ($(interpolative.north) + (-0.8em,-1.1em)$) to [bend left=12.5]  ($(cr.west) + (2.5em,0.3em)$);

\begin{pgfonlayer}{background}
\draw [join=round,cyan,dotted,fill={coloruppermiddle}] ($(rank.north west  -| interpolative.north west) + (-0.8em, +0.8em)$) rectangle ($(skeleton.south east) + (0.8em, -0.8em)$);
\end{pgfonlayer}

\end{tikzpicture}
}
\end{widepage}
\caption{Data interpretation relationship. Refer also to its position in the matrix decomposition landscape depicted in Figure~\ref{fig:matrix-decom-world-picture}.}
\label{fig:data-interpretation-world-picture}
\end{figure}

\subsection*{\textbf{A Brief Introduction to SVD}}
For the analysis of the	 decomposition in this part, particular the interpolative decomposition, we need to introduce the singular value decomposition (SVD) and the Eckart-Young-Mirsky theorem.
Further details are available in Chapter~\ref{chapter:SVD}. For any matrix $\bA\in \real^{m\times n}$ with rank $r$, the SVD is given by:
$$
\bA = \bU \bSigma \bV^\top = \sum_{i=1}^r \sigma_i \bu_i \bv_i^\top,
$$ 
where $\bu_i$ and $\bv_i$ represent the left and right singular vectors of $\bA$, respectively.  
The matrix $\bSigma\in \real^{m\times n}$ is a rectangular diagonal matrix whose entries are the singular values arranged in decreasing order, 
$\sigma_1\geq \sigma_2\geq \ldots \geq \sigma_r \textcolor{mylightbluetext}{>} \sigma_{r+1}=\ldots=0$, along the main diagonal, containing only $r$ nonzero singular values.

For a given integer $k$ satisfying $1\leq k\leq r$, we define $\bA_k$ as the \textit{truncated SVD (TSVD)} of $\bA$, consisting of the largest $k$ terms, i.e., $\bA_k = \sum_{i=1}^{k} \sigma_i\bu_i\bv_i^\top$ from the SVD of $\bA=\sum_{i=1}^{r} \sigma_i\bu_i\bv_i^\top$ by zeroing out the $r-k$ trailing singular values of $\bA$. 
It's important to note that $\bA_k$ provides the best rank-$k$ approximation to $\bA$ in terms of the spectral norm or the Frobenius norm (Section~\ref{section:svd-low-rank-approxi}). 
That is, for any rank-$k$ matrix $\bB$, we can establish that $\norm{\bA-\bA_k}_2 \leq \norm{\bA-\bB}_2$ or $\norm{\bA-\bA_k}_F \leq \norm{\bA-\bB}_F$, where the norms are the spectral norm and Frobenius norm, respectively.
And the distance between $\bA$ and $\bA_k$ can be quantified as follows:
$$
\norm{\bA-\bA_k}_F = \sqrt{\sigma_{k+1}^2+\sigma_{k+2}^2+\ldots +\sigma_{r}^2} \leq \norm{\bA-\bB}_F,
$$
or 
$$
\norm{\bA-\bA_k}_2 = \sigma_{k+1} \leq \norm{\bA-\bB}_2.
$$
\subsection*{\textbf{Subset Selection Problem}}

Subset selection is a technique employed to choose a subset of columns from a real matrix $\bA\in\real^{m\times n}$ such that the subset accurately depicts the complete matrix and is far from being rank-deficient.
(DI1) above serves as an illustrative instance of column selection with full rank $r$ columns. 
However, for any integer $1 \leq k \leq r$, we always attempt to find $k$ linearly independent columns that best capture the information in the matrix. 
More formally, the subset selection problem entails discovering a permutation matrix $\bP$ in such a way that
$$
\bA\bP = 
\begin{bmatrix}
\underset{m\times k}{\bA_1}& \underset{m\times (n-k)}{\bA_2}
\end{bmatrix},
$$ 
where $\bA_1\in \real^{m\times k}$ contains $k$ linearly independent columns of $\bA$ and
\begin{enumerate}
\item The smallest singular value of $\bA_1$ is as large as possible (similarly to SVD!). That is, there exists a constant value $\eta$ such that the $k$-th singular value of $\bA_1$ (i.e., the smallest singular value of $\bA_1$) is bounded by that of $\bA$:
$$
\frac{\sigma_k(\bA)}{\eta} \leq \sigma_k(\bA_1) \leq \sigma_k(\bA).
$$
\item The remaining $n-k$ redundant columns of $\bA_2$ can be effectively approximated using the $k$ columns from $\bA_1$. 
That is, there exists a matrix $\bW$ capable of reconstructing $\bA_2$ from $\bA_1$, or the value of
$$
\mathop{\min}_{\bW \in \real^{k\times (n-k)}} \norm{\bA_1\bW-\bA_2}_2
$$
is small, i.e., there exists a constant $\eta$ such that the distance is bounded by the $(k+1)$-th singular value of $\bA$:
$$
\sigma_{k+1}(\bA) \leq \mathop{\min}_{\bW \in \real^{k\times (n-k)}} \norm{\bA_1\bW-\bA_2}_2 \leq \eta \sigma_{k+1}(\bA).
$$
\end{enumerate}

\newpage
\chapter{CR Decomposition}\label{section:cr-decomposition}
\begingroup
\hypersetup{
	linkcolor=structurecolor,
	linktoc=page,  
}
\minitoc \newpage
\endgroup

\section{CR Decomposition}
\lettrine{\color{caligraphcolor}T}
The CR decomposition of a matrix is introduced in \citet{strang2021every, stranglu} to reveal its rank and the relationships between its columns and rows. As usual,  we shall first present the result, reserving the discussion of its existence and derivation for subsequent sections.

\index{Decomposition: CR}
\begin{theoremHigh}[CR Decomposition]\label{theorem:cr-decomposition}
Let $\bA \in \real^{m \times n}$ be any matrix with  rank-$r$. Then it can be factored as 
$$
\underset{m\times n}{\bA} = \underset{m\times r}{\bC} \gapthree \underset{r\times n}{\bR},
$$
where $\bC$ contains the first $r$ linearly independent columns of $\bA$, and $\bR$ is an $r\times n$ matrix utilized  to reconstruct the columns of $\bA$ based on the columns of $\bC$. In particular, $\bR$ corresponds to the \textit{reduced row echelon form (RREF)} of $\bA$ without the zero rows.

The storage required for the decomposition is subsequently reduced or, in some cases, increased from $mn$ floating-point numbers to $r(m+n)$ floating-point numbers.
\end{theoremHigh}

The advantages of the CR decomposition are immediately apparent: the CR decomposition helps in identifying the independent columns and rows of a matrix, which are crucial for understanding the matrix's rank and its column and row spaces.
For large matrices, it can be used to approximate the matrix by retaining only the most significant columns and rows, which can be useful in reducing computational complexity.
The decomposition can be used to solve linear systems $\bA\bx = \bb $ by transforming the system into a simpler form using the matrices $ \bC $ and $ \bR $.
Since both $\bC$ and $\bR$ have full rank, the CR decomposition provides a way to compute the Moore-Penrose pseudo-inverse of the matrix, which is given by $ \bA^+ = \bR^+\bC^+ $, where $ \bR^+ $ and $ \bC^+ $ are the pseudo-inverses of $ \bR $ and $ \bC $, respectively (see Problem~\ref{prob:cr_pseudo}).
In the context of graphs, the CR decomposition can be applied to the incidence matrix of a graph to analyze the flow of currents or other conserved quantities, as seen in Kirchhoff's Current Law \citep{strang2021three}.

\section{Existence of  CR Decomposition}

Given that the matrix $\bA$ has a rank of $r$, it contains $r$ linearly independent columns. Consequently, we select these linearly independent columns from $\bA$ and incorporate them into $\bC$:

\begin{tcolorbox}[title={Find $r$ Linearly Independent Columns from $\bA$},colback=\mdframecolorTheorem]
1. If column 1 of $\bA$ is nonzero, include it in the columns of $\bC$.

2. If column 2 of $\bA$ is not a scalar multiple of column 1, include it in the columns of $\bC$.

3. If column 3 of $\bA$ is not a linear combination of columns 1 and 2, include it in the columns of $\bC$.

4. Continue this process until we find $r$ linearly independent columns (or all the linearly independent columns if the rank $r$ is not known in advance).
\end{tcolorbox}

Once we have obtained the $r$ linearly independent columns from $\bA$, we can establish the existence of the CR decomposition through the perspective of matrix multiplication within the column space.

\paragraph{Perspective on matrix multiplication in column space.} The multiplication of two matrices, $\bD\in \real^{m\times k}$ and $\bE\in \real^{k\times n}$, results in the matrix $\bA=\bD\bE=\bD[\be_1, \be_2, \ldots, \be_n] = [\bD\be_1, \bD\be_2, \ldots, \bD\be_n]$. In other words, each column in  $\bA$ is a linear combination of the columns of $\bD$. 

\begin{proof}[of Theorem~\ref{theorem:cr-decomposition}]
Given that the rank of matrix $\bA$ is $r$ and $\bC$ comprises $r$ linearly independent columns from $\bA$, the column space of $\bC$ is equivalent to that of $\bA$. 
Any other column $\ba_i$ of $\bA$  can be represented as a linear combination of the columns of $\bC$, i.e., there exists a vector $\br_i$ such that $\ba_i = \bC \br_i$, $\forall i\in \{1, 2, \ldots, n\}$. By arranging  these vectors $\br_i$'s as the columns of a matrix $\bR$, we obtain 
$$
\bA = [\ba_1, \ba_2, \ldots, \ba_n] = [\bC \br_1, \bC \br_2, \ldots, \bC \br_n]= \bC \bR,
$$
from which the result follows.
\end{proof}

\section{Reduced Row Echelon Form (RREF)}\label{section:rref-cr}
In Section~\ref{section:gaussian-elimination} on Gaussian elimination, we introduced the elimination matrix (a lower triangular matrix; see \eqref{equation:elimination_mat}) and the permutation matrix to facilitate the transformation of $\bA$ into an upper triangular form. 
Let's now revisit the Gaussian elimination process applied to a $4\times 4$ square matrix, where $\boxtimes$ represents a value that is not necessarily zero, and \textbf{boldface} indicates the value has just been changed:
\begin{tcolorbox}[title={Gaussian Elimination for a Square Matrix},colback=\mdframecolorTheorem]
$$
\footnotesize
\begin{sbmatrix}{\bA}
\boxtimes & \boxtimes & \boxtimes & \boxtimes \\
\boxtimes & \boxtimes & \boxtimes & \boxtimes \\
\boxtimes & \boxtimes & \boxtimes & \boxtimes \\
\boxtimes & \boxtimes & \boxtimes & \boxtimes
\end{sbmatrix}
\stackrel{\bE_1}{\longrightarrow}
\begin{sbmatrix}{\bE_1\bA}
\boxtimes & \boxtimes & \boxtimes & \boxtimes \\
\bm{0} & \bm{0} & \bm{\boxtimes} & \bm{\boxtimes} \\
\bm{0} & \bm{\boxtimes} & \bm{\boxtimes} & \bm{\boxtimes} \\
\bm{0} & \bm{\boxtimes} & \bm{\boxtimes} & \bm{\boxtimes}
\end{sbmatrix}
\stackrel{\bP_1}{\longrightarrow}
\begin{sbmatrix}{\bP_1\bE_1\bA}
\boxtimes & \boxtimes & \boxtimes & \boxtimes \\
\bm{0} & \textcolor{mylightbluetext}{\bm{\boxtimes}} & \bm{\boxtimes} & \bm{\boxtimes} \\
\bm{0}  & \bm{0} & \textcolor{mylightbluetext}{\bm{\boxtimes}} & \bm{\boxtimes} \\
0 & \boxtimes & \boxtimes & \boxtimes
\end{sbmatrix}
\stackrel{\bE_2}{\longrightarrow}
\begin{sbmatrix}{\bE_2\bP_1\bE_1\bA}
\boxtimes & \boxtimes & \boxtimes & \boxtimes \\
0 &  \textcolor{mylightbluetext}{\boxtimes} & \boxtimes & \boxtimes \\
0 & 0  & \textcolor{mylightbluetext}{\boxtimes} & \boxtimes \\
0 & \bm{0}  & \bm{0} & \textcolor{mylightbluetext}{\bm{\boxtimes}}
\end{sbmatrix}.
$$
\end{tcolorbox}

Moreover, Gaussian elimination can be extended to rectangular matrices. We illustrate this with an example involving a $4\times 5$ matrix:
\begin{tcolorbox}[title={Gaussian Elimination for a Rectangular Matrix},colback=\mdframecolorTheorem]
$$
\footnotesize
\begin{sbmatrix}{\bA}
\textcolor{mylightbluetext}{2} & \boxtimes & 10 & 9 & \boxtimes\\
\boxtimes & \boxtimes & \boxtimes & \boxtimes & \boxtimes\\
\boxtimes & \boxtimes & \boxtimes & \boxtimes & \boxtimes\\
\boxtimes & \boxtimes & \boxtimes & \boxtimes & \boxtimes\\
\end{sbmatrix}
\stackrel{\bE_1}{\longrightarrow}
\begin{sbmatrix}{\bE_1\bA}
\textcolor{mylightbluetext}{2} & \boxtimes & 10 & 9 & \boxtimes\\
\bm{0} & \bm{0} & \textcolor{mylightbluetext}{\bm{5}} & \bm{6} & \bm{\boxtimes}\\
\bm{0} & \bm{0} & \bm{2} & \bm{\boxtimes} & \bm{\boxtimes}\\
\bm{0} & \bm{0} & \bm{\boxtimes} & \bm{\boxtimes} & \bm{\boxtimes}\\
\end{sbmatrix}
\stackrel{\bE_2}{\longrightarrow}
\begin{sbmatrix}{\bE_2\bE_1\bA}
\textcolor{mylightbluetext}{2} & \boxtimes & 10 & 9 & \boxtimes\\
0 & 0 & \textcolor{mylightbluetext}{5} & 6 & \boxtimes\\
0 & 0 & \bm{0} & \textcolor{mylightbluetext}{\bm{3}} & \bm{\boxtimes}\\
0 & 0 & \bm{0} & \bm{0} & \bm{0}\\
\end{sbmatrix},
$$
\end{tcolorbox}
\noindent where the numbers highlighted in \textcolor{mylightbluetext}{blue} are \textit{pivots}, as defined previously (Definition~\ref{definition:pivot}. In the context of Gaussian elimination, a pivot element is the first nonzero element in a row when performing row operations to transform a matrix into its row echelon form or reduced row echelon form.). The resulting matrix is referred to as the \textit{row echelon form} of  $\bA$. 
Note that, in this specific example,  the fourth row becomes a zero row. 
To proceed, we can perform row operations by subtracting each row by a scalar multiple of the next row to ensure that the entries above the pivots are all zero:
\begin{tcolorbox}[title={Reduced Row Echelon Form: Get Zero Above Pivots},colback=\mdframecolorTheorem]
$$
\footnotesize
\begin{sbmatrix}{\bE_2\bE_1\bA}
\textcolor{mylightbluetext}{2} & \boxtimes & 10 & 9 & \boxtimes\\
0 & 0 & \textcolor{mylightbluetext}{5} & 6 & \boxtimes\\
0 & 0 & 0 & \textcolor{mylightbluetext}{3} & \boxtimes\\
0 & 0 & 0 & 0 & 0\\
\end{sbmatrix}
\stackrel{\bE_3}{\longrightarrow}
\begin{sbmatrix}{\bE_3\bE_2\bE_1\bA}
\textcolor{mylightbluetext}{2} & \boxtimes & \bm{0}  & \bm{-3}  & \bm{\boxtimes} \\
0 & 0 & \textcolor{mylightbluetext}{5} & 6 & \boxtimes\\
0 & 0 & 0 & \textcolor{mylightbluetext}{3} & \boxtimes\\
0 & 0 & 0 & 0 & 0\\
\end{sbmatrix}
\stackrel{\bE_4}{\longrightarrow}
\begin{sbmatrix}{\bE_4\bE_3\bE_2\bE_1\bA}
\textcolor{mylightbluetext}{2} & \boxtimes & 0 & \bm{0} & \bm{\boxtimes}\\
0 & 0 & \textcolor{mylightbluetext}{5} & \bm{0} & \bm{\boxtimes}\\
0 & 0 & 0 & \textcolor{mylightbluetext}{3} & \boxtimes\\
0 & 0 & 0 & 0 & 0\\
\end{sbmatrix},
$$
\end{tcolorbox}
\noindent where $\bE_3$ subtracts twice the second row from the first row, while $\bE_4$ adds the third row to the first row and subtracts twice the third row from the second row. Ultimately, achieving the   full \textit{reduced row echelon form} involves ensuring that the pivots are set to 1:
\begin{tcolorbox}[title={Reduced Row Echelon Form: Make The Pivots To Be 1},colback=\mdframecolorTheorem]
\begin{equation}\label{equation:cr_exp1}
\footnotesize
\begin{sbmatrix}{\bE_4\bE_3\bE_2\bE_1\bA}
	\textcolor{mylightbluetext}{2} & \boxtimes & 0 & 0 & \boxtimes\\
	0 & 0 & \textcolor{mylightbluetext}{5} & 0 & \boxtimes\\
	0 & 0 & 0 & \textcolor{mylightbluetext}{3} & \boxtimes\\
	0 & 0 & 0 & 0 & 0\\
\end{sbmatrix}
\stackrel{\bE_5}{\longrightarrow}
\begin{sbmatrix}{\bE_5\bE_4\bE_3\bE_2\bE_1\bA}
	\textcolor{mylightbluetext}{\bm{1}} & \bm{\boxtimes} & \bm{0} & \bm{0} & \bm{\boxtimes}\\
	\bm{0} & \bm{0} & \textcolor{mylightbluetext}{\bm{1}} & \bm{0} & \bm{\boxtimes}\\
	\bm{0} & \bm{0} & \bm{0} & \textcolor{mylightbluetext}{\bm{1}} & \bm{\boxtimes}\\
	0 & 0 & 0 & 0 & 0\\
\end{sbmatrix},
\end{equation}
\end{tcolorbox}
\noindent where $\bE_5$ is employed  to set the pivots to 1. Note here, it is not necessary for 
the transformation matrices $\bE_1, \bE_2, \ldots, \bE_5$ to be lower triangular matrices, as typically required in LU decomposition; they can also take the form of permutation matrices or other types.  The resultant matrix is denoted as the \textit{reduced row echelon form (RREF)} of $\bA$, characterized by containing 1's as pivots and zeros above the pivots.

More generally, suppose the row echelon form (with zero above the pivots) of $\bA\in\real^{m\times n}$ is $\bF$. To obtain the reduced row echelon form of $\bA$, we apply a set of transformations $\bE_1, \bE_2, \ldots, \bE_r$, where 
$$
\bE_i\triangleq
\bE_{i;\eta} \triangleq
\begin{bmatrix}
1&&&&&&\\
&\ddots&&&&&\\
&& 1&&&&\\
&&& \eta &&&\\
&&&& 1 &&\\
&&&&& \ddots &\\
&&&&&&1\\
\end{bmatrix}
=
\bI + (\eta-1)\be_i\be_i^\top,
\gap 
\text{with $1\leq i\leq r$},
$$
where $\be_i$ represents the $i$-th unit basis vector, $r$ denotes the rank of $\bA$ (the number of nonzero rows in $\bF$), and $\eta$ is the inverse of the $i$-th pivot of $\bF$. Each $\bE_i$ is invertible with $\bE_{i;\eta}^{-1} = \bE_{i;\eta^{-1}} = \bI + (\frac{1}{\eta}-1)\be_i\be_i^\top$.
The reduced row echelon form of $\bA$  can then be obtained by $(\bE_{r}\bE_{r-1}\ldots\bE_1\bF)$.

Rigorously, we present the definition of the RREF as follows.
\begin{definition}[Reduced Row Echelon Form, RREF]
Let $\bA\in\real^{m\times n}$. The row echelon form of $\bA$ is characterized by the following conditions:
\begin{enumerate}
\item[1.] The leading nonzero entry (pivot, Definition~\ref{definition:pivot}) of the $(i+1)$-th row occurs  to the right of the leading nonzero entry of the $i$-th row.
\item[2.] All entries in a column below a  leading nonzero entry are zeros.
\item[3.] All zero rows are placed at the bottom of the matrix.
\end{enumerate} 
Property 2 is a simple consequence of property 1, but we include it for emphasis.
If the matrix in row echelon form satisfies the following additional conditions, then it is
in reduced row echelon form:
\begin{enumerate}
\item[4.] The leading nonzero entry in each row is 1. That is, the pivots are 1.
\item[5.] The entries above a pivot are zero.
\end{enumerate}
\end{definition}
\begin{exercise} Let $\bA\in\real^{m\times  n}$ be any matrix. Show that there exists a sequence of row transformations $\bE_1, \bE_2, \ldots, \bE_k$ such that $\bB = \bE_k\bE_{k-1}\ldots\bE_1\bA$ is in reduced row echelon form. \textit{Hint: Use induction}.
\end{exercise}

\begin{lemma}[Rank and Pivots]\label{lemma:rank-is-pivots}
The rank of a matrix $\bA$ is equal to the number of pivots (in its reduced row echelown form).
Thus, the rank of a matrix is equal to the rank of its RREF.
\end{lemma}

\begin{proof}[of Lemma~\ref{lemma:rank-is-pivots}]
According to Proposition~\ref{proposition:rowspa_rowele}, the row space of $\bA$ is identical to that of its RREF. The rank of the RREF is equal to the number of pivots. Therefore, the rank of $\bA$ is equal to the number of pivots as well.
\end{proof}

\begin{lemma}[RREF in CR\index{Reduced row echelon form}]\label{lemma:r-in-cr-decomposition}
The reduced row echelon form of the matrix $\bA$, excluding zero rows, corresponds to the matrix $\bR$ in the CR decomposition.
\end{lemma}
\begin{proof}[Informal Proof of Lemma~\ref{lemma:r-in-cr-decomposition}]
Informally, following the example provided earlier in \eqref{equation:cr_exp1}, we can express the matrix $\bA$ as follows
$$
\bE_5\bE_4\bE_3\bE_2\bE_1\bA \triangleq \bR_0 
\quad \implies \quad  
\bA = (\bE_5\bE_4\bE_3\bE_2\bE_1)^{-1}\bR_0.
$$
We observe that  columns 1, 3, and 4 of $\bR_0$ contain only a single nonzero element, equal to 1. This observation allows us to construct a matrix $\bC$ (identical to the  ``column matrix" in the CR decomposition) whose first three columns are equal to columns 1, 3, and 4 of  $\bA$, i.e., $\bC=[\ba_1, \ba_3, \ba_4]$. 
Additionally, since the elements of last row of $\bR_0$ are all zeros, we can safely disregard the last row of $\bR_0$ in our computation. 
Notably, this matrix $\bC$ is unique in its ability to reconstruct columns 1, 3, and 4 of $\bA$, as the pivots of $\bR_0$ are all 1. We obtain
$
\bA=\bC\bR,
$
the CR decomposition.
\end{proof}

We provide a rigorous proof for the uniqueness of the RREF of a matrix.
\begin{theorem}[Uniqueness of RREF]\label{theorem:unique_rref}
Let $\bA\in\real^{m\times n}$ be any matrix with rank $r$. If $\bX$ and $\bY$ are two reduced row echelon forms of $\bA$, obtained by applying two sequences of elementary row operations $\bE_1, \bE_2, \ldots, \bE_p$ and $\bF_1, \bF_2, \ldots, \bF_q$, respectively, where 
$$
\bX = \underbrace{\bE_p\ldots \bE_2\bE_1}_{\bE}\bA 
\gap\text{and}\gap
\bY = \underbrace{\bF_q\ldots \bF_2\bF_1}_{\bF}\bA.
$$
Then, the two reduced row echelon forms are the same, i.e., $\bX=\bY$ and $\bE=\bF$.

\end{theorem}
\begin{proof}[of Theorem~\ref{theorem:unique_rref}]
Let $\bB \triangleq \bE\bF^{-1}=\bE_p\ldots\bE_2\bE_1\bF_1^{-1}\bF_2^{-1}\ldots \bF_q^{-1}$. 
Then, we have $\bX = \bB\bY$ and $\bY = \bB^{-1}\bX$.
Suppose $\be_i$ represents the $i$-th unit basis vector in $\real^n$. Then the  $i$-th column of $\bX$ and $\bY$ can be denoted by $\bx_i =\bX\be_i$ and $\by_i = \bY\be_i$, respectively.
\paragraph{Zero columns match.} We note that if $\bx_i=\bzero$, then $\by_i =\bB^{-1}\bx_i =\bzero$; and if  $\by_i=\bzero$, then $\bx_i =\bB\by_i =\bzero$. Thus, the zero columns in $\bX$ and $\bY$ are aligned.
Without loss of generality, we assume that $\bX$ and $\bY$ do not contain zero columns for the following analysis.
\paragraph{First column.} Since we assume $\bX$ and $\bY$ do not contain zero columns, it follows that $\bx_1=\by_1=\be_1$. This also implies the first column of $\bB$ is $\be_1$.
We refer to columns in $\bX$ or $\bY$ that contain pivots as \textit{pivot columns}, and those do not  as \textit{non-pivot columns}.

\index{Pivot columns}
\index{Non-pivot columns}

\paragraph{Non-pivot columns between the first and second pivot columns.}
Suppose the indices of the pivot columns in $\bX$ are $\{i_1, i_2, \ldots, i_r\}$, and the indices of the pivot columns in $\bY$ are $\{j_1, j_2, \ldots, j_r\}$. According to Lemma~\ref{lemma:rank-is-pivots}, there are $r$ pivot columns. 
And we have already shown that $i_1=j_1=1$. Then for $k\in\{2,3,\ldots, j_2-1\}$, we suppose $\by_k=\lambda \be_1$ for some nonzero $\lambda$. We have
$$
\bx_k = \bB\by_k = \bB\lambda \be_1=\lambda\bb_1  =\lambda\be_1
\quad\implies\quad \bx_k = \by_k, \gap \forall k\in\{2,3,\ldots, j_2-1\}.
$$
Conversely, suppose $\bx_k=\lambda \be_1$ for some nonzero $\lambda$ with $k\in\{2,3,\ldots, i_2-1\}$. In this case, we also find that $\bx_k=\by_k$. This implies the non-pivot columns $k=2, 3, \ldots, j_2-1$ of $\bX$ and $\bY$ are the same, and also $i_2=j_2$ (the second index for the pivot column of $\bX$ and $\bY$ are the same as well, and $\bx_{j_2} = \by_{j_2}=\be_2$).

To conclude, we have demonstrated that the first $j_2$ columns of  $\bX$ and $\bY$ are identical.
\paragraph{Non-pivot columns between the second and third pivot columns.} Since $\bx_{j_2} = \by_{j_2} = \be_2$, we have $\bx_{j_2}=\bB\by_{j_2} = \bB\be_2 = \bb_2=\be_2$, i.e., the second column of $\bB$ is $\be_2$. On the other hand, considering the non-pivot columns of $k\in\{j_2+1, j_2+2, \ldots, j_3-1\}$ in $\bY$, we assume $\by_k = \lambda_1\be_1+\lambda_2\be_2$. Then we have 
$$
\begin{aligned}
\bx_k = \bB\by_k =\bB(\lambda_1\be_1+&\lambda_2\be_2) =\lambda_1\bb_1+\lambda_2\bb_2 =\lambda_1\be_1+\lambda_2\be_2 \\
&\quad\implies\quad \bx_k=\by_k,\gap \forall 
k\in\{j_2+1, j_2+2, \ldots, j_3-1\}.
\end{aligned}
$$
Conversely, suppose $\bx_k = \lambda_1\be_1+\lambda_2\be_2$ with $k\in\{i_2+1, i_2+2, \ldots, i_3-1\}$, we also have $\bx_k=\by_k$. 
This implies the non-pivot columns $k=j_2+1, j_2+2, \ldots, j_3-1$ of $\bX$ and $\bY$ are the same, and also $j_3=i_3$ (the third index for the pivot column of $\bX$ and $\bY$ are the same, and $\bx_{j_3} = \by_{j_3}=\be_3$).

To conclude, for now we have shown that the first $j_3$ columns of $\bX$ and $\bY$ are identical.
Continuing with this process, we conclude the claim.
\end{proof}

Readers interested in an inductive proof of the above theorem can also refer to \citet{gallier2017fundamentals}.

\begin{exercise}[Determinant of RREF]
Show that the determinant of a matrix $\bA$ is nonzero if and only if its RREF is $\bI$.
Therefore,  from the multiplicative  property of determinants (Remark~\ref{remark:determinant-intermezzo}), the determinant of the matrix can be calculated  by recording the operations performed during the process leading to its RREF.
\end{exercise}

In summary, we begin by calculating the reduced row echelon form of matrix $\bA$ using $rref(\bA)$. Then, the matrix $\bC$ in the CR decomposition is derived from $\bA$ by removing all non-pivot columns (identified by examining columns in $rref(\bA)$ that do not contain a pivot).  
Simultaneously, the factor $\bR$ is constructed by eliminating zero rows from $rref(\bA)$. 
This process is essentially a specialized case of the \textit{rank decomposition} (Theorem~\ref{theorem:rank-decomposition}), but it stands out because it explicitly involves the RREF. Thus, we introduce it here for its specific relevance.

An important characteristic of $\bR$  is that a subset of its  $r$ columns, each containing a pivot, together form an $r\times r$ identity matrix. It's worth reiterating that we can obtain this matrix $\bR$ simply by eliminating the zero rows from the RREF. As noted in \citet{strang2021every}, a notation for the RREF that retains the zero rows is denoted by $\bR_0$:
$$
\bR_0 = rref(\bA)=
\begin{bmatrix}
	\bR \\
	\bzero
\end{bmatrix}=
\begin{bmatrix}
	\bI_r & \bF \\
	\bzero & \bzero
\end{bmatrix}\bP,~\footnote{Permutation matrix $\bP$ on the right side of a matrix is used to interchange the columns of that matrix. 
}
$$
where the $n\times n$ permutation matrix $\bP$ arranges the columns of the $r\times r$ identity matrix $\bI_r$ into their correct positions, aligning them with the first $r$ linearly independent columns of the original matrix $\bA$.

Previously, we established a fundamental theorem in linear algebra using the UTV framework, demonstrating the equality of row rank and column rank for any matrix (Theorem~\ref{theorem:equal-dimension-rank}).
The CR decomposition also elucidates this theorem.
\begin{proof}[{of Theorem~\ref{theorem:equal-dimension-rank}, the Second Way}]
In the context of the CR decomposition of matrix $\bA = \bC\bR$, we can express $\bR$ as $\bR = [\bI_r, \bF]\bP$, where $\bP$ represents an $n\times n$ permutation matrix used to arrange the columns of the $r\times r$ identity matrix $\bI_r$ in their appropriate positions, as discussed above. 
It can be easily verified that the $r$ rows of $\bR$ are linearly independent (since the submatrix $\bI_r$ within $\bR$ is nonsingular), so  the row rank of $\bR$ is $r$. 

First, according to the definition of the CR decomposition, the $r$ columns of $\bC$ are derived from $r$ linearly independent columns of $\bA$, and the column rank of $\bA$ is $r$. Furthermore,
\begin{itemize}
\item Since $\bA=\bC\bR$, all rows of $\bA$ are linear combinations of the rows of $\bR$. In other words, the row rank of $\bA$ is not greater than the row rank of $\bR$: $\cspace(\bA^\top)\subseteq \cspace(\bR^\top)$.

\item From  $\bA=\bC\bR$, we  have $(\bC^\top\bC)^{-1}\bC^\top\bC\bR = (\bC^\top\bC)^{-1}\bC^\top\bA$, which simplifies to $\bR = (\bC^\top\bC)^{-1}\bC^\top\bA$. 
Here, $\bC^\top\bC$ is nonsingular since $\bC$ has full column rank $r$. 
Consequently, all rows of $\bR$ are also linear combinations of the rows of $\bA$. That is, the row rank of $\bR$ is not greater than the row rank of $\bA$: $\cspace(\bR^\top)\subseteq \cspace(\bA^\top)$.

\item By ``sandwiching," the row rank of $\bA$ is equal to the row rank of $\bR$, which is $r$ (and also $\cspace(\bA^\top)= \cspace(\bR^\top)$).
\end{itemize}

Therefore, both the row rank and column rank of $\bA$ are equal to $r$, from which the result follows.
\end{proof}

In the proof above, we leverage the CR decomposition to demonstrate the equivalence of a matrix's row rank and column rank. 
Alternatively, an elementary proof without using CR decomposition or Gaussian elimination is provided in Appendix~\ref{append:row-equal-column}. 
Furthermore, we also discuss the special form of the pseudo-inverse derived from the CR decomposition in Appendix~\ref{appendix:pseudo-inverse}.

\section{Computing  CR Decomposition via  Gaussian Elimination}
The central step in computing the CR decomposition is to determine the reduced row  echelon form of  matrix $\bA$. Suppose $\bA$ is of size $m\times n$, we will consider the worst-case scenario in which $\bA$ is of full rank.

$\bullet$ Obtain the row echelon form (REF)\index{Row echelon form}:

A. Initialize $\bR$ as $\bA$. Utilize the first row of $\bA$ to nullify values below $r_{11}$ in the first column, with the consideration of permutation if $r_{11}=0$. Store the result in $\bR$. (This step requires  $2(m-1)n+(m-1)$ flops, as discussed in the subsequent proof);

B. Employ the second row of the resultant matrix $\bR$ to zero out the values below $r_{22}$, with the consideration of permutation if $r_{22}=0$. (This step requires $2(m-2)(n-1)+(m-2)$ flops, as discussed in the subsequent proof);

C. Continue this process until the row echelon form is obtained.

$\bullet$ Obtain the reduced row  echelon form (RREF)\index{Reduced row echelon form}:
	
D. Utilize the last row to make the values above the last pivot to be zero (ignore if the row consists entirely of zeros); 

E. Employ the penultimate row to set the values above the second-to-last pivot to  zero (ignore if the row consists entirely of zeros); 

F. Continue this process until we obtained the reduced row  echelon form, and then divide each row by a scalar to set the pivots to  1. Note that there are $m-1$ such steps if $m\leq n$ and $n-1$ such steps if $m>n$ (further details are discussed in the subsequent proof).

This process is described in Algorithm~\ref{alg:cr-decomposition-alg}.
\begin{algorithm}[H] 
\caption{CR Decomposition via the Gaussian Elimination} 
\label{alg:cr-decomposition-alg} 
\begin{algorithmic}[1] 
\Require 
Rank-$r$ matrix $\bA=[\ba_1, \ba_2, \ldots, \ba_n]$ with size $m\times n $;
\Statex \textbf{Stage A: Obtain the upper triangular matrix} 
\State Initially set $\bR=\bA$; 
\For{$k=1$ to $m-1$}  \Comment{Operate on the $k$-th column}
\State \algoalign{Use row $k$ of $\bR$ to make the values below $r_{kk}$ to be 0 (permutation involved if $r_{kk}=0$);}
\EndFor
\Statex \textbf{Stage B: Zero out the entries above the pivots} 
\For{($k=m$ to $2$ if $m\leq n$)  or ($k=n$ to $2$ if $m > n$)} 
\State \algoalign{Use row $k$ to make the values above the $k$-th pivot to be zero (ignore if the row consists entirely of zeros);}
\EndFor
\State Set the pivot values to 1 by dividing specific scalar values;
\State Store the columns of $\bA$ correspond to the indices of pivot columns into matrix $\bC$.
\end{algorithmic} 
\end{algorithm}

\begin{theorem}[Algorithm Complexity: CR via the Gaussian Elimination\index{Gaussian elimination complexity}]\label{theorem:cr-decomposition-alg} 
Algorithm~\ref{alg:cr-decomposition-alg} requires 	
$$
\text{cost}=\left\{
\begin{aligned}
	&\sim 2m^2n-m^3 \text{ flops}, \qquad  &\text{if }&m\leq n; \\
	&\sim mn^2 \text{ flops},\qquad &\text{if }&m>n.
\end{aligned}
\right.
$$
to compute the CR decomposition of an $m\times n$ matrix.
\end{theorem}

\begin{proof}[of Theorem~\ref{theorem:cr-decomposition-alg} ]
Procedure (A) discussed above requires $(m-1)$ divisions to obtain the \textit{multipliers}, followed by  $(m-1)\times n$ multiplications with the multipliers (we operate on $(m-1)$ rows of length-$n$ vector). To zero out the values below $r_{11}$, it involves $(m-1)\times $ $n$ subtractions. Consequently, procedure (A) costs  \underline{$2(m-1)n+(m-1)$} flops.

Procedure (B) requires $(m-2)$ divisions to calculate  the multipliers and it involves $(m-2)\times (n-1)$ multiplications with the multipliers (we operate on $(m-2)$ rows of length-$(n-1)$ vector). To zero out the values below $r_{22}$, it involves $(m-2)\times (n-1)$ subtractions. As a result, procedure (B) costs \underline{$2(m-2)(n-1)+(m-2)$} flops.

The procedure can continue, and we can summarize the cost for each loop of step 3 in the following table:
\begin{table}[H]
\setlength{\tabcolsep}{7.6pt}
\begin{tabular}{l|l|l|l}
$k$      & Get multipliers & Multipliers multiply each row & Rows subtraction \\ \hline
1        & $2:m=m-1$ rows      & $(m-1)(n)$                      & $(m-1)(n)$         \\
2        & $3:m=m-2$ rows      & $(m-2)(n-1)$                  & $(m-2)(n-1)$     \\
3        & $4:m=m-3$ rows      & $(m-3)(n-2)$                  & $(m-3)(n-2)$     \\
$\vdots$ & $\vdots$            & $\vdots$                      & $\vdots$         \\
$k$      & $k+1:m=m-k$ rows    & $(m-k)(n-k+1)$                & $(m-k)(n-k+1)$   \\
$\vdots$ & $\vdots$            & $\vdots$                      & $\vdots$         \\
$m-1$    & $m:m=1$ row        & $(1)(n-m+2)$                  & $(1)(n-m+2)$    
\end{tabular}
\end{table}
We notice that the value $n-m+2$ in the last row of the  table above may not always be positive. 
Therefore, we will discuss the complexity of computing the row echelon form (REF) separately for two cases.

\paragraph{Get REF, case 1: $m\leq n$.}  The procedure can go on until  loop $m-1$.
Thus, in step 3 to get the row echelon form, we need 
$$\sum_{i=1}^{m-1}\left[ 2(m-i)(n-i+1) +(m-i)\right] =\sum_{i=1}^{m-1} \left[2i^2-(2m+2n+3)i+(2mn+3m)\right] 
$$
flops
or \underline{$m^2n-\frac{1}{3}m^3$} flops if keep only the leading terms.

\paragraph{Get REF, case 2: $m>n$.} The procedure should stop at loop $n$, otherwise, $n-k+1$ will be zero if $k=n+1$.  Thus, in step 3 to get the row echelon form, we need 
$$
\sum_{i=1}^{n} \left[2(m-i)(n-i+1) +(m-i)\right]=\sum_{i=1}^{n} \left[2i^2-(2m+2n+3)i+(2mn+3m)\right]
$$
flops or \underline{$mn^2-\frac{1}{3}n^3$} flops if keep only the leading terms.

\begin{figure}[H]
\centering  
\vspace{-0.35cm} 
\subfigtopskip=2pt 
\subfigbottomskip=2pt 
\subfigcapskip=-5pt 
\subfigure[RREF for $m\leq n$.]{\label{fig:cr-crcomputation1}
\includegraphics[width=0.4\linewidth]{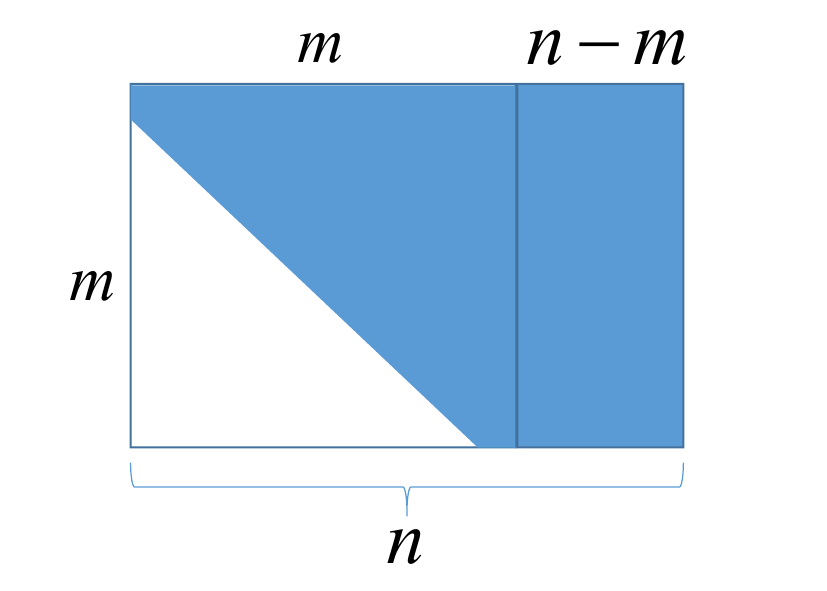}}
\quad 
\subfigure[RREF for $m> n$.]{\label{fig:cr-crcomputation2}
\includegraphics[width=0.4\linewidth]{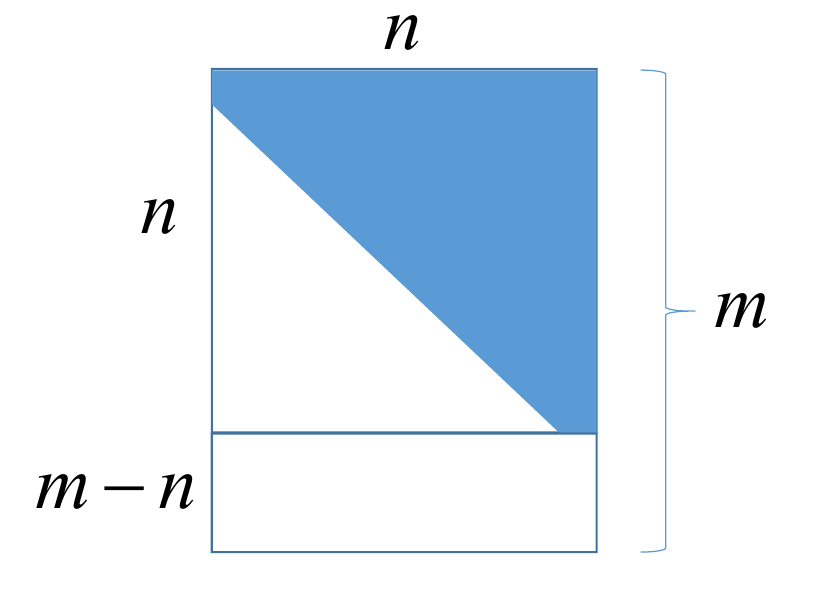}}
\caption{Get RREF from REF, where blank entries indicate zeros, blue entries indicate values that are not necessarily zero, pivots are positioned on the diagonal in the ``worst" case. }
\label{fig:cr-crcomputation12}
\end{figure}
To obtain the reduced row  echelon form (RREF) from the row echelon form (REF), we can also consider two cases: $m\leq n$ and $m>n$, as shown in Figure~\ref{fig:cr-crcomputation12}, where blank entries indicate zeros, blue entries indicate values that are not necessarily zero, pivots are positioned on the diagonal in the ``worst" case.  

\paragraph{Get RREF, case 1: $m\leq n$.} 

Procedure (D) requires $(m-1)$ divisions to obtain the multipliers, followed by $(m-1)\times (n-m+1)$ multiplications with the multipliers (we operate on $(m-1)$ length-$(n-m+1)$ vectors). To zero out the values  above the pivot, $(m-1)\times (n-m+1)$ subtractions are involved. As a result, procedure (D) incurs a cost of \underline{$2(m-1)(n-m+1)+m-1$} flops.

Similarly, procedure (E) requires $(m-2)$ divisions to compute the multipliers, followed by $(m-2)\times (n-m+2)$ multiplications with these multipliers (we operate on $(m-2)$ length-$(n-m+2)$ vectors). 
To zero out the values above the pivot, it involves $(m-2)\times (n-m+2)$ subtractions. Consequently, procedure (E) costs \underline{$2(m-2)(n-m+2)+m-2$} flops.

The procedure can continue, and we can  again summarize the cost for each loop in the following table:
\begin{table}[H]
\setlength{\tabcolsep}{7.6pt}
\begin{tabular}{l|l|l|l}
$k$      & Get multipliers & Multipliers multiply each row & Rows subtraction   \\ \hline
$m$      & $1:m-1=m-1$ rows      & $(m-1)(n-m+1)$                & $(m-1)(n-m+1)$     \\
$m-1$    & $1:m-2=m-2$ rows      & $(m-2)(n-m+2)$                & $(m-2)(n-m+2)$     \\
$m-2$    & $1:m-3=m-3$ rows      & $(m-3)(n-m+3)$                & $(m-3)(n-m+3)$     \\
$\vdots$ & $\vdots$            & $\vdots$                      & $\vdots$           \\
$k$      & $1:k-1=k-1$ rows    & $(k-1)(n-k+1)$            & $(k-1)(n-k+1)$ \\
$\vdots$ & $\vdots$            & $\vdots$                      & $\vdots$           \\
$2$      & $1:1=1$ row        & $(1)(n-2+1)$                  & $(1)(n-2+1)$    
\end{tabular}
\end{table}


Therefore, when $m\leq n$, in the second loop of the algorithm to get the reduced row  echelon form from the row echelon form, we require 
$$
\begin{aligned}
\sum_{k=2}^{m} \left[2(k-1)(n-k+1)+k-1\right]
&\stackrel{(i\triangleq k-1)}{=}\sum_{i=1}^{m-1} \left[2(i)(n-i)+i\right]
= \sum_{i=1}^{m-1} \left[-2i^2+(2n+1)i\right]
\end{aligned}
$$ 
flops
or \underline{$-\frac{2}{3}m^3+m^2n$} flops if we keep only the leading terms.

\paragraph{Get RREF, case 2: $m>n$.} 
Procedure (D) requires $(n-1)$ divisions to get the multipliers, followed by $(n-1)\times 1$ multiplications with the multipliers (we operate on $(n-1)$ length-1 vectors). To zero out the values above the pivot, it involves $(n-1)\times 1$ subtractions. As a result, procedure (D) costs \underline{$2(n-1)(1)+n-1$} flops.

Procedure (E) needs $(n-2)$ divisions to get the multipliers, followed by $(n-2)\times 2$ multiplications with the multipliers (we operate on $(n-2)$ length-2 vectors). To zero out the values above the pivot, it involves $(n-2)\times 2$ subtractions. As a result, procedure (E) costs \underline{$2(n-2)(2)+n-2$} flops.

The procedure can go on, and we can  again summarize the cost for each loop in the following table:
\begin{table}[H]
\setlength{\tabcolsep}{8.9pt}
\begin{tabular}{l|l|l|l}
$k$      & Get multipliers & Multipliers multiply each row & Rows subtraction   \\ \hline
$n$      & $1:n-1=n-1$ rows      & $(n-1)(1)$                & $(n-1)(1)$      \\
$n-1$    & $1:n-2=n-2$ rows      & $(n-2)(2)$                & $(n-2)(2)$        \\
$n-2$    & $1:n-3=n-3$ rows      & $(n-3)(3)$                & $(n-3)(3)$           \\
$\vdots$ & $\vdots$            & $\vdots$                      & $\vdots$           \\
$k$      & $1:k-1=k-1$ rows    & $(k-1)(n-k+1)$            & $(k-1)(n-k+1)$ \\
$\vdots$ & $\vdots$            & $\vdots$                      & $\vdots$           \\
$2$      & $1:1=1$ row        & $(1)(n-2+1)$                  & $(1)(n-2+1)$    
\end{tabular}
\end{table}

Therefore, when $m>n$, during the second loop to obtain the reduced row  echelon form from the row echelon form, we require
$$
\begin{aligned}
\sum_{k=2}^{n} \left[2(k-1)(n-k+1)+k-1\right]
&\stackrel{i\triangleq k-1}{=}\sum_{i=1}^{n-1} \left[2(i)(n-i)+i\right]
= \sum_{i=1}^{n-1} \left[-2i^2+(2n+1)i\right]
\end{aligned}
$$ 
flops 
or \underline{$\frac{1}{3}n^3$} flops if we keep only the leading terms.

\paragraph{Total cost.}
In step 8, there are a maximum of $mn$ flops (divisions) required to set the pivots to 1. Thus, the total cost can be summarized as follows:
\begin{itemize}
\item The total cost for $m\leq n$ is  $m^2n-\frac{1}{3}m^3-\frac{2}{3}m^3+m^2n = $\underline{$2m^2n-m^3$} flops if we keep only the leading terms;

\item The total cost for $m> n$ is  $mn^2-\frac{1}{3}n^3+\frac{1}{3}n^3 = $\underline{$mn^2$} flops if we keep only the leading term.
\end{itemize}
If $m=n$, the cost is $m^3$ flops.
And this completes the proof.
\end{proof}

In the skeleton decomposition (Section~\ref{section:skeleton_decom}), we will use the  Gram-Schmidt process to find the linearly independent columns of $\bA$, which has a similar complexity to the Gaussian elimination used for finding the reduced row  echelon form. However, the Gram-Schmidt process is clearer from the perspective of the equations to be computed, meaning we can easily write out specific mathematical forms for the entries of the matrices.

\section{Rank Decomposition}
We previously mentioned that the CR decomposition is a special case of rank decomposition. Formally, we prove the existence of the rank decomposition rigorously in the following theorem.

\begin{theoremHigh}[Rank Decomposition]\label{theorem:rank-decomposition}
Let $\bA \in \real^{m \times n}$ be any rank-$r$ matrix. Then it can be factored as the following \textit{rank decomposition} (a.k.a., \textit{full-rank decomposition}):
$$
\underset{m\times n}{\bA }= \underset{m\times r}{\bD}\gapthree  \underset{r\times n}{\bF},
$$
where both $\bD \in \real^{m\times r}$ and $\bF \in \real^{r\times n}$  have (full) rank $r$. 
The storage for the decomposition is then reduced or potentially increased from $mn$ floating-point numbers to $r(m+n)$ floating-point numbers.
\end{theoremHigh}
\begin{proof}[of Theorem~\ref{theorem:rank-decomposition}]
According to the ULV decomposition in Theorem~\ref{theorem:ulv-decomposition}, we can decompose $\bA$ as follows: 
$$
\bA = \bU \begin{bmatrix}
	\bL & \bzero \\
	\bzero & \bzero 
\end{bmatrix}\bV.
$$
Let $\bU_0 \triangleq \bU[:,1:r]$ and $\bV_0 \triangleq \bV[1:r,:]$, i.e., $\bU_0$ comprises only the first $r$ columns of $\bU$, and $\bV_0$ includes only the first $r$ rows of $\bV$. 
Thus, we can still express $\bA$ as  $\bA = \bU_0 \bL\bV_0$, where $\bU_0 \in \real^{m\times r}$ and $\bV_0\in \real^{r\times n}$. This is also referred to as the reduced ULV decomposition, as shown in Figure~\ref{fig:ulv-comparison}. Let \{$\bD \triangleq \bU_0\bL$ and $\bF \triangleq\bV_0$\}  or \{$\bD \triangleq \bU_0$ and $\bF \triangleq\bL\bV_0$\}, we find such  a rank decomposition.
\end{proof}

The rank decomposition is not unique. Even by elementary transformations, we have 
$$
\bA = 
\bE_1
\begin{bmatrix}
	\bZ & \bzero \\
	\bzero & \bzero 
\end{bmatrix}
\bE_2,
$$
where $\bE_1 \in \real^{m\times m}$ and $\bE_2\in \real^{n\times n}$ represent products of nonsingular elementary row and column operations, and $\bZ\in \real^{r\times r}$. The transformation is rather general, and  there are numerous such $\bE_1,\bE_2,$ and $\bZ$ matrices. 
When $\bZ=\bI_r$, where $r$ is the rank of matrix $\bA$, the decomposition is known as the \textit{Smith decomposition or Smith form} of $\bA$ \citep{bernstein2009matrix}.
With a construction similar to that demonstrated in the previous proof, we can derive an alternative rank decomposition.
\index{Decomposition: Smith}
\index{Smith decomposition}

Analogously, we can obtain such $\bD$ and $\bF$ using SVD, URV, CR, CUR, and numerous other decompositional algorithms. However, we may establish a connection between these various rank decompositions through the following lemma.
\begin{corollary}[Connection Between Rank Decompositions]\label{corollary:connection-rank-decom}
Let $\bA=\bD_1\bF_1=\bD_2\bF_2\in\real^{m\times n}$ be  two rank decompositions of $\bA$. Then there exists a nonsingular matrix $\bP$ such that 
$$
\bD_1 = \bD_2\bP
\qquad
\text{and}
\qquad 
\bF_1 = \bP^{-1}\bF_2.
$$
More generally, given $\bA,\bB\in\real^{m\times n}$, $\bA$ and $\bB$ are biequivalent (Definition~\ref{definition:biequivalent}) if and only if $\bA$ and $\bB$ have the same Smith form.
\end{corollary}
\begin{proof}[of Corollary~\ref{corollary:connection-rank-decom}]
Since $\bD_1\bF_1=\bD_2\bF_2$, we have $\bD_1\bF_1\bF_1^\top=\bD_2\bF_2\bF_1^\top$. It is trivial that $\rank(\bF_1\bF_1^\top)=\rank(\bF_1)=r$ such that $\bF_1\bF_1^\top$ is a square matrix with full rank and thus is nonsingular. This implies $\bD_1=\bD_2\bF_2\bF_1^\top(\bF_1\bF_1^\top)^{-1}$. Let $\bP\triangleq\bF_2\bF_1^\top(\bF_1\bF_1^\top)^{-1}$, we have $\bD_1=\bD_2\bP$ and $\bF_1 = \bP^{-1}\bF_2$. The second part can be proved similarly.
\end{proof}

To conclude this section, we present an alternative form of the rank decomposition that will facilitate connections to tensor decomposition (Theorem~\ref{theorem:cp-decomp}).
\begin{theoremHigh}[Rank Decomposition, An Alternative Form]\label{theorem:rank-decomposition-alternative}
Let $\bA \in \real^{m \times n}$  be any rank-$r$ matrix. Then it can be factored as 
$$
\underset{m\times n}{\bA }= \underset{m\times r}{\bD}\gapthree  \underset{n\times r}{\bE^\top} = \sum_{i=1}^{r} \bd_i \be_i^\top,
$$
where $\bD=[\bd_1, \bd_2, \ldots, \bd_r] \in \real^{m\times r}$ is of rank $r$, and $\bE=[\be_1, \be_2, \ldots, \be_r] \in \real^{n\times r}$ also has rank $r$, i.e., both $\bD$ and $\bE$ have full rank $r$. 
\end{theoremHigh}

\section{Application: Rank and Trace of an Idempotent Matrix}
The CR decomposition serves as a powerful tool for understanding the rank characteristics of idempotent matrices. 
Additionally, its utility in orthogonal projections is discussed in Appendix~\ref{appendix:orthogonal}.
\begin{proposition}[Rank and Trace of an Idempotent Matrix\index{Trace}]\label{proposition:rank-of-symmetric-idempotent2_tmp}
Let $\bA$ be idempotent (i.e., $\bA^2 = \bA$). Then its rank is equal to its trace.
\end{proposition}
\begin{proof}[of Proposition~\ref{proposition:rank-of-symmetric-idempotent2_tmp}]
For any $n\times n$ matrix $\bA$ with rank $r$, the CR decomposition exists as $\bA = \bC\bR$, where $\bC\in\real^{n\times r}$ and $\bR\in \real^{r\times n}$, both with a full rank of $r$. Consequently,
$$
\begin{aligned}
\bA^2 = \bA
\quad\implies \quad
\bC\bR\bC\bR = \bC\bR
\quad\implies \quad 
\bR\bC\bR =\bR
\quad\implies \quad
\bR\bC =\bI_r,
\end{aligned}
$$ 
where $\bI_r$ denotes the  $r\times r$ identity matrix. Therefore,
$$
\trace(\bA) = \trace(\bC\bR) =\trace(\bR\bC) = \trace(\bI_r) = r, 
$$
which corresponds to the rank of $\bA$. This equality holds due to the invariance of the trace under cyclic permutations.
\end{proof}

\index{Fundamental theorem}
\section{Other Applications}
The CR decomposition or rank decomposition plays a pivotal role in proving numerous essential theorems.
For instance, it is instrumental in establishing the existence of the pseudo-inverse, as demonstrated in Theorem~\ref{theorem:existence-of-pseudo-inverse}.
Furthermore, it plays a key role in determining the basis of the four fundamental subspaces in linear algebra, as discussed in Appendix~\ref{appendix:cr-decomposition-four-basis}. 

Moreover, the CR factorization finds applications in data interpretation and solving computational problems. For instance, it proves valuable in solving  least squares problems, enabling the removal of redundant variables by considering a reduced linear system. Throughout this book, readers will encounter various applications of the CR decomposition  across different contexts.

\begin{problemset}
\item Discuss what rank number $r$ in Theorem~\ref{theorem:cr-decomposition} and Theorem~\ref{theorem:rank-decomposition} can reduce storage requirements.
\item Determine the reduced row echelon form and the CR decomposition for the matrix 
$$
\bA = 
\begin{bmatrix}
1 & 3 & 2 \\
3 & 7 & 6 \\
4 & 5 & 8
\end{bmatrix}.
$$

\item Apply the RREF to the matrix
$$
\bA = \begin{bmatrix}
	1 & 2 & 1 & 1 \\
	1 & 4 & 2 & 3 \\
	1 & 1 & 2 & -1 \\
	-3 & -1 & 4 & 0
\end{bmatrix}.
$$

\item \label{prob:cr_pseudo} Find the pseudo-inverse of a matrix $\bA$ using its CR decomposition (refer to Appendix~\ref{section:pseudo_cr}).

\item \label{problem:aug_lin} Show that the solution of the linear system $\bA\bx=\bb$ remains unchanged if the same sequence of elementary row transformations  is applied to both $\bA$ and $\bb$. Thus, the solution is revealed by the RREF of the augmented matrix $[\bA,\bb]$.

\item Following Problem~\ref{problem:aug_lin}, show that the two linear systems $\bA_1\bx=\bb_1$ and $\bA_2\bx=\bb_2$ have the same set of solutions if and only if $[\bA_1,\bb_1]$ and $[\bA_2,\bb_2]$ have the same RREF.

\item Show that if a system of linear equations has two distinct solutions, then it must have infinitely many solutions.
\item Show that if a linear system $\bA\bx=\bb$ has more than one solution, then system $\bA\bx=\bzero$ holds. 
\item A system of linear equations with fewer equations than unknowns is sometimes referred to as an \textit{underdetermined} system. Provide an example of an inconsistent underdetermined system of two equations in three unknowns  (Remark~\ref{remark:into_lin_syssol}).

\item Suppose an underdetermined system is consistent (Remark~\ref{remark:into_lin_syssol}). Explain why such a system must have an infinite number of solutions.

\item A system of linear equations with more equations than unknowns is sometimes referred to as an \textit{overdetermined} system. Discuss the conditions under which such a system can be consistent.

\item Two matrices are called \textit{row equivalent} if there is a sequence of elementary row operations that transforms one matrix into the other. Show that if matrices $\bA$ and $\bB$ are row equivalent, they have the same RREF.

\item Let $\bA\bx=\bb$ be a consistent system with $\bA\in\real^{m\times n}$. Show that $\bA$ has $m$ pivot columns. Let further $m=n$, show that the RREF of $\bA$ is the identity matrix.

\item 
Let $ \bA=\footnotesize\begin{bmatrix}
	a & b \\
	c & d
\end{bmatrix}$ be any $ 2 \times 2 $  nonsingular matrix. 
Show that there exists an nonsingular matrix $\bS$ such that
$$
\bS\bA = \begin{bmatrix}
	1 & 0 \\
	0 & ad - bc
\end{bmatrix},
$$
where $\bS$ is the product of at most four elementary matrices of the form $ \bE_{i,j;\alpha} $ in \eqref{equation:ele_exam_eijal}.
More generally, let $ \bA\in\real^{n\times n} $ be nonsingular. Show that there is a matrix $\bS$ such that
$$
\bS\bA = \begin{bmatrix}
	\bI_{n-1} & 0 \\
	0 & d
\end{bmatrix},
$$
where $d = \det(\bA) $, and  $\bS $ is again a product of elementary matrices of the form $ \bE_{i,j;\alpha} $.

\end{problemset}

%% file: chapter-skeleton.tex
\newpage
\chapter{Skeleton/CUR Decomposition}
\begingroup
\hypersetup{
	linkcolor=structurecolor,
	linktoc=page,  
}
\minitoc \newpage
\endgroup

\section{Skeleton Decomposition}\label{section:skeleton_decom}
\lettrine{\color{caligraphcolor}T}
The CR decomposition reuses the actual columns of the underlying matrix, whereas the skeleton decomposition extends this  by using both actual columns and rows of the matrix.

\index{Decomposition: Skeleton}
\index{Decomposition: CUR}
\begin{theoremHigh}[Skeleton Decomposition]\label{theorem:skeleton-decomposition}
Let $\bA \in \real^{m \times n}$ be any rank-$r$ matrix. Then it can be factored as 
$$
\underset{m\times n}{\bA }= 
\underset{m\times r}{\bC} \gapthree \underset{r\times r}{\bU^{-1} }\gapthree \underset{r\times n}{\bR},~\footnote{The nonsingular $\bU$ extends a terminology called \textit{rank principal} for a square matrix $\bB\in\real^{n\times n}$. $\bB$, with rank $r$, is called rank principal if it has a nonsingular $r\times r$ principal submatrix (Definition~\ref{definition:principle-minors}). 
That is, there exists an index set $I\subseteq \{1,2,\ldots,n\}$ with cardinality $r$ such that $\rank(\bB)=\rank(\bB[I,:])=\rank(\bB[:,J])=\rank(\bB[I,J])=r$.
See Problem~\ref{problem:rank_pcin}.
}
$$
where $\bC$ contains some $r$ linearly independent columns of $\bA$, $\bR$ contains some $r$ linearly independent rows of $\bA$, and $\bU$ is the \textbf{nonsingular} submatrix at the intersection of $\bC$ and $\bR$. 
\begin{itemize}
\item The storage required for the decomposition is  reduced or, in some cases, potentially increased from $mn$ floating-point numbers to $r(m+n)+r^2$ floating-point numbers. 
\item Alternatively, if we only record the positions of the indices, it would require $mr$ and $nr$ floating-point numbers for storing $\bC$ and $\bR$, respectively. Additionally, it would necessitate an  extra $2r$ integers to remember the position of each column of $\bC$ in that of $\bA$ and each row of $\bR$ in that of $\bA$ (i.e., to construct $\bU$ from $\bC$ and $\bR$).
\end{itemize}
\end{theoremHigh}

\begin{figure}[H]
\centering
\includegraphics[width=0.7\textwidth]{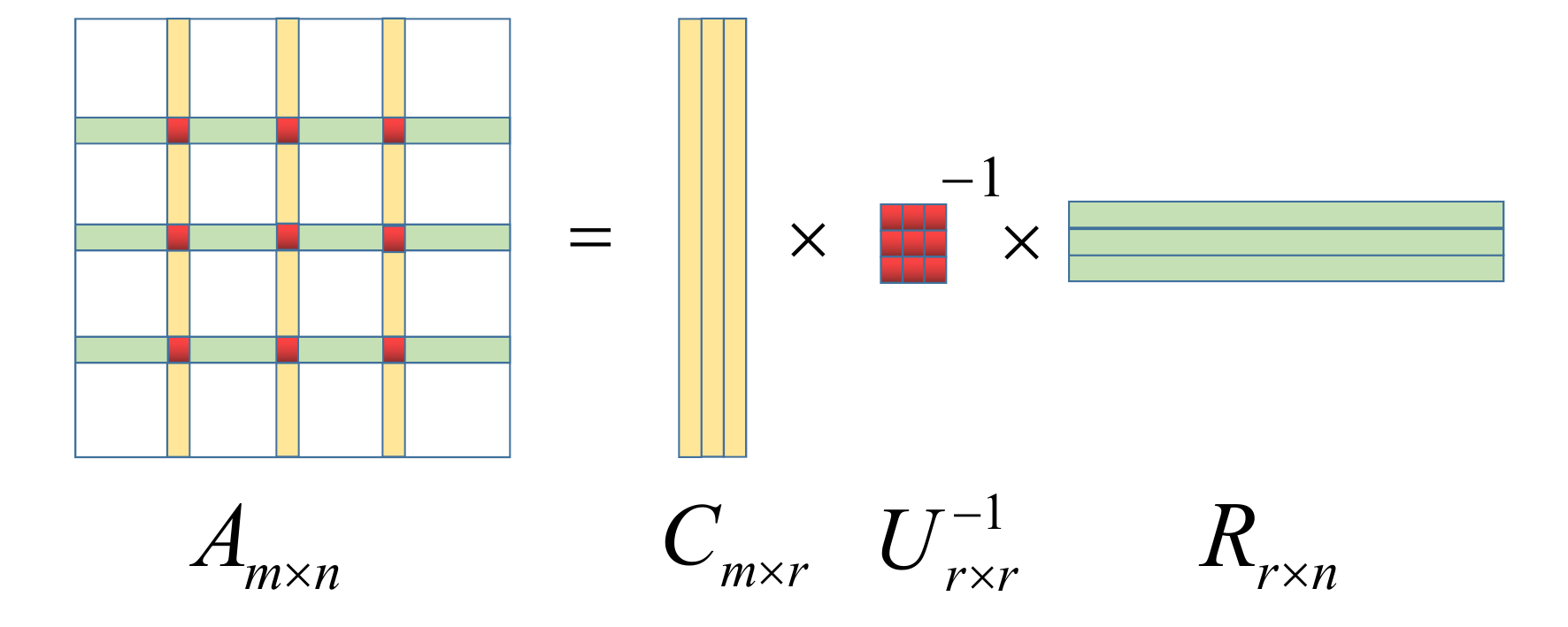}
\caption{Illustration of the skeleton decomposition of a matrix, where the \textcolor{mydarkyellow}{yellow} vectors denote the linearly independent columns of $\bA$, and \textcolor{mydarkgreen}{green} vectors denote the linearly independent rows of $\bA$.}
\label{fig:skeleton}
\end{figure}

The skeleton decomposition is also  referred to as the \textit{CUR decomposition} following from its notation in the decomposition. 
Compared to SVD, CUR offers advantages in terms of reification and interpretability  issues since it uses the actual columns (rows) of the matrix, whereas SVD uses some artificial singular vectors that may not accurately represent physical reality \citep{mahoney2009cur}. 
Moreover,  CUR exhibits sparsity preservation when the underlying data is sparse, making it particularly suitable for applications involving sparse matrices.

Conversely, similar to SVD, the CUR decomposition serves as a valuable tool in various application domains for tasks such as data compression, feature extraction, and data analysis and is a versatile tool for approximating matrices with reduced computational complexity, making it suitable for handling large-scale data in various domains \citep{mahoney2009cur, an2012large, lee2008cur+}. 
For example, CUR is used to compress large datasets by retaining only a fraction of the original data (columns and rows) while maintaining the essential information, providing a low-rank approximation useful in various numerical linear algebra applications, such as solving linear systems, eigenvalue problems, and matrix inversion. It also provides tools for image compression and analysis, where the original image matrix can be approximated with a lower-dimensional representation.
In machine learning, CUR can be used for dimensionality reduction, feature extraction, and data representation, which can improve the efficiency of algorithms and reduce computational costs.
In collaborative filtering for recommendation systems, CUR can be employed to approximate large user-item interaction matrices, making the system more scalable.
CUR decomposition can also be extended to higher-dimensional arrays (tensors) for applications in multi-linear algebra and data analysis, enabling the processing of complex datasets \citep{kishore2017literature}.

An illustration of the skeleton decomposition is shown in Figure~\ref{fig:skeleton}, where  \textcolor{mydarkyellow}{yellow} vectors denote the linearly independent columns of $\bA$, and \textcolor{mydarkgreen}{green} vectors denote the linearly independent rows of $\bA$. Specifically, given  index vectors $I$ and $J$, both of size $r$, containing the indices of rows and columns selected from $\bA$ to form $\bR$ and $\bC$, respectively, the submatrix $\bU$ can be expressed as $\bU=\bA[I,J]$ (refer to Definition~\ref{definition:matlabnotation}).

\section{Existence of  Skeleton Decomposition}
In Theorem~\ref{theorem:equal-dimension-rank}, we have demonstrated that the row rank and the column rank of a matrix are equal.
In other words, we can state that the dimensions of the column space and the row space are identical. This property is crucial for establishing the existence of the skeleton decomposition.

We are now prepared to demonstrate the existence of the skeleton decomposition. The proof is rather elementary.
\begin{proof}[of Theorem~\ref{theorem:skeleton-decomposition}]
The proof hinges on the existence of such a nonsingular matrix $\bU$, which is central to this decomposition method. 

\paragraph{Existence of such a nonsingular matrix $\bU$.} Since the matrix $\bA$ is of rank $r$, we can select $r$ linearly independent columns from $\bA$.
Suppose we place these  $r$ linearly independent columns $\ba_{i1}, \ba_{i2}, \ldots, \ba_{ir}$ into the columns of an $m\times r$ matrix $\bN=[\ba_{i1}, \ba_{i2}, \ldots, \ba_{ir}] \in \real^{m\times r}$. By Theorem~\ref{theorem:equal-dimension-rank}, the dimension of the column space of $\bN$ is $r$, and thus the dimension of the row space of $\bN$ is also $r$. Consequently, we can choose $r$ linearly independent rows $\bn_{j1}^\top,\bn_{j2}^\top, \ldots, \bn_{jr}^\top $ from $\bN$ and place the specific $r$ rows into the rows of an $r\times r$ matrix $\bU = [\bn_{j1}^\top; \bn_{j2}^\top; \ldots; \bn_{jr}^\top]\in \real^{r\times r}$.
Applying Theorem~\ref{theorem:equal-dimension-rank} once again, the dimension of the column space of $\bU$ is also $r$, which means there are  $r$ linearly independent columns within $\bU$.
Therefore, $\bU$ is such a nonsingular matrix of size $r\times r$. (see Problem~\ref{problem:nons_sub} for more insights of the rank of a submatrix.)

\paragraph{Main proof.}
Given the presence of a nonsingular $r\times r$ matrix $\bU$ within $\bA$, the existence of the skeleton decomposition can be established as follows.
Suppose $\bU=\bA[I,J]$, where $I$ and $J$ are index vectors of size $r$. Since $\bU$ is a nonsingular matrix, the columns of $\bU$ are linearly independent. Thus, the columns of matrix $\bC$, based on the columns of $\bU$, are also linearly independent (i.e., select the $r$ columns of $\bA$ that correspond to the entries of    $\bU$. Here $\bC$ is equivalent to the matrix $\bN$ we construct previously, denoted by $\bC=\bA[:,J]$).

Given that the rank of  $\bA$ is $r$, for any  column $\ba_i$ in $\bA$, $\ba_i$ can be represented as a linear combination of the columns of $\bC$. In other words, there exists a vector $\bx$ such that  $\ba_i = \bC \bx$ for all $ i\in \{1, 2, \ldots, n\}$. 
Let the $r$ rows (entries) of $\ba_i\in\real^n$ corresponding to the row entries of $\bU$ be $\br_i \in \real^r$ for all $i\in \{1, 2, \ldots, n\}$ (i.e., $\br_i$ contains $r$ entries of $\ba_i$). That is, select the $r$ entries of $\ba_i$'s corresponding to the entries of $\bU$ as follows:
$$
\bA = [\ba_1,\ba_2, \ldots, \ba_n]\in \real^{m\times n} \qquad \implies \qquad
\bA[I,:]=[\br_1, \br_2, \ldots, \br_n] \in \real^{r\times n}.
$$
Since $\ba_i = \bC\bx$, $\bU$ is a submatrix inside $\bC$, and $\br_i$ is a subvector inside $\ba_i$, we have $\br_i = \bU \bx$, which states that $\bx = \bU^{-1} \br_i$. 
Therefore, for each $i\in\{1,2,\ldots,n\}$, we can express $\ba_i = \bC \bU^{-1} \br_i$. 
By combining the $n$ columns of such $\br_i$ into $\bR=[\br_1, \br_2, \ldots, \br_n]$, we obtain
$$
\bA = [\ba_1, \ba_2, \ldots, \ba_n] = \bC \bU^{-1} \bR,
$$
from which the result follows.

In summary, we start by identifying $r$ linearly independent columns of $\bA$ and placing them into $\bC\in \real^{m\times r}$. Within $\bC$, we find an $r\times r$ nonsingular submatrix $\bU$. 
The $r$ rows of $\bA$ corresponding to the entries of $\bU$ can help to reconstruct the columns of $\bA$. Again, this process is visually illustrated in Figure~\ref{fig:skeleton}.
\end{proof}

In the case where $\bA$ is square and invertible,  the skeleton decomposition simplifies to $\bA=\bC\bU^{-1} \bR$, where $\bC=\bR=\bU=\bA$ such that the decomposition essentially becomes $\bA = \bA\bA^{-1}\bA$.

\paragraph{CR decomposition vs skeleton decomposition.} We note that the CR decomposition and the skeleton decomposition share a similar form, even in terms of the symbols employed--$\bA=\bC\bR$ for the CR decomposition and $\bA=\bC\bU^{-1}\bR$ for the skeleton decomposition. 

In both the CR decomposition and the skeleton decomposition, we \textbf{can}~\footnote{Here, we emphasize that we ``can," but not necessarily, as we will see in the randomized algorithm for the skeleton decomposition.} select the first $r$ linearly independent columns to form the matrix $\bC$ (the symbol for both the CR decomposition and the skeleton decomposition). Thus, in this scenario, the symbols of $\bC$ in the CR decomposition and  skeleton decompositions are exactly the same. 

However, there is a distinction between the two methods when it comes to $\bR$. 
In the CR decomposition, $\bR$ represents the reduced row echelon form without zero rows, whereas in the skeleton decomposition, $\bR$ consists of specific rows from $\bA$, giving $\bR$ different meanings in these two approaches. 
While $\bU^{-1}\bR$ in the skeleton decomposition is the reduced row  echelon form without zero rows, provided that we select the first $r$ linearly independent columns in the skeleton algorithm (see Algorithm~\ref{alg:skeleton-decomposition}, Algorithm~\ref{alg:skeleton-decomposition-modified}, or Algorithm~\ref{alg:skeleton-by-row-reduced-echelon}).

\paragraph{A word on the uniqueness of CR decomposition and skeleton decomposition.} As previously mentioned, in both the CR and  skeleton decompositions, we can select the first $r$ linearly independent columns to obtain the matrix $\bC$. From this perspective, the CR and skeleton decompositions exhibit a unique form. 
However, if we select the last $r$ linearly independent columns, we will obtain a different CR decomposition or skeleton decomposition.

To reiterate, in the proof outlined above for the existence of the skeleton decomposition, we initially identify  $r$ linearly independent columns of $\bA$ and store them in the matrix $\bC\in\real^{m\times r}$.  
From $\bC$, we extract an $r\times r$ nonsingular submatrix $\bU$. 
Finally, we identify the $r$ rows of $\bA$ that correspond to the entries of $\bU$ to form the row matrix $\bR\in\real^{r\times n}$. 
A further question can be posed that if matrix $\bA$ has rank $r$, matrix $\bC$ contains $r$ linearly independent columns of $\bA$, and matrix $\bR$ contains $r$ linearly independent rows of $\bA$, is the $r\times r$ ``intersection" of $\bC$ and $\bR$ invertible or not \footnote{We express our gratitude to Gilbert Strang for posing this  question.}?

\begin{corollary}[Nonsingular Intersection]\label{corollary:invertible-intersection}
Let $\bA \in \real^{m\times n}$ be any matrix that has rank $r$, and let matrix $\bC$ contain $r$ linearly independent columns of $\bA$, and matrix $\bR$ contain $r$ linearly independent rows of $\bA$. Then, the $r\times r$ ``intersection" matrix $\bU$ of $\bC$ and $\bR$ is invertible.
\end{corollary}
\begin{proof}[of Corollary~\ref{corollary:invertible-intersection}]
Let $I$ and $J$ represent the indices of the rows and columns selected from $\bA$ to form $\bR$ and $\bC$, respectively, then $\bR$ can be expressed as $\bR=\bA[I, :]$, $\bC$ can be represented as $\bC = \bA[:,J]$, and $\bU$ can be denoted by $\bU=\bA[I,J]$.

Since $\bC$ comprises $r$ linearly independent columns of $\bA$, any column $\ba_i$ from $\bA$ can be expressed as $\ba_i = \bC\bx_i = \bA[:,J]\bx_i$ for all $i \in \{1,2,\ldots, n\}$. 
This implies that the $r$ entries of $\ba_i$ corresponding to the $I$ indices can be represented by the columns of $\bU$ such that $\ba_i[I] = \bU\bx_i \in \real^{r}$ for all $i \in \{1,2,\ldots, n\}$, i.e.,
$$
\ba_i = \bC\bx_i = \bA[:,J]\bx_i \in \real^{m} \qquad \implies  \qquad
\ba_i[I] =\bA[I,J]\bx_i= \bU\bx_i \in \real^{r}.
$$ 
Given that $\bR$ comprises $r$ linearly independent rows of $\bA$, both the row rank and column rank of $\bR$ are equal to $r$. Combining these facts, the $r$ columns of $\bR$ corresponding to the indices in $J$ (i.e., the $r$ columns of $\bU$) are also linearly independent. 

By applying Theorem~\ref{theorem:equal-dimension-rank} again, the dimension of the row space of $\bU$ is also equal to $r$, which means there are  $r$ linearly independent rows within $\bU$, and thus $\bU$ is invertible.
\end{proof}

\begin{algorithm}[h] 
\caption{Skeleton Decomposition via Gram-Schmidt Process} 
\label{alg:skeleton-decomposition} 
\begin{algorithmic}[1] 
\Require 
Rank-$r$ matrix $\bA=[\ba_1, \ba_2, \ldots, \ba_n]$ with size $m\times n $; 
\State Initially set column count $ck=0$ and row count $rk=0$; 
\Statex \textbf{Stage A: Obtain independent columns into $\bC$}
\State Set $\bq_1   = \ba_1/r_{11},  r_{11}=\normtwo{\ba_1}$; \Comment{Suppose the first column is nonzero and $r_{11}\neq 0$}
\For{$k=2$ to $n$} 
\State \algoalign{$\bq_k = \frac{(\ba_k-\sum_{i=1}^{k-1}r_{ik}\bq_i)}{r_{kk}},\,
r_{ik}=\bq_i^\top\ba_k, \,
r_{kk}=\normtwo{\ba_k-\sum_{i=1}^{k-1}r_{ik}\bq_i}$, $  \forall i \in \{1,\ldots, k-1\}$;}
\If{$r_{kk}\neq 0$}
\State 	Select the $k$-th column of $\bA$ into $ck$-th column of $\bC$;
\State  $ck=ck+1$;
\EndIf
\EndFor
\Statex \textbf{Stage B: Obtain independent rows into $\bU$ (and $\bR$)}
\State Suppose $\bC=[\bc_1^\top; \bc_2^\top; \ldots; \bc_m^\top]$, where $\bc_i^\top$ is the $i$-th row of $\bC$;
\State Set $\bq_1   = \bc_1/r_{11},  r_{11}=\normtwo{\bc_1}$;\Comment{Suppose the first row is nonzero and $r_{11}\neq 0$}
\For{$k=2$ to $m$} 
\State \algoalign{$\bq_k = \frac{(\bc_k-\sum_{i=1}^{k-1}r_{ik}\bq_i)}{r_{kk}}, \,
r_{ik}=\bq_i^\top\bc_k,\,  r_{kk}=\normtwo{\bc_k-\sum_{i=1}^{k-1}r_{ik}\bq_i} $, $ 
\forall i \in \{1,\ldots, k-1\}$;}
\If{$r_{kk}\neq 0$}
\State 	Select the $k$-th row of $\bC$ into $rk$-th row of $\bU$;
\State  $rk=rk+1$;
\EndIf
\EndFor
\State Select the rows from $\bA$ corresponding to the rows of $\bU$ into $\bR$.
\end{algorithmic} 
\end{algorithm}

\index{Gram–Schmidt}
\section{Computing  Skeleton Decomposition via  Gram-Schmidt Process}

In Section~\ref{section:dependent-gram-schmidt-process}, we  discussed how to handle linearly dependent columns when orthogonalizing a matrix. 
We can thus apply the Gram-Schmidt process to extract $r$ linearly independent columns from $\bA$, forming the matrix $\bC$. Furthermore, we can use the Gram-Schmidt process to select $r$ linearly independent columns from $\bC^\top$, resulting in the matrices $\bU$ and $\bR$. The process is outlined in Algorithm~\ref{alg:skeleton-decomposition}.
While the Gram-Schmidt process aims to find orthonormal vectors, the orthogonal property is not essential for the skeleton decomposition. Instead, we use the orthogonal property to determine the linear independence of the vectors.

\begin{theorem}[Algorithm Complexity: Skeleton via Gram-Schmidt]\label{theorem:comp-skeleton-decom1}
Algorithm~\ref{alg:skeleton-decomposition} requires $\sim 2(mn^2+rm^2+r^3)$ flops to compute a skeleton decomposition of an $m\times n$ matrix. We will show that this algorithm is not the most efficient way to obtain the skeleton decomposition.
\end{theorem}

\begin{proof}[of Theorem~\ref{theorem:comp-skeleton-decom1}]
This involves applying Theorem~\ref{theorem:qr-reduced} twice to obtain matrix $\bU$: once to get $\bC$ (which costs $2mn^2$ flops) and once to get $\bU$ (which costs $2rm^2$ flops). Additionally, we apply Theorem~\ref{theorem:inverse-by-lu2} to compute the inverse of $\bU$ (which costs $2r^3$ flops).
\end{proof}

It should be noted that selecting the first $r$ linearly independent columns from $\bA$ can lead to instability in the inverse of $\bU$. 
To address this issue, a \textit{maxvol procedure}, as described in \citet{goreinov2010find}, can be  employed. 
Towards the end of this section, we will introduce a \textit{pseudoskeleton decomposition} to mitigate this problem.

\index{Gram–Schmidt}
\section{Computing  Skeleton Decomposition via Modified Gram-Schmidt}
In the Gram-Schmidt process, the goal is to determine orthonormal vectors that span the same column space as the original matrix.  
However, this orthogonal property is not necessary  for finding the skeleton decomposition. Our primary objective is simply to determine whether the vectors are linearly independent.
To achieve this, we can adapt the Gram-Schmidt process slightly to select $r$ linearly independent columns from $\bA$ to form $\bC$.  
Similarly, we can employ the adapted Gram-Schmidt process to select $r$ linearly independent columns from $\bC^\top$ to obtain $\bU$ and $\bR$.

For a set of $n$ vectors $\ba_1, \ba_2, \ldots, \ba_n$ from the columns of $\bA$, when we intend to project vector $\ba_k$ ($k\in\{2,3,\ldots,n\}$) onto the space that is orthogonal to the subspace formed by $\ba_1, \ba_2, \ba_{k-1}$ (assuming the $k-1$ vectors are linearly independent), we can express this projection using Equation~\eqref{equation:gram-schdt-eq2}:
$$
\begin{aligned}
\bb_k &= \ba_k -  \left(\frac{ \bb_1^\top\ba_k}{\bb_1^\top\bb_1} \bb_1 + \frac{ \bb_2^\top\ba_k}{\bb_2^\top\bb_2} \bb_2+\ldots + \frac{ \bb_{k-1}^\top\ba_k}{\bb_{k-1}^\top\bb_{k-1}} \bb_{k-1} \right)
=\ba_k -  \sum_{i=1}^{k-1}\frac{ \bb_i^\top\ba_k}{\bb_i^\top\bb_i} \bb_i,
\end{aligned}
$$
where $\bb_k$ represents the projection of $\ba_k$ onto the subspace that is orthogonal to the one spanned by  $\{\ba_1, \ba_2, \ba_{k-1}\}$, $\bb_{k-1}$ is the projection of $\ba_{k-1}$ onto the subspace that is orthogonal to the one spanned by $\{\ba_1, \ba_2, \ba_{k-2}\}$, and so on.

If $\bb_k$ is nonzero, it indicates that $\ba_k$ is linearly independent from $\{\ba_1, \ba_2, \ldots, \ba_{k-1}\}$ for any $k\in\{2,3,\ldots,n\}$, and therefore, it will be included in the columns of $\bC$.
Similarly, we can select  linearly independent rows from $\bC$ to construct $\bU$. The procedure is illustrated  in Algorithm~\ref{alg:skeleton-decomposition-modified}.

\begin{algorithm}[h] 
\caption{Skeleton Decomposition via \textbf{Adapted} Gram-Schmidt Process} 
\label{alg:skeleton-decomposition-modified} 
\begin{algorithmic}[1] 
\Require 
Rank-$r$ matrix $\bA=[\ba_1, \ba_2, \ldots, \ba_n]$ with size $m\times n $; 
\State Initially set column count $ck=0$ and row count $rk=0$; 
\Statex \textbf{Stage A: Obtain independent columns into $\bC$}
\State Set $\bb_1= \ba_1$; \Comment{Suppose the first column is nonzero, $0$ flops}
\For{$k=2$ to $n$} 
\State $\bb_k=\ba_k -  \sum_{i=1}^{k-1}\frac{ \bb_i^\top\ba_k}{\bb_i^\top\bb_i} \bb_i$, (Skip the zero terms $\bb_i$); 
\If{$\bb_k\neq \bzero$}
\State 	Select the $k$-th column of $\bA$ into the $ck$-th column of $\bC$;
\State  $ck=ck+1$;
\EndIf
\EndFor
\Statex \textbf{Stage B: Obtain independent rows into $\bU$ (and $\bR$)}
\State Suppose $\bC=[\bc_1^\top; \bc_2^\top; \ldots; \bc_m^\top]$, where $\bc_i^\top$ is the $i$-th row of $\bC$;
\State Set $\bb_1   = \bc_1$; \Comment{Suppose the first row is nonzero, 0 flops}
\For{$k=2$ to $m$} 
\State $\bb_k = \bc_k -  \sum_{i=1}^{k-1}\frac{ \bb_i^\top\bc_k}{\bb_i^\top\bb_i} \bb_i$, (Skip the zero terms $\bb_i$); 
\If{$\bb_k\neq \bzero$}
\State 	Select the $k$-th row of $\bC$ into the $rk$-th row of $\bU$;
\State  $rk=rk+1$;
\EndIf
\EndFor
\State Select the rows from $\bA$ corresponding to the rows of $\bU$ into $\bR$.
\end{algorithmic} 
\end{algorithm}

\begin{theorem}[Algorithm Complexity: Skeleton via Adapted Gram-Schmidt]\label{theorem:comp-skeleton-decom2}
Algorithm~\ref{alg:skeleton-decomposition-modified} requires $\sim 2(mn^2+rm^2+ r^3)$ flops to compute a skeleton decomposition of an $m\times n$ matrix.
\end{theorem}

\begin{proof}[of Theorem~\ref{theorem:comp-skeleton-decom2}]
For the last loop $n$ in step 4,  the computation involving  $\frac{ \bb_i^\top\ba_n}{\bb_i^\top\bb_i} \bb_i$ consists of the following:

a. $m$ multiplications and $m-1$ additions for the calculation of $\bb_i^\top\ba_n$;  

b. $m$ multiplications for computing $\bb_i^\top\ba_n \times \bb_i $;

c. $m$ multiplications and $m-1$ additions for determining $\bb_i^\top\bb_i$;

d. 1 division; 

\noindent
There are $n-1$ such terms to compute procedure $(a), (b),$ and $(d)$ above. However, for procedure $(c)$, we've already calculated $\bb_i^\top\bb_i$ for $i \in \{1, 2, \ldots, n-2\}$ in the previous loops. As a result, we require a total of \underline{$3m(n-1)+2m-1$} flops to compute $\frac{ \bb_i^\top\ba_n}{\bb_i^\top\bb_i} \bb_i$ for $i\in \{1, 2, \ldots, n-1\}$, where the last $2m-1$ flops are to compute $\bb_{n-1}^\top\bb_{n-1}$.
The final expression $\ba_k -  \sum_{i=1}^{k-1}\frac{ \bb_i^\top\ba_k}{\bb_i^\top\bb_i} \bb_i$ then takes $(n-2)m$ additions and $m$ subtractions, which cost \underline{$(n-1)m$} flops. 
Thus, the overall cost of step 4 in the final iteration $n$ is $4m(n-1)+2m-1$ flops. Let $f(i)\triangleq 4m(i-1)+2m-1$, the total cost for step 3 to step 9 can be obtained by 
$$
\mathrm{cost}=f(2)+f(3)+\ldots+ f(n).
$$
Simple calculations can show this loop requires $4m \frac{n^2-n}{2} + (2m-1)(n-1)$ flops, or \underline{$2mn^2$} flops if we keep only the leading term.

Similarly, for step 13, the total cost for the second loop is $4r \frac{m^2-m}{2} + (2r-1)(m-1)$ flops, or \underline{$2rm^2$} flops if we keep only the leading term.

As a result, the total cost amounts to $2(mn^2+rm^2+r^3)$ flops if we keep only the leading terms, where $2r^3$ stems from the computation of the inverse of $\bU$ as per Theorem~\ref{theorem:inverse-by-lu2}. 
\end{proof}

The complexity of Algorithm~\ref{alg:skeleton-decomposition-modified}  is identical to that of Algorithm~\ref{alg:skeleton-decomposition} since they are mathematically equivalent. However, Algorithm~\ref{alg:skeleton-decomposition-modified} provides greater clarity and eliminates minor computations that do not contribute significantly to the primary complexity.

\section{Computing  Skeleton Decomposition via  Gaussian Elimination}
In Algorithm~\ref{alg:cr-decomposition-alg}, we  discussed the procedure  to obtain the reduced row  echelon form through Gaussian elimination. And the columns containing the pivots are the ones that are linearly independent.  Consequently, we can employ this reduced row echelon form to derive the skeleton decomposition.

\begin{algorithm}[H] 
\caption{Skeleton Decomposition via Gaussian Elimination} 
\label{alg:skeleton-by-row-reduced-echelon} 
\begin{algorithmic}[1] 
\Require 
Rank-$r$ matrix $\bA=[\ba_1, \ba_2, \ldots, \ba_n]$ with size $m\times n $; 
\State Use Algorithm~\ref{alg:cr-decomposition-alg} to obtain the reduced row  echelon form of $\bA$. Subsequently, extract  the columns containing pivots from $\bA$  and store them in matrix $\bC$;
\State Use Algorithm~\ref{alg:cr-decomposition-alg} to obtain the reduced row  echelon form of $\bC^\top$. Subsequently, extract  the columns containing pivots from $\bC^\top$and store them in matrix $\bU^\top$; 
\State Select the rows from $\bA$ corresponding to the rows of $\bU$ into $\bR$.
\end{algorithmic} 
\end{algorithm}

\begin{theorem}[Algorithm Complexity: Skeleton via Guassian Elimination]\label{theorem:comp-skeleton-decom-row-reduced}
Algorithm~\ref{alg:skeleton-by-row-reduced-echelon} requires 
$$
\text{cost}=\left\{
\begin{aligned}
	&(2m^2n-m^3) + (2r^2m-r^3) + 2r^3=2m^2n-m^3 + 2r^2m+r^3 \text{}, \qquad  &\text{if }m\leq n; \\
	&(mn^2)+ (2r^2m-r^3)+2r^3 = mn^2+ 2r^2m+r^3\text{},\qquad &\text{if }m>n.
\end{aligned}
\right.
$$
flops to compute a skeleton decomposition of an $m\times n$ matrix. The cost of $2r^3$ above is again from the calculation of the inverse of $\bU$  (Theorem~\ref{theorem:inverse-by-lu2}).
\end{theorem}

The proof of the above theorem is straightforward  since $r\leq m$ by Theorem~\ref{theorem:cr-decomposition-alg}.

\section{Randomized Algorithms}
The numerical methods for computing skeleton decomposition described here involve selecting the initial linearly independent columns from $\bA$ and placing them into $\bC$, same as that in the CR decomposition. However, this step is not strictly necessary. Randomized algorithms, as discussed in  \citet{mahoney2009cur, boutsidis2009improved, drineas2012fast, kishore2017literature}, do not require the selection of linearly independent columns in a left-to-right manner.

In \citet{goreinov1997pseudo, goreinov2001maximal}, an approximation method for skeleton decomposition known as the \textit{pseudoskeleton approximation} was introduced. 
In this method, the $k$ columns in $\bC$ and $k$ rows in $\bR$ with $k<r$ are chosen such that their intersection $\bU_{k\times k}$ has the maximum volume.
Here, the term ``volume" refers to the determinant of $\bU_{k\times k}$, which should be maximized among all possible $k\times k$ submatrices of $\bA$.

\section{Pseudoskeleton Decomposition via the SVD}\label{section:pseudoskeleton}
We will extensively cover singular value decomposition (SVD) in Chapter~\ref{chapter:SVD}.
For now, let's assume a prior understanding of SVD, and we will demonstrate how to utilize it to approximate the skeleton decomposition. 
Feel free to skip this section on a first reading.

For a  matrix $\bA\in\real^{m\times n}$, we aim to construct an approximation of $\bA$ with a rank $\gamma \leq \min(m,n)$ in the form of the skeleton decomposition. That is, $\bA$ is approximated as $\bA\approx \bC\bU^{-1}\bR$, where $\bC$ and $\bR$ contain $\gamma$ selected columns and rows, respectively, and $\bU$ is the intersection of $\bC$ and $\bR$. 
More precisely, if $I$ and $J$ are the indices of the selected rows and columns, respectively, $\bU$ can be denoted by $\bU=\bA[I,J]$. Note that $\gamma$ is not necessarily equal to the rank $r$ of $\bA$. 

Instead of selecting the $r$ linearly independent columns from $\bA$ (as seen in skeleton decomposition), we opt for $k$ random columns to form matrix $\bC$ with $k>r$ or even $k=\min\{m,n\}$, given by the indices $J$ from  $\bA$. That is, $\bC = \bA[:,J]\in \real^{m\times k}$ in mathematical terms. 
Simultaneously, we select $k$ rows from $\bA$ by the indices $I$ from   $\bA$ into  $\bR=\bA[I,:]$ such that the volume of the intersection matrix $\bU=\bA[I,J]$ is maximized, i.e., $\det(\bU)$ is maximized given the matrix $\bC$ chosen randomly. While $\bC$ is selected randomly, the choice of $\bR$ is deterministic. This leads to the following decomposition:
$$
\bA = \bC_{m\times k}\bU_{k\times k}^{-1}\bR_{k\times n}.
$$
Once more, it's important to note that the inverse of $\bU_{k\times k}$ lacks stability due to the random selection process. 
To address this issue, we perform a full SVD on $\bU_{k\times k}$ (refer to Chapter~\ref{chapter:SVD} for the difference between the reduced SVD and full SVD):
$$
\bU_{k\times k} = \bU_k\bSigma_k\bV_k^\top,
$$
where $\bU_k, \bV_k\in \real^{k\times k}$ are orthogonal matrices, $\bSigma_k$ is a diagonal matrix containing $k$ singular values $\sigma_1 \geq \sigma_2 \geq \ldots \geq \sigma_k$, some of which may be zero.
Subsequently, we choose $\gamma$ singular values that are greater than some value $\epsilon$ and truncate  $\bU_k, \bV_k,$ and $ \bSigma_k$ according to the $\gamma$ selected singular values such that $\bU_{k\times k}$ is approximated by a rank-$\gamma$ matrix $\bU_{k\times k} \approx \bU_\gamma\bSigma_\gamma\bV_\gamma^\top$, where $\bU_\gamma,\bV_\gamma\in \real^{k\times \gamma}$, and $\bSigma_\gamma \in \real^{\gamma\times \gamma}$. Therefore, the pseudo-inverse of $\bU_{k\times k}$ is \footnote{See Section~\ref{section:application-ls-svd} and Appendix~\ref{appendix:pseudo-inverse}  for a detailed discussion of the pseudo-inverse.}
$$
\bU^+ = (\bU_\gamma\bSigma_\gamma\bV_\gamma^\top)^{+} =\bV_\gamma \bSigma_\gamma^{-1}\bU_\gamma^\top.
$$
As a result, the matrix $\bA$ is approximated as a rank-$\gamma$ matrix
\begin{equation}\label{equation:skeleton-low-rank}
\begin{aligned}
\bA 
&\approx \bC\bV_\gamma \bSigma_\gamma^{-1}\bU_\gamma^\top \bR\\
&=\bC_2 \bR_2,  \qquad (\text{let $\bC_2\triangleq\bC\bV_\gamma \bSigma_\gamma^{-1/2}$ and $\bR_2\triangleq\bSigma_\gamma^{-1/2}\bU_\gamma^\top \bR$})
\end{aligned}
\end{equation}
where $\bC_2$ and $\bR_2$ are rank-$\gamma$ matrices. 
For a systematic approach to selecting $\epsilon$, see \citet{goreinov1997pseudo} and \citet{kishore2017literature}. 
In the method described above, we only choose $\bC$ randomly, while $\bR$ is determined. Algorithms presented in \citet{zhu2011randomised} choose both $\bC$ and $\bR$ randomly, yielding more stable results.

In the introduction section of this Part~\ref{part:data-interation}, we claimed that  data interpretation decomposition methods, such as  CR and skeleton decomposition, facilitate identifying a subset of columns or rows that capture  the matrix's essential information by retaining the original columns or rows of the underlying matrix.
A cursory examination of the pseudoskeleton decomposition in Equation~\eqref{equation:skeleton-low-rank} might not seem to exhibit this capability.
However, given that $\bC_2=\bC(\bV_\gamma \bSigma_\gamma^{-1/2})$ and  columns of $\bV_\gamma \bSigma_\gamma^{-1/2}$ are orthogonal, one can infer that the columns of $\bC_2$ represent the columns of $\bC$ in a different coordinate system. 
The same concept applies to the rows of $\bR_2$ (refer to Section~\ref{section:coordinate-transformation} for further insights into matrix decomposition transformations).

\begin{problemset}
\item Find the CUR decomposition for the matrix 
$$
\bA = 
\begin{bmatrix}
	1 & 3 & 2 \\
	3 & 7 & 6 \\
	4 & 5 & 8
\end{bmatrix}.
$$

\item Find the pseudo-inverse of a matrix $\bA$ using its CUR decomposition.

\item \label{problem:nons_sub} Consider a matrix $\bA\in\real^{n\times n}$ in block form $\bA=\scriptsize\begin{bmatrix}
	\bK & \bL\\
	\bM & \bN
\end{bmatrix}$, where $\bK\in\real^{r\times r}$ and $\bN\in\real^{(n-r)\times (n-r)}$. Show that 
\begin{itemize}
	\item If $\bK$ is nonsingular, then $[\bK, \bL]$ and $\scriptsize\begin{bmatrix}
		\bK\\
		\bM 
	\end{bmatrix}$ has full ranks.
	\item If $\rank(\bA) = \rank([\bK, \bL])=\rank(\scriptsize\begin{bmatrix}
		\bK\\
		\bM 
	\end{bmatrix}\normalsize)$, then $\bK$ is nonsingular.
\end{itemize}

\item \label{problem:rank_pcin} Let $\bA\in\real^{n\times n}$ be  a symmetric or skew-symmetric matrix (Definition~\ref{definition:speci_mat}). Show that 
\begin{itemize}
\item $\rank(\bA[I, :])=\rank(\bA[:,I])$ for any index set $I\subseteq \{1,2,\ldots,n\}$.
\item $\bA$ is rank principal (if it has a nonsingular $r\times r$ principal submatrix; Definition~\ref{definition:principle-minors}).
\end{itemize}

\item Investigate how different column and row selection strategies affect the accuracy of the CUR decomposition. Generate a matrix $\bA\in\real^{100\times 100}$, and implement different selection strategies for columns and rows (e.g., random, based on column norms, leverage scores).
For each strategy, perform CUR decomposition/approximation with $r=10$. Evaluate the approximation error for each strategy. Discuss which strategy provides the best approximation.

\item Compare the CUR decomposition with other matrix decomposition methods like SVD and QR decomposition in terms of accuracy and computational efficiency.
Generate a matrix $\bA\in\real^{100\times 100}$. Perform CUR decomposition, SVD, and QR decomposition on $\bA$. For each method, compute the approximation error using a suitable norm. Measure the computational time for each method. Discuss the trade-offs between accuracy and computational efficiency for each method.
\end{problemset}

%% file: chapter-id.tex
\newpage
\chapter{Interpolative Decomposition (ID)}\label{chapter:id}
\begingroup
\hypersetup{
	linkcolor=structurecolor,
	linktoc=page,  
}
\minitoc \newpage
\endgroup

\section{Interpolative Decomposition}
\lettrine{\color{caligraphcolor}T}
The column interpolative decomposition (ID) factors a matrix as the product of two matrices, one of which contains selected columns from the original matrix, and the other of which includes a subset of columns consisting of an identity matrix, with all its values being no greater than 1 in absolute value.  Formally, the specific details of the column ID are outlined in the following theorem.
Without special mention, we will refer to the columns ID simply as ID throughout this discussion, 

\index{Decomposition: ID}
\begin{theoremHigh}[Column Interpolative Decomposition]\label{theorem:interpolative-decomposition}
Let $\bA \in \real^{m \times n}$ be any rank-$r$ matrix. Then it can be factored as 
$$
\underset{m \times n}{\bA} = \underset{m\times r}{\bC} \gapthree  \underset{r\times n}{\bW},~\footnote{See Theorem~\ref{theorem:interpolative-decomposition-row} for details of row and two-sided IDs.}
$$
where $\bC\in \real^{m\times r}$ contains $r$ linearly independent columns of $\bA$, and $\bW\in \real^{r\times n}$ represents the matrix used to reconstruct $\bA$. The matrix $\bW$ contains an $r\times r$ identity submatrix (under a mild column permutation). 
Specifically, the entries in $\bW$ do not exceed 1 in magnitude:
$$
\max |w_{ij}|\leq 1, \,\, \forall \gapforall i\in [1,r], j\in [1,n].
$$

The storage for the decomposition is then reduced or potentially increased from $mn$ floating-point numbers to $mr$ and $(n-r)r$ floating-point numbers for storing $\bC$ and $\bW$, respectively.
Additionally, an extra $r$ integers are required to record the position of each column of $\bC$ in that of $\bA$.
\end{theoremHigh}

\begin{figure}[H]
\centering
\includegraphics[width=0.7\textwidth]{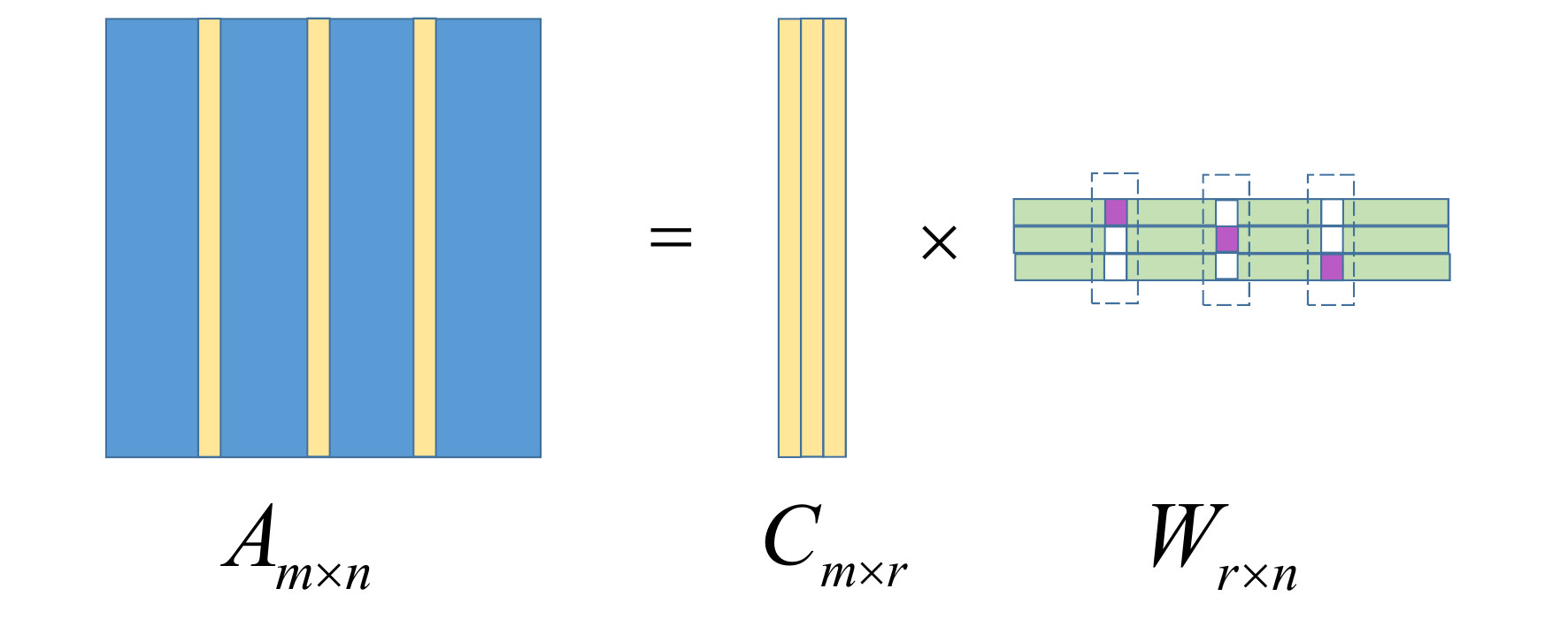}
\caption{Demonstration of the column ID of a matrix, where the \textcolor{mydarkyellow}{yellow} vectors represent the linearly independent columns of $\bA$, white entries denote zeros, and \textcolor{mydarkpurple}{purple} entries signify  ones.}
\label{fig:column-id}
\end{figure}

The illustration of the column ID, which is referred to as ID when there is no ambiguity (see Section~\ref{section:row_twoside_ID} for details on row and two-sided ID), is shown in Figure~\ref{fig:column-id}, where the \textcolor{mydarkyellow}{yellow} vectors denote the linearly independent columns of $\bA$, and the \textcolor{mydarkpurple}{purple} vectors within $\bW$ form an $r\times r$ identity submatrix. 
Notably, the positions of the \textcolor{mydarkpurple}{purple} vectors inside $\bW$ match  those of  the \textcolor{mydarkyellow}{yellow} vectors within $\bA$. The benefits and applications of ID include data compression, storage efficiency, data interpretation, and many more. 
ID reduces the storage needs for large matrices by identifying a subset of columns (or rows) that encapsulate the essential information, providing a compact representation that is significantly smaller than the original matrix, especially for sparse matrices. 
ID can also be used to precondition linear systems, making them easier to solve by approximating the matrix with a low-rank structure, which is particularly advantageous for large and sparse systems.
When the original matrix $\bA$ has certain properties, such as sparsity or nonnegativity, the matrices $\bC$ and $\bR$ (see Section~\ref{section:row_twoside_ID}) derived from the ID retain these properties, which can aid in interpretation and further analysis  \citep{liberty2007randomized, martinsson2011randomized, halko2011finding}.

The concept of column ID closely resembles that of the CR decomposition (Theorem~\ref{theorem:cr-decomposition}), both methods select $r$ linearly independent columns and arrange these columns into the first factor, with the second factor containing an $r\times r$ identity submatrix. 
However, there are distinctions. 
In the CR decomposition, the method precisely chooses the first $r$ linearly independent columns into the first factor, and the identity submatrix appears in the pivot positions (Definition~\ref{definition:pivot}). And more importantly, the second factor in the CR decomposition comes from the RREF (Lemma~\ref{lemma:r-in-cr-decomposition}). 
As a result, the column ID can be applied in contexts where CR decomposition is used, such as proving that the rank of an idempotent matrix equals its trace (Proposition~\ref{proposition:rank-of-symmetric-idempotent2_tmp}) and demonstrating that the column rank equals the row rank of a matrix (Theorem~\ref{theorem:equal-dimension-rank}).
Moreover, the column ID is also a special case of the rank decomposition (Theorem~\ref{theorem:rank-decomposition}) and is apparently not unique. 
The relationship between different column IDs is explained in  Corollary~\ref{corollary:connection-rank-decom}.

\paragraph{Notations that will be extensively used in the sequel.} Following again the Matlab-style notation, if $J_s$ is an index vector of size $r$ that contains the indices of columns selected from $\bA$ to form $\bC$, then $\bC$ can be denoted as $\bC=\bA[:,J_s]$ (Definition~\ref{definition:matlabnotation}).
The matrix $\bC$ consists of the ``skeleton" columns of $\bA$, hence the subscript $s$ in $J_s$. From the ``skeleton" index vector $J_s$, the $r\times r$ identity submatrix inside $\bW$ can be recovered by
$$
\bW[:,J_s] = \bI_r \in \real^{r\times r}.
$$ 
Suppose we proceed by placing the remaining indices of matrix $\bA$ into another index vector denoted by $J_r$, where
$$
J_s\cap J_r=\varnothing \qquad \text{and}\qquad J_s\cup J_r = \{1,2,\ldots, n\}.
$$
The remaining $n-r$ columns in $\bW$ consist of an $r\times (n-r)$ \textit{expansion matrix}:
$$
\bE = \bW[:,J_r] \in \real^{r\times (n-r)},
$$
where the entries of $\bE$ are called \textit{expansion coefficients}  because it contains elements for reconstructing the columns of $\bA$ from $\bC$. Moreover, let $\bP\in \real^{n\times n}$ be a (column) permutation matrix (Definition~\ref{definition:permutation-matrix}) defined by $\bP=\bI_n[:,(J_s, J_r)]$ so that
$$
\bA\bP = \bA[:,(J_s, J_r)] = \left[\bC, \bA[:,J_r]\right],
$$
and 
\begin{equation}\label{equation:interpolatibve-w-ep}
	\bW\bP = \bW[:,(J_s, J_r)] =\left[\bI_r, \bE \right] \quad\implies\quad  \bW = \left[\bI_r, \bE \right] \bP^\top.
\end{equation}
More insights about this notation can be found in the discussion in Problem~\ref{problem:mtb_dis}.

\section{Existence of  Column Interpolative Decomposition}\label{section:proof-column-id}
\index{Cramer's rule}
\paragraph{Cramer's rule.}
The proof of the existence of the column ID depends on Cramer's rule, which we will  discuss briefly here; see Problem~\ref{prob:cramer_adj_1}$\sim$\ref{prob:cramer_adj_4} for more details. 
Cramer's rule provides a formula for each entry of the solution to a linear system.
Consider a system of $n$ linear equations in $n$ unknowns, represented as:
$$
\bM \bx = \bl,
$$
where $\bM\in \real^{n\times n}$ is nonsingular, and $\bx,\bl \in \real^n$. Then the theorem states that, in this case, the system possesses a unique solution, whose individual values for the unknowns are given by:
$$
x_i = \frac{\det(\bM_{\bl}(i))}{\det(\bM)}, \qquad \text{for all}\gapthree i\in \{1,2,\ldots, n\},
$$
where $\bM_{\bl}(i)$ represents the matrix formed by replacing the $i$-th column of $\bM$ with the column vector $\bl$. In full generality, Cramer's rule applies to  the matrix equation
$$
\bM\bX = \bL,
$$
where $\bM\in \real^{n\times n}$ is nonsingular, and $\bX,\bL\in \real^{n\times m}$. Let $I\triangleq[i_1, i_2, \ldots, i_k]$ and $J\triangleq[j_1,j_2,\ldots, j_k]$ be two index vectors, where $1\leq i_1\leq i_2\leq \ldots\leq i_k\leq n$ and $1\leq j_1\leq j_2\leq \ldots\leq j_k\leq m$. 
Then, $\bX[I,J]$ represents a $k\times k$ submatrix of $\bX$. 
Furthermore,  let $\bM_{\bL}(I,J)$ be the $n\times n$ matrix formed by replacing the $(i_s)$-th column of $\bM$ with $(j_s)$-th column of $\bL$ for all $s\in \{1,2,\ldots, k\}$. 
Cramer's rule then states:
$$
\det(\bX[I,J]) = \frac{\det\left(\bM_{\bL}(I,J)\right)}{\det(\bM)}.
$$ 
When $I$ and $J$ are of size 1,  this simplifies to:
\begin{equation}\label{equation:cramer-rule-general}
	x_{ij} = \frac{\det\left(\bM_{\bL}(i,j)\right)}{\det(\bM)}.
\end{equation}
We are now prepared to demonstrate the existence of the column ID.
\begin{proof}[{of Theorem~\ref{theorem:interpolative-decomposition}: Existence of Column ID}]
As previously mentioned, our proof relies on Cramer's rule. 
To show the existence of the column ID, we need to demonstrate that the entries of $\bW$  do not exceed 1 in absolute value.
If we can show the entries of $\bW$ can be denoted by  Cramer's rule equality in Equation~\eqref{equation:cramer-rule-general}, and furthermore, establish that the numerator is smaller than the denominator, then we can complete the proof. However, we notice that the matrix in the denominator of Equation~\eqref{equation:cramer-rule-general} is a square matrix. 
This is where a clever technique comes into play.

\paragraph{Step 1: column ID for full row rank matrix.}
For a start, we first consider the full row rank matrix $\bA$ (which implies $r=m$, $m\leq n$, and $\bA\in \real^{r\times n}$ such that the matrix $\bC\in \real^{r\times r}$ is a square matrix in the column ID of $\bA=\bC\bW$ that we want to obtain). We can determine the ``skeleton" index vector $J_s$ by 
\begin{equation}\label{equation:interpolative-choose-js}
	\boxed{	J_s = \mathop{\arg\max}_{J} \left\{|\det(\bA[:,J])|: \text{$J$ is a subset of $\{1,2,\ldots, n\}$ with size $r=m$} \right\}.}
\end{equation}
In other words, $J_s$ is the index vector that is determined by maximizing the magnitude of the determinant of $\bA[:,J_s]$. As explained in the previous section, there exists a (column) permutation matrix $\bP$ such that
$$
\bA\bP = 
\begin{bmatrix}
	\bA[:,J_s]&\bA[:,J_r]
\end{bmatrix}.
$$
Since $\bC=\bA[:,J_s]$ has a full column rank of $r=m$, it is  nonsingular. The  equation above can be reformulated as follows:
$$
\begin{aligned}
\bA
&=\begin{bmatrix}
\bA[:,J_s]&\bA[:,J_r]
\end{bmatrix}\bP^\top
= 
\bA[:,J_s]
\bigg[
\bI_r \gap \bA[:,J_s]^{-1}\bA[:,J_r]
\bigg]
\bP^\top\\
&= \bC 
\underbrace{\begin{bmatrix}
	\bI_r & \bC^{-1}\bA[:,J_r]
\end{bmatrix}
\bP^\top}_{\triangleq\bW},
\end{aligned}
$$
where the matrix $\bW$ is given by 
$
\begin{bmatrix}
	\bI_r & \bC^{-1}\bA[:,J_r]
\end{bmatrix}\bP^\top
\triangleq
\begin{bmatrix}
	\bI_r & \bE
\end{bmatrix}\bP^\top
$, according to Equation~\eqref{equation:interpolatibve-w-ep}. 
To prove the claim that the magnitudes of the entries in $\bW$ do not exceed 1, it suffices to show that the entries in $\bE\triangleq\bC^{-1}\bA[:,J_r]\in \real^{r\times (n-r)}$ are no greater than 1 in absolute value.

Define the index vector $[j_1,j_2,\ldots, j_n]$ as a permutation of $[1,2,\ldots, n]$ in a way that:
$$
[j_1,j_2,\ldots, j_n] = [1,2,\ldots, n] \bP = [J_s, J_r].~\footnote{Note here $[j_1,j_2,\ldots, j_n] $, $[1,2,\ldots, n]$, $J_s$, and $J_r$ are row vectors.}
$$
Thus, it follows from $\bC\bE=\bA[:,J_r]$ that 
$$
\begin{aligned}
	\underbrace{	[\ba_{j_1}, \ba_{j_2}, \ldots, \ba_{j_r}]}_{=\bC=\bA[:,J_s]} \bE &=
	\underbrace{[\ba_{j_{r+1}}, \ba_{j_{r+2}}, \ldots, \ba_{j_n}]}_{=\bA[:,J_r]\triangleq \bB},
\end{aligned}
$$ 
where $\ba_i$ represents the $i$-th column of $\bA$, and we let $\bB\triangleq\bA[:,J_r]$.
Therefore, according to Cramer's rule in Equation~\eqref{equation:cramer-rule-general}, we have 
\begin{equation}\label{equation:column-id-expansionmatrix}
	e_{kl} = 
	\frac{\det\left(\bC_{\bB}(k,l)\right)}
	{\det\left(\bC\right)}, \qquad \text{for all}\gapthree k\in \{1,2,\ldots, r\}, \,\,\,l\in \{1,2,\ldots, n-r\},
\end{equation}
where $e_{kl}$ represents the entry ($k,l$) of $\bE$, and $\bC_{\bB}(k,l)$ denotes the $r\times r$ matrix formed by replacing the $k$-th column of $\bC$ with the $l$-th column of $\bB$. For example, 
$$
\begin{aligned}
	e_{11} &= 
	\frac{\det\left([\textcolor{mylightbluetext}{\ba_{j_{r+1}}}, \ba_{j_2}, \ldots, \ba_{j_r}]\right)}
	{\det\left([\ba_{j_1}, \ba_{j_2}, \ldots, \ba_{j_r}]\right)},
	\qquad 
	&e_{12} &=
	\frac{\det\left([\textcolor{mylightbluetext}{\ba_{j_{r+2}}}, \ba_{j_2},\ldots, \ba_{j_r}]\right)}
	{\det\left([\ba_{j_1}, \ba_{j_2}, \ldots, \ba_{j_r}]\right)},\\
	e_{21} &= 
	\frac{\det\left([\ba_{j_1},\textcolor{mylightbluetext}{\ba_{j_{r+1}}}, \ldots, \ba_{j_r}]\right)}
	{\det\left([\ba_{j_1}, \ba_{j_2}, \ldots, \ba_{j_r}]\right)},
	\qquad 
	&e_{22} &= 
	\frac{\det\left([\ba_{j_1},\textcolor{mylightbluetext}{\ba_{j_{r+2}}}, \ldots, \ba_{j_r}]\right)}
	{\det\left([\ba_{j_1}, \ba_{j_2}, \ldots, \ba_{j_r}]\right)}.
\end{aligned}
$$
Since $J_s$ is chosen to maximize the magnitude of $\det(\bC)$ in Equation~\eqref{equation:interpolative-choose-js}, it follows that 
$$
|e_{kl}|\leq 1, \qquad \text{for all}\gapthree k\in \{1,2,\ldots, r\}, \,\,\,l\in \{1,2,\ldots, n-r\}.
$$

\paragraph{Step 2: apply to general matrices.}
To summarize what we have proved above and to abuse the notation, we have demonstrated that for any matrix $\bF\in \real^{r\times n}$ with \textbf{full} rank $r\leq n$, the column ID exists such that $\bF=\bC_0\bW$, where the entries in $\bW$ are not greater than 1 in absolute value.

Applying the finding to the  general matrix $\bA\in \real^{m\times n}$ with rank $r\leq \min\{m,n\}$, it is evident that the matrix $\bA$ admits a rank decomposition (Theorem~\ref{theorem:rank-decomposition}):
$$
\underset{m\times n}{\bA} = \underset{m\times r}{\bD}\gapthree \underset{r\times n}{\bF},
$$
where $\bD$ and $\bF$ have full column rank $r$ and full row rank $r$, respectively. For the column ID of $\bF=\bC_0\bW$, where $\bC_0=\bF[:,J_s]$ contains $r$ linearly independent columns of $\bF$. We notice from $\bA=\bD\bF$ such that 
$$
\bA[:,J_s]=\bD\bF[:,J_s],
$$
i.e., the columns indexed by $J_s$ of $(\bD\bF)$ can be obtained by $\bD\bF[:,J_s]$, which in turn are the columns of $\bA$ indexed by $J_s$. This makes
$$
\underbrace{\bA[:,J_s]}_{\triangleq\bC}= \underbrace{\bD\bF[:,J_s]}_{\bD\bC_0},
$$
and 
$$
\bA = \bD\bF =\bD\bC_0\bW = \underbrace{\bD\bF[:,J_s]}_{\bC}\bW=\bC\bW.
$$
This completes the proof.
\end{proof}

The above proof reveals an intuitive way to compute the ``optimal" column ID of a matrix $\bA$, as illustrated in Algorithm~\ref{alg:column-id-intuitive}. However, any algorithm that is guaranteed to find such an optimally-conditioned factorization must have combinatorial complexity \citep{martinsson2019randomized}. In the subsequent sections, we will consider alternative approaches to find a relatively well-conditioned factorization. 
\begin{algorithm}[h] 
	\caption{An \textcolor{mylightbluetext}{Intuitive} Method to Compute the Column ID} 
	\label{alg:column-id-intuitive} 
	\begin{algorithmic}[1] 
		\Require 
		Rank-$r$ matrix $\bA$ with size $m\times n $; 
		\State Compute the rank decomposition $\underset{m\times n}{\bA} = \underset{m\times r}{\bD}\gapthree \underset{r\times n}{\bF}$ (Theorem~\ref{theorem:rank-decomposition}) such as from UTV (Section~\ref{section:ulv-urv-decomposition});
		\State Compute column ID of $\bF$: $\bF=\bF[:,J_s]\bW = \widetildebC\bW$:
		$$
		\begin{aligned}
			2.1.\,\,\, &\left\{
			\begin{aligned}
				J_s &= \mathop{\arg\max}_{J} \left\{|\det(\bF[:,J])|: \text{$J$ is a subset of $\{1,2,\ldots, n\}$ with size $r$} \right\};&\\
				J_r &= \{1,2,\ldots, n\} \backslash J_s;&\\
			\end{aligned}
			\right.\\
			2.2.\,\,\,&\left\{
			\begin{aligned}
				\widetildebC &= \bF[:,J_s]; \\
				\bM &= \bF[:,J_r];
			\end{aligned}
			\right.\\
			2.3. \,\,\,&\bF\bP = \bF[:,(J_s,J_r)] \text{ to obtain permutation matrix $\bP$};\\
			2.4. \,\,\,&e_{kl} = 
			\frac{\det\left(\widetildebC_{\bM}(k,l)\right)}
			{\det\left(\widetildebC\right)},  \qquad \text{for all}\gapforall k\in [1, r], l\in [1,n-r] \text{  ~(Equation~\eqref{equation:column-id-expansionmatrix})};\\
			2.5. \,\,\,&\bW= [\bI_r, \bE]\bP^\top \text{ ~(Equation~\eqref{equation:interpolatibve-w-ep})}.
		\end{aligned}
		$$
		\State 	$\bC=\bA[:,J_s]$;
		\State Output the column ID $\bA=\bC\bW$;
	\end{algorithmic} 
\end{algorithm}

\begin{example}[Compute the Column ID]\label{example:column-id-a}
Let
$
\bA=
\scriptsize
\begin{bmatrix}
56 & 41 & 30\\
32 & 23 & 18\\
80 & 59 & 42
\end{bmatrix}
$
be given with rank 2. The trivial process for computing the column ID of $\bA$ is shown as follows. 
First, we find a rank decomposition:
$$
\bA = \bD\bF=
\begin{bmatrix}
1 & 0 \\
0 & 1 \\
2 &-1
\end{bmatrix}
\begin{bmatrix}
56 & 41 & 30 \\
32 & 23 & 18 
\end{bmatrix}.
$$
Since the rank $r=2$, the index vector $J_s$ can take one of the following values: $[1,2], [0,2], [0,1]$, where the absolute determinants of $\bF[:,J_s]$ are $48, 48$, and $ 24$, respectively. We proceed by choosing $J_s=[0,2]$:
$$
\begin{aligned}
\widetildebC &\triangleq \bF[:,J_s]=
\begin{bmatrix}
56 & 30 \\
32 & 18 
\end{bmatrix},\qquad 
\bM \triangleq \bF[:,J_r]=\begin{bmatrix}
41 \\
23
\end{bmatrix}.
\end{aligned}
$$
Thus, 
$$
\bF\bP = \bF[:(J_s,J_r)] = \bF[:,(0,2,1)]
\quad\implies\quad 
\bP = 
\begin{bmatrix}
1 &  & \\
& &1\\
& 1 & 
\end{bmatrix}.
$$
In this example, $\bE\in \real^{2\times 1}$:
$$
\begin{aligned}
e_{11} &=
\det\left(
\begin{bmatrix}
41 & 30 \\
23 & 18
\end{bmatrix}\right)\bigg/
\det\left(
\begin{bmatrix}
56 & 30 \\
32 & 18
\end{bmatrix}\right)=1;\\
e_{21} &=
\det\left(
\begin{bmatrix}
56 & 41 \\
32 & 23
\end{bmatrix}\right)\bigg/
\det\left(
\begin{bmatrix}
56 & 30 \\
32 & 18
\end{bmatrix}\right)=-\frac{1}{2}.
\end{aligned}
$$
This makes 
$$
\bE = 
\begin{bmatrix}
1\\-\frac{1}{2}
\end{bmatrix}
\quad \implies\quad  
\bW = [\bI_2, \bE]\bP^\top =
\begin{bmatrix}
1 & 1 & 0\\
0 & -\frac{1}{2} & 1
\end{bmatrix}.
$$
The final selected columns  are
$$
\bC = \bA[:,J_s] = 
\begin{bmatrix}
56 & 30\\
32 & 18\\
80 & 42
\end{bmatrix}.
$$
And the net result is given by
$$
\bA=\bC\bW =
\begin{bmatrix}
56 & 30\\
32 & 18\\
80 & 42
\end{bmatrix}
\begin{bmatrix}
1 & 1 & 0\\
0 & -\frac{1}{2} & 1
\end{bmatrix},
$$
where entries of $\bW$ are no greater than 1 in absolute value, as desired.
\end{example}

To conclude this section, we discuss  the non-uniqueness of the column ID.
\begin{remark}[Non-uniqueness of the Column ID]
In the  specific Example~\ref{example:column-id-a} provided above, we notice that both $\bF[:,(1,2)]$ and $\bF[:,(0,2)]$ yield   the  maximum absolute determinant. 
Therefore, both of them can result in a column ID of $\bA$. Whilst, we only select one $J_s$ from $[1,2], [0,2]$, and $ [0,1]$. 
Additionally, once we fix the index set $J_s$ based on the maximum absolute determinant search, any permutation of it can also be selected, e.g., $J_s=[0,2]$ or $J_s=[2,0]$ are both good. The two choices on the selection of the column index search yield the non-uniqueness of the column ID.
\end{remark}

\section{Row ID and Two-Sided ID}\label{section:row_twoside_ID}

We refer to the decomposition described above as the column ID. This terminology is not arbitrary since it has its siblings:
\begin{theoremHigh}[The Complete  Interpolative Decomposition]\label{theorem:interpolative-decomposition-row}
Let $\bA \in \real^{m \times n}$ be any rank-$r$ matrix, Then it can be factored as 
$$
\begin{aligned}
\text{Column ID: }&\gap \underset{m \times n}{\bA} &=& \boxed{\underset{m\times r}{\bC}} \gap  \underset{r\times n}{\bW} ; \\
\text{Row ID: } &\gap &=&\underset{m\times r}{\bZ} \gap  \boxed{\underset{r\times n}{\bR}}; \\
\text{Two-Sided ID: } &\gap &=&\underset{m\times r}{\bZ} \gap \boxed{\underset{r\times r}{\bU}}  \gap  \underset{r\times n}{\bW}, \\
\end{aligned}
$$
where
\begin{itemize}
\item $\bC=\bA[:,J_s]\in \real^{m\times r}$ represents a selection of $r$  linearly independent columns from $\bA$, and $\bW\in \real^{r\times n}$ is the matrix used to reconstruct $\bA$ and it contains an $r\times r$ identity submatrix (under a mild column permutation): $\bW[:,J_s]=\bI_r$;
\item $\bR=\bA[I_s,:]\in \real^{r\times n}$ represents a selection of $r$ linearly independent rows from $\bA$, and $\bZ\in \real^{m\times r}$ is the matrix used to reconstruct $\bA$ and it contains an $r\times r$ identity submatrix (under a mild row permutation): $\bZ[I_s,:]=\bI_r$;
\item The entries in $\bW$ and $\bZ$ have values no larger than 1 in magnitude: $\max |w_{ij}|\leq 1$ and $\max |z_{ij}|\leq 1$ for all $i,j$;
\item $\bU=\bA[I_s,J_s] \in \real^{r\times r}$ is the nonsingular submatrix at the intersection of $\bC$ and $\bR$;
\item \textbf{Skeleton decomposition:} the three matrices $\bC,\bR$, and $\bU$ in the $\boxed{\text{boxed}}$ texts share same notations as the skeleton decomposition (Theorem~\ref{theorem:skeleton-decomposition}), where they even have the same meanings such that the three matrices make the skeleton decomposition of $\bA$: $\bA=\bC\bU^{-1}\bR$.
\end{itemize}
\end{theoremHigh}
The proof of the row ID is analogous to that of the column ID. Suppose the column ID of $\bA^\top$ is given by $\bA^\top=\bC_0\bW_0$, where $\bC_0$ contains $r$ linearly independent columns of $\bA^\top$ (i.e., $r$ linearly independent rows of $\bA$). Let $\bR\triangleq\bC_0^\top$ and $\bZ\triangleq\bW_0^\top$. Then, the row ID is obtained by $\bA=\bZ\bR$. 

For the two-sided ID, recalling from the skeleton decomposition (Theorem~\ref{theorem:skeleton-decomposition}), if $\bU$ is the intersection of $\bC$ and $\bR$, it follows that $\bA=\bC\bU^{-1}\bR$. 
Consequently, we can deduce $\bC\bU^{-1}=\bZ$ through the row ID, which, in turn, implies $\bC=\bZ\bU$. Leveraging  column ID, we can further infer  $\bA=\bC\bW=\bZ\bU\bW$, thereby confirming the existence of the two-sided ID.

\paragraph{Data storage.} Regarding the data storage requirements for each ID, we can summarize as follows:
\begin{itemize}
\item \textit{Column ID.} It requires $mr$ and $(n-r)r$ floating-point numbers to store $\bC$ and $\bW$, respectively, and an extra $r$ integers to store the indices of the selected columns in $\bA$;
\item \textit{Row ID.} It requires $nr$ and $(m-r)r$ floating-point numbers to store $\bR$ and $\bZ$, respectively, and an extra $r$ integers to store the indices of the selected rows in $\bA$;
\item \textit{Two-Sided ID.} It requires $(m-r)r$, $(n-r)r$, and $r^2$ floating-point numbers to store $\bZ,\bW$, and $\bU$, respectively. And an extra $2r$ integers are required to store the indices of the selected rows and columns in $\bA$.
\end{itemize}

\paragraph{Storage reduction for sparse matrices.} 
Suppose the column ID of $\bA$ is given by $\bA=\bC\bW$, where $\bC=\bA[:,J_s]$, and a good spanning row index  set $I_s$ of $\bC$ could be found such that:
$$
\bA[I_s,:] = \bC[I_s,:]\bW = \bA[I_s,J_s]\bW.
$$ 
We observe that $\bC[I_s,:] = \bA[I_s,J_s]\in \real^{r\times r}$, which is nonsingular (since it has full rank $r$ in the sense of both row rank and column rank). Consequently, it follows that 
$$
\bW = (\bA[I_s,J_s])^{-1} \bA[I_s,:].
$$
Therefore, there is no necessity to store the matrix $\bW$ explicitly in the two-sided ID. We only need to store $\bA[I_s,:]$ and $(\bA[I_s,J_s])^{-1}$. 
Alternatively, when we are able to  compute the inverse of $\bA[I_s,J_s]$ on the fly, it only requires $r$ integers to store $J_s$, and we can recover $\bA[I_s,J_s]$ from $\bA[I_s,:]$. The storage of $\bA[I_s,:]$ is cheap if $\bA$ is sparse.

\section{Computing  Column ID via the CPQR}
The method used in the previous section can be applied to compute the "optimal" column ID.
However, any algorithm that is guaranteed to find such an optimally-conditioned factorization must have combinatorial complexity. An inexpensive alternative is to ensure that the factor $\bW$ is small in norm rather than bounding each entry in modulus by one.
Recall the column-pivoted QR decomposition (CPQR, Theorem~\ref{theorem:rank-revealing-qr-general}). For a matrix $\bA\in \real^{m\times n}$, the \textit{reduced} CPQR is given by:
$$
\bA\bP = \bQ_r
\underbrace{\begin{bmatrix}
		\bR_{11} & \bR_{12} \\ 
\end{bmatrix}}_{\bR}
=
\begin{bmatrix}
	\bQ_r\bR_{11} & \bQ_r\bR_{12} \\ 
\end{bmatrix},
$$
where $\bR_{11} \in \real^{r\times r}$ is upper triangular, $\bR_{12} \in \real^{r\times (n-r)}$, $\bQ_r\in \real^{m\times r}$ contains orthonormal columns, and $\bP$ is a permutation matrix. 
The computational complexity of the CRPQ is $\mathcalO(mnr)$ flops (Theorem~\ref{theorem:qr-reduced-rank-revealing}). And the operation  $\bA\bP$ permutes the $r$ linearly independent columns into the first $r$ columns of $\bA\bP$:
$$
\bA\bP = \bA[:, (J_s, J_r)]=
\begin{bmatrix}
	\bQ_r\bR_{11} & \bQ_r\bR_{12} \\ 
\end{bmatrix}.
$$
In the ``practical" CPQR via CGS  introduced in Section~\ref{section:practical-cpqr-cgs}, $\bQ_r\bR_{11}$ contains the $r$ linearly independent columns of $\bA$ with the largest norm and $\bQ_r\bR_{12}$ is small in norm. This is important to our aim that the column ID is \textit{well-conditioned} in that entries of $\bW$ are small in magnitude.  
Let $\bC\triangleq\bQ_r\bR_{11}$, and solve the linear equation $\bC\bE=\bQ_r\bR_{12}$, the column ID then follows:
$$
\bA = \bC \underbrace{[\bI_r, \bE]\bP^\top}_{\bW}.
$$
\paragraph{Calculation of the linear system.}
The linear system $\bC\bE=\bQ_r\bR_{12}$ is well defined and does not involve least squares, as was shown in Section~\ref{section:application-ls-qr}, since every column of $\bQ_r\bR_{12}$ lies within the column space of $\bC$ (i.e., the column space of $\bA$). The solution can be obtained via the \textit{normal equation} \footnote{We shall briefly discuss in Definition~\ref{definition:normal-equation-als}.}: $\bE=(\bC^\top\bC)^{-1} \bC^\top\bQ_r\bR_{12}$.
An extra cost for computing the columnd ID thus comes from the calculation of the matrix $\bE$. The calculation of $(\bC^\top\bC)^{-1}$ requires $\underline{r^2(2m-1)+2r^3}$ flops, where $2r^3$ comes from the computation of the inverse of an $r\times r$ matrix (Theorem~\ref{theorem:inverse-by-lu2}). 
Let's break down the remaining steps:
$$
\text{step 2: }\qquad 
\bE=\underbrace{(\bC^\top\bC)^{-1} }_{r\times r}
\gap
\underbrace{ \bC^\top}_{r\times m} 
\gap 
\underbrace{\bQ_r}_{m\times r}
\gap 
\underbrace{ \bR_{12}}_{r\times (n-r)}.
$$
Since $r<m$, the product $\bC^\top\bQ_r$ should be considered as the  next step to make a smaller $r\times r$ matrix, requiring $\underline{(2m-1)r^2}$ flops:
$$
\text{step 3: }\qquad 
\bE=\underbrace{(\bC^\top\bC)^{-1} }_{r\times r}
\gap
\underbrace{ (\bC^\top\bQ_r)}_{r\times r} 
\gap 
\underbrace{ \bR_{12}}_{r\times (n-r)}.
$$
When $r<(n-r)$, the calculation of $(\bC^\top\bC)^{-1} (\bC^\top\bQ_r)$ in the  equation above should be performed first, resulting in a remaining complexity of $\underline{(2r-1)r^2+ (2r-1)r(n-r)}$ flops. Otherwise, $(\bC^\top\bQ_r)\bR_{12}$ should be performed first, leading to a remaining complexity of $\underline{2(2r-1)r(n-r)}$ flops. To summarize, the final complexity for the calculation of the normal equation is:
$$
\text{cost=}
\left\{
\begin{aligned}
	&\{r^2(2m-1)+2r^3 + (2m-1)r^2\}+(2r-1)r^2+ (2r-1)r(n-r), \gap &r&<(n-r); \\
	&\{r^2(2m-1)+2r^3 + (2m-1)r^2\}+ 2(2r-1)r(n-r), \gap &r&\geq (n-r).\\
\end{aligned}
\right.
$$
The normal equation provides a straightforward approach to solving the linear equation mentioned above. Alternatively, iterative methods like gradient descent (Section~\ref{section:als-gradie-descent}) can also be utilized.

\subsubsection*{\textbf{Partial factorization via the CPQR}}
The computational complexity of  algorithms for computing the column ID of a matrix via the CPQR depends on the complexity of the CPQR decomposition. 
For a matrix $\bA\in \real^{m\times n}$ with an ``exact" rank of $r$, the complexity of the CPQR is $\mathcalO(mr^2)$ flops. 
However, when dealing with a rank-deficient matrix $\bA$, the complexity increases to $\mathcalO(mn^2)$. 
The partial factorization CPQR via MGS (Section~\ref{section:partial-cpqr-mgs}) can be employed to address this issue, where the partial CPQR decomposition is given by 
$$
\bA\bP = 
\bQ\bR=
\bQ
\begin{bmatrix}
	\bR_{11} & \bR_{12} \\
	\bzero & \bR_{22}
\end{bmatrix},
$$
where $\bR_{22}$ is small in norm. Such a partial factorization can either take a predetermined rank $k$, or a tolerance $\delta$ such that whenever $r_{kk}<\delta$ in the upper triangular matrix $\bR$ of the CPQR decomposition, we terminate. 
This approach is similar to what we will introduce in the low-rank column ID via the rank-revealing QR (RRQR) decomposition. 
However,  in the partial CPQR, there is no guarantee that the norm of $\bR_{22}$ will be minimized, which is assured in RRQR.

\section{Low-Rank Column ID via  RRQR}

An approximate rank-$\gamma$ ID of a matrix $\bA\in \real^{m\times n}$ is the approximate factorization:
$$
\begin{aligned}
\underset{m\times n}{\bA} \approx \underset{m\times \gamma}{\bC} \gap\underset{ \gamma \times n}{ \bW}
\quad\text{and}\quad
\bA \bP \approx \bC \gap[\bI, \bE] ,\\
\end{aligned}
$$
where  the partial column skeleton $\bC\in \real^{m\times \gamma}$
is given by a subset of the columns of $\bA$, $\gamma$ is known as the \textit{numerical rank}, the entries of $\bE$ are known as the \textit{expansion coefficients} as we have shown previously,  and $\bW$
is well-conditioned in a sense that we will make precise shortly.

The low-rank approximation of column ID has been  extensively studied in the context of rank-revealing QR decomposition
\citep{voronin2017efficient, martinsson2019randomized, martinsson2002randomized, halko2011finding}.
In the context of rank-revealing QR (RRQR) (Equation~\eqref{equation:rankr-reval-qr}\index{Rank-revealing QR}), there exists a permutation $\bP$ such that the linearly independent columns of $\bA$ can be permuted to the left and 
$$
\begin{aligned}
\bA\bP = 
\bQ\bR
&=
\begin{bmatrix}
\bQ_1 & \bQ_2
\end{bmatrix}
\begin{bmatrix}
\bL & \bM \\
\bzero & \bN
\end{bmatrix}
\triangleq
\begin{bmatrix}
\bQ_1\bL & \bQ_1\bM+\bQ_2\bN 
\end{bmatrix}
\triangleq
\bQ_1\bL
\begin{bmatrix}
\bI_{\gamma} & \bY
\end{bmatrix},
\end{aligned}
$$
where $\bN \in \real^{(n-\gamma)\times (n-\gamma)}$, and $\norm{\bN}$ is small in some norm. Here, $\bY$ is  the solution of the linear system $(\bQ_1\bL) \bY = (\bQ_1\bM+\bQ_2\bN)$.
We observe that $\bQ_1 \bL$ represents the first $\gamma$ columns of $\bA\bP$. 
Let $\bC\triangleq\bQ_1\bL$ and $\bE\triangleq\bY$. Then, the low-rank column ID approximation can be obtained by solving the linear system 
\begin{equation}\label{equation:rrqr-columnid-linear}
	(\bQ_1\bL) \bE = (\bQ_1\bM+\bQ_2\bN).
\end{equation}
Since $\bN$ is small in norm, we can approximate it to determine the solution of
$$
\begin{aligned}
	\underset{\gamma\times \gamma}{\bL} & &\underset{\gamma \times (n-\gamma)}{\bE}= &\gap \underset{(n-\gamma)  \times (n-\gamma)  }{\bM} .\\
\end{aligned}
$$
The problem is well defined since $\bL$ is nonsingular (upper triangular with nonzero diagonals). 
If, by chance, $\bL$ were singular, it would imply that $\bA$ has a rank, denoted by $\gamma^\prime$, which is less than $\gamma$.
In such a scenario, the lower $\gamma-\gamma^\prime$
rows in the  linear system above would consist of all zeros, so there exists a solution in this case as well. 

The approximation error of the column ID obtained via RRQR with pivoting is identical to  that of
the standard RRQR method:
$$
\norm{\bA-\bC\bW} = \norm{\bQ_2\bN},
$$
where $\bW=[\bI_{\gamma}, \bE]\bP^\top$.

\begin{remark}[Condition of the Linear System]
Unfortunately, $\bL$ in the linear Equation~\eqref{equation:rrqr-columnid-linear} is typically quite ill-conditioned. 
However, there is still a solution $\bE$ to the linear system in Equation~\eqref{equation:rrqr-columnid-linear} with entries of moderate size.
Informally, the directions
where $\bL$ and $\bM$ point in are ``lined up."
The details can be found in \citet{cheng2005compression}.
\end{remark}

\section{Computing  ID via Randomized Algorithm}\label{section:randomi-id}
Suppose the matrix $\bA$ admits a rank decomposition:
$$
\underset{m\times n}{\bA} = \underset{m\times r}{\bD}\gapthree \underset{r\times n}{\bF}.
$$
Upon finding the rank decomposition, suppose further that the row ID of $\bD$ is obtained by 
$$
\text{row ID of $\bD$: }\qquad \underset{m\times r}{\bD} = \underset{m\times r}{\bZ}\gapthree \underset{r\times r}{\bR_0}
=\bZ\bD[I_s,:].
$$
Similar to the proof of the column ID in Section~\ref{section:proof-column-id}, we observe that $\bZ \bA[I_s,:]$  automatically forms a row ID of $\bA$: $\bA=\bZ\bA[I_s,:]$. This can be demonstrated as follows:
\begin{equation}\label{equation:row-id-sub-d}
\begin{aligned}
	\bZ\bA[I_s,:] &= \bZ (\bD[I_s,:]\bF) &\qquad &\text{(since $\bA=\bD\bF$)}\\
	&=\bD \bF=\bA. &\qquad &\text{(since $\bD=\bZ\bD[I_s,:]$)}
\end{aligned}
\end{equation}
We notice that $\bD$ spans the same column space as $\bA$: $\cspace(\bD)=\cspace(\bA)$. 
The finding above indicates that as long as we find an $m\times r$ matrix $\bD$ spanning the same column space as $\bA$, a row ID can be applied to $\bD$ to find the row ID of $\bA$. This reveals the randomized algorithm for computing the row ID of $\bA$. 
To show this, we require a few relevant facts.
\begin{lemma}[Subspace of $\bA^\top \bA$ and $\bA\bA^\top$]\label{lemma:rank-of-ttt}
Let $\bA\in \real^{m\times n}$ be given. Then,
\begin{itemize}
	\item The column space of $\bA^\top \bA$ is identical  to the column space of $\bA^\top$ (i.e., row space of $\bA$): $\cspace(\bA^\top\bA)=\cspace(\bA^\top)$;
	\item The column space of $\bA\bA^\top$ is identical  to the column space of $\bA$: $\cspace(\bA\bA^\top)=\cspace(\bA)$.
\end{itemize}
\end{lemma}
\begin{proof}[of Lemma~\ref{lemma:rank-of-ttt}]
Let $\bx\in \nspace(\bA)$, we have 
$
\bA\bx  = \bzero \implies \bA^\top\bA \bx =\bzero, 
$
i.e., $\bx\in \nspace(\bA) \implies \bx \in \nspace(\bA^\top \bA)$. Therefore, $\nspace(\bA) \subseteq \nspace(\bA^\top\bA)$. 
Furthermore, let $\bx \in \nspace(\bA^\top\bA)$, we have 
$$
\bA^\top \bA\bx = \bzero\implies \bx^\top \bA^\top \bA\bx = 0\implies \normtwo{\bA\bx}^2 = 0 \implies \bA\bx=\bzero, 
$$
i.e., $\bx\in \nspace(\bA^\top \bA) \implies \bx\in \nspace(\bA)$. Therefore, $\nspace(\bA^\top\bA) \subseteq\nspace(\bA) $. 
As a result, by ``sandwiching," it follows that  
$
\nspace(\bA) = \nspace(\bA^\top\bA).
$
According to the fundamental theorem of linear algebra in Appendix~\ref{appendix:fundamental-rank-nullity}, we have 
$$
\cspace(\bA^\top)=\cspace(\bA^\top\bA).
$$
Applying the same process to $\bA^\top$ leads to the second part of the lemma.
\end{proof}
\begin{exercise}
Let $\bA\in \real^{m\times n}$ be given. Show that $\cspace(\bA\bA^\top\bA) = \cspace(\bA)$.
\end{exercise}

The  lemma above indicates that if we have a matrix $\bA\in \real^{m\times n}$ with rank $r$ and a matrix $\bB\in \real^{m\times k}$ with $k>r$ such that $\cspace(\bB)=\cspace(\bA)$,  a moment of reflexion reveals that $\bA\bA^\top\bB$ spans the same column space as $\bA$: 
\begin{equation}\label{equation:random-row-id-column-space0}
	\cspace(\bA\bA^\top\bB)=\cspace(\bA), \qquad \text{if }\gap  \cspace(\bB)=\cspace(\bA).
\end{equation}
Furthermore, according to Lemma~\ref{lemma:column-basis-from-row-basis}, for any matrix $\bA\in \real^{m\times n}$, suppose that $\{\bg_1, \bg_2, \ldots, \bg_r\}$ is a set of vectors in $\real^n$ that forms a basis for the row space of $\bA$. Then, $\{\bA\bg_1, \bA\bg_2, \ldots, \bA\bg_r\}$ is a basis for the column space of $\bA$:
\begin{equation}\label{equation:random-row-id-column-space-eq2}
	\cspace(\bA\bG) =\cspace(\bA),
	\quad\text{where $\bG = [\bg_1, \bg_2, \ldots, \bg_r]$ and $\cspace(\bG^\top)=\cspace(\bA^\top)$.}
\end{equation}
Thus, the  matrix $\bB$ in Equation~\eqref{equation:random-row-id-column-space0} can be constructed as $\bA\bG$, where the columns of $\bG$ comprise the row basis of $\bA$:
\begin{equation}\label{equation:random-row-id-column-space2}
	\cspace(\bA\bA^\top \bA\bG) =\cspace(\bA). 
\end{equation}
\index{Fundamental theorem}

\begin{algorithm}[h] 
\caption{A \textcolor{mylightbluetext}{Randomized} Method to Compute the \textcolor{mylightbluetext}{Row} ID} 
\label{alg:row-id-randomize1} 
\begin{algorithmic}[1] 
\Require 
Rank-$r$ matrix $\bA$ with size $m\times n $; 
\State Decide the over-sampling parameter $k$ (e.g., $k=10$)$\rightarrow$ let $z=r+k$; 
\State Decide the iteration number $\eta$ (e.g., $\eta=0,1$ or $2$);
\State Generate $r+k$ Gaussian random vectors in $\real^n$ into columns of matrix $\bG\in \real^{n\times (r+k)}$;\gap\Comment{i.e., probably contain the row basis of $\bA$}
\State Initialize $\bD\leftarrow\bA\bG\in \real^{m\times (r+k)}$; \Comment{i.e., probably $\cspace(\bD)=\cspace(\bA)$, $(2n-1)mz$ flops}
\For{$i=1$ to $\eta$}
\State $\bD \leftarrow \bA\bA^\top\bD$; \Comment{$(2m-1)nz+(2n-1)mz$ flops}
\EndFor
\State Calc. the row ID of small matrix $\bD$: $\bD=\bZ\bR_0=\bZ\bD[I_s,:]$;\Comment{$\mathcalO(mz^2)$ flops by CPQR}
\State Output the row ID of $\bA$: $\bA=\bZ\bA[I_s,:]$; 
\end{algorithmic} 
\end{algorithm}

Normally, we could simply stop at $\bA\bG$ since $\cspace(\bA\bG) =\cspace(\bA)$, and the row ID of $\bA\bG$ would reveal the row ID of $\bA$ by Equation~\eqref{equation:row-id-sub-d}. Computing the row ID of $\bA\bG$ is a relatively simpler task compared to that of $\bA$ directly, since $\bA\bG\in \real^{m\times r}$ and $r\ll n$. 

However, in some situations, when the matrix $\bA$ is not exactly rank-$r$, i.e., some singular values $\sigma_k$ (with $k>r$) of $\bA$ may be small enough to be regarded as zero (truncated). We will introduce the SVD in Chapter~\ref{chapter:SVD} such that any matrix $\bA$ admits the factorization $\bA=\bU\bSigma\bV^\top$, where $\bU$ and $\bV$ are orthogonal, and roughly speaking, $\bSigma$ is a diagonal matrix containing singular values, with the number of nonzero singular values being the rank of $\bA$. The columns of $\bU$ contain the column basis of $\bA$ and the columns of $\bV$ contain the row basis of $\bA$ (Lemma~\ref{section:four-space-svd}). Therefore, $\bA\bA^\top \bA\bG$ results in
$$
\bA\bA^\top \bA\bG = (\bU\bSigma\bV^\top) (\bV\bSigma\bU^\top) (\bU\bSigma\bV^\top)\bG=\bU\bSigma^3\bV^\top\bG,
$$
where $\bSigma^3$ will \textbf{diminish the significance of small singular values, effectively pushing them closer to zero}. 
This is reasonable in that the small singular values do not contribute significantly to the \textit{numerical rank} (Definition~\ref{definition:effective-rank-in-svd}). 
In situations where the singular values of matrices decay slowly, it may be necessary to repeat this process two or even three times:
$$
\bA\bA^\top(\bA\bA^\top \bA\bG ) =\bU\bSigma^5\bV^\top\bG.
$$
The procedure for computing the row ID of a matrix is presented in Algorithm~\ref{alg:row-id-randomize1}, where a parameter of $k$ is utilized to ensure  a high probability that $\bG$=$[\bg_1, \ldots,\bg_r, \bg_{r+1}, \ldots, \bg_{r+k}]$ contains the row basis of $\bA$ (a choice of $k=10$ is often good) and each $\bg_i$ is generated as a \textit{random Gaussian vector}. The iteration parameter $\eta$ is usually set to 1 or 2 to enhance accuracy (particularly in reducing the influence of small singular values). 
The complexity of the algorithm is outlined in the comment of each step, making it
\begin{equation}\label{equation:rnadom-rowid-complexity}
	\underline{\mathcalO(mn(r+k))} \,\, \text{flops.}
\end{equation}

\begin{remark}[A Word on the Source of Row Basis]\label{remark:source-row-basis}
In step 4 of Algorithm~\ref{alg:row-id-randomize1}, we aim for $\bD$ to contain linearly independent vectors that can span the column space of $\bA$, i.e., find a column basis. 
The method we present here for achieving this is known as the \textit{random projection method}, in which case we transform row basis matrix $\bG$ to the column basis matrix $\bD=\bA\bG$. 
The matrix $\bD$ formed by these columns is expected to be very close
to $\bA$ in a sense that the basis of the range of $\bD$ covers the range of $\bA$ well. The probability of failure is negligible.
However, further question can be posed, what if $\bG$ does not contain enough row basis vectors? 
In such cases, the column basis matrix $\bD$ can be chosen in different ways: by \textit{subsampling of the input matrix} directly (in the sense of row basis or column basis that we will show in next paragraph). 
The orthonormal basis consisting of $r$ linearly independent vectors can be obtained using exact methods since the size of $\bD$ or $\bG$ is very small. These techniques are relatively insensitive to the quality of randomness and yield highly accurate results. 
\end{remark}

\index{Gram–Schmidt}
\paragraph{Subsampling method for constructing basis.} 
We have shown in the Gram-Schmidt process (Section~\ref{section:project-onto-a-vector}), the projection of a vector $\ba$ onto $\bb$ results in the projection $\hat{\ba} = \frac{\ba^\top\bb}{\bb^\top\bb}\bb$, which lies in the direction of $\bb$ since it is a scalar multiple of $\bb$. 
The component of $\ba$ along $\bb$ thus can be obtained by $\ba-\hat{\ba} = \ba-\frac{\ba^\top\bb}{\bb^\top\bb}\bb$, which is $\bzero$ if $\ba$ and $\bb$ are linearly dependent (i.e., in the same direction). When $\ba,\bb\in \real^n$, the computational complexity is $\mathcalO(n)$ flops. 
Therefore, when  generating the $\bg_i$'s vectors in Algorithm~\ref{alg:row-id-randomize1}, we can project it onto each row of $\bA$ to check if it lies within the row space of $\bA$, and this operation costs $\mathcalO(mn)$ flops. 
And since there are $\sim (r+k)$ such vector generations, the overall cost of this check operation becomes $\underline{\mathcalO(mn(r+k))}$ flops to assure $\{\bg_1, \bg_2, \ldots, \bg_{r+k}\}$ can span the row space of $\bA$. 
This complexity is of the same order as the total algorithmic cost in  Equation~\eqref{equation:rnadom-rowid-complexity} when $n\leq m$. 
Analogously, we can directly generate the columns of $\bD=[\bd_1, \bd_2, \ldots, \bd_{r+k}]$ using Gaussian random vectors and check if each $\bd_i$ resides in the column space of $\bA$. 
This costs $\underline{\mathcalO(mn(r+k))}$, which again matches the cost of the algorithm in Equation~\eqref{equation:rnadom-rowid-complexity}. Therefore, the complexity of the algorithm would not blow up.

\begin{algorithm}[h] 
\caption{A \textcolor{mylightbluetext}{Randomized} Method to Compute the \textcolor{mylightbluetext}{Column} ID} 
\label{alg:column-id-randomize1} 
\begin{algorithmic}[1] 
\Require 
Rank-$r$ matrix $\bA$ with size $m\times n $; 
\State Decide the over-sampling parameter $k$ (e.g., $k=10$)$\rightarrow$ let $z=r+k$;
\State Decide the number of iterations $\eta$ (e.g., $\eta=0,1$ or $2$);
\State Generate $r+k$ Gaussian random vectors in $\textcolor{mylightbluetext}{\real^m}$ into columns of matrix $\bG\in \real^{\textcolor{mylightbluetext}{m}\times (r+k)}$;\Comment{i.e., probably contain the column basis of $\bA$}
\State Initialize $\bF\leftarrow\textcolor{mylightbluetext}{\bA^\top}\bG\in \real^{n\times (r+k)}$; \Comment{i.e., probably $\cspace(\bF)=\cspace(\bA^\top)$, $\textcolor{mylightbluetext}{(2m-1)nz}$ flops}
\For{$i=1$ to $\eta$}
\State $\bF \leftarrow \textcolor{mylightbluetext}{\bA^\top\bA}\bF$;  \Comment{$(2m-1)nz+(2n-1)mz$ flops}
\EndFor
\State Calc. the \textcolor{mylightbluetext}{column} ID of small $\bF$: $\bF=\bC_0\bW=\bF[:,J_s]\bW$;\Comment{$O\textcolor{mylightbluetext}{(mz^2)}$ flops by CPQR}
\State Output column ID of $\bA$: $\bA=\bA[:,J_s]\bW$;
\end{algorithmic} 
\end{algorithm}

\paragraph{Compute the column ID by randomized method.}
Analogously, the  procedure described in Algorithm~\ref{alg:row-id-randomize1}  for computing the row ID can be readily applied to compute the column ID of the matrix $\bA\in \real^{m\times n}$ by calculating the row ID of $\bA^\top$. 
Alternatively, we can derive it from the proof of the column ID. 
If $\bA$ admits a rank decomposition (Theorem~\ref{theorem:rank-decomposition}):
$$
\underset{m\times n}{\bA} = \underset{m\times r}{\bD}\gapthree \underset{r\times n}{\bF},
$$
where $\bD$ and $\bF$ have full column rank $r$ and full row rank $r$, respectively. 
Upon finding the rank decomposition, suppose further that the column ID of $\bF$ is obtained by 
$$
\text{column ID of $\bF$: }\qquad \underset{r\times n}{\bF} = \underset{r\times r}{\bC_0}\gapthree \underset{r\times n}{\bW}
=\bF[:,J_s]\bW.
$$
We observe that $ \bA[:,J_s]\bW$  automatically forms a column ID of $\bA$: $\bA=\bA[:,J_s]\bW$. This can be shown as follows:
\begin{equation}\label{equation:column-id-sub-d}
\begin{aligned}
\bA[:,J_s]\bW &= (\bD\bF[:,J_s])\bW &\qquad &\text{(since $\bA=\bD\bF$)}\\
&=\bD \bF=\bA. &\qquad &\text{(since $\bF=\bF[:,J_s]\bW$)}
\end{aligned}
\end{equation}
We notice that $\bF$ spans the same row space as $\bA$: $\cspace(\bF^\top)=\cspace(\bA^\top)$. 
The finding above indicates that as long as we find an $r\times n$ matrix $\bF$ spanning the same row space as $\bA$, a column ID can be applied to $\bF$ to find the column ID of $\bA$.
This reveals the randomized algorithm for computing the column ID of $\bA$. 
The underlying idea corresponds to the second part of the column ID proof in Section~\ref{section:proof-column-id}.

We may further notice that if the columns of $\bG$ comprise the column basis of $\bA$, then $\bA^\top\bG$ will contain the row basis of $\bA$, as indicated by Lemma~\ref{lemma:column-basis-from-row-basis} again. Moreover, $\bA^\top\bA\bA^\top\bG$ will span the row space of $\bA$, where the left multiplication by $\bA^\top\bA$ can help shrink small singular values:
$$
\bA^\top\bA\bA^\top\bG = (\bV\bSigma\bU^\top) (\bU\bSigma\bV^\top)(\bV\bSigma\bU^\top)\bG=\bV\bSigma^3\bU^\top\bG.
$$
The randomized procedure for computing the column ID is formulated in Algorithm~\ref{alg:column-id-randomize1}, with distinctions from the row ID highlighted in \textcolor{mylightbluetext}{blue} text.

%% file: chapter-id-bid.tex
\section{Computing  ID via Bayesian Approach}\label{section:bayes-id}

The Bayesian approach for computing the interpolative decomposition was proposed in \citet{lu2022bayesian, lu2022comparative}, demonstrating promising results in terms of low overfitting when computing low-rank approximation of column IDs. 
In this section, we will provide an overview of the derivation of this approach, and we will refer to the column ID as ID for simplicity without special mention.
The principle forms the basis for widely used learning algorithms such as those used for feature selection, finding quantitative strategies~\citep{lu2022feature}.
For a background on the Bayesian inference approach, including concepts such as conjugacy, Gibbs sampling, and prior distribution,  readers are encouraged to consult \citet{hoff2009first}

\index{Gamma distribution}
\index{Inverse-Gamma distribution}
\index{Gibbs sampling}
\index{Gibbs sampler}
\index{Truncated-normal distribution}
\index{General-truncated-normal distribution}
\subsection*{Probability Distributions}

To make the reading experience self-contained, we will first devote some time developing intuitions and notations in probability distributions in this section.

$\normal(x\mid \mu, \tau^{-1}) =\sqrt{\frac{ \tau}{2\pi}}\exp\{-\frac{\tau}{2} (x-\mu)^2\}$ is a Gaussian distribution with mean $\mu$ and precision $\tau$ (variance $\sigma^2=\tau^{-1}$).

\begin{figure}[h]
\centering  
\vspace{-0.15cm} 
\subfigtopskip=2pt 
\subfigbottomskip=2pt 
\subfigcapskip=-5pt 
\subfigure[Gamma distribution.]{\label{fig:dists_gamma}
\includegraphics[width=0.482\linewidth]{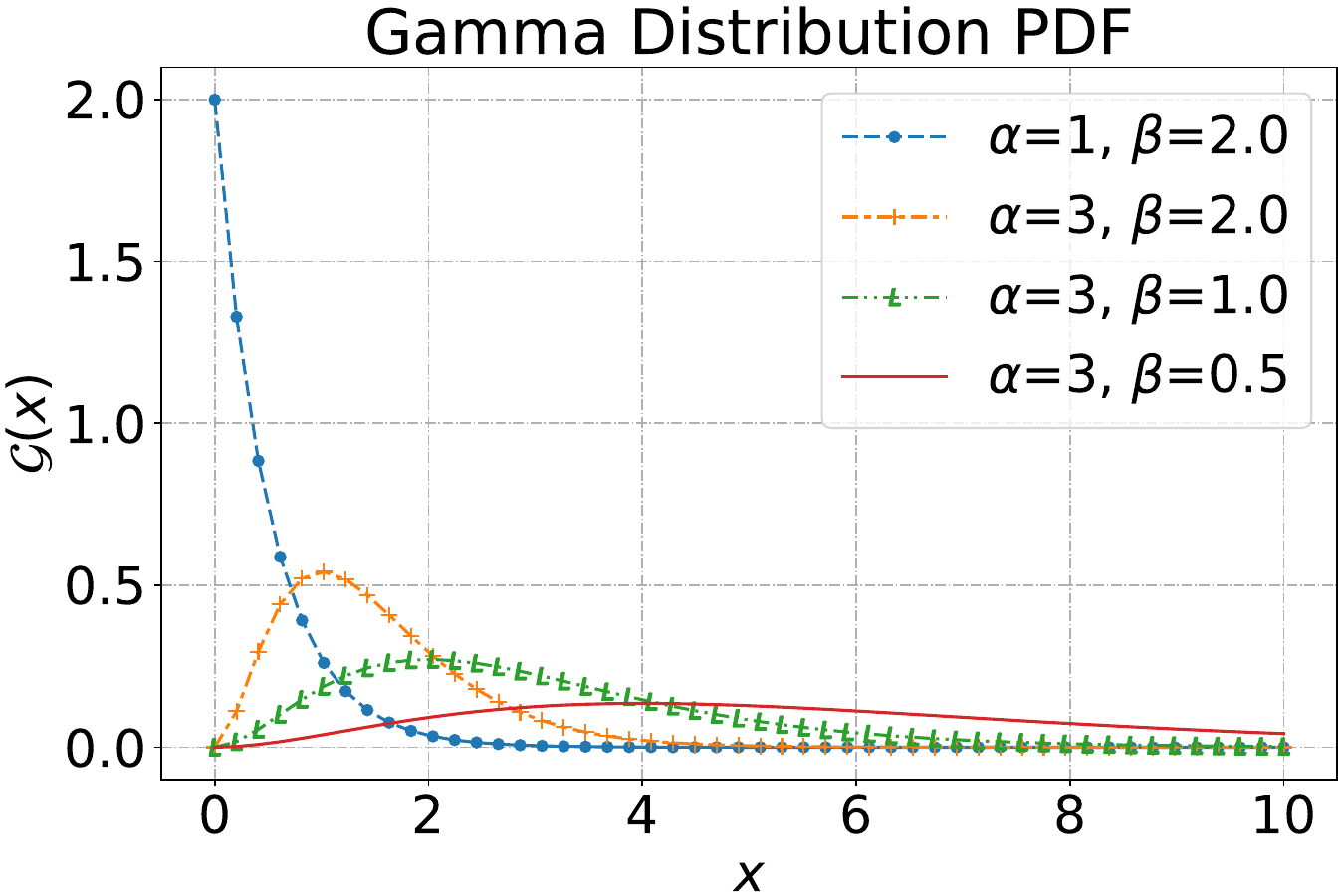}}
\subfigure[Inverse-Gamma distribution.]{\label{fig:dists_inversegamma}
\includegraphics[width=0.482\linewidth]{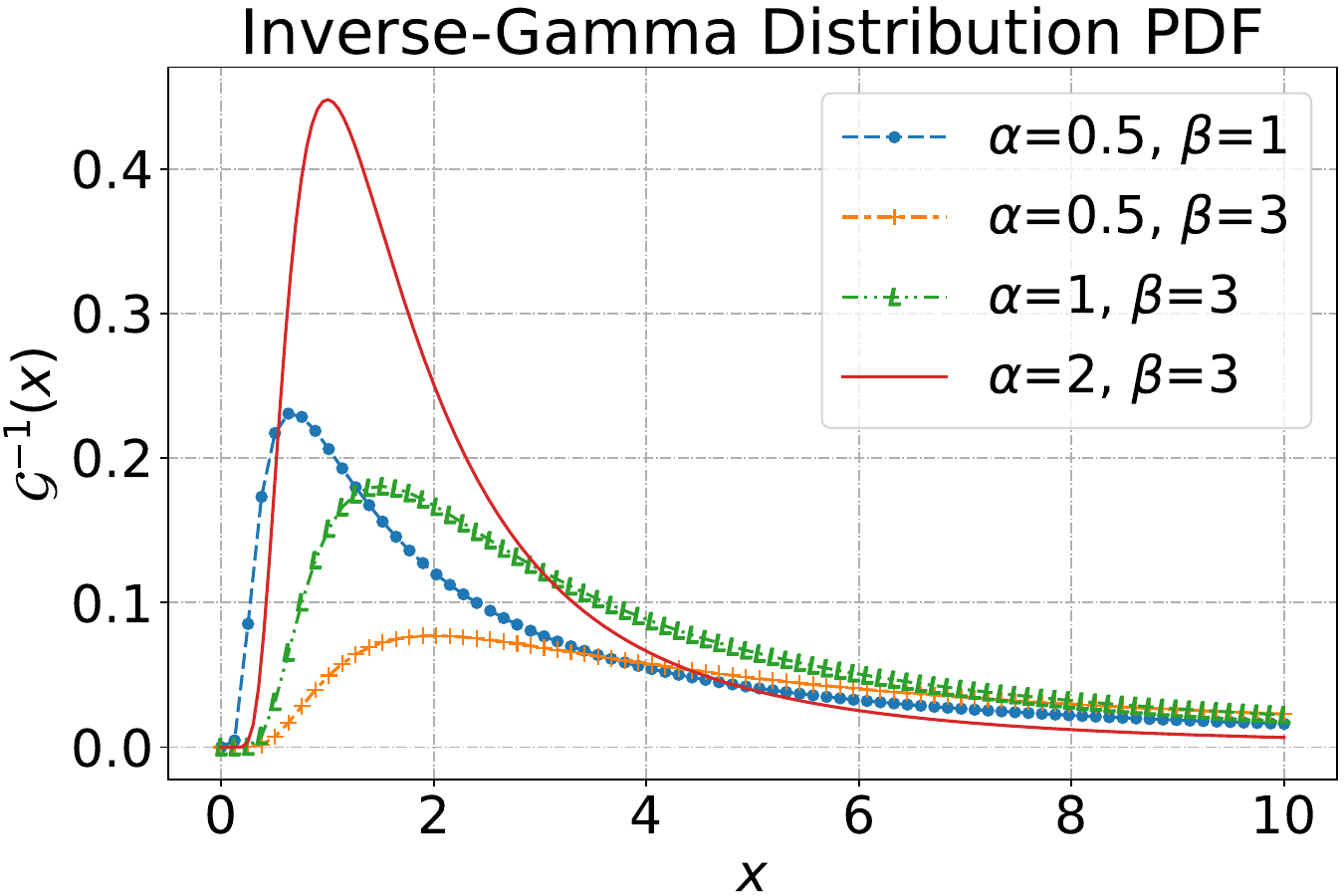}}
\caption{Gamma and inverse-Gamma probability density functions for different values of the parameters $\alpha$ and $\beta$.}
\label{fig:dists_gamma_inversegamma}
\end{figure}

$\gammadist(x\mid\alpha, \beta)= \frac{\beta^\alpha}{\Gamma(\alpha)} x^{\alpha-1}\exp\{-\beta x\}u(x)$ is a Gamma distribution, where $\Gamma(\cdot)$ is the Gamma function, and $u(x)$ represents the unit step function, which equals 1   when $x\geq0$ and 0 otherwise. Figure~\ref{fig:dists_gamma} compares different parameters $\alpha$ and $\beta$ for the Gamma distribution.

$\inversegammadist(x\mid\alpha, \beta)= \frac{\beta^\alpha}{\Gamma(\alpha)} x^{-\alpha-1}\exp\{-\frac{\beta}{x}\}u(x)$ is an inverse-Gamma distribution. 
Figure~\ref{fig:dists_inversegamma} compares different parameters $\alpha$ and $\beta$ for the inverse-Gamma distribution.

\begin{figure}[h]
\centering  
\vspace{-0.35cm} 
\subfigtopskip=2pt 
\subfigbottomskip=2pt 
\subfigcapskip=-5pt 
\subfigure[Truncated-normal.]{\label{fig:dists_truncatednorml}
\includegraphics[width=0.482\linewidth]{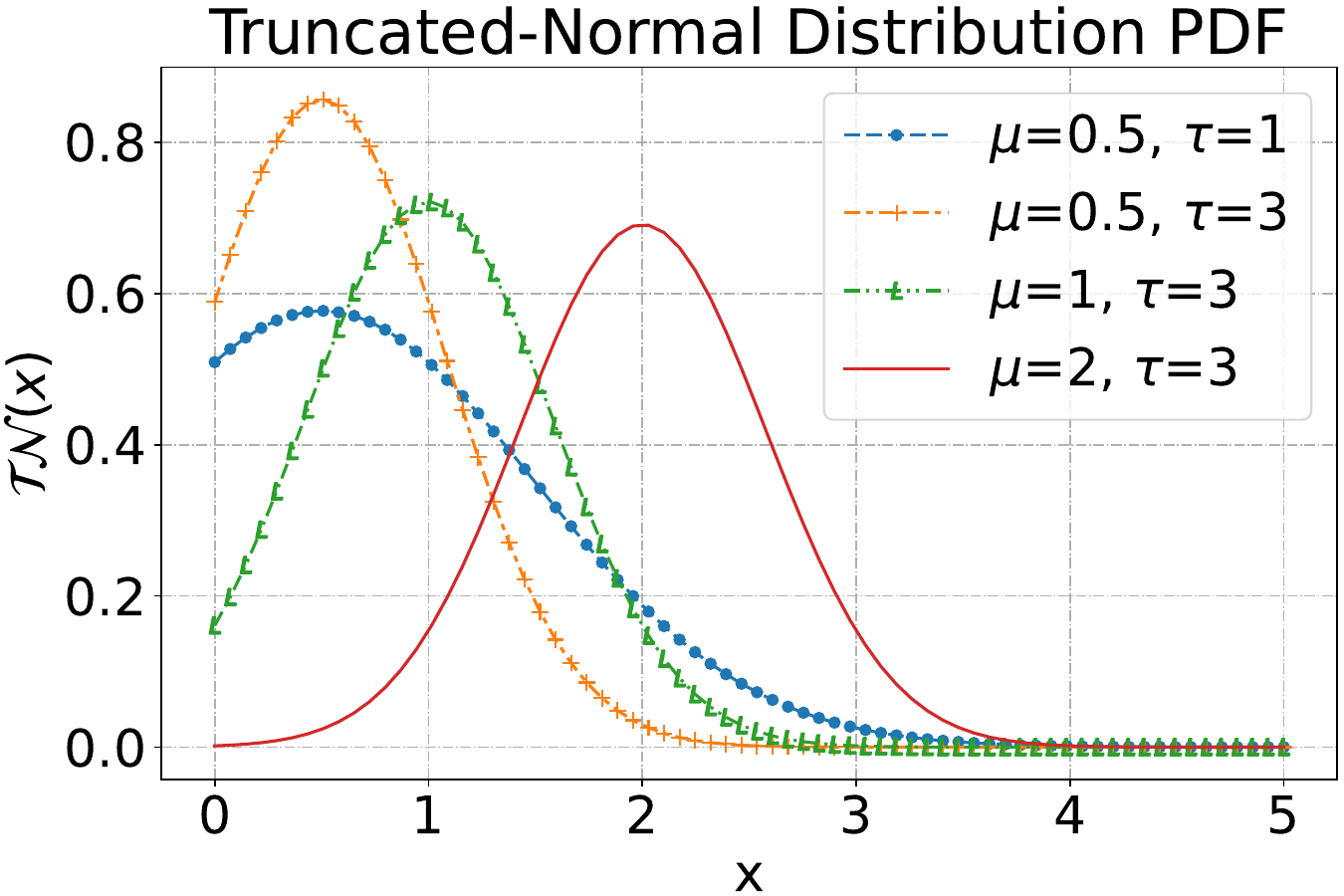}}
\subfigure[General-truncated-normal.]{\label{fig:dists_generaltruncatednorml}
\includegraphics[width=0.482\linewidth]{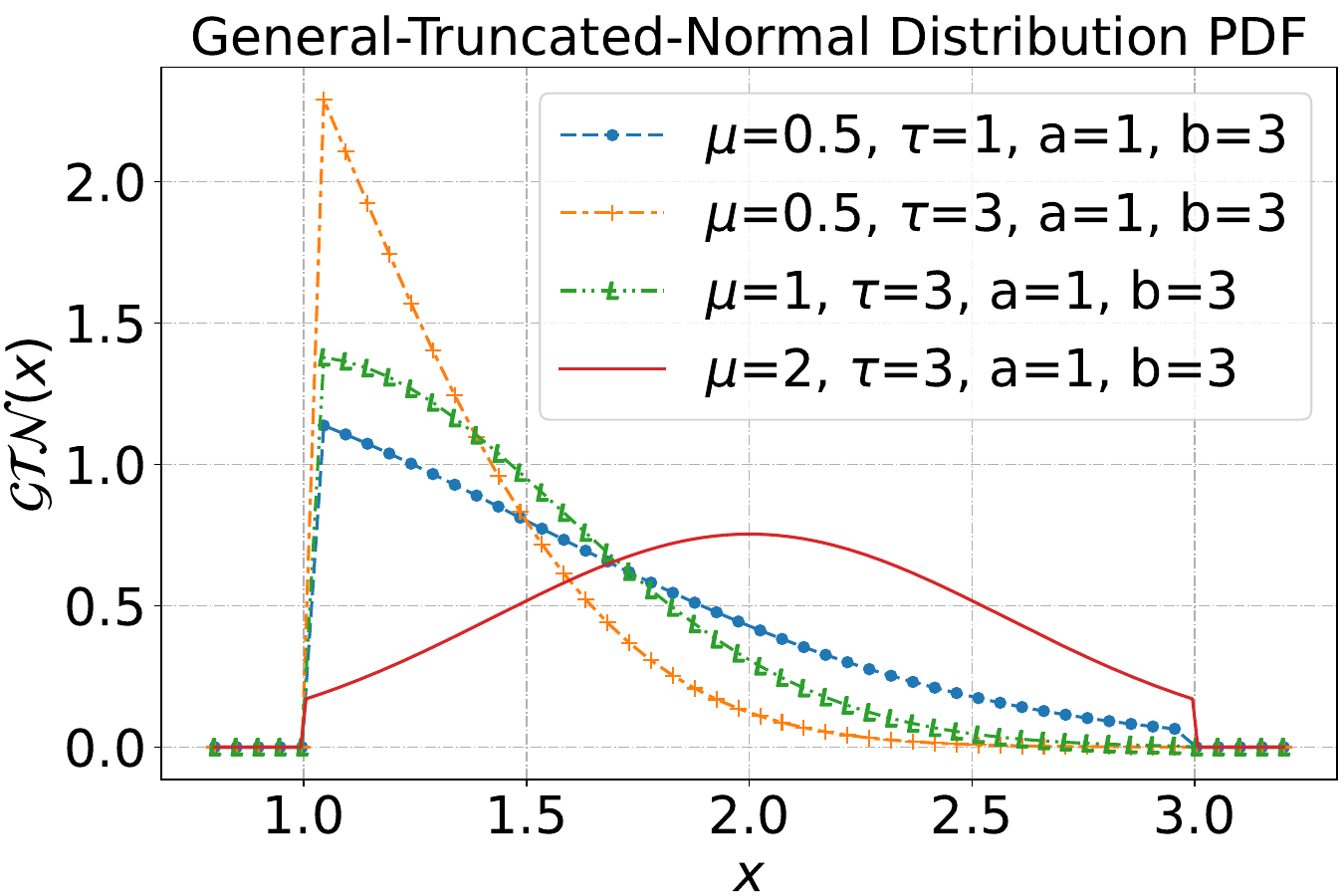}}
\caption{Truncated-normal and general-truncated-normal probability density functions for different values of the parameters $\mu$ and $\tau$.}
\label{fig:dists_truncatednorml_and_general}
\end{figure}

$\truncatednormal(x\mid\mu,\tau^{-1}) =\frac{\sqrt{\frac{\tau}{2\pi}} \exp\{-\frac{\tau}{2} (x-\mu)^2 \} } 
{1-\Phi(-\mu\sqrt{\tau})} u(x)$
is a truncated-normal (TN) with zero density below $x=0$ and being renormalized to integrate to one. Parameters $\mu$ and $\tau$ are known as the ``parent mean" and ``parent precision." $\Phi(\cdot)$ is the cumulative distribution function of the standard normal density $\normal(0,1)$. Figure~\ref{fig:dists_truncatednorml} compares different parameters $\mu$ and $\tau$ for the TN distribution.

$\generaltruncatednormal(x\mid\mu, \frac{1}{\tau}, a, b)=\frac{\sqrt{\frac{\tau}{2\pi}} \exp \{-\frac{\tau}{2}(x-\mu)^2  \}  }{\Phi((b-\mu)\cdot \sqrt{\tau})-\Phi((a-\mu)\cdot \sqrt{\tau})}$$u(x\mid a,b)$
is a general-truncated-normal (GTN) with zero density below $x=a$ or above $x=b$ and being renormalized to integrate to one, where $u(x\mid a,b)$ is a step function that has a value of 1 when $a\leq x\leq b$ and 0 otherwise. Similarly, parameters $\mu$ and $\tau$ are known as the ``parent mean" and ``parent precision" of the normal distribution. When $a=0$ and $b=\infty$, the GTN distribution reduces to a TN density. Figure~\ref{fig:dists_generaltruncatednorml} compares different parameters $\mu$ and $\tau$ for the GTN distribution.


\subsection*{Bayesian GBT and GBTN Models for ID}

\begin{figure}[h]
\centering  
\vspace{-0.35cm} 
\subfigtopskip=2pt 
\subfigbottomskip=6pt 
\subfigcapskip=-2pt 
\subfigure[GBT model.]{\includegraphics[width=0.3\textwidth]{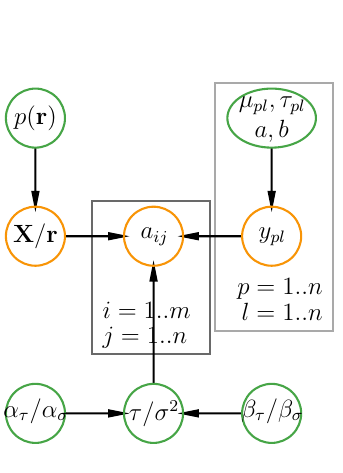} \label{fig:bmf_bid_GBT}}
\subfigure[GBTN model.]{\includegraphics[width=0.3\textwidth]{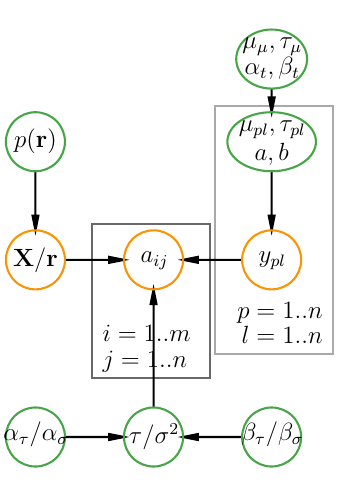} \label{fig:bmf_bid_GBTN}}
\caption{Graphical representation of GBT and GBTN models. Orange circles represent observed and latent variables, green circles denote prior variables, and plates represent repeated variables. The symbol ``/" in the variable represents ``or," and comma ``," in the variable represents ``and." Parameters $a$ and $b$ are fixed with $a=-1$ and $b=1$ in our case; while a weaker construction can set them to $a=-2$ and $b=2$.}
\label{fig:bmf_bids}
\end{figure}

To reiterate, the low-rank ID problem of an observed matrix $\bA$ can be stated as $\bA=\bC\bW+\bE$, where $\bA= [\ba_1, \ba_2, \ldots, \ba_n]\in \real^{m\times n}$ is approximately factorized into two matrices: an $m\times k$ matrix $\bC\in \real^{m\times k}$ containing $k$ basis columns of $\bA$, and a $k\times n$ matrix $\bW\in \real^{k\times n}$ with entries no larger than 1 in magnitude \footnote{A weaker construction is to assume no entry of $\bW$ has an absolute value greater than 2.}; the noise is characterized by the matrix $\bE\in \real^{m\times n}$. 
The task of training such models amounts to finding the best rank-$k$ approximation to the observed $m\times n$ target matrix $\bA$ under  a predefined loss function. 
To facilitate this, we introduce a binary \textit{state vector} $\br\in \{0,1\}^n$, where each element indicates the type of the corresponding column, basis column or interpolated (remaining) column: if $r_i=1$, then the $i$-th column $\ba_i$ is a basis column; if $r_i=0$, then $\ba_i$ is interpolated using the basis columns plus some error term. Suppose further $J_s$ is the set of  indices for the selected basis columns, $J_r$ is the set of  indices for the interpolated columns such that $J_s \cup J_r =\{1,2,\ldots, n\}$, and 
$$
J_s=J_s(\br)=\{i\mid r_i=1\}_{i=1}^n, \qquad 
J_r=J_r(\br)=\{i\mid r_i=0\}_{i=1}^n.
$$ 
Then, $\bC$ can be described as $\bC=\bA[:,J_s]$. 
The approximation $\bA\approx \bC\bW$ can be expressed equivalently as $\bA\approx\bC\bW\triangleq\bX\bY$, where $\bX\in \real^{m\times n}$ and $\bY\in \real^{n\times n}$, such that 
$$
\begin{aligned}
\bX[:,J_s]&=\bC\in \real^{m\times k}, \qquad \bX[:,J_r] = \bzero\in \real^{m\times (n-k)}, 
\qquad 
\bW = \bY[J_s,:]\in \real^{k\times n}.
\end{aligned}
$$ 
We also observe the presence of an identity  matrix $\bI\in \real^{k\times k}$ in both $\bW$ and $\bY$:
\begin{equation}\label{equation:submatrix_bid_identity}
\bI = \bW[:,J_s] = \bY[J_s,J_s].
\end{equation}
\begin{figure*}[h]
\centering  
\subfigtopskip=2pt 
\subfigbottomskip=9pt 
\subfigcapskip=-5pt 
\includegraphics[width=0.95\textwidth]{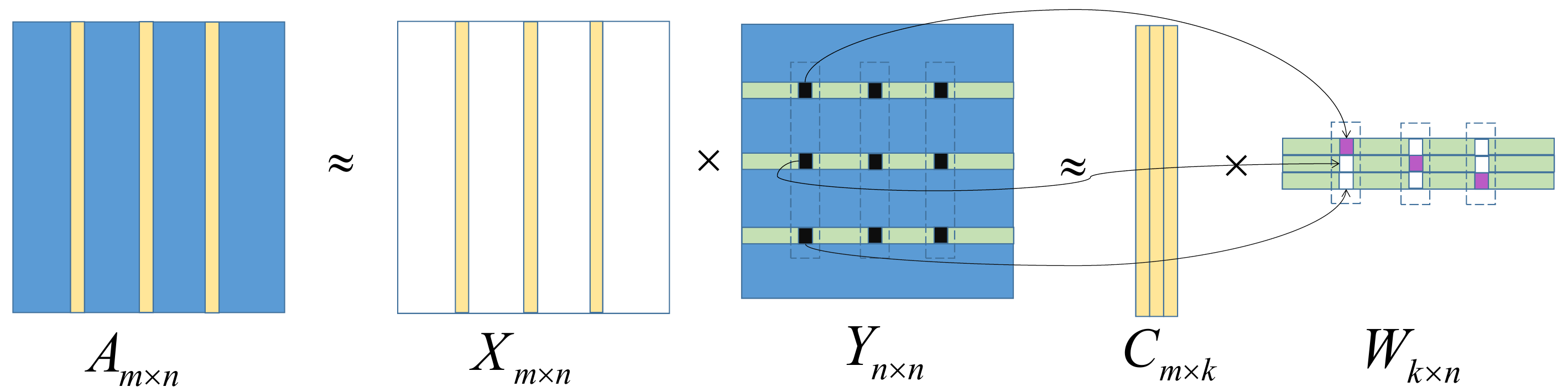}
\caption{Demonstration of the interpolative decomposition of a matrix, where  \textcolor{mydarkyellow}{yellow} vectors denote the basis columns of matrix $\bA$, white entries denote zero, \textcolor{mydarkpurple}{purple} entries denote one, \textcolor{mylightbluetext}{blue} and black entries denote elements that are not necessarily zero. The Bayesian ID models find the approximation $\bA\approx\bX\bY$, while the post-processing procedure finds the approximation $\bA\approx\bC\bW$.}
\label{fig:id-column}
\end{figure*}
To determine the low-rank ID of $\bA\approx\bC\bW$, we can transform the problem into finding the approximation  $\bA\approx\bX\bY$, with the state vector $\br$ recovering the submatrix $\bC$ (see Figure~\ref{fig:id-column}).
To assess the quality of this approximation, \textit{reconstruction error} measured by mean squared error (MSE or Frobenius norm) is minimized (assuming $k$ is known):
\begin{equation}\label{equation:idbid-per-example-loss}
\mathop{\min}_{\bW,\bZ} \,\, \frac{1}{mn}\sum_{j=1}^n \sum_{i=1}^{m} \left(a_{ij} - \bx_i^\top\by_j\right)^2,
\end{equation}
where $\bx_i$ and $\by_j$ represent the $i$-th \textbf{row} of $\bX$ and $j$-th \textbf{column} $\bY$, respectively (with $i\in\{1,2,\ldots,m\}, j\in\{1,2,\ldots,n\}$).
In the Bayesian method, we approach the magnitude constraint in $\bW$ and $\bY$ by considering the Bayesian ID model as a latent factor model. 
We establish a fully specified graphical model for the problem and employ Bayesian learning methods to infer the latent factors. In this regard, explicit magnitude constraints  on the latent factors are not necessary, as they are  naturally taken care of by the appropriate choice of the prior distribution; here we use the general-truncated-normal prior.

To be more concrete, we view the data $\bA$ as being produced according to the probabilistic generative process shown in Figure~\ref{fig:bmf_bids}. The observed $(i,j)$-th entry $a_{ij}$ of the matrix $\bA$ is modeled using a Gaussian likelihood function with variance $\sigma^2$ and a  mean given by the latent decomposition $\bx_i^\top\by_j$ (Equation~\eqref{equation:idbid-per-example-loss}):
\begin{equation}\label{equation:grrn_data_entry_likelihood}
\begin{aligned}
p(a_{ij} \mid  \bx_i^\top\by_j, \sigma^2) &= \normal(a_{ij}\mid \bx_i^\top\by_j, \sigma^2);\\
p(\bA\mid  \bX,\bY,\sigma^2) = \prod_{i,j=1}^{m,n}\normal \left(a_{ij}\mid  (\bX\bY)_{ij}, \sigma^2 \right) &
= \prod_{i,j=1}^{m,n} \normal \left(a_{ij}\mid  (\bX\bY)_{ij}, \tau^{-1} \right),
\end{aligned}
\end{equation}
where $\tau = (\sigma^2)^{-1}$ represents the precision parameter.
We choose a conjugate prior over the data variance, an inverse-Gamma distribution with shape $\alpha_\sigma$ and scale $\beta_\sigma$:
\begin{equation}\label{equation:prior_grrn_gamma_on_variance}
p(\sigma^2 \mid  \alpha_\sigma, \beta_\sigma) = \inversegammadist(\sigma^2 \mid  \alpha_\sigma, \beta_\sigma).
\end{equation}
Alternatively,  a conjugate Gamma prior (i.e., $\gammadist(\tau\mid \alpha_\tau, \beta_\tau)$) over the precision can be used; but we shall not repeat the details.

We treat the latent variables $y_{pl}$'s as random variables ($p,l\in\{1,2,\ldots, n\}$, see Figure~\ref{fig:bmf_bids}). And we need prior densities over these latent variables to express our beliefs about their values, e.g., constraint with magnitude smaller than 1 in this context though there are many other constraints (e.g., nonnegativity in \citet{lu2022flexible, lu2022robust}, semi-nonnegativity in \citet{ding2008convex}, and discreteness in \citet{gopalan2014bayesian}).
Furthermore, we make the additional assumption  that the latent variable $y_{pl}$'s are independently drawn from a general-truncated-normal prior:
\begin{equation}\label{equation:rn_prior_bidd}
p(y_{pl} \mid  \cdot ) = \generaltruncatednormal(y_{pl} \mid  \mu_{pl}, (\tau_{pl})^{-1}, a=-1, b=1).
\end{equation}
This prior serves to enforce the constraint on the component $\bY$ with no entry of $\bY$ having an absolute value greater than 1, and is conjugate to the Gaussian likelihood. The resulting posterior density  also follows a general-truncated-normal distribution.
In a less strict variant of the interpolative decomposition, the magnitude constraint can be relaxed to 2.
In such cases,  the prior can be adjusted flexibly by setting the parameters to $a=-2$ and $b=2$ accordingly.

\paragraph{Hierarchical prior.}
To further enhance flexibility, we choose a convenient joint hyperprior density over the parameters $\{\mu_{pl}, \tau_{pl}\}$ of the GTN prior in Equation~\eqref{equation:rn_prior_bidd}, namely, the \textit{GTN-scaled-normal-Gamma (GTNSNG) prior},
\begin{equation}
\begin{aligned}
&\gap p(\mu_{pl}, \tau_{pl} \mid \cdot) 
=\gtnsng(\mu_{pl}, \tau_{pl}\mid  \mu_\mu, (\tau_\mu)^{-1},\alpha_t, \beta_t)\\
&=	\big\{\Phi((b-\mu_\mu)\cdot \sqrt{\tau_\mu})-\Phi((a-\mu_\mu)\cdot \sqrt{\tau_\mu})\big\}\cdot 
\normal(\mu_{pl}\mid  \mu_\mu, (\tau_\mu)^{-1}) \cdot \gammadist(\tau_{pl} \mid  a_t, b_t).
\end{aligned}
\end{equation}
This prior can decouple the parameters $\{\mu_{pl}\}$ and $\{\tau_{pl}\}$, making their posterior conditional densities    normal, and Gamma, respectively, due to this convenient scale.

\paragraph{Terminology.} There are three key choices that determine the specific type of matrix decomposition model we use: the likelihood function, the priors placedover the factor matrices $\bW$ and $\bZ$, and whether we use any further hierarchical priors. 
Thus, we name the model based on the density functions in the order of the likelihood and priors.
For instance, if we select a Gaussian density as the likelihood function and opt for exponential and Gaussian density functions as the two prior density functions, respectively, 
then the model will be denoted as a Gaussian Exponential-Gaussian (GEG) model. 
At times, we introduce a hyperprior over the parameters of the prior density functions. 
For instance, if we place a Gamma prior over the Gaussian density, then it will further be termed as a Gaussian Exponential-Gaussian Gamma (GEGA) model. In this sense, the introduced hierarchical model for interpolative decomposition is named the \textit{GBT} and \textit{GBTN} model, where $B$ stands for Beta-Bernoulli density intrinsically.

\subsection*{Gibbs Sampler}

In this context, we employ Gibbs sampling, as it tends to be very accurate in finding the true posterior. 
While variational Bayesian inference is an alternative method, we will not delve into its details here; see \citet{lu2023bayesian} for more details. 
We shortly describe the posterior conditional densities in this section.

\paragraph{Conditional density of $y_{pl}$.}
The conditional density of $y_{pl}$ follows a general-truncated-normal density. Denoting all elements of $\bY$ except $y_{pl}$ by $\bY_{-pl}$, and 
following the graphical representation of the GBT (or the GBTN) model in Figure~\ref{fig:bmf_bids}, 
we can derive the conditional density of $y_{pl}$ as follows:
\begin{equation}\label{equation:posterior_gbt_ypl}
\begin{aligned}
&\gap p(y_{pl} \mid  \bA, \bX, \bY_{-pl}, \mu_{pl}, \tau_{pl}, \sigma^2) \propto  p(\bA\mid \bX,\bY, \sigma^2) \cdot p(y_{pl}\mid \mu_{pl}, \tau_{pl} )\\
&=\prod_{i,j=1}^{m,n} \normal \left(a_{ij}\mid  \bx_i^\top\by_j, \sigma^2 \right)\times
\generaltruncatednormal(y_{pl} \mid  \mu_{pl}, (\tau_{pl})^{-1},a=-1,b=1) \\
&\propto
\exp\bigg\{  
-\frac{1}{2\sigma^2} \sum_{i,j=1}^{m,n} (a_{ij} - \bx_i^\top \by_j)^2
\bigg\}
\exp \{-\frac{\tau_{pl}}{2}(y_{pl}-\mu_{pl})^2  \}
u(y_{pl} \mid  a,b)\\
&\propto
\exp\bigg\{  
-\frac{1}{2\sigma^2} \sum_{i}^{m} (a_{il} - \bx_i^\top \by_l)^2
\bigg\}
\exp \{-\frac{\tau_{pl}}{2}(y_{pl}-\mu_{pl})^2  \}
u(y_{pl} \mid  a,b)\\
&\propto
\exp\bigg\{  
-\frac{1}{2\sigma^2} \sum_{i}^{m} \bigg( x_{ip} ^2y_{pl }^2  + 2x_{ip} y_{pl } \big(\sum_{j\neq p}^{n}x_{ij} y_{jl}-a_{il}\big)\bigg)
\bigg\}
\exp \{-\frac{\tau_{pl}}{2}(y_{pl}-\mu_{pl})^2  \}
u(y_{pl} \mid  a,b)\\
&\propto
\exp\bigg\{  
-y_{pl }^2
\underbrace{\bigg(\frac{\sum_{i}^{m}  x_{ip} ^2}{2\sigma^2}+\textcolor{black}{\frac{\tau_{pl}}{2}} \bigg)}_{
\textcolor{mylightbluetext}{\frac{1}{2} \widetilde{\tau}}
}
+y_{pl } 
\underbrace{\bigg(\frac{1}{\sigma^2}  \sum_{i}^{m} x_{ip}  \big(a_{il}-\sum_{j\neq p}^{n}x_{ij}
y_{jl}\big)
+\textcolor{black}{\tau_{pl}\mu_{pl}}
\bigg)}_{\textcolor{mylightbluetext}{\widetilde{\tau} \widetilde{\mu}}}
\bigg\}
u(y_{pl} \mid  a,b)\\
&\propto \normal(y_{pl}\mid  \widetilde{\mu},( \widetilde{\tau})^{-1})u(y_{pl} \mid  a,b) 
\propto \generaltruncatednormal(y_{pl}\mid  \widetilde{\mu},( \widetilde{\tau})^{-1}, a=-1,b=1),
\end{aligned}
\end{equation}
where again, for simplicity, we assume the rows of $\bX$ are denoted by $\bx_i$'s and columns of $\bY$ are denoted by $\by_j$'s, $\widetilde{\tau} ={\sum_{i}^{m}  x_{ip} ^2}/{\sigma^2} +\tau_{pl}$ is the posterior ``parent precision" of the general-truncated-normal distribution, and 
$
\widetilde{\mu} = \big(\frac{1}{\sigma^2}  \sum_{i}^{m} x_{ip}  \big(a_{il}-\sum_{j\neq p}^{n}x_{ij}
y_{jl}\big)
+\textcolor{black}{\tau_{pl}\mu_{pl}}
\big) \big/ \widetilde{\tau}
$
is the posterior ``parent mean" of the general-truncated-normal distribution.

\paragraph{Conditional density of  state vector $\br$.}
Given the state vector $\br=[r_1,r_2, \ldots, r_n]^\top\in \real^n$, the relationship between $\br$
and the index sets $J_s$ and $J_r$ is simple: $J_s = J_s(\br) = \{p\mid r_p = 1\}_{p=1}^n$ and $J_r = J_r(\br) = \{p\mid r_p = 0\}_{p=1}^n$. 
When updating the state vector $\br$ to a new value, we select one index $j$ from the index set $J_s$ and another index $i$ from the index set $J_r$ (we note that $r_j=1$ and $r_i=0$ for the old state vector $\br$) such that 
\begin{equation}\label{equation:postrerior_gbt_rvector}
\begin{aligned}
o_j &= 
\frac{p(r_j=0, r_i=1\mid \bA,\sigma^2, \bY, \br_{-ji})}
{p(r_j=1, r_i=0\mid \bA,\sigma^2, \bY, \br_{-ji})}\\
&=
\frac{p(r_j=0, r_i=1)}{p(r_j=1, r_i=0)}\times
\frac{p(\bA\mid \sigma^2, \bY, \br_{-ji}, r_j=0, r_i=1)}{p(\bA\mid \sigma^2, \bY, \br_{-ji}, r_j=1, r_i=0)},
\end{aligned}
\end{equation}
where $\br_{-ji}$ denotes all elements of $\br$ except the $j$-th and $i$-th entries.
Trivially, we can set $p(r_j=0, r_i=1)=p(r_j=1, r_i=0)$. Then the full conditionally probability of $p(r_j=0, r_i=1\mid \bA,\sigma^2, \bY, \br_{-ji})$ can be calculated by 
\begin{equation}\label{equation:postrerior_gbt_rvector222}
p(r_j=0, r_i=1\mid \bA,\sigma^2, \bY, \br_{-ji}) = \frac{o_j}{1+o_j}.
\end{equation}

\paragraph{Conditional density of  variance $\sigma^2$.}
Finally, the conditional density of $\sigma^2$ follows an inverse-Gamma distribution due to conjugacy,
\begin{equation}\label{equation:posterior_gnt_sigma2}
\begin{aligned}
&\gap p(\sigma^2 \mid  \bX, \bY, \bA)
= \inversegammadist(\sigma^2 \mid  \widetilde{\alpha_\sigma}, \widetilde{\beta_\sigma}),
\end{aligned}
\end{equation}
where $\widetilde{\alpha_\sigma} = \frac{MN}{2}+\alpha_\sigma$, 
$\widetilde{\beta_\sigma}=\frac{1}{2} \sum_{i,j=1}^{M,N}(a_{ij}-\bx_i^\top\by_j)^2+\beta_\sigma$.

\paragraph{Extra update for GBTN model.}
Following the graphical representation of the GBTN model in Figure~\ref{fig:bmf_bids}, the conditional density of $\mu_{pl}$ can be derived as follows:
\begin{equation}\label{equation:posterior_gbt_mupl}
\begin{aligned}
&\gap p(\mu_{pl} \mid  \tau_{pl}, \mu_\mu, \tau_\mu, a_t, b_t, y_{pl})\\
&\propto \generaltruncatednormal(y_{pl} \mid  \mu_{pl}, (\tau_{pl})^{-1}, a=-1, b=1)
\cdot \gtnsng(\mu_{pl}, \tau_{pl}\mid  \mu_\mu, (\tau_\mu)^{-1},\alpha_t, \beta_t)\\
&\propto\generaltruncatednormal(y_{pl} \mid  \mu_{pl}, (\tau_{pl})^{-1}, a=-1, b=1)
\cdot  \\
&\gap \gap \big\{\Phi((b-\mu_\mu)\cdot \sqrt{\tau_\mu})-\Phi((a-\mu_\mu)\cdot \sqrt{\tau_\mu})\big\}
\cdot 
{\normal(\mu_{pl}\mid  \mu_\mu, (\tau_\mu)^{-1})} \cdot 
\cancel{\gammadist(\tau_{pl} \mid  a_t, b_t)}\\
&\propto 
\sqrt{\tau_{pl}}\cdot \exp\left\{ -\frac{\tau_{pl}}{2} (y_{pl}-\mu_{pl})^2\right\}
\cdot \exp\left\{ -\frac{\tau_\mu}{2}(\mu_\mu - \mu_{pl})^2  \right\}\\
&\propto \exp\bigg\{  
-\underbrace{\frac{\tau_{pl}+\tau_\mu}{2} }_{
\textcolor{mylightbluetext}{\frac{1}{2} \widetilde{t}}
}
\mu_{pl}^2 + \mu_{pl}
\underbrace{(\tau_{pl}y_{pl}+\tau_\mu\mu_\mu)}_{
\textcolor{mylightbluetext}{\widetilde{m}\cdot \widetilde{t}}
}  \bigg\}\propto 
\normal(\mu_{pl}\mid  \widetilde{m},(\,\widetilde{t}\,)^{-1}),
\end{aligned}
\end{equation}
where $\widetilde{t} = \tau_{pl}+\tau_\mu$ and $\widetilde{m} =(\tau_{pl}y_{pl}+\tau_\mu\mu_\mu)/\widetilde{t}$ are the posterior precision and mean of the normal density, respectively. 
Similarly, the conditional density of $\tau_{pl}$ is given by:
\begin{equation}\label{equation:posterior_gbt_taupl}
\begin{aligned}
&\gap p(\tau_{pl} \mid  \mu_{pl}, \mu_\mu, \tau_\mu, a_t, b_t, y_{pl})\\
&\propto \generaltruncatednormal(y_{pl} \mid  \mu_{pl}, (\tau_{pl})^{-1}, a=-1, b=1)
\cdot \gtnsng(\mu_{pl}, \tau_{pl}\mid  \mu_\mu, (\tau_\mu)^{-1},\alpha_t, \beta_t)\\
&\propto\generaltruncatednormal(y_{pl} \mid  \mu_{pl}, (\tau_{pl})^{-1}, a=-1, b=1)
\cdot \\
&\gap\gap \big\{\Phi((b-\mu_\mu)\cdot \sqrt{\tau_\mu})-\Phi((a-\mu_\mu)\cdot \sqrt{\tau_\mu})\big\}
\cdot 
\cancel{\normal(\mu_{pl}\mid  \mu_\mu, (\tau_\mu)^{-1})} \cdot 
{\gammadist(\tau_{pl} \mid  a_t, b_t)}\\
&\propto \exp\left\{  -\tau_{pl}  \frac{(y_{pl}- \mu_{pl})^2}{2}  \right\}
\tau_{pl}^{1/2} \tau_{pl}^{a_t-1} \exp\left\{  -b_t \tau_{pl} \right\}\\
&\propto \exp\left\{   -\tau_{pl}\left[ b_t +  \frac{(y_{pl}- \mu_{pl})^2}{2}  \right] \right\}
\cdot \tau_{pl}^{(a_t+1/2)-1}
\propto \gammadist(\tau_{pl} \mid  \widetilde{a}, \widetilde{b}),
\end{aligned}
\end{equation}
where $\widetilde{a} = a_t+1/2$ and $\widetilde{b}=b_t +  \frac{(y_{pl}- \mu_{pl})^2}{2}$ are the posterior parameters of the Gamma density.
The complete procedure is then presented in Algorithm~\ref{alg:gbtn_gibbs_sampler} \footnote{For a tutorial on the Gibbs sampling procedure, we refer to \citet{hoff2009first}.}.

\begin{algorithm}[tb] 
\caption{Gibbs sampler for GBT and GBTN ID models. The procedure presented here may not be efficient but is explanatory. A more efficient one can be implemented in a vectorized manner. By default, uninformative priors are $a=-1, b=1,\alpha_\sigma=0.1, \beta_\sigma=1$, ($\{\mu_{pl}\}=0, \{\tau_{pl}\}=1$) for GBT, ($\mu_\mu =0$, $\tau_\mu=0.1, \alpha_t=\beta_t=1$) for GBTN.} 
\label{alg:gbtn_gibbs_sampler}  
\begin{algorithmic}[1] 
\For{$t=1$ to $T$}
\State Sample state vector $\br$ from Equation~\eqref{equation:postrerior_gbt_rvector222};
\State \algoalign{Update matrix $\bX$ by $\bA[:,J_s]$, where the index vector $J_s$ is the index of $\br$ with value 1, and set $\bX[:,J_r]=\bzero$, where the index vector $J_r$ is the index of $\br$ with value 0;
 }
\State Sample $\sigma^2$ from $p(\sigma^2 \mid  \bX,\bY, \bA)$ in Equation~\eqref{equation:posterior_gnt_sigma2};
\For{$p=1$ to $n$} 
\For{$l=1$ to $n$} 
\State Sample $y_{pl}$ from Equation~\eqref{equation:posterior_gbt_ypl};
\State (GBTN only) Sample $\mu_{pl}$ from Equation~\eqref{equation:posterior_gbt_mupl};
\State (GBTN only) Sample $\tau_{pl}$ from Equation~\eqref{equation:posterior_gbt_taupl};
\EndFor
\EndFor
\State Report loss in Equation~\eqref{equation:idbid-per-example-loss}, stop if it converges.
\EndFor
\State Report mean loss in Equation~\eqref{equation:idbid-per-example-loss} after burn-in iterations.
\end{algorithmic} 
\end{algorithm}

\begin{algorithm}[h] 
\caption{\textit{Aggressive} Gibbs sampler for GBT ID model. The procedure presented here may not be efficient but is explanatory. A more efficient one can be implemented in a vectorized manner. By default, uninformative priors are $a=-1, b=1,\alpha_\sigma=0.1, \beta_\sigma=1$, ($\{\mu_{pl}\}=0, \{\tau_{pl}\}=1$) for GBT. 
} 
\label{alg:gbtn_gibbs_sampler_aggressive}  
\begin{algorithmic}[1] 
\For{$t=1$ to $T$}
\State Sample state vector $\br$ from $\{\br_1, \br_2\}$ by Equation~\eqref{equation:postrerior_gbt_rvector222};
\State Decide $\bY$: $\bY=\bY_1$ if $\br$ is $\br_1$; $\bY=\bY_2$ if $\br$ is $\br_2$;
\State Update state vector $\br_1=\br$; 
\State Sample proposal state vector $\br_2$ based on $\br$;
\State Update matrix $\bX$ by $\br=\br_1$; 
\State Update proposal $\bX_2$ by $\br_2$;
\State Sample $\sigma^2$ from $p(\sigma^2 \mid  \bX,\bY, \bA)$ in Equation~\eqref{equation:posterior_gnt_sigma2}; 
\State Sample $\bY_1 = \{y_{pl}\}$ using $\bX$;
\State Sample $\bY_2 = \{y_{pl}\}$ using $\bX_2$;
\State Report loss in Equation~\eqref{equation:idbid-per-example-loss}, stop if it converges.
\EndFor
\State Report mean loss in Equation~\eqref{equation:idbid-per-example-loss} after burn-in iterations.
\end{algorithmic} 
\end{algorithm}

\subsection*{Aggressive Update}
In Algorithm~\ref{alg:gbtn_gibbs_sampler}, we notice that we set $\bX[:,J_r]=\bzero$ when the new state vector $\br$ is sampled. However, in the subsequent iteration, the index set $J_r$ may be updated, in which case one entry $i$ of $J_r$ may be altered to have a value of 1: 
$$
r_i=0 \rightarrow r_i=1.
$$
This could potentially cause problems in the update of $y_{pl}$ in Equation~\eqref{equation:posterior_gbt_ypl}, since a zero $i$-th column in $\bX$ cannot effectively update $\bY$. One solution is to record a proposal state vector $\br_2$ and a corresponding proposal factor matrix $\bX_2$ based on this $\br_2$ vector. 
When the update in the next iteration selects the old state vector $\br$, the factor matrix $\bX$ is adopted to finish the updates; while the algorithm chooses the proposal state vector $\br_2$, the proposal factor matrix $\bX_2$ is applied instead to do the updates. We call this procedure the \textit{aggressive} update.
The \textit{aggressive} sampler for the GBT model is formulated in Algorithm~\ref{alg:gbtn_gibbs_sampler_aggressive}. 
For simplicity, we omit the sampler for the GBTN model, as it can be implemented similarly.

\begin{algorithm}[h] 
\caption{Gibbs sampler for GBT and GBTN ID with \textbf{ARD} models.  The procedure presented here can be inefficient but is explanatory. While a vectorized manner can be implemented to find a more efficient algorithm. By default, weak priors are $a=-1, b=1,\alpha_\sigma=0.1, \beta_\sigma=1$, ($\{\mu_{pl}\}=0, \{\tau_{pl}\}=1$) for GBT, ($\mu_\mu =0$, $\tau_\mu=0.1, \alpha_t=\beta_t=1$) for GBTN. \textcolor{mylightbluetext}{Number of critical steps: $\nu$}.} 
\label{alg:gbtn_gibbs_sampler_withard}  
\begin{algorithmic}[1] 
\For{$t=1$ to $T$}
\For{\textcolor{mylightbluetext}{$j=1$ to $n$}}
\State \textcolor{mylightbluetext}{Sample state vector element $r_j$ from Equation~\eqref{equation:postrerior_gbt_rvector222_ard}};
\EndFor
\State \algoalign{Update matrix $\bX$ by $\bA[:,J_s]$, where the index vector $J_s$ is the index of $\br$ with value 1, and set $\bX[:,J_r]=\bzero$, where the index vector $J_r$ is the index of $\br$ with value 0;}
\State Sample $\sigma^2$ from $p(\sigma^2 \mid  \bX,\bY, \bA)$ in Equation~\eqref{equation:posterior_gnt_sigma2}; 
\For{\textcolor{mylightbluetext}{$n=1$ to $\nu$}}
\For{$p=1$ to $n$} 
\For{$l=1$ to $n$} 
\State Sample $y_{pl}$ from Equation~\eqref{equation:posterior_gbt_ypl};
\State (GBTN only) Sample $\mu_{pl}$ from Equation~\eqref{equation:posterior_gbt_mupl};
\State (GBTN only) Sample $\tau_{pl}$ from Equation~\eqref{equation:posterior_gbt_taupl};
\EndFor
\EndFor
\EndFor
\State Output loss in Equation~\eqref{equation:idbid-per-example-loss}, stop iteration if it converges.
\EndFor
\State Output averaged loss in Equation~\eqref{equation:idbid-per-example-loss} for evaluation after burn-in iterations.
\end{algorithmic} 
\end{algorithm}

\subsection*{Post-Processing}
The Gibbs sampling algorithm finds the approximation $\bA\approx \bX\bY$, where $\bX\in \real^{m\times n}$ and $\bY\in \real^{n\times n}$. As stated above, the redundant columns in $\bX$ and redundant rows in $\bY$ can be eliminated  using the index vector $J_s$: 
$$
\begin{aligned}
\bC&=\bX[:,J_s]=\bA[:,J_s]
\quad \text{and}\quad 
\bW= \bY[J_s,:].
\end{aligned}
$$
Since the submatrix $\bY[J_s,J_s]=\bW[:,J_s]$ (see Equation~\eqref{equation:submatrix_bid_identity}), generated during the Gibbs sampling procedure, is not inherently constrained to be an identity matrix (however, this is required in the interpolative decomposition), we must manually set it to an identity matrix. 
This will marginally reduce the reconstruction error (see Problem~\ref{problem:id1}).
The post-processing procedure is illustrated in Figure~\ref{fig:id-column}.

\index{ARD}
\index{Bayesian inference}
\section{Computing  ID via Bayesian Approach with ARD}
In the GBT and GBTN models, it is necessary to predefine the latent dimension $k$. 
However, the \textit{automatic relevance determination (ARD)} methods  extend the Bayesian models to eliminate the need for model selection.

Given the state vector $\br=[r_1,r_2, \ldots, r_n]^\top\in \real^n$ whose index sets are $J_s = J_s(\br) = \{p\mid r_p = 1\}_{p=1}^n$ and $J_r = J_r(\br) = \{p\mid r_p = 0\}_{p=1}^n$. A new value of the state vector $\br$ is to select one index $j$ from either the index set $J_s$ or the index set $J_r$, i.e., $j\in J_s\cup J_r$,
such that 
\begin{equation}\label{equation:postrerior_gbt_rvector_ard}
\begin{aligned}
o_j &= 
\frac{p(r_j=0\mid \bA,\sigma^2, \bY, \br_{-j})}
{p(r_j=1\mid \bA,\sigma^2, \bY, \br_{-j})}=
\frac{p(r_j=0)}{p(r_j=1)} \times
\frac{p(\bA\mid \sigma^2, \bY, \br_{-j}, r_j=0)}{p(\bA\mid \sigma^2, \bY, \br_{-j}, r_j=1)},
\end{aligned}
\end{equation}
where $\br_{-j}$ denotes all elements of $\br$ except the $j$-th element.
Comparing Equation~\eqref{equation:postrerior_gbt_rvector_ard} with Equation~\eqref{equation:postrerior_gbt_rvector}, we may find that in the former equation, the number of selected columns is no longer fixed. 
Therefore, we let the inference determine the number of columns in the basis matrix $\bC$ for interpolative decomposition.
Once again, we can set $p(r_j=0)=p(r_j=1)=0.5$. Then the full conditionally probability of $p(r_j=0, r_i=1\mid \bA,\sigma^2, \bY, \br_{-ji})$ can be calculated as:
\begin{equation}\label{equation:postrerior_gbt_rvector222_ard}
p(r_j=0\mid \bA,\sigma^2, \bY, \br_{-j}) = \frac{o_j}{1+o_j}.
\end{equation}
The complete algorithm for GBT and GBTN with ARD is outlined in Algorithm~\ref{alg:gbtn_gibbs_sampler_withard}, where the difference lies in that we need to iterate over all elements of the state vector rather than just one or two elements. We aware that many elements in the state vector can change their signs, making the update of matrix $\bY$ (see Figure~\ref{fig:id-column}) unstable. 
To address this, we also define a value of \textit{critical steps $\nu$}: after sampling the entire state vector $\br$, we update  the matrix $\bY$ and its related parameters multiple times (here we repeat $\nu$ times). 
The differences are highlighted in blue in Algorithm~\ref{alg:gbtn_gibbs_sampler_withard}.

\begin{problemset}
\index{Cramer's rule}
\item \label{prob:cramer_adj_1} \textbf{Cramer's rule.}	Consider the linear system $\bM\bx=\bl$, where $\bM\in\real^{n\times n}$, and $\bx,\bl\in\real^{n}$. Let $\bM_{\bl}(i)$ be the matrix formed by replacing the $i$-th column of $\bM$ with $\bl$. Show that the $i$-th element of $\adjugate(\bM)\bl\in\real^n$ ($\adjugate(\bM)$ is the adjugate of $\bM$; see Definition~\ref{definition:adjugate}) is 
\begin{equation}\label{equation:cramer_adj_1}
\big(\adjugate(\bM)\bl\big)_i = \det(\bM_{\bl}(i)), \gap i\in\{1,2,\ldots,n\}.
\end{equation}
Now consider the linear system $\bM\bX=\bL$, where $\bM\in\real^{n\times n}$, and $\bX,\bL\in\real^{n\times m}$. Let $\bM_{\bL}(i,j)$ be the matrix formed by replacing the $i$-th column of $\bM$ with the $j$-th column $\bl_j$ of $\bL$. Show that the $(i,j)$-th element of $\adjugate(\bM)\bL\in\real^{n\times m}$ is 
\begin{equation}\label{equation:cramer_adj_2}
\big(\adjugate(\bM)\bL\big)_{ij} = \det(\bM_{\bL}(i,j)), \gap i\in\{1,2,\ldots,n\}, j\in\{1,2,\ldots,m\}.
\end{equation}
\textit{Hint: Examine the definitions of determinant and adjugate (Definition~\ref{definition:determinant},~\ref{definition:adjugate}).}

\item \label{prob:cramer_adj_2} \textbf{Cramer's rule.} Consider the same setting as Problem~\ref{prob:cramer_adj_1} and denote the vector  $\adjugate(\bM)\bl\in\real^n$ and  the matrix 
$\adjugate(\bM)\bL\in\real^{n\times m}$ by
$$
\adjugate(\bM)\bl = \big[\det(\bM_{\bl}(i))\big]_{i=1}^n\in\real^n,
\gap 
\adjugate(\bM)\bL = \big[\det(\bM_{\bl}(i,j))\big]_{i,j=1}^{n,m}\in\real^{n\times m},
$$
respectively; i.e., the $i$-th element of the vector is $\det(\bM_{\bl}(i))$,  the  $(i,j)$-th element of the matrix is $\det(\bM_{\bL}(i,j))$. 
Show that
\begin{equation}\label{equation:cramer_adj_3}
\begin{aligned}
\bM \big[\det(\bM_{\bl}(i))\big]_{i=1}^n &= \bM \adjugate(\bM) \bl=\det(\bM)\bl;\\
\bM \big[\det(\bM_{\bL}(i,j))\big]_{i,j=1}^{n,m} &= \bM \adjugate(\bM) \bL=\det(\bM)\bL.
\end{aligned}
\end{equation}
\textit{Hint: Examine the definition of adjugate (Definition~\ref{definition:adjugate}) and Exercise~\ref{exercise:pro_adjug}.}

\item \label{prob:cramer_adj_3} \textbf{Cramer's rule.} Consider the same setting as Problem~\ref{prob:cramer_adj_1} and assume further that $\bM$ is nonsingular. Show that the $i$-th element of the solution $\bx$ is 
\begin{equation}\label{equation:cramer_adj_res1}
x_i = \frac{\det(\bM_{\bl}(i))}{\det(\bM)}, \gap \forall i\in\{1,2,\ldots,n\}.
\end{equation}
Similarly, show that the $(i,j)$-th element of the solution $\bX$ is 
\begin{equation}\label{equation:cramer_adj_res2}
x_{ij} = \frac{\det\left(\bM_{\bL}(i,j)\right)}{\det(\bM)}, \gap \forall i\in\{1,2,\ldots,n\}, j\in\{1,2,\ldots,m\}.
\end{equation}
That is,  Cramer's rule.

\item \label{prob:cramer_adj_4} \textbf{Cramer's rule: the simple way.}
Consider the same setting as Problem~\ref{prob:cramer_adj_1} and assume further $\bM$ is nonsingular. We notice that 
\begin{equation}\label{equation:cramer_adj_4}
\bM \bI_{\bl}(i) = \bM_{\bl}(i),  \gap \forall i\in\{1,2,\ldots,n\},
\end{equation}
where $\bI_{\bl}(i)$ represents the identity matrix with  the $i$-th column replaced by $\bl$. 
Taking the determinant yields:
\begin{equation}\label{equation:cramer_adj_5}
	\det(\bM) \det(\bI_{\bl}(i)) = \det(\bM_{\bl}(i)).
\end{equation}
Show that  $\det(\bI_{\bl}(i))=x_i$, thereby verifying the result in \eqref{equation:cramer_adj_res1}.

\item \label{prob:cramer_adj_5} \textbf{Determinant of inverses for subsets, Jacobi's equality.} Consider $\bM\in\real^{n\times n}$ and two index sets $I,J\subseteq\{1,2,\ldots,n\}$ (their complementary sets are $I^C$ and $J^C$, respectively). 
Show that 
\begin{equation}
\det\big(\bM^{-1}[I^C,J^C]\big) = (-1)^{\gamma} \frac{\det(\bM[J,I])}{\det(\bM)},
\end{equation}
where $\gamma=\sum_{i\in I} i +\sum_{j\in J}j$ is the sum of indices.
When $I=J$, this also indicates
\begin{equation}
	\det\big(\bM^{-1}[I^C,I^C]\big) = \frac{\det(\bM[I,I])}{\det(\bM)},
\end{equation}
which is known as \textit{Jacobi's equality}.
\textit{Hint: Examine the definitions of determinant and adjugate (Definition~\ref{definition:determinant},~\ref{definition:adjugate}). We may also prove the latter using Schur complement; Exercise~\ref{exercise:jacob_eq}.}

\item Determine the column ID for the matrix 
$$
\bA = 
\begin{bmatrix}
1 & 3 & 2 \\
3 & 7 & 6 \\
4 & 5 & 8
\end{bmatrix}.
$$

\item \label{problem:mtb_dis} \textbf{Matlab-style notation.} Consider a rectangular matrix $\bA\in\real^{m\times n}$ with rank $r$, which admits a  rank decomposition $\bA=\bD\bF$; and consider index sets $I,J\subseteq \{1,2,\ldots,m\}$ and $K,L\subseteq\{1,2,\ldots,n\}$ with cardinality $r$. Then, we have  $\bA[I,K]=\bD[I,:]\bF[:,K]$.
Show that 
\begin{itemize}
	\item $\bA[I,K]$ is nonsingular if and only if $\rank(\bD[I,:])=\rank(\bF[:,K])=r$.
	\item $\det(\bA[I,K]) \det(\bA[J,L])=\det(\bA[I,L])\det(\bA[J,K])$.
\end{itemize}

\item \label{problem:id1} Identify an example that demonstrates the effectiveness of the post-processing method described in Section~\ref{section:bayes-id} in reducing reconstruction errors.

\item \textbf{Conjugacy, Beta-Bernoulli model \citep{lu2023bayesian}.} Given the Beta distribution $\betadist(\theta \mid a, b)=\frac{1}{B(a,b)} \theta^{a-1} (1-\theta)^{b-1}\indicator(0\leq \theta\leq 1)$ and the Bernoulli distribution $\bernoulli(x\mid \theta)=\theta^x(1-\theta)^x\indicator(x\in\{0,1\})$. Discuss the form of the posterior distribution of $\theta$ if we assume the prior is $\betadist(\theta \mid a,b)$ and the likelihood is $\bernoulli(x\mid \theta)$.

\end{problemset}

%% file: chapter-hessenberg.tex
\part{Reduction to Hessenberg, Tridiagonal, and Bidiagonal Form}
\index{Similarity transformation}
\section*{Introduction}
\lettrine{\color{caligraphcolor}I}
In practical applications, it is often necessary to decompose a matrix $\bA$ into two orthogonal matrices, denoted by $\bA = \bQ\bLambda\bQ^\top$, where $\bLambda$ is either diagonal or upper triangular, and $\bQ$ is orthogonal.
This type of decomposition is widely used in various fields,  including eigenanalysis through  Schur decomposition and principal component analysis (PCA) via spectral decomposition.
The process of computing this decomposition involves a series of \textit{orthogonal similarity transformations}:
$$
\underbrace{\bQ_k^\top\ldots \bQ_2^\top \bQ_1^\top}_{\triangleq\bQ^\top} \bA \underbrace{\bQ_1\bQ_2\ldots\bQ_k}_{\triangleq\bQ},
$$
which ultimately  converges to $\bLambda$ when $k\rightarrow\infty$ or $k$ becomes largely enough. However, performing these transformations can be challenging in practice, especially when using methods like Householder reflectors. 
For instance, after applying the QR decomposition using Householder reflectors (where $\boxtimes$ represents a value that is not necessarily zero, and \textbf{boldface} indicates a value has just been changed), the sequence of orthogonal similarity transformations can be constructed through the use of Householder reflectors: 
$$
\begin{aligned}
	\begin{sbmatrix}{\bA}
		\boxtimes & \boxtimes & \boxtimes & \boxtimes& \boxtimes \\
		\boxtimes & \boxtimes & \boxtimes & \boxtimes& \boxtimes \\
		\boxtimes & \boxtimes & \boxtimes & \boxtimes& \boxtimes \\
		\boxtimes & \boxtimes & \boxtimes & \boxtimes & \boxtimes\\
		\boxtimes & \boxtimes & \boxtimes & \boxtimes& \boxtimes
	\end{sbmatrix}
	&\stackrel{\bH_1\times }{\rightarrow}
	\begin{sbmatrix}{\bH_1\bA}
		\bm{\boxtimes} & \bm{\boxtimes} & \bm{\boxtimes} & \bm{\boxtimes}& \bm{\boxtimes} \\
		\bm{0} & \bm{\boxtimes} & \bm{\boxtimes} & \bm{\boxtimes}& \bm{\boxtimes} \\
		\bm{0} & \bm{\boxtimes} & \bm{\boxtimes} & \bm{\boxtimes}& \bm{\boxtimes} \\
		\bm{0} & \bm{\boxtimes} & \bm{\boxtimes} & \bm{\boxtimes} & \bm{\boxtimes}\\
		\bm{0} & \bm{\boxtimes} & \bm{\boxtimes} & \bm{\boxtimes}& \bm{\boxtimes}
	\end{sbmatrix}
	\stackrel{\times \bH_1^\top }{\rightarrow}
	\begin{sbmatrix}{\bH_1\bA\bH_1^\top}
		\bm{\boxtimes} & \bm{\boxtimes} & \bm{\boxtimes} & \bm{\boxtimes}& \bm{\boxtimes} \\
		\bm{\boxtimes} & \bm{\boxtimes} & \bm{\boxtimes} & \bm{\boxtimes}& \bm{\boxtimes} \\
	\bm{\boxtimes} & \bm{\boxtimes} & \bm{\boxtimes} & \bm{\boxtimes} & \bm{\boxtimes}\\
		\bm{\boxtimes} & \bm{\boxtimes} & \bm{\boxtimes} & \bm{\boxtimes} & \bm{\boxtimes}\\
		\bm{\boxtimes} & \bm{\boxtimes} & \bm{\boxtimes} & \bm{\boxtimes} & \bm{\boxtimes}
	\end{sbmatrix},
\end{aligned}
$$
where the left Householder ($\bH_1\times $) introduces zeros below the main diagonal  in the first column (see Section~\ref{section:qr-via-householder}).  Unfortunately, the right Householder ($\times \bH_1^\top$) will destroy the zeros introduced by the left Householder.

To address this issue, a more practical approach involves dividing the algorithm into two phases. 
In the first phase, we aim to transform the matrix into either a Hessenberg matrix (Definition~\ref{definition:upper-hessenbert}) or a tridiagonal matrix (Definition~\ref{definition:tridiagonal-hessenbert}). Subsequently, in the second phase, we seek another algorithm to further process the results obtained in the first phase and achieve the desired goal, thus completing the algorithm:
$$
\begin{aligned}
	\begin{sbmatrix}{\bA}
		\boxtimes & \boxtimes & \boxtimes & \boxtimes& \boxtimes \\
		\boxtimes & \boxtimes & \boxtimes & \boxtimes& \boxtimes \\
		\boxtimes & \boxtimes & \boxtimes & \boxtimes& \boxtimes \\
		\boxtimes & \boxtimes & \boxtimes & \boxtimes & \boxtimes\\
		\boxtimes & \boxtimes & \boxtimes & \boxtimes& \boxtimes
	\end{sbmatrix}
	&\stackrel{\bH_1\times }{\rightarrow}
	\begin{sbmatrix}{\bH_1\bA}
		\boxtimes & \boxtimes & \boxtimes & \boxtimes& \boxtimes \\
		\bm{\boxtimes} & \bm{\boxtimes} & \bm{\boxtimes} & \bm{\boxtimes}& \bm{\boxtimes} \\
		\bm{0} & \bm{\boxtimes} & \bm{\boxtimes} & \bm{\boxtimes} & \bm{\boxtimes}\\
		\bm{0} & \bm{\boxtimes} & \bm{\boxtimes} & \bm{\boxtimes} & \bm{\boxtimes}\\
		\bm{0} & \bm{\boxtimes} & \bm{\boxtimes} & \bm{\boxtimes}& \bm{\boxtimes}
	\end{sbmatrix}
	\stackrel{\times \bH_1^\top }{\rightarrow}
	\begin{sbmatrix}{\bH_1\bA\bH_1^\top}
		\boxtimes & \bm{\boxtimes} & \bm{\boxtimes} & \bm{\boxtimes}& \bm{\boxtimes} \\
		\boxtimes & \bm{\boxtimes} & \bm{\boxtimes} & \bm{\boxtimes}& \bm{\boxtimes} \\
		0 & \bm{\boxtimes} & \bm{\boxtimes} & \bm{\boxtimes} & \bm{\boxtimes}\\
		0 & \bm{\boxtimes} & \bm{\boxtimes} & \bm{\boxtimes} & \bm{\boxtimes}\\
		0 & \bm{\boxtimes} & \bm{\boxtimes} & \bm{\boxtimes} & \bm{\boxtimes}
	\end{sbmatrix},
\end{aligned}
$$
where the left Householder ($\bH_1\times $) will not influence the first row, and the right Householder ($\times \bH_1^\top$) will not influence the first column. 
A phase two algorithm, which is typically an iterative algorithm, for finding the triangular matrix is presented as follows:
$$
\begin{aligned}
	\begin{sbmatrix}{\bH_3\bH_2\bH_1\bA\bH_1^\top\bH_2^\top \bH_3^\top}
		\boxtimes & \boxtimes  & \boxtimes  & \boxtimes & \boxtimes  \\
		\boxtimes & \boxtimes  & \boxtimes  & \boxtimes & \boxtimes  \\
		0 & \boxtimes  & \boxtimes  & \boxtimes & \boxtimes  \\
		0 & 0 & \boxtimes  & \boxtimes & \boxtimes  \\
		0 & 0 & 0 & \boxtimes  & \boxtimes 
	\end{sbmatrix}
	\stackrel{\text{Phase 2} }{\longrightarrow}
	\begin{sbmatrix}{\bLambda}
	\bm{\boxtimes} & \bm{\boxtimes} & \bm{\boxtimes} & \bm{\boxtimes}& \bm{\boxtimes} \\
		\bm{0} & \bm{\boxtimes} & \bm{\boxtimes} & \bm{\boxtimes}& \bm{\boxtimes} \\
		0 & \bm{0}  & \bm{\boxtimes} & \bm{\boxtimes}& \bm{\boxtimes} \\
		0 & 0 & \bm{0}  & \bm{\boxtimes}& \bm{\boxtimes} \\
		0 & 0 & 0 & \bm{0}  & \bm{\boxtimes}
	\end{sbmatrix}
\end{aligned}
$$

Based on the preceding discussion, when computing the spectral decomposition, Schur decomposition, or singular value decomposition (SVD), it is customary to make a compromise. 
This typically involves calculating the Hessenberg, tridiagonal, or bidiagonal form in the first phase and reserving the second phase for completing the computation \citep{van2012families, van2014restructuring, trefethen1997numerical}.  
This part then surveys the \textit{Hessenberg decomposition}, \textit{tridiagonal decomposition}, and \textit{bidiagonal decomposition}, which often serve as the first phase of these more complex algorithms.

\newpage
\chapter{Hessenberg Decomposition and Similarity Transformations}\label{chapter:hessenberg}
\begingroup
\hypersetup{
	linkcolor=structurecolor,
	linktoc=page,  
}
\minitoc \newpage
\endgroup

\section{Hessenberg Decomposition}
\lettrine{\color{caligraphcolor}T}
The Hessenberg decomposition is a procedure employed to convert a matrix into an upper Hessenberg form.
This transformation simplifies the matrix, making it suitable for use as a first phase in various algorithms, thereby reducing the computational complexity.
Let us commence by providing a precise definition of upper Hessenberg matrices.

\begin{definition}[Upper Hessenberg Matrix\index{Hessenbert matrix}]\label{definition:upper-hessenbert}
An \textit{upper Hessenberg matrix} (simply called Hessenberg matrix when there is no confusion) is a square matrix in which all the entries below the   subdiagonal (i.e., the entries below the \textit{main diagonal}) are zeros.
Likewise, a lower Hessenberg matrix is a square matrix in which all the entries above the  superdiagonal (i.e., the entries above the main diagonal) are zeros.

The definition of an upper Hessenberg matrix can also be extended to rectangular matrices, and the specific form can be inferred from the context.

In matrix language, for any matrix $\bH\in \real^{n\times n}$ with each entry represented as $h_{ij}$ for all $i,j\in \{1,2,\ldots, n\}$, $\bH$ is considered an upper Hessenberg matrix if $h_{ij}=0$ for all $i\geq j+2$.

Let $i$ denote the smallest positive integer such that $h_{i+1, i}=0$ for $i\in \{1,2,\ldots, n-1\}$. By convention, the upper Hessenberg matrix $\bH$ is called \textbf{unreduced} if $i=n$.
\end{definition}

Taking a $5\times 5$ matrix as an example, the lower triangular below the  subdiagonal are zero in the upper Hessenberg matrix:
$$
\begin{sbmatrix}{possibly\,\, unreduced}
\boxtimes & \boxtimes & \boxtimes & \boxtimes & \boxtimes\\
\boxtimes & \boxtimes & \boxtimes & \boxtimes & \boxtimes\\
0 & \boxtimes & \boxtimes & \boxtimes & \boxtimes\\
0 & 0 & \boxtimes & \boxtimes & \boxtimes\\
0 & 0 & 0 & \boxtimes & \boxtimes
\end{sbmatrix}
\qquad 
\text{or}
\qquad 
\begin{sbmatrix}{reduced}
	\boxtimes & \boxtimes & \boxtimes & \boxtimes & \boxtimes\\
	\boxtimes & \boxtimes & \boxtimes & \boxtimes & \boxtimes\\
	0 & \boxtimes & \boxtimes & \boxtimes & \boxtimes\\
	0 & 0 & \textcolor{mylightbluetext}{0} & \boxtimes & \boxtimes\\
	0 & 0 & 0 & \boxtimes & \boxtimes
\end{sbmatrix}.
$$
We now present the Hessenberg decomposition:
\index{Decomposition: Hessenberg}
\begin{theoremHigh}[Hessenberg Decomposition]\label{theorem:hessenberg-decom}
Let $\bA$ be any $n\times n$ square matrix. Then it can be factored as 
$$
\bA = \bQ\bH\bQ^\top \qquad \text{or} \qquad \bH = \bQ^\top \bA\bQ, 
$$
where $\bH$ is an upper Hessenberg matrix, and $\bQ$ is an orthogonal matrix.
\end{theoremHigh}
It's not hard to find that a lower Hessenberg decomposition of $\bA^\top$ is given by $\bA^\top = \bQ\bH^\top\bQ^\top$ if $\bA$ admits the (upper) Hessenberg decomposition $\bA = \bQ\bH\bQ^\top$. 
The Hessenberg decomposition closely resembles the QR decomposition, as both transformations result in a matrix with zero entries in its lower sections, thereby achieving sparsity.
\begin{remark}[Why Hessenberg Decomposition]
We will observe that the introduction of zeros into $\bH$ from $\bA$ is accomplished by the left orthogonal matrix $\bQ$ (akin to the QR decomposition). 
In this context, the right orthogonal matrix $\bQ^\top$ does not result in any significant or simplified transformation of the matrix. 
So why do we opt for the Hessenberg decomposition instead of a QR decomposition, which exhibits a simpler structure with zeros even in the  subdiagonal? 

The answer lies in the preceding section, where we elucidated that the Hessenberg decomposition primarily serves as a phase one step for other algorithms.
It facilitates the factorization of the matrix into two orthogonal matrices, as seen in methods like SVD, UTV, and more. If we employ an aggressive algorithm that favors zeros in the  subdiagonal, such as in the QR decomposition, the right orthogonal transformation $\bQ^\top$ would rapidly eliminate the zeros, as we will soon observe.


On the other hand, expressing $\bA$ as $\bQ\bH\bQ^\top$ with respect to $\bH$ corresponds to an \textit{orthogonal similarity transformation} (Definition~\ref{definition:similar-matrices}) applied to $\bA$. 
This transformation ensures that  the eigenvalues, rank, and trace of both $\bA$ and $\bH$ are identical (Proposition~\ref{proposition:eigenvalue-similar-matrices}). Consequently, if we want to study the properties of $\bA$, exploration on $\bH$ can be a relatively simpler task than that on the original matrix $\bA$.

Furthermore, let $\bA=\bQ\bH\bQ^\top\in\real^{n\times n}$ be given. In some situations, we may want to solve the linear system $(\bA+\gamma\bI)\bx=\bb$ for different values of $\gamma\in\real$ and $\bb\in\real^n$.
The linear system can be equivalently expressed as $(\bH+\gamma\bI)\bQ^\top\bx=\bQ^\top\bb$. Since $\bH$ is upper Hessenberg, the system can be solved efficiently using methods like forward and backward substitution.
\end{remark}

\index{Similar matrices}
\index{Similarity transformation}
\index{Orthogonal similarity transformation}
\section{(Orthogonal) Similarity Transformation}
As previously mentioned, the Hessenberg decomposition presented in this section, the tridiagonal decomposition discussed in the next section, the Schur decomposition (Theorem~\ref{theorem:schur-decomposition}), and the spectral decomposition (Theorem~\ref{theorem:spectral_theorem}) all exhibit a common structure that converts the matrix into a \textit{similar matrix}. Herein, we provide a precise definition of similar matrices and similarity transformations.
\begin{definition}[Similar Matrices and Similarity Transformation]\label{definition:similar-matrices}
Let two matrices $\bA,\bB\in\real^{n\times n}$ be given.
Then  $\bA$ and $\bB$ are called \textit{similar matrices} if there exists a nonsingular matrix $\bP$ such that $\bB = \bP\bA\bP^{-1}$. 
	
In words, for any nonsingular matrix $\bP$, the matrices $\bA$ and $\bP\bA\bP^{-1}$ are similar. 
In this context, given the nonsingular matrix $\bP$, the operation $\bP\bA\bP^{-1}$
is referred to as a \textit{similarity transformation} applied to the matrix $\bA$.

Moreover, when $\bP$ is orthogonal, then $\bP\bA\bP^\top$ is also known as the \textit{orthogonal similarity transformation} of $\bA$. The orthogonal similarity transformation is significant  in the sense that the condition number of the transformed matrix $\bP\bA\bP^\top$ is not worse than that of the original matrix $\bA$. \footnote{Note that  $\bA$ and $\bB$ are called \textit{congruent} if $\bB = \bS\bA\bS^\top$ for some nonsingular matrix $\bS$. In this sense, an orthogonal similarity transformation is both a similarity transformation and a congruence transformation.}
\end{definition}

The difference between  similarity transformations and  orthogonal similarity transformations is partly explained in the context of coordinate transformations (Section~\ref{section:coordinate-transformation}).
Next, we prove some important properties of similar matrices in the following proposition: the eigenvalue, trace and rank of similar matrices.
\begin{proposition}[Eigenvalue, Trace, and Rank Invariant of Similar Matrices\index{Trace}\index{Similar matrices}]\label{proposition:eigenvalue-similar-matrices}
Let $\bA,\bP\in\real^{n\times n}$ be given, where $\bP$ is nonsingular. Then any eigenvalue of $\bA$ is also an eigenvalue of $\bP\bA\bP^{-1}$. 
The converse is also true that any eigenvalue of $\bP\bA\bP^{-1}$ is likewise an eigenvalue of $\bA$. 
In other words, $\Lambda(\bA) = \Lambda(\bP\bA\bP^{-1})$, where $\Lambda(\cdot)$ represents the spectrum of a matrix  (Definition~\ref{definition:spectrum}).

Furthermore, the trace, rank, and determinant of $\bA$ are equal to those of  $\bP\bA\bP^{-1}$ for any nonsingular matrix $\bP$ (both rank and trace are \textit{similarity invariant}): 
$$
\rank(\bA)=\rank(\bP\bA\bP^{-1}), \quad \trace(\bA)=\trace(\bP\bA\bP^{-1}),\quad \det(\bA)=\det(\bP\bA\bP^{-1}).
$$
\end{proposition}
\begin{proof}[of Proposition~\ref{proposition:eigenvalue-similar-matrices}]
For any eigenpair $(\lambda, \bx)$ of $\bA$, we have $\bA\bx =\lambda \bx$. Consequently, $\lambda \bP\bx = \bP\bA\bP^{-1} \bP\bx$, demonstrating that $\bP\bx$ is an eigenvector of $\bP\bA\bP^{-1}$ associated with $\lambda$.
Similarly, for any eigenpair $(\lambda, \by)$ of $\bP\bA\bP^{-1}$, we have $\bP\bA\bP^{-1} \by = \lambda \by$. Then, $\bA\bP^{-1} \by = \lambda \bP^{-1}\by$, indicating that $\bP^{-1}\by$ is an eigenvector of $\bA$ corresponding to $\lambda$. 

Regarding  the trace of $\bP\bA\bP^{-1}$, we can establish that $\trace(\bP\bA\bP^{-1}) = \trace(\bA\bP^{-1}\bP) = \trace(\bA)$, where the first equality comes from the fact that the
trace of a product is invariant under cyclical permutations of the factors:
$
\trace(\bA\bB\bC) = \trace(\bB\bC\bA) = \trace(\bC\bA\bB),
$
if all $\bA\bB\bC$, $\bB\bC\bA$, and $\bC\bA\bB$ exist.
The determinant follows that $\det(\bP\bA\bP^{-1})=\det(\bP)\det(\bA)\det(\bP^{-1})=\det(\bA)$.

Regarding the  rank of $\bP\bA\bP^{-1}$, we separate it into two claims as follows.
\paragraph{Rank claim 1: $\rank(\bZ\bA)=\rank(\bA)$ if $\bZ$ is nonsingular.}
We will begin by demonstrating that $\rank(\bZ\bA)=\rank(\bA)$ if $\bZ$ is nonsingular. 
Consider any vector $\bn$ in the null space of $\bA$, that is, $\bA\bn = \bzero$. 
Consequently, $\bZ\bA\bn = \bzero$, that is, $\bn$ also resides in the null space of $\bZ\bA$. 
This implies  $\nspace(\bA)\subseteq \nspace(\bZ\bA)$.

Conversely, for any vector $\bmm$ in the null space of $\bZ\bA$, i.e., $\bZ\bA\bmm = \bzero$, we can deduce that $\bA\bmm = \bZ^{-1} \bzero=\bzero$. That is, $\bmm$ is also in the null space of $\bA$. And this indicates $\nspace(\bZ\bA)\subseteq \nspace(\bA)$.

By combining these two arguments, we conclude:
$$
\nspace(\bA) = \nspace(\bZ\bA)\quad  \implies \quad \rank(\bZ\bA)=\rank(\bA).
$$

\paragraph{Rank claim 2: $\rank(\bA\bZ)=\rank(\bA)$ if $\bZ$ is nonsingular.}
We observe that the row rank of a matrix is equivalent to its column rank (Theorem~\ref{theorem:equal-dimension-rank}).
Therefore, $\rank(\bA\bZ) = \rank(\bZ^\top\bA^\top)$. Since $\bZ^\top$ is nonsingular, as per claim 1, we can conclude that $\rank(\bZ^\top\bA^\top) = \rank(\bA^\top) = \rank(\bA)$, where the last equality follows again from the fact that the row rank is equal to the column rank for any matrix. This establishes that $\rank(\bA\bZ)=\rank(\bA)$, as claimed.

Since $\bP$ and $\bP^{-1}$ are nonsingular, we can deduce that $\rank(\bP\bA\bP^{-1}) = \rank(\bA\bP^{-1}) = \rank(\bA)$. This completes the proof.
\end{proof}
The proposition  above will prove to be highly valuable in the subsequent discussion (refer to Lemma~\ref{lemma:rank-of-symmetric-idempotent} for its application in establishing the trace and rank of symmetric idempotent matrices).
Moreover, if we can find a set of similarity transformations on symmetric $\bA$ such that $\bP_T\ldots\bP_2\bP_1\bA\bP_1^{-1}\bP_2^{-1}\ldots\bP_T^{-1}$ is a diagonal matrix, then the eigenvalues of $\bA$ can be obtained. This is the main idea behind the methods discussed in Chapter~\ref{section:eigenvalue-problem} for calculating eigenvalues.

\begin{remark}
Having the same eigenvalues is a necessary condition for similarity transformations, but it is not a sufficient condition. 
For example, the matrices $\footnotesize\begin{bmatrix}
1&0\\
0&0
\end{bmatrix}$
and 
$\footnotesize\begin{bmatrix}
	0&0\\
	0&0
\end{bmatrix}$
both have the same eigenvalues $\{1,0\}$, but they are not similar.
\end{remark}

\begin{exercise}[Similarity Transformation]\label{exercise:sim_trans_inhess}
Let $\bA$ and $\bB$ be similar matrices. Show that 
\begin{itemize}
\item  $\det(\bA)=\det(\bB)$.
\item $\bA^k$ and $\bB^k$ are similar for any integer $k>1$.
\item $\bA$ is nonsingular if and only if $\bB$ is nonsingular.
\item  $(\bA-\lambda\bI)^k$ and $(\bB-\lambda\bI)^k$ are similar for any integer  $k>1$ and any $\lambda\in\real$; and their ranks are equal.
\item $\bB=\bzero$ if and only if $\bA=\bzero$.
\item $\bB=\bI$ if and only if $\bA=\bI$.
\end{itemize}
\end{exercise}

\section{Existence of  Hessenberg Decomposition}
We will demonstrate that any $n\times n$ matrix can be reduced into Hessenberg form through a sequence of Householder transformations applied alternately from the left and the right to the matrix  in an interleaved manner.  
In our earlier work, we employed Householder reflectors to triangularize matrices and introduce zeros below the \textbf{main diagonal}, resulting in the QR decomposition.  
A similar approach can be applied to introduce zeros below the \textbf{subdiagonal}.  
To comprehend this, a revisit to the fundamental concepts of Householder reflectors in Definition~\ref{definition:householder-reflector}  is necessary.

Before delving into the mathematical construction of such a decomposition, we would like to highlight the following remark, which will prove to be highly beneficial in the derivation of the decomposition.
\begin{remark}[Left and Right Multiplied by a Matrix with Block Identity]\label{remark:left-right-identity}
Let  $\bA\in \real^{n\times n}$ be a square matrix, and given the   matrix 
$
\bH = 
\footnotesize
\begin{bmatrix}
\bI_k &\bzero \\
\bzero & \bH_{n-k}
\end{bmatrix},
$
where $\bI_k$ is the $k\times k$ identity matrix. Then, $\bH\bA$ will not change the first $k$ rows of $\bA$, and $\bA\bH$ will not change the first $k$ columns of $\bA$.
\end{remark}
The proof of this remark is straightforward. 

\subsubsection*{\textbf{First Step: Introducing Zeros in the First Column}}	
Let $\bA=[\ba_1, \ba_2, \ldots, \ba_n]$ represent the column partition of matrix $\bA$, where each column $\ba_i$ belongs to $\real^{n}$. 
Suppose we obtain vectors $\bar{\ba}_1, \bar{\ba}_2, \ldots, \bar{\ba}_n \in \real^{n-1}$ by striking out the first component in each $\ba_i$. 
We define the following quantities:
$$
r_1 \triangleq \normtwo{\bar{\ba}_1}, \qquad \bu_1 \triangleq \frac{\bar{\ba}_1 - r_1 \be_1}{\normtwo{\bar{\ba}_1 - r_1 \be_1}}, \qquad \text{and}\qquad \widetilde{\bH}_1 \triangleq \bI - 2\bu_1\bu_1^\top \in \real^{(n-1)\times (n-1)},
$$
where $\be_1$ represents the first unit basis vector in $\real^{n-1}$, i.e., $\be_1=[1;0;0;\ldots;0]\in \real^{n-1}$. To introduce zeros below the subdiagonal and operate on the submatrix $\bA_{2:n,1:n}$, we append the Householder reflector into $\bH_1$, defined as 
$
\bH_1 \triangleq
\footnotesize
\begin{bmatrix}
1 &\bzero \\
\bzero & \widetilde{\bH}_1
\end{bmatrix},
$
in which case, $\bH_1\bA$ will introduce zeros in the first column of $\bA$ below entry (2,1). The first row of $\bA$ will remain unchanged, as noted in Remark~\ref{remark:left-right-identity}. And we can easily verify that both $\bH_1$ and $\widetilde{\bH}_1$ are orthogonal and symmetric matrices (from the definition of the Householder reflector, Remark~\ref{remark:householder-propes}). 
To achieve the form outlined in Theorem~\ref{theorem:hessenberg-decom}, we multiply $\bH_1\bA$ on the right by $\bH_1^\top$,  resulting in $\bH_1\bA\bH_1^\top$. Importantly, the multiplication by $\bH_1^\top$ on the right preserves the first column of $\bH_1\bA$, thereby retaining the zeros introduced in the first column.

An example of a $5\times 5$ matrix is shown as follows, where $\boxtimes$ represents a value that is not necessarily zero, and \textbf{boldface} indicates the value has just been changed:
$$
\begin{aligned}
	\begin{sbmatrix}{\bA}
		\boxtimes & \boxtimes & \boxtimes & \boxtimes & \boxtimes \\
		\boxtimes & \boxtimes & \boxtimes & \boxtimes & \boxtimes\\
		\boxtimes & \boxtimes & \boxtimes & \boxtimes & \boxtimes\\
		\boxtimes & \boxtimes & \boxtimes & \boxtimes & \boxtimes\\
		\boxtimes & \boxtimes & \boxtimes & \boxtimes & \boxtimes
	\end{sbmatrix}
	\stackrel{\bH_1\times}{\rightarrow}
	&\begin{sbmatrix}{\bH_1\bA}
		\boxtimes & \boxtimes & \boxtimes & \boxtimes & \boxtimes \\
		\bm{\boxtimes} & \bm{\boxtimes} & \bm{\boxtimes} & \bm{\boxtimes} & \bm{\boxtimes}\\
		\bm{0} & \bm{\boxtimes} & \bm{\boxtimes} & \bm{\boxtimes} & \bm{\boxtimes}\\
		\bm{0} & \bm{\boxtimes} & \bm{\boxtimes} & \bm{\boxtimes} & \bm{\boxtimes}\\
		\bm{0} & \bm{\boxtimes} & \bm{\boxtimes} & \bm{\boxtimes} & \bm{\boxtimes}
	\end{sbmatrix}
	\stackrel{\times\bH_1^\top}{\rightarrow}
	\begin{sbmatrix}{\bH_1\bA\bH_1^\top}
		\boxtimes & \bm{\boxtimes} & \bm{\boxtimes} & \bm{\boxtimes} & \bm{\boxtimes} \\
		\boxtimes & \bm{\boxtimes} & \bm{\boxtimes} & \bm{\boxtimes} & \bm{\boxtimes}\\
		0 & \bm{\boxtimes} & \bm{\boxtimes} & \bm{\boxtimes} & \bm{\boxtimes}\\
		0 & \bm{\boxtimes} & \bm{\boxtimes} & \bm{\boxtimes} & \bm{\boxtimes}\\
		0 & \bm{\boxtimes} & \bm{\boxtimes} & \bm{\boxtimes} & \bm{\boxtimes}
	\end{sbmatrix}.\\
\end{aligned}
$$

\subsubsection*{\textbf{Second Step: Introducing Zeros in the Second Column}}	
Let $\bB \triangleq \bH_1\bA\bH_1^\top$, where the entries in the first column below entry (2,1) are all zeros. The objective is to introduce zeros in the second column below entry (3,2). 
Let $\bB_2 \triangleq \bB_{2:n,2:n}=[\bb_1, \bb_2, \ldots, \bb_{n-1}]\in\real^{(n-1)\times(n-1)}$. 
Consider vectors $\bar{\bb}_1, \bar{\bb}_2, \ldots, \bar{\bb}_{n-1} \in \real^{n-2}$ obtained by striking out the first component in each $\bb_i$. We can  construct another Householder reflector as follows:
\begin{equation}\label{equation:householder-qr-lengthr}
r_1 \triangleq \normtwo{\bar{\bb}_1},  \qquad \bu_2 \triangleq \frac{\bar{\bb}_1 - r_1 \be_1}{\normtwo{\bar{\bb}_1 - r_1 \be_1}}, \qquad \text{and}\qquad \widetilde{\bH}_2 \triangleq \bI - 2\bu_2\bu_2^\top\in \real^{(n-2)\times (n-2)},
\end{equation}
where $\be_1$ now represents the first unit basis vector in $\real^{n-2}$. To introduce zeros below the subdiagonal and operate on the submatrix $\bB_{3:n,1:n}$, we append the Householder reflector into $\bH_2$, defined as 
$
\bH_2 \triangleq 
\footnotesize
\begin{bmatrix}
\bI_2 &\bzero \\
\bzero & \widetilde{\bH}_2
\end{bmatrix},
$
where $\bI_2$ represents the $2\times 2$ identity matrix. Notably, $\bH_2\bH_1\bA\bH_1^\top$ preserves the first two rows of $\bH_1\bA\bH_1^\top$; and since the Householder cannot reflect a zero vector, the zeros in the first column will be maintained. 
Furthermore, applying  $\bH_2^\top$ to the right of $\bH_2\bH_1\bA\bH_1^\top$ leaves the first two columns unchanged, ensuring that the zeros remain unaltered.

Following the example of a $5\times 5$ matrix, the second step is illustrated as follows, where $\boxtimes$ represents a value that is not necessarily zero, and \textbf{boldface} indicates the value has just been changed:
$$
\begin{aligned}
	\begin{sbmatrix}{\bH_1\bA\bH_1^\top}
		\boxtimes & \boxtimes & \boxtimes & \boxtimes &  \boxtimes  \\
		\boxtimes & \boxtimes & \boxtimes & \boxtimes &  \boxtimes  \\
		0 & \boxtimes & \boxtimes & \boxtimes &  \boxtimes  \\
		0 & \boxtimes & \boxtimes & \boxtimes &  \boxtimes  \\
		0 & \boxtimes & \boxtimes & \boxtimes &  \boxtimes 
	\end{sbmatrix}
	\stackrel{\bH_2\times}{\rightarrow}
	\begin{sbmatrix}{\bH_2\bH_1\bA\bH_1^\top}
		\boxtimes & \boxtimes & \boxtimes & \boxtimes & \boxtimes \\
		\boxtimes & \boxtimes & \boxtimes & \boxtimes & \boxtimes \\
		0 & \bm{\boxtimes} & \bm{\boxtimes} & \bm{\boxtimes} & \bm{\boxtimes}\\
		0 & \bm{0} & \bm{\boxtimes} & \bm{\boxtimes} & \bm{\boxtimes}\\
		0 & \bm{0} & \bm{\boxtimes} & \bm{\boxtimes} & \bm{\boxtimes}
	\end{sbmatrix}
	\stackrel{\times\bH_2^\top}{\rightarrow}
	\begin{sbmatrix}{\bH_2\bH_1\bA\bH_1^\top\bH_2^\top}
		\boxtimes & \boxtimes & \bm{\boxtimes} & \bm{\boxtimes} & \bm{\boxtimes} \\
		\boxtimes & \boxtimes & \bm{\boxtimes} & \bm{\boxtimes} & \bm{\boxtimes} \\
		0 & \boxtimes & \bm{\boxtimes} & \bm{\boxtimes} & \bm{\boxtimes}\\
		0 & 0 & \bm{\boxtimes} & \bm{\boxtimes} & \bm{\boxtimes}\\
		0 & 0 & \bm{\boxtimes} & \bm{\boxtimes} & \bm{\boxtimes}
	\end{sbmatrix}.
\end{aligned}
$$

The same process can go on, and there are a total of $n-2$ such steps. Ultimately, we will obtain the Hessenberg form by (given that $\bH_i$'s are both symmetric and orthogonal):
$$
\bH = \bH_{n-2} \bH_{n-3}\ldots\bH_1 \bA\bH_1^\top\bH_2^\top\ldots\bH_{n-2}^\top
= \bH_{n-2} \bH_{n-3}\ldots\bH_1 \bA\bH_1\bH_2\ldots\bH_{n-2}.
$$
Note that only $n-2$ such stages are involved, rather than $n-1$ or $n$. 
We will confirm this with an example below.
The example of a $5\times 5$ matrix in its entirety is shown as follows, where again $\boxtimes$ represents a value that is not necessarily zero, and \textbf{boldface} indicates the value has just been changed. 
\begin{mdframed}[hidealllines=\mdframehidelineNote,backgroundcolor=\mdframecolor, frametitle={A Complete Example of Hessenberg Decomposition}]
$$
\begin{aligned}
\begin{sbmatrix}{\bA}
	\boxtimes & \boxtimes & \boxtimes & \boxtimes & \boxtimes \\
	\boxtimes & \boxtimes & \boxtimes & \boxtimes & \boxtimes\\
	\boxtimes & \boxtimes & \boxtimes & \boxtimes & \boxtimes\\
	\boxtimes & \boxtimes & \boxtimes & \boxtimes & \boxtimes\\
	\boxtimes & \boxtimes & \boxtimes & \boxtimes & \boxtimes
\end{sbmatrix}
\stackrel{\bH_1\times}{\rightarrow}
&\begin{sbmatrix}{\bH_1\bA}
	\boxtimes & \boxtimes & \boxtimes & \boxtimes & \boxtimes \\
	\bm{\boxtimes} & \bm{\boxtimes} & \bm{\boxtimes} & \bm{\boxtimes} & \bm{\boxtimes}\\
	\bm{0} & \bm{\boxtimes} & \bm{\boxtimes} & \bm{\boxtimes} & \bm{\boxtimes}\\
	\bm{0} & \bm{\boxtimes} & \bm{\boxtimes} & \bm{\boxtimes} & \bm{\boxtimes}\\
	\bm{0} & \bm{\boxtimes} & \bm{\boxtimes} & \bm{\boxtimes} & \bm{\boxtimes}
\end{sbmatrix}
\stackrel{\times\bH_1^\top}{\rightarrow}
\begin{sbmatrix}{\bH_1\bA\bH_1^\top}
	\boxtimes & \bm{\boxtimes} & \bm{\boxtimes} & \bm{\boxtimes} & \bm{\boxtimes} \\
	\boxtimes & \bm{\boxtimes} & \bm{\boxtimes} & \bm{\boxtimes} & \bm{\boxtimes}\\
	0 & \bm{\boxtimes} & \bm{\boxtimes} & \bm{\boxtimes} & \bm{\boxtimes}\\
	0 & \bm{\boxtimes} & \bm{\boxtimes} & \bm{\boxtimes} & \bm{\boxtimes}\\
	0 & \bm{\boxtimes} & \bm{\boxtimes} & \bm{\boxtimes} & \bm{\boxtimes}
\end{sbmatrix}\\
\end{aligned}
$$

$$
\begin{aligned}
\gap\gap\gap\gap\gap\gap\gap
\stackrel{\bH_2\times}{\rightarrow}
&\begin{sbmatrix}{\bH_2\bH_1\bA\bH_1^\top}
	\boxtimes & \boxtimes & \boxtimes & \boxtimes & \boxtimes \\
	\boxtimes & \boxtimes & \boxtimes & \boxtimes & \boxtimes \\
	0 & \bm{\boxtimes} & \bm{\boxtimes} & \bm{\boxtimes} & \bm{\boxtimes}\\
	0 & \bm{0} & \bm{\boxtimes} & \bm{\boxtimes} & \bm{\boxtimes}\\
	0 & \bm{0} & \bm{\boxtimes} & \bm{\boxtimes} & \bm{\boxtimes}
\end{sbmatrix}
\stackrel{\times\bH_2^\top}{\rightarrow}
\begin{sbmatrix}{\bH_2\bH_1\bA\bH_1^\top\bH_2^\top}
	\boxtimes & \boxtimes & \bm{\boxtimes} & \bm{\boxtimes} & \bm{\boxtimes} \\
	\boxtimes & \boxtimes & \bm{\boxtimes} & \bm{\boxtimes} & \bm{\boxtimes} \\
	0 & \boxtimes & \bm{\boxtimes} & \bm{\boxtimes} & \bm{\boxtimes}\\
	0 & 0 & \bm{\boxtimes} & \bm{\boxtimes} & \bm{\boxtimes}\\
	0 & 0 & \bm{\boxtimes} & \bm{\boxtimes} & \bm{\boxtimes}
\end{sbmatrix}\\
\end{aligned}
$$

$$
\begin{aligned}
\gap\gap\gap\gap\gap\gap\gap\gap
\stackrel{\bH_3\times}{\rightarrow}
&\begin{sbmatrix}{\bH_3\bH_2\bH_1\bA\bH_1^\top\bH_2^\top}
	\boxtimes & \boxtimes & \boxtimes & \boxtimes & \boxtimes \\
	\boxtimes & \boxtimes & \boxtimes & \boxtimes & \boxtimes \\
	0 & \boxtimes & \boxtimes & \boxtimes & \boxtimes \\
	0 & 0 & \bm{\boxtimes} & \bm{\boxtimes} & \bm{\boxtimes}\\
	0 & 0 & \bm{0} & \bm{\boxtimes} & \bm{\boxtimes}
\end{sbmatrix}
\stackrel{\times\bH_3^\top}{\rightarrow}
\begin{sbmatrix}{\bH_3\bH_2\bH_1\bA\bH_1^\top\bH_2^\top\bH_3^\top}
	\boxtimes & \boxtimes & \boxtimes & \bm{\boxtimes} & \bm{\boxtimes} \\
	\boxtimes & \boxtimes & \boxtimes & \bm{\boxtimes} & \bm{\boxtimes} \\
	0 & \boxtimes & \boxtimes & \bm{\boxtimes} & \bm{\boxtimes}\\
	0 & 0 & \boxtimes & \bm{\boxtimes} & \bm{\boxtimes}\\
	0 & 0 & 0 & \bm{\boxtimes} & \bm{\boxtimes}
\end{sbmatrix},
\end{aligned}
$$
\end{mdframed}
where we find 
\begin{itemize}
\item when multiplying by $\bH_1$ on the left, we operate on a $(5-1)\times (5-1+1)$ submatrix; \item when multiplying by $\bH_2$ on the left, we operate on a $(5-2)\times (5-2+1)$ submatrix;  \item and when multiplying by $\bH_3$ on the left, we operate on a $(5-3)\times(5-3+1)$ submatrix. 
\end{itemize}
Similarly, 
\begin{itemize}
\item when multiplying by $\bH_1^\top$ on the right, we operate on a $5\times (5-1)$ submatrix; \item when multiplying by $\bH_2^\top$ on the right, we operate on a $5\times (5-2)$ submatrix; 
\item and when multiplying by $\bH_3^\top$ on the right, we operate on a $5\times (5-3)$ submatrix. 
\end{itemize}
These observations are crucial for efficiently computing the Hessenberg decomposition and reducing its computational complexity.

\index{Householder}
\section{Computing  Hessenberg Decomposition}\label{section:householder_hessenger}
For a given matrix $\bA\in \real^{n\times n}$,when we left-multiply it by $\bH_i$ for each step $i\in \{1, 2, \ldots, n-2\}$, we  perform operations on a submatrix of size $(n-i)\times(n-i+1)$. 
Furthermore, if we consider  the first column of this $(n-i)\times(n-i+1)$ matrix, we can set the values for the first column directly by the reflected vector used to obtain the Householder reflector (i.e., the first column is a scalar multiple of $\be_1\in\real^{n-i}$), allowing us to operate on an $(n-i)\times(n-i)$ submatrix instead.

When we right-multiply by $\bH_i^\top$ in each step $i\in \{1, 2, \ldots, n-2\}$, we perform operations on a submatrix of size $n\times(n-i)$.
After applying the Householder transformations, we obtain the final Hessenberg matrix $\bH$. The full process is outlined in Algorithm~\ref{alg:hessenbert-decomposition-householder}.

Similar to the computation of the QR decomposition using Householder reflectors (Section~\ref{section:householder_qr_cp}), in Algorithm~\ref{alg:hessenbert-decomposition-householder}, when calculating $\bH = \bH_{n-2} \bH_{n-3}\ldots\bH_1 \bA\bH_1^\top\bH_2^\top\ldots\bH_{n-2}^\top$, we express the equation as follows:
$$
\begin{aligned}
\bH &= \cleft[black](\bH_{n-2} \ldots\cleft[green](\bH_3\cleft[red](\bH_2\cleft[blue](\bH_1 \bA\bH_1\cright[blue])\bH_2\cright[red])\bH_3\cright[green])\ldots\bH_{n-2}\cright[black]) \\
&= 
\begin{bmatrix}
	\bI_{n-2}& \bzero \\
	\bzero& \bI - 2\bu_{n-2}\bu_{n-2}^\top 
\end{bmatrix}
\ldots
\begin{bmatrix}
	\bI_3& \bzero \\
	\bzero& \bI - 2\bu_3\bu_3^\top 
\end{bmatrix}
\begin{bmatrix}
	\bI_2& \bzero \\
	\bzero& \bI - 2\bu_2\bu_2^\top 
\end{bmatrix}
\begin{bmatrix}
	\bI_1 &\bzero \\
	\bzero & \bI - 2\bu_1\bu_1^\top
\end{bmatrix}\\
&\bA
\begin{bmatrix}
		\bI_1 &\bzero \\
		\bzero & \bI - 2\bu_1\bu_1^\top
	\end{bmatrix}
	\begin{bmatrix}
		\bI_2& \bzero \\
		\bzero& \bI - 2\bu_2\bu_2^\top 
	\end{bmatrix}
	\begin{bmatrix}
		\bI_3& \bzero \\
		\bzero& \bI - 2\bu_3\bu_3^\top 
	\end{bmatrix}
	\ldots
	\begin{bmatrix}
		\bI_{n-2}& \bzero \\
		\bzero& \bI - 2\bu_{n-2}\bu_{n-2}^\top 
	\end{bmatrix},
\end{aligned}
$$
where the different colors of the parentheses indicate the computational order. Specifically:
\begin{itemize}
\item The upper-left of $\bH_1$ is a $1\times 1$ identity matrix. Multiplying on the left will not change the \textbf{first row} of $\bA$, and multiplying on the right will not change the \textbf{first column} of $\bH_1\bA$.
\item The upper-left of $\bH_2$ is a $2\times 2$ identity matrix. Multiplying on the left will not change the \textbf{first 2 rows} of $\bH_1\bA\bH_1$, and multiplying on the right will not change the \textbf{first 2 columns} of $\bH_2\bH_1\bA\bH_1$.
\item The upper-left of $\bH_3$ is a $3\times 3$ identity matrix. Multiplying on the left will not change the \textbf{first 3 rows} of $\bH_2\bH_1\bA\bH_1\bH_2$, and multiplying on the right will not change the \textbf{first 3 columns} of $\bH_3\bH_2\bH_1\bA\bH_1\bH_2$.
\item The process can go on, and this property yields the steps 8 and 9 in the algorithm.
\end{itemize}

Similarly, to obtain the final orthogonal matrix $\bQ=\bH_1 \bH_2\ldots\bH_{n-2}$, we write out the equation:
$$
\begin{aligned}
	\bQ
	&=\bH_1 \bH_2\bH_3\ldots\bH_{n-2} \\
	&= \begin{bmatrix}
		\bI_1 &\bzero \\
		\bzero & \bI - 2\bu_1\bu_1^\top
	\end{bmatrix}
	\begin{bmatrix}
		\bI_2& \bzero \\
		\bzero& \bI - 2\bu_2\bu_2^\top 
	\end{bmatrix}
	\begin{bmatrix}
		\bI_3& \bzero \\
		\bzero& \bI - 2\bu_3\bu_3^\top 
	\end{bmatrix}
	\ldots
	\begin{bmatrix}
		\bI_{n-2}& \bzero \\
		\bzero& \bI - 2\bu_{n-2}\bu_{n-2}^\top 
	\end{bmatrix},
\end{aligned}
$$
where the upper-left portion of $\bH_1$ is a $1\times 1$ identity matrix; the upper-left portion of $\bH_2$ is a $2\times 2$ identity matrix, and it will not change the first 2 columns of $\bH_1$; and the upper-left portion of $\bH_3$ is a $3\times 3$ identity matrix, which will not modify the first 3 columns of $\bH_1\bH_2$; $\ldots$. This property yields the step 15 in the algorithm.

\begin{algorithm}[h] 
\caption{Hessenberg Decomposition via  Householder Reflectors} 
\label{alg:hessenbert-decomposition-householder} 
\begin{algorithmic}[1] 
\Require Matrix $\bA$ with size $n\times n $; 
\Statex \textbf{Stage A: Obtain the Hessenberg matrix}
\State Initially set $\bH \leftarrow \bA$; (\textcolor{mylightbluetext}{note $\bH$ is the Hessenberg, $\bH_i$'s are Householders})
\For{$i=1$ to $n-2$} 
\State $\ba \leftarrow \bH_{i+1:n,i}$, i.e., first column of $\bH_{i+1:n,i:n}\in \real^{(n-i)\times(n-i+1)}$;
\State $r \leftarrow \normtwo{\ba}$;\Comment{$2(n-i)$ flops}
\State $\bu_i \leftarrow \ba-r\be_1$; \Comment{1 flop}
\State $\bu_i \leftarrow \bu_i / \normtwo{\bu_i}$ (Reflector is $ \widetilde{\bH}_i = \bI - 2\bu_i\bu_i^\top \in \real^{(n-i)\times (n-i)}$); \Comment{$3(n-i)$ flops}
\State $\bH_{i+1,i}\leftarrow r$, $\bH_{i+2:n,i}\leftarrow\bzero$, i.e., obtain the first column of $\bH_{i+1:n,i:n}$; \Comment{0 flops}
\State \algoalign{Left: set the value of the second   to the last column of $\bH_{i+1:n,i:n}$, i.e., operate on $\bH_{i+1:n,i\textcolor{mylightbluetext}{+1}:n}\in\real^{(n-i)\times (n-i)}$:}
$$
\begin{aligned}
&\bH_{i+1:n,i+1:n} \leftarrow (\bI-2\bu_i\bu_i^\top)\bH_{i+1:n,i+1:n} \\
&= \bH_{i+1:n,i+1:n} - 2\bu_i\bu_i^\top\bH_{i+1:n,i+1:n} \in \real^{(n-i)\times(n-i)},\quad
\text{($4(n-i)^2$ flops);}
\end{aligned}
$$
\State Right: operate on $\bH_{1:n,i+1:n}\in\real^{n\times (n-i)}$:
$$
\begin{aligned}
&\bH_{1:n,i+1:n} \leftarrow \bH_{1:n,i+1:n}(\bI-2\bu_i\bu_i^\top) \\
&=\bH_{1:n,i+1:n} - \bH_{1:n,i+1:n} 2\bu_i\bu_i^\top\in \real^{n\times(n-i)},
\quad
\text{($4n(n-i)-n$ flops);}
\end{aligned}
$$
\EndFor
\State Output $\bH$ as the upper Hessenberg matrix;
\Statex \textbf{Stage B: Obtain the orthogonal matrix}
\State Get $\bQ\leftarrow\bH_1\bH_2\ldots\bH_{n-2}$, where $\bH_i$'s are Householder reflectors: $\bH_i=\diag(\bI_i, \widetildebH_i)$.
\State Initially set $\bQ \leftarrow \bH_1$;
\For{$i=1$ to $n-3$} 
\State Compute $\bQ$:
$$\begin{aligned}
&\bQ_{1:n,i+2:n} \leftarrow \bQ_{1:n,i+2:n}\widetildebH_{i+1}
= \bQ_{1:n,i+2:n}(\bI - 2\bu_{i+1}\bu_{i+1}^\top)\\
&=\bQ_{1:n,i+2:n}-\bQ_{1:n,i+2:n}2\bu_{i+1}\bu_{i+1}^\top \in \real^{n\times (n-i-1)},
\quad
\text{($4n(n-i-1) - n$ flops);}
\end{aligned}$$
\EndFor
\State Output $\bQ$ as the orthogonal matrix.
\end{algorithmic} 
\end{algorithm}

\begin{theorem}[Algorithm Complexity: Hessenberg via Householder]\label{theorem:hessenbert-householder}
Algorithm~\ref{alg:hessenbert-decomposition-householder} requires $\sim \frac{10}{3}n^3$ flops to compute the Hessenberg decomposition of an $n\times n$ square matrix. Furthermore, if $\bQ$ is needed explicitly, an additional $\sim 2n^3$ flops are required.
\end{theorem}

\begin{proof}[of Theorem~\ref{theorem:hessenbert-householder}] 
We divide the proof into two parts: obtaining the Hessenberg matrix and the orthogonal matrix.
\paragraph{Obtaining the Hessenberg matrix.}
For each loop iteration $i$, the submatrix $\bH_{i+1:n,i:n}\in \real^{(n-i)\times(n-i+1)}$. Therefore, $\ba$ in step 3 is in $\real^{n-i}$.

In step 4,  computing $r = \normtwo{\ba}$ involves $n-i$ multiplications, $n-i-1$ additions, and 1 square root operation,  totaling \underline{$2(n-i)$} flops.

In step 5, calculating $\bu_i = \ba-r\be_1$ requires 1 subtraction, which is \underline{$1$} flop due to the special structure of $\be_1$;

In step 6, similar to step 4, it requires $2(n-i)$ flops to compute the norm and $(n-i)$ additional divisions, resulting in \underline{$3(n-i)$} flops.


In step 8, suppose in loop $i$, computing $\bu_i^\top \bH_{i+1:n,i+1:n}$ requires $(n-i)\times$($n-i$ multiplications and $n-i-1$ additions), which is \underline{$(n-i)(2(n-i)-1)$} flops. Calculating $2\bu_i$ requires \underline{$n-i$} multiplications. 
Furthermore, the product $2\bu_i (\bu_i^\top \bH_{i+1:n,i+1:n})$ entails \underline{$(n-i)^2$} multiplications to form an $(n-i)\times (n-i)$ matrix. The final matrix subtraction involves \underline{$(n-i)^2$} subtractions. Thus, the overall complexity for step 8  is \underline{$4(n-i)^2$} flops for each iteration $i$. 

In step 9, the computation of $2\bu_i$ needs 0 flops since it has already been calculated in step 8. Similarly, the product $\bH_{1:n,i+1:n} 2\bu_i$ involves $n\times$($n-i$ multiplications and $n-i-1$ additions), which is \underline{$n(2(n-i)-1)$} flops. $\bH_{1:n,i+1:n} 2\bu_i\bu_i^\top$ takes \underline{$n(n-i)$} multiplications to generate an $n\times (n-i)$ matrix. The final matrix subtraction requires additional \underline{$n(n-i)$} subtractions. Consequently, the computational complexity of step 9 amounts to \underline{$4n(n-i)-n$} flops for each iteration $i$.

For each iteration $i$ in the loop, the total computational complexity is given by $4i^2 - (12n+4)i + (8n^2+3n+2)$ flops. We can define a function $f(i)\triangleq 4i^2 - (12n+4)i + (8n^2+3n+2)$  to represent this complexity, the complexity for step 2 to step 10 can be obtained by 
$
\mathrm{cost}=f(1)+f(2)+\ldots+f(n-2).
$
The sum of these $n-2$ loops converges to $\frac{10}{3}n^3$ flops if we keep only the leading term.

\paragraph{Obtaining the orthogonal matrix.}
For the additional computation of $\bQ$ in step 15, the situation is similar to step 9. The computation of $2\bu_{i+1}$ needs 0 flops since it has already been calculated in step 8. The product of $\bQ_{1:n,i+2:n}2\bu_{i+1}$ involves $n\times$($n-i-1$ multiplications and $n-i-2$ additions), which is \underline{$n(2(n-i-1)-1)$} flops. The calculation of $\bQ_{1:n,i+2:n}2\bu_{i+1}\bu_{i+1}^\top$ takes \underline{$n(n-i-1)$} flops to make an $n\times (n-i-1)$ matrix. The final matrix subtraction needs additional \underline{$n(n-i-1)$} subtractions. 
Consequently, the computational complexity of step 15 amounts to \underline{$4n(n-i-1) - n$} flops for each iteration $i$. Let $g(i)\triangleq 4n(n-i-1) - n$, the overall complexity for step 14 to step 16 can be expressed as 
$
\mathrm{cost}=g(1)+g(2)+\ldots +g(n-3),
$
which shows  that the  complexity is $2n^3$ flops if we keep only the leading term.
\end{proof}

\section{Properties of  Hessenberg Decomposition}\label{section:hessenberg-decomposition}

\begin{theorem}[Shift Hessenberg]
Let $\bH\in\real^{n\times n}$ be an unreduced Hessenberg matrix and $\lambda$ be its eigenvalue. Let further $\bH-\lambda\bI=\bQ\bR$ be the QR decomposition of the shifted matrix $\bH-\lambda\bI$. Let $\bA\triangleq \bR\bQ+\lambda\bI$, then $a_{n,n-1=0}$ and $a_{nn}=\lambda$.
\end{theorem}
The proof is straightforward since $\bH$ is unreduced such that the first $n-1$ columns are linearly independent and $r_{ii}\neq 0$ for $i\in\{1,2,\ldots,n-1\}$.
However, $\bH-\lambda\bI$ is singular, implying $r_{nn}=0$ and $a_{nn}=\lambda$.

The Hessenberg decomposition is not unique due to the different ways to construct the Householder reflectors (see Equation~\eqref{equation:householder-qr-lengthr}). However, under certain mild conditions, we can observe a similar structure in different decompositions.

\begin{theorem}[Implicit Q Theorem for Hessenberg Decomposition\index{Implicit Q theorem}]\label{theorem:implicit-q-hessenberg}
Let $\bA\in \real^{n\times n}$ be given and suppose two Hessenberg decompositions of $\bA$  are given by $\bA=\bU\bH\bU^\top=\bV\bG\bV^\top$, where $\bU=[\bu_1, \bu_2, \ldots, \bu_n]$ and $\bV=[\bv_1, \bv_2, \ldots, \bv_n]$ are the column partitions of $\bU$ and $\bV$, respectively. Suppose further that $k$ is the smallest positive integer for which $h_{k+1,k}=0$, where $h_{ij}$ denotes the entry $(i,j)$ of $\bH$. Then,
\begin{itemize}
\item If $\bu_1=\bv_1$, then $\bu_i = \pm \bv_i$ and $|h_{i,i-1}| = |g_{i,i-1}|$ for $i\in \{2,3,\ldots,k\}$. 
\item When $k=n$, the Hessenberg matrix $\bH$ is known as \textit{unreduced} (Definition~\ref{definition:upper-hessenbert}). However, if $k<n$, then $g_{k+1,k}=0$.
\end{itemize}
\end{theorem}
\begin{proof}[of Theorem~\ref{theorem:implicit-q-hessenberg}]
Define the orthogonal matrix $\bZ\triangleq\bV^\top\bU$. We have:
$$
\left. 
\begin{aligned}
\bG\bZ &= \bV^\top\bA\bV \bV^\top\bU = \bV^\top\bA\bU \\
\bZ\bH &= \bV^\top\bU \bU^\top\bA\bU = \bV^\top\bA\bU 
\end{aligned}
\right\}
\quad\implies\quad  
\bG\bZ = \bZ\bH.
$$
The column indexed by $(i-1)$ in each of them can be represented as:
$
\bG\bz_{i-1} = \bZ\bh_{i-1},
$
where $\bz_{i-1}$ and $\bh_{i-1}$ are the $(i-1)$-th columns of $\bZ$ and $\bH$, respectively. Since $h_{l,i-1}=0$ for $l\geq i+1$ (as per the definition of upper Hessenberg matrices), $\bZ\bh_{i-1}$ can be represented as 
$$
\bZ\bh_{i-1} = \sum_{j=1}^{i} h_{j,i-1} \bz_j = h_{i,i-1}\bz_i + \sum_{j=1}^{i-1} h_{j,i-1} \bz_j.
$$
Combining the two findings above, it follows that 
$$
h_{i,i-1}\bz_i  = \bG\bz_{i-1} - \sum_{j=1}^{i-1} h_{j,i-1} \bz_j.
$$
A moment of reflexion reveals that $[\bz_1, \bz_2, \ldots,\bz_k]$ is upper triangular. And since $\bZ$ is orthogonal, it must be diagonal, and each value on the diagonal is in $\{-1, 1\}$ for $i\in \{2,\ldots, k\}$. Then $\bz_1=\be_1$ and $\bz_i = \pm \be_i$ for $i\in \{2,\ldots, k\}$. Additionally,  $\bz_i =\bV^\top\bu_i$ and $h_{i,i-1}=\bz_i^\top (\bG\bz_{i-1} - \sum_{j=1}^{i-1} h_{j,i-1} \bz_j)=\bz_i^\top \bG\bz_{i-1}$. Therefore, for $i\in \{2,\ldots,k\}$, $\bz_i^\top \bG\bz_{i-1}$ is just $\pm g_{i,i-1}$. It follows that 
$
\begin{aligned}
|h_{i,i-1}| &= |g_{i,i-1}| \,\, \text{and}\,\,   \bu_i = \pm\bv_i,\,\, \forall i\in \{2,3,\ldots,k\}.
\end{aligned}
$ 
This proves the first part. For the second part, if $k<n$, 
$$
\begin{aligned}
g_{k+1,k} &= \be_{k+1}^\top\bG\be_{k} = \pm \be_{k+1}^\top\underbrace{\bG\bZ}_{\bZ\bH} \be_{k} =  \pm \be_{k+1}^\top\underbrace{\bZ\bH \be_{k}}_{\text{$k$-th column of $\bZ\bH$}} \\
&= \pm \be_{k+1}^\top \bZ\bh_k = \pm \be_{k+1}^\top \sum_{j=1}^{k+1} h_{jk}\bz_j 
=\pm \be_{k+1}^\top \sum_{j=1}^{\textcolor{mylightbluetext}{k}} h_{jk}\bz_j=0,
\end{aligned}
$$
where the penultimate equality is derived from the assumption that $h_{k+1,k}=0$. This completes the proof.
\end{proof}
\paragraph{Key observations.}
From the theorem mentioned above, when two Hessenberg decompositions of matrix $\bA$ are both unreduced and have the same first column in the orthogonal matrices, then the Hessenberg matrices $\bH, \bG$ are similar matrices such that $\bH = \bD\bG\bD^{-1}$, where $\bD=\diag(\pm 1, \pm 1, \ldots, \pm 1)$. \textit{Moreover, and most importantly, if we impose the constraint that  the elements in the  subdiagonal of the Hessenberg matrix $\bH$ are positive (if possible), then the Hessenberg decomposition $\bA=\bQ\bH\bQ^\top$ is uniquely determined by $\bA$ and the first column of $\bQ$.} This is similar to what we have claimed about the uniqueness of the QR decomposition (Corollary~\ref{corollary:unique-qr}).

The next finding involves a Krylov matrix defined as follows:
\begin{definition}[Krylov Matrix\index{Krylov matrix}]\label{definition:krylov-matrix}
Let $\bA\in \real^{n\times n}$, $\bq\in \real^n$, and $k$ be given, where $k$ is a positive integer. Then the \textit{Krylov matrix} is defined as follows:
$$
\bK(\bA, \bq, k) = 
\begin{bmatrix}
	\bq, & \bA\bq, & \ldots, & \bA^{k-1}\bq 
\end{bmatrix}
\in \real^{n\times k}.
$$
\end{definition}

\begin{theorem}[Unreduced Hessenberg]\label{theorem:implicit-q-hessenberg-v2}
Let $\bA\in \real^{n\times n}$ be given, and suppose there exists an orthogonal matrix $\bQ$ such that $\bA$ can be factored as $\bA = \bQ\bH\bQ^\top$. Then, $\bQ^\top\bA\bQ=\bH$ is an unreduced upper Hessenberg matrix if and only if $\bR=\bQ^\top \bK(\bA, \bq_1, n)$ is both nonsingular and upper triangular, where $\bq_1$ is the first column of $\bQ$.  

Moreover, if $\bR$ is singular and $k$ is the smallest index such that $r_{kk}=0$, then $k$ is also the smallest index for which $h_{k,k-1}=0$.
\end{theorem}
\begin{proof}[of Theorem~\ref{theorem:implicit-q-hessenberg-v2}] 
First, assume $\bH$ is unreduced, we can explicitly express the following matrix:
$$
\bR = \bQ^\top \bK(\bA, \bq_1, n) = [\be_1, \bH\be_1, \ldots, \bH^{n-1}\be_1],
$$
where, obviously,  $\bR$ is upper triangular with $r_{11}=1$. Observe that $r_{ii} = h_{21}h_{32}\ldots h_{i,i-1}$ for $i\in \{2,3,\ldots, n\}$. Since $\bH$ is unreduced, $\bR$ is nonsingular as well.

Conversely, assume $\bR$ is upper triangular and nonsingular, we observe that $\br_{k+1} = \bH\br_{k}$, implying that the $(k+2:n)$-th rows of $\bH[:,1:k]$ are zero and $h_{k+1,k}\neq 0$ for $k\in \{1,2,\ldots, n-1\}$. Therefore, $\bH$ is unreduced.

If $\bR$ is singular and $k$ is the smallest index satisfying $r_{kk}=0$, then 
$$
\left. 
\begin{aligned}
r_{k-1,k-1}&=h_{21}h_{32}\ldots h_{k-1,k-2}&\neq 0; \\
r_{kk}&=h_{21}h_{32}\ldots h_{k-1,k-2} h_{k,k-1}&= 0,
\end{aligned}
\right\}
\quad\implies\quad  
h_{k,k-1} =0,
$$
from which the result follows.
\end{proof}

\begin{algorithm}
\caption{Moler and Stewart's HT reduction \citep{moler1973algorithm}}
\label{alg:hess_tri_form} 
\begin{algorithmic}[1] 
\Require A general matrix $\bA \in \real^{n \times n}$ and an upper triangular matrix $\bB \in \real^{n \times n}$;
\Ensure Orthogonal matrices $\bQ, \bZ \in \real^{n \times n}$ such that $(\bH,\bT) = (\bQ^\top \bA \bZ, \bQ^\top \bB \bZ)$ is in HT form;
\State \textbf{Remark:} $\bL_{i-1,i}, \bR_{i,i-1} \in \real^{n \times n}$ denote  Givens rotations (Section~\ref{section:qr-givens}) acting on rows/columns $i-1$ and $i$. 
\State Initially set $\bQ \leftarrow \bI_n$, $\bZ \leftarrow \bI_n$, $\bH\leftarrow \bA$, and $\bT\leftarrow \bB$;
\For{$j = 1, 2, \ldots, n-2$} \Comment{Introduce zeros in the $j$-th column of $\bA$}
\For{$i = n, n-1, \ldots, j+2$}
\State Construct $\bL_{i-1,i}$ such that the $(i,j)$-th entry of $\bL_{i-1,i}^\top \bH$ is zero.
\State Update $\bH \leftarrow \bL_{i-1,i}^\top \bH$, $\bT \leftarrow \bL_{i-1,i}^\top \bT$, $\bQ \leftarrow \bQ \bL_{i-1,i}$.
\State Construct $\bR_{i,i-1}$ such that the fill-in $(i,i-1)$ entry of $\bT \bR_{i,i-1}$ is zero.
\State Update $\bH \leftarrow \bH \bR_{i,i-1}$, $\bT \leftarrow \bT \bR_{i,i-1}$, $\bZ \leftarrow \bZ \bR_{i,i-1}$.
\EndFor
\EndFor
\State Output $(\bH,\bT) = (\bQ^\top \bA \bZ, \bQ^\top \bB \bZ)$.
\end{algorithmic}
\end{algorithm}

\section{Hessenberg-Triangular Decomposition}\label{section:ht_form}
Given a matrix pair $(\bA, \bB)$, where $\bA, \bB \in \real^{n \times n}$, a preprocessing step of the \textit{QZ decomposition or generalized Schur decomposition} \citep{moler1973algorithm} for solving the regular \textit{generalized eigenvalue problem}~\footnote{See Problem~\ref{prob:geneig1}$\sim$\ref{prob:geneig2} for more details.} $(\bA - \lambda \bB)\bx = \bzero$ involves  computing orthogonal matrices $\bQ, \bZ \in \real^{n \times n}$ such that $\bQ^\top \bA \bZ$ is upper Hessenberg while $\bQ^\top \bB \bZ$ is upper triangular. This so-called \textit{Hessenberg-triangular (HT) form} of the matrix pair $(\bA, \bB)$  significantly reduces the computational cost during the iterative part of the QZ algorithm, which in turn plays a crucial role in the computation of quadratic eigenvalue problems \citep{zhang2017matrix}.

The HT reduction begins by computing a QR decomposition $\bB = \bQ_0 \bB_0$, where $\bQ_0$ is orthogonal and $\bB_0$ is upper triangular. The matrices $\bA$ and $\bB$ are then overwritten by $\bQ_0 \bA$ and $\bQ_0 \bB = \bB_0$, respectively. 
Thus, for the rest of this section, we assume that the matrix $\bB$ in the pair $(\bA, \bB)$ is already in upper triangular form. 
In the HZ algorithm, matrix $\bA$ is then reduced to Hessenberg form by applying a sequence of Givens rotations; see Algorithm~\ref{alg:hess_tri_form}.
The goal is to reduce $\bA$ to Hessenberg form while maintaining the triangular form of $\bB$. 
This is achieved by premultiplying $\bA$ with Householder reflections or Givens rotations to annihilate elements below the first subdiagonal, and postmultiplying $\bB$ with a different set of Householder reflections or Givens rotations to preserve its triangular form (we use Givens rotations in Algorithm~\ref{alg:hess_tri_form}).

An example of a $7\times 7$ matrix is shown as follows at $i=5$ and $j=2$, where $\boxtimes$ represents a value that is not necessarily zero, and \textbf{boldface} indicates the value has just been changed. The \textcolor{mylightbluetext}{blue} elements are introduced to zero from a nonzero value; while the \textcolor{winestain}{red} elements are modified to nonzero from a zero value:
$$
\setlength{\arraycolsep}{2pt}
\begin{sbmatrix}{\bL_{4,5}^\top\bH}
\boxtimes & \boxtimes & \boxtimes & \boxtimes & \boxtimes & \boxtimes & \boxtimes \\
\boxtimes & \boxtimes & \boxtimes & \boxtimes & \boxtimes & \boxtimes & \boxtimes \\
0 & \boxtimes & \boxtimes & \boxtimes & \boxtimes & \boxtimes & \boxtimes \\
0 & \bm{\boxtimes} & \bm{\boxtimes} & \bm{\boxtimes} & \bm{\boxtimes} & \bm{\boxtimes} & \bm{\boxtimes} \\
0 & \textcolor{mylightbluetext}{\bm{0}} & \bm{\boxtimes} & \bm{\boxtimes} & \bm{\boxtimes} & \bm{\boxtimes} & \bm{\boxtimes} \\
0 & 0 & \boxtimes & \boxtimes & \boxtimes & \boxtimes & \boxtimes \\
0 & 0 & \boxtimes & \boxtimes & \boxtimes & \boxtimes & \boxtimes \\
\end{sbmatrix}
\begin{sbmatrix}{\bL_{4,5}^\top\bT}
\boxtimes & \boxtimes & \boxtimes & \boxtimes & \boxtimes & \boxtimes & \boxtimes \\
0 & \boxtimes & \boxtimes & \boxtimes & \boxtimes & \boxtimes & \boxtimes \\
0 & 0 & \boxtimes & \boxtimes & \boxtimes & \boxtimes & \boxtimes \\
0 & 0 & 0 & \bm{\boxtimes} & \bm{\boxtimes} & \bm{\boxtimes} & \bm{\boxtimes} \\
0 & 0 & 0 & \textcolor{winestain}{\bm{\boxtimes}} & \bm{\boxtimes} & \bm{\boxtimes} & \bm{\boxtimes} \\
0 & 0 & 0 & 0 & 0 & \boxtimes & \boxtimes \\
0 & 0 & 0 & 0 & 0 & 0 & \boxtimes \\
\end{sbmatrix}
\rightarrow
\begin{sbmatrix}{\bL_{4,5}^\top\bH\bR_{5,4}}
\boxtimes & \boxtimes & \boxtimes & \bm{\boxtimes} & \bm{\boxtimes} & {\boxtimes} & \boxtimes \\
\boxtimes & \boxtimes & \boxtimes & \bm{\boxtimes} & \bm{\boxtimes} & {\boxtimes} & \boxtimes \\
\boxtimes & \boxtimes & \boxtimes & \bm{\boxtimes} & \bm{\boxtimes} & {\boxtimes} & \boxtimes \\
0 & \boxtimes & \boxtimes & \bm{\boxtimes} & \bm{\boxtimes} & {\boxtimes} & \boxtimes \\
0 & 0 & \boxtimes & \bm{\boxtimes} & \bm{\boxtimes} & {\boxtimes} & \boxtimes \\
0 & 0 & \boxtimes & \bm{\boxtimes} & \bm{\boxtimes} & {\boxtimes} & \boxtimes \\
0 & 0 & \boxtimes & \bm{\boxtimes} & \bm{\boxtimes} & {\boxtimes} & \boxtimes \\
\end{sbmatrix}
\begin{sbmatrix}{\bL_{4,5}^\top\bT\bR_{5,4}}
\boxtimes & \boxtimes & \boxtimes & \bm{\boxtimes} & \bm{\boxtimes} & \boxtimes & \boxtimes \\
0 & \boxtimes & \boxtimes & \bm{\boxtimes} & \bm{\boxtimes} & \boxtimes & \boxtimes \\
0 & 0 & \boxtimes & \bm{\boxtimes} & \bm{\boxtimes} & \boxtimes & \boxtimes \\
0 & 0 & 0 & \bm{\boxtimes} & \bm{\boxtimes} & \boxtimes & \boxtimes \\
0 & 0 & 0 & \textcolor{mylightbluetext}{\bm{0}} & \bm{\boxtimes} & \boxtimes & \boxtimes \\
0 & 0 & 0 & 0 & 0 & \boxtimes & \boxtimes \\
0 & 0 & 0 & 0 & 0 & 0 & \boxtimes \\
\end{sbmatrix}.
$$

\begin{problemset}
\item Let $\bA,\bB\in\real^{n\times n}$ be similar. Show that $\adjugate(\bA)$ and $\adjugate(\bB)$ are similar.
\item Let $\bA$ be given and $\bP$ be nonsingular. Show that if $\bP\bA\bP^{-1}$ is upper triangular, then the diagonal entries $\bP\bA\bP^{-1}$ are the eigenvalues of $\bA$.
\item \textbf{Power property of similar matrices.} Let $\bB=\bP\bA\bP^{-1}$. Show that $\bB^k = \bP\bA^k\bP^{-1}$ for $k=1,2,\ldots$, that is, $\bB^k$ and $\bA^k$ are similar if $\bB$ and $\bA$ are similar. If one of $\bA$ and $\bB$ is nonsingular, show that $\bB^{-1} = \bP\bA^{-1}\bP^{-1}$ also holds.
\item In the main section, we transform the given matrix $\bA\in\real^{n\times n}$ into its orthogonal similarity transformation. Use Gaussian elimination matrices (see \eqref{equation:elimination_mat}) to transform into its similarity transformation. Discuss the complexity of your algorithm.	
	
\item Show  that if $\bA=\bE\bC\bE^{-1}$ and $\bB=\bF\bC\bF^{-1}$, then $\bA$ and $\bB$ are similar matrices.

\item Show that the matrices $\footnotesize\begin{bmatrix}
4 & 1\\
-1 & 0
\end{bmatrix}$
and 
$\footnotesize\begin{bmatrix}
	1 & 1\\
	0 & 3
\end{bmatrix}$
are similar.

\item \label{prob:hess_poly1} \textbf{Polynomial.} Let $\bA$ and $\bB$ be similar, and consider a polynomial $p(\bC)=\gamma_n\bC^n+\gamma_{n-1}\bC^{n-1}+\ldots+\gamma_0$. Show that 
$p(\bA)$ and $p(\bB)$ are also similar.
\item \label{prob:hess_poly2} \textbf{Polynomial.}  Let $\bA$ and a nonsingular $\bP$ be given, and consider a polynomial $p(\bC)=\gamma_n\bC^n+\gamma_{n-1}\bC^{n-1}+\ldots+\gamma_0$. Show that $p(\bP\bA\bP^{-1}) = \bP p(\bA)\bP^{-1}$.

\item Let  $\bA, \bB\in\real^{n\times n}$ be given. Show that $\bA\bB$ and $\bB\bA$ are similar if at least one of them is nonsingular.

\item \textbf{Similarity transformation.} Let $\bA\in\real^{n\times n}$ and let $\bP\in\real^{n\times n}$ be nonsingular. Show that $\det(\bP^{-1}\bA\bP - \lambda \bI)=\det(\bA-\lambda\bI)$. This again demonstrates that the eigenvalues remain unchanged under similarity transformations (see Definition~\ref{definition:characteristic_polynomial}). 

\item  Let $\bH\in\real^{n\times n}$ be an unreduced upper Hessenberg matrix. Show that $\rank(\bH-\lambda\bI)\geq n-1$ for any $\lambda\in\real$.
\item Let $\bH\in\real^{n\times n}$ be an unreduced upper Hessenberg matrix. Show that its geometric multiplicity is 1 for any eigenvalue (Definition~\ref{definition:eigen_multipli}). 

\item Let $\bA\in\real^{n\times n}$ be given with a lower bandwidth of $p$ (Definition~\ref{defin:matrix-bandwidth}). Provide an algorithm that computes the Hessenberg decomposition of $\bA$ using Householder reflectors or Givens rotations.

\item \textbf{Hessenberg LU.} Let $\bH\in\real^{n\times n}$ be upper Hessenberg. Show that there exists a set of Gaussian elimination matrices $\bE_1, \bE_2, \ldots, \bE_{n-1}$ with entries bounded by unity (see Equation~\eqref{equation:elimination_mat}) and a set of permutation matrices $\bP_1, \bP_2, \bP_{n-1}$ such that $\bE_{n-1}\bP_{n-1}\ldots \bE_2\bP_2\bE_1\bP_1\bH$ is upper triangular. Discuss the complexity of your algorithm.

\item \textbf{Hessenberg QR.} Let $\bH\in\real^{n\times n}$ be upper Hessenberg. Provide an algorithm that computes the QR decomposition of $\bH$ using Givens rotations with a complexity of $\mathcalO(n^2)$ flops.

\item Let $\bH \in \real^{n \times n}$ be upper Hessenberg with an eigenpair  $(\lambda, \bv)$. Provide an algorithm that computes an orthogonal matrix $\bQ$ such that
$ \bQ^\top\bH \bQ = \scriptsize\begin{bmatrix} \lambda & \bu^\top \\ \bzero & \bH_1 \end{bmatrix}, $
where $\bH_1 \in \real^{(n-1) \times (n-1)}$ is also upper Hessenberg. \textit{Hint: Consider $\bQ$ as a product of Givens rotations.}

\item  (Read Chapter~\ref{chapter:spectral-decomposition} first) Let the following matrix be a $4\times 4$ Hessenberg matrix:
$$
\bH =
\begin{bmatrix}
b_1 & c_1 & d_1 & e_1 \\
a_1 & b_2 & c_2 & d_2 \\
0   & a_2 & b_3 & c_3 \\
0   & 0   & a_3 & b_4
\end{bmatrix}.
$$
Show that
\begin{itemize}
\item If $a_1, a_2, a_3$ are all nonzero, and any eigenvalue $\lambda$ of $\bH$ is a real number, then the geometric multiplicity of $\lambda$ must be equal to one.
\item If $\bH$ is similar to a symmetric matrix $\bA$, and the algebraic multiplicity of some eigenvalue $\lambda$ of $\bA$ is greater than 1, then at least one of  $a_1, a_2, a_3$  must be zero.
\end{itemize}

\item Consult \citet{kagstrom2008blocked, bujanovic2018householder} and derive the complexity of the Hessenberg-triangular decomposition.

\item Let $\bA\in\real^{n\times n}$. Show that $\bA$ is idempotent if and only if there exists an orthogonal matrix $\bB\in\real^{n\times n}$ such that $\bA$ and $\bB$ are similar.

\item Show that if $\bA\in\real^{n\times n}$ is similar to an orthogonal matrix, then $\bA^{-1}$ is similar to $\bA^\top$.


\end{problemset}

%% file: chapter-tridiagonal.tex
\newpage
\chapter{Tridiagonal Decomposition}\label{chapter:tridiagonal_decom}
\begingroup
\hypersetup{
	linkcolor=structurecolor,
	linktoc=page,  
}
\minitoc \newpage
\endgroup

\section{Tridiagonal Decomposition}\label{section:tridiagonal_decom}
\lettrine{\color{caligraphcolor}S}
Similar to the Hessenberg decomposition, the tridiagonal decomposition can also simplify a matrix and serve as a preliminary step for other algorithms, thereby reducing their computational burden.
We will begin by providing a formal definition of tridiagonal matrices.

\begin{definition}[Tridiagonal Matrix\index{Tridiagonal matrix}]\label{definition:tridiagonal-hessenbert}
A tridiagonal matrix is a square matrix, where all the entries below the  subdiagonal and the entries above the superdiagonal are zeros. In other words,  a tridiagonal matrix is a special type of \textit{band matrix}.

The definition of a tridiagonal matrix can also be extended to rectangular matrices, and the form can be inferred from the context.

In matrix language, consider any matrix $\bT\in \real^{n\times n}$ with the entry ($i,j$) denoted by $t_{ij}$ for all $i,j\in \{1,2,\ldots, n\}$. Then $\bT$ is a tridiagonal matrix if  $t_{ij}=0$ for all $i\geq j+2$ and $i \leq j-2$.

Let $i$ denote the smallest positive integer such that $t_{i+1, i}=0$ for $i\in \{1,2,\ldots, n-1\}$. By convention, the tridiagonal matrix $\bT$ is called  \textbf{unreduced} if $i=n$.
\end{definition}

Taking a $5\times 5$ matrix as an example, the lower triangular below the  subdiagonal and upper triangular above the superdiagonal are zero in the tridiagonal matrix:
$$
\begin{sbmatrix}{possibly\,\, unreduced}
\boxtimes & \boxtimes & 0 & 0 & 0\\
\boxtimes & \boxtimes & \boxtimes & 0 & 0\\
0 & \boxtimes & \boxtimes & \boxtimes & 0\\
0 & 0 & \boxtimes & \boxtimes & \boxtimes\\
0 & 0 & 0 & \boxtimes & \boxtimes
\end{sbmatrix}
\qquad 
\text{or}
\qquad 
\begin{sbmatrix}{reduced}
\boxtimes & \boxtimes & 0 & 0 & 0\\
\boxtimes & \boxtimes & \boxtimes & 0 & 0\\
0 & \boxtimes & \boxtimes & \boxtimes & 0\\
0 & 0 & \textcolor{mylightbluetext}{0} & \boxtimes & \boxtimes\\
0 & 0 & 0 & \boxtimes & \boxtimes
\end{sbmatrix}.
$$
Obviously, a tridiagonal matrix is a special case of an upper Hessenberg matrix. 
Subsequently, we can formulate the following tridiagonal decomposition:
\index{Decomposition: Tridiagonal}
\begin{theoremHigh}[Tridiagonal Decomposition]\label{theorem:tridiagonal-decom}
Let $\bA$ be any $n\times n$ symmetric matrix. Then it can be factored as 
$$
\bA = \bQ\bT\bQ^\top \qquad \text{or} \qquad \bT = \bQ^\top \bA\bQ, 
$$
where $\bT$ is a \textit{symmetric} tridiagonal matrix, and $\bQ$ is an orthogonal matrix.
\end{theoremHigh}
The existence of the tridiagonal matrix can be readily demonstrated by applying the Hessenberg decomposition to the symmetric matrix $\bA$.  
Nevertheless, we will observe that the computational complexity is significantly reduced compared to that of the Hessenberg decomposition, due to the symmetric property of the matrix.

\section{Computing  Tridiagonal Decomposition}\label{section:compute-tridiagonal}

As symmetry introduces zeros in both rows and columns, we can save computational effort by disregarding these extra zeros.

An example of a $5\times 5$ matrix is provided below, where $\boxtimes$ or a letter represents a value that is not necessarily zero, and \textbf{boldface} indicates the value has just been changed: 
$$
\begin{aligned}
\begin{sbmatrix}{\bA}
\boxtimes & \boxtimes & \boxtimes & \boxtimes & \boxtimes \\
\boxtimes & \boxtimes & \boxtimes & \boxtimes & \boxtimes\\
\boxtimes & \boxtimes & \boxtimes & \boxtimes & \boxtimes\\
\boxtimes & \boxtimes & \boxtimes & \boxtimes & \boxtimes\\
\boxtimes & \boxtimes & \boxtimes & \boxtimes & \boxtimes
\end{sbmatrix}
\stackrel{\bH_1\times}{\rightarrow}
&\begin{sbmatrix}{\bH_1\bA}
\boxtimes & \boxtimes & \boxtimes & \boxtimes & \boxtimes \\
\bm{a} & \bm{\boxtimes} & \bm{\boxtimes} & \bm{\boxtimes} & \bm{\boxtimes}\\
\bm{0} & \bm{\boxtimes} & \bm{\boxtimes} & \bm{\boxtimes} & \bm{\boxtimes}\\
\bm{0} & \bm{\boxtimes} & \bm{\boxtimes} & \bm{\boxtimes} & \bm{\boxtimes}\\
\bm{0} & \bm{\boxtimes} & \bm{\boxtimes} & \bm{\boxtimes} & \bm{\boxtimes}
\end{sbmatrix}
\stackrel{\times\bH_1^\top}{\rightarrow}
\begin{sbmatrix}{\bH_1\bA\bH_1^\top}
\boxtimes & \bm{a} & \bm{0} & \bm{0} & \bm{0} \\
a & \bm{\boxtimes} & \bm{\boxtimes} & \bm{\boxtimes} & \bm{\boxtimes}\\
0 & \bm{\boxtimes} & \bm{\boxtimes} & \bm{\boxtimes} & \bm{\boxtimes}\\
0 & \bm{\boxtimes} & \bm{\boxtimes} & \bm{\boxtimes} & \bm{\boxtimes}\\
0 & \bm{\boxtimes} & \bm{\boxtimes} & \bm{\boxtimes} & \bm{\boxtimes}
\end{sbmatrix}\\
\stackrel{\bH_2\times}{\rightarrow}
&\begin{sbmatrix}{\bH_2\bH_1\bA\bH_1^\top}
\boxtimes & a & 0 & 0 & 0 \\
a & \boxtimes & \boxtimes & \boxtimes & \boxtimes \\
0 & \bm{b} & \bm{\boxtimes} & \bm{\boxtimes} & \bm{\boxtimes}\\
0 & \bm{0} & \bm{\boxtimes} & \bm{\boxtimes} & \bm{\boxtimes}\\
0 & \bm{0} & \bm{\boxtimes} & \bm{\boxtimes} & \bm{\boxtimes}
\end{sbmatrix}
\stackrel{\times\bH_2^\top}{\rightarrow}
\begin{sbmatrix}{\bH_2\bH_1\bA\bH_1^\top\bH_2^\top}
\boxtimes & a & 0 & 0 & 0 \\
a & \boxtimes & \bm{b} & \bm{0} & \bm{0} \\
0 & b & \bm{\boxtimes} & \bm{\boxtimes} & \bm{\boxtimes}\\
0 & 0 & \bm{\boxtimes} & \bm{\boxtimes} & \bm{\boxtimes}\\
0 & 0 & \bm{\boxtimes} & \bm{\boxtimes} & \bm{\boxtimes}
\end{sbmatrix}\\
\stackrel{\bH_3\times}{\rightarrow}
&\begin{sbmatrix}{\bH_3\bH_2\bH_1\bA\bH_1^\top\bH_2^\top}
\boxtimes & a & 0 & 0 & 0 \\
a & \boxtimes & b & 0 & 0 \\
0 & b & \boxtimes & \boxtimes & \boxtimes \\
0 & 0 & \bm{c} & \bm{\boxtimes} & \bm{\boxtimes}\\
0 & 0 & \bm{0} & \bm{\boxtimes} & \bm{\boxtimes}
\end{sbmatrix}
\stackrel{\times\bH_3^\top}{\rightarrow}
\begin{sbmatrix}{\bH_3\bH_2\bH_1\bA\bH_1^\top\bH_2^\top\bH_3^\top}
\boxtimes & a & 0 & 0 & 0 \\
a & \boxtimes & b & 0 & 0 \\
0 & b & \boxtimes & \bm{c} & \bm{0}\\
0 & 0 & c & \bm{\boxtimes} & \bm{\boxtimes}\\
0 & 0 & 0 & \bm{\boxtimes} & \bm{\boxtimes}
\end{sbmatrix}.
\end{aligned}
$$

The complete procedure for computing the tridiagonal decomposition of a symmetric matrix is presented in Algorithm~\ref{alg:tridiagonal-decomposition-householder}.
\begin{theorem}[Algorithm Complexity: Tridiagonalization via Householder]\label{theorem:tridiagonal-householder}
Algorithm~\ref{alg:tridiagonal-decomposition-householder} requires $\sim \frac{4}{3}n^3$ flops to compute a tridiagonal decomposition of an $n\times n$ symmetric matrix. Furthermore, if the orthogonal matrix $\bQ$ is needed explicitly, an additional $\sim 2n^3$ flops are required.
\end{theorem}

\begin{proof}[of Theorem~\ref{theorem:tridiagonal-householder}]
The complexity of step 9 warrants  close scrutiny. It arises from both the left  and right Householder updates. 
Given that we have computed the first column of $\bT_{\textcolor{mylightbluetext}{i+1}:n,i:n}$ and the first row of $\bT_{i:n,\textcolor{mylightbluetext}{i+1}:n}$ explicitly, the left and right Householder updates take the following forms (where the \textcolor{winestain}{red}-colored text indicates differences from the process of computing the Hessenberg decomposition for non-symmetric matrices):
\begin{mdframed}[hidealllines=\mdframehidelineNote,backgroundcolor=\mdframecolor,frametitle={Where Does the Step 9 Come From: A Splitting Way}]
Left: update the second to the last  columns  of $\bT_{\textcolor{mylightbluetext}{i+1}:n,i:n}$, i.e., working on $\bT_{i+1:n,i+1:n}\in \real^{(n-i)\times (n-i)}$:
$$
\begin{aligned}
	\bT_{i+1:n,i+1:n} &\leftarrow (\bI-2\bu_i\bu_i^\top)\bT_{i+1:n,i+1:n} \\
	&= \bT_{i+1:n,i+1:n} - 2\bu_i\bu_i^\top\bT_{i+1:n,i+1:n} \in \real^{(n-i)\times(n-i)}, \qquad \text{($4(n-i)^2$ flops);}
\end{aligned}
$$
Right: update the second to the last  rows  of $\bT_{i:n,\textcolor{mylightbluetext}{i+1}:n}$, i.e., working on $\bT_{i+1:n,i+1:n}\in \real^{(n-i)\times (n-i)}$:
$$
\begin{aligned}
\bT_{\textcolor{winestain}{i+1}:n,i+1:n} &\leftarrow \bT_{\textcolor{winestain}{i+1}:n,i+1:n}(\bI-2\bu_i\bu_i^\top) \\
&=\bT_{\textcolor{winestain}{i+1}:n,i+1:n} - \bT_{\textcolor{winestain}{i+1}:n,i+1:n} 2\bu_i\bu_i^\top\in \real^{\textcolor{winestain}{(n-i)}\times(n-i)}, \qquad\text{($4(n-i)^2$ flops).}
\end{aligned}
$$
\end{mdframed}
As both updates now operate on the same submatrix  $\bT_{i+1:n,i+1:n}$, we can combine them into step 9 of the algorithm.  Let $\bZ \triangleq \bT_{i+1:n,i+1:n} \in \real^{(n-i)\times (n-i)}$. A delve into the reduction of the complexity in step 9 proceeds  as follows:
$$
\begin{aligned}
\bZ &\leftarrow  \bH_i \bZ \bH_i   \qquad  &\text{($\bH_i$ is the $i$-th Householder reflector)}\\ 
&= (\bI-2\bu_i\bu_i^\top)\bZ (\bI-2\bu_i\bu_i^\top)\\
&= (\bZ - 2\bu_i\underbrace{\bu_i^\top \bZ}_{\by^\top })(\bI-2\bu_i\bu_i^\top)  \\
&= (\bZ - 2\bu_i\by^\top )(\bI-2\bu_i\bu_i^\top)  & \text{(let $\by^\top \triangleq \bu_i^\top \bZ\rightarrow \by \triangleq \bZ\bu_i$)}\\
&= \bZ - 2\bu_i\by^\top - 2\bZ\bu_i\bu_i^\top + 4\beta \bu_i\bu_i^\top & \text{(let $\beta \triangleq \by^\top\bu_i$)}\\
&= \bZ - (2\bu_i\by^\top-2\beta \bu_i\bu_i^\top)  - (2\underbrace{\bZ\bu_i}_{\by}\bu_i^\top -2\beta \bu_i\bu_i^\top )  \\
&= \bZ - (2\bu_i\by^\top-2\beta \bu_i\bu_i^\top)  - (2\by\bu_i^\top -2\beta \bu_i\bu_i^\top )  &\text{($\by = \bZ\bu_i$)}\\
&= \bZ - \{\bu_i\underbrace{(2\by^\top- 2\beta \bu_i^\top)}_{\bx^\top }  +\underbrace{ (2\by^\top-2\beta \bu_i^\top)^\top}_{\bx } \bu_i^\top   \}\\
&= \bZ - \{\bu_i\bx^\top  +( \bu_i\bx^\top)^\top   \}.   &(\text{let $\bx\triangleq 2(\by- \beta \bu_i)$})
\end{aligned}
$$
The costs come from 
\begin{itemize}
\item $\by \leftarrow \bZ\bu_i$: the computational cost is $[2(n-i)-1](n-i) = 2(n-i)^2 - (n-i)$ flops from Lemma~\ref{lemma:matrix-multi-complexity} (on the complexity of a matrix multiplication);
\item $\beta \leftarrow \by^\top\bu_i$: $2(n-i)-1$ flops;
\item  $\bx\leftarrow 2(\by- \beta \bu_i)$: $3(n-i)$ flops;
\item $\bu_i\bx^\top$: $(n-i)^2$ flops;
\item $( \bu_i\bx^\top)^\top $: 0 flops;
\item $\bu_i\bx^\top  +( \bu_i\bx^\top)^\top$: it involves $1+2+\ldots+(n-i)= \frac{(n-i)^2+(n-i)}{2}$ additions since it results in a symmetric matrix;
\item $\underbrace{\bZ}_{\text{symmetric}} - \underbrace{\{\bu_i\bx^\top  +( \bu_i\bx^\top)^\top   \}}_{\text{symmetric}} $: it involves $1+2+\ldots+(n-i)= \frac{(n-i)^2+(n-i)}{2}$ subtractions, as both matrices are symmetric.
\end{itemize}
If we keep only the leading terms from step 4 to step 9, the total complexity for loop $i$ is given by 
$f(i) = 4(n-i)^2$ flops. Summing over all iterations, we obtain the final cost:
$$
\text{cost=} f(1)+f(2)+\ldots +f(n-2)=~ \frac{4}{3}n^3 \text{ flops}.
$$
To determine the computational  complexity of calculating the orthogonal matrix $\bQ$, it is the same as that in Theorem~\ref{theorem:hessenbert-householder}.
This completes the proof.
\end{proof}

\index{Householder}
\begin{algorithm}[h] 
\caption{Tridiagonal Decomposition via the Householder Reflector} 
\label{alg:tridiagonal-decomposition-householder} 
\begin{algorithmic}[1] 
\Require Matrix $\bA$ with size $n\times n $; 
\Statex \textbf{Stage A: Obtain the tridiagonal matrix}
\State Initially set $\bT \leftarrow \bA$;
\For{$i=1$ to $n-2$} 
\State $\ba \leftarrow \bT_{i+1:n,i}$, i.e., first column of $\bT_{i+1:n,i:n}\in \real^{(n-i)\times(n-i+1)}$;
\State $r \leftarrow \normtwo{\ba}$; \Comment{$2(n-i)$ flops}
\State $\bu_i \leftarrow \ba-r\be_1$; \Comment{1 flop}
\State $\bu_i \leftarrow \bu_i / \normtwo{\bu_i}$; \Comment{$3(n-i)$ flops}
\State $\bT_{i+1,i}\leftarrow r$, $\bT_{i+2:n,i}\leftarrow\bzero$, i.e., set the value of first column of $\bT_{\textcolor{mylightbluetext}{i+1}:n,i:n}$; \Comment{0 flops}
\State  $\bT_{i,i+1}\leftarrow r$, $\bT_{i,i+2:n}\leftarrow\bzero$, i.e., set the value of first row of $\bT_{i:n,\textcolor{mylightbluetext}{i+1}:n}$; \Comment{0 flops}

\State Left and Right: let $\bZ = \bT_{i+1:n,i+1:n} \in \real^{(n-i)\times (n-i)}$, 
$$
\begin{aligned}
\bZ &\leftarrow \bH_i \bZ \bH_i   \qquad \text{($\bH_i$ is the $i$-th Householder reflector)}\\ 
&= (\bI-2\bu_i\bu_i^\top)\bZ (\bI-2\bu_i\bu_i^\top)\\
&= \bZ - \bZ\cdot 2\bu_i\bu_i^\top  - 2\bu_i\bu_i^\top\bZ+ 2\bu_i\bu_i^\top\bZ\cdot 2\bu_i\bu_i^\top,\qquad (\text{$4(n-i)^2$ flops}); \\
\end{aligned}
$$
\EndFor
\State Output $\bT$ as the tridiagonal matrix;
\Statex \textbf{Stage B: Obtain the orthogonal matrix}
\State Get $\bQ=\bH_1\bH_2\ldots\bH_{n-2}$, where $\bH_i$'s are Householder reflectors.
\State Initially set $\bQ \leftarrow \bH_1$;
\For{$i=1$ to $n-3$} 
\State Compute $\bQ$:
$$\begin{aligned}
&\bQ_{1:n,i+2:n} \leftarrow \bQ_{1:n,i+2:n}\bH_{i+1}= \bQ_{1:n,i+2:n}(\bI - 2\bu_{i+1}\bu_{i+1}^\top)\\
&=\bQ_{1:n,i+2:n}-\bQ_{1:n,i+2:n}2\bu_{i+1}\bu_{i+1}^\top \in \real^{n\times (n-i-1)}, \qquad \text{($4n(n-i-1) - n$ flops)};\\
\end{aligned}$$
\EndFor
\State Output $\bQ$ as the orthogonal matrix.
\end{algorithmic} 
\end{algorithm}

\section{Properties of  Tridiagonal Decomposition}\label{section:tridiagonal-decomposition}
Similarly to the Hessenberg decomposition, the tridiagonal decomposition is not unique. However, and most importantly, if we restrict the elements in the  subdiagonal of the tridiagonal matrix $\bT$ to be positive (if possible), then the tridiagonal decomposition $\bA=\bQ\bT\bQ^\top$ is uniquely determined by $\bA$ and the first column of $\bQ$ under mild conditions.
\begin{theorem}[Implicit Q Theorem for Tridiagonal\index{Implicit Q theorem}]\label{theorem:implicit-q-tridiagonal}
Let $\bA\in\real^{n\times n}$ be given and 
suppose two tridiagonal decompositions of a symmetric matrix $\bA$ are given by $\bA=\bU\bT\bU^\top=\bV\bG\bV^\top$, where $\bU=[\bu_1, \bu_2, \ldots, \bu_n]$ and $\bV=[\bv_1, \bv_2, \ldots, \bv_n]$ are the column partitions of $\bU$ and $\bV$, respectively. Suppose further that $k$ is the smallest positive integer for which $t_{k+1,k}=0$, where $t_{ij}$ is the entry $(i,j)$ of $\bT$. Then, 
\begin{itemize}
\item If $\bu_1=\bv_1$, then $\bu_i = \pm \bv_i$ and $|t_{i,i-1}| = |g_{i,i-1}|$ for $i\in \{2,3,\ldots,k\}$. 
\item When $k=n$, the tridiagonal matrix $\bT$ is known as unreduced. However, if $k<n$, then $g_{k+1,k}=0$.
\end{itemize}
\end{theorem}
From the above theorem, we observe that  if we restrict the elements in the lower sub-diagonal of the tridiagonal matrix $\bT$ to be positive (if possible), i.e., \textit{unreduced}, then the tridiagonal decomposition $\bA=\bQ\bT\bQ^\top$ is uniquely determined by $\bA$ and the first column of $\bQ$. This again is similar to what we have claimed on the uniqueness of the QR decomposition (Corollary~\ref{corollary:unique-qr}).

Similarly, an unreduced tridiagonal decomposition can be obtained from the implication of the Krylov matrix (Definition~\ref{definition:krylov-matrix}).
\begin{theorem}[Unreduced Tridiagonal]\label{theorem:implicit-q-tridiagonal-v2}
Let $\bA\in\real^{n\times n}$ be given and suppose there exists an orthogonal matrix $\bQ$ such that $\bA\in \real^{n\times n}$ can be factored as $\bA = \bQ\bT\bQ^\top$. Then, $\bQ^\top\bA\bQ=\bT$ is an unreduced tridiagonal matrix if and only if $\bR=\bQ^\top \bK(\bA, \bq_1, n)$ is nonsingular and upper triangular, where $\bq_1$ is the first column of $\bQ$.  

If $\bR$ is singular and $k$ is the smallest index satisfying $r_{kk}=0$, then $k$ is also the smallest index for which $t_{k,k-1}=0$.
\end{theorem}

\begin{problemset}
\item Let $\bH\in\real^{n\times n}$ be upper Hessenberg. Provide an algorithm that computes the decomposition $\bH\bR=\bR\bT$, where $\bR$ is unit upper triangular, and $\bT$ is tridiagonal.
\item Based on the proofs of Theorems~\ref{theorem:implicit-q-hessenberg} and~\ref{theorem:implicit-q-hessenberg-v2}, prove Theorems~\ref{theorem:implicit-q-tridiagonal} and \ref{theorem:implicit-q-tridiagonal-v2}.

\item \citep{golub2013matrix} Let $\bA=\bS+\sigma\bu\bu^\top\in\real^{n\times n}$, where $\bS\in\real^{n\times n}$ is skew-symmetric (see Definition~\ref{definition:speci_mat}), $\bu\in\real^n$, and $\sigma\in\real$. Show that there exists an orthogonal matrix $\bQ$ such that $\bQ^\top\bA\bQ = \bT+\sigma\be_1\be_1^\top$, where $\bT$ is tridiagonal and skew-symmetric.


\item Let $\gamma_0, \gamma_1, \ldots, \gamma_n>0$. Show that the following $n\times n$ tridiagonal matrix is positive definite:
$$
\scriptsize
\begin{bmatrix}
\gamma_0+\gamma_1 & -\gamma_1 & 0 & \ldots & 0\\
-\gamma_1 & \gamma_1+\gamma_2 & -\gamma_2 & \ldots & 0\\
0 & -\gamma_2 & \gamma_2+\gamma_3 & \ldots & 0\\
\vdots & \vdots & \vdots & \ddots & \vdots\\
0 & 0 & 0 & \ldots & \gamma_{n-1}+\gamma_n
\end{bmatrix}.
$$
\textit{Hint: Consider the leading principal minors in Section~\ref{section:band_lu_wop}.}

\index{Toeplitz matrix}
\item \citep{higham2002accuracy} Let  $\bT_n(a,b,c)\in\real^{n\times n}$ be a tridiagonal matrix of the following form:
$$
\textbf{(toeplitz tridiagonal matrix)}:
\qquad 
\bT_n(a,b,c)=
\scriptsize
\begin{bmatrix}
b & c & 0 & \ldots & 0\\
a & b & c & \ldots & 0\\
0 & a & b & \ldots & 0\\
\vdots & \vdots & \vdots & \ddots & \vdots\\
0 & 0 & 0 & \ldots & b
\end{bmatrix}.
$$
Show that the eigenvalues of $\bT_n(a,b,c)$ are $b+2\sqrt{ac}\cos(\frac{k\pi}{n+1})$ for $k\in\{1,2,\ldots,n\}$.

\item \citep{noschese2013tridiagonal} Show that $\bT_n(a,b,c)$ is normal if and only if $\abs{a}=\abs{c}$.

\item  Let $\bT\in\real^{n\times n}$ be an unreduced tridiagonal matrix. Show that $\rank(\bT-\lambda\bI)\geq n-1$ for any $\lambda\in\real$.
\item Let $\bT\in\real^{n\times n}$ be an unreduced tridiagonal matrix. Show that its geometric multiplicity is 1 for any eigenvalue (Definition~\ref{definition:eigen_multipli}). 

\item Let $\bA\in\real^{n\times n}$ be tridiagonal. Show that if $a_{i,i+1}a_{i+1,i} > 0$ for all $i\in\{1,2,\ldots, n - 1\}$, then $\bA$ has $n$ distinct real eigenvalues. Moreover, show that if $a_{i,i+1}a_{i+1,i} \geq 0$ for all  $i\in\{1,2,\ldots, n - 1\}$, then all \ eigenvalues of $\bA$ are real. \textit{Hint: Use the Jordan decomposition discussed in Chapter~\ref{chapter:eig_jordan}.}
\end{problemset}

%% file: chapter-bidiagonal.tex
\newpage
\chapter{Bidiagonal Decomposition}\label{section:bidiagonal-decompo}
\begingroup
\hypersetup{
	linkcolor=structurecolor,
	linktoc=page,  
}
\minitoc \newpage
\endgroup

\section{Bidiagonal Decomposition}
\lettrine{\color{caligraphcolor}W}
When the matrix is not square and symmetric, obtaining a tridiagonal form becomes challenging. However, we can take a step back and seek a decomposition involving two different orthogonal matrices. We start by providing a rigorous definition of upper bidiagonal matrices.

\begin{definition}[Upper Bidiagonal Matrix\index{Bidiagonal matrix}]\label{definition:bidiagonal-matrix}	
An \textit{upper bidiagonal} matrix is a square matrix characterized by a banded structure with nonzero entries along both the \textit{main diagonal} and the \textit{superdiagonal} (i.e., the elements directly above the main diagonal). This means the matrix has exactly two nonzero diagonals.

Conversely, if the diagonal below the main diagonal contains  nonzero entries and the superdiagonal contains zero entries, the matrix is referred to as \textit{lower bidiagonal}.

The concept of bidiagonal matrices can also be extended to rectangular matrices, and the specific form can be inferred from the context.
\end{definition}

Taking a $7\times 5$ matrix as an example, the lower triangular below the main diagonal and the upper triangular above the superdiagonal are zero in the upper bidiagonal matrix:
$$
\footnotesize
\begin{bmatrix}
\boxtimes & \boxtimes & 0 & 0 & 0\\
0 & \boxtimes & \boxtimes & 0 & 0\\
0 & 0 & \boxtimes & \boxtimes & 0\\
0 & 0 & 0 & \boxtimes & \boxtimes\\
0 & 0 & 0 & 0 & \boxtimes\\
0 & 0 & 0 & 0 & 0\\
0 & 0 & 0 & 0 & 0
\end{bmatrix}.
$$

Then we have the following bidiagonal decomposition:
\index{Decomposition: Bidiagonal}
\begin{theoremHigh}[Bidiagonal Decomposition]\label{theorem:Golub-Kahan-Bidiagonalization-decom}
Let $\bA\in\real^{m\times n}$ be any $m\times n$ matrix. Then it  can be factored as 
$$
\bA = \bU\bB\bV^\top \qquad \text{or} \qquad \bB= \bU^\top \bA\bV, 
$$
where $\bB$ is an upper bidiagonal matrix, and $\bU\in\real^{m\times m}$ and $\bV\in\real^{n\times n}$ are orthogonal matrices.
Specifically, $\bV$ has the following structure 
$$
\bV = 
\begin{bmatrix}
1 & \bzero \\
\bzero & \bQ 
\end{bmatrix},
$$
where $\bQ\in\real^{(n-1)\times (n-1)}$ is orthogonal.
\end{theoremHigh}
We will observe that the bidiagonalization  bears a resemblance to the structure of a singular value decomposition. The primary difference lies in the values of $\bB$ in bidiagonal form, which contain nonzero entries on the superdiagonal such that it will be shown to play an important role in the calculation of the singular value decomposition (Chapter~\ref{section:eigenvalue-problem}).

\section{Existence of  Bidiagonal Decomposition: Golub-Kahan Approach}

In previous work, Householder reflectors were employed to triangularize matrices by introducing zeros below the main diagonal, facilitating the QR decomposition, and introducing zeros below the subdiagonal to achieve the Hessenberg decomposition.
A similar approach can be employed to find the bidiagonal decomposition.

\index{Householder}
\subsubsection*{\textbf{First Step 1.1: Introduce Zeros for the First Column}}	
For brevity, we assume $m\geq n$; see Problem~\ref{prob:bidia_mln} for further considerations.
Let $\bA=[\ba_1, \ba_2, \ldots, \ba_n]$ be the column partition of $\bA$, where each $\ba_i \in \real^{m}$. 
We can construct the Householder reflector as follows:
$$
r_1 = \normtwo{\ba_1}, \qquad \bu_1 = \frac{\ba_1 - r_1 \be_1}{\normtwo{\ba_1 - r_1 \be_1}} ,\qquad \text{and}\qquad \bH_1 = \bI - 2\bu_1\bu_1^\top \in \textcolor{mylightbluetext}{\real^{m\times m}},
$$
where $\be_1$ represents the first unit basis vector in $\textcolor{mylightbluetext}{\real^{m}}$, i.e., $\be_1=[1;0;0;\ldots;0]\in \textcolor{mylightbluetext}{\real^{m}}$.
The matrix $\bH_1$ is  symmetric and orthogonal  (from the definition of Householder reflectors).
Applying $\bH_1$ to $\bA$ introduces zeros in the first column of $\bA$ below the $(1,1)$ entry, effectively reflecting $\ba_1$ to $r_1 \be_1$.

An example of a $7\times 5$ matrix is shown as follows, where $\boxtimes$ represents a value that is not necessarily zero, and \textbf{boldface} indicates the value has just been changed:
$$
\begin{aligned}
	\footnotesize
\begin{sbmatrix}{\bA}
\boxtimes & \boxtimes & \boxtimes & \boxtimes & \boxtimes \\
\boxtimes & \boxtimes & \boxtimes & \boxtimes & \boxtimes\\
\boxtimes & \boxtimes & \boxtimes & \boxtimes & \boxtimes\\
\boxtimes & \boxtimes & \boxtimes & \boxtimes & \boxtimes\\
\boxtimes & \boxtimes & \boxtimes & \boxtimes & \boxtimes\\
\boxtimes & \boxtimes & \boxtimes & \boxtimes & \boxtimes\\
\boxtimes & \boxtimes & \boxtimes & \boxtimes & \boxtimes
\end{sbmatrix}
\stackrel{\bH_1\times}{\rightarrow}
&\footnotesize\begin{sbmatrix}{\bH_1\bA}
\bm{\boxtimes} & \bm{\boxtimes} & \bm{\boxtimes} & \bm{\boxtimes} & \bm{\boxtimes}\\
\bm{0} & \bm{\boxtimes} & \bm{\boxtimes} & \bm{\boxtimes} & \bm{\boxtimes}\\
\bm{0} & \bm{\boxtimes} & \bm{\boxtimes} & \bm{\boxtimes} & \bm{\boxtimes}\\
\bm{0} & \bm{\boxtimes} & \bm{\boxtimes} & \bm{\boxtimes} & \bm{\boxtimes}\\
\bm{0} & \bm{\boxtimes} & \bm{\boxtimes} & \bm{\boxtimes} & \bm{\boxtimes}\\
\bm{0} & \bm{\boxtimes} & \bm{\boxtimes} & \bm{\boxtimes} & \bm{\boxtimes}\\
\bm{0} & \bm{\boxtimes} & \bm{\boxtimes} & \bm{\boxtimes} & \bm{\boxtimes}\\
\end{sbmatrix}.
\end{aligned}
$$
Up to this point, the process is identical to the QR decomposition using Householder reflectors, as detailed in Section~\ref{section:qr-via-householder}.
To proceed  further, the act of introducing zeros above the superdiagonal of $\bH_1\bA$ is equivalent to introducing zeros below the subdiagonal of $(\bH_1\bA)^\top$.

\subsubsection*{\textbf{First Step 1.2: Introduce Zeros for the First Row}}	

Now consider the \textit{transpose} of $\bH_1\bA$, that is $(\bH_1\bA)^\top =\bA^\top\bH_1^\top \in \real^{n\times m}$. And let the column partition of  $\bA^\top\bH_1^\top$ be given by $\bA^\top\bH_1^\top \triangleq [\bz_1, \bz_2, \ldots, \bz_m]$, where each $\bz_i \in \real^n$. 
Suppose $\bar{\bz}_1, \bar{\bz}_2, \ldots, \bar{\bz}_m \in \real^{n-1}$ are vectors obtained by removing the first component from $\bz_i$'s. 
We can  construct another Householder reflector as follows:
$$
r_1 = \normtwo{\bar{\bz}_1}, \qquad \bv_1 = \frac{\bar{\bz}_1 - r_1 \be_1}{\normtwo{\bar{\bz}_1 - r_1 \be_1}},  \qquad \text{and}\qquad \widetilde{\bL}_1 = \bI - 2\bv_1\bv_1^\top \in\textcolor{mylightbluetext}{\real^{(n-1)\times (n-1)}},
$$
where $\be_1$ now represents the first unit basis vector in $\textcolor{mylightbluetext}{\real^{n-1}}$, i.e., $\be_1=[1;0;0;\ldots;0]\in \textcolor{mylightbluetext}{\real^{n-1}}$. To introduce zeros below the subdiagonal and operate on the submatrix $(\bA^\top\bH_1^\top)[2:n,1:m]$, we append the Householder reflector into $\bL_1$, defined as
$
\bL_1 = \footnotesize\begin{bmatrix}
1 &\bzero \\
\bzero & \widetilde{\bL}_1
\end{bmatrix},
$
where both $\bL_1$ and $\widetilde{\bL}_1$ are  orthogonal  and symmetric (from the definition of  Householder reflectors).
In this case, $\bL_1(\bA^\top\bH_1^\top)$  introduces zeros in the first column of $(\bA^\top\bH_1^\top)$ below entry (2,1), i.e., reflects $\bar{\bz}_1$ to $r_1\be_1$. The first row of $(\bA^\top\bH_1^\top)$ remains unaffected and kept unchanged by Remark~\ref{remark:left-right-identity}, preserving the zeros introduced in step (1.1). 

Returning to the original \textit{untransposed} matrix $\bH_1\bA$, postmultiplying by $\bL_1^\top$ introduces zeros in the first row to the right of the entry (1,2).
Once again, using the example provided earlier, the procedure on the $7\times 5$ matrix is shown as follows, where $\boxtimes$ represents a value that is not necessarily zero, and \textbf{boldface} indicates the value has just been changed:
$$
\begin{aligned}
	\footnotesize
\begin{sbmatrix}{\bA^\top\bH_1^\top}
\boxtimes & 0 & 0 & 0 & 0 & 0 & 0 \\
\boxtimes & \boxtimes & \boxtimes & \boxtimes & \boxtimes & \boxtimes & \boxtimes\\
\boxtimes & \boxtimes & \boxtimes & \boxtimes & \boxtimes & \boxtimes & \boxtimes\\
\boxtimes & \boxtimes & \boxtimes & \boxtimes & \boxtimes & \boxtimes & \boxtimes\\
\boxtimes & \boxtimes & \boxtimes & \boxtimes & \boxtimes & \boxtimes & \boxtimes\\
\end{sbmatrix}
\stackrel{\bL_1\times}{\rightarrow}
\footnotesize
\begin{sbmatrix}{\bL_1 \bA^\top\bH_1^\top}
\boxtimes & 0 & 0 & 0 & 0 & 0 & 0 \\
\bm{\boxtimes} & \bm{\boxtimes} & \bm{\boxtimes} & \bm{\boxtimes} & \bm{\boxtimes} & \bm{\boxtimes} & \bm{\boxtimes}\\
\bm{0} & \bm{\boxtimes} & \bm{\boxtimes} & \bm{\boxtimes} & \bm{\boxtimes} & \bm{\boxtimes} & \bm{\boxtimes}\\
\bm{0} & \bm{\boxtimes} & \bm{\boxtimes} & \bm{\boxtimes} & \bm{\boxtimes} & \bm{\boxtimes} & \bm{\boxtimes}\\
\bm{0} & \bm{\boxtimes} & \bm{\boxtimes} & \bm{\boxtimes} & \bm{\boxtimes} & \bm{\boxtimes} & \bm{\boxtimes}\\
\end{sbmatrix}
\stackrel{(\cdot)^\top}{\rightarrow}
\footnotesize
\begin{sbmatrix}{\bH_1\bA\bL_1^\top }
\boxtimes & \bm{\boxtimes} & \bm{0} & \bm{0} & \bm{0} \\
0 & \bm{\boxtimes} & \bm{\boxtimes} & \bm{\boxtimes} & \bm{\boxtimes}\\
0 & \bm{\boxtimes} & \bm{\boxtimes} & \bm{\boxtimes} & \bm{\boxtimes}\\
0 & \bm{\boxtimes} & \bm{\boxtimes} & \bm{\boxtimes} & \bm{\boxtimes}\\
0 & \bm{\boxtimes} & \bm{\boxtimes} & \bm{\boxtimes} & \bm{\boxtimes}\\
0 & \bm{\boxtimes} & \bm{\boxtimes} & \bm{\boxtimes} & \bm{\boxtimes}\\
0 & \bm{\boxtimes} & \bm{\boxtimes} & \bm{\boxtimes} & \bm{\boxtimes}\\
\end{sbmatrix}.
\end{aligned}
$$

In short, $\bH_1\bA\bL_1^\top$ finishes the first step of introducing zeros to the first column and the first row of $\bA$.

\subsubsection*{\textbf{Second Step 2.1: Introduce Zeros for the Second Column}}	
Let $\bB \triangleq \bH_1\bA\bL_1^\top$, where all entries  in the first column  below entry (1,1) are zero, and similarly, all entries in the first row to the right of entry $(1,2)$ are also zero.
In this step, the objective is to introduce zeros in the second column below entry (2,2). 
To achieve this, we define $\bB_2 \triangleq \bB[2:m,2:n]\triangleq[\bb_1, \bb_2, \ldots, \bb_{n-1}] \in \real^{(m-1)\times (n-1)}$. 
We can once again construct a Householder reflector as follows:
$$
r_1 = \normtwo{\bb_1},\qquad	\bu_2 = \frac{\bb_1 - r_1 \be_1}{\normtwo{\bb_1 - r_1 \be_1}}, \qquad  \text{and}\qquad \widetilde{\bH}_2 = \bI - 2\bu_2\bu_2^\top\in \textcolor{mylightbluetext}{\real^{(m-1)\times (m-1)}},
$$
where $\be_1$ now represents the first unit basis vector in $\textcolor{mylightbluetext}{\real^{m-1}}$ i.e., $\be_1=[1;0;0;\ldots;0]\in \textcolor{mylightbluetext}{\real^{m-1}}$. To introduce zeros below the main diagonal and perform operations on the submatrix $\bB[2:m,2:n]$, we append the Householder reflector into $\bH_2$, defined as 
$
\bH_2 = \footnotesize\begin{bmatrix}
1 &\bzero \\
\bzero & \widetilde{\bH}_2
\end{bmatrix}.
$
In this scenario, it becomes evident that  $\bH_2(\bH_1\bA\bL_1^\top)$ will not change the first  row of $(\bH_1\bA\bL_1^\top)$, as indicated in Remark~\ref{remark:left-right-identity}. Additionally,  since the Householder reflector cannot reflect a zero vector, the zeros in the first column will remain unaffected as well. 

Continuing from the previous example, the process applied to a $7\times 5$ matrix is shown as follows, where $\boxtimes$ represents a value that is not necessarily zero, and \textbf{boldface} indicates the value has just been changed:
$$
\begin{aligned}
	\footnotesize
\begin{sbmatrix}{\bH_1\bA\bL_1^\top }
\boxtimes & \boxtimes & 0 & 0 & 0 \\
0 & \boxtimes & \boxtimes & \boxtimes & \boxtimes\\
0 & \boxtimes & \boxtimes & \boxtimes & \boxtimes\\
0 & \boxtimes & \boxtimes & \boxtimes & \boxtimes\\
0 & \boxtimes & \boxtimes & \boxtimes & \boxtimes\\
0 & \boxtimes & \boxtimes & \boxtimes & \boxtimes\\
0 & \boxtimes & \boxtimes & \boxtimes & \boxtimes\\
\end{sbmatrix}
\stackrel{\bH_2\times }{\rightarrow}
\footnotesize
\begin{sbmatrix}{\bH_2\bH_1\bA\bL_1^\top }
\boxtimes & \boxtimes & 0& 0 & 0 \\
0 & \bm{\boxtimes} & \bm{\boxtimes} & \bm{\boxtimes} & \bm{\boxtimes}\\
0 & \bm{0} & \bm{\boxtimes} & \bm{\boxtimes} & \bm{\boxtimes}\\
0 & \bm{0} & \bm{\boxtimes} & \bm{\boxtimes} & \bm{\boxtimes}\\
0 & \bm{0} & \bm{\boxtimes} & \bm{\boxtimes} & \bm{\boxtimes}\\
0 & \bm{0} & \bm{\boxtimes} & \bm{\boxtimes} & \bm{\boxtimes}\\
0 & \bm{0} & \bm{\boxtimes} & \bm{\boxtimes} & \bm{\boxtimes}\\
\end{sbmatrix}.
\end{aligned}
$$

\subsubsection*{\textbf{Second Step 2.2: Introduce Zeros for the Second Row}}	
Similar to step (1.2), we now consider the  \textit{transpose} of $\bH_2\bH_1\bA\bL_1^\top$, given by $(\bH_2\bH_1\bA\bL_1^\top)^\top =\bL_1\bA^\top\bH_1^\top\bH_2^\top \in \real^{n\times m}$. 
Let the column partition of this transposed matrix be given by $\bL_1\bA^\top\bH_1^\top\bH_2^\top \triangleq [\bx_1, \bx_2, \ldots, \bx_m]$, where each $\bx_i \in \real^n$. 
Suppose $\bar{\bx}_1, \bar{\bx}_2, \ldots, \bar{\bx}_m \in \real^{n-2}$ are vectors obtained by removing the first two components from $\bx_i$'s. 
We can now proceed to construct the Householder reflector as follows:
$$
r_1 = \normtwo{\bar{\bx}_1},\qquad \bv_2 = \frac{\bar{\bx}_1 - r_1 \be_1}{\normtwo{\bar{\bx}_1 - r_1 \be_1}}, \qquad \text{and}\qquad \widetilde{\bL}_2 = \bI - 2\bv_2\bv_2^\top \in \textcolor{mylightbluetext}{\real^{(n-2)\times (n-2)}},
$$
where $\be_1$ now represents the first unit basis vector in $\textcolor{mylightbluetext}{\real^{n-2}}$, i.e., $\be_1=[1;0;0;\ldots;0]\in \textcolor{mylightbluetext}{\real^{n-2}}$. To introduce zeros below the subdiagonal and operate on the submatrix $(\bL_1\bA^\top\bH_1\bH_2)[3:n,1:m]$, we append the Householder reflector into $\bL_2$, defined as
$
\bL_2 = \footnotesize\begin{bmatrix}
\bI_2 &\bzero \\
\bzero & \widetilde{\bL}_2
\end{bmatrix},
$
where $\bI_2$ is the $2\times 2$ identity matrix. 
By definition, $\bL_2$ and $\widetilde{\bL}_2$ are both symmetric and orthogonal matrices (from the properties of Householder reflectors).
In this case, $\bL_2(\bL_1\bA^\top\bH_1^\top\bH_2^\top)$ will introduce zeros in the second column of $(\bL_1\bA^\top\bH_1^\top\bH_2^\top)$ below entry (3,2). 
Notably, the first two rows of $(\bL_1\bA^\top\bH_1^\top\bH_2^\top)$ remain unaffected and kept unchanged, as demonstrated by Remark~\ref{remark:left-right-identity}. 
{Furthermore, the first column of the matrix also remains unaltered}. 

Returning to the original \textit{untransposed} matrix $\bH_2\bH_1\bA\bL_1^\top$, multiplying on the right by $\bL_2^\top$ serves the purpose of introducing zeros in the second row to the right of entry (2,3).  

Following the above example, the procedure on the $7\times 5$ matrix is shown as follows, where $\boxtimes$ represents a value that is not necessarily zero, and \textbf{boldface} indicates the value has just been changed:
$$
\begin{aligned}
	\footnotesize
\begin{sbmatrix}{\bL_1\bA^\top\bH_1^\top\bH_2^\top }
\boxtimes & 0 & 0 & 0 & 0 & 0 & 0 \\
\boxtimes & \boxtimes & 0 & 0 & 0 & 0 & 0\\
0 & \boxtimes & \boxtimes & \boxtimes & \boxtimes & \boxtimes & \boxtimes\\
0 & \boxtimes & \boxtimes & \boxtimes & \boxtimes & \boxtimes & \boxtimes\\
0 & \boxtimes & \boxtimes & \boxtimes & \boxtimes & \boxtimes & \boxtimes\\
\end{sbmatrix}
\stackrel{\bL_2\times }{\rightarrow}
\footnotesize
\begin{sbmatrix}{\bL_2\bL_1\bA^\top\bH_1^\top\bH_2^\top }
\boxtimes & 0 & 0 & 0 & 0 & 0 & 0 \\
\boxtimes & \boxtimes & 0 & 0 & 0 & 0 & 0\\
0 & \bm{\boxtimes} & \bm{\boxtimes} & \bm{\boxtimes} & \bm{\boxtimes} & \bm{\boxtimes} & \bm{\boxtimes}\\
0 & \bm{0} & \bm{\boxtimes} & \bm{\boxtimes} & \bm{\boxtimes} & \bm{\boxtimes} & \bm{\boxtimes}\\
0 & \bm{0} & \bm{\boxtimes} & \bm{\boxtimes} & \bm{\boxtimes} & \bm{\boxtimes} & \bm{\boxtimes}\\
\end{sbmatrix}
\stackrel{(\cdot)^\top }{\rightarrow}
\footnotesize
\begin{sbmatrix}{\bH_2\bH_1\bA\bL_1^\top\bL_2^\top}
\boxtimes & \boxtimes & 0 & 0 & 0\\
0 & \boxtimes & \bm{\boxtimes} & \bm{0} & \bm{0} \\
0 & 0 & \bm{\boxtimes} & \bm{\boxtimes} & \bm{\boxtimes}\\
0 & 0 & \bm{\boxtimes} & \bm{\boxtimes} & \bm{\boxtimes}\\
0 & 0 & \bm{\boxtimes} & \bm{\boxtimes} & \bm{\boxtimes}\\
0 & 0 & \bm{\boxtimes} & \bm{\boxtimes} & \bm{\boxtimes}\\
0 & 0 & \bm{\boxtimes} & \bm{\boxtimes} & \bm{\boxtimes}
\end{sbmatrix}.\\
\end{aligned}
$$

\index{Lanczos-Golub-Kahan (LGK) bidiagonalization}
\index{Lanczos bidiagonalization}
\index{Golub-Kahan bidiagonalization}
\index{Bidiagonalization}

In short, the operation $\bH_2(\cdot)\bL_2^\top$ (specifically, the operation is $\bH_2(\bH_1\bA\bL_1^\top)\bL_2^\top$) accomplishes  the second step of introducing zeros into the second column and the second row of $\bA$.

The same process can continue, and we shall notice that there are $n$ such $\bH_i$ Householder reflectors on the left and $n-2$ such $\bL_i$ Householder reflectors on the right (assuming $m>n$ for simplicity; the number of reflectors will be clear from the complete example in the sequel). 
This interleaved Householder factorization is commonly referred to as the \textit{Golub-Kahan bidiagonalization} \citep{golub1965calculating} \footnote{The Golub-Kahan bidiagonalization is also known as \textit{Lanczos-Golub-Kahan (LGK) bidiagonalization} or simply \textit{Lanczos bidiagonalization}. }. Eventually, we achieve the  bidiagonalization by 
$$
\bB = \underbrace{\bH_{n} \bH_{n-1}\ldots\bH_1}_{\text{$n$ stages}} \bA\underbrace{\bL_1^\top\bL_2^\top\ldots\bL_{n-2}^\top}_{\text{$n-2$ stages}}.
$$
Since the $n$ left Householder reflectors $\bH_i$'s and $n-2$ right Householder reflectors $\bL_i$'s are symmetric and orthogonal, we can express the bidiagonalization as
$$
\bB =\bH_{n} \bH_{n-1}\ldots\bH_1 \bA\bL_1\bL_2\ldots\bL_{n-2}.
$$
Let $\bU\triangleq\bH_1\bH_2\ldots\bH_n$ and $\bV = \bL_1\bL_2\ldots\bL_{n-2}$. The bidiagonal decomposition of $\bA$ can then be expressed as $\bA \triangleq \bU\bB\bV^\top$. Since the upper-left portion of each $\bL_i$ contains an $i\times i$ identity matrix, the matrix $\bV$ exhibits the  structure
$
\bV = \footnotesize
\begin{bmatrix}
1 & \bzero \\
\bzero & \bQ
\end{bmatrix},
$
where $\bQ$ is an orthogonal matrix.

A complete example of this Golub-Kahan bidiagonalization for a $7\times 5$ matrix is presented as follows, where again $\boxtimes$ represents a value that is not necessarily zero, and \textbf{boldface} indicates the value has just been changed:

\begin{mdframed}[hidealllines=\mdframehidelineNote,backgroundcolor=\mdframecolor,frametitle={A Complete Example of Golub-Kahan Bidiagonalization}]
$$
\begin{aligned}
	\footnotesize
\begin{sbmatrix}{\bA}
\boxtimes & \boxtimes & \boxtimes & \boxtimes & \boxtimes \\
\boxtimes & \boxtimes & \boxtimes & \boxtimes & \boxtimes\\
\boxtimes & \boxtimes & \boxtimes & \boxtimes & \boxtimes\\
\boxtimes & \boxtimes & \boxtimes & \boxtimes & \boxtimes\\
\boxtimes & \boxtimes & \boxtimes & \boxtimes & \boxtimes\\
\boxtimes & \boxtimes & \boxtimes & \boxtimes & \boxtimes\\
\boxtimes & \boxtimes & \boxtimes & \boxtimes & \boxtimes
\end{sbmatrix}
&\stackrel{\bH_1\times}{\rightarrow}
\footnotesize
\begin{sbmatrix}{\bH_1\bA}
\bm{\boxtimes} & \bm{\boxtimes} & \bm{\boxtimes} & \bm{\boxtimes} & \bm{\boxtimes}\\
\bm{0} & \bm{\boxtimes} & \bm{\boxtimes} & \bm{\boxtimes} & \bm{\boxtimes}\\
\bm{0} & \bm{\boxtimes} & \bm{\boxtimes} & \bm{\boxtimes} & \bm{\boxtimes}\\
\bm{0} & \bm{\boxtimes} & \bm{\boxtimes} & \bm{\boxtimes} & \bm{\boxtimes}\\
\bm{0} & \bm{\boxtimes} & \bm{\boxtimes} & \bm{\boxtimes} & \bm{\boxtimes}\\
\bm{0} & \bm{\boxtimes} & \bm{\boxtimes} & \bm{\boxtimes} & \bm{\boxtimes}\\
\bm{0} & \bm{\boxtimes} & \bm{\boxtimes} & \bm{\boxtimes} & \bm{\boxtimes}\\
\end{sbmatrix}
\stackrel{\times\bL_1^\top}{\rightarrow}
\footnotesize
\begin{sbmatrix}{\bH_1\bA\bL_1^\top}
\boxtimes & \bm{\boxtimes} & \bm{0} & \bm{0} & \bm{0} \\
0 & \bm{\boxtimes} & \bm{\boxtimes} & \bm{\boxtimes} & \bm{\boxtimes}\\
0 & \bm{\boxtimes} & \bm{\boxtimes} & \bm{\boxtimes} & \bm{\boxtimes}\\
0 & \bm{\boxtimes} & \bm{\boxtimes} & \bm{\boxtimes} & \bm{\boxtimes}\\
0 & \bm{\boxtimes} & \bm{\boxtimes} & \bm{\boxtimes} & \bm{\boxtimes}\\
0 & \bm{\boxtimes} & \bm{\boxtimes} & \bm{\boxtimes} & \bm{\boxtimes}\\
0 & \bm{\boxtimes} & \bm{\boxtimes} & \bm{\boxtimes} & \bm{\boxtimes}\\
\end{sbmatrix}\\
&\stackrel{\bH_2\times}{\rightarrow}
\footnotesize
\begin{sbmatrix}{\bH_2\bH_1\bA\bL_1^\top}
\boxtimes & \boxtimes & \boxtimes & \boxtimes & \boxtimes \\
0 & \bm{\boxtimes} & \bm{\boxtimes} & \bm{\boxtimes} & \bm{\boxtimes}\\
0 & \bm{0} & \bm{\boxtimes} & \bm{\boxtimes} & \bm{\boxtimes}\\
0 & \bm{0} & \bm{\boxtimes} & \bm{\boxtimes} & \bm{\boxtimes}\\
0 & \bm{0} & \bm{\boxtimes} & \bm{\boxtimes} & \bm{\boxtimes}\\
0 & \bm{0} & \bm{\boxtimes} & \bm{\boxtimes} & \bm{\boxtimes}\\
0 & \bm{0} & \bm{\boxtimes} & \bm{\boxtimes} & \bm{\boxtimes}
\end{sbmatrix}
\stackrel{\times\bL_2^\top}{\rightarrow}
\footnotesize
\begin{sbmatrix}{\bH_2\bH_1\bA\bL_1^\top\bL_2^\top}
\boxtimes & \boxtimes & 0 & 0 & 0\\
0 & \boxtimes & \bm{\boxtimes} & \bm{0} & \bm{0} \\
0 & 0 & \bm{\boxtimes} & \bm{\boxtimes} & \bm{\boxtimes}\\
0 & 0 & \bm{\boxtimes} & \bm{\boxtimes} & \bm{\boxtimes}\\
0 & 0 & \bm{\boxtimes} & \bm{\boxtimes} & \bm{\boxtimes}\\
0 & 0 & \bm{\boxtimes} & \bm{\boxtimes} & \bm{\boxtimes}\\
0 & 0 & \bm{\boxtimes} & \bm{\boxtimes} & \bm{\boxtimes}
\end{sbmatrix}\\
&\stackrel{\bH_3\times}{\rightarrow}
\footnotesize
\begin{sbmatrix}{\bH_3\bH_2\bH_1\bA\bL_1^\top\bL_2^\top}
\boxtimes & \boxtimes & \boxtimes & \boxtimes & \boxtimes \\
0 & \boxtimes & \boxtimes & \boxtimes & \boxtimes \\
0 & 0 & \bm{\boxtimes} & \bm{\boxtimes} & \bm{\boxtimes}\\
0 & 0 & \bm{0} & \bm{\boxtimes} & \bm{\boxtimes}\\
0 & 0 & \bm{0} & \bm{\boxtimes} & \bm{\boxtimes}\\
0 & 0 & \bm{0} & \bm{\boxtimes} & \bm{\boxtimes}\\
0 & 0 & \bm{0} & \bm{\boxtimes} & \bm{\boxtimes}
\end{sbmatrix}
\stackrel{\times\bL_3^\top}{\rightarrow}
\footnotesize
\begin{sbmatrix}{\bH_3\bH_2\bH_1\bA\bL_1^\top\bL_2^\top\bL_3^\top}
\boxtimes & \boxtimes & 0 & 0 & 0\\
\boxtimes & \boxtimes & \boxtimes & 0 & 0\\
0 & \boxtimes & \boxtimes & \bm{\boxtimes} & \bm{0}\\
0 & 0 & \boxtimes & \bm{\boxtimes} & \bm{\boxtimes}\\
0 & 0 & 0 & \bm{\boxtimes} & \bm{\boxtimes}\\
0 & 0 & 0 & \bm{\boxtimes} & \bm{\boxtimes}\\
0 & 0 & 0 & \bm{\boxtimes} & \bm{\boxtimes}
\end{sbmatrix}\\
&\stackrel{\bH_4\times}{\rightarrow}
\footnotesize
\begin{sbmatrix}{{\scriptsize\bH_4\bH_3\bH_2\bH_1\bA\bL_1^\top\bL_2\bL_3^\top}}
\boxtimes & \boxtimes & 0 & 0 & 0 \\
0 & \boxtimes & \boxtimes & 0 & 0 \\
0 & 0 & \boxtimes& \boxtimes & 0\\
0 & 0 & 0 & \bm{\boxtimes} & \bm{\boxtimes}\\
0 & 0 & 0 & \bm{0} & \bm{\boxtimes}\\
0 & 0 & 0 & \bm{0} & \bm{\boxtimes}\\
0 & 0 & 0 & \bm{0} & \bm{\boxtimes}
\end{sbmatrix}
\stackrel{\bH_5\times}{\rightarrow}
\footnotesize
\begin{sbmatrix}{\bH_5\bH_4\bH_3\bH_2\bH_1\bA\bL_1^\top\bL_2\bL_3^\top}
\boxtimes & \boxtimes & 0 & 0 & 0 \\
0 & \boxtimes & \boxtimes & 0 & 0 \\
0 & 0 & \boxtimes& \boxtimes & 0\\
0 & 0 & 0 & \boxtimes & \boxtimes\\
0 & 0 & 0 & 0 & \bm{\boxtimes}\\
0 & 0 & 0 & 0 & \bm{0}\\
0 & 0 & 0 & 0 & \bm{0}
\end{sbmatrix}.
\end{aligned}
$$
\end{mdframed}
We present the method where a right Householder reflector $\bL_i$ follows from a left one $\bH_i$. However, a common mistake that might be employed is to perform the left reflectors altogether, followed by all the right ones. 
In this approach, a bidiagonal decomposition combines elements of both QR decomposition and Hessenberg decomposition.
Nevertheless, this is problematic because the right Householder reflector $\bL_1$ would destroy the zeros introduced by the left ones. Therefore, the left and right reflectors must be applied in an interleaved manner to properly introduce the zeros.

\section{Computing  Bidiagonal Decomposition: Golub-Kahan Approach}

Building upon the computation of QR decomposition and Hessenberg decomposition using the Householder reflector (Section~\ref{section:householder_qr_cp}; Section~\ref{section:householder_hessenger}), it is straightforward to derive the procedure to compute the bidiagonalization, as outlined in Algorithm~\ref{alg:bidiagonal-decomposition-householder}, where the \textcolor{winestain}{red}-colored text highlights the  differences between the right Householder reflectors
in bidiagonal decomposition and  those in Hessenberg decomposition (Algorithm~\ref{alg:hessenbert-decomposition-householder}).
\begin{theorem}[Algorithm Complexity: Golub-Kahan Bidiagonalization]\label{theorem:bidiagonal-full-householder}
Algorithm~\ref{alg:bidiagonal-decomposition-householder} requires $\sim 4mn^2-\frac{4}{3}n^3$ flops to compute a bidiagonal decomposition of an $m\times n$ matrix with $m>n$. 
Moreover, if $\bU$ and $\bV$ are needed explicitly, an additional $\sim 4m^2n-2mn^2 + 2n^3$ flops are required.
\end{theorem}
The proof is straightforward: the procedure shown above requires roughly twice the complexity of QR decomposition via the Householder reflector. This is because it involves two interleaved Householder QR decompositions: one operating on the  $m\times n$ matrix $\bA$ and the other on the $n\times m$ matrix $\bA^\top$. 
The computational complexity to obtain the orthogonal matrix $\bU$ is $4m^2n-2mn^2$ flops, which is the same as that in the QR decomposition (Theorem~\ref{theorem:qr-full-householder}). And similarly, the complexity to obtain the orthogonal matrix $\bV$ is $2n^3$ flops, which is the same as that in the Hessenberg decomposition (Theorem~\ref{theorem:hessenbert-householder}).

\begin{algorithm}[htp] 
\caption{Golub-Kahan Bidiagonal Decomposition} 
\label{alg:bidiagonal-decomposition-householder} 
\begin{algorithmic}[1] 
\Require Matrix $\bA$ with size $m\times n $ with $m>n$; 
\State Initially set $\bB = \bA^\top$; \Comment{Initially, $\bB\in\real^{n\times m}$}
\For{$i=1$ to $n$} 
\State \textcolor{winestain}{// do the left Householder reflector;}
\If{$i<=n-2$}
\State \textcolor{black}{$\bB \leftarrow \bB^\top$ such that $\bB\in \textcolor{mylightbluetext}{\real^{m\times n}}$;}
\EndIf
\State \textcolor{black}{$\ba \leftarrow \bB_{i:m,i}$, i.e., first column of $\bB_{i:m,i:n} \in \real^{(m-i+1)\times(n-i+1)}$;}
\State $r \leftarrow \normtwo{\ba}$; \Comment{$2(m-i+1)$ flops}
\State $\bu_i \leftarrow \ba-r\be_1 \in \real^{m-i+1}$; \Comment{1 flop}
\State $\bu_i \leftarrow \bu_i / \normtwo{\bu_i}$; \Comment{$3(m-i+1)$ flops}
\State \algoalign{$\bB_{i,i} \leftarrow r$, $\bB_{i+1:m,i}\leftarrow\bzero$; \Comment{update the first column of $\bB_{i:m,i:n}$, 0 flops}}
\State Update  second to  last columns  of $\bB_{i:m,i:n}$:\Comment{$4(m-i+1)(n-i) +(m-n+1)$ flops}
$$
\bB_{i:m,i+1:n} \leftarrow \bB_{i:m,i+1:n} - 2\bu_i (\bu_i^\top \bB_{i:m,i+1:n});
$$ 
\If{$i<=n-2$} 
\State \textcolor{winestain}{// do the right Householder reflector;}
\State $\bB \leftarrow \bB^\top$ such that $\bB\in \textcolor{mylightbluetext}{\real^{n\times m}}$;
\State $\bz \leftarrow \bB_{i+1:n,i}$, i.e., first column of $\bB_{i+1:n,i:\textcolor{winestain}{m}}\in \real^{(n-i)\times(\textcolor{winestain}{m}-i+1)}$;
\State $s \leftarrow \normtwo{\bz}$; \Comment{$2(n-i)$ flops}
\State $\bv_i \leftarrow \bz-s\be_1 \in \real^{n-i}$; \Comment{1 flop}
\State $\bv_i \leftarrow \bv_i / \normtwo{\bv_i}$; \Comment{$3(n-i)$ flops}
\State $\bB_{i+1,i}\leftarrow r$, $\bB_{i+2:n,i}\leftarrow\bzero$; \Comment{update the first column of $\bB_{i+1:n,i:\textcolor{winestain}{m}}$, 0 flops}
\State Update second to last columns of $\bB_{i+1:n,i:\textcolor{winestain}{m}}$: \Comment{$4(n-i)(m-i) +(n-m)$ flops}
$$
\begin{aligned}
\bB_{i+1:n,i+1:\textcolor{winestain}{m}} &\leftarrow (\bI-2\bv_i\bv_i^\top)\bB_{i+1:n,i+1:\textcolor{winestain}{m}}\\
&= \bB_{i+1:n,i+1:\textcolor{winestain}{m}} - 2\bv_i(\bv_i^\top\bB_{i+1:n,i+1:\textcolor{winestain}{m}}) \in \real^{(n-i)\times(\textcolor{winestain}{m}-i)}\\
\end{aligned}
$$
\EndIf
\EndFor
\State Output $\bB$ as the bidiagonal matrix;
\State Get $\bU=\bH_1 \bH_2\ldots\bH_n$ and $\bV = \bL_1\bL_2\ldots\bL_{n-2}$;
\State Initially set $\bU \leftarrow \bH_1$;
\For{$i=1$ to $n-1$} 
\State $\bU_{1:m,i+1:m} \leftarrow \bU_{1:m,i+1:m}(\bI - 2\bu_{i+1}\bu_{i+1}^\top)=\bU_{1:m,i+1:m}-\bU_{1:m,i+1:m}2\bu_{i+1}\bu_{i+1}^\top$.
\EndFor
\State Initially set $\bV \leftarrow \bL_1$;
\For{$i=1$ to $n-3$} 
\State $\bV_{1:n,i+2:n} \leftarrow \bV_{1:n,i+2:n}(\bI - 2\bv_{i+1}\bv_{i+1}^\top)=\bV_{1:n,i+2:n}-\bV_{1:n,i+2:n}2\bv_{i+1}\bv_{i+1}^\top \in \real^{n\times (n-i-1)}$;
\EndFor
\State Output $\bU$ and $\bV$ as the orthogonal matrix;
\end{algorithmic} 
\end{algorithm}

\section{Computing  Bidiagonal Decomposition: LHC Approach}
\begin{figure}[H]
\centering
\includegraphics[width=0.9\textwidth]{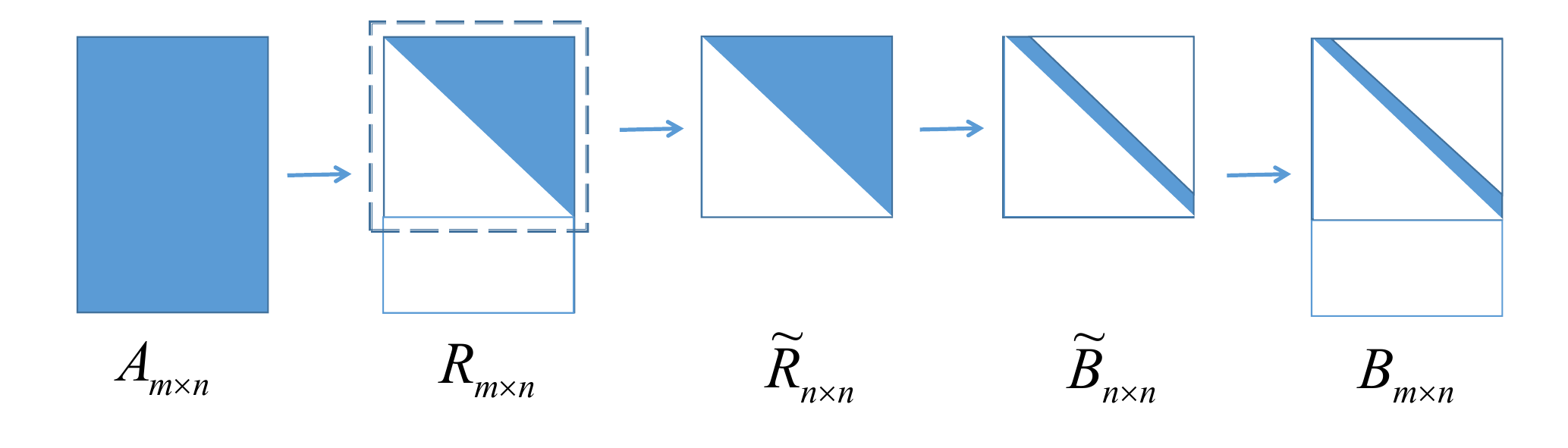}
\caption{Illustration of the LHC bidiagonalization process for  a matrix: first using QR decomposition, then using Golub-Kahan bidiagonalization.}
\label{fig:lhc-bidiagonal}
\end{figure}

In the previous section, we discussed that left   and right  Householder reflectors must be applied in an interleaved manner to avoid destroying the zeros introduced by the left reflectors. However, when $m\gg n$, we can optimize the process by extracting a square triangular matrix using QR decomposition and applying the Golub-Kahan bidiagonalization to the $n\times n$ square triangular matrix instead. 
This technique is known as the \textit{Lawson-Hanson-Chan (LHC) bidiagonalization} \citep{lawson1995solving, chan1982improved}, and the procedure is shown in Figure~\ref{fig:lhc-bidiagonal}.

\index{LHC bidiagonalization}

The LHC bidiagonalization process commences with the computation of the QR decomposition, $\bA = \bQ\bR$.  
It is then followed by applying the Golub-Kahan process to $\widetilde{\bR}$, where $\widetilde{\bR}=\bR[1:n,1:n]$ is the square $n\times n$ triangular submatrix inside $\bR$ (see Figure~\ref{fig:lhc-bidiagonal}), resulting in $\widetilde{\bR} = \widetilde{\bU} \widetilde{\bB} \bV^\top$, where $\widetilde{\bB}$ is bidiagonal, and $\widetilde{\bU}$ and $\bV$ are orthogonal. Subsequently,  $\widetilde{\bU}$ is appended to form an $m\times m$  orthogonal matrix, $\widetilde{\bB}$ is appended to form a thin matrix:
$$
\bU_0 = 
\begin{bmatrix}
\widetilde{\bU} & \bzero \\
\bzero & \bI_{m-n}
\end{bmatrix} \in\real^{m\times m}
\qquad \text{and}\qquad 
\bB = \begin{bmatrix}
\widetilde{\bB} \\ 
\bzero_{(m-n)\times n}	
\end{bmatrix}
\in\real^{m\times n},
$$
which results in $\bR=\bU_0\bB \bV^\top$ and $\bA = \bQ\bU_0\bB \bV^\top$. Let $\bU\triangleq\bQ\bU_0$, we obtain the bidiagonal decomposition. The QR decomposition requires $2mn^2-\frac{2}{3}n^3$ flops and the Golub-Kahan process now requires $\frac{8}{3}n^3$ (operating on an $n\times n$ submatrix). Thus, the total complexity to obtain the bidiagonal matrix $\bB$ is 
$$
\text{LHC bidiagonalization:   } \sim 2mn^2 + 2n^3 \text{  flops}.
$$
The LHC process creates zeros and subsequently  destroys them again in the lower triangle portion of the upper $n\times n $ square matrix of $\bR$ (i.e., zeros introduced to the lower triangle portion of $\widetilde{\bR}$ by the QR decomposition  will be destroyed by the second phase of the Golub-Kahan bidiagonalization).
However, the zeros in the lower $(m-n)\times n$ rectangular matrix of $\bR$ remain unaffected.
Therefore, when $m-n$ is sufficiently large (or $m\gg n$), there is a net gain in efficiency. Simple calculations show that  the LHC bidiagonalization costs less when $m>\frac{5}{3}n$ compared to the standard Golub-Kahan bidiagonalization.

\section{Computing  Bidiagonal Decomposition: Three-Step Approach}
\begin{figure}[H]
\centering
\includegraphics[width=0.9\textwidth]{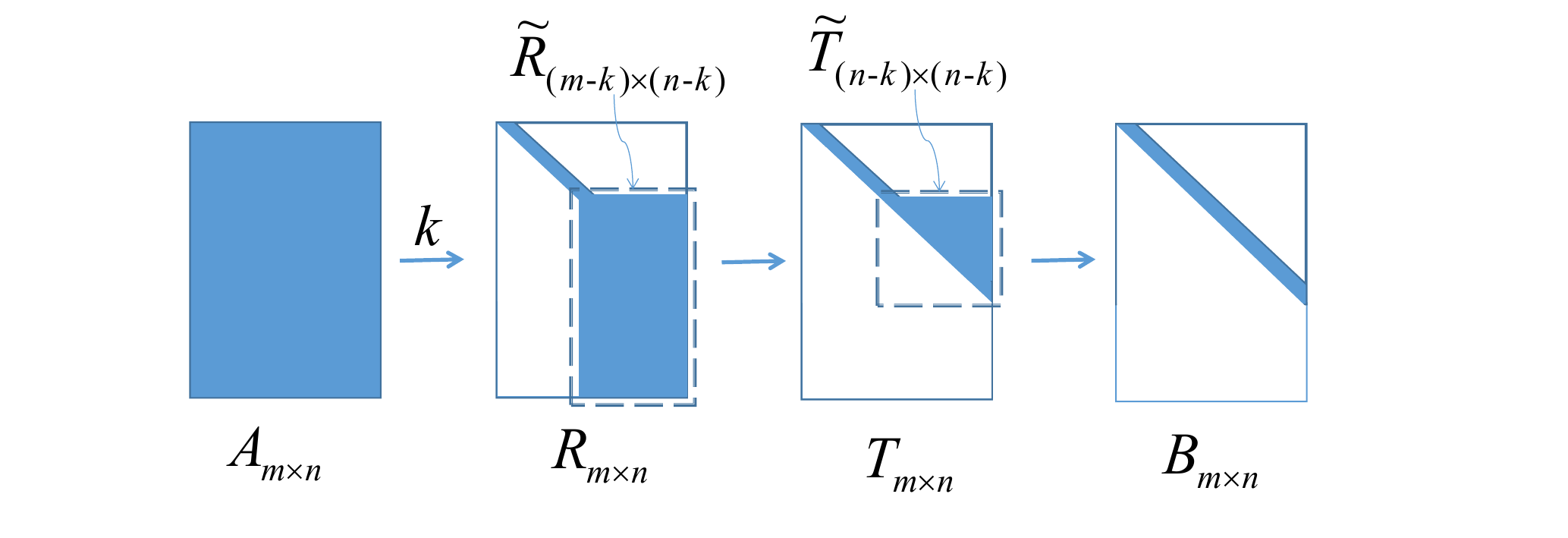}
\caption{Illustration of the Three-Step bidiagonalization process for a matrix: first applying the Golub-Kahan process for $k$ steps, and applying the LHC process for the remaining.}
\label{fig:lhc-bidiagonal2}
\end{figure}

The LHC procedure is advantageous only when $m>\frac{5}{3}n$. A further trick  involves applying the QR decomposition not at the beginning of the computation, but at a suitable point in the middle \citep{trefethen1997numerical}. 
Specifically, the procedure is outlined in Figure~\ref{fig:lhc-bidiagonal2}, where we apply the first $k$ steps of left and right Householder reflectors as those in the Golub-Kahan process, while leaving the lower-right $(m-k)\times(n-k)$ submatrix ``unreflected." 
Then, the LHC process is applied to this submatrix to obtain the final bidiagonal decomposition. This approach reduces complexity when  $n<m<2n$.

The complexity of the Three-Step bidiagonalization can be broken down into three parts. 
The complexity of the first $k$ loops in Algorithm~\ref{alg:bidiagonal-decomposition-householder} can be demonstrated as follows:
$$
\textbf{(Step 1)}: \qquad f_1 = 8mnk - (4m+4n)k^2 + \frac{8}{3}k^3 \text{  flops},
$$
which makes $4mn^2 - \frac{4}{3}n^3$ flops when $k=n$.
The complexity of the QR decomposition using the Householder reflector for $\widetilde{\bR} \in \real^{(m-k)\times (n-k)}$ is 
$$
\textbf{(Step 2)}:  \qquad f_2 = 2(m-k)(n-k)^2 - \frac{2}{3} (n-k)^3\text{  flops}.
$$
And the complexity of the Golub-Kahan bidiagonalization for $\widetilde{\bT} \in \real^{(n-k)\times (n-k)}$ is 
$$
\textbf{(Step 3)}: \qquad f_3 = \frac{8}{3}(n-k)^3\text{  flops}.
$$
Thus, the overall complexity of the three steps is given by 
$$
g(k) =f_1 +f_2+f_3 = -\frac{4}{3}k^3 + (6n-2m)k^2  +(4mn-8n^2)k + 2mn^2+2n^3  .
$$
The problem now becomes to finding a value of $k<n$ such that $g(k)$ is minimized. 
By taking the derivative of the function above, we can derive the following:
$$
g^\prime(k) = -4k^2 +(12n-4m)k + (4mn-8n^2),
$$
the root of which  is $k=n$ or $2n-m$. We notice that $0<k<n$ such that when $2n-m>0$, the optimal value appears in one of $\{0, n, 2n-m\}$. Trivial calculation shows that
the optimal value for $k$ is $k=2n-m$ if $0\leq 2n-m<n$ (or $k=n$ if $n\geq 2n-m$, which is equivalent to the Golub-Kahan process) and the final complexity is  reduced to 
$$
\begin{aligned}
g(2n-m) &=  2mn^2 + 2m^2n -\frac{2}{3}m^3 -\frac{2}{3}n^3 \,\, \text{ flops}, \qquad  \text{when $n\leq m<2n$}.
\end{aligned}
$$

\begin{figure}[H]
\centering
\vspace{-0.35cm}
\subfigtopskip=2pt
\subfigbottomskip=2pt
\subfigcapskip=-5pt
\subfigure[$g(k)$.]{\label{fig:bidiagonal-gk-sample1}
\includegraphics[width=0.475\linewidth]{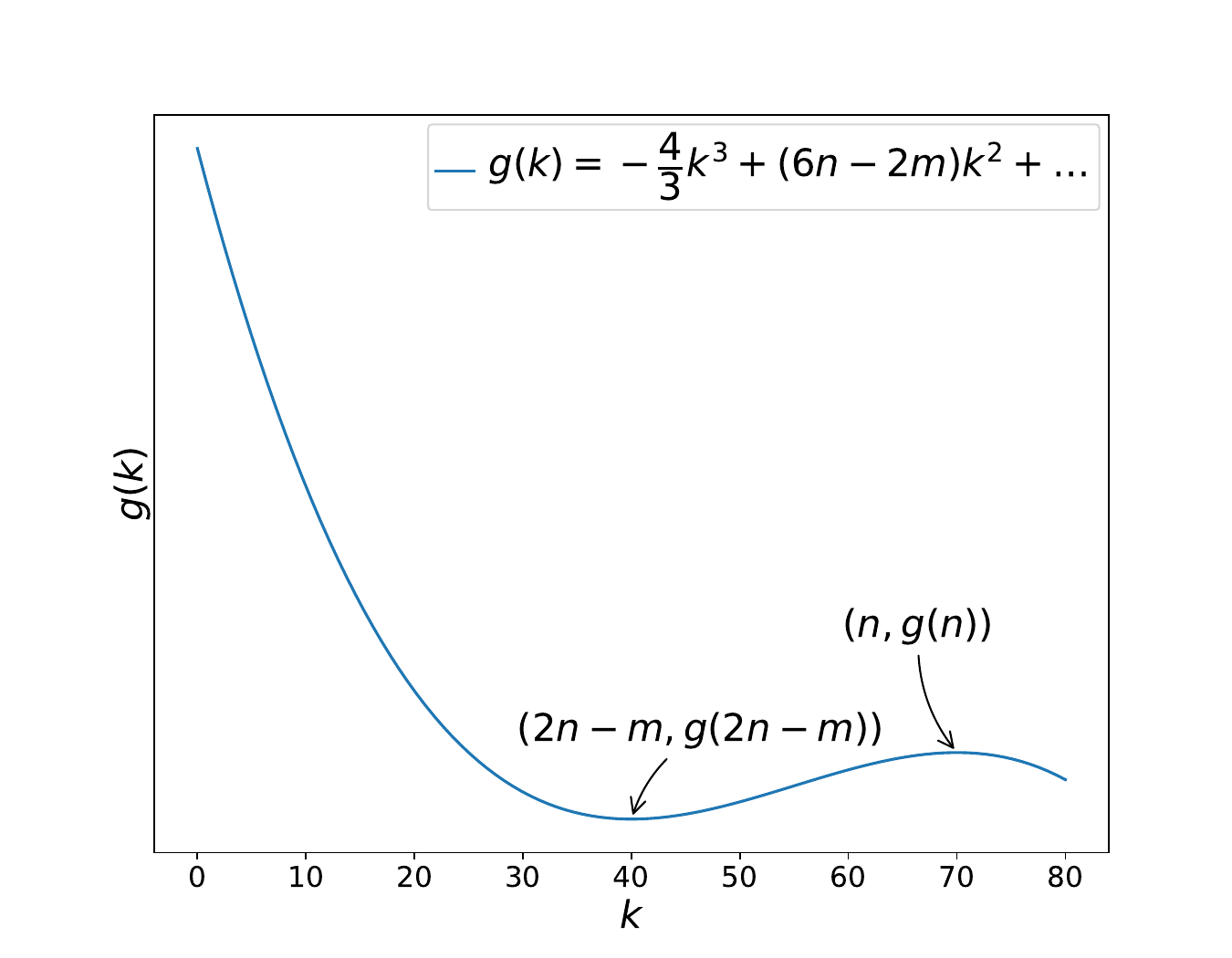}}
\quad 
\subfigure[$g^\prime(k)$.]{\label{fig:bidiagonal-gk-sample2}
\includegraphics[width=0.475\linewidth]{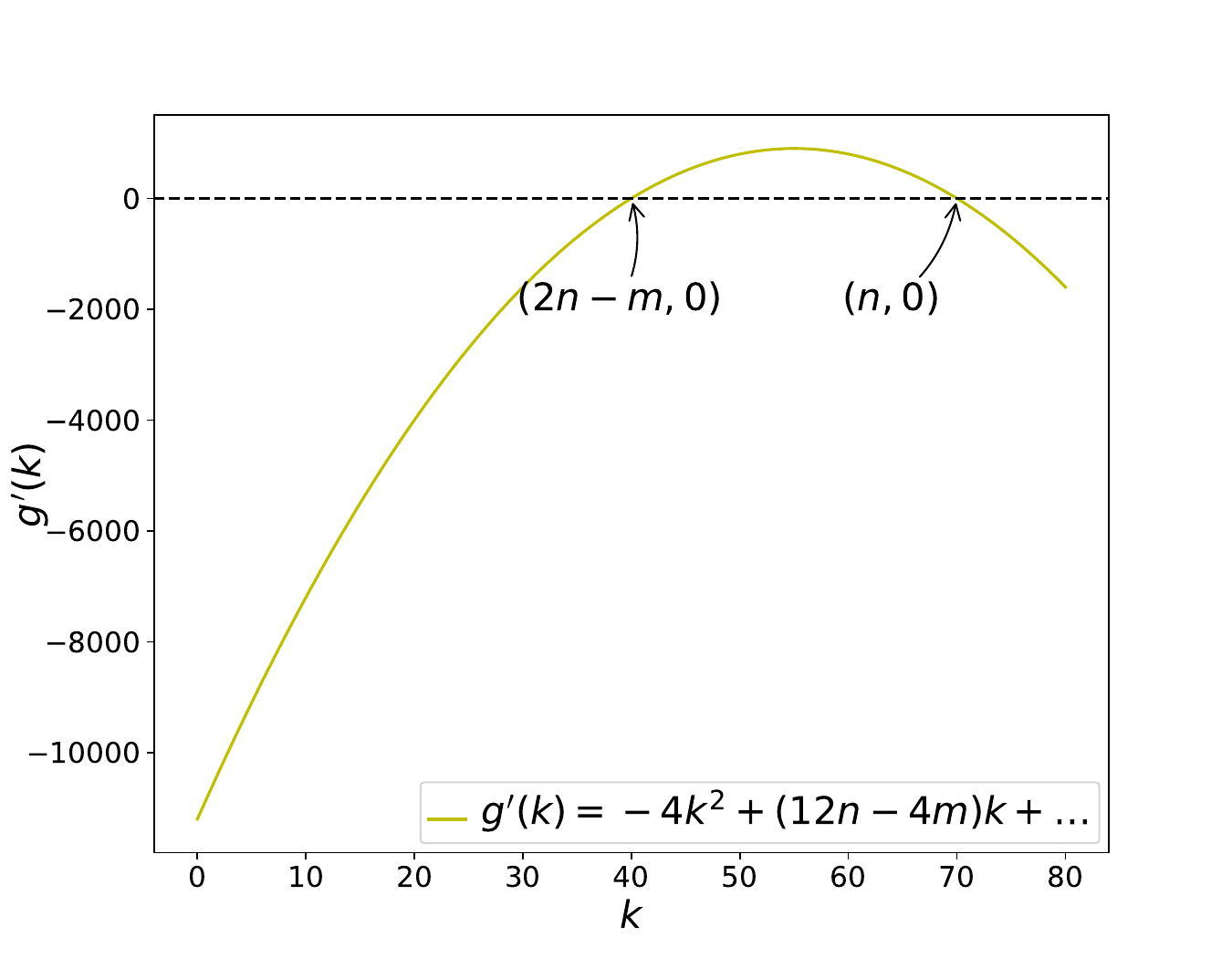}}
\caption{An example of the complexity when $n=70, m=100$.}
\label{fig:bidiagonal-gk-sample}
\end{figure}

An example illustrating the complexity of the Three-Step method when $n=70, m=100$ is shown in Figure~\ref{fig:bidiagonal-gk-sample}, where the roots of the derivative are found to be at $2n-m=40$ and $n=70$, such that $g^\prime(40)=g^\prime(70)=0$. In this specific case, the function $g(k)$ is decreasing when $k\in (0, 2n-m]$ and increasing when $k\in (2n-m, n]$.

\noindent

\begin{SCfigure}
\centering
\includegraphics[width=0.6\textwidth]{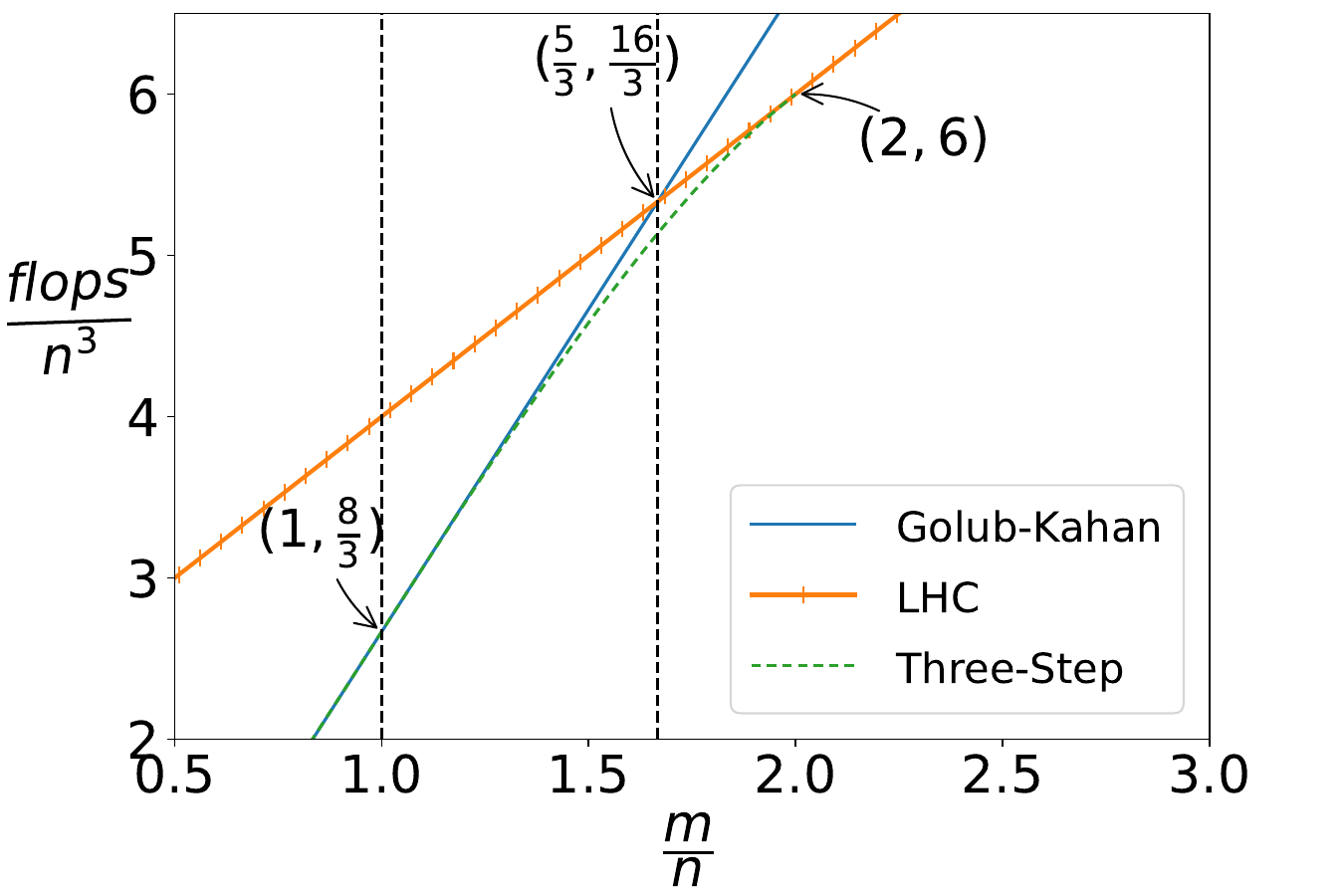}
\caption{Comparison of the computational complexities of the three bidiagonalization methods. When $m>2n$, the LHC method is preferred; when $n\leq m<2n$, the Three-Step method is preferred, though the improvement is relatively modest.}
\label{fig:bidiagonal-loss-compare}
\end{SCfigure}

To conclude, the costs of the three methods are shown as follows:
$$ 
\left\{
\begin{aligned}
&\text{Golub-Kahan: }  \sim 4mn^2-\frac{4}{3}n^3 \,\, \text{ flops},   \\
&\text{LHC: } \sim 2mn^2 + 2n^3  \,\, \text{ flops},  \\
&\text{Three-Step: } \sim 2mn^2 + 2m^2n -\frac{2}{3}m^3 -\frac{2}{3}n^3 \,\, \text{ flops}, \,\,\text{if $n\leq m<2n$}.
\end{aligned}
\right.
$$
When $m>2n$, the LHC method is the preferred choice for computing the bidiagonal decomposition; when $n\leq m<2n$, the Three-Step method is preferred, although the improvement is relatively modest, as shown in Figure~\ref{fig:bidiagonal-loss-compare}, where the operation counts for the three methods are plotted as a function of $\frac{m}{n}$.

Notice that the complexity discussed here does not involve the additional computations required for matrices $\bU$ and  $\bV$. We shall not discuss this issue for simplicity.

\section{Connection to Tridiagonal Decomposition}

We first illustrate the connection between  tridiagonal and bidiagonal decompositions with the following lemma, which  reveals how to construct a tridiagonal matrix from a bidiagonal one.
\begin{lemma}[Construct Tridiagonal From Bidiagonal]\label{lemma:construct-triangular-from-bidia}
Let $\bB\in \real^{n\times n}$ be upper bidiagonal. Then, $\bT_1=\bB^\top\bB$ and $\bT_2=\bB\bB^\top$ are \textit{symmetric} tridiagonal matrices.
\end{lemma}
\begin{proof}[of Lemma~\ref{lemma:construct-triangular-from-bidia}]
Suppose $\bB$ has the following form
$$
\bB=
\begin{bmatrix}
b_{11} & b_{12} & 0 & 0 &\ldots \\
0      & b_{22} & b_{23} & 0 & \ldots \\
0      & 0      & b_{33} & b_{34} & \ldots \\
\vdots & \vdots & \vdots & \ddots & \ddots \\
\ldots & \ldots & \ldots & \ldots & b_{nn} 
\end{bmatrix}.
$$
Then, $ \bT_1=\bB^\top \bB$ is computed as:
$$
\begin{aligned}
&\gap \bT_1=\bB^\top\bB =\\
&\begin{bmatrix}
b_{11} & 0 & 0 & 0 &\ldots \\
b_{12}     & b_{22} & 0  & 0 & \ldots \\
0      & b_{23}     & b_{33} & 0 & \ldots \\
0 & 0 & b_{34} & \vdots & \ddots \\
\vdots & \vdots & \vdots & \vdots & \vdots 
\end{bmatrix}
\begin{bmatrix}
b_{11} & b_{12} & 0 & 0 &\ldots \\
0      & b_{22} & b_{23} & 0 & \ldots \\
0      & 0      & b_{33} & b_{34} & \ldots \\
\vdots & \vdots & \vdots & \ddots & \ddots \\
\ldots & \ldots & \ldots & \ldots & b_{nn} 
\end{bmatrix}
=
\begin{bmatrix}
b_{11}^2      & b_{11}b_{12} & 0  &\ldots \\
b_{11}b_{12} &b_{12}^2+ b_{22}^2 & b_{22}b_{23} & \ldots \\
0      & b_{22}b_{23}      & b_{23}^2+b_{33}^2 &  \ldots \\
\vdots & \vdots & \ddots & \ddots  \\
\end{bmatrix},
\end{aligned}
$$
which is symmetric and tridiagonal as claimed. Similarly, we can demonstrate that $\bT_2=\bB\bB^\top$ is also symmetric and tridiagonal:
$$
\begin{aligned}
&\gap \bT_2=\bB\bB^\top =\\
&
\begin{bmatrix}
b_{11} & b_{12} & 0 & 0 &\ldots \\
0      & b_{22} & b_{23} & 0 & \ldots \\
0      & 0      & b_{33} & b_{34} & \ldots \\
\vdots & \vdots & \vdots & \ddots & \ddots \\
\ldots & \ldots & \ldots & \ldots & b_{nn} 
\end{bmatrix}
\begin{bmatrix}
b_{11} & 0 & 0 & 0 &\ldots \\
b_{12}     & b_{22} & 0  & 0 & \ldots \\
0      & b_{23}     & b_{33} & 0 & \ldots \\
0 & 0 & b_{34} & \vdots & \ddots \\
\vdots & \vdots & \vdots & \vdots & \vdots 
\end{bmatrix}
=
\begin{bmatrix}
b_{11}^2 +b_{12}^2     & b_{12}b_{22} & 0  &\ldots \\
b_{12}b_{22} &b_{22}^2+ b_{23}^2 & b_{23}b_{33} & \ldots \\
0      & b_{23}b_{33}      & b_{33}^2+b_{34}^2 &  \ldots \\
\vdots & \vdots & \ddots & \ddots  \\
\end{bmatrix}.
\end{aligned}
$$
\end{proof}
The lemma above reveals an important property. Suppose $\bA=\bU\bB\bV^\top$ is the bidiagonal decomposition of $\bA$, then the symmetric matrix $\bA\bA^\top$ admits the following tridiagonal decomposition:
$$
\bA\bA^\top=\bU\bB\bV^\top \bV\bB^\top\bU^\top = \bU\bB\bB^\top\bU^\top.
$$
Similarly, the symmetric matrix $\bA^\top\bA$ has the following tridiagonal decomposition:
$$
\bA^\top\bA=\bV\bB^\top\bU^\top \bU\bB\bV^\top=\bV\bB^\top\bB\bV^\top.
$$ 

As a final result in this section, we state a theorem giving the tridiagonal decomposition of a symmetric matrix with special eigenvalues. 
\begin{theoremHigh}[Tridiagonal Decomposition for Nonnegative Eigenvalues]\label{theorem:tri-nonnegative-eigen}
Let $\bA$ be   an $n\times n$ symmetric matrix with nonnegative eigenvalues. Then, there exists a matrix $\bZ$ such that 
$$
\bA=\bZ\bZ^\top.
$$
Moreover, the tridiagonal decomposition of $\bA$ can be reduced to a problem of finding the bidiagonal decomposition of $\bZ =\bU\bB\bV^\top$ such that the tridiagonal decomposition of $\bA$ is given by 
$$
\bA = \bZ\bZ^\top = \bU\bB\bB^\top\bU^\top.
$$
\end{theoremHigh}
\begin{proof}[of Theorem~\ref{theorem:tri-nonnegative-eigen}]
The eigenvectors of symmetric matrices can be chosen to be orthogonal (Proposition~\ref{proposition:orthogonal-eigenvectors}) such that the symmetric matrix $\bA$ can be decomposed as $\bA=\bQ\bLambda\bQ^\top$ (spectral theorem~\ref{theorem:spectral_theorem}), where $\bLambda$ is a diagonal matrix containing the eigenvalues of $\bA$. When the eigenvalues are nonnegative, $\bLambda$ can be factored as $\bLambda=\bLambda^{1/2} \bLambda^{1/2}$. Let $\bZ \triangleq \bQ\bLambda^{1/2}$. Then, $\bA$ can be factored as $\bA=\bZ\bZ^\top$. Combining these findings yields the desired result.
\end{proof}

\index{Least squares}
\section{Bidiagonal Least Squares and LGK Bidiagonalization}
We consider an  overdetermined linear system $\bC\bx=\bb$, where $\bC\in\real^{m\times (n-1)}$ and $m> (n-1)$. 
\footnote{Note we set the matrix dimensions to $m\times (n-1)$ since we consider the bidiagonal decomposition of an $m\times n$ augmented matrix. Generally, we can also consider $\bC\in\real^{m\times n}$ with $m\geq n$.}
We are interested in the bidiagonalization of the augmented matrix $\bA \triangleq [\bb, \bC]\in\real^{m\times n}$, which admits the following bidiagonal decomposition:
\begin{equation}\label{equation:lgk_bidi_ls}
\begin{aligned}
\bA = &\bU\bB\bV^\top 
\triangleq
\bU\bB
\begin{bmatrix}
1 & \bzero \\
\bzero & \bQ
\end{bmatrix}^\top\\ 
&\implies  
\bB=
\bU^\top [\bb, \bC] 
\begin{bmatrix}
1 & \bzero \\
\bzero & \bQ
\end{bmatrix}
=
[\bU^\top\bb, \bU^\top\bC\bQ]
=
\begin{bmatrix}
b_{11}\be_1& \bB_2\\
\bzero & \bzero
\end{bmatrix},
\end{aligned}
\end{equation}
where $b_{11}$ represents the (1,1) entry of matrix $\bB$, $\bB_2 \triangleq \bB[1:n,2:n]\in\real^{n\times (n-1)}$, and $\bQ\triangleq\bV[2:n,2:n]\in\real^{(n-1)\times (n-1)}$ is orthogonal.
We then have 
$$
\begin{aligned}
\normtwo{\bb-\bC\bx} =
\normtwo{
[\bb, \bC]
\begin{bmatrix}
1 \\
-\bx 
\end{bmatrix}
}
=
\normtwo{
\bU^\top
[\bb, \bC]\bV\bV^\top
\begin{bmatrix}
1 \\
-\bx 
\end{bmatrix}
}.
\end{aligned}
$$
Let $\bd \triangleq \bQ^\top\bx$, we have 
$$
\normtwo{\bb-\bC\bx}
=
\normtwo{\bB \bV^\top
\begin{bmatrix}
1 \\
-\bx 
\end{bmatrix}
}
=
\normtwo{b_{11} \be_1- \bB_2\bd}.
$$
The least squares problem of $\normtwo{\bb-\bC\bx}$ then can be equivalently recovered by finding the least squares solution of $\normtwo{b_{11} \be_1- \bB_2\bd}$ in terms of the variable $\bd$.

\subsection*{LGK Bidiagonalization}

Express $\bB$ as follows:
$$
\bB = 
\footnotesize
\begin{bmatrix}
b_{11} & b_{12} &   &  \ldots &  \\
 & b_{22} & b_{23} & \ldots  &   \\
&&\ddots &\ddots & \vdots \\
&&&b_{n-1,n-1}& b_{n-1,n} \\
&&&& b_{nn}\\
\bzero&\bzero&\bzero&\bzero& \bzero
\end{bmatrix}\in\real^{m\times n}.
$$
We have 
$$
\bB_2 = 
\footnotesize
\begin{bmatrix}
b_{12} &  &  \ldots &   \\
b_{22} & b_{23} & \ldots  &  \\
&\ddots &\ddots & \vdots \\
&&b_{n-1,n-1}& b_{n-1,n} \\
&&& b_{nn}
\end{bmatrix}
\normalsize
\bB_2^\top = 
\footnotesize
\begin{bmatrix}
	\footnotesize
b_{12} & b_{22} &  \ldots &   &  \\
0 & b_{23} & b_{33}  &   &  \\
& &\ddots &b_{n-1,n-1} &  \\
&& & b_{n-1,n} & b_{nn} \\
\end{bmatrix}
\in\real^{(n-1)\times n}.
$$
\paragraph{First step.}
From Equation~\eqref{equation:lgk_bidi_ls}, we find that 
${b_{11}  = \normtwo{\bb} \triangleq \normtwo{\ba_1}}$, where $\bb\triangleq\ba_1$ is the first column of $\bA$. 
Additionally, we have:
$$
\bC\bQ = 
\bU
\begin{bmatrix}
	\bB_2\\
	\bzero 
\end{bmatrix}
\implies
\bC^\top 
\begin{bmatrix}
\bu_1 & \bu_2&\ldots &\bu_{n}
\end{bmatrix}
=
\bQ\bB_2^\top, 
\quad \text{where}\quad \bU=[\bu_1,\bu_2,\ldots,\bu_m].
$$
Let $\bQ\triangleq[\bq_1,\bq_2,\ldots,\bq_{n-1}]$ be the column partition of $\bQ$, and let $\bq_0\triangleq\bzero$. We then have 
$$
\left\{
\begin{aligned}
\bC^\top\bu_i &= b_{ii}\bq_{i-1} + b_{i,i+1}\bq_i \,\,\implies\,\, b_{i,i+1}\bq_i = \bC^\top\bu_i - b_{ii}\bq_{i-1}, \gap \forall i\in\{1,2,\ldots, n-1\};\\
\bC^\top\bu_n &= b_{nn}\bq_{n-1}.
\end{aligned}
\right.
$$
If $b_{ii}$ and $\bq_{i-1}$ are known, $b_{i,i+1}$ can be determined as the norm of the right-hand side equation in the above equality: 
\begin{equation}\label{equation:lgk_bidi_ls_ste11}
{b_{i,i+1}=\pm\normtwo{\bC^\top\bu_i - b_{ii}\bq_{i-1}}}, \gap \forall i\in\{1,2,\ldots, n-1\}.
\end{equation}
and 
\begin{equation}\label{equation:lgk_bidi_ls_ste12}
{\bq_i = \frac{\bC^\top\bu_i - b_{ii}\bq_{i-1}}{b_{i,i+1}}}, \gap \text{if }b_{i,i+1}\neq 0, \,\,\forall i\in\{1,2\ldots, n-1\}.
\end{equation}
and 
\begin{equation}\label{equation:lgk_bidi_ls_ste13}
{b_{ii}=\pm\normtwo{\bC^\top\bu_i - b_{i,i+1}\bq_i}}, \gap \forall i\in\{2,3,\ldots, n-1\}.
\end{equation}

\paragraph{Second step.}
Similarly, from  Equation~\eqref{equation:lgk_bidi_ls}, we have ${\bu_1 = \bb/b_{11}=\ba_1/b_{11}}$ and 
$
\bC\bQ =
\footnotesize 
\bU
\begin{bmatrix}
\bB_2\\
\bzero 
\end{bmatrix}.
$
This leads to:
$$
\bC \bq_i = b_{i,i+1}\bu_i + b_{i+1,i+1}\bu_{i+1}, \gap \forall i \in\{1,2,\ldots,n-1\}.
$$
This implies
\begin{equation}\label{equation:lgk_bidi_ls_eq2}
{\bu_{i+1} = \frac{\bC \bq_i -b_{i,i+1}\bu_i}{b_{i+1,i+1}}} , \gap \text{if }b_{i+1,i+1}\neq 0,\,\, \forall i \in\{1,2,\ldots,n-1\}.
\end{equation}

The two steps described above form a recursive algorithm for computing the bidiagonal decomposition of the matrix  $\bA$, and is known as the \textit{LGK bidiagonalization}. 
The derivation above is valid  when $m\geq n$. A similar approach can be applied when $n\geq m$. 
Simple calculations can show the complexity is $\sim 4mn^2$ flops to obtain all $\bB, \bU$, and $\bV$, which  reduces the complexity of Theorem~\ref{theorem:bidiagonal-full-householder}.

\begin{algorithm}[H] 
\caption{LGK Bidiagonal Decomposition} 
\label{alg:lgk_bidiagonal} 
\begin{algorithmic}[1] 
\Require Matrix $\bA$ with size $m\times n $ and $m\geq n$; 
\State Initially set $b_{11}  \leftarrow  \normtwo{\ba_1}$, $\bu_1 \leftarrow \ba_1/b_{11}$, $\bq_0\leftarrow\bzero$; \Comment{$3m$ flops}
\For{$i=1$ to $n-1$} 
\State $b_{i,i+1}\leftarrow \pm\normtwo{\bC^\top\bu_i - b_{ii}\bq_{i-1}}$ by Equation~\eqref{equation:lgk_bidi_ls_ste11}; \Comment{$(2m+3)(n-1)$ flops}
\State $\bq_i \leftarrow \frac{\bC^\top\bu_i - b_{ii}\bq_{i-1}}{b_{i,i+1}}$ by 
Equation~\eqref{equation:lgk_bidi_ls_ste12}; \Comment{$n-1$ flops}
\State $\bu_{i+1} \leftarrow \frac{\bC \bq_i -b_{i,i+1}\bu_i}{b_{i+1,i+1}}$ by 
Equation~\eqref{equation:lgk_bidi_ls_eq2}; \Comment{$(2n-1)m$ flops}
\State $b_{jj}\leftarrow\pm\normtwo{\bC^\top\bu_j - b_{j,j+1}\bq_j}$ by 
Equation~\eqref{equation:lgk_bidi_ls_ste13}, where $j=i+1$; \Comment{$4(n-1)$ flops}
\EndFor
\State Output $\bB, \bU$, and $\bV$;
\end{algorithmic} 
\end{algorithm}

The algorithm breaks down if any $b_{i,i+1}$ or $b_{jj}$ is equal to zero. It can be shown in the context of solving least squares problems, these zero divisions can be handled with special cases \citep{bjorck2004acta}.
Another issue arises is that, in floating-point arithmetic, the columns in $\bU$ and $\bV$ can lose orthogonality as the recursion proceeds (similar to the CGS and MGS for computing the QR decomposition, Section~\ref{section:qr-gram-compute}).
Additionally, a variant of LGK bidiagonalization has been proposed to avoid storing all the vectors $\bu_i$ and $\bq_i$ by \citet{paige1982lsqr}.

\index{Approximate least squares}
\subsection*{Approximate Least Squares}
We further explore the approximation of the least squares problem $\mathop{\min}_{\bx}\normtwo{\bC\bx-\bb}$. Denote $\bQ_k\triangleq[\bq_1,\bq_2,\ldots,\bq_k]$, $\bU_k\triangleq[\bu_1,\bu_2,\ldots,\bu_k]$, and $\bU=[\bU_{k+1}, \bU_\perp]$. 
Additionally, let $\bB_k$ be the upper-left $k\times (k-1)$ submatrix of $\bB_2$. 
Once again, referring to Equation~\eqref{equation:lgk_bidi_ls}, we have:
$$
\bC\bQ_k = \bU_{k+1}\bB_{k+1}.
$$
Note that the variable $\bx$ lies in $\real^{n-1}$, and the vectors $\{\bq_1,\bq_2,\ldots,\bq_k\}$ are orthogonal in $\real^{n-1}$.
Approximately, we can estimate $\bx$ using a linear combination of the $k$ vectors, i.e., there exists a vector $\by$ such that $\bx\approx\bQ_k\by$.
Assume we want to find the optimal approximate solution within the subspace spanned by the $k$ vectors $\{\bq_1,\bq_2,\ldots,\bq_k\}$, i.e., solving the following problem in terms of $\by$:
\begin{equation}\label{equation:reduced_rank_bidiagonal}
\mathop{\min}_{\by} \normtwo{\bC\bQ_k\by - \bb},
\end{equation}
where $\by\in\real^k$ (it can be shown that $\by=\bQ_k^\top\bx \in\real^k$).
Based on the preceding discussion, the optimization problem is equivalent to:
$$
\begin{aligned}
\mathop{\min}_{\by}
\normtwo{\bU_{k+1}\bB_{k+1}\by - \bb}
&=
\mathop{\min}_{\by}
\normtwo{\bU^\top(\bU_{k+1}\bB_{k+1}\by - \bb)}\\
&=
\mathop{\min}_{\by}
\normtwo{\begin{bmatrix}
\bB_{k+1}\by \\
\bzero 
\end{bmatrix}
-
\begin{bmatrix}
b_{11}\be_1\\
\bzero 
\end{bmatrix}
}
=
\mathop{\min}_{\by}
\normtwo{\bB_{k+1}\by
-
b_{11}\be_1
}.
\end{aligned}
$$
Thus, the approximate least squares problem becomes 
$
\mathop{\min}_{\by}
\normtwo{
\bB_{k+1}\by
-
b_{11}\be_1
},
$
where $b_{11} = \normtwo{\bb}$.
Due to the bidiagonal structure, the problem can be solved in $\sim n$ flops \citep{elden2007matrix}.

\paragraph{Reduced-rank model.}
The problem in Equation~\eqref{equation:reduced_rank_bidiagonal} is known as  the least squares problem corresponding to the \textit{reduced-rank model}. Instead of considering the full model 
$\mathop{\min}_{\by} \normtwo{\bC\bx - \bb}$, we introduce an approximate orthogonal basis of low dimension in $\real^{n-1}$ where the solution $\bx$ lies (i.e.,  $\{\bq_1,\bq_2,\ldots, \bq_k\}$). 
This approach helps reduce the ill-conditioning of the problem and makes the solution less sensitive to perturbations in the data \citep{elden2007matrix}. 

\begin{problemset}
\item \label{prob:bidia_mln} We discussed the bidiagonalization for a matrix $\bA\in\real^{m\times n}$ with $m\geq n$ in the main section. Provide an algorithm that computes the bidiagonalization with $m<n$; discuss the complexity. Alternatively, discuss the algorithm for computing $\bA=\bU\bB\bV^\top$ with orthogonal $\bU,\bV$ and lower bidiagonal $\bB$ when $m\geq n$.
\item Prove in details that the LHC bidiagonalization method costs less when $m>\frac{5}{3}n$ compared to the Golub-Kahan bidiagonalization, as illustrated in Figure~\ref{fig:bidiagonal-gk-sample}.

\item Prove in details that the Three-Step bidiagonalization method costs less when $n\leq m<2n$ compared to the Golub-Kahan and LHC bidiagonalization methods.

\item (Read Chapter~\ref{chapter:SVD} first) Let  $\bA\in\real^{n\times n}$ be upper bidiagonal with  a repeated singular value. Show that $\bA$ must have a zero on its diagonal or superdiagonal.

\item \textbf{Singular values of bidiagonal (read Chapter~\ref{chapter:SVD} first \citep{bernstein2008matrix, mathias2014singular}).}
Let $\bA\in\real^{n\times n}$ be upper bidiagonal with the main diagonal values $\{a_1, a_2, \ldots, a_n\}$ and the superdiagonal values $\{b_1,b_2, \ldots, b_{n-1}\}$, and let $\bB\in\real^{n\times n}$ be bidiagonal. Show that
\begin{enumerate}
\item The singular values of $\bA$ are distinct.
\item If $\abs{\bB} = \abs{\bA}$, where $\abs{\cdot}$ returns the absolute values of a matrix, then $ \bA $ and $ \bB $ have the same singular values.
\item If $ \abs{\bA} \preceq \abs{\bB} $ (i.e., $\abs{\bB}-\abs{\bA}$ is PSD)  and $ \abs{\bA} \neq \abs{\bB} $, then $ \sigma_{\max}(\bA) < \sigma_{\max}(\bB) $.
\item If $ \abs{\bI \hadaprod \bA} \preceq \abs{\bI \hadaprod \bB} $ and $ \abs{\bI \hadaprod \bA} \neq \abs{\bI \hadaprod \bB} $, where $\hadaprod$ denotes the Hadamard product (Definition~\ref{definition:hada_prod}), then $ \sigma_{\min}(\bA) < \sigma_{\min}(\bB) $.
\item If $ \abs{\bI_{\text{up}} \hadaprod \bA} \preceq \abs{\bI_{\text{up}} \hadaprod \bB} $, and $ \abs{\bI_{\text{up}} \hadaprod \bA} \neq \abs{\bI_{\text{up}} \hadaprod \bB} $, where $ \bI_{\text{up}} $ denotes the matrix with all entries on the superdiagonal equal to 1 and all other entries equal to 0, then $ \sigma_{\min}(\bB) < \sigma_{\min}(\bA) $.
\end{enumerate}

\item Explore the process of bidiagonalization utilizing Givens rotations. What happens if the matrix is upper triangular or tridiagonal?

\item Let  $\bA\in\real^{n\times n}$ be upper bidiagonal with $a_{nn} = 0$. Show how to construct orthogonal matrices $\bU$ 
and $\bV$ (as products of Givens rotations) such that $\bU^\top\bA\bV$ is upper bidiagonal with the $n$-th column being zero.

\end{problemset}

%% file: chapter-eigenvalue.tex
\newpage 
\part{Eigenvalue Problem}
\section*{Introduction}
\lettrine{\color{caligraphcolor}T}
This part will introduce decompositional approaches  related to eigenvalues and eigenvectors.
\textit{Leonhard  Euler} revealed the existence of directions that are preserved by linear transformations (referred to as eigenvectors) in his investigation of the motions of rigid bodies.
\textit{Lagrange, Cauchy, Fourier, and Hermite} subsequently continued this line of research. 
The study of eigenvectors and eigenvalues has become increasingly  important due to their applications in stability theory and heat conduction. 
\textit{Hilbert}  later extended the research on eigenvalues to functional analysis, particularly in the theory of integral operators. The terms ``eigenvalue" and ``eigenvector" were defined by him. The word ``eigenvalue" is a combination of the  German word ``eigen" (which means ``own") with the English word ``value." 

To be more specific, the part will deal with decompositions related to orthogonal equivalence and orthogonal similarity transformations (Definition~\ref{definition:biequivalent}, Definition~\ref{definition:simiar_congru}).
The \textit{orthogonal similarity transformation} $\bA\in\real^{n\times n}\rightarrow \bQ^\top\bA\bQ$ corresponds to changing the basis from the given one to another orthonormal basis, where the orthogonal matrix $\bQ$ serves as the change of basis matrix.
While the \textit{orthogonal equivalence transformation} $\bA\in\real^{m\times n}\rightarrow \bU^\top\bA\bV$ involves changing the bases of the spaces 
$\real^m$ and $\real^n$ from the given orthonormal bases $\widetildebU$ and $\widetildebV$ to other orthonormal bases. Here, $\bU$ and $\bV$ are the orthogonal matrices representing the new bases.

\newpage
\chapter{Eigenvalue, Jordan Decomposition, Characteristic Polynomial}\label{chapter:eig_jordan}
\begingroup
\hypersetup{
	linkcolor=structurecolor,
	linktoc=page,  
}
\minitoc \newpage
\endgroup

\section{Eigenvalue Decomposition}
\index{Decomposition: EVD}
\begin{theoremHigh}[Eigenvalue Decomposition]\label{theorem:eigenvalue-decomposition}
Let $\bA\in \real^{n\times n}$ be any square matrix with a set of linearly independent eigenvectors. Then it can be factored as 
$$
\bA = \bX\bLambda\bX^{-1},
$$
where $\bX\in\real^{n\times n}$ is a matrix whose columns are the eigenvectors of $\bA$, and $\bLambda\in\real^{n\times n}$ is a diagonal matrix    with  eigenvalues $\lambda_1, \lambda_2, \ldots, \lambda_n$ on its diagonal, i.e., $\bLambda=\diag(\lambda_1, \lambda_2, \ldots, \lambda_n)$.
\end{theoremHigh}

\index{EVD}
\index{Eigenvalue decomposition}
Eigenvalue decomposition (EVD), also known as matrix diagonalization (see Definition~\ref{definition:diagonalizable}), is the process of transforming a matrix $\bA$ into a diagonal form.  
When all eigenvalues of $\bA$ are distinct, the associated  eigenvectors are guaranteed to be linearly independent, which enables the diagonalization of $\bA$. 
It is crucial to note that diagonalization is only possible if there exists a full set of $n$ linearly independent eigenvectors. 
In Section~\ref{section:otherform-spectral}, we will explore in more detail the conditions required for a matrix to have a complete set of linearly independent eigenvectors.

\begin{remark}[Left/Right Eigenvectors]
Suppose $\bA\in\real^{n\times n}$ admits the eigenvalue decomposition $\bA=\bX\bLambda\bX^{-1}$. The eigenvectors of $\bA^\top$ correspond to the columns of $(\bX^{-1})^\top=[\by_1, \by_2 , \ldots, \by_n]$. Therefore, the columns of $(\bX^{-1})^\top$ are often referred to as the \textbf{left eigenvectors} of $\bA$. Consequently, the matrix $\bA$ can also be expressed as:
$
\bA = \bX\bLambda\bX^{-1} = \lambda_1\bx_1\by_1^\top
+\lambda_2\bx_2\by_2^\top
+\ldots +
\lambda_n\bx_n\by_n^\top,
$
where $\bX=[\bx_1, \bx_2 , \ldots, \bx_n]$.
\end{remark}

\begin{definition}[Diagonalizable]\label{definition:diagonalizable}
A matrix $\bA$  is said to be diagonalizable if there exists a nonsingular matrix $\bP$ and a diagonal matrix $\bD$ such that $\bA = \bP\bD\bP^{-1}$.
In other words, $\bA$ is similar to a diagonal matrix. {In this case, the rows of $\bP^{-1}$ are left eigenvectors of $\bA$; and the columns of $\bP$ are right eigenvectors of $\bA$.}
\end{definition}

\section{Existence of  Eigenvalue Decomposition}
\begin{proof}[of Theorem~\ref{theorem:eigenvalue-decomposition}]
Let $\bX=[\bx_1, \bx_2, \ldots, \bx_n]$ denote the linearly independent eigenvectors of $\bA$. Clearly, we have
$$
\bA\bx_1=\lambda_1\bx_1,\qquad \bA\bx_2=\lambda_2\bx_2, \qquad \ldots, \qquad\bA\bx_n=\lambda_n\bx_n.
$$
In  matrix form, this can be written as:
$$
\bA\bX = [\bA\bx_1, \bA\bx_2, \ldots, \bA\bx_n] = [\lambda_1\bx_1, \lambda_2\bx_2, \ldots, \lambda_n\bx_n] = \bX\bLambda.
$$
Since we assume the eigenvectors are linearly independent,  $\bX$ has full rank and is invertible. Consequently, we can deduce:
$
\bA = \bX\bLambda \bX^{-1}.
$
This completes the proof.
\end{proof}

We will discuss some variants and similar forms of eigenvalue decomposition in the spectral decomposition section. 
In one scenario, matrix $\bA$ is required to be symmetric, and the factor $\bX$ is not only nonsingular but also orthogonal. 
Alternatively, if the matrix $\bA$  is a \textit{simple matrix} (i.e., the algebraic and geometric multiplicities of its eigenvalues are the same), $\bA$ can also be diagonalized. The decomposition also carries  a geometric interpretation, which we will explore in Section~\ref{section:coordinate-transformation}. 

An essential property of the matrix decomposition in the form of $\bA =\bX\bLambda\bX^{-1}$ is its computational efficiency in computing the $m$-th power.

\begin{remark}[$m$-th Power]\label{remark:power-eigenvalue-decom}
The $m$-th power of $\bA$ is $\bA^m = \bX\bLambda^m\bX^{-1}$ if  $\bA$ can be factored as $\bA=\bX\bLambda\bX^{-1}$.
Moreover, if $\bA=\bX\bLambda\bX^{-1}$, then the eigenvalues of $\bA^m$ are the $m$-th powers of the eigenvalues of $\bA$ (the result can be generalized using the Schur decomposition; see Remark~\ref{remark:mth_schur}).
\end{remark}
We note that computing $\bLambda^m$ is easy because we can apply this operation individually to each diagonal element.

It is important to note that the existence of the eigenvalue decomposition requires $\bA$ to have linearly independent eigenvectors. Fortunately, under certain conditions, this requirement is naturally satisfied.
\begin{proposition}[Different Eigenvalues]\label{proposition:diff-eigenvec-decompo}
Let the eigenvalues $\lambda_1, \lambda_2, \ldots, \lambda_n$ of $\bA\in \real^{n\times n}$ be all distinct. Then the corresponding eigenvectors are inherently linearly independent. 
In other words, any square matrix with distinct eigenvalues can be readily  diagonalized. 
\end{proposition}
\begin{proof}[of Proposition~\ref{proposition:diff-eigenvec-decompo}]
Suppose the eigenvalues $\lambda_1, \lambda_2, \ldots, \lambda_n$ are distinct, and the eigenvectors $\bx_1,\bx_2, \ldots, \bx_n$ are linearly dependent. That is, without loss of generality, we assume there exists a  nonzero vector $\bc = [c_1,c_2,\ldots,c_{n-1}]^\top$ such that 
$
\bx_n = \sum_{i=1}^{n-1} c_i\bx_{i}. 
$
Consequently, we can express
$$
\begin{aligned}
\bA \bx_n &= \bA \left(\sum_{i=1}^{n-1} c_i\bx_{i}\right) 
=c_1\lambda_1 \bx_1 + c_2\lambda_2 \bx_2 + \ldots + c_{n-1}\lambda_{n-1}\bx_{n-1}.
\end{aligned}
$$
and 
$$
\begin{aligned}
	\bA \bx_n &= \lambda_n\bx_n
	=\lambda_n (c_1\bx_1 +c_2\bx_2+\ldots +c_{n-1} \bx_{n-1}).
\end{aligned}
$$
Equating these two expressions gives
$
\sum_{i=1}^{n-1} (\lambda_n - \lambda_i)c_i \bx_i = \bzero .
$
This leads to a contradiction since $\lambda_n \neq \lambda_i$ for all $i\in \{1,2,\ldots,n-1\}$, from which the result follows.
\end{proof}

\begin{remark}[Limitations of Eigenvalue Decomposition]
Eigenvalue decomposition has certain limitations, which will be  addressed in the following chapters: 
\begin{itemize}
\item The eigenvectors in $\bX$ are often not orthogonal, and there might not always be a sufficient number of distinct eigenvectors (some eigenvalues may be identical).
	
\item To compute the eigenvalues and eigenvectors from the equation $\bA\bx = \lambda\bx$, the matrix   $\bA$ must be square. Rectangular matrices cannot undergo eigenvalue decomposition.
\end{itemize}
\end{remark}

\index{Rayleigh quotient}
The process of computing eigenvalues and eigenvectors involves solving a polynomial equation. Algorithms like the Rayleigh quotient iteration can effectively handle this task. We will delve deeper into various algorithms for determining the eigenvalues and eigenvectors of a matrix in Section~\ref{section:eigenvalue-problem}. Once the eigenvalues and eigenvectors are known, the next step is to compute the inverse of the nonsingular matrix $\bX$. 
As demonstrated in Theorem~\ref{theorem:inverse-by-lu2}, this operation typically requires a complexity of  $2n^3$ flops.

\section{Characteristic Polynomial and Multiplicity}
Let $\bA\in\complex^{n\times n}$ be given. Then the matrix $\lambda\bI-\bA$ is singular if and only if its determinant $\det(\lambda\bI-\bA)=0$, and if and only if $\lambda$ is an eigenvalue of $\bA$.
Therefore, $\lambda\bI-\bA$ is called the \textit{characteristic matrix} of $\bA$.
We also provide a definition for the characteristic polynomial:
\begin{definition}[Characteristic Polynomial\index{Characteristic polynomial}]\label{definition:characteristic_polynomial}
Let $\bA \in \complex^{n\times n}$ be any square matrix. Then, the \textit{characteristic polynomial}, denoted by  $p_{\bA}(\lambda)=\det(\lambda \bI- \bA)$, is given by 
\begin{equation}\label{equation:char_po_main}
\begin{aligned}
	p_{\bA}(\lambda)=\det(\lambda\bI-\bA ) &=\lambda^n + \gamma_{n-1} \lambda^{n-1} + \ldots + \gamma_1 \lambda  + \gamma_0\\
	&=(\lambda-\lambda_1)^{k_1} (\lambda-\lambda_2)^{k_2} \ldots (\lambda-\lambda_m)^{k_m},~\footnote{Note that $\det(\bA-\lambda\bI) = (-1)^n\det(\lambda\bI-\bA )$. See Problem~\ref{prob:consmatfrom_pol} for how to construct a matrix from a general polynomial.}
\end{aligned}
\end{equation}
where $\lambda_1, \lambda_2, \ldots, \lambda_m$ are the distinct roots of $\det( \lambda\bI-\bA)=0$, which are also the eigenvalues of $\bA$, and $k_1+k_2+\ldots +k_m=n$. 
Thus, the equation $\det( \lambda\bI-\bA)=0$ is called the \textit{characteristic equation} of $\bA$.
In other words, $\det(\lambda\bI-\bA)$ is a polynomial of degree $n$ for any matrix $\bA\in \complex^{n\times n}$ (see proof of Proposition~\ref{proposition:eigen-multiplicity}). It can be observed that (see Remark~\ref{remark:det_altsumper}):
$$
p_{\bA}(\lambda)= \lambda^n - \trace(\bA)t^{n-1} +\ldots+(-1)^n\det(\bA).
$$

\end{definition}
From the definition of the characteristic polynomial, if  $\bA\in\complex^{n\times n}$, then $\bA$ has $n$ eigenvalues (counting multiplicities). 
Although the eigenvalues do not have to be real, it can be observed that the coefficients of   $p_{\bA}(\lambda)$ are all real if $\bA$ is real.

\begin{proposition}[Characteristic Polynomial]
Let $\bA \in \complex^{n\times n}$ be any square matrix. Then, its characteristic polynomial can be expressed as 
$$
\det(\lambda \bI- \bA) = 
\sum_{k=0}^{n} (-1)^k a_k \lambda^{n-k},
$$
where $a_k$ is the sum of the principal minors (Definition~\ref{definition:principle-minors}) of order $k$ for matrix $\bA$. By convention, we set $a_0=1$.
\end{proposition}

\begin{exercise}[Characteristic Polynomial of Similar Matrices]
Let $\bA$ and $\bB$ be similar (Definition~\ref{definition:similar-matrices}). Show that they have the same characteristic polynomial. 
Show that the reverse is not. \textit{Hint: Consider matrices
$\footnotesize\begin{bmatrix}
0&0\\
0&0 
\end{bmatrix}$ 
and 
$\footnotesize\begin{bmatrix}
0&1\\
0&0 
\end{bmatrix}$.}
\end{exercise}

The characteristic polynomial of a matrix is utilized  to define  important multiplicities  as follows. 
While we define these concepts in the complex domain, their applications in matrix decomposition often involve real matrices in our case.
\begin{definition}[Algebraic Multiplicity and Geometric Multiplicity\index{Algebraic multiplicity}\index{Geometric multiplicity}]\label{definition:eigen_multipli}
Given the characteristic polynomial of a matrix $\bA\in \complex^{n\times n}$:
$$
\begin{aligned}
\det(\lambda\bI-\bA ) =(\lambda-\lambda_1)^{k_1} (\lambda-\lambda_2)^{k_2} \ldots (\lambda-\lambda_m)^{k_m},
\end{aligned}
$$
the integer $k_i$ is called the \textit{algebraic multiplicity} of the eigenvalue $\lambda_i$, i.e., the algebraic multiplicity of an eigenvalue $\lambda_i$ is equal to the multiplicity of the corresponding root of the characteristic polynomial.
Therefore, the set of distinct eigenvalues is called the spectrum of $\bA$; while the set including eigenvalues with their multiplicities is called the multispectrum (see Definition~\ref{definition:spectrum}).

\index{Eigenspace}
The \textit{eigenspace associated with the eigenvalue $\lambda_i$} is defined as the null space of $(\bA - \lambda_i\bI)$, i.e., $\nspace(\bA - \lambda_i\bI)$.
And the dimension of the eigenspace associated with $\lambda_i$, $\nspace(\bA - \lambda_i\bI)$, is referred to as the \textit{geometric multiplicity} of $\lambda_i$.

In summary, we denote the algebraic multiplicity of $\lambda_i$ by $alg(\lambda_i)$ and its geometric multiplicity by $geo(\lambda_i)$.
\end{definition}

\begin{exercise}
Suppose all the eigenvalues of $\bA\in\real^{n\times n}$ are zero. Show that $\bA^n=\bzero$. 
\textit{Hint: Use Cayley-Hamilton Theorem~\ref{theorem:cayley_hami}.}
\end{exercise}

\index{Multiplicity}
\begin{remark}[Geometric Multiplicity]\label{remark:geometric-mul-meaning}
Note that for a matrix $\bA$ and its eigenspace $\nspace(\bA-\lambda_i\bI)$ associated with the eigenvalue $\lambda_i$, the dimension of the eigenspace is equivalent to the number of linearly independent eigenvectors of $\bA$ associated with $\lambda_i$, which form a basis for the  eigenspace. 
This implies that,
although there are infinitely many  eigenvectors associated with each eigenvalue $\lambda_i$,  they form a subspace (provided the zero vector is added) that can be described by a finite number of vectors.
\end{remark}

Proposition~\ref{proposition:eigenvalue-similar-matrices} shows that similar matrices have the same eigenvalues. The definition of multiplicities also indicates their equivalence.
\begin{corollary}[Multiplicity in Similar Matrices\index{Similar matrices}]\label{corollary:multipli-similar-matrix}
Similar matrices share the same algebraic  and geometric multiplicities.
\end{corollary}
\begin{proof}[of Corollary~\ref{corollary:multipli-similar-matrix}]
In Proposition~\ref{proposition:eigenvalue-similar-matrices}, we proved that the eigenvalues of similar matrices are the same. Therefore, the algebraic multiplicities of similar matrices are also the same.

Suppose $\bA$ and $\bB= \bP\bA\bP^{-1}$ are similar matrices, where $\bP$ is nonsingular. 
Let the geometric multiplicity of an eigenvalue of $\bA$, denoted by $\lambda$, be $k$. 
Then, there exists a set of orthogonal (linearly independent) vectors $\bv_1, \bv_2, \ldots, \bv_k$ that constitute a basis for the eigenspace $\nspace(\bA-\lambda\bI)$, satisfying $\bA\bv_i = \lambda \bv_i$ for all $i\in \{1, 2, \ldots, k\}$. 
Define $\bw_i \triangleq \bP\bv_i$ for each $i$. These $\bw_i$'s are eigenvectors of $\bB$ associated with the eigenvalue $\lambda$. 
Furthermore, these $\bw_i$'s are linearly independent since $\bP$ is nonsingular. Thus, the dimension of the eigenspace $\nspace(\bB-\lambda\bI)$ is at least $k$. That is, $\dim(\nspace(\bA-\lambda\bI)) \leq \dim(\nspace(\bB-\lambda\bI)) $. 

Similarly, there exists  a set of orthogonal vectors $\bw_1, \bw_2, \ldots, \bw_k$ that form a basis for the eigenspace $\nspace(\bB-\lambda\bI)$. Then, $\bv_i \triangleq \bP^{-1}\bw_i$ for all $i \in \{1, 2, \ldots, k\}$ are  eigenvectors of $\bA$ associated with the eigenvalue $\lambda$. This  results in $\dim(\nspace(\bB-\lambda\bI)) \leq \dim(\nspace(\bA-\lambda\bI)) $. 

By combining both inequalities, we conclude that $\dim(\nspace(\bA-\lambda\bI)) = \dim(\nspace(\bB-\lambda\bI)) $, establishing the equality of geometric multiplicities for similar matrices.
\end{proof}

According to the definition, the sum of the algebraic multiplicities of all eigenvalues is equal to $n$, while the sum of the geometric multiplicities can be strictly smaller. For each eigenvalue, this inequality also holds true.
\begin{theorem}[Bounded Geometric Multiplicity]\label{theorem:bounded-geometri}
Let  $\bA\in \real^{n\times n}$ be  any matrix. Then, its geometric multiplicity of  any eigenvalue $\lambda_i$ is bounded by its algebraic multiplicity:
$$
geo(\lambda_i) \leq alg(\lambda_i).
$$
This also indicates $\rank(\bA-\lambda\bI)\geq n-k$ if $alg(\lambda)=k$, with equality for $k=1$.
Alternatively, the  $\rank(\bA-\lambda\bI)=n-r$ if $geo(\lambda)=r$.
\end{theorem}
\begin{proof}[of Theorem~\ref{theorem:bounded-geometri}]
If we can find a similar matrix $\bB$ of $\bA$ that has a specific form of the characteristic polynomial, then we complete the proof. 

Suppose $\bP_1 = [\bv_1, \bv_2, \ldots, \bv_k]$ contains a set of  linearly independent eigenvectors of $\bA$ associated with $\lambda_i$. 
That is, these $k$ vectors form a basis for the eigenspace $\nspace(\bA-\lambda_i\bI)$, and the geometric multiplicity associated with $\lambda_i$ is $k$. We can extend  this set to a basis of  $\real^n$ by adding more $n-k$ linearly independent vectors, yielding:
$$
\bP = [\bP_1, \bP_2] = [\bv_1, \bv_2, \ldots, \bv_k, \bv_{k+1}, \ldots, \bv_n],
$$
where $\bP$ is nonsingular. Consequently, $\bA\bP = [\lambda_i\bP_1, \bA\bP_2]$.

Construct a matrix 
$\bB = 
\scriptsize
\begin{bmatrix}
\lambda_i \bI_k & \bC \\
\bzero  & \bD
\end{bmatrix}$, where $\bA\bP_2 = \bP_1\bC + \bP_2\bD$. 
Then, $\bP^{-1}\bA\bP = \bB$ such that $\bA$ and $\bB$ are similar matrices. We can always find such matrices $\bC$ and $\bD$ satisfying the above condition, since $\bv_i$'s are linearly independent  and span the entire space $\real^n$, ensuring that any column of $\bA\bP_2$ belongs to the column space of $\bP=[\bP_1,\bP_2]$.
Therefore, 
$$
\begin{aligned}
\det(\bA-\lambda\bI) &\stackrel{*}{=} \det(\bP^{-1})\det(\bA-\lambda\bI)\det(\bP)  
\stackrel{+}{=}\det(\bP^{-1}(\bA-\lambda\bI)\bP)  
=  \det(\bB-\lambda\bI) \\
&= 
\det\left(
\footnotesize\begin{bmatrix}
(\lambda_i-\lambda) \bI_k  & \bC \\
\bzero  & \bD - \lambda \bI
\end{bmatrix}
\right)
\normalsize= (\lambda_i-\lambda)^k \det(\bD-\lambda\bI),
\end{aligned}
$$ 
where equality ($*$) follows from $(\text{$\det(\bP^{-1}) = 1/\det(\bP)$})$, equality $(+)$ follows from the fact that $(\det(\bA)\det(\bB) = \det(\bA\bB))$, and the last equality follows from the determinant of block matrices (Remark~\ref{remark:determinant-intermezzo}).
This implies 
$
geo(\lambda_i) \leq alg(\lambda_i).
$
And we complete the proof.
\end{proof}

\section{Jordan Decomposition}
The \textit{canonical form} of a matrix, also known as its \textit{normal or standard form}, is a method of representing the matrix in a simplified  and standardized way. 
One of the simplest forms of a matrix's canonical representation is its diagonal form. However, not every matrix can be diagonalized.
In eigenvalue decomposition, we assume the matrix $\bA$ has $n$ linearly independent eigenvectors. This assumption, however, does not apply to all square matrices.  
In numerous applications, it is crucial to simplify a given matrix to its most basic form, which is where matrix similarity reduction comes into play (for example, the Hessenberg decomposition in Theorem~\ref{theorem:hessenberg-decom}).
To address limitations in eigenvalue decomposition, a more generalized version, called the \textit{Jordan decomposition}, is introduced. Named after Camille Jordan \textit{Camille Jordan} \citep{jordan1870traite}, the Jordan canonical form serves as the standard for similarity reduction.

Before delving into Jordan decomposition, let us define \textit{Jordan blocks} and \textit{Jordan form}.
\begin{definition}[Jordan Block]\label{definition:jordan_block}
An $m\times m$ upper triangular matrix $\bJ_m(\lambda)$ is referred to as a Jordan block provided all $m$ diagonal elements are identical and equal to the value $\lambda$, and all superdiagonal elements are set to one:
$$
\bJ_m(\lambda)=
\begin{bmatrix}
\lambda & 1        & 0 & \ldots & 0 & 0 & 0 \\
0       & \lambda  & 1 & \ldots & 0 & 0 & 0 \\
0       & 0  & \lambda & \ldots & 0 & 0 & 0 \\
\vdots  & \vdots   & \vdots & \ddots & \vdots & \vdots\\
0 & 0 & 0 & \ldots & \lambda & 1& 0\\
0 & 0 & 0 & \ldots & 0 &\lambda & 1\\
0 & 0 & 0 & \ldots & 0 & 0 & \lambda
\end{bmatrix}_{m\times m}.
$$
\end{definition}

\begin{exercise}[Properties of Jordan Block]\label{exercise:nil_jor}
For the Jordan block $\bJ_m(\lambda)$, show that 
\begin{itemize}
\item $\bJ_{m}(0)$ is nilpotent, and the nilpotency (Definition~\ref{definition:niopotent_mat}) is $m$. 
\item $\rank(\bJ_m(0)^i)=m-i$, for $i\in\{1,2,\ldots,m\}$.
\item $\bJ_{m}(0) \be_{i+1}=\be_i$ for all $i\in\{1,2,\ldots,m-1\}$.
\item $(\bI_m - \bJ_{m}(0)^\top \bJ_{m}(0))\bx = \bx^\top\be_1\be_1$.
\item $\bJ_{m}(0)^\top \bJ_{m}(0)=\begin{bmatrix}
	0 & \bzero\\
	\bzero & \bI_{m-1}
\end{bmatrix}$.
\end{itemize}

\end{exercise}

\begin{exercise}[Jordan Block]
For a Jordan block $\bJ_m(\lambda)$, show that 
\begin{itemize}
\item The eigenspace associated to $\lambda$ has a dimension of one, i.e., $geo(\lambda)=1$.
\item The algebraic multiplicity of $\lambda$ is $m$.
\end{itemize}
\end{exercise}

\begin{definition}[Jordan Form\index{Jordan block}]
Let $\bA$ be  an $n\times n$ matrix. Then, a Jordan form $\bJ$ for $\bA$ is a block-diagonal matrix defined as
$$
\bJ=\diag\big(\bJ_{m_1}(\lambda_1), \bJ_{m_2}(\lambda_2), \ldots \bJ_{m_k}(\lambda_k)\big),~\footnote{In many texts, $\bJ$ is also denoted as
$
\bJ= \bJ_{m_1}(\lambda_1)\oplus \bJ_{m_2}(\lambda_2)\oplus \ldots \oplus \bJ_{m_k}(\lambda_k) 
$, i.e., the direct sum notation.}
$$
where $\lambda_1, \lambda_2, \ldots, \lambda_k$ are eigenvalues of $\bA$ (duplicates possible), and $m_1+m_2+\ldots+m_k=n$.
\end{definition}

Although not all matrices admit the eigenvalue decomposition, they can be factored using the Jordan decomposition.
A nondiagonalizable matrix $\bA$ with multiple eigenvalues can be reduced to its Jordan canonical form through a similarity transformation.
\begin{theoremHigh}[Jordan Decomposition]\label{theorem:jordan-decomposition}
Let $\bA\in \real^{n\times n}$ be any  square matrix  with real eigenvalues. Then, it  can be factored as 
$$
\bA = \bX\bJ\bX^{-1}, ~\footnote{If $\bA\in\complex^{n\times n}$ or eigenvalues are not real, then $\bX,\bJ\in\complex^{n\times n}$.}
$$
where $\bX\in\real^{n\times n}$ is a nonsingular matrix containing the \textit{generalized eigenvectors} of $\bA$ as its columns, and $\bJ\in\real^{n\times n}$ is a Jordan form matrix $\diag(\bJ_{m_1}(\lambda_1), \bJ_{m_2}(\lambda_2), \ldots, \bJ_{m_k}(\lambda_k))$, with 
$$
\bJ_{m_i}(\lambda_i) = 
\begin{bmatrix}
\lambda_i & 1         & \ldots & 0 & 0  \\
0       & \lambda_i   & \ldots & 0 & 0  \\
\vdots  & \vdots   & \vdots & \vdots & \vdots\\
0 & 0  & \ldots & \lambda_i & 1\\
0 & 0  & \ldots & 0 &\lambda_i \\
\end{bmatrix}_{m_i\times m_i}, 
\gap i\in\{1,2,\ldots,k\},
$$
an $m_i\times m_i$ square matrix with $m_i$ being the number of repetitions of the eigenvalue $\lambda_i$,  $m_1+m_2+\ldots +m_k = n$, and $\{\lambda_1, \lambda_2, \ldots, \lambda_k\}$ are eigenvalues of $\bA$ with duplicates possible. 
The Jordan form $\bJ$ is uniquely determined by $\bA$ up to permutation of the Jordan blocks, and is called the \textit{Jordan  form} or \textit{Jordan canonical form} of $\bA$.
Furthermore, the nonsingular matrix $\bX$ is referred to as  the \textit{matrix of generalized eigenvectors} of $\bA$.
\end{theoremHigh}

As an example, a Jordan form can have the following structure:
$$
\begin{aligned}
\bJ&=\diag(\bJ_{m_1}(\lambda_1), \bJ_{m_2}(\lambda_2), \ldots, \bJ_{m_k}(\lambda_k))\\
&=
\footnotesize
\begin{bmatrix}
\begin{bmatrix}
\lambda_1 & 1        & 0 \\
0       & \lambda_1  & 1 \\
0       & 0  & \lambda_1
\end{bmatrix} &  & & &\\  
& \begin{bmatrix}
\lambda_2
\end{bmatrix} & & &\\  
&  &\begin{bmatrix}
\lambda_3 & 1 \\
0 & \lambda_3
\end{bmatrix} & &\\  
&  &   & \ddots & &\\  
&  &   &  & &\begin{bmatrix}
\lambda_k & 1 \\
0 & \lambda_k
\end{bmatrix}\\  
\end{bmatrix}.
\end{aligned}
$$
Note that zeros can appear on the superdiagonal of $\bJ$, and the first column is always a vector containing only eigenvalues of $\bA$ in each block. 

\subsection{Properties of Jordan Decomposition}

\paragraph{Jordan decomposition of a block-diagonal matrix.}
Let $\bA\triangleq\diag(\bA_1, \bA_2, \ldots,\bA_k)$ be a block-diagonal matrix, where each block has a Jordan decomposition $\bA_i = \bX_i\bJ_i\bX_i^{-1}$ for all $i\in\{1,2,\ldots,k\}$.
Let $\bX\triangleq\diag(\bX_1,\bX_2,\ldots,\bX_k)$ and $\bJ\triangleq\diag(\bJ_1,\bJ_2,\ldots,\bJ_k)$. Therefore, the Jordan decomposition of $\bA$ is given by 
$
\bA = \bX\bJ\bX^{-1}.
$

\paragraph{Nilpotency of Jordan blocks.}
The Jordan block $\bJ_m(0)$ with zeros on the diagonal is nilpotent such that $\bJ_m(0)^p=\bzero$ with $p\geq m$, and 
\begin{equation}\label{equation:npjd1}
\rank(\bJ_m(0)^{p-1}) - \rank(\bJ_m(0)^{p})
=
\begin{cases}
	1, \quad \text{if $p\leq m$};\\
	0, \quad \text{if $p>m$}.
\end{cases}
\end{equation}

\paragraph{Commuting decomposition.}
For any Jordan block, we have the identity $\bJ_m(\lambda) = \lambda \bI_m + \bJ_m (0)$, where $\bJ_m (0)$ is nipotent such that  $\bJ_m (0)^m = \bzero$. Thus, any Jordan
block is the sum of a diagonal matrix and a nilpotent matrix.
More generally, the Jordan canonical form $\bJ=\bX^{-1}\bA\bX$ of any matrix $\bA$ can be written as $\bJ = \bD + \bN$, where $\bD$ is a diagonal matrix whose main diagonal is the same as that of $\bJ$, and $\bN = \bJ -\bD$.
The matrix $\bN$ is also nilpotent, and $\bN^m = \bzero$ if $m$ is the size of the largest Jordan block in $\bJ$.
Since both $\bD$ and $\bN$ has conformal block-diagonal matrices, they commute such that $\bD\bN=\bN\bD$.

In this sense, $\bA$ can be decomposed as $\bA=\bX\bJ\bX^{-1} =\bX\bD\bX^{-1}+\bX\bN\bX^{-1} \triangleq \bA_d+\bA_n$, where $\bA_d\triangleq\bX\bD\bX^{-1}$ is the eigenvalue decomposition of $\bA_d$ and $\bA_n\triangleq\bX\bN\bX^{-1}$ is nilpotent. 
Moreover, since $\bD$ and $\bN$ commute,  it also follows that $\bA_d \bA_n = \bA_n \bA_d$. This demonstrates  that any square matrix can be decomposed into the sum of two commuting matrices.

\paragraph{Rank of power: $\rank(\bA-\lambda\bI)^p$.} The Jordan decomposition is unique up to permutation of eigenvalues.
Without loss of generalization, we assume 
$$
\bX^{-1}\bA\bX = \bJ = \diag(\bJ_{m_1}(\lambda), \bJ_{m_2}(\lambda), \ldots, \bJ_{m_z}(\lambda), \widetildebJ), 
$$
where $\widetildebJ\in\real^{m\times m}$ and $m\triangleq n-m_1-m_2-\ldots -m_z$.
That is, the eigenvalue  $\lambda$ appears in the upper left blocks.
Exercise~\ref{exercise:sim_trans_inhess} shows $(\bA-\lambda\bI)^p$ and $(\bJ-\lambda\bI)^p$ are similar for $p=1,2,\ldots$, and they share the same rank (Proposition~\ref{proposition:eigenvalue-similar-matrices}). Since $\rank((\widetildebJ-\lambda\bI)^p)=m$, we have:
$$
\begin{aligned}
\rank&\left((\bA-\lambda\bI)^p\right) 
=\rank\left((\bJ-\lambda\bI)^p\right)\\
&=
\rank(\bJ_{m_1}(0)^p)+\rank(\bJ_{m_2}(0)^p) +\ldots+\rank(\bJ_{m_z}(0)^p) +m.
\end{aligned}
$$
When $p\geq \max\{m_1,m_2,\ldots,m_z\}$, it follows that $\rank\left((\bA-\lambda\bI)^p\right) =m$.
Define 
\begin{equation}
r_p(\bA,\lambda) \triangleq \rank\left((\bA-\lambda\bI)^p\right)
\quad\text{and}\quad
d_p(\bA,\lambda)\triangleq r_{k-1}(\bA,\lambda) - r_p(\bA,\lambda).
\end{equation}
When $\lambda $ is not an eigenvalue of $\bA$, then $r_{k-1}(\bA,\lambda) = r_p(\bA,\lambda)=n$ and therefore, $d_p(\bA,\lambda)=0$.
However, when $\lambda$ is an eigenvalue of $\bA$, the above discussions and \eqref{equation:npjd1} show that the value $d_p(\bA,\lambda)$ has the following interpretation:
$$
\begin{aligned}
&d_p(\bA,\lambda) = (1 \text{ if $m_1\geq p$}) + (1 \text{ if $m_z\geq p$}) +\ldots+(1 \text{ if $m_z\geq p$})\\
&= \text{number of blocks in the Jordan form with the  eigenvalue $\lambda$ that have size at least $p$}.
\end{aligned}
$$ 
This leads to the  algorithm to determine the number of Jordan blocks for each eigenvalue in the next paragraph.

\index{Index of an eigenvalue}
\begin{definition}[Index of an Eigenvalue]\label{definition:index_eig_jor}
Let $\bA\in\real^{n\times n}$ with a real eigenvalue $\lambda$. If $d_p(\bA,\lambda) = 0$ for all $q>m$, then the value $m$ is called the index of $\lambda$ as an eigenvalue of $\bA$, i.e., the size of the largest Jordan block of $\bA$ associated with the eigenvalue $\lambda$.
If $\lambda$ is not an eigenvalue of $\bA$, then the index of $\lambda$ is 0.
\end{definition}

\begin{exercise}[Index of Block-Diagonal and Block Upper Triangular]
Let the index of $\lambda$ in $\bA_{1}$ be $d_1$ and the index of $\lambda$ in $\bA_2$ be $d_2$,  where $\bA_1\in\real^{n_1\times n_1}$ and $\bA_2\in\real^{n_2\times n_2}$ are given. Show that the index of $\lambda$ in $\diag(\bA_1,\bA_2)\in\real^{(n_1+n_2)\times(n_1+n_2)}$ is $\max\{d_1,d_2\}$.

More generally, let $\bA=[\bA_{ij}]_{i,j=1}^{k,k}$ be a block upper triangular containing $k\times k$ blocks, where $\bA_{ij}=\bzero$ if $i>j$.
Let further the index of  $\lambda$ in each diagonal block $\bA_{ii}$ be $d_i$ for $i\in\{1,2,\ldots,k\}$. Show that the index of $\lambda$ as an eigenvalue of $\bA$ is at most $d_1 +d_2+\ldots+ d_k$.
\end{exercise}

\paragraph{Number of Jordan blocks.}
We previously noted that the $k$ eigenvalues ($\lambda_1,\lambda_2,\ldots,\lambda_k$) mentioned in the theorem might include duplicates.
A repeated eigenvalue $\lambda$ with algebraic multiplicity $m_i$ can correspond to one or more Jordan blocks, depending on the geometric multiplicity of $\lambda$ (see discussion in Definition~\ref{definition:eigen_multipli}). For example, consider a $3\times 3$ matrix $\bB$ with an eigenvalue $\lambda$ having an algebraic multiplicity of 3. The Jordan  form of the matrix $\bB$ could take on any of the following three configurations:
\begin{align*}
&\bJ_1=\begin{bmatrix}\bJ_1(\lambda)& 0& 0\\ 0& \bJ_1(\lambda)& 0\\ 0& 0& \bJ_1(\lambda)\end{bmatrix}=\begin{bmatrix}\lambda& 0& 0\\ 0&\lambda& 0\\ 0& 0&\lambda\end{bmatrix}&\quad&\text{(geometric multiplicity=3)};\\
&\bJ_2=\begin{bmatrix}\bJ_1(\lambda)& 0\\ 0& \bJ_2(\lambda)\end{bmatrix}=\begin{bmatrix}\lambda& 0& 0\\ 0&\lambda& 1\\ 0& 0&\lambda\end{bmatrix}&\quad&\text{(geometric multiplicity=2)};\\
&\bJ_3=\bJ_3(\lambda)=\begin{bmatrix}\lambda& 1& 0\\ 0&\lambda& 1\\ 0& 0&\lambda\end{bmatrix}&\quad&\text{(geometric multiplicity=1)},
\end{align*}
because the meaning of the geometric multiplicity $p$ of the eigenvalue $\lambda$ is that the number of linearly independent eigenvectors corresponding to $\lambda$ is $p$ (Remark~\ref{remark:geometric-mul-meaning}).

\begin{proposition}[Examination of Number of Jordan Blocks]
Let $\bA\in\real^{n\times n}$. Then the number of Jordan blocks associated with the eigenvalue $\lambda_i$ is given by:
\begin{itemize}
\item The number of Jordan blocks with order greater than $1$ (i.e., $\geq 2$) is determined by
$
\alpha_{1} = \rank( \bA - \lambda_{i} \bI) - \rank(\bA - \lambda_{i} \bI)^{2}.
$

\item The number of Jordan blocks with order greater than $2$ (i.e., $\geq 3$) is determined by
$
\alpha_{2} = \rank(\bA - \lambda_{i} \bI)^{2} - \rank(\bA - \lambda_{i} \bI)^{3}.
$

\item More generally, the number of Jordan blocks with the order greater than $p-1$ (i.e., $\geq p$) is determined by
$
\alpha_{p-1} = \rank(\bA - \lambda_{i} \bI)^{p-1} - \rank(\bA - \lambda_{i} \bI)^{p}.
$

\item The sum of the orders of the Jordan blocks corresponding to an eigenvalue $\lambda_{i}$ is equal to its algebraic multiplicity.
\end{itemize}
\end{proposition}

It is worth noting that Jordan canonical forms that differ only in the permutation order of their Jordan blocks are considered equivalent. For example, the second-order Jordan canonical form $\bJ_2$ can also be represented as:
\begin{align*}
\bJ_2=\begin{bmatrix}\bJ_2 (\lambda)& 0\\ 0& \bJ_1(\lambda)\end{bmatrix}=\begin{bmatrix}\lambda& 1& 0\\ 0&\lambda& 0\\ 0& 0&\lambda\end{bmatrix}.
\end{align*}

\paragraph{Similarity of a matrix and its transpose.}
Consider a Jordan block $\bJ_m(\lambda)$ and a reversal matrix $\bP$ (satisfying $\bP^2=\bI$ and $\bP=\bP^\top$; see Definition~\ref{definition:permutation-matrix}). We can verify that 
$$
\bP\bJ_m(\lambda) = \bJ_m(\lambda)^\top\bP
\quad\implies\quad 
\bJ_m(\lambda) = \bP\bJ_m(\lambda)^\top\bP.
$$
Therefore, a Jordan form $\bJ$ containing several Jordan blocks on its main diagonal also has this property: $\bJ=\bP\bJ^\top\bP$.
Given the Jordan decomposition of $\bA=\bX\bJ\bX^{-1}$. It then follows that
\begin{equation}\label{equation:simi_trans}
\bA^\top = \bX^{-\top}\bJ^\top\bX^\top = \bX^{-\top}\bP\bJ\bP\bX^\top 
=
\bX^{-\top}\bP\bX^{-1}\bA\bX\bP\bX^\top 
\triangleq
\bS^{-1}\bA\bS,
\end{equation}
where $\bS\triangleq\bX\bP\bX^\top $. Therefore, $\bA$ and  $\bA^\top$ are similar.

\subsection{Generalized Jordan Decomposition}
\index{Decomposition: Generalized Jordan}
\index{Generalized Jordan decomposition}
Let $\bJ_m(\lambda, \epsilon)$ be the \textit{generalized Jordan block} replacing positive ones of a (standard) Jordan block $\bJ_m(\lambda)$ with $\epsilon$, i.e., with $\epsilon$'s on the superdiagonal instead of ones. Then, we can find similarity transformations to a generalized Jordan form for a matrix. 
\begin{corollaryHigh}[Generalized Jordan Decomposition]\label{corollary:gen_jordan_decom}
Let $\bA\in\real^{n\times n}$ be any square matrix  with real eigenvalues. Then, it can be factored as 
$$
\bA = \bX_\epsilon\bJ_\epsilon\bX_\epsilon^{-1}, ~\footnote{If $\bA\in\complex^{n\times n}$ or eigenvalues are not real, then $\bX_\epsilon,\bJ_\epsilon\in\complex^{n\times n}$.}
$$
where $\bX_\epsilon\in\real^{n\times n}$ is a nonsingular matrix, and $\bJ_\epsilon\in\real^{n\times n}$ is a generalized Jordan form matrix $\diag(\bJ_{m_1}(\lambda_1, \epsilon), \bJ_{m_2}(\lambda_2, \epsilon), \ldots, \bJ_{m_k}(\lambda_k, \epsilon))$,
with
$$
\bJ_{m_i}(\lambda_i) = 
\begin{bmatrix}
	\lambda_i & \epsilon         & \ldots & 0 & 0  \\
	0       & \lambda_i   & \ldots & 0 & 0  \\
	\vdots  & \vdots   & \vdots & \vdots & \vdots\\
	0 & 0  & \ldots & \lambda_i & \epsilon\\
	0 & 0  & \ldots & 0 &\lambda_i \\
\end{bmatrix}_{m_i\times m_i}
\gap i\in\{1,2,\ldots,k\},
$$
an $m_i\times m_i$ square matrix with $m_i$ being the number of repetitions of the eigenvalue $\lambda_i$,  $m_1+m_2+\ldots +m_k = n$, and $\{\lambda_1, \lambda_2, \ldots, \lambda_k\}$ are eigenvalues of $\bA$ with duplicates possible. 
The matrix $\bJ_\epsilon$ is uniquely determined by $\bA$ up to permutation of the Jordan blocks.

\end{corollaryHigh}
The generalized Jordan decomposition can be obtained by applying a further similarity transformation on the (standard) Jordan decomposition: 
$$
\bJ=\bX^{-1}\bA\bX
\implies  
\bJ_\epsilon=\bS_\epsilon^{-1}\bX^{-1}\bA\underbrace{\bX\bS_\epsilon}_{\bX_\epsilon},
$$  
where $\bS_\epsilon = \diag(\bS_{\epsilon, m_1}, \bS_{\epsilon, m_2}, \ldots, \bS_{\epsilon, m_k})$, and $\bS_{\epsilon, m_i}=\diag(1, \epsilon,\epsilon^2, \ldots, \epsilon^{m_i-1})$.

\section{Existence of Jordan Decomposition}
There are two primary  approaches to obtaining the Jordan canonical form of a matrix:
\begin{itemize}
\item \textit{Direct similarity reduction.} This method involves directly reducing a constant matrix to its simplest form through similarity transformations.
\item \textit{Balanced reduction and Smith normal form.} This approach starts with a balanced reduction of the polynomial matrix associated with the original constant matrix. Subsequently, the Smith normal form derived from the balanced reduction is transformed into the Jordan canonical form for the similarity reduction of the original matrix.
\end{itemize}
This section focuses on the implementation of the first method.
The Jordan decomposition is not particularly interesting in practice as it is extremely sensitive to perturbations. Even with the smallest random change to a matrix, the matrix can be made diagonalizable \citep{van2020advanced}. As a result, there are no practical mathematical software libraries or tools that compute it. 
To prove the existence of the Jordan decomposition, we use the following lemma, which identifies the similarity transformation for upper triangular matrices:
\begin{lemma}[Jordan Form of Strictly Upper Triangular Matrices]\label{lemma:upptri_jordan}
Let $\bA\in\real^{n\times n}$ be any square strictly upper triangular. Then, it can be factored as
$$
\bA=\bX\bJ\bX^{-1},~\footnote{When $\bA\in\complex^{n\times n}$ is complex, $\bX$ can be chosen complex as well.}
$$
where $\bX\in\real^{n\times n}$ is nonsingular, $\bJ=\diag\big(\bJ_{m_1}(0), \bJ_{m_2}(0), \ldots, \bJ_{m_k}(0)\big)$ is strictly upper triangular, $m_1\geq m_2\geq \ldots m_k\geq 1$, and $m_1+m_2+\ldots+m_k=n$.
\end{lemma}
\begin{proof}[of Lemma~\ref{lemma:upptri_jordan}]
We prove the existence of the  decomposition using induction. The result is self-evident for the case of  $n=1$. Assume it holds for all $n-1$. If we prove the decomposition also exists for $n$, then we complete the proof.
For brevity,  let $\bJ_{m}\triangleq\bJ_m(0)$ denote the Jordan block with $\lambda=0$.
We consider any strictly upper triangular $\bA=\footnotesize\begin{bmatrix}
0 & \ba^\top \\
\bzero & \bA_1 
\end{bmatrix}$, where $\bA_1\in\real^{(n-1)\times (n-1)}$.
By the induction hypothesis, there is a nonsingular matrix $\bX_1\in\real^{(n-1)\times (n-1)}$ such that 
$\bA_1=\bX_1\bJ_1\bX_1^{-1}\in\real^{(n-1)\times (n-1)}$, where 
$$
\bJ_1=\diag\big(\bJ_{n_1}, \bJ_{n_2}, \ldots,\bJ_{n_p} \big)
=
\begin{bmatrix}
\bJ_{n_1} & \bzero \\
\bzero & \bJ
\end{bmatrix},
$$
$n_1\geq n_2\geq \ldots n_p\geq 1$, and $\sum_{j=1}^{p}n_j=n-1$. This indicates $\bJ^{n_1}=\bzero $ by the nilpotency property (Exercise~\ref{exercise:nil_jor}).
Thus,
$$
\begin{bmatrix}
1 & \bzero \\
\bzero & \bX_1 
\end{bmatrix}^{-1}
\begin{bmatrix}
	0 & \ba^\top \\
	\bzero & \bA_1 
\end{bmatrix}
\begin{bmatrix}
	1 & \bzero \\
	\bzero & \bX_1 
\end{bmatrix}
=
\begin{bmatrix}
0 & \ba^\top \bX_1 \\
\bzero & \bX_1^{-1}\bA_1\bX_1 
\end{bmatrix}
\stackrel{[\bc_1^\top,\bc_2^\top]\triangleq\ba^\top \bX_1}{=}
\begin{bmatrix}
0 & \bc_1^\top & \bc_2^\top \\
\bzero & \bJ_{n_1} &  \bzero \\
\bzero & \bzero & \bJ
\end{bmatrix}.
$$
Next, consider a further similarity transformation on the above equation (and use the property in Exercise~\ref{exercise:nil_jor}):
\begin{equation}\label{equation:jd_q_s1}
\begin{bmatrix}
1 & -\bc_1^\top \bJ_{n_1}^\top & \bzero \\
\bzero & \bI & \bzero \\
\bzero & \bzero & \bI
\end{bmatrix}
\begin{bmatrix}
0 & \bc_1^\top & \bc_2^\top \\
\bzero & \bJ_{n_1} &  \bzero \\
\bzero & \bzero & \bJ
\end{bmatrix}
\begin{bmatrix}
1 & \bc_1^\top \bJ_{n_1}^\top & \bzero \\
\bzero & \bI & \bzero \\
\bzero & \bzero & \bI
\end{bmatrix}
=
\begin{bmatrix}
0 & (\bc_1^\top\be_1) \be_1^\top  & \bc_2^\top \\
\bzero & \bJ_{n_1} & \bzero \\
\bzero & \bzero & \bJ
\end{bmatrix}.
\end{equation}
\paragraph{Case 1: $\gamma\triangleq\bc_1^\top\be_1\neq 0$. }
Then, a similarity transformation on \eqref{equation:jd_q_s1} can be constructed: 
$$
\begin{bmatrix}
\frac{1}{\gamma}  & \bzero   & \bzero  \\
\bzero & \bI & \bzero \\
\bzero & \bzero & \frac{1}{\gamma}\bI
\end{bmatrix}
\begin{bmatrix}
0 & \gamma \be_1^\top  & \bc_2^\top \\
\bzero & \bJ_{n_1} & \bzero \\
\bzero & \bzero & \bJ
\end{bmatrix}
\begin{bmatrix}
\gamma & \bzero   & \bzero  \\
\bzero & \bI & \bzero \\
\bzero & \bzero &  \gamma\bI
\end{bmatrix}
\triangleq
\begin{bmatrix}
0 & \be_1^\top & \bc_2^\top\\
\bzero & \bJ_{n_1} & \bzero \\
\bzero & \bzero & \bJ
\end{bmatrix}=
\begin{bmatrix}
\widehat{\bJ} & \be_1\bc_2^\top \\
\bzero &\bJ 
\end{bmatrix},
$$
where 
$\bJ_{n_1+1}(0)\triangleq\widehat{\bJ}=\footnotesize\begin{bmatrix}
0 & \be_1^\top \\
\bzero & \bJ_{n_1} 
\end{bmatrix}$ is a Jordan block of size $n_1+1$ with $\lambda=0$
and 
$ \be_1\bc_2^\top=\footnotesize\begin{bmatrix}
\bc_2^\top \\
\bzero 
\end{bmatrix}$ ($\be_1$ is the first standard basis vector in $\real^{n_1+1}$).
We note by Exercise~\ref{exercise:nil_jor} that $\widehat{\bJ}\be_{i+1}=\be_i$ for $i\in\{1,2, \ldots, n_1\}$.
Therefore, we can proceed by another similarity transformation to obtain the desired strictly upper triangular Jordan form:
$$
\begin{bmatrix}
\bI& \be_2\bc_2^\top\\
\bzero & \bI 
\end{bmatrix}
\begin{bmatrix}
\widehat{\bJ} & \be_1\bc_2^\top \\
\bzero &\bJ 
\end{bmatrix}
\begin{bmatrix}
\bI& -\be_2\bc_2^\top\\
\bzero & \bI 
\end{bmatrix}
=
\begin{bmatrix}
\widehat{\bJ} & \be_2\bc_2^\top\bJ \\
\bzero & \bJ 
\end{bmatrix},
$$
and for $i\in\{2,3,\ldots, n_1\}$:
$$
\implies 
\begin{bmatrix}
\bI& \be_{i+1}\bc_2^\top\bJ^{i-1}\\
\bzero & \bI 
\end{bmatrix}
\begin{bmatrix}
\widehat{\bJ} & \be_i\bc_2^\top\bJ^{i-1} \\
\bzero &\bJ 
\end{bmatrix}
\begin{bmatrix}
\bI& -\be_{i+1}\bc_2^\top\bJ^{i-1}\\
\bzero & \bI 
\end{bmatrix}
=
\begin{bmatrix}
\widehat{\bJ} & \be_{i+1}\bc_2^\top\bJ^{i} \\
\bzero & \bJ 
\end{bmatrix}.
$$
Exercise~\ref{exercise:nil_jor} show that  $\bJ^{k}=\bzero$ for $k\geq n_1$. After at most $n_1$ steps, we can obtain  the desired strictly upper triangular Jordan form: $\bA$ is similar to
$$
\begin{bmatrix}
\widehat{\bJ} & \bzero \\
\bzero & \bJ 
\end{bmatrix}.
$$

\paragraph{Case 2: $\gamma\triangleq\bc_1^\top\be_1= 0$. }
We can apply a permutation similarity transformation such that 
$$
\bP^{-1}
\begin{bmatrix}
0 & \bzero   & \bc_2^\top \\
\bzero & \bJ_{n_1} & \bzero \\
\bzero & \bzero & \bJ
\end{bmatrix}
\bP
=
\begin{bmatrix}
\bJ_{n_1} & \bzero   & \bzero  \\
\bzero & 0 & \bc_2^\top  \\
\bzero & \bzero & \bJ
\end{bmatrix}.
$$
The bottom right matrix $\footnotesize\begin{bmatrix}
0 & \bc_2^\top  \\
\bzero & \bJ
\end{bmatrix}$ is similar to a strictly upper triangular Jordan form by the induction hypothesis. Therefore, $\bA$ is also similar to a strictly upper triangular Jordan form. 
This completes the proof.
\end{proof}

A direct consequence of Lemma~\ref{lemma:upptri_jordan} is that  the upper triangular matrix with all diagonals equal to $\lambda$ also admits this similarity transformation that transforms the matrix into a Jordan form, with all diagonals equal to $\lambda$ as well. 
\begin{lemma}[Jordan Form of Upper Triangular Matrices]\label{lemma:jd_spup}
Let $\bA\in\real^{n\times n}$ be any square  upper triangular \textit{with all diagonals equal to $\lambda$}. Then, it can be factored as
$$
\bA=\bX\bJ\bX^{-1},~\footnote{When $\bA\in\complex^{n\times n}$ is complex, $\bX$ can be chosen complex as well.}
$$
where $\bX\in\real^{n\times n}$ is nonsingular,  $\bJ=\diag\big(\bJ_{m_1}(\lambda), \bJ_{m_2}(\lambda), \ldots, \bJ_{m_k}(\lambda)\big)$ is  upper triangular \textit{with all diagonals equal to $\lambda$}, $m_1\geq m_2\geq \ldots m_k\geq 1$, and $m_1+m_2+\ldots+m_k=n$.
\end{lemma}
The proof is straightforward. Let $\bA_0 \triangleq sut(\bA)$ be the strictly upper triangular part of $\bA$. Then, $\bA_0=\bX\diag\big(\bJ_{m_1}(0), \bJ_{m_2}(0), \ldots, \bJ_{m_k}(0)\big)\bX^{-1}$ and $\lambda\bI=\lambda\bX\bX^{-1}$. 
Therefore, 
$$
\bA=\bX\diag\big(\bJ_{m_1}(\lambda), \bJ_{m_2}(\lambda), \ldots, \bJ_{m_k}(\lambda)\big)\bX^{-1}.
$$

Therefore, if we can apply a similarity transformation to a general matrix $\bA$ to convert it into a block upper triangular matrix, where each block has the same diagonal values. Then, applying Lemma~\ref{lemma:jd_spup} to each of these upper triangular blocks will yield  the Jordan decomposition (Theorem~\ref{theorem:jordan-decomposition}). 
Theorem~\ref{theorem:block_real_eigens} provides such a block upper triangularization, which is a consequence of the Schur decomposition.

\section{Application: Computing Fibonacci Numbers}\label{section:fibonacci}
Eigenvalue decomposition can be used to compute Fibonacci numbers. 
This particular example is derived and adapted from \citet{strang1993introduction}. 
Each subsequent Fibonacci number, denoted by $F_{k+2}$, is obtained by adding the two preceding Fibonacci numbers, namely $F_{k+1}+F_{k}$. The sequence begins as $0, 1, 1, 2, 3, 5, 8, \ldots$. 
A natural question arises: what is the value of $F_{100}$?
Eigenvalue decomposition facilitates the derivation of the general formula for the sequence.
\index{Fibonacci number}
\index{General formula of a sequence}

Let, $\bu_{k}\triangleq\footnotesize\begin{bmatrix}
F_{k+1}\\
F_k
\end{bmatrix}$. 
Then it follows that $\bu_{k+1}=\footnotesize\begin{bmatrix}
	F_{k+2}\\
	F_{k+1}
\end{bmatrix}=
\footnotesize\begin{bmatrix}
	1&1\\
	1&0
\end{bmatrix}
\bu_k
$ by the rule that $F_{k+2}=F_{k+1}+F_k$ and $F_{k+1}=F_{k+1}$. 
Let $\bA\triangleq\footnotesize\begin{bmatrix}
	1&1\\
	1&0
\end{bmatrix}$. 
We then have the general formula $\bu_{100} = \bA^{100}\bu_0$, where $\bu_0=\footnotesize
\begin{bmatrix}
	1\\
	0
\end{bmatrix}
$.

We have demonstrated in Remark~\ref{remark:determinant-intermezzo} that $\det(\bA-\lambda\bI)=0$ when $\lambda$ is an eigenvalue of $\bA$. Simple calculations reveals that $\det(\bA-\lambda\bI) = \lambda^2-\lambda+1=0$, yielding the eigenvalues and the corresponding eigenvectors as follows: 
$$
(\lambda_1, \bx_1) = \left(\frac{1+\sqrt{5}}{2},  
\,\,\,\,
\begin{bmatrix}
\lambda_1\\
1
\end{bmatrix}\right)
\quad\text{and}\quad
(\lambda_2, \bx_2)
= 
\left(\frac{1-\sqrt{5}}{2}, 
\,\,\,\,
\begin{bmatrix}
\lambda_2\\
1
\end{bmatrix}\right).
$$
As per Remark~\ref{remark:power-eigenvalue-decom},
we can express $\bA^{100} = \bX\bLambda^{100}\bX^{-1} = \bX
\footnotesize
\begin{bmatrix}
	\lambda_1^{100}&0\\
	0&\lambda_2^{100}
\end{bmatrix}\bX^{-1}$, where $\bX^{-1}$ can be easily calculated as $\bX^{-1} =
\footnotesize
\begin{bmatrix}
	\frac{1}{\lambda_1-\lambda_2} & \frac{-\lambda_2}{\lambda_1-\lambda_2} \\
	-\frac{1}{\lambda_1-\lambda_2} & \frac{\lambda_1}{\lambda_1-\lambda_2}
\end{bmatrix} 
=
\footnotesize
\begin{bmatrix}
\frac{\sqrt{5}}{5} & \frac{5-\sqrt{5}}{10} \\
-\frac{\sqrt{5}}{5} & \frac{5+\sqrt{5}}{10}
\end{bmatrix}$. We notice that $\bu_{100} = \bA^{100}\bu_0$ corresponds to  the first column of $\bA^{100}$, which can be represented as:
 $$
\bu_{100} = 
\begin{bmatrix}
F_{101}\\
F_{100}
\end{bmatrix}=
\begin{bmatrix}
	\frac{\lambda_1^{101}-\lambda_2^{101}}{\lambda_1-\lambda_2}\\
	\frac{\lambda_1^{100}-\lambda_2^{100}}{\lambda_1-\lambda_2}
\end{bmatrix}.
$$
Upon a simple check of the calculation,  we have $F_{100}=3.542248481792631e+20$. Or more generally, we can express $\bu_{K}$ as follows:
$$
\bu_{K} = 
\begin{bmatrix}
	F_{K+1}\\
	F_{K}
\end{bmatrix}=
\begin{bmatrix}
	\frac{\lambda_1^{K+1}-\lambda_2^{K+1}}{\lambda_1-\lambda_2}\\
	\frac{\lambda_1^{K}-\lambda_2^{K}}{\lambda_1-\lambda_2}
\end{bmatrix},
$$
where the general form of $F_K$ is given by $F_K=\frac{\lambda_1^{K}-\lambda_2^{K}}{\lambda_1-\lambda_2}$.

\paragraph{Extension to general 2 by 2 matrices.}
The  matrix $\footnotesize\begin{bmatrix}
	1 & 1 \\ 1 & 0
\end{bmatrix}$ presented above is special in that it has two linearly independent eigenvectors so that it admits eigenvalue decomposition. However, it is also interesting to compute the $n$-th power of general 2 by 2 matrices, say
$$
\bA=
\begin{bmatrix}
	a & b \\ c & d
\end{bmatrix}.
$$
Suppose $\alpha$ and $\beta$ are eigenvalues of $\bA$, which are two roots of the characteristic polynomial 
$$
\det(\lambda\bI-\bA)= \lambda^2 + a_1 \lambda + a_0 = 0
\quad\implies\quad 
 \alpha, \beta = \frac{-a_1 \pm \sqrt{a_1^2 - 4a_0}}{2}.
$$
This reveals that $a_1 = -(\alpha+\beta)$ and $a_0=\alpha\beta$. According to the \textit{Cayley-Hamilton theorem} (Theorem~\ref{theorem:cayley_hami}), a matrix satisfies its own characteristic equation:
$
\bA^2 + a_1 \bA + a_0 \bI = \bzero,
$
so that 
\begin{equation}\label{eqiation:2by2-cayley}
\bA^2 - (\alpha+\beta)\bA + \alpha\beta\bI = \bzero.
\end{equation}
Following the tricks in \citet{Williams1992npower}, we define matrices $\bX,\bY,\bZ$ such that 
$$
\left\{
\begin{aligned}
\bX &\triangleq \frac{\bA-\beta\bI}{\alpha-\beta}, \gap \bY \triangleq \frac{\bA-\alpha\bI}{\beta-\alpha}, \gap &\text{if $\alpha\neq \beta$;}\\
\bZ &\triangleq \bA-\alpha\bI, &\text{if $\alpha=\beta$,}
\end{aligned}
\right.
$$
where $\bA = \alpha\bX + \beta\bY$. From Equation~\eqref{eqiation:2by2-cayley}, for any $k\geq 2$, we have:
$$
\left\{
\begin{aligned}
	\bX^k &= \bX, \gap \bX\bY=\bY\bX=\bzero, \gap \bY^k=\bY, \gap  &\text{if $\alpha\neq \beta$;}\\
	\bZ^k &= \bzero, &\text{if $\alpha=\beta$.}
\end{aligned}
\right.
$$
Therefore, the $n$-th power of $\bA$ can be expressed as:
\begin{equation}\label{equation:genera2by2power}
\bA^n=
\left\{
\begin{aligned}
	(\alpha\bX+\beta\bY)^n &= \alpha^n\bX+ \beta^n\bY= \alpha^n \left(\frac{\bA-\beta\bI}{\alpha-\beta}\right)+\beta^n\left(\frac{\bA-\alpha\bI}{\beta-\alpha}\right), \gap  &\text{if $\alpha\neq \beta$;}\\
	(\bZ+\alpha\bI)^n &= n\alpha^{n-1}\bZ + \alpha^n\bI=\alpha^{n-1}\left( n\bA - \alpha(n-1)\bI \right), &\text{if $\alpha=\beta$.}
\end{aligned}
\right.
\end{equation}
Hence, if matrix $\bA$ is nonsingular, $\alpha\neq0$, and $\beta\neq 0$, Equation~\eqref{equation:genera2by2power} holds for all integer values of $n$; if matrix $\bA$ is real, but the eigenvalues are complex with some power of them being real, e.g., $\alpha^m=\beta^m=f$ is real, then we have $\bA^m = f\bI$.  
The $n$-th power for general matrices results in the general form of any sequence, e.g., 
$$
\bu_{k+1} = 
\begin{bmatrix}
2 & 3 \\2 & 0
\end{bmatrix}
\bu_k.
$$
And we shall not repeat the details.

\index{Matrix polynomial}
\section{Application: Matrix Polynomials}
We have shown in Problems~\ref{prob:hess_poly1} and \ref{prob:hess_poly2} that $f(\bP\bA\bP^{-1}) =\bP f(\bA) \bP^{-1}$ if $f(\bC) = \gamma_m\bC^m+\gamma_{m-1}\bC^{m-1}+\ldots+\gamma_0$ is a polynomial.
Let 
$
\bA = \bX\bJ\bX^{-1}
$
be the Jordan decomposition of $\bA\in\real^{n\times n}$,
where $\bX\in\real^{n\times n}$ is a nonsingular matrix containing the \textit{generalized eigenvectors} of $\bA$ as its columns, and $\bJ\in\real^{n\times n}$ is a Jordan form matrix $\diag(\bJ_{m_1}(\lambda_1), \bJ_{m_2}(\lambda_2), \ldots, \bJ_{m_k}(\lambda_k)) \triangleq\diag(\bJ_1, \bJ_{2}, \ldots, \bJ_{k})$, where $\bJ_{m_i}(\lambda_i)\in\real^{m_i\times m_i}$ and $\sum_{i=1}^{k}m_i  = n$.
Therefore, it follows that 
\begin{equation}
f(\bA) = \bX f(\bJ) \bX^{-1}
=
\bX \diag(f(\bJ_1), f(\bJ_2), \ldots, f(\bJ_k)) \bX^{-1},
\end{equation}
where 
\begin{equation}
f(\bJ_1) = 
\begin{bmatrix}
f(\lambda_i) & f'(\lambda_i) & \frac{1}{2!} f''(\lambda_i) &  \ldots & \frac{1}{(m_i-1)!} f^{(m_i-1)}(\lambda_i)\\
0  & f(\lambda_i) & f'(\lambda_i) &  \ldots & \frac{1}{(m_i-2)!} f^{(m_i-2)}(\lambda_i)\\
0  & 0 & \ddots &  \ddots & \vdots\\
0  & 0 & 0 &  f(\lambda_i) &f'(\lambda_i)\\
0  & 0 & 0 &  \ldots & f(\lambda_i)\\
\end{bmatrix},
\end{equation}
and $f^{(k)}(x)$ denotes the $k$-th  derivative of $f(x)$.
Therefore, it is easy to extend this to common matrix functions that can be expressed as polynomials as well \citep{zhang2017matrix}:

\paragraph{Powers of a matrix.}
$\bA^m\triangleq  \bX\bJ^m\bX^{-1} = \bX f(\bJ)\bX^{-1}$, where $f(x) = x^m$ for all $m=1,2,\ldots$.

\paragraph{Matrix logarithm.} Let $f(x) = \ln(1 + x)$. Then,
\begin{equation}
\ln(\bI + \bA) \triangleq \sum_{i=1}^{\infty} \frac{(-1)^{i-1}}{i} \bA^i = \bX \left( \sum_{i=1}^{\infty} \frac{(-1)^{i-1}}{i} \bA^i \right) \bX^{-1} = \bX f(\bJ)\bX^{-1}.
\end{equation}

\paragraph{Sine and cosine functions.} Let  $f_1(x) = \sin(x)$ and $f_2(x) = \cos (x)$. Then,
\begin{align}
\sin(\bA) &\triangleq \sum_{i=0}^{\infty} \frac{(-1)^i}{(2i+1)!} \bA^{2i+1} = \bX \left( \sum_{i=0}^{\infty} \frac{(-1)^i}{(2i+1)!} \bJ^{2i+1} \right) \bX^{-1} = \bX f_1(\bJ)\bX^{-1};\\
\cos (\bA)&\triangleq \sum_{i=0}^{\infty} \frac{(-1)^i}{(2i)!} \bA^{2i} = \bX \left( \sum_{i=0}^{\infty} \frac{(-1)^i}{(2i)!} \bJ^{2i} \right) \bX^{-1} = \bX f_2(\bJ)\bX^{-1}.
\end{align}

\paragraph{Matrix exponentials.} Let $f_1(x) = e^x$ and $f_2(x) = e^{-x}$. Then,
\begin{align}
	e^{\bA} &\triangleq \sum_{i=0}^{\infty} \frac{1}{i!} \bA^i = \bX \left( \sum_{i=0}^{\infty} \frac{1}{i!} \bJ^i \right) \bX^{-1} = \bX f_1(\bJ)\bX^{-1}; \label{equation:max_expo}\\
	e^{-\bA} &\triangleq \sum_{i=0}^{\infty} \frac{1}{i!} (-1)^i \bA^i = \bX \left( \sum_{i=0}^{\infty} \frac{1}{i!} (-1)^i \bJ^i \right) \bX^{-1} = \bX f_2(\bJ)\bX^{-1}.
\end{align}

\paragraph{Matrix exponential functions.} Let $f_1(x) = e^{xt}$ and $f_2(x) = e^{-xt}$. Then,
\begin{align}
e^{\bA t} &\triangleq \sum_{i=0}^{\infty} \frac{1}{i!} \bA^i t^i = \bX \left( \sum_{i=0}^{\infty} \frac{1}{i!} \bJ^i t^i \right) \bX^{-1} = \bX f_1(\bJ)\bX^{-1};\label{equation:max_expo3}\\
e^{-\bA t} &\triangleq \sum_{i=0}^{\infty} \frac{1}{i!} (-1)^i \bA^i t^i = \bX \left( \sum_{i=0}^{\infty} \frac{1}{i!} (-1)^i \bJ^i t^i \right) \bX^{-1} = \bX f_2(\bJ)\bX^{-1}.
\end{align}

\begin{problemset}
\item \textbf{Polynomial of upper triangular.} Let $\bR=[r_{ij}]_{i,j=1}^{n,n}\in\real^{n\times n}$ be upper triangular. Show that its characteristic polynomial is $p_{\bR}(\lambda)=(\lambda-r_{11})(\lambda-r_{22})\ldots(\lambda-r_{nn})$. This concludes that the diagonal elements of an upper triangular matrix are its eigenvalues.
Moreover, let $\bA=[\bA_{ij}]_{i,j=1}^{n,n}$ be a upper triangular block matrix with $\bA_{ij}=\bzero$ if $i>j$. Show that its characteristic polynomial is $p_{\bA}(\lambda) = \prod_{i=1}^{n}p_{\bA_{ii}}(\lambda)$.

\item Show that $\bA,\bB\in\real^{n\times n}$ have the same characteristic polynomials (the same eigenvalues) if and only if $\trace(\bA^k)=\trace(\bB^k)$ for all $k=\{1,2,\ldots,n\}$.

\item Show that a matrix $\bA\in\real^{n\times n}$ and its transpose $\bA^\top$ have the same eigenvalues. \textit{Hint: Use  the characteristic polynomial of $\bA$; or consider Equation~\eqref{equation:simi_trans}.}

\item Let  $\lambda$ be a  complex eigenvalue of $\bA$. Use the characteristic polynomial to show that its conjugate $\overline{\lambda}$ is also an eigenvalue of $\bA$ with the same algebraic multiplicity.

\item \label{prob:consmatfrom_pol} We have shown that each $n$-th order square matrix has its $n$-th order polynomial in Definition~\ref{definition:characteristic_polynomial}. On the other hand, each  $n$-th order polynomial can also be written as the characteristic polynomial of an $n \times  n$ matrix.
Given the polynomial $p(\lambda) =  \lambda^n + \gamma_{n-1} \lambda^{n-1} + \ldots + \gamma_1 \lambda  + \gamma_0$, show that $p(\lambda) = \det(\lambda\bI-\bA)$, where 
$$
\bA=
\footnotesize
\begin{bmatrix}
-\gamma_{n-1} & -\gamma_{n-2} & \ldots & -\gamma_1 & \gamma_0 \\
1 & 0 & \ldots & 0 & 0 \\
0 & 1 & \ldots & 0 & 0 \\
\vdots  & \vdots & \ddots & \vdots & \vdots \\
0 & 0 & \ldots & 1 & 0 \\
\end{bmatrix}.
$$
	
\item \label{prob:speci_eign} Consider the  matrix $\bA=\bone\bone^\top\in\real^{n\times n}$, where all entries equal to 1, and $\bone\in\real^n$ is the vector of all ones. Find $n$ linearly independent eigenvectors of $\bA$, and determine the corresponding eigenvalues. \textit{Hint: Consider $\bx_i=\bone-n\be_i$ and $\bone$.} 

\item  What are the eigenvalues of the matrix $\bA=\scriptsize\begin{bmatrix}
5 & -1 & -1\\
-1 & 5 & -1\\
-1& -1 & 5
\end{bmatrix}$?
\textit{Hint: $\bA=6\bI-\bone\bone^\top$ and use Problem~\ref{prob:speci_eign}.}

\item Given an idempotent matrix $\bA\in\real^{n\times n}$ (Definition~\ref{definition:Idempotent_mat}), show that $\bB\bA$ and $\bA\bB\bA$ share the same eigenvalues.

\item Given a Householder transformation matrix $\bH = \bI - 2\bu\bu^\top\in\real^{n\times n}$, where $\normtwo{\bu}=1$, show that $\bu$ is an eigenvector of $\bH$ and determine its corresponding eigenvalue. Provide a geometric interpretation of the eigenvalues of $\bH$.
Suppose further that $\bv^\top\bu=0$, where $\bv$ is a nonzero vector. Show that $\bv$ is also an eigenvector of $\bH$ and find its corresponding eigenvalue.

\item Suppose $\lambda$ is an eigenvalue of $\bA\in\real^{n\times n}$. Show that $\lambda-\mu$ is an eigenvalue of $\bA-\mu\bI$.

\item Suppose $\lambda$ is an eigenvalue of a nonsingular matrix $\bA\in\real^{n\times n}$. Show that $\lambda^{-1}$ is an eigenvalue of $\bA^{-1}$.

\item \label{prob:geneig1} \textbf{Generalized eigenproblem.} Many scientific packages also deals with the generalized problem $\bA\bx=\lambda\bB\bx$, where $\bB$ is nonsingular. If $\bA$ is symmetric and $\bB$ is PD with Cholesky decomposition $\bB=\bR^\top\bR$, show that the eigenvalue $\lambda$ is a (standard) eigenvalue of $\bC=(\bR^{-1})^\top \bA\bR^{-1}$, corresponding to the eigenvector $\bR\bx$.

\item \label{prob:geneig2} \textbf{Generalized eigenproblem \citep{teukolsky1992numerical}.} Suppose $\bA\lambda^2 +\bB\lambda +\bC=\bzero$. Show that $\lambda$  can be solved by a standard eigenproblem. \textit{Hint: Let $\by=\lambda\bx$, and consider the matrix $\footnotesize\begin{bmatrix}
	\bzero & \bI\\
	-\bA^{-1}\bC & -\bA^{-1}\bB 
	\end{bmatrix}$.}

\item Let $\bA\in\real^{n\times n}$ and suppose $\lambda$ is not an eigenvalue of $\bA$. Show that $\rank((\bA-\lambda\bI)^{k})=\rank((\bA-\lambda\bI)^{k-1})$ for all $k\geq  2$. \textit{Hint: Use Jordan decomposition and Exercise~\ref{exercise:sim_trans_inhess}.}

\item Given an integer $n\geq 1$, show that 
\begin{itemize}
\item If $m=2n$, then the Jordan canonical form of $\bJ_m(0)^2$ is $\diag(\bJ_n(0), \bJ_n(0))$.
\item If $m=2n+1$, then the Jordan canonical form of $\bJ_m(0)^2$ is $\diag(\bJ_{n+1}(0), \bJ_n(0))$. 
\end{itemize}

\item \textbf{Matrix exponentials.} Given the definition of matrix exponentials in \eqref{equation:max_expo}, let  $\bA$ and $\bB$ commute, i.e., $\bA\bB=\bB\bA$. Show that $e^{\bA+\bB} = e^{\bA}\cdot e^{\bB}$.

\item \label{prob:mtexp2} \textbf{Matrix exponentials.} Let $f(t):\real\rightarrow \real^n$ be a function satisfying $f'(t)=\bA f(t)$ with $f(0)=\bx_0\in\real^n$ and $\bA\in\real^{n\times n}$. Show that the unique solution is $f(t)=e^{\bA t}\bx_0$.

\item \textbf{Matrix exponentials.} Consider the matrix exponential function in \eqref{equation:max_expo3}. Let $\bX^{-1}\bA\bX=\bJ=\diag(\bJ_1,\bJ_2,\ldots,\bJ_k)$ be the Jordan form of $\bA$, where $\bJ_i\in\real^{m_i\times m_i}$ contains the eigenvalue $\lambda_i$ on the diagonal. Show that 
$$
e^{\bJ_i t} =
e^{\lambda_i\cdot t}
\begin{bmatrix}
1 & t  & \frac{t^2}{2!} & \ldots & \frac{t^{m_i-1}}{(m_i-1)!} \\
0 &  1 & t & \ldots & \frac{t^{m_i-2}}{(m_i-2)!} \\
\vdots & \vdots & \ddots &\vdots  & \vdots \\
0 &  0 & 0 & \ddots & t \\
0 &  0 & 0 & \ldots & 1 \\
\end{bmatrix} .
$$  
Rewrite the matrix equation in Problem~\ref{prob:mtexp2} under the form of $g'(t)=\bJ g(t)$ and determine the explicit form of $g(t)$.
\textit{Hint: Decompose every Jordan block as a sum of a diagonal matrix and a nilpotent matrix.}

\end{problemset}

%% file: chapter-schur.tex
\newpage
\chapter{Schur Decomposition}
\begingroup
\hypersetup{
linkcolor=structurecolor,
linktoc=page,  
}
\minitoc 
\newpage
\endgroup

\section{Schur Decomposition}
\index{Decomposition: Schur}
\lettrine{\color{caligraphcolor}T}
The Schur decomposition uses an orthogonal similarity transformation to transform an arbitrary square matrix into an upper triangular matrix. 
This allows many properties of the original matrix to be explored using the simpler form of the upper triangular matrix.

\begin{theoremHigh}[Schur Decomposition]\label{theorem:schur-decomposition}
Let $\bA\in \real^{n\times n}$ be any \textbf{real} square matrix with \textbf{real} eigenvalues. Then, it can be factored as~\footnote{The complex Schur decomposition is discussed in Theorem~\ref{theorem:schur-decomposition_complex}.}
$$
\bA = \bQ\bU\bQ^\top,
$$
where $\bQ$ is a (real) orthogonal matrix, and $\bU$ is a (real) upper triangular matrix. 
In other words, any real square matrix $\bA$ with real eigenvalues can be triangularized.
\end{theoremHigh}

 The first column of $\bA\bQ$ and $\bQ\bU$ are $\bA\bq_1$ and  $u_{11}\bq_1$, respectively.
Consequently,  $u_{11}$corresponds to an eigenvalue of $\bA$, while  $\bq_1$ serves as its associated eigenvector. 
However, the other columns of $\bQ$ may not necessarily be eigenvectors of $\bA$. 

\paragraph{Schur decomposition for symmetric matrices.} When dealing with a symmetric matrix $\bA=\bA^\top$, we find that $\bQ\bU\bQ^\top = \bQ\bU^\top\bQ^\top$. 
Consequently, $\bU$ must be a diagonal matrix, and this diagonal matrix  contains the eigenvalues of $\bA$. Furthermore, all the columns of $\bQ$ are eigenvectors of $\bA$. 
This leads us to the conclusion that symmetric matrices are inherently diagonalizable, even in the presence of repeated eigenvalues.

\section{Existence of  Schur Decomposition}
In order to establish the validity of Theorem~\ref{theorem:schur-decomposition}, we must utilize the following lemma.
\begin{lemma}[Submatrix with Same Eigenvalues]\label{lemma:submatrix-same-eigenvalue}
Let $\bA_{k+1}\in \real^{(k+1)\times (k+1)}$ be  any square matrixhas with real eigenvalues $\lambda_1, \lambda_2, \ldots, \lambda_{k+1}$. Then, we can construct a $k\times k$ matrix $\bA_{k}$ with eigenvalues $\lambda_2, \lambda_3, \ldots, \lambda_{k+1}$ as follows: 
$$
\bA_{k} = 
\begin{bmatrix}
-\bp_2^\top- \\
-\bp_3^\top- \\
\vdots \\
-\bp_{k+1}^\top-
\end{bmatrix}
\bA_{k+1}
\begin{bmatrix}
\bp_2 & \bp_3 &\ldots &\bp_{k+1}
\end{bmatrix},
$$
where $\bp_1$ is an eigenvector of $\bA_{k+1}$ with a norm of 1 corresponding to the eigenvalue $\lambda_1$, and $\bp_2, \bp_3, \ldots, \bp_{k+1}$ denote any mutually orthonormal vectors orthogonal to $\bp_1$. 
\end{lemma}

\begin{proof}[of Lemma~\ref{lemma:submatrix-same-eigenvalue}]
Let $\bP_{k+1} \triangleq [\bp_1, \bp_2, \ldots, \bp_{k+1}]$. It follows that $\bP_{k+1}^\top\bP_{k+1}=\bI$, and 
$$
\bP_{k+1}^\top \bA_{k+1} \bP_{k+1} =
\begin{bmatrix}
\lambda_1 & \bzero \\
\bzero  & \bA_{k}
\end{bmatrix}. 
$$
For any eigenvalue $\lambda\in\{\lambda_2, \lambda_3, \ldots, \lambda_{k+1}\}$, as per Remark~\ref{remark:determinant-intermezzo}, we obtain
$$
\begin{aligned}
\det(\bA_{k+1} -\lambda\bI) 
&= \det(\bP_{k+1}^\top (\bA_{k+1}-\lambda\bI)  \bP_{k+1}) 
=\det(\bP_{k+1}^\top \bA_{k+1}\bP_{k+1} - \lambda\bP_{k+1}^\top\bP_{k+1}) \\
&= \det\left(
\begin{bmatrix}
\lambda_1-\lambda & \bzero \\
\bzero & \bA_k - \lambda\bI
\end{bmatrix}
\right)
=(\lambda_1-\lambda)\det(\bA_k-\lambda\bI).
\end{aligned}
$$
The last equality follows from the determinant of a block matrix (Remark~\ref{remark:determinant-intermezzo}).
Since $\lambda$ is an eigenvalue of $\bA$ and $\lambda \neq \lambda_1$ (assume $\lambda_1\neq \lambda_2,\ldots,\lambda_{k+1}$), then $\det(\bA_{k+1} -\lambda\bI) = (\lambda_1-\lambda)\det(\bA_{k}-\lambda\bI)=0$ means $\lambda$ is also an eigenvalue of $\bA_{k}$.
\end{proof}

We proceed to establish the existence of the Schur decomposition through an inductive proof.
\begin{proof}[of Theorem~\ref{theorem:schur-decomposition}: Existence of Schur Decomposition]
We note that the theorem holds trivially when $n=1$ by simply setting $Q=1$ and $U=A$. Suppose the theorem holds true for $n=k$ given $k> 1$. Now, if we can demonstrate that the theorem also holds true for $n=k+1$, then we complete the proof.
Assume that for $n=k$, the theorem is valid for $\bA_k =\bQ_k \bU_k \bQ_k^\top$. 

Suppose further that \textit{$\bP_{k+1}$ contains orthogonal vectors $\bP_{k+1} = [\bp_1, \bp_2, \ldots, \bp_{k+1}]$ as constructed in Lemma~\ref{lemma:submatrix-same-eigenvalue}, where $\bp_1$ is an eigenvector of $\bA_{k+1}\in\real^{(k+1)\times(k+1)}$ corresponding to eigenvalue $\lambda_1$ and its norm is 1;  and $\bp_2, \ldots, \bp_{k+1}$ are orthonormal to $\bp_1$}. Let the other $k$ eigenvalues of $\bA_{k+1}$ be $\lambda_2, \lambda_3, \ldots, \lambda_{k+1}$. 
Since we assume  the theorem is true for $n=k$, we can find a matrix $\bA_{k}$ with eigenvalues $\lambda_2, \lambda_3, \ldots, \lambda_{k+1}$ and it admits the Schur decomposition $\bA_k =\bQ_k \bU_k \bQ_k^\top$. 
Thus, we can deduce the following property based on Lemma~\ref{lemma:submatrix-same-eigenvalue}:
$$
\bP_{k+1}^\top \bA_{k+1} \bP_{k+1} = 
\begin{bmatrix}
\lambda_1 &\bzero \\
\bzero & \bA_k
\end{bmatrix} 
\qquad\implies \qquad 
\bA_{k+1} \bP_{k+1} =
\bP_{k+1} 
\begin{bmatrix}
\lambda_1 &\bzero \\
\bzero & \bA_k
\end{bmatrix}.
$$
Let 
$
\bQ_{k+1} \triangleq
\bP_{k+1}
$
$
\footnotesize
\begin{bmatrix}
1 &\bzero \\
\bzero & \bQ_k
\end{bmatrix}.
$
Then,  it follows that
$$
\begin{aligned}
\bA_{k+1} \bQ_{k+1} &= 
\bA_{k+1}
\bP_{k+1}
\begin{bmatrix}
1 &\bzero \\
\bzero & \bQ_k
\end{bmatrix}
=
\bP_{k+1} 
\begin{bmatrix}
\lambda_1 &\bzero \\
\bzero & \bA_k
\end{bmatrix}
\begin{bmatrix}
1 &\bzero \\
\bzero & \bQ_k
\end{bmatrix} 
=
\bP_{k+1}
\begin{bmatrix}
\lambda_1 & \bzero \\
\bzero & \bA_k\bQ_k
\end{bmatrix}\\
&=
\bP_{k+1}
\begin{bmatrix}
\lambda_1 & \bzero \\
\bzero & \bQ_k \bU_k   
\end{bmatrix} 
=
\underbrace{\bP_{k+1}
\begin{bmatrix}
	1 &\bzero \\
	\bzero & \bQ_k
\end{bmatrix}}_{=\bQ_{k+1}}
\underbrace{\begin{bmatrix}
	\lambda_1 &\bzero \\
	\bzero & \bU_k
\end{bmatrix}}_{\triangleq\bU_{k+1}}  
=\bQ_{k+1}\bU_{k+1}.
\end{aligned}
$$
We then have $\bA_{k+1} = \bQ_{k+1}\bU_{k+1}\bQ_{k+1}^\top$, where $\bU_{k+1}$ is an upper triangular matrix, and $\bQ_{k+1}$ is an orthogonal matrix since $\bP_{k+1}$ and 
$\scriptsize\begin{bmatrix}
1 &\bzero \\
\bzero & \bQ_k
\end{bmatrix}$ are both orthogonal matrices.
\end{proof}

\section{Other Forms of  Schur Decomposition}\label{section:other-form-schur-decom}
In the proof of the Schur decomposition, we obtain the upper triangular matrix $\bU_{k+1}$ by appending the eigenvalue $\lambda_1$ to $\bU_k$. As a result, the values on the diagonal are always eigenvalues. Therefore, we can decompose the upper triangular matrix into two parts.  

\begin{corollaryHigh}[Form 2 of Schur Decomposition]\label{corollary:schur-second-form}
Let $\bA\in \real^{n\times n}$ be any real square matrix with real eigenvalues. Then, it can be factored as 
$$
\bQ^\top\bA\bQ = \bLambda +\bT \qquad \textbf{or} \qquad \bA = \bQ(\bLambda +\bT)\bQ^\top,
$$
where $\bQ$ is an orthogonal matrix, $\bLambda=\diag(\lambda_1, \lambda_2, \ldots, \lambda_n)$ is a diagonal matrix containing the eigenvalues of $\bA$ (the eigenvalues can be arranged in any prescribed order), and $\bT$ is a \textit{strictly upper triangular} matrix (with zeros on the diagonal).
\end{corollaryHigh}
A strictly upper triangular matrix is an upper triangular matrix having 0's along the diagonal as well as the lower portion. Another proof of this decomposition is that $\bA$ and $\bU$ (where $\bU = \bQ^\top\bA\bQ$) are similar matrices, implying they share the same eigenvalues (Proposition~\ref{proposition:eigenvalue-similar-matrices}).
Furthermore, the eigenvalues of any upper triangular matrices are located on the diagonal. 
To illustrate this, for any upper triangular matrix $\bR \in \real^{n\times n}$, where the diagonal values are $r_{ii}$ for all $i\in \{1,2,\ldots,n\}$, we have
$$
\bR \be_i = r_{ii}\be_i,
$$
where $\be_i$ is the $i$-th unit basis vector in $\real^n$ (i.e., $\be_i$ corresponds to the $i$-th column of the $n\times n$ identity matrix $\bI_n$).
Thus, we can decompose $\bU$ into the sum of $\bLambda$ and $\bT$.

\begin{remark}[$m$-th Power]\label{remark:mth_schur}
The above observation shows that the eigenvalues of the $m$-th power $\bA^m$ are the  $m$-th powers of the eigenvalues of $\bA$.
\end{remark}

A final observation regarding the second form of the Schur decomposition is presented as follows. 
Given $\bA\bQ = \bQ(\bLambda +\bT)$, it follows that 
$$
\bA \bq_k = \lambda_k\bq_k + \sum_{i=1}^{k-1}t_{ik}\bq_i,
$$
where $t_{ik}$ is the ($i,k$)-th entry of $\bT$. 
The form is quite similar to the eigenvalue decomposition. However, instead of being eigenvectors, the columns of $\bQ$ are orthonormal bases that are correlated with each other.

\section{Complex and Quasi-Triangular Schur Decomposition*}
In the primary result of Theorem~\ref{theorem:schur-decomposition}, we focus exclusively on real matrices with real eigenvalues. However, this restriction does not always apply in practical situations. A more general scenario is addressed in the subsequent theorem.
This decomposition is attributed to \textit{Issai Schur} (1875{\textendash}1941), a Russian mathematician, who studied and worked in Germany for most of his life.
\index{Decomposition: Complex Shur}
\begin{theoremHigh}[Complex Schur Decomposition]\label{theorem:schur-decomposition_complex}
Let  $\bA\in \complex^{n\times n}$ be any \textbf{complex} square matrix. Then, it  can be factored as 
$$
\bA = \bU\bT\bU^\ast,
$$
where $\bU\in\complex^{n\times n}$ is a unitary matrix, and $\bT\in\complex^{n\times n}$ is an upper triangular matrix (not necessarily real). 
\end{theoremHigh}
\begin{proof}[of Theorem~\ref{theorem:schur-decomposition_complex}]
We establish this through recursive computation and the principle of similarity transformations.
The proof draws from \citet{bernstein2009matrix, horn2012matrix}, with minor adaptations.
Suppose $\bA$ has a complex eigenpair $(\lambda_1\in\complex, \bx_1\in\complex^n)$ with $\bx_1^*\bx_1=1$ (of unit length).
Extend $\bx_1$ to form a unitary matrix $\bU_1=[\bx_1, \widehatbU_1]\in\complex^{n\times n}$, where $\widehatbU_1^*\widehatbU_1=\bI_{n-1}$ and $\bx_1^*\widehatbU_1=\bzero$.
Therefore, the first column of $\bU_1^{-1}\bA\bU_1$ is 
$$
(\bU_1^{-1}\bA\bU_1)[:,1] = \bU_1^{-1}\bA\bx_1 = \lambda_1\bU_1^{-1}\bx_1 = \lambda_1\be_1,
\qquad  (\text{since $\bU_1\be_1=\bx_1$}).
$$
This indicates $\bA$ admits the following decomposition:
$$
\bA =
\bU_1
\begin{bmatrix}
\lambda_1 & \bc_1^\top \\
\bzero & \bA_1 
\end{bmatrix}
\bU_1^{-1},
\gap 
\text{where $\bc_1\in\complex^{n-1}$, $\bA_1\in\complex^{(n-1)\times (n-1)}$}.
$$
Since $\bA$ and $\footnotesize\begin{bmatrix}
\lambda_1 & \bc_1^\top \\
\bzero & \bA_1 
\end{bmatrix}$ are unitarily  similar, they have the same eigenvalues (see Proposition~\ref{proposition:eigenvalue-similar-matrices} for a hint on the proof).
If $\bA$ has eigenvalues $\{\lambda_1,\lambda_2\ldots,\lambda_n\}$, then $\bA_1$ has eigenvalues $\{\lambda_2, \lambda_3,\ldots,\lambda_n\}$.
This, in return, implies $\bA_1$ admits the following decomposition:
$$
\bA_1 =
\widehatbU_2
\begin{bmatrix}
\lambda_2 & \bc_2^\top \\
\bzero & \bA_2 
\end{bmatrix}
\widehatbU_2^{-1},
\gap 
\text{where $\bc_2\in\complex^{n-2}$, $\bA_2\in\complex^{(n-2)\times (n-2)}$}.
$$
Therefore, $\bA$ admits 
$$
\bA=
\bU_1\bU_2
\begin{bmatrix}
\lambda_1 & \widehat{c}_{11} & \widehat{\bc}_{12}^\top \\
0 & \lambda_2 & \bc_2^\top \\
0 & 0 & \bA_2
\end{bmatrix}
\bU_2^{-1}\bU_1^{-1},
$$
where $\widehat{\bc}_1=[\widehat{c}_{11};\widehat{\bc}_{12}]=\widehatbU_2^{-1}\bc_1$, and $\footnotesize\bU_2=\begin{bmatrix}
1 & \bzero \\
\bzero & \widehatbU_2
\end{bmatrix}$
is unitary.
By continuing this iterative process, we can construct the unitary matrix $\bU=\bU_1\bU_2\ldots\bU_n$ to achieve the Schur decomposition.
We may also note that the first column of $\bU$ is $\bx_1$, an eigenvector of $\bA$ associated to the eigenvalue $\lambda_1$.
\end{proof}

\begin{remark}[Conjugate Complex Eigenvalues]\label{remark:conjug_pari}
Theorem~\ref{theorem:schur-decomposition_complex} deals with complex matrices, while Theorem~\ref{theorem:schur-decomposition} focuses on real matrices with real eigenvalues. 
However, when $\bA\in\real^{n\times n}$ is real with complex eigenvalues, we may obtain the following results.
Suppose $\bA$ admits a complex Schur decomposition $\bT=\bU^*\bA\bU$. Then, $\overline{\bT}=\bU^\top\bA\overline{\bU}$ is the complex conjugate of $\bT$ and is unitarily similar to $\bA$ and $\bT$.
Therefore, $\bT$, $\overline{\bT}$, and $\bA$ have the same eigenvalues. This implies if $a+ib$ is an eigenvalue of $\bT$, then $a-ib$ is also an eigenvalue of $\bT$. That is, non-real eigenvalues of $\bA$ appear in complex conjugate pairs.
\end{remark}

\paragraph{Quasi-triangular.} 
In the context of the above theorem, if we restrict $\bU$ and $\bT$ to be \textbf{real matrices}, and $\bA$ may contain complex eigenvalues, then $\bT$ takes a \textit{quasi-triangular form}, i.e., it has diagonal blocks of order 1 and 2.
To see this, we need the following lemmas.

For a complex eigenvector of a general matrix, we can find the eigenvalue of its conjugate easily.
\begin{lemma}[Conjugate Eigenpair]\label{lemma:conjugate_eig_pair}
Let $\bw=\bx+i\by$ be an eigenvector of $\bA$ corresponding to the eigenvalue $z=\lambda+i\mu$ ($\lambda, \mu, \bx, \by$ are real scalars or real vectors). 
Then, $\overline{\bw}=\bx-i\by$ is an eigenvector of $\bA$ corresponding to the eigenvalue $\overline{z}=\lambda-i\mu$.
\end{lemma}

For any real matrix with complex eigenvalues, we can find the following block Schur decomposition.
\begin{lemma}[Real Block Schur Decomposition]\label{lemma:real_block_schur}
Let $\bA\in \real^{n\times n}$ be any \textbf{real} square matrix with \textbf{complex} eigenvalues. 
Then, it can be factored as
$$
\bA = \bS\bU\bS^\top,
$$
where $\bS$ is nonsingular, and $\bU$ is a  (real) block upper triangular matrix, where the upper-left block is a 2-by2 block. 
\end{lemma}
\begin{proof}[of Lemma~\ref{lemma:real_block_schur}]
Suppose $\bw=\bx+i\by$ is an eigenvector of $\bA$ corresponding to the eigenvalue $z=\lambda+i\mu$ ($\lambda, \mu, \bx, \by$ are real scalars or real vectors). 
Then, $\overline{\bw}=\bx-i\by$ is an eigenvector of $\bA$ corresponding to the eigenvalue $\overline{z}=\lambda-i\mu$.
It then can be shown that $\bw$ and $\overline{\bw}$ are linearly independent, and $\bx$ and $\by$ are also linearly independent; otherwise, the eigenvalue and eigenvector are real.
Let 
$
\bB\triangleq
\footnotesize
\begin{bmatrix}
\lambda & \mu\\
-\mu & \lambda
\end{bmatrix}
$
and construct the  nonsingular matrix $\bS$:
$
\bS=
\begin{bmatrix}
\bx & \by & \bS_1 
\end{bmatrix}\in\real^{n\times n}.
$
Since $\bS^{-1}[\bx,\by] = \footnotesize\begin{bmatrix}
\bI_2\\
\bzero
\end{bmatrix}$, it follows that 
$$
\bS^{-1}\bA\bS = 
\bS^{-1} 
\begin{bmatrix}
\bA[\bx, \by] & \bA\bS_1 
\end{bmatrix}= 
\bS^{-1} 
\begin{bmatrix}
[\bx, \by]\bB & \bA\bS_1 
\end{bmatrix}
=
\begin{bmatrix}
\bB & \bC \\
\bzero & \bA_1 
\end{bmatrix},
$$
where $\bC\in\real^{2\times (n-2)}$ and $\bA_1\in\real^{(n-1)\times (n-2)}$.
This completes the proof.
\end{proof}

\index{Quasi-triangular}
\begin{theoremHigh}[Quasi-Triangular Schur Decomposition]\label{theorem:quasi_schur}
Let $\bA\in \real^{n\times n}$ be any \textbf{real} square matrix with possible \textbf{complex} eigenvalues. Then, it can be factored as
\begin{equation}\label{equation:quasi_schur_form1}
\bA = \bS\bD\bS^{-1},
\end{equation}
where $\bS\in\real^{n\times n}$ is an $n\times n$ nonsingular matrix, and $\bD\in\real^{n\times n}$ is a real upper quasi-triangular matrix with the following form
$$
\bD = 
\begin{bmatrix}
\bA_1 &\boxtimes&\boxtimes&\boxtimes \\
& \bA_2 &\boxtimes&\boxtimes \\
&&\ddots &\boxtimes \\
&&&\bA_p
\end{bmatrix},
$$
where $\boxtimes$ denotes entries that are not necessarily zero, and each $\bA_i$ is either a real scalar (corresponding to a real eigenvalue of $\bA$) or a $2\times 2$ real matrix (corresponding a complex eigenvalue of $\bA$) of the following form 
\begin{equation}\label{equation:quasi_schur}
\bA_i=
\begin{bmatrix}
\lambda_i & \mu_i \\
-\mu_i & \lambda_i
\end{bmatrix}, 
\gap 
\lambda_i,\mu_i\in\real, \mu_i>0, \text{ and $\lambda_i\pm \mu_i$ are eigenvalues of $\bA$}.
\end{equation}
In the meantime, $\bA$ also can be factored as 
\begin{equation}\label{equation:quasi_schur_form2}
\bA = \bU\widetildebD\bU^\top,
\end{equation}
where $\bU\in\real^{n\times n}$ is real orthogonal, and $\widetildebD\in\real^{n\times n}$ is also real upper quasi-triangular.
\end{theoremHigh}
\begin{proof}[of Theorem~\ref{theorem:quasi_schur}]
The first part can be proved recursively using Lemma~\ref{lemma:real_block_schur}.
For the second part, suppose $\bS$ in \eqref{equation:quasi_schur_form1} admits QR decomposition $\bS=\bU\bR$. Then, 
$$
\bU^\top\bA\bU = 
\bR\bD\bR^{-1}
=
\begin{bmatrix}
\bR_{11}\bA_1\bR_{11}^{-1} &\boxtimes&\boxtimes&\boxtimes \\
& \bR_{22}\bA_2\bR_{22}^{-1} &\boxtimes&\boxtimes \\
&&\ddots &\boxtimes \\
&&&\bR_{pp}\bA_p\bR_{pp}^{-1}
\end{bmatrix}
$$
is upper quasi-triangular.
For $i\in\{1,2,\ldots,p\}$, if $\bA_i$ is real scalar, then it is identical to that in \eqref{equation:quasi_schur_form1}; if $\bA_i$ is a $2\times 2$ real matrix, it is similar to the corresponding block of \eqref{equation:quasi_schur_form1}.
\end{proof}

\section{Applications and Properties}
A few results can be easily proved using Schur decomposition. For example, to prove the trace of a matrix is equal to the sum of eigenvalues (Theorem~\ref{theorem:eigen_trace2}), to prove the existence of the block-diagonalization (Theorem~\ref{theorem:block_real_eigens}), to prove the Schur inequality (Theorem~\ref{theorem:schur_inequality}), and to prove the existence of the spectral decomposition (Theorem~\ref{theorem:spectral_theorem}).
We provide more results derived from the Schur decomposition in this section.


\index{Rank-one perturbation}
\paragraph{Rank-one perturbation.}
The Schur decomposition of the rank-one perturbed matrix can be obtained efficiently.
\begin{theorem}[Rank-One Perturbed Schur]\label{theorem:rno_pertsch}
Suppose $\bA\in\complex^{n\times n}$ has eigenvalues $\lambda_1, \lambda_2,\ldots,\lambda_n\in\complex$, where $\bA\bx=\lambda_1\bx$. 
Then, for any vector $\bv\in\complex^n$, the eigenvalues of $\bA+\bx\bv^*$ are $\lambda_1+\bv^*\bx, \lambda_2, \lambda_3, \ldots,\lambda_n$. The Schur decomposition of $\bA+\bx\bv^*$ can be obtained efficiently if the Schur decomposition of $\bA$ is known.
\end{theorem}
\begin{proof}[of Theorem~\ref{theorem:rno_pertsch}]
Following the proof of Theorem~\ref{theorem:schur-decomposition_complex},  we have 
$$
\bU_1^*\bA\bU_1 =
\begin{bmatrix}
\lambda_1 & \bc_1^\top \\
\bzero & \bA_1 
\end{bmatrix}
,
$$
where $\bc_1\in\complex^{n-1}$, $\bA_1\in\complex^{(n-1)\times (n-1)}$,  the first column of $\bU_1$ is $\bu_{11}\triangleq\bx_1=\frac{\bx}{\normtwo{\bx}}$, and $\bA_1\in\complex^{(n-1)\times(n-1)}$ has eigenvalues $\lambda_2, \lambda_3, \ldots,\lambda_n$.
Therefore, let the $i$-th column of $\bU_1$ be $\bu_{1i}$, we have 
$$
\bU_1^*\bx\bv^*\bU_1 
=
\begin{bmatrix}
\bu_{11}^*\bx \\
\bu_{12}^*\bx\\
\ldots\\
\bu_{1n}^*\bx
\end{bmatrix}
\begin{bmatrix}
\bv^*\bu_{11} & \bv^*\bu_{12} & \ldots  &  \bv^*\bu_{1n} 
\end{bmatrix}
=
\begin{bmatrix}
\bu_{11}^*\bx \\
0\\
\ldots\\
0
\end{bmatrix}
\begin{bmatrix}
\bv^*\bu_{11} & \bv^*\bu_{12} & \ldots  &  \bv^*\bu_{1n} 
\end{bmatrix}.
$$
Only the first row is nonzero. 
Therefore, 
$$
\bU_1^*(\bA+ \bx\bv^*)\bU_1 
=
\begin{bmatrix}
\lambda_1+\bv^*\bx & \boxtimes \\
\bzero & \bA_1 
\end{bmatrix}, 
$$
where $\boxtimes$ denotes a vector that is not necessarily zero.
This implies the eigenvalues of $\bA+\bx\bv^*$ are $\lambda_1+\bv^*\bx, \lambda_2, \lambda_3, \ldots,\lambda_n$.
\end{proof}

\index{Cayley-Hamilton theorem}
\paragraph{Cayley-Hamilton theorem.}
We use the Cayley-Hamilton theorem to demonstrate the application of eigenvalue decomposition in Section~\ref{section:fibonacci}. Here, we provide a rigorous proof of the theorem.
\begin{theorem}[Cayley-Hamilton Theorem]\label{theorem:cayley_hami}
A matrix satisfies its own characteristic equation (Definition~\ref{definition:characteristic_polynomial}). That is, given a matrix $\bA\in\real^{n\times n}$, $p_{\bA}(\lambda)=\det(\lambda\bI-\bA)=\prod_{i=1}^{n}(\lambda-\lambda_i)$ and $p_{\bA}(\lambda)=0$ if $\lambda$ is an eigenvalue of $\bA$.
Then, $\bA$ also satisfies this characteristic equation: $p_{\bA}(\bA)=\prod_{i=1}^{n}(\bA-\lambda_i\bI)=\bzero$.
\end{theorem}
\begin{proof}[of Theorem~\ref{theorem:cayley_hami}]
Suppose $\bA$ admits the Schur decomposition $\bA=\bQ\bU\bQ^\top$. Then, 
$$
p_{\bA}(\bA)=\prod_{i=1}^{n}(\bQ\bU\bQ^\top-\lambda_i\bI) =\bQ \cdot p_{\bA}(\bU) \cdot\bQ^\top.
$$
Therefore, it suffices to show $p_{\bA}(\bU)=\prod_{i=1}^{n}(\bU-\lambda_i\bI)=\bzero$.
We observe that the upper left 2-by-2 block of $(\bU-\lambda_1\bI)(\bU-\lambda_2\bI)$ is zero. This again invokes the upper left 3-by-3 block $(\bU-\lambda_1\bI)(\bU-\lambda_2\bI)(\bU-\lambda_3\bI)$ to be zero. Continuing this process, the result follows.
\end{proof}

We have shown in Remark~\ref{remark:power-eigenvalue-decom} that  eigenvalue decomposition can help identify the $m$-th power of a matrix easily. 
The Cayley-Hamilton theorem can be used to express the $m$-th power of a square matrix $\bA\in\real^{n\times n}$ as a linear combination of $\bI, \bA, \bA^2, \ldots, \bA^{m-1}$, i.e., the linear combination of lower power values.
Suppose  $p_{\bA}(\lambda)=\det(\lambda\bI-\bA ) =\lambda^n + \gamma_{n-1} \lambda^{n-1} + \ldots + \gamma_1 \lambda  + \gamma_0$. Then, 
\begin{equation}\label{equation:calhm_eq1}
\bA^n =- \gamma_{n-1} \bA^{n-1} - \ldots - \gamma_1 \bA  - \gamma_0\bI.
\end{equation}
This also implies 
$$
\bI = -\frac{1}{\gamma_0} (\bA^{n-1}+\gamma_{n-1}\bA^{n-2}+\ldots+\gamma_1)\bA.
$$

\paragraph{Computation of inverses.}
If $\bA$ is nonsingular, then multiplying \eqref{equation:calhm_eq1} by $\bA^{-1}$ yields
$$
\bA^{-1} = -\left(\frac{1}{\gamma_0}\bA^{n-1} +\frac{\gamma_{n-1}}{\gamma_0} \bA^{n-2} + \ldots + \frac{\gamma_1}{\gamma_0} \bI\right).
$$

\index{Sylvester's theorem}
\paragraph{Sylvester's theorem.}
Given $\bA,\bX\in\real^{n\times n}$, $\bA$ and $\bX$ are said to \textit{commute} if $\bA\bX=\bX\bA$. 
More generally, consider the  equation $\bA\bX=\bX\bB$, where $\bA\in\real^{n\times n}$, $\bB\in\real^{m\times m}$, and $\bX\in\real^{n\times m}$. 
The Cayley-Hamilton theorem indicates (see Problem~\ref{problem:poly_cay}):
\begin{equation}\label{equation:poly_cay_gen}
p(\bA)\bX=\bX p(\bB), 
\gap 
\text{for any polynomial $p(\lambda)$.}
\end{equation}
This reveals  Sylvester's theorem.

\begin{theorem}[Sylvester's Theorem]\label{theorem:sylvesters_theorem}
Let $\bA\in\real^{n\times n}$ and $\bB\in\real^{m\times m}$. 
\begin{itemize}
\item If $\Lambda(\bA)\cup \Lambda(\bB)=\varnothing$ (i.e., the intersection of the spectrum sets is empty),  the equation $\bA\bX-\bX\bB=\bzero$ is satisfied only when $\bX=\bzero\in\real^{n\times m}$.
\item More generally,  \textit{Sylvester's equation} $\bA\bX-\bX\bB=\bC$ has a unique solution $\bX\in\real^{n\times m}$ for each $\bC\in\real^{n\times m}$ if and only if $\Lambda(\bA)\cup \Lambda(\bB)=\varnothing$.~\footnote{If $\bA$ and $\bB$ are complex, then there is a unique complex solution $\bX$ for each $\bC\in\complex^{n\times m}$.}
\end{itemize}
\end{theorem}
\begin{proof}[of Theorem~\ref{theorem:sylvesters_theorem}]
The second part is a direct result of the first part; so we only prove the first part. 
For the first part, it suffices to show that $p_{\bB}(\bA)\bX=\bX p_{\bB}(\bB)=\bzero$ due to \eqref{equation:poly_cay_gen}.
Suppose $\bB$ has eigenvalues $\lambda_1, \lambda_2, \ldots,\lambda_n$ and admits the characteristic polynomial $p_{\bB}(\lambda)=\prod_{i=1}^{n}(\lambda-\lambda_i)$ and $p_{\bB}(\bA)=\prod_{i=1}^{n}(\bA-\lambda_i\bI)$.
If  $\Lambda(\bA)\cup \Lambda(\bB)=\varnothing$, then each component $(\bA-\lambda_i\bI)$ is nonsingular, and $p_{\bB}(\bA)$ is nonsingular.
Therefore, $p_{\bB}(\bA)\bX=\bzero$ if and only if $\bX=\bzero$.
Conversely, if $p_{\bB}(\bA)\bX=\bzero$ has nontrivial solution, then there must be at least one component $(\bA-\lambda_i\bI)$ is singular. Thus, $\Lambda(\bA)\cup \Lambda(\bB)\neq \varnothing$.
\end{proof}

The existence of the Schur decomposition reveals the eigenvalues of $\bB^{-1}\bA$ (if $\bB$ is nonsingular) from the upper triangular matrices.
\begin{corollary}[Eigenvalues from Schur]\label{corollary:eig_utv}
Suppose $\bA,\bB\in\real^{n\times n}$ admit  decompositions $\bA=\bQ\bT_A\bV^\top$ and $\bB=\bQ\bT_B\bV^\top$, respectively, where $\bQ, \bV$ are orthogonal and $\bT_A, \bT_B$ are upper triangular.
Then, the diagonal elements of $\bT_B^{-1}\bT_A$ are the eigenvalues of $\bB^{-1}\bA$ (we assumes all the eigenvalues are real).
\end{corollary}
\begin{proof}[of Corollary~\ref{corollary:eig_utv}]
The proof relies on the Schur decomposition (Theorem~\ref{theorem:schur-decomposition}) and we assume all the eigenvalues discussed are real for simplicity.
Suppose $\bB^{-1}\bA$ admits  a Schur decomposition $\bB^{-1}\bA=\bV\bU\bV^\top$ ($\bV$ is orthogonal, $\bU$ is upper triangular)~\footnote{In the corollary, if we don't assume real eigenvalues, then $\bU$ can be upper quasi-triangular. And $\bT_A$ shown below is also upper quasi-triangular.}, and $\bB\bV$ admits a QR decomposition $\bB\bV=\bQ\bT_B\implies \bB=\bQ\bT_B\bV^\top$ ($\bQ$ is orthogonal, $\bT_B$ is upper triangular).
Then, $\bA=\bB\bV\bU\bV^\top=\bQ\underbrace{(\bT_B\bU)}_{\triangleq\bT_A}\bV^\top$, where $\bT_A\triangleq \bT_B\bU$ is upper triangular.
This completes the proof.
\end{proof}

\begin{exercise}
Discuss the connection between the decompositions in Corollary~\ref{corollary:eig_utv} and the UTV decomposition (Section~\ref{section:ulv-urv-decomposition}).
\end{exercise}

\begin{problemset}


\item Given $\bA\in\real^{n\times n}$, for $n=2$, show that 
\begin{itemize}
\item $\det(\bI+\bA) = 1+\det(\bA) + \tr(\bA)$.
\end{itemize}
For $n=3$, show that
\begin{itemize}
\item $\det(\bI+\bA) = 1+ \det(\bA) +\tr(\bA) + \frac{1}{2} \tr(\bA)^2 - \frac{1}{2}\tr(\bA^2)$.
\end{itemize}

\item Given $\bA\in\real^{n\times n}$, and let $\bB$ be the matrix obtained by multiplying a row of $\bA$ by nonnegative scalar $\gamma$. Prove that $\det(\bB) = \gamma\det(\bA)$.  

\item \label{problem:poly_cay} Given any polynomial $p(\lambda)$, show that $p(\bA)\bX=\bX p(\bB)$ if $\bA\bX=\bX\bB$.

\item Given any polynomial $p(\lambda)$, show that $\bA\bB p(\bA\bB)=\bA p(\bB\bA)\bB$ if $\bA\in\real^{m\times n}$ and $\bB\in\real^{n\times m}$.

\item \label{prob:diag_uppt} \textbf{Diagonalization of upper triangular matrices.} Let  $\bU\in\real^{n\times n}$ be an upper triangular matrix whose $(i,j)$-th entry is denoted by $u_{ij}$, and let $\bD_t=\diag(t, t^2, \ldots,t^n)$ be a diagonal matrix. Show that the similarity transformation on $\bU$ takes the following form:
$$
\bD_t\bU\bD_t^{-1}
=
\begin{bmatrix}
	u_{11} & t^{-1} u_{12} & t^{-2} u_{13} & \ldots & t^{-n+1}u_{1n}\\
	0   &  u_{22} & t^{-1} u_{23} & \ldots & t^{-n+2}u_{2n}\\
	0   &  0  &  u_{33} & \ldots & t^{-n+3}u_{3n}\\
	0   &  \vdots & \vdots & \ddots & \vdots\\
	0   &  0 & 0 & \ldots & u_{nn}\\
\end{bmatrix}.
$$
Thus, when $t$ is sufficiently large, the off-diagonal values can be made arbitrarily  small.

\item \label{prob:bounded_spec} \textbf{Bounded norm using spectral radius.}
Let $\bA\in\real^{n\times n}$ and $\epsilon>0$. Show that there exists a matrix norm $\norm{\cdot}$ on $\real^{n\times n}$ such that $\rho(\bA)\leq \rho(\bA)+\epsilon$, where $\rho(\bA)$ is the spectral radius of $\bA$. \textit{Hint: Use Problem~\ref{prob:diag_uppt}, Schur decomposition $\bA=\bU\bT\bU^\top$, and consider the matrix-1-norm (Equation~\eqref{equation:mat_one_norm}) by $\norm{\bA}=\norm{\bD_t\bU^\top\bA\bU\bD_t^{-1}}_1$.}

\item \label{prob:schur_range_hermi} \textbf{Complex Schur for range-Hermitian \citep{bernstein2009matrix}.} Let $\bA\in\complex^{n\times n}$ with $\rank(\bA)=r$. Show that there exist a unitary $\bU\in\complex^{n\times n}$ and a nonsingular $\bB\in\complex^{r\times r}$ such that 
$$
\bA = 
\bU
\begin{bmatrix}
\bB & \bzero \\
\bzero & \bzero 
\end{bmatrix}
\bU^*.
$$
In addition, if $\bA$ is normal, then $\bB$ can be chosen to be diagonal (Theorem~\ref{theorem:normal_Complex_spectral_theorem}).

\item Use the Cayley-Hamilton theorem to compute the inverse of
$
\footnotesize
\begin{bmatrix}
3 & 0 & 1\\
0 & 2 & 0\\
0 & 0 & 3
\end{bmatrix}$.

\item Show that $\bA\in\real^{n\times n}$ is nilpotent if and only if $\trace(\bA^k)=0$ for all $k\in\{1,2,\ldots,n\}$. \textit{Hint: Use the $m$-th power eigenvalues, Remark~\ref{remark:mth_schur}.}
\end{problemset}

%% file: chapter-spectral.tex
\newpage
\chapter{Spectral Decomposition and Diagonalization}\label{chapter:spectral-decomposition}
\begingroup
\hypersetup{
	linkcolor=structurecolor,
	linktoc=page,  
}
\minitoc \newpage
\endgroup
\section{Spectral Decomposition}
\lettrine{\color{caligraphcolor}T}
The spectral theorem, also referred to as  the spectral decomposition for symmetric matrices, states that symmetric matrices possess  real eigenvalues and  can be diagonalized using an (real) orthonormal basis \footnote{Note that the spectral decomposition for \textit{Hermitian matrices} states that Hermitian matrices also have real eigenvalues and  can be diagonalized using a complex orthonormal basis.}. 
In the subsequent sections, we will outline the primary outcome and postpone more detailed explorations.

\index{Decomposition: Spectral}
\begin{theoremHigh}[Spectral Decomposition]\label{theorem:spectral_theorem}
A real matrix $\bA \in \real^{n\times n}$ is symmetric if and only if there exist an orthogonal matrix $\bQ$ and a diagonal matrix $\bLambda$ such that~\footnote{The spectral decomposition for real normal matrices is discussed in Theorem~\ref{theorem:normal_real_spectral_theorem} and for general complex normal matrices is discussed in Theorem~\ref{theorem:normal_Complex_spectral_theorem}. We will primarily consider the spectral theorem for real symmetric matrices; while some properties of complex Hermitian matrices are discussed in Problems~\ref{problem:real_herm}$\sim$\ref{problem:real_skewhe_final}.}
\begin{equation*}
	\bA = \bQ \bLambda \bQ^\top,
\end{equation*}
where the columns of $\bQ = [\bq_1, \bq_2, \ldots, \bq_n]$ are eigenvectors of $\bA$ and are mutually orthonormal, and the diagonal entries of $\bLambda=\diag(\lambda_1, \lambda_2, \ldots, \lambda_n)$ are the corresponding eigenvalues of $\bA$, which are real. 
The matrix $\bLambda$ is \textbf{uniquely} determined  by $\bA$ up to permutation of eigenvalues.
 ~This is known as the \textit{spectral decomposition} or \textit{spectral theorem} for a real symmetric matrix $\bA$. Specifically, the following properties hold:
\begin{enumerate}
\item A symmetric matrix exclusively has \textit{real eigenvalues}.

\item Its eigenvectors can be normalized to form an orthonormal set.

\item The rank of $\bA$ matches the count of its nonzero eigenvalues.

\item When the eigenvalues are distinct, the eigenvectors are automatically linearly independent (Theorem~\ref{theorem:independent-eigenvector-theorem}).
\end{enumerate}

\end{theoremHigh}

In contrast to the eigenvalue decomposition  (Theorem~\ref{theorem:eigenvalue-decomposition}), which necessitates a square matrix $\bA$ with linearly independent eigenvectors, the spectral theorem applies to any symmetric matrix, ensuring it can be diagonalized with orthonormal eigenvectors.

According to Proposition~\ref{proposition:eigenvalue-similar-matrices}, similar matrices have identical eigenvalues.
From the spectral decomposition, it follows  that $\bA$ and $\bLambda$ are similar matrices, thus sharing the same eigenvalues.  
For any diagonal matrix, the eigenvalues are simply the diagonal entries. \footnote{Indeed, as demonstrated in the previous section, the diagonal elements of a triangular matrix represent its eigenvalues.} This relationship is evident from
$
\bLambda \be_i = \lambda_i \be_i,
$
where $\be_i$ denotes  the $i$-th standard  basis vector. Hence, the matrix $\bLambda$ encapsulates the eigenvalues of $\bA$.

The spectral decomposition has a wide range of applications across various fields. In machine learning, spectral decomposition is primarily used in principal component analysis (PCA), a statistical procedure that uses an orthogonal transformation to convert a set of observations of possibly correlated variables into a set of values of linearly uncorrelated variables called principal components (see Section~\ref{section:pca_isvd}). 
In quantum mechanics, observables are represented by Hermitian operators, which generalize symmetric matrices. The spectral theorem ensures that these operators have real eigenvalues, corresponding to the possible outcomes of measurements. The eigenvectors represent the states in which the observable has a definite value.
In signal processing, the spectral decomposition is used to analyze and process signals. The Fourier transform, which is a special case of spectral decomposition, decomposes a signal into its frequency components. This is essential for filtering, compression, and noise reduction.

\section{Existence of  Spectral Decomposition}\label{section:existence-of-spectral}

We will prove the theorem in several steps. 

\begin{tcolorbox}[title={Symmetric Matrix Property 1 of 4},colback=\mdframecolorTheorem]
\begin{proposition}[Real Eigenvalues]\label{proposition:real-eigenvalues-spectral}
The eigenvalues of any real symmetric matrix are all real. \footnote{In fact, it can be shown that any complex Hermitian matrix also has only real eigenvalues; see Problem~\ref{problem:real_herm}.}
\end{proposition}
\end{tcolorbox}
\begin{proof}[of Proposition~\ref{proposition:real-eigenvalues-spectral}]
Consider an eigenvalue $\lambda$ of the real symmetric matrix $\bA$, represented as a complex number $\lambda=a+ib$, where $a$ and $b$ are real. Its complex conjugate is $\bar{\lambda}=a-ib$. 
Similarly, for the corresponding complex eigenvector $\bx = \bc+i\bd$, its complex conjugate is $\bar{\bx}=\bc-i\bd$, where $\bc$ and $\bd$ are real vectors. We then have the following property:
$$
\bA \bx = \lambda \bx\qquad   \underrightarrow{\text{ leads to }}\qquad  \bA \bar{\bx} = \bar{\lambda} \bar{\bx}\qquad   \underrightarrow{\text{ transpose to }}\qquad  \bar{\bx}^\top \bA =\bar{\lambda} \bar{\bx}^\top.
$$
Taking the dot product of the first equation with $\bar{\bx}$ and the last equation with $\bx$:
$$
\bar{\bx}^\top \bA \bx = \lambda \bar{\bx}^\top \bx \qquad \text{and } \qquad \bar{\bx}^\top \bA \bx = \bar{\lambda}\bar{\bx}^\top \bx.
$$
From these equations, we arrive at the equality $\lambda\bar{\bx}^\top \bx = \bar{\lambda} \bar{\bx}^\top\bx$, where $\bar{\bx}^\top\bx = (\bc-i\bd)^\top(\bc+i\bd) = \bc^\top\bc+\bd^\top\bd$ is a real number. Therefore, the imaginary part of $\lambda$ must be zero, implying that $\lambda$ is real.
\end{proof}

\begin{exercise}[Normal and Symmetric]
Show that a real square matrix $\bA\in\real^{n\times n}$ is normal and all its eigenvalues are real if and only if $\bA$ is symmetric.
\end{exercise}

\begin{exercise}[Conjugate]
Given a complex square matrix $\bA\in\complex^{n\times n}$ and its eigenpair $(\lambda,\bu)$, show that $(\overline{\lambda}, \overline{\bu})$ is an eigenpair of its conjugate $\overline{\bA}$.
This indicates that when $\bA$ is real, its eigenvalues come in conjugate pairs (Remark~\ref{remark:conjug_pari}).
\end{exercise}

\begin{tcolorbox}[title={Symmetric Matrix Property 2 of 4},colback=\mdframecolorTheorem]
\begin{proposition}[Orthogonal Eigenvectors]\label{proposition:orthogonal-eigenvectors}
The eigenvectors  corresponding to distinct eigenvalues of any real symmetric matrix are orthogonal.
Therefore, we can normalize these  eigenvectors to make them orthonormal since $\bA\bx = \lambda \bx \implies \bA\frac{\bx}{\normtwo{\bx}} = \lambda \frac{\bx}{\normtwo{\bx}}$, which corresponds to the same eigenvalue.
\footnote{In fact, the eigenvectors of a complex Hermitian matrix corresponding to distinct eigenvalues are also orthogonal to each other.}
\end{proposition}
\end{tcolorbox}
\begin{proof}[of Proposition~\ref{proposition:orthogonal-eigenvectors}]
Suppose  $\lambda_1$ and $\lambda_2$ are distinct eigenvalues of the real symmetric matrix $\bA$, with corresponding  eigenvectors $\bx_1$ and $\bx_2$, respectively, satisfying $\bA\bx_1=\lambda \bx_1$ and $\bA\bx_2 = \lambda_2\bx_2$. We have the following equality:
$$
\bA\bx_1=\lambda_1 \bx_1 \quad\implies\quad \bx_1^\top \bA =\lambda_1 \bx_1^\top \quad\implies\quad \bx_1^\top \bA \bx_2 =\lambda_1 \bx_1^\top\bx_2,
$$
and 
$$
\bA\bx_2 = \lambda_2\bx_2 \quad\implies\quad  \bx_1^\top\bA\bx_2 = \lambda_2\bx_1^\top\bx_2.
$$
Thus, we have  $\lambda_1 \bx_1^\top\bx_2=\lambda_2\bx_1^\top\bx_2$. Since  $\lambda_1\neq \lambda_2$, the eigenvectors are orthogonal.
\end{proof}

In the above Proposition~\ref{proposition:orthogonal-eigenvectors}, we demonstrate that the eigenvectors corresponding to distinct eigenvalues  of symmetric matrices are orthogonal. More generally, we prove the crucial theorem that eigenvectors corresponding to distinct eigenvalues of any matrix are linearly independent.
\begin{theorem}[Independent Eigenvector Theorem]\label{theorem:independent-eigenvector-theorem}
Let $\bA\in \real^{n\times n}$ be any matrix with $k$ distinct eigenvalues, then any set of $k$ corresponding (nonzero) eigenvectors are linearly independent.
\end{theorem}
\begin{proof}[of Theorem~\ref{theorem:independent-eigenvector-theorem}]
We will prove this by induction. First, we will prove that any two eigenvectors corresponding to distinct eigenvalues are linearly independent. Suppose  $\bv_1$ and $\bv_2$ are eigenvectors corresponding to distinct eigenvalues $\lambda_1$ and $\lambda_2$, respectively. Assume there exists a nonzero vector $\bx=[x_1,x_2] \neq \bzero $ such that  
\begin{equation}\label{equation:independent-eigenvector-eq1}
x_1\bv_1+x_2\bv_2=\bzero.
\end{equation}
That is, we assume $\bv_1$ and $\bv_2$ are linearly dependent.
Premultiply Equation~\eqref{equation:independent-eigenvector-eq1}  by $\bA$, we obtain
\begin{equation}\label{equation:independent-eigenvector-eq2}
x_1 \lambda_1\bv_1 + x_2\lambda_2\bv_2 = \bzero.
\end{equation}
Premultiply Equation~\eqref{equation:independent-eigenvector-eq1}  by $\lambda_2$, we get 
\begin{equation}\label{equation:independent-eigenvector-eq3}
x_1\lambda_2\bv_1 + x_2\lambda_2\bv_2 = \bzero.
\end{equation}
Subtract Equation~\eqref{equation:independent-eigenvector-eq2} from Equation~\eqref{equation:independent-eigenvector-eq3}, we find
$$
x_1(\lambda_2-\lambda_1)\bv_1 = \bzero.
$$
Since $\lambda_2\neq \lambda_1$ and $\bv_1\neq \bzero$, we must have $x_1=0$. From Equation~\eqref{equation:independent-eigenvector-eq1} and $\bv_2\neq \bzero$, we must also have $x_2=0$, which arrives at a contradiction. Therefore, $\bv_1$ and $\bv_2$ are linearly independent.

Now, assume that any $j<k$ eigenvectors corresponding to distinct eigenvalues are linearly independent. We aim to show that any $j+1$ eigenvectors are also linearly independent. Suppose $\bv_1, \bv_2, \ldots, \bv_j$ are linearly independent, and $\bv_{j+1}$ is dependent on the first $j$ eigenvectors. That is, there exists a nonzero vector $\bx=[x_1,x_2,\ldots, x_{j}]\neq \bzero$ satisfying  
\begin{equation}\label{equation:independent-eigenvector-zero}
\bv_{j+1}=	x_1\bv_1+x_2\bv_2+\ldots+x_j\bv_j .
\end{equation}
Suppose the $j+1$ eigenvectors correspond to distinct eigenvalues $\lambda_1,\lambda_2,\ldots,\lambda_j,\lambda_{j+1}$.
Premultiply Equation~\eqref{equation:independent-eigenvector-zero} by $\bA$, we obtain
\begin{equation}\label{equation:independent-eigenvector-zero2}
\lambda_{j+1} \bv_{j+1} = x_1\lambda_1\bv_1+x_2\lambda_2\bv_2+\ldots+x_j \lambda_j\bv_j .
\end{equation}
Premultiply Equation~\eqref{equation:independent-eigenvector-zero} by $\lambda_{j+1}$, we get
\begin{equation}\label{equation:independent-eigenvector-zero3}
	\lambda_{j+1} \bv_{j+1} = x_1\lambda_{j+1}\bv_1+x_2\lambda_{j+1}\bv_2+\ldots+x_j \lambda_{j+1}\bv_j .
\end{equation}
Subtract Equation~\eqref{equation:independent-eigenvector-zero3} from Equation~\eqref{equation:independent-eigenvector-zero2}, we find
$$
x_1(\lambda_{j+1}-\lambda_1)\bv_1+x_2(\lambda_{j+1}-\lambda_2)\bv_2+\ldots+x_j (\lambda_{j+1}-\lambda_j)\bv_j = \bzero. 
$$
From the assumption, $\lambda_{j+1} \neq \lambda_i$ for all $i\in \{1,2,\ldots,j\}$, and $\bv_i\neq \bzero$ for all $i\in \{1,2,\ldots,j\}$. We must have $x_1=x_2=\ldots=x_j=0$, which leads to a contradiction. Therefore, the eigenvectors $\bv_1,\bv_2,\ldots,\bv_j,\bv_{j+1}$ are linearly independent. This completes the proof.
\end{proof}

An immediate consequence of the above theorem is as follows:
\begin{corollary}[Independent Eigenvector Theorem, CNT.]\label{theorem:independent-eigenvector-theorem-basis}
If a matrix $\bA\in \real^{n\times n}$ has $n$ distinct eigenvalues, then any set of $n$ corresponding eigenvectors form a basis for the space $\real^n$.
\end{corollary}

\begin{tcolorbox}[title={Symmetric Matrix Property 3 of 4},colback=\mdframecolorTheorem]
\begin{proposition}[Orthonormal Eigenvectors for Duplicate Eigenvalue]\label{proposition:eigen-multiplicity}
Let $\bA\in\real^{n\times n}$ be symmetric.
If $\bA$ has a duplicate eigenvalue $\lambda_i$ with multiplicity (Definition~\ref{definition:eigen_multipli}) $k\geq 2$, then there exist $k$ orthonormal eigenvectors corresponding to $\lambda_i$.
\end{proposition}
\end{tcolorbox}
\begin{proof}[of Proposition~\ref{proposition:eigen-multiplicity}]
We start by noting that there exists at least one unit-length eigenvector $\bx_{i1}$ corresponding to $\lambda_i$. 
Furthermore, for such an eigenvector $\bx_{i1}$, we can consistently find  $n-1$ additional orthonormal vectors $\by_2, \by_3, \ldots, \by_n$ such that $\{\bx_{i1}, \by_2, \by_3, \ldots, \by_n\}$ constitutes an orthonormal basis of $\real^n$. 
Define the matrices $\bY_1$ and $\bP_1$ as follows:
$$
\bY_1\triangleq[\by_2, \by_3, \ldots, \by_n] \qquad \text{and} \qquad \bP_1\triangleq[\bx_{i1}, \bY_1].
$$
Since $\bA$ is symmetric, we then have
$$
\bP_1^\top\bA\bP_1 = \begin{bmatrix}
	\lambda_i &\bzero \\
	\bzero & \bY_1^\top \bA\bY_1
\end{bmatrix}.
$$
Since $\bP_1$ is nonsingular and orthogonal, $\bA$ and $\bP_1^\top\bA\bP_1$ are similar matrices such that they have the same eigenvalues  (see Proposition~\ref{proposition:eigenvalue-similar-matrices}).
Using the determinant of block matrices  (Remark~\ref{remark:determinant-intermezzo}), we get: 
$$
\det(\bP_1^\top\bA\bP_1 - \lambda\bI_n) =
(\lambda_i - \lambda )\det(\bY_1^\top \bA\bY_1 - \lambda\bI_{n-1}).
$$
If $\lambda_i$ has a multiplicity of $k\geq 2$, then the term $(\lambda_i-\lambda)$ appears $k$ times in the characteristic  polynomial resulting from the determinant $\det(\bP_1^\top\bA\bP_1 - \lambda\bI_n)$, i.e., this term appears  $k-1$ times in the characteristic  polynomial from $\det(\bY_1^\top \bA\bY_1 - \lambda\bI_{n-1})$. In other words, $\det(\bY_1^\top \bA\bY_1 - \lambda_i\bI_{n-1})=0$, and $\lambda_i$ is an eigenvalue of $\bY_1^\top \bA\bY_1$ with multiplicity $k-1$. 

Let $\bB\triangleq\bY_1^\top \bA\bY_1$. Since $\det(\bB-\lambda_i\bI_{n-1})=0$, the null space of $\bB-\lambda_i\bI_{n-1}$ is non-empty. Suppose $(\bB-\lambda_i\bI_{n-1})\bn = \bzero$, i.e., $\bB\bn=\lambda_i\bn$, where $\bn$ is an eigenvector of $\bB$. 

From $
\bP_1^\top\bA\bP_1 = 
\footnotesize
\begin{bmatrix}
\lambda_i &\bzero \\
\bzero & \bB
\end{bmatrix},
$
we have $
\bA\bP_1 
\footnotesize
\begin{bmatrix}
z \\
\bn 
\end{bmatrix} 
= 
\bP_1
\footnotesize
\begin{bmatrix}
\lambda_i &\bzero \\
\bzero & \bB
\end{bmatrix}
\begin{bmatrix}
z \\
\bn 
\end{bmatrix}$, where $z$ is any scalar. 
From the left side of this equation:
\begin{equation}\label{equation:spectral-pro4-right}
\begin{aligned}
\bA\bP_1 
\begin{bmatrix}
z \\
\bn 
\end{bmatrix} 
&=
\begin{bmatrix}
\lambda_i\bx_{i1}, \bA\bY_1
\end{bmatrix}
\begin{bmatrix}
z \\
\bn 
\end{bmatrix} 
=\lambda_iz\bx_{i1} + \bA\bY_1\bn.
\end{aligned}
\end{equation}
From the right side of the equation:
\begin{equation}\label{equation:spectral-pro4-left}
\begin{aligned}
\bP_1
\begin{bmatrix}
	\lambda_i &\bzero \\
	\bzero & \bB
\end{bmatrix}
\begin{bmatrix}
	z \\
	\bn 
\end{bmatrix}
&=
\begin{bmatrix}
	\bx_{i1} & \bY_1
\end{bmatrix}
\begin{bmatrix}
	\lambda_i &\bzero \\
	\bzero & \bB
\end{bmatrix}
\begin{bmatrix}
	z \\
	\bn 
\end{bmatrix}
=
\begin{bmatrix}
\lambda_i\bx_{i1} & \bY_1\bB 
\end{bmatrix}
\begin{bmatrix}
	z \\
	\bn 
\end{bmatrix}\\
&= \lambda_i z \bx_{i1} + \bY_1\bB \bn 
=\lambda_i z \bx_{i1} + \lambda_i \bY_1 \bn,
\end{aligned}
\end{equation}
where the last equality is due to  $\bB \bn=\lambda_i\bn$.
Combining Equations~\eqref{equation:spectral-pro4-left} and \eqref{equation:spectral-pro4-right}, we obtain 
$$
\bA\bY_1\bn = \lambda_i\bY_1 \bn,
$$
which means $\bY_1\bn$ is an eigenvector of $\bA$ corresponding to the eigenvalue $\lambda_i$ (the same eigenvalue corresponding to $\bx_{i1}$). Since $\bY_1\bn$ is a linear combination of $\by_2, \by_3, \ldots, \by_n$, which are orthonormal to $\bx_{i1}$, it can be chosen to be orthonormal to $\bx_{i1}$.

To conclude, if there exists an  eigenvector, $\bx_{i1}$, corresponding to the eigenvalue  $\lambda_i$ with a multiplicity  $k\geq 2$, we can construct a second eigenvector by choosing a vector from the null space of $(\bB-\lambda_i\bI_{n-1})$, as constructed above. 

Now, suppose  we have constructed the second eigenvector $\bx_{i2}$, which is orthonormal to $\bx_{i1}$.  
For such eigenvectors $\bx_{i1}$ and $\bx_{i2}$, we can always find  $n-2$ additional orthonormal vectors $\by_3, \by_4, \ldots, \by_n$ such  that $\{\bx_{i1},\bx_{i2}, \by_3, \by_4, \ldots, \by_n\}$ forms an orthonormal basis for $\real^n$. 
Place these vectors  $\by_3, \by_4, \ldots, \by_n$ into matrix $\bY_2$ and $\{\bx_{i1},\bx_{i2},  \by_3, \by_4, \ldots, \by_n\}$ into matrix $\bP_2$:
$$
\bY_2\triangleq[\by_3, \by_4, \ldots, \by_n] \qquad \text{and} \qquad \bP_2\triangleq[\bx_{i1}, \bx_{i2},\bY_1].
$$
Since $\bA$ is symmetric, we then have
$$
\bP_2^\top\bA\bP_2 = 
\begin{bmatrix}
	\lambda_i & 0 &\bzero \\
	 0& \lambda_i &\bzero \\
	\bzero & \bzero & \bY_2^\top \bA\bY_2
\end{bmatrix}
\triangleq
\begin{bmatrix}
	\lambda_i & 0 &\bzero \\
	0& \lambda_i &\bzero \\
	\bzero & \bzero & \bC
\end{bmatrix},
$$
where $\bC\triangleq\bY_2^\top \bA\bY_2$ such that $\det(\bP_2^\top\bA\bP_2 - \lambda\bI_n) = (\lambda_i-\lambda)^2 \det(\bC - \lambda\bI_{n-2})$. If the multiplicity of $\lambda_i$ is $k\geq 3$, then $\det(\bC - \lambda_i\bI_{n-2})=0$, and the null space of $\bC - \lambda_i\bI_{n-2}$ is not empty. Thus, we can still find a vector $\bn$ from the null space of $\bC - \lambda_i\bI_{n-2}$ such that  $\bC\bn = \lambda_i \bn$. Now we can construct a vector $\footnotesize\begin{bmatrix}
	z_1 \\
	z_2\\
	\bn
\end{bmatrix}\in \real^n $, where $z_1$ and $ z_2$ are any scalar values, such that 
$$
\bA\bP_2\begin{bmatrix}
	z_1 \\
	z_2\\
	\bn
\end{bmatrix} = \bP_2 
\begin{bmatrix}
	\lambda_i & 0 &\bzero \\
	0& \lambda_i &\bzero \\
	\bzero & \bzero & \bC
\end{bmatrix}
\begin{bmatrix}
	z_1 \\
	z_2\\
	\bn
\end{bmatrix}.
$$
Similarly, from the left side of the above equation, we will get $\lambda_iz_1\bx_{i1} +\lambda_iz_2\bx_{i2}+\bA\bY_2\bn$. From the right side of the above equation, we will get $\lambda_iz_1\bx_{i1} +\lambda_i z_2\bx_{i2}+\lambda_i\bY_2\bn$. As a result, 
$$
\bA\bY_2\bn = \lambda_i\bY_2\bn,
$$
where $\bY_2\bn$ is an eigenvector of $\bA$, orthogonal to both $\bx_{i1}$ and $\bx_{i2}$. This eigenvector can also be normalized to ensure orthonormality with the first two eigenvectors.

The process can continue, ultimately yielding a set of $k$ orthonormal eigenvectors corresponding to $\lambda_i$.

In fact, the dimension of the null space of $\bP_1^\top\bA\bP_1 -\lambda_i\bI_n$ is equal to the multiplicity $k$. It also follows that if the multiplicity of $\lambda_i$ is $k$, there cannot be more than $k$ orthogonal eigenvectors corresponding to $\lambda_i$. 
If there were more than $k$, it would lead to the conclusion that there are more than $n$  orthogonal eigenvectors in $\real^n$, which is a contradiction.
\end{proof}

The existence of the spectral decomposition can be readily established based on the propositions presented above. Alternatively, we can employ the Schur decomposition to demonstrate its existence.
\begin{proof}[{of Theorem~\ref{theorem:spectral_theorem}: Existence of Spectral Decomposition}, the Second Way]
From the Schur decomposition in Theorem~\ref{theorem:schur-decomposition}, for a symmetric matrix $\bA=\bA^\top$, we have $\bQ\bU\bQ^\top = \bQ\bU^\top\bQ^\top$. 
This implies that $\bU$ is a diagonal matrix containing the eigenvalues of $\bA$. Consequently,  the columns of $\bQ$ are eigenvectors of $\bA$. We thus conclude that every symmetric matrix can be orthogonally diagonalized, even in cases of repeated eigenvalues.
\end{proof}

More compactly, we can establish the proof using the Gram-Schmidt process and induction.
\begin{proof}[{of Theorem~\ref{theorem:spectral_theorem}: Existence of Spectral Decomposition, the Third Way}]
We prove the existence of the spectral decomposition through induction. The result is self-evident for the case of  $n=1$. Assume the decomposition exists for all $n-1$. If we prove the decomposition also exists for $n$, then we complete the proof.
Suppose $\bA\in \real^{n\times n}$ has eigenvalues $\lambda_1, \lambda_2, \ldots, \lambda_n$ and $\bp_1$ is the eigenvalue corresponding to $\lambda_1$: $\bA\bp_1 = \lambda_1\bp_1$. Let $\bP$ be an orthogonal matrix with $\bp_1$ as the first column:
$$
\bP = 
[\bp_1, \bP_{n-1}],
$$
where $\bP_{n-1}\in \real^{n\times (n-1)}$ can be constructed using the Gram-Schmidt process such that $\bP_{n-1}^\top\bP_{n-1}=\bI_{n-1}$. Then we have 
$$
\bP^\top \bA\bP\triangleq
\begin{bmatrix}
\lambda_1 & 0 & \ldots & 0 \\
0 & && \\
\vdots &&\bB&\\
0 &&&
\end{bmatrix},
$$
where $\bB \in \real^{(n-1)\times (n-1)}$.
Since $\bP^\top\bA\bP$ are $\bA$ are similar matrices, they share the same eigenvalues (Proposition~\ref{proposition:eigenvalue-similar-matrices}). Therefore, $\bB$ has eigenvalues $\lambda_2, \ldots, \lambda_n$. By the induction assumption, the symmetric matrix $\bB$ (since $\bP^\top\bA\bP$ is symmetric) admits a spectral decomposition: $\bR^\top\bB\bR \triangleq \diag(\lambda_2, \ldots, \lambda_n)$. Set
$$
\bQ \triangleq
\bP
\begin{bmatrix}
	\lambda_1 & 0 & \ldots & 0 \\
	0 & && \\
	\vdots &&\bR&\\
	0 &&&
\end{bmatrix},
$$ 
then $\bQ^\top\bA\bQ = \diag(\lambda_1, \lambda_2, \ldots, \lambda_n)$.
\end{proof}

For a symmetric matrix $\bA^\top \bA$, the rank remains unchanged from that of $\bA$, a property we will utilize in proving the singular value decomposition in the following chapter.
However, it's important to note that, generally, the rank of the product of two matrices does not exceed the rank of either individual matrix.
\begin{lemma}[Rank of $\bA\bB$]\label{lemma:rankAB}
Let $\bA\in \real^{m\times n}$ and  $\bB\in \real^{n\times k}$. Then, the matrix multiplication $\bA\bB\in \real^{m\times k}$ satisfies $\rank$($\bA\bB$)$\leq\min\{\rank(\bA), \rank(\bB)\}$.
\end{lemma}
\begin{proof}[of Lemma~\ref{lemma:rankAB}]
Considering the matrix multiplication $\bA\bB$, the following observations hold:
\begin{itemize}
\item All rows of $\bA\bB$ are linear combinations of the rows of $\bB$. Therefore, the row space of $\bA\bB$ is a subset of the row space of $\bB$, implying $\rank$($\bA\bB$)$\leq$$\rank$($\bB$).

\item All columns of $\bA\bB$ are linear combinations of the columns of $\bA$. Therefore,  the column space of $\bA\bB$ is a subset of the column space of $\bA$, implying $\rank$($\bA\bB$)$\leq$$\rank$($\bA$).
\end{itemize}
Therefore, $\rank$($\bA\bB$)$\leq\min\{\rank(\bA), \rank(\bB)\}$.
\end{proof}

\begin{tcolorbox}[title={Symmetric Matrix Property 4 of 4},colback=\mdframecolorTheorem]
\begin{proposition}[Rank of Symmetric Matrices]\label{proposition:rank-of-symmetric}
Let $\bA\in\real^{n\times n}$ be an $n\times n$ real symmetric matrix. Then, $\rank(\bA)$ equals
the total number of nonzero eigenvalues of $\bA$. 
In particular, $\bA$ has full rank if and only if $\bA$ is nonsingular. Furthermore, $\cspace(\bA)$ can be equivalently represented as the linear space spanned by the eigenvectors of $\bA$ that correspond to nonzero eigenvalues.
\end{proposition}
\end{tcolorbox}
\begin{proof}[of Proposition~\ref{proposition:rank-of-symmetric}]
For any symmetric matrix $\bA$, we have $\bA$, in spectral form, as $\bA = \bQ \bLambda\bQ^\top$ and also $\bLambda = \bQ^\top\bA\bQ$. 
Using Lemma~\ref{lemma:rankAB}, we have:
\begin{itemize}
\item From $\bA = \bQ \bLambda\bQ^\top$, we have $\rank(\bA) \leq \rank(\bQ \bLambda) \leq \rank(\bLambda)$;
\item From $\bLambda = \bQ^\top\bA\bQ$, we have $\rank(\bLambda) \leq \rank(\bQ^\top\bA) \leq \rank(\bA)$, 
\end{itemize}
The inequalities  imply that  $\rank(\bA) = \rank(\bLambda)$, which corresponds to the total number of nonzero eigenvalues.

Since $\bA$ is nonsingular if and only if all  its eigenvalues are nonzero, we can affirm that $\bA$ has full rank if and only if $\bA$ is nonsingular.
\end{proof}

Similar to  eigenvalue decomposition, we can compute the $m$-th power of a matrix $\bA$ more efficiently using the spectral decomposition.
\begin{remark}[$m$-th Power]\label{remark:power-spectral}
The $m$-th power of $\bA$ is $\bA^m = \bQ\bLambda^m\bQ^\top$ if the matrix $\bA$ can be factored as the spectral decomposition $\bA=\bQ\bLambda\bQ^\top$.
\end{remark}

\section{Uniqueness of  Spectral Decomposition}\label{section:uniqueness-spectral-decomposition}
Clearly, the spectral decomposition is not unique, primarily due to the multiplicity of eigenvalues. 
When eigenvalues $\lambda_i$ and $\lambda_j$ coincide for certain $1\leq i,j\leq n$, interchanging the corresponding eigenvectors in $\bQ$ yields equivalent results, yet the decompositions themselves differ. 
However, the \textit{eigenspaces} (i.e., the null space $\nspace(\bA - \lambda_i\bI)$ corresponding to eigenvalue $\lambda_i$) remain constant for each eigenvalue. 
Therefore, there is a unique decomposition in terms of eigenspaces. Any orthonormal basis for these eigenspaces can be selected without altering the overall spectral decomposition.

\section{Other Forms, Connecting Eigenvalue Decomposition}\label{section:otherform-spectral}
In this section, we delve into various forms of spectral decomposition under different conditions. To begin, we rigorously define simple and semisimple eigenvalues.
From the proof of Proposition~\ref{proposition:eigen-multiplicity}, we observe that for symmetric matrices, the algebraic multiplicity equals the geometric multiplicity. Such matrices are referred to as simple matrices.

\begin{definition}[Simple, Semisimple Eigenvalue]\label{definition:simple_eig}
An eigenvalue $\lambda$ of a matrix $\bA$ is termed \textit{simple} if its algebraic multiplicity, $alg(\lambda)$, is 1. 
If the geometric multiplicity, $geo(\lambda)$, equals the algebraic multiplicity, $alg(\lambda)$, then $\lambda$  is referred to as a semisimple eigenvalue.
\end{definition}
The left and right eigenvectors associated with the same simple eigenvalue cannot be orthogonal; this follows from the block-diagonalization theorem (see Exercise~\ref{exercise:nonsing_eigv}).

\begin{definition}[Simple, Defective Matrix]
A matrix is termed   \textit{simple}  (or \textit{non-defective}) when the algebraic  and geometric multiplicities coincide for all its eigenvalues (i.e., all the eigenvalues are semisimple).
Conversely, a matrix is \textit{defective} if there is at least one eigenvalue that is not semisimple.
\footnote{Note that  a matrix $\bA$ is a simple matrix if $geo(\lambda)=1$ for every eigenvalue $\lambda$.
Note further that if $\lambda$ is a simple eigenvalue of $\bA$, then 	$geo(\lambda)=alg(\lambda)=1$.}
\end{definition}

\begin{exercise}
Show that the matrix $\scriptsize\begin{bmatrix}
0 & 1 \\
0 & 0
\end{bmatrix}$ is not a simple matrix.
\end{exercise}

\begin{exercise}
Proposition~\ref{proposition:diff-eigenvec-decompo} indicates any $n$-by-$n$ square matrix is diagonalizable if it has $n$ distinct eigenvalues.
Show that this is a sufficient condition for diagonalization but not a necessary condition: provide an example of a diagonalizable matrix that does not have distinct eigenvalues.
\end{exercise}
\begin{exercise}[Diagonalization of Block-Diagonal Matrices]\label{exercise:block_diagon}
Consider a block-diagonal matrix $\bA=\diag(\bA_1, \bA_2, \ldots, \bA_n)$. Show that $\bA$ is diagonalizable if and only if all of the blocks $\bA_1, \bA_2, \ldots, \bA_n$ are diagonalizable. \textit{Hint: Use induction to prove the forward implication.}
\end{exercise}

\begin{exercise}[Simultaneously Diagonalization \citep{horn2012matrix}]\label{exercise:simu_diag}
Two matrices
$\bA,\bB\in\real^{n\times n}$ are said to be  \textit{simultaneously diagonalizable} if there exists a nonsingular $\bS$ such that $\bS\bA\bS^{-1}$ and $\bS\bB\bS^{-1}$ are both diagonal. 
Given that  $\bA, \bB\in\real^{n\times n}$ be diagonalizable, show that $\bA$ and $\bB$ are simultaneously diagonalizable if and only if they commute ($\bA\bB=\bB\bA$).
This indicates that if $\bA$ and $\bB$ are symmetric, then there exists an orthogonal $\bU$ such that $\bU\bA\bU^\top$ and $\bU\bB\bU^\top$ are both diagonal if and only if $\bA\bB=\bB\bA$.
\textit{Hint: Examine the diagonalization of $\bA$ on both $\bA$ and $\bB$, show that $\bP^{-1}\bA\bP$ and $\bP^{-1}\bB\bP$ also commute, use Remark~\ref{remark:commute} and Exercise~\ref{exercise:block_diagon}.}
\end{exercise}

Diagonal matrices exhibit a straightforward structure, which simplifies the computation of determinants and inverses. 
Both the eigenvalue decomposition presented in Theorem~\ref{theorem:eigenvalue-decomposition} and the spectral decomposition in Theorem~\ref{theorem:spectral_theorem} involves matrices that are diagonalizable.
In fact, they both involve a simple matrix.

\begin{theorem}[Simple Matrices are Diagonalizable]\label{theorem:simple-diagonalizable}
A matrix is  a simple matrix if and only if it is diagonalizable.
\end{theorem}
\begin{proof}[of Theorem~\ref{theorem:simple-diagonalizable}]
Suppose that $\bA\in \real^{n\times n}$ is a simple matrix, meaning the algebraic and geometric multiplicities for each eigenvalue are equal. 
For a specific eigenvalue $\lambda_i$, let $\{\bv_1^i, \bv_2^i, \ldots, \bv_{k_i}^i\}$ be a basis for the eigenspace $\nspace(\bA - \lambda_i\bI)$. 
In other words, $\{\bv_1^i, \bv_2^i, \ldots, \bv_{k_i}^i\}$ is a set of linearly independent eigenvectors of $\bA$ associated with $\lambda_i$, where ${k_i}$ is the algebraic or geometric multiplicity of $\lambda_i$: $alg(\lambda_i)=geo(\lambda_i)=k_i$. Suppose there are $m$ distinct eigenvalues. 
Since $k_1+k_2+\ldots +k_m = n$, the set of eigenvectors consists of the union of $n$ vectors. 
Consider a  linear combination of these eigenvectors:
\begin{equation}\label{equation:proof-simple-diagonalize}
	\bz \triangleq \sum_{j=1}^{k_1} x_j^1 \bv_j^1+ \sum_{j=1}^{k_2} x_j^2 \bv_j^2 + \ldots+ \sum_{j=1}^{k_m} x_j^m \bv_j^m = \bzero.
\end{equation}
Let $\bw^i \triangleq \sum_{j=1}^{k_i} x_j^i \bv_j^i$. Then $\bw^i$ is either an eigenvector associated with $\lambda_i$ or the zero vector. Therefore, $\bz = \sum_{i=1}^{m} \bw^i$ is a sum of either zero vectors or  eigenvectors associated with different eigenvalues of $\bA$. Since eigenvectors associated with different eigenvalues are linearly independent. We must have $\bw^i =\bzero$ for all $i\in \{1, 2, \ldots, m\}$. That is,
$$
\bw^i = \sum_{j=1}^{k_i} x_j^i \bv_j^i = \bzero, \qquad \text{for all $i\in \{1, 2, \ldots, m\}$}.
$$
Since we assume the eigenvectors $\bv_j^i$'s associated with $\lambda_i$ are linearly independent, we must have $x_j^i=0$ for all $i \in \{1,2,\ldots, m\}, j\in \{1,2,\ldots,k_i\}$. Thus, the $n$ vectors are linearly independent:
$$
\{\bv_1^1, \bv_2^1, \ldots, \bv_{k_i}^1\},\{\bv_1^2, \bv_2^2, \ldots, \bv_{k_i}^2\},\ldots,\{\bv_1^m, \bv_2^m, \ldots, \bv_{k_i}^m\}.
$$
According to the eigenvalue decomposition presented in Theorem~\ref{theorem:eigenvalue-decomposition}, matrix $\bA$ can be diagonalizable.

Conversely, suppose $\bA$ is diagonalizable. That is, there exists a nonsingular matrix $\bP$ and a diagonal matrix $\bD$ such that $\bA =\bP\bD\bP^{-1} $. $\bA$ and $\bD$ are similar matrices such that they have the same eigenvalues (Proposition~\ref{proposition:eigenvalue-similar-matrices}), same algebraic multiplicities, and geometric multiplicities (Corollary~\ref{corollary:multipli-similar-matrix}). It can be readily verified that a diagonal matrix has equal algebraic  and geometric multiplicities. Therefore, $\bA$ is a simple matrix.
\end{proof}

\begin{remark}[Equivalence on Diagonalization]
Based on Theorem~\ref{theorem:independent-eigenvector-theorem}, which states that any eigenvectors corresponding to different eigenvalues are linearly independent; and Remark~\ref{remark:geometric-mul-meaning}, which states that the geometric multiplicity is the dimension of the eigenspace. We realize that if the geometric multiplicity is equal to the algebraic multiplicity (for all eigenvalues), the eigenspaces can span the entire space $\real^n$ for a  matrix $\bA\in \real^{n\times n}$. 
Therefore, the theorem mentioned above is equivalent to stating that if the eigenspaces can span the entire space $\real^n$, then $\bA$ can be diagonalizable.
\end{remark}
 
\begin{corollary}
A square matrix $\bA$ is considered simple if it has a complete set of linearly independent eigenvectors.
Additionally, if $\bA$ is symmetric, it is also a simple matrix.
\end{corollary}
The corollary presented can be easily proven through the eigenvalue decomposition outlined in Theorem~\ref{theorem:eigenvalue-decomposition} and the spectral decomposition detailed in Theorem~\ref{theorem:spectral_theorem}.

\begin{exercise}\label{exercise:diagonaliza_indep}
Show that a matrix $\bA\in\real^{n\times n}$ is diagonalizable ($\bA = \bP\bD\bP^{-1}$) if and only if there exists a set of $n$ linearly independent  eigenvectors of $\bA$. In this sense, the columns of any matrix $\bP$ are eigenvectors of $\bA$, and the diagonal entries of $\bD$ are the eigenvalues corresponding to the columns of $\bP$. (Note we have already shown the forward implication in Theorem~\ref{theorem:eigenvalue-decomposition}.)
\end{exercise}

\subsection*{Block-Diagonal Matrices}
We provide the diagonalization result for a block-diagonal matrix.
\begin{theoremHigh}[Block-Diagonal Matrices]\label{theorem:block_diagonal_spec}
Given a matrix $\bA\in \real^{n\times n}$, it takes the form of a  block-diagonal (quasi-diagonal) matrix as shown below:
$$
\bA=
\begin{bmatrix}
\bA_1 &&& \\
&\bA_2&& \\
&&\ddots& \\
&&&\bA_p 
\end{bmatrix},
$$
where $\bA_j\in\real^{n_j\times n_j}$ and $\sum_{j=1}^{p}n_j = n$. Then $\bA$ is diagonalizable if and only if every matrix $\bA_j$ is diagonalizable for all $j\in\{1,2,\ldots, p\}$.
\end{theoremHigh}
\begin{proof}[of Theorem~\ref{theorem:block_diagonal_spec}]
Suppose the block-diagonal matrix $\bA$ is diagonalizable. Then, there exists a nonsingular matrix $\bP\in \real^{n\times n}$ such that $\bP^{-1}\bA \bP = \bD$, where $\bD$ is a  diagonal matrix: $\bD\in \diag(\lambda_1, \lambda_2, \ldots, \lambda_n)$. 
Suppose further each column of $\bP$ is partitioned into $p$ blocks with each block $\bp_i^j$ of size $\real^{n_j}$. Then, the $\bP$ has the following form:
$$
\bP = [\bp_1, \bp_2, \ldots, \bp_n]
=
\begin{bmatrix}
\bp_1^1 & \bp_2^1 & \ldots & \bp_n^1 \\
\bp_1^2 & \bp_2^2 & \ldots & \bp_n^2 \\
\vdots & \vdots & \ddots & \vdots \\
\bp_1^p & \bp_2^p & \ldots & \bp_n^p \\
\end{bmatrix},
\gap 
\text{with  }
\bp_i^j \in \real^{n_j}.
$$
Therefore, from $\bA\bp_i = \lambda_i\bp_i$, we have 
$$
\begin{bmatrix}
	\bA_1 &&& \\
	&\bA_2&& \\
	&&\ddots& \\
	&&&\bA_p 
\end{bmatrix}
\begin{bmatrix}
\bp_i^1 \\
\bp_i^2 \\
\vdots \\
\bp_i^p
\end{bmatrix}
=
\lambda_i
\begin{bmatrix}
	\bp_i^1 \\
	\bp_i^2 \\
	\vdots \\
	\bp_i^p
\end{bmatrix},
\gap 
\text{for all } i\in \{1,2,\ldots,n\}.
$$
This implies that
$\bA_j \bp_i^j = \lambda_i\bp_i^j$ for all $i\in\{1,2,\ldots,n\}, j\in\{1,2,\ldots,p\}$.
Suppose further that $\bB^j = [\bp_1^j, \bp_2^j, \ldots , \bp_n^j ] \in \real^{n_j\times n}$ such that 
$$
\bP = 
\footnotesize
\begin{bmatrix}
	\bB^1 \\
\bB^2 \\
\vdots \\
\bB^p
\end{bmatrix}.
$$
The rank of $\bB^j$ is $n_j$ since $\bP$ is nonsingular. Therefore, $\bA_j$ is diagonalizable by Exercise~\ref{exercise:diagonaliza_indep}.

For the converse implication, suppose all $\bA_j$ for $j\in\{1,2,\ldots,p\}$ are diagonalizable. That is, there exists a nonsingular $\bQ_j$ such that $\bQ_j^{-1}\bA_j\bQ_j$ is a diagonal matrix. Let 
$$
\bQ\triangleq
\begin{bmatrix}
	\bQ1 &&& \\
	&\bQ_2&& \\
	&&\ddots& \\
	&&&\bQ_p 
\end{bmatrix}.
$$
It then follows that $\bQ^{-1}\bA\bQ$ is a diagonal matrix.
\end{proof}

\subsection*{Spectrum of Special Matrices}
For special matrices, the spectrum (Definition~\ref{definition:spectrum}) can  lie either on or inside the unit circle \citep{simovici2008mathematical}. To see this, we will rely on the following result, which is left as an exercise.
\begin{exercise}\label{exercise:block_trian_det}
Consider the block upper triangular matrix $\bA\in\real^{n\times n}$ and the block lower triangular matrix $\bB\in\real^{n\times n}$ with the following forms:
$$
\bA=
\begin{bmatrix}
	\bA_{11} & \bA_{12}&\ldots &\bA_{1p} \\
	&\bA_{22}& \ldots & \bA_{2p} \\
	&&\ddots& \\
	&&&\bA_{pp} 
\end{bmatrix},
\gap 
\bB=
\begin{bmatrix}
	\bA_{11} & &\ldots & \\
\bA_{21}	&\bA_{22}& \ldots &  \\
\vdots	&\vdots &\ddots& \\
\bA_{p1}	&\bA_{p2}&&\bA_{pp} 
\end{bmatrix}.
$$
Show that the determinant $\det(\bA)=\det(\bB) = \det(\bA_{11})\det(\bA_{22})\ldots\det(\bA_{pp})$.
\textit{Hint: Use induction}.
\end{exercise}

\begin{theorem}[Block Upper Triangular Matrices]\label{theorem:block_upper_spectrum}
Given the block upper triangular matrix $\bA\in\real^{n\times n}$ with the following form:
$$
\bA=
\begin{bmatrix}
\bA_{11} & \bA_{12}&\ldots &\bA_{1p} \\
&\bA_{22}& \ldots & \bA_{2p} \\
&&\ddots& \\
&&&\bA_{pp} 
\end{bmatrix},
$$
where $\bA_{jj}\in \real^{n_j\times n_j}$ for all $j\in\{1,2,\ldots, p\}$. Then the spectrum  $\Lambda(\bA) = \cup_{j=1}^p \Lambda(\bA_{jj})$. The same result holds for block lower triangular matrices.
\end{theorem}
\begin{proof}[of Theorem~\ref{theorem:block_upper_spectrum}]
It can be shown that (by Exercise~\ref{exercise:block_trian_det})
$$
\det(\lambda\bI-\bA) = 
\prod_{j=1}^{p}\det(\lambda \bI_{n_j}-\bA_{jj}).
$$
This completes the proof.
\end{proof}
From the theorem above, it becomes evident that the spectrum of any upper (lower) triangular matrix is the set of its diagonal elements.

Consider a matrix $\bA$ and the eigenpair $(\lambda, \bx)$: $\bA\bx = \lambda\bx$. We have $\bA^k\bx = \lambda^k\bx$ for $k\geq 1$. Therefore, if $\lambda\in\Lambda(\bA)$, then $\lambda^k\in \Lambda(\bA^k)$ for $k\geq 1$.

\begin{theorem}[Nilpotent Matrices]\label{theorem:spec_nilp}
Let $\bA\in\real^{n\times n}$ be a nilpotent matrix (Definition~\ref{definition:niopotent_mat}). Then, $\Lambda(\bA)=\{0\}$.
\end{theorem}
\begin{proof}[of Theorem~\ref{theorem:spec_nilp}]
Suppose the nilpotency is $\nilp(\bA)=k$ and $\bA$ has an eigenvalue $\lambda$. We have $\lambda^k\in \Lambda(\bA^k)=\Lambda(\bzero) = \{0\}$. Therefore, any eigenvalue $\lambda$ of $\bA$ must be 0.
\end{proof}

\begin{theorem}[Idempotent Matrices]\label{theorem:spec_idem}
Let $\bA\in\real^{n\times n}$ be a idempotent matrix (Definition~\ref{definition:Idempotent_mat}). Then, $\Lambda(\bA)=\{0,1\}$.
\end{theorem}
\begin{proof}[of Theorem~\ref{theorem:spec_idem}]
Given the eigenpair $(\lambda, \bx)$: $\bA\bx = \lambda\bx$. We have $\bA^2\bx = \bA\bx = \lambda\bx$. And it also follows that $\bA^2\bx = \bA(\bA\bx) = \lambda\bA\bx = \lambda^2\bx$. This implies $\lambda^2=\lambda$, and thus $\lambda\in\{0,1\}$.
\end{proof}

\begin{theorem}[Stochastic Matrices]\label{theorem:spec_stoc}
Let $\bA\in\real^{n\times n}$ be a stochastic matrix (Definition~\ref{definition:stochastic_mat}). Then, at least one eigenvalue of $\bA$ is equal to 1, and all eigenvalues lie on or inside the unit circle.
\end{theorem}
\begin{proof}[of Theorem~\ref{theorem:spec_stoc}]
Since all rows of $\bA$ sum to 1 (with nonnegative elements), it can be shown that the eigenvector $\bone$  corresponds to the eigenvalue $\lambda=1$. Given the eigenpair $\bA\bx = \lambda\bx$, we have 
$$
\lambda x_i = \sum_{j=1}^{n} a_{ij} x_j, \gap 
\forall i \in\{1,2,\ldots, n\}
\gap\implies \gap 
|\lambda| |x_i| = \left\vert \sum_{j=1}^{n} a_{ij} x_j\right\vert
\leq 
\sum_{j=1}^{n} a_{ij} |x_j|,
$$
since the elements $\{a_{ij}\}$ are nonnegative.
Suppose $x_p$ has the largest magnitude for some $p\in\{1,2,\ldots,n\}$: $|x_p| = \max\{|x_i| : 1\leq i\leq n\}$. We get 
$
|\lambda| \leq 
\sum_{j=1}^{n} a_{pj} \frac{|x_j|}{|x_p|}
\leq \sum_{j=1}^{n} a_{pj} =1.
$
Thus, all eigenvalues of $\bA$ lie on or inside the unit cycle.
\end{proof}

\begin{exercise}[Orthogonal Matrices]\label{exercise:eig_orthos}
Show that all eigenvalues of an orthogonal matrix are located on the unit cycle.
\end{exercise}

%

\subsection*{Other Forms}
We are now prepared to present the second form of the spectral decomposition.
\begin{theoremHigh}[Spectral Decomposition: The Second Form]\label{theorem:spectral_theorem_secondForm}
Let $\bA \in \real^{n\times n}$ be a \textbf{simple matrix}. Then it  can be factored as a sum of a set of idempotent matrices
$$
\bA = \sum_{i=1}^{n} \lambda_i \bA_i,
$$
where each $\lambda_i$, for  $i\in \{1,2,\ldots, n\}$, represents an eigenvalue of $\bA$ (duplicates are possible). The idempotent matrices $\bA_i$ satisfy the following properties:
\begin{enumerate}
\item  \textit{Idempotent.} $\bA_i^2 = \bA_i$ for all $i\in \{1,2,\ldots, n\}$;
\item \textit{Orthogonal.} $\bA_i\bA_j = \bzero$ for all $i \neq j$;
\item \textit{Additivity.} $\sum_{i=1}^{n} \bA_i = \bI_n$;
\item \textit{Rank-Additivity.} $\rank(\bA_1) + \rank(\bA_2) + \ldots + \rank(\bA_n) = n$.
\end{enumerate}

\end{theoremHigh}

\begin{proof}[of Theorem~\ref{theorem:spectral_theorem_secondForm}]
Since $\bA$ is a simple matrix, according to Theorem~\ref{theorem:simple-diagonalizable}, there exists a nonsingular matrix $\bP$ and a diagonal matrix $\bLambda$ such that $\bA=\bP\bLambda\bP^{-1}$, where $\bLambda=\diag(\lambda_1, \lambda_2, \ldots, \lambda_n)$,  $\lambda_i$'s are eigenvalues of $\bA$, and the columns of $\bP$ consist of the corresponding eigenvectors. Suppose
$$
\bP \triangleq \begin{bmatrix}
\bv_1 & \bv_2&\ldots & \bv_n
\end{bmatrix}
\qquad
\text{and }
\qquad
\bP^{-1} \triangleq 
\footnotesize
\begin{bmatrix}
\bw_1^\top \\
\bw_2^\top \\
\vdots \\
\bw_n^\top
\end{bmatrix}
$$
are the column and row partitions of $\bP$ and $\bP^{-1}$, respectively. Then, we can rewrite $\bA$ as
$$
\bA= \bP\bLambda\bP^{-1} = 
\begin{bmatrix}
	\bv_1 & \bv_2&\ldots & \bv_n
\end{bmatrix}
\bLambda
\begin{bmatrix}
	\bw_1^\top \\
	\bw_2^\top \\
	\vdots \\
	\bw_n^\top
\end{bmatrix}=
\sum_{i=1}^{n}\lambda_i \bv_i\bw_i^\top.
$$
By defining $\bA_i \triangleq \bv_i\bw_i^\top$, we obtain $\bA = \sum_{i=1}^{n} \lambda_i \bA_i$.
It follows from $\bP^{-1}\bP = \bI$  that 
$$ 
\left\{
\begin{aligned}
	&\bw_i^\top\bv_j = 1 ,& \mathrm{\,\,if\,\,} i = j;  \\
	&\bw_i^\top\bv_j = 0 ,& \mathrm{\,\,if\,\,} i \neq j. 
\end{aligned}
\right.
$$
Therefore, 
$$ 
\bA_i\bA_j =\bv_i\bw_i^\top\bv_j\bw_j^\top = \left\{
\begin{aligned}
	&\bv_i\bw_i^\top = \bA_i ,& \mathrm{\,\,if\,\,} i = j;  \\
	& \bzero ,& \mathrm{\,\,if\,\,} i \neq j. 
\end{aligned}
\right.
$$
This confirms  the idempotency and orthogonality of $\bA_i$'s. 
Additionally, we observe that $\sum_{i=1}^{n}\bA_i = \bP\bP^{-1}=\bI$, implying the additivity of $\bA_i$'s. 
The rank-additivity property holds trivially because  $\rank(\bA_i)=1$ for all $i\in \{1,2,\ldots, n\}$.
\end{proof}
The decomposition is closely related to  \textit{Cochran's theorem}, which is detailed in  Appendix~\ref{appendix:cochran-theorem}. It plays a significant role in the distribution theory of linear models, as discussed in \citet{lu2021rigorous}.

Continuing, if we assume the presence of $k$ distinct eigenvalues for a simple matrix $\bA$, we can derive the following result.
\begin{theoremHigh}[Spectral Decomposition: The Third Form]\label{Corollary:spectral_theorem_3Form}
Let $\bA \in \real^{n\times n}$ be a  \textbf{simple matrix}  \textcolor{mylightbluetext}{with $k$ distinct eigenvalues}. Then, it can be factored as a sum of a set of idempotent matrices:
$$
\bA = \sum_{i=1}^{\textcolor{mylightbluetext}{k}} \lambda_i \bA_i,
$$
where each $\lambda_i$, for  $i\in \{1,2,\ldots, \textcolor{mylightbluetext}{k}\}$,  represents one of the distinct eigenvalues of $\bA$. 
The idempotent matrices $\bA_i$ satisfy the following properties:
\begin{enumerate}
\item \textit{Idempotent.} $\bA_i^2 = \bA_i$ for all $i\in \{1,2,\ldots, \textcolor{mylightbluetext}{k}\}$;
\item \textit{Orthogonal.} $\bA_i\bA_j = \bzero$ for all $i \neq j$;
\item \textit{Additivity.} $\sum_{i=1}^{\textcolor{mylightbluetext}{k}} \bA_i = \bI_n$;
\item \textit{Rank-Additivity.} $\rank(\bA_1) + \rank(\bA_2) + \ldots + \rank(\bA_{\textcolor{mylightbluetext}{k}}) = n$.
\end{enumerate}

\end{theoremHigh}
\begin{proof}[of Theorem~\ref{Corollary:spectral_theorem_3Form}]
According to Theorem~\ref{theorem:spectral_theorem_secondForm}, we can decompose $\bA$ as $\bA =\sum_{j=1}^{n} \beta_j \bB_j$, where $\beta_j$ are the eigenvalues and $\bB_j$ are the corresponding idempotent matrices from the second form of the spectral decomposition. 
Assume without loss of generality that the eigenvalues are ordered such that $\beta_1 \leq \beta_2 \leq \ldots \leq\beta_n$, where duplicates are possible. Let $\{\lambda_1, \lambda_2, \ldots, \lambda_k\}$ denote the set of $k$ distinct eigenvalues, and let  $\bA_i$ represent the sum of the $\bB_j$ matrices associated with $\lambda_i$. 
Suppose the multiplicity of $\lambda_i$ is $m_i$ for $i\in\{1,2,\ldots,k\}$, and the $\bB_j$'s associated with $\lambda_i$ can be denoted by $\{\bB_{1}^i, \bB_{2}^i, \ldots, \bB_{m_i}^i\}$. Thus, $\bA_i$ can be defined as $\bA_i = \sum_{j=1}^{m_i} \bB_{j}^i$. Consequently, $\bA = \sum_{i=1}^{k} \lambda_i \bA_i$.

\paragraph{Idempotency.} $\bA_i^2 = (\bB_1^i + \bB_2^i+\ldots \bB_{m_i}^i)(\bB_1^i + \bB_2^i+\ldots \bB_{m_i}^i)= \bB_1^i + \bB_2^i+\ldots \bB_{m_i}^i = \bA_i$  due to the idempotency and orthogonality of the $\bB_j^i$ matrices.

\paragraph{Orthogonality.} $\bA_i\bA_j = (\bB_1^i + \bB_2^i+\ldots \bB_{m_i}^i)(\bB_1^j + \bB_2^j+\ldots \bB_{m_j}^j)=\bzero$  due to the orthogonality of the $\bB_j^i$ matrices.

\paragraph{Additivity.} It is evident that  $\sum_{i=1}^{k} \bA_i = \bI_n$.

\paragraph{Rank-Additivity.} $\rank(\bA_i ) = \rank(\sum_{j=1}^{m_i} \bB_{j}^i) = m_i$ such that $\rank(\bA_1) + \rank(\bA_2) + \ldots + \rank(\bA_{k}) = m_1+m_2+\ldots+m_k=n$.
\end{proof}

\begin{theoremHigh}[Spectral Decomposition: Backward Implication]\label{Corollary:spectral_theorem_4Form}
If a matrix $\bA \in \real^{n\times n}$ with $k$ distinct eigenvalues can be factored as a sum of a set of idempotent matrices
$
\bA = \sum_{i=1}^{k} \lambda_i \bA_i,
$
where each $\lambda_i$, for  $i\in \{1,2,\ldots, k\}$, represents one of   the distinct eigenvalues of $\bA$, and 
$\bA_i$'s satisfy the four conditions outlined in Theorem~\ref{Corollary:spectral_theorem_3Form}, then the matrix $\bA$ is a simple matrix.
\end{theoremHigh}
\begin{proof}[of Theorem~\ref{Corollary:spectral_theorem_4Form}]
Assume that $\rank(\bA_i) = r_i$ for $i \in \{1,2,\ldots, k\}$. By  ULV decomposition in Theorem~\ref{theorem:ulv-decomposition}, $\bA_i$ can be factored as 
$
\bA_i = 
\scriptsize
\bU_i \begin{bmatrix}
\bL_i & \bzero \\
\bzero & \bzero 
\end{bmatrix}\bV_i,
$
where $\bL_i \in \real^{r_i \times r_i}$ is lower triangular, and $\bU_i \in \real^{n \times n}$ and $\bV_i\in \real^{n \times n}$ are orthogonal. Define 
$$
\bX_i \triangleq
\bU_i \begin{bmatrix}
\bL_i  \\
\bzero  
\end{bmatrix}
\qquad 
\text{and}
\qquad 
\bV_i \triangleq 
\begin{bmatrix}
\bY_i \\
\bZ_i
\end{bmatrix},
$$
where $\bX_i$ is of size $\real^{n\times r_i}$, and $\bY_i \in \real^{r_i \times n}$ consists of the first $r_i$ rows of $\bV_i$. Consequently, we have
$$
\bA_i = \bX_i \bY_i.
$$
This can be viewed as a \textit{reduced} ULV decomposition of $\bA_i$. By appending the $\bX_i$'s and $\bY_i$'s into matrices $\bX$ and $\bY$, we obtain:
$$
\bX \triangleq [\bX_1, \bX_2, \ldots, \bX_k]
\qquad
\text{and}\qquad
\bY \triangleq
\footnotesize
\begin{bmatrix}
\bY_1\\
\bY_2\\
\vdots \\
\bY_k
\end{bmatrix},
$$
where $\bX\in \real^{n\times n}$ and $\bY\in \real^{n\times n}$ (due to the rank-additivity property). Utilizing  block matrix multiplication and the additivity of $\bA_i$'s, we have 
$
\bX\bY = \sum_{i=1}^{k} \bX_i\bY_i = \sum_{i=1}^{k} \bA_i = \bI.
$
Therefore, $\bY$ is the inverse of $\bX$, and 
$$
\bY\bX = 
\begin{bmatrix}
	\bY_1\\
	\bY_2\\
	\vdots \\
	\bY_k
\end{bmatrix}
[\bX_1, \bX_2, \ldots, \bX_k]
=
\begin{bmatrix}
\bY_1\bX_1 & \bY_1\bX_2 & \ldots & \bY_1\bX_k\\
\bY_2\bX_1 & \bY_2\bX_2 & \ldots & \bY_2\bX_k\\
\vdots & \vdots & \ddots & \vdots\\
\bY_k\bX_1 & \bY_k\bX_2 & \ldots & \bY_k\bX_k\\
\end{bmatrix}
=\bI,
$$
such that 
$$ 
\bY_i\bX_j = \left\{
\begin{aligned}
	&\bI_{r_i}  ,& \mathrm{\,\,if\,\,} i = j;  \\
	& \bzero ,& \mathrm{\,\,if\,\,} i \neq j. 
\end{aligned}
\right.
$$
This implies 
$$ 
\bA_i\bX_j = \left\{
\begin{aligned}
	&\bX_i  ,& \mathrm{\,\,if\,\,} i = j;  \\
	& \bzero ,& \mathrm{\,\,if\,\,} i \neq j, 
\end{aligned}
\right.
\qquad 
\text{and}
\qquad
\bA \bX_i = \lambda_i\bX_i.
$$
Finally, we have 
$$
\begin{aligned}
\bA\bX &= \bA[\bX_1, \bX_2, \ldots, \bX_k] = [\lambda_1\bX_1, \lambda_2\bX_2, \ldots, \lambda_k\bX_k] = \bX\bLambda,
\end{aligned}
$$
where 
$$\bLambda =
\begin{bmatrix}
\lambda_1 \bI_{r_1} & \bzero & \ldots & \bzero \\
\bzero &  \lambda_2 \bI_{r_2} & \ldots & \bzero \\
\vdots &  \vdots & \ddots & \vdots \\
\bzero & \bzero & \ldots &  \lambda_k \bI_{r_k} \\
\end{bmatrix}
$$
is a diagonal matrix. This implies that $\bA$ can be diagonalized, and by Theorem~\ref{theorem:simple-diagonalizable}, $\bA$ is indeed a simple matrix.
\end{proof}

\begin{corollary}[Forward and Backward Spectral]
Combining Theorem~\ref{Corollary:spectral_theorem_3Form} and Theorem~\ref{Corollary:spectral_theorem_4Form}, we can claim that a matrix $\bA \in \real^{n\times n}$ is a simple matrix with $k$ distinct eigenvalues if and only if it can be factored as a sum of a set of idempotent matrices
$
\bA = \sum_{i=1}^{k} \lambda_i \bA_i,
$
where each $\lambda_i$, for  $i\in \{1,2,\ldots, k\}$, represents one of the distinct eigenvalues of $\bA$, and $\bA_i$'s satisfy the four conditions outlined in Theorem~\ref{Corollary:spectral_theorem_3Form}.
\end{corollary}

\section{Rank-Revealing Orthogonal Similar Decomposition}
The UTV decomposition discussed in Chapter~\ref{section:ulv-urv-decomposition} provides an orthogonal biequivalent decomposition  for rectangular matrices that reveals the rank of the matrix. Analogously, the rank-revealing orthogonal similar (RROS) decomposition provides such a decomposition for square matrices.
\begin{theoremHigh}[Rank-Revealing Orthogonal Similar (RROS) Decomposition]\label{theorem:rrorthosimi}
Let $\bA\in\real^{n\times n}$ be any square matrix  with rank $r$ and real eigenvalues. Then, it can be factored as 
$$
\bQ\bT\bQ^{\top},~\footnote{When $\bA$ is complex or the eigenvalues are complex, then $\bQ$ is unitary and $\bT$ is also complex.}
$$
where $\bQ\in\real^{n\times n}$ is orthogonal, and $\bT\in\real^{n\times n}$ is upper triangular with the first $r$ rows being linearly independent and the last $n-r$ rows being zero.
\end{theoremHigh}
\begin{proof}[of Theorem~\ref{theorem:rrorthosimi}]
Consider $\bU=[\bU_1,\bU_2]\in\real^{n\times n}$, where the columns of $\bU_1\in\real^{n\times r}$ form an orthonormal basis for the column space of $\bA$, and the columns of $\bU_2\in\real^{n\times (n-r)}$ form an orthonormal basis for the null space of $\bA^\top$. Therefore, 
$$
\bU^\top\bA\bU=
\begin{bmatrix}
\bU_1^\top\bA\bU_1 & \bU_1^\top\bA\bU_2 \\
\bzero & \bzero 
\end{bmatrix}.
$$
Suppose $\bU^\top\bA\bU$ admits the Schur decomposition $\bU_1^\top\bA\bU_1=\bV\bR\bV^\top$ (when $\bU_1^\top\bA\bU_1$ has real eigenvalues, $\bV$ and $\bR$ are real; otherwise, $\bV$ and $\bR$ are complex). 
Let $\bP\triangleq\diag(\bV, \bI_{n-r})$, it follows that 
$$
(\bU\bP)^\top \bA \underbrace{(\bU\bP)}_{\bQ}
=
\begin{bmatrix}
\bV^\top & \bzero \\
\bzero & \bI 
\end{bmatrix}
\bU^\top \bA\bU 
\begin{bmatrix}
	\bV & \bzero \\
	\bzero & \bI 
\end{bmatrix}
=
\begin{bmatrix}
\bR & \bV^\top\bU_1^\top\bA\bU_2\\
\bzero & \bzero 
\end{bmatrix}.
$$
This completes the proof.
\end{proof}

\begin{exercise}
Find the decomposition $\bA=\bQ\bT\bQ^\top$ where $\bT$ is \textbf{lower} triangular with the  first $r$ columns being linearly independent and the last $n-r$ columns being zero. \textit{Hint: Find the Schur decomposition of $\bU_1^\top\bA^\top\bU_1$ in the above proof.}
\end{exercise}

\section{Block-Triangularization of Matrices with (Real) Eigenvalues}
The spectral decomposition applies to symmetric matrices, which have real eigenvalues. For general matrices with real eigenvalues, a block-triangularization can be achieved, a result derived from the Schur decomposition.

\begin{theoremHigh}[Block-Triangularization of Matrices with (Real) Eigenvalues]\label{theorem:block_real_eigens}
Let $\bA \in \real^{n\times n}$ be a {real}  matrix  with real eigenvalues. Then, it  can be factored as
\begin{equation*}
	\bA = \bP \bT \bP^{-1},~\footnote{When $\bA$ is complex or the eigenvalues are not real, then $\bP$ and $\bT$ can be complex as well.}
\end{equation*}
where $\bP\in \real^{n\times n}$ is an $n\times n$ nonsingular matrix, and $\bT\in \real^{n\times n}$ is a special upper triangular (a quasi-diagonal)  matrix with the following form
$$
\bT = 
\diag\left( \bT_{11}, \bT_{22}, \ldots, \bT_{pp}\right)
=
\begin{bmatrix}
	\bT_{11} &&& \\
	& \bT_{22} && \\
	&&\ddots & \\
	&&&\bT_{pp}
\end{bmatrix},
$$
where each $\bT_{ii}\in\real^{n_i\times n_i}$ is an upper triangular matrix whose diagonal values are $\lambda_i$ (the algebraic multiplicity of $\lambda_i$ is $n_i$, i.e., $\bT_{ii}=\lambda_i\bI_{n_i}+\bR_{ii}$, with $\bR_{ii}$ being strictly upper triangular), and $\sum_{i=1}^{p}n_i=n$.
(Without loss of generality, we assume the eigenvalues are ordered in descending order $\lambda_1>\lambda_2>\ldots>\lambda_p$).
\end{theoremHigh}

\begin{proof}
From the Schur decomposition (Theorem~\ref{theorem:schur-decomposition} and Corollary~\ref{corollary:schur-second-form}), we reorder the eigenvalues in descending order: $\bA=\bQ\bU\bQ^{-1}$, where $\bU$ is upper triangular with the block formulation:
$$
\bU=
[\bT_{ij}]_{i,j=1}^{p}
=
\begin{bmatrix}
\bT_{11} & \bT_{12} & \ldots & \bT_{1p}\\
\bzero  & \bT_{22} & \ldots & \bT_{2p}\\
\bzero & \bzero & \ddots & \bT_{p-1,p}\\
\bzero & \bzero & \ldots & \bT_{pp}\\
\end{bmatrix}=
\begin{bmatrix}
\bT_{11} & \bB \\
\bzero & \bC
\end{bmatrix},
$$
where each $\bT_{ii}$ is upper triangular with the same eigenvalue $\lambda_i$ along the diagonal.
Sylvester's theorem (Theorem~\ref{theorem:sylvesters_theorem}) indicates that $\bT_{11}\bX-\bX\bC=-\bB$ has a unique solution. Therefore, we may construct $\bP_1$
$$
\bP_1=
\begin{bmatrix}
\bI & \bX \\
\bzero & \bI
\end{bmatrix}
\gap 
\text{and}
\gap 
\bP_1^{-1}=
\begin{bmatrix}
\bI & -\bX \\
\bzero & \bI
\end{bmatrix}
\gap \implies \gap
\bP_1^{-1} \bU\bP_1
=
\begin{bmatrix}
\bT_{11} & \bzero \\
\bzero & \bC
\end{bmatrix}.
$$
Continuing this process with $p-1$ reductions, we can find $p-1$ nonsingular matrices $\bP_1, \bP_2, \ldots, \bP_{p-1}$ and let $\bP\triangleq\bQ\bP_1\bP_2\ldots\bP_{p-1}$. This completes the proof.
\end{proof}

\begin{exercise}[Non-Orthogonality of Eigenvectors]\label{exercise:nonsing_eigv}
Let $\bA\in\real^{n\times n}$ be given, and let $\bu$ be a right eigenvector  and $\bv$ be a left eigenvector of $\bA\in\real^{n\times n}$, respectively, corresponding to the same simple eigenvalue $\lambda$ (Definition~\ref{definition:simple_eig}). Show that $\bu^\top\bv\neq 0$.
\textit{Hint: Construct the block-diagonalization $\bA=\bP\bT\bP^{-1}$, show that the first column of $\bP$ is $\bu$, and the first row of $\bP^{-1}$ is $\bv^\top$. Examine the value of $\bu^\top\bv$ in terms of $\bP$ and $\bP^{-1}$.}
\end{exercise}

\begin{exercise}
Let $\bA\in\real^{n\times n}$ be given, and let $\bu$ be a right eigenvector  and $\bv$ be a left eigenvector of $\bA\in\real^{n\times n}$, respectively, corresponding to the same simple eigenvalue $\lambda$. Show that $\bA-\lambda\bI+\kappa \bu\bv^\top$ is nonsingular for all $\kappa\neq 0$.
\textit{Hint: Use Exercise~\ref{exercise:nonsing_eigv}; consider  $\bx$ satisfying $(\bA-\lambda\bI+\kappa \bu\bv^\top)\bx=\bzero$, show that $\bx=c\bu$, and show that $c$ must be zero.}
\end{exercise}

\section{Block-Diagonalization of Skew-Symmetric Matrices*}
We have introduced the spectral decomposition for symmetric matrices. Another type of matrix closely related to symmetric matrices is the \textit{skew-symmetric matrix}: $\bA^\top = -\bA$ (Definition~\ref{definition:speci_mat}). Notably, under this definition,  the diagonal values $a_{ii}$ for all $i \in \{1,2,\ldots, n\}$ satisfy the condition $a_{ii} = -a_{ii}$, implying that all the diagonal components are zero.

\begin{exercise}[Quadratic of Skew-Symmetric]
	Let $\bA\in\real^{n\times n}$ be  skew-symmetric. Show that
	$
	\bx^\top\bA\bx = 0,  \forall\, \bx\in\real^n.
	$
\end{exercise}

We have established in Proposition~\ref{proposition:real-eigenvalues-spectral}  that all  eigenvalues of symmetric matrices are real. Similarly, it can be demonstrated that all  eigenvalues of skew-symmetric matrices are imaginary or zero.

\begin{lemma}[Imaginary Eigenvalues]\label{lemma:real-eigenvalues-spectral-skew}
All eigenvalues of any skew-symmetric matrix are either purely imaginary or zero.~\footnote{The result can also be extended  to skew-Hermitian matrices; see Problem~\ref{problem:real_skewherm}.}
\end{lemma}
\begin{proof}[of Lemma~\ref{lemma:real-eigenvalues-spectral-skew}]
Suppose   $\lambda=a+ib$ is a complex eigenvalue of the skew-symmetric matrix $\bA$ corresponding to the complex eigenvector $\bx = \bc+i\bd$, where $a$, $b$, $\bc$, and $\bd$ are real numbers or vectors. The complex conjugate of $\lambda$ is $\overline{\lambda}=a-ib$, and the complex conjugate of $\bx$ is $\overline{\bx}=\bc-i\bd$. The eigenvalue equation and its conjugate can be written as:
$$
\bA \bx = \lambda \bx\qquad   \underrightarrow{\text{ leads to }}\qquad  \bA \overline{\bx} = \overline{\lambda} \overline{\bx}\qquad   \underrightarrow{\text{ transpose to }}\qquad  \overline{\bx}^\top \bA^\top =\overline{\lambda} \overline{\bx}^\top.
$$
Taking the dot product of the first equation with $\overline{\bx}$ and the last equation with $\bx$, we get:
$$
\overline{\bx}^\top \bA \bx = \lambda \overline{\bx}^\top \bx, \qquad \text{and } \qquad \overline{\bx}^\top \bA^\top \bx = \overline{\lambda}\overline{\bx}^\top \bx.
$$
Then we have the equality $-\lambda\overline{\bx}^\top \bx = \overline{\lambda} \overline{\bx}^\top\bx$ (since $\bA^\top=-\bA$). Since $\overline{\bx}^\top\bx = (\bc-i\bd)^\top(\bc+i\bd) = \bc^\top\bc+\bd^\top\bd$ is a real number, the real part of $\lambda$ must be zero. Therefore, $\lambda$ is  either purely imaginary or zero.
\end{proof}

\begin{proposition}[Odd Skew-Symmetric Determinant]\label{proposition:skew-symmetric-determinant}
Let $\bA\in \real^{n\times n}$ be a skew-symmetric matrix. If $n$ is odd, then $\det(\bA)=0$.
\end{proposition}
\begin{proof}[of Proposition~\ref{proposition:skew-symmetric-determinant}]
When $n$ is odd, we have 
$$
\det(\bA) = \det(\bA^\top) = \det(-\bA) = (-1)^n \det(\bA) = -\det(\bA),
$$
which can only be true if $\det(\bA)=0$.
\end{proof}

\begin{theoremHigh}[Block-Diagonalization of Skew-Symmetric Matrices]\label{theorem:skew-block-diagonalization_theorem}
Let  $\bA \in \real^{n\times n}$ be a real skew-symmetric matrix. Then, it can be factored as
\begin{equation*}
\bA = \bZ \bD \bZ^\top,
\end{equation*}
where $\bZ\in \real^{n\times n}$ is an $n\times n$ nonsingular matrix, and $\bD\in \real^{n\times n}$ is a block-diagonal (quasi-diagonal) matrix with the following form
$$
\bD = 
\diag\left(\begin{bmatrix}
	0 & 1 \\
	-1 & 0
\end{bmatrix}, 
\ldots, 
\begin{bmatrix}
	0 & 1 \\
	-1 & 0
\end{bmatrix}, 
0, \ldots, 0\right).
$$
\paragraph{Note.} Under certain mild conditions, the block-diagonalization can be reformulated into a spectral decomposition for skew-symmetric matrices making $\bZ$   orthogonal. See Corollary~\ref{corollary:spectral_normal_skew}.
\end{theoremHigh}
\begin{proof}[of Theorem~\ref{theorem:skew-block-diagonalization_theorem}]
We will prove this through a  recursive calculation. 
As customary, we will denote the entry ($i,j$) of a matrix $\bA$ as $a_{ij}$. 

\paragraph{Case 1).} Suppose the first row of $\bA$ is nonzero. Note that $\bE\bA\bE^\top$ is skew-symmetric for any matrix $\bE$ if $\bA$ is skew-symmetric. This ensures that  both the diagonals of $\bA$ and $\bE\bA\bE^\top$ are zero. The upper-left $2\times 2$ submatrix of $\bE\bA\bE^\top$ takes the following form
$$
(\bE\bA\bE^\top)[1:2,1:2]=
\begin{bmatrix}
0 & x \\
-x & 0
\end{bmatrix}.
$$
Since we suppose the first row of $\bA$ is nonzero, there exists a permutation matrix $\bP$ (Definition~\ref{definition:permutation-matrix}) such that we will exchange any nonzero value, say $a$, in the first row to the second column of $\bP\bA\bP^\top$. And as discussed above, the upper-left $2\times 2$ submatrix of $\bP\bA\bP^\top$ has the following form
$$
(\bP\bA\bP^\top)[1:2,1:2]=
\begin{bmatrix}
	0 & a \\
	-a & 0
\end{bmatrix}.
$$
Construct a nonsingular matrix $\bM = 
\footnotesize\begin{bmatrix}
1/a & \bzero \\
\bzero   & \bI_{n-1}
\end{bmatrix}$ such that the upper-left $2\times 2$ submatrix of $\bM\bP\bA\bP^\top\bM^\top$ becomes
$$
(\bM\bP\bA\bP^\top\bM^\top)[1:2,1:2]=
\begin{bmatrix}
0 & 1 \\
-1 & 0
\end{bmatrix}.
$$
Now we finish block-diagonalizing the upper-left $2\times 2$ block. Suppose now $(\bM\bP\bA\bP^\top\bM^\top)$  has a nonzero value, say $b$, in the first row at entry $(1,j)$ for some $j>2$. We can construct a nonsingular matrix $\bL = \bI - b\cdot\bE_{j2}$, where $\bE_{j2}$ is an all-zero matrix except the entry ($j, 2$) is 1. Then, $\bL(\bM\bP\bA\bP^\top\bM^\top)\bL^\top$  introduces a zero for the entry with value $b$.

\begin{mdframed}[hidealllines=\mdframehideline,backgroundcolor=\mdframecolor,frametitle={A Trivial Example}]
For example, suppose $\bM\bP\bA\bP^\top\bM^\top$ is a $3\times 3$ matrix with the following value 
$$
\bM\bP\bA\bP^\top\bM^\top = 
\begin{bmatrix}
0 & 1 & b \\
-1 & 0 & \times \\
\times & \times & 0
\end{bmatrix}, 
\qquad \text{and}\qquad
\bL =\bI - b\cdot\bE_{j2}=
\begin{bmatrix}
	1 & 0 & 0 \\
	0 & 1 & 0 \\
	0 &-b & 1
\end{bmatrix},
$$
where $j=3$ in this particular example. This results in 
$$
\bL\bM\bP\bA\bP^\top\bM^\top\bL^\top =
\begin{bmatrix}
	1 & 0 & 0 \\
	0 & 1 & 0 \\
	0 &-b & 1
\end{bmatrix}
\begin{bmatrix}
	0 & 1 & \textcolor{mylightbluetext}{b} \\
	-1 & 0 & \times \\
	\times & \times & 0
\end{bmatrix}
\begin{bmatrix}
1 & 0 & 0 \\
0 & 1 & -b \\
0 & 0 & 1
\end{bmatrix} 
= 
\begin{bmatrix}
	0 & 1 & \textcolor{mylightbluetext}{0} \\
	-1 & 0 & \times \\
	\times & \times & 0
\end{bmatrix}.
$$
\end{mdframed}
Similarly, if the second row of $\bL\bM\bP\bA\bP^\top\bM^\top\bL^\top$ contains a nonzero value, say $c$, we can construct a nonsingular matrix $\bK = \bI+c\cdot \bE_{j1}$ such that $\bK(\bL\bM\bP\bA\bP^\top\bM^\top\bL^\top)\bK^\top$  introduces a zero for the entry with value $c$.
\begin{mdframed}[hidealllines=\mdframehideline,backgroundcolor=\mdframecolor,frametitle={A Trivial Example}]
	For example, suppose $\bL\bM\bP\bA\bP^\top\bM^\top\bL^\top$ is a $3\times 3$ matrix with the following value 
	$$
	\bL\bM\bP\bA\bP^\top\bM^\top\bL^\top = 
	\begin{bmatrix}
		0 & 1 & 0 \\
		-1 & 0 & c \\
		\times & \times & 0
	\end{bmatrix}, 
	\qquad \text{and}\qquad
	\bK =\bI + c\cdot\bE_{j1}=
	\begin{bmatrix}
		1 & 0 & 0 \\
		0 & 1 & 0 \\
		c &0 & 1
	\end{bmatrix},
	$$
	where $j=3$ for this specific example. This results in 
	$$
	\bK\bL\bM\bP\bA\bP^\top\bM^\top\bL^\top\bK^\top =
	\begin{bmatrix}
		1 & 0 & 0 \\
		0 & 1 & 0 \\
		c & 0 & 1
	\end{bmatrix}
	\begin{bmatrix}
		0 & 1 & 0 \\
		-1 & 0 & \textcolor{mylightbluetext}{c} \\
		\times & \times & 0
	\end{bmatrix}
	\begin{bmatrix}
		1 & 0 & c \\
		0 & 1 & 0 \\
		0 & 0 & 1
	\end{bmatrix} 
	= 
	\begin{bmatrix}
		0 & 1 & 0 \\
		-1 & 0 & \textcolor{mylightbluetext}{0} \\
		\times & \times & 0
	\end{bmatrix}.
	$$
Since we have shown that $\bK\bL\bM\bP\bA\bP^\top\bM^\top\bL^\top\bK^\top$ is also skew-symmetric, then it simplifies to 
$$
\bK\bL\bM\bP\bA\bP^\top\bM^\top\bL^\top\bK^\top= 
\begin{bmatrix}
	0 & 1 & 0 \\
	-1 & 0 & \textcolor{mylightbluetext}{0} \\
	\textcolor{winestain}{0} & \textcolor{winestain}{0} & 0
\end{bmatrix},
$$
so we do not need to address the first two columns further.
\end{mdframed}
By applying the same process to the bottom-right $(n-2)\times(n-2)$ submatrix, we complete the proof for this case.

\paragraph{Case 2).} If the first row of $\bA$ is zero, we can use a permutation matrix to move the first row to the last row and then apply the process described in case 1 to complete the proof.
\end{proof}

The block-diagonalization of skew-symmetric matrices, as discussed above, reveals that the rank of a skew-symmetric matrix is always even.   
Moreover, we can establish that the determinant of skew-symmetric matrices of even order is nonnegative as follows.
\begin{corollary}[Even Skew-Symmetric Determinant]\label{corollary:skew-symmetric-determinant-even}
Let $\bA\in \real^{n\times n}$ be a skew-symmetric matrix. If $n$ is even, then $\det(\bA)\geq 0$.
\end{corollary}
\begin{proof}[of Corollary~\ref{corollary:skew-symmetric-determinant-even}]
According to Theorem~\ref{theorem:skew-block-diagonalization_theorem}, we can block-diagonalize $\bA = \bZ\bD\bZ^\top$ such that 
$
\det(\bA) = \det(\bZ\bD\bZ^\top) = \det(\bZ)^2 \det(\bD) \geq 0.
$
This completes the proof.
\end{proof}

\section{Spectral Decomposition of Normal Matrices*}
More generally, we consider the spectral decomposition of a normal matrix $\bA$: $\bA^\top\bA=\bA\bA^\top$ (Definition~\ref{definition:complex_special}).
\begin{theoremHigh}[Spectral Decomposition of Real Normal Matrices]\label{theorem:normal_real_spectral_theorem}
Let $\bA \in \real^{n\times n}$ be a \textbf{real} normal matrix. Then, it can be factored as
\begin{equation*}
	\bA = \bQ \bD \bQ^\top,
\end{equation*}
where $\bQ\in \real^{n\times n}$ is an $n\times n$ (real) orthogonal matrix (with possible complex eigenvalues), and $\bD\in \real^{n\times n}$ is a block-diagonal (\textit{quasi-diagonal}) matrix with the following form
$$
\bD = 
\begin{bmatrix}
\bA_1 &&& \\
& \bA_2 && \\
&&\ddots & \\
&&&\bA_p
\end{bmatrix},
$$
where each  $\bA_i$ is either a real scalar (corresponding to a real eigenvalue of $\bA$) or a $2\times 2$ real matrix of the following form 
\begin{equation}\label{equation:spec_normal_equ}
\bA_i=
\begin{bmatrix}
	\lambda_i & \mu_i \\
	-\mu_i & \lambda_i
\end{bmatrix},
\,
\text{corresponding to a complex eigenvalue $\lambda_i\pm i\mu_i $ of $\bA$ ($\mu_i>0$)}.
\end{equation}
\noindent
\paragraph{Note.} This theorem generalizes the main theorem (Theorem~\ref{theorem:spectral_theorem}) in terms of the matrix type: the generalization extends from symmetric matrices to normal matrices.
The  blocks $\{\bA_i\}_{i=1}^p$ may appear in any prescribed order.
Remark~\ref{remark:conjug_pari} shows that  non-real eigenvalues appear in complex conjugate pairs; and this theorem shows that the complex conjugate eigenvalue pairs correspond to the $2\times 2$ blocks in the block-diagonal matrix $\bD$.
\end{theoremHigh}


To see this, we first examine the properties of normal matrices. We provide general properties of \textbf{complex} normal matrices, which can be easily adapted to the context of real normal matrices.

\subsection*{Properties of Normal Matrices}
\begin{lemma}[Null Space of Complex Normal Matrices]\label{lemma:null_normal_space}
Let $\bA\in \complex^{n\times n}$ be  a complex  normal matrix. Then, $\bA$ and its Hermitian transpose, $\bA^\ast$, share the same null space: $\nspace(\bA)= \nspace(\bA^\ast)$.
\end{lemma}
\begin{proof}[of Lemma~\ref{lemma:null_normal_space}]
If $\bA$ is normal, we have 
$$
\begin{aligned}
\bA^\ast \bx = \bzero 
\implies 
\langle \bA^\ast \bx, \bA^\ast \bx\rangle = 0 
\implies
\langle \bA \bx, \bA \bx\rangle = 0
\implies
\bA \bx = \bzero.
\end{aligned}
$$
\footnote{The dot product notation, $\langle \bx,\bx\rangle = \bx^\top\bx$ if $\bx\in \real^n$ or $=\bx^\ast\bx$ if $\bx\in \complex^n$.}
Therefore, $\nspace(\bA^\ast)\subseteq \nspace(\bA)$.
Similarly, 
$$
\bA \bx=\bzero 
\implies 
\bA^\ast \bx=\bzero
$$
This implies $\nspace(\bA)\subseteq \nspace(\bA^\ast)$. 
By ``sandwiching," we have 
$\nspace(\bA)= \nspace(\bA^\ast)$.
\end{proof}

\begin{lemma}[Conjugate Eigenpair of Normal Matrices]\label{lemma:conjugate_eigen_pair}
Let $\bA\in \complex^{n\times n}$ be  a complex  normal matrix. Then, the eigenpair $\bA\bx=\lambda \bx$ implies 
$\bA^\ast\bx=\overline{\lambda}\bx$.
\end{lemma}
\begin{proof}[of Lemma~\ref{lemma:conjugate_eigen_pair}]
We note that $\bx$ is a vector in the null space of $\bB \triangleq \bA-\lambda\bI$, and $\bB$ is also normal. 
According to Lemma~\ref{lemma:null_normal_space}, $\bx$ is also in the null space of $\bB^\ast=\bA^\ast - \overline{\lambda}\bI$ and this implies $\bA^\ast\bx=\overline{\lambda}\bx$.
\end{proof}

\begin{lemma}[Distinct Eigenvectors of Normal Matrices]\label{lemma:distinct_normal_vec}
Let $\bA\in \complex^{n\times n}$ be a complex  normal matrix, and let $\bx, \by$ be eigenvectors of $\bA$ corresponding to distinct eigenvalues. 
Then, $\langle \bx, \by\rangle =0$, i.e., $\bx,\by$ are orthogonal.
The also indicates the eigenspaces (i.e., the null spaces $\nspace(\bA - \lambda_i\bI)$ corresponding to eigenvalue $\lambda_i$) of a normal matrix $\bA\in\complex^{n\times n}$ are mutually orthogonal.
\footnote{If $\bx$ is complex, then the real parts and imaginary parts are orthogonal respectively.}
\footnote{Note we have shown in Proposition~\ref{proposition:orthogonal-eigenvectors} that any eigenvectors corresponding to distinct eigenvalues of symmetric matrices are orthogonal.}
\end{lemma}
\begin{proof}[of Lemma~\ref{lemma:distinct_normal_vec}]
Suppose $\bA\bx=\lambda\bx$ and $\bA\by=\mu\by$ with $\lambda\neq \mu$. 
By Lemma~\ref{lemma:conjugate_eigen_pair}, we have 
$$
\left\{
\begin{aligned}
\bA\bx &= \lambda\bx \quad\implies\quad \bA^\ast\bx = \overline{\lambda}\bx;  \\
\bA\by &= \mu\by \quad\implies\quad \bA^\ast\by = \overline{\mu}\by.
\end{aligned}
\right.
$$
Therefore, we get 
$$
\left\{
\begin{aligned}
\langle \bx, \bA\by\rangle 
&= \langle \bx, \mu\by\rangle = \overline{\mu} \langle\bx,\by\rangle;\\ 
\langle \bx, \bA\by\rangle 
&=\langle \bA^\ast \bx, \by\rangle = 
\langle \overline{\lambda}\bx, \by\rangle=
\overline{\lambda} \langle \bx, \by\rangle.
\end{aligned}
\right.
$$
Since we assume $\lambda\neq \mu$, we get $\langle \bx, \by\rangle =0$.
\end{proof}
The above lemma also shows that we can pick an orthonormal basis for each eigenspace of a normal matrix. The union of these bases forms an orthonormal basis for $\complex^n$.

\subsection*{Main Proof}
We now provide the proof of the spectral theorem for normal matrices.
\begin{proof}[{of Theorem~\ref{theorem:normal_real_spectral_theorem}: Spectral Decomposition of Real Normal Matrices}]
Suppose the real normal matrix $\bA$ has an eigenpair $(z, \bw)$, where $z=\lambda+i\mu$ and $\bw=\bx+i\by$ ($\lambda, \mu, \bx, \by$ are real scalars or real vectors).
Then,
$$
\begin{aligned}
\bA\bw = z\bw &\quad\implies\quad 
\bA(\bx+i\by) = (\lambda+i\mu)(\bx+i\by)\\
&\quad\implies\quad \left\{
\begin{aligned}
\bA\bx &= \lambda\bx - \mu\by; \\
\bA\by&= \lambda\by +\mu\bx.
\end{aligned}
\right.
\end{aligned}
$$
We consider the proof in two cases: $\mu=0$ and $\mu\neq 0$.

\paragraph{Case 1).} Suppose $\mu\neq 0$.
According to Lemma~\ref{lemma:conjugate_eig_pair}, we also have $\bA(\bx-i\by) = (\lambda-i\mu)(\bx-i\by)$. 
Therefore, $\bw$ and $\overline{\bw}$ are two eigenvectors of $\bA$ corresponding to different eigenvalues. By Lemma~\ref{lemma:distinct_normal_vec}, $\langle \bw, \overline{\bw}\rangle=0$:
$$
\langle \bx+i\by, \bx-i\by\rangle = 
(\bx^\top-i\by^\top)(\bx-i\by) = \bx^\top\bx -\by^\top\by - 2i \bx^\top\by = 0.
$$
That is, $\bx^\top\by=0$ and $\bx^\top\bx = \by^\top\by$. Without loss of generality, we assume $\normtwo{\bx} = \normtwo{\by} = 1$ such that $\bx$ and $\by$ are orthonormal.

Let $\bP$ be an orthogonal matrix with $\bx$ as the first column and $\by$ as the second column:
$
\bP \triangleq
[\bx, \by, \bP_{n-2}],
$
where $\bP_{n-2}\in \real^{n\times (n-2)}$ can be constructed using the Gram-Schmidt process such that $\bP_{n-2}^\top\bP_{n-2}=\bI_{n-2}$. 
Let further 
$$
\bA_1 \triangleq 
\begin{bmatrix}
	\lambda_1 & \mu_1 \\
	-\mu_1 & \lambda_1
\end{bmatrix}\triangleq
\begin{bmatrix}
	\lambda & \mu \\
	-\mu & \lambda
\end{bmatrix}
\gap \text{such that}\gap
\bP^\top \bA\bP=
\begin{bmatrix}
	\lambda_1 & \mu_1 & \ldots & 0 \\
-\mu_1 & \lambda_1 && \\
	\vdots &&\bB&\\
	0 &&&
\end{bmatrix},
$$
where $\bB \in \real^{(n-2)\times (n-2)}$.
Since $\bA^\top\bA=\bA\bA^\top$ and $\bA_1^\top\bA_1=\bA_1\bA_1^\top$, it can be shown that $\bB$ is also normal.
Therefore, $\bB$ can be factored as 
$$
\widetildebR^\top \bB\widetildebR=
\begin{bmatrix}
\lambda_2 & \mu_2 & \ldots & 0 \\
-\mu_2 & \lambda_2 && \\
\vdots &&\bC&\\
0 &&&
\end{bmatrix}
\quad 
\text{with}\quad
\bA_2 = 
\begin{bmatrix}
\lambda_2 & \mu_2 \\
	-\mu_2 & \lambda_2
\end{bmatrix},
$$
where $\bC \in \real^{(n-4)\times (n-4)}$, and $\lambda_2+i\mu_2$ is an eigenvalue of both $\bA$ and $\bB$.
Let 
$$
\bR\triangleq
\begin{bmatrix}
\bI_2 & \bzero \\
\bzero & \widetildebR 
\end{bmatrix}
\gap \text{such that}\gap 
\bP^\top 
\bR^\top 
\begin{bmatrix}
\bA_1 & &\\
& \bA_2 & \\
&& \bC
\end{bmatrix}
\bR\bP,
$$
where $\bR\bP$ is also orthogonal and $\bC$ is also normal.
By continuing this recursive process, we conclude the proof.

\paragraph{Case 2).} Suppose $\mu=0$. Then $\bx$ and $\by$ are eigenvectors of $\bA$ corresponding to the eigenvalue $\lambda$. The proof then proceeds similarly, with $\bA_i$ being a scalar.
\end{proof}

A direct consequence of the spectral decomposition for real normal matrices is discussed below.
The proof is left as an exercise.
\begin{corollary}[Orthogonally Similar Real Normal Matrices]\label{corollary:ort_sim_real}
Two real orthogonal matrix are (real) orthogonally similar to each other if and only if they have the same (complex or real) eigenvalues.
\end{corollary}


\subsection*{Theorem~\ref{theorem:normal_real_spectral_theorem} for Symmetric Matrices}
We note that symmetric matrices belong to the category of special normal matrices. 
Subsequently, we will elucidate the  connection  between the spectral theorem for normal matrices and  that for symmetric matrices.
\begin{corollaryHigh}[Spectral Decomposition of Symmetric Matrices]\label{corollary:spectral_normal_symmetricx}
When dealing with a real symmetric $\bA\in\real^{n\times n}$ in Theorem~\ref{theorem:normal_real_spectral_theorem}, the $\mu_i$'s in Equation~\eqref{equation:spec_normal_equ} are equal to zero. 
This implies that a real symmetric matrix $\bA$ is diagonalizable, confirming  the main result  stated in Theorem~\ref{theorem:spectral_theorem}.
\end{corollaryHigh}
The proof is evident since real symmetric matrices have only real eigenvalues. Alternatively, we provide another proof below.
\begin{proof}[of Corollary~\ref{corollary:spectral_normal_symmetricx}]
Suppose the real symmetric matrix $\bA$ has an eigenpair $(z, \bw)$, where $z=\lambda+i\mu$ and $\bw=\bx+i\by$ ($\lambda, \mu, \bx, \by$ are real scalars or real vectors).
Then,
$$
\begin{aligned}
\bA\bw = z\bw 
\implies \left\{
\begin{aligned}
\bA\bx &= \lambda\bx - \mu\by; \\
\bA\by&= \lambda\by +\mu\bx.
\end{aligned}
\right.
\end{aligned}
\implies
\left\{
\begin{aligned}
	\langle \bA\bx, \by \rangle 
	&= \langle \lambda\bx - \mu\by, \by \rangle 
	= \lambda\langle \bx , \by \rangle - \mu\langle\by, \by \rangle;\\
	\langle \bx, \bA\by \rangle 
	&= \langle \bx, \lambda\by +\mu\bx \rangle 
	= \lambda\langle \bx , \by \rangle+ \mu\langle\bx, \bx \rangle.
\end{aligned}
\right.
$$
Therefore, $\mu \left(\langle\bx, \bx \rangle + \langle\by, \by \rangle\right)=0$. Since either $\bx\neq \bzero $ or $\by\neq \bzero$, we have $\mu=0$.
\end{proof}

\subsection*{Theorem~\ref{theorem:normal_real_spectral_theorem} for Skew-Symmetric Matrices}
We further note that skew-symmetric matrices are also a subset of special normal matrices. 
We will now show the connection between the spectral theorem for normal matrices and that for skew-symmetric matrices.
\begin{corollaryHigh}[Spectral Decomposition of Skew-Symmetric Matrices]\label{corollary:spectral_normal_skew}
When dealing with a real skew-symmetric $\bA\in\real^{n\times n}$ in Theorem~\ref{theorem:normal_real_spectral_theorem}, the $\lambda_i$'s in Equation~\eqref{equation:spec_normal_equ} are zeros.
\end{corollaryHigh}
Again, the proof is evident since real skew-symmetric matrices have only complex eigenvalues, which are either purely imaginary or zero. Alternatively, we provide another proof below.
\begin{proof}[of Corollary~\ref{corollary:spectral_normal_skew}]
Suppose the real skew-symmetric matrix $\bA$ has an eigenpair $(z, \bw)$, where $z=\lambda+i\mu$ and $\bw=\bx+i\by$ ($\lambda, \mu, \bx, \by$ are real scalars or real vectors).
Then,
$$
\begin{aligned}
\bA\bw= z\bw  
&\quad\implies\quad \left\{
\begin{aligned}
\bA\bx &= \lambda\bx - \mu\by; \\
\bA\by&= \lambda\by +\mu\bx.
\end{aligned}
\right.
\end{aligned}
$$
Since $\bA = -\bA^\top$ and the inner product $\langle \cdot, \cdot \rangle$ is symmetric, it follows that 
$$
\langle \bA\bw, \bw \rangle 
=
\langle \bw, \bA^\top \bw \rangle
=
\langle \bw, -\bA\bw \rangle 
= 
- \langle \bw, \bA\bw \rangle
=
-\langle \bA\bw, \bw \rangle.
$$
This implies $\langle \bA\bw, \bw \rangle =0$.
We then have 
$$
\begin{aligned}
\langle \bA\bx, \bx \rangle 
&= \langle \lambda\bx - \mu\by, \bx \rangle 
= \lambda\langle \bx , \bx \rangle - \mu\langle\by, \bx \rangle =0,\\
\langle \bA\by, \by \rangle
&= \langle \lambda\by +\mu\bx, \by \rangle 
= \lambda\langle \by , \by \rangle + \mu\langle\by, \bx \rangle =0.\\
\end{aligned}
$$
Therefore, $\lambda \left(\langle\bx, \bx \rangle + \langle\by, \by \rangle\right)=0$. Since either $\bx\neq \bzero $ or $\by\neq \bzero$, we have $\lambda=0$.
\end{proof}

\begin{exercise}
Use the  above result to prove the block-diagonalization for skew-symmetric matrices  in Theorem~\ref{theorem:skew-block-diagonalization_theorem}.
\end{exercise}
\begin{exercise}
Prove the converse of Corollary~\ref{corollary:spectral_normal_skew} is also true. Specifically, show that any matrix of the form described in Corollary~\ref{corollary:spectral_normal_skew} is skew-symmetric. (Note that any matrix of the form in Corollary~\ref{corollary:spectral_normal_symmetricx} is symmetric.)
\end{exercise}

\subsection*{Theorem~\ref{theorem:normal_real_spectral_theorem} for Orthogonal Matrices}
We also observe that orthogonal matrices belong to the class of special normal matrices. 
We then show the connection of spectral theorem between normal matrices and orthogonal matrices in the following corollary.
\begin{corollaryHigh}[Spectral Decomposition of Orthogonal Matrices]\label{corollary:spectral_normal_Orthogonal}
When dealing with a real orthogonal $\bA\in\real^{n\times n}$ in Theorem~\ref{theorem:normal_real_spectral_theorem}, then each block in Equation~\eqref{equation:spec_normal_equ} is either 1, $-1$, or a two-dimensional matrix of the following form 
$$
\bA_i
=
\begin{bmatrix}
\cos \theta_i & \sin\theta_i \\
-\sin \theta_i & \cos\theta_i 
\end{bmatrix},
\quad 
\text{where$\quad 0<\theta_i<\pi$.}
$$
\end{corollaryHigh}
\begin{proof}[of Corollary~\ref{corollary:spectral_normal_Orthogonal}]
The result follows immediate from the fact that the eigenvalues of orthogonal matrices have an absolute value of 1.
\end{proof}

For the spectral decomposition of an orthogonal matrix, we can also reorder the orthonormal basis  such that the orthogonal matrix $\bA$ can be factored as
$
\bA = \bQ \bD \bQ^\top,
$
where $\bQ$ is an $n\times n$ orthogonal matrix, and $\bD$ is a block-diagonal matrix with the following form
$$
\bD = 
\footnotesize
\begin{bmatrix}
\bA_1 &\ldots&&  && \\
\vdots & \bA_2 && &&\\
&&\ddots &\vdots && \\
&\ldots&\ldots&\bA_k&&\\
&&&& -\bI_p& \\
&&&&& \bI_q \\
\end{bmatrix}
\normalsize
\gap 
\text{with\,\, }
\bA_i
=
\begin{bmatrix}
	\cos \theta_i & \sin\theta_i \\
	-\sin \theta_i & \cos\theta_i 
\end{bmatrix}.
$$

In conclusion, when applying the spectral theorem to different types of real normal matrices, we obtain various forms as shown in Figure~\ref{fig:spec_normal-summary}.

\begin{figure}[htbp]
\centering
\centering
\resizebox{0.9\textwidth}{!}{%
\begin{tikzpicture}[>=latex]

\tikzstyle{state} = [draw, very thick, fill=white, rectangle, minimum height=3em, minimum width=6em, node distance=8em, font={\sffamily\bfseries}]
\tikzstyle{stateEdgePortion} = [black,thick];
\tikzstyle{stateEdge} = [stateEdgePortion,->];
\tikzstyle{stateEdge2} = [stateEdgePortion,<->];
\tikzstyle{edgeLabel} = [pos=0.5, text centered, font={\sffamily\small}];

\node[ellipse, name=pdmatrix, draw,font={\sffamily\bfseries},  node distance=7em, xshift=-9em, yshift=-1em,fill={colorals}]  {Normal $\bA$};
\node[state, name=bsqure, below of=pdmatrix, xshift=0em, yshift=1em, fill={colorlu}] {$\bQ\bD\bQ^\top$};
\node[state, name=ptp, right of=bsqure, xshift=3em, fill={colorlu}] {$\bQ\bD\bQ^\top$};
\node[state, name=rsqure, left of=bsqure, xshift=-3em, fill={colorlu}] {$\bQ\bD\bQ^\top$};
\node[ellipse, name=utv, below of=pdmatrix,draw,  node distance=7em, xshift=0em, yshift=-4em,font={\tiny},fill={coloruppermiddle}]  {\parbox{9em}{Diagonal of Order 2,\\i.e., $\lambda_i=0$}};
\node[ellipse, name=upperr, left of=utv, draw, node distance=8em, xshift=-3em,font={\tiny},fill={coloruppermiddle}]  {\parbox{9em}{Diagonal of Order 1,\\i.e., $\mu_i=0$}};
\node[ellipse, name=nonp, right of=utv,draw,  node distance=8em, xshift=3em, font={\tiny},fill={coloruppermiddle}]  {\parbox{9.9em}{Diagonal of Order 2, \\ Columns have length 1}};

\coordinate (lq2inter3) at ($(pdmatrix.east -| ptp.north) + (-0em,0em)$);
\draw (pdmatrix.east) edge[stateEdgePortion] (lq2inter3);
\draw (lq2inter3) edge[stateEdge] 
node[edgeLabel, text width=7.25em, yshift=0.8em]{{Orthogonal}} (ptp.north);

\coordinate (rqr2inter1) at ($(pdmatrix.west) + (0,0em)$);
\coordinate (rqr2inter3) at ($(rqr2inter1-| rsqure.north) + (-0em,0em)$);
\draw (rqr2inter1) edge[stateEdgePortion] (rqr2inter3);
\draw (rqr2inter3) edge[stateEdge] 
node[edgeLabel, text width=8em, yshift=0.8em]{{Symmetric}} (rsqure.north);

\draw (pdmatrix.south)
edge[stateEdge] node[edgeLabel,yshift=0.5em]{{Skew-Symmetric} } 
(bsqure.north);

\draw (upperr.north)
edge[stateEdge] node[edgeLabel,yshift=0.5em]{} 
(rsqure.south);

\draw (utv.north)
edge[stateEdge] node[edgeLabel,yshift=0.5em]{} 
(bsqure.south);

\draw (nonp.north)
edge[stateEdge] node[edgeLabel,yshift=0.5em]{} 
(ptp.south);

\begin{pgfonlayer}{background}
\draw [join=round,cyan,dotted,fill={colormiddle}] ($(upperr.south west) + (-1.6em, -1em)$) rectangle ($( pdmatrix.east-|ptp.north east) + (2.3em, +1.8em)$);
\end{pgfonlayer}

\end{tikzpicture}
}
\caption{Demonstration of the spectral theorem on a normal matrix $\bA$.}
\label{fig:spec_normal-summary}
\end{figure}

\subsection*{Spectral Decomposition for Normal by Parts}
A real normal matrix $\bA$ can be decomposed into $\bA = (\bA+\bA^\top)+(\bA-\bA^\top)$, where $\bA+\bA^\top$ is symmetric and $\bA-\bA^\top$ is skew-symmetric. Therefore, the decomposition of a general real normal matrix can be done in two parts.
It  can then be shown that the real parts of $\bA$'s eigenvalues are the eigenvalues of the symmetric matrix $(\bA+\bA^\top)$, and the imaginary parts of its eigenvalues are the eigenvalues of the skew-symmetric matrix $(\bA-\bA^\top)$.

\subsection{Spectral Decomposition for Complex Normal Matrices}

In Theorem~\ref{theorem:normal_real_spectral_theorem}, we focus on the case of \textbf{real} normal matrices, where the factored components contain only real elements. 
However, when we consider \textbf{complex} normal matrices, we have  the following result. Instead of a block-diagonalization, we achieve a full diagonalization (which may include complex elements).

\begin{theoremHigh}[Spectral Decomposition of Complex Normal Matrices]\label{theorem:normal_Complex_spectral_theorem}
Let $\bA \in \complex^{n\times n}$ be a \textbf{complex} normal matrix. Then, it can be factored as
\begin{equation*}
\bA = \bU \bD \bU^\ast,
\end{equation*}
where $\bU$ is an $n\times n$ unitary matrix (Definition~\ref{definition:complex_special}) containing unit eigenvectors of $\bA$, and $\bD$ is a diagonal matrix containing the corresponding eigenvalues (\textbf{not necessarily real}).
\footnote{This theorem is a generalization of the main theorem (Theorem~\ref{theorem:spectral_theorem}) in terms of the matrix domain type: from real symmetric matrix to complex normal matrix; and it is also a special case of the Schur decomposition for range-Hermitian matrices (Problem~\ref{prob:schur_range_hermi}).}
\footnote{\textbf{Important Takeaway:} Note that this claim also holds for \textcolor{mylightbluetext}{real} normal matrix $\bA\in\real^{n\times n}$. For real symmetric matrices or complex Hermitian matrices, the eigenvalues are real. However, even for real normal matrices, the eigenvalues are not necessarily real. Therefore, for real normal matrices, the diagonal matrix $\bD$ may contain complex values on the diagonal:
\begin{enumerate}
\item \underline{Real symmetric $\bA$}: $\bA = \bQ\bLambda\bQ^\top$, where $\bQ$ is real orthogonal and $\bLambda$ is \textbf{real} diagonal;
\item \underline{Complex Hermitian $\bA$}: $\bA=\bU\bLambda\bU^\ast$, where $\bU$ is complex unitary and $\bLambda$ is \textbf{real} diagonal (Prob~\ref{problem:real_herm});
\item \underline{Real or complex normal $\bA$}: $\bA = \bU\bD\bU^\ast$, where $\bU$ is complex unitary and $\bD$ is \textbf{complex} diagonal.
\end{enumerate}
}
\end{theoremHigh}

To illustrate this, we require an exercise as follows.
\begin{exercise}
Show that a complex matrix $\bA\in \complex^{n\times n}$ is normal and upper triangular (or lower triangular) if and only if $\bA$ is a diagonal matrix.

The result can be generalized to a block normal matrix. Suppose $\bA=\footnotesize\begin{bmatrix}
\bA_{11}&\bA_{12}\\
\bA_{21}&\bA_{22}
\end{bmatrix}$. Show that $\bA$ is normal if and only if $\bA_{11}$ and $\bA_{22}$ are normal and $\bA_{12}=\bzero$. 
This can be generalized further that 
a block upper triangular matrix is normal if and only if each of its off-diagonal blocks is zero and each of its diagonal blocks is normal.
\end{exercise}

\begin{proof}[of Theorem~\ref{theorem:normal_Complex_spectral_theorem}]
For simplicity, we employ the result of complex Schur decomposition (Theorem~\ref{theorem:schur-decomposition_complex}) to prove this theorem.
For a complex normal matrix $\bA$, it admits a complex Schur decomposition
$$
\begin{aligned}
\bA &= \bU\bD\bU^\ast
\gap&\text{and}&\gap 
&\bA^\ast &= \bU\bD^\ast\bU^\ast\\
\implies\gap 
\bA\bA^\ast &= \bU\bD\bD^\ast\bU^\ast
\gap&\text{and}&\gap 
&\bA^\ast\bA &=\bU\bD^\ast \bD\bU^\ast.
\end{aligned}
$$
Since $\bA$ is normal, we get $\bD\bD^\ast = \bD^\ast \bD$. Because $\bD$ is normal and upper triangular,  it must be diagonal, as demonstrated in the preceding exercise. We complete the proof.
\end{proof}

\begin{exercise}
There is a converse statement to Theorem~\ref{theorem:normal_Complex_spectral_theorem}: if a matrix can be diagonalized (with possible complex elements), then the matrix is normal.
\end{exercise}

\begin{exercise}
Consider a complex square matrix $\bA\in\complex^{n\times n}$, show that $\bu\in\complex^n$ is a right eigenvector of $\bA$ corresponding to eigenvalue $\lambda$ if and only if $\bu^*$ is a left eigenvector of $\bA$ corresponding to the same eigenvalue $\lambda$. \textit{Hint: Consider the spectral decomposition.}
\end{exercise}

To summarize this section, we examine the case of a \textbf{complex}  matrix $\bA$ for Corollary~\ref{corollary:spectral_normal_symmetricx}, Corollary~\ref{corollary:spectral_normal_skew}, and Corollary~\ref{corollary:spectral_normal_Orthogonal}.
\begin{theoremHigh}[Spectral Decomposition of Complex Normal Matrices]
Given a complex normal matrix $\bA\in \complex^{n\times n}$, there exists a unitary matrix $\bU$ and a diagonal matrix $\bD$ such that 
$$
\bA=\bU\bD\bU^\ast.
$$ 
\begin{itemize}
\item If $\bA$ is Hermitian, then $\bD$ is a real matrix. 
\item If $\bA$ is skew-Hermitian, then the entries in $\bD$ are purely imaginary or zero. 
\item If $\bA$ is unitary, then the entries in $\bD$ have absolute value 1.
\end{itemize}
\end{theoremHigh}

\subsection{Properties of (Complex) Normal Matrices}

Theorem~\ref{theorem:normal_Complex_spectral_theorem} demonstrates that any complex normal matrix can be \textit{unitarily diagonalizable}.
For a complex normal matrix, the equivalent conditions can be summarized in the following theorem.
\begin{theorem}[Equivalence of Complex Normal Matrices]\label{theorem:equiv_cmnormal}
Suppose $\bA\in\complex^{n\times n}$ has eigenvalues $\lambda_1,\lambda_2,\ldots,\lambda_n\in\complex$. The following four statements are equivalent:
\begin{enumerate}
\item $\bA$ is normal;
\item $\bA$ is unitarily diagonalizable;
\item $\sum_{i,j=1}^{n}\abs{a_{ij}}^2 = \sum_{i=1}^{n}\abs{\lambda_i}^2$;
\item $\bA$ has $n$ orthonormal eigenvectors.
\end{enumerate}
\end{theorem}
\begin{proof}[of Theorem~\ref{theorem:equiv_cmnormal}]
The indication from (1) to (2) is already proved in Theorem~\ref{theorem:normal_Complex_spectral_theorem}.
Suppose (2) holds $\bA=\bQ\bD\bQ^{*}$. Then $\sum_{i,j=1}^{n}\abs{a_{ij}}^2=\trace(\bA^*\bA)=\trace(\bD^*\bD)$, which proves (3).
The remaining implications are straightforward and are left as exercises.
\end{proof}

We have shown in Corollary~\ref{corollary:ort_sim_real} that two real normal matrices are orthogonally similar to each other if and only if they have the same (complex or real) eigenvalues.
A similar result can be observed in the complex case.
\begin{theorem}[Unitarily Similar Complex Normal Matrices]\label{theorem:uni_sim_nmat}
Two (complex)  normal matrices are unitarily similar if and only if they have the same eigenvalues.
\end{theorem}
\begin{proof}[of Theorem~\ref{theorem:uni_sim_nmat}]
Consider normal matrices $\bA,\bB\in\complex^{n\times n}$ and suppose $\bA$ admits the spectral decomposition $\bA=\bU\bD\bU^*$.
If $\bA$ and $\bB$ are unitarily similar, then $\bB=\bQ\bA\bQ^*= \bQ\bU\bD\bU^*\bQ^*$. Therefore, $\bA$ and $\bB$ have the same eigenvalues.

Conversely, assume $\bA$ and $\bB$ have the same eigenvalues. Since $\bA$ admits the spectral decomposition $\bA=\bU\bD\bU^*$, where the diagonals of $\bD$ are the eigenvalues of $\bA$ and $\bB$. $\bB$ also admits the spectral decomposition $\bB=\bV\bD\bV^*$. Therefore, 
$\bA=\bU\bV^*(\bV\bD\bV^*)\bV\bU^*$ such that $\bA$ and $\bB$ are unitarily similar.
\end{proof}


\section{Properties  of Eigenvalues of Symmetric Matrices}\label{section:propert_eig_sym}
Eigenvalues are crucial components in various matrix decompositions, including eigenvalue decomposition, Schur decomposition, and spectral decomposition. In the following sections, we will explore several key  facts, inequalities, and computational approaches related to eigenvalues.
Many of the results presented in this section can be generalized from real matrices $\bA\in\real^{n\times n}$ to matrices $\bA\in\sF^{n\times n}$, where $\sF^{n \times n}$ denotes either $\real^{n \times n}$  (real matrices) or $\complex^{n \times n}$ (complex matrices).
Additionally, these results often apply to Hermitian matrices $\bA$ when considering symmetric matrices $\bA$ (with mild adjustments).
\begin{theorem}[Product of Eigenvalues]\label{theorem:eigen_trace}
Let $\bA\in \real^{n\times n}$ be  a square matrix. Then, the product of its $n$ eigenvalues is equal to  the determinant of $\bA$.
\end{theorem}
\begin{proof}[of Theorem~\ref{theorem:eigen_trace}]
Suppose the characteristic polynomial (Definition~\ref{definition:characteristic_polynomial}) is
$$
\begin{aligned}
p_{\bA}(\lambda)=\det(\lambda\bI-\bA ) &=\lambda^n - \gamma_{n-1} \lambda^{n-1} + \ldots + \gamma_1 \lambda  + \gamma_0\\
&=(\lambda-\lambda_1) (\lambda-\lambda_2) \ldots (\lambda-\lambda_n),
\end{aligned}
$$
where $\lambda_1,\lambda_2,\ldots,\lambda_n$ are eigenvalues of $\bA$, and some eigenvalues may be repeated.
We realize that the constant term of $\det(\lambda\bI-\bA ) $ is given by 
$$
\begin{aligned}
f(0) &= \det(-\bA) = (-1)^n\det(\bA)
= (-1)^n \lambda_1 \lambda_2\ldots \lambda_n.
\end{aligned}
$$
Therefore, the product of the $n$ eigenvalues is equal to the determinant of $\bA$.
\end{proof}

\begin{theorem}[Sum of Eigenvalues]\label{theorem:eigen_trace2}
Let $\bA\in \real^{n\times n}$ be  a square matrix. 
Then, the trace of $\bA$ is equal to the sum of its  $n$ eigenvalues.
\end{theorem}
\begin{proof}[of Theorem~\ref{theorem:eigen_trace2}]
By Proposition~\ref{proposition:eigenvalue-similar-matrices}, the traces and eigenvalues of similar matrices are equal. Using the Schur decomposition, we can write   $\bA=\bQ\bU\bQ^\top$, where $\bA$ and $\bU$ are similar matrices. Therefore, the trace of $\bA$ is equal to  the trace of $\bU$, which is the sum of the eigenvalues of $\bA$, as stated in Corollary~\ref{corollary:schur-second-form}.
\end{proof}

\index{Eigenspace}
\index{Interlacing property}
\index{Rayleigh-Ritz}
\index{Poincar{\'e} Separation}
\subsection{Rayleigh-Ritz Theorem and Eigenvalue Interlacing}
The following \textit{Rayleigh-Ritz theorem} relates the Rayleigh quotient and eigenvalues of a symmetric matrix. 
Furthermore, it frequently serves as a tool to verify the validity of different optimization problems, such as PCA and graph drawing  \citep{gallier2017fundamentals}.

\noindent
\begin{theorem}[(Restricted) Rayleigh-Ritz Theorem]\label{theorem:rayleigh_v1}
Let $\bA\in \real^{n\times n}$ be a symmetric matrix, and let $\lambda_1\leq  \lambda_2\leq \ldots \leq \lambda_n$ be its eigenvalues corresponding to mutually orthonormal eigenvectors $\bq_1, \bq_2, \ldots, \bq_n$, respectively.
Suppose the subspace $\mathcalV$ is spanned by $\{\bq_p, \bq_{p+1}, \ldots, \bq_q\}$:
$$
\mathcalV =\spn\left(\bq_p, \bq_{p+1}, \ldots, \bq_q \right),
$$
where $1\leq p\leq q\leq n$. 
For a  unit vector $\bx\in \mathcalV$ satisfying $\normtwo{\bx}=1$, we have 
$$
\lambda_p \leq 
\bx^\top\bA\bx
\leq \lambda_q,
\gap \forall\bx\in\mathcalV \text{ with } \normtwo{\bx}=1
.
\footnote{Similar results can be obtained for complex Hermitian matrices. If we assume $\lambda_1\geq \lambda_2\geq \ldots \geq \lambda_n$, 
then we have 
$
\lambda_q \leq 
\bx^\top\bA\bx
\leq \lambda_p
$.
Similar adjustments apply to the following theorems.
}
$$
\end{theorem}
\begin{proof}[of Theorem~\ref{theorem:rayleigh_v1}]
Since $\bx\in\mathcalV$ is a unit vector, it can be decomposed as $\bx = a_p\bq_p + a_{p+1}\bq_{p+1} +\ldots + a_q\bq_q$ such that $\bx^\top\bq_i = a_i$ for $i\in \{p,p+1,\ldots, q\}$ and $|a_p|^2 +|a_{p+1}|^2 +\ldots +|a_q|^2=1$. Then we have 
$$
\begin{aligned}
\bx^\top\bA\bx &= \bx^\top \bA (a_p\bq_p +\ldots + a_q\bq_q) 
= \bx^\top  (a_p\lambda_p \bq_p  +\ldots + a_q\lambda_q\bq_q) \\
&= a_p^2\lambda_p +\ldots + a_q^2\lambda_q.
\end{aligned}
$$
Since $|a_p|^2 +|a_{p+1}|^2 +\ldots +|a_q|^2=1$, the result  follows immediately.
\end{proof}

The above result can be extended to any vector $\bx$ (not necessarily belonging to $\mathcalV$)
\begin{theorem}[Rayleigh-Ritz Theorem]\label{theorem:rayleigh_v2}
Let $\bA\in \real^{n\times n}$ be  a symmetric matrix, and let $\lambda_1\leq \lambda_2\leq \ldots \leq \lambda_n$ be its eigenvalues corresponding to mutually orthonormal eigenvectors $\bq_1, \bq_2, \ldots, \bq_n$, respectively. If $\bx\in \real^n$ is a unit vector, we have 
$$
\lambda_n \geq 
\bx^\top\bA\bx
\geq \lambda_1, 
\gap \normtwo{\bx}=1.~\footnote{Similar results can be obtained for complex Hermitian matrices.}
$$
\end{theorem}

The Rayleigh-Ritz theorem states that if a symmetric matrix $\bA\in\real^{n\times n}$ has eigenvalues $\lambda_1\leq \lambda_2\leq \ldots \leq \lambda_n$ and orthonormal eigenvectors $\bq_1, \bq_2, \ldots, \bq_n$, then the \textit{Rayleigh quotient} \footnote{We will see its role in computing the eigenvalues of a matrix in Equation~\eqref{equation:raylei_ratio}.}
$r(\bx) = \frac{\bx^\top\bA\bx}{\bx^\top\bx}$ satisfies the following optimization (for general vector $\bx$ that may have any length):
\begin{equation}\label{equation:rr_rr_1}
\begin{aligned}
	&\mathop{\max}_{\bx\neq \bzero}
	\frac{\bx^\top\bA\bx}{\bx^\top\bx}
	=
	\lambda_n
	\qquad\text{and}\qquad
	\mathop{\min}_{\bx\neq \bzero}
	\frac{\bx^\top\bA\bx}{\bx^\top\bx}
	=
	\lambda_1,
\end{aligned}
\end{equation}
with the maximum value $\lambda_n$ attained at $\bx=\alpha \bq_n$ (for nonzero $\alpha$), and the minimum value $\lambda_1$ attained at $\bx=\alpha\bq_1$. Or if $\mathcalV$ is the subspace spanned by  $\{\bq_p, \bq_{p+1}, \ldots, \bq_q\}$, then 
\begin{equation}\label{equation:rr_rr_2}
\begin{aligned}
&\mathop{\max}_{\bx\neq \bzero, \bx\in\mathcalV}
\frac{\bx^\top\bA\bx}{\bx^\top\bx}
=
\lambda_q
\qquad\text{and}\qquad
\mathop{\min}_{\bx\neq \bzero, \bx\in\mathcalV}
\frac{\bx^\top\bA\bx}{\bx^\top\bx}
=
\lambda_p,
\end{aligned}
\end{equation}
with the maximum value $\lambda_q$ achieved at $\bx=\alpha\bq_q$, and the minimum value $\lambda_p$ achieved at $\bx=\alpha\bq_p$.
Let $\bx\triangleq\be_i$ for $i\in\{1,2,\ldots,n\}$. This implies
\begin{equation}\label{equation:symm_ei_dia_ineq}
\textbf{(Symmetric $\bA$): }\qquad \lambda_{\min}(\bA) \leq d_{\min}(\bA) \leq d_{\max}(\bA) \leq \lambda_{\max}(\bA),
\end{equation} 
where $d_{\min}(\bA)$ and $d_{\max}(\bA)$ are the smallest and largest diagonal values of $\bA$, respectively.

\begin{remark}[Eigenvalues as Optimization]
\index{Lagrange multiplier}
Equation~\eqref{equation:rr_rr_1} can also be demonstrated using the following constrained optimization problem:
\begin{equation}\label{equation:eign_opt}
	\mathop{\max}_{\bx \in \real^n} \bx^\top\bA\bx, \qquad \text{s.t.,} \gapthree \norm{\bx}_2=1.
\end{equation}
By forming the Lagrangian, the optimization can be transformed into 
$$
\mathcal{L}(\bx, \lambda)= \bx^\top\bA\bx - \lambda(\bx^\top\bx-1),
$$
where $\lambda$ is called the \textit{Lagrange multiplier}. To find the solution, it is necessary for the gradient of the Lagrangian to be zero at the optimal solution $\bx^*$:
$
\nabla_{\bx} \mathcal{L}(\bx, \lambda) = 2\bA\bx-2\lambda\bx = \bzero.
$
This implies $\bA\bx = \lambda\bx$  and  the optimal points indicate $\lambda$ and $\bx^*$ are an eigenvalue and the corresponding eigenvector of $\bA$, respectively. \eqref{equation:rr_rr_1} shows that the solution to \eqref{equation:eign_opt} corresponds to the largest real eigenvalue of $\bA$.
\end{remark}

\begin{remark}[Positive Definiteness from Rayleigh-Ritz Theorem]
The Rayleigh-Ritz theorem immediately implies that  a symmetric matrix $\bA$ is positive definite if and only if all its eigenvalues are positive, and is positive semidefinite if and only if all its eigenvalues are nonnegative.
\end{remark}

\index{Bendixon's theorem}
\index{Hirsch's theorem}

Bendixon's and Hirsch's theorem follow directly from the Rayleigh-Ritz theorem by noticing that $\text{Re}(\lambda) = \frac{\lambda+\lambda^*}{2}$ and $\text{Im}(\lambda) = \frac{\lambda-\lambda^*}{2i}$ if the eigenvalue $\lambda$ is complex.
\begin{exercise}[Bendixson’s and Hirsch’s theorem]\label{exercise:bendixson_hirch}
Let $\bA \in \sF^{n \times n}$, and let $\lambda \in \Lambda(\bA)$. Show that following inequalities hold:~\footnote{$\sF^{n \times n}$ denotes either $\real^{n \times n}$ or $\complex^{n \times n}$; $\text{Re}(x)$ and $\text{Im}(x)$ denote the real and imaginary parts of $x$, respectively.}
\begin{itemize}
\item $\lambda_{\min}\left[\frac{1}{2}(\bA + \bA^*)\right] \leq \text{Re}(\lambda) \leq \lambda_{\max}\left[\frac{1}{2}(\bA + \bA^*)\right]$.
\item $\lambda_{\min}\left[\frac{1}{2i}(\bA - \bA^*)\right] \leq \text{Im}(\lambda) \leq \lambda_{\max}\left[\frac{1}{2i}(\bA - \bA^*)\right]$.
\end{itemize}
\end{exercise}

The Rayleigh-Ritz theorem can also be applied non-symmetric matrices.
\begin{exercise}[Eigenvalues of Non-Symmetric Matrices]
Let $\lambda_{1}, \lambda_{2}, \ldots, \lambda_{n}$ be the eigenvalues of $\bA \in \real^{n\times n}$ (not necessarily symmetric). Show that
$$
\min_{\bx \neq \bzero} \left| \frac{\bx^\top \bA \bx}{\bx^\top \bx} \right| \leq \left| \lambda_i \right| \leq \max_{\bx \neq \bzero} \left| \frac{\bx^\top \bA \bx}{\bx^\top \bx} \right|, \gap i \in \{1, 2, \ldots, n\},
$$
and that strict inequality is possible in either inequality.
\textit{Hint: Use Problem~\ref{prob:sam_quad}.}
\end{exercise}

\index{Eigenvalue interlace}
\index{Eigenvalue interlacing}
\index{Interlace}

The Rayleigh-Ritz theorem establishes a connection between the eigenvalues of two symmetric matrices, $\bA$ and $\bB=\bZ^\top\bA\bZ$, even when they have different sizes. Here, $\bZ$ is a rectangular matrix with the property $\bZ^\top\bZ=\bI$ (meaning $\bZ$ is semi-orthogonal). To understand this relationship, we introduce the concept of \textit{interlace}.
\begin{definition}[Eigenvalue Interlacing]
Let $\bA\in \real^{n\times n}$ and $\bB\in\real^{m\times m}$ be  symmetric matrices with $m\leq n$. Let $\lambda_1 \leq \lambda_2\leq \ldots \leq \lambda_n$ be the eigenvalues of $\bA$ and $\mu_1\leq \mu_2\leq \ldots \leq\mu_m$ be the eigenvalues of $\bB$. Then, we say that the eigenvalues of $\bB$ interlace the eigenvalues of $\bA$ if 
$$
\lambda_i \leq \mu_i \leq \lambda_{n-m+i}, \gap \forall \, i\in\{1,2,\ldots,m\}.
$$
For example, 
$$
\begin{aligned}
\lambda_1 &\leq \mu_1 \leq \lambda_4, 
\gap
\lambda_2  \leq \mu_2 \leq \lambda_5,
\end{aligned}
$$
if $n=5, m=2$.
\end{definition}

\begin{theorem}[Eigenvalue Interlacing]\label{theorem:eigen_interlac}
Let $\bA\in \real^{n\times n}$ be a symmetric matrix, and let  $\bZ\in\real^{n\times m}$ be a semi-orthogonal matrix such that $\bZ^\top\bZ=\bI$ with $m\leq n$. Define $\bB\triangleq\bZ^\top\bA\bZ\in\real^{m\times m}$. Suppose $\lambda_1\leq \lambda_2\leq \ldots \leq \lambda_n$ are the eigenvalues of $\bA$, 
and $\mu_1\leq \mu_2\leq \ldots \leq\mu_m$ are the eigenvalues of $\bB$. Then we have:
\begin{enumerate}
\item The eigenvalues of $\bB$ interlace the eigenvalues of $\bA$.
\item If $\lambda_i=\mu_i$, then there is an eigenvector $\by$ of $\bB$ associated with the eigenvalue $\mu_i$ such that $\bZ\by$ is an eigenvector of $\bA$ associated with the eigenvalue $\lambda_i$.
\end{enumerate}
\end{theorem}
\begin{proof}[of Theorem~\ref{theorem:eigen_interlac}]
1). 
Suppose $\{\bu_1, \bu_2, \ldots, \bu_n\}\in\real^n$ and $\{\bv_1, \bv_2, \ldots, \bv_m\}\in\real^m$ are orthonormal eigenvector bases for $\bA$ and $\bB$, respectively. 
Let further $\mathcalU_j$ be the subspace spanned by $\{\bu_1, \bu_2, \ldots, \bu_j\}$, and  $\mathcalV_j$ be the subspace spanned by $\{\bv_1, \bv_2, \ldots, \bv_j\}$.
For any $i\in\{1,2\ldots, m\}$, 
$$
\begin{aligned}
&\text{the subspace $\mathcalV_i  \subseteq \real^m$ has dimension $i$};\\
&\text{the subspace $\bZ^\top\mathcalU_{i-1}\subseteq \real^m$ has dimension at most $i-1$};\\
&\text{the subspace $(\bZ^\top\mathcalU_{i-1})^\perp \subseteq \real^m$ has dimension at least $m-i+1$}.
\end{aligned}
$$
Since $\dim(\mathcalV_i)+ \dim((\bZ^\top\mathcalU_{i-1})^\perp)>m$, there exists a nonzero vector $\by\in \mathcalV_i \cap (\bZ^\top\mathcalU_{i-1})^\perp$.
It follows that 
$$
\by^\top \bZ^\top \bu_k = (\bZ\by)^\top\bu_k=0, \gap k\in\{1,2,\ldots, i-1\}.
$$
That is, $\bZ\by \in (\mathcalU_{i-1})^\perp$.
By the Rayleigh-Ritz theorem and the fact that $\bZ^\top\bZ=\bI$, we can establish the following inequality:
$$
\lambda_i \leq 
\frac{(\bZ\by)^\top\bA(\bZ\by)}{(\bZ\by)^\top(\bZ\by)}
=
\frac{\by^\top\bB\by}{\by^\top\by}.
$$
Moreover, applying the Rayleigh-Ritz theorem again, we get 
$$
\mu_i = 
\mathop{\max}_{\bx\neq 0, \bx\in\{\bv_{i+1}, \ldots, \bv_m\}^\perp}
\frac{\bx^\top\bB\bx}{\bx^\top\bx}
=
\mathop{\max}_{\bx\neq 0, \bx\in\{\bv_1, \ldots, \bv_i\}}
\frac{\bx^\top\bB\bx}{\bx^\top\bx}.
$$
This implies 
$
\frac{\bw^\top\bB\bw}{\bw^\top\bw}
\leq \mu_i,
\forall\, \bw\in \mathcalV_i.
$
Since $\by\in \mathcalV_i$, we have 
$$
\lambda_i\leq 
\frac{\by^\top\bB\by}{\by^\top\by}
\leq \mu_i , 
\gap
\forall \,
i\in\{1,2,\ldots,m\}.
$$
Apply this finding to the symmetric matrices $-\bA$ and $-\bB$, we conclude that 
$$
-\lambda_{n-m+i} \leq -\mu_i.
$$
We conclude the eigenvalues of $\bB$ interlace those of $\bA$.

2). If $\lambda_i=\mu_i$, we have 
$$
\lambda_i = \frac{(\bZ\by)^\top\bA(\bZ\by)}{(\bZ\by)^\top(\bZ\by)}
=
\frac{\by^\top\bB\by}{\by^\top\by}
=
\mu_i.
$$
Therefore, $\by$ is an eigenvector of $\bB$, and $\bZ\by$ is an eigenvector of $\bA$ corresponding to the same eigenvalue $\lambda_i=\mu_i$.
\end{proof}
\begin{remark}[Sum of Eigenvalues]
Let $\bA\in\real^{n\times n}$ be symmetric with eigenvalues $\lambda_1\leq  \lambda_2\leq \ldots \leq \lambda_n$. Given the fact that the trace of a symmetric matrix is the sum of the eigenvalues, we can apply the eigenvalue interlacing theorem to observe:
$$
\lambda_1 + \lambda_2+\ldots +\lambda_m 
\leq \trace(\bZ^\top\bA\bZ)
\leq 
\lambda_{n-m+1} +\lambda_{n-m+2} +\ldots + \lambda_{n}.
$$
\end{remark}

The \textit{Poincar{\'e} separation theorem} follows immediately from the eigenvalue interlacing theorem (they represent the same fundamental concept in different forms). The Poincar{\'e} separation theorem is widely used in multivariate analysis and sample surveys \citep{rao1979separation, scott1985separation}.
\begin{theorem}[Poincar{\'e} Separation Theorem]
Let $\bA\in\real^{n\times n}$ be a symmetric matrix with eigenvalues $\lambda_1(\bA)\leq \lambda_2(\bA)\leq\ldots \leq \lambda_n(\bA)$. Suppose $\bU=[\bu_1, \bu_2, \ldots, \bu_m]\in \real^{n\times m}$ contains $m$ orthonormal vectors with $1\leq m\leq n$.
Let $\bB\triangleq\bU^\top\bA\bU \in \real^{m\times m}$  whose eigenvalues satisfy
$\lambda_1(\bB)\leq \lambda_2(\bB)\leq\ldots \leq \lambda_m(\bB)$. Then we have 
$$
\lambda_k(\bA) \leq \lambda_k(\bB) \leq \lambda_{n-m+k}(\bA), \gap \forall k\in\{1,2,\ldots,m\}.
$$
\end{theorem}

\index{Cauchy interlacing theorem}
We can construct the semi-orthogonal matrix $\bU$  using special matrices to obtain specific cases of the eigenvalue interlacing property.
\begin{remark}[Special Case and Cauchy Interlacing Theorem]\label{remark:poincare_sep_eign}
Let $\bA\in\real^{n\times n}$ be a given symmetric matrix.
Let further $\bP \triangleq \bI[1:n,1:n-1]$ be the $n\times (n-1)$ matrix obtained by dropping out the last column of an $n\times n$ identity matrix. We have $\bP^\top\bP=\bI_{n-1}$, and $\bB=\bP^\top\bA\bP$ is the matrix obtained from $\bA\in\real^{n\times n}$ by deleting its last row and  last column, i.e., $\bB=\bA[1:n-1,1:n-1]$ and $\bA\triangleq\scriptsize\begin{bmatrix}
\bB & \bc\\
\bc^\top & d
\end{bmatrix}$. The Poincar{\'e} separation theorem then implies that if $\bA$ is symmetric, 
\begin{equation}\label{equation:cauch_ineq_spe}
\lambda_1 \leq \mu_1 \leq 
\lambda_2 \leq \mu_2 \leq 
\ldots 
\leq 
\mu_{n-2} \leq \lambda_{n-1}
\leq 
\mu_{n-1} \leq \lambda_n,
\end{equation}
where $\lambda_k$ is the $k$-th smallest eigenvalue of $\bA$ and $\mu_k$ is the $k$-th smallest eigenvalue of $\bB$.
In this context, for some $k\in\{1,2,\ldots,n-1\}$, 
\begin{itemize}
\item  $\lambda_k = \mu_k$ if and only if there exists a nonzero $\bx\in\real^{n-1}$ such that $\bB\bx = \mu_k\bx$,
$\bc^\top \bx = 0$, and $\bB\bx = \lambda_k\bx$.
\item $\mu_k = \lambda_{k+1}$ if and only if there exists a nonzero $\bx\in\real^{n-1}$ such that $\bB\bx = \mu_k\bx$, $\bc^\top\bx = 0$, and $\bB\bx = \lambda_{k+1}\bx$. 
\item If no eigenvector of $\bB$ is orthogonal to $\bx$, then every inequality in \eqref{equation:cauch_ineq_spe} is strict.
\end{itemize}

If we further set $\bP\triangleq\bI[1:n,1:m]$, i.e., dropping out the last $n-m$ columns of an identity matrix, then  $\bB=\bP^\top\bA\bP$ is the matrix obtained from $\bA$ by deleting its last $n-m$ rows and its last $n-m$ columns, i.e., $\bB=\bA[1:m,1:m]$ (the $m\times m$ leading principal submatrix of $\bA$). We obtain the \textit{Cauchy interlacing theorem}:
$$
\lambda_k \leq \mu_k \leq \lambda_{n-m+k}, \gap k\in\{1,2,\ldots, m\},
$$
where $\lambda_k$ is the $k$-th smallest eigenvalue of $\bA$ and $\mu_k$ is the $k$-th smallest eigenvalue of $\bB$.
Apparently, these results  can be  extended to any $m\times m$ principal submatrix of $\bA$ (not necessarily its leading principal submatrix).
\end{remark}

\index{Courant-Fischer}
\subsection{Courant-Fischer, Fan, Weyl, and Eigenvalue Interlacing}
The \textit{Courant-Fischer theorem}, also known as the \textit{Min-max theorem}, can be employed to establish other significant theorems.
\begin{theorem}[Courant-Fischer Theorem]\label{theorem:courant_fischer}
Let $\bA\in \real^{n\times n}$ be  a symmetric matrix with eigenvalues $\lambda_1\leq  \lambda_2\leq \ldots \leq \lambda_n$.
Let further $\sV_k$ denote the set of subspaces of $\real^n$ with dimension $k$ with $k\in\{1,2,\ldots,n\}$. Then,
$$
\begin{aligned}
\lambda_k &= 
\mathop{\max}_{\mathcalS \in \sV_{n-k+1}}
\mathop{\min}_{\bx\neq \bzero,\bx\in\mathcalS }
\frac{\bx^\top\bA\bx}{\bx^\top\bx}
\gap \text{and}\gap 
\lambda_k =
\mathop{\min}_{\mathcalS \in \sV_{k}}
\mathop{\max}_{\bx\neq \bzero,\bx\in\mathcalS }
\frac{\bx^\top\bA\bx}{\bx^\top\bx}.~\footnote{The result can be generalized to the complex Hermitian matrix $\bA$, and $\sV_k$ then denotes the set of subspace of $\complex^n$.}
\end{aligned}
$$
\end{theorem}
\begin{proof}[of Theorem~\ref{theorem:courant_fischer}]
For simplicity, we will only prove the second equality here, as the first one can be demonstrated similarly.
Suppose $\{\bu_1, \bu_2, \ldots, \bu_n\}$ are the orthonormal eigenvectors of $\bA$, where $\bu_i$ is a unit eigenvector corresponding to $\lambda_i$. Then the subspace $\mathcalV_k$ spanned by 
$\{\bu_1, \bu_2, \ldots, \bu_k\}$ has dimension $k$, and $\mathcalV_k \in \sV_k$.
By the Rayleigh-Ritz theorem, we have 
$$
\lambda_k = 
\mathop{\max}_{\bx\neq\bzero, \bx\in\mathcalV_k}
\frac{\bx^\top\bA\bx}{\bx^\top\bx}
\geq 
\mathop{\inf}_{\mathcalS\in \sV_k}
\mathop{\max}_{ {\bx\neq \bzero,\\\bx\in\mathcalS} }
\frac{\bx^\top\bA\bx}{\bx^\top\bx}.
$$
Therefore, we need to show $\lambda_k
\leq 
\mathop{\max}_{\bx\neq \bzero,\bx\in\mathcalS }
\frac{\bx^\top\bA\bx}{\bx^\top\bx}$ for all $\mathcalS \in \sV_k$ to prove the second equality in Theorem~\ref{theorem:courant_fischer}.

To see this,
we note that $\dim(\mathcalV_{k-1}^\perp) = n - \dim(\mathcalV_{k-1}) = n-k+1$, and $\dim(\mathcalS)=k$ since $\mathcalS\in \sV_k$. Consequently, this implies $\dim(\mathcalS) + \dim(\mathcalV_{k-1}^\perp) = n+1$, leading to $\dim(\mathcalS\cap \mathcalV_{k-1}^\perp)> 0$.
Then for any nonzero $\bv \in \mathcalS\cap \mathcalV_{k-1}^\perp$, by the Rayleigh-Ritz theorem, we have 
$$
\lambda_k 
= \mathop{\min}_{\bx\neq\bzero, \bx\in\mathcalV_{k-1}^\perp} \frac{\bx^\top\bA\bx}{\bx^\top\bx}
\leq 
\frac{\bv^\top\bA\bv}{\bv^\top\bv}
\leq 
\mathop{\max}_{\bx\neq\bzero, \bx\in\mathcalS} \frac{\bx^\top\bA\bx}{\bx^\top\bx}.
$$
The result then follows:
$
\lambda_k =
\mathop{\min}_{\mathcalS \in \sV_{k}}
\mathop{\max}_{\bx\neq \bzero,\bx\in\mathcalS }
\frac{\bx^\top\bA\bx}{\bx^\top\bx}.
$
\end{proof}
When $k=n$ in the second equality of Theorem~\ref{theorem:courant_fischer} (or $k=1$ in its first equality), we may omit the outer optimization and set $\mathcalS=\real^n$. This reduces to the Rayleigh-Ritz Theorem~\ref{theorem:rayleigh_v2}.
\begin{proposition}[Courant-Fischer Theorem, CNT.]\label{proposition:cour_fish}
Let $\bA\in\real^{n\times n}$ be  a symmetric matrix with eigenvalues $\lambda_1\leq \lambda_2\leq \ldots\leq \lambda_n$. Let $\mathcalV_k$ be a $k$-dimensional subspace of $\real^n$: $\dim(\mathcalV_k)=k$. Then, 
\begin{itemize}
\item If $\frac{\bx^\top\bA\bx}{\bx^\top\bx}\geq \gamma$ (resp. $>\gamma$) for any vector $\bx\in\mathcalV_k$, then $\lambda_{n-k+1}\geq \gamma$ (resp. $>\gamma$).
\item If $\frac{\bx^\top\bA\bx}{\bx^\top\bx}\leq  \gamma$ (resp. $<\gamma$) for any vector $\bx\in\mathcalV_k$, then $\lambda_{k}\leq  \gamma$ (resp. $<\gamma$).
\end{itemize}
\end{proposition}
\begin{proof}[of Proposition~\ref{proposition:cour_fish}]
Suppose $\{\bu_1, \bu_2, \ldots, \bu_n\}$ are the orthonormal eigenvectors of $\bA$, where $\bu_i$ is a unit eigenvector corresponding to $\lambda_i$. 
Let $\mathcalW \triangleq  \spn\{\bu_1, \bu_2, \ldots, \bu_{n-k+1}\}$. This indicates $\dim(\mathcalV_k)+\dim(\mathcalW)=n+1$. Therefore, there exist a nonzero vector $\bv\in \mathcalV_k\cap \mathcalW$.
By the Rayleigh-Ritz Theorem~\ref{theorem:rayleigh_v1} and the assumption, we have  $\gamma \leq \frac{\bv^\top\bA\bv}{\bv^\top\bv}\leq \lambda_{n-k+1} $. 
The other results follow similarly by applying the observation to $\bA$ or $-\bA$.
\end{proof}

The Courant-Fischer theorem finds the optimization for individual eigenvalues. The Fan theorem, on the other hand, finds the optimization for the sum of largest and smallest $k$ eigenvalues.
\begin{theorem}[Ky Fan Eigenvalue Theorem \citep{fan1949theorem, fan1950theorem}\index{Ky Fan eigenvalue theorem}]\label{theorem:kyfan_eig_th}
Let $\bA\in\real^{n\times n}$ be symmetric with eigenvalues $\lambda_1(\bA) \geq  \lambda_2(\bA)\geq \ldots \geq \lambda_n(\bA)$. The following inequalities hold:
\begin{equation}\label{equation:kyfan_eig_th1}
\max_{\bU^\top\bU = \bI_k} \trace(\bU^\top\bA\bU) = \sum_{i=1}^{k} \lambda_i(\bA)
\quad \text{and}\quad
\min_{\bU^\top\bU = \bI_k} \trace(\bU^\top\bA\bU) = \sum_{i=n-k+1}^{n} \lambda_i(\bA),
\end{equation}
where $\bU$ is semi-orthogonal  satisfying $\bU^\top\bU = \bI_k$.
Let $\bA,\bB\in\real^{n\times n}$ be symmetric, \eqref{equation:kyfan_eig_th1} indicates
\begin{equation}
\sum_{i=1}^{k} \lambda_i(\bA+\bB) \leq \big[\sum_{i=1}^{k} \lambda_i(\bA) +\lambda_i(\bB)\big].
\end{equation}
\end{theorem}
\begin{proof}[of Theorem~\ref{theorem:kyfan_eig_th}]
The full proof relies on the majorization principle, we only sketch the main steps here; see \citet{fan1949theorem, fan1950theorem, marshall1979inequalities} for a detailed proof.
Let $\bA = \bQ\bLambda\bQ^\top$ be the spectral decomposition of $\bA$,  the optimizations in \eqref{equation:kyfan_eig_th1} become
$$
\max_{\bV^\top\bV = \bI_k} \trace(\bV^\top\bLambda\bV)
\quad\text{and}\quad
\min_{\bV^\top\bV = \bI_k} \trace(\bV^\top\bLambda\bV),
$$
where $\bV$ is also semi-orthogonal  satisfying $\bV^\top\bV = \bI_k$.
The maximum and minimum values are attained at setting the columns of $\bV$ to the corresponding eigenvectors of $\bA$.
\end{proof}

The Courant-Fischer theorem is a valuable tool for proving various theorems. For instance, it can be applied to establish both lower and upper bounds for the $k$-th eigenvalue of the sum of matrices $\bA+\bB$.
\begin{theorem}[Weyl's Theorem\index{Weyl's theorem}]\label{theorem:weyl}
Let $\bA,\bB\in\real^{n\times n}$ be symmetric. Let further $\{\lambda_j(\bA)\}_{j=1}^n$, $\{\lambda_j(\bB)\}_{j=1}^n$, and $\{\lambda_j(\bA+\bB)\}_{j=1}^n$ denote the set of eigenvalues of $\bA, \bB,$ and $\bA+\bB$ in nondecreasing order, respectively.
Then, for any $k\in\{1,2,\ldots, n\}$, we have 
\begin{equation}\label{equation:Weyl_main1}
\lambda_{k-i+1}(\bA) +\lambda_i(\bB)
\leq 
\lambda_k(\bA+\bB)
\leq 
\lambda_{k+j}(\bA) +\lambda_{n-j}(\bB), 
\end{equation}
where $i\in\{1,2,\ldots,k\},j\in\{0,1,\ldots, n-k\}$. 
The equality in the upper bound holds  for some pair $k, j$ if and only if there is a nonzero
vector $\bx$ such that $\bA\bx =\lambda_{k+j}(\bA)\bx$, $\bB\bx = \lambda_{n-j}(\bB)\bx$, and $(\bA + \bB)\bx = \lambda_k(\bA + \bB)\bx$;
and the equality in the lower bound holds for some pair $k,i$ if and only if there exists a nonzero vector $\bx$ such that $\bA\bx = \lambda_{k- i+1}(\bA)\bx$, $\bB\bx = \lambda_i(\bB)\bx$, and $(\bA + \bB)\bx = \lambda_k(\bA + \bB)x$.
If $\bA$ and $\bB$ have no common eigenvector, then the inequalities are strict.

\paragraph{Other implications.}
The result for $i=1,j=0$ is frequently used, for example,  in low-rank approximation results:
\begin{equation}
\lambda_k(\bA) +\lambda_1(\bB)
\leq 
\lambda_k(\bA+\bB)
\leq 
\lambda_k(\bA) +\lambda_n(\bB).
\end{equation}
Due to symmetry, we also have 
\begin{equation}
\lambda_k(\bB) +\lambda_1(\bA)
\leq 
\lambda_k(\bA+\bB)
\leq 
\lambda_k(\bB) +\lambda_n(\bA).
\end{equation}
When $\{k=i=1, j=0\}$ and $\{k=n,i=1, j=0\}$, we have, respectively, 
\begin{equation}
\begin{aligned}
\lambda_1(\bB) +\lambda_1(\bA)
&\leq 
\lambda_1(\bA+\bB)
\leq 
\lambda_1(\bB) +\lambda_n(\bA);\\
\lambda_n(\bB) +\lambda_1(\bA)
&\leq 
\lambda_n(\bA+\bB)
\leq 
\lambda_n(\bB) +\lambda_n(\bA).
\end{aligned}
\end{equation}
Proposition~\ref{propo:sing_eig_ineq1} also shows 
\begin{equation}
\lambda_{1}(\bA) + \lambda_{1}(\bB) \leq \lambda_{1}(\bA+\bB) \leq \lambda_{n}(\bA+\bB) \leq \lambda_{n}(\bA) + \lambda_{n}(\bB).
\end{equation}
\end{theorem}
\begin{proof}[of Theorem~\ref{theorem:weyl}]
For brevity, we only prove the result when $i=1,j=0$. The general case can be derived similarly.
According to the Rayleigh-Ritz theorem, we have 
$$
\lambda_1(\bB) \leq 
\frac{\bx^\top\bB\bx}{\bx^\top\bx}
\leq \lambda_n(\bB), \gap 
\text{for all }\bx\neq \bzero.
$$
According to the  Courant-Fischer theorem, given the set $\sV_k$ of subspaces  of $\real^n$ with dimension $k$ and for $k\in\{1,2,\ldots,n\}$, we have
$$
\begin{aligned}
\lambda_k(\bA+\bB)
&=
\mathop{\min}_{\mathcalS \in \sV_{k}}
\mathop{\max}_{\bx\neq \bzero,\bx\in\mathcalS }
\frac{\bx^\top (\bA+\bB)\bx}{\bx^\top\bx}
=
\mathop{\min}_{\mathcalS \in \sV_{k}}
\mathop{\max}_{\bx\neq \bzero,\bx\in\mathcalS }
\big(\frac{\bx^\top \bA\bx}{\bx^\top\bx} +\frac{\bx^\top \bB\bx}{\bx^\top\bx}\big)\\
&\geq 
\mathop{\min}_{\mathcalS \in \sV_{k}}
\mathop{\max}_{\bx\neq \bzero,\bx\in\mathcalS }
\frac{\bx^\top \bA\bx}{\bx^\top\bx} +\lambda_1(\bB)
=
\lambda_k(\bA)+\lambda_1(\bB).
\end{aligned}
$$
Similarly, we have $\lambda_k(\bA+\bB)
\leq 
\lambda_k(\bA) +\lambda_n(\bB)$. We complete the proof.
\end{proof}

Weyl's theorem describes how the eigenvalues of a symmetric matrix  $\bA\in\real^{n\times n}$ can change when it is perturbed by another symmetric matrix.
Specifically, if a symmetric matrix $\bE$ is a perturbation of a symmetric matrix $\bA$, Weyl's theorem reveals the important inequality
\begin{equation}
\abs{\lambda_k(\bA) - \lambda_k(\bA+\bE)} \leq \normtwo{\bE},\gap k\in\{1,2,\ldots, n\}.
\end{equation}
That is, the eigenvalues of a real symmetric matrix  are stable under  small perturbations.
Let $\bLambda_A$ and $\bLambda_{A+E}$ denote the diagonal matrices containing eigenvalues  of $\bA$ and $\bA+\bE$, respectively. This also implies, by the definition of the spectral norm (Definition~\ref{definition:spectral_norm}):
\begin{equation}
	\normtwo{\bLambda_A-\bLambda_{A+E}}\leq \normtwo{\bE}.
\end{equation}
A variant of this result for unitarily invariant norms (Definition~\ref{definition:unitarily_invaria}), known as \textit{Mirsky's theorem},   is discussed in Problem~\ref{prob:pert_eig_uniinvar}.
\index{Mirsky's theorem}

\begin{exercise}[Perturbation on Singular Values]\label{exercise:perturb_sings}
As a hindsight, we will introduce singular values in the next section. 
Show that when a  matrix $\bE\in\real^{m\times n}$ is a perturbation on a  matrix $\bA\in\real^{m\times n}$, Weyl's theorem reveals the inequality
\begin{equation}
\abs{\sigma_k(\bA) - \sigma_k(\bA+\bE)} \leq \normtwo{\bE}, 
\gap \forall k\in\{1,2,\ldots,\min\{m,n\}\},
\end{equation}
where $\sigma_k(\cdot)$ is the $k$-th largest singular value of a matrix.
Let $\bSigma_A$ and $\bSigma_{A+E}$ denote the diagonal matrices containing singular values of $\bA$ and $\bA+\bE$, respectively. This also implies, by the definition of the spectral norm (Definition~\ref{definition:spectral_norm}):
\begin{equation}
\normtwo{\bSigma_A-\bSigma_{A+E}}\leq \normtwo{\bE}.~\footnote{The bound in terms of the Frobenius norm is discussed in Theorem~\ref{theorem:vonneu_ineq2_sings}.}
\end{equation}
This indicates that the singular values of a rectangular matrix are also stable under small perturbations.
\end{exercise}

If a given symmetric matrix is modified  by bordering, then the new and old eigenvalues must interlace (Remark~\ref{remark:poincare_sep_eign}). 
Weyl's theorem reveals that if a given symmetric matrix is modified  by adding a rank-one symmetric matrix, the new and old eigenvalues also interlace.
To see this, we need the following  results.

\begin{corollary}\label{corollary:weyl_gens}
Let $ \bA, \bB \in\real^{n\times n} $ be symmetric. Suppose that $ \bB $ has exactly $\beta$ positive eigenvalues and exactly $\alpha $ negative eigenvalues. 
Let further $\{\lambda_j(\bA)\}_{j=1}^n$, $\{\lambda_j(\bB)\}_{j=1}^n$, and $\{\lambda_j(\bA+\bB)\}_{j=1}^n$ denote the set of eigenvalues of $\bA, \bB,$ and $\bA+\bB$ in nondecreasing order, respectively.
Then,
\begin{equation}\label{equation:weyl_gens1}
\lambda_k(\bA+\bB) \leq \lambda_{k+\beta}(\bA), \quad k\in\{1,2\ldots, n-\beta\} ,
\end{equation}
with equality for some $ k $ if and only if $ \bB $ is singular and there exists a nonzero vector $ \bx $ such that $ \bA \bx = \lambda_{k+\beta}(\bA) \bx $, $ \bB \bx = \bzero $, and $ (\bA+\bB) \bx = \lambda_k(\bA+\bB) \bx $. Similarly, we have 
\begin{equation}\label{equation:weyl_gens2}
\lambda_{k-\alpha}(\bA) \leq \lambda_k(\bA+\bB), \quad k\in\{\alpha+1,\alpha+2\ldots, n\},
\end{equation}
with equality for some $ k $ if and only if $ \bB $ is singular and there exists a nonzero vector $ \bx $ such that $ \bA \bx = \lambda_{k-\alpha}(\bA) \bx $, $ \bB \bx = \bzero $, and $ (\bA+\bB) \bx = \lambda_k(\bA+\bB) \bx $. 
Every inequality in \eqref{equation:weyl_gens1} and \eqref{equation:weyl_gens2} is strict  if either (a). $ \bB $ is nonsingular or (b). $ \bB \bx \neq \bzero $ for every eigenvector $ \bx $ of $ \bA $.
\end{corollary}
\begin{proof}[of Corollary~\ref{corollary:weyl_gens}]
We first note that if the symmetric matrix $\bB $ has exactly $ \beta $ positive eigenvalues and exactly $\alpha$ negative eigenvalues, 
then $ \lambda_{n-\beta}(\bB) \leq 0 $ and $ \lambda_{\alpha+1}(\bB) \geq 0 $, with equality if and only if $ n > \beta + \alpha $, that is, if and only if $ \bB $ is singular.
Let $ j = \beta $ in \eqref{equation:Weyl_main1}, we can obtain $ \lambda_k(\bA+\bB) \leq \lambda_{k+\beta}(\bA) + \lambda_{n-\beta}(\bB) \leq \lambda_{k+\beta}(\bA) $ with equality if and only if $ \bB $ is singular and there exists a nonzero vector $ \bx $ such that $ \bA \bx = \lambda_{k+\beta}(\bA) \bx $, $ \bB \bx = \bzero $, and $ (\bA+\bB) \bx = \lambda_k(\bA+\bB) \bx $. A similar argument shows that \eqref{equation:weyl_gens2} follows  by letting $i = \alpha + 1 $ in \eqref{equation:Weyl_main1}.
\end{proof}
\begin{exercise}
Discuss the results when $\bB$ has exactly one positive eigenvalue and exactly one negative eigenvalues in Corollary~\ref{corollary:weyl_gens}.
\end{exercise}

Similar results can be obtained when a symmetric matrix $\bA$ is perturbed by a singular symmetric matrix $\bB$.
\begin{corollary}[Weyl's Theorem for Singular Symmetric]\label{corollary:sing_weyl_theo}
Let $\bA, \bB \in \real^{n\times n}$ be symmetric. Suppose that $\bB$ is singular with $\rank(\bB) = r$. 
Let further $\{\lambda_j(\bA)\}_{j=1}^n$, $\{\lambda_j(\bB)\}_{j=1}^n$, and $\{\lambda_j(\bA+\bB)\}_{j=1}^n$ denote the set of eigenvalues of $\bA, \bB,$ and $\bA+\bB$ in nondecreasing order, respectively.
Then, 
\begin{equation}\label{equation:sing_sym_int1}
\lambda_k(\bA+\bB) \leq \lambda_{k+r}(\bA), \quad k \in\{1, \ldots, n-r\},
\end{equation}
with equality for some $k$ if and only if $\lambda_{n-r}(\bB) = 0$ and there exists a nonzero vector $\bx$ such that $\bA\bx = \lambda_{k+r}(\bA)\bx$, $\bB\bx = \bzero$, and \((\bA+\bB)\bx = \lambda_k(\bA+\bB)\bx\). Similarly, we have
\begin{equation}\label{equation:sing_sym_int2}
\lambda_{k-r}(\bA) \leq \lambda_k(\bA+\bB), \quad k\in\{r+1, \ldots, n\},
\end{equation}
with equality for some $k$ if and only if $\lambda_{r+1}(\bB) = 0$ and there exists a nonzero vector $\bx$ such that $\bA\bx = \lambda_{k-r}(\bA)\bx$, $\bB\bx = \bzero$, and $(\bA+\bB)\bx = \lambda_k(\bA+\bB)\bx$. If $\bB\bx \neq \bzero$ for every eigenvector $\bx$ of $\bA$, then every inequality in \eqref{equation:sing_sym_int1} and \eqref{equation:sing_sym_int2} is a strict inequality.

\end{corollary}
\begin{proof}[of Corollary~\ref{corollary:sing_weyl_theo}]
We first note that if symmetric $\bB$ is singular with $\rank(\bB)=r$, then $\lambda_{n-r}(\bB) \leq 0$ and $\lambda_{r+1}(\bB) \geq 0$ (since only $r$ eigenvalues are nonzero).
To verify \eqref{equation:sing_sym_int1}, set $j = r$ in \eqref{equation:Weyl_main1}, we obtain
$
\lambda_k(\bA+\bB) \leq \lambda_{k+r}(\bA) + \lambda_{n-r}(\bB) \leq \lambda_{k+r}(\bA)
$
with equality if and only if \(\lambda_{n-r}(\bB) = 0\) and there exists equality in \eqref{equation:Weyl_main1} with $j = r$. A similar argument shows that \eqref{equation:sing_sym_int2} follows from  \eqref{equation:Weyl_main1} with $i = r+1$. 
\end{proof}

%
%

The eigenvalue interlacing property results from Weyl's theorem is then provided below.
\begin{theorem}[Rank-One Eigenvalue Interlacing Theorem]\label{theorem:rkon_eg_inter}
Let $\bA\in\real^{n\times n}$ be symmetric with $n\geq 2$, and let $\bx\in\real^n$ be nonzero. 
Let $\{\lambda_j(\bA)\}_{j=1}^n$  and $\{\lambda_j(\bA+\bx\bx^\top)\}_{j=1}^n$ denote the set of eigenvalues of $\bA$ and $\bA+\bx\bx^\top$ in nondecreasing order, respectively.
Then, it follows that 
\begin{equation}\label{equation:rankone_inter1}
\begin{aligned}
\lambda_{k}(\bA) &\leq \lambda_{k}(\bA + \bx\bx^\top) \leq \lambda_{k+1}(\bA), \quad k\in\{1,2, \ldots, n-1\}, \\
\lambda_{n}(\bA) &\leq \lambda_{n}(\bA + \bx\bx^\top).
\end{aligned}
\end{equation}
The cases of equality in \eqref{equation:rankone_inter1} are as described in Corollary~\ref{corollary:weyl_gens} with $\beta = 1$ and $\alpha = 0$. 
For example, $\lambda_{k}(\bA + \bx\bx^\top) = \lambda_{k+1}(\bA)$ if and only if there exists a nonzero vector $\bv$ such that $\bA \bv = \lambda_{k+1}(\bA) \bv$, $\bx^\top \bv = 0$, and $(\bA + \bx\bx^\top) \bv = \lambda_{k}(\bA + \bx\bx^\top) \bv$. Similarly, we have

\begin{equation}\label{equation:rankone_inter2}
	\begin{aligned}
\lambda_1(\bA - \bx\bx^\top) &\leq \lambda_1(\bA),\\
\lambda_{k-1}(\bA) &\leq \lambda_{k}(\bA - \bx\bx^\top) \leq \lambda_{k}(\bA), \quad k\in\{2, 3,\ldots, n\}.
\end{aligned}
\end{equation}
The cases of equality in \eqref{equation:rankone_inter2} are as described in Corollary~\ref{corollary:weyl_gens} with $\beta = 0$ and $\alpha = 1$. If no eigenvector of $\bA$ is orthogonal to $\bx$, then every inequality in \eqref{equation:rankone_inter1} and \eqref{equation:rankone_inter2} is  strict.
\end{theorem}
\begin{proof}[of Theorem~\ref{theorem:rkon_eg_inter}]
We first note that if  $\bx\in\real^n$ is nonzero and $n\geq 2$, then $\lambda_{n-1}(\bx\bx^\top) = 0 = \lambda_2(\bx\bx^\top)$ since $\rank(\bx\bx^\top)=1$.
Since $\bx\bx^\top$ has the eigenpair $(\normtwo{\bx}^2, \bx)$ and $-\bx\bx^\top$ has the eigenpair $(\normtwo{\bx}^2, -\bx)$, the results follow
by letting $\beta = 1$ and $\alpha = 0$ in \eqref{equation:weyl_gens1}; 
and by letting $i=1$ in \eqref{equation:Weyl_main1} (since the eigenvalue $\lambda_1(\bx\bx^\top)=\normtwo{\bx}^2$ is positive).
The second part can be proved similarly.
\end{proof}

\begin{exercise}[Weyl's Theorem for PSD/PD]
Let $\bA, \bB \in\real^{n\times n}$ be symmetric. Show that if $\bB$ is positive semidefinite, then 
$$
\lambda_k(\bA) \leq \lambda_k(\bA+\bB), \quad k\in\{1,2, \ldots, n\}.
$$
Discuss the conditions under which the equality holds.
Show that if $\bB$ is positive definite, then the inequality becomes strict.
\end{exercise}

\subsection{Eigenvalue Inequalities for PSD/PD}\label{section:eig_ine_pdpsd}

\begin{lemma}[PSD Difference, Multispectrum]\label{lemma:psd_diff_mus}
Let $\bA,\bB\in\real^{n\times n}$ be symmetric with $\bB-\bA\succeq \bzero$ and $\widehat{\Lambda}(\bA)=\widehat{\Lambda}(\bB)$ (multispectrum; see Definition~\ref{definition:spectrum}). Then, $\bA=\bB$.
\end{lemma}
\begin{proof}[of Lemma~\ref{lemma:psd_diff_mus}]
Let $\gamma\geq 0$, $\widetildebA\triangleq\bA+\gamma\bI$, and $\widetildebB\triangleq\bB+\gamma\bI$. We have $\widehat{\Lambda}(\widetildebA)=\widehat{\Lambda}(\widetildebB)$ and $\det(\widetildebA)=\det(\widetildebB)$ by definition.
Since $\widetildebA$ and $\widetildebB$ are symmetric, their spectral decompositions show that $\widetildebB\succeq\widetildebA\implies\widetildebA^{-1/2}\widetildebB\widetildebA^{-1/2}\succeq\bI$.
Theorem~\ref{theorem:eigen_charac} shows $\lambda_{\min}(\widetildebA^{-1/2}\widetildebB\widetildebA^{-1/2}) \geq 1$. 
We also have $\det(\widetildebA^{-1/2}\widetildebB\widetildebA^{-1/2}) = \det(\widetildebB)/\det(\widetildebA)=1$.
Therefore, $\lambda_i(\widetildebA^{-1/2}\widetildebB\widetildebA^{-1/2})=1$ for all $i\in\{1,2,\ldots,n\}$. 
Hence, $\widetildebA^{-1/2}\widetildebB\widetildebA^{-1/2}=\bI$, implying $\bA=\bB$.
\end{proof}

\index{Weyl's inequality}
\begin{theorem}[Weyl's Theorem for PSD Difference]\label{theorem:psd_diff}
Let $\bA,\bB\in\real^{n\times n}$ be symmetric with $\bB-\bA\succeq \bzero$. Then, we have 
$$
\lambda_i(\bA) \leq \lambda_i(\bB), \quad \forall i\in\{1,2,\ldots,n\}.
$$
If $\bA\neq \bB$, then there exists $i\in\{1,2,\ldots,n\}$ such that the inequality holds $\lambda_i(\bA) < \lambda_i(\bB)$.
If $\bB-\bA\succ \bzero$, then all the inequalities hold: $\lambda_i(\bA) < \lambda_i(\bB)$, $\forall i\in\{1,2,\ldots,n\}$.
\end{theorem}
\begin{proof}[Theorem~\ref{theorem:psd_diff}]
Theorem~\ref{theorem:eigen_charac} shows $\lambda_{\min}(\bA)\bI\preceq \bA\preceq \bB\preceq \lambda_{\max}(\bB)\bI$, $\lambda_{\min}(\bA)\leq\lambda_{\min}(\bB)$, and $\lambda_{\max}(\bA)\leq \lambda_{\max}(\bB)$.
Let $\bQ^\top\bA\bQ=\bLambda=\diag(\lambda_1(\bA), \lambda_2(\bA), \ldots,\lambda_n(\bA))$ be the spectral decomposition of $\bA$.
Let further $\bA_i\triangleq\diag(\lambda_1(\bA), \lambda_2(\bA), \ldots,\lambda_i(\bA))$ and $\bB_i$ be the leading principal submatrix of $\bQ^\top\bB\bQ$ for all $i\in\{2,3,\ldots, n-1\}$.
Theorem~\ref{theorem:quadratic-pd} shows $\bQ^\top(\bB-\bA)\bQ$ is also PSD if $\bB-\bA$ is PSD.
This implies $\bB_i-\bA_i$ is also PSD and $\lambda_{\min}(\bA_i)\leq \lambda_{\min}(\bB_i)$.
Equation~\eqref{equation:symm_ei_dia_ineq} shows 
$$
\lambda_i(\bA)=\lambda_{\min}(\bA_i)\leq \lambda_{\min}(\bB_i) \leq \lambda_i(\bQ^\top\bB\bQ) =\lambda_i(\bB), \quad \forall i\in\{2,3,\ldots, n-1\}.
$$
If $\bA\neq \bB$, then Lemma~\ref{lemma:psd_diff_mus} shows $\widehat{\Lambda}(\bA)\neq \widehat{\Lambda}(\bB)$, and thus there exists $i\in\{1,2,\ldots,n\}$ such that the inequality holds $\lambda_i(\bA) < \lambda_i(\bB)$.
Finally, if $\bB-\bA\succ \bzero$, we have $\lambda_{\min}(\bA_i)< \lambda_{\min}(\bB_i)$. This completes the proof.
\end{proof}

A direct consequence of  Weyl's theorem is the relationship between the determinants of two matrices.
\begin{corollary}[Trace Results from Weyl's Inequality]
Let $\bA,\bB\in\real^{n\times n}$ be symmetric. Then,
\begin{itemize}
\item If $ \bA \preceq \bB$, then $\trace(\bA) \leq \trace(\bB)$.
\item If $\bA \preceq \bB$ and $\trace(\bA) = \trace(\bB)$, then $\bA = \bB$.
\item If $\bA \prec \bB$, then $\trace(\bA) < \trace(\bB)$.
\item If $\bzero \preceq \bA \preceq \bB$, then $0 \leq \det(\bA) \leq \det \bB$.
\item If $\bzero \preceq \bA \prec \bB$, then $0 \leq \det(\bA) < \det(\bB)$.
\item If $\bzero \prec \bA\preceq \bB$ and $\det(\bA)=\det(\bB)$, then $\bA=\bB$.
\end{itemize}
\end{corollary}

\index{Orthogonal projection}
\index{Projection matrix}
\index{Projector}
\section{Application: Eigenvalue of Projection Matrix}\label{section:spec_app_eigproj}
In Section~\ref{section:application-ls-qr}, we introduced the application of QR decomposition to solve the least squares problem, where we consider the overdetermined system $\bA\bx = \bb$ with $\bA\in \real^{m\times n}$ being the data matrix, and $\bb\in \real^m$ being the observation vector ($m\geq n$). 
Normally, $\bA$ exhibits full column rank since the data from real work has a large chance to be unrelated (or becomes full rank after preprocessing). And the least squares solution is given by $\bx_{LS} = (\bA^\top\bA)^{-1}\bA^\top\bb$ for minimizing $\normtwo{\bA\bx-\bb}^2$, where $\bA^\top\bA$ is invertible since $\bA$ has full column rank and $\rank(\bA^\top\bA)=\rank(\bA)$.
Subsequently, the recovered observation vector is  $\hat{\bb} = \bA\bx_{LS} = \bA(\bA^\top\bA)^{-1}\bA^\top\bb$. 
While $\bb$ may not lie within the column space of $\bA$,  the recovered vector $\hat{\bb}$ does. 
Mathematically, this process can be described using the matrix $\bH=\bA(\bA^\top\bA)^{-1}\bA^\top$, known as the (orthogonal) projection matrix (Appendix~\ref{section:by-geometry-hat-matrix}),
i.e., projecting vector $\bb$ onto the column space of $\bA$. 
Alternatively, $\bH$ is also referred to as the \textit{hat matrix} because it ``put a hat" on  $\bb$. It can be shown that the projection matrix is \textit{symmetric and idempotent} (i.e., $\bH=\bH^\top$ and $\bH^2=\bH$). 
\index{Projection matrix}
\index{Idempotent}
\index{Idempotency}
\begin{remark}[Column Space of Projection Matrix]
The hat matrix $\bH = \bA(\bA^\top\bA)^{-1}\bA^\top$ projects any vector in $\real^m$ onto the column space of $\bA$. That is, $\bH\by \in \cspace(\bA)$ for any $\by\in\real^m$. Since $\bH\by$ is a linear combination of the columns of $\bH$, it implies that $\cspace(\bH) = \cspace(\bA)$. 
In general, for any projection matrix $\bH$ that projects vectors onto a subspace $\mathcalV$, it holds that  $\cspace(\bH) = \mathcalV$. 
This property can be formally proven using singular value decomposition (Chapter~\ref{chapter:SVD}).
\end{remark}

We will now demonstrate that the eigenvalues of any projection matrix are limited to 0 and 1. For a deeper exploration of orthogonal projections, see Appendix~\ref{appendix:orthogonal}.
\begin{proposition}[Eigenvalue of Projection Matrix]\label{proposition:eigen-of-projection-matrix}
The only possible eigenvalues of a projection matrix are   0 and 1.
\footnote{Actually, in Theorem~\ref{theorem:spec_idem}, we proved the only possible eigenvalues of idempotent matrices are 0 and 1. This result demonstrates how to apply spectral decomposition in real-world applications.}
\end{proposition}
\begin{proof}[of Proposition~\ref{proposition:eigen-of-projection-matrix}]
Since the projection matrix $\bH$ is symmetric, we have the spectral decomposition $\bH =\bQ\bLambda\bQ^\top$. Leveraging  the idempotent property of $\bH$, we have
$$
\begin{aligned}
(\bQ\bLambda\bQ^\top)^2 &= \bQ\bLambda\bQ^\top
\quad\implies \quad
\bQ\bLambda^2\bQ^\top = \bQ\bLambda\bQ^\top 
\quad\implies \quad
\bLambda^2 =\bLambda 
\quad\implies \quad
\lambda_i^2 =\lambda_i.
\end{aligned}
$$
Therefore, the only possible eigenvalues for $\bH$ are 0 and 1.
\end{proof}

This property of  projection matrices holds significance in the analysis of distribution theory for linear models \citep{lu2021rigorous}. 
Building on the eigenvalues of the projection matrix, we can also define the perpendicular projection $\bI-\bH$.
\begin{proposition}[Project onto $\mathcalV^\perp$]\label{proposition:orthogonal-projection_tmp}
Let $\mathcalV$ be a subspace, and $\bH$ be the projection matrix onto $\mathcalV$. Then, $\bI-\bH$ serves as the projection matrix onto the orthogonal complement $\mathcalV^\perp$.
\end{proposition}

\begin{proof}[of Proposition~\ref{proposition:orthogonal-projection_tmp}]
First, $(\bI-\bH)$ is symmetric, $(\bI-\bH)^\top = \bI - \bH^\top = \bI-\bH$ since $\bH$ is symmetric. Additionally, we have
$
(\bI-\bH)^2 = \bI^2 -\bI\bH -\bH\bI +\bH^2 = \bI-\bH,
$
which implies that $\bI-\bH$ is idempotent.
Thus, $\bI-\bH$ qualifies as a projection matrix. 
Applying the spectral theorem again, let $\bH =\bQ\bLambda\bQ^\top$. 
Then, we can express $\bI-\bH$ as  $\bI-\bH = \bQ\bQ^\top - \bQ\bLambda\bQ^\top = \bQ(\bI-\bLambda)\bQ^\top$. Hence, the column space of $\bI-\bH$ is spanned by the eigenvectors of $\bH$ corresponding to the zero eigenvalues of $\bH$ (by Proposition~\ref{proposition:eigen-of-projection-matrix}), which coincides with $\mathcalV^\perp$.
\end{proof}

For a detailed examination of the derivation and implications of the projection matrix, readers are encouraged to consult Appendix~\ref{section:by-geometry-hat-matrix}, although it may not be the primary focus of matrix decomposition results.

\section{Application: Alternative Definition of PD and PSD}\label{section:equivalent-pd-psd}
In Definition~\ref{definition:psd-pd-defini}, we introduced the definitions of positive definite (PD) and positive semidefinite (PSD) matrices based on their quadratic forms. Here, we prove that a symmetric matrix is positive definite if and only if all of its eigenvalues are positive, i.e., Theorem~\ref{theorem:eigen_charac}.


\begin{proof}[of Theorem~\ref{theorem:eigen_charac}] 
Suppose $\bA$ is PD. For any eigenvalue $\lambda$ and its corresponding eigenvector $\bv$ of $\bA$, we have $\bA\bv = \lambda\bv$. Consequently, 
$
\bv^\top \bA\bv = \lambda\normtwo{\bv}^2 > 0.
$
This implies $\lambda>0$.
Conversely, suppose all eigenvalues of $\bA$ are positive, and consider the spectral decomposition of $\bA =\bQ\bLambda \bQ^\top$. Let $\bx$ be any nonzero vector, and let $\by\triangleq\bQ^\top\bx$. We have
$$
\bx^\top \bA \bx = \bx^\top (\bQ\bLambda \bQ^\top) \bx = (\bx^\top \bQ) \bLambda (\bQ^\top\bx) = \by^\top\bLambda\by = \sum_{i=1}^{n} \lambda_i y_i^2>0.
$$
That is, $\bA$ is PD. 
Similarly, we can establish the proof for the second part of the claim regarding PSD matrices.
\end{proof}

The same technique can be used to prove properties of indefinite matrices.
\begin{exercise}
Show that a symmetric matrix $\bA\in\real^{n\times n}$ is indefinite if and only if it has at least one positive eigenvalue and at least one negative eigenvalue.
\end{exercise}

\subsection{Nonsingular Factor of PSD and PD Matrices}\label{section:nonsing_fac_pdpsd}
Theorem~\ref{theorem:eigen_charac} can be applied to prove the nonsingular factorization  of PD  matrices. See Figure~\ref{fig:pd-summary} for its role in the factorizations of PD matrices.
\begin{theoremHigh}[Nonsingular Factor of PSD and PD Matrices]\label{theorem:nonsingular-factor-of-PD}
Let $\bA\in\real^{n\times n}$ be a real symmetric matrix. Then, $\bA$ is PSD if and only if $\bA$ can be factored as $\bA=\bP^\top\bP$, where $\bP\in\real^{n\times n}$; and it is PD if and only if $\bP$ is nonsingular. 
The result can be applied to an $m\times n$  matrix $\bS$ with $m\geq n$ such that $\bA=\bS^\top\bS$ for both PD and PSD $\bA$ ($\bS$ has full column rank if $\bA$ is PD). In all cases, we have
\begin{itemize}
\item $\rank(\bA)=\rank(\bP)=\rank(\bS)$.
\item Null space: $\nspace(\bA)=\nspace(\bP)=\nspace(\bS)$.
\end{itemize}
\end{theoremHigh}
\begin{proof}[of Theorem~\ref{theorem:nonsingular-factor-of-PD}]
For the first part, suppose $\bA$ is PSD. Its spectral decomposition is given by $\bA = \bQ\bLambda\bQ^\top$. Since eigenvalues of PSD matrices are nonnegative, we can decompose $\bLambda$ as $\bLambda=\bLambda^{1/2}\bLambda^{1/2}$~\footnote{When $\bLambda$ is positive definite,  $\bLambda^{1/2}$ is called the \textit{positive definite square root}; while $\bLambda$ is positive semidefinite,  $\bLambda^{1/2}$ is called the \textit{positive semidefinite square root}. In such cases, $\bA^{1/2}=\bQ\bLambda^{1/2}\bQ^\top$ is called the \textit{square root} of $\bA$.}. Let $\bP \triangleq \bLambda^{1/2}\bQ^\top$. Then, $\bA$ can be factored  as $\bA=\bP^\top\bP$.

Conversely, if $\bA$ can be factored as $\bA=\bP^\top\bP$, then all eigenvalues of $\bA$ are nonnegative. To see this, for any eigenvalues $\lambda$ and its corresponding eigenvector $\bv$ of $\bA$, we have 
$$
\lambda = \frac{\bv^\top\bA\bv}{\bv^\top\bv} = \frac{\bv^\top\bP^\top\bP\bv}{\bv^\top\bv}=\frac{\normtwo{\bP\bv}^2}{\normtwo{\bv}^2} \geq 0.
$$
This implies $\bA$ is PSD according to Theorem~\ref{theorem:eigen_charac}.

Similarly, we can prove the second part of the argument for PD matrices, where the positive definiteness will result in the nonsingular $\bP$, and the nonsingularity of  $\bP$ implies  the positivity  of the eigenvalues. 

We can of course append zero rows into $\bP\rightarrow \bS\in\real^{m\times n}$ such that $\bA=\bS^\top\bS$. In this case, $\bS$ has full column rank $n$. Conversely, given a matrix $\bS\in\real^{m\times n}$ with full column rank $n$ (i.e., $m\geq n$) such that $\bA=\bS^\top\bS$, the null space of $\bS$ is of dimension zero. Therefore, $\bx^\top\bA\bx=\bx^\top\bS^\top\bS\bx > 0$ for all nonzero $\bx$.

To see the null space, if $\bA\bx=\bzero$, then $\bx^\top\bA\bx=\normtwo{\bP\bx}^2=\normtwo{\bS\bx}^2=0$ such that $\nspace(\bA)\subseteq\nspace(\bS)=\nspace(\bP)$. On the other hand, if $\bP\bx=\bS\bx=\bzero$, we also have $\bA\bx=\bzero$ such that $\nspace(\bS)=\nspace(\bP)\subseteq \nspace(\bA)$. This shows $\nspace(\bA)=\nspace(\bP)=\nspace(\bS)$. The fundamental theorem of linear algebra~\ref{theorem:fundamental-linear-algebra} implies $\rank(\bA)=\rank(\bP)=\rank(\bS)$.
\end{proof}

By the definition of congruence (Definition~\ref{definition:simiar_congru}), Theorem~\ref{theorem:nonsingular-factor-of-PD} implies that a matrix $\bA$ is PD if and only if it is congruent to the identity matrix.

\subsection{Proof for Semidefinite Rank-Revealing Decomposition}\label{section:semi-rank-reveal-proof}


In this section, we provide a proof for Theorem~\ref{theorem:semidefinite-factor-rank-reveal},  which establishes the existence of the rank-revealing decomposition for positive semidefinite matrices.
\index{Rank-revealing}\index{Semidefinite rank-revealing}
\begin{proof}[of Theorem~\ref{theorem:semidefinite-factor-rank-reveal}]
The proof is a consequence of the ``nonsingular" factor of PSD matrices (Theorem~\ref{theorem:nonsingular-factor-of-PD}) and the existence of the column-pivoted QR decomposition (Theorem~\ref{theorem:rank-revealing-qr-general}).
	
According to Theorem~\ref{theorem:nonsingular-factor-of-PD}, the ``nonsingular" factor of the PSD matrix $\bA$ is given by $\bA = \bZ^\top\bZ$, where $\bZ=\bLambda^{1/2}\bQ^\top$, and $\bA=\bQ\bLambda\bQ^\top$ is the spectral decomposition of $\bA$.
As per Proposition~\ref{proposition:rank-of-symmetric}, the rank of matrix $\bA$ is the number of nonzero eigenvalues (here, the number of positive eigenvalues since $\bA$ is PSD). Therefore, only $r$ components in $\bLambda^{1/2}$ are nonzero, and $\bZ=\bLambda^{1/2}\bQ^\top$ contains only $r$ independent columns, i.e., $\bZ$ is of rank $r$. By utilizing the column-pivoted QR decomposition, we have
$$
\bZ\bP = \bQ
\begin{bmatrix}
	\bR_{11} & \bR_{12} \\
	\bzero   & \bzero 
\end{bmatrix},
$$
where $\bP$ is a permutation matrix, $\bR_{11}\in \real^{r\times r}$ is upper triangular with positive diagonals, and $\bR_{12}\in \real^{r\times (n-r)}$. Therefore,
$$
\bP^\top\bA\bP  = 
\bP^\top\bZ^\top\bZ\bP = 
\begin{bmatrix}
\bR_{11}^\top & \bzero \\
\bR_{12}^\top & \bzero 
\end{bmatrix}
\begin{bmatrix}
	\bR_{11} & \bR_{12} \\
	\bzero   & \bzero 
\end{bmatrix},
\gap \text{with}\gap 
\bR \triangleq \begin{bmatrix}
	\bR_{11} & \bR_{12} \\
	\bzero   & \bzero 
\end{bmatrix}.
$$
Thus, we find the rank-revealing decomposition for the semidefinite matrix: $\bP^\top\bA\bP = \bR^\top\bR$.
\end{proof}

This decomposition is derived  using complete pivoting,  at each stage of which, we permute the largest diagonal element
in the active submatrix into the pivot position. The procedure is similar to the partial pivoting technique discussed in Section~\ref{section:partial-pivot-lu}.

\section{Application: Cholesky via  QR and  Spectral Decompositions}\label{section:cholesky-by-qr-spectral}
In this section, we present an alternative proof for the existence of the Cholesky decomposition using the QR and spectral decompositions.
%
\begin{proof}[of Theorem~\ref{theorem:cholesky-factor-exist}]
According to Theorem~\ref{theorem:nonsingular-factor-of-PD}, the PD matrix $\bA$ can be factored as $\bA=\bP^\top\bP$, where $\bP$ is a nonsingular matrix. Applying the QR decomposition to $\bP$, we get $\bP = \bQ\bR$, where $\bQ$ is orthogonal and $\bR$ is upper triangular. This implies
$$
\bA = \bP^\top\bP = \bR^\top\bQ^\top\bQ\bR = \bR^\top\bR,
$$
where we notice that the form $\bR^\top\bR$ resembles the Cholesky decomposition, except that we have not yet established that $\bR$ has positive diagonal entries. However, from Algorithm~\ref{alg:reduced-qr}, the existence of QR decomposition via the Gram-Schmidt process (Section~\ref{section:qr-gram-compute}), we realize that the diagonals of $\bR$ are nonnegative; and if $\bP$ is nonsingular, the diagonals of $\bR$ are indeed positive. 
\end{proof}
The proof for the above theorem is a direct result of the existence of both the QR decomposition and the spectral decomposition. Thus, the existence of the Cholesky decomposition can be demonstrated through the QR decomposition and the spectral decomposition in this context.

\index{Gram-Schmidt}
\section{Application: Unique Power Decomposition of PD Matrices}\label{section:unique-posere-pd}

\begin{theoremHigh}[Unique Power Decomposition of PD/PSD Matrices]\label{theorem:unique-factor-pd}
Any $n\times n$ PSD (resp. PD) matrix $\bA$ can be \textbf{uniquely} factored as a power of another PSD (resp. PD) matrix $\bB$ such that $\bA =\bB^k$ with $k=\{1,2,\ldots\}$, where $\rank(\bB)=\rank(\bA)$. 
\end{theoremHigh}
\begin{proof}[of Theorem~\ref{theorem:unique-factor-pd}]
We first show the existence of such a positive semidefinite matrix $\bB$ that satisfies $\bA = \bB^k$. 
\paragraph{Existence.} Since $\bA$ is symmetric and positive semidefinite, its spectral decomposition  is given by $\bA = \bQ\bLambda\bQ^\top$, where $\bQ$ is orthogonal and $\bLambda$ is diagonal containing the  eigenvalues of $\bA$. Since eigenvalues of PSD matrices are nonnegative by Theorem~\ref{theorem:eigen_charac}, the $k$-th square root of $\bLambda$ exists. We can define $\bB \triangleq\bQ\bLambda^{1/k}\bQ^\top$ such that $\bA = \bB^k$, where $\bB$ is apparently PSD.

\paragraph{Uniqueness.} 
Suppose the factorization is not unique. Then, there exist two positive definite matrices $\bB_1$ and $\bB_2$ such that
$$
\bA = \bB_1^k = \bB_2^k.
$$
Their spectral decompositions are given by 
$$
\bB_1 = \bQ_1 \bLambda_1\bQ_1^\top \qquad \text{and} \qquad \bB_2 = \bQ_2 \bLambda_2\bQ_2^\top.
$$
We notice that $\bLambda_1^k$ and $\bLambda_2^k$ contain the eigenvalues of $\bA$, and both eigenvalues of $\bB_1$ and $\bB_2$ contained in $\bLambda_1$ and $\bLambda_2$ are nonnegative (since both $\bB_1$ and $\bB_2$ are  PSD). Without loss of generality, we suppose $\bLambda_1=\bLambda_2=\bLambda^{1/k}$, and $\bLambda=\diag(\lambda_1,\lambda_2, \ldots, \lambda_n)$ such that $\lambda_1\geq \lambda_2 \geq \ldots \geq \lambda_n$. Utilizing the equation  $\bB_1^k = \bB_2^k$, we have
$$
\bQ_1 \bLambda \bQ_1^\top = \bQ_2 \bLambda \bQ_2^\top  \quad\implies\quad \bQ_2^\top\bQ_1 \bLambda = \bLambda \bQ_2^\top\bQ_1.
$$
Let $\bZ \triangleq \bQ_2^\top\bQ_1 $ (which is orthogonal), this implies that $\bLambda$ and $\bZ$ commute, and $\bZ$ must be a block-diagonal matrix whose partitioning conforms to the block structure of $\bLambda$ (Remark~\ref{remark:commute}). This results in $\bLambda^{1/k} = \bZ\bLambda^{1/k}\bZ^\top$ and
$$
\bB_2 = \bQ_2 \bLambda^{1/k}\bQ_2^\top = \bQ_2 \bQ_2^\top\bQ_1\bLambda^{1/k} \bQ_1^\top\bQ_2 \bQ_2^\top=\bB_1.
$$
Thus, the decomposition is unique. In a similar manner, we can establish the unique decomposition of a PD matrix $\bA =  \bB^k$, where $\bB$ is also PD. 
For a more detailed discussion, see \citet{koeber2006unique, horn2012matrix}, which provides an alternative proof using polynomials.
\end{proof}

\paragraph{Decompositions for PD matrices.}  In summary, for a PD matrix $\bA$, we can factor it into $\bA=\bR^\top\bR$, where $\bR$ is an upper triangular matrix with positive diagonals as shown in Theorem~\ref{theorem:cholesky-factor-exist} by Cholesky decomposition; 
$\bA = \bP^\top\bP$, where $\bP$ is nonsingular in Theorem~\ref{theorem:nonsingular-factor-of-PD};
and $\bA = \bB^k$ for $k=\{1,2,\ldots\}$, where $\bB$ is PD in Theorem~\ref{theorem:unique-factor-pd}.
An overview of these decompositions is provided in Figure~\ref{fig:pd-summary}.

\begin{problemset}
	
\item \label{problem:real_herm} \textbf{Hermitian properties.} Let $\bA\in\complex^{n\times n}$ be complex Hermitian (Definition~\ref{definition:complex_special}). Show that 
\begin{itemize}
\item The eigenvalues of $\bA$  are all real.
\item The main diagonals of $\bA$ are all real.
\item $\bx^*\bA\bx$ is real for all $\bx\in\complex^n$.
\item $\bS^*\bA\bS$ is complex Hermitian for all $\bS\in\complex^{n\times n}$.
\end{itemize}

\item \textbf{Hermitian observations.}  Let $\bA=\bB+i\bC\in\complex^{n\times n}$ be complex Hermitian, where $\bB,\bC\in\real^{n\times n}$. Show that 
\begin{itemize}
\item If $\bA$ is Hermitian, then $\bB$ is symmetric and $\bC$ is skew-symmetric (Definition~\ref{definition:speci_mat}).
\item $\bA+\bA^*$ is Hermitian; $\bA-\bA^*$ is skew-Hermitian (Definition~\ref{definition:complex_special}).
\item If $\bA$ is Hermitian, then $i\bA$ is skew-Hermitian; if $\bA$ is skew-Hermitian, then $i\bA$ is Hermitian. How about the reverse claims?
\end{itemize}

\index{Toeplitz decomposition}
\item \textbf{Toeplitz decomposition.} Let $\bA\in\complex^{n\times n}$. Show that $\bA$ can be \textbf{uniquely} factored as $\bA=\bB+i\bC$, where both $\bB\in\complex^{n\times n}$ and $\bC\in\complex^{n\times n}$ are Hermitian. \textit{Hint: Consider $\frac{1}{2}(\bA+\bA^*) + i [\frac{1}{2i}(\bA-\bA^*)]$.}

\item \textbf{Hermitian decomposition.} Suppose the complex Hermitian $\bA\in\complex^{n\times n}$ admits a spectral decomposition $\bA=\bU\bLambda\bU^*$, where $\bU$ is unitary and $\bLambda=\diag(\lambda_1, \lambda_2, \ldots, \lambda_n)$ is real diagonal with $\lambda_1\geq \lambda_2\geq\ldots\geq \lambda_n$ (Theorem~\ref{theorem:normal_Complex_spectral_theorem}).
Define $\lambda_i^+\triangleq\max\{\lambda_i, 0\}$ and $\lambda_i^-\triangleq\min\{\lambda_i, 0\}$. Let  $\bLambda_+ \triangleq \diag(\lambda_1^+, \lambda_2^+, \ldots, \lambda_n^+)$ and $\bLambda_- \triangleq \diag(\lambda_1^-, \lambda_2^-, \ldots, \lambda_n^-)$; and $\bA_+\triangleq\bU\bLambda_+\bU^*$ and $\bA_-\triangleq-\bU\bLambda_-\bU^*$. Therefore, $\bA=\bA_+-\bA_-$, and $\bA_+$ is called the PSD part of $\bA$. Show that
\begin{itemize}
\item $\bA_+$ and $\bA_-$ are PSD (i.e., $\bx^*\bA_+\bx\geq 0$ and $\bx^*\bA_-\bx\geq 0$ for all nonzero $\bx\in\complex^n$).
\item $\bA_+$ and $\bA_-$ commute.
\item $\rank(\bA)=\rank(\bA_+)+\rank(\bA_-)$.
\item $\bA_+\bA_- = \bA_-\bA_+=\bzero$.
\end{itemize}

\item Show that $\bA$ is complex Hermitian if and only if $\bA^2=\bA^*\bA$.
\item Let $\bA, \bS\in\complex^{n\times n}$ be given, and $\bA$ be Hermitian. Show that $\bS\bA\bS^*$ is Hermitian.
\item Let $\bA,\bB\in\complex^{n\times n}$ be Hermitian. Show that $\bA$ and $\bB$ commute if and only if $\bA\bB$ is Hermitian.

\item \label{problem:real_skewherm} \textbf{Skew-Hermitian.}  Show that the eigenvalues of a skew-Hermitian matrix $\bA\in\complex^{n\times n}$ (Definition~\ref{definition:complex_special})  are either  purely imaginary or zero.

\item \textbf{Skew-Hermitian.} Show that $\bx^*\bA\bx$ is purely imaginary for all $\bx\in\complex^n$ if and only if $\bA$ is skew-Hermitian.

\index{Sylvester's law of inertia}
\item \label{problem:real_skewhe_final} \textbf{Inertia of Hermitian: Sylvester's law of inertia.} Let $\bA,\bB\in\complex^{n\times n}$ be Hermitian. Show that $\bA$ and $\bB$ are congruent (Definition~\ref{definition:simiar_congru}) if and only if $\inertia(\bA) = \inertia(\bB)$, where the inertia of a matrix $\inertia(\cdot)$ is a vector of three elements containing the numbers of negative eigenvalues, zero eigenvalues, and positive eigenvalues. \textit{Hint: Use Theorem~\ref{theorem:normal_Complex_spectral_theorem} for Hermitian matrices.}

\item Show that two normal matrices are similar if and only if they have the same characteristic polynomial.

\item  \label{problem:symm_square} \textbf{Square of symmetric matrices.} Let $\bA\in\real^{n\times n}$ be symmetric. Show that the eigenvalues of $\bA^2$ (can only be nonnegative) are the squares  of the eigenvalues of $\bA\in\real^{n\times n}$. What is the relationship between the eigenvectors?

\item \textbf{Symmetric idempotent.} Let $\bA=\bA^\top=\bA^2\in\real^{n\times n}$ with $\rank(\bA)=r$. Show that there exists an orthogonal matrix $\bQ$ such that $\bQ^\top\bA\bQ=\diag(\bI_r, \bzero)$.

\item \label{problem:skewsymm_square} \textbf{Square of skew-symmetric matrices.} We have shown in Lemma~\ref{lemma:real-eigenvalues-spectral-skew} that all eigenvalues of any skew-symmetric matrix $\bA\in\real^{n\times n}$ are either purely imaginary or zero and have shown in Remark~\ref{remark:conjug_pari} that the non-real eigenvalues appear in complex conjugate pairs. Show that the eigenvalue of $\bA^2$ (can only be real and non-positive) are the negatives of the squares of the magnitudes of the eigenvalues of $\bA$. \textit{Hint: Use Corollary~\ref{corollary:spectral_normal_skew}.}

\item \label{problem:skewsymm_square2}  \textbf{Square of skew-symmetric matrices.} Following Problem~\ref{problem:skewsymm_square}, discuss  the eigenvalues of $-\bA^2$.

\item \textbf{Skew-symmetric.} Let $\bA\in\real^{n\times n}$ be skew-symmetric. Show that $\trace(\bA)=0$. Additionally, if $\bB\in \real^{n\times n}$ be symmetric, show that $\trace(\bA\bB)=0$.

\item Show that the null space of $\bA$ is the same as that of $\bA^*$ if $\bA\in\complex^{n\times n}$ is normal. \textit{Hint: Use Exercise~\ref{exercise:left_right_normal}.}

\item Let $\bA$ and $\bB$ be symmetric positive definite matrices. Show that the product $\bA\bB$ may not retain  symmetry, but its eigenvalues remain positive. \textit{Hint: Take the product of $\bB\bx$ and $\bA\bB\bx=\lambda\bx$.}


\item Given a symmetric positive definite matrix $\bA\in\real^{n\times n}$ with eigenvalues arranged in descending order as  $\lambda_1\geq \lambda_2 \geq \ldots \geq \lambda_n$.
\begin{itemize}
\item What are the eigenvalues of $\lambda\bI-\bA$? 
\item Prove the positive semidefinteness of $\lambda_1\bI-\bA$.
\item Show that $\lambda_1\bx^\top\bx\geq \bx^\top\bA\bx$ for every $\bx\in\real^n$.
\item Determine  the maximum value of $\bx^\top\bA\bx/\bx^\top\bx$.
\end{itemize}

\item Given a matrix $\bA\in\real^{n\times n}$ satisfying $\bA^2-\bA=2\bI$, show that $\bA$ is diagonalizable.




\item Show that when $\bA$ and $\bB$ are positive semidefintie, then  the condition $\trace(\bA\bB) = 0$ is equivalent to $\bA\bB = \bzero$. \textit{Hint: The trace is invariant under cyclic permutations, and write out the trace using the elements of matrices from the spectral  decomposition.}

\index{Fan's inequality}
\index{Hardy-Littlewood-P\'olya inequality}
\item \label{prob:fans_ineq} \textbf{Fan's inequality \citep{fan1949theorem, borwein2006convex}.} 
Let $\bA,\bB$ be real symmetric, and let $\blambda^{\downarrow}(\bA)$ be the vector containing the eigenvalues of $\bA$ in nonincreasing order. Show that  $\trace(\bA\bB) \leq \blambda^{\downarrow}(\bA)^{\top}\blambda^{\downarrow}(\bB)$.
The equality holds if and only if $\bA$ and $\bB$ admit spectral decompositions $\bA=\bQ\blambda^{\downarrow}(\bA)\bQ^\top$ and $\bB=\bQ\blambda^{\downarrow}(\bB)\bQ^\top$ (called \textit{simultaneous ordered spectral decomposition}). \footnote{Fan's inequality is a refinement of the Cauchy-Schwarz inequality for symmetric matrices, see Equation~\eqref{equation:vector_form_cauchyschwarz}.}

\item \label{prob:hardy_ineq} \textbf{Hardy-Littlewood-P\'olya inequality  \citep{borwein2006convex}.} Let $[\bx]^{\downarrow}$ denote the vector with the same components of $\bx$ permuted into nonincreasing order. Show that  $\bx^\top\by \leq [\bx]^{\downarrow \top}[\by]^{\downarrow}$. \textit{Hint: Apply  Fan's inequality to diagonal matrices.}

\index{Mirsky's theorem}
\item \label{prob:pert_eig_uniinvar} \textbf{Mirsky's theorem.} Let $\bA,\bB\in\real^{n\times n}$ be symmetric. Let further $\blambda^{\downarrow}(\bA)$ and $\blambda^{\downarrow}(\bB)$ be the vectors containing the eigenvalues of $\bA$ and $\bB$, respectively, in nonincreasing order. Show that $\norm{\diag(\blambda^{\downarrow}(\bA)) - \diag(\blambda^{\downarrow}(\bB))} \leq \norm{\bA-\bB}$ if $\norm{\cdot}$ is unitarily/orthogonally invariant (Definition~\ref{definition:unitarily_invaria}).

\item \label{problem:eig_det_house} Let $\bH=\bI - 2  \frac{\bu\bu^\top}{\bu^\top\bu}\in\real^{n\times n}$ be  a Householder reflector. Show  that the  eigenvalue $\lambda_1=1$ has multiplicity  $n-1$; and the  eigenvalue $\lambda_2=-1$ has multiplicity  $1$.
This explains why the determinant of a Householder reflector is $-1$ (by Theorem~\ref{theorem:eigen_trace}).

\item \label{problem:inverse_block_upp} \textbf{Inverse of block upper triangular.} Given a matrix $\bA\in\real^{m\times n}$, show that the inverse of 
$\scriptsize
\begin{bmatrix}
\bI_m & \bA \\
\bzero & \bI_n 
\end{bmatrix}$ 
is 
$\scriptsize
\begin{bmatrix}
\bI_m & -\bA \\
\bzero & \bI_n 
\end{bmatrix}$.

\item \label{problem:eig_reverse} \textbf{Eigenvalues of reverse product.} Let $\bA\in\real^{m\times n}$ and $\bB\in\real^{n\times m}$ with $m\leq n$. Show that the $n$ eigenvalues $\bB\bA$  are the eigenvalues of $\bA\bB$ together with $n-m$ zeros. \textit{Hint: Show that
$\scriptsize
\begin{bmatrix}
\bA\bB&\bzero \\
\bB& \bzero 
\end{bmatrix}$
and 
$\scriptsize
\begin{bmatrix}
\bzero &\bzero \\
\bB& \bA\bB
\end{bmatrix}$
are similar using Problem~\ref{problem:inverse_block_upp}, and use Proposition~\ref{proposition:eigenvalue-similar-matrices} to discuss the eigenvalues of the two matrices.
}

\item \textbf{Eigenvalues of rank decomposition.} Consider the rank decomposition of $\bA=\bD\bF\in\real^{n\times n}$ with rank $r$ (Theorem~\ref{theorem:rank-decomposition}). Show that the eigenvalues of $\bA$ are the same as those of $\bF\bD$ together with $n-r$ zeros. 
How does this result change if the decomposition satisfies $\bD\in\real^{n\times k}$ and $\bF\in\real^{k\times n}$ with $k>r$?
\textit{Hint: Use Problem~\ref{problem:eig_reverse}.}

\index{Courant-Fischer theorem}
\item \textbf{Courant-Fischer.} Let $\bA \in \real^{n\times n}$ be symmetric. Show that if $\bx^\top \bA\bx \geq 0$ for all $\bx$ in a $k$-dimensional subspace, then $\bA$ has at least $k$ nonnegative eigenvalues; if $\bx^\top \bA\bx > 0$ for all nonzero $\bx$ in a $k$-dimensional subspace, then $\bA$ has at least $k$ positive eigenvalues. \textit{Hint: Use Proposition~\ref{proposition:cour_fish}.}

\index{Kantorovich's inequality}
\index{Wielandt's inequality}
\item \textbf{Kantorovich's and Wielandt's inequalities.} Let $\bA\in\real^{n\times n}$ be PD. Show that 
$$
\begin{aligned}
&\qquad\qquad\qquad \qquad (\bx^\top\bx)^2 \leq (\bx^\top\bA\bx)(\bx^\top\bA^{-1}\bx) \quad \text{for all $\bx\in\real^n$};\\
&\textbf{Kantorovich's: }(\bx^\top\bA\bx)(\bx^\top\bA^{-1}\bx) \leq \frac{(\lambda_1+\lambda_n)^2}{4\lambda_1\lambda_n}\normtwo{\bx}^4, \quad \text{for all $\bx\in\real^n$};\\
&\textbf{Wielandt's: }{(\bx^\top\bA\by)}^2 \leq \frac{(\lambda_1-\lambda_n)^2}{(\lambda_1+\lambda_n)^2} (\bx^\top\bA\bx)(\by^\top\bA\by), \quad \text{for all orthogonal $\bx,\by\in\real^n$},
\end{aligned}
$$
where $\lambda_1$ and $\lambda_n$ are the smallest and largest eigenvalues of $\bA$, respectively. 

\item \textbf{Subspace in symmetric.} Let $\bA\in\real^{n\times n}$ be symmetric. Show that $\cspace(\bA)=\cspace(\bA^k)$ and $\nspace(\bA)=\nspace(\bA^k)$ for $k=2,3,\ldots$.

\item Show that $\bA$ is symmetric $\iff$ 
$\scriptsize\begin{bmatrix}
\bzero & \bA\\
\bA& \bzero
\end{bmatrix}$ is symmetric 
$\iff$
$\scriptsize\begin{bmatrix}
\bzero & \bA\\
-\bA& \bzero
\end{bmatrix}$ 
is skew-symmetric.

\index{Decell's theorem}
\item \textbf{Decell's theorem \citep{decell1965application}.} Let the characteristic polynomial of $\bA\bA^\top$ ($\bA\in\real^{m\times n}$) be $\det(\bA\bA^\top-\lambda\bI)=  (-1)^m (a_0\lambda^m+a_1\lambda^{m-1}+\ldots+ a_{m-1}\lambda +a_m)$, where $a_0=1$. Let $k$ be the largest integer such that $a_k\neq 0$. Show that the generalized inverse $\bA^-$ (satisfying $\bA\bA^-\bA=\bA$; see Section~\ref{section:g_inv}) is 
$$
\bA^-=
-a_k^{-1}\bA^\top\left( (\bA\bA^\top)^{k-1}+ a_1(\bA\bA^\top)^{k-2}+\ldots+a_{k-2}(\bA\bA^\top) + a_{k-1} \bI\right).
$$
When $k=0$,  $\bA^-=\bzero$.
\textit{Hint: Use the Cayley-Hamilton theorem and the spectral theorem.}
\end{problemset}

%% file: chapter-svd.tex
\newpage
\chapter{Singular Value Decomposition (SVD)}\label{chapter:SVD}
\begingroup
\hypersetup{
	linkcolor=structurecolor,
	linktoc=page,  
}
\minitoc \newpage
\endgroup

\section{Singular Value Decomposition}
\lettrine{\color{caligraphcolor}I}
In eigenvalue decomposition, we factor the matrix into a diagonal matrix. However, this is not always possible. 
When $\bA$ lacks linearly independent eigenvectors, such diagonalization cannot be achieved.  
The \textit{singular value decomposition (SVD)} fills this gap. 
Instead of using an eigenvector matrix, SVD decomposes a matrix into two orthogonal matrices. Below, we state the theorem for SVD and discuss its existence in the following sections.

\index{Karhunen-Loewe expansion}
\index{Decomposition: SVD}
\begin{theoremHigh}[Reduced SVD for Rectangular Matrices \citep{eckart1939principal, klema1980singular}]\label{theorem:reduced_svd_rectangular}
Let $\bA\in\real^{m\times n}$ be any real $m\times n$ matrix with rank $r$. Then, it can be factored as
$$
\bA = \bU \bSigma \bV^\top,
$$ 
where $\bSigma\in \real^{r\times r}$ is a real diagonal matrix, $\bSigma=\diag(\sigma_1, \sigma_2 \ldots, \sigma_r)$ with $\sigma_1 \geq \sigma_2 \geq \ldots \geq \sigma_r>0$ and 
\footnote{Note when $\bA$ is complex,  both $\bU$ and $\bV$ become semi-unitary. However, $\bSigma$ remains real and diagonal because the eigenvalues of Hermitian matrices $\bA\bA^*$ or $\bA^*\bA$ are real.}
\begin{itemize}
\item The $\sigma_i$'s are the nonzero \textit{singular values} of $\bA$, in the meantime, which are the positive square roots of the nonzero \textit{eigenvalues} of $\trans{\bA} \bA$ and $ \bA \trans{\bA}$.

\item The columns of $\bU\in \real^{m\times r}$ contain the $r$ eigenvectors of $\bA\bA^\top$ corresponding to the $r$ nonzero eigenvalues of $\bA\bA^\top$ (semi-orthogonal). 

\item The columns of $\bV\in \real^{n\times r}$ contain the $r$ eigenvectors of $\bA^\top\bA$ corresponding to the $r$ nonzero eigenvalues of $\bA^\top\bA$ (semi-orthogonal). 


\item The columns of $\bU$ and $\bV$ are called the \textit{left and right singular vectors} of $\bA$, respectively. 

\item Furthermore, the columns of both $\bU$ and $\bV$ are mutually orthonormal (by spectral theorem~\ref{theorem:spectral_theorem}). 
\end{itemize}

In particular, the matrix decomposition can be expressed as the sum of outer products of vectors:
$$
\bA = \bU \bSigma \bV^\top = \sum_{i=1}^r \sigma_i \bu_i \bv_i^\top,
$$ which is a sum of $r$ rank-one matrices.

\end{theoremHigh}

By appending $m-r$ additional orthonormal columns to the  $r$  eigenvectors of   $\bA\bA^\top$,  similar to the silent columns in the QR decomposition (Section~\ref{section:silentcolu_qrdecomp}), we obtain an orthogonal matrix $\bU\in \real^{m\times m}$. 
The same process applies to the columns of $\bV$.
We then present the full SVD for matrices in the following theorem. 
Differences between the reduced and full SVD are highlighted in \textcolor{mylightbluetext}{blue} text.
\begin{theoremHigh}[Full SVD for Rectangular Matrices]\label{theorem:full_svd_rectangular}
Let $\bA\in\real^{m\times n}$ be any real $m\times n$ matrix with rank $r$. Then, it can be factored as
$$
\bA = \bU \bSigma \bV^\top,
$$ 
where the  top-left part of $\bSigma\in $\textcolor{mylightbluetext}{$\real^{m\times n}$} is a real diagonal matrix, that is $\bSigma=\footnotesize\begin{bmatrix}
	\bSigma_1 & \bzero \\
	\bzero & \bzero
\end{bmatrix}$, where $\bSigma_1=\diag(\sigma_1, \sigma_2 \ldots, \sigma_r)\in \real^{r\times r}$ with $\sigma_1 \geq \sigma_2 \geq \ldots \geq \sigma_r>0$ and 
\begin{itemize}
\item The $\sigma_i$'s are the nonzero \textit{singular values} of matrix $\bA$, which are the positive square roots of the nonzero \textit{eigenvalues} of $\trans{\bA} \bA$ and $ \bA \trans{\bA}$. 

\item  The columns of $\bU\in \textcolor{mylightbluetext}{\real^{m\times m}}$ contain the $r$ eigenvectors of $\bA\bA^\top$ corresponding to the $r$ nonzero eigenvalues of $\bA\bA^\top$ \textcolor{mylightbluetext}{and $m-r$ extra orthonormal vectors from $\nspace(\bA^\top)$}. 

\item  The columns of $\bV\in \textcolor{mylightbluetext}{\real^{n\times n}}$ contain the $r$ eigenvectors of $\bA^\top\bA$ corresponding to the $r$ nonzero eigenvalues of $\bA^\top\bA$ \textcolor{mylightbluetext}{and $n-r$ extra orthonormal vectors from $\nspace(\bA)$}.

\item  The columns of  $\bU$ and $\bV$ are called the \textit{left and right singular vectors} of $\bA$, respectively. 

\item Furthermore, the columns of $\bU$ and $\bV$ are mutually orthonormal (by spectral theorem~\ref{theorem:spectral_theorem}), and \textcolor{mylightbluetext}{$\bU$ and $\bV$ are orthogonal matrices}. 
\end{itemize}

In particular, the matrix decomposition can be expressed as the sum of outer products of vectors: $ \bA = \bU \bSigma \bV^\top = \sum_{i=1}^r \sigma_i \bu_i \bv_i^\top$, which is a sum of $r$ rank-one matrices.
\end{theoremHigh}

\begin{figure}[H]
\centering  
\vspace{-0.35cm} 
\subfigtopskip=2pt 
\subfigbottomskip=2pt 
\subfigcapskip=-5pt 
\subfigure[Reduced SVD decomposition.]{\label{fig:svdhalf}
\includegraphics[width=0.47\linewidth]{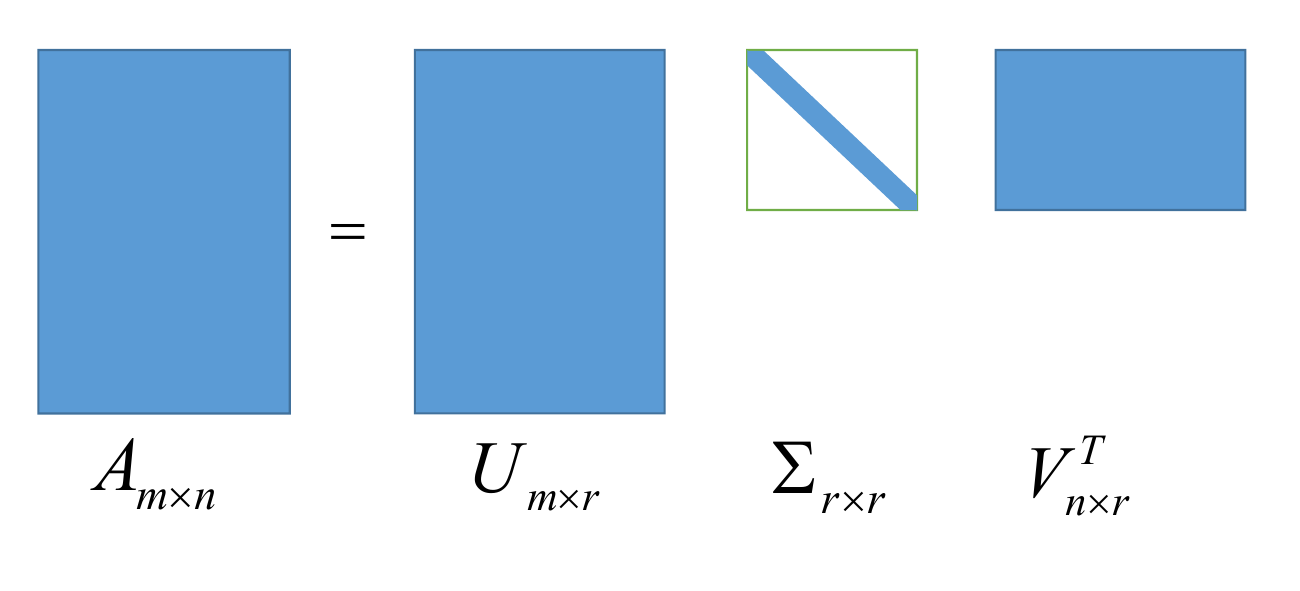}}
\quad 
\subfigure[Full SVD decomposition.]{\label{fig:svdall}
\includegraphics[width=0.47\linewidth]{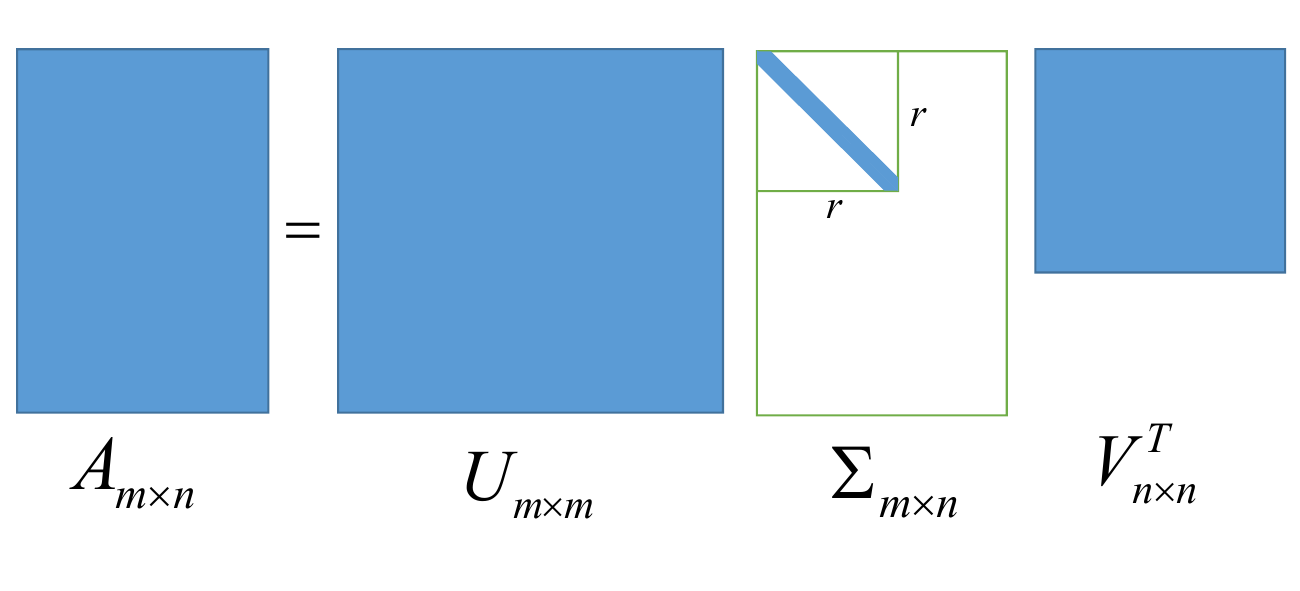}}
\caption{Comparison between the reduced and full SVD.}
\label{fig:svd-comparison}
\end{figure}

In image processing, the SVD is also known as the \textit{Karhunen-Loewe expansion}.
The comparison between the reduced and the full SVD is shown in Figure~\ref{fig:svd-comparison}, where white entries are zero, and \textcolor{mylightbluetext}{blue} entries are not necessarily zero.

\begin{remark}
We clarify some terminologies and observations in the following.
\begin{itemize}
\item The \textit{multiplicity} of a singular value $\sigma_i$ of $\bA$ is the multiplicity of $\sigma_i^2$ as an eigenvalue of $\bA^\top\bA$ or, equivalently, of $\bA\bA^\top$. Therefore, a singular value $\sigma_i$ of $\bA$ is said to be \textbf{simple} if $\sigma_i^2$ is a simple eigenvalue (Definition~\ref{definition:simple_eig}) of $\bA^\top\bA$ or, equivalently, of $\bA\bA^\top$.

\item Complex Hermitian matrices has only real eigenvalues (Problem~\ref{problem:real_herm}). Therefore, for a complex matrix $\bA\in\complex^{m\times n}$, the singular matrix $\bSigma$ remains real, while the orthogonal matrices $\bU$ and $\bV$ become unitary (resp. semi-unitary) for the full SVD (resp. reduced SVD).

\item We have shown that product and sum of eigenvalues are equal to the determinant and trace of the given matrix, respectively (Theorem~\ref{theorem:eigen_trace}, \ref{theorem:eigen_trace2}). Therefore, the product of singular values of $\bA$ is equal to $\abs{\det(\bA)}$; and the sum of squared singular values is equal to $\trace(\bA^\top\bA)=\trace(\bA\bA^\top)$.

\item The SVD of a matrix is not unique. For example, replacing $\bU$ by $-\bU$ and $\bV$ by $-\bV$ yields another valid SVD.

\item When $\bA$ is symmetric, then $\bU=\bV$. The singular values are equal to the eigenvalues of $\bA$ (Theorem~\ref{theorem:spectral_theorem}, Corollary~\ref{corollary:spectral_normal_symmetricx}, Problem~\ref{problem:symm_square}).

\item When $\bA$ is skew-symmetric, then $\bU=\bV$. The singular values are equal to the square roots of the eigenvalues of $-\bA^2$ (Corollary~\ref{corollary:spectral_normal_skew}, Problem~\ref{problem:skewsymm_square}, Problem~\ref{problem:skewsymm_square2}).
\end{itemize}
\end{remark}

\section{Existence of SVD via Spectral Decomp. of Symmetric Matrices}
To establish the existence of the SVD, we  rely on the  following lemmas. 
As previously mentioned, the singular values correspond to the square roots of the eigenvalues of $\bA^\top\bA$. 
Since negative values lack real square roots, it follows that the eigenvalues must be nonnegative.
\begin{lemma}[Nonnegative Eigenvalues of $\bA^\top \bA$]\label{lemma:nonneg-eigen-ata}
For any matrix $\bA\in \real^{m\times n}$, the  $\bA^\top \bA$ possesses nonnegative eigenvalues.
\end{lemma}
\begin{proof}[of Lemma~\ref{lemma:nonneg-eigen-ata}]
For an eigenvalue $\lambda$ and its corresponding eigenvector $\bx$ of $\bA^\top \bA$, we have
$
\bA^\top \bA \bx = \lambda \bx \implies \bx^\top \bA^\top \bA \bx = \lambda \bx^\top\bx. 
$
Since $\bx^\top \bA^\top \bA \bx  = \normtwo{\bA \bx}^2 \geq 0$ and $\bx^\top\bx> 0$, it follows that $\lambda \geq 0$.
\end{proof}
It can also be shown that $\bA^\top \bA$ is positive semidefinite. Thus, the eigenvalue characterization theorem also guarantees the nonnegative eigenvalues.
Given that $\bA^\top\bA$ has nonnegative eigenvalues, we  can define the \textit{singular value} $\sigma\geq 0$ of $\bA$ such that $\sigma^2$ is an eigenvalue of $\bA^\top\bA$, i.e., {$\bA^\top\bA \bv = \sigma^2 \bv$}. This definition is crucial for establishing the existence of the SVD.

We have shown in Lemma~\ref{lemma:rankAB} that $\rank$($\bA\bB$)$\leq$min$\{\rank$($\bA$), $\rank$($\bB$)\}.
However, the symmetric matrix $\bA^\top \bA$ is  special in that the rank of $\bA^\top \bA$ is equal to $\rank(\bA)$. We now prove this.
\begin{lemma}[Rank of $\bA^\top \bA$]\label{lemma:rank-of-ata}
The matrices $\bA^\top \bA$ and $\bA$ share the same rank.
Extending this observation to $\bA^\top$, we can also prove that $\bA\bA^\top$ and $\bA$ share the same rank. 
\end{lemma}
\begin{proof}[of Lemma~\ref{lemma:rank-of-ata}]
Let $\bx\in \nspace(\bA)$. Then,
$
\bA\bx  = \bzero \implies \bA^\top\bA \bx =\bzero, 
$
i.e., $\bx\in \nspace(\bA) \implies \bx \in \nspace(\bA^\top \bA)$. Therefore, $\nspace(\bA) \subseteq \nspace(\bA^\top\bA)$. 

Conversely, let $\bx \in \nspace(\bA^\top\bA)$. Then,
$
\bA^\top \bA\bx = \bzero\implies \bx^\top \bA^\top \bA\bx = 0\implies \normtwo{\bA\bx}^2 = 0 \implies \bA\bx=\bzero, 
$
i.e., $\bx\in \nspace(\bA^\top \bA) \implies \bx\in \nspace(\bA)$. Therefore, $\nspace(\bA^\top\bA) \subseteq\nspace(\bA) $. 
As a result, by ``sandwiching," it follows that  
$$\nspace(\bA) = \nspace(\bA^\top\bA) \qquad
\text{and} \qquad 
\dim(\nspace(\bA)) = \dim(\nspace(\bA^\top\bA)).
$$
According to the fundamental theorem of linear algebra in Appendix~\ref{appendix:fundamental-rank-nullity}, $\bA^\top \bA$ and $\bA$ have the same rank.
\end{proof}

Applying the same argument to $\bA^\top$, we can also prove that $\bA\bA^\top$ and $\bA$ have the same rank:
$$
\rank(\bA) = \rank(\bA^\top \bA) = \rank(\bA\bA^\top).
$$ 
\index{Fundamental theorem}


In the context of the SVD, we claim that the matrix $\bA$ is a sum of $r$ rank-one matrices, where $r$ is the number of nonzero singular values. And the number of nonzero singular values actually aligns with the rank of the matrix.

\begin{lemma}[Number of Nonzero Singular Values]\label{lemma:rank-equal-singular}
The number of nonzero singular values of a matrix $\bA$ is equal to its rank.
\end{lemma}
\begin{proof}[of Lemma~\ref{lemma:rank-equal-singular}]
The rank of any symmetric matrix (here $\bA^\top\bA$) equals the number of nonzero eigenvalues (with repetitions) by Proposition~\ref{proposition:rank-of-symmetric}. Therefore, the number of nonzero singular values equals the rank of $\bA^\top \bA$. By Lemma~\ref{lemma:rank-of-ata}, the number of nonzero singular values equals the rank of $\bA$.
\end{proof}

We are now prepared to demonstrate the existence of the SVD.
\begin{proof}[{of Theorem~\ref{theorem:reduced_svd_rectangular}: Existence of  Reduced SVD}]
Since $\bA^\top \bA$ is a symmetric matrix, according to the spectral theorem~\ref{theorem:spectral_theorem} and Lemma~\ref{lemma:nonneg-eigen-ata}, there exists a semi-orthogonal matrix $\bV\in\real^{n\times r}$ such that
$
{\bA^\top \bA = \bV \bSigma^2 \bV^\top},
$
where $\bV=[\bv_1,\bv_2,\ldots,\bv_r]$,
and $\bSigma$ is a diagonal matrix containing the nonzero singular values of $\bA$, i.e., $\bSigma^2$ contains the nonzero eigenvalues of $\bA^\top \bA$.
Specifically, $\bSigma=\diag(\sigma_1, \sigma_2, \ldots, \sigma_r)$, where $\{\sigma_1^2, \sigma_2^2, \ldots, \sigma_r^2\}$ are the nonzero eigenvalues of $\bA^\top \bA$, and $r$ is the rank of $\bA$. 
Now, we delve into the central part of the proof.
\begin{mdframed}[hidealllines=\mdframehidelineNote,backgroundcolor=\mdframecolor]
Starting from \fbox{$\bA^\top\bA \bv_i = \sigma_i^2 \bv_i$}, $\forall i \in \{1, 2, \ldots, r\}$, i.e., the eigenvector $\bv_i$ of $\bA^\top\bA$ is corresponding to the eigenvalue $\sigma_i^2$:

1. Multiplying both sides by $\bv_i^\top$:
$$
\bv_i^\top\bA^\top\bA \bv_i = \sigma_i^2 \bv_i^\top \bv_i 
\quad\implies\quad
\normtwo{\bA\bv_i}^2 = \sigma_i^2 
\quad\implies\quad
\normtwo{\bA\bv_i}=\sigma_i.
$$

2. Multiplying both sides by $\bA$:
$$
\bA\bA^\top\bA \bv_i = \sigma_i^2 \bA \bv_i 
\quad\implies\quad
\bA\bA^\top \frac{\bA \bv_i }{\sigma_i}= \sigma_i^2 \frac{\bA \bv_i }{\sigma_i} 
\quad\implies\quad
\bA\bA^\top \bu_i = \sigma_i^2 \bu_i,
$$
where we notice that this form allows us to identify the eigenvector of $\bA\bA^\top$ corresponding to $\sigma_i^2$ as $\bA \bv_i$. Since the length of $\bA \bv_i$ is $\sigma_i$, we can introduce the definition of $\bu_i \triangleq \frac{\bA \bv_i }{\sigma_i}$ with a unit norm.
\end{mdframed}
These $\bu_i$'s are orthogonal because $(\bA\bv_i)^\top(\bA\bv_j)=\bv_i^\top\bA^\top\bA\bv_j=\sigma_j^2 \bv_i^\top\bv_j=0$ if $i\neq j$.
Therefore, the spectral decomposition of $\bA\bA^\top$ can be obtained by
$
{\bA \bA^\top = \bU \bSigma^2 \bU^\top},
$
where $\bU=[\bu_1, \bu_2, \ldots,\bu_r]$.
Since {$\bA\bv_i = \sigma_i\bu_i$}, we have 
\begin{equation}\label{equation:svd_eqprof1}
[\bA\bv_1, \bA\bv_2, \ldots, \bA\bv_r] = [ \sigma_1\bu_1,  \sigma_2\bu_2, \ldots,  \sigma_r\bu_r]\quad\implies\quad
\bA\bV = \bU\bSigma.
\end{equation}
Since $\bV\bV^\top \neq \bI$, we cannot obtain the reduced SVD directly. 
Suppose we append the semi-orthogonal matrix  $\bV$ into an orthogonal matrix $\widetilde{\bV}=[\bV, \bV_2]$, and append the semi-orthogonal matrix $\bU$ into an orthogonal matrix $\widetilde{\bU}=[\bU, \bU_2]$. We then obtain
$$
\bA \widetilde{\bV}= \widetilde{\bU} \widetilde{\bSigma} ,\gap \text{where}\gap 
\widetilde{\bSigma} 
= \begin{bmatrix}
	\bSigma & \bzero  \\
	\bzero &\bzero
\end{bmatrix}
\quad\implies\quad
\bA = \widetilde{\bU} \widetilde{\bSigma}\widetilde{\bV}^\top,
$$
i.e., the full SVD (since $\widetilde{\bV}\widetilde{\bV}^\top=\bI$).
Simplifying the product, we get:
$
\bA =\bU \bSigma \bV^\top + \bU_2 \cdot \bzero \cdot \bV_2^\top = \bU \bSigma \bV^\top,
$
i.e., the reduced SVD,
which completes the proof.
\end{proof}

\index{Orthogonal projection}
\index{Projection matrix}
\index{Projector}
The proof also demonstrates that if $\bA=\bU\bSigma\bV^\top$ is the reduced SVD of $\bA$, it follows from \eqref{equation:svd_eqprof1} that $\bA\bV\bV^\top=\bA$. This means that $\bV\bV^\top$ is an (orthogonal) projection matrix that maps every row of $\bA$ onto itself (a projection matrix onto the row space of $\bA$; Definition~\ref{definition:projection_matrix_intro}).
For more information on alternative projectors derived from the SVD, see Section~\ref{section:svd_ortho_proj}.

An additional outcome of the aforementioned proof is that the spectral decomposition of $\bA^\top\bA = \bV \bSigma^2 \bV^\top$ will result in the spectral decomposition of $\bA \bA^\top = \bU \bSigma^2 \bU^\top$ with identical eigenvalues.

\begin{corollary}[Eigenvalues of $\bA^\top\bA$ and $\bA\bA^\top$]
The nonzero eigenvalues of $\bA^\top\bA$ and $\bA\bA^\top$ are the same.
\end{corollary}
We have shown in Lemma~\ref{lemma:nonneg-eigen-ata}  that the eigenvalues of $\bA^\top \bA$ are nonnegative. Consequently, the eigenvalues of $\bA\bA^\top$ are nonnegative as well.
\begin{corollary}[Nonnegative Eigenvalues of $\bA^\top\bA$ and $\bA\bA^\top$]
The eigenvalues of both  $\bA^\top\bA$ and $\bA\bA^\top$ are nonnegative.
\end{corollary}

\begin{exercise}[Eigenvalues and Singular Values]
Let $\bA\in\real^{m\times n}$ with $\rank(\bA)=r$. Show that
\begin{itemize}
\item $\lambda_i(\bA\bA^\top)=\lambda_i(\bA^\top\bA)=\sigma_i(\bA\bA^\top)=\sigma_i(\bA^\top\bA)=\sigma_i(\bA)^2$ for $i\in\{1,2,\ldots,r\}$.
\item $\sigma_i(\bA\bA^\top\bA)=\sigma_i(\bA)^3$ for $i\in\{1,2,\ldots,r\}$.
\item $\sigma_{\max}(\bA)\leq 1$ if and only if $\bI-\bA^\top\bA$ is PSD.
\end{itemize}
\end{exercise}

\section{Existence of  SVD via Spectral Decomp. of Normal Matrices}
The provided proof of the SVD's existence relies on the spectral decomposition of the symmetric matrices $\bA^\top\bA$ or $\bA\bA^\top$. 
Employing this approach facilitates the verification of various SVD properties, including the existence of four orthonormal bases in the SVD (see Section~\ref{section:four-space-svd}). 
However,
if we restrict our discussion to  \textbf{square normal matrices}, an alternative way to prove the existence of the SVD emerges through the spectral decomposition of normal matrices (Theorem~\ref{theorem:normal_real_spectral_theorem}).

As a recap, in Theorem~\ref{theorem:normal_real_spectral_theorem}, we claimed that a real normal matrix $\bA\in\real^{n\times n}$ (with $\bA^\top\bA=\bA\bA^\top$) admits the following spectral decomposition 
$$
\bA = \bQ\bD\bQ^\top,
$$
where $\bQ$ is an $n\times n$ orthogonal matrix, and $\bD$ is a block-diagonal matrix $\bD=\diag(\bA_1, \bA_2, $ $\ldots, \bA_p)$ with $p\leq n$. Each block within $\bD$ can either be one-dimensional, such as $\bA_i = \lambda_i$, or two-dimensional blocks of the following form 
$$
\bA_i =
\begin{bmatrix}
	\lambda_i & \mu_i \\
	-\mu_i & \lambda_i
\end{bmatrix},
\gap 
\text{with\,\, } \mu_i\neq 0.
$$

\begin{proof}[{of Theorem~\ref{theorem:reduced_svd_rectangular}: Existence of  SVD for Square Normal Matrices}]
Without loss of generality,
we assume $\bA_1, \bA_2, \ldots,  \bA_k$ are two-dimensional blocks and $\lambda_{2k+1}, \lambda_{2k+2}, \ldots, \lambda_n$ are the scalar entries. 
From the proof of Theorem~\ref{theorem:normal_real_spectral_theorem}, we know that the $\lambda_i\pm i \mu_i$ ($i\in \{1,2,\ldots, k\}$) are the complex eigenvalues of $\bA$, and $\lambda_{2k+1}, \lambda_{2k+1}, \ldots, \lambda_n$ are the real eigenvalues of $\bA$.
Let 
$$
\begin{aligned}
\rho_{2i-1}&\triangleq\rho_{2i} \triangleq \sqrt{\lambda_i^2+\mu_i^2} = \sqrt{\det(\bA_i)},\gap 
&\forall& i\in\{1,2,\ldots, k\}; \\
\rho_{i} &\triangleq |\lambda_i|, \gap &\forall& i \in\{2k+1, 2k+2, \ldots, \textcolor{mylightbluetext}{r}\},
\end{aligned}
$$
where we will soon see that the value of $r$ is the rank of the matrix $\bA$. Additionally, for $r<i\leq n$, it holds that $\lambda_i=0$.
Without loss of generality, we assume $\rho_1\geq \rho_2 \geq \ldots \geq \rho_r >0$ (this can be achieved via a permutation matrix). Then we have 
$$
\bA\bA^\top = \bA^\top\bA = \bQ\bD \bQ^\top \bQ\bD^\top\bQ^\top = \bQ\bD\bD^\top\bQ^\top,
$$
with 
$$
\bD\bD^\top = \diag\left(\rho_1^2,\rho_2^2, \ldots, \rho_r^2, 0, \ldots, 0   \right)\in\real^{n\times n}.
$$
Therefore, $\bQ\bD\bD^\top\bQ^\top$ is the spectral decomposition of matrix $\bA\bA^\top$, and $\{\rho_1^2, \ldots, \rho_r^2\}$ are the eigenvalues of $\bA\bA^\top$. This confirms that $r$ corresponds to the rank of the matrix  $\bA$.
Thus, let $\sigma_i\triangleq\rho_i$ for all $i\in\{1,2,\ldots, r\}$. Then, $\{\sigma_1, \sigma_2, \ldots, \sigma_r\}$ will be the (nonzero) singular values of $\bA$.  
We further define the block-diagonal matrix $\bZ$ (which is also an orthogonal matrix):
$$
\bZ \triangleq
\footnotesize
\begin{bmatrix}
	\frac{1}{\sigma_1} \bA_1 & &&&&& \\
	&\ddots &&&&& \\
	&& \frac{1}{\sigma_{2k}}\bA_k&&&& \\
	&&&  \frac{\lambda_{2k+1}}{|\lambda_{2k+1}|}&&& \\
	&&&&  \ddots && \\
	&&&&& \frac{\lambda_r}{|\lambda_r|}& \\
	&&&&&& \bI_{n-r} 
\end{bmatrix}
\normalsize
\in \real^{n\times n}, 
\gap 
\text{with }\bZ^\top\bZ = \bZ\bZ^\top = \bI_n.
$$
We also define the diagonal matrix $\bSigma \triangleq \diag(\sigma_1, \sigma_2, \ldots, \sigma_r, 0, \ldots, 0)\in\real^{n\times n}$.
We have 
$$
\bD = \bZ\bSigma.
$$
Let $\bQ\bZ \triangleq \bU$ and $ \bQ\triangleq\bV$. It  follows that
$$
\bA = \bQ\bD\bQ^\top = \bQ\bZ\bSigma\bQ^\top =\bU\bSigma\bV^\top.
$$
That is, $\bA$ admits the SVD $\bA=\bU\bSigma\bV^\top$.
\end{proof}

This proof will be proved valuable in finding the pseudo-inverse of normal matrices (Appendix~\ref{appendix:pseudo_inv_normal}).

\index{Gram–Schmidt}
\section{Existence of  SVD via Inductive Optimization}
To provide an alternative proof of the SVD,
without loss of generality, we assume $\bA\in\real^{m\times n}$ with $m\geq n$. Then $\bA$ admits the SVD $\bA=\bU\footnotesize\begin{bmatrix}
\bSigma \\ \bzero 
\end{bmatrix}\normalsize\bV^\top$, where $\bSigma\in\real^{n\times n}$ contains the singular values along its diagonal. If we could prove the existence of the SVD for this matrix $\bA$ with $m\geq n$, we can prove the existence for general matrices by applying the theorem to $\bA^\top$. The following proof is adapted from \citet{elden2007matrix, golub2013matrix}.

\begin{proof}[{of Theorem~\ref{theorem:reduced_svd_rectangular}: Existence of the SVD}]
Without loss of generality,  assume $m\geq n$. Consider the following optimization problem:
$$
\mathop{\sup}_{\normtwo{\bx}=1}\normtwo{\bA\bx},
$$
where the maximum is attained since we are seeking the supremum of a continuous function over a closed set. Once we have determined the optimal $\bx$, we set $\bA\bx = \sigma_1\by$, with the constraint  $\normtwo{\by}=1$. By definition, we have $\normtwo{\bA} = \sigma_1$. Employing the Gram-Schmidt process, we can extend $\bx$ and $\by$ to form orthogonal matrices:
$$
\bY_1 =
\begin{bmatrix}
	\by & \overline{\bY}_1
\end{bmatrix}
\in\real^{m\times m},
\gap 
\bX_1 =
\begin{bmatrix}
	\bx & \overline{\bX}_1
\end{bmatrix}
\in\real^{n\times n}.
$$
Then we have 
$$
\bY_1^\top\bA\bX_1=
\begin{bmatrix}
\sigma_1 & \by^\top \bA\overline{\bX}_1 \\
\overline{\bY}_1^\top \bA\bx & \overline{\bY}_1^\top \bA\overline{\bX}_1
\end{bmatrix}
=
\begin{bmatrix}
	\sigma_1 & \bw^\top  \\
\bzero & \overline{\bY}_1^\top \bA\overline{\bX}_1
\end{bmatrix}\triangleq
\bA_1,
\gap (\bw^\top\triangleq\by^\top \bA\overline{\bX}_1)
$$
where the penultimate equality arises  from the fact that $\bA\bx=\sigma_1\by$.
Therefore, by defining $\bB\triangleq\overline{\bY}_1^\top \bA\overline{\bX}_1$, we obtain
$$
\frac{1}{\sigma_1^2+\bw^\top\bw} \normtwo{\bA_1
\begin{bmatrix}
	\sigma_1 \\ \bw
\end{bmatrix}
}^2
=
\frac{1}{\sigma_1^2+\bw^\top\bw} \normtwo{
	\begin{bmatrix}
		\sigma_1^2 + \bw^\top\bw\\
		\bB\bw
	\end{bmatrix}
}^2
\geq \sigma_1^2+\bw^\top\bw.
$$
Since $\normtwo{\bA}^2=\normtwo{\bA_1}^2= \sigma_1^2$ \footnote{The matrix norm  used here is the spectral norm (Definition~\ref{definition:spectral_norm}).}, according to the subordinate property of matrix norm (Definition~\ref{definition:subordinate_matrix_norm}), we must have 
$$
\normtwo{\bA_1\begin{bmatrix}
		\sigma_1 \\ \bw
\end{bmatrix}}^2
\leq 
\normtwo{\bA_1}^2
\normtwo{\begin{bmatrix}
		\sigma_1 \\ \bw
\end{bmatrix}}^2
=
\sigma_1^2 \left( \sigma_1^2 +\bw^\top\bw \right).
$$
Combining the two inequalities, this implies $\bw^\top\bw=0$ and $\bw=0$.  Therefore, $\bY_1^\top\bA\bX_1$ finishes the first step of diagonalization for $\bA$:
$$
\bA
=
\bY_1
\begin{bmatrix}
\sigma_1 &\bzero \\
\bzero & \bB
\end{bmatrix}
\bX_1^\top.
$$
Suppose that any $(m-1)\times (n-1)$ matrix $\bB$ admits an SVD:
$$
\bB = 
\bU_2 
\begin{bmatrix}
\bSigma_2  \\
\bzero  
\end{bmatrix}
\bV_2^\top 
\,\, \implies \,\, 
\bA = 
\bY_1
\begin{bmatrix}
	\sigma_1 &\bzero \\
	\bzero & \bB
\end{bmatrix}
\bX_1^\top
=
\bY_1\begin{bmatrix}
1 & \bzero \\
\bzero & \bU_2 
\end{bmatrix}
\begin{bmatrix}
\sigma_1 & \bzero  \\
\bzero & \bSigma_2 \\
\bzero & \bzero 
\end{bmatrix}
\begin{bmatrix}
1 & \bzero  \\
\bzero & \bV_2^\top 
\end{bmatrix}
\bX_1^\top.
$$
Let 
$$
\bU \triangleq \bY_1\begin{bmatrix}
	1 & \bzero \\
	\bzero & \bU_2 
\end{bmatrix}, 
\gap 
\bSigma \triangleq 
\begin{bmatrix}
	\sigma_1 & \bzero  \\
	\bzero & \bSigma_2 \\
\end{bmatrix}, 
\gap 
\text{and }\quad
\bV\triangleq
\bX_1
\begin{bmatrix}
	1 & \bzero  \\
	\bzero & \bV_2
\end{bmatrix},
$$
we find the SVD of $\bA$ by induction.
\end{proof}

\section{Computing the SVD}\label{section:comput-svd-in-svd}
Suppose again we have an oracle algorithm capable of computing the eigenvalues and eigenvectors of $\bA^\top\bA$ at a cost of $f(m,n)$ flops. Then the computation of SVD is straightforward based on the previously outlined steps. 
The algorithm is shown in Algorithm~\ref{alg:svd-oracle}.

\begin{algorithm}[H] 
\caption{A Simple SVD} 
\label{alg:svd-oracle} 
\begin{algorithmic}[1] 
\Require 
Rank-$r$ matrix $\bA=[\ba_1, \ba_2, \ldots, \ba_n]$ with size $m\times n $; 
\State Initially get $\bA^\top\bA \bx_i = \sigma_i^2 \bx_i$ $\forall i\in \{1, 2, \ldots, r\}$; \Comment{ $f(m,n)$ flops}
\State Normalize each eigenvectors $\bv_i = \frac{\bx_i}{\normtwo{\bx_i}}$;\Comment{$r\times 3n = 3nr$ flops}
\State Normalize each eigenvectors $\bu_i = \frac{\bA \bv_i}{\sigma_i}$; \Comment{$r(m(2n-1)+m)$ flops}
\end{algorithmic} 
\end{algorithm}

By completing $\{\bv_{r+1}, \bv_{r+2}, \ldots, \bv_n\}$ into a full orthogonal basis via the Gram-Schmidt process, we obtain the orthogonal $\bV$.
And similarly, completing $\{\bu_{r+1}, \bu_{r+2}, \ldots, \bu_m\}$ into a full orthogonal basis using the Gram-Schmidt process, we obtain the orthogonal $\bU$.

\index{Gram–Schmidt}
The detailed analysis of the computational complexity for SVD is complicated. It is typically computed numerically using a two-step procedure. In the first step, the matrix is reduced to a bidiagonal matrix. This takes $\mathcalO(mn^2)$ flops (Section~\ref{section:bidiagonal-decompo}). 
In the second step, the SVD of the bidiagonal matrix is computed, requiring $\mathcalO(n)$ iterations with each involving $\mathcalO(n)$ flops. Thus, the overall cost is $\mathcalO(mn^2)$ flops. For those who are interested in delving deeper into the computation of SVD, please refer to Section~\ref{section:eigenvalue-problem} or \citet{trefethen1997numerical, golub2013matrix, kishore2017literature} for more details. 

\subsection*{Randomized Method for Computing the SVD Approximately}
Suppose  the matrix $\bA$ admits a rank decomposition (Theorem~\ref{theorem:rank-decomposition}):
$$
\underset{m\times n}{\bA} = \underset{m\times r}{\bD}\gapthree \underset{r\times n}{\bF},
$$
where the columns of $\bD$ span the same column space as $\bA$: $\cspace(\bD)=\cspace(\bA)$. If we orthogonalize the columns of $\bD=[\bd_1, \bd_2, \ldots, \bd_r]$ into $\bQ_r=[\bq_1, \bq_2, \ldots, \bq_r]$ such that 
$$
\spn([\bq_1, \bq_2, \ldots, \bq_k]) = \spn([\bd_1, \bd_2, \ldots, \bd_k]), \gap \text{for $k\in \{1,2,\ldots,r\}$}.
$$
This can be achieved using  the \textit{reduced} QR decomposition via the Gram-Schmidt process,
with a complexity of $\mathcalO(mr^2)$ if $\bD\in \real^{m\times r}$ (Section~\ref{section:qr-gram-compute}; where we complete the $\bQ_r$ into an orthogonal one $\bQ=[\bQ_r, \bQ_2]$ such that $\bQ\bQ^\top=\bI$ for further analysis, but note that $\bQ_2$ is not needed for the final algorithm!). One can show that if the SVD of the $m\times n$ matrix $\widetildebE\triangleq\bQ^\top \bA \in \real^{m\times n}$ is given by $\widetildebE = \widetildebU\bSigma\bV^\top$, the SVD of $\bA$ can be obtained by 
\begin{equation}\label{equation:svd-random-au}
\bA=\bQ\widetildebE=\underbrace{(\bQ\widetildebU)}_{\triangleq\bU}\bSigma\bV^\top.
\end{equation}
Breaking down the above equation, we have:
$$
\begin{aligned}
\widetildebE = 
\begin{bmatrix}
\bQ_r^\top \bA \\
\bQ_2^\top \bA
\end{bmatrix}
&=
\widetildebU\bSigma\bV^\top
=
\begin{bmatrix}
\widetildebU_r & \widetildebU_{12} \\
\widetildebU_{21} & \widetildebU_{22}
\end{bmatrix}
\begin{bmatrix}
\bSigma_r & \bzero \\
\bzero & \bzero 
\end{bmatrix}
\begin{bmatrix}
\bV_r^\top\\
\bV_2^\top
\end{bmatrix}
= 
\begin{bmatrix}
\widetildebU_r \bSigma_r  \\
\widetildebU_{21}\bSigma_r 
\end{bmatrix}
\bV_r^\top \\
&\implies 
\left\{
\begin{aligned}
\bQ_r^\top \bA &= \widetildebU_r \bSigma_r\bV_r^\top;\\	
\bQ_2^\top \bA &= \widetildebU_{21}\bSigma_r \bV_r^\top.\\
\end{aligned}\right.
\end{aligned}
$$
By Equation~\eqref{equation:svd-random-au}, we have 
$$
\widetildebU = \bQ^\top \bU =
\begin{bmatrix}
\bQ_r^\top \\
\bQ_2^\top 
\end{bmatrix}
\begin{bmatrix}
\bU_r & \bU_2 
\end{bmatrix}
=
\begin{bmatrix}
\bQ_r^\top \bU_r & \bQ_r^\top\bU_2 \\
\bQ_2^\top \bU_r & \bQ_2^\top\bU_2 \\ 
\end{bmatrix}
=
\begin{bmatrix}
\widetildebU_r & \widetildebU_{12} \\
\widetildebU_{21} & \widetildebU_{22}
\end{bmatrix}.
$$
Since $\cspace(\bQ_r) = \cspace(\bA)=\cspace(\bU_r)$, it follows that $\bQ_2^\top \bU_r=\bzero $ given that $\bQ_2$ lies in the orthogonal complement of $\bQ_r$ (which is also the orthogonal complement of $\bU_r$). Therefore, 
$$
\begin{aligned}
\widetildebE &= 
\begin{bmatrix}
\bQ_r^\top \bA \\
\bzero
\end{bmatrix}
\end{aligned}
\gap \implies \gap
\bA=\bQ\widetildebE=
\begin{bmatrix}
\bQ_r & \bQ_2 
\end{bmatrix}
\begin{bmatrix}
\bQ_r^\top \bA \\
\bzero
\end{bmatrix}
=
\bQ_r\bQ_r^\top \bA.
$$
That is, $\bQ_r\bQ_r^\top$ is an orthogonal projection matrix that maps every column of $\bA$ onto itself (a projection matrix onto the column space of $\bA$; Definition~\ref{definition:projection_matrix_intro}).
We observe that $\bQ_r^\top \bA = \widetildebU_r \bSigma_r\bV_r^\top$ above is the SVD of $\bQ_r^\top \bA \in\real^{r\times n}$.

\begin{mdframed}[hidealllines=\mdframehidelineNote,backgroundcolor=\mdframecolor]
To conclude the observation, suppose we find a matrix $\bD\in \real^{m\times r}$ that spans the column space of $\bA\in \real^{m\times n}$. 
We perform the \textit{reduced} QR decomposition of $\bD=\bQ_r\bR$, where $\bQ_r$ also spans the column space of $\bA$. 
Then we compute the SVD of the \textbf{small} matrix $\bE=\bQ_r^\top\bA \in\real^{r\times n} $, which costs $\mathcalO(nr^2)$ flops: $\bE=\bQ_r^\top\bA=\widetildebU_r\bSigma_r\bV_r^\top$. The \textit{reduced} SVD of the large matrix can be obtained by 
$$
\bA = \bQ_r(\bQ_r^\top\bA) = \underbrace{\bQ_r\widetildebU_r}_{\triangleq\bU_r}\bSigma_r\bV_r^\top.
$$
By completing the matrices $\bU_r\in \real^{m\times r}$ and $\bV_r\in \real^{n\times r}$ into full orthogonal matrices, we find the \textit{full} SVD of $\bA$.
\end{mdframed}
Further, similar to the randomized algorithm for interpolative decomposition in Section~\ref{section:randomi-id}, by Lemma~\ref{lemma:column-basis-from-row-basis}, for any matrix $\bA\in \real^{m\times n}$, suppose that $\{\bg_1, \bg_2, \ldots, \bg_r\}$ is a set of vectors in $\real^n$, which forms a basis for the row space, then $\{\bA\bg_1, \bA\bg_2, \ldots, \bA\bg_r\}$ is a basis for the column space of $\bA$:
\begin{equation}\label{equation:random-svd-column-space1}
\cspace(\bA\bG) =\cspace(\bA)\gap \text{with} \gap \bG=[\bg_1, \bg_2, \ldots, \bg_r]
\end{equation}
And a small integer $k$ (say $k=10$) should be selected to over-sample such that there is a high probability that $\bG$=$[\bg_1, \ldots,\bg_r, \bg_{r+1}, \ldots, \bg_{r+k}]$ contains the row basis of $\bA$. 
Other methods to guarantee that $\bG$ can span the row basis of $\bA$ or that $\bA\bG$ can span the column space of $\bA$ are discussed in Remark~\ref{remark:source-row-basis}. 
Again, The choice of  $k=10$ is often suitable. 

The procedure is formulated in Algorithm~\ref{alg:svd-randomized}, where the end result costs $\mathcalO(mn(r+k))$ flops compared to the original $\mathcalO(mn^2)$ flops. 
We notice that the leading term in the computational cost of $\mathcalO(mn(r+k))$ flops in the algorithm comes from step 5, involving the matrix multiplication $\bQ_r^\top\bA$. A structured choice of the random matrix $\bG$ can further reduce this to $\mathcalO(mn\log(r+k))$ flops \citep{martinsson2019randomized, ailon2006approximate}. 


\begin{algorithm}[h] 
\caption{A Randomized Method to Compute the SVD} 
\label{alg:svd-randomized} 
\begin{algorithmic}[1] 
\Require 
Rank-$r$ matrix $\bA$ with size $m\times n $; 
\State Decide the over-sampling parameter $k$ (e.g., $k=10$), and let $z=r+k$; 
\State Generate $r+k$ Gaussian random vectors in $\real^n$ into columns of matrix $\bG\in \real^{n\times (r+k)}$;\Comment{i.e., probably contain the row basis of $\bA$}
\State Initialize $\bD=\bA\bG \in \real^{m\times (r+k)}$; \Comment{probably $\cspace(\bD)=\cspace(\bA)$, $m(2n-1)(r+k)$ flops}
\State Compute full QR decomposition $\bD = \underbrace{\bQ_r}_{m\times (r+k)}\bR$;\Comment{$\mathcalO(mz^2)$ flops}
\State Form matrix $\bE=\bQ_r^\top \bA \in\real^{(r+k)\times n}$; \Comment{$nz(2m-1)$ flops}
\State Compute the SVD of the \textbf{small} matrix $\bE$: $\bE = \bU_0\bSigma\bV^\top$; \Comment{$\mathcalO(nz^2)$ flops}
\State Form $\bU=\bQ_r\bU_0$ such that the \textit{reduced} SVD of $\bA=\bU\bSigma\bV^\top$; \Comment{$mz(2z-1)$ flops}
\end{algorithmic} 
\end{algorithm}

\section{Polar Decomposition}\label{section:polar_decom}

A decomposition closely related to the SVD is the \textit{polar form} of a matrix. 
In the context of continuum mechanics, it is imperative to distinguish between stretching and rotation.
The polar decomposition factors any matrix into an orthogonal matrix (corresponds to rotation) and a symmetric matrix (corresponds to stretching or compression, see Section~\ref{section:coordinate-transformation}).

\begin{theoremHigh}[Polar Decomposition]\label{theorem:polar-decomposition}
Let $\bA\in\real^{m\times n}$. Then $\bA$ can be factored as 
\begin{itemize}
\item \textbf{Case $m>n$: left polar decomposition.} $\bA=\bQ_l\bS_l$, where $\bS_l^2=\bA^\top\bA$ is PSD and is \textbf{uniquely} determined. 
The factor $\bQ_l$ has orthonormal columns, and it is \textbf{uniquely} determined if $\rank(\bA)=n$.
\item \textbf{Case $m<n$: right polar decomposition.} $\bA=\bS_r\bQ_r$, where $\bS_r^2=\bA\bA^\top$ is PSD and is \textbf{uniquely} determined.
The factor $\bQ_r$ has orthonormal rows, and it is \textbf{uniquely} determined if $\rank(\bA)=m$.
\item \textbf{Case $m=n$: left/right polar decomposition.} $\bA=\bQ\bS_l=\bS_r\bQ$, where $\bS_l^2=\bA^\top\bA$ and $\bS_r^2=\bA\bA^\top$ are PSD and are \textbf{uniquely} determined. The factor $\bQ$ is orthonoal, and it is the same for both the left and right polar decompositions. $\bQ$ is \textbf{uniquely} determined if $\bA$ is nonsingular ($\rank(\bA)=n$).
\end{itemize}
\noindent
Note in all cases, the PSD factors are uniquely determined, and become PD if $\bA$ has full rank  (full row rank or column rank). The semi-orthogonal factors $\bQ_l, \bQ_r$, and $\bQ$ are uniquely determined only when $\bA$ has full rank.~\footnote{When $\bA$ is complex, then the orthogonal (resp. semi-orthogonal) matrices become unitary (resp. semi-unitary) matrices, and the PSD matrices become  complex Hermitian and PSD matrices.}
\end{theoremHigh}
\begin{proof}[of Theorem~\ref{theorem:polar-decomposition}]
Let the SVD of $\bA$ be $\bA=\bU\bSigma\bV^\top={(\bU\bV^\top)}{(\bV\bSigma\bV^\top)}\triangleq\bQ_l\bS_l$ such that $\bS_l^2=\bV\bSigma^2\bV^\top=\bA^\top\bA$. Since $\bA^\top\bA$ is PSD, $\bS_l$ is uniquely determined (Theorem~\ref{theorem:unique-factor-pd}). 
If further $\rank(\bA)=n$, i.e., $\bA$ has full (column) rank, $\bA^\top\bA$ is PD and $\bS_l$ has full rank (Theorem~\ref{theorem:unique-factor-pd}) such that $\bQ_l=\bA\bS_l^{-1}$, implying $\bQ_l$ is uniquely determined.

The second case can be similarly proved such that $\bA=\bU\bSigma\bV^\top=(\bU\bSigma\bU^\top)(\bU\bV^\top)\triangleq\bS_r\bQ_r$.
Since $\bA\bA^\top$ is PSD, $\bS_r$ is uniquely determined. If further $\rank(\bA)=m$, $\bA\bA^\top$ is PD and $\bS_r$ has full rank such that $\bQ_r=\bS_r^{-1}\bA$ is uniquely determined.

The third case is a combination of the previous cases. This completes the proof.
\end{proof}
\begin{exercise}[Trace of PSD in Polar Decomposition]\label{exercise:trace_polar}
Show that the PSD matrices in the polar decomposition satisfy  $\trace(\bS_r)$ or $\trace(\bS_l)$ is equal to the sum of the singular values of $\bA$.
\end{exercise}
\begin{exercise}[Normal from Polar]
Let $\bA\in\real^{n\times n}$ be nonsingular, and suppose it admits the polar decomposition $\bA=\bS_r\bQ$, where $\bS_r$ is PD and $\bQ$ is orthogonal. Show that $\bA$ is normal if and only if $\bS_r\bQ=\bQ\bS_r$.
\end{exercise}
\begin{exercise}
Let $\bA,\bB\in\real^{n\times n}$ be orthogonal, and let $\bA+\bB$ be nonsingular. Show that the orthogonal factor in the polar decomposition of $\bA+\bB$ is $\bA(\bA^\top\bB)^{1/2}$. 
\end{exercise}

For the other way around, given the left polar decomposition of $\bA = \bQ_l\bS_l$. Since $\bS_l$ is symmetric, it admits the spectral decomposition $\bS_l = \bV\bSigma\bV^\top$. Therefore, the SVD is obtained as $\bA = \bQ_l(\bV\bSigma\bV^\top) \triangleq \bU\bSigma\bV^\top$ by letting $\bU\triangleq\bQ_l\bV$.

\begin{example}[Polar Decomposition]
Consider the matrix
$
\scriptsize
\bA = \begin{bmatrix}
1 & 1 \\
1 & 1\\
\end{bmatrix}.
$
We can find two left polar decompositions:
$$
\bA = 
\begin{bmatrix}
1 & 0 \\
0 & 1\\
\end{bmatrix}
\begin{bmatrix}
1 & 1 \\
1 & 1\\
\end{bmatrix}
= 
\begin{bmatrix}
0 & 1 \\
1 & 0\\
\end{bmatrix}
\begin{bmatrix}
1 & 1 \\
1 & 1\\
\end{bmatrix}.
$$
And two right polar decompositions:
$$
\bA = 
\begin{bmatrix}
1 & 1 \\
1 & 1\\
\end{bmatrix}
\begin{bmatrix}
1 & 0 \\
0 & 1\\
\end{bmatrix}
= 
\begin{bmatrix}
1 & 1 \\
1 & 1\\
\end{bmatrix}
\begin{bmatrix}
0 & 1 \\
1 & 0\\
\end{bmatrix}.
$$
Thus, the polar decomposition is not unique if $\bA$ is singular.
\end{example}

The polar decomposition reveals the closest orthogonal matrix and the closest rotation matrix to a given matrix, which we will discuss below. The problems of finding the closest scalar multiple of an orthogonal matrix, one-sided orthogonal transformations, and two-sided orthogonal transformations are addressed in Section~\ref{section:svd_von_cou_etc}, as a result of von Neumann's trace theorem and polar decomposition.

\paragraph{Complex scalar analogy.}
In complex scalar case,
each complex scalar $x$  can be factored as $x = re^{i\theta}$, where $r$ is a real and nonnegative number (which can be thought of as a $1 \times 1$ positive semidefinite matrix) and $e^{i\theta}$ has a modulus of one (which can be thought of as a $1 \times 1$ orthogonal matrix). The factor $r = \abs{x}$ is always uniquely determined, but the factor $e^{i\theta}$ is uniquely determined only if $x$ is nonzero. Thus, the polar decomposition is a matrix analog of this scalar factorization.

\paragraph{Rotation or reflection.}
Consider a real and nonsingular matrix $\bA\in\real^{n\times n}$ and its left polar decomposition $\bA=\bQ_l\bS_l$ (the following result applies to the right polar decomposition analogously), it can be shown that 
\begin{equation}\label{equation:det_polar}
\det(\bQ_l) = \sign(\det(\bA)).
\end{equation}
We know that $\det(\bA)=\det(\bQ_l)\det(\bS_l)$ (Remark~\ref{remark:determinant-intermezzo}). 
And since $\bS_l$ is PD and the determinant is equal to the product of its eigenvalues (Theorem~\ref{theorem:eigen_trace}), then $\det(\bS_l)> 0$.
As the determinant of an orthogonal matrix  is either 1 or $-1$ (Exercise~\ref{exercise:eig_orthos}), we conclude the result in Equation~\eqref{equation:det_polar}.
This shows that $\bQ_l$ is a \textit{rotation} when $\det(\bA)$ is positive and a \textit{reflection} when it is negative (see examples in Figure~\ref{fig:rotation}).

We prove that the polar decomposition gives us the closest orthogonal matrix to $\bA$ in terms of the Frobenius norm.

\begin{theorem}[Closest Orthogonal Matrix \citep{julian2015polar}]\label{theorem:closes_ortho}
Given the left polar decomposition of matrix $\bA\in\real^{n\times n}$: $\bA=\bQ_l\bS_l$. Then $\bQ_l$ is the closest orthogonal matrix to $\bA$ in the Frobenius norm sense:~\footnote{When $\bA$ is complex,  we can similarly find the closest unitary matrix.}
$$
\normf{\bA - \bQ_l}^2 < \normf{\bA - \bP}^2, \quad \forall\, \bP \text{ is orthogonal}, \bP \neq \bQ_l.
$$
Furthermore, this shortest distance is 
$$
\normf{\bA - \bQ_l}^2 = \sum_{i=1}^n (\sigma_i - 1)^2,
$$
where $\sigma_i$'s are the singular values of $\bA$.
\end{theorem}
\begin{proof}[of Theorem~\ref{theorem:closes_ortho}]
To see this, we have 
$$
\begin{aligned}
\normf{\bA - \bQ_l}^2 &= \trace((\bA - \bQ_l)^\top (\bA - \bQ_l))= \normf{\bA}^2 - 2 \trace(\bA^\top \bQ_l) + \trace(\bI) \\
& = \normf{\bA}^2 - 2 \trace(\bS_l) + \trace(\bI).
\end{aligned}
$$
Similarly, for any orthogonal matrix $\bP$, we have 
$
\normf{\bA - \bP}^2 = \normf{\bA}^2 - 2 \trace(\bA^\top \bP) + \trace(\bI).
$
Subtracting the two equations:
$$
\normf{\bA - \bQ_l}^2 - \normf{\bA - \bP}^2 = 2 \trace(\bA^\top \bP - \bS_l).
$$
Since $\bS_l$ is PSD, it admits the spectral decomposition $\bS_l = \bV \bLambda \bV^\top$. We have 
$$
\begin{aligned}
\normf{\bA - \bQ_l}^2 - \normf{\bA - \bP}^2 
&= 2 \trace(\bV \bLambda \bV^\top \bQ_l^\top \bP - \bV \bLambda \bV^\top)\\
&= 2 \trace( \bLambda \underbrace{\bV^\top \bQ_l^\top \bP\bV}_{\triangleq\bY} -  \bLambda ) = 2 \sum_i \lambda_i (y_{ii} - 1) \\
\end{aligned}
$$
where the second equality follows from the trace invariance under cyclic permutations, 
and  $\bY \triangleq \bV^\top \bQ_l^\top \bP \bV$ is an orthogonal matrix, implying $\abs{y_{ii}} \leq 1$. Since $\bS_l$ is PSD ($\lambda_i \geq 0$),
this sum  achieves its maximum (equal to zero) only when $y_{ii} = 1$:
$$
\normf{\bA - \bQ_l}^2 - \normf{\bA - \bP}^2 \leq 0
\gap\implies\gap 
\normf{\bA - \bQ_l}^2 \leq \normf{\bA - \bP}^2.
$$
Moreover, equality holds only when
$$
\bY = \bV^\top \bQ_l^\top \bP \bV = \bI \gap\implies\gap  \bQ_l^\top \bP = \bV \bV^\top = \bI \gap\implies\gap  \bP = \bQ_l,
$$
which proves the inequality when $\bP \neq  \bQ_l$.

For the second part, let $\bA=\bU\bSigma\bV^\top$ be the SVD of $\bA$.  We have 
$
\normf{\bA - \bQ_l}^2 = \normf{\bA}^2 - 2 \trace(\bS_l) + \trace(\bI) = \sum_{i=1}^n (\sigma_i^2 - 2 \lambda_i + 1).
$
Notice that 
$$
\bA^\top \bA = \bS_l^\top \bS_l \gap\implies\gap  \bV \bSigma^2 \bV^\top = \bV \bLambda^2 \bV^\top \gap\implies\gap  \lambda_i = \pm \sigma_i.
$$
Since $\bS_l$ is PSD, $\lambda_i = \sigma_i$, and we have:
$
\normf{\bA - \bQ_l}^2 = \sum_{i=1}^n (\sigma_i^2 - 2 \sigma_i + 1) = \sum_{i=1}^n (\sigma_i - 1)^2.
$
This completes the proof.
\end{proof}

From the above theorem, in the case that $\bA$ is real and $\det(\bA) > 0$, the left polar decomposition $\bQ_l$ is a rotation, and thus it is  the closest rotation to $\bA$. 
When $\det(\bA) = 0$, an SVD can be computed $\bA=\bU\bSigma\bV^\top$ such that $\bU \bV^\top$ is a rotation (since $\bQ_l=\bU\bV^\top$; and the sign for the left and right null space basis vectors can be chosen arbitrarily). 
However, when $\det(\bA) < 0$, $\bQ_l$ must be a reflection. To make $\bQ_l$ a rotation when $\det(\bA) < 0$, $\bS_l$ would need to have an odd number of negative eigenvalues. 
This can be achieved by flipping some signs in the SVD used to derive $\bA = \bQ_l\bS_l$.
However, Theorem~\ref{theorem:closes_ortho} no longer holds in this case. 


\begin{theorem}[Closest Rotation Matrix]\label{theorem:closes_rotat}
Let $\bA\in\real^{n\times n}$ be nonsingular  with $\det(\bA)<0$, and let its singular value decomposition be given by  $\bA=\bU\bSigma\bV^\top$. Then the closest rotation to $\bA$ is given by 
$$
\bU\bD\bV^\top,
\quad \text{where}\quad
\text{$\bD$ is diagonal with }\,\,
d_{ii} = \left\{ \begin{array}{cc}
	1 & \text{if} \ i < n; \\
	-1 & \text{if} \ i = n,
\end{array} \right.
$$
\end{theorem}
\begin{proof}[of Theorem~\ref{theorem:closes_rotat}]
From the discussion above, to find the closest rotation to $\bA\in\real^{n\times n}$ when $\det(\bA)<0$, we turn to a constrained minimization problem:
$$
\min_{\bP \in \real^{n \times n}} \normf{\bA - \bP}^2 \quad \text{s.t.} \quad \bP^\top \bP = \bI, \, \det(\bP) = 1.
$$
The corresponding Lagrangian (see \citet{boyd2004convex} for more details) is:
$$
L(\bP, \bM, m) = \normf{\bA - \bP}^2 + \trace((\bP^\top \bP - \bI)^\top \bM) + m(\det(\bP) - 1),~\footnote{$\trace((\bP^\top \bP - \bI)^\top \bM)$ follows from the Frobenius matrix inner product: $\langle \bA, \bB\rangle=\sum_{i,j=1}^{m,n}a_{ij}b_{ij}=\trace(\bA\bB^\top)=\trace(\bA^\top\bB)$.}
$$
where $m$ is a single Lagrange multiplier, and $\bM$ is a symmetric matrix of Lagrange multipliers (with one multiplier per unique equation in the symmetric constraint $\bP^\top \bP = \bI$). Differentiating with respect to $\bM$ and $m$  recovers the constraints, while differentiating with respect to $\bP$ and setting to zero yields:
$$
\frac{\partial L(\bP, \bM, m)}{\partial \bP} 
= -2\bA -2\bP+ 2\bP\bM + m \det(\bP)(\bP^{-1})^\top = \bzero.
~\footnote{From the fact that $\frac{\partial \trace(\bX\bA)}{\partial \bX}=\bA^\top$, and $\frac{\partial \det(\bX)}{\partial \bX}=\det(\bX)(\bX^{-1})^\top$.}
$$
Since the constraints must be satisfied at the optimum, we can substitute  $\det(\bP) = 1$ and $[\bP^{-1}]^\top = \bP$:
$$
2\bA = 2\bP\bM + m\bP -2\bP \gap\implies\gap  \bA = \bP \bigg( \bM + \frac{m-2}{2} \bI \bigg).
$$
This means that an optimal rotation $\bP$ puts $\bA$ in a very similar form to the polar decomposition: $\bA = \bP\bS_l$, where $\bS_l \triangleq \left( \bM + \frac{m-2}{2} \bI \right)$ is a symmetric matrix. Using the  spectral decomposition $\bS_l = \bQ \bLambda \bQ^\top$, we have 
$$
\bA^\top \bA = \bS_l^\top \bS_l = \bQ \bLambda^2 \bQ^\top.
$$
As shown above, $\lambda_i = \pm \sigma_i$ ($\sigma_i$ is the $i$-th singular value of $\bA$), with the specific signs determined by $\bP$. 
By encoding the signs in a diagonal matrix $\bD$ so that $\lambda_i = d_{ii} \sigma_i$ with $ d_{ii}\in\{-1,1\}$ (i.e., $\bLambda = \bD\bSigma$), we see that the optimal distance must be of the form:
\begin{equation}\label{equation:polar_minrot1}
	\normf{\bA - \bP}^2 = \normf{\bP\bQ \bLambda \bQ^\top - \bP}^2 = \trace(\bLambda^2 - 2 \bLambda + \bI) = \sum_i (\lambda_i - 1)^2 = \sum_i (d_{ii} \sigma_i - 1)^2. 
\end{equation}
Furthermore, any assignment of signs multiplying to $-1$ is obtained by some rotation $\bP$. To see this, consider arbitrary signs, $\bLambda = \bD\bSigma$, $\det(\bD) = -1$:
\begin{equation}\label{equation:polar_minrot2}
	\bA = \bU \bSigma \bV^\top = \bU (\bD \bV^\top \bV \bD) \bSigma \bV^\top = \underbrace{(\bU \bD \bV^\top)}_{\bP} \underbrace{(\bV \bLambda \bV^\top)}_{\bS_l, \, \bV\triangleq\bQ}, 
\end{equation}
and $\det(\bP) = \frac{\det(\bA)}{\det(\bS_l)} = 1$. 
In other words, we can find an optimal $\bP$ by optimizing over $d_{ii} = \pm 1$. Since any $\sigma_i > 0$ when $\det(\bA)<0$, and 
$$
\det(\bA)=\det(\bU)\det(\bV)\prod_{i=1}^{n}\sigma_i \gap \text{and}\gap \det(\bS_l)=\det(\bV)^2\prod_{i=1}^{n}\lambda_i,
$$
each negative sign increases distance \eqref{equation:polar_minrot1} by $4\sigma_i$.
Therefore, the optimal sign assignment must negate only the single smallest singular value (i.e., we optimize Equation~\eqref{equation:polar_minrot1} from \eqref{equation:polar_minrot2}). Choosing these optimal signs, $d_{ii} = \left\{ \begin{array}{cc}
	1 & \text{if} \ i < n; \\
	-1 & \text{if} \ i = n,
\end{array} \right.$ (assuming singular values are sorted in decreasing order), Equation~\eqref{equation:polar_minrot2} yields the closest rotation $\bP = \bU \bD \bV^\top$.
\end{proof}

\section{Properties of the SVD}\label{section:property-svd}
\subsection{Four Subspaces in SVD}\label{section:four-space-svd}
\index{Fundamental theorem}
For any matrix $\bA\in  \real^{m\times n}$, the following properties hold:
\begin{itemize}
\item $\nspace(\bA)$ is the orthogonal complement of the row space $\cspace(\bA^\top)$ in $\real^n$: $\dim(\nspace(\bA))+\dim(\cspace(\bA^\top))=n$;

\item $\nspace(\bA^\top)$ is the orthogonal complement of the column space $\cspace(\bA)$ in $\real^m$: $\dim(\nspace(\bA^\top))+\dim(\cspace(\bA))=m$.
\end{itemize}
This is known as  the fundamental theorem of linear algebra, also referred to as the rank-nullity theorem. And the proof can be found in Appendix~\ref{appendix:fundamental-rank-nullity}. Additionally, we can find the bases for the four subspaces using the CR decomposition in Appendix~\ref{appendix:cr-decomposition-four-basis}.
Specially, the SVD allows us to identify an orthonormal basis for each subspace. 

\begin{figure}[h!]
	\centering
	\includegraphics[width=0.98\textwidth]{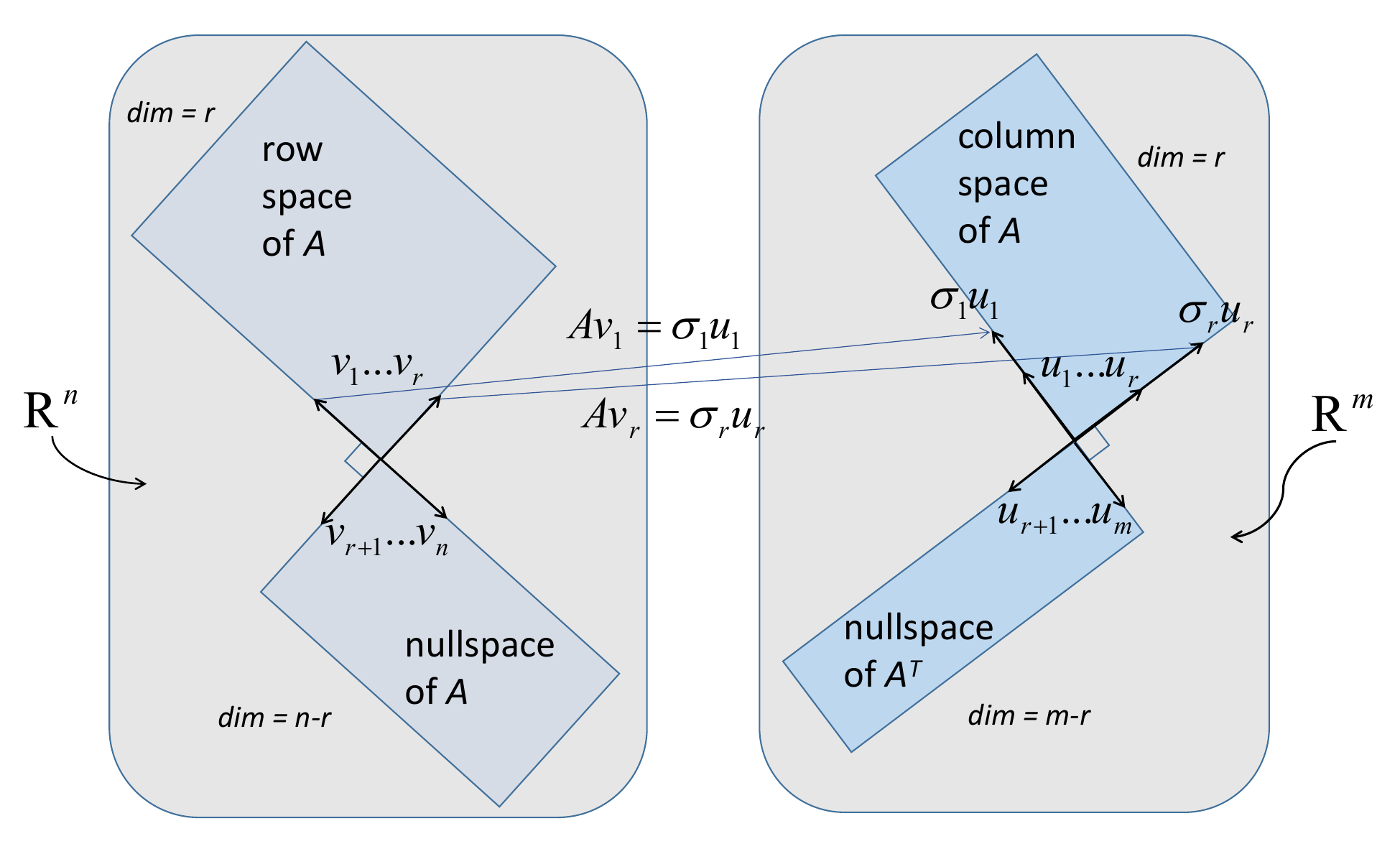}
	\caption{Orthonormal bases that diagonalize $\bA$ using the SVD.}
	\label{fig:lafundamental3-SVD}
\end{figure}

\begin{proposition}[Four Orthonormal Bases]\label{proposition:svd-four-orthonormal-Basis}
Given the full SVD of a matrix $\bA = \bU \bSigma \bV^\top$, where $\bU=[\bu_1, \bu_2, \ldots,\bu_m]$ and $\bV=[\bv_1, \bv_2, \ldots, \bv_n]$ represent the column partitions of $\bU$ and $\bV$, respectively. Then  we have the following property:
\begin{itemize}
\item $\{\bv_1, \bv_2, \ldots, \bv_r\} $ is an orthonormal basis of $\cspace(\bA^\top)$;

\item $\{\bv_{r+1},\bv_{r+2}, \ldots, \bv_n\}$ is an orthonormal basis of $\nspace(\bA)$;

\item $\{\bu_1,\bu_2, \ldots,\bu_r\}$ is an orthonormal basis of $\cspace(\bA)$;

\item $\{\bu_{r+1}, \bu_{r+2},\ldots,\bu_m\}$ is an orthonormal basis of $\nspace(\bA^\top)$. 
\end{itemize}

The relationships among the four subspaces are illustrated in Figure~\ref{fig:lafundamental3-SVD}, where $\bA$ transforms the row basis $\bv_i$ into the column basis $\bu_i$ by $\sigma_i\bu_i=\bA\bv_i$ for all $i\in \{1, 2, \ldots, r\}$.
\end{proposition}
\begin{proof}[of Proposition~\ref{proposition:svd-four-orthonormal-Basis}]
According to Proposition~\ref{proposition:rank-of-symmetric},  in the case of a  symmetric matrix $\bA^\top\bA$, the column space $\cspace(\bA^\top\bA)$ is spanned by the eigenvectors, thus $\{\bv_1,\bv_2, \ldots, \bv_r\}$ constitutes  an orthonormal basis for $\cspace(\bA^\top\bA)$.

Considering the following:
\begin{enumerate}
\item Since $\bA^\top\bA$ is symmetric,  the row space of $\bA^\top\bA$ is equivalent to its column space.

\item All rows of $\bA^\top\bA$ are  linear combinations of the rows of $\bA$, meaning the row space of $\bA^\top\bA$ is a subset of the row space of $\bA$, i.e., $\cspace(\bA^\top\bA) \subseteq \cspace(\bA^\top)$. 

\item Since $\rank(\bA^\top\bA) = \rank(\bA)$ by Lemma~\ref{lemma:rank-of-ata}, we can conclude that:

The row space of $\bA^\top\bA$ = the column space of $\bA^\top\bA$ =  the row space of $\bA$, i.e., $\cspace(\bA^\top\bA) = \cspace(\bA^\top)$. Thus, $\{\bv_1, \bv_2,\ldots, \bv_r\}$ forms an orthonormal basis for $\cspace(\bA^\top)$. 
\end{enumerate}

Moreover, the subspace spanned by $\{\bv_{r+1}, \bv_{r+2},\ldots, \bv_n\}$ serves as the orthogonal complement to the subspace spanned by $\{\bv_1,\bv_2, \ldots, \bv_r\}$. Therefore, $\{\bv_{r+1},\bv_{r+2}, \ldots, \bv_n\}$ is an orthonormal basis for $\nspace(\bA)$. 

Applying this process to $\bA\bA^\top$ allows us to establish the remaining claims in the lemma.  

Alternatively, since we have established that $\{\bv_1, \bv_2,\ldots, \bv_r\}$ forms an orthonormal basis for the row space of $\bA$, i.e., the basis for $\cspace(\bA^\top)$. 
According to Lemma~\ref{lemma:column-basis-from-row-basis}~\footnote{For any matrix $\bA$, let $\{\br_1, \br_2, \ldots, \br_r\}$ be a set of vectors in $\real^n$, which forms a basis for the row space, then $\{\bA\br_1, \bA\br_2, \ldots, \bA\br_r\}$ is a basis for the column space of $\bA$.} 
and  given that  $\bu_i = \frac{\bA\bv_i}{\sigma_i},\, \forall i \in\{1, 2, \ldots, r\}$, we observe that $\{\bu_1,\bu_2, \ldots,\bu_r\}$ forms an orthonormal basis for the column space of $\bA$. 
\end{proof}

\index{Orthogonal projection}
\index{Projection matrix}
\index{Projector}
\subsection{SVD-Related Orthogonal Projections}\label{section:svd_ortho_proj}
In the context of SVD, several key orthogonal projections arise from the four fundamental subspaces.
A comprehensive  discussion of the orthogonal projection is available in Appendix~\ref{section:by-geometry-hat-matrix} and Definition~\ref{definition:projection_matrix_intro}. In words,  orthogonal projection matrices are matrices that are symmetric and idempotent. The projection matrix will project any vector onto its column space. 
The idempotency has a geometrical meaning such that projecting a vector twice is equivalent to projecting it once. 
The symmetry also has a geometrical meaning such that the distance between the original vector and the projected vector (which lies within the column space of the projection matrix) is minimized. 
Suppose $\bA=\bU\bSigma\bV^\top$ is the full SVD of $\bA$ with rank $r$. If we have the following column partitions
\[
\begin{blockarray}{ccc}
	\begin{block}{c[cc]}
\bU=&	\bU_r & \bU_m   \\
	\end{block}
	&m\times r & m\times (m-r)   \\
\end{blockarray}
,\qquad
\begin{blockarray}{ccc}
	\begin{block}{c[cc]}
\bV=	&	\bV_r & \bV_n   \\
	\end{block}
	&n\times r & n\times (n-r)   \\
\end{blockarray},
\]
where $\bU_r$ and $\bV_r$ consist of the first $r$ columns of $\bU$ and $\bV$, respectively.
Then  the four orthogonal projections can be obtained by 
$$
\begin{aligned}
\bV_r\bV_r^\top &= \text{projection onto $\cspace(\bA^\top)$},
\quad
&\bV_n\bV_n^\top &=\text{projection onto $\nspace(\bA)$},\\
\bU_r\bU_r^\top &= \text{projection onto $\cspace(\bA)$},
\quad
&\bU_m\bU_m^\top &= \text{projection onto $\nspace(\bA^\top)$}.\\
\end{aligned}
$$

\index{Singular}
\index{Nonsingular}
\subsection{Relationship between Singular Values and Determinant}\label{section:nonsingul_inver_svd}
For a square matrix $\bA\in \real^{n\times n}$, its singular value decomposition  is given by  $\bA = \bU\bSigma\bV^\top$. It follows that 
$$
|\det(\bA)| = |\det(\bU\bSigma\bV^\top)| = |\det(\bSigma)| = \sigma_1 \sigma_2\ldots \sigma_n.
$$
If all the singular values $\sigma_i$ are nonzero, then $\det(\bA)\neq 0$. That is, $\bA$ is \textbf{nonsingular}. 
Conversely, if there is at least one singular value such that $\sigma_i =0$, then $\det(\bA)=0$, signifying that $\bA$ lacks  full rank and is not invertible. Then the matrix is called \textbf{singular}. This explains why the $\sigma_i$ values are referred to as the singular values.

\subsection{Orthogonal Equivalence}

We previously defined that for any nonsingular matrix $\bP$ and square matrix $\bA$, the matrices $\bA$ and $\bP\bA\bP^{-1}$ are similar (Definition~\ref{definition:similar-matrices}).
The concept of \textit{orthogonal equivalence} extends this idea to rectangular matrices.

\index{Orthogonal equivalence}
\begin{definition}[Orthogonal Equivalent Matrices\index{Orthogonal equivalent matrices}]
Given two nonsingular matrices $\bM\in\real^{m\times m}$ and $\bS\in\real^{n\times n}$, 
the matrices $\bA\in\real^{m\times n}$ and $\bM\bA\bS$ are called \textit{equivalent matrices}.

Similarly, given any orthogonal matrices $\bU\in\real^{m\times m}$ and $\bV\in\real^{n\times n}$, the matrices $\bA\in\real^{m\times n}$ and $\bU\bA\bV$ are called \textit{orthogonally equivalent matrices}.
In the complex domain, if $\bU$ and $\bV$ are unitary, they are called \textit{unitarily equivalent}.
\end{definition}

Next, we encounter the following property concerning orthogonally equivalent matrices.

\begin{lemma}[Orthogonal Equivalent Matrices]\label{lemma:orthogonal-equivalent-matrix}
Given any orthogonally equivalent matrices $\bA$ and $\bB$, their singular values are identical.
\end{lemma}
\begin{proof}[of Lemma~\ref{lemma:orthogonal-equivalent-matrix}]
Since $\bA$ and $\bB$ are orthogonally equivalent, there exist orthogonal matrices $\bU$ and $\bV$ such that $\bB = \bU\bA\bV$. Thus,
$$
\bB\bB^\top = (\bU\bA\bV)(\bV^\top\bA^\top\bU^\top) = \bU\bA\bA^\top\bU^\top.
$$
This shows that  $\bB\bB^\top$ and $\bA\bA^\top$ are similar matrices. 
According to Proposition~\ref{proposition:eigenvalue-similar-matrices}, the eigenvalues of similar matrices are the same, thereby proving the singular values of $\bA$ and $\bB$ are the same.
\end{proof}

\subsection{SVD for QR}
\begin{lemma}[SVD for QR]\label{lemma:svd-for-qr}
Let $\bA\in \real^{m\times n}$ with $m\geq n$, and let its   full QR decomposition be given by $\bA=\bQ\bR$, where $\bQ\in \real^{m\times m}$ is orthogonal and $\bR\in \real^{m\times n}$ is upper triangular. Then, $\bA$ and $\bR$ have the same singular values and right singular vectors.
\end{lemma}
\begin{proof}[of Lemma~\ref{lemma:svd-for-qr}]
We notice that $\bA^\top\bA = \bR^\top\bR$ such that $\bA^\top\bA$ and $\bR^\top\bR$ share the same eigenvalues and eigenvectors. 
Consequently, $\bA$ and $\bR$ have the same singular values and right singular vectors (i.e., the eigenvectors of $\bA^\top\bA$ or $\bR^\top\bR$).
\end{proof}

The  lemma above plays a crucial role in demonstrating the existence and characteristics of the rank-revealing QR decomposition (Section~\ref{section:rank-one-qr-revealing}). 
Moreover, it also implies that the SVD of a matrix can be constructed by performing the QR decomposition on the same matrix. 
Specifically, if the QR decomposition of $\bA$ is given by $\bA=\bQ\bR$ and the SVD of $\bR$ is given by $\bR=\bU_0 \bSigma\bV^\top$, then the SVD of $\bA$ can be obtained by 
$
\bA =\underbrace{ \bQ\bU_0}_{\triangleq\bU} \bSigma\bV^\top.
$

\subsection{SVD for Square Matrices}
When the underlying matrix $\bA\in\real^{n\times n}$ is square,  the SVD can be expressed in a special form involving orthogonal similarity transformations  (Definition~\ref{definition:simiar_congru}).
Given the full SVD  $\bA=\bU\bSigma\bV^\top$, where $\bSigma=\footnotesize\begin{bmatrix}
\bSigma_r & \bzero \\
\bzero & \bzero 
\end{bmatrix}$ and $\bSigma_r\in\real^{r\times r}$ is the diagonal singular matrix, 
we have 
$$
\bA=\bU\bSigma\bV^\top=
\bU\bSigma\bV^\top\bU\bU^\top,
$$
where $\bQ\triangleq\bV^\top\bU$ is an orthogonal matrix with the block form 
$\bQ=
\footnotesize\begin{bmatrix}
\bQ_{11} & \bQ_{12}\\
\bQ_{21} & \bQ_{22}
\end{bmatrix}$,  $\bQ_{11}\in\real^{r\times r}$,  $\bQ_{22}\in\real^{(n-r)\times (n-r)}$, and $\bQ_{11}\bQ_{11}^\top+\bQ_{12}\bQ_{12}^\top=\bI_r$.
Therefore, 
\begin{equation}
\bA=
\bU
\begin{bmatrix}
	\bSigma_r\bQ_{11} & \bSigma_r\bQ_{12}\\
	\bzero & \bzero 
\end{bmatrix}
\bU^\top,
\gap 
\text{with $\bQ_{11}\bQ_{11}^\top+\bQ_{12}\bQ_{12}^\top=\bI_r$}.
\end{equation}
This shows $\bA$ and $\footnotesize\begin{bmatrix}
	\bSigma_r\bQ_{11} & \bSigma_r\bQ_{12}\\
	\bzero & \bzero 
\end{bmatrix}$ are orthogonally similar.

\subsection{Positive Singular Values}
Suppose $\bA\in\real^{m\times n}$ with $m\geq n$ admits all positive singular values, then there exists a positive value $\eta\geq 0$ such that 
\begin{equation}
\normtwo{\bA\bz} \geq \eta \normtwo{\bz}.
\end{equation}
To see this, suppose $\bA$ admits full SVD $\bA=\bU\bSigma\bV^\top$. We have, by the invariance under orthogonal transformations of $\ell_2$ vector norms, 
$$
\normtwo{\bA\bz}=\normtwo{\bU\bSigma\bV^\top\bz}=\normtwo{\bSigma\bV^\top\bz} \geq 
\eta \normtwo{\bV^\top\bz}=
\eta \normtwo{\bz}.
$$

\index{Machine precision}
\subsection{Numerical Rank}
The existence of the SVD is essential for defining the effective rank of a matrix.
\begin{definition}[Effective Rank vs Exact Rank]\label{definition:effective-rank-in-svd}
The \textit{effective rank},  also known as the \textit{numerical rank},  is defined as follows:
Following Lemma~\ref{lemma:rank-equal-singular}, the number of nonzero singular values is equal to the rank of a matrix.
Let the $i$-th largest singular value of $\bA$ be denoted as $\sigma_i(\bA)$.
Then if $\sigma_r(\bA)\gg \sigma_{r+1}(\bA)\approx 0$, we refer to $r$ as the numerical rank of $\bA$ 
\footnote{or if $\sigma_r(\bA)\gg \sigma_{r+1}(\bA)= \mathcalO(\mu)$, where $\mu$ is the \textit{machine precision}.}
\footnote{In this context, we use the notation ``$a\gg b$" to indicate that ``$a$ is much larger than $b$."}. Whereas, when $\sigma_i(\bA)>\sigma_{r+1}(\bA)=0$, it is known as having \textit{exact rank} $r$ (also known as the \textit{ordinary rank}), as extensively employed in our prior discussions. 
More formally, the numerical rank $\nr(\bA, d)$ of a matrix $\bA$, given a nonnegative value $d\geq 0$, is defined as follows
\begin{equation}
\nr(\bA, d) = \min \{\rank(\bB) : \normtwo{\bA-\bB} \leq d\},
\end{equation}
with respect to the spectral norm (Definition~\ref{definition:spectral_norm}).
Alternatively, we can define the numerical rank $\nr(\bA, d)$ as the unique integer $k$ such that both of the following two conditions hold:
\begin{enumerate}
\item There exists a matrix $\bB$ of precise rank $k$ such that $\Vert\bA-\bB\Vert_2\leq d$;
\item There is no matrix $\bB$ with a rank less than $k$ such that (1) holds.
\end{enumerate}
\end{definition}

\index{Numerical rank}
\index{Numerical rank}
\index{Truncated}
\index{Truncated SVD}
\index{Eckart-Young-Mirsky theorem}
Given the nonnegative tolerance value $d$, the numerical rank of a matrix then can be determined by the following theorem. 
\begin{theorem}[Numerical Rank]\label{theorem:nmum_rank_d}
Let $\bA\in\real^{m\times n}$ be any matrix with nonzero singular values $\sigma_1\geq \sigma_2 \geq  \ldots \geq \sigma_r$. Then, given a nonnegative value $d\geq 0$, the numerical rank $\nr(\bA, d)=k <r$ if and only if $\sigma_k > d \geq \sigma_{k+1}$.
\end{theorem}
\begin{proof}[of Theorem~\ref{theorem:nmum_rank_d}]
The proof relies on the  \textit{Eckart-Young-Mirsky theorem} (Section~\ref{section:svd-low-rank-approxi}); feel free to skip for a first reading. 
Suppose $\sigma_k >d \geq \sigma_{k+1}$ and let $\bA_k$ be the truncated SVD of $\bA$ with the largest $k$ terms, i.e., $\bA_k = \sum_{i=1}^{k} \sigma_i\bu_i\bv_i^\top$ from the SVD of $\bA=\sum_{i=1}^{r} \sigma_i\bu_i\bv_i^\top$ by zeroing out the $r-k$ trailing singular values of $\bA$. 
By assumption, we have 
$$
\Vert\bA-\bA_{k-1}\Vert_2 >d \geq \Vert\bA-\bA_k\Vert_2.
$$
According to Theorem~\ref{theorem:young-theorem-spectral}, we also have 
$$
\Vert\bA-\bA_{k-1}\Vert_2 = \min \{\Vert\bA-\bB \Vert_2 : \rank(\bB)=k-1\} >d.
$$
This implies $\min\{\rank(\bB) : \Vert\bA-\bB\Vert_2 \leq d\}=k$ and $\nr(\bA,d)=k$.

For the converse implication, suppose $\nr(\bA,d)=k$. This means $\min\{\rank(\bB) : \Vert\bA-\bB\Vert_2 \leq d\}=k$. Therefore, $\sigma_k=\Vert\bA-\bA_{k-1}\Vert_2>d$.
Moreover, the assumption also indicates that $d\geq \Vert\bA-\bB\Vert_2$ if $\bB$ is of rank $k$. 
By Theorem~\ref{theorem:young-theorem-spectral}, we also have 
$$
\Vert\bA-\bA_{k}\Vert_2 = \min \{\Vert\bA-\bB \Vert_2 : \rank(\bB)=k\}.
$$
Therefore, $d\geq \Vert\bA-\bA_{k}\Vert_2 = \sigma_{k+1}$. This completes the proof.
\end{proof}

In addition to determining the numerical rank via SVD, it can also be determined using QR decomposition.

\begin{theorem}[Numerical Rank via QR]\label{theorem:nrank_qr}
Let $\bA\in\real^{m\times n}$ be any matrix with $m\geq n$, and let its  QR decomposition be given by $\bA=\bQ\bR$, where 
$\bR=
\footnotesize\begin{bmatrix}
\bR_{11} & \bR_{12}\\
\bzero & \bR_{22}
\end{bmatrix}$ and $\bR_{11}\in\real^{r\times r}$.
If $\sigma_{\min}(\bR_{11})\gg \normtwo{\bR_{22}} = \mathcalO(\mu)$, where $\mu$ is the machine precision, then $\bA$ has numerical rank $r$.
\end{theorem}
\begin{proof}[of Theorem~\ref{theorem:nrank_qr}]
Suppose $\bA$ has singular values $\sigma_1\geq\sigma_2\ldots \geq \sigma_n$.
Based on Remark~\ref{remark:qr_inter_svd}, we have $\sigma_r\geq\sigma_{r}(\bR_{11})=\sigma_{\min}(\bR_{11})$
and $\normtwo{\bR_{22}} =\sigma_{\max}(\bR_{22}) \geq \sigma_{r+1}$.
Therefore, we have 
$$
\sigma_r\geq\sigma_{\min}(\bR_{11}) \gg 
\normtwo{\bR_{22}} \geq \sigma_{r+1}= \mathcalO(\mu),
$$
from which the result follows.
\end{proof}

We have briefly introduced the notion of rank-revealing QR decomposition in Section~\ref{section:partial-cpqr-mgs}.
Its rigorous definition can be given using the numerical rank.
\begin{definition}[Rank-Revealing QR Decomposition]\label{definition:rrqr_in_svd}
Let $\bA\in\real^{m\times n}$ be any matrix with numerical rank $r$. If there exists a permutation matrix $\bP$ such that the QR decomposition of $\bA\bP$ is
$$
\bA\bP = 
\bQ\bR=
\bQ
\begin{bmatrix}
\bR_{11} & \bR_{12} \\
\bzero & \bR_{22}
\end{bmatrix},
$$
where $\bR_{11}\in\real^{r\times r}$, and $\sigma_{\min}(\bR_{11}) \gg \normtwo{\bR_{22}} = \mathcalO(\mu)$. Then the decomposition $\bA\bP = 
\bQ\bR$ is called a rank-revealing QR decomposition (RRQR) of $\bA$.
\end{definition}

\section{Properties of Singular Values}
Section~\ref{section:propert_eig_sym} describes the properties and eigenvalues of symmetric matrices. 
This section builds on that understanding to delve into the properties of singular values, exploring several important theorems and inequalities related to singular values, such as the Courant-Fischer Theorem,  Weyl's Theorem, and others like the von Neumann and interlacing inequalities. These theorems provide a framework for understanding the behavior of singular values under various operations and conditions, including matrix multiplication and addition.
\subsection{Courant-Fischer, Von Neumann, Interlacing, and other Inequalities}\label{section:svd_von_cou_etc}
\subsection*{Courant-Fischer and Weyl Theorem for Singular Values}
The Courant-Fischer Theorem~\ref{theorem:courant_fischer} describes how to obtain the eigenvalues of a symmetric matrix. Similarly, we can use min-max optimizations to obtain the singular values of a matrix.
\begin{theorem}[Courant-Fischer Theorem for Singular Values]\label{theorem:courant_fischer_sing}
Let $\bA\in \real^{m\times n}$ be any matrix with singular values $\sigma_1\geq  \sigma_2\geq \ldots \geq \sigma_r$, where $r$ is the rank of $\bA$.
Let $\sV_k$ denote the set of subspaces of $\real^n$ with dimension $k$  with $k\in\{1,2,\ldots,n\}$. Then,
$$
\begin{aligned}
\sigma_k &= 
\mathop{\min}_{\mathcalS \in \sV_{n-k+1}}
\mathop{\max}_{\bx\neq \bzero,\bx\in\mathcalS }
\frac{\normtwo{\bA\bx}}{\normtwo{\bx}}
\gap \text{and}\gap 
\sigma_k =
\mathop{\max}_{\mathcalS \in \sV_{k}}
\mathop{\min}_{\bx\neq \bzero,\bx\in\mathcalS }
\frac{\normtwo{\bA\bx}}{\normtwo{\bx}}.~\footnote{The result can be generalized to the complex  matrix $\bA$, and $\sV_k$ then denotes the set of subspace of $\complex^n$.}
\end{aligned}
$$
\end{theorem}
\begin{proof}[of Theorem~\ref{theorem:courant_fischer_sing}]
Suppose $\lambda_1\leq\lambda_2 \leq \ldots\leq \lambda_n$ are the ordered eigenvalues of the PSD matrix $\bA^\top\bA$. Then $\sigma_k^2 = \lambda_{n-k+1}$.
Therefore, 
$$
\begin{aligned}
\sigma_k^2=\lambda_{n-k+1} &= 
\mathop{\max}_{\mathcalS \in \sV_{k}}
\mathop{\min}_{\bx\neq \bzero,\bx\in\mathcalS }
\frac{\bx^\top\bA^\top\bA\bx}{\bx^\top\bx} 
=
\mathop{\min}_{\mathcalS \in \sV_{n-k+1}}
\mathop{\max}_{\bx\neq \bzero,\bx\in\mathcalS }
\frac{\bx^\top\bA^\top\bA\bx}{\bx^\top\bx}.
\end{aligned}
$$
This yields the result.
\end{proof}

Moreover, using a similar argument, for any matrix $\bA\in\real^{m\times n}$, the Rayleigh-Ritz Theorem~\ref{theorem:rayleigh_v2} for eigenvalues also shows
\begin{equation}\label{equation:rayleig_ritz_sing}
\begin{aligned}
&\mathop{\max}_{\bx\neq \bzero}
\frac{\normtwo{\bA\bx}}{\normtwo{\bx}}
=
\sigma_{\max}(\bA)
\qquad\text{and}\qquad
\mathop{\min}_{\bx\neq \bzero}
\frac{\normtwo{\bA\bx}}{\normtwo{\bx}}
=
\sigma_{\min}(\bA).
\end{aligned}
\end{equation}
This confirms the definition of the spectral norm of a matrix such that $\normtwo{\bA}=\sigma_{\max}(\bA)$ for any matrix $\bA$  and the matrix-vector inequality such that $\normtwo{\bA\bx}\leq \sigma_{\max}(\bA)\normtwo{\bx}$ for any matrix $\bA$ and vector $\bx$ (Definition~\ref{definition:spectral_norm}).

Alternatively, \eqref{equation:rayleig_ritz_sing} can be proved using the invariance of the $\ell_2$ vector norm under orthogonal transformations.
\begin{theorem}\label{theorem:svd_opt_app}
Let $\bA = \bU\bSigma\bV^\top\in \real^{m\times n}$ be the SVD of matrix $\bA$ with singular values $\sigma_1\geq \sigma_2\geq \ldots\geq \sigma_q$, where $q=\min\{m,n\}$. Then, the problem 
$$
\begin{aligned}
&\gap \min \normtwo{\bA\bx} \quad \text{s.t. }\quad \normtwo{\bx}=1
\end{aligned}
$$
has the solution $\bx=\bv_q$ (the $q$-th column of $\bV$) and the value of the minimum is 
$$
\mathop{\min}_{\normtwo{\bx}=1} \normtwo{\bA\bx} = \sigma_q.
$$
Similarly, $
\begin{aligned}
	\max \normtwo{\bA\bx}  \text{ s.t. } \normtwo{\bx}=1
\end{aligned}
$ has the solution $\bx=\bv_1$   (the first column of $\bV$), and the value of the maximum is 
$
\mathop{\max}_{\normtwo{\bx}=1} \normtwo{\bA\bx} = \sigma_{1}.
$
\end{theorem}

\begin{proof}[of Theorem~\ref{theorem:svd_opt_app}]
Using the invariance of the $\ell_2$ vector norm under orthogonal transformations, we have 
$$
\begin{aligned}
&\mathop{\min}_{\normtwo{\bx}=1} \normtwo{\bA\bx}^2 = \mathop{\min}_{\normtwo{\bV^\top\bx}=1} \normtwo{\bU\bSigma\bV^\top \bx}^2
\xlongequal{\by=\bV^\top\bx}
\mathop{\min}_{\normtwo{\by}=1} \normtwo{\bU\bSigma\by}^2\\
&=\mathop{\min}_{\normtwo{\by}=1} \normtwo{\bSigma\by}^2
=\mathop{\min}_{\normtwo{\by}=1} (\sigma_1^2y_1^2 + \sigma_2^2y_2^2 +\ldots +\sigma_q^2y_q^2 )
\geq \sigma_q^2.
\end{aligned}
$$
Therefore, the minimum is achieved when $\by=\be_q$ and $\bx=\bV\by=\bv_q$.
The second part can be proved similarly.
\end{proof}

The Courant-Fischer theorem for singular values can be used to bound the singular values of a product of matrices. 
Let $\sigma_k(\cdot)$ denote the $k$-th largest singular value of a matrix.
For any $\bA,\bB\in\real^{n\times n}$, we want to show that 
\begin{equation}
	\sigma_k(\bA\bB) \leq \sigma_1(\bA)\sigma_k(\bB).
\end{equation}
To see this, \eqref{equation:rayleig_ritz_sing} shows $\normtwo{\bA\bB\bx}\leq \sigma_1(\bA)\normtwo{\bB\bx}$.
By the Courant-Fischer theorem for $ \sigma_k(\bB)$, we have:
$
\sigma_k(\bB) = 
\mathop{\max}_{\mathcalS \in \sV_{k}}
\mathop{\min}_{\bx\neq \bzero,\bx\in\mathcalS }
\frac{\normtwo{\bB\bx}}{\normtwo{\bx}}.
$
Thus, for any subspace $\mathcalV$ of dimension $k$:
$$
\min_{\substack{\bx\neq\bzero, \bx \in \mathcalV }} \frac{\normtwo{\bB\bx}}{\normtwo{\bx}} \leq \sigma_k(\bB).
$$
Using the Courant-Fischer theorem for $\sigma_k(\bA\bB)$, we obtain the result:
$$
\begin{aligned}
	\sigma_k(\bA\bB) 
	&= 
	\mathop{\max}_{\mathcalS \in \sV_{k}}
	\mathop{\min}_{\bx\neq \bzero,\bx\in\mathcalS }
	\frac{\normtwo{\bA\bB\bx}}{\normtwo{\bx}}
	\leq 
	\mathop{\max}_{\mathcalS \in \sV_{k}}
	\mathop{\min}_{\bx\neq \bzero,\bx\in\mathcalS }
	\frac{\sigma_1(\bA)\normtwo{\bB\bx}}{\normtwo{\bx}}\\
	&\leq  
	\mathop{\max}_{\mathcalS \in \sV_{k}}
	\mathop{\min}_{\bx\neq \bzero,\bx\in\mathcalS }
	\sigma_1(\bA)\sigma_k(\bB)
	=\sigma_1(\bA)\sigma_k(\bB).
\end{aligned}
$$

Theorem~\ref{theorem:courant_fischer_sing} provides an max-min optimization method for finding the singular values of a matrix. 
Nevertheless, if the right singular vectors are given, the optimization method can be simplified to a single maximization.
\begin{theorem}[Singular Value Optimization]\label{theorem:sig_va_op2}
Let $\bA\in\real^{m\times n}$ with ordered singular values $\sigma_1\geq \sigma_2\geq \ldots\geq\sigma_q$ and $q=\min\{m,n\}$, and let $\bA=\bU\bSigma\bV^\top$ be the SVD of $\bA$ where $\bV=[\bv_1, \bv_2, \ldots,\bv_n]$ are the right singular vectors of $\bA$. Let further $\mathcalV=\spn\{\bv_k, \ldots,\bv_n\}$.
Then $\max \{\normtwo{\bA\bx}: \bx\in\mathcalV, \normtwo{\bx}=1\} = \sigma_k$, for all $k\in\{1,2,\ldots,q\}$.
\end{theorem}
\begin{proof}[of Theorem~\ref{theorem:sig_va_op2}]
Let $\bx\triangleq a_k\bv_k +a_{k+1}\bv_{k+1}+\ldots+a_n\bv_n$ with $\sum_{i=k}^{n} a_i^2=1$.
Then, $\normtwo{\bA\bx}^2 = \normtwo{\bU\bSigma\bV^\top\bx}^2 = \sum_{i=k}^{q} a_i^2 \sigma_i^2\leq \sigma_k^2$, and the equality is attained at $\bx=\bv_k$.
\end{proof}

\index{Weyl's theorem}
Theorem~\ref{theorem:weyl} provides the relationship between the eigenvalues of a sum of two matrices based on the Courant-Fisher theorem for eigenvalues. Similar results can be observed for singular values.
\begin{theorem}[Weyl's Theorem for Singular Values]\label{theorem:weyl_sing}
Let $\bA, \bB \in \real^{m\times n}$ be given and let $q = \min\{m,n\}$. The following inequalities hold for the nonincreasing ordered singular values of $\bA$, $\bB$, $\bA + \bB$, and $\bA\bB^\top$:
\begin{align}
\sigma_{i+j-1}(\bA + \bB) &\leq \sigma_i(\bA) + \sigma_j(\bB); \label{equation:weyl_sin1}  \\
\sigma_{i+j-1}(\bA\bB^\top) &\leq \sigma_i(\bA) \sigma_j(\bB), \label{equation:weyl_sin2}
\end{align}
for $1 \leq i, j \leq q$ and $i + j \leq q + 1$. In particular,
\begin{align}
|\sigma_i(\bA + \bB) - \sigma_i(\bA)| &\leq \sigma_1(\bB) \quad \text{for } i \in\{1, \ldots, q\}; \label{equation:weyl_sin3}\\
\sigma_i(\bA\bB^\top) &\leq \sigma_i(\bA) \sigma_1(\bB) \quad \text{for } i \in\{1, \ldots, q\}. \label{equation:weyl_sin4}
\end{align}

\end{theorem}
\begin{proof}[of Theorem~\ref{theorem:weyl_sing}]
The proof is adapted from \citet{horn1994topics}.
Let $\bA = \bU \bSigma_A \bV^\top$ and $\bB = \bW \bSigma_B \bZ^\top$ be singular value decompositions of $\bA$ and $\bB$ with orthogonal $\bU = [\bu_1,\bu_2, \ldots, \bu_m]$, $\bW = [\bw_1,\bw_2, \ldots, \bw_m] \in \real^m$ and orthogonal $\bV = [\bv_1,\bv_2, \ldots, \bv_n]$, $\bZ = [\bz_1,\bz_2, \ldots, \bz_n] \in \real^n$. Let $i$ and $j$ be positive integers with $1 \leq i, j \leq q$ and $i + j \leq q + 1$.

First, consider the sum inequalities \eqref{equation:weyl_sin1}. Let 
$$
\mathcalV \triangleq\spn\{\bv_{i}, \ldots, \bv_n\}
\quad \text{and} \quad
\mathcalW \triangleq \spn\{\bz_{j}, \ldots, \bz_n\}.
$$
Note that $\dim(\mathcalV) = n-i+1$ and $\dim(\mathcalW) = n-j+1$. Then,
\begin{align*}
\alpha &\triangleq \dim(\mathcalV \cap \mathcalW) = \dim(\mathcalV) + \dim (\mathcalW) - \dim(\mathcalV + \mathcalW) \\
&> (n-i+1) + (n-j+1) - n = n - (i+j-1) + 1 \geq 1
\end{align*}
because of the bounds assumed for $i$ and $j$. Thus, the subspace $\mathcalV \cap \mathcalW$ has positive dimension $\alpha$, $n-\alpha+1 \leq i+j-1$, and we can use the Courant-Fischer Theorem~\ref{theorem:courant_fischer_sing} for singular values and Theorem~\ref{theorem:sig_va_op2} to compute
$$
\begin{aligned}
\sigma_{i+j-1}&(\bA + \bB) \leq \sigma_{n-\alpha+1}(\bA + \bB)
= \min_{\substack{\mathcalS \subset \real^n \\ \dim(\mathcalS) = \alpha}} \max_{\substack{\bx \in \mathcalS \\ \normtwo{\bx} = 1}} \normtwo{(\bA + \bB)\bx}
\leq \max_{\substack{\bx \in \mathcalV \cap \mathcalW \\ \normtwo{\bx} = 1}} \normtwo{(\bA + \bB)\bx}\\
&\leq \max_{\substack{\bx \in \mathcalV \cap \mathcalW \\ \normtwo{\bx} = 1}} (\normtwo{\bA \bx} + \normtwo{\bB \bx})
\leq \max_{\substack{\bx \in \mathcalV \\ \normtwo{\bx} = 1}} \normtwo{\bA \bx} + \max_{\substack{\bx \in \mathcalW \\ \normtwo{\bx} = 1}} \normtwo{\bB \bx} = \sigma_i(\bA) + \sigma_j(\bB).
\end{aligned}
$$

Next, consider the product inequalities \eqref{equation:weyl_sin2}. Let the left polar decomposition (Theorem~\ref{theorem:polar-decomposition}) of $\bA\bB^\top$ be $\bA \bB^\top = \bQ_l\bS_l$ (the case for $m\leq n$ can be proved similarly using right polar decomposition), where $\bQ_l \in \real^{m\times m}$ is orthogonal and $\bS_l \in \real^{m\times m}$ is positive semidefinite and has the same singular values (which are also its eigenvalues) as $\bA \bB^\top$. Let
$$
\mathcalV \triangleq \spn\{\bQ_l^\top \bu_i, \ldots, \bQ_l^\top \bu_n\} \quad \text{and} \quad \mathcalW \triangleq \spn\{\bz_j, \ldots, \bz_n\}.
$$
so $\alpha \triangleq \dim(\mathcalV \cap \mathcalW) \geq n - (i + j - 1) + 1 \geq 1$. Since $\bS_l = \bQ_l^\top\bA \bB^\top$, we have
$$
\bx^\top \bS_l \bx = \bx^\top \bQ_l^\top \bA \bB^\top \bx = (\bA^\top \bQ_l \bx)^\top(\bB^\top \bx) \leq \|\bA^\top \bQ_l\bx\|_2 \|\bB^\top \bx\|_2 \quad \text{for any } \bx \in \real^n.
$$
Therefore, we can use the Courant-Fischer Theorem~\ref{theorem:courant_fischer_sing} for singular values, the Courant-Fischer Theorem~\ref{theorem:courant_fischer}, and Theorem~\ref{theorem:sig_va_op2}  to compute
$$
\begin{aligned}
\sigma_{i+j-1}&(\bA\bB^\top) = \sigma_{i+j-1}(\bS_l) \leq \sigma_{n-\alpha+1}(\bS_l)
= \min_{\substack{\mathcalS \subset \real^n \\ \dim(\mathcalS) = \alpha}} \max_{\substack{\bx \in \mathcalS \\ \normtwo{\bx} = 1}} \bx^\top\bS_l\bx\\
&\leq \max_{\substack{\bx \in \mathcalV \cap \mathcalW \\ \normtwo{\bx} = 1}} \bx^\top \bS_l \bx
\leq \max_{\substack{\bx \in \mathcalV \cap \mathcalW \\ \normtwo{\bx} = 1}} (\|\bA^\top \bQ_l\bx\|_2 \|\bB^\top \bx\|_2)
\leq \max_{\substack{\bx \in \mathcalV \cap \mathcalW \\ \normtwo{\bx} = 1}} \|\bA^\top \bQ_l\bx\|_2 \max_{\substack{\bx \in \mathcalV \cap \mathcalW \\ \normtwo{\bx} = 1}} \|\bB^\top \bx\|_2\\
&\leq \max_{\substack{\bx \in \mathcalV \\ \normtwo{\bx} = 1}} \|\bA^\top \bQ_l \bx\|_2 \max_{\substack{\bx \in \mathcalW \\ \normtwo{\bx} = 1}} \|\bB^\top \bx\|_2 = \sigma_i(\bA) \sigma_j(\bB).
\end{aligned}
$$
Notice that the nondecreasing ordered eigenvalues $\lambda_1(\bS_l) \leq \cdots \leq \lambda_n(\bS_l)$ are related to the nonincreasing ordered singular values $\sigma_1(\bS_l) \geq \cdots \geq \sigma_n(\bS_l)$ by $\lambda_i(\bS_l) = \sigma_{n-i+1}(\bS_l)$.

The inequalities $\sigma_i(\bA + \bB) \leq \sigma_i(\bA) + \sigma_1(\bB)$ and $\sigma_i(\bA \bB^\top) \leq \sigma_i(\bA) \sigma_1(\bB)$ follow from setting $j = 1$ in \eqref{equation:weyl_sin1} and \eqref{equation:weyl_sin2}. The two-sided bound in the additive case now follows from observing that $\sigma_i(\bA) = \sigma_i((\bA + \bB) - \bB) \leq \sigma_i(\bA + \bB) + \sigma_1(-\bB) = \sigma_i(\bA + \bB) + \sigma_1(\bB)$.
\end{proof}

\subsection*{Von Neumann}
\textit{Von Neumann's trace theorem}, also known as the \textit{von Neumann trace inequality}, is a fundamental result in matrix theory and functional analysis. It establishes  an important relationship between the traces of the product of two matrices and the singular values of those matrices.
\begin{theorem}[Von Neumann's Trace Theorem]\label{theorem:vonneumann1}
Let $\bA,\bB\in\real^{m\times n}$ be given with  ordered singular values $\sigma_1(\bA) \geq \ldots \geq \sigma_p(\bA)$ and $\sigma_1(\bB) \geq \ldots \geq \sigma_p(\bB)$, respectively, where $p=\min\{m,n\}$. Then, 
$$
\trace(\bA\bB^\top) \leq \sum_{i=1}^{p} \sigma_i(\bA) \sigma_i(\bB).~\footnote{When $\bA,\bB\in\complex^{n\times n}$, $\trace(\bA\bB)$ becomes $\operatorname{Re} \trace(\bA\bB)$.}
$$
When $m=n$, we also have $\trace(\bA\bB) \leq \sum_{i=1}^{n} \sigma_i(\bA) \sigma_i(\bB)$.
\end{theorem}
\begin{proof}[of Theorem~\ref{theorem:vonneumann1}]
We first discuss the case $m=n$.
Let $ \bA = \bU_1 \bSigma_1 \bV_1^\top $ and $ \bB^\top = \bU_2 \bSigma_2 \bV_2^\top $ be the SVDs of $\bA$ and $\bB^\top$, respectively, where $ \bU_1, \bV_1, \bU_2, \bV_2 \in \real^{n\times n}$ are orthogonal, $\bSigma_1 = \diag(\sigma_1(\bA), \ldots, \sigma_n(\bA))$, and $\bSigma_2 = \diag(\sigma_1(\bB), \ldots, \sigma_n(\bB))$. Let $ \bQ \triangleq \bV_1^\top \bU_2$ and $ \bP \triangleq \bV_2^\top \bU_1$. Then,
$$
\begin{aligned}
	\trace(\bA\bB^\top) 
	&=  \trace(\bU_1 \bSigma_1 \bV_1^\top \bU_2 \bSigma_2 \bV_2^\top) 
	=  \trace(\bSigma_1 \bV_1^\top \bU_2 \bSigma_2 \bV_2^\top \bU_1) =  \trace(\bSigma_1 \bQ \bSigma_2 \bP) \\
	&=  \sum_{i, j=1}^{n,n} \sigma_i(\bA) \sigma_j(\bB) q_{ij} p_{ji} 
	\leq \sum_{i, j=1}^{n,n} \sigma_i(\bA) \sigma_j(\bB) \abs{q_{ij} p_{ji}}.
\end{aligned}
$$
The matrix $\bZ$ with $(i,j)$-th element given by $z_{ij}=\abs{q_{ij} p_{ji}}$ is doubly substochastic (Exercise~\ref{exercise:dou_substo_int}), and there exists a doubly stochastic matrix $\bM$ such that $\bM-\bZ$ is nonnegative (Exercise~\ref{exercise:dou_substo_int2}). This results in
$$
\begin{aligned}
	\trace(\bA\bB^\top) &\leq \sum_{i, j=1}^{n,n} \sigma_i(\bA) \sigma_j(\bB) m_{ij} 
	\leq \max \left\{ \sum_{i, j=1}^{n,n} \sigma_i(\bA) \sigma_j(\bB) s_{ij} : \bS\text{ is doubly stochastic} \right\}.
\end{aligned}
$$
The function $f(\bS) = \sum_{i, j=1}^{n, n} \sigma_{i}(\bA) \sigma_{j}(\bB) s_{ij} $ is a linear (and a convex) function on the set of doubly stochastic matrices, and  it attains its maximum at a permutation matrix. Let $I$ be a permutation of $\{1, 2, \ldots, n\}$ such that $s_{ij} = 1 $ if and only if $ j = I(i) $, we have
$$
\begin{aligned}
	\trace(\bA\bB^\top) &\leq \sum_{i, j=1}^{n,n} \sigma_i(\bA) \sigma_j(\bB) s_{ij} = \sum_{i=1}^n \sigma_i(\bA) \sigma_{I(i)}(\bB) 
	\leq \sum_{i=1}^n \sigma_i(\bA) \sigma_i(\bB).
\end{aligned}
$$
The above result can be applied to $\bB$ rather than on $\bB^\top$. Therefore, we also have $\trace(\bA\bB) \leq \sum_{i=1}^{n} \sigma_i(\bA) \sigma_i(\bB)$.

When $m>n$, we can extend $\bA$ and $\bB$ into $\widetildebA=[\bA,\bzero]\in\real^{m\times m}$ and $\widetildebB=[\bB,\bzero]\in\real^{m\times m}$, respectively. 
Thus, $\bA\bB^\top=\widetildebA\widetildebB^\top$. The above discussion shows $\tr(\widetildebA\widetildebB^\top)\leq \sum_{i=1}^m \sigma_i(\widetildebA) \sigma_{i}(\widetildebB) =\sum_{i=1}^n \sigma_i(\bA) \sigma_{i}(\bB)$.
The case of $m<n$ can be proved analogously.
\end{proof}

Let $\bB=\bI$ in von Neumann's trace theorem, we get $\trace(\bA) \leq \sum_{i=1}^{p} \sigma_i(\bA)$. The following theorem addresses the conditions under which this inequality holds with equality.
\begin{theorem}[Equality in Von Neumann's Theorem]\label{theorem:equa_vonneu}
Let $ \bA\in \real^{m\times n} $ with ordered singular values $\sigma_1\geq \sigma_2 \geq \ldots \geq \sigma_p$, where $ p = \min\{m, n\} $ and $ q = \max\{m, n\} $. Let further $I = \{1, \ldots, p\}$. Then, $ \trace(\bA) \leq \sum_{i=1}^{p} \sigma_i$, with equality if and only if the leading principal submatrix $ \bA[I] $ is positive semidefinite and $ \bA $ has no nonzero entries outside this leading principal submatrix.
\end{theorem}
\begin{proof}[of Theorem~\ref{theorem:equa_vonneu}]
To prove the sufficiency of this equality, observe that if the leading principal submatrix $ \bA[I] $ is positive semidefinite and other entries of $ \bA $ are all zero, then its eigenvalues are its singular values, which are also the singular values of $ \bA $; the trace of $ \bA[I] $ is the sum of its eigenvalues, which is the sum of the singular values of $ \bA $.

Conversely, suppose that $\sum_{i=1}^{p} a_{ii} = \sum_{i=1}^{p} \sigma_i$. If $ \bA = \bzero $ there is nothing to prove. Otherwise, let $\rank(\bA) = r \geq 1$. Complete $ \bA $ with zero blocks into a square matrix
$
\widetildebA = 
\scriptsize\begin{bmatrix}
\bA & \bzero_{m, q-n} \\
\bzero_{q-m, n} & \bzero_{q-m, q-n}
\end{bmatrix} 
\normalsize\in \real^{q\times q},
$
which has the same trace and singular values as $ \bA $. Let $\widetildebA = \bU \bSigma_r \bV^\top$ be a reduced SVD of $\widetildebA$, where
$
\bU = 
[\bu_1,\bu_2,\ldots, \bu_r] \in \real^{q\times r} 
\quad \text{and} \quad 
\bV = 
[\bv_1, \bv_2,\ldots,\bv_r]\in \real^{q\times r}
$
are the column partitions of $\bU$ and $\bV$, respectively, and $\bSigma_r = \diag(\sigma_1, \sigma_2,\ldots, \sigma_r)$. Then, to obtain the equality, we need to have
$$
\begin{aligned}
\trace(\widetildebA) &= \sum_{i=1}^{p} a_{ii} = \sum_{i=1}^{q} \sum_{j=1}^{r} u_{ij} \sigma_j {v}_{ij} 
= \sum_{j=1}^{r} \sigma_j \sum_{i=1}^{q} u_{ij} {v}_{ij} 
= \sum_{j=1}^{r} \sigma_j (\bv_j^\top \bu_j) = \sum_{j=1}^{r} \sigma_j = \sum_{j=1}^{p} \sigma_j.
\end{aligned}
$$
That is, we need to ensure that $(\bv_j^\top \bu_j) = 1$ for each $ j\in\{1,2, \ldots, r \}$ to obtain the equality. Since
$$
1 = (\bv_j^\top \bu_j) \stackrel{(*)}{\leq} |\bv_j^\top \bu_j| \stackrel{(\dag)}{\leq} \|\bu_j\|_2^2 \|\bv_j\|_2^2 = 1,
$$
the equality at $(\dag)$ and the equality case of the Cauchy-Schwarz inequality ensure that there are scalars $ \gamma_j $ such that $ \bu_j = \gamma_j \bv_j $ for each $ j\in\{1,2, \ldots, r \}$; the equality at $(*)$ ensures that each $ \gamma_j = 1 $. Therefore, $ \bU = \bV $ and $\widetildebA = \bU \bSigma_r \bU^\top$ is positive semidefinite. It follows that its leading principal submatrix $ \bA[I] $ is positive semidefinite, and no other entries of $\widetildebA$ (and hence of $ \bA $) are nonzero. 
\end{proof}

Exercise~\ref{exercise:perturb_sings} shows that the difference between two corresponding singular values of two matrices is bounded by the spectral norm of the difference between the two matrices: $\normtwo{\bSigma_A-\bSigma_B} \leq \normtwo{\bA-\bB}$, where $\bSigma_A$ and $\bSigma_B$ are diagonal matrices containing singular values of $\bA$ and $\bB$, respectively.
The von Neumann theorem can be applied to prove, for example, the following inequalities that are frequently used in matrix approximation problems.
The proof relies on the existence of the polar decomposition (Theorem~\ref{theorem:polar-decomposition}).
Theorem~\ref{theorem:weyl_ineq2_sings} based on Weyl's inequality will discuss a more general case for these inequalities.
\begin{theorem}[Singular Value Inequalities]\label{theorem:vonneu_ineq2_sings}
Let $\bA,\bB\in\real^{m\times n}$ be given with  ordered singular values $\sigma_1(\bA) \geq \ldots \geq \sigma_p(\bA)$ and $\sigma_1(\bB) \geq \ldots \geq \sigma_p(\bB)$, respectively, where $p=\min\{m,n\}$. 
Let further  $\sigma_1(\bA\bB^\top) \geq \ldots \geq \sigma_p(\bA\bB^\top)$ be the ordered singular values of $\bA\bB^\top$.
Then, it follows that 
\begin{enumerate}
\item \textbf{Frobenius bounding.} $\normf{\bA-\bB}^2 \geq \sum_{i=1}^{p}(\sigma_i(\bA)-\sigma_i(\bB))^2=\normf{\bSigma_A-\bSigma_B}^2$~\footnote{It can actually be shown that this inequality holds for any unitarily invariant norm, a category that includes both the spectral and Frobenius norms: $\norm{\bA-\bB}^2 \geq \norm{\bSigma_A-\bSigma_B}^2$ if $\norm{\cdot}$ is unitarily invariant \citep{mirsky1960symmetric}.}, where $\bSigma_A$ and $\bSigma_B$ are diagonal matrices containing singular values of $\bA$ and $\bB$, respectively:\begin{itemize}
\item  equality is attained if and only if $\trace(\bA\bB^\top)=\sum_{i=1}^{p}\sigma_i(\bA)\sigma_i(\bB)$;
\item if the equality is attained, then both $\bA\bB^\top$ and $\bB^\top\bA$ are PSD such that $\trace(\bA\bB^\top)$ is nonnegative.
\end{itemize}
\item $\sum_{i=1}^{p} \sigma_i(\bA\bB^\top)\leq \sum_{i=1}^{p}\sigma_i(\bA)\sigma_i(\bB)$.
\item $\sum_{i=1}^{p} \sigma_i(\bA)\sigma_i(\bB)=\max\{\trace(\bA\bP\bB^\top\bQ): \bP\in\real^{n\times n} \text{ and }\bQ\in\real^{m\times m} \text{ are orthogonal}\}$.
\item $\sum_{i=1}^{p} \sigma_i(\bA)=
\left\{
\begin{aligned}
\max_{\text{orthogonal }\bQ} \trace(\bA\bQ), \gap \text{if $m\leq n$};\\
\max_{\text{orthogonal }\bQ} \trace(\bQ\bA), \gap \text{if $m\geq n$}.
\end{aligned}
\right.
$
\end{enumerate}
\end{theorem}
\begin{proof}[of Theorem~\ref{theorem:vonneu_ineq2_sings}]
For (1), use Theorem~\ref{theorem:vonneumann1}, we have 
$
\normf{\bA-\bB}^2
=\normf{\bA}^2+\normf{\bB}^2-2\trace(\bA\bB^\top)
\geq \normf{\bA}^2+\normf{\bB}^2-2 \sum_{i=1}^{p}\sigma_i(\bA)\sigma_i(\bB) = \sum_{i=1}^{p}(\sigma_i(\bA)-\sigma_i(\bB))^2.
$
If the equality is attained, then Theorem~\ref{theorem:vonneumann1} and (2) indicate that
$$
\trace(\bA\bB^\top) \leq \sum_{i=1}^{p} \sigma_i(\bA\bB^\top)
\leq \sum_{i=1}^{p} \sigma_i(\bA) \sigma_i(\bB)
=
\trace(\bA\bB^\top).
$$
Therefore, $\trace(\bA\bB^\top) = \sum_{i=1}^{p} \sigma_i(\bA\bB^\top)$, where Theorem~\ref{theorem:equa_vonneu} shows that $\bA\bB^\top$ is PSD. Applying this process to $\normf{\bB^\top-\bA^\top}^2$ shows that $\bB^\top\bA$ is also PSD.

For (2), using the right polar decomposition of $\bA\bB^\top=\bS_r\bQ_r$, where $\bS_r$ is PSD and $\bQ_r$ is orthogonal (see Theorem~\ref{theorem:polar-decomposition}, the case for $m>n$ can be proved similarly using the left polar decomposition), we have:
$$
\sum_{i=1}^{p} \sigma_i(\bA\bB^\top) = \trace(\bS_r)
=
\trace(\bA\bB^\top\bQ_r^\top)
\leq \sum_{i=1}^{p} \sigma_i(\bA)\sigma_i(\bQ_r\bB)= \sum_{i=1}^{p} \sigma_i(\bA)\sigma_i(\bB),
$$
where the inequality follows from Theorem~\ref{theorem:vonneumann1}.

For (3),  Theorem~\ref{theorem:vonneumann1} shows 
$
\trace(\bA\bP\bB^\top\bQ)
\leq \sum_{i=1}^{p} \sigma_i(\bA\bP)\sigma_i(\bQ^\top\bB)
=\sum_{i=1}^{p} \sigma_i(\bA)\sigma_i(\bB)
$.
Let $ \bA = \bU_1 \bSigma_1 \bV_1^\top $ and $ \bB = \bU_2 \bSigma_2 \bV_2^\top $ be the SVDs of $\bA$ and $\bB$, respectively, where $ \bU_1, \bU_2\in\real^{m\times m}$, $\bV_1, \bV_2 \in \real^{n\times n}$ are orthogonal, $\bSigma_1 = \diag(\sigma_1(\bA), \ldots, \sigma_p(\bA), 0,\ldots)$, and $\bSigma_2 = \diag(\sigma_1(\bB), \ldots, \sigma_p(\bB), 0,\ldots)$.
Then, the upper bound can be achieved at $\bP=\bV_1\bV_2^\top$ and $\bQ=\bU_2\bU_1^\top$.

For (4), when $m\geq n$, complete $\bA$ into $\widetildebA=[\bA,\bzero]\in\real^{m\times m}$. Theorem~\ref{theorem:vonneumann1} shows 
$
\trace(\bQ\bA)=\trace(\bQ\widetildebA)
\leq \sum_{i=1}^{m} \sigma_i(\bQ)\sigma_i(\widetildebA^\top)
=
\sum_{i=1}^{m} \sigma_i(\widetildebA)
=
\sum_{i=1}^{p} \sigma_i(\bA),
$
where we use the fact that the singular values of an orthogonal matrix are all one.
Let $\bA=\bQ_l\bS_l$ be the left polar decomposition of $\bA$, the upper bound is achieved at $\trace(\bQ_l^\top\bA)=\trace(\bS_l)=\sum_{i=1}^{p} \sigma_i(\bA)$.
The case of $m\leq n$ can be proved similarly by completing $\bA$ into $\scriptsize\begin{bmatrix}
\bA\\
\bzero 
\end{bmatrix}\in \real^{n\times n}$.
\end{proof}
The Frobenius bounding inequality in Theorem~\ref{theorem:vonneu_ineq2_sings} is often used to derive the truncated SVD (TSVD) that approximates a matrix with a low-rank matrix; see Section~\ref{section:svd-low-rank-approxi} for a detailed discussion.
\index{Truncated}
\index{Truncated SVD}

\begin{exercise}[Singular Value Inequalities]
Let $\bA \in \real^{m \times n}$ and $\bB \in \mathbb{F}^{n \times l}$. If $r \geq 0$, show that
$$
\sum_{i=1}^{k} \sigma_{i}^{r}(\bA\bB) \leq \sum_{i=1}^{k} \sigma_{i}^{r}(\bA) \sigma_{i}^{r}(\bB).
\quad k\in\{1, \ldots, \min\{ m, n,l\}\},
$$
If $r < 0$, $n = m = l$, and $\bA$ and $\bB$ are nonsingular, show that
$$
\sum_{i=1}^{n} \sigma_{i}^{r}(\bA\bB) \leq \sum_{i=1}^{n} \sigma_{i}^{r}(\bA) \sigma_{i}^{r}(\bB).
$$
\end{exercise}

\subsection*{Orthogonal Approximation Problems}
Section~\ref{section:polar_decom} discussed the closest orthogonal and rotation matrices  to a (square) matrix in terms of the Frobenius norm using polar decomposition. 
The Frobenius bounding inequality can be applied to determine when a scalar multiple of an orthogonal matrix is the closest to a given matrix.
\begin{theorem}[Closest (Scalar Multiple) Orthogonal]\label{theorem:closes_scal_ortho}
Let $\bA\in\real^{n\times n}$ be given with ordered singular values $\sigma_1\geq \sigma_2 \geq \ldots \geq \sigma_n$, and  let $\bQ\in\real^{n\times n}$ be an  orthogonal matrix.
Then, $\normf{\bA-\gamma\bQ}^2$ is minimized when 
$$
\gamma =\pm \frac{1}{n}\sum_{i=1}^{n}\sigma_i\triangleq\beta,
\quad\text{s.t.}\quad
\normf{\bA-\gamma\bQ}^2 \geq \sum_{i=1}^{n}(\sigma_i-\beta)^2.
$$
The equality is attained at $\bQ=\bQ_r$, where $\bS_r\bQ_r$ is the right polar decomposition of $\bA$.
\end{theorem}
\begin{proof}[of Theorem~\ref{theorem:closes_scal_ortho}]
The Frobenius bounding inequality shows that 
$$
\normf{\bA-\gamma\bQ}^2 \geq \sum_{i=1}^{n}(\sigma_i-\sigma_i(\gamma\bQ))^2
=
\sum_{i=1}^{n}\sigma_i^2- 2\abs{\gamma}\sum_{i=1}^{n} \sigma_i + n\abs{\gamma}^2,
$$
which is minimized at $\gamma =\pm \frac{1}{n}\sum_{i=1}^{n}\sigma_i$.
To see the equality, let $\bA=\bS_r\bQ_r$ be the right polar decomposition of $\bA$, where $\bS_r$ is PSD and $\bQ_r$ is orthogonal.
Since $\trace(\bS_r)=\sum_{i}^n \sigma_i =n\beta$ (Exercise~\ref{exercise:trace_polar}), we have  
$$
\normf{\bA-\beta\bQ_r}^2
=
\normf{\bS_r\bQ_r - \beta\bQ_r}^2
=
\normf{\bS_r-\beta\bI}^2
=
\trace(\bS_r^2) - 2\beta\trace(\bS_r)+n\beta^2
=
\normf{\bA}^2 - n\beta^2.
$$
Therefore, $\beta\bQ_r$ is the best approximation to $\bA$ by a scalar multiple of an orthogonal matrix.
\end{proof}

Given two matrices \( \bA \) and \( \bB \) of the same dimensions \( m \times n \), the \textit{unitary Procrustes problem} is to find a unitary matrix \( \bQ \) that minimizes the Frobenius norm of the difference between \( \bA \) and \( \bQ\bB \):
$
\min_{\bQ} \normf{\bA - \bQ\bB}^2,
$
where \( \bQ \) is an \( m \times m \) unitary matrix, i.e., \( \bQ^* \bQ = \bQ \bQ^* = \bI \), and \( \bQ^* \) denotes the conjugate transpose of \( \bQ \).
For simplicity, we consider orthogonal $\bQ$, though the term ``unitary Procrustes problem" is  more commonly used.
\begin{theorem}[Unitary/Orthogonal Procrustes]\label{theorem:ortho_procruste}
Let $\bA,\bB\in\real^{ m \times n }$  be given, and let $\bQ\in\real^{m\times m}$ be an  orthogonal matrix. Then,
$$
\textbf{(Orthogonal $\bQ$): }\quad \normf{\bA-\bQ\bB}^2 
\geq \normf{\bA}^2 - 2\sum_{i=1}^{m}\sigma_i(\bA\bB^\top) + \normf{\bB}^2.
$$ 
The equality is attained at $\bQ=\bQ_r$, where $\bS_r\bQ_r$ is the right polar decomposition of $\bA\bB^\top$.
Similarly, let further $\bP\in\real^{n\times n}$ be an orthogonal matrix. Then,
$$
\textbf{(Orthogonal $\bP$): }\quad \normf{\bA-\bB\bP}^2 
\geq \normf{\bA}^2 - 2\sum_{i=1}^{m}\sigma_i(\bA^\top\bB) + \normf{\bB}^2.
$$
The equality is attained at $\bP^\top=\widetilde{\bQ_r}$, where $\widetilde{\bS_r}\widetilde{\bQ_r}$ is the right polar decomposition of $\bB\bA^\top$.
\end{theorem}
\begin{proof}[of Theorem~\ref{theorem:ortho_procruste}]
Using von Neumann's trace Theorem~\ref{theorem:vonneumann1}, we have
$$
\normf{\bA-\bQ\bB}^2 = \normf{\bA}^2 - 2\trace(\bA\bB^\top\bQ^\top) + \normf{\bB}^2
\geq \normf{\bA}^2 - 2\sum_{i=1}^{m}\sigma_i(\bA\bB^\top) + \normf{\bB}^2.
$$	
The Frobenius bounding inequality (Theorem~\ref{theorem:vonneu_ineq2_sings}) shows that the equality is attained only if  $\bA\bB^\top\bQ^\top$ is PSD.
Let $\bA\bB^\top=\bS_r\bQ_r$ be the right polar decomposition, where $\bS_r$ is PSD and $\bQ_r$ is orthogonal. 
Thus, $\bA\bB^\top\bQ_r^\top=\bS_r$ is PSD and $\trace(\bA\bB^\top\bQ_r^\top)=\trace(\bS_r)=\sum_{i=1}^{m} \sigma_i(\bA\bB^\top)$, implying the lower bound is achieved.
The second part can be proved by noticing $\normf{\bA-\bB\bP}^2=\normf{\bB-\bP^\top\bA}^2$.
\end{proof}

The unitary Procrustes problem can be considered as a one-sided orthogonal optimization problem. The following generalization involves two orthogonal matrices.
\begin{theorem}[Two-Sided Orthogonal]\label{theorem:twos_si_pros}
Let $\bA,\bB\in\real^{ m \times n }$  be given with $p=\min\{m,n\}$, and let $\bQ\in\real^{m\times m}$ and $\bP\in\real^{n\times n}$ be  orthogonal matrices. Then, 
$$
\normf{\bA-\bQ\bB\bP}^2 \geq \sum_{i=1}^{p} (\sigma_i(\bA)-\sigma_i(\bB))^2.
$$
The equality is attained at $\bQ=\bU_1\bU_2^\top$ and $\bP=\bV_2\bV_1^\top$, where $ \bA = \bU_1 \bSigma_1 \bV_1^\top $ and $ \bB = \bU_2 \bSigma_2 \bV_2^\top $ are the SVDs of $\bA$ and $\bB$, respectively.
\end{theorem}
\begin{proof}[of Theorem~\ref{theorem:twos_si_pros}]
The Frobenius bounding inequality (Theorem~\ref{theorem:vonneu_ineq2_sings}) shows that 
$$
\normf{\bA-\bQ\bB\bP}^2 \geq 
\sum_{i=1}^{p} (\sigma_i(\bA)-\sigma_i(\bQ\bB\bP))^2
=
\sum_{i=1}^{p} (\sigma_i(\bA)-\sigma_i(\bB))^2.
$$
Let $ \bA = \bU_1 \bSigma_1 \bV_1^\top $ and $ \bB = \bU_2 \bSigma_2 \bV_2^\top $ be the SVDs of $\bA$ and $\bB$, respectively, where  $ \bU_1, \bU_2\in\real^{m\times m}$, $\bV_1, \bV_2 \in \real^{n\times n}$ are orthogonal, $\bSigma_1 = \diag(\sigma_1(\bA), \ldots, \sigma_p(\bA), 0,\ldots)$, and $\bSigma_2 = \diag(\sigma_1(\bB), \ldots, \sigma_p(\bB), 0,\ldots)$.
Let $\bQ\triangleq\bU_1\bU_2^\top$ and $\bP\triangleq\bV_2\bV_1^\top$, we have $\normf{\bA-\bQ\bB\bP}^2=\normf{\bU_1\bSigma_1\bV_1^\top - \bU_1\bSigma_2\bV_1^\top}^2=\sum_{i=1}^{p} (\sigma_i(\bA)-\sigma_i(\bB))^2$, achieving the lower bound.
\end{proof}

\subsection*{Interlacing Property for Singular Values}
The interlacing property of SVD is an extension of the interlacing property observed in symmetric matrices. 
We present the theorem directly, with the proof available in \citet{wilkinson1971algebraic} and further discussed in \citet{golub2013matrix} (Theorem 8.1.7, p. 443). \index{Interlacing property}
Alternatively, one can easily derive the result by applying Remark~\ref{remark:poincare_sep_eign}.
\begin{theorem}[Interlacing Property for Symmetric Matrices]
Let $\bA\in \real^{n\times n}$ be symmetric, and let $\bA_r$ be the upper-left $r\times r$ submatrix of $\bA$ such that $\bA_r=\bA[1:r, 1:r]$. 
Let $\lambda_i(\cdot)$ denote the $i$-th largest eigenvalue of the given matrix. 
For $k=r+1$,  we have 
$$
\lambda_{r+1}(\bA_{k})\leq \lambda_r(\bA_r) \leq \lambda_r(\bA_{k})\leq \ldots \leq \lambda_2(\bA_{k})\leq \lambda_1 (\bA_r)\leq \lambda_1(\bA_{k}).
$$
\end{theorem}

The interlacing property of singular values can be derived immediately from that of the symmetric matrices:
\begin{theorem}[Interlacing Property for Singular Values]\label{theorem:interlacing-singular}
Let $\bA=[\ba_1, \ba_2, \ldots, \ba_n]\in \real^{m\times n}$ be any matrix with $m\geq n$, and let $\bA_r=[\ba_1, \ba_2, \ldots, \ba_r]$ be the matrix containing the first $r$ columns of $\bA$. 
Define $\sigma_i(\cdot)$ as the $i$-th largest singular value of the given matrix. 
For $k=r+1$,  we have 
$$
\sigma_{k}(\bA_{k})\leq \sigma_r(\bA_r) \leq \sigma_r(\bA_{k})\leq \ldots \leq \sigma_2(\bA_{k})\leq \sigma_1 (\bA_r)\leq \sigma_1(\bA_{k}).
$$
\end{theorem}
The above result can be easily verified since $\bA_r^\top\bA_r$ can be obtained by deleting the last row and the last column of $\bA_k^\top\bA_k$.

To see the next result, we require the following lemmas.
\begin{lemma}[Interlacing Property for Singular Values]\label{lemma:interlacing-singular2}
Let  $\bA\in \real^{m\times n}$ be any matrix, and let $\bA_r$ denote the submatrix of $\bA$ after dropping out a total of $r$ columns of $\bA$. Then,
$$
\sigma_{k}(\bA) \geq 
\sigma_k(\bA_r)
\geq \sigma_{k+r}(\bA), 
\gap \forall k\in\{1,2,\ldots, \min\{m,n\}\},
$$
where $\sigma_k(\cdot)$ is the $k$-th largest singular value of the given matrix, and we set $\sigma_j(\bA)=0$ if $j>\min\{m,n\}$.
\end{lemma}
The lemma above can be derived from the Cauchy interlacing theorem (Remark~\ref{remark:poincare_sep_eign}).

\index{Gram–Schmidt}
\begin{lemma}[Submatrix Singular Value Inequality]\label{lemma:interl_sing_uk}
Let $\bA\in \real^{m\times n}$ be any matrix, and let $\bU_k\in\real^{m\times k}$ and $\bV_k\in\real^{n\times k}$ contain orthonormal columns: $\bU_k^\top\bU_k = \bV_k^\top\bV_k = \bI_k$, where $k\leq \min\{m,n\}$. 
Let $\sigma_i(\cdot)$ denote the $i$-th largest singular value of the given matrix. 
Then,
\begin{enumerate}
\item For $i\in\{1,2,\ldots,k\}$, we have $\sigma_{i}(\bU_k^\top \bA\bV_k) \leq \sigma_i(\bA)$.
\item $|\det(\bU_k^\top \bA\bV_k)| \leq \sigma_1(\bA) \sigma_2(\bA) \ldots \sigma_k(\bA) $.
\end{enumerate} 
\end{lemma}
\begin{proof}[of Lemma~\ref{lemma:interl_sing_uk}]
For (1), we can extend the semi-orthogonal matrices $\bU_k$ and $\bV_k$ to orthogonal matrices $\bU$ and $\bV$ using the Gram-Schmidt process. Therefore, $\bU_k^\top\bA\bV_k$ is the $k\times k$ upper-left submatrix of $\bU^\top\bA\bV$.
According to Lemma~\ref{lemma:interlacing-singular2}, we have $\sigma_i(\bA)=\sigma_i(\bU^\top\bA\bV) \geq \sigma_i(\bU_k^\top\bA\bV_k)$ for $i\in\{1,2,\ldots,k\}$.

For (2), by the fact that the magnitude of the determinant of a matrix is equal to the product of its singular values (Section~\ref{section:nonsingul_inver_svd}), we get 
$$
\begin{aligned}
|\det(\bU_k^\top \bA\bV_k)| &= \sigma_1(\bU_k^\top \bA\bV_k) \sigma_2(\bU_k^\top \bA\bV_k)\ldots \sigma_k(\bU_k^\top \bA\bV_k)
\leq \sigma_1(\bA) \sigma_2(\bA)\ldots \sigma_k(\bA).
\end{aligned}
$$
This completes the proof.
\end{proof}

\begin{remark}[Submatrix Singular Values Inequality in QR Decomposition]\label{remark:qr_inter_svd}
 Lemma~\ref{lemma:interlacing-singular2} and Lemma~\ref{lemma:interl_sing_uk} presented here are employed in the demonstration of \textit{Weyl's inequality}. 
In the context of low-rank QR decomposition, these lemmas can also be applied to establish  bounds for low-rank matrix approximation \citep{gu1996efficient}. 
To see this, we consider the rank-revealing QR decomposition (Section~\ref{section:rank-r-qr}) of $\bA\in\real^{m\times n}$ with $m\geq n $:
$$
\bA\bP = \bQ\bR = 
\bQ
\begin{bmatrix}
\bA_k & \bB_k\\
\bzero & \bC_k
\end{bmatrix}
\quad\implies\quad  
\begin{bmatrix}
	\bA_k & \bB_k\\
	\bzero & \bC_k
\end{bmatrix}
=
\bQ^\top\bA\bP,
$$
where $\bQ\in\real^{m\times m}$ is orthogonal, $\bA_k\in\real^{k\times k}$ is upper triangular with nonnegative diagonal elements, $\bB_k\in\real^{k\times (n-k)}$, $\bC_k\in\real^{(m-k)\times (n-k)}$, and $\bP\in\real^{n\times n}$ is a permutation matrix used to reveal linear dependencies among the columns of $\bA$. 
The choice of the parameter $k$ is typically made to be the smallest integer that yields a sufficiently small value of  $\normtwo{\bC_k}$.

Let $\bQ_k$ denote the first $k$ columns of $\bQ$, and $\bP_k$ denote the first $k$ columns of $\bP$.
This yields
$
\bA_k = \bQ_k^\top\bA\bP_k.
$
Based on Lemma~\ref{lemma:interl_sing_uk}, we obtain 
$$
\text{(QR1)}\gap \sigma_i(\bA_k)\leq \sigma_i(\bA), \gap \forall i\in\{1,2,\ldots, k\}.
$$
We can further construct a row permutation matrix $\bP_1$ and a column permutation matrix $\bP_2$ such that:
$$
\begin{bmatrix}
	\bC_k^\top & \bB_k^\top\\
	\bzero & \bA_k^\top 
\end{bmatrix}
=
\underbrace{\bP_1\bP^\top}_{\triangleq\bU} \bA^\top\underbrace{\bQ\bP_2 }_{\triangleq\bV}.
$$
By Lemma~\ref{lemma:interlacing-singular2}, we have 
$$
\text{(QR2)}\gap \sigma_j(\bC_k) \geq \sigma_{j+k}(\bA), \gap \forall j\in\{1,2,\ldots, n-k\},
$$
where $\sigma_j(\cdot)$ represents the $j$-th largest singular value of the given matrix.
The inequalities (QR1) and (QR2) also imply
$$
\sigma_{\min}(\bA_k)\leq \sigma_k(\bA) \gap 
\text{and}
\gap
\sigma_{\max}(\bC_k) \geq \sigma_{k+1}(\bA),
$$
where $\sigma_{\min}(\bA)$ and $\sigma_{\max}(\bA)$ denote the smallest and largest singular values of $\bA$, respectively.
Therefore, if we find an algorithm that can make $\sigma_{\max}(\bC_k)$ small in norm, it inherently reveals the rank-revealing QR decomposition.
\end{remark}

\index{Weyl’s inequality}
We are now ready to demonstrate \textit{Weyl's inequality}.
\begin{theorem}[Weyl's Inequality]\label{theorem:weyl_ineq}
Let $\bA\in\real^{n\times n}$ be any matrix,  with  eigenvalues $\lambda_1, \lambda_2, \ldots, \lambda_n\in \complex $ and singular values $\sigma_1, \sigma_2, \ldots, \sigma_n\in \real_+$ such that 
$
|\lambda_1| \geq |\lambda_2| \geq \ldots \geq |\lambda_n|
$
and 
$
\sigma_1 \geq \sigma_2 \geq \ldots \geq \sigma_n
$. 
Then,
$$
\begin{aligned}
&|\lambda_1||\lambda_2| \ldots  |\lambda_n| = \sigma_1  \sigma_2  \ldots  \sigma_n;\\
&|\lambda_1||\lambda_2| \ldots  |\lambda_k| \leq  \sigma_1  \sigma_2  \ldots  \sigma_k, \gap \forall k\in\{1,2,\ldots, n-1\}.
\end{aligned}
$$
Since all the components are nonnegative, we also have 
$$
\begin{aligned}
&|\lambda_1|+|\lambda_2|+ \ldots  +|\lambda_k| \leq  \sigma_1 + \sigma_2+  \ldots  +\sigma_k, \gap \forall k\in\{1,2,\ldots, n\}.
\end{aligned}
$$
\end{theorem}
\begin{proof}[of Theorem~\ref{theorem:weyl_ineq}]
For simplicity, we only show the case for $\bA$ with real eigenvalues. The result  can be easily extended to complex cases.
Suppose $\bA$ admits Schur decomposition $\bA = \bU\bT\bU^\top$ (Theorem~\ref{theorem:schur-decomposition}).
Let $\bU_k$ contain the first $k$ columns of $\bU$, then $\bT_k\triangleq\bU_k^\top\bA\bU_k$ is the $k\times k$ upper-left submatrix of $\bT=\bU^\top\bA\bU$. 
By the fact that the determinant of an upper triangular matrix is equal to the product of its main diagonal entries (which is also the product of its eigenvalues), we have 
$$
|\lambda_1 \lambda_2 \ldots \lambda_k| = |\det(\bT_k)|
=|\det(\bU_k^\top\bA\bU_k)|
\leq \sigma_1\sigma_2\ldots \sigma_k,
$$
where the last inequality is derived from Lemma~\ref{lemma:interl_sing_uk}. 
When $k=n$, the equality is attained because $|\det(\bA)| = \sigma_1\sigma_2\ldots\sigma_n$ and $\det(\bA)=\lambda_1\lambda_2\ldots\lambda_n$.
This completes the proof.
\end{proof}

Weyl's inequality relates the eigenvalues and singular values of a matrix. The following result then establishes a connection between the singular values of the product of matrices.
\begin{theorem}[Singular Value Inequalities]\label{theorem:weyl_ineq2_sings}
Let $\bA\in \real^{m\times r}, \bB\in\real^{r\times n}$ be given with $q=\min\{n,r,m\}$. Denote the ordered singular values of $\bA,\bB,$ and $\bA\bB$ by 
$$
\begin{aligned}
\sigma_1(\bA)&\geq \sigma_2(\bA) \geq \ldots \geq \sigma_{\min\{m,r\}}(\bA) \geq 0, \\
\sigma_1(\bB)&\geq \sigma_2(\bB) \geq \ldots \geq \sigma_{\min\{r,n\}}(\bB) \geq 0,\\
\sigma_1(\bA+\bB)&\geq \sigma_2(\bA+\bB) \geq \ldots \geq \sigma_{\min\{m,n\}}(\bA+\bB) \geq 0,
\end{aligned}
$$
respectively. 
Then,
$$
\prod_{i=1}^{k}\sigma_i(\bA\bB) 
\leq \prod_{i=1}^{k} \sigma_i(\bA)\sigma_i(\bB), \gap k\in\{1,2,\ldots,q\}.
$$
If $n=r=m$, then the equality holds for $k=n$. Since all the components are nonnegative, we also have 
$$
\sum_{i=1}^{k}\sigma_i(\bA\bB) 
\leq \sum_{i=1}^{k} \sigma_i(\bA)\sigma_i(\bB), \gap k\in\{1,2,\ldots,q\}.
$$
\end{theorem}
\begin{proof}[of Theorem~\ref{theorem:weyl_ineq2_sings}]
Suppose the matrix $\bA\bB$ admits SVD $\bA\bB=\bU\bSigma\bV^\top$. Let $\bU_k$ and $\bV_k$ denote the first $k$ columns of $\bU$ and $\bV$, respectively. 
Consequently, $\bU_k^\top(\bA\bB)\bV_k = \diag(\sigma_1(\bA\bB), \sigma_2(\bA\bB),$ $\ldots,\sigma_k(\bA\bB))$, as $\bU_k^\top(\bA\bB)\bV_k$ represents the $k\times k$ upper-left submatrix of $\bU^\top(\bA\bB)\bV$.

Suppose  the matrix $\bB\bV_k\in\real^{r\times k}$ admits reduced SVD $\bB\bV_k=\bQ_1\bD\bQ_2^\top$, where $\bQ_1\in \real^{r\times k}, \bD\in\real^{k\times k}$, and $\bQ_2\in\real^{k\times k}$ since $r\geq k$. 
Thus, we have $\bQ_2\bD^2\bQ_2^\top = (\bB\bV_k)^\top(\bB\bV_k) = \bV_k^\top\bB^\top\bB\bV_k$. 
As per Lemma~\ref{lemma:interl_sing_uk}, it follows that 
$$
\begin{aligned}
\left\vert \det\left(\bQ_2\bD^2\bQ_2^\top\right) \right\vert&=
\left\vert \det\left(\bV_k^\top\bB^\top\bB\bV_k  \right)\right\vert
\leq \sigma_1(\bB^\top\bB) \sigma_2(\bB^\top\bB)\ldots  \sigma_k(\bB^\top\bB)\\
&= \sigma_1(\bB)^2 \sigma_2(\bB)^2 \ldots  \sigma_k(\bB)^2.
\end{aligned}
$$
By the fact that the determinant of a matrix is equal to the product of eigenvalues and  referring back to Lemma~\ref{lemma:interl_sing_uk}, we have 
$$
\begin{aligned}
\sigma_1&(\bA\bB) \sigma_2(\bA\bB)\ldots\sigma_k(\bA\bB) =\left\vert\det\left(\bU_k^\top(\bA\bB)\bV_k\right)\right\vert
= \left\vert\det\left(\bU_k^\top\bA \bQ_1\bD\bQ_2^\top  \right)\right\vert\\
&=\left\vert\det\left(\bU_k^\top\bA \bQ_1  \right)\right\vert
\left\vert\det\left(\bD\bQ_2^\top \right)\right\vert
\leq 
\sigma_1(\bA) \sigma_2(\bA) \ldots  \sigma_k(\bA)
\sigma_1(\bB) \sigma_2(\bB) \ldots  \sigma_k(\bB).
\end{aligned}
$$
If $n=r=m$, then 
$$
\sigma_1(\bA\bB)\sigma_2(\bA\bB)\ldots \sigma_n(\bA\bB)
=
\left\vert\det (\bA\bB)\right\vert
=\left\vert\det (\bA)\right\vert\left\vert\det (\bB)\right\vert.
$$
Applying Weyl's inequality, we complete the proof.
\end{proof}

For further inequalities connecting the eigenvalues and singular values of a matrix, additional references can be found in \citet{horn1994matrix, horn2012matrix}.

\subsection*{Other Singular Value Inequalities}

\index{Browne’s theorem}
\begin{theorem}[Browne’s Theorem]\label{theorem:eig_sig_bound}
Let $\bA \in \complex^{n \times n}$, and let $\lambda \in \Lambda(\bA)$. Then, $\sigma_{\min}(\bA) \leq |\lambda| \leq \sigma_{\max}(\bA)$, which means the magnitude of an eigenvalue of a complex matrix is bounded by its smallest and largest singular values.
\end{theorem}
\begin{proof}
Let $(\lambda, \bu)$ be an eigenpair of $\bA$.
The inequality in Definition~\ref{definition:spectral_norm} shows $\abs{\lambda}\cdot\normtwo{\bu}=\normtwo{\bA\bu}\leq \sigma_{\max}(\bA)\normtwo{\bu}$ for any eigenpair of $\bA$. This proves $\abs{\lambda} \leq \sigma_{\max}(\bA)$. 
Moreover, Equation~\ref{equation:rayleig_ritz_sing} shows
$
\sigma_{\min}(\bA)=
\mathop{\min}_{\bx\neq \bzero}
\frac{\normtwo{\bA\bx}}{\normtwo{\bx}}
\leq
\frac{\normtwo{\bA\bu}}{\normtwo{\bu}} = \abs{\lambda}.
$
This proves the result.
\end{proof}

\begin{proposition}[Eigenvalue and Singular Value Inequality]\label{propo:sing_eig_ineq1}
Let $\bA, \bB \in \real^{n \times n}$. Then, we have 
$$
\begin{aligned}
\rho(\bA+\bB) &\leq \sigma_{\max}(\bA+\bB) \leq \sigma_{\max}(\bA) + \sigma_{\max}(\bB).
\end{aligned}
$$
Suppose further $\bA$ and $\bB$ are symmetric, then
$$
\rho(\bA+\bB) = \sigma_{\max}(\bA+\bB) \leq \sigma_{\max}(\bA) + \sigma_{\max}(\bB) = \rho(\bA) + \rho(\bB)
$$
and
$$
\lambda_{\min}(\bA) + \lambda_{\min}(\bB) \leq \lambda_{\min}(\bA+\bB) \leq \lambda_{\max}(\bA+\bB) \leq \lambda_{\max}(\bA) + \lambda_{\max}(\bB).
$$
\end{proposition}
\begin{proof}[of Proposition~\ref{propo:sing_eig_ineq1}]
Theorem~\ref{theorem:eig_sig_bound} shows the spectral radius of $\bA+\bB$ satisfies $\rho(\bA+\bB) \leq \sigma_{\max}(\bA+\bB)$.
The inequality $\sigma_{\max}(\bA+\bB) \leq \sigma_{\max}(\bA) + \sigma_{\max}(\bB)$ follows from the property of the spectral norm:
$$
\sigma_{\max}(\bA+\bB)= \normtwo{\bA+\bB} \leq \normtwo{\bA}+\normtwo{\bB}=\sigma_{\max}(\bA) + \sigma_{\max}(\bB).
$$
When  $\bA$ is symmetric, all the eigenvalues are real and $\sigma_{\max}(\bA)=\rho(\bA)=\lambda_{\max}(\bA)$.
This completes the proof.
\end{proof}

\index{Matrix norm}
\subsection{Ky Fan $k$-Norms}\label{section:kyfan_knorm}
A norm on a vector (resp. matrix) space is a function that assigns a nonnegative value to each vector (resp. matrix), representing its ``length" or ``size." To ensure that this function behaves in a way that is consistent with our intuitive understanding of length and distance, a norm must satisfy three key properties: positive homogeneity, the triangle inequality, and nonnegativity.
In Appendix~\ref{appendix:matrix-norm}, we introduce various matrix norms and their properties. In short, a matrix norm $\norm{\cdot}$ for a matrix $\bA\in \real^{m\times n}$ should adhere to the following criteria:
\begin{itemize}
	\item \textit{Nonnegativity}. $\norm{\bA} \geq 0$, and the equality is obtained if and only if $\bA=\bzero $. 
	\item \textit{Positive homogeneity}. $\norm{\lambda \bA} = |\lambda| \cdot \norm{\bA}$ for any $\lambda \in \real$.
	\item \textit{Triangle inequality}. $\norm{\bA+\bB} \leq \norm{\bA}+\norm{\bB}$ for any matrices $\bA, \bB\in \real^{m\times n}$.
\end{itemize}

The spectral norm introduced in Definition~\ref{definition:spectral_norm} satisfies the above criteria. 
More generally, the singular values of a matrix can be employed to define a general matrix norm \citep{horn1994matrix, gallier2017fundamentals}.
\begin{definition}[Ky Fan $k$-Norm and Schatten Norms]\label{definition:ky_fan_norm}
\index{\kyfannorm}
\index{Schatten norm}
Let $\bA\in\real^{m\times n}$ be any matrix with $q=\min\{m,n\}$ and singular values $\sigma_1\geq  \sigma_2 \geq \ldots \geq \sigma_q\geq 0$. Then, for any $k\in\{1,2,\ldots, q\}$, the following quantity is called the \textit{\kyfannorm} of $\bA$:
$$
\norm{\bA}_{[k]} = \sigma_1+\sigma_2+\ldots +\sigma_k.
$$
More generally, given $p\geq 1$ and $1\leq k\leq q$, 
$$
\norm{\bA}_{[k,p]} = \left( \sigma_1^p+ \sigma_2^p+\ldots +\sigma_k^p\right)^{1/p},
$$
is called the \textit{\kyfanpknorm} of $\bA$. Specially, when $k=q$, $\norm{\bA}_{[q,p]}$ is called the \textit{Schatten $p$-norm} of $\bA$, which is the $\ell_p$ vector norm on the vector formed by singular values:
\begin{itemize}
\item When $p=2$, the Schatten 2-norm $\norm{\bA}_{[q,2]}$ reduces to the Frobenius norm.
\item When $p=\infty$, the Schatten $\infty$-norm $\norm{\bA}_{[q,\infty]}$ reduces to the spectral norm.
\item When $p=1$, the Schatten 1-norm $\norm{\bA}_{[q,1]}$ is also known as  the \textit{nuclear norm} or \textit{trace norm}, which is defined as the sum of its singular values. It is used as the measurement to find low-rank approximations of matrices and is a convex relaxation of the rank of matrices \citep{jain2017non}.
\item Let $\infty \geq z \geq p\geq 2\geq 1$, Proposition~\ref{prop:lp_norm_ineqs} shows $
	\norm{\bA}_{[q,\infty]}\leq \norm{\bA}_{[q,z]} \leq \norm{\bA}_{[q,p]}\leq \norm{\bA}_{[q,2]} \leq \norm{\bA}_{[q,1]}.
$
\end{itemize}

\end{definition}

The nonnegative and positive homogeneity properties of the \kyfannorm can be  easily verified. The triangle inequality can be shown through the following steps.

\index{Partial isometry}
\index{Rank $k$ partial isometry}
\begin{definition}[Rank $k$ Partial Isometry]
A matrix $\bA\in\real^{m\times n}$ is said to be a rank $k$ partial isometry if $\sigma_1(\bA)=\sigma_2(\bA) = \ldots =\sigma_k(\bA)=1$ and 
$\sigma_{k+1}(\bA)=\sigma_{k+2}(\bA)=\ldots=\sigma_q(\bA)=0$, where $q=\min\{m,n\}$ and $\sigma_i(\bA)$ is the $i$-th largest singular value of $\bA$. 
In the meantime, two partial isometries $\bA, \bB\in\real^{m\times n}$ are said to be \textit{orthogonal} if $\bA^\top\bB=\bzero$ and $\bA\bB^\top=\bzero$.
\end{definition}

\begin{lemma}\label{lemma:kytriang_lem1}
Let $\bA\in\real^{m\times n}$ with  ordered singular values $\sigma_1 \geq \sigma_2 \geq \ldots \geq \sigma_q\geq 0$ and $q=\min\{m,n\}$. Then, for $k\in\{1,2,\ldots,q\}$, we have 
$$
\begin{aligned}
\sum_{i=1}^{k} \sigma_i 
&=
\max\left\{ |\trace(\bA\bB)|: \bB \in\real^{n\times m} \text{ is a rank $k$ partial isometry}  \right\}
\end{aligned}
$$
\end{lemma}
\begin{proof}[of Lemma~\ref{lemma:kytriang_lem1}]
Using the fact that the determinant of a matrix is equal to the sum of its eigenvalues, along with the results from Theorem~\ref{theorem:weyl_ineq} and Theorem~\ref{theorem:weyl_ineq2_sings}, we have  
$$
\begin{aligned}
\left\vert \trace(\bA\bB )\right\vert
&=
\left\vert \sum_{i=1}^{m} \lambda_i(\bA\bB) \right\vert
\leq 
\sum_{i=1}^{m} \left\vert \lambda_i(\bA\bB) \right\vert
\leq 
\sum_{i=1}^{m} \sigma_i(\bA\bB)\\
&\leq 
\sum_{i=1}^{m} \sigma_i(\bA)\sigma_i(\bB)
\leq \sum_{i=1}^{q} \sigma_i(\bA)\sigma_i(\bB)
=\sum_{i=1}^{k} \sigma_i(\bA).
\end{aligned}
$$
Suppose $\bA$ admits SVD $\bA=\bU\bSigma\bV^\top$ and let $\bB = \bV\bD\bU^\top$, where 
$
\bD = 
\scriptsize
\begin{bmatrix}
	\bI_k & \bzero  \\
	\bzero & \bzero 
\end{bmatrix}.
$
Then, $\trace(\bA\bB) = \trace(\bU\bSigma\bD\bU^\top ) = \trace(\bSigma\bD) = 
\sum_{i=1}^{k} \sigma_i(\bA),
$
where the second equality follows from the invariance of the trace under cyclic permutations.
Therefore, the equality can be achieved.
\end{proof}

\begin{theorem}[Ky Fan Triangle Inequality]\label{theorem:kyfan_triangle}
Let $\bA,\bB\in\real^{m\times n}$ with ordered singular values $\sigma_1(\bA)\geq \sigma_2(\bA)\geq \ldots \sigma_q(\bA)\geq 0$ and $\sigma_1(\bB)\geq \sigma_2(\bB)\geq \ldots \sigma_q(\bB)\geq 0$, where $q=\min\{m,n\}$. Let further $\sigma_1(\bA+\bB)\geq \sigma_2(\bA+\bB)\geq \ldots \sigma_q(\bA+\bB)\geq 0$ be the ordered singular values of $\bA+\bB$. Then,
$$
\norm{\bA+\bB}_{[k]} \leq \norm{\bA}_{[k]} + \norm{\bB}_{[k]}.
$$
\end{theorem}
\begin{proof}[of Theorem~\ref{theorem:kyfan_triangle}]
Let $\sS_{k}$ denote the set of rank $k$ partial isometries in $\real^{n\times m}$. 
According to Lemma~\ref{lemma:kytriang_lem1}, we have
$$
\begin{aligned}
\sum_{i=1}^{k} &\sigma_i(\bA+\bB)
= \max\{|\trace((\bA+\bB)\bC)|: \bC\in \sS_k\}
= \max\{|\trace(\bA\bC+\bB\bC)|: \bC\in \sS_k\}\\
&\leq  \max\{|\trace(\bB\bC)|+|\trace(\bB\bC)|: \bC\in \sS_k\}\\
&\leq \max\{|\trace(\bB\bC)|: \bC\in \sS_k\}+
\max\{|\trace(\bB\bC)|: \bC\in \sS_k\}
=\sum_{i=1}^{k} \sigma_i(\bA) +\sum_{i=1}^{k} \sigma_i(\bB).
\end{aligned}
$$
This completes the proof.
Alternatively, we may construct $\widetildebA\triangleq\scriptsize\begin{bmatrix}
\bzero & \bA\\
\bA^\top &\bzero
\end{bmatrix}$
and 
$\widetildebB\triangleq\scriptsize\begin{bmatrix}
	\bzero & \bB\\
	\bB^\top &\bzero
\end{bmatrix}$, which are symmetric.
Problem~\ref{prob:svd_4symm} shows $\lambda_i(\widetildebA)=\sigma_i(\bA)$ and $\lambda_i(\widetildebB)=\sigma_i(\bB)$ for $i\in\{1,2,\ldots,\min\{m,n\}\}$ (assume the eigenvalues and singular values are all in nonincreasing order).
Theorem~\ref{theorem:kyfan_eig_th} shows $\sum_{i=1}^{k} \lambda_i(\widetildebA+\widetildebB) \leq \big[\sum_{i=1}^{k} \lambda_i(\widetildebA) +\lambda_i(\widetildebB)\big]$. This again proves the result.
\end{proof}
Therefore, we verify the \kyfannorm is a valid matrix norm.

\begin{exercise}[Schatten Norm]
Verify that the \kyfanpknorm and the Scatten $p$-norm are valid matrix norms.
\end{exercise}

\section{Low-Rank Approximation via SVD}\label{section:svd-low-rank-approxi}
In the context of low-rank approximation problems, there are primarily two types, which arise from the interplay between rank and error:   the \textit{fixed-precision approximation problem} and the \textit{fixed-rank approximation problem}. For the fixed-precision approximation problem,  given a matrix $\bA$ and a  tolerance $\epsilon$, one wants to find a matrix $\bB$ with rank $r = r(\epsilon)$ such that $\norm{\bA-\bB} \leq \epsilon$ under an appropriate matrix norm. On the contrary, in the fixed-rank approximation problem, one seeks a matrix $\bB$ with a fixed rank $k$ and  aims to minimize the error $\norm{\bA-\bB}$. This section will concentrate on the latter problem. Some excellent examples of this can also be found in the literature by \citet{stewart1993early, kishore2017literature, martinsson2019randomized}.

Suppose we wish to approximate a matrix $\bA\in \real^{m\times n}$ of rank $r$ by a lower-rank matrix $\bB$ of rank $k$ ($k<r$). The quality of the approximation can be  measured using the spectral norm (Definition~\ref{definition:spectral_norm}):\index{Low-rank approximation}
$$
\bB = \mathop{\arg\min}_{\bB} \, \normtwo{\bA - \bB}.
$$
The optimal rank-$k$ approximation can be obtained using the following theorem, attributed to Eckart-Young-Mirsky.

\index{Truncated}
\index{Truncated SVD}
\begin{theorem}[Eckart-Young-Mirsky Theorem w.r.t. Spectral Norm\index{Eckart-Young-Mirsky theorem}]\label{theorem:young-theorem-spectral}
Let $\bA\in \real^{m\times n}$ be given with ordered singular values $\sigma_1\geq \sigma_2\geq \ldots \geq \sigma_r$ and $1\leq k\leq \rank(\bA)=r$, and let $\bA_k$ denote the \textit{truncated SVD (TSVD)} of $\bA$ with the largest $k$ singular terms, i.e., $\bA_k = \sum_{i=1}^{k} \sigma_i\bu_i\bv_i^\top$ derived from the SVD of $\bA=\sum_{i=1}^{r} \sigma_i\bu_i\bv_i^\top$ by zeroing out the $r-k$ trailing singular values of $\bA$. Then, $\bA_k$ is the optimal rank-$k$ approximation to $\bA$ in terms of the spectral norm.~\footnote{Note that $\bA_k$ can be stored using $(m+n)k + k$ entries, as opposed to $mn$ entries for the full matrix. }
\end{theorem}

\begin{proof}[of Theorem~\ref{theorem:young-theorem-spectral}]
To prove the theorem, we need to show that for any matrix $\bB$, if $\rank(\bB)=k$, then the inequality $\normtwo{\bA-\bB} \geq \normtwo{\bA-\bA_k}$ holds.

Since $\rank(\bB)=k$, then $\dim (\nspace(\bB))=n-k$. As a result, any set of $k+1$ basis vectors in $\real^n$ must intersect $\nspace(\bB)$. As shown in Proposition~\ref{proposition:svd-four-orthonormal-Basis}, $\{\bv_1,\bv_2, \ldots, \bv_r\}$ is an orthonormal basis for $\cspace(\bA^\top)\subset \real^n$; so  we can choose the first $k+1$ $\bv_i$'s as a $k+1$ basis for $\real^n$. Let $\bV_{k+1} \triangleq [\bv_1, \bv_2, \ldots, \bv_{k+1}]$. Then there exist a unit vector $\bx$ such that
$$
\bx \in \nspace(\bB) \cap \cspace(\bV_{k+1}),\qquad \text{s.t.}\,\,\,\, \normtwo{\bx}=1.
$$
That is, $\bx$ can be expressed as $\bx = \sum_{i=1}^{k+1} a_i \bv_i$, and $\normtwo{\sum_{i=1}^{k+1} a_i \bv_i} = \sum_{i=1}^{k+1}a_i^2=1$. Additionally, $\bB\bx=\bzero$.
Thus, using the definition of the spectral norm (Definition~\ref{definition:spectral_norm}), we have
$$
\begin{aligned}
\normtwo{\bA-\bB}^2 &\geq \normtwo{(\bA-\bB)\bx}^2\big/ \normtwo{\bx}^2\stackrel{*}{=} \normtwo{\bA\bx}^2
\stackrel{+}{=}\sum_{i=1}^{k+1} \sigma_i^2 (\bv_i^\top \bx)^2\\
&\stackrel{\dag}{\geq} \sigma_{k+1}^2\sum_{i=1}^{k+1}  (\bv_i^\top \bx)^2
\stackrel{\star}{=} \sigma_{k+1}^2\sum_{i=1}^{k+1} a_i^2=\sigma_{k+1}^2,
\end{aligned}
$$
where equality ($*$) follows because  $\bx$ is in the null space of $\bB$, equality (+) is due to $\bx$ being orthogonal to $\bv_{k+2}, \ldots, \bv_r$, inequality ($\dag$) follows from $\sigma_{k+1}\leq \sigma_{k}\leq\ldots\leq \sigma_{1}$, and equality ($\star$) follows from $\bv_i^\top \bx = a_i$.
On the other hand, it is evident that $\normtwo{\bA-\bA_k}^2 = \normtwo{\sum_{i=k+1}^{r}\sigma_i\bu_i\bv_i^\top}^2=\sigma_{k+1}^2$. Thus, $\normtwo{\bA-\bA_k} \leq \normtwo{\bA-\bB}$, which completes the proof. 
\end{proof}

\paragraph{Uniqueness.} If $\sigma_{k+1}=\sigma_k$, then another rank-$k$ approximation $\widetildebA_k=\sum_{i=1}^{k-1} \sigma_i\bu_i\bv_i^\top + \sigma_{k+1}\bu_{k+1}\bv_{k+1}^\top\neq\bA_k$ also achieves the minimum error  in the spectral norm. Therefore, the rank-$k$ approximation is usually not unique when there are repeated singular values. However, if all singular values are distinct, the approximation is unique.

Moreover, we can also prove that $\bA_k$ is the optimal rank-$k$ approximation to $\bA$ in terms of the Frobenius norm (Definition~\ref{definition:frobernius-in-svd}). 
The minimal error is given by the Euclidean norm of the singular values that have been zeroed out in the process: $\normf{\bA-\bA_k} =\sqrt{\sigma_{k+1}^2 +\sigma_{k+2}^2+\ldots +\sigma_{r}^2}$. The Frobenius bounding in Theorem~\ref{theorem:vonneu_ineq2_sings} can be applied directly to show this. 
Below, we provide an alternative proof derived from more elementary discussions.
\begin{theorem}[Eckart-Young-Mirsky Theorem w.r.t. Frobenius Norm\index{Eckart-Young-Mirsky theorem}]\label{theorem:young-theorem_frob}
Let $\bA\in \real^{m\times n}$ be given with ordered singular values $\sigma_1\geq \sigma_2\geq \ldots \geq \sigma_r$ and $1\leq k\leq \rank(\bA)=r$, and let $\bA_k$ denote the truncated SVD (TSVD) of $\bA$ with the largest $k$ singular terms, i.e., $\bA_k = \sum_{i=1}^{k} \sigma_i\bu_i\bv_i^\top$ derived from the SVD of $\bA=\sum_{i=1}^{r} \sigma_i\bu_i\bv_i^\top$ by zeroing out the $r-k$ trailing singular values of $\bA$. Then, $\bA_k$ is the optimal rank-$k$ approximation to $\bA$ in terms of the Frobenius norm.
\end{theorem}
\begin{proof}[of Theorem~\ref{theorem:young-theorem_frob}]
As per Definition~\ref{definition:frobernius-in-svd}, we note that $\normf{\bA - \bA_k}^2 = \normf{\bA}^2 - \sum_{i=1}^{k}\sigma_i^2$. We need to show that for any matrix $\bB$, if $\rank(\bB) = k$, it holds that $\normf{\bA-\bB}^2 \geq \normf{\bA-\bA_k}^2$.

The matrix $\bB$ can be written as a sum of $k$ rank-one matrices: $\bB = \sum_{i=1}^{k}\bx_i\by_i^\top$. Without loss of generality, we assume that the vectors $\{\bx_1, \bx_2, \ldots, \bx_k\}$ are mutually orthonormal; if not, we can use a linear combination  of the orthonormal set to express each vector. With this assumption,  we then have 
$$
\begin{aligned}
\normf{\bA-\bB}^2 &=\trace\bigg(  
\big(\bA - \sum_{i=1}^{k}\bx_i\by_i^\top\big)^\top 
\big(\bA - \sum_{i=1}^{k}\bx_i\by_i^\top\big)
\bigg)\\
&= \trace\bigg(
\bA^\top\bA + \sum_{i=1}^{k}(\by_i - \bA^\top\bx_i)(\by_i - \bA^\top\bx_i)^\top 
- \sum_{i=1}^{k}\bA^\top\bx_i\bx_i^\top \bA
\bigg).
\end{aligned}
$$
Since $\sum_{i=1}^{k}(\by_i - \bA^\top\bx_i)(\by_i - \bA^\top\bx_i)^\top $ is positive semidefinite, its diagonal values are nonnegative. Therefore, 
\begin{equation}\label{equation:ecth_eq1}
\begin{aligned}
\Vert\bA-\bB\Vert_F^2 &\geq 
\trace\bigg(
\bA^\top\bA 
- \sum_{i=1}^{k}\bA^\top\bx_i\bx_i^\top \bA
\bigg)
=
\Vert\bA\Vert_F^2 - \trace\bigg(\sum_{i=1}^{k}\bA^\top\bx_i\bx_i^\top \bA\bigg).
\end{aligned}
\end{equation}
Let $\bA = \bU\bSigma\bV^\top$ be the SVD of $\bA$, where $\bU = [\bU_1, \bU_2]$, $\bV=[\bV_1, \bV_2]$ with $\bU_1, \bV_2$ being the first $k$ columns of $\bU, \bV$ respectively. Let further $\bSigma_1 = \diag(\sigma_1, \sigma_2, \ldots, \sigma_k)$ and $\bSigma_2 = \diag(\sigma_{k+1}, \sigma_{k+2}, \ldots)$. Then we have 
$$
\bA\bA^\top = \bU\bSigma^2\bU^\top = [\bU_1, \bU_2] 
\begin{bmatrix}
\bSigma_1^2 & \bzero\\
\bzero & \bSigma_2^2
\end{bmatrix} 
\begin{bmatrix}
\bU_1^\top \\
\bU_2^\top
\end{bmatrix}
=
\bU_1 \bSigma_1^2 \bU_1^\top +
\bU_2 \bSigma_2^2 \bU_2^\top.
$$ 
Let $\bz_i \triangleq \bA^\top\bx_i\in \real^n$, we have
$$
\begin{aligned}
\trace\bigg(\sum_{i=1}^{k}\bA^\top\bx_i\bx_i^\top \bA\bigg) 
&=  \trace \bigg(\sum_{i=1}^{k} \bz_i \bz_i^\top  \bigg) 
=  \sum_{i=1}^{k} \sum_{j=1}^{n} z_{ij}^2
= \sum_{i=1}^{k} \bz_i^\top\bz_i
= \sum_{i=1}^{k} \bx_i^\top\bA\bA^\top\bx_i,
\end{aligned}
$$
which is the sum of $k$ scalars (invariance under cyclic permutations).
For each component of the above equation, we decompose
$$
\begin{aligned}
\bx_i^\top&\bA\bA^\top\bx_i = \bx_i^\top \bU_1\bSigma_1^2\bU_1^\top \bx_i 
+ \bx_i^\top \bU_2\bSigma_2^2\bU_2^\top \bx_i 
= \Vert\bSigma_1\bU_1^\top\bx_i\Vert_2^2+ \Vert\bSigma_2\bU_2^\top\bx_i\Vert_2^2\\
&= \sigma_k^2 + \Vert\bSigma_1\bU_1^\top\bx_i\Vert_2^2 - \sigma_k^2 \Vert\bU_1^\top\bx_i\Vert_2^2 
-\big( \sigma_k^2\Vert\bU_2^\top\bx_i\Vert_2^2- \Vert\bSigma_2\bU_2^\top\bx_i\Vert_2^2  \big)
-\sigma_k^2(1-\Vert\bU^\top\bx_i\Vert_2^2).
\end{aligned}
$$
Given that $\bU$ is orthogonal and $\bx_i$ has a unit length, we have $\Vert\bU^\top\bx_i\Vert_2^2=1$. And it can be easily shown that $\left( \sigma_k^2\Vert\bU_2^\top\bx_i\Vert_2^2- \Vert\bSigma_2\bU_2^\top\bx_i\Vert_2^2  \right)\geq 0$. Hence, we can deduce: 
\begin{equation}\label{equation:ecth_eq2}
\begin{aligned}
\sum_{i=1}^{k} \bx_i^\top\bA\bA^\top\bx_i &\leq 
k\sigma_k^2 +\sum_{i=1}^{k}\left(
\Vert\bSigma_1\bU_1^\top\bx_i\Vert_2^2 - \sigma_k^2 \Vert\bU_1^\top\bx_i\Vert_2^2
\right)
=
k\sigma_k^2 +\sum_{i=1}^{k} \sum_{j=1}^{k} (\sigma_j^2-\sigma_k^2)(\bv_j^\top\bx_i)^2\\
&=  k\sigma_k^2 +\sum_{j=1}^{k} (\sigma_j^2-\sigma_k^2)\bigg(\sum_{i=1}^{k} (\bv_j^\top\bx_i)^2\bigg)
\leq \sum_{j=1}^{k}\sigma_j^2.
\end{aligned}
\end{equation}
Combining Equation~\eqref{equation:ecth_eq1} and Equation~\eqref{equation:ecth_eq2}, we conclude that
$$
\Vert\bA-\bB\Vert_F^2 \geq 
\Vert\bA\Vert_F^2 - \sum_{j=1}^{k}\sigma_j^2,
$$
from which the result follows.
\end{proof}

While our exposition has centered around the optimal low-rank approximation using the Frobenius and spectral norms, it's worth highlighting that the Eckart-Young-Mirsky theorem extends its reach much further. Indeed, this theorem applies to any norm that exhibits unitary invariance, a category that includes not only the aforementioned norms but also the entire spectrum of Schatten norms and the operator norm. Demonstrating these broader applications would involve additional proofs, which we invite interested readers to explore as exercises. See the footnote of Theorem~\ref{theorem:vonneu_ineq2_sings} for a sketch of this topic.


\begin{SCfigure}
\centering
\includegraphics[width=0.6\textwidth]{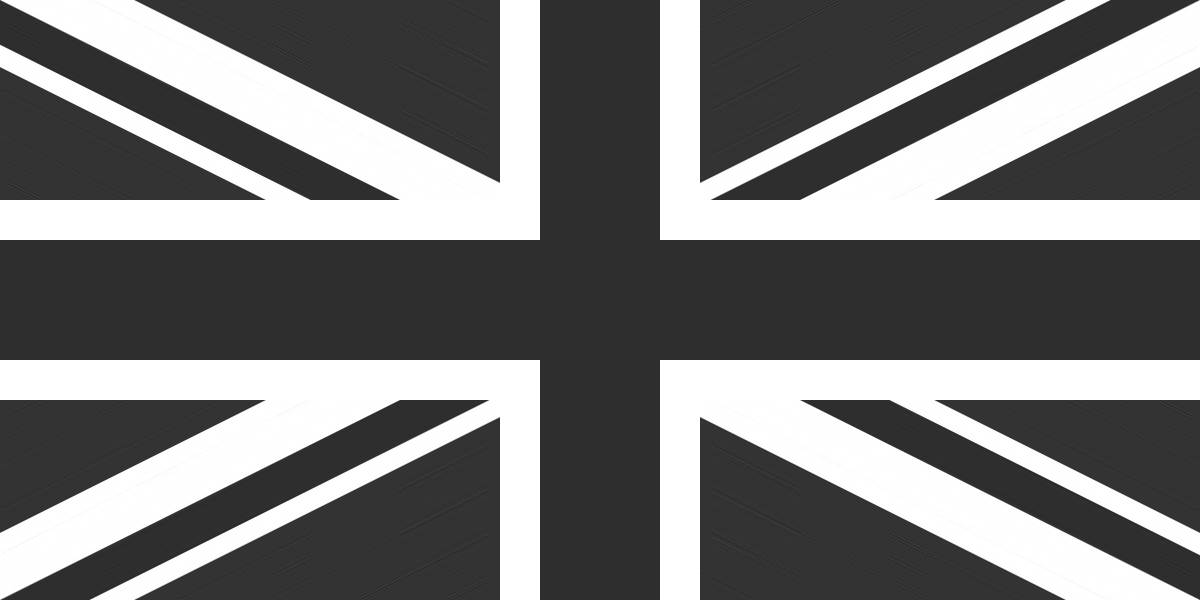}
\caption{The input is a grayscale flag image intended for compression. The dimension of the image is $600\times 1200$ with a rank of 402. }
\label{fig:eng300}
\end{SCfigure}

\begin{figure}[h]
\centering   
\vspace{-0.35cm}  
\subfigtopskip=2pt  
\subfigbottomskip=2pt  
\subfigcapskip=-5pt  
\subfigure[$\sigma_1\bu_1\bv_1^\top$, $F_1 = 60,217$.]{\label{fig:svd1}%
	\includegraphics[width=0.32\linewidth]{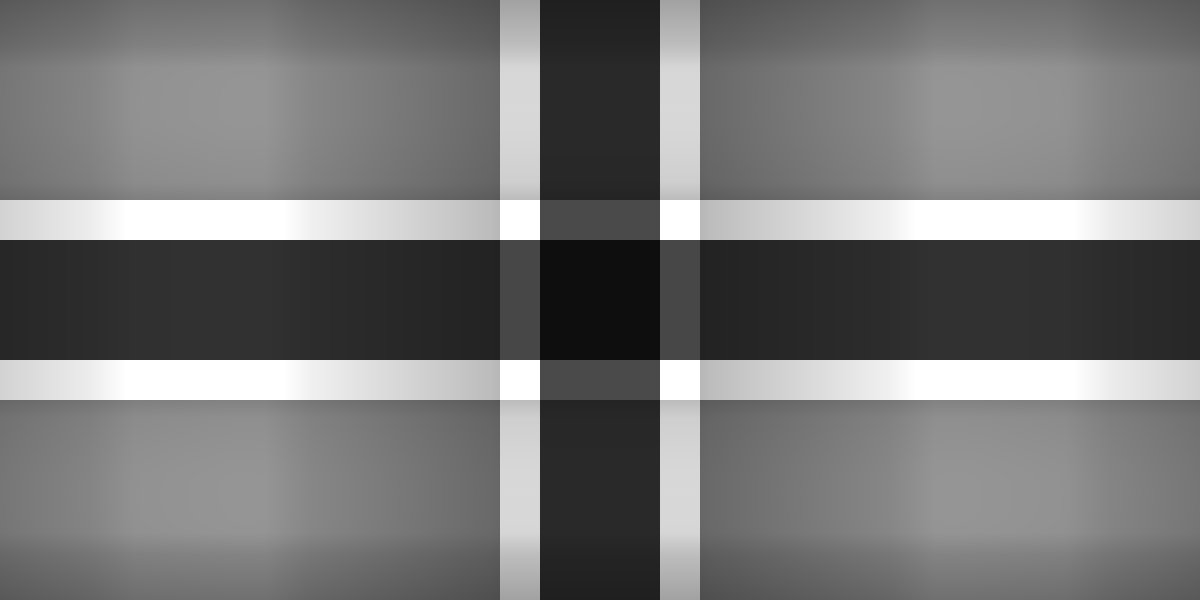}}
\subfigure[$\sigma_2\bu_2\bv_2^\top$, $F_2 = 120,150$.]{\label{fig:svd2}%
	\includegraphics[width=0.32\linewidth]{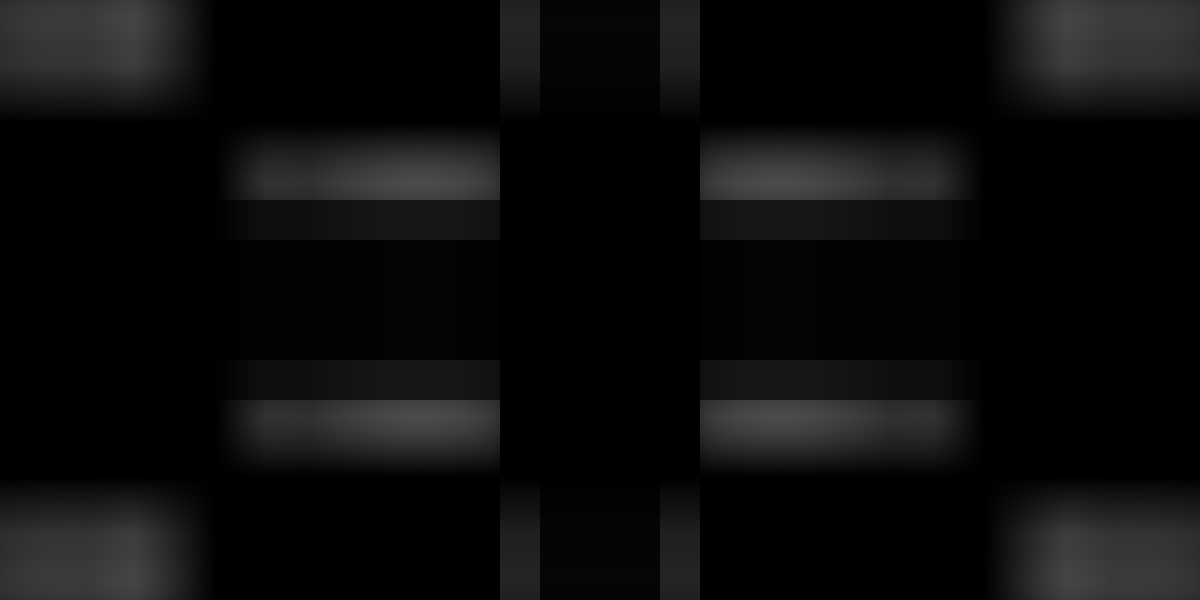}}
\subfigure[$\sigma_3\bu_3\bv_3^\top$, $F_3 = 124,141$.]{\label{fig:svd3}%
	\includegraphics[width=0.32\linewidth]{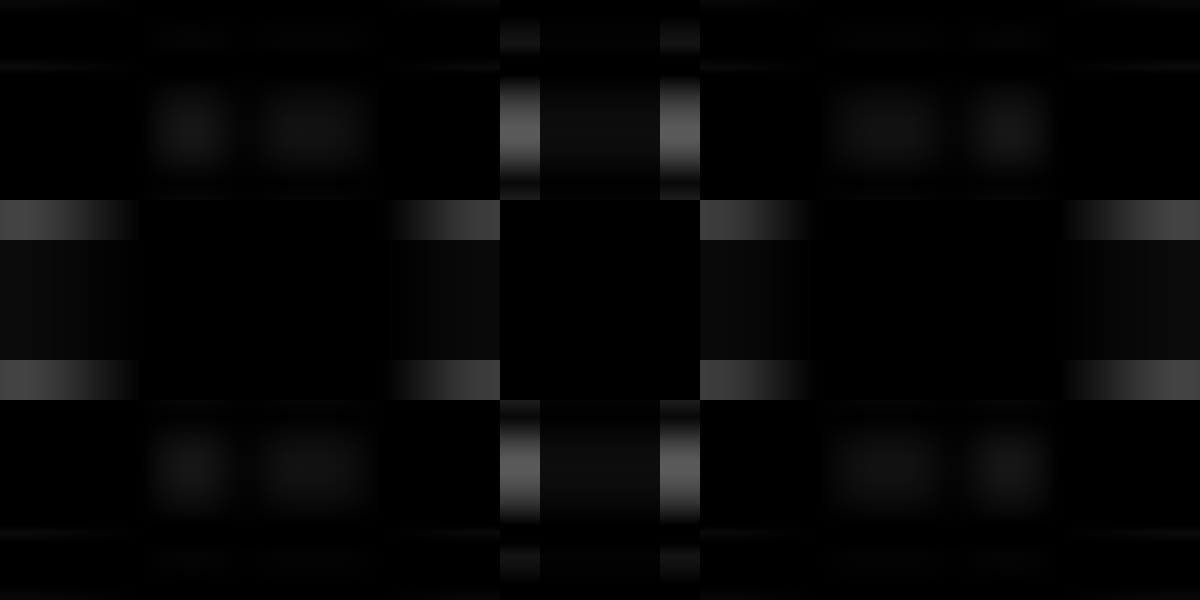}}\\
\subfigure[$\sigma_4\bu_4\bv_4^\top$, $F_4 = 125,937$.]{\label{fig:svd4}%
	\includegraphics[width=0.32\linewidth]{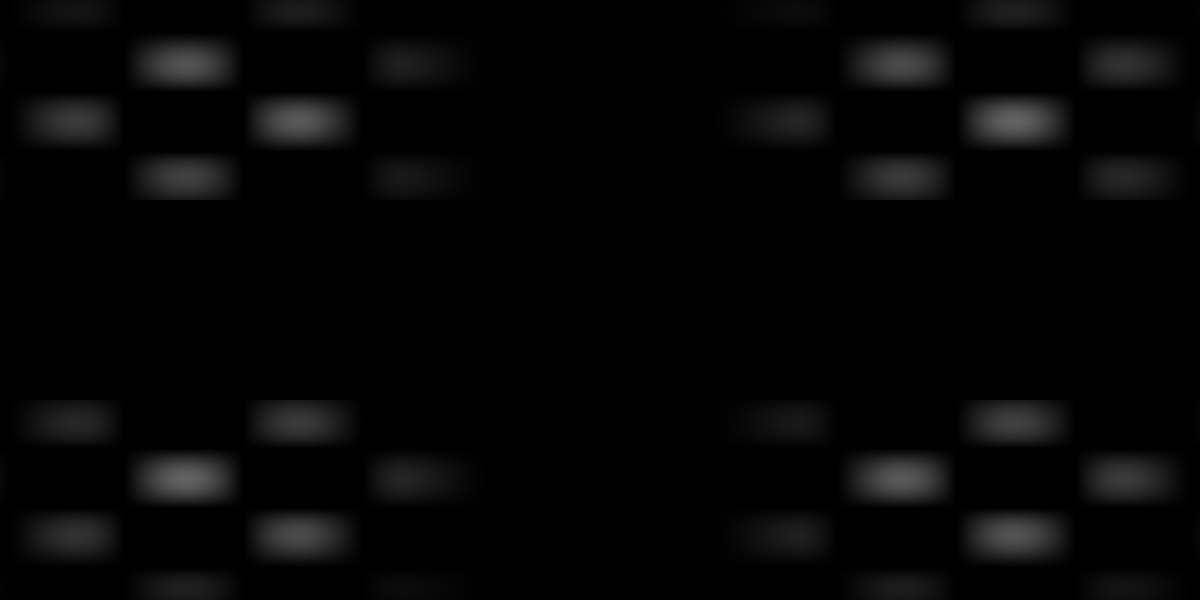}}
\subfigure[$\sigma_5\bu_5\bv_5^\top$, $F_5 = 126,127$.]{\label{fig:svd5}%
	\includegraphics[width=0.32\linewidth]{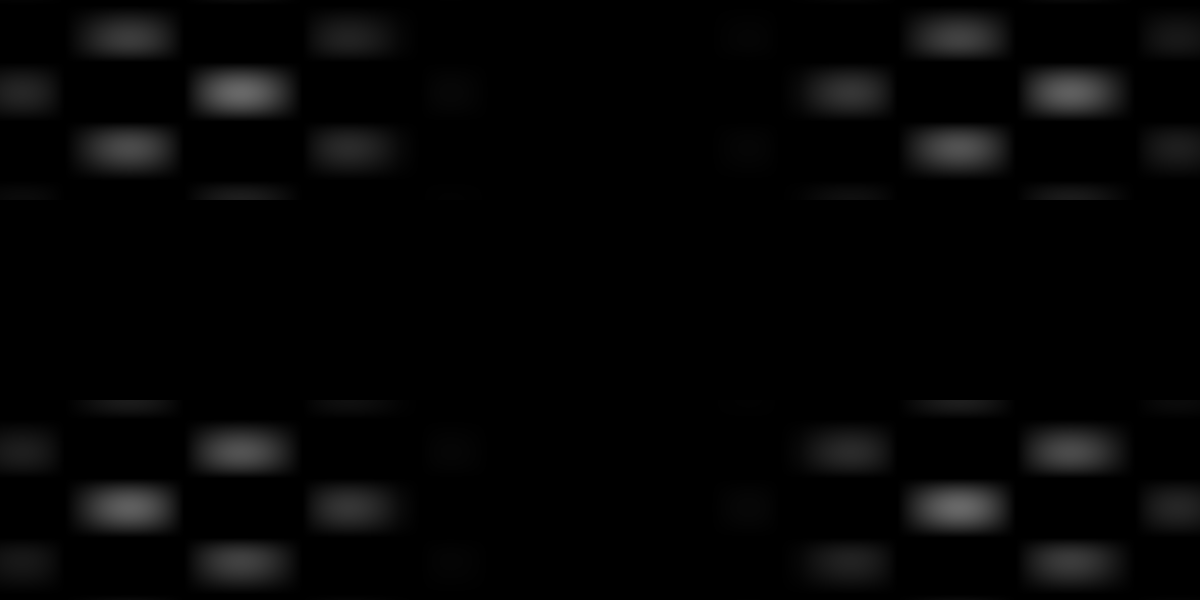}}
\subfigure[All 5 singular values: $\sum_{i=1}^{5}\sigma_i\bu_i\bv_i^\top$, \protect $F=44,379$.]{\label{fig:svd6}%
	\includegraphics[width=0.32\linewidth]{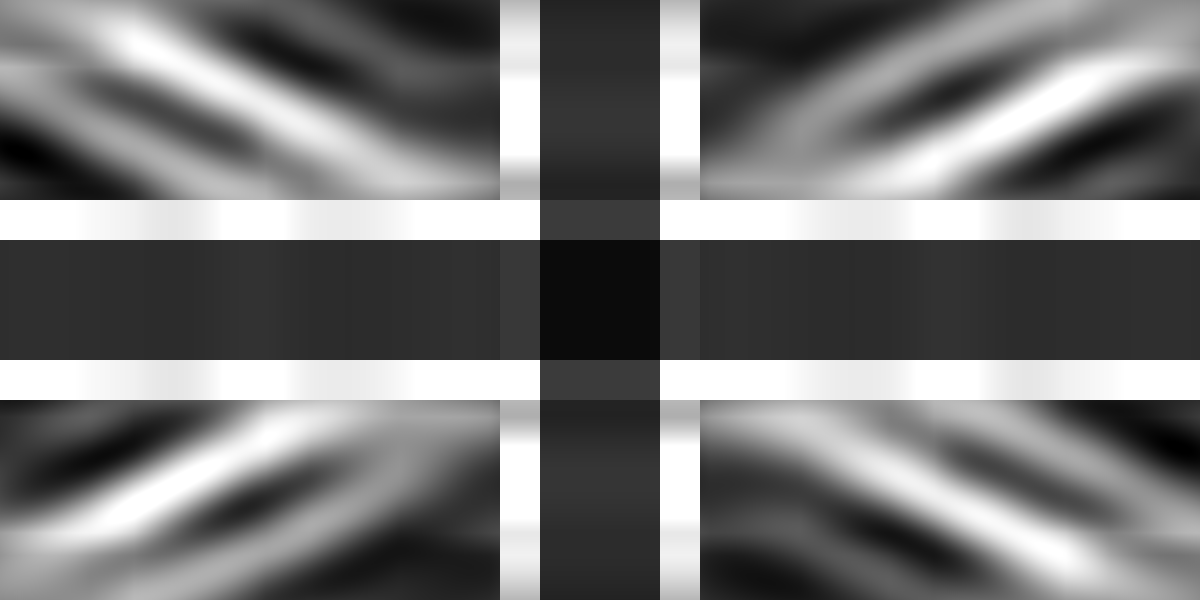}}
\caption{Image compression for a gray flag image into a rank-5 matrix via the SVD, and decompose into 5 parts, where $\sigma_1 \geq \sigma_2 \geq \ldots \geq \sigma_{5}$, i.e., $F_1\leq F_2\leq \ldots \leq F_5$ with $F_i \triangleq \normf{\sigma_i\bu_i\bv^\top - \bA}$ for $i\in \{1,2,\ldots, 5\}$. And reconstruct images using a single singular value and its corresponding left and right singular vectors.}
\label{fig:svdd-by-parts}
\end{figure}

\begin{figure}[h]
\centering  
\vspace{-0.35cm} 
\subfigtopskip=2pt 
\subfigbottomskip=2pt 
\subfigcapskip=-5pt 
\subfigure[$\bc_1\br_1^\top$, $G_1 = 60,464$.]{\label{fig:skeleton1}%
\includegraphics[width=0.32\linewidth]{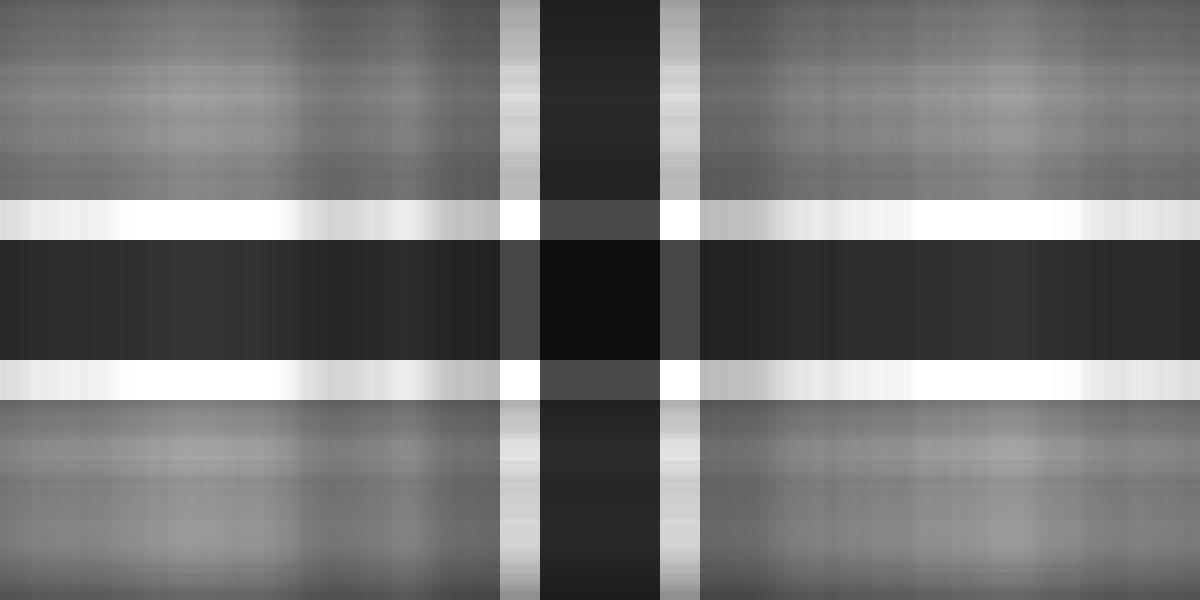}}
\subfigure[$\bc_2\br_2^\top$, $G_2 = 122,142$.]{\label{fig:skeleton2}%
\includegraphics[width=0.32\linewidth]{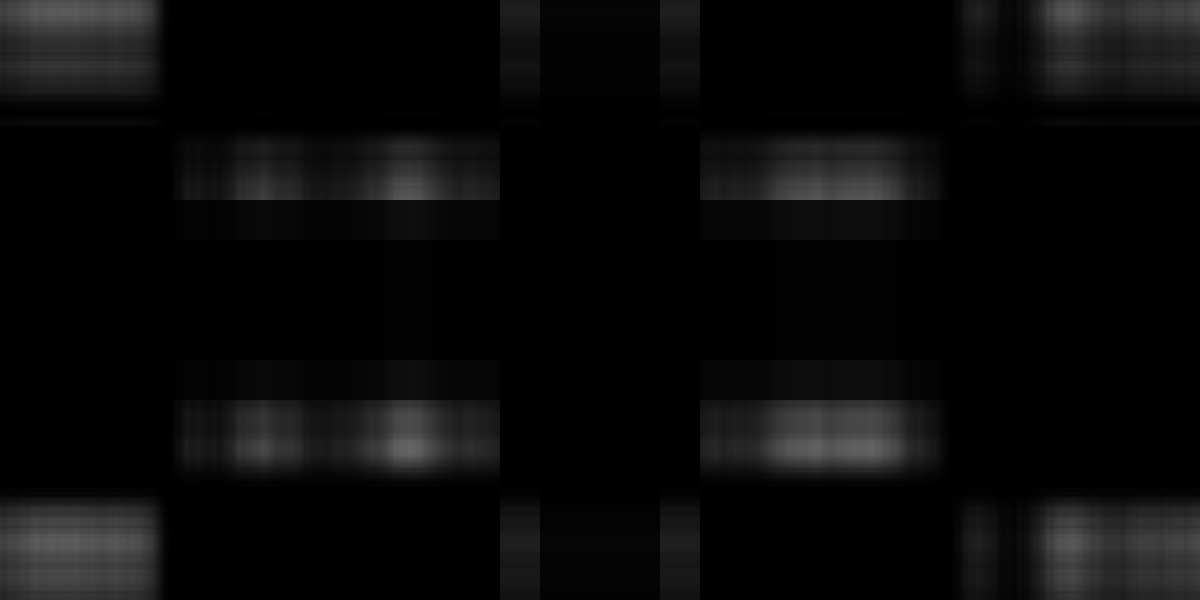}}
\subfigure[$\bc_3\br_3^\top$, $G_3 = 123,450$.]{\label{fig:skeleton3}%
\includegraphics[width=0.32\linewidth]{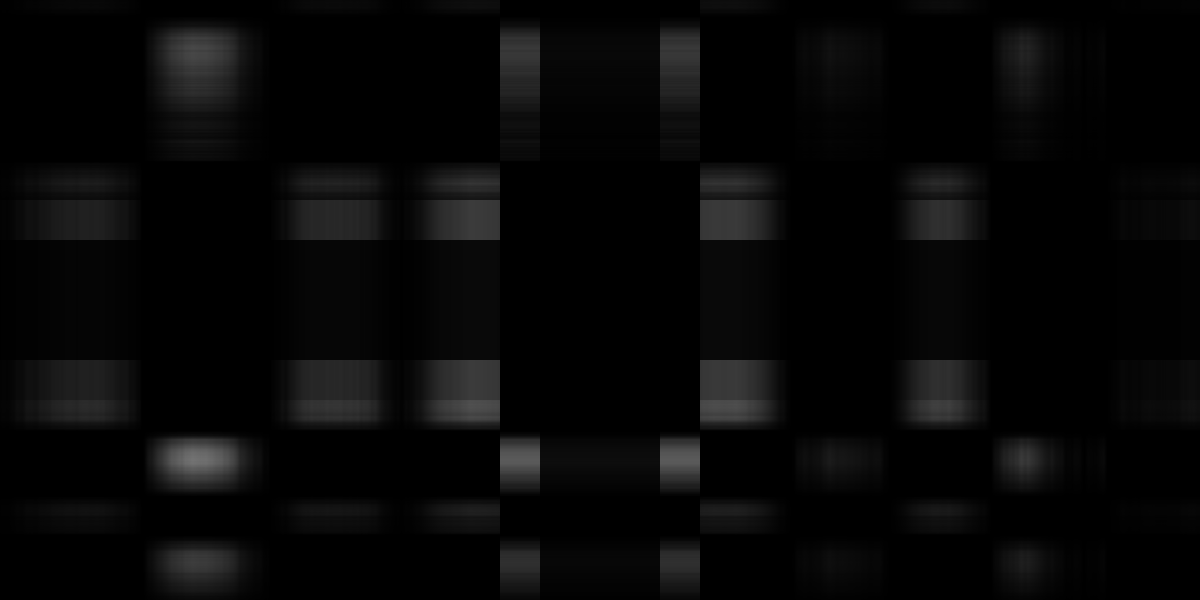}}\\
\subfigure[$\bc_4\br_4^\top$, $G_4 = 125,975$.]{\label{fig:skeleton4}%
\includegraphics[width=0.32\linewidth]{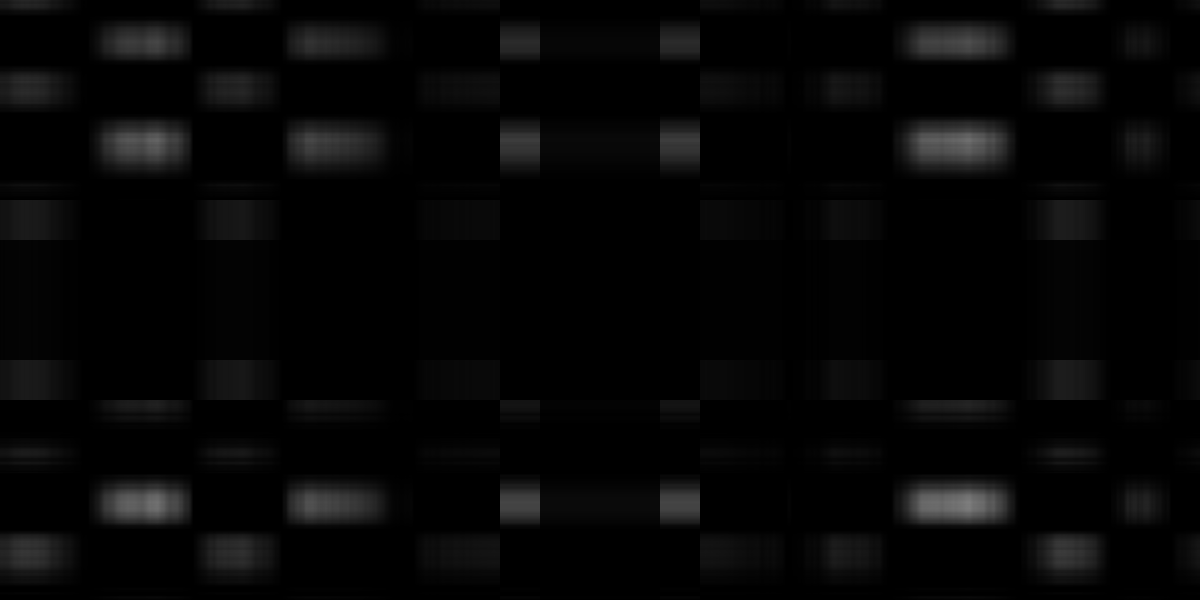}}
\subfigure[$\bc_5\br_5^\top$, $G_5 = 124,794$.]{\label{fig:skeleton5}%
\includegraphics[width=0.32\linewidth]{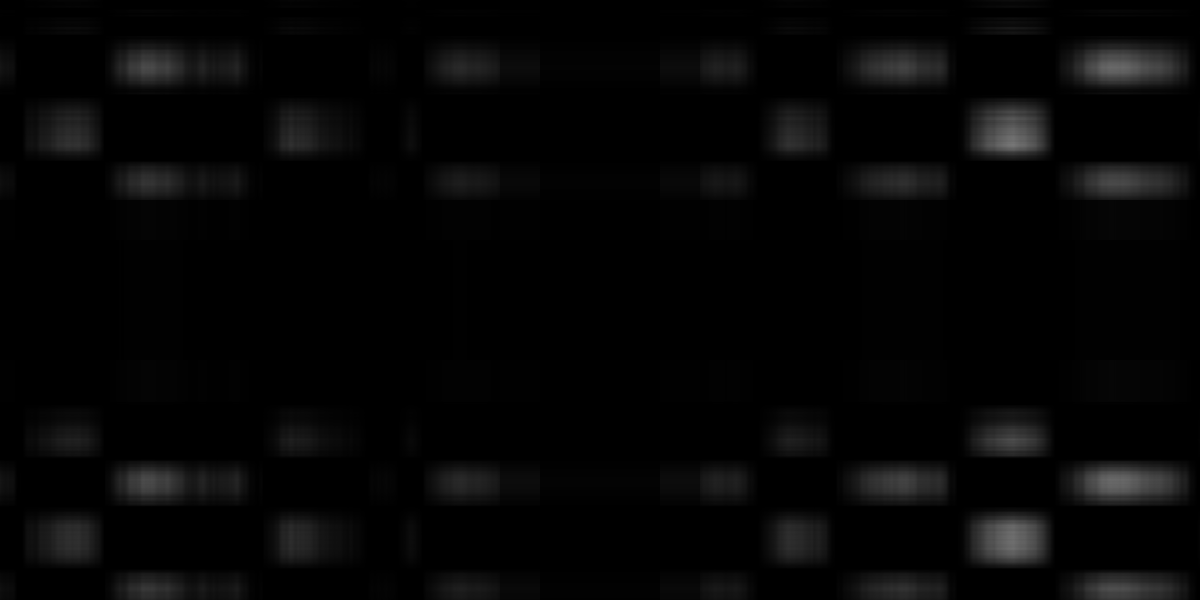}}
\subfigure[All 5 parts: $\sum_{i=1}^{5}\bc_i\br_i^\top$, \protect\newline $G=45,905$.]{\label{fig:skeleton6}%
\includegraphics[width=0.32\linewidth]{./imgs/svd_pic6_all.png}}
\caption{Image compression for a gray flag image into a rank-5 matrix via the pseudoskeleton decomposition, and decompose into 5 parts, where $G_i\triangleq\normf{\bc_i\br_i^\top-\bA}$ for $i\in \{1,2,\ldots, 5\}$ and $G_1\leq G_2 \leq \ldots \leq G_5$. And reconstruct images by $\bc_i\br_i^\top$.}
\label{fig:skeletond-by-parts}
\end{figure}

SVD thus provides the optimal approximation of a matrix in terms of various matrix norms. As noted by \citet{stewart1998matrix, kishore2017literature}, \emph{the singular value decomposition is the creme de la creme of rank-reducing decompositions---the decomposition that all others try to beat}. And also, \citet{strang1993introduction} states that \emph{the SVD is the climax of this linear algebra course}.

\index{Truncated}
\index{Truncated SVD}
Figure~\ref{fig:eng300} presents an example of a grayscale  image to be compressed. The dimension of the image is $600\times 1200$ with a rank of 402. 
In Figure~\ref{fig:svdd-by-parts}, we approximate the image into a rank-5 matrix using truncated SVD: $\bA\approx \sum_{i=1}^{5}\sigma_i\bu_i\bv_i^\top$. It is known that the singular values contain the spectrum information with higher singular values containing lower-frequency information. 
And low-frequency components often carry more pertinent information  \citep{leondes1995multidimensional}. 
We find that the image, $\sigma_1\bu_1\bv_1^\top$, reconstructed using the largest singular value closely resembles the original flag image. 
The subsequent images, reconstructed using the second to fifth singular values and their corresponding singular vectors, capture finer details to reconstruct the full image.

Similar results can be observed for the low-rank approximation through the pseudoskeleton decomposition (Section~\ref{section:pseudoskeleton}). In Equation~\eqref{equation:skeleton-low-rank}, we derived the low-rank approximation as $\bA\approx \bC_2\bR_2$, where $\bC_2\in \real^{m\times \gamma}$ and $\bR_2\in \real^{\gamma\times n}$ if $\bA\in \real^{m\times n}$ such that $\bC_2$ and $\bR_2$ are rank-$\gamma$ matrices. For $\gamma=5$: 
$$
\bC_2=[\bc_1, \bc_2, \ldots, \bc_5]
\qquad
\text{and}
\qquad 
\bR_2 = 
\footnotesize
\begin{bmatrix}
\br_1^\top \\
\br_2^\top \\
\vdots \\
\br_5^\top 
\end{bmatrix},
$$ 
are the column and row partitions of $\bC_2$ and $\bR_2$, respectively. Then $\bA$ can be approximated by $\sum_{i=1}^{5}\bc_i\br_i^\top$. The partitions are ordered such that 
$$
\underbrace{\normf{\bc_1\br_1^\top-\bA}}_{G_1} \leq 
\underbrace{\normf{\bc_2\br_2^\top-\bA}}_{G_2}
\leq \ldots \leq 
\underbrace{\normf{\bc_5\br_5^\top-\bA}}_{G_5}.
$$
We observe that $\bc_1\br_1^\top$ behaves similarly to $\sigma_1\bu_1\bv^\top$,  with reconstruction errors measured by the Frobenius norm being very close (60,464 for the pseudoskeleton case compared to  60,217 for the SVD case). This similarity is due to the pseudoskeleton decomposition relying on the SVD (Section~\ref{section:pseudoskeleton}), making $\bc_1\br_1^\top$ internally has the largest ``singular value" meaning in this sense.

\begin{SCfigure}
\centering
\includegraphics[width=0.53\textwidth]{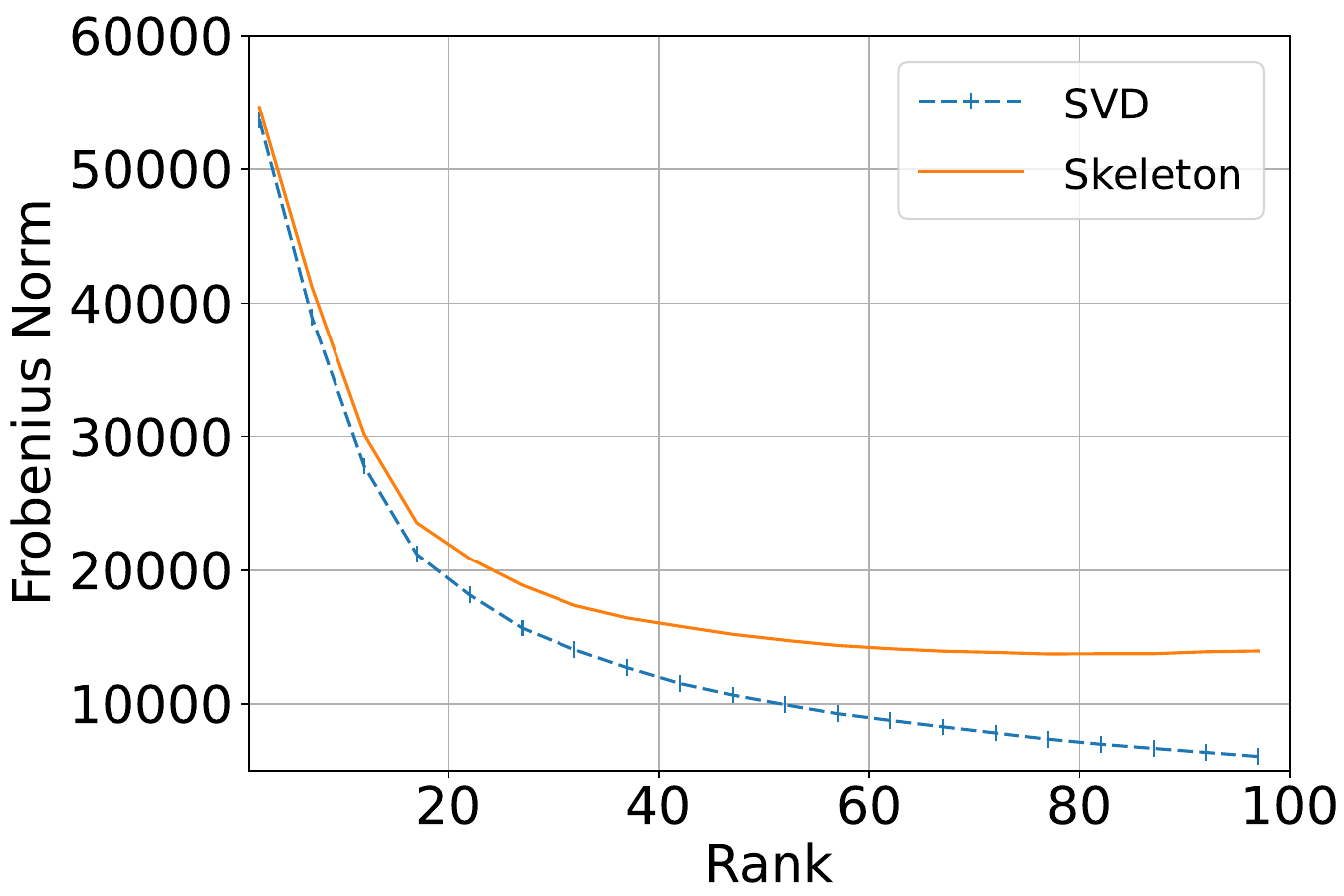}
\caption{Comparison of reconstruction errors measured by the Frobenius norm between the SVD and the pseudoskeleton approximation.}
\label{fig:svd_skeleton_fnorm}
\end{SCfigure}

\begin{figure}[h]
\centering 
\vspace{-0.35cm} 
\subfigtopskip=2pt 
\subfigbottomskip=2pt 
\subfigcapskip=-5pt 
\subfigure[SVD with rank 90,\gap\gap  \protect\newline Frobenius norm=\textbf{6,498}.]{\label{fig:svd90}%
\includegraphics[width=0.32\linewidth]{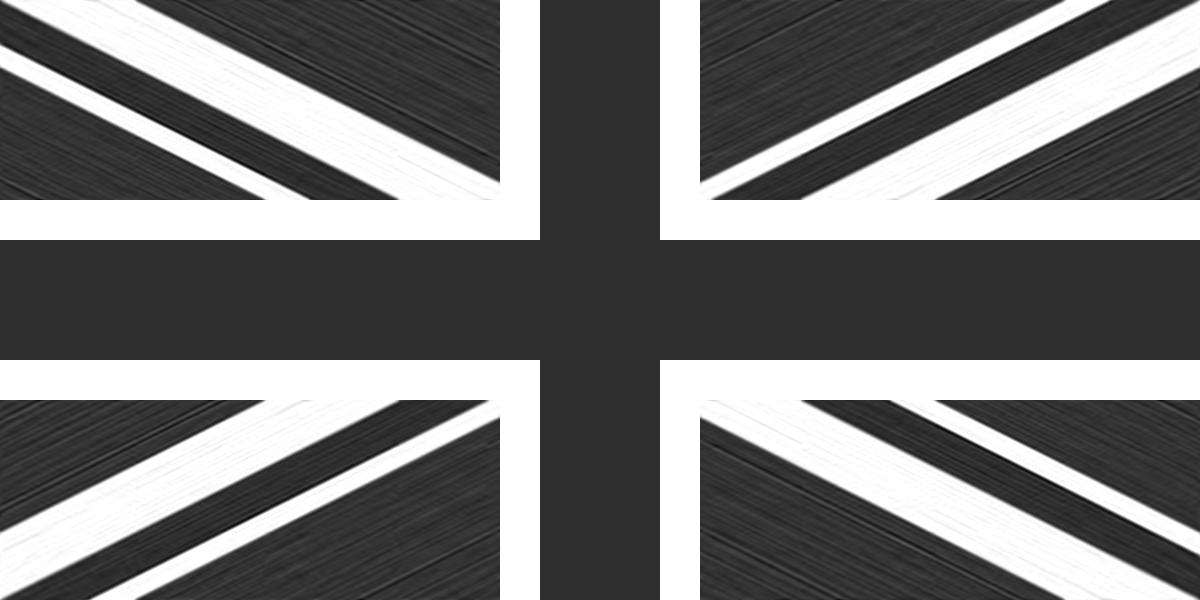}}
\subfigure[SVD with rank 30,\gap\gap \protect\newline Frobenius norm=\textbf{14,586}.]{\label{fig:svd30}%
\includegraphics[width=0.32\linewidth]{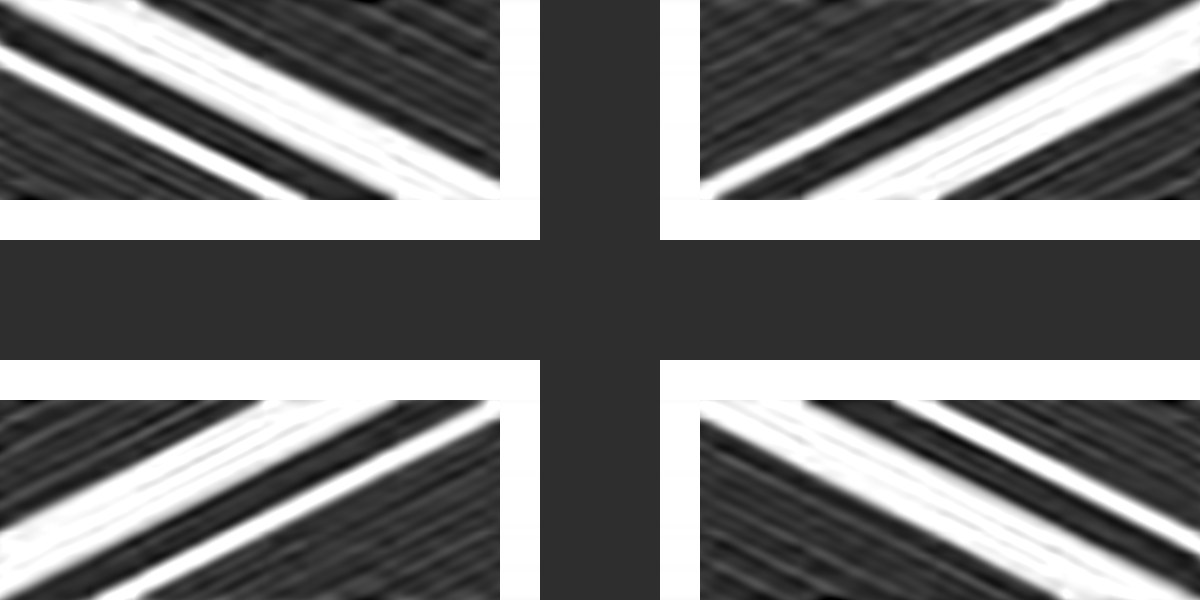}}
\subfigure[SVD with rank 10,\gap\gap \protect\newline Frobenius norm=\textbf{31,402}.]{\label{fig:svd10}%
\includegraphics[width=0.32\linewidth]{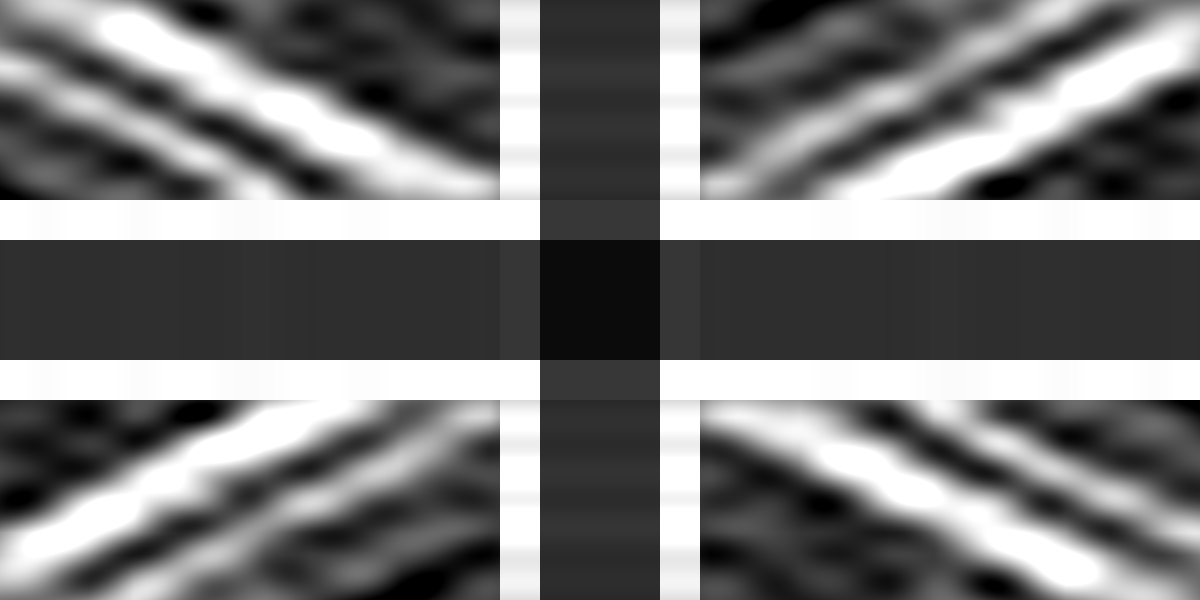}}\\
\subfigure[Pseudoskeleton with rank 90,\protect\newline Frobenius norm=13,751.]{\label{fig:skeleton90}%
\includegraphics[width=0.32\linewidth]{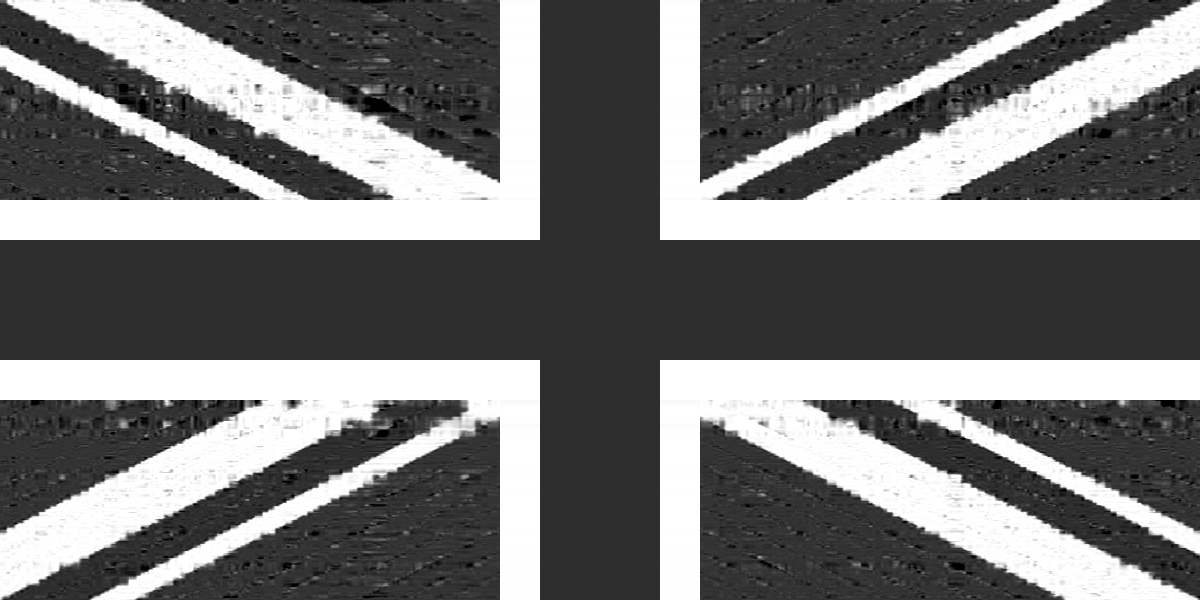}}
\subfigure[Pseudoskeleton with rank 30,\protect\newline Frobenius norm=17,853.]{\label{fig:skeleton30}%
\includegraphics[width=0.32\linewidth]{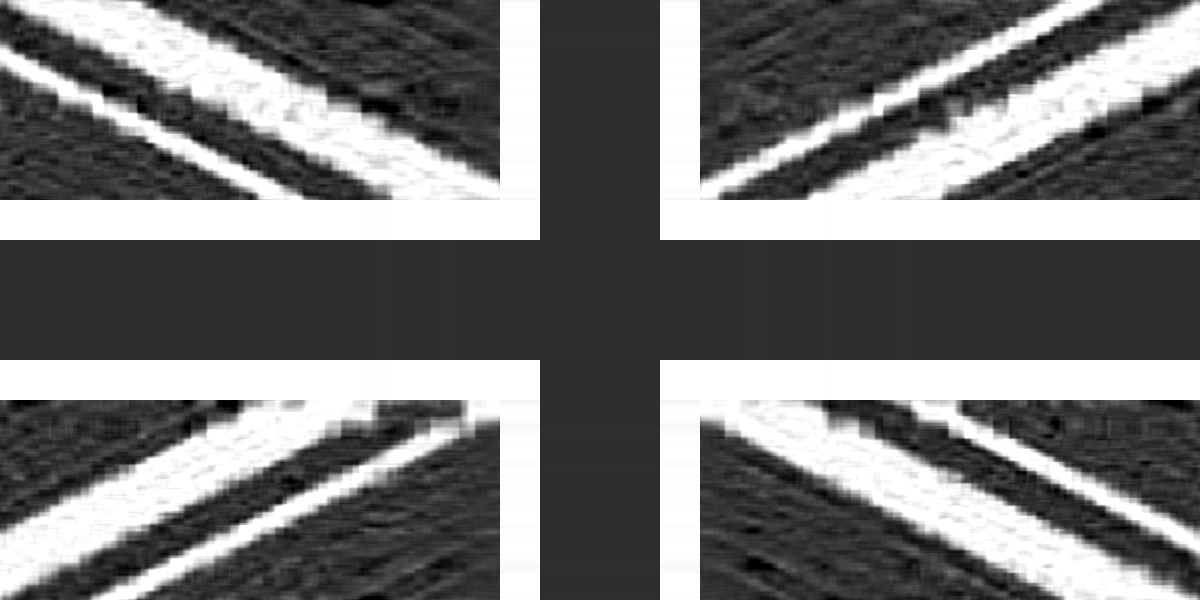}}
\subfigure[Pseudoskeleton with rank 10,\protect\newline Frobenius norm=33,797.]{\label{fig:skeleton10}%
\includegraphics[width=0.32\linewidth]{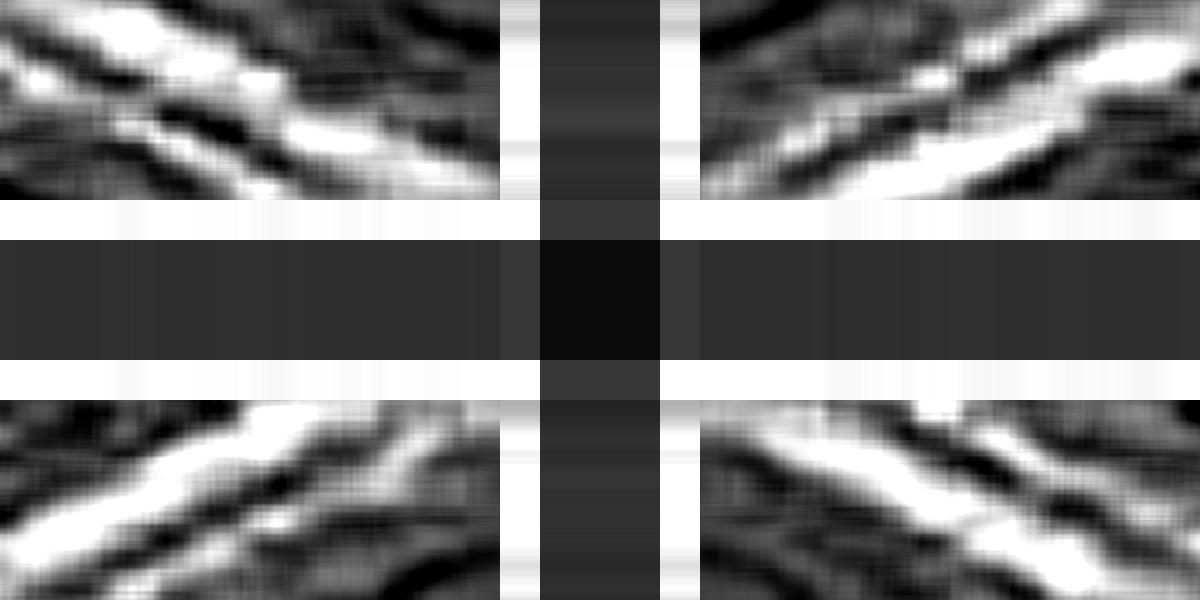}}
\caption{Image compression for grayscale flag image with varying ranks.}
\label{fig:svdd-pseudoskeleton}
\end{figure}

We finally compare the low-rank approximations between the SVD and the pseudoskeleton with varying ranks. 
Figure~\ref{fig:svdd-pseudoskeleton} shows the differences in compression at ranks of 90, 30, and 10, respectively. We observe that the SVD performs well at ranks of 90 and 30. 
In contrast, the pseudoskeleton-approximation effectively compresses  the black horizontal and vertical lines in the image but struggles to preserve the finer details of the flag. Figure~\ref{fig:svd_skeleton_fnorm} presents a comparison of the reconstruction errors between the SVD and the pseudoskeleton-approximation, measured by the Frobenius norm ranging from rank $1$ to $100$, where we find in all cases, the truncated SVD outperforms the pseudoskeleton in terms of the Frobenius norm. 
Similar trends can be observed when applied to the spectral norm.

\index{Matrix multiplication}
\index{Coordinate transformation}

\section{Coordinate Transformation in Matrix Decomposition}\label{section:coordinate-transformation}

Consider a vector $\bv\in \real^3$ with elements $\bv = [3 ,7, 2]^\top$. 
It is essential to clarify the significance of these values: In the Cartesian coordinate system, they represent a component of 3 along the $x$-axis, a component of 7 along the $y$-axis, and a component of 2 along the $z$-axis.
These scalar values are the \textit{coordinates} of $\bv$ with respect to the basis of the Cartesian system.
Matrix multiplication, on the other hand, gains significance when applied in high-dimensional spaces.

\paragraph{Coordinate defined by a nonsingular matrix.} Suppose we have a $3\times 3$ nonsingular matrix $\bB$, which  is invertible and possesses linearly independent columns. 
Consequently, the three columns of $\bB$ collectively form a basis for the $\real^{3}$ space. 
Taking a step further, the three columns of $\bB$ can serve as the basis for a \textcolor{black}{\textbf{new coordinate system}}, referred to as the \textcolor{black}{\textbf{$B$ coordinate system}}.

Returning to the Cartesian coordinate system, we also have a set of three vectors forming a basis, denoted by $\{\be_1, \be_2, \be_3\}$. 
If we arrange the three vectors as columns in a matrix, this matrix will be the identity matrix. 
Therefore, when we multiply a vector $\bv$ by the identity matrix, denoted by $\bI\bv$, we are essentially performing a coordinate transformation that leaves $\bv$ in the same coordinate system.
In other words, $\bI\bv = \bv$ means \textcolor{black}{\textbf{transferring $\bv$ from the Cartesian coordinate system into the Cartesian coordinate system}}, the same coordinate.

Similarly, when we multiply a vector $\bv$ by the matrix $\bB$, denoted by $\bB\bv$, \textcolor{black}{\textbf{we are transforming $\bv$ from the Cartesian coordinate system into the $B$ coordinate system}}. 
To illustrate this with a specific example, consider $\bv = [3, 7, 2]^\top$ and $\bB=[\bb_1, \bb_2, \bb_3]$. In this case, we have $\bu=\bB\bv = 3\bb_1+7\bb_2+2\bb_3$, i.e., vector $\bu$ contains 3 units of the first basis $\bb_1$ of $\bB$, 7 units of the second basis $\bb_2$ of $\bB$, and 2 units of the third basis $\bb_3$ of $\bB$. 
Now, if we wish to transform the vector $\bu$ from the $B$ coordinate system back to the Cartesian coordinate system, we can achieve this by multiplying $\bu$ by the inverse of $\bB$, denoted by $\bB^{-1}$. This operation results in $\bB^{-1}\bu = \bv$.

\index{Geometrical interpretation}
\paragraph{Coordinate defined by an orthogonal matrix.} A $3\times 3$ orthogonal matrix $\bQ$ defines a ``better" coordinate system since its three columns, forming the basis, are mutually orthonormal  (same as those in the Cartesian coordinate system). 
The operation $\bQ\bv$ facilitates the transition of $\bv$ from the Cartesian coordinate system to the one defined by the orthogonal matrix. 
Since the basis vectors from the orthogonal matrix exhibit orthonormality, just like the three vectors $\be_1, \be_2, \be_3$ in the Cartesian coordinate system, the transformation induced by the orthogonal matrix involves rotating or reflecting the Cartesian system.
To revert to the Cartesian coordinate system, one can utilize $\bQ^{-1}=\bQ^\top$.

\begin{figure}[h]
\centering
\includegraphics[width=0.99\textwidth]{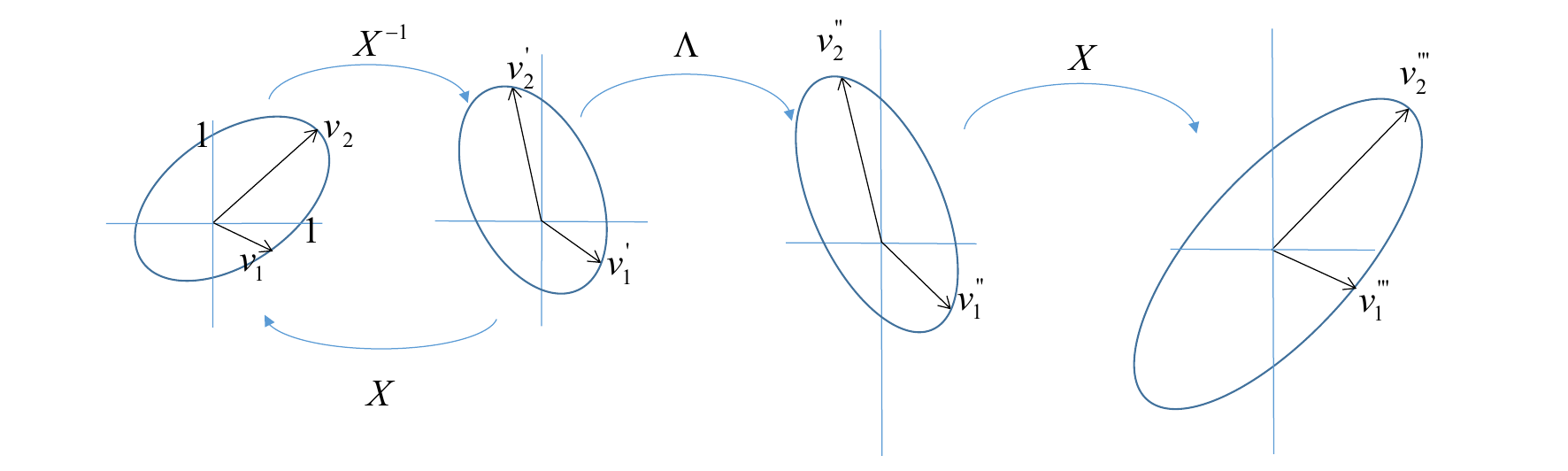}
\caption{Eigenvalue decomposition $\bA = \bX\bLambda\bX^{-1}$: $\bX^{-1}$ undergoes a transformation into a different coordinate system, followed by stretching with $\bLambda$, and then transforming back with $\bX$. 
	$\bX^{-1}$ and $\bX$ are nonsingular, which will change the basis of the system, and the angle between the vectors $\bv_1$ and $\bv_2$ will \textbf{not} be preserved. In other words, the angle between $\bv_1$ and $\bv_2$ is \textbf{different} from the angle between $\bv_1^\prime$ and $\bv_2^\prime$. The lengths of $\bv_1$ and $\bv_2$ are also \textbf{not} preserved; that is, $\normtwo{\bv_1} \neq \normtwo{\bv_1^\prime}$ and $\normtwo{\bv_2} \neq \normtwo{\bv_2^\prime}$.}
\label{fig:eigen-rotate}
\end{figure}

\section*{Eigenvalue Decomposition}
A square matrix $\bA$ with linearly independent eigenvectors can be factored as $\bA = \bX\bLambda\bX^{-1}$, where $\bX$ and $\bX^{-1}$ are nonsingular so that they define a system transformation inherently. 
The operation $\bA\bu = \bX\bLambda\bX^{-1}\bu$ firstly transfers $\bu$ into the coordinate system defined by $\bX^{-1}$, which we shall refer to as the \textit{eigen coordinate system}. 
Subsequently, the operation $\bLambda(\cdot)$  stretches each component of the vector in the eigen system by the length of the corresponding eigenvalue. 
Finally, $\bX$ facilitates the transformation of the resultant vector back to the Cartesian coordinate system. 
A visual representation of the coordinate system transformation via eigenvalue decomposition is presented in Figure~\ref{fig:eigen-rotate}, where $\bv_1$ and $\bv_2$ are two linearly independent eigenvectors of $\bA$ such that they form a basis for $\real^2$.

\begin{figure}[h]
\centering
\includegraphics[width=0.99\textwidth]{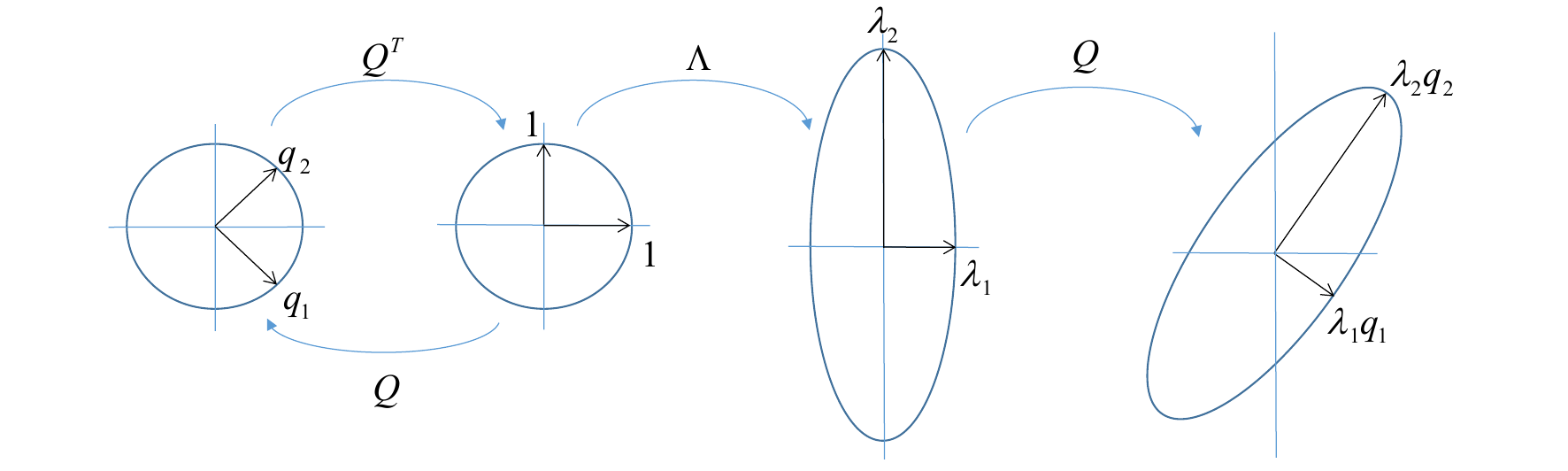}
\caption{Spectral decomposition $\bQ\bLambda \bQ^\top$: $\bQ^\top$ rotates or reflects, $\bLambda$ stretches the cycle to an ellipse, and $\bQ$ rotates or reflects back. Orthogonal matrices $\bQ^\top$ and $\bQ$ only change the basis of the system. However, they preserve both the angle between the vectors $\bq_1$ and $\bq_2$, and their lengths.}
\label{fig:spectral-rotate}
\end{figure}

\section*{Spectral Decomposition}
A symmetric matrix $\bA$ can be decomposed as $\bA = \bQ\bLambda\bQ^\top$, where $\bQ$ and $\bQ^\top$ are orthogonal matrices so that they define a coordinate system transformation inherently as well. The operation $\bA\bu = \bQ\bLambda\bQ^\top\bu$ firstly rotates or reflects $\bu$ into the coordinate system defined by $\bQ^\top$, which we shall refer to as the \textit{spectral coordinate system}. 
The operation $\bLambda(\cdot)$ stretches each component of the vector in the spectral system by the length of the corresponding eigenvalue. 
Subsequently, $\bQ$ facilitates the rotation or reflection of the resultant vector back to the original coordinate system. 
A demonstration of how the spectral decomposition transforms between coordinate systems in $\real^2$ space is shown in Figure~\ref{fig:spectral-rotate}, where $\bq_1$ and $\bq_2$ represent two linearly independent eigenvectors of $\bA$ such that they form a basis for $\real^2$. The coordinate transformation in the spectral decomposition is similar to that in the eigenvalue decomposition, with the distinction that in the spectral decomposition,  orthogonal vectors transformed by $\bQ^\top$ remain orthogonal. 
This is also a property of orthogonal matrices. That is, orthogonal matrices can be viewed as matrices, which change the basis of other matrices while preserving the angle (inner product) between  vectors:
$$
\bu^\top \bv = (\bQ\bu)^\top(\bQ\bv).
$$
The invariance of the angle between vectors also relies on the invariance of their lengths:
$$
\normtwo{\bQ\bu}= \normtwo{\bu}.
$$

\section*{SVD}
\begin{figure}[h]
\centering
\includegraphics[width=0.99\textwidth]{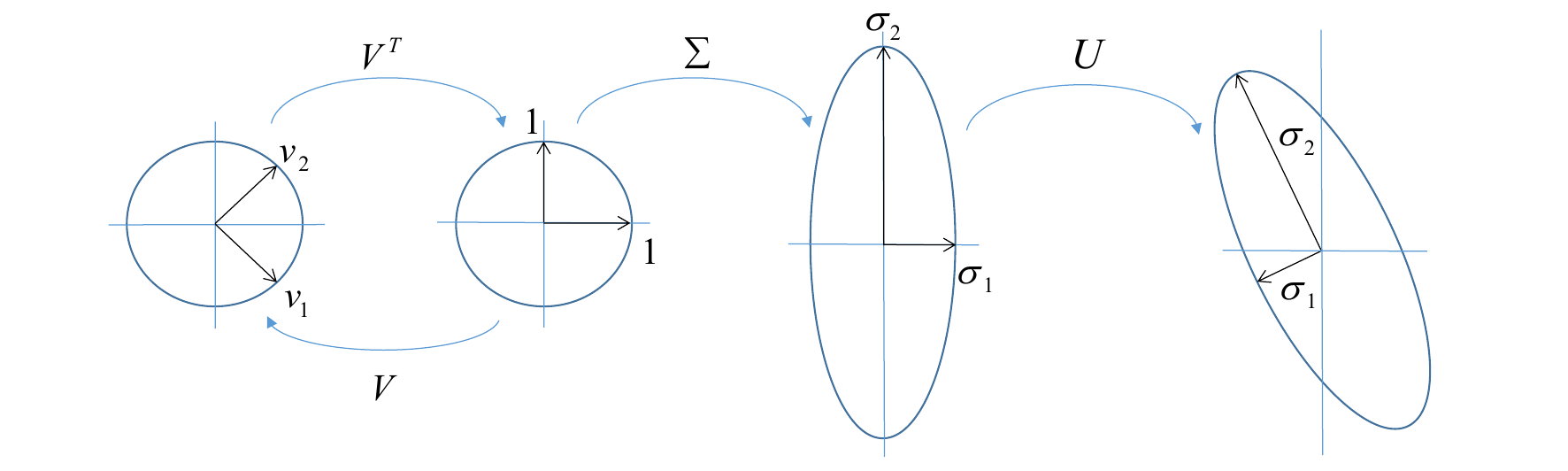}
\caption{SVD $\bA=\bU\bSigma\bV^\top$: $\bV^\top$ and $\bU$ rotate or reflect, $\bSigma$ stretches the circle to an ellipse. Orthogonal matrices $\bV^\top$ and $\bU$ only change the basis of the system. However, they preserve both the angle between the vectors $\bv_1$ and $\bv_2$, and their lengths.}
\label{fig:svd-rotate}
\end{figure}
\begin{figure}[h]
\centering
\includegraphics[width=0.99\textwidth]{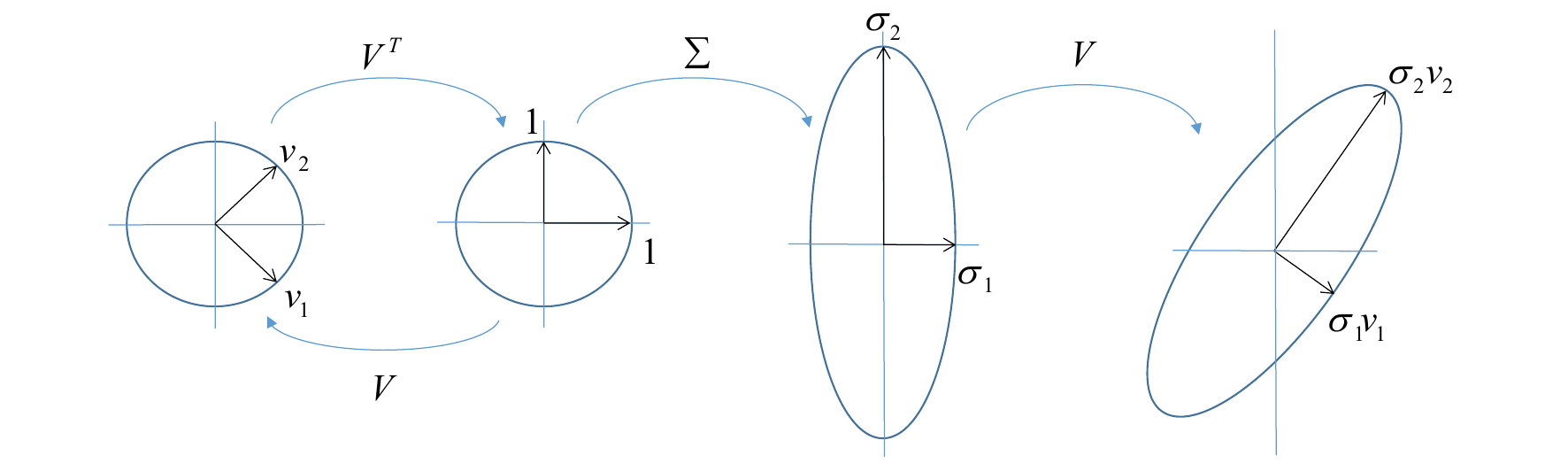}
\caption{$\bV\bSigma \bV^\top$ from SVD or polar decomposition: $\bV^\top$ rotates or reflects, $\bSigma$ stretches the cycle to an ellipse, and $\bV$ rotates or reflects back. Orthogonal matrices $\bV^\top$ and $\bV$ only change the basis of the system. However, they preserve both the angle between the vectors $\bv_1$ and $\bv_2$, and their lengths.}
\label{fig:polar-rotate}
\end{figure}
Any $m\times n$ matrix with rank $r$ can be factored as $\bA=\bU\bSigma\bV^\top$, which represents the SVD. The operation $\bA\bu=\bU\bSigma\bV^\top\bu$ then firstly rotates or reflects vector $\bu$ into the system defined by $\bV^\top$, which we refer to as the \textit{$V$ coordinate system}. $\bSigma$ stretches the first $r$ components of the resulting vector in the $V$ system by the lengths of the singular values. If $n\geq m$, then $\bSigma$ only keeps  $m-r$ additional components, which are scaled to zero, while removing the final $n-m$ components. If $m>n$, then $\bSigma$ scales $n-r$ components to zero and also adds  $m-n$  additional zero components. Finally, $\bU$ rotates or reflects the resulting vector into the \textit{$U$ coordinate system} defined by $\bU$. 
A visual demonstration of how the SVD transforms in a $2\times 2$ example is shown in Figure~\ref{fig:svd-rotate}. Further, Figure~\ref{fig:polar-rotate} demonstrates the transformation of $\bV\bSigma \bV^\top$ by a $2\times 2$ example. 
Similar to the spectral decomposition, orthogonal matrices $\bV^\top$ and $\bU$ only change the basis of the system but preserve the angle between vectors $\bv_1$ and $\bv_2$.

\begin{figure}[h]
\centering
\includegraphics[width=0.99\textwidth]{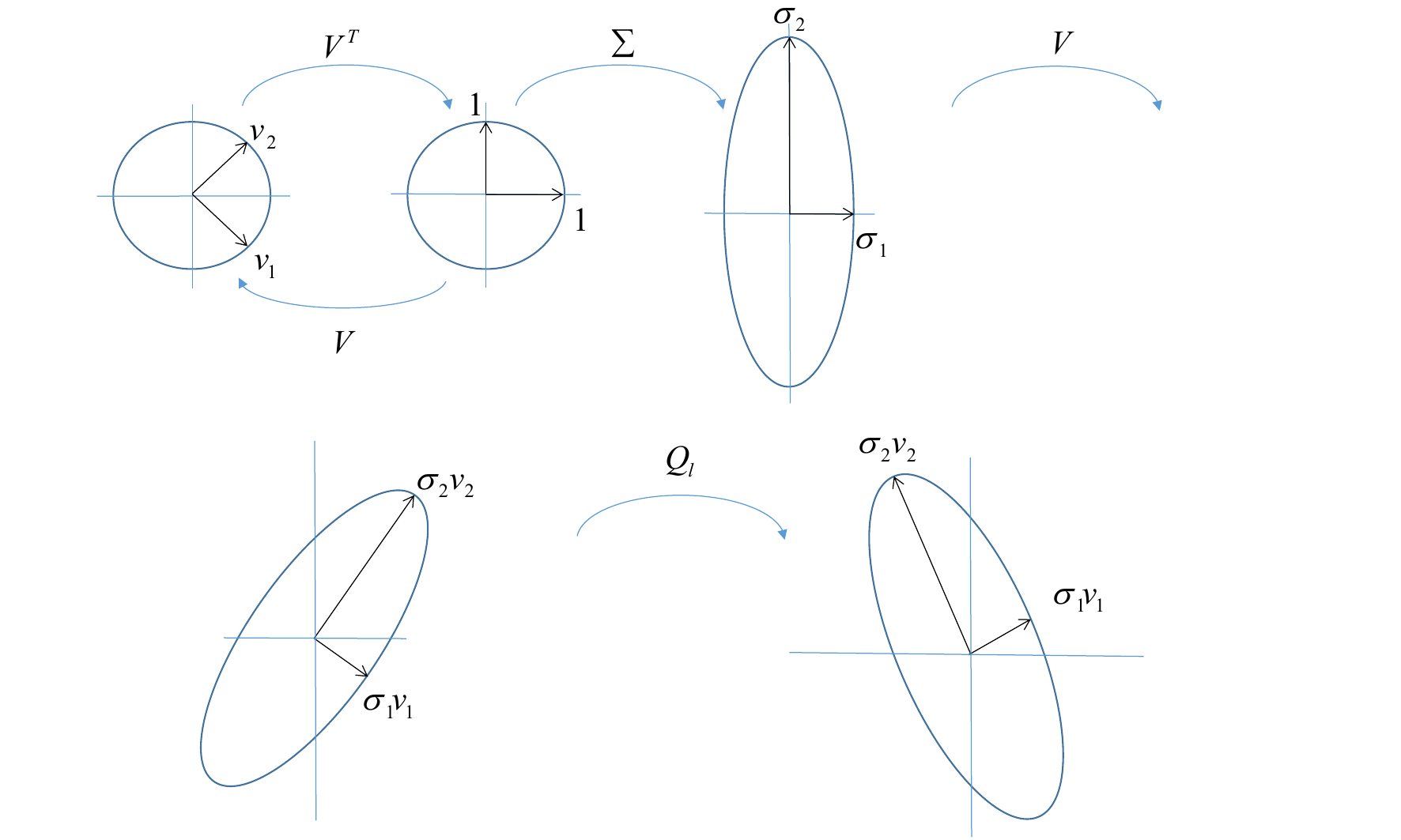}
\caption{Polar decomposition $\bA=\bQ_l\bS$: $\bV^\top$ rotates or reflects, $\bSigma$ stretches the cycle to an ellipse, and $\bV$ rotates or reflects back. Orthogonal matrices $\bV^\top$, $\bV$, and $\bQ_l$ only change the basis of the system. However, they preserve both the angle between the vectors $\bv_1$ and $\bv_2$, and their lengths.}
\label{fig:polar-rotate2}
\end{figure}

\section*{Polar Decomposition}
Any square matrix $\bA\in\real^{n\times n}$ can be factored as the left polar decomposition $\bA = (\bU\bV^\top)( \bV\bSigma \bV^\top) = \bQ_l\bS$. Similarly, the operation $\bA\bu = \bQ_l( \bV\bSigma \bV^\top)\bu$ transforms $\bu$ into the system defined by $\bV^\top$, and stretch each component by the lengths of the corresponding singular values. 
Subsequently, the resulting vector is transferred back into the Cartesian coordinate system by $\bV$. Finally, $\bQ_l$ will rotate or reflect the resulting vector from the Cartesian coordinate system into the $Q$ system defined by $\bQ_l$.
The right polar decomposition carries a similar interpretation. 
Similar to the spectral decomposition, orthogonal matrices $\bV^\top$ and $\bV$ only change the basis of the system but preserve the angle between the vectors $\bv_1$ and $\bv_2$.

\index{Least squares}
\section{Application: Least Squares via SVD for Rank-Deficient Matrices}\label{section:application-ls-svd}
The least squares problem is discussed in Section~\ref{section:application-ls-qr}. As a recap, 
let's consider the overdetermined system $\bA\bx = \bb$, where $\bA\in \real^{m\times n}$ is the data matrix, and $\bb\in \real^m$ with $m\geq n$ is the observation matrix. 
Normally, the data matrix $\bA$ will have full column rank  since real-world data often has little correlation. 
And the least squares (LS) solution is given by $\bx_{LS} = (\bA^\top\bA)^{-1}\bA^\top\bb$ for minimizing $\normtwo{\bA\bx-\bb}^2$, where $\bA^\top\bA$ is invertible since $\bA$ has full column rank and $\rank(\bA^\top\bA)=\rank(\bA)$.

However, when $\bA$ does not possess full column rank, $\bA^\top\bA$ becomes non-invertible. 
In such cases, we can employ the singular value decomposition of $\bA$ to solve the least squares problem. This approach is elucidated in the ensuing theorem.
\begin{theorem}[LS via SVD for Rank-Deficient Matrix]\label{theorem:svd-deficient-rank}
Let $\bA\in \real^{m\times n}$, and let $\bA=\bU\bSigma\bV^\top$ represent its full SVD decomposition, where $\bU\in\real^{m\times m}$ and $\bV\in \real^{n\times n}$ are orthogonal matrices, and $\rank(\bA)=r$. 
Let $\bU=[\bu_1, \bu_2, \ldots, \bu_m]$ and $\bV=[\bv_1, \bv_2, \ldots, \bv_n]$ be the column partitions of $\bU$ and $\bV$, respectively, and let  $\bb\in \real^m$ be the observation vector. 
Then the least squares solution with the minimal $\ell_2$ norm to $\bA\bx=\bb$ is given by 
\begin{equation}\label{equation:svd-ls-solution}
\bx_{LS} = \sum_{i=1}^{r}\frac{\bu_i^\top \bb}{\sigma_i}\bv_i = \bV\bSigma^+\bU^\top \bb, 
\end{equation}
where the upper-left side of $\bSigma^+ \in \real^{n\times m}$ is a diagonal matrix, and 
$\bSigma^+ = 
\footnotesize
\begin{bmatrix}
\bSigma_1^+ & \bzero \\
\bzero & \bzero
\end{bmatrix}$ with $\bSigma_1^+=\diag(\frac{1}{\sigma_1}, \frac{1}{\sigma_2}, \ldots, \frac{1}{\sigma_r})$.
\end{theorem}

\begin{proof}[of Theorem~\ref{theorem:svd-deficient-rank}]
Write out the loss to be minimized:
$$
\begin{aligned}
\normtwo{\bA\bx-\bb}^2 &= (\bA\bx-\bb)^\top(\bA\bx-\bb)
\stackrel{*}{=}\normtwo{\bU^\top \bA \bx-\bU^\top\bb}^2 
\stackrel{\dag}{=}\normtwo{\bU^\top \bA \bV\bV^\top \bx-\bU^\top\bb}^2 \\
&=\normtwo{\bSigma\balpha - \bU^\top\bb}^2 
\stackrel{+}{=}\sum_{i=1}^{r}(\sigma_i\alpha_i - \bu_i^\top\bb)^2 +\sum_{i=r+1}^{m}(\bu_i^\top \bb)^2,
\end{aligned}
$$
where equality ($*$) follows from the invariance under orthogonal transformations, equality ($\dag$) is due to orthogonality of $\bV$, $\balpha\triangleq\bV^\top \bx$, and equality (+) is due to $\sigma_{r+1}=\sigma_{r+2}= \ldots= \sigma_m=0$.
Since $\bx$ only appears in $\balpha$, we minimize the loss by setting $\alpha_i = \frac{\bu_i^\top\bb}{\sigma_i}$ for all $i\in \{1, 2, \ldots, r\}$. 
Any values assigned to $\alpha_{r+1}, \alpha_{r+2}, \ldots, \alpha_{n}$ will not affect the result. From the regularization perspective (or to achieve the minimal $\ell_2$ norm), these values can be set to 0 (see  Problem~\ref{prob:als_pseudo1}$\sim$\ref{prob:als_pseudon} for more insights). This approach yields the least squares solution via SVD:
$$
\bx_{LS} = \sum_{i=1}^{r}\frac{\bu_i^\top \bb}{\sigma_i}\bv_i=\bV\bSigma^+\bU^\top \bb \triangleq \bA^+\bb,
$$
where $\bA^+\triangleq\bV\bSigma^+\bU^\top\in \real^{n\times m}$ is known as the \textit{pseudo-inverse} of $\bA$. 
For more details on the pseudo-inverse, including the proof that the column space of $\bA^+$ is equal to the row space of $\bA$, and the row space of $\bA^+$ is equal to the column space of $\bA$, please refer to Appendix~\ref{appendix:pseudo-inverse}.\index{Pseudo-inverse}
\end{proof}

\begin{figure}[h!]
	\centering
	\includegraphics[width=0.98\textwidth]{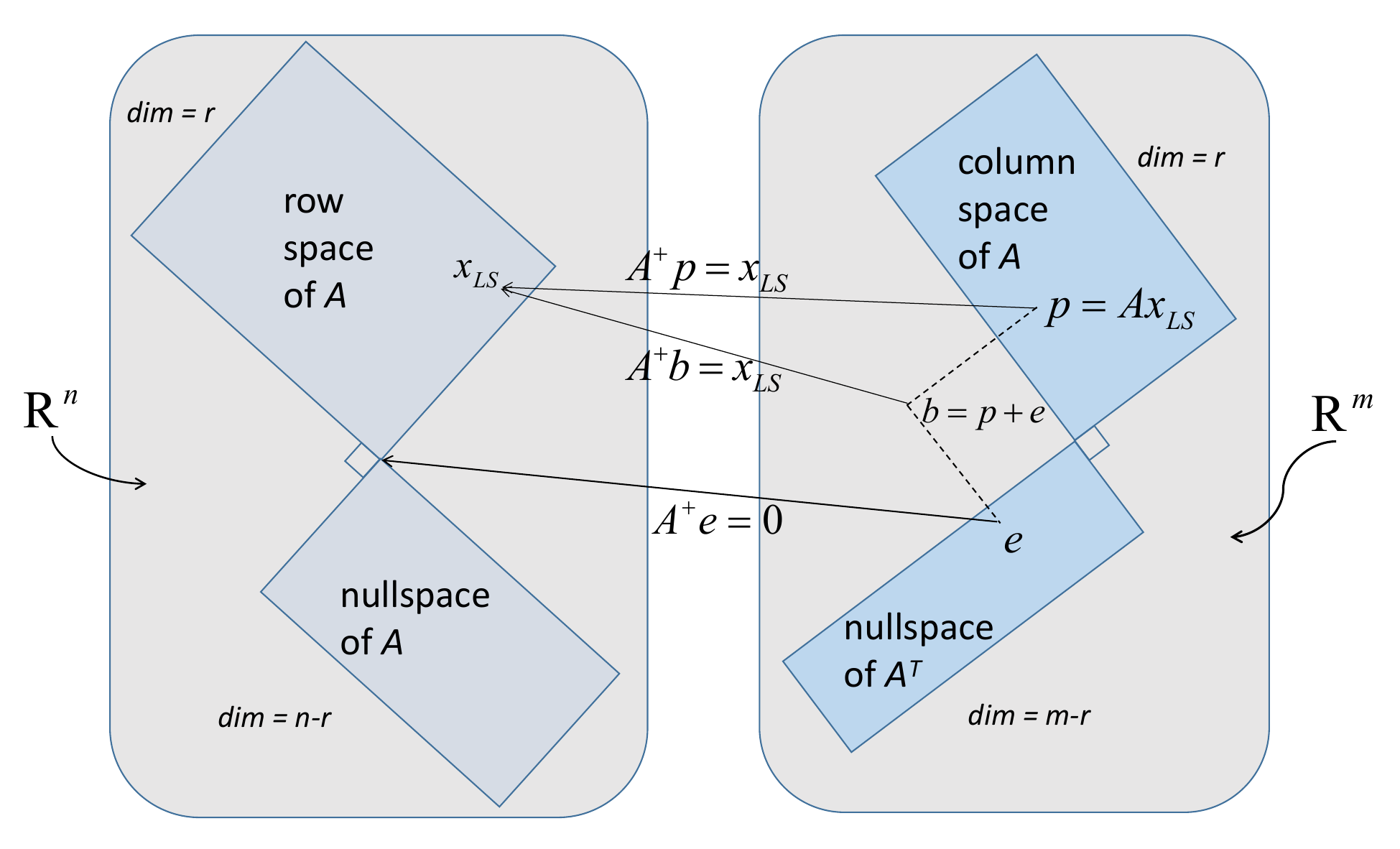}
	\caption{$\bA^+$: Pseudo-inverse of $\bA$.}
	\label{fig:lafundamental4-LS-SVD}
\end{figure}

Let $\bA^+ \triangleq \bV\bSigma^+\bU^\top$ be the pseudo-inverse of $\bA$. The pseudo-inverse $\bA^+$ agrees with $\bA^{-1}$ when $\bA$ is invertible. 
In the context of the least squares solution, the goal is to minimize the mean squared error by making the error $\bb-\bA\bx$ as small as possible.
Since $\bA\bx$ is a linear combination of the columns of $\bA$ and it  always resides within the column space of $\bA$, 
the optimal choice is the point  $\bp$ in the column space of $\bA$ that is closest to $\bb$.
This point is the projection $\bp$ of $\bb$ onto the column space of $\bA$. 
Then the error vector $\be=\bb-\bp$ has the minimal length. In other words, the optimal combination $\bp = \bA\bx_{LS}$ is the projection of $\bb$ onto the column space, while the error $\be$ is orthogonal to the column space of $\bA$ (see Appendix~\ref{section:by-geometry-hat-matrix} and Appendix~\ref{appendix:pseudo-inverse} for more details). Therefore, according to the fundamental theorem of linear algebra, the error $\be=\bb-\bA\bx_{LS}$ lies within the null space of $\bA^\top$:
$$
\bA^\top(\bb-\bA\bx_{LS}) = \bzero  \qquad \text{or} \qquad \bA^\top\bb=\bA^\top\bA\bx_{LS},
$$
which is also referred to  as the \textit{normal equation} for least squares. The relationship between $\be$ and $\bp$ is shown in Figure~\ref{fig:lafundamental4-LS-SVD}, where $\bb$ is decomposed into $\bp+\be$. Since $\be$ is in $\nspace(\bA^\top)$ and is orthogonal to $\cspace(\bA)$, and as shown in Section~\ref{section:property-svd}, $\{\bu_1,\bu_2, \ldots,\bu_r\}$ is an orthonormal basis of $\cspace(\bA)$, then the first $r$ components of $\bU^\top\be$ are all zeros. Therefore, $\bA^+\be = \bV\bSigma^+\bU^\top\be=\bzero$. Moreover, $\bx_{LS}=\bA^+\bb = \bA^+(\bp+\be) = \bA^+\bp$.

Additionally, it was demonstrated in Section~\ref{section:property-svd} that $\{\bv_1, \bv_2, \ldots, \bv_r\} $ forms an orthonormal basis of $\cspace(\bA^\top)$. 
Consequently, $\bx_{LS} = \sum_{i=1}^{r}\frac{\bu_i^\top \bb}{\sigma_i}\bv_i$ resides within the row space of $\bA$. In other words, $\bx_{LS}$ cannot be split into a combination of two components that are in the row space  and  null space of $\bA$, respectively.

Apart from this least squares solution obtained from SVD, in practice, a direct solution of the normal equations can lead to numerical difficulties when $\bA^\top\bA$ is nearly singular. In particular, when two or more  columns in $\bA^\top\bA$ are co-linear, the resulting parameter values can have substantial magnitudes. Such near-degeneracies will not be uncommon when dealing with real data sets. The resulting numerical difficulties can be addressed using  SVD as well; see Section~\ref{section:pre_ls} or  \citet{bishop2006pattern} for more information.

\index{Least squares}
\index{Truncated}
\index{Truncated SVD}
\index{Norm ratio}
\section{Application: Least Squares with Norm Ratio Method}\label{section:application-ls-svd-norm-ratio}

Continuing from the previous section's setup, let $\bA_k \in \real^{m\times n}$ be the rank-$k$ approximation to the original $m\times n$ matrix $\bA$. Define the \textit{Frobenius norm ratio} \citep{zhang2017matrix} as 
$$
\nu(k) = \frac{\normf{\bA_k}}{\normf{\bA}} = \frac{\sqrt{\sigma_1^2+\sigma_2^2+\ldots +\sigma_k^2}}{\sqrt{\sigma_1^2+\sigma_2^2+\ldots +\sigma_p^2}}, \quad p = \min\{m,n\},
$$
where 
$\bA_k$ is the truncated SVD of $\bA$ containing the largest $k$ terms, i.e., $\bA_k = \sum_{i=1}^{k} \sigma_i\bu_i\bv_i^\top$ from the SVD of $\bA=\sum_{i=1}^{p} \sigma_i\bu_i\bv_i^\top$. 
We determine the minimum integer $k$ satisfying 
$$
\nu(k) \geq \alpha
$$
as the \textit{effective rank estimate} $\hat{r}$, 
with $\alpha$ being the threshold, typically set at $\alpha=0.997$ (with a maximum value of 1). 
Once the effective rank $\hat{r}$ is established, it is substituted into Equation~\eqref{equation:svd-ls-solution}, yielding:
$$
\hat{\bx}_{LS} = \sum_{i=1}^{\textcolor{mylightbluetext}{\hat{r}}}\frac{\bu_i^\top \bb}{\sigma_i}\bv_i , 
$$
which can be regarded as an approximation to the LS solution $\bx_{LS}$. 
This solution corresponds to the LS solution of the linear equation $\bA_{\hat{r}}\bx = \bb$, where 
$
\bA_{\hat{r}} = \sum_{i=1}^{\hat{r}} \sigma_i \bu_i\bv_i^\top.
$
This filtering method is particularly useful when working with noisy matrices, such as in the processing of financial data for quantitative strategies.

\index{Least squares}
\section{Application: Perturbation Theory for  Least Squares Problems}
In Equation~\eqref{equation:qr_condition_num}, we derived the condition number for the linear system problem $\bA\bx=\bb$  as $\kappa(\bA) = \normtwo{\bA} \normtwo{\bA^{-1}}$. 
Using the SVD and the definition of the spectral norm, if $\bA\in\real^{m\times n}$ has full rank with $m\geq n$, we have 
$$
\kappa(\bA) = \frac{\sigma_1(\bA)}{\sigma_n(\bA)},
$$
where $\sigma_i(\bA)$ represents the $i$-th largest singular value of $\bA$.
The following perturbation result holds:
\begin{theorem}[Perturbation of Least Squares \citep{wedin1973perturbation}]
Let $\bA\in \real^{m\times n}$ be any matrix with full rank and $m\geq n$. Consider the least squares problem $\mathop{\min}_{\bx} \normtwo{\bA\bx-\bb}^2$.
Suppose that  perturbations $\delta\bA$ and $\delta\bb$ are such that 
$$
\epsilon = \frac{\normtwo{\delta\bA}}{\normtwo{\bA}}, \gap 
\eta = \frac{\normtwo{\delta\bA}}{\sigma_n} = \kappa(\bA)\epsilon,
$$
where $\sigma_n$ denotes the smallest nonzero singular value of $\bA$. 
Consequently, the perturbed matrix $\bA+\delta\bA$ has full rank, and the absolute error bound of $\delta\bx$ is 
$$
\normtwo{\delta\bx} \leq 
\frac{\kappa(\bA)}{1-\eta} 
\left( \epsilon \normtwo{\bx}
+
\frac{\normtwo{\delta\bb}}{\normtwo{\bA}}
+
\epsilon \kappa(\bA)\frac{\normtwo{\bb-\bA\bx}}{\normtwo{\bA}}
\right).
$$
\end{theorem}
The proof is available in  \citet{wedin1973perturbation}. We find that the condition number defined from the linear system problem $\bA\bx=\bb$ (Section~\ref{section:qr_condition}) can be applied to bound the least squares problem. 
This underscores the significance of the condition number's definition.
We observe that
\begin{itemize}
\item If $m=n$, having $\bA$ with full rank implies that $\bx$ has an exact solution $\bx=\bA^{-1}\bb$. 
In this scenario, it reduces to the linear system problem (Section~\ref{section:qr_condition}).
\item Furthermore, if we suppose that $\delta\bA=\bzero$, the absolute error bound reduces to Equation~\eqref{equation:qr_well_cond}.
\end{itemize}

\index{PCA}
\index{Principal component analysis}
\section{Application: PCA via  Spectral Decomposition and  SVD}\label{section:pca_isvd}
An important application of SVD is to apply principal component analysis (PCA). 
PCA is frequently employed  to identify patterns within data and to uncover the variance-covariance structure of the data, thereby achieving the following objectives:
\begin{enumerate}
\item \textit{Data reduction}. We reduce the dimensionality of the data using a smaller number of \textit{principal components}.
\item \textit{Interpretation}.  PCA can help reveal relationships that were not suspected previously.
\end{enumerate}

Dimensionality reduction is also beneficial in scenarios requiring lower-dimensional data, such as data visualization, storage, and computationally intensive tasks. 
Given a data set consisting of $n$ observations $\{\bx_1,\bx_2,\ldots,\bx_n\}$, where $\bx_i\in \real^p$ for all $i\in \{1,2,\ldots,n\}$, our goal is to project the data into a lower-dimensional space, say $m<p$. Define the sample mean vector and sample covariance matrix as follows:
$$
\overline{\bx} \triangleq \frac{1}{n}\sum_{i=1}^{n}\bx_i
\qquad 
\text{and}
\qquad 
\bS \triangleq \frac{1}{n-1}\sum_{i=1}^{n} (\bx_i - \overline{\bx})(\bx_i-\overline{\bx})^\top,
$$
where the divisor $n-1$  serves the purpose of making it to be an \textit{unbiased consistent estimator} of the  covariance matrix  \citep{lu2021rigorous}. 
Alternatively, the covariance matrix can also be defined as $\bS \triangleq \frac{1}{\textcolor{mylightbluetext}{n}}\sum_{i=1}^{n} (\bx_i - \overline{\bx})(\bx_i-\overline{\bx})^\top$, which is also a \textit{consistent estimator} of the covariance matrix \footnote{Consistency: An estimator $\theta_n $ of $\theta$ constructed on the basis of a sample of size $n$ is said to be consistent if $\theta_n\stackrel{p}{\rightarrow} \theta$   as $n \rightarrow \infty $.}.

Each data point $\bx_i$ is then projected onto a scalar value using a vector $\bu_1$ (see discussion below) such that the projection is given by  $\bu_1^\top\bx_i$. The mean of the projected data is obtained by $\Exp[\bu_1^\top\bx_i] = \bu_1^\top \overline{\bx}$, and the variance of the projected data is given by 
$$
\begin{aligned}
\Cov[\bu_1^\top\bx_i] &= \frac{1}{n-1} \sum_{i=1}^{n}( \bu_1^\top \bx_i - \bu_1^\top\overline{\bx})^2=
\frac{1}{n-1} \sum_{i=1}^{n}\bu_1^\top  ( \bx_i -\overline{\bx})( \bx_i -\overline{\bx})^\top\bu_1
=\bu_1^\top\bS\bu_1.
\end{aligned}
$$

\begin{figure}[h]
\centering   
\vspace{-0.35cm}  
\subfigtopskip=2pt  
\subfigbottomskip=2pt 
\subfigcapskip=-5pt 
\subfigure[Project onto y-axis.]{\label{fig:pca_cluster1}
\includegraphics[width=0.31\linewidth]{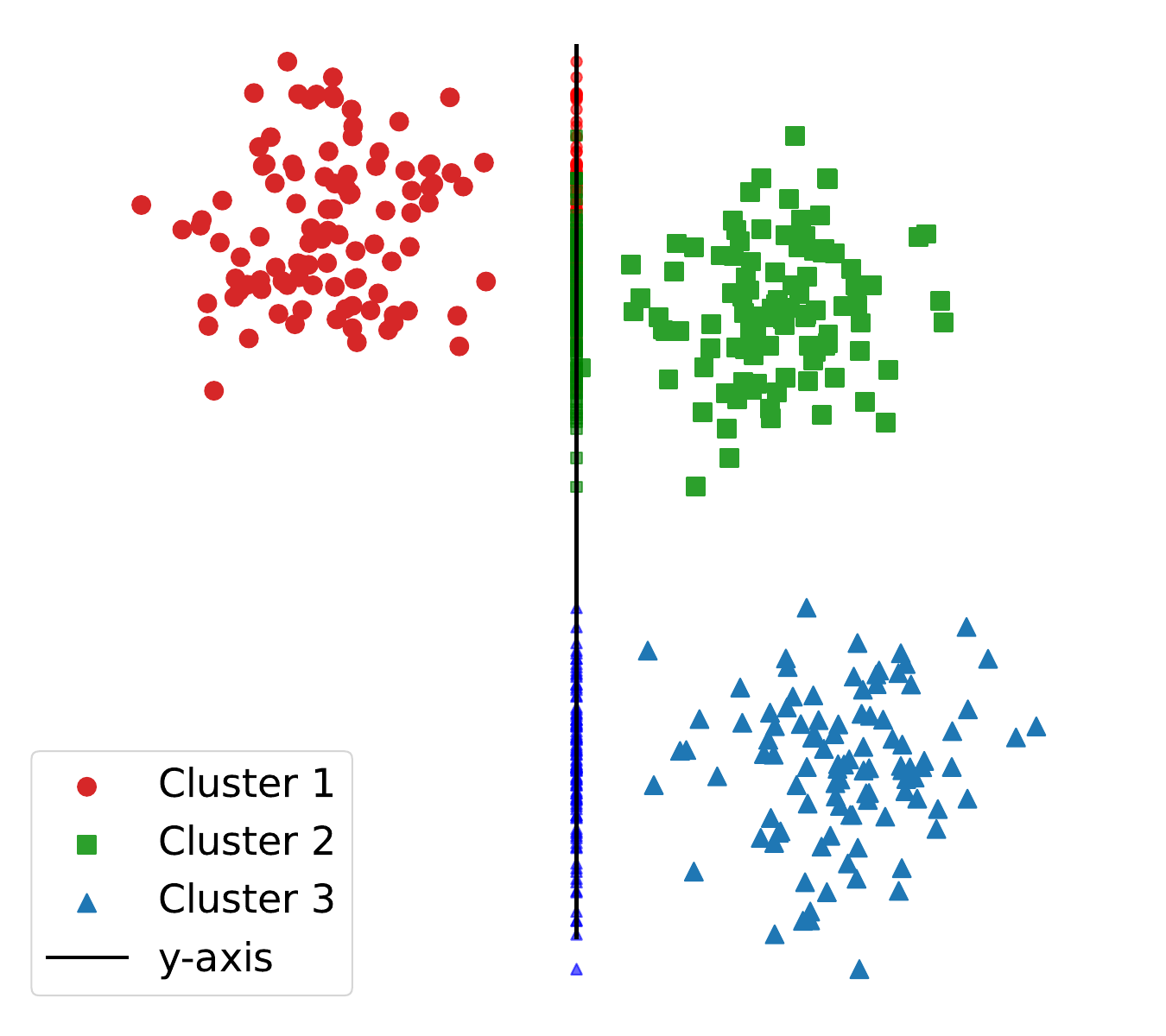}}
\subfigure[Project onto x-axis.]{\label{fig:pca_cluster2}
\includegraphics[width=0.31\linewidth]{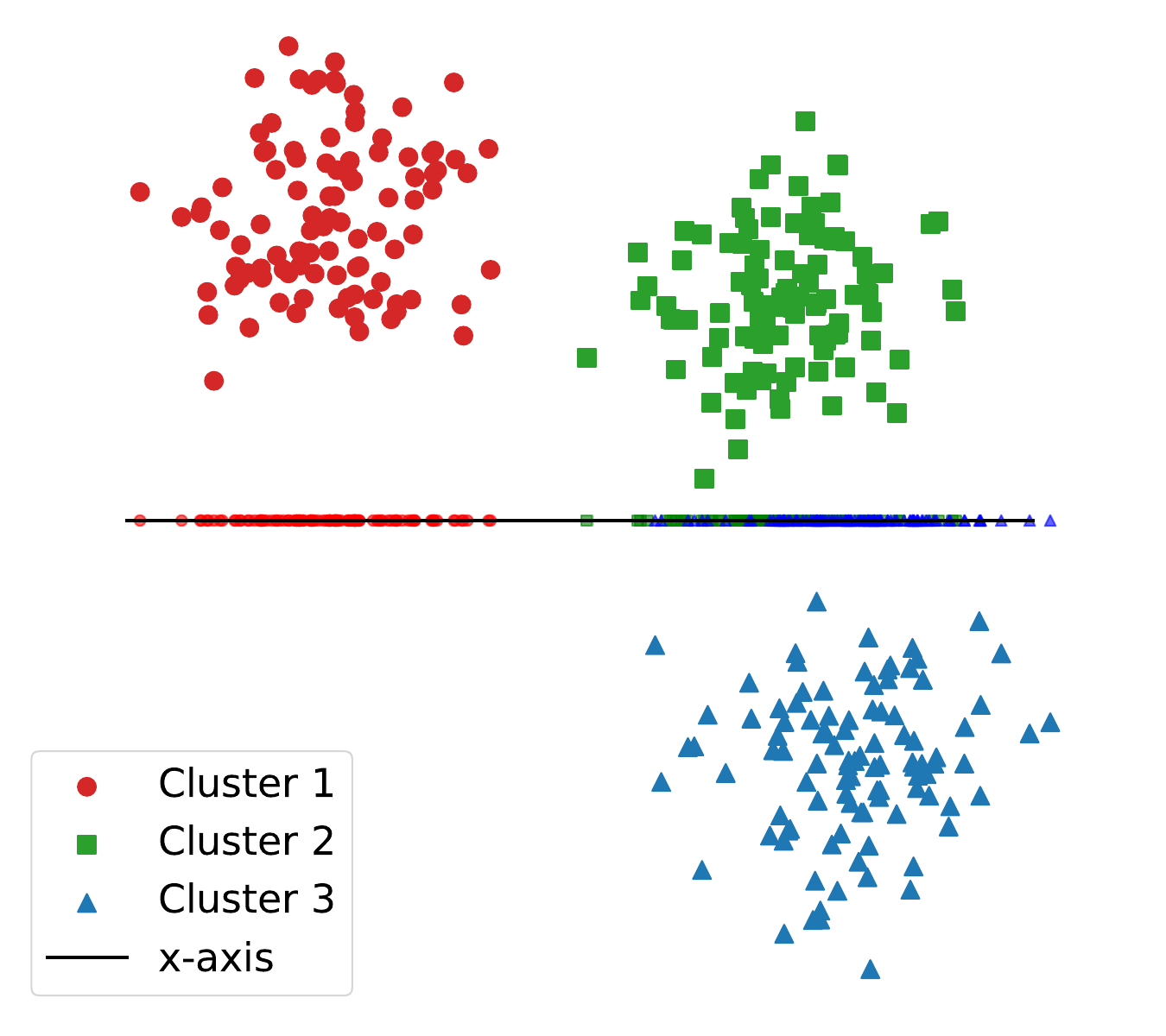}}
\subfigure[Project onto the first principal axis.]{\label{fig:pca_cluster3}
	\includegraphics[width=0.31\linewidth]{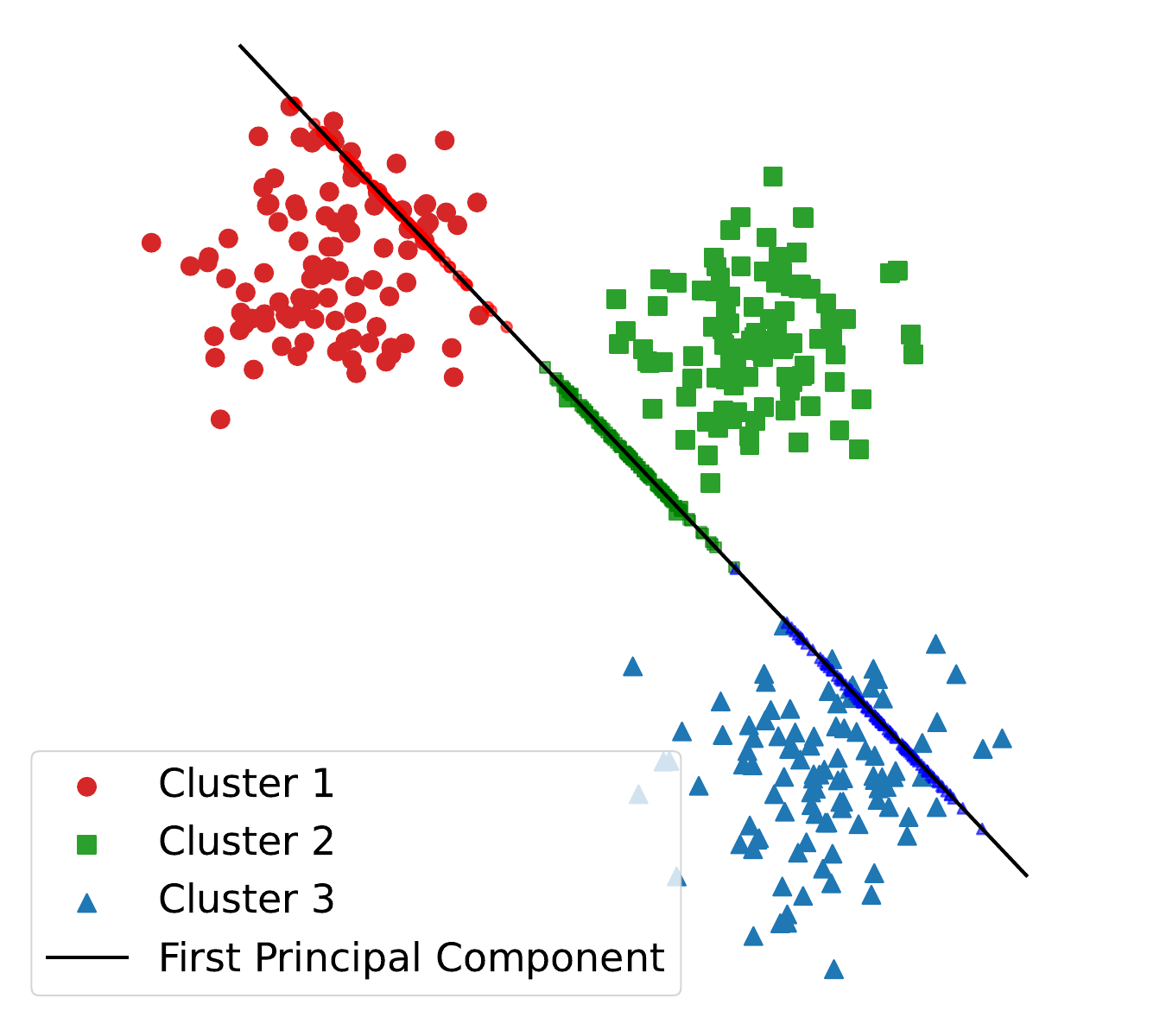}}
\caption{Dimension reduction of a two-dimensional data set that contains three clusters can lead to significant information loss when projecting onto either the x-axis or the y-axis. In contrast, projecting the data onto the first principal axis with maximal projected variance allows for a substantial preservation of the differentiation among the clusters.}
\label{fig:pca_clusters}
\end{figure}

The objective of PCA is to maximize the projected variance (see visual description in Figure~\ref{fig:pca_clusters}) $\bu_1^\top\bS\bu_1$ concerning $\bu_1$ to retain more information from the projected data. To achieve this, we need to impose a constraint on $\normtwo{\bu_1}$ to prevent it from approaching infinity, which is achieved by setting $\bu_1^\top\bu_1=1$. 
Employing Lagrange multipliers (see, for example,  \citet{bishop2006pattern, boyd2004convex}), we can maximize the following objective: 
\begin{equation}\label{equation:pca_raw1}
\bu_1^\top\bS\bu_1 + \lambda_1 (1 - \bu_1^\top\bu_1).
\end{equation}
This yields
$$
\bS\bu_1 = \lambda_1\bu_1 \implies \bu_1^\top\bS\bu_1 = \lambda_1.
$$
That is, $\bu_1$ is an eigenvector of $\bS$ corresponding to the eigenvalue $\lambda_1$. And the maximum variance projection $\bu_1$ corresponds to the largest eigenvalue of $\bS$. This eigenvector is known as the \textit{first principal axis}.

Define the other principal axes by the decremental eigenvalues until we have $m$ such principal components. The procedure brings about the desired dimension reduction and is known as the \textit{maximum-variance formulation} of PCA \citep{hotelling1933analysis, bishop2006pattern, shlens2014tutorial}.
An alternative approach, known as the \textit{minimum-error formulation} of PCA, is discussed in \citet{pearson1901liii, bishop2006pattern}.

\paragraph{Second perspective.} 
Alternatively put, 
now let's assume the data are centered such that the sample mean vector $\overline{\bx}$ is zero, or we can set $\bx_i \leftarrow \bx_i-\overline{\bx}$ to centralize the data, i.e., achieving data centering by subtracting the mean from each data point.
Our aim is to project data points $\{\bx_1, \bx_2, \ldots, \bx_n\}$ into a lower-dimensional space $\real^p\rightarrow \real^m$. Suppose $\bP\in\real^{p\times m}$ is a semi-orthogonal projection satisfying $\bP^\top\bP=\bI_m$, i.e., $\bP$ forms an orthogonormal basis for a low-dimensional subspace of $\real^p$.
Then $\bH=\bP\bP^\top$ defines an orthogonal projection (Definition~\ref{definition:orthogonal-projection-matrix}) onto the low-dimensional subspace defined by the column space $\mathcalV$ of $\bP$ (see Proposition~\ref{proposition:projection-from-matrix}). Then the projection of any data point $\bx_i$ onto the
subspace $\mathcalV$ is $\mathcal{P}_{\bP}(\bx_i) =\bP\bP^\top\bx_i$.
The principal component analysis framework finds the $\bP$ that   maximizes the variance.
It can be shown that the covariance matrix after projection is:
$$
\frac{1}{n-1} \sum_{i=1}^{n}\bP\bP^\top \bx_i  (\bP\bP^\top \bx_i  )^\top 
=\frac{1}{n-1} \bP\bP^\top \bX^\top\bX \bP\bP^\top,
$$
where $\bX\in\real^{n\times p}$ contains the data observations as rows. 
Since the variance is the trace of the covariance matrix, we can use the following optimization to find  $\bP$:
\begin{equation}\label{equation:pca_second}
\mathop{\max}_{\bP}\,\, \trace(\bP\bP^\top \bX^\top\bX \bP\bP^\top), \gap \text{s.t.} \gap \bP^\top\bP=\bI_m.
\end{equation}
By applying the cyclic property of traces, it follows that 
$$
\trace(\bP\bP^\top \bX^\top\bX \bP\bP^\top)
=\trace(\bP^\top \bX^\top\bX \bP).
$$
It  can then be shown that the columns of $\bP$ contain the eigenvalues of $\bX^\top\bX$ corresponding to the largest $m$ eigenvalues.

\begin{SCfigure}
\centering
\includegraphics[width=0.4\textwidth]{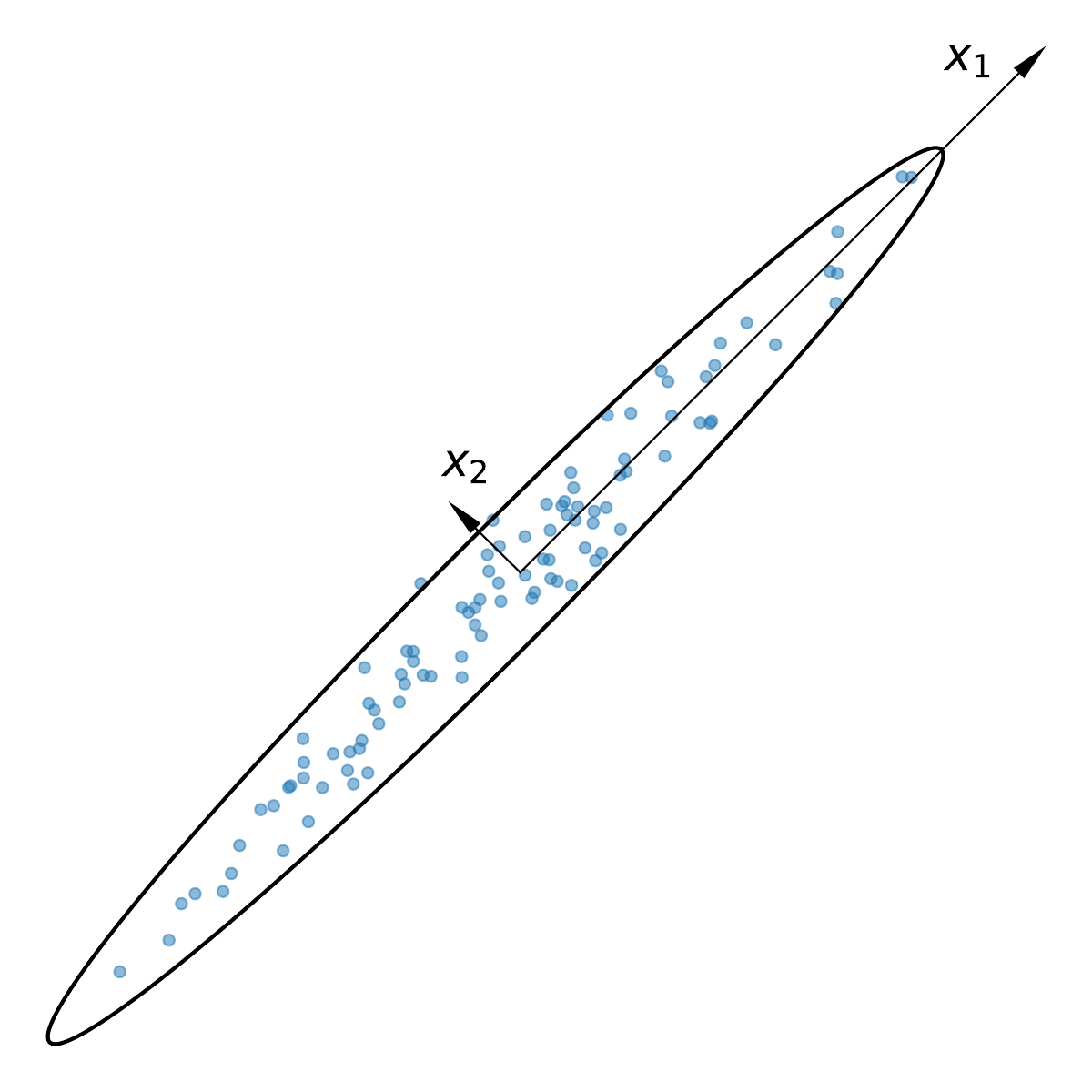}
\caption{Description of PCA in a two-dimensional case. $\bx_1$ and $\bx_2$ are the directions of corresponding eigenvectors of the covariance matrix. Therefore, $\bx_1$ encodes the first principal axis, and $\bx_2$ is the second principal axis.}
\label{fig:pca_axis_2d}
\end{SCfigure}

\paragraph{Third perspective: reconstruction.}
As shown above, the projection of data point $\bx_i$ onto the subspace $\mathcalV$ (the column space of $\bP$) is $\bP\bP^\top \bx_i$. Then the sum of squared loss is 
$$
\sum_{i=1}^{n} \normtwo{\bP\bP^\top \bx_i-\bx_i}^2
=\normf{\bP\bP^\top \bX^\top - \bX^\top}^2
=-\trace(\bP^\top\bX^\top\bX\bP)+\trace(\bX\bX^\top).
$$
Then the problem of minimizing the sum of squared loss is
\begin{equation}\label{equation:pca_third}
\mathop{\min}_{\bP}\,\, -\trace(\bP^\top\bX^\top\bX\bP)+\trace(\bX\bX^\top)\gap \text{s.t.} \gap \bP^\top\bP=\bI_m.
\end{equation}
It can be shown that problem~\eqref{equation:pca_second} is equivalent to problem~\eqref{equation:pca_third}. 
Therefore, we can conclude that PCA seeks to find the projection onto the subspace with maximum variance, which is essentially minimizing the reconstruction error of the projected points.
Figure~\ref{fig:pca_axis_2d} shows an example of a set of data points in two-dimensional space, where $\bx_1$ encodes the first principal axis, and $\bx_2$ represents the second principal axis.

\paragraph{PCA via the spectral decomposition.}
We again assume that  the data are centered such that the sample mean vector $\overline{\bx}$ is zero. 
Let the data matrix $\bX \in \real^{n\times p}$ contain the data observations as rows. The covariance matrix is given by 
\begin{equation}\label{equation:pca-equ0}
\bS= \frac{\bX^\top\bX}{n-1},
\end{equation}
which is a symmetric matrix, and its spectral decomposition is given by 
\begin{equation}\label{equation:pca-equ1}
\bS = \bU\bLambda\bU^\top,
\end{equation}
where $\bU$ is an orthogonal matrix of eigenvectors (the columns of $\bU$ are eigenvectors of $\bS$), and $\bLambda=\diag(\lambda_1, \lambda_2,\ldots, \lambda_p)$ is a diagonal matrix containing the eigenvalues (ordered such that $\lambda_1 \geq \lambda_2 \geq \ldots \geq \lambda_p$). The eigenvectors are called the \textit{principal axes} of the data, and they \textit{decorrelate} the  covariance matrix. Projections of the data onto the principal axes are called the \textit{principal components}. The $i$-th principal component is given by the $i$-th column of $\bX\bU$. 
If our objective is to reduce the dimension from $p$ to $m$, we simply select the first $m$ columns of $\bX\bU$.

\paragraph{PCA via the SVD.}
If the SVD of $\bX$ is given by $\bX = \bP\bSigma\bQ^\top$, then the covariance matrix can be expressed as 
\begin{equation}\label{equation:pca-equ2}
\bS= \frac{\bX^\top\bX}{n-1} = \bQ \frac{\bSigma^2}{n-1}\bQ^\top,
\end{equation}
where $\bQ\in \real^{p\times p}$ is an orthogonal matrix containing the right singular vectors of $\bX$, and the upper-left part of $\bSigma$ is a diagonal matrix containing the singular values $\diag(\sigma_1,\sigma_2,\ldots)$, with $\sigma_1\geq \sigma_2\geq \ldots$. The number of singular values is equal to $\min\{n,p\}$, which will not exceed $p$, and some  may be zero. 

The comparison between Equation~\eqref{equation:pca-equ2} and Equation~\eqref{equation:pca-equ1} suggests that Equation~\eqref{equation:pca-equ2} also represents a spectral decomposition of $\bS$.
This is due to the fact that both the eigenvalues in $\bLambda$ and singular values in $\bSigma$ are ordered in a descending order, and the spectral decomposition in terms of the eigenspaces is unique (Section~\ref{section:uniqueness-spectral-decomposition}). 

In other words, the right singular vectors $\bQ$ can also serve as the principal axes, which decorrelate the covariance matrix,
and the singular values are related to the eigenvalues of the covariance matrix through the relationship: $\lambda_i = \frac{\sigma_i^2 }{n-1}$. To reduce the dimensionality of the data from $p$ to $m$, we  select the largest $m$ singular values and their corresponding right singular vectors. This is also related to the truncated SVD (TSVD) $\bX_m = \sum_{i=1}^{m}\sigma_i \bp_i\bq_i^\top$, where $\bp_i$'s and $\bq_i$'s are the columns of $\bP$ and $\bQ$, respectively.

\index{Truncated}
\index{Truncated SVD}

\paragraph{A byproduct of PCA via the SVD for high-dimensional data.} For a principal axis $\bu_i$ of $\bS = \frac{\bX^\top\bX}{n-1}$, we have 
$$
\frac{\bX^\top\bX}{n-1} \bu_i = \lambda_i \bu_i.
$$
Premultiplying by $\bX$, we obtain 
$$
\frac{\bX\bX^\top}{n-1} (\bX\bu_i) = \lambda_i (\bX\bu_i),
$$
which implies that $\lambda_i$ is also an eigenvalue of $\frac{\bX\bX^\top}{n-1} \in \real^{n\times n}$, and the corresponding eigenvector is $\bX\bu_i$. This relationship is also stated in the proof of Theorem~\ref{theorem:reduced_svd_rectangular}, which establishes the existence of the SVD. If $p \gg n$, instead of finding the eigenvectors of $\bS=\frac{\bX^\top\bX}{n-1}$, i.e., the principal axes of $\bS=\frac{\bX^\top\bX}{n-1}$, we can find the eigenvectors of $\frac{\bX\bX^\top}{n-1}$. This reduces the computational complexity from $\mathcalO(p^3)$ to $\mathcalO(n^3)$, which is more efficient when $p \gg n$. 

Now, returning to the principal axes of $\bS=\frac{\bX^\top\bX}{n-1}$, suppose  the eigenvector of $\frac{\bX\bX^\top}{n-1}$ is $\bv_i$, corresponding to a nonzero eigenvalue $\lambda_i$:
$$
\frac{\bX\bX^\top}{n-1} \bv_i = \lambda_i \bv_i.
$$
Premultiplying by $\bX^\top$, we obtain 
$$
\frac{\bX^\top\bX}{n-1} (\bX^\top\bv_i) = \bS(\bX^\top\bv_i)   = \lambda_i (\bX^\top\bv_i),
$$
indicating that the eigenvector $\bu_i$ of $\bS$ is proportional to $\bX^\top\bv_i$, where $\bv_i$ is the eigenvector of $\frac{\bX\bX^\top}{n-1}$ corresponding to the same eigenvalue $\lambda_i$. Note that a further normalization step is needed to ensure $\normtwo{\bu_i}=1$.
Thus, when $p\gg n$,  we can also find the principal axes using the SVD of $\frac{\bX\bX^\top}{n-1}$.

\begin{figure}[H]
	\centering
	\includegraphics[width=0.99\textwidth]{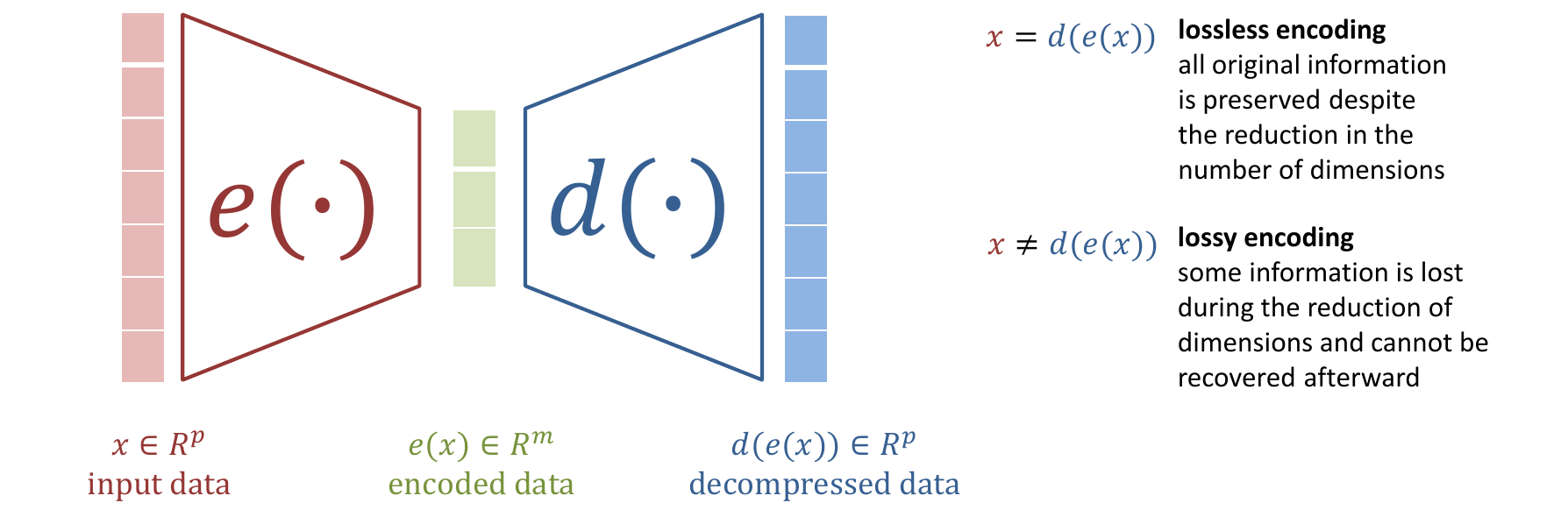}
	\caption{Description of an autoencoder.}
	\label{fig:autoencoder}
\end{figure}
\index{Autoencoder}
\subsection*{Autoencoder}
In machine learning, an autoencoder performs dimensionality reduction, which decreases the number of features describing a data set ($\bx$ in Figure~\ref{fig:autoencoder}). 
This process, known as the \textit{encoding process} ($e(\bx)$ in Figure~\ref{fig:autoencoder}), involves either selecting a subset of the existing features or extracting a smaller set of new features derived from the originals.
The \textit{decoding process}, denoted as $d(e(\bx))$, reverses the encoding to reconstruct the original data from its compressed form. Depending on the initial data distribution, the dimensionality of the latent space, and the specifics of the encoder, this compression might be lossy---meaning some information is irretrievably lost during encoding and thus cannot be fully restored during decoding.

Thus, the primary goal of an autoencoding method is to identify the optimal encoder/decoder pair within a specified set of candidates. In essence, given a set of potential encoders and decoders, the aim is to find the pair that retains the maximum amount of information during the encoding phase, thereby minimizing the reconstruction error during the decoding phase.  If we denote  $\gE$ and $\gD$ as the families of encoders and decoders we are considering, respectively, the dimensionality reduction problem can be formulated as:
$$
(e, d) = \mathop{\argmin}_{(e,d)\in(\gE,\gD)} f\left(\bx, d(e(\bx))\right),
$$
where $f(\cdot, \cdot)$ measures the loss.
The families of encoders and decoders can include any functions that are suitable for the task, such as multilayer perceptrons or neural networks.
However, when $\gE$ and $\gD$ are linear transformations (a.k.a., a \textit{linear autoencoder}), the autoencoder aligns with the SVD (or PCA) when $f(\cdot, \cdot)$ is the spectral norm loss or the Frobenius norm loss (or others; see Section~\ref{section:svd-low-rank-approxi})~\footnote{The SVD or PCA is actually a special case of the linear autoencoder, in which case the basis vectors are orthonormal. However, the vectors in a linear autoencoder need not necessarily to be orthonormal.}.

Again, we assume the data are centered such that the sample mean vector $\overline{\bx}$ is zero. 
Let the data matrix $\bX \in \real^{n\times p}$ contain the data observations as rows.
As demonstrated in Section~\ref{section:svd-low-rank-approxi}, the truncated SVD provides the optimal low-rank approximation of a matrix.
Let the truncated SVD of $\bX^\top$ be $\widetilde{\bX}^\top=\bU_m\bSigma_m\bV_m^\top$, where $\bU_m\in\real^{p\times m}$, $\bV_m\in\real^{n\times m}$, and $\bSigma_m\in\real^{m\times m}$. Therefore, from the description of PCA, the encoding process can be formulated as 
$$
\text{encoder: }\gap e(\bx_i) = \bU_m^\top \bx_i, \gap \bx_i\in\real^{p\times 1}, \forall i\in\{1,2,\ldots, n\}.
$$
That is, the encoder represents orthonormal eigenvectors corresponding to the $m$ greatest eigenvalues of the covariance  matrix.
This implies  $\bU_m^\top\widetilde{\bX}^\top=\bSigma_m\bV_m^\top$ since $\bU_m^\top\bU_m=\bI_m$. 
Therefore, the decoding process is 
$$
\text{decoder: }\gap d(e(\bx_i)) = \bU_m e(\bx_i), \gap  \forall i\in\{1,2,\ldots, n\}.
$$
This leads to $d(e({\bX}^\top))=\bU_m\bSigma_m\bV_m^\top$, which is the truncated SVD of $\bX^\top$, confirming that when $\gE$ and $\gD$ represent linear transformations, the autoencoder aligns with the SVD (or PCA).

\index{Orthogonal matrix factorization}
\index{Nonnegative PCA}
\subsection*{PCA as Orthogonal Matrix Factorization and Nonnegative PCA}
Consider the same setting of a matrix $\bX$ as described in \eqref{equation:pca-equ0}, and the following matrix factorization problem (see Chapter~\ref{chapter:als} and \ref{chapter:nmf} for more details):
\begin{equation}\label{equation:or_pca1}
\mathopmin{\bW,\bZ} \normf{\bX-\bW\bZ}^2, 
\gap 
\text{with}\gap
\bZ\bZ^\top =\bI_k,
\end{equation} 
where $\bX\in\real^{n\times p}, \bW\in\real^{n\times k}$, and $\bZ\in\real^{k\times p}$ with $k\leq \min\{n,p\}$. 
Given $\bZ$,  the optimal value of $\bW$ is $\bX\bZ^\top$ (see Problem~\ref{prob:ortho_mf}, a consequence of root finding for the gradient).
Therefore, \eqref{equation:or_pca1} can be reformulated as 
\begin{equation}\label{equation:or_pca2}
\mathopmin{\bZ\bZ^\top =\bI_k} \normf{\bX-\bX\bZ^\top\bZ}^2
=
\mathopmin{\bZ\bZ^\top =\bI_k}\normf{\bX}^2 - \normf{\bX\bZ^\top}^2
=
\mathopmax{\bZ\bZ^\top =\bI_k}\normf{\bX\bZ^\top}^2.
\end{equation} 
When $k=1$, this reduces to $\mathopmax{\normtwo{\bz} =1}\bz^\top(\bX^\top\bX)\bz$, which is the same as the PCA problem in \eqref{equation:pca_raw1}.
Therefore,  nonnegative PCA extends this by considering the following problem 
$$
\mathopmax{\bz\in\real_+^p, \normtwo{\bz} =1}\bz^\top(\bX^\top\bX)\bz,
$$
which finds applications in the analysis of gene expression data to identify sets of genes that are co-expressed across different samples and in extracting features from images in a way that preserves the nonnegative structure of pixel values \citep{montanari2015non}.
If $\bX\in\real_+^{n\times p}$, the two problems are equivalent.

\section{Application: Low-Rank Neural Networks}\label{section:low-rank-neural}
\index{Low-rank neural networks}
\subsection*{Low-Rank Neural Networks}
We start with the basic LeNet5 neural network \citep{lecun2015lenet} to illustrate the concept of the low-rank neural networks. 
We have made modifications to the fully connected layers accordingly, which we now refer to as $LenetModified$, as shown in the gray box below. 
Furthermore, we recognize that a layer with 120 input features  and 100 output features  can be represented as a matrix of size $120\times 100$.  
To introduce a form of ``regularization" on the matrix, for instance, we can limit the rank of the matrix to 50. 
This enables us to decompose the $120\times 100$ matrix into two smaller matrices through matrix multiplication.
In essence, whether we have the original $120\times 100$ matrix or  the product of two matrices with dimensions $120\times 50$ and $50\times 100$, there is a functional equivalence in certain contexts. 
The resulting matrix from the multiplication operation also has a size of $120\times 100$. 
Consequently, we obtain a low-rank version of the fully connected layer, termed $LeNetDecom$ structure, where  each fully connected layer in $LenetModified$ is split into two low-rank layers. The fully connected layers of $LenetModified$ and $LenetDecom$ are displayed below:
\begin{svgraybox}
\begin{lstlisting}
LenetModified{
	Convolutional Layers Omitted;
	(Fully Connected Layers): 
	(0): Linear(in_features=120, out_features=100)
	(1): Tanh()
	(2): Linear(in_features=100, out_features=80)
	(3): Tanh()
	(4): Linear(in_features=80, out_features=60)
	(5): Tanh()
	(6): Linear(in_features=60, out_features=40)
	(7): Tanh()
	(8): Linear(in_features=40, out_features=10)
	(9): LogSoftmax()
}
\end{lstlisting}
\end{svgraybox}

\begin{svgraybox}
\begin{lstlisting}
LenetDecom{
	Convolutional Layers Omitted;
	(Fully Connected Layers):
	(0): Linear(in_features=120, out_features=50)
	(1): Linear(in_features=50, out_features=100)
	(2): Tanh()
	(3): Linear(in_features=100, out_features=40)
	(4): Linear(in_features=40, out_features=80)
	(5): Tanh()
	(6): Linear(in_features=80, out_features=30)
	(7): Linear(in_features=30, out_features=60)
	(8): Tanh()
	(9): Linear(in_features=60, out_features=20)
	(10): Linear(in_features=20, out_features=40)
	(11): Tanh()
	(12): Linear(in_features=40, out_features=10)
	(13): LogSoftmax()
}
\end{lstlisting}
\end{svgraybox} 

We realize that by reducing a fully connected layer with size $(m\times n)$ to a low-rank layer $(m\times r) \times (r\times n)$, we can decrease the storage requirement for the model from $mn$ floating-point numbers to $r(m+n)$ floating-point numbers. 
In this specific example mentioned above, we reduce the storage from $120\times 100=12,000$ to $50(120+100)=11,000$ floating-point numbers. 
This reduction also involves the computational complexity from matrix multiplication operations. 

\begin{figure}[h]
\centering  
\vspace{-0.35cm} 
\subfigtopskip=2pt 
\subfigbottomskip=2pt 
\subfigcapskip=-5pt 
\subfigure[Training loss.]{\label{fig:lenetLoss_train}
\includegraphics[width=0.475\linewidth]{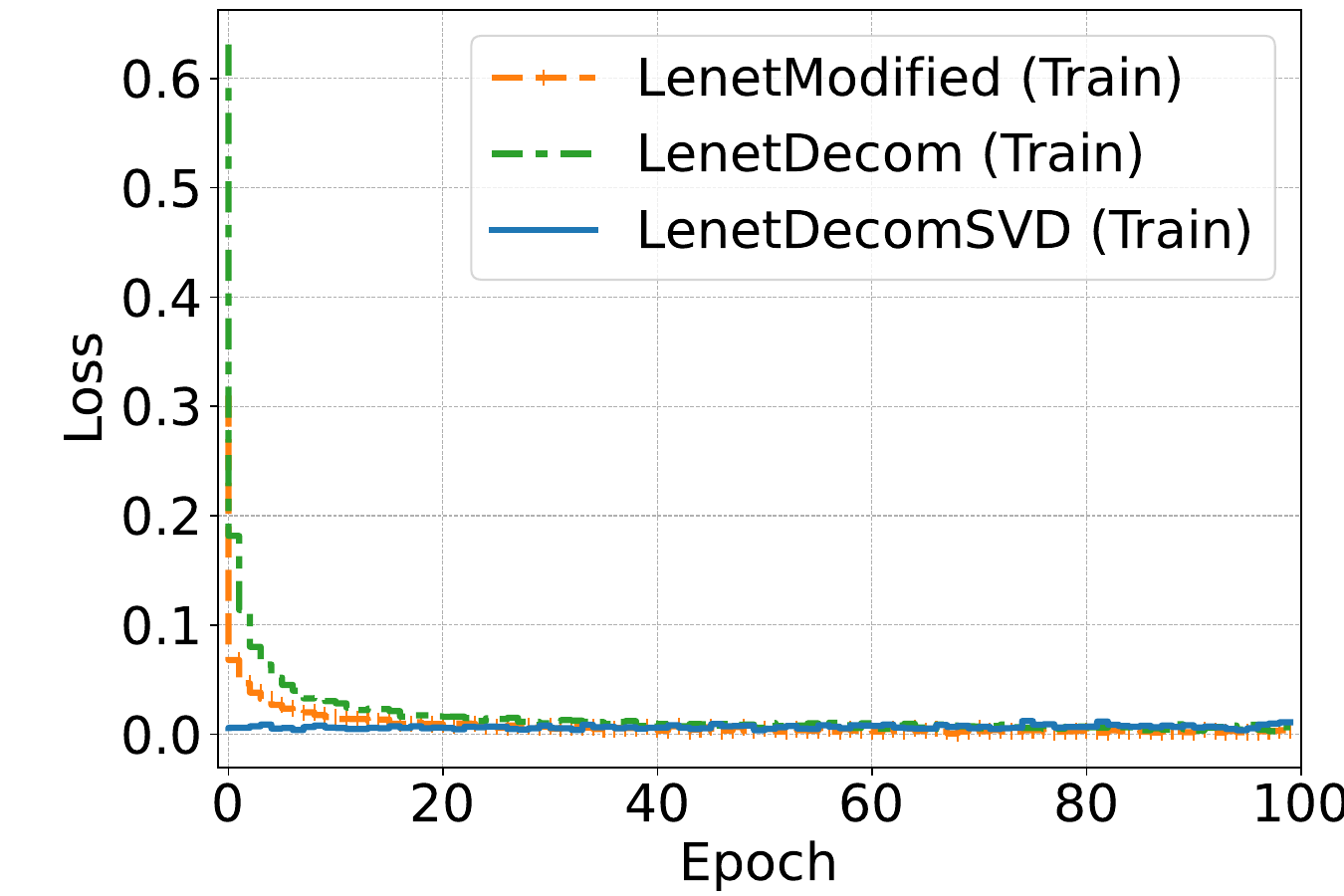}}
\quad 
\subfigure[Validation loss.]{\label{fig:lenetLoss_val}
\includegraphics[width=0.475\linewidth]{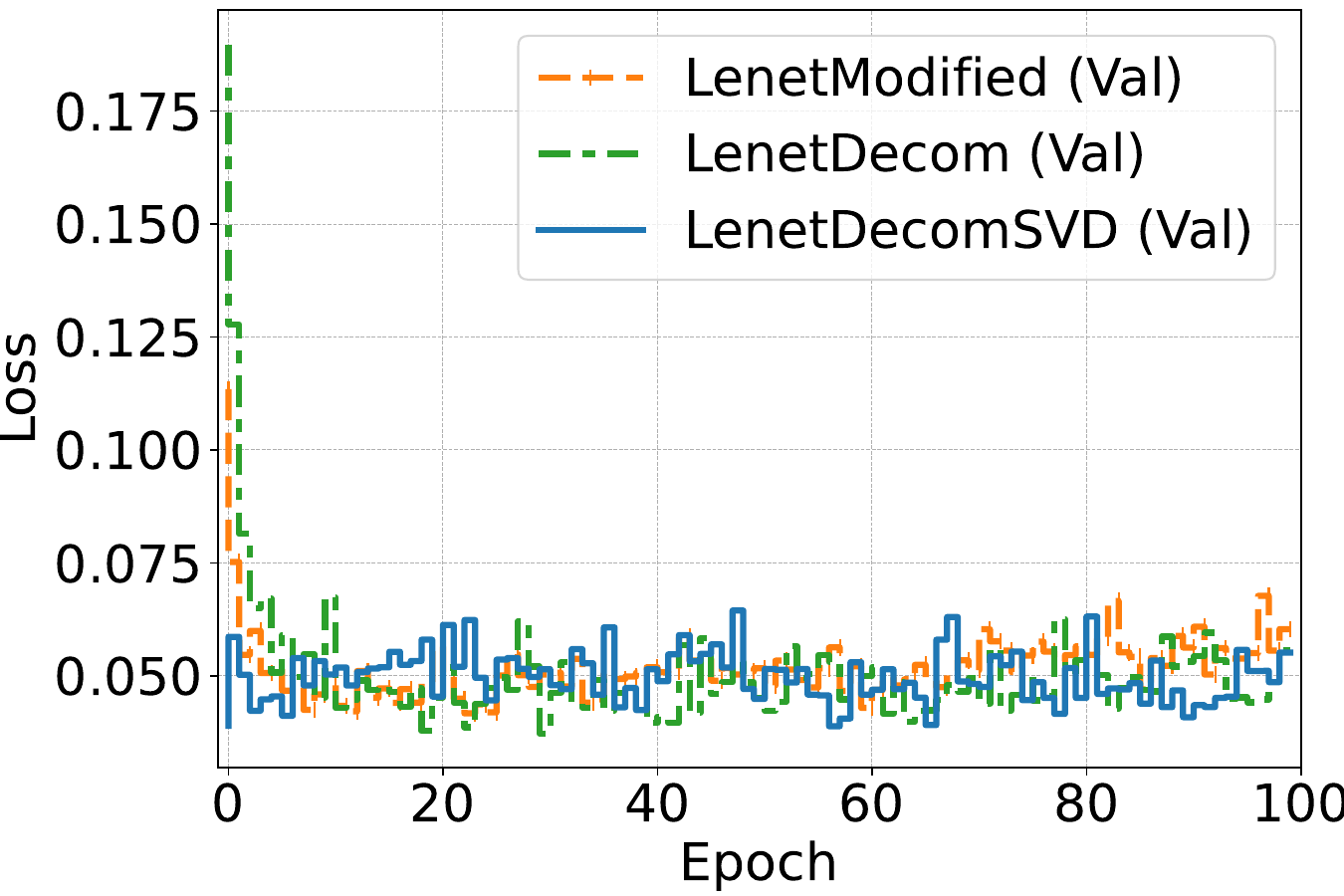}}
\caption{Comparison of full neural networks and low-rank neural networks.}
\label{fig:lenetLoss}
\end{figure}

After training for 100 epochs on the MNIST data set \citep{lecun1998mnist}, we observe that  the minimal training loss of $LenetModified$ is lower than that of $LenetDecom$, as shown in Figure~\ref{fig:lenetLoss_train}. However, the optimal validation loss of $LenetModified$ (greater than 0.04) is higher than that of $LenetDecom$ (less than 0.04), as illustrated in Figure~\ref{fig:lenetLoss_val}. 
This demonstrates the ``regularization" property of low-rank neural networks through this simple example.

Furthermore, suppose we have already trained $LenetModified$, and we wish to utilize the structure of $LenetDecom$ instead, as it saves space and matrix operations. To avoid training from scratch, we could load the weights from $LenetModified$ to $LenetDecom$. 
One approach could be employed is to apply the reduced SVD to decompose the weight matrix and keep only the first $r$ singular values, i.e., zeroing out the trailing  singular values $\sigma_{r+1}= \sigma_{r+2}= \ldots=0$ ($r=50$ in the example mentioned above). That is, a weight matrix $\bW = \bU_r\bSigma_r\bV_r \triangleq \bW_1\bW_2$, where we can set $\bW_1\triangleq\bU_r\bSigma_r$ and $\bW_2\triangleq\bV_r$. 
Then we load the convolutional layers of $LenetDecom$ directly from the convolutional layers of $LenetModified$, and load the fully connected layers of $LenetDecom$ via the decomposition of $LenetModified$. 
Fine-tuning the network will improve performance.  The training and validation losses of this method, denoted as $LenetDecomSVD$, are shown in Figure~\ref{fig:lenetLoss}, where the training loss of $LenetDecomSVD$ approaches 0 even in the first epoch (since we load the weights from a trained $LenetModified$), and the optimal validation loss of $LenetDecomSVD$ is also smaller than 0.04, surpassing the original $LenetModified$. 

This serves as a simple example  of how the low-rank neural network works. 
It is important to note that low-rank neural networks may not necessarily outperform other approaches in all scenarios.  
Furthermore, the matrix decomposition approach can be extended to \textit{convolutional layers} based on the equivalence between convolutional and fully connected layers, as  described in \citet{ma2017equivalence}. 
For a more comprehensive exploration of deep neural networks, refer to the studies conducted by \citet{chen2015net2net, wei2016network, lu2018compnet}

\subsection*{One More Step: Adding a Nonlinear Function Layer}
In the previous section, we approximated a fully connected layer using the product of two low-rank matrices. 
The next step involves introducing a nonlinear function between these two low-rank matrices, resulting in a final layer of the form $f(\bW_1) \bW_2$, where $f(\cdot)$ represents a nonlinear function. 
We call this approach in the above example   $LenetDecomNonlinear$.

\index{Nonlinear function layer}

Furthermore, if we have already trained $LenetModified$, we can modify the structure of previously trained models and transfer to a new structure, thus avoiding the need to train from scratch every time. We can then once again factor the fully connected networks through matrix decomposition.

\begin{SCfigure}
	\centering
	\includegraphics[width=0.45\textwidth]{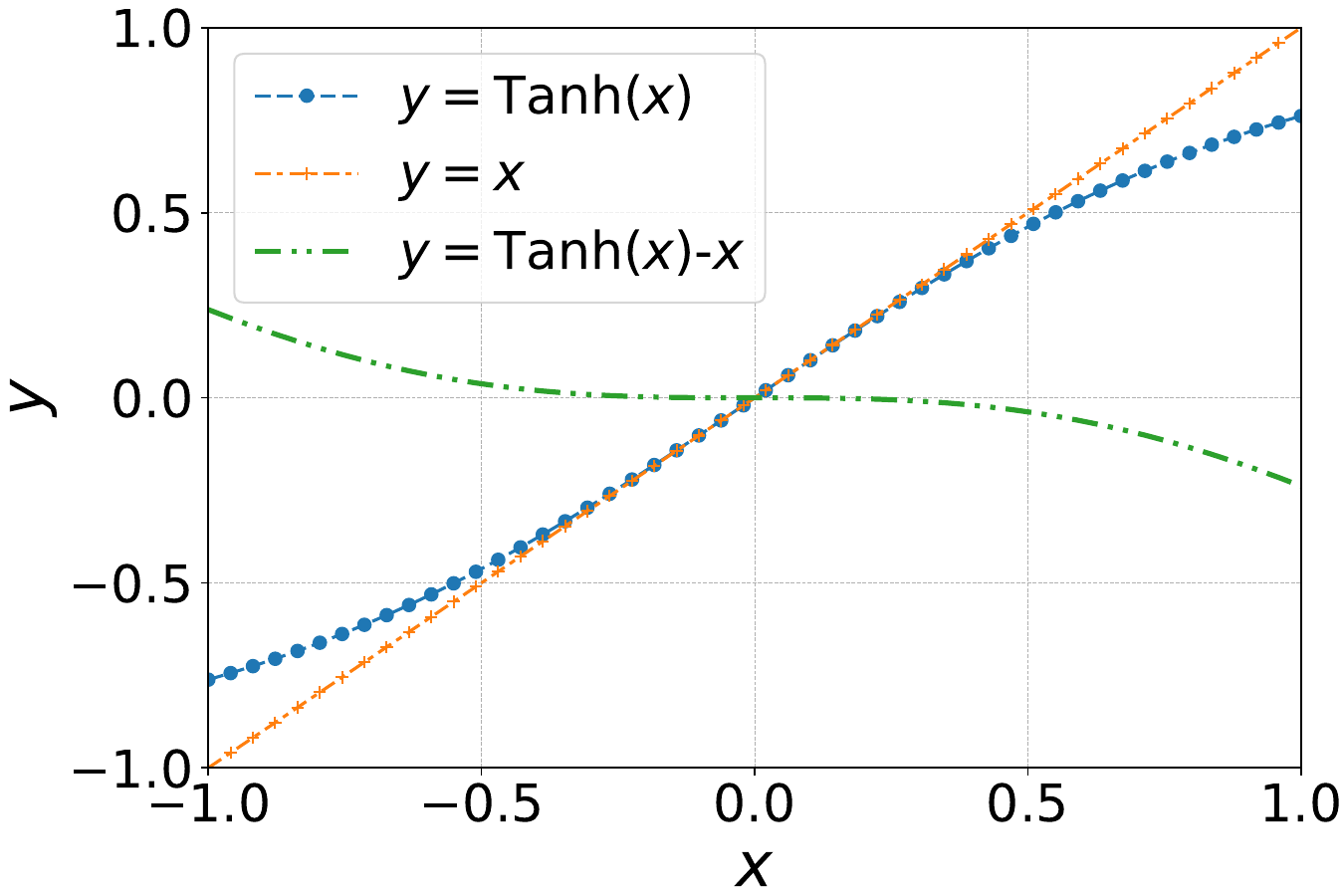}
	\caption{Demonstration of $y=Tanh(x)$ vs $y=x$. We notice that in the input field between $(-0.25, 0.25)$, $Tanh$ is almost close to a linear function.}
	\label{fig:tanh}
\end{SCfigure}

Specifically, when we incorporate the $Tanh$ function in neural networks as the nonlinear function, where $Tanh(x) = \frac{e^x-e^{-x}}{e^x + e^{-x}}$, we observe that  when the input values fall within the range of  $(-0.25, 0.25)$, $Tanh$ behaves almost like a linear function, as shown in Figure~\ref{fig:tanh}, where the (\textcolor{mydarkgreen}{green}) dashed line represents the difference between $y=Tanh(x)$ and $y=x$ and the difference is nearly 0 when $x\in (-0.25, 0.25)$.

After the matrix multiplication of $\bW=\bW_1\bW_2$, if we can ensure that all the values $\bA_{in}\bW_1$ fall within the range of $(-0.25, 0.25)$, where $\bA_{in}$ is the output from the previous layer, then we can add the $Tanh$ function/layer without any suffer, as it operates within the linear space of $Tanh$. 
To address this issue, we observe that matrix decomposition has an equivalence: $\bW = \bW_1*\bW_2 = (\sigma \bW_1)*(\frac{1}{\sigma}\bW_2)$, where $\sigma$ is any nonzero scalar, i.e., \textit{scaling invariance}. 
By adjusting $\sigma$, we can control the maximal absolute value of $\bA_{in}(\sigma \bW_1)$ to be 0.25, thereby adding the nonlinear function without significant impact. 
Following the specific example in the previous section, we denote this method as $LenetDecomNonlinear\_SVD$, where we introduce an extra $Tanh$ function  between the factored matrices, as shown follows:
\begin{svgraybox} 
\begin{lstlisting}
LenetDecomNonlinear{
	Convolutional Layers Omitted;
	(Fully Connected Layers):
	(0): Linear(in_features=120, out_features=50)
	(1): Tanh()  [Differece!]
	(2): Linear(in_features=50, out_features=100)
	(3): Tanh()
	(4): Linear(in_features=100, out_features=40)
	(5): Tanh() [Differece!]
	(6): Linear(in_features=40, out_features=80)
	(7): Tanh()
	(8): Linear(in_features=80, out_features=30)
	(9): Tanh() [Differece!]
	(10): Linear(in_features=30, out_features=60)
	(11): Tanh()
	(12): Linear(in_features=60, out_features=20)
	(13): Tanh() [Differece!]
	(14): Linear(in_features=20, out_features=40)
	(15): Tanh()
	(16): Linear(in_features=40, out_features=10)
	(17): LogSoftmax()
}
\end{lstlisting}
\end{svgraybox} 
\begin{figure}[h]
\centering  
\vspace{-0.35cm} 
\subfigtopskip=2pt 
\subfigbottomskip=2pt 
\subfigcapskip=-5pt 
\subfigure[Training loss.]{\label{fig:nonlinear_lenetLoss_train}
\includegraphics[width=0.475\linewidth]{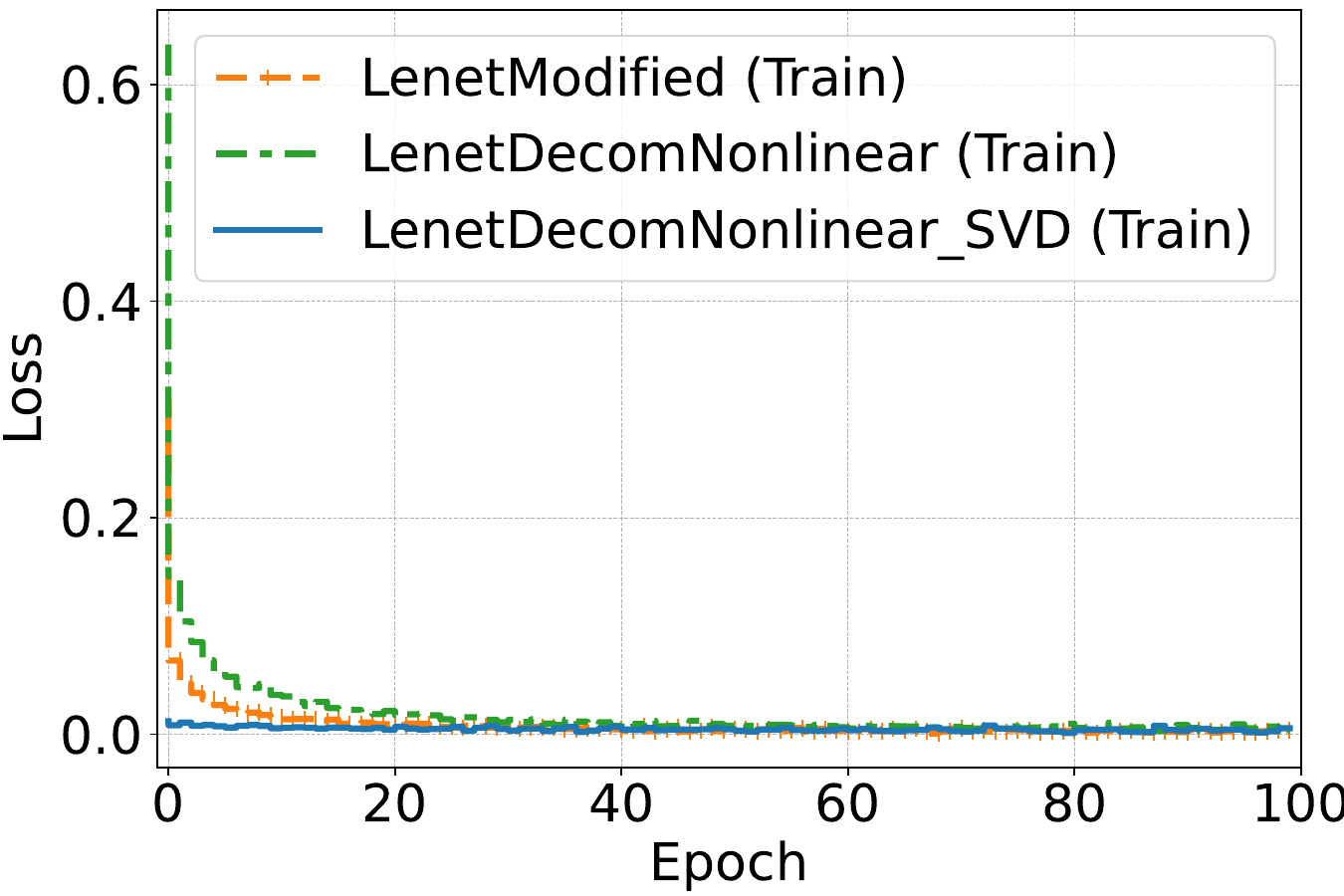}}
\quad 
\subfigure[Validation loss.]{\label{fig:nonlinear_lenetLoss_val}
\includegraphics[width=0.475\linewidth]{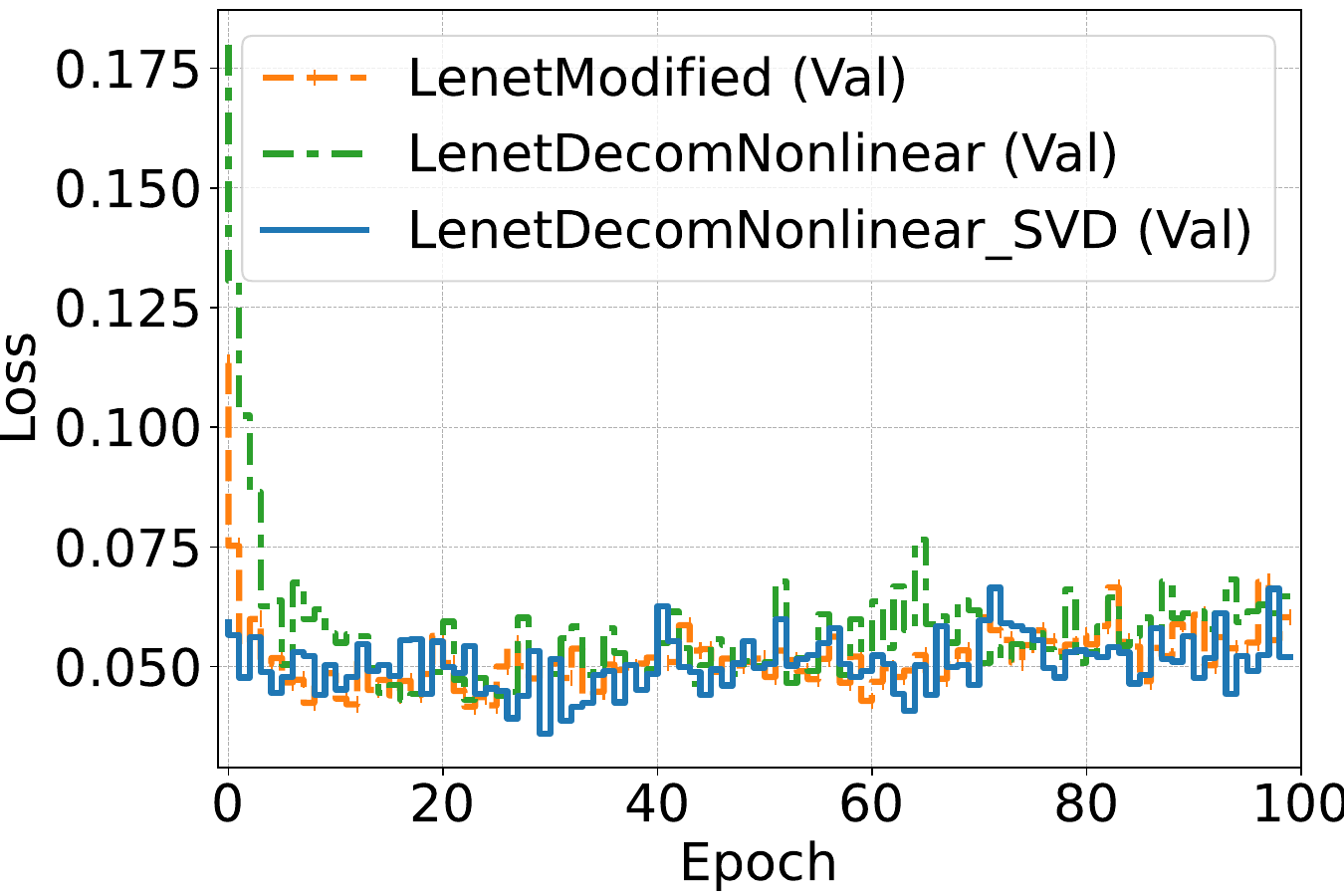}}
\caption{Comparison of full neural networks and $Tanh$ in low-rank neural networks.}
\label{fig:nonlinear_lenetLoss}
\end{figure}

After training for 100 epochs, we find that the minimal training losses for $LenetModified$,\\ $LenetDecomNonlinear$, and $LenetDecomNonlinear\_SVD$ are similar, as shown in Figure~\ref{fig:nonlinear_lenetLoss_train}. 
Additionally, the optimal/minimal validation losses for $LenetModified$ and $LenetDecomNonlinear$ are also similar, both exceeding 0.04. However, the minimal validation loss for  $LenetDecomNonlinear\_SVD$ falls below 0.04, presenting a promising outcome for the application of this matrix decomposition method in neural networks, as illustrated in Figure~\ref{fig:nonlinear_lenetLoss_val}. 
For further insights into these \textit{neural architecture search (NAS)} methods, readers are encouraged to explore the studies by \citet{chen2015net2net, wei2016network, lu2018compnet}.

\index{Symmetric dilation}
\begin{problemset}
	
\item Show that $(\bA\bA^\top)^{1/2}\bA = \bA(\bA^\top\bA)^{1/2}$.
\item Let $\bA\in\real^{m\times n}$ be any matrix. Show that the trace of $\bA^\top\bA$ is equal to the sum of all $a_{ij}^2$, i.e., $\trace(\bA^\top\bA) = \sum_{i,j=1}^{m,n}a_{ij}^2$.
\item Let $\bA\in\real^{m\times n}$ be any matrix with rank $r$. Suppose  the columns of $\bB\in\real^{m\times r}$ span the column space of  $\bA$, and the columns of $\bC\in\real^{n\times r}$ span the row space of  $\bA$. Show that the matrix $\bA$ admits the factorization $\bA=\bB\bE\bC^\top$, where $\bE$ is an $r$ by $r$ nonsingular matrix.
\item Let $\bA\in\real^{n\times n}$  be any square matrix with rank $r$. Consider the $(2n)\times (2n)$ symmetric matrix
$$
\bB=
\begin{bmatrix}
\bzero & \bA \\
\bA^\top & \bzero 
\end{bmatrix}.
$$
Suppose $\bA$ admits the full SVD $\bA=\bU\bSigma\bV^\top$, where $\bSigma=\diag(\sigma_1,\sigma_2,\ldots,\sigma_n)$. 
\begin{itemize}
\item Show that $\sigma_k$ is an eigenvalue of $\bB$ corresponding to the eigenvector $\footnotesize\begin{bmatrix}
	\bv_k \\
	\bu_k
\end{bmatrix}$ for any $k\in\{1,2,\ldots,n\}$,
and that $-\sigma_k$ is an eigenvalue of $\bB$ corresponding to the eigenvector $\footnotesize\begin{bmatrix}
	\bv_k \\
	-\bu_k
\end{bmatrix}$  for any $k\in\{1,2,\ldots,n\}$.
\item Show that the $2n$ eigenvectors are pairwise orthogonal.
\end{itemize}

\item \label{prob:svd_4symm} Let $\bA\in\real^{m\times n}$ be any  rectangular matrix with rank $r$. Consider the $(m+n)\times (m+n)$ symmetric matrix (usually called the \textit{symmetric dilation of $\bA$})
$$
\bB=
\begin{bmatrix}
	\bzero & \bA \\
	\bA^\top & \bzero 
\end{bmatrix}.
$$
Suppose $\bA$ admits the full SVD $\bA=\bU\bSigma\bV^\top$, where $\bSigma=\diag(\sigma_1,\sigma_2,\ldots,\sigma_n)$. 
\begin{itemize}
\item Show that $\sigma_k$ is an eigenvalue of $\bB$ corresponding to the eigenvector $\footnotesize\begin{bmatrix}
	\bv_k \\
	\bu_k
\end{bmatrix}$ for any $k\in\{1,2,\ldots,r\}$,
and that $-\sigma_k$ is an eigenvalue of $\bB$ corresponding to the eigenvector $\footnotesize\begin{bmatrix}
	\bv_k \\
	-\bu_k
\end{bmatrix}$  for any $k\in\{1,2,\ldots,r\}$.
\item Show that the remaining $m+n-2r$ eigenvectors of $\bB$ correspond to the eigenvalue 0.
\item Show that the $m+n$ eigenvectors  are pairwise orthogonal.
\end{itemize}

\item Consider the symmetric dilation of $\bA\in\real^{m\times n}$:
$
\bB=
\footnotesize
\begin{bmatrix}
	\bzero & \bA \\
	\bA^\top & \bzero 
\end{bmatrix}.
$
Show that $\normtwo{\bA}=\normtwo{\bB}$.


\item We have shown in Lemma~\ref{lemma:orthogonal-equivalent-matrix}  that orthogonally equivalent matrices have the same singular values. 
Prove the reverse implication that if two matrices have the same singular values, then they are orthogonally equivalent.

\item Given two nonzero vectors $\bu,\bv\in\real^n$, and let $\bA=\bu\bv^\top$. Show that the nonzero singular value of $\bA$ is $\normtwo{\bu}\cdot \normtwo{\bv}$.

\item Given a square matrix $\bA\in\real^{n\times n}$ with singular values $\sigma_1\geq \sigma_2\geq\ldots \geq \sigma_n$. Show that $\sigma_1^3, \sigma_2^3,\ldots , \sigma_n^3$ are the singular values of $\bA\bA^\top\bA$.

\item Given a rectangular matrix $\bA\in\real^{m\times n}$, and let $\bB\in\real^{\widehat{m}\times \widehat{n}}$ be a submatrix of $\bA$, where $\widehat{m}\leq m$ and $\widehat{n}\leq n$. Show that the largest singular value of $\bB$ is not greater than the largest singular value of $\bA$.

\item Given a  matrix $\bA\in\real^{n\times n}$ and a positive definite matrix $\bB\in\real^{n\times n}$, prove that the singular values of $\bB\bA$ and $\bA$ are the same. Additionally, discuss the relationship between the left and right singular vectors of $\bB\bA$ and $\bA$.

\item We only discuss the SVD of real matrices in this chapter: $\bA=\bU\bSigma\bV^\top$. Show that if $\bA$ is real, then both $\bU$ and $\bV$ are also real.

\item Given a Householder transformation matrix $\bH = \bI - 2\bu\bu^\top\in\real^{n\times n}$, where $\normtwo{\bu}=1$, determine the eigenvalues, determinant, and singular values of this Householder transformation matrix.

\item Given the nonzero singular values $\sigma_1, \sigma_2, \ldots, \sigma_r$ of $\bA$, discuss the singular values of $\bA^\top$, $\gamma\bA$ with $\gamma>0$, and $\bA^{-1}$ (if $\bA$ is nonsingular).

\item Given a square and real matrix $\bA\in\real^{n\times n}$, show that $\bA=\bzero$ if and only if all eigenvalues of $\bA$ are zero.

\item Given a square matrix $\bA\in\real^{n\times n}$, show that $\bA^\top\bA$ and $\bA\bA^\top$ are similar (Definition~\ref{definition:similar-matrices}). \textit{Hint: Proceeding with the SVD of $\bA$}.

\item Consider the definition of a general norm (not necessarily an $\ell_p$ norm) in Section~\ref{section:kyfan_knorm}, show that $\norm{\bx} \leq \sum_{i=1}^{n} \abs{\bx_i}\norm{\be_i}$, where $\bx\in\real^n$ and $\be_i$ is the $i$-th standard basis vector in $\real^n$. \textit{Hint: Examine the properties of norms  on $\bx=\sum_{i=1}^{n}x_i\be_i$.}

\item \citep{higham1994matrix, higham2002accuracy}$^\star$  Given the left polar decomposition of $\bA\in\real^{m\times n}$ with $m\geq n$: $\bA=\bQ\bS$, show that 
$$
\frac{\normtwo{\bA^\top\bA-\bI}}{1+\normtwo{\bA}}
\leq \normtwo{\bA-\bQ}
\leq 
\frac{\normtwo{\bA^\top\bA-\bI}}{1+\sigma_{\min}(\bA)}.
$$ 
This result shows that the two measures of $\normtwo{\bA^\top\bA-\bI}$ and $\normtwo{\bA-\bQ}$ are essentially equivalent.

\item Show that all eigenvalues of a square matrix are less than or equal to its maximal singular value $\sigma_1$.

\item Suppose $\bx$ is an eigenvector of $\bA^\top\bA$ corresponding to a nonzero eigenvalue. Discuss the corresponding eigenvector of $\bA\bA^\top$. \textit{Hint: Premultiply by $\bA$}.

\item Given the SVD of a nonsingular square matrix $\bA=\bU\bSigma\bV^\top\in\real^{n\times n}$, determine the singular values of $\bA^\top\bA$.

\item Find the optimal rank-1 approximation in terms of the spectral norm for the following matrix:
$$\bA=\begin{bmatrix}
\cos\theta & -\sin\theta \\
\sin\theta & \cos\theta
\end{bmatrix}.
$$

\item \textbf{Skew-symmetric.} Given a skew-symmetric and tridiagonal matrix $\bS\in\real^{n\times n}$, show that it admits the decomposition 
$$
\bP^\top\bS\bP 
=
\begin{bmatrix}
\bzero & \bB^\top \\
\bB & \bzero
\end{bmatrix},
$$
where $\bB\in\real^{m\times m}$, $n=2m$, and $\bP$ is a permutation matrix.
Given further the SVD of $\bB=\bU\bSigma\bV^\top$, find the eigenvalues and eigenvectors of $\bS$.

\item Discuss the uniqueness of the polar decomposition for the  matrix
$
\bA = 
\footnotesize\begin{bmatrix}
	1 & 0 \\
	0 & 0\\
\end{bmatrix}.
$

\item Show that the eigenvalues of a PSD matrix are equal to its singular values.

\item Show that the trace of a PSD matrix is equal to the its nuclear norm.

\item Let $\bx\in\real^m$ and $\by\in\real^n$. Show that $\sigma_{\max}(\bx\by^\top) = \sqrt{\bx^\top\bx\by^\top\by}$. \textit{Hint: Consider the definition of the Frobenius norm (Definition~\ref{definition:frobernius-in-svd}).}

\item Let $\bA\in\real^{n\times n}$ with $\rank(\bA)=1$. Show that $\sigma_{max}(\bA) = (\trace(\bA\bA^\top))^{1/2}$.

\item Let $\bA\in\real^{m\times n}$, $\gamma\in\real^n$. Show that $\sigma_i(\gamma\bA)=\abs{\gamma}\sigma_i(\bA)$ for all $i\in\{1,2,\ldots,\min\{m,n\}\}$.

\item \textbf{Orthogonal decomposition \citep{zhan2001span}.} Show that every square matrix $\bA\in\real^{n\times n}$ (with SVD $\bA=\bU\bSigma\bV^\top$ and singular values $\sigma_1\geq \ldots\geq \sigma_n$) is a linear combination  of $n$ orthogonal matrices.
\textit{Hint: Consider $\bA=\sum_{i=1}^{n}a_i \bU\bQ_i\bV$, where $\bQ_i=\diag(-\bI_{i-1}, \bI_{n-i+1})$, and
$$
\footnotesize
\begin{bmatrix}
1	& -1 & -1 & \ldots &  -1\\
1 & 1 & -1 & \ldots  &-1 \\
1 & 1& 1 & \ldots  &-1 \\
\vdots &\vdots & \vdots  & \ddots&\vdots \\
1  & 1 & 1&  \ldots &1 \\
\end{bmatrix}
\begin{bmatrix}
a_1\\a_2\\a_3\\ \vdots\\ a_n
\end{bmatrix}
=
\begin{bmatrix}
	\sigma_1\\\sigma_2\\\sigma_3\\ \vdots\\ \sigma_n
\end{bmatrix}.
$$
}

\item Is there any coordinate transformation involved  in the QR or LQ decomposition?

\item Given the SVD of the matrix 
$$
\bA =
\begin{bmatrix}
	2 & 2 \\
	-1 & 1 
\end{bmatrix}
=
\begin{bmatrix}
	1 & 0 \\
	0 & 1 
\end{bmatrix}
\begin{bmatrix}
	2\sqrt{2} & 0 \\
	0 & \sqrt{2}
\end{bmatrix}
\begin{bmatrix}
	1/\sqrt{2} & 1/\sqrt{2} \\
	-1/\sqrt{2} & 1/\sqrt{2}
\end{bmatrix}
=\bU\bSigma\bV^\top,
$$
illustrate the coordinate transformation of this decomposition in a two-dimensional figure.

\item \citep{horn2012matrix} Let $\bQ\in\real^{n\times n}$ be an orthogonal matrix. Show that $\bQ$ can be decomposed as $\bQ=\bU_1\bU_2\ldots\bU_N\bD$, where $\bD=\diag(1,1,\ldots,1,\det(\bQ))$, each $\bU_i$ represents a plane rotation (Definition~\ref{definition:givens-rotation-in-qr}), and $N=n(n-1)/2$. \textit{Hint: Use the result in Problem~\ref{problem:part_ortho}.}
\end{problemset}

%% file: chapter-eigencalc.tex
\newpage
\chapter{Eigenvalue Problem}\label{section:eigenvalue-problem}
\begingroup
\hypersetup{
	linkcolor=structurecolor,
	linktoc=page,  
}
\minitoc \newpage
\endgroup

\section{Background}
\lettrine{\color{caligraphcolor}T}
The decompositional methods discussed in the previous sections are all related to the eigenvalues and eigenvectors of matrices, making the eigenvalue problem worthy of independent consideration. 
The eigenvalues of a matrix $\bA$ can be found by solving the equation $\det(\bA-\lambda\bI)=0$. However, root-searching in the characteristic equation is usually a very poor computational method for finding eigenvalues.
Following Proposition~\ref{proposition:eigenvalue-similar-matrices}, the ``grand strategy" of virtually all modern eigensystem routines is to nudge
the matrix $\bA$ towards diagonal form by a sequence of similarity transformations such that $\ldots\bP_2\bP_2\bA\bP_1^{-1}\bP_2^{-1}\ldots$ is diagonalized.
There are two distinct sets of techniques for implementing this grand strategy. These techniques complement each other well, and most modern eigensystem routines incorporate both.

The first set of techniques involves constructing individual transformations, known as ``atomic" transformations, which are designed to perform specific tasks. For example, these include:
\begin{itemize}
	\item Jacobi transformations to zero out a particular off-diagonal element (Section~\ref{section:jacobi-spectral}).
	\item Householder transformations to zero out an entire row or column (Chapter~\ref{chapter:hessenberg}).
\end{itemize}
Although a finite sequence of these simple transformations cannot fully diagonalize a matrix, they can reduce it to a special form, such as tridiagonal or Hessenberg (see Chapter~\ref{chapter:tridiagonal_decom} and Chapter~\ref{chapter:hessenberg}). Alternatively, the sequence can be iterated repeatedly until the matrix is nearly diagonal, at which point the second set of techniques can be applied to further refine the result.

The second set of techniques, referred to as factorization methods, is more sophisticated. Suppose a matrix $\bA$ can be factored into a left factor \(\bP_l\) and a right factor \(\bP_r\). Then,
$$
\bA = \bP_l  \bP_r
\quad \implies \quad 
\bP_l^{-1}  \bA = \bP_r.
$$
Multiplying the factors in reverse order yields
$
\bP_r \cdot \bP_l = \bP_l^{-1} \cdot \bA \cdot \bP_l,
$
which represents a similarity transformation on $\bA$ with the transformation matrix \(\bP_l\). In Section~\ref{section:qr_algo} and \ref{section:qr_totrid}, we will examine the QR method, which leverages this idea.
Factorization methods do not converge exactly in a finite number of transformations. However, this set of methods converge rapidly and reliably. When combined with an appropriate initial reduction using simple similarity transformations, these methods are the preferred choice for practical applications.

In this section, we present some classical eigenvalue algorithms that will be often useful for computing eigenvalue-related decompositions. 
To simplify the discussion, we consider only real and symmetric matrices. For a real and symmetric matrix $\bA\in\real^{n\times n}$,  the spectral theorem (Theorem~\ref{theorem:spectral_theorem}) allows us to factor it as
\begin{equation*}
	\bA = \bQ \bLambda \bQ^\top = \bQ \bLambda \bQ^{-1},
\end{equation*}
where the columns of $\bQ = [\bq_1, \bq_2, \ldots, \bq_n]$ are eigenvectors of $\bA$, which are  mutually orthonormal, and the entries of $\bLambda=\diag(\lambda_1, \lambda_2, \ldots, \lambda_n)$ are the corresponding real eigenvalues of $\bA$. We assume further that the eigenvalues are ordered by magnitude such that
$$
|\lambda_1| \geq |\lambda_2| \geq \ldots \geq |\lambda_n| \geq 0,
$$ 
with the understanding that in certain discussions, the equality may not hold due to specific considerations. One thing to note further is that since $\bA$ is diagonalizable, the set of eigenvectors $\{\bq_1, \bq_2, \ldots, \bq_n\}$ spans the entire space $\real^n$, allowing any vector $\bv$ in $\real^n$ to be expressed as a linear combination of the eigenvectors:
$$
\bv = x_1\bq_1 +x_2\bq_2+\ldots+x_n\bq_n, \qquad \text{for all}\gap \bv\in\real^n.
$$

\section{Rate of Convergence}
Before we delve into specific algorithms for computing eigenvalues, it is essential to establish criteria for evaluating the speed at which these algorithms converge, as most of them are iterative methods. We define the convergence of a sequence as follows. Note that the $k$-th element in a sequence is denoted by a superscript in parentheses. 
For instance, $\bA^{(k)}$ denotes the $k$-th matrix in a sequence, and $\ba^{(k)}$ denote the $k$-th vector in a sequence. \index{Rate of convergence}
\begin{definition}[Convergence of a Sequence]
Let $\alpha^{(1)}, \alpha^{(2)}, \ldots \in \real$ be an infinite sequence of scalars. The sequence $\alpha^{(k)}$ is said to converge to $\alpha^\star$ if 
$$
\mathop{\lim}_{k\rightarrow \infty}  \abs{\alpha^{(k)} - \alpha^\star} = 0.
$$
Similarly, let $\balpha^{(1)}, \balpha^{(2)}, \ldots \in \real^n$ be an infinite sequence of vectors. The sequence $\balpha^{(k)}$ is said to converge to $\balpha^\star$ if 
$$
\mathop{\lim}_{k\rightarrow \infty} \norm{\balpha^{(k)} - \balpha^\star} = 0.
$$
\end{definition}
The convergence of a sequence of vectors or matrices depends on the norm used. It's important to note that, according to the equivalence of vector or matrix norms (Theorem~\ref{theorem:equivalence-vector-norm}, Theorem~\ref{theorem:equiv_mat_norm}), if a sequence of vectors or matrices converges in one norm, it will also converge in all other norms. This observation is crucial for the analysis of eigenvector convergence.

\begin{definition}[Linear Convergence]\label{definition:linear-convergence}
A sequence $\alpha^{(k)}$ with limit $\alpha^\star$ is called \textit{linearly convergence} if there exists a constant \textcolor{mylightbluetext}{$c\in (0,1)$} such that 
$$
\abs{\alpha^{(k+1)} - \alpha^\star} \leq c \abs{\alpha^{(k)} - \alpha^\star}.
$$
In other words, the \textit{linearly convergent sequence} has the following property:
$$
\mathop{\lim}_{k\rightarrow \infty}  \frac{\abs{\alpha^{(k+1)} - \alpha^\star}}{\abs{\alpha^{(k)} - \alpha^\star}} = c \in (0,1).
$$
\end{definition}
For example, the sequence $\alpha^{(k)} = 4+ (1/4)^k$ converges linearly to $\alpha^\star = 4$ since 
$$
\mathop{\lim}_{k\rightarrow \infty}  \frac{\abs{\alpha^{(k+1)} - \alpha^\star}}{\abs{\alpha^{(k)} - \alpha^\star}} = \frac{1}{4} \in (0,1).
$$
\begin{definition}[Superlinear Convergence]
	A sequence $\alpha^{(k)}$ with limit $\alpha^\star$ is called \textit{superlinearly convergence} if there exists a constant \textcolor{mylightbluetext}{$c_k >0$ with $c_k \rightarrow 0$}  such that 
	$$
	\abs{\alpha^{(k+1)} - \alpha^\star} \leq c_k \abs{\alpha^{(k)} - \alpha^\star}.
	$$
	In other words, the \textit{superlinearly convergent sequence} has the following property:
	$$
	\mathop{\lim}_{k\rightarrow \infty}  \frac{\abs{\alpha^{(k+1)} - \alpha^\star}}{\abs{\alpha^{(k)} - \alpha^\star}} =0.
	$$
\end{definition}
For example, the sequence $\alpha^{(k)} = 4+\left(\frac{1}{k+4}\right)^{k+3}$ converges superlinearly to $\alpha^\star = 4$ since 
$$
\mathop{\lim}_{k\rightarrow \infty}  \frac{\abs{\alpha^{(k+1)} - \alpha^\star}}{\abs{\alpha^{(k)} - \alpha^\star}}= 
\left(\frac{k+4}{k+5}\right)^{k+3}  \frac{1}{k+5}
= 0.
$$

\begin{definition}[Quadratic Convergence]\label{definition:quadratic-convergence}
	A sequence $\alpha^{(k)}$ with limit $\alpha^\star$ is called \textit{quadratically convergence} if there exists a constant \textcolor{mylightbluetext}{$c>0$} such that 
	$$
	\abs{\alpha^{(k+1)} - \alpha^\star} \leq c \abs{\alpha^{(k)} - \alpha^\star}^2.
	$$
	In other words, the \textit{quadratically convergent sequence} has the following property:
	$$
	\mathop{\lim}_{k\rightarrow \infty}  \frac{\abs{\alpha^{(k+1)} - \alpha^\star}}{\abs{\alpha^{(k)} - \alpha^\star}^2} = c .
	$$
\end{definition}
For example, the sequence $\alpha^{(k)} = 4+ (1/4)^{2^k}$ converges quaeratically to $\alpha^\star = 4$ since 
$$
\mathop{\lim}_{k\rightarrow \infty}  \frac{\abs{\alpha^{(k+1)} - \alpha^\star}}{\abs{\alpha^{(k)} - \alpha^\star}^2} = 1.
$$

\begin{definition}[Cubic Convergence]
	A sequence $\alpha^{(k)}$ with limit $\alpha^\star$ is called \textit{cubically convergence} if there exists a constant \textcolor{mylightbluetext}{$c>0$} such that 
	$$
	\abs{\alpha^{(k+1)} - \alpha^\star} \leq c \abs{\alpha^{(k)} - \alpha^\star}^3.
	$$
	In other words, the \textit{cubically convergent sequence} has the following property:
	$$
	\mathop{\lim}_{k\rightarrow \infty}  \frac{\abs{\alpha^{(k+1)} - \alpha^\star}}{\abs{\alpha^{(k)} - \alpha^\star}^3} = c .
	$$
\end{definition}
For example, the sequence $\alpha^{(k)} = 4+ (1/4)^{3^k}$ converges cubically to $\alpha^\star = 4$ since 
$$
\mathop{\lim}_{k\rightarrow \infty}  \frac{\abs{\alpha^{(k+1)} - \alpha^\star}}{\abs{\alpha^{(k)} - \alpha^\star}^3} = 1.
$$

\section{Rayleigh Quotient}\index{Rayleigh quotient}
The \textit{Rayleigh quotient} of a vector $\bx\in \real^n$ associated with the matrix $\bA$ is the scalar given by the quadratic form:
\begin{equation}\label{equation:raylei_ratio}
r(\bx) = \frac{\bx^\top\bA\bx}{\bx^\top\bx} = \left(\frac{\bx}{\norm{\bx}}\right)^\top \bA \left(\frac{\bx}{\norm{\bx}}\right),
\end{equation}
where $\left(\bx/\norm{\bx}\right)$ is a unit-length vector. When $\bx$ is an eigenvector of $\bA$, $r(\bx)$ corresponds to the eigenvalue $\lambda$ associated with $\bx$. To see this, suppose $\bA\bx=\lambda\bx$, it follows that 
$$\bx^\top\bA\bx = \lambda\bx^\top\bx 
\gap
\implies
\gap 
\lambda = \frac{\bx^\top\bA\bx}{\bx^\top\bx}=r(\bx).
$$ 

\begin{figure}[h]
\centering  
\vspace{-0.35cm} 
\subfigtopskip=2pt 
\subfigbottomskip=2pt 
\subfigcapskip=-5pt 
\subfigure[Surface and contour plots of $r(\bx)$.]{\label{fig:project-rayleigh1}
\includegraphics[width=0.47\linewidth]{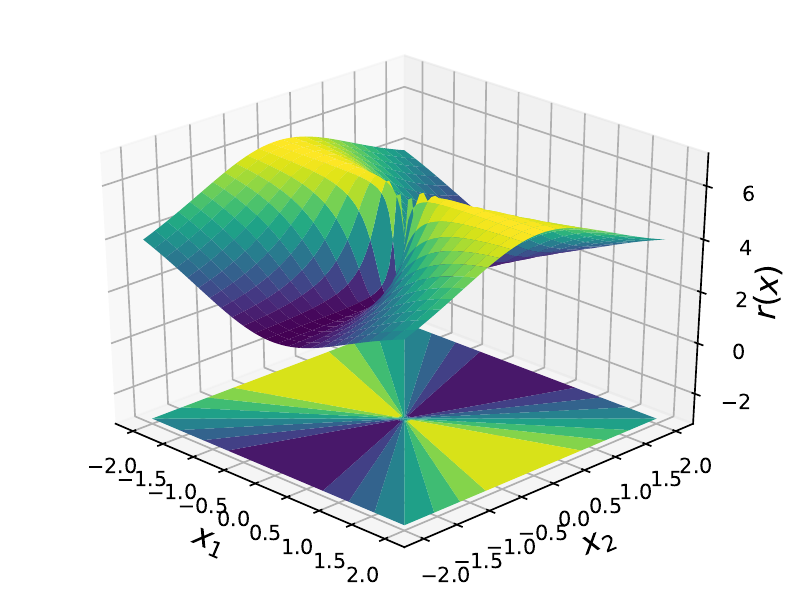}}
\quad 
\subfigure[Surface and contour plots of $\nabla r(\bx)$.]{\label{fig:project-rayleigh2}
\includegraphics[width=0.47\linewidth]{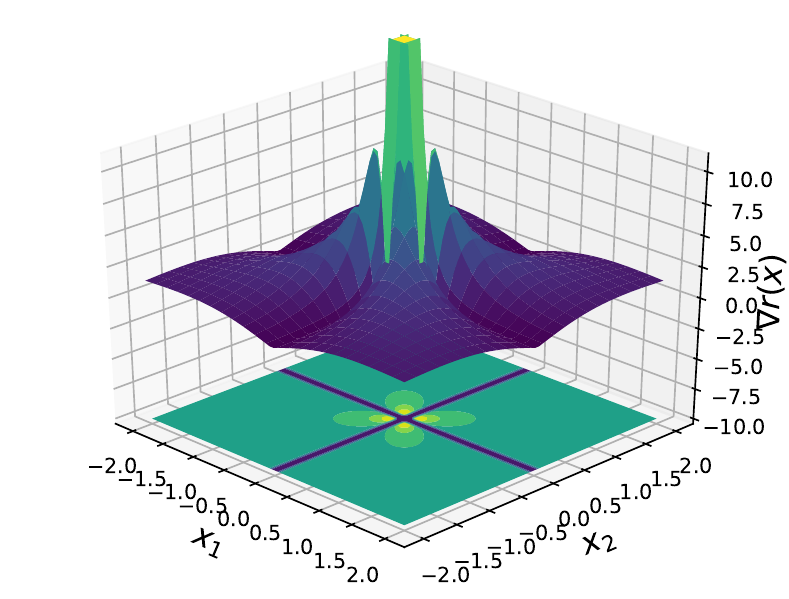}}
\caption{Surface and contour plots (\textcolor{mylightbluetext}{blue}=low, \textcolor{mydarkyellow}{yellow}=high), where in Figure~\ref{fig:project-rayleigh1}, the upper graph is the surface plot, and the lower one is its projection (i.e., contour). The example is drawn by setting $\bA=\scriptsize\begin{bmatrix}
	6 & 0\\ 0 & 2
\end{bmatrix}$, where the input vector $\bx=[x_1,x_2]^\top$ lies in $\real^2$. The eigenvalues are $6$ and $2$, with the corresponding eigenvectors lying on the lines of
$z\cdot \scriptsize\begin{bmatrix}
	1 \\ 0
\end{bmatrix}$ and 
$z\cdot \scriptsize\begin{bmatrix}
	0 \\ 1
\end{bmatrix}$, where $z\in (-\infty, \infty)$ is a scalar.}
\label{fig:rayleigh-projection}
\end{figure}

One of the most significant properties of the Rayleigh quotient is that when $\bx$ is not an eigenvector of $\bA$, the scalar $r(\bx)$ acts most like an eigenvalue in the sense that the squared norm $\norm{\bA\bx-r(\bx)\bx}^2$ is minimized. To see this, suppose we want to find $\lambda$ such that $\norm{\bA\bx-\lambda\bx}^2$ is minimized. Write out the equation with respect to $\lambda$:
$$
\norm{\bA\bx-\lambda\bx}^2 = \bx^\top\bx \lambda^2 - 2\bx^\top\bA\bx \lambda + \bx^\top\bA^\top\bA\bx.
$$
Since $\bx^\top\bx \geq 0$ and assuming $\bx^\top\bx\neq 0$, the  equation above is minimized by setting the derivative to  zero, resulting in $\lambda=r(\bx)$, i.e., the Rayleigh quotient of $\bx$. Therefore, the Rayleigh quotient is a natural  estimate of the eigenvalue   if $\bx$ is close to, but not necessarily equal to, an eigenvector. To see this, it is reasonable to consider the vector $\bx\in \real^n$ as an input variable and $r(\bx)$ as the output. Let $a \triangleq \bx^\top\bA\bx$ and $ b\triangleq\bx^\top\bx$, the gradient of $r(\bx)$ with respect to $\bx$ is given by 
$$
\nabla r(\bx) = \frac{a^\prime b - ab^\prime}{b^2} = \frac{\bx^\top\bx (\bA+\bA^\top)\bx - 2\bx^\top\bA\bx \bx}{\norm{\bx}^4}. \qquad (\text{$\bA$ is any square matrix})
$$
If we further restrict $\bA$ to be symmetric, as we will mostly consider in this section, the gradient of the Rayleigh quotient reduces to 
$$
\nabla r(\bx)  = \frac{2\bx^\top\bx \bA\bx - 2\bx^\top\bA\bx \bx}{\norm{\bx}^4}
=\frac{2}{\norm{\bx}^2}\left(\bA\bx - r(\bx)\bx\right).
 \qquad (\text{$\bA$ is symmetric})
$$
We observe that the gradient is a zero vector if and only if $\bx$ is an eigenvector of $\bA$, in which case $r(\bx)$ is the corresponding eigenvalue. 
This observation has a geometric interpretation: when viewing $\bx$ as the input variable, the \textit{stationary points} of the function $r(\bx)$ correspond to the eigenvectors of $\bA$, and the function outputs are the corresponding eigenvalues. 
An illustrative  example of these stationary points is shown in Figure~\ref{fig:rayleigh-projection}, where  $r(\bx)$ and $\nabla r(\bx)$ are plotted from the matrix $\bA=
\scriptsize
\begin{bmatrix}
6 & 0\\ 0 & 2
\end{bmatrix}$.

\section{Power, Inverse Power, and Rayleigh Quotient Methods}
We start by computing the eigenvalues using several partial methods, which determine the extremal eigenvalues of the symmetric matrix $\bA$, i.e., those with the maximum or minimum magnitude.
\subsection{The Power Method}\label{section:power_method}
The power method  generates a sequence $\bv^{(k)}$ that converges linearly to an eigenvector corresponding to the largest eigenvalue (in magnitude) of $\bA$. The corresponding eigenvalue can then be obtained using the Rayleigh quotient.
Our first attempt to find the eigenvalue of $\bA$ will focus on the value of the largest eigenvalue (in magnitude), denoted by  $\lambda_1$. 
As we progress through this section and  subsequent ones, these ``theoretical" algorithms will shed light on the convergence analysis of the methods.

\begin{algorithm}[H] 
\caption{Power Iteration (Theoretical but Impossible: For Convergence Analysis Only)} 
\label{alg:power-iteration-theoretical}
\begin{algorithmic}[1] 
\Require A real and symmetric $\bA\in \real^{n\times n}$;
\State Initialize $\bv^{(0)}\leftarrow$ some vector with $\norm{\bv^{(0)}}=1$;
\For{$k=1,2,\ldots$} 
\State $\bw \leftarrow \bA \bv^{(k-1)}$;
\State $\bv^{(k)} \leftarrow \bw / \textcolor{mylightbluetext}{\lambda_1}$;
\State $\lambda^{(k)} \leftarrow \frac{(\bv^{(k)})^\top\bA\bv^{(k)}}{(\bv^{(k)})^\top\bv^{(k)}}$; \Comment{i.e., Rayleigh quotient}
\EndFor
\end{algorithmic} 
\end{algorithm}
Suppose the eigenvalue $\lambda_1$ with the largest magnitude is given in advance. Algorithm~\ref{alg:power-iteration-theoretical} provides an iterative approach to finding the corresponding eigenvector. 
Express $\bv^{(0)}$ as a linear combination of the orthonormal eigenvectors $\{\bq_1, \bq_2, \ldots, \bq_n\}$ (since they span the entire space $\real^n$, see Spectral Theorem~\ref{theorem:spectral_theorem}):
\begin{equation}\label{equation:power-inprac-eq1}
\begin{aligned}
\bv^{(0)} = x_1 \bq_1 + x_2\bq_2 + \ldots + x_n\bq_n = \bQ\bx
\quad\implies\quad \bx = \bQ^{-1}\bv^{(0)},
\end{aligned}
\end{equation}
where $\bx = [x_1, x_2, \ldots, x_n]^\top$.
Since the matrix $\bA\in \real^{n\times n}$ is real and symmetric, the $k$-th power of $\bA$ is given by 
\begin{equation}\label{equation:power-inprac-eq2}
\bA^k = \bQ \bLambda^k \bQ^{-1} 
\quad\implies\quad 
	\bA^k \bQ = \bQ\bLambda^k 
	\quad\implies\quad 
\bA^k \bq_i = \lambda_i^k \bq_i,\forall\, i
\end{equation}
where $\bA = \bQ \bLambda \bQ^{-1}$ represents the spectral decomposition of $\bA$ (Theorem~\ref{theorem:spectral_theorem}; Remark~\ref{remark:power-spectral}).
Clearly, we can obtain the $k$-th element $\bv^{(k)}$ in the sequence from Algorithm~\ref{alg:power-iteration-theoretical} by 
$$
\begin{aligned}
\bv^{(k)} &= \bA\bv^{(k-1)} /\lambda_1 = \bA^k\bv^{(0)} / \lambda_1^k
=1/ \lambda_1^k(x_1 \lambda_1^k\bq_1 + x_2\lambda_2^k\bq_2 + \ldots + x_n\lambda_n^k\bq_n)\\
&=x_1 \bq_1 + x_2\left(\frac{\lambda_2}{\lambda_1}\right)^k\bq_2 + \ldots + x_n\left(\frac{\lambda_n}{\lambda_1}\right)^k\bq_n.
\end{aligned}
$$ 
Then,  it follows that
$$
\begin{aligned}
\bv^{(k)} - x_1\bq_1 &= \bA^k\bv^{(0)} / \lambda_1^k- x_1\bq_1= (\bQ \bLambda^k \bQ^{-1}) \bv^{(0)}/\lambda_1^k - x_1\bq_1 \\
\stackrel{\bx = \bQ^{-1}\bv^{(0)}}{\longeqf} &\bQ 
\footnotesize
\begin{bmatrix}
1 & & & \\
& \left(\frac{\lambda_2}{\lambda_1}\right)^k & & \\
& & \ddots & \\
& & & \left(\frac{\lambda_n}{\lambda_1}\right)^k
\end{bmatrix}
\normalsize
\bx- x_1\bq_1  
= \bQ 
\footnotesize
\begin{bmatrix}
0 & & & \\
& \left(\frac{\lambda_2}{\lambda_1}\right)^k & & \\
& & \ddots & \\
& & & \left(\frac{\lambda_n}{\lambda_1}\right)^k
\end{bmatrix}
\normalsize
\bx.
\end{aligned}
$$
Therefore, we have the $k$-th and $(k+1)$-th elements of the  sequence as follow: 
$$
\bQ^{-1}(\bv^{(k)} - x_1\bq_1) = 
\footnotesize
\begin{bmatrix}
	0 & & & \\
	& \left(\frac{\lambda_2}{\lambda_1}\right)^k & & \\
	& & \ddots & \\
	& & & \left(\frac{\lambda_n}{\lambda_1}\right)^k
\end{bmatrix}
\normalsize
\bx,
$$
and 
$$
\bQ^{-1}(\bv^{(k+1)} - x_1\bq_1) = 
\footnotesize
\begin{bmatrix}
	0 & & & \\
	& \left(\frac{\lambda_2}{\lambda_1}\right)^{k+1} & & \\
	& & \ddots & \\
	& & & \left(\frac{\lambda_n}{\lambda_1}\right)^{k+1} 
\end{bmatrix}
\normalsize
\bx
=
\footnotesize
\begin{bmatrix}
	0 & & & \\
	& \left(\frac{\lambda_2}{\lambda_1}\right) & & \\
	& & \ddots & \\
	& & & \left(\frac{\lambda_n}{\lambda_1}\right)
\end{bmatrix}
\normalsize
\bQ^{-1}(\bv^{(k)} - x_1\bq_1).
$$
Define a matrix norm $\norm{\bb}_{\bQ^{-1}}$ as $\norm{\bQ^{-1} \bb}_2$ for any vector $\bb\in\real^n$.
This results in 
\begin{equation}\label{equation:power-iteration-analysis}
\begin{aligned}
\norm{\bv^{(k+1)} - x_1\bq_1}_{\bQ^{-1}} &= 
\norm{\bQ^{-1}(\bv^{(k+1)} - x_1\bq_1)}_2\\
&\leq \left|\frac{\lambda_2}{\lambda_1}\right| \cdot 
\norm{\bQ^{-1}(\bv^{(k)} - x_1\bq_1)}_2 
= \left|\frac{\lambda_2}{\lambda_1}\right| \cdot 
\norm{\bv^{(k)} - x_1\bq_1}_{\bQ^{-1}},
\end{aligned}
\end{equation}
where the inequality is derived from the matrix-vector inequality (Definition~\ref{definition:spectral_norm}).
The above deduction shows that if we can establish that $\norm{\bv}_{\bX^{-1}} = \norm{\bX^{-1}\bv}_2$ is a valid vector norm satisfying the three criteria of a vector norm (Definition~\ref{definition:matrix-norm}), then, according to
the equivalence of vector norms (Theorem~\ref{theorem:equivalence-vector-norm}), we can shown that  $\bv^{(k)}$ in Algorithm~\ref{alg:power-iteration-theoretical} converges linearly to $x_1\bq_1$  (Definition~\ref{definition:linear-convergence}) since $\left|\frac{\lambda_2}{\lambda_1}\right|\leq 1$. This is indeed  the case.

\begin{exercise}[The Q-Norm]
Show that the vector norm $\norm{\bb}_{\bQ^{-1}}$ for any vector $\bb$ with an orthogonal $\bQ$ is a valid vector norm that satisfies the nonnegativity, positive homogeneity, and the triangle inequality properties (Section~\ref{section:kyfan_knorm} or Definition~\ref{definition:matrix-norm}).
The Q-norm is indeed equivalent to the $\ell_2$ vector norm due to the orthogonal invariance of the $\ell_2$ norm (Exercise~\ref{exercise:orthogo_ell2}).
\end{exercise}

Nonetheless, a limitation of Algorithm~\ref{alg:power-iteration-theoretical} lies in our lack of knowledge about $\lambda_1$, rendering it impractical for computing the eigenvector.
To address this issue, we scale it to have unit length at each iteration. The modified algorithm is  presented in Algorithm~\ref{alg:power-iteration}, 
with changes highlighted in \textcolor{mylightbluetext}{blue} text.

\begin{algorithm}[H] 
\caption{Power Iteration (A Practical One, Compare to Algorithm~\ref{alg:power-iteration-theoretical})} 
\label{alg:power-iteration}
\begin{algorithmic}[1] 
\Require A real and symmetric $\bA\in \real^{n\times n}$;
\State Initialize $\bv^{(0)}\leftarrow$ some vector with $\norm{\bv^{(0)}}=1$;
\For{$k=1,2,\ldots$} 
\State $\bw \leftarrow \bA \bv^{(k-1)}$;
\State $\bv^{(k)} \leftarrow \bw / \textcolor{mylightbluetext}{\norm{\bw}}$;
\State $\lambda^{(k)} \leftarrow (\bv^{(k)})^\top\bA\bv^{(k)}$;  \Comment{i.e., Rayleigh quotient}
\EndFor
\end{algorithmic} 
\end{algorithm}
\paragraph{Convergence analysis.}
Clearly, $\bv^{(k)}$ is still a scalar  multiple of $\bA^k \bv^{(0)}$, such that  $\bv^{(k)} = c_k \bA^k \bv^{(0)} = \frac{\bA^k \bv^{(0)}}{\norm{\bA^k \bv^{(0)}}}$ for some constant $c_k$. We have 
\begin{equation}\label{equation:converge-equa-of-power-ite}
\begin{aligned}
\bv^{(k)} &=  c_k \bA^k\bv^{(0)}=c_k (x_1 \lambda_1^k\bq_1 + x_2\lambda_2^k\bq_2 + \ldots + x_n\lambda_n^k\bq_n) \\
&= c_k \lambda_1^k\bigg(x_1 \bq_1 + \underbrace{x_2\left(\frac{\lambda_2}{\lambda_1}\right)^k\bq_2 + \ldots + x_n\left(\frac{\lambda_n}{\lambda_1}\right)^k\bq_n}_{\triangleq\by^{(k)}}\bigg).
\end{aligned}
\end{equation}
The sequence $\by^{(k)}$ in the above equation vanishes as $k\rightarrow \infty$ if $|\lambda_1| \textcolor{mylightbluetext}{>} |\lambda_2| \geq |\lambda_3| \geq \ldots \geq |\lambda_n| \geq 0$. Furthermore, if $x_1\neq 0$, then as $k\rightarrow  \infty$, the vector sequence $\bv^{(k)}$ converges to $\pm \bq_1$, where the sign is decided by $\lambda_1^k x_1$. 

Given that $\bv^{(k)}$ converges to $\pm \bq_1$ with $c_k\lambda_1^k x_1 \rightarrow \pm 1$, following the discussion in \citet{quarteroni2010numerical}, we  provide an improved result on the convergence of the eigenvector $\bv^{(k)}$ with a lower bound such that a form in the linear convergence can be obtained (Definition~\ref{definition:linear-convergence}).
This observation proceeds the following theorem.
\begin{theorem}[Convergence of Power Iteration: Eigenvector]\label{theorem:power-convergence}
Let $\bA\in \real^{n\times n}$ be real and symmetric. 
Moreover, assume the following conditions hold:
\begin{itemize}
\item  $|\lambda_1| \textcolor{mylightbluetext}{>} |\lambda_2| \geq |\lambda_3| \geq \ldots \geq |\lambda_n| \geq 0$, i.e., $\lambda_1$ has an algebraic multiplicity  of 1; 
\item  $\bq_1^\top \bv^{(0)} \neq 0$, i.e., the initial guess $\bv^{(0)}$ has a component in the direction of the eigenvector $\bq_1$ associated with the eigenvalue $\lambda_1$.
\end{itemize}
\item
Then, for Algorithm~\ref{alg:power-iteration},  there exists a constant $c$ such that 
$$
\norm{\widetildebv^{(k)}-\bq_1} \leq c\cdot \left|\frac{\lambda_2}{\lambda_1}\right|^k,
$$
with the sequence $\widetildebv^{(k)} \triangleq \frac{\bv^{(k)}}{c_k\lambda_1^k x_1} = \bq_1+ \sum_{i=2}^n \frac{x_i}{x_1}\left(\frac{\lambda_i}{\lambda_1}\right)^k \bq_i $.
\end{theorem} 
\begin{proof}[of Theorem~\ref{theorem:power-convergence}]
We notice that $\norm{\widetildebv^{(k)}-\bq_1} = 
\norm{\sum_{i=2}^n \frac{x_i}{x_1}\left(\frac{\lambda_i}{\lambda_1}\right)^k \bq_i}$. 
Let $\bz\triangleq\left[0, \frac{x_2}{x_1}\left(\frac{\lambda_2}{\lambda_1}\right)^k, \frac{x_3}{x_1}\left(\frac{\lambda_3}{\lambda_1}\right)^k, \ldots, \frac{x_n}{x_1}\left(\frac{\lambda_n}{\lambda_1}\right)^k\right]$, $\bQ=[\bq_1,\bq_2, \bq_3, \ldots, \bq_n]$. It follows that 
$$
\begin{aligned}
\norm{\widetildebv^{(k)}-\bq_1} &=\norm{\bQ\bz}  \stackrel{\star}{=}  \norm{\bz} 
=\left(\sum_{i=2}^{n}\left(\frac{x_i}{x_1}\right)^2\left(\frac{\lambda_i}{\lambda_1}\right)^{2k} \right)^{1/2}
\leq \left|\frac{\lambda_2}{\lambda_1}\right|^k
\left(\sum_{i=2}^{n}\left(\frac{x_i}{x_1}\right)^2 \right)^{1/2},
\end{aligned}
$$
where  the equality ($\star$) is derived from the invariance under orthogonal transformations (Exercise~\ref{exercise:orthogo_ell2}) and the inequality comes from the assumption $|\lambda_2| \geq |\lambda_3| \geq \ldots \geq |\lambda_n|$. Let $c\triangleq\bigg(\sum_{i=2}^{n}\left(\frac{x_i}{x_1}\right)^2 \bigg)^{1/2}$, the result follows.


The above equation also indicates  (based on  the assumption $|\lambda_2| \geq |\lambda_3| \geq \ldots \geq |\lambda_n|$):
$$
\begin{aligned}
\norm{\widetildebv^{(k)}-\bq_1} &= \norm{\bQ\bz}  \stackrel{\star}{=}  \norm{\bz} 
=\left(\sum_{i=2}^{n}\left(\frac{x_i}{x_1}\right)^2\left(\frac{\lambda_i}{\lambda_1}\right)^{2k} \right)^{1/2}\geq  \left|\frac{\lambda_n}{\lambda_1}\right|^k
\left(\sum_{i=2}^{n}\left(\frac{x_i}{x_1}\right)^2 \right)^{1/2}.
\end{aligned}
$$
Therefore, a bound on the convergence is given by
$$
\left.
\begin{aligned}
\norm{\widetildebv^{(k+1)}-\bq_1}
&\leq c\cdot \left|\frac{\lambda_2}{\lambda_1}\right|^{k+1}\\
\norm{\widetildebv^{(k)}-\bq_1} 
&\geq c\cdot \left|\frac{\lambda_n}{\lambda_1}\right|^k 
\end{aligned}
\right\}
\quad\implies\quad  
\frac{\norm{\widetildebv^{(k+1)}-\bq_1}}{\norm{\widetildebv^{(k)}-\bq_1}}
\leq 
\left| 
\frac{\lambda_2^{k+1}}{\lambda_1 \cdot \lambda_n^k}
\right|.
$$
However, since the values of  $\lambda_n, \lambda_2,$ and $\lambda_1$
are unknown, the above bound may be uninformative. It becomes tight only when $|\lambda_1| \gg |\lambda_2|$.
\end{proof}

Going further and following from the discussion in \citet{golub2013matrix}, we  provide a result on the convergence of the eigenvalue $\lambda^{(k)}$ with a  bound such that the sequence $\lambda^{(k)}$ converges to $\lambda_1$ \textbf{quadratically} with respect to the ratio $\left|\frac{\lambda_2}{\lambda_1}\right|$.
\begin{theorem}[Convergence of Power Iteration: Eigenvalue]\label{Theorem:power-convergence-eigenvalue}
(Under the same conditions as Theorem~\ref{theorem:power-convergence}) Let $\bA\in \real^{n\times n}$ be real and symmetric. Moreover, assume the following conditions hold:
\begin{itemize}
\item  $|\lambda_1| \textcolor{mylightbluetext}{>} |\lambda_2| \geq |\lambda_3| \geq \ldots \geq |\lambda_n| \geq 0$, i.e., $\lambda_1$ has an algebraic multiplicity of 1; 
\item  $\bq_1^\top \bv^{(0)} \neq 0$, i.e., the initial guess $\bv^{(0)}$ has a component in the direction of the eigenvector $\bq_1$ associated with the eigenvalue $\lambda_1$.
\end{itemize}
\item
Then, for Algorithm~\ref{alg:power-iteration}, define $\theta_k \in [0, \pi/2]$ as 
$$
c_k = \cos\theta_k = |\bq_1^\top \bv^{(k)}|.
$$
The quantity $c_k$ is well defined since $\norm{\bq_1}=\norm{\bv^{(k)}}=1 \rightarrow 0< |\bq_1^\top \bv^{(k)}|\leq1$.
Then the sequences $s_k=\sin \theta_k$ and $t_k = \tan \theta_k$ follow
\begin{itemize}
\item Convergence of $s_k$: $ |s_k| \leq \textcolor{black}{t_0} \left|\frac{\lambda_2}{\lambda_1}\right|^k$;
\item Convergence of $\lambda^{(k)}$: $\begin{aligned}
	|\lambda^{(k)} - \lambda_1|  \leq  \mathop{\max}_{2\leq i \leq n} |\lambda_1-\lambda_i| \cdot  \textcolor{black}{t_0^2} \left( \frac{\lambda_i}{\lambda_1}\right)^{2k}.
\end{aligned}
$
\end{itemize}
\end{theorem} 
\begin{proof}[of Theorem~\ref{Theorem:power-convergence-eigenvalue}]
Since $s_k^2 = 1-c_k^2 = 1- (\bq_1^\top \bv^{(k)})^2 = 1-  \left(\frac{\bq_1^\top \bA^k \bv^{(0)}}{\norm{\bA^k \bv^{(0)}}} \right)^2$ by Equation~\eqref{equation:converge-equa-of-power-ite}, where $\bv^{(0)}$ can be expressed as the linear combination of the orthonormal eigenvectors:
$$
\bv^{(0)} = x_1 \bq_1 + x_2\bq_2 + \ldots + x_n\bq_n = \bQ\bx  \quad\implies\quad  \bx = \bQ^{-1}\bv^{(0)}.
$$
Since $\norm{\bx}^2 = \norm{\bQ^{-1}\bv^{(0)}}^2$ and $\bv^{(0)}$ is of length 1, this indicates $x_1^2+x_2^2+\ldots+x_n^2=\norm{\bx}^2=1$. 
Furthermore, as demonstrated  in Equation~\eqref{equation:converge-equa-of-power-ite} that
$$
\bA^k \bv^{(0)} = x_1 \lambda_1^k\bq_1 + x_2\lambda_2^k\bq_2 + \ldots + x_n\lambda_n^k\bq_n.
$$
These findings imply that 
\begin{equation}\label{equation:power-lambda-1}
\begin{aligned}
s_k^2 &=1-  \left(\frac{\bq_1^\top \bA^k \bv^{(0)}}{\norm{\bA^k \bv^{(0)}}} \right)^2 
=1-  \left(\frac{x_1^2 \lambda_1^{2k}}{\sum_{i=1}^{n} x_i^2 \lambda_i^{2k}} \right)
= \frac{\sum_{i=2}^{n} x_i^2 \lambda_i^{2k}}{\sum_{i=1}^{n} x_i^2 \lambda_i^{2k}} 
\leq \textcolor{black}{\frac{\sum_{i=2}^{n} x_i^2 \lambda_i^{2k}}{x_1^2 \lambda_1^{2k}}} \\
&=\frac{1}{x_1^2} \sum_{i=2}^{n} x_i^2 \left( \frac{\lambda_i}{\lambda_1}\right)^{2k} 
\leq \frac{1}{x_1^2} \left( \sum_{i=2}^{n} x_i^2 \right) \left( \frac{\lambda_2}{\lambda_1}\right)^{2k}
= \frac{1-x_1^2}{x_1^2}\left( \frac{\lambda_2}{\lambda_1}\right)^{2k} = t_0^2 \left( \frac{\lambda_2}{\lambda_1}\right)^{2k},
\end{aligned}
\end{equation}
where the last equality follows from the definition of $c_k$; when $k=0$, we have $c_0 = |\bq_1^\top\bv^{(0)}|=|x_1| \neq 0$.
Therefore, the first result follows:
$
\abs{s_k}\leq t_0 \left|\frac{\lambda_2}{\lambda_1}\right|^k.
$
The above result  aligns with the work of \citet{golub2013matrix}. 
Furthermore, since $\bv^{(k)}=\frac{\bA^k \bv^{(0)}}{\norm{\bA^k \bv^{(0)}}}$ as per Equation~\eqref{equation:converge-equa-of-power-ite}, we have 
$$
\begin{aligned}
\lambda^{(k)} &= \bv^{(k)\top}\bA\bv^{(k)} = 
\bigg(\frac{\bA^k \bv^{(0)}}{\norm{\bA^k \bv^{(0)}}}\bigg)^\top 
\bA 
\bigg(\frac{\bA^k \bv^{(0)}}{\norm{\bA^k \bv^{(0)}}}\bigg)
=
\frac{\bv^{(0)\top} \bA^{2k+1}\bv^{(0)}}
{\bv^{(0)\top} \bA^{2k}\bv^{(0)}}
=
\frac{\sum_{i=1}^{n}x_i^2 \lambda_i^{2k+1}}{\sum_{i=1}^{n}x_i^2 \lambda_i^{2k}}.
\end{aligned}
$$
Therefore, 
$$
\begin{aligned}
|\lambda^{(k)} - \lambda_1| 
&= 
\left|
\frac{\sum_{i=\textcolor{mylightbluetext}{2}}^{n}x_i^2 \lambda_i^{2k}(\lambda_i-\lambda_1)}
{\sum_{i=1}^{n}x_i^2 \lambda_i^{2k}}
\right|
\leq \mathop{\max}_{2\leq i \leq n} |\lambda_1-\lambda_i|
\left(
\frac{\sum_{i=\textcolor{mylightbluetext}{2}}^{n}x_i^2 \lambda_i^{2k}}
{\sum_{i=1}^{n}x_i^2 \lambda_i^{2k}}
\right)\\
&\leq \mathop{\max}_{2\leq i \leq n} |\lambda_1-\lambda_i| \cdot  t_0^2 \left( \frac{\lambda_i}{\lambda_1}\right)^{2k},
\end{aligned}
$$
where the last inequality comes from Equation~\eqref{equation:power-lambda-1} shown above.
\end{proof}

From the above deduction, we conclude the following convergence result:
\begin{theorem}[Convergence of Power Iteration]\label{theorem:convergence-power-iteration}
(Under the same conditions as Theorem~\ref{theorem:power-convergence}) Let $\bA\in \real^{n\times n}$ be real and symmetric. Moreover, assume the following conditions hold:
\begin{itemize}
\item  $|\lambda_1| \textcolor{mylightbluetext}{>} |\lambda_2| \geq |\lambda_3| \geq \ldots \geq |\lambda_n| \geq 0$, i.e., $\lambda_1$ has an algebraic multiplicity  of 1; 
\item  $\bq_1^\top \bv^{(0)} \neq 0$, i.e., the initial guess $\bv^{(0)}$ has a component in the direction of the eigenvector $\bq_1$ associated with the eigenvalue $\lambda_1$.
\end{itemize}
\item 
Then the iterates of Algorithm~\ref{alg:power-iteration} satisfy
$$
\norm{\bv^{(k)} - (\pm \bq_1)}= \mathcalO\left(\left|\frac{\lambda_2}{\lambda_1}\right|^k\right)
\qquad\text{and} \qquad 
\abs{\lambda^{(k)}-\lambda_1}= \mathcalO\left(\left|\frac{\lambda_2}{\lambda_1}\right|^{2k}\right)
$$
as $k\rightarrow \infty$.
\end{theorem}
The $\pm$ sign associated with $\bq_1$ signifies that either $\bq_1$ or $-\bq_1$ should be included in the result. 
This implies that  $\bv^{(k)}$ converges linearly to $\pm\bq_1$ in Algorithm~\ref{alg:power-iteration}.


\paragraph{Initial eigenvector assumptions.}
The assumption $\bq_1^\top \bv^{(0)} \neq 0$ in the above algorithm presupposes that the initial guess of the eigenvector possesses a component in the direction of the eigenvector $\bq_1$ that we want to find. If this condition is not met, convergence to  $\pm \bq_1$ cannot be guaranteed.

\paragraph{Dominant eigenvalue assumptions.}
We also  note that the assumption $|\lambda_1|>|\lambda_2|$ does not include an equality sign. In this scenario, $\lambda_1$ is referred to as the \textit{dominant} eigenvalue of matrix $\bA$. Nevertheless, if \{$|\lambda_1|\textcolor{mylightbluetext}{=|}\lambda_2| \textcolor{mylightbluetext}{>} |\lambda_3| \geq |\lambda_4| \geq \ldots $, $\bq_1^\top \bv^{(0)} \neq 0$, and $\bq_2^\top \bv^{(0)} \neq 0 $\}, then $\bv^{(k)}$ will converge to a scalar multiple of $x_1\bq_1 \pm x_2\bq_2$, i.e., lies in the subspace spanned by $\{\bq_1, \bq_2\}$. To see this, 
\begin{enumerate}
\item $\lambda_1=\lambda_2$, i.e., the two dominant eigenvalues are coincident. Equation~\eqref{equation:converge-equa-of-power-ite} shows that  $\bv^{(k)}$ and $\lambda^{(k)}$ converge to 
$$
\left\{
\begin{aligned}
\bv^{(k)} &\stackrel{k\rightarrow \infty}{\longrightarrow} \bbeta_1\triangleq \frac{ x_1 \bq_1 +  x_2\bq_2}{\norm{x_1 \bq_1 +  x_2\bq_2}}\in \spn\{\bq_1, \bq_2\}
\qquad \text{i.e.,} \gap
\bA\bbeta_1 =  \lambda_1 \bbeta_1,\\
\lambda^{(k)} &\stackrel{k\rightarrow \infty}{\longrightarrow} \bbeta_1^\top \bA\bbeta_1 = \lambda_1=\lambda_2.
\end{aligned}
\right.
$$
Therefore, the vector sequence  $\bv^{(k)}$ still converges to an eigenvector of $\bA$, which lies in the space spanned by $\{\bq_1,\bq_2\}$, and the scalar sequence $\lambda^{(k)}$ still converges to $\lambda_1=\lambda_2$. We  notice that when $\lambda_1= \lambda_2$, any vector in $\spn\{\bq_1, \bq_2\}$ qualifies as an eigenvector of $\bA$.

\item $\lambda_1=-\lambda_2$, i.e., the two dominant eigenvalues are of opposite signs. Equation~\eqref{equation:converge-equa-of-power-ite} shows that $\bv^{(k)}$ and $\lambda^{(k)}$ converge to 
$$
\begin{aligned}
\bv^{(k)} &\stackrel{k\rightarrow \infty}{\longrightarrow} \bbeta_2\triangleq \frac{ \lambda_1^k x_1 \bq_1 +  \lambda_2^k x_2\bq_2}{\norm{\lambda_1^k x_1 \bq_1 +  \lambda_2^k x_2\bq_2}}\in \spn\{\bq_1, \bq_2\}.
\end{aligned}
$$
Then  we have
$$
\left\{
\begin{aligned}
\bA\bbeta_2 &=  \lambda_1 \frac{x_1\bq_1+x_2\bq_2}{\norm{x_1\bq_1-x_2\bq_2}}, \qquad \text{when $k$ is odd;}\\
\bA\bbeta_2 &=  \lambda_1 \frac{x_1\bq_1-x_2\bq_2}{\norm{x_1\bq_1+x_2\bq_2}}, \qquad \text{when $k$ is even.}
\end{aligned}
\right.
$$
Therefore, $\bv^{(k)}$ will not converge to an eigenvector of $\bA$, nor will $\lambda^{(k)}$. In this case, we observe that $\bA\bx=\lambda\bx \rightarrow \bA^2\bx=\lambda^2\bx$ such that the eigenvalues of $\bA^2$ will be nonnegative, and $\lambda_1^2=\lambda_2^2$ are the same eigenvalue of $\bA^2$ if $\lambda_1$ and $-\lambda_2$ are eigenvalues of $\bA$. Applying the power method to $\bA^2$ will lead to the convergence to the eigenvector and eigenvalue of $\bA^2$, which is the same situation as discussed in case (1).

\item $\lambda_1=\bar{\lambda}_2$, i.e., the two dominant eigenvalues are complex conjugate. The power method is not convergent  \citep{wilkinson1971algebraic, quarteroni2010numerical}, and we shall not delve into the details.
\end{enumerate}
\paragraph{What if $\bq_1^\top \bv^{(0)} = 0$?} All the convergence results are derived  under the assumption that $\bq_1^\top \bv^{(0)} \neq 0$. 
However, since $\bq_1$ is not known in advance, this requirement cannot always be met. 
In cases where  $\bq_1^\top \bv^{(0)} = 0$, as deduced from  Equation~\eqref{equation:converge-equa-of-power-ite} again, the vector sequence $\bv^{(k)}$ converges to $\bq_2$, leading to $\lambda^{(k)}\rightarrow \lambda_2$. Therefore, the requirement on the initial guess will not impede  the convergence of the power method, although the ultimate outcome differs.

\paragraph{Why do we start from the ``theoretical" one?} In the theoretical power method Algorithm~\ref{alg:power-iteration-theoretical}, we assume that $\lambda_1$ is known and show that the vector sequence converges linearly to the eigenvector in Equation~\eqref{equation:power-iteration-analysis}. This result matches the convergence of the ``practical" power method Algorithm~\ref{alg:power-iteration} as shown in Theorem~\ref{theorem:convergence-power-iteration}. The ``theoretical" one thus can be employed to find an initial analysis on the power method. In the subsequent section, we will explore its counterpart in the  \textit{inverse power method}.

\index{Conditional gradient method}
\index{Frank-Wolfe method}
\subsubsection*{Connection to Conditional Gradient Method$^\ast$}

Consider the problem 
$$
\mathop{\min}_{\bx} f(\bx),\gap \text{s.t.}\gap \bx\in\sS.
$$
The \textit{conditional gradient method (a.k.a., the Frank-Wolfe method)} computes the next step as a convex combination of the current iterate and a minimizer of a linearized version of the objective function over the feasible set $\sS$.

\begin{algorithm}[h] 
\caption{Conditional Gradient Method} 
\label{alg:frank-wolfe}
\begin{algorithmic}[1] 
\Require Consider problem $\mathop{\min}_{\bx} f(\bx),\,\, \text{s.t. } \bx\in\sS$; 
\State Initialize $\bx^{(0)}\in\sS$;
\For{$k=1,2,\ldots$} 
\State Compute $\bt^{(k)} \leftarrow \argmin_{\bt\in\sS} \langle\bt, \nabla f(\bx^{(k)})\rangle$;
\State Choose $s_k\in[0,1]$ and set $\bx^{(k+1)} \leftarrow \bx^{(k)}+s_k(\bt^{(k)} - \bx^{(k)})$;
\EndFor
\end{algorithmic} 
\end{algorithm}

The conditional gradient method achieves an $\mathcalO(1/\sqrt{k})$ rate of convergence when $f(x)$ is non-convex, and an $\mathcalO(1/k)$ rate of convergence when it is convex (for the conditional gradient norm) \citep{luss2013conditional, beck2017first}.

Consider the specific problem 
$$
\mathop{\min}_{\bx} f(\bx)=-\frac{1}{2} \bx^\top\bA\bx,\gap \text{s.t.}\gap \sS=\{\norm{\bx}_2\leq 1\},
$$
where $\bA$ is positive semidefinite. 
Then, $\bt^{(k)} = \argmin_{\bt\in\sS} \langle\bt, \nabla f(\bx^{(k)})\rangle$ is $\frac{\bA\bx^{(k)}}{\norm{\bA\bx^{(k)}}_2}$. And the conditional gradient update is
$$
\bx^{(k+1)} = (1-s_k)\bx^{(k)} + s_k \frac{\bA\bx^{(k)}}{\norm{\bA\bx^{(k)}}_2}.
$$
If we set the step sizes using an exact line search strategy, then 
$$
s_k =\mathop{\argmin}_{{s\in[0,1]}} f(\bx^{(k)}+s_k(\bt^{(k)} - \bx^{(k)})).
$$
When $\bx^{(k)}$ is not the optimal solution of the problem, and since $f(\bx)$ is concave, we can choose $s_k=1$. Therefore, the update becomes 
$$
\bx^{(k+1)} = \frac{\bA\bx^{(k)}}{\norm{\bA\bx^{(k)}}_2}.
$$
This is equivalent to the power iteration in Algorithm~\ref{alg:power-iteration} for finding the eigenvector corresponding to the maximal eigenvalue.

\subsection{The Inverse Power Method}
The \textit{power method} focuses on determining an eigenvector associated with the largest eigenvalue  (in magnitude). 
In contrast, the \textit{inverse power method} aims to find an eigenvector associated with the smallest
eigenvalue (in magnitude). To demonstrate this, we first introduce  a lemma concerning the eigenpair of the inverse of matrices as follows.
\begin{lemma}[Eigenpair of Inverse Matrix]\label{lemma:inverse-eigenpair}
Let $\bA^{n\times n}$ be nonsingular, and let $(\lambda, \bx)$ be an eigenpair of $\bA$. Then, $(1/\lambda, \bx)$ is an eigenpair of $\bA^{-1}$.
\end{lemma}

Furthermore, we previously assumed that the matrix $\bA$ is real and symmetric. We have shown in Chapter~\ref{chapter:spectral-decomposition} that the rank of a real and symmetric matrix is equal to the number of nonzero eigenvalues. The inverse power method involves the inverses of the eigenvalues, so we  further assume that the matrix involved in this section is nonsingular, i.e., all the eigenvalues are nonzero, and the matrix is invertible.

\begin{algorithm}[h] 
\caption{Inverse Power Iteration (Compare to Algorithm~\ref{alg:power-iteration})} 
\label{alg:inverse-power-iteration}
\begin{algorithmic}[1] 
\Require \textcolor{mylightbluetext}{Nonsingular} matrix $\bA\in \real^{n\times n}$ that is real and symmetric; 
\State Initialize $\bv^{(0)}\leftarrow$ some vector with $\norm{\bv^{(0)}}=1$;
\For{$k=1,2,\ldots$} 
\State $\bw \leftarrow \textcolor{mylightbluetext}{\bA^{-1}} \bv^{(k-1)}$;
\State $\bv^{(k)} \leftarrow \bw / \norm{\bw}$;
\State $\lambda^{(k)} \leftarrow (\bv^{(k)})^\top\bA\bv^{(k)}$;  \Comment{i.e., Rayleigh quotient}
\EndFor
\end{algorithmic} 
\end{algorithm}
The inverse power method is formulated in Algorithm~\ref{alg:inverse-power-iteration}. The underlying concept of the algorithm is that an eigenvector associated with the smallest eigenvalue (in magnitude) of $\bA$ corresponds to an eigenvector
associated with the largest eigenvalue (in magnitude)  of $\bA^{-1}$, 
as indicated by Lemma~\ref{lemma:inverse-eigenpair}.
As a result of the power method, it is evident that the sequence $\bv^{(k)}$ in Algorithm~\ref{alg:inverse-power-iteration} converges linearly to the eigenvector of $\bA$ corresponding to the smallest eigenvalue $\lambda_n$ in magnitude (or equivalently, converges linearly to the eigenvector of $\bA^{-1}$ corresponding to the largest eigenvalue $1/\lambda_n$ in magnitude).

To facilitate  the analysis of the convergence of the inverse power method, we can carefully employ a theoretical (but impossible) way  as shown in Algorithm~\ref{alg:power-iteration-theoretical}, where we assume $\lambda_1$ is known and the convergence result is given in Equation~\eqref{equation:power-iteration-analysis}. 
Now, returning to the inverse power method, suppose $\lambda_n$ is known and a ``theoretical" inverse power method is induced, the convergence is thus given by (analogous to Equation~\eqref{equation:power-iteration-analysis}):
\begin{equation}\label{equation:-inverse-power-iteration-analysis1}
\begin{aligned}
\norm{\bv^{(k+1)} - x_n \bq_n}_{\bQ^{-1}} &= \norm{\bQ^{-1}(\bv^{(k+1)} - x_n\bq_n)}_2\\
&\leq \left|\frac{\lambda_n}{\lambda_{n-1}}\right| \cdot \norm{\bQ^{-1}(\bv^{(k)} - x_n\bq_n)}_2 
= \left|\frac{\lambda_n}{\lambda_{n-1}}\right| \cdot \norm{\bv^{(k)} - x_n\bq_n}_{\bQ^{-1}},
\end{aligned}
\end{equation}
where we use the fact that if the spectral decomposition of $\bA$ is $\bA=\bQ\bLambda\bQ^{-1}$, then the spectral decomposition of $\bA^{-1}$ is $\bA^{-1} =\bQ\bLambda^{-1} \bQ^{-1}$. Here, $\lambda_n$ and $\lambda_{n-1}$ are the two smallest  eigenvalues of $\bA$ (in magnitude). However, it is possible (and not uncommon) that the two eigenvalues are very close to each other, leading to $\left|\frac{\lambda_n}{\lambda_{n-1}}\right| $ approaching 1, which results in slow convergence.

\subsection{The Shifted Inverse Power Method}
Now, let's consider a scenario where we have prior knowledge that one of the eigenvalues of matrix $\bA$ is in close proximity to a specific value $\mu\in \real$. 
\textbf{\textcolor{black}{For any value $\mu$ that is not an eigenvalue of $\bA$, the matrix $\bA-\mu\bI$ is nonsingular, even if $\bA$ is singular}}. The matrix $\bA-\mu\bI$ is referred to as \textit{the matrix $\bA$ that has been ``shifted" by $\mu$}, and $\mu$ is called a \textit{shift}.

The following lemma reveals that a shifted version of the inverse power method can be employed to find the eigenvector associated with the eigenvalue  closest to $\mu$.

\begin{lemma}[Eigenpair of Shifted Matrix]\label{lemma:shifted-eigenpair}
Let $(\lambda, \bx)$ be an eigenpair of $\bA\in \real^{n\times n}$, and suppose that $\mu\in \real $ is not an eigenvalue of $\bA$. Then, $(\lambda-\mu, \bx)$ is an eigenpair of $\bA-\mu\bI$.
Note that the eigenvectors of $\bA-\mu\bI$ remain identical to those of $\bA$ since $\bA\bx =\lambda\bx \rightarrow (\bA-\mu\bI)\bx=(\lambda-\mu)\bx$.
\end{lemma}

\begin{algorithm}[h] 
\caption{Shifted Inverse Power Iteration (Compare to Algorithm~\ref{alg:inverse-power-iteration})} 
\label{alg:inverse-power-iteration-shifted}
\begin{algorithmic}[1] 
\Require Matrix $\bA\in \real^{n\times n}$ that is real and symmetric, \textcolor{mylightbluetext}{$\mu$ is not an eigenvalue of $\bA$}; 
\State Initialize $\bv^{(0)}\leftarrow$ some vector with $\norm{\bv^{(0)}}=1$;
\For{$k=1,2,\ldots$} 
\State $\bw \leftarrow \textcolor{mylightbluetext}{(\bA-\mu\bI)^{-1}} \bv^{(k-1)}$;
\State $\bv^{(k)} \leftarrow \bw / \norm{\bw}$;
\State $\lambda^{(k)} \leftarrow (\bv^{(k)})^\top\bA\bv^{(k)}$;  \Comment{i.e., Rayleigh quotient}
\EndFor
\end{algorithmic} 
\end{algorithm}

The procedure is outlined in Algorithm~\ref{alg:inverse-power-iteration-shifted}. To facilitate  the analysis of convergence again, suppose $\mu$ is closest to the smallest eigenvalue $\lambda_n$ in magnitude, the convergence is (analogous to Equation~\eqref{equation:power-iteration-analysis}):
\begin{equation}\label{equation:-inverse-power-iteration-analysis-shifted}
\begin{aligned}
&\norm{\bv^{(k+1)} - x_n \bq_n}_{\bQ^{-1}} = \norm{\bQ^{-1}(\bv^{(k+1)} - x_n\bq_n)}_2\\
&\gap\gap\gap\gap\leq \left|\frac{\lambda_n\textcolor{mylightbluetext}{-\mu}}{\lambda_{n-1}\textcolor{mylightbluetext}{-\mu}}\right| \cdot \norm{\bQ^{-1}(\bv^{(k)} - x_n\bq_n)}_2 
= \left|\frac{\lambda_n \textcolor{mylightbluetext}{-\mu}}{\lambda_{n-1}\textcolor{mylightbluetext}{-\mu}}\right| \cdot \norm{\bv^{(k)} - x_n\bq_n}_{\bQ^{-1}}.
\end{aligned}
\end{equation}
When $\mu$ is close to $\lambda_n$, the bound $\left|\frac{\lambda_n \textcolor{mylightbluetext}{-\mu}}{\lambda_{n-1}\textcolor{mylightbluetext}{-\mu}}\right|$ becomes small such that the convergence is faster compared to the (naive) inverse power method (see Equation~\eqref{equation:-inverse-power-iteration-analysis1}), although it remains linearly convergent.

The formal convergence result is presented in the following theorem (one can follow the deduction as demonstrated in Theorem~\ref{theorem:convergence-power-iteration} to obtain the result).

\begin{theorem}[Convergence of Shifted Inverse Power Iteration]
Suppose $\lambda_J$ is the eigenvalue closest  to $\mu$ and $\lambda_K$ is the second closest. Moreover, assume $\bq_J^\top \bv^{(0)}\neq 0$. Then, the iterates of Algorithm~\ref{alg:inverse-power-iteration-shifted} satisfy
$$
\norm{\bv^{(k)} - (\pm \bq_J)}= \mathcalO\left(\left|\frac{\lambda_J-\mu}{\lambda_K-\mu}\right|^k\right) 
\qquad\text{and}\qquad
\abs{\lambda^{(k)}-\lambda_J} =\mathcalO\left(\left|\frac{\lambda_J-\mu}{\lambda_K-\mu}\right|^{2k}\right)
$$
as $k\rightarrow \infty$.
\end{theorem}
And this indicates $\bv^{(k)}$ converges to $\pm\bq_J$ linearly for Algorithm~\ref{alg:inverse-power-iteration-shifted}.


\subsection{The Rayleigh Quotient Method}
We know that the \textit{inverse power iteration} converges to the eigenvector corresponding to the smallest eigenvalue of $\bA$ in magnitude, and the \textit{shifted inverse power iteration} converges to the eigenvector corresponding to the eigenvalue closest to $\mu$ with  potentially faster convergence. Both methods are, in a sense,  \textit{inverse} power iteration. 
By combining the concepts from both algorithms, i.e., using the Rayleigh quotient of the estimated eigenvector in each iteration as the estimate of the eigenvalue, we can develop a faster algorithm; see Algorithm~\ref{alg:rayleigh-quotient-iteration2}.

\begin{algorithm}[h] 
\caption{Rayleigh Quotient Iteration (Compare to Algorithm~\ref{alg:inverse-power-iteration-shifted})} 
\label{alg:rayleigh-quotient-iteration2}
\begin{algorithmic}[1] 
	\Require A real and symmetric $\bA\in \real^{n\times n}$;
	\State Initialize $\bv^{(0)}\leftarrow$ some vector with $\norm{\bv^{(0)}}=1$;
	\State $\lambda^{(0)}\leftarrow(\bv^{(0)})^\top\bA\bv^{(0)}$; \Comment{i.e., Rayleigh quotient}  
	\For{$k=1,2,\ldots$} 
	\State $\bw \leftarrow \textcolor{mylightbluetext}{(\bA-\lambda^{(k-1)}\bI)^{-1}} \bv^{(k-1)}$;
	\State $\bv^{(k)} \leftarrow \bw / \norm{\bw}$;
	\State $\lambda^{(k)} \leftarrow (\bv^{(k)})^\top\bA\bv^{(k)}$;  \Comment{i.e., Rayleigh quotient} 
	\EndFor
\end{algorithmic} 
\end{algorithm}

The Rayleigh quotient iteration finds an eigenvector; however, it is not initially clear which eigenvalue it is associated with.

\section{QR Algorithm}\label{section:qr_algo}
The QR algorithm for computing the eigenvalues and eigenvectors of matrices has been recognized as one of the ten most important algorithms of the twentieth century \citep{dongarra2000guest, cipra2000best}.
This algorithm  was published by \textit{John G. F. Francis} in the works of \citet{francis1961qr, francis1962qr} that is quoted as one of the \textit{jewels of numerical analysis} \citep{trefethen1997numerical}. 
In this discussion, we will introduce the QR algorithm, starting from its simplest form and progressing to a shifted version with implicit calculations.

The QR algorithm advances the process of \textit{simultaneously} calculating the eigenvalues of a given matrix $\bA$. The central idea is to reduce the matrix $\bA$ through a series of \textit{similarity transformations} (Definition~\ref{definition:similar-matrices}) that preserve the eigenvalues (Proposition~\ref{proposition:eigenvalue-similar-matrices}), making them easier to compute. 
The resulting procedure is straightforward: the QR algorithm involves taking a QR decomposition, multiplying the computed factors $\bQ$ and $\bR$ together in reverse order as $\bR\bQ$, and repeating the process. Hence the name \textit{QR algorithm}.

\subsection{Preliminary: Power Iteration with Known Eigenvector}
We will first illustrate how the QR algorithm evolves from the power iteration algorithm.
We will use subscripts for all vectors $\bv^{(k)}$'s to emphasize their convergence to the eigenvectors associated with $\lambda_i$'s: for example, $\bv^{(k)}_{\textcolor{mylightbluetext}{i}}$ will be demonstrated to converge to the eigenvector corresponding to the eigenvalue $\lambda_{\textcolor{mylightbluetext}{i}}$. 
Additionally, suppose we have prior knowledge that the normalized eigenvector $\bq_1$ is associated with $\lambda_1$. 
For a second initial vector, $\bv_2^{(0)}$, it is essential that it does not contain a component in the direction of $\bq_1$. 
This can be achieved by making it  orthogonal to  $\bq_1$, which can be accomplished through the projection method introduced in the Gram-Schmidt process (Section~\ref{section:project-onto-a-vector}): 
$$
\bv_2^{(k+1)} \leftarrow \bv_2^{(k)} -\bq_1^\top \bv_2^{(k)}\bq_1
$$ 
such that $\bq_1^\top \bv_2^{(k+1)}=0$ (recall that $\bq_1^\top \bv_2^{(k)}\bq_1$ represents the component of $ \bv_2^{(k)}$ in the direction of $\bq_1$, given that $\bq_1$ has  unit length). With $\bq_1$ known beforehand, we consider the method outlined in Algorithm~\ref{alg:power-iteration-in-qr-preliminary1}.

\begin{algorithm}[H] 
\caption{Power Iteration (\textbf{assume $\bq_1$ is Known Up Front})} 
\label{alg:power-iteration-in-qr-preliminary1}
\begin{algorithmic}[1] 
\Require A real and symmetric $\bA\in \real^{n\times n}$;
\State The normalized eigenvector \textcolor{mylightbluetext}{$\bq_1$ is known up front}, corresponding to eigenvalue $\lambda_1$;
\State Initialize $\bv_2^{(0)}\leftarrow$ some vector in $\real^n$;
\State $\bv_2^{(0)} \leftarrow \bv_2^{(0)} -\textcolor{mylightbluetext}{\bq_1}^\top \bv_2^{(0)} \textcolor{mylightbluetext}{\bq_1}$; \Comment{i.e., project along $\bq_1$: $\bq_1^\top \bv_2^{(0)}=0$}
\State $\bv_2^{(0)} \leftarrow \bv_2^{(0)} \big/\norm{\bv_2^{(0)}}$;  \Comment{normalize to have length one} 
\For{$k=1,2,\ldots$} 
\State $\bv_2^{(k)}\leftarrow \bA \bv_2^{(k-1)}$;
\State $\bv_2^{(k)} \leftarrow \bv_2^{(k)} -\textcolor{mylightbluetext}{\bq_1}^\top \bv_2^{(k)} \textcolor{mylightbluetext}{\bq_1}$; \Comment{make sure $\bv_2^{(k)}$ is orthogonal to $\bq_1$}
\State $\bv_2^{(k)} \leftarrow \bv_2^{(k)}\big/ \norm{\bv_2^{(k)}}$;
\State $\lambda_2^{(k)} \leftarrow (\bv_2^{(k)})^\top\bA\bv_2^{(k)}$;  \Comment{i.e., Rayleigh quotient}
\EndFor
\end{algorithmic} 
\end{algorithm}

Once again, we express $\bv_2^{(0)}$ as a linear combination of the orthonormal eigenvectors $\{\bq_i\}$:
\begin{equation}\label{equation:power-qr-premi-1}
\bv_2^{(0)} = x_1 \bq_1 + x_2\bq_2 + \ldots + x_n\bq_n.
\end{equation}
Since in step 3 of the Algorithm~\ref{alg:power-iteration-in-qr-preliminary1}, we project  the vector $\bv_2^{(0)}$ along $\bq_1$, then the component $x_1$ of Equation~\eqref{equation:power-qr-premi-1} is equal to zero.
Similarly, $\bv_2^{(k)}$ can be expressed as a scalar multiple of $\bA^k \bv_2^{(0)}$ such that  $\bv_2^{(k)} = c_{2k} \bA^k \bv_2^{(0)} = \frac{\bA^k \bv_2^{(0)}}{\norm{\bA^k \bv_2^{(0)}}}$ (see Problem~\ref{prob:power_nai1}). We have 
\begin{equation}\label{equation:power_nai1}
\begin{aligned}
	\bv_2^{(k)} &=  c_{2k} \bA^k\bv_2^{(0)} 
	=c_{2k} (x_1 \lambda_1^k\bq_1 + x_2\lambda_2^k\bq_2+x_3\lambda_3^k\bq_3 + \ldots + x_n\lambda_n^k\bq_n) \\
	&=c_{2k} (x_2\lambda_2^k\bq_2 + x_3\lambda_3^k\bq_3+ \ldots + x_n\lambda_n^k\bq_n) \\
	&= c_{2k} \lambda_2^k\left(x_2\bq_2 + x_3\left(\frac{\lambda_3}{\lambda_2}\right)^k\bq_3 + \ldots + x_n\left(\frac{\lambda_n}{\lambda_1}\right)^k\bq_n\right).
\end{aligned}
\end{equation}
Therefore, following from Theorem~\ref{theorem:convergence-power-iteration}, if we assume $|\lambda_2| \textcolor{mylightbluetext}{>} |\lambda_3| \geq |\lambda_3| \geq \ldots \geq |\lambda_n| \geq 0$ and $\bq_2^\top \bv_2^{(0)}\neq 0$, $\bv_2^{(k)}$ will converge linearly to $\pm \bq_2$, and consequently, $\lambda_2^{(k)}$ will converge to $\lambda_2$ using Algorithm~\ref{alg:power-iteration-in-qr-preliminary1}.


Analogously, when \{$|\lambda_2|\textcolor{mylightbluetext}{=|}\lambda_3| \textcolor{mylightbluetext}{>} |\lambda_4| \geq |\lambda_5| \ldots $, $\bq_2^\top \bv_2^{(0)} \neq 0$, and $\bq_3^\top \bv_2^{(0)} \neq 0$\}, $\bv_2^{(k)}$ will converge to a scalar multiple of $x_2\bq_2 \pm x_3\bq_3$, i.e., lies in the space spanned by $\{\bq_2, \bq_3\}$ (see discussion in Section~\ref{section:power_method}).

\subsection{Preliminary: Power Iteration with Unknown Eigenvector}
However, the approach described in Algorithm~\ref{alg:power-iteration-in-qr-preliminary1} is not practical because we typically do not know $\bq_1$ up front. 
However, since a simple power iteration can compute $\bv_1^{(0)}$ that converges to $\pm \bq_1$ (Algorithm~\ref{alg:power-iteration}), we can construct a simultaneous algorithm. 
Instead of projecting $\bv_2^{(k)}$ along $\bq_1$, we can project it along $\bv_1^{(k)}$ in Algorithm~\ref{alg:power-iteration}. 
Thus, a method to find both $\bq_1$ and $\bq_2$ simultaneously can be constructed as outlined in Algorithm~\ref{alg:power-iteration-in-qr-preliminary2-unknown-q1}.

\begin{algorithm}[htp] 
\caption{Power Iteration ($\bq_1$ is Unknown Up Front, Compare to Algorithm~\ref{alg:power-iteration-in-qr-preliminary1})} 
\label{alg:power-iteration-in-qr-preliminary2-unknown-q1}
\begin{algorithmic}[1] 
\Require A real and symmetric $\bA\in \real^{n\times n}$;
\State Initialize $\bv_1^{(0)}, \bv_2^{(0)}\leftarrow$ two random vectors in $\real^n$;
\State $\bv_2^{(0)} \leftarrow \bv_2^{(0)} -\textcolor{mylightbluetext}{\bv_1^{(0)}}^\top \bv_2^{(0)} \textcolor{mylightbluetext}{\bv_1^{(k)}}$; \Comment{i.e., project along $\bv_1^{(0)}$: \textcolor{mylightbluetext}{$\bv_1^{(0)\top} \bv_2^{(0)}=0$}}
\State $\bv_1^{(0)} \leftarrow \bv_1^{(0)} \big/\norm{\bv_1^{(0)}}$;  \Comment{normalize to have length one} 
\State $\bv_2^{(0)} \leftarrow \bv_2^{(0)} \big/\norm{\bv_2^{(0)}}$;  \Comment{normalize to have length one} 
\For{$k=1,2,\ldots$} 
\State // update $\bv_1^{(k)}$
\State $\bv_1^{(k)}  \leftarrow \bA\bv_1^{(k-1)}$;
\State $\bv_1^{(k)} \leftarrow \bv_1^{(k)} \big/ \norm{\bv_1^{(k)}}$;
\State // update $\bv_2^{(k)}$
\State $\bv_2^{(k)} \leftarrow \bA \bv_2^{(k-1)}$;
\State $\bv_2^{(k)} \leftarrow \bv_2^{(k)} -\textcolor{mylightbluetext}{\bv_1^{(k)}}^\top \bv_2^{(k)} \textcolor{mylightbluetext}{\bv_1^{(k)}}$; \Comment{make sure $\bv_2^{(k)}$ is orthogonal to \textcolor{mylightbluetext}{$\bv_1^{(k)}$}}
\State $\bv_2^{(k)} \leftarrow \bv_2^{(k)}\big/ \norm{\bv_2^{(k)}}$;
\State // compute the corresponding eigenvalues
\State $\lambda_1^{(k)} \leftarrow (\bv_1^{(k)})^\top\bA\bv_1^{(k)}$;  \Comment{i.e., Rayleigh quotient}
\State $\lambda_2^{(k)} \leftarrow (\bv_2^{(k)})^\top\bA\bv_2^{(k)}$;  \Comment{i.e., Rayleigh quotient}
\EndFor
\end{algorithmic} 
\end{algorithm}
Combining the findings above, the iterates of Algorithm~\ref{alg:power-iteration-in-qr-preliminary2-unknown-q1} satisfy that (when $\bq_1^\top \bv_1^{(0)} \neq 0$ and $\bq_2^\top \bv_2^{(0)}\neq 0$)
\begin{itemize}
\item If $|\lambda_1| > |\lambda_2|$, the vector sequence $\bv_1^{(k)}$ will converge linearly to $\pm \bq_1$ at a rate of $\left|\frac{\lambda_2}{\lambda_1}\right|$;
\item If $|\lambda_1| > |\lambda_2|> |\lambda_3|$, the vector sequence $\bv_2^{(k)}$ will converge linearly to $\pm \bq_2$  at a rate of $\left|\frac{\lambda_3}{\lambda_2}\right|$;
\item If $|\lambda_1| = |\lambda_2| > |\lambda_3|$, the vector sequence $\bv_1^{(k)} $ will converge to a scalar multiple of $x_1\bq_1 \pm x_2\bq_2$, and the vector sequence $\bv_2^{(k)}$ will converge linearly to $\pm \bq_2$. In other words, the span of $\{\bq_1, \bq_2\}$ can be approximated by the span of $\{\bv_1^{(k)}, \bv_2^{(k)}\}$.
\end{itemize}

\subsection{Preliminary: Power Iteration with Unknown Eigenvector and QR}\label{section:power-eigen-unknown-prelimise}
We observe that  steps 2 to  4 of Algorithm~\ref{alg:power-iteration-in-qr-preliminary2-unknown-q1} are equivalent to applying a QR decomposition on an $\real^{n\times 2}$ matrix $[\bv_1^{(0)}, \bv_2^{(0)}]$, which can be written as:
$$
\underbrace{[\bv_1^{(0)}, \bv_2^{(0)}] }_{\widehat{\bV}^{(0)}}, \bR \leftarrow QR([\bv_1^{(0)}, \bv_2^{(0)}]),
$$
where $QR(\cdot)$ denotes the function to obtain the QR decomposition of a matrix.
Moreover, steps 6 to  12 of Algorithm~\ref{alg:power-iteration-in-qr-preliminary2-unknown-q1} can also be restated as performing a QR decomposition on the $n\times 2$ matrix $\widehat{\bV}^{(k)},\bR\leftarrow QR(\bA[\bv_1^{(k-1)}, \bv_2^{(k-1)}])\triangleq QR(\bA \widehat{\bV}^{(k-1)})$. A further simplification on the form is to obtain the eigenvalues $\lambda_1^{(k)}, \lambda_2^{(k)} $ as the diagonals of the following matrix $\widehat{\bA}^{(k)}$:
$$
\underbrace{\begin{bmatrix}
\lambda_1^{(k)}  & \boxtimes\\
\boxtimes & \lambda_2^{(k)} 
\end{bmatrix}}_{\triangleq\widehat{\bA}^{(k)}}
\leftarrow 
\underbrace{\begin{bmatrix}
\bv_1^{(k)} & \bv_1^{(k)} 
\end{bmatrix}^\top }_{\widehat{\bV}^{(k)\top}}
\bA
\underbrace{\begin{bmatrix}
\bv_1^{(k)} & \bv_1^{(k)} 
\end{bmatrix}}_{\widehat{\bV}^{(k)}}
$$
We then reformulate  Algorithm~\ref{alg:power-iteration-in-qr-preliminary2-unknown-q1} using this QR decomposition as shown in  Algorithm~\ref{alg:power-iteration-in-qr-preliminary2-unknown-qr-function}.

\begin{algorithm}[H] 
\caption{Power Iteration (On 2 Vectors, Equivalent to Algorithm~\ref{alg:power-iteration-in-qr-preliminary2-unknown-q1})} 
\label{alg:power-iteration-in-qr-preliminary2-unknown-qr-function}
\begin{algorithmic}[1] 
\Require A real and symmetric $\bA\in \real^{n\times n}$;
\State Initialize  $\widehat{\bV}^{(0)}\leftarrow [\bv_1^{(0)}, \bv_2^{(0)} ]\in \real^{n\times 2}\leftarrow$ two random vectors in $\real^n$;
\State $\widehat{\bV}^{(0)}, \cancel{\bR}\leftarrow QR(\widehat{\bV}^{(0)})$;
\For{$k=1,2,\ldots$} 
\State $\widehat{\bV}^{(k)}, \cancel{\bR}\leftarrow QR(\bA\widehat{\bV}^{(k-1)})$;
\State $\widehat{\bA}^{(k)}\leftarrow\widehat{\bV}^{(k)\top} \bA \widehat{\bV}^{(k)}$; \Comment{ compute the corresponding eigenvalues}
\EndFor
\end{algorithmic} 
\end{algorithm}

Note that the $\widehat{\text{widehat}}$ above the matrices in Algorithm~\ref{alg:power-iteration-in-qr-preliminary2-unknown-qr-function} helps differentiate them from those in QR algorithms (see Section~\ref{section:qralg-without-shifts}).
With hindsight, it is natural to extend this algorithm to handle any number of vectors $p\leq n$. The complete algorithm is formulated in Algorithm~\ref{alg:power-iteration-in-qr-preliminary2-unknown-qr-function-mn}.
\begin{algorithm}[H] 
\caption{Power Iteration (On $p$ Vectors, Compare to Algorithm~\ref{alg:power-iteration-in-qr-preliminary2-unknown-qr-function})}
\label{alg:power-iteration-in-qr-preliminary2-unknown-qr-function-mn}
\begin{algorithmic}[1] 
\Require A real and symmetric $\bA\in \real^{n\times n}$; 
\State Initialize  $\widehat{\bV}^{(0)}\leftarrow$ random matrix in $\real^{n\times \textcolor{mylightbluetext}{p}}$;
\State $\widehat{\bV}^{(0)}, \cancel{\bR}\leftarrow QR(\widehat{\bV}^{(0)})$;
\For{$k=1,2,\ldots$} 
\State $\widehat{\bV}^{(k)}, \cancel{\bR}\leftarrow QR(\bA\widehat{\bV}^{(k-1)})$;
\State $\widehat{\bA}^{(k)}\leftarrow \widehat{\bV}^{(k)\top} \bA \widehat{\bV}^{(k)}$; \Comment{ compute the corresponding eigenvalues}
\EndFor
\end{algorithmic} 
\end{algorithm}
Again, we observe that, when $\bq_i^\top \bv_i^{(0)} \neq 0$ for all $i \in \{1,2,\ldots, p\}$, \footnote{That is, the initial guess $\bv_i^{(0)}$ is not orthogonal to eigenvector $\bq_i$, and contains a component in the direction of the eigenvector.} the iterates of Algorithm~\ref{alg:power-iteration-in-qr-preliminary2-unknown-qr-function-mn} satisfy that 
\begin{enumerate}
\item If $|\lambda_1| > |\lambda_2| > \ldots > |\lambda_p| > |\lambda_{p+1}| \geq |\lambda_{p+2}| \geq \ldots $, then each column $i$ of $\widehat{\bV}^{(k)}$ (i.e., $\bv_i^{(k)}$) will converge linearly to $\pm \bq_i$, where the rate of removing from the component in the direction of $\bq_j$ is recorded as $\left|\frac{\lambda_i}{\lambda_j}\right|$ ($0<i\leq p $ and $p<j\leq n$);
\item If some of the eigenvalues have an equal magnitude, then the subspace spanned by the corresponding columns of $\widehat{\bV}^{(k)}$ will approximate the subspace spanned by the corresponding eigenvectors associated with those eigenvalues;
\item If $p=n$, and $|\lambda_1| > |\lambda_2| > \ldots >|\lambda_n|$, then we will find all the eigenvectors of $\bA$.
\end{enumerate}

\subsection{A Simple QR Algorithm from Power Iteration: without Shifts}\label{section:qralg-without-shifts}
Up to this point, we have developed the power iteration into an algorithm that is capable of finding all the eigenvectors (under mild conditions) by employing QR decomposition.

Now, let's consider a slight modification to the power iteration, as shown in Algorithm~\ref{alg:qr-algorithm-simple1}.
This variant treats the $\bR$ matrix from the QR decomposition as a sequence of matrices, and initializes $\widehat{\bA}^{(0)}$  to be the original matrix $\bA$ (where we ignored the value indexed by 0 in Algorithm~\ref{alg:power-iteration-in-qr-preliminary2-unknown-qr-function-mn}). Compare the  Algorithms~\ref{alg:qr-algorithm-simple1} and \ref{alg:qr-algorithm-simple2} below.

\noindent
\begin{minipage}[t]{0.495\linewidth}
\begin{algorithm}[H] 
\caption{Power Iteration}
\label{alg:qr-algorithm-simple1}
\begin{algorithmic}[1] 
\Require $\bA\in \real^{n\times n}$ is real and symmetric;
\State Same as Algorithm~\ref{alg:power-iteration-in-qr-preliminary2-unknown-qr-function-mn};
\State $\widehat{\bA}^{(0)} \leftarrow \bA$; 
\State $\widehat{\bV}^{(0)}\leftarrow\bI_n$; \Comment{initial eigenvector guess}
\State $\widehat{\bR}^{(0)} \leftarrow \bI_n$; \Comment{compensate the sequence}
\For{$k=1,2,\ldots$} 
\State $\widehat{\bV}^{(k)}, \widehat{\bR}^{(k)}\leftarrow QR(\bA\widehat{\bV}^{(k-1)})$;
\State $\widehat{\bA}^{(k)}\leftarrow\widehat{\bV}^{(k)\top} \bA \widehat{\bV}^{(k)}$;
\EndFor
\end{algorithmic} 
\end{algorithm}
\end{minipage}%
\hfil 
\begin{minipage}[t]{0.495\linewidth}
\begin{algorithm}[H] 
	\caption{Simple QR Algorithm}
	\label{alg:qr-algorithm-simple2}
	\begin{algorithmic}[1] 
		\Require $\bA\in \real^{n\times n}$ is real and symmetric; 
		\State $\bA^{(0)} \leftarrow \bA$; 
		\State $\bV^{(0)}\leftarrow\bI_n$;  \Comment{initial eigenvector guess}
		\State $\bR^{(0)} \leftarrow \bI_n$; \Comment{compensate the sequence}
		\State $\bQ^{(0)} \leftarrow \bI_n$; \Comment{compensate the sequence}
		\For{$k=1,2,\ldots$} 
		\State $\bQ^{(k)}, \bR^{(k)}\leftarrow QR(\bA^{(k-1)})$;
		\State $\bA^{(k)} \leftarrow \bR^{(k)} \bQ^{(k)}$; 
		\State $\bV^{(k)} \leftarrow \bV^{(k-1)}\bQ^{(k)} $;
		\EndFor
	\end{algorithmic} 
\end{algorithm}
\end{minipage}

We will first demonstrate the equivalence of the two algorithms through the following lemma.
\begin{lemma}[QR Algorithm from Power Iteration]\label{lemma:qr-algo-from-power}
The two algorithms, Algorithm~\ref{alg:qr-algorithm-simple1} and Algorithm~\ref{alg:qr-algorithm-simple2}, are equivalent in the sense that for all iterations $k\in \{0,1,2,\ldots\}$, we have 
$$
\left\{
\begin{aligned}
	\widehat{\bA}^{(k)} &= \bA^{(k)};  \qquad \text{(diagonals converge to eigenvalues)}\\
	\widehat{\bR}^{(k)} &= \bR^{(k)}; \\
	\widehat{\bV}^{(k)} &= \bV^{(k)}. \qquad \text{(columns converge to eigenvectors)}
\end{aligned}
\right.
$$
\end{lemma}
For clarity, the proof is postponed to Section~\ref{section:proofs-qralgorithms}.
You might wonder how we arrived at the simplification from Algorithm~\ref{alg:qr-algorithm-simple1} to Algorithm~\ref{alg:qr-algorithm-simple2}. The answer is unglamorous! It was by trial and error. 

\index{Similarity transformation}
\paragraph{What's in the QR algorithm?}
All we do in Algorithm~\ref{alg:qr-algorithm-simple2} is to perform a QR decomposition, 
followed by a multiplication of the computed factors $\bQ$ and $\bR$ in  reverse order, denoted by $\bR\bQ$, and then repeat this process.  
For the convergence to a diagonal form, which allows us to find the eigenvalues, the operations from $\bA^{(k-1)}$ to $\bA^{(k)}$ should be similarity transformations (Definition~\ref{definition:similar-matrices}).
According to Algorithm~\ref{alg:qr-algorithm-simple2}, it can be shown that 
\begin{tcolorbox}[title={Simple QR Algorithm Property 1},colback=\mdframecolorTheorem]
\begin{equation}\label{equation:simple-qr-find}
\text{(SQR 1)} \qquad \bA^{(k)} = \bQ^{(k)\top } \bA^{(k-1)} \bQ^{(k)},
\end{equation}
\end{tcolorbox}
\noindent since $\bA^{(k)}=\bR^{(k)} \bQ^{(k)} = (\bQ^{(k)\top} \bA^{(k-1)}) \bQ^{(k)} $. Therefore, $\bA^{(k)}$ is a \textit{similarity transformation}  \footnote{Here, when the nonsingular matrix is orthogonal, it is also known as an \textit{orthogonal similarity transformation}.} of $\bA^{(k-1)}$, which is also a similar transformation from $\bA$. 
Consequently, $\bA^{(k)}$ and $\bA$ share the \textit{same eigenvalues, trace, and rank} (Proposition~\ref{proposition:eigenvalue-similar-matrices}). By tracing back to the initial induction step in  Equation~\eqref{equation:simple-qr-find}, it follows that
\begin{tcolorbox}[title={Simple QR Algorithm Property 2},colback=\mdframecolorTheorem]
\begin{equation}\label{equation:simple-qr-find-root}
\begin{aligned}
\text{(SQR 2)} \qquad \bA^{(k)} &= \bQ^{(k)\top } \bA^{(k-1)} \bQ^{(k)} \\
&=\bQ^{(k)\top } \left(\bQ^{(k-1)\top } \bA^{(k-2)} \bQ^{(k-1)} \right) \bQ^{(k)} \\
&=\ldots \\
&=
\underbrace{\bQ^{(k)\top } \bQ^{(k-1)\top }\ldots \bQ^{(0)\top } }_{\bV^{(k)\top}}
   \bA
\underbrace{\bQ^{(0)} \ldots  \bQ^{(k-1)}  \bQ^{(k)} }_{\bV^{(k)}}
\\
&=\bV^{(k)\top}\bA \bV^{(k)}.
\end{aligned}
\end{equation}
\end{tcolorbox}

The  lemma above states that the simple QR algorithm is equivalent to the power iteration algorithm, under which the diagonal elements of $\bA^{(k)}$ converge to the eigenvalues of $\bA$ (under mild conditions). 
The expression $\bA^{(k)} =\bV^{(k)\top}\bA \bV^{(k)} $   not only represents a similarity transformation but also an \textit{orthogonal similarity transformation}, due to the orthogonality of $\bV^{(k)}$. 
This fact holds particular significance for the stability of the iterative method, as the condition of $\bA^{(k)}$ is not worse than that of the original matrix $\bA$. 
Subsequently,  the equivalence indicates that the $i$-th column of $\bV^{(k)}$  converges linearly to the $i$-th eigenvector of $\bA$, i.e., $\pm \bq_i$. 
Upon a closer examination of the process, it shows by Algorithms~\ref{alg:qr-algorithm-simple1} and~\ref{alg:qr-algorithm-simple2} that 
\begin{tcolorbox}[title={Simple QR Algorithm Property 3},colback=\mdframecolorTheorem]
\begin{equation}\label{equation:simple-qr-find2}
\begin{aligned}
&\left.
\begin{aligned}
	\bV^{(k)} &= \bV^{(0)}\bQ^{(1)}\bQ^{(2)}\ldots \bQ^{(k)}    \\
	\bA^{\textcolor{mylightbluetext}{k}}&=\bV^{(k)} \bR^{(k)} \bR^{(k-1)} \ldots \bR^{(0)} 
\end{aligned}
\right\} 
\implies  \\
&\text{(SQR 3)} \qquad \bA^k= \underbrace{\bQ^{(0)}\bQ^{(1)}\bQ^{(2)}\ldots \bQ^{(k)}}_{=\bV^{(k)}, \text{orthogonal}}
\underbrace{\bR^{(k)} \bR^{(k-1)} \ldots \bR^{(0)}}_{\triangleq\bU^{(k)}, \text{upper triangular}}
,
\end{aligned}
\end{equation}
\end{tcolorbox}
\noindent where the left-hand side of the above equation is delayed to be proved in Section~\ref{section:proofs-qralgorithms}.
Additionally, we set $\bQ^{(0)} = \bV^{(0)} = \bI$ for simplicity.
Thus, $\bA^{k}$ can be expressed as a QR decomposition by $\bA^{k}=\bV^{(k)}\bU^{(k)}$.\footnote{Note the difference between the notation $\bA^k$ and $\bA^{(k)}$, where $\bA^k$ represents the $k$-th power of $\bA$.} 

The equivalence between  the simultaneous power method and the simple QR algorithm (Algorithms~\ref{alg:qr-algorithm-simple1} and \ref{alg:qr-algorithm-simple2}) provides valuable insights into their convergence behavior.
Similar to the discussion in Section~\ref{section:power-eigen-unknown-prelimise}, we can conclude that
\begin{enumerate}
\item If $|\lambda_1| > |\lambda_2| > \ldots > |\lambda_p| > |\lambda_{p+1}| \geq |\lambda_{p+2}| \geq \ldots \geq |\lambda_n|$, then each column $i$ of $\bV^{(k)}$ (i.e., $\bv_i^{(k)}$)  converges linearly to $\pm \bq_i$, where the rate of removing from the component in the direction of $\bq_j$ is recorded as $\left|\frac{\lambda_i}{\lambda_j}\right|$ ($0<i\leq p $ and $p<j\leq n$);
\item If some  eigenvalues have an equal magnitude,  the subspace spanned by the corresponding columns of $\bV^{(k)}$ will approximate the subspace spanned by the  associated   eigenvalues;
\item If $p=n$ and $|\lambda_1| > |\lambda_2| > \ldots >|\lambda_n|$ (all eigenvalues are distinct), then we will find all the eigenvectors of $\bA$. Specifically, if 
$$
\bA^{(k)} = 
\begin{bmatrix}
\lambda_1^{(k)} & a_{12} & a_{13} & \ldots & a_{1n}\\
a_{12} & \lambda_2^{(k)} & a_{23} & \ldots & a_{2n}\\
\vdots & \vdots & \ddots & \ddots & \vdots\\
a_{1,n-1} &  a_{2,n-1} &a_{3,n-1}  & \ddots & a_{n-1,n}\\
a_{1n} &  a_{2n} &a_{3n}  & \ldots & \lambda_n^{(k)}\\
\end{bmatrix},
$$
which is symmetric,  the  rate of convergence  is
$
|a_{i,i-1}| = \mathcalO\left(\left|\frac{\lambda_i}{\lambda_{i-1}}\right|^k\right).
$
\end{enumerate}
Under different conditions, say, $|\lambda_1|\geq |\lambda_2| \geq \ldots \geq |\lambda_p| > |\lambda_{p+1}| > \ldots |\lambda_{n}|$, then
\begin{enumerate}
\item The span of the first $p$ columns of $\bV^{(k)}$  converges to the subspace spanned by the first $p$ orthonormal eigenvectors, $\spn\{\bq_1, \bq_2, \ldots, \bq_p\}$;
\item The span of the last $n-p$ columns of $\bV^{(k)}$  converges to the subspace spanned by the last $n-p$ orthonormal eigenvectors, $\spn\{\bq_{p+1}, \bq_{p+2}, \ldots, \bq_n\}$;
\item Since  $\bA$ is assumed to be  real and symmetric, it follows that
 $$\spn\{\bq_1, \bq_2, \ldots, \bq_p\} \perp \spn\{\bq_{p+1}, \bq_{p+2}, \ldots, \bq_n\} .$$
Therefore, the first $p$ columns and last $n-p$ columns of $\bV^{(k)}$ reside within two mutually orthogonal complement subspaces. And the rate of convergence at which the two subspaces become orthogonal to each other is linear, characterized by the constant $|\lambda_{p+1}/\lambda_p|$;
\item Furthermore, in the scenario where $p=n$, it follows that the last column of $\bV^{(k)}$ converges linearly to $\pm \bq_n$, featuring a constant rate denoted by $\mathcalO(|\lambda_n/\lambda_{n-1}|)$.
\end{enumerate}
Based on the discoveries made through the shifted inverse power method (Algorithm~\ref{alg:inverse-power-iteration-shifted}), it has been established that introducing a shift equal to the estimated smallest eigenvalue (in magnitude) during each iteration can expedite the convergence of the algorithm. 
This insight reveals a ``practical" QR algorithm, which will be introduced in the next section.

\begin{remark}[Asymmetric Matrix $\bA$]
In the preceding  discussions, we assume that $\bA$ is real and symmetric. For an asymmetric matrix $\bA$ with distinct real eigenvalues, it can be shown that the $\bA^{(k)}$ in the QR algorithm converges to an upper triangular matrix, with the eigenvalues  on the diagonal (see the second form of Schur decomposition, Corollary~\ref{corollary:schur-second-form}). Fore a more general matrix $\bA$, the sequence converges to an upper \textit{quasi-triangular} matrix. See \citet{quarteroni2010numerical, golub2013matrix} for more details.
\end{remark}

\paragraph{LU algorithm.} 
We  provide a brief overview of the LU algorithm. 
Instead of performing a QR decomposition of the matrix in each iteration, we can also apply the LU decomposition (Theorem~\ref{theorem:lu-factorization-with-permutation}). 
To illustrate this, let's examine the steps outlined in Algorithm~\ref{alg:lu-algorithm-simple2}, where the differences from the simple QR algorithm are highlighted in \textcolor{mylightbluetext}{blue} text. 
Additionally, we  assume that $\bA$ is nonsingular, ensuring that its LU decomposition consists of nonsingular factors.
In such cases, we obtain the following equations:
\begin{equation}\label{equation:simple-qr-find-lu-simple}
\boxed{\begin{aligned}
		\text{(SLU 1)} &\qquad \bA^{(k)} = \bU^{(k)}\bL^{(k)}
		=((\bL^{(k)})^{-1} \bL^{(k)})\bU^{(k)}\bL^{(k)}=(\bL^{(k)})^{-1}\bA^{(k-1)}\bL^{(k)}, \\
		\text{(SLU 2)} &\qquad \bA^{(k)} = 
		(\bL^{(k) } )^{-1} (\bL^{(k-1) })^{-1}\ldots (\bL^{(0) } )^{-1}
		\bA
		\bL^{(0)} \ldots  \bL^{(k-1)}  \bL^{(k)} .
\end{aligned}}
\end{equation}
Based on the equation (SLU 1), the update from $\bA^{(k-1)}$ to $\bA^{(k)}$ involves a similarity transformation (not an orthogonal similarity transformation now). 
Consequently, the accuracy of the method  depends on the condition of each $\bL^{(k) }$, which can sometimes be difficult to control. As a result, the LU algorithm is less commonly employed. 
For more comprehensive information, please refer to \citet{rutishauser1958solution, francis1961qr}, where this algorithm is originally referred to as the \textit{LR algorithm}.

\index{Similarity transformation}

\begin{algorithm}[H] 
\caption{Simple LU Algorithm (Compare to Algorithm~\ref{alg:qr-algorithm-simple2})}
\label{alg:lu-algorithm-simple2}
\begin{algorithmic}[1] 
\Require A real and symmetric $\bA\in \real^{n\times n}$;  
\State $\bA^{(0)} \leftarrow \bA$; 
\State $\bV^{(0)}\leftarrow \bI_n$;  \Comment{initial eigenvector guess}
\State $\textcolor{mylightbluetext}{\bU}^{(0)} \leftarrow \bI_n$; \Comment{compensate the sequence}
\State $\textcolor{mylightbluetext}{\bL}^{(0)} \leftarrow \bI_n$; \Comment{compensate the sequence}
\For{$k=1,2,\ldots$} 
\State $\textcolor{mylightbluetext}{\bL}^{(k)}, \textcolor{mylightbluetext}{\bU}^{(k)}\leftarrow \textcolor{mylightbluetext}{LU}(\bA^{(k-1)})$;
\State $\bA^{(k)} \leftarrow \textcolor{mylightbluetext}{\bU}^{(k)} \textcolor{mylightbluetext}{\bL}^{(k)}$; 
\State $\bV^{(k)} \leftarrow \bV^{(k-1)}\textcolor{mylightbluetext}{\bL}^{(k)} $;
\EndFor
\end{algorithmic} 
\end{algorithm}

\subsection{A Practical QR Algorithm: with Shifts}
Based on the preceding discussion and the shifted inverse power method (Algorithm~\ref{alg:inverse-power-iteration-shifted}), 
we can compute the Rayleigh quotient of the last column of $\bV^{(k-1)}$ during each iteration to obtain an estimate for the \textit{shift value}. 
Alternatively, since the matrix $\bA^{(k-1)}$ approximates the eigenvalues on its  diagonal, $a_{nn}^{(k-1)}$ (the last element on the diagonal of $\bA^{(k-1)}$) can be used  as a shift. The procedure is detailed in  Algorithm~\ref{alg:qr-algorithm-practical}, with distinctions from Algorithm~\ref{alg:qr-algorithm-simple2} highlighted in \textcolor{mylightbluetext}{blue} text.


\begin{algorithm}[H] 
\caption{Practical QR Algorithm (The Final Algorithm! Compare to Algorithm~\ref{alg:qr-algorithm-simple2})}
\label{alg:qr-algorithm-practical}
\begin{algorithmic}[1] 
\Require A real and symmetric $\bA\in \real^{n\times n}$; 
\State $\bA^{(0)} \leftarrow \bA$; 
\State $\bV^{(0)}\leftarrow\bI_n$;\Comment{initial eigenvector guess}
\State $\bR^{(0)} \leftarrow \bI_n$;  \Comment{compensate the sequence}
\State $\bQ^{(0)} \leftarrow \bI_n$; \Comment{compensate the sequence}
\For{$k=1,2,\ldots$} 
\State Pick a shift $\mu^{(k)}$; \Comment{e.g., $\mu^{(k)}=a_{nn}^{(k-1)}$}
\State $\bQ^{(k)}, \bR^{(k)}\leftarrow QR\left(\bA^{(k-1)} - \textcolor{mylightbluetext}{\mu^{(k)}\bI}\right)$; \Comment{QR decomposition}
\State $\bA^{(k)}\leftarrow\bR^{(k)} \bQ^{(k)} + \textcolor{mylightbluetext}{\mu^{(k)}\bI}$;  \Comment{diagonals converge to eigenvalues}
\State $\bV^{(k)} \leftarrow \bV^{(k-1)}\bQ^{(k)} $; \Comment{columns converge to eigenvectors}
\EndFor
\end{algorithmic} 
\end{algorithm}

Similar observations can be made:
\begin{tcolorbox}[title={Practical QR Algorithm Property 1},colback=\mdframecolorTheorem]
\begin{equation}\label{equation:practical-qr-find1}
\left.
\begin{aligned}
\bR^{(k)}&= \bQ^{(k)\top} \left(\bA^{(k-1)} - \mu^{(k)}\bI\right)\\
\bA^{(k)}&=\bR^{(k)} \bQ^{(k)} + \mu^{(k)}\bI  
\end{aligned}
\right\} \implies
\underbrace{\bA^{(k)} = \bQ^{(k)\top } \bA^{(k-1)} \bQ^{(k)}}_{\text{(PQR 1)=(SQR 1)}},
\end{equation}
\end{tcolorbox}
\noindent which is the same \textit{similar transformation} as in the simple QR algorithm, as shown in Equation~\eqref{equation:simple-qr-find}.
Similarly, for Equation~\eqref{equation:simple-qr-find-root}, we can observe:
\begin{tcolorbox}[title={Practical QR Algorithm Property 2},colback=\mdframecolorTheorem]
\begin{equation}\label{equation:practical-qr-find-root}
\begin{aligned}
\bA^{(k)} =\bV^{(k)\top}\bA \bV^{(k)} 
\quad \implies\quad 
\underbrace{\bA =\bV^{(k)}  \bA^{(k)} \bV^{(k)\top}}_{\text{(PQR 2)=(SQR 2)}},
\end{aligned}
\end{equation}
\end{tcolorbox}
\noindent where $\bV^{(k)}=\bQ^{(0)}\bQ^{(1)}\bQ^{(2)}\ldots \bQ^{(k)}$ (same as (SQR 2)).
Once again, $\bA^{(k)} =\bV^{(k)\top}\bA \bV^{(k)} $ is an \textit{orthogonal similarity transformation} because $\bV^{(k)}$ is orthogonal. The condition of $\bA^{(k)}$ remains no worse than that of the original matrix $\bA$. 
Furthermore, the third property of the practical QR algorithm differs slightly  from that of the simple QR algorithm:
\begin{tcolorbox}[title={Practical QR Algorithm Property 3},colback=\mdframecolorTheorem]
\begin{equation}\label{equation:practical-qr-find2}
\begin{aligned}
&\text{(PQR 3)} \neq  \text{(SQR 3)}: \\
&(\bA-\mu^{(k)}\bI)(\bA-\mu^{(k-1)}\bI)\ldots (\bA-\mu^{(1)}\bI)= \underbrace{\bQ^{(0)}\bQ^{(1)}\bQ^{(2)}\ldots \bQ^{(k)}}_{=\bV^{(k)}, \text{orthogonal}}
\underbrace{\bR^{(k)} \bR^{(k-1)} \ldots \bR^{(0)}}_{\triangleq\bU^{(k)}, \text{upper triangular}}.
\end{aligned}
\end{equation}
\end{tcolorbox}
\noindent Again, for clarity, the proof of Equation~\eqref{equation:practical-qr-find2} will be deferred to  Section~\ref{section:proofs-qralgorithms}. Similar to the analysis in the simple QR algorithm (Section~\ref{section:qralg-without-shifts}), assume first $\mu^{(k)}=\mu$ is fixed in Algorithm~\ref{alg:qr-algorithm-practical}, and the eigenvalues are ordered such that $|\lambda_1-\mu| > |\lambda_2-\mu| > \ldots >|\lambda_n-\mu|$, then the $(i,i-1)$-th entry in $\bA^{(k)}$ converges linearly to zero with a constant rate of
$
\mathcalO\big(\left|\frac{\lambda_i-\mu}{\lambda_{i-1}-\mu}\right|^k\big).
$
This implies that if $\mu^{(k)} = a_{nn}^{(k-1)}$, $|\lambda_n - \mu^{(k)}|$ tends to be much smaller than $|\lambda_i - \mu^{(k)}|$ for $i\in \{1,2,\ldots, n-1\}$:
$$
|\lambda_n - \mu^{(k)}| \ll |\lambda_i - \mu^{(k)}|, \gap i\in \{1,2,\ldots, n-1\}.
$$
This  leads to rapid convergence of the last column of  $\bA^{(k)}$ to the corresponding eigenvalue of $\bA$  (i.e., $a_{nn}^{(k-1)}$ converges rapidly to the eigenvalue).

\section{Apply the Practical QR Algorithm to Tridiagonal Matrices}\label{section:qr_totrid}
We observe that the practical QR algorithm  results in the eigengectors and eigenvalues through a sequence of \textit{orthogonal similarity transformations}, as described by the property (PQR 1) in Equation~\eqref{equation:practical-qr-find1}: $\bA^{(k)} = \bQ^{(k)\top } \bA^{(k-1)} \bQ^{(k)}$. When $\bA$ is symmetric, there exists a tridiagonal decomposition $\bA=\bQ^{(0)\top}\bT^{(0)}\bQ^{(0)}$ (Theorem~\ref{theorem:tridiagonal-decom}), which is an orthogonal similarity transformation as well. 
This decomposition can be employed as the first phase  of the QR algorithm because the tridiagonal matrix $\bT^{(0)}$ is closer to a diagonal form, making it  easier to converge to the diagonal eigenvalue matrix. The procedure is outlined in Algorithm~\ref{alg:qr-algorithm-practical-tridiagonal}.
\subsection{Explicit Shifted QR Algorithm}
\begin{algorithm}[h] 
\caption{Practical QR Algorithm (Two Phases, Compare to Algorithm~\ref{alg:qr-algorithm-practical})}
\label{alg:qr-algorithm-practical-tridiagonal}
\begin{algorithmic}[1] 
\Require A real and symmetric $\bA\in \real^{n\times n}$; 
\State \textcolor{mylightbluetext}{$\bA \leftarrow \bQ^{(0)}\bT^{(0)}\bQ^{(0)\top}$}; \Comment{tridiagonal decomposition of $\bA$}
\State $\bV^{(0)} \leftarrow \textcolor{mylightbluetext}{\bQ^{(0)}}$;    \Comment{previously $\bV^{(0)}=\bI_n$}
\State $\bR^{(0)} \leftarrow \bI_n$;
\For{$k=1,2,\ldots$} 
\State Pick a shift $\mu^{(k)}$; \Comment{e.g., $\mu^{(k)}=\textcolor{mylightbluetext}{t}_{nn}^{(k-1)}$}
\State $\bQ^{(k)}, \bR^{(k)}\leftarrow QR\left(\textcolor{mylightbluetext}{\bT}^{(k-1)} - \mu^{(k)}\bI\right)$;
\State $\textcolor{mylightbluetext}{\bT}^{(k)} \leftarrow \bR^{(k)} \bQ^{(k)} + \mu^{(k)}\bI$; 
\State $\bV^{(k)} \leftarrow \bV^{(k-1)}\bQ^{(k)} $;
\EndFor
\end{algorithmic} 
\end{algorithm}

The properties of the practical QR algorithm still hold that 
\begin{equation}\label{equation:tpqr-1-2}
{
\begin{aligned}
\text{(TPQR 1)} \qquad 	\bT^{(k)} &= \bQ^{(k)\top } \bT^{(k-1)} \bQ^{(k)};\\
\text{(TPQR $2^\prime$)} \qquad 	\bT^{(k)} 
&= \bQ^{(k)\top } \bQ^{(k-1)\top }\ldots \bQ^{(1)\top } \bT^{(0)}\bQ^{(1)} \ldots  \bQ^{(0)}  \bQ^{(k)} \bQ^{(k)};\\
&= \underbrace{\bQ^{(k)\top } \bQ^{(k-1)\top }\ldots \bQ^{(1)\top }\bQ^{(0)\top } }_{\bV^{(k)\top}}
\bA
\underbrace{\bQ^{(0) }\bQ^{(1)} \ldots  \bQ^{(0)}  \bQ^{(k)} \bQ^{(k)}}_{\bV^{(k)}}.
\end{aligned}}
\end{equation}
Except now that the $\bT^{(k)}$'s are tridiagonal matrices (see discussion in the next section).
Suppose $\bT^{(k-1)}$ has the following form, where we want to decide the shift value $\mu^{(k)}$:
\begin{equation}\label{equation:implicit-matrix-qr}
\bT^{(k-1)} = 
\footnotesize
\begin{bmatrix}
a_1 & b_1 & & \ldots & 0 \\
b_1 & a_2 & \ddots & & \vdots \\
& \ddots & \ddots & \ddots & \\
\vdots & & \ddots & a_{n-1} & b_{n-1} \\
0 & \ldots & & b_{n-1} & a_n 
\end{bmatrix}.
\end{equation}
As discussed, one reasonable choice for the shift is $\mu^{(k)} = a_n$. 
Additionally, \citet{wilkinson1968global} shows that a more effective choice is to shift by the eigenvalue of 
$$
\bT^{(k-1)}[n-1:n, n-1:n]=
\begin{bmatrix}
a_{n-1} & b_{n-1}\\
b_{n-1} & a_n
\end{bmatrix}.
$$
The  matrix above has two possible eigenvalues, where the one closer to $a_n$ is chosen. It is given by 
$$
\mu^{(k)} = \frac{a_{n-1}+a_n   }{2} + \text{sign}(d) \frac{1}{2}\sqrt{(a_{n-1}-a_n)^2+4b_{n-1}^2},
$$
where $d=a_n-a_{n-1}$. This is known as the \textit{Wilkinson shift}. \citet{wilkinson1968global} also shows the algorithm converges \textit{cubically} to the eigenvectors of $\bA$ with any  of the shift strategies, and the Wilkinson shift is preferred for some heuristic reasons.

\subsection{Implicit Shift QR Algorithm}\label{section:implifit-shift-qr}
In Algorithm~\ref{alg:qr-algorithm-practical-tridiagonal}, there can be instances where the shift $\mu^{(k)}$ is significantly larger than some of the diagonal values $a_i$ in $\bT^{(k-1)}$. 
Thus, it is more reasonable to prioritize  an implicit update on the $\bT^{(k)}$ from $\bT^{(k-1)}$, in which case there is no such a sequence of $\bV^{(k)}$ converging to the eigenvectors of $\bA$; and therefore, it requires an extra ``explicit" computation on the eigenvectors. To see this, we present some observations that are necessary for the \textit{implicit shift QR algorithm} (please review the properties of the tridiagonal decomposition in Section~\ref{section:tridiagonal-decomposition}  before proceeding).
\begin{itemize}
\item \textit{Preservation of Form in Simple QR Algorithm.} We note that the tridiagonal decomposition of a tridiagonal matrix may not be a tridiagonal matrix, i.e., $\bT_+ = \bQ\bT\bQ^\top$ may not be tridiagonal even if $\bT$ is tridiagonal. However, in our case, if $\bT=\bQ\bR$ is the QR decomposition of a symmetric and tridiagonal matrix $\bT\in \real^{n\times n}$, then $\bQ$ has a lower bandwidth of 1 and $\bR$ has an upper bandwidth of 2 (Definition~\ref{defin:matrix-bandwidth}). Consequently, the reverse QR update $\bT_+ = \bR\bQ = \bQ^\top \bQ(\bR\bQ) = \bQ^\top \bT \bQ$ is also symmetric and tridiagonal; see Problem~\ref{prob:reverse_qr_trid}.

\item \textit{Positive Lower Subdiagonals of the Tridiagonal Matrix.} The QR decomposition is not unique (Section~\ref{section:nonunique-qr}). However, when the diagonals of $\bR$ are restricted to be positive, the QR decomposition is unique (Corollary~\ref{corollary:unique-qr}). Then, if $\bT$ has \textit{positive lower subdiagonals} (which implies $\bT$ is \textit{unreduced}, Definition~\ref{definition:tridiagonal-hessenbert}, i.e., nonzero lower subdiagonals, and further the lower subdiagonals are \textcolor{mylightbluetext}{positive}),  $\bT_+=\bQ^\top\bT\bQ$ also has positive lower subdiagonals.~\footnote{This is an important claim often overlooked in many texts.}

\item \textit{Preservation of Form in Practical QR Algorithm.} If $\mu\in \real$, and $\bT-\mu\bI=\bQ\bR$ is the QR decomposition of the shifted $\bT$, then $\bT_+=\bR\bQ +\mu\bI$ is also symmetric and tridiagonal.

\item \textit{Implicit Q Theorem.} We observe that if we restrict the elements in the lower subdiagonal of the tridiagonal matrix $\bT$ to be \textcolor{mylightbluetext}{positive} (if possible), i.e., \textit{unreduced} with positive lower subdiagonals, then the tridiagonal decomposition $\bT_+=\bQ\bT\bQ^\top$ is uniquely determined by $\bT$ and the first column of $\bQ$ (Theorem~\ref{theorem:implicit-q-tridiagonal}). This claim is significant; see the next point.
 
\index{Implicit Q theorem}
\item \textit{Tridiagonal Update.} According to condition (TPQR 1) in Equation~\eqref{equation:tpqr-1-2}, when the QR algorithm is applied to a tridiagonal matrix, the update of the ``eigenvalue matrix" forms a tridiagonal decomposition implicitly: $\bT^{(k)} = \bQ^{(k)\top } \bT^{(k-1)} \bQ^{(k)}\rightarrow  \bT^{(k-1)}=\bQ^{(k) }\bT^{(k)}\bQ^{(k)\top } $, i.e., the tridiagonal decomposition of $\bT^{(k-1)}$ is given by $\bT^{(k-1)}=\bQ^{(k) }\bT^{(k)}\bQ^{(k)\top } $. If we assume $\bT^{(k-1)}$ is \textit{unreduced with positive lower subdiagonals}, then $\bT^{(k)}$ is also \textit{unreduced with positive lower subdiagonals}. 
As per the implicit Q theorem mentioned above, the tridiagonal decomposition is \textbf{uniquely} determined by $\bT^{(k-1)}$ itself and the first column of $\bQ^{(k) }$.

\item \textit{Convergence of Lower Subdiagonals.} If $\bT^{(k-1)}$ has positive lower subdiagonals, $\bT^{(k)}$ will also have positive lower subdiagonals, which in turn results in the positive lower subdiagonals in $\bT^{(k+1)}$. However, at a certain stage, the lower subdiagonals (which are equal to the superdiagonals) will converge to zero. 

\item \textit{Connection to the ``simple" QR algorithm.} We notice that the ``simple" QR algorithm without shifts (Algorithm~\ref{alg:qr-algorithm-simple2})   applied to the tridiagonal matrix also has this tridiagonal update (as that in the practical QR algorithm). However, the difference lies in the first column of $\bQ^{(k) }$. 
\end{itemize}

\index{Bulge}
We then introduce the implicit update for the practical QR algorithm.
\paragraph{Step 1: introducing the bulge.}
To illustrate the implicit shift algorithm, we will walk through a step-by-step example using a $5\times 5$ matrix. 
Suppose at the $(k-1)$-th iteration, the elements in $\bT^{(k-1)}$ are given by Equation~\eqref{equation:implicit-matrix-qr}. We find the $2$ by 2 Givens rotation $\widetildebG_1^\top$, where $c=\cos(\theta)$ and $s=\sin(\theta)$ are computed such that
$$
\underbrace{\begin{bmatrix}
		c & s \\
		-s & c
\end{bmatrix}}_{\triangleq\widetildebG_1^\top}
\begin{bmatrix}
a_1 - \textcolor{mylightbluetext}{\mu^{(k)}}\\
b_1
\end{bmatrix}
= 
\begin{bmatrix}
\boxtimes \\
0
\end{bmatrix}.	
$$
And an $n$ by $n$ Givens rotation can be constructed by 
$$
\bG_{12}^\top = 
\begin{bmatrix}
\widetildebG_1^\top & \\
 & \bI_{n-2}
\end{bmatrix},
$$
where $n$ denotes the matrix size, and the subscript ``$12$" of $\bG_{12}^\top$ indicates the position where the rotation occurs (Definition~\ref{definition:givens-rotation-in-qr}). Up to this point, this is equivalent to what we have done in the first step of the QR decomposition of ($\bT^{(k-1)}-\mu^{(k)}\bI$) via the Givens rotation (Section~\ref{section:qr-givens}). \textbf{This is also what we do in step 6 of Algorithm~\ref{alg:qr-algorithm-practical-tridiagonal}.} 
For a $5\times 5$ example, we observe that $\bG_{12}^\top$ applied to $(\bT^{(k-1)}-\textcolor{mylightbluetext}{\mu^{(k)}\bI})$ will introduce a zero in entry (2,1) and destroy the zero in entry $(1,3)$. The process is shown as follows, where $\boxtimes$ represents a value that is not necessarily zero, and \textbf{boldface} indicates the value has just been changed:
\begin{equation}\label{equation:bulge_tri1}
\begin{aligned}
	\begin{sbmatrix}{(\bT^{(k-1)}-\textcolor{mylightbluetext}{\mu^{(k)}\bI})}
		\boxtimes & \boxtimes & 0 & 0 & 0\\
		\boxtimes & \boxtimes & \boxtimes & 0 & 0\\
		0 & \boxtimes & \boxtimes & \boxtimes & 0\\
		0 & 0 & \boxtimes & \boxtimes & \boxtimes\\
		0 & 0 & 0 & \boxtimes& \boxtimes
	\end{sbmatrix}
	&\stackrel{\bG_{12}^\top\times }{\rightarrow}
	\begin{sbmatrix}{\bG_{12}^\top(\bT^{(k-1)}-\textcolor{mylightbluetext}{\mu^{(k)}\bI})}
		\bm{\boxtimes} & \bm{\boxtimes} & \textcolor{mylightbluetext}{\bm{\boxtimes}} & \bm{0} & \bm{0}\\
		\textcolor{mylightbluetext}{\bm{0}} & \bm{\boxtimes} & \bm{\boxtimes} & \bm{0} & \bm{0}\\
		0 & \boxtimes & \boxtimes & \boxtimes & 0 \\
		0 & 0 & \boxtimes & \boxtimes & \boxtimes\\
		0 & 0 & 0 & \boxtimes & \boxtimes
	\end{sbmatrix}.
\end{aligned}
\end{equation}
That is, the Givens rotation $\bG_{12}^\top$ introduces a zero into the (2,1)-th entry of $(\bT^{(k-1)}-\textcolor{black}{\mu^{(k)}\bI})$.
For the following steps, instead of operating on $(\bT^{(k-1)}-\textcolor{black}{\mu^{(k)}\bI})$, we operate directly on $\bT^{(k-1)}$:
\begin{enumerate}
\item To conform to the (TPQR 1) property outlined in Equation~\eqref{equation:tpqr-1-2}, i.e., a tridiagonal decomposition. 
\item We  will see that the first column of $\bG_{12}$ is identical  to the first column of $\bQ^{(k)}$ in Equation~\eqref{equation:tpqr-1-2}.
\end{enumerate}


Additionally, it's important to note that the operation of $\bG_{12}^\top$ on $\bT^{(k-1)}$ will eliminate the zeros in entries (1,3) and (3,1):
$$
\begin{aligned}
	\begin{sbmatrix}{\bT^{(k-1)}}
		\boxtimes & \boxtimes & 0 & 0 & 0\\
		\boxtimes & \boxtimes & \boxtimes & 0 & 0\\
		0 & \boxtimes & \boxtimes & \boxtimes & 0\\
		0 & 0 & \boxtimes & \boxtimes & \boxtimes\\
		0 & 0 & 0 & \boxtimes& \boxtimes
	\end{sbmatrix}
	&\stackrel{\bG_{12}^\top\times }{\rightarrow}
	\begin{sbmatrix}{\bG_{12}^\top\bT^{(k-1)}}
		\bm{\boxtimes} & \bm{\boxtimes} & \textcolor{mylightbluetext}{\bm{\boxtimes}}  & \bm{0} & \bm{0}\\
		\bm{\boxtimes} & \bm{\boxtimes} & \bm{\boxtimes} & \bm{0} & \bm{0}\\
		0 & \boxtimes & \boxtimes & \boxtimes & 0 \\
		0 & 0 & \boxtimes & \boxtimes & \boxtimes\\
		0 & 0 & 0 & \boxtimes & \boxtimes
	\end{sbmatrix}
&\stackrel{\times\bG_{12}}{\rightarrow}
\mathop{\left[\begin{array}{ccccc}
		\bm{\boxtimes} & \bm{\boxtimes} & \textcolor{mylightbluetext}{\bm{\boxtimes}}  & \bm{0} & \bm{0}\\\cline{1-1}
		\bmboxtimes & \bm{\boxtimes} & \bm{\boxtimes} & \bm{0} & \bm{0}\\
		\colorbmboxtimes & \bm{\boxtimes} & \boxtimes & \boxtimes & 0 \\\cline{1-1}
		\bm{0} & \bm{0} & \boxtimes & \boxtimes & \boxtimes\\
		\bm{0} & \bm{0} & 0 & \boxtimes & \boxtimes
	\end{array}\right]}_{\textstyle\mathstrut{\bG_{12}^\top\bT^{(k-1)}\bG_{12}}},
\end{aligned}
$$
where we find the Givens rotation, $\bG_{12}^\top$, preultiplying  a matrix will modify the first two rows of the underlying matrix; and the transpose of the Givens rotation, $\bG_{12}$, postmultiplying a matrix will modify the first two columns of the underlying matrix. The blue $\textcolor{mylightbluetext}{\bm{\boxtimes}}$'s indicate the zero entries in the tridiagonal matrix that are destroyed by the orthogonal similarity transformation, which  is known as ``\textit{introducing the bulge}."

\index{Similarity transformation}

\paragraph{Step 2: chasing the bulge.}
Now, the problem becomes ``\textit{chasing the bulge}", which involves  transferring the ``bulge" back to zero. 
In the meantime, since the final tridiagonal decomposition is determined by the first column. The second Givens rotation can be constructed by calculating $c=\cos(\theta), s=\sin(\theta)$ such that 
$$
\underbrace{\begin{bmatrix}
		c & s \\
		-s & c
\end{bmatrix}}_{\widetildebG_2^\top}
\underbrace{(\bG_{12}^\top\bT^{(k-1)}\bG_{12})_{2:3,1}}_{\text{the vector in the \fbox{box} of above matrix}}
= 
\begin{bmatrix}
	\boxtimes \\
	0
\end{bmatrix}.	
$$
And the $n\times n$ Givens rotation is constructed by 
$$
\bG_{23}^\top = 
\begin{bmatrix}
	1 & & \\
	& \widetildebG_2^\top & \\
& 	& \bI_{n-3}
\end{bmatrix}.
$$
Following the above example, we have
$$
\begin{aligned}
	\begin{sbmatrix}{\bG_{12}^\top\bT^{(k-1)}\bG_{12}}
		\boxtimes & \boxtimes & \textcolor{mylightbluetext}{\boxtimes}  & 0 & 0\\
		\boxtimes & \boxtimes & \boxtimes & 0 & 0\\
		\textcolor{mylightbluetext}{\boxtimes} & \boxtimes & \boxtimes & \boxtimes & 0 \\
		0 & 0 & \boxtimes & \boxtimes & \boxtimes\\
		0 & 0 & 0 & \boxtimes & \boxtimes
	\end{sbmatrix}
	&\stackrel{\bG_{23}^\top\times }{\rightarrow}
	\begin{sbmatrix}{\bG_{23}^\top(\bG_{12}^\top\bT^{(k-1)}\bG_{12})}
		\boxtimes & \boxtimes & \textcolor{mylightbluetext}{\boxtimes}  & 0 & 0\\
		\bm{\boxtimes} & \bm{\boxtimes} &\bm{\boxtimes}  & \textcolor{mylightbluetext}{\bm{\boxtimes} } & \bm{0}\\
		\textcolor{brown}{\bm{0}} & \bm{\boxtimes} & \bm{\boxtimes} & \bm{\boxtimes} & \bm{0} \\
		0 & 0 & \boxtimes & \boxtimes & \boxtimes\\
		0 & 0 & 0 & \boxtimes & \boxtimes
	\end{sbmatrix}
	&\stackrel{\times\bG_{23} }{\rightarrow}
	\mathop{\left[\begin{array}{ccccc}
		\boxtimes & \bm{\boxtimes} & \textcolor{brown}{\bm{0}}  & 0 & 0\\
		\bm{\boxtimes} & \bm{\boxtimes} &\bm{\boxtimes}  & \textcolor{mylightbluetext}{\bm{\boxtimes} } & \bm{0}\\\cline{2-2}
		\textcolor{brown}{\bm{0}} & \bmboxtimes & \bm{\boxtimes} & \bm{\boxtimes} & \bm{0} \\
		0 & \colorbmboxtimes & \bm{\boxtimes} & \boxtimes & \boxtimes\\\cline{2-2}
		0 & \bm{0} & \bm{0} & \boxtimes & \boxtimes 
\end{array}\right]}_{\textstyle\mathstrut{\bG_{23}^\top(\bG_{12}^\top\bT^{(k-1)}\bG_{12})\bG_{23}}},
\end{aligned}
$$
where we find the Givens rotation, $\bG_{23}^\top$, premultiplying  a matrix will modify its rows $2$ and $3$; and the transpose, $\bG_{23}$, postmultiplying a matrix will modify its columns $2$ and $3$. This process reintroduces zeros for entries (3,1) and (1,3) from step 1; however, it also introduces new ``bulges" for entries (4,2) and (2,4). Thus, we can apply another Givens rotation $\bG_{34}$ to introduce zeros back. The same process can go on, and is known as ``chasing the bulge".

The entire example can be summarized as follows, where the blue $\textcolor{mylightbluetext}{\bm{\boxtimes}}$ indicates the introduced bulge, the brown $\textcolor{brown}{\bm{0}}$ indicates the chased-out bulge, and \textbf{boldface} indicates the value has just been changed:
\begin{mdframed}[hidealllines=\mdframehidelineNote,backgroundcolor=\mdframecolor, frametitle={A Complete Example of Implicit QR Algorithm}]
\begin{equation}\label{equation:tridia-update-implicit}
\begin{aligned}
\begin{sbmatrix}{\bT^{(k-1)}}
		\boxtimes & \boxtimes & 0 & 0 & 0\\
		\boxtimes & \boxtimes & \boxtimes & 0 & 0\\
		0 & \boxtimes & \boxtimes & \boxtimes & 0\\
		0 & 0 & \boxtimes & \boxtimes & \boxtimes\\
		0 & 0 & 0 & \boxtimes& \boxtimes
	\end{sbmatrix}
	&\stackrel{\bG_{12} }{\rightarrow}
	\mathop{\left[\begin{array}{ccccc}
			\bm{\boxtimes} & \bm{\boxtimes} & \textcolor{mylightbluetext}{\bm{\boxtimes}}  & \bm{0} & \bm{0}\\\cline{1-1}
			\bmboxtimes & \bm{\boxtimes} & \bm{\boxtimes} & \bm{0} & \bm{0}\\
			\colorbmboxtimes & \bm{\boxtimes} & \boxtimes & \boxtimes & 0 \\\cline{1-1}
			\bm{0} & \bm{0} & \boxtimes & \boxtimes & \boxtimes\\
			\bm{0} & \bm{0} & 0 & \boxtimes & \boxtimes
		\end{array}\right]}_{\textstyle\mathstrut{\bG_{12}^\top(\cdot)\bG_{12}}}
	\stackrel{\bG_{23} }{\rightarrow}
	\mathop{\left[\begin{array}{ccccc}
			\boxtimes & \bm{\boxtimes} & \textcolor{brown}{\bm{0}}  & 0 & 0\\
			\bm{\boxtimes} & \bm{\boxtimes} &\bm{\boxtimes}  & \textcolor{mylightbluetext}{\bm{\boxtimes} } & \bm{0}\\\cline{2-2}
			\textcolor{brown}{\bm{0}} & \bmboxtimes & \bm{\boxtimes} & \bm{\boxtimes} & \bm{0} \\
			0 & \colorbmboxtimes & \bm{\boxtimes} & \boxtimes & \boxtimes\\\cline{2-2}
			0 & \bm{0} & \bm{0} & \boxtimes & \boxtimes 
		\end{array}\right]}_{\textstyle\mathstrut{\bG_{23}^\top(\cdot)\bG_{23}}}\\
\end{aligned}
\end{equation}
$$
\begin{aligned}
\qquad \qquad \qquad  \gap &\stackrel{\bG_{34} }{\rightarrow}
\mathop{\left[\begin{array}{ccccc}
		\boxtimes & \boxtimes & \bm{0}   & \bm{0} & 0\\
		\boxtimes & \boxtimes &\bm{\boxtimes}  & \textcolor{brown}{\bm{0} }& 0\\
		\bm{0} & \bm{\boxtimes} & \bm{\boxtimes} &\bm{\boxtimes} & \textcolor{mylightbluetext}{\bm{\boxtimes} } \\\cline{3-3}
		\bm{0} & \textcolor{brown}{\bm{0} } & \bmboxtimes & \bm{\boxtimes} & \bm{\boxtimes}\\
		0 & 0 & \colorbmboxtimes & \bm{\boxtimes} & \boxtimes  \\\cline{3-3}
	\end{array}\right]}_{\textstyle\mathstrut{\bG_{34}^\top(\cdot)\bG_{34}}}
\stackrel{\bG_{45} }{\rightarrow}
\mathop{\left[\begin{array}{ccccc}
		\boxtimes & \boxtimes & 0 & \bm{0} &\bm{0}\\
		\boxtimes & \boxtimes &\boxtimes  & \bm{0} & \bm{0}\\
		0 & \boxtimes & \boxtimes & \bm{\boxtimes} & \textcolor{brown}{\bm{0} } \\
		\bm{0} & \bm{0} & \bm{\boxtimes} & \bm{\boxtimes} & \bm{\boxtimes}\\
		\bm{0} & \bm{0} & \textcolor{brown}{\bm{0} } & \bm{\boxtimes} & \bm{\boxtimes}  \\
	\end{array}\right]}_{\textstyle\mathstrut{\bT^{(k)}=\bG_{45}^\top(\cdot)\bG_{45}}}.
\end{aligned}
$$
\end{mdframed}
For a general $n\times n$ matrix,  one can compute rotations $\bG_{12}, \bG_{23}, \ldots, \bG_{n-1,n}$ with the property that if $\bZ \triangleq \bG_{12}\bG_{23} \ldots\bG_{n-1,n}$, then $\bT^{(k-1)} = \bZ \bT^{(k)}\bZ^\top$ represents the tridiagonal decomposition, where the first column of $\bZ$ is given by $\bZ\be_1 = \bG_{12}\be_1 = \bQ^{(k)}\be_1$. We realize that the first column of $\bZ$ and $\bQ^{(k)}$ are identical (both are the first column of orthogonal matrix from the QR decomposition of $(\bT^{(k-1)}-\textcolor{black}{\mu^{(k)}\bI})$), and therefore $\bZ=\bQ^{(k)}$ by implicit Q theorem under the following conditions:
\begin{enumerate}
\item $\bT^{(k-1)}$ is unreduced with positive  subdiagonals;
\item The QR decomposition method employed in the QR algorithm is unique, with the diagonals of the upper triangular matrix having positive diagonal values (Corollary~\ref{corollary:unique-qr}).
\end{enumerate}
\index{Implicit Q theorem}

\paragraph{The complete algorithm.}
For simplicity, we represent the construction of the $2\times 2$ Givens rotation $\widetildebG_i$ such that $\widetildebG_i^\top \bx =\widetildebG_i^\top 
\footnotesize
\begin{bmatrix}
	x_1\\x_2
\end{bmatrix}= \begin{bmatrix}
	\boxtimes \\ 0
\end{bmatrix}$ (i.e., the second element is 0 after applying the Givens rotation) using the following function:
$$
\widetildebG_i^\top =\text{givens}(x_1,x_2).
$$ 
In all iterations, $\widetildebG_i^\top $ will be of size $2\times 2$. And we denote the construction of the $n\times n$ Givens rotation $\bG_{i,i+1}^\top$ by the following function:
$$
\bG_{i,i+1}^\top=
G(\widetildebG_i^\top )=
\begin{bmatrix}
\bI_{i-1} & & \\
& \widetildebG_i^\top& \\
& & \bI_{n-i-1} 
\end{bmatrix}.
$$
For further simplification, we will denote $\bG_{i,i+1}^\top$ by $\bG_i^\top$, which performs a left multiplication on another matrix, implicitly modifying the $i$-th and $(i+1)$-th rows of the matrix.
(And performing a right multiplication with  $\bG_i$  will modify the $i$-th and $(i+1)$-th columns of the matrix). 
The complete procedure is  formulated in Algorithm~\ref{alg:qr-algorithm-practical-tridiagonal-implicit-shift}, where $t_{ij}^{(k-1)}$ is the $(i,j)$-th entry of $\bT^{(k-1)}$.

\begin{algorithm}[h] 
\caption{Practical QR Algorithm with Implicit Shift}
\label{alg:qr-algorithm-practical-tridiagonal-implicit-shift}
\begin{algorithmic}[1] 
\Require A real and symmetric $\bA\in \real^{n\times n}$; 
\State $\bA=\bQ^{(0)}\bT^{(0)}\bQ^{(0)\top}$; \Comment{tridiagonal decomposition of $\bA$ with positive subdiagonals}
\For{$k=1,2,\ldots$} 
\State Pick a shift $\mu^{(k)}$; \Comment{e.g., $\mu^{(k)}=t_{nn}^{(k-1)}$}
\State $x_1 \leftarrow t_{11}-\mu^{(k)}, x_2\leftarrow t_{21}$; \Comment{$t_{ij} = t^{(k-1)}_{ij}$}
\State $\bT^{(k)} \leftarrow \bT^{(k-1)}$; \Comment{initialize $\bT^{(k)} $}
\For{$i=1:n-1$}
\State $\widetildebG_i^\top \leftarrow \text{givens}(x_1,x_2)$;
\State $\bG_{i}^\top \leftarrow G(\widetildebG_i^\top)$; \Comment{It modifies the $i$-th and $(i+1)$-th rows of the matrix}
\State $\bT^{(k)} \leftarrow \bG_{i}^\top\bT^{(k)}\bG_{i}$;
\If{$i<n-1$}
\State $x_1 \leftarrow t_{i+1,i}, x_2\leftarrow t_{i+2,i}$;\Comment{$t_{ij} \leftarrow t^{(k-1)}_{ij}$}
\EndIf
\EndFor 
\State $\bQ^{(k)\top} \leftarrow \bG_{n-1}^\top\ldots\bG_1^\top$; \Comment{this results in $\bT^{(k)} =\bQ^{(k)\top} \bT^{(k-1)}\bQ^{(k)} $}

\EndFor
\end{algorithmic} 
\end{algorithm}
Suppose, for iteration $p$, $\bT^{(p)}$ converges to a diagonal matrix (within  machine precision). 
Then, let us  document the updates for each iteration:
$$
\left.
\begin{aligned}
	\bT^{(p)}&=\bQ^{(p)\top}\bT^{(p-1)} \bQ^{(p)}\\
	\bT^{(p-1)}&=\bQ^{(p-1)\top}\bT^{(p-2)} \bQ^{(p-1)}\\
	\vdots &= \vdots \\
	\bT^{(1)}&=\bQ^{(1)\top}\bT^{(0)} \bQ^{(1)}\\
 \bT^{(0)}&=\bQ^{(0)\top}\bA\bQ^{(0)}
\end{aligned}
\right\}
\implies 
\bA = 
\underbrace{\bQ^{(0)}\ldots \bQ^{(p)}}_{\triangleq\bQ}
\bT^{(p)}
\underbrace{(\bQ^{(0)}\ldots \bQ^{(p)})^\top}_{\triangleq\bQ^\top}
$$
represents the approximated spectral decomposition of the real symmetric $\bA$, where $\bQ$ is orthogonal and contains the eigenvectors of $\bA$, and $\bT^{(p)}$ is diagonal and contains the eigenvalues of $\bA$ (Theorem~\ref{theorem:spectral_theorem}).

\paragraph{Decouple.}
We can ensure that the $\bT^{(k-1)}$ possesses nonnegative  subdiagonals (the specific QR decomposition favored as mentioned above). 
However, it may also occur that $\bT^{(k-1)}$ is \textit{reduced}, i.e., some  subdiagonals are equal to zero. 
In this case, the eigenproblem splits into a pair of smaller problems. For example, when $\bT^{(k-1)}_{k+1,k}=0$, then, a ``practical" QR algorithm can be applied to the submatrices:
$$
\bT^{(k-1)}_{1:k,1:k}   \qquad \text{and} \qquad \bT^{(k-1)}_{k+1:n,k+1:n}.
$$
And the eigenvalues can be obtained by 
$$
\Lambda\left(\bT^{(k-1)} \right) = \Lambda\left(\bT^{(k-1)}_{1:k,1:k}  \right) \cup \Lambda\left(\bT^{(k-1)}_{k+1:n,k+1:n}\right),
$$
where $\Lambda(\cdot)$ represents the spectrum of a matrix (Definition~\ref{definition:spectrum}).


\index{Jacobi's rotation}\index{Off-diagonal norm}
\section{Jacobi's Method}\label{section:jacobi-spectral}
Jacobi's method, introduced in 1846 by \citet{jacobi1846theory},  is one of the oldest techniques for computing  the eigenvalues of a matrix. The idea is to iteratively  diagonalize a small submatrix at each time such that the entire matrix can be diagonalized eventually. 
The mathematical measure used to quantify the reduction is known as the \textit{off-diagonal norm}:
$$
\off(\bA) = \sqrt{\sum_{i=1}^{n} \sum_{j=1, j\neq i}^{n} a_{ij}^2 },
$$
i.e., the Frobenius norm of the off-diagonal entries.
The method aims to iteratively reduce this off-diagonal quantity using Jacobi's rotations:
$$
\footnotesize
\begin{blockarray}{cccccccccccc}
	\begin{block}{c[cccccccccc]c}
		&1 &          &   &  &   &   & &  & & &\\
		&& \ddots  &  &  &  & && & & &\\
		&&      & 1 &  & & &  && & & \\
		&&      &  & c &  &  &  & s & & &k\\
		&&& &   & 1 & & && & &\\
		\bJ_{kl}=&&& &   &   &\ddots &  && & &\\
		&&& &  &   &  & 1&& & &\\
		&&& & -s &  &  & &c& & & l\\
		&&& & &  &  & & &1 & &\\
		&&& & &  &  & & & &\ddots &\\
	\end{block}
	&&      &  & k &  &  &  & l & & &\\
\end{blockarray},
$$
which are the same as the Givens rotations (Definition~\ref{definition:givens-rotation-in-qr}). 
However, the distinction lies in how they are employed. 
In Jacobi's rotation, the choice of the angle $\theta$ for $s = \cos \theta$ and $c = \cos \theta$ is made to ensure that the submatrix of $\bJ^\top \bA\bJ$ becomes diagonal.

\subsection{The 2 by 2 Case}
To understand the operation of Jacobi's rotation, let's consider a 2 by 2 submatrix within a symmetric matrix:
$$
\bA(k,l)\triangleq
\begin{bmatrix}
	a_{kk} & a_{kl} \\
	a_{lk} & a_{ll}
\end{bmatrix}
=
\begin{bmatrix}
	a_{kk} & a_{kl} \\
	a_{kl} & a_{ll}
\end{bmatrix}.
$$
Then, the angle $\theta$ can be computed such that 
$$
\begin{aligned}
\bJ^\top \bA(k,l)\bJ
&=
\begin{bmatrix}
	c & -s \\
	s & c
\end{bmatrix}
\begin{bmatrix}
	a_{kk} & a_{kl} \\
	a_{kl} & a_{ll}
\end{bmatrix}
\begin{bmatrix}
	c & s \\
	-s & c
\end{bmatrix}\\
&=
\begin{bmatrix}
	(c^2-cs) a_{kk}+s^2 a_{ll} -cs \cdot a_{kl}  & (c^2-s^2)a_{kl}+cs (a_{kk}-a_{ll}) \\
	(c^2-s^2)a_{kl} +cs (a_{kk}-a_{ll}) & s^2 a_{kk} +c^2 a_{ll}+2cs \cdot a_{kl}
\end{bmatrix}
=
\begin{bmatrix}
	\neq 0 & 0 \\
	0 & \neq 0
\end{bmatrix}.
\end{aligned}
$$
If $a_{kl}=0$, we simply set $c=1$ and $s=0$, leaving the submatrix $\bA(k,l)$ unchanged. 
Otherwise, it is evident that both $c\neq 0$ and $s\neq 0$. Dividing the off-diagonal element $\{(c^2-s^2)a_{kl} +cs (a_{kk}-a_{ll}) \}$ by $c^2$ and $a_{kl}$, it follows that 
$
(1- \frac{s^2}{c^2} + \frac{s}{c} \frac{a_{kk}-a_{ll}}{a_{kl}}) = 0.
$
Let 
$$
\tan \theta \triangleq t \triangleq \frac{\sin \theta}{\cos \theta} = \frac{s}{c}
\qquad 
\text{and}
\qquad 
\tau \triangleq \frac{a_{ll}-a_{kk}}{2a_{kl}},
$$
it suffices to solve the equation
$
t^2+2\tau t -1=0.
$
This gives
$
t = -\tau \pm \sqrt{\tau^2+1}.
$
We notice that
$$
\norm{\bA - \bJ_{kl}^\top \bA\bJ_{kl}}_F^2 = 4(1-c) \sum_{i\neq k,l}(a_{ik}^2+a_{il}^2) + 2a_{kl}^2 /c^2
\qquad\text{and}\qquad
c = \frac{1}{\sqrt{1+t^2}}.
$$
Therefore, the larger the $c$, the smaller the $\norm{\bA - \bJ_{kl}^\top \bA\bJ_{kl}}_F^2$. This implies that we should choose a $t$ with a smaller magnitude:
$
t_{min}=
\left\{
\begin{aligned}
&-\tau + \sqrt{\tau^2+1}, \qquad \text{if $\tau\geq 0$};\\
&-\tau - \sqrt{\tau^2+1}, \qquad \text{if $\tau< 0$}.\\
\end{aligned}
\right.
$
Based on the preceding discussion, we can define the \textit{ComputeJacobiRotation} function, which computes the Jacobi rotation from the submatrix of $\bA$.
\begin{algorithm}[H] 
	\caption{Compute Jacobi's Rotation Given the Submatrix}
	\label{alg:jacobi-submatrix}
	\begin{algorithmic}[1] 
		\Require Matrix $\bA\in \real^{n\times n}$ is real and symmetric, ($k,l$) such that $1\leq k< l\leq n$; 
		\Function{computeJacobiRotation}{$\bA$, $k,l$}
		\If{$A_{kl}\neq 0$}
		\State $\tau \leftarrow (a_{ll} - a_{kk})/(2a_{kl})$
		\If{$\tau \geq 0$}
			\State $t\leftarrow -\tau + \sqrt{\tau^2+1}$;
		\Else
			\State $t\leftarrow -\tau - \sqrt{\tau^2+1}$;
		\EndIf
		\State $c\leftarrow 1/\sqrt{1+t^2}, s\leftarrow tc$;
		\Else
		\State $c\leftarrow 1,s\leftarrow 0$;
		\EndIf
		\State Output $c,s$;
		\EndFunction
	\end{algorithmic} 
\end{algorithm}

\subsection{The Complete Jacobi's Method}
The term ``complete" is derived from the concept of  complete pivoting, which involves searching the largest element in the matrix to pivot (Section~\ref{section:complete-pivoting}). While in the complete Jacobi's method, at each iteration, we need to decide which submatrix $\bA(k,l)$ to diagonalize. 
Specifically, in the \textit{complete Jacobi's method}, we select the indices $(k,l)$ such that the $(k,l)$-th entry,  $a_{kl}^2$, is maximal, with the goal of ``maximizing"  the reduction of the off-diagonal quantity. The complete search for the largest magnitude defines  Algorithm~\ref{alg:jacobi-complete-alg}.
\begin{algorithm}[H] 
	\caption{Complete Jacobi's Method}
	\label{alg:jacobi-complete-alg}
	\begin{algorithmic}[1] 
		\Require $\bA\in \real^{n\times n}$ is real and symmetric, a positive tolerance $tol$ such that $\delta = tol \cdot \norm{\bA}_F$;
		\State $\bQ \leftarrow  \bI_n$;
		\State $i\leftarrow 0$; \Comment{Count the number of iterations}
		\State $\bLambda \leftarrow \bA$;
		\While{$\off(\bLambda)> \delta$}
		\State Choose $(k,l)$ so that $a_{kl}^2 \leftarrow \arg\max \bLambda_{ij}^2,\, \forall 1\leq i,j \leq n$;
		\State Compute $c,s$ from \textit{computeJacobiRotation}($\bLambda, k,l$);
		\State Decide the $n\times n$ Jacobi's rotation $\bJ_{kl}$ by $c,s$;
		\State $\bLambda  \leftarrow \bJ_{kl}^\top \bLambda \bJ_{kl}$;
		\State $\bQ \leftarrow \bQ\bJ_{kl}$;
		\State Compute $\off(\bLambda^{(i)})^2 $; \Comment{which, we will show, converges linearly to 0.}
		\State $i \leftarrow i+1$;
		\EndWhile
		\State Output the sequence $\off(\bLambda^{(i)})^2 $,  approximated diagonal $\bLambda$, and orthogonal $\bQ$;
	\end{algorithmic} 
\end{algorithm}
The algorithm computes the spectral decomposition by producing $\bA = \bQ\bLambda\bQ^\top$.
It is not hard to see that, at each iteration, it requires $\mathcalO(n^2)$ flops to search for the largest magnitude entry $(k,l)$ and $\mathcalO(n)$ flops to perform the update of the iteration. We should note that although a symmetric pair of zeros is introduced into the matrix at each iteration,
the previously introduced zeros may be destroyed. However, what matters is that the off-diagonal quantity reduces steadily.

Suppose we select the  index $(k,l)$ and the Jacobi rotation $\bJ_{kl}$ is constructed such that 
$$
\bLambda_+ = \bJ_{kl}^\top\bLambda \bJ_{kl}.
$$
Since the Jacobi rotation is orthogonal, it follows that 
$$
\norm{\bLambda_+}_F^2 = \norm{\bLambda}_F^2
\qquad\text{and}\qquad
\off(\bLambda_+)^2 = \off(\bLambda)^2 -2a_{kl}^2.
$$
As $a_{kl}$ has the largest magnitude in $\bLambda$, we also have 
$$
\off(\bLambda)^2 \leq n(n-1)a_{kl}^2.
$$
Therefore, 
$$
\off(\bLambda_+)^2= \off(\bLambda)^2 -2a_{kl}^2 \leq \off(\bLambda)^2 - \frac{2}{n(n-1)} \off(\bLambda)^2 = \left(1-\frac{2}{n(n-1)}\right)\off(\bLambda)^2.
$$
By induction, this implies that  Jacobi's method reduces the off-diagonal quantity after $k$ iterations by 
$$
\off(\bLambda^{(i)})^2 \leq \left(1-\frac{2}{n(n-1)}\right)^2 \off(\bLambda^{(0)}).
$$
and 
$$
\mathop{\lim}_{k\rightarrow \infty}  \frac{|\off(\bLambda^{(i)})^2 - 0|}{|\off(\bLambda^{(i-1)})^2 -0|} = 1-\frac{2}{n(n-1)} \in (0,1),
$$
such that the sequence $\off(\bLambda^{(i)})^2 $ converges linearly to 0 (Definition~\ref{definition:linear-convergence}). However, in practice, \citet{henrici1958speed, schonhage1964quadratic, van1966quadratic} show that a quadratic convergence (Definition~\ref{definition:quadratic-convergence}) can be obtained such that 
$$
\mathop{\lim}_{k\rightarrow \infty}  \frac{|\off(\bLambda^{(i+\frac{n(n-1)}{2})})^2|}{|\off(\bLambda^{(i)})^2|^2} = c,
$$
where $c$ is a constant. We will not provide further details here.

\subsection{The Cyclic-by-Row Jacobi's Method}

In the previous section, we discussed that the majority of the computational cost associated with the complete Jacobi's method arises from the search for the element with the largest magnitude, a task that consumes $\mathcalO(n^2)$ flops.
To mitigate this computational burden, it is advisable to adopt a row-by-row approach for updates, iteratively addressing the $n(n-1)$ upper triangular entries.
This approach is known as the \textit{cyclic-by-row} algorithm, and the procedure is shown in Algorithm~\ref{alg:jacobi-cyclicrow-alg}. 

Similarly, \citet{wilkinson1962note, van1966quadratic} show that the cyclic-by-row algorithm converges \textit{quadratically}. 
Moreover, due to its avoidance of element comparisons within the matrix, this approach effectively reduces computational complexity.

\begin{algorithm}[H] 
	\caption{Cyclic-by-Row Jacobi's Method}
	\label{alg:jacobi-cyclicrow-alg}
	\begin{algorithmic}[1] 
		\Require $\bA\in \real^{n\times n}$ is real and symmetric, a positive tolerance $tol$ such that $\delta = tol \cdot \normf{\bA}$;
		\State $\bQ \leftarrow \bI_n$;
		\State $i\leftarrow 0$;  \Comment{Count the number of iterations}
		\State $\bLambda\leftarrow \bA$;
		\While{$\off(\bLambda)> \delta$}
		\For{$k=1:n-1$}
		\For{$l=k+1:n$}
		\State Compute $c,s$ from \textit{computeJacobiRotation}($\bLambda, k,l$);
		\State Decide the $n\times n$ Jacobi's rotation $\bJ_{kl}$ by $c,s$;
		\State $\bLambda  \leftarrow \bJ_{kl}^\top \bLambda \bJ_{kl}$;
		\State $\bQ \leftarrow \bQ\bJ_{kl}$;
		\State Compute $\off(\bLambda^{(i)})^2 $; 
		\State $i\leftarrow i+1$;
		\EndFor
		\EndFor
		\EndWhile
		\State Output the sequence $\off(\bLambda^{(i)})^2 $,  approximated diagonal $\bLambda$, and orthogonal $\bQ$;
	\end{algorithmic} 
\end{algorithm}

\subsection{Other Issues}
In practice, when computing on a $p$-processor computer, it is reasonable to implement  Jacobi's algorithm in a block-wise manner, thereby facilitating a parallelizable approach. 
We refer the issue to \citet{bischof1986two, shroff1989convergence, golub2013matrix}, and its counterpart in SVD computation has been explored by \citet{van1985block}.

\section{Computing the SVD}
\subsection{Implicit Shift QR Algorithm}
We  have previously shown that any symmetric matrix can be reduced to a tridiagonal form through a sequence of Householder reflectors that are applied  alternately from the left and right. 
This process is a special case of the Hessenberg decomposition (Section~\ref{section:compute-tridiagonal}) and  can reduce the cost of the QR algorithm when  computing the spectral decomposition of a matrix. This two-phase computation is not unique and has its counterpart. 
Due to the work of Golub, Kahan, and others on the bidiagonalization in the 1960s (Theorem~\ref{theorem:Golub-Kahan-Bidiagonalization-decom}), an analogous two-phase approach has become the standard for computing the SVD. 
In  SVD computation, the matrix is first reduced to a bidiagonal form \footnote{In this section, when we refer to a bidiagonal matrix, we implicitly mean an upper bidiagonal matrix.} and then the bidiagonal matrix is diagonalized:
$$
\begin{bmatrix}
	\boxtimes & \boxtimes & \boxtimes & \boxtimes & \boxtimes\\
	\boxtimes & \boxtimes & \boxtimes & \boxtimes & \boxtimes\\
	\boxtimes & \boxtimes & \boxtimes & \boxtimes & \boxtimes\\
	\boxtimes & \boxtimes & \boxtimes & \boxtimes & \boxtimes\\
	\boxtimes & \boxtimes & \boxtimes & \boxtimes & \boxtimes\\
	\boxtimes & \boxtimes & \boxtimes & \boxtimes & \boxtimes\\
	\boxtimes & \boxtimes & \boxtimes & \boxtimes & \boxtimes\\
\end{bmatrix}
\mathop{\longrightarrow}^{\text{phase 1} }_{\text{bidiagonalize} }
\begin{bmatrix}
	\boxtimes & \boxtimes & 0 & 0 & 0\\
	0 & \boxtimes & \boxtimes & 0 & 0\\
	0 & 0 & \boxtimes & \boxtimes & 0\\
	0 & 0 & 0 & \boxtimes & \boxtimes\\
	0 & 0 & 0 & 0 & \boxtimes\\
	0 & 0 & 0 & 0 & 0\\
	0 & 0 & 0 & 0 & 0
\end{bmatrix}
\mathop{\longrightarrow}^{\text{phase 2} }_{\text{diagonalize} }
\begin{bmatrix}
	\boxtimes & 0 & 0 & 0 & 0\\
	0 & \boxtimes & 0 & 0 & 0\\
	0 & 0 & \boxtimes & 0 & 0\\
	0 & 0 & 0 & \boxtimes & 0\\
	0 & 0 & 0 & 0 & \boxtimes\\
	0 & 0 & 0 & 0 & 0\\
	0 & 0 & 0 & 0 & 0
\end{bmatrix}.
$$
Bidiagonalization differs from  tridiagonalization in that it does not require the two orthogonal matrices on the left and on the right to be the transpose of the same matrix, and their shapes can be different. This is precisely what is needed for the SVD. For the second phase of the SVD computation, we will proceed directly to the solution using the implicit shift QR algorithm. 
The development to this final procedure is similar to what we have developed for the spectral decomposition.

\paragraph{Phase 2 of SVD.}
We  note the fact that $\bT=\bB^\top\bB$ is tridiagonal if $\bB$ is bidiagonal (Lemma~\ref{lemma:construct-triangular-from-bidia}). Following further the tridiagonal update in the QR algorithm (Section~\ref{section:implifit-shift-qr}), suppose in the $k$-th iteration, we have the bidiagonal matrix $\bB^{(k-1)}$ and its tridiagonal companion \textcolor{mylightbluetext}{$\bT^{(k-1)} = \bB^{(k-1)\top}\bB^{(k-1)}$}, for $k=1,2,\ldots$.
The SVD of the bidiagonal matrix $\bB^{(k-1)}$ can be obtained by finding the spectral decomposition of the tridiagonal matrix $\bT^{(k-1)}$.
Therefore, we consider the tridiagonal update (i.e., introducing and chasing the bulge) applied to $\bT^{(k-1)} = \bB^{(k-1)\top}\bB^{(k-1)}$  (Section~\ref{section:implifit-shift-qr}).
The tridiagonal update is followed by a set of Givens rotations on the left and right iteratively, as shown in Equation~\eqref{equation:tridia-update-implicit}  with a $5\times 5$ example, where the blue $\textcolor{mylightbluetext}{\bm{\boxtimes}}$ indicates the bulge introduced, and the \fbox{boxed} vector indicates how the Givens matrix is constructed:
$$
\begin{aligned}
	\bT^{(k-1)}=	\begin{sbmatrix}{\bT^{(k-1)}}
		\boxtimes & \boxtimes & 0 & 0 & 0\\
		\boxtimes & \boxtimes & \boxtimes & 0 & 0\\
		0 & \boxtimes & \boxtimes & \boxtimes & 0\\
		0 & 0 & \boxtimes & \boxtimes & \boxtimes\\
		0 & 0 & 0 & \boxtimes& \boxtimes
	\end{sbmatrix}
	&\stackrel{\bG_{12} }{\rightarrow}
	\mathop{\left[\begin{array}{ccccc}
			\bm{\boxtimes} & \bm{\boxtimes} & \textcolor{mylightbluetext}{\bm{\boxtimes}}  & \bm{0} & \bm{0}\\\cline{1-1}
			\bmboxtimes & \bm{\boxtimes} & \bm{\boxtimes} & \bm{0} & \bm{0}\\
			\colorbmboxtimes & \bm{\boxtimes} & \boxtimes & \boxtimes & 0 \\\cline{1-1}
			\bm{0} & \bm{0} & \boxtimes & \boxtimes & \boxtimes\\
			\bm{0} & \bm{0} & 0 & \boxtimes & \boxtimes
		\end{array}\right]}_{\textstyle\mathstrut{\bG_{12}^\top(\cdot)\bG_{12}}}
	\stackrel{\bG_{23} }{\rightarrow}
	\mathop{\left[\begin{array}{ccccc}
			\boxtimes & \bm{\boxtimes} & \textcolor{brown}{\bm{0}}  & 0 & 0\\
			\bm{\boxtimes} & \bm{\boxtimes} &\bm{\boxtimes}  & \textcolor{mylightbluetext}{\bm{\boxtimes} } & \bm{0}\\\cline{2-2}
			\textcolor{brown}{\bm{0}} & \bmboxtimes & \bm{\boxtimes} & \bm{\boxtimes} & \bm{0} \\
			0 & \colorbmboxtimes & \bm{\boxtimes} & \boxtimes & \boxtimes\\\cline{2-2}
			0 & \bm{0} & \bm{0} & \boxtimes & \boxtimes 
		\end{array}\right]}_{\textstyle\mathstrut{\bG_{23}^\top(\cdot)\bG_{23}}}\\
	&\stackrel{\bG_{34} }{\rightarrow}
	\mathop{\left[\begin{array}{ccccc}
			\boxtimes & \boxtimes & \bm{0}   & \bm{0} & 0\\
			\boxtimes & \boxtimes &\bm{\boxtimes}  & \textcolor{brown}{\bm{0} }& 0\\
			\bm{0} & \bm{\boxtimes} & \bm{\boxtimes} &\bm{\boxtimes} & \textcolor{mylightbluetext}{\bm{\boxtimes} } \\\cline{3-3}
			\bm{0} & \textcolor{brown}{\bm{0} } & \bmboxtimes & \bm{\boxtimes} & \bm{\boxtimes}\\
			0 & 0 & \colorbmboxtimes & \bm{\boxtimes} & \boxtimes  \\\cline{3-3}
		\end{array}\right]}_{\textstyle\mathstrut{\bG_{34}^\top(\cdot)\bG_{34}}}
	\stackrel{\bG_{45} }{\rightarrow}
	\mathop{\left[\begin{array}{ccccc}
			\boxtimes & \boxtimes & 0 & \bm{0} &\bm{0}\\
			\boxtimes & \boxtimes &\boxtimes  & \bm{0} & \bm{0}\\
			0 & \boxtimes & \boxtimes & \bm{\boxtimes} & \textcolor{brown}{\bm{0} } \\
			\bm{0} & \bm{0} & \bm{\boxtimes} & \bm{\boxtimes} & \bm{\boxtimes}\\
			\bm{0} & \bm{0} & \textcolor{brown}{\bm{0} } & \bm{\boxtimes} & \bm{\boxtimes}  \\
		\end{array}\right]}_{\textstyle\mathstrut{\bG_{45}^\top(\cdot)\bG_{45}}}=\bT^{(k)}.
\end{aligned}
$$
That is,
$$
\begin{aligned}
\bT^{(k)} &= 
(\bG_{45}^\top\bG_{34}^\top\bG_{23}^\top\bG_{12}^\top)
\textcolor{mylightbluetext}{\bT^{(k-1)}}
(\bG_{12}\bG_{23}\bG_{34}\bG_{45})\\
&=
(\bG_{45}^\top\bG_{34}^\top\bG_{23}^\top\bG_{12}^\top)
\textcolor{mylightbluetext}{(\bB^{(k-1)\top}\bB^{(k-1)})}
(\bG_{12}\bG_{23}\bG_{34}\bG_{45}).
\end{aligned}
$$
If we can show that $\bT^{(k)}$ can also be decomposed into a product of bidiagonal matrices, then the ``tridiagonal update" can be replaced with  a ``bidiagonal update." To achieve this goal, two issues must be addressed:
\begin{itemize}
\item Decompose $\bT^{(k)}$ into the bidiagonal form if $\bT^{(k-1)}$ has the bidiagonal form  $\bT^{(k-1)} = \bB^{(k-1)\top}\bB^{(k-1)}$;
\item Find the corresponding Givens rotations using the bidiagonal matrix $\bB^{(k-1)}$ instead of the tridiagonal one $\bT^{(k-1)}$. 
\end{itemize}
What we aim to demonstrate is not difficult to ascertain. This can be either 
\begin{equation}\label{equation:svd_chase_choce}
\begin{aligned}
	&\textbf{(Choice 1)}: \qquad \bB^{(k)\top} = (\bG_{45}^\top\bG_{34}^\top\bG_{23}^\top\bG_{12}^\top)\bB^{(k-1)\top} \qquad \text{or}\\
	&\textbf{(Choice 2)}: \qquad \bB^{(k)\top} = (\bG_{45}^\top\bG_{34}^\top\bG_{23}^\top\bG_{12}^\top)\bB^{(k-1)\top}
	\textcolor{mylightbluetext}{\bV_{12}\bV_{23}\bV_{34}\bV_{45}}.
\end{aligned}
\end{equation}
That is, out goal is to ensure that one of the aforementioned choices remains in a lower bidiagonal form if  $ \bB^{(k-1)\top}$ is lower bidiagonal. The $\bV_{i,i+1}$ matrices are orthogonal such that they will cancel out in the product of $\bB^{(k)\top} \bB^{(k)} $. With hindsight, we will see that the second form of $\bB^{(k)\top} $ will be employed (to chase the bulge).
For clarity, let's assume that $\bT^{(k-1)}$ and $\bB^{(k-1)}$ take the following forms (taking a $5\times 5$ matrix as an example):
$$
\bT^{(k-1)} = 
\begin{bmatrix}
	a_1 & b_1 & & \ldots & 0 \\
	b_1 & a_2 & \ddots & & \vdots \\
	& \ddots & \ddots & \ddots & \\
	\vdots & & \ddots & \ddots & b_{4} \\
	0 & \ldots & & b_{4} & a_5
\end{bmatrix}=
\begin{sbmatrix}{=\bB^{(k-1)\top}}
	c_1 & 0 & & \ldots & 0 \\
	d_1 & c_2 & \ddots & & \vdots \\
	& \ddots & \ddots & \ddots & \\
	\vdots & & \ddots & \ddots & 0 \\
	0 & \ldots & & d_{4} & c_5
\end{sbmatrix}
\begin{sbmatrix}{=\bB^{(k-1)} }
	c_1 & d_1 & & \ldots & 0 \\
	0 & c_2 & \ddots & & \vdots \\
	& \ddots & \ddots & \ddots & \\
	\vdots & & \ddots & \ddots & d_{4} \\
	0 & \ldots & & 0 & c_5
\end{sbmatrix}.
$$

\paragraph{Step 1: introducing and chasing the bulge in the bidiagonal update.}
Following the Givens rotations constructed in Section~\ref{section:implifit-shift-qr}, it follows that 
we can obtain a Givens rotation from the bidiagonal matrix as follows:
$$
\bG_{12}^\top = 
\begin{bmatrix}
	\widetildebG_1^\top & \\
	& \bI_{n-2}
\end{bmatrix},
\qquad 
\text{with}
\gap
\underbrace{\begin{bmatrix}
		c & s \\
		-s & c
\end{bmatrix}}_{\widetildebG_1^\top }
\begin{bmatrix}
	c_1^2 - \mu^{(k)}\\
	c_1d_1
\end{bmatrix}
= 
\begin{bmatrix}
	\boxtimes \\
	0
\end{bmatrix},
$$
where $\mu^{(k)}$ is the shift in the $k$-th iteration, $n=5$ and $a_1 = c_1^2, b_1=c_1d_1$ for our example. This operation will \textit{introduce a bulge} in $\bB^{(k-1)\top}$ (same as the bulge in the tridiagonal form of $\bB^{(k-1)\top}\bB^{(k-1)}$ in \eqref{equation:bulge_tri1}):
$$
\begin{aligned}
	\begin{sbmatrix}{\bB^{(k-1)\top}}
		\boxtimes & 0 & 0 & 0 & 0\\
		\boxtimes & \boxtimes & 0 & 0 & 0\\
		0 & \boxtimes & \boxtimes & 0 & 0\\
		0 & 0 & \boxtimes & \boxtimes & 0\\
		0 & 0 & 0 & \boxtimes& \boxtimes
	\end{sbmatrix}
	&\stackrel{\bG_{12}^\top\times }{\rightarrow}
	\begin{sbmatrix}{\bG_{12}^\top\bB^{(k-1)\top}}
		\bm{\boxtimes} & \textcolor{mylightbluetext}{\bm{\boxtimes}}  &\bm{0} & \bm{0} & \bm{0}\\
		\bm{\boxtimes} & \bm{\boxtimes} & \bm{0} & \bm{0} & \bm{0}\\
		0 & \boxtimes & \boxtimes & 0 & 0 \\
		0 & 0 & \boxtimes & \boxtimes & 0\\
		0 & 0 & 0 & \boxtimes & \boxtimes
	\end{sbmatrix}.
\end{aligned}
$$
We observe that applying a left rotation $\bG_{23}^\top$ will not help chase out the bulge \footnote{Applying $\bG_{23}^\top$ to the left of $\bG_{12}^\top\bB^{(k-1)\top}$ will modify its rows 2 and 3, which will not introduce zeros back in the first row.}. 
To address this, we need to construct a right Givens rotation in such a way that (this is the choice 2 in \eqref{equation:svd_chase_choce}):
$$
\bV_{12}^\top = 
\begin{bmatrix}
	\widetildebV_1^\top & \\
	& \bI_{n-2}
\end{bmatrix},
\qquad 
\text{where}
\gap
\underbrace{\begin{bmatrix}
		c & s \\
		-s & c
\end{bmatrix}}_{\triangleq\widetildebV_1^\top}
\begin{bmatrix}
	\Big(\bG_{12}^\top\bB^{(k-1)\top}\Big)_{11}\\
\Big(\bG_{12}^\top\bB^{(k-1)\top}\Big)_{12}
\end{bmatrix}
= 
\begin{bmatrix}
	\boxtimes \\
	0
\end{bmatrix}.	
$$
This results in 
$$
\begin{aligned}
	\begin{sbmatrix}{\bB^{(k-1)\top}}
		\boxtimes & 0 & 0 & 0 & 0\\
		\boxtimes & \boxtimes & 0 & 0 & 0\\
		0 & \boxtimes & \boxtimes & 0 & 0\\
		0 & 0 & \boxtimes & \boxtimes & 0\\
		0 & 0 & 0 & \boxtimes& \boxtimes
	\end{sbmatrix}
	&\stackrel{\bG_{12}^\top\times }{\rightarrow}
	\begin{sbmatrix}{\bG_{12}^\top\bB^{(k-1)\top}}
		\bm{\boxtimes} & \textcolor{mylightbluetext}{\bm{\boxtimes}}  &\bm{0} & \bm{0} & \bm{0}\\
		\bm{\boxtimes} & \bm{\boxtimes} & \bm{0} & \bm{0} & \bm{0}\\
		0 & \boxtimes & \boxtimes & 0 & 0 \\
		0 & 0 & \boxtimes & \boxtimes & 0\\
		0 & 0 & 0 & \boxtimes & \boxtimes
	\end{sbmatrix}
&\stackrel{\times\bV_{12} }{\rightarrow}
\begin{sbmatrix}{\bG_{12}^\top\bB^{(k-1)\top} \bV_{12}}
	\bm{\boxtimes} & \textcolor{brown}{\bm{0}}  &\bm{0} & \bm{0} & \bm{0}\\
	\bm{\boxtimes} & \bm{\boxtimes} & \bm{0} & \bm{0} & \bm{0}\\
	\textcolor{mylightbluetext}{\bm{\boxtimes}} & \bm{\boxtimes} & \boxtimes & 0 & 0 \\
	\bm{0} & \bm{0} & \boxtimes & \boxtimes & 0\\
	\bm{0} & \bm{0}& 0 & \boxtimes & \boxtimes
\end{sbmatrix}.
\end{aligned}
$$
As long as the $\bG_{12}^\top$ is constructed in the same way as  in the ``tridiagonal update," then 
$\bG_{12}^\top \bT^{{(k-1)}}\bG_{12}$
is equal to 
$\bG_{12}^\top \bB^{{(k-1)}\top}{\bV_{12}\bV_{12}^\top} \bB^{(k-1)}\bG_{12}$ since $\bV_{12}\bV_{12}^\top=\bI$.

\paragraph{Step 2.} Step 2 differs from what we did in the ``tridiagonal update." 
In the ``tridiagonal update," we construct a Givens matrix $\bG_{23}$ to chase the bulge in $\bT^{(k-1)}$. 
By applying the implicit Q theorem again (Theorem~\ref{theorem:implicit-q-tridiagonal}), we can determine the Givens rotation directly using $\bB^{(k-1)}$  (\textit{without explicitly constructing $\bT^{(k-1)}=\bB^{(k-1)\top}\bB^{(k-1)}$}). We denote this new Givens rotation as  $\bU_{23}$ to distinguish it from $\bG_{23}$:
$$
\bU_{23}^\top = 
\begin{bmatrix}
	1 & & \\
	& \widetildebU_2^\top & \\
	& & \bI_{n-3}
\end{bmatrix},
\qquad 
\text{where}
\gap
\underbrace{\begin{bmatrix}
		c & s \\
		-s & c
\end{bmatrix}}_{\triangleq\widetildebU_2^\top }
\begin{bmatrix}
	\Big(\bG_{12}^\top\bB^{(k-1)\top} \bV_{12}\Big)_{21}\\
	\Big(\bG_{12}^\top\bB^{(k-1)\top} \bV_{12}\Big)_{31}
\end{bmatrix}
= 
\begin{bmatrix}
	\boxtimes \\
	0
\end{bmatrix}.	
$$

The process can go on, and the complete example for chasing the bulge can be shown as follows, where the blue $\textcolor{mylightbluetext}{\bm{\boxtimes}}$ indicates the bulge introduced, the brown $\textcolor{brown}{\bm{0}}$ indicates the bulge has been chased out, and \textbf{boldface} indicates the value has just been changed:
\begin{mdframed}[hidealllines=\mdframehidelineNote,backgroundcolor=\mdframecolor, frametitle={A Complete Example of Implicit QR Algorithm for SVD}]
$$
\begin{aligned}
	\begin{sbmatrix}{\bB^{(k-1)\top}}
		\boxtimes & 0 & 0 & 0 & 0\\
		\boxtimes & \boxtimes & 0 & 0 & 0\\
		0 & \boxtimes & \boxtimes & 0 & 0\\
		0 & 0 & \boxtimes & \boxtimes & 0\\
		0 & 0 & 0 & \boxtimes& \boxtimes
	\end{sbmatrix}
	&\stackrel{\bG_{12}^\top\times }{\rightarrow}
	\begin{sbmatrix}{\bG_{12}^\top\bB^{(k-1)\top}}
		\bm{\boxtimes} & \textcolor{mylightbluetext}{\bm{\boxtimes}}  &\bm{0} & \bm{0} & \bm{0}\\
		\bm{\boxtimes} & \bm{\boxtimes} & \bm{0} & \bm{0} & \bm{0}\\
		0 & \boxtimes & \boxtimes & 0 & 0 \\
		0 & 0 & \boxtimes & \boxtimes & 0\\
		0 & 0 & 0 & \boxtimes & \boxtimes
	\end{sbmatrix}
	&\stackrel{\times\bV_{12} }{\rightarrow}
	\begin{sbmatrix}{\bG_{12}^\top\bB^{(k-1)\top} \bV_{12}}
		\bm{\boxtimes} & \textcolor{brown}{\bm{0}}  &\bm{0} & \bm{0} & \bm{0}\\
		\bm{\boxtimes} & \bm{\boxtimes} & \bm{0} & \bm{0} & \bm{0}\\
		\textcolor{mylightbluetext}{\bm{\boxtimes}} & \bm{\boxtimes} & \boxtimes & 0 & 0 \\
		\bm{0} & \bm{0} & \boxtimes & \boxtimes & 0\\
		\bm{0} & \bm{0}& 0 & \boxtimes & \boxtimes
	\end{sbmatrix}\\
&\stackrel{\bU_{23}^\top \times}{\rightarrow}
\begin{sbmatrix}{\bU_{23}^\top (\cdot) }
	\boxtimes & 0  & 0 & 0 & 0\\
	\bm{\boxtimes} & \bm{\boxtimes} & \textcolor{mylightbluetext}{\bm{\boxtimes}}  & \bm{0} & \bm{0}\\
	\textcolor{brown}{\bm{0}} & \bm{\boxtimes} & \bm{\boxtimes} & \bm{0} & \bm{0} \\
	0 & 0 & \boxtimes & \boxtimes & 0\\
	0 & 0& 0 & \boxtimes & \boxtimes
\end{sbmatrix}
&\stackrel{\times\bV_{23} }{\rightarrow}
\begin{sbmatrix}{\bU_{23}^\top(\cdot ) \bV_{23}}
	\boxtimes & \bm{0}  & \bm{0} & 0 & 0\\
	\bm{\boxtimes} & \bm{\boxtimes} & \textcolor{brown}{\bm{0}}  & \bm{0} & \bm{0}\\
	\bm{0} & \bm{\boxtimes} & \bm{\boxtimes} & \bm{0} & \bm{0} \\
	0 & \textcolor{mylightbluetext}{\bm{\boxtimes}} & \bm{\boxtimes} & \boxtimes & 0\\
	0 &  \bm{0}&  \bm{0}& \boxtimes & \boxtimes
\end{sbmatrix}\\
\end{aligned}
$$
$$
\begin{aligned}
\qquad \qquad \qquad\qquad\,\,\,\,\,&\stackrel{\bU_{34}^\top \times}{\rightarrow}
\begin{sbmatrix}{\bU_{34}^\top(\cdot)}
	\boxtimes & 0  & 0 & 0 & 0\\
	\boxtimes & \boxtimes & 0  & 0 & 0\\
	\bm{0} & \bm{\boxtimes} & \bm{\boxtimes} & \textcolor{mylightbluetext}{\bm{\boxtimes}}  & \bm{0} \\
	\bm{0} & \textcolor{brown}{\bm{0}} & \bm{\boxtimes} & \bm{\boxtimes} & \bm{0}\\
	0 & 0&  0& \boxtimes & \boxtimes
\end{sbmatrix}
&\stackrel{\times\bV_{34} }{\rightarrow}
\begin{sbmatrix}{\bU_{34}^\top(\cdot )\bV_{34}}
	\boxtimes & 0  & \bm{0} & \bm{0} & 0\\
	\boxtimes & \boxtimes & \bm{0}  & \bm{0} & 0\\
	\bm{0} & \bm{\boxtimes} & \bm{\boxtimes} & \textcolor{brown}{\bm{0}} & \bm{0} \\
	\bm{0} & \bm{0} & \bm{\boxtimes} & \bm{\boxtimes} & \bm{0}\\
	0 & 0&  \textcolor{mylightbluetext}{\bm{\boxtimes}}& \bm{\boxtimes}  & \boxtimes
\end{sbmatrix}\\
&\stackrel{\bU_{45}^\top \times}{\rightarrow}
\begin{sbmatrix}{\bU_{45}^\top(\cdot)}
	\boxtimes & 0  & 0 & 0& 0\\
	\boxtimes & \boxtimes & 0  & 0 & 0\\
	0 & \boxtimes & \boxtimes& 0 & 0\\
	\bm{0} & \bm{0} & \bm{\boxtimes} & \bm{\boxtimes} & \textcolor{mylightbluetext}{\bm{\boxtimes}}\\
	\bm{0} & \bm{0} &  \textcolor{brown}{\bm{0}} & \bm{\boxtimes}  & \bm{\boxtimes}
\end{sbmatrix}
&\stackrel{\times\bV_{45} }{\rightarrow}
\begin{sbmatrix}{\bU_{45}^\top(\cdot) \bV_{45}}
	\boxtimes & 0  & 0  & \bm{0}& \bm{0}\\
	\boxtimes & \boxtimes & 0 & \bm{0}& \bm{0}\\
	0 & \boxtimes & \boxtimes& \bm{0}& \bm{0}\\
	\bm{0} & \bm{0} & \bm{\boxtimes} & \bm{\boxtimes} & \textcolor{brown}{\bm{0}}\\
	\bm{0} & \bm{0} &  \bm{0} & \bm{\boxtimes}  & \bm{\boxtimes}
\end{sbmatrix},
\end{aligned}
$$
where $\bG_{12}$ is employed to represent the \textbf{same} construction as in the ``tridiagonal update," while $\bU_{i,i+1}$ for $i\geq 2$ is employed to indicate the \textbf{different} construction based on the implicit Q theorem.  
For clarity, we will denote $\bG_{12}$ as $\bU_{12}$ and $\widetildebG_1$ as $\widetildebU_1$ in the following algorithm.
\end{mdframed}

\paragraph{The complete algorithm.}
Once again, for simplicity, we denote the construction of the $2\times 2$ Givens rotation $\widetildebU_i$ such that $\widetildebU_i^\top \bx =\widetildebU_i^\top 
\footnotesize
\begin{bmatrix}
	x_1\\x_2
\end{bmatrix}= \begin{bmatrix}
	\boxtimes \\ 0
\end{bmatrix}$ (i.e., the second element is 0 after applying the Givens rotation) by the following function:
$$
\widetildebU_i^\top =\text{givens}(x_1,x_2).
$$ 
In all iterations, $\widetildebU_i^\top $ will be of size $2\times 2$. And we denote the construction of the $n\times n$ Givens rotation $\bG_{i,i+1}^\top$ by the following function:
$$
\bU_{i,i+1}^\top=
G(\widetildebU_i^\top )=
\begin{bmatrix}
	\bI_{i-1} & & \\
	& \widetildebU_i^\top& \\
	& & \bI_{n-i-1} 
\end{bmatrix}.
$$
For further simplicity, we will denote $\bU_{i,i+1}^\top$ as $\bU_i^\top$,  which, when  multiplied  on the left of another matrix, modifies the $i$-th and $(i+1)$-th rows; and $\bV_{i,i+1}$ as $\bV_i$. The procedure is outlined in Algorithm~\ref{alg:qr-svd-implifit}, where the shift $\mu^{(k)} = t_{nn}^{(k-1)}$ (i.e., the last diagonal element of $\bT^{(k-1)}$) can be obtained from the last diagonal of the matrix multiplication  $\bB^{(k-1)\top}\bB^{(k-1)}$, i.e., $\mu^{(k)} =b_{n,n-1}^2+b_{nn}^2$.



\begin{algorithm}[H] 
	\caption{Golub-Kahan SVD: Practical QR Algorithm with Implicit Shift (Compare to Algorithm~\ref{alg:qr-algorithm-practical-tridiagonal-implicit-shift})}
	\label{alg:qr-svd-implifit}
\begin{algorithmic}[1] 
\Require Matrix $\bA\in \real^{m\times n}$; 
\State \textcolor{mylightbluetext}{$\bA^\top= \bV^{(0)}\bB^{(0)}\bU^{(0)\top}$}; \Comment{bidiagonal decom. of $\bA^\top$: $\bA\bA^\top=\bU^{(0)}\bB^{(0)\top} \bB^{(0)} \bU^{(0)\top}$}
\For{$k=1,2,\ldots$} 
\State Pick a shift $\mu^{(k)}$; \Comment{e.g., $\mu^{(k)}=t_{nn}^{(k-1)}=b_{n,n-1}^2+b_{nn}^2$}
\State $x_1 \leftarrow t_{11}-\mu^{(k)}, x_2\leftarrow t_{21}$; \Comment{$t_{ij} = t^{(k-1)}_{ij}$}
\State $\bB^{(k)} \leftarrow \bB^{(k-1)}$; \Comment{initialize $\bB^{(k)} $}
\For{$i=1:n-1$}
\State $\bU_{i}^\top \leftarrow G(\widetildebU_i^\top)$, where $\widetildebU_i^\top \leftarrow \text{givens}(x_1,x_2)$;
\State $\bB^{(k)\top} \leftarrow \bU_{i}^\top\bB^{(k)\top}$; \Comment{left update}
\State $x_1 \leftarrow b_{ii}, x_2\leftarrow b_{i,i+1}$; \Comment{$b_{ij} = (\bB^{(k)\top})_{ij}$}
\State $\bV_{i}\textcolor{mylightbluetext}{^\top} \leftarrow G(\widetildebV_i^\top)$, where $\widetildebV_i\textcolor{mylightbluetext}{^\top} \leftarrow \text{givens}(x_1,x_2) \implies  \bV_i = (\bV_i^\top)^\top$;
\State $\bB^{(k)\top} \leftarrow \bU_{i}^\top\bB^{(k)\top}\bV_i$; \Comment{right update}
\If{$i<n-1$}
\State $x_1\leftarrow b_{i+1,i}, x_2\leftarrow b_{i+2,i}$; \Comment{$b_{ij} = (\bB^{(k)\top})_{ij}$}
\EndIf
\EndFor 
\State $\bU^{(k)\top} \leftarrow \bU_{n-1}^\top\ldots \bU_1^\top$;
\State  $\bV^{(k)} \leftarrow \bV_{1}\ldots\bV_{n-1}$; \Comment{this results in $\bB^{(k)^\top} =\bU^{(k)\top} \bB^{(k-1)^\top}\bV^{(k)} $}
\EndFor
\end{algorithmic} 
\end{algorithm}

Again, suppose for iteration $p$, $\bB^{(p)}$ converges to a diagonal matrix (within  machine precision). Then, let us  document the updates for each iteration:
$$
\left.
\begin{aligned}
\bB^{(p)\top}&=\bU^{(p)\top}\bB^{(p-1)\top} \bV^{(p)}\\
\bB^{(p-1)\top}&=\bU^{(p-1)\top}\bB^{(p-2)\top} \bV^{(p-1)}\\
\vdots &= \vdots \\
\bB^{(1)\top}&=\bU^{(1)\top}\bB^{(0)\top} \bV^{(1)}\\
\bB^{(0)\top}&=\bU^{(0)\top}\bA\bV^{(0)}
\end{aligned}
\right\}
\longleftrightarrow
\left\{
\begin{aligned}
\bT^{(p)}&=\bB^{(p)\top}\bB^{(p)}&=&\bU^{(p)\top}\bT^{(p-1)} \bU^{(p)}\\
\bT^{(p-1)}&=\bB^{(p-1)\top}\bB^{(p-1)}&=&\bU^{(p-1)\top}\bT^{(p-2)} \bU^{(p-1)}\\
\vdots &= \vdots \\
\bT^{(1)}&=\bB^{(1)\top}\bB^{(1)}&=&\bU^{(1)\top}\bT^{(0)} \bU^{(1)}\\
\bT^{(0)}&=\bB^{(0)\top}\bB^{(0)}&=&\bU^{(0)\top}\textcolor{mylightbluetext}{\bA\bA^\top} \bU^{(0)}
\end{aligned}
\right.
$$
This yields
$$
\bB^{(p)\top} \approx \bB^{(p)} = 
\underbrace{\bU^{(p)\top} \ldots \bU^{(0)\top}}_{\triangleq\bU^\top}
\bA 
\underbrace{\bV^{(0)}\ldots\bV^{(p)}}_{\triangleq\bV}
\quad\implies\quad
\left\{
\begin{aligned}
	\bA &=\bU \bB^{(p)}\bV^\top\\
	\bA\bA^\top &= \bU \bB^{(p)}\bB^{(p)}\bU^\top
\end{aligned}
\right.,
$$
which approximates the SVD of $\bA$.

\subsection{Jacobi's SVD Method}
We have previously explored   Jacobi's method for computing the spectral decomposition of a matrix, where a pair of orthogonal matrices are applied on the left and right iteratively to reduce the off-diagonal quantity. The orthogonal matrices applied on the left and right  are transposes of each other ($\bQ$ and $\bQ^\top$). It is then straightforward to apply  Jacobi's method using two different sequences of orthogonal matrices  in an attempt to further reduce the off-diagonal quantity. The problem remains
$$
\begin{aligned}
	\bJ_1^\top \bA(k,l)\bJ_2
	&=
	\begin{bmatrix}
		c_1 & -s_1 \\
		s_1 & c_1
	\end{bmatrix}
	\begin{bmatrix}
		a_{kk} & a_{kl} \\
		a_{kl} & a_{ll}
	\end{bmatrix}
	\begin{bmatrix}
		c_2 & s_2 \\
		-s_2 & c_2
	\end{bmatrix}
	=
	\begin{bmatrix}
		\neq 0 & 0 \\
		0 & \neq 0
	\end{bmatrix}.
\end{aligned}
$$
Such a problem can be decomposed into two components:
\begin{enumerate}
\item Find a Jacobi's rotation such that $\widetilde{\bJ}_1\bA(k,l)$ is symmetric;
\item Apply  Jacobi's method to the symmetric $\widetilde{\bJ}_1\bA(k,l)$ so that $\bJ_2^{\top}(\widetilde{\bJ}_1\bA(k,l)) \bJ_2$ becomes diagonal.
\end{enumerate}
Therefore, let $\bJ_1^\top \triangleq \bJ_2^\top\widetilde{\bJ}_1 $, we can diagonalize the submatrix within $\bA$. We will not delve into the details here.

\section{Proof of Results}\label{section:proofs-qralgorithms}
We present proofs for several results in this section.
\begin{proof}[of Lemma~\ref{lemma:qr-algo-from-power}]
The equality of the three sequences at iteration $0$ is trivial due to the same initialization. 
Assume  it holds for the $(k-1)$-th iteration that \{$\widehat{\bA}^{(k-1)} = \bA^{(k-1)}, 
\widehat{\bR}^{(k-1)} = \bR^{(k-1)}, 
\widehat{\bV}^{(k-1)} = \bV^{(k-1)}$\}, we will show these equalities also hold true for the $k$-th iteration that 
\{$\widehat{\bA}^{(k)} = \bA^{(k)}, 
\widehat{\bR}^{(k)} = \bR^{(k)}, 
\widehat{\bV}^{(k)} = \bV^{(k)}$\}, thereby completing the proof by induction.
According to Algorithm~\ref{alg:qr-algorithm-simple1}, 
$$
\left\{
\begin{aligned}
\widehat{\bV}^{(k)}, \widehat{\bR}^{(k)}= QR(\bA\widehat{\bV}^{(k-1)}) &\quad\implies\quad
\widehat{\bV}^{(k)}\widehat{\bR}^{(k)} = \bA\widehat{\bV}^{(k-1)}; \\
\widehat{\bA}^{(k-1)}=\widehat{\bV}^{(k-1)\top} \underbrace{\bA \widehat{\bV}^{(k-1)}}_{\widehat{\bV}^{(k)}\widehat{\bR}^{(k)}} &\quad\implies\quad  
\widehat{\bA}^{(k-1)}= \underbrace{\widehat{\bV}^{(k-1)\top} \widehat{\bV}^{(k)}}_{\text{orthogonal}}\widehat{\bR}^{(k)}.
\end{aligned}
\right.
$$
Since $\widehat{\bV}^{(k-1)\top} \widehat{\bV}^{(k)}$ is a product of two orthogonal matrices, it is also an orthogonal matrix. Therefore, \textcolor{mylightbluetext}{$\widehat{\bA}^{(k-1)}= \widehat{\bV}^{(k-1)\top} \widehat{\bV}^{(k)}\widehat{\bR}^{(k)}$} represents the QR decomposition of $\widehat{\bA}^{(k-1)}$.

According to Algorithm~\ref{alg:qr-algorithm-simple2},
$$
\bQ^{(k)}, \bR^{(k)}= QR(\bA^{(k-1)}) \quad\implies\quad  \bA^{(k-1)}=\bQ^{(k)} \bR^{(k)}.
$$ 
Therefore, another QR decomposition of $\bA^{(k-1)}$ is given by \textcolor{mylightbluetext}{$\bA^{(k-1)}=\bQ^{(k)} \bR^{(k)}$}. Although, the QR decomposition is not unique, if we ensure that the diagonal elements of the upper triangular matrix are positive, the QR decomposition is unique (Corollary~\ref{corollary:unique-qr}). Given that $\widehat{\bA}^{(k-1)} = \bA^{(k-1)}$ (from the assumption),  it follows that
$
\widehat{\bR}^{(k)} = \bR^{(k)} 
$
and 
$$
\begin{aligned}
\bQ^{(k)} = \widehat{\bV}^{(k-1)\top} \widehat{\bV}^{(k)} \quad\implies\quad \widehat{\bV}^{(k-1)}\bQ^{(k)} &=  \widehat{\bV}^{(k)} 
\quad\implies\quad
\underbrace{\bV^{(k-1)}\bQ^{(k)}}_{\bV^{(k)}} =  \widehat{\bV}^{(k)} .
\end{aligned}
$$
Therefore, it follows that
$
\widehat{\bV}^{(k)} = \bV^{(k)}.
$
To establish that  $\widehat{\bA}^{(k)} = \bA^{(k)} $, we proceed as follows:
$$
\begin{aligned}
\widehat{\bA}^{(k)} &=\widehat{\bV}^{(k)\top} \bA \widehat{\bV}^{(k)}  &\text{(By Algorithm~\ref{alg:qr-algorithm-simple1})}\\
&=\bV^{(k)\top} \bA {\bV}^{(k)} &\text{($\widehat{\bV}^{(k)} = \bV^{(k)}$)} \\
&=(\bV^{(k-1)}\bQ^{(k)} )^\top \bA (\bV^{(k-1)}\bQ^{(k)} ) &\text{(By Algorithm~\ref{alg:qr-algorithm-simple2})}\\
&= \bQ^{(k)\top} \widehat{\bA}^{(k-1)} \bQ^{(k)}=\bQ^{(k)\top} \bA^{(k-1)} \bQ^{(k)}  &\text{(By Algorithm~\ref{alg:qr-algorithm-simple1})}\\
&=\bA^{(k)}.    &\text{(By Equation~\eqref{equation:simple-qr-find})}
\end{aligned}
$$
This completes the proof.
\end{proof}

\begin{proof}[of Equation~\eqref{equation:simple-qr-find2}]
The first equality on the left-hand side of Equation~\eqref{equation:simple-qr-find2} is a result of  step 8 in Algorithm~\ref{alg:qr-algorithm-simple2}. 
The second equation follows from the iterative process described in Algorithm~\ref{alg:qr-algorithm-simple1}, which can be expressed as:
$$
\left.
\begin{aligned}
\bA &= \widehat{\bV}^{(k)} \widehat{\bR}^{(k)} \widehat{\bV}^{(k-1)\top} 
=\bV^{(k)} \bR^{(k)} \bV^{(k-1)\top} \\
&= \widehat{\bV}^{(k-1)} \widehat{\bR}^{(k-1)} \widehat{\bV}^{(k-2)\top} 
=\bV^{(k-1)} \bR^{(k-1)} \bV^{(k-2)\top} \\
& \ldots 
\end{aligned}
\right\} \rightarrow
\bA^{k}=\bV^{(k)} \bR^{(k)} \bR^{(k-1)} \ldots \bR^{(0)} .
$$
This completes the proof.
\end{proof}

\begin{proof}[of Equation~\eqref{equation:practical-qr-find2}]
For $k=1$, it is evident that  $(\bA-\mu^{(1)}\bI) = \bQ^{(1)}\bR^{(1)} = \bQ^{(0)}\bQ^{(1)}\bR^{(1)} \bR^{(0)}$. Suppose the statement holds for $k-1$, such that 
$$
(\bA-\mu^{(k-1)}\bI)(\bA-\mu^{(k-2)}\bI)\ldots (\bA-\mu^{(1)}\bI)= \underbrace{\bQ^{(0)}\bQ^{(1)}\bQ^{(2)}\ldots \bQ^{(k-1)}}_{=\bV^{(k-1)}, \text{orthogonal}}
\underbrace{\bR^{(k-1)} \bR^{(k-2)} \ldots \bR^{(0)}}_{\triangleq\bU^{(k-1)}, \text{upper triangular}}.
$$
If we can prove it also holds for $k$, then we complete the proof by induction. To see this, we have 
$$
\begin{aligned}
&\gap (\bA-\mu^{(k)}\bI)(\bA-\mu^{(k-1)}\bI)\ldots (\bA-\mu^{(1)}\bI)
=(\bA-\mu^{(k)}\bI)\bV^{(k-1)}\left(\bR^{(k-1)} \bR^{(k-2)} \ldots \bR^{(0)}\right)\\
&=\bigg(\underbrace{\bV^{(k) } \bA^{(k)} \bV^{(k)\top}}_{\text{Equation}~\eqref{equation:practical-qr-find-root}}    -\mu^{(k)}\underbrace{\bV^{(k)}\bV^{(k)\top} }_{=\bI}\bigg) \bV^{(k-1)}\left(\bR^{(k-1)} \bR^{(k-2)} \ldots \bR^{(0)}\right)\\
&=\bV^{(k) }\left( \bA^{(k)} -\mu^{(k)}  \bI   \right)\bV^{(k)\top}  \bV^{(k-1)}\left(\bR^{(k-1)} \bR^{(k-2)} \ldots \bR^{(0)}\right)\\ 
&=\bV^{(k)}\bigg( \underbrace{\bR^{(k)}\bQ^{(k)}  }_{\text{step 8}}\bigg)
\underbrace{( \bV^{(k-1)}\bQ^{(k)} )^\top}_{\text{step 9}}
\bV^{(k-1)}\left(\bR^{(k-1)} \bR^{(k-2)} \ldots \bR^{(0)}\right)\\
&=\bV^{(k)}\left(\bR^{(k)}\bR^{(k-1)} \bR^{(k-2)} \ldots \bR^{(0)}\right).
\end{aligned}	
$$
This completes the proof.
\end{proof}

\begin{problemset}
\item \label{prob:power_nai1} Prove that $\bv_2^{(k)}$ in \eqref{equation:power_nai1} is a scalar multiple of $\bA^k \bv_2^{(0)}$: $\bv_2^{(k)} = c_{2k} \bA^k \bv_2^{(0)} = \frac{\bA^k \bv_2^{(0)}}{\norm{\bA^k \bv_2^{(0)}}}$. \textit{Hint: Prove by induction.}	

\item \label{prob:reverse_qr_trid} Let $\bT\in \real^{n\times n}$ be symmetric and tridiagonal with a QR decomposition $\bT=\bQ\bR$. Show that $\bT_+ = \bR\bQ $ is also symmetric and tridiagonal.

\item \textbf{Rayleigh quotient.}	 Let $\bA \in\real^{n\times n}$ be symmetric, and suppose that at least one eigenvalue of $\bA$ is positive. Show that 
$\lambda_{\max}(\bA) = \max \left\{ \frac{\bx^\top \bA \bx}{\bx^\top \bx} : \bx^\top \bA \bx = 1 \right\}$. \textit{Hint: The maximum value of the Rayleigh quotient over all nonzero vectors $\bx$ is the largest eigenvalue $\lambda_{\max}(\bA)$ of $\bA$. This maximum is achieved when $\bx$ is an eigenvector corresponding to $\lambda_{\max}(\bA)$. Consider $\bx^\top\bA\bx=\gamma^2\by^\top\bA\by=1$ with $\bx=\gamma\by$.}
	
\item Find the conditional gradient (i.e., the Frank-Wolfe method) update for the inverse power method.
\item Find the eigenvalues of the following matrices using the power method and the QR algorithm:
$$
\begin{bmatrix}
1 & 1 \\
1 & 0
\end{bmatrix}
\qquad 
\text{and}
\qquad 
\begin{bmatrix}
2 & -4\\
-3 & 3
\end{bmatrix}.
$$

\item Find the SVD for $\bA=\footnotesize\begin{bmatrix}
0 & 1 \\
1 & 0
\end{bmatrix}$.

\item \textbf{Givens rotation.} Given the 2-nd order Givens rotation $\bG=\footnotesize\begin{bmatrix}
c & s \\
-s & c
\end{bmatrix}$, and two vectors $\bx, \by\in\real^n$, find the values of $c$ and $s$ such that $[\bx, \by]\bQ$ contains orthogonal columns.

\item \textbf{Implicit shift QR algorithm.} Find the Givens rotation to chase the $3$ in the $(3,1)$-th entry of the following matrix:
$$
\bA=
\begin{bmatrix}
5 & 1 & 0 & 0 \\
2 & 2 & 5 & 0 \\
3 & 4 & 5 & 1\\
0 & 0 & 2 & 5 
\end{bmatrix}.
$$

\item Use the Weierstrass theorem (Theorem~\ref{theorem:weierstrass_them}) to demonstrate that the optimization problem in \eqref{equation:eign_opt} has a solution.
\end{problemset}

%% file: chapter-als.tex
\newpage 
\part{Special Topics}

\newpage
\chapter{Alternating Least Squares (ALS)}\label{chapter:als}
\begingroup
\hypersetup{
	linkcolor=structurecolor,
	linktoc=page,  
}
\minitoc \newpage
\endgroup

\index{Least squares}
\index{Linear models}
\index{Regression analysis}
\section{Preliminary: Least Squares Approximations}\label{section:pre_ls}
\lettrine{\color{caligraphcolor}T}
The linear model stands as the primary technique in regression analysis, utilizing the least squares approximation as the fundamental tool to minimize the sum of squared errors (refer to Section~\ref{section:application-ls-qr}).
This approach is a natural choice when seeking the regression function that minimizes the corresponding expected squared error. 
Over the past few decades, linear models have been widely applied in diverse areas, including decision-making  \citep{dawes1974linear}, time series analysis \citep{christensen1991linear, lu2017machine}, quantitative finance \citep{menchero2011barra}, and numerous other disciplines like production science, social science, and soil science \citep{fox1997applied, lane2002generalized, schaeffer2004application, mrode2014linear}.

To be more concrete, let's consider an overdetermined system represented by $\bb = \bA\bx $, where $\bA\in \real^{m\times n}$ represents the \textit{input data matrix} (also known as the \textit{predictor variables}), $\bb\in \real^m$ is the \textit{observation vector} (or \textit{target/response vector}), and the  number of samples $m$ exceeds  the dimensionality $n$. 
The vector $\bx$ represents a vector of weights of the linear model. 
Typically, $\bA$ will have full column rank since real-world data is often uncorrelated or can be made so through preprocessing.
In practical scenarios, a \textit{bias term} (a.k.a., an \textit{intercept}) is added to the first column of $\bA$ to enable the least squares method to solve for:
\begin{equation}\label{equation:ls-bias}
	\widetildebA \widetildebx = 
[\bm{1} ,\bA ] 
\begin{bmatrix}
	x_0\\
	\bx
\end{bmatrix}
 = \bb .
\end{equation}


However, it is common for the equation $\bb = \bA\bx$ to have  no exact solution (the system is \textit{inconsistent}) due to the  being overdetermined----that is, there are more equations than unknowns.
Define the column space of $\bA$ as $\{\bA\bgamma: \,\, \forall \bgamma \in \real^n\}$, denoted by $\cspace(\bA)$.
In essence, when we say $\bb = \bA\bx$ has no solution, it implies that $\bb$ lies outside the column space of $\bA$. 
In other words, the error $\be = \bb -\bA\bx$ cannot be reduced to zero. 
In such cases, the goal shifts to minimizing the error, typically measured by the mean squared error.
The resulting solution $\bx_{LS}$, which minimizes $\normtwo{\bb-\bA\bx_{LS}}^2$, is known as the \textit{least squares solution}. The least squares method is a cornerstone of mathematical sciences, and there are numerous resources dedicated entirely to this topic, including works by   \citet{trefethen1997numerical, strang2019linear, strang2021every,  lu2021rigorous}.


\paragraph{Least squares by calculus.}
When $\normtwo{\bb-\bA\bx}^2$ is differentiable and the parameter space of $\bx$ spans the entire space $\real^n$ (i.e., unconstrained optimization)
\footnote{In other words, the \textit{domain} of the optimization problem $\mathop{\min}_{\bx} \normtwo{\bb-\bA\bx}^2$ is the entire space  $\real^n$.}, 
the least squares estimate corresponds to the root of the gradient of $\normtwo{\bb-\bA\bx}^2$. 
This leads us to the following lemma.~\footnote{Variants of the least squares problem are explored in Problems~\ref{problem:rls}$\sim$\ref{problem:twls2}.}

\index{Fermat's theorem}
\begin{lemma}[Least Squares by Calculus]\label{lemma:ols}
Let  $\bA \in \real^{m\times n}$ be a  fixed data matrix with full rank  and $m\geq n$ (i.e., its columns  are linearly independent)~\footnote{See Problem~\ref{prob:als_pseudo1}$\sim$\ref{prob:als_pseudon} for a relaxation of this condition using the pseudo-inverse.}. 
For the overdetermined system $\bb = \bA\bx$, the least squares solution, obtained by employing calculus and  setting the partial derivatives in every direction of $\normtwo{\bb-\bA\bx}^2$ to  zero (i.e., the gradient vanishes), is given by $\bx_{LS} = (\bA^\top\bA)^{-1}\bA^\top\bb$ \footnote{This is known as the \textit{first-order optimality condition} for local optima points. Note that the proof of the first-order optimality condition for multivariate functions strongly relies on the first-order optimality conditions for univariate functions, which is also known as  \textit{Fermat's theorem}. See Problem~\ref{problem:fist_opt}.}. The value, $\bx_{LS} = (\bA^\top\bA)^{-1}\bA^\top\bb$, is commonly referred to as the \textit{ordinary least squares (OLS)} estimate or simply the \textit{least squares (LS)} estimate of $\bx$.
\end{lemma}

To prove the lemma above, we need to demonstrate the invertibility of $\bA^\top\bA$. Since we assume $\bA$ has full rank and $m\geq n$, $\bA^\top\bA \in \real^{n\times n}$ is invertible if it has a rank of $n$, which is equivalent to the rank of $\bA$. This assertion has been validated in Lemma~\ref{lemma:rank-of-ata}.
\begin{proof}[of Lemma \ref{lemma:ols}]
Recalling from calculus that  a function $f(\bx)$ attains a minimum  at  $\bx_{LS}$ when the gradient $\nabla f(\bx)=\bzero$. The gradient of $\normtwo{\bb-\bA\bx}^2$ is given by $2\bA^\top\bA\bx -2\bA^\top\bb$. $\bA^\top\bA$ is invertible since we assume $\bA$ is fixed and has full rank with $m\geq n$ (Lemma~\ref{lemma:rank-of-ata}). 
Consequently, the OLS solution for $\bx$ is $\bx_{LS} = (\bA^\top\bA)^{-1}\bA^\top\bb$, which completes the proof.
\end{proof}

\begin{definition}[Normal Equation]\label{definition:normal-equation-als}
The condition for the gradient of $\normtwo{\bb-\bA\bx}^2$ to be zero can be expressed as $\bA^\top\bA \bx_{LS} = \bA^\top\bb$. The equation is also known as the \textit{normal equation}. 
Under the assumption that $\bA$ has full rank with $m\geq n$, it follows that $\bA^\top\bA$ is invertible, implying $\bx_{LS} = (\bA^\top\bA)^{-1}\bA^\top\bb$.
\end{definition}

\index{Convex functions}
\begin{figure}[h!]
\centering  
\vspace{-0.35cm} 
\subfigtopskip=2pt 
\subfigbottomskip=2pt 
\subfigcapskip=-5pt 
\subfigure[A convex function.]{\label{fig:convex-1}
\includegraphics[width=0.26\linewidth]{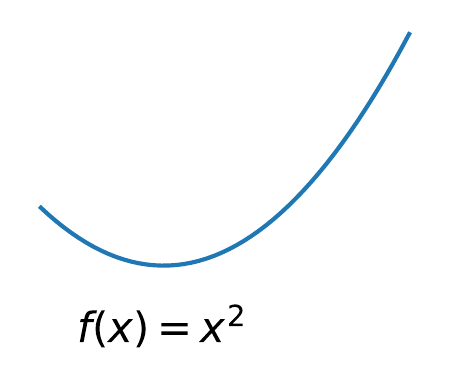}}
\subfigure[A concave function.]{\label{fig:convex-2}
\includegraphics[width=0.26\linewidth]{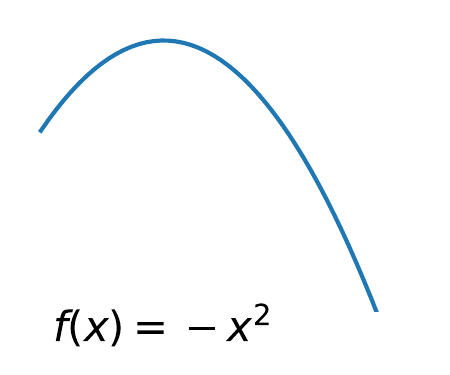}}
\subfigure[A random function.]{\label{fig:convex-3}
\includegraphics[width=0.26\linewidth]{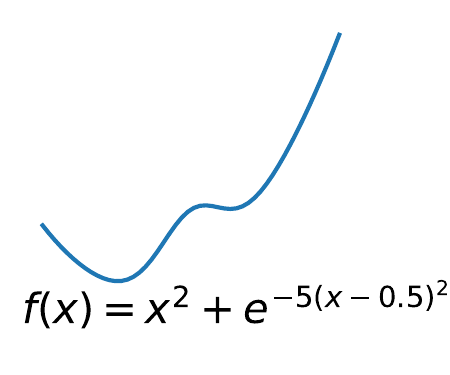}}
\caption{Three functions.}
\label{fig:convex-concave-none}
\end{figure}
However, it is not certain whether the least squares estimate obtained in Lemma~\ref{lemma:ols} is the smallest, largest, or neither.
An example illustrating this ambiguity is presented in Figure~\ref{fig:convex-concave-none}. 
What we can confidently assert is the existence of at least one root for the gradient of the function  $f(\bx)=\normtwo{\bb-\bA\bx}^2$, and this root is a necessary condition for the minimum point rather than a sufficient condition.
The following remark provides clarification on this matter.
\begin{remark}[Verification of Least Squares Solution]
Why does a zero gradient imply the least mean squared error? 
The typical explanation stems from convex analysis, as we will see shortly.  
However, here we directly confirm that the OLS solution indeed minimizes the mean squared error.  
For any $\bx \neq \bx_{LS}$, we have 
\begin{equation}
\begin{aligned}
\normtwo{\bb - \bA\bx}^2 &= \normtwo{\bb - \bA\bx_{LS} + \bA\bx_{LS} - \bA\bx}^2
= \normtwo{\bb-\bA\bx_{LS} + \bA (\bx_{LS} - \bx)}^2 \\
&=\normtwo{\bb-\bA\bx_{LS}}^2 + \normtwo{\bA(\bx_{LS} - \bx)}^2 + 2\left(\bA(\bx_{LS} - \bx)\right)^\top(\bb-\bA\bx_{LS}) \\ 
&=\normtwo{\bb-\bA\bx_{LS}}^2 + \normtwo{\bA(\bx_{LS} - \bx)}^2 + 2(\bx_{LS} - \bx)^\top(\bA^\top\bb - \bA^\top\bA\bx_{LS}), \nonumber
\end{aligned} 
\end{equation}
where the third term is zero due to the normal equation, and $\normtwo{\bA(\bx_{LS} - \bx)}^2 \geq 0$. Therefore,
$
\normtwo{\bb - \bA\bx}^2 \geq \normtwo{\bb-\bA\bx_{LS}}^2. 
$
Thus, we demonstrate that the OLS estimate indeed corresponds to the minimum, rather than the maximum or a saddle point~\footnote{
A \textit{saddle point} is a point at which the gradient vanishes (a \textit{stationary point}), and  there exists a direction where the objective function decreases and another direction where it increases.
}, using the calculus approach.
As a matter of fact, this condition from the least squares estimate is also known as the \textit{sufficiency of stationarity under convexity}. When $\bx$ is defined over the entire space $\real^n$, this condition is also known as the \textit{necessity of stationarity under convexity}.
\end{remark}
\index{Saddle point}

Another question that arises is: Why does the normal equation seemingly ``magically" yield solutions for  $\bx$?
A simple example can help illustrate the answer. The equation $x^2=-1$ has no real solution. 
However, $x\cdot x^2 = x\cdot (-1)$ does have a real solution $\hat{x} = 0$, in which case, $\hat{x}$ minimizes the difference between $x^2$ and $-1$, making them as close as possible.

\begin{example}[Altering the Solution Set by Left Multiplication]
Consider the  data matrix and target vector:
$\tiny
\bA=
\begin{bmatrix}
-3 & -4 \\
4 & 6  \\
1 & 1
\end{bmatrix}
$
and
$
\bb=
\tiny
\begin{bmatrix}
1  \\
-1   \\
0
\end{bmatrix}
.
$
It can be easily verified that the system $\bA\bx = \bb$ has no solution for $\bx$. 
However, if we multiply both sides on the left by
$
\bB=
\scriptsize
\begin{bmatrix}
0 & -1 & 6\\
0 & 1  & -4
\end{bmatrix},
$
then  $\bx_{LS} = [1/2, -1/2]^\top$ becomes the solution to $\bB\bA\bx= \bB\bb$. 
This example illustrates why the normal equation can lead to the least squares solution. Multiplying a linear system on the left alters the solution set, effectively projecting the problem into a different space where the least squares solution can be computed.
\end{example}

\paragraph{Rank-deficiency.}
In this discussion, we assume $\bA\in \real^{m\times n}$ has full rank with $m\geq n$, ensuring that $\bA^\top\bA$ is invertible. 
However, if two or more columns of $\bA$ are perfectly correlated, the matrix $\bA$ becomes deficient, and $\bA^\top\bA$ becomes singular. 
To address this issue, one can choose $\bx$ that minimizes $\normtwo{\bx_{LS}}^2$ while satisfying the normal equation. 
That is, we select the least squares solution with the smallest magnitude. 
In Sections~\ref{section:ls-utv} and~\ref{section:application-ls-svd}, we briefly discuss how to use UTV decomposition and SVD to address this rank-deficient least squares problems.
See Problems~\ref{prob:als_pseudo1}$\sim$\ref{prob:als_pseudon} or the following paragraph for more insights.

\index{Contion number}
\index{Tikhonov regularization}
\index{$\ell_2$ regularization}
\paragraph{Regularizations and stability.}
A common  problem that arise in the ordinary least square solution is the near-singularity of $\bA$.
Let the SVD of $\bA$ be $\bA=\bU\bSigma\bV^\top\in\real^{m\times n}$, where $\bU\in\real^{m\times m}$ and $\bV\in\real^{n\times n}$ are orthogonal, and the main diagonal of $\bSigma\in\real^{m\times n}$ contains the singular values. Consequently, $\bA^\top\bA = \bV(\bSigma^\top\bSigma)\bV^\top \triangleq \bV\bS\bV^\top$, where $\bS\triangleq \bSigma^\top\bSigma  = \diag(\sigma_1^2, \sigma_2^2, \ldots,\sigma_n^2)\in\real^{n\times n}$ contains the squared singular values of $\bA$. When $\bA$ is nearly singular, $\sigma_n^2\approx 0$, making the inverse operation $(\bA^\top\bA)^{-1} = \bV\bS^{-1}\bV^\top$ numerically unstable. 
As a result, the solution $\bx_{LS} =(\bA^\top\bA)^{-1}\bA^\top\bb $ may diverge.
To address this issue, we typically  add an $\ell_2$ regularization term to obtain the solution for the following optimization problem:
\begin{equation}
\bx_{Tik} = \mathop{\argmin}_{\bx} \normtwo{\bb-\bA\bx}^2 +\lambda\normtwo{\bx}^2.
\end{equation}
This method is known as the  \textit{Tikhonov regularization method} (or simply the $\ell_2$ regularized method) \citep{tikhonov1963solution}.
The gradient of the problem is $2(\bA^\top\bA+\lambda\bI)\bx-2\bA^\top\bb$. Thus, the least squares solution is given by 
$
\bx_{Tik} = (\bA^\top\bA+\lambda\bI)^{-1}\bA^\top\bb.
$
The inverse operation becomes $(\bA^\top\bA+\lambda\bI)^{-1} = \bV(\bS+\lambda\bI)^{-1}\bV^\top$, where $\widetildebS\triangleq(\bS+\lambda\bI)=\diag(\sigma_1^2+\lambda, \sigma_2^2+\lambda, \ldots,\sigma_n^2+\lambda)$. 
The solutions for OLS and Tikhonov regularized LS are given, respectively, by 
\begin{equation}
\begin{aligned}
\bx_{LS} &= (\bA^\top\bA)^{-1}\bA^\top\bb = \bV\left(\bS^{-1}\bSigma\right)\bU^\top\bb;\\
\bx_{Tik} &= (\bA^\top\bA+\lambda\bI)^{-1}\bA^\top\bb = \bV\left((\bS+\lambda\bI)^{-1}\bSigma\right)\bU^\top\bb,\\
\end{aligned}
\end{equation}
where the main diagonals of $\left(\bS^{-1}\bSigma\right)$ are $\diag(\frac{1}{\sigma_1}, \frac{1}{\sigma_2}, \ldots, \frac{1}{\sigma_n})$; and the main diagonals of $\left((\bS+\lambda\bI)^{-1}\bSigma\right)$ are $\diag(\frac{\sigma_1}{\sigma_1^2+\lambda}, \frac{\sigma_2}{\sigma_2^2+\lambda}, \ldots, \frac{\sigma_n}{\sigma_n^2+\lambda})$. The latter solution is more stable if $\lambda$ is greater than the   smallest nonzero squared singular value.
The condition number becomes smaller if  the smallest singular value $\sigma_n$ is close to zero (see Equation~\eqref{equation:qr_condition_num} or Appendix~\ref{appendix:condition_number}):
$$
\kappa(\bA^\top\bA) = \frac{\sigma_1^2}{\sigma_n^2}
\qquad \rightarrow \qquad
\kappa(\bA^\top\bA+\lambda\bI) = \frac{\lambda+\sigma_1^2}{\lambda+\sigma_n^2}.
$$
Tikhonov regularization effectively prevents  divergence  in the least squares solution 
$\bx_{LS} = (\bA^\top\bA)^{-1} \bA^\top \bb$ when the matrix $\bA$ is nearly singular or even rank-deficient. This improvement enhances the convergence properties of both the LS algorithm and its variants, such as alternating least squares,  while addressing identifiability issues  in various settings (see Section~\ref{section:regularization-extention-general}). As a result, Tikhonov regularization has become a widely applied technique.

\index{Data least squares}
\index{Total least squares}
\paragraph{Data least squares.}
While the OLS method accounts for errors in the response variable $\bb$, the \textit{data least sqaures (DLS)} method considers errors in the predictor variables:
\begin{equation}
\bx_{DLS} = \mathop{\argmin}_{\bx, \widetildebA} \normf{\widetildebA}^2, \quad \text{s.t.}\quad \bb\in\cspace(\bA+\widetildebA),
\end{equation}
where $\widetildebA$ represents a perturbation of $\bA$ (i.e., a noise in the predictor variables).
That is, $(\bA+\widetildebA) \bx_{DLS} = \bb$, assuming the measured response $\bb$ is noise-free.
The Lagrangian function and its gradient w.r.t. $\bx$ are, respectively, given by
$$
\begin{aligned}
L(\bx, \widetildebA, \blambda) &= \trace(\widetildebA\widetildebA^\top) +\blambda^\top (\bA\bx+\widetildebA\bx-\bb);\\
\nabla_{\widetildebA} L(\bx, \widetildebA,\blambda) &= \widetildebA+\blambda\bx^\top = \bzero \quad\implies\quad \widetildebA=-\blambda\bx^\top,
\end{aligned}
$$
where $\blambda\in\real^m$  is a vector of Lagrange multipliers.
Substituting the value of the vanishing gradient into $(\bA+\widetildebA) \bx = \bb$ yields $\blambda = \frac{\bA\bx-\bb}{\bx^\top\bx}$ and $\widetildebA=-\frac{(\bA\bx-\bb)\bx^\top}{\bx^\top\bx} $.
Therefore, using the invariance of cyclic permutation of factors in trace,  the objective function becomes 
$$
\mathop{\argmin}_{\bx}
\frac{(\bA\bx-\bb)^\top (\bA\bx-\bb)}{\bx^\top\bx} .
$$

\paragraph{Total least squares.} Similar to  data least squares, the \textit{total least squares (TLS)} method considers errors in both the predictor variables and the response variables. The TLS problem can be formulated as:
\begin{equation}
\bx_{TLS} = \mathop{\argmin}_{\bx, \widetildebA, \widetildebb} \normf{[\widetildebA, \widetildebb]}^2, 
\quad \text{s.t.}\quad (\bb+\widetildebb)\in\cspace(\bA+\widetildebA), 
\end{equation}
where $\widetilde{\bA}$ and $\widetilde{\bb}$ are perturbations in the predictor variables and the response variable, respectively.
Let $\bC\triangleq[\bA,\bb]\in\real^{m\times (n+1)}$,  $\bD\triangleq[\widetildebA, \widetildebb]\in\real^{m\times (n+1)}$, and $\by\in\footnotesize\begin{bmatrix}
\bx\\
-1
\end{bmatrix}$, the problem can be equivalently stated as
\begin{equation}
\bx_{TLS} = \mathop{\argmin}_{\by, \bD} \normf{\bD}^2, 
\quad \text{s.t.}\quad \bD\by = -\bC\by, 
\end{equation}

\index{Decomposition: ALS}
\index{Netflix}
\section{Netflix Recommender and Matrix Factorization}\label{section:als-netflix}
The explosion of data from advancements in sensor technology and computer hardware has presented new challenges for data analysis. 
The large volume of data often contains noise and other distortions, requiring preprocessing for deductive science to be applied. 
For example, signals received by antenna arrays often are contaminated by noise and other degradations. 
To effectively analyze the data, it is necessary to reconstruct or represent it in a way that reduces inaccuracies while maintaining certain feasibility conditions.

Additionally, in many cases, the data collected from complex systems is the result of multiple interrelated variables acting in unison. When these variables are not well defined, the information contained in the original data can be overlapping and unclear. By creating a reduced system model, we can achieve a level of accuracy that is close to the original system. The common approach in removing noise, reducing the model, compressing the data, and reconstructing feasibility is to replace the original data with a lower-dimensional representation obtained through subspace approximation.
Therefore, \textit{low-rank matrix approximations (LRMA) or low-rank matrix decompositions}  play a crucial role in 
a wide range of applications such as compression, feature selection and noise filtering.~\footnote{More accurately, the term ``approximation" usually refers to a situation where a matrix $\bA$ is expressed as $\bA\approx\bW\bZ$. In this case, 
$\bW$ and $\bZ$ are matrices whose product is not exactly equal to $\bA$, but rather a close representation or estimate of $\bA$. On the other hand, the term ``decomposition" typically implies that 
$\bA$ can be exactly represented by the product of $\bW$ and $\bZ$, i.e., $\bA=\bW\bZ$.
However, in this context, we use the terms ``approximation" and ``decomposition" interchangeably, even though they technically have distinct meanings. This means that when discussing matrix factorizations, we may refer to both exact  and approximate representations using either term.
}

Low-rank matrix decomposition is a powerful technique   used in machine learning and data mining  to represent a given matrix as the product of two or more matrices with lower dimensions.
It is used to capture the essential structure of a matrix while ignoring noise and redundancies. The most common methods for low-rank matrix decomposition include singular value decomposition (SVD), principal component analysis (PCA),  multiplicative update nonnegative matrix factorization (NMF), and the alternating least squares (ALS) approach, which we will introduce in this section.

For example, in the Netflix Prize competition \citep{bennett2007netflix}, the goal is to predict the ratings of users for different movies, given the existing ratings (resp. interaction) of those users for other movies (resp. items).
We index $M$ movies with $m= 1, 2,\ldots,M$ and $N$ users
with $n = 1, 2,\ldots,N$. (In the matrix approximation context, lowercase letters e.g., $m,n,k$, are used for the subscripts in running indices, while  uppercase letters $M, N, K$ denote the upper bound of an index.) We denote the rating of the $n$-th user for the $m$-th movie by $a_{mn}$. 
Define $\bA$ as an $M \times N$ rating matrix (a \textit{movie-by-user matrix}) with columns $\{\ba_n\} \in \real^M$ containing ratings provided by the $n$-th user (also known as the \textit{preference matrix}). Note that many ratings $\{a_{mn}\}$ are missing, and our goal is to predict these missing ratings accurately, i.e., to complete the matrix.

It is  clear that without some inherent structure in the matrix, and consequently in the way users rate items, there would be no relationship between the observed and unobserved entries. This would mean there is no unique method to complete the matrix. Therefore, it is crucial to impose some structure on the matrix. A common structural assumption is that of low rank: we aim to fill in the missing entries of matrix $\bA$, assuming $\bA$ is a low-rank matrix. This assumption makes the problem well-posed and allows for a unique solution, as the low-rank structure connects the entries of the matrix (i.e., a \textit{matrix completion} problem). Consequently, the unobserved entries can no longer be independent of the observed values. \footnote{However, this is actually a strong assumption. Consider the matrix $\bA=\sum_{i=1}^{r} \be_i\widetilde{\be}_j^\top$, where $\be_i$ is the standard basis for $\real^M$ and $\widetilde{\be}_j$ is the standard basis for $\real^N$. $\bA$ is a rank-$r$ matrix and contains only $r$ nonzero entries.  
In recommendation settings, we typically observe only a few random entries of the matrix. As a result, there is a high possibility  that some entries may never be observed. This poses a significant challenge for the matrix completion problem.
However, we will not cover this topic in this book.}
It is important to note that, except for very special data structures, a matrix cannot be compressed/decomposed without incurring some compression error, since a low-rank matrix representation is only an approximation of the original matrix.
This procedure, often known as \textit{collaborative filtering}, seeks to exploit co-occurring patterns in the observed behaviors across users in order to predict future  behaviors of users.

Consider the mask matrix $\bM\in \{0,1\}^{M\times N}$, where $m_{mn}\in \{0,1\}$ indicates whether  user $n$ has rated  movie $m$ or not.
Then the low-rank matrix completion problem can be formulated as 
\begin{equation}
\widetildebA = \mathop{\argmin}_{\bX\in\real^{M\times N}} \sum_{m,n=1}^{M,N} (x_{mn} - a_{mn})^2\cdot m_{mn}\gap \text{s.t.} \gap \rank(\bX)\leq K.
\end{equation}
However, this problem is NP-hard (non-deterministic polynomial) \citep{hardt2014computational}. While it can be equivalently written (proof from singular value decomposition) in the following  unconstrained form:
\begin{equation}
\widetildebA =\widetildebW\widetildebZ = \mathop{\argmin}_{\substack{\bW\in\real^{M\times K}\\ \bZ\in\real^{K\times N}}} \sum_{m,n=1}^{M,N} ((\bW\bZ)_{mn}- a_{mn})^2\cdot m_{mn},
\end{equation}
which allows for indirect solution or approximation using alternate algorithms.

We then formally consider algorithms for solving the following problem: The matrix $\bA$ is approximately factorized into an $M\times K$ matrix $\bW$ and a $K \times  N$ matrix $\bZ$. 
Typically, $K$ is selected to be smaller than both $M$ and $N$, ensuring that $\bW$ and $\bZ$ have reduced dimensions compared to the original  matrix $\bA$. 
This reduction in dimensionality results in a compressed version of the original data matrix. 
An appropriate decision on the value of $K$ is critical in practice; but the choice of $K$ is very
often problem-dependent.
The factorization is significant in the sense that if $\bA=[\ba_1, \ba_2, \ldots, \ba_N]$ and $\bZ=[\bz_1, \bz_2, \ldots, \bz_N]$ are the column partitions of $\bA$ and $\bZ$, respectively, then we have $\ba_n = \bW\bz_n$. This means each column $\ba_n$ is approximated by a linear combination of the columns of $\bW$, weighted by the components in $\bz_n$. 
Therefore, the columns of $\bW$ can be thought of as containing the column basis (\textit{template columns}, or the approximation of the column basis) of $\bA$; and $\bz_n$ indicates the coordinates (or \textit{activations}) of $\ba_n$ in the basis $\bW$. 
This concept is similar to the factorization methods discussed in the data interpretation part (Part~\ref{part:data-interation}). 
The key difference is that we do not restrict  $\bW$ to consist of exact columns from $\bA$.

\index{Cross-validation}

\index{Two-block coordinate descent}
\begin{algorithm}[H] 
\caption{2-Block Coordinate Descent: Framework of Most ALS and NMF Algorithms}
\label{alg:two_bcd_gen_inals}
\begin{algorithmic}[1] 
\Require A loss function for a variable with two blocks $\bX=(\bW,\bZ)$: $f(\bX)=f(\bW,\bZ)$, and data matrix $\bA$;
\Ensure Constraint on $\bW$ and $\bZ$;
\State Generate some initial matrices $\bW^{(0)}$ and $\bZ^{(0)}$;
\For{$t = 1, 2, \ldots$}
\State $\bW^{(t)} \leftarrow \text{update}\big(\bA, \bZ^{(t-1)}, \bW^{(t-1)}\big)$;
\State $\bZ^{(t)} \leftarrow \text{update}\big(\bA, \bW^{(t)T}, \bZ^{(t-1)T}\big)$;
\EndFor
\end{algorithmic} 
\end{algorithm}
However, in most cases, the resulting factorization problem has no exact solution, thus requiring optimization procedures to find suitable numerical approximations. The problem is usually solved using a \textit{two-block coordinate descent (2-BCD)} approach (see Algorithm~\ref{alg:two_bcd_gen_inals} for a general illustration).
In order to obtain the approximation $\bA\approx\bW\bZ$, we must establish a loss function such that the distance between $\bA$ and $\bW\bZ$ can be measured. 
In our discussion, the chosen loss function is the Frobenius norm  (a.k.a., the Euclidean distance, Definition~\ref{definition:frobernius-in-svd}) between two matrices, which vanishes to zero if $\bA=\bW\bZ$, and its advantages will become evident shortly.

To simplify the problem, let's first assume that there are no missing ratings. 
We project data vectors $\ba_n\in\real^M$ into a lower dimension $\bz_n \in \real^K$  with $K<\min\{M, N\}$
in a way that the \textit{reconstruction error} (a.k.a., \textit{criterion function, objective function, cost function, or loss function}) as measured by the Frobenius norm (a.k.a., sum of squared loss) is minimized (assume $K$ is known):
\begin{equation}\label{equation:als-per-example-loss2}
L(\bW,\bZ)\triangleq D(\bA,\bW\bZ) \triangleq \frac{1}{2}\sum_{n=1}^N \sum_{m=1}^{M} \left(a_{mn} - \bw_m^\top\bz_n\right)^2 
=\frac{1}{2} \normf{\bW\bZ-\bA}^2,~\footnote{
Note that we include a scaling factor of $\frac{1}{2}$ for easier discussion of gradients. Minimizing over $\frac{1}{2}\normf{\bW\bZ-\bA}^2$ is equivalent to minimizing over $\normf{\bW\bZ-\bA}^2$ or $\normf{\bW\bZ-\bA}$.
The choice of the Frobenius norm assumes i.i.d. Gaussian noise on the data ($\bA=\bW\bZ+\bN$, where each entry of $\bN$ follows i.i.d. Gaussian noise) and leads to a smooth optimization via least squares. When the loss is measured by the $\ell_1$ matrix norm (Definition~\ref{definition:lp-matrix_norm_app}), one obtains a robust low-rank matrix factorization; and the noise is assumed i.i.d. Laplace.}
\end{equation}
where $\bW=[\bw_1^\top; \bw_2^\top; \ldots; \bw_M^\top]\in \real^{M\times K}$ and $\bZ=[\bz_1, \bz_2, \ldots, \bz_N] \in \real^{K\times N}$ contain $\bw_m$'s and $\bz_n$'s as \textbf{rows and columns}, respectively. 
In \eqref{equation:als-per-example-loss2},  $L(\bW,\bZ)$ indicates it is a loss function w.r.t. $\bW$ and $\bZ$, and $D(\bA, \bW\bZ)$ implies it is a distance/divergence~\footnote{In words, the \textit{distance} $D(\bE,\bF)$ indicates $D(\bE,\bF)=D(\bF,\bE)\geq 0$ and the equality holds if and only if $\bE=\bF$; while the \textit{divergence} holds that  $D(\bE,\bF)\neq D(\bF,\bE)\geq 0$ and the equality holds if and only if $\bE=\bF$.} between $\bA$ and $\bW\bZ$ (we will use the two terms interchangeably  when necessary).

Moreover, the loss function $L(\bW,\bZ)=\frac{1}{2} \normf{\bW\bZ-\bA}^2$ demonstrates convexity concerning $\bZ$ when $\bW$ is held constant, and analogously, with respect to $\bW$ when $\bZ$ is fixed. 
This characteristic  motivates the alternate algorithm that alternately fixes one of the variables and
optimizes over the other.
Therefore, we can first minimize the loss with respect to $\bZ$
while keeping $\bW$ fixed, and subsequently minimize it with respect to $\bW$ with $\bZ$ fixed.
This leads to two optimization problems, denoted by ALS1 and ALS2, respectively:
$$
\left\{
\begin{aligned}
	\bZ &\leftarrow \mathop{\arg \min}_{\bZ} L(\bW,\bZ);    \qquad \text{(ALS1)} \\ 
	\bW &\leftarrow \mathop{\arg \min}_{\bW} L(\bW,\bZ). \qquad \text{(ALS2)}
\end{aligned}
\right.
$$
This is known as the \textit{two-block coordinate descent (2-BCD)  algorithm} as mentioned previously, where we alternate between optimizing the least squares with respect to $\bW$ and $\bZ$. 
Hence, it is also called the \textit{alternating least squares (ALS)} algorithm \citep{comon2009tensor, takacs2012alternating, giampouras2018alternating}. Convergence is guaranteed if the loss function $L(\bW,\bZ)$ decreases at each iteration, and we shall discuss this further  in the sequel.

\index{Coordinate descent algorithm}
\index{Convexity}
\index{Global minimum}
\index{ALS}
\begin{remark}[Convexity and Global Minimum]
Although the loss function defined by the Frobenius norm $\frac{1}{2}\normf{\bW\bZ-\bA}^2$ is convex either with respect to  $\bW$ when $\bZ$ is fixed or vice versa (called \textit{marginally convex}, Definition~\ref{definition:marginal_convexfuncs}), it is not \textit{jointly convex} in both variables simultaneously. Therefore, locating the global minimum is infeasible. 
However, the convergence is assured to find local minima.

More generally, let $D(\bA, \bB)$ be convex in the second argument $\bB$. Then, $D(\bA,\bW\bZ)$ is convex in $\bW$ for $\bZ$ fixed and vice versa; see Problem~\ref{prob:sep_conv}.
\end{remark}

\subsection*{Given $\bW$, Optimizing $\bZ$}

Now, let's examine the problem of $\bZ \leftarrow \mathop{\argmin}_{\bZ} L(\bW,\bZ)$. When there exists a unique minimum of the loss function $L(\bW,\bZ)$ with respect to $\bZ$, we refer to it as the \textit{least squares} minimizer of $\mathop{\argmin}_{\bZ} L(\bW,\bZ)$. 
With $\bW$ fixed, $L(\bW,\bZ)$  can be represented as $L(\bZ|\bW)$ (or more compactly, as $L(\bZ)$) to emphasize  the dependence on $\bZ$:
$$
\begin{aligned}
2L(\bZ|\bW) &= \normf{\bW\bZ-\bA}^2= \left\Vert\bW[\bz_1,\bz_2,\ldots, \bz_N]-[\ba_1,\ba_2,\ldots,\ba_N]\right\Vert^2=
\normf{
\scriptsize
\begin{bmatrix}
\bW\bz_1 - \ba_1 \\
\bW\bz_2 - \ba_2\\
\vdots \\
\bW\bz_N - \ba_N
\end{bmatrix}
}^2. 
\end{aligned}
$$
Now, if we define 
$$
\footnotesize
\widetildebW \triangleq
\begin{bmatrix}
	\bW & \bzero & \ldots & \bzero\\
	\bzero & \bW & \ldots & \bzero\\
	\vdots & \vdots & \ddots & \vdots \\
	\bzero & \bzero & \ldots & \bW
\end{bmatrix}
\in \real^{MN\times KN}, 
\gap 
\widetildebz\triangleq
\begin{bmatrix}
	\bz_1 \\ \bz_2 \\ \vdots \\ \bz_N
\end{bmatrix}
\in \real^{KN},
\gap 
\widetildeba\triangleq
\begin{bmatrix}
	\ba_1 \\ \ba_2 \\ \vdots \\ \ba_N
\end{bmatrix}
\in \real^{MN},
$$
then the (ALS1) problem can be reduced to the ordinary least squares problem for minimizing  $\big\Vert{\widetildebW \widetildebz - \widetildeba}\big\Vert_2^2$ with respect to $\widetildebz$. And the solution is given by 
$
\widetildebz = (\widetildebW^\top\widetildebW)^{-1} \widetildebW^\top\widetildeba.
$
However, it is not advisable  to obtain the result using this approach, as computing the inverse of  $\widetildebW^\top\widetildebW$ requires $2(KN)^3$ flops (Theorem~\ref{theorem:inverse-by-lu2}).
Alternatively, a more direct way to solve the  (ALS1)  problem  is to find the gradient of $L(\bZ|\bW)$ with respect to $\bZ$ (assuming all  partial derivatives of this function exist): 
\begin{equation}\label{equation:givenw-update-z-allgd}
\begin{aligned}
\nabla_{\bZ} L(\bZ|\bW) &= \frac{1}{2}
\frac{\partial \,\,\trace\left((\bW\bZ-\bA)(\bW\bZ-\bA)^\top\right)}{\partial \bZ}
\stackrel{\star}{=}  \bW^\top(\bW\bZ-\bA) \in \real^{K\times N},
\end{aligned}
\end{equation}
where the first equality arises from the definition of the Frobenius norm (Definition~\ref{definition:frobernius-in-svd}) such that $\normf{\bA} = \sqrt{\sum_{m=1,n=1}^{M,N} (a_{mn})^2}=\sqrt{\trace(\bA\bA^\top)}$, and equality ($\star$)  is a consequence of the fact that $\frac{\partial \trace(\bA\bA^\top)}{\partial \bA} = 2\bA$. When the loss function is a differentiable function of $\bZ$, we can determine the least squares solution using differential calculus. 
Since we optimize over an open set $\real^{K\times N}$,  any minimum of the function 
$L(\bZ|\bW)$ must be a root of the equation:
$$
\nabla_{\bZ} L(\bZ|\bW)  = \bzero.
$$
By finding the root of the  equation above, we obtain the ``candidate" update for $\bZ$ that finds the minimizer of $L(\bZ|\bW)$:
\begin{equation}\label{equation:als-z-update}
\textbf{(``Candidate" update for $\bZ$)}: \gap {\bZ = (\bW^\top\bW)^{-1} \bW^\top \bA  \leftarrow \mathop{\arg \min}_{\bZ} L(\bZ|\bW).}
\end{equation}
This computation requires $2K^3$ flops to compute the inverse of $\bW^\top\bW$, compared to $2(KN)^3$ flops to get the inverse of $\widetildebW^\top\widetildebW$.
Prior to confirming that a root of the equation above is indeed a minimizer (as opposed to a maximizer, hence the term ``candidate" update), it is imperative to establish the convexity of the function. 
For a twice continuously differentiable function, this verification  can be equivalently achieved by confirming (see Problem~\ref{problem:pos_hessian} for more details): 
$$
\nabla^2_{\bZ} L(\bZ|\bW) \succ 0.~
\footnote{In short, a twice continuously differentiable function $f$ over an open convex set $\sS$ is called \textit{convex} if and only if $\nabla^2f(\bx)\succeq  \bzero $ for any $\bx\in \sS$ (sufficient and necessary for convex); and called \textit{strictly convex} if $\nabla^2f(\bx)\succ \bzero$ for any $\bx\in \sS$ (only sufficient for strictly convex, e.g., $f(x)=x^6$ is strictly convex, but $f^{\prime\prime}(x)=30x^4$ is equal to zero at $x=0$.). 
And when the convex function $f$ is a continuously differentiable function over a convex set $\sS$, the stationary point $\nabla f(\bx^\star)=\bzero$ of $\bx^\star\in\sS$ is  a \textit{global minimizer} of $f$ over $\sS$.
In our context, when given $\bW$ and updating $\bZ$, the function is defined over the entire space $\real^{K\times N}$.
}
$$
That is, the Hessian matrix is positive definite (Definition~\ref{definition:psd-pd-defini}; see, for example, \citet{beck2014introduction}). To demonstrate this, we explicitly express the Hessian matrix as
\begin{equation}\label{equation:als-z-update_hessian}
\nabla^2_{\bZ} L(\bZ|\bW)= \widetildebW^\top\widetildebW \in \real^{KN\times KN},
~\footnote{
A block-diagonal matrix whose block matrix on the diagonal is $\bW^\top\bW$. And it can be equivalently denoted as $\nabla^2_{\bZ} L(\bZ|\bW) = \diag(\bW,\bW,\ldots,\bW)^\top\diag(\bW,\bW,\ldots,\bW)$.
Using the Kronecker product defined in Definition~\ref{definition:kronecker-product}, this can be equivalently written as $\nabla^2_{\bZ} L(\bZ|\bW) = \bI_{N} \kronecker (\bW^\top\bW)$, where $\bI_N$ is the $N\times N$ identity matrix.
}
\end{equation}
which maintains full rank if $\bW\in \real^{M\times K}$ has full rank  and $K<M$ (Lemma~\ref{lemma:rank-of-ata}).

\begin{remark}[Positive Definite Hessian if $\bW$ Has Full Rank]
We  claim that if $\bW\in\real^{M\times K}$ has full rank $K$ with $K<M$, then $\nabla_{\bZ}^2 L(\bZ|\bW)$ is positive definite. This can be demonstrated by confirming that when $\bW$ has full rank, the equation $\bW\bx=\bzero$  holds true only when $\bx=\bzero$, since the null space of $\bW$ has a  dimension of 0. Therefore, 
$$
\bx^\top (\bW^\top\bW)\bx >0, \qquad \text{for any nonzero vector $\bx\in \real^K$}.
$$ 
And this in turn implies $\widetildebW^\top\widetildebW \succ \bzero $.
\end{remark}
Now, the problem becomes  \textcolor{black}{\textbf{whether $\bW$ has full rank so that the Hessian of $L(\bZ|\bW)$ is positive definite}}; otherwise, we cannot claim the update of $\bZ$ in Equation~\eqref{equation:als-z-update} reduces the loss (due to convexity) so that the matrix decomposition progressively improves the approximation of the original matrix $\bA$ by $\bW\bZ$ in each iteration.
We will address the positive definiteness of the Hessian matrix shortly, relying on the following lemma.
\begin{lemma}[Rank of $\bZ$ after Updating]\label{lemma:als-update-z-rank}
	Suppose $\bA\in \real^{M\times N}$ has full rank with \textcolor{mylightbluetext}{$M\leq N$} and $\bW\in \real^{M\times K}$ has full rank with $K<M$ (i.e., $K<M\leq N$). Then the update of $\bZ=(\bW^\top\bW)^{-1} \bW^\top \bA \in \real^{K\times N}$ in Equation~\eqref{equation:als-z-update} has full rank.
\end{lemma}
\begin{proof}[of Lemma~\ref{lemma:als-update-z-rank}]
Since $\bW^\top\bW\in \real^{K\times K}$ has full rank if $\bW$ has full rank (Lemma~\ref{lemma:rank-of-ata}), it follows that $(\bW^\top\bW)^{-1} $ has full rank. 

Suppose $\bW^\top\bx=\bzero$. This implies that $(\bW^\top\bW)^{-1} \bW^\top\bx=\bzero$. Thus, the following two null spaces satisfy:
$
\nspace(\bW^\top) \subseteq \nspace\left((\bW^\top\bW)^{-1} \bW^\top\right).
$
Moreover, suppose $\bx$ lies in the null space of $(\bW^\top\bW)^{-1} \bW^\top$ such that $(\bW^\top\bW)^{-1} \bW^\top\bx=\bzero$. And since $(\bW^\top\bW)^{-1} $ is invertible, it implies $ \bW^\top\bx=(\bW^\top\bW)\bzero=\bzero$, leading to
$
\nspace\left((\bW^\top\bW)^{-1} \bW^\top\right)\subseteq \nspace(\bW^\top).
$
Consequently, through ``sandwiching," it follows that 
\begin{equation}\label{equation:als-z-sandiwch1}
\nspace(\bW^\top) = \nspace\left((\bW^\top\bW)^{-1} \bW^\top\right).
\end{equation}
Therefore, $(\bW^\top\bW)^{-1} \bW^\top$ has full rank $K$. Let $\bT\triangleq(\bW^\top\bW)^{-1} \bW^\top\in \real^{K\times M}$, and suppose $\bT^\top\bx=\bzero$. This implies $\bA^\top\bT^\top\bx=\bzero$, yielding 
$
\nspace(\bT^\top) \subseteq \nspace(\bA^\top\bT^\top).
$
Similarly, suppose $\bA^\top(\bT^\top\bx)=\bzero$. Since $\bA$ has full rank with the dimension of the null space being 0: $\dim\left(\nspace(\bA^\top)\right)=0$, $(\bT^\top\bx)$ must be zero. The claim follows  since $\bA$ has full rank $M$ with the row space of $\bA^\top$ being equal to the column space of $\bA$, where $\dim\left(\cspace(\bA)\right)=M$ and  $\dim\left(\nspace(\bA^\top)\right) = M-\dim\left(\cspace(\bA)\right)=0$. 
Consequently, $\bx$ is in the null space of $\bT^\top$ if $\bx$ is in the null space of $\bA^\top\bT^\top$:
$
\nspace(\bA^\top\bT^\top)\subseteq \nspace(\bT^\top).
$
By ``sandwiching" again, we obtain
\begin{equation}\label{equation:als-z-sandiwch2}
\nspace(\bT^\top) = \nspace(\bA^\top\bT^\top).
\end{equation}
Since $\bT^\top$ has full rank $K<M\leq N$, it follows that $\dim\left(\nspace(\bT^\top) \right) = \dim\left(\nspace(\bA^\top\bT^\top)\right)=0$.
Therefore,
$\bZ^\top=\bA^\top\bT^\top$ has full rank $K$.
We complete the proof.
\end{proof}

\subsection*{Given $\bZ$, Optimizing $\bW$}
The matrix factorization problem exhibits symmetry: $\bA=\bW\bZ$ if and only if $\bA^\top=\bZ^\top\bW^\top$ such that $D(\bA, \bW\bZ)=D(\bA^\top, \bZ^\top\bW^\top)$. The analysis of optimizing $\bW$ given $\bZ$ directly follows  from the previously discussed methodology.
Below, we provide a brief outline of the results.
With $\bZ$ fixed, $L(\bW,\bZ)$ can be expressed as $L(\bW|\bZ)$ (or more compactly, as $L(\bW)$)  to emphasize  the dependence on $\bW$:
$
\begin{aligned}
	L(\bW|\bZ) &= \frac{1}{2}\normf{\bW\bZ-\bA}^2.
\end{aligned}
$
To solve the optimization problem (ALS2) directly, we need to compute the gradient of  $L(\bW|\bZ)$ with respect to $\bW$:
$$
\begin{aligned}
\nabla_{\bW} L(\bW|\bZ) &= \frac{1}{2}
\frac{\partial \trace\left((\bW\bZ-\bA)(\bW\bZ-\bA)^\top\right)}{\partial \bW}
= (\bW\bZ-\bA)\bZ^\top \in \real^{M\times K}.
\end{aligned}
$$
Similarly, the ``candidate" update for  $\bW$ can be obtained by identifying the root of the gradient $\nabla_{\bW} L(\bW|\bZ)$:
\begin{equation}\label{equation:als-w-update}
\textbf{(``Candidate" update for $\bW$)}:\gap{\bW^\top = (\bZ\bZ^\top)^{-1}\bZ\bA^\top  \leftarrow \mathop{\arg\min}_{\bW} L(\bW|\bZ).}
\end{equation}
Once more, we emphasize that the update is merely a ``candidate" update. 
Further validation is necessary  to ascertain the positive definiteness of the Hessian matrix.
The Hessian matrix is obtained as follows:
\begin{equation}\label{equation:als-w-update_hessian}
\begin{aligned}
\nabla_{\bW}^2 L(\bW|\bZ) =\widetildebZ\widetildebZ^\top \in \real^{KM\times KM},
\end{aligned}
\end{equation}
where $\widetildebZ\triangleq\diag(\bZ,\bZ,\ldots,\bZ)\in\real^{KM\times NM}$ is defined analogously to $\widetildebW$ in \eqref{equation:als-z-update_hessian}.
Therefore, by similar analysis, if $\bZ$ has full rank with $K<N$, the Hessian matrix is positive definite.

In Lemma~\ref{lemma:als-update-z-rank}, we proved that $\bZ$ has full rank under certain conditions, ensuring that  the Hessian matrix in Equation~\eqref{equation:als-w-update_hessian} is positive definite, and the update in Equation~\eqref{equation:als-w-update} exists.
We now prove that $\bW$ also has full rank under certain conditions, such that the Hessian in Equation~\eqref{equation:als-z-update_hessian} is positive definite, and the update in  Equation~\eqref{equation:als-z-update} exists.
\begin{lemma}[Rank of $\bW$ after Updating]\label{lemma:als-update-w-rank}
Suppose $\bA\in \real^{M\times N}$ has full rank with \textcolor{mylightbluetext}{$M\geq N$} and $\bZ\in \real^{K\times N}$ has full rank with $K<N$ (i.e., $K<N\leq M$). Then the update of $\bW^\top = (\bZ\bZ^\top)^{-1}\bZ\bA^\top$ in Equation~\eqref{equation:als-w-update} has full rank.
\end{lemma}
The proof of Lemma~\ref{lemma:als-update-w-rank} is similar to that of  Lemma~\ref{lemma:als-update-z-rank}, and we shall not repeat the details.

\paragraph{Key observation.}
Combining the observations in Lemma~\ref{lemma:als-update-z-rank} and Lemma~\ref{lemma:als-update-w-rank}, as long as we \textcolor{mylightbluetext}{initialize $\bZ$ and $\bW$ to have full rank}, the updates in Equation~\eqref{equation:als-z-update} and Equation~\eqref{equation:als-w-update} are reasonable \textbf{since the Hessians in Equation~\eqref{equation:als-z-update_hessian} and \eqref{equation:als-w-update_hessian} are positive definite}. 
\textbf{
Note that we need an additional condition to satisfy  
both Lemma~\ref{lemma:als-update-z-rank} 
and Lemma~\ref{lemma:als-update-w-rank}: $M=N$, meaning there must be an equal number of movies and   users.
} 
We will relax this condition in the next section through regularization.
(Alternatively,  Problems~\ref{prob:als_pseudo1}$\sim$\ref{prob:als_pseudon} relax this condition using pseudo-inverse.)
We summarize the process in Algorithm~\ref{alg:als}.
Since the loss $\frac{1}{2}\normf{\bA-\bW\bZ}^2$ in each iteration is monotonically nonincreasing and bounded below, it converges (Section~\ref{section:als_convergence}).
In particular, $\nabla_{\bZ} L(\bZ|\bW)=\bzero$ and $\nabla_{\bW} L(\bW|\bZ)=\bzero$ when $iter\rightarrow \infty$.

\begin{algorithm}[H] 
\caption{Alternating Least Squares}
\label{alg:als}
\begin{algorithmic}[1] 
\Require Matrix $\bA\in \real^{M\times N}$ \textcolor{mylightbluetext}{with $M= N$};
\State Initialize $\bW\in \real^{M\times K}$, $\bZ\in \real^{K\times N}$ \textcolor{mylightbluetext}{with full rank and $K<M= N$}; 
\State Choose a stoping criterion on the approximation error $\delta$;
\State Choose the maximum number of iterations $C$;
\State $iter=0$; \Comment{Count for the number of iterations}
\While{$\normf{\bA-\bW\bZ}>\delta $ and $iter<C$} 
\State $iter=iter+1$;
\State $\bZ \leftarrow (\bW^\top\bW)^{-1} \bW^\top \bA  \leftarrow \mathop{\arg \min}_{\bZ} L(\bZ|\bW)$;
\State $\bW^\top \leftarrow (\bZ\bZ^\top)^{-1}\bZ\bA^\top  \leftarrow \mathop{\arg\min}_{\bW} L(\bW|\bZ)$;
\EndWhile
\State Output $\bW,\bZ$;
\end{algorithmic} 
\end{algorithm}


\section{More on the Error Measure and Statistical Interpretation$^*$}\label{section:more_err_sta_als}
We briefly introduce alternative error measures and their statistical interpretation for matrix factorization/approximation problems in this section.
\paragraph{Frobenius norm and Gaussian noise.}
As mentioned previously, the choice of the Frobenius norm assumes i.i.d. Gaussian noise on the entries of the underlying matrix, leading  to a smooth optimization via least squares.
To see this, we denote $\bB\triangleq\bW\bZ$ and assume the noise $\epsilon$ is i.i.d. Gaussian with mean zero and variance $\sigma^2$:
\begin{equation}\label{equation:gau_noise}
a_{mn}= b_{mn}+ \epsilon_{mn}, \gap 
\epsilon_{mn}\sim \normal(0, \sigma^2), 
\gap \forall m,n.
\end{equation}
Under this i.i.d. Gaussian  assumption, the log-likelihood function is 
$$
\ln p(\bA\mid \bB, \sigma^2) 
=
-\frac{1}{2\sigma^2}\sum_{m,n=1}^{M,N} (a_{mn}-b_{mn})^2 + C(\sigma)
=
-\frac{1}{2\sigma^2}\normf{\bA-\bW\bZ}^2 +C(\sigma),
$$
where $C(\sigma^2)$ is a term dependent on  $\sigma$.
Therefore, given the observed values in $\bA$, the unknown parameters $\bW,\bZ$, and $\sigma^2$ can be estimated by maximizing the log-likelihood function (maximum likelihood estimator, MLE):
$$
\mathopmax{\bW,\bZ,\sigma^2} -\frac{1}{2\sigma^2}\normf{\bA-\bW\bZ}^2
=
\mathopmin{\bW,\bZ,\sigma^2}\frac{1}{2\sigma^2}\normf{\bA-\bW\bZ}^2.
$$
The  reduces to \eqref{equation:als-per-example-loss2} if the noise level is constant. 
When the noise is not i.i.d.,  the loss function becomes a weighted $\ell_2$ norm between $\bA$ and $\bW\bZ$; see Problem~\ref{prob:non_iid_gaus}. 

\paragraph{Matrix $\ell_1$ norm and Laplace noise.}
Similarly, if we assume the noise is i.i.d. Laplace with location $\mu=0$ and scale $\sigma$, the log-likelihood function becomes 
$$
\ln p(\bA\mid \bB, \sigma) = -\frac{1}{2\sigma}\sum_{m,n=1}^{M,N}\abs{a_{mn}-b_{mn}} +C(\sigma)
=
-\frac{1}{2\sigma} \normmone{\bA-\bW\bZ} +C(\sigma),
$$
where $\normmone{\cdot}$ denotes the matrix $\ell_1$ norm (Definition~\ref{definition:lp-matrix_norm_app}).
The MLE is equivalent to the matrix factorization with the loss function given by the matrix $\ell_1$ norm if the scale is held constant.
Same as the Gaussian noise case, the Laplace noise is also additive:
$$
a_{mn}= b_{mn}+ \epsilon_{mn}, \gap 
\epsilon_{mn}\sim \text{Laplace}(0, \sigma), 
\gap \forall m,n.
$$
In practice, the (matrix) $\ell_1$ norm is more robust than the Frobenius norm in the presence of outliers.

\paragraph{Matrix $\ell_\infty$ norm and uniform noise.}
Alternative, an i.i.d. uniform noise can be  assumed:
$$
a_{mn} \sim \text{Uniform}(b_{mn} -\sigma, b_{mn}+\sigma), 
\gap \forall m,n,
$$
where $\text{Uniform}(a, b)$ is the uniform distribution over the interval $[a,b]$.
The uniform noise is additive such that 
$$
a_{mn}= b_{mn}+ \epsilon_{mn}, \gap 
\epsilon_{mn}\sim \text{Uniform}(-\sigma, \sigma), 
\gap \forall m,n.
$$
The MLE corresponds to minimizing the matrix $\ell_\infty$  norm (Definition~\ref{definition:lp-matrix_norm_app}):
$$
L(\bW,\bZ) = \normminf{\bA-\bW\bZ} = \mathopmax{m,n} \abs{a_{mn}-b_{mn}}.
$$
\paragraph{Poisson and nonnegative integers.}
When the entries of $\bA$ are nonnegative integers (count data, such as vector of word or counts in text mining), the Poisson distribution can be assumed \citep{lu2023bayesian}:
$$
p(a_{mn}=x\mid b_{mn}) = \frac{b_{mn}^{x}}{x!} \exp(-b_{mn})
\,\,\text{ with }\,\,
\Exp[a_{mn}] = b_{mn}, \,\, \Var[a_{mn}] = b_{mn},
\gap \forall m,n.
$$
This indicates $a_{mn}=0$ if $b_{mn}=0$.
The MLE corresponds to minimizing the Kullback-Leibler (KL) divergence between $\bA$ and $\bW\bZ$:
\begin{equation}\label{equation:als_poi_los}
L(\bW,\bZ) = \sum_{m,n=1}^{M,N} l(a_{mn}, b_{mn}) = \sum_{m,n=1}^{M,N}\left(a_{mn} \ln(\frac{a_{mn}}{b_{mn}}) -a_{mn}+b_{mn}\right),
\end{equation}
where $l(a_{mn}, b_{mn})$ is assumed to be 0 if $b_{mn}=0$. 
The term $2l(a_{mn}, b_{mn}) =2( \ln p(a_{mn} \mid a_{mn}) - \ln p(a_{mn} \mid b_{mn}))$ is called the \textit{deviance (goodness-of-fit statistic)}, which measures the goodness of fit of the model by  comparing the log-likelihood difference between the \textit{saturated model} and the \textit{current model}.
The Poisson noise is not additive.

\paragraph{Multiplicative Gamma noise.}
The Gaussian, Laplace, and uniform noises are all additive noise. The Gamma noise, on the contrary, assumes multiplicative noise:
$$
a_{mn}= b_{mn} \cdot  \epsilon_{mn}, \gap 
\epsilon_{mn}\sim \text{Gamma}(r, \lambda)=\frac{\lambda^r}{\Gamma(r)} \epsilon_{mn}^{r-1}\exp(-\lambda \epsilon_{mn}), 
\gap \forall m,n.
$$
where $\Gamma(r)$ is the Gamma function.
When the mean $\Exp[a_{mn}]=\frac{r}{\lambda}=1$, 
the MLE corresponds to minimizing the following loss (\textit{Itakura-Saito divergence, IS divergence} \citep{itakura1968analysis}):
\begin{equation}\label{equation:als_gama_los}
L(\bW,\bZ) = D(\bA,\bB) = \sum_{m,n=1}^{M,N} \frac{a_{mn}}{b_{mn}} - \ln\left(\frac{a_{mn}}{b_{mn}}\right)-1.
\end{equation}
\begin{exercise}[Scale Invariant of IS]
Show that the IS divergence is scale invariant: $D(\bA,\bB) = D(\gamma\bA, \gamma\bB)$.
\end{exercise}
The scale invariant property is particularly useful in audio source separation problems since low-power frequency bands can perceptually contribute as much as high-power frequency bands \citep{gillis2020nonnegative}.

\index{Regularization}
\section{Regularization and Identifiability: Extension to General Matrices}\label{section:regularization-extention-general}

\textit{Regularization} is a machine learning technique employed to prevent overfitting and improve model generalization. Overfitting occurs when a model is overly complex and fits the training data too closely, resulting in poor performance on  unseen data. 
To mitigate this issue, regularization introduces a constraint or a penalty term into the loss function used for model optimization, discouraging the development of overly complex models. 
This creates  a trade-off between having a simple, generalizable model and fitting the training data well. 
Common types of regularization include $\ell_1$ regularization, $\ell_2$ regularization (Tikhonov regularization), and elastic net regularization (a combination of $\ell_1$ and $\ell_2$ regularizations). 
Regularization finds extensive applications in machine learning algorithms such as linear regression, logistic regression, and neural networks.

In the context of the alternating least squares problem, we can add an $\ell_2$ regularization term  to minimize the following loss:
\begin{equation}\label{equation:als-regularion-full-matrix}
L(\bW,\bZ)  = \frac{1}{2}\normf{\bW\bZ-\bA}^2 +\frac{1}{2}\lambda_w \normf{\bW}^2 + \frac{1}{2}\lambda_z \normf{\bZ}^2, \qquad \lambda_w>0, \lambda_z>0,
\end{equation}
where the gradient with respect to $\bZ$ and $\bW$ are given, respectively, by 
\begin{equation}\label{equation:als-regulari-gradien}
\left\{
\begin{aligned}
\nabla_{\bZ} L(\bZ|\bW) &= \bW^\top(\bW\bZ-\bA) + \lambda_z\bZ \in \real^{K\times N};\\
\nabla_{\bW} L(\bW|\bZ)  &= (\bW\bZ-\bA)\bZ^\top + \lambda_w\bW \in \real^{M\times K}.
\end{aligned}
\right.
\end{equation}
The Hessian matrices are given, respectively, by 
$$
\left\{
\begin{aligned}
\nabla^2_{\bZ} L(\bZ|\bW) &= \widetildebW^\top\widetildebW+ \lambda_z\bI \in \real^{KN\times KN};\\
\nabla^2_{\bW} L(\bW|\bZ)  &= \widetildebZ\widetildebZ^\top + \lambda_w\bI \in \real^{KM\times KM}, \\
\end{aligned}
\right.
$$
which are positive definite due to the perturbation by the regularization:
$$
\left\{
\begin{aligned}
\bx^\top (\widetildebW^\top\widetildebW +\lambda_z\bI)\bx 
&= \underbrace{\bx^\top\widetildebW^\top\widetildebW\bx}_{\geq 0} + \lambda_z \normtwo{\bx}^2>0, \gap \text{for nonzero $\bx$};\\
\bx^\top (\widetildebZ\widetildebZ^\top +\lambda_w\bI)\bx 
&= \underbrace{\bx^\top\widetildebZ\widetildebZ^\top\bx}_{\geq 0} + \lambda_w \normtwo{\bx}^2>0,\gap \text{for nonzero $\bx$}.
\end{aligned}
\right.
$$
\textbf{The regularization ensures that the Hessian matrices remain positive definite, even if $\bW$ and $\bZ$ are rank-deficient}. 
Consequently, matrix decomposition can be extended to any matrix, regardless of whether $M>N$ or $M<N$. In rare cases, $K$ even can be chosen as $K>\max\{M, N\}$ to obtain a \textit{high-rank approximation} of $\bA$. However, in most scenarios, we aim to find a \textit{low-rank approximation} of $\bA$ with $K<\min\{M, N\}$. For example,  ALS can be employed to discover low-rank neural networks, reducing  memory usage of  neural networks while improving performance (Section~\ref{section:low-rank-neural}).
Therefore, the minimizers can be determined by identifying the roots of the gradients:
\begin{equation}\label{equation:als-regular-final-all}
\left.
\begin{aligned}
\bZ &= (\bW^\top\bW+ \lambda_z\bI)^{-1} \bW^\top \bA 
\gap \text{and}\gap 
\bW^\top = (\bZ\bZ^\top+\lambda_w\bI)^{-1}\bZ\bA^\top .
\end{aligned}
\right.
\end{equation}
The regularization parameters $\lambda_z, \lambda_w\in \real_{++}$ are used to balance the trade-off
between the accuracy of the approximation and the smoothness of the computed solution. The selection of these parameters is typically problem-dependent and can be determined through \textit{cross-validation} (CS). Again, we summarize the process in Algorithm~\ref{alg:als-regularizer}.
We will also introduce the \textit{alternating direction methods of multipliers  (ADMM)} for solving  matrix factorization problems with $\ell_2$ or $\ell_1$ regularization in Section~\ref{section:nmf_admm_all}, where the method can be extended to other regularizations and constraints, such as  nonnegative constraints.

The $\ell_2$ (or $\ell_1$ ) regularizations can be applied to generalize the ALS problem to general matrices.
However, we will consider the missing entries in the following sections, thus the problem becomes a matrix completion formulation.
In this sense, the $\ell_1$ and $\ell_2$ regularizations  are not the only applicable regularizations; for example, the \textit{nuclear norm} \footnote{Also called the \textit{Schatten 1-norm} or \textit{trace norm}; see Definition~\ref{definition:ky_fan_norm}.} of $\bW\bZ$ (the sum of singular values of the matrix) can be applied, for which the \textit{Soft-Impute for matrix completion} algorithm guarantees the recovery of the matrix when the number of observed entries $z$ satisfies
$
z\geq C r n \log n,
$
where the underlying matrix $\bA$ is of size $\real^{n\times n}$ and $C > 0$ is a fixed universal constant \citep{gross2011recovering, hastie2015statistical}. 
However, the $\ell_2$ regularization on $\bW$ and $\bZ$ can somehow be reformulated into the nuclear norm form (see Problem~\ref{problem:nuclear_equi}).

\index{Cross-validation}
\begin{algorithm}[H] 
\caption{Alternating Least Squares with Regularization}
\label{alg:als-regularizer}
\begin{algorithmic}[1] 
\Require Matrix $\bA\in \real^{M\times N}$;
\State Initialize $\bW\in \real^{M\times K}$, $\bZ\in \real^{K\times N}$ \textcolor{mylightbluetext}{randomly without condition on the rank and the relationship between $M, N, K$}; 
\State Choose a stoping criterion on the approximation error $\delta$;
\State Choose regularization parameters $\lambda_w, \lambda_z$;
\State Choose the  maximal number of iterations $C$;
\State $iter=0$; \Comment{Count for the number of iterations}
\While{$\normf{\bA-\bW\bZ}>\delta $ and $iter<C$}
\State $iter=iter+1$; 
\State $\bZ \leftarrow (\bW^\top\bW+ \lambda_z\bI)^{-1} \bW^\top \bA  \leftarrow \mathop{\arg \min}_{\bZ} L(\bZ|\bW)$;
\State $\bW^\top \leftarrow (\bZ\bZ^\top+\lambda_w\bI)^{-1}\bZ\bA^\top  \leftarrow \mathop{\arg\min}_{\bW} L(\bW|\bZ)$;
\EndWhile
\State Output $\bW,\bZ$;
\end{algorithmic} 
\end{algorithm}

\index{K-means problem}
\index{Regularization}
\index{Constraint}
\index{Identifiability}
\index{Sparsity}
\index{Group sparsity}
\paragraph{Regularization as constraints and identifiability.} 
Regularization terms, such as $\lambda_w\normf{\bW}^2$ in \eqref{equation:als-regularion-full-matrix}, can be interpreted as  constraints $\normf{\bW}\leq C$, where $C$ is a constant, via Lagrangian multipliers (see, for example, \citet{boyd2004convex}). 
Different constraints can be placed on the factors $\bW$ and $\bZ$. For example, the nonnegative constraint discussed in Chapter~\ref{chapter:nmf} and the sparsity constraint discussed in Section~\ref{section:reg_geom_inter}.
Moreover, the two matrices $\bW\in\real^{M\times K}$ and $\bZ\in\real^{K\times N}$ have $(M+N)K$ degrees of freedom.
However, due to the scaling degree of freedom of the columns of $\bW$ and rows of $\bZ$ in $\bA=\bW\bZ$, the factorization $\bW\bZ$ has $(M+N-1)K$ degrees of freedom: $\bW[:, k]\bZ[k,:] = (\gamma\bW[:, k])(\frac{1}{\gamma}\bZ[k,:])$ for any scalar $\gamma\neq 0$ and $k\in\{1,2,\ldots,K\}$.
Therefore, the factorization is not identifiable.  Regularization helps reduce overfitting and addresses the issue of identifiability by incorporating prior information through regularizations/constraints.
Besides the $\ell_2$ and $\ell_1$ regularizations mentioned above, we may also impose alternatives \citep{lee2009semi, bach2011convex, cai2010graph, iordache2012total, gillis2020nonnegative}:
\begin{itemize}
\item \textbf{Minimum-volume.} Impose the regularizer $\lambda_w \det(\bW^\top\bW)$ such that the factor has a small volume and is as close as possible to the data points.

\item \textbf{K-means constraint (vector quantization).}
In the context of the \textit{K-means problem}, where each column of $\bA$ is a data point, the goal is to determine a set of $K$ centroids $\bw_k, k\in\{1,2,\ldots,K\}$ (i.e., the columns of $\bW$) such that the sum of distances between each data point and its nearest centroid is minimized. 
This setup is equivalent to the low-rank matrix factorization problem where the second factor $\bZ$ must have exactly one nonzero entry per column, which is set to one, indicating the assignment of data points to their respective centroids: $\bZ\bZ^\top$ is diagonal and $\bZ\in\{0,1\}^{K\times N}$, where the diagonal values of $\bZ\bZ^\top$ indicate the number of data points associated with each cluster. The columns of $\bW$ then represent the cluster centroids \citep{zhang2017matrix, gillis2020nonnegative}.

\item \textbf{Sparsity.} 
The $\ell_0$ norm (the term is used even it is not a norm), defined as the number of nonzero elements in a vector, is a direct measure of sparsity. However, it is not commonly used for sparsity optimization due to its  non-convex, discontinuous, and non-differentiable nature.
Instead of the $\ell_0$ norm, the $\ell_1$ norm is often used as a convex relaxation of $\ell_0$. The $\ell_1$ norm promotes sparsity and is convex, which allows for the use of efficient optimization algorithms. Additionally, the $\ell_1$ norm is continuous and differentiable almost everywhere, making it more amenable to gradient descent and related techniques.
The $\ell_1$ matrix norm on $\bZ$ can induce sparsity among the elements of $\bZ$. 
Sparsity is useful, for example, in facial feature extraction, which leads to more localized features, meaning that fewer features are used to reconstruct each input image.
The \textit{spar} operator, introduced in \citet{hoyer2004non}, is continuous and promotes sparsity based on the ratio of $\ell_1$ to $\ell_2$ norm: for $\bx\neq\bzero\in\real^n$, $spar(\bx)=\frac{\sqrt{n}-\normone{\bx}/\normtwo{\bx}}{\sqrt{n}-1}\in[0,1]$.
$spar(\bx)=0$ if and only if $\normone{\bx}=\sqrt{n}\normtwo{\bx}$ and all entries of $\bx$ are equal. 
And $spar(\bx)=1$ if and only if $\normone{\bx}=\normtwo{\bx}$, in which case $\normzero{\bx}=1$.
The higher the value, the sparser; for example $spar([1,0,0]) > spar([1,1,1])$.
In the matrix factorization context, such a sparsity constraint can be applied by ensuring $spar(\bW[:,k]) \geq r_w$ and $spar(\bZ[k,:])\geq r_z$ for all $k\in\{1,2,\ldots,K\}$, where $r_w$ and $r_z$ are constants in the interval $[0,1]$ that impose a minimal sparsity level on the columns of $\bW$ and rows of $\bZ$.

\item \textbf{Group sparsity.}
When we believe that the columns of $\bZ$ can be grouped into subsets $\sG$, where each subset represents a coherent group of columns, the goal is to enforce sparsity at the group level rather than at the individual element level. Specifically, the union of all groups covers all columns of $\bZ$, i.e.,  $\cup_{\scriptsize g\in\sG}g=\{1,2,\ldots,N\}$. 
To achieve group sparsity, we can apply the  $\ell_1$ vector norm to a vector in $\real^{\abs{\sG}}$, where each component of this vector is the  $\ell_1$ or $\ell_2$ matrix norm of the columns within each group $g\in\sG$ (the $\ell_2$ matrix norm within each group does not promote sparsity). This means that instead of penalizing individual elements directly, we penalize the aggregated norms of the columns within each group.
This approach is particularly useful when the underlying structure of the data suggests that features or variables naturally cluster into meaningful groups, and we wish to select entire groups rather than individual features.

\index{Orthogonal matrix factorization}
\item \textbf{Orthogonality.} To promote minimal overlap between the features $\bW$ and/or the activations $\bZ$, orthogonality can be utilized. Specifically, one can enforce the columns of $\bW$ or the rows of $\bZ$ to be orthogonal through the penalty terms 
$\lambda_w\normf{\bW^\top \bW - \bI_K}^2$
and 
$\lambda_z\normf{\bZ \bZ^\top - \bI_K}^2$,
respectively. This regularization encourages the columns of $\bW$ and/or the rows of $\bZ$ to be as distinct as possible, which is the opposite goal of minimum-volume regularizers. When the constraint $\bZ\bZ^\top = \bI_K$ is enforced, matrix factorization effectively transforms into a \textit{soft clustering problem}: each data point is approximated by multiples of the columns of $\bW$.
When we instead impose an orthogonal constraint such that  $\bZ\bZ^\top=\bI_K$ rather than an orthogonal regularization, the problem becomes \textit{orthogonal matrix factorization}; see Problem~\ref{prob:ortho_mf}. In this case, $\bZ$ has orthonormal rows.


\item \textbf{Spatial smoothness.} When matrix factorization is used for feature extraction among a set of images, these images are typically first vectorized and stacked as the columns (resp. rows) of  $\bA$. Consequently, the spatial information inherent in the images is lost. To reintroduce this spatial coherence, regularizers can be employed. Since the columns of $\bW$ (resp. the rows of  $\bZ$) should represent basis images, they should exhibit some spatial consistency. To preserve the sharp edges present in the images, $\ell_1$-based regularizers are more suitable; for example, \textit{total variation} is often used and has the form:
$
\lambda_w \big(\sum_{k=1}^{K} \sum_{(i_1, i_2) \in \sS} \abs{w_{i_1,k}-w_{i_2,k}}\big),
$
where $\sS$ indicates the set of neighboring pixels in the image. See also the denoising least squares problem in Problem~\ref{prob:denoise_rls}.

\item  \textbf{Graph regularization.} Assume we want to preserve the geometric relationships among the input data points when they are projected into the low-dimensional subspace spanned by the columns of $\bW$. In other words, we aim to ensure that the distances between the columns of $\bZ$ closely match the distances between the columns of $\bA$. Specifically, if the Euclidean distance between two columns of $\bA$ is small, we expect the distance between the corresponding columns of $\bZ$ to be similarly small, and vice versa. This can be achieved using the regularizer
$\lambda_z \sum_{i,j} r_{ij} \normtwo{\bz_i-\bz_j}^2$,
where $r_{ij}$ is inversely proportional to the distance between the data points $i$ and $j$ (the $(i,j)$-th entry of a matrix $\bR$). The matrix $\bR$ can be computed in various ways; a common choice is $r_{ij} = \exp\{-\gamma\normtwo{\ba_i-\ba_j}\}$
for some parameter $\gamma > 0$. The entries in the matrix $\bR$ can be interpreted as weights in a graph connecting the data points, hence the term ``\textit{graph-regularized matrix factorization}".
Another application of this regularizer is in scenarios where partial label information is available. For example, in a dataset of facial images, we might know that certain subsets of images contain the same person. In such cases, one can construct $\bR$ as an $N\times N$ binary matrix where $r_{ij} = 1$ if and only if the data points $i$ and $j$ share the same label. Approaches that utilize this partial label information are known as \textit{semi-supervised matrix factorization}, in contrast to \textit{unsupervised matrix factorization}.
\end{itemize}

\index{Missing entries}
\index{Hadamard product}
\index{Netflix recommender}
\section{Missing Entries and Rank-One Update}\label{section:alt-columb-by-column}
Since  matrix decomposition via  ALS is extensively used in the context of Netflix recommender data, where a substantial number of entries are missing due to users not having watched certain movies or choosing not to rate them for various reasons.
In this scenario, the low-rank matrix decomposition problem is also known as \textit{matrix completion} that can help recover unobserved entries \citep{jain2017non}.
To address this, we can introduce an additional mask matrix $\bM\in \{0,1\}^{M\times N}$, where $m_{mn}\in \{0,1\}$ indicates whether  user $n$ has rated  movie $m$ or not. Therefore, the loss function can be defined as 
$$
L(\bW,\bZ) = \frac{1}{2}\normf{\bM\hadaprod  \bA- \bM\hadaprod (\bW\bZ)}^2,
$$
where $\hadaprod$ represents the \textit{Hadamard product} between matrices~\footnote{Note that in most textbooks, the Hadamard product is denoted by ``$\circ$";  however, to avoid conflict with the tensor outer product discussed in Part~\ref{part:tensor_decom}, we use ``$\hadaprod$" here.}. For example, the Hadamard product of a $3 \times 3$ matrix $\bA$ and a $3\times 3$ matrix $\bB$ is
$$
\bA\hadaprod \bB = 
\begin{bmatrix}
	a_{11} & a_{12} & a_{13} \\
	a_{21} & a_{22} & a_{23} \\
	a_{31} & a_{32} & a_{33}
\end{bmatrix}
\hadaprod
\begin{bmatrix}
	b_{11} & b_{12} & b_{13} \\
	b_{21} & b_{22} & b_{23} \\
	b_{31} & b_{32} & b_{33}
\end{bmatrix}
=
\begin{bmatrix}
	a_{11}b_{11} &  a_{12}b_{12} & a_{13}b_{13} \\
	a_{21}b_{21} & a_{22}b_{22} & a_{23}b_{23} \\
	a_{31}b_{31} & a_{32}b_{32} & a_{33}b_{33}
\end{bmatrix}.
$$
The above formulation concisely expresses our goal of finding a completion of the ratings matrix that is both of low rank and consistent with observed user ratings.
To find the solution to this problem, we decompose the updates in Equation~\eqref{equation:als-regular-final-all} into:
\begin{equation}\label{equation:als-ori-all-wz}
	\left\{
	\begin{aligned}
		\bz_n &= (\bW^\top\bW+ \lambda_z\bI)^{-1} \bW^\top \ba_n, &\gap& \text{for $n\in \{1,2,\ldots, N\}$}  ;\\
		\bw_m &= (\bZ\bZ^\top+\lambda_w\bI)^{-1}\bZ\bb_m,  &\gap& \text{for $m\in \{1,2,\ldots, M\}$} ,
	\end{aligned}
	\right.
\end{equation}
where $\bZ=[\bz_1, \bz_2, \ldots, \bz_N]$ and $\bA=[\ba_1,\ba_2, \ldots, \ba_N]$ represent the column partitions of $\bZ$ and $\bA$, respectively. Similarly, $\bW^\top=[\bw_1, \bw_2, \ldots, \bw_M]$ and $\bA^\top=[\bb_1,\bb_2, \ldots, \bb_M]$ are the column partitions of $\bW^\top$ and $\bA^\top$, respectively. This decomposition of the updates indicates the updates can be performed in a column-by-column fashion (the rank-one updates).

\paragraph{Given $\bW$.}
Let $\bo_n\in \{0,1\}^M$ represent the movies rated by user $n$, where $o_{nm}=1$ if user $n$ has rated movie $m$, and $o_{nm}=0$ otherwise. Then the $n$-th column of $\bA$ without missing entries can be denoted using the Matlab-style notation as $\ba_n[\bo_n]$. 
And we want to approximate the existing entries of the $n$-th column by $\ba_n[\bo_n] \approx \bW[\bo_n, :]\bz_n$, which is indeed a rank-one least squares problem:
\begin{equation}\label{equation:als-ori-all-wz-modif-z}
\begin{aligned}
\bz_n &= \left(\bW[\bo_n, :]^\top\bW[\bo_n, :]+ \lambda_z\bI\right)^{-1} \bW[\bo_n, :]^\top \ba_n[\bo_n], \quad \text{for $n\in \{1,2,\ldots, N\}$} .
\end{aligned}
\end{equation}
Moreover, the loss function with respect to $\bz_n$ and $\bZ$  can be described, respectively, by
$$
\begin{aligned}
L(\bz_n|\bW) &=\sum_{m\in \bo_n} \left(a_{mn} - \bw_m^\top\bz_n\right)^2
\gap \text{and}\gap
L(\bZ|\bW) =\sum_{n=1}^N\ \sum_{m\in \bo_n} \left(a_{mn} - \bw_m^\top\bz_n\right)^2.
\end{aligned}
$$

\paragraph{Given $\bZ$.}
Similarly, if $\bp_m \in\{0,1\}^{N}$ denotes the users who have rated  movie $m$, with $p_{mn}=1$ if  movie $m$ has been rated by user $n$, and $p_{mn}=0$ otherwise. Then the $m$-th row of $\bA$ without missing entries can be denoted by the Matlab-style notation as $\bb_m[\bp_m]$. We want to approximate the existing entries of the $m$-th row by $\bb_m[\bp_m] \approx \bZ[:, \bp_m]^\top\bw_m$, 
\footnote{Note that $\bZ[:, \bp_m]^\top$ is the transpose of $\bZ[:, \bp_m]$, which is equal to $\bZ^\top[\bp_m,:]$, i.e., transposing first and then selecting.}
which  is again a rank-one least squares problem:
\begin{equation}\label{equation:als-ori-all-wz-modif-w}
\begin{aligned}
\bw_m &= (\bZ[:, \bp_m]\bZ[:, \bp_m]^\top+\lambda_w\bI)^{-1}\bZ[:, \bp_m]\bb_m[\bp_m],  \quad \text{for $m\in \{1,2,\ldots, M\}$} .
\end{aligned}
\end{equation}
Similarly, the loss function with respect to $\bw_m$ and $\bW$ can be described, respectively,  by
$$
\begin{aligned}
L(\bw_m|\bZ) &=\sum_{n\in \bp_m} \left(a_{mn} - \bw_m^\top\bz_n\right)^2 
\gap \text{and}\gap
L(\bW|\bZ) =\sum_{m=1}^M  \sum_{n\in \bp_m} \left(a_{mn} - \bw_m^\top\bz_n\right)^2 .
\end{aligned}
$$
The procedure is once again presented in Algorithm~\ref{alg:als-regularizer-missing-entries}.
Other approaches, such as \textit{singular value projection (SVP)}, also exist to address the matrix completion problem. At a high level, SVP is a type of projected gradient descent (PGD) method that updates iteratively via gradient descent, projecting the updated matrix into a low-rank form through singular value decomposition at each step. However, the alternating least squares approach generally  outperforms SVP in the context of matrix completion, so we will not delve into SVP here. For more details, refer to \citet{jain2017non}.

\begin{algorithm}[h] 
\caption{Alternating Least Squares with Missing Entries and Regularization}
\label{alg:als-regularizer-missing-entries}
\begin{algorithmic}[1] 
\Require Matrix $\bA\in \real^{M\times N}$;
\State Initialize $\bW\in \real^{M\times K}$, $\bZ\in \real^{K\times N}$ \textcolor{mylightbluetext}{randomly without condition on the rank and the relationship between $M, N, K$}; 
\State Choose a stoping criterion on the approximation error $\delta$;
\State Choose regularization parameters $\lambda_w, \lambda_z$;
\State Compute the mask matrix $\bM$ from $\bA$;
\State Choose the maximum number of iterations $C$;
\State $iter=0$; \Comment{Count for the number of iterations}
\While{\textcolor{mylightbluetext}{$\normf{\bM\hadaprod  \bA- \bM\hadaprod (\bW\bZ)}^2>\delta $} and $iter<C$}
\State $iter=iter+1$; 
\For{$n=1,2,\ldots, N$}
\State $\bz_n \leftarrow \left(\bW[\bo_n, :]^\top\bW[\bo_n, :]+ \lambda_z\bI\right)^{-1} \bW[\bo_n, :]^\top \ba_n[\bo_n]$; \Comment{$n$-th column of $\bZ$}
\EndFor

\For{$m=1,2,\ldots, M$}
\State $\bw_m \leftarrow (\bZ[:, \bp_m]\bZ[:, \bp_m]^\top+\lambda_w\bI)^{-1}\bZ[:, \bp_m]\bb_m[\bp_m]$;\Comment{$m$-th column of $\bW^\top$}
\EndFor
\EndWhile
\State Output $\bW^\top=[\bw_1, \bw_2, \ldots, \bw_M],\bZ=[\bz_1, \bz_2, \ldots, \bz_N]$;
\end{algorithmic} 
\end{algorithm}

\index{Hidden features}
\index{Inner product}
\section{Vector Inner Product and Hidden Vectors}\label{section:als-vector-product}
We have observed that the ALS algorithm seeks to find lower-dimensional matrices $\bW$ and $\bZ$ such that their product $\bW\bZ$ can approximate $\bA\approx \bW\bZ$ in terms of the  squared loss:
$
\mathop{\min}_{\bW,\bZ}  \sum_{n=1}^N \sum_{m=1}^{M} \left(a_{mn} - \bw_m^\top\bz_n\right)^2,
$
that is, each entry $a_{mn}$ in $\bA$ can be approximated as the inner product of  two vectors $\bw_m^\top\bz_n$. The geometric definition of the vector inner product is given by 
$$
\bw_m^\top\bz_n = \normtwo{\bw_m}\cdot \normtwo{\bz_n} \cos \theta,
$$
where $\theta$ represents the angle between vectors $\bw_m$ and $\bz_n$ (Definition~\ref{definition:angle_bet_vec_ineq}). Thus, if the vector norms of $\bw_m$ and $\bz_n$ are determined, the smaller the angle, the larger the inner product.

In the context of Netflix data,   movie ratings range from 0 to 5,
with higher ratings indicating a stronger user preference for the movie. 
If $\bw_m$ and $\bz_n$ fall sufficiently ``close," the value $\bw_m^\top\bz_n$ becomes larger. 
This concept elucidates the essence of ALS, where $\bw_m$ represents the features or attributes of movie $m$, while $\bz_n$ encapsulates  the features or preferences of user $n$. 
In other words,  ALS associates each user with a \textit{latent vector of preference} and each movie with a \textit{latent vector of attributes}.
Furthermore, each element in $\bw_m$ and $\bz_n$ signifies a specific feature. For example, it could be that the second feature $w_{m2}$ ($w_{m2}$ denotes the second element of vector $\bw_{m}$) represents whether the movie is an action movie or not, and $z_{n2}$ might denote whether  user $n$ has a preference for action movies. When this holds true, then $\bw_m^\top\bz_n$ becomes large and provides a good approximation of  $a_{mn}$.

In the decomposition $\bA\approx \bW\bZ$, it is established that the rows of $\bW$ contain the hidden features of the movies, and the columns of $\bZ$ contain the hidden features of the users. 
Nevertheless, the explicit meanings of the rows in $\bW$ or the columns in $\bZ$ remain undisclosed.
Although they might correspond to categories or genres of the movies, fostering underlying connections between users and movies, their precise nature remains uncertain.
It is precisely this ambiguity that gives rise to the terminology ``latent" or ``hidden."

\index{Stochastic gradient descent}
\index{Gradient descent}
\index{Matrix inverse}
\index{LU decomposition}
\section{Gradient Descent}\label{section:als-gradie-descent}
In Algorithm~\ref{alg:als}, \ref{alg:als-regularizer}, and \ref{alg:als-regularizer-missing-entries}, we reduce the loss through the inversion of matrices (e.g., using LU decomposition as described in Theorem~\ref{theorem:inverse-by-lu2}). 
The reality, however, is frequently far from straightforward, particularly in the big data era of today. As data volumes explode, the size of the inversion matrix will grow at a pace proportional to the cube of the number of samples,  which poses a great challenge to the storage and computational resources.
On the other hand, this leads to the creation of an ongoing development of the gradient-based optimization technique.
The \textit{gradient descent (GD)} method and its variant, the \textit{stochastic gradient descent (SGD)} method, are among them the simplest, fastest, and most efficient methods \citep{lu2022gradient}. Convex loss function optimization problems are frequently solved using this type of approach. We now go into more details about its principle.

In Equation~\eqref{equation:als-ori-all-wz}, we derive the column-by-column update rules directly from the full matrix approach outlined in Equation~\eqref{equation:als-regular-final-all} (with regularization taken into account). 
To understand the underlying concept, consider the loss function with regularization, as given by Equation~\eqref{equation:als-regularion-full-matrix}.
When minimizing the  loss in \eqref{equation:als-regularion-full-matrix} with respect to $\bz_n$, we can break down the loss as follows:
\begin{equation}\label{als:gradient-regularization-zn}
\footnotesize
\begin{aligned}
L(\bz_n)  &=\frac{1}{2}\normf{\bW\bZ-\bA}^2 +\frac{1}{2}\lambda_w \normf{\bW}^2 + \frac{1}{2}\lambda_z \normf{\bZ}^2
= \frac{1}{2}\normtwo{\bW\bz_n-\ba_n}^2 + \frac{1}{2}\lambda_z \normtwo{\bz_n}^2 + C_{z_n},
\end{aligned}
\end{equation}
where $C_{z_n}$ is a constant with respect to $\bz_n$, and $\bZ=[\bz_1, \bz_2, \ldots, \bz_N]$ and $\bA=[\ba_1,\ba_2, \ldots, \ba_N]$ represent the column partitions of $\bZ$ and $\bA$, respectively. 
The gradient and the root are given, respectively, by 
$$
\begin{aligned}
\nabla_{\bz_n} L(\bz_n) = \bW^\top\bW\bz_n - \bW^\top\ba_n + \lambda_z\bz_n
\,\,\implies \,\,
\bz_n = (\bW^\top\bW+ \lambda_z\bI)^{-1} \bW^\top \ba_n, \,\,  \forall\,n.
\end{aligned}
$$
This solution corresponds to the first update rule in the column-wise update in Equation~\eqref{equation:als-ori-all-wz}.
Similarly, when minimizing the loss with respect to $\bw_m$, we have:
\begin{equation}\label{als:gradient-regularization-wd}
\footnotesize
\begin{aligned}
L(\bw_m )  
&
=\frac{1}{2}\normf{\bZ^\top\bW-\bA^\top}^2 +\frac{1}{2}\lambda_w \normf{\bW^\top}^2 + \frac{1}{2}\lambda_z \normf{\bZ}^2
= \frac{1}{2}\normtwo{\bZ^\top\bw_m-\bb_n}^2 + \frac{1}{2}\lambda_w \normtwo{\bw_m}^2 + C_{w_m},
\end{aligned}
\end{equation}
where $C_{w_m}$ is a constant with respect to $\bw_m$, and $\bW^\top=[\bw_1, \bw_2, \ldots, \bw_M]$ and $\bA^\top=[\bb_1,\bb_2, \ldots,$ $\bb_M]$ represent the column partitions of $\bW^\top$ and $\bA^\top$, respectively. 
Analogously, taking the gradient with respect to $\bw_m$, it follows that
$$
\begin{aligned}
\nabla_{\bw_m} L(\bw_m) = \bZ\bZ^\top\bw_m - \bZ\bb_n + \lambda_w\bw_m
\,\,\implies\,\,
\bw_m = (\bZ\bZ^\top+\lambda_w\bI)^{-1}\bZ\bb_m, \,\, \forall \, m.
\end{aligned}
$$
This solution corresponds to the second update rule in the column-wise update in Equation~\eqref{equation:als-ori-all-wz}:

Now suppose we express the iteration number ($t=1,2,\ldots$)  as the superscript, and we want to find the updates $\{\bz^{(t+1)}_n, \bw^{(t+1)}_m\}$  in the $(t+1)$-th iteration  base on $\{\bZ^{(t)}, \bW^{(t)}\}$  in the $t$-th iteration:
$$
\left.
\begin{aligned}
\bz^{(t+1)}_n    &\leftarrow \mathop{\arg \min}_{\bz_n^{(t)}} L(\bz_n^{(t)})
\qquad\text{and}\qquad
\bw_m^{(t+1)}    \leftarrow \mathop{\arg\min}_{\bw_m^{(t)}} L(\bw_m^{(t)}).
\end{aligned}
\right.
$$
For simplicity, we will only derive for $\bz^{(t+1)}_n    \leftarrow \mathop{\arg \min}_{\bz_n^{(t)}} L(\bz_n^{(t)})$, and the derivation for the update on $\bw_m^{(t+1)}$ will follow a similar process.

\index{Linear approximation}
\index{Linear update}
\index{Greedy search}
\index{Gradient descent}
\paragraph{Approximation by linear update.} 
Suppose we want to approximate $\bz^{(t+1)}_n$ using a \textit{linear update} on $\bz^{(t)}_n$:
$$
\textbf{(Linear Update)}: \qquad {\bz^{(t+1)}_n = \bz^{(t)}_n + \eta \bv.}
$$
The problem now reduces to finding the solution of $\bv$ such that
$$
\bv=\mathop{\arg \min}_{\bv} L(\bz^{(t)}_n + \eta \bv) .
$$
By Taylor's formula (Appendix~\ref{appendix:taylor-expansion}), $L(\bz^{(t)}_n + \eta \bv)$ can be approximated by 
$$
L(\bz^{(t)}_n + \eta \bv) \approx L(\bz^{(t)}_n ) + \eta \bv^\top \nabla  L(\bz^{(t)}_n ),
$$
where $\eta$ is a small value, and $\nabla  L(\bz^{(t)}_n )$ represents the gradient of $L(\bz)$ at $\bz^{(t)}_n$. 
To find $\bv$ under the constraint $\normtwo{\bv}=1$ for a positive $\eta$, we perform the following minimization:
$$
\bv=\mathop{\argmin}_{\normtwo{\bv}=1} L(\bz^{(t)}_n + \eta \bv) 
\approx\mathop{\argmin}_{\normtwo{\bv}=1}
 \left\{L(\bz^{(t)}_n ) + \eta \bv^\top \nabla  L(\bz^{(t)}_n )\right\}.
$$
This approach is known as the \textit{greedy search}. The optimal $\bv$ can be obtained by 
$$
\bv = -\nabla L(\bz^{(t)}_n )\big/{\big\Vert{\nabla L(\bz^{(t)}_n )}\big\Vert_2},
$$
which means $\bv$ points in the opposite direction of $\nabla L(\bz^{(t)}_n )$. Therefore, it is reasonable to update  $\bz_n^{(t+1)}$ as follows:
$$
\bz^{(t+1)}_n =\bz^{(t)}_n + \eta \bv = \bz^{(t)}_n - \eta {\nabla L(\bz^{(t)}_n )}\big/{\big\Vert{\nabla L(\bz^{(t)}_n )}\big\Vert_2},
$$
which is commonly referred to as   the \textit{gradient descent} (GD). Similarly, the gradient descent for $\bw_m^{(t+1)}$ is given by
$$
\bw^{(t+1)}_m =\bw^{(t)}_m + \eta \bv = \bw^{(t)}_m - \eta {\nabla L(\bw^{(t)}_m )}\big/{\big\Vert{\nabla L(\bw^{(t)}_m )}\big\Vert_2}.
$$
The revised procedure for Algorithm~\ref{alg:als-regularizer} employing a gradient descent approach is presented in Algorithm~\ref{alg:als-regularizer-missing-stochas-gradient}.

It's noteworthy that the ALS without GD (Algorithm~\ref{alg:als-regularizer})  lacks explicit parameters like step size. 
This characteristic can be both advantageous and disadvantageous. 
On one hand, it absolves the  user from the time-consuming task of fine-tuning parameters, making the method more accessible and less demanding. 
On the other hand, this absence of adjustable parameters also restricts the user's control to directly influence the progression of the algorithm, leaving the convergence of ALS entirely contingent upon the inherent structure of the optimization problem at hand.

In practical applications, it is customary to alternate between the pure ALS iterations outlined in Algorithm~\ref{alg:als-regularizer} and the modified, gradient-descent variants  mentioned in this section. These descent adaptations offer the user a degree of control through a tunable step length parameter, allowing for a more customized approach to the optimization process.

\begin{algorithm}[h] 
\caption{Alternating Least Squares with Full Entries and Gradient Descent}
\label{alg:als-regularizer-missing-stochas-gradient}
\begin{algorithmic}[1] 
\Require Matrix $\bA\in \real^{M\times N}$;
\State Initialize $\bW\in \real^{M\times K}$, $\bZ\in \real^{K\times N}$ \textcolor{mylightbluetext}{randomly without condition on the rank and the relationship between $M, N, K$}; 
\State Choose a stoping criterion on the approximation error $\delta$;
\State Choose regularization parameters $\lambda_w, \lambda_z$, and step sizes $\eta_w, \eta_z$;
\State Choose the maximum number of iterations $C$;
\State $iter=0$; \Comment{Count for the number of iterations}
\While{$\normf{\bA- (\bW\bZ)}^2>\delta $ and $iter<C$}
\State $iter=iter+1$; 
\For{$n=1,2,\ldots, N$}
\State $\bz^{(t+1)}_n \leftarrow\bz^{(t)}_n - \eta_z {\nabla L(\bz^{(t)}_n )}\big/{\big\Vert{\nabla L(\bz^{(t)}_n )}\big\Vert_2}$; \Comment{$n$-th column of $\bZ$}
\EndFor

\For{$m=1,2,\ldots, M$}
\State $\bw^{(t+1)}_m  \leftarrow \bw^{(t)}_m - \eta_w {\nabla L(\bw^{(t)}_m )}\big/{\big\Vert{\nabla L(\bw^{(t)}_m )}\big\Vert_2}$;\Comment{$m$-th column of $\bW^\top$}
\EndFor
\EndWhile
\State Output $\bW^\top=[\bw_1, \bw_2, \ldots, \bw_M],\bZ=[\bz_1, \bz_2, \ldots, \bz_N]$;
\end{algorithmic} 
\end{algorithm}

\begin{figure}[h]
\centering  
\vspace{-0.15cm}    
\subfigtopskip=2pt  
\subfigbottomskip=2pt 
\subfigcapskip=-5pt  
\subfigure[A two-dimensional convex function $L(\bz)$.]{\label{fig:alsgd1}
\includegraphics[width=0.47\linewidth]{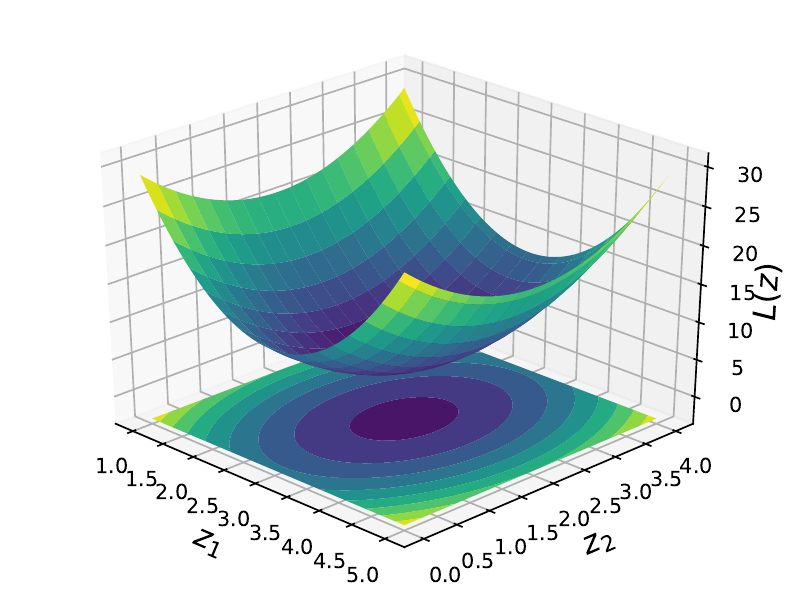}}
\subfigure[$L(\bz)=c$ is a constant.]{\label{fig:alsgd2}
\includegraphics[width=0.44\linewidth]{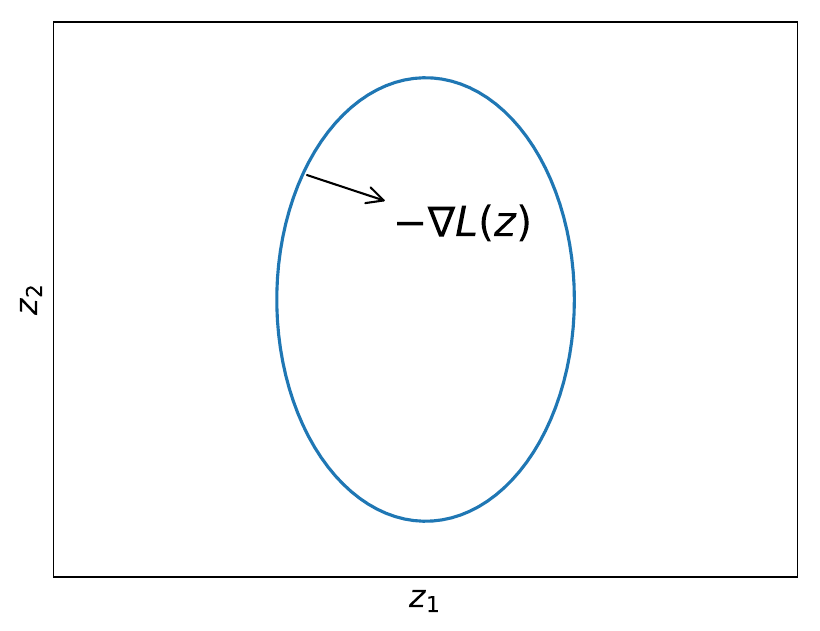}}
\caption{Figure~\ref{fig:alsgd1} shows  surface and  contour plots for a specific function (\textcolor{mydarkblue}{blue}=low, \textcolor{mydarkyellow}{yellow}=high), where the upper graph is the surface plot, and the lower one is its projection  (i.e., contour). Figure~\ref{fig:alsgd2}: $-\nabla L(\bz)$ pushes the loss to decrease for the convex function $L(\bz)$.}
\label{fig:alsgd-geometrical}
\end{figure}

\index{Level curves}
\index{Level surfaces}
\paragraph{Geometrical interpretation of gradient descent.} 
\begin{lemma}[Direction of Gradients]\label{lemm:direction-gradients}
Gradients of variables given a loss function are perpendicular to the level curves (a.k.a., level surfaces).
\end{lemma}
\begin{proof}[of Lemma~\ref{lemm:direction-gradients}, the informal proof]
This proof involves showing that the gradient is orthogonal to the tangent of the level curve.  For simplicity, let's start with the two-dimensional case.
Suppose the level curve takes the form $f(x,y)=c$. 
This implicitly establishes a relationship between $x$ and $y$ such that $y=y(x)$, where $y$ can be regarded as a function of $x$~\footnote{This is known as the implicit function  theorem, under the conditions of the nonzero partial derivative and smoothness.}. Therefore, the level curve can be expressed as 
$
f(x, y(x)) = c.
$
Applying the chain rule, we get:
$$
\frac{\partial f}{\partial x} \underbrace{\frac{dx}{dx}}_{=1} + \frac{\partial f}{\partial y} \frac{dy}{dx}=0
\gap \implies \gap 
\left\langle \frac{\partial f}{\partial x}, \frac{\partial f}{\partial y}\right\rangle
\cdot 
\left\langle \frac{dx}{dx}, \frac{dy}{dx}\right\rangle=0.
$$
That is, the gradient is perpendicular to the tangent.

In full generality, consider the level curve of a vector $\bx\in \real^n$: $f(\bx) = f(x_1, x_2, \ldots, x_n)=c$. Each variable $x_i$ can be regarded as a function of a parameter $t$ on the level curve $f(\bx)=c$: $f(x_1(t), x_2(t), \ldots, x_n(t))=c$. Differentiating the equation with respect to $t$ using the chain rule:
$$
\frac{\partial f}{\partial x_1} \frac{dx_1}{dt} + \frac{\partial f}{\partial x_2} \frac{dx_2}{dt}
+\ldots + \frac{\partial f}{\partial x_n} \frac{dx_n}{dt}
=0.
$$
Thus, the gradient is perpendicular to the tangent in the $n$-dimensional case:
$$
\left\langle \frac{\partial f}{\partial x_1}, \frac{\partial f}{\partial x_2}, \ldots, \frac{\partial f}{\partial x_n}\right\rangle
\cdot 
\left\langle \frac{dx_1}{dt}, \frac{dx_2}{dt}, \ldots \frac{dx_n}{dt}\right\rangle=0.
$$
This completes the proof.
\end{proof}
The lemma  reveals the geometrical interpretation of gradient descent. When seeking a solution to minimize a convex function $L(\bz)$, gradient descent proceeds in the direction opposite to the gradient, which leads to a reduction in the loss. Figure~\ref{fig:alsgd-geometrical} depicts a two-dimensional case, where $-\nabla L(\bz)$ pushes the loss to decrease for the convex function $L(\bz)$. 
\index{Convex function}
\index{Contour plot}
\index{Contour plot}
\index{Regularization}

\index{Geometrical interpretation}
\index{Projection gradient descent}
\index{Overfitting}
\section{Regularization: A Geometrical Interpretation}\label{section:reg_geom_inter}
\begin{figure}[h]
\centering
\includegraphics[width=0.95\textwidth]{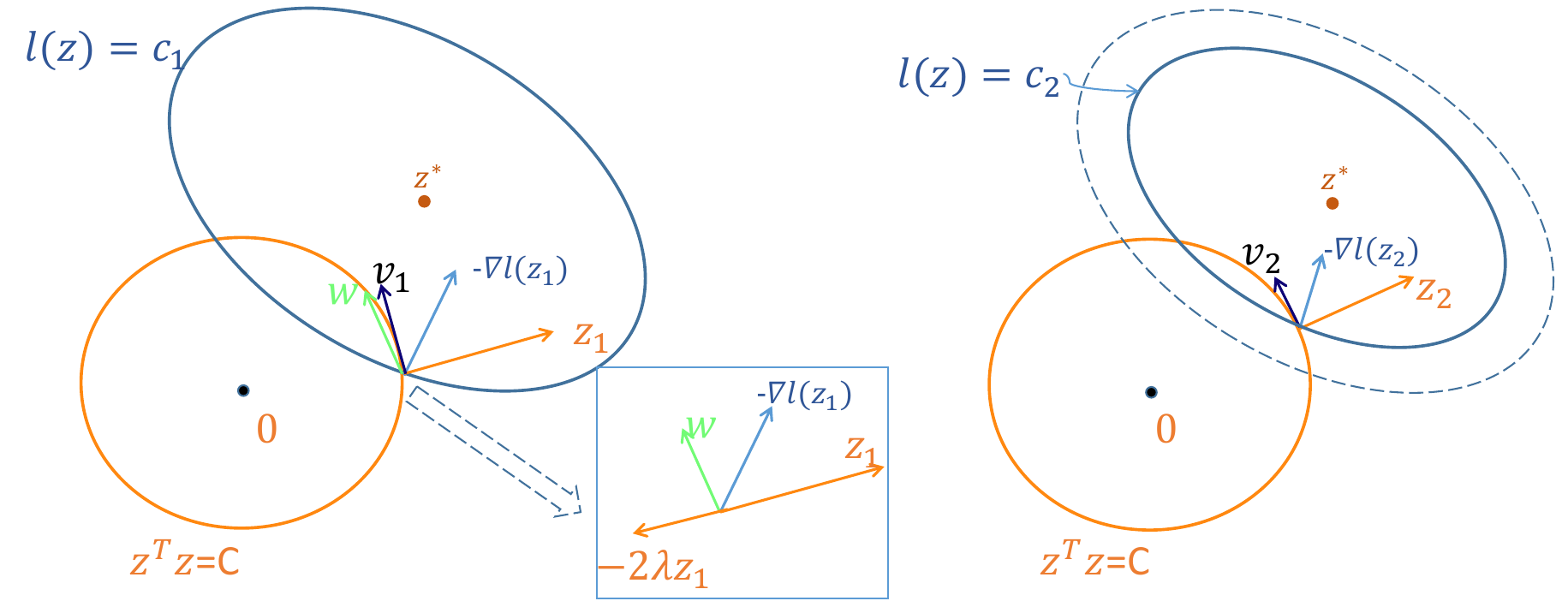}
\caption{Constrained gradient descent with $\bz^\top\bz\leq C$. The \textcolor{mydarkgreen}{green} vector $\bw$ represents the projection of $\bv_1$ onto the set $\bz^\top\bz\leq C$, where $\bv_1$ is the component of $-\nabla l(\bz)$ that is perpendicular to $\bz_1$. 
The image on the right illustrates the next step after the update in the left picture. $\bz^\star$ denotes the optimal solution of \{$\min l(\bz)$\}.}
\label{fig:alsgd3}
\end{figure}
In Section~\ref{section:regularization-extention-general}, we discussed how  regularization can extend the ALS algorithm to general matrices.
Gradient descent can provide a geometric interpretation of regularization. 
To avoid confusion, we denote the loss function without regularization as $l(\bz)$ and the loss function with regularization as $L(\bz) = \l(\bz)+\lambda_z \normtwo{\bz}^2$, where $l(\bz): \real^n \rightarrow \real$. When minimizing $l(\bz)$, a descent method searches for a solution in $\real^n$. However, in machine learning, searching across the entire space $\real^n$ can lead to overfitting. 
A partial solution is to search within a subset of the vector space, e.g., searching in $\bz^\top\bz < C$ for some constant $C$. 
This can be formulated as:
$$
\mathop{\arg\min}_{\bz} \,\, l(\bz), \gap \text{s.t.,} \gap \bz^\top\bz\leq C.
$$
As demonstrated above, a basic gradient descent method will go further in the direction of $-\nabla l(\bz)$, i.e., update $\bz$ as $\bz\leftarrow \bz-\eta \nabla l(\bz)$ for a small step size $\eta$. When the level curve is $l(\bz)=c_1$ and the current position of parameter $\bz$ is $\bz=\bz_1$, where $\bz_1$ is the intersection of $\bz^\top\bz=C$ and $l(\bz)=c_1$, the descent direction $-\nabla l(\bz_1)$ will be perpendicular to the level curve of $l(\bz_1)=c_1$, as shown in the left picture of Figure~\ref{fig:alsgd3} (by Lemma~\ref{lemm:direction-gradients}). However, if we further restrict that the optimal value must lie  within $\bz^\top\bz\leq C$, the trivial descent direction $-\nabla l(\bz_1)$ will lead the update $\bz_2=\bz_1-\eta \nabla l(\bz_1)$ beyond the boundary of $\bz^\top\bz\leq C$. 
One solution is to decompose the step $-\nabla l(\bz_1)$ into 
$$
-\nabla l(\bz_1) = a\bz_1 + \bv_1,
$$ 
where $a\bz_1$ represents the component perpendicular to the curve of $\bz^\top\bz=C$, and $\bv_1$ is the component parallel to the curve of $\bz^\top\bz=C$. By keeping only the step $\bv_1$, the update becomes
$$
\bz_2 = \text{project}(\bz_1+\eta \bv_1) = \text{project}\bigg(\bz_1 + \eta
\underbrace{(-\nabla l(\bz_1) -a\bz_1)}_{\bv_1}\bigg),~\footnote{where the operation project($\bx$) 
will project the vector $\bx$ to the closest point inside $\bz^\top\bz\leq C$. Notice here the unprojected update $\bz_2 = \bz_1+\eta \bv_1$ can still make $\bz_2$ fall outside the curve of $\bz^\top\bz\leq C$.}
$$ 
which will lead to a smaller loss from $l(\bz_1)$ to $l(\bz_2)$ while maintaining the constraint $\bz^\top\bz\leq C$. This approach is known as  \textit{projection gradient descent (PGD)}. 
It is not hard to see that the update $\bz_2 = \text{project}(\bz_1+\eta \bv_1)$ can be understood as finding a vector $\bw$ (represented by the \textcolor{mydarkgreen}{green} vector in the left picture of Figure~\ref{fig:alsgd3}) such that $\bz_2=\bz_1+\bw$ lies inside the curve of $\bz^\top\bz\leq C$. Mathematically, the vector $\bw$ can be determined as $-\nabla l(\bz_1) -2\lambda \bz_1$ for some $\lambda$, as shown in the middle picture of Figure~\ref{fig:alsgd3}. This is precisely the negative gradient of $L(\bz)=l(\bz)+\lambda\normtwo{\bz}^2$ such that 
$$
\begin{aligned}
\bw=-\nabla L(\bz) &= -\nabla l(\bz) - 2\lambda \bz 
\gap \implies \gap
\bz_2 = \bz_1+ \bw =\bz_1 -  \nabla L(\bz).
\end{aligned}
$$
And in practice, a small step size $\eta$ prevents the trajectory from moving outside the curve of $\bz^\top\bz\leq C$:
$$
\bz_2  =\bz_1 -  \eta\nabla L(\bz),
$$
which aligns with the regularization term discussed in Section~\ref{section:regularization-extention-general}.

\index{Sparsity}
\index{$\ell_1$ norm}
\begin{figure}[h]
\centering
\includegraphics[width=0.95\textwidth]{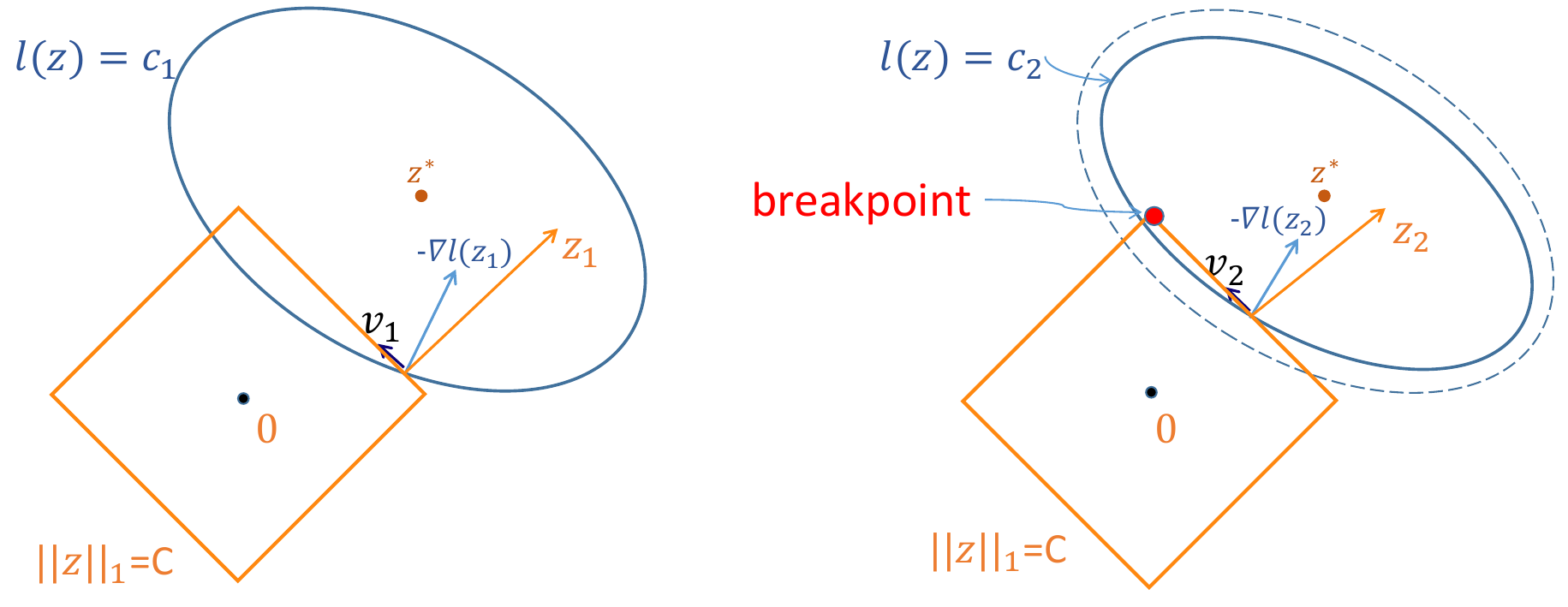}
\caption{Constrained gradient descent with $\norm{\bz}_1\leq C$, where the \textcolor{winestain}{red} dot denotes the breakpoint in the $\ell_1$ norm. The right picture illustrates the next step after the update in the left picture. $\bz^\star$ denotes the optimal solution of \{$\min l(\bz)$\}.}
\label{fig:alsgd4}
\end{figure}
\paragraph{Sparsity.}
In certain scenarios, we seek to identify a sparse solution $\bz$ such that $l(\bz)$ is minimized. 
For example, in facial feature extraction, sparsity leads to more localized features, meaning that fewer features are used to reconstruct each input image.
Regularization to be constrained in $\norm{\bz}_1 \leq C$ exists to this purpose, where $\norm{\cdot}_1$ is the $\ell_1$ norm of a vector or a matrix. 
Visualizations of the $\ell_1$ norm in two- and three-dimensional spaces are shown in Figure~\ref{fig:p-norm-2d} and \ref{fig:p-norm-comparison-3d}, respectively. Similar to the previous case, the $\ell_1$ constrained optimization pushes the gradient descent towards the border of the level set $\norm{\bz}_1=C$. The situation in the two-dimensional case is shown in Figure~\ref{fig:alsgd4}. 
In  high-dimensional cases, many elements in $\bz$ will be pushed towards the breakpoint of $\norm{\bz}_1=C$ so that the solution favors sparsity, as shown in the right picture of Figure~\ref{fig:alsgd4}.

\index{Stochastic gradient descent}
\index{Stochastic coordinate descent}
\section{Stochastic Gradient Descent}
The gradient descent method is a valuable optimization algorithm; however, it exhibits certain limitations in practical applications. 
To comprehend the issues associated with the gradient descent method, we consider the mean squared error (MSE) derived from  Equation~\eqref{equation:als-per-example-loss2}:
\begin{equation}\label{equation:als-per-example-loss_mse}
\frac{1}{MN}\mathop{\min}_{\bW,\bZ}  \sum_{n=1}^N \sum_{m=1}^{M} \left(a_{mn} - \bw_m^\top\bz_n\right)^2.
\end{equation}
The MSE requires calculating the residual $e_{mn} = (a_{mn} - \bw_m^\top\bz_n)^2$ for each observed entry $a_{mn}$, representing the squared difference between predicted and actual values.  The total sum of residual squares is denoted by $e = \sum_{m,n=1}^{MN}e_{mn}$.
In cases with a substantial number of training entries (i.e., large $MN$), the entire computation process becomes notably slow. 
Additionally, the gradients from different input samples may cancel out, resulting in small changes in the final update.
To address these challenges, researchers have enhanced the gradient descent method with the \textit{stochastic gradient descent (SGD)} method (see, for example, \citet{lu2022gradient}). 
In the SGD algorithm, instead of calculating the full gradient of the objective function with respect to the parameters  across all samples in the data set, which can be computationally expensive, the algorithm takes a more efficient approach. It randomly chooses one sample and calculates the gradient of the objective function with respect to the parameters using only  this single sample. 
This gradient estimate is then used to update the parameters in the direction that minimizes the objective function. 
By using a single sample at each iteration, the SGD algorithm provides a fast and often sufficient approximation of the full gradient, making it particularly useful for large data sets.

In particular, we consider again the per-example loss:
$$
L(\bW,\bZ)= \frac{1}{2} \sum_{n=1}^N \sum_{m=1}^{M} \left(a_{mn} - \bw_m^\top\bz_n\right)^2 +\frac{1}{2} \lambda_w\sum_{m=1}^{M}\normtwo{\bw_m}^2 +\frac{1}{2}\lambda_z\sum_{n=1}^{N}\normtwo{\bz_n}^2.
$$
As we iteratively reduce the  loss term $l(\bw_m, \bz_n)=\frac{1}{2}\left(a_{mn} - \bw_m^\top\bz_n\right)^2+\frac{1}{2} \lambda_w\normtwo{\bw_m}^2 +\frac{1}{2}\lambda_z\normtwo{\bz_n}^2$ for all $m\in \{1,2,\ldots, M\}, n\in\{1,2,\ldots,N\}$ (referred to as the per-example loss term), the overall loss $L(\bW,\bZ)$ decreases accordingly.
This process is also known as  \textit{stochastic coordinate descent}. The gradients  with respect to $\bw_m$ and $\bz_n$, and their roots are given, respectively, by 
$$
\left\{
\begin{aligned}
\nabla_{\bz_n} l(\bz_n) &= \bw_m\bw_m^\top \bz_n +\lambda_z\bz_n  -a_{mn} \bw_m 
 &\implies&\,\, \bz_n= a_{mn}(\bw_m\bw_m^\top+\lambda_z\bI)^{-1}\bw_m;\\
\nabla_{\bw_m} l(\bw_m) &= \bz_n\bz_n^\top\bw_m +\lambda_w\bw_m - a_{mn}\bz_n &\implies&\,\, \bw_m= a_{mn}(\bz_n\bz_n^\top+\lambda_w\bI)^{-1}\bz_n.
\end{aligned}
\right.
$$
Alternatively, the update can be performed using gradient descent. Since we update based on the per-example loss, this approach is also known as the \textit{stochastic gradient descent (SGD)}:
$$
\left.
\begin{aligned}
\bz_n&\leftarrow \bz_n - \eta_z \frac{\nabla_{\bz_n} l(\bz_n)}{\normtwo{\nabla_{\bz_n} l(\bz_n)}}
\gap \text{and}\gap 
\bw_m\leftarrow \bw_m - \eta_w \frac{\nabla_{\bw_m} l(\bw_m)}{\normtwo{\nabla_{\bw_m} l(\bw_m)}}.
\end{aligned}
\right.
$$
The stochastic gradient descent update for ALS is formulated in Algorithm~\ref{alg:als-regularizer-missing-stochas-gradient-realstoch}. 
It is possible that the gradient descent or stochastic gradient descent algorithm may fail to converge. In such cases, it is appropriate to re-run the algorithm using a smaller learning rate.
And in practice, the values of $m$ and $n$ in the algorithm can be randomly generated, which is why the method is termed ``stochastic." ~\footnote{When we iteratively choose the values of $m$ and $n$ from $\{1,2,\ldots, M\}$ and $\{1,2,\ldots, N\}$ in a deterministic cyclic order, respectively, the stochastic method can be referred to as  ``\textit{incremental gradient descent}."}

\begin{algorithm}[h] 
\caption{Alternating Least Squares with Full Entries and SGD}
\label{alg:als-regularizer-missing-stochas-gradient-realstoch}
\begin{algorithmic}[1] 
\Require  Matrix $\bA\in \real^{M\times N}$;
\State Initialize $\bW\in \real^{M\times K}$, $\bZ\in \real^{K\times N}$ \textcolor{mylightbluetext}{randomly without condition on the rank and the relationship between $M, N, K$}; 
\State Choose a stoping criterion on the approximation error $\delta$;
\State Choose regularization parameters $\lambda_w, \lambda_z$, and step size $\eta_w, \eta_z$;
\State Choose the maximum number of iterations $C$;
\State $iter=0$; \Comment{Count for the number of iterations}
\While{$\normf{  \bA- (\bW\bZ)}^2>\delta $ and $iter<C$}
\State $iter=iter+1$; 
\For{$n=1,2,\ldots, N$}
\For{$m=1,2,\ldots, M$} \Comment{in practice, $m,n$ can be randomly produced}
\State $\bz_n\leftarrow \bz_n - \eta_z {\nabla l(\bz_n)}/{\normtwo{\nabla l(\bz_n)}}$;\Comment{$n$-th column of $\bZ$}
\State $\bw_m\leftarrow \bw_m - \eta_w {\nabla l(\bw_m)}/{\normtwo{\nabla l(\bw_m)}}$;\Comment{$m$-th column of $\bW^\top$}
\EndFor
\EndFor

\EndWhile
\State Output $\bW^\top=[\bw_1, \bw_2, \ldots, \bw_M],\bZ=[\bz_1, \bz_2, \ldots, \bz_N]$;
\end{algorithmic} 
\end{algorithm}

\section{Bias Term}

\begin{figure}[htp]
\centering
\includegraphics[width=0.95\textwidth]{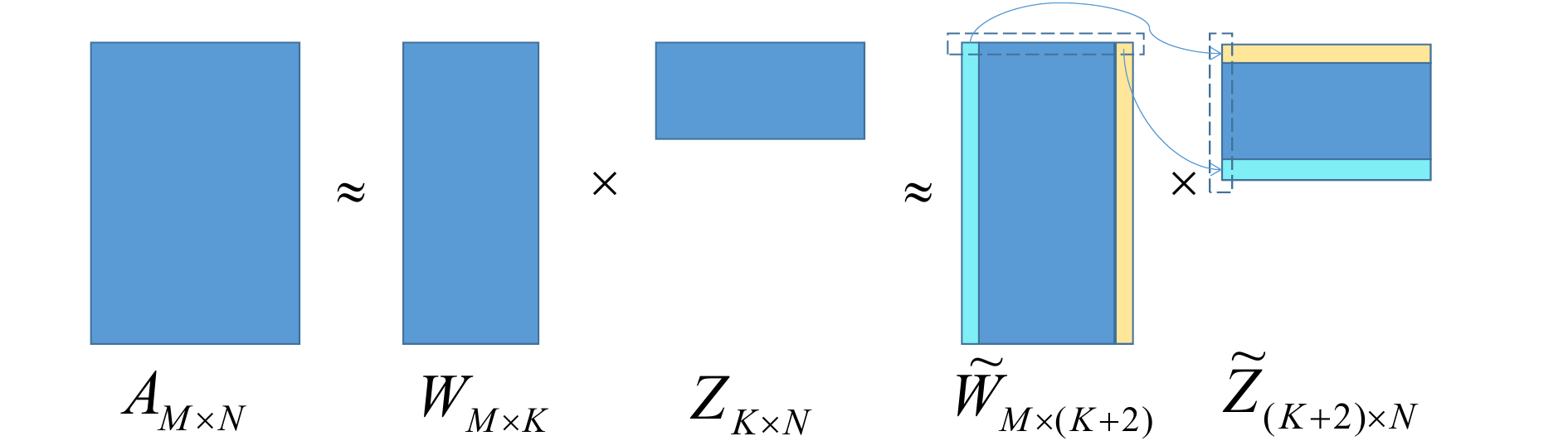}
\caption{Bias terms in alternating least squares, where the \textcolor{mydarkyellow}{yellow} entries denote ones (which are fixed), and the \textcolor{cyan}{cyan} entries denote the added features to fit the bias terms. The dotted boxes provide an example of how the bias terms work.}
\label{fig:als-bias}
\end{figure}
In the context of ordinary least squares, a bias term is usually incorporated into the raw matrix to improve model performance, as illustrated in Equation~\eqref{equation:ls-bias}. 
A similar approach can be applied to the ALS problem. 
Specifically, a fixed column filled with all ones can be appended to the \textbf{last column} of matrix $\bW$, thus an extra row should be added to the last row of matrix $\bZ$ to fit the features introduced by the bias term in $\bW$. 
Analogously, a fixed row with all ones can be added to the \textbf{first row} of $\bZ$, and an extra column in the first column of $\bW$ can be added to fit the features. The setup is depicted in Figure~\ref{fig:als-bias}.

Given the loss function with respect to the columns of $\bZ$ in Equation~\eqref{als:gradient-regularization-zn}, let 
$
\widetildebz_n
\triangleq
\scriptsize
\begin{bmatrix}
	1\\
	\bz_n
\end{bmatrix}
\in\real^{K+2}
$ be the $n$-th column of $\widetildebZ$, we have
\begin{equation}
\footnotesize
\begin{aligned}
2L(\bz_n) 
&=\normf{\widetildebW\widetildebZ-\bA}^2 +\lambda_w \normf{\widetildebW}^2 + \lambda_z \normf{\widetildebZ}^2
= 
\left\Vert
\widetildebW
\begin{bmatrix}
	1 \\
	\bz_n
\end{bmatrix}-\ba_n
\right\Vert_2^2 
+ 
\underbrace{\lambda_z \normtwo{\widetildebz_n}^2}_{=\lambda_z \normtwo{\bz_n}^2+\lambda_z}
+ 
C_{z_n}\\
&= 
\left\Vert
\begin{bmatrix}
	\widebarbw_0 & \widebarbW
\end{bmatrix}
\begin{bmatrix}
	1 \\
	\bz_n
\end{bmatrix}-\ba_n
\right\Vert_2^2 
+ \lambda_z \normtwo{\bz_n}^2 + C_{z_n}
= 
\bigg\Vert
\widebarbW \bz_n - 
\underbrace{(\ba_n-\widebarbw_0)}_{\triangleq\widebarba_n}
\bigg\Vert_2^2 
+ \lambda_z \normtwo{\bz_n}^2 + C_{z_n},
\end{aligned}
\end{equation}
where $\widebarbw_0$ represents the first column of $\widetildebW$, $\widebarbW$ denotes the remaining $K+1$ columns of $\widetildebW$ (i.e., $\widetildebW\triangleq[\widebarbw_0, \widebarbW]$), and $C_{z_n}$ is a constant with respect to $\bz_n$. Let $\widebarba_n \triangleq \ba_n-\widebarbw_0$, the update for $\bz_n$ is just similar to the one in Equation~\eqref{als:gradient-regularization-zn}, with the gradient given by
$$
\nabla_{\bz_n} L(\bz_n) = \widebarbW^\top\widebarbW\bz_n - \widebarbW^\top\widebarba_n + \lambda_z\bz_n.
$$
Therefore, the update for $\bz_n$ is given by determining the root of the  gradient above:
$$
\textbf{(update for $\widetildebz_n$)}: \quad \bz_n = (\widebarbW^\top\widebarbW+ \lambda_z\bI)^{-1} \widebarbW^\top \widebarba_n
\gap 
\implies 
\gap 
\widetildebz_n = \begin{bmatrix}
	1\\\bz_n 
\end{bmatrix},
\,\forall n.
$$
Similarly, following the loss with respect to each row of $\bW$ in Equation~\eqref{als:gradient-regularization-wd}, let $\widetildebw_m \triangleq
\scriptsize
\begin{bmatrix}
	\bw_m \\
	1
\end{bmatrix} \in\real^{K+2}$ be the $m$-th row of $\widetildebW$ (or $m$-th column of $\widetildebW^\top$), we have 
\begin{equation}
\footnotesize
\begin{aligned}
&2L(\bw_m ) 
=\normf{\widetildebZ^\top\widetildebW^\top-\bA^\top}^2 +\lambda_w \normf{\widetildebW^\top}^2 + \lambda_z \normf{\widetildebZ}^2
= 
\normtwo{\widetildebZ^\top\widetildebw_m-\bb_m}^2 + 
\underbrace{\lambda_w \normtwo{\widetildebw_m}^2}_{=\lambda_w \normtwo{\bw_m}^2+\lambda_w}
+ 
C_{w_m} \\
&= 
\bigg\Vert
\begin{bmatrix}
\widebarbZ^\top&
\widebarbz_0
\end{bmatrix}
\scriptsize
\begin{bmatrix}
\bw_m \\
1
\end{bmatrix}
\footnotesize
-\bb_m
\bigg\Vert_2^2 
+ 
\lambda_w \normtwo{\bw_m}^2
+ 
C_{w_m}
= 
\left\Vert
\widebarbZ^\top\bw_m 
-(\bb_m-\widebarbz_0) 
\right\Vert_2^2+ 
\lambda_w \normtwo{\bw_m}^2
+ 
C_{w_m},
\end{aligned}
\end{equation}
where $\widebarbz_0$ represents the last column of $\widetildebZ^\top$, $\widebarbZ^\top$ contains the remaining $K+1$ columns of $\widetildebZ^\top$ (i.e., $\widetildebZ^\top\triangleq[\widebarbZ^\top, \widebarbz_0]$),
$C_{w_m}$ is a constant with respect to $\bw_m$.  
$\bW^\top=[\bw_1, \bw_2, \ldots, \bw_M]$ and $\bA^\top=[\bb_1,\bb_2, \ldots, \bb_M]$ are the column partitions of $\bW^\top$ and $\bA^\top$, respectively. Let $\widebarbb_m \triangleq \bb_m-\widebarbz_0$, the update for $\bw_m$ is again just similar to  the one in Equation~\eqref{als:gradient-regularization-wd}, with the gradient given by
$$
\nabla_{\bw_m} L(\bw_m) = \widebarbZ\cdot \widebarbZ^\top\bw_m - \widebarbZ\cdot \widebarbb_m + \lambda_w\bw_m.
$$
Therefore, the update for $\bw_m$ is given by the root of the  gradient above:
$$
\textbf{(update for $\widetildebw_m$)}:\quad 
\bw_m=(\widebarbZ\cdot \widebarbZ^\top+\lambda_w\bI)^{-1}\widebarbZ\cdot \widebarbb_m
\gap \implies \gap 
\widetildebw_m = \begin{bmatrix}
	\bw_m \\ 1
\end{bmatrix}, 
\forall m.
$$
Similar updates can be derived using gradient descent, taking into account the bias terms and handling missing entries (see Section~\ref{section:als-gradie-descent} and \ref{section:alt-columb-by-column} for a reference).

\section{Newton's Method$^*$}\label{section:als_quasi_newton}
The ALS algorithms, whether using GD or SGD, fall under the category of first-order optimization algorithms.
However, we can also apply  Newton's method, as introduced in Section~\ref{section:app_cho_md_newton}, to derive a second-order update method:
\begin{equation}
\begin{aligned}
\bZ &\leftarrow \bZ - \big(\nabla_{\bZ}^2L(\bZ|\bW)\big)^{-1}  \nabla_{\bZ} L(\bZ|\bW) = \bZ - (\bW^\top\bW)^{-1} \bW^\top (\bW\bZ-\bA);\\
\bW &\leftarrow \bW -   \nabla_{\bW} L(\bW|\bZ) \big(\nabla_{\bW}^2L(\bW|\bZ)\big)^{-1} = \bW -  (\bW\bZ-\bA)\bZ^\top (\bZ\bZ^\top)^{-1}.
\end{aligned}
\end{equation}
In order to prevent the Hessian matrix from being singular or having a large condition number, one can still incorporate the $\ell_2$ regularization:
\begin{equation}
\begin{aligned}
\bZ &\leftarrow \bZ - \big(\nabla_{\bZ}^2L(\bZ|\bW)\big)^{-1}  \nabla_{\bZ} L(\bZ|\bW) = \bZ - (\bW^\top\bW+\lambda\bI)^{-1} \bW^\top (\bW\bZ-\bA);\\
\bW &\leftarrow \bW -   \nabla_{\bW} L(\bW|\bZ) \big(\nabla_{\bW}^2L(\bW|\bZ)\big)^{-1} = \bW -  (\bW\bZ-\bA)\bZ^\top (\bZ\bZ^\top+\lambda\bI)^{-1}.
\end{aligned}
\end{equation}

\section{Convergence$^*$}\label{section:als_convergence}

We mentioned previously that the ALS Algorithm~\ref{alg:als} belongs to a class of methods known as block coordinate descent (BCD) methods (Algorithm~\ref{alg:two_bcd_gen_inals}), where we   iteratively optimize over subsets of the variables (blocks) rather than over all variables simultaneously. This approach is particularly useful when dealing with large-scale optimization problems, where optimizing over all variables at once can be computationally expensive or impractical, though the solution may not always converge to the global optimum, especially for non-convex problems.
In some cases, multiple blocks can be updated independently, allowing for parallel execution.
In the context of the  ALS problem, $(\bW,\bZ)$ can be considered as two blocks of a single variable, hence the name 2-BCD.
The convergence to a stationary point (where the gradient vanishes) for the 2-BCD approach is discussed below.
\begin{theorem}[Convergence of 2-BCD Method \citep{grippo2000convergence, beck2017first, jain2017non, gillis2020nonnegative}]
The limit points of the iterates of a 2-BCD algorithm are stationary points, provided that the following two conditions are met:
\begin{enumerate}
\item  The objective function is continuously differentiable.
\item Each block of variables is required to belong to a closed convex set.
\end{enumerate}
\end{theorem}
For the ALS algorithms (and most NMF approaches introduced in the next chapter), both  conditions hold. Therefore, a guarantee of convergence to stationary points of the algorithms is achieved.
More generally, the convergence conditions for a general BCD approach are discussed below.

\begin{theorem}[Convergence of BCD Method \citep{bertsekas1997nonlinear, gillis2020nonnegative}]
The limit points of the iterates of a BCD algorithm are stationary points, provided that the following conditions hold:
\begin{enumerate}
\item  The objective function is continuously differentiable.
\item  Each block of variables is required to belong to a closed convex set.
\item  The minimum computed at each iteration for a given block of variables is uniquely attained.
\item  The objective function values are monotonically decreasing between all iterates and the next (which is obtained by updating a single block of variables).
\end{enumerate}
\end{theorem}

\section{Applications}

\subsection{Low-Rank Approximation}\label{section:als-low-flag}
We have discussed and compared the effects of using SVD and pseudoskeleton methods for low-rank approximation in Section~\ref{section:svd-low-rank-approxi}. The image intended for compression is shown in Figure~\ref{fig:eng300} with a size of $600\times 1200$ and rank  402. 
Figure~\ref{fig:svdd-by-parts} shows the image reconstructed using only the first singular value already closely approximates the original image. 
Additionally, in Figure~\ref{fig:svdd-pseudoskeleton}, we observe the differences in compression quality for various ranks, specifically 90, 60, 30, and 10. 
Pseudoskeleton performs well in the black horizontal and vertical lines within the image, but struggles to maintain the fine details in the flag.

\noindent
\begin{figure}[h]
\centering  
\vspace{-0.35cm}  
\subfigtopskip=2pt 
\subfigbottomskip=2pt 
\subfigcapskip=-5pt  
\subfigure[$\sigma_1\bu_1\bv_1^\top$\protect\newline$F_1=60217$]{\label{fig:svd12}
\includegraphics[width=0.15\linewidth]{./imgs/svd_pic1.png}}
\subfigure[$\sigma_2\bu_2\bv_2^\top$\protect\newline$F_2=120150$]{\label{fig:svd22}
\includegraphics[width=0.15\linewidth]{./imgs/svd_pic2.png}}
\subfigure[$\sigma_3\bu_3\bv_3^\top$\protect\newline$F_3=124141$]{\label{fig:svd32}
\includegraphics[width=0.15\linewidth]{./imgs/svd_pic3.png}}
\subfigure[$\sigma_4\bu_4\bv_4^\top$\protect\newline$F_4=125937$]{\label{fig:svd42}
\includegraphics[width=0.15\linewidth]{./imgs/svd_pic4.png}}
\subfigure[$\sigma_5\bu_5\bv_5^\top$\protect\newline$F_5=126127$]{\label{fig:svd52}
\includegraphics[width=0.15\linewidth]{./imgs/svd_pic5.png}}
\subfigure[{ All 5 singular values: $\sum_{i=1}^{5}\sigma_i\bu_i\bv_i^\top$,\protect\newline$F=$\textbf{44379}.}]{\label{fig:svd62}
\includegraphics[width=0.15\linewidth]{./imgs/svd_pic6_all.png}}
\quad
\subfigure[$\bc_1\br_1^\top$\protect\newline$G_1=60464$]{\label{fig:skeleton51}
\includegraphics[width=0.15\linewidth]{./imgs/skeleton_5_1.png}}
\subfigure[$\bc_2\br_2^\top$\protect\newline$G_2=122142$]{\label{fig:skeleton52}
\includegraphics[width=0.15\linewidth]{./imgs/skeleton_5_2.png}}
\subfigure[$\bc_3\br_3^\top$\protect\newline$G_3=123450$]{\label{fig:skeleton53}
\includegraphics[width=0.15\linewidth]{./imgs/skeleton_5_3.png}}
\subfigure[$\bc_4\br_4^\top$\protect\newline$G_4=125975$]{\label{fig:skeleton54}
\includegraphics[width=0.15\linewidth]{./imgs/skeleton_5_5.png}}
\subfigure[$\bc_5\br_5^\top$\protect\newline$G_5=124794$]{\label{fig:skeleton55}
\includegraphics[width=0.15\linewidth]{./imgs/skeleton_5_4.png}}
\subfigure[Pseudoskeleton Rank 5 $\sum_{i=1}^{5}\bc_i\br_i^\top$,\protect\newline$G=45905$.]{\label{fig:skeleton5_all}
\includegraphics[width=0.15\linewidth]{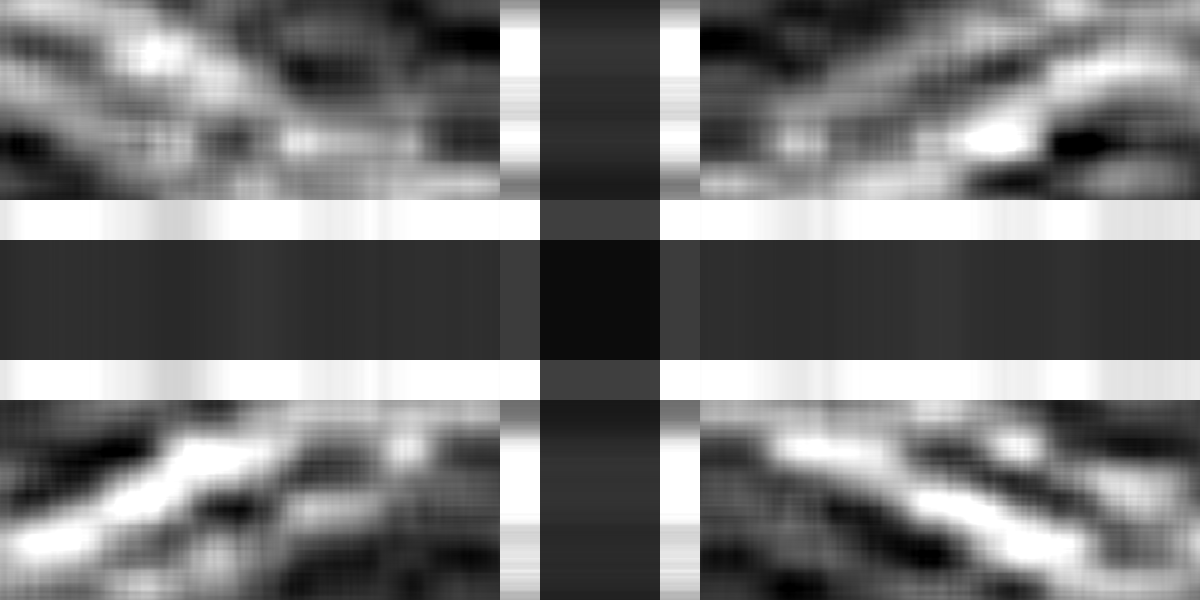}}
\quad
\subfigure[$\bw_1\bz_1^\top$\protect\newline$S_1=82727$]{\label{fig:als51}
\includegraphics[width=0.15\linewidth]{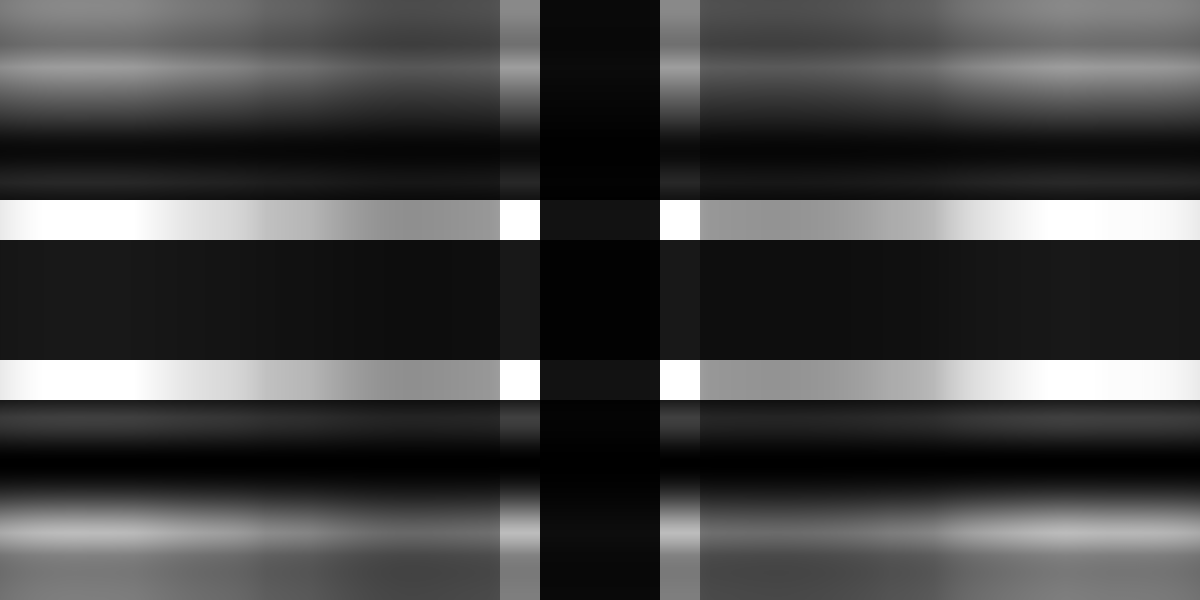}}
\subfigure[$\bw_2\bz_2^\top$\protect\newline$S_2=107355$]{\label{fig:als52}
\includegraphics[width=0.15\linewidth]{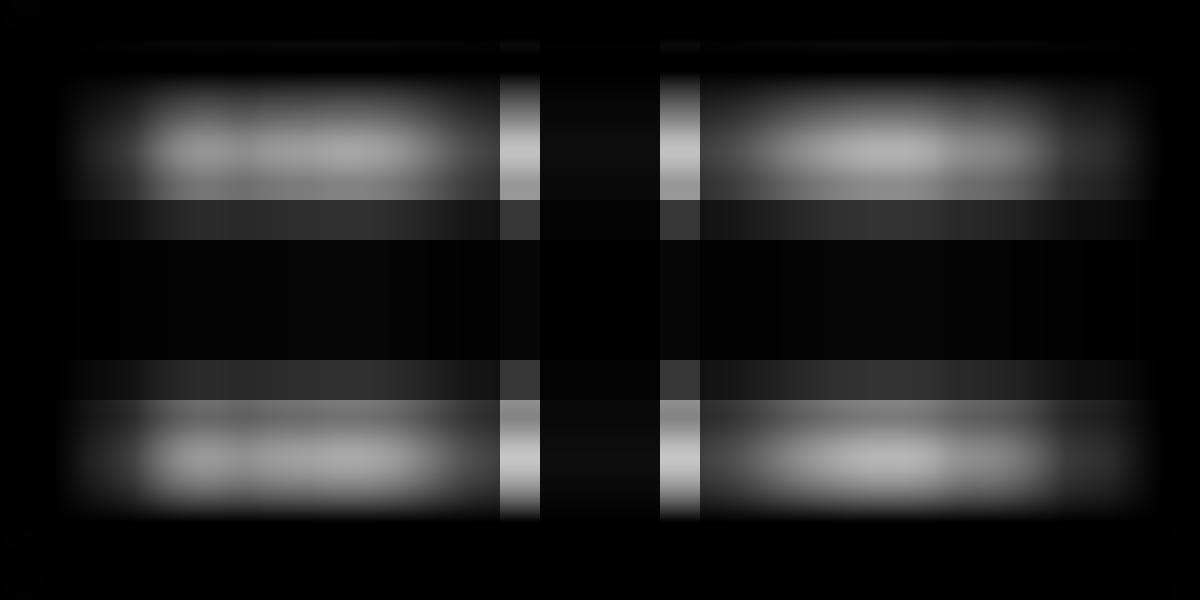}}
\subfigure[$\bw_3\bz_3^\top$\protect\newline$S_3=119138$]{\label{fig:als53}
\includegraphics[width=0.15\linewidth]{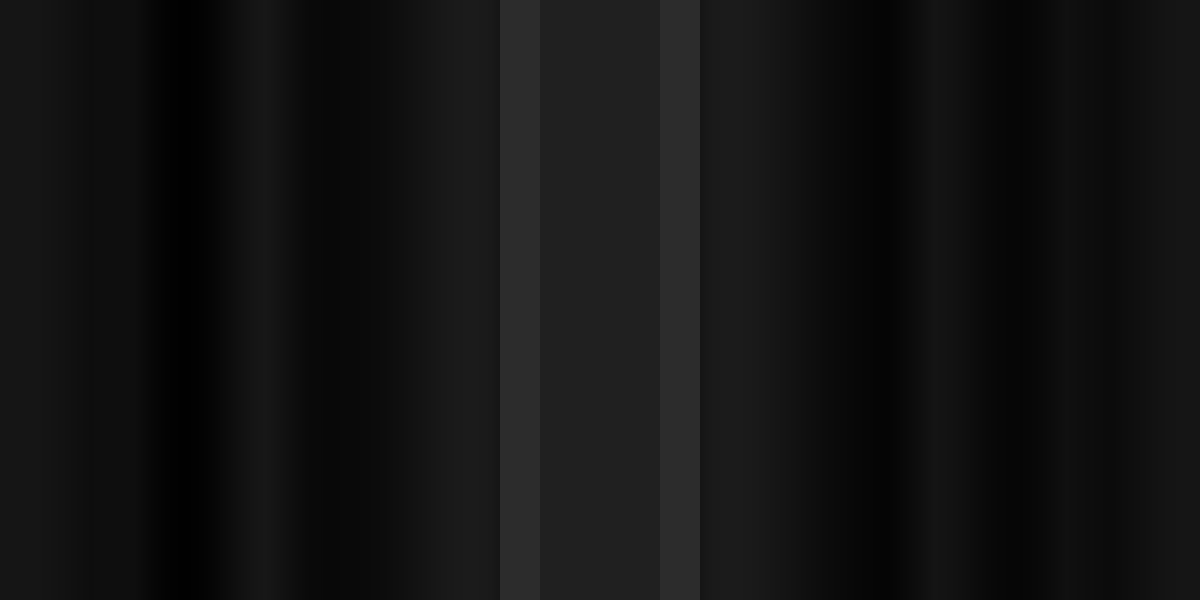}}
\subfigure[$\bw_4\bz_4^\top$\protect\newline$S_4=120022$]{\label{fig:als54}
\includegraphics[width=0.15\linewidth]{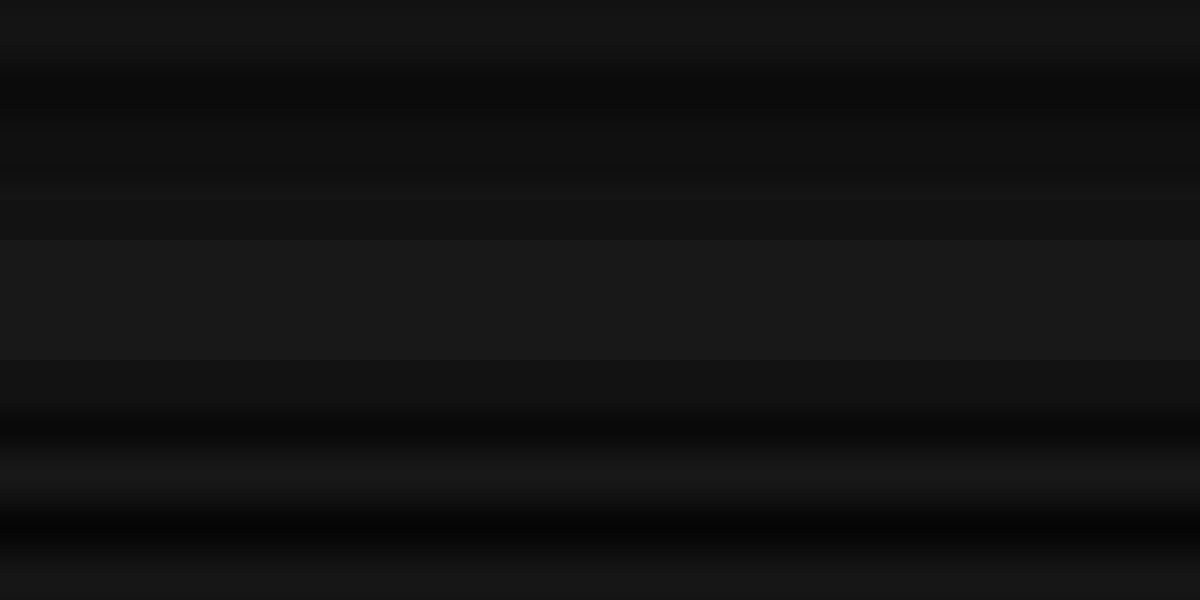}}
\subfigure[$\bw_5\bz_5^\top$\protect\newline$S_5=120280$]{\label{fig:als55}
\includegraphics[width=0.15\linewidth]{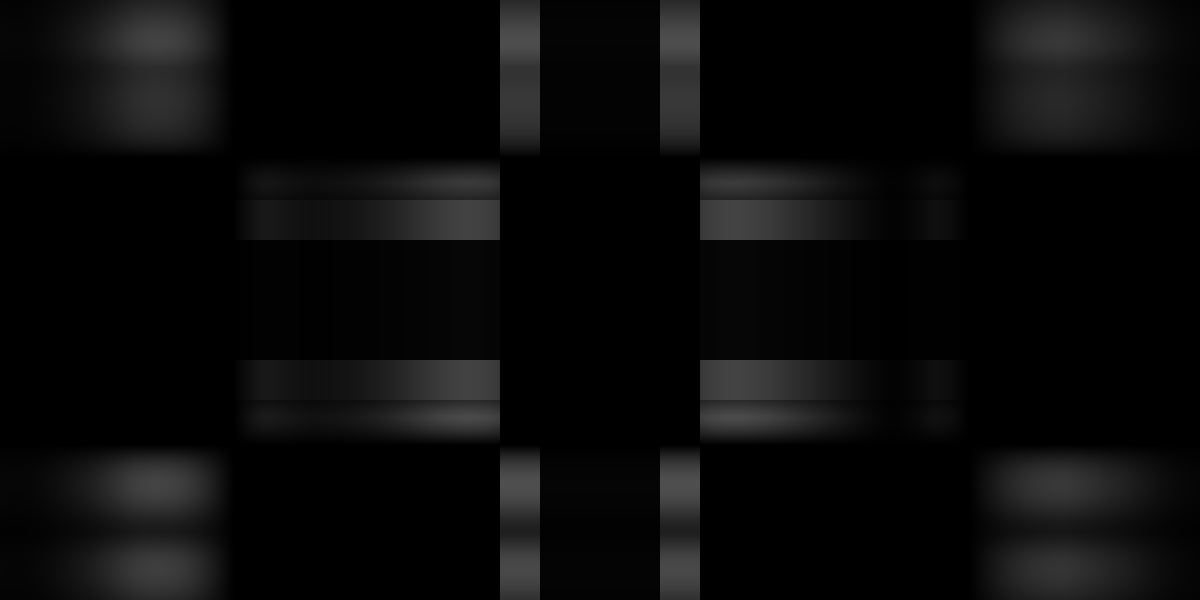}}
\subfigure[ALS Rank 5 $\sum_{i=1}^{5}\bw_i\bz_i^\top$,\protect\newline$S=52157$.]{\label{fig:als5_all}
\includegraphics[width=0.15\linewidth]{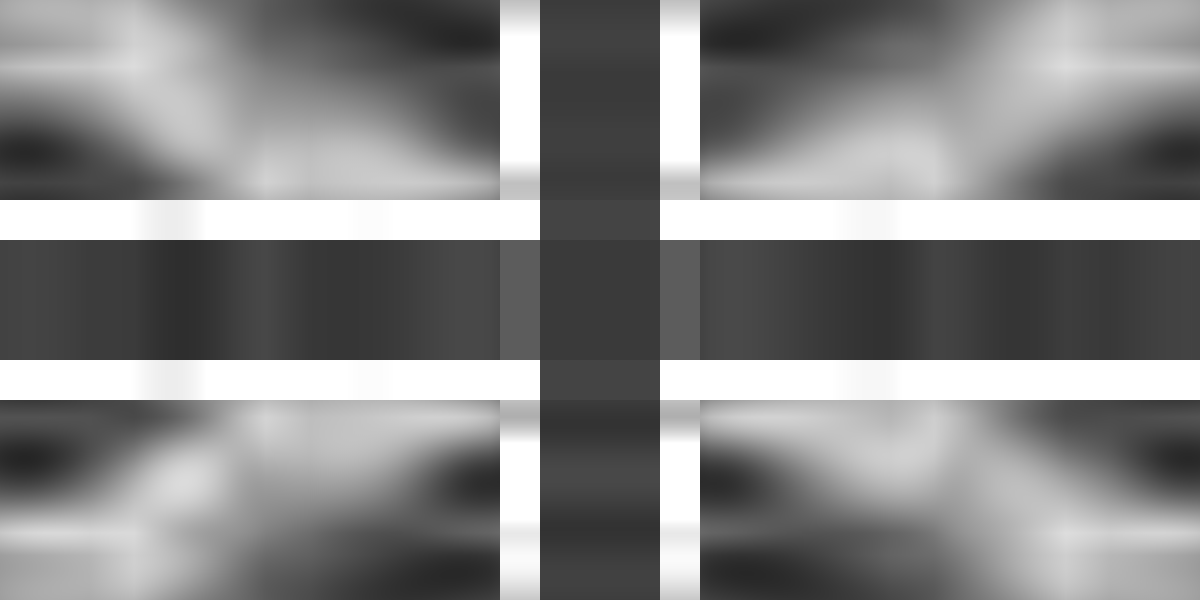}}
\caption{Image compression on a grayscale flag image, reducing it into a rank-5 matrix via the SVD, and decomposing into 5 parts, where $\sigma_1 \geq \sigma_2 \geq \ldots \geq \sigma_{5}$, i.e., $F_1\leq F_2\leq \ldots \leq F_5$ with $F_i \triangleq\normf{\sigma_i\bu_i\bv^\top - \bA}$ for $i\in \{1,2,\ldots, 5\}$. And reconstruct images by single singular value and its corresponding left and right singular vectors. Similar approaches are applied to the pseudoskeleton and ALS algorithms ($\bc_i\br_i^\top$ and $\bw_i\bz_i^\top$, respectively). \textbf{Upper:} SVD; \textbf{Middle:} Pseudoskeleton; \textbf{Lower:} ALS.}
\label{fig:svdd-by-parts-als}
\end{figure}

Similar results can be observed when employing the ALS decomposition for low-rank approximation. 
The ALS approximation is given by $\bA\approx\bW\bZ$, where $\bW\in \real^{m\times \gamma}$ and $\bZ\in \real^{\gamma\times n}$ if $\bA\in \real^{m\times n}$ such that $\bW$ and $\bZ$ are rank-$\gamma$ matrices. Suppose $\gamma=5$, and 
$$
\bW=[\bw_1, \bw_2, \ldots, \bw_5]
\qquad
\text{and}
\qquad 
\bZ = 
\footnotesize
\begin{bmatrix}
	\bz_1^\top; 
	\bz_2^\top;
	\ldots;
	\bz_5^\top
\end{bmatrix}
$$ 
are the column and row partitions of $\bW$ and $\bZ$, respectively \footnote{For simplicity, note that this definition differs from the one  provided in Section~\ref{section:als-netflix}, where   $\bw_i$ was defined as the rows of $\bW$. }. Then $\bA$ can be approximated as $\sum_{i=1}^{5}\bw_i\bz_i^\top$. The partitions are ordered such that 
$$
\underbrace{\normf{\bw_1\bz_1^\top-\bA}}_{S_1} \leq 
\underbrace{\normf{\bw_2\bz_2^\top-\bA}}_{S_2}
\leq \ldots \leq 
\underbrace{\normf{\bw_5\bz_5^\top-\bA}}_{S_5}.
$$
It is observed that $\bw_1\bz_1^\top$ behaves slightly \textbf{differently} compared to $\sigma_1\bu_1\bv^\top$, where the reconstruction errors   measured by the Frobenius norm are not as close  (82,727 in the ALS case compared to that of 60,217 in the SVD case). As  mentioned previously, $\bc_1\br_1^\top$ works similarly to  $\sigma_1\bu_1\bv^\top$ because the pseudoskeleton method relies on the SVD. However, in ALS, the reconstruction is based on least squares optimization. 
The key difference between ALS and SVD lies in how they prioritize the vectors in the basis.  
In SVD, the significance of each vector is determined by the corresponding singular value. This usually means that the first vector of the basis dominates and is the most used vector to reconstruct data; followed by the second vector, and so on. 
Thus, SVD implicitly establishes a hierarchy among the basis vectors, which is not present in ALS.
In Figure~\ref{fig:als52}, we find the second component $\bw_2\bz_2^\top$ obtained via ALS  plays an important role in the reconstruction of the original figure. In contrast, the second component $\sigma_2\bu_2\bv_2^\top$ obtained via SVD in Figure~\ref{fig:svd22}  contributes less significantly to the reconstruction.
\begin{SCfigure}
\caption{Comparison of reconstruction errors measured using Frobenius norm, among the SVD, pseudoskeleton, and ALS, where the approximated rank ranges from $3$ to 100. It is worth noting that ALS with well-selected parameters performs similarly to SVD.}
\includegraphics[width=0.5\textwidth]{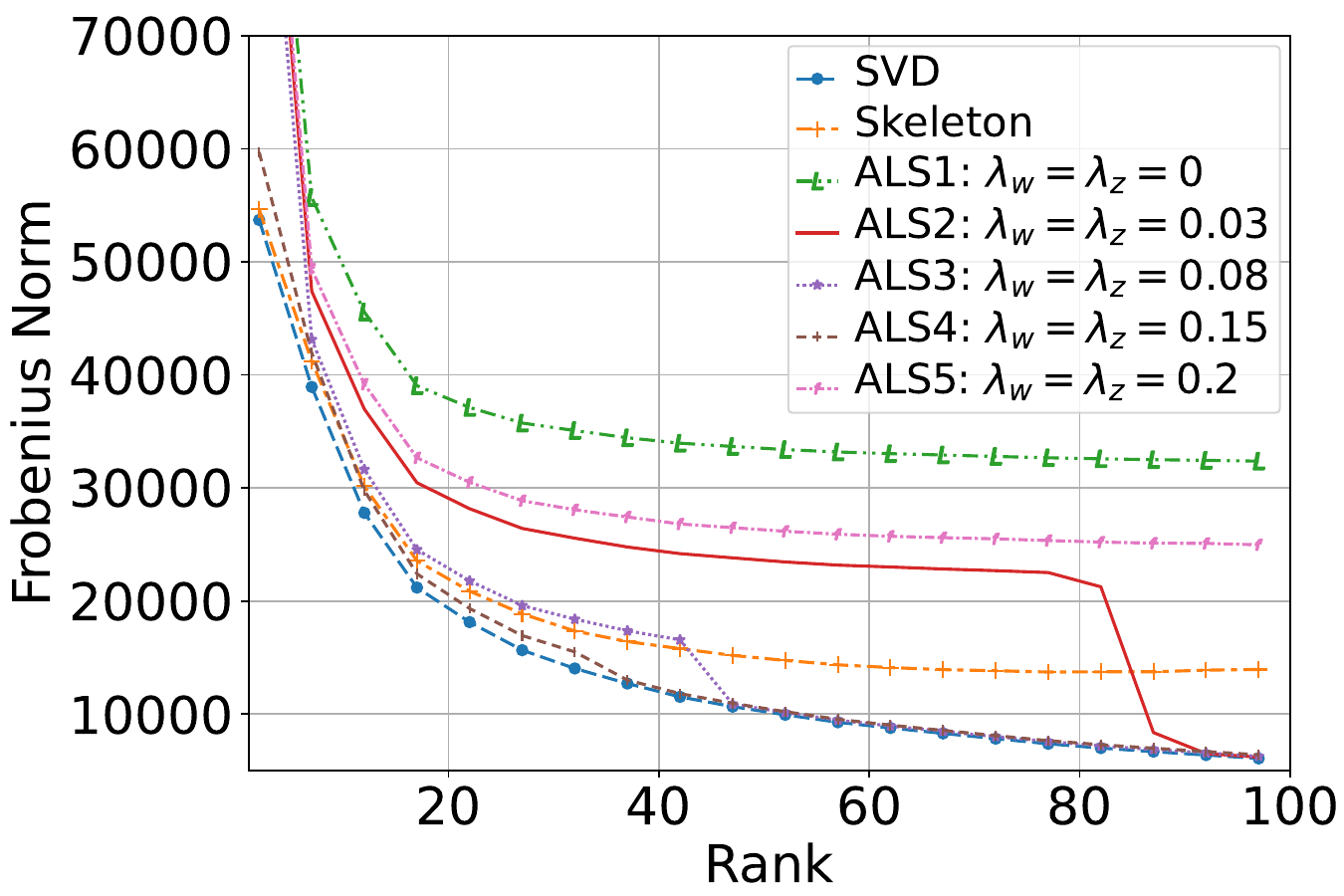}
\label{fig:svd_skeleton_als_fnorm}
\end{SCfigure}

\begin{figure}[h]
	\centering
	\vspace{-0.35cm}
	\subfigtopskip=2pt
	\subfigbottomskip=2pt
	\subfigcapskip=-5pt
	\subfigure[SVD with rank 90\protect\newline Frobenius norm=\textbf{6,498}]{\label{fig:svd902}
		\includegraphics[width=0.3\linewidth]{./imgs/svd90.png}}
	\quad 
	\subfigure[Pseudoskeleton with rank 90\protect\newline Frobenius norm=13,751]{\label{fig:skeleton902}
		\includegraphics[width=0.3\linewidth]{./imgs/skeleton90.png}}
	\subfigure[ALS with rank 90\protect\newline Frobenius norm=6,622]{\label{fig:als_90}
		\includegraphics[width=0.3\linewidth]{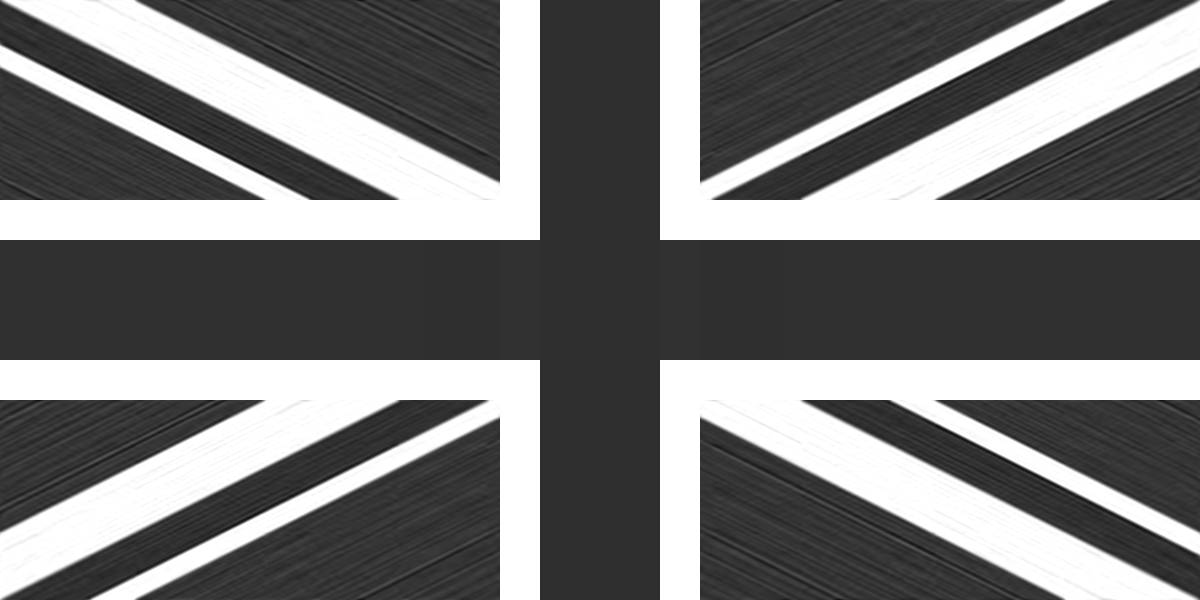}}\\
	\subfigure[SVD with rank 60\protect\newline Frobenius norm=\textbf{8,956}]{\label{fig:svd602}
		\includegraphics[width=0.3\linewidth]{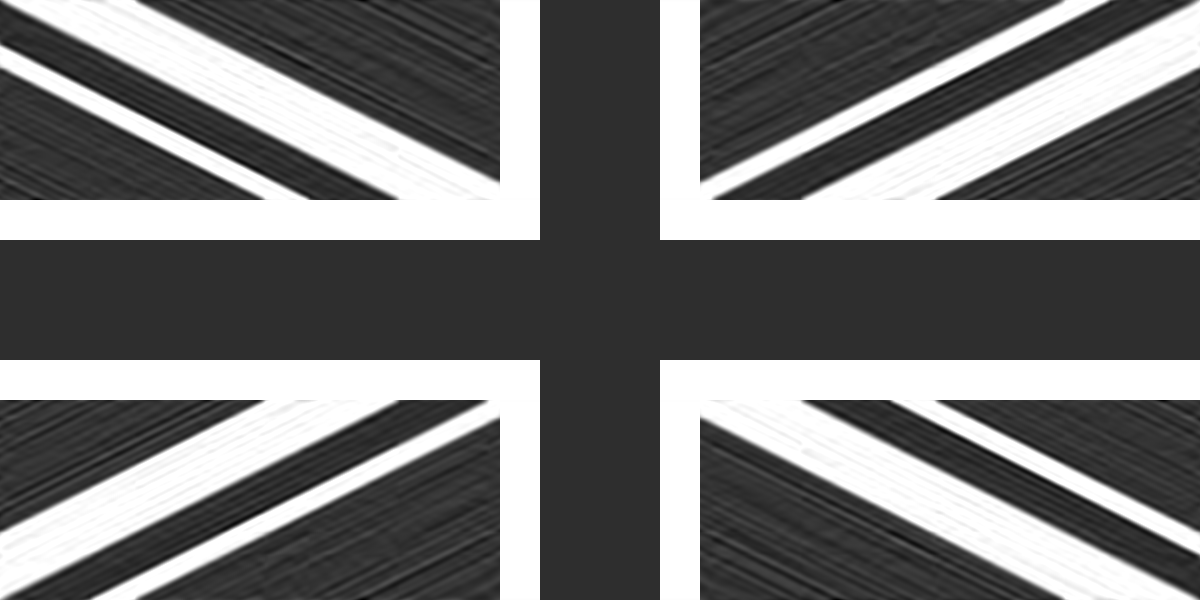}}
	\quad 
	\subfigure[Pseudoskeleton with rank 60\protect\newline Frobenius norm=14,217]{\label{fig:skeleton602}
		\includegraphics[width=0.3\linewidth]{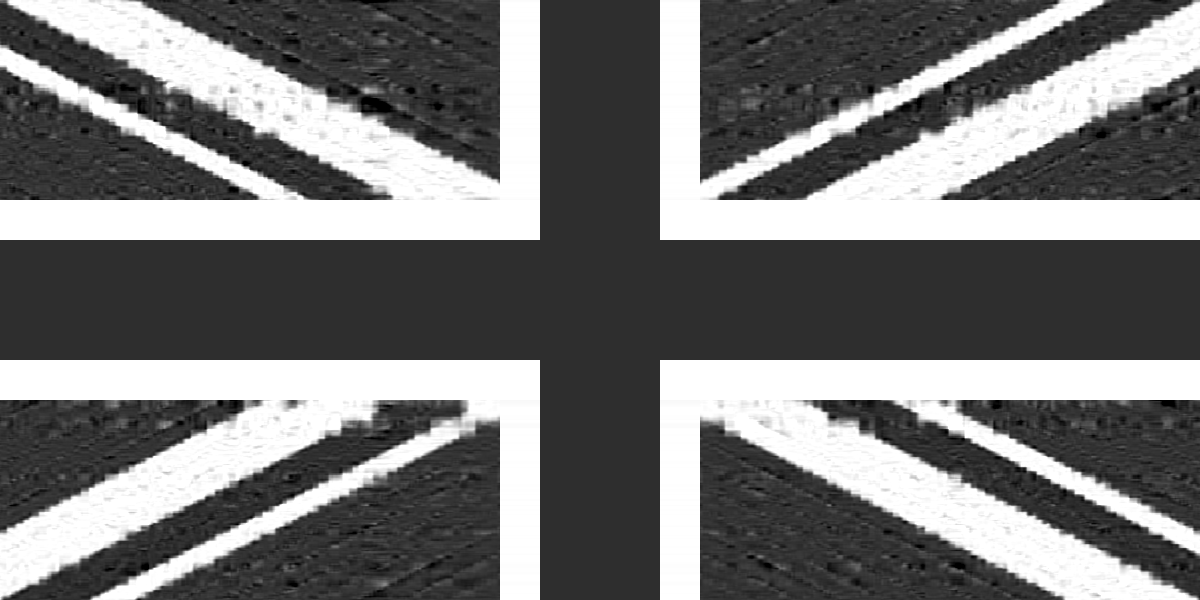}}
	\subfigure[ALS with rank 60\protect\newline Frobenius norm=9,028]{\label{fig:als_60}
		\includegraphics[width=0.3\linewidth]{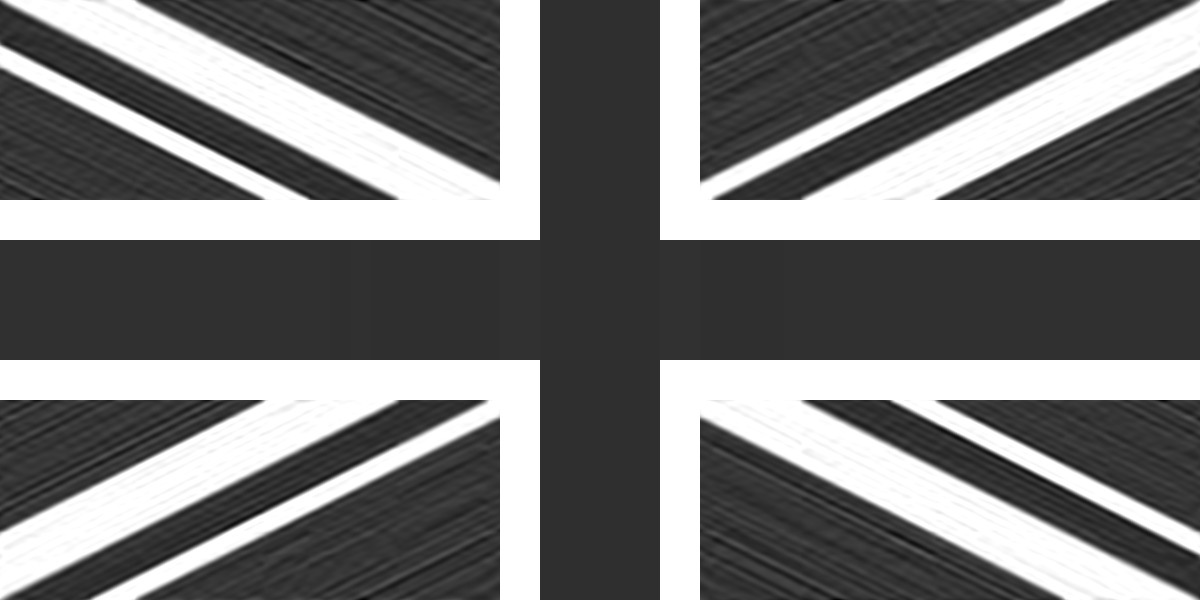}}\\
	\subfigure[SVD with rank 30\protect\newline Frobenius norm=\textbf{14,586}]{\label{fig:svd302}
		\includegraphics[width=0.3\linewidth]{./imgs/svd30.png}}
	\quad 
	\subfigure[Pseudoskeleton with rank 30\protect\newline Frobenius norm=17,853]{\label{fig:skeleton302}
		\includegraphics[width=0.3\linewidth]{./imgs/skeleton30.png}}
	\subfigure[ALS with rank 30\protect\newline Frobenius norm=18,624]{\label{fig:als_30}
		\includegraphics[width=0.3\linewidth]{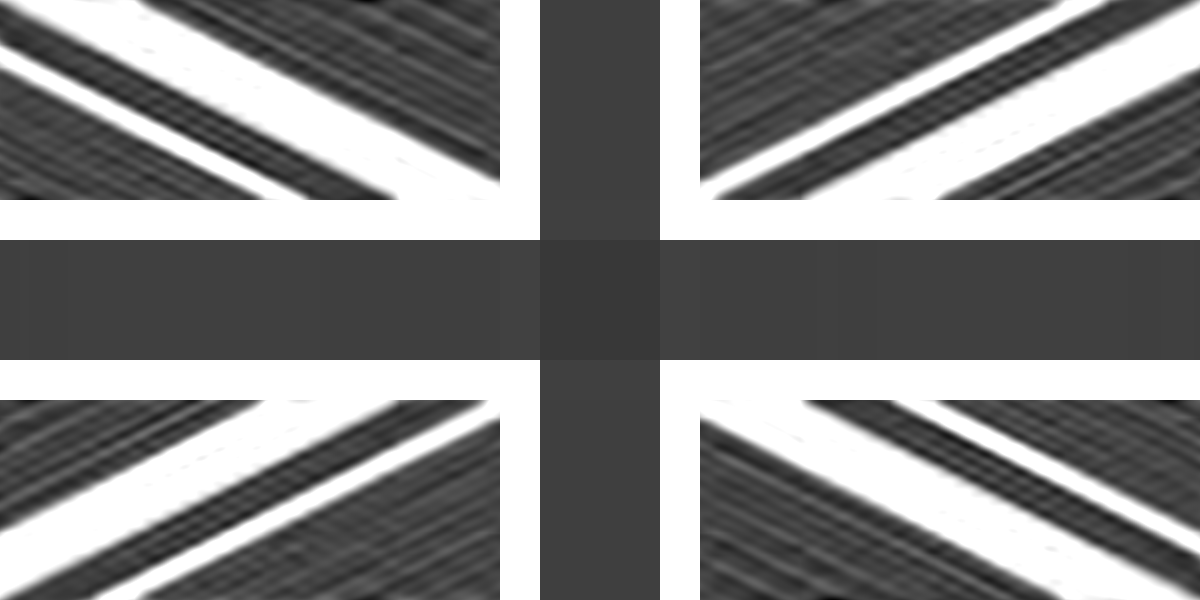}}
	\subfigure[SVD with rank 10\protect\newline Frobenius norm=\textbf{31,402}]{\label{fig:svd102}
		\includegraphics[width=0.3\linewidth]{./imgs/svd10.png}}
	\quad 
	\subfigure[Pseudoskeleton with rank 10\protect\newline Frobenius norm=33,797]{\label{fig:skeleton102}
		\includegraphics[width=0.3\linewidth]{./imgs/skeleton10.png}}
	\subfigure[ALS with rank 10\protect\newline Frobenius norm=33,449]{\label{fig:als_10}
		\includegraphics[width=0.3\linewidth]{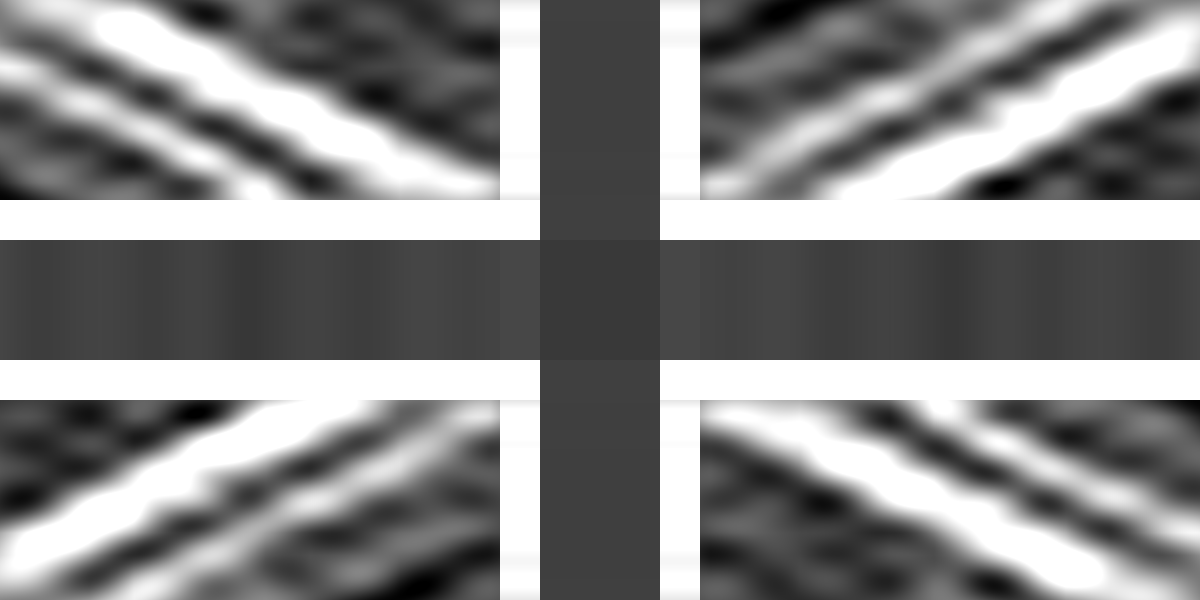}}
	\caption{Image compression for a grayscale flag image with different ranks.}
	\label{fig:svdd-pseudoskeleton-als}
\end{figure}

We finally compare low-rank approximations using   SVD, pseudoskeleton, and ALS with varying ranks (range from 3 to 100). Figure~\ref{fig:svdd-pseudoskeleton-als} demonstrates the differences in compression quality at ranks 90, 60, 30, and 10. 
We observe that  SVD reconstructs well at ranks 90, 60, and 30. 
The pseudoskeleton-approximation compresses well in the black horizontal and vertical lines within the image but performs poorly in capturing  the finer details of the flag. 
ALS works similarly to  SVD in terms of visual expression and reconstruction errors measured by the Frobenius norm.
Figure~\ref{fig:svd_skeleton_als_fnorm} compares the reconstruction errors measured by the Frobenius norm among the (truncated) SVD,  pseudoskeleton, and  ALS approximations, with ranks   ranging from  $3$ to $100$. 
In all cases, truncated SVD consistently outperforms the others in terms of the Frobenius norm. 
Similar trends are observed when evaluating the spectral norm. 
Notably, ALS outperforms  the pseudoskeleton decomposition when $\lambda_w=\lambda_z=0.15$. An interesting cutoff happens at $\lambda_w=\lambda_z=\{0.03, 0.08, 0.15\}$: as the rank increases,  ALS approaches the performance of  SVD closely in  low-rank approximation.

\index{Truncated}
\index{Truncated SVD}

\subsection{Movie Recommender}\label{section:movie_rec_als}
The ALS algorithm has been extensively developed for movie recommendation systems. 
To illustrate this, we use  the ``MovieLens 100K" data set from MovieLens \citep{harper2015movielens}~\footnote{http://grouplens.org}. 
This data set is widely recognized and used in the field of recommender systems research due to its comprehensive set of user ratings for movies. 
It consists of 100,000 ratings from 943 users for 1,682 movies, with rating values ranging from 0 to 5. The data was collected through the MovieLens website 
over a seven-month period from September 19th, 
1997 to April 22nd, 1998. This data has been cleaned up---users
who had less than 20 ratings or did not have complete demographic
information were removed from this data set such that simple demographic info for the users (age, gender, occupation, zip) can be obtained. However, our focus will solely be on the raw rating matrix  to evaluate how well the low-rank ALS approach can capture the underlying structure of the data, leading to accurate and meaningful recommendations.

\begin{figure}[h]
\centering  
\vspace{-0.35cm} 
\subfigtopskip=2pt 
\subfigbottomskip=2pt
\subfigcapskip=-5pt
\subfigure[Training set.]{\label{fig:movie100k1}
\includegraphics[width=0.47\linewidth]{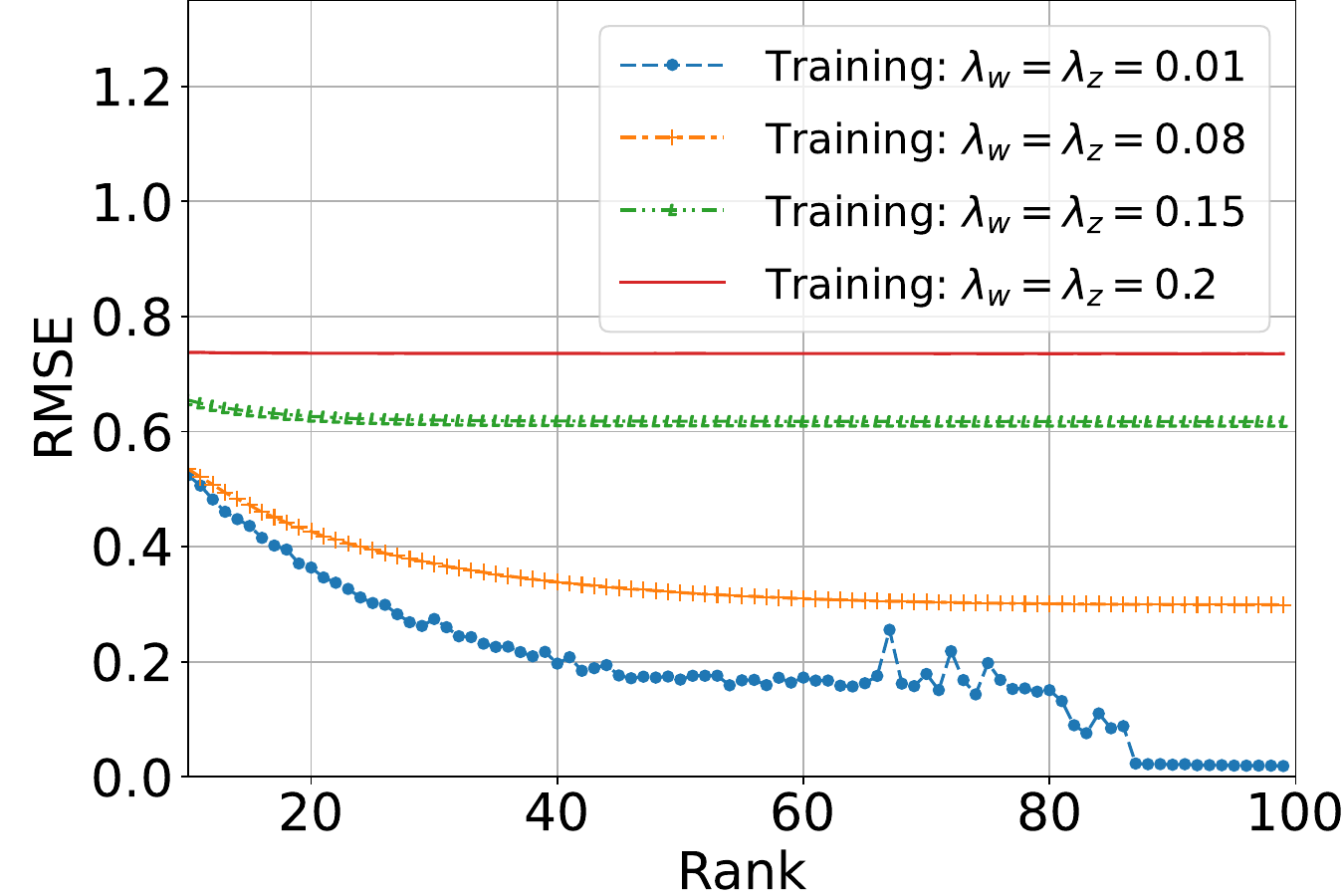}}
\quad 
\subfigure[Validation set.]{\label{fig:movie100k2}
\includegraphics[width=0.47\linewidth]{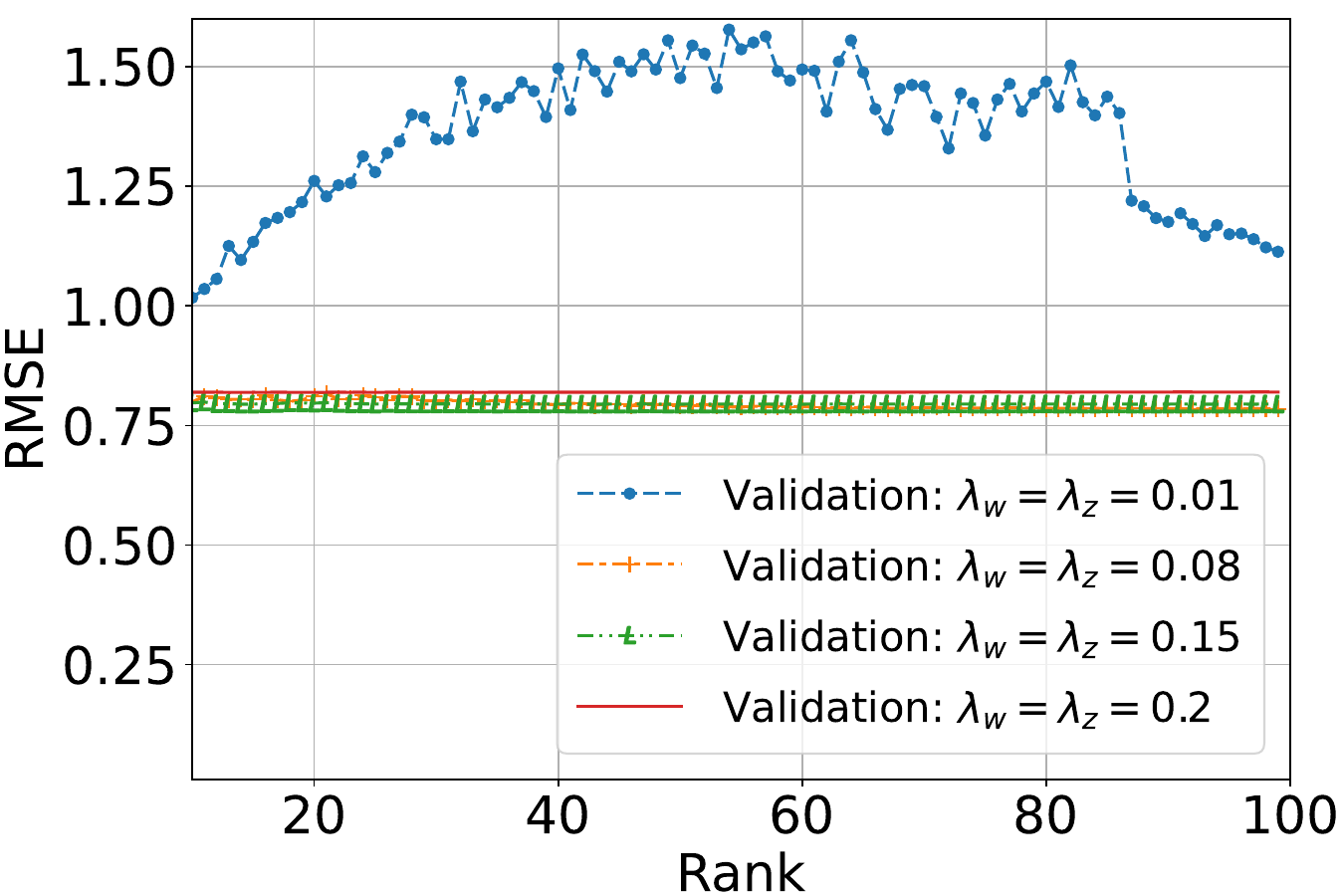}}
\caption{Comparison of training and validation error for the ``MovieLens 100K" data set with different reduction dimensions and regularization parameters.}
\label{fig:movie100k}
\end{figure}

The data set is split into training and validation set, comprising approximately 95,015 and 4,985 ratings, respectively, to fit by the ALS algorithm.
The error is quantified using the \textit{root mean squared error (RMSE)}. The RMSE is a common measure of the difference between actual and predicted values. For a set of values $\{x_1, x_2, \ldots, x_n\}$ and their predictions $\{\hat{x}_1, \hat{x}_2, \ldots, \hat{x}_n\}$, the RMSE can be described as 
$
\text{RMSE}(\bx, \hat{\bx}) = \sqrt{\frac{1}{n} \sum_{i=1}^{n}(x_i-\hat{x}_i)^2}.
$
For evaluating the ALS algorithm, the minimum RMSE for the validation set is achieved with $K=62$ and $\lambda_w=\lambda_z=0.15$, resulting in an RMSE of $0.806$ (less than 1), as shown in Figure~\ref{fig:movie100k}. 
Given that ratings range from 0 to 5, the ALS algorithm can predict whether a user is likely to enjoy a movie (e.g., ratings of 4 to 5) or not  (e.g., ratings of 0 to 2). 

\paragraph{Recommender 1.}
A simple recommender system suggests movie $m$ to user $n$ if $a_{mn}\geq4$ and user $n$ has not yet rated movie $m$. 
\paragraph{Recommender 2.}
Alternatively, we can recommend similar movies to those highly rated by the user. 
Suppose user $n$ has rated movie $m$ with a 5 ($a_{mn}=5$). Under the ALS approximation $\bA=\bW\bZ$, where each row of $\bW$ represents the hidden features of each movie (see Section~\ref{section:als-vector-product} on  vector inner products), the solution involves identifying the most similar movies that user  $n$ has not rated (or watched), to movie $m$. Mathematically, this is expressed as:
$$
\mathop{\argmax}_{\bw_i} \gap \text{similarity}(\bw_i, \bw_m), \qquad \text{for all} \gap i \notin \bo_n,
$$
where $\bw_i$'s are the rows of $\bW$, each representing the hidden features of movie $i$, and $\bo_n$ represents a mask vector, indicating the movies that user $n$ has rated.

The method described above relies on a similarity function applied to two vectors. The \textit{cosine similarity} is the most commonly used measure. It is defined as the cosine of the angle between the two vectors:
$$
\cos(\bx, \by) = \frac{\bx^\top\by}{\normtwo{\bx}\cdot \normtwo{\by}},
$$
where the value ranges from $-1$ to 1, with $-1$ representing perfectly dissimilar and 1 being perfectly similar. Based on this definition, it follows that the cosine similarity depends only on the angle between the two nonzero vectors, but not on their magnitudes since it can be regarded as the inner product between the normalized  versions of these vectors. 
Another measure for calculating similarity is the \textit{Pearson similarity}:
$$
\text{Pearson}(\bx,\by) =\frac{\Cov(\bx,\by)}{\sigma_x \cdot \sigma_y}
= \frac{\sum_{i=1}^{n} (x_i - \bar{x} ) (y_i -\bar{y})}{ \sqrt{\sum_{i=1}^{n} (x_i-\bar{x})^2}\sqrt{ \sum_{i=1}^{n} (y_i-\bar{y})^2 }}.
$$
It is calculated as the ratio between the covariance of two variables and the product of their standard deviations, whose range varies between $-1$ and 1, where $-1$ is perfectly dissimilar, 1 is perfectly similar, and 0 indicates no linear relationship.  Pearson similarity is commonly used to measure the linear correlation between two sets of data. 

Both Pearson correlation and cosine similarity are widely used in various fields, including machine learning and data analysis. Pearson correlation is often used in regression analysis, while cosine similarity is commonly used in recommendation systems and information retrieval. 
In our context, cosine similarity performs better in \textit{precision-recall (PR) curve} analysis. 
\begin{figure}[h]
\centering
\vspace{-0.35cm}
\subfigtopskip=2pt
\subfigbottomskip=2pt
\subfigcapskip=-5pt
\subfigure[Cosine Bin Plot.]{\label{fig:als-cosine}%
\includegraphics[width=0.32\linewidth]{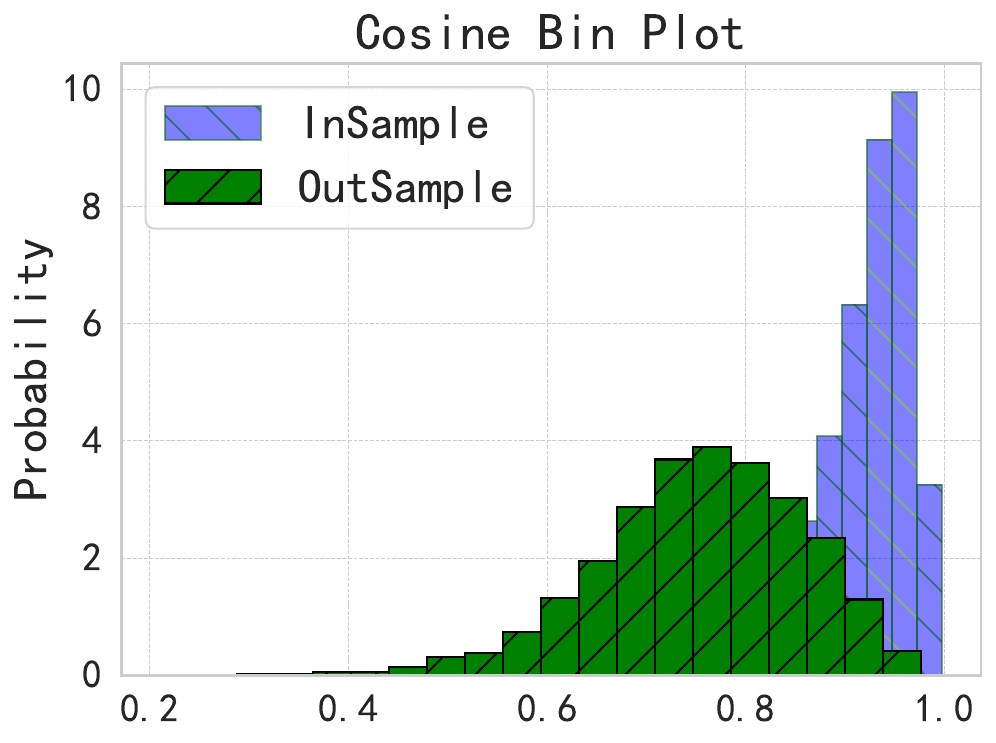}}%
\subfigure[Pearson Bin Plot.]{\label{fig:als-pearson}%
\includegraphics[width=0.32\linewidth]{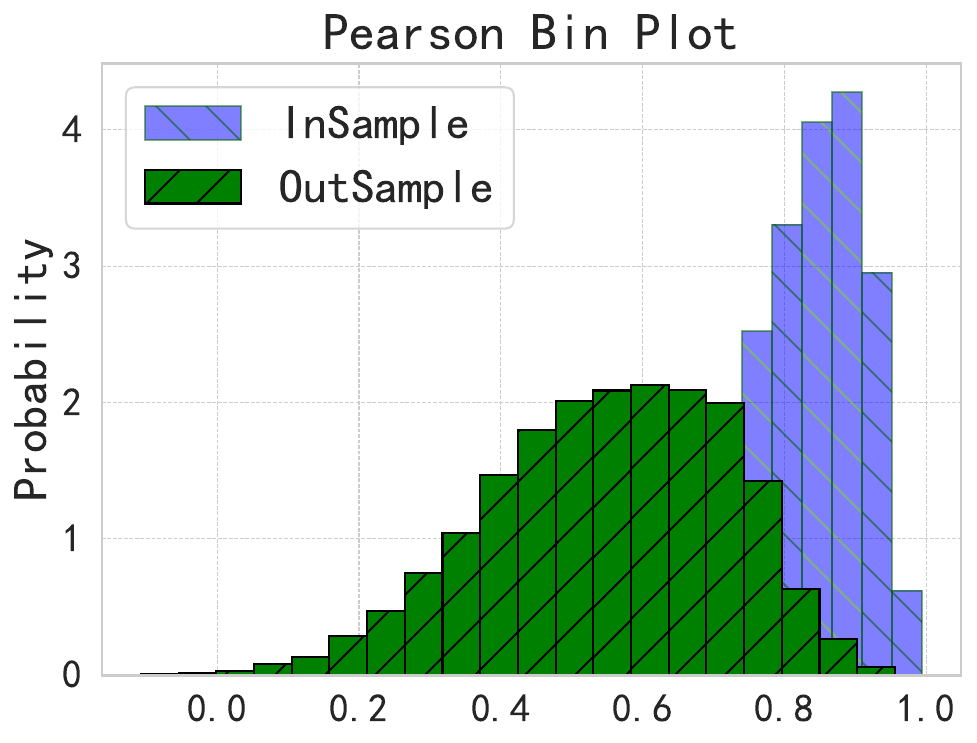}}%
\subfigure[PR Curve.]{\label{fig:als-prcurve}%
\includegraphics[width=0.35\linewidth]{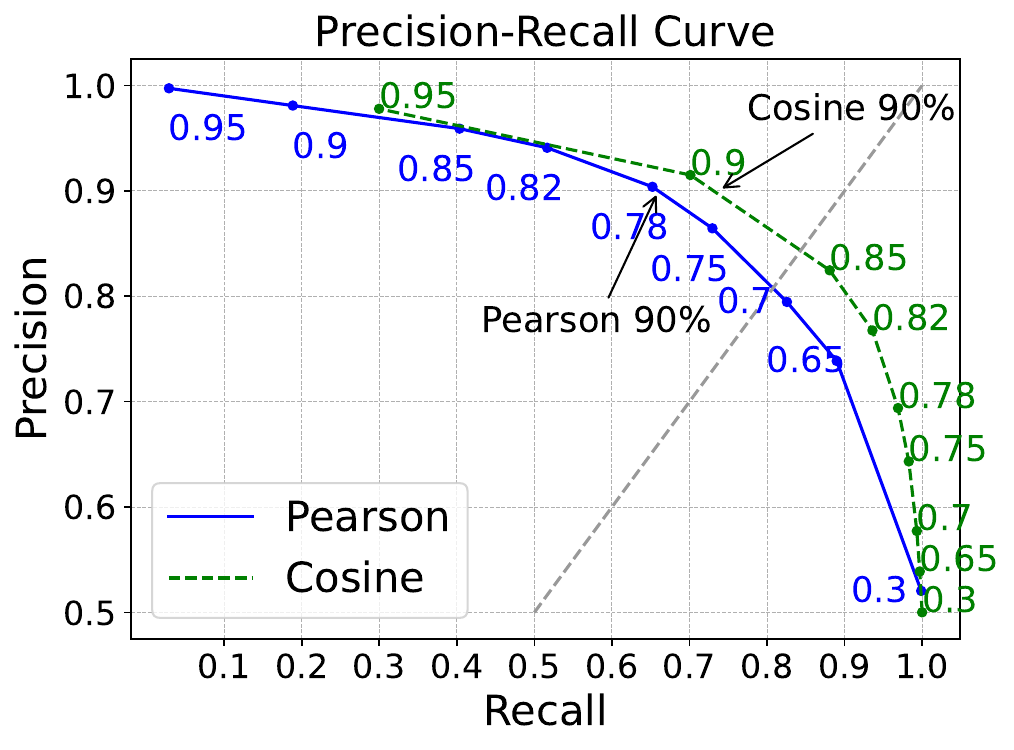}}%
\caption{Distribution of the insample and outsample using cosine and Pearson similarities, and the Precision-Recall curves for both.}
\label{fig:als-prcurive-bin}
\end{figure}

Building upon the previous example  using the MovieLens 100K data set, we set $\lambda_w=\lambda_z=0.15$ for  regularization and a rank of $62$ to minimize  RMSE. 
We aim to analyze the similarity between different movie hidden vectors, and the goal of Recommender 2 is to see whether the matrix factorization can help differentiate high-rated from low-rated movies, thereby recommending movies correlated with the user's high-rated ones.
Define further the term ``insample" as the similarity between the movies having rates $5$ for each user, and ``outsample" as the similarity between the movies having rates $5$ and $1$ for each user. Figure~\ref{fig:als-cosine} and \ref{fig:als-pearson} depict the bin plots of the distributions of insample and outsample under cosine and Pearson similarities, respectively.
In both scenarios, a clear distinction is observed between the distributions of the ``insample" and ``outsample" data, indicating that  ALS decomposition can actually find the hidden features of different movies for each user.
Figure~\ref{fig:als-prcurve} displays the \textit{precision-recall (PR) curve} for these scenarios, where we find  cosine similarity outperforms Pearson similarity, achieving over $73\%$ recall with  $90\%$ precision. However, Pearson similarity can  identify only about $64\%$ of the high-rated movies with the same precision. 
In practice, other measures, such as  \textit{negative Euclidean distance}, can also be explored. The Euclidean distance  measures the ``dissimilarity" between two vectors; and a negative value thus represents their similarity.

The ALS method for recommendation discussed here is designed for \textit{explicit data}, where the ratings provided by each user have a clear hierarchical meaning. In contrast, there are also recommendation systems for \textit{implicit data}, where the system automatically infers users' preferences by tracking their actions, such as which items they viewed, where they clicked, which products they purchased, or how long they spent on a web page.
In such cases,  ALS can be extended to more complex models, such as using a dictionary matrix to transform  the explicit data into user and item latent vectors \citep{he2017neural}, incorporating multinomial prior into a variational auto-encoder,  and enhancing the model's ability to handle implicit feedback by leveraging probabilistic modeling techniques \citep{liang2018variational}.

\begin{problemset}
	

\item \label{prob:als_pseudo1} \textbf{Least squares for rank-deficiency.} Let $\bA\in\real^{m\times n}$ and $\bb\in\real^m$. Show that the least squares problem $L(\bx)=\normtwo{\bA\bx-\bb}^2$ has a minimizer $\bx^*\in\real^n$ if and only if there exists a vector $\by\in\real^n$ such that $\bx^*=\bA^+\bb+(\bI-\bA^+\bA)\by$, where $\bA^+$ is the pseudo-inverse of $\bA$ (Appendix~\ref{appendix:pseudo-inverse_main}). 
\begin{itemize}
\item This shows that the least squares has a \textbf{unique} minimizer of $\bx^*=\bA^+\bb$ only when $\bA^+$ is a left inverse of $\bA$ (Definition~\ref{definition:one_side_inverse}). The solution in Lemma~\ref{lemma:ols} is a special case.
\item The optimal value is $L(\bx^*)=\bb^\top(\bI-\bA\bA^+)\bb$.
\item If $\by\neq \bzero$: $\normtwo{\bA^+\bb}\leq \normtwo{\bA^+\bb+(\bI-\bA^+\bA)\by}$.
\end{itemize}
\textit{Hint: See LS via SVD in Section~\ref{section:application-ls-svd}.}

\item \label{prob:als_pseudo2} \textbf{Least squares for rank-deficiency.} Let  $\bA\in\real^{m\times n}$ and $\bB\in\real^{m\times p}$. Show that the least squares problem $L(\bX) = \normf{\bA\bX-\bB}^2$ has a minimizer $\bX^*=\bA^+\bB\in\real^{n\times p}$. Determine all the minimizers using Problem~\ref{prob:als_pseudo1}. (Notice that  problems involving the regularization of   $\bX$ to be orthogonal are discussed in Theorem~\ref{theorem:ortho_procruste} and \ref{theorem:twos_si_pros}.)

\item \label{prob:als_pseudon} \textbf{Least squares for rank-deficiency.}  Let  $\bA\in\real^{m\times n}$ and $\bB\in\real^{p\times n}$. Show that the least squares problem $L(\bX) = \normf{\bX\bA-\bB}^2$ has a minimizer $\bX^*=\bB\bA^+\in\real^{p\times m}$.
	
\item Prove Lemma~\ref{lemma:als-update-w-rank}.

\item \label{prob:sep_conv} \textbf{Marginally convex.} Let $D(\bA, \bB)$ be convex in the second argument $\bB$. Show that $D(\bA,\bW\bZ)$ is convex in $\bW$ for a  fixed $\bZ$, and vice versa.

\item Show that a function that is jointly convex is necessarily marginally convex as well.

\item Derive the update of column-by-column fashion for Algorithm~\ref{alg:als-regularizer}.

\item \label{prob:non_iid_gaus} \textbf{Weighted $\ell_2$ loss from non i.i.d. Gaussian noise.} Suppose the Gaussian noise in \eqref{equation:gau_noise} is not i.i.d.. Discuss the likelihood function for the problem $\bA=\bW\bZ$. Show that the loss function is a \textit{weighted $\ell_2$ norm (or a weighted Frobenius norm)}: $L(\bW,\bZ)=\normf{\bW\hadaprod(\bA-\bW\bZ)}^2=\sum_{m,n=1}^{M,N} w_{mn}(a_{mn}-b_{mn})^2$ if $\bB\triangleq\bW\bZ=\{b_{mn}\}\in\real^{M\times N}$. 
Explain the relationship between the noise variance and the weight $w_{mn}$ for each $(m,n)$-th entry.

\item Show that the loss function in \eqref{equation:als_gama_los} is derived using the deviance discussed in \eqref{equation:als_poi_los}.

\item \label{prob:ortho_mf} \textbf{Orthogonal and projective  matrix factorization.} Consider the optimization $\mathopmin{\bW}\normf{\bA-\bW\bZ}^2$ such that $\bZ\bZ^\top=\bI_K$, where $\bA\in\real^{M\times N}, \bW\in\real^{M\times K}, \bZ\in\real^{K\times N}$, and $K\leq \min\{M,N\}$. Show that the optimal value $\bW^*$ given $\bZ$ is $\bA\bZ^\top$.
This indicates that the matrix factorization optimization can be equivalently stated as $\mathopmin{\bZ\bZ^\top=\bI_K}\normf{\bA-\bA\bZ^\top\bZ}^2$. And the relaxed problem is called the \textit{projective matrix factorization} \citep{yuan2005projective, yang2010linear}:
$$
\mathopmin{\bZ}\normf{\bA-\bA\bZ^\top\bZ}^2,
$$
where each row of $\bA$ is projected onto a $K$-dimensional subspace (Appendix~\ref{section:by-geometry-hat-matrix}), hence the name.
The interpretations of orthogonal and projective matrix factorizations are discussed in Problem~\ref{prob:ortho_nmf}.

\item \label{problem:rls} \textbf{Regularized least squares (RLS).} 
Given $\bA\in\real^{m\times n}, \bb\in\real^{m}, \bB\in\real^{p\times n}$, and $\lambda\in\real_{++}$, we consider the constrained least squares problem:
$$
\mathop{\min}_{\bx\in\real^n} \normtwo{\bA\bx-\bb}^2 + \lambda\normtwo{\bB\bx}^2.
$$
Show that the constrained least squares (RLS) problem has a unique solution if and only if $\nspace(\bA)\cap \nspace(\bB) = \{\bzero\}$.

\item \label{prob:denoise_rls} \textbf{Denoising via RLS.} Consider a noisy measurement of a signal $\bx\in\real^n$:
$
\by = \bx+\be,
$
where $\by$ is the observed measurement, and $\be$ is the noise vector. We want to find an estimate $\bx$ of the observed measurement $\by$ such that $\bx \approx \by$:
$
\min \normtwo{\bx-\by}^2.
$
Apparently, the optimal solution of this optimization is given by $\bx=\by$; however, it is meaningless.
To improve the estimate,  we can add a penalty term for the differences between  consecutive observations:
$
R(\bx) = \sum_{i=1}^{n-1} (x_i - x_{i+1})^2.
$
Then, 
\begin{itemize}
\item Find the constrained least squares representation for this problem and derive the constrained least squares solution.
\item Find some applications of this denoising problem. For example, when we model the profit and loss signal of a financial asset, the two observations over consecutive days of the underlying asset should exhibit smooth transitions rather than abrupt changes.
\end{itemize}


\item \textbf{Weighted least squares (WLS).}
Building upon the assumptions in Lemma~\ref{lemma:ols}, we consider further that each data point $i\in\{1,2,\ldots, m\}$ (i.e., each row of $\bA$) has a weight $w_i$. 
This means some  data points may carry greater significance than others, and we can produce approximate minimizers that reflect this.
Show that the value $\bx_{WLS} = (\bA^\top\bW^2\bA)^{-1}\bA^\top\bW^2\bb$ serves as the \textit{weighted least squares (WLS)}  estimate of $\bx$, where $\bW=\diag(w_1, w_2, \ldots, w_m)\in\real^{m\times m}$. \textit{Hint: Find the normal equation for this problem.}

\item \textbf{Positive definite  weighted least squares (PDWLS).}
Building upon the assumptions in Lemma~\ref{lemma:ols}, we consider further  the matrix equation $\bA\bx + \be =\bb$, where $\be$ is an error vector. Define the weighted error squared sum $E_w = \be^\top \bW \be$, where the weighting matrix $\bW$ is  positive definite. 
Show that the positive definite weighted least squares solution is $\bx^* = (\bA^\top\bW\bA)^{-1}\bA^\top\bW\bb$. \textit{Hint: Compute the gradient of $E_w = (\bb-\bA\bx)^\top\bW(\bb-\bA\bx)$.}

\item \textbf{Weighted color noise least squares.}
Building upon the assumptions in Lemma~\ref{lemma:ols}, we consider  the matrix equation $ \bA\bx + \be = \bb $, where $\be$ is an additive color noise vector satisfying the conditions $\Exp[\be] = \bzero$ and $\Exp[\be\be^\top] = \bSigma$, where $\bSigma$ is known. Use the weighting error function $E_w = \be^\top \bW \be$ as the cost function for finding the optimal estimate $\bx^*$. Show that
$\bx^* = (\bA^\top \bW \bA)^{-1} \bA^\top \bW \bb$,
where the optimal choice of the weighting matrix $\bW$ is $\bW^* = \bSigma^{-1}$.
\textit{Hint: Compute the gradient of $E_w = (\bb-\bA\bx)^\top\bW(\bb-\bA\bx)$.}

\item \label{problem:tls} \textbf{Transformed least squares (TLS).}
Building upon the assumptions in Lemma~\ref{lemma:ols}, we consider further the restriction $\bx=\bC\bgamma+\bc$, where $\bC\in\real^{n\times k}$ is a known matrix such that $\bA\bC$ has full rank, $\bc$ is a known vector, and $\bgamma$ is an unknown vector.
Show that the value $\bx_{TLS}=\bC(\bC^\top\bA^\top\bA\bC)^{-1}(\bC^\top\bA^\top)(\bb-\bA\bc) +\bc$ serves as the \textit{transformed least squares (TLS)} estimate of $\bx$.


\item \label{problem:twls2} Find the transformed weighted least squares estimate.

\index{First-order optimality condition}
\index{Fermat's theorem}
\item \label{problem:fist_opt} \textbf{First-order optimality condition for local optima points.} 
Consider  \textit{Fermat's theorem}: for a one-dimensional function $g(\cdot)$ defined and differentiable over an interval ($a, b$), if a point $x^*\in(a,b)$ is a local maximum or minimum, then $g^\prime(x^*)=0$. 
Prove the first-order optimality conditions for multivariate functions based on  Fermat's theorem for one-dimensional functions.
That is, let  $f: \sS\rightarrow \real$ be a function defined on a set $\sS\subseteq \real^n$. Suppose that $\bx^*\in\text{int}(\sS)$, i.e., in the interior point of the set, is a local optimum point and that all the partial derivatives (Definition~\ref{definition:partial_deri}) of $f$ exist at $\bx^*$. Then $\nabla f(\bx^*)=\bzero$, i.e., the gradient vanishes at all local optimum points. (Note that, this optimality condition is a necessary condition but not sufficient; however, there could be vanished points which are not local maximum or minimum point.)
\textit{Hint: Consider the one-dimensional function $g(t)=f(\bx^* + t\be_i)$ for $i\in\{1,2,\ldots,n\}$.}

\item \label{problem:pos_hessian} \textbf{Global minimum point of convex functions.} Let function $f$ be a twice continuously differentiable function defined over $\real^n$. Suppose that the Hessian $\nabla^2f(\bx) \succeq 0$ for any $\bx\in\real^n$ (i.e., the Hessian is always positive semidefinite~\footnote{Instead, if we assume the Hessian is positive semidefinite at a given point, then the point is a local minimum point.}).
This property is also referred to as the \textit{convexity}.
Show that $\bx^*$ is a global minimum point of $f$ if $\nabla f(\bx^*)=\bzero$. \textit{Hint: Use the linear approximation theorem in Theorem~\ref{theorem:linear_approx}.}

\item \textbf{Two-sided matrix least squares \citep{friedland2007generalized, aggarwal2020linear}.} Let $\bB$ be an $M\times K$ matrix and $\bC$ be a $P\times N$ matrix. Find the $K\times P$ matrix $\bX$ such that $L(\bX)=\norm{\bA - \bB\bX\bC}_F^2$ is minimized, where $\bA\in\real^{M\times N}$ is known. 
\begin{itemize}
\item Derive the derivative of $L$ with respect to $\bX$ and the optimality conditions. 
\item Show that one possible solution to the optimality conditions is $\bX^*=\bB^+\bA\bC^+$, where $\bB^+$ and $\bC^+$ are the pseudo-inverses of $\bB$ and $\bC$, respectively (see Appendix~\ref{appendix:pseudo-inverse_main}).
\end{itemize}
Similarly, consider the optimization with $\rank(\bX)\leq p$:
$
L(\bX)=\norm{\bA - \bB\bX\bC}_F^2$, s.t. 
$\rank(\bX)\leq p$.
Show that 
\begin{itemize}
\item One possible solution to this is $\bX^*=\bB^+\bA_p\bC^+$, where $\bA_p$ a truncated SVD of $\bB\bB^+\bA\bC^+\bC$ by replacing all but the $p$ largest singular values by zero.
\item  $\bX^*$ also minimizes $\normf{\bX}$, i.e., has the smallest magnitude among all  solutions.
\item $\bX^*$ is the \textbf{unique} solution if and only if either $\rank(\bB\bB^+\bA\bC^+\bC)\leq p$ or both $\rank(\bB\bB^+\bA\bC^+\bC)\geq p$ and $\sigma_{p+1}(\bB\bB^+\bA\bC^+\bC) < \sigma_{p}(\bB\bB^+\bA\bC^+\bC)$.
\end{itemize}

\item \label{problem:mono_gd} \textbf{Monotonic progress of gradient descent.} Consider the gradient descent for a differentiable function $f(\bx):\real^n\rightarrow \real$ that is $L$-strongly smooth (Definition~\ref{definition:scss_func}). Suppose the iterate $\bx^{(t+1)}$ is obtained from iterate $\bx^{(t)}$ by 
$
\bx^{(t+1)} = \bx^{(t)} - \eta\nabla f(\bx^{(t)}).
$
Show that 
\begin{itemize}
\item If the step size $\eta \leq \frac{2}{L}$,  the function value $f$ is nonincreasing: $f(\bx^{(t+1)})\leq f(\bx^{(t)})$.
\item If the step size $\eta \in [\frac{1}{2L}, \frac{1}{L}]$,  the gradient satisfies $\normtwo{\nabla f(\bx^{(t)})}\leq \epsilon$ after $T=\mathcalO(\frac{1}{\epsilon^2})$ steps.
\textit{Hint: Use the descent lemma for $L$-strongly smooth functions (Definition~\ref{definition:scss_func}).}
\end{itemize}

\item \label{problem:nuclear_equi} \citep{rennie2005fast, mazumder2010spectral} Consider the nuclear norm~\footnote{The nuclear norm is defined as the sum of singular values of a matrix and provides the tightest convex envelope of the rank function of a matrix (see Definition~\ref{definition:ky_fan_norm} and \citet{jain2017non}).} $\norm{\bA}_n$ of any matrix $\bA\in\real^{m\times n}$ with rank $r$. Show that 
$$
\norm{\bA}_n = \mathop{\min}_{\substack{\bW\in\real^{m\times r} \\ 
\bZ\in\real^{r\times n} \\
}}
\frac{1}{2} (\normf{\bW}^2 + \normf{\bZ}^2)
\gap \text{s.t.} \gap 
\bA=\bW\bZ
$$

\item Discuss the gradient descent updates for different regularizations in \S~\ref{section:regularization-extention-general}.

\end{problemset}

%% file: chapter-NMF.tex
\newpage
\chapter{Nonnegative Matrix Factorization (NMF)}\index{NMF}\label{chapter:nmf}
\begingroup
\hypersetup{
	linkcolor=structurecolor,
	linktoc=page,  
}
\minitoc \newpage
\endgroup

\index{Decomposition: NMF}
\index{Sparsity}
\index{Nonnegativity constraint}
\section{Nonnegative Matrix Factorization}
\lettrine{\color{caligraphcolor}I}
In the era of big data, extracting meaningful patterns and latent structures from high-dimensional data sets has become a central challenge in various scientific and technological fields. 
Singular value decomposition (SVD) is supported by robust theoretical foundations and is applicable in a wide range of contexts. However, it does have limitations; for example, when applied to a nonnegative matrix~\footnote{Nonnegative matrices possess unique properties in linear algebra and are crucial for theoretical analysis; see Problems~\ref{prob:nonn_lin_1}$\sim$\ref{prob:nonn_lin_12}.},  SVD can produce negative values, which can be challenging to interpret.
To overcome this limitation, \textit{nonnegative matrix factorization (NMF)} has emerged as a powerful and interpretable tool for dimensionality reduction, feature extraction, and discovering latent structures within complex data.
Early consideration of the NMF problem was due to \citet{paatero1994positive, cohen1993nonnegative}, who referred it as \textit{positive matrix factorization}. 
Later, \citet{lee2001algorithms} popularized the problem with the introduction of the \textit{multiplicative update} rule.

Following  the matrix factorization using the alternating least squares (ALS) method, we now consider algorithms for solving the NMF problem: 
\begin{itemize}
\item Given a nonnegative matrix $\bA\in \real_+^{M\times N}$ with rank $r$, find nonnegative matrix factors $\bW\in \real_+^{M\times K}$ and $\bZ\in \real_+^{K\times N}$ such that: 
$
\bA\approx\bW\bZ.
$
\end{itemize}

As mentioned in the ALS section, the fundamental challenge in linear data analysis involves transforming or decomposing a high-dimensional data vector into a linear combination of lower-dimensional vectors. This transformation captures the essential characteristics of the original data, making it suitable for pattern recognition tasks. Consequently, these lower-dimensional vectors are often referred to as \textit{``hidden vectors," ``pattern vectors," or ``feature vectors."}
When conducting data analysis, building models, and processing information, it is crucial to meet two primary requirements for a pattern vector:
\begin{itemize}
\item \textit{Interpretability}. Each component of a pattern vector should possess clear physical or physiological significance, allowing for a meaningful understanding of the data.
\item \textit{Statistical fidelity}. In cases where the data are reliable and contain minimal error or noise, the components of a pattern vector should effectively capture the variability within the data, reflecting its primary distribution of information.
\end{itemize}
The NMF approach addresses these issues in various applications. For example:
\begin{itemize}
\item In document collections, documents are represented as vectors, with each vector element indicating the frequency (often weighted) of a specific term within the document. Arranging these document vectors sequentially forms a nonnegative term-by-document matrix, which provides a numerical representation of the entire document collection.
\item In image collections, each image is depicted by a vector, where each vector element represents a pixel. The value of each element, a nonnegative number, reflects the intensity and color of the corresponding pixel, leading to a nonnegative pixel-by-image matrix.

\item In gene expression analysis, observations from gene sequences under different experimental conditions are compiled into gene-by-experiment matrices. These matrices encapsulate the variations in gene expression across experiments.

\item For item sets or recommendation systems, customer purchase histories or ratings for a selection of items are recorded in a nonnegative sparse matrix. This matrix efficiently captures the sparse nature of user interactions with a large number of potential items.
\end{itemize}

Unlike arbitrary linear combinations, the linear combinations in the NMF context involve only nonnegative weights of nonnegative template vectors (the columns of $\bW$). As a result, there are no effects such as destructive interference, where a positive component could be canceled out by adding a negative component.  Instead, the data vectors must be explained using purely constructive methods, involving only positive components.
The nonnegativity constraint  inherently imposes   sparsity, enabling the factorization to capture additive features, which is especially advantageous in applications where parts-based representations are meaningful. 
This constraint make it find applications in diverse fields such as text mining, image processing, document analysis, and bioinformatics, where the identified components often correspond to distinct parts or features.
For example, in image processing, NMF has proven valuable for tasks such as  object detection,  image segmentation, and facial recognition \citep{lee2001algorithms, gillis2014and, gillis2020nonnegative}. The decomposition into nonnegative components aligns with the intuitive notion that images are composed of identifiable parts.
In the topic recovery problem, each column of $\bA$ denotes a document; the NMF aligns with a soft clustering approach where each column of $\bW$ represents a topic, and the positive entries of each column of $\bZ$ represent the  positive weights of each document for those topics \citep{shahnaz2006document}.
On the other hand, a nonnegative matrix factorization $\bA\approx \bW\bZ$ can be applied for clustering algorithm. 
Specifically, the data vector $\ba_j$ is assigned to cluster $i$ if $z_{ij}$ is the largest element in column $j$ of $\bZ$ \citep{brunet2004metagenes, gao2005improving}.
See also the applications of the NMF in the survey paper \citet{berry2007algorithms}.
In conclusion, the popularity of NMF stems from its ability to automatically extract sparse and easily interpretable factors.

To measure the quality of the  approximation, we evaluate the loss by computing  the Frobenius norm of the difference between the original matrix and the approximation:
\begin{equation}\label{equation:frob_nmf}
L(\bW,\bZ) \triangleq D(\bA, \bW\bZ) = \frac{1}{2}\normf{\bW\bZ-\bA}^2,~\footnote{Note that we include a scaling factor of $\frac{1}{2}$ for easier discussion.}
\end{equation}
where $L(\bW,\bZ)$ indicates it is a loss function w.r.t. $\bW$ and $\bZ$, and $D(\bA, \bW\bZ)$ implies it is a distance/divergence between $\bA$ and $\bW\bZ$ (we will use the two notations interchangeably  when  necessary).
The Frobenius norm is arguably the most widely used norm for NMF because it corresponds to Gaussian additive noise, which is reasonable in many situations and allows for the design of particularly efficient algorithms; see Section~\ref{section:more_err_sta_als}.
For nonnegative data, Gaussian noise can be interpreted as a truncated version of Gaussian noise \citep{lu2023bayesian}.
In later sections, we will extend this approach to include more general  $\beta$-divergences (Section~\ref{section:beta_div_altmu}).


 When we want to find two nonnegative matrices $\bW\in\real^{M\times r}_+$ and $\bZ\in\real_+^{r\times N}$ such that $\bA=\bW\bZ$, the problem is known as the \textit{Exact NMF} of $\bA$ of size $r$. Exact NMF is NP-hard \citep{vavasis2010complexity, gillis2020nonnegative}.
Thus, we only consider the approximation of NMF here.

\index{Overfitting}
In the context of collaborative filtering, it is recognized  that  NMF via multiplicative updates can result in overfitting  despite favorable convergence results.
The overfitting can be partially mitigated through regularization, but its out-of-sample performance remains low. 
Bayesian optimization through the use of generative models, on the other hand, can effectively prevent overfitting in nonnegative matrix factorization \citep{brouwer2017comparative, lu2022flexible, lu2023bayesian}.
In the following sections, we introduce several methods for solving NMF problems and provide a brief discussion of their applications.

\index{Bayesian inference}
\index{Bayesian optimization}
\index{Bayesian matrix decomposition}

\begin{algorithm}[h] 
\caption{Projected Gradient Descent Method}
\label{alg:pgd_gen}
\begin{algorithmic}[1] 
\Require A function $f(\bx)$ and a set $\sS$; 
\For{$t=1,2,\ldots$}
\State Pick a step size $\eta_t$;
\State Set $\bx^{(t+1)} \leftarrow \mathcalP_{\sS}(\bx^{(t)} - \eta_t \nabla f(\bx^{(t)}))$;
\EndFor
\State Output final  $\bx$;
\end{algorithmic} 
\end{algorithm}
\section{NMF via Alternating Projected Gradient Descent (APGD)}\label{section:nmf_apgd}
The projected gradient descent (PGD, Algorithm~\ref{alg:pgd_gen}) addresses the minimization of a function over a set $\sS$:
$$
\mathopmin{\bx\in\sS} f(\bx).
$$
The orthogonal projection is defined as
$
\mathcalP_{\sS} (\bx) \triangleq \mathop{\argmin}_{\by\in\sS} \normtwo{\by-\bx}.
$
When $\sS$ is the nonnegative orthant, $\mathcalP_{\sS} (\bx)$ simplifies to $\mathcalP_{\sS} (\bx) = \max\{\bzero, \bx\}$, where the max operator is applied componentwise.

Therefore, the \textit{alternating PGD (APGD)} approach for NMF updates the factored components iteratively by 
$$
\bZ\leftarrow \max\bigg\{\bzero, \mathop{\argmin}_{\bZ\in\real^{K\times N}}\normf{\bW\bZ-\bA} \bigg\}
\gap\text{and}\gap
\bW\leftarrow \max\bigg\{\bzero, \mathop{\argmin}_{\bW\in\real^{M\times K}}\normf{\bW\bZ-\bA} \bigg\},
$$
where each update can be solved using a least squares method followed by projection onto the nonnegative orthant.
Due to the projection, the solution may not be properly scaled. A closed-form scaling factor $\gamma$ can improve the solution at each iteration:
$$
\gamma^* = \mathop{\argmin}_{\gamma\geq 0} \normf{\gamma\bW\bZ-\bA} 
=
\frac{\langle \bA, \bW\bZ\rangle}{\langle \bW\bZ, \bW\bZ\rangle}
=
\frac{\langle \bA\bZ^\top, \bW\rangle}{\langle \bW^\top\bW, \bZ\bZ^\top\rangle}.
$$
While it is generally not advised to use APGD due to its convergence challenges, APGD can be quite effective as an initialization method. This approach involves running a few iterations of APGD before switching to a different NMF algorithm, which is particularly beneficial for sparse matrices \citep{gillis2014and}.

\index{Nonnegative least squares}
\index{NNLS|see {Nonnegative least squares}}
\index{ANLS|see {Nonnegative least squares}}
\section{NMF via Alternating Nonnegative Least Squares (ANLS)}\label{section:nmf_anls}
The fundamental element of  the ALS approach is the least squares problem (Lemma~\ref{lemma:ols}). 
For NMF, we focus on the \textit{nonnegative least squares (NNLS)} problem:
\begin{equation}
\mathopmin{\bx\geq \bzero } f(\bx) = \mathopmin{\bx\geq \bzero } \frac{1}{2}\normtwo{\bb-\bM\bx}^2
\gap
\text{with }\bM\in\real^{m\times n}, \bb\in\real^m, \bx\in\real_+^n. 
\end{equation}
The KKT conditions (Appendix~\ref{appendix:KKT_proj}) imply the complementary slackness condition $\lambda_ix_i^*=0, \forall i$, where $\lambda_i$ is the Lagrangian multiplier; and the optimal condition $\nabla f(\bx^*) -\sum_{i}\lambda_i \be_i=\bzero$, where $\bx^*$ denotes the optimal solution of the NNLS problem.
Together, the complementary slackness and the optimal condition indicate that:
$$
\nabla f(\bx^*)
=
\sum_{i:x_i^*=0} \lambda_i \be_i.
$$
Therefore, we conclude that 
\begin{equation}\label{equation:kkv_nnn_raw}
(\textbf{KKT of NNLS})\gap 
	\bx^*\geq \bzero, 
	\gap
	\nabla f(\bx^*)\geq 0, 
	\gap
	\text{and}
	\gap
	x_i^* (\nabla f(\bx^*))_i=0,\, \forall i.
\end{equation}
These conditions imply sparsity when the nonnegative constraint     is applied, meaning the NNLS or NMF problem inherently  imposes a \textbf{sparsity constraint}.

Assume we are given the inactive set $I\subseteq \{1,2,\ldots,n\}$:
$$
I = \left\{ i \mid x_{i}^{*} > 0, \,\forall i \in\{1,2,\ldots,n\} \right\}.
$$
The complement of $I$ is the so-called \textit{active set}, where the corresponding constraints are active. That is, the active set contains  indices $i$ such that $x_{i}^{*} = 0$. The nonzero entries of $\bx^{*}$ can be determined by solving the following linear system:
$$
\begin{aligned}
[\nabla_{\bx} f(\bx)]_{I} = \bzero 
\gapthree\Longleftrightarrow\gapthree [\bM^\top(\bM \bx - \bb)]_{I} = \bzero 
\gapthree\Longleftrightarrow\gapthree \bM[:, I]^\top \bM[:, I] \bx[I] = \bM[:, I]^\top \bb.
\end{aligned}
$$
These are the normal equations of the unconstrained least squares problem for $\bx[I]$, that is,
\[
\min_{\bx[I]} \frac{1}{2}\normtwo{\bb -\bM[:, I] \bx[I]}^2.
\]
This leads to  the \textit{active-set method}, which iteratively updates the active set through pivoting (that is, entering and removing variables from the active set) to ensure the objective function  decreases \citep{lawson1995solving}; see Algorithm~\ref{alg:nmf_anls}. 

\paragraph{Alternating nonnegative least squares (ANLS).} 
Once we have the active-set method for NNLS problems,  NMF can be achieved by replacing OLS in ALS algorithms with NNLS, known as \textit{alternating nonnegative least squares (ANLS)} \citep{kim2011fast}.
Given a fixed $\bW$, the NMF can be solved for each column of $\bZ$ separately:
$$
\frac{1}{2}\normf{\bA-\bW\bZ}^2 =
\frac{1}{2}\sum_{n=1}^{N}\normtwo{\ba_n - \bW\bz_n}^2,
$$
where each subproblem $\mathopmin{\bz_n\geq \bzero}\normtwo{\ba_n - \bW\bz_n}^2$ can be solved using NNLS.
Since the NMF problem is symmetric: $\bA=\bW\bZ$ if and only if $\bA^\top=\bZ^\top\bW^\top$ such that $D(\bA, \bW\bZ)=D(\bA^\top, \bZ^\top\bW^\top)$. The analysis of optimizing $\bW$ given $\bZ$  follows directly from the previous methodology.
We should also note that  since the initial guess of $\bW$ and $\bZ$ typically offers a poor approximation of  $\bA$, solving the NNLS subproblems exactly in the early stages of the alternating algorithms is often unnecessary. 
Instead, it can be more efficient to use ANLS as a refinement step within a less computationally expensive NMF algorithm, such as  APGD or MU (discussed in later sections)

\begin{algorithm}[h] 
\caption{Nonnegative Least Squares (NNLS) via Active-Set Method}
\label{alg:nmf_anls}
\begin{algorithmic}[1] 
\Require A real-valued matrix $\bM\in\real^{m \times n}$, a real-valued vector $\bb\in\real^m$;
\State Initialize index sets $I = \emptyset$ and $J = \{1, \ldots, n\}$;
\State Initialize unknown $\bx\in\real^n$ to an all-zero vector and let $\bw \leftarrow \bM^\top(\bb - \bM\bx)$;
\State Let $\bw[J]$ denote the sub-vector with indices from $J$;
\State Choose a stopping criterion on the approximation error $\delta$;
\State Choose the maximum number of iterations $C$;
\State $iter=0$; \Comment{Count for the number of iterations}
\While{$J \neq \emptyset$ and $\max(\bw[J]) > \delta$ and $iter<C$}
\State $iter=iter+1$; 
\State Let $j$ in $J$ be the index of $\max(\bw[J])$ in $\bw$: $j=\mathop{\argmax}_{j\in J} w_j$;
\State Add $j$ to $I$  and remove $j$ from $J$ such that $I\cup J = \{1,2,\ldots,n\}$;
\State Let $\bM[:,I]$ be $\bM$ restricted to the variables/columns included in $I$;
\State \parbox[t]{\dimexpr\linewidth-\algorithmicindent}{Let $\bs$ be vector of same length as $\bx$;
Let $\bs[I]$ denote the sub-vector with indices from $I$, and let $\bs[J]$ denote the sub-vector with indices from $J$;}
\State Set $\bs[I] \leftarrow ((\bM[:,I])^\top \bM[:,I])^{-1} (\bM[:,I])^\top \bb$ and $\bs[J]$ to zero;
\While{$\min(\bs[I]) \leq 0$}
\State Let $\alpha \leftarrow \min \frac{x_i}{x_i - s_i}$ for $i$ in $I$ where $s_i \leq 0$;
\State Set $\bx\leftarrow \bx + \alpha(\bs - \bx)$;
\State Move to $J$ all indices $j$ in $I$ such that $x_j \leq 0$;
\State Set $\bs[I] \leftarrow ((\bM[:,I])^\top \bM[:,I])^{-1} (\bM[:,I])^\top \bb$;
\EndWhile
\State Set $\bs[J]$ to zero;
\State Set $\bx\leftarrow \bs$;
\State Set $\bw\leftarrow \bM^\top(\bb - \bM\bx)$;
\EndWhile
\State Output $\bx$;
\end{algorithmic} 
\end{algorithm}

\index{Hierarchical ANLS}
\section{NMF via Hierarchical Alternating Nonnegative Least Squares}
Let $\ba, \bb\in\real_+^n$ be two nonnegative vectors, then the \textit{univariate NNLS} problem can be formulated as 
$$
\mathopmin{x\geq 0} \normtwo{\ba-x\bb}^2.
$$
The solution of $x$ can be obtained in closed-form: $x=\max\big\{0, \frac{\bb^\top\ba}{\normtwo{\bb}^2}\big\}$ if $\normtwo{\bb}\neq 0$.
Having this univariate NNLS solution in mind, we consider the $k$-th row of $\bZ$ for $k\in\{1,2,\ldots,K\}$, the subproblem in NMF is
\begin{equation}\label{equation:llipschi_hianls}
\mathopmin{\bZ[k,:]\geq \bzero} \bigg\Vert\underbrace{\big(\bA-\sum_{p\neq k}^{K} \bW[:,p]\bZ[p,:]\big)}_{\triangleq\bA_k} - \bW[:,k]\bZ[k,:]\bigg\Vert_F^2, 
\gap \forall k,~\footnote{This subproblem is convex and is $L$-Lipschitz gradient continuous/$L$-strongly smooth (Definition~\ref{definition:scss_func}); see Problem~\ref{prob:llipschi_hianls}.}
\end{equation}
which indicates the entries in a row of $\bZ$ do not interact (similarly, entries in a column of $\bW$ do not interact). Therefore, the optimization of each entry in a row of $\bZ$ can be decoupled. 
Let $\bA_k\triangleq\big(\bA-\sum_{p\neq k}^{K} \bW[:,p]\bZ[p,:]\big)$, the NMF problem becomes a set of rank-one updates on $\bA_k$, for $k\in\{1,2,\ldots,K\}$.
The solution is 
$$
\bZ^*[k,:]=\mathop{\argmin}_{\bZ[k,:]\geq \bzero} \normf{\bA_k -\bW[:,k]\bZ[k,:]}^2
=
\max\left(
\bzero, 
\frac{\bW[:,k]^\top\bA_k}{\normtwo{\bW[:,k]}^2}
\right),
\gap \forall k,
$$
where the max operator is applied componentwise.
This leads to the \textit{hierarchical ANLS (Hi-ANLS)} solution for NMF problems, which iteratively solves a univariate NNLS problem.
The procedure is described in Algorithm~\ref{alg:hie_anls}, where we note that $\bZ[k,:]^\top = \bZ^\top[:,k]$.
In the algorithm, we update the $k$-th row of $\bZ$ and $k$-th column of $\bW$ in an interleaved manner. \citet{gillis2012accelerated} show that updating $\bZ$ several times before updating $\bW$  can significantly improve the performance  since this reuses the results of $\bW^\top\bA$ and $\bW^\top\bW$.

\begin{algorithm}[h] 
\caption{NMF via Hierarchical Alternating Nonnegative Least Squares (Hi-ANLS)}
\label{alg:hie_anls}
\begin{algorithmic}[1] 
\Require Matrix $\bA\in \real_+^{M\times N}$;
\State Initialize $\bW\in \real_{++}^{M\times K}$, $\bZ\in \real_{++}^{K\times N}$ randomly with positive entries;
\State Choose a stoping criterion on the approximation error $\delta$;
\State Choose maximal number of iterations $C$;
\State $iter=0$; \Comment{Count for the number of iterations}
\While{$\normf{\bA- (\bW\bZ)}^2>\delta $ and $iter<C$}
\State $iter=iter+1$;  
\For{$k=1$ to $K$}
\State $
\bZ[k,:]\leftarrow 
\max\left(
\bzero, 
\frac{\bW[:,k]^\top\bA_k}{\normtwo{\bW[:,k]}^2}
\right)$; \Comment{$\bA_k\triangleq\big(\bA-\sum_{p\neq k}^{K} \bW[:,p]\bZ[p,:]\big)$}

\State $
\bW[:, k]\leftarrow 
\max\left(
\bzero, 
\frac{\bA_k\bZ[k,:]^\top}{\normtwo{\bZ[k,:]}^2}
\right)$;
\EndFor
\EndWhile
\State Output $\bW,\bZ$;
\end{algorithmic} 
\end{algorithm}

\section{NMF via Alternating Direction Methods of Multipliers (ADMM)}\label{section:nmf_admm_all}
We briefly introduce the \textit{alternating direction methods of multipliers (ADMM)} and then  discuss its applications in matrix factorization and NMF.
\paragraph{ADMM.}
ADMM addresses the following convex optimization problem:
\begin{equation}\label{equation:admm_prob}
\mathopmin{\bx, \bz} f(\bx)+g(\bz), \gap \text{s.t.}\gap \bD\bx+\bE\bz=\bff.
\end{equation}
Given a parameter $\rho>0$, the \textit{augmented Lagrangian} of \eqref{equation:admm_prob} is 
\begin{equation}
L_\rho (\bx, \bz, \bl) = f(\bx)+g(\bz) +\langle \bl, \bD\bx+\bE\bz-\bff\rangle + \frac{\rho}{2}\normtwo{\bD\bx+\bE\bz-\bff}^2.
\end{equation}
When $\rho=0$, the augmented Lagrangian reduces to the Lagrangian function; when $\rho>0$, the augmented Lagrangian function  acts as  a penalized version of the Lagrangian function.
The \textit{augmented Lagrangian method} solves the problem by performing the following steps  (at the $(t+1)$-th iteration):
$$
\text{augmented Lagrangian:}
\gap 
\left\{
\begin{aligned}
(\bx^{(t+1)}, \bz^{(t+1)}) &\in \mathop{\argmin}_{\bx, \bz} L_\rho (\bx, \bz, \bl);\\
\bl^{(t+1)}&=\bl^{(t)} + \rho(\bD\bx^{(t+1)}+\bE\bz^{(t+1)} -\bff ),
\end{aligned}
\right.
$$
where the update on $\bl^{(t+1)}$ is derived from the \textit{conjugate subgradient theorem} (see, for example, \citet{bach2011convex}), and 
the symbol `$\in$' indicates that the minimum points may not be uniquely determined.
One source of difficulty is the coupling term between the $\bx$ and  $\bz$ variables, which is of the form $\rho(\bx^\top\bD^\top\bE\bz)$.
ADMM tackles this difficulty by replacing the exact minimization of $(\bx,\bz)$ with one iteration of the alternating minimization method.
To be more specific, for the  $(t+1)$-iteration, the solution of ADMM takes the following form:
\begin{equation}
\text{ADMM:}\gap
\left\{
\begin{aligned}
	\bx^{(t+1)}&\in\mathop{\argmin}_{\bx} \left\{ f(\bx)+ \frac{\rho}{2}\normtwo{\bD\bx+\bE\bz^{(t)} -\bff +\frac{1}{\rho}\bl^{(t)}}^2 \right\};\\
	\bz^{(t+1)}&\in\mathop{\argmin}_{\bz} \left\{ g(\bz)+ \frac{\rho}{2}\normtwo{\bD\bx^{(t+1)}+\bE\bz -\bff +\frac{1}{\rho}\bl^{(t)}}^2 \right\};\\
	\bl^{(t+1)}&=\bl^{(t)} + \rho(\bD\bx^{(t+1)}+\bE\bz^{(t+1)} -\bff ).
\end{aligned}
\right.
\end{equation}
Let $\widetildebl\triangleq\frac{1}{\rho}\bl$, this can be equivalently stated as (the one we will use in the sequel):
\begin{equation}\label{equation:admm_gen_up}
\text{ADMM:}\gap
\left\{
\begin{aligned}
\bx^{(t+1)}&\in\mathop{\argmin}_{\bx} \left\{ f(\bx)+ \frac{\rho}{2}\normtwo{\bD\bx+\bE\bz^{(t)} -\bff +\widetildebl^{(t)}}^2 \right\};\\
\bz^{(t+1)}&\in\mathop{\argmin}_{\bz} \left\{ g(\bz)+ \frac{\rho}{2}\normtwo{\bD\bx^{(t+1)}+\bE\bz -\bff +\widetildebl^{(t)}}^2 \right\};\\
\widetildebl^{(t+1)}&=\widetildebl^{(t)} + (\bD\bx^{(t+1)}+\bE\bz^{(t+1)} -\bff ).
\end{aligned}
\right.
\end{equation}
That is, iteratively updating $\bx, \bz$, and $\bl$.

\paragraph{ADMM applied to matrix factorization.}
We return to the problem discussed in ALS (Equation~\eqref{equation:als-per-example-loss2}, i.e., matrix factorization with Frobenius norm; not necessarily a NMF problem) together with a regularization function $r(\bZ)$:
$$
\mathopmin{\bZ} \frac{1}{2}\normf{\bA-\bW\bZ}^2+r(\bZ).
$$
The problem can be equivalently stated with an auxiliary variable $\widetildebZ\in\real^{K\times N}$:
\begin{equation}\label{equation:mf_admm_prob1}
\mathopmin{\bZ} \frac{1}{2}\normf{\bA-\bW\bZ}^2+r(\widetildebZ), 
\gap 
\text{s.t.}
\gap 
\bZ=\widetildebZ.
\end{equation}
Following \eqref{equation:admm_gen_up}, let a. \{$\bx\leftarrow \bZ$, $\bz\leftarrow \widetildebZ$, $\widetildebl\leftarrow \bL$, $\bD=-\bI$,  $\bE=\bI$\} or b. \{$\bx\leftarrow \bZ$, $\bz\leftarrow \widetildebZ$, $\widetildebl\leftarrow \bL$, $\bD=\bI$,  $\bE=-\bI$\}, 
the ADMM update  for \eqref{equation:mf_admm_prob1} is 
\begin{equation}\label{equation:admm_gen_als}
\left\{
\begin{aligned}
\bZ 
&\stackrel{(a)}{\leftarrow} (\bW^\top\bW+\rho \bI)^{-1} \left[ \bW^\top\bA +\rho(\widetildebZ+\bL) \right]
&\stackrel{(b)}{\leftarrow}& (\bW^\top\bW+\rho \bI)^{-1} \left[ \bW^\top\bA +\rho(\widetildebZ-\bL) \right];\\
\widetildebZ
&\stackrel{(a)}{\leftarrow}\mathop{\argmin}_{\widetildebZ} r(\widetildebZ) + \frac{\rho}{2}\normf{-\bZ+\widetildebZ + \bL}^2
&\stackrel{(b)}{\leftarrow}&\mathop{\argmin}_{\widetildebZ} r(\widetildebZ) + \frac{\rho}{2}\normf{\bZ-\widetildebZ + \bL}^2\\
\bL&\stackrel{(a)}{\leftarrow}\bL -\bZ+\widetildebZ &\stackrel{(b)}{\leftarrow}& \bL +\bZ-\widetildebZ.
\end{aligned}
\right.
\end{equation}
In practice, the Cholesky decomposition of $(\bW^\top\bW+\rho \bI)$ can be calculated such that the update can be obtained by forward and backward substitutions.
The update on $\bW$ can be obtained similarly due to symmetry.
We will consider the setting (a) in \eqref{equation:admm_gen_als} for the following discussions.

%

\paragraph{ADMM applied to $\ell_1$ regularization.}
We may also consider the $\ell_1$ regularization (see Section~\ref{section:regularization-extention-general}): $r(\widetildebZ)=\lambda \Vert\widetildebZ\Vert_1$. The update for each element $(k,n)$ of $\widetildebZ$ is $\widetildez_{kn}\leftarrow \max(0, 1-\frac{\lambda}{\rho} \abs{h_{kn}}^{-1}) h_{kn}$ for all $k\in\{1,2,\ldots,K\}$ and $n\in\{1,2,\ldots,N\}$, where $h_{kn} = z_{kn}-l_{kn}$ (i.e., the elements of $\bH=\bZ-\bL$). 

\paragraph{ADMM applied to smoothness/denoising regularization.}
The smoothness regularization on $\bZ$ can be defined as $r(\widetildebZ)=\frac{\lambda}{2} \Vert\bT\widetildebZ^\top\Vert_F^2$, where $\bT$ is an $N\times N $ tridiagonal matrix with 2 on the main diagonal and $-1$ on the superdiagonal and subdiagonal. This regularization ensures the proximal components in each row of $\widetildebZ$ is smooth (see Problem~\ref{prob:denoise_rls}). The update  of $\widetildebZ$ becomes $\widetildebZ\leftarrow \rho\bZ(\lambda \bT^\top\bT +\rho\bI)^{-1}$ \citep{huang2016flexible}.

\paragraph{ADMM applied to NMF.}
The NMF with ADMM is achieved simply by replacing $r(\bZ)$ with an indicator function.
The update on $\widetildebZ$ becomes $\max\left(\bzero, \bZ-\bL\right)$, where the max operator is applied  componentwise.
However, unlike the approaches introduced above and the MU approach in the next section, the updates from the ADMM approach are usually not monotonically nonincreasing.


\index{Alternating update}
\index{Kullback-Leibler divergence}
\index{Multiplicative update}
\section{NMF via Multiplicative Update (MU)}\label{section:nmf_frob_mu}
We further consider  NMF via an alternative alternating update. The hidden features in $\bW$ and $\bZ$ are modeled as nonnegative vectors in low-dimensional space. These latent vectors are randomly initialized and iteratively updated via an alternating multiplicative update rule to minimize the Frobenius norm distance between the observed and modeled matrices. 
Following Section~\ref{section:als-netflix}, we consider the low-rank with $K$ components; given $\bW\in \real_+^{M\times K}$, we aim to update $\bZ\in \real_+^{K\times N}$. The gradient with respect to $\bZ$ is given by Equation~\eqref{equation:givenw-update-z-allgd}:
$
\begin{aligned}
\nabla_{\bZ} L(\bW, \bZ) =\bW^\top(\bW\bZ-\bA) \in \real^{K\times N}.
\end{aligned}
$
Applying the gradient descent idea discussed in Section~\ref{section:als-gradie-descent}, the straightforward update for $\bZ$ can be achieved by 
$$
(\text{GD on $\bZ$})\gap \bZ \leftarrow \bZ - \eta \left(\nabla_{\bZ} L(\bW, \bZ)\right)=\bZ - \eta \nabla_{\bZ} L(\bW, \bZ),
$$
where $\eta$ represents a small positive step size. 
\paragraph{Multiplicative update (MU).}
If we allow a different step size for each entry of $\bZ$, the update can be written as:
$$
(\text{GD$^\prime$ on $\bZ$})\gap 
\begin{aligned}
	z_{kn} &\leftarrow z_{kn} - \frac{\eta_{kn}}{2} \left(\nabla_{\bZ} L(\bW, \bZ)\right)_{kn}
	=z_{kn} - \eta_{kn}(\bW^\top\bW\bZ-\bW^\top\bA)_{kn}, \,\, \forall k,n,
\end{aligned}
$$
where $z_{kn}$ denotes the $(k,n)$-th entry of $\bZ$. 
To proceed, we further rescale the step size:
$$
\eta_{kn} = \frac{z_{kn}}{(\bW^\top\bW\bZ)_{kn}}.
$$
Then we obtain the update rule:
\begin{equation}\label{equation:multi-update-z}
(\text{MU on $\bZ$})\gap 
\bZ \leftarrow \bZ\hadaprod \frac{[\bW^\top\bA]}{[\bW^\top\bW\bZ]}
\stackrel{*}{=}
\bZ - \frac{[\bZ]}{[\bW^\top\bW\bZ]}\hadaprod \nabla_{\bZ} L(\bW, \bZ) 
,
\end{equation}
where $\frac{[\cdot]}{[\cdot]}$ represents the componentwise division between two matrices. This is known as the \textit{multiplicative update (MU)}, and is first developed in \citet{lee2001algorithms} for NMF problems. 
Analogously, the multiplicative update for $\bW$ can be obtained by 
\begin{equation}\label{equation:multi-update-w}
(\text{MU on $\bW$})\gap
\bW \leftarrow \bW \hadaprod \frac{[\bA\bZ^\top]}{[\bW\bZ\bZ^\top]} 
\stackrel{*}{=} 
\bW - \frac{[\bW]}{[\bW\bZ\bZ^\top]}\hadaprod \nabla_{\bW} L(\bW, \bZ) .
\end{equation}
The factors $\frac{(\bW^\top\bA)_{kn}}{(\bW^\top\bW\bZ)_{kn}}$ and $\frac{(\bA\bZ^\top)_{mk}}{(\bW\bZ\bZ^\top)_{mk}}$ for all $m,k,n$ in \eqref{equation:multi-update-z} and \eqref{equation:multi-update-w} are called \textit{multiplicative factors}.
When $\bA=\bW\bZ$, these multiplicative factors reduce to one, indicating that the corresponding gradients vanish.

\paragraph{MU vs gradient descent.}
From the above construction, it is revealed that multiplication update algorithms are fundamentally similar to gradient descent algorithms, but they differ in their step size selection. With an appropriate choice of step size, a multiplication algorithm can transform the subtraction update rule of the standard gradient descent method into a multiplicative update rule.

In the gradient descent algorithm, a fixed or adaptive step length is typically used, and this step length is independent of the specific variable being updated. In other words, the step size may vary over time, but at any given update step, all entries of the variable matrix are updated using the same step size.
In contrast, the multiplication algorithm uses different step sizes ($\eta_{kn}$ above) for different entries of the factor matrix. This means that the step length is adaptive to each matrix entry. This adaptability is a key reason why the multiplication algorithm can outperform the gradient descent algorithm.

\paragraph{KKT conditions for NMF with Frobenius norm.}
The KKT conditions (Appendix~\ref{appendix:KKT_proj}) indicate that (see derivation in \eqref{equation:kkv_nnn_raw}):
\begin{equation}\label{equation:nmf_fro_kkt1}
\begin{aligned}
\bZ\geq \bzero,\gap&\nabla_{\bZ} L(\bW,\bZ)&\geq& \bzero, \gap \langle \bZ, \nabla_{\bZ} L(\bW,\bZ)\rangle &=&\bzero_{K\times N}; \\
\bW\geq \bzero,\gap&\nabla_{\bW} L(\bW,\bZ)&\geq& \bzero, \gap \langle \bW, \nabla_{\bW} L(\bW,\bZ)\rangle &=&\bzero_{M\times K}.
\end{aligned}
\end{equation}
This also implies
\begin{equation}
\begin{aligned}
\min\{\bZ, \nabla_{\bZ} L(\bW,\bZ)\} = \bzero_{K\times N}
\gap \text{and}\gap 
\min\{\bW, \nabla_{\bW} L(\bW,\bZ)\} = \bzero_{M\times K},
\end{aligned}
\end{equation}
where the min operator $\min\{\cdot, \cdot \}$ is applied componentwise. Any pair $(\bW,\bZ)$ satisfying the KKT conditions is a stationary point of the NMF problem in \eqref{equation:frob_nmf}.

\paragraph{Problems in MU.}
Equality ($*$) in \eqref{equation:multi-update-z} indicates a rescaled gradient descent update in the MU rules, which also implies
$$
\frac{[\bW^\top\bA]_{kn}}{[\bW^\top\bW\bZ]_{kn}}\geq 1 \gap\Longleftrightarrow\gap (\nabla_{\bZ} L(\bW, \bZ))_{kn}\leq 0, \gap  \forall k, n.
$$
Therefore, the MU algorithm induces three-fold rules: (i) Increase if its partial derivative is negative; (ii) Decrease it if its partial derivative is positive; (iii) Keep it unchanged if its partial derivative is zero.
However, if an element of $\bZ$ is equal to zero, the MU rule cannot modify it. Therefore, it is possible for an entry of $\bZ$  to be zero while its partial derivative is negative, which would not satisfy the KKT conditions in \eqref{equation:nmf_fro_kkt1}.
As a result, the iterates from the MU rule are not guaranteed to converge to a stationary point. There are several ways to address this issue, such as: (i) Using a small positive lower bound for the entries of $\bZ$ and  $\bW$\citep{gillis2012accelerated}.
(ii) Using the MU rule while reinitializing zero entries of $\bZ$ and  $\bW$ to a small positive constant when their partial derivatives become negative \citep{chi2012tensors}.

We now prove the monotonicity of the MU rule for NMF problems.
\begin{theorem}[Monotonically Nonincreasing of Multiplicative Update]\label{theorem:conv_mu_fro}
The loss $L(\bW,\bZ)=\frac{1}{2}\normf{\bW\bZ-\bA}^2$ remains nonincreasing under the following multiplicative update rules:~\footnote{More general results for $\beta$-divergences are discussed in Theorem~\ref{theorem:conv_mu_beta}.}
$$
\begin{aligned}
\bZ &\leftarrow \bZ\hadaprod \frac{[\bW^\top\bA]}{[\bW^\top\bW\bZ]}
\gap\text{and}\gap 
\bW &\leftarrow \bW \hadaprod \frac{[\bA\bZ^\top]}{[\bW\bZ\bZ^\top]},
\end{aligned}
$$
where $\bA\in\real_+^{M\times N}, \bW\in\real_+^{M\times K}$, and $\bZ\in\real_+^{K\times N}$. 
The operator $\frac{[\cdot]}{[\cdot]}$ represents the componentwise division between two matrices, and $\hadaprod$ denotes the Hadamard product between two matrices.

The MU update requires that $\bZ$ and $\bW$ should be initialized with positive (nonzero) entries; otherwise, the MU will not modify any entries due to the Hadamard product.
\end{theorem}

The MU approach for NMF sets a spark on NMF that has ignited a fire that has been growing steadily since then, becoming a cornerstone in the NMF literature for several reasons: (i) The update rules are extremely easy to implement; (ii) In practice, the convergence is relatively fast compared to many other methods; (iii) One major benefit of using MU rules is that the nonnegativity constraints are enforced automatically.
To prove the monotonicity of the MU rules, we need the following results.
\begin{definition}[Auxiliary Function (Majorizer)]\label{definition:aux_func}
$G(\bx, \widetildebx)$ is called an auxiliary function for $F(\bx)$ (or a majorizer of $F$ at $\widetildebx$) if the conditions~\footnote{$\bx$ can be scalars, vectors, or matrices.}
$$
G(\bx, \widetildebx) \geq F(\bx)
\gap 
\text{and}
\gap 
G(\bx, \bx) = F(\bx)
$$
are satisfied.
In other words, the auxiliary function $G(\bx, \widetildebx)$ is an upper bound of $F(\bx)$ and the bound is tight when $\widetildebx=\bx$.
\end{definition}

\begin{lemma}[Nonincreasing in Auxiliary Functions]\label{lemma:noninmuaux}
If $ G $ is an auxiliary function for $F$, then $F$ is nonincreasing under the update
\begin{equation}\label{equation:aux_update}
\bx^{(t+1)} = \mathop{\argmin}_{\bx} \, G(\bx, \bx^{(t)}).
\end{equation}
\end{lemma}
\begin{proof}[of Lemma~\ref{lemma:noninmuaux}]
The definition of the auxiliary function indicates that $ F(\bx^{(t+1)}) \leq G(\bx^{(t+1)}, \bx^{(t)}) \leq G(\bx^{(t)}, \bx^{(t)}) = F(\bx^{(t)})$.
\end{proof}

Note that $ F(\bx^{(t+1)}) = F(\bx^{(t)})$ only if $ \bx^{(t)}$ is a local minimum of $ G(\bx, \bx^{(t)})$ w.r.t. $\bx$. If the partial derivatives of $ F$ exist and are continuous in a small neighborhood of $ \bx^{(t)}$, this also implies that the gradient $ \nabla F(\bx^{(t)}) = \bzero$. Thus, by iterating the update in \eqref{equation:aux_update}, we obtain a sequence of estimates that converge to a local minimum $ \bx_{\min} = \argmin_{\bx} F(\bx)$ of the objective function:
\begin{equation}
 F(\bx^{(0)}) \geq 	 F(\bx^{(1)}) \geq  F(\bx^{(2)})\geq \ldots  \geq  F(\bx^{(t)}) \geq  F(\bx^{(t+1)})\geq \ldots \geq F(\bx_{\min}).
\end{equation}
Definition~\ref{definition:aux_func} finds a majorizer $G$ of $F$, and Lemma~\ref{lemma:noninmuaux} shows the minimization property in $G$, hence the algorithm is often referred to as the \textit{majorization-minimization (MM) framework}.
The update benefits when the global minimizer of $G$ has a closed-form solution or can be computed efficiently.

Therefore, if we can find an appropriate auxiliary function $ G(\bx, \bx^{(t)})$ for both variables in $\normf{\bA-\bW\bZ}$, the update rules in  Theorem~\ref{theorem:conv_mu_fro} follow from  \eqref{equation:aux_update}.
To apply the auxiliary function to the NMF problem, we consider a column in $\bA$ or $\bZ$: $\ba\triangleq\ba_n$ and $\bz\triangleq\bz_n$  in the following lemma, where $n\in\{1,2,\ldots,N\}$.

\begin{lemma}[Auxiliary Function for NMF]\label{lemma:aux_nmf}
Let $\bW\in\real^{K\times N}, \ba\in\real^M$, and $\bz\in\real^{K}$.
Let further $\bD\in\real^{K\times K}$ be a diagonal matrix with the $(k,k)$-th entry being $d_{kk}=\frac{(\bW^\top\bW \bz)_k}{z_k}=\frac{\bw_k^\top\bW\bz}{z_k} = \frac{\sum_{j=1}^{K} (\bW^\top\bW)_{kj}z_j}{z_k}, \, \,\forall k\in\{1,2,\ldots,K\}$, where $\bw_k$ is the $k$-th column of $\bW$ and $z_k$ is the $k$-th component of $\bz$.
Then, the following function is an auxiliary function for  $F(\bz)=\frac{1}{2}\normtwo{\ba-\bW\bz}^2$:
\[ 
G(\bz, \bz^{(t)}) = F(\bz^{(t)}) + (\bz-\bz^{(t)})^\top \nabla F(\bz^{(t)}) +\frac{1}{2}(\bz-\bz^{(t)})^\top \bD (\bz-\bz^{(t)}).
\]
\end{lemma}
\begin{proof}[of Lemma~\ref{lemma:aux_nmf}]
Since the third-order partial derivatives of $F(\bz)$ vanish (see Problem~\ref{prob:third_order_nmf}), $F(\bz)$ can be factored as 
$$
F(\bz) = F(\bz^{(t)})+(\bz-\bz^{(t)})^\top \nabla F(\bz^{(t)}) +\frac{1}{2}(\bz-\bz^{(t)})^\top \bW^\top\bW (\bz-\bz^{(t)}).
$$
Apparently, $G(\bz, \bz)=F(\bz) $. To complete the proof, we need to show that $G(\bz, \bz^{(t)}) \geq F(\bz)$; that is, $\bD-\bW^\top\bW$ is positive semidefinite.
To prove this, consider the matrix $\bM\in\real^{K\times K}$ whose entries are $m_{ij}=z_i (\bD-\bW^\top\bW)_{ij}z_j$ for all $i,j\in\{1,2,\ldots,K\}$, which is a rescaling of the components of $\bD-\bW^\top\bW$. Then $\bD-\bW^\top\bW$ is positive semidefinite if and only if $\bM$ is: 
$$
\begin{aligned}
\bx^\top& \bM\bx
= \sum_{i,j=1}^{K,K} x_i m_{ij}x_j
\stackrel{*}{=}\sum_{i,j=1}^{K,K} \left\{(\bW^\top\bW)_{ij} z_i z_j  x_i^2 - (\bW^\top\bW)_{ij}z_i z_j  x_i x_j\right\}\\
&\stackrel{\dag}{=}\sum_{i,j=1}^{K,K} (\bW^\top\bW)_{ij}z_i z_j \left( \frac{1}{2}x_i^2 + \frac{1}{2}x_j^2 - x_ix_j\right)
=\sum_{i,j=1}^{K,K} (\bW^\top\bW)_{ij}z_i z_j\frac{1}{2} \left( x_i-x_j\right)^2 \geq 0,
\end{aligned}
$$
where equality $(\dag)$ follows from the symmetry of $\bM$, and equality ($*$) follows from the diagonality of $\bD$:
$$
\sum_{i,j=1}^{K,K} x_iz_i d_{ij}z_jx_j
=
\sum_{i=1}^{K} x_iz_i d_{ii}z_ix_i
=
\sum_{i=1}^{K} x_i^2 z_i^2 \frac{\sum_{j=1}^{K}(\bW^\top\bW)_{ij}z_j}{z_i}
=
\sum_{i,j=1}^{K,K} (\bW^\top\bW)_{ij} z_i z_j  x_i^2.
$$
This completes the proof.
\end{proof}

The proof of Theorem~\ref{theorem:conv_mu_fro} follows directly from the above lemmas.
Clearly, the approximations $\bW$ and $\bZ$ remain nonnegative during the updates.
It is generally better to update $\bW$ and $\bZ$ ``simultaneously” rather than ``sequentially," i.e., updating each matrix completely before the other. In this case, after updating a row of $\bZ$, we update the corresponding column of $\bW$.  
In the implementation, it is advisable to introduce a small positive quantity, say the square root of the machine precision, to the denominators in the approximations of $\bW$ and $\bZ$
at each iteration. 
And a trivial value like $\epsilon=10^{-9}$  suffices. The full procedure is shown in Algorithm~\ref{alg:nmf-multiplicative}.
In practice, the algorithm can also be accelerated by updating $\bW$ several times before updating $\bZ$, during which process we can reuse the result of $\bA\bZ^\top$ and $\bZ\bZ^\top$, and vice versa.

\index{Machine precision}
\begin{algorithm}[h] 
\caption{NMF via Multiplicative Updates}
\label{alg:nmf-multiplicative}
\begin{algorithmic}[1] 
\Require Matrix $\bA\in \real_+^{M\times N}$;
\State Initialize $\bW\in \real_{++}^{M\times K}$, $\bZ\in \real_{++}^{K\times N}$ randomly with positive entries;
\State Choose a stoping criterion on the approximation error $\delta$;
\State Choose maximal number of iterations $C$;
\State $iter=0$; \Comment{Count for the number of iterations}
\While{$\normf{\bA- (\bW\bZ)}^2>\delta $ and $iter<C$}
\State $iter=iter+1$; 
\State $\bZ \leftarrow \bZ\hadaprod \frac{[\bW^\top\bA]}{[\bW^\top\bW\bZ]+\epsilon}$;
\State $\bW \leftarrow \bW \hadaprod \frac{[\bA\bZ^\top]}{[\bW\bZ\bZ^\top]+\epsilon}$;
\EndWhile
\State Output $\bW,\bZ$;
\end{algorithmic} 
\end{algorithm}


\index{Regularization}
\subsection{Regularization}
As mentioned in \eqref{equation:kkv_nnn_raw},  the NNLS or NMF problem implicitly imposes a \textbf{sparsity constraint}. 
On the other hand, similar to the ALS method with regularization discussed in Section~\ref{section:regularization-extention-general}, recall that the regularization can help extend the applicability of ALS to general matrices. 
Additionally, a regularization term can be incorporated into the NMF framework to enhance its performance:
$$
L(\bW,\bZ)  =\frac{1}{2}\normf{\bW\bZ-\bA}^2 +\frac{1}{2}\lambda_w \normf{\bW}^2 + \frac{1}{2}\lambda_z \normf{\bZ}^2, \qquad \lambda_w>0, \lambda_z>0,
$$
where the employed matrix norm is still the Frobenius norm. The gradient with respect to $\bZ$ given $\bW$ is the same as that in Equation~\eqref{equation:als-regulari-gradien}:
$$
\begin{aligned}
	\frac{\partial L(\bZ|\bW)}{\partial \bZ} =\bW^\top(\bW\bZ-\bA) + \textcolor{mylightbluetext}{\lambda_z\bZ}  \in \real^{K\times N}.
\end{aligned}
$$
The gradient descent update can be obtained by 
$$
(\text{GD on }\bZ) \gap \bZ \leftarrow \bZ - \eta \left(\frac{\partial L(\bZ|\bW)}{\partial \bZ}\right)=\bZ - \eta \left(  \bW^\top\bW\bZ-\bW^\top\bA+\textcolor{mylightbluetext}{\lambda_z\bZ}\right),
$$
Analogously, if we assume a different step size for each entry of $\bZ$, the update can be obtained by
$$
(\text{GD$^\prime$ on $\bZ$})\gap 
\begin{aligned}
z_{kn} 
&=z_{kn} - \eta_{kn}(\bW^\top\bW\bZ-\bW^\top\bA+\textcolor{mylightbluetext}{\lambda_z\bZ})_{kn}, 
 \,\, \forall k,n.
\end{aligned}
$$
We again rescale the step size:
$
\eta_{kn} = \frac{z_{kn}}{(\bW^\top\bW\bZ)_{kn}}.
$
Then we obtain the MU rules for $\bZ$ and $\bW$ (due to symmetry):
$$
\begin{aligned}
\bZ \leftarrow  \bZ\hadaprod \frac{[\bW^\top\bA-\lambda_z\bZ]}{[\bW^\top\bW\bZ]}
\gap\text{and}\gap 
\bW \leftarrow \bW \hadaprod \frac{[\bA\bZ^\top-\lambda_w\bW]}{[\bW\bZ\bZ^\top]}.
\end{aligned}
$$

\paragraph{Modified MU.}
Since the update for regularized NMF can result in negative values, a modified MU can be applied such that
$$
\textbf{(MMU1):}\quad  \bZ \leftarrow \left[\bZ\hadaprod \frac{[\bW^\top\bA-\lambda_z\bZ]}{[\bW^\top\bW\bZ]}\right]_+;
\gap\text{and}\gap 
\bW \leftarrow \left[\bW \hadaprod \frac{[\bA\bZ^\top-\lambda_w\bW]}{[\bW\bZ\bZ^\top]}\right]_+,
$$
where $[x]_+ = \max\{x, \epsilon\}$. The parameter $\epsilon$ is usually a very small positive number that prevents the emergence of negative update.
That is, we add a small lower bound for entries of $\bW$ and $\bZ$. Alternatively, another approach is: 
$$
\begin{aligned}
\textbf{(MMU2):}\quad 
\bZ \leftarrow \bZ\hadaprod \frac{[\bW^\top\bA-\lambda_z\bZ]_+}{[\bW^\top\bW\bZ]}
\gap\text{and}\gap 
\bW \leftarrow\bW \hadaprod \frac{[\bA\bZ^\top-\lambda_w\bW]_+}{[\bW\bZ\bZ^\top]}.
\end{aligned}
$$

\section{NMF with Three Factors}
The NMF method, when extended to incorporate three factor matrices, is referred to as  \textit{nonnegative matrix trifactorization (tri-NMF)}. 
This approach introduces an additional factor:
\begin{equation}\label{equation:tri_nmv}
\bA \approx \bW\bU\bZ,
\end{equation}
where $\bW\in\real_+^{M\times K}$, $\bU\in\real_+^{K\times J}$, and $\bZ\in\real_+^{J\times N}$.
Consider  the matrix-by-user matrix $\bA\in\real^{M\times N}$, where each element is a binary number $\{0,1\}$ (called \textit{implicit data}; in contrast to  \textit{explicit data} used previously).~\footnote{For example, in datasets like Netflix or MovieLens, ratings above 4 can be mapped to 1, while ratings below 1 can be mapped to 0 to obtain an implicit data set.}
NMF on this matrix provides a sum of  $K$ rank-one matrices $\bA\approx\sum_{k=1}^{K} \bW[:,k]\bZ[k,:]$.
In the context of implicit data, each rank-one matrix can be interpreted as finding a subset of users and a subset of items (e.g., movies)  that interact strongly with each other.
The tri-NMF, on the other hand, can be approximated by 
$$
\bA\approx\sum_{k=1}^{K}\sum_{j=1}^{J} \bW[:,k]\bU[k,j]\bZ[j,:].
$$
This formulation can be interpreted as identifying  separately $J$ subsets of movies that are watched together (the rows of $\bZ$) and $K$ subset of users that behave similarly (the columns of $\bW$); while the matrix $\bU$ tells us how these subsets interact together.
If $u_{kj}>0$, then the $k$-th subset of users (corresponding to the positive entries of $\bW[:,k]$) watches the movies from the $j$-th subset of movies (corresponding to the positive entries of $\bZ[j,:]$).

In other words, tri-NMF identifies groups of users who exhibit similar behavior (by watching the same movies) and groups of movies that are similar (because they are watched by the same users), while connecting these groups through the nonnegative matrix $\bU$. Another application of tri-NMF is in text mining, where it identifies subsets of documents that contain the same words (the columns of $\bW$) and subsets of words that appear in the same documents (the rows of $\bZ$), linking these subsets via the matrix $\bU$ \citep{brouwer2017comparative, gillis2020nonnegative, lu2023bayesian}.

\section{$\beta$-Divergence, Alternative Perspectives of MU}\label{section:beta_div_altmu}

We have introduced several alternative error measures for matrix factorization problems in Section~\ref{section:more_err_sta_als}. 
As mentioned previously, the sum of squared loss, as given in  \eqref{equation:als-per-example-loss2} or \eqref{equation:frob_nmf}, is convex when one of the factors is held constant, leading to a smooth optimization process.
This type of loss function falls under a broader class of distance/divergence estimators known as \textit{$\beta$-divergence} in the context of NMF.
Given two nonnegative scalars $x$ and $y$, the $\beta$-divergence between $x$ and $y$ is defined as follows:
\begin{equation}\label{equation:scalar_beta_div}
d_{\beta}(x, y)
=
\left\{
\begin{aligned}
&\frac{x}{y}-\ln \frac{x}{y}-1, &\text{if } \beta=0;\\
&x\ln\frac{x}{y} -x+y , &\text{if } \beta=1;\\
&\frac{1}{\beta^2-\beta}(x^\beta+(\beta-1)y^\beta - \beta xy^{\beta-1} ) , &\text{otherwise}.\\
\end{aligned}
\right.
\end{equation}
The $\beta$-divergence is continuous in $\beta$ since $\mathoplim{\beta\rightarrow 0}(x^\beta -y^\beta)/\beta=\ln(x/y)$. 
When $\beta=0, 1,$ and $2$, the $\beta$-divergences are also known as the \textit{Itakura-Saito (IS), KL, and Frobenius/Euclidean distances/divergences}, respectively.
The $\beta$-divergence between two matrices $\bB$ and $\bC$ is 
\begin{equation}\label{equation:beta_div_mat_def}
D_{\beta}(\bB,\bC)=
\sum_{j}d_{\beta}(\bb_j, \bc_j)
=
\sum_{i,j} d_{\beta}({b_{ij}, c_{ij}}).
\end{equation}

\begin{figure}[h]
\centering  
\vspace{-0.15cm}    
\subfigtopskip=2pt  
\subfigbottomskip=2pt 
\subfigcapskip=-5pt  
\subfigure[$\beta$-divergence for $d_{\beta}(1,y)$.]{\label{fig:beta_divergence_1}
\includegraphics[width=0.47\linewidth]{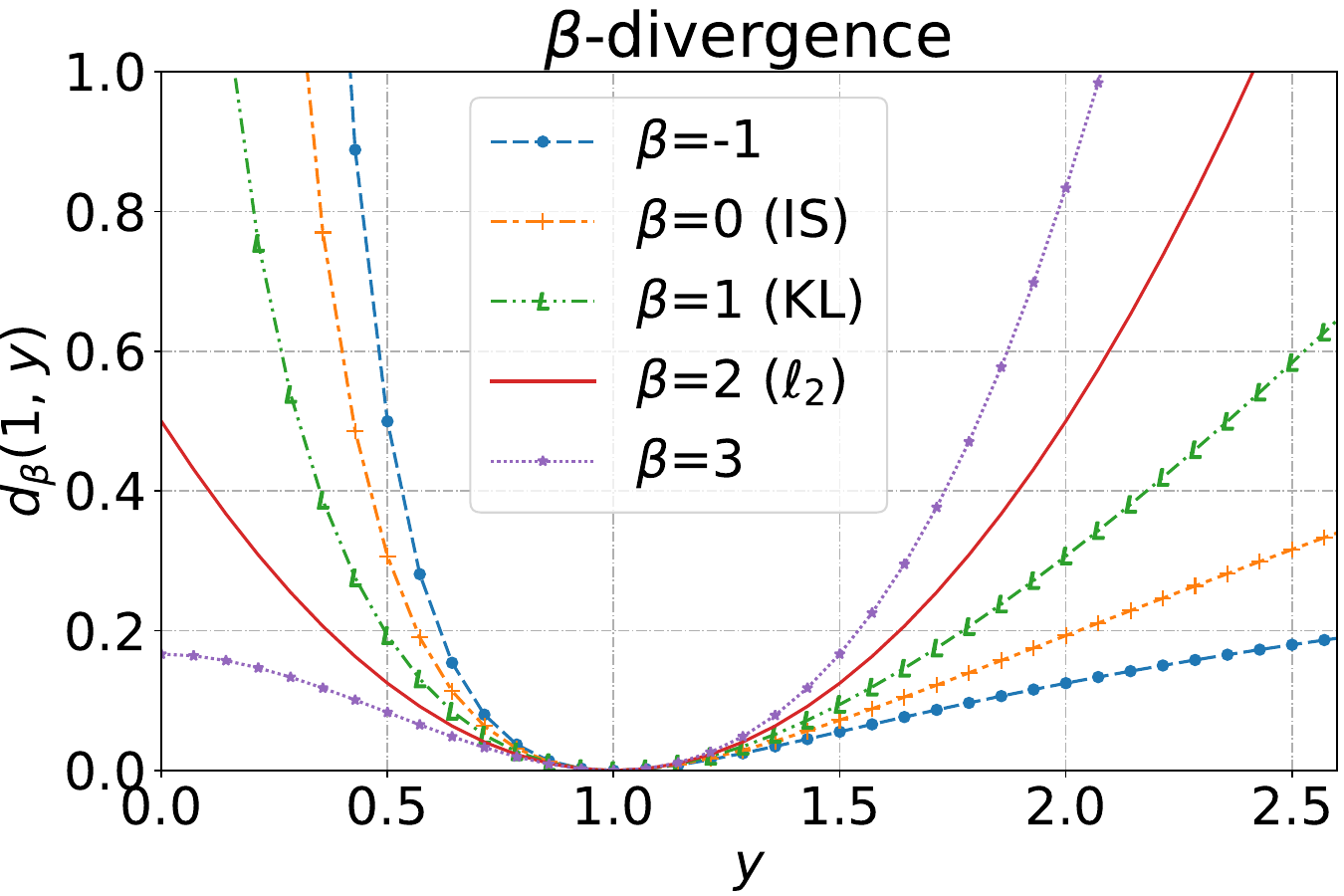}}
\subfigure[$\beta$-divergence for $d_{\beta}(2,y)$.]{\label{fig:beta_divergence_2}
\includegraphics[width=0.47\linewidth]{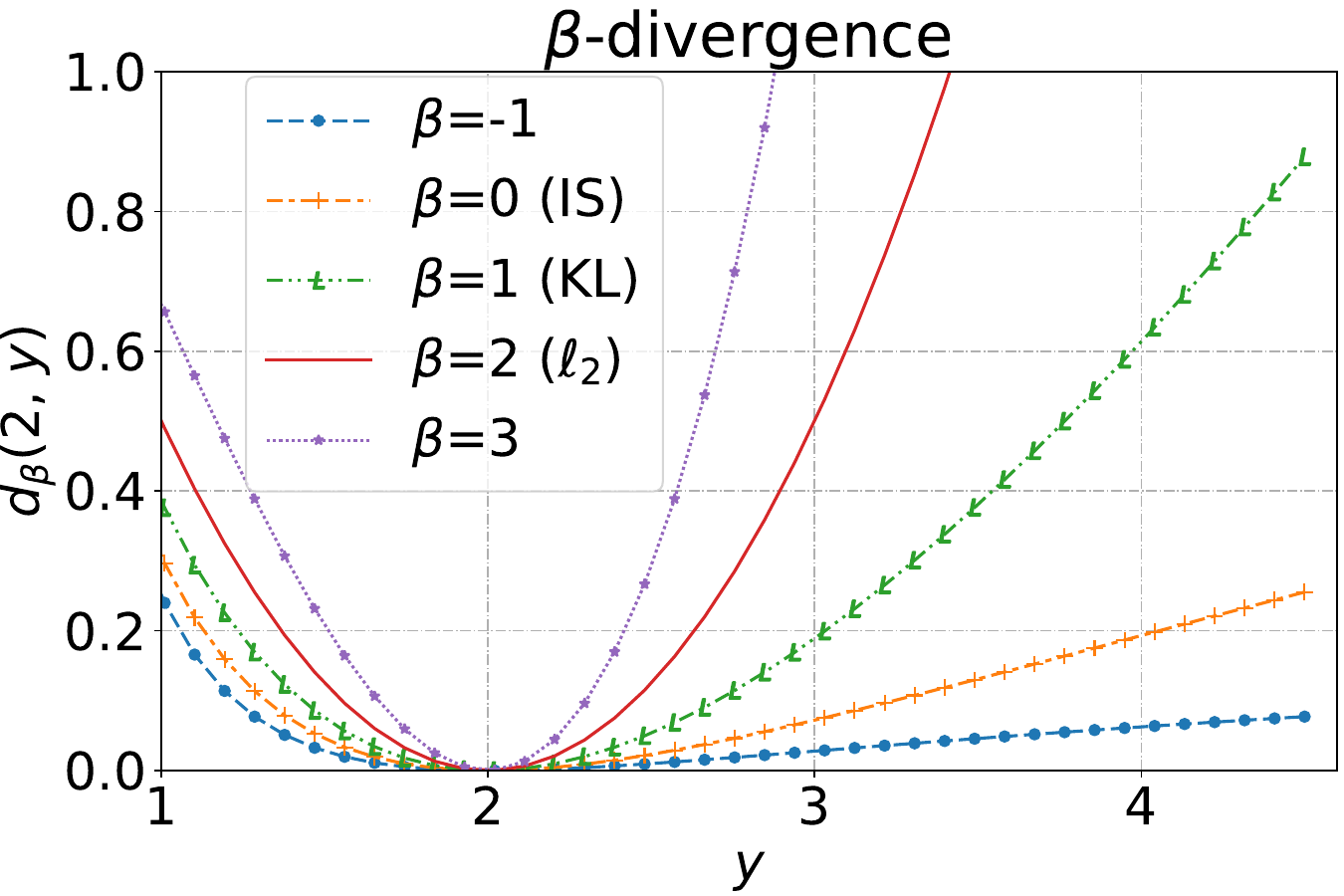}}
\caption{
The analysis of $\beta$-divergence is complex. When the first argument is fixed at 1,  smaller values are less penalized as  the $\beta$ value increases; however, when the first argument is $2$,  smaller values  are more penalized as  the $\beta$ value increases.
In both cases, larger values are more heavily penalized as the $\beta$ value increases.
}
\label{fig:beta_divergence_all}
\end{figure}

The analysis of $\beta$-divergence is complex. When the first argument is fixed at 1,  smaller values are less penalized as  the $\beta$ value increases; however, when the first argument is $2$,  smaller values  are more penalized as  the $\beta$ value increases.
In both cases, larger values are more heavily penalized as the $\beta$ value increases. See Figure~\ref{fig:beta_divergence_all}.

\paragraph{Convexity of $\beta$-divergence.}
When $\beta\in[1,2]$, the function $d_{\beta}(x,y)$ is convex in the second argument $y$. This implies $D_{\beta}(\bA,\bW\bZ)$ is convex in $\bW$ when fixing $\bZ$, and vice versa (Problem~\ref{prob:sep_conv}).  
Therefore, coordinate descent algorithms can be effectively applied for NMF using the $\beta$-divergence.

\paragraph{Scaling in arguments.} Let $\gamma>0$ be a scale factor, then 
\begin{equation}
d_{\beta}(\gamma x, \gamma y) = \gamma^\beta d_{\beta}(x,y).
\end{equation}
This indicates that the larger the $\beta$, the more sensitive the $\beta$-divergence is to large values of $x$ or $y$; on the contrary, $\beta$-divergence with small $\beta<0$ values relies more heavily on the smallest data values.
However, when $\beta=0$ (called \textit{Itakura-Saito divergence, IS divergence}, see \eqref{equation:als_gama_los}), the $\beta$-divergence is not related to the $\beta$ value, and thus it is invariant to scaling. 
What matters is the ratio between $x$ and $y$; see \eqref{equation:scalar_beta_div}.

\paragraph{Gradient.} 
Since we consider a nonnegative matrix $\bA$ for NMF problems, we should note that $d_{\beta}(x, \cdot)$ for $x=0$ is not defined for all values of $\beta$:
$$
d_{\beta}(0, y)=
\left\{
\begin{aligned}
\text{not defined}, \gap &\text{if $\beta\leq 0$};\\
\frac{1}{\beta} y^\beta, \gap \gap &\text{if $\beta> 0$},
\end{aligned}
\right.
\gapthree \implies \gapthree
d'_{\beta}(0, y)=
\left\{
\begin{aligned}
\text{not defined}, \gap &\text{if $\beta\leq 0$};\\
 y^{\beta-1}, \gap \gap &\text{if $\beta> 0$},
\end{aligned}
\right.
$$
where the derivative $d'_{\beta}(0, y)$ corresponds to the second variable $y$.
Therefore, the algorithm developed in the following sections can only be applied to positive matrices when $\beta\leq 0$.
Table~\ref{table:dom_beta_1} and Table~\ref{table:dom_beta_2} present the domains of $d_\beta(x, \cdot)$ and $d'_{\beta}(x, \cdot)$, respectively, for different values of $\beta$ and $x$.
\noindent
\begin{table}[h]
\centering
\begin{minipage}{0.37\textwidth}
\centering
\caption{Domain of $d_\beta(x, \cdot)$.}
\label{table:dom_beta_1}
\begin{tabular}{c|c|c|c}
\hline
& $\beta \leq 0$ & $\beta \in (0,1]$ & $\beta > 1$ \\
\hline
$x = 0$ & $\varnothing$ & $\real_+$ & $\real_+$ \\
$x > 0$ & $\real_{++}$ & $\real_{++}$ & $\real_+$ \\
\hline
\end{tabular}
\end{minipage}
\quad
\begin{minipage}{0.54\textwidth}
\centering
\caption{Domain of $d'_{\beta}(x, \cdot)$.}
\label{table:dom_beta_2}
\begin{tabular}{c|c|c|c|c}
\hline
& $\beta \leq 0$ & $\beta \in (0,1)$ & $\beta \in [1,2)$ & $\beta \geq 2$ \\
\hline
$x = 0$ & $\varnothing$ & $\real_{++}$ & $\real_+$ & $\real_+$ \\
$x > 0$ & $\real_{++}$ & $\real_{++}$ & $\real_{++}$ & $\real_+$ \\
\hline
\end{tabular}
\end{minipage}
\end{table}

On the other hand, the gradients of  $D_{\beta}(\bA,\bW\bZ)$ w.r.t. $\bZ$ and $\bW$ are given by (if exist) 
\begin{equation}
\begin{aligned}
\nabla_{\bZ} D_{\beta}(\bA,\bW\bZ) &= \bW^\top \big( (\bW\bZ)^{\beta-2} \hadaprod (\bW\bZ-\bA) \big);\\
\nabla_{\bW} D_{\beta}(\bA,\bW\bZ) &=  \big( (\bW\bZ)^{\beta-2} \hadaprod (\bW\bZ-\bA) \big)\bZ^\top, 
\end{aligned}
\end{equation}
where $(\bW\bZ)^{\beta-2}$ denotes the componentwise exponent.
When $\beta=2$, the gradient reduces to the one in \eqref{equation:givenw-update-z-allgd}.

\paragraph{Decomposition of $\beta$-divergence.}
The $\beta$-divergence can be divided into three parts: convex, concave, and constant terms. 
We should note that this decomposition is not unique since any affine term is both convex and concave. 
We follow the convention in  \citet{fevotte2011algorithms}:
\begin{equation}\label{equation:beta_decom_defi}
d_{\beta}(x,y) = \convd_{\beta}(x,y)+\concd_{\beta}(x,y)+\cnstd_{\beta}(x,y),
\end{equation} 
where $\convd_{\beta}(x,y)$ is convex in $y$, $\concd_{\beta}(x,y)$ is concave in $y$, and $\cnstd_{\beta}(x,y)$ is constant in $y$; see Table~\ref{table:beta_decom} for different $\beta$ values.

\begin{table}[h]
\setlength{\tabcolsep}{6.pt}    
\centering
\begin{tabular}{c|c|c|c}
\hline
& $\convd_{\beta}(x,y)/\convd'_{\beta}(x,y)$, convex & $\concd_{\beta}(x,y)/\concd'_{\beta}(x,y)$, concave & $\cnstd_{\beta}(x,y)$, constant \\ \hline\hline
$\beta < 1, \beta \neq 0$ & $-\frac{1}{\beta-1}xy^{\beta-1}/ -xy^{\beta-2}$ & $\frac{1}{\beta}y^{\beta}/y^{\beta-1}$ & $\frac{1}{\beta(\beta-1)}x^{\beta}$ \\ \hline
$\beta = 0$ & $xy^{-1}/-xy^{-2}$ & $\ln y/y^{-1}$ & $x(\ln x - 1)$ \\ \hline
$1 \leq \beta \leq 2$ & $d_{\beta}(x,y)/d'_{\beta}(x,y)$ & 0/0 & 0 \\ \hline
$\beta > 2$ & $\frac{1}{\beta}y^{\beta}/y^{\beta-1}$ & $-\frac{1}{\beta-1}xy^{\beta-1}$ & $\frac{1}{\beta(\beta-1)}x^{\beta}$ \\ \hline
\end{tabular}
\caption{Scalar convex-concave-constant decomposition  of $d_{\beta}(x,y)$ with respect to the second variable $y$, and the corresponding derivatives with respect to the second variable $y$.}
\label{table:beta_decom}
\end{table}

\paragraph{KKT conditions for NMF with $\beta$-divergence.}

The KKT conditions (Appendix~\ref{appendix:KKT_proj}) indicate that (see derivation in \eqref{equation:kkv_nnn_raw}):
\begin{equation}\label{equation:nmf_beta_kkt1}
\begin{aligned}
\bZ\geq \bzero,\gap&\nabla_{\bZ} D_{\beta}(\bA, \bW\bZ) &\geq& \bzero, \gap \langle \bZ, \nabla_{\bZ} D_{\beta}(\bA, \bW\bZ)\rangle &=&\bzero_{K\times N}; \\
\bW\geq \bzero,\gap&\nabla_{\bW} D_{\beta}(\bA, \bW\bZ) &\geq& \bzero, \gap \langle \bW, \nabla_{\bW} D_{\beta}(\bA, \bW\bZ)\rangle &=&\bzero_{M\times K}.
\end{aligned}
\end{equation}
This also implies
\begin{equation}
\begin{aligned}
\min\{\bZ, \nabla_{\bZ} D_{\beta}(\bA, \bW\bZ)\} = \bzero_{K\times N}
\gap \text{and}\gap 
\min\{\bW, \nabla_{\bW} D_{\beta}(\bA, \bW\bZ)\} = \bzero_{M\times K},
\end{aligned}
\end{equation}
where the min operator $\min\{\cdot, \cdot \}$ is applied componentwise.

\subsection{MU for $\beta$-Divergence Obtained by Gradient Ratio Heuristic}\label{section:mu_gd_ratio}
We have shown that the MU update with Frobeius norm can be derived from rescaled gradient descent. 
For brevity, let $\nabla_{\bZ}\triangleq\nabla_{\bZ} D_{\beta}(\bA,\bW\bZ) = \nabla_{\bZ}^+ - \nabla_{\bZ}^-$, where 
\begin{equation}\label{equation:mu_grati_decom}
\nabla_{\bZ}^+ = \bW^\top \big( (\bW\bZ)^{\beta-1} \big)
\gap \text{and}\gap 
\nabla_{\bZ}^- = \bW^\top \big( (\bW\bZ)^{\beta-2} \hadaprod \bA \big).
\end{equation}
When $z_{kn}>0, \forall k,n$, the KKT conditions show that $(\nabla_{\bZ}^+)_{kn} = (\nabla_{\bZ}^-)_{kn}$. The rule from gradient descent (i.e., $\bZ^{(t+1)}=\bZ^{(t)} - \eta\nabla_{\bZ}$) indicates a small decrease (resp. increase) of $z_{kn}$ will lead to a  decrease of the loss function if $(\nabla_{\bZ})_{kn}>0$ (resp. $<0$).
Therefore, it is reasonable to update $z_{kn}$ using the componentwise ratio between $\nabla_{\bZ}^-$ and $\nabla_{\bZ}^+$:
\begin{equation}\label{equation:mu_grati_z}
\bZ\leftarrow  \bZ \hadaprod \frac{[\nabla_{\bZ}^-]}{[\nabla_{\bZ}^+]},
\end{equation}
where $\frac{[\cdot]}{[\cdot]}$ represents the componentwise ratio of two matrices.
When $\beta=2$, the MU reduces to the case in Theorem~\ref{theorem:conv_mu_fro}.
When $\beta=1$, the loss function becomes the KL divergence, and the update for $\bZ$ is 
$$
\textbf{($\beta=1$)}:\gap \bZ\leftarrow \bZ \hadaprod \frac{[\bW^\top \frac{[\bA]}{[\bW\bZ]}]}{[\bW^\top \bone_{M\times N}]}.
$$
It can be shown that when $\beta\in[1,2]$, the MU algorithms will monotonically decrease $D_{\beta}(\bA, \bW\bZ)$.

\subsection{MU for $\beta$-Divergence  Obtained by Rescaled PGD}
As discussed in Section~\ref{section:nmf_apgd},
the PGD approach involves projecting the gradient descent update onto the feasible set. Consider a standard GD update on $f(\bx)$: $\bx^{(t+1)}=\bx^{(t)}-\eta\nabla f(\bx^{(t)})$, where $\eta$ is a step size and  $-\nabla f(\bx^{(t)})$ is a \textit{descent direction} ($\bg$ is a descent direction if $\bg^\top\nabla f(\bx^{(t)})<0$.)
Consider further  a diagonal $\bD$ such that $-\eta \nabla f(\bx^{(t)}) \rightarrow -\bD \nabla f(\bx^{(t)})$ is also a descent direction (replacing the step size by a diagonal matrix)~\footnote{$\bD$ can be relaxed to any positive definite matrices.}.
In this case, if the feasible set of $\bx$ is nonnegative, then the PGD is useful: $\bx^{(t+1)}=\mathcalP(\bx^{(t)}-\bD\nabla f(\bx^{(t)}))$, where $\mathcalP(x)=\max\{x, 0\}$~\footnote{see, for example, \citet{beck2017first}.}.
If we further decompose the gradient into positive and negative parts: $\nabla f(\bx^{(t)}) = \nabla^+f(\bx^{(t)}) - \nabla^-f(\bx^{(t)})$ with $\nabla^+f(\bx^{(t)})>0 $ and $\nabla^-f(\bx^{(t)})>0$, taking $\bD=\diag\big( \frac{[\bx^{(t)}]}{[\nabla^+f(\bx^{(t)}))]}  \big)$, the rescaled PGD update becomes a MU rule:
\begin{equation}
\bx^{(t+1)}=\mathcalP\bigg(\bx^{(t)}-\diag\big( \frac{[\bx^{(t)}]}{[\nabla^+f(\bx^{(t)})]}  \big)\nabla f(\bx^{(t)})\bigg)
=
\mathcalP\bigg(\bx^{(t)} \hadaprod \frac{[\nabla^-f(\bx^{(t)})]}{[\nabla^+f(\bx^{(t)})]}\bigg).
\end{equation}
If we use the decomposition of gradient in \eqref{equation:mu_grati_decom}, the rescaled PGD becomes the MU update for NMF in \eqref{equation:mu_grati_z}.
If we further incorporate a step size $\eta$ in the rescaled PGD update, it becomes 
\begin{equation}
\bx^{(t+1)}
=
\mathcalP\bigg((1-\eta)\bx^{(t)} + \eta\bx^{(t)} \hadaprod \frac{[\nabla^-f(\bx^{(t)})]}{[\nabla^+f(\bx^{(t)})]}\bigg).
\end{equation}
Since $-\bD \nabla f(\bx^{(t)})$ is a descent direction, the step size $\eta \in(0,1)$ can ensure the update is monotonically nonincreasing.
Note that the projection operator can be omitted since all  updates are nonnegative.

\subsection{MU for $\beta$-Divergence Obtained by MM Framework}
The $\beta$-divergence between two matrices can be defined columnwise (Equation~\eqref{equation:beta_div_mat_def}), and the $\beta$-divergence can be divided into three parts (convex, concave, and constant, Equation~\eqref{equation:beta_decom_defi}). Thus, the loss function in NMF can be decomposed into (note the loss function can be further divided componentwise):
$$
D_{\beta}(\bA, \bW\bZ)
=
\sum_{n=1}^{N}d_{\beta}(\ba_n, \bW\bz_n)
=
\sum_{n=1}^{N}\left(\convd_{\beta}(\ba_n, \bW\bz_n) + \concd_{\beta}(\ba_n, \bW\bz_n) + \cnstd_{\beta}(\ba_n, \bW\bz_n)\right).
$$
For each column $n$, the MM framework involves finding auxiliary functions for the three components separately.
To see this, we need the following lemma:
\begin{lemma}[Auxiliary Function By Parts]
Let $F(\bx)=\sum_{i=1}^{n} F_i(\bx)$, and  let $G_i(\bx, \widetildebx)$ be an auxiliary function for $F_i(\bx)$ at $\widetildebx$ for all $i$. 
Then, $G(\bx, \widetildebx)=\sum_{i=1}^{n}G_i(\bx, \widetildebx)$ is an auxiliary function for $F(\bx)$ at $\widetildebx$.
\end{lemma}
This lemma indicates that if the auxiliary function is constructed separately for each component, it allows us to decouple the optimization.

\paragraph{Constant part.}
There is no need to find an auxiliary function for the constant term  $\cnstd_{\beta}(\ba_n, \bW\bz_n)$, since it does not influence the minimization of $d_{\beta}(\ba_n, \bW\bz_n)$ with respect to $\bz_n$.
\paragraph{Concave part.}
Any concave function can be upper-bounded using  linearization (the tangent plane, Theorem~\ref{theorem:conv_gradient_ineq}):
$$
\concd_{\beta}(x, y) 
\leq 
\concd_{\beta}(x, \widetildey)
+
(y-\widetildey) \concd_{\beta}'(x, \widetildey), 
$$
where $\concd_{\beta}'(x, \widetildey)$ denotes the gradient of $\concd(x, \widetildey)$ with respect to the second component $\widetildey$.
Therefore, for any $\widetildebz_n\in\real^{K}$, the auxiliary function for the concave component $\concd_{\beta}(\ba_n, \bW\bz_n)$ can be constructed by  
$$
\concG(\bz_n, \widetildebz_n) = \concd_{\beta}(\ba_n, \bW\widetildebz_n) + (\bW\bz_n-\bW\widetildebz_n) \hadaprod \concd_{\beta}' (\ba_n, \bW\widetildebz_n).
$$

\paragraph{Convex part.}
The auxiliary function for the convex part follows from the convexity inequality (Exercise~\ref{exercise:convexfuncs}).
Construct a matrix $\bP\in\real^{M\times K}$ as follows:
\begin{equation}
p_{mk} = \frac{w_{mk}{\widetildez}_{kn}}{\sum_{j}w_{mj}{\widetildez}_{jn}}
=
\frac{w_{mk}{\widetildez}_{kn}}{\bW[m,:]\widetildebz_n}
\gap \implies\gap
\bP\geq\bzero \text{ and } \bP\bone = \bone.
\end{equation}
That is, each row of $\bP$ belongs to the unit simplex in $\real^K$.
Therefore, we have
$$
\begin{aligned}
\convd_{\beta}(a_{mn}, \bW[m,:]\bz_n ) 
&=
\convd_{\beta}\big(a_{mn}, \sum_{k=1}^{K}w_{mk}z_{kn} \big) 
=
\convd_{\beta}\big(a_{mn}, \sum_{k=1}^{K} p_{mk}\frac{w_{mk}z_{kn}}{p_{mk}} \big) \\
&\leq 
\sum_{k=1}^{K} p_{mk}\convd_{\beta}\big(a_{mn}, \frac{w_{mk}z_{kn}}{p_{mk}} \big). 
\end{aligned}
$$
The decomposition finds the auxiliary function for $D_{\beta}(\bA,\bW\bZ)$ w.r.t. $\bZ$.
\begin{theorem}[Auxiliary Function for $D_{\beta}(\bA,\bW\bZ)$ w.r.t. $\bZ$]\label{theorem:aux_db_z}
Let $\widetildeba_n=\bW\widetildebz_n$ with $\widetildea_{mn}\triangleq\bW[m,:]\widetildebz_n$ for all $m, n$, where $\widetildebz_n$ is any vector in $\real^K$. Then, $G(\bZ,\widetildebZ)=\sum_{n=1}^{N} G_{n}(\bz_n, \widetildebz_n) = \sum_{n=1}^{N} \sum_{m=1}^{M} G_{mn}$ is an auxiliary function for $D_{\beta}(\bA,\bW\bZ)$ w.r.t. $\bZ$, where 
$$
\footnotesize
\begin{aligned}
G_{mn}
&=
\cnstd_{\beta}(a_{mn}, \widetildea_{mn})
+
\concd_{\beta}(a_{mn}, \widetildea_{mn}) + \sum_{k=1}^{K}w_{mk}(z_{kn}-\widetildez_{kn}) \concd_{\beta}'(a_{mn}, \widetildea_{mn})
+ \sum_{k=1}^{K} \frac{w_{mk}\widetildez_{kn}}{\widetildea_{mn}}\convd_{\beta}\big(a_{mn},  \frac{\widetildea_{mn}z_{kn}}{\widetildez_{kn}} \big).
\end{aligned}
$$
\end{theorem}

\begin{exercise}[Gradient and Hessian of Auxiliary Functions]\label{exercise:gra_hes_aux}
Consider the setting and notations in Theorem~\ref{theorem:aux_db_z}. 
Let $G_n(\bz_n, \widetildebz_n)= \sum_{k=1}^{K} G_k (\bz_{kn}, \widetildebz_n) + C(\bz_n)$ where $C(\bz_n)$ is a constant w.r.t. $\bz_n$. 
That is, 
$$
G_k (\bz_{kn}, \widetildebz_n) 
= 
\sum_{m=1}^{M}w_{mk} z_{kn}\concd_{\beta}'(a_{mn}, \widetildea_{mn})
+ 
\sum_{m=1}^{M} \frac{w_{mk}\widetildez_{kn}}{\widetildea_{mn}}\convd_{\beta}\big(a_{mn},  \frac{\widetildea_{mn}z_{kn}}{\widetildez_{kn}} \big).
$$
Show that gradient of the auxiliary function is 
$$
\nabla_{z_{kn}}G_n(\bz_n, \widetildebz_n) 
=
\sum_{m=1}^{M} 
w_{mk}
\bigg(
\concd_{\beta}'(a_{mn}, \widetildea_{mn})
+
\convd_{\beta}'\big(a_{mn},  \frac{\widetildea_{mn}z_{kn}}{\widetildez_{kn}} \big)
\bigg), 
$$
and the Hessian matrix is diagonal with entries
$$
\nabla^2_{z_{kn}}G_n(\bz_n, \widetildebz_n) 
=
\sum_{m=1}^{M} 
w_{mk}
\frac{\widetildea_{mn}}{\widetildez_{kn}}
\bigg(
\convd_{\beta}{''}\big(a_{mn},  \frac{\widetildea_{mn}z_{kn}}{\widetildez_{kn}} \big)
\bigg).
$$
Note in all cases, the first-order derivative or the second-order derivative corresponds to the second argument of $d_{\beta}(\cdot, \cdot)$.
\end{exercise}

Since $\convd_{\beta}(\cdot, \cdot)$ is convex in the second argument, the Hessian is positive definite. Thus, the auxiliary function is convex.
These constructions result in the following theorem by minimizing the  auxiliary function obtained in Theorem~\ref{theorem:aux_db_z}.
\begin{theorem}[Nonincreasing of MU for $\beta$-Divergence \citep{fevotte2011algorithms, gillis2020nonnegative}]\label{theorem:conv_mu_beta}
Let $\bA\in\real_+^{M\times N}$, $\bW\in\real_{++}^{M\times K}$, and $\bZ\in\real_{++}^{K\times N}$. 
The  loss $D_{\beta}(\bA,\bW\bZ)$ remains nonincreasing under the following multiplicative update rules:
$$
\footnotesize
\begin{aligned}
\bZ\leftarrow \bZ\hadaprod \left(\frac{\left[\bW^\top \left\{ (\bW\bZ)^{(\beta-2)} \hadaprod \bA \right\}\right]}{[\bW^\top (\bW\bZ)^{(\beta-1)}]}\right)^{m(\beta)},
\gapthree \text{and}\gapthree
\bW\leftarrow \bW\hadaprod \left(\frac{\left[\left\{ (\bW\bZ)^{(\beta-2)} \hadaprod \bA \right\} \bZ^\top \right]}{[ (\bW\bZ)^{(\beta-1)}\bZ^\top]}\right)^{m(\beta)},
\end{aligned}
$$
where 
$$
m(\beta)=\left\{
\begin{aligned}
&\frac{1}{2-\beta}, \gap &\textit{if }&\beta<1; \\
&1,                  &\textit{if }& 1\leq\beta\leq 2; \\
&\frac{1}{\beta-1}, &\textit{if }& \beta>1.
\end{aligned}
\right.
$$
When $\beta=2$, the result reduces to Theorem~\ref{theorem:conv_mu_fro}.
When $1\leq \beta\leq 2$, the MU obtained via the MM framework coincides with heuristic described in Section~\ref{section:mu_gd_ratio}.
\end{theorem}
The update in Theorem~\ref{theorem:conv_mu_beta} ensures nonnegativity of the parameter updates, provided they are initialized with positive values.

\paragraph{Choice of $\beta$ for NMF.}
The choice of $\beta$-divergence for NMF is problem-dependent.
\citet{fevotte2009nonnegative} present  results of decomposing a piano power spectrogram using $\beta=0$ and demonstrate that components corresponding to very low residual noise and hammer strikes on the strings are extracted with great accuracy; these components are either ignored or severely degraded when using Euclidean or KL distances/divergences. 
\citet{fitzgerald2009use} show that  $\beta=0.5$ is optimal for  music source separation problems.

\paragraph{Convergence.}
An algorithm is said to be  \textit{convergent} if it produces a sequence of iterates $\{\bZ^{(t)}\}_{t\geq 1}$ or $\{\bW^{(t)}\}_{t\geq 1}$ that converges to a limit point $\bW^*$ or $\bZ^*$ satisfying the KKT conditions in \eqref{equation:nmf_beta_kkt1}. Monotonic nonincreasingness does not imply convergence in general, and neither is monotonicity necessary for convergence. Proving convergence of the MU methods is beyond the scope of this book; we refer the readers to \citet{gillis2020nonnegative, fevotte2011algorithms} and references therein for more details.

\section{Initialization}
A significant challenge in NMF is the absence of guaranteed convergence to a global minimum.  
Often, the convergence process is slow, and the algorithm may reach a suboptimal approximation.
In the preceding discussion, we initialized $\bW$ and $\bZ$ randomly. 
To mitigate this issue, there are also alternative strategies designed to obtain better initial estimates in the hope of converging more rapidly to a good solution \citep{boutsidis2008svd, gillis2014and}. We sketch the methods as follows for reference:
\begin{itemize}
\item \textit{Clustering techniques.} Apply some clustering methods to the columns of $\bA$, set the cluster means of the top $K$ clusters as the columns of $\bW$, and initialize $\bZ$ as a proper scaling of
the cluster indicator matrix (that is, $z_{kn}\neq 0$ indicates  that $\ba_n$ belongs to the $k$-th cluster);
\item \textit{Subset selection.} Pick $K$ columns of $\bA$, and set those as the initial columns for $\bW$. And analogously, $K$ rows of $\bA$ are selected to form the rows of $\bZ$;
\item \textit{SVD-based approach.} Suppose the optimal rank-$K$ approximation of $\bA$ is $\bA=\sum_{i=1}^{K}\sigma_i\bu_i\bv_i^\top$, where each factor $\sigma_i\bu_i\bv_i^\top$ is a rank-one matrix with possible negative values in $\bu_i$ and $\bv_i$, and nonnegative $\sigma_i$. Denote $[x]_+\triangleq\max(x, 0)$, we notice 
$$
\bu_i\bv_i^\top = [\bu_i]_+[\bv_i]_+^\top+[-\bu_i]_+[-\bv_i]_+^\top-[-\bu_i]_+[\bv_i]_+^\top-[\bu_i]_+[-\bv_i]_+^\top,
$$
where the first two rank-one factors in this decomposition are nonnegative.
Then, either $[\bu_i]_+[\bv_i]_+^\top$ or $[-\bu_i]_+[-\bv_i]_+^\top$ can be selected to replace the factor $\bu_i\bv_i^\top$. \citet{boutsidis2008svd} suggest to replace each rank-one factor in $\sum_{i=1}^{K}\sigma_i\bu_i\bv_i^\top$ with  either $[\bu_i]_+[\bv_i]_+^\top$ or $[-\bu_i]_+[-\bv_i]_+^\top$, selecting the one with the larger norm and scaling it properly.
In other words, if we select $[\bu_i]_+[\bv_i]_+^\top$, then $\sigma_i\cdot [\bu_i]_+  $ can be initialized as the $i$-th column of $\bW$, and $[\bv_i]_+^\top$ can be chosen as the $i$-th row of $\bZ$.
\end{itemize}
However, these techniques are not guaranteed to yield better performance theoretically. 
We recommend referring to the aforementioned papers for more detailed information.

\index{Implicit hierarchy}
\section{Movie Recommender Context}
Both  NMF and  ALS methods approximate a matrix and reconstruct its entries using a set of basis/template vectors. 
The key difference lies in the nature of these basis vectors and how the approximation is carried out.
The basis  in  NMF is composed of vectors with nonnegative elements while the basis vectors in  ALS can have positive or negative values.
In  NMF, each vector is reconstructed  as a nonnegative summation of the basis vectors with  ``relatively" small components in the direction of each basis vector.
In contrast, in the ALS approximation, the data is modeled as a linear combination of the basis vector such that we can add or subtract vectors as needed; and the components in the direction of each basis vector can be large positive values or negative values. Therefore, depending on the application, one or the other factorization can be utilized to describe the data with different meanings.

\paragraph{Movie recommender context.}
In the context of a movie recommender system,  the rows of $\bW$ represent the hidden features of  movies, while the columns of $\bZ$ represent the hidden features of users. 
For example, in NMF, a movie might be described as 0.5 comedy, 0.002 action, and 0.09 romantic. However, in the ALS approach, we can get combinations such as 4 comedy, $-0.05$ action, and $-3$ drama, indicating positive or negative contributions to each feature.

\paragraph{Implicit hierarchy.}
Both ALS and NMF do not rank the importance of each basis vector hierarchically. In contrast, singular value decomposition (SVD) ranks the importance of each basis vector based on the corresponding singular value.
In the SVD representation of $\bA=\sum_{i=1}^{r}\sigma_i\bu_i\bv_i^\top$, 
this usually means that the reconstruction $\sigma_1\bu_1\bv_1^\top$ via the first set of basis vectors dominates and is the most used set to reconstruct data, followed by the second set, and so on. 
This creates an implicit hierarchy in the SVD basis that doesn't happen in the ALS or the NMF approach. 
For instance, recall the low-rank approximation of a flag image in Section~\ref{section:als-low-flag}, where we find the second component $\bw_2\bz_2^\top$ via  ALS in Figure~\ref{fig:als52} plays an important role in the reconstruction of the original figure. 
Conversely, the second component $\sigma_2\bu_2\bv_2^\top$ via  SVD in Figure~\ref{fig:svd22} plays a less important role in the reconstruction.  

\paragraph{Interpretability of basis vectors.}
In  SVD, the basis vectors can be statistically interpreted as the directions of maximum variance, but many of these directions lack a clear visual or intuitive interpretation due to the presence of zero, positive, and negative entries. 
When these basis vectors are used in a linear combination, the combination involves complex cancellations of positive and negative values, which can obscure the intuitive physical meaning of individual basis vectors. As a result, many basis vectors do not provide a meaningful explanation for nonnegative data, such as  pixel values in a color image.
On one hand, the entries of a nonnegative pattern vector should all be nonnegative values. On the other hand, mutually orthogonal singular vectors must contain negative entries. For example, if all  entries of the singular vector  $\bu_1$ corresponding to the maximum singular value are nonnegative, then any other singular vector orthogonal to $\bu_1$ must contain at least one negative entry; otherwise, the orthogonality condition $\bu_1^\top\bu_j=0$ for $j\neq 1$ cannot be satisfied. This indicates that mutually orthogonal singular vectors are not suitable as pattern vectors or basis vectors in nonnegative data analysis.

\section{Other Applications}

\begin{figure}[h]
\centering       
\subfigtopskip=2pt               
\subfigbottomskip=-2pt         
\subfigcapskip=-30pt      
\includegraphics[width=0.98\textwidth]{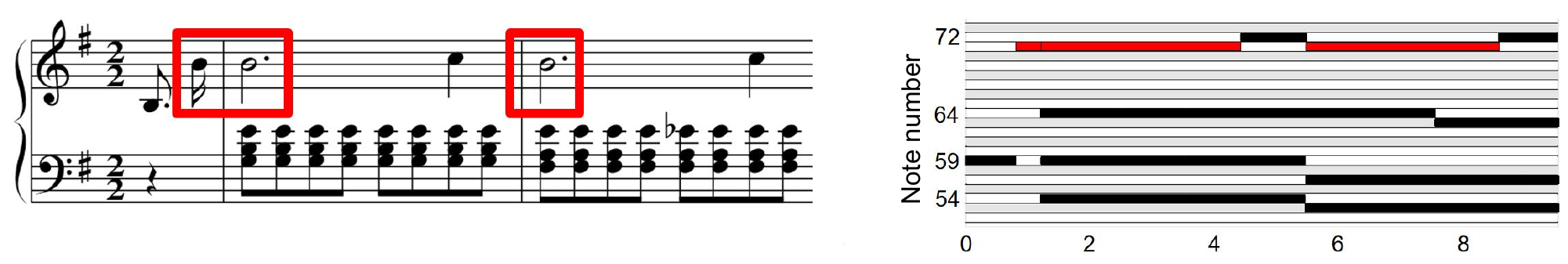}
\caption{Musical score and piano-roll representation. Figure is adapted from \citet{muller2015fundamentals}.}
\label{fig:nmf_music_note}
\end{figure}

\paragraph{Music spectral reconstruction.}
To illustrate the application of NMF, we demonstrate how this technique can decompose a spectrogram of a music recording into components that carry musical significance \citep{muller2015fundamentals}. As an example, let's examine the opening measures of \textit{Frédéric Chopin's Prélude Op. 28, No. 4}. Figure~\ref{fig:nmf_music_note} presents the musical notation alongside a piano-roll visualization that is synchronized with an audio recording of the piece. For clarity, all information pertaining to the note numbered $p=71$ are emphasized with red rectangular frames.

\begin{figure}[h]
	\centering       
	\vspace{-0.35cm}                 
	\subfigtopskip=2pt               
	\subfigbottomskip=-2pt         
	\subfigcapskip=-10pt      
	\includegraphics[width=0.98\textwidth]{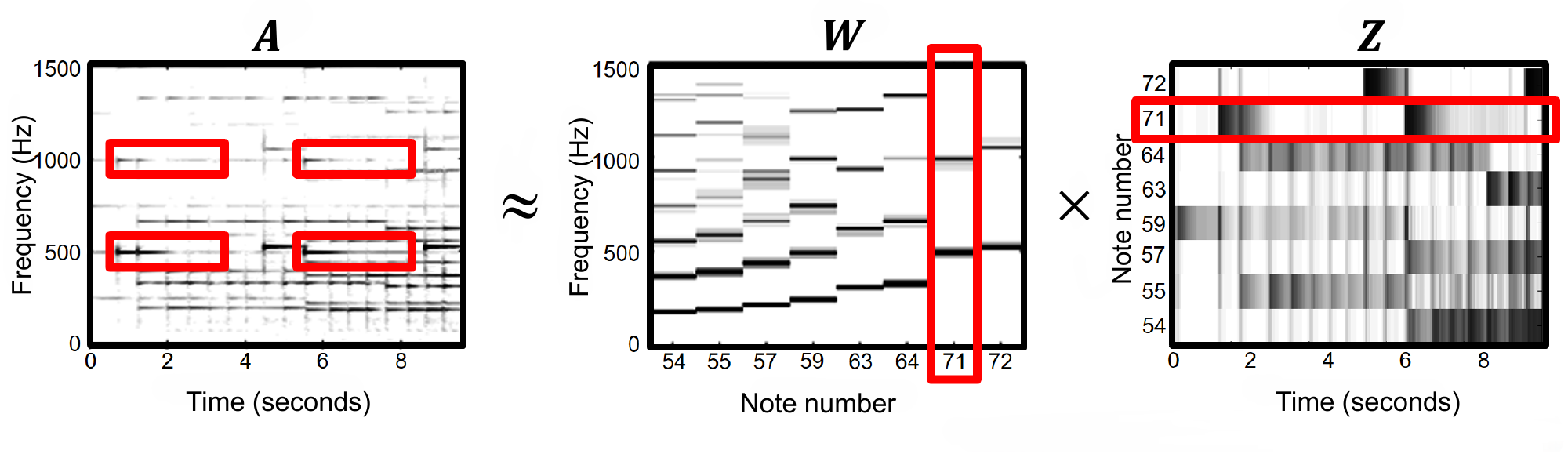}
	\caption{Ideal NMF of the spectrogram using NMF. Figure is adapted from \citet{muller2015fundamentals}.}
	\label{fig:nmf_music_decom}
\end{figure}

Regarding the original data matrix $\bA$, we utilize the magnitude STFT (see, for example, \citet{lopez2019nmf}), which consists of a series of spectral vectors. By applying NMF, this matrix can be decomposed into two nonnegative matrices, $\bW$ and $\bZ$. Ideally,  $\bW$ encapsulates the spectral patterns corresponding to the pitches of the notes present in the musical piece, whereas $\bZ$ indicates the temporal points at which these patterns appear in the audio recording. Figure~\ref{fig:nmf_music_decom} illustrates such a decomposition applied to the Chopin prelude.

In this scenario, each template represented by the matrix  $\bW$ corresponds to the spectral manifestation of a specific pitch within  $\bA$, and the activation matrix $\bZ$ resembles the piano-roll representation of the musical score. 
Therefore, the advantages of NMF over general matrix factorization are evident:
\begin{itemize}
\item \textbf{Nonnegativity constraint.} NMF enforces nonnegativity on both the factorization matrices $\bW$ and $\bZ$. This constraint makes the resulting matrices more interpretable because they can be directly related to physical or perceptual quantities in the domain of interest. In the case of music, the nonnegative factors correspond to meaningful musical elements like notes or chords.
\item \textbf{Interpretability.} In NMF, the matrix $\bW$ represents the spectral patterns (timbres) of the notes present in the music piece, and $\bZ$ indicates the temporal activations of these patterns. This leads to a more interpretable decomposition compared to unconstrained matrix factorization methods, where the factors might not have a clear physical or musical interpretation.
\end{itemize}

\begin{figure}[h]
\centering       
\vspace{-0.35cm}                 
\subfigtopskip=2pt               
\subfigbottomskip=-2pt         
\subfigcapskip=-10pt      
\includegraphics[width=0.98\textwidth]{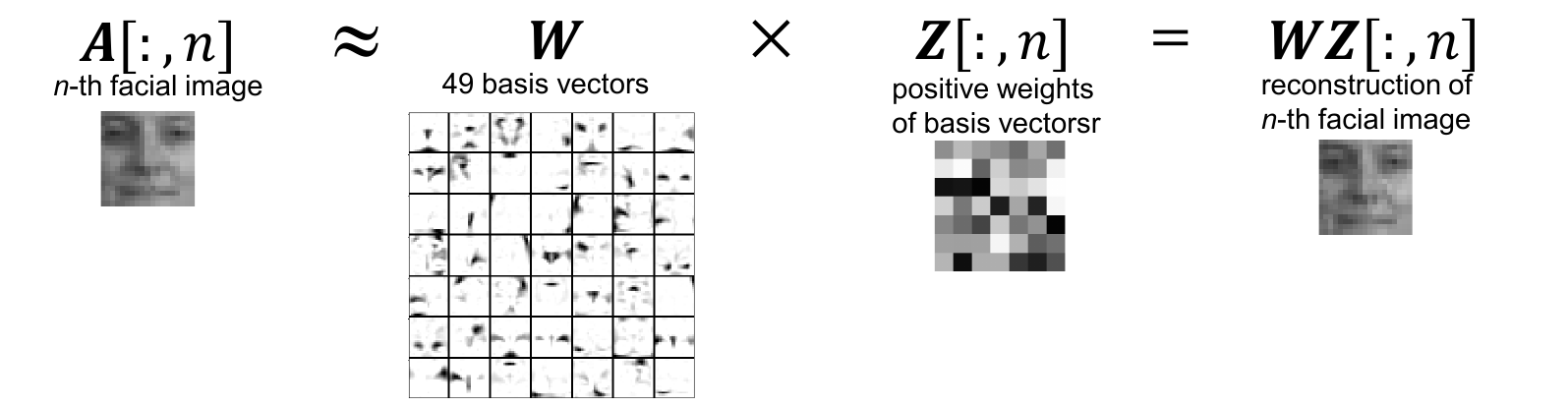}
\caption{NMF of the CBCL face database with $K=49$. The basis vectors in $\bW$ are reshaped into $19\times 19$ images. Facial features can be observed from these reshaped basis vectors, e.g., eyes, noses, nasolabial folds, and lips. Figure is adapted from \citet{lee1999learning, gillis2014and}.}
\label{fig:nmf_face_decom}
\end{figure}

\paragraph{Facial feature extraction and reconstruction.}
Suppose each column of the data matrix $\bA \in \real_+^{M\times N}$ represents a vectorized grayscale image of a face, where the $(m,n)$-th entry of $\bA$ corresponds to the intensity of the $m$-th pixel in the $n$-th face. NMF decomposes $\bA$ into two nonnegative matrices $\bW$ and $\bZ$ such that each image $\ba_n$ can be approximated by a linear combination of the columns of $\bW$. Since $\bW$ is nonnegative, its columns can be interpreted as images, referred to as template images, which are vectors of pixel intensities. The nonnegative weights in $\bZ$ ensure that these template images are combined additively to reconstruct each original face image. Given that the number of basis images $K$ is typically much smaller than the number of images $N$, the basis images should capture localized and sparse features that are common across multiple images. For facial images, these basis images often represent features such as eyes, noses, nasolabial folds, and lips (see Figure~\ref{fig:nmf_face_decom}, NMF for  the CBCL face data~\footnote{http://cbcl.mit.edu/software-datasets/FaceData2.html}), while the columns of $\bZ$ indicate the presence of these features in each image \citep{lee1999learning, gillis2014and}.

On the other hand, if each column of $\bA$ indicates a facial image of a single person, the  NMF approach can be utilized for face recognition. NMF has been shown to be more robust to occlusions compared to PCA or ALS, which generates dense factors. For instance, if a new face with occlusions (e.g., wearing sunglasses or distortions) needs to be mapped into the NMF basis, the non-occluded parts (e.g., the lips or the noise) can still be accurately approximated \citep{jain2017non}.

\paragraph{Topic recovery.} 
As mentioned at the very beginning of this chapter, NMF can be effectively applied to topic recovery problems. 
Typically, this application involves constructing a term-document matrix  $\bA$, where the rows correspond to terms (words or phrases) and the columns correspond to documents. Each entry  $a_{mn}$ in the matrix represents the weight or importance of term  $m$ in document $n$. 
This weight could be binary (presence/absence), \textit{term frequency (TF)}, or \textit{term frequency-inverse document frequency (TF-IDF)} \citep{shahnaz2006document}.
In this framework, each column of $\bW$ can be interpreted as a topic, while each column of $\bZ$ indicates the positive weight of each topic for the given document;  NMF in this context aligns with a soft clustering approach where each data point can belong to multiple clusters. 
NMF is particularly well-suited for topic recovery because it captures the additive nature of topics in text data and produces interpretable results. However, the choice of the number of topics $K$  and the initialization of $\bW$  and  $\bZ$ can significantly affect the quality of the results. Additionally, NMF is sensitive to the scaling of the input data, so it's important to preprocess the term-document matrix appropriately.

\index{$L$-strongly smoothness}
\begin{problemset}
	

\item \label{prob:llipschi_hianls} \textbf{$L$-strongly smooth and PGD in Hi-ANLS problems.} A function $f(\bx): \real^n\rightarrow \real$ is called an $L$-Lipschitz gradient continuous function (a.k.a., a $L$-strongly smooth function) if it satisfies that $\normtwo{\nabla f(\bx) - \nabla f(\by)} \leq L\normtwo{\bx-\by}$ for all $\bx,\by$. Show that the subproblem~\eqref{equation:llipschi_hianls} in Hi-ANLS is $L$-strongly smooth with constant $L=\normtwo{\bW[:,k]}^2$. Therefore, the subproblem can be treated as a \textit{projected gradient descent (PGD)} update with a step size $\eta=\frac{1}{L}$, i.e., using gradient descent update with a step size $\eta=\frac{1}{L}$ first and projecting the update onto the nonnegative orthant afterwards. 

\item \label{prob:lsmooth_quadra} \textbf{$L$-strongly smoothness of quadratic functions.} Let $f(\bx)=\frac{1}{2}\bx^\top\bA\bx + \bb^\top\bx+c$, where $\bA\in\real^{n\times n}, \bb\in\real^n$.  Show that $\norm{\nabla f(\bx) -\nabla f(\by)}_a = \norm{\bA\bx-\bA\by}_a  \leq \norm{\bA}_{b,a}\norm{\bx-\by}_{b}$ for any $\bx,\by\in\real^n$, where $\norm{\bA}_{b,a}$ is the induced norm (Definition~\ref{definition:induced_norm_app}) and $\frac{1}{a}+\frac{1}{b}=1$ with $a\in[1,\infty]$. 
This implies that affine functions $g(\bx)=\bb^\top\bx+c$ are $0$-strongly smooth.
How is this related to Problem~\ref{prob:llipschi_hianls}?
\item \label{prob:lsmooth_dslemma} \textbf{Descent lemma for $L$-strongly smooth functions.} Let $f:\sS\rightarrow (-\infty, \infty]$ be a function defined over a convex set $\sS$ such that $\normtwo{\nabla f(\bx) - \nabla f(\by)} \leq L\normtwo{\bx-\by}$ for all $\bx$ and $\by$.
Show that 
$
f(\by)\leq f(\bx)+ \nabla f(\bx)^\top (\by-\bx) + \frac{L}{2}\normtwo{\bx-\by}^2.
$
\textit{Hint: Use fundamental theorem of calculus (Appendix~\ref{appendix:taylor-expansion}). This is the  initial definition in Definition~\ref{definition:scss_func}.}

\item \label{prob:third_order_nmf} Let $\ba\in\real^M$, $\bz\in\real^K$, and $\bW\in\real^{K\times L}$. Show that the third-order partial derivatives of  $F(\bz)=\frac{1}{2}\normtwo{\ba-\bW\bz}^2$ vanish.

\item \textbf{MM applied to $L$-strongly smooth functions.} Let $f(\bx): \real^n\rightarrow \real$ be a $L$-strongly smooth function such that $\normtwo{\nabla f(\bx) - \nabla f(\by)} \leq L\normtwo{\bx-\by}$ for all $\bx,\by$. Show that $g(\bx, \widetildebx) = f(\widetildebx)+\nabla f(\widetildebx)^\top (\bx-\widetildebx) +\frac{L}{2} \normtwo{\bx-\widetildebx}^2$ is an auxiliary function for $f(\bx)$.
Find the update rule for this problem.

\item Derive the gradients and gradient descent updates for the tri-NMF problem in \eqref{equation:tri_nmv}.

\item \label{prob:projpro0} \textbf{Projection property-O.} Let $\sS\subset \real^n$ be \textbf{any set} and $\by\in\real^n$ such that $\widetildeby\triangleq\mathcalP_{\sS}(\by)$ is the projection of $\by$ onto set $\sS$. Show that for all $\bx\in\sS$, we have $\normtwo{\widetildeby - \by}\leq \normtwo{\bx-\by}$.

\item \label{prob:projpro1} \textbf{Projection property-I.} Let $\sS\subset \real^n$ be a \textbf{convex set} and $\by\in\real^n$ such that $\widetildeby\triangleq\mathcalP_{\sS}(\by)$. Show that for all $\bx\in\sS$, we have $\langle \bx-\widetildeby, \by-\widetildeby \rangle\leq 0$, i.e., the angle between the two vectors is greater than 90\textdegree.

\item \label{prob:projpro2} \textbf{Projection property-II.} Let $\sS\subset \real^n$ be a \textbf{convex set} and $\by\in\real^n$ such that $\widetildeby\triangleq\mathcalP_{\sS}(\by)$. Show that for all $\bx\in\sS$, we have $\normtwo{\widetildeby - \bx} \leq \normtwo{\by-\bx}$ and $\normtwo{\widetildeby - \bx}^2\leq \normtwo{\by-\bx}^2 - \normtwo{\by-\widetildeby}^2$ (the latter is related to the Pythagorean theorem). \textit{Hint: Examine $\normtwo{\by-\bx}^2=\normtwo{(\widetildeby-\bx)-(\widetildeby-\by)}^2$ and Problem~\ref{prob:projpro1}.}

\item \textbf{Linear feasibility.} Let $\sS=\{\bx\in\real^n: \bA\bx=\bb\}$ with full row rank $\bA$. Show that $\mathcalP_{\sS}(\bx) = \bx-\bA^\top(\bA\bA^\top)^{-1}(\bA\bx-\bb)$.

\item \label{prob:ortho_nmf} \textbf{Orthogonal and projective NMF, and clustering.} Consider the same setting as the  orthogonal or projective matrix factorization in Problem~\ref{prob:ortho_mf},  and consider further that $\bA,\bW$, and $\bZ$ are nonnegative. 
Show that there is only one positive entry in each column of $\bZ$ in this case. 
How  is this related to the K-means problem discussed in Section~\ref{section:regularization-extention-general}?
When each column of $\bA$ represents a data point, discuss the interpretation of $z_{kn}$ (the $(k,n)$-th entry of $\bZ$) as the importance of the $k$-th cluster to the $n$-th data point in the projective NMF case; that is, each data point can belong to several clusters.

\item Show that the Poisson loss in \eqref{equation:als_poi_los} is equivalent to minimizing the $\beta$-divergence between $\bA$ and $\bW\bZ$ with $\beta=1$.

\item Show that the Gamma loss in \eqref{equation:als_gama_los} is equivalent to minimizing the $\beta$-divergence between $\bA$ and $\bW\bZ$ with $\beta=0$.

\item \textbf{Regular and irreducible.} Let $\bA\in\real^{n\times n}$ be nonnegative and irreducible (Definition~\ref{definition:reduc_irreduc}). Show that $(\bI+\bA)^{n-1}> \bzero$, i.e., $\bI+\bA$ is regular (Definition~\ref{definition:pox_nonnegamat}). \textit{Hint: Find the relationship in the sequence $\bx^{(t+1)}=(\bI+\bA)\bx^{(t)}$ and show that $(\bI+\bA)\bx>\bzero$ for any nonzero $\bx\geq \bzero$.}

\item \label{prob:ab_diverg} \textbf{AB divergence \citep{amari2000methods}.} Let the \textit{$\alpha$-$\beta$ (AB) divergence} be given as follows:
\begin{align*}
d_{\alpha,\beta}(x,y)&=\begin{cases}
-\frac{1}{\alpha\beta}(x^\alpha y^\beta-\frac{\alpha}{\alpha+\beta}x^{\alpha+\beta}-\frac{\beta}{\alpha+\beta}y^{\alpha+\beta}),&\alpha,\beta,\alpha+\beta\neq 0;\\
\frac{1}{\alpha^2}(x^\alpha\ln(\frac{x^\alpha}{y^\alpha})-x^\alpha+y^\alpha),&\alpha\neq 0,\beta=0;\\
\frac{1}{\alpha^2}(\ln(\frac{y^\alpha}{x^\alpha})+(\frac{y^\alpha}{x^\alpha})^{-1}-1),&\alpha=-\beta\neq 0;\\
\frac{1}{\beta^2}(y^\beta\ln(\frac{y^\beta}{x^\beta})-y^\beta+x^\beta),&\alpha=0,\beta\neq 0;\\
\frac{1}{2}(\ln (x)-\ln (y))^2,&\alpha=0,\beta=0.
\end{cases}
\end{align*}
When $\alpha+\beta=1$, it is called the \textit{$\alpha$-divergence}. Discuss when  it reduces to the $\beta$-divergence.
Show that $d_{\alpha,\beta}(x,y)\geq 0$ and the equality holds if and only if $x=y$.

\item \label{prob:nonn_lin_1} Suppose $\bC=\abs{\bA}$ is defined as the matrix obtained by setting each entry of $\bC$ as the absolute value of $\bA\in\real^{n\times n}$. Show that $\abs{\bA}\geq \bzero_{m,n}$ and $\abs{\bA}=\bzero$ if and only if $\bA=\bzero$.

\item Suppose $\bA \geq \bzero_{n}$ is nonnegative and has a positive row. If $\abs{\bA\bx}=\bA\abs{\bx}$, where $\bx\in\complex^n$, then there exists a real $\theta\in[0, 2\pi)$ such that $e^{-i\theta}\bx = \abs{\bx}$, where $e^{-i\theta}\bx$ indicates a complex vector with $j$-th element being $e^{-i\theta} x_j$. \textit{Hint: Use triangle inequality $\abs{\bA\bx}\leq \abs{\bA}\abs{\bx}$, and examine the positive row. In the polar coordinate notation, $e^{i\theta}=\cos\theta+i\sin\theta$ and $\abs{e^{i\theta}x}=\abs{x} \implies \normtwo{e^{i\theta}\bx}=\normtwo{\bx}$.}

\item \label{prob:nnga_algebra} \textbf{Nonnegative algebra.} A bounty of results can be harvested from nonnegative conditions. We investigate several of them in this problem. Given square matrices $\bA,\bB,\bC,\bD\in\real^{n\times n}$, show that 
\begin{itemize}
\item \textbf{Triangle inequality.} $\abs{\bA\bB}\leq \abs{\bA}\abs{\bB}$.
\item \textbf{Nonexpansiveness.} $\abs{\bA^k}\leq \abs{\bA}^k$, for all $m=\{1,2,\ldots\}$.
\item \textbf{Equal norm.} $\normf{\bA}=\normf{\abs{\bA}}$.
\item If $\abs{\bB}\geq \abs{\bA}$, then $\normf{\bB}\geq \normf{\bA}$.
\item If $ \bB \geq \bA\geq \bzero$ and $ \bD \geq \bC\geq \bzero$, then $ \bB\bD \geq \bA\bC\geq \bzero$.
\item If $ \bB \geq \bA\geq \bzero$, then $\bB^k\geq \bA^k\geq \bzero$, for all $m=\{1,2,\ldots\}$.
\end{itemize}
Given rectangular matrices $\bA,\bB\in\real^{m\times n}$, show that 
\begin{itemize}
\item $\abs{\bA+\bB}\leq \abs{\bA}+\abs{\bB}$.
\end{itemize}

\item$^\ast$ \label{prob:nnga_algebra2} \textbf{Eigenvalue interlacing in nonnegative matrices.}  Let   $\bB - \abs{\bA}\in\real_+^{n\times n}$. Show  that 
$$
\rho(\bA) \leq \rho(\abs{\bA}) \leq \rho(\bB),
$$
where $\rho(\bX)$ represents the spectral radius of matrix $\bX$ (Definition~\ref{definition:spectrum}).
\textit{Hint: Use Problem~\ref{prob:nnga_algebra} and Gelfand formula (Exercise~\ref{exercise:gelfand_formula}); show that $\normf{\bA^k}\leq \normf{\abs{\bA}^k} \leq \normf{\bB^k}$.}

\item \label{prob:nnga_algebra3} Use Problem~\ref{prob:nnga_algebra2} to show that $\rho(\bB)\geq \rho(\bA)$ if $\bB\geq \bA\geq \bzero$.

\item Let $\bA\in\real_+^{n\times n}$, $\bB=\bA[1:k,1:k], \,\forall k\in\{1,2,\ldots,n\}$  (any leading principal submatrix of $\bA$, Definition~\ref{definition:leading-principle-minors}), and $\bC\in\real^{k\times k}, \,\forall k\in\{1,2,\ldots,n\}$ be any principal submatrix (Definition~\ref{definition:principle-minors}). Show that 
\begin{itemize}
\item $\rho(\scriptsize\begin{bmatrix}
\bB & \bzero \\
\bzero & \bzero 
\end{bmatrix}
\normalsize
) 
\leq 
\rho (\bA)
$
$\implies\rho(\bB)\leq \rho(\bA)$.
\item Use the first result to prove $\rho(\bC)\leq \rho(\bA)$. \textit{Hint: Use permutation transformations.}
\item $\mathopmax{i=1,2,\ldots,n}a_{ii} \leq \rho(\bA)$.
\end{itemize}

\item$^\ast$  \label{prob:nonn_lin_12} Let $\bA\in\real_+^{n\times n}$. Show that 
$$
\begin{aligned}
\text{Row sum: }\gap  \mathop{\min}_{1\leq i \leq n} \sum_{j=1}^{n} a_{ij} 
&\leq \rho(\bA)
\leq 
\mathop{\max}_{1\leq i \leq n} \sum_{j=1}^{n} a_{ij}; \\
\text{Column sum: }\gap \mathop{\min}_{1\leq j \leq n} \sum_{i=1}^{n} a_{ij} 
&\leq \rho(\bA)
\leq 
\mathop{\max}_{1\leq j \leq n} \sum_{i=1}^{n} a_{ij}. \\
\end{aligned}
$$

\end{problemset}

%% file: chapter-special.tex
\newpage 
\chapter{Decomposition with Special Matrix Products}
\begingroup
\hypersetup{
	linkcolor=structurecolor,
	linktoc=page,  
}
\minitoc \newpage
\endgroup

\section{Special Matrix Products} \index{Matrix products}

In the forthcoming sections, several key matrix products will play a significant role in elucidating the algorithms under discussion. These matrix operations not only serve as fundamental components in the development and understanding of these algorithms but also find extensive use in the context of tensor decomposition, which is explored in Part \ref{part:tensor_decom}. 
These special matrix products have been extensively studied in the field of tensor approximation and decomposition, particularly for enhancing the rank and descriptive power of low-rank tensors. With the rise of low-rank adaptations in large language models, these matrix products have also gained popularity in matrix analysis.
The importance of these matrix products cannot be overstated, as they form the basis for many computational techniques and methodologies that are central to both algorithmic illustration and tensor-based analyses.
\subsection{Hadamard Product/Power}
The Hadamard product, also known as the element-wise product, is defined for two matrices of identical dimensions. It involves multiplying each element from one matrix by the corresponding element in the other matrix, resulting in a new matrix of the same size. 
This operation is particularly useful in various mathematical and computational contexts (as we have already seen in the ALS Section~\ref{section:alt-columb-by-column}), where element-by-element multiplication can simplify complex operations or facilitate the implementation of certain algorithms.
\begin{definition}[Hadamard Product/Power]\label{definition:hada_prod}
The \textit{Hadamard product} (a.k.a., \textit{Schur product or element-wise product}) of matrices $\bA, \bB\in \real^{I\times J}$, and the Hadamard power (a.k.a., Schur power) of $\bA$, denoted by $\bA\hadaprod\bB$ and $\bA^{\hadaprod\gamma}$, respectively:
$$
\bA\hadaprod\bB = 
\begin{bmatrix}
a_{11}b_{11} & a_{12}b_{12} & \ldots & a_{1J}b_{1J}\\
a_{21}b_{21} & a_{22}b_{22} & \ldots & a_{2J}b_{2J}\\
\vdots  & \vdots & \ddots & \vdots\\
a_{I1}b_{I1} & a_{I2}b_{I2} & \ldots & a_{IJ}b_{IJ}\\
\end{bmatrix},
\gap 
\bA^{\hadaprod\gamma}=
\begin{bmatrix}
	a_{11}^{\gamma} & a_{12}^{\gamma} & \ldots & a_{1J}^{\gamma}\\
	a_{21}^{\gamma} & a_{22}^{\gamma}& \ldots & a_{2J}^{\gamma}\\
	\vdots  & \vdots & \ddots & \vdots\\
	a_{I1}^{\gamma} & a_{I2}^{\gamma} & \ldots & a_{IJ}^{\gamma}\\
\end{bmatrix},
$$
which are  matrices of size $I\times J$. That is, each entry in $\bA\hadaprod\bB$ is the product of the corresponding elements in $\bA$ and $\bB$, and each entry in $\bA^{\hadaprod\gamma}$ is the $\gamma$-th power of the corresponding element in $\bA$. Note that $\bA^{\hadaprod 2}=\bA\hadaprod \bA$ and $\bA^{\hadaprod 0}=\bone_{I\times J}$. Additionally, $\bI\hadaprod\bA$ denotes the diagonal part of $\bA$.~\footnote{In most textbooks, the Hadamard product is denoted by ``$\circ$"; not to being conflicted with the tensor outer product discussed in Part~\ref{part:tensor_decom}, we use ``$\hadaprod$" here.}
\end{definition}

\begin{exercise}[Schur's Product Theorem \citep{schur1911bemerkungen, horn1994topics}]\label{exercise:schur_product}
Let $\bA,\bB\in\real^{n\times n}$. Show that 
\begin{itemize}
\item $\bA\hadaprod \bB$ is PSD if $\bA$ and $\bB$ are PSD.
\item $\bA\hadaprod \bB$ is PD if $\bA$ and $\bB$ are PD.
\item $\bA\hadaprod \bB$ is PD if $\bA$ is PD and the diagonals of $\bB$ are all positive.
\end{itemize}
\end{exercise}
\begin{exercise}[Eigenvalue Inequality \citep{bapat1997nonnegative}]
Let $\bA,\bB\in\real^{n\times n}$ be PSD, and let $\lambda_1,\lambda_2,\ldots,\lambda_n$ be the eigenvalues of $\bA\hadaprod\bB$ and $\widetilde{\lambda}_1,\widetilde{\lambda}_2 \ldots,\widetilde{\lambda}_n$ be the eigenvalues of $\bA\bB$. Show that 
$\prod_{i=k}^{n}\lambda_i \geq \prod_{i=k}^{n}\widetilde{\lambda}_i $ for all $k\in\{1,2,\ldots,n\}$.
\end{exercise}

Further properties of the Hadamard product are discussed in Problems~\ref{prob:hadama_1}$\sim$\ref{prob:hadama_n}.
\subsection{Kronecker Product}
The \textit{Kronecker product} of two vectors $\ba \in \real^{I}$ and $\bb\in \real^{K}$, denoted by $\ba\kronecker \bb$, is defined as follows:
$$
\ba\kronecker \bb=
\begin{bmatrix}
a_1\bb \\
a_2\bb \\
\vdots \\
a_I\bb
\end{bmatrix}
\in\real^{IK},
$$
which is a column vector of size $(IK)$. It can be easily verified that if $\normtwo{\ba}=\normtwo{\bb}=1$, then $\normtwo{\ba\kronecker \bb}=1$.
In such a case, it follows that 
\begin{equation}
\ba\bb^\top = \ba\kronecker\bb^\top=\bb^\top\kronecker\ba \in \real^{I\times K}.
\end{equation}
Note that $\real^{IK}$ denotes a vector, while $\real^{I\times K}$ denotes a $I$-by-$K$ matrix.

The Kronecker product of two matrices allows for the construction of large matrices from smaller ones, which can be useful in creating block matrices or structured matrices. This property is particularly valuable in areas like system theory or control theory, where systems are often represented as large matrices.
Formally, we define the Kronecker product as follows.
\begin{definition}[Matrix Kronecker Product and Power]\label{definition:kronecker-product}
Similarly, the \textit{Kronecker product} of two matrices $\bA \in \real^{I\times J}$ and $\bB\in \real^{K\times L}$, denoted by $\bA\kronecker \bB$, is defined as follows:
$$
\begin{aligned}
\bA\kronecker \bB &= 
\begin{bmatrix}
a_{11} \bB & a_{12}\bB & \ldots & a_{1J}\bB \\
a_{21} \bB & a_{22}\bB & \ldots & a_{2J}\bB \\
\vdots  & \vdots  & \ddots & \vdots \\
a_{I1} \bB & a_{I2}\bB & \ldots & a_{IJ}\bB \\
\end{bmatrix}\\
&=
\begin{bmatrix}
\ba_1 \kronecker \bb_1 \,\,\ldots & \ba_1\kronecker \ba_L \mid 
\ba_2 \kronecker \bb_1   \,\,\ldots &\ba_2\kronecker \ba_L \mid 
\ba_J \kronecker \bb_1   \,\,\ldots &\ba_J\kronecker \ba_L\,
\end{bmatrix},
\end{aligned}
$$
which is a matrix of size $(IK)\times (JL)$. 
That is, the Kronecker product $\bA\kronecker \bB$ can be divided into $I\times J$ blocks; for each block $(i,j)$, it is a $K\times L$ matrix recorded by $a_{ij}\bB$.
When $\bA$ and $\bB$ are of the same shape, then it can be shown that $\bA\hadaprod \bB$ is a principal submatrix of $\bA\kronecker \bB$.
The \textit{Kronecker power} of a matrix $\bA$ is defined as $\bA^{\kronecker k}=\bA\kronecker\bA\kronecker \ldots\kronecker\bA$ with $\bA$ appearing $k$ times for $k=1,2,\ldots$.
\end{definition}
\begin{exercise}[Hadamard as Submatrix of Kronecker]\label{exercise:hada_sub_kro}
Let $\bA, \bB\in\real^{m\times n}$. Show that $\bA\hadaprod\bB=(\bA\kronecker\bB)[I,J]$, where $I=\{1,m+2, 2m+3, \ldots, m^2\}$ and $J=\{1,n+2,2n+3,\ldots,n^2\}$. That is, $\bA\hadaprod\bB$ is a  submatrix of $\bA\kronecker\bB$. When $m=n$, $\bA\hadaprod\bB$ is a principal submatrix of $\bA\kronecker\bB$.
\end{exercise}

The Kronecker product maintains specific matrix properties that can simplify computations and ensure stability.
\begin{lemma}[Kronecker of Orthogonal, Triangular, Diagonal, (Semi)definite, Nonsingular]\label{lemma:krokecker_keep_special}
The Kronecker products of two orthogonal, two triangular, two diagonal,  two (semi)definite,  two nonsingular, or two permutation matrices are also orthogonal, triangular,  diagonal, (semi)definite,  nonsingular, or permutation matrices, respectively. 
\end{lemma}

Specifically,  given four vectors \{$\ba \in \real^{I}$ and $\bb\in \real^{K}$\} and \{$\bc \in \real^{I}$ and $\bd\in \real^{K}\}$, then 
\begin{equation}\label{equation:kronecker-vector-find2}
(\ba\kronecker \bb)^\top(\bc\kronecker \bd)=
\begin{bmatrix}
a_1\bb^\top & 
a_2\bb^\top &
\ldots &
a_I\bb^\top  
\end{bmatrix}
\begin{bmatrix}
c_1\bd \\
c_2\bd \\
\vdots \\
c_I\bd
\end{bmatrix}
=\sum_{i=1}^{I}a_ic_i \bb^\top\bd=(\ba^\top\bc)(\bb^\top\bd).
\end{equation}
Particularly, when $\bc=\ba$ and $\bd=\bb$, it follows that 
\begin{equation}
(\ba\kronecker \bb)^\top(\ba\kronecker \bb) =\normtwo{\ba}^2 \cdot \normtwo{\bb}^2.
\end{equation}
Similarly, given four matrices $\bA,\bC\in \real^{I\times J}$ and $\bB,\bD\in \real^{K\times L}$, it follows that 
\begin{equation}\label{equation:kron_pro1_trans}
(\bA\kronecker \bB)^\top (\bC\kronecker \bD) = (\bA^\top\bC) \kronecker (\bB^\top\bD)\in\real^{JL\times JL}.
\end{equation}
Note also, for $\bA \in \real^{I\times J}$, $\bB\in \real^{K\times L}$, $\bC\in \real^{J\times I}$, and $\bD\in \real^{L\times K}$, the  equation above reduces to
\begin{equation}\label{equation:kron_pro1_notrans}
(\bA\kronecker \bB) (\bC\kronecker \bD) = (\bA\bC) \kronecker (\bB\bD) \in\real^{IK\times IK}.
\end{equation}
When $\bA \in \real^{I\times J}$, $\bB\in \real^{K\times L}$, $\bc\in \real^{J}$, and $\bd\in \real^{L}$, the equality becomes
\begin{equation}\label{equation:kronecker_eq3}
(\bA\kronecker \bB) (\bc \kronecker \bd) = (\bA\bc) \kronecker (\bB\bd)\in\real^{IK}.
\end{equation}
More generally, when $\bA \in \real^{I\times J}$, $\bB\in \real^{K\times L}$, $\bC\in \real^{J\times P}$, and $\bD\in \real^{L\times Q}$, we have (see Problem~\ref{problem:gene_kro_mat}) 
\begin{equation}
	(\bA\kronecker \bB) (\bC\kronecker \bD) = (\bA\bC) \kronecker (\bB\bD)\in\real^{IK\times PQ}.
\end{equation}

From the definition of the Kronecker product, it can be readily verified that any eigenvalue of $\bA\kronecker \bB$ arises as a product of the eigenvalues of $\bA$ and $\bB$.
\begin{lemma}[Eigenvalue, Determinant, and Trace of Kronecker Product]
Suppose $\bA\in\real^{m\times m}$ has an eigenpair $(\lambda, \bx)$, and $\bB\in\real^{n\times n}$ has an eigenpair $(\mu, \by)$. 
Then, $\lambda\mu$ is an eigenvalue of $\bA\kronecker \bB$ corresponding to the eigenvector $\bx\kronecker \by$.
This also indicates $\det(\bA\kronecker\bB)=\det(\bB\kronecker\bA)$ and $\trace(\bA\kronecker\bB)=\trace(\bB\kronecker\bA)$.
\end{lemma}
\begin{exercise}[Determinant and Trace of Kronecker Product]
Let $\bA\in\real^{n\times n}$ and $\bB\in\real^{m\times m}$. Show that 
$$
\det(\bA\kronecker\bB)=\det(\bB\kronecker\bA)=(\det(\bA))^m(\det(\bB))^n
$$ 
and 
$$
\trace(\bA\kronecker\bB)=\trace(\bB\kronecker\bA) =\trace(\bA)\trace(\bB).
$$
\end{exercise}

Before delving into the discussion of  the pseudo-inverse and orthogonality in the context of the Kronecker product (which will be proved useful in tensor decomposition methods, as covered in Part~\ref{part:tensor_decom}), we introduce the associative and distributive properties of Kronecker products without a proof. Detailed proofs can be found in \citet{horn1994topics}.
\begin{remark}[Properties of Kronecker Products]
The Kronecker product is \textit{scalar homogeneous}:
\begin{equation}
\bA\kronecker (\gamma\bB) = (\gamma\bA)\kronecker\bB = \gamma(\bA\kronecker\bB).
\end{equation}
The Kronecker product is \textit{associative}, i.e.,
\begin{equation}\label{equation:kro_associ}
(\bA\kronecker \bB) \kronecker \bC =  \bA\kronecker (\bB\kronecker \bC).
\end{equation}
The Kronecker product is \textit{right-distributive}, i.e.,
$$
(\bA\pm \bB) \kronecker \bC = \bA \kronecker \bC \pm \bB \kronecker \bC.
$$
The Kronecker product is \textit{left-distributive}, i.e.,
$$
\bA\kronecker (\bB\pm \bC)  = \bA \kronecker \bB \pm \bA \kronecker \bC.
$$
Taking the transpose before carrying out the Kronecker product yields the same result as doing so afterwards, i.e.,
\begin{equation}\label{equation:kro_trans_pro}
(\bA\kronecker \bB)^\top=\bA^\top \kronecker \bB^\top.
\end{equation}
Taking the inverse before carrying out the Kronecker product yields the same result as doing so afterwards, i.e.,
$$
(\bA\kronecker \bB)^{-1}=\bA^{-1} \kronecker \bB^{-1}.
$$
The Kronecker product of an $m\times m$ identity matrix and an $n\times n$ identity matrix is equal to an $mn \times mn$ identity matrix, i.e., $\bI_m\kronecker \bI_n = \bI_{mn}$.
\end{remark}
Using equality~\eqref{equation:kro_trans_pro}, we can thus prove Equation~\eqref{equation:kron_pro1_trans} from Equation~\eqref{equation:kron_pro1_notrans}, or vice versa. Additional properties of the Kronecker product are discussed in Problems~\ref{prob:kronec_1}$\sim$\ref{prob:kronec_n}.

\subsection*{Pseudo-Inverse in Kronecker Product}
Moreover,  the pseudo-inverse of $(\bA\kronecker \bB)$ is given by
\begin{equation}\label{equation:otimes-psesudo}
(\bA\kronecker \bB)^+ = \bA^+ \kronecker \bB^+,
\end{equation}
where $\bA^+$ represents the pseudo-inverse of matrix $\bA$ \citep{van2000ubiquitous}. Recall that the pseudo-inverse of a  matrix $\bA$ with full column rank is simply $\bA^+=(\bA^\top\bA)^{-1}\bA^\top$. 
When both $\bA$ and $\bB$ are  semi-orthogonal (a matrix is semi-orthogonal if its columns are mutually orthogonal and have unit length, see definition in Section~\ref{section:orthogonal-orthonormal-qr}), their pseudo-inverses are $\bA^+=\bA^\top$ and $\bB^+=\bB^\top$, respectively. And it follows that 
\begin{equation}\label{equation:otimes-psesudo-semi}
(\bA\kronecker \bB)^+ = \bA^\top \kronecker \bB^\top,
\quad\text{with semi-orthogonal }\bA,\bB.
\end{equation}
Analogously, the  pseudo-inverse can be applied to the Kronecker product of a sequence of matrices:
\begin{equation}\label{equation:otimes-psesudo-semi22}
(\bA\kronecker \bB \kronecker\bC\ldots )^+ = \bA^\top \kronecker \bB^\top\kronecker \bC^\top \ldots.
\quad\text{with semi-orthogonal }\bA,\bB, \bC,\dots.
\end{equation}

\subsection*{Orthogonality in Kronecker Product}
Suppose further $\ba\in \real^I$, and $\bb_1, \bb_2\in \real^K$ with $\bb_1^\top\bb_2=0$, then 
\begin{equation}\label{equation:orthogonal-in-kronecker1}
(\ba\kronecker \bb_1)  \perp (\ba\kronecker \bb_2).
\end{equation}
Or suppose $\ba_1,\ba_2\in \real^I$ with $\ba_1^\top\ba_2=0$, and $\bb\in \real^K$, then 
\begin{equation}\label{equation:orthogonal-in-kronecker2}
(\ba_1\kronecker\bb) \perp (\ba_2 \kronecker \bb).
\end{equation}
These two findings imply that $\bA\kronecker \bB$ contains mutually orthogonal (resp. orthonormal) columns if both $\bA \in \real^{I\times J}$ and $\bB\in \real^{K\times L}$ contain mutually orthogonal (resp. orthonormal) columns:
\begin{equation}
\text{semi-orthogonal } \bA\text{ and }\bB 
\quad  \implies \quad 
\text{semi-orthogonal } \bA\kronecker \bB.
\end{equation}

\subsection{Khatri-Rao Product}

The Khatri-Rao product, often denoted by $\khatrirao$, is a specialized matrix product that is closely related to the Kronecker product but is typically used in settings where the matrices involved are column-wise partitioned.
It  is a key component in tensor decompositions, especially in the CP decomposition (Section~\ref{section:cp_decom}). In CP decomposition, a tensor is approximated as a sum of rank-one tensors, and the factors of these rank-one components can be represented using the Khatri-Rao product of matrices. This makes it a powerful tool for analyzing multi-way data. 
When dealing with multiple data sets or signals that need to be combined, the Khatri-Rao product can be used to integrate information from different sources. For example, in sensor networks, where each sensor provides a different type of measurement, the Khatri-Rao product allows for the creation of a comprehensive model that takes into account all available data.
In statistical analysis, especially in the context of multivariate analysis, the Khatri-Rao product is used to manipulate and analyze covariance matrices or correlation structures. It can help in understanding the relationships between different variables across multiple dimensions. This is particularly relevant in fields like genomics, where interactions between genetic markers are studied.

\begin{definition}[Khatri-Rao Product]\label{definition:khatri-rao-product}
The \textit{Khatri-Rao product} of two matrices  $\bA \in \real^{I\times K}$ and $\bB\in \real^{J\times K}$,  denoted by $\bA\khatrirao \bB$, is defined as follows:
$$
\bA\khatrirao\bB =
\begin{bmatrix}
\ba_1\kronecker \bb_1 & \ba_2\kronecker \bb_2 & \ldots & \ba_K\kronecker \bb_K
\end{bmatrix},
$$
which is a matrix of size $(IJ)\times K$. And it is known as the ``matching column-wise" Kronecker product.

\end{definition}

For partitioned matrices, there is a partition-wise counterpart of the Khatri-Rao product, denoted by $\khatrirao_b$.
\begin{definition}[Partition-Wise Khatri-Rao Product]\label{definition:partition_khatri_prod}
The \textit{partition-wise Khatri-Rao product} of two matrices  
$$\bA=[\underbrace{\bA_1}_{\real^{I\times M_1}}, \underbrace{\bA_2}_{\real^{I\times M_2}}, \ldots, \underbrace{\bA_R}_{\real^{I\times M_R}}]
\gap \text{and}\gap 
\bB=[\underbrace{\bB_1}_{\real^{J\times N_1}}, \underbrace{\bB_2}_{\real^{J\times N_2}}, \ldots, \underbrace{\bB_R}_{\real^{J\times N_R}}],$$  
denoted by $\bA\khatrirao_b \bB$, is defined as follows:
$$
\bA\khatrirao_b\bB =
\begin{bmatrix}
	\bA_1\kronecker \bB_1 & \bA_2\kronecker \bB_2 & \ldots & \bA_R\kronecker \bB_R
\end{bmatrix},
$$
which is a matrix of size $(IJ)\times (\sum_{i=1}^{R} M_iN_i)$.
\end{definition}

The Khatri-Rao product shares some properties with the Kronecker product, but it is tailored for operations involving matrices with a common number of columns.
Based on the  definition of the Khatri-Rao product, for two vectors $\ba$ and $\bb$, we can observe that the Khatri-Rao product of these two vectors is equivalent to their Kronecker product:
\begin{equation}
\ba \khatrirao \bb = \ba\kronecker\bb.
\end{equation}
Given three matrices $\bA, \bB$, and $\bC$, the ``distributive law" and ``associative law" for the Khatri-Rao product follow that 
$$
\begin{aligned}
\textit{distributive law}:&\quad 	(\bA+ \bB)\khatrirao \bC = (\bA\khatrirao \bC)+ (\bB\khatrirao\bC) ;\\
\textit{associative law}: &\quad \bA\khatrirao \bB\khatrirao \bC = (\bA\khatrirao \bB)\khatrirao \bC = \bA\khatrirao (\bB\khatrirao \bC).
\end{aligned}
$$
Additionally, when $\bA, \bB\in\real^{I\times K}$ are semi-orthogonal matrices, the Khatri-Rao product $\bA\khatrirao \bB$ is also semi-orthogonal (see Problem~\ref{problem:semi_khatri}):
\begin{equation}
\text{semi-orthogonal } \bA\text{ and }\bB 
\quad  \implies \quad 
\text{semi-orthogonal } \bA\khatrirao \bB.
\end{equation}

\subsection{More Properties of Special Matrix Products}
\paragraph{\colorbox{\mdframecolorBox}{$(\bA\khatrirao \bB )^\top (\bA\khatrirao \bB )=
(\bA^\top\bA) \hadaprod (\bB^\top\bB)$}.}
Moreover, we observe that for two matrices   $\bA \in \real^{I\times K}$ and $\bB\in \real^{J\times K}$, the following relationship holds:
\begin{equation}
\bZ \triangleq (\bA\khatrirao \bB )^\top (\bA\khatrirao \bB )=
\begin{bmatrix}
(\ba_1\kronecker \bb_1)^\top  \\
(\ba_2\kronecker \bb_2)^\top  \\
\vdots \\ 
(\ba_K\kronecker \bb_K)^\top
\end{bmatrix}
\begin{bmatrix}
\ba_1\kronecker \bb_1 & \ba_2\kronecker \bb_2 & \ldots & \ba_K\kronecker \bb_K
\end{bmatrix},
\end{equation}
where $\bZ\in \real^{K\times K}$, and each entry $(i,j)$, denoted by $z_{ij}$, is given by 
$$
z_{ij} = (\ba_i\kronecker \bb_i)^\top (\ba_j\kronecker \bb_j) = (\ba_i^\top\ba_j)(\bb_i^\top\bb_j),
$$
where the last equality is derived from Equation~\eqref{equation:kronecker-vector-find2}. Therefore, $\bZ$ can be  equivalently  expressed as 
\begin{equation}\label{equation:two-khatri-rao-pro-equi}
\bZ = (\bA\khatrirao \bB )^\top (\bA\khatrirao \bB ) =
(\bA^\top\bA) \hadaprod (\bB^\top\bB).
\end{equation}
More generally, this property also holds for three or more matrices; see Equation~\eqref{equ:aa_hada_bb2}.

\paragraph{\colorbox{\mdframecolorBox}{$(\bA\khatrirao \bB )^\top (\bC\khatrirao \bD )=
(\bA^\top\bC) \hadaprod (\bB^\top\bD)$}.}
Similarly, given  $\bA,\bC \in \real^{I\times K}$ and $\bB, \bD\in \real^{J\times K}$, it follows that 
\begin{equation}
(\bA\khatrirao \bB )^\top (\bC\khatrirao \bD ) =
(\bA^\top\bC) \hadaprod (\bB^\top\bD).
\end{equation}

\paragraph{\colorbox{\mdframecolorBox}{$(\bA\kronecker\bB)(\bC\khatrirao \bD)=(\bA\bC)\khatrirao(\bB\bD)$}.}
To see this, consider $\bA\in\real^{I\times J}$, $\bB\in\real^{K\times L}$, $\bC\in\real^{J\times P}$, and $\bD\in\real^{L\times P}$. 
The matrix  $(\bC\khatrirao \bD)$ has a shape of ${JL\times P}$. Each column of $(\bC\khatrirao \bD)$ is the Kronecker product of the corresponding columns of $\bC$ and $\bD$. Specifically, the $p$-th column of $(\bC\khatrirao \bD)$ is $\bc_p\kronecker \bd_p$, $p\in\{1,2,\ldots, P\}$.
Therefore, the $p$-th column of the left-hand side, by Equation~\eqref{equation:kronecker_eq3}, is
$$
(\bA\kronecker\bB)(\bc_p\kronecker \bd_p)=(\bA\bc_p)\kronecker(\bB \bd_p).
$$
For the right-hand side, the $p$-th columns of the matrices $(\bA\bC)$ and $(\bB\bD)$ are $\bA\bc_p$ and $\bB\bd_p$, respectively. Thus, the $p$-th column of the Khatri-Rao product  $(\bA\bC)\khatrirao(\bB\bD)$ is
$$
(\bA\bc_p)\kronecker (\bB\bd_p).
$$
Since the $p$-th columns of both sides of the equation are identical for all  $p\in\{1,2,\ldots, P\}$, this concludes the proof that:
$$
(\bA\kronecker\bB)(\bC\khatrirao \bD)=(\bA\bC)\khatrirao(\bB\bD).
$$

\paragraph{\colorbox{\mdframecolorBox}{$(\bA\khatrirao\bB)\hadaprod(\bC\khatrirao \bD)=(\bA\hadaprod\bC)\khatrirao(\bB\hadaprod\bD)$}.}
To demonstrate this, given $\bA, \bC\in\real^{I\times K}$ and $\bB, \bD\in\real^{J\times K}$.
The $k$-th column of the left-hand side is $(\ba_k\kronecker\bb_k)\hadaprod (\bc_k\kronecker\bd_k)$, and  the $k$-th column of the right-hand side is $(\ba_k\hadaprod\bc_k)\kronecker (\bb_k\hadaprod\bd_k)$. These two expressions are equivalent, which proves the result.

\paragraph{\colorbox{\mdframecolorBox}{$(\bA\khatrirao\bB)^+=\big((\bA^\top\bA)\hadaprod(\bB^\top\bB)\big)^+ (\bA\khatrirao\bB)^\top$}.}
This result follows directly from the definition of the pseudo-inverse and Equation~\eqref{equation:two-khatri-rao-pro-equi}.
More generally, it can be shown that $(\bA\khatrirao\bB \khatrirao\bC)^+=\big((\bA^\top\bA)\hadaprod(\bB^\top\bB)\hadaprod(\bC^\top\bC)\big)^+ (\bA\khatrirao\bB\khatrirao\bC)^\top$

To summarize, given the definitions of the Hadamard, Kronecker, and Khatri-Rao products, we have (some results are provided without proof):
\begin{mdframed}[hidealllines=\mdframehidelineNote,backgroundcolor=\mdframecolor,frametitle={Properties of Special Matrix Products}]
\begin{align}
&	\bb\kronecker \bc^\top &=&\bc^\top \kronecker \bb;\label{equ:kro_vec1}\\
&\ba\kronecker 	\bb\kronecker \bc^\top &=&\ba\kronecker\bc^\top \kronecker \bb = \bc^\top \kronecker\ba\kronecker \bb;\label{equ:kro_vec2}\\
&	(\ba\kronecker \bb)^\top(\bc\kronecker \bd)&=&(\ba^\top\bc)(\bb^\top\bd);\\
&	(\ba\kronecker \bb)^\top(\ba\kronecker \bb) &=&\normtwo{\ba}^2 \cdot \normtwo{\bb}^2;\\
&	(\bA\kronecker \bB)^\top (\bC\kronecker \bD) &=& (\bA^\top\bC) \kronecker (\bB^\top\bD), \gap \left(\text{\parbox{11em}{with $\bA,\bC$ same shape, \\$\bB,\bD$ same shape}}\right);\\
&(\bA\kronecker \bB) (\bC\kronecker \bD) &=& (\bA\bC) \kronecker (\bB\bD);\label{equ:mat_kro_prod}\\
&(\bA\kronecker \bB) (\bc \kronecker \bd) &=& (\bA\bc) \kronecker (\bB\bd);\\
& (\bA\kronecker \bB)^+ &=& \bA^+ \kronecker \bB^+; \qquad (\bA\kronecker \bB)^- = \bA^- \kronecker \bB^-;\\
&	\bA\khatrirao \bB\khatrirao \bC &=& (\bA\khatrirao \bB)\khatrirao \bC 	= \bA\khatrirao (\bB\khatrirao \bC);\\
&	(\bA\khatrirao \bB )^\top (\bA\khatrirao \bB ) &=& 	(\bA^\top\bA) \hadaprod (\bB^\top\bB);\label{equ:aa_hada_bb}\\
&	(\bA\khatrirao \bB\khatrirao\bC )^\top (\bA\khatrirao \bB\khatrirao\bC ) &=& 	(\bA^\top\bA) \hadaprod (\bB^\top\bB)\hadaprod (\bC^\top\bC);\label{equ:aa_hada_bb2}\\
&(\bA\khatrirao \bB )^\top (\bC\khatrirao \bD ) &=&  (\bA^\top\bC) \hadaprod (\bB^\top\bD);\label{equ:ab_top_cd}\\
&(\bA\kronecker\bB)(\bC\khatrirao \bD) &=& (\bA\bC)\khatrirao(\bB\bD);\label{equ:atimesb_codotb}\\
& (\bA\khatrirao\bB)\hadaprod(\bC\khatrirao \bD) &=& (\bA\hadaprod\bC)\khatrirao(\bB\hadaprod\bD);\\
&(\bA\khatrirao\bB)^+ &=&\big((\bA^\top\bA)\hadaprod(\bB^\top\bB)\big)^+ (\bA\khatrirao\bB)^\top;\\
&(\bA\khatrirao\bB \khatrirao\bC)^+&=&\big((\bA^\top\bA)\hadaprod(\bB^\top\bB)\hadaprod(\bC^\top\bC)\big)^+ (\bA\khatrirao\bB\khatrirao\bC)^\top,
\end{align}
\end{mdframed}
where \eqref{equ:kro_vec2} follows from \eqref{equ:kro_vec1} and \eqref{equation:kro_associ}.

\begin{figure}[h]
	\centering
	\includegraphics[width=1\textwidth]{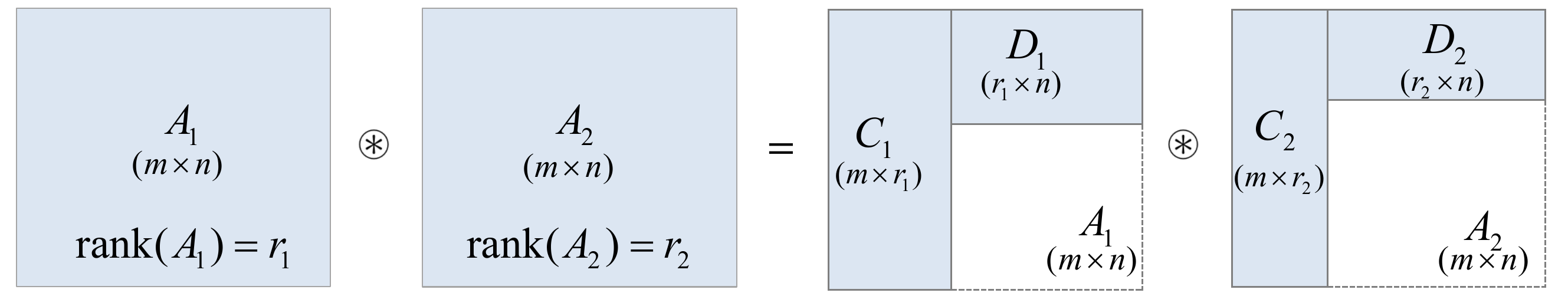}
	\caption{Diagram illustrating the rank in a Hadamard product.}
	\label{fig:hada_rank}
\end{figure}

We conclude this section by describing that the ranks for Hadamard products and Kronecker products are ``multiplicative", and the rank for Khatri-Rao products is nondecreasing.
\begin{theorem}[Rank of Hadamard Products, \citep{horn1990hadamard}]\label{theorem:rank_hada_prod}
Let $\bA_1,\bA_2\in\real^{m\times n}$ be any $m\times n$ matrices with rank $r_1$ and rank $r_2$, respectively.
Then their Hadamard product $\bA_1\hadaprod\bA_2$ has rank at most $r_1\cdot r_2$: $\rank(\bA_1\hadaprod\bA_2)\leq \rank(\bA_1)\rank(\bA_2)$.
\end{theorem}
\begin{proof}[of Theorem~\ref{theorem:rank_hada_prod}]
For any rank-$r_1$ matrix $\bA_1$ and rank-$r_2$ matrix $\bA_2$, they can be expressed as sums of rank-1 matrices:
$$
\begin{aligned}
\bA_1&=\bC_1\bD_1^\top = \sum_{i=1}^{r_1} \bc_{1i}\bd_{1i}^\top
\qquad\text{and}\qquad
\bA_2=\bC_2\bD_2^\top = \sum_{j=1}^{r_2} \bc_{2j}\bd_{2j}^\top,
\end{aligned}
$$
where $\bC_1\in\real^{m\times r_1}$ and $\bD_1\in\real^{n\times r_1}$ are rank-$r_1$ matrices,  $\bC_2\in\real^{m\times r_2}$ and $\bD_2\in\real^{n\times r_2}$ are rank-$r_2$ matrices, and $\bc_{1i}, \bd_{1i}$ are the $i$-th column of $\bC_1, \bD_1$ (same for $\bc_{2j},\bd_{2j}$, see Figure~\ref{fig:hada_rank}).
Therefore, 
$$
\bA_1\hadaprod\bA_2 = \left(\sum_{i=1}^{r_1} \bc_{1i}\bd_{1i}^\top\right)\hadaprod  \left(\sum_{j=1}^{r_2} \bc_{2j}\bd_{2j}^\top\right)
=\sum_{i=1}^{r_1}\sum_{j=1}^{r_2} \left(\bc_{1i}\bd_{1i}^\top \right) \hadaprod \left(\bc_{2j}\bd_{2j}^\top\right).
$$
The Hadamard product is the sum of $r_1\cdot r_2$ rank-1 matrices, and thus has rank at most $r_1\cdot r_2$.
\end{proof}

\begin{theorem}[Rank, Trace of Kronecker Products]\label{theorem:rank_kronec_prod}
Let $\bA \in \real^{I\times J}$ and $\bB\in \real^{K\times L}$. Then, 
\begin{equation}
\rank(\bA\kronecker \bB) = \rank(\bA)\rank(\bB).
\end{equation}
This indicates that the rank is multiplicative under the Kronecker product. Consequently, if either  $\bA=\bzero$ or $\bB=\bzero$, then $\bA\kronecker\bB=\bzero $.
For square and symmetric matrices $\bA$ and $\bB$, we also have 
\begin{flalign}
&\textbf{(Square and symmetric $\bA,\bB$): }\quad \quad\quad
\trace(\bA\kronecker \bB) = \trace(\bA)\trace(\bB).&
\quad\quad\quad
\end{flalign}
When $\bA\in\real^{I\times J}$ and $\bB\in\real^{I\times J}$ have the same dimension, Exercise~\ref{exercise:hada_sub_kro} also implies
\begin{flalign}
&\textbf{(Same shape $\bA,\bB$): }\quad\quad
\rank(\bA\hadaprod \bB) \leq \rank(\bA\kronecker \bB) = \rank(\bA)\rank(\bB).&
\end{flalign}
If $\bA$ is PD and $\bB$ is PSD, and they have the same dimension, then 
\noindent
\begin{flalign}\label{equation:rkjkro_pdpsd}
&\textbf{(PD $\bA$, PSD $\bB$): }\, \rank(\bB)\leq \rank(\bA\hadaprod \bB) \leq \rank(\bA\kronecker \bB) = \rank(\bA)\rank(\bB).&
\end{flalign}
\end{theorem}
\begin{proof}[of Theorem~\ref{theorem:rank_kronec_prod}]
Suppose $\bA$ and $\bB$ admit full SVD as $\bU_1^\top\bA\bV_1=\bSigma_1$ and $\bU_2^\top\bB\bV_2=\bSigma_2$, respectively. 
Using Equality~\eqref{equ:mat_kro_prod}, we have 
$$
\begin{aligned}
\bSigma_1\kronecker \bSigma_2 &= (\bU_1^\top\bA\bV_1)\kronecker (\bU_2^\top\bB\bV_2)\\
&= (\bU_1^\top\kronecker \bU_2^\top) (\bA\bV_1\kronecker \bB\bV_2)
=(\bU_1^\top\kronecker \bU_2^\top)(\bA\kronecker \bB)(\bV_1\kronecker \bV_2).
\end{aligned}
$$
Since the Kronecker product of orthogonal matrices is also orthogonal (Lemma~\ref{lemma:krokecker_keep_special}), we have, by Proposition~\ref{proposition:eigenvalue-similar-matrices},
$$
\rank(\bSigma_1\kronecker \bSigma_2) = \rank(\bA\kronecker \bB),
$$
where $\bSigma_1\kronecker \bSigma_2$ has $\rank(\bSigma_1)\rank(\bSigma_2)$ nonzero entries, indicating that $\rank(\bSigma_1\kronecker \bSigma_2)=\rank(\bSigma_1)\rank(\bSigma_2)$.

For the second part, since $\bA$ and $\bB$ are symmetric, $\bSigma_1\kronecker \bSigma_2$ and $\bA\kronecker \bB$ are similar matrices,  their traces are the same (Proposition~\ref{proposition:eigenvalue-similar-matrices}).

To prove \eqref{equation:rkjkro_pdpsd} (see \citet{djokovic1965hadamard}). If $\bB=\bzero$, the proof is trivial. Suppose $\rank(\bB)=r$. Since $\bB$ is PSD, there exists a principal submatrix $\bB[I, I]\in\real^{r\times r}$ such that $\bB[I, I]$ is nonsingular and PD.
Exercise~\ref{exercise:schur_product} shows that $\bA[I,I]\hadaprod\bB[I,I]$ is also nonsingular PD.
This completes the proof.
\end{proof}

\begin{theorem}[Rank of Khatri-Rao Products]\label{theorem:rank_khatri_prod}
Let $\bA \in \real^{I\times K}$ and $\bB\in \real^{J\times K}$. Then, 
$$
\rank(\bA\khatrirao  \bB) \geq  \max\{\rank(\bA),\rank(\bB)\}.
$$
\end{theorem}
\begin{proof}[of Theorem~\ref{theorem:rank_khatri_prod}]
For any matrix $\bC$, we can find a nonsingular $\bS$ such that there exists a row of $\bS\bC$ with all nonzeros.
Therefore, assume there exist nonsingular matrices $\bP$ and $\bT$ such that $\bP\bA$ and $\bT\bB$ each have at least one row containing all nonzero elements.
Since $\rank(\bP\bA)=\rank(\bA)$, $\rank(\bT\bB)=\rank(\bB)$, and $(\bP\bA)\khatrirao (\bT\bB) = (\bP\kronecker \bT)(\bA\khatrirao \bB)$ with $\rank(\bA\khatrirao \bB) = \rank((\bP\kronecker \bT)(\bA\khatrirao \bB))$ (Lemma~\ref{lemma:left_mul_krank}, Equation~\eqref{equ:atimesb_codotb}, and $(\bP\kronecker \bT)$ is also nonsingular from Lemma~\ref{lemma:krokecker_keep_special}), it suffices to show that $\rank(\bP\bA\khatrirao \bT\bB) \geq \max\{\rank(\bP\bA),\rank(\bT\bB)\}$.

Given that $\bT\bB$ contains a row that has all nonzero elements, then $(\bP\bA\khatrirao \bT\bB)$ contains a submatrix equal to $\bP\bA$, rescaled columnwise by the elements of that row of $\bT\bB$. Therefore, $\rank(\bP\bA\khatrirao \bT\bB) \geq \rank(\bP\bA)$. Similarly, we can prove $\rank(\bP\bA\khatrirao \bT\bB) \geq \rank(\bT\bB)$. This completes the proof.
\end{proof}

\begin{theorem}[$k$-Rank of Khatri-Rao Products]\label{theorem:kkrank_khatri_prod}
Let $\bA \in \real^{I\times K}$ and $\bB\in \real^{J\times K}$. Then,
$$
\rank_k(\bA\khatrirao  \bB) \geq  \min\{\rank_k(\bA)+\rank_k(\bB)-1, K\},
$$
where $\rank_k(\cdot)$ denotes the Kruskal rank (Definition~\ref{definition:krus_rk}).
\end{theorem}
\begin{proof}[of Theorem~\ref{theorem:kkrank_khatri_prod}]
The proof is based on \citet{ten2000k}.
Without loss of generality, we assume the $k$-rank of $\bA\khatrirao \bB$ is less than $K$; otherwise the proof is trivial.
Let $S$ be the smallest number of linearly dependent columns of $\bA\khatrirao \bB$. Therefore, $\rank_k(\bA\khatrirao  \bB)  = S-1$.
We collect some set of $S$ linearly dependent columns of $\bA\khatrirao  \bB$ into $\bA_S\khatrirao \bB_S$, where $\bA_S$ and $\bB_S$ contain the corresponding columns of $\bA$ and $\bB$, respectively.
Thus, there exists a vector $\bc_S$ with nonzero elements satisfying $(\bA_S\khatrirao \bB_S)\bc_S=\bzero$; otherwise we can reduce the number $S$.
This implies $\bA_S\bC_S\bB_S^\top=\bzero $, where $\bC_S=\diag(\bc_S)$ is nonsingular.
Therefore, we have 
$$
0=\rank(\bA_S\bC_S\bB_S^\top)\geq \rank(\bA_S)+\rank(\bB_S)-S.~\footnote{Follows from the Frobenius inequality: $\rank(\bA\bB\bC)\geq \rank(\bA\bB)+\rank(\bB\bC)-\rank(\bB)$.}
$$
Since $\rank(\bA_S)\geq \rank_k(\bA_S)\geq \rank_k(\bA)$ and $\rank(\bB_S)\geq \rank_k(\bB_S)\geq \rank_k(\bB)$, this completes the proof.
\end{proof}

\section{Hadamard and Kronecker Representations of  Matrix Decomposition}

The Hadamard/Kronecker product  of matrices that have been decomposed into their constituent parts (decompositional forms) can be computed by directly multiplying these factored components. When matrices are expressed in terms of their decompositions, such as through singular value decomposition, eigenvalue decomposition, or other factorizations, the Hadamard/Kronecker product can often be simplified and  more clearly  expressed in terms of these factors.

By examining the Hadamard/Kronecker product through the lens of the decompositions, we might uncover patterns or structures within the product that would not be as apparent if we were to simply compute the Hadamard/Kronecker product of the original matrices. For example, if two matrices are decomposed into their singular values and orthogonal matrices, the Hadamard/Kronecker product of these decomposed forms can sometimes reveal insights about the singular values or the orthogonality of the resulting matrix.

This approach can be particularly useful in various applications, such as signal processing, machine learning, and numerical linear algebra, where understanding the structure of the Hadamard/Kronecker product can lead to more efficient algorithms or better comprehension of the underlying data.
\begin{theoremHigh}[Hadamard of Spectral]
Let $\bA,\bB\in\real^{n\times n}$ be two symmetric matrices with spectral decompositions $\bA=\bU\bLambda\bU^\top=\sum_{i=1}^{n}\lambda_i\bu_i\bu_i^\top$ and $\bB=\bV\bD\bV^\top=\sum_{i=1}^{n}d_i\bv_i\bv_i^\top$, respectively.
Then, their Hadamard product $\bA\hadaprod\bB$ is also symmetric with the following decomposition 
\begin{equation}\label{equation:hada_spec_e1}
\bA\hadaprod\bB= \sum_{i=1}^n\sum_{j=1}^n \lambda_i d_j (\bu_i\bu_i^\top)\hadaprod (\bv_j\bv_j^\top)
=\sum_{i=1}^n\sum_{j=1}^n \lambda_i d_j (\bu_i\hadaprod \bv_j) (\bu_i\hadaprod \bv_j)^\top.
\end{equation}
Using Equation~\eqref{equ:aa_hada_bb} on vectors, we also have 
$$
\bA\hadaprod\bB= \sum_{i=1}^n\sum_{j=1}^n \lambda_i d_j (\bu_i\bu_i^\top)\hadaprod (\bv_j\bv_j^\top)
=\sum_{i=1}^n\sum_{j=1}^n \lambda_i d_j   (\bu_i^\top \khatrirao \bv_j^\top )^\top (\bu_i^\top\khatrirao \bv_j^\top )
$$
When $\bA,\bB$ are positive semidefinite, we have $\bA\triangleq\bC^\top\bC$ and $\bB\triangleq\bD^\top\bD$ with $\bC\triangleq\bLambda^{1/2}\bU^\top$ and $\bD\triangleq\bD^{1/2}\bV$, respectively (since the eigenvalues are nonnegative for positive semidefinite matrices).
Using Equation~\eqref{equ:aa_hada_bb}, we have 
$$
\bA\hadaprod\bB= (\bC^\top\bC) \hadaprod (\bD^\top\bD)=(\bC\khatrirao \bD )^\top (\bC\khatrirao \bD ).
$$
When $\bA,\bB$ are positive definite (resp. positive semidefinite), then $\bA\hadaprod \bB$ is also positive definite (resp. positive semidefinite). This is known as  \textit{Schur's product theorem} (Exercise~\ref{exercise:schur_product}).
\end{theoremHigh}
\begin{proof}[of Schur's product theorem]
From the above discussion, it can be readily shown that if $\bA, \bB$ are positive semidefinite, then $\bA\hadaprod\bB$ is also positive semideifnite.
For positive definite matrices, we have $\lambda_id_j>0, i,j\in\{1,2,\ldots,n\}$ in Equation~\eqref{equation:hada_spec_e1}. 
And for any $\bx\neq \bzero$, we have $\bx^\top (\bu_i\hadaprod \bv_j) (\bu_i\hadaprod \bv_j)^\top \bx\geq 0$. It suffices to show that there exist  $i,j$ such that $\bx^\top (\bu_i\hadaprod \bv_j) (\bu_i\hadaprod \bv_j)^\top \bx\neq  0$.
To see this, we have 
$$
\bx^\top (\bu_i\hadaprod \bv_j) (\bu_i\hadaprod \bv_j)^\top \bx 
=\sum_{k=1}^{n} (x_k u_{i,k} v_{j,k})^2.
$$
Since $\bA$ is positive definite, there exist an $i$ such that $\bx\hadaprod \bu_i \neq \bzero$. Otherwise, $\bx^\top\bu_i=0$ for all $i\in\{1,2,\ldots,n\}$; and thus $\bx^\top\bA\bx=\sum_{i=1}^{n}\lambda_i \bx^\top\bu_i\bu_i^\top\bx=0$.
Likewise, for $\bx\hadaprod \bu_i\neq \bzero$, there exist a $j$ such that $\bv_j \hadaprod (\bx\hadaprod \bu_i)\neq \bzero$.
This completes the proof.~\footnote{The proof can also be carried out using the trace norm or through Gaussian integration.}
\end{proof}

\begin{theoremHigh}[Hadamard of SVD]\label{theorem:hada_svd}
Let $\bA,\bB\in\real^{m\times n}$ be  two matrices with full SVD $\bA=\bU_1\bSigma_1\bV_1^\top$ and $\bB=\bU_2\bSigma_2\bV_2^\top$, respectively. Then, the Hadamard product $\bA\hadaprod \bB$ admits the following decomposition: 
$$
\bA\hadaprod\bB = (\bU_1^\top\khatrirao\bU_2^\top)^\top (\bSigma_1 \kronecker \bSigma_2) (\bV_1^\top \khatrirao \bV_2^\top).
$$
\end{theoremHigh}
\begin{proof}[of Theorem~\ref{theorem:hada_svd}]
Using Equations~\eqref{equ:ab_top_cd} and~\eqref{equ:atimesb_codotb}, we have 
$$
\begin{aligned}
\bA\hadaprod\bB &= (\bU_1\bSigma_1\bV_1^\top)\hadaprod (\bU_2\bSigma_2\bV_2^\top)\\
&=(\bU_1^\top \khatrirao \bU_2^\top)^\top (\bSigma_1\bV_1^\top \khatrirao \bSigma_2\bV_2^\top)\\
&= (\bU_1^\top \khatrirao \bU_2^\top)^\top(\bSigma_1\kronecker \bSigma_2)(\bV_1^\top \khatrirao \bV_2^\top).
\end{aligned}
$$
Since $\bU_1, \bV_1, \bU_2, \bV_2$ are orthogonal matrices, $(\bU_1^\top \khatrirao \bU_2^\top)$ and $(\bV_1^\top \khatrirao \bV_2^\top)$ are semi-orthogonal.
And $(\bSigma_1\kronecker \bSigma_2)$ is diagonal.~\footnote{However, the form does not constitute a SVD of the Hadamard product.}
This completes the proof.
\end{proof}

Lemma~\ref{lemma:krokecker_keep_special} shows that the Kronecker product of two permutation, two triangular, or two orthogonal matrices results in permutation, triangular, or orthogonal matrices, respectively. This implies that the Kronecker product of two matrices with known decompositions also has a known decomposition form, due to \eqref{equ:mat_kro_prod}.
\begin{theoremHigh}[Kronecker Decomposition\, \citep{van1993approximation, schacke2004kronecker}]
Given the LU decomposition of $\bA$ and $\bB$, it follows that 
$$
\text{LU:}\gap 
\left.
\begin{aligned}
\bA&= \bP_1\bL_1\bU_1 \\
\bB&= \bP_2\bL_2\bU_2
\end{aligned}
\right\}
\gap 
\implies 
\gap 
\bA\kronecker \bB = (\bP_1\kronecker \bP_2) (\bL_1\kronecker \bL_2)(\bU_1\kronecker \bU_2).
$$
Given the Cholesky decomposition of $\bA$ and $\bB$, it follows that 
$$
\text{Cholesky:}\gap 
\left.
\begin{aligned}
\bA&= \bR_1^\top \bR_1 \\
\bB&= \bR_2^\top \bR_2
\end{aligned}
\right\}
\gap 
\implies 
\gap 
\bA\kronecker \bB = (\bR_1\kronecker \bR_2)^\top (\bR_1\kronecker \bR_2).
$$
Given the QR decomposition of $\bA$ and $\bB$, it follows that 
$$
\text{QR:}\gap 
\left.
\begin{aligned}
\bA&= \bQ_1 \bR_1 \\
\bB&= \bQ_2 \bR_2
\end{aligned}
\right\}
\gap 
\implies 
\gap 
\bA\kronecker \bB = (\bQ_1\kronecker \bQ_2) (\bR_1\kronecker \bR_2).
$$
Given the SVD of $\bA$ and $\bB$, it follows that 
$$
\text{SVD:}\gap 
\left.
\begin{aligned}
\bA&= \bU_1\bSigma_1 \bV_1^\top \\
\bB&= \bU_2\bSigma_2 \bV_2^\top
\end{aligned}
\right\}
\gap 
\implies 
\gap 
\bA\kronecker \bB = (\bU_1\kronecker \bU_2) (\bSigma_1\kronecker \bSigma_2) (\bV_1\kronecker \bV_2)^\top.
$$
\end{theoremHigh}


\section{Low-Rank Hadamard Decomposition}\label{section:low_rank_hadamard}
In the field of linear algebra and data analysis, matrix decomposition techniques are essential for extracting meaningful information from complex datasets. One common objective is to approximate a given matrix using a lower-rank representation, which simplifies the data while preserving its key characteristics. The Hadamard product, also known as the element-wise product, provides an alternative approach to matrix decomposition compared to traditional matrix multiplication.

As discussed  in Chapter~\ref{chapter:als}, the alternating least squares (ALS) algorithm is an iterative method used to find the best low-rank approximation of a matrix by decomposing it into two or more matrices. ALS is particularly advantageous for large-scale problems, such as those found in recommender systems, where the goal is to predict missing entries in a user-item interaction matrix. During each iteration, the ALS algorithm alternates between updating one matrix while keeping the other fixed, thereby minimizing the reconstruction error at every step.
Nonnegative matrix factorization (NMF), introduced in Chapter~\ref{chapter:nmf}, is a variant of matrix factorization where both the original matrix and the resulting factorized matrices have nonnegative entries. This constraint makes NMF especially suitable for applications where the data represents quantities that cannot be negative, such as images, audio signals, or document-term matrices in text mining.

Ws further explore  the  Hadamard decomposition of a matrix $\bA$, where $\bA$ can be expressed as the Hadamard product of two low-rank matrices: $\bA=\bA_1\hadaprod \bA_2$. This type of decomposition is advantageous when the data exhibits multiplicative relationships, and a low-rank approximation is desired to reduce complexity or enhance interpretability.
\paragraph{Non-Factorizability  Issue.}
Although when $\bA_1\in\real^{n^2\times n^2}$ and $\bA_2\in\real^{n^2\times n^2}$ share the same rank $n$, the Hadamard product $\bA_1\hadaprod \bA_2$ can achieve a maximum rank of $n^2$ (Theorem~\ref{theorem:rank_hada_prod}), not all matrices $\bA\in\real^{n^2\times n^2}$ of rank $n^2$  can be  represented as the Hadamard product of two lower-rank matrices:
\begin{itemize}
\item The Hadamard decomposition $\bA = \bA_1 \hadaprod  \bA_2$, where $\bA_1$ and $\bA_2$ are rank-$n$ factors, encodes a system of nonlinear equations.
\item This system comprises $n^2 \times n^2 = n^4$ equations (one per entry of $ \bA$) and, due to the low-rank constraint on the two Hadamard factors $\bA_1$ and $ \bA_2$, only $(n^2 n + n n^2) = 2n^3$ variables exist.
\item For $n > 2$, there are more equations than variables, suggesting that all the equations will be simultaneously satisfied only in special cases.
For example, if the matrix $\bA$ includes a row or a column with all but a single entry being zero, then not all the equations in the system can be satisfied \citep{ciaperoni2024hadamard}.

\end{itemize}

Therefore, we focus on solving the low-rank reconstruction problem for the Hadamard decomposition.
In Theorem~\ref{theorem:rank_hada_prod}, assuming that  $\bA_1$ and $\bA_2$ share the same rank $K$, our aim is to reconstruct the design matrix $\bA$ through the Hadamard product $\bA_1\hadaprod\bA_2$. 
Building upon the matrix factorization method used in   alternating least squares (Section~\ref{section:als-netflix}),
we now concentrate on algorithms for solving the \textit{low-rank Hadamard decomposition} problem  (we may use the term ``Hadamard decomposition" to refer to ``low-rank Hadamard decomposition" for brevity, when it is clear from the context; this applies to the low-rank Kronecker and Khatri-Rao decompositions in the sequel):
\begin{itemize}
\item Given a real matrix $\bA\in \real^{M\times N}$, find  matrix factors $\bA_1\in \real^{M\times N}$ and $\bA_2\in \real^{M\times N}$ such that: 
\begin{equation}
\min\,\,L(\bC_1, \bD_1, \bC_2,\bD_2) = \normf{\bA_1\hadaprod \bA_2-\bA}^2\triangleq \normf{(\bC_1\bD_1)\hadaprod (\bC_2\bD_2) -\bA}^2,
\end{equation}
where $\bC_1, \bC_2\in\real^{M\times K}$, and $\bD_1, \bD_2\in\real^{K\times N}$: $\bA_1\triangleq\bC_1\bD_1$ and $\bA_2\triangleq\bC_2\bD_2$ such that $\bA_1$ and $\bA_2$ are rank-$K$ matrices.~\footnote{In the proof of Theorem~\ref{theorem:rank_hada_prod}, we consider the matrices $\bD_1, \bD_2\in\real^{N\times K}$. To abuse notations, we use $\bD_1, \bD_2\in\real^{K\times N}$ here for  ease of deriving gradients in the sequel.}
\end{itemize}
Low-rank (Hadamard) decomposition is often  necessary because many natural phenomena exhibit multiplicative or conjunctive relationships \citep{ciaperoni2024hadamard}.
For instance, consider a study on risk factors for a disease with two predictors: smoking status (yes/no) and alcohol consumption (yes/no). The multiplicative model would account for not only the main effects of smoking and alcohol consumption but also their interaction.
The (low-rank) Hadamard decomposition offers an alternative approach.

Following the alternating descent framework (Section~\ref{section:als-netflix}), in each iteration, the matrices $\bC_1, \bD_1, \bC_2$, and $\bD_2$ are updated sequentially by taking a step in the direction opposite to the gradient of the objective function (a.k.a., the reconstruction error in this context).
It then can be shown that 
$$
\nabla L(\bC_1)=\nabla L(\bC_1|\bD_1, \bC_2,\bD_2)=2\big(\left((\bC_1\bD_1)\hadaprod (\bC_2\bD_2) -\bA\right)\hadaprod (\bC_2\bD_2)\big)\bD_1^\top.
$$
\begin{proof}
For simplicity, we derive the gradient of $\bE$ for $f(\bE) = \normf{\bE\bF\hadaprod \bC - \bD}^2$.
We have 
$$
\begin{aligned}
f(\bE)&=\normf{\bE\bF\hadaprod \bC - \bD}^2 
= \trace\left((\bE\bF\hadaprod \bC - \bD)^\top(\bE\bF\hadaprod \bC - \bD) \right)\\
&= \trace\left((\bE\bF\hadaprod \bC)^\top (\bE\bF\hadaprod \bC)\right)-2\trace\left((\bE\bF\hadaprod \bC)^\top\bD\right) + \trace(\bD^\top\bD).
\end{aligned}
$$
Considering the first term, we get 
$$
\frac{\partial \trace\left((\bE\bF\hadaprod \bC)^\top (\bE\bF\hadaprod \bC)\right)}{\partial \bE}
=2(\bE\bF)\hadaprod \bC\hadaprod \bC \cdot \bF^\top.
~\footnote{Use the fact that $\frac{\partial \trace\left((\bE\hadaprod \bC)^\top(\bE\hadaprod \bC) \right)}{\partial \bE}=2\bE\hadaprod \bC\hadaprod \bC$, which can be derived element-wise.}
$$
For the second term, it follows that 
$$
-2\frac{\partial \trace\left((\bE\bF\hadaprod \bC)^\top\bD\right)}{\partial \bE}
=-2\bD\hadaprod \bC \cdot \frac{\partial \bE\bF}{\partial \bE}
=-2\bD\hadaprod \bC \cdot\bF^\top.
~\footnote{Use the fact that $\frac{\partial\trace( (\bE\hadaprod \bC)^\top\bD )}{\partial \bE} = \bD\hadaprod \bC$, which can be derived element-wise. Since $\trace( (\bE\hadaprod \bC)^\top\bD)=\sum_{i,j} d_{ij}a_{ij}c_{ij}$ and thus $\frac{\partial \trace( (\bE\hadaprod \bC)^\top\bD)}{\partial a_{ij}}=d_{ij}c_{ij}$.}
$$
The third term is a constant w.r.t. to $\bE$. 
Therefore, $\frac{\partial f(\bE)}{\partial \bE} = 2(\bE\bF)\hadaprod \bC\hadaprod \bC \cdot \bF^\top-2\bD\hadaprod \bC \cdot\bF^\top
=2\big((\bE\bF)\hadaprod \bC -\bD \big)\hadaprod \bC\cdot\bF^\top
.$
Substituting  $\bE\triangleq\bC_1$, $\bF\triangleq\bD_1$, $\bC\triangleq\bC_2\bD_2$, and $\bD\triangleq\bA$ completes the proof.
\end{proof}
The gradients with respect to $\bD_1, \bC_2$, and $\bD_2$ can be derived analogously.
Thus, the alternating descent method for obtaining the low-rank approximation of Hadamard decomposition can be described by  Algorithm~\ref{alg:ad_hadamad_svd}.

\begin{algorithm}[h] 
\caption{Alternating Descent with Gradient Descent for Low-Rank Hadamard Decomposition: A regularization can also be added into the gradient descent update (see Section~\ref{section:regularization-extention-general}).}
\label{alg:ad_hadamad_svd}
\begin{algorithmic}[1] 
\Require Matrix $\bA\in \real^{M\times N}$;
\State Initialize $\bC_1,\bC_2\in \real^{M\times K}$, and $\bD_1,\bD_2\in \real^{K\times N}$; 
\State Choose a stoping criterion on the approximation error $\delta$;
\State Choose  step size $\eta$;
\State Choose the maximum number of iterations $C$;
\State $iter=0$; \Comment{Count for the number of iterations}
\While{$\normf{(\bC_1\bD_1)\hadaprod (\bC_2\bD_2) -\bA}^2>\delta $ and $iter<C$}
\State $iter=iter+1$; 
\State $\Delta \leftarrow \left((\bC_1\bD_1)\hadaprod (\bC_2\bD_2) -\bA\right)$;
\State $\bC_1 \leftarrow \bC_1-\eta \nabla L(\bC_1)=\bC_1-\eta\cdot 2\left(\Delta \hadaprod (\bC_2\bD_2)\right)\bD_1^\top$;
\State $\bD_1 \leftarrow \bD_1-\eta \nabla L(\bD_1)=\bD_1-\eta\cdot 2 \left\{\left(\Delta^\top \hadaprod (\bC_2\bD_2)^\top\right)\bC_1\right\}^\top$;
\State $\bC_2 \leftarrow \bC_2-\eta \nabla L(\bC_2)=\bC_2 - \eta \cdot 2 \left(\Delta \hadaprod (\bC_1\bD_1)\right)\bD_2^\top$;
\State $\bD_2 \leftarrow \bD_2-\eta \nabla L(\bD_2)=\bD_2-\eta \cdot 2 \left\{\left( \Delta^\top \hadaprod (\bC_1\bD_1)^\top \right)\bC_2\right\}^\top$;

\EndWhile
\State Output $\bC_1,\bD_1, \bC_2,\bD_2$;
\end{algorithmic} 
\end{algorithm}

\subsection{Rank-One Update }
Following the rank-one update approach in  ALS (Section~\ref{section:alt-columb-by-column}), we consider updating the $n$-th column $\bd_{1,n}$ of $\bD_1$, $n\in\{1,2,\ldots, N\}$.
Analogously, we can obtain the gradient of $\bd_{1,n}$:
\begin{equation}\label{equation:hada_rkone1}
\begin{aligned}
\nabla L(\bd_{1,n}) 
=\frac{\partial L(\bd_{1,n})}{\partial \bd_{1,n}}
&= 2\bC_1^\top \left((\bC_1\bd_{1,n}) \hadaprod \ba_{2,n} \hadaprod \ba_{2,n}\right) - 2\bC_1^\top (\ba_n \hadaprod \ba_{2,n})\\
&=2\bC_1^\top \left(\left[(\bC_1\bd_{1,n}) \hadaprod \ba_{2,n}-\ba_n \right] \hadaprod\ba_{2,n}\right), \gap n\in\{1,2,\ldots, N\}, 
\end{aligned}
\end{equation}
where $\ba_{2,n}$ is the $n$-th column of $\bA_2\triangleq\bC_2\bD_2$, $n\in\{1,2,\ldots, N\}$. The gradient of the columns of $\bD_2$ can be calculated similarly.

Suppose further that $\bC_1^\top=[\bc_{1,1}, \bc_{1,2}, \ldots, \bc_{1,M}]\in\real^{K\times M}$,  $\bB\triangleq\bA^\top=[\bb_1, \bb_2, \ldots, \bb_M]\in\real^{N\times M}$, and $\bB_2\triangleq\bA_2^\top\triangleq(\bC_2\bD_2)^\top=[\bb_{2,1}, \bb_{2,2}, \ldots, \bb_{2,M}]\in\real^{N\times M}$, i.e., the row partitions of $\bC_1$, $\bA$, and $\bA_2=(\bC_2\bD_2)$, respectively.
Then, the gradient of $\bc_{1,m}$ is 
\begin{equation}\label{equation:hada_rkone2}
\nabla L(\bc_{1,m}) 
=\frac{\partial L(\bc_{1,m})}{\partial \bc_{1,m}}
= 2\bD_1 \left([(\bD_1^\top\bc_{1,m}) \hadaprod \bb_{2,m}-\bb_m ] \hadaprod\bb_{2,m}\right), \, m\in\{1,2,\ldots, M\}.
\end{equation}
The gradient of the rows of $\bC_2$ can be obtained analogously.
Therefore,  Algorithm~\ref{alg:ad_hadamad_svd} can be modified to update the columns of $\bD_1, \bD_2$ and the rows of $\bC_1, \bC_2$ iteratively (referred to as  rank-one updates).

\subsection{Missing Entries}
The rank-one update can be extended to the Netflix context, in which case many entries of $\bA$ are missing. 
Assuming $\bA$ is a low-rank matrix, we aim to fill in the missing entries of matrix $\bA$ (where $M$ represents the number of movies, and $N$ represents the number of users).

Let $\bo_n\in \{0,1\}^M, n\in\{1,2,\ldots, N\}$, represent the movies rated by user $n$, where $o_{nm}=1$ if user $n$ has rated movie $m$, and $o_{nm}=0$ otherwise.
Similarly, let $\bp_m \in\{0,1\}^{N}, m\in\{1,2,\ldots,M\}$ denote the users who have rated  movie $m$, with $p_{mn}=1$ if the movie $m$ has been rated by user $n$, and $p_{mn}=0$ otherwise.
Then, Equations~\eqref{equation:hada_rkone1} and~\eqref{equation:hada_rkone2} become
\begin{align}
\nabla L(\bd_{1,n}) 
&=2\bC_1[\bo_n,:]^\top \left(\left[(\bC_1[\bo_n,:]\bd_{1,n}) \hadaprod \ba_{2,n}[\bo_n]-\ba_n[\bo_n] \right] \hadaprod\ba_{2,n}[\bo_n]\right),\nonumber \\
&\gap\gap\gap\gap\gap\gap\gap\gap\gap\gap\gap\gap\gap\gap\gap\gap n\in\{1,2,\ldots, N\};\label{equation:hada_rkone3} \\
\nabla L(\bc_{1,m}) 
&= 2\bD_1[:,\bp_m] \left(\left[(\bD_1[:,\bp_m]^\top\bc_{1,m}) \hadaprod \bb_{2,m}[\bp_m]-\bb_m[\bp_m] \right] \hadaprod\bb_{2,m}[\bp_m]\right),\nonumber\\
&\gap\gap\gap\gap\gap\gap\gap\gap\gap\gap\gap\gap\gap\gap\gap\gap  m\in\{1,2,\ldots, M\} \label{equation:hada_rkone4}.
\end{align}
Since the Hadamard product commutes,  the gradient of $L(\bd_{2,n}), n\{1,2,\ldots, N\}$ and $\bc_{2,m},  m\in\{1,2,\ldots, M\}$ can be obtained similarly due to symmetry. 
The process for predicting the missing entries in $\bA$ is then formulated in Algorithm~\ref{alg:ad_hadamad_missen}.

\begin{algorithm}[h] 
\caption{Alternating Descent with Gradient Descent for Hadamard Decomposition with Missing Entries: A regularization can also be added into the gradient descent update (see Section~\ref{section:regularization-extention-general}).}
\label{alg:ad_hadamad_missen}
\begin{algorithmic}[1] 
\Require Matrix $\bA\in \real^{M\times N}$;
\State Initialize $\bC_1,\bC_2\in \real^{M\times K}$, and $\bD_1,\bD_2\in \real^{K\times N}$; 
\State Choose a stoping criterion on the approximation error $\delta$;
\State Choose  step size $\eta$;
\State Choose the maximum number of iterations $C$;
\State $iter=0$; \Comment{Count for the number of iterations}
\While{$\normf{(\bC_1\bD_1)\hadaprod (\bC_2\bD_2) -\bA}^2>\delta $ and $iter<C$}
\State $iter\leftarrow iter+1$; 
\For{$n=1,2,\ldots, N$}
\State $\bd_{1,n}\leftarrow\bd_{1,n}-\eta \nabla L(\bd_{1,n}) $;  \Comment{Equation~\eqref{equation:hada_rkone3}}
\State $\bd_{2,n}\leftarrow\bd_{2,n}-\eta \nabla L(\bd_{2,n}) $;
\EndFor
\For{$m=1,2,\ldots, M$}
\State $\bc_{1,m}\leftarrow\bc_{1,m}-\eta \nabla L(\bc_{1,m})  $; \Comment{Equation~\eqref{equation:hada_rkone4}}
\State $\bc_{2,m}\leftarrow\bc_{2,m}-\eta \nabla L(\bc_{2,m})  $;
\EndFor
\EndWhile
\State Output $\bC_1,\bD_1, \bC_2,\bD_2$;
\end{algorithmic} 
\end{algorithm}

\subsection*{Example}
\begin{figure}[h]
\centering  
\vspace{-0.35cm} 
\subfigtopskip=2pt 
\subfigbottomskip=2pt
\subfigcapskip=-5pt
\subfigure[Training.]{\label{fig:movie100k1_hada}
\includegraphics[width=0.301\linewidth]{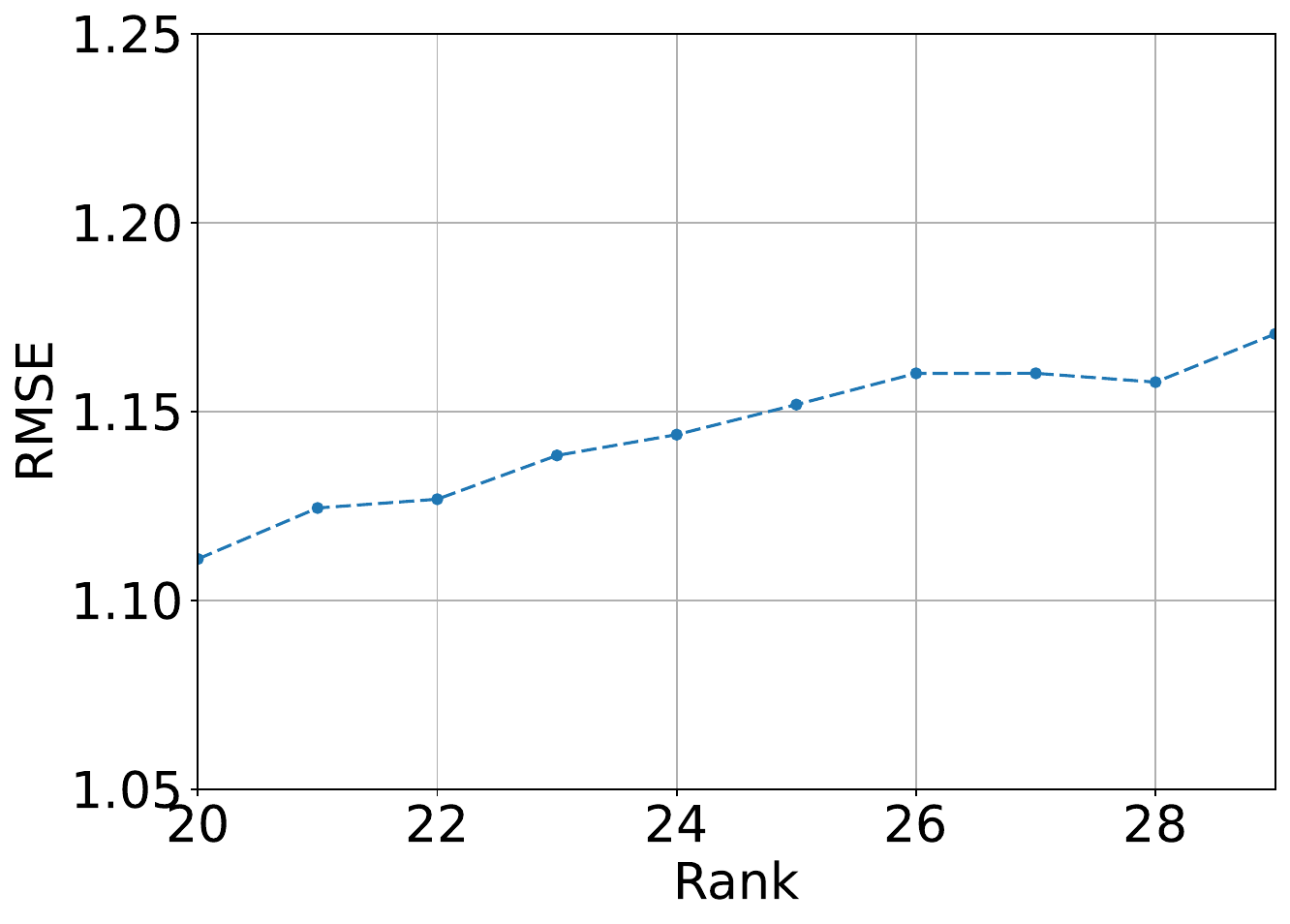}}
\subfigure[Validation.]{\label{fig:movie100k_hada2}
\includegraphics[width=0.301\linewidth]{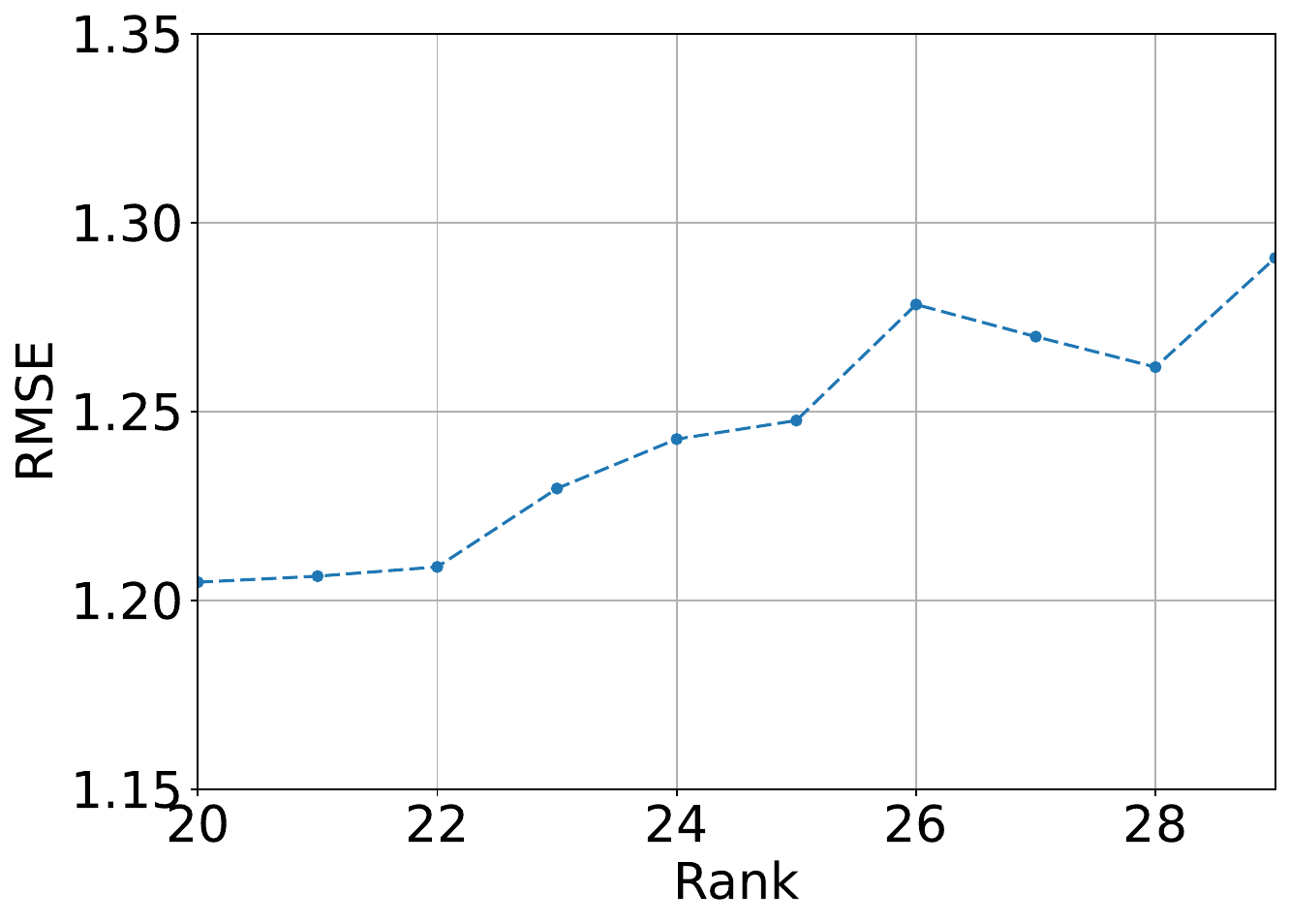}}
\subfigure[Iterations.]{\label{fig:movie100k_hada3}
\includegraphics[width=0.301\linewidth]{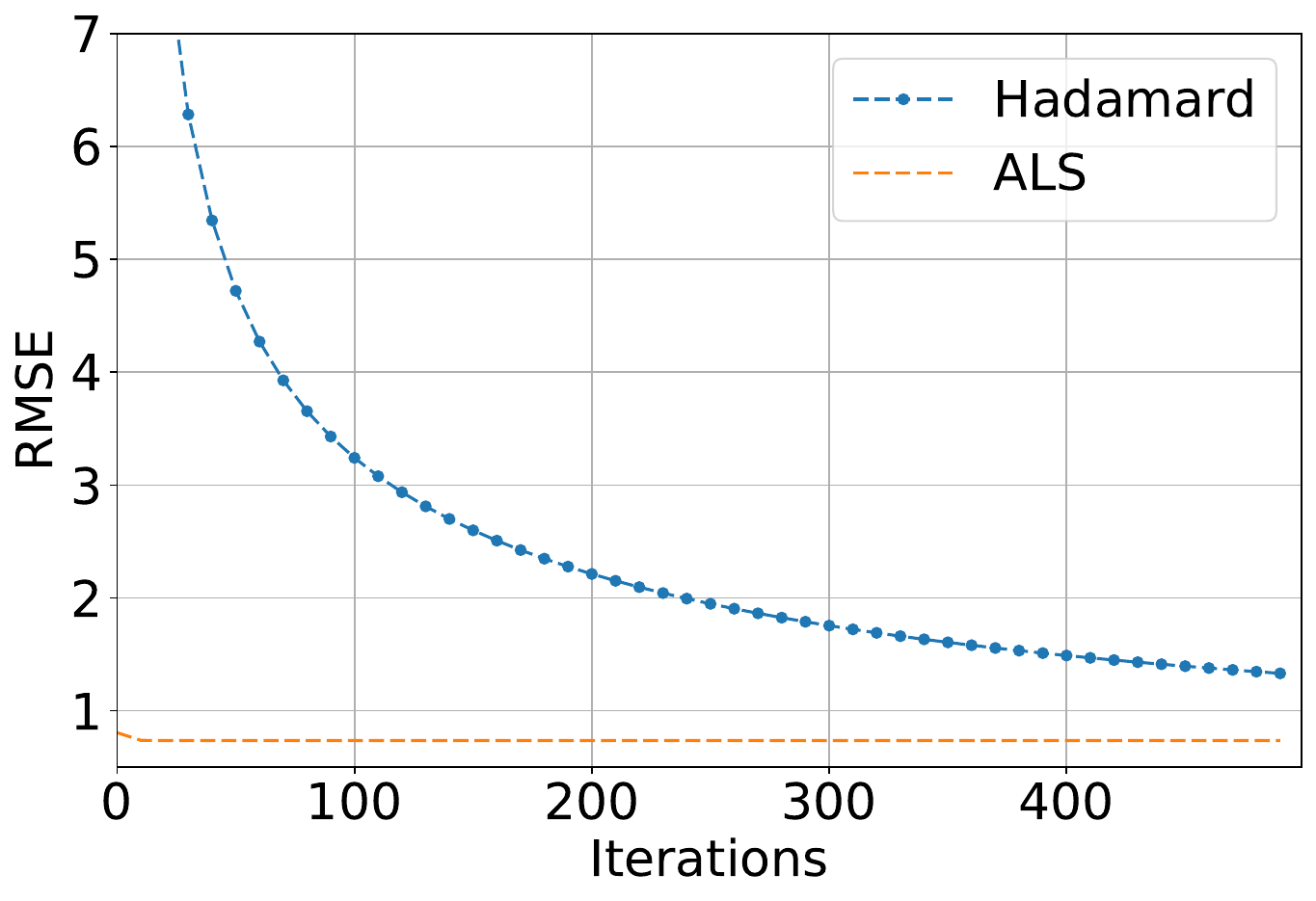}}
\caption{Comparison of training and validation error for the ``MovieLens 100K" data set with different reduction dimensions $K$ in low-rank Hadamard decomposition.}
\label{fig:movie100k_hada}
\end{figure}

We utilize the ``MovieLens 100K" dataset, as introduced in Section~\ref{section:movie_rec_als}, to investigate the impact of the low-rank Hadamard decomposition on the reconstruction of missing entries.
As a reminder,  the ALS method achieved  an RMSE of 0.806 on the  validation set (Figure~\ref{fig:movie100k}); and the validation error generally decreased  as the rank dimension increased.
Figure~\ref{fig:movie100k_hada} displays  the training and validation errors for various values of the reduction dimension $K$ in the low-rank Hadamard decomposition.
Notably, both the training and validation errors  consistently remain above 1.
In this particular scenario, increasing the dimension $K$ does not result in a reduction of error.
This behavior can be attributed to the relatively small size of the data set; we anticipate better performance with larger data sets. Additionally, since the objective function in the low-rank Hadamard decomposition is not jointly convex and is more complex than the problem addressed by ALS, the gradient descent algorithm may easily become trapped in local minima.

Another challenge observed during the experiments was the tendency of the low-rank Hadamard decomposition to diverge, which requires a very small step size (learning rate) of $1\times 10^{-5}$. Figure~\ref{fig:movie100k_hada3} compares the number of iterations needed for convergence between ALS and the low-rank Hadamard decomposition, with the dimension $K$ set to 20. While ALS converges in fewer than 20 iterations, the low-rank Hadamard decomposition needs over 500 iterations to reach convergence. Accelerated update methods, such as the \textit{Barzilai and Borwein gradient method, FISTA, and Nesterov's accelerated gradient method} \citep{barzilai1988two, nesterov2005smooth, beck2009fast}, could potentially expedite the convergence. However, a detailed discussion of these advanced optimization techniques is outside the scope of this book, and we will not discuss them in detail.


\section{Low-Rank Kronecker Decomposition}\label{section:low_rank_kronecker}

\subsection{Low-Rank Kronecker Decomposition}
Similar to the low-rank Hadamard decomposition, we can also consider the low-rank Kronecker decomposition.
Suppose the design matrix $\bA\in\real^{M\times N}$ has dimensions  $M=m_1m_2$ and $N=n_1n_2$. 
Let us consider two matrices $\bB\in\real^{m_1\times n_1}$ and $\bC\in\real^{m_2\times n_2}$ such that $\bB\kronecker \bC\in\real^{M\times N}$.
Theorem~\ref{theorem:rank_kronec_prod} states that $\rank(\bB\kronecker \bC)=\rank(\bB)\rank(\bC)$.
The goal then becomes 
\begin{itemize}
\item Given a real matrix $\bA\in \real^{M\times N}$, find  matrix factors $\bB\in\real^{m_1\times n_1}$ and $\bC\in\real^{m_2\times n_2}$ such that: 
$$
\min\,\,L(\bB, \bC) = \normf{\bB\kronecker \bC-\bA}^2.
$$
\end{itemize}

Given the definition of the Kronecker product $\bB\kronecker \bC$ (Definition~\ref{definition:kronecker-product}), we may consider the \textit{uniform blocking} of matrix $\bA$:
$$
\bA=
\begin{bmatrix}
\bA_{11} &\bA_{12} & \ldots &\bA_{1,n_1} \\
\bA_{21} &\bA_{22} & \ldots &\bA_{2,n_1} \\
\vdots & \vdots & \ddots & \vdots \\
\bA_{m_1,1} &\bA_{m_1,2} & \ldots &\bA_{m_1,n_1} \\
\end{bmatrix}, 
\gap 
\bA_{ij}\in\real^{m_2\times n_2},
$$
where the $(i,j)$-th block is an $m_2\times n_2$ matrix with $\bA_{ij}=\bA[(i-1)m_2+1:i\cdot m_2, (j-1)n_2+1:j\cdot n_2]\in\real^{m_2\times n_2}$, for all  $i\in\{1,2,\ldots,m_1\}$ and $j\in\{1,2,\ldots,n_1\}$.
While the $(i,j)$-th block of $\bB\kronecker \bC$ is $b_{ij}\bC$.
Therefore, when $\bC$ is held constant,  the objective function with respect to $\bB$ (or $b_{ij}$ for $i\in\{1,2,\ldots,m_1\}, j\in\{1,2,\ldots,n_1\}$) is 
$$
L(\bB) = \sum_{i=1}^{m_1}\sum_{j=1}^{n_1} \normf{\bA_{ij} - b_{ij}\bC}^2.
$$
Analogously, when  $\bB$ is kept fixed, the objective function with respect to $\bC$ is 
$$
L(\bC) =\sum_{i=1}^{m_2}\sum_{j=1}^{n_2} \normf{\widetildebA_{ij}- c_{ij}\bB}^2,
$$
where $\widetildebA_{ij}=\bA[i:m_2:M, j:n_2:N]\in\real^{m_1\times n_1}$, i.e., slicing the rows every $m_2$ indices and the columns every $n_2$ indices; the row indices of $\widetildebA_{ij}$ are $i, i+m_2, i+2m_2, \ldots$, and the column indices are $j, j+n_2, j+2n_2,\ldots$.

\subsection{Kronecker Decomposition via Alternating Least Squares}
Thinking at the block level for matrices can enable us to find  least squares solutions separately.
When dealing with large matrices, particularly those that can be naturally partitioned into smaller blocks, it is often beneficial to think at the block level. This approach allows us to find least squares solutions separately for each block independently, which can significantly simplify the overall optimization process.
\begin{theorem}[Kronecker Decomposition via Least Squares \citep{van1993approximation}]\label{theorem:lrank_kro_opt}
Let $\bA\in \real^{M\times N}$ be any matrix  with $M=m_1m_2$ and $N=n_1n_2$. If $\bC\in\real^{m_2\times n_2}$ is fixed, then the matrix $\bB$ defined by 
\begin{equation}\label{equation:kro_ls_1}
b_{ij} =\frac{\trace(\bA_{ij}^\top \bC)}{\trace(\bC^\top\bC)}, \gap i\in\{1,2,\ldots, m_1\}, j\in\{1,2,\ldots, n_1\}
\footnote{The trace operator takes into account only the diagonal elements of a matrix. Therefore, in practice, it's more appropriate to compute the Hadamard product of the two matrices and then sum all the squared elements.}
\end{equation}
minimizes $\normf{\bB\kronecker \bC-\bA}$, where $\bA_{ij}=\bA[(i-1)m_2+1:i\cdot m_2, (j-1)n_2+1:j\cdot n_2]$ is $(i,j)$-th block of $\bA$ with size $m_2\times n_2$.
Analogously, if $\bB\in\real^{m_1\times n_1}$ is fixed, then the matrix $\bC$ defined by 
\begin{equation}\label{equation:kro_ls_2}
c_{ij} =\frac{\trace(\widetilde{\bA}_{ij}^\top\bB)}{\trace(\bB^\top\bB)}, \gap i\in\{1,2,\ldots, m_2\}, j\in\{1,2,\ldots, n_2\}
\end{equation}
minimizes $\normf{\bB\kronecker \bC-\bA}$, where $\widetilde{\bA}_{ij}=\bA[i:m_2:M, j:n_2:N]\in\real^{m_1\times n_1}$.
\end{theorem}
\begin{proof}[of Theorem~\ref{theorem:lrank_kro_opt}]
Consider the  loss of the $(i,j)$-th block:
$$
\normf{\bA_{ij} - b_{ij}\bC}^2 = \trace\left((\bA_{ij} - b_{ij}\bC)^\top(\bA_{ij} - b_{ij}\bC)\right)
=\normf{\bA_{ij}}^2 - 2b_{ij}\trace(\bC^\top\bA_{ij}) + b_{ij}^2 \normf{\bC}^2.
$$
Therefore, it follows that 
$$
\frac{\partial L(\bB)}{\partial b_{ij}} 
= -2\trace(\bC^\top\bA_{ij})+2b_{ij}\normf{\bC}^2.
$$
Setting the partial derivatives to zero yields the result.
The second part follows similarly.
\end{proof}
This approach leverages the fact that the least squares solution can be found separately for each block, which simplifies the computation and can lead to faster convergence. It is particularly useful in scenarios where the matrix $\bA$ exhibits a structured pattern that can be exploited through block-wise processing.
The lemma suggests an alternating descent update for obtaining the decomposition, which alternatively improve $\bB$ and $\bC$.
The low-rank Kronecker decomposition has a closed-form for each update, and the procedure is formulated in Algorithm~\ref{alg:ad_Kronecker_zerograd}.

\begin{algorithm}[h] 
\caption{Alternating Descent for Low-Rank Kronecker Decomposition}
\label{alg:ad_Kronecker_zerograd}
\begin{algorithmic}[1] 
\Require Matrix $\bA\in \real^{M\times N}$ with $M=m_1m_2$ and $N=n_1n_2$;
\State Initialize $\bB\in\real^{m_1\times n_1}$ and $\bC\in\real^{m_2\times n_2}$; 
\State Choose a stoping criterion on the approximation error $\delta$;
\State Choose the maximum number of iterations $C$;
\State $iter=0$; \Comment{Count for the number of iterations}
\While{$\normf{\bB\kronecker \bC -\bA}^2>\delta $ and $iter<C$}
\State $iter=iter+1$; 
\State $C_1 \leftarrow \trace(\bC^\top\bC)$
\For{$i=1,2,\ldots,m_1, j=1,2,\ldots,n_1$}
\State $b_{ij} \leftarrow \frac{\trace(\bA_{ij}^\top \bC)}{C_1}$;
\EndFor
\State $C_2 \leftarrow \trace(\bB^\top\bB)$
\For{$i=1,2,\ldots,m_2, j=1,2,\ldots,n_2$}
\State $c_{ij} \leftarrow \frac{\trace(\widetilde{\bA}_{ij}^\top\bB)}{C_2}$;
\EndFor
\EndWhile
\State Output $\bB, \bC$;
\end{algorithmic} 
\end{algorithm}

\subsection{Missing Entries}
The least squares framework for low-rank Kronecker decomposition can be readily applied to the ``Netflix" context.
We introduce an additional mask matrix $\bM\in\{0,1\}^{M\times N}$, where $m_{mn}\in\{0,1\}$ indicates whether user $n$ has rated movie $m$ or not. 
Consequently, the update of $b_{ij}$ in Equation~\eqref{equation:kro_ls_1} can be obtained by 
$$
b_{ij} =\frac{\normf{\bA_{ij} \hadaprod \bC \hadaprod \bM_{ij}}^2}{\normf{\bC \hadaprod \bC \hadaprod \bM_{ij}}^2}, \gap i\in\{1,2,\ldots, m_1\}, j\in\{1,2,\ldots, n_1\},
$$
where $\bM_{ij}=\bM[(i-1)m_2+1:i\cdot m_2, (j-1)n_2+1:j\cdot n_2]\in\real^{m_2\times n_2}$.
And the update of $c_{ij}$ in Equation~\eqref{equation:kro_ls_2} can be derived as:
$$
c_{ij} =\frac{ \normf{\widetilde{\bA}_{ij}\hadaprod\bB\hadaprod \widetildebM_{ij}}^2 }{\normf{\bB\hadaprod \bB \hadaprod \widetildebM_{ij}}^2}, \gap i\in\{1,2,\ldots, m_2\}, j\in\{1,2,\ldots, n_2\},
$$
where $\widetildebM_{ij}=\bM[i:m_2:M, j:n_2:N]\in\real^{m_1\times n_1}$.

\section{Low-Rank Khatri-Rao Decomposition}\label{section:lrank_khatri_decom}
\subsection{Low-Rank Khatri-Rao Decomposition}

As  previously mentioned, in models where interactions between different sets of features are crucial, the Khatri-Rao product can be employed to represent these interactions in a compact and meaningful way. This is especially relevant in fields such as genomics, where the interactions between genetic markers are a subject of study.

Similar to the low-rank Hadamard decomposition, we may also consider the low-rank Kronecker decomposition.
Suppose the design matrix $\bA\in\real^{M\times N}$ has dimensions $M=m_1m_2$. 
Consider two matrices $\bB\in\real^{m_1\times N}$ and $\bC\in\real^{m_2\times N}$ such that their Khatri-Rao product is $\bB\khatrirao  \bC\in\real^{M\times N}$.
Theorems~\ref{theorem:rank_khatri_prod} and \ref{theorem:kkrank_khatri_prod} demonstrate that $\rank(\bB\khatrirao  \bC)\geq \max\{\rank(\bB), \rank(\bC)\}$ and $\rank_k(\bB\khatrirao  \bC)\geq \min\{\rank_k(\bB)+\rank_k(\bC)-1, N\}$. In this sense, the expressive power of the Khatri-Rao product is at least as strong as the maximum rank of the individual factors, and it is more powerful than the sum of the Kruskal ranks of the two factors minus one, up to the dimension 
$N$.
The goal then becomes 
\begin{itemize}
	\item Given a real matrix $\bA\in \real^{M\times N}$, find  matrix factors $\bB\in\real^{m_1\times N}$ and $\bC\in\real^{m_2\times N}$ such that: 
	$$
	\min\,\,L(\bB, \bC) = \normf{\bB\khatrirao  \bC-\bA}^2.
	$$
\end{itemize}

\subsection{Khatri-Rao Decomposition via Alternating Least Squares}
Following Section~\ref{section:project-onto-a-vector} about the projection of a vector onto another vector, the alternating least squares algorithm can be obtained as follows.
\begin{theorem}[Khatri-Rao Decomposition via Least Squares]\label{theorem:lrank_khatri_opt}
Let $\bA\in \real^{M\times N}$ be any matrix with $M=m_1m_2$. If $\bC\in\real^{m_2\times N}$ is fixed, then the matrix $\bB\in\real^{m_1\times N}$ defined by 
\begin{equation}\label{equation:khatri_opt_1}
b_{ij} =\frac{\bc_j^\top\ba_{ij}}{\bc_j^\top \bc_j}, \gap i\in\{1,2,\ldots, m_1\}, j\in\{1,2,\ldots, N\}
\end{equation}
minimizes $\normf{\bB\khatrirao  \bC-\bA}$, where $\ba_{ij}=\bA[(i-1)m_2+1:i\cdot m_2, j]\in\real^{m_2\times 1}$ is the $(i,j)$-th block of $\bA$ with size $m_2\times 1$, and $\bc_j$ is the $j$-th column of $\bC$.
Analogously, if $\bB\in\real^{m_1\times N}$ is fixed, then the matrix $\bC\in\real^{m_2\times N}$ defined by 
\begin{equation}\label{equation:khatri_opt_2}
c_{ij} =\frac{\bb_j^\top\widetilde{\ba}_{ij} }{\bb_j^\top\bb_j}, \gap i\in\{1,2,\ldots, m_2\}, j\in\{1,2,\ldots, N\}
\end{equation}
minimizes $\normf{\bB\khatrirao \bC-\bA}$, where $\widetilde{\ba}_{ij}=\bA[i:m_2:M, j]\in\real^{m_1\times 1}$, and $\bb_j$ is the $j$-th column of $\bB$.
\end{theorem}
\begin{proof}[of Theorem~\ref{theorem:lrank_khatri_opt}]
Consider the  loss for the $(i,j)$-th block:
$$
\normtwo{\ba_{ij} - b_{ij}\bc_j}^2 = b_{ij}^2\bc_j^\top\bc_j - 2\bc_j^\top\ba_{ij} b_{ij} + \ba_{ij}^\top\ba_{ij}.
$$
Therefore, $b_{ij} =\frac{\bc_j^\top\ba_{ij}}{\bc_j^\top \bc_j}$ obtains the minimum.
The second part follows similarly.
\end{proof}

The Khatri-Rao decomposition is highly flexible, allowing for the addition of more components in the Khatri-Rao product. This flexibility enables us to extend the basic Khatri-Rao decomposition to incorporate multiple factors, which can be useful in various applications where the data exhibit complex interactions, making it suitable for scenarios where multiple sources of variation need to be accounted for. The \textit{cascaded Khatri-Rao decomposition} can be formulated as follows:
\begin{itemize}
\item Given a real matrix $\bA\in\real^{M\times N}$,  the goal is to find matrix factors $\bB\in\real^{m_1\times N}$, $\bW\in\real^{n\times N}$, $\bC\in\real^{m_2\times N}$, and potentially additional factors, such that: 
$$
\min\,\,L(\bB, \bW, \bC, \ldots) = \normf{\bB\khatrirao  \bW\khatrirao \bC-\bA}^2.
$$
\end{itemize}
\begin{theorem}[Cascaded Khatri-Rao Decomposition via Least Squares]\label{theorem:casc_lrank_khatri_opt}
Let $\bA\in \real^{M\times N}$ be any matrix with $M=m_1\cdot\textcolor{mylightbluetext}{n}\cdot m_2$. If $\textcolor{mylightbluetext}{\bW\in\real^{n\times N}}$ and $\bC\in\real^{m_2\times N}$ are fixed, then the matrix $\bB\in\real^{m_1\times N}$ defined by 
\begin{equation}\label{equation:cas_khatri_opt_1}
b_{ij} =\frac{\textcolor{mylightbluetext}{\widetildebc_j}^\top\ba_{ij}}{\textcolor{mylightbluetext}{\widetildebc_j}^\top \textcolor{mylightbluetext}{\widetildebc_j}}, \gap i\in\{1,2,\ldots, m_1\}, j\in\{1,2,\ldots, N\}
\end{equation}
minimizes $\normf{\bB\khatrirao \bW\khatrirao  \bC-\bA}$, where $\ba_{ij}=\bA[(i-1)\textcolor{mylightbluetext}{\cdot (n\cdot m_2)}+1:i \textcolor{mylightbluetext}{\cdot (n\cdot m_2)}, j]$ is the $(i,j)$-th block of $\bA$ with size $\textcolor{mylightbluetext}{ (n\cdot m_2)}\times 1$, and $\textcolor{mylightbluetext}{\widetildebc_j}$ is the $j$-th column of $\bW\khatrirao \bC$.

Analogously, if $\textcolor{mylightbluetext}{\bW\in\real^{n\times N}}$ and $\bB\in\real^{m_1\times N}$ are fixed, then the matrix $\bC$ defined by 
\begin{equation}\label{equation:cas_khatri_opt_2}
c_{ij} =\frac{\textcolor{mylightbluetext}{\widetildebb_j}^\top\widetilde{\ba}_{ij} }{\textcolor{mylightbluetext}{\widetildebb_j}^\top\textcolor{mylightbluetext}{\widetildebb_j}}, \gap i\in\{1,2,\ldots, m_2\}, j\in\{1,2,\ldots, N\}
\end{equation}
minimizes $\normf{\bB\khatrirao\bW\khatrirao \bC-\bA}$, where $\widetilde{\ba}_{ij}=\bA[i:\textcolor{black}{ m_2}:M, j]\in\real^{\textcolor{mylightbluetext}{n\cdot m_1}\times 1}$, and $\textcolor{mylightbluetext}{\widetildebb_j}$ is the $j$-th column of $\bB\khatrirao \bW$.

\paragraph{Obtain $\bW$. }
Up to this point,  the result is similar to Theorem~\ref{theorem:lrank_khatri_opt} with careful attention to the dimensions.
If $\bB\in\real^{m_1\times N}$ and $\bC\in\real^{m_2\times N}$ are fixed, then the matrix $\bW$ defined by 
\begin{equation}\label{equation:cas_khatri_opt_3}
w_{ij} 	=\frac{\widetildebw_{ij}^\top\widehat{\ba}_{ij} }{\widetildebw_{ij}^\top\widetildebw_{ij}}, \gap i\in\{1,2,\ldots, n\}, j\in\{1,2,\ldots, N\}
\end{equation}	
minimizes $\normf{\bB\khatrirao\bW\khatrirao \bC-\bA}$, where
$\widehat{\ba}_{ij} =\bA[I_s, j] \in\real^{m_1\cdot m_2\times 1}$,  
$$
\widetildebw_{ij}=
[b_{1j}  \bc_j, b_{2j}  \bc_j, \dots b_{m_1,j}  \bc_j]^\top \in\real^{m_1\cdot m_2\times 1},
$$
$[b_{1j}  , b_{2j}  , \dots b_{m_1,j}  ]^\top $ is the $j$-th column of $\bB$, and $\bc_j$ is the $j$-th column of $\bC$.
The  index set $I_s$ contains $m_1$ blocks, each with $m_2$ components:
$$
\footnotesize
\begin{aligned}
 I_s = \{(j:j+m_2-1), (j+n:j+n+m_2-1), \ldots  (j+(m_1-1)n:j+(m_1-1)n+m_2-1)\}, 
\end{aligned}
$$
\end{theorem}
\begin{proof}[of Theorem~\ref{theorem:casc_lrank_khatri_opt}]
The $j$-th column of $\bB\khatrirao\bW\khatrirao \bC$ is
$$
\bb_i \kronecker \bw_i \kronecker \bc_i
= 
[
(b_{1j} w_{1j} \bc_j, b_{1j} w_{2j} \bc_j, \dots b_{1j} w_{nj} \bc_j)\mid 
(b_{2j} w_{1j} \bc_j, b_{2j} w_{2j} \bc_j, \dots b_{2j} w_{nj} \bc_j)\mid
\ldots 
]^\top.
$$
The result follows by solving  least squares problem for each $w_{ij}$ for $i\in\{1,2,\ldots, n\}, j\in\{1,2,\ldots, N\}$.
\end{proof}

With this proof, it is straightforward to extend the cascaded Khatri-Rao decomposition to more than three matrices: : $\bA=\bA_1\khatrirao \bA_2\khatrirao \bA_3\khatrirao \bA_4\khatrirao\ldots$.
The update for $\bA_2$ can be obtained by setting $\bA_1\triangleq\bB$,  $\bA_2\triangleq\bW$, $\bA_3\khatrirao \bA_4\khatrirao\ldots\triangleq\bC$ in Theorem~\ref{theorem:casc_lrank_khatri_opt}. 
Similarly,  the update for $\bA_3$ can be obtained by setting $\bA_1\khatrirao \bA_2\triangleq\bB$,  $\bA_3\triangleq\bW$, $ \bA_4\khatrirao \ldots\triangleq\bC$ in Theorem~\ref{theorem:casc_lrank_khatri_opt}.

\section{Low-Rank Adaptation in Large Language Models}

Large language models (LLMs) have become a cornerstone of modern natural language processing (NLP) and have driven significant advancements in the field. 
This section provides a brief look at the background of large language models, including their development, key milestones, and their impact on NLP.
The concept of using neural networks for language modeling dates back to the early 2000s. However, the rise of LLMs really began with the advent of deep learning and the transformer architecture. Prior to this, recurrent neural networks (RNNs) and long short-term memory (LSTM) networks were the dominant architectures for sequence modeling.

The transformer architecture, introduced in 2017 by \citet{vaswani2017attention}, marked a significant departure from RNNs. Transformers use self-attention mechanisms to process input sequences in parallel, making them highly efficient for large-scale training. This architecture became the foundation  for many state-of-the-art models.

The pre-training paradigm involves training a model on a large corpus of text to learn general linguistic patterns before fine-tuning it for specific tasks. This approach was popularized by BERT (bidirectional encoder representations from transformers) in 2018 \citep{devlin2018bert}, which achieved state-of-the-art results on a variety of NLP benchmarks. BERT and similar models leverage the transformer architecture to capture contextual representations of words.

Following the success of BERT, researchers began to scale up the size of pre-trained models. GPT-2 (generative pre-trained transformer 2) and RoBERTa (robustly optimized BERT pre-training approach) were among the first models to demonstrate that increasing model size could lead to improved performance. GPT-3, released in 2020, took this trend to the extreme, with 175 billion parameters \citep{brown2020language}, and showed remarkable abilities in few-shot learning, where the model can perform tasks with minimal additional training.

Recent developments in large language models have focused on improving efficiency and scalability.
As models grew larger, the computational and storage costs of fine-tuning them became prohibitive. This led to the development of parameter-efficient tuning (PET) techniques, which aim to adapt pre-trained models to new tasks with minimal changes to the model's parameters. Techniques like LoRA (low-rank adaptation) and adapters are examples of PET methods \citep{hu2021lora, houlsby2019parameter}.

Adapters  are small neural networks inserted between layers of the pre-trained model, which are fine-tuned for specific tasks while keeping the rest of the model frozen. Adapters have been shown to be effective in transferring knowledge from pre-training to downstream tasks.

\subsection{LoRA and Other Tales}
LoRA  is a technique that freezes the pre-trained model weights and injects trainable rank decomposition matrices into each layer of the transformer architecture. LoRA significantly reduces the number of trainable parameters and the GPU memory requirement, making it more efficient than full fine-tuning.
LoRA does not introduce any additional latency during inference. After fine-tuning, the trainable matrices can be merged with the frozen weights, eliminating the need for extra computations during inference.

Full fine-tuning, which involves retraining all the parameters of a pre-trained model, becomes impractical as the size of these models grows. For example, GPT-3 175B has 175 billion parameters , and deploying multiple instances of such a model, each fine-tuned for a different task, is prohibitively expensive in terms of both computational resources and storage.
Moreover, storing an entire model checkpoint for each downstream application makes deployment and switching between different tasks extremely inefficient. This is particularly problematic when dealing with large models.

LoRA reduces the number of trainable parameters significantly, which in turn reduces the computational and memory overhead during training and inference. For instance, compared to GPT-3 175B  fine-tuned with Adam, LoRA can reduce the number of trainable parameters by 10,000 times and the GPU memory requirement by 3 times.

\begin{figure}[h]
\centering
\includegraphics[width=1.02\textwidth]{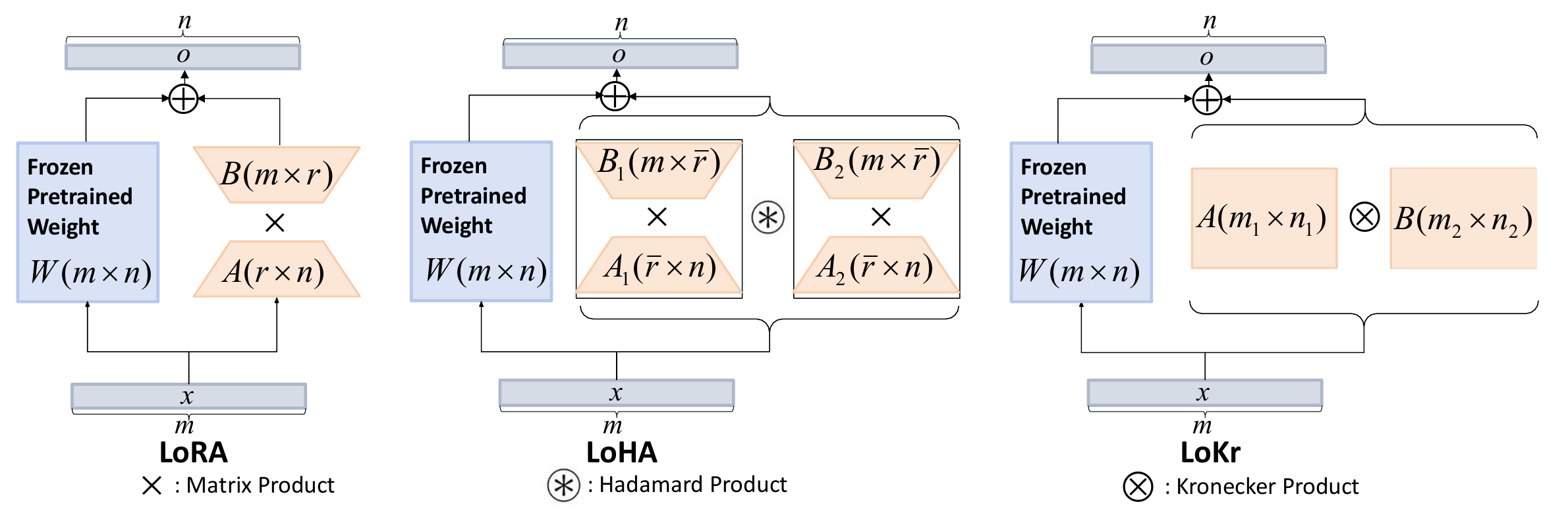}
\caption{Diagram illustrating LoRA, LoHA, and LoKr (KronA).}
\label{fig:lora_loha_lokr}
\end{figure}

LoRA also enables task-switching flexibility.
A pre-trained model can be shared and used to build many small LoRA modules for different tasks. By freezing the shared model and efficiently switching tasks by replacing the small LoRA matrices (the matrices $\bA$ and $\bB$ in Figure~\ref{fig:lora_loha_lokr}), the storage requirement and task-switching overhead are significantly reduced.

\index{LoRA}
\index{LoHA}
\index{LoKr (KronA)}
\paragraph{LoRA.}  To be more specific, the weight update $\Delta \bW \in \real^{m\times n}$ is decomposed into two low-rank matrices $\bB \in\real^{m\times r}$ and $\bA\in\real^{r\times n}$, where $m$ and $n$ are the dimensions of the original model parameters, and $r$, the rank of the decomposition, satisfies $r \leq  \min\{m,n\}$. 
During the fine-tuning phase, the pre-trained model parameter $\bW$ remains unchanged (often referred to as ``frozen"), while only the low-rank matrices $\bB$ and $\bA$ are updated.
The forward pass, originally defined as $\bo = \bW \bx + \bb$, is adjusted to:
\begin{equation}\label{equation:lora_eq}
\bo = \bW\bx+\bb +\alpha \Delta \bW\bx
=\bW\bx+\bb +\alpha \bB\bA\bx,
\end{equation}
where $\bb$ is the bias term, and $\alpha$ is called the \textit{merging ratio}, which controls the balance between retaining the pre-trained model's information and adapting it to  new target concepts. Hence the name \textit{low-rank adaptation (LoRA)} \citep{hu2021lora}.

\paragraph{LoHA.} 
In particular, it is widely acknowledged that methods based on matrix factorization face limitations due to the low-rank constraint. Within the LoRA framework, weight updates are confined to a low-rank space, which can impact the performance of the fine-tuned model. For improved fine-tuning, we hypothesize that a relatively higher rank might be necessary, especially when dealing with larger fine-tuning datasets or when the data distribution of downstream tasks significantly differs from the pretraining data. However, this could lead to increased memory usage and greater storage requirements.
We have demonstrated in Theorem~\ref{theorem:rank_hada_prod} that the Hadamard product of two rank-$\overline{r}$ matrices $\bW_1=\bB_1\bA_1$ and $\bW_2=\bB_2\bA_2$ has a rank of at most $\overline{r}^2$.
Given that  $\bA$ and $\bB$ in Equation~\eqref{equation:lora_eq} require $(m+n)r$ floating-points numbers, if we set $\overline{r} = \frac{r}{2}$, the number of floating-point numbers remains $(m+n)r$ for the Hadamard product of $(\bB_1\bA_1)\hadaprod (\bB_2\bA_2)$, while the rank of the low-rank approximation becomes $\frac{r^2}{4}$ at most. 
This technique is frequently referred to as  \textit{LoHA (low-rank adaptation with Hadamard product)} in the field (originally proposed to address  low-rank constraints  in federated learning problems  \citep{hyeon2021fedpara}). 
When $r>4$, LoHA can represent more complex models compared to LoRA.
The forward pass then becomes:
$$
\bo = \bW\bx+\bb +\alpha \Delta \bW\bx
=\bW\bx+\bb +\alpha \left((\bB_1\bA_1)\hadaprod (\bB_2\bA_2)\right) \bx.
$$
\paragraph{Warranty.} We have demonstrated in Section~\ref{section:low_rank_hadamard} that not all matrices can be factored by the Hadamard product of low-rank matrices. 
Therefore, when $\overline{r}=\frac{r}{2}$,  LoRA can search all the low-rank matrices within the space with rank $r$, whereas   LoHA  searches within an unspecified subset of the space with rank $\frac{r^2}{4}$. The advantages of LoHA in the field of federated learning or fine-tuning LLMs still need to be elucidated.

\paragraph{LoKr.} 
Similarly, the low-rank adaptation can be represented by the Kronecker product of two low-rank matrices $\bA\in\real^{m_1\times n_1}$ and $\bB\in\real^{m_2\times n_2}$ with $m=m_1m_2$ and $n=n_1n_2$. 
This approach is commonly known as LoKr (low-rank adaptation with Kronecker product, a.k.a., Kronecker adapter) in the literature \citep{edalati2022krona, yeh2023navigating}.
The forward pass is then given by 
$$
\bo = \bW\bx+\bb +\alpha \Delta \bW\bx
=\bW\bx+\bb +\alpha \left(\bA\kronecker \bB\right) \bx.
$$
Alternatively, one of the factors can be further  approximated by a low-rank matrix decomposition. For instance, consider  $\bB=\bC\bD$ with $\bC\in\real^{m_2\times k}$ and $\bD\in\real^{k\times n_2}$. In this case,  LoKr can be further represented as 
$$
\bo = \bW\bx+\bb +\alpha \Delta \bW\bx
=\bW\bx+\bb +\alpha \left(\bA\kronecker (\bC\bD)\right) \bx.
$$
The comparison among the three adaptation methods is shown in Figure~\ref{fig:lora_loha_lokr}.

\begin{figure}[h]
	\centering
	\includegraphics[width=1.02\textwidth]{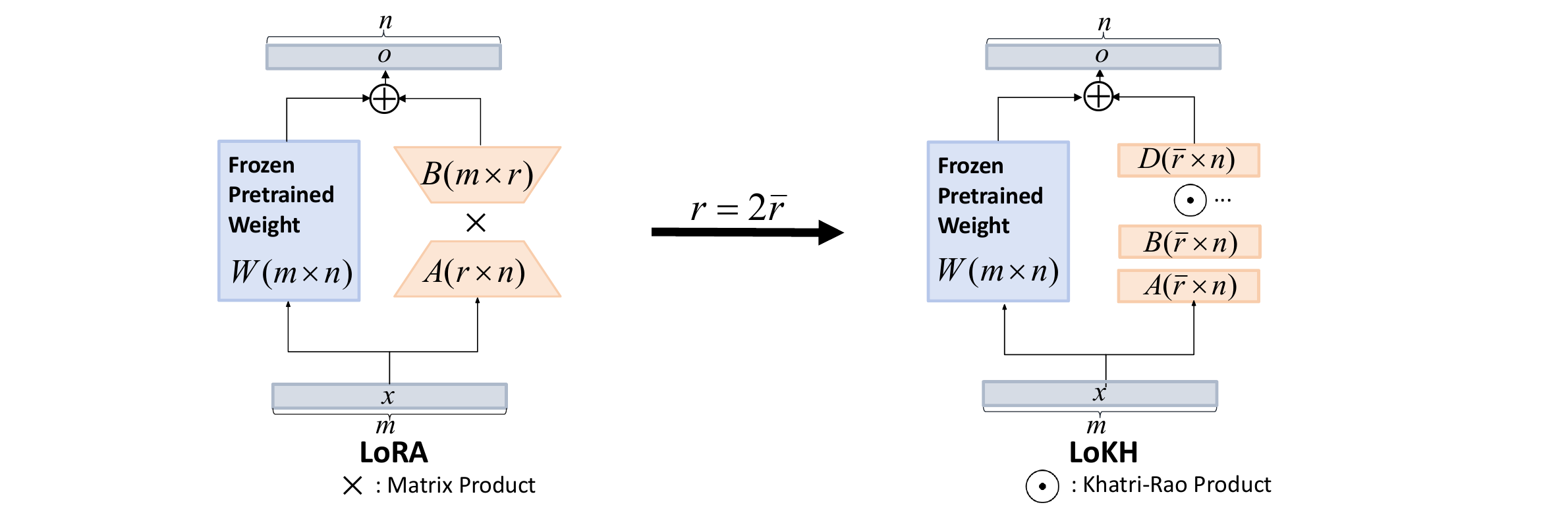}
	\caption{Diagram illustrating LoRA and LoKH.}
	\label{fig:lora_loKH}
\end{figure}

\index{LoKH}
\paragraph{Insights from Khatri-Rao Products: LoKH.}
We consider a special case where $\bW\in\real^{n\times n}$ and each low-rank matrix has the same rank and $k$-rank. This is not uncommon since the matrix structures in language models are extremely complex and usually have full rank. In the context of transformer architectures, the dimension $n$ is typically a power of $2$. As a result, there exists an integer $\overline{r}$ such that $\overline{r}^4 = n$. Let $\bA,\bB, \bC, \bD\in\real^{\overline{r}\times n}$ represent the low-rank matrices, the \textit{low-rank adaptation with Khatri-Rao product (LoKH)} can be expressed as (see Figure~\ref{fig:lora_loKH}):
$$
\bo = \bW\bx+\bb +\alpha \Delta \bW\bx
=\bW\bx+\bb +\alpha \left(\bA\khatrirao  \bB \khatrirao \bC\khatrirao \bD\right) \bx.
$$
In this case,  LoRA with rank $r=2\overline{r}$ can represent models with a rank of at most $r$ (when the matrix $\bA\in\real^{n\times n}$ is square  with $m=n$). However, if we explore the space where the rank is equal to the  $k$-rank, the LoKH method can find models with a rank greater than  $4\overline{r}-3$ (Theorem~\ref{theorem:kkrank_khatri_prod}).
When $r>3$, LoKH can represent more complex models compared to LoRA. 
Therefore, LoKH (rank $2r-3$) strikes a balance between LoRA (rank $r$) and LoHA (rank $\frac{r^2}{4}$).
LoHA represents a larger range of models due to its higher rank, but in a small subset of this high-rank space; LoRA offers a simpler model.
To the best of our knowledge, the LoKH method has not been explored in the fields of text-to-image models or large language models. This presents an opportunity for future research to investigate the potential benefits of LoKH in these domains, potentially leading to improved performance and efficiency in transformer-based models.
Furthermore, the LoKH method is straightforward to extend by cascading more components, allowing for the representation of more complex models: $\Delta\bW= \bA\khatrirao \bB\khatrirao \bC\khatrirao \bD\khatrirao \bE\ldots$.

Moreover, the Khatri-Rao decomposition has the property of matching the columns of the weight $\bW$. In the transformer structure, the attention mechanism  computes
$$
\text{Attention} (\bQ, \bK, \bV) = \text{softmax}\left( \frac{\bQ^\top\bK}{\sqrt{d_k}}\right)\bV^\top,
$$
where $\bQ, \bK \in\real^{d_k\times n}$ are the query and key matrices, $\bV\in\real^{ d_v\times n}$ is the value matrix, $n$ is the sequence length   (tokens for texts or patches for images), and $d_k, d_v$ are the dimension lengths. Specifically, the attention  captures the interactions between the query, key, and value matrices in a way that emphasizes the matching of features across different tokens/patches.
The LoKH method on these matrices, on the other hand, can be seen as a technique that facilitates the matching of features for each token/patch in the sequence before sending them into the attention mechanism.
By leveraging the Khatri-Rao product, LoKH enables the model to focus on the specific features of individual tokens/patches instead of interleaving the features of different tokens/patches (before sending them into the attention mechanism), which can be beneficial for tasks requiring fine-grained attention mechanisms.

Finally, the Khatri-Rao method is flexible, allowing us to explore  partition-wise decomposition (Definition~\ref{definition:partition_khatri_prod}). In this case, the low-rank decomposition separates the tokens/patches into several blocks, where  local information is exchanged more extensively within each block.
The method can be made even more flexible to identify various shifted windows across the blocks and then combine two (or more) partition-wise decompositions using the Hadamard product. 
The low-rank decomposition then assumes local/group connections for some tokens/patches, while each group is not connected (akin to the group Lasso formulation; see Section~\ref{section:regularization-extention-general} and \citet{yuan2006model, bach2011convex}).
This approach can be tailored to transformer architectures for computer vision, akin to the shifted windows and patch merging techniques used in the Swin Transformer \citep{liu2021swin}.

\begin{figure}[h]
\centering
\includegraphics[width=1.02\textwidth]{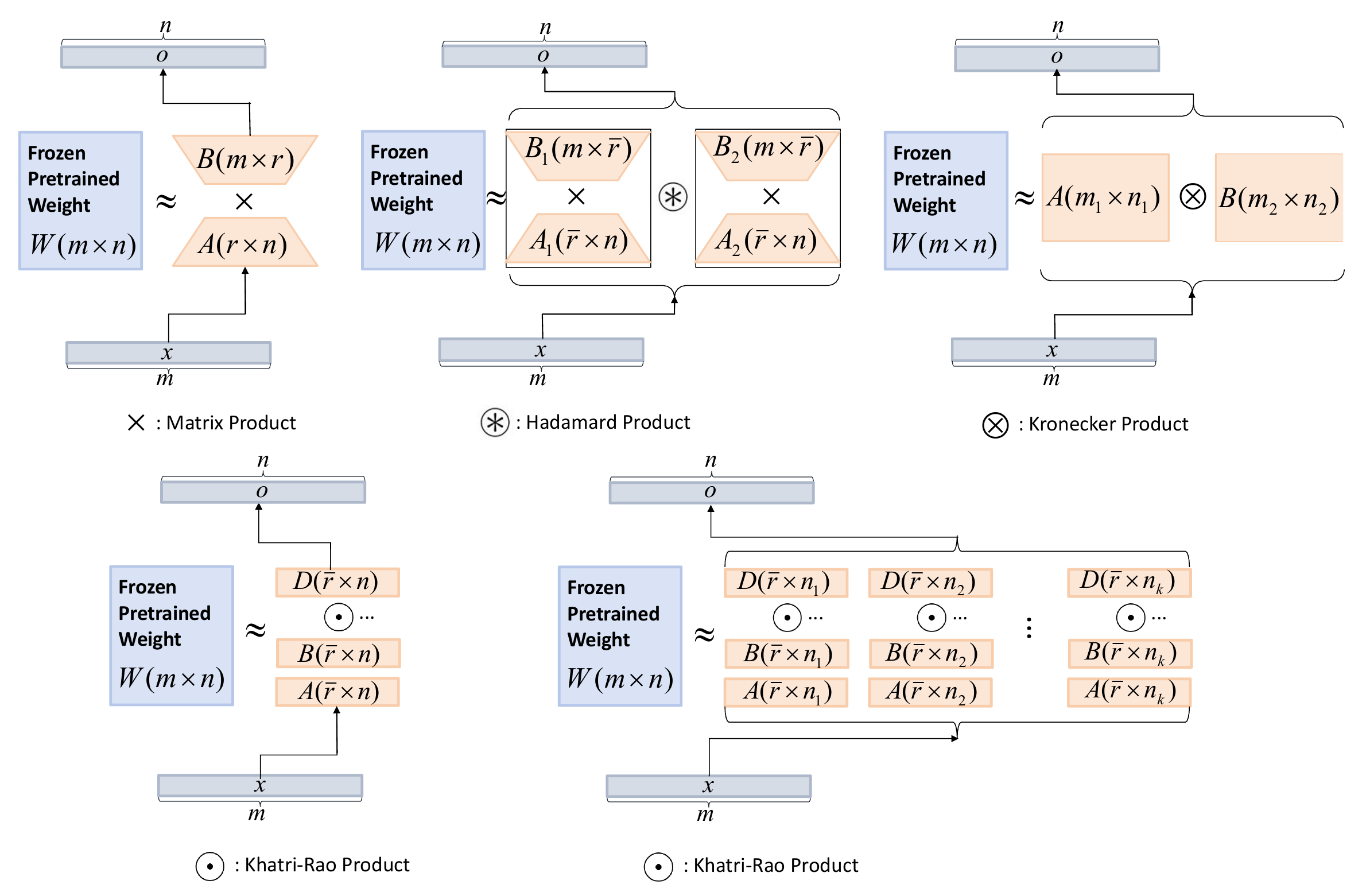}
\caption{Diagram illustrating low-rank approximation for transformer architectures. The bottom-right figure illustrates the cascading of Khatri-Rao products such that $n=n_1n_2\ldots n_k$.}
\label{fig:lora_loha_lokr_approx}
\end{figure}

\subsection{Low-Rank Approximation of Transformers}

We have demonstrated that low-rank matrix decomposition methods can be effectively applied to find low-rank approximations of neural networks in Section~\ref{section:low-rank-neural}. In this context, each layer of the neural network is factored into low-rank expressions, and the resulting network is then fine-tuned. This process enables the compression of deep neural networks, making them suitable for deployment on embedded devices with limited computational resources.

With the success of methods like LoRA, LoHA, and LoKr in the field of fine-tuning pre-trained large language models, it remains worthwhile to explore the application of low-rank decomposition techniques on pre-trained transformer architectures. By compressing these transformers using low-rank approximations (e.g., ALS in Section~\ref{section:als-netflix}, low-rank Hadamard decomposition in Section~\ref{section:low_rank_hadamard},  low-rank Kronecker decomposition in Section~\ref{section:low_rank_kronecker}, and low-rank Khatri-Rao decomposition in Section~\ref{section:lrank_khatri_decom},  see Figure~\ref{fig:lora_loha_lokr_approx}) and then fine-tuning the compressed models for downstream tasks, we can potentially achieve significant reductions in model size without sacrificing performance.

\begin{problemset}
\item \label{prob:hadama_1} Let $\ba,\bb,\bc\in\real^n$. Show that $\ba^\top(\bb\hadaprod\bc)=\bb^\top(\ba\hadaprod\bc)=\bc^\top(\ba\hadaprod\bb) = \sum_{i=1}^{n}a_ib_ic_i$.

\item Let $\bA\in\real^{n\times n}$ and $\bb\in\real^n$. Show that $\diag(\bb)\bA = \bA\hadaprod (\bb\bone_{n}^\top)$ and $\bA\diag(\bb) = \bA\hadaprod (\bone_{n}\bb^\top)$, where $\bone_{ n}\in\real^{n}$ is a   vector of all ones. 

\item Let $\ba,\bb\in\real^m$ and $\bc,\bd\in\real^n$. Show that $(\ba\bc^\top)\hadaprod(\bb\bd^\top)=(\ba\hadaprod\bb)(\bc\hadaprod\bd)^\top$.

\item Let $\bA,\bB\in\real^{m\times n}$ be given, and let $\bD_1\in\real^{m\times m}$ and $\bD_2\in\real^{n\times n}$ be diagonal. Show that $(\bD_1\bA)\hadaprod(\bB\bD_2)=\bD_1(\bA\hadaprod\bB)\bD_2$.

\item Let $\bA,\bB\in\real^{m\times n}$. Show that $\trace\big((\bA\hadaprod\bB)(\bA\hadaprod\bB)^\top\big) = \trace\big((\bA\hadaprod\bA)(\bB\hadaprod\bB)^\top\big)$.

\item Let $\bA\in\real^{n\times n}$ be nonsingular. Show that 
$(\bA\hadaprod\bA^{-\top})\bone_n=\bone_n$.

\item Let $\bA,\bB\in\real^{m\times n}$ and $\bC\in\real^{n\times m}$. Show that 
$
\bI\hadaprod(\bA(\bB^\top\hadaprod\bC))=\bI\hadaprod((\bA\hadaprod\bB)\bC)
=
\bI\hadaprod((\bA\hadaprod\bC^\top)\bB^\top)
$
and 
$
\trace(\bA(\bB^\top\hadaprod\bC))=\trace((\bA\hadaprod\bB)\bC)
=
\trace((\bA\hadaprod\bC^\top)\bB^\top).
$

\item  Let $\bA, \bB\in\real^{n\times n}$ and $\bx,\by\in\real^n$ be given. Show that 
$$
\begin{aligned}
\bx^\top(\bA\hadaprod\bB)\by &= \trace(\diag(\bx)\bA\diag(\by)\bB^\top);\\
(\bx\by^\top)\hadaprod\bA &= \diag(\bx)\bA\diag(\by).
\end{aligned}
$$
\textit{Hint: Use the implication from $\normf{\bA}^2 = \trace(\bA\bA^\top)$.}


\index{Fejer theorem}
\item \textbf{Hadamard PD/PSD, Fejer theorem \citep{horn2012matrix}.} Let $\bA\in\real^{n\times n}$ be given. Show that $\bA$ is PSD if and only if $\sum_{i,j=1}^{n,n}a_{ij}b_{ij}\geq 0$ holds for  all $n\times n$ PSD matrices $\bB=[b_{ij}]\in\real^{n\times n}$.

\item \textbf{Oppenheim’s inequality \citep{million2007hadamard}.} Let $\bA,\bB\in\real^{n\times n}$ be PSD. Show that $\det(\bA\hadaprod\bB)\geq a_{11}a_{22}\ldots a_{nn}\det(\bB)$ and $\det(\bA\hadaprod\bB) \geq \det(\bA\bB)$. \textit{Hint: Theorem~\ref{theorem:hadamard_ineq_psd}.}

\item \label{prob:hadama_n} \textbf{Hadamard PD/PSD.} Let $\bA,\bB\in\real^{n\times n}$ be PSD. Show that $\lambda_{\min}(\bA\bB) \leq \lambda_{\min}(\bA\hadaprod\bB)$.

\item \label{problem:gene_kro_mat}  Let $\bA \in \real^{I\times J}$, $\bB\in \real^{K\times L}$, $\bC\in \real^{J\times P}$, and $\bD\in \real^{L\times Q}$. Show that $(\bA\kronecker \bB) (\bC\kronecker \bD) = (\bA\bC) \kronecker (\bB\bD)$. \textit{Hint: Express each element of the equality explicitly.}

\item \label{prob:kronec_1} Let $\ba,\bb\in\real^n$. Show that $\ba\kronecker\bb = (\ba\kronecker\bI)\bb=(\bI\kronecker \bb)\ba$.

\item Let $\bA\in\real^{m\times m}$ and $\bB\in\real^{n\times n}$. Show that $\bA\kronecker\bB = (\bA\kronecker\bI_n)(\bI_m\kronecker \bB)$. \textit{Hint:use Equation~\ref{equ:mat_kro_prod}.}

\item  Let $\ba,\bb,\bc,\bd\in\real^n$. Show that $\ba^\top\bc\bb^\top\bd=(\ba^\top\kronecker\bb^\top)(\bc\kronecker\bd)=(\ba\kronecker\bb)^\top(\bc\kronecker\bd)$.

\item Let $\bA\in\real^{n\times n}, \bB\in\real^{m\times m}$. Show that $(\bA\kronecker \bB)^k=(\bA^k\kronecker\bB^k)$ for $k=1,2,\ldots$.

\item Let $\bA\in\real^{m\times \gamma}, \bB\in\real^{\gamma\times n}$. Show that $(\bA\bB)^{\kronecker k}=\bA^{\kronecker k}\bB^{\kronecker k}$ for $k=1,2,\dots$.

\item Let $\bA\in\real^{n\times n}$. Show that $\rho(\bA\kronecker\bA)=\rho(\bA)^2$ (spectral radius, Definition~\ref{definition:spectrum}).

\item \label{prob:kronec_n} Let $\bA_1,\bA_2,\ldots,\bA_k\in\real^{n\times n}$ be skew-symmetric. Show that $\bA_1\kronecker\bA_2\kronecker\ldots\kronecker\bA_k$ is symmetric (resp. skew-symmetric) if $k$ is even (resp. odd).

\item Prove Lemma~\ref{lemma:krokecker_keep_special}.
\item \label{problem:semi_khatri} Let $\bA, \bB\in\real^{I\times K}$ be semi-orthogonal. Show that $\bA\khatrirao \bB$ is also semi-orthogonal.
\item Conduct experiments on the  ``MovieLens 100K" data set  for the low-rank Kronecker decomposition.

\end{problemset}

%% file: chapter-biconjugate.tex
\newpage
\chapter{Biconjugate Decomposition}
\begingroup
\hypersetup{
	linkcolor=structurecolor,
	linktoc=page,  
}
\minitoc \newpage
\endgroup
\section{Existence of the Biconjugate Decomposition}
\lettrine{\color{caligraphcolor}T}
The concept of biconjugate decomposition was introduced by \citet{chu1995rank}, while the rank-diminishing operator, on the other hand, has roots in the work of  \citet{egervary1960rank, householder1964theory, stewart1973conjugate}.
A variety of matrix decomposition methods can be unified through  this biconjugate decomposition. 
In Section~\ref{section:conn_bi}, biconjugate decomposition is put into perspective by providing connections with standard decompositional methods, namely LDU, Cholesky, QR, and SVD decompositions.
The existence of the biconjugate decomposition is supported by the rank-one reduction theorem, as presented below.

\begin{theorem}[(Wedderburn's) Rank-One Reduction\index{Rank-one reduction}]\label{theorem:rank-1-reduction}
Let $\bA\in \real^{m\times n}$  be an $m\times n$ matrix with rank $r$, and let  $\bx\in \real^n$ and $\by\in \real^m$ be a pair of vectors such that $w=\by^\top\bA\bx \neq 0$. Then the matrix
\begin{equation}\label{equation:rank_reduc_operat}
\bB=\bA-w^{-1}\bA\bx\by^\top\bA
\end{equation}
has a rank of $r-1$, which is exactly one less than the rank of $\bA$: $\rank(\bB)=\rank(\bA)-1$.
A generalization of the rank-one reduction is discussed in Problem~\ref{problem:rk_reduc}.
\end{theorem}
\begin{proof}[of Theorem~\ref{theorem:rank-1-reduction}]
To prove the theorem, it suffices to show that the dimension of the null space of $\bB$ is one greater than that of $\bA$, indicating that $\bB$ has a rank exactly one less than the rank of $\bA$.

For any vector $\bn \in \nspace(\bA)$, i.e., $\bA\bn=\bzero$, we  have $\bB\bn =\bA\bn-w^{-1}\bA\bx\by^\top\bA\bn=\bzero $, implying that $\nspace(\bA)\subseteq \nspace(\bB)$.

Now, consider any vector $\bmm \in \nspace(\bB)$, i.e., $\bB\bmm=\bzero$. We have $\bB\bmm = \bA\bmm-w^{-1}\bA\bx\by^\top\bA\bmm =\bzero$.

Let $k\triangleq w^{-1}\by^\top\bA\bmm$, which is a scalar. 
Therefore, $\bB\bmm=\bA(\bmm - k\bx)=\bzero$, i.e., for any vector $\bn\in \nspace(\bA)$, we could find a vector $\bmm\in \nspace(\bB)$ such that $\bn=(\bmm - k\bx)\in \nspace(\bA)$. 
Note that $\bA\bx\neq \bzero$ based on the definition of $w$. Thus, the null space of $\bB$ is therefore obtained from the null space of $\bA$ by adding $\bx$ to its basis, which will increase the order of the space by one. 
Consequently, the dimension of $\nspace(\bA)$ is smaller than the dimension of $\nspace(\bB)$ by one, which completes the proof.
\end{proof}

The converse of the above theorem is also true, as shown in the following corollary.
\begin{corollary}[Rank-One Reduction, \citep{egervary1960rank}]\label{corollary:rk_redu_one}
Let  $\bA\in \real^{m\times n}$ be any $m\times n$ matrix, and let $\bu\in\real^m$ and $\bv\in\real^n$ be two vectors.
Then, the rank of matrix $\bB = \bA-\sigma^{-1}\bu\bv^\top$ is less than that of $\bA$ if and only if there exist vectors $\bx\in\real^n$ and $\by\in\real^m$ such that $\bu=\bA\bx, \bv=\bA^\top\by$, and $\sigma=\by^\top\bA\bx\neq 0$. In this case, it holds that $\rank(\bB)=\rank(\bA)-1$.
\end{corollary}

More generally, the rank-one reduction can be extended to the case where a matrix of rank possibly greater than one is subtracted.

\begin{corollary}[Rank-$k$ Reduction, \citep{cline1979rank}]\label{corollary:rk_K_redu_one}
Let $\bA\in \real^{m\times n}$ be any $m\times n$ matrix. Let further $\bP\in\real^{m\times k}$,  $\bU\in\real^{k\times k}$ be nonsingular, and $\bQ\in\real^{n\times k}$. Then,
$$
\rank(\bA-\bP\bU^{-1}\bQ^\top) = \rank(\bA)-\rank(\bP\bU^{-1}\bQ^\top)
$$
if and  only if there exist $\bX\in\real^{n\times k}$ and $\bY\in\real^{m\times k}$ such that 
$$
\bP=\bA\bX, \gap \bQ = \bA^\top\bY, \gap \text{and} \gap\bU=\bY^\top\bA\bX.
$$
\end{corollary}

Suppose a matrix $\bA\in \real^{m\times n}$ has rank $r$. We can define a rank-reducing process to generate a sequence of \textit{Wedderburn matrices} $\{\bA_k\}$:
\begin{equation}\label{equation:rank_reduc_proces}
\bA_1 = \bA \gap \text{and}\gap \bA_{k+1} = \bA_k-w_k^{-1}\bA_k\bx_k\by_k^\top\bA_k, \quad \forall k\in\{1,2,\ldots, r\},
\end{equation}
where $\bx_k \in \real^n$ and $\by_k\in \real^m$ are any vectors satisfying $w_k = \by_k^\top\bA\bx_k \neq 0$.
The operator in Equation~\eqref{equation:rank_reduc_operat} is known as a \textit{rank-diminishing operator}, and the process described by Equation~\eqref{equation:rank_reduc_proces} is referred to the \textit{rank-reducing process}.
And the sets $\{\bx_1,\bx_2,\ldots, \bx_r\}$ and $\{\by_1,\by_2,\ldots, \by_r\}$ are called the \textit{vectors associated with the rank-reducing process}. Alternatively, if we let $\bX=[\bx_1,\bx_2,\ldots, \bx_r]$ and $\bY=[\by_1,\by_2,\ldots, \by_r]$, then ($\bX,\bY$) is said to  \textit{effect a rank-reducing process} for $\bA$. 
The sequence will terminate after $r$ steps since the rank of $\bA_k$ decreases by exactly one at each step. 
The sequence can be written out as follows:
$$
\begin{aligned}
	\bA_1 &= \bA, \\
\bA_1-\bA_{2} &= w_1^{-1}\bA_1\bx_1\by_1^\top\bA_1,\\
\bA_2-\bA_3 &=w_2^{-1}\bA_2\bx_2\by_2^\top\bA_2, \\
\bA_3-\bA_4 &=w_3^{-1}\bA_3\bx_3\by_3^\top\bA_3, \\
\vdots &=\vdots\\
\bA_{r-1}-\bA_{r} &=w_{r-1}^{-1}\bA_{r-1}\bx_{r-1}\by_{r-1}^\top\bA_{r-1}, \\
\bA_r-\bzero &=w_{r}^{-1}\bA_{r}\bx_{r}\by_{r}^\top\bA_{r}.
\end{aligned}
$$
By summing up the sequence, we obtain
\begin{equation}\label{equation:bi_rreduc_equa}
\begin{aligned}
\textbf{(Rank-reducing)}: &\quad (\bA_1-\bA_2)+(\bA_2-\bA_3)+\ldots+(\bA_{r-1}-\bA_{r})+(\bA_r-\bzero ) \\
&=\bA_1=\bA= \sum_{i=1}^{r}w_i^{-1}\bA_i\bx_i\by_i^\top\bA_i.
\end{aligned}
\end{equation}
\index{Decomposition: Biconjugate}
Therefore, we can derive the following decomposition directly from this rank-reducing process.
\begin{theoremHigh}[Biconjugate Decomposition: Form 1]\label{theorem:biconjugate-form1}
Let $\bA\in\real^{m\times n}$ any matrix with rank $r$.
This equality~\eqref{equation:bi_rreduc_equa}, which results from the rank-reducing process, implies the following matrix decomposition
$$
\bA = \bPhi \bOmega^{-1} \bPsi^\top,
$$
where $\bOmega=\diag(w_1, w_2, \ldots, w_r)$, $\bPhi=[\bphi_1,\bphi_2, \ldots, \bphi_r]\in \real^{m\times r}$, and $\bPsi=[\bpsi_1, \bpsi_2, \ldots, \bpsi_r]$ with 
$$
\bphi_k = \bA_k\bx_k 
\qquad \text{and}\qquad 
\bpsi_k=\bA_k^\top \by_k, 
\qquad \forall k\in\{1,2,\ldots,r\}.
$$
\end{theoremHigh}
Thus, different choices of $\bx_k$'s and $\by_k$'s will result in different biconjugate factorizations, making this factorization  quite general and  versatile. In the following sections, we will demonstrate its connection with some well-known decomposition methods.
\index{Wedderburn sequence}

\begin{remark}\label{remark:weded_space}
Regarding the vectors $\bx_k$ and $\by_k$ in the Wedderburn sequence, the following property holds:
$$
\begin{aligned}
\bx_k &\in \nspace(\bA_{k+1})      \quad \implies \quad   \bx_k\bot \cspace(\bA_{k+1}^\top), \\
\by_k &\in \nspace(\bA_{k+1}^\top) \quad \implies \quad   \by_k\bot \cspace(\bA_{k+1}).
\end{aligned}
$$
To see why, consider the following:
$$
\begin{aligned}
\bA_{k+1} \bx_k 
&=  (\bA_k-w_k^{-1}\bA_k\bx_k\by_k^\top\bA_k) \bx_k\\
&=  \bA_k(\bx_k - \frac{\by_k^\top\bA_k\bx_k}{\by_k^\top\bA_k\bx_k}\bx_k)=\bzero.
\end{aligned}
$$
Similarly, it can be shown that $ \bA_{k+1}^\top\by_k=\bzero$.
\end{remark}

\begin{lemma}[General Term Formula of Wedderburn Sequence: V1]\label{lemma:wedderburn-sequence-general}
Let $\bA\in\real^{m\times n}$ be any matrix with rank $r$, and let $\bA_1\triangleq\bA$.
For each matrix in the sequence defined by $\bA_{k+1} = \bA_k -w_k^{-1} \bA_k\bx_k\by_k^\top \bA_k$ ($k\in\{1,2,\ldots,r-1\}$), the matrix $\bA_{k+1} $ can be expressed as 
$$
\bA_{k+1} = \bA - \sum_{i=1}^{k}w_{i}^{-1} \bA\bu_i \bv_i^\top \bA,
\qquad \forall k\in\{1,2,\ldots,r-1\},
$$
where 
$$
\bu_k=\bx_k -\sum_{i=1}^{k-1}\frac{\bv_i^\top \bA\bx_k}{w_i}\bu_i
\qquad \text{and}\qquad 
\bv_k=\by_k -\sum_{i=1}^{k-1}\frac{\by_k^\top \bA\bu_i}{w_i}\bv_i.
$$
Let $\bX=[\bx_1,\bx_2,\ldots,\bx_r]$, $\bY=[\by_1,\by_2,\ldots,\by_r]$, $\bU=[\bu_1,\bu_2,\ldots,\bu_r]$, and $\bV=[\bv_1,\bv_2,\ldots,\bv_r]$ be the column partitions for each set of vectors.
Then, the rank-reducing process can be viewed as transforming the  matrix pair $(\bX,\bY)$ into $(\bU,\bV)$.~\footnote{It can be shown that if $\bA$ is symmetric and $\bX=\bY$, then $\bU=\bV$.}
\end{lemma}

\begin{figure}[H]
\centering  
\vspace{-0.35cm} 
\subfigtopskip=2pt 
\subfigbottomskip=2pt 
\subfigcapskip=-5pt 
\subfigure[``Project" onto a line.]{\label{fig:project-line-wedd}
\includegraphics[width=0.47\linewidth]{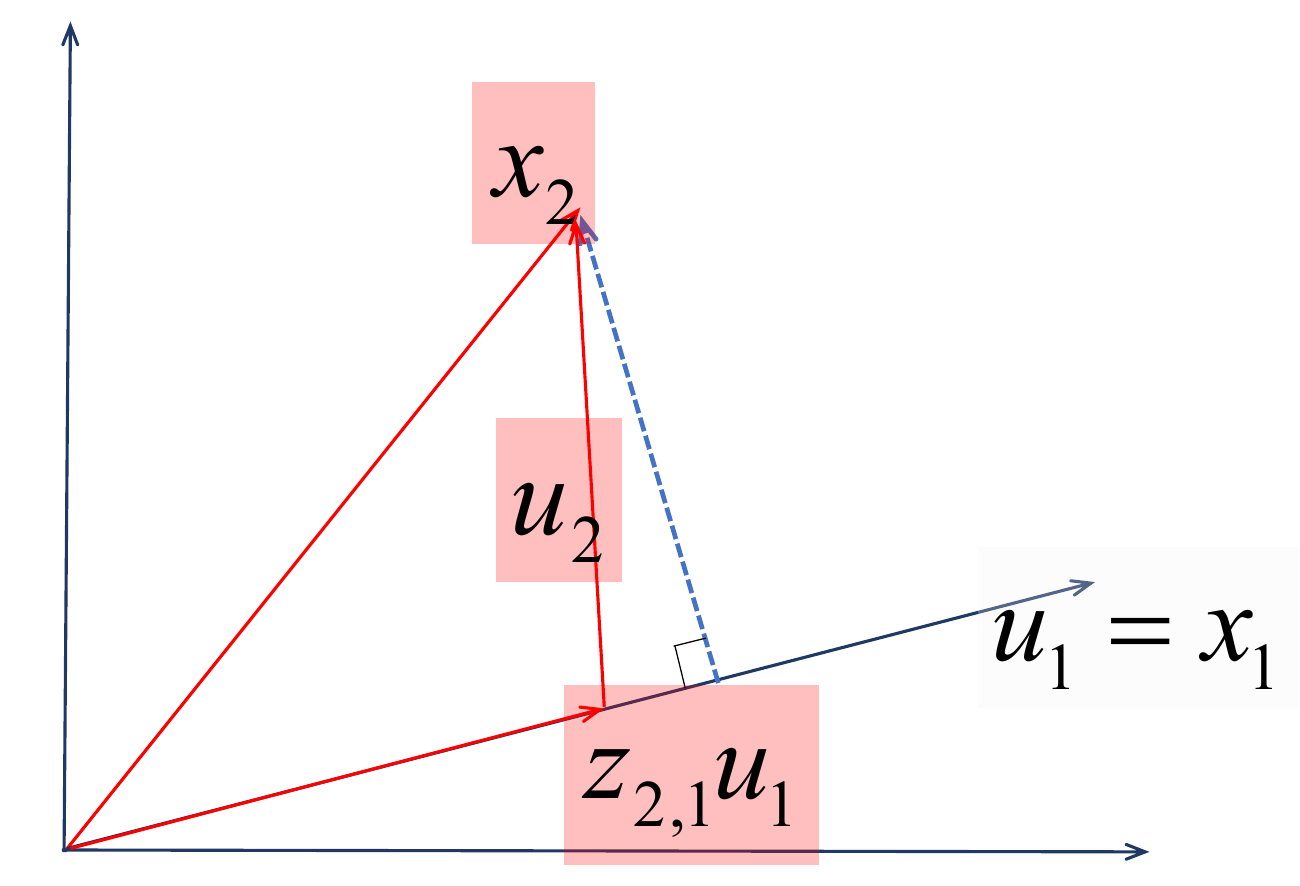}}
\quad 
\subfigure[``Project" onto a space.]{\label{fig:project-space-wedd}
\includegraphics[width=0.47\linewidth]{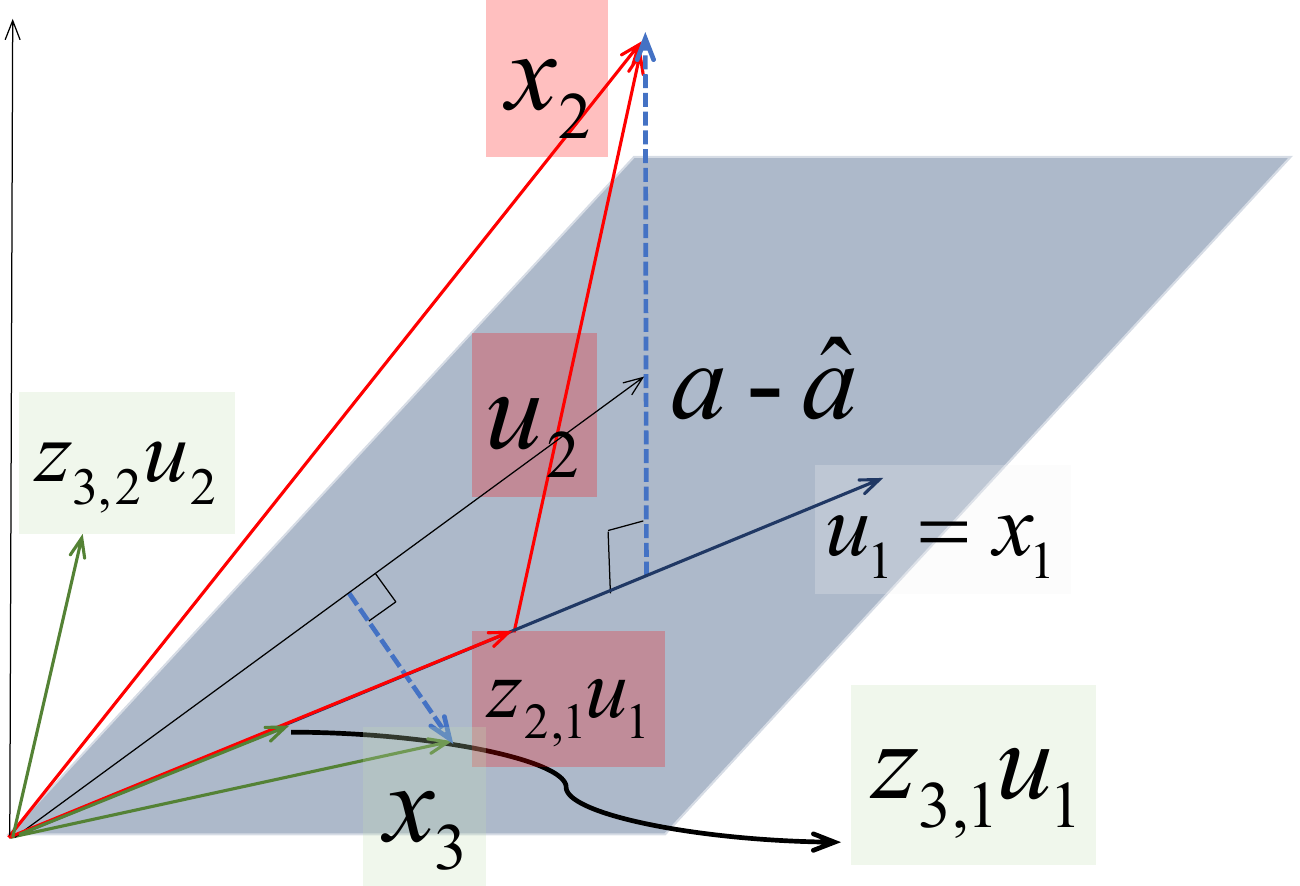}}
\caption{``Project" a vector onto a line and onto a space. Compare to the Gram-Schmidt process in Figure~\ref{fig:gram-schmidt-12}.}
\label{fig:projection-wedd}
\end{figure}

The proof of this lemma is deferred to Section~\ref{section:wedderburn-general-term}. We notice that $w_i =\by_i^\top\bA_i\bx_i$ in the general term formula is related to $\bA_i$, which means it is not the true general term formula. 
We will later reformulate $w_i$ to be related to $\bA$ rather than $\bA_i$.
From the general term formula of the  Wedderburn sequence, we have: 
$$
\begin{aligned}
\bA_{k+1} &= \bA - \sum_{i=1}^{k}w_{i}^{-1} \bA\bu_i \bv_i^\top \bA, \\
\bA_{k} &= \bA - \sum_{i=1}^{k-1}w_{i}^{-1} \bA\bu_i \bv_i^\top \bA.
\end{aligned}
$$
Thus, we obtain $\bA_{k+1} - \bA_{k} = -w_{k}^{-1} \bA\bu_k \bv_k^\top \bA$. 
Since  the sequence is defined as $\bA_{k+1} = \bA_k -w_k^{-1} \bA_k\bx_k\by_k^\top \bA_k$, we can deduce that $w_{k}^{-1} \bA\bu_k \bv_k^\top \bA = w_k^{-1} \bA_k\bx_k\by_k^\top \bA_k$. 
Consequently, it follows that
\begin{equation}\label{equation:wedderburn-au-akxk}
\begin{aligned}
\bA\bu_k &=\bA_k\bx_k,\\
\bv_k^\top \bA&=\by_k^\top \bA_k.
\end{aligned}
\end{equation}
Let $z_{k,i} \triangleq \frac{\bv_i^\top \bA\bx_k}{w_i}$, which is a scalar. Referring to the definitions of $\bu_k$ and $\bv_k$ in the  lemma above, we can express them explicitly as follows:
\begin{itemize}
\item $\bu_1=\bx_1$; 

\item $\bu_2 = \bx_2 - z_{2,1}\bu_1 \implies $ $\bx_2$ is a linear combination of $\bu_1$ and $\bu_2$;

\item $\bu_3 = \bx_3 - z_{3,1}\bu_1-z_{3,2}\bu_2 \implies$ $\bx_3$ is a linear combination of $\bu_1,\bu_2$, and $\bu_3$;

\item $\ldots$.
\end{itemize}
Each $z_{k,i}$ ($i< k$) encodes the component of $\bx_k$ in that of $\bu_i$.
This process bears resemblance to the Gram-Schmidt process (Section~\ref{section:gram-schmidt-process}). 
However, in this process,  we do not perform an \textbf{orthogonal projection} of $\bx_2$ onto $\bx_1$ to find the vector component of $\bx_2$ along $\bx_1$, as we would do in an orthogonal projection (Definition~\ref{definition:orthogonal-projection-matrix}). 
Instead, the vector of $\bx_2$ along $\bx_1$ is now defined by $z_{2,1}$ (i.e., an \textbf{oblique projection}, Figure~\ref{fig:ls-geometric1-compare}). 
This process is illustrated in Figure~\ref{fig:projection-wedd}. 
In Figure~\ref{fig:project-line-wedd}, $\bu_2$ is not perpendicular to $\bu_1$ (in the Gram-Schmidt process, $\bu_2$ would be perpendicular to $\bu_1$ via orthogonal projections). 
Nevertheless, $\bu_2$ does not lie on the same line as $\bu_1$, so $\{\bu_1, \bu_2\}$  can still  span a $\real^2$ subspace. Similarly, in Figure~\ref{fig:project-space-wedd}, $\bu_3= \bx_3 - z_{3,1}\bu_1-z_{3,2}\bu_2$ does not lie in the subspace spanned by $\{\bu_1, \bu_2\}$, allowing $\{\bu_1, \bu_2, \bu_3\}$  to span a $\real^3$ subspace.
\index{Gram–Schmidt}

A moment of reflexion would reveal that the span of $\{\bx_2, \bx_1\}$ is the same as the span of $\{\bu_2, \bu_1\}$. This equivalence extends to the $\bv_i$ vectors and $\by_i$ vectors as well.  
We can express this property as follows:
\begin{equation}\label{equation:wedderburn-span-same}
\left\{
\begin{aligned}
\spn\{\bx_1, \bx_2, \ldots, \bx_j\} &= \spn\{\bu_1, \bu_2, \ldots, \bu_j\},\gap \forall j\in\{1,2,\ldots, r\};\\
\spn\{\by_1, \by_2, \ldots, \by_j\} &= \spn\{\bv_1, \bv_2, \ldots, \bv_j\},\gap \forall j\in\{1,2,\ldots, r\}.\\
\end{aligned}
\right.
\end{equation}

Furthermore, based on the rank-reducing property in the Wedderburn sequence, we have 
$$
\left\{
\begin{aligned}
	\cspace(\bA_1) &\supset \cspace(\bA_2) \supset \cspace(\bA_3) \supset \ldots;\\
	\nspace(\bA_1^\top) &\subset \nspace(\bA_2^\top) \subset \nspace(\bA_3^\top) \subset \ldots.
\end{aligned}
\right.
$$
Since $\by_k \in \nspace(\bA_{k+1}^\top)$, it then follows that $\by_j \in \nspace(\bA_{k+1}^\top)$ for all $j<k+1$, i.e., $\bA_{k+1}^\top \by_j=\bzero$ for all $j<k+1$. 
This also holds true for $\bx_{k+1}^\top \bA_{k+1}^\top \by_j=0$ for all $j<k+1$. From Equation~\eqref{equation:wedderburn-au-akxk}, we also have $\bu_{k+1}^\top \bA^\top \by_j=0$ for all $j<k+1$. Following Equation~\eqref{equation:wedderburn-span-same}, we obtain 
\begin{equation}
\bx_{k+1}^\top \bA_{k+1}^\top \by_j=0
\,\,\stackrel{\eqref{equation:wedderburn-au-akxk}}{\implies}\,\,
\bu_{k+1}^\top \bA^\top \by_j=0
\,\,\stackrel{\eqref{equation:wedderburn-span-same}}{\implies}\,\,
\bv_j^\top\bA\bu_{k+1}=0\quad  \text{ for all } j<k+1.
\end{equation}
Similarly, we can prove 
\begin{equation}
\bv_{k+1}^\top\bA\bu_{j}=0\quad \text{ for all } j<k+1.
\end{equation}
Moreover, since $w_k = \by_k^\top \bA_k\bx_k$, according to Equation~\eqref{equation:wedderburn-au-akxk}, we can express  $w_k$ as:
$$
\begin{aligned}
w_k &= \by_k^\top \bA_k\bx_k
=\bv_k^\top \bA\bx_k \\
&=\bv_k^\top\bA (\bu_k +\sum_{i=1}^{k-1}\frac{\bv_i^\top \bA\bx_k}{w_i}\bu_i) \qquad &(\text{by the definition of }\bu_k \text{ in Lemma~\ref{lemma:wedderburn-sequence-general}})\\
&=\bv_k^\top\bA \bu_k,\qquad &(\text{by } \bv_{k}^\top\bA\bu_{j}=0 \text{ for all } j<k)
\end{aligned}
$$
which can be utilized to substitute the value of $w_k$ in Lemma~\ref{lemma:wedderburn-sequence-general}. We then have the full version of the general term formula of the Wedderburn sequence. 
In this form, the formula no longer depends on $\bA_k$'s (in the form of $w_k$'s): 
\begin{equation}\label{equation:uk-vk-to-mimic-gram-process}
\bu_k=\bx_k -\sum_{i=1}^{k-1}\frac{\bv_i^\top \bA\bx_k}{\textcolor{mylightbluetext}{\bv_i^\top\bA\bu_i}}\bu_i
\qquad \text{and}\qquad 
\bv_k=\by_k -\sum_{i=1}^{k-1}\frac{\by_k^\top \bA\bu_i}{\textcolor{mylightbluetext}{\bv_i^\top\bA\bu_i}}\bv_i.
\end{equation}

\index{Gram–Schmidt}
\paragraph{Gram-Schmidt process from Wedderburn sequence.} 
Suppose the matrices $\bX=[\bx_1, 
\bx_2, \ldots, \bx_r]\in \real^{n\times r}$ and $\bY=[\by_1, \by_2, \ldots, \by_r]\in \real^{n\times r}$ effect a rank-reducing process for $\bA$. 
If $\bA=\bI\in\real^{n\times n}$ is the identity matrix and $(\bX=\bY)$ are identical, containing the vectors for which an orthogonal basis is desired, i.e., we aim to obtain an orthogonal basis spanning the same column space of $\bX=\bY$, then  the result of the rank-reducing process, $(\bU = \bV)$, gives the resultant orthogonal basis (but not an orthonormal basis).

To see this, we follow the computation of the Gram-Schmidt process as given in Equation~\eqref{equation:qr-gsp-equation}, where we replace the data matrix by $\bX=\bY\in\real^{m\times n}$ (assuming $\bX$ has full rank for simplicity), and we use matrix $\bQ=[\bq_1,\bq_2, \ldots, \bq_n]$ to denote the orthonormal basis.
Then, for any $k\in\{1,2,\ldots,n\}$, the Gram-Schmidt process admits 
\begin{equation}\label{equation:bicon_cgs}
\text{Gram-Schmidt process}=\left\{
\begin{aligned}
\bx_k^\perp&= \left(\bI-\sum_{i=1}^{k-1}\bq_i\bq_i^\top\right)\bx_k;\\
\bq_k &= \frac{\bx_k^\perp}{\norm{\bx_k^\perp}}.
\end{aligned}
\right.
\end{equation}
Considering the rank-reducing process with $\bA=\bI$, $\bX=\bY$. Since $\bA$ is symmetric and $\bX=\bY$, we obtain $\bU=\bV$. We can thus focus on the analysis of $\bU$. For any $k\in\{1,2,\ldots,n\}$, the rank-reducing process yields 
\begin{equation}\label{equation:bicon_wedd}
\text{Rank-reducing process}=\left\{
\begin{aligned}
\bu_k&=\bx_k -\sum_{i=1}^{k-1}\frac{\bv_i^\top \bA\bx_k}{\bv_i^\top\bA\bu_i}\bu_i
=\left(\bI -\sum_{i=1}^{k-1}\frac{\bu_i\bu_i^\top }{\bu_i^\top\bu_i}\bu_i\right)\bx_k ;\\
\widetilde{\bu}_k &= \frac{\bu_k}{\norm{\bu_k}}.
\end{aligned}
\right.
\end{equation}
Therefore, by comparing Equations~\eqref{equation:bicon_cgs} and \eqref{equation:bicon_wedd}, we can find the equivalence between $\bq_k$ and $\widetilde{\bu}_k$ for $k\in\{1,2,\ldots,n\}$. Thus, the Wedderburn sequence effects a Gram-Schmidt process when $\bX=\bY$ is the data matrix and $\bA=\bI$.


The expressions for $\bu_k$ and $\bv_k$ in Equation~\eqref{equation:uk-vk-to-mimic-gram-process} closely resemble the projection onto the perpendicular space in the Gram-Schmidt process, as shown in Equation~\eqref{equation:gram-schdt-eq2}. We then define the bilinear form \fbox{$<\bx, \by>\triangleq\by^\top\bA\bx$} to explicitly mimic the form of projection in Equation~\eqref{equation:gram-schdt-eq2}. 	

We consolidate the results obtained so far into the following lemma, which provides a concise overview of what we have been working on. These results will be extensively utilized  in the sequel.
\begin{lemma}[Properties of Wedderburn Sequence]\label{lemma:wedderburn-sequence-general-v2}
Let $\bA\in\real^{m\times n}$ be any matrix with rank $r$, and let $\bA_1\triangleq\bA$.
For each matrix in the sequence defined by $\bA_{k+1} = \bA_k -w_k^{-1} \bA_k\bx_k\by_k^\top \bA_k$ ($k\in\{1,2,\ldots,r-1\}$), the matrix $\bA_{k+1} $ can be expressed as 
$$
\bA_{k+1} = \bA - \sum_{i=1}^{k}w_{i}^{-1} \bA\bu_i \bv_i^\top \bA,
$$
where 
\begin{equation}\label{equation:properties-of-wedderburn-ukvk}
\bu_k=\bx_k -\sum_{i=1}^{k-1}\frac{\textcolor{mylightbluetext}{<\bx_k, \bv_i>}}{\textcolor{mylightbluetext}{<\bu_i,\bv_i>}}\bu_i
\qquad \text{and}\qquad 
\bv_k=\by_k -\sum_{i=1}^{k-1}\frac{\textcolor{mylightbluetext}{<\bu_i,\by_k>}}{\textcolor{mylightbluetext}{<\bu_i,\bv_i>}}\bv_i.
\end{equation}
Furthermore, we can observe the following properties:
\begin{equation}\label{equation:wedderburn-au-akxk-2}
\begin{aligned}
	\bA\bu_k &=\bA_k\bx_k;\\
	\bv_k^\top \bA&=\by_k^\top \bA_k;
\end{aligned}
\end{equation}
\begin{equation}\label{equation:wedderburn-au-akxk-33}
<\bu_k, \bv_j>=<\bu_j, \bv_k>=0 \text{ for all } j<k;
\end{equation}
\begin{equation}\label{equation:wk-by-ukvk}
w_k = \by_k^\top\bA_k\bx_k = <\bu_k, \bv_k>.
\end{equation}

\end{lemma}

By substituting Equation~\eqref{equation:wedderburn-au-akxk-2} into Form 1 of the biconjugate decomposition and using Equation~\eqref{equation:wk-by-ukvk}, which implies $w_k = \bv_k^\top\bA\bu_k$, we obtain the Form 2 and Form 3 of this decomposition:
\begin{theoremHigh}[Biconjugate Decomposition: Form 2 and Form 3]\label{theorem:biconjugate-form2}
Let $\bA\in\real^{m\times n}$ be any matrix with rank $r$.
This equality~\eqref{equation:bi_rreduc_equa}, which results from the rank-reducing process, implies the following matrix decomposition:
\begin{equation}\label{equation:bicon_form2_1}
\bA = \bA\bU_r \bOmega_r^{-1} \bV_r^\top\bA,
\end{equation}
where $\bOmega_r=\diag(w_1, w_2, \ldots, w_r)$, $\bU_r=[\bu_1,\bu_2, \ldots, \bu_r]\in \real^{n\times r}$, and  $\bV_r=[\bv_1, \bv_2,$ $\ldots,$ $\bv_r] \in \real^{m\times r}$ with  
\begin{equation}\label{equation:properties-of-wedderburn-ukvk2-inform2}
\bu_k=\bx_k -\sum_{i=1}^{k-1}\frac{<\bx_k, \bv_i>}{<\bu_i,\bv_i>}\bu_i
\qquad \text{and}\qquad 
\bv_k=\by_k -\sum_{i=1}^{k-1}\frac{<\bu_i,\by_k>}{<\bu_i,\bv_i>}\bv_i.
\end{equation}
Additionally, the following decomposition holds:
\begin{equation}\label{equation:wedderburn-vgamma-ugamma}
\bV_\gamma^\top \bA \bU_\gamma = \bOmega_\gamma,
\end{equation}
where $\bOmega_\gamma=\diag(w_1, w_2, \ldots, w_\gamma)$, $\bU_\gamma=[\bu_1,\bu_2, \ldots, \bu_\gamma]\in \real^{n\times \gamma}$, and $\bV_\gamma=[\bv_1, \bv_2,$ $\ldots,$ $\bv_\gamma]\in \real^{m\times \gamma}$. Note the difference between the subscripts $r$ and $\gamma$  employed here, with $\gamma \leq r$.
\end{theoremHigh}
Note that Equation~\eqref{equation:bicon_form2_1} is derived from  \eqref{equation:wedderburn-au-akxk-2}, 
and Equation~\eqref{equation:wedderburn-vgamma-ugamma} comes from \eqref{equation:wedderburn-au-akxk-33}. 
We observe that, in the above two forms of the biconjugate decomposition, they no longer depend on the Wedderburn matrices  $\{\bA_k\}$. 

\paragraph{A word on the notation.} We will employ subscripts to indicate the dimension of the matrix to prevent confusion in the following discussion, e.g., the $r, \gamma$ notations used in the preceding theorem.

\section{Properties of the Biconjugate Decomposition}
\begin{corollary}[Connection of $\bU_\gamma$ and $\bX_\gamma$]\label{corollary:biconjugate-connection-u-x}
Let $\bA\in\real^{m\times n}$ be any matrix with rank $r\geq \gamma$.
If $(\bX_\gamma, \bY_\gamma) \in \real^{n\times \gamma}\times \real^{m\times \gamma}$ effects a rank-reducing process for $\bA$, then there exist unique unit upper triangular matrices $\bR_\gamma^{(x)}\in\real^{\gamma\times \gamma}$ and $\bR_\gamma^{(y)}\in\real^{\gamma\times \gamma}$ such that 
$$
\bX_\gamma = \bU_\gamma \bR_\gamma^{(x)}
\qquad \text{and} \qquad 
\bY_\gamma=\bV_\gamma\bR_\gamma^{(y)},
$$
where $\bU_\gamma$ and $\bV_\gamma$ are matrices with columns derived from the Wedderburn sequence, as described in Equation~\eqref{equation:wedderburn-vgamma-ugamma}.
\end{corollary} 
\begin{proof}[of Corollary~\ref{corollary:biconjugate-connection-u-x}]
The proof follows directly from the definitions of $\bu_k$ and $\bv_k$ in Equation~\eqref{equation:properties-of-wedderburn-ukvk} or Equation~\eqref{equation:properties-of-wedderburn-ukvk2-inform2}. We set  the $j$-th column of $\bR_\gamma^{(x)}$ and $\bR_\gamma^{(y)}$ as follows:
$$
\left[\frac{<\bx_j,\bv_1>}{<\bu_1,\bv_1>}, \frac{<\bx_j,\bv_2>}{<\bu_2,\bv_2>}, \ldots, \frac{<\bx_j,\bv_{j-1}>}{<\bu_{j-1},\bv_{j-1}>}, 1, 0, 0, \ldots, 0  \right]^\top,
$$
and 
$$
\left[\frac{<\bu_1, \by_j>}{<\bu_1,\bv_1>},\frac{<\bu_2, \by_j>}{<\bu_2,\bv_2>}, \ldots, \frac{<\bu_{j-1}, \by_j>}{<\bu_{j-1},\bv_{j-1}>}, 1, 0, 0, \ldots, 0  \right]^\top.
$$
And the uniqueness stems  from the fact that the matrices $\bU_\gamma$ and $\bV_\gamma$  have independent columns from the rank-reducing process.
This completes the proof.
\end{proof}

The pair $(\bU_\gamma, \bV_\gamma) \in \real^{m\times \gamma}\times \real^{n\times \gamma}$ in Theorem~\ref{theorem:biconjugate-form2} is called a \textbf{biconjugate pair} with respect to $\bA$ if $\bOmega_\gamma$ is nonsingular and diagonal. 
Furthermore, suppose the pair $(\bX_\gamma, \bY_\gamma) \in \real^{n\times \gamma}\times \real^{m\times \gamma}$ effects a rank-reducing process for $\bA$. 
Then, the pair $(\bX_\gamma, \bY_\gamma)$ is said to be \textbf{biconjugatable} and can be \textbf{biconjugated into a biconjugate pair} of matrices $(\bU_\gamma, \bV_\gamma)$, if there exist unit upper triangular matrices $\bR_\gamma^{(x)}$ and $\bR_\gamma^{(y)}$ such that $\bX_\gamma = \bU_\gamma \bR_\gamma^{(x)}$ and $\bY_\gamma=\bV_\gamma\bR_\gamma^{(y)}$.


\section{Connection to Well-Known Decomposition Methods}\label{section:conn_bi}
In this section, we demonstrate how biconjugate decomposition can be related to well-known decompositional methods.

\subsection{LDU Decomposition}
\begin{theorem}[LDU, \cite{chu1995rank} Theorem 2.4]\label{theorem:biconjugate-ldu}
Let $\bA\in\real^{m\times n}$ be any matrix with rank $r\geq \gamma$.
Let further $(\bX_\gamma, \bY_\gamma) \in \real^{n\times \gamma}\times \real^{m\times \gamma}$ 
 with $\gamma$ in $\{1, 2, \ldots, r\}$. Then, the pair  $(\bX_\gamma, \bY_\gamma)$ can be biconjugated if and only if $\bY_\gamma^\top\bA\bX_\gamma$ has an LDU decomposition. 
\end{theorem}
\begin{proof}[of Theorem~\ref{theorem:biconjugate-ldu}]
Suppose $\bX_\gamma$ and $\bY_\gamma$ are biconjugatable; then there exist unit upper triangular matrices $\brx$ and $\bry$ such that $\bxgamma = \bugamma\brx$, $\bygamma = \bvgamma\bry$, and $\bvgamma^\top\bA\bugamma = \bomegagamma$ is a nonsingular diagonal matrix. Then  it follows that 
$$
\bygamma^\top\bA\bxgamma = \bryt \bvgamma^\top \bA \bugamma\brx = \bryt \bomegagamma \brx
$$
is the unique  LDU decomposition of $\bygamma^\top\bA\bxgamma$. This expression can be regarded as the \textbf{fourth form of biconjugate decomposition}.

Conversely, suppose $\bygamma^\top\bA\bxgamma = \bR_2^\top \bD\bR_1$ is an LDU decomposition, with both $\bR_1$ and $\bR_2$ being unit upper triangular matrices. Since $\bR_1^{-1}$ and $\bR_2^{-1}$ are also unit upper triangular matrices,  $(\bX_\gamma, \bY_\gamma)$ biconjugates into $(\bX_\gamma\bR_1^{-1}, \bY_\gamma\bR_2^{-1})$.
\end{proof}

\begin{corollary}[Determinant]\label{corollary:lu-determinant}
Let $\bA\in\real^{m\times n}$ be any matrix with rank $r\geq \gamma$.
Suppose the pair $(\bX_\gamma, \bY_\gamma) \in \real^{n\times \gamma}\times \real^{m\times \gamma}$ can be biconjugated into $(\bU_\gamma, \bV_\gamma)$ such that $\bvgamma^\top\bA\bugamma = \bomegagamma=\diag(w_1,w_2,\ldots,w_\gamma)$ is a nonsingular diagonal matrix. 
Then it follows that
$$
\det(\bygamma^\top \bA\bxgamma) = \prod_{i=1}^{\gamma} w_i.
$$
\end{corollary}
\begin{proof}[of Corollary~\ref{corollary:lu-determinant}]
By Theorem~\ref{theorem:biconjugate-ldu}, since $(\bX_\gamma, \bY_\gamma)$ are biconjugatable, then there exist unit upper triangular matrices $\brx$ and $\bry$ such that $\bygamma^\top\bA\bxgamma = \bryt \bomegagamma \brx$. The determinant is simply the product of the diagonal elements.
\end{proof}

\begin{lemma}[Biconjugatable in Principal Minors]\label{lemma:Biconjugatable-in-Principal-Minors}
Let $\bA\in\real^{m\times n}$ be any matrix with rank $r\geq \gamma$.
In the Wedderburn sequence, we choose $\bx_i$ as the $i$-th standard basis in $\real^n$ for $i \in \{1, 2, \ldots, \gamma\}$ (i.e., $\bx_i = \be_i \in \real^n$), and $\by_i$ as the $i$-th standard basis in $\real^m$ for $i \in \{1, 2, \ldots, \gamma\}$ (i.e., $\by_i=\be_i \in 
\real^m$). 
That is, $\bygamma^\top \bA\bxgamma$ corresponds to the leading principal submatrix of $\bA$, i.e., $\bygamma^\top \bA\bxgamma = \bA[1:\gamma, 1:\gamma]$. 
Then, $(\bX_\gamma, \bY_\gamma)$ is biconjugatable  into $(\bU_\gamma, \bV_\gamma)$ such that $\bvgamma^\top\bA\bugamma = \bomegagamma=\diag(w_1,w_2,\ldots,w_\gamma)$ is a nonsingular diagonal matrix
if and only if the $\gamma$-th leading principal minor of $\bA$ is nonzero: $\det(\bA[1:\gamma, 1:\gamma])\neq 0$. In this case, the $\gamma$-th leading principal minor of $\bA$ is given by $\prod_{i=1}^{\gamma} w_i$.\index{Leading principal minors}
\end{lemma}
\begin{proof}[of Lemma~\ref{lemma:Biconjugatable-in-Principal-Minors}]
The proof is straightforward that the $\gamma$-th leading principal minor of $\bA$ being nonzero will imply that $w_i \neq 0$ for all $i\leq \gamma$. Thus, the Wedderburn sequence can be successfully obtained. The converse holds because Corollary~\ref{corollary:lu-determinant} implies that $\det(\bygamma^\top \bA\bxgamma)$ is nonzero.
\end{proof}

We have now arrived at the LDU decomposition for square matrices. 
\begin{theorem}[LDU: Biconjugate Decomposition for Square Matrices]\label{theorem:biconjugate-square-ldu}
For any matrix $\bA\in \real^{n\times n}$, the pair $(\bI_n, \bI_n)$ is biconjugatable if and only if all the leading principal minors of $\bA$ are nonzero. In this case, $\bA$ can be factored as 
$$
\bA = \bV_n^{-\top} \bOmega_n \bU_n^{-1} = \bL\bD\bU,
$$
where $\bOmega_n = \bD$ is a diagonal matrix with nonzero values along its diagonal, $\bV_n^{-\top} = \bL$ is a unit lower triangular matrix, and $\bU_n^{-1} = \bU$ is a unit upper triangular matrix.
\end{theorem}
\begin{proof}[of Theorem~\ref{theorem:biconjugate-square-ldu}]
As per Lemma~\ref{lemma:Biconjugatable-in-Principal-Minors}, it is evident that the pair $(\bI_n, \bI_n)$ is biconjugatable. 
Based on Corollary~\ref{corollary:biconjugate-connection-u-x}, we have $\bU_n \bR_n^{(x)} = \bI_n$ and $\bI_n=\bV_n\bR_n^{(y)}$. Thus, $\bR_n^{(x)} = \bU_n^{-1}$ and $\bR_n^{(y)} = \bV_n^{-1}$ are well defined. This completes the proof.
\end{proof}

\subsection{Cholesky Decomposition}
For symmetric and positive definite matrices, the leading principal minors are always positive. 
The proof for this statement can be found in Section~\ref{appendix:leading-minors-pd}.
The following theorem describes the biconjugate decomposition for positive definite matrices.
\begin{theorem}[Cholesky: Biconjugate Decomposition for PD Matrices]\label{theorem:biconjugate-square-cholesky}
For any symmetric and positive definite matrix $\bA\in \real^{n\times n}$, the Cholesky decomposition of $\bA$ can be derived from the Wedderburn sequence by setting $(\bX_n, \bY_n)$ as the pair $(\bI_n, \bI_n)$.
In this case, $\bA$ can be factored as 
$$
\bA = \bU_n^{-\top} \bOmega_n \bU_n^{-1} = (\bU_n^{-\top} \bOmega_n^{1/2})( \bOmega_n^{1/2} \bU_n^{-1}) =\bR^\top\bR,
$$
where $\bOmega_n$ is a diagonal matrix with positive values along the diagonal, and $\bU_n^{-1}$ is a unit upper triangular matrix.
\end{theorem}
\begin{proof}[of Theorem~\ref{theorem:biconjugate-square-cholesky}]
Given that the leading principal minors of positive definite matrices are positive, $w_i>0$ for all $i\in \{1, 2,\ldots, n\}$. 
It follows from the LDU factorization via biconjugation and the symmetry of $\bA$ that $\bA = \bU_n^{-\top} \bOmega_n \bU_n^{-1}$.  
Since $w_i$'s are positive,  $\bOmega_n$ is positive definite and  can be factored as $\bOmega_n = \bOmega_n^{1/2}\bOmega_n^{1/2}$. This implies that $\bOmega_n^{1/2} \bU_n^{-1}$ is the Cholesky factor. 
\end{proof} 


\subsection{QR Decomposition}

 
Without loss of generality, we  assume that $\bA\in \real^{n\times n}$ has full rank, which allows for the QR decomposition: $\bA=\bQ\bR$, where $\bQ\in\real^{n\times n}$ is orthogonal, and $\bR \in \real^{n\times n}$ is upper triangular with  full rank and positive diagonal values. We then have the following result.
\begin{theorem}[QR: Biconjugate Decomposition for Nonsingular Matrices]\label{theorem:biconjugate-square-qr}
For any nonsingular matrix $\bA\in \real^{n\times n}$, the QR decomposition of $\bA$ can be obtained from the Wedderburn sequence by setting $(\bX_n, \bY_n)$ as $(\bI_n, \bA)$.
Thus, $\bA$ can be factored as 
$$
\bA = \bQ\bR,
$$
where $\bQ=\bV_n \bOmega_n^{-1/2}$ is an orthogonal matrix and $\bR = \bOmega_n^{1/2}\bR_n^{(x)}$ is an upper triangular matrix, according to the \textcolor{black}{\textbf{Form 4}} in Theorem~\ref{theorem:biconjugate-ldu}, with $\gamma=n$:
$$
\bY_n^\top\bA\bX_n = \bR_n^{(y)\top} \bV_n^\top \bA \bU_n\bR_n^{(x)}=\bR_n^{(y)\top} \bOmega_n \bR_n^{(x)},
$$
where we set $\gamma=n$ because $\gamma$ is any value such that $\gamma\leq r$ and the rank $r=n$.
\end{theorem}

\begin{proof}[of Theorem~\ref{theorem:biconjugate-square-qr}]
Since $(\bX_n, \bY_n) = (\bI_n, \bA)$, applying Theorem~\ref{theorem:biconjugate-ldu}, we have the decomposition
$$
\bY_n^\top\bA\bX_n = \bR_n^{(y)\top} \bV_n^\top \bA \bU_n\bR_n^{(x)}=\bR_n^{(y)\top} \bOmega_n \bR_n^{(x)}.
$$
Substituting $(\bX_n,\bY_n)=(\bI_n, \bA)$ into the  decomposition above, we obtain: 
\begin{equation}\label{equation:biconjugate-qr-ata1}
\begin{aligned}
	\bY_n^\top\bA\bX_n &= \bR_n^{(y)\top} \bV_n^\top \bA \bU_n\bR_n^{(x)} = \bR_n^{(y)\top} \bOmega_n \bR_n^{(x)};\\
	\bA^\top \bA &= \bR_n^{(y)\top} \bOmega_n \bR_n^{(x)}; \\
	\bA^\top \bA &= \bR_1^\top \bOmega_n \bR_1; \qquad (\text{$\bA^\top\bA$ is symmetric and let $\bR_1=\bR_n^{(x)}=\bR_n^{(y)}$})\\
	\bA^\top \bA &= (\bR_1^\top \bOmega_n^{1/2\top}) (\bOmega_n^{1/2}\bR_1);\\
	\bA^\top \bA &= \bR^\top\bR. \qquad \qquad\qquad\qquad(\text{let $\bR \triangleq \bOmega_n^{1/2}\bR_1$})
\end{aligned}
\end{equation}
To see why $\bOmega_n$ can be factored as $\bOmega_n = \bOmega_n^{1/2\top}\bOmega_n^{1/2}$, we consider the following steps. 
Suppose $\bA=[\ba_1, \ba_2, \ldots, \ba_n]$ is the column partition of $\bA$. We obtain $w_i = \by_i^\top\bA\bx_i=\ba_i^\top \ba_i>0$ since $\bA$ is nonsingular. Therefore, $\bOmega_n=\diag(w_1, w_2, \ldots, w_n)$ is positive definite and it can be factored as 
\begin{equation}\label{equation:omega-half-qr}
\bOmega_n = \bOmega_n^{1/2}\bOmega_n^{1/2}= \bOmega_n^{1/2\top}\bOmega_n^{1/2}.
\end{equation}
By $\bxgamma = \bugamma\brx$ in Theorem~\ref{theorem:biconjugate-ldu} for all $\gamma\in \{1, 2, \ldots, n\}$, we have 
$$
\begin{aligned}
\bX_n &= \bU_n\bR_1; \\
\bI_n &= \bU_n\bR_1; \qquad (\text{Since $\bX_n = \bI_n$}) \\
\bU_n &= \bR_1^{-1}.
\end{aligned}
$$
By $\bygamma = \bvgamma\bry$ in Theorem~\ref{theorem:biconjugate-ldu}  for all $\gamma\in \{1, 2, \ldots, n\}$, we have
\begin{equation}\label{equation:biconjugate-qr-ata2}
\begin{aligned}
\bY_n &= \bV_n\bR_1;\\
\bA &= \bV_n\bR_1; \qquad &(\text{$\bA=\bY_n$}) \\
\bA^\top \bA &= \bR_1^\top\bV_n^\top \bV_n\bR_1; \\
\bR_1^\top \bOmega_n \bR_1&=\bR_1^\top\bV_n^\top \bV_n\bR_1; \qquad &(\text{Equation~\eqref{equation:biconjugate-qr-ata1}})\\
(\bR_1^\top \bOmega_n^{1/2\top}) (\bOmega_n^{1/2}\bR_1) &= (\bR_1^\top \bOmega_n^{1/2\top} \bOmega_n^{-1/2\top})\bV_n^\top \bV_n (\bOmega_n^{-1/2}\bOmega_n^{1/2}\bR_1); \qquad &\text{(Equation~\eqref{equation:omega-half-qr})} \\
\bR^\top\bR &= \bR^\top (\bOmega_n^{-1/2\top} \bV_n^\top) (\bV_n \bOmega_n^{-1/2})  \bR.
\end{aligned}
\end{equation}
Thus, $\bQ=\bV_n \bOmega_n^{-1/2}$ is an orthogonal matrix.
\end{proof}

\subsection{SVD}\label{section:bicon_SVD}
To explore  the SVD of a square matrix $\bA\in\real^{n\times n}$ within the biconjugation decomposition, we introduce the following notation: let $\bA=\bU^\svd \bSigma^\svd \bV^{\svd\top}$ be the SVD of $\bA$, where $\bU^\svd = [\bu_1^\svd, \bu_2^\svd, \ldots, \bu_n^\svd]$ is orthogonal, $\bV^\svd = [\bv_1^\svd, \bv_2^\svd, \ldots, \bv_n^\svd]$ is orthogonal, and $\bSigma^\svd = \diag(\sigma_1, \sigma_2, \ldots, \sigma_n)$ is diagonal. Without loss of generality, we assume $\bA\in \real^{n\times n}$ and $\rank(\bA)=n$. Readers can verify the equivalence for a general matrix $\bA\in \real^{m\times n}$.

If the pair ($\bX_n=\bV^\svd$, $\bY_n=\bU^\svd$) effects a rank-reducing process for $\bA$.
From the definitions of $\bu_k$ and $\bv_k$ in Equation~\eqref{equation:properties-of-wedderburn-ukvk} or Equation~\eqref{equation:properties-of-wedderburn-ukvk2-inform2}, we have
$$
\bu_k = \bv_k^\svd \qquad \text{and} \qquad \bv_k = \bu_k^\svd \qquad \text{and} \qquad w_k = \by_k^\top \bA \bx_k=\sigma_k.
$$
This implies $\bV_n = \bU^\svd$, $\bU_n = \bV^\svd$, and $\bOmega_n = \bSigma^\svd$, where we set $\gamma=n$ because $\gamma$ can be any value such that $\gamma\leq r$ and the rank $r=n$.

By $\bX_n= \bU_n\bR_n^{(x)}$ in Theorem~\ref{theorem:biconjugate-ldu}, we have 
$$
\bX_n = \bU_n\bR_n^{(x)}   \quad\implies\quad
\bV^\svd = \bV^\svd\bR_n^{(x)}	 \quad\implies\quad
\bI_n = \bR_n^{(x)}.
$$
By $\bY_n = \bV_n\bR_n^{(y)}$ in Theorem~\ref{theorem:biconjugate-ldu}, we have
$$
\bY_n = \bV_n\bR_n^{(y)}    \quad\implies\quad
\bU^\svd = \bU^\svd\bR_n^{(y)} \quad\implies\quad
\bI_n=\bR_n^{(y)}.
$$
Applying Theorem~\ref{theorem:biconjugate-ldu} again and setting $\gamma=n$, we have
$$
\bY_n^\top\bA\bX_n = \bR_n^{(y)\top} \bV_n^\top \bA \bU_n\bR_n^{(x)}=\bR_n^{(y)\top} \bOmega_n \bR_n^{(x)}.
$$
This simplifies to
$$
\bU^{\svd\top}\bA \bV^\svd = \bSigma^\svd,
$$
which corresponds precisely to the form of a SVD. This demonstrates the equivalence between the SVD and the biconjugate decomposition when the Wedderburn sequence is applied with $(\bV^\svd, \bU^\svd)$ as $(\bX_n, \bY_n)$.

\section{Proof: General Term Formula of Wedderburn Sequence}\label{section:wedderburn-general-term}
In Lemma~\ref{lemma:wedderburn-sequence-general}, we present the general term formula of the Wedderburn sequence.
Given any matrix  $\bA\in\real^{m\times n}$, the Wedderburn sequence of $\bA$ is defined by $\bA_{k+1} = \bA_k -w_k^{-1} \bA_k\bx_k\by_k^\top \bA_k$ with $\bA_1 = \bA$. The proof of the general term formula for this sequence is as follows:
\begin{proof}[of Lemma~\ref{lemma:wedderburn-sequence-general}]
For $\bA_2$, let $\bu_1\triangleq\bx_1$ and $\bv_1 \triangleq\by_1$. We have:
$$
\begin{aligned}
\bA_2 &=\bA_1 -w_1^{-1} \bA_1\bx_1\by_1^\top \bA_1=\bA -w_1^{-1} \bA\bu_1\bv_1^\top \bA.
\end{aligned}
$$
For $\bA_3$, we can write out the equation as:
$$
\footnotesize
\begin{aligned}
&\bA_3 = \bA_2 -w_2^{-1} \bA_2\bx_2\by_2^\top \bA_2 \\
&=(\bA -w_1^{-1} \bA\bu_1\bv_1^\top \bA) - w_2^{-1}(\bA -w_1^{-1} \bA\bu_1\bv_1^\top \bA)\bx_2\by_2^\top(\bA -w_1^{-1} \bA\bu_1\bv_1^\top \bA) \qquad &\text{(substitute $\bA_2$)}\\
&=(\bA -w_1^{-1} \bA\bu_1\bv_1^\top \bA) 
- w_2^{-1}\textcolor{mylightbluetext}{\bA}(\textcolor{mylightbluetext}{\bx_2 }-w_1^{-1} \bu_1\bv_1^\top \bA\textcolor{mylightbluetext}{\bx_2 })(\textcolor{mylightbluetext}{\by_2^\top} -w_1^{-1}\textcolor{mylightbluetext}{\by_2^\top} \bA\bu_1\bv_1^\top )\textcolor{mylightbluetext}{\bA} \qquad &\text{(factor out $\bA$)}\\
&=\bA -w_1^{-1} \bA\bu_1\bv_1^\top \bA - w_2^{-1}\bA\bu_2\bv_2^\top\bA
=\bA -\sum_{i=1}^{2}w_i^{-1} \bA\bu_i\bv_i^\top \bA,
\end{aligned}
$$
where $\bu_2=\bx_2 -w_1^{-1} \bu_1\bv_1^\top \bA\bx_2=\bx_2 -\frac{\bv_1^\top \bA\bx_2}{w_1}\bu_1$, and $\bv_2=\by_2 -w_1^{-1}\by_2^\top \bA\bu_1\bv_1=\by_2 -\frac{\by_2^\top \bA\bu_1}{w_1}\bv_1$.
Similarly, we can find the expression of $\bA_4$ by $\bA$:
$$
\footnotesize
\begin{aligned}
&\bA_4  = \bA_3 -w_3^{-1} \bA_3\bx_3\by_3^\top \bA_3 \\
&=\bA -\sum_{i=1}^{2}w_i^{-1} \bA\bu_i\bv_i^\top \bA 
 - w_3^{-1} \big(\bA -\sum_{i=1}^{2}w_i^{-1} \bA\bu_i\bv_i^\top \bA\big)\bx_3\by_3^\top \big(\bA -\sum_{i=1}^{2}w_i^{-1} \bA\bu_i\bv_i^\top \bA\big) \gap &\text{{(substitute $\bA_3$)}}\\
&=\bA -\sum_{i=1}^{2}w_i^{-1} \bA\bu_i\bv_i^\top \bA 
- w_3^{-1}\textcolor{mylightbluetext}{\bA} \big(\textcolor{mylightbluetext}{\bx_3} -\sum_{i=1}^{2}w_i^{-1} \bu_i\bv_i^\top \bA\textcolor{mylightbluetext}{\bx_3}\big) \big(\textcolor{mylightbluetext}{\by_3^\top} -\sum_{i=1}^{2}w_i^{-1}\textcolor{mylightbluetext}{\by_3^\top} \bA\bu_i\bv_i^\top \big)\textcolor{mylightbluetext}{\bA} \gap &\text{(factor out $\bA$)} \\
&=\bA -\sum_{i=1}^{2}w_i^{-1} \bA\bu_i\bv_i^\top \bA - w_3^{-1}\bA\bu_3\bv_3^\top\bA
=\bA -\sum_{i=1}^{3}w_i^{-1} \bA\bu_i\bv_i^\top \bA,
\end{aligned}
$$	
where $\bu_3=\bx_3 -\sum_{i=1}^{2}\frac{\bv_i^\top \bA\bx_3}{w_i}\bu_i$, and $\bv_3=\by_3 -\sum_{i=1}^{2}\frac{\by_3^\top \bA\bu_i}{w_i}\bv_i$.
Continuing this process, we can define 
$$
\bu_k=\bx_k -\sum_{i=1}^{k-1}\frac{\bv_i^\top \bA\bx_k}{w_i}\bu_i
\qquad \text{and} \qquad
\bv_k=\by_k -\sum_{i=1}^{k-1}\frac{\by_k^\top \bA\bu_i}{w_i}\bv_i,
$$ 
and  the general term of the Wedderburn sequence can be proved by induction.
\end{proof}

\begin{problemset}
\item Following the proof of Theorem~\ref{theorem:rank-1-reduction}, prove Corollarys~\ref{corollary:rk_redu_one} and \ref{corollary:rk_K_redu_one}.

\item Discuss the Wedderburn sequence of $(\bX_n,\bY_n)=(\bI_n, \bA)$ in Theorem~\ref{theorem:biconjugate-square-qr} if $\bA$ is singular.

\item Following Section~\ref{section:bicon_SVD}, verify the equivalence between the SVD and the biconjugate decomposition for a general matrix $\bA$ of size $m\times n$.

\item \label{problem:rk_reduc} \textbf{Rank reduction theorem.} Let $\bA\in\real^{m\times n}$, $\bX\in\real^{n\times k}$, and $\bY\in\real^{m\times k}$ such that $\bW=\bY^\top\bA\bX$ is nonsingular. Show that
$$
\rank(\bA-\bA\bX\bW^{-1}\bY^\top\bA) = \rank(\bA)-\rank(\bA\bX\bW^{-1}\bY^\top\bA).
$$
When $k=1$, this is the rank-one reduction (Theorem~\ref{theorem:rank-1-reduction}).
Discuss its connection to Corollary~\ref{corollary:rk_K_redu_one}.

\item Show that if $\bA$ is symmetric and $\bX=\bY$, then $\bU=\bV$ in Lemma~\ref{lemma:wedderburn-sequence-general}.
\end{problemset}

%% file: chapter-tensornotation.tex
\newpage 
\part{Tensor Decomposition}\label{part:tensor_decom}

\chapter{Background and Notations}
\begingroup
\hypersetup{
	linkcolor=structurecolor,
	linktoc=page,  
}
\minitoc \newpage
\endgroup

\section{Matrices to Tensors}
\lettrine{\color{caligraphcolor}W}
We have previously explored matrix decomposition and approximation, where the primary objects of interest were vectors and matrices. 
The utilization  of tensors in data science was pioneered by researchers in psychometrics and chemometrics during the 1960s \citep{smilde2005multi, kolda2006multilinear}. These early contributions laid the foundation for the modern application of tensors in a wide range of fields.

In today's era, the volume and structural complexity of modern datasets are increasing to a point where conventional  analysis methods and algorithms are no longer sufficient. 
Such massive datasets can contain billions of entries and are often represented as large block matrices or tensors. This has sparked renewed interest in developing matrix and tensor algorithms that are suitable for very large-scale datasets.
To address this challenge, tensor decomposition offers a natural, sparse, and distributed representation for big data. Tensor decomposition encompasses both established and emerging methodologies for tensor-based representations and optimization.
It provides a theoretical and computational framework for analyzing computationally prohibitive large volumes of data by separating them into relevant and irrelevant information, both of which have lower and reduced dimensionality. This approach often allows for the super-compression of datasets with up to $10^{50}$ entries down to more manageable sizes of $10^{7}$ or fewer entries \citep{cichocki2016tensor}.
Tensor approximation can be visualized as special graph structures that decompose high-order tensors into a set of sparsely interconnected low-order core tensors or matrices. This decomposition improves both interpretability and computational efficiency, making tensor decomposition a potent tool for managing the complexity and scale of modern datasets.

\index{Curse of dimensionality}
On the other hand, the term \textit{``curse of dimensionality"} was originally coined to describe the exponential increase in the number of samples required to estimate an arbitrary function with a specific accuracy as the number of variables (or dimensions) grows. In the broader context of machine learning and optimization, this term also captures the exponential rise in parameters required to model data or systems, along with the overwhelming number of degrees of freedom that can arise in high-dimensional spaces \citep{bishop2006pattern}.
In tensor contexts, the \textit{``curse of dimensionality"} refers specifically to the rapid growth in the number of elements, $I^N$, in an $N$-order tensor of size $(I \times I \times \cdots \times I)$. As the tensor order $N$ increases, this growth exponentially rendering the volume of multiway arrays impractically large, placing significant demands on computational resources and memory. (Even for $I=2$, storing a tensor of order $N=50$ would require 9 petabytes of storage). Processing high-dimensional tensors involves addressing dependencies across an excessive number of degrees of freedom, which presents both substantial challenges and intriguing opportunities. However, such processing often necessitates approximations that may compromise  accuracy \citep{oseledets2011tensor, grasedyck2013literature}.

To mitigate these issues, \textit{low-rank tensor approximations (LRTAs)} offer a promising approach. These approximations help manage the scale, accuracy, and complexity of high-dimensional data by leveraging efficient tensor decomposition algorithms that rely on relatively simple mathematical operations. A notable advantage of LRTAs is their use of compact, low-order core tensors or matrices, which act as manageable ``building blocks" of the tensor. These core tensors facilitate  super-compression of large-scale, incomplete, and noisy datasets, making the data more manageable and easier to visualize.

In this part, our focus is on LRTAs, which enables the compression of large data tensors using interconnected low-order core tensors or matrices so as to offer a powerful method for addressing the curse of dimensionality, provide a feasible solution to the broader challenge of managing exceedingly large datasets within practical computation limits,  and enable significant data compression and efficient processing.
In conclusion, the advantages of LRTAs are twofold:
(i) They provide a tensor approximation of a specific cost (objective) function while preserving the required level of accuracy; (ii) They derive physically significant latent variables from data using a method that is both sufficiently precise and computationally efficient.
In this section, we will provide a brief overview of the foundational concepts for tensor analysis, following the notational conventions largely from the reviews of \citet{kolda2006multilinear, kolda2009tensor, cichocki2016tensor, golub2013matrix} and \citet{favier2014overview}.
A tensor is essentially a multidimensional array; when the dimension is $N$, we call the tensor an \textit{Nth-order tensor}. For example, a first-order tensor corresponds to a vector, while a second-order tensor corresponds to a matrix. 
In this sense, tensors are multidimensional extensions of matrices, which are used to represent ubiquitous multidimensional data, such as RGB images, hyperspectral images, and videos.

Previously, we use \textbf{boldface} lowercase letters possibly with subscripts to denote vectors (e.g., $\bmu$, $\bx$, $\bx_n$, $\bz$) and
\textbf{boldface} uppercase letters possibly with subscripts to denote matrices (e.g., $\bA$, $\bL_j$). To avoid confusion, we will use a \textbf{boldface} Euler script letter to denote a tensor with an order larger than 2, e.g., $\eX, \eY,\eA,\eB_i,\eG_j \in \real^{I\times J\times K}$ will represent third-order tensors.

\begin{figure}[h]
	\centering  
	\vspace{-0.35cm}
	\subfigtopskip=2pt 
	\subfigbottomskip=2pt 
	\subfigcapskip=-5pt 
	\subfigure[2nd-order tensor.]{\label{fig:tensor-2dd}%
		\includegraphics[width=0.33\linewidth]{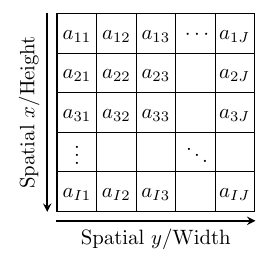}}%
	\subfigure[3rd-order tensor.]{\label{fig:tensor-3dd}%
		\includegraphics[width=0.5\linewidth]{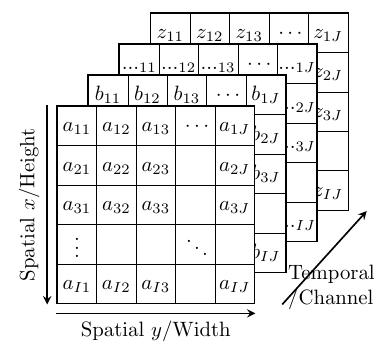}}%
	\caption{Transition from a 2nd-order tensor to a 3rd-order tensor.}
	\label{fig:tensor-2d-2-3d}
\end{figure}

In a second-order tensor (i.e., a matrix), the first axis is typically taken as the \textbf{coordinate of spatial $x$} or \textbf{height}, and the second axis is then taken as the \textbf{coordinate of spatial $y$} or \textbf{width}. 
For example, Figure~\ref{fig:tensor-2dd} is an example of a second-order tensor $\bX\in \real^{I\times J}$, in which case its column vectors are called the \textit{mode-1 vectors} or \textit{mode-1 fibers}, and its rows are called the \textit{mode-2 vectors} or \textit{mode-2 fibers}. 
A third-order tensor extends this concept by adding a third dimension, which is often  referred to as a \textbf{``temporal" coordinate} (or a \textbf{``channel" coordinate} in an RGB picture). The arrangement of a third-order tensor $\eX\in \real^{I\times J\times K}$ is shown in Figure~\ref{fig:tensor-3dd}, in which case $K=3$ in the case of an RGB image. 
Similarly, the vectors along the channel direction of $\eX$ are called the \textit{mode-3 vectors} or \textit{mode-3 fibers}.

More generally, an \textit{Nth-order tensor} is denoted by $\eX\in \real^{I_1\times I_2\times \ldots \times I_N}$, where the $n$-th size of the tensor is $I_n$. The distinction  between the \textit{orders} (or \textit{ways}) and the \textit{sizes} (or \textit{dimensions}) of the tensor is that $N$ represents the order of the tensor, while $I_n$ represents the size of the \textit{mode-$n$ fibers}. That is, the term ``dimensions/sizes" stands for
the number of values that an index can take in a particular mode.

\section{Tensor Indexing}\index{Tensor indexing}
\textit{Slices} (\textit{matrix slices}) refer to two-dimensional sections of a tensor, defined by fixing all but two indices. Figure~\ref{fig:thifd-slices} shows the slices of a third-order tensor $\eX\in \real^{I\times J\times K}$, where the \textit{horizontal ones} are the slices by fixing the first index and varying the last two: $\eX_{i,:,:}$ for all $i\in \{1,2,\ldots, I\}$ \footnote{where the indices typically range from 1 to their capital version, here, $i\in \{1,2,\ldots, I\}$. That is, lower-case letters e.g, $i, j, k$, are used for the subscripts in running indices, while capital letters $I, J, K$ denote the upper bound of an index.}; the \textit{lateral ones} are the slices by fixing the second index: $\eX_{:,j,:}$ for all $j\in \{1,2,\ldots, J\}$; and the \textit{frontal ones} are the slices by fixing the third index: $\eX_{:,:,k}$ for all $k\in \{1,2,\ldots, K\}$.

\index{Tensor slices}
\index{Tensor fibers}

Similarly, \textit{fibers} are the higher-order analogues of matrix rows and columns. A fiber is defined by fixing all indices except one. Figure~\ref{fig:thifd-fibers} shows the fibers of a third-order tensor $\eX\in \real^{I\times J\times K}$, where the \textit{column fibers} (\textit{mode-1 fibers}) are the ones by varying only the first index: $\eX_{:,j,k}$; the \textit{row fibers} (\textit{mode-2 fibers}) are the ones by varying only the second index: $\eX_{i,:,k}$; and the \textit{tube fibers} (\textit{mode-3 fibers}) are the ones by varying only the third index: $\eX_{i,j,:}$.

\begin{figure}[h]
	\centering
	\vspace{-0.35cm}
	\subfigtopskip=2pt
	\subfigbottomskip=2pt
	\subfigcapskip=-5pt
	\subfigure[Horizontal slices $\eX_{i,:,:}$]{\label{fig:third-slice1}%
		\centering
		\begin{minipage}[b]{0.32\linewidth}
			\centering
   
		\end{minipage}%
	}
\caption{Slices of a third-order tensor $\eX\in \real^{I\times J\times K}$. The \textbf{darker}, the larger the index. Clearly, the slices of an Nth-order tensor are $(N-1)$-th order tensors, and the slices of a third-order tensor are matrices.}
\label{fig:thifd-slices}
\end{figure}

\begin{figure}[h]
	\centering
	\vspace{-0.35cm}
	\subfigtopskip=2pt
	\subfigbottomskip=2pt
	\subfigcapskip=-5pt
	\subfigure[Mode-1 fibers: columns, $\eX_{:,j,k}$]{\label{fig:third-fibdr1}%
		\centering
		\begin{minipage}[b]{0.32\linewidth}
			\centering
			%
     
		\end{minipage}%
	}
	
	\caption{Fibers of a third-order tensor $\eX\in \real^{I\times J\times K}$.}
	\label{fig:thifd-fibers}
\end{figure}

More generally, the $(i,j,k)$-th element of a third-order tensor $\eX\in \real^{I\times J\times K}$ is indexed by $x_{ijk}$. And the general index for an Nth-order tensor can be inferred  from the context.

\section{Inner Product and Frobernius Norm}
The \textit{inner product} of two Nth-order tensors with the same size, denoted by $\eX,\eY\in \real^{I_1\times I_2\times \ldots \times I_N}$, is the sum of the element-wise products of their entries:
$$
\langle \eX,\eY\rangle = \sum_{i_1}^{I_1}\sum_{i_2}^{I_2}\ldots \sum_{i_N}^{I_N} (x_{i_1,i_2,\ldots, i_N})\cdot (y_{i_1,i_2,\ldots, i_N}).
$$
Similarly, the Frobenius norm of an Nth-order tensor $\eX\in \real^{I_1\times I_2\times \ldots \times I_N}$ is given by $\sqrt{\langle \eX,\eX\rangle}$:
$$
\norm{\eX}_F =\sqrt{\langle \eX,\eX\rangle}
= \sqrt{\sum_{i_1}^{I_1}\sum_{i_2}^{I_2}\ldots \sum_{i_N}^{I_N} (x_{i_1,i_2,\ldots, i_N})^2}.
$$
Again, we will denote $\norm{\eX}_F$ simply as $\norm{\eX}$ without any further specification.

\section{Outer Product and Rank-One Tensor}\label{section:rank-one-tensor}\index{Rank-one tensor}
Given two vectors $\ba\in \real^I$ and $\bb\in \real^J$, the outer product is $\ba\bb^\top \in \real^{I\times J}$ in a basic  vector context for two vectors. Analogously, in the tensor language, we use the symbol $``\circ"$ to denote the outer product, e.g., $\ba\bb^\top = \ba \circ\bb$ denotes the outer product of two vectors.
\begin{definition}[Nth-Order Vector Outer Product, Decomposed Tensor Set]
In higher dimensions, the outer product of $N$ vectors $\ba^{(1)}\in \real^{I_1}, \ba^{(2)}\in \real^{I_2}, \ldots, \ba^{(N)}\in \real^{I_N}$ is given by 
$$
\ba^{(1)}\circ \ba^{(2)}\circ \ldots \circ \ba^{(N)}
\in\real^{I_1\times I_2\times \ldots \times I_N},
$$
where the $(i_1,i_2,\ldots, i_N)$-th element can be obtained by 
$$
(\ba^{(1)}\circ \ba^{(2)}\circ \ldots \circ \ba^{(N)})_{i_1,i_2,\ldots, i_N}=
a^{(1)}_{i_1}a^{(2)}_{i_2}\ldots a^{(N)}_{i_N}.
$$
Specifically, any tensor $\eX \in  \real^{I_1\times I_2\times \ldots \times I_N}$ that can be expressed as the outer product of $N$ vectors is called \textit{decomposed}.
And the set of all $I_1\times I_2\times \ldots \times I_N$ decomposed tensors is called the \textit{decomposed tensor set} and is denoted by $\sD(I_1, I_2, \ldots , I_N)$ or simply $\sD$.
\end{definition}

In the meantime, an Nth-order tensor $\eX \in  \real^{I_1\times I_2\times \ldots \times I_N}$ is of \textit{rank-one} if it
can be expressed as the outer product of $N$ vectors, i.e.,
\begin{equation}\label{equation:tensor_outer_decom}
\eX = \ba^{(1)}\circ \ba^{(2)}\circ \ldots \circ \ba^{(N)},
\end{equation}
The $(i_1,i_2,\ldots, i_N)$-th element of $\eX$ is thus given by 
$
\eX_{i_1,i_2,\ldots, i_N} = x_{i_1,i_2,\ldots, i_N} = a_{i_1}^{(1)}a_{i_2}^{(2)}\ldots a_{i_N}^{(N)}.
$
And the representation for a third-order rank-one tensor is illustrated in Figure~\ref{fig:rank-one-tensor}.
\paragraph{Expressing a tensor using outer products.}
For a matrix $\bA\in\real^{m\times n}$, it can be equivalently expressed as 
$
\bA = \sum_{i=1}^{m}\sum_{j=1}^{n} a_{ij}\be_i\be_j^\top,
$
where $\be_i$ and $\be_j$ are the standard basis vectors in $\real^m$ and $\real^n$, respectively.
Similarly, a general tensor $\eX \in  \real^{I_1\times I_2\times \ldots \times I_N}$ can be equivalently expressed using high-order vector outer products:
\begin{equation}
\eX =\sum_{i_1=1}^{I_1}\ldots \sum_{i_N=1}^{I_N}x_{i_1,i_2,\ldots,i_N} 
\be_{i_1}^{(I_1)}\circ \be_{i_2}^{(I_2)}\circ \ldots \circ \be_{i_N}^{(I_N)},
\end{equation}
where $\be_{i_n}^{(I_n)}\in\real^{I_n}$ is the $i_n$-th basis vector in $\real^{I_n}$.

\paragraph{Inner product and norm of outer product.}
Given a set of vectors $\ba^{(1)}, \bb^{(1)}\in \real^{I_1}, \ba^{(2)}, \bb^{(2)}\in \real^{I_2}, \ldots, \ba^{(N)}, \bb^{(N)}\in \real^{I_N}$. 
If $\eX=\ba^{(1)}\circ \ba^{(2)}\circ \ldots \circ \ba^{(N)}$ 
and 
$\eY=\bb^{(1)}\circ \bb^{(2)}\circ \ldots \circ \bb^{(N)}$, then it follows that 
\begin{itemize}
\item 	$\langle \eX, \eY \rangle
=
\prod_{n=1}^{N} 
\langle \ba^{(n)}, \bb^{(n)} \rangle$.

\item  $\norm{\eX}_F = \prod_{n=1}^{N} \norm{\ba^{(n)}}_2$ 
and $\norm{\eY}_F = \prod_{n=1}^{N} \norm{\bb^{(n)}}_2$.

\item \textbf{Orthogonality.} Assume further $\norm{\eX}_F=\norm{\eY}_F=1$, then $\eX$ and $\eY$ are said to be orthogonal if their inner product is zero: $\langle \eX, \eY\rangle =
\prod_{n=1}^{N} 
\langle \ba^{(n)}, \bb^{(n)} \rangle = 0$.
\end{itemize}

\begin{figure}[htbp]
\centering
\resizebox{0.65\textwidth}{!}{%
\begin{tikzpicture}
\draw [very thick] (5.8+0.5,-3.8+0.5) rectangle (7+0.5,-2.6+0.5);
\filldraw [fill=gray!60!white,draw=green!40!black] (5.8+0.5,-3.8+0.5) rectangle (7+0.5,-2.6+0.5);
\draw [step=0.4/2, very thin, color=gray] (5.8+0.5,-3.8+0.5) grid (7+0.5,-2.6+0.5);

\draw [very thick] (5.8+0.4,-3.8+0.4) rectangle (7+0.4,-2.6+0.4);
\filldraw [fill=gray!50!white,draw=green!40!black] (5.8+0.4,-3.8+0.4) rectangle (7+0.4,-2.6+0.4);
\draw [step=0.4/2, very thin, color=gray] (5.8+0.4,-3.8+0.4) grid (7+0.4,-2.6+0.4);

\draw [very thick] (5.8+0.3,-3.8+0.3) rectangle (7+0.3,-2.6+0.3);
\filldraw [fill=gray!40!white,draw=green!40!black] (5.8+0.3,-3.8+0.3) rectangle (7+0.3,-2.6+0.3);
\draw [step=0.4/2, very thin, color=gray] (5.8+0.3,-3.8+0.3) grid (7+0.3,-2.6+0.3);

\draw [very thick] (5.8+0.2,-3.8+0.2) rectangle (7+0.2,-2.6+0.2);
\filldraw [fill=gray!30!white,draw=green!40!black] (5.8+0.2,-3.8+0.2) rectangle (7+0.2,-2.6+0.2);
\draw [step=0.4/2, very thin, color=gray] (5.8+0.2,-3.8+0.2) grid (7+0.2,-2.6+0.2);

\draw [very thick] (5.8+0.1,-3.8+0.1) rectangle (7+0.1,-2.6+0.1);
\filldraw [fill=gray!20!white,draw=green!40!black] (5.8+0.1,-3.8+0.1) rectangle (7+0.1,-2.6+0.1);
\draw [step=0.4/2, very thin, color=gray] (5.8+0.1,-3.8+0.1) grid (7+0.1,-2.6+0.1);

\draw [very thick] (5.8,-3.8) rectangle (7,-2.6);
\filldraw [fill=gray!10!white,draw=green!40!black] (5.8,-3.8) rectangle (7,-2.6);
\draw [step=0.4/2, very thin, color=gray] (5.8,-3.8) grid (7,-2.6);

\draw [very thick] (5.8-0.1,-3.8-0.1) rectangle (7-0.1,-2.6-0.1);
\filldraw [fill=gray!6!white,draw=green!40!black] (5.8-0.1,-3.8-0.1) rectangle (7-0.1,-2.6-0.1);
\draw [step=0.4/2, very thin, color=gray] (5.8-0.1,-3.8-0.1) grid (7-0.1,-2.6-0.1);

\draw[->] (5.6,-2.6) -- ++(0.6,0.6); 
\draw[->] (5.55,-2.7) -- ++(0,-1.2);  
\draw[->] (5.66,-4.04) -- ++(1.25,0);  

\draw (6.35,-4.25) node {{\color{black}\scriptsize{spatial $y$}}};
\draw (5.2-0.3,-3.2) node {{\color{black}\scriptsize{spatial $x$}}};
\draw (5.3,-2.2) node {{\color{black}\scriptsize{temporal}}};

\draw (6.2,-4.6) node {{\color{black}\scriptsize{$\eX\in\real^{I\times J \times K}$}}};

\draw (8,-3.2) node {{\color{black}\large{$=$}}};

\draw [very thick] (8.6,-4.4+0.2) rectangle (8.8,-3.2+0.2);
\filldraw [fill=WildStrawberry!40!white,draw=green!40!black] (8.6,-4.4+0.2) rectangle (8.8,-3.2+0.2);
\draw (8.7,-4.8+0.2) node {{\color{black}\scriptsize{$\ba \in \real^I$}}};

\draw [very thick] (9,-3-0.2+0.2) rectangle (10.2,-2.8-0.2+0.2);
\filldraw [fill=RubineRed!60!white,draw=green!40!black] (9,-3-0.2+0.2) rectangle (10.2,-2.8-0.2+0.2);
\draw (9.6,-3.3-0.2+0.2) node {{\color{black}\scriptsize{$\bb \in \real^J$}}};

\draw[fill=RedOrange!40!white, line width=0.8pt] (9.2,-2.4+0.2) -- (9.4,-2.4+0.2) -- (8.8,-3+0.2) -- (8.6,-3+0.2) -- cycle;
\draw (10-0.2,-1.9) node {{\color{black}\scriptsize{$\bc\in \real^K$}}};
\end{tikzpicture}
}
\caption{A third-order tensor with rank-one, $\eX=\ba\circ\bb\circ\bc$, where the $(i, j, k)$-th element of $\eX$ is given by $x_{ijk} = a_ib_j c_k$.}
\label{fig:rank-one-tensor}
\end{figure}

The outer product of tensors is defined in a similar way.
\begin{definition}[Nth-Order Tensor Outer Product]
Given  tensors $\eX\in\real^{I_1\times I_2\times \ldots \times I_P}$ and $\eY\in\real^{J_1\times J_2\times \ldots \times J_Q}$, the tensor outer product is given by 
$$
\eX\circ \eY \in \real^{I_1\times I_2\times \ldots \times I_P\times J_1\times J_2\times \ldots \times J_Q},
$$
where the $(i_1, i_2, \ldots , i_P, j_1, j_2, \ldots , j_Q)$-th element is given by 
$$
(\eX\circ \eY)_{i_1, i_2, \ldots , i_P, j_1, j_2, \ldots , j_Q}
=
x_{i_1, i_2, \ldots , i_P } 
\cdot 
y_{j_1, j_2, \ldots , j_Q}.
$$
\end{definition}

\section{Diagonal and Identity Tensors}
Diagonal matrices and identity matrices have their counterparts in tensor language.
\begin{definition}[Diagonal and Identity Tensors]\label{definition:identity-tensors}
A tensor $\eX \in  \real^{I_1\times I_2\times \ldots \times I_N}$ is called a \textit{diagonal tensor} if $x_{i_1i_2\ldots i_N}\neq 0$ if and only if $i_1=i_2=\ldots=i_N$. And when $x_{i_1i_2\ldots i_N}=1$ if and only if $i_1=i_2=\ldots=i_N$, the tensor is known as the \textit{Nth-order identity tensor}: 
$\eX=\sum_{i=1}^{I}\be_{i}^{(I)}\circ \be_{i}^{(I)}\circ \ldots \circ \be_{i}^{(I)}$ is the identity tensor if all the dimensions are equal to $I$, where $\be_i^{(I)}$ is the $i$-th basis vector in $\real^I$.
\end{definition}

\begin{figure}[h]
	\centering
	\includegraphics[width=0.5\textwidth]{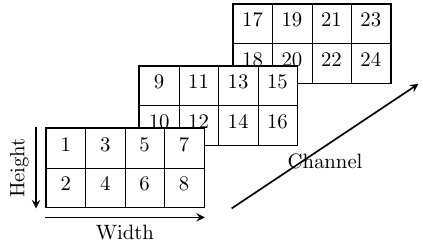}
	\caption{An example of a 3rd-order tensor: $\eY\in \real^{2\times 4\times 3}$, where $y_{ijk}=1+(i-1)+2(j-1)+8(k-1)$.}
	\label{fig:3rd-example}
\end{figure}

\section{Mode-$n$ Rank and Tensor Rank}
In second-order tensors, i.e., matrices, the column rank is equal to the row rank \footnote{The column rank of a matrix is equal to the dimension of the column space, the row rank of a matrix is equal to the dimension of the row space.}. 
However, this equivalence does not necessarily hold for higher-order tensors.
\begin{definition}[Mode-$n$ Rank of Tensor]
Given a tensor $\eX\in \real^{I_1\times I_2\times \ldots \times I_N}$, the \textit{mode-$n$ rank} of the tensor is the dimension 
$$
\rank_n(\eX) = \dim\left(
\spn\left\{\eX_{i_1, \ldots, i_{n-1},:, i_{n+1},\ldots, i_N}|
i_1\in\{1,2,\ldots,I_1\}, \ldots, i_N\in\{1,2,\ldots,I_N\} \right\}
\right),
$$
where $\rank_n(\eX)\leq I_n$.
\end{definition}
In tensors, the rank can vary along different dimensions, i.e., the mode-$m$ and the mode-$n$ ranks of a tensor are generally not the same if $m\neq n$.

However, many researchers attempt  to define the unique rank for a tensor, i.e., the \textit{tensor rank} \citep{kruskal1977three}.
\begin{definition}[Tensor Rank]
Given a tensor $\eX\in \real^{I_1\times I_2\times \ldots \times I_N}$, then it has a tensor rank of $R$ if $R$ is the smallest  value such that $\eX$ can be decomposed into a weighted sum of $R$ decomposed tensors
$$
\eX = \sum_{r=1}^{R} \gamma_r \eA_r,
$$
where $\gamma_r>0$, $\norm{\eA_r}_F=1$, and $\eA_r \in\sD(I_1, I_2, \ldots , I_N)$ (i.e., $\eA_r$ can be decomposed into the outer product of $N$ vectors). And the rank is denoted by 
$$
\rank(\eX) = R.
$$

\end{definition}

\section{Matricization: Matrix Representation of High-Order Tensors}\label{section:matricization-tensor-original}
In tensor analysis, it is desirable to represent a tensor as a matrix,  in which case the tensor is rearranged into a matrix; the process is known as  \textit{tensor matricization}.
There are several types of matrix representations of an Nth-order tensor $\eX\in \real^{I_1\times I_2\times \ldots \times I_N}$. 
However, for the matrix representation along the $n$-th mode, the size will always be $I_n\times (I_1\ldots I_{n-1} I_{n+1} \ldots  I_N)$, which is also called the \textbf{mode-$n$ matricization of
the tensor} $\eX$ and is denoted by 
$$
\mathcalM_{(n)}(\eX)\triangleq
\bX_{(n)}\in\real^{I_n\times (I_1\ldots I_{n-1} I_{n+1} \ldots  I_N)}.~\footnote{Note that we use the ``$\times$" symbol to denote a delimiter of a tensor; a scalar product that is not represented by a ``$\times$" symbol  indicates that the dimension of the mode is the product of a set of scalars.}
$$ 

To illustrate, we consider a specific example, $\eY\in \real^{2\times 4\times 3}$, where the entries are filled with values from $\{1,2,\ldots, 24\}$. 
In this arrangement, each number is stored in ascending order, starting from the first index (height), then the second index (width), and finally the third index (channel), as shown in Figure~\ref{fig:3rd-example}.
Mathematically, each entry of $\eY$ can be described by 
$$
y_{ijk} = 1+(i-1)+2(j-1)+8(k-1),
$$
where $i\in\{1,2\}$, $j\in\{1,2,3,4\}$, and $k\in\{1,2,3\}$.

\subsection{Kolda Matricization}\label{section:kolda_matric}
The primary matricization used in this book is called the \textit{Kolda matricization method} \citep{kolda2006multilinear}.
\paragraph{Mode-1 matricization.} Suppose now  we want to obtain the matricized form along the 1st-mode, $\bY_{(1)}\in \real^{2\times 12}$, the elements will be \textbf{\textcolor{mylightbluetext}{fetched}} first from the 1-st index of $\eY$ (height), followed by the 2-nd index (width), and finally the 3-rd index (channel). And these elements are then \textbf{\textcolor{mylightbluetext}{stored}} into $\bY_{(1)}$ starting from the first index, and then the second:
$$
\bY_{(1)}=
\begin{bmatrix}
1&3 & 5& 7&9 & 11& 13& 15& 17& 19 \,\,\,\, 21\,\,\,\,23\\
2&4 & 6& 8&10 & 12& 14& 16& 18& 20\,\,\,\, 22\,\,\,\,24\\
\end{bmatrix}.
$$
For a general tensor $\eX\in\real^{I\times J\times K}$, the mode-1 matricization of $\eX$ can be represented as 
\begin{equation}\label{equation:kold_mode1}
\eX  \xrightarrow{\text{matricization}}\bX_{(1)} = 
\begin{bmatrix}
\underset{I\times J}{\eX_{:,:,1}}, & \underset{I\times J}{\eX_{:,:,2}}, &\ldots, & \underset{I\times J}{\eX_{:,:,K}}
\end{bmatrix}
\in\real^{I\times (JK)}
\end{equation}
\paragraph{Mode-2 matricization.} If we want to obtain the mode-2 matricization of tensor $\eY$, i.e., $\bY_{(2)}\in \real^{4\times 6}$, the elements will be \textbf{\textcolor{mylightbluetext}{fetched}} first from the 2-nd index (width), followed by the 1-st index (height), and finally the 3-rd index (channel). \textbf{The storing of the elements will always be the same:} the elements are \textbf{\textcolor{mylightbluetext}{stored}} into $\bY_{(2)}$ starting from the first index, and then the second:
$$
\bY_{(2)}=
\begin{bmatrix}
1& 2& 9& 10& 17& 18 \\
3& 4& 11& 12& 19& 20\\
5& 6& 13& 14& 21& 22\\
7& 8& 15& 16& 23& 24\\
\end{bmatrix}.
$$
For a general tensor $\eX\in\real^{I\times J\times K}$, the mode-2 matricization of $\eX$ can be represented as 
\begin{equation}\label{equation:kold_mode2}
\eX \xrightarrow{\text{matricization}} \bX_{(2)} = 
\begin{bmatrix}
\underset{J\times I}{\eX_{:,:,1}^{\textcolor{black}{\top}}}, 
& \underset{J\times I}{\eX_{:,:,2}^{\textcolor{black}{\top}}}, 
&\ldots, 
& \underset{J\times I}{\eX_{:,:,K}^{\textcolor{black}{\top}}}
\end{bmatrix}
\in\real^{J\times (IK)}
\end{equation}
\paragraph{Mode-3 matricization.} If we want to obtain the mode-3 matricization of tensor $\eY$, i.e., $\bY_{(3)}\in \real^{3\times 8}$, the elements will be \textbf{\textcolor{mylightbluetext}{fetched}} first from the  3-rd index (channel), followed by the 1-st index (height), and finally the 2-nd index (width). \textbf{The storing of the elements is still the same}: the elements are stored into $\bY_{(3)}$ starting from the first index, and then the second:
$$
\bY_{(3)}=
\begin{bmatrix}
1& 2& 3&    4&  5&  6&  7& 8\\
9& 10& 11&  12& 13& 14& 15& 16 \\
17& 18& 19& 20& 21& 22& 23&  24\\
\end{bmatrix}.
$$
For a general tensor $\eX\in\real^{I\times J\times K}$, the mode-3 matricization of $\eX$ can be represented as 
\begin{equation}\label{equation:kold_mode3}
\eX \xrightarrow{\text{matricization}} \bX_{(3)} = 
\begin{bmatrix}
\underset{K\times I}{\eX^\top_{:,1,:}}, & \underset{K\times I}{\eX^\top_{:,2,:}}, &\ldots, & \underset{K\times I}{\eX^\top_{:,J,:}}
\end{bmatrix}
\in\real^{K\times (IJ)}
\end{equation}

To summarize,  the following properties are derived from the matricization process (for a third-order tensor):
\begin{itemize}
\item The column vectors of $\eY$ are column vectors of $\bY_{(1)}$;
\item The row vectors  of $\eY$ are column vectors of $\bY_{(2)}$;
\item The tube vectors  of $\eY$ are column vectors of $\bY_{(3)}$.
\end{itemize}

More generally, for an Nth-order tensor $\eX\in \real^{I_1\times I_2\times \ldots \times I_N}$, in the mode-$n$ matricization, the  $(i_1,i_2,\ldots, i_N)$-th tensor element is mapped into the following  $(i_n, j)$-th matrix element:
\begin{equation}
j= 1+\sum_{k=1,k\neq n}^{N} (i_k-1) J_k, \qquad \text{where }\gap J_k = \prod_{m=1,m\neq n}^{k-1} I_m.
\end{equation}

\subsection{Kiers Matricization}
The Kolda matricization method has a straightforward interpretation: for finding the mode-$n$ matricization, the method first fetches elements from the $n$-th order, and then fetches from other orders in  ascending order $\{1,2,(n-1), (n+1), \ldots, N\}$.
However, the closed-form representations of the matricization are not consistent, e.g., the fixing indexing in Equation~\eqref{equation:kold_mode1} and Equation~\eqref{equation:kold_mode2} is the 3-rd mode, whereas the fixing indexing   Equation~\eqref{equation:kold_mode3} is the 2-nd mode.
The \textit{Kiers matricization method} addresses this inconsistency \citep{kiers2000towards}.

For a general tensor $\eX\in\real^{I\times J\times K}$, the mode-1, mode-2, and mode-3  Kiers matricizations of $\eX$ can be represented as:
\begin{align}\label{equation:kieres_mode1}
\eX  \xrightarrow{\text{matricization}}
&\bX_{(1)} = 
\begin{bmatrix}
	\underset{I\times J}{\eX_{:,:,1}}, & \underset{I\times J}{\eX_{:,:,2}}, &\ldots, & \underset{I\times J}{\eX_{:,:,K}}
\end{bmatrix}
\in\real^{I\times (JK)};\\
&\bX_{(2)} = 
\begin{bmatrix}
\underset{J\times K}{\eX_{1,:,:}}, & \underset{J\times K}{\eX_{2,:,:}}, &\ldots, & \underset{J\times K}{\eX_{I,:,:}}
\end{bmatrix}
\in\real^{J\times (IK)};\\
&\bX_{(3)} = 
\begin{bmatrix}
	\underset{K\times I}{\eX^\top_{:,1,:}}, & \underset{K\times I}{\eX^\top_{:,2,:}}, &\ldots, & \underset{K\times I}{\eX^\top_{:,J,:}}
\end{bmatrix}
\in\real^{K\times (IJ)},
\end{align}
where we find the mode-1 and mode-3 Kiers matricizations are the same as the Kolda matricizations.

More generally, given the Nth-order tensor $\eX\in \real^{I_1\times I_2\times \ldots \times I_N}$, in the mode-$n$ matricization, the  $(i_1,i_2,\ldots, i_N)$-th tensor element is mapped into the following   $(i_n, j)$-th matrix element:
\begin{equation}
j= \sum_{p=1}^{N-2}
\left(
(i_{N+n-p}-1) \prod_{q=n+1}^{N+n-p-1} I_q
\right) 
+i_{n+1}, \gaps \text{where }\gaps I_{N+m}=I_m \text{ and } i_{N+m}=i_m. 
\end{equation}

\subsection{LMV Matricization}
The \textit{LMV matricization method}  works similarly to the Kiers matricization method \citep{de2000multilinear}.
For a general tensor $\eX\in\real^{I\times J\times K}$, the mode-1, mode-2, and mode-3  Kiers matricizations of $\eX$ can be represented as:
\begin{align}\label{equation:LMV_mode1}
\eX  \xrightarrow{\text{matricization}}
&\bX_{(1)} = 
\begin{bmatrix}
\underset{I\times K}{\eX_{:,1,:}}, 
& \underset{I\times K}{\eX_{:,2,:}}, 
&\ldots, 
& \underset{I\times K}{\eX_{:,J,:}}
\end{bmatrix}
\in\real^{I\times (JK)};\\
&\bX_{(2)} = 
\begin{bmatrix}
\underset{J\times I}{\eX_{:,:,1}^{\textcolor{black}{\top}}}, 
& \underset{J\times I}{\eX_{:,:,2}^{\textcolor{black}{\top}}}, 
&\ldots, 
& \underset{J\times I}{\eX_{:,:,K}^{\textcolor{black}{\top}}}
\end{bmatrix}
\in\real^{J\times (IK)};\\
&\bX_{(3)} = 
\begin{bmatrix}
\underset{K\times J}{\eX_{1,:,:}^{\textcolor{black}{\top}}}, 
& \underset{K\times J}{\eX_{2,:,:}^{\textcolor{black}{\top}}}, 
&\ldots, 
& \underset{K\times J}{\eX_{I,:,:}^{\textcolor{black}{\top}}}
\end{bmatrix}
\in\real^{K\times (IJ)},
\end{align}
where we find the mode-2 LMV matricization is the same as  the Kolda matricization.

To summarize the three matricizations,  the equivalence of matricizations for a third-order tensor is provided in ~\ref{table:different_matric}.

\begin{table}[]
\setlength{\tabcolsep}{3.4pt}
\begin{tabular}{l|l|l}
\hline
 & Method 1    & Method 2 \\ \hline
mode-$1$ 
& \shortstack{$\bX_{(1)}^{Kolda} =\bX_{(1)}^{Kiers}=$ \\ $\begin{bmatrix}
		\underset{I\times J}{\eX_{:,:,1}}, & \underset{I\times J}{\eX_{:,:,2}}, &\ldots, & \underset{I\times J}{\eX_{:,:,K}}
	\end{bmatrix}
	\in\real^{I\times (JK)}$}
& \shortstack{$\bX_{(1)}^{LMV}=$ \\ $\begin{bmatrix}
		\underset{I\times K}{\eX_{:,1,:}}, 
		& \underset{I\times K}{\eX_{:,2,:}}, 
		&\ldots, 
		& \underset{I\times K}{\eX_{:,J,:}}
	\end{bmatrix}
	\in\real^{I\times (JK)}$}      \\\hline
mode-$2$ 
& \shortstack{$\bX_{(2)}^{Kolda} =\bX_{(2)}^{LMV}=$ \\ $\begin{bmatrix}
		\underset{J\times I}{\eX_{:,:,1}^{\textcolor{black}{\top}}}, 
		& \underset{J\times I}{\eX_{:,:,2}^{\textcolor{black}{\top}}}, 
		&\ldots, 
		& \underset{J\times I}{\eX_{:,:,K}^{\textcolor{black}{\top}}}
	\end{bmatrix}
	\in\real^{J\times (IK)}$}  
& \shortstack{$\bX_{(2)}^{Kiers}=$  \\ $\begin{bmatrix}
		\underset{J\times K}{\eX_{1,:,:}}, & \underset{J\times K}{\eX_{2,:,:}}, &\ldots, & \underset{J\times K}{\eX_{I,:,:}}
	\end{bmatrix}
	\in\real^{J\times (IK)}$}    \\\hline
mode-$3$ 
& \shortstack{$\bX_{(3)}^{Kolda} =\bX_{(3)}^{Kiers}=$ \\ $\begin{bmatrix}
		\underset{K\times I}{\eX^\top_{:,1,:}}, & \underset{K\times I}{\eX^\top_{:,2,:}}, &\ldots, & \underset{K\times I}{\eX^\top_{:,J,:}}
	\end{bmatrix}
	\in\real^{K\times (IJ)}$} 
& \shortstack{$\bX_{(3)}^{LMV} =$ \\ $\begin{bmatrix}
		\underset{K\times J}{\eX_{1,:,:}^{\textcolor{black}{\top}}}, 
		& \underset{K\times J}{\eX_{2,:,:}^{\textcolor{black}{\top}}}, 
		&\ldots, 
		& \underset{K\times J}{\eX_{I,:,:}^{\textcolor{black}{\top}}}
	\end{bmatrix}
	\in\real^{K\times (IJ)} $}    \\ \hline
\end{tabular}
\caption{Comparison of different matricizations for a third-order tensor  $\eX\in\real^{I\times J\times K}$.}
\label{table:different_matric}
\end{table}

\subsection{Vectorization}
\index{Vectorization}
Specifically, the \textit{vectorization} of a tensor is defined in a similar way, 
where the entries are fetched in an ordered manner starting from the 1-st index, then the 2-nd, and so on, until the last index. This process converts the tensor into a long one-dimensional vector:
$$
vec(\eY) = 
\begin{bmatrix}
1,2, \ldots ,24
\end{bmatrix}^\top.
$$

\begin{exercise}[Vectorization of Rank-One Tensors]\label{exer:vec_rkoneten}
Let $\eX=\ba^{(1)}\circ \ba^{(2)}\circ \ldots \circ \ba^{(N)}$ be a rank-one tensor in $\real^{I_1\times I_2\times \ldots \times I_N}$,
where $\ba^{(1)}\in \real^{I_1}, \ba^{(2)}\in \real^{I_2}, \ldots, \ba^{(N)}\in \real^{I_N}$. 
Show that $vec(\eX) = \ba^{(N)} \kronecker  \ldots \kronecker \ba^{(2)}\kronecker \ba^{(1)}$, where $\kronecker$ denotes the Kronecker product (Definition~\ref{definition:kronecker-product}).
\end{exercise}
Specifically, for three vectors $\ba\in\real^I$, $\bb\in\real^J$, and $\bc\in\real^K$, it follows that
$$
\begin{aligned}
\eX &= \ba\circ\bb\circ \bc \in\real^{I\times J\times K} \iff x_{ijk} = a_i b_j c_k;\\
\bx &= \ba\kronecker \bb\kronecker \bc = \real^{I\cdot J\cdot K} \iff x_{k+(j-1)K+(i-1)JK} = a_i b_j c_k.
\end{aligned}
$$
This  demonstrates  the connection between  high-order vector outer products and Kronecker products.

\begin{exercise}[Vectorization of Matrix Products]
Show that 
\begin{itemize}
\item $vec(\ba\bb^\top) = \bb\kronecker \ba$.
\item $vec(\bA\bB\bC) = (\bC^\top\kronecker \bA)vec(\bB)$ with appropriate matrices $\bA,\bB$, and $\bC$.
\item $vec(\bA\diag(\bb)\bC) = (\bC^\top\khatrirao\bA)\bb$ with appropriate matrices $\bA,\bC$, and vector $\bb$.
\end{itemize}
\end{exercise}

As we can see,  both the matricization and vectorization of a tensor are not unique, as long as we maintain consistency in the entries, they will produce the same results in the analysis. 
For additional insights, refer to the works of \citet{de2000multilinear, kiers2000towards, kolda2009tensor, zhang2017matrix}.

\paragraph{Norm of the difference of two tensors.} The matricization or vectorization can aid in deriving the properties of  tensors. For instance, given $\eX,\eY\in \real^{I_1\times I_2\times \ldots \times I_N}$, it follows that 
\begin{equation}\label{equation:differe-tensor-norm}
\begin{aligned}
\norm{\eX-\eY}^2 &= \norm{vec(\eX)-vec(\eY)}^2 = \norm{vec(\eX)}^2 - 2vec(\eX)^\top vec(\eY) + \norm{vec(\eY)}^2\\
	&=\norm{\eX}^2 -2\langle \eX,\eY\rangle + \norm{\eY}^2.
\end{aligned}
\end{equation}

\subsection{Tuple Matricization}\label{section:tuple_matricization}
The mode-$n$ matricization is obtained along a single mode. Similarly, for the set of all modes $\overline{\mathcalN}=\{1,2,\ldots, N\}$ for $\eX\in \real^{I_1\times I_2\times \ldots \times I_N}$, we consider a subset $\overline{\mathcalI}\subset\overline{\mathcalN}$ and its complementary set $\overline{\mathcalJ} = \overline{\mathcalN}\backslash \overline{\mathcalI}$. 
The the matricization along $\overline{\mathcalI}$ follows similarly:
\begin{equation}
\textbf{(Tuple matricization)}: \qquad  
\mathcalM_{(\overline{\mathcalI})}(\eX) \triangleq \eX_{(\overline{\mathcalI})}\in\real^{\overline{\mathcalI}\times \overline{\mathcalJ}}.
\end{equation}
For example, if $\overline{\mathcalI}=\{2,3\}$, then $\mathcalM_{(\overline{\mathcalI})}(\eX)\in\real^{(I_2I_3)\times (I_1I_4\ldots I_N)}$.

\begin{exercise}[Tuple Matricization]\label{exercise:tuple_matci}
Let $\eX=\ba\circ \bb\circ\bc\circ\bd \in\real^{I_1\times I_2\times I_3\times I_4}$. 
Show that 
$$
\begin{aligned}
\eX_{(1,2)} &= (\bb\kronecker \ba)(\bd\kronecker \bc)^\top \in\real^{(I_1I_2)\times (I_3I_4)};\\
\eX_{(2,3)} &=(\bc\kronecker \bb)(\bd\kronecker \ba)^\top \in\real^{(I_2I_3)\times (I_1I_4)}.
\end{aligned}
$$
\end{exercise}

\section{Tensor Multiplication} 
We consider the \textit{mode-$n$ tensor multiplication}: multiplying the tensor in the $n$-th mode  by a matrix \footnote{This operation was originally defined in \citet{tucker1964extension, tucker1966some}, and is sometimes called the \textit{Tucker product} or \textit{Tucker mode-$n$ product}.}. Given an Nth-order tensor $\eX \in  \real^{I_1\times I_2\times \ldots \times I_N}$ and a matrix $\bA\in \real^{M\times I_n}$, the \textit{mode-$n$ tensor multiplication} of $\eX$ and $\bA$ transforms the $n$-th dimension of $\eX$ from $I_n$ to $M$. The mode-$n$ tensor multiplication is denoted by $\eY\triangleq\eX \times_n \bA \in \real^{I_1\times \ldots \times I_{n-1}\times \textcolor{mylightbluetext}{M}\times I_{n+1}\times \ldots \times I_N}$, and its element is given by
\begin{equation}\label{equation:moden-tensor-multi}
\begin{aligned}
\textbf{(Tensor multip.): }	
(\eX \times_n \bA)_{i_1\times \ldots \times i_{n-1}\times \textcolor{mylightbluetext}{m}\times i_{n+1}\times \ldots \times i_N}=
\sum_{i_n=1}^{I_n} (x_{i_1i_2\ldots i_N}) (a_{\textcolor{mylightbluetext}{m} i_n})\\
\end{aligned}
\end{equation}
The tensor multiplication can also be expressed in matrix forms (Lemma~\ref{lemma:tensor-multi-matriciz1}).
When $M=1$, the result $\eY=\eX \times_n \bA$ is an $(N-1)$th-order tensor.
Going further, the \textit{model-$n$ multilinear multiplication} matricize the result of mode-$n$ tensor multiplication:
\begin{equation}\label{equation:moden_tensor_multilinear}
\textbf{(Multilinear multip.): }\gap
\eX \circ_n \bA = \mathcalM_{(n)}(\eX \times_n \bA) \in\real^{M\times (I_1 \ldots  I_{n-1} I_{n+1} \ldots  I_N)}.
\end{equation}
That is, $(\eX \times_n \bA)$ results in a tensor (which modifies the dimension of model-$n$); $(\eX \circ_n \bA)$ results in a matrix.

\paragraph{Matrix multiplication as tensor multiplication.} In matrix language, given two matrices $\bA\in \real^{I\times K}$ and $\bB \in \real^{K\times J}$, then the matrix multiplication can be equivalently denoted by 
$$
\bA\bB = \bB\times_1 \bA.
$$

Suppose now, $\bA\in \real^{I\times K}$ and $\bB \in \real^{J\times K}$, then the matrix multiplication can be equivalently denoted by
$$
\bA\bB^\top = \bA\times_2 \bB.
$$

\paragraph{SVD and tensor multiplication.}
Supposing $\bA\in \real^{M\times N}$ has a  reduced SVD (Theorem~\ref{theorem:reduced_svd_rectangular}) given by 
$$
\underset{M\times N}{\bA} = \underset{M\times R}{\bU}
\gap
\underset{R\times R}{\bSigma}
\gap 
\underset{R\times N}{\bV^\top}.
$$
Then the reduced SVD of $\bA$ can be equivalently denoted by 
\begin{equation}\label{equation:svd-by-tensor-multi0}
	\bA = \bU\bSigma\bV^\top = \bSigma\times_1 \bU \times_2 \bV.
\end{equation}
The full SVD of $\bA$ can be represented in a similar way.

The matricization of a tensor can elucidate  the processes involved in  tensor multiplication.
\begin{lemma}[Tensor Multiplication in Matricization]\label{lemma:tensor-multi-matriciz1}
Given an Nth-order tensor $\eX \in  \real^{I_1\times I_2\times \ldots \times I_N}$ and a matrix $\bA\in \real^{M\times I_n}$, it follows that 
$$
\begin{aligned}
\eY &= \eX\times_n\bA\in \real^{I_1\times \ldots \times I_{n-1}\times \textcolor{mylightbluetext}{M}\times I_{n+1}\times \ldots \times I_N} \\  
&\gap\gap \implies\quad   \bY_{(n)} = \bA\bX_{(n)} \in \real^{M\times I_{-n}},
\end{aligned}
$$
where $I_{-n}=I_1 I_2\ldots I_{n-1}I_{n+1}\ldots I_N$. 
\end{lemma}
From the  lemma above, suppose the columns of $\bA$ are mutually orthonormal with $I_n\leq M$ (semi-orthogonal, see definition in Section~\ref{section:orthogonal-orthonormal-qr}). Then it follows that 
$$
\bA^\top \bY_{(n)} = \underbrace{\bA^\top \bA}_{\bI}\bX_{(n)}  = \bX_{(n)}.
$$
This reveals a significant property of tensor multiplication: $\eX =  \eY \times_n \bA^\top $. That is, \textbf{if $\bA$ is semi-orthogonal}, we have 
\begin{equation}\label{equation:semiorthogonal-in-tensor-multi}
\eY = \eX\times_n\bA 
\qquad\implies\qquad 
\eX =  \eY \times_n \bA^\top.
\end{equation}

We summarize the properties of tensor multiplication in the following lemma.
\begin{lemma}[Tensor Multiplication]\label{lemma:tensor-multi2}
Given an Nth-order tensor $\eX \in  \real^{I_1\times I_2\times \ldots \times I_N}$, we have 
\begin{enumerate}
\item \textit{Distributive law.} Given further the matrices $\bA\in \real^{J_m\times I_m}$ and $\bB\in \real^{J_n\times I_n}$ with $m\neq n$, it follows that 
$$
\eX \times_m \bA \times_n \bB =  (\eX \times_m \bA) \times_n \bB = (\eX  \times_n \bB) \times_m \bA.
$$
The result can be extended to a set of indices $\{n_1, n_2, \ldots, n_p\}\subseteq \{1,2,\ldots,N\}$, as long as the elements in the set are distinct.
That is, the order of mode-$n$ multiplications is irrelevant when the indices are distinct.
\item Given the matrices $\bA\in \real^{P\times I_n}$ and $\bB\in \real^{Q\times P}$, it follows that
$$
\eX \times_n\bA \times_n \bB = \eX\times_n (\bB\bA)  \in \real^{I_1\times \ldots \times I_{n-1}\times \textcolor{mylightbluetext}{Q}\times I_{n+1}\times \ldots \times I_N}.
$$
\item \label{equationbracket} Given $\eX\in \real^{I_1\times \ldots \times I_{n-1}\times \textcolor{mylightbluetext}{M}\times I_{n+1}\times \ldots \times I_N} $, $\eY\in \real^{I_1\times \ldots \times I_{n-1}\times \textcolor{mylightbluetext}{K}\times I_{n+1}\times \ldots \times I_N} $, and $\bA\in \real^{M\times K}$, it follows that 
$$
\langle\eX, \eY\times_n\bA \rangle = \langle \eX\times_n\bA^\top, \eY\rangle .
$$
\item \label{equationbracketlast} Given a semi-orthogonal matrix $\bA\in \real^{P\times I_n}$ (i.e., $\bA^\top\bA=\bI$), 
and an Nth-order tensor $\eX \in  \real^{I_1\times I_2\times \ldots \times I_N}$, it follows that 
$$
	\eY = \eX\times_n\bA \leadto \eX =  \eY \times_n \bA^\top.
$$
and 
$$
\norm{\eY} = \norm{\eX},
$$
i.e., the length is  preserved under semi-orthogonal transformations.

\item \label{equationbracket5} Given a matrix $\bA\in \real^{P\times I_n}$ with full column rank, 
and an Nth-order tensor $\eX \in  \real^{I_1\times I_2\times \ldots \times I_N}$, it follows that 
$$
\eY = \eX\times_n\bA \leadto  \eX = \eY \times_n \bA^+,
$$
where $\bA^+$ is the pseudo-inverse of $\bA$ (Appendix~\ref{appendix:pseudo-inverse_main}).
\end{enumerate}
\end{lemma}

\begin{problemset}
\item Find the closed-form representation of the mode-$n$ matricization of the  Nth-order tensor $\eX\in \real^{I_1\times I_2\times \ldots \times I_N}$ with the LMV matricization method.

\item Given the tensor $\eX\in\real^{3\times 4\times 2}$, where 
$$
\eX_{:,:,1} = 
\begin{bmatrix}
1 & 4 & 7 & 10 \\
2 & 5 & 8 & 11  \\
3 & 6 & 9 & 12
\end{bmatrix},
\qquad 
\eX_{:,:,2} = 
\begin{bmatrix}
	13 & 16 & 19 & 22 \\
	14 & 17 & 20 & 23  \\
	15 & 18 & 21 & 24
\end{bmatrix},
$$
Find the mode-1 product with $\bA=\begin{bmatrix}
1 & 3 & 5\\
2 & 4 & 6
\end{bmatrix}$
and the mode-3 product with 
$\bB=\begin{bmatrix}
	1 & 4 \\
	2 & 5 \\
	3 & 6\\
\end{bmatrix}$.

\item Prove the properties stated in Lemma~\ref{lemma:tensor-multi2}.

\end{problemset}

%% file: chapter-tensorcp.tex
\newpage
\chapter{CP, Tucker, HOSVD, TT Decompositions}
\begingroup
\hypersetup{
	linkcolor=structurecolor,
	linktoc=page,  
}
\minitoc \newpage
\endgroup

\index{Decomposition: CP}
\section{CP Decomposition}\label{section:cp_decom}

Similar to the assumption of low-rank structure in low-rank matrix approximations (LRMA) or low-rank matrix decompositions, the low-rank assumption in tensor decomposition posits that a high-dimensional tensor can be approximated by a combination of lower-dimensional tensors.
This assumption is crucial because it allows for the efficient representation and manipulation of large and complex datasets.
However, unlike matrices, tensors often exhibit intricate dependencies across multiple modes, which makes their low-rank approximation inherently challenging. In practice, except in cases of very structured data, tensor decomposition will involve some compression error because a low-rank tensor approximation generally cannot fully represent the original high-dimensional tensor. 
The quality of this approximation depends on how well the data adheres to a low-rank structure, and it may require a trade-off between approximation accuracy and computational efficiency.
Different types of tensor decompositions leverage the low-rank assumption in various ways. The \textit{CP (CANDECOMP-PARAFAC) decomposition}, for instance, assumes that the tensor can be represented as a sum of rank-one tensors---essentially outer products of vectors corresponding to each mode of the tensor. This decomposition simplifies the tensor by expressing it through a series of component vectors, providing an efficient yet interpretable structure that captures the core patterns across all modes.

\begin{theoremHigh}[CP Decomposition \citep{harshman1970foundations, carroll1970analysis}]\label{theorem:cp-decomp}
The CP decomposition factors a tensor into a sum of  rank-one
tensors. For a general Nth-order tensor $\eX \in \real^{I_1\times I_2\times \ldots \times I_N}$, it admits the CP decomposition
\begin{equation}\label{equation:cp_form1}
\eX \approx \llbracket \bA^{(1)}, \bA^{(2)}, \ldots, \bA^{(N)}   \rrbracket  
=
\sum_{r=1}^{R} \ba_r^{(1)} \circ \ba_r^{(2)} \circ \ldots \circ \ba_r^{(N)} ,
\end{equation}
where $\bA^{(n)}=[\ba_1^{(n)}, \ba_2^{(n)}, \ldots, \ba_R^{(n)}] \in \real^{I_n \times R}$ for all $n\in \{1,2,\ldots, N\}$ is the column partition of the matrix $\bA^{(n)}\in \real^{I_n \times R}$. The memory requirement to store such a full tensor is then reduced from $\prod_{n=1}^{N}I_n$ to $R\cdot(\sum_{n=1}^{N}I_n)$, which \textit{scales linearly} with the tensor order $N$ and the maximum dimension size $\max\{I_1, I_2, \ldots,I_N\}$.

\paragraph{Form 2.}
Alternatively, it is often useful to assume that the columns of $\bA^{(1)}, \bA^{(2)}, \ldots, \bA^{(N)}$ are normalized to a length of one, with the weights absorbed into the vector $\blambda\in \real^R$ so that
$$
\eX \approx \llbracket \blambda; \bA^{(1)}, \bA^{(2)}, \ldots, \bA^{(N)}   \rrbracket  
=\sum_{r=1}^{R} \lambda_r \cdot  \ba_r^{(1)} \circ \ba_r^{(2)} \circ \ldots \circ \ba_r^{(N)}.
$$

\paragraph{Form 3.}
Taking a step beyond the normalized version, we can construct a diagonal Nth-order tensor $\eL\in\real^{R\times R\times\ldots\times R}$ such that
$$
\left\{
\begin{aligned}
l_{i,i,\ldots,i} &= \lambda_i, \gap \text{if }i\in \{1,2,\ldots, R\};\\
l_{i_1,i_2,\ldots,i_N}&= 0, \gap \text{if not \{$i_1=i_2=\ldots=i_N$\}}.
\end{aligned}
\right.
$$
Then, the CP decomposition of the Nth-order tensor $\eX \in \real^{I_1\times I_2\times \ldots \times I_N}$ can be expressed as 
$$
\eX \approx \llbracket \eL; \bA^{(1)}, \bA^{(2)}, \ldots, \bA^{(N)}   \rrbracket  
=
\sum_{r_1=1}^{R} \sum_{r_2=1}^{R} \ldots \sum_{r_N=1}^{R}
l_{r_1r_2\ldots r_N}
\ba_{\textcolor{mylightbluetext}{r_1}}^{(1)} \circ \ba_{\textcolor{mylightbluetext}{r_2}}^{(2)} \circ \ldots \circ \ba_{\textcolor{mylightbluetext}{r_N}}^{(N)}.
$$
Here, we follow the notation $\llbracket \ldots  \rrbracket$ introduced by \citet{kruskal1977three}, which is also known as the \textit{Kruskal operator form}.
\end{theoremHigh}

The underlying concept of the CP decomposition is to represent the tensor as a sum of rank-one tensors, which are essentially the outer products of vectors (Section~\ref{section:rank-one-tensor}). The CP decomposition is also known as the \textit{Canonical Polyadic Decomposition}, \textit{CANDECOMP-PARAFAC}, or simply \textit{PARAFAC decomposition}, and the concept of the CP decomposition stems from two papers of Hitchcock in 1927 and 1928, respectively \citep{hitchcock1927expression, hitchcock1928multiple}.

\begin{figure}[h]
\centering
\resizebox{1.0\textwidth}{!}{%
\begin{tikzpicture}
\draw [very thick] (5.8+0.5,-3.8+0.5) rectangle (7+0.5,-2.6+0.5);
\filldraw [fill=gray!60!white,draw=green!40!black] (5.8+0.5,-3.8+0.5) rectangle (7+0.5,-2.6+0.5);
\draw [step=0.4/2, very thin, color=gray] (5.8+0.5,-3.8+0.5) grid (7+0.5,-2.6+0.5);

\draw [very thick] (5.8+0.4,-3.8+0.4) rectangle (7+0.4,-2.6+0.4);
\filldraw [fill=gray!50!white,draw=green!40!black] (5.8+0.4,-3.8+0.4) rectangle (7+0.4,-2.6+0.4);
\draw [step=0.4/2, very thin, color=gray] (5.8+0.4,-3.8+0.4) grid (7+0.4,-2.6+0.4);

\draw [very thick] (5.8+0.3,-3.8+0.3) rectangle (7+0.3,-2.6+0.3);
\filldraw [fill=gray!40!white,draw=green!40!black] (5.8+0.3,-3.8+0.3) rectangle (7+0.3,-2.6+0.3);
\draw [step=0.4/2, very thin, color=gray] (5.8+0.3,-3.8+0.3) grid (7+0.3,-2.6+0.3);

\draw [very thick] (5.8+0.2,-3.8+0.2) rectangle (7+0.2,-2.6+0.2);
\filldraw [fill=gray!30!white,draw=green!40!black] (5.8+0.2,-3.8+0.2) rectangle (7+0.2,-2.6+0.2);
\draw [step=0.4/2, very thin, color=gray] (5.8+0.2,-3.8+0.2) grid (7+0.2,-2.6+0.2);

\draw [very thick] (5.8+0.1,-3.8+0.1) rectangle (7+0.1,-2.6+0.1);
\filldraw [fill=gray!20!white,draw=green!40!black] (5.8+0.1,-3.8+0.1) rectangle (7+0.1,-2.6+0.1);
\draw [step=0.4/2, very thin, color=gray] (5.8+0.1,-3.8+0.1) grid (7+0.1,-2.6+0.1);

\draw [very thick] (5.8,-3.8) rectangle (7,-2.6);
\filldraw [fill=gray!10!white,draw=green!40!black] (5.8,-3.8) rectangle (7,-2.6);
\draw [step=0.4/2, very thin, color=gray] (5.8,-3.8) grid (7,-2.6);

\draw [very thick] (5.8-0.1,-3.8-0.1) rectangle (7-0.1,-2.6-0.1);
\filldraw [fill=gray!6!white,draw=green!40!black] (5.8-0.1,-3.8-0.1) rectangle (7-0.1,-2.6-0.1);
\draw [step=0.4/2, very thin, color=gray] (5.8-0.1,-3.8-0.1) grid (7-0.1,-2.6-0.1);
	
\draw[->] (5.6,-2.6) -- ++(0.6,0.6); 
\draw[->] (5.55,-2.7) -- ++(0,-1.2);  
\draw[->] (5.66,-4.04) -- ++(1.25,0);  

\draw (6.35,-4.25) node {{\color{black}\scriptsize{spatial $y$}}};
\draw (5.2-0.3,-3.2) node {{\color{black}\scriptsize{spatial $x$}}};
\draw (5.3,-2.2) node {{\color{black}\scriptsize{temporal}}};
	
	\draw (6.2,-4.6) node {{\color{black}\scriptsize{$\eX\in\real^{I\times J \times K}$}}};
	
	\draw (8,-3.2) node {	{\color{black}\large{$\approx$}}};

	\draw [very thick] (8.6,-4.4+0.2) rectangle (8.8,-3.2+0.2);
	\filldraw [fill=WildStrawberry!40!white,draw=green!40!black] (8.6,-4.4+0.2) rectangle (8.8,-3.2+0.2);
	\draw (8.7,-4.8+0.2) node {{\color{black}\scriptsize{$\ba_1 \in \real^I$}}};
	
	\draw [very thick] (9,-3-0.2+0.2) rectangle (10.2,-2.8-0.2+0.2);
	\filldraw [fill=RubineRed!60!white,draw=green!40!black] (9,-3-0.2+0.2) rectangle (10.2,-2.8-0.2+0.2);
	\draw (9.6,-3.3-0.2+0.2) node {{\color{black}\scriptsize{$\bb_1 \in \real^J$}}};

	\draw[fill=RedOrange!40!white, line width=0.8pt] (9.2,-2.4+0.2) -- (9.4,-2.4+0.2) -- (8.8,-3+0.2) -- (8.6,-3+0.2) -- cycle;
	\draw (10-0.2,-1.9) node {{\color{black}\scriptsize{$\bc_1\in \real^K$}}};
	
	\draw (10.8,-3.2) node {{\color{black}\large{$+$}}};
	\draw (11.6,-3.2) node {{\color{black}\large{$\dots$}}};
	\draw (12.4,-3.2) node {{\color{black}\large{$+$}}};
	
	\draw [very thick] (8.6+4.4,-4.4+0.2) rectangle (8.8+4.4,-3.2+0.2);
	\filldraw [fill=WildStrawberry!40!white,draw=green!40!black] (8.6+4.4,-4.4+0.2) rectangle (8.8+4.4,-3.2+0.2);
	\draw (8.7+4.4,-4.8+0.2) node {{\color{black}\scriptsize{$\ba_R\in \mathbb{R}^I $}}};
	
	\draw [very thick] (9+4.4,-3-0.2+0.2) rectangle (10.2+4.4,-2.8-0.2+0.2);
	\filldraw [fill=RubineRed!60!white,draw=green!40!black] (9+4.4,-3-0.2+0.2) rectangle (10.2+4.4,-2.8-0.2+0.2);
	\draw (9.6+4.4,-3.3-0.2+0.2) node {{\color{black}\scriptsize{$\bb_R\in \real^J$}}};
	
	\draw[fill=RedOrange!40!white, line width=0.8pt] (9.2+4.4,-2.4+0.2) -- (9.4+4.4,-2.4+0.2) -- (8.8+4.4,-3+0.2) -- (8.6+4.4,-3+0.2) -- cycle;
	\draw (10.1+4,-1.9) node {{\color{black}\scriptsize{$\bc_R\in \real^{K}$}}};
	
\end{tikzpicture}
}
\caption{The CP decomposition of a third-order tensor: $\eX \approx \llbracket \bA, \bB, \bC   \rrbracket  
	=
	\sum_{r=1}^{R} \ba_r\circ \bb_r \circ \bc_r$, where $\eX\in \real^{I\times J\times K}$, $\bA\in \real^{I\times R}, \bB\in \real^{J\times R}$, and $\bC\in \real^{K\times R}$. Compare to the third-order rank-one tensor in Figure~\ref{fig:rank-one-tensor}.}
\label{fig:cp-decom-third}
\end{figure}

\begin{figure}[h]
	\centering
	\resizebox{0.7\textwidth}{!}{%
\begin{tikzpicture}
	\coordinate (O) at (0,0,0);
	\coordinate (A) at (0,\Width,0);
	\coordinate (B) at (0,\Width,\Height);
	\coordinate (C) at (0,0,\Height);
	\coordinate (D) at (\Depth,0,0);
	\coordinate (E) at (\Depth,\Width,0);
	\coordinate (F) at (\Depth,\Width,\Height);
	\coordinate (G) at (\Depth,0,\Height);
	\draw[red!60!black,fill=red!5] (O) -- (C) -- (G) -- (D) -- cycle;
	\draw[red!60!black,fill=red!5] (O) -- (A) -- (E) -- (D) -- cycle;
	\draw[red!60!black,fill=red!5] (O) -- (A) -- (B) -- (C) -- cycle;
	\draw[red!60!black,fill=red!5,opacity=0.8] (D) -- (E) -- (F) -- (G) -- cycle;
	\draw[red!60!black,fill=red!5,opacity=0.6] (C) -- (B) -- (F) -- (G) -- cycle;
	\draw[red!60!black,fill=red!5,opacity=0.8] (A) -- (B) -- (F) -- (E) -- cycle;
	
	\coordinate (O) at (0+\xx,0+\yy,0+\zz);
	\coordinate (A) at (0+\xx,0.25\Width+\yy,0+\zz);
	\coordinate (B) at (0+\xx,0.25\Width+\yy,0.25\Height+\zz);
	\coordinate (C) at (0+\xx,0+\yy,0.25\Height+\zz);
	\coordinate (D) at (0.25\Depth+\xx,0+\yy,0+\zz);
	\coordinate (E) at (0.25\Depth+\xx,0.25\Width+\yy,0+\zz);
	\coordinate (F) at (0.25\Depth+\xx,0.25\Width+\yy,0.25\Height+\zz);
	\coordinate (G) at (0.25\Depth+\xx,0+\yy,0.25\Height+\zz);
	\draw[green!80!black,fill=green!10] (O) -- (C) -- (G) -- (D) -- cycle;
	\draw[green!80!black,fill=green!10] (O) -- (A) -- (E) -- (D) -- cycle;
	\draw[green!80!black,fill=green!10] (O) -- (A) -- (B) -- (C) -- cycle;
	\draw[green!40!black,fill=green!10,opacity=0.8] (D) -- (E) -- (F) -- (G) -- cycle;
	\draw[green!40!black,fill=green!10,opacity=0.6] (C) -- (B) -- (F) -- (G) -- cycle;
	\draw[green!40!black,fill=green!10,opacity=0.8] (A) -- (B) -- (F) -- (E) -- cycle;
	\draw (0.2,-1.2,0) node {\scriptsize{\color{gray}$J$}};
	\draw (0.2,-0.95,0) node[rotate = 0] {{\color{gray!65}$\underbrace{\hspace{2cm}}$}};
	\draw (-0.4,1,2) node[rotate = 90] {\scriptsize{\color{gray}$I$}};
	\draw (-0.15,1,2) node[rotate = 270] {{\color{gray!65}$\underbrace{\hspace{2cm}}$}};
	\draw (2.5,-0.2,1.4) node[rotate = 45] {\scriptsize{\color{gray}$K$}};
	\draw (2.2,0,1.2) node[rotate = 45] {{\color{gray!65}$\underbrace{\hspace{1.1cm}}$}};
	\draw [draw=gray!75,thick,->] (0.8,0.1,1) -- (1,0.8,1) node [right] {{\color{green!50!black}$x_{ijk}$}};
	\draw (0.8,0,1) node {\scriptsize{\color{gray}$(i,j,k)$-th}};
	\draw (1,-1.2,1) node {\color{black}$\eX\in\real^{I\times J\times K}$};
	
	\draw (3,0.5,0) node[rotate = 0] {{\color{black}\LARGE{$\approx$}}};
	
	\coordinate (O) at (0+\xd,0+\yd,0+\zd);
	\coordinate (A) at (0+\xd,\Width+\yd,0+\zd);
	\coordinate (B) at (0+\xd,\Width+\yd,0.25\Height+\zd);
	\coordinate (C) at (0+\xd,0+\yd,0.25\Height+\zd);
	\coordinate (D) at (0.5\Depth+\xd,0+\yd,0+\zd);
	\coordinate (E) at (0.5\Depth+\xd,\Width+\yd,0+\zd);
	\coordinate (F) at (0.5\Depth+\xd,\Width+\yd,0.25\Height+\zd);
	\coordinate (G) at (0.5\Depth+\xd,0+\yd,0.25\Height+\zd);
	\draw[red!60,fill=yellow!25] (O) -- (C) -- (G) -- (D) -- cycle;
	\draw[red!60,fill=yellow!25] (O) -- (A) -- (E) -- (D) -- cycle;
	\draw[red!60,fill=yellow!25] (O) -- (A) -- (B) -- (C) -- cycle;
	\draw[red!60,fill=yellow!25,opacity=0.8] (D) -- (E) -- (F) -- (G) -- cycle;
	\draw[red!60,fill=yellow!25,opacity=0.6] (C) -- (B) -- (F) -- (G) -- cycle;
	\draw[red!60,fill=yellow!25,opacity=0.8] (A) -- (B) -- (F) -- (E) -- cycle;
	
	\coordinate (O) at (0+\xd,0+\yy-0.1,0+\zd);
	\coordinate (A) at (0+\xd,0.25\Width+\yy-0.1,0+\zd);
	\coordinate (B) at (0+\xd,0.25\Width+\yy-0.1,0.25\Height+\zd);
	\coordinate (C) at (0+\xd,0+\yy-0.1,0.25\Height+\zd);
	\coordinate (D) at (0.5\Depth+\xd,0+\yy-0.1,0+\zd);
	\coordinate (E) at (0.5\Depth+\xd,0.25\Width+\yy-0.1,0+\zd);
	\coordinate (F) at (0.5\Depth+\xd,0.25\Width+\yy-0.1,0.25\Height+\zd);
	\coordinate (G) at (0.5\Depth+\xd,0+\yy-0.1,0.25\Height+\zd);
	\draw[green!80!black,fill=green!10] (O) -- (C) -- (G) -- (D) -- cycle;
	\draw[green!80!black,fill=green!10] (O) -- (A) -- (E) -- (D) -- cycle;
	\draw[green!80!black,fill=green!10] (O) -- (A) -- (B) -- (C) -- cycle;
	\draw[green!40!black,fill=green!10,opacity=0.8] (D) -- (E) -- (F) -- (G) -- cycle;
	\draw[green!40!black,fill=green!10,opacity=0.6] (C) -- (B) -- (F) -- (G) -- cycle;
	\draw[green!40!black,fill=green!10,opacity=0.8] (A) -- (B) -- (F) -- (E) -- cycle;
	
	\draw[red!60] (0.5\Depth+\xd,\Width+\yd,0.25\Height+\zd) -- (0.5\Depth+\xd,0+\yd,0.25\Height+\zd);
	\draw (\xd+0.2,\yy-1.9,\zz) node {\color{black}$\bA\in\real^{I\times R}$};
	\draw (\xd-0.3,\yy-0.1,\zz) node {\color{green!50!black}$\boldsymbol{a}_{i}$};
	
	\coordinate (O) at (0+\xd+1.2,0+\ye,0+\zd);
	\coordinate (A) at (0+\xd+1.2,0.5\Width+\ye,0+\zd);
	\coordinate (B) at (0+\xd+1.2,0.5\Width+\ye,0.25\Height+\zd);
	\coordinate (C) at (0+\xd+1.2,0+\ye,0.25\Height+\zd);
	\coordinate (D) at (\Depth+\xd+1.2,0+\ye,0+\zd);
	\coordinate (E) at (\Depth+\xd+1.2,0.5\Width+\ye,0+\zd);
	\coordinate (F) at (\Depth+\xd+1.2,0.5\Width+\ye,0.25\Height+\zd);
	\coordinate (G) at (\Depth+\xd+1.2,0+\ye,0.25\Height+\zd);
	\draw[red!60,fill=blue!15] (O) -- (C) -- (G) -- (D) -- cycle;
	\draw[red!60,fill=blue!15] (O) -- (A) -- (E) -- (D) -- cycle;
	\draw[red!60,fill=blue!15] (O) -- (A) -- (B) -- (C) -- cycle;
	\draw[red!60,fill=blue!15,opacity=0.8] (D) -- (E) -- (F) -- (G) -- cycle;
	\draw[red!60,fill=blue!15,opacity=0.6] (C) -- (B) -- (F) -- (G) -- cycle;
	\draw[red!60,fill=blue!15,opacity=0.8] (A) -- (B) -- (F) -- (E) -- cycle;
	
	\coordinate (O) at (0+\xd+2.5,0+\ye,0+\zd);
	\coordinate (A) at (0+\xd+2.5,0.5\Width+\ye,0+\zd);
	\coordinate (B) at (0+\xd+2.5,0.5\Width+\ye,0.25\Height+\zd);
	\coordinate (C) at (0+\xd+2.5,0+\ye,0.25\Height+\zd);
	\coordinate (D) at (0.25\Depth+\xd+2.5,0+\ye,0+\zd);
	\coordinate (E) at (0.25\Depth+\xd+2.5,0.5\Width+\ye,0+\zd);
	\coordinate (F) at (0.25\Depth+\xd+2.5,0.5\Width+\ye,0.25\Height+\zd);
	\coordinate (G) at (0.25\Depth+\xd+2.5,0+\ye,0.25\Height+\zd);
	\draw[green!80!black,fill=green!10] (O) -- (C) -- (G) -- (D) -- cycle;
	\draw[green!80!black,fill=green!10] (O) -- (A) -- (E) -- (D) -- cycle;
	\draw[green!80!black,fill=green!10] (O) -- (A) -- (B) -- (C) -- cycle;
	\draw[green!40!black,fill=green!10,opacity=0.8] (D) -- (E) -- (F) -- (G) -- cycle;
	\draw[green!40!black,fill=green!10,opacity=0.6] (C) -- (B) -- (F) -- (G) -- cycle;
	\draw[green!40!black,fill=green!10,opacity=0.8] (A) -- (B) -- (F) -- (E) -- cycle;
	
	\draw[red!60] (0+\xd+1.2,0.5\Width+\ye,0.25\Height+\zd) -- (\Depth+\xd+1.2,0.5\Width+\ye,0.25\Height+\zd);
	\draw (\xd+2.0,\yy-0.8,\zz) node {\color{black}$\bB\in\real^{J\times R}$};
	\draw (\xd+2.6,\yy-0.3,\zz) node {\color{green!50!black}$\bb_{j}$};
	
	\coordinate (O) at (0+\xd+1.2,0+\ye+1.5,0+\zd);
	\coordinate (A) at (0+\xd+1.2,0.25\Width+\ye+1.5,0+\zd);
	\coordinate (B) at (0+\xd+1.2,0.25\Width+\ye+1.5,\Height+\zd);
	\coordinate (C) at (0+\xd+1.2,0+\ye+1.5,\Height+\zd);
	\coordinate (D) at (0.5\Depth+\xd+1.2,0+\ye+1.5,0+\zd);
	\coordinate (E) at (0.5\Depth+\xd+1.2,0.25\Width+\ye+1.5,0+\zd);
	\coordinate (F) at (0.5\Depth+\xd+1.2,0.25\Width+\ye+1.5,\Height+\zd);
	\coordinate (G) at (0.5\Depth+\xd+1.2,0+\ye+1.5,\Height+\zd);
	\draw[red!60,fill=red!15] (O) -- (C) -- (G) -- (D) -- cycle;
	\draw[red!60,fill=red!15] (O) -- (A) -- (E) -- (D) -- cycle;
	\draw[red!60,fill=red!15] (O) -- (A) -- (B) -- (C) -- cycle;
	\draw[red!60,fill=red!15,opacity=0.8] (D) -- (E) -- (F) -- (G) -- cycle;
	\draw[red!60,fill=red!15,opacity=0.6] (C) -- (B) -- (F) -- (G) -- cycle;
	\draw[red!60,fill=red!15,opacity=0.8] (A) -- (B) -- (F) -- (E) -- cycle;
	
	\coordinate (O) at (0+\xd+0.8,0+\yy+1.2,0+\zd);
	\coordinate (A) at (0+\xd+0.8,0.25\Width+\yy+1.2,0+\zd);
	\coordinate (B) at (0+\xd+0.8,0.25\Width+\yy+1.2,0.25\Height+\zd);
	\coordinate (C) at (0+\xd+0.8,0+\yy+1.2,0.25\Height+\zd);
	\coordinate (D) at (0.5\Depth+\xd+0.8,0+\yy+1.2,0+\zd);
	\coordinate (E) at (0.5\Depth+\xd+0.8,0.25\Width+\yy+1.2,0+\zd);
	\coordinate (F) at (0.5\Depth+\xd+0.8,0.25\Width+\yy+1.2,0.25\Height+\zd);
	\coordinate (G) at (0.5\Depth+\xd+0.8,0+\yy+1.2,0.25\Height+\zd);
	\draw[green!80!black,fill=green!10] (O) -- (C) -- (G) -- (D) -- cycle;
	\draw[green!80!black,fill=green!10] (O) -- (A) -- (E) -- (D) -- cycle;
	\draw[green!80!black,fill=green!10] (O) -- (A) -- (B) -- (C) -- cycle;
	\draw[green!40!black,fill=green!10,opacity=0.8] (D) -- (E) -- (F) -- (G) -- cycle;
	\draw[green!40!black,fill=green!10,opacity=0.6] (C) -- (B) -- (F) -- (G) -- cycle;
	\draw[green!40!black,fill=green!10,opacity=0.8] (A) -- (B) -- (F) -- (E) -- cycle;
	
	\draw[red!60] (0.5\Depth+\xd+1.2,0.25\Width+\ye+1.5,0+\zd) -- (0.5\Depth+\xd+1.2,0.25\Width+\ye+1.5,\Height+\zd);
	\draw (\xd+2.5,\yy+1.3,\zz) node {\color{black}$\bC\in\real^{K\times R}$};
	\draw (\xd+0.5,\yy+1.4,\zz) node {\color{green!50!black}$\bc_{k}$};

\end{tikzpicture}
}
\caption{Index of the CP decomposition of a third-order tensor: $\eX\approx\llbracket \bA,\bB,\bC\rrbracket$.}
\label{fig:cp-decom-third-index}
\end{figure}

\begin{figure}[h]
	\centering
	\resizebox{0.7\textwidth}{!}{%
		\begin{tikzpicture}
\pgfmathsetmacro{\bxx}{0.3}
\pgfmathsetmacro{\begupwhite}{0.2}

\pgfmathsetmacro{\bcwd}{0.3}
\pgfmathsetmacro{\bbwd}{0.3}
\pgfmathsetmacro{\bbright}{0.9}
\pgfmathsetmacro{\bbup}{0.225}

\pgfmathsetmacro{\distrightdown}{4.2}
\pgfmathsetmacro{\distl}{1.2}
\pgfmathsetmacro{\aprimeright}{0.4}
			
			\draw [very thick] (5.8+0.5,-3.8+0.5) rectangle (7+0.5,-2.6+0.5);
			\filldraw [fill=gray!60!white,draw=green!40!black] (5.8+0.5,-3.8+0.5) rectangle (7+0.5,-2.6+0.5);
			\draw [step=0.4/2, very thin, color=gray] (5.8+0.5,-3.8+0.5) grid (7+0.5,-2.6+0.5);
			
			\draw [very thick] (5.8+0.4,-3.8+0.4) rectangle (7+0.4,-2.6+0.4);
			\filldraw [fill=gray!50!white,draw=green!40!black] (5.8+0.4,-3.8+0.4) rectangle (7+0.4,-2.6+0.4);
			\draw [step=0.4/2, very thin, color=gray] (5.8+0.4,-3.8+0.4) grid (7+0.4,-2.6+0.4);
			
			\draw [very thick] (5.8+0.3,-3.8+0.3) rectangle (7+0.3,-2.6+0.3);
			\filldraw [fill=gray!40!white,draw=green!40!black] (5.8+0.3,-3.8+0.3) rectangle (7+0.3,-2.6+0.3);
			\draw [step=0.4/2, very thin, color=gray] (5.8+0.3,-3.8+0.3) grid (7+0.3,-2.6+0.3);
			
			\draw [very thick] (5.8+0.2,-3.8+0.2) rectangle (7+0.2,-2.6+0.2);
			\filldraw [fill=gray!30!white,draw=green!40!black] (5.8+0.2,-3.8+0.2) rectangle (7+0.2,-2.6+0.2);
			\draw [step=0.4/2, very thin, color=gray] (5.8+0.2,-3.8+0.2) grid (7+0.2,-2.6+0.2);
			
			\draw [very thick] (5.8+0.1,-3.8+0.1) rectangle (7+0.1,-2.6+0.1);
			\filldraw [fill=gray!20!white,draw=green!40!black] (5.8+0.1,-3.8+0.1) rectangle (7+0.1,-2.6+0.1);
			\draw [step=0.4/2, very thin, color=gray] (5.8+0.1,-3.8+0.1) grid (7+0.1,-2.6+0.1);
			
			\draw [very thick] (5.8,-3.8) rectangle (7,-2.6);
			\filldraw [fill=gray!10!white,draw=green!40!black] (5.8,-3.8) rectangle (7,-2.6);
			\draw [step=0.4/2, very thin, color=gray] (5.8,-3.8) grid (7,-2.6);

			\draw [very thick] (5.8-0.1,-3.8-0.1) rectangle (7-0.1,-2.6-0.1);
			\filldraw [fill=gray!6!white,draw=green!40!black] (5.8-0.1,-3.8-0.1) rectangle (7-0.1,-2.6-0.1);
			\draw [step=0.4/2, very thin, color=gray] (5.8-0.1,-3.8-0.1) grid (7-0.1,-2.6-0.1);
			
			\draw[->] (5.6,-2.6) -- ++(0.6,0.6); 
			\draw[->] (5.55,-2.7) -- ++(0,-1.2);  
			\draw[->] (5.66,-4.04) -- ++(1.25,0);  
			
			\draw (6.35,-4.25) node {{\color{black}\scriptsize{spatial $y$}}};
			\draw (5.2-0.3,-3.2) node {{\color{black}\scriptsize{spatial $x$}}};
			\draw (5.3,-2.2) node {{\color{black}\scriptsize{temporal}}};
			
			\draw (6.2,-4.6) node {{\color{black}\scriptsize{$\eX\in\real^{I\times J \times K}$}}};
			
			\draw (8,-3.2) node {	{\color{black}\large{$\approx$}}};

\draw [very thick] (8.6+\distl,-5.0+0.2) rectangle (9.2+\distl,-3.2+0.2);
\filldraw [fill=RubineRed!60!white,draw=green!40!black] (8.6+\distl,-5.0+0.2) rectangle (9.2+\distl,-3.2+0.2);
\draw (8.7+\distl,-5.3+0.2) node {{\color{black}\scriptsize{$\bB \in \real^{J\times R}$}}};
			
\draw [very thick] 
(8.8+\bbright+\distl,-2.9+\bbup) rectangle (10.7+\bbright+\distl,-2.3+\bbup);
\filldraw  [fill=WildStrawberry!40!white,draw=green!40!black]
(8.8+\bbright+\distl,-2.9+\bbup) rectangle (10.7+\bbright+\distl,-2.3+\bbup);
\draw (9.6+\bbright+\distl,-3.2+\bbup) node {{\color{black}\scriptsize{$\bA^\top \in \real^{R\times I}$}}};
			
\draw[fill=RedOrange!40!white, line width=0.8pt] 
(8.45+\distl,-2.0-\bxx) -- (8.8+\distl,-1.8-\bxx) -- (8.25+\distl-0.7,-1.6-\bxx+0.6) -- (7.9+\distl-0.7,-1.8-\bxx+0.6) -- cycle;
\draw (8.8+\distl,-1.5) node {{\color{black}\scriptsize{$\bC\in \real^{K\times R}$}}};
			
\draw[fill=RedOrange!10!white, line width=0.5pt] 
(9.2+\distl,-2.0-\bxx) -- (9.3+\bcwd+\distl-0.02,-1.8+0.01-\bxx) -- (8.7+\bcwd+\distl,-1.8-\bxx) -- (8.6+\distl+0.02,-2.0+0.01-\bxx) -- cycle;
\draw[fill=RedOrange!10!white, line width=0.8pt]
(8.6+\distl,-2.0-\bxx) -- (9.2+\distl,-2.0-\bxx) -- (9.2+\distl,-2.6-\bxx) -- (8.6+\distl,-2.6-\bxx) -- cycle;
\draw[fill=RedOrange!10!white, line width=0.8pt]
(9.3+\bcwd+\distl,-1.8-\bxx) -- (9.2+\distl,-2.1+0.1-\bxx) -- (9.2+\distl,-3+0.4-\bxx) -- (9.3+\bcwd+\distl,-2.4-\bxx) -- cycle;
\draw (10-0.2+\distl,-1.5-\bxx) node {{\color{black}\scriptsize{$\eL\in \real^{R\times R\times R}$}}};

\end{tikzpicture}
}
\caption{The CP decomposition of a third-order tensor: $\eX \approx  \llbracket\eL; \bA^{(1)}, \bA^{(2)}, \ldots, \bA^{(N)}   \rrbracket   
=
\sum_{r=1}^{R} \lambda_r\cdot \ba_r\circ \bb_r \circ \bc_r$, where $\eX\in \real^{I\times J\times K}$, $\eL\in \real^{R\times R\times R}$, and where the columns of $\bA\in \real^{I\times R}, \bB\in \real^{J\times R}$ and $\bC\in \real^{K\times R}$ are of \textcolor{mylightbluetext}{unit length}. 
Compare to the first form of the CP decomposition in Figure~\ref{fig:cp-decom-third}.
}
\label{fig:cp-decom-thirr-2ndform}
\end{figure}

\paragraph{Rank of a tensor.} From the CP decomposition, the rank of a tensor $\eX$ can be defined as the smallest value of $R$ for which the CP decomposition holds \textbf{exactly}.

\paragraph{Matricization.}
When $\eX$ is a $2$nd-order tensor (i.e., a matrix), the CP decomposition simplifies to the rank decomposition, in a sense of ignoring the noise term (Theorem~\ref{theorem:rank-decomposition-alternative}). For simplicity, we will consider the CP decomposition for a third-order tensor $\eX\in \real^{I\times J\times K}$, as depicted in  Figure~\ref{fig:cp-decom-third}:
\begin{equation}\label{equation:third-fact-1}
\eX \approx \llbracket \bA, \bB, \bC   \rrbracket  
=
\sum_{r=1}^{R} \ba_r\circ \bb_r \circ \bc_r.
\end{equation}
In terms of matrix dimensions, we have $\bA\in \real^{I\times R}, \bB\in \real^{J\times R}$, and $\bC\in \real^{K\times R}$.
Expressing this CP decomposition element-wise, as shown in Figure~\ref{fig:cp-decom-third-index}, the  Equation~\eqref{equation:third-fact-1} can be equivalently rephrased as 
$$
x_{ijk} = \sum_{r=1}^{R}a_{ir}b_{jr}c_{kr}.
$$
Therefore, the matricized form of the third-order CP decomposition can be expressed as 
\begin{equation}\label{equation:cp-matriciza-third-1}
\left\{
\begin{aligned}
	\bX_{(1)}  &\approx \bA (\bC\khatrirao \bB )^\top \in \real^{I\times (JK)}; \\
	\bX_{(2)} &\approx \bB (\bC\khatrirao \bA )^\top \in \real^{J\times (IK)};\\
	\bX_{(3)} &\approx \bC(\bB\khatrirao \bA)^\top \in \real^{K\times (IJ)},
\end{aligned}
\right.
\footnote{Unless otherwise specified, the three matricizations are derived using the Kolda matricization method (Section~\ref{section:kolda_matric}). 
In a comprehensive representation, the six matricizations can be obtained by 
$$
\begin{aligned}
\bX_{(1)}^{Kolda} =\bX_{(1)}^{Kiers}
&\approx \bA (\bC\khatrirao \bB )^\top \in \real^{I\times (JK)},
\qquad \bX_{(1)}^{LMV}
\approx \bA(\bB\khatrirao \bC)^\top \in \real^{I\times (JK)},\\
\bX_{(2)}^{Kolda} =\bX_{(2)}^{LMV}
&\approx \bB (\bC\khatrirao \bA )^\top \in \real^{J\times (IK)},
\qquad\bX_{(2)}^{Kiers}
\approx \bB(\bA\khatrirao \bC)^\top \in \real^{J\times (IK)},\\
\bX_{(3)}^{Kolda} =\bX_{(3)}^{Kiers}
&\approx \bC(\bB\khatrirao \bA)^\top \in \real^{K\times (IJ)},
\qquad\bX_{(3)}^{LMV}
\approx \bC(\bA\khatrirao \bB)^\top\in \real^{K\times (IJ)},
\end{aligned}
$$
respectively. See also Table~\ref{table:different_matric}.
}
\end{equation}
where $\bX_{(n)}\in \real^{I\times (JK)}$ is the mode-$n$ matricization of the tensor $\eX$ for $n\in \{1,2,3\}$, and ``$\khatrirao$" denotes the Khatri-Rao product of matrices (Definition~\ref{definition:khatri-rao-product}) such that
$$
\left\{
\begin{aligned}
\bC\khatrirao \bB&=[\bc_1\kronecker \bb_1, \bc_2\kronecker \bb_2, \ldots, \bc_R\kronecker \bb_R]\in \real^{JK\times R};\\
\bC\khatrirao \bA&=[\bc_1\kronecker \ba_1, \bc_2\kronecker \ba_2, \ldots, \bc_R\kronecker \ba_R]\in \real^{IK\times R};\\
\bB\khatrirao \bA&=[\bb_1\kronecker \ba_1, \bb_2\kronecker \ba_2, \ldots, \bb_R\kronecker \ba_R]\in \real^{IJ\times R}.\\
\end{aligned}
\right.
$$

\paragraph{Matricization of general Nth-order tensors.}
In full generality, consider the Nth-order tensor $\eX \in \real^{I_1\times I_2\times \ldots \times I_N}$, the mode-$n$ matricized form of the CP decomposition $\eX \approx \llbracket \bA^{(1)}, \bA^{(2)}, \ldots, \bA^{(N)}   \rrbracket  $ is given by 
$$
\boxed{\underbrace{\bX_{(n)}}_{I_n\times (I_{-n})}\approx
	\underbrace{\bA^{(n)}}_{{I_n\times R}} \underbrace{\left(\bA^{(N)} \khatrirao \bA^{(N-1)} \khatrirao \ldots \khatrirao \bA^{(n+1)}\khatrirao \bA^{(n-1)}\khatrirao \ldots \khatrirao \bA^{(2)}\khatrirao \bA^{(1)}  \right)^\top}_{R\times (I_{-n})}},
$$
where $I_{-n}=I_1 I_2\ldots I_{n-1}I_{n+1}\ldots I_N$.

\paragraph{Vectorization of general Nth-Order tensors.} Similarly, the vectorization of the CP decomposition with the form $\eX \approx \llbracket \bA^{(1)}, \bA^{(2)}, \ldots, \bA^{(N)}   \rrbracket  $ is given by 
$$
\boxed{
\begin{aligned}
\underbrace{vec(\eX)}_{(I_1 I_2\ldots I_N)\times 1}
 \approx 
\underbrace{\left(\bA^{(N)} \khatrirao \bA^{(N-1)} \khatrirao \ldots  \khatrirao \bA^{(2)}\khatrirao \bA^{(1)}  \right)}_{(I_1 I_2\ldots I_N)\times R} \cdot \underbrace{\bm{1}}_{R\times 1}
\end{aligned}
},
$$
where $\bm{1}=[1,1,\ldots, 1]^\top \in \real^{R\times 1}$ denotes the all-ones vector.
Alternatively, the vectorization of CP forms can be represented via a summation of $R$ vectors; see Problem~\ref{prob:vec_cpform}.

Table~\ref{table:different_cpforms} presents a summary of equivalent  representations for the CP decomposition in Equation~\eqref{equation:cp_form1}.

\begin{table}[]
\setlength{\tabcolsep}{8pt}    
\begin{tabular}{l|l}
\hline\hline
Form                  & Formula                                                                                                                                                                                                                                                               \\ \hline
Operator      & $\eX \approx \llbracket \bA^{(1)}, \bA^{(2)}, \ldots, \bA^{(N)}   \rrbracket$                                                                                                                                                                                         \\\hline
Outer prod.         & $\eX \approx\sum_{r=1}^{R} \ba_r^{(1)} \circ \ba_r^{(2)} \circ \ldots \circ \ba_r^{(N)}$                                                                                                                                                                              \\\hline
Matricization & $\footnotesize\underbrace{\bX_{(n)}}_{I_n\times (I_{-n})}\approx\underbrace{\bA^{(n)}}_{{I_n\times R}} \underbrace{\left(\bA^{(N)} \khatrirao \bA^{(N-1)} \khatrirao \ldots \khatrirao \bA^{(n+1)}\khatrirao \bA^{(n-1)}\khatrirao \ldots \khatrirao \bA^{(2)}\khatrirao \bA^{(1)}  \right)^\top}_{R\times (I_{-n})}$ \\\hline
Vectorization         & $\underbrace{vec(\eX)}_{(I_1 I_2\ldots I_N)\times 1}\approx \underbrace{\left(\bA^{(N)} \khatrirao \bA^{(N-1)} \khatrirao \ldots  \khatrirao \bA^{(2)}\khatrirao \bA^{(1)}  \right)}_{(I_1 I_2\ldots I_N)\times R} \cdot \underbrace{\bm{1}}_{R\times 1}$\\
\hline  \hline                            
\end{tabular}
\caption{Equivalent forms for the CP decomposition in Equation~\eqref{equation:cp_form1}.}
\label{table:different_cpforms}
\end{table}

\paragraph{Equivalent forms of the CP decomposition.}
As stated in  Theorem~\ref{theorem:cp-decomp}, when the columns of $\bA^{(n)}$'s are of unit length with the weights absorbed into the vector $\blambda\in \real^R$, then the CP decomposition of the Nth-order tensor $\eX \in \real^{I_1\times I_2\times \ldots \times I_N}$ can be expressed as 
$$
\eX \approx \llbracket \blambda; \bA^{(1)}, \bA^{(2)}, \ldots, \bA^{(N)}   \rrbracket  
=\sum_{r=1}^{R} \lambda_r \cdot  \ba_r^{(1)} \circ \ba_r^{(2)} \circ \ldots \circ \ba_r^{(N)}.
$$
Suppose further the Nth-order tensor $\eL\in \real^{R\times R\times \ldots \times R}$ is a diagonal tensor (Definition~\ref{definition:identity-tensors}) with 
\begin{equation}\label{equation:cp-oth3er-fuorm-nth1}
\left\{
\begin{aligned}
	l_{i,i,\ldots,i} &= \lambda_i, \gap \text{if }i\in \{1,2,\ldots, R\};\\
	l_{i_1,i_2,\ldots,i_N}&= 0, \gap \text{if not \{$i_1=i_2=\ldots=i_N$\}}.
\end{aligned}
\right.
\end{equation}
Then, the CP decomposition of the Nth-order tensor $\eX \in \real^{I_1\times I_2\times \ldots \times I_N}$ can be expressed as 
\begin{equation}\label{equation:cp-oth3er-fuorm-nth2}
\eX \approx \llbracket \eL; \bA^{(1)}, \bA^{(2)}, \ldots, \bA^{(N)}   \rrbracket  
=
\sum_{r_1=1}^{R} \sum_{r_2=1}^{R} \ldots \sum_{r_N=1}^{R}
l_{r_1r_2\ldots r_N}
\ba_{\textcolor{mylightbluetext}{r_1}}^{(1)} \circ \ba_{\textcolor{mylightbluetext}{r_2}}^{(2)} \circ \ldots \circ \ba_{\textcolor{mylightbluetext}{r_N}}^{(N)}.
\end{equation}
Utilizing  the mode-$n$ tensor multiplication as defined in Equation~\eqref{equation:moden-tensor-multi}, the CP decomposition in Equation~\eqref{equation:cp-oth3er-fuorm-nth2} can  be alternatively  expressed as:
\begin{equation}\label{equation:cp-oth3er-fuorm-nth3}
	\eX \approx \llbracket\eL; \bA^{(1)}, \bA^{(2)}, \ldots, \bA^{(N)}   \rrbracket  
	=
	\eL \times_1 \bA^{(1)} \times_2 \bA^{(2)} \ldots \times_N \bA^{(N)}.
\end{equation}
The CP decomposition in the form of Equation~\eqref{equation:cp-oth3er-fuorm-nth3} for a third-order tensor is then illustrated in Figure~\ref{fig:cp-decom-thirr-2ndform}.

\paragraph{Storage Advantage of the CP decomposition.} When considering the first equality in Equation~\eqref{equation:cp-matriciza-third-1}, $\bX_{(1)}$ is of size $I\times (JK)$, and it admits the reduced SVD $\bX_{(1)}=\bU\bSigma\bV^\top$ the form of which requires $R(I+JK+1)$ parameters if $\bX_{(1)}$ is of rank $R$.
In contrast, the  CP decomposition in Equation~\eqref{equation:third-fact-1} requires only $R(I+J+K) \ll R(I+JK+1) $ in most cases (since it scales linearly in the mode dimensions).

\paragraph{Non-uniqueness.} The CP decomposition is not unique. Consider the form of the CP decomposition for $\eX \in \real^{I_1\times I_2\times \ldots \times I_N}$ in Equation~\eqref{equation:cp_form1}, it can be observed that, given any column permutation matrix $\bP\in\real^{R\times R}$, it follows that
$$
\eX \approx \llbracket \bA^{(1)}, \bA^{(2)}, \ldots, \bA^{(N)}\rrbracket
 = \llbracket \bA^{(1)}\bP, \bA^{(2)}\bP, \ldots, \bA^{(N)}\bP\rrbracket.
$$
This is known as the \textit{permutation inderterminancy} \citep{zhang2017matrix}.

Moreover, it also follows that 
$$
\eX \approx \sum_{r=1}^{R} \ba_r^{(1)} \circ \ba_r^{(2)} \circ \ldots \circ \ba_r^{(N)}
=\sum_{r=1}^{R} (s_r\ba_r^{(1)}) \circ (t_r\ba_r^{(2)}) \circ \ldots \circ( p_r\ba_r^{(N)}),
$$
for any $s_rt_r \ldots p_r=1$ for all $r\in\{1,2,\ldots,R\}$. This is known as the \textit{scaling indeterminancy}.

In conclusion, we define the set of c-rank of a tensor using the terms from the CP decomposition.
\begin{definition}[C-Rank of a Tensor]\label{definition:c_rank}
The tensor c-rank of $\eX$ is defined as the minimal $R$ such that $\eX$ has a CP decomposition with $R$ terms.
The set $CP(R)$ of tensors is the set of tensors that have at most  c-rank of $R$.
\end{definition}
The set $CP(R)$ is not closed, which motivates researchers to develop more consistent algorithms; see Problem~\ref{prob:close_cp_tucset}.

\section{Computing the CP Decomposition}\index{ALS}
\begin{algorithm}[h] 
\caption{CP Decomposition via ALS}
\label{alg:cp-decomposition-full-gene}
\begin{algorithmic}[1] 
\Require Nth-order tensor $\eX \in \real^{I_1\times I_2\times \ldots \times I_N}$;
\State Pick a rank $R$;
\State Initialize $\bA^{(n)}\in \real^{I_n\times R}$ for all $n\in \{1,2,\ldots, N\}$ randomly; 
\State Choose maximal number of iterations $C$;
\State $iter=0$; \Comment{Count for the number of iterations}
\While{$iter<C$}
\State $iter=iter+1$; 
\For{$n=1,2,\ldots, N$}
\State \algoalign{$\bV= \bA^{(N)\top} \bA^{(N)}\circledast \ldots \circledast \bA^{(n+1)\top} \bA^{(n+1)}\circledast \bA^{(n-1)\top} \bA^{(n-1)}\circledast \ldots \circledast \bA^{(1)\top} \bA^{(1)} \in \real^{R\times R}$;}
\State $\bW=\left(\bA^{(N)} \khatrirao \bA^{(N-1)} \khatrirao \ldots \khatrirao \bA^{(n+1)}\khatrirao \bA^{(n-1)}\khatrirao \ldots \khatrirao \bA^{(2)}\khatrirao \bA^{(1)}  \right)\in \real^{ I_{-n}\times R}$;
\State $\bX_{(n)} \in \real^{I_n\times I_{-n}}$; \Comment{Calculate the mode-$n$ matricization of tensor $\eX$}
\State $\bA^{(n)} = \bX_{(n)} \bW\bV^+\in \real^{I_n\times R}$;
\EndFor
\EndWhile
\State Output $\bA^{(1)}, \bA^{(2)}, \ldots, \bA^{(N)}$;
\end{algorithmic} 
\end{algorithm}

There is no
finite algorithm for determining the rank of a tensor, and the complexity for it is an NP-hard problem \citep{haastad1989tensor}. Consequently, the rank $R$ in the CP decomposition cannot be predetermined. However, the rank $R$ can be treated as a hyperparameter, and the CP decomposition can be approximated through low-rank tensor decomposition.

In the third-order case, $
\eX \approx \llbracket \bA, \bB, \bC   \rrbracket  
=
\sum_{r=1}^{R} \ba_r\circ \bb_r \circ \bc_r 
$, suppose we fix the tensor rank $R$, we want to find a rank-$R$ tensor $\weX$ such that
$$
\mathop{\min }_{\weX} \normf{\eX - \weX}, \qquad \text{with $\weX=\sum_{r=1}^{R} \ba_r\circ \bb_r \circ \bc_r $}.
$$
When all but one of the matrices are fixed, the problem simplifies to a least squares problem in its matricized form.  
Therefore, the alternating least squares (ALS) method can be employed.
The ALS approach proceeds as follows: first fixes $\bB$ and $\bC$ to update  $\bA$; then fixes $\bA$ and $\bC$ to update  $\bB$;
followed by fixing $\bA$ and $\bB$ to update  $\bC$; and continues to repeat the entire procedure until
a convergence criterion is met or a sufficient number of iterations  has been performed.

\paragraph{Given $\bB$ and $\bC$, update for $\bA$.}
For example, suppose that $\bB$ and $\bC$ are fixed \footnote{Again, we consider the Kolda matricization by default.}. Then we want to solve 
$$
\mathop{\min }_{\widehatbA} \normf{\bX_{(1)} - \widehatbA (\bC\khatrirao\bB)^\top}.
$$
For simplicity, we only consider the ALS without regularization and bias terms. 
Referring to the update of the ALS used in the Netflix recommender (Chapter~\ref{chapter:als}), the update for $\widehatbA$ is just the same update as that of $\bW$ in Equation~\eqref{equation:als-w-update}: given fixed matrices $\bZ$ and $\bD$, we seek $\bW$ that minimizes $\normf{\bD-\bW\bZ}^2$. If $\bZ\bZ^\top$ is invertible, the update is given by
$$
\bW = \bD\bZ^\top (\bZ\bZ^\top)^{-1}  \leftarrow \mathop{\arg\min}_{\bW} \normf{\bD-\bW\bZ}^2.
$$
Therefore, returning to the update in the third-order CP decomposition, it follows that
$$
\begin{aligned}
\widehatbA &\leftarrow \bX_{(1)} (\bC\khatrirao\bB) \left( (\bC\khatrirao\bB)^\top (\bC\khatrirao\bB)\right)^{-1}\\
&=\bX_{(1)} (\bC\khatrirao\bB)\left( (\bC^\top\bC)\circledast (\bB^\top\bB)\right)^{-1},
\end{aligned}
$$
where the last equality is derived from Equation~\eqref{equation:two-khatri-rao-pro-equi}. 
When $\left( (\bC^\top\bC)\circledast (\bB^\top\bB)\right)$ is not invertible, a pseudo-inverse can be applied instead (Appendix~\ref{appendix:pseudo-inverse}):
$$
\begin{aligned}
	\widehatbA\leftarrow \bX_{(1)} (\bC\khatrirao\bB)\left( (\bC^\top\bC)\circledast (\bB^\top\bB)\right)^{+},
\end{aligned}
$$
where $\bA^+$ represents the pseudo-inverse of  matrix $\bA$. 
The complete procedure for computing the CP decomposition of an Nth-order tensor is outlined in Algorithm~\ref{alg:cp-decomposition-full-gene}. 
The update for the \textit{nonnegative CP decomposition} follows similarly from the similarity between ALS and nonnegative matrix factorization (NMF, Chapter~\ref{chapter:nmf}). 
The ALS update can be replaced by the multiplicative update rule (see Equation~\eqref{equation:multi-update-w}).

The ALS algorithm offers primary advantages in terms of simplicity and ease of implementation.  
However, its notable drawback lies in its slow convergence, which can sometimes result in divergence. 
To address these limitations, it is essential to regularize the cost function of the  ALS algorithm, e.g., by $\ell_2$ norm.

%% file: chapter-tensortucker.tex

\index{Decomposition: Tucker}
\section{Tucker Decomposition}

Although the CP format in Theorem~\ref{theorem:cp-decomp} effectively bypasses the curse of dimensionality (the memory requirement to store such a full tensor is  reduced from $\prod_{n=1}^{N}I_n$ to $R\cdot(\sum_{n=1}^{N}I_n)$, which \textit{scales linearly} with the tensor order $N$ and size $\max\{I_1, I_2, \ldots,I_N\}$), the CP approximation can suffer from numerical instability and expressiveness challenges for high-order tensors (i.e., difficulty to arrive at a CP format). 
The corresponding algorithms for CP decompositions are often ill-posed \citep{de2008tensor}.
To address these limitations of the CP decomposition, the Tucker decomposition offers a more flexible and robust alternative. While CP decomposition represents a tensor as a sum of rank-one tensors (outer products of vectors), Tucker decomposition extends this by introducing a core tensor that captures complex interactions across modes, along with factor matrices for each mode that help capture structure within each dimension.

\begin{theoremHigh}[Tucker Decomposition \citep{tucker1966some}]\label{theorem:tucker-decomp}
The Tucker decomposition factors a tensor into a sum of rank-one
tensors. For a general Nth-order tensor, $\eX \in \real^{I_1\times I_2\times \ldots \times I_N}$, it admits the Tucker decomposition
\begin{equation}\label{equation:tucker_decom}
\eX \approx \llbracket\eG; \bA^{(1)}, \bA^{(2)}, \ldots, \bA^{(N)}   \rrbracket  
=
\sum_{r_1=1}^{R_1} \sum_{r_2=1}^{R_2} \ldots \sum_{r_N=1}^{R_N}
g_{r_1r_2\ldots r_N}
\ba_{r_1}^{(1)} \circ \ba_{r_2}^{(2)} \circ \ldots \circ \ba_{r_N}^{(N)} ,~
\footnote{With the definition of tensor multiplication in Equation~\eqref{equation:moden-tensor-multi}, the Tucker decomposition can be denoted by 
$\eX\approx \eG \times_1 \bA^{(1)} \times_2\bA^{(2)}\ldots \times_N \bA^{(N)}$.}
\end{equation}
where
\begin{itemize}
\item $R_1< I_1, R_2<I_2, \ldots, R_N<I_N$;
\item  $\eG$ of size ${R_1\times R_2\times \ldots \times R_N}$ is called the \textit{core tensor} so that $\eG$ can be thought of as a compressed version of $\eX$;
\item $\bA^{(n)}=[\ba_1^{(n)}, \ba_2^{(n)}, \ldots, \ba_{R_n}^{(n)}] \in \real^{I_n \times R_n}$ for all $n\in \{1,2,\ldots, N\}$ is the column partition of the matrix $\bA^{(n)}\in \real^{I_n \times R_n}$; $\bA^{(n)}$'s are also called \textit{mode frames} for the Tucker tensor representation;
\item The $\bA^{(n)}$'s usually have mutually orthonormal columns and can be thought of as the principal components of each mode. In this sense, the $\bA^{(n)}$'s are \textit{semi-orthogonal matrices} (see the definition in Section~\ref{section:orthogonal-orthonormal-qr}); 
\item We can complete the semi-orthogonal matrices into \textit{full orthogonal matrices} by adding \textit{silent columns} into $\bA^{(n)}$'s so that $\bA^{(n)}\in \real^{I_n\times I_n}$ becomes an orthogonal matrix. In this case, $\eG$ will be expanded to a tensor of size $\real^{I_1\times I_2\times \ldots \times I_N}$, where $g_{r_1r_2\ldots r_N} =0$ when any $r_n>R_n$ for $n\in \{1,2,\ldots, N\}$. This is known as the \textit{full Tucker decomposition}, while the previous one is also called the \textit{reduced Tucker decomposition}  to avoid confusion. We will primarily focus on the reduced case in our discussions.
\end{itemize}

\end{theoremHigh}

The advantages of Tucker decomposition over  CP decomposition are then revealed:
\begin{itemize}
\item CP decomposition can be numerically unstable for very high-order tensors due to its reliance on summing rank-one components. In contrast, Tucker decomposition avoids this instability by directly capturing relationships in the core tensor and factor matrices, making it more suitable for handling large and noisy data.
\item Although Tucker decomposition is more complex than CP, it still achieves significant memory savings. The original tensor requires storing $\prod_{n=1}^{N} I_n$ elements, while Tucker decomposition reduces this to the storage of the core tensor $\eG$ and factor matrices $\bA^{(n)}$, with a total storage cost of $\prod_{n=1}^{N} R_n + \sum_{n=1}^{N} I_n R_n$. This efficiency makes Tucker decomposition feasible for high-dimensional data, especially when the core tensor has significantly lower ranks than the original dimensions.
\item By assigning different ranks $R_n$ to each mode, Tucker decomposition enables mode-specific compression, which can adapt to the structure of the data in each dimension. This flexibility allows for more targeted compression and can be particularly useful for data with varying degrees of redundancy across modes.
\end{itemize}
While Tucker decomposition offers greater expressiveness and numerical stability, it can be more computationally demanding compared to CP decomposition. To address this issue, a hierarchical Tucker decomposition, which organizes a tensor into a binary tree structure by recursively splitting its modes, is proposed to allow for a hierarchical and efficient low-rank representation \citep{hackbusch2009new, grasedyck2010hierarchical}.

In practice, Tucker decomposition is widely used in applications where richer representations are needed, such as in image and video compression, multiway data analysis, and feature extraction for machine learning. The decomposition’s ability to maintain essential relationships across modes makes it particularly valuable for tasks requiring interpretability and robustness, such as chemometrics, neuroscience, and psychometrics \citep{bro1998multi, smilde2005multi}.

\paragraph{Equivalent forms of the Tucker decomposition.} 
The representation in \eqref{equation:tucker_decom} is known as the \textit{Tucker format}.
When the columns of $\bA^{(n)}$ has mutually orthonormal columns, it is also called the \textit{orthogonal Tucker format}.
Using the mode-$n$ tensor multiplication defined in Equation~\eqref{equation:moden-tensor-multi}, the Tucker decomposition can also be expressed as:
\begin{equation}\label{equation:tucker-decom-in-tuckeropera}
\eX \approx \llbracket\eG; \bA^{(1)}, \bA^{(2)}, \ldots, \bA^{(N)}   \rrbracket  
=
\eG \times_1 \bA^{(1)} \times_2 \bA^{(2)} \ldots \times_N \bA^{(N)}.
\end{equation}
Based on the result in Lemma~\ref{lemma:tensor-multi2}~\eqref{equationbracketlast}, since $\bA^{(n)}$'s are semi-orthogonal, it also follows that 
\begin{equation}\label{equation:tucker-decom-in-tuckeropera2}
\eG = \llbracket\eX; \bA^{(1)\top}, \bA^{(2)\top}, \ldots, \bA^{(N)\top}\rrbracket   = \eX \times_1 \bA^{(1)\top} \times_2 \bA^{(2)\top} \ldots \times_N \bA^{(N)\top}.
\end{equation}
That is, for fixed \textbf{semi-orthogonal mode frames} $\bA^{(n)}$, $\forall n\in\{1,2,\ldots,N\}$, the unique core tensor $\eG$ minimizing  $\normf{\eX - \llbracket\eG; \bA^{(1)}, \bA^{(2)}, \ldots, \bA^{(N)}   \rrbracket  }$ is $\llbracket\eX; \bA^{(1)\top}, \bA^{(2)\top}, \ldots, \bA^{(N)\top}\rrbracket$.

The operator defined in Equation~\eqref{equation:tucker-decom-in-tuckeropera} is sometimes referred to as the \textit{Tucker operator} \citep{kolda2006multilinear}.
Expressing this Tucker decomposition element-wise, the $(i_1,i_2,\ldots,i_N)$-th element of $\eX$ can be obtained by 
$$
x_{i_1,i_2,\ldots,i_N} = \sum_{r_1=1}^{R_1} \sum_{r_2=1}^{R_2} \ldots \sum_{r_N=1}^{R_N}
g_{r_1r_2\ldots r_N}
a_{i_1r_1}^{(1)} a_{i_2r_2}^{(2)} \ldots a_{i_Nr_N}^{(N)}.
$$

\paragraph{Properties of the Tucker operator.}
For an Nth-order tensor $\eX\in\real^{I_1\times I_2\times \ldots \times I_N}$,
the Tucker operator exhibits the following properties:
\begin{itemize}
\item \textit{Distributive law.} Given further two sequences of matrices $\bA^{(n)}\in\real^{ J_n \times I_n}$ and
$\bB^{(n)}\in\real^{K_n\times J_n}$ where $n\in\{1,2,\ldots, N\}$, it follows that 
$$
\begin{aligned}
 & \left\llbracket
\llbracket \eX; \bA^{(1)}, \bA^{(2)}, \ldots, \bA^{(N)}\rrbracket;
\bB^{(1)}, \bB^{(2)}, \ldots, \bB^{(N)}
\right\rrbracket\\
&\gap\gap\gap=
\llbracket \eX; \bB^{(1)}\bA^{(1)}, \bB^{(2)}\bA^{(2)}, \ldots,\bB^{(N)} \bA^{(N)}\rrbracket
\in\real^{K_1\times K_2\times \ldots \times K_N}
\end{aligned}
$$

\item \textit{Pseudo-inverse.}  Given further a sequence of matrices $\bA^{(n)}\in\real^{ J_n \times I_n}$ with full column rank, as indicated by Lemma~\ref{lemma:tensor-multi2}~\eqref{equationbracket5}, we have 
$$
\eY=\llbracket \eX; \bA^{(1)}, \bA^{(2)}, \ldots, \bA^{(N)}\rrbracket
\leadto
\eX = \llbracket \eY; \bA^{(1)+}, \bA^{(2)+}, \ldots, \bA^{(N)+}\rrbracket,
$$
where $\bA^{(n)+}$ is the pseudo-inverse of $\bA^{(n)}$.

\item \textit{Semi-orthogonal.} As mentioned earlier, given further a sequence of semi-orthogonal matrices $\bU^{(n)}\in\real^{ J_n \times I_n}$ with $\bU^{(n)\top}\bU^{(n)}=\bI$, by Lemma~\ref{lemma:tensor-multi2}~\eqref{equationbracketlast}, it follows that 
$$
\eY=\llbracket \eX; \bU^{(1)}, \bU^{(2)}, \ldots, \bU^{(N)}\rrbracket
\leadto
\eX = \llbracket \eY; \bU^{(1)\top}, \bU^{(2)\top}, \ldots, \bU^{(N)\top}\rrbracket.
$$
\end{itemize}

\begin{figure}[htp]
\centering
\resizebox{1.0\textwidth}{!}{%
\begin{tikzpicture}
\pgfmathsetmacro{\bxx}{0.3}
\pgfmathsetmacro{\begupwhite}{0.2}

\pgfmathsetmacro{\bcwd}{0.3}
\pgfmathsetmacro{\bbwd}{0.3}
\pgfmathsetmacro{\bbright}{0.9}
\pgfmathsetmacro{\bbup}{0.225}

\pgfmathsetmacro{\distrightdown}{4.2}
\pgfmathsetmacro{\dist}{4}
\pgfmathsetmacro{\distl}{0}
\pgfmathsetmacro{\aprimeright}{0.4}

\draw [very thick] (5.8+0.5,-3.8+0.5) rectangle (7+0.5,-2.6+0.5);
\filldraw [fill=gray!60!white,draw=green!40!black] (5.8+0.5,-3.8+0.5) rectangle (7+0.5,-2.6+0.5);
\draw [step=0.4/2, very thin, color=gray] (5.8+0.5,-3.8+0.5) grid (7+0.5,-2.6+0.5);

\draw [very thick] (5.8+0.4,-3.8+0.4) rectangle (7+0.4,-2.6+0.4);
\filldraw [fill=gray!50!white,draw=green!40!black] (5.8+0.4,-3.8+0.4) rectangle (7+0.4,-2.6+0.4);
\draw [step=0.4/2, very thin, color=gray] (5.8+0.4,-3.8+0.4) grid (7+0.4,-2.6+0.4);

\draw [very thick] (5.8+0.3,-3.8+0.3) rectangle (7+0.3,-2.6+0.3);
\filldraw [fill=gray!40!white,draw=green!40!black] (5.8+0.3,-3.8+0.3) rectangle (7+0.3,-2.6+0.3);
\draw [step=0.4/2, very thin, color=gray] (5.8+0.3,-3.8+0.3) grid (7+0.3,-2.6+0.3);

\draw [very thick] (5.8+0.2,-3.8+0.2) rectangle (7+0.2,-2.6+0.2);
\filldraw [fill=gray!30!white,draw=green!40!black] (5.8+0.2,-3.8+0.2) rectangle (7+0.2,-2.6+0.2);
\draw [step=0.4/2, very thin, color=gray] (5.8+0.2,-3.8+0.2) grid (7+0.2,-2.6+0.2);

\draw [very thick] (5.8+0.1,-3.8+0.1) rectangle (7+0.1,-2.6+0.1);
\filldraw [fill=gray!20!white,draw=green!40!black] (5.8+0.1,-3.8+0.1) rectangle (7+0.1,-2.6+0.1);
\draw [step=0.4/2, very thin, color=gray] (5.8+0.1,-3.8+0.1) grid (7+0.1,-2.6+0.1);

\draw [very thick] (5.8,-3.8) rectangle (7,-2.6);
\filldraw [fill=gray!10!white,draw=green!40!black] (5.8,-3.8) rectangle (7,-2.6);
\draw [step=0.4/2, very thin, color=gray] (5.8,-3.8) grid (7,-2.6);

\draw [very thick] (5.8-0.1,-3.8-0.1) rectangle (7-0.1,-2.6-0.1);
\filldraw [fill=gray!6!white,draw=green!40!black] (5.8-0.1,-3.8-0.1) rectangle (7-0.1,-2.6-0.1);
\draw [step=0.4/2, very thin, color=gray] (5.8-0.1,-3.8-0.1) grid (7-0.1,-2.6-0.1);

\draw[->] (5.6,-2.6) -- ++(0.6,0.6); 
\draw[->] (5.55,-2.7) -- ++(0,-1.2);  
\draw[->] (5.66,-4.04) -- ++(1.25,0);  

\draw (6.35,-4.25) node {{\color{black}\scriptsize{spatial $y$}}};
\draw (5.2-0.3,-3.2) node {{\color{black}\scriptsize{spatial $x$}}};
\draw (5.3,-2.2) node {{\color{black}\scriptsize{temporal}}};

\draw (6.2,-4.6) node {{\color{black}\scriptsize{$\eX\in\real^{I\times J \times K}$}}};

\draw (8,-3.2) node {	{\color{black}\large{$\approx$}}};

\draw [very thick] (8.6+\distl,-4+0.2) rectangle (9.2+\distl,-3.2+0.2);
\filldraw [fill=RubineRed!60!white,draw=green!40!black]  (8.6+\distl,-4+0.2) rectangle (9.2+\distl,-3.2+0.2);
\draw (8.7+\distl,-4.4+0.2) node {{\color{black}\scriptsize{$\bB \in \real^{J\times Q}$}}};

\draw [very thick] 
(9+\bbright+\distl,-2.9+\bbup) rectangle (9.5+\bbright+\distl,-2.6+\bbup);
\filldraw [fill=WildStrawberry!40!white,draw=green!40!black]
(9+\bbright+\distl,-2.9+\bbup) rectangle (9.5+\bbright+\distl,-2.6+\bbup);
\draw (9.6+\bbright+\distl,-3.2+\bbup) node {{\color{black}\scriptsize{$\bA^\top \in \real^{P\times I}$}}};

\draw[fill=RedOrange!40!white, line width=0.8pt] 
(8.45+\distl,-2.3-\bxx) -- (8.9+\distl,-2.1-\bxx) -- (8.35+\distl,-1.9-\bxx) -- (7.9+\distl,-2.1-\bxx) -- cycle;
\draw (8.4+\distl,-1.8) node {{\color{black}\scriptsize{$\bC\in \real^{K\times R}$}}};

\draw[fill=RedOrange!10!white, line width=0.5pt] 
(9.2+\distl,-2.4+0.1-\bxx) -- (9.4+\bcwd+\distl-0.02,-2.1+0.01-\bxx) -- (8.8+\bcwd+\distl,-2.1-\bxx) -- (8.6+\distl+0.02,-2.4+0.1+0.01-\bxx) -- cycle;
\draw[fill=RedOrange!10!white, line width=0.8pt]
(8.6+\distl,-2.4+0.1-\bxx) -- (9.2+\distl,-2.4+0.1-\bxx) -- (9.2+\distl,-3+0.4-\bxx) -- (8.6+\distl,-3+0.4-\bxx) -- cycle;
\draw[fill=RedOrange!10!white, line width=0.8pt]
(9.4+\bcwd+\distl,-2.1-\bxx) -- (9.2+\distl,-2.4+0.1-\bxx) -- (9.2+\distl,-3+0.4-\bxx) -- (9.4+\bcwd+\distl,-2.4-\bxx) -- cycle;
\draw (10-0.2+\distl,-1.8-\bxx) node {{\color{black}\scriptsize{$\eG\in \real^{P\times Q\times R}$}}};

\draw (11.6,-3.2) node {{\color{black}\large{$=$}}};

\draw [very thick] 
(8.6+\dist,-4+0.2) rectangle (9.4+\dist,-3.2+0.2);
\filldraw [fill=blue!40,draw=green!40!black] 
(8.6+\dist,-4+0.2) rectangle (9.4+\dist,-3.2+0.2);
\filldraw [fill=blue!40,draw=green!40!black, pattern=north east lines] 
(8.6+\dist,-4+0.2) rectangle (9.4+\dist,-3.2+0.2);

\draw [very thick] (8.6+\dist,-4+0.2) rectangle (9.2+\dist,-3.2+0.2);
\filldraw [fill=RubineRed!60!white,draw=green!40!black]  (8.6+\dist,-4+0.2) rectangle (9.2+\dist,-3.2+0.2);
\draw (8.7+\dist,-4.4+0.2) node {{\color{black}\scriptsize{$\bB^\prime \in \real^{J\times \textcolor{mylightbluetext}{J}}$}}};

\draw [very thick] 
(9+\bbright+\dist+\aprimeright,-2.8+\bbup) rectangle (9.5+\bbright+\dist+\aprimeright,-2.3+\bbup);
\filldraw [fill=blue!40,draw=green!40!black] 
(9+\bbright+\dist+\aprimeright,-2.8+\bbup) rectangle (9.5+\bbright+\dist+\aprimeright,-2.3+\bbup);
\filldraw [fill=blue!40,draw=green!40!black, pattern=north east lines] 
(9+\bbright+\dist+\aprimeright,-2.8+\bbup) rectangle (9.5+\bbright+\dist+\aprimeright,-2.3+\bbup);

\draw [very thick] 
(9+\bbright+\dist+\aprimeright,-2.8+\bbup) rectangle (9.5+\bbright+\dist+\aprimeright,-2.5+\bbup);
\filldraw [fill=WildStrawberry!40!white,draw=green!40!black]
(9+\bbright+\dist+\aprimeright,-2.8+\bbup) rectangle (9.5+\bbright+\dist+\aprimeright,-2.5+\bbup);
\draw (9.6+\bbright+\dist+\aprimeright,-3.1+\bbup) node {{\color{black}\scriptsize{$\bA^{\prime\top }\in \real^{ \textcolor{mylightbluetext}{I}\times I}$}}};

\draw[fill=blue!40, line width=0.8pt] 
(8.45+\dist,-2.3-\bxx+\begupwhite) -- (8.9+\dist+0.3,-2.1-\bxx+\begupwhite+0.14) -- (8.35+\dist+0.3,-1.9-\bxx+\begupwhite+0.15) -- (7.9+\dist,-2.1-\bxx+\begupwhite) -- cycle;
\draw[fill=blue!40, line width=0.8pt, pattern=north east lines] 
(8.45+\dist,-2.3-\bxx+\begupwhite) -- (8.9+\dist+0.3,-2.1-\bxx+\begupwhite+0.14) -- (8.35+\dist+0.3,-1.9-\bxx+\begupwhite+0.15) -- (7.9+\dist,-2.1-\bxx+\begupwhite) -- cycle;

\draw[fill=RedOrange!40!white, line width=0.8pt] 
(8.45+\dist,-2.3-\bxx+\begupwhite) -- (8.9+\dist,-2.1-\bxx+\begupwhite) -- (8.35+\dist,-1.9-\bxx+\begupwhite) -- (7.9+\dist,-2.1-\bxx+\begupwhite) -- cycle;
\draw (8.4+\dist,-1.8+\begupwhite) node {{\color{black}\scriptsize{$\bC^\prime\in \real^{K\times \textcolor{mylightbluetext}{K}}$}}};

\draw[fill=RedOrange!1!white, line width=0.8pt,opacity=0.5] 
(8.6+\dist,-2.4+0.1-\bxx) -- (9.4+\dist,-2.4+0.1-\bxx) -- (9.4+\dist,-3+0.4-\bxx) -- (8.6+\dist,-3+0.4-\bxx) -- cycle;

\draw[fill=RedOrange!10!white, line width=0.8pt] 
(9.2+\dist,-2.4+0.1-\bxx+\begupwhite) -- (9.4+\bcwd+\dist,-2.1-\bxx+\begupwhite) -- (8.8+\bcwd+\dist+0.04,-2.1-\bxx+\begupwhite+0.01) -- (8.6+\dist,-2.4+0.1-\bxx+\begupwhite) -- cycle;
\draw[fill=RedOrange!10!white, line width=0.8pt]
(8.6+\dist,-2.4+0.1-\bxx+\begupwhite) -- (9.2+\dist,-2.4+0.1-\bxx+\begupwhite) -- (9.2+\dist,-3+0.4-\bxx+\begupwhite) -- (8.6+\dist,-3+0.4-\bxx+\begupwhite) -- cycle;
\draw[fill=RedOrange!10!white, line width=0.8pt]
(9.4+\bcwd+\dist,-2.1-\bxx+\begupwhite) -- (9.2+\dist,-2.4+0.1-\bxx+\begupwhite) -- (9.2+\dist,-3+0.4-\bxx+\begupwhite) -- (9.4+\bcwd+\dist,-2.4-\bxx+\begupwhite) -- cycle;
\draw (10+\dist,-1.7-\bxx+\begupwhite) node {{\color{black}\scriptsize{$\eG^\prime\in \real^{\textcolor{mylightbluetext}{I\times J\times K}}$}}};


\draw[fill=RedOrange!1!white, line width=0.8pt,opacity=0.5] 
(9.4+\dist,-2.4+0.1-\bxx+\begupwhite) -- (9.6+\bcwd+\dist+0.24,-2.1-\bxx+\begupwhite+0.12) -- (8.8+\bcwd+\dist+0.26,-2.1-\bxx+\begupwhite+0.12) -- (8.6+\dist,-2.4+0.1-\bxx+\begupwhite) -- cycle;
\draw[fill=RedOrange!1!white, line width=0.8pt,opacity=0.5] 
(9.6+\bcwd+\dist+0.25,-2.1-\bxx+\begupwhite+0.12) -- (9.4+\dist,-2.4+0.1-\bxx+\begupwhite) -- (9.4+\dist,-3+0.4-\bxx) -- (9.6+\bcwd+\dist+0.25,-2.4-\bxx+0.12) -- cycle;

\end{tikzpicture}
}
\caption{The Tucker decomposition of a third-order tensor: $\eX \approx \llbracket \eG; \bA, \bB, \bC   \rrbracket  
	=
	\sum_{p=1}^{P}\sum_{q=1}^{Q}\sum_{r=1}^{R}
	g_{pqr} \cdot \ba_p\circ \bb_q \circ \bc_r$. Middle: \textbf{reduced Tucker decomposition}. Right: \textbf{full Tucker decomposition}, where the \textcolor{mylightbluetext}{blue hatched} entries represent  \textit{silent orthogonal} columns, and the white entries of $\eG^\prime$ are zero.}
\label{fig:tucker-decom-third}
\end{figure}

\paragraph{Matricization.}
For simplicity, we consider the Tucker decomposition for the third-order tensor $\eX\in \real^{I\times J\times K}$, as illustrated in Figure~\ref{fig:tucker-decom-third}:
\begin{equation}\label{equation:tucker-third-fact-1}
\eX \approx \llbracket \eG; \bA, \bB, \bC   \rrbracket  
=
\sum_{p=1}^{P}\sum_{q=1}^{Q}\sum_{r=1}^{R}
g_{pqr} \cdot \ba_p\circ \bb_q \circ \bc_r,
\end{equation}
where $\bA\in \real^{I\times P}, \bB\in \real^{J\times Q}$, $\bC\in \real^{K\times R}$, and $\eG\in \real^{P\times Q\times R}$. Analogously, following the (Kolda) matricization of the third-order CP decomposition in Equation~\eqref{equation:cp-matriciza-third-1}, we have 
\begin{equation}\label{equation:tucker-matriciza-third-1}
	\left\{
	\begin{aligned}
		\bX_{(1)}  &\approx \bA\bG_{(1)} (\bC\kronecker \bB )^\top \in \real^{I\times (JK)}; \\
		\bX_{(2)} &\approx \bB\bG_{(2)} (\bC\kronecker \bA )^\top \in \real^{J\times (IK)};\\
		\bX_{(3)} &\approx \bC\bG_{(3)}(\bB\kronecker \bA)^\top \in \real^{K\times (IJ)},
	\end{aligned}
	\right.
\end{equation}
where the symbol ``$\kronecker$" represents the Kronecker product (Definition~\ref{definition:kronecker-product}).

\paragraph{Matricization of general Nth-order tensors.}
In full generality, for an Nth-order tensor $\eX \in \real^{I_1\times I_2\times \ldots \times I_N}$, the mode-$n$ matricized form of the Tucker decomposition is given by 
$$
\boxed{\underbrace{\bX_{(n)}}_{I_n\times (I_{-n})}\approx
	\underbrace{\bA^{(n)}}_{{I_n\times R_n}} 
	\underbrace{\bG_{(n)}}_{R_n\times (R_{-n})}
	\underbrace{\left(\bA^{(N)} \kronecker \bA^{(N-1)} \kronecker \ldots \kronecker \bA^{(n+1)}\kronecker \bA^{(n-1)}\kronecker \ldots \kronecker \bA^{(2)}\kronecker \bA^{(1)}  \right)^\top}_{(R_{-n})\times (I_{-n})}}
$$
where $I_{-n}=I_1 I_2\ldots I_{n-1}I_{n+1}\ldots I_N$ and $R_{-n}=R_1 R_2\ldots R_{n-1}R_{n+1}\ldots R_N$.
It then follows that $\rank(\bX_{(n)})\leq R_n$.
This motivates the definition of the \textit{multilinear rank} (t-rank) of a tensor $\eX$ as the tuple
\begin{equation}\label{equation:t_rank}
\textbf{Multilinear rank:} 
(\overline{R}_1, \overline{R}_2, \ldots, \overline{R}_N),
\,\, \text{where} \,\, 
\overline{R}_n = \rank(\bX_{(n)}), \forall n\in\{1,2,\ldots,N\}.
\end{equation}

\begin{definition}[Multilinear Rank (T-Rank) of a Tensor]\label{definition:t_rank}
The tensor t-rank of $\eX \in \real^{I_1\times I_2\times \ldots \times I_N}$ is defined as the tuple in \eqref{equation:t_rank}.
The set $Tucker(\overline{R}_1, \overline{R}_2, \ldots, \overline{R}_N)$ of tensors is the set of tensors that have at most t-rank of $(\overline{R}_1, \overline{R}_2, \ldots, \overline{R}_N)$.
\end{definition}

\paragraph{Tucker truncation.} The definition of multilinear rank motivates the Tucker truncation \citep{de2000multilinear}. Consider the Tucker format in Theorem~\ref{theorem:tucker-decomp} and the SVD of mode-$n$ matricization:
$$
\bX_{(n)} = \bU_n\bSigma_n\bV_n^\top, 
\gap\text{where}\gap 
\bU_n\in\real^{I_n\times I_n}.
$$ 
Then, the truncation of $\eX$ to Tucker rank $(\overline{R}_1, \overline{R}_2, \ldots, \overline{R}_N)$ (\textit{Tucker truncation}) is given by 
\begin{equation}\label{equation:tucker_truncation}
\mathcalT_{\overline{R}_1, \overline{R}_2, \ldots, \overline{R}_N}(\eX) 
=
\llbracket\eX; \widetildebA_1\widetildebA_1^\top, \widetildebA_2\widetildebA_2^\top, \ldots, \widetildebA_N\widetildebA_N^\top   \rrbracket 
 \in \real^{I_1\times I_2\times \ldots \times I_N}
\end{equation}
where $\widetildebA_n\in\real^{I_n\times \overline{R}_n}$ is the matrix of the first $\overline{R}_n$ columns of $\bU_n$: $\widetildebA_n \triangleq \bU_n[:,:\overline{R}_n]$.
Tucker truncation transforms a (non-orthogonal) Tucker format into an orthogonal Tucker format. The core tensor is uniquely determined after this truncation, and this form is related to the high-order SVD (HOSVD, Section~\ref{section:hosvd}).

\paragraph{Vectorization.} Going further from the matricization of the Nth-order tensor $\eX$, the vectorization of the Tucker decomposition is given by 
\begin{equation}\label{equation:tucker-vec-in-theorem}
\boxed{\underbrace{vec(\eX) }_{(I_1\ldots I_N)\times 1}
	\approx 
	\underbrace{\left(\bA^{(N)} \kronecker \bA^{(N-1)} \kronecker \ldots\kronecker \bA^{(1)}  \right)}_{(I_1\ldots I_N)\times (R_1\ldots R_N)}
	\underbrace{vec(\eG)}_{(R_1\ldots R_N)\times 1}
}.
\end{equation}

\paragraph{Counterpart in $\eG$.} For the problem in Equation~\eqref{equation:tucker-decom-in-tuckeropera2}, it admits the matricization and vectorization in a similar way:
$$
\boxed{
\begin{aligned}
\bG_{(n)}\approx
\bA^{(n)\top}
\bG_{(n)}\left(\bA^{(N)\top} \kronecker \bA^{(N-1)\top} \kronecker \ldots \kronecker \bA^{(n+1)\top}\kronecker \bA^{(n-1)\top}\kronecker \ldots \kronecker \bA^{(2)\top}\kronecker \bA^{(1)\top}  \right)^\top
\end{aligned}
}
$$
and 
\begin{equation}\label{equation:tucker-vec-in-theorem-eg}
\boxed{\underbrace{vec(\eG) }_{(R_1\ldots R_N)\times 1 }
\approx 
\underbrace{\left(\bA^{(N)\top} \kronecker \bA^{(N-1)\top} \kronecker \ldots\kronecker \bA^{(1)\top}  \right)}_{ (R_1\ldots R_N)  \times (I_1\ldots I_N)}
\underbrace{vec(\eX)}_{(I_1\ldots I_N)\times 1}
}.
\end{equation}


\section{Computing the Tucker Decomposition}

To compute the Tucker decomposition of $\eX \in \real^{I_1\times I_2\times \ldots \times I_N}$, we  consider algorithms for solving the problem:\index{ALS}
$$
\boxed{\{{\eG,\bA^{(1)}, \bA^{(2)}, \ldots, \bA^{(N)}}\}
	=
	\mathop{\arg \gap \min}_{\eG,\bA^{(1)}, \bA^{(2)}, \ldots, \bA^{(N)}}
	\left\Vert\eX - \llbracket\eG; \bA^{(1)}, \bA^{(2)}, \ldots, \bA^{(N)}   \rrbracket  \right\Vert_F^2},
$$
where $\eG$ is a tensor of size ${R_1\times R_2\times \ldots \times R_N}$, and $\bA^{(n)}$'s are semi-orthogonal of size ${I_n\times R_n}$ for $n\in \{1,2,\ldots, N\}$ (we focus on the orthogonal Tucker format for simplicity). Similar to the CP decomposition, an \textit{alternating descent} algorithm can be employed to find the solution approximately.

\subsubsection*{\textbf{Given $\bA^{(1)}, \bA^{(2)}, \ldots, \bA^{(N)}$, Update $\eG$:}} 
The update for $\eG$ is straightforward, as derived from Equation~\eqref{equation:tucker-decom-in-tuckeropera2} that 
$$
\eG = \llbracket\eX; \bA^{(1)\top}, \bA^{(2)\top}, \ldots, \bA^{(N)\top}\rrbracket   = \eX \times_1 \bA^{(1)\top} \times_2 \bA^{(2)\top} \ldots \times_N \bA^{(N)\top}.
$$
To simplify matters, using the vectorized form from Equation~\eqref{equation:tucker-vec-in-theorem}, we have:
$$
\begin{aligned}
&\left\Vert\eX - \llbracket\eG; \bA^{(1)}, \bA^{(2)}, \ldots, \bA^{(N)}   \rrbracket  \right\Vert_F^2
=
\left\Vert 
vec(\eX)
-
\left(\bA^{(N)} \kronecker \bA^{(N-1)} \kronecker \ldots\kronecker \bA^{(1)}  \right)
vec(\eG)
\right\Vert_F^2,
\end{aligned}
$$
which is essentially a least squares problem. The solution is given by
\begin{equation}\label{equation:tucker-eg-vector-undo}
	\begin{aligned}
		vec(\eG) &\leftarrow \left(\bA^{(N)} \kronecker \bA^{(N-1)} \kronecker \ldots\kronecker \bA^{(1)}  \right)^+ vec(\eX)\\
		&=\left(\bA^{(N)\top} \kronecker \bA^{(N-1)\top} \kronecker \ldots\kronecker \bA^{(1)\top}  \right)vec(\eX),
	\end{aligned}
\end{equation}
where the last equality follows from Equation~\eqref{equation:otimes-psesudo-semi} since $\bA^{(n)}$'s are semi-orthogonal. 
Once we have completed the update for the vectorized version, an \textit{un-vectorization operation} can be applied to find the updated tensor $\eG$.

\subsubsection*{\textbf{Given $\bA^{(1)},\ldots, \bA^{(n-1)}, \bA^{(n+1)}, \ldots, \bA^{(N)}$ and $\eG$, Update $\bA^{(n)}$:}}
When all but $\bA^{(n)}$ are fixed, as indicated by Equation~\eqref{equation:differe-tensor-norm}, the Frobenius norm of the difference is given by 
$$
\begin{aligned}
&\gap \left\Vert\eX - \llbracket\eG; \bA^{(1)}, \bA^{(2)}, \ldots, \bA^{(N)}   \rrbracket  \right\Vert_F^2\\
&=\left\Vert\eX \right\Vert_F^2
-2\big\langle \eX,\llbracket\eG; \bA^{(1)}, \bA^{(2)}, \ldots, \bA^{(N)}   \rrbracket\big\rangle 
+ \left\Vert\llbracket\eG; \bA^{(1)}, \bA^{(2)}, \ldots, \bA^{(N)}   \rrbracket\right\Vert_F^2. 
\end{aligned}
$$
Using Lemma~\ref{lemma:tensor-multi2} (\ref{equationbracket}), Equation~\eqref{equation:tucker-decom-in-tuckeropera}, and Equation~\eqref{equation:tucker-decom-in-tuckeropera2}, we can deduce that
$$
\begin{aligned}
&\gap\big\langle \eX,\llbracket\eG; \bA^{(1)}, \bA^{(2)}, \ldots, \bA^{(N)}   \rrbracket\big\rangle\\
&=
\big\langle \llbracket\eX; \bA^{(1)\top}, \bA^{(2)\top}, \ldots, \bA^{(N)\top}   \rrbracket, \eG\big\rangle \gap &\text{(Lemma~\ref{lemma:tensor-multi2} (\ref{equationbracket}), Equation~\eqref{equation:tucker-decom-in-tuckeropera})}\\
&= \langle \eG, \eG \rangle=\left\Vert\eG \right\Vert_F^2,  &\gap \text{(Equation~\eqref{equation:tucker-decom-in-tuckeropera2})}
\end{aligned}
$$
Also, from Equation~\eqref{equation:tucker-decom-in-tuckeropera} and Lemma~\ref{lemma:tensor-multi2}~\eqref{equationbracketlast}, and from the fact that the length is preserved under semi-orthogonal, we obtain 
$$
\left\Vert\llbracket\eG; \bA^{(1)}, \bA^{(2)}, \ldots, \bA^{(N)},   \rrbracket\right\Vert_F^2 = \norm{\eG}_F^2.
$$
Combining all the findings, we have
$$
\begin{aligned}
\left\Vert\eX - \llbracket\eG; \bA^{(1)}, \bA^{(2)}, \ldots, \bA^{(N)}   \rrbracket  \right\Vert_F^2
=\left\Vert\eX \right\Vert_F^2- \left\Vert\eG \right\Vert_F^2.
\end{aligned}
$$
Hence, minimizing the left-hand of the above equation is equivalent to maximizing $\left\Vert\eG \right\Vert_F^2$. Therefore, to update $\bA^{(n)}$, the problem becomes 
$$
\boxed{\begin{aligned}
		&\mathop{\max}_{\bA^{(n)}} \,\,\,\,\eG = \mathop{\max}_{\bA^{(n)}} \,\,\,\, \llbracket\eX; \bA^{(1)\top}, \bA^{(2)\top}, \ldots,\bA^{(n-1)\top},\bA^{(n)\top},\bA^{(n+1)\top},\ldots,  \bA^{(N)\top}\rrbracket\\
		&\text{subject to \gap  $\bA^{(n)} \in \real^{I_n\times R_n}$ is semi-orthogonal.} 
\end{aligned} }
$$
Define $\eY\triangleq \llbracket\eX; \bA^{(1)\top}, \bA^{(2)\top}, \ldots,\bA^{(n-1)}, \textcolor{mylightbluetext}{\bI},\bA^{(n+1)},\ldots,  \bA^{(N)\top}\rrbracket$. Using Lemma~\ref{lemma:tensor-multi-matriciz1}, this,  once again, is equivalent to finding the solution of 
$$
\boxed{
\begin{aligned}
&\mathop{\max}_{\bA^{(n)}} \,\,\,\,\normf{\eY \times_n \bA^{(n)}}= \mathop{\max}_{\bA^{(n)}} \,\,\,\, \normf{\bA^{(n)}\bY_{(n)}}\\
&\text{subject to \gap  $\bA^{(n)} \in \real^{I_n\times R_n}$ is semi-orthogonal,} 
\end{aligned} 
}
$$
where $\bY_{(n)} \in \real^{R_n\times R_{1}\ldots R_{n-1}R_{n+1}\ldots R_N}$ is the matricization along the mode-$n$ of $\eY$. The solution can be obtained by setting the rows of $\bA^{(n)}$ to the first $R_n$ leading left singular vectors of $\bY_{(n)}$ \citep{kolda2006multilinear, kolda2009tensor}. 
We notice that the  update on $\bA^{(n)}$ mentioned above is independent of $\eG$, allowing us to update the tensor $\eG$ once after the convergence of $\bA^{(n)}$'s. The procedure is outlined in Algorithm~\ref{alg:tucker-decomposition-full-gene}.

\begin{algorithm}[h] 
	\caption{Tucker Decomposition via ALS}
	\label{alg:tucker-decomposition-full-gene}
	\begin{algorithmic}[1] 
		\Require Nth-order tensor $\eX \in \real^{I_1\times I_2\times \ldots \times I_N}$;
		\State Pick a rank $R_1, R_2, \ldots, R_N$;
		\State Initialize $\bA^{(n)}\in \real^{I_n\times R_n}$ for all $n\in \{1,2,\ldots, N\}$ randomly; 
		\State Choose maximal number of iterations $C$;
		\State $iter=0$; \Comment{Count for the number of iterations}
		\While{$iter<C$}
		\State $iter=iter+1$; 
		\For{$n=1,2,\ldots, N$}
		\State Set $\eY= \llbracket\eX; \bA^{(1)\top}, \bA^{(2)\top}, \ldots,\bA^{(n-1)}, \textcolor{mylightbluetext}{\bI},\bA^{(n+1)},\ldots,  \bA^{(N)\top}\rrbracket$;
		\State Find the matricization along mode-$n$: $\bY_{(n)}$;
		\State Set the rows of $\bA^{(n)}$ to first $R_n$ leading left singular vectors of $\bY_{(n)}$;
		\EndFor
		\EndWhile
		\State $\eG = \llbracket\eX; \bA^{(1)\top}, \bA^{(2)\top}, \ldots, \bA^{(N)\top}\rrbracket$;\Comment{by vectorize and un-vectorize, Eq.~\eqref{equation:tucker-eg-vector-undo}}
		\State Output $\bA^{(1)}, \bA^{(2)}, \ldots, \bA^{(N)}, \eG$; 
	\end{algorithmic} 
\end{algorithm}

\paragraph{Initialization by High-Order SVD (HOSVD).}\index{HOSVD}
In Algorithm~\ref{alg:tucker-decomposition-full-gene}, the matrices $\bA^{(n)}$'s are initialized randomly. 
However, the high-order SVD can be utilized as an initial point. The problem can be formulated as follows: 
$$
\boxed{\begin{aligned}
&\mathop{\max}_{\bA^{(n)}} \,\,\,\, \llbracket\eX; \bI, \bI, \ldots,\bI,\bA^{(n)\top},\bI,\ldots,  \bI\rrbracket\\
&\text{subject to \gap  $\bA^{(n)} \in \real^{I_n\times R_n}$ is semi-orthogonal.} 
\end{aligned} }
$$
Similarly, the problem can be equivalently formulated as finding the solution of:
\begin{equation}\label{equation:hosvd-initia-for-tucker-0}
\boxed{
\begin{aligned}
&\,\mathop{\max}_{\bA^{(n)}}\,\,\,\, \normf{\eX\times_n \bA^{(n)\top}}=\normf{\bA^{(n)\top} \bX_{(n)}}\\
&\,\text{subject to \gap  $\bA^{(n)} \in \real^{I_n\times R_n}$ is semi-orthogonal,} 
\end{aligned} 
}
\end{equation}
which is solved by setting the columns of $\bA^{(n)}$ to the first $R_n$ leading left singular vectors of $\bX_{(n)}$. 
An equivalent derivation using the matricization of the HOSVD will be presented in the next section; refer to Equation~\eqref{equation:hosvd-initia-for-tucker-1}.

%% file: chapter-tensorhosvd.tex

\index{Decomposition: HOSVD}
\section{High-Order SVD (HOSVD)}\label{section:hosvd}
We previously discussed the utilization of high-order SVD (HOSVD) as an initialization step  for calculating the Tucker decomposition. We now explore the properties of  HOSVD.

\begin{theoremHigh}[High-Order SVD (HOSVD) \citep{de2000multilinear}]\label{theorem:hosvd-decomp}
The HOSVD factors a tensor into a sum of rank-one
tensors. For a general Nth-order tensor, $\eX \in \real^{I_1\times I_2\times \ldots \times I_N}$, it admits the HOSVD
\begin{equation}\label{equation:hosvd_decom}
\eX \approx \llbracket\eG; \bA^{(1)}, \bA^{(2)}, \ldots, \bA^{(N)}   \rrbracket  
=
\sum_{r_1=1}^{R_1} \sum_{r_2=1}^{R_2} \ldots \sum_{r_N=1}^{R_N}
g_{r_1r_2\ldots r_N}
\ba_{r_1}^{(1)} \circ \ba_{r_2}^{(2)} \circ \ldots \circ \ba_{r_N}^{(N)},~
\footnote{With the definition of tensor multiplication in Equation~\eqref{equation:moden-tensor-multi}, the HOSVD can be denoted by 
	$\eX\approx \eG \times_1 \bA^{(1)} \times_2\bA^{(2)}\ldots \times_N \bA^{(N)}$.}
\end{equation}
where
\begin{itemize}
\item $R_1<  I_1, R_2< I_2, \ldots, R_N< I_N$;
\item  $\eG$ of size ${R_1\times R_2\times \ldots \times R_N}$ is called the \textit{core tensor} so that $\eG$ can be thought of as a compressed version of $\eX$;
\item $\bA^{(n)}=[\ba_1^{(n)}, \ba_2^{(n)}, \ldots, \ba_{R_n}^{(n)}] \in \real^{I_n \times R_n}$ for all $n\in \{1,2,\ldots, N\}$ is the column partition of the matrix $\bA^{(n)}\in \real^{I_n \times R_n}$;  $\bA^{(n)}$'s are also called \textit{mode frames}.
\item The $\bA^{(n)}$'s have mutually orthonormal columns and can be thought of as the principal component of each mode. In this sense, the $\bA^{(n)}$'s are \textit{semi-orthogonal matrices} (see the definition in Section~\ref{section:orthogonal-orthonormal-qr});
\item We can complete the semi-orthogonal matrices into \textit{full orthogonal matrices} by adding \textit{silent columns} into $\bA^{(n)}$'s so that $\bA^{(n)}\in \real^{I_n\times I_n}$ becomes an orthogonal matrix. 
In this case, $\eG$ will be expanded to a tensor of size $\real^{I_1\times I_2\times \ldots \times I_N}$, where $g_{r_1r_2\ldots r_N} =0$ when any  $r_n>R_n$ for $n\in \{1,2,\ldots, N\}$. This is known as the \textit{full HOSVD}, while the previous one is also called the \textit{reduced HOSVD} to avoid confusion; and we shall only consider the reduced case in most of our discussions.
\end{itemize}
\noindent
\textbf{Up to this point, HOSVD is essentially equivalent to the Tucker decomposition}. However, they differ in the following aspects (we will use the full HOSVD as an example, and the reduced case can be inferred from the context):
\begin{itemize}
\item \textit{All orthogonality.} The slices in each mode are mutually orthogonal. Let \textit{subtensor} $\eG_{r_n=\alpha}$ denote the slice of $\eG$, where the $n$-th index is set to $\alpha$, i.e., an $(N-1)$th-order subtensor: $\eG_{r_n=\alpha}\triangleq\eG_{:,\ldots,:,r_n=\alpha,:, \ldots,:}$, then it follows that 
\begin{equation}\label{equation:hosvd-conditio1}
\langle  \eG_{r_n=\alpha},  \eG_{r_n=\beta}\rangle =0, \gap \alpha\neq \bbeta \in \{1,2,\ldots, I_n\},
\end{equation}
for all possible values of $n\in \{1,2,\ldots, N\}$.
\item \textit{Ordering.} The Frobenius norms of slices in each mode  decrease as the running index increases:
\begin{equation}\label{equation:hosvd-conditio2}
\underbrace{\normf{\eG_{r_n=1}}}_{\triangleq\sigma_i^{(n)}}
\geq \underbrace{\normf{\eG_{r_n=2}} }_{\triangleq\sigma_2^{(n)}}
\geq \ldots 
\geq \underbrace{\normf{\eG_{r_n=I_n}} }_{\triangleq\sigma_{I_n}^{(n)}}
\geq 0,
\end{equation}
for all possible values of $n\in \{1,2,\ldots, N\}$.
\end{itemize}
\end{theoremHigh}
The HOSVD is extensively used in applications of face recognition. While SVD is a popular technique in recognizing faces, referred to as ``EigenFaces" \citep{chellappa1995human, elden2007matrix},
it may not perform well when several environmental factors vary.
More recently, the HOSVD is applied via a tensor model, known as the ``TensorFaces" approach \citep{vasilescu2002human, vasilescu2002multilinear, elden2007matrix}, which improves the precision of the recognition.
\index{Eigenface}
\index{Tensorface}

The Frobenius norm $\normf{\eG_{r_n=i}}$ is usually denoted as $\normf{\eG_{r_n=i}} \triangleq \sigma_i^{(n)}$, and is known as the \textit{mode-$n$ singular values of $\eX$}. And the vector $\ba_i^{(n)}$ is the \textit{$i$-th mode-$n$ singular vector}. 
The comparison between the matrix and high-order (tensor) SVD reveals a clear analogy between
the two cases. In matrix language, consider the reduced SVD of $\bA=\bU\bSigma\bV^\top\in \real^{M\times N}$, the singular values in matrix $\bSigma\in \real^{R\times R}$ correspond to  the norms of the rows or columns of $\bSigma$ (since $\bSigma$ is diagonal). 
The left and right singular vectors in the columns of $\bU$ and $\bV$, respectively, can now be generalized to the mode-$n$ singular vectors of tensors.

\paragraph{Equivalent forms of the HOSVD.} 
The equivalent forms of the HOSVD are identical to those of the Tucker decomposition. 
Given an Nth-order tensor $\eX \in \real^{I_1\times I_2\times \ldots \times I_N}$, using the mode-$n$ tensor multiplication in Equation~\eqref{equation:moden-tensor-multi}, the HOSVD can also be expressed as:
\begin{equation}\label{equation:tucker-decom-in-hosvd}
	\eX \approx \llbracket\eG; \bA^{(1)}, \bA^{(2)}, \ldots, \bA^{(N)}   \rrbracket  
	=
	\eG \times_1 \bA^{(1)} \times_2 \bA^{(2)} \ldots \times_N \bA^{(N)}.
\end{equation}
Note the analogy with the SVD for a matrix in Equation~\eqref{equation:svd-by-tensor-multi0}:
\begin{equation}\label{equation:svd-by-tensor-multi1}
	\bA = \bU\bSigma\bV^\top = \bSigma\times_1 \bU \times_2 \bV\in \real^{M\times N}.
\end{equation}
That is, if the diagonal singular value matrix $\bSigma$ is regarded as a second-order tensor, then $\bSigma$ serves as the core tensor of the second-order tensor $\eA=\bA$.

Based on the result in Lemma~\ref{lemma:tensor-multi2}~\eqref{equationbracketlast}, since $\bA^{(n)}$'s are semi-orthogonal such that $\bA^{(n)\top}\bA^{(n)}=\bI$, it also follows that 
\begin{equation}\label{equation:tucker-decom-in-hosvd2}
	\eG = \llbracket\eX; \bA^{(1)\top}, \bA^{(2)\top}, \ldots, \bA^{(N)\top}\rrbracket   = \eX \times_1 \bA^{(1)\top} \times_2 \bA^{(2)\top} \ldots \times_N \bA^{(N)\top}.
\end{equation}

Expressing the HOSVD element-wise, the $(i_1,i_2,\ldots,i_N)$-th element of $\eX$ can be obtained by 
$$
x_{i_1,i_2,\ldots,i_N} = \sum_{r_1=1}^{R_1} \sum_{r_2=1}^{R_2} \ldots \sum_{r_N=1}^{R_N}
g_{r_1r_2\ldots r_N}
a_{i_1r_1}^{(1)} a_{i_2r_2}^{(2)} \ldots a_{i_Nr_N}^{(N)}.
$$

\paragraph{Matricization.}
Similarly, for an Nth-order tensor $\eX \in \real^{I_1\times I_2\times \ldots \times I_N}$, the mode-$n$ (Kolda) matricized form of the HOSVD is given by 
$$
\boxed{
\underbrace{\bX_{(n)}}_{I_n\times (I_{-n})}
\approx
\underbrace{\bA^{(n)}}_{{I_n\times R_n}} 
\underbrace{\bG_{(n)}}_{R_n\times (R_{-n})}
\underbrace{\left(\bA^{(N)} \kronecker \bA^{(N-1)} \kronecker \ldots \kronecker \bA^{(n+1)}\kronecker \bA^{(n-1)}\kronecker \ldots \kronecker \bA^{(2)}\kronecker \bA^{(1)}  \right)^\top}_{(R_{-n})\times (I_{-n})}
}
$$
where $I_{-n}=I_1 I_2\ldots I_{n-1}I_{n+1}\ldots I_N$ and $R_{-n}=R_1 R_2\ldots R_{n-1}R_{n+1}\ldots R_N$. 
Moreover, the property of \textbf{all orthogonality} becomes apparent when examining the matricization form.
Following the conditions in Equation~\eqref{equation:hosvd-conditio1} and \eqref{equation:hosvd-conditio2}, this implies $\bG_{(n)}$ has mutually orthogonal rows (not necessarily orthonormal), with vector Frobenius norms equal to $\sigma_1^{(n)}, \sigma_2^{(n)},\ldots, \sigma_{R_n}^{(n)}$, respectively (in the reduced HOSVD scenario). 

\paragraph{Matrix SVD in HOSVD.}
Construct x the diagonal matrix (in the reduced HOSVD scenario)
$$
\begin{aligned}
\bSigma^{(n)} \triangleq \diag\left(\sigma_1^{(n)}, \sigma_2^{(n)}, \ldots, \sigma_{R_n}^{(n)}\right) \in \real^{R_n\times R_n},
\end{aligned}
$$
where $\sigma_i^{(n)} = \norm{\eG_{r_n=i}}_F$. This implies, for the row normalized version $\widetildebG_{(n)}$ of $\bG_{(n)}$, we have 
$$
\underbrace{\bG_{(n)}}_{R_n\times R_{-n}}= 
\underbrace{\bSigma^{(n)}}_{R_n\times R_n} \underbrace{\widetildebG_{(n)}}_{R_n\times R_{-n}}.
$$
Let further, 
$$
\underbrace{\bV^{(n)\top}}_{R_n \times I_{-n}}
 \triangleq
\underbrace{\widetildebG_{(n)}}_{R_n\times R_{-n}}
\underbrace{\left(\bA^{(N)} \kronecker \bA^{(N-1)} \kronecker \ldots \kronecker \bA^{(n+1)}\kronecker \bA^{(n-1)}\kronecker \ldots \kronecker \bA^{(2)}\kronecker \bA^{(1)}  \right)^\top}_{(R_{-n})\times (I_{-n})},
$$
where the columns of $\bV^{(n)}$ are mutually orthonormal , as shown by Equation~\eqref{equation:orthogonal-in-kronecker1} and \eqref{equation:orthogonal-in-kronecker2}.
This reveals the (matrix) reduced SVD of $\bX_{(n)}$:
\begin{equation}\label{equation:hosvd-initia-for-tucker-1}
\boxed{	
	\underbrace{\bX_{(n)}}_{I_n\times (I_{-n})} \approx
	\underbrace{\bA^{(n)}}_{{I_n\times R_n}} 
	\underbrace{\bSigma^{(n)}}_{R_n\times R_n}
	\underbrace{\bV^{(n)\top}}_{R_n \times I_{-n}}.
}
\end{equation}
And if \{$\eG, \bA^{(N)}, \ldots, \bA^{(n+1)}, \bA^{(n-1)}, \ldots, \bA^{(1)}$\} are fixed, the update for $\bA^{(n)}$ can be obtained by setting its columns  as the first $R_n$ left singular vectors of $\bX_{(n)}$. This matches the update in the subproblem outlined in Equation~\eqref{equation:hosvd-initia-for-tucker-0}. The \textbf{uniqueness} of the mode-$n$ singular values thus stems from the uniqueness of the (matrix) reduced SVD.

\paragraph{Vectorization.} Expanding upon the matricization process for the Nth-order tensor $\eX$, the vectorization is given by 
\begin{equation}\label{equation:hosvd-vec-in-theorem}
	\boxed{\underbrace{vec(\eX) }_{(I_1\ldots I_N)\times 1}
		\approx 
		\underbrace{\left(\bA^{(N)} \kronecker \bA^{(N-1)} \kronecker \ldots\kronecker \bA^{(1)}  \right)}_{(I_1\ldots I_N)\times (R_1\ldots R_N)}
		\underbrace{vec(\eG)}_{(R_1\ldots R_N)\times 1}.
	}
\end{equation}

\paragraph{Third-order case.}
For the sake of simplicity and enhancing comprehension, we  consider the reduced HOSVD applied to a third-order tensor $\eX\in \real^{I\times J\times K}$:
\begin{equation}\label{equation:hosvd-third-fact-1}
	\eX \approx \llbracket \eG; \bA, \bB, \bC   \rrbracket  
	=
	\sum_{p=1}^{P}\sum_{q=1}^{Q}\sum_{r=1}^{R}
	g_{pqr} \cdot \ba_p\circ \bb_q \circ \bc_r,
\end{equation}
where $\bA\in \real^{I\times P}, \bB\in \real^{J\times Q}$, $\bC\in \real^{K\times R}$, and $\eG\in \real^{P\times Q\times R}$. Then,  it follows that
\begin{itemize}
\item \textit{All orthogonality.}  Any subtensors $\eG_{:,\alpha,:}$ and $\eG_{:,\beta,:}$ for $\alpha\neq\beta$  are orthogonal:
$$
\langle  \eG_{:,\alpha,:},  \eG_{:,\beta,:}\rangle =0, \gap \alpha\neq \beta \in \{1,2,\ldots, J\}.
$$
\item \textit{Ordering.} The Frobenius norms of slices in each mode decreases as the running index increases:
$$
\norm{\eG_{:,1,:}}_F \geq \norm{\eG_{:,2,:}}_F \geq \ldots \geq \norm{\eG_{:,J,:}}_F \geq 0.
$$
\end{itemize}
The illustration of this third-order HOSVD (reduced and full versions) is similar to that of the Tucker decomposition and is shown in Figure~\ref{fig:tucker-decom-third}.

\section{Computing the HOSVD}

The calculation of the HOSVD is essentially equivalent to performing a (matrix) reduced SVD on the matricized form of the HOSVD, as described in Equation~\eqref{equation:hosvd-initia-for-tucker-1}. The procedure is presented in Algorithm~\ref{alg:hosvd-decomposition-full-gene}, and we shall not repeat the details. 
In certain constructions, we compute the tensor $\eG$  in each iteration, and terminate  if the norm between the updated $\eG$ and the old $\eX$  falls below a predefined constant threshold $\epsilon$ (see Proposition~\ref{proposition:frobenius-hosvd}).
See also \citet{de2000multilinear} for a more detailed discussion.
\begin{algorithm}[h] 
\caption{HOSVD via matrix SVD}
\label{alg:hosvd-decomposition-full-gene}
\begin{algorithmic}[1] 
\Require Nth-order tensor $\eX \in \real^{I_1\times I_2\times \ldots \times I_N}$;
\State Pick a rank $R_1, R_2, \ldots, R_N$ (or determine these values by the effective rank via matrix SVD);
\State Initialize $\bA^{(n)}\in \real^{I_n\times R_n}$ for all $n\in \{1,2,\ldots, N\}$ randomly; 
\State Choose maximal number of iterations $C$;
\State $iter=0$; \Comment{Count for the number of iterations}
\While{$iter<C$}
\State $iter=iter+1$; 
\For{$n=1,2,\ldots, N$} 
\State Find the matricization along mode-$n$: $\bX_{(n)}$;
\State Set the columns of $\bA^{(n)}$ to first $R_n$ leading left singular vectors of $\bX_{(n)}$;
\EndFor
\EndWhile
\State $\eG = \llbracket\eX; \bA^{(1)\top}, \bA^{(2)\top}, \ldots, \bA^{(N)\top}\rrbracket$;
\State Output $\bA^{(1)}, \bA^{(2)}, \ldots, \bA^{(N)}, \eG$; 
\end{algorithmic} 
\end{algorithm}

\section{Properties of the HOSVD}

\subsection{Frobenius Norm}
\begin{proposition}[Frobenius Norm of a Tensor]\label{proposition:frobenius-hosvd}
Let the HOSVD of an Nth-order tensor $\eX \in \real^{I_1\times I_2\times \ldots \times I_N}$ be given by $ \eX \approx \llbracket\eG; \bA^{(1)}, \bA^{(2)}, \ldots, \bA^{(N)}   \rrbracket  $. Then, it follows that 
$$
\begin{aligned}
\normf{\eX}^2 &= \normf{\eG}^2 
=\sum_{i=1}^{R_1}  (\sigma_i^{(1)})^2 
=\sum_{i=1}^{R_2}  (\sigma_i^{(2)})^2 
=\ldots
=\sum_{i=1}^{R_N}  (\sigma_i^{(N)})^2 .
\end{aligned}
$$
\end{proposition}
The squared Frobenius norm of a matrix is defined as  the sum of squares of its singular values  (see Definition~\ref{definition:frobernius-in-svd}). 
The  lemma  above indicates that the singular values within each slice of a tensor exhibit a similar characteristic. The proof of this lemma follows immediately from the (matrix) reduced SVD of the matricized form of the HOSVD, as presented Equation~\eqref{equation:hosvd-initia-for-tucker-1}.

\index{Truncated}
\index{Truncated SVD}
\subsection{Low-Rank Approximation}
Let $\bA\in \real^{M\times N}$. We have previously discussed that the best rank-$k$ approximation $\bA_k$ of the matrix can be obtained through the truncated SVD (Theorem~\ref{theorem:young-theorem-spectral}), where 
$$
\normf{\bA - \bA_k}^2 = \sigma_{k+1}^2+\sigma_{k+2}^2 +\ldots+\sigma_{\min\{M,N\}}^2,
$$
A similar result holds for low-rank tensor approximation.
\begin{proposition}[Low-Rank Approximation]\label{proposition:low-rank-hosvd}
Let the HOSVD of an Nth-order tensor $\eX \in \real^{I_1\times I_2\times \ldots \times I_N}$ be given by $ \eX \approx \llbracket\eG; \bA^{(1)}, \bA^{(2)}, \ldots, \bA^{(N)}   \rrbracket  $. Define a tensor $\weX$ by discarding the smallest mode-$n$ singular values $\sigma_{K_n+1}^{(n)}, \sigma_{K_n+2}^{(n)}, \ldots, \sigma_{R_n}^{(n)}$ for given values of $K_n\leq R_n$ (for all $n\in \{1,2,\ldots, N\}$). Then  it follows that 
$$
\normf{\eX-\weX}^2 \leq 
\sum_{r_1=K_1+1}^{R_1} (\sigma_{r_1}^{(1)})^2+
\sum_{r_2=K_2+1}^{R_2} (\sigma_{r_2}^{(2)})^2+
\ldots 
+
\sum_{r_N=K_N+1}^{R_N} (\sigma_{r_N}^{(N)})^2.
$$
\end{proposition}
\begin{proof}[of Proposition~\ref{proposition:low-rank-hosvd}]
According to Proposition~\ref{proposition:frobenius-hosvd}, we obtain the following
$$
\footnotesize
\begin{aligned}
&\normf{\eX-\weX}^2= 
\sum_{r_1=1}^{R_1} \sum_{r_2=1}^{R_2} \ldots \sum_{r_N=1}^{R_N}
g_{r_1r_2\ldots r_N}^2
-
\sum_{r_1=1}^{K_1} \sum_{r_2=1}^{K_2} \ldots \sum_{r_N=1}^{K_N}
g_{r_1r_2\ldots r_N}^2 \\
&=\sum_{r_1=K_1+1}^{R_1} \sum_{r_2=K_2+1}^{R_2} \ldots \sum_{r_N=K_N+1}^{R_N}
g_{r_1r_2\ldots r_N}^2 \\
&\leq 
\sum_{\textcolor{mylightbluetext}{r_1=K_1+1}}^{R_1} \sum_{r_2=1}^{R_2} \ldots \sum_{r_N=1}^{R_N}
g_{r_1r_2\ldots r_N}^2 
+
\sum_{r_1=1}^{R_1} \sum_{\textcolor{mylightbluetext}{r_2=K_2+1}}^{R_2} \ldots \sum_{r_N=1}^{R_N}
g_{r_1r_2\ldots r_N}^2 
+
\ldots 
+
\sum_{r_1=1}^{R_1} \sum_{r_2=1}^{R_2} \ldots \sum_{\textcolor{mylightbluetext}{r_N=K_N+1}}^{R_N}
g_{r_1r_2\ldots r_N}^2 \\
&=\sum_{r_1=K_1+1}^{R_1} (\sigma_{r_1}^{(1)})^2+
\sum_{r_2=K_2+1}^{R_2} (\sigma_{r_2}^{(2)})^2+
\ldots 
+
\sum_{r_N=K_N+1}^{R_N} (\sigma_{r_N}^{(N)})^2.
\end{aligned}
$$
This completes the proof.
\end{proof}

%% file: chapter-tensorttdecom.tex

\index{Decomposition: TT}
\section{Tensor-Train (TT) Decomposition}

\begin{theoremHigh}[Tensor-Train (TT) Decomposition \citep{oseledets2011tensor}]\label{theorem:ttrain-decomp}
For a general Nth-order tensor, $\eX \in \real^{I_1\times I_2\times \ldots \times I_N}$, it admits the tensor-train decomposition
$$
\eX \approx \eG^{(1)} \boxtimes \eG^{(2)} \boxtimes \ldots \boxtimes \eG^{(N)},
$$
where the symbol ``$\boxtimes$" means that the elements of $\eX$ can be obtained by a product of $N$ matrices
$$
x_{i_1,i_2,\ldots,i_N} \triangleq\eX_{i_1,i_2,\ldots,i_N}= \eG^{(1)}_{:,:,i_1} \eG^{(2)}_{:,:,i_2}\ldots \eG^{(N)}_{:,:,i_N}.
$$ 
Note here
\begin{itemize}
\item Each $\eG^{(n)} \in \real^{R_{n-1}\times R_{n}\times I_n}$ for all $n\in \{1,2,\ldots, N\}$  is a third-order tensor, and is referred to as a \textit{TT core};
\item Each $\eG^{(n)}_{:,:,i_n} \in \real^{R_{n-1}\times R_n}$ for all $i_n\in \{1,2,\ldots, I_n\}$ is a frontal slice of $\eG^{(n)}$ (see Figure~\ref{fig:thifd-slices});
\item $R_0, R_N$ are imposed to be 1 for boundary conditions; 
\item $R_1,R_2,\ldots, R_{N-1}$ are known as the \textit{tensor ranks (TT-ranks)} of the corresponding dimensions;
\end{itemize}
\noindent
The illustration of how to extract each element of the decomposition is shown in Figure~\ref{fig:tensortrain-decom}.
\end{theoremHigh}

Tensor Train (TT) decomposition offers distinct advantages over CP or Tucker decomposition, primarily due to its stability, flexibility, and control over approximation accuracy. While CP decomposition represents a tensor as a sum of rank-one components, TT decomposition arranges tensors in a sequential chain of smaller core tensors, connected in a structured format. This difference provides several advantages:
\begin{itemize}
\item TT decomposition is built upon well-established techniques, such as truncated singular value decomposition (tSVD) and adaptive cross-approximation, which enable robust approximation and reduce instability issues (see below).
\item TT decomposition provides explicit control over approximation accuracy by setting a threshold for the singular values retained in each core tensor. This allows for quasi-optimal rank reduction, meaning the decomposition can maintain a specific level of accuracy with minimal computational overhead. In contrast, CP or Tucker decomposition lacks this direct control, often requiring heuristic approaches to manage rank selection and approximation error, which can lead to inefficiencies or less reliable results.
\item TT decomposition allows for structured low-rank representation by sequentially compressing each mode using low-rank matrices. This rank reduction approach provides a more refined and  hierarchical compression of data, preserving essential patterns with less information loss. TT can achieve compact representations with fewer components than CP or Tucker, especially for high-dimensional data, which reduces computational and memory costs while maintaining accuracy.
\item By organizing core tensors in a sequence, TT can capture complex interactions among latent variables across modes more effectively than CP, while avoiding the exponentially growing memory requirements associated with the core tensor in the Tucker decomposition (Tucker decomposition becomes  unfavourable for large order $N$). The tensor train structure inherently separates latent variables, allowing for a more interpretable and sophisticated representation of the data’s underlying structure. This separation is particularly beneficial in applications that involve hierarchical or sequential dependencies, such as time series or natural language processing tasks.
\item TT decomposition algorithms are simpler and often faster to compute than CP algorithms, particularly in high dimensions. Since each core tensor is computed independently, TT algorithms can leverage matrix-based techniques like SVD for efficient processing. This simplicity enables TT to be applied more easily to large-scale data, with predictable computational requirements and low-rank approximation control that CP decomposition lacks.
\end{itemize}

\begin{figure}[h]
\centering
\resizebox{1.0\textwidth}{!}{%
\begin{tikzpicture}
\pgfmathsetmacro{\bxx}{0.3}
\pgfmathsetmacro{\begupwhite}{0.2}

\pgfmathsetmacro{\bcwd}{0.3}
\pgfmathsetmacro{\bbwd}{0.3}
\pgfmathsetmacro{\bbright}{0.9}
\pgfmathsetmacro{\bbup}{0.225}

\pgfmathsetmacro{\distrightdown}{4.2}
\pgfmathsetmacro{\distl}{1.2}
\pgfmathsetmacro{\aprimeright}{0.4}

\pgfmathsetmacro{\abvd}{3.3}
\pgfmathsetmacro{\abvdd}{3.9}
\draw (6.6,-3.2+\abvd) node {{\color{black}\large{$\eX_{i_1, i_2,\ldots, i_N}$}}};
\draw (8,-3.2+\abvd) node {	{\color{black}\large{$\approx$}}};

\draw [very thick] (5.8+0.5,-3.8+0.5) rectangle (7+0.5,-2.6+0.5);
\filldraw [fill=gray!60!white,draw=green!40!black] (5.8+0.5,-3.8+0.5) rectangle (7+0.5,-2.6+0.5);
\draw [step=0.4/2, very thin, color=gray] (5.8+0.5,-3.8+0.5) grid (7+0.5,-2.6+0.5);

\draw [very thick] (5.8+0.4,-3.8+0.4) rectangle (7+0.4,-2.6+0.4);
\filldraw [fill=gray!50!white,draw=green!40!black] (5.8+0.4,-3.8+0.4) rectangle (7+0.4,-2.6+0.4);
\draw [step=0.4/2, very thin, color=gray] (5.8+0.4,-3.8+0.4) grid (7+0.4,-2.6+0.4);

\draw [very thick] (5.8+0.3,-3.8+0.3) rectangle (7+0.3,-2.6+0.3);
\filldraw [fill=gray!40!white,draw=green!40!black] (5.8+0.3,-3.8+0.3) rectangle (7+0.3,-2.6+0.3);
\draw [step=0.4/2, very thin, color=gray] (5.8+0.3,-3.8+0.3) grid (7+0.3,-2.6+0.3);

\draw [very thick] (5.8+0.2,-3.8+0.2) rectangle (7+0.2,-2.6+0.2);
\filldraw [fill=gray!30!white,draw=green!40!black] (5.8+0.2,-3.8+0.2) rectangle (7+0.2,-2.6+0.2);
\draw [step=0.4/2, very thin, color=gray] (5.8+0.2,-3.8+0.2) grid (7+0.2,-2.6+0.2);

\draw [very thick] (5.8+0.1,-3.8+0.1) rectangle (7+0.1,-2.6+0.1);
\filldraw [fill=gray!20!white,draw=green!40!black] (5.8+0.1,-3.8+0.1) rectangle (7+0.1,-2.6+0.1);
\draw [step=0.4/2, very thin, color=gray] (5.8+0.1,-3.8+0.1) grid (7+0.1,-2.6+0.1);

\draw [very thick] (5.8,-3.8) rectangle (7,-2.6);
\filldraw [fill=gray!10!white,draw=green!40!black] (5.8,-3.8) rectangle (7,-2.6);
\draw [step=0.4/2, very thin, color=gray] (5.8,-3.8) grid (7,-2.6);

\draw [very thick] (5.8-0.1,-3.8-0.1) rectangle (7-0.1,-2.6-0.1);
\filldraw [fill=gray!6!white,draw=green!40!black] (5.8-0.1,-3.8-0.1) rectangle (7-0.1,-2.6-0.1);
\draw [step=0.4/2, very thin, color=gray] (5.8-0.1,-3.8-0.1) grid (7-0.1,-2.6-0.1);



\draw (6.2,-4.6) node {{\color{black}\scriptsize{$\eX\in\real^{I_1\times I_2 \times \ldots \times I_N}$}}};

\draw (8,-3.2) node {	{\color{black}\large{$\approx$}}};

\pgfmathsetmacro{\cubex}{1}
\pgfmathsetmacro{\cubey}{1}
\pgfmathsetmacro{\cubez}{0.1}
\pgfmathsetmacro{\slicedist}{0.3}
\pgfmathsetmacro{\xa}{3.7}
\pgfmathsetmacro{\ya}{3.1}
\draw[black,fill=gray!70!white] (6+\xa,0-\ya,-6*\slicedist) -- ++(-\cubex,0,0) -- ++(0,-\cubey,0) -- ++(\cubex,0,0) -- cycle;
\draw[black,fill=gray!70!white] (6+\xa,0-\ya,-6*\slicedist) -- ++(0,0,-\cubez) -- ++(0,-\cubey,0) -- ++(0,0,\cubez) -- cycle;
\draw[black,fill=gray!70!white] (6+\xa,0-\ya,-6*\slicedist) -- ++(-\cubex,0,0) -- ++(0,0,-\cubez) -- ++(\cubex,0,0) -- cycle;

\draw[black,fill=gray!60!white] (6+\xa,0-\ya,-5*\slicedist) -- ++(-\cubex,0,0) -- ++(0,-\cubey,0) -- ++(\cubex,0,0) -- cycle;
\draw[black,fill=gray!60!white] (6+\xa,0-\ya,-5*\slicedist) -- ++(0,0,-\cubez) -- ++(0,-\cubey,0) -- ++(0,0,\cubez) -- cycle;
\draw[black,fill=gray!60!white] (6+\xa,0-\ya,-5*\slicedist) -- ++(-\cubex,0,0) -- ++(0,0,-\cubez) -- ++(\cubex,0,0) -- cycle;

\draw[black,fill=gray!50!white] (6+\xa,0-\ya,-4*\slicedist) -- ++(-\cubex,0,0) -- ++(0,-\cubey,0) -- ++(\cubex,0,0) -- cycle;
\draw[black,fill=gray!50!white] (6+\xa,0-\ya,-4*\slicedist) -- ++(0,0,-\cubez) -- ++(0,-\cubey,0) -- ++(0,0,\cubez) -- cycle;
\draw[black,fill=gray!50!white] (6+\xa,0-\ya,-4*\slicedist) -- ++(-\cubex,0,0) -- ++(0,0,-\cubez) -- ++(\cubex,0,0) -- cycle;
\draw[black,fill=gray!40!white] (6+\xa,0-\ya,-3*\slicedist) -- ++(-\cubex,0,0) -- ++(0,-\cubey,0) -- ++(\cubex,0,0) -- cycle;
\draw[black,fill=gray!40!white] (6+\xa,0-\ya,-3*\slicedist) -- ++(0,0,-\cubez) -- ++(0,-\cubey,0) -- ++(0,0,\cubez) -- cycle;
\draw[black,fill=gray!40!white] (6+\xa,0-\ya,-3*\slicedist) -- ++(-\cubex,0,0) -- ++(0,0,-\cubez) -- ++(\cubex,0,0) -- cycle;
\draw[black,fill=gray!30!white] (6+\xa,0-\ya,-2*\slicedist) -- ++(-\cubex,0,0) -- ++(0,-\cubey,0) -- ++(\cubex,0,0) -- cycle;
\draw[black,fill=gray!30!white] (6+\xa,0-\ya,-2*\slicedist) -- ++(0,0,-\cubez) -- ++(0,-\cubey,0) -- ++(0,0,\cubez) -- cycle;
\draw[black,fill=gray!30!white] (6+\xa,0-\ya,-2*\slicedist) -- ++(-\cubex,0,0) -- ++(0,0,-\cubez) -- ++(\cubex,0,0) -- cycle;
\draw[black,fill=gray!20!white] (6+\xa,0-\ya,-1*\slicedist) -- ++(-\cubex,0,0) -- ++(0,-\cubey,0) -- ++(\cubex,0,0) -- cycle;
\draw[black,fill=gray!20!white] (6+\xa,0-\ya,-1*\slicedist) -- ++(0,0,-\cubez) -- ++(0,-\cubey,0) -- ++(0,0,\cubez) -- cycle;
\draw[black,fill=gray!20!white] (6+\xa,0-\ya,-1*\slicedist) -- ++(-\cubex,0,0) -- ++(0,0,-\cubez) -- ++(\cubex,0,0) -- cycle;
\draw[black,fill=gray!10!white] (6+\xa,0-\ya,0) -- ++(-\cubex,0,0) -- ++(0,-\cubey,0) -- ++(\cubex,0,0) -- cycle;
\draw[black,fill=gray!10!white] (6+\xa,0-\ya,0) -- ++(0,0,-\cubez) -- ++(0,-\cubey,0) -- ++(0,0,\cubez) -- cycle;
\draw[black,fill=gray!10!white] (6+\xa,0-\ya,0) -- ++(-\cubex,0,0) -- ++(0,0,-\cubez) -- ++(\cubex,0,0) -- cycle;
\draw (7+\xa,0-\ya) node {	{\color{black}\large{$\boxtimes$}}};
\draw (5.7+\xa,-4.8) node {{\color{black}\scriptsize{$\eG^{(1)}\in\real^{R_0\times R_1 \times  I_1}$}}};

\draw[->] (4.9+\xa,-3) -- ++(0.8,0.8); 
\draw[->] (4.85+\xa,-3.1) -- ++(0,-1.);  
\draw[->] (4.95+\xa,-4.2) -- ++(1.25,0);  
\draw (5.5+\xa,-4.45) node {{\color{black}\scriptsize{$R_1$}}};
\draw (4.6+\xa,-3.5) node {{\color{black}\scriptsize{$R_0$}}};
\draw (5.2+\xa,-2.35) node {{\color{black}\scriptsize{$I_1$}}};

\draw[black,fill=gray!40!white] (6+\xa,0-\ya+\abvd,-3*\slicedist) -- ++(-\cubex,0,0) -- ++(0,-\cubey,0) -- ++(\cubex,0,0) -- cycle;
\draw[black,fill=gray!40!white] (6+\xa,0-\ya+\abvd,-3*\slicedist) -- ++(0,0,-\cubez) -- ++(0,-\cubey,0) -- ++(0,0,\cubez) -- cycle;
\draw[black,fill=gray!40!white] (6+\xa,0-\ya+\abvd,-3*\slicedist) -- ++(-\cubex,0,0) -- ++(0,0,-\cubez) -- ++(\cubex,0,0) -- cycle;
\draw (5.8+\xa,-4.7+\abvdd) node {{\color{black}\scriptsize{$\eG^{(1)}_{:,:,i_1}\in\real^{R_0\times R_1}$}}};
\draw (7+\xa,0-\ya+\abvd) node {	{\color{black}\large{$\times$}}};
\draw[-{Computer Modern Rightarrow}] (5.8+\xa,-2.) -- ++(0,0.8); 
\draw (5.8+\xa,-1.65) node {{\color{black}\scriptsize{\textcolor{mylightbluetext}{$i_1$-th slice}}}};
\pgfmathsetmacro{\xa}{6.2}
\draw[black,fill=gray!70!white] (6+\xa,0-\ya,-8*\slicedist) -- ++(-\cubex,0,0) -- ++(0,-\cubey,0) -- ++(\cubex,0,0) -- cycle;
\draw[black,fill=gray!70!white] (6+\xa,0-\ya,-8*\slicedist) -- ++(0,0,-\cubez) -- ++(0,-\cubey,0) -- ++(0,0,\cubez) -- cycle;
\draw[black,fill=gray!70!white] (6+\xa,0-\ya,-8*\slicedist) -- ++(-\cubex,0,0) -- ++(0,0,-\cubez) -- ++(\cubex,0,0) -- cycle;
\draw[black,fill=gray!70!white] (6+\xa,0-\ya,-7*\slicedist) -- ++(-\cubex,0,0) -- ++(0,-\cubey,0) -- ++(\cubex,0,0) -- cycle;
\draw[black,fill=gray!70!white] (6+\xa,0-\ya,-7*\slicedist) -- ++(0,0,-\cubez) -- ++(0,-\cubey,0) -- ++(0,0,\cubez) -- cycle;
\draw[black,fill=gray!70!white] (6+\xa,0-\ya,-7*\slicedist) -- ++(-\cubex,0,0) -- ++(0,0,-\cubez) -- ++(\cubex,0,0) -- cycle;
\draw[black,fill=gray!70!white] (6+\xa,0-\ya,-6*\slicedist) -- ++(-\cubex,0,0) -- ++(0,-\cubey,0) -- ++(\cubex,0,0) -- cycle;
\draw[black,fill=gray!70!white] (6+\xa,0-\ya,-6*\slicedist) -- ++(0,0,-\cubez) -- ++(0,-\cubey,0) -- ++(0,0,\cubez) -- cycle;
\draw[black,fill=gray!70!white] (6+\xa,0-\ya,-6*\slicedist) -- ++(-\cubex,0,0) -- ++(0,0,-\cubez) -- ++(\cubex,0,0) -- cycle;

\draw[black,fill=gray!60!white] (6+\xa,0-\ya,-5*\slicedist) -- ++(-\cubex,0,0) -- ++(0,-\cubey,0) -- ++(\cubex,0,0) -- cycle;
\draw[black,fill=gray!60!white] (6+\xa,0-\ya,-5*\slicedist) -- ++(0,0,-\cubez) -- ++(0,-\cubey,0) -- ++(0,0,\cubez) -- cycle;
\draw[black,fill=gray!60!white] (6+\xa,0-\ya,-5*\slicedist) -- ++(-\cubex,0,0) -- ++(0,0,-\cubez) -- ++(\cubex,0,0) -- cycle;

\draw[black,fill=gray!50!white] (6+\xa,0-\ya,-4*\slicedist) -- ++(-\cubex,0,0) -- ++(0,-\cubey,0) -- ++(\cubex,0,0) -- cycle;
\draw[black,fill=gray!50!white] (6+\xa,0-\ya,-4*\slicedist) -- ++(0,0,-\cubez) -- ++(0,-\cubey,0) -- ++(0,0,\cubez) -- cycle;
\draw[black,fill=gray!50!white] (6+\xa,0-\ya,-4*\slicedist) -- ++(-\cubex,0,0) -- ++(0,0,-\cubez) -- ++(\cubex,0,0) -- cycle;
\draw[black,fill=gray!40!white] (6+\xa,0-\ya,-3*\slicedist) -- ++(-\cubex,0,0) -- ++(0,-\cubey,0) -- ++(\cubex,0,0) -- cycle;
\draw[black,fill=gray!40!white] (6+\xa,0-\ya,-3*\slicedist) -- ++(0,0,-\cubez) -- ++(0,-\cubey,0) -- ++(0,0,\cubez) -- cycle;
\draw[black,fill=gray!40!white] (6+\xa,0-\ya,-3*\slicedist) -- ++(-\cubex,0,0) -- ++(0,0,-\cubez) -- ++(\cubex,0,0) -- cycle;
\draw[black,fill=gray!30!white] (6+\xa,0-\ya,-2*\slicedist) -- ++(-\cubex,0,0) -- ++(0,-\cubey,0) -- ++(\cubex,0,0) -- cycle;
\draw[black,fill=gray!30!white] (6+\xa,0-\ya,-2*\slicedist) -- ++(0,0,-\cubez) -- ++(0,-\cubey,0) -- ++(0,0,\cubez) -- cycle;
\draw[black,fill=gray!30!white] (6+\xa,0-\ya,-2*\slicedist) -- ++(-\cubex,0,0) -- ++(0,0,-\cubez) -- ++(\cubex,0,0) -- cycle;
\draw[black,fill=gray!20!white] (6+\xa,0-\ya,-1*\slicedist) -- ++(-\cubex,0,0) -- ++(0,-\cubey,0) -- ++(\cubex,0,0) -- cycle;
\draw[black,fill=gray!20!white] (6+\xa,0-\ya,-1*\slicedist) -- ++(0,0,-\cubez) -- ++(0,-\cubey,0) -- ++(0,0,\cubez) -- cycle;
\draw[black,fill=gray!20!white] (6+\xa,0-\ya,-1*\slicedist) -- ++(-\cubex,0,0) -- ++(0,0,-\cubez) -- ++(\cubex,0,0) -- cycle;
\draw[black,fill=gray!10!white] (6+\xa,0-\ya,0) -- ++(-\cubex,0,0) -- ++(0,-\cubey,0) -- ++(\cubex,0,0) -- cycle;
\draw[black,fill=gray!10!white] (6+\xa,0-\ya,0) -- ++(0,0,-\cubez) -- ++(0,-\cubey,0) -- ++(0,0,\cubez) -- cycle;
\draw[black,fill=gray!10!white] (6+\xa,0-\ya,0) -- ++(-\cubex,0,0) -- ++(0,0,-\cubez) -- ++(\cubex,0,0) -- cycle;
\draw (7.2+\xa,0-\ya) node {	{\color{black}\large{$\boxtimes$}}};
\draw (5.9+\xa,-4.8) node {{\color{black}\scriptsize{$\eG^{(2)}\in\real^{R_1\times R_2 \times  I_2}$}}};
\draw[->] (4.9+\xa,-3) -- ++(0.85,0.85); 
\draw[->] (4.85+\xa,-3.1) -- ++(0,-1.);  
\draw[->] (4.95+\xa,-4.2) -- ++(1.25,0);  
\draw (5.5+\xa,-4.45) node {{\color{black}\scriptsize{$R_2$}}};
\draw (4.6+\xa,-3.6) node {{\color{black}\scriptsize{$R_1$}}};
\draw (5.2+\xa,-2.35) node {{\color{black}\scriptsize{$I_2$}}};

\draw[black,fill=gray!40!white] (6+\xa,0-\ya+\abvd,-3*\slicedist) -- ++(-\cubex,0,0) -- ++(0,-\cubey,0) -- ++(\cubex,0,0) -- cycle;
\draw[black,fill=gray!40!white] (6+\xa,0-\ya+\abvd,-3*\slicedist) -- ++(0,0,-\cubez) -- ++(0,-\cubey,0) -- ++(0,0,\cubez) -- cycle;
\draw[black,fill=gray!40!white] (6+\xa,0-\ya+\abvd,-3*\slicedist) -- ++(-\cubex,0,0) -- ++(0,0,-\cubez) -- ++(\cubex,0,0) -- cycle;
\draw (5.9+\xa,-4.7+\abvdd) node {{\color{black}\scriptsize{$\eG^{(2)}_{:,:,i_2}\in\real^{R_1\times R_2}$}}};
\draw (7+\xa,0-\ya+\abvd) node {	{\color{black}\large{$\times$}}};
\draw[-{Computer Modern Rightarrow}] (5.8+\xa,-2.) -- ++(0,0.8); 
\draw (5.8+\xa,-1.65) node {{\color{black}\scriptsize{\textcolor{mylightbluetext}{$i_2$-th slice}}}};

\draw (7.9+\xa,0-\ya) node {	{\color{black}\large{$\ldots \boxtimes$}}};
\draw (7.9+\xa,0-\ya+\abvd) node {	{\color{black}\large{$\ldots \times$}}};
\pgfmathsetmacro{\xa}{9.9}
\draw[black,fill=gray!70!white] (6+\xa,0-\ya,-7*\slicedist) -- ++(-\cubex,0,0) -- ++(0,-\cubey,0) -- ++(\cubex,0,0) -- cycle;
\draw[black,fill=gray!70!white] (6+\xa,0-\ya,-7*\slicedist) -- ++(0,0,-\cubez) -- ++(0,-\cubey,0) -- ++(0,0,\cubez) -- cycle;
\draw[black,fill=gray!70!white] (6+\xa,0-\ya,-7*\slicedist) -- ++(-\cubex,0,0) -- ++(0,0,-\cubez) -- ++(\cubex,0,0) -- cycle;
\draw[black,fill=gray!70!white] (6+\xa,0-\ya,-6*\slicedist) -- ++(-\cubex,0,0) -- ++(0,-\cubey,0) -- ++(\cubex,0,0) -- cycle;
\draw[black,fill=gray!70!white] (6+\xa,0-\ya,-6*\slicedist) -- ++(0,0,-\cubez) -- ++(0,-\cubey,0) -- ++(0,0,\cubez) -- cycle;
\draw[black,fill=gray!70!white] (6+\xa,0-\ya,-6*\slicedist) -- ++(-\cubex,0,0) -- ++(0,0,-\cubez) -- ++(\cubex,0,0) -- cycle;

\draw[black,fill=gray!60!white] (6+\xa,0-\ya,-5*\slicedist) -- ++(-\cubex,0,0) -- ++(0,-\cubey,0) -- ++(\cubex,0,0) -- cycle;
\draw[black,fill=gray!60!white] (6+\xa,0-\ya,-5*\slicedist) -- ++(0,0,-\cubez) -- ++(0,-\cubey,0) -- ++(0,0,\cubez) -- cycle;
\draw[black,fill=gray!60!white] (6+\xa,0-\ya,-5*\slicedist) -- ++(-\cubex,0,0) -- ++(0,0,-\cubez) -- ++(\cubex,0,0) -- cycle;

\draw[black,fill=gray!50!white] (6+\xa,0-\ya,-4*\slicedist) -- ++(-\cubex,0,0) -- ++(0,-\cubey,0) -- ++(\cubex,0,0) -- cycle;
\draw[black,fill=gray!50!white] (6+\xa,0-\ya,-4*\slicedist) -- ++(0,0,-\cubez) -- ++(0,-\cubey,0) -- ++(0,0,\cubez) -- cycle;
\draw[black,fill=gray!50!white] (6+\xa,0-\ya,-4*\slicedist) -- ++(-\cubex,0,0) -- ++(0,0,-\cubez) -- ++(\cubex,0,0) -- cycle;
\draw[black,fill=gray!40!white] (6+\xa,0-\ya,-3*\slicedist) -- ++(-\cubex,0,0) -- ++(0,-\cubey,0) -- ++(\cubex,0,0) -- cycle;
\draw[black,fill=gray!40!white] (6+\xa,0-\ya,-3*\slicedist) -- ++(0,0,-\cubez) -- ++(0,-\cubey,0) -- ++(0,0,\cubez) -- cycle;
\draw[black,fill=gray!40!white] (6+\xa,0-\ya,-3*\slicedist) -- ++(-\cubex,0,0) -- ++(0,0,-\cubez) -- ++(\cubex,0,0) -- cycle;
\draw[black,fill=gray!30!white] (6+\xa,0-\ya,-2*\slicedist) -- ++(-\cubex,0,0) -- ++(0,-\cubey,0) -- ++(\cubex,0,0) -- cycle;
\draw[black,fill=gray!30!white] (6+\xa,0-\ya,-2*\slicedist) -- ++(0,0,-\cubez) -- ++(0,-\cubey,0) -- ++(0,0,\cubez) -- cycle;
\draw[black,fill=gray!30!white] (6+\xa,0-\ya,-2*\slicedist) -- ++(-\cubex,0,0) -- ++(0,0,-\cubez) -- ++(\cubex,0,0) -- cycle;
\draw[black,fill=gray!20!white] (6+\xa,0-\ya,-1*\slicedist) -- ++(-\cubex,0,0) -- ++(0,-\cubey,0) -- ++(\cubex,0,0) -- cycle;
\draw[black,fill=gray!20!white] (6+\xa,0-\ya,-1*\slicedist) -- ++(0,0,-\cubez) -- ++(0,-\cubey,0) -- ++(0,0,\cubez) -- cycle;
\draw[black,fill=gray!20!white] (6+\xa,0-\ya,-1*\slicedist) -- ++(-\cubex,0,0) -- ++(0,0,-\cubez) -- ++(\cubex,0,0) -- cycle;
\draw[black,fill=gray!10!white] (6+\xa,0-\ya,0) -- ++(-\cubex,0,0) -- ++(0,-\cubey,0) -- ++(\cubex,0,0) -- cycle;
\draw[black,fill=gray!10!white] (6+\xa,0-\ya,0) -- ++(0,0,-\cubez) -- ++(0,-\cubey,0) -- ++(0,0,\cubez) -- cycle;
\draw[black,fill=gray!10!white] (6+\xa,0-\ya,0) -- ++(-\cubex,0,0) -- ++(0,0,-\cubez) -- ++(\cubex,0,0) -- cycle;
\draw (5.9+\xa,-4.8) node {{\color{black}\scriptsize{$\eG^{(N)}\in\real^{R_{N-1}\times R_N \times  I_N}$}}};
\draw[->] (4.9+\xa,-3) -- ++(0.8,0.8); 
\draw[->] (4.85+\xa,-3.1) -- ++(0,-1.);  
\draw[->] (4.95+\xa,-4.2) -- ++(1.25,0);  
\draw (5.5+\xa,-4.45) node {{\color{black}\scriptsize{$R_N$}}};
\draw (4.48+\xa,-3.7) node {{\color{black}\scriptsize{$R_{N-1}$}}};
\draw (5.2+\xa,-2.35) node {{\color{black}\scriptsize{$I_N$}}};

\draw[black,fill=gray!40!white] (6+\xa,0-\ya+\abvd,-3*\slicedist) -- ++(-\cubex,0,0) -- ++(0,-\cubey,0) -- ++(\cubex,0,0) -- cycle;
\draw[black,fill=gray!40!white] (6+\xa,0-\ya+\abvd,-3*\slicedist) -- ++(0,0,-\cubez) -- ++(0,-\cubey,0) -- ++(0,0,\cubez) -- cycle;
\draw[black,fill=gray!40!white] (6+\xa,0-\ya+\abvd,-3*\slicedist) -- ++(-\cubex,0,0) -- ++(0,0,-\cubez) -- ++(\cubex,0,0) -- cycle;
\draw (5.9+\xa,-4.7+\abvdd) node {{\color{black}\scriptsize{$\eG^{(N)}_{:,:,i_N}\in\real^{R_{N-1}\times R_N}$}}};
\draw[-{Computer Modern Rightarrow}] (5.8+\xa,-2.) -- ++(0,0.8); 
\draw (5.8+\xa,-1.65) node {{\color{black}\scriptsize{\textcolor{mylightbluetext}{$i_N$-th slice}}}};
\end{tikzpicture}
}
\caption{Tensor-train decomposition of an Nth-order tensor: $\eX \in \real^{I_1\times I_2\times \ldots \times I_N}\approx \eG^{(1)} \boxtimes \eG^{(2)} \boxtimes \ldots \boxtimes \eG^{(N)}$, where each $\eG^{(n)} \in \real^{R_{n-1}\times R_{n}\times I_n}$ for $n\in \{1,2,\ldots,N\}$. Each $(i_1,i_2,\ldots,i_N)$-th element is obtained by the matrix multiplication of the corresponding frontal slices of $\eG^{(n)}$'s.
Note here $R_0=R_N=1$.
}
\label{fig:tensortrain-decom}
\end{figure}

The illustration of the TT decomposition for an Nth order tensor is shown in Figure~\ref{fig:tensortrain-decom}.
In other words, the TT format approximates each entry of the tensor $\eX$ with a product of $N$ matrices, in particular with a sequence of $R_n \times R_{n+1}$ matrices, each indexed by the parameter $i_{n+1}$. The figure resembles  a
train with links connecting them, hence the name ``train." Suppose that the TT-ranks are all equal, $R_1=R_2=\ldots=R_{N-1} = R$, and that $I_1=I_2=\ldots=I_N = I$, then the TT decomposition
requires the storage of $\mathcalO(N I R^2)$ floats. Thus, the memory complexity of the TT decomposition scales linearly with dimension.

\paragraph{The ``$\boxtimes$" notation.} The symbol ``$\boxtimes $" is defined as a special tensor product. Let $\eA\in \real^{R_1\times R_2 \times I_1\times I_2\times \ldots\times I_M}$ and $\eB\in \real^{R_2\times R_3\times J_1\times J_2\times \ldots \times J_N}$ (i.e., the second dimension of $\eA$ and the first dimension of $\eB$ align), then it follows that 
$$
\eA\boxtimes \eB \in \real^{R_1\times R_3 \times  I_1\times I_2\times \ldots\times I_M \times J_1\times J_2\times \ldots \times J_N},
$$
where each element is given by
$$
\begin{aligned}
(\eA\boxtimes \eB)_{r_1, r_3 ,  i_1, i_2, \ldots, i_M , j_1, j_2, \ldots , j_N}
=\sum_{r_2=1}^{R_2}
(a_{r_1,\textcolor{mylightbluetext}{r_2},i_1, i_2, \ldots, i_M})
(b_{\textcolor{mylightbluetext}{r_2},r_3,j_1, j_2, \ldots , j_N}).
\end{aligned}
$$

\paragraph{TT-rank.}
Define the \textit{tensor unfolding} for an Nth-order tensor $\eX \in \real^{I_1\times I_2\times \ldots \times I_N}$ in the following way
\begin{equation}\label{equation:tt_unfond1}
\abovebX{n} \triangleq \mathcalM_{\{1,2,\ldots,n\}} (\eX) \in \real^{(I_1I_2\ldots I_n) \times (I_{n+1}I_{n+2}\ldots I_N)},
\end{equation}
i.e., a tuple matricization (Section~\ref{section:tuple_matricization}), where the $\big(\{i_1\ldots i_n\},\{i_{n+1}\ldots i_N\}\big)$-th element of $\abovebX{n}$ is given by $\abovebX{n}_{i_1\ldots i_n,i_{n+1}\ldots i_N} = \eX_{i_1,i_2,\ldots,i_N}$. 
This specific tensor unfolding can also be expressed using a reshape operation:
\begin{equation}\label{equation:tt_unfond2}
\abovebX{n} = \text{reshape}\bigg(\eX, (I_1I_2\ldots I_n), (I_{n+1}I_{n+2}\ldots I_N) \bigg).
\end{equation}
This unfolding implies $\rank(\abovebX{n})\leq R_n$ for all $n\in\{1,2,\ldots,N\}$. 

\begin{definition}[TT-Rank of a Tensor]\label{definition:ttde_rank}
The full tensor TT-ranks of  $\eX \in \real^{I_1\times I_2\times \ldots \times I_N}$ are defined as the tuple in $(R_1, R_2, \ldots, R_{N-1})$ derived from the TT decomposition (Theorem~\ref{theorem:ttrain-decomp}).
The set $TT(R_1, R_2, \ldots, R_{N-1})$ consists of tensors whose TT-ranks do not exceed  $(R_1, R_2, \ldots, R_{N-1})$.
\end{definition}
It has been shown that the set $TT(R_1, R_2, \ldots, R_{N-1})$ of tensors with TT-ranks bounded by $(R_1, R_2, \ldots, R_{N-1})$ is closed, and under the condition of full rank, this set constitutes a smooth manifold \citep{holtz2012manifolds, uschmajew2013geometry}.

\index{Truncated}
\index{Truncated SVD}
\section{Computing the TT Decomposition}
\begin{algorithm}[h] 
	\caption{TT-SVD}
	\label{alg:tt-decomposi}
	\begin{algorithmic}[1] 
		\Require Nth-order tensor $\eX \in \real^{I_1\times I_2\times \ldots \times I_N}$;
		\State Set the initial matrix $\bC \leftarrow \abovebX{1} \in \real^{I_1 \times I_2I_3\ldots I_N}$;
		\State Set $R_0=R_N=1$;
		\For{$n=1,2,\ldots, N-1$} 
		\State $\bC=\text{reshape}\bigg(\bC, (I_nR_{n-1}), (I_{n+1}I_{n+2}\ldots I_N)\bigg)$;
		\State Compute the $\delta$-truncated SVD of $\bC = \bU\bSigma\bV^\top$, numerical rank $R_n=\rank(\bC)$;
		\State Set the rows of $\bG^{(n)}_{(2)}\in \real^{R_n \times I_nR_{n-1}}$ by first $R_n$ leading left singular vectors of $\bC$;
		\State Un-matricization: $\eG^{(n)}=\text{reshape}(\bG^{(n)}_{(2)})\in \real^{R_{n-1}\times  R_n\times  I_n}$;
		\State Get the new matrix $\bC = \bSigma\bV^\top\in \real^{R_n\times I_{n+1}I_{n+2}\ldots I_N}$;
		\EndFor
		\State Get the last core tensor: $\eG^{(N)}=\text{reshape}(\bC)\in \real^{R_{N-1}\times  R_N\times  I_N}$;
		\State Output $\eG^{(1)}, \eG^{(2)}, \ldots, \eG^{(N)}$; 
	\end{algorithmic} 
\end{algorithm}

The unfolding in \eqref{equation:tt_unfond1} and \eqref{equation:tt_unfond2} reveals a recursive algorithm to calculate the TT decomposition of $\eX$. 
For this tensor unfolding along mode-$1$:
$$
\abovebX{1} \in \real^{I_1 \times I_2I_3\ldots I_N}.
$$ 
And recall the matricization in mode-$2$ of $\eG^{(1)}$ is given by (Section~\ref{section:matricization-tensor-original}):
$$
\bG^{(1)}_{(2)} \in \real^{ R_1\times R_0I_1 }=\real^{ R_1 \times I_1 }. 
$$
Then the matrix  $\bG^{(1)}_{(2)}$ can be though of as the data-distilled version of $\abovebX{1}$, i.e., the row space of $\bG^{(1)}_{(2)}$ spans the same column space of $\abovebX{1}$. And therefore, the value of  $R_1$ can be determined by the numerical rank (Definition~\ref{definition:effective-rank-in-svd}) of the SVD of $\abovebX{1}$, i.e., obtained by discarding the singular values of $\abovebX{1}$ smaller than $\delta$. Similar to the Tucker and HOSVD methods, the rows of $\bG^{(1)}_{(2)}$ can be obtained as the first $R_1$ left singular vectors of $\abovebX{1}$. For the (truncated) SVD of $\abovebX{1}$:
$$
\begin{aligned}
\abovebX{1} &= \bU_1\bSigma_1\bV_1^\top
\gap\implies\gap
\boxed{
\bG^{(1)}_{(2)} = \bU_1^\top\in \real^{R_1\times I_1}
}.
\end{aligned}
$$
Now, what remains for decomposition is $\bSigma_1\bV_1^\top \in \real^{R_1\times I_2I_3\ldots I_N}$. Using a similar ``matrix unfolding" approach, suppose we reshape $\bSigma_1\bV_1^\top$ into a $I_2R_1 \times I_3I_4\ldots I_N$ matrix $\bC$:
$$
\bC \in \real^{I_2R_1 \times I_3I_4\ldots I_N} = \text{reshape}\bigg(\bSigma_1\bV_1^\top, (I_2R_1 ), (I_3I_4\ldots I_N) \bigg).
$$
The second core tensor $\eG^{(2)}$ can be obtained in a similar manner through the (truncated) SVD of $\bC$:
$$
\begin{aligned}
\bC &= \bU_2 \bSigma_2\bV_2^\top 
\gap\implies\gap
\boxed{
	\bG^{(2)}_{(2)} = \bU_2^\top\in \real^{ R_2\times I_2R_1}
}.
\end{aligned}
$$
The same process can continue, and eventually, all the $N$ core tensors $\{\eG^{(1)}, \eG^{(2)}, \ldots, \eG^{(N)}\}$ will be obtained via the set of (truncated) SVDs. This algorithm is known as the TT-SVD algorithm, as described in \citet{oseledets2011tensor}. The complete procedure is outlined in Algorithm~\ref{alg:tt-decomposi}. 

Furthermore, if the rank $R_n\leq \rank(\bC)$ in each iteration,  there always exists a best low-rank TT  approximation to $\eX$ in terms of Frobenius norm, denoted by $\eX_{best}$.
Additionally,  if  the truncation tolerance for the SVD of each unfolding is set to $\delta = \epsilon/\sqrt{N-1} \normf{\eX}$, then the TT-SVD is able to construct the quasi-optimal approximation, denoted by  $\eX_{SVD}$, satisfying
$$
\normf{\eX-\eX_{SVD}} \leq \sqrt{N-1} \normf{\eX-\eX_{best}}.
$$

\paragraph{Complexity and curse of dimensionality.} Let's consider again a scenario where all the TT-ranks are equal, i.e.,  $R_1=R_2=\ldots=R_{N-1} = R$, and where all the dimensions are the same, i.e.,  $I_1=I_2=\ldots=I_N = I$.  The computational complexity of the SVD calculation for  a matrix  of size $M\times N$  is $\mathcalO(MN^2)$ if (Section~\ref{section:comput-svd-in-svd}). 
Then we can demonstrate that the complexity of each step $n$ in Algorithm~\ref{alg:tt-decomposi} can be expressed as 
$
f(n) = (R I) (I^{(N-n)})^2.
$
Summing these complexities gives the total cost of the TT-SVD algorithm:
$$
\text{cost}= f(1)+f(2)+\ldots+f(N-1)= RI (I^{2(N-1)} + I^{2(N-2)}+ I^{2}). 
$$
Therefore, the complexity grows exponentially with the dimension $I$, and as a result,  the curse of dimensionality for the calculation is still involved.
However, the number of parameters is reduced from $\prod_{n=1}^{N}I_n$ to $\sum_{n=1}^{N} R_{n-1}R_nI_n$. More details are discussed in \citet{oseledets2011tensor}.

\paragraph{Alternative calculation methods.} We notice that the calculation of the TT decomposition relies on  rank-revealing decomposition, such as the SVD. 
Alternative methods, such as the rank-revealing QR (Section~\ref{section:rank-r-qr}), column-pivoted QR (Theorem~\ref{theorem:rank-revealing-qr-general}), CUR (Theorem~\ref{theorem:skeleton-decomposition}), UTV (Theorem~\ref{theorem:ulv-decomposition}), Column ID (Theorem~\ref{theorem:interpolative-decomposition}) can be applied in each iteration to perform the matrix decomposition for finding the spanning columns. 
We won't delve into the details here. For more information, refer to \citet{oseledets2011tensor, bigoni2016spectral}.

\begin{problemset}

\item Prove the matricization of the third-order tensor in Equation~\eqref{equation:cp-matriciza-third-1}.
\item \label{prob:vec_cpform} Show that the CP decomposition in Theorem~\ref{theorem:cp-decomp} can be represented as 
$vec(\eX)\approx \sum_{r=1}^{R}  \ba_r^{(N)}\kronecker  \ldots \kronecker  \ba_r^{(2)} \kronecker \ba_r^{(1)}$.
\textit{Hint: Use Exercise~\ref{exer:vec_rkoneten}.}

\item Delve into the footnote of Equation~\eqref{equation:cp-matriciza-third-1} for the Kolda, Kiers, and LMV matricizations related to the CP decomposition. Then, find all  six matricitizations for the Tucker decomposition in Equation~\eqref{equation:tucker-matriciza-third-1}. Similarly  find all the six matricitizations for the HOSVD.

\item Prove the general matricization forms of the Tucker decomposition and the HOSVD.

\item \label{prob:close_cp_tucset} Show that the set $CP(R)$ is not closed, and the set $Tucker(\overline{R}_1, \overline{R}_2, \ldots, \overline{R}_N)$ is closed (Definition~\ref{definition:c_rank},~\ref{definition:t_rank}).



\item Discuss the complexity of using  RRQR (Section~\ref{section:rank-r-qr}) or CPQR (Theorem~\ref{theorem:rank-revealing-qr-general}) for calculating the TT decomposition.

\item Consider a tensor $\eX\in\real^{3\times 4\times 5}$, where 
$$
\eX_{i, j, k} = i+j+k.
$$
Show that its TT cores are 
$$
\begin{aligned}
\eG^{(1)}
&=
\left(
\begin{bmatrix}
	1 & 1
\end{bmatrix},
\begin{bmatrix}
	2 & 1
\end{bmatrix},
\begin{bmatrix}
	3 & 1
\end{bmatrix}
\right), 
\gap
\eG^{(2)}
=
\left(
\begin{bmatrix}
	1 & 0 \\
	1 & 1
\end{bmatrix}, 
\begin{bmatrix}
	1 & 0 \\
	2 & 1
\end{bmatrix},
\begin{bmatrix}
	1 & 0 \\
	3 & 1
\end{bmatrix}, 
\begin{bmatrix}
	1 & 0 \\
	4 & 1
\end{bmatrix}
\right), \\
\gap
\eG^{(3)}
&=
\left(
\begin{bmatrix}
	1\\1
\end{bmatrix},
\begin{bmatrix}
	1\\2
\end{bmatrix},
\begin{bmatrix}
	1\\3
\end{bmatrix},
\begin{bmatrix}
	1\\4
\end{bmatrix},
\begin{bmatrix}
	1\\5
\end{bmatrix}
\right), 
\gap
\end{aligned}
$$
where the TT format uses 32 elements to describe the original 60 elements.

\item Use the \textit{Tensor Toolbox in Python: Bigoni D. Tensor toolbox, 2015} to find  the tensor-train decomposition of  the following tensor $\eX\in\real^{4\times 5 \times 3}$:
$$
\footnotesize
\begin{aligned}
\eX_{:,:,1}&=
\begin{bmatrix}
	1 & 5 & 9 & 13& 17\\
	2 & 6 & 10 & 14& 18\\
	3 & 7 & 11 & 15& 19\\
	4 & 8 & 12 & 16& 20
\end{bmatrix},\, 
\eX_{:,:,2}=
\begin{bmatrix}
	21 & 25 & 29 & 33& 37\\
	22 & 26 & 30 & 34& 38\\
	23 & 27 & 31 & 35& 39\\
	24 & 28 & 32 & 36& 40
\end{bmatrix},\, 
\eX_{:,:,3}&=
\begin{bmatrix}
	41 & 45 & 49 & 53& 57\\
	42 & 46 & 50 & 54& 58\\
	43 & 47 & 51 & 55& 59\\
	44 & 48 & 52 & 56& 60
\end{bmatrix}.
\end{aligned}
$$ 
\end{problemset}

%% file: chapter-app_proofs.tex
\part{Appendix}
\appendix

\newpage
\chapter{Proofs}
\begingroup
\hypersetup{
	linkcolor=structurecolor,
	linktoc=page,  
}
\minitoc \newpage
\endgroup

\section{The Fundamental Theorem of Linear Algebra}
\subsection{Dimension of Column Space and Row Space}\label{append:row-equal-column}
\label{app:theorem}

In Theorem~\ref{theorem:equal-dimension-rank}, we established that the row rank and column rank of any matrix $\bA\in \real^{m\times n}$ are equal, meaning the dimension of the column space of $\bA$ is the same as the dimension of its row space, using UTV or CR decomposition. In this appendix, we provide another proof using elementary methods, which also highlights the fundamental theorem of linear algebra. \index{Matrix rank}

\begin{proof}[{of Theorem~\ref{theorem:equal-dimension-rank}, the Third Way}]
We first notice that the null space of $\bA$ is orthogonal complementary to the row space of $\bA$: $\nspace(\bA) \bot \cspace(\bA^\top)$ (where the row space of $\bA$ is equivalent to the column space of $\bA^\top$), that is, vectors in the null space of $\bA$ are orthogonal to vectors in the row space of $\bA$. To see this, suppose $\bA$ has rows $\ba_1^\top, \ba_2^\top, \ldots, \ba_m^\top$ and $\bA=[\ba_1^\top; \ba_2^\top; \ldots; \ba_m^\top]$. For any vector $\bx\in \nspace(\bA)$, we have $\bA\bx = \bzero$, that is, $[\ba_1^\top\bx; \ba_2^\top\bx; \ldots; \ba_m^\top\bx]=\bzero$. And since the row space of $\bA$ is spanned by $\ba_1^\top, \ba_2^\top, \ldots, \ba_m^\top$, $\bx$ is perpendicular to all vectors in $\cspace(\bA^\top)$, which means $\nspace(\bA) \bot \cspace(\bA^\top)$.

Next, suppose the dimension of the row space of $\bA$ is $r$. \textcolor{mylightbluetext}{Let $\br_1, \br_2, \ldots, \br_r$ be a set of vectors in $\real^n$ and form a basis for the row space}. Then the $r$ vectors $\bA\br_1, \bA\br_2, \ldots, \bA\br_r$ are in the column space of $\bA$ and are linearly independent. To see this, suppose we have a linear combination of the $r$ vectors: $x_1\bA\br_1 + x_2\bA\br_2+ \ldots+ x_r\bA\br_r=\bzero$, that is, $\bA(x_1\br_1 + x_2\br_2+ \ldots+ x_r\br_r)=\bzero$, and the vector $\bv\triangleq x_1\br_1 + x_2\br_2+ \ldots+ x_r\br_r$ is in null space of $\bA$. But since $\{\br_1, \br_2, \ldots, \br_r\}$ is a basis for the row space of $\bA$, $\bv$ is thus also in the row space of $\bA$. We have shown that vectors from the null space of $\bA$ is perpendicular to vectors from the row space of $\bA$; thus, it holds that $\bv^\top\bv=0$ and $x_1=x_2=\ldots=x_r=0$. Therefore, \textcolor{mylightbluetext}{$\bA\br_1, \bA\br_2, \ldots, \bA\br_r$ are in the column space of $\bA$, and they are linearly independent}. This implies that the dimension of the column space of $\bA$ is larger than $r$. This result shows that \textbf{row rank of $\bA\leq $ column rank of $\bA$}. 

If we apply this process again to $\bA^\top$, we will establish that \textbf{column rank of $\bA\leq $ row rank of $\bA$}. This completes the proof that the row rank and column rank of $\bA$ are equal.
\end{proof}

From this proof, we can infer that if $\br_1, \br_2, \ldots, \br_r$ is a set of vectors in $\real^n$ that forms a basis for the row space, then \textcolor{mylightbluetext}{$\bA\br_1, \bA\br_2, \ldots, \bA\br_r$ form a basis for the column space of $\bA$}. 
This observation is formalized  in the following lemma.

\begin{lemma}[Column Basis from Row Basis]\label{lemma:column-basis-from-row-basis}
Let $\bA\in \real^{m\times n}$ be any matrix, and suppose that $\{\br_1, \br_2, \ldots, \br_r\}$ is a set of vectors in $\real^n$, which forms a basis for the row space. Then, $\{\bA\br_1, \bA\br_2, \ldots, \bA\br_r\}$ is a basis for the column space of $\bA$.
\end{lemma}

\index{Fundamental theorem}
\subsection{The Fundamental Theorem of Linear Algebra}\label{appendix:fundamental-rank-nullity}
For any matrix $\bA\in \real^{m\times n}$, it is straightforward to verify that any vector in the row space of $\bA$ is orthogonal to any vector in the null space of $\bA$. Suppose $\bx_n \in \nspace(\bA)$, then $\bA\bx_n = \bzero$, indicating that $\bx_n$ is orthogonal to every row of $\bA$, thus confirming our assertion.

Similarly, we can also demonstrate that any vector in the column space of $\bA$ is perpendicular to any vector in the null space of $\bA^\top$. 
Furthermore, the column space of $\bA$ together with the null space of $\bA^\top$ span the entire space $\real^m$, which is a key part of the fundamental theorem of linear algebra.

The fundamental theorem consists of two aspects: the dimensions of the subspaces and their orthogonality. 
The orthogonality can be easily verified as shown above. 
Additionally, when the row space has dimension $r$, the null space has dimension $n-r$. 
This is not immediately obvious and is proven in the following theorem.

\begin{figure}[h!]
	\centering
	\includegraphics[width=0.98\textwidth]{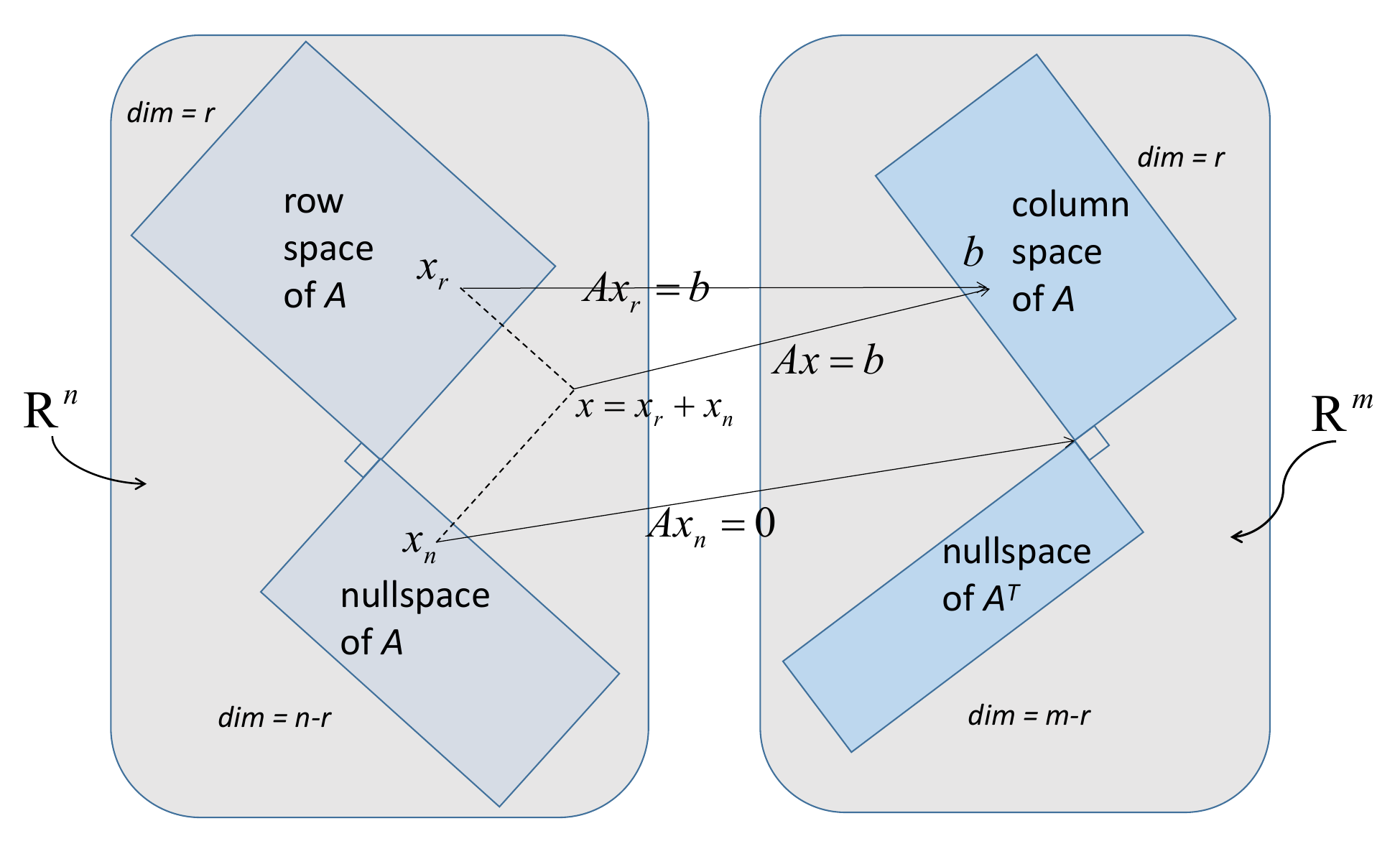}
	\caption{Two pairs of orthogonal subspaces in $\real^n$ and $\real^m$. $\dim(\cspace(\bA^\top)) + \dim(\nspace(\bA))=n$ and $\dim(\nspace(\bA^\top)) + \dim(\cspace(\bA))=m$. The null space component maps to zero as $\bA\bx_n = \bzero \in \real^m$. The row space component maps to the column space as $\bA\bx_r = \bA(\bx_r+\bx_n)=\bb \in \cspace(\bA)$.}
	\label{fig:lafundamental}
\end{figure}
\begin{theorem}[The Fundamental Theorem of Linear Algebra]\label{theorem:fundamental-linear-algebra}
Orthogonal Complement and Rank-Nullity Theorem: Let $\bA\in \real^{m\times n}$ be any matrix. Then, 
\begin{itemize}
\item $\nspace(\bA)$ is orthogonal complement to the row space $\cspace(\bA^\top)$ in $\real^n$: $\dim(\nspace(\bA))+\dim(\cspace(\bA^\top))=n$;
\item $\nspace(\bA^\top)$ is orthogonal complement to the column space $\cspace(\bA)$ in $\real^m$: $\dim(\nspace(\bA^\top))+\dim(\cspace(\bA))=m$;
\item For a rank-$r$ matrix $\bA$, $\dim(\cspace(\bA^\top)) = \dim(\cspace(\bA)) = r$, hence, $\dim(\nspace(\bA)) = n-r$ and $\dim(\nspace(\bA^\top))=m-r$.
\end{itemize}
\end{theorem}

\begin{proof}[of Theorem~\ref{theorem:fundamental-linear-algebra}]
Building upon the proof of Theorem~\ref{theorem:equal-dimension-rank} in Appendix~\ref{append:row-equal-column}, consider a set of vectors $\br_1, \br_2, \ldots, \br_r$  in $\real^n$ forming a basis for the row space.
Consequently, \textcolor{black}{$\{\bA\br_1, \bA\br_2, \ldots, \bA\br_r\}$ constitutes a basis for the column space of $\bA$}. 
Let $\bn_1, \bn_2, \ldots, \bn_k \in \real^n$ form a basis for the null space of $\bA$. Following again from the proof of Theorem~\ref{theorem:equal-dimension-rank} in Appendix~\ref{append:row-equal-column}, it follows that $\nspace(\bA) \bot \cspace(\bA^\top)$,  indicating orthogonality between $\br_1, \br_2, \ldots, \br_r$ and $\bn_1, \bn_2, \ldots, \bn_k$. Consequently, the set   $\{\br_1, \br_2, \ldots, \br_r, \bn_1, \bn_2, \ldots, \bn_k\}$ is linearly independent in $\real^n$.

For any vector $\bx\in \real^n $, $\bA\bx$ is in the column space of $\bA$. Then it can be expressed as a linear combination of $\bA\br_1, \bA\br_2, \ldots, \bA\br_r$: $\bA\bx \triangleq \sum_{i=1}^{r}a_i\bA\br_i$, which states that $\bA(\bx-\sum_{i=1}^{r}a_i\br_i) = \bzero$ and $\bx-\sum_{i=1}^{r}a_i\br_i$ is thus in $\nspace(\bA)$. Since $\{\bn_1, \bn_2, \ldots, \bn_k\}$ is a basis for the null space of $\bA$, $\bx-\sum_{i=1}^{r}a_i\br_i$ can be expressed as a linear combination of $\bn_1, \bn_2, \ldots, \bn_k$: $\bx-\sum_{i=1}^{r}a_i\br_i = \sum_{j=1}^{k}b_j \bn_j$, i.e., $\bx=\sum_{i=1}^{r}a_i\br_i + \sum_{j=1}^{k}b_j \bn_j$. That is, any vector $\bx\in \real^n$ can be expressed as $\{\br_1, \br_2, \ldots, \br_r, \bn_1, \bn_2, \ldots, \bn_k\}$ and the set forms a basis for $\real^n$. Thus, the dimensions satisfy: $r+k=n$, i.e., $\dim(\nspace(\bA))+\dim(\cspace(\bA^\top))=n$. Similarly, we can prove that $\dim(\nspace(\bA^\top))+\dim(\cspace(\bA))=m$.
\end{proof}

Figure~\ref{fig:lafundamental} demonstrates two pairs of such orthogonal subspaces and shows how $\bA$ maps a vector $\bx$ into the column space. The dimensions of the row space of $\bA$ and the null space of $\bA$ add to $n$. And the dimensions of the column space of $\bA$ and the null space of $\bA^\top$ add to $m$. The null space component maps to zero as $\bA\bx_{\bn} = \bzero \in \real^m$, which is the intersection of the column space of $\bA$ and the null space of $\bA^\top$. 
Conversely, the row space component transforms into the column space as $\bA\bx_{\br} = \bA(\bx_{\br} + \bx_{\bn})=\bb\in \real^m$. 

\subsection{Find the Basis of the Four Subspaces via  CR Decomposition}\label{appendix:cr-decomposition-four-basis}
Given the CR decomposition of a matrix $\bA=\bC\bR \in \real^{m\times n}$ (Section~\ref{section:cr-decomposition}), we have $\bR = [\bI_r, \bF ]\bP$, where $\bP$ is an $n\times n$ permutation matrix that rearranges the columns of the $r\times r$ identity matrix $\bI_r$ into the correct positions as shown in Section~\ref{section:rref-cr}. We can thus use the $r$ linearly independent columns of $\bC$ as the basis for the column space of $\bA$, i.e., basis for $\cspace(\bA)$.

Further, let $\bN
\triangleq
\bP^\top
\footnotesize
\begin{bmatrix}
-\bF \\
\bI_{n-r}
\end{bmatrix}$, 
where $\bI_{n-r}$ is an $(n-r)\times (n-r)$ identity matrix. 
Then,
$$
\begin{aligned}
\bA  \bN 
&=\begin{bmatrix}
\bC & \bC\bF
\end{bmatrix}\bP \bP^\top 
\begin{bmatrix}
-\bF \\
\bI_{n-r}
\end{bmatrix}=\begin{bmatrix}
\bC & \bC\bF
\end{bmatrix}
\begin{bmatrix}
-\bF \\
\bI_{n-r}
\end{bmatrix}=-\bC\bF+\bC\bF=\bzero.
\end{aligned}
$$
Moreover, the $n\times (n-r)$ matrix $\footnotesize\begin{bmatrix}
-\bF \\
	\bI_{n-r}
\end{bmatrix}$ has $n-r$ independent columns since $\bI_{n-r}$ is an identity matrix. 
And multiplying from the left by a permutation matrix $\bP^\top$ does not affect the linear independence of these columns. Thus, $\bN$ has $n-r$ linearly independent columns that belong to the null space of $\bA$. Furthermore, we have shown that the dimension of $\nspace(\bA)$ is $n-r$, thus the columns of $\bN$ form a basis for $\nspace(\bA)$.

Similarly, by applying the CR decomposition to the transpose $\bA^\top$, we can determine the bases for the row space of $\bA$ and the  null space of $\bA^\top$.

In Section~\ref{section:property-svd}, we further find the orthonormal basis for the four subspaces via SVD. 
For a more detailed discussion of this fundamental theorem of linear algebra, see \citet{lu2021revisit}, which provides seven figures illustrating the theorem from various perspectives.

\index{Fundamental theorem}
\subsection{The Fundamental Theorem of Linear Algebra: Least Squares View}\label{appendix:ls-fundation-theorem}
\begin{figure}[h!]
\centering
\includegraphics[width=0.999\textwidth]{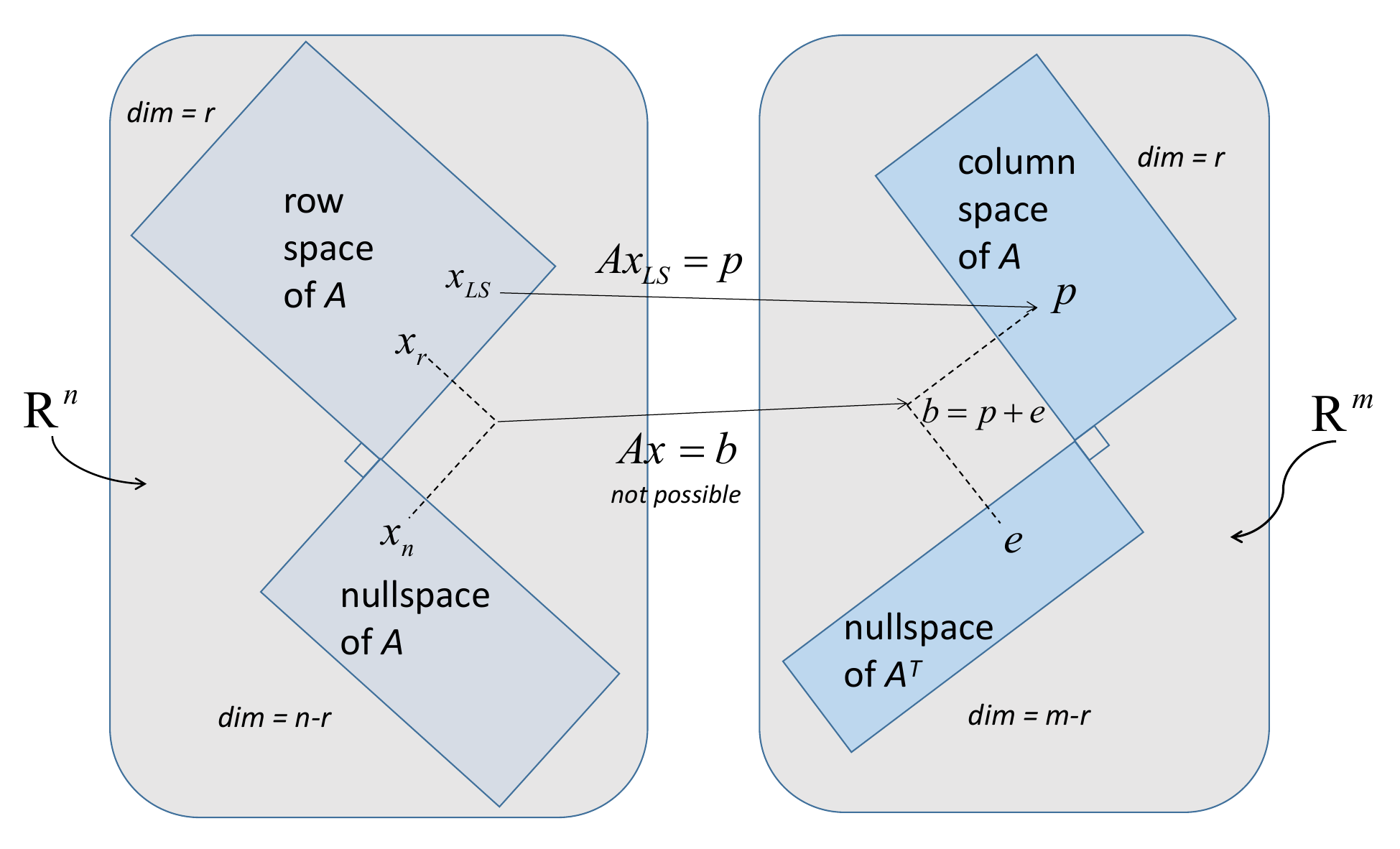}
\caption{Least squares: a row space to column space view. Transferring from the row space of $\bA$ to the column space of $\bA$. The least squares solution $\bx_{LS}$ minimizes the distance  $\normtwo{\bA\bx-\bb}^2$.}
\label{fig:lafundamental2-LS}
\end{figure}

The least squares problem is described in Section~\ref{section:application-ls-qr}, Section~\ref{section:ls-utv}, and Section~\ref{section:application-ls-svd}  using different matrix decompositions. As a recap, 
let's consider the overdetermined system $\bA\bx = \bb$, where $\bA\in \real^{m\times n}$ is the data matrix, and $\bb\in \real^m$ with $m\geq n$ is the observation vector. Typically, $\bA$ (or after pre-processing) has full column rank since the data from real work has a large chance to be unrelated. And the least squares (LS) solution is given by $\bx_{LS} = (\bA^\top\bA)^{-1}\bA^\top\bb$ for minimizing $\normtwo{\bA\bx-\bb}^2$, where $\bA^\top\bA$ is invertible since $\bA$ has full column rank and $\rank(\bA^\top\bA)=\rank(\bA)$.
\index{Least squares}

The goal of the least squares problem is to minimize the error $\bb-\bA\bx$ concerning the mean squared error (MSE, or $\ell_2$ norm).  
Since $\bA\bx$ is a linear combination of the columns of $\bA$, it always lies within the column space of $\bA$. 
Therefore, the best approximation of $\bb$ in the column space of $\bA$ is the projection $\bp$ of $\bb$ onto this space. 
Then the error vector $\be=\bb-\bp=\bb-\bA\bx_{LS}$ has minimal length and the error $\be$ is perpendicular to the column space of $\bA$. 
Therefore, $\be=\bb-\bA\bx_{LS}$ is in the null space of $\bA^\top$:
$$
\bA^\top(\bb-\bA\bx_{LS}) = \bzero  \qquad \implies \qquad \bA^\top\bb=\bA^\top\bA\bx_{LS},
$$
which is also known as the \textit{normal equation} of least squares. The relationship between $\be$ and $\bp$ is shown in Figure~\ref{fig:lafundamental2-LS}, where $\bb$ is decomposed into $\bb=\bp+\be$ ($\bp$ lies in the column space of $\bA$,  $\be$ lies in the null space of $\bA^\top$, and $\bp^\top\be=0$). Furthermore, it can be shown that $\bx_{LS}$ lies in the row space of $\bA$, i.e., it cannot be split into a combination of two components that are in the row space of $\bA$ and null space of $\bA$, respectively (see $\bx_{LS}$ via the pseudo-inverse of $\bA$ in Section~\ref{section:application-ls-svd}, $\bx_{LS}$ is a linear combination of the orthonormal basis of the row space).

\section{The Fundamental Theorem of Algebra}\label{appendix:fund_plain_alge}

The fundamental theorem of algebra states that every nonzero, single-variable polynomial of degree $n$ with complex coefficients has, counted with multiplicity, exactly $n$ complex roots.
More formally, consider a polynomial of degree $n$:
$$
p(z) = a_n z^n + a_{n-1} z^{n-1} + \cdots + a_1 z + a_0,
$$
where $a_n\neq 0$ and the coefficients $a_i$ are complex numbers. Then there exist exactly $n$  complex numbers $z_1, z_2, \ldots, z_n$
(which may include repetitions) such that:
$$
p(z_i) = 0 \text{ for } i = 1, 2, \ldots, n.
$$
The roots $z_i$ are the solutions to the equation  $p(z) = 0$. They can be real numbers, which are special cases of complex numbers with an imaginary part of zero.

The theorem guarantees that any polynomial equation will have a solution in the complex plane, which means that the field of complex numbers is algebraically closed. This is a powerful result that underpins much of modern mathematics and its applications.

It's worth noting that the theorem applies to polynomials over the complex numbers but not necessarily to polynomials over other fields, such as the real numbers, where a polynomial of degree $n$ may not have $n$ real roots. However, considering the polynomial over the complex numbers ensures that all roots can be found.

\section{The Fundamental Theorem of Calculus and Taylor's Expansion}\label{appendix:taylor-expansion}
Let $f(\cdot):\real^n\rightarrow \real$. The fundamental theorem of calculus states that
\begin{align}
\nabla f(\by) -\nabla f(\bx) &= \left( \int_{0}^{1} \nabla^2f(\bx+t(\by-\bx))dt \right) \cdot (\by-\bx);\\
\nabla f(\bx+\alpha\bd) - \nabla f(\bx) &= \int_{0}^{\alpha} \nabla^2 f(\bx+t\bd)\bd dt;\\
f(\by) - f(\bx) &= \int_{0}^{1} \langle \nabla f(\bx+t(\by-\bx)), \by-\bx  \rangle dt.
\end{align}
The last equality can be shown by letting $g(t)\triangleq f(\bx+t(\by-\bx))$. Then we can obtain $f(\by)$ by 
$
f(\by) = g(1) = g(0) + \int_{0}^{1} g'(t)dt.
$
The other equalities  can be demonstrated similarly. Moreover, use directional derivative for twice continuously differentiable functions, we also have
\begin{equation}
f(\by) = f(\bx) + \nabla f(\bx)^\top (\by-\bx) + \int_{0}^{1} (1-t) \frac{\partial^2 f(\bx+t(\by-\bx))}{\partial t^2} dt.
\end{equation}

\begin{theorem}[Taylor’s Expansion with Lagrange Remainder]
Let $f(x): \real\rightarrow \real$ be $k$-times continuously differentiable on the closed interval $I$ with endpoints $x$ and $y$, for some $k\geq 0$. If $f^{(k+1)}$ exists on the interval $I$, then there exists a $x^\star \in (x,y)$ such that 
$$
\begin{aligned}
f(x)& = f(y) + f^\prime(y)(x-y) +\ldots + \frac{f^{(k)}(y)}{k!}(x-y)^k
+ \frac{f^{(k+1)}(x^\star)}{(k+1)!}(x-y)^{k+1}\\
&=\sum_{i=0}^{k} \frac{f^{(i)}(y)}{i!} (x-y)^i + \frac{f^{(k+1)}(x^\star)}{(k+1)!}(x-y)^{k+1}.
\end{aligned}
$$ 
Taylor's expansion can be extended to a function of vector $f(\bx):\real^n\rightarrow \real$ or a function of matrix $f(\bX): \real^{m\times n}\rightarrow \real$.
\end{theorem}
Taylor's expansion, or also known as  \textit{Taylor's series}, approximates the function $f(x)$ around the value of $y$ by a polynomial in a single indeterminate $x$. To understand the origin of this series, we recall from the elementary calculus course that the approximated function around $\theta=0$ for $\cos (\theta)$ is given by 
$
\cos (\theta) \approx 1-\frac{\theta^2}{2}.
$
That is, the $\cos(\theta)$ function is approximated by a polynomial of degree 2. Suppose we want to approximate $\cos(\theta)$ more generally by a polynomial of degree 2: $ f(\theta) = c_1+c_2 \theta+ c_3 \theta^2$. An intuitive idea is to match  the function and its derivatives at $\theta=0$. That is,
$$\left\{
\begin{aligned}
\cos(0) &= f(0); \\
\cos^\prime(0) &= f^\prime(0);\\
\cos^{\prime\prime}(0) &= f^{\prime\prime}(0);\\
\end{aligned}
\right.
\quad\implies\quad 
\left\{
\begin{aligned}
1 &= c_1; \\
-\sin(0) &=0= c_2;\\
-\cos(0) &=-1= 2c_3.\\
\end{aligned}
\right.
$$
Solving these equations, we find: $f(\theta) = c_1+c_2 \theta+ c_3 \theta^2 = 1-\frac{\theta^2}{2}$, which matches our initial approximation $\cos (\theta) \approx 1-\frac{\theta^2}{2}$ around  $\theta=0$.

\section{Cochran's Theorem}\label{appendix:cochran-theorem}
In this appendix, we provide a proof for Theorem~\ref{theorem:cochran-theorem},  Cochran's theorem \citep{james1952notes, tan1975some, anderson1980cochran, gut2009multivariate}. 
To see this, we first give the definition of the Chi-square distribution.

\begin{SCfigure}[0.8][h]
\includegraphics[width=0.5\textwidth]{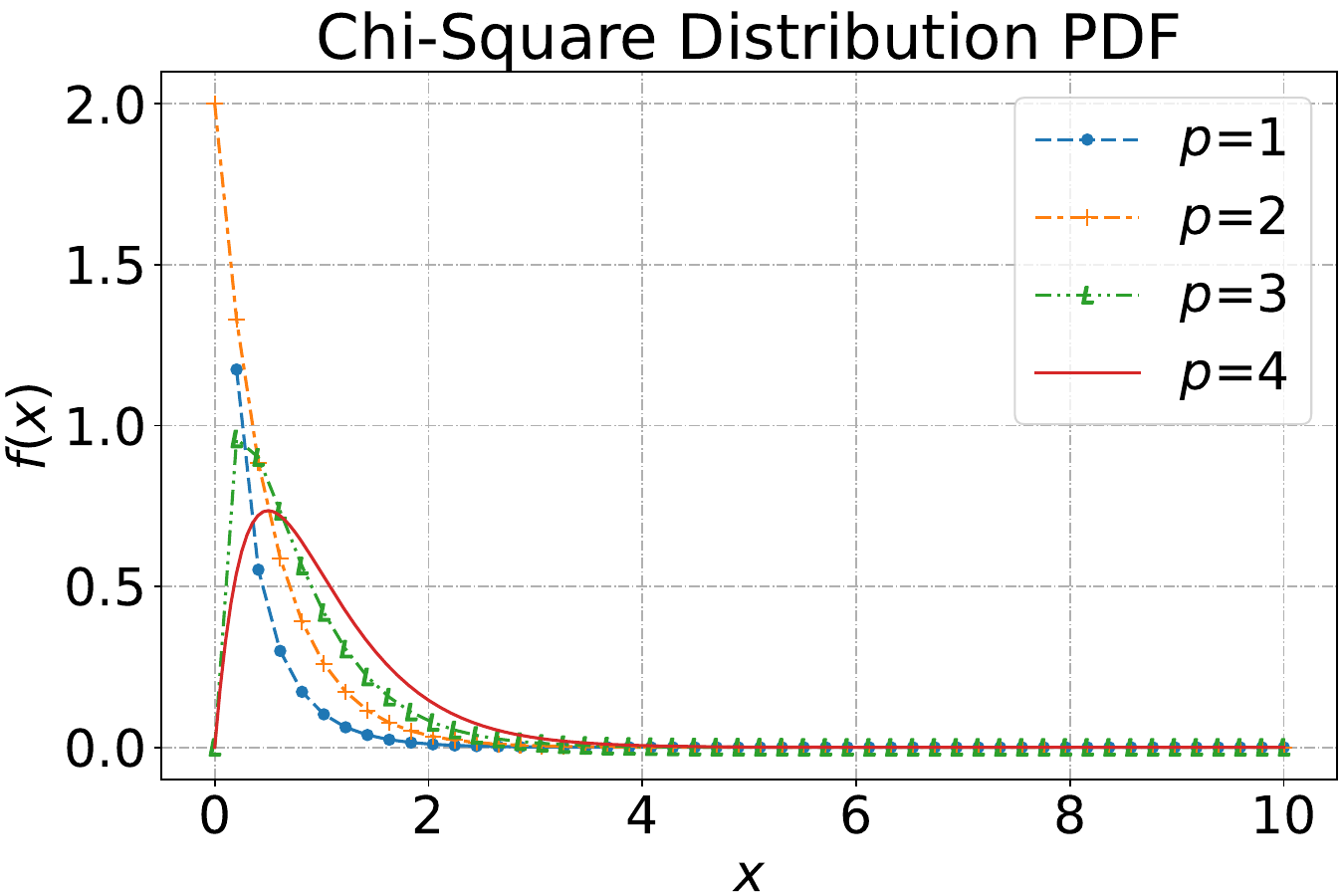}
\caption{Chi-square probability density functions for different values of the parameter $p$.}
\label{fig:chi_square_dist}
\end{SCfigure}

\begin{definition}[Chi-Square Distribution]\index{Chi-square distribution}
Let $\rva \sim \normal(\bzero, \bI_{p\times p})$, where $\normal(\bmu, \bSigma)$ denotes the multivariate normal distribution with mean $\bmu$ and variance $\bSigma$. Then, $\rx=\sum_i^p \ra_{i}^2$ follows a Chi-square distribution with $p$ degrees of freedom, denoted by $\rx \sim \chi_{(p)}^2$:
$$ 
f(x; p)=\left\{
\begin{aligned}
&\frac{1}{2^{p/2}\Gamma(\frac{p}{2})} x^{\frac{p}{2}-1} \exp\{-\frac{x}{2}\},& \mathrm{\,\,if\,\,} x \geq 0,  \\
&0 , &\mathrm{\,\,if\,\,} x <0.
\end{aligned}
\right.
$$
The mean and  variance of $\rx\sim \chi_{(p)}^2$ are given by $\Exp[\rx]=p$ and $\Var[\rx]=2p$, respectively.

The function $\Gamma(\alpha) = \int_{0}^{\infty}  t^{\alpha-1} e^{-t} dt $ is the Gamma function, which normalizes the distribution so that it integrates to 1.  In special cases when $y$ is a positive integer, $\Gamma(y) = (y-1)!$. Figure~\ref{fig:chi_square_dist} compares different parameters for the Chi-square distribution.
\end{definition}

\noindent

\begin{theorem}[Cochran's Theorem\index{Cochran's theorem}]\label{theorem:cochran-theorem}
Let $\rvy$ \footnote{Note that we use normal fonts of boldface lowercase letter $\rvy$ to denote a random vector variable and use italic letter $\by$ to denote a vector.} be random variable in $\real^n$. Then  $\rvy^\top\rvy$ can be factored into $k>0$ terms:
$$
\rvy^\top\rvy = \rvy^\top\bA_1\rvy+ \rvy^\top\bA_2\rvy +\ldots + \rvy^\top\bA_k\rvy. 
$$
The matrices $\bA_i$'s must satisfy   the following conditions:
\begin{enumerate}[(i).]
\item  $\bA_1, \bA_2, \ldots, \bA_k$ are positive semidefinite (PSD);
\item  $\bA_1+\bA_2+ \ldots+ \bA_k = \bI_n$ is the $n\times n$ identity matrix;
\item  Let $r_i \triangleq \rank(\bA_i)$, and $r_1 + r_2+\ldots +r_k=n$. 
\end{enumerate}

\item Let $\rq_i \triangleq \rvy^\top\bA_i\rvy$ for all $i\in\{1,2,\ldots,k\}$. Then, the following results hold:
\begin{enumerate}
\item  If $\rvy \sim \normal(\textcolor{mylightbluetext}{\bzero}, \sigma^2\bI)$, then $\rq_i \sim \sigma^2\chi_{(r_i)}^2$;
\item  If $\rvy \sim \normal(\textcolor{mylightbluetext}{\bmu}, \sigma^2\bI)$ and $\bmu^\top\bA_i\bmu=0$, then $\rq_i \sim \sigma^2\chi_{(r_i)}^2$;
\item  $\rq_i$'s are independent to each other.
\end{enumerate}

\end{theorem}

To prove  Cochran's theorem, we need the following lemma.
\begin{lemma}[Idempotent Decomposition: Rank-Additivity\index{Idempotent decomposition}]\label{lemma:cochran-lemma-symmetric-decomposition}
Given $n\times n$ square matrices $\bA_1, \bA_2, \ldots, \bA_k$ satisfying $\bA_1+\bA_2+\ldots+\bA_k = \bI_n$,  the following three conditions are equivalent:
\begin{enumerate}[(i).]
\item  $\bA_i^2 = \bA_i$, for all $i\in \{1, 2, \ldots, k\}$, i.e., $\bA_i$'s are idempotent;
\item  \textbf{Rank is additive:} $\rank(\bA_1)+\rank(\bA_2)+\ldots+\rank(\bA_k)=n$;
\item  \textbf{Orthogonality:} $\bA_i\bA_j = \bzero$ for all $i\neq j$ and $i,j \in \{1, 2, \ldots, k\}$.
\end{enumerate}
\end{lemma}
\begin{proof}[of Lemma~\ref{lemma:cochran-lemma-symmetric-decomposition}]
From (i) to (ii), by Lemma~\ref{lemma:rank-of-symmetric-idempotent2}, the trace and rank of any idempotent matrix are the same. Then, we have   
$
\sum_{i=1}^{k} \rank(\bA_i) = \sum_{i=1}^{k} \trace(\bA_i) = \trace(\bI_n) = n.
$
	
From (ii) to (iii), consider the following block Gaussian elimination for a $(k+1)\times(k+1)$ block matrix (that is, a $(k+1)n\times(k+1)n$ matrix), where the $(2,2), (3,3), \ldots, (k+1, k+1)$ blocks are $\bA_1, \bA_2, \ldots, \bA_k$, respectively:
$$
\footnotesize
\begin{aligned}
	\bX=&\begin{bmatrix}
		\bzero_n & & &\\  
		& \bA_1 & & \\  
		&  &  \ddots  &\\  
		&  &   &  \bA_k\\  
	\end{bmatrix}\stackrel{\bE_1\times}{\rightarrow}
	\begin{bmatrix}
		\bzero_n & \bA_1  & \ldots& \bA_k\\  
		& \bA_1 &\\  
		&  &   \ddots  &\\  
		&  &    &  \bA_k\\  
	\end{bmatrix}\stackrel{\bE_2\times}{\rightarrow}
	\begin{bmatrix}
		\bI_n & \bA_1 & \ldots& \bA_k\\  
		\bA_1& \bA_1 &\\  
		\vdots	&  &   \ddots  &\\  
		\bA_k	&  &   &  \bA_k\\  
	\end{bmatrix}\stackrel{\bE_3\times}{\rightarrow}\\
	&\begin{bmatrix}
		\bI_n & \bA_1 & \ldots& \bA_k\\  
		& \bA_1-\bA_1^2  & &-\bA_1\bA_k\\  
		&  \vdots &    \ddots  & \vdots\\  
		&  -\bA_k\bA_1&     &  \bA_k-\bA_k^2\\  
	\end{bmatrix}\stackrel{\bE_4\times}{\rightarrow}
	\begin{bmatrix}
		\bI_n &  & & \\  
		& \bA_1-\bA_1^2  & &-\bA_1\bA_k\\  
		&  \vdots &    \ddots  & \vdots\\  
		&  -\bA_k\bA_1&     &  \bA_k-\bA_k^2\\  
	\end{bmatrix}=\bY,
\end{aligned}
$$
where blank entries indicate zeros. And 
\begin{itemize}
\item $\bE_1$ is used to add row-2, row-3, $\ldots$, row-(k+1) to row-1;

\item $\bE_2$ is used to add column-2, column-3, $\ldots$, column-(k+1) to column-1; 

\item $\bE_3$ is used to subtract the row-2 by $\bA_1$*(row-1), subtract the row-3 by $\bA_2$*(row-1), $\ldots$;

\item $\bE_4$ is used to subtract the column-2 by $\bA_1$*(column-1), subtract column-3 by $\bA_2$*(column-1), $\ldots$.
\end{itemize}
We notice that elementary operations/transformations do not alter the rank of a matrix. Since $\bX$ is of rank $\sum_{i=1}^{k}r_k \triangleq n$ (where $r_i=\rank(\bA_i)$), and $\bI_n$ in $\bY$ is of rank $n$ as well. We must have
$$
\begin{bmatrix}
\bA_1-\bA_1^2  & &-\bA_1\bA_k\\  
\vdots &    \ddots  & \vdots\\  
-\bA_k\bA_1&     &  \bA_k-\bA_k^2\\  
\end{bmatrix}=\bzero,
$$
which implies $\bA_i\bA_j = \bzero$ for all $i\neq j$ and $i,j \in \{1, 2, \ldots, k\}$.

From (iii) to (i), we have 
$
\begin{aligned}
\bA_i &= \bA_i \bI_n
=\bA_i(\bA_1+\bA_2+\ldots+\bA_k)
=\bA_i^2.
\end{aligned}
$
This completes the proof.
\end{proof}

Now we are ready to prove  Cochran's Theorem as follows:
\begin{proof}[of Theorem~\ref{theorem:cochran-theorem}]
\paragraph{Case 1.} If $\rvy \sim \normal(\bzero, \sigma^2\bI)$: 
From Lemma~\ref{lemma:cochran-lemma-symmetric-decomposition}, $\bA_i$'s are idempotent.
By spectral theorem~\ref{theorem:spectral_theorem} and Lemma~\ref{proposition:eigenvalues-of-projection2} (the only possible eigenvalues of idempotent matrices are 0 and 1), we can rewrite  $\rq_i$ as $\rq_i = \rvy^\top\bA_i\rvy = \rvy^\top(\bQ \bLambda\bQ^\top)\rvy$, where $\bA_i=\bQ \bLambda\bQ^\top$ is the spectral decomposition of $\bA_i$ ($\bLambda$ contains only 1 and 0 on the diagonal). 
Since rotations on the normal distribution do not effect the distribution~\footnote{Rotations on the Gaussian distribution do not effect the distribution. That is for any orthogonal matrix $\bQ$ with $\bQ\bQ^\top=\bQ^\top\bQ=\bI$, if $\rvv\sim \normal(\bzero, \sigma^2\bI)$, then $\bQ\rvv\sim \normal(\bzero, \sigma^2\bI)$.}, we can define $\boldeta$ and obtain:
\begin{equation*}
\begin{aligned}
\boldeta &\triangleq \bQ^\top\rvy \sim \normal(\bzero, \sigma^2\bI)
\implies
\rq_i =  \boldeta^\top \bLambda \boldeta \sim \sigma^2 \chi^2_{\mathrm{rank(\bA_i)}}=\sigma^2 \chi^2_{(r_i)},
\quad\forall i.
\end{aligned}
\end{equation*}

\paragraph{Case 2.}
If $\rvy \sim \normal(\bmu, \sigma^2\bI)$, and $\bmu^\top \bA_i \bmu = 0$:
Let $\rp_i \triangleq (\rvy-\bmu)^\top\bA_i(\rvy-\bmu)$. 
Similarly, we can rewrite  $\rp_i$ as $\rp_i = (\rvy-\bmu)^\top\bA_i(\rvy-\bmu) = (\rvy-\bmu)^\top(\bQ \bLambda\bQ^\top)(\rvy-\bmu)$, where $\bA_i=\bQ \bLambda\bQ^\top$ is the spectral decomposition of $\bA_i$. Using the fact that rotations on the normal distribution do not effect the distribution, we can define $\boldeta$ and obtain:
\begin{equation*}
\begin{aligned}
\boldeta &\triangleq \bQ^\top(\rvy-\bmu) \sim \normal(\bzero, \sigma^2\bI)
\implies
\rp_i =  \boldeta^\top \bLambda \boldeta \sim \sigma^2 \chi^2_{\mathrm{rank(\bA_i)}}=\sigma^2 \chi^2_{(r_i)},
\quad\forall i.
\end{aligned}
\end{equation*}
Next,  we decompose  $\rp_i$ into
$$
\begin{aligned}
\rp_i &= \rvy^\top\bA_i \rvy - 2\rvy^\top\bA_i\bmu + \bmu^\top \bA_i\bmu
= \rq_i - 2\rvy^\top(\bQ \bLambda\bQ^\top)\bmu + \bmu^\top (\bQ \bLambda\bQ^\top)\bmu.
\end{aligned}
$$
Since we assume $\bmu^\top \bA_i\bmu=\bmu^\top (\bQ \bLambda\bQ^\top)\bmu=0$, and $\bLambda$ contains only 1 and 0 on the diagonal, we have $\bLambda =\bLambda \bLambda^\top$. That is,
$$
\bmu^\top (\bQ \bLambda\bQ^\top)\bmu = \bmu^\top (\bQ \bLambda\bLambda^\top\bQ^\top)\bmu = \normtwo{\bLambda^\top\bQ^\top\bmu}^2 = 0,
$$
which implies that $\rq_i\stackrel{d}{=}\rp_i\sim\sigma^2 \chi^2_{(r_i)}$.

\paragraph{Case 3.}
From Lemma~\ref{lemma:cochran-lemma-symmetric-decomposition}, $\bA_i \bA_j = \bzero$ if $i\neq j$.
Let $\bA_i = \bQ_i\bLambda_i \bQ_i^\top$ and $\bA_j = \bQ_j\bLambda_j \bQ_j^\top$ be the spectral decomposition of $\bA_i$ and $\bA_j$, respectively. Then,  we have 
\begin{equation}\label{equation:cochran-proof-zero}
\begin{aligned}
\bA_i \bA_j &= \bQ_i\bLambda_i \bQ_i^\top \bQ_j\bLambda_j \bQ_j^\top = \bzero
\gap \implies\gap
\bQ_i^\top\bA_i \bA_j\bQ_j =\bLambda_i \bQ_i^\top \bQ_j\bLambda_j  = \bzero.
\end{aligned}
\end{equation}
Write out $\rq_i$ and $\rq_j$: 
$$
\begin{aligned}
\rq_i &= \rvy^\top \bQ_i\bLambda_i \bQ_i^\top \rvy = \rvy^\top \bQ_i\bLambda_i\bLambda_i^\top \bQ_i^\top \rvy
\quad\text{and}\quad 
\rq_j = \rvy^\top \bQ_j\bLambda_j \bQ_j^\top \rvy = \rvy^\top \bQ_j\bLambda_j\bLambda_j^\top \bQ_j^\top \rvy.
\end{aligned}
$$
Let $\rva_i \triangleq \bLambda_i^\top \bQ_i^\top \rvy$ and $\rva_j \triangleq \bLambda_j^\top \bQ_j^\top \rvy$. Then,
$$
\begin{aligned}
	\Cov[\rva_i, \rva_j] &= \bLambda_i^\top \bQ_i^\top \Cov[\rvy, \rvy] \bQ_j\bLambda_j = \sigma^2\bLambda_i^\top \bQ_i^\top  \bQ_j\bLambda_j = \bzero,
\end{aligned}
$$
where the last equality follows from Equation~\eqref{equation:cochran-proof-zero}. This implies $\Cov[\rq_i,\rq_j]=0$ since $\rq_i = \rva_i^\top \rva_i$ and $\rq_j = \rva_j^\top \rva_j$.
\end{proof}

An alternative proof is provided in \citet{gut2009multivariate}. 
However, the author does not provide an inductive case for $k>2$.  Interested readers can refer to it.

\index{KKT conditions}
\section{KKT Conditions and Weierstrass Theorem}\label{appendix:KKT_proj}

We illustrate the \textit{Karush-Kuhn-Tucker (KKT) conditions}  for linearly constrained problems in the following theorem without a proof. A detailed discussion can be found in \citet{beck2014introduction}, and the proof relies on  \textit{Farkas' lemma}.

\begin{theorem}[KKT Conditions for Linearly Constrained Problems, Theorem 10.7 in \citet{beck2014introduction}]\label{theorem:kktcond}
Consider the minimization problem
$$
(\text{P})\qquad
\begin{aligned}
\min  f(\bx), \qquad 
\text{s.t.}\gapthree &g_i(\bx) =\ba_i^\top \bx - b_i\leq 0, \gapthree i\in\{1,2,\ldots, m\},\\
&h_j(\bx)=\bc_j^\top\bx-d_j=0, \gapthree j\in\{1,2,\ldots, p\},
\end{aligned}
$$
where $f(\cdot)$ is  continuously differentiable over $\real^n$, $\ba_1,\ba_2,\ldots, \ba_m,\bc_1,\bc_2,\ldots, \bc_p\in\real^n$, and $b_1,b_2,\ldots,b_m, d_1,d_2,\ldots,d_p\in\real$.
Then, 
\begin{itemize}[(a)]
\item \textbf{Necessity of the KKT conditions.}
Let $\bx^*$ be a local minimum point of (P). Then, there exist $\lambda_1, \lambda_2, \ldots,\lambda_m\geq 0$ and $\mu_1, \mu_2, \ldots, \mu_p\in\real$ such that 
\begin{align}
&\nabla f(\bx^*) + \sum_{i=1}^{m} \lambda_i\ba_i +\sum_{j=1}^{p}\mu_j\bc_j=\bzero ,\label{equa:kkt1}\\
&\lambda_i (\ba_i^\top\bx^*-b_i)=0, \gapthree i=\{1,2,\ldots, m\}. \label{equa:kkt2}
\end{align}

\item \textbf{Sufficiency in the convex case.} Suppose further that the function $f$ is convex over $\real^n$, and $\bx^*$ is a feasible solution of (P) for which there exist $\lambda_1, \lambda_2, \ldots,\lambda_m\geq 0$ and $\mu_1, \mu_2, \ldots, \mu_p\in\real$ such that \eqref{equa:kkt1} and \eqref{equa:kkt2} are satisfied. Then, $\bx^*$ is an optimal solution of (P).
\end{itemize}

Moreover, we can define the \textit{Lagrangian function}:
$$
L(\bx, \blambda, \bmu) = f(\bx) + \sum_{i=1}^{m} \lambda_ig_i(\bx) + \sum_{j=1}^{p} \mu_jh_j(\bx).
$$
The KKT condition \eqref{equa:kkt1} can be derived from the root of the gradient of the Lagrangian function:
$
\nabla_{\bx} L(\bx, \blambda, \bmu) 
= \nabla f(\bx) + \sum_{i=1}^{m} \lambda_i \nabla g_i(\bx) + \sum_{j=1}^{p} \mu_j \nabla h_j(\bx)
=\bzero.
$
\end{theorem}

The KKT conditions also hold under the generalized Slater's condition.
\begin{theorem}[KKT Conditions under the Generalized Slater's Condition]
Consider the minimization problem
$$
(\text{P})\qquad
\begin{aligned}
\min f(\bx),\qquad 
\text{s.t.}\gapthree 
&g_i(\bx) \leq 0, \gapthree i\in\{1,2,\ldots, m\}, \gapthree \text{(convex)}\\
&h_j(\bx)\leq 0, \gapthree j\in\{1,2,\ldots, p\}, \gapthree \text{(affine)}\\
&s_k(\bx)= 0, \gapthree k\in\{1,2,\ldots, q\},  \gapthree\text{(affine)}\\
\end{aligned}
$$
where $f, g_1, \ldots, g_m$ are a continuously differentiable over $\real^n$, $h_j, s_k$ for all $j,k$ are affine. 
Then, we have the following results:
\begin{itemize}[(a)]
\item \textbf{Necessity of the KKT conditions.}
Let $\bx^*$ be a local minimum point of (P). 
Suppose that there exists $\widehat{\bx}\in\real^n$ such that the following generalized Slater's condition holds:
$$
\text{generalized Slater's condition: }
\left\{
\begin{aligned}
g_i(\widehat{\bx}) &\textcolor{mylightbluetext}{<}0, \gapthree i\in\{1,2,\ldots, m\}, \\
h_j(\widehat{\bx})&\leq 0, \gapthree j\in\{1,2,\ldots, p\},\\
s_k(\widehat{\bx})&= 0, \gapthree k\in\{1,2,\ldots, q\},
\end{aligned}
\right.
$$
Then, there exist $\lambda_1, \lambda_2, \ldots,\lambda_m, \eta_1, \eta_2,\ldots,\eta_p\geq 0$ and $\mu_1, \mu_2, \ldots, \mu_q\in\real$ such that 
\begin{align}
&\nabla f(\bx^*) + \sum_{i=1}^{m} \lambda_i\nabla g_{i}(\bx^*) +\sum_{j=1}^{p}\eta_j \nabla h_j(\bx^*) +\sum_{k=1}^{q}\mu_k \nabla s_k(\bx^*) =\bzero,\label{equa:kkt1_gen}\\
&\lambda_i g_i(\bx^*)=0, \gapthree i=\{1,2,\ldots, m\}; 
\quad 
\eta_j h_j(\bx^*)=0, \gapthree j=\{1,2,\ldots, p\}.\label{equa:kkt2_gen}
\end{align}

\item \textbf{Sufficiency in the convex case.} Suppose further that the function $f$ is convex over $\real^n$, and $\bx^*$ is a feasible solution of (P) for which there exist $\lambda_1, \lambda_2, \ldots,\lambda_m\geq 0$ and $\mu_1, \mu_2, \ldots, \mu_p\in\real$ such that \eqref{equa:kkt1_gen} and \eqref{equa:kkt2_gen} are satisfied. Then, $\bx^*$ is an optimal solution of (P).
\end{itemize}

\end{theorem}

The KKT conditions can be applied to find the orthogonal projection onto an affine space:
$$
\begin{aligned}
\mathop{\min}_{\bx\in\real^n}\gapthree &f(\bx) =\normtwo{\bx-\by}^2 \qquad \text{s.t.}\gapthree & \bA\bx=\bb,
\end{aligned}
$$
where $\bA\in\real^{m\times n}$ and $\bb\in\real^m$.
This problem aims to find the closest point to $\by$ in the \textit{affine space} $\mathcalV=\{\bx\in\real^n: \bA\bx=\bb\}$.
Since the norm function is convex, the KKT conditions are necessary and sufficient for obtaining the optimal solution of the problem.
Therefore, the KKT conditions are 
$
2(\bx-\by) + 2\bA^\top\bmu = \bzero
$
and 
$
\bA\bx = \bb.
$
Thus, we obtain 
$
\bA(\by-\bA^\top\bmu) = \bb
\implies 
\bA\by - \bb = \bA\bA^\top\bmu.
$
If $\bA$ has full row rank, then
$
\bmu = (\bA\bA^\top)^{-1} (\bA\by - \bb).
$
Therefore, the optimal solution $\bx^*$ can be obtained by 
$
\begin{aligned}
\bx^* &= \by - \bA^\top \bmu
= \by - \bA^\top  (\bA\bA^\top)^{-1} (\bA\by - \bb).
\end{aligned}
$

When $\bA$ is a vector $\ba^\top$, the problem reduces to finding the orthogonal projection onto the \textit{hyperplane}:
$
\mathcalV' = \{\bx\in\real^n: \ba^\top\bx = b\}.
$
Then, the optimal solution is 
$
\bx^* =\by-\frac{ (\ba^\top\by - b)}{\ba^\top\ba} \ba.
$

We  provide the Weierstrass for optimization problems involving  continuous or proper closed functions.
\begin{theorem}[Weierstrass Theorem and Variants]\label{theorem:weierstrass_them}
We consider the Weierstrass theorem and its variants for different types of functions and sets:
\begin{enumerate}
\item Let $f:\sS\rightarrow (-\infty, \infty]$ be a \textbf{continuous function} defined over a \textit{nonempty and compact (closed and bounded) set} $\sS\subseteq\real^n$.
Then there exists a global minimum point of $f$ over $\sS$ and a global maximum point of $f$ over $\sS$.
\item Let $f:\sS\rightarrow  (-\infty, \infty]$ be a \textbf{continuous and coercive function} and let $\sS\subseteq\real^n$ be a \textit{nonempty and closed set}. Then $f$ has a global minimum point over $\sS$. (\text{The coerciveness ensures the function is bounded over a subset.})

\item Let $f:\sS\rightarrow (-\infty, \infty]$ be a \textbf{proper closed function}, and one of the following is satisfied:
\begin{enumerate}
\item   $f$ is defined over a \textit{nonempty and compact set} $\sS$. 
\item   assume there exists a \textit{nonempty and bounded level set} $\text{Lev}(f, \alpha)=\{\bx\in\sS: f(\bx)\leq \alpha\}$. 
\item  $f$ is coercive, i.e., $\mathop{\lim}_{\normtwo{\bx}\rightarrow\infty}=\infty$.
\end{enumerate}
Then there exists a global minimum point of $f$ over $\sS$. And the set of minimizers $\{\bx\in\sS:  f(\bx)\leq  f(\by), \forall \by\in \sS\}$  of $\mathopmin{\sS} f(\bx)$ is non-empty and compact.
\end{enumerate}

\end{theorem}
Note the coerciveness in (2) ensures the function is bounded over a subset. The three conditions of (3.a)$\sim$(3.c) essentially ensure that the minimum value of $f(\bx)$ cannot be attained at infinity.

\section{Generalized Singular Value Decomposition (GSVD)*}
Following from \citet{golub2013matrix}, we provide a concise overview of the generalized singular value decomposition (GSVD). See also \citet{zhang2017matrix, bai1993computing, zha1989restricted, golub2013matrix, paige1981towards} for a more detailed discussion on GSVD.

\begin{theoremHigh}[CS Decomposition]
Let 
$$
\begin{blockarray}{cccc}
\begin{block}{c[cc]c}
&	\bQ_{11} & \bQ_{12} & m_1  \\
\bQ=	&	\bQ_{21} & \bQ_{22}  & m_2 \\
\end{block}
& n_1 & n_2 &  \\
\end{blockarray},
$$
be an orthogonal matrix with $m_1 \geq n_1$ and $m_1 \geq m_2$. Define the nonnegative integers $p$ and $q$ by $p\triangleq\max\{0, n_1-m_2\}$ and $q\triangleq\max\{0, m_2-n_1\}$, respectively. There exist orthogonal matrices $\bU_1 \in \real^{m_1\times m_1}$, $\bU_2\in \real^{m_2\times m_2}$, $\bV_1\in \real^{n_1\times n_1}$, and $\bV_2\in \real^{n_2\times n_2}$ such that if 
$$
\bU \triangleq \begin{bmatrix}
\bU_1 & \bzero \\
\bzero & \bU_2
\end{bmatrix}
\qquad 
\text{and}
\qquad 
\bV \triangleq \begin{bmatrix}
\bV_1 & \bzero \\
\bzero & \bV_2
\end{bmatrix},
$$
then 
$$
\begin{blockarray}{ccccccc}
\begin{block}{c[ccccc]c}
&	\bI & 0 & 0 & 0 & 0 & p  \\
& 0 & \bC & \bS & 0	 & 0  & n_1-p \\
\bU^\top\bQ\bV=	& 0 & 0 & 0 &0	 &  \bI  & m_1-n_1 \\
& 0 & \bS & -\bC & 0	 & 0  & n_1-p \\
& 0 & 0 & 0 & \bI  & 0  & q \\
\end{block}
& p & n_1-p &  n_1-p& q &  m_1-n_1&  \\
\end{blockarray},
$$
where 
$$
\begin{aligned}
\bC &= \diag(\cos(\theta_{p+1}), \ldots, \cos(\theta_{n_1})) = \diag(c_{p+1}, \ldots, c_{n_1}),\\
\bS &= \diag(\sin(\theta_{p+1}), \ldots, \sin(\theta_{n_1})) = \diag(s_{p+1}, \ldots, s_{n_1}),
\end{aligned}
$$
and $0\leq \theta_{p+1} \leq \ldots \leq \theta_{n_1} \leq \pi/2$.
\end{theoremHigh}

The proof can be found in \citet{paige1981towards} and a thin version of the CS decomposition is provided in \citet{golub2013matrix}.


\begin{theoremHigh}[Generalized Singular Value Decomposition]\label{theorem:gsvd-main}
Let $\bA \in \real^{m_1\times n}$ and $\bB \in \real^{m_2\times n}$ with $m_1 \geq n$, and let 
$$
r = \rank\left(\begin{bmatrix}
\bA \\
\bB
\end{bmatrix}\right).
$$
Then there exist orthogonal matrices $\bU_1 \in \real^{m_1\times m_1}$ and $\bU_2\in \real^{m_2\times m_2}$, and invertible matrix $\bX \in \real^{n\times n}$ such that 
$$
\begin{aligned}
&\begin{blockarray}{ccccc}
\begin{block}{c[ccc]c}
&	\bI & \bzero & \bzero & p  \\
\bU_1^\top \bA \bX = \bD_A = 	&  \bzero   &	\bS_A & \bzero  & r-p \\
&	\bzero  &\bzero& \bzero & m_1-r  \\
\end{block}
& p & r-p &  n-r & \\
\end{blockarray}\\
&\begin{blockarray}{ccccc}
\begin{block}{c[ccc]c}
&	\bzero & \bzero & \bzero & p  \\
\bU_2^\top \bB \bX = \bD_B = 	&   \bzero  &	\bS_B & \bzero  & r-p \\
&	\bzero  &\bzero& \bzero & m_2-r  \\
\end{block}
& p & r-p &  n-r & \\
\end{blockarray},
\end{aligned}
$$
where $p = \max\{r-m_2, 0\}$, $\bS_A = \diag(\alpha_{p+1}, \ldots, \alpha_r)$,  
$\bS_B = \diag(\beta_{p+1}, \ldots, \beta_r)$, and 
$$
\alpha_i^2 + \beta_i^2 = 1, \qquad \forall i \in \{p+1, \ldots, r\}.
$$
\end{theoremHigh}
\begin{proof}[of Theorem~\ref{theorem:gsvd-main}]
Suppose the SVD of $\footnotesize\begin{bmatrix}
\bA \\
\bB
\end{bmatrix}$ is given by 
$$
\begin{bmatrix}
\bA \\
\bB
\end{bmatrix}
=
\begin{bmatrix}
\bQ_{11} & \bQ_{12} \\
\bQ_{21} &\bQ_{22}
\end{bmatrix}
\begin{bmatrix}
\bSigma_r & \bzero \\
\bzero & \bzero 
\end{bmatrix}
\bP^\top,
$$
where $\bSigma_r \in \real^{r\times r}$ is nonsingular, $\bQ_{11}\in \real^{m_1\times r}$, and $\bQ_{21}\in \real^{m_2\times r}$. Apply the CS decomposition to $\bQ$, there exist orthogonal matrices $\bU_1\in \real^{m_1\times m_1}$, $\bU_2\in \real^{m_2\times m_2}$, and $\bV_1 \in \real^{r\times r}$ such that 
$$
\begin{bmatrix}
\bU_1 & \bzero \\
\bzero & \bU_2
\end{bmatrix}^\top 
\begin{bmatrix}
\bQ_{11}\\
\bQ_{21} 
\end{bmatrix}\bV_1 
=
\begin{bmatrix}
\bD_A^r \\
\bD_B^r
\end{bmatrix},
$$
where $\bD_A^r$ and $\bD_B^r$ are the first $r$ columns of $\bD_A$ and $\bD_B$, respectively (note that row permutation needed here). It then follows that 
$$
\scriptsize
\begin{aligned}
\begin{bmatrix}
\bU_1 & \bzero \\
\bzero & \bU_2
\end{bmatrix}^\top 
\begin{bmatrix}
\bA\\
\bB 
\end{bmatrix}\bP
&= 
\begin{bmatrix}
\bD_A^r & \bU_1\bQ_{12}\\
\bD_B^r & \bU_2\bQ_{22}
\end{bmatrix}
\begin{bmatrix}
\bV_1^\top \bSigma_r & \bzero \\
\bzero  & \bzero 
\end{bmatrix} 
= \begin{bmatrix}
\bD_A^r & \bzero\\
\bD_B^r & \bzero
\end{bmatrix}
\begin{bmatrix}
\bV_1^\top \bSigma_r & \bzero \\
\bzero  & \bI_{n_1-r} 
\end{bmatrix} 
= \begin{bmatrix}
\bD_A\\
\bD_B
\end{bmatrix}
\begin{bmatrix}
\bV_1^\top \bSigma_r & \bzero \\
\bzero  & \bI_{n_1-r} 
\end{bmatrix}. \\
\end{aligned}
$$
By setting 
$$
\bX = \bP \begin{bmatrix}
\bV_1^\top \bSigma_r & \bzero \\
\bzero  & \bI_{n_1-r} 
\end{bmatrix}^{-1},
$$
we complete the proof.
\end{proof}
Note that if $\bB =\bI_{n_1}$ and we set $\bX = \bU_2$, then we obtain the SVD of $\bA$.

%% file: chapter-app_inequality.tex
\newpage
\chapter{Selected Topics}
\begingroup
\hypersetup{
linkcolor=structurecolor,
linktoc=page,  
}
\minitoc \newpage
\endgroup

\newpage


\section{Famous Inequalities}\label{appendix:inequalities}
In this section, we present several well-known inequalities that will be frequently utilized.
Additional inequalities are discussed throughout the main text; for example, the Schur inequality in Theorem~\ref{theorem:schur_inequality}, Hadamard inequalities in \S~\ref{section:hadama_ineq}, Sylvester's inequality and Frobenius' inequality in Remark~\ref{remark:rank_prop}, Fan's inequality in Problem~\ref{prob:fans_ineq}, and Hardy-Littlewood-P\'olya inequality in Problem~\ref{prob:hardy_ineq}.

We  state Jensen's inequality for a convex function without providing a proof;  for a detailed discussion, refer to  \citet{beck2014introduction}.

\begin{theorem}[Jensen's Inequality]
Let $f: \sS \rightarrow \real$ be  a convex function, where $\sS \subseteq \real^n$ is a convex set. Then for any $\bx_1, \bx_2, \ldots, \bx_k\in\sS$ and $\blambda \in \Delta_k$ (Definition~\ref{definition:simplex}),  the following inequality holds:
\begin{equation*}
f\left( \sum_{i=1}^{k} \lambda_i \bx_i\right) \leq \sum_{i=1}^{k} \lambda_i f(\bx_i).
\end{equation*}
\end{theorem}

The AM-GM inequality  is a fundamental tool in competitive mathematics, particularly useful for determining the maximum or minimum values of multivariable functions or expressions.  It establishes a relationship between the \textit{arithmetic mean (AM)} to the \textit{geometric mean (GM)}. 
\begin{theorem}[AM-GM Inequality]
For  any nonnegative real numbers $x_1,x_2,\ldots,x_n$,  the following inequality holds:
$$
\frac{\sum_{i=1}^{n} x_i}{n} \geq \sqrt[n]{\prod_{i=1}^{n}x_i}.
$$ 
That is, the \textit{geometric mean} of a set of nonnegative numbers does not exceed the \textit{arithmetic mean}, and equality holds if and only if all the numbers are equal.
\end{theorem}

When $n=2$, the AM-GM inequality can be expressed as:
\begin{equation}\label{equation:amgm_ineq}
a^2+b^2\geq 2\sqrt{a^2b^2} = 2\abs{ab}\geq ab.
\end{equation}

\begin{proposition}[Weighted AM-GM Inequality]\label{proposition:weighted_amgm}
For any nonnegative real numbers $x_1,x_2,\ldots,x_n$ and nonnegative weights $w_1,w_2,\ldots, w_n$,  the following inequality holds:
$$
\frac{\sum_{i=1}^{n} w_i x_i}{\sum_{i=1}^{n}w_i} \geq
\sqrt[(\sum w_i)]{\prod_{i=1}^{n}x_i^{w_i}}.
$$ 
When $w_1=w_2=\ldots=w_n=1$, this reduces to  the standard AM-GM inequality. Alternatively, when $\bw\in \Delta_k$, i.e., $\bw$ is a unit-simplex with $\sum_{i=1}^{n}w_i=1$, it follows that 
$$
\sum_{i=1}^{n} w_i x_i \geq
\prod_{i=1}^{n}x_i^{w_i},
\quad\text{where}\quad \bw\in \Delta_k.
$$
\end{proposition}
\begin{proof}[of Proposition~\ref{proposition:weighted_amgm}]
For simplicity, we will only prove the second part, the first one can be established similarly.
Applying  Jensen's inequality to the convex function $f(x)=-\ln (x)$ with $x_1, x_2,  \ldots, x_n>0$ and $\bw \in \Delta_n$, we have 
$
-\ln\left( \sum_{i=1}^{k} w_i x_i\right) \leq -\sum_{i=1}^{k} w_i \ln(\bx_i).
$
Taking the exponent of both sides, we get 
$
\sum_{i=1}^{n}w_i x_i \geq \exp^{\sum_{i=1}^{n} w_i \ln(x_i)}
=\prod_{i=1}^{n}x_i^{w_i}.
$
This completes the proof.
\end{proof}

\begin{exercise}
Given nonnegative real numbers $x_1, x_2, \ldots, x_n\in \real_+$, show that 
$$
\frac{\sum_{i=1}^{n }x_i}{n} \leq \sqrt{\frac{\sum_{i=1}^{n}x_i^2}{n}}.
$$
\end{exercise}

\begin{exercise}
Given positive real numbers $x_1, x_2, \ldots, x_n\in \real_{++}$, show that 
$$
\frac{\sum_{i=1}^{n} x_i^2}{\sum_{i=1}^{n} x_i} 
\leq \sqrt{\frac{\sum_{i=1}^{n} x_i^3}{\sum_{i=1}^{n} x_i}}.
$$
\end{exercise}

\begin{exercise}
Given $\bw\in \Delta_n$, show that 
$$
\sum_{i=1}^{n} \frac{x_i}{\sqrt{1-x_i}} \geq \sqrt{\frac{n}{n-1}}.
$$
\end{exercise}

\begin{theorem}[Markov's Inequality]\label{theorem:markov-inequality}
Let $\rx$ be a nonnegative random variable. Then, for any $\epsilon >0$, we have 
$$
\prob[\rx \geq \epsilon] \leq \frac{\Exp[\rx]}{\epsilon}.
$$
\end{theorem}
\begin{proof}[of Theorem~\ref{theorem:markov-inequality}]
We notice the trick that $0\leq \epsilon \indicator\{\rx \geq \epsilon\} \leq \rx$ since $\rx$ is nonnegative. 
Taking expectations on both sides, we get $\Exp[\epsilon \indicator\{\rx\geq \epsilon\}] \leq \Exp[\rx]$. This simplifies to:
$$
\Exp\left[\epsilon \indicator\{\rx\geq \epsilon\}\right] = \epsilon\Exp\left[ \indicator\{\rx\geq \epsilon\}\right] = \epsilon \left(1\cdot \prob[\rx\geq \epsilon] + 0\cdot \prob[\rx< \epsilon] \right) = \epsilon \cdot \prob[\rx\geq \epsilon] \leq \Exp[\rx].
$$
This completes the proof.
\end{proof}

\begin{theorem}[Chebyshev's Inequality]
Let $\rx$ be a random variable with finite mean $\Exp[\rx]< \infty$. Then, for any $\epsilon >0$, we have 
$$
\prob[|\rx - \Exp[\rx]| \geq \epsilon ] \leq \frac{\Var[\rx]}{\epsilon^2}.
$$
\end{theorem}
Chebyshev's inequality can be readily verified by defining $\ry \triangleq (\rx - \Exp[\rx])^2$ (which is nonnegative) and applying Markov's inequality to $\ry$.

\subsection{Cauchy-Schwarz Inequality}
\index{Cauchy–Schwarz inequality}

The \textit{Cauchy–Schwarz inequality}  is considered one of the most important and widely used inequalities in mathematics.
\begin{theorem}[Cauchy-Schwarz Inequality]
For any random $m\times n$ matrices $\rmX$ and $\rmY$, we have 
$$
\Exp\left[\norm{\rmX^\top \rmY}\right] \leq \Exp\left[\norm{\rmX}^2\right]^{1/2} \Exp\left[\norm{\rmY}^2\right]^{1/2},
$$
where the inner product is defined as $\langle\rmX,\rmY\rangle = \Exp\left[\norm{\rmX^\top \rmY}\right]$.
\end{theorem}

The result can be applied to non-random matrices and vectors.
\begin{proposition}[Cauchy-Schwarz Matrix (Vector) Inequality]\label{proposition:cauchy-schwarz-inequ}
For any $m\times n$ matrices $\bX$ and $\bY$, we have 
$$
\norm{\bX^\top \bY} \leq \Vert\bX\Vert \cdot \Vert\bY\Vert.
$$
This is a special form of the Cauchy-Schwarz inequality, where the inner product is defined as $\langle\bX,\bY\rangle = \norm{\bX^\top \bY}$.
Similarly, for any vectors $\bu$ and $\bv$, we have 
\begin{equation}\label{equation:vector_form_cauchyschwarz}
\abs{\bu^\top \bv} \leq \normtwo{\bu} \cdot  \normtwo{\bv},
\end{equation}
where the equality holds if and only if $\bu$ and $\bv$ are linearly dependent.
In two-dimensional case, it becomes
$$
(ac+bd)^2 \leq (a^2 +b^2)(c^2+d^2).
$$
\end{proposition}
The vector form of the Cauchy-Schwarz inequality plays an important role in different branches of modern mathematics, including Hilbert space theory and numerical analysis \citep{wu2009various}. Here, we  only provide the proof for the vector form of the Cauchy-Schwarz inequality for simplicity. To see this, given two vectors $\bu,\bv\in \real^n$, we have 
$$
\begin{aligned}
&0\leq \sum_{i=1}^{n}\sum_{j=1}^{n} (u_i v_j - u_j v_i)^2 = 
\sum_{i=1}^{n}\sum_{j=1}^{n} u_i^2v_j^2 + \sum_{i=1}^{n}\sum_{j=1}^{n} v_i^2 u_j^2 - 2\sum_{i=1}^{n}\sum_{j=1}^{n} u_iu_j v_iv_j\\
&=\left(\sum_{i=1}^{n} u_i^2\right) \left(\sum_{j=1}^{n} v_j^2\right) +
\left(\sum_{i=1}^{n} v_i^2\right) \left(\sum_{j=1}^{n} u_j^2\right) - 
2\left(\sum_{i=1}^{n} u_iv_i \right)^2
=2 \normtwo{\bu}^2 \cdot \normtwo{\bv}^2 -2 \abs{\bu^\top\bv}^2,
\end{aligned}
$$
from which the result follows.
The equality holds if and only if $\bu = k\bv$ for some constant $k\in \real$, i.e., $\bu$ and $\bv$ are linearly dependent.

The Cauchy-Schwarz inequality can be applied to prove  \textit{Titu's lemma}.
\begin{proposition}[Titu's Lemma]\label{proposition:titus_lemma}
Let $x_1, x_2, \ldots, x_n$ and $y_1, y_2, \ldots, y_n$ be positive reals. Then,
$$
\frac{x_1^2}{y_1}+\frac{x_2^2}{y_2}+\ldots + \frac{x_n^2}{y_n}
\geq \frac{(x_1+x_2+\ldots +x_n)^2}{ y_1+y_2+\ldots+y_n}.
$$
\end{proposition}
\begin{proof}[of Proposition~\ref{proposition:titus_lemma}]
Suppose $a_i$ = $\frac{x_i}{\sqrt{y_i}}$ and $b_i = \sqrt{y_i}$. By the Cauchy-Schwarz inequality, we have 
$$
\begin{aligned}
|\ba^\top\bb|^2 &\leq \Vert\ba\Vert^2 \cdot  \Vert\bb\Vert^2
\quad\implies\quad 
\left(\sum_{i=1}^{n} x_i\right)^2	\leq \left( \sum_{i=1}^{n} \frac{x_i^2}{y_i} \right)
\left( \sum_{i=1}^{n} y_i\right).
\end{aligned}
$$
This completes the proof.
\end{proof}

\begin{exercise}[Nesbitt's Inequality]
Let $x,y,z\in \real^+$. Show that 
$$
\frac{x}{y+z} + \frac{y}{x+z} + \frac{z}{x+y} \geq \frac{3}{2}.
$$
\textit{Hint: Apply Titu's lemma to $
\frac{x^2}{xy+xz} + \frac{y^2}{xy+yz} + \frac{z^2}{xz+yz}.
$}
\end{exercise}

\paragraph{Angle between two vectors.} From Equation~\eqref{equation:vector_form_cauchyschwarz}, given two vectors $\bx,\by$, we note that 
$$
-1\leq 
\frac{\bx^\top\by}{\Vert\bx\Vert_2 \Vert\by\Vert_2} 
\leq 1.
$$
This two-side inequality illustrates the concept of the angle between two vectors.
\begin{definition}[Angle Between Vectors]\label{definition:angle_bet_vec_ineq}
The angle between two vectors $\bx$ and $\by$ is the number $\theta\in [0,\pi]$ such that 
$$
\cos \theta = \frac{\bx^\top\by}{\Vert\bx\Vert_2 \Vert\by\Vert_2} .
$$
\end{definition}
The definition of the angle between vectors will be useful in the discussion of transformation in matrix decompositions (Section~\ref{section:coordinate-transformation}).

\subsection{Generalized Cauchy-Schwarz}
From the vector form of the Cauchy-Schwarz inequality, let $x_i\triangleq u_i^2$, $y_i\triangleq v_i^2$ for all $i=\{1,2,\ldots,n\}$, where $\bu,\bv\in\real^n$ in Equation~\eqref{equation:vector_form_cauchyschwarz}, we have 
\begin{equation}\label{equation:gen_ca_sc_1}
\left( \sum_{i=1}^{n} x_i \right)^{1/2} \left( \sum_{i=1}^{n} y_i \right)^{1/2}
\geq  \sum_{i=1}^{n} (x_iy_i)^{1/2} .
\end{equation}
More generally, we have the generalized Cauchy-Schwarz inequality.
\begin{theorem}[Generalized Cauchy-Schwarz Inequality]
Given a set of vectors $\ba,\bb,\ldots, \bz\in \real^n$, and weights $\lambda_1, \lambda_2, \ldots, \lambda_z$ with $\lambda_1+\lambda_2+\ldots+\lambda_z=1$, it follows that 
$$
\begin{aligned}
& \left( \sum_{i=1}^{n} a_i \right)^{\lambda_a}
\left( \sum_{i=1}^{n} b_i \right)^{\lambda_b}
\ldots
\left( \sum_{i=1}^{n} z_i \right)^{\lambda_z}\\
&\gap
\geq 
(a_1^{\lambda_a} b_1^{\lambda_b} \ldots z_1^{\lambda_z})+
(a_2^{\lambda_a} b_2^{\lambda_b} \ldots z_2^{\lambda_z})+
\ldots +
(a_n^{\lambda_a} b_n^{\lambda_b} \ldots z_n^{\lambda_z})
\end{aligned}.
$$
The equality holds if $a_i=b_i=\ldots=z_i$ for all $i\in\{1,2,\ldots,n\}$.
\end{theorem}

\index{Schur's inequality}
\subsection{Schur's Inequality}
 Schur's inequality relates three nonnegative real numbers.
\begin{theorem}[Schur's Inequality]\label{theorem:schurs_real_inequality}
Given nonnegative real numbers $x, y, z$, and a positive real number $t$, it follows that 
\begin{equation}\label{equation:sschurs_real_inequality}
x^t(x-y)(x-z) + y^t(y-x)(y-z) +z^t(z-x)(z-y) \geq 0.
\end{equation}
The equality holds if and only if $x=y=z$, or if two of $x,y,z$ are equal and the third is 0.
When $t=1$, it becomes 
$$
x^3+y^3+z^3 + 3xyz \geq xy(x+y) + xz(x+z) + yz(y+z).
$$
\end{theorem}
\begin{proof}[of Theorem~\ref{theorem:schurs_real_inequality}]
We note the inequality is symmetric in the variables $x,y,z$. Without loss of generality, suppose $x\geq y\geq z$. We factor the left-hand side of Equation~\eqref{equation:sschurs_real_inequality} as 
$
(x-y)\bigg(x^t (x-z) - y^t(y-z) \bigg) + z^t (z-x)(z-y).
$
We note that $(x-y)\geq 0 $, $\big(x^t (x-z) - y^t(y-z) \big) \geq 0$, and $z^t (z-x)(z-y)\geq 0$. This completes the proof.
\end{proof}

\subsection{Young's Inequality}
Young's inequality is a special case of the weighted AM-GM inequality (Proposition~\ref{proposition:weighted_amgm}), while it has a wider range of applications.
\begin{theorem}[Young’s Inequality]\label{theorem:holder-inequality1}
For nonnegative numbers $x,y \geq 0$, and positive real numbers $p,q>1$ with $\frac{1}{p}+\frac{1}{q} = 1$, it follows that
\begin{equation}\label{equation:holder-v1}
xy \leq  \frac{1}{p} x^p +\frac{1}{q} y^q.
\end{equation}
\end{theorem}

\begin{figure}[h]
\centering
\vspace{-0.35cm}
\subfigtopskip=2pt
\subfigbottomskip=2pt
\subfigcapskip=-5pt
\subfigure[Case 1: $x< y^{1/(p-1)}$ and $p\geq 2$.]{\label{fig:holder1}
\includegraphics[width=0.4\linewidth]{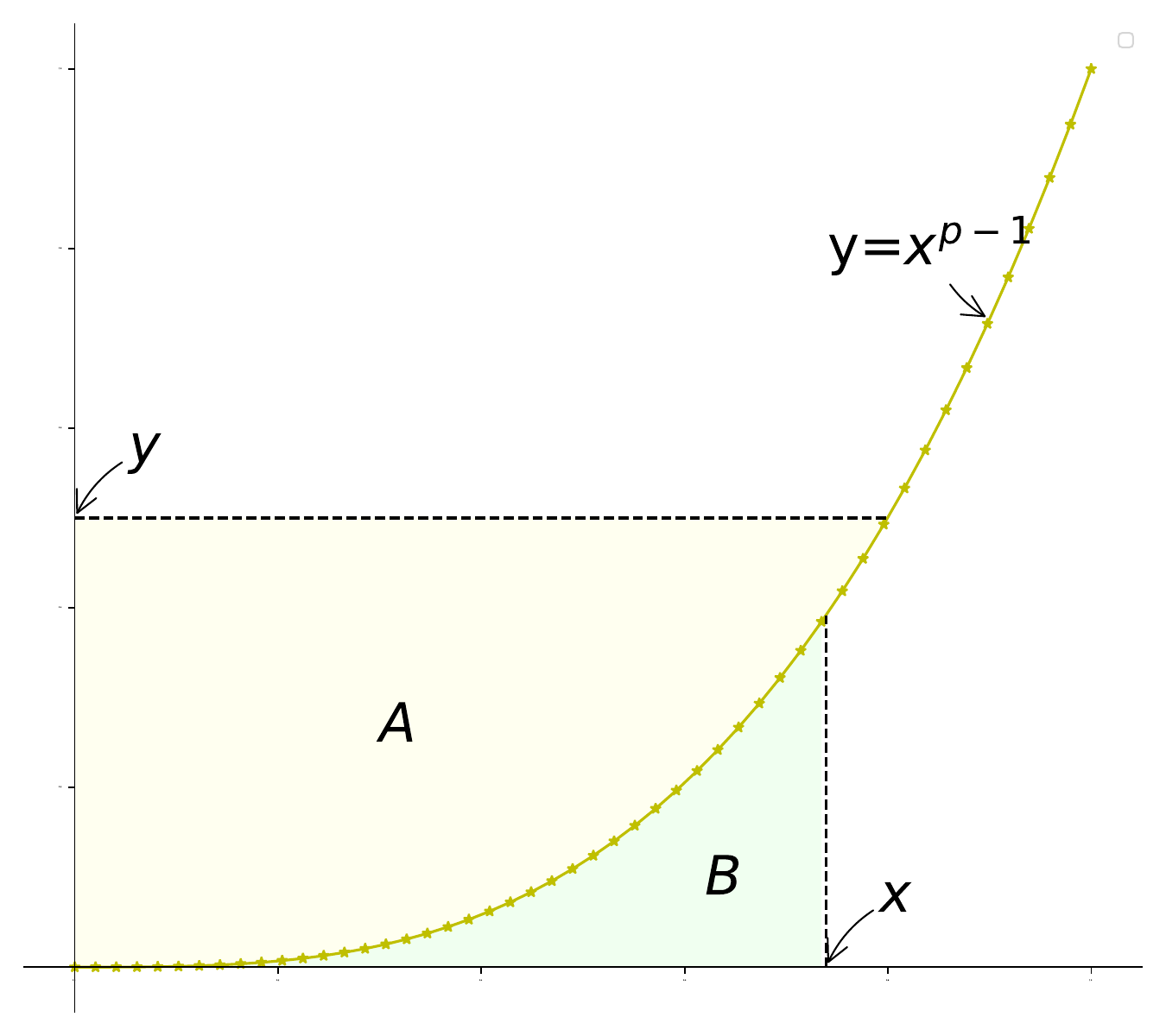}}
\quad 
\subfigure[Case 2: $x\geq  y^{1/(p-1)}$ and $p\geq 2$.]{\label{fig:holder2}
\includegraphics[width=0.4\linewidth]{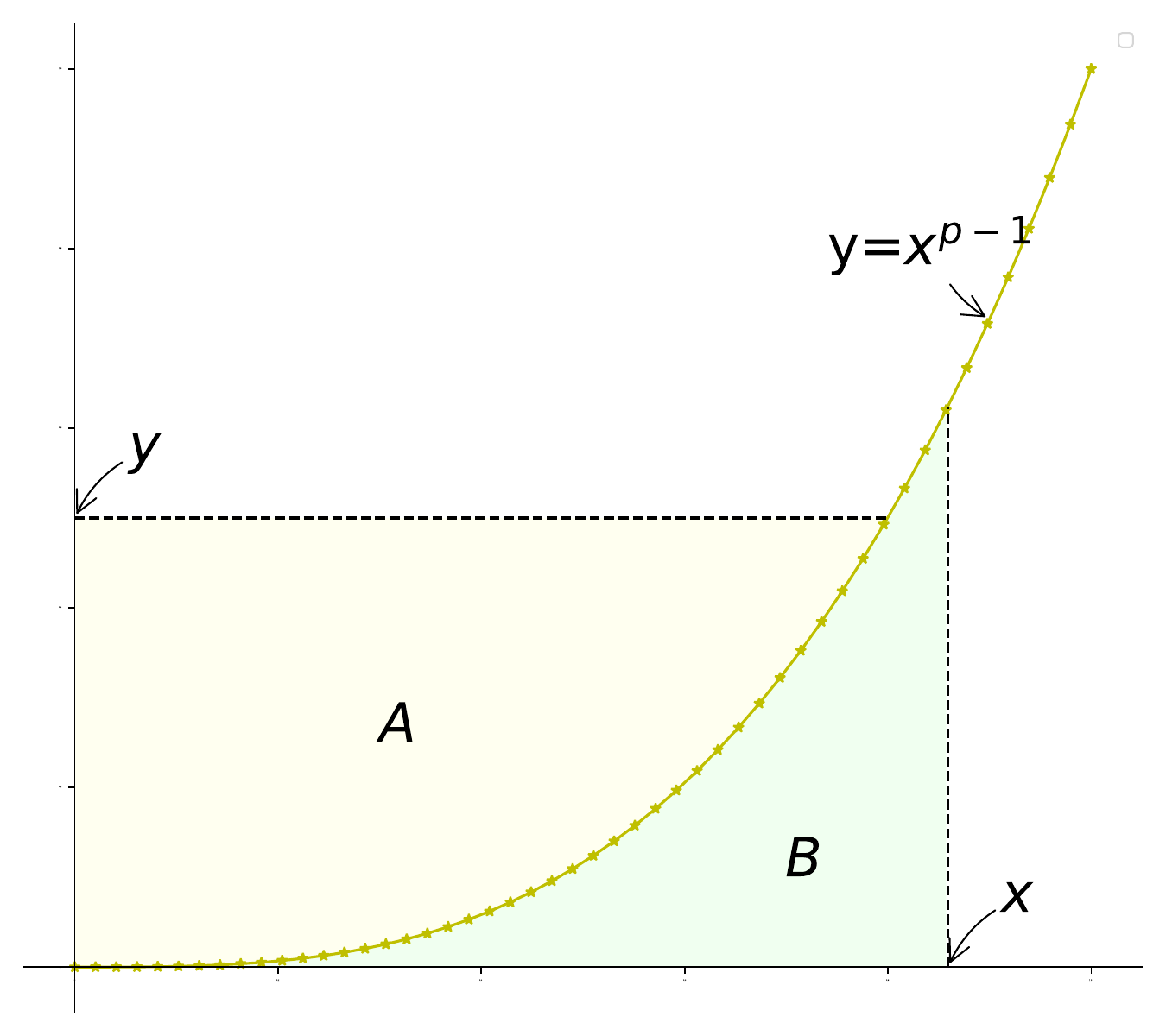}}
\subfigure[Case 3: $x< y^{1/(p-1)}$ and $1<p< 2$.]{\label{fig:holder3}
\includegraphics[width=0.4\linewidth]{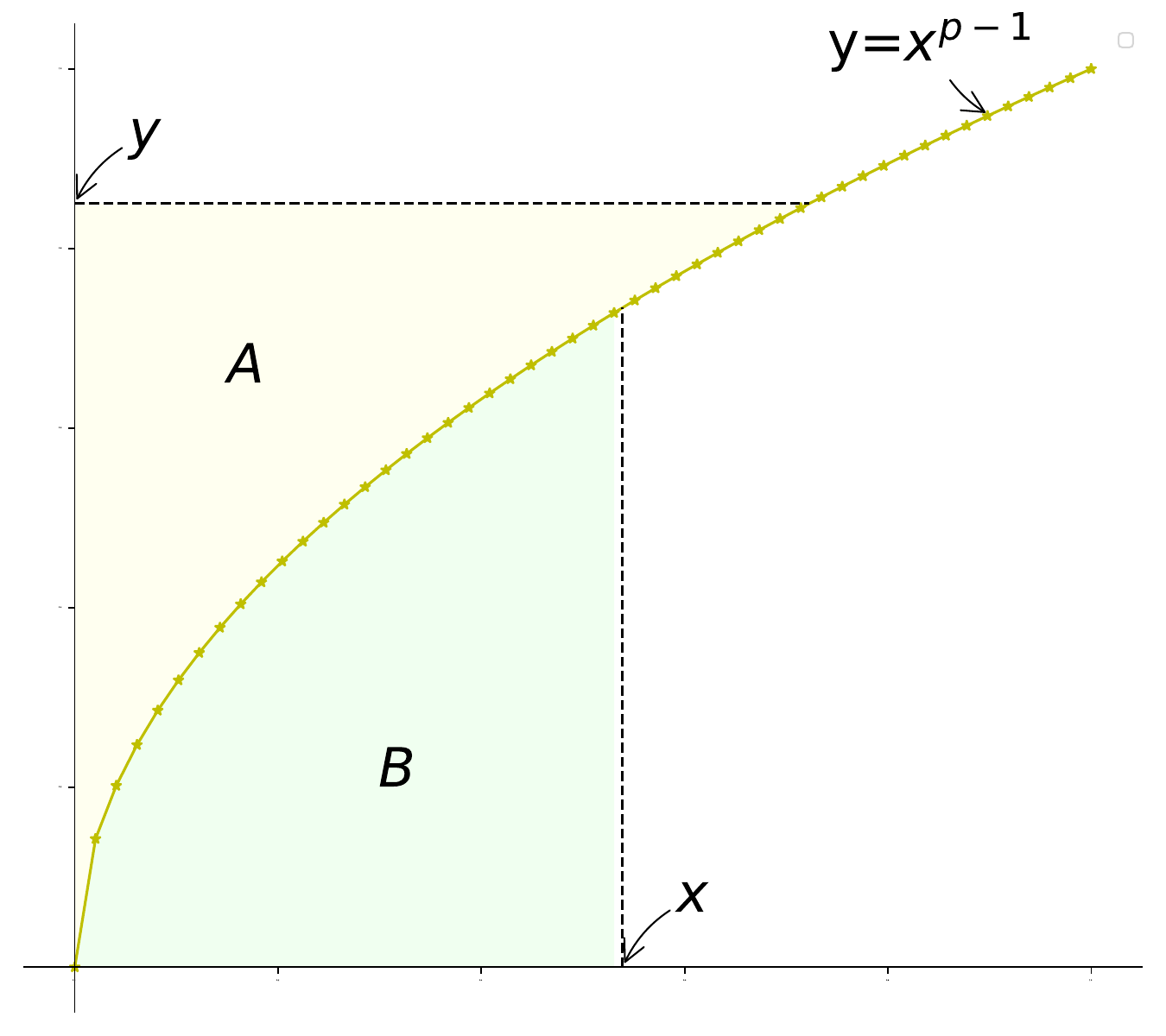}}
\quad 
\subfigure[Case 4: $x\geq  y^{1/(p-1)}$ and $1<p< 2$.]{\label{fig:holder4}
\includegraphics[width=0.4\linewidth]{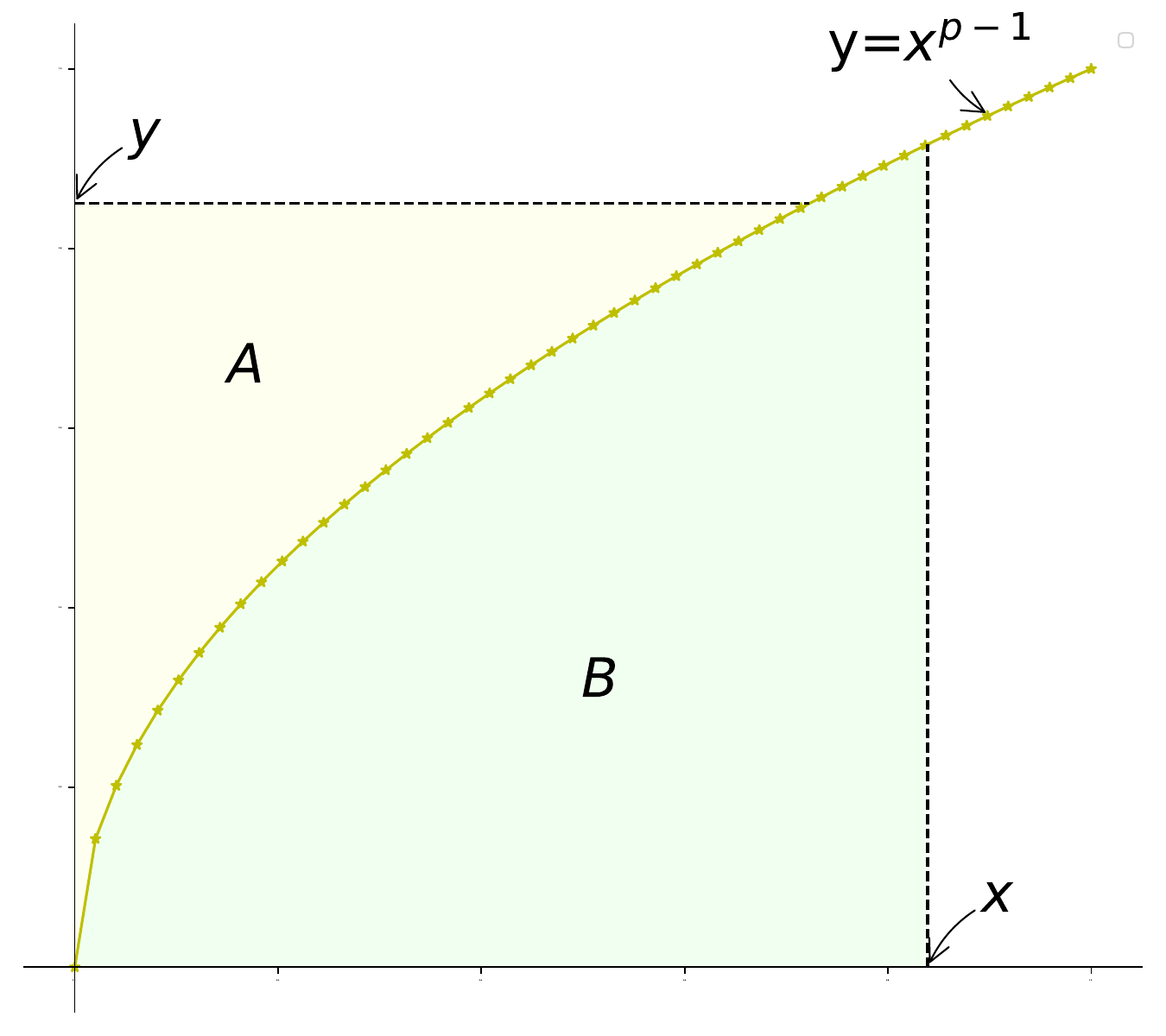}}
\caption{Demonstration of Young’s inequality for different cases.}
\label{fig:holderineqd}
\end{figure}
\begin{proof}[of Theorem~\ref{theorem:holder-inequality1}]
The area $xy$ is smaller than the areas of the sum of the two trapezoids with curved edges (the colored areas $A$ and $B$) as shown in Figure~\ref{fig:holderineqd}:
$$
\mathrm{area } \,\,\,B =\int_0^x x^{p-1} dx, \qquad \mathrm{area } \,\,\,A=\int_{0}^y y^{1/(p-1)}dy.
$$
That is, 
$$
\begin{aligned}
xy &\leq \int_0^x x^{p-1} dx + \int_{0}^y y^{1/(p-1)}dy 
= \frac{1}{p} x^p + \int_{0}^y y^{q/p}dy \\
&= \frac{1}{p} x^p +  (\frac{q}{p}+1)^{-1} y^{q/p+1} 
= \frac{1}{p} x^p +\frac{1}{q} y^q.
\end{aligned}
$$
This completes the proof.
\end{proof}

\begin{SCfigure}
\centering
\includegraphics[width=0.53\textwidth]{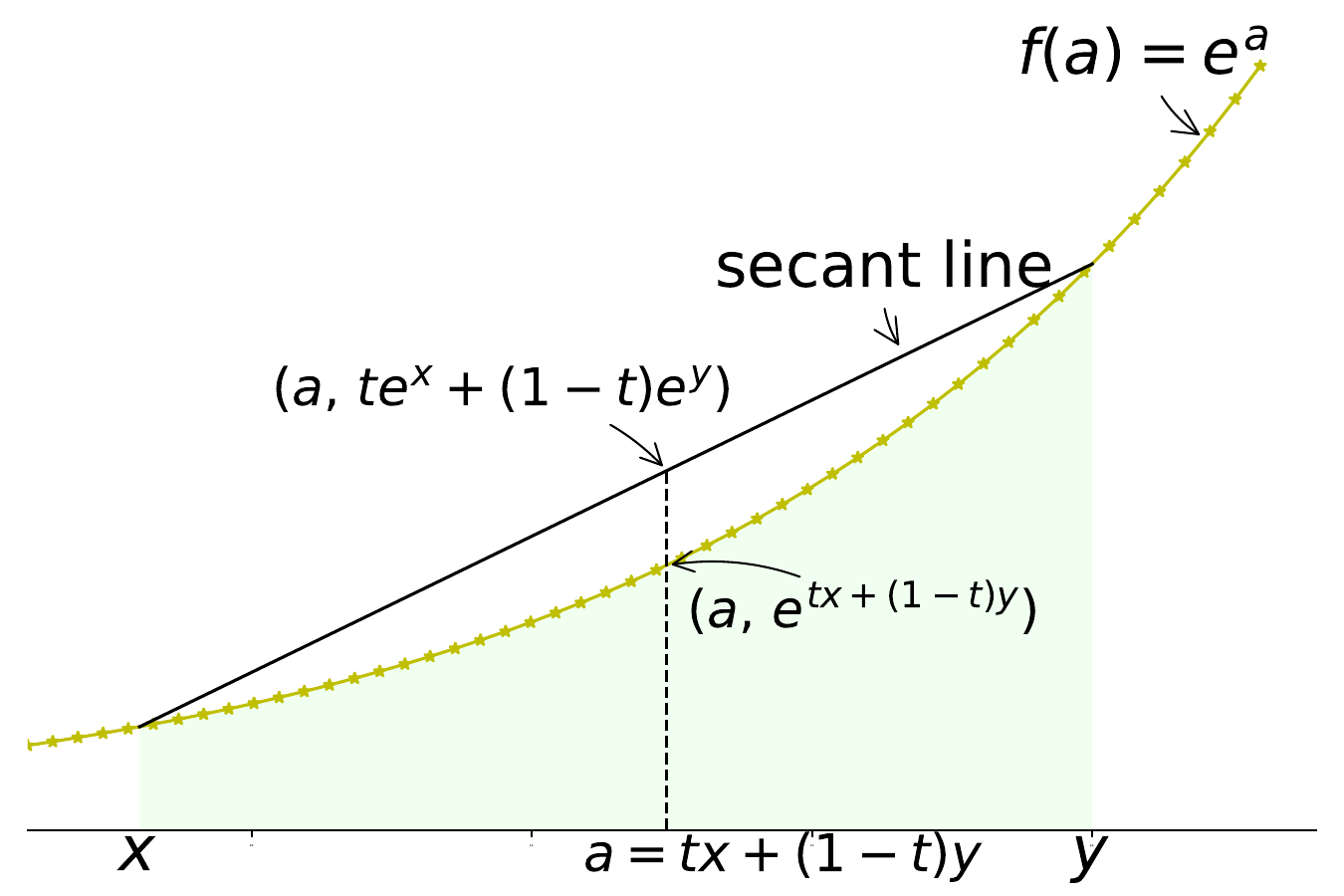}
\caption{Demonstration of interpolation inequality for $e^x$.}
\label{fig:older_interpolation}
\end{SCfigure}

There are several conceptually different ways to prove Young's inequality. To provide an alternative, we can use the \textit{interpolation inequality} as follows.
\begin{lemma}[Interpolation Inequality for Exponential function]\label{lemma:interpolation_inequality}
Given $t\in [0,1]$, we have 
$$
e^{tx+(1-t)y} \leq t e^x + (1-t)e^y. 
$$
\end{lemma}
\begin{proof}[of Lemma~\ref{lemma:interpolation_inequality}]
The function of the scant line through the points $(x,e^x)$ and $(y, e^y)$ on the graph of $f(a)=e^a$ is (see Figure~\ref{fig:older_interpolation}):
$$
f(t) = (tx+(1-t)y, te^x+(1-t)e^y).
$$
Since the function $f(a)=e^a$ is convex, we have 
$$
e^{tx+(1-t)y}\leq t e^x + (1-t)e^y. 
$$
This completes the proof.
\end{proof}


By using the interpolation inequality for the exponential function, we provide an alternative proof for  Theorem~\ref{theorem:holder-inequality1}.
\begin{proof}[of Theorem~\ref{theorem:holder-inequality1}, an alternative way]
We observe that
$$
\begin{aligned}
xy &= \exp\{\log x+\log y\} 
= \exp\left \{ \frac{q-1}{q} \frac{q}{q-1} \log x + \frac{1}{q} q \log y \right\}.
\end{aligned}
$$
Then, by Lemma~\ref{lemma:interpolation_inequality}, we have 
$$
\begin{aligned}
xy&\leq \frac{q-1}{q} e^{\frac{q}{q-1} \log x} + \frac{1}{q} e^{q\log y} 
= \frac{q-1}{q} x^{\frac{q}{q-1}} + \frac{1}{q} y^q
= \frac{1}{p} x^p +\frac{1}{q}y^q,
\end{aligned}
$$
from which the result follows.
\end{proof}


\subsection{H{\"o}lder's Inequality}

\holders inequality, named after \textit{Otto  H{\"o}lder}, is widely used in optimization, machine learning, and many other topics. 
It is also a generalization of the (vector) Cauchy-Schwarz inequality.

\begin{theorem}[\holders  Inequality]\label{theorem:holder-inequality}
Suppose $p,q>1$ such that $\frac{1}{p}+\frac{1}{q} = 1$. Then  for any vector $\bx,\by\in \real^n$, we have
$$
\sum_{i=1}^{n}x_i y_i
\leq 
\abs{\sum_{i=1}^{n}x_i y_i}
\leq \sum_{i=1}^{n}|x_i| |y_i| \leq \left(\sum_{i=1}^{n}  |x_i|^p\right)^{1/p}  \left(\sum_{i=1}^{n} |y_i|^q\right)^{1/q}=\norm{\bx}_p\norm{\by}_q,
$$
where $\norm{\bx}_p = \left(\sum_{i=1}^{n}  |x_i|^p\right)^{1/p} $ is known as the \textbf{$\ell_p$ norm} or \textbf{$p$ norm} of the vector $\bx$ (see Section~\ref{section:vector-norm}). The equality holds if the two sequences $\{\abs{x_i}^p\}$ and $\{|y_i|^q\}$ are linearly dependent~\footnote{To be more concrete, the equality attains if and only if $\abs{\bx^\top\by}=\abs{\bx}^\top\abs{\by}$ and 
$$
\left\{
\begin{aligned}
&\abs{\bx}\hadaprod\abs{\by}=\norminf{\by}\abs{\bx}, &\text{if }& p=1; \\ 
&\abs{\bx}\hadaprod\abs{\by}=\norminf{\bx}\abs{\by}, &\text{if }& p=\infty;\\ 
&\norm{\by}_q^{1/p} \abs{\bx}^{\hadaprod1/q} = \norm{\bx}_p^{1/q}\abs{\by}^{\hadaprod1/p}, &\text{if }&1<p<\infty, \\ 
\end{aligned}
\right.
$$
where $\hadaprod$ denotes Hadamard product/power (Definition~\ref{definition:hada_prod}); see \citet{bernstein2008matrix} and the references therein.
}.
When $p=q=2$, this reduces to the vector Cauchy-Schwarz inequality \eqref{equation:vector_form_cauchyschwarz}.
\end{theorem}
\begin{proof}[of Theorem~\ref{theorem:holder-inequality}]
Let $ u \triangleq \frac{|x_i|}{\Vert\bx\Vert_p}$ and $ v \triangleq \frac{|y_i|}{\Vert\by\Vert_q}$. From Equation~\eqref{equation:holder-v1}, it follows that
$$
uv = \frac{|x_i||y_i|}{\norm{\bx}_p\norm{\by}_q} \leq \frac{1}{p} \frac{|x_i|^p}{\norm{\bx}_p^p} + \frac{1}{q}\frac{|y_i|^q}{\norm{\by}_q^q}, \qquad \forall i \in \{1,2,\ldots,n\}.
$$
Therefore,
$$
\sum_{i=1}^{n} \frac{|x_i||y_i|}{\norm{\bx}_p\norm{\by}_q} \leq \frac{1}{p\norm{\bx}_p^p} \sum_{i=1}^{n}|x_i|^p  +
\frac{1}{q\norm{\by}_q^q}   \sum_{i=1}^{n}|y_i|^q = \frac{1}{p}+\frac{1}{q} = 1.
$$
That is, 
$
\sum_{i=1}^{n} |x_i||y_i|  \leq  \norm{\bx}_p\norm{\by}_q.
$
It is trivial that $\sum_{i=1}^{n}x_i y_i\leq \sum_{i=1}^{n}|x_i| |y_i| $,
from which the result follows.
\end{proof}

\begin{corollary}[Cauchy Inequality]
When take $p=q=2$ in \holders inequality, we recover the \textit{Cauchy inequality}:
$$
\sum_{i=1}^{n} x_i y_i \leq
\sum_{i=1}^{n} |x_i| |y_i| \leq 
\sqrt{\sum_{i=1}^{n} |x_i|^2}
\sqrt{\sum_{i=1}^{n} |y_i|^2}.
$$
\end{corollary}

Minkowski's inequality follows immediately from the \holders inequality.
\begin{theorem}[Minkowski's Inequality]\label{theorem:minkoski}
Given nonnegative reals $x_1, x_2, \ldots, x_n\geq 0$ and $y_1, y_2, \ldots, y_n\geq 0$, and $p\geq 1$, we have 
$$
\left( \sum_{i=1}^{n} (x_i+y_i)^p  \right)^{1/p} \leq  
\left( \sum_{i=1}^{n}x_i^p \right)^{1/p}+
\left( \sum_{i=1}^{n}y_i^p \right)^{1/p}.
$$
The equality holds if and only if the sequence $x_1, x_2, \ldots, x_n$ and $y_1, y_2, \ldots, y_n$ are proportional. If $p<1$, the inequality sign is reversed.
\end{theorem}
Minkowski's inequality is essentially  the triangle inequality of $\ell_p$ norms (Section~\ref{section:vector-norm}).
\begin{proof}[of Theorem~\ref{theorem:minkoski}]
If $p=1$, the inequality is immediate. We then assume $p>1$.
We note that 
$$
\sum_{i=1}^{n} (x_i+y_i)^p = 
\sum_{i=1}^{n} x_i (x_i+y_i)^{p-1}
+ 
\sum_{i=1}^{n} y_i(x_i+y_i)^{p-1}.
$$
Using \holders inequality, given $p,q$ such that $p>1$ and $\frac{1}{p}+\frac{1}{q}=1$, we have 
$$
\footnotesize
\begin{aligned}
\sum_{i=1}^{n} x_i (x_i+y_i)^{p-1}
&\leq 
\bigg( \sum_{i=1}^{n} x_i^p \bigg)^{1/p}
\bigg( \sum_{i=1}^{n} (x_i+y_i)^{(p-1)q} \bigg)^{1/q}
=\bigg( \sum_{i=1}^{n} x_i^p \bigg)^{1/p}
\bigg( \sum_{i=1}^{n} (x_i+y_i)^{p} \bigg)^{1/q},\\
\sum_{i=1}^{n} y_i (x_i+y_i)^{p-1}
&\leq 
\bigg( \sum_{i=1}^{n} y_i^p \bigg)^{1/p}
\bigg( \sum_{i=1}^{n} (x_i+y_i)^{(p-1)q} \bigg)^{1/q}
=\bigg( \sum_{i=1}^{n} y_i^p \bigg)^{1/p}
\bigg( \sum_{i=1}^{n} (x_i+y_i)^{p} \bigg)^{1/q}.
\end{aligned}
$$
Combining these results, we get
$$
\sum_{i=1}^{n} (x_i+y_i)^p
\leq 
\bigg(
\bigg( \sum_{i=1}^{n} x_i^p \bigg)^{1/p}
+
\bigg( \sum_{i=1}^{n} y_i^p \bigg)^{1/p}
\bigg)
\bigg( \sum_{i=1}^{n} (x_i+y_i)^{p} \bigg)^{1/q}.
$$
We conclude the claimed result.
\end{proof}

%% file: chapter-app_norm.tex
\newpage
\section{Norm}\label{appendix:matrix-norm}
The concept of a norm is essential for evaluating the magnitude of vectors and, consequently, allows for the definition of certain metrics on linear spaces equipped with norms. Norms provide a measure of the magnitude of a vector or matrix, which is useful in many applications, such as determining the length of a vector in Euclidean space or the size of a matrix in a multidimensional setting.
Additionally, norms enable us to define distances between vectors or matrices. 
The distance between two vectors $\bu$ and $\bv$ can be computed using the norm of their difference $\norm{\bu-\bv}$. This is critical for tasks involving proximity measures, such as clustering algorithms in machine learning.~\footnote{We only discuss the norms (and inner products) for real vector or matrix spaces. Most results can be applied directly to complex cases.}
\begin{definition}[Vector Norm and Matrix Nrom\index{Matrix norm}\index{Vector norm}]\label{definition:matrix-norm}
Given a norm $\norm{\cdot}: \real^n\rightarrow \real$ on vectors or a norm $\norm{\cdot}: \real^{m\times n}\rightarrow \real$ on matrices, for any vector $\bx \in \real^{n}$ and  any matrix $\bA \in \real^{m\times n}$, we have~\footnote{When $\norm{\bx}=\bzero$ for some nonzero vector $\bx$, the norm is called a \textit{semi-norm}.}~\footnote{When the vector has a single element, the norm can be understood as the absolute value operation.}~\footnote{When $\norm{\bI}=1$ for a matrix norm, the matrix norm is said to be \textit{normalized}.}
\begin{enumerate}
\item \textit{Positive homogeneity}. $\norm{\lambda \bA} = |\lambda| \cdot \norm{\bA}$ or $\norm{\lambda \bx} = |\lambda| \cdot \norm{\bx}$ for any $\lambda \in \real$.
\item \textit{Triangle inequality, a.k.a., subadditivity}. $\norm{\bA+\bB} \leq \norm{\bA}+\norm{\bB}$, or $\norm{\bx+\by} \leq \norm{\bx}+\norm{\by}$ for any matrices $\bA, \bB\in \real^{m\times n}$ or vectors $\bx,\by\in \real^n$.
\begin{enumerate}
\item The triangle inequality also indicates~\footnote{
$\big| \, \norm{\bx} - \norm{\by}\, \big| \leq \norm{\bx-\by}$ since $$
\left.
\begin{aligned}
\norm{\bx} &= \norm{\bx-\by+\by} \leq \norm{\bx-\by} + \norm{\by};\\
\norm{\by} &= \norm{\by-\bx+\bx} \leq \norm{\by-\bx} + \norm{\bx},\\
\end{aligned}
\right\}
\implies 
\begin{aligned}
\norm{\bx}-\norm{\by}  &\leq \norm{\bx-\by}; \\
\norm{\by} -\norm{\bx} &\leq  \norm{\by-\bx},
\end{aligned}
$$
}
\begin{equation}
\big| \, \norm{\bA} - \norm{\bB}\, \big| \leq \norm{\bA-\bB}
\quad\text{or}\quad
\big| \, \norm{\bx} - \norm{\by}\, \big| \leq \norm{\bx-\by}.
\end{equation}

\end{enumerate}
\item \textit{Nonnegativity}. $\norm{\bA} \geq 0$ or $\norm{\bx}\geq 0$.
\begin{enumerate}
\item  the equality holds if and only if $\bA=\bzero $ or $\bx=\bzero$. 
\end{enumerate}
\end{enumerate}
\end{definition}
\index{Vector norm}
\index{Matrix norm}
The vector space $\real^n$, together with a given norm $\norm{\cdot}$, is called a \textit{normed vector space}.
On the other hand, one way to define norms for matrices is by viewing a matrix $\bA\in\real^{m\times n}$ as a vector in $\real^{mn}$, e.g., using the vectorization of the matrix.
What distinguishes a matrix norm is a property called \textit{submultiplicativity}: $\norm{\bA\bB}\leq \norm{\bA}\norm{\bB}$ if $\norm{\cdot}$ is a submultiplicative matrix norm (see discussions below). Almost all of the matrix norms we discuss are submultiplicative (Frobenius norm in Proposition~\ref{propo:submul_frob} and spectral norm in Proposition~\ref{propo:submul_spec}, both of which are special cases of the \textit{Schatten $p$-norm} discussed in Section~\ref{section:kyfan_knorm}, as a consequence of the singular value decomposition).

\index{Orthogonally invariant}
\index{Initarily invariant}
\begin{definition}[Orthogonally Invariant Norms]\label{definition:unitarily_invaria}
A matrix norm on $\real^{m\times n}$ is \textit{orthogonally invariant} if $\norm{\bU\bA\bV} =\norm{\bA}$ for all orthogonal $\bU\in\real^{m\times m}$ and $\bV\in\real^{n\times n}$ and for all $\bA\in\real^{m \times n}$; and  it is \textit{weakly orthogonally invariant} if  $\norm{\bU\bA\bU^\top} =\norm{\bA}$ for all orthogonal $\bU\in\real^{m\times m}$ and $\bA\in\real^{m\times m}$ is square.
Similarly, a vector norm on $\real^{n}$ is \textit{orthogonally invariant} if $\norm{\bQ\bx}=\norm{\bx}$ for all orthogonal $\bQ\in\real^{n\times n}$ and for all $\bx\in\real^n$.~\footnote{The term \textit{unitarily invariant} is used more frequently in the literature for complex matrices.}
\end{definition}

\begin{exercise}[Sum and Max or Two Norms]
Show that the sum or max  of two norms is a norm.
\end{exercise}


\index{Semi-norm}
\index{Semi-inner product}
\index{Inner product}
\begin{definition}[Inner Product]\label{definition:inner_prod}
In most cases, the norm can be derived from the vector \textit{inner product} $\langle\cdot, \cdot\rangle: \real^n\times \real^n\rightarrow \real$ (the inner product of vectors $\bx,\by\in\real^n$ is given by $\langle\bx,\by\rangle$), which satisfies the following three axioms:~\footnote{When $\inner{\bx}{\bx}=0$ for some nonzero $\bx$, the inner product is called a \textit{semi-inner product}.}
\begin{enumerate}
\item \textit{Commutativity}. $\langle\bx,\by\rangle = \langle\by,\bx\rangle$ for any $\bx,\by\in \real^n$. 
\item \textit{Linearity}. $\langle\lambda_1 \bx_1+\lambda_2\bx_2, \by\rangle = \lambda_1\langle\bx,\by\rangle+\lambda_2\langle\bx_2,\by\rangle$ for any $\lambda_1,\lambda_2 \in \real$ and $\bx,\by\in \real^n$.
\item \textit{Nonnegativity}. $\langle \bx,\bx \rangle \geq 0$ for any $\bx \in \real^n$.
\begin{enumerate}
\item $\langle \bx,\bx \rangle = 0$ if and only if $\bx=\bzero$.
\end{enumerate}
\end{enumerate}
The matrix inner product can be defined similarly as a function of $\langle\cdot, \cdot\rangle: \real^{m\times n}\times \real^{m\times n}\rightarrow \real$.
\end{definition}

For example, the \textit{Euclidean inner product}, defined as $\inner{\bu}{\bv}=\bu^\top\bv, \forall\bu,\bv\in\real^n$, is an inner product on vectors; the \textit{Frobenius inner product}, defined as $\inner{\bA}{\bB} = \trace(\bA^\top\bB)$, is an inner product on matrices;

\begin{exercise}[Cauchy-Schwarz Inequality for Inner Product]
Given an inner product $\langle\cdot, \cdot\rangle: \real^n\times \real^n\rightarrow \real$, show that 
$$
\abs{\langle \bu, \bv\rangle}^2 \leq \inner{\bu}{\bu}\inner{\bv}{\bv} 
\gap 
\text{for all }
\gap 
\bu,\bv\in\real^n.
$$
\textit{Hint: Consider $\bx \triangleq \inner{\bv}{\bv}\bu - \inner{\bu}{\bv}\bv$ and examine $\inner{\bx}{\bx}\geq 0$.}
\end{exercise}

\begin{exercise}[Norm Derived from Inner Product, and its Property]\label{exercise:norm_innr_pro}
Let $\langle\cdot, \cdot \rangle$ be an inner product on $\real^n\times \real^n$. Then the function $\norm{\cdot}=\langle\cdot, \cdot\rangle^{1/2}:  \real^n\rightarrow [0, \infty)$ is a norm on $\real^n$.
Show that the function satisfies 
\begin{equation}\label{equation:norm_inn_pro}
\textbf{(Parallelogram identity):}\,\, \norm{\bu+\bv}^2+\norm{\bu-\bv}^2 = 2(\norm{\bu}^2+\norm{\bv}^2), \,\, \forall \bu,\bv\in\real^n.
\end{equation} 
(This is known as the \textit{parallelogram identity}: consider a parallelogram in a vector space where the sides are represented by the vectors $\bu$ and $\bv$, the diagonals of this parallelogram are represented by the vectors $\bu+\bv$ and $\bu-\bv$.)
Prove the following \textit{polarization identity}:
\begin{equation}\label{equation:norm_pol_pro}
\textbf{(Polarization identity):} \gap \inner{\bu}{\bv} = \frac{1}{4}(\norm{\bu+\bv}^2-\norm{\bu-\bv}^2).
\end{equation}
\end{exercise}
Therefore, the vector space $\real^n$, together with a given inner product $\inner{\cdot}{\cdot}$, is called an \textit{inner product space}, which is also a normed vector space with the derived norm.

\index{$\ell_p$ norm}
\index{$\ell_2$ norm}
\index{$\ell_1$ norm}
\index{$\ell_\infty$ norm}
\subsection{Vector Norm}\label{section:vector-norm}
The vector norm is derived from the definition of the inner product. In most cases, the inner product in $\real^n$ is implicitly  referred to as  the  \textit{dot product (a.k.a., Euclidean inner product)}, defined by 
$
\langle\bx,\by\rangle = \sum_{i=1}^{n} x_i y_i.
$
The $\ell_2$ \textit{norm} (a.k.a., $\ell_2$ vector norm) is induced from the dot product
$
\norm{\bx}_2 = \sqrt{\langle\bx,\bx\rangle} = \sqrt{\sum_{i=1}^{n}x_i^2}.
$
More generally,  given a $p\geq 1$ \footnote{When $p<1$, the $\ell_p$ does not meet the third axiom in Definition~\ref{definition:matrix-norm}.}, the $\ell_p$ \textit{norm} (a.k.a., $\ell_p$ vector norm, \holders norm) is given by 
\begin{equation}
\norm{\bx}_p = \sqrt[p]{ \sum_{i=1}^{n}|x_i|^p  },
\quad 
p\geq 1,
\end{equation}
where the positive homogeneity and nonnegativity are apparent, and the triangle inequality follows from Minkowski's inequality (Theorem~\ref{theorem:minkoski}; alternatively, we can prove the triangle inequality of the $\ell_p$ norm using H{\"o}lder's inequality; see below).
From this,  the $\ell_1$ and $\ell_\infty$ \textit{norms}  can be obtained by, respectively, 
$$
\norm{\bx}_1 = \sum_{i=1}^{n} |x_i|
\gap \text{and}\gap
\norm{\bx}_\infty = \mathop{\max}_{i=1,2,\ldots,n} |x_i| .
$$

Following the definition of the $\ell_p$ norm, we can obtain the famous H{\"o}lder's inequality in Theorem~\ref{theorem:holder-inequality}. Thus, the $\ell_p$ norm is sometimes referred to as  \textit{H{\"o}lder's norm}. 
Conversely, H{\"o}lder's inequality can be used to prove the validity of the  $\ell_p$norm.
\begin{proof}[Triangle Inequality of $\ell_p$ Norm] 
\textbf{Case $p>1$ for $\norm{\bx}_p$.}
From  H{\"o}lder's inequality (requires $p,q>1$; Theorem~\ref{theorem:holder-inequality}), we can prove the validity of the $\ell_p$ norm. Let $\frac{1}{p}+\frac{1}{q} = 1$, it follows that 
$$
\begin{aligned}
\sum_{i=1}^{n}& (|x_i| + |y_i|)^p = \sum_{i=1}^{n}|x_i| (|x_i|+|y_i|)^{p-1} + \sum_{i=1}^{n}|y_i| (|x_i|+|y_i|)^{p-1} \\
&\leq \left(\sum_{i=1}^{n}  |x_i|^p\right)^{1/p} \left( \sum_{i=1}^{n} (|x_i| + |y_i|)^{(p-1)q}  \right)^{1/q}
\left(\sum_{i=1}^{n} |y_i|^p\right)^{1/p}  \left( \sum_{i=1}^{n} (|x_i| + |y_i|)^{(p-1)q}  \right)^{1/q} \\
&= \left[\left(\sum_{i=1}^{n}  |x_i|^p   \right)^{1/p} +\left(\sum_{i=1}^{n} |y_i|^p\right)^{1/p} \right]  \left( \sum_{i=1}^{n} (|x_i| + |y_i|)^{(p-1)q}  \right)^{1/q} ,\\
\end{aligned}
$$
Since $(p-1)q=p$, the above inequality implies 
$$
\left( \sum_{i=1}^{n} (|x_i| + |y_i|)^{p}  \right)^{1/p} \leq  \left(\sum_{i=1}^{n}  |x_i|^p   \right)^{1/p} +\left(\sum_{i=1}^{n} |y_i|^p\right)^{1/p} = \norm{\bx}_p + \norm{\by}_p.
$$
Therefore, it follows that 
$$
\norm{\bx+\by}_p = \left( \sum_{i=1}^{n} (|x_i + y_i|)^{p}  \right)^{1/p}  \leq  \left( \sum_{i=1}^{n} (|x_i| + |y_i|)^{p}  \right)^{1/p} 
\leq \norm{\bx}_p + \norm{\by}_p.
$$

\paragraph{Case $p=1$ for $\norm{\bx}_1$.}
It is straightforward that $\norm{\bx+\by}_1 \leq \norm{\bx}_1 + \norm{\by}_1$.

\paragraph{Case $p<1$ for $\norm{\bx}_p$.}
However, when $p<1$, the triangle inequality condition for vector norms is not satisfied. For example, when $p=1/3$, and $\bx=[0,1]^\top$ and $\by=[1,0]^\top$. It follows that 
$$
\begin{aligned}
\norm{\bx+\by}_{1/3}= 8, \gap \norm{\bx}_{1/3} = 1, \gap \norm{\by}_{1/3}=1
\gap \implies 
\norm{\bx+\by}_{1/3}> \norm{\bx}_{1/3}+\norm{\by}_{1/3}.
\end{aligned}
$$
This proves the validity of the $\ell_p$ norm.
\end{proof}


Alternatively, the triangle inequality for $\ell_2$ norm can be proved using the Cauchy-Schwarz inequality.
\begin{proof}[Triangle Inequality of $\ell_2$  Norm by Cauchy-Schwarz Inequality]
By the Cauchy-Schwarz inequality (Proposition~\ref{proposition:cauchy-schwarz-inequ}), for any vectors $\bx$ and $\by$, we have $\bx^\top \by \leq \normtwo{\bx} \cdot \normtwo{\by}$. Therefore,
$$
\begin{aligned}
\normtwo{\bx+\by}^2 &= \normtwo{\bx}^2 + \normtwo{\by}^2 + 2\bx^\top\by
\leq  \normtwo{\bx}^2 + \normtwo{\by}^2 + 2 \normtwo{\bx} \cdot \normtwo{\by} 
= (\normtwo{\bx} + \normtwo{\by})^2.
\end{aligned}
$$
This results in $\normtwo{\bx+\by} \leq \normtwo{\bx}+\normtwo{\by}$.
\end{proof}

\begin{exercise}[Orthogonal Invariance of $\ell_2$]\label{exercise:orthogo_ell2}
Show that the $\ell_2$ norm is orthogonally invariant: $\norm{\bQ\bx}_2=\norm{\bx}_2$ for all $\bx\in\real^n$ if $\bQ$ is orthogonal.
How about other $\ell_p$ norms?
\end{exercise}

\begin{exercise}[$\ell_1$, $\ell_\infty$ norm]
Show that $\ell_1$, $\ell_\infty$ are not norms derived from an inner product 
\textit{Hint: Examine $\be_1$ and $\be_2$, and the parallelogram identity in \eqref{equation:norm_inn_pro} or the polarization identity in \eqref{equation:norm_pol_pro}.}
\end{exercise}

\begin{exercise}[$k$-Norm]
Show that the $k$-norm $\norm{\bx}^{[k]}$ on $\bx\in\real^m$ defined as the sum of the largest $k$ absolute elements of a vector is a valid norm. 
Specifically, $\norm{\bx}_\infty = \norm{\bx}^{[1]}$ and $\norm{\bx}^{[n]}=\norm{\bx}_1$. 
Show that
\begin{equation}
\norm{\bx}_\infty = \norm{\bx}^{[1]} 
\leq 
\norm{\bx}^{[2]}
\leq\ldots 
\leq 
\norm{\bx}^{[n]}=\norm{\bx}_1.
\end{equation}
\end{exercise}

\section*{Properties}
For any vector norm, we also have the following properties.
\begin{proposition}[$\ell_p$ Norm Inequalities]\label{prop:lp_norm_ineqs}
Given any $\bx\in\real^n$, the function $f(p)=\norm{\bx}_p$ is nonincreasing on $p=[1,\infty)$. Therefore, 
\begin{equation}
\norm{\bx}_\infty\leq \norm{\bx}_q \leq \norm{\bx}_p\leq \normtwo{\bx} \leq \normone{\bx},
\quad \text{if }  \infty \geq q \geq p\geq 2\geq 1.
\end{equation}
If $\bx$ has at least two nonzero elements, the strict inequalities can be attained.
This can also be observed in the figures of unit balls (Figure~\ref{fig:p-norm-2d} and Figure~\ref{fig:p-norm-comparison-3d}).
Moreover, for any $\infty\geq p \geq 1$, we also have 
\begin{equation}
\norm{\bx}_\infty \leq \norm{\bx}_p \leq  n^{1/p} \norm{\bx}_\infty,
\quad \text{if }  \infty \geq  p \geq 1.
\end{equation}
\end{proposition}
\begin{proof}[of Proposition~\ref{prop:lp_norm_ineqs}]
Let $S(p) \triangleq \norm{\bx}_p^p$ and take the logarithm of $f(p)$, we have
$\ln f(p) = \frac{1}{p} \ln S(p). $
Differentiate both sides with respect to $p$:
$$\frac{f'(p)}{f(p)} = \frac{d}{dp} \left( \frac{1}{p} \ln S(p) \right)
= -\frac{1}{p^2} \ln S(p) + \frac{1}{p} \cdot \frac{S'(p)}{S(p)}
. 
$$
This shows
$$f'(p) = f(p) \left( -\frac{1}{p^2} \ln S(p) + \frac{1}{p} \cdot \frac{S'(p)}{S(p)} \right)
=
\frac{1}{p}\norm{\bx}_p^{1-p}\sum_{i=1}^{n} \alpha_i,
$$
where 
$$
S'(p)=\sum_{i=1}^{n}\abs{x_i}^p \ln\abs{x_i}
\quad \text{and} \quad 
\alpha_i =
\left\{
\begin{aligned}
&\abs{x_i}^p (\ln \abs{x_i} - \ln \norm{\bx}_p), &\text{if } x_i\neq 0;\\
&0, &\text{if } x_i=0.
\end{aligned}
\right.
$$
This shows $f'(p)\leq 0$ since $\ln \abs{x_i} - \ln \norm{\bx}_p\leq 0$.
If $\bx$ has at least two nonzero elements, $f'(p)< 0$ and $f(p)$ is decreasing.

For the second part, following the definition of the $\ell_p$ norm, we notice that 
$$
(\norm{\bx}_\infty )^p = \left(\mathop{\max}_{i=1,2,\ldots,n} |x_i|\right)^p \leq \sum_{i=1}^{n} |x_i|^p \leq n \mathop{\max}_{i=1,2,\ldots,n} |x_i|^p = n (\norm{\bx}_\infty )^p,
$$
from which the result follows.
\end{proof}



\begin{proposition}[Orthogonally Invariant Vector Norm]\label{proposition:orinv_vecnorm}
Let $\norm{\cdot}$ be an orthogonally invariant norm on $\real^n$ and $\norm{\cdot}_2$ be the $\ell_2$ norm. Then, $\norm{\bx}=\norm{\bx}_2\norm{\be_1}$ for any $\bx\in\real^n$.
\end{proposition}
\begin{proof}[of Proposition~\ref{proposition:orinv_vecnorm}]
Consider $\gamma=\norm{\bx}_2$ and $\widehat{\bx}=\bx/\gamma$, so $\bx=\gamma \widehat{\bx}$. 
Let $\bU$ be an orthogonal matrix whose first column is $\widehat{\bx}$.
From the definition of a norm and the orthogonally invariant property, we have $\norm{\bx}=\gamma\norm{\widehat{\bx}}=\gamma\norm{\bU\be_1}=\gamma\norm{\be_1}$. 
\end{proof}
The result also indicates that the $\ell_2$ norm is the only orthogonally invariant norm that has $\norm{\be_1}=1$.

\paragraph{Dual norm.}
Consider the $\ell_p$ vector norm. From \holders inequality (Theorem~\ref{theorem:holder-inequality}), we have 
$
\bx^\top\by \leq \norm{\bx}_p \norm{\by}_q,
$
where $p,q>1$, $\frac{1}{p}+\frac{1}{q}=1$, and $\bx,\by\in \real^n$. The equality holds if the two sequences $\{|x_i|^p\}$ and $\{|y_i|^q\}$ are linearly dependent. This implies
\begin{equation}\label{equation:dual_norm_equa}
\mathop{\max}_{\norm{\by}_q=1} \bx^\top\by = \norm{\bx}_p.
\end{equation}
For this reason, we call $\norm{\cdot}_q$ the \textit{dual norm} of $\norm{\cdot}_p$.
On the other hand, for each $\bx\in \real^n$ with $\norm{\bx}_p=1$, there exists a vector $\by\in \real^n$ such that $\norm{\by}_q=1$ and $\bx^\top\by=1$.
Notably, the $\ell_2$ norm is dual to itself, and the $\ell_1$ and $\ell_\infty$ norms are dual to each other.

We conclude this section by introducing an important property of vector norms that will be often useful. The following equivalence of vector norms states that if a vector is small in one norm, it will also be small in another norm as well, and vice versa.
\begin{theorem}[Equivalence of Vector Norms]\label{theorem:equivalence-vector-norm}
Let $\norm{\cdot}_a$ and $\norm{\cdot}_b$ be two different vector norms: $\real^n\rightarrow \real$. Then there exist positive scalars $\alpha$ and $\beta$ such that for all $\bx \in \real^n$, the following inequality holds:
$$
\alpha\norm{\bx}_a \leq \norm{\bx}_b \leq \beta\norm{\bx}_a.
$$
This implies 
$$
\frac{1}{\beta}\norm{\bx}_b  \leq \norm{\bx}_a \leq \frac{1}{\alpha}\norm{\bx}_b.
$$
Hence, the name ``equivalence" is justified
\end{theorem}
\begin{proof}[of Theorem~\ref{theorem:equivalence-vector-norm}]
In advanced calculus, it is stated that $\mathop{\sup}_{\bx\in \sS} f(\bx)$ is attained for some vector $\bx\in \sS$ as long as $f(\cdot)$ is continuous and $\sS$ is a compact set (closed and bounded); see Theorem~\ref{theorem:weierstrass_them}. When the supremum is an element in the set $\sS$, this supremum is known as the maximum such that $\mathop{\sup}_{\bx\in\sS} f(\bx) =\mathop{\max}_{\bx\in\sS} f(\bx)$.
Without loss of generality, we assume $\bx\neq 0$. Then we have:
$$
\begin{aligned}
\norm{\bx}_b &= \frac{\norm{\bx}_b}{\norm{\bx}_a} \norm{\bx}_a 
\leq \mathop{\sup}_{\bz\neq 0} \frac{\norm{\bz}_b}{\norm{\bz}_a} \norm{\bx}_a
= \mathop{\sup}_{\norm{\by}_a=1} \norm{\by}_b \norm{\bx}_a
=\norm{\bx}_a \mathop{\max}_{\norm{\by}_a=1} \norm{\by}_b.
\end{aligned}
$$
The last equality holds since $\{\by: \norm{\by}_a=1\}$ is a compact set.
By setting $\beta\triangleq \mathop{\max}_{\norm{\by}_a=1} \norm{\by}_b  $, we have $\norm{\bx}_b \leq \beta\norm{\bx}_a$.
From the above argument, there exists a $\tau$ such that 
$
\norm{\bx}_a \leq \tau \norm{\bx}_b.
$
Let $\alpha = \frac{1}{\tau}$, we have $\alpha\norm{\bx}_a \leq \norm{\bx}_b$, from which the result follows.
\end{proof}

\begin{example}[Equivalence of Vector Norms]
The following inequalities hold for all $\bx\in\real^n$:
$$
\begin{aligned}
\norm{\bx}_\infty &\leq \norm{\bx}_1 \leq n\norm{\bx}_\infty; \\
\norm{\bx}_\infty &\leq \norm{\bx}_2 \leq \sqrt{n}\norm{\bx}_\infty; \\
\norm{\bx}_2 &\leq \norm{\bx}_1 \leq \sqrt{n}\norm{\bx}_2. \\
\end{aligned}
$$
This example shows the equivalence of the $\ell_1, \ell_2$, and $\ell_\infty$ vector norms.
\end{example}

\index{Unit ball}
\section*{Unit Ball}
\begin{SCfigure}
	\centering
	\includegraphics[width=0.5\textwidth]{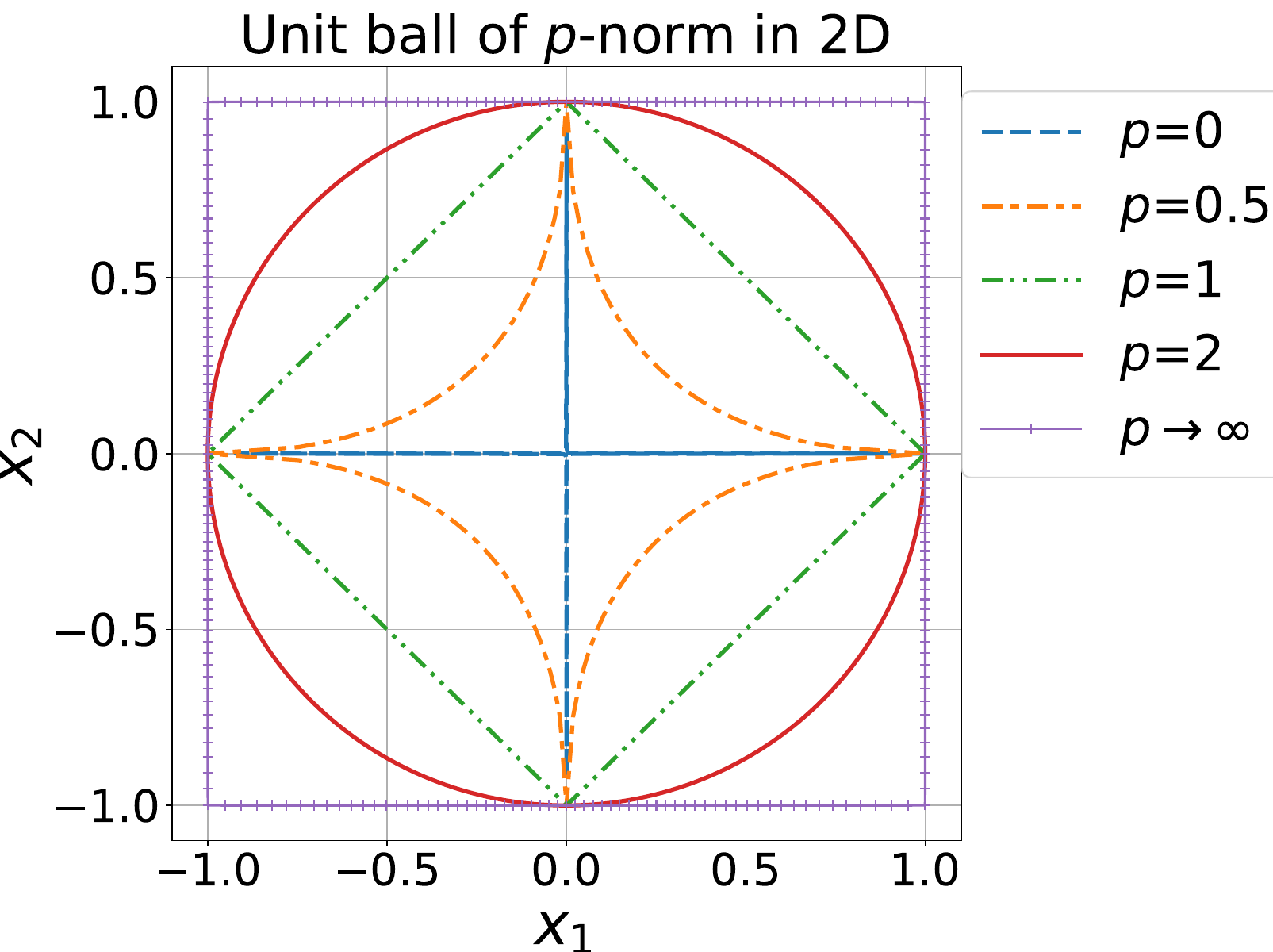}
	\caption{Unit ball of $\ell_p$ norm in two-dimensional space. When $p<1$, the metric is not a norm since it does not meet the third axiom of the norm in Definition~\ref{definition:matrix-norm}.}
	\label{fig:p-norm-2d}
\end{SCfigure}

\begin{figure}[H]
\centering  
\vspace{-0.35cm} 
\subfigtopskip=2pt
\subfigbottomskip=2pt 
\subfigcapskip=-5pt 
\subfigure[$p=\infty$.]{\label{fig:p-norm-3d1}
\includegraphics[width=0.18\linewidth]{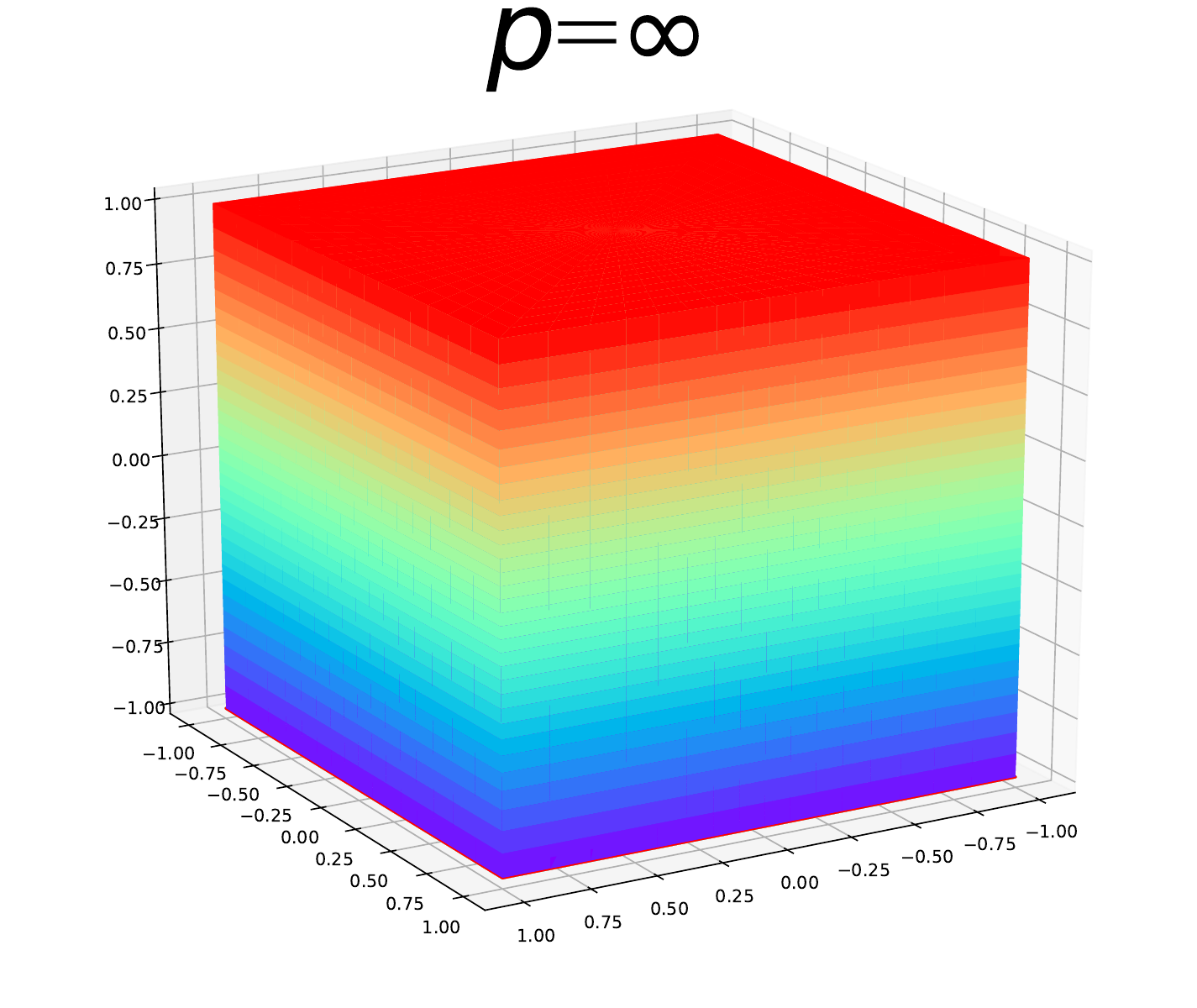}}
\subfigure[$p=2$.]{\label{fig:p-norm-3d2}
\includegraphics[width=0.18\linewidth]{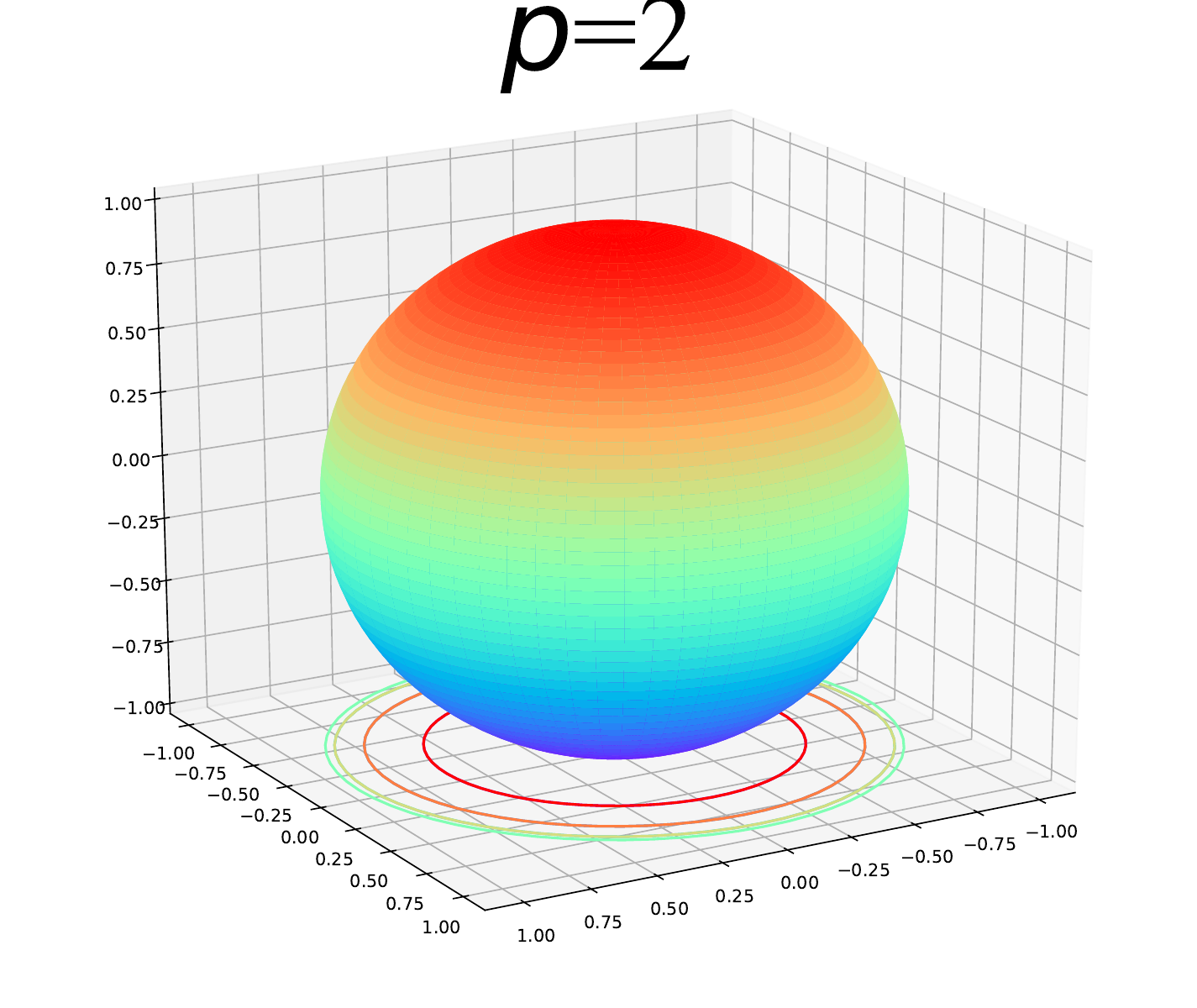}}
\subfigure[$p=1$.]{\label{fig:p-norm-3d3}
\includegraphics[width=0.18\linewidth]{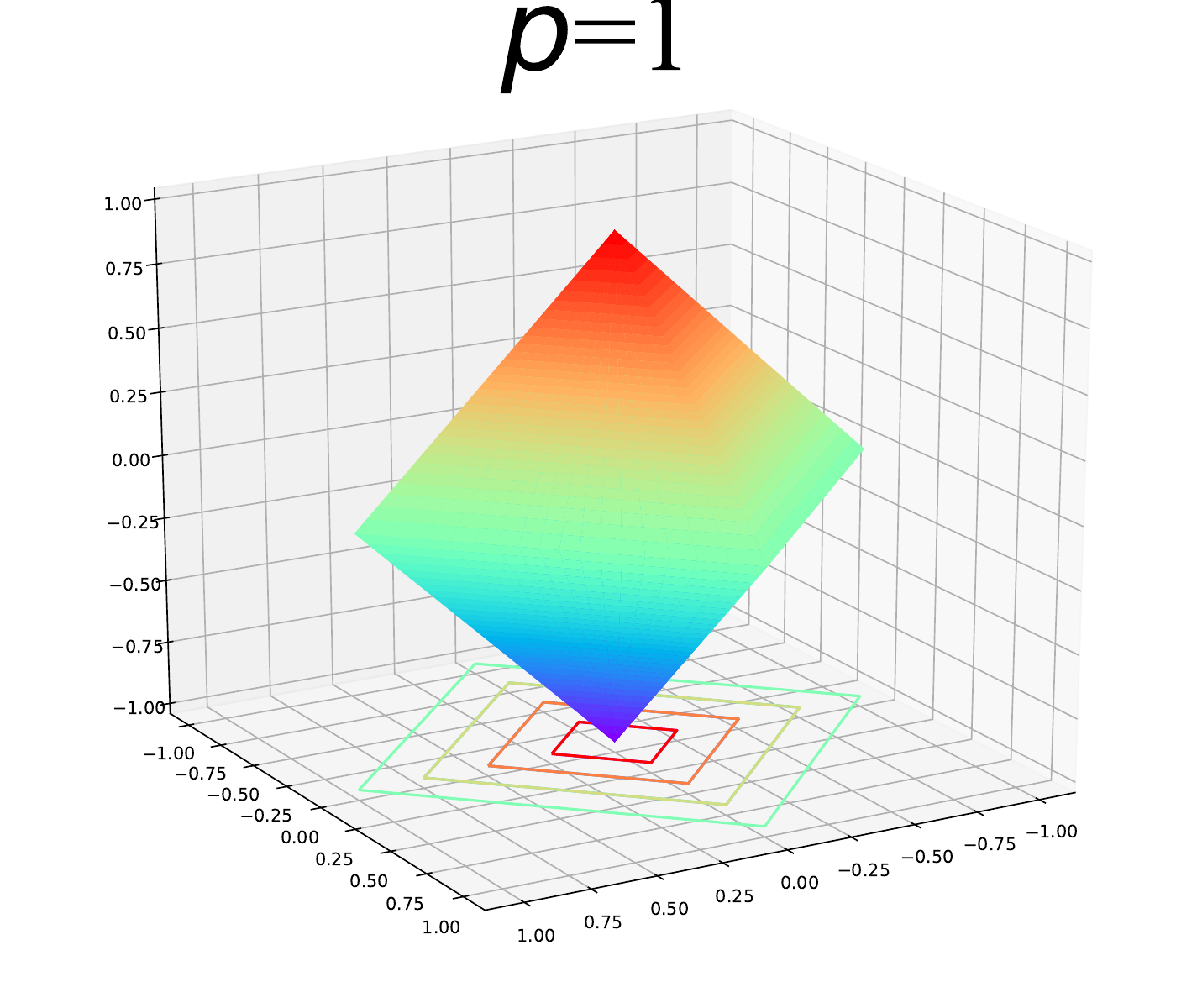}}
\subfigure[$p=0.5$.]{\label{fig:p-norm-3d4}
\includegraphics[width=0.18\linewidth]{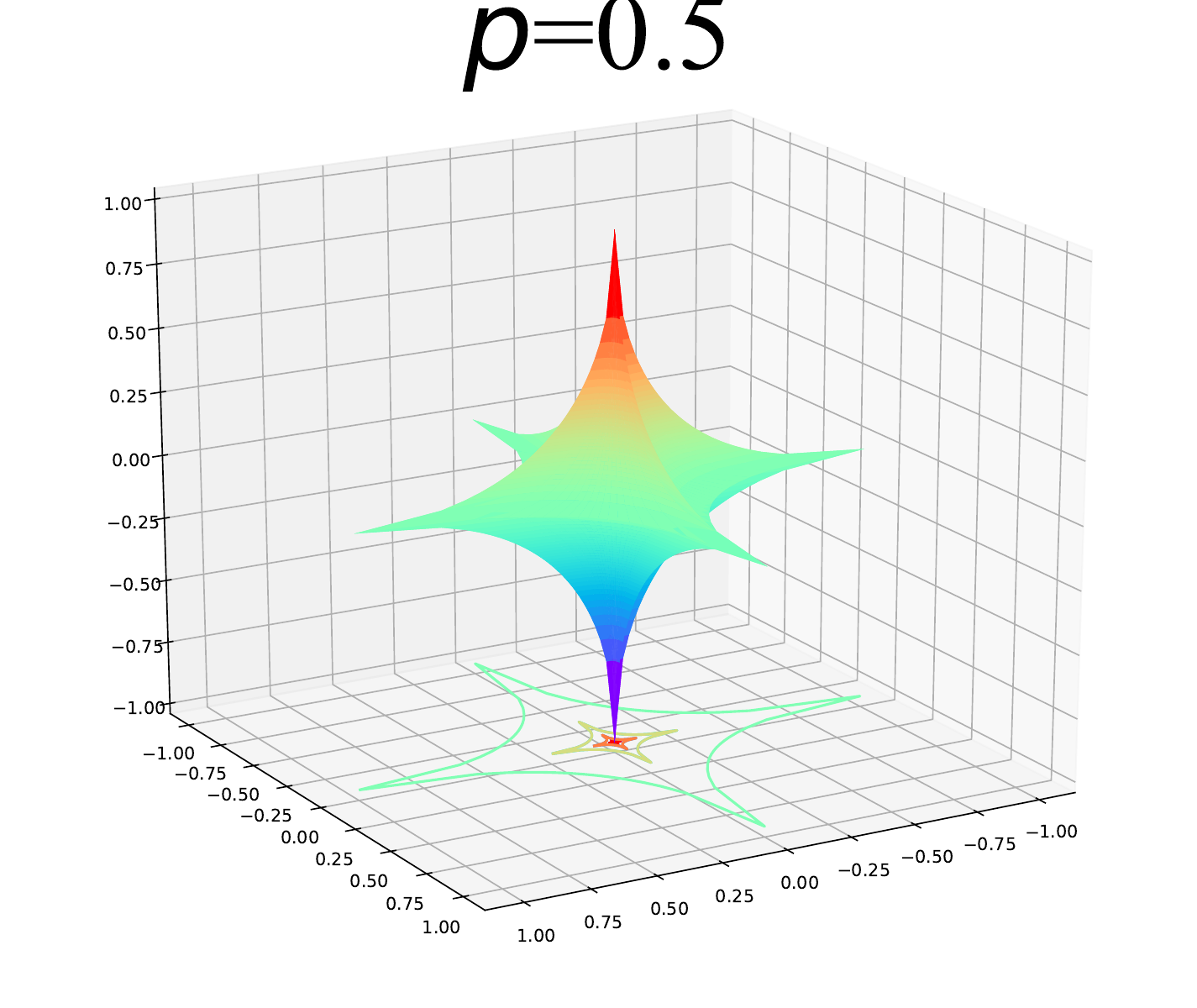}}
\subfigure[$p=0$.]{\label{fig:p-norm-3d5}
\includegraphics[width=0.18\linewidth]{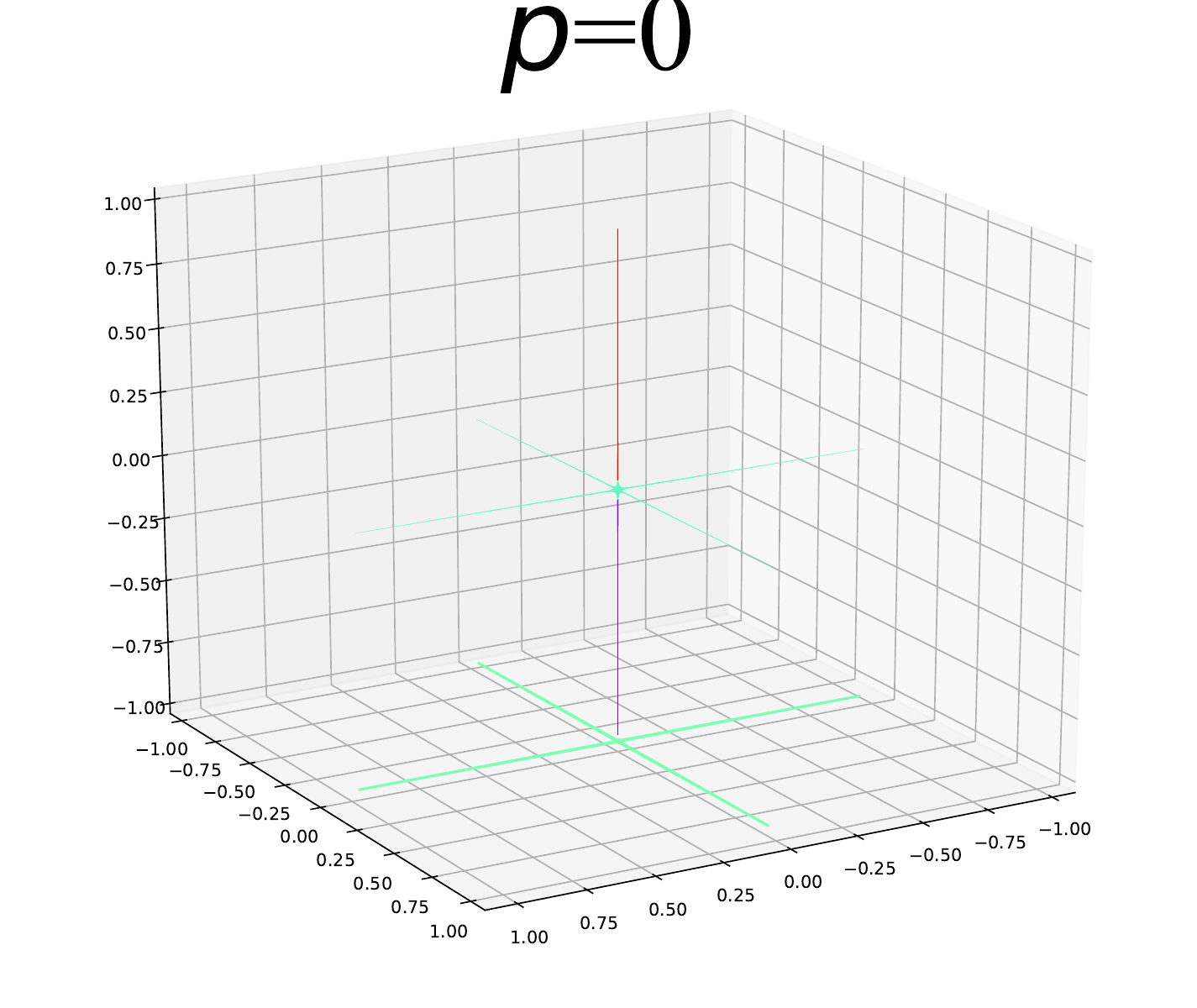}}
\caption{Unit ball of $\ell_p$ norm in three-dimensional space. When $p<1$, the metric is not a norm since it does not meet the third axiom of the norm in Definition~\ref{definition:matrix-norm}.}
\label{fig:p-norm-comparison-3d}
\end{figure}

The \textit{unit balls} of the norms represent  the set of all points whose distances from the origin (i.e., the zero vector) equal 1.
If the distance is defined by the $\ell_p$ norm,  the unit ball is the collection of 
$$
\sB = \{\bx: \norm{\bx}_p=1\}.
$$
The comparison of the  $\ell_p$ norm in two-dimensional and three-dimensional spaces with different values of $p$ is depicted in Figure~\ref{fig:p-norm-2d} and Figure~\ref{fig:p-norm-comparison-3d}.

Vector norms play a crucial role in machine learning. In Section~\ref{section:application-ls-qr}, we mentioned that the least squares aims to minimize the squared distance between the observation $\bb$ and the expected observation $\bA\bx$: $\norm{\bA\bx-\bb}_2$, i.e., the $\ell_2$ norm of $\bA\bx-\bb$.
Alternatively, minimizing the $\ell_1$ norm of the difference between the observed and predicted values can lead to a more robust estimation of $\bx$, particularly in the presence of outliers  \citep{zoubir2012robust}.

\section*{$\bQ$-Inner Product, $\bQ$-Norm, $\bS$-Norm}
The dot product is not the only possible inner product that can be defined over $\real^n$. Given a positive definite $n\times n$ matrix $\bQ$, a $\bQ$\textit{-dot product} can be defined as~\footnote{When $\bQ$ is positive semidefinite, the $\bQ$-dot product is a semi-inner product.}
$$
\langle \bx,\by \rangle_{\bQ} = \bx^\top\bQ\by.
$$
One can verify that the $\bQ$-dot product defined above satisfies the three axioms for the inner product discussed at the beginning of this section. When $\bQ=\bI$, we recover the dot product (Euclidean inner product). From the $\bQ$-dot product, the \textit{$\bQ$-$\ell_2$-norm} can be defined as 
$
\norm{\bx}_{\bQ,2} = \sqrt{\bx^\top\bQ\bx}.
$
More generally, given a positive definite matrix $\bQ\in\real^{n\times n}$ and  a general norm $\norm{\cdot}$ on $\real^n$, one may consider the \textit{$\bQ$-norm} as:
\begin{equation}
\textbf{$\bQ$-norm:} \gap\norm{\bx}_{\bQ} = \norm{\sqrt{\bQ}\bx}
\end{equation}	
Alternatively, given a matrix $\bS\in\real^{m\times n}$ with full column rank and a general norm $\norm{\cdot}$ on $\real^m$,  the \textit{$\bS$-norm} is defined as:
\begin{equation}\label{equation:s_norm}
\textbf{$\bS$-norm:} \gap \norm{\bx}_{\bS}=\norm{\bS\bx},
\end{equation}	
which is a norm on $\real^n$.
To see this, we introduce the following lemma.
\begin{lemma}[Construct Norms from Other Norms]\label{lemma:construct_norm}
Let $\norm{\cdot}$ be a norm over $\real^m$. Given a matrix $\bS\in \real^{m\times n}$, then $\norm{\bS(\cdot)}$ is a semi-norm~\footnote{Might be equal to zero for nonzero vectors in $\real^n$.} over $\real^n$. 
If $\bS$ has full column rank $n$, then $\norm{\bS(\cdot)}$ is a norm over $\real^m$.
\end{lemma}
\begin{proof}[of Lemma~\ref{lemma:construct_norm}]
Given any vectors $\bx, \by\in \real^m$, following from Definition~\ref{definition:matrix-norm}, we have 
$$
\norm{\bx}\geq 0, \gap \norm{\lambda \bx} = |\lambda|\cdot  \norm{\bx}, \gap \norm{\bx+\by} \leq \norm{\bx} + \norm{\by}.
$$
Suppose $\bc, \bd\in \real^n$.
Therefore, we have  
$$
\norm{\bS\bc}\geq 0, \gap \norm{\lambda \bS\bc} = |\lambda|\cdot  \norm{\bS\bc}, \gap \norm{\bS\bc+\bS\bd} \leq \norm{\bS\bc} + \norm{\bS\bb}.
$$
Therefore, $\norm{\bS(\cdot)}$ is a semi-norm.
If $\bS$ has full column rank, $\bS\bc=\bzero$ only if $\bc=\bzero$. The semi-norm becomes a norm.
This completes the proof.
\end{proof}

\begin{exercise}[$\bQ$-Norm, $\bS$-Norm]\label{exercise:qnorm}
Show that the $\bQ$-norm and $\bS$-norm defined above satisfy the three axioms of norms in Definition~\ref{definition:matrix-norm}.
Show that when $\norm{\cdot}$ is the $\ell_2$ norm, the $\bS$-norm reduces to the $\bQ$-$\ell_2$-norm. 
\textit{Hint: Use Cholesky and spectral decompositions (\S~\ref{section:conc_pd}, Theorem~\ref{theorem:nonsingular-factor-of-PD}).}
\end{exercise}

\begin{exercise}[Weighted $\ell_p$ Norm]
Let $w_1, w_2, \ldots,w_p$ be positive real numbers and $p\geq 1$. Show that  the weighted $\ell_p$ norm $\norm{\bx}=\sqrt[p]{ \sum_{i=1}^{n}w_i|x_i|^p  }$ on $\real^n$ is a valid norm. \textit{Hint: Show that it is a norm of the form $\norm{\bS\bx}_p$.} 
\end{exercise}


\subsection{Matrix Norm}\label{appendix:matrix-norm-sect2}

\index{Submultiplicativity}
\index{Submultiplicative matrix norms}
\subsection*{Submultiplicativity of Matrix Norms}
In some texts, a matrix norm that is not \textit{submultiplicative} is termed as a \textit{vector norm on matrices} or a \textit{generalized matrix norm}.
The submultiplicativity of a matrix norm is important for the analysis of square matrices, although the definition of a matrix norm applies to both square and rectangular matrices.
For a submultiplicative matrix norm $\norm{\cdot}$ that satisfies $\norm{\bA\bB}\leq \norm{\bA}\norm{\bB}$, consider $\bA\in\real^{n\times n}$, it follows that 
\begin{equation}\label{equation:power_subm}
\norm{\bA^2} \leq \norm{\bA}^2
\quad\implies\quad
\norm{\bA^k} \leq \norm{\bA}^k, \forall k\in\{1,2,\ldots,\}.
\end{equation}
Therefore, if the matrix is idempotent, i.e., $\bA^2=\bA$, we have $\norm{\bA}\geq 1$, which also indicates 
\begin{equation}\label{equation:power_subm2}
\norm{\bI}\geq 1,
\quad\text{if}\quad
\norm{\cdot} \text{ is submultiplicative}.
\end{equation}
On the other hand, if $\bA$ is nonsingular, we have the inequality for submultiplicative norms: 
$$
1\leq \norm{\bI}=\norm{\bA\bA^{-1}}\leq \norm{\bA}\norm{\bA^{-1}}.
$$  
That is, a submultiplicative norm has $\norm{\bI}\geq 1$ and is \textit{normalized} if and only if $\norm{\bI}=1$.

\paragraph{Bounds on the spectral radius.}
The submultiplicativity of a matrix norm can also be used to find a bound on the spectral radius of a matrix (Definition~\ref{definition:spectrum}).
Given an eigenpair $(\lambda,\bx)$ of a matrix $\bA\in\real^{n\times n}$, and consider the matrix $\bX=\bx\bone^\top=[\bx,\bx,\ldots,\bx]$. 
For a submultiplicative matrix norm $\norm{\cdot}$ on $\real^{n\times n}$, it follows that 
$
\abs{\lambda}\norm{\bX}
=
\norm{\lambda\bX}
=
\norm{\bA\bX}
\leq \norm{\bA}\norm{\bX}.
$
Thus, the spectral radius is bounded by the matrix norm $\norm{\bA}$.
Similarly, we can prove that $\abs{\lambda^{-1}}\leq \norm{\bA^{-1}}$.
Therefore, we  obtain the lower and upper bounds on the spectral radius $\rho(\bA)$:
\begin{equation}\label{equation:bounds_submulnorm}
\textbf{Bounds on spectral radius:} \gap \frac{1}{\norm{\bA^{-1}}}\leq \rho(\bA)\leq \norm{\bA}.
\end{equation}
The property of submultiplicativity further helps us understand the bounds on the spectral radius of the product of matrices. For instance, given matrices $\bA$ and $\bB$ (provided  the matrix product $\bA\bB$ is defined), we have
\begin{equation}
	\rho(\bA\bB)\leq \norm{\bA\bB} \leq \norm{\bA}\norm{\bB}.
\end{equation}

\index{Convergent matrices}
\paragraph{Power of square matrices.}
When $\norm{\bA}\leq 1$ for a square matrix $\bA\in\real^{n\times n}$, \eqref{equation:power_subm} shows that $\norm{\bA^k} \leq \norm{\bA}^k\stackrel{k\rightarrow \infty}{=}0$.~\footnote{Regardless which norm is used, as long as the norm is submultiplicative.}
This means $\bA^k$ tends to  zero matrix when $k\rightarrow \infty$ if $\norm{\bA}\leq 1$.
A matrix with this property is called \textit{convergent}.
The convergence of a matrix can also be characterized by its spectral radius.
\begin{exercise}[Convergence Matrices]\label{exer:conv_mat}
Let $\bA\in\real^{n\times n}$. Show that $\mathoplim{k\rightarrow \infty}\bA^k=\bzero$ if and only if the spectral radius $\rho(\bA)<1$.
\textit{Hint: Examine $\bA^k\bx=\lambda^k\bx$ and Problem~\ref{prob:bounded_spec}.}
\end{exercise}

\begin{exercise}[Convergence Matrices]
Let $\bA\in\real^{n\times n}$ and  $\epsilon>0$. Show that there is a constant $C=C(\bA,\epsilon)$ such that $\abs{(\bA^k)_{ij}} \leq C(\rho(\bA)+\epsilon)^k$ for $k\in\{1,2,\ldots\}, i,j\in\{1,2,\ldots,n\}$.
\textit{Hint: Use Exercise~\ref{exer:conv_mat} and examine $\bB=\frac{1}{\rho(\bA)+\epsilon}\bA$, whose spectral radius is strictly less than 1.}
\end{exercise}

\index{Gelfand formula}
The Gelfand formula provides a way to estimate the spectral radius of a matrix using submultiplicative norms.
\begin{exercise}[Gelfand Formula]\label{exercise:gelfand_formula}
Let $\bA\in\real^{n\times n}$ and let $\norm{\cdot}$ be a submultiplicative matrix norm on $\real^{n\times n}$. 
Show that $\rho(\bA)=\mathoplim{k\rightarrow \infty}\norm{\bA^k}^{1/k}$.
\textit{Hint: Consider $\rho(\bA^k)=\rho(\bA)^k\leq \norm{\bA}$ and examine $\bB=\frac{1}{\rho(\bA)+\epsilon}\bA$, whose spectral radius is strictly less than 1.}
\end{exercise}

\index{Frobenius}
\subsection*{Frobenius Norm}
The norm of a matrix serves the same purpose as the norm of a vector.
 One of the important matrix norms is defined as follows, the matrix counterpart of the vector $\ell_2$-norm.
\begin{definition}[Frobenius Norm]\label{definition:frobenius}
The Frobenius norm of a matrix $\bA\in \real^{m\times n}$ is defined as 
$$
\norm{\bA}_F = \sqrt{\sum_{i=1,j=1}^{m,n} (a_{ij})^2}=\sqrt{\trace(\bA\bA^\top)}=\sqrt{\trace(\bA^\top\bA)} = \sqrt{\sigma_1^2+\sigma_2^2+\ldots+\sigma_r^2},
$$
i.e., the square root of the sum of the squares of the elements of $\bA$. The values of $\sigma_i$'s are the singular values of $\bA$, and $r$ is the rank of $\bA$.
\end{definition}
The equivalence of $\sqrt{\sum_{i=1,j=1}^{m,n} (a_{ij})^2}$, $\sqrt{\trace(\bA\bA^\top)}$, and $\sqrt{\trace(\bA^\top\bA)}$ is straightforward. 
The equivalence between $\sqrt{\trace(\bA\bA^\top)}$ and $\sqrt{\sigma_1^2+\sigma_2^2+\ldots+\sigma_r^2}$ can be shown using the singular value decomposition (SVD). Suppose $\bA$ admits the SVD $\bA = \bU\bSigma\bV^\top$, then:
$$
\sqrt{\trace(\bA\bA^\top)} = \sqrt{\trace(\bU\bSigma\bV^\top \bV\bSigma\bU^\top)} = \sqrt{\trace(\bSigma^2)}=\sqrt{\sigma_1^2+\sigma_2^2+\ldots+\sigma_r^2}.
$$
Apparently, the Frobenius norm can  also be defined using the  vector $\ell_2$ norm such that $\norm{\bA}_F = \sqrt{\sum_{i=1}^{n} \norm{\ba_i}^2}$, where $\ba_i$ for all $i \in \{1,2,\ldots, n\}$ are the columns of $\bA$.

\index{Submultiplicativity}
\begin{proposition}[Submultiplicativity of Frobenius]\label{propo:submul_frob}
The Frobnenius norm is submultiplicative. That is, $\norm{\bA\bB}_F\leq \norm{\bA}_F\norm{\bB}_F$.
\end{proposition}
\begin{proof}[of Proposition~\ref{propo:submul_frob}]
Suppose $\bA\in\real^{m\times n}$ and $\bB\in\real^{n\times p}$. We have 
$$
\norm{\bA\bB}_F = 
\bigg(\sum_{i,j=1}^{m,p} \big(\sum_{k=1}^{n}a_{ik}b_{kj} \big)^2 \bigg)^{1/2}
\leq 
\bigg(\sum_{i,j=1}^{m,p} \big(\sum_{k=1}^{n}a_{ik}^2 \big)\big(\sum_{k=1}^{n}b_{kj}^2 \big) \bigg)^{1/2}
=
\norm{\bA}_F\norm{\bB}_F.
$$
This completes the proof.
\end{proof}

\index{Orthogonally invariance}
\begin{proposition}[Orthogonally Invariance of Frobenius]\label{proposition:frobenius-orthogonal-equi}
Let $\bA\in \real^{m\times n}$ be given, and let $\bU\in \real^{m\times m}$ and $\bV\in \real^{n\times n}$ be orthogonal matrices. Then,
$
\norm{\bA }_F =  \norm{\bU\bA\bV }_F.
$
\end{proposition}
\begin{proof}[of Proposition~\ref{proposition:frobenius-orthogonal-equi}]
We observe that 
$$
\begin{aligned}
\norm{\bU\bA\bV}_F  &=\sqrt{\trace((\bU\bA\bV)(\bU\bA\bV)^\top)} = \sqrt{\trace(\bU\bA\bA^\top\bU^\top)}\\
&= \sqrt{\trace(\bA\bA^\top\bU^\top\bU)} = \sqrt{\trace(\bA\bA^\top)}=
\norm{\bA}_F ,
\end{aligned}
$$
where the third equality holds because the trace is invariant under cyclic permutations.
\end{proof}

\index{Schur inequality}
The Frobenius norm is defined as the square root of the sum of the squares of the matrix elements. Additionally, we have the well-known \textit{Schur inequality} (to differentiate from  Schur's inequality in Theorem~\ref{theorem:schurs_real_inequality}) from the definition of the Frobenius norm.
\begin{theorem}[Schur Inequality]\label{theorem:schur_inequality}
Let $\lambda_1, \lambda_2, \ldots, \lambda_n$ be  real eigenvalues of matrix $\bA\in \real^{n\times n}$. Then,
$
\sum_{i=1}^{n} |\lambda_i|^2 \leq \sum_{i=1,j=1}^{n,n} |a_{ij}|^2 = \norm{\bA}_F^2,
$
i.e., the sum of absolute eigenvalues is bounded by the Frobenius norm.
\end{theorem}
\begin{proof}[of Theorem~\ref{theorem:schur_inequality}]
Suppose the Schur decomposition of $\bA$ is given by $\bA=\bQ\bU\bQ^\top$ (Theorem~\ref{theorem:schur-decomposition}), where the diagonal of $\bU$ contains the eigenvalues of $\bA$ (Corollary~\ref{corollary:schur-second-form}). By the orthogonal invariance Proposition~\ref{proposition:frobenius-orthogonal-equi}, we have
$
\norm{\bU}_F = \norm{\bQ\bU\bQ^\top}_F = \norm{\bA}_F. 
$
Therefore,
$$
\sum_{i=1}^{n} |\lambda_i|^2 = \sum_{i=1}^{n} u_{ii}^2 \leq \sum_{i=1}^{n} u_{ii}^2 + \sum_{ i\neq j} u_{ij}^2 = \norm{\bU}_F^2
\quad\implies \quad
\sum_{i=1}^{n} |\lambda_i|^2  \leq \norm{\bA}_F^2.
$$
This completes the proof.
\end{proof}
\begin{exercise}[Schur Inequality]
Prove the Schur inequality for general matrices that do not necessarily have real eigenvalues. Under what conditions does equality hold?
\end{exercise}

\subsection*{The $\ell_p$ Matrix Norm: Generalizing Frobenius Norm}

Similar to the $\ell_p$ vector norm, we can  define the $\ell_p$ matrix norm as follows.
\begin{definition}[$\ell_p$ Matrix Norm]\label{definition:lp-matrix_norm_app}
Given a  matrix $\bA\in \real^{m\times n}$, for  $p\geq 1$, the $\ell_p$ \textit{matrix norm} is defined as:
$$
\norm{\bA}_{m_p} = \sqrt[p]{ \sum_{j=1}^{n} \sum_{i=1}^{m}|a_{ij}|^p  }.
$$	
From this definition, we have the specific cases of the $\ell_1, \ell_2$, and $\ell_\infty$ matrix norms:
$$
\begin{aligned}
\norm{\bA}_{m_1} &= \sum_{j=1}^{n} \sum_{i=1}^{m} |a_{ij}|; \qquad
\norm{\bA}_{m_2} = \left(\sum_{j=1}^{n} \sum_{i=1}^{m} |a_{ij}|^2 \right)^{1/2}; \\
\norm{\bA}_{m_\infty} &= \mathop{\max}_{i,j} |a_{ij}|, \gap \forall i\in\{1,2,\ldots,m\}, j\in \{1,2,\ldots,n\},\\
\end{aligned}
$$
where the subscript $m_p$ is used to distinguish the  $\ell_p$  matrix norm  from the  $\ell_p$ vector norm and the induced matrix norms (discussed later). When $p=2$, the $\ell_2$ matrix norm reduces to the Frobenius norm.
Since $\norm{\bA}_{m_p} = \norm{vec(\bA)}_p$, it follows from Propositrion~\ref{prop:lp_norm_ineqs} that
\begin{equation}
\norm{\bA}_{m_\infty}\leq \norm{\bA}_{m_q} \leq \norm{\bA}_{m_p}\leq \normtwo{\bA} \leq \norm{\bA}_{m_1},
\quad \text{if }  \infty \geq q \geq p\geq 2\geq 1;
\end{equation}
\begin{equation}
\norm{\bA}_{m_\infty} \leq \norm{\bA}_{m_p} \leq  n^{1/p} \norm{\bA}_{m_\infty},
\quad \text{if }  \infty \geq  p \geq 1.
\end{equation}
\end{definition}

\begin{exercise}[$\ell_1$, $\ell_\infty$ Matrix Norm]
Show that the $\ell_1$ matrix norm is submultiplicative and the $\ell_\infty$ matrix norm is not submultiplicative.
How about the matrix norm $\norm{\bA}=n\norm{\bA}_{m_\infty}$ if $\bA\in\real^{n\times n}$?
\end{exercise}

\subsection*{Spectral Norm}
Another important matrix norm that is extensively used is the \textit{spectral norm}.
\begin{definition}[Spectral Norm]\label{definition:spectral_norm_app}
The spectral norm of a matrix $\bA\in \real^{m\times n}$ is defined as 
$$
\norm{\bA}_2 = \mathop{\max}_{\bx\neq\bzero} \frac{\norm{\bA\bx}_2}{\norm{\bx}_2}  =\mathop{\max}_{\bu\in \real^n: \norm{\bu}_2=1}  \norm{\bA\bu}_2 ,
$$
which is also the largest singular value of $\bA$, i.e., $\norm{\bA}_2 = \sigma_{\max}(\bA)$~\footnote{See Equation~\eqref{equation:rayleig_ritz_sing}. When $\bA$ is an $n\times n$  positive semidefinite matrix, $\norm{\bA}_2=\sqrt{\lambda_{\max}(\bA^2)}=\lambda_{\max}(\bA)$, i.e., the largest eigenvalue of $\bA$.}. The second equality holds because scaling $\bx$ by a nonzero scalar does not change the norm:
$
\frac{\norm{\lambda \cdot \bA\by}_2 }{\norm{\lambda \cdot \by}_2 }  =  \norm{\bA\by}_2.
$
The definition also indicates the matrix-vector inequality: $\norm{\bA\bx}_2 \leq \norm{\bA}_2\norm{\bx}_2$ for all vectors $\bx\in\real^n$.
\end{definition}
To see why the spectral norm of a matrix is equal to its largest singular value, consider the singular value decomposition $\bA=\bU\bSigma\bV^\top$. We have:
$$
\begin{aligned}
\norm{\bU\bSigma\bV^\top}_2 &= \max_{\bx \neq \bzero} \frac{\norm{\bU\bSigma\bV^\top \bx}_2}{\norm{\bx}_2} \stackrel{*}{=} \max_{\bV^\top\bx \neq \bzero} \frac{\norm{\bSigma\bV^\top \bx}_2}{\norm{\bx}_2}\\
&= \max_{\bV^\top\bx \neq \bzero} \frac{\norm{\bSigma\bV^\top \bx}_2}{\norm{\bV^\top\bx}_2} \frac{\norm{\bV^\top\bx}_2}{\norm{\bx}_2} 
=\max_{\by \neq \bzero} \frac{\norm{\bSigma \by}_2}{\norm{\by}_2}
\leq \sigma_{\max}(\bA),
\end{aligned}
$$
where the equality ($*$) holds since $\bV$ is orthogonal and the $\ell_2$ vector norm is orthogonally invariant. The inequality holds because the largest singular value of $\bA$ is the maximum value of $\frac{\norm{\bSigma \by}_2}{\norm{\by}_2}$ over all nonzero vectors $\by$.
Alternatively, this can be shown by noting that: $\norm{\bA\bx}_2^2=\abs{\bx^\top\bA^\top\bA\bx}\leq \sigma_{\max}^2(\bA)$.

\index{Submultiplicativity}
\begin{proposition}[Submultiplicativity of Spectral]\label{propo:submul_spec}
The spectral norm is submultiplicative. That is, $\norm{\bA\bB}_2\leq \norm{\bA}_2\norm{\bB}_2$.
\end{proposition}
\begin{proof}[of Proposition~\ref{propo:submul_spec}]
The definition of the spectral norm shows that
$$
\begin{aligned}
\norm{\bA\bB}_2 &= \max_{\bx \neq \bzero} \frac{\norm{\bA\bB \bx}_2}{\norm{\bx}_2} \stackrel{*}{=} \max_{\bB\bx \neq \bzero} \frac{\norm{\bA\bB \bx}_2}{\norm{\bx}_2} 
= \max_{\bB\bx \neq \bzero} \frac{\norm{\bA\bB \bx}_2}{\norm{\bB\bx}_2} \frac{\norm{\bB\bx}_2}{\norm{\bx}_2} \\
&\leq \max_{\by \neq \bzero} \frac{\norm{\bA \by}_2}{\norm{\by}_2} \max_{\bx \neq \bzero} \frac{\norm{\bB\bx}_2}{\norm{\bx}_2}
=\norm{\bA}_2\norm{\bB}_2,
\end{aligned}
$$
where the equality ($*$) holds because if the maximum is obtained when $\bx\neq \bzero$ and $\bB\bx=\bzero$, the norm is $\norm{\bA\bB}=\bzero$. This holds only when $\bA\bB=\bzero$, and the submultiplicativity holds obviously.

Alternatively, we have 
$
\norm{\bA\bB}_2 = \max_{\norm{\bu}=1}\norm{\bA\bB \bu}_2,
$
where $\norm{\bA\bB \bu}_2\leq \norm{\bA}_2\norm{\bB\bu}_2\leq\norm{\bA}_2\norm{\bB}_2 $. This again indicates $\norm{\bA\bB}_2\leq \norm{\bA}_2\norm{\bB}_2 $ implicitly.
\end{proof}

\index{Orthogonally invariance}
\begin{proposition}[Orthogonally Invariance of Spectral]\label{proposition:submul_spec}
Let $\bA\in \real^{m\times n}$ be given, and let $\bU\in \real^{m\times m}$ and $\bV\in \real^{n\times n}$ be orthogonal matrices. Then,
$
\norm{\bA }_2 =  \norm{\bU\bA\bV }_2.
$
\end{proposition}
\begin{proof}[of Proposition~\ref{proposition:submul_spec}]
We notice that 
$$
\begin{aligned}
\normtwo{\bU\bA\bV}  &=\mathop{\max}_{\bx\neq\bzero} \frac{\norm{\bU\bA\bV\bx}_2}{\norm{\bx}_2}
=\mathop{\max}_{\bx\neq\bzero} \frac{\norm{\bA\bV\bx}_2}{\norm{\bx}_2}
=\mathop{\max}_{\bx\neq\bzero} \frac{\norm{\bA\bV\bx}_2}{\norm{\bV\bx}_2}
=\mathop{\max}_{\by\neq\bzero} \frac{\norm{\bA\by}_2}{\norm{\by}_2}
=\norm{\bA}_2 ,
\end{aligned}
$$
This completes the proof.
\end{proof}

We conclude that the spectral norm is a normalized matrix norm.
\begin{proposition}[Normalized Spectral]\label{proposition:norma_spec}
The spectral norm is normalized such that:
\begin{enumerate}
\item \textit{Normalization}. When $\bA\neq \bzero$, we have $\norm{\frac{1}{\norm{\bA}_2} \bA}_2=1$.
\item \textit{Normalized}. $\norm{\bI}_2=1$.
\end{enumerate}
\end{proposition}

\subsection*{Induced Matrix Norm: Generalizing Spectral Norm}
More generally, many matrix norms can be generated  using the concept of induced norms.

\index{Induced norm}
\begin{definition}[Induced Matrix Norm: General Matrix Norm]\label{definition:induced_norm_app}
 Given a matrix $\bA\in \real^{m\times n}$, and two vector norms $\norm{\cdot}_a$ and $\norm{\cdot}_b$ on $\real^n$ and $\real^m$ respectively, the \textit{induced matrix norm} $\norm{\bA}_{a,b}$ is defined by 
$$
\norm{\bA}_{a,b} =  \mathop{\max}_{\norm{\bx}_a \neq  0 } \frac{\norm{\bA\bx}_b}{\norm{\bx}_a}  
=  \mathop{\max}_{\norm{\bx}_a = 1} \norm{\bA\bx}_b .~\footnote{Although the norms  $\norm{\cdot}_a$ and $\norm{\cdot}_b$ can be chosen from a variety of options, $\ell_p$ norms are most commonly discussed in this context. Therefore, we may also refer to it as the  \textit{\holders induced matrix norm}.}
$$
The \textit{matrix-vector product inequality} can be easily derived from the above definition:
\begin{equation}\label{equation:induced_ineqy_intern}
\norm{\bA\bx}_b \leq \norm{\bA}_{a,b} \cdot \norm{\bx}_a.
\end{equation}
The induced matrix norm can also be referred to as the $(a,b)$-norm. When $a=b$, we simply call it an \textit{$a$-norm} and use the notation $\norm{\bA}_{a}$ instead of $\norm{\bA}_{a,a}$. When $a=b=2$, we have the classic spectral norm. The induced norm is defined as the maximum factor by which $\bA$ magnifies the length of nonzero vectors, where the length of the vector $\bx$ is measured with norm $\norm{\cdot}_a$ and the length of the transformed vector $\bA\bx$ is measured with norm $\norm{\cdot}_b$. 
\end{definition}

From the definition of the induced norm, we can see that the spectral norm is a special induced norm when $a=b=2$ such that $\norm{\bA}_2 = \norm{\bA}_{2,2}$. Similarly, \textit{matrix 1-norm} can be obtained by
\begin{equation}\label{equation:mat_one_norm}
\textbf{Matrix 1-norm: }\norm{\bA}_1 = \mathop{\max}_{j=1,2,\ldots,n} \sum_{i=1}^{m}|a_{ij}| 
=\mathop{\max}_{j=1,2,\ldots,n} \norm{\ba_j}_1
,
\end{equation}
which is also called the \textit{maximum absolute column sum norm}. 
To see this, given $\bA\in\real^{m\times n}$, we have 
$$
\norm{\bA\bx}_1=\norm{\sum_{i=1}^{n} x_i\ba_i}_1
\leq 
\sum_{i=1}^{n} \abs{x_i}\norm{\ba_i}_1
\leq 
\sum_{i=1}^{n} \abs{x_i} \left(\mathop{\max}_{j=1,2,\ldots,n}\norm{\ba_j}_1 \right)
=
\norm{\bx}_1 \norm{\bA}_1.
$$
Therefore, $ \mathop{\max}_{\norm{\bx}_1 = 1} \norm{\bA\bx}_1 \leq \norm{\bA}_1$.
On the other hand, consider $\mathop{\max}_{\norm{\bx}_1 = 1}\norm{\bA\bx}_1$ for $\bx=\be_j,\,\forall j\in\{1,2,\ldots,n\}$, we have
$
\mathop{\max}_{\norm{\bx}_1 = 1}\norm{\bA\bx}_1\geq \mathopmax{j=1,2,\ldots,n}\norm{\ba_j}_1=\norm{\bA}_1.
$
By ``sandwiching", we have $\mathop{\max}_{\norm{\bx}_1 = 1}\norm{\bA\bx}_1=\norm{\bA}_1$.
Similarly, the \textit{matrix $\infty$-norm} can be obtained by 
\begin{equation}\label{equation:mat_inf_norm}
\textbf{Matrix $\infty$-norm: }\norm{\bA}_\infty = \mathop{\max}_{i=1,2,\ldots,m} \sum_{j=1}^{n}|a_{ij}| ,
\end{equation}
which is also called the \textit{maximum absolute row sum norm}.
\begin{exercise}[Maximum Absolute Row Sum Norm]
Prove that $\norm{\bA}_\infty$ is induced from the $\ell_\infty$ vector norm.
\end{exercise}
\index{Equi-induced matrix norm}
\index{Normalized matrix norm}
\begin{exercise}[Equi-Induced and Normalized Norm]
Let $a=b$ and let $\bA\in\real^{n\times n}$ be square in Definition~\ref{definition:induced_norm_app}. Then the induced norm is also called a \textbf{equi-induced norm}. Show that the equi-induced norm is normalized: $\norm{\bI_n}=1$ and $\norm{\bA}\geq 1$.
\end{exercise}

\index{Submultiplicativity}
We have shown that the spectral norm, which is a special induced norm, is submultiplicative.
In fact, all induced norms satisfy the submultiplicative property.
\begin{theorem}[Submultiplicativity of Induced Norms]\label{theorem:Submultiplicativity_induced}
For $a,b,c\geq 1$ and matrices $\bA,\bB$ (provided the matrix product $\bA\bB$ is defined), we have 
$$
\norm{\bA\bB}_{a,b} \leq \norm{\bA}_{c,b}\norm{\bB}_{a,c}.
$$
That is, the norm of a
product is bounded by the product of the norms.
\end{theorem}
\begin{proof}[of Theorem~\ref{theorem:Submultiplicativity_induced}]
From the definition of the induced norm, we have
$$
\begin{aligned}
\norm{\bA\bB}_{a,b} &= \mathop{\max}_{\norm{\bx}_a=1} \norm{\bA\bB\bx}_b
\leq \mathop{\max}_{\norm{\bx}_a=1}  \norm{\bA}_{c, b} \norm{\bB\bc}_{c} 
\leq \mathop{\max}_{\norm{\bx}_a=1}\norm{\bA}_{c,b} \norm{\bB}_{a,c} \norm{\bx}_a .
\end{aligned}
$$
This completes the proof.
\end{proof}

\begin{remark}[Condition Number]
For a nonsingular $\bA\in\real^{n\times n}$,
we note that $1 = \norm{\bI}_{b,b} \leq \norm{\bA}_{a,b} \norm{\bA^{-1}}_{b,a}$. The value $\kappa_{a,b}=\norm{\bA}_{a,b} \norm{\bA^{-1}}_{b,a}$ is normally known as the \textit{condition number} with respect to the induced vector norms $\norm{\cdot}_a$ and $\norm{\cdot}_b$. See Appendix~\ref{appendix:condition_number}  for a further discussion.
\end{remark}

The following properties of the induced norm are useful for evaluating the theory of condition numbers  (Appendix~\ref{appendix:condition_number}).
\begin{lemma}[Induced Norm Properties]\label{lemma:induced_norm_property}
\begin{enumerate}
\item
Given $\by\in \real^m$, $\bx\in \real^n$, and $a,b\geq 1$ with $\frac{1}{a}+\frac{1}{a^\star}=1$ and $\frac{1}{b}+\frac{1}{b^\star}=1$, it follows that 
$$
\norm{\by\bx^\top}_{a,b} = \norm{\by}_b \norm{\bx}_{a^\star}.
$$
\item 
Given $\by\in \real^m$ and $\bx\in \real^n$ with $\norm{\bx}_a =\norm{\by}_b=1$, then there exists a matrix $\bA$ such that $\norm{\bA}_{a,b}=1$ and $\bA\bx=\by$.
\item 
Particularly, we have $\norm{\bA^\top}_{a,b} = \norm{\bA}_{b^\star, a^\star}$.
\end{enumerate}
\end{lemma}
\begin{proof}[of Lemma~\ref{lemma:induced_norm_property}]
For (1), from the definition of the induced norm, we have 
$$
\norm{\by\bx^\top}_{a,b} =\mathop{\max}_{\norm{\bz}_a= 1} \norm{\by\bx^\top\bz}_b
= \norm{\by}_b \mathop{\max}_{\norm{\bz}_a= 1} |\bx^\top\bz| =  \norm{\by}_b  \norm{\bx}_{a^\star},
$$
where the last equality follows from Equation~\eqref{equation:dual_norm_equa}.

For (2), again from Equation~\eqref{equation:dual_norm_equa}, there exists  a vector $\bz\in \real^n$ such that $\norm{\bz}_{a^\star}=1$ and $\bz^\top\bx=1$. Therefore, $\bA=\by\bz^\top$ satisfies $\bA\bx=\by$. And from (1), we have $\norm{\bA}_{a,b}=\norm{\by}_b \norm{\bz}_{a^\star}=1.$ 

For (3), again from the definition of the induced norm, we have 
$$
\begin{aligned}
\norm{\bA^\top}_{a,b} &= \mathop{\max}_{\norm{\bx}_a=1} \norm{\bA^\top\bx}_b
= \mathop{\max}_{\norm{\bx}_a=1} \mathop{\max}_{\norm{\bz}_{b^\star}=1} \bx^\top \bA\bz
=\mathop{\max}_{\norm{\bz}_{b^\star}} \norm{\bA\bz}_{a^\star} 
= \norm{\bA}_{b^\star, a^\star},
\end{aligned}
$$
where the second and third equalities follow from Equation~\eqref{equation:dual_norm_equa}.
This completes the proof.
\end{proof}

Similarly to the vector norms in Theorem~\ref{theorem:equivalence-vector-norm},  matrix norms also admit an equivalence statement.
\begin{theorem}[Equivalence of Matrix Norms]\label{theorem:equiv_mat_norm}
Let $\norm{\cdot}_a$ and $\norm{\cdot}_b$ be two different matrix norms: $\real^{m\times n}\rightarrow \real$. Then there exist positive scalars $\alpha$ and $\beta$ such that for all $\bA \in \real^{m\times n}$,
$$
\alpha\norm{\bA}_a \leq \norm{\bA}_b \leq \beta\norm{\bA}_a.
$$
\end{theorem}
The proof is similar to that of Theorem~\ref{theorem:equivalence-vector-norm}, and we shall not repeat the details. 
The equivalence theorem again states that if a matrix is small in one norm, it is also small in other norms, and vice versa.

We have shown in Lemma~\ref{lemma:construct_norm} how to construct a norm from another norm.
The induced norm can also be constructed using another induced norm (when the underlying matrix is square).
\begin{proposition}[Construct Induced Norms from Induced Norms]\label{propo:cons_ind}
Let $\norm{\cdot}$ be a matrix norm (resp. a submultiplicative matrix norm) on $\real^{n\times n}$ and $\bS\in\real^{n\times n}$ be nonsingular. Then, the following function on matrix $\bA\in\real^{n\times n}$ is also a matrix norm (resp. a submultiplicative matrix norm):
$$
\norm{\bA}_{\bS} = \norm{\bS\bA\bS^{-1}}.
$$
Moreover, if $\norm{\cdot }$ is induced from the vector norm $\norm{\cdot}_v$ on $\real^n$, the matrix norm $\norm{\bA}_{\bS}$ is induced from the vector norm $\norm{\bx}_{\bS}=\norm{\bS\bx}_v$ on $\real^n$ (the $\bS$-norm in \eqref{equation:s_norm} with nonsingular $\bS$).
\end{proposition}
\begin{proof}[of Proposition~\ref{propo:cons_ind}]
The proof for the first part is left as an exercise.
And we note that the submultiplicativity of $\norm{\bA}_{\bS}$ follows from the submultiplicativity of $\norm{\bA}$.
For the second part, we have
$$
\mathopmax{\norm{\bx}_{\bS}=1} \norm{\bA\bx}_{\bS}
=
\mathopmax{\norm{{\bS}\bx}=1} \norm{\bS\bA\bx}_v
=
\mathopmax{\norm{\by}=1} \norm{\bS\bA\bS^{-1}\by}_v
=
\norm{\bS\bA\bS^{-1}},
$$
from which the result follows.
\end{proof}

\index{Subordinate}
\subsection*{Subordinate Property}

We note in the induced norm that $	\norm{\bA\bx}_b \leq \norm{\bB}_{a,b} \cdot \norm{\bx}_a$ (Equation~\eqref{equation:induced_ineqy_intern}). 
Recall that the matrix norm measures how much a vector is ``stretched" in this context. 
This property is known as the \textit{subordinate property} of a matrix norm with respect to given vector norms.
\begin{definition}[Subordinate Matrix Norm]\label{definition:subordinate_matrix_norm}
A matrix norm $\norm{\cdot}: \real^{m\times n} \rightarrow \real$ is said to be subordinate to vector norms $\norm{\cdot}_a: \real^n\rightarrow \real$ and  $\norm{\cdot}_b: \real^m\rightarrow \real$ if for all vectors $\bx\in\real^n$ and matrices $\bA\in \real^{m\times n}$,
$$
	\norm{\bA\bx}_b \leq \norm{\bA} \cdot \norm{\bx}_a.
$$
\end{definition}

The following corollary  follows immediately from the definition of the induced norm.
\begin{corollary}[Subordinate Induced Norm]
The induced matrix norm $\norm{\cdot}_{a,b}$ is subordinate to the vector norms $\norm{\cdot}_a$ and $\norm{\cdot}_b$ that induce it.
\end{corollary}

\begin{corollary}[Subordinate Frobenius Norm]\label{corollary:subordinate_frobenius}
The Frobenius norm is subordinate to the $\ell_2$ vector norm.
\end{corollary}
\begin{proof}[of Corollary~\ref{corollary:subordinate_frobenius}]
Without loss of generality, we assume $\bx\neq \bzero \in\real^n $. 
Given a matrix $\bA\in \real^{m\times n}$, we have 
$$
\begin{aligned}
\norm{\bA\bx}_2 = \frac{\norm{\bA\bx}_2}{\norm{\bx}_2} \norm{\bx}_2 \leq
\mathop{\max }_{\by\neq \bzero } \frac{\norm{\bA\by}_2}{\norm{\by}_2} \norm{\bx}_2 
= \mathop{\max }_{\norm{\by}_2=1 }  \norm{\bA\by}_2 \norm{\bx}_2   = 
\norm{\bA}_2 \norm{\bx}_2 .
\end{aligned}
$$
Since $\norm{\bA}_F =\sqrt{\sigma_1^2+\sigma_2^2+\ldots+\sigma_r^2}$ and $\norm{\bA}_2 = \sigma_1$ ($\sigma_1$ is the largest singular value), we also have $\norm{\bA}_2 \leq \norm{\bA}_F$. Combining the two results, the claim follows.
\end{proof}


We provide an inequality involving subordinate norms, which will be useful in the analysis of condition numbers.
To see this, we  need to use the spectral lemma of norms (given submultiplicativity). We prove in Theorem~\ref{theorem:schur_inequality} that the sum of squared eigenvalues is bounded by the Frobenius norm. The spectral radius (Definition~\ref{definition:spectrum}) of any square matrix is bounded for any (submultiplicative) matrix norm (Equation~\eqref{equation:bounds_submulnorm}).

\begin{proposition}[Subordinate Inequality]\label{proposition:subordinate_ineq1}
Let $\norm{\cdot}:\real^{n\times n}\rightarrow \real$ be any submultiplicative induced matrix norm. Then,
\begin{enumerate}
\item If $\norm{\cdot}$ is a subordinate matrix norm,  and  $\bA\in\real^{n\times n}$ is any square matrix satisfying $\norm{\bA}<1$, then the matrix $\bI+\bA$ is invertible and 
$$
\frac{1}{1+\norm{\bA}}
\leq 
\norm{(\bI+\bA)^{-1}}
\leq \frac{\norm{\bI}}{1-\norm{\bA}}.
$$
When $\norm{\cdot}$ is normalized, $\norm{\bI}=1$ in the right inequality~\footnote{Actually, the only subordinate unitarily invariant norm is the spectral norm, which is normalized: $\normtwo{\bI}=1$.
}.
\item If  the matrix $\bI+\bA$ is singular, then $\norm{\bA}\geq 1$ for every matrix norm (nor necessarily subordinate).
\end{enumerate}
\end{proposition}
\begin{proof}[of Proposition~\ref{proposition:subordinate_ineq1}]
For (1), since $(\bI+\bA)\bx = \bzero$ implies $\bA\bx = -\bx$, we have $\Vert\bx\Vert = \Vert\bA\bx\Vert$.
By the subordinate property, if $(\bI+\bA)\bx = \bzero$, we get 
$$
\Vert\bx\Vert = \Vert\bA\bx\Vert \leq 
\norm{\bA} \norm{\bx}.
$$
Since $\Vert\bA\Vert<1$, if $\bx\neq \bzero$, we have 
$
\Vert\bA\bx\Vert
< \Vert\bx\Vert.
$
This leads to a contradiction to $\Vert\bx\Vert = \Vert\bA\bx\Vert$.
Therefore, $\bx$ must be $\bzero$. Hence, the null space of  $\bI+\bA$ is of dimension 0, and  $\bI+\bA$ is nonsingular. Similarly, we can also show that $\bI-\bA$ is nonsingular.
Then we get
$$
(\bI+\bA)^{-1} +\bA(\bI+\bA)^{-1} = (\bI+\bA)(\bI+\bA)^{-1} = \bI
\quad\implies\quad 
(\bI+\bA)^{-1} = \bI - \bA(\bI+\bA)^{-1}.
$$
By the triangle inequality and submultiplicative property, we have 
$$
\Vert (\bI+\bA)^{-1}\Vert 
\leq 
\norm{\bI}+ \Vert\bA\Vert \left\Vert (\bI+\bA)^{-1}\right\Vert
\quad\implies\quad
\norm{(\bI+\bA)^{-1}}
\leq \frac{\norm{\bI}}{1-\Vert\bA\Vert}.
$$
On the other hand, since $\norm{\cdot}$ is submultiplicative, we have $\norm{\bI}\geq 1$ (Equation~\eqref{equation:power_subm2}). Since $\bI-\bA$ is nonsingular, we obtain the left inequality: 
$$
1\leq \norm{\bI}
=\norm{(\bI-\bA)(\bI-\bA)^{-1}}
\leq 
\norm{(\bI-\bA)}\norm{(\bI-\bA)^{-1}}
\leq 
(1+\norm{\bA})\norm{(\bI-\bA)^{-1}}.
$$

For (2), if the norm is not subordinate and $\bI+\bA$ is singular, then $-1$ is an eigenvalue of $\bA$. By Equation~\eqref{equation:bounds_submulnorm}, it follows that 
$$
|-1|\leq \rho(\bA) \leq \Vert\bA\Vert,
$$
from which the result follows.
\end{proof}

%% file: chapter-app_proj.tex
\newpage
\section{(Orthogonal) Projection and  Geometric Interpretation of LS}\label{section:by-geometry-hat-matrix}
The least squares problem is addressed in Section~\ref{section:application-ls-qr}, Section~\ref{section:ls-utv}, and Section~\ref{section:application-ls-svd}  using various matrix decompositions.
As discussed previously, the \textit{ordinary least squares (OLS) estimate} or least squares (LS) estimate aims to minimize $\normtwo{\bb-\bA\bx}^2$, which searches for an estimate $\bx_{LS}$ such that $\bA\bx_{LS}$ is in the column space of $\bA$, $\cspace(\bA)$, thereby minimizing  the distance between $\bA\bx_{LS}$ and $\bb$. 
The predicted value $\bp = \bA\bx_{LS}$ is the (orthogonal) projection of $\bb$ onto column space $\cspace(\bA)$, achieved through a \textit{projection matrix} $\bH = \bA(\bA^\top\bA)^{-1}\bA^\top$:\index{Projection matrix}\index{Orthogonal projection}
\begin{equation}
\bp = \bA\bx_{LS} = \bH\bb, \nonumber
\end{equation}
where the matrix $\bH$ is also called the \textit{hat matrix}  due to its role in ``putting a hat" on the output ($\bp$ is sometimes denoted as $\hat{\bb}$ in literature). 

However, what exactly is a projection matrix? Simply stating that $\bH = \bA(\bA^\top\bA)^{-1}\bA^\top$ is a projection matrix is insufficient without explanation. 
Before exploring the properties of the projection matrix, let us review key characteristics of symmetric and idempotent matrices, which will be crucial in the subsequent sections.

\subsection{Properties of Symmetric and Idempotent Matrices}

Symmetric idempotent matrices have specific eigenvalues, which will be often used.
\begin{lemma}[Eigenvalue of Symmetric Idempotent Matrices]\label{lemma:eigenvalues-of-projection}
The only possible eigenvalues of any symmetric idempotent matrix are 0 and 1.
\end{lemma}
In Lemma~\ref{proposition:eigenvalues-of-projection2}, we show that the eigenvalues of idempotent matrices are 1 or 0 as well, which relaxes the conditions required here (both idempotent and symmetric). However, the method used in the proof is quite useful,  so we present both results. To prove the lemma above, we need to use the result of the spectral theorem. 

\begin{proof}[of Lemma~\ref{lemma:eigenvalues-of-projection}]
Suppose matrix $\bA$ is symmetric idempotent. By the spectral theorem (Theorem~\ref{theorem:spectral_theorem}), we can decompose $\bA = \bQ \bLambda\bQ^\top$, where $\bQ$ is an orthogonal matrix and $\bLambda$ is a diagonal matrix. Therefore,
\begin{equation}
\begin{aligned}
(\bQ\bLambda\bQ^\top)^2 = \bQ\bLambda\bQ^\top
\quad\implies \quad
\bQ\bLambda^2\bQ^\top = \bQ\bLambda\bQ^\top 
\quad\implies \quad
\bLambda^2 = \bLambda 
\quad\implies \quad
\lambda_i^2 = \lambda_i \nonumber.
\end{aligned}
\end{equation}
Thus, the eigenvalues of $\bA$ must satisfy  $\lambda_i \in \{0,1\},\,\, \forall i$. This completes the proof.
\end{proof}

In the above lemma, we use the spectral theorem to demonstrate that the only eigenvalues of any symmetric idempotent matrix are 0 and 1. 
This technique is commonly employed in mathematical proofs, particularly in distribution theory (see \citet{lu2021rigorous}). 
With simple modifications, we can extend this result to idempotent matrices that are not necessarily symmetric.

\begin{lemma}[Eigenvalue of Idempotent Matrices\index{Idempotent matrices}]\label{proposition:eigenvalues-of-projection2}
The only possible eigenvalues of any idempotent matrix are 0 and 1.
\end{lemma}
\begin{proof}[of Lemma~\ref{proposition:eigenvalues-of-projection2}]
Let $\bx$ be an eigenvector of the idempotent matrix $\bA$ corresponding to eigenvalue $\lambda$. That is,
$
\bA \bx = \lambda \bx.
$
Furthermore, we observe that
$$
\begin{aligned}
\bA^2\bx &=(\bA^2)\bx =\bA\bx = \lambda\bx
=\bA(\bA\bx) = \bA(\lambda\bx)=\lambda\bA\bx=\lambda^2\bx.
\end{aligned}
$$
Since $\bx\neq \bzero$, this implies $\lambda^2 = \lambda$ and $\lambda$ is either 0 or 1.
\end{proof}

We also prove the rank of a symmetric idempotent matrix is equal to its trace, which will be extremely useful in the following sections.
\begin{lemma}[Rank and Trace of Symmetric Idempotent Matrices\index{Trace}]\label{lemma:rank-of-symmetric-idempotent}
Let $\bA$ be  any $n\times n$ symmetric idempotent matrix. Then, the rank of $\bA$ equals the trace of $\bA$.
\end{lemma}
\begin{proof}[of Lemma~\ref{lemma:rank-of-symmetric-idempotent}]
By the spectral theorem~\ref{theorem:spectral_theorem}, we have spectral decomposition of $\bA$ as $\bA = \bQ\bLambda\bQ^\top$. 
Since $\bA$ and $\bLambda$ are similar matrices, they share the same rank and trace (Proposition~\ref{proposition:eigenvalue-similar-matrices}). That is,
$$
\begin{aligned}
\rank(\bA) &= \rank(\diag(\lambda_1, \lambda_2, \ldots, \lambda_n))
\quad\text{and}\quad
\trace(\bA) = \trace(\diag(\lambda_1, \lambda_2, \ldots, \lambda_n)).\\
\end{aligned}
$$
By Lemma~\ref{lemma:eigenvalues-of-projection}, the only eigenvalues of $\bA$ are 0 and 1. Therefore,
$\rank(\bA) = \trace(\bA)$.
\end{proof}

In the above lemma, we prove the rank and trace of any symmetric idempotent matrix are the same. 
However, this condition is somewhat restrictive. We now show that the same result holds for any idempotent matrix, not just symmetric and idempotent ones. Although this lemma is a more general version, we provide both proofs because the methods used in the proof  are valuable for subsequent discussions.

\begin{lemma}[Rank and Trace of an Idempotent Matrices]\label{lemma:rank-of-symmetric-idempotent2}
Let $\bA\in\real^{n\times n}$ be idempotent, the rank of $\bA$ equals its trace: $\rank(\bA)=\rank(\bA^2)=\trace(\bA)$.
\end{lemma}
\begin{proof}[of Lemma~\ref{lemma:rank-of-symmetric-idempotent2}]
Any $n\times n$  rank-$r$ matrix $\bA$ admits the CR decomposition $\bA = \bC\bR$ (Theorem~\ref{theorem:cr-decomposition}), where $\bC\in\real^{n\times r}$ and $\bR\in \real^{r\times n}$ have full rank $r$.
Then, 
$$
\begin{aligned}
\bA^2 &= \bA 
\quad\implies \quad
\bC\bR\bC\bR = \bC\bR 
\quad\implies \quad
\bR\bC\bR =\bR 
\quad\implies \quad
\bR\bC =\bI_r,
\end{aligned}
$$ 
where $\bI_r$ is an $r\times r$ identity matrix. Thus,
$
\trace(\bA) = \trace(\bC\bR) =\trace(\bR\bC) = \trace(\bI_r) = r, 
$
which equals the rank of $\bA$. 
The second equality follows from the fact that the trace is invariant under cyclic permutations.
\end{proof}

\subsection{Orthogonal Projection and Geometric Interpretation for LS}\label{appendix:orthogonal}

Formally, we define the projection matrix as follows:

\begin{definition}[Projection Matrix]\label{definition:projection-matrix}
A (resp. symmetric) matrix $\bH\in\real^{m\times m}$ is called a (resp. orthogonal) projection matrix (a.k.a., a projector) onto a subspace $\mathcalV \in \real^m$ if and only if $\bH$ satisfies the following properties:

\item (P1). $\bH\bb \in \mathcalV$ for all $\bb \in \real^m$: any vector can be projected onto  the subspace $\mathcalV$;

\item (P2). $\bH\bb =\bb$ for all $\bb \in \mathcalV$: projecting a vector already in that subspace has no effect;

\item (P3). $\bH^2 = \bH$, i.e., projecting twice is the same as projecting once because we are already in that subspace, i.e., $\bH$ is idempotent.
\end{definition}

Since we project a vector in $\real^m$ onto a subspace of $\real^m$,  any projection matrix is a square matrix. Otherwise, we would project onto a subspace of $\real^k$ instead of  $\real^m$.
We realize that $\bH\bb$ is always in the column space of $\bH$, and we are interested in the relationship between $\mathcalV$ and $\cspace(\bH)$ (the column space of $\bH$). 
In fact, the column space of $\bH$ is equal to the subspace $\mathcalV$ we want to project onto. Suppose $\mathcalV=\cspace(\bH)$, and suppose further that $\bb$ is already in the subspace $\mathcalV=\cspace(\bH)$, i.e., there exists a vector $\balpha$ such that $\bb=\bH\balpha$. Given only the  condition $(P3)$ above, we have,
$$
\bH\bb = \bH\bH\balpha = \bH\balpha = \bb.
$$
That is, the condition $(P3)$ implies conditions $(P1)$ and $(P2)$. 
Therefore, the definition of a projection matrix can be simplified to only require the idempotent property $(P3)$.

Intuitively, we also want the projection $\bp=\bH\bb$ of any vector $\bb$ to be perpendicular to $\bb - \bp$ such that the distance between $\bp$ and $\bb$ is minimized, resulting in the least squared error (in the least squares problem). 
This is called the \textit{orthogonal projection}.

\begin{definition}[Orthogonal Projection and Oblique Projection Matrix]\label{definition:orthogonal-projection-matrix}
A matrix $\bH$ is called an \textbf{orthogonal projection} matrix onto a subspace $\mathcalV \in \real^m$ if and only if $\bH$ is a projection matrix, and the projection $\bp$ of any vector $\bb\in \real^m$ is perpendicular to $\bb - \bp$, i.e., $\bH$ projects onto $\mathcalV$ and along $\mathcalV^\perp$.
Otherwise, if $\bp$ is not perpendicular to $\bb-\bp$, the projection matrix is called an \textbf{oblique projection} matrix.
\end{definition}

\begin{figure}[h!]
\centering
\vspace{-0.35cm}
\subfigtopskip=2pt
\subfigbottomskip=2pt
\subfigcapskip=-5pt
\subfigure[Orthogonal projection: project $\bb$ to $\bp$.]{\label{fig:ls-geometric1}
\includegraphics[width=0.47\linewidth]{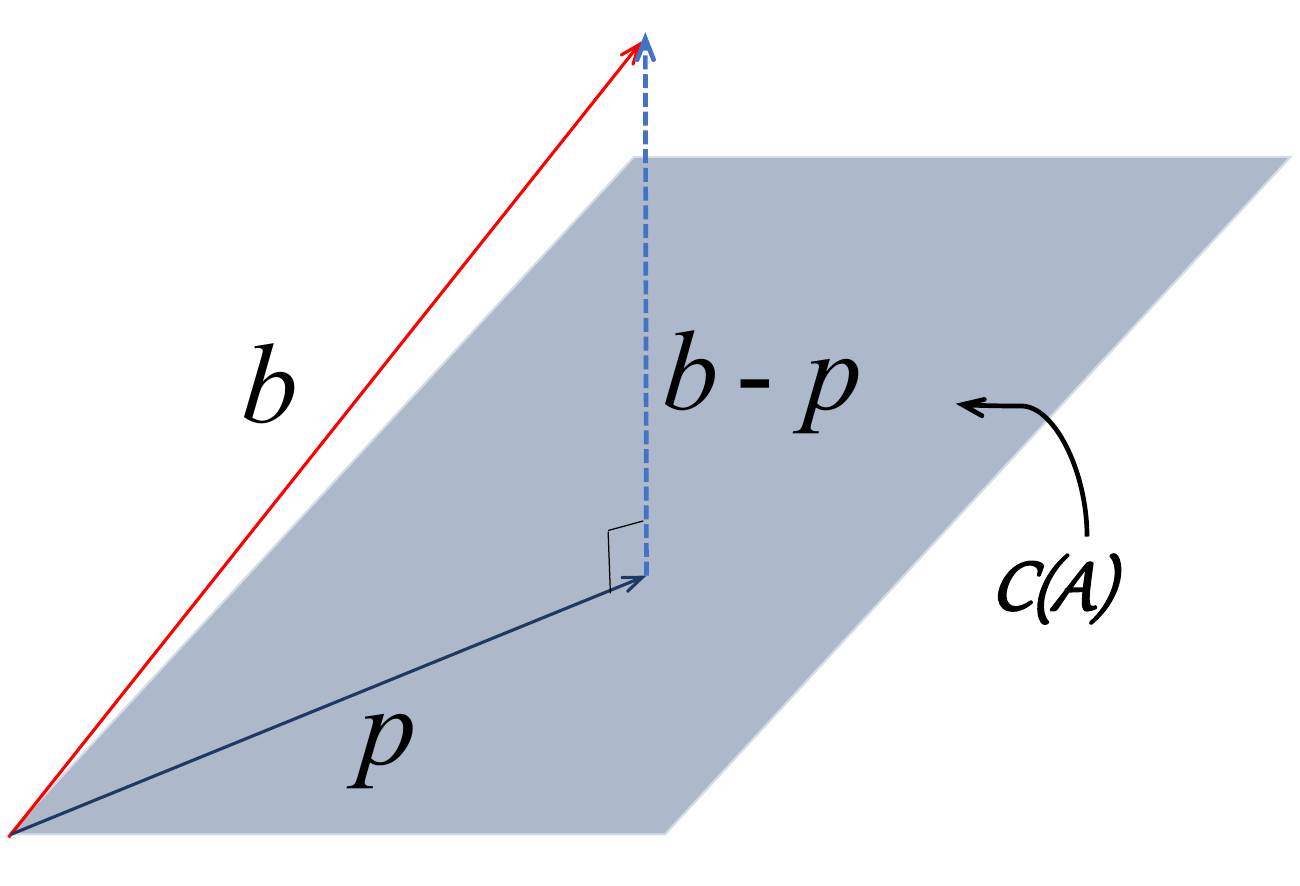}}
\quad 
\subfigure[Oblique projection: project $\bb$ to $\bp_1$ or $\bp_2$.]{\label{fig:ls-geometric1-oblique}
\includegraphics[width=0.47\linewidth]{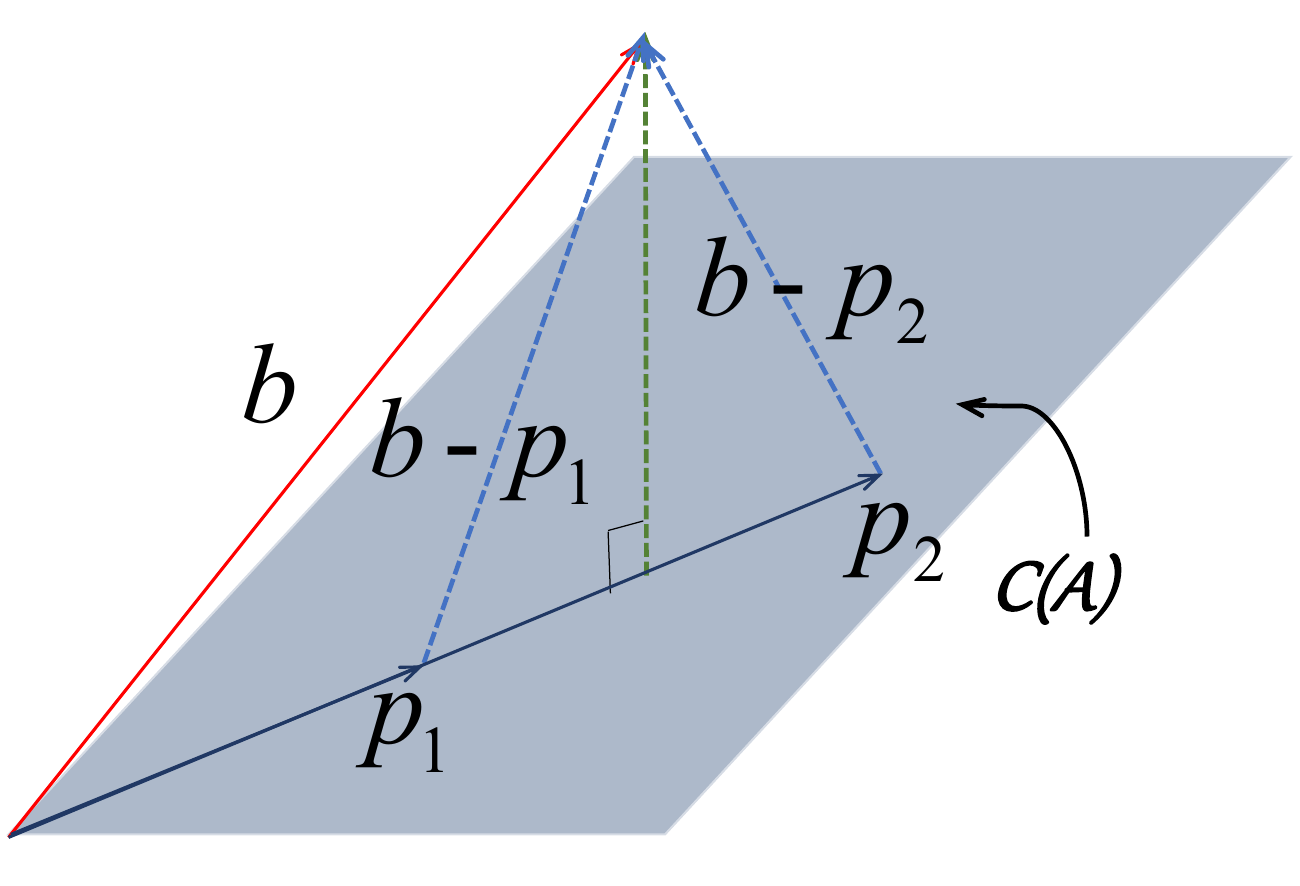}}
\caption{Projection onto the hyperplane of $\cspace(\bA)$, demonstrating orthogonal and oblique projections.}
\label{fig:ls-geometric1-compare}
\end{figure}

\begin{remark}[Orthogonal Projections vs Orthogonal Matrices]
Note that the \textit{orthogonal projection matrix} does not mean that the projection matrix is an orthogonal matrix (discussed in Section~\ref{section:orthogonal-orthonormal-qr}). 
Instead, it means that the projection $\bp$ is perpendicular to $\bb- \bp$. 
This type of projection matrix is so common that we will implicitly assume a projection is orthogonal unless otherwise specified. To avoid confusion, the term \textit{oblique projection matrix} is used for non-orthogonal cases, as illustrated in Figure~\ref{fig:ls-geometric1-compare}.
\end{remark}

\begin{proposition}[Symmetric Orthogonal Projection Matrix]\label{proposition:symmetric-projection-matrix}
Let $\bH$ be any (oblique) projection matrix. Then, $\bH$ is an orthogonal projection matrix if and only if  $\bH$ is symmetric.
\end{proposition}
\begin{proof}[of Proposition~\ref{proposition:symmetric-projection-matrix}]
First, suppose $\bH$ is an orthogonal projection matrix, which projects vectors onto a subspace $\mathcalV$. Then, any vectors $\bv$ and $\bw$ can be decomposed into components lying  in $\mathcalV$ ($\bv_p$ and $\bw_p$) and components lying in $\mathcalV^\perp$ ($\bv_n$ and $\bw_n$) such that
$$
\begin{aligned}
\bv &= \bv_p + \bv_n
\quad \text{and}\quad
\bw = \bw_p + \bw_n.
\end{aligned}
$$
Since the projection matrix $\bH$ projects vectors onto $\mathcalV$, we have $\bH\bv = \bv_p$ and $\bH\bw = \bw_p$. Therefore,
$$
\begin{aligned}
(\bH\bv)^\top \bw&=\bv_p^\top \bw=\bv_p^\top (\bw_p + \bw_n)  &\qquad  \bv^\top (\bH\bw) &=\bv^\top \bw_p=(\bv_p+\bv_n)^\top \bw_p \\
&=\bv_p^\top \bw_p, &\qquad &=\bv_p^\top \bw_p, \\
\end{aligned}
$$
where the last equalities follow from the fact that $\bv_p$ is perpendicular to $\bw_n$, and $\bv_n$ is perpendicular to $\bw_p$. Thus, we have 
$$
\begin{aligned}
(\bH\bv)^\top \bw = \bv^\top (\bH\bw) 
\quad\implies\quad
 \bv^\top \bH^\top \bw = \bv^\top \bH\bw
 \quad\implies\quad
 \bH^\top = \bH.
\end{aligned}
$$

Conversely, if a projection matrix $\bH$ (not necessarily an orthogonal projection) is symmetric, then any vector $\bv$ can be decomposed into $\bv = \bH\bv + (\bI-\bH)\bv$. If we can prove $\bH\bv$ is perpendicular to $(\bI-\bH)\bv$, then we complete the proof. To see this, we have 
$$
\begin{aligned}
(\bH\bv)^\top (\bI-\bH)\bv &= \bv^\top \bH^\top (\bI-\bH)\bv
=\bv^\top  (\bH^\top-\bH^\top\bH)\bv \\
&=\bv^\top  (\bH-\bH\bH)\bv
=\bv^\top  (\bH-\bH)\bv = \bzero,
\end{aligned}
$$
which completes the proof.
\end{proof}

We claimed that the orthogonal projection has a minimum length, i.e., the distance between $\bb$ and $\bp$ is minimized. We rigorously prove this property in the following proposition.
\begin{proposition}[Minimum Distance in Orthogonal Projection]\label{proposition:minimal-distance-orthogonal-projection}
Let $\mathcalV$ be a subspace of $\real^m$ and $\bH$ be an orthogonal projection onto $\mathcalV$. Then, for any vector $\bb\in\real^n$, we have 
$$
\normtwo{\bb - \bH\bb}^2 \leq \normtwo{\bb - \bv}^2, \qquad \forall \bv \in \mathcalV.
$$
\end{proposition}
\begin{proof}[of Proposition~\ref{proposition:minimal-distance-orthogonal-projection}]
Let $\bH = \bQ\bLambda\bQ^\top \in \real^{m\times m}$ be the spectral decomposition of the orthogonal projection $\bH$, where $\bQ=[\bq_1, \bq_2, \ldots, \bq_m]$ is the column partition of $\bQ$, and $\bLambda=\diag(\lambda_1, \lambda_2, \ldots, \lambda_m)$. Let $\dim(\mathcalV) = r$. 
From Lemma~\ref{lemma:eigenvalues-of-projection}, the only possible eigenvalues of an orthogonal projection matrix are 1 and 0. Without loss of generality, let $\lambda_1=\lambda_2=\ldots=\lambda_r=1$ and $\lambda_{r+1}=\lambda_{r+2}=\ldots=\lambda_m=0$. Then, it follows that
\begin{itemize}
\item $\{\bq_1, \bq_2, \ldots, \bq_m\}$ is an orthonormal basis of $\real^m$;

\item $\{\bq_1, \bq_2, \ldots, \bq_r\}$ is an orthonormal basis of $\mathcalV$. So for any vector $\bv\in \mathcalV$, we have $\bv^\top\bq_i=0$ for $i\in \{r+1, r+2, \ldots, m\}$. 
\end{itemize}
Then we have,
$$
\begin{aligned}
\normtwo{\bb - \bH\bb}^2 &\stackrel{*}{=} \normtwo{\bQ^\top\bb - \bQ^\top\bH\bb}^2
= \sum_{i=1}^{m}(\bb^\top \bq_i - (\bH\bb)^\top \bq_i)^2 
\stackrel{\dag}{=}\sum_{i=1}^{m}(\bb^\top \bq_i - \bb^\top\bH \bq_i)^2\\
&\stackrel{+}{=}\sum_{i=1}^{m}(\bb^\top \bq_i - \lambda_i\bb^\top \bq_i)^2 
\stackrel{\star}{=}0+\sum_{i=r+1}^{m}(\bb^\top \bq_i)^2  
\leq \sum_{i=1}^{r}(\bb^\top\bq_i - \bv^\top\bq_i)^2+ \sum_{i=r+1}^{m}(\bb^\top \bq_i)^2 \\
&\stackrel{\ddag}{=}\normtwo{\bQ^\top \bb - \bQ^\top \bv}^2 
=\normtwo{\bb - \bv}^2,
\end{aligned}
$$
where equality ($*$) follows from the invariance under orthogonal transformations, equality ($\dag$) follows from the symmetry of $\bH$, equality (+) follows from $\bH\bQ =  \bQ\bLambda$, equality ($\star$) follows from the eigenvalues of $\bH$ being 1 or 0, and equality ($\ddag$) follows from $\bv^\top\bq_i=0$ for $i>r$.
This completes the proof.
\end{proof}

\begin{proposition}[Angle between the Original and Projected Vectors]\label{proposition:angle-orthogonal-projection}
Let $\bH$ be an orthogonal projection onto $\mathcalV$. Then,
\begin{enumerate}
\item $\bb^\top (\bH\bb) \geq 0$, i.e., angle between $\bb$ and $\bH\bb$ is smaller than $90^\circ$;
\item $\normtwo{\bH\bb}^2 \leq \normtwo{\bb}^2$.
\end{enumerate}
\end{proposition}
\begin{proof}[of Proposition~\ref{proposition:angle-orthogonal-projection}]
By the definition of the orthogonal projection, we have $\bb^\top (\bH\bb) $$=\bb^\top \bH (\bH\bb)$$=\bb^\top \bH^\top (\bH\bb) =\normtwo{\bH\bb}^2 \geq 0$.
And we can decompose $\bb$ as 
$$
\begin{aligned}
\normtwo{\bb}^2 = \normtwo{(\bI-\bH+\bH)\bb}^2 
&= \normtwo{(\bI-\bH)\bb}^2 + \normtwo{\bH\bb}^2 + 2 \bb^\top (\bI-\bH)^\top\bH\bb\\
&= \normtwo{(\bI-\bH)\bb}^2 + \normtwo{\bH\bb}^2 \geq \normtwo{\bH\bb}^2.
\end{aligned}
$$
This completes the proof.
\end{proof}

In conclusion, regarding the properties of projection and orthogonal projection matrices, we define a projection matrix to be idempotent, and it is symmetric when restricted to the orthogonal projection. We have proven that the distance between the original vector and the projected vector is minimal in the context of orthogonal projections.
In the following proposition, we construct an orthogonal projection matrix from a set of vectors.

\begin{proposition}[Projection Matrix from a Set of Vectors]\label{proposition:projection-from-matrix}
Let $\ba_1, \ba_2, \ldots, \ba_n \in\real^m$ be linearly independent and span the subspace $\mathcalV$: $\cspace([\ba_1, \ba_2, \ldots, \ba_n]) = \mathcalV$, where $m\geq n$. Then, the orthogonal projection onto $\mathcalV$ can be represented as 
$$
\bH = \bA(\bA^\top\bA)^{-1}\bA^\top,
$$
where $\bA=[\ba_1, \ba_2, \ldots, \ba_n]\in \real^{m\times n}$.
\end{proposition}
\begin{proof}[of Proposition~\ref{proposition:projection-from-matrix}]
It can be easily verified $\bH$ is symmetric and idempotent. By the SVD of $\bA=\bU\bSigma\bV^\top$, we have $\bH = \bA(\bA^\top\bA)^{-1}\bA^\top = \bU\bSigma(\bSigma^\top\bSigma)^{-1} \bSigma^\top \bU^\top$. Let $\bU=[\bu_1, \bu_2, \ldots, \bu_m]$ be the column partition of $\bU$. From Proposition~\ref{proposition:svd-four-orthonormal-Basis}, $\{\bu_1, \bu_2, \ldots, \bu_n\}$ is an orthonormal basis of $\cspace(\bA)$.
And the matrix $\bSigma(\bSigma^\top\bSigma)^{-1} \bSigma^\top$ in $\bH$ is an $m\times m$ matrix with the upper-left part being an $n\times n$ identity matrix and the other parts being zero. 
Applying this observation of $\bH$ to the spectral theorem, $\{\bu_1, \bu_2, \ldots, \bu_n\}$ is also an orthonormal basis of $\cspace(\bH)$. 
Thus, we have $\cspace(\bH)=\cspace(\bA)$, and the orthogonal projection $\bH$ projects onto $\cspace(\bA)$, which completes the proof.
\end{proof}


\paragraph{Geometric interpretation of LS.}
The proposition above brings us back to the result we demonstrated at the beginning of this section. 
For the LS estimate to minimize $\normtwo{\bb-\bA\bx}^2$, we seek an estimate $\bx_{LS}$ such that $\bp = \bA\bx_{LS}$ lies in $\cspace(\bA)$ and  minimizes the distance between $\bA\bx_{LS}$ and $\bb$. An orthogonal projection matrix $\bH = \bA(\bA^\top\bA)^{-1}\bA^\top$ can project $\bb$ onto the column space of $\bA$, resulting in the  projected vector  $\bp = \bH\bb$. By Proposition~\ref{proposition:minimal-distance-orthogonal-projection},  the squared distance between $\bp$ and $\bb$ is minimal.

To reiterate, the hat matrix $\bH$ has a geometric interpretation. $\bH$ drops a perpendicular from $\bb$ to the hyperplane defined by the column space of $\bA$. 
Specifically, $\bH = \bA(\bA^\top\bA)^{-1}\bA^\top$ projects $\bb$ onto the column space of $\bA$: $\bp = \bH\bb$.
The idempotency of $\bH$ also has a geometric interpretation: applying  $\bH$ multiple times has no additional effect because the vector is already on the hyperplane. 
Therefore, $\bH^2 \bb = \bH\bb$. This situation is illustrated in Figure~\ref{fig:ls-geometric1}. The sum of squared error is then equal to the squared Euclidean distance between $\bb$ and $\bp$. 
Consequently, the least squares solution for $\bx$ corresponds to the orthogonal projection of $\bb$ onto the column space of $\bA$.

\begin{proposition}[Column Space of Projection Matrices]\label{proposition:column-space-of-projection}
We notice that the hat matrix $\bH = \bA(\bA^\top\bA)^{-1}\bA^\top$ is employed to project any vector in $\real^m$ onto the column space of $\bA \in \real^{m\times n}$ (with full column rank). That is, $\bH\bb \in \cspace(\bA)$. Notice again that  $\bH\bb$ is nothing but a linear combination of the columns of $\bH$, hence $\cspace(\bH) = \cspace(\bA)$. 

In general, for any projection matrix $\bH$ that  projects a vector onto a subspace $\mathcalV$, it follows that $\cspace(\bH) = \mathcalV$.
\end{proposition}
While it is straightforward to prove this using the proof of Proposition~\ref{proposition:projection-from-matrix}, we will demonstrate it using a different approach.
\begin{proof}[of Proposition~\ref{proposition:column-space-of-projection}]
Since $\bH = \bA(\bA^\top\bA)^{-1}\bA^\top \triangleq \bA\bC$ (let $\bC\triangleq(\bA^\top\bA)^{-1}\bA^\top$), the columns of $\bH$ are linear combinations of the columns of $\bA$. Thus, we have $\cspace(\bH) \subseteq \cspace(\bA)$.	
By Lemma~\ref{lemma:rank-of-symmetric-idempotent2}, we have 
$$
\rank(\bH) = \trace(\bH) =\trace(\bA(\bA^\top\bA)^{-1}\bA^\top) = \trace((\bA^\top\bA)^{-1}\bA^\top\bA) =\trace(\bI_n)=n,
$$
where the third equality follows from the cyclic property of the trace: $\trace(\bX\bY\bZ) = \trace(\bY\bZ\bX) = \trace(\bZ\bX\bY)$.
Thus, the rank of $\bH$ equals the rank of $\bA$, implying that  $\cspace(\bH) = \cspace(\bA)$.
\end{proof}

\subsection{Properties of Orthogonal Projection Matrices}
In fact, $(\bI-\bH)$ is also a symmetric idempotent matrix if $\bH$ is  symmetric and idempotent. 
Specifically, when $\bH$ projects onto a subspace $\mathcalV$, $\bI-\bH$ projects onto the perpendicular subspace $\mathcalV^\perp$.

\begin{proposition}[Project onto $\mathcalV^\perp$]\label{proposition:orthogonal-projection}
Let $\mathcalV$ be a subspace and $\bH$ be an orthogonal projection onto $\mathcalV$. Then, $\bI-\bH$ is the orthogonal projection matrix onto $\mathcalV^\perp$.

The claim can be extended further that suppose $\mathcalV_1 \subseteq \mathcalV$ and $\mathcalV_2 \subseteq \mathcalV^\perp$. Then, $\bH_1$ is the orthogonal projection that projects onto $\mathcalV_1$ and $\bH_2$ is the orthogonal projection that projects onto $\mathcalV_2$ if and only if $\bH_1\bH_2 = \bzero$.
\end{proposition}

\begin{proof}[of Proposition~\ref{proposition:orthogonal-projection}]
First, $(\bI-\bH)$ is symmetric, $(\bI-\bH)^\top = \bI - \bH^\top = \bI-\bH$ since $\bH$ is symmatric. And 
$
(\bI-\bH)^2 = \bI^2 -\bI\bH -\bH\bI +\bH^2 = \bI-\bH.
$
Thus, $\bI-\bH$ is an orthogonal projection matrix. By the spectral theorem again, let $\bH =\bQ\bLambda\bQ^\top$. Then, $\bI-\bH = \bQ\bQ^\top - \bQ\bLambda\bQ^\top = \bQ(\bI-\bLambda)\bQ^\top$. Hence the column space of $\bI-\bH$ is spanned by the eigenvectors of $\bH$ corresponding to the zero eigenvalues of $\bH$ (by Lemma~\ref{lemma:eigenvalues-of-projection}), which coincides with $\mathcalV^\perp$.

For the second part, since $\cspace(\bH_1) = \mathcalV_1$ and $\cspace(\bH_2) = \mathcalV_2$, every column of $\bH_1$ is perpendicular to columns of $\bH_2$. Therefore, it follows that $\bH_1\bH_2=\bzero$. Conversely, if $\bH_1\bH_2=\bzero$, then $\bH_1(\bH_2\bb) = 0$ for all $\bb$. Thus, $\mathcalV_1\perp \mathcalV_2$.
Moreover, it can be easily verified when $\mathcalV_1 = \mathcalV$ and $\mathcalV_2 = \mathcalV^\perp$, then $\bH(\bI-\bH) = \bzero$. 
\end{proof}

A projection matrix that can project any vector onto a subspace is not unique. However, when restricted to the orthogonal projection,  the projection is unique.

\begin{proposition}[Uniqueness of Orthogonal Projection]\label{proposition:unique-projection-orthogonal}
If $\bP$ and $\bH$ are orthogonal projection matrices onto the same subspace $\mathcalV$, then $\bP = \bH$.
\end{proposition}
\begin{proof}[of Proposition~\ref{proposition:unique-projection-orthogonal}]
For any vector $\bv$ in $\real^m$, it can be decomposed into a vector $\bv_p$ in $\mathcalV$ and a vector $\bv_n$ in $\mathcalV^\perp$ such that $\bv = \bv_p+\bv_n$ and $\bv_p^\top \bv_n=0$. Then,
$
\bP \bv = \bv_p = \bH\bv,
$
which implies  $(\bP-\bH)\bv = \bzero$. Since any vector $\bv\in \real^m$ is in the null space of $\bP-\bH$, then $\bP-\bH$ must be of rank 0. 
Therefore, $\bP=\bH$.
\end{proof}

\begin{proposition}[Nested Projection\index{Nested projection}]\label{proposition:nested-projection}
Let $\mathcalV_1 \subseteq \mathcalV_2 \subseteq \real^m$ be two nested linear subspaces. If $\bH_1$ is the orthogonal projection onto $\mathcalV_1$, and $\bH_2$ is the orthogonal projection onto $\mathcalV_2$, then the following hold: 
\begin{enumerate}
\item  $\bH_2\bH_1 = \bH_1 = \bH_1\bH_2$;
\item  $\bH_2-\bH_1$ is also an orthogonal projection.
\end{enumerate}
\end{proposition}
\begin{proof}[of Proposition~\ref{proposition:nested-projection}]
For all $\bb \in \real^m$, we have $\bH_1\bb \in \mathcalV_1$. This implies $\bH_1\bb\in \mathcalV_1\subseteq\mathcalV_2$. Thus,
$$
\bH_2(\bH_1\bb) = \bH_1\bb. \qquad (\text{from Definition~\ref{definition:projection-matrix}})
$$
Therefore, $(\bH_2\bH_1 -\bH_1)\bb = \bzero$ for all $\bb\in \real^m$. 
This means the dimension of the null space $\nspace(\bH_2\bH_1 -\bH_1)=m$ and the rank of $\bH_2\bH_1 -\bH_1$ is 0, which results in $\bH_2\bH_1 =\bH_1$.

For $\bH_1\bH_2$, since both $\bH_1$ and $\bH_2$ are symmetric, we have $\bH_1\bH_2 = \bH_1^\top\bH_2^\top = (\bH_2\bH_1)^\top = \bH_1^\top = \bH_1$ which completes the proof of part 1.

To prove the second part, we notice that $(\bH_2-\bH_1)^\top = \bH_2-\bH_1$ and 
$$
\begin{aligned}
(\bH_2-\bH_1)^2 &= \bH_2^2-\bH_2\bH_1 - \bH_1\bH_2 + \bH_1^2 \\
&=\bH_2-\bH_1 - \bH_1 + \bH_1 
=\bH_2 - \bH_1,
\end{aligned}
$$
which states that $\bH_2-\bH_1$ is both symmetric and idempotent. This completes the proof.
\end{proof}

From the proposition above, we can conclude that orthogonal projection matrices are positive semidefinite (PSD).
\begin{proposition}[PSD in Projection]\label{proposition:symmetric-projection-psd}
Any orthogonal projection matrix $\bH$ is positive semidefinite (PSD).
\end{proposition}
\begin{proof}[of Proposition~\ref{proposition:symmetric-projection-psd}]
Since $\bH$ is symmetric and idempotent. For any vector $\bx$, we have 
$
\bx^\top \bH \bx = \bx^\top \bH\bH \bx = \bx^\top \bH^\top\bH \bx=\normtwo{\bH \bx}^2 \geq 0.
$
Thus, $\bH$ is PSD.
\end{proof}

To summarize, we list some key properties of the orthogonal projection matrix in the following remark.
\begin{remark}[Important Facts About Hat Matrix (Part 1)]
\item 1. As we assume $\bA\in \real^{m\times n}$ is fixed and has full rank with $m\geq n$. It is known that the rank of $\bA$ is equal to the rank of its \textit{Gram matrix}, defined as $\bA^\top\bA$, such that 
$
\rank(\bA) = \rank(\bA^\top\bA).
$

\item 2. The rank of an orthogonal projection matrix is the dimension of the subspace onto which it projects. Hence, the rank of $\bH$ is $n$ when $\bA$ has full rank and $m\geq n$:
\begin{equation}
\rank(\bH) = \rank\left(\bA(\bA^\top\bA)^{-1}\bA^\top\right) = \rank\left((\bA^\top\bA)^{-1}\bA^\top\bA\right) = \rank(\bI_n) = n;\nonumber
\end{equation}

\item 3. The column space of $\bH$ is identical to the column space of $\bA$.
\end{remark}

\subsection{Distance Between Subspaces}
Suppose $\mathcalV_1$ and $\mathcalV_2$ are subspaces of $\real^n$ with $\dim(\mathcalV_1)=\dim(\mathcalV_2)$. We define the \textit{distance} between the two spaces by 
$$
\dist(\mathcalV_1,\mathcalV_2) = \normtwo{\bH_1 - \bH_2},~\footnote{Alternative definitions of the distance between two spaces are discussed in \citet{chen2021spectral}.}
$$
the spectral norm between the matrices,
where $\bH_1$ is the orthogonal projection onto $\mathcalV_1$, and $\bH_2$ is the orthogonal projection onto $\mathcalV_2$.
\begin{theorem}[Subspace Distance]\label{theorem:subspace-distance}
Let $\bX, \bY\in \real^{n\times n}$ be orthogonal matrices. Suppose their column partitions are
\[
\begin{blockarray}{ccc}
\begin{block}{c[cc]}
\bX=&	\bX_1 & \bX_2  \\
\end{block}
&n\times r & n\times (n-r)   \\
\end{blockarray}
,\qquad
\begin{blockarray}{ccc}
\begin{block}{c[cc]}
\bY=	&	\bY_1 & \bY_2   \\
\end{block}
&n\times r & n\times (n-r)   \\
\end{blockarray},
\]
where $\bX_r$ and $\bY_r$ are the first $r$ columns of $\bX$ and $\bY$, respective. Suppose further that $\mathcalV_1 = \cspace(\bX_1)$ and $\mathcalV_2 = \cspace(\bY_1)$. Then,
$$
\dist(\mathcalV_1, \mathcalV_2) = \normtwo{\bX_1^\top \bY_2} = \normtwo{\bY_1^\top\bX_2}.
$$
\end{theorem}

\begin{proof}[of Theorem~\ref{theorem:subspace-distance}]
It can be easily verified $\bX_1\bX_1^\top$ and $\bY_1\bY_1^\top$ are the orthogonal projections onto $\mathcalV_1$ and $\mathcalV_2$, respectively.
Write out the equation
$$
\begin{aligned}
\dist(\mathcalV_1, \mathcalV_2) &= \normtwo{\bX_1\bX_1^\top - \bY_1\bY_1^\top}
\stackrel{\dag}{=} \normtwo{\bX^\top (\bX_1\bX_1^\top - \bY_1\bY_1^\top)\bY}   \\
&= 
\normtwo{\begin{bmatrix}
		\bX_1^\top \\
		\bX_2^\top
	\end{bmatrix} (\bX_1\bX_1^\top - \bY_1\bY_1^\top)[\bY_1, \bY_2]}
\stackrel{+}{=} 
\normtwo{\begin{bmatrix}
		\bI_r \\
		\bzero 
	\end{bmatrix} \bX_1^\top[\bY_1, \bY_2] -\begin{bmatrix}
		\bX_1^\top \\
		\bX_2^\top
	\end{bmatrix}  \bY_1[\bI_r, \bzero]}\\
&= \left\|\begin{bmatrix}
\bzero & \bX_1^\top \bY_2 \\
-\bX_2^\top\bY_1 & \bzero
\end{bmatrix}\right\|_2,
\end{aligned}
$$
where equality $(\dag)$ follows from the invariance under orthogonal transformations, and equality (+) follows from the fact that $\bX_2^\top\bX_1 =\bzero$ and $ \bY_1^\top\bY_2=\bzero$.
We note that $\bX_1^\top \bY_2$ and $\bX_2^\top\bY_1$ are submatrices of the following orthogonal matrix
$$
\bQ = \bX^\top\bY = \begin{bmatrix}
\bQ_{11} & \bQ_{12} \\
\bQ_{21} & \bQ_{22}
\end{bmatrix}
=\begin{bmatrix}
\bX_1^\top \bY_1 & \bX_1^\top\bY_2 \\
\bX_2^\top \bY_1 & \bX_2^\top\bY_2
\end{bmatrix}.
$$
For any unit vector $\bx\in \real^{r}$ ($\normtwo{\bx}=1$), it follows that 
$$
\bQ\begin{bmatrix}
\bx \\
\bzero 
\end{bmatrix}
=
\begin{bmatrix}
\bQ_{11} \bx \\
\bQ_{21}\bx
\end{bmatrix}.
$$
By invariant under orthogonal matrices \footnote{$\normtwo{\bx} = \normtwo{\bQ\bx}$ given the orthogonal matrix $\bQ$.}, we have $\normtwo{\bx}^2 = \normtwo{\bQ_{11} \bx}^2 + \normtwo{\bQ_{21}\bx}^2=1$. Thus, by the Definition~\ref{definition:spectral_norm_app} of the spectral norm, we have 
$$
\normtwo{\bQ_{21}}^2  =\mathop{\max}_{\bx\in \real^r: \normtwo{\bx}=1}  \normtwo{\bQ_{21}\bx}^2 = 
1-\mathop{\min}_{\bx\in \real^r: \normtwo{\bx}=1}  \normtwo{\bQ_{11}\bx}^2
= 1 - \sigma_{min}(\bQ_{11})^2,
$$
where $\sigma_{min}(\bQ_{11})$ is the smallest signular value of $\bQ_{11}$.

Analogously, by applying the process above to $\bQ^\top$, we will obtain 
$$
\normtwo{\bQ^\top_{12}}^2 = 1- \sigma_{min}(\bQ_{11}^\top)^2
\quad\implies \quad 
\normtwo{\bQ_{12}}^2 = 1- \sigma_{min}(\bQ_{11})^2
\quad\implies \quad 
\normtwo{\bQ_{21}} = \normtwo{\bQ_{12}}.
$$
Write out the distance again by 
$$
\begin{aligned}
\dist(\mathcalV_1, \mathcalV_2) 
&= \left\|\begin{bmatrix}
\bzero & \bQ_{12} \\
-\bQ_{21} & \bzero
\end{bmatrix}\right\|_2 \triangleq\left\|\bZ\right\|_2 .
\end{aligned}
$$
Let $\bw \in \real^{n}$ with $\normtwo{\bw}\triangleq\normtwo{\footnotesize\begin{bmatrix}
\bu\\
\bv
\end{bmatrix}}=1$, we have 
$$
\begin{aligned}
\normtwo{\bZ\bw}
&= \left\|\begin{bmatrix}
\bzero & \bQ_{12} \\
-\bQ_{21} & \bzero
\end{bmatrix} \begin{bmatrix}
\bu \\
\bv
\end{bmatrix}\right\| 
= \left\|\begin{bmatrix}
\bQ_{12}\bv \\
-\bQ_{21}\bu
\end{bmatrix}\right\|
=\sqrt{\normtwo{\bQ_{12}\bv}^2 + \normtwo{\bQ_{21}\bu}^2} \\
&\stackrel{*}{\leq} \normtwo{\bQ_{12}} \sqrt{\normtwo{\bv}^2 +\normtwo{\bu}^2} =\normtwo{\bQ_{21}} \sqrt{\normtwo{\bv}^2 +\normtwo{\bu}^2}  
= \normtwo{\bQ_{12}} \normtwo{\bw} =\normtwo{\bQ_{21}} \normtwo{\bw} \\
&= \normtwo{\bQ_{12}}  =\normtwo{\bQ_{21}},
\end{aligned}
$$
where the inequality ($*$) follows from the matrix-vector product that $\normtwo{\bA\bx} \leq \normtwo{\bA} \normtwo{\bx}$ for any matrix $\bA$ and vector $\bx$.
The upper bound of $\normtwo{\bZ\bw} $ is the spectral norm $\dist(\mathcalV_1, \mathcalV_2)=\normtwo{\bZ} $, and this completes the proof.
\end{proof}
From the theorem above, we realize that 
$
0 \leq \dist(\mathcalV_1, \mathcalV_2) \leq 1.
$
And it is easy to show that 
$$
\left\{
\begin{aligned}
\dist(\mathcalV_1, \mathcalV_2) &= 0 \qquad \rightarrow \qquad \mathcalV_1=\mathcalV_2;\\
\dist(\mathcalV_1, \mathcalV_2) &= 1 \qquad\rightarrow \qquad \mathcalV_1 \cap \mathcalV_2^\perp \neq \{0\}.
\end{aligned}
\right.
$$


\subsection{Projection for LS with Noise Disturbance\index{Noise disturbance}}\label{sec:geometry-noise-disturbance}
\begin{figure}[h!]
\centering
\includegraphics[width=0.5\textwidth]{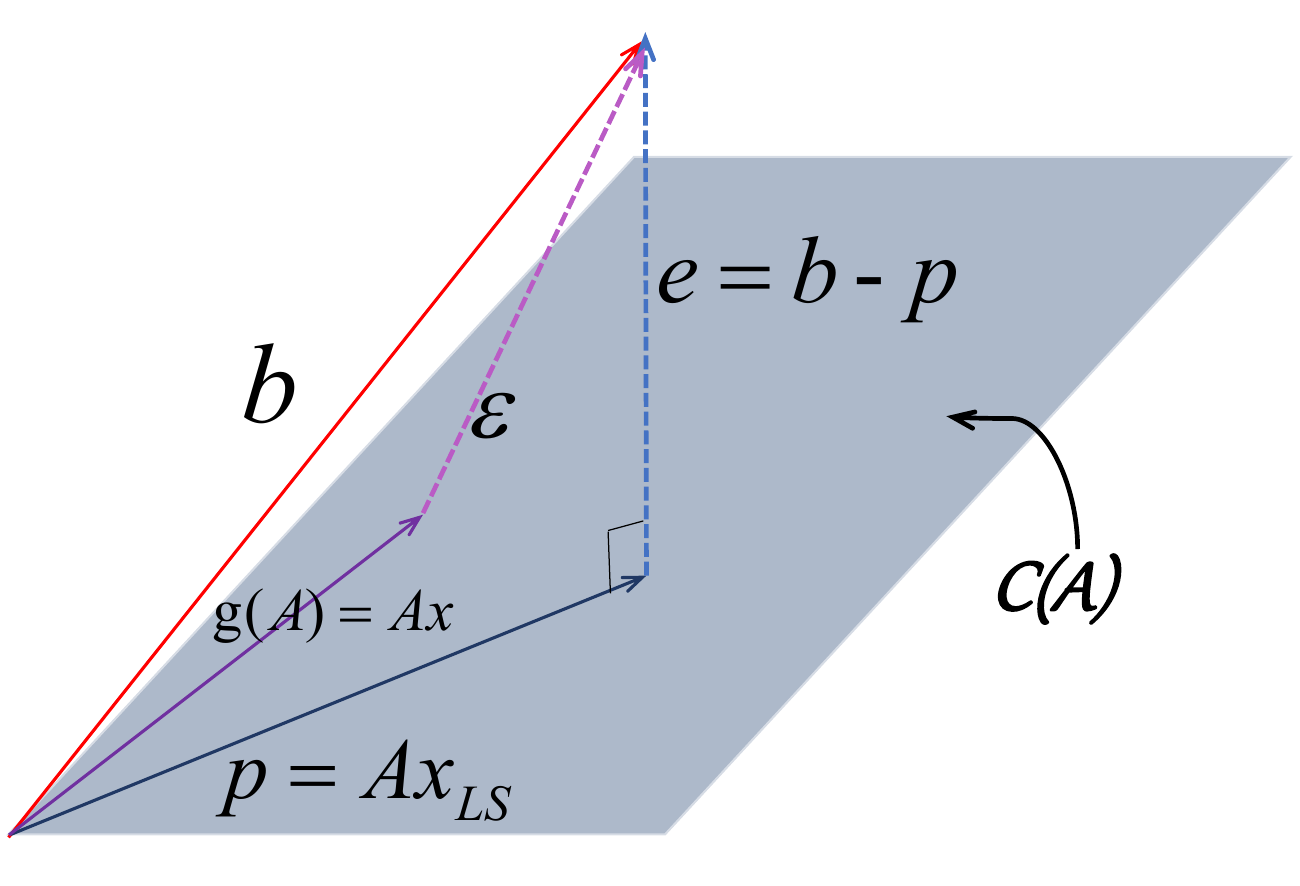}
\caption{Projection onto the hyperplane of $\cspace(\bA)$ and disturbed by noise.}
\label{fig:ls-geometric2}
\end{figure}

We   further consider that $\bb$ comes from some ideal function $g(\bA) \in \cspace(\bA) $ and $\bb =g(\bA)+  \bepsilon$, where $\bepsilon$ is the noise. 
That is, the actual observation $\bb$ is perturbed by some noise.
In this case, we assume that the observed value  $\bb$ differs from the true function $g(\bA)=\bA\bx$ by an additive noise. 
This scenario is illustrated in Figure~\ref{fig:ls-geometric2}, providing a comprehensive view of the linear problem.

\begin{mdframed}[hidealllines=\mdframehidelineNote,backgroundcolor=\mdframecolor]
\item 1. Vector of outputs (responses): $\bb  \in \mathbb{R}^m$ is an $m\times 1$ vector of observations of the output variable, where $m$ is the sample size.

\item 2. Design matrix (data matrix): $\bA$ is an $m\times n$ matrix of inputs, and $n$ is the dimension of the inputs for each observation.~\footnote{The term of \textit{design matrix} is usually used to mean that the matrix possesses  some special properties, such as the \textit{restricted isometry property (RIP)}. This property signifies that the matrix maintains the geometric structure of vectors subsequent to linear transformation, thereby ensuring the convergence of diverse algorithms \citep{jain2017non}.}

\item 3. Vector of parameters: $\bx \in \mathbb{R}^n$ is an $n\times 1$ vector of regression coefficients.

\item 4. Vector of noises: $\bepsilon \in \mathbb{R}^m$ is an $m\times 1$ vector of noises.
\item 5. Vector of errors: $\be \in \mathbb{R}^m$ is an $m\times 1$ vector of errors. For predicted outputs $\bp$, $\be = \bb - \bp$. Thus, $\be$ is different from $\bepsilon$. The latter is unobservable, while the former is a byproduct of our linear model. In some texts, $\be$ is denoted by $\hat{\bepsilon}$ to establish a connection to $\bepsilon$.
\end{mdframed}

By introducing the noise vector, we can provide additional insights about the hat matrix:
\begin{remark}[Important Facts About Hat Matrix (Part 2)]
\item 4. \textbf{Error vector.} $\be = \bb-\bp = \underline{(\bI-\bH)\bb} = (\bI-\bH)(\bA\bx+\bepsilon) = (\bI-\bH)\bA\bx + (\bI-\bH)\bepsilon=\underline{(\bI-\bH)\bepsilon}$: projecting $\bb$ onto the perpendicular space is equivalent to projecting $\bepsilon$ onto the perpendicular space. This property can be easily verified from the geometric meaning of $\bb$ and $\bepsilon$, as shown in Figure~\ref{fig:ls-geometric2};

\item 5. \textbf{Orthogonality.} $\bp$ and $\be$ are orthogonal, $\bH \bepsilon$ and $\be$ are orthogonal;

\item 6. \textbf{Pythagoras.} $\normtwo{\bb}^2 = \normtwo{\bp}^2 + \normtwo{\be}^2$ and  $\normtwo{\bepsilon}^2 = \normtwo{\bH\bepsilon}^2 + \normtwo{\be}^2$;

\item 7. \textbf{Pythagoras in general.} for any orthogonal projection matrix $\bP$, we have $\normtwo{\bx}^2 = \normtwo{\bP\bx}^2 + \normtwo{(\bI-\bP)\bx}^2$.
\end{remark}
The general Pythagorean theorem can be verified as follows:
$$
\begin{aligned}
\normtwo{\bP\bx}^2 &+ \normtwo{(\bI-\bP)\bx}^2 = \bx^\top\bP^\top\bP\bx + \bx^\top (\bI-\bP)^\top(\bI-\bP)\bx\\
&=\bx^\top\bP\bx + \bx^\top (\bI-\bP)\bx 
=\bx^\top[\bP\bx+(\bI-\bP)\bx]
=\normtwo{\bx}^2.
\end{aligned}
$$
A more detailed analysis of this noise-perturbed  linear model can be found in \citet{lu2021rigorous}. This additive noise framework also serves as the foundation for linear models using Bayesian approaches.

%% file: chapter-app_pseudo.tex
\newpage
\section{Pseudo-Inverse}\label{appendix:pseudo-inverse}
If the matrix $\bA$ is nonsingular, then the linear system $\bb = \bA\bx$ can be easily solved using the inverse of $\bA$ such that $\widehat{\bx} = \bA^{-1}\bb$. However, the inverse of an $m\times n$ matrix $\bA$ does not exist if $\bA$ is not square or if $\bA$ is singular. 
In such cases, we can still find its \textit{pseudo-inverse}, an $n\times m$ matrix denoted by $\bA^+$ (if $\bA\in\real^{m\times n}$).
The pseudo-inverse plays a crucial role in solving linear systems by the method of least squares.
Before the discussion of pseudo-inverse, we will introduce the concepts of \textit{one-sided inverse}, \textit{generalized inverse}, and \textit{reflexive generalized inverse}, which are essential prerequisites. However, readers can skip the three sections and still understand the main ideas about the pseudo-inverse.

\index{Pseudo-inverse}

\subsection{One-Sided Inverse}
The one-sided inverse is defined as follows:
\begin{definition}[One-Sided Inverse\index{One-sided inverse}]\label{definition:one_side_inverse}
Let   $\bA\in \real^{m\times n}$ be any matrix. If there is a matrix $\bA_L^{-1} \in \real^{n\times m}$ such that
$$
\bA_L^{-1} \bA = \bI_n,
$$
then $\bA_L^{-1}$ is a \textbf{left inverse} of $\bA$, and $\bA$ is said to be  \textbf{left-invertible}.
Similarly, if there is a matrix $\bA_R^{-1} \in \real^{n\times m}$ such that 
$$
\bA \bA_R^{-1}= \bI_m,~\footnote{Note here the superscript $-1$ in $\bA_L^{-1}$ and $\bA_R^{-1}$ does not denote the inverse of $\bA_L$ or $\bA_R$ but rather the one-sided inverse of $\bA$.}
$$
then $\bA_R^{-1}$ is a \textbf{right inverse} of $\bA$, and $\bA$ is said to be \textbf{right-invertible}.
\end{definition}

\begin{proposition}[One-Sided Invertible]\label{proposition:one-sided-invertible}
Let $\bA\in \real^{m\times n}$ be any matrix. Then,
\begin{itemize}
\item   $\bA$ is left-invertible if and only if $\bA$ has full column rank (which implies $m\geq n$);
\item   $\bA$ is right-invertible if and only if $\bA$ has full row rank (which implies $m\leq n$).
\end{itemize}
\end{proposition}
\begin{proof}[of Proposition~\ref{proposition:one-sided-invertible}]
Suppose $\bA$ has full column rank, then $\bA^\top \bA \in \real^{n\times n}$ has full rank (Lemma~\ref{lemma:rank-of-ata}). Therefore, $(\bA^\top \bA)^{-1}(\bA^\top \bA) = \bI_n$. That implies \textcolor{mylightbluetext}{$(\bA^\top \bA)^{-1}\bA^\top$ is a left inverse of $\bA$}. 

Conversely,  suppose $\bA$ is left-invertible and there exists an $\bA_L^{-1}$ such that $\bA_L^{-1} \bA = \bI_n$. Since each  row of $\bA_L^{-1} \bA$ is a linear  combination of the rows of $\bA$,  the row space of $\bA_L^{-1} \bA$ is a subset of the row space of $\bA$. Thus, $\rank(\bA) \geq \rank(\bA_L^{-1} \bA) =\rank(\bI_n) = n$, which implies $\rank(\bA)=n$, and $\bA$ has full column rank. 

Similarly, we can show $\bA$ is right-invertible if and only if $\bA$ has full row rank, and in this case, \textcolor{mylightbluetext}{$\bA^\top (\bA\bA^\top)^{-1}$ is a right inverse of $\bA$}.
\end{proof}

We have shown in the above proof that $(\bA^\top \bA)^{-1}\bA^\top$ is a specific left inverse of $\bA$ if $\bA$ has full column rank, and $\bA^\top (\bA\bA^\top)^{-1}$ is a specific right inverse of $\bA$ if $\bA$ has full row rank.
However, the inverse of a $k\times k$ nonsingular matrix requires $2k^3$ floating-points operations (flops), which is a computationally intensive process (Section~\ref{section:inverse-by-lu}).
In our case, the inverse of $\bA^\top\bA$ requires $2n^3$ flops, and the inverse of $\bA\bA^\top$ requires $2m^3$ flops. A simpler approach to finding a one-sided inverse is through elementary operations.

Suppose $\bA\in \real^{m\times n}$ has full column rank. By applying \textit{row elementary operations} represented by a matrix $\bE\in \real^{m\times m}$ to $[\bA, \bI_m]$, we obtain:
\begin{equation}\label{equation:onesided-1}
\bE \begin{bmatrix}
\bA & \bI_m
\end{bmatrix}=
\begin{bmatrix}
\bI_n & \bG \\
\bzero & \bZ
\end{bmatrix},
\end{equation}
where $\bG \in \real^{n\times m}$, $\bI_m$ is the $m\times m$ identity matrix, $\bI_n$ is the $n\times n$ identity matrix, and $\bZ$ is an $(m-n)\times m$ matrix.
Subsequently, it can be easily verified that $\bG\bA = \bI_n$, establishing  $\bG$ as a left inverse of $\bA$.

Similarly, suppose $\bA\in \real^{m\times n}$ has full row rank. By applying \textit{column elementary operations} represented by a matrix  $\bE \in \real^{n\times n}$ to $[\bA^\top, \bI_n]^\top$, we obtain:
\begin{equation}\label{equation:onesided-2}
\begin{bmatrix}
\bA \\
\bI_n
\end{bmatrix} \bE=
\begin{bmatrix}
\bI_m & \bzero\\
\bG & \bZ
\end{bmatrix},
\end{equation}
where $\bZ$ is an $n\times (n-m)$ matrix.
Consequently,  $\bG \in \real^{n\times m}$ serves as a right inverse of $\bA$.

More generally, the following two propositions demonstrate   how to find more left inverses or right inverses of a matrix.

\begin{proposition}[Finding Left Inverses]
Let $\bA\in \real^{m\times n}$ be left-invertible (i.e., $m\geq n$ with full column rank). Then,
$$
\bA_L^{-1} = \begin{bmatrix}
(\bA_1^{-1} -\bY \bA_2 \bA_1^{-1}) & \bY
\end{bmatrix}\bE
$$
is a left inverse of $\bA$, where $\bY\in \real^{n\times (m-n)}$ can be any matrix, and 
$\bE\bA = 
\footnotesize\begin{bmatrix}
\bA_1 \\
\bA_2
\end{bmatrix}$ is the result of row elementary transformations applied to $\bA$ such that $\bA_1 \in \real^{n\times n}$ is invertible (since $\bA$ has full column rank $n$) and $\bE\in \real^{m\times m}$.
\end{proposition}
One can verify that $\bG$ in Equation~\eqref{equation:onesided-1} is a specific left inverse of $\bA$ by setting $\bY=\bzero$. 
Since 
$\bE = 
\footnotesize
\begin{bmatrix}
\bG \\
\boxtimes
\end{bmatrix}$, $\bA_1 = \bI_n$, and $\bA_2 = \bzero$, we have 
$
\bA_L^{-1} = 
[(\bA_1^{-1} -\bY \bA_2 \bA_1^{-1}),\, \bY]\bE = \bG + \bY\bZ = \bG,
$
where the last equation follows from the assumption that $\bY=\bzero$.

\begin{proposition}[Finding Right Inverses]
Let $\bA\in \real^{m\times n}$ be right-invertible (i.e., $m\leq n$ with full row rank). Then,
$$
\bA_R^{-1} = \bE\begin{bmatrix}
(\bA_1^{-1}-\bA_1^{-1}\bA_2\bY) \\
\bY
\end{bmatrix}
$$
is a right inverse of $\bA$, where $\bY\in \real^{(n-m)\times m}$ can be any matrix, and $\bA\bE = [\bA_1 , \bA_2]$ is the  result of column elementary transformations applied to  $\bA$ such that $\bA_1 \in \real^{m\times m}$ is invertible (since $\bA$ has full row rank $m$) and $\bE\in \real^{n\times n}$.
\end{proposition}
Similarly, one can verify that $\bG$ in Equation~\eqref{equation:onesided-2} is a specific right inverse of $\bA$ by setting $\bY=\bzero$. Since $\bE = [\bG, \bZ]$, $\bA_1=\bI_m$, and $\bA_2=\bzero$, we have 
$
\bA_R^{-1} = \bE 
\footnotesize\begin{bmatrix}
(\bA_1^{-1}-\bA_1^{-1}\bA_2\bY) \\
\bY
\end{bmatrix} \normalsize= \bG+\bZ\bY = \bG,
$
where again the last equality follows from the assumption that $\bY=\bzero$.

Under specific conditions, the linear system $\bA\bx=\bb$ has a unique solution.
\begin{theorem}[Unique Linear System Solution]\label{theorem:unique-linear-system-solution}
Let $\bA\in \real^{m\times n}$ be left-invertible (i.e., $m\geq n$ with full column rank), and let $\bA_L^{-1}$ be a left inverse of $\bA$. Then, the linear system $\bA\bx = \bb$ has a \textbf{unique} solution if and only if 
$$
(\bI_m - \bA \bA_L^{-1})\bb = \bzero.
$$
And the unique solution is given by 
$
\widehat{\bx} = (\bA^\top\bA)^{-1}\bA^\top \bb.
$
\end{theorem}
\begin{proof}[of Theorem~\ref{theorem:unique-linear-system-solution}]
Let $\bx_0$ be a solution of $\bA\bx =\bb$. Then,
$
\bA \bA_L^{-1} (\bA \bx_0) = \bA\bA_L^{-1} \bb
$
implies 
$
\bA (\bA_L^{-1} \bA) \bx_0 = \bA \bx_0 = \bb.
$
That is, $\bA\bA_L^{-1} \bb=\bb$ and $(\bI_m-\bA\bA_L^{-1})\bb = \bzero$. 

Conversely, suppose $(\bI_m - \bA \bA_L^{-1})\bb = \bzero$, and let $\bx_0 = \bA_L^{-1} \bb$. Substituting $\bx_0 = \bA_L^{-1} \bb$ into $(\bI_m - \bA \bA_L^{-1})\bb = \bzero$, we obtain 
$
\bA\bx_0 = \bb,
$
which implies $\bx_0 = \bA_L^{-1} \bb$ is a solution of $\bA\bx = \bb$ if $(\bI_m - \bA \bA_L^{-1})\bb = \bzero$.

To prove the uniqueness, suppose $\bx_0$ and $\bx_1$ are two solutions of $\bA\bx = \bb$. We have $\bA\bx_0=\bA\bx_1 = \bb$ such that $\bA(\bx_0-\bx_1) = \bzero$. Since $\bA$ is left-invertible so that $\bA$ has full column rank $n$, the dimension of the row space of $\bA$ is $n$ as well such that the null space of $\bA$ is of dimension 0 (i.e., $\dim(\cspace(\bA^\top)) + \dim(\nspace(\bA))=n$ by the fundamental theorem of linear algebra, see Theorem~\ref{theorem:fundamental-linear-algebra}). Therefore, $\bx_0=\bx_1$, which completes the proof.
\end{proof}

\index{Fundamental theorem}
In the fundamental theorem of linear algebra, Figure~\ref{fig:lafundamental}, \textbf{the row space of $\bA$ is the entire space of $\real^n$ if $\bA$ is left-invertible} (i.e., $\bA$ has full column rank $n$). If the condition $(\bI_m - \bA \bA_L^{-1})\bb = \bzero$ is satisfied, then $\bA\bA_L^{-1} \bb=\bb$, it implies that $\bb$ is in the column space of $\bA$. 
Therefore, the linear system $\bA\bx=\bb$ has at least one solution, and the above theorem shows that this solution is unique.

\begin{theorem}[Always Have Solution]\label{theorem:always-have-solution-right-inverse}
Let $\bA\in \real^{m\times n}$ be right-invertible (i.e., $m\leq n$ with full row rank), and let $\bA_R^{-1}$ be a right inverse of $\bA$. Then for any $\bb\in \real^m$, the linear system $\bA\bx = \bb$ has solutions, and the solution is given by
$$
\widehat{\bx} = \bA_R^{-1}\bb, 
$$
where the right inverse $\bA_R^{-1}$ is \textbf{not} necessarily unique.
\end{theorem}
\begin{proof}[of Theorem~\ref{theorem:always-have-solution-right-inverse}]
It can be readily verified that 
$
(\bA \bA_R^{-1}) \bb = \bI_m \bb =\bb,
$
so that $\bA_R^{-1}\bb$ is a solution of $\bA\bx=\bb$.
\end{proof}
We notice that if $\bA$ is right-invertible, then it has full row rank $m$. In the fundamental theorem of linear algebra, Figure~\ref{fig:lafundamental}, \textbf{the column space of $\bA$ is the entire space of $\real^m$ if $\bA$ is right-invertible}. Consequently, any vector $\bb\in \real^m$ is in the column space of $\bA$ such that $\bA\bx=\bb$ has at least one solution. 

\index{Fundamental theorem}
\subsection{Generalized Inverse (g-inverse)}\label{section:g_inv}
We  previously mentioned that if the matrix $\bA$ is nonsingular,  the linear system $\bb = \bA\bx$ can be readily solved using the inverse of $\bA$, such that $\widehat{\bx} = \bA^{-1}\bb$. However, the inverse of an $m\times n$ matrix $\bA$ does not exist if $\bA$ is either not square or  singular. 
Nevertheless, when $\bb$ is in the column space of $\bA$, we can still find the solution to the linear system. 
The relationship between the solution $\widehat{\bx}$ and $\bb$ is given by the generalized inverse of $\bA$: $\widehat{\bx} = \bA^-\bb$.
\begin{definition}[Generalized Inverse\index{Generalized inverse}]
Let $\bA \in \real^{m\times n}$ be any matrix  with rank $r$ where $r \leq \min\{n, m\}$, i.e., could be full-rank or non-full rank. 
Then, a generalized inverse $\bA^- \in\real^{n\times m}$ of $\bA$ is a matrix satisfying
$$
(C1). \qquad  \bA\bA^-\bA = \bA.
$$
That is, $\bA^-\bA $ maps the columns of $\bA$ to themselves.
Or equivalently, for any $\bb \in \cspace(\bA)$, it follows that
$$
(C1.1). \qquad \bA\bA^-\bb = \bb, 
\quad \text{for any $\bb \in \cspace(\bA)$.}
$$
\end{definition}
To see the equivalence between $(C1)$ and $(C1.1)$, that is, we need to show that $\bA$ satisfies $(C1)$ if and only if it satisfies $(C1.1)$.
\paragraph{Forward implication.}
For any $\bb \in \cspace(\bA)$, we can find an $\bx\in\real^n$ such that $\bA\bx = \bb$. If $\bA$ and $\bA^-$ satisfy $(C1)$, then
$$
\bA\bA^-\bA \bx = \bA\bx 
\quad\implies\quad
\bA\bA^-\bb = \bb,
$$
which implies $\bA$ and $\bA^-$ also satisfy $(C1.1)$. 
\paragraph{Reverse implication.}
For the reverse, suppose $\bA$ and $\bA^-$ satisfy $(C1.1)$, then 
$$
\bA\bA^-\bb = \bb 
\quad\implies\quad
\bA\bA^-\bA \bx = \bA\bx ,
$$
which implies $\bA$ and $\bA^-$ also satisfy $(C1)$.

Multiplying both sides of $(C1)$ on the left  by $\bA^-$ and using the definition of the projection matrix in Definition~\ref{definition:projection-matrix}, we obtain $\bA^-\bA\bA^-\bA = \bA^-\bA$ such that $ \bA^-\bA$ is idempotent, which implies $\bA^-\bA$ is a projection matrix (not necessarily  an orthogonal projection).
\begin{theorem}[Projection Matrix from Generalized Inverse]\label{theorem:idempotent-of-ginverse}
Given any matrix $\bA$ and its generalized inverse $\bA^-$, $\bA^-\bA$ is a projection matrix but not necessarily an orthogonal projection. 
The same claim applies to $\bA\bA^-$ as well.
\end{theorem}

\begin{proposition}[Rank of Generalized Inverse]\label{proposition:rank-of-ginverse}
Given any matrix $\bA \in \real^{m\times n}$ and its generalized inverse $\bA^- \in\real^{n\times m}$, then 
$$
\rank(\bA^-) \geq \rank(\bA).
$$
Moreover, we have $\rank(\bA)=\rank(\bA\bA^-)=\rank(\bA^-\bA)$.
\end{proposition}
\begin{proof}[of Proposition~\ref{proposition:rank-of-ginverse}]
From $(C1)$, we notice that $\rank(\bA) = \rank( \bA\bA^-\bA)$. Additionally,
$
\rank( \bA\bA^-\bA) \leq \rank(\bA\bA^-)\leq \rank(\bA^-),
$
where the first inequality holds because the columns of $\bA\bA^-\bA$ are linear combinations of columns of $\bA\bA^-$, and the second inequality holds because the rows of $\bA\bA^-$ are linear combinations of rows of $\bA^-$. This implies that $\rank(\bA^-) \geq \rank(\bA)$.

For the second part, we have 
$
\rank(\bA) \geq \rank(\bA\bA^-) \geq \rank(\bA\bA^-\bA),
$
where the first inequality holds because the columns of $\bA\bA^-$ are linear combinations of the columns of $\bA$, and the second inequality holds because the columns of $\bA\bA^-\bA$ are linear combinations of the columns of $\bA\bA^-$. Utilizing $(C1)$ again, we have $\rank(\bA) = \rank( \bA\bA^-\bA)$, which implies by ``sandwiching" that
$$
\rank(\bA) = \rank(\bA\bA^-) = \rank(\bA\bA^-\bA).
$$ 
Similarly, we also have 
$
\rank(\bA) \geq \rank(\bA^-\bA) \geq \rank(\bA\bA^-\bA),
$
where the first inequality holds because the rows of $\bA^-\bA$ are linear combinations of the rows of $\bA$, and the second inequality holds because the rows of $\bA\bA^-\bA$ are linear combinations of the rows of $\bA^-\bA$. By ``sandwiching" again, we have 
$$
\rank(\bA) = \rank(\bA^-\bA) = \rank(\bA\bA^-\bA),
$$ 
which completes the proof.
\end{proof}

In Proposition~\ref{proposition:one-sided-invertible}, we showed that the left inverse exists if and only if $\bA$ has full column rank, and the right inverse exists if and only if $\bA$ has full row rank.
However,  these conditions are not necessary for generalized inverses. 
When this full-rank condition is satisfied, we further have the following properties for generalized inverses.
\begin{proposition}[Full-Rank Generalized Inverse]\label{proposition:fullrank-ginverse}
Given any matrix $\bA \in \real^{m\times n}$ and its generalized inverse $\bA^- \in\real^{n\times m}$, then,
\begin{enumerate}
\item   $\bA$ has full column rank if and only if $\bA^-\bA = \bI_n$;
\item   $\bA$ has full row rank if and only if $\bA\bA^- = \bI_m$.
\end{enumerate}
\end{proposition}
\begin{proof}[of Proposition~\ref{proposition:fullrank-ginverse}]
Suppose $\bA$ has full column rank. From Proposition~\ref{proposition:rank-of-ginverse}, we know that $\rank(\bA)=\rank(\bA\bA^-)=\rank(\bA^-\bA)$. Therefore, $\rank(\bA^-\bA)=\rank(\bA)=n$, and $\bA^-\bA\in \real^{n\times n}$ is nonsingular. We obtain
$$
\bI_n= (\bA^-\bA)(\bA^-\bA)^{-1} = \bA^-(\bA\bA^-\bA)(\bA^-\bA)^{-1} = \bA^-\bA.
$$
Conversely, suppose $\bA^-\bA = \bI_n$. This implies $\rank(\bA^-\bA)=n$. Since $\rank(\bA^-\bA)=\rank(\bA)$, we have $\rank(\bA)=n$, meaning $\bA$ has full column rank.

Similarly, we can show $\bA$ has full row rank if and only if $\bA\bA^- = \bI_m$.
\end{proof}

\begin{proposition}[Constructing Generalized Inverse]\label{proposition:construct-ginverse}
Given any matrix $\bA \in \real^{m\times n}$  and its generalized inverse $\bA^- \in\real^{n\times m}$, there exists an $n\times m$ matrix $\bL$ such that 
\begin{equation}\label{equation:constructing-ginverse}
\overline{\bA} \triangleq \bA^-+ \bL-\bA^-\bA\bL\bA\bA^-
\end{equation}
is also a generalized inverse of $\bA$. 
Furthermore, for any generalized inverse $\overline{\bA}$, there exists a matrix $\bL$ such that Equation~\eqref{equation:constructing-ginverse} is satisfied.
\end{proposition}
\begin{proof}[of Proposition~\ref{proposition:construct-ginverse}]
Write out the equation
$$
\begin{aligned}
\bA\overline{\bA}\bA &= \bA(\bA^-+ \bL-\bA^-\bA\bL\bA\bA^-)\bA = \bA\bA^-\bA + \bA\bL\bA  -(\bA\bA^-\bA)\bL(\bA\bA^-\bA)\\
&= \bA\bA^-\bA + \bA\bL\bA  -\bA\bL\bA = \bA.
\end{aligned}
$$
Thus, $\overline{\bA}$ is a generalized inverse of $\bA$. 

Now, suppose  that $\bM$ is any generalized inverse of $\bA$, and define $\bL \triangleq \bM-\bA^-$. Recall that $\bA\bM\bA = \bA$. We have: 
$$
\begin{aligned}
\bA^-+ &\bL-\bA^-\bA\bL\bA\bA^- = \bA^-+ (\bM-\bA^-)-\bA^-\bA(\bM-\bA^-)\bA\bA^- \\
&= \bM - \bA^-(\bA\bM\bA)\bA^- + \bA^-(\bA\bA^-\bA)\bA^-
= \bM - \bA^-\bA\bA^- + \bA^-\bA\bA^- 
= \bM,
\end{aligned}
$$
which implies $\bL$ can be constructed for any generalized inverse $\bM$.
\end{proof}

\begin{proposition}[Generalized Inverse Properties]\label{proposition:ginverse-properties}
Given any matrix $\bA \in \real^{m\times n}$ and its generalized inverse $\bA^- \in\real^{n\times m}$, then,
\begin{enumerate}
\item   $(\bA^\top)^- = (\bA^-)^\top$, i.e., $(\bA^-)^\top$ is a generalized inverse of $\bA^\top $;
\item   For any $a\neq0$, $\frac{1}{a} \bA^-$ is a generalized inverse of $a\bA$;
\item   Suppose $\bL\in \real^{m\times m}$ and $\bM\in \real^{n\times n}$ are both invertible. Then, $\bM^{-1}\bA^- \bL^{-1}$ is a generalized inverse of $\bL\bA\bM$;
\item   $\cspace(\bA\bA^-) = \cspace(\bA)$ and $\nspace(\bA^-\bA) = \nspace(\bA)$.
\end{enumerate}
\end{proposition}
\begin{proof}[of Proposition~\ref{proposition:ginverse-properties}]
For (1), from  $(C1)$, $\bA\bA^-\bA = \bA$, we have $\bA^\top (\bA^-)^\top\bA^\top =\bA^\top $ such that $(\bA^-)^\top$ is a generalized inverse of $\bA^\top $.

For (2), it can be easily verified that $(a\bA) (\frac{1}{a} \bA^-) (a\bA) = (a\bA)$ such that $\frac{1}{a} \bA^-$ is a generalized inverse of $a\bA$ for any $a\neq 0$.

For (3), we observe that $(\bL\bA\bM) (\bM^{-1}\bA^- \bL^{-1})(\bL\bA\bM)= \bL\bA\bA^- \bA\bM=\bL\bA\bM$, which implies $\bM^{-1}\bA^- \bL^{-1}$ is a generalized inverse of $\bL\bA\bM$.

For (4), since the columns of $\bA\bA^-$ are linear combinations of the columns of $\bA$, then $\cspace(\bA\bA^-) \subseteq \cspace(\bA)$. And we proved that $\rank(\bA)=\rank(\bA\bA^-)$ in Proposition~\ref{proposition:rank-of-ginverse}. Thus, $\cspace(\bA\bA^-) = \cspace(\bA)$. Similarly, we could prove $\nspace(\bA^-\bA) = \nspace(\bA)$.
\end{proof}

\subsection{Reflexive Generalized Inverse (rg-inverse)}

\begin{figure}[h!]
\centering
\includegraphics[width=0.8\textwidth]{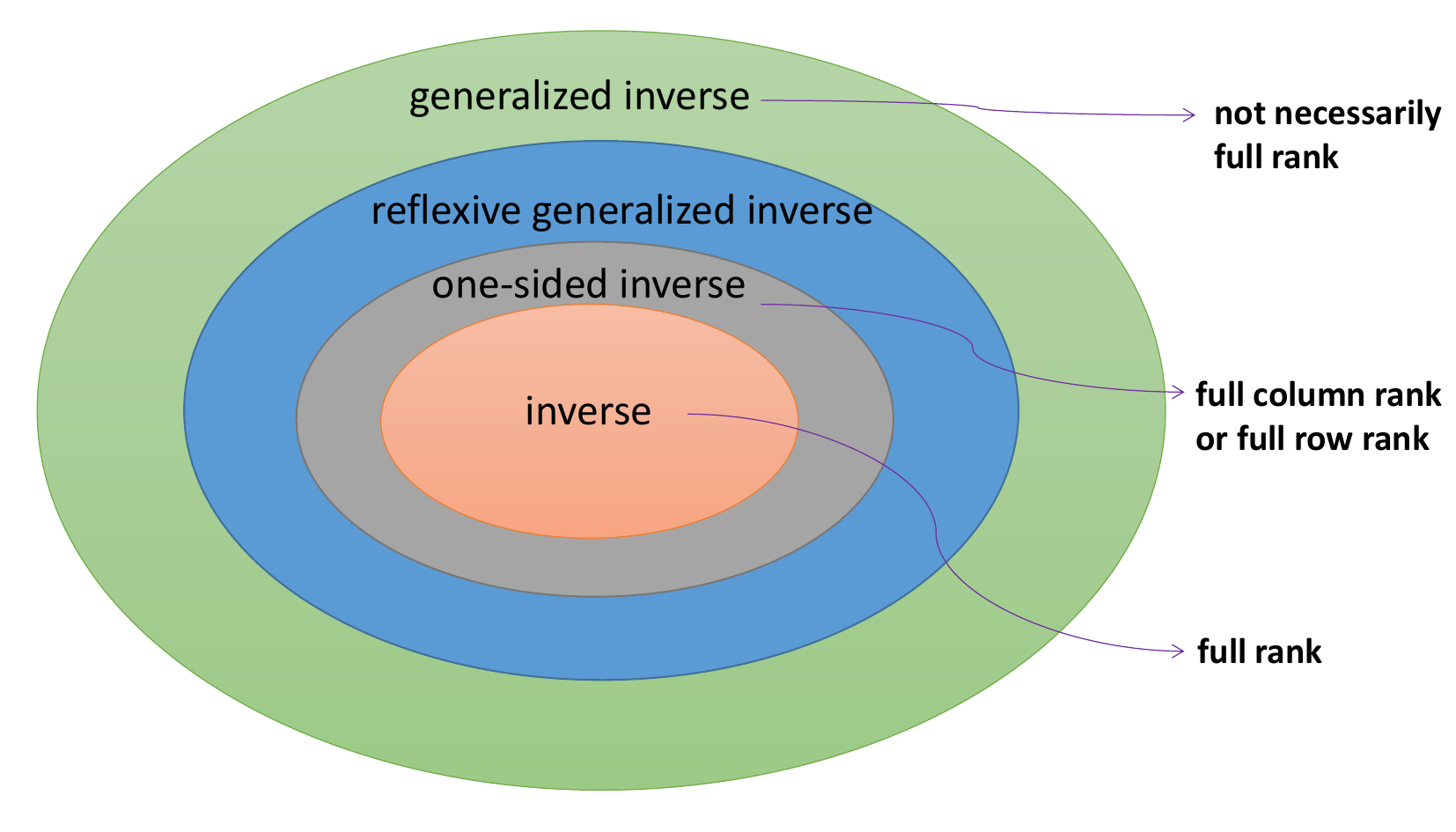}
\caption{Relationship of different inverses: inverse, one-sided inverse, reflexive generalized inverse, and generalized inverse.}
\label{fig:pseudo-inverse-comparison}
\end{figure}
\begin{definition}[Reflexive Generalized Inverse (a.k.a., Middle Inverse)\index{Reflexive generalized inverse}]
Let $\bA \in \real^{m\times n}$ be any matrix with rank $r$ where $r \leq \min\{n, m\}$, i.e., could be full-rank or non-full ranked. Then, a \textbf{reflexive generalized inverse (rg-inverse)} or a \textbf{middle inverse} $\bA_r^- \in\real^{n\times m}$ of $\bA$ is a matrix that satisfies
$$
\begin{aligned}
(C1). \qquad  &\bA\bA_r^-\bA = \bA;\\
(C2). \qquad &\bA_r^- \bA \bA_r^- = \bA_r^-.
\end{aligned}
$$ 
In other words, $\bA_r^-$ is a g-inverse of $\bA$, and $\bA$ is a g-inverse of $\bA_r^-$.
\end{definition}

Suppose $\bA$ has rank $r$. Then, it can be factored as $\bA = \bE_1 
\footnotesize
\begin{bmatrix}
\bI_r & \bzero\\
\bzero & \bzero
\end{bmatrix}
\normalsize\bE_2$, where $\bE_1\in \real^{m\times m}$ and $\bE_2\in \real^{n\times n}$ are elementary transformations on $\bA$.   We can then construct such a reflexive generalized inverse of $\bA$ as 
$$
\bA_r^- = \bE_2^{-1} 
\begin{bmatrix}
\bI_r & \bL \\
\bM  & \bM\bL
\end{bmatrix}
\bE_1^{-1},
$$
where $\bL\in \real^{r\times (m-r)}$, $\bM\in \real^{(n-r)\times r}$ can be any matrices, making the reflexive generalized inverse  \textbf{not unique}. This construction of the reflexive generalized inverse also shows that a reflexive generalized inverse always exists for any matrix, implying that the reflexive generalized inverse is a more general inverse of $\bA$ compared to the one-sided inverse, which may not always exist.
\begin{proposition}[Reflexive Generalized Inverse from G-Inverse]\label{proposition:generalized-for-reflexive-two}
Given any matrix $\bA \in \real^{m\times n}$, and $\bL$ and $\bM$ are two g-inverses of $\bA$, then,
$$
\bZ = \bL\bA\bM 
$$
is a reflexive generalized inverse of $\bA$.
\end{proposition}
It can be easily verified that $\bA\bZ\bA = \bA$ and $\bZ\bA\bZ = \bZ$ for the  above proposition.

\begin{proposition}[Reflexive Generalized Inverse from G-Inverse]\label{proposition:generalized-for-reflexive-two22}
Given any matrix $\bA \in \real^{m\times n}$, then the following two matrices $\bL$ and $\bM$ are two reflexive generalized inverses of $\bA$:
$$
\begin{aligned}
\bL &=(\bA^\top\bA)^-\bA^\top
\quad \text{and}\quad 
\bM =\bA^\top (\bA\bA^\top)^-,
\end{aligned}
$$
where $(\bA^\top\bA)^-$ is a g-inverse of $(\bA^\top\bA)$, and $(\bA\bA^\top)^-$ is a g-inverse of $(\bA\bA^\top)$.
\end{proposition}
\begin{proof}[of Proposition~\ref{proposition:generalized-for-reflexive-two22}]
To prove the result,  we first establish the following result.
\paragraph{Claim: $\cspace(\bA^\top \bA) = \cspace(\bA^\top)$  {and}  $\nspace(\bA^\top \bA)=\nspace(\bA)$.}
Since the columns of $\bA^\top \bA$ are linear combinations of the columns of $\bA^\top$, we have $\cspace(\bA^\top \bA) \subseteq \cspace(\bA^\top)$.
In Lemma~\ref{lemma:rank-of-ata}, we proved that $\rank(\bA^\top \bA)=\rank(\bA)$. This implies $\rank(\bA^\top \bA)=\rank(\bA^\top)$ and $\cspace(\bA^\top \bA) = \cspace(\bA^\top)$. Furthermore, the orthogonal complement of $\cspace(\bA^\top)$ is $\nspace(\bA)$, and the orthogonal complement of $\cspace(\bA^\top \bA)$ is $\nspace(\bA^\top \bA)$. Therefore, by the fundamental theorem of linear algebra in Appendix~\ref{appendix:fundamental-rank-nullity}, we have
$$
\cspace(\bA^\top \bA) = \cspace(\bA^\top)
\quad \text{and}\quad
\nspace(\bA^\top \bA)=\nspace(\bA).
$$
Next, there exists  a set of vectors $\bz_1, \bz_2, \ldots, \bz_n \in \real^n$ such that the  column-$i$ of $\bA^\top$ can be expressed as $\bA^\top\bA\bz_i$. That is, for $\bZ=[\bz_1, \bz_2, \ldots, \bz_n]$, we have
$$
\bA^\top = \bA^\top\bA\bZ
\quad\implies\quad
\bA\bL\bA = (\bA^\top\bA\bZ)^\top  (\bA^\top\bA)^-\bA^\top \bA = \bZ^\top \bA^\top\bA (\bA^\top\bA)^-\bA^\top\bA.
$$
By the condition $(C1.1)$ of the g-inverse, we have $\bA^\top\bA (\bA^\top\bA)^- \bb = \bb $ for any $\bb\in\cspace(\bA^\top\bA)$. This implies  $\bA^\top\bA (\bA^\top\bA)^-\bA^\top\bA = \bA^\top\bA$ and 
\begin{equation}\label{equation:reflex-eq1}
\bA\bL\bA = (\bA^\top\bA\bZ)^\top  (\bA^\top\bA)^-\bA^\top \bA = \bZ^\top \bA^\top\bA=\bA.
\end{equation}
Next, we consider $\bL\bA\bL$:
$
\bL\bA\bL = (\bA^\top\bA)^-\bA^\top \bA(\bA^\top\bA)^-\bA^\top.
$
Same argument can be applied to $\bA^\top \bA(\bA^\top\bA)^-\bA^\top = \bA^\top$. Thus, 
\begin{equation}\label{equation:reflex-eq2}
\bL\bA\bL = (\bA^\top\bA)^-\bA^\top = \bL.
\end{equation}
Combining Equation~\eqref{equation:reflex-eq1} and Equation~\eqref{equation:reflex-eq2}, we conclude that $\bL$ is a reflexive generalized inverse of $\bA$. Similarly, the same process can be applied to show that $\bM$ is also a reflexive generalized inverse of $\bA$.
\end{proof}
\index{Fundamental theorem}

From the definition, we realize that the reflexive generalized inverse is a special type of generalized inverse. Under certain conditions, the two types of inverses become equivalent.
\begin{proposition}[Reflexive Generalized Inverse in G-Inverse]\label{proposition:reflexive-from-ginverse}
Given any matrix $\bA \in \real^{m\times n}$  and its generalized inverse $\bA^-\in \real^{n\times m}$, then $\bA^-$ is a reflexive generalized inverse of $\bA$ if and only if $\rank(\bA)=\rank(\bA^-)$.
\end{proposition}
\begin{proof}[of Proposition~\ref{proposition:reflexive-from-ginverse}]
Suppose $\bA^-$ is a generalized inverse of $\bA$, then $\bA\bA^-\bA=\bA$. Suppose further that $\bA^-$ is also a reflexive generalized inverse, then $\bA^-\bA\bA^- = \bA^-$. We want to show $\rank(\bA)=\rank(\bA^-)$. Firstly, we have 
$$
\begin{aligned}
\rank(\bA) &= \rank(\bA\bA^-\bA)\leq \rank(\bA^-) = \rank(\bA^-\bA\bA^-)\leq \rank(\bA),\\
\end{aligned}
$$
where the two inequalities are from Proposition~\ref{proposition:rank-of-ginverse}. This implies $\rank(\bA)=\rank(\bA^-)$ by ``sandwiching."

Conversely, suppose $\bA^-$ is a generalized inverse of $\bA$ and $\rank(\bA)=\rank(\bA^-)$, then $\bA\bA^-\bA=\bA$. We want to show $\bA^-$ is also a reflexive generalized inverse. First, we have
$$
\rank(\bA) = \rank(\bA\bA^-\bA) \leq \rank(\bA^-\bA) \leq \rank(\bA^-) = \rank(\bA),
$$
where the first inequality holds because the rows of $\bA\bA^-\bA$ are linear combinations of the rows of $\bA^-\bA$, and the second inequality holds because  the columns of $\bA^-\bA$ are linear combinations of the columns of $\bA^-$. Again by ``sandwiching," we have
$$
\rank(\bA) = \rank(\bA\bA^-\bA) = \rank(\bA^-\bA) = \rank(\bA^-) = \rank(\bA),
$$
This equality $\rank(\bA^-\bA) = \rank(\bA^-)$ and the subspace inequality $\cspace(\bA^-\bA)\subseteq \cspace(\bA^-)$ imply $\cspace(\bA^-\bA)=\cspace(\bA^-)$. 
Consequently. there exists  a set of vectors $\balpha_1, \balpha_2, \ldots, \balpha_m \in \real^n$ such that the column-$i$ of $\bA^-$ can be expressed as $\bA^-\bA\balpha_i$. That is, for $\bL\triangleq[\balpha_1, \balpha_2, \ldots, \balpha_m]$, we have
$
\bA^- = \bA^-\bA\bL.
$
Since $\bA=\bA\bA^-\bA$, it follows that 
$
\bA=\bA\bA^-\bA = \bA(\bA^-\bA\bL)\bA = \bA\bL\bA,
$
where the last equality follows form the condition $(C1)$ of the g-inverse; and thus, $\bL$ is a generalized inverse of $\bA$. From Proposition~\ref{proposition:generalized-for-reflexive-two}, $\bA^-=\bA^-\bA\bL$ is a reflexive generalized inverse of $\bA$, which completes the proof.
\end{proof}

From the above results, we can conclude the rank in reflexive generalized inverses.
\begin{proposition}[Rank of Reflexive Generalized Inverse]
Given any matrix $\bA \in \real^{m\times n}$ and its reflexive generalized inverse $\bA_r^- \in\real^{n\times m}$.
Combining the result in Proposition~\ref{proposition:reflexive-from-ginverse} and the rank properties  of g-inverses in Proposition~\ref{proposition:rank-of-ginverse}, we have
$$
\rank(\bA_r^-)=\rank(\bA)=\rank(\bA\bA_r^-)=\rank(\bA_r^-\bA).
$$
\end{proposition}

\begin{proposition}[Reflexive Generalized Inverse Properties]\label{proposition:refle-ginverse-properties}
Given any matrix $\bA \in \real^{m\times n}$ and its reflexive generalized inverse $\bA_r^- \in\real^{n\times m}$, the following properties hold:
\begin{enumerate}
	\item   $\cspace(\bA\bA_r^-) = \cspace(\bA)$ and $\nspace(\bA_r^-\bA) = \nspace(\bA)$.
	\item   $\cspace(\bA_r^-\bA) = \cspace(\bA_r^-)$ and $\nspace(\bA\bA_r^-) = \nspace(\bA_r^-)$. 
\end{enumerate}
\end{proposition}

\begin{proof}[of Proposition~\ref{proposition:refle-ginverse-properties}]
Suppose $\bA^-$ is a g-inverse of $\bA$,
we proved in Proposition~\ref{proposition:ginverse-properties} that $\cspace(\bA\bA^-) = \cspace(\bA)$ and $\nspace(\bA^-\bA) = \nspace(\bA)$. Since $\bA^-_r$ is a g-inverse of $\bA$, and $\bA$ is a g-inverse of $\bA^-_r$, the properties follow directly.
\end{proof}

\index{Pseudo-inverse}
\index{Moore-Penrose pseudo-inverse}
\subsection{Pseudo-Inverse}\label{appendix:pseudo-inverse_main}
\begin{figure}[h!]
\centering
\includegraphics[width=0.98\textwidth]{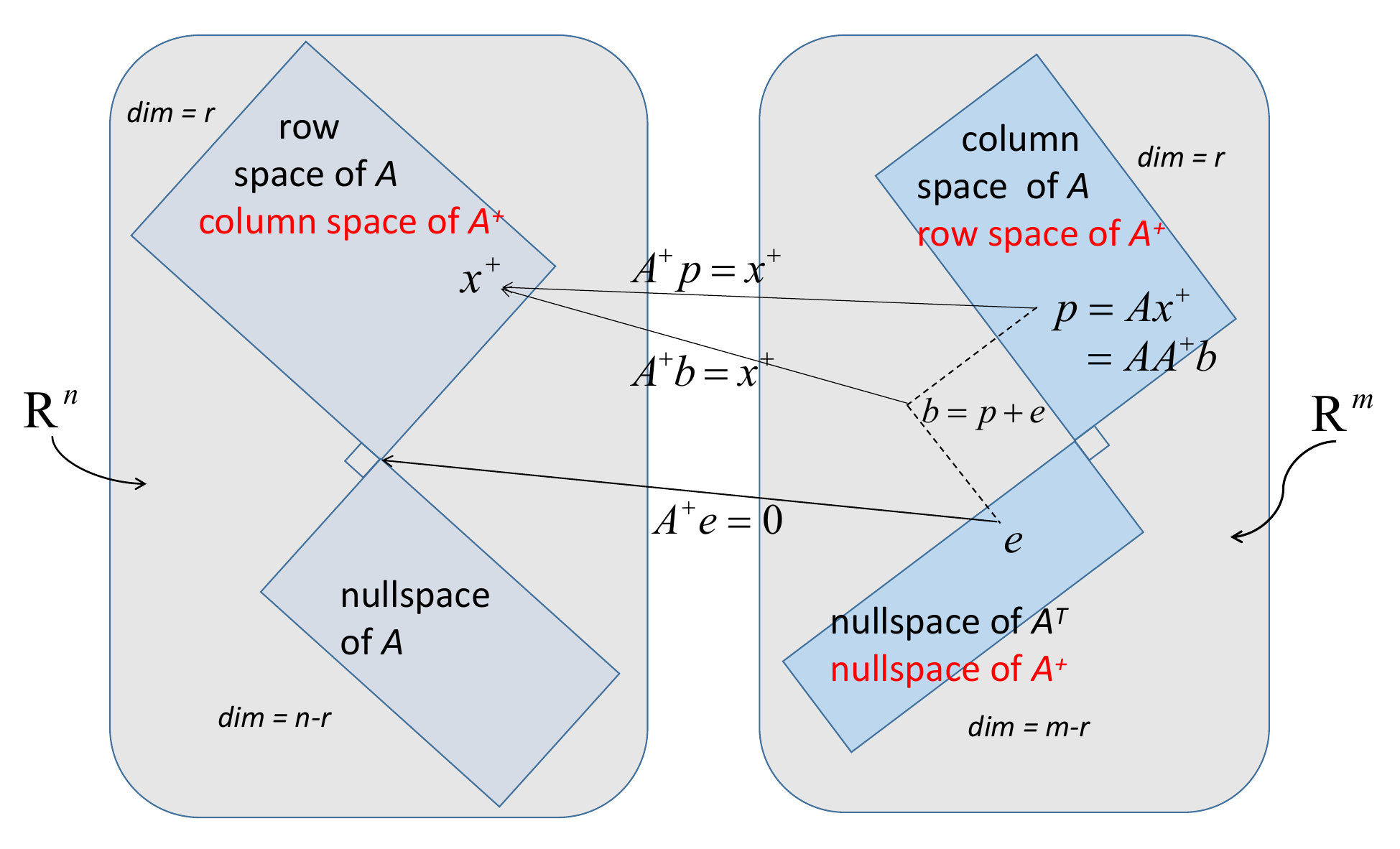}
\caption{Column space and row space of the pseudo-inverse $\bA^+$. $\bA$ maps vectors from the row space to the column space. 
Conversely, $\bA^+$ maps vectors from the column space to the row space. The decomposition of $\bb$ into $\bp+\be$ and the transformation to $\bx^+$ are discussed in Section~\ref{section:pseudo-in-svd}. \textbf{This is a more detailed picture of the pseudo-inverse compared to Figure~\ref{fig:lafundamental4-LS-SVD}}.}
\label{fig:lafundamental5-pseudo}
\end{figure}

As we mentioned previously, for a matrix $\bA\in \real^{m\times n}$, we can find its \textit{pseudo-inverse}, an $n\times m$ matrix denoted by $\bA^+$.
In simple terms, when $\bA$ multiplies a vector $\bx$ in its row space, this produces $\bA\bx$ in the column space (see Figure~\ref{fig:lafundamental}). Both these two spaces have the same dimension $r$, which is the rank of $\bA$. 
$\bA$ is always ``invertible" when restricted to these spaces, and $\bA^+$ serves as the ``inverse" of $\bA$. \textbf{That is, $\bA^+\bA\bx = \bx$ when $\bx$ is in the row space of $\bA$.} And \textbf{$\bA\bA^+\bb = \bb$ when $\bb$ is in the column space of $\bA$} (see Figure~\ref{fig:lafundamental5-pseudo}).

\textbf{The null space of $\bA^+$ is the null space of $\bA^\top$}. It consists of  vectors $\bb$ in $\real^m$ such that $\bA^\top\bb = \bzero$. 
These vectors $\bb$ are orthogonal to every $\bA\bx$ in the column space. 
The proof of this property is deferred to Theorem~\ref{theorem:pseudo-four-basis-space}.

More formally, the pseudo-inverse, or also known as the \textit{Moore-Penrose pseudo-inverse}, $\bA^+$, is defined as the \textbf{unique} $n\times m$ matrix satisfying the following four criteria:
\begin{equation}\label{equation:pseudi-four-equations}
\boxed{
\begin{aligned}
&(C1). \qquad  \bA\bA^+\bA &=& \bA; \qquad &(\bA^+\text{ is a g-inverse of }\bA) \\
&(C2). \qquad  \bA^+\bA\bA^+ &=&\bA^+ ;\qquad &(\bA\text{ is a g-inverse of }\bA^+)\\
&(C3). \qquad  (\bA\bA^+)^\top &=&\bA\bA^+;\\
&(C4). \qquad  	(\bA^+\bA)^\top &=& \bA^+\bA.
\end{aligned}
}
\end{equation}

In Theorem~\ref{theorem:idempotent-of-ginverse}, we established that $\bA\bA^+$ and $\bA^+\bA$ are idempotent if $\bA^+$ is a g-inverse of $\bA$, and thus they are both projection matrices. 
Since $\bA^+$ is the pseudo-inverse of $\bA$ \footnote{We use ``the pseudo-inverse" rather than ``a pseudo-inverse" since the pseudo-inverse is unique, as we will prove shortly.},  conditions $(C3)$ and $(C4)$ ensure that  $\bA\bA^+$ and $\bA^+\bA$ are symmetric such that they are \textbf{orthogonal} projections as well (see Proposition~\ref{proposition:symmetric-projection-matrix}).

The existence of the pseudo-inverse for any matrix is further supported by the pseudo-inverse obtained from the CR decomposition of the matrix.
\begin{theorem}[Existence of Pseudo-Inverse]\label{theorem:existence-of-pseudo-inverse}
Every matrix $\bA \in \real^{m\times n}$ has a pseudo-inverse.
\end{theorem}
\begin{proof}[of Theorem~\ref{theorem:existence-of-pseudo-inverse}]
Given the CR decomposition of $\bA=\bC\bR$ (Theorem~\ref{theorem:cr-decomposition}), let 
$$
\bA^+ = \bR^+\bC^+ = \bR^\top (\bR\bR^\top)^{-1} (\bC^\top\bC)^{-1}\bC^\top,
$$
where $\bR^+=\bR^\top (\bR\bR^\top)^{-1}$ and $\bC^+ =(\bC^\top\bC)^{-1}\bC^\top$. \footnote{It can be easily verified that $\bR^+$ is the pseudo-inverse of $\bR$ and $\bC^+$ is the pseudo-inverse of $\bC$.} $\bR\bR^\top$ and $\bC^\top\bC$ are invertible since $\bC\in \real^{m\times r}$ and $\bR\in \real^{r\times n}$ have full rank $r$ from the property of the CR decomposition. 
We can verify that $\bA^+$ satisfies the four criteria for the pseudo-inverse:
$$
\begin{aligned}
&(C1). \qquad \bA\bA^+\bA &=& \bC\bR\left(\bR^\top (\bR\bR^\top)^{-1} (\bC^\top\bC)^{-1}\bC^\top\right)\bC\bR = \bC\bR = \bA, \\
&(C2). \qquad \bA^+\bA\bA^+ &=&\left(\bR^\top (\bR\bR^\top)^{-1} (\bC^\top\bC)^{-1}\bC^\top\right) \bC\bR\left(\bR^\top (\bR\bR^\top)^{-1} (\bC^\top\bC)^{-1}\bC^\top\right)\\
&&=&\bR^\top(\bR\bR^\top)^{-1}(\bC^\top\bC)^{-1}\bC^\top = \bA^+,\\
&(C3). \qquad (\bA\bA^+)^\top &=& \bC(\bC^\top\bC)^{-1} (\bR\bR^\top)^{-1} \bR\bR^\top\bC^\top = \bC(\bC^\top\bC)^{-1} \bC^\top\\
&&=&\bC\bR\bR^\top (\bR\bR^\top)^{-1} (\bC^\top\bC)^{-1}\bC^\top = \bA\bA^+,\\
&(C4). \qquad (\bA^+\bA)^\top &=& \bR^\top\bC^\top \bC(\bC^\top\bC)^{-1} (\bR\bR^\top)^{-1} \bR= \bR^\top (\bR\bR^\top)^{-1} \bR\\
&&=& \bR^\top (\bR\bR^\top)^{-1} (\bC^\top\bC)^{-1}\bC^\top\bC\bR = \bA^+\bA.
\end{aligned}
$$
Thus, $\bA^+$ is the pseudo-inverse of $\bA$, proving the existence of the pseudo-inverse.
\end{proof}

\begin{theorem}[Uniqueness of Pseudo-Inverse]\label{theorem:uniqueness-of-pseudo-inverse}
Every matrix $\bA$ has a unique pseudo-inverse.
\end{theorem}
\begin{proof}[of Theorem~\ref{theorem:uniqueness-of-pseudo-inverse}]
Suppose $\bA_1^+$ and $\bA_2^+$ are two pseudo-inverses of $\bA$. Then
$$
\begin{aligned}
\bA_1^+ &= \bA_1^+\bA\bA_1^+ = \bA_1^+(\bA\bA_2^+\bA)\bA_1^+ = \bA_1^+ (\bA\bA_2^+)(\bA\bA_1^+) \qquad &(\text{by $(C2), (C1)$})\\
&=\bA_1^+ (\bA\bA_2^+)^\top(\bA\bA_1^+)^\top =  \bA_1^+\bA_2^{+\top}\bA^\top \bA_1^{+\top}\bA^\top \qquad &(\text{by $(C3)$})\\
&=\bA_1^+\bA_2^{+\top}(\bA \bA_1^{+}\bA)^\top= \bA_1^+\bA_2^{+\top}\bA^\top \qquad &(\text{by $(C1)$})\\
&= \bA_1^+(\bA \bA_2^{+})^\top =\bA_1^+ \bA \bA_2^{+} = \bA_1^+ (\bA \bA_2^+\bA) \bA_2^{+}  \qquad &(\text{by $(C3),(C1)$})\\
&= (\bA_1^+ \bA) (\bA_2^+\bA) \bA_2^{+} =  (\bA_1^+ \bA)^\top (\bA_2^+\bA)^\top \bA_2^{+}  \qquad &(\text{by $(C4)$})\\
&=(\bA \bA_1^{+} \bA)^\top \bA_2^{+\top} \bA_2^{+} = \bA^\top \bA_2^{+\top} \bA_2^{+}   \qquad &(\text{by $(C1)$})\\
&=  (\bA_2^{+}\bA)^\top \bA_2^{+} =\bA_2^{+}\bA \bA_2^{+} = \bA_2^{+}, \qquad &(\text{by $(C4),(C2)$})\\
\end{aligned}
$$
which implies the uniqueness of the pseudo-inverse.
\end{proof}


We are now ready to discuss the four subspaces associated with the pseudo-inverse.
\begin{theorem}[Four Subspaces in Pseudo-Inverse]\label{theorem:pseudo-four-basis-space} 
Given the pseudo-inverse $\bA^+$ of $\bA$, we have the following properties:
\begin{enumerate}
\item The column space of $\bA^+$ is  identical to the row space of $\bA$;
\item The row space of $\bA^+$ is  identical to the column space of $\bA$; 
\item The null space of $\bA^+$ is  identical to the null space of $\bA^\top$;
\item The null space of $\bA^{+\top}$ is  identical to the null space of $\bA$.
\end{enumerate}
\item The relationships among these subspaces are illustrated in Figure~\ref{fig:lafundamental5-pseudo}.
\end{theorem}
\begin{proof}[of Theorem~\ref{theorem:pseudo-four-basis-space}]
Since $\bA^+$ is a special type of rg-inverse, by Proposition~\ref{proposition:refle-ginverse-properties}, we have
$$
\begin{aligned}
\cspace(\bA\bA^+) &= \cspace(\bA) \qquad \text{and} \qquad \nspace(\bA^+\bA) = \nspace(\bA);\\
\cspace(\bA^+\bA) &= \cspace(\bA^+)\qquad \text{and} \qquad \nspace(\bA\bA^+) = \nspace(\bA^+).
\end{aligned}
$$
By conditions $(C3)$ and $(C4)$ of the definition of pseudo-inverses, we also have 
$$
(\bA^+\bA)^\top = \bA^+\bA \qquad \text{and} \qquad (\bA\bA^+)^\top =\bA\bA^+.
$$
According to the fundamental theorem of linear algebra, we realize that $\cspace(\bA\bA^+)$ is the orthogonal complement to $\nspace((\bA\bA^+)^\top)$, and $\cspace(\bA^+\bA)$ is the orthogonal complement to $\nspace((\bA^+\bA)^\top)$: 
$$
\begin{aligned}
\cspace(\bA\bA^+) \perp \nspace((\bA\bA^+)^\top) &\quad\implies\quad \cspace(\bA\bA^+) \perp \nspace(\bA\bA^+); \\
\cspace(\bA^+\bA) \perp \nspace((\bA^+\bA)^\top) &\quad\implies\quad \cspace(\bA^+\bA) \perp \nspace(\bA^+\bA).
\end{aligned}
$$
This implies 
$$
\begin{aligned}
\cspace(\bA) \perp \nspace(\bA^+) \qquad \text{and} \qquad \cspace(\bA^+)\perp \nspace(\bA).
\end{aligned}
$$
That is, $\nspace(\bA^+) = \nspace(\bA^\top)$ and $\cspace(\bA^+) = \cspace(\bA^\top)$. According to  the fundamental theorem of linear algebra, this also implies $\cspace(\bA^{+\top})=\cspace(\bA)$ and $\nspace(\bA^{+\top})=\nspace(\bA)$. This completes the proof.
\end{proof}
\index{Fundamental theorem}

To conclude, we compare the properties of different inverses of $\bA$ in Table~\ref{table:different-inverses}.
\begin{table}[h!]
\setlength{\tabcolsep}{7pt}
\begin{tabular}{c|l|l|l}
\hline
& \multicolumn{1}{c|}{g-inverse}                                                                                   & \multicolumn{1}{c|}{rg-inverse}                                                                                                                                                                             & pseudo-inverse                                                                                                                                                                                                                                                                \\ \hline
subspaces & \begin{tabular}[c]{@{}l@{}}$\cspace(\bA\bA^-) = \cspace(\bA)$\\ $\nspace(\bA^-\bA) = \nspace(\bA)$\end{tabular}  & \begin{tabular}[c]{@{}l@{}}$\cspace(\bA\bA_r^-) = \cspace(\bA)$ \\ $\nspace(\bA_r^-\bA) = \nspace(\bA)$\\ $\cspace(\bA_r^-\bA) = \cspace(\bA_r^-)$ \\ $\nspace(\bA\bA_r^-) = \nspace(\bA_r^-)$\end{tabular} & \begin{tabular}[c]{@{}l@{}}$\cspace(\bA\bA^+) = \cspace(\bA)=\cspace(\bA^{+\top})$ \\ $\nspace(\bA^+\bA) = \nspace(\bA)=\nspace(\bA^{+\top})$\\ $\cspace(\bA^+\bA) = \cspace(\bA^+)=\cspace(\bA^\top)$ \\ $\nspace(\bA\bA^+) = \nspace(\bA^+)=\nspace(\bA^\top)$\end{tabular} \\ \hline \hline
rank      & \begin{tabular}[c]{@{}l@{}}$\rank(\bA\bA^-)$\\ $=\rank(\bA^-\bA)$\\ $=\rank(\bA)$\\ $\leq \rank(\bA^-)$\end{tabular} & \begin{tabular}[c]{@{}l@{}}$\rank(\bA_r^-)$\\ $=\rank(\bA)$\\ $=\rank(\bA\bA_r^-)$\\ $=\rank(\bA_r^-\bA)$\end{tabular}                                                                                          & \begin{tabular}[c]{@{}l@{}}$\rank(\bA^+)$\\ $=\rank(\bA)$\\ $=\rank(\bA\bA^+)$\\ $=\rank(\bA^+\bA)$\end{tabular}                                                                                                                                                                  \\ \hline
\end{tabular}\caption{Comparison of different inverses in terms of subspace and rank.}\label{table:different-inverses}
\end{table}

In Theorem~\ref{theorem:idempotent-of-ginverse}, we stated that $\bA^-\bA$ and $\bA\bA^-$ are projection matrices, where $\bA^-$ is a generalized inverse of $\bA$ (but not necessarily orthogonal projections). 
However, when considering the pseudo-inverse, we find orthogonal projections.
\begin{theorem}[Orthogonal Projection onto Column Space and Row Space]\label{theorem:orthogonal-from-pseudo-inverse}
Given any matrix $\bA \in \real^{m\times n}$ and its pseudo-inverse $\bA^+ \in \real^{n\times m}$, $\bH = \bA\bA^+$ is the \textit{orthogonal projection} onto the column space of $\bA$. Similarly, $\bP= \bA^+\bA$ is the \textit{orthogonal projection} onto the row space of $\bA$.
\end{theorem}
\begin{proof}[of Theorem~\ref{theorem:orthogonal-from-pseudo-inverse}]
Since  $\bH^\top = (\bA\bA^+)^\top=\bA\bA^+=\bH$ from the definition of the pseudo-inverse, and $\bH$ is idempotent such that $\bH$ is an {orthogonal projection}.
From Table~\ref{table:different-inverses}, we conclude that $\cspace(\bH)=\cspace(\bA\bA^+)=\cspace(\bA)$. This implies $\bH$ is the {orthogonal projection} onto the column space of $\bA$. Similarly, we can prove $\bP= \bA^+\bA$ is the {orthogonal projection} onto the row space of $\bA$.
\end{proof}

\subsection{Pseudo-Inverse with Varying Dimensions}
In conclusion, let us examine the pseudo-inverse of a matrix $\bA$ in $\real^{m\times n}$ with varying dimensions.
Specifically, we suppose $\bA$ has rank $r$. And we define the pseudo-inverse in either of the following ways:
\paragraph{Case $m>n=r$.} That is, matrix $\bA$ has linearly independent columns: 
$\bA^\top\bA$ is an $n\times n$ invertible matrix, and we define the left-pseudo-inverse:
$$
{\text{left-pseudo-inverse} = \bA^+ = (\bA^\top\bA)^{-1}\bA^\top},
$$
which satisfies
\begin{equation}
\bA^+ \bA =  (\bA^\top\bA)^{-1}\bA^\top\bA = \bI_n
\quad\text{and}\quad
\bA\bA^+ = \bA(\bA^\top\bA)^{-1}\bA^\top \neq\bI.
\end{equation}

\paragraph{Case $n>m=r$.} That is, matrix $\bA$ has linearly independent rows: 
$\bA\bA^\top$ is an $m\times m$ invertible matrix, and we define the right-pseudo-inverse:
$$
{\text{right-pseudo-inverse} = \bA^+ = \bA^\top(\bA\bA^\top)^{-1}},
$$
which satisfies  
\begin{equation}\label{equation:right-inverse-ba-appendix}
\bA\bA^+ = \bA\bA^\top(\bA\bA^\top)^{-1} = \bI_m
\quad\text{and}\quad
\bA^+ \bA =  \bA^\top(\bA\bA^\top)^{-1}\bA \neq \bI.
\end{equation}

\paragraph{Case rank-deficient.} We delay the pseudo-inverse for rank-deficient matrices in the next section using singular value decomposition.

\paragraph{Pseudo-inverse of pseudo-inverse.}
Additionally, we can also show that $(\bA^+)^+ = \bA$. If $m>n=r$, we have 
$$
\begin{aligned}
(\bA^+)^+ &= [(\bA^\top\bA)^{-1}\bA^\top]^+
= \bA^{+\top}(\bA^+\bA^{+\top})^{-1} &(m>n=r) \\
&= \left[(\bA^\top\bA)^{-1}\bA^\top\right]^{\top}\left\{\left[(\bA^\top\bA)^{-1}\bA^\top\right]\left[(\bA^\top\bA)^{-1}\bA^\top\right]^{\top}\right\}^{-1}\\
&=\bA (\bA^\top\bA)^{-1}\left\{(\bA^\top\bA)^{-1}\bA^\top\bA(\bA^\top\bA)^{-1}   \right\}^{-1}
=\bA.
\end{aligned}
$$
Similarly, we can show $(\bA^+)^+ = \bA$ if $n>m=r$.

In particular, when $m=n$, $\bA$ is a square invertible matrix, then both left- and right-pseudo-inverse are the inverse of $\bA$:
$$
\begin{aligned}
\text{left-pseudo-inverse} &= \bA^+ = (\bA^\top\bA)^{-1}\bA^\top=\bA^{-1} \bA^{-\top} \bA^\top = \bA^{-1}, \\
\text{right-pseudo-inverse} &= \bA^+ = \bA^\top(\bA\bA^\top)^{-1} = \bA^\top\bA^{-\top}\bA^{-1} = \bA^{-1}.
\end{aligned}
$$

\begin{exercise}
Given a matrix $\bA\in\real^{m\times n}$ and a positive constant $\gamma>0$, it is trivial $\bA^\top\bA+\gamma\bI$ is nonsingular. Show that the limit of the matrix $(\bA^\top\bA+\gamma\bI)^{-1}\bA^\top$ is equal to the pseudo-inverse $\bA^+$ of $\bA$ as $\gamma$ approaches to zero.
\end{exercise}

\subsection{Pseudo-Inverse in SVD}\label{section:pseudo-in-svd}
When considering the full SVD $\bA = \bU\bSigma\bV^\top$ of a matrix $\bA\in\real^{m\times n}$, where $\bU\in\real^{m\times m}$ and $\bV\in\real^{n\times n}$ are orthogonal matrices, and diagonal of $\bSigma\in\real^{m\times n}$ contains the singular values of $\bA$, we discuss the pseudo-inverse of the matrix in different cases:
\paragraph{Case $m>n=r$.} Let $\bSigma^+\triangleq(\bSigma^\top\bSigma)^{-1}\bSigma^\top$. The left-pseudo-inverse can be obtained by 
$$
\begin{aligned}
\text{left-pseudo-inverse} 
&= \bA^+ = (\bA^\top\bA)^{-1}\bA^\top 
= (\bV \bSigma^\top \bU^\top\bU\bSigma\bV^\top)^{-1}\bV \bSigma^\top \bU^\top\\
&= \bV (\bSigma^\top\bSigma)^{-1}\bV^\top \bV \bSigma^\top \bU^\top 
= \bV[(\bSigma^\top\bSigma)^{-1}\bSigma^\top] \bU^\top 
= \bV \bSigma^+ \bU^\top.  
\end{aligned}
$$

\paragraph{Case $n>m=r$.} Let $\bSigma^+\triangleq\bSigma^\top(\bSigma\bSigma^\top)^{-1}$. The right-pseudo-inverse can be obtained by 
$$
\begin{aligned}
\text{right-pseudo-inverse} 
&= \bA^+ = \bA^\top(\bA\bA^\top)^{-1}
= (\bU\bSigma\bV^\top)^\top[(\bU\bSigma\bV^\top)(\bU\bSigma\bV^\top)^\top]^{-1}\\
&= \bV\bSigma^\top\bU^\top(\bU\bSigma\bV^\top\bV\bSigma^\top\bU^\top)^{-1} 
= \bV\bSigma^\top\bU^\top\bU^{-\top}(\bSigma\bSigma^\top)^{-1}\bU^{-1} \\
&= \bV\bSigma^\top(\bSigma\bSigma^\top)^{-1}\bU^{-1} 
=\bV \bSigma^+ \bU^\top. 
\end{aligned}
$$

\paragraph{Case rank-deficient.} $\bA^+=\bV \bSigma^+ \bU^\top$,
where the upper-left side of $\bSigma^+ \in \real^{n\times m}$ is a diagonal matrix $\diag(\frac{1}{\sigma_1}, \frac{1}{\sigma_2}, \ldots, \frac{1}{\sigma_r})$. It can be easily verified that this definition of the pseudo-inverse satisfies the four criteria in Equation~\eqref{equation:pseudi-four-equations}.

In all cases, $\bSigma^+$ is  the pseudo-inverse of $\bSigma$ with $1/\sigma_1, 1/\sigma_2, \ldots, 1/\sigma_r$ on its diagonal. We thus conclude the pseudo-inverse from SVD as summarized in  Table~\ref{table:pseudo-inverse-svd}.
\begin{table}[h!]\centering
\setlength{\tabcolsep}{11.4pt}
\begin{tabular}{c|c|c|c|c|c}
\hline
\addlinespace[0.03cm]
& $\bA$                & $\bA^\top$           & $\bA^+$                & $\bA^{+\top}$  &   $(\bA^{+})^+$      \\ \hline\hline
\addlinespace[0.03cm]
SVD & $\bU\bSigma\bV^\top$ & $\bV\bSigma\bU^\top$ & $\bV\bSigma^+\bU^\top$ & $\bU\bSigma^+\bV^\top$& $\bU\bSigma\bV^\top=\bA$ \\ \hline
\end{tabular}
\caption{SVD for the transpose and pseudo-inverse of a matrix.}
\label{table:pseudo-inverse-svd}
\end{table}

From the pseudo-inverse via the SVD, we can provide an alternative way to understand the orthogonal projection in pseudo-inverses.
\begin{mdframed}[hidealllines=\mdframehidelineNote,backgroundcolor=\mdframecolor,frametitle={Another Way to See the Orthogonal Projection in Pseudo-Inverse via SVD}]

Given $\bA\in\real^{m\times n}$, we have  previously shown that $\bH=\bA\bA^+\in\real^{m\times m}$ is an orthogonal projection. We now need to show that it projects onto $\cspace(\bA)$. For any vector $\bb\in \real^m$, we have 
$$
\bH\bb= \bA{\bA^+\bb} \triangleq \bA\bx^+, 
\quad \text{let } \bx^+\triangleq \bA^+\bb,
$$
which is a linear combination of the columns of $\bA$. Thus, $\cspace(\bH)\subseteq \cspace(\bA)$.

Moreover, since $\bH$ is an orthogonal projection, by Lemma~\ref{lemma:rank-of-symmetric-idempotent2}, we have $\rank(\bH)=\trace(\bH)=\trace(\bA\bA^+) = \trace(\bU\bSigma\bV^\top \bV\bSigma^+\bU^\top)=\trace(\bU\bSigma\bSigma^+\bU^\top)=r$, where $\bU\bSigma\bV^\top$ is the SVD of $\bA$. That is, the rank of $\bH$ equals the rank of $\bA$. This implies $\cspace(\bH)= \cspace(\bA)$.

Similarly, we can prove $\bP= \bA^+\bA$ is the orthogonal projection onto the row space of $\bA$ using  the pseudo-inverse from the SVD.
\end{mdframed}

We finally provide another proof of  the four subspaces in the pseudo-inverse $\bA^+$ using its SVD. First, we need to show the following lemma that $\bA^+ \bA^{+\top}$ and $\bA^+$ have the same rank.
\begin{lemma}[Rank of $\bA^+ \bA^{+\top}$]\label{lemma:apapt-rank}
Let  $\bA^+$ be the pseudo-inverse of the matrix $\bA$. Then, $\bA^+ \bA^{+\top}$ and $\bA^+$ have same rank.
\end{lemma}
The proof follows directly from Lemma~\ref{lemma:rank-of-ata}.
%
Using this lemma, we can provide an alternative way to show the four subspaces in pseudo-inverses.
\begin{mdframed}[hidealllines=\mdframehidelineNote,backgroundcolor=\mdframecolor,frametitle={Another Way to See the Subspaces in Pseudo-Inverse via the SVD}]
Let the SVD of $\bA$ be $\bA = \bU\bSigma\bV^\top$ ($\bA$ has rank $r$), and let its pseudo-inverse  be $\bA^+ = \bV\bSigma^+\bU^\top$.
From Proposition~\ref{proposition:rank-of-symmetric}, consider the symmetric matrix and its spectral decomposition:  $\bA^+ \bA^{+\top} = \bV(\bSigma^\top\bSigma)^{-1}\bV^\top$ (this is the spectral decomposition of $\bA^+ \bA^{+\top} $, and $\{\bv_1,\bv_2 \ldots, \bv_r\}$ are the eigenvectors of $\bA^+ \bA^{+\top} $), $\cspace(\bA^+ \bA^{+\top})$ is spanned by the eigenvectors.
Therefore, $\{\bv_1,\bv_2 \ldots, \bv_r\}$ forms an orthonormal basis for $\cspace(\bA^+ \bA^{+\top})$.

We can easily derive the following results:
\begin{itemize}
\item $\bA^+ \bA^{+\top}$ is symmetric, so the row space of $\bA^+ \bA^{+\top}$ equals the column space of $\bA^+ \bA^{+\top}$. 

\item All columns of $\bA^+ \bA^{+\top}$ are linear combinations of the columns of $\bA^+$, so the column space of $\bA^+ \bA^{+\top}$ $\subseteq$ the column space of $\bA^+$, i.e., $\cspace(\bA^+ \bA^{+\top}) \subseteq \cspace(\bA^+)$. 

\item $\rank(\bA^+ \bA^{+\top}) = \rank(\bA^+)$ by Lemma~\ref{lemma:apapt-rank}.
\end{itemize}
Therefore, the row space of $\bA^+ \bA^{+\top}$ = the column space of $\bA^+ \bA^{+\top}$ =  the column space of $\bA^+$, i.e., $\cspace(\bA^+ \bA^{+\top}) = \cspace(\bA^+ )$. \textcolor{mylightbluetext}{Thus, $\{\bv_1, \bv_2,\ldots, \bv_r\}$ is an orthonormal basis of $\cspace(\bA^+)$}. We also proved in Proposition~\ref{proposition:svd-four-orthonormal-Basis} that $\{\bv_1, \bv_2,\ldots, \bv_r\}$ is an orthonormal basis of the row space of $\bA$ (i.e., basis of $\cspace(\bA^\top)$).
Therefore, $\cspace(\bA^+)=\cspace(\bA^\top)$, as shown in Figure~\ref{fig:lafundamental5-pseudo}.

Similarly, if we apply this process to $ \bA^{+\top}\bA^+$, we can show that the row space of $\bA^+$ is equal to the column space of $\bA$, and the null space of $\bA^+$ is equal to the null space of $\bA^\top$.
%
\end{mdframed}
We provide additional  descriptions for Figure~\ref{fig:lafundamental5-pseudo}. 
\paragraph{Comment 1). $\bA^+\bp =\bx^+$ with $\bp\in \cspace(\bA)$ and $\bx^+\in \cspace(\bA^\top)$.}
Furthermore, for $\bx^+$ in row space of $\bA$, we can express $\bx^+$ as  $\bx^+ \triangleq \sum_{i=1}^{r} x_i \bv_i$ since $\{\bv_1, \bv_2, \ldots, \bv_r\}$ is an orthonormal basis for the row space of $\bA$ (Proposition~\ref{proposition:svd-four-orthonormal-Basis}). For the vector $\bp = \bA\bx^+$ in the column space of $\bA$, we have $\bp = \bU\bSigma\bV^\top \bx^+$ and
$$
\begin{aligned}
\bA^+\bp &= \bV\bSigma^+\bU^\top \bU\bSigma\bV^\top \bx^+ = (\bV\bSigma^+\bSigma\bV^\top) (\bx^+)\\ 
&=\left(\sum_{i=1}^{r} \bv_i \bv_i^\top\right) \left(\sum_{i=1}^{r} x_i \bv_i\right) = \sum_{i=1}^{r} x_i \bv_i = \bx^+.
\end{aligned}
$$

\paragraph{Comment 2). $\bA^+\bb = \bx^+$ with $\bb\in \real^m$ and $\bx^+\in \cspace(\bA^\top)$.}
For any vector $\be$ in the null space of $\bA^\top$, we have $\bA^\top\be=\bA^+\be=\bzero$ since $\nspace(\bA^+)=\nspace(\bA^\top)$. Any  vector $\bb \in \real^m$ can be decomposed into $\bb = \bp+\be$, where $\bp$ is a vector lies in the column space of $\bA$ and $\be$ is a vector resides in the null space of $\bA^\top$. That is 
$$
\bA^+ \bp = \bA^+\bb = \bx^+,
$$
where $\bx^+$ is in the row space of $\bA$.

To conclude, for any vector $\bx^+$ in row space of $\bA$, there exists a vector $\bp$ in the column space of $\bA$ such that 
$$
\bp = \bA\bx^+ \leadto \bA^+\bp=\bx^+,
$$
and the relationship is depicted in Figure~\ref{fig:lafundamental5-pseudo}.

\paragraph{Linear systems in pseudo-inverses.}
In Section~\ref{section:application-ls-svd}, we introduced how the SVD can be applied to solve the least squares problem with rank-deficient matrices, i.e., solving the overdetermined linear system $\bA\bx=\bb$ where $\bA\in\real^{m\times n}$, $m\geq n$, and rank $r\leq n$. 
From Theorem~\ref{theorem:orthogonal-from-pseudo-inverse}, we can further conclude the following results related to the linear system problem.
\begin{proposition}[Target Vector of LS in Column Space]\label{proposition:targv_in_col}
Let $\bA\in\real^{m\times n}$ be any rank-$r$ matrix with SVD $\bA=\bU\bSigma\bV^\top$.
Then, the subspace $\cspace(\bA)=\cspace(\bA\bA^+)$ consists of all \textit{target vector} $\bb\in\real^m$  (satisfying $\bA\bx=\bb$) such that 
$$
\bU^\top \bb = 
\begin{bmatrix}
\bz \\
\bzero 
\end{bmatrix}, 
\gap 
\text{with }\bz\in\real^r.
$$
\end{proposition}
\begin{proof}[of Proposition~\ref{proposition:targv_in_col}]
Consider the linear system $\bA\bx=\bb$, we have 
$$
\bU^\top\bb
= \bU^\top\bA\bx 
= \bU^\top \bU\bSigma\bV^\top\bx
= \bSigma\bV^\top\bx 
= \begin{bmatrix}
\bSigma_r & \bzero \\
\bzero & \bzero_{m-r} 
\end{bmatrix}
\bV^\top \bx
=
\begin{bmatrix}
\bz \\
\bzero 
\end{bmatrix},
$$
where $\bSigma_r =\diag(\sigma_1, \sigma_2, \ldots, \sigma_r) \in\real^{r\times r}$ and $\bz\in\real^r$. 
Conversely,  if 
$\bU^\top \bb =
\footnotesize
\begin{bmatrix}
\bz \\
\bzero 
\end{bmatrix}$, then 
$
\bb
=\bU
\footnotesize
\begin{bmatrix}
\bz \\
\bzero 
\end{bmatrix}
$. And we have 
$$
\begin{aligned}
\bA\bA^+ \bb 
&= 
\bU\bSigma\bV^\top  \bV\bSigma^+\bU^\top \bb 
= \bU 
\begin{bmatrix}
\bI_r & \bzero \\
\bzero & \bzero_{m-r} 
\end{bmatrix}
\begin{bmatrix}
\bz \\
\bzero 
\end{bmatrix}
=\bb.
\end{aligned}
$$
That is, $\bb$ lies in the column space of $\bA$ and $\bA\bA^+$.
\end{proof}

Similarly, we have the following result.
\begin{proposition}[Target Vector of LS in Row Space]\label{proposition:targv_in_row}
Let $\bA\in\real^{m\times n}$ be any rank-$r$ matrix with SVD $\bA=\bU\bSigma\bV^\top$.
Then, the subspace $\cspace(\bA^\top)=\cspace(\bA^+\bA)$ (i.e., the row space of $\bA$) consists of all \textit{target vector} $\bb\in\real^n$ (satisfying $\bA^\top\bx=\bb$) such that 
$$
\bV^\top \bb = 
\begin{bmatrix}
\bz \\
\bzero 
\end{bmatrix}, 
\gap 
\text{with }\bz\in\real^r.
$$
\end{proposition}
\begin{proof}[of Proposition~\ref{proposition:targv_in_row}]
Consider the linear system $\bb = \bA^+\bA \bx$ (or $\bA^\top\bx=\bb$), we have 
$$
\bb = \bA^+\bA\bx = \bV\bSigma^+\bU^\top \bU\bSigma\bV^\top \bx 
=
\bV
\begin{bmatrix}
\bI_r & \bzero \\
\bzero & \bzero_{m-r} 
\end{bmatrix}
\bV^\top 
\bx
= 
\bV
\begin{bmatrix}
\bz \\
\bzero 
\end{bmatrix},
$$
where $\bz\in\real^r$. 
Conversely, if $\bV^\top \bb =
\footnotesize
\begin{bmatrix}
\bz \\
\bzero 
\end{bmatrix}$, then 
$ \bb =\bV
\footnotesize
\begin{bmatrix}
\bz \\
\bzero 
\end{bmatrix}$. 
And we have 
$$
\bA^+\bA\bb
=
\bV\bSigma^+\bU^\top \bU\bSigma\bV^\top\bV
\begin{bmatrix}
\bz \\
\bzero 
\end{bmatrix}
=
\bV
\begin{bmatrix}
\bI_r & \bzero \\
\bzero & \bzero_{m-r} 
\end{bmatrix}
\begin{bmatrix}
\bz \\
\bzero 
\end{bmatrix}
=
\bV\begin{bmatrix}
\bz \\
\bzero 
\end{bmatrix}=\bb.
$$
That is, $\bb$ lies in the subspace $\cspace(\bA^\top)=\cspace(\bA^+\bA)$.
\end{proof}

\subsection{Pseudo-Inverse in Spectral Theorem for Normal Matrices}\label{appendix:pseudo_inv_normal}
As a recap, in Theorem~\ref{theorem:normal_real_spectral_theorem}, we stated that a real normal matrix $\bA\in\real^{n\times n}$ (with $\bA^\top\bA=\bA\bA^\top$) admits a spectral decomposition 
$
\bA = \bQ\bD\bQ^\top,
$
where $\bQ$ is an $n\times n$ orthogonal matrix, and $\bD$ is a block-diagonal matrix $\bD=\diag(\bA_1, \bA_2, $ $\ldots, \bA_p)$ with $p\leq n$ consisting of one-dimensional blocks $\bA_i = \lambda_i$, or of two-dimensional blocks of the form 
$$
\bA_i =
\begin{bmatrix}
\lambda_i & \mu_i \\
-\mu_i & \lambda_i
\end{bmatrix},
\gap 
\text{with } \mu_i\neq 0.
$$
The pseudo-inverses of  normal matrices also have the spectral decomposition.
\begin{theoremHigh}[Spectral Decomposition in Pseudo-Inverse]\label{theorem:norm_pseudo_spec}
Let $\bA\in\real^{n\times n}$ be  a real normal matrix  that admits the  block-diagonalization $\bA=\bQ\bD\bQ^\top$ (discussed above). Then the pseudo-inverse of $\bA$ can be factored as 
$$
\bA^+ = \bQ\bD^+\bQ^\top,
$$
where $\bD^+$ is the pseudo-inverse of $\bD$.
Specifically, if 
$$
\bD = \begin{bmatrix}
\bD_r & \bzero \\
\bzero & \bzero 
\end{bmatrix},
\gap 
\text{with }\bD_r\in\real^{r\times r},
\qquad\underrightarrow{ \text{then} }\qquad
\bD^+ = \begin{bmatrix}
\bD_r^{-1} &\bzero \\
\bzero & \bzero 
\end{bmatrix}.
$$
\end{theoremHigh}
\begin{proof}[of Theorem~\ref{theorem:norm_pseudo_spec}]
Without loss of generality,
we assume $\bA_1, \bA_2, \ldots,  \bA_k$ are two-dimensional blocks and $\lambda_{2k+1}, \lambda_{2k+2}, \ldots, \lambda_n$ are the scalar entries. 
From the proof of Theorem~\ref{theorem:normal_real_spectral_theorem}, we know that the $\lambda_i\pm i \mu_i$ ($i\in \{1,2,\ldots, k\}$) are complex eigenvalues of $\bA$, and $\lambda_{2k+1}, \lambda_{2k+1}, \ldots, \lambda_n$ are real eigenvalues of $\bA$.
Let 
$$
\begin{aligned}
\rho_{2i-1}&\triangleq\rho_{2i} \triangleq \sqrt{\lambda_i^2+\mu_i^2} = \sqrt{\det(\bA_i)},\gap 
&\forall& i\in\{1,2,\ldots, k\}; \\
\rho_{i} &\triangleq |\lambda_i|, \gap &\forall& i \in\{2k+1, 2k+2, \ldots, \textcolor{mylightbluetext}{r}\},
\end{aligned}
$$
where we shall shortly see that the scalar $r$ is the rank of the matrix $\bA$. For $r<i\leq n$, $\lambda_i=0$.
Without loss of generality, we assume $\rho_1\geq \rho_2 \geq \ldots \geq \rho_r >0$ (this can also be achieved via a permutation matrix). Then we have 
$$
\bA\bA^\top = \bA^\top\bA = \bQ\bD \bQ^\top \bQ\bD^\top\bQ^\top = \bQ\bD\bD^\top\bQ^\top,
$$
with 
$$
\bD\bD^\top = \diag\left(\rho_1^2,\rho_2^2, \ldots, \rho_r^2, 0, \ldots, 0   \right)\in\real^{n\times n}.
$$
Therefore, $\bQ\bD\bD^\top\bQ^\top$ is the spectral decomposition of the matrix $\bA\bA^\top$, and $\{\rho_1^2, \ldots, \rho_r^2\}$ are the eigenvalues of $\bA\bA^\top$. Thus, let $\sigma_i=\rho_i$ for all $i\in\{1,2,\ldots, r\}$, $\{\sigma_1, \sigma_2, \ldots, \sigma_r\}$ are the (nonzero) singular values of $\bA$ (see the proof of Theorem~\ref{theorem:reduced_svd_rectangular}). This confirms that $r$ is the rank of the matrix $\bA$. We further define the block-diagonal matrix $\bZ$ (in the meantime, an orthogonal matrix):
$$
\bZ = 
\begin{bmatrix}
\frac{1}{\sigma_1} \bA_1 & &&&&& \\
&\ddots &&&&& \\
&& \frac{1}{\sigma_{2k}}\bA_k&&&& \\
&&&  \frac{\lambda_{2k+1}}{|\lambda_{2k+1}|}&&& \\
&&&&  \ddots && \\
&&&&& \frac{\lambda_r}{|\lambda_r|}& \\
&&&&&& \bI_{n-r} 
\end{bmatrix}\in \real^{n\times n}, 
\gap 
\text{with }\bZ^\top\bZ = \bZ\bZ^\top = \bI_n.
$$
Define the diagonal matrix $\bSigma \triangleq \diag(\sigma_1, \sigma_2, \ldots, \sigma_r, 0, \ldots, 0)\in\real^{n\times n}$.
We have 
$$
\bD = \bZ\bSigma.
$$
Let $ \bU \triangleq \bQ\bZ $ and $ \bV \triangleq \bQ$, it then follows that
$$
\bA = \bQ\bD\bQ^\top = \bQ\bZ\bSigma\bQ^\top =\bU\bSigma\bV^\top.
$$
Thus, $\bA$ admits the SVD $\bA=\bU\bSigma\bV^\top$.
The pseudo-inverse of $\bA$ is (from Table~\ref{table:pseudo-inverse-svd})
$$
\bA^+ = \bV\bSigma^+\bU^\top = \bV\bSigma^+\bZ^\top\bQ^\top,
$$
where $\bSigma^+\bZ^\top = \bSigma^+\bZ^{-1} = \bD^+$. This implies the spectral decomposition for normal matrices
$$
\bA = \bQ\bD\bQ^\top 
\quad\implies\quad
\bA^+ = \bQ \bD^+\bQ^\top.
$$
Specifically, if 
$
\bD = 
\footnotesize
\begin{bmatrix}
\bD_r & \bzero \\
\bzero & \bzero 
\end{bmatrix},
$
with
$ 
\bD_r\in\real^{r\times r},
$
then $\bD_r$ is invertible. Therefore,
$
\bD^+ = 
\footnotesize
\begin{bmatrix}
\bD_r^{-1} &\bzero \\
\bzero & \bzero 
\end{bmatrix}.
$
This completes the proof.
\end{proof}

\subsection{Pseudo-Inverse in CR Decomposition and Skeleton Decomposition}\label{section:pseudo_cr}
Suppose the CR decomposition (Theorem~\ref{theorem:cr-decomposition}) of a matrix $\bA \in \real^{m\times n}$ is given by $\bA=\bC\bR$, where $\bC\in\real^{m\times r}$ and $\bR\in \real^{r\times  n}$ has full ranks. And since in the CR decomposition, the row number of $\bC$  is  greater than or equal to the column number, and the row number of $\bR$ is less than or equal to the column number: $r\leq \min\{m,n\}$. Then we have the pseudo-inverse of $\bC$ (left-pseudo-inverse) and $\bR$ (right-pseudo-inverse):

$$
\bC^+ = (\bC^\top\bC)^{-1}\bC^\top \qquad \text{and} \qquad \bR^+ = \bR^\top (\bR\bR^\top)^{-1}.
$$
Consider the following two cases:
\paragraph{Case $m>n=r$.} We can obtain the left-pseudo-inverse by 
$$
\begin{aligned}
\text{left-pseudo-inverse} &= \bA^+ = (\bA^\top\bA)^{-1}\bA^\top 
= (\bR^\top\bC^\top\bC\bR)^{-1}\bR^\top\bC^\top.
\end{aligned}
$$
We then have 
$
(\bR^\top\bC^\top\bC\bR) \bR^+\bC^+ = (\bR^\top\bC^\top\bC\bR)\bR^\top (\bR\bR^\top)^{-1}(\bC^\top\bC)^{-1}\bC^\top = \bR^\top\bC^\top$. 
Thus, 
$
\bA^+ = \bR^+\bC^+.
$
Moreover, since $m>n$, we have $\bA^+\bA = \bI_n$, and
\begin{equation}\label{equation:cr-right-pseudo-inverse-appendix}
\begin{aligned}
\bA^+\bA &=\bR^+\bC^+ \bC\bR 
=\bR^\top (\bR\bR^\top)^{-1}(\bC^\top\bC)^{-1}\bC^\top \bC\bR \\
&= \bR^\top (\bR\bR^\top)^{-1}\bR 
=\bR^+\bR = \bI_n.
\end{aligned}
\end{equation}
This might seem confusing since we claimed in Equation~\eqref{equation:right-inverse-ba-appendix} that $\bR^+\bR \neq \bI_n$. However, this is not true. We assume $\bA$ has full rank  and $m>n$; this implies $\bR$ is a square invertible matrix. Therefore, $\bR^+=\bR^{-1}$ and the result in Equation~\eqref{equation:cr-right-pseudo-inverse-appendix}  does not apply to the left-pseudo-inverses of other matrices.

\paragraph{Case $n>m=r$.} We can obtain the right-pseudo-inverse by 
$$
\begin{aligned}
\text{right-pseudo-inverse} &= \bA^+ = \bA^\top(\bA\bA^\top)^{-1}
= \bR^\top\bC^\top(\bC\bR\bR^\top\bC^\top)^{-1}.
\end{aligned}
$$
Again, we have $\bR^+\bC^+(\bC\bR\bR^\top\bC^\top) = \bR^\top (\bR\bR^\top)^{-1}(\bC^\top\bC)^{-1}\bC^\top\bC\bR\bR^\top\bC^\top=\bR^\top\bC^\top$. Thus,
$
\bA^+ = \bR^+\bC^+.
$

From the two cases above, we can then conclude that for any full-rank matrix $\bA$ (either $m>n$ or $n\geq m$; see the proof of Theorem~\ref{theorem:existence-of-pseudo-inverse} for more general cases), the pseudo-inverse of $\bA$ is 
$$
\bA^+ = \bR^+\bC^+.
$$

The pseudo-inverse of $\bA$ for skeleton decomposition and interpolative decomposition is defined similarly. For the skeleton decomposition of $\bA=\bC\bU^{-1}\bR$,  the pseudo-inverse of $\bA$ is $\bA^+ = \bR^+\bU\bC^+$, where 
$
\bC^+ = (\bC^\top\bC)^{-1}\bC^\top $ and $\bR^+ = \bR^\top (\bR\bR^\top)^{-1}.
$

%% file: chapter-app_condition.tex
\newpage
\index{Perturbation}
\index{Condition number}
\section{Condition Number}\label{appendix:condition_number}
Let $\sX, \sY$ be vector spaces, and let $f:\sX\rightarrow \sY$ be a function. 
For any input $\bx\in \sX\mybackslash  \{\bzero\}$ and a perturbation $\widetildebx\in \sX$ of $\bx$, we define the \textit{absolute error} and \textit{relative error} of $\widetildebx$ with respect to $\bx$ and the norm $\norm{\cdot}$ by, respectively,
$$
\abserror (\widetildebx, \bx) = \norm{\widetildebx - \bx}
\quad\text{and}\quad
\relerror (\widetildebx, \bx) =\frac{\norm{\widetildebx - \bx}}{\norm{\bx}}.
$$
When performing computations, we are typically more interested in relative errors than absolute errors because the magnitude of the error by itself does not provide much useful information. The concept of the condition number is closely tied to the relative erro \citep{burgisser2013condition}.
\begin{definition}[(Strict) Condition Number]\label{definition:strict_condition_number}

The \textbf{condition number} of an input $\bx\in \sX$ (with respect to problem $f$) is, roughly speaking, the worst possible amplification  of the output error with respect to a small input
perturbation (i.e., small $\relerror\left(\widetildebx, \bx\right)$). More formally, the condition is defined as 
$$
\cond^f(\bx) 
= \mathop{\lim}_{\delta\rightarrow 0}
\mathop{\sup}_{\relerror(\widetildebx, \bx)\leq \delta}
\frac{\relerror\left(f(\widetildebx), f(\bx)\right) }{\relerror\left(\widetildebx, \bx\right) }.
$$
\end{definition}

\subsection{Condition Number for Matrix Inversion Problem}
When working with matrix problems,
for the induced matrix norm $\norm{\bA}_{a,b}$ from vector norms $\norm{\cdot}_a$ and $\norm{\cdot}_b$, and a matrix $\bA\in \real^{n\times n}$, define 
$$
\kappa_{ab}(\bA) \triangleq \norm{\bA}_{a,b} \norm{\bA^{-1}}_{b,a}.
$$
When $a=b=2$, we obtain the classic condition number with respect to the spectral norm (Equation~\eqref{equation:qr_condition_num}). More generally, we can prove the condition number with respect to the general induced norms in the following theorem.

\begin{theorem}[Condition Number of Matrix Inversion with Induced Norm]\label{theorem:condition_induced}
Given a nonsingular matrix $\bA\in \real^{n\times n}$ and the function to compute the inverse $g(\bA)=\bA^{-1}$. 
We measure the relative error on the data space and solution space with respect to $\norm{\cdot}_{a,b}$ and $\norm{\cdot}_{b,a}$, respectively. Then,
 we have $\cond^g(\bA) = \kappa_{ab}(\bA)$.
\end{theorem}

\begin{proof}[of Theorem~\ref{theorem:condition_induced}]
Suppose $\widetildebA$ is a perturbation of $\bA$ such that $\widetildebA = \bA-\bE$, where $\bE\in \real^{n\times n}$ is the error matrix with $\Vert\bE\Vert_{a,b}\rightarrow 0$.
The relative error of the output (with respective to the induced norm $\norm{\cdot}_{a,b}$) is 
$$
\relerror\left(g(\widetildebA), g(\bA)\right) 
= 
\frac{\norm{(\bA-\bE)^{-1} - \bA^{-1}}_{b,a}}
{\norm{\bA^{-1}}_{b,a}}
=\frac{\norm{\bA^{-1}(\bI-\bE\bA^{-1})^{-1} - \bA^{-1}}_{b,a}}
{\norm{\bA^{-1}}_{b,a}}
$$
Since we are considering the perturbation $\Vert\bE\Vert_{a,b}\rightarrow 0$, it follows that the derivative of the function $(\bI-\bE\bA^{-1})^{-1}$ has a Taylor series expansion (Appendix~\ref{appendix:taylor-expansion}):
$$
(\bI-\bE\bA^{-1})^{-1} = \bI + (\bE\bA^{-1})+ (\bE\bA^{-1})^2 + (\bE\bA^{-1})^3+\ldots 
= \sum_{i=0}^{\infty}(\bE\bA^{-1})^i.
$$
Therefore, we have
\begin{equation}\label{equation:conditon_inverse_relerrlreq}
\begin{aligned}
\relerror&\left(g(\widetildebA), g(\bA)\right)
=\frac{\norm{\bA^{-1}\sum_{i=0}^{\infty}(\bE\bA^{-1})^i- \bA^{-1}}_{b,a}}
{\Vert\bA^{-1}\Vert_{b,a}}\\
&= \frac{\norm{\bA^{-1}\bE\bA^{-1}+\bA^{-1}\sum_{i=2}^{\infty}(\bE\bA^{-1})^i}_{b,a}}
{\Vert\bA^{-1}\Vert_{b,a}} \\
&\leq 
\frac{\norm{\bA^{-1}\bE\bA^{-1}}_{b,a}+
	\norm{\bA^{-1}\sum_{i=2}^{\infty}(\bE\bA^{-1})^i}_{b,a}}
{\Vert\bA^{-1}\Vert_{b,a}} \\
&\leq
\frac{\norm{\bA^{-1}\bE\bA^{-1}}_{b,a}+
\norm{\bA^{-1}}_{b,a}\norm{\sum_{i=2}^{\infty}(\bE\bA^{-1})^i}_{b,b}}
{\Vert\bA^{-1}\Vert_{b,a}}\\
&\leq
\frac{\norm{\bA^{-1}\bE\bA^{-1}}_{b,a}+
\norm{\bA^{-1}}_{b,a} \sum_{i=2}^{\infty} \norm{\bE}_{a,b}^i \norm{\bA^{-1}}_{b,a}^i}
{\Vert\bA^{-1}\Vert_{b,a}},
\end{aligned}
\end{equation}
from the triangle inequality and 
submultiplicative property of matrix norms (Theorem~\ref{theorem:Submultiplicativity_induced}).
Hence,  the condition number can be derived,
$$
\begin{aligned}
\cond^g(\bA) 
&= \mathop{\lim}_{\delta\rightarrow 0}
\mathop{\sup}_{\relerror(\widetildebA, \bA)\leq \delta}
= \frac{\relerror\left(g(\widetildebA), g(\bA)\right) }{\relerror\left(\widetildebA, \bA\right) }\\
&=\mathop{\lim}_{\delta\rightarrow 0}
\mathop{\sup}_{\frac{\Vert\bE\Vert_{a,b}}{\norm{\bA}_{a,b}}\leq \delta}
\frac{\Vert\bA^{-1}\bE\bA^{-1}+\bA^{-1}\sum_{i=2}^{\infty}(\bE\bA^{-1})^i\Vert_{b,a}}
{\Vert\bA^{-1}\Vert_{b,a}}
\bigg/
\frac{\Vert\bE\Vert_{a,b}}{\norm{\bA}_{a,b}} \\
&\stackrel{*}{\leq}
 \mathop{\lim}_{\delta\rightarrow 0}
\mathop{\sup}_{\frac{\Vert\bE\Vert_{a,b}}{\norm{\bA}_{a,b}}\leq \delta}
\frac{\norm{\bA^{-1}}_{b,a}^2  \Vert\bE\Vert_{a,b}
	+
\norm{\bA^{-1}}_{b,a} \sum_{i=2}^{\infty} \Vert\bE\Vert_{a,b}^i \norm{\bA^{-1}}_{b,a}^i
 }{\norm{\bA^{-1}}_{b,a}} 
\cdot 
\frac{\norm{\bA}_{a,b}}{\Vert\bE\Vert_{a,b}}\\
&=
\norm{\bA}_{a,b} \cdot 
\mathop{\lim}_{\delta\rightarrow 0}
\mathop{\sup}_{\frac{\Vert\bE\Vert_{a,b}}{\norm{\bA}_{a,b}}\leq \delta}
\left( \norm{\bA^{-1}}_{b,a} 
+ \norm{\bA^{-1}}_{b,a}\sum_{i=2}^{\infty} 
\left(\norm{\bE}_{a,b} \cdot  \norm{\bA^{-1}}_{b,a}\right)^{i-1}
\right) \\
&=\norm{\bA}_{a,b} \norm{\bA^{-1}}_{b,a} 
= \kappa_{ab}(\bA),
\end{aligned}
$$
where the inequality ($*$) is again from the triangle inequality and submultiplicative property of matrix norms  (Theorem~\ref{theorem:Submultiplicativity_induced}).
Therefore, we prove $\cond^g(\bA) \leq \kappa_{ab}(\bA)$.

To show the reverse inequality (i.e., $\kappa_{ab}(\bA)\leq \cond^g(\bA)$),  given $\by\in\real^n$ such that $\norm{\bA^{-1}\by}_a = \norm{\bA^{-1}}_{b,a}$ and $\Vert\by\Vert_{b}=1$.
Define $\bx = \frac{\bA^{-1}\by}{\norm{\bA^{-1}\by}_a}$ such that $\norm{\bx}_a=1$.
By Lemma~\ref{lemma:induced_norm_property},
there exists a matrix $\bB\in \real^{n\times n}$ such that $\norm{\bB}_{a,b}=1$ and $\bB\bx=\by$.
Since we consider $\Vert\bE\Vert_{a,b}\rightarrow 0$, we can define $\bE=\delta\bB$ for small $\delta>0$ (with $\Vert\bE\Vert_{a,b} = \delta$).
Then we have 
$$
\begin{aligned}
\kappa_{ab}(\bA) &= \norm{\bA}_{a,b} \norm{\bA^{-1}}_{b,a}
= \norm{\bA}_{a,b} \norm{\bA^{-1}\by}_{a}
= \norm{\bA}_{a,b} \norm{\bA^{-1}\bB\bx}_{a}\\
&=  \norm{\bA}_{a,b}  \left\Vert \bA^{-1}\bB \frac{\bA^{-1}\by}{\Vert\bA^{-1}\by\Vert_a} \right\Vert_a
=\norm{\bA}_{a,b}   \frac{\norm{\bA^{-1}\bB\bA^{-1}\by}_a}{\norm{\bA^{-1}}_{b,a}} \\
&\stackrel{*}{\leq }
\norm{\bA}_{a,b}
\frac{\norm{\bA^{-1}\bB\bA^{-1}}_{b,a} \Vert\by\Vert_b}
{\norm{\bA^{-1}}_{b,a}} 
=\norm{\bA}_{a,b}
\frac{\norm{\bA^{-1}\bE\bA^{-1}}_{b,a}}
{\delta\cdot  \norm{\bA^{-1}}_{b,a}},
\end{aligned}
$$
where the inequality ($*$) follows from Equation~\eqref{equation:induced_ineqy_intern}.
On the other hand, following from Equation~\eqref{equation:conditon_inverse_relerrlreq}, we note that 
$$
\mathop{\lim}_{\Vert\bE\Vert_{a,b}\rightarrow 0}\norm{\bA^{-1}}_{b,a} \sum_{i=2}^{\infty} \Vert\bE\Vert_{a,b}^i \norm{\bA^{-1}}_{b,a}^i=
o(\Vert\bE\Vert_{a,b})
$$
and
$$
\begin{aligned}
\relerror\left(g(\widetildebA), g(\bA)\right)
&=\frac{\norm{\bA^{-1}\bE\bA^{-1}}_{b,a}
}
{\norm{\bA^{-1}}_{b,a}}
+o(\Vert\bE\Vert_{a,b}).
\end{aligned}
$$
This implies 
$$
\begin{aligned}
\kappa_{ab}(\bA) &=  
\frac{\norm{\bA}_{a,b}}{\delta }
\left(\relerror\left(g(\widetildebA), g(\bA)\right) - o(\Vert\bE\Vert_{a,b})\right)\\
&=\frac{\norm{\bA}_{a,b}}{\Vert\bE\Vert_{a,b} }
\left(\relerror\left(g(\widetildebA), g(\bA)\right) - o(\Vert\bE\Vert_{a,b})\right)\\
&\leq \frac{\relerror\left(g(\widetildebA), g(\bA)\right)}
{\relerror(\widetildebA, \bA)}
\leq 
\mathop{\sup}_{\Vert\bE\Vert_{a,b} \leq \delta }
\frac{\relerror\left(g(\widetildebA), g(\bA)\right)}
{\relerror(\widetildebA, \bA)}\\
&=\mathop{\sup}_{\relerror(\widetildebA, \bA) \leq \frac{\delta}{\norm{\bA}_{a,b}}}
\frac{\relerror\left(g(\widetildebA), g(\bA)\right)}
{\relerror(\widetildebA, \bA)}.
\end{aligned}
$$
Let $\delta\rightarrow 0$, it follows that $\kappa_{ab}(\bA)\leq \cond^g(\bA)$.
This completes the proof.
\end{proof}

\subsection{Condition Number for Linear System}
In last section, we prove that the condition number of the matrix inversion problem is equal to $\kappa_{ab}(\bA)=\norm{\bA}_{a,b} \norm{\bA^{-1}}_{b,a}$.
We consider further the condition number for linear system problem.
Note in this section, we consider general subordinate matrix norms $\norm{\cdot}:\real^{n\times n}\rightarrow \real$ (Definition~\ref{definition:subordinate_matrix_norm}).

Instead of considering the worst possible magnification of the output error with respect to a small input perturbation in Definition~\ref{definition:strict_condition_number}, in the linear system problem, we define the condition number as $\kappa(\bA)$ with respect to the subordinate matrix norm as follows. The general matrix condition number provides similar measure as the (strict) condition number in the context of linear systems. The larger the (strict or general) condition number, the more illness of the problem.
The reason will be clear in the sequel. 
\begin{definition}[General Matrix Condition Number]\label{definition:general_matrix_condition_num}
Given any nonsingular matrix $\bA$ and a subordinate matrix norm $\norm{\cdot}$, we define the general condition number 
$$
\cond(\bA)=\kappa(\bA) = \norm{\bA}\norm{\bA^{-1}}.
$$
\end{definition}
We then consider the condition number in linear systems with different perturbations.
\begin{theorem}[Perturb The Target]\label{theorem:perturb_targer_ols}
Let $\bA$ be a nonsingular matrix, and let $\bx, \bx+\delta\bx$ be the solutions of the following linear systems, respectively,
$$
\begin{aligned}
\bA\bx&=\bb
\quad \text{and}\quad 
\bA(\bx+\delta\bx) = \bb+\delta\bb.
\end{aligned}
$$ 
If $\bb\neq \bzero$, then the  following inequality holds:
$$
\frac{\Vert\delta \bx\Vert}{\norm{\bx}}
\leq \cond(\bA)
\frac{\Vert\delta\bb\Vert}{\norm{\bb}}.
$$

Note in this case, the condition number is the same as the (strict) condition number defined in Definition~\ref{definition:strict_condition_number}, 
i.e., it is the  worst possible magnification of the output error with respect to a small input perturbation. 
\end{theorem}
\begin{proof}[of Theorem~\ref{theorem:perturb_targer_ols}]
From the linear system equations, we get 
$
\delta\bx = \bA^{-1}\delta \bb.
$
Since we are considering subordinate matrix norms, it follows that 
$$
\begin{aligned}
\Vert\delta \bx\Vert &\leq \left\Vert\bA^{-1}\right\Vert \Vert\delta \bb\Vert
\quad\text{and}\quad
\Vert \bb \Vert &\leq \left\Vert \bA\right \Vert \norm{\bx}.
\end{aligned}
$$
Therefore, the relative error $\frac{\Vert\delta \bx\Vert}{\norm{\bx}}$ is bounded in terms of the relative error $\frac{\Vert\delta\bb\Vert}{\norm{\bb}}$:
$$
\frac{\Vert\delta \bx\Vert}{\norm{\bx}}
\leq 
\left\Vert\bA\right\Vert  \left\Vert\bA^{-1}\right\Vert 
\frac{\Vert\delta\bb\Vert}{\norm{\bb}}.
$$
This  completes the proof.
\end{proof}

\begin{theorem}[Perturb The Matrix]\label{theorem:perturb_in_matrix_ols}
Let $\bA$ be a nonsingular matrix, and let $\bx, \bx+\delta\bx$ be the solutions of the following linear systems, respectively,
$$
\begin{aligned}
\bA\bx&=\bb
\quad \text{and}\quad 
(\bA+\delta\bA)(\bx+\delta\bx) = \bb.
\end{aligned}
$$ 
If $\bb\neq \bzero$, then the  following inequality holds:
$$
\frac{\Vert\delta \bx\Vert}{\Vert\bx +\delta\bx\Vert}
\leq \cond(\bA)
\frac{\Vert\delta\bA\Vert}{\Vert\bA\Vert}.
$$

Note in this case, the condition number is \textbf{not} exactly the same as the (strict) condition number defined in Definition~\ref{definition:strict_condition_number}  since we modify the relative error from $\frac{\Vert\delta \bx\Vert}{\norm{\bx}}$ 
to 
the (modified) relative error $\frac{\Vert\delta \bx\Vert}{\Vert\bx+\delta\bx\Vert}$.
That's why we refer to the condition number in Definition~\ref{definition:general_matrix_condition_num} as the \textit{general matrix condition number}.

Furthermore, when considering  small $\Vert\delta\bA\Vert$ with $\Vert\delta\bA\Vert <\frac{1}{\norm{\bA^{-1}}} $, we have 
\begin{equation}\label{equation:perturb_final_matrix_cond}
\begin{aligned}
\frac{\Vert\delta\bx\Vert}{\norm{\bx}}
&\leq \cond(\bA)
\left(  
\frac{1}{1- \norm{\bA^{-1}}\cdot \Vert\delta\bA\Vert}
\right)
\frac{\Vert\delta\bA\Vert}{\Vert\bA\Vert}\\
&=
 \cond(\bA)
\left(  
1+O\left(\Vert\delta\bA\Vert\right)
\right)
\frac{\Vert\delta\bA\Vert}{\Vert\bA\Vert}
\end{aligned}
\end{equation}
\end{theorem}
\begin{proof}[of Theorem~\ref{theorem:perturb_in_matrix_ols}]
From the linear system equations, we obtain 
\begin{equation}\label{equation:ols_cond_perturb_A1}
\bA\bx +\bA\delta\bx + \delta\bA\bx + \delta\bA\delta\bx = \bb
\quad\implies\quad 
\delta\bx = -\bA^{-1} \delta\bA(\bx+\delta\bx).
\end{equation}
Since we are applying subordinate matrix norms, it follows that 
$$
\begin{aligned}
\Vert\delta \bx\Vert &\leq 
\left\Vert\bA^{-1}\right\Vert 
\left\Vert\delta\bA\right\Vert
\Vert\bx+\delta\bx \Vert.
\end{aligned}
$$
Therefore, the (modified) relative error $\frac{\Vert\delta \bx\Vert}{\Vert\bx+\delta\bx\Vert}$ is bounded in terms of the relative error $\frac{\Vert\delta\bA\Vert}{\Vert\bA\Vert}$:
$$
\frac{\Vert\delta \bx\Vert}{\Vert\bx +\delta\bx\Vert}
\leq 
\kappa(\bA)
\frac{\Vert\delta\bA\Vert}{\Vert\bA\Vert}.
$$
Since $\delta\bA < \frac{1}{\norm{\bA^{-1}}}$,  we have 
$$
\left\Vert\bA^{-1} \delta\bA\right\Vert 
\leq 
\left\Vert\bA^{-1}\right\Vert 
\Vert\delta\bA\Vert 
<1.
$$
by Proposition~\ref{proposition:subordinate_ineq1} (assume the norm is normalized such that $\norm{\bI}=1$ for simplicity) and the matrix $\bI+\bA^{-1}\delta\bA$ is invertible, we get
\begin{equation}\label{equation:ols_cond_perturb_A2}
\left\Vert (\bI+\bA^{-1}\delta\bA)^{-1} \right\Vert 
\leq 
\frac{1}{1- \Vert\bA^{-1}\delta\bA\Vert}
\leq
\frac{1}{1- \Vert\bA^{-1}\Vert\cdot \Vert\delta\bA\Vert}.
\end{equation}
By adding $\bx$ to Equation~\eqref{equation:ols_cond_perturb_A1}, we obtain 
$$
\begin{aligned}
&(\bI+\bA^{-1}\delta\bA)(\bx+\delta\bx) =\bx
\quad\implies\quad 
\bx+\delta\bx =(\bI+\bA^{-1}\delta\bA)^{-1}\bx\quad\implies \\
&\gap\gap
\delta\bx = \left((\bI+\bA^{-1}\delta\bA)^{-1}-\bI\right)\bx
=
(\bI+\bA^{-1}\delta\bA)^{-1}\left( \bI-(\bI+\bA^{-1}\delta\bA)  \right)\bx\\
&\gap\gap\gap =
-(\bI+\bA^{-1}\delta\bA)^{-1}\left( \bA^{-1}\delta\bA  \right)\bx.
\end{aligned}
$$
Again by the subordinate property, we get 
\begin{equation}\label{equation:ols_cond_perturb_A3}
\Vert\delta\bx\Vert
\leq 
\left\Vert(\bI+\bA^{-1}\delta\bA)^{-1}\right\Vert
\left\Vert  \bA^{-1}\right\Vert
\left\Vert  \delta\bA\right\Vert
\left\Vert    \bx\right\Vert.
\end{equation}
Combine Equation~\eqref{equation:ols_cond_perturb_A2} and \eqref{equation:ols_cond_perturb_A3}, we get 
$$
\frac{\Vert\delta\bx\Vert}{\norm{\bx}}
\leq \cond(\bA)
\left(  
\frac{1}{1- \Vert\bA^{-1}\Vert\cdot \Vert\delta\bA\Vert}
\right)
\frac{\Vert\delta\bA\Vert}{\Vert\bA\Vert}.
$$
This completes the proof.
\end{proof}

Same setting as Theorem~\ref{theorem:perturb_in_matrix_ols},
suppose we consider the function $f(\bA,\bb) = \bA^{-1}\bb$. Then the strict condition number in Definition~\ref{definition:strict_condition_number} for this problem is 
$$
\cond^f(\bA,\bb) = \cond(\bA)
\left(  
\frac{1}{1- \Vert\bA^{-1}\Vert\cdot \Vert\delta\bA\Vert}
\right).
$$
The strict condition number is proportional to the general condition number up to a constant. That's why we relax the definition of the strict condition number to the general condition in the context of linear systems. This modified definition can simplify the matters.

\begin{corollary}[Perturb Matrix and Target]
When the matrix $\bA$ and the target $\bb$ are perturbed simultaneously, we have the following linear system
$$
(\bA+\delta\bA)(\bx+\delta\bx) = \bb+\delta\bb.
$$
Then it can be shown when $\Vert\delta\bA\Vert <\frac{1}{\Vert\bA^{-1}\Vert}$, it follows that 
$$
\frac{\Vert\delta\bx\Vert}{\norm{\bx}}
\leq \cond(\bA)
\left(  
\frac{1}{1- \norm{\bA^{-1}}\cdot \Vert\delta\bA\Vert}
\right)
\left(
\frac{\Vert\delta\bA\Vert}{\Vert\bA\Vert}
+
\frac{\Vert\delta\bb\Vert}{\norm{\bb}}
\right).
$$
\end{corollary}

\subsection*{Properties of Condition Number}
To conclude this section, we provide some properties of the general matrix condition number.
\begin{proposition}[Condition Number Properties]\label{proposition:condi_properties}
Given a nonsingular matrix $\bA\in \real^{n\times n}$, and the condition number $\cond(\bA)=\norm{\bA} \norm{\bA^{-1}}$ with respect to a normalized subordinate norm $\Vert\cdot\Vert:\real^{n\times n}\rightarrow \real$, we have following properties:
\begin{enumerate}
\item A normalized matrix norm indicates $\norm{\bI}=1$. Then,
$$
\begin{aligned}
\cond(\bA) &\geq 1;\\
\cond(\bA) &=\cond(\bA^{-1});\\
\cond(\alpha\bA)&=\cond(\bA), \gap \text{with }\alpha\in \real\mybackslash\{0\}.
\end{aligned}
$$
\item The condition number $\cond_2(\bA)$ with respect to the spectral norm is 
$$
\cond_2(\bA) = \frac{\sigma_1}{\sigma_n},
$$
where $\sigma_1$ is the largest singular value, and $\sigma_n$ is the smallest singular value of $\bA$.
\item If $\bA$ is orthogonal, then $\cond_2(\bA)=1$.
\item The condition number $\cond_2(\bA)$ is invariant under orthogonal transformation, 
$$
\begin{aligned}
\cond_2(\bA) &= \cond_2(\bU^\top\bA\bV)= \cond_2(\bU\bA\bV^\top)=\\
\cond_2(\bU\bA) &= \cond_2(\bA\bV) =\cond_2(\bU\bA\bV),
\end{aligned}
$$
where $\bU$ and $\bV$ are orthogonal.
\end{enumerate}
\end{proposition}
\begin{proof}[of Proposition~\ref{proposition:condi_properties}]
For (1). Since the norm is normalized, we have $1=\Vert\bI\Vert \leq \Vert\bA\Vert\left\Vert\bA^{-1}\right\Vert=\cond(\bA)$ from the submultiplicativity property. $\cond(\bA) =\cond(\bA^{-1})$ follows immediately from symmetry. And $\cond(\alpha\bA)=\cond(\bA)$ follows from the definition of the condition number.

For (2). We know that $\Vert\bA\Vert_2=\sigma_1$, and by SVD (Chapter~\ref{chapter:SVD}), it can be easily shown that $\left\Vert\bA^{-1}\right\Vert = \frac{1}{\sigma_n}$. This proves the result.

To see (3), for an orthogonal matrix, we have  $\bA^\top\bA=\bI$. Its eigenpari shows $\bA^\top\bA\bx=\lambda\bx$. Therefore, $\lambda=1$. Thus, it follows that $\cond_2(\bA)=1$.

Part (4) follows from the invariance of the spectral norm under orthogonal transformations.
\end{proof}